%% file: CamposCiclotomicos.tex
\documentclass[envcountsame,envcountsect]{CIEA}

\usepackage{makeidx}
\usepackage{graphicx}
\usepackage{multicol}
\usepackage[spanish]{babel}
\usepackage{amsmath}
\usepackage{amssymb}
\usepackage{amscd}
\usepackage{latexsym}
\usepackage{mathrsfs}
\usepackage{picinpar}
\usepackage{epsfig}
\usepackage[all]{xy}
\usepackage{pifont}
\usepackage{mathtools}
\usepackage{rotating}
\usepackage{xcolor}

\def\ed{\end{document}}
\def\s{\bigskip}

\font\letritas=pcrro7t.tfm at 48pt

\renewcommand\abstractname{Resumen}
\renewcommand\chaptername{Cap{\'\i}tulo}

\renewcommand\contentsname{Contenido}
\renewcommand\appendixname{Ap\'{e}ndice}
\renewcommand\bibname{Bibliograf\'{\i}a} 
\renewcommand\indexname{\'Indice alfab\'etico}

\renewcommand\gcd{\operatorname{mcd}}

\newcommand{\carlitzbinom}{\genfrac[]{0pt}{0}}
\newcommand{\xbinom}{\genfrac(){.5pt}0}
\newcommand{\frobeniusbinom}{\genfrac[]{1pt}0}
\newcommand{\frobenius}{\genfrac[]{0.5pt}0}
\newcommand{\artin}{\genfrac(){.5pt}1}
\newcommand{\artinp}{\genfrac(){0.5pt}0}

\def\eu{\mathfrak}
\def\ma{\mathbb}
\def\mc{\mathcal}

\def\SE{\begin{list}
	{{\rm {(\roman{otrobean})}}}{\usecounter{otrobean}}
	\setlength{\rightmargin}{\leftmargin}}
	
\DeclareMathOperator{\mcd}{mcd}
\DeclareMathOperator{\mcm}{mcm}

\newcommand{\ucong}{\mathbin{\rotatebox[origin=c]{-90}{$\cong$}}}

\DeclareRobustCommand\longtwoheadrightarrow
     {\relbar\joinrel\twoheadrightarrow}

\newcommand{\Frobenius}{\operatorname{Fr}}
\newcommand{\sop}{\operatorname{sop}}
\newcommand{\B}{\operatorname{Br}}
\newcommand{\Ver}{\operatorname{Ver}}
\newcommand{\Drin}{\operatorname{Drin}}
\newcommand{\Red}{\operatorname{Red}}
\newcommand{\inmenosnoventa}{\mathbin{\rotatebox[origin=c]{-90}{$\in$}}}
\newcommand{\incientoochenta}{\mathbin{\rotatebox[origin=c]{180}{$\in$}}}
\newcommand{\Res}{\operatornamewithlimits{Res}}
\renewcommand{\ker}{\operatorname{n{\text{\rm{\'u}}}c}}
\renewcommand{\deg}{\operatorname{gr}}
\newcommand{\card}{\operatorname{card}}
\newcommand{\equivalente}{\operatornamewithlimits{\equiv}}
\newcommand{\igual}{\operatornamewithlimits{=}}
\newcommand{\equis}{\operatornamewithlimits{\hbox{$x$}}}
\newcommand{\des}{\operatornamewithlimits{\leq}}

\newcommand{\tr}{\operatorname{tr}}
\newcommand{\modulo}{\operatorname{mod}}
\newcommand{\Adj}{\operatorname{Adj}}
\newcommand{\im}{\operatorname{im}}
\newcommand{\cog}{\operatorname{cog}}
\newcommand{\Irr}{\operatorname{Irr}}
\newcommand{\MC}{\operatorname{MC}}
\newcommand{\isomo}{\operatornamewithlimits{\cong}}
\newcommand{\longto}{\operatorname{\longrightarrow}}
\newcommand{\Id}{\operatorname{Id}}
\newcommand{\tor}{\operatorname{tor}}
\newcommand{\rest}{\operatorname{rest}}
\newcommand{\id}{\operatorname{Id}}

\newcommand{\Gal}{\operatornamewithlimits{Gal}}
\newcommand{\biy}{\operatornamewithlimits{\rightleftarrows}}

\newcommand{\Tr}{\operatorname{Tr}}

\newcommand{\N}{\operatorname{N}}

\newcommand{\ran}{\operatorname{rango}}
\newcommand{\Aut}{\operatorname{Aut}}
\newcommand{\con}{\operatorname{con}}
\newcommand{\coc}{\operatorname{coc}}
\newcommand{\cotr}{\operatorname{cotr}}

\newcommand{\ab}{\operatorname{ab}}
\newcommand{\res}{\operatorname{res}}
\newcommand{\cores}{\operatorname{cor}}
\newcommand{\car}{\operatorname{car}}
\newcommand{\ord}{\operatorname{ord}}
\newcommand{\Hom}{\operatorname{Hom}}

\newcommand{\coker}{\operatorname{con{\text{\rm{\'u}}}cleo}}

\newcommand{\ex}{\operatorname{ex}}
\newcommand{\an}{\operatorname{an}}
\newcommand{\pic}{\operatorname{Pic}}
\newcommand{\sgn}{\operatorname{sgn}}
\newcommand{\End}{\operatorname{End}}
\newcommand{\partereal}{\operatorname{{\ma R}e}}

\newcommand{\sen}{\operatorname{sen}}

\newcommand{\gr}{\operatorname{gr}}

\newcommand{\GL}{\operatorname{GL}}
\newcommand{\codim}{\operatorname{codim}}

\newcommand{\masmenos}{\operatornamewithlimits{\pm}}

\newcommand{\integer}{\genfrac[]{0.5pt}0}
\newcommand{\integerchico}{\genfrac[]{0.5pt}1}
\newcommand{\mayor}{\genfrac\lceil\rceil{0.5pt}0}
\newcommand{\mayorchico}{\genfrac\lceil\rceil{0.5pt}1}
\newcommand{\hooklongrightarrow}{\lhook\joinrel\longrightarrow}
\newcommand{\uigual}{\mathbin{\rotatebox[origin=c]{90}{$=$}}}
\renewcommand{\o}{\mathcal O}
\newcommand{\ca}{\operatorname{car}}
\newcommand{\inv}{\operatorname{inv}}
\newcommand{\Spl}{\operatorname{Spl}}
\newcommand{\Isog}{\operatorname{Isog}}
\newcommand{\isomorfo}{\operatornamewithlimits{\cong}}
\newcommand{\infla}{\operatorname{inf}}

\newcommand{\lcm}{\operatorname{mcm}}
\newcommand{\Sg}{\operatorname{Sig}}
\newcommand{\sg}{\operatorname{sig}}
\newcommand{\exte}{\operatorname{ext}}
\newcommand{\cupdot}{\mathbin{\mathaccent\cdot\cup}}

\def\dps{\displaystyle}

\def\Br#1{\B(#1)}
\def\pKK{\pK\in{\ma P}_K}
\def\NN#1{{\mathcal N}_{#1}}

\def\abe#1{#1^{\rm{ab}}}
\def\norma{\N_{L_{\pL}/K_{\pK}}}
\def\norm{{\mc N}}
\def\unidades#1{U_K^{(#1)}}
\def\unidadess#1{U_\pK^{(#1_\pK)}}
\def\d{\displaystyle}
\def\e{\quad}
\def\normal{\triangleleft}

\def\Lra{\Longrightarrow}
\def\lra{\longrightarrow}

\def\iso{\simeq}
\def\pK{{\eu p}}
\def\pL{{\eu P}}
\def\p{{\eu p}_{\infty}}
\def\G#1{\big(R_T/\langle#1\rangle\big)^{\ast}}
\def\raizm#1{\sqrt[l]{(-1)^{\deg P_{#1}}P_{#1}}}

\def\raiz#1{\sqrt[l]{P_{#1}}}
\def\Witt#1{\stackrel{_{\bullet}}{#1}}
\def\vWitt#1#2{\big(#1^{(1)},\ldots,#1^{(#2)}\big)}
\def\S#1{S_{\infty}(#1)}
\def\F{{\ma F}_q}
\def\K{\mathscr K}
\def\LL{\mathscr L}
\def\nr#1{#1^{\rm{nr}}}
\def\mo#1{#1^{\rm{mod}}}
\def\pinfty{{\eu p}_{1,{\infty}}^{e_{1,{\infty}}}
\cdots {\eu p}_{r_{\infty},\infty}^{e_{r_{\infty},\infty}}}
\def\elemental#1#2#3{#1_{#3-1}#2^{p^{#3-1}}+#1_{#3-2}#2^{p^{#3-2}}+
\cdots +#1_2#2^{p^2} + #1_1#2^p+#1_0#2}
\def\elementalito#1#2#3{#1_{#3-1}#2^{p^{#3-1}}+
\cdots + #1_1#2^p+#1_0#2}
\def\va#1#2{\big({#1}_1,\ldots,{#1}_{#2}}

\def\L{L^{\ast}}
\def\cond#1{{\eu f}_{#1}}
\def\aditivos#1#2{#1\langle\tau_{#2}\rangle}
\def\aditivo#1{#1\langle\tau\rangle}
\def\serie#1#2#3{#1(\tau)=\sum_{i=0}^{#2}#3_i \tau^i}
\def\Ki{\K_{\infty}}
\def\Kii{K_{\infty}}
\def\vi{v_{\infty}}
\def\Ci{{\ma C}_{\infty}}
\def\Fi{{\ma F}_{\infty}}
\def\di{d_{\infty}}
\def\exponencial#1#2#3{#1\prod_{#2\in #3\setminus \{0\}}\big(1-\frac{#1}{#2}\big)}
\def\torcido{\langle\negthinspace\langle\tau\rangle\negthinspace\rangle}
\def\mK{{\eu m}}
\def\JKm{J_K^\mK}
\def\CKm{C_K^\mK}
\def\DKm{D_K^\mK}
\def\PKm{P_K^\mK}
\def\JKSm{J_{K,S}^\mK}
\def\PKSm{P_{K,S}^\mK}
\def\CKSm{C_{K,S}^\mK}
\def\ClKSm{Cl_{K}^\mK(\o_S)}
\def\KSm{K_S^\mK}
\def\enc#1#2{\lceil#1\rceil_#2}
\def\ence#1#2{\lceil#1_#2\rceil_#2}
\def\lubintatemas#1{\oplus_{#1}}
\def\lubintatepor#1{\odot_{#1}}
\def\Km{K^\mK}
\def\Jm{J_K^{\langle\mK\rangle}}
\def\poly{P_1^{\alpha_1}\cdots P_r^{\alpha_r}}
\def\polyn#1{P_1^{\alpha_{1,#1}}\cdots P_r^{\alpha_{r,#1}}}
\def\Ku#1#2{K\big(\sqrt[#1]{#2}\big)}
\def\zg{{\ma Z}[G]}
\def\co#1#2#3{H^#1(#2,#3)}
\def\cob#1#2#3{H^#1(#2|#3)}
\def\cot#1#2#3{\widehat H^#1(#2,#3)}
\def\hot#1#2#3{\widehat H_#1(#2,#3)}
\def\ho#1#2#3{H_#1(#2,#3)}
\def\suma#1#2#3{\sum_{#1\in #2} #3_#1#1}
\def\too#1#2{\xrightarrow[]{\partial^{#2}_#1}}
\def\tooi#1#2{\xleftarrow[]{d^{#2}_#1}}
\def\tooo#1#2{\xrightarrow[]{d^{#2}_#1}}
\def\toooh#1{\xrightarrow[]{\widehat d_#1}}
\def\Fr#1#2{\Frobenius_{#1|#2}}
\def\n#1#2{\eta_{#1|#2}}
\def\Fro#1{\Frobenius_{#1}}

\def\P{{\mathcal P}}
\def\Pi{{\mc P}_{\infty}}
\def\PP#1{{\ma P}_#1}

\def\sep#1{#1^{{\rm sep}}}
\def\me{{\eu m}=\prod_{\pK\in {\ma P}_K}\pK^{n_\pK}}

\def\AE#1{{\cal O}_#1}
\def\cic#1#2{{\ma Q}(\zeta_{#1^{#2}})}
\def\f#1{{\eu f}_{#1}}
\def\ext#1#2#3{[K(\Lambda_{P^{#1}}): {#2}K(\Lambda_{P^{#3}})]}
\def\lam#1#2{\lambda_{#1^{#2}}}
\def\lama#1{K(\Lambda_{P^#1})}
\def\Lam#1#2{\Lambda_{#1^{#2}}}
\def\cicl#1#2{K(\Lambda_{#1^{#2}})}
\def\units #1#2{\big(R_T/\langle #1^{#2} \rangle\big)^{\ast}}
\def\m{\xrightarrow[m\to\infty]{}}
\def\g#1{#1_{\eu{ge}}}
\def\ge#1#2{#1_{\eu{ge},#2}}
\def\gge#1#2{(#1_#2)_{\eu{ge}}}
\def\ggge#1{(#1_1#1_2)_{\eu{ge}}}
\def\gex#1{#1_{\eu{gex}}}
\def\ggex#1#2{#1_{\eu{gex},#2}}
\def\dif{\eu{Dif}}
\def\*#1{#1^{\ast}}
\def\producto{{\eu p}_1^{e_1}\cdots {\eu p}_r^{e_r}}
\def\ff{{\ma F}_q^{\ast}}
\def\finite#1{{\ma F}_{q^{#1}}}
\def\matriz#1#2{\left[\begin{array}{ccccc}
#1& #1^p&\cdots&#1^{p^{n-2}}&#1^{p^{n-1}}\\
#2_2&#2_2^p&\cdots&#2_2^{p^{n-2}}&#2_2^{p^{n-1}}\\
\vdots&\vdots&\ddots&\vdots&\vdots\\
#2_{n-1}&#2_{n-1}^p&\cdots&#2_{n-1}^{p^{n-2}}&#2_{n-1}^{p^{n-1}}\\
#2_n&#2_n^p&\cdots&#2_n^{p^{n-2}}&#2_n^{p^{n-1}}
\end{array}\right]}
\def\vmatriz#1{\left[\begin{array}{ccccc}
\vec#1_1& \vec#1^p_1&\cdots&\vec#1^{p^{n-2}}_1&\vec #1^{p^{n-1}}_1\\
\vec#1_2&\vec#1_2^p&\cdots&\vec#1_2^{p^{n-2}}&\vec#1_2^{p^{n-1}}\\
\vdots&\vdots&\ddots&\vdots&\vdots\\
\vec#1_{n-1}&\vec#1_{n-1}^p&\cdots&\vec#1_{n-1}^{p^{n-2}}&
\vec#1_{n-1}^{p^{n-1}}\\
\vec#1_n&\vec#1_n^p&\cdots&\vec#1_n^{p^{n-2}}&\vec#1_n^{p^{n-1}}
\end{array}\right]}
\def\menoszeta#1{\zeta_{2^{#1}}-\zeta_{2^{#1}}^{-1}}

\def\idel#1{\tilde{\vec{#1}}}

\def\l{\begin{list}
	{{\rm {(\roman{bean})}}}{\usecounter{bean}
	\setlength{\rightmargin}{\leftmargin}}}
	
\def\las{\begin{list}
	{{\rm {(\arabic{2bean})}}}{\usecounter{2bean}
	\setlength{\rightmargin}{\leftmargin}}}

\def\lasa{\begin{list}
	{{\rm {(\alph{3bean})}}}{\usecounter{3bean}
	\setlength{\rightmargin}{\leftmargin}}}
	
\def\fin{\hfill\qed\medskip}

\newcounter{bean}
\newcounter{otrobean}
\newcounter{2bean}
\newcounter{3bean}

\spdefaulttheorem{notacion}{Notaci\'on}{\bfseries}{\rmfamily}
\spdefaulttheorem{ejercicio}{Ejercicio}{\bfseries}{\rmfamily}
\spdefaulttheorem{teorema}{Teorema}{\bfseries}{\itshape}
\spdefaulttheorem{proposicion}{Proposici\'on}{\bfseries}{\itshape}
\spdefaulttheorem{definicion}{Definici\'on}{\bfseries}{\rmfamily}
\spdefaulttheorem{corolario}{Corolario}{\bfseries}{\itshape}
\spdefaulttheorem{lema}{Lema}{\bfseries}{\itshape}
\spdefaulttheorem{observacion}{Observaci\'on}{\bfseries}{\rmfamily}
\spdefaulttheorem{ejemplo}{Ejemplo}{\bfseries}{\rmfamily}
\spdefaulttheorem{ejemplos}{Ejemplos}{\bfseries}{\rmfamily}
\spdefaulttheorem{pregunta}{Pregunta}{\bfseries}{\rmfamily}
\spdefaulttheorem{conjetura}{Conjetura}{\bfseries}{\rmfamily}
\spdefaulttheorem {teo} {Teorema}{\bfseries}{\itshape}
\spdefaulttheorem {prop} {Proposici\'on}{\bfseries}{\itshape}
\spdefaulttheorem {cor} {Corolario}{\bfseries}{\itshape}
\spdefaulttheorem {lem}  {Lema}{\bfseries}{\itshape}
\spdefaulttheorem {ej} {Ejemplo}{\bfseries}{\rmfamily}
\spdefaulttheorem{dfn}{Definici\'on}{\bfseries}{\rmfamily}
\spdefaulttheorem{notacion y definicion}{Notaci\'on y definici\'on}{\bfseries}{\rmfamily}
\spdefaulttheorem{consecuencia}{Consecuencia}{\bfseries}{\rmfamily}

\makeindex

\begin{document}

\author{{\bf{Martha Rzedowski Calder\'on}}\\
\vspace{-.5cm}
{\scriptsize{Departamento de Control Autom\'atico,}}\\
\vspace{-.5cm}
{\scriptsize{Centro de Investigaci\'on y de Estudios Avanzados del I.P.N.}}\\
\vspace{-.5cm}
{\scriptsize{\tt{mrzedowski@ctrl.cinvestav.mx}}}\\
\vspace{-.5cm}
\\
\\
{\bf{Gabriel Villa Salvador}}\\
\vspace{-.5cm}
{\scriptsize{Departamento de Control Autom\'atico,}}\\
\vspace{-.5cm}
{\scriptsize{Centro de Investigaci\'on y de Estudios Avanzados del I.P.N.}}\\
\vspace{-.5cm}
{\scriptsize{\tt{gvillasalvador@gmail.com, gvilla@ctrl.cinvestav.mx}}}}
\title{{\bf{Campos ciclot\'omicos}}}
\subtitle{Num\'ericos y de funciones (tercera versi\'on)\\
Con una introducci\'on a la teor\'ia de campos de clase}

\date{\huge\fontfamily{pzc}\selectfont 14 de diciembre de 2022}

\maketitle

\frontmatter

\include{dedic}

\include{Introduccion}

\tableofcontents

\mainmatter

\include{Capitulo1} 
\include{Capitulo2} 
\include{Capitulo3}

\include{Capitulo4}

\include{Capitulo5}

\include{Capitulo6}

\include{Capitulo7} 
\include{Capitulo8} 
\include{Capitulo9} 
\include{Capitulo10} 
\include{Capitulo11}
\include{Capitulo12}
\include{Capitulo13} 
\include{Capitulo14}
\include{Capitulo15}

\include{Capitulo16}
\include{Capitulo17}

\backmatter

\include{Notaciones}
\include{Bibliografia}
\printindex

\end{document}

%% file: dedic.tex
\thispagestyle{empty}
\vspace*{2.5cm}
\begin{flushright}

{\large A Sof{\'\i}a, a mi padre y a la memoria de mi madre}

\vspace{2.5cm}

{\large A la memoria de mis padres}

\vspace{2.5 cm}

{\large A los que se doblan pero no se doblegan}

\vspace{2.5cm}

{\large{\em{Si nunca se pide un imposible a un alumno,
nunca ense\~nar\'a de lo que es capaz.}}}

\vspace{1cm}

{\large{\em{El arte de la ense\~naza es dar la impresi\'on
de haber sabido toda la vida lo que aprendimos en la
noche anterior.}}}

\vspace{1cm}

{\large{\em{Las matem\'aticas son lo que hacen los matem\'aticos,
cuando ellos dicen que est\'an haciendo matem\'aticas.}}}

\end{flushright}

%% file: Introduccion.tex
\chapter*{Introducci\'on}

Los campos ciclot\'omicos tuvieron su origen en el trabajo de Kummer
en 1847 sobre el \'Ultimo Teorema de Fermat. Un campo ciclot\'omico
$\cic n{}$
es simplemente el campo que se obtiene al
agregar al campo de los n\'umeros racionales ${\ma 
Q}$ las ra{\'\i}ces de la ecuaci\'on $x^n-1$ donde $n$ es un n\'umero
natural.

Debemos hacer notar que este trabajo toca \'unicamente la superficie
de la teor{\'\i}a de campos ciclot\'omicos y no trata sobre muchos de los
temas presentados en libros como el de L. Washington \cite{Was97} o
el de S. Lang \cite{Lan90}. Por otro lado hemos incluido el tema de los
campos de funciones ciclot\'omicos que muestra nuevamente y de manera
clara la similitud que existe entre los campos num\'ericos y los campos
de funciones.

La importancia de los campos ciclot\'omicos es su simplicidad y su utilidad
y relevancia para diversos objetivos, como son el estudio de las extensiones
abelianas, su uso para el estudio del \'Ultimo Teorema de Fermat, 
ejemplifican de diversas formas conceptos que usualmente son de dif{\'\i}cil
acceso como son la ramificaci\'on, los discriminantes, las bases enteras
y as{\'\i} sucesivamente.

Suponemos que el lector est\'a familiarizado con un curso general de
teor{\'\i}a de n\'umeros y uno sobre teor{\'\i}a de Galois
y que conoce, aunque sea de manera superficial,
los conceptos de dominios de Dedekind, anillos de enteros, discriminante,
diferente, bases enteras, etc.

Parte de lo que presentamos en este libro est\'a basado en el libro
de Washington mencionado anteriormente. Nuestro primer cap{\'\i}tulo
es una compilaci\'on de varios conceptos y resultados que usaremos 
a lo largo de este trabajo: discriminante, diferente, grupos de descomposici\'on
e inercia, as{\'\i} como los grupos de ramificaci\'on para el estudio de la
ramificaci\'on salvaje.

El Cap{\'\i}tulo \ref{ch2} presenta la teor{\'\i}a de Galois de extensiones infinitas
con el objetivo de estudiar la teor{\'\i}a de Iwasawa la cual veremos
en el Cap{\'\i}tulo \ref{Ch7}. El Cap{\'\i}tulo \ref{Ch3} es la introducci\'on
a los campos ciclot\'omicos. En el Cap{\'\i}tulo \ref{Ch1} damos una 
demostraci\'on del teorema cl\'asico de Kronecker--Weber
basada en la demostraci\'on
dada por D. Hilbert \cite{Hil98} usando los grupos de ramificaci\'on.

El Cap{\'\i}tulo \ref{ch4} trata sobre un caso especial del teorema de
Dirichlet sobre la infinidad de n\'umeros primos en progresiones aritm\'eticas
y el estudio de los subcampos de los campos ciclot\'omicos $\cic pm$
donde $p$ es un n\'umero primo.

En el Cap{\'\i}tulo \ref{Ch12} usamos los caracteres de Dirichlet para el 
estudio aritm\'etico de las extensiones abelianas de ${\ma Q}$, tal como
fue introducido por Leopoldt \cite{Leo62}. En el Cap{\'\i}tulo \ref{Ch13}
probamos el caso general del teorema de Dirichlet antes mencionado.

El Cap{\'\i}tulo \ref{Ch5} da una breve introducci\'on a los campos de 
funciones y el Cap{\'\i}tulo \ref{Ch6} desarrolla los campos de funciones
ciclot\'omicos basados en los trabajos de Carlitz \cite{Car35,Car38} y de 
Hayes \cite{Hay74} los cuales son an\'alogos a los campos 
ciclot\'omicos num\'ericos. 

El Cap\'itulo \ref{ChRam} consta de dos partes. En la primera desarrollamos
los grupos de clases de divisores y los grupos de clases de ideales para
campos de funciones, con \'enfasis especial en los campos de
funciones congruentes. En la segunda parte presentamos diversos
resultados sobre ramificaci\'on en campos de funciones. Muchos de
estos resultados son completamente generales y otros son aplicables
\'unicamente a campos de funciones congruentes. Un lugar especial
en esta parte lo ocupa el estudio de las extensiones c\'iclicas de grado primo sobre
un campo de funciones racionales.

En el Cap{\'\i}tulo \ref{Ch9*} presentamos
una teor{\'\i}a de extensiones de campos de funciones con campo
de constantes un campo finito, los cuales est\'an generados
por elementos de torsi\'on bajo la acci\'on de Carlitz--Hayes. Esta
teor{\'\i}a puede ser considerada como una teor{\'\i}a de tipo de
Kummer para extensiones del tipo de Carlitz--Hayes.

En el Cap{\'\i}tulo \ref{Ch9'} estudiamos las $p$--extensiones para
campos de funciones de caracter\'istica $p\geq 0$. 
Damos una breve introducci\'on a los
vectores de Witt los cuales estudian las $p$--extensiones
c{\'\i}clicas en caracter{\'\i}stica $p$. Otro tipo de extensiones
que vemos en el Cap\'itulo \ref{Ch9'} son las $p$--extensiones
elementales abelianas. Tambi\'en damos las
diversas definiciones de conductor y establecemos algunas de sus
relaciones. En el Cap{\'\i}tulo \ref{Ch11} aplicamos lo
desarrollado anteriormente y damos una demostraci\'on combinatoria
del an\'alogo del teorema de Kronecker--Weber para campos
de funciones racionales congruentes.

El Cap\'itulo \ref{Ch12*} trata sobre los campos de g\'eneros de
campos de funciones congruentes. Presentamos los resultados
generales de esta teor\'ia basados principalmente en los caracteres
de Dirichlet. En el Cap\'itulo \ref{DrinfeldCh15} damos una breve 
introducci\'on a los m\'odulos de Drinfeld, los cuales generalizan
al m\'odulo de Carlitz el cual es la base de los campos de funciones
ciclot\'omicos. En este cap\'itulo presentamos, usando campos
de clase, la m\'axima extensi\'on abeliana de un campo de
funciones congruente.

El Cap\'itulo \ref{Ch7} es una introducci\'on a la 
teor{\'\i}a de Iwasawa. En esta
\'ultima parte, nos basamos en los Cap{\'\i}tulos 7 y 13 de \cite{Was97}. Para 
poder dar una introducci\'on a la teor{\'\i}a de Iwasawa hemos 
necesitado usar varios teoremas de teor{\'\i}a de campos de clases
los cuales escapan a los alcances de este escrito y por tanto los hemos
usado sin demostraci\'on en varios casos. Algunos de los resultados
necesarios pueden ser consultados en el Cap\'itulo \ref{CClaseC17}.
El lector puede, sin perder
continuidad, esquivar las partes t\'ecnicas de estos teoremas.

Originalmente, el \'ultimo cap{\'\i}tulo, que
pretend\'ia m\'as bien ser un ap\'endice del libro,
trata sobre la teor\'ia de campos de clases, tanto de campos globales como
de campos locales y considerando tambi\'en los campos 
de funciones. El cap\'itulo creci\'o, de manera desmesurada, despu\'es de
que consideramos que val\'ia la pena presentar las bases de la 
teor\'ia de campos de clase con pr\'acticamente todas las demostraciones
y como consecuencia, en el cap\'itulo se estudian la teor\'ia de
Kummer, campos locales, cohomolog\'ia de grupos, campos 
globales y la teor\'ia de campos de clase, tanto locales como
globales. Hacemos notar que decidimos incluir, completo, el caso
de caracter\'istica positiva, incluyendo las demostraciones 
respectivas. Este cap\'itulo
justifica varias afirmaciones hechas a lo largo del libro.

\vfill

\begin{flushright}

Martha Rzedowski Calder\'on,

Gabriel D. Villa Salvador.

Ciudad de M\'exico, diciembre de 2022.

\end{flushright}

%% file: Capitulo1.tex
\chapter{Teor{\'\i}a algebraica de n\'umeros}\label{Ch0}

En este cap{\'\i}tulo introductorio, pretendemos recordar r\'apidamente
varios conceptos que, en diversa medida, nos servir\'an para nuestros
objetivos en este libro. Todos los resultados pueden ser consultados en
\cite{Hil98, Jan73, Lan86, MurEsm2005, Nar90, Ser}. Hacemos notar
que, aunque en este cap\'itulo consideramos campos num\'ericos,
pr\'acticamente todos los resultados son aplicables a los campos
de funciones \cite{Art,Che51,Deu73,Iwa95,Ros2002,Ser,Sti2009,
Vil2003,Vil2006}

\section{Discriminante y diferente}\label{Sec0.1}

Por un campo num\'erico $K$ entenderemos una extensi\'on 
finita del campo de
los n\'umeros racionales ${\ma Q}$ y ${\cal O}_K$ 
denota el anillo de los enteros de
$K$, es decir, ${\cal O}_K:=\{\alpha \in K\mid \Irr(\alpha, x, {\ma Q})
\in {\ma Z}[x]\}$, donde $\Irr(\beta,t,E)$ denota al polinomio
m\'onico de m{\'\i}nimo grado $f(t)$ en la variable $t$ sobre el campo $E$,
$f(t)\in E[t]$, tal que $f(\beta)=0$.

Resulta ser que ${\cal O}_K$ es un ${\ma Z}$--m\'odulo libre de rango
$[K:{\ma Q}]$. En general, si $L/K$ es una extensi\'on de campos 
num\'ericos, ${\cal O}_L$ no necesariamente es un ${\cal O}_K$--m\'odulo
libre.
Lo que s{\'\i} se tiene es que 
\[
{\cal O}_L\cong {\cal O}_K^{[L:K]-1}\oplus I
\]
como ${\cal O}_K$--m\'odulos, donde $I$ es un ideal de ${\cal O}_K$
(ver Teorema \ref{DrinfeldT1.3.1}).

Empezamos con el resultado de Dirichlet sobre las unidades.

\begin{teorema}[Teorema de las unidades de
Dirichlet\index{Dirichlet!teorema de las unidades de $\sim$}\index{teorema
de las unidades de Dirichlet}]\label{T11.3}
Sea $K$ cualquier campo num\'erico, $[K:{\ma Q}]= r_1+2r_2=n<\infty$,
donde $r_1$ es el n\'umero de encajes reales de $K$ y
$2r_2$ es el n\'umero de encajes complejos de $K$.
Sea $U_K$ el grupo de unidades de $K$, es decir del anillo de enteros
${\cal O}_K$ de $K$: $U_K:={\cal O}_K^{\ast}$. Entonces, como grupos,
tenemos $U_K\cong W_K\times {\ma Z}^{r_1+r_2-1}$ donde
$W_K$ son las ra{\'\i}ces de unidad que hay en ${\cal O}_K$. En otras
palabras, existen unidades $w_1,\ldots, w_{r_1},w_{r_1+1},\ldots,
w_{r_1+r_2-1}$ tales que todo elemento $u$ de $U_K$ se escribe de
manera \'unica como
\[
u=\zeta w_1^{\alpha_1}\cdots w_{r_1+r_2-1}^{\alpha_{r_1+r_2-1}}
\]
con $\zeta$ una ra{\'\i}z de unidad en $K$ y $\alpha_1,\ldots,\alpha_{
r_1+r_2+1}\in{\ma Z}$. En particular, el grupo de torsi\'on de $U_K$ es
$W_K$. $\fin$
\end{teorema}

Sea $L/K$ una extensi\'on finita de campo num\'ericos.

Para cualquier subconjunto $\{\xi_1,\ldots, \xi_n\}\subseteq {\cal O}_L$ se define
el {\em discriminante\index{discriminante}} de $\{\xi_1,\ldots, \xi_n\}$ por
\[
{\eu d}_{L/K}(\xi_1,\ldots, \xi_n)=\big(\det[\sigma_j \xi_i]_{1\leq i,j\leq n}\big)^2=
\det(\Tr_{L/K} \xi_i \xi_j)
\]
donde $\sigma_1,\ldots,\sigma_n$ son los $n=[L:K]$ diferentes $K$--encajes
de $L$ en una cerradura algebraica $\bar{K}$ de $K$ y $\Tr_{L/K}$ denota
la traza de $L$ a $K$.

Si $L=K(a)$ y $p(x)=\Irr(a,x,K)$ entonces
\begin{align*}
{\eu d}_{L/K}(a):&={\eu d}_{L/K}(1,a,a^2,\ldots,a^{n-1})=\prod_{i<j}(\sigma_i (a)-
\sigma_j(a))^2\\
&=(-1)^{n(n-1)/2}N_{L/K}(p'(a))
\end{align*}
donde $N_{L/K}(a)$ denota la norma de $a$.

Sea $\{\alpha_1,\ldots,\alpha_n\}$ una base entera de $K/{\ma Q}$,
es decir, ${\cal O}_K={\ma Z}\alpha_1\oplus \cdots\oplus{\ma Z}\alpha_n$, y
 $\{\beta_1,
\ldots,\beta_n\}\subset {\cal O}_K$. Entonces si $\beta_i=\sum_{j=1}^n a_{ij}
\alpha_j$, $n=[K:{\ma Q}]$ con $a_{ij}\in{\ma Z}$, se tiene
\[
{\eu d}_{K/{\ma Q}}(\beta_1,\ldots,\beta_n)= m^2{\eu d}_{L/K}(\alpha_1,
\ldots,\alpha_n), \quad \text{donde}\quad m=\det[a_{ij}]\in{\ma Z}
\]
y en particular, si $\{\beta_1,\ldots,\beta_n\}$ es tambi\'en 
una base entera, entonces
$m=\pm 1$ y 
\[
\delta_K:={\eu d}_{K/{\ma Q}}(\beta_1,\ldots,\beta_n)=
{\eu d}_{K/{\ma Q}}(\alpha_1,\ldots,\alpha_n)
\]
se llama el {\em discriminante del campo\index{discriminante del campo}} y
es independiente de la base entera.

Se tiene que el signo de $\delta_K$ es igual a $(-1)^{r_2}$ donde $[K:{\ma Q}]=
r_1+2r_2$, $r_1$ es el n\'umero de encajes reales de $K$ y $2r_2$ es el n\'umero
de encajes complejos, esto es, no reales, de $K$. 
Para una demostraci\'on ver Teorema \ref{T1.2.1.4}. Adem\'as se tiene que
\[
\delta_K\equiv 0 \text{\ \'o\ } 1\bmod 4.
\]

Para una extensi\'on finita $L/K$ cualquiera de campos num\'ericos, sea 
$M\subseteq L$ un subgrupo aditivo de $L$. El {\em m\'odulo
complementario\index{modulo complementario@m\'odulo complementario}} 
a $M$ se define por:
\[
M':=\{x\in L\mid \Tr_{L/K}(xm)\in {\cal O}_K \ \forall\ m\in M\}.
\]
Notemos que si $M$ es un ${\cal O}_K$--m\'odulo, entonces
$M'$ tambi\'en es un ${\cal O}_K$--m\'odulo. Ahora bien, si
$\{\alpha_1,\ldots,\alpha_n\}$ es una base de $L/K$, la {\em base
dual\index{base dual}} de $\{\alpha_1,\ldots,\alpha_n\}$ se define
por $\{\alpha_1',\ldots,\alpha_n'\}$ donde 
\[
\Tr_{L/K}(\alpha_i\alpha_j')=\delta_{ij}=
\left\{
\begin{array}{ccc}
1&{\text si}&i=j\\0&{\text si}&i\neq j
\end{array}
\right.,
\quad 1\leq i,j\leq n.
\]
Si $M={\cal O}_K\alpha_1+\cdots+{\cal O}_K\alpha_n$, entonces 
se tiene que $M'={\cal O}_K\alpha_1'+\cdots+{\cal O}_K\alpha_n'$.

El {\em diferente\index{diferente}} ${\eu D}_{L/K}$ de $L/K$ se define por
\[
{\eu D}_{L/K}^{-1}:={\cal O}_L'
\]
es decir, como el inverso del m\'odulo complementario de ${\cal O}_L$.

Si ${\cal P}$ es un ideal primo de ${\cal O}_K$, se define la
{\em conorma\index{conorma}} de ${\cal P}$ por
\[
\con_{K/L}{\pK}={\pK}{\cal O}_L={\pL}_1^{e_1}\cdots {\pL}_g^{e_g}.
\]
Si $K\subseteq L\subseteq F$ es una torre de campos num\'ericos, entonces
se tiene que
\begin{gather}\label{Ec1.0}
{\eu D}_{M/K}={\eu D}_{M/L} \con_{L/M}{\eu D}_{L/K}.
\end{gather}

 En general, para una extensi\'on $E/F$, si $\pK$ es un
primo en $F$ y $\pL$ es un primo en $E$ sobre $\pK$, denotamos
por $a_{E/F}(\pL|\pK)$ al exponente de $\pL$ en el diferente ${\eu D}_{
E/F}$ y por $e_{E/F}(\pL|\pK)$ al \'indice de ramificaci\'on de $\pL/\pK$.

Sea $K\subseteq L\subseteq M$ una torre ya sea de campos num\'ericos
o de campos de funciones. Sea $\P$ un primo en $K$, $\pK$ un primo
en $L$ sobre $\P$ y $\pL$ un primo en $M$ sobre $\P$ y $\pK$.
\[
\xymatrix{
M\ar@{--}[r]\ar@{-}[d]&\pL\ar@{--}[d]
\ar@/^1pc/@{-}[dd]^{{\eu D}_{M/K}={\eu D}_{M/L}\con_{L/M}{\eu D}_{L/K}}\\ 
L\ar@{-}[d]\ar@{--}[r]&\pK\ar@{--}[d]\\
K\ar@{--}[r]&\P
}
\]

Como consecuencia inmediata de (\ref{Ec1.0}), se tiene
\begin{gather}\label{Ec2.0}
a_{M/K}(\pL|\P)=a_{M/L}(\pL|\pK)+e_{M/L}(\pL|\pK)a_{L/K}(\pK|\P).
\end{gather}

El {\em discriminante\index{discriminante}} ${\eu d}_{L/K}$ de la extensi\'on
$L/K$ se define por
\[
{\eu d}_{L/K}:= N_{L/K} {\eu D}_{L/K}.
\]
Se tiene que ${\eu d}_{L/{\ma Q}}=\langle\delta_L\rangle$.

Un resultado muy \'util para el c\'alculo del diferente, es el siguiente:

\begin{teorema}\label{T1.3.7} Se tiene que ${\eu D}_{L/K}$ es el m\'aximo
com\'un divisor del conjunto
\begin{gather*}
\{f'(\alpha)\mid \alpha\in \AE L, L=K(\alpha), f(x)=\Irr(\alpha,x,K)\}\\
= \langle f'(\alpha)\mid \alpha\in \AE L, L=K(\alpha), f(x)=\Irr(\alpha,x,K)\rangle.
\tag*{$\fin$}
\end{gather*}
\end{teorema}

\begin{corolario}\label{C1.3.7'} Si $\AE L=\AE K[\alpha]$ 
entonces ${\eu D}_{L/K}
=\langle f'(\alpha)\rangle$ donde $f(x):=\Irr(\alpha,x,K)$. $\fin$
\end{corolario}

Es decir, en general tenemos para una extensi\'on $L/K$:
\begin{gather*}
{\eu D}_{L/K}=\langle f'(\alpha)\mid \alpha\in {\cal O}_L, L=K(\alpha), 
f(x)=\Irr(\alpha,x,K)\rangle\\
\intertext{y}
{\eu d}_{L/K}=\langle N_{L/K}f'(\alpha)\mid \alpha\in {\cal O}_L, L=K(\alpha), 
f(x)=\Irr(\alpha,x,K)\rangle.
\end{gather*}
En particular, si ${\cal O}_K={\ma Z}[\alpha]$, ${\eu D}_{L/{\ma Q}}=\langle
f'(\alpha)\rangle$ y 
\[
\delta_K=(-1)^{n(n-1)/2}N_{K/{\ma Q}}(f'(\alpha))=\prod_{i<j}(\alpha^{(i)}
-\alpha^{(j)})^2
\]
donde $\alpha=\alpha^{(1)},\ldots, \alpha^{(n)}$ son los conjugados
de $\alpha$ y $n=[K:{\ma Q}]$.

\begin{teorema}[Kummer\index{teorema de Kummer}\index{Kummer!teorema
de $\sim$}]\label{T8.12}
Sean $A$ un dominio Dedekind, $K=\coc A$ el campo de cocientes
de $A$. Sean $E/K$ una extensi\'on finita y separable y $B$
la cerradura entera de $A$ en $E$. Supongamos que
existe $\alpha\in B$ tal que $B=A[\alpha]$. Sea
$f(x):=\Irr(\alpha,x,{\ma Q})$. Fijemos un primo $\pK$ en $A$.
Sea $\overline{f(x)}:=f(x)\bmod \pK$ y sea $\overline{f(x)}=
\overline{P_1(x)}^{e_1}\cdots \overline{P_g(x)}^{e_g}$ la 
factorizaci\'on de $\overline{f(x)}$ en factores irreducibles
m\'onicos en $(A/\pK)[x]$ y $P_i(x)\in A[x]$ m\'onico.
Entonces
\[
\pK B=\pL_1^{e_1}\cdots \pL_g^{e_g}\qquad \qquad
\xymatrix{
\pL_1,\cdots, \pL_g\ar@{--}[r]\ar@{-}[d]
& B\ar@{-}[r]\ar@{-}[d]&E\ar@{-}[d]\\
\pK\ar@{--}[r]&A\ar@{-}[r]&K
}
\]
donde $\pL_i=\pK+\langle P_i(\alpha)\rangle$, es su
descomposici\'on como producto de ideales primos de $B$.
$\fin$
\end{teorema}

\begin{theorem}[Dedekind]\label{T0.1.1}
Sean $L/K$ una extensi\'on finita 
de campos num\'ericos y ${\cal P}$ un ideal primo
no cero de ${\cal O}_K$. Entonces si
\[
\con_{L/K}{\cal P}={\eu p}_1^{e_1}\ldots {\eu p}_g^{e_g}
\]
se tiene que alg\'un $e_i>1$, es decir, ${\eu p}_i$ es ramificado y ${\cal P}$
es ramificado si y solamente si ${\eu p}_i |{\eu D}_{L/K}$ y ${\cal P}|
{\eu d}_{L/K}$. 
\end{theorem}

\begin{proof} Presentamos un esquema de  una parte de
la demostraci\'on para $K/{\ma Q}$.

Si $p\in {\ma Z}$ es un n\'umero primo, ${\eu p}\subseteq {\cal O}_K$ es un
primo encima de $p$ y ${\eu D}_{K/{\ma Q}}$ es el diferente de $K/{\ma Q}$,
entonces veremos que
si ${\eu p}^e|\langle p\rangle $, se tiene que ${\eu p}^{e-1}|{\eu D}_{K/{\ma Q}}$ lo
cual implicar\'a en particular que si ${\eu p}$ es ramificado, entonces
$e\geq 2$ y  por lo tanto $e-1\geq 1$ y en particular ${\eu p}|{\eu D}$.

Para probar la afirmaci\'on anterior
pongamos $(p)=p{\cal O}_K={\eu p}^m {\eu a}$ con $\mcd({\eu a},{\eu p})=1$
y $m\geq e$. Entonces si $x\in {\eu p}{\eu a}$, $x=\sum_{i=1}^n \alpha_i
\beta_i$ con $\alpha_i \in {\eu p}, \beta_i \in {\eu a}$. Por tanto
$x^{p^t}\equiv \sum_{i=1}^n\alpha_i^{p^t} \beta_i^{p^t}\bmod p$. Para $t$ 
suficientemente grande se tiene que $\alpha_i^{p^t}\in{\eu p}^m$ y por tanto
$x^{p^t}\in {\eu p}^m {\eu a}=\langle p\rangle$. En particular obtenemos que
$\Tr_{K/{\ma Q}} x^{p^t}\in p{\ma Z}$.

Se tiene que $\Tr_{K/{\ma Q}} x^{p^t}\in p{\ma Z}
\Rightarrow (\Tr_{K/{\ma Q}}x)^{p^t}\in p{\ma Z}
\Rightarrow \Tr_{K/{\ma Q}} x\in p{\ma Z}\Rightarrow \Tr_{K/{\ma Q}}(p^{-1}
{\eu p}{\eu a})\subseteq {\ma Z}\Rightarrow p^{-1}{\eu p}{\eu a}\subseteq
{\eu D}_{K/{\ma Q}}^{-1}\Rightarrow {\eu D}_{K/{\ma Q}} p^{-1}{\eu p}{\eu a}
\subseteq {\cal O}_K\Rightarrow {\eu D}_{K/{\ma Q}}\subseteq p{\eu p}^{-1}
{\eu a}^{-1}={\eu p}^m{\eu a}{\eu p}^{-1}{\eu a}^{-1}= {\eu p}^{m-1}
\Rightarrow {\eu p}^{m-1} |{\eu D}_{K/{\ma Q}}$. \fin
\end{proof}

\section{Ramificaci\'on en campos num\'ericos}\label{Sec0.2}

Los resultados que aqu{\'\i} presentamos pueden ser consultados en
\cite{Lan86} y en \cite[Cap{\'\i}tulo 4]{Nar90}.

Uno de los problemas centrales que se presentan en la Teor{\'\i}a de N\'umeros
es el problema de la {\em ramificaci\'on}. Recordaremos a continuaci\'on algunos
resultados generales que aplicaremos a nuestro caso particular de los campos
ciclot\'omicos.

Como hemos recordado anteriormente, el diferente inverso  ${\eu D}_{L/K}^{-1}$
de una extensi\'on de campos num\'ericos est\'a dado como el m\'odulo
complementario del anillo de enteros ${\cal O}_L$ mediante la traza. M\'as 
precisamente
\begin{equation}\label{Eq1.3.1'}
{\eu D}_{L/K}^{-1}:= \{x\in L\mid \Tr_{L/K}(x{\cal O}_L)\subseteq {\cal O}_K\}=:
{\cal O}'_L.
\end{equation}

Se tiene que ${\eu D}_{L/K}^{-1}\supseteq {\cal O}_L$ y es un ideal fraccionario.

Si ${\eu p}$ es un ideal primo no cero de $\AE K$, entonces el ideal extendido
de ${\eu p}$ en $\AE L$ es $\pK \AE L= \pL_1^{e_1}\cdots \pL_g^{e_g}$ donde 
\[
[L:K]=\sum_{i=1}^g e_i f_i\quad\text{y}\quad f_i=[\AE L/\pL_i:\AE K/\pK].
\]

Si $e_i>1$ decimos que $\pL_i$ es {\em ramificado\index{ramificado!primo}} y que
$\pK$ es {\em ramificado}. Adem\'as tenemos que $N_{L/K}\pL_i=\pK^{f_i}$.

La conexi\'on entre el diferente y la ramificaci\'on est\'a expl{\'\i}citamente
dada en el siguiente resultado, el cual precisa el Teorema \ref{T0.1.1}.

\begin{teorema}\label{T1.3.1} Si $\pK$ es un ideal no primo de $\AE K$ y si
$\pK \AE L= \pL_1^{e_1}\cdots \pL_g^{e_g}$, entonces $\pL_i^{e_i-1}|{\eu D}_{L/K}$.
M\'as a\'un, $\pL_i^{e_i}|{\eu D}_{L/K}$ si y solamente si $p|e_i$ donde $p$
es la caracter{\'\i}stica del campo residual $\AE K/\pK$. $\fin$
\end{teorema}

\begin{corolario}\label{C1.3.2} Se tiene que $\pL|{\eu D}_{L/K}$ si y solamente si
$\pL$ es ramificado. $\fin$
\end{corolario}

\begin{corolario}\label{C1.3.3} El n\'umero de primos de $\AE L$ ramificados
en $L/K$ es finito y los primos ramificados
son exactamente los divisores de ${\eu D}_{L/K}$. $\fin$
\end{corolario}

\begin{definicion}\label{D1.3.4} Con las notaciones anteriores, decimos que
$\pK$ (o que $\pL_i$) es {\em salvajemente ramificado\index{ramificaci\'on
salvaje}} si $p|e_i$ y {\em moderadamente ramificado\index{ramificaci\'on
moderada}} si $p\nmid e_i$.
\end{definicion}

Con respecto al discriminante, tenemos:

\begin{corolario}\label{C1.3.6} El n\'umero de primos de $\AE K$ ramificados
en $L$ es finito y los primos ramificados
son exactamente los divisores de ${\eu d}_{L/K}$. $\fin$
\end{corolario}

\begin{notacion}\label{N1.3.7} Cuando $K={\ma Q}$, ponemos simplemente
${\eu d}_L:={\eu d}_{L/{\ma Q}}$ y ${\eu D}_L:={\eu D}_{L/{\ma Q}}$.
\end{notacion}

\begin{ejemplo}\label{E1.3.8} Consideremos una extensi\'on cuadr\'atica de
${\ma Q}$: $L:={\ma Q}(\sqrt{d})$ donde $d\in {\ma Z}$ es libre de cuadrados,
esto es, $d=p_1\cdots p_r$ donde $p_1,\ldots, p_r$ son primos distintos.

Se tiene que $\sqrt{d}\in\AE L$ y que $f(x)=
\Irr(\sqrt{d},x,{\ma Q})=x^2-d$. Por tanto
$f'(\sqrt{d})=2\sqrt{d}$. Se sigue que ${\eu D}_L|\langle2\sqrt{d}\rangle$ y que
${\eu d}_L|\langle N_{L/{\ma Q}}(2\sqrt{d})\rangle = 
\langle -4d\rangle =\langle 4d\rangle$.
En particular, los primos ramificados se encuentran en $\{2,p_1,\ldots,p_r\}$.

Sea $\alpha=a+b\sqrt{d}\in \AE L$. Si $b=0$, entonces $\alpha=a\in {\ma Z}$ e
$\Irr(\alpha, x,{\ma Q})=x-\alpha=:p_{\alpha}(x)$. 
Se tiene $p_{\alpha}'(x)=1$. Si $b\neq 0$,
entonces $\Irr(\alpha,x, {\ma Q})=(x-a)^2-b^2d=p_{\alpha}(x)=x^2-2ax+a^2-b^2d\in
{\ma Z}[x]$. Obtenemos que $p_{\alpha}'(\alpha)=2(\alpha-a)=2b\sqrt{d}$ y
$N_{L/{\ma Q}}(p_{\alpha}'(\alpha))=-4b^2 d=(2b)^2(-d)\in{\ma Z}$.

Puesto que $b\in{\ma Q}$ y $(2b)^2d\in{\ma Z}$, 
si escribimos $b=\frac{\gamma}{\beta}$
con $\gamma,\beta\in{\ma Z}$ y primos relativos 
se tiene que $\beta|2$ pues si alg\'un
n\'umero primo $p$ divide a $\beta$, entonces tenemos que $(2b)^2=\frac{
4\gamma^2}{p^2 \beta_1}$ lo cual implica que $p^2|4$. 
En particular obtenemos que
\[
\langle N_{L/{\ma Q}}(p_{\alpha}'(\alpha))\mid \alpha\in{\ma Q}\rangle =
\left\{
\begin{array}{ccl}
\langle d \rangle&\text{si}&\text {existe } \alpha=a+b\sqrt{d}\in \AE L, b\not\in {\ma Z}\\
\langle 4d\rangle&& \text{en otro caso}.
\end{array}
\right.
\]

Ahora si existe $\alpha=a+b\sqrt{d}\in \AE L$ con $b\not\in{\ma Z}$, entonces
$b=\frac{b_1}{2}$ con $b_1\in{\ma Z}$ impar. Tenemos que si
 $a\in {\ma Z}$, entonces $b\sqrt{d}\in\AE L$ pero $\Irr(b\sqrt{d},x,{\ma Q})=
 x^2-b^2d\not\in {\ma Z}[x]$ lo cual es absurdo. Se sigue que $a\not\in{\ma Z}$ pero
 $2a\in{\ma Z}$, esto es, $a=\frac{a_1}{2}$ con $a_1$ impar.
 
 Obtenemos que $a^2-b^2d=\frac{a_1^2}{4}-\frac{b_1^2}{4}d=\frac{a_1^2-b_1^2d}{4}
 \in{\ma Z}$ lo cual implica que $a_1^2-b_1^2d\equiv 0\bmod 4$. Por lo tanto
 $1\equiv a_1^2\equiv b_1^2 d\bmod 4\equiv d\bmod 4$. De esto se obtiene que
 $d\equiv 1\bmod 4$ y ${\eu d}_{{\ma Q}(\sqrt{d})}=\langle d\rangle$.
 
 En otro caso, esto es, si $d\not\equiv 1\bmod 4$, entonces $d\equiv 2, 3\bmod 4$
 y ${\eu d}_{{\ma Q}(\sqrt{d})} =\langle 4d \rangle$. Escribiendo $p_i\AE L {{\ma Q}(
 \sqrt{d})}=\pL_i^2$, $1\leq i\leq r$, se obtiene finalmente,
 \[
 {\eu D}_{{\ma Q}(\sqrt{d})}=\left\{
 \begin{array}{cl}
 \pL_1\cdots \pL_r,& d\equiv 1\bmod 4,\\
 \pL_0^2\pL_1\cdots \pL_r,& d\equiv 3\bmod 4,\\
 \pL_1^3\pL_2\cdots \pL_r,& d\equiv 2\bmod 4, \quad p_1=2,
 \end{array}
 \right.
 \]
 donde $2{\cal O}_{{\ma Q}(\sqrt{d})}=\pL_0^2$ cuando $d\equiv 3\bmod 4$.
\end{ejemplo}

\section{Grupos de inercia, de descomposici\'on y de ramificaci\'on}\label{Sec0.3}

Para estudiar ramificaci\'on en extensiones 
de Galois tenemos a nuestra disposici\'on
los grupos de inercia, de descomposici\'on y de ramificaci\'on. M\'as precisamente,
sea $L/K$ una extensi\'on finita
de Galois de campos num\'ericos. Sea $\pK$ un primo
en $\AE K$ y sea $\pL$ un primo en $\AE L$ que divide a $\pK$, esto es,
$\pK\AE L=\pL{\eu a}$ con ${\eu a}$ un ideal de $\AE L$. Sea $G:=\Gal(L/K)$.
Se define:

\begin{definicion}\label{D1.3.9} El {\em grupo de descomposici\'on\index{grupo
de descomposici\'on}} $D(\pL|\pK)$ se define por
\[
D(\pL|\pK):=\{\sigma\in G\mid \sigma\pL=\pL\}.
\]
\end{definicion}

Notemos que si $\sigma\in D(\pL|\pK)$ entonces $\sigma(\AE L)=\AE L$ y por
tanto $\sigma$ induce un automorfismo
\[
\tilde{\sigma}\colon \AE L/\pL\to \AE L/\pL
\]
tal que $\tilde{\sigma}|_{\AE K/\pK}=\Id_{\AE K/\pK}$. En otras palabras
$\tilde{\sigma}\in\Gal\big((\AE L/\pL)/(\AE K/\pK)\big)$. El mapeo 
\begin{eqnarray*}
\theta\colon D(\pL|\pK)&\lra&\Gal((\AE L/\pL)/(\AE K/\pK))\\
\sigma&\longmapsto&\tilde{\sigma}
\end{eqnarray*}
es un epimorfismo. El n\'ucleo de $\theta$ es el {\em grupo de inercia} de
$\pL$ sobre $\pK$. M\'as precisamente

\begin{definicion}\label{D1.3.10} El {\em grupo de inercia\index{grupo de
inercia}} $I(\pL|\pK)$ de $\pL$ sobre $\pK$ se define por
\[
I(\pL|\pK):=\{\sigma\in D(\pL|\pK)\mid \tilde{\sigma}=\Id_{\AE L/\pL}\}=
\{\sigma\in G\mid \sigma x - x \in \pL\ \forall\  x\in \AE L\}.
\]
\end{definicion}

Se tiene que
\begin{gather*}
\frac{D(\pL|\pK)}{I(\pL|\pK)}\cong\Gal(\AE L/\pL:\AE K/\pK).\\
\intertext{En particular tenemos que el {\em grado relativo\index{grado relativo}} es}
f(\pL|\pK):=[\AE L/\pL : \AE K/\pK]=\frac{|D(\pL|\pK)|}{|I(\pL|\pK)|}.
\end{gather*}
M\'as a\'un, el {\em {\'\i}ndice de ramificaci\'on\index{indice@{\'\i}ndice de ramificaci\'on}}
$e(\pL|\pK)= e$ de $\pL$ sobre $\pK$, es decir $\pK\AE L=\pL^e{\eu a}$ con
${\eu a}$ de $\AE L$ y $\pL$ y ${\eu a}$ primos relativos, es igual a la cardinalidad
de $I(\pL|\pK)$:
\[
e(\pL|\pK)=|I(\pL|\pK)|.
\]
En consecuencia tenemos que $|D(\pL|\pK)|= e(\pL|\pK)f(\pL|\pK)$.

Como $L/K$ es de Galois, tenemos que $e(\sigma \pL|\pK)= e(\pL|\pK)$ y 
$f(\sigma \pL|\pK)=f(\pL|\pK)$ para toda $\sigma \in G$ y si ponemos $e:=e(\pL|\pK)$
y $f:= f(\pL|\pK)$, entonces 
\begin{gather*}
[L:K]=efg\\
\intertext{donde $g$ es el n\'umero de ideales primos
de $\AE L$ que dividen a $\pK$:}
\pK\AE L=(\pL_1\cdots\pL_g)^e \quad {\rm con}\quad f=[\AE L/
\pL_i:\AE K/\pK].
\end{gather*}

Ahora consideremos una extensi\'on finita de Galois $L/K$ y un subcampo
$K\subseteq E \subseteq L$ tal que $E/K$ tambi\'en es Galois. En particular
tenemos que $\Gal(L/E)$ es un subgrupo normal de $\Gal(L/K)$,
$\Gal(L/E)\lhd \Gal(L/K)$. Fijemos un primo $\pL$ de $L$ y sean
$\pK:=\pL\cap E$ y ${\mc P}:=\pL\cap K=\pK\cap K$ las restricciones
de $\pL$ a $E$ y $K$ respectivamente.

\begin{teorema}\label{P3.1.1.1} Se tiene
\las
\item $D(\pL|\pK)=D(\pL|{\mc P})\cap \Gal(L/E)$ y 
$I(\pL|\pK)=I(\pL|{\mc P})\cap \Gal(L/E)$.

En este inciso, no es necesario pedir que $E/K$ sea Galois.

\item El siguiente diagrama es conmutativo y tanto sus columnas como
sus filas son exactas:
\[
\xymatrix{
&1\ar@{->}[d]&1\ar@{->}[d]&1\ar@{->}[d]\\
1\ar@{->}[r]&I(\pL|\pK)\ar@{->}[r]\ar@{->}[d]&
I(\pL|{\mc P})\ar@{->}[r]\ar@{->}[d]&
I(\pK|{\mc P})\ar@{->}[r]\ar@{->}[d]&1\\
1\ar@{->}[r]&D(\pL|\pK)\ar@{->}[r]\ar@{->}[d]&
D(\pL|{\mc P})\ar@{->}[r]\ar@{->}[d]&
D(\pK|{\mc P})\ar@{->}[r]\ar@{->}[d]&1\\
1\ar@{->}[r]&\Gal(\bar{L}/\bar{E})\ar@{->}[r]\ar@{->}[d]&
\Gal(\bar{L}/\bar{K})\ar@{->}[r]\ar@{->}[d]&
\Gal(\bar{E}/\bar{K})\ar@{->}[r]\ar@{->}[d]&1\\
&1&1&1
}
\]
donde $\bar{L}={\mc O}_L/\pL$, $\bar{E}={\mc O}_E/\pK$
y $\bar{K}={\mc O}_K/{\mc P}$.

En particular se tiene que
\begin{gather*}
D(\pK|{\mc P})\cong \frac{D(\pL|{\mc P})}{D(\pL|\pK)}=
\frac{D(\pL|{\mc P})}{D(\pL|{\mc P})\cap \Gal(L/E)}
\cong \frac{D(\pL|{\mc P}) \Gal(L/E)}{\Gal(L/E)}
\intertext{y}
I(\pK|{\mc P})\cong \frac{I(\pL|{\mc P})}{I(\pL|\pK)}=
\frac{I(\pL|{\mc P})}{I(\pL|{\mc P})\cap \Gal(L/E)}
\cong \frac{I(\pL|{\mc P}) \Gal(L/E)}{\Gal(L/E)}.
\end{gather*}
\end{list}
\end{teorema}

\begin{proof} \cite[Proposition 22, Chapter 1]{Ser}. \fin
\end{proof}

Notemos que, como $\AE K/\pK$ es un campo finito, digamos que $\AE K/\pK
\cong {\ma F}_q$, $\AE L/\pL\cong {\ma F}_{q^f}$, entonces $\Gal(\AE L/\pL:\AE K/\pK)$
es un grupo c{\'\i}clico de orden $f$ generado por el {\em automorfismo de
Frobenius\index{automorfismo de Frobenius}}:
\[
\varphi_{\pK}\colon {\ma F}_{q^f}\to {\ma F}_{q^f},\quad \varphi_{\pK}(x)= x^q.
\]

Cuando $\pK$ no es ramificado se tiene que $I(\pL|\pK)=\{\Id\}$ y por lo tanto
existe un \'unico $\theta\in G$ tal que $\tilde{\theta}=\varphi_{\pK}$. 
Este $\theta$ es el 
automorfismo de Frobenius y se denota por: $\theta=\left[\frac{L/K}{\pL}\right]$.

Se tiene que $\left[\frac{L/K}{\sigma \pL}\right]=\sigma\theta\sigma^{-1}=
\sigma \left[\frac{L/K}{\pL}\right]\sigma^{-1}$. Notemos que
\[
D(\sigma \pL|\pK)=\sigma D(\pL|\pK)\sigma^{-1}\quad {\text{y}}
\quad I(\sigma \pL|\pK)=\sigma
I(\pL|\pK)\sigma^{-1}.
\]

El {\em s{\'\i}mbolo de Artin\index{simbolo de Artin@s{\'\i}mbolo de Artin}} $\left(\frac{L/K}{\pK}\right)$
est\'a definido por la clase de conjugaci\'on
\[
\left(\frac{L/K}{\pK}\right)=\left\{\sigma\left[\frac{L/K}{\pL}\right]\sigma^{-1}\mid
\sigma\in G\right\}.
\]

Cuando $L/K$ es una extensi\'on abeliana y $\pK$ es no ramificado se tiene que
el s{\'\i}mbolo de Artin consta de un elemento, esto es, $\left(\frac{L/K}{\pK}\right)
\in G$ y satisface 
\[
\left(\frac{L/K}{\pK}\right)(x)\equiv x^q \bmod \pL \ \forall\ x\in \AE L.
\]

\begin{definicion}\label{D1.3.10bis} Con la notaci\'on anterior, ponemos
$G_{-2}:= G$, $G_{-1}:=D(\pL|\pK)$, $G_0:=I(\pL|\pK)$ y en general para
$i\geq -1$, $i\in {\ma Z}$, se define el {\em $i$--\'esimo grupo de 
ramificaci\'on\index{grupo de ramificaci\'on}} $G_i$ por:
\[
G_i:=\{\sigma\in G_{-1}\mid \sigma a-a\in \pL^{i+1}\ \forall\ a\in\AE L\}.
\]
\end{definicion}

\begin{proposicion}\label{P1.3.10-2bis} Las siguientes condiciones son equivalentes
\l
\item $\sigma\in G_i$, $i\geq -1$, es decir, $\sigma a-a\in \pL^{i+1}\ \forall\ a\in \AE L$.
\item $\sigma \pi-\pi\in \pL^{i+1}$ para un elemento $\pi\in\AE L$ tal que 
$v_{\pL}(\pi)=1$.
\end{list}
\end{proposicion}

\begin{proof} {\underline{(I) $\Rightarrow$ (II)}}: Es inmediata.

\noindent
{\underline{(II) $\Rightarrow$ (I)}}: Sea $E$ el {\em campo de inercia\index{campo de
inercia}} de $\pL$, esto es, $E:=L^{G_0}$. Entonces si ${\eu q}:=\pL\cap \AE E$,
se tiene que ${\eu q}$ es totalmente ramificado en la extensi\'on $L/E$. Tomamos
las localizaciones $B:=\AE {{L,\pL}} \supseteq \AE L$ y $A:= \AE {E,{\eu q}}$. Entonces
$A$ y $B$ son anillos de valuaci\'on discreta y $B$ es un $A$--m\'odulo libre de
rango 
$|G_0|$. Entonces $B=A[\pi]$. Si $a\in B$, entonces se tiene
\[
a=\sum_{i=0}^{e-1}\alpha_i\pi^i,\quad e:=|G_0|, \quad {\text{y}} \quad\alpha_i\in A.
\]
Por lo tanto
\begin{align*}
\sigma a-a &=\sum_{i=0}^{e-1}\alpha_i(\sigma\pi^i-\pi^i)\\
&=
\sum_{i=1}^{e-1}\alpha_i(\sigma\pi-\pi)(\sigma\pi^{i-1}+(\sigma\pi^{i-2})\cdot \pi+
\cdots+ (\sigma\pi)\cdot\pi^{i-2}+\pi^{i-1})\\
&\in \pL^{i+1}B.
\end{align*}
En particular, si $a\in \AE L$, $\sigma a-a \in \pL^{i+1}B\cap \AE L=\pL^{i+1}$. $\fin$
\end{proof}

Se tiene que $G_i$ es un subgrupo normal de $G_{-1}=D(\pL|\pK)$, $G_{i+1}
\subseteq G_i$. Adem\'as para $i$ suficientemente grande tenemos que $G_i=
\{\Id\}$.

Para $\sigma\in G_{-1}$, $\sigma\neq \Id$, existe $i$ tal que $\sigma\in
G_i\setminus G_{i+1}$. Se define $i_{G_{-1}}(\sigma):=i$. Si $\sigma = \Id$ definimos
$i_{G_{-1}}(\sigma)=\infty$. Notemos que $i_{G_{-1}}(\sigma)\geq i+1$ si y s\'olo si
$\sigma\in G_i$. Adem\'as se tiene 
\[
\sum_{\sigma\neq \Id} i_{G_{-1}}(\sigma)=\sum_{i=0}^{\infty}(|G_i|-1).
\]

La conexi\'on con el diferente la obtenemos del siguiente resultado:

\begin{teorema}\label{T1.3.11} Sean $\pL$ y $\pK$ como antes. Sea
$s\geq 0$ la potencia de $\pL$ que aparece en ${\eu D}_{L/K}$. Entonces
\[
s=\sum_{\sigma\neq \Id} i_{G_{-1}}(\sigma)=\sum_{i=0}^{\infty}(|G_i|-1).
\tag*{$\fin$}
\]
\end{teorema}

\begin{corolario}\label{C1.3.11} Se tiene que $\pL$ es salvajemente
ramificado si y solamente si $G_1\neq\{\Id\}$. $\fin$
\end{corolario}

Ahora consideremos un subgrupo $H$ de $G_{-1}$. Sean $E:=L^H$,
el campo fijo de $H$, y ${\eu q}:={\eu P}\cap E$. Entonces
los grupos de ramificaci\'on de $H$ satisfacen:

\begin{proposicion}\label{P1.3.11'}
Para $\sigma \in H$ se tiene $i_H(\sigma)=i_{G_{-1}}(\sigma)$ y
$H_i=G_i\cap H$ para toda $i\geq -1$. $\fin$
\end{proposicion}

Ahora veamos algunas otras propiedades de los grupos de inercia que
usaremos m\'as adelante, en particular en la demostraci\'on del
Teorema de Kronecker--Weber.

Sea $E$ un campo num\'erico, $\AE E$ su anillo de enteros. Sea $\pL$
un ideal primo no cero de $\AE E$. Sea $S:=\AE E\setminus \pL$ el 
complemento de $\pL$ y $\hat{\cal O}_E:= S^{-1}\AE E= (\AE E)_{\pL}$ la
localizaci\'on de de $\AE E$ en $\pL$. Se tiene que  $\hat{\cal O}_E =
\big\{\frac{\alpha}{\beta}\in E\mid \alpha,\beta \in \AE E, \beta \notin \pL\big\}$.
Entonces $\hat{\cal O}_E$ es un anillo local con ideal m\'aximo $\hat{\pL}:=
\pL \hat{\cal O}_E=\big\{\frac{\alpha}{\beta}\in E\mid \alpha\in\pL,\beta\in
\AE E\setminus \pL\big\}$ y las unidades de $\hat{\cal O}_E$ son 
$\hat{\cal O}_E^{\ast}=\hat{\cal O}_E\setminus \hat{\pL}$. Sean
\[
U^{(0)}=\hat{\cal O}_E^{\ast} \quad \text{y} \quad U^{(n)}=1+\hat{\pL}^n
\subseteq \hat{\cal O}_E^{\ast}, \quad n\geq 1.
\]
Entonces $U^{(i)}$ es cerrado con la multiplicaci\'on para $i\geq 0$ y de
hecho si $v_{\pL}$ es la valuaci\'on $\hat{\pL}$--\'adica $U^{(i)}=\{x\in
\hat{\cal O}_E^{\ast}\mid v_{\pL}(x-1)\geq i\}$. Notemos que para
$i\geq 1$, $U^{(i)}/U^{(i+1)}$, el cual entenderemos como las clases de 
equivalencia 
\[
1+x\sim 1+y\iff x-y\in U^{(i+1)},
\]
 es un grupo pues
el inverso de la clase $(\overline{1+x})$ es $(\overline{1+x})^{-1} =
\overline{1-x}$, lo cual se sigue del hecho de que $(1+x)(1-x)=1-x^2\in
U^{(i+1)}$.

\begin{proposicion}\label{P1.3.12} Sea $E$ cualquier campo num\'erico.
Entonces
\las
\item $U^{(0)}/U^{(1)}\cong (\AE E/\pL)^{\ast}$.

\item $U^{(i)}/U^{(i+1)}\cong \hat{\pL}^i/\hat{\pL}^{i+1}\cong \pL^i/\pL^{i+1}
\cong \AE E/\pL$ para $i\geq 1$.
\end{list}
\end{proposicion}

\begin{proof}

\las
\item Sea $\varphi\colon U^{(0)}\to (\hat{\cal O}_E/\hat{\pL})^{\ast}$, 
$\varphi(x)=x\bmod \hat{\pL}$. Entonces $\varphi$ es un epimorfismo de
grupos abelianos multiplicativos. Adem\'as $\ker \varphi = \{x\in \hat{\cal O}_E^{\ast}
\mid \varphi(x)=x\bmod \hat{\pL}\equiv 1\bmod \hat{\pL}\}= U^{(1)}$. Por tanto
$U^{(0)}/U^{(1)}\cong (\AE E/\pL)^{\ast}$.

\item Sea $\psi\colon \hat{\pL}^i\lra U^{(i)}$, $\psi_i(x)=1+x$. Entonces $\psi$ es un
mapeo biyectivo, que no es homomorfismo, tal que
 compuesto con la proyecci\'on natural
$U^{(i)}\to U^{(i)}/U^{(i+1)}$ nos da una funci\'on suprayectiva $\tilde{\psi}_i\colon
\hat{\pL}^i\to U^{(i)}/U^{(i+1)}$ donde tanto $\hat{\pL}^i$ como $U^{(i)}/U^{(i+1)}$ son
grupos abelianos. Sean $x,y\in\hat{\pL}^i$, entonces $\tilde{\psi}_i(x+y)=1+(x+y)\bmod
U^{(i+1)}$ y
\[
\tilde{\psi}_i(x)\tilde{\psi}_i(y) =(1+x)(1+y)\bmod U^{(i+1)}= 1+(x+y)+(xy)\bmod U^{(i+1)}.
\]
Ahora bien, puesto que $x,y\in\hat{\pL}^i$ con $i\geq 1$ y en particular $i+1\leq 2 i$,
entonces $xy\in\hat{\pL}^{2i}\subseteq \hat{\pL}^{i+1}$ y por tanto $1+xy\equiv 1
\bmod U^{(i+1)}$. Se sigue que $\tilde{\psi}_i(xy)=\tilde{\psi}_i(x)\tilde{\psi}_i(y)$.

Entonces $\hat{\psi}_i$ es un epimorfismo de grupos con $\ker \tilde{\psi}_i=
\hat{\pL}^{i+1}$ de donde obtenemos que $U^{(i)}/U^{(i+1)}\cong \hat{\pL}^i/
\hat{\pL}^{i+1}$.

Como siguiente paso probemos que $\pL^i/\pL^{i+1}\cong \hat{\pL}^i/\hat{\pL}^{i+1}$
para $i\geq 1$. Sea $\pL^i\xrightarrow[]{\alpha}\hat{\pL}^i\xrightarrow[]{\beta}
\hat{\pL}^i/\hat{\pL}^{i+1}$ los homomorfismos naturales $\alpha(x)=\frac{x}{1}$
y $\beta(y)=y\bmod \hat{\pL}^{i+1}$. Puesto que claramente tenemos que
$\ker (\beta\circ \alpha)=\pL^{i+1}$, basta probar que $\beta\circ \alpha$ es 
suprayectiva. Sea $\frac{x}{t}\in\hat{\pL}^i$, es decir $x\in\hat{\pL}^i$, $t
\notin \hat{\pL}$. Entonces puesto que $\pL$ es maximal y $t\notin \pL$, se tiene
que $(t)+\pL=\AE E$ y por tanto $\pL^i(t)+\pL^{i+1}=\pL^i$. En particular existen
$a\in\pL^i$ y $z\in\pL^{i+1}$ tales que $at+z=x$ y $\frac{x}{t}-\frac{a}{1}=\frac{z}{t}
\in\hat{\pL}^{i+1}$.

Por tanto $(\beta\circ \alpha)(a)=\frac{a}{1}\bmod \hat{\pL}^{i+1}=\frac{x}{t}
\bmod \hat{\pL}^{i+1}$ y $\beta\circ\alpha$ es suprayectiva probando que
$\pL^i/\pL^{i+1}\cong \hat{\pL}^i/\hat{\pL}^{i+1}$ para $i\geq 0$.

Para $i=0$, $\hat{\cal O}_E/\hat{\pL}\cong \AE E/\pL$. Finalmente probaremos que
$\hat{\pL}^i/\hat{\pL}^{i+1}\cong \hat{\cal O}_E/\hat{\pL}$ y con esto terminaremos la
demostraci\'on. Puesto que $\hat{\cal O}_E$ es un anillo de valuaci\'on y 
$\hat{\pL}=(\pi)$ es principal, el mapeo 
\[
\begin{array}{ccccl}
\hat{\cal O}_E&\to&\hat{\pL}^i&\to&\hat{\pL}^i/\hat{\pL}^{i+1}\\
x&\mapsto&x\pi^i&\mapsto&x\pi^i \bmod \hat{\pL}^{i+1}
\end{array}
\]
es un epimorfismo de grupos con n\'ucleo $\hat{\pL}^i$. $\fin$

\end{list}
\end{proof}

La conexi\'on que tenemos entre los grupos de ramificaci\'on y los grupos
$U^{(i)}/U^{(i+1)}$ nos lo da el siguiente resultado.

\begin{proposicion}\label{P1.3.13} Sea $K/{\ma Q}$ una extensi\'on finita
de Galois con grupo de Galois $G$. Sean $p$ un n\'umero primo racional y
$\pL$ un primo en $\AE K$ sobre $p$. Sea $G_{-1}=D(\pL|p)$ el grupo
de descomposici\'on, $G_0=I(\pL|p)$ el grupo de inercia y $G_i$, $i\geq 1$
los grupos de ramificaci\'on. Sean $U^{(i)}=1+\hat{\pL}^i$, $i\geq 0$, con
$U^{(0)}=\hat{\cal O}_K^{\ast}$. Sea $\pi\in \AE K$ un elemento primo de
$\hat{\pL}$, es decir, $(\pi)=\hat{\pL}$ para lo cual es suficiente
seleccionar $\pi\in\pL\setminus \pL^2$. Entonces $\sigma\in G_i\iff \sigma(\pi)/
\pi\in U^{(i)}$, $i\geq 0$.
\end{proposicion}

\begin{proof}
Recordemos que $\sigma\in G_i\iff \sigma\pi-\pi
\in\pL^{i+1}$ (Proposici\'on \ref{P1.3.10-2bis}).
Por lo tanto para $\sigma\in G_i\iff \sigma\pi-\pi \in \hat{\pL}^{i+1}\iff
\frac{\sigma\pi}{\pi}-1\in \hat{\pL}^i$ siendo esto \'ultimo equivalente a que
$\frac{\sigma\pi}{\pi}\in 1+\hat{\pL}^i=U^{(i)}$. $\fin$
\end{proof}

\begin{proposicion}\label{P1.3.14} Sea $L/K$ una extensi\'on de Galois
con grupo $G$.
Se tiene para $i\geq 0$ que el mapeo
\[
\varphi\colon G_i/G_{i+1}\to U^{(i)}/U^{(i+1)}, \quad \varphi(\sigma \bmod
G_{i+1}):=\sigma\pi/\pi\bmod U^{(i+1)},
\]
donde $\pi$ es un elemento primo de $\pL$, $v_{\pL}(\pi)=1$, es un monomorfismo
de grupos que es independiente de $\pi$. En particular $G_i/G_{i+1}$ es un
$p$--grupo elemental abeliano donde $p$ es la caracter{\'\i}stica de $\AE K/\pL$.
\end{proposicion}

\begin{proof}
Sea $\pi_1$ otro elemento primo, por lo que $\pi_1=\alpha \pi$ con $\alpha$
una unidad. Por lo tanto
\[
\frac{\sigma \pi_1}{\pi_1}=\frac{\sigma\pi}{\pi}\cdot\frac{\sigma\alpha}{\alpha}.
\]
Si $\sigma\in G_i$, $\sigma\alpha-\alpha\in\pL^{i+1}$ lo cual implica que
$\frac{\sigma\alpha}{\alpha}-1\in U^{(i+1)}$ lo cual demuestra que $\varphi$
no depende del elemento primo $\pi$.

Sean ahora $\sigma,\theta\in G_i$, entonces
\[
\frac{\sigma\theta(\pi)}{\pi}=\frac{\sigma\pi}{\pi}\cdot\frac{\theta\pi}{\pi}\cdot\frac
{\sigma u}{u},\quad u=\frac{\theta\pi}{\pi}.
\]
Puesto que $u$ es una unidad, la observaci\'on anterior muestra que $\frac{
\sigma u}{u}\in U^{(i+1)}$ y por lo tanto $\varphi(\sigma\theta)=
\varphi(\sigma)\varphi(\theta)$ y $\varphi$ es un homomorfismo de grupos.
Finalmente, si $\varphi(\sigma)=1$, entonces $\frac{\sigma\pi}{\pi}\in U^{(i+1)}$.
Por la Proposici\'on \ref{P1.3.13} se sigue que $\sigma\in G_{i+1}$ por lo que
$\varphi$ es inyectiva. $\fin$
\end{proof}

Notemos que la Proposici\'on \ref{P1.3.14} prueba que $G_i/G_{i+1}$ es
elemental abeliano para $i\geq 0$ y por tanto $G_1$ es un $p$--grupo, donde
$p$ es la caracter{\'\i}stica de los campos residuales $\AE L/\pL$ y $\AE K/\pK$,
y que $G_0/G_1\subseteq (\AE L/\pL)^{\ast}$, es decir, $G_0/G_1$ es un grupo
c{\'\i}clico de orden primo relativo a $|({\mc O}_L/\pL)^{\ast}|$
y por tanto de orden primo relativo a $p$. En particular

\begin{corolario}\label{C1.3.14'} Sea $p$ la caracter\'istica del campo residual
${\mc O}_L/\pL$. Se tiene que $G_{-1}=D(\pL|\pK)$ 
es un grupo soluble, $G_0/G_1$
es un grupo c{\'\i}clico de orden primo relativo a $p$ 
y $G_1$ es un $p$--grupo. Si $\pL|\pK$ es moderadamente
ramificado, entonces $G_0$ es c{\'\i}clico. Finalmente $G_i/G_{i+1}$ es un
$p$--grupo elemental abeliano para $i\geq 1$. $\fin$.
\end{corolario}

Un resultado que necesitaremos para demostrar el Teorema de Kronecker--Weber es
el siguiente.

\begin{proposicion}\label{P1.3.16} Sea $K/{\ma Q}$ una extensi\'on finita de
Galois con grupo $G$.
Supongamos que $G_{-1}/G_1$ es un grupo
abeliano. Entonces si $\varphi\colon G_0/G_1\to U^{(0)}/U^{(1)}\cong (\AE K/\pK)^{\ast}$
es el mapeo dado en la Proposici\'on \rm{\ref{P1.3.14}}, se tiene que $\im \varphi \subseteq
({\ma Z}/p{\ma Z})^{\ast}$, donde $\pK\cap {\ma Z}=(p)$.
\end{proposicion}

\begin{proof}
Sean $\sigma\in G_0$ y 
$\varphi(\bar{\sigma})=\alpha\in(\AE K/\pK)^{\ast}$, 
es decir, $\varphi(\bar{\sigma})=
\frac{\sigma\pi}{\pi}\bmod \pK = \alpha$, esto es,
$\sigma \pi=\alpha\pi\bmod \pK$. Sea $\theta\in G_{-1}$ arbitrario y sea
$\pi_1:=\theta^{-1}(\pi)$, el cual es un 
elemento primo para $\pK$, entonces se tiene
\begin{gather*}
\sigma\theta^{-1}(\pi)\equiv \alpha\theta^{-1}(\pi) \bmod \pK.\\
\intertext{Por ser $G_{-1}/G_1$ abeliano tenemos}
(\theta\sigma\theta^{-1})(\pi)\equiv \sigma\pi\equiv \theta(\alpha)\pi\bmod \pK
\equiv \alpha\pi\bmod \pK.
\end{gather*}
Por tanto $\theta(\alpha)\equiv \alpha\bmod\pK$ para toda $\theta\in G_{-1}/G_1$.
Se sigue que $\alpha$ es invariante bajo $G_{-1}/G_0\cong \Gal
\big((\AE K/\pK):({\ma Z}/p{\ma Z})\big)$
y por tanto $\alpha \in {\ma Z}/p{\ma Z}$. $\fin$
\end{proof}

Regresaremos a los grupos de ramificaci\'on en la Subseci\'on
\ref{CClaseS3.2.4} cuando estudemos teor\'ia de campos de clase locales.

\section{Primos infinitos, bases normales y campo de clase de
Hilbert}\label{S0.4}

En esta secci\'on recordamos algunos resultados no relacionados
entre si, pero que encontraremos a lo largo del libro.

Sea $K$ un campo num\'erico, $[K:{\ma Q}]=n=r_1+2r_2$. Sean
\begin{gather*}
\sigma_1,\ldots, \sigma_{r_1}, \sigma_{r_1+1},\ldots, \sigma_{
r_1+r_2},\overline{\sigma}_{r_1+1},\ldots, \overline{\sigma}_{
r_1+r_2}\\
\intertext{los $r_1$ encajes reales de $K$ y los $2r_2$ encajes
complejos:}
\begin{alignat*}{3}
\sigma_i&\colon K \longto &\ {\ma R}, &\quad 1\leq i\leq r_1,\\
\sigma_{r_1+j}&\colon K\longto &\ {\ma C}, &\quad 1\leq j\leq r_2,
\quad \sigma_{r_1+j}(K) \nsubseteq {\ma R}.
\end{alignat*}
\end{gather*}

Sea $|\ |$ el valor absoluto usual de ${\ma C}$. Definimos
los siguientes valores absolutos definidos sobre $K$:
\[
|x|_{\sigma_i}:= |\sigma_i x|,\quad 1\leq i\leq r_1+r_2, \quad
\text{donde}\quad |\sigma_{r_1+j} x|=|\overline{\sigma}_{r_1+j} x|.
\]

\begin{definicion}\label{D9.1} Los valores absolutos
$\big\{|\ |_{\sigma_i}\big\}_{1\leq i\leq r_1+r_2}$ son los
{\em primos infinitos\index{primos infinitos}} de $K$. Los
valores absolutos $\big\{|\ |_{\sigma_i}\big\}_{1\leq i\leq r_1}$
son los {\em primos infinitos reales\index{primos infinitos reales}}
y $\big\{|\ |_{\sigma_1+j}\big\}_{1\leq j\leq r_2}$ son los
{\em primos infinitos complejos\index{primos infinitos complejos}}.
\end{definicion}

En ${\ma Q}$ \'unicamente existe un primo infinito, el cual es real, y
corresponde al valor absoluto usual. Para $[K:{\ma Q}]=n=
r_1+2 r_2$ se tienen $r_1+r_2$ primos infinitos, $r_1$ reales
y $r_2$ complejos.

Notemos que esta definici\'on de hecho generaliza el concepto
de primo, pues si $\pK$ es un primo de ${\cal O}_K$, a $\pK$
le podemos asociar su valuaci\'on: si $x\in K^{\ast}$, $x{\cal O}_K
=\langle x\rangle = \pK^n {\eu a}$ con $(\pK,{\eu a})=1$, $n\in{\ma
Z}$. Entonces $v_{\pK}(x)=n$ y definimos el {\em valor
absoluto ${\eu p}$--\'adico\index{valor absoluto ${\eu p}$--\'adico}}:
\[
|x|_{\pK}:=p^{-fn} \quad \text{donde}\quad \pK\cap {\ma Q}=
\langle p\rangle, \quad n_{L/{\ma Q}}\pK=p^f,
\]
es decir, donde $f$ es el grado de inercia.

En otras palabras, $|x|_{\pK}=p^{-fn}=(N_{K/{\ma Q}}\pK)^{
v_{\pK}(x)}$. Se define $|0|_{\pK}=0$. Este valor absoluto es no
{\em arquimediano\index{valor absoluto no arquimediano}}:
\[
|x+y|_{\pK}\leq \max\{|x|_{\pK},|y|_{\pK}\}
\]
y los valores absolutos $|\ |_{\sigma_i}$ son arquimedianos.
Resumiendo, podemos pensar en``primos'' de $K$ como un
valor absoluto de $K$.

Ahora consideremos $L/K$ una extensi\'on de campos num\'ericos.
Sea $\omega\colon L\to {\ma C}$ un encaje de $L$ y sea $
\omega|_K=
\sigma$, $\sigma\colon K\to {\ma C}$ es un encaje de $K$. 
Notemos que si $\omega$ es real, $\sigma$ necesariamente
es real, pero si $\omega$ es complejo, entonces $\sigma$ puede
ser real o complejo.

\begin{definicion}\label{D9.2}
Con la notaci\'on anterior, decimos que $\omega$ (o $\sigma$)
es {\em ramificado} si $\sigma$ es real y $\omega$ es complejo
y definimos que el {\'\i}ndice de ramificaci\'on como $2$: $e(
\omega|\sigma)=2$. Esto se hace pensando en que $[{\ma C}:
{\ma R}]=2$ o en que $\omega|_K=\overline{\omega}|_K=\sigma$.

En el caso de que $\omega$ y $\sigma$ sean ambos reales o
ambos complejos, entonces definimos el {\'\i}ndice de
ramificaci\'on como $1$: $e(\omega|\sigma)=1$.

En cualquier caso, el grado de inercia lo definimos como $1$:
$f(\omega|\sigma)=1$ siempre.
\end{definicion}

Se tiene $e(\omega|\sigma)f(\omega|\sigma)=
\begin{cases} 2& \text{$2$ si $\omega$ es complejo y $\sigma$
es real}\\ $1$& \text{en otro caso}.\end{cases}$.

En particular, si fijamos $\sigma\colon K\to {\ma C}$ un encaje.
$\sigma$ tiene $[L:K]$ extensiones a encajes $\omega\colon
L\to {\ma C}$. Si $\omega$ y $\overline{\omega}$ son dos
complejos conjugados de estos encajes, $\omega$ y $\overline{\omega}$
los consideramos los mismos y denotamos a cualquier extensi\'on
$\omega$ por $\omega | \sigma$.

Entonces 
\begin{gather*}
\sum_{\omega|\sigma} e(\omega|\sigma) f(\omega|
\sigma)=[L:K].\\
\intertext{Esta f\'ormula es exactamente la misma f\'ormula
que para los primos finitos:}
\sum_{\pL|\pK} e(\pL|\pK)f(\pL|\pK)=
[L:K]
\end{gather*}
 donde $\pK$ es un primo de ${\cal O}_K$ y $\pL|\pK$ 
recorre los primos de ${\cal O}_L$ sobre $\pK$.

A continuaci\'on, estudiamos bases normales de una
extensi\'on finita de Galois de campos arbitrarios. La 
demostraci\'on que presentamos del teorema de la base
normal se debe a \cite{Kim}.

\begin{teorema}[de la base normal]\label{T1.4.3} 
Sea $L/K$ una extensi\'on finita de Galois de campos
arbitrarios. Sea $G=\Gal(L/K)$. Entonces existe $\alpha\in
L$ tal que $\{\alpha^{\sigma}\}_{\sigma\in G}$ es
base de $L/K$. En particular $L\cong K[G]$ como
$G$--m\'odulos.
\end{teorema}

\begin{proof}
(1) Sea $K$ finito. Entonces $G=\langle\sigma\rangle\cong
C_n$. Consideremos a $\sigma$ como un endomorfismo 
$K$--lineal del $K$--espacio vectorial $L$. M\'as precisamente,
$\sigma\colon L\lra L$ es una transformaci\'on lineal del
$K$--espacio vectorial $L$.

Ahora $\{1,\sigma,\ldots, \sigma^{n-1}\}$
es un conjunto linealmente independientes sobre $L$
por el teorema
de la independencia de Artin (ver Teorema \ref{CClaseT1.1.1}). 
En particular, el 
polinomio m\'inimo $m(x)$ de $\sigma$, considerado como transformaci\'on
lineal, es de grado $n$. Entonces, como $L$ es un 
$K[x]$--m\'odulo finitamente generado con 
acci\'on $x\circ a :=\sigma a$,
$a\in L$ y $K[x]$ es un dominio de ideales principales,
se tiene una descomposici\'on
\begin{gather*}
L\cong M_1\oplus \cdots \oplus M_r\cong
K[x]\alpha_1\oplus \cdots \oplus K[x]\alpha_r,\\
\text{con}\quad \an(\alpha_i)|\an(\alpha_{i+1}),
\quad 1\leq i\leq r-1\quad\text{y}\quad \an(\alpha_r)
=m(x),
\end{gather*}
donde $\an(\alpha)$ denota al anulador de $\alpha$.
Por lo tanto se tiene $L\cong M\alpha_r=M\alpha_1\cong
K[x]/\an(\alpha_1)=K[x]/m(x)$ y $\theta =x+m(x)$
satisface $L=K[\theta]$.

Notemos que esta demostraci\'on es aplicable a cualquier
extensi\'on c\'iclica $L/K$, independientemente de que $K$
sea finito o no.

(2) Ahora sea $K$ infinito. Por el teorema del elemento primitivo,
$L=K(a)$ para alg\'un $a\in L$. Sean $f(x)=\Irr (a,x, K)$ y
$G=\{\sigma_1,\ldots,\sigma_n\}$, $\deg f(x)=n$. Sean
$a_i=\sigma_i(a)$, $g(x)=\frac{f(x)}{(x-a)f'(a)}$ y $
\sigma_i(g(x))=\frac{f(x)}{(x-a_i)f'(a_i)}=g_i(x)$.

Para $i\neq j$ se tiene $g_i(x)g_j(x)\equiv 0\bmod f(x)$
pues $x-a_i\neq x-a_j$. Sea $f(x)=(x-a_i)h(x)$,
$f'(x)=h(x)+(x-a_i)h'(x)$, $f'(a_i)=h(a_i)$. Se sigue que
$g_i(a_i)=1$. Ahora bien, para $i\neq j$, $g_i(a_j)=
\frac{f(a_j)}{(a_j-a_i) f'(a_i)}=0$. En resumen
$g_i(a_j)=\delta_{ij}$ es la delta de Kronecker. Sea
$t(x)=g_1(x)+\cdots+g_n(x)-1$. Entonces
$\deg(t(x))\leq n-1$ pues $\deg g_i(x)=n-1$ para
$1\leq i\leq n$ y $t(a_1)=\cdots=t(a_n)=0$. Se sigue
que $t(x)=0$, es decir, $g_1(x)+\cdots+g_n(x)-1=0$.

Ahora $g_i(x)=g_i(x)(g_1(x)+\cdots+g_n(x))=g_i(x)^2+
\sum_{j\neq i} g_i(x)g_j(x)\equiv g_i^2(x) + 0 \bmod f(x)$.
Por tanto $g_i^2(x)\equiv g_i(x) \bmod f(x)$.

Sean $A(x)=(b_{ij}(x))_{\substack{1\leq i\leq n\\ 1\leq j\leq n}}$
y $D(x)=\det(A(x))$ donde $b_{ij}(x)=
\sigma_i\sigma_j(g(x))$. Sea $A^2(x)=(c_{ij}(x)
)_{\substack{1\leq i\leq n\\ 1\leq j\leq n}}$. Entonces
$c_{ij}(x)=\sum_{k=1}^n b_{ik}(x)b_{kj}(x)$. Ahora bien
\begin{gather*}
b_{ij}(x)=\sigma_i\sigma_j(g(x))=\sigma_i(g_j(x))=
\sigma_i\Big(\frac{f(x)}{(x-a_j)f'(a_j)}\Big)=
\frac{f(x)}{(x-\sigma_i a_j) f'(\sigma_i a_j)},\\
b_{ij}(a_t)=\begin{cases}
0&\text{si $a_t\neq \sigma_i a_j$}\\
1&\text{si $a_t=\sigma_i a_j$}
\end{cases}.
\end{gather*}

Sea $\sigma_i a_j:=a_{m_i(j)}$. Entonces $b_{ij}(a_t)=\delta_{
t, m_i(j)}$. As\'i
\[
c_{ij}(a_t)=\sum_{k=1}^n b_{ik}(a_t)b_{kj}(a_t)=
\sum_{k=1}^n \delta_{t,m_i(k)}\delta_{t,m_k(j)}=
\begin{cases}
0&\text{si $m_i(k)\neq m_k(j)$}\\
1&\text{si $m_i(k)=m_k(j)$}
\end{cases}.
\]
Ahora bien $m_i(k)=m_k(j)$ si y s\'olo si $\sigma_i(a_k)=
\sigma_k(a_j)$ si y s\'olo si $\sigma_k^{-1}\sigma_i(a_k)=
a_j$.

Regresando a $A(x)=(b_{ij}(x))_{\substack{1\leq i\leq n\\ 1\leq j\leq n}}$,
se tiene $b_{ij}(a_t)=\delta_{t,m_i(j)}$. Esto es, cada fila y cada columna
de $A(a_t)$ tiene exactamente un $1$ y todos los dem\'as son
$0$. Fijando $j$, se tiene que si $\delta_{t,m_i(j)}=\delta_{t,m_{i'}(j)}=1$ 
entonces $t=m_i(j)=m_{i'}(j)$ y $\sigma_i(a_j)=m_i(j)=m_{i'}(j)=\sigma_{
i'}(a_j)$ lo cual implica que $\sigma_i^{-1}\sigma_{i'}(a_j)=a_j$
estos es, $\sigma_i^{-1}\sigma_{i'}$ deja fijo a $K(a_j)=
K(a)=L$. Por tanto $\sigma_{i}^{-1}\sigma_{i'}=\Id_L$, esto es,
$\sigma_{i'}=\sigma_i$ y finalmente $i=i'$.

De lo anterior obtenemos $\det A(a_t)=D(a_t)=\pm 1$ y por tanto
$ D(a_t)^2=1$ para toda $t\in\{1,\ldots,n\}$. Se
sigue que $D(x)^2-1$ tiene como ra\'ices a $a_1,\ldots,a_n$
y en particular $D(x)\neq 0$. Puesto que $K$ es infinito, existe
$\alpha\in K$ tal que $D(\alpha)\neq 0$. Sea $\theta=g(\alpha)$.
Entonces $D(\alpha)=\det(\sigma_i\sigma_j(g(\alpha)))=
\det (\sigma_i\sigma_j(\theta))\neq 0$.

Consideremos cualquier relaci\'on lineal $x_1
\sigma_1(\theta)+\cdots+x_n\sigma_n(\theta)=0$ con
$x_1,\ldots,x_n\in K$. Aplicando $\sigma_j$ para 
$j=1,\ldots,n$ se obtiene
\[
(\sigma_i\sigma_j(\theta))_{\substack{1\leq i\leq n\\ 1\leq j\leq n}}
\left(
\begin{array}{c} x_1\\ \vdots \\ x_n
\end{array}
\right)
=0=\left(
\begin{array}{c}
0\\ \vdots \\ 0
\end{array}
\right).
\]
Puesto que $\det(\sigma_i\sigma_j(\theta))_{\substack{
1\leq i\leq n\\ 1\leq j\leq n}}\neq 0$, $C=
(\sigma_i\sigma_j(\theta))_{\substack{1\leq i\leq n\\ 1\leq j\leq n}}$
es invertible y $\left(
\begin{array}{c} x_1\\ \vdots \\ x_n
\end{array}
\right)
= C^{-1}\left(
\begin{array}{c}
0\\ \vdots \\ 0
\end{array}
\right)
=\left(
\begin{array}{c}
0\\ \vdots \\ 0
\end{array}
\right)$ por lo cual $x_1=\cdots=x_n=0$.

Por tanto $\{\sigma_i(\theta)\}_{i=1}^n$ es un conjunto linealmente
independiente sobre $K$ y por tanto es base de $L/K$. $\fin$
\end{proof}

Terminamos esta secci\'on mencionando el Teorema de
Hilbert sobre la m\'axima extensi\'on abeliana no ramificada de 
un campo. Ver Corolario \ref {CClaseC4.5.13} y
Proposici\'on \ref{CClaseP4.7.10}

\begin{teorema}[Campo de clase de
Hilbert\index{Hilbert!campo de clase de $\sim$}]\label{T9.3}
Sea $K$ una extensi\'on finita de ${\ma Q}$ y sea $I_K$ su grupo de clases de
ideales de $K$. Sea $H_K$ la m\'axima extensi\'on abeliana
de $K$ no ramificada en ning\'un primo, finito o infinito. Entonces
$H_K$ es una extensi\'on finita y de Galois de $K$ con
grupo de Galois isomorfo a $I_K$: $\Gal(H_K/K)\cong I_K$
(ver Corolario {\rm{\ref{CClaseC4.5.13}}}).
$\fin$
\end{teorema}

\section{Aritm\'etica de extensiones separables finitas}\label{S1.5}

En esta secci\'on, todas las extensiones son separables y finitas. Los resultados
son v\'alidos tanto para campos num\'ericos como para campos de
funciones sobre un campo de constantes $k$ perfecto. Algunos de estos
resultados ser\'an probados en el Cap\'itulo \ref{ChRam} en el caso de
campos de funciones congruentes, esto es, con $k$ finito. La demostraci\'on
para $k$ perfecto y para campos num\'ericos es totalmente an\'aloga.
Hacemos notar que la mayor\'ia de los resultados siguen siendo v\'alidos
para $k$ arbitrario.

\begin{teorema}\label{T1.5.1} Sea $E/F$ una extensi\'on y 
sea $\pK$ un primo de $F$
ramificado en $E/F$, es decir $e_{E/F}({\mc P}|\pK) >1$ con $\pK={\mc P} \cap
F$ para alg\'un ${\mc P}$ un primo de $E$ sobre $\pK$. 
Sea $E=E_1E_2$ con $F\subseteq E_i\subseteq E$, $i=1,2$. Entonces
si ${\mc P}_i={\mc P}\cap E_i$, $i=1,2$,
se tiene $e_{E_1/F}({\mc P}_1|\pK) >1$ o
$e_{E_2/F}({\mc P}_2|\pK) >1$.

Equivalentemente, si $e_{E_i/F}({\mc P}_i|\pK) =1$ para todo 
${\mc P}_i$ en $E_i$ con ${\mc P}_i\cap F=\pK$ para $i=1,2$, entonces
$e_{E/F}({\mc P}|\pK) =1$ para todo ${\mc P}$ en $E$ con ${\mc P}
\cap F=\pK$.
En otras palabras, si un primo de $F$ es no ramificado en $E_i/F$,
$i=1,2$, entonces es no ramificado en $E_1E_2/F$.
\end{teorema}

\begin{proof}
Sea $\tilde E$ la cerradura de Galois de $E/F$. Supongamos que 
$e_{E_i/F}({\mc P}_i|\pK) =1$ para todo primo ${\mc P}_i$ en
$E_i$ sobre $\pK$, $i=1,2$. Sea $I=I(\tilde {\mc P}|\pK)$
el grupo de inercia de $\tilde{\mc P}$ en $\tilde E/F$ con $\tilde{\mc P}
\cap F=\pK$. Queremos probar que $e_{\tilde E/F}(\tilde {\mc P}|\pK) =1$.

Sean $M:=\tilde E^I$, $\pL=\tilde{\mc P}\cap M$, ${\mc P}=\tilde{\mc P}
\cap E$ y ${\mc P}_i=\tilde{\mc P}\cap E_i$, $i=1,2$.
\[
\xymatrix{
&&\tilde E\ar@{--}[r]\ar@{-}[dl]\ar@{-}[dd]&\tilde{\mc P}
\ar@{--}[dddlll]\ar@{--}[dddr]\ar@/^6pc/@{--}[dddd]\ar@/_2pc/@{--}[dlll]\\
\pL\ar@{--}[r]&M\\
&&E\ar@{-}[dl]\ar@{-}[dr]\ar@{--}[r]&{\mc P}\ar@{--}[uu]\\
{\mc P}_1\ar@{--}[r]&E_1\ar@{-}[dr]&&E_2\ar@{--}[r]\ar@{-}[dl]&{\mc P}_2\\
&&F\ar@{--}[r]&\pK
}
\]

Se tiene $e_{\tilde E/M}(\tilde{\mc P}|\pL)=[\tilde E:M]=|I|$ y $e_{M/F}
(\pL|\pK)=1$. Sea $R=ME_1$. Se tiene $M\subseteq R\subseteq \tilde E$.
Entonces $I_{\tilde E/E_1}(\tilde{\mc P}|{\mc P}_1)\subseteq I_{\tilde E/F}(
\tilde{\mc P}|\pK)=I$. Por tanto $\tilde E^I=M\supseteq \tilde E^{I_{\tilde E/E_1}
(\tilde{\mc P}|{\mc P}_1)}$ y $I_{\tilde E/E_1}(\tilde{\mc P}|{\mc P}_1)\subseteq
\Gal(\tilde E/E_1)$. Por tanto $\tilde E^{I_{\tilde E/E_1}(\tilde{\mc P}|{\mc P}_1)}
\supseteq \tilde E^{\Gal(\tilde E/E_1)}=E_1$ de donde se sigue que
$E_1\subseteq M$. Similarmente $E_2\subseteq M$ y por tanto $E=E_1E_2
\subseteq M$. Se sigue que $e_{E_1E_2/F}({\mc P}|\pK)|e_{M/F}(
\pL|\pK)=1$.
$\fin$
\end{proof}

\begin{corolario}\label{C1.5.2} Con las notaciones del Teorema {\rm{\ref{T1.5.1}}},
si $F\subseteq S\subseteq E$ y $e_{S/F}(\tilde{\mc P}\cap S|\pK)=1$,
entonces $S\subseteq M$. $\fin$
\end{corolario}

\begin{teorema}\label{T1.5.3} Con las notaciones del Teorema {\rm{\ref{T1.5.1}}},
si $e_{E/F}({\mc P}|\pK)=1$ para todo primo ${\mc P}$ de $E$ sobre $\pK$,
entonces $e_{\tilde E/F}(\tilde{\mc P}|\pK)=1$. En otras palabras, si $\pK$ es
no ramificado en $E/F$, entonces $\pK$ es no ramificado en $\tilde E/F$.
\end{teorema}

\begin{proof}
Ver la demostraci\'on del Teorema \ref{T.Ram10} (1).
$\fin$
\end{proof}

\begin{teorema}\label{T1.5.4} Con las notaciones del Teorema {\rm{\ref{T1.5.1}}}
se tiene que $\pK$ es totalmente descompuesto en $E/F$ si y solamente si
$\pK$ es totalmente descompuesto en $\tilde E/F$.
\end{teorema}

\begin{proof}
Ver la demostraci\'on del Teorema \ref{T.Ram10} (2).
$\fin$
\end{proof}

\begin{proposicion}\label{P1.5.5} 
Sean $E/F$ y $L/F$ dos extensiones separables finitas
donde el primo $\pK$ de $F$ es totalmente descompuesto.
Entonces $\pK$ es totalmente descompuesto en $EL/F$.
\end{proposicion}

\begin{proof}
Ver la demostraci\'on de la Proposici\'on \ref{P.Ram10(1)}.
$\fin$
\end{proof}

\begin{proposicion}\label{P1.5.6}
Sea $E/F$ una extensi\'on finita de Galois.
Sea $K/F$ una extensi\'on
arbitraria finita y sea $L=EK$. Sea $\pK_L$ un divisor
primo en $L$ y sean $\pK_E, \pK_F$ y $\pK_K$
sus restricciones a $E, F$ y $K$ respectivamente.
Sea $\rest\colon \Gal(L/K)\lra \Gal(E/F)$ el monomorfismo
de restricci\'on. Sean
$D$ e $I$ los grupos de descomposici\'on
y de inercia respectivamente. Entonces
\las
\item $I(\pK_L|\pK_K)|_E \subseteq I(\pK_E|\pK_F)$,
\item $D(\pK_L|\pK_K)|_E \subseteq D(\pK_E|\pK_F)$.
\end{list}
\end{proposicion}

\begin{proof}
Ver la demostraci\'on de la Proposici\'on \ref{P10.4.1.Ram}.
\end{proof}

\begin{corolario}\label{C1.5.7}
Con las hip\'otesis de la Proposici\'on {\rm{\ref{P1.5.6}}},
se tiene:
\las
\item Si $\pK_F$ se descompone totalmente en $E/F$,
entonces $\pK_K$ se descompone totalmente en $L/K$.

\item Si $\pK_F$ es no ramificado en $E/F$, entonces 
$\pK_K$ es no ramificado en $L/K$. Equivalentemente
si $\pK_K$ es ramificado en $L/K$, entonces $\pK_F$
es ramificado en $E/F$.
\end{list}
\end{corolario}

\begin{proof}
Ver el Corolario \ref{C.10.4.2.Ram}.
$\fin$
\end{proof}

\begin{corolario}\label{C1.5.8}
Sean $E/F$ y $L/F$ dos extensiones finitas y separables.
Si el primo ${\mc P}$ de $L$ es ramificado en $EL/L$,
entonces $\pK:={\mc P}\cap F$ es ramificado en $E/F$.
\end{corolario}

\begin{proof}
\[
\xymatrix{
&E\ar@{-}[r]\ar@{-}[d]&EL\ar@{-}[d]\\
&F\ar@{-}[r]&L\ar@{--}[rd]\\
\pK={\mc P}\cap F\ar@{--}[ru]\ar@{--}[rrr]&&&{\mc P}
}
\]
Sea $\tilde E$ la cerradura de Galois $E/F$. Entonces $\tilde EL$
es Galois sobre $L$ y ${\mc P}$ es ramificado en $\tilde EL/L$ lo
que implica que $\pK$ es ramificado en $\tilde E/F$. Por tanto
$\pK$ es ramificado en $E/F$.
$\fin$
\end{proof}

\begin{corolario}\label{C1.5.9}
Con las notaciones del Corolario {\rm{\ref{C1.5.8}}}, si $\pK$
es totalmente descompuesto en $E/F$, entonces ${\mc P}$
es totalmente descompuesto en $EL/L$.
\end{corolario}

\begin{proof}
Por el Teorema \ref{T1.5.4}, $\pK$ es totalmente descompuesto
en $\tilde E/F$. Por el Corolario \ref{C1.5.7}, ${\mc P}$ es totalmente
descompuesto en $\tilde EL/L$ por lo que ${\mc P}$ es totalmente
descompuesto en $EL/L$.
$\fin$
\end{proof}

%% file: Capitulo2.tex
\chapter{Teor{\'\i}a de Galois infinita}\label{ch2}

\section{L{\'\i}mites directos y l{\'\i}mites inversos}\label{S2.1}

Primero recordemos algunos conceptos generales.
Sea $A$ un anillo conmutativo
con unidad. 

\begin{definicion}\label{D5.1} Un {\em conjunto dirigido\index{conjunto dirigido}}
$I$ es un conjunto no vac{\'\i}o parcialmente ordenado tal que
para cualesquiera $i,j\in I$ existe $k\in I$ tal que $i\leq k$ y 
$j\leq k$.
\end{definicion}

Sean $I$ un conjunto dirigido y $\{M_i\}_{i\in I}$ un conjunto de
$A$--m\'odulos. Si para cualesquiera $i,j\in I$ con $i\leq j$ existe
$\varphi_{i,j}\colon M_i\to M_j$ un $A$--homomorfismo tal que
\l
\item $\varphi_{ii}=\Id_{M_i}$ para toda $i\in I$,
\item $\varphi_{ik}=\varphi_{jk}\circ \varphi_{ij}$ para cualesquiera
$i\leq j\leq k$,
\end{list}
entonces decimos que $\{M_i,\varphi_{ij}, I\}_{\substack{i,j\in I\\
i\leq j}}$ forman un {\em sistema directo\index{sistema directo}}
sobre $I$.

\begin{definicion}\label{D5.2} Si $\{M_i,
\varphi_{ij}, I\}_{\substack{i,j\in I\\ i\leq j}}$ es un sistema directo,
el {\em l{\'\i}mite directo\index{l{\'\i}mite directo}} se define como
el $A$--m\'odulo 
\[
\lim_{\substack{\longrightarrow\\ i\in I}}M_i=:M
\]
definido por $M=P/N$ donde $P=\oplus_{i\in I}M_i$ y $N=\langle
m_i-\varphi_{ij}(m_i)\mid i\in I, m_i\in M_i, i\leq j\rangle$.
\end{definicion}

Sea $h_i\colon M_i\to M$ dado por $h_i=
\pi\circ \mu_i$ donde $\mu_i\colon M_i\to P$ es el encaje
natural y $\pi\colon P\to P/N$ es la proyecci\'on natural.
En otras palabras
\[
h_i(x):=(\xi_j)_{j\in I}\bmod N \quad{\text{donde}}\quad
\xi_j=\left\{
\begin{array}{ccc}0&{\text{si}}&j\neq i\\
x&{\text{si}}&j=i.
\end{array}
\right.
\]

Se tiene que $h_i=h_j\circ \varphi_{ij}$ para $i\leq j$, $i,j\in I$
pues $m_i\equiv \varphi_{ij}(m_i)\bmod N$ para toda $m_i\in M_i$.

Enunciamos en el siguiente resultado todas las propiedades
que necesitaremos a lo largo de este libro.

\begin{teorema}\label{T5.3}
Sea $I$ un conjunto dirigido, $\{M_i,
\varphi_{ij}, I\}_{\substack{i,j\in I\\ i\leq j}}$ un sistema directo
y $M=\lim\limits_{\longrightarrow} M_i$, el l{\'\i}mite directo.
Sean $h_i\colon M_i\to M$ los homomorfismos naturales.
Entonces
\las
\item Dado $m\in M$, existen $i\in I$ y $x\in M_i$ tales que
$m=h_i(x)$.

\item Todo elemento que es $0$ en $M$, es que era 
eventualmente $0$ en los $M_i$'s. Esto es, si $h_i(x)=0$
con $x\in M_i$, entonces existe $j\geq i$ tal que $\varphi_{
ij}(x)=0\in M_j$.

\item $M$ satisface la siguiente propiedad universal: Sea $N$\
un $A$--m\'odulo tal que para toda $i\in I$ existe $\theta_i\colon
M_i\to N$ un homomorfismo de $A$--m\'odulos tal que
\begin{gather*}
\xymatrix{
M_i\ar[rr]^{\varphi_{ij}}\ar[dr]_{\theta_i}&&M_j\ar[ld]^{\theta_j}\\
&N}\quad \theta_j\circ \varphi_{ij}=\theta_i \quad {\text{para toda}}
\quad i\leq j.\\
\intertext{Entonces existe un \'unico homomorfismo $\theta\colon M\to N$
tal que}
\xymatrix{
M_i\ar[rr]^{h_i}\ar[rd]_{\theta_i}&&M\ar[dl]^{\theta}\\ &N}
\quad \theta\circ h_i=\theta_i \quad {\text{para toda}}\quad i\in I.
\end{gather*}

\item $M$ est\'a caracterizado por la propiedad universal
de {\rm (3)}.

\item Si $\{M_i\}_{i\in I}$ es una familia de subm\'odulos de un
$A$--m\'odulo $N$ tal que para cualesquiera $i,j\in I$ existe 
$k\in I$ tal que $M_i+M_j\subseteq M_k$. Se define $i\leq j$
en $I$ si y s\'olo si $M_i\subseteq M_j$. Entonces $I$ es un
sistema dirigido y se define $\varphi_{ij}\colon M_i\to M_j$ como
el encaje natural cuando $i\leq j$. Entonces
\[
M=\lim_{\substack{\longrightarrow\\ i\in I}} M_i=
\sum_{i\in I}M_i = \bigcup_{i\in I} M_i
\]
donde $h_i:M_i\to M$ es el encaje natural.

\end{list}
\end{teorema}

\begin{proof}
\cite[pags. 32--33]{AtiMac69}. $\fin$
\end{proof}

Estamos interesados en el caso especial de una extensi\'on de
Galois de campos. Podemos considerar a los campos como
${\ma Z}$--m\'odulos. Sea $L/K$ una extensi\'on cualquiera de
Galois, no necesariamente finita. Para cada $\alpha\in L$,
$K(\alpha)/K$ es una extensi\'on finita. Entonces
$L=\bigcup\limits_{\alpha
\in L}K(\alpha)=\lim\limits_{\substack{\longrightarrow\\ \alpha\in L}}
K(\alpha)$, donde los homomorfismos est\'an dados por 
el encaje natural $K(\alpha)\to K(\beta)$ para $K(\alpha)\subseteq
K(\beta)$.

Notemos que con esta notaci\'on, el conjunto dirigido es
$I:=\{K(\alpha)\mid \alpha\in L\}$ donde se define $K(\alpha)\leq
K(\beta)\iff K(\alpha)\subseteq K(\beta)$. Veamos cual es la
conexi\'on con los grupos de Galois.

Para este fin, consideremos ahora {\em l{\'\i}mites
inversos\index{l{\'\i}mites inversos}}. Sea $I$ un conjunto dirigido.
Para cada $i\in I$ consideremos un objeto $A_i$, el cual puede
ser un grupo, un espacio topol\'ogico, un anillo, un m\'odulo, un conjunto,
etc. Nuestro inter\'es ser\'a fundamentalmente en considerar
$A_i$ un grupo finito con la topolog{\'\i}a discreta.

\begin{definicion}\label{D5.2'}
Un {\em sistema inverso} $\{A_i,\phi_{ji},I\}_{\substack{i,j\in I\\ i\leq j}}$
es un sistema tal que $\phi_{ji}\colon A_j\to A_i$ son morfismos para $i\leq j$
tales que:
\las
\item $\phi_{ii}=\Id_{A_i}$ para toda $i\in I$.

\item $\phi_{ji}\circ \phi_{kj}=\phi_{ki}\quad {\text{para}}\quad i\leq j\leq k$,
\[
\xymatrix{
A_k\ar[rr]^{\phi_{kj}}\ar[dr]_{\phi_{ki}}&&A_j\ar[dl]^{\phi_{ji}}\\ & A_i}
\]
(por {\em morfismo\index{morfismo}} entendemos homomorfismo 
si los objetos son grupos, anillos, campos, m\'odulos o cualquier
otra estructura algebraica; funci\'on continua si los objetos son
espacios topol\'ogicos, homomorfismo continuo si los objetos son
grupos o anillos topol\'ogicos y simples funciones si los $A_i$'s
son simplemente conjuntos).
\end{list}
\end{definicion}

Si $\{A_i,\phi_{ji},I\}_{\substack{i,j\in I\\ i\leq j}}$ es un sistema
inverso, entonces decimos que $(X,\varphi_i)_{i\in I}$ es un
{\em l{\'\i}mite inverso\index{l{\'\i}mite inverso}} del sistema si
$\varphi_i\colon X\to A_i$ son morfismos tales que
\[
\phi_{ji}\circ \varphi_j=\varphi_i\quad{\text{para}}\quad i\leq j
\quad
\xymatrix{X\ar[rr]^{\varphi_j}\ar[rd]_{\varphi_i}&&A_j\ar[dl]^{\phi_{ji}}\\
&A_i}
\]
y si $(Y,\mu_i)_{i\in I}$ es otro sistema tal que $\mu_i\colon
Y\to A_i$ son morfismos tales que
\begin{gather*}
\xymatrix{
Y\ar[rr]^{\mu_j}\ar[rd]_{\mu_i}&&A_j\ar[ld]^{\varphi_{ij}}\\&A_i}
\quad \varphi_{ij}\circ \mu_j=\mu_i\quad{\text{para toda}} \quad i\leq j\\
\intertext{entonces existe un \'unico morfismo $\xi\colon Y\to X$ tal que}
\xymatrix{
Y\ar[rr]^{\xi}\ar[dr]_{\mu_i}&&X\ar[dl]^{\varphi_i}\\&A_i}
\quad \varphi_i\circ \xi =\mu_i \quad {\text{para toda}}\quad i\in I.\\
\intertext{Escribimos}
X=\lim_{\substack{\longleftarrow\\ i\in I}}A_i=\lim_{\substack{
\longleftarrow\\ i}}A_i=\lim_{\longleftarrow} A_i.
\end{gather*}

\begin{teorema}\label{T5.4}
Dado un sistema inverso $\{A_i,\phi_{ji},I\}_{\substack{i,j\in I\\ i\leq j}}$
existe el l{\'\i}mite inverso $(X,\varphi_i)_{i\in I}$, $X=\lim\limits_{
\substack{\longleftarrow\\ i\in I}} A_i$. Se tiene que $(X,\varphi_i)_{
i\in I}$ es \'unico salvo isomorfismo.

M\'as a\'un $(X,\varphi_i)_{i\in I}$ se puede realizar como las
``{\em sucesiones coherentes\index{sucesiones coherentes}}'' de
$B:=\prod_{i\in I}A_i$, el producto directo, es decir
\[
X=\{(a_i)_{i\in I}\in B\mid a_i=\phi_{ji}(a_j){\text{\ para toda\ }}
i,j\in I, i\leq j\}
\]
donde $\phi_i\colon X\to A_i$ es la $i$--\'esima proyecci\'on.
\end{teorema}

\begin{proof} Se puede consultar la demostraci\'on en
\cite[Cap{\'\i}tulo 11]{Vil2006} o en \cite{RibZal2000}. 
Aqu{\'\i} la volvemos a presentar.

Primero veamos la unicidad. Si $(Z,\theta_i)_{i\in J}$ es otro l{\'\i}mite
inverso, entonces existe mapeos \'unicos $\alpha\colon X\to Z$,
$\beta\colon Z\to X$ tales que los siguientes diagramas 
conmutan:
\[
\theta_i\circ \alpha =\varphi_i\qquad
\xymatrix{
X\ar[r]^{\alpha}\ar[dr]_{\varphi_i}&Z\ar[r]^{\beta}\ar[d]^{\theta_i}
&X\ar[dl]^{\varphi_i}\\&A_i
}
\qquad \varphi_i\circ \beta=\theta_i
\]

Entonces $\beta\circ \alpha$ y $\Id_X$ satisfacen que $\varphi_i\circ
(\beta\circ \alpha)=\varphi_i=\varphi_i\circ(\Id_X)$. Por la unicidad
tenemos que $\beta\circ \alpha=\Id_X$. An\'alogamente tenemos
$\alpha\circ \beta=\Id_Z$. Esto prueba que $\alpha$ y $\beta$ son
isomorfismos inversos uno del otro entre $X$ y $Z$.

Para ver la existencia, sea $B:=\prod\limits_{i\in I}A_i$. Se define
$X:=\{(a_i)_{i\in I}\in B\mid \phi_{kj}(a_k)=a_j\text{\ si\ }j\leq k\}$ ($B$
se considera con las operaciones entrada por entrada en el caso
de grupos, anillos, campos, m\'odulos, etc. o con la topolog{\'\i}a
producto en el caso de espacios topol\'ogicos). Sea $\pi_i\colon B
\to A_i$ la proyecci\'on y sea $\varphi_i\colon X\to A_i$, $\varphi_i:
=\pi_i|_X$. Se tiene que
\[
\big(\phi_{ji}\circ \varphi_j\big)\big((a_k)_{k\in I}\big)=\phi_{ji}
(a_j)=a_i=\varphi_i\big((a_k)_{k\in I}\big)
\]
para todo $(a_k)_{k\in I}\in X$. Si $(Y,\xi_i)_{i\in I}$ es tal que $\xi_i
\colon Y\to A_i$ satisface $\phi_{ji}\circ \xi_j=\xi_i$ para $i\leq j$, sea
$\xi\colon Y\to X$ dada por $\xi(y):=\big(\xi_i(y)\big)_{i\in I}$.
Entonces $\xi$ est\'a bien definido puesto que $(\varphi_i\circ \xi)(y)=
\varphi_i\big((\xi_k(y)_{k\in I}\big)$ de tal forma que $\xi(y)\in X$. 
pues $\phi_{ji}(\xi_j(x))=\xi_i(x)$. Por
tanto $X$ es un l{\'\i}mite inverso de $\big\{A_i,\phi_{ji}, I\big\}$.
$\fin$
\end{proof}

\begin{observacion}\label{O7.1.5} Si para cada $i\in I$, $A_i$
es un espacio topol\'ogico Hausdorff, damos a $A:=\prod\limits_{i\in I}
A_i$ la topolog{\'\i}a producto y $\lim\limits_{\leftarrow} A_i$ es un
espacio topol\'ogico con la topolog{\'\i}a inducida. Siempre
supondremos que los mapeos $\phi_{ji}$ son continuos. Ahora, las
funciones $\phi_i$ son siempre continuas pues si $U$ es un 
abierto de $A_i$, tenemos que $\phi_i^{-1}(U)=\pi_i^{-1}(U)\cap
\lim\limits_{\leftarrow}A_i$ y $\pi_i^{-1}(U)$ es un conjunto abierto en
$A$ por la definici\'on de la topolog{\'\i}a producto.

Se tiene mucho m\'as. Si $V$ es un conjunto abierto de $X=\lim
\limits_{\leftarrow}A_i$, veremos que $V$ contiene a alg\'un conjunto
de la forma $\phi_k^{-1}(U_k)$ para alg\'un conjunto abierto $U_k$
de $A_k$ y alg\'un $k\in I$. Se tiene que $V$ est\'a generado
por uniones e intersecciones finitas de conjuntos de la forma
$\pi_j^{-1}(U_j)\cap X$, por lo que basta verificar que $\phi_i^{-1}
(U_i)\cap \phi_j^{-1}(U_j)=\phi_k^{-1}(U_k)$ para alg\'un $k$.

Sea $k\geq i, j$ y sea $U_k:= \phi^{-1}_{kj}(U_j)\cap \phi^{-1}_{ki}(U_i)
$. Entonces
\[
\phi_k^{-1}(U_k)=\phi_k^{-1}(\phi_{kj}^{-1}(U_j))\cap
\phi_k^{-1}(\phi_{ki}^{-1}(U_i))=\phi_j^{-1}(U_j)\cap\phi_i^{-1}(U_i).
\]
\end{observacion}

\begin{observacion}\label{O5.5}
Se puede probar que si $\{A_i,\phi_{ji},I\}_{\substack{i,j\in I\\ i\leq j}}$
es un sistema inverso de espacios compactos Hausdorff no vac{\'\i}os,
entonces $X\neq \emptyset$.
\end{observacion}

Tambi\'en en general, si los $A_i$ son espacios topol\'ogicos
(y posiblemente algo m\'as) consideramos a $B=\prod_{
i\in I}A_i$ con la topolog{\'\i}a producto. Entonces $X$ es un subespacio cerrado de $B$ (Proposici\'on \ref{P7.1.6}).
En particular, si cada, $A_i$ es compacto, entonces
$B$ es compacto y por lo tanto $X$ es compacto.

\begin{proposicion}\label{P7.1.6} $X$ es cerrado en $B$.
\end{proposicion}

\begin{proof} Sea $\big(a_i\big)_{i\in I}\in B\setminus X$. Entonces
existen $i\leq j$ tales que $\phi_{ji}(a_j)\neq a_i$. Por ser $A_i$
Hausdorff existen vecindades abiertas $U$ de $\phi_{ji}(a_j)$ y $V$
de $a_i$ tales que $U\cap V=\emptyset$. Sea $W:=\phi^{-1}_{ji}(U)$,
el cual es un abierto de $A_j$. El conjunto $\tilde{U}=V\times W\times
\prod\limits_{k\neq i, j}A_k\subseteq B$ es un abierto y $(a_i)_{i\in I}
\in \tilde{U}$. Puesto que $\phi_{ji}(W)\subseteq U$ y $U\cap V=
\emptyset$ se tiene que $\tilde{U}\cap X=\emptyset$ lo cual 
prueba que $B\setminus X$ es abierto y que $X$ es cerrado. $\fin$
\end{proof}

Como mencionamos anteriormente, estamos interesados en grupos
de Galois, por ello recordamos la siguiente definici\'on.

\begin{definicion}\label{D5.6} Un grupo de $G$ se llama 
{\em grupo topol\'ogico\index{grupo topol\'ogico}} si $G$ es un
espacio topol\'ogico tal que las operaciones de grupo
\begin{alignat*}{3}
\circ\colon G\times G&\to G\qquad {\text{y}}\qquad & i\colon G&\to G\\
(x,y)&\mapsto xy &x&\mapsto x^{-1}
\end{alignat*}
son funciones continuas.
\end{definicion}

Cuando $G$ sea un grupo finito, siempre le daremos a $G$ la
topolog{\'\i}a discreta. En general tenemos:

\begin{proposicion}\label{P5.7}
Un grupo topol\'ogico $G$ es Hausdorff si y solamente si la
identidad de $G$, $\{e\}$ es cerrada en $G$.
\end{proposicion}

\begin{proof}
Si $G$ es Hausdorff, los puntos son cerrados.

Rec{\'\i}procamente, si $\{e\}$ es un conjunto cerrado en $G$, entonces
$\mu^{-1}(\{e\})\subseteq G\times G$ es un conjunto cerrado
donde $\mu\colon G\times G\to G$ est\'a dada por
$\mu(x,y)=xy^{-1}$. Ahora bien, $\mu$ es una funci\'on continua por
ser la composici\'on de las funciones continuas $(\Id, i)\colon
G\times G\to G\times G$, $(\Id,i)(x,y)=(x,y^{-1})$ y la multiplicaci\'on.
Ahora, $\mu^{-1}(\{e\})=\{(x,x)\mid x\in G\}=\Delta$. Sabemos en
general que un espacio topol\'ogico $X$ es Hausdorff si y solamente
si $\Delta=\{(x,x)\mid x\in X\}$ es cerrado en $X\times X$. Por lo
tanto $G$ es Hausdorff. $\fin$.
\end{proof}

Otra observaci\'on es que no solamente $\{e\}$ caracteriza si
$G$ es Hausdorff o no, sino que la topolog{\'\i}a misma de $G$
est\'a determinada por las vecindades de $\{e\}$. Esto se sigue
de que si $g\in G$ est\'a fijo, entonces $\xi_g\colon G\to G$ dada
por $\xi_g(h)=gh$ es un homeomorfismo de espacios topol\'ogicos
pues $\xi_g$ es continua y $\xi_g^{-1}=\xi_{g^{-1}}$. Adem\'as
$\xi_g(e)=g$. Por lo tanto $W$ es un vecindad de $g$
si y solamente si $g^{-1}W=\xi_{g^{-1}}(W)$ es una vecindad
de $\{e\}$.

Una pregunta que contestaremos a continuaci\'on es: ?`Qu\'e
grupos $G$ pueden ser grupos de Galois de alguna extensi\'on
de campos? Sabemos que si $G$ es finito, entonces $G$
es el grupo de Galois de una extensi\'on de campos $L/K$.
Recordemos r\'apidamente su demostraci\'on para ver m\'as
adelante su contraparte infinita.

Sea $k$ cualquier campo y sea $n\in{\ma N}$ tal que $G$ es
subgrupo del grupo sim\'etrico $S_n$. Sean $x_1,\ldots, x_n$
variables independientes y $L:=k(x_1,x_2,\ldots,x_n)=\coc
k[x_1,x_2,\ldots,x_n]$ el campo de las funciones racionales en
$n$ variables, esto es, $k[x_1,\ldots,x_n]$ es el anillo de 
polinomios en $n$ variables y $k(x_1,\ldots, x_n)$ el campo
de cocientes.

Hacemos actuar $S_n$ sobre $L$ de la siguiente forma. Si
$\sigma\in S_n$ y  si $f(x_1,\ldots,x_n)\in L$, entonces
\[
(\sigma\circ f)(x_1,\ldots,x_n):=f(x_{\sigma(1)},\ldots, x_{\sigma(n)}).
\]
Al considerar $G$ como subgrupo de $S_n$, se sigue del
Teorema de Artin\index{teorema!Artin}\index{Artin!teorema de $\sim$},
que $L/L^G:=K$ es una extensi\'on de Galois con grupo de Galois
$G$.

Resulta ser que no cualquier grupo infinito puede ser el grupo de
Galois de alguna extensi\'on. Por ejemplo, veremos que ${\ma Z}$
no puede ser grupo de Galois de ninguna extensi\'on de campos.

\begin{definicion}\label{D5.8} Un {\em grupo
profinito\index{grupo!profinito}} $G$ es un grupo topol\'ogico
Hausdorff compacto que contiene una base de vecindades abiertas
de $\{e\}$ que consiste de subgrupos normales de $G$.
\end{definicion}

\begin{observacion}\label{O7.1.12}
Si $G$ es un grupo finito, a $G$ se le da la topolog{\'\i}a discreta
y con esta topolog{\'\i}a $G$ es un grupo profinito.
\end{observacion}

La raz\'on por la cual los grupos de la Definici\'on \ref{D5.8}
se llaman profinitos es debido a que son l{\'\i}mites inversos
de grupos finitos. De hecho tenemos los siguientes resultados.

\begin{teorema}\label{CClaseT1.4.4} Si $G$ es profinito y si $N$ var{\'\i}a
a trav\'es de los subgrupos abiertos normales de $G$, entonces
tanto algebraica como topol\'ogicamente, se tiene
$G\cong \lim\limits_{\substack{\leftarrow\\ N}}G/N$.

Rec{\'\i}procamente, si $\{G_i,f_i\}_{i\in I}$ es un sistema
proyectivo consistente de grupos finitos $G_i$, se tiene que
$G:=\lim\limits_{\substack{\leftarrow\\ i}}G_i$ es un 
grupo profinito.
\end{teorema}
\begin{proof} \cite[Theorem 11.3.15]{Vil2006}. $\fin$
\end{proof}

\begin{observacion}[Teorema \ref{T6.2}]\label{CClaseO1.4.5}
Los grupos de Galois son grupos profinitos pues si $\{K\mid
K/k\text{\ es Galois y\ } [K:k]<\infty\}$, entonces
\[
\lim\limits_{\substack{\leftarrow\\ K}}\Gal(K/k)\cong\Gal(\Omega/k).
\]

Se puede probar que todo grupo profinito es grupo de Galois
de alguna extensi\'on de campos $\Omega/k$ (ver 
Teorema \ref{T6.3'}).
\end{observacion}

\begin{teorema}\label{T5.9} Sea $G$ un grupo topol\'ogico.
Las siguientes condiciones son equivalentes.
\l
\item $G$ es un grupo profinito.
\item $G$ es el l{\'\i}mite inverso de grupos finitos.
\item $G$ es un grupo topol\'ogico Hausdorff compacto totalmente
disconexo, es decir, las componentes conexas de $G$ son los
puntos.
\item $G$ es un grupo topol\'ogico Hausdorff compacto que
tiene una base de vecindades de $\{e\}$ que consiste de 
subgrupos normales de $G$.
\end{list}
\end{teorema}

\begin{proof} Ver \cite[Theorem 11.3.16, p\'agina 398]{Vil2006}.
$\fin$.
\end{proof}

\begin{observacion}\label{O5.10}
Notemos que un grupo profinito $G$ es {\em completo}, es
decir, toda sucesi\'on de Cauchy en $G$ converge en $G$.
\end{observacion}

\begin{ejemplos}\label{Ej5.10'}
\las
\item Si $G$ es finito, entonces $G$ es profinito.
\item Sea $I:={\ma N}$ con orden definido por $n\leq m \iff
n|m$. Sea
\begin{align*}
f_{m.n}\colon {\ma Z}/m{\ma Z}&\to {\ma Z}/n{\ma Z}\\
a\bmod m&\mapsto a\bmod n.
\end{align*}
Entonces el {\em anillo de Pr\"ufer\index{Pr\"ufer!anillo de
$\sim$}\index{anillo de Pr\"ufer}} $\hat{{\ma Z}}$ se define
por
\[
\hat{{\ma Z}}:=\lim_{\substack{\longleftarrow\\ n\in{\ma N}}}
{\ma Z}/n{\ma Z}\subseteq \prod_{n\in {\ma N}} {\ma Z}/
n{\ma Z}.
\]
$\hat{{\ma Z}}$ se llama {\em proc{\'\i}clico\index{grupo
proc{\'\i}clico}} por ser l{\'\i}mite directo de grupos c{\'\i}clicos finitos.

Sea $\varphi\colon {\ma Z}\mapsto \hat{{\ma Z}}$, $\varphi(x):=
(x\bmod n)_{n\in{\ma N}}\in\hat{{\ma Z}}$. Entonces $\varphi$ es
inyectivo y adem\'as $\varphi({\ma Z})$ es denso en $\hat{{\ma Z}}$.

Por otro lado $\begin{array}{ccc}
\hat{{\ma Z}}&\to&n\hat{{\ma Z}}\\ x&\mapsto&nx\end{array}$
es un isomorfismo de grupos y un homeomorfismo de espacios
topol\'ogicos.

\item Sea $p$ un n\'umero primo. Para $n\in {\ma N}\cup\{0\}$ y
$m\leq n$, la proyecci\'on natural $\varphi_{n,m}\colon
{\ma Z}/p^n{\ma Z}\to {\ma Z}/p^m{\ma Z}$ define un sistema 
inverso.

Sea $X:=\lim\limits_{\substack{\longleftarrow\\ n}}{\ma Z}/
p^n{\ma Z}\subseteq \prod\limits_{n=0}^{\infty}{\ma Z}/p^n{\ma Z}$.
Definimos ${\ma Z}_p:=\{\sum_{n=0}^{\infty} a_np^n\mid
a_n\in\{0,1,\ldots,p-1\}\}$. ${\ma Z}_p$ es la completaci\'on
de ${\ma Z}$ con la topolog{\'\i}a $p$--\'adica, es decir, se define
$|x|_p:=p^{-v_p(n)}$ para $x\in{\ma Z}$, donde $v_p$ es la 
valuaci\'on $p$--\'adica la cual est\'a definida como sigue: si
$x\in{\ma Z}$, $x\neq 0$, digamos $x=p^mb$ con $b$ primo
relativo a $p$ y entonces $v_p(x):=m$. Tambi\'en se define
$v_p(0):=\infty$ y $|0|_p:= 0$.

Definimos 
\[
\begin{array}{rcl}\mu\colon {\ma Z}_p&\to& X\\
\sum\limits_{n=0}^{\infty}a_np^n&\mapsto& \big(\sum\limits_{n=0}^i a_n p^n
\big)_{i\in {\ma N}\cup\{0\}}\end{array}.
\]
Entonces $\mu$ es un
isomorfismo de anillos y por tanto ${\ma Z}_p\cong \lim\limits_{
\substack{\longleftarrow\\ n\in{\ma N}}}{\ma Z}/p^n{\ma Z}$ que
tambi\'en es un homeomorfismo de espacios topol\'ogicos.

\end{list}
\end{ejemplos}

\begin{observacion}\label{O5.14} Sea $G$ un grupo profinito
y sea $N$ un subgrupo abierto y normal. Entonces $G:=\cup_{
x\in G}xN$. Puesto que $N$ es abierto, $xN$ es abierto y puesto
que $G$ es compacto, entonces la cubierta abierta $\{xN\}_{x\in G}$
tiene una subcubierta finita, es decir, existen $x_1,\ldots, x_r\in G$
tales que $G=\cup_{i=1}^r x_iN$. En particular $[G:N]\leq r<
\infty$ y por lo tanto $N$ es de {\'\i}ndice finito. Digamos que
$[G:N]=t$ y $y_1=e,\ldots, y_t$ es un conjunto completo de
representantes de las clases m\'odulo $N$: $G=\biguplus_{i=1}^ty_i N$.
En particular $N=G\setminus\big(\biguplus_{i=2}^t y_iN\big)$ es
abierto, por lo que $N$ es cerrado.

M\'as generalmente, si $H$ es un subgrupo abierto de $G$, 
$\cup_{x\not\in H}xH$ es abierto y por tanto $H=G\setminus
\cup_{x\not\in H} xH$ es cerrado. El rec{\'\i}proco no se
cumple: $H=\{e\}$ es cerrado pero si $G$ no es finito, $H$
no es abierto pues $G$ es compacto.

En resumen tenemos que si $N$ es un subgrupo abierto,
entonces $N$ es cerrado y de {\'\i}ndice finito.
\end{observacion}

\section{Teor{\'\i}a de Galois infinita}\label{S2.2}

Antes de presentar un resumen de la teor\'ia infinita de Galois,
recordamos un resultado b\'asico de teor\'ia de Galois finita.
Este resultado se presenta pues ser\'a de uso intensivo en el
Cap\'itulo \ref{Ch12*} cuando estudiemos los campos de g\'eneros
de campos de funciones.

\begin{teorema}\label{T2.2.Gal1}
Sea $F$ cualquier campo y sea $E/F$ una extensi\'on de Galois
finita. Sea $K$ cualquier extensi\'on de $F$. Entonces $KE/K$
es una extensi\'on de Galois y el mapeo de restricci\'on
\[
\rest\colon\Gal(KE/K)\lra\Gal(E/F),\quad \sigma\longmapsto
\sigma|_E
\]
es un monomorfismo que define un isomorfismo $\Gal(KE/K)\cong
\Gal(E/E\cap K)\subseteq \Gal(E/F)$. $\fin$
\end{teorema}

La correspondencia que usaremos en el Cap\'itulo \ref{Ch12*}
es la siguiente consecuencia del Teorema \ref{T2.2.Gal1}.

\begin{corolario}\label{C2.2.Gal2}
Sean $E/F$ una extensi\'on de Galois 
y $K/F$ una extensi\'on arbitraria tal que $E\cap K=F$. Entonces
existe una correspondencia biyectiva entre las redes de
subcampos ${\mc A}=\{M\mid F\subseteq M\subseteq E\}$
y los subcampos ${\mc B}=\{R\mid K\subseteq  R\subseteq
KE\}$ dada por
\[
\Psi\colon{\mc A}\lra{\mc B}, \quad M\xmapsto{\ \Psi\ } MK.
\]

La biyecci\'on inversa est\'a dada por
\[
\Phi\colon{\mc B}\lra {\mc A}, \quad R\xmapsto{\ \Phi\ } R\cap E.
\]

En particular tenemos para cualquier subcampo $F\subseteq 
M\subseteq E$ y para cualquier subcampo $K\subseteq R
\subseteq KE$
\begin{gather*}
M=MK\cap E\quad\text{y}\quad R=(R\cap E)K. \tag*{\fin}
\end{gather*}
\end{corolario}
\[
\xymatrix{
E\ar@{-}[rr]\ar@{-}[d]&& EK\ar@{-}[d]\\
R\cap E\ar@{<->}[rr]\ar@{-}[d]&&R=(R\cap E)K\ar@{-}[d]\\
M=MK\cap E\ar@{<->}[rr]\ar@{-}[d]&&MK\ar@{-}[d]\\
F=E\cap K\ar@{-}[rr]&&K}
\]

\begin{observacion}\label{O6.0}
Si $A$ y $B$ son subgrupos de $\Gal(E/F)$ y si $E^C$
denota al campo fijo de $E$ bajo el subgrupo $C$, entonces
se tiene
\begin{gather*}
E^A E^B=E^{A\cap B}\quad\text{y}\quad E^A\cap E^B=
E^{\langle A,B\rangle}.
\end{gather*}
\end{observacion}

\begin{definicion}\label{D6.1} Dado un campo $k$, se
denota por $G_k:=\Gal(\bar{k}/k)$ al grupo de Galois
de $\bar{k}/k$ donde $\bar{k}$ es una cerradura 
separable de $k$ y $G_k$ es el {\em grupo absoluto de
Galois de $k$\index{grupo absoluto de Galois de
un campo}}.
\end{definicion}

En general $G_k$ es un grupo infinito y el Teorema de
Correspondencia de Galois ya no se cumple en este
caso.

\begin{ejemplo}\label{E6.2} Si ${\ma F}_p$ es campo
finito de $p$ elementos y si $G:=G_{{\ma F}_p}=\Gal(
\bar{{\ma F}}_p/{\ma F}_p)$, entonces si $H=(\varphi)$
es el grupo generado por el automorfismo de Frobenius,
$\varphi\colon \bar{{\ma F}}_p\to\bar{{\ma F}}_p$, $
\varphi(x)=x^p$ satisface que ${\ma F}_p=\bar{{\ma F}}_p^H
=\bar{{\ma F}}_p^G$ pero $G\neq H$.
\end{ejemplo}

Para establecer la correspondencia de Galois debemos
de proveer a $G$ de una topolog{\'\i}a.

En general, sea $L/K$ una extensi\'on algebraica normal
y separable de campos, es decir, una extensi\'on de
Galois. Sea ${\cal K}:=\{K_i\mid i\in I\}$ la colecci\'on
de todos los campos $K_i$ tales que $K\subseteq K_i
\subseteq L$ y $K_i/K$ es una extensi\'on finita de 
Galois. Entonces $L=\cup_{i\in I}K_i$.

Sea $G:=\Gal(L/K)$, $N_i:=\Gal(L/K_i)$, $i\in I$. Entonces
$K_i=L^{N_i}=\{\alpha\in L\mid \sigma(\alpha)=\alpha\ 
\forall\ \sigma\in N_i\}$ y se tiene $N_i\lhd G$, $G/N_i\cong
\Gal(K_i/K)$ es un grupo finito.

Se define en $G$ la {\em topolog{\'\i}a de
Krull\index{Krull!topolog{\'\i}a de $\sim$}\index{topolog{\'\i}a
de Krull}} definiendo $\{\sigma N_i\mid i\in I\}$ como
un sistema de vecindades abiertas de $\sigma\in G$. Se 
tiene que la multiplicaci\'on y la inversi\'on
\begin{alignat*}{8}
\varphi\colon G\times G &\to&&\ 
G\quad &, \quad &i\colon &G
&\to &&\ G\\
(\sigma,\psi)&\mapsto&&\ 
\sigma\psi &&&\sigma&\mapsto &&\ 
\sigma^{-1}
\end{alignat*}
son mapeos continuos pues $\varphi^{-1}(\sigma \psi
N_i)\supseteq \sigma N_i \times \psi N_j$ e $i^{-1}
(\sigma^{-1} N_j)=\sigma N_j$ y por lo tanto $G$
con esta topolog{\'\i}a es un grupo topol\'ogico que
adem\'as es Hausdorff ya que $\cap_{i\in I}N_i=\{e\}$.

Se tiene:

\begin{teorema}\label{T6.2} El grupo $G=\Gal(L/K)$ con la
topolog{\'\i}a de Krull es un grupo profinito y $G\cong
\lim\limits_{\substack{\longleftarrow\\ i\in I}} G/N_i \cong
\lim\limits_{\substack{\longleftarrow\\ i\in I}} \Gal(K_i/K)$
tanto algebraica como topol\'ogicamente. En otras palabras,
$\Gal\Big(\bigcup\limits_i K_i/K\Big)=\lim\limits_{\substack{
\leftarrow\\ i}}\Gal(K_i/K)$. Equivalentemente
\[
\Gal\Big(\lim_{\substack{\rightarrow\\ i}}K_i/K\Big)\cong
\lim_{\substack{\leftarrow\\ i}}\Gal(K_i/K).
\]
\end{teorema}

\begin{proof} \cite[Theorem 11.4.5, p\'agina 402]{Vil2006}.
$\fin$
\end{proof}

Con esta topolog{\'\i}a tenemos:

\begin{teorema}[Teorema Fundamental de la Teor{\'\i}a
de Galois\index{teorema!fundamental de la teor{\'\i}a
de Galois}]\label{T6.3}
Sea $K/F$ una extensi\'on de Galois con grupo $G=\Gal(
K/F)$. Sean 
\begin{gather*}
{\cal F}(K/F)=\{L\mid L {\text{\ es un campo
tal que\ }} F\subseteq L\subseteq K\}\quad
{\text{y}} \\
S(G)=\{H\mid H {\text{\ es un subgrupo cerrado de\ }} G\}.\\
\intertext{Sean}
\Phi\colon {\cal F}(K/F)\longrightarrow S(G)\quad {\text{y}}
\quad \Psi\colon S(G)\longrightarrow {\cal F}(K/F)\\
\intertext{dadas por:}
\Phi(L):=\{\sigma\in G\mid \sigma|_L=\Id_L\}=
\Gal(K/L),\\
\Psi(H):=\{\alpha\in K\mid \sigma \alpha=\alpha
\ \forall\ \sigma \in H\}=K^H.
\end{gather*}

Entonces $\Phi$ y $\Psi$ son biyecciones mutuamente 
inversas. M\'as a\'un,  $L_1\subseteq L_2$ si y solamente
si $\Gal(K/L_1)\supseteq \Gal(K/L_2)$ y $H_1\subseteq
H_2$ si y solamente si $K^{H_1}\supseteq K^{H_2}$.

Si $\sigma\in G$ y $L\in {\cal F}(K/F)$, entonces
\[
\Gal(K/\sigma L)=\sigma \Gal(K/L)\sigma^{-1}
\]
y en particular $L\in {\cal F}(K/F)$ es una extensi\'on normal
si y solamente si $\Gal(K/L)$ es normal en $G$ y en este
caso
\[
\Gal(L/F)\cong \frac{\Gal(K/F)}{\Gal(K/L)}.
\]
Finalmente, los subgrupos abiertos de $G$ corresponden
a subextensiones finitas de $K/F$.
\end{teorema}

\begin{proof}
\cite[Theorem 11.4.9, p\'agina 405]{Vil2006}. $\fin$
\end{proof}

Se tiene que si $G$ es un grupo de Galois, entonces $G$
necesariamente $G$ es un grupo profinito. El rec{\'\i}proco
tambi\'en se cumple.

\begin{teorema}[Leptin\index{teorema!Leptin}]\label{T6.3'}
Sea $G$ un grupo profinito cualquiera. Entonces existe
una extensi\'on de campos $K/F$ tal que $G\cong
\Gal(K/F)$ tanto algebraica como topol\'ogicamente donde
$\Gal(K/F)$ tiene la topolog{\'\i}a de Krull.
\end{teorema}

\begin{proof}
\cite[Theorem 11.4.10, p\'agina 407]{Vil2006}. $\fin$
\end{proof}

%% file: Capitulo3.tex
\chapter{Campos ciclot\'omicos}\label{Ch3}

\section{La funci\'on exponencial
y el n\'umero $\pi$}\label{S3.1}

La funci\'on exponencial $\exp(z)$ ha jugado un
papel central en diversas \'areas de las
Matem\'aticas y de otras disciplinas: Ingenier{\'\i}a,
F{\'\i}sica, etc. La Teor{\'\i}a de N\'umeros y en
particular, la Teor{\'\i}a de Campos de Clase no es
la excepci\'on. Por esto damos un muy breve 
repaso a las propiedades b\'asicas de esta
funci\'on.

\begin{definicion}\label{D1.1.1} La {\em funci\'on
exponencial\index{funci\'on exponencial}} $\exp
\colon {\ma C}\to {\ma C}$ se define por
$\exp(z)=e^z:=\sum_{n=0}^{\infty}\frac{z^n}{n!}$.
\end{definicion}

Es un ejercicio elemental probar que la serie 
converge absoluta y uniformemente por compactos
en ${\ma C}$. En particular $f(z):=e^z$ es una
funci\'on holomorfa en ${\ma C}$ y se tiene
\[
f'(z)=\sum_{n=0}^{\infty}\frac{nz^{n-1}}{n!}=
\sum_{n=1}^{\infty}\frac{z^{n-1}}{(n-1)!}=
\sum_{n=0}^{\infty}\frac{z^n}{n!}=e^z=f(z).
\]

Usando el producto de Cauchy para series y el
Binomio de Newton, se tiene que para todo $z,
w\in{\ma C}$,
\begin{align}\label{Eq1.1}
e^z \cdot e^w&=\sum_{n=0}^{\infty}\frac{z^n}{n!}\cdot
\sum_{n=0}^{\infty}\frac{w^n}{n!} =
\sum_{n=0}^{\infty}\Big(\sum_{k=0}^n \frac{z^k w^{n
-k}}{k! (n-k)!}\Big)=\nonumber \\
&=\sum_{n=0}^{\infty}\frac{1}{n!}\Big(
\sum_{k=0}^n \binom{n}{k} z^k w^{n-k}\Big) =
\sum_{n=0}^{\infty}\frac{(z+w)^n}{n!}=e^{z+w}.
\end{align}

\begin{definicion}\label{D1.1.2} Se definen las
funciones {\em seno} y {\em coseno}: $\sen
\colon {\ma C}\to {\ma C}$, $\cos\colon{\ma C}\to
{\ma C}$ como
\begin{equation}\label{Eq1.2}
\sen(z)=\sum_{n=0}^{\infty}\frac{(-1)^n z^{2n+1}}
{(2n+1)!},\qquad \cos(z)=\sum_{n=0}^{\infty}
\frac{(-1)^n z^{2n}}{(2n)!}.
\end{equation}
\end{definicion}

Es f\'acil verificar las siguientes propiedades para
cualesquiera $z,w\in{\ma C}$:

\begin{lema}\label{L1.1.2'}{\ }
\lasa
\item $e^{iz}=\cos z+ i\sen z$,
\item $\cos^2 z+\sen^2 z= 1$,
\item $\sen (z+w)=\sen z\cos w+\cos z\sen w$,
\item $\cos(z+w)=\cos z\cos w-\sen z \sen w$,
\item $\cos z=\frac{e^{iz}+e^{-iz}}{2}, \quad
\sen z=\frac{e^{iz}-e^{-iz}}{2i}$. $\fin$
\end{list}
\end{lema}

En particular, si $y\in {\ma R}$, entonces
$\cos^2 y+\sen^2 y = 1$ y por tanto $|\cos y|\leq 1$,
$|\sen y|\leq 1$ para $y\in{\ma R}$.

De esta forma tenemos la expresi\'on debida a
Euler: Para $z=x+iy\in {\ma C}$:
\[
e^z=e^{x+iy}=e^xe^{iy}=e^x(\cos y+i\sen y).
\]
Adem\'as, $e^ze^{-z}=e^{z-z}=e^0=1$ lo cual en particular
implica que $e^z\neq 0$ para todo $z\in{\ma C}$. Por otro
lado, para $x\in{\ma R}$, $e^x=e^{x/2+x/2}=(e^{x/2})^2>0$.

Tambi\'en se tiene que $|e^z|=|e^x||e^{iy}|=e^x\sqrt{\cos^2 y
+\sen^2 y}=e^x$ y $(\sen z)'=\cos z$, $(\cos z)'=-\sen z$.
Para $y\in{\ma R}$, 
\begin{align*}
\sen y &=\sum_{n=0}^{\infty}\frac{(-1)^n y^{2n+1}}{(2n+1)!}=
\sum_{m=0}^{\infty}\Big(\frac{(-1)^{2m}y^{4m+1}}{(4m+1)!}+
\frac{(-1)^{2m+1}y^{4m+2}}{(4m+3)!}\Big)\\
&=\sum_{m=0}^{\infty}\frac{y^{4m+1}}{(4m+3)!}((4m+3)(4m+2)
-y^2).
\end{align*}
En particular, para $0<\sqrt{y}<\sqrt{6}$ se tiene que 
$\sen y>0$ y $\cos y$ es una funci\'on decreciente en el intervalo
$(0,\sqrt{6})$. Puesto que $\cos 0=1>0$ y $\cos 2<0$, existe
un \'unico $\xi_0\in (0,2)$ tal que $\cos \xi_0=0$.

\begin{definicion}\label{D1.1.3}
Se define el n\'umero $\pi$ como el \'unico elemento $\pi\in(0,4)$
tal que $\cos \pi/2=0$.
\end{definicion}

Notemos que $\sen^2 \pi/2=1$ y $0<\pi/2<\sqrt{6}$, por lo que
$\sen \pi/2 >0$ de donde se sigue que $\sen \pi/2=1$.
Del Lema \ref{L1.1.2'} es f\'acil verificar que $\sen \pi =0$,
$\cos 3\pi/2=0$ y que $\sen (z+2\pi)=\sen z$, $\cos(z+2\pi)=
\cos z$ para toda $z\in {\ma C}$. En particular $\sen$ y $\cos$
son funciones peri\'odicas. Sea $t_0\in{\ma R}$, $t_0>0$ m{\'\i}nimo
tal que $\sen(y+t_0)=\sen y$ para toda $y\in{\ma R}$. Se tiene
que $0<t_0\leq 2\pi$.

Sea $n\in{\ma N}$ tal que $nt_0\leq 2\pi<(n+1)t_0$. Entonces
se tiene $t_0=(n+1)t_0-nt_0>2\pi-nt_0$ y puesto que
$\sen(z+(2\pi-nt_0))=\sen z$, y $2\pi-nt_0<t_0$, se sigue
que $2\pi=nt_0$. Ahora bien, $\sen y>0=\sen 0$ para $y\in(0,
\pi/2]$ lo cual implica que $t_0>\pi/2$ y en particular
$2\pi=nt_0>n\pi/2$, esto es, $n<4$. Se tiene que $n\neq 2$ pues
si $2\pi/2=\pi$ satisface $\cos \pi =-1\neq 1=\cos 0$. Similarmente
$n\neq 3$ pues $2\pi/3$ satisface $\sen(2\pi/3)\neq 0=\sen 0$. 
Se sigue que $n=1$, esto es, $t_0=2\pi$.

Similarmente para la funci\'on coseno.

\begin{teorema}\label{T1.1.4} El m{\'\i}nimo per{\'\i}odo para
las funciones seno y coseno es $t_0=2\pi$. $\fin$
\end{teorema}

\begin{corolario}\label{C1.1.5} Se tiene $e^z=1\iff z=2n\pi i$, $n
\in{\ma Z}$.
\end{corolario}

\begin{proof} {\ }

\noindent
$\Leftarrow)$ $e^{2n\pi i}=\cos 2n\pi+i \sen 2n\pi=\cos 0+i
\sen 0=1$.

\noindent
$\Rightarrow)$ Si $e^z=\cos z+i\sen z=1$, entonces $|e^z|=e^x=1$
esto es $x=0$ pues $e^x$ es una funci\'on creciente y $e^0=1$. 
Por lo tanto $z=iy$, $y\in{\ma R}$. Tenemos entonces que
$\cos y+i \sen y=1$ de donde $\cos y=1$ y $\sen y=0$. Por la
discusi\'on anterior, $\cos 0=\cos y =1 $ y $\sen 0=\sen y =0$ de
donde se sigue el resultado. $\fin$
\end{proof}

El desarrollo anterior nos conduce a nuestra \'area de inter\'es,
esto es, el c\'alculo de las ra{\'\i}ces $n$--\'esimas de la unidad,
$n\in{\ma N}$. M\'as generalmente, tenemos:

\begin{proposicion}[F\'ormula de De Moivre\index{Moivre!f\'ormula
de $\sim$}]\label{P1.1.6}
Sea $z_0\in{\ma C}$, $z_0\neq 0$. Entonces $z_0$ se puede
escribir como $z_0=\rho e^{i\alpha}$, $\alpha\in{\ma R}$,
$\rho\in{\ma R}$, $\rho=|z_0|>0$. Adem\'as, para $n\in{\ma N}$, existen
exactamente $n$ n\'umeros complejos $\omega_k$, $k=0,1,\ldots,
n-1$ tales que $\omega_k^n=z_0$. Los elementos $\omega_k$
est\'an dados por
\[
\omega_k=\rho^{1/n}e^{((\alpha+2n\pi)/k)i}, k=0,1,\ldots, n-1.
\]
\end{proposicion}

\begin{proof}
Notemos que la funci\'on $g\colon{\ma R}\to S^1$,  $g(y):=e^{iy}$ donde
$S^1=\{\xi\in{\ma C}\mid |\xi|=1\}$, es
suprayectiva. Adem\'as $g\colon [0,2\pi)\to S^1$ es una
funci\'on biyectiva. Todo lo anterior es consecuencia de la
discusi\'on anterior sobre las funciones seno y coseno y no
presentamos los detalles. Por tanto, dado $z_0\in{\ma C}$,
$z_0\neq 0$, entonces $z_1=\frac{z_0}{|z_0|}$ satisface que
$|z_1|=1$ y por tanto existe un \'unico $\alpha\in[0,2\pi)$ tal que
$e^{i\alpha}=z_1$. Por tanto $z_0=\rho e^{i\alpha}$, $\rho=|z_0|$.

Sea $\omega\in{\ma C}$ tal que $\omega^n=z_0$. Escribamos
$\omega=\mu e^{i\beta}$, $\mu=|\omega|>0$, $\beta
\in{\ma R}$. Entonces $\omega^n=\mu^n e^{in\beta}=
\rho e^{i\alpha}=z_0$. Por lo tanto $\mu^n=|\omega^n|=
|z_0|=\rho$, esto es, $\mu=\rho^{1/n}$. Adem\'as $e^{in\beta}=
e{i\alpha}$ lo cual equivale a $e^{i(n\beta-\alpha)}=1$.

Por el Corolario \ref{C1.1.5} se tiene que $n\beta-\alpha=2m\pi$
para alg\'un $m\in{\ma Z}$. Se sigue que $\beta=\frac{\alpha+
2m\pi}{n}$. Sea $\omega_m:=\rho^{1/n}e^{((\alpha+2m\pi)/n)i}$,
$m\in{\ma Z}$. Es inmediato que $\omega_m^n=z_0$ y que
$\omega_m=\omega_{m'}\iff m\equiv m\bmod n$. Por lo tanto
hay exactamente $n$ ra{\'\i}ces: $\omega_0,\ldots, \omega_{n-1}$.
$\fin$

\end{proof}

\begin{definicion}\label{D1.1.7} Se define $\zeta_n$ por
$\zeta_n=\exp \big(\frac{2\pi i}{n}\big)$.
\end{definicion}

Notemos que $\zeta_n^n=1$ y que $\zeta_n^m\neq 1$ para
$1\leq m\leq n-1$. Adem\'as $\{\zeta_0^n=1, \zeta_n, \zeta_n^2,
\ldots, \zeta_n^{n-1}=\zeta_n^{-1}\}=W_n$ son las ra{\'\i}ces
del polinomio $p(z)=z^n-1\in {\ma C}[z]$. Notemos que $W_n$
es un grupo c{\'\i}clico de orden $n$. Los generadores de $W_n$
son los elementos $\zeta_n^a$ con $\mcd(a,n)=1$.

\begin{observacion}\label{O1.1.8} Si $n|m$ entonces
$\zeta_m^n=\exp\big(\frac{2\pi i}{m}\cdot n\big)=\exp\big(\frac{
2\pi i}{m/n}\big)=\zeta_{m/n}$ y m\'as generalmente, si
$n=xt$ con $t|m$, $\zeta_m^n=\zeta_m^{xt}=\zeta_{m/t}^x$.
\end{observacion}

\section{Campos ciclot\'omicos}\label{S3.2}

\begin{definicion}\label{D1.2.1} Para $n\in{\ma N}$ se define
el $n$--\'esimo campo ciclot\'omico\index{campo ciclot\'omico}
 por $\cic n{}$.
\end{definicion}

Notemos que $\cic n{}/{\ma Q}$ es una extensi\'on de Galois
pues al ser ${\ma Q}$ de caracter{\'\i}stica $0$, la extensi\'on
es separable y $\cic n{}$ es el campo de descomposici\'on del
polinomio $x^n-1$ sobre ${\ma Q}$. Sea $G_n:=\Gal(\cic n{}/
{\ma Q})$. Entonces $\sigma\in G_n$, $\sigma$ est\'a 
determinado por su acci\'on en $\zeta_n$ y $\sigma \zeta_n$
debe ser una ra{\'\i}z de $x^n-1$, por lo tanto $\sigma \zeta_n
=\zeta_n^a$. Denotamos a este elemento por $\sigma=\sigma_a$.
Ahora bien, si $\sigma^{-1}\zeta_n=\zeta_n^b$, se tiene
$\zeta_n=\sigma^{-1}\sigma \zeta_n=\zeta_n^{ab}$, esto es
$ab\equiv 1\bmod n$ y en particular $a\in U_n=\{t\in{\ma Z}/
n{\ma Z}\mid (t,n)=1\}$ donde $t\in{\ma Z}$, $\bar{t}=t\bmod n$.
Es claro que la funci\'on $\varphi\colon G_n\to U_n$, $\varphi(
\sigma_a)=a$ es un monomorfismo de grupos. En particular
$G_n$ es un grupo abeliano.

\begin{definicion}\label{D1.2.2} Para $n\in{\ma N}$ se define
el $n$--\'esimo polinomio ciclot\'omico\index{polinomio
ciclot\'omico}\index{ciclot\'omico!polinomio} por 
\[
\psi_n(x)=\prod_{\substack{(i,n)=1\\ 0\leq i<n}}(x-\zeta_n^i).
\]
\end{definicion}

Se tiene que $\gr \psi_n(x)=|U_n|=|\{a\in{\ma Z}\mid 0\leq a
<n, \mcd(a,n)=1\}|=\varphi(n)$ donde $\varphi$ es la
funci\'on fi de Euler\index{funci\'on fi de Euler}\index{Euler!funci\'on
fi de $\sim$}.

\begin{proposicion}\label{P1.2.3} Para $n\in {\ma N}$ se tiene
\[
x^n-1=\prod_{d|n}\psi_d(x).
\]
\end{proposicion}

\begin{proof}
Se tiene $x^n-1=\prod_{i=0}^{n-1}(x-\zeta_n^i)$, de donde
se sigue que $\psi_d(x)|x^n-1$ para toda $d|n$ pues $
\psi_d(x)=\prod_{(i,d)=1}(x-\zeta_d^i)=\prod_{(i,d)=1}(x-
\zeta_n^{in/d})$, esto es, 
$\zeta_n^{in/d}=\zeta_d^i$.

Ahora bien, veamos que si $d_1|n$, $d_2|n$ y $d_1\neq d_2$
entonces $\mcd (\psi_{d_1},\psi_{d_2})=1$. En efecto, si $\delta$ fuese
una ra{\'\i}z com\'un de $\psi_{d_1}$ y $\psi_{d_2}$, entonces
$\delta=\zeta_{d_1}^{i_1}=\zeta_{d_2}^{i_2}$ para algunos
$i_1, i_2$ tales que $\mcd(i_j,d_j)=1$ para $j=1,2$. Entonces
$\delta=\zeta_{d_1d_2}^{i_1d_2}=\zeta_{d_1d_2}^{i_2d_1}$ de
donde se seguir{\'\i}a que $i_1d_2=i_2d_1$. Puesto que $\mcd(i_1,
d_1)=1$, se tiene que $i_1|i_2$ y viceversa por lo que $i_1=i_2$
y $d_1=d_2$ contrario a lo supuesto. Por lo tanto
\[
\prod_{d|n}\psi_d(x)|x^n-1.
\]
La igualdad se sigue de que ambos polinomios son m\'onicos
y de que 
\begin{gather*}
\gr(\prod_{d|n}\psi_d(x))=\sum_{d|n}\varphi(d)
=n=\gr(x^n-1). \tag*{$\fin$}
\end{gather*}
\end{proof}

La igualdad $\sum_{d|n}\varphi(d)=n$ la probamos a continuaci\'on.

\begin{proposicion}\label{P1.2.4} Sea $n\in{\ma N}$ y sea
$\varphi$ la funci\'on fi de Euler. Entonces $\sum_{d|n}
\varphi(d)=n$.
\end{proposicion}

\begin{proof}

Damos dos demostraciones. Para la primera, consideremos $C_n$
un grupo c{\'\i}clico de $n$ elementos. Sea $A_t:=\{x\in C_n\mid o(x)
=t\}$ donde $o(x)$ denota el orden del elemento $x$. Si $t\nmid n$,
se tiene que $A_t=\emptyset$. Si $t\mid n$, entonces $C_n$ tiene
un \'unico subgrupo $H_t$ de orden $t$ y puesto que $C_n$
es c{\'\i}clico, este subgrupo es a su vez c{\'\i}clico. Los elementos
de orden $t$ de $C_n$ son precisamente los generadores de 
$H_t$ y por tanto $|A_t|=\varphi(t)$.

Se tiene que si $d_1\neq d_2$, $A_{d_1}\cap A_{d_2}=\emptyset$
y cada $x\in C_n$ est\'a en alg\'un $A_t$, de donde:
\[
n=|C_n|=\sum_{t=1}^n |A_t|=\sum_{t|n}A_t =\sum_{t|d} \varphi(t).
\]
Esto termina la primera demostraci\'on.

Presentamos una segunda demostraci\'on m\'as directa. Primero,
si $p$ es un n\'umero primo y $\alpha\in{\ma N}$, entonces
$\varphi(p^{\alpha})=p^{\alpha}-p^{\alpha-1}$. Por tanto
\[
p^{\alpha}=\sum_{t=1}^{\alpha}(p^t-p^{t-1})+1=
\sum_{t=1}^{\alpha} \varphi(p^t) +1=\sum_{t=0}^{\alpha}
\varphi(p^t),
\]
En general, puesto que $\varphi$ es una funci\'on multiplicativa,
es decir, si $\mcd (n,m)=1$, $\varphi(nm)=\varphi(n)\varphi(m)$,
se tiene en general que si $n=p_1^{\alpha_1}\cdots p_r^{\alpha_r}$
con $p_1,\ldots, p_r$ primos distintos y $\alpha_i\geq 1$, $1\leq
i\leq r$, entonces
\begin{align*}
n&=p_1^{\alpha_1}\cdots p_r^{\alpha_r}=\prod_{i=1}^{r} p_i^{\alpha_i}
=\prod_{i=1}^r\big(\sum_{j_i=0}^{\alpha_i}\varphi(p_i^{j_i})\big)=
\sum_{0\leq j_i\leq \alpha_i} \prod_{i=1}^r \varphi(p_i^{j_i})\\
&=\sum_{0\leq j_i\leq \alpha_i}\varphi\big(\prod_{i=1}^r p_i^{j_i}\big)=
\sum_{\substack{0\leq \beta_i\leq \alpha_i \\ 0\leq i\leq r}}
\varphi(p_1^{\beta_1}\cdots p_r^{\beta_r})=\sum_{d|n}\varphi(d).
\end{align*}
Esto termina la segunda demostraci\'on. $\fin$

\end{proof}

Ahora bien como consecuencia de la Proposici\'on \ref{P1.2.3}
tenemos:

\begin{corolario}\label{C1.2.5} $\psi_n(x)\in{\ma Z}[x]$ para toda
$n\in {\ma N}$.
\end{corolario}

\begin{proof}
Lo hacemos por inducci\'on en $n$. Para $n=1$, se tiene
que $\psi_i(x)=x-1\in{\ma Z}[x]$. Sea $n>1$ y suponemos que
$\psi_d(x)\in{\ma Z}[x]$ para toda $d<n$. Entonces
$x^n-1=\prod_{d|n}\psi_d(x)=\psi_n(x)\cdot\prod_{\substack{d|n\\
d<n}}\psi_d(x)$.

Ahora $\prod_{\substack{d|n\\ d<n}} \psi_d(x)=h(x)\in{\ma Z}[x]$.
Por lo tanto $\psi_n(x)=\frac{x^n-1}{h(x)}\in{\ma Q}[x]$ de donde
$x^n-1=h(x)\psi_n(x)$. Ahora bien, usando ya sea el Lema de
Gauss o el algoritmo de la divisi\'on para dominios enteros y que
$h(x)$ es un polinomio m\'onico, se sigue que $\psi_n(x)\in
{\ma Z}[x]$. $\fin$

\end{proof}

\begin{ejemplos}\label{Ej1.2.6}
{\ }

\las
\item $\psi_1(x)=x-1$,

$\psi_2(x)=x+1=\frac{x^2-1}{x-1}$,

$\psi_3(x) =x^2+x+1=\frac{x^3-1}{x-1}$,

$\psi_4(x)=x^2+1=\frac{x^4-1}{\psi_1(x)\psi_2(x)}=\frac{x^4-1}{x^2-1}$,

$\psi_5(x)=x^4+x^3+x^2+x+1=\frac{x^5-1}{x-1}$.

\item Si $p$ es un n\'umero primo, 
\[
\psi_p(x)=\frac{x^p-1}{\psi_1(x)}=
\frac{x^p-1}{x-1}=x^{p-1}+x^{p-2}+\cdots+x+1.
\]

\item Si $p$ es un n\'umero primo entonces
 \begin{align*}
 \psi_{p^n}(x)&=\frac{x^{p^n}-1}{
\prod_{i=0}^{n-1}\psi_{p^i}(x)}=\frac{x^{p^n-1}}{x^{p^{n-1}}-1}\\
&=x^{p^{n-1}(p-1)}+x^{p^{n-1}(p-2)}+\cdots+x^{p^{n-1}}+1=
\psi_p(x^{p^{n-1}}).
\end{align*}

En particular $\psi_{p^n}(1)=\underbrace{1+1+\cdots+1}_{p\text{\ 
veces}}=p$ y 
\[
p=\prod_{\substack{i=0\\ \mcd(i,p)=1}}^{p^n-1}(1-\zeta_{p^n}^i).
\]

Notemos adem\'as que por el criterio de Eiseinstein, $\psi_p(x)\in
{\ma Z}[x]$ es irreducible. Esto no es casualidad como veremos
a continuaci\'on.
\end{list}
\end{ejemplos}

\begin{proposicion}\label{P1.2.6'} Si $n$ y $m$ son primos relativos,
entonces $\cic n{}\cic m{}=\cic {nm}{}$. En consecuencia si la
descomposici\'on en primos de $n$ est\'a dado por $n=p_1^{\alpha_1}
\cdots p_r^{\alpha_r}$ entonces $\cic n{} =\prod_{i=1}^r \cic {p_i}{
\alpha_i}$.
\end{proposicion}

\begin{proof}

Se tiene $\zeta_n=\zeta_{nm}^m$ y $\zeta_m=\zeta_{nm}^n$ por
tanto $\cic n{} \cic m{}\subseteq \cic {nm}{}$ (de hecho esto
se cumple para todas $n,m\in{\ma N}$).

Ahora, sean $\alpha,\beta\in{\ma Z}$ tales que $\alpha n+\beta m=1$.
Por tanto
\[
\zeta_{nm}=\zeta_{nm}^{\alpha n+\beta m}=\zeta_{nm}^{\alpha n}
\zeta_{nm}^{\beta m}=\zeta_m^{\alpha}\zeta_n^{\beta}\in
\cic n{} \cic m{}.
\]
Por tanto $\cic {nm}{}\subseteq \cic n{}\cic m{}$ de donde se sigue
la igualdad. $\fin$
\end{proof}

\begin{teorema}\label{T1.2.7}
Para cualquier $n\in{\ma N}$, $\psi_n(x)\in{\ma Z}[x]$ es irreducible
sobre ${\ma Q}$.
\end{teorema}

\begin{proof}

Sea $f(x):=\Irr(\zeta_n,x,{\ma Q})$ el polinomio irreducible de 
$\zeta_n$ sobre ${\ma Q}$. Puesto que $\psi_n(\zeta_n)=0$, se tiene
que $f(x)|x^n-1$. Ahora bien, sea $x^n-1=f(x)g(x)$ con $f(x)$ y
$g(x)$ con coeficiente l{\'\i}der igual a $1$. Por el Lema de Gauss
se sigue que $f(x),g(x)\in{\ma Z}[x]$. Notemos que las ra{\'\i}ces
de $\psi_n(x)$ son $\{\zeta_n^d\}_{(d,n)=1}$. En particular cualquier
ra{\'\i}z de $\psi_n(x)$ es de la forma $\zeta_n^{p_1 \cdots p_r}$
con $p_1,\ldots,p_r$ n\'umeros primos, no necesariamente
distintos, tales que $p_i\nmid n$. Ahora bien, si probamos que dada
cualquier ra{\'\i}z $\lambda$ de $f(x)$, entonces $\lambda^p$, con
$p$ un n\'umero primo tal que $p\nmid n$, es ra{\'\i}z de $f(x)$,
entonces se tendr\'a que toda ra{\'\i}z de $\psi_n(x)$ ser\'a
tambi\'en ra{\'\i}z de $f(x)$ y en particular se seguir\'a que
$\psi_n(x)|f(x)$ de donde se obtendr\'a la igualdad $\psi_n(x)=
f(x)$ y que $\psi_n(x)$ es irreducible.

En resumen, vamos a probar que si $\lambda$ es cualquier ra{\'\i}z
de $f(x)$, entonces $\lambda^p$ es tambi\'en de $f(x)$ con $p$
es un n\'umero primo tal que $p\nmid n$.

Supongamos que $\lambda$ es ra{\'\i}z de $f(x)$ pero que 
$\lambda^p$ no lo es. Puesto que $\lambda^p$ 
es ra{\'\i}z de $x^n-1$, entonces $g(\lambda^p)=0$. Puesto que 
$f(\lambda)=0$ y $f(x)$ es irreducible y $\lambda$ es ra{\'\i}z de
$g(x^p)$, se sigue que $f(x)|g(x^p)$. Pongamos $g(x^p)=f(x)h(x)$ 
con $h(x)\in{\ma Z}[x]$ por el Lema de Gauss.

Por otro lado, si $g(x)=x^m+b_{m-1}x^{m-1}+\cdots +b_1x +b_0\in
{\ma Z}[x]$, entonces
\begin{align*}
g(x)^p&\equiv (x^m + b_{m-1}x^{m-1}+\cdots +b_1 x+b_0)^p
\bmod p\\ 
&\equiv
x^{pm}+b_{m-1}^p x^{p(m-1)}+\cdots+b_1^px^p+b_0^p\bmod p\\
&\equiv (x^p)^m+b_{m-1}(x^p)^{m-1}+\cdots + b_1(x^p)+b_0
\bmod p \equiv g(x^p)\bmod p.
\end{align*}
Esto es, m\'odulo $p$, $g(x)^p\equiv f(x)h(x)\bmod p$. En particular
$\overline{g(x)}:=g(x)\bmod p\in{\ma F}_p[x]$ y $\overline{f(x)}$
no son primos relativos en ${\ma F}_p[x]$ y puesto que 
$\overline{x^n-1}=\overline{f(x)}\overline{g(x)}$, se tiene que
$\overline{x^n-1}\in{\ma F}_p[x]$ tiene ra{\'\i}ces m\'ultiples. Sin
embargo la derivada de $\overline{x^n-1}$ es $\overline{nx^{n-1}}\not\equiv 0
\bmod p$ pues $p\nmid n$. La \'unica ra{\'\i}z de la derivada
de $\overline{x^n-1}$ es $\overline{0}$ la cual no es ra{\'\i}z de
$\overline{x^n-1}$ de donde se sigue que $\overline{x^n-1}$ no 
tiene ra{\'\i}ces m\'ultiples. Esta contradicci\'on prueba que 
$\lambda^p$ es ra{\'\i}z de $f(x)$ y termina la demostraci\'on del
teorema. $\fin$

\end{proof}

\begin{observacion}\label{O1.2.8} Se podr{\'\i}a dar otra 
demostraci\'on de que $\psi_n(x)$ es irreducible probando que si
$p$ es un n\'umero primo y $m\in{\ma N}$, entonces $\psi_{p^m}
(x)$ es irreducible por medio del c\'alculo del {\'\i}ndice de
ramificaci\'on de $p$ en $\cic pm$ probando que $e\geq \varphi(
p^m)$. De esta forma, y viendo que no hay m\'as ramificaci\'on,
se seguir{\'\i}a que si $\mcd(m,n)=1$, $\cic n{}\cap \cic m{}={\ma Q}$
y en consecuencia, usando que la funci\'on fi de Euler es 
multiplicativa se deducir{\'\i}a que $[\cic n{}:{\ma Q}]=\varphi(n)$
y que $\psi_n(x)$ es irreducible. Claramente esta demostraci\'on
es mucho m\'as complicada que la presentada, sin embargo
basta hallar una demostraci\'on independiente de que $\cic n{}
\cap \cic m{}={\ma Q}$ para $n, m$ primos relativos.
\end{observacion}

\begin{corolario}\label{C1.2.9} Para $n\in{\ma N}$, $[\cic n{}:{\ma Q}]=
\varphi(n)=\gr \psi_n(x)$ y $\cic n{}/{\ma Q}$ es una extensi\'on
de Galois con grupo de Galois isomorfo a $U_n:=({\ma Z}/n{\ma Z}
)^{\ast}$.
\end{corolario}

\begin{proof}

Se tiene $G_n=\Gal(\cic n{}/{\ma Q})\subseteq U_n$
y $|G_n|=[\cic n{}:{\ma Q}]=\varphi(n)=\gr \psi_n=|
U_n|$ de donde se sigue que son iguales. $\fin$.

\end{proof}

\begin{corolario}\label{C1.2.10}
Si $m$ y $n$ son primos relativos, entonces
$\cic n{}\cap \cic m{}={\ma Q}$. En particular
$\Gal(\cic {nm}{}/{\ma Q})\cong \Gal(\cic n{}/{\ma Q})
\times \Gal(\cic m{}/{\ma Q})$.
\end{corolario}.

\begin{proof} Sea $K:=\cic n{}\cap \cic m{}$. Se tiene
el diagrama
\[
\xymatrix{
&\cic n{}\ar@{-}[r]\ar@{-}[d]&\cic n{}\cic m{}=\cic {nm}{}
\ar@{-}[d]\\
&K\ar@{-}[r]\ar@{-}[dl]&\cic m{}\\{\ma Q}
}
\]
Entonces
\begin{align*}
\varphi(nm)=\varphi(n)\varphi(m)=[\cic {nm}{}:
{\ma Q}]&=[\cic {nm}{}\cic n{}][\cic n{}:{\ma Q}]\\
&=[\cic {nm}{}\cic m{}][\cic m{}:{\ma Q}].
\end{align*}
Es decir
\[
\varphi(n)\varphi(m)=[\cic {nm}{}\cic n{}]\varphi(n)
=[\cic {nm}{}\cic m{}]\varphi(m).
\]

Se sigue que $\varphi(m)=[\cic {nm}{}\cic n{}]$ y
$\varphi(n)=[\cic {nm}{}\cic m{}]$.
En particular tenemos
\begin{align*}
\varphi(n)&=[\cic {nm}{}\cic m{}]\leq 
[\cic n{}:K]\\
&\leq [\cic n{}:K][K:{\ma Q}]=
[\cic n{}:{\ma Q}]=\varphi(n)
\end{align*}
lo cual implica que $[\cic n{}:K]=\varphi(n)$ y por lo
tanto $K={\ma Q}$.
La igualdad $\Gal(\cic {nm}{}/{\ma Q})\cong
\Gal(\cic n{}/{\ma Q})\times \Gal(\cic m{})/{\ma Q})$
se sigue de la Teor{\'\i}a de Galois. $\fin$

\end{proof}

El siguiente resultado es una consecuencia inmediata
del Teorema Chino del Residuo. Aqu{\'\i} presentamos
otra demostraci\'on usando los resultados hasta
ahora obtenidos en campos ciclot\'omicos.

\begin{corolario}\label{C1.2.11}
Si $n=p_1^{\alpha_1}\cdots p_r^{\alpha_r}$ con
$p_1,\ldots, p_r$ primos distintos, entonces
$U_n\cong U_{p_1^{\alpha_1}}\times \cdots \times
U_{p_r^{\alpha_r}}$.
\end{corolario}

\begin{proof}
Se tiene el diagrama
\begin{gather*}
\xymatrix{
&\cic n{}\ar@{-}[ld]\ar@{-}[d]\ar@{-}[rd]\\
\cic {p_1}{\alpha_1}\ar@{--}[r]\ar@{-}[dr]_{
U_{p_1^{\alpha_1}}}
&\cic {p_i}{\alpha_i}\ar@{-}[d]^{U_{p_i^{
\alpha_i}}}\ar@{--}[r]&\cic {p_r}{\alpha_r}
\ar@{-}[dl]^{U_{p_r^{\alpha_r}}}\\ &{\ma Q}
}\\
\intertext{Por lo tanto}
U_n\cong \Gal(\cic n{}/{\ma Q})\cong
\prod_{i=1}^r\Gal(\cic {p_i}{\alpha_i}/{\ma Q})
\cong \prod_{i=1}^r U_{p_i^{\alpha_i}}.
\tag*{$\fin$}
\end{gather*}

\end{proof}

Como hicimos notar antes, se tiene que si $p$
es un n\'umero primo y $n\in{\ma N}$, entonces
$\psi_{p^n}(x)=\psi_p(x^{p^{n-1}})$ (Ejemplo
\ref{Ej1.2.6} (3)). M\'as generalmente, tenemos

\begin{proposicion}\label{P1.2.12}
Si $n=p_1^{\alpha_1}\cdots p_r^{\alpha_r}$ es
la descomposici\'on en primos, entonces
$\psi_n(x)=\psi_{p_1\cdots p_r}\big(
x^{p_1^{\alpha_1-1}\cdots p_r^{\alpha_r-1}}\big)$.
\end{proposicion}

\begin{proof} Primero notemos que
\begin{align*}
\gr \big(\psi_{p_1\cdots p_r}\big(
x^{p_1^{\alpha_1-1}\cdots p_r^{
\alpha_r-1}}\big)\big)&=
x^{p_1^{\alpha_1-1}\cdots p_r^{\alpha_r-1}}
\gr \psi_{p_1\cdots p_r}(x)\\
&=p_1^{\alpha_1-1}\cdots p_r^{\alpha_r-1}
\varphi(p_1\cdots p_r)=\varphi(n)=
\gr \psi_n(x).
\end{align*}

Ahora $\zeta_n^{p_1^{\alpha_1-1}\cdots
p_r^{\alpha_r-1}}=\zeta_{p_1\cdots p_r}$ lo
cual implica que $\zeta_n$ es ra{\'\i}z del
polinomio
$\psi_{p_1\cdots p_r}\big(
x^{p_1^{\alpha_1-1}\cdots p_r^{
\alpha_r-1}}\big)$ de donde $\psi_n(x)|
\psi_{p_1\cdots p_r}\big(
x^{p_1^{\alpha_1-1}\cdots p_r^{
\alpha_r-1}}\big)$ lo cual implica que ambos
son iguales. $\fin$.

\end{proof}

Recordemos la f\'ormula de inversi\'on de
M\"obius\index{f\'ormula de inversi\'on de
M\"obius}\index{M\"obius!f\'ormula de 
inversi\'on de $\sim$}. Consideremos las
funciones $\mu\colon{\ma N}\to {\ma Q}$,
$\varepsilon\colon {\ma N}\to {\ma Q}$ dadas
por
\begin{align*}
\mu(n)&=\begin{cases}
1&\text{si $n=1$;}\\
(-1)^r &\text{si $n=p_1\cdots p_r$ con $p_1,\ldots,
p_r$ son primos distintos;}\\
0 &\text{en otro caso, esto es, si existe $d>1$,
$d^2|n$.}
\end{cases}\\
\varepsilon(n)&=\begin{cases}
1&\text{si $n=1$;}\\ 0&\text{si $n>1$.}
\end{cases}
\end{align*}

Entonces se tiene

\begin{lema}\label{L1.2.13}
$\sum_{d|n}\mu(d)=\varepsilon(n)$.
\end{lema}

\begin{proof}
Sea $n=p_1^{\alpha_1}\cdots p_r^{
\alpha_r}$ la descomposici\'on en primos de
$n$. Entonces
\begin{align*}
\sum_{d|n}\mu(d)&= \sum_{i_1<\cdots <i_t}
\mu(p_{i_1}\cdots p_{i_t})=\sum_{t=0}^r
\binom{r}{t}(-1)^t = (1-1)^r=0^r\\
&=\begin{cases}
1 &\text{si $r=0$}\\0&\text{si $r>0$}\end{cases}
=\begin{cases} 1&\text{si $n=1$}\\ 0&\text{si
$n>1$}\end{cases} = \varepsilon(n).
\tag*{$\fin$}
\end{align*}
\end{proof}

\begin{corolario}[F\'ormula de Inversi\'on de
M\"obius\index{f\'ormula de inversi\'on de
M\"obius}\index{M\"obius!f\'ormula de 
inversi\'on de $\sim$}]\label{C1.2.14}
Sea $k$ cualquier campo y sean $f.g\colon
{\ma N}\to k$ dos funciones tales que:
\las
\item $f(n)=\sum_{d|n}g(d)$. Entonces
\[
g(n)=\sum_{d|n}\mu(n/d)f(d)=\sum_{d|n}
\mu(d)f(n/d).
\]

\item Si $f(n), g(n)\neq 0$ para toda $n\in
{\ma N}$ y $f(n)=\prod_{d|n} g(d)$. Entonces
\[
g(n)=\prod_{d|n}f(d)^{\mu(n/d)}=\prod_{d|n}
f(n/d)^{\mu(d)}.
\]
\end{list}
\end{corolario}

\begin{proof}
{\ } 
\las 
\item Se tiene que
$\sum_{d|n}\mu(d)f(n/d)=\sum_{d|n}\mu(d)
\big(\sum_{s|\frac{n}{d}} g(s)\big)$. Ahora bien,
si $s|\frac{n}{d}$, entonces $d|\frac{n}{s}$. Por
tanto, la \'ultima suma es del tipo $\sum_{t|n}
a_t g(t)$ para algunos $a_t\in{\ma Z}$.

Obtenemos,  por 
el Lema \ref{L1.2.13}, que
$a_t=\sum_{d|\frac{n}{t}}\mu(d)
=\varepsilon(n/t)=
\begin{cases} 1&\text{si $n=t$}\\0&\text{si
$n\neq t$}\end{cases}$. Por lo tanto, se tiene
\[
\sum_{d|n}\mu(d) f(n/d)=\sum_{t|n}
a_t g(t)=g(n).
\]

\item Tenemos
\begin{gather*}
\prod_{d|n} f(d)^{\mu(n/d)} =\prod_{d|n}\Big[
\prod_{t|d}g(t)\Big]^{\mu(n/d)} = \prod_{a|t}
g(a)^{s(a)}\\
\intertext{donde}
s(a)=\sum_{a|d}\mu(n/d)=\sum_{t|\frac{n}{a}}
\mu(t)=\varepsilon(n/a)
\end{gather*}
de donde se sigue el resultado. $\fin$

\end{list}

\end{proof}

La f\'ormula de inversi\'on de M\"obius nos da una
expresi\'on para el polinomio ciclot\'omico
en t\'erminos de los polinomios $x^n-1$.

\begin{proposicion}\label{P1.2.15}
Para $n\in{\ma N}$ se tiene 
\[
\psi_n(x)=\prod_{d|n}(x^d-1)^{\mu(n/d)}.
\]
\end{proposicion}

\begin{proof}

Sean $f,g\colon {\ma N}\to {\ma Q}(x)$, donde
${\ma Q}(x)$ es el campo de las funciones 
racionales sobre ${\ma Q}$ dadas por
\[
f(n):=x^n-1,\qquad g(m):=\psi_m(x).
\]

Por la Proposici\'on \ref{P1.2.3} se tiene que
$f(n)=\prod_{d|n}g(d)$. Del Corolario \ref{C1.2.14}
se sigue que $g(n)=\prod_{d|n}f(d)^{\mu(n/d)}$
que es el resultado enunciado. $\fin$

\end{proof}

\begin{ejemplo}\label{Ej1.2.15'}
Se tiene
\begin{align*}
\psi_{12}(x)&= \prod_{d|12} (x^d-1)^{\mu(12/d)}
\igual_{\substack{\uparrow\\ d\in\{1,2,3,4,6,12\}}}\\
&=
(x-1)^{\mu(12)}(x^2-1)^{\mu(6)}(x^3-1)^{\mu(4)}\\
&\hspace{2cm}
(x^4-1)^{\mu(3)} (x^6-1)^{\mu(2)}(x^{12}-1)^{
\mu(1)}\\
&=(x-1)^0(x^2-1)^1(x^3-1)^0(x^4-1)^{-1}
(x^6-1)^{-1}(x^{12}-1)\\
&=\frac{x^{12}-1}{x^6-1}
\cdot \frac{1}{\Big(\frac{x^4-1}{x^2-1}\Big)}
=\frac{x^6+1}{x^2+1}= x^4-x^2+1.
\end{align*}

Por lo tanto $\psi_{12}(x)=x^4-x^2+1=
(x-\zeta_{12})(x-\zeta_{12}^5)(x-\zeta_{12}^7)
(x-\zeta_{12}^{11})$.
\end{ejemplo}

\begin{observacion}\label{O1.2.16}
Si $n$ es impar, entonces $\varphi(2n)=
\varphi(n)$ y $\cic n{}\subseteq \cic {2n}{}$
de donde se sigue que $\cic n{}=\cic {2n}{}$.
De hecho se tiene 
\[
\cic {2n}{}\igual\limits_{
\substack{\uparrow\\ \mcd(2,n)=1}}\cic 2{}\cic n{}=
{\ma Q}(-1)\cic n{}={\ma Q}\cic n{}=\cic n{}.
\]
Por lo tanto, siempre que consideremos un
campo ciclot\'omico $\cic m{}$, supondremos que
$m\not\equiv 2\bmod 4$.
\end{observacion}

\subsection{Estructura de $U_n$}\label{S1.2.1}

Puesto que $\Gal(\cic n{}/{\ma Q})\cong U_n$
es importante determinar la estructura de este
\'ultimo grupo.

\begin{definicion}\label{D1.2.1.1}
Sea $p\in{\ma N}$ un n\'umero primo. Entonces
definimos
la {\em valuaci\'on $p$--\'adica\index{valuaci\'on
$p$--\'adica}} $v_p$ de ${\ma Q}^{\ast}$ por
$v_p\colon {\ma Q}^{\ast}\to {\ma Z}$ dada de la
siguiente forma. Para $a\in {\ma Z}$, podemos
escribir $a=\pm p^n b$ con $n\geq 0$, $b\in
{\ma N}$ y $p\nmid b$. Entonces $v_p(a):=n$.
\end{definicion}

Si $\alpha=\frac{a}{b}\in{\ma Q}^{\ast}$, $a,b\in
{\ma Z}\setminus \{0\}$, se define
\[
v_p(\alpha):=v_p(a)-v_p(b).
\]
Notemos que si $x,y\in{\ma Z}\setminus \{0\}$,
entonces $v_p(xy)=v_p(x)+v_p(y)$ por lo tanto si
$\alpha=\frac{a}{b}=\frac{c}{d}\in{\ma Q}^{\ast}$,
entonces $ad=bc$ y $v_p(ad)=v_p(a)+v_p(d)=
v_p(b)+v_p(c)=v_p(bc)$ de donde se sigue que
$v_p(a)-v_p(b)=v_p(c)-v_p(d)$ y por lo tanto
la definici\'on de $v_p(\alpha)$ no depende de la 
representaci\'on de $\alpha$ como cociente
de dos enteros.

De la misma forma, se sigue que si $\alpha,\beta
\in{\ma Q}^{\ast}$, entonces $v_p(\alpha \beta)=
v_p(\alpha)+v_p(\beta)$; $v(\alpha^{-1})=
-v_p(\alpha)$ y $v_p\big(\frac{\alpha}{\beta}\big)=
v_p(\alpha)-v_p(\beta)$. M\'as a\'un $\alpha=\frac{a}
{b}$ puede escribirse de manera \'unica como
$\alpha=p^m \frac{c}{d}$ con $p\nmid cd$ y $m
\in{\ma Z}$. Entonces se tiene $v_p(\alpha)=m$.

Notemos que $v_p(-\alpha)=v_p(\alpha)$ y que
$v_p(1)=0$. Se define $v_p(0):=\infty$ donde
$\infty$ es cualquier s{\'\i}mbolo al que supondremos
sujeto a las siguientes reglas:
\begin{enumerate}
\item Para toda $a\in{\ma Z}$, se tiene $a<\infty$;
\item $\infty+\infty=\infty\cdot \infty=\infty$;
\item Si $a\in{\ma Z}\setminus \{0\}$, $a\cdot
\infty =\infty$;
\item El s{\'\i}mbolo $0\cdot \infty$ no se define.
\end{enumerate}

Con esta convenci\'on se tiene:

\begin{teorema}\label{T12.1.2}
Para $\alpha,\beta\in{\ma Q}$ se tiene que
\begin{gather*}
v_p(\alpha+\beta)\geq \min \{v_p(\alpha),v_p(\beta)\}\\
\intertext{y si $v_p(\alpha)\neq v_p(\beta)$ entonces}
v_p(\alpha+\beta)= \min \{v_p(\alpha),v_p(\beta)\}.
\end{gather*}
\end{teorema}

\begin{proof}

Si $\alpha$ o $\beta=0$ no hay nada que probar.
Sean $\alpha,\beta\neq 0$. Escribamos
$\alpha=p^n\frac{c}{d}$, $\beta=p^m\frac{e}{f}$
con $p\nmid cdef$, $n,m\in {\ma Z}$. Entonces
\begin{gather*}
\alpha+\beta=p^n\frac{c}{d}+p^m\frac{e}{f}=
\frac{p^ncf+p^med}{df}=\frac{p^rg}{df}\\
\intertext{donde $r:=\min\{n,m\}$ y $g\in{\ma Z}$.
Por lo tanto}
v_p(\alpha+\beta)\geq r =\min \{v_p(\alpha),v_p(\beta)\}
\end{gather*}
lo cual prueba nuestra primera afirmaci\'on

Ahora, si $v_p(\alpha)\neq v_p(\beta)$, es decir
$n\neq m$, se tiene que $p\nmid g$ y por lo tanto
$v_p(\alpha+\beta)= \min \{v_p(\alpha),v_p(\beta)\}$.

Alternativamente, digamos $n<m$. Entonces
$v_p(\alpha+\beta)\geq \min \{v_p(\alpha),v_p(\beta)\}$
y adem\'as 
\begin{align*}
v_p(\alpha)&=v_p(\alpha+\beta-\beta)\geq \min\{
v_p(\alpha+\beta),v_p(-\beta)\}\\
&= \min\{
v_p(\alpha+\beta),v_p(\beta)\} \geq v_p(\alpha).
\end{align*}

Por lo tanto $v_p(\alpha)=\min\{
v_p(\alpha+\beta),v_p(\beta)\}$ y $v_p(\alpha)<
v_p(\beta)$ lo cual implica que $v_p(\alpha)=
v_p(\alpha+\beta)$. $\fin$.
\end{proof}

Consideremos ahora el grupo $U_{p^n}$ con $p$ un
n\'umero primo. Si $n=1$, entonces $U_p=({\ma Z}/p{\ma
Z})^{\ast}= {\ma F}_p^{\ast}$ con ${\ma F}_p$ el
campo finito de $p$ elementos. Se tiene que el grupo
multiplicativo de un campo finito es c{\'\i}clico. Por lo
tanto $U_p\cong {\ma Z}/(p-1){\ma Z}$.

Ahora supongamos que $p>2$ y que $n\geq 1$. Entonces
\begin{gather*}
|U_{p^n}|=\varphi(p^n)=p^{n-1}(p-1).\\
\intertext{Sea $x:=1+p$. Entonces}
x^{p^k}=(1+p)^{p^k}=1+\sum_{i=1}^{p^k}\binom{p^k}{i}
p^i= 1+p^{k+1}+\sum_{i=2}^{p^k}\binom{p^k}{i}p^i.
\end{gather*}

Veamos que $v_p\Big(\binom{p^k}{i}p^i\Big)>v_p\Big(
\binom{p^k}{1} p\Big)= p^{k+1}$ para $2\leq i\leq p^k$.
Sea 
\begin{align*}
A:&=v_p\Big(\binom{p^k}{i}p^i\Big)-
v_p\Big(\binom{p^k}{1}p\Big)=v_p\Big(\frac{1}{p^k}
\binom{p^k}{i}p^{i-1}\Big)\\
&= v_p\Big(\frac{1}{i} \binom{p^k-1}{i-1}\Big)+(i-1)\geq
i-1-v_p(i).
\end{align*}

Se tiene que para $a\geq 1$ y $p$ un n\'umero
primo $a\leq p^a-(p-1)$ y la desigualdad es estricta para
$a\geq 2$. Si $\mcd(i,p)=1$ entonces $v_p(i)=0$ y
por tanto $A\geq i-1\geq 1>0$, $2\leq i\leq p^k$. Si
$i=p^a b$, $\mcd(p,b)=1$,  entonces $v_p(i) =a\leq
p^a b-(p-1)=i-(p-1)$, de donde $A\geq (i-1)-i+(p-1)\geq
p-2>0$.

En resumen, si $x=1+p$, entonces $x^{p^k}=1+p^{k+1}+
sp^{k+2}$ para alg\'un $s\in {\ma Z}$ y en particular
$x^{p^k}\equiv 1\bmod p^n\iff k+1\geq n\iff k\geq n-1$.
Se sigue que el orden de $x \bmod p^n$ es $p^{n-1}$.
Por otro lado tenemos el epimorfismo natural
\begin{eqnarray*}
U_{p^n}&\to&U_p\\
\xi\bmod p^n&\mapsto&\xi\bmod p\\
\end{eqnarray*}
de donde tenemos que existe $y\in U_{p^n}$ de orden
$(p-1)$ y por tanto $xy$ es de orden $p^{n-1}(p-1)=\varphi(
p^n)=|U_{p^n}|$ probando que $U_{p^n}$ es un grupo
c{\'\i}clico para $p>2$, con $p$ primo y $n\in{\ma N}$.

Ahora consideremos el caso $2^n$. Se tiene 
$U_2=\{1\}$; $U_4=({\ma Z}/4{\ma Z})^{\ast}\cong\{\pm 1\}
\cong C_2$. Notemos que $U_8$ no es c{\'\i}clico:
$U_8=\{1,3,5,7\}$ y todos sus elementos son de orden $2$:
$1^2\equiv 3^2\equiv 5^2\equiv 7^2\equiv 1\bmod 8$, es decir,
$U_8\cong C_2\times C_2$.

Para $n\geq 3$ se tiene la sucesi\'on exacta
\begin{equation}\label{Eq1.2.1.3'}
1\longrightarrow D_{2^n,4}\longto
 U_{2^n}\stackrel{\varphi}{\longrightarrow} U_4\longrightarrow 1
\end{equation}
donde $\varphi$ es el epimorfismo natural y $D_{2^n,4}:=
\ker \varphi =\{x\bmod 2^n\mid x\equiv 1\bmod 4\}$.

Se tiene en particular que $5=1+2^2\in D_{2^n,4}$ y de
manera similar como antes, es decir considerando las potencias
$(1+2^2)^{2^k}$, se tiene que $o(5\bmod 2^n)=2^{n-2}$ y
en particular $D_{2^n,4}$ es un grupo c{\'\i}clico de orden $2^{n-2}$.

Ahora, para $x\in U_{2^n}$, si $x\in D_{2^n,4}=\langle 5\rangle$
se tiene que $o(x)|2^{n-2}$. Si $x\notin D_{2^n,4}$ entonces
$x\equiv 3\bmod 4$. Escribamos $x=3+2^2 a$. Entonces
$x^2=9+24 a+2^4 a^2=1+2^3 t$ de donde obtenemos, como en
la primera parte, que $o(x^2)|2^{n-3}$. De aqu{\'\i} se sigue
que $o(x)|2^{n-2}$ y que todo elemento $x\in U_{2^n}$
tiene orden menor o igual a $2^{n-2}$ por lo que para $n\geq 3$,
$U_{2^n}$ no es un grupo c{\'\i}clico lo cual tambi\'en se sigue
del hecho de que existe un epimorfismo natural $U_{2^n}\to
U_8$ y de que $U_8$ no es c{\'\i}clico, pero de esta forma
obtuvimos un elemento de orden exactamente $2^{n-2}$.

Puesto que $|U_{2^n}|=2^{n-1}$ se sigue que $U_{2^n}\cong
C_{2^{n-2}}\times C_2$. Este isomorfismo tambi\'en se sigue
de que la sucesi\'on (\ref{Eq1.2.1.3'}) se escinde: sea
$\psi:U_4\to U_{2^n}$, $\psi(3)=\psi(-1)=2^n-1$ y $(\varphi
\circ \psi)(3)=\varphi(2^n-1)=\varphi(3+2^n-4)=3$ y en particular
\[
U_{2^n}\cong D_{2^n,4}\times U_4.
\]

Finalmente, para $n\in{\ma N}$, $n\geq 3$, se tiene que si la
descomposici\'on en primos de $n$ es $n=2^mp_1^{\alpha_1}
\cdots p_r^{\alpha_r}$ con $p_1,\ldots, p_r$ primos
impares distintos, $m\geq 0$, $\alpha_i\geq0$, $r\geq 0$, 
entonces por el Teorema Chino del Residuo se tiene
\[
U_n\cong U_{2^m}\times U_{p_1^{\alpha_1}}\times \cdots
\times U_{p_r^{\alpha_r}}.
\]

	Resumimos nuestra discusi\'on anterior en el siguiente
resultado.

\begin{teorema}\label{T1.2.1.3} Sea $n\geq 3$, $n=2^mp_1^{\alpha_1}
\cdots p_r^{\alpha_r}$ en su descomposici\'on en primos. 
Entonces, si $m\geq 2$
\begin{equation}\label{Eq1.2.1}
U_n\cong C_2\times C_{2^m-2}\times C_{p_1-1}\times
\cdots \times C_{p_r-1}\times C_{p_1^{\alpha_1-1}}
\times\cdots\times
C_{p_r^{\alpha_r-1}}.
\end{equation}

Si $m=0,1$, entonces
\begin{equation}\label{Eq1.2.2}
U_n\cong C_{p_1-1}\times
\cdots \times C_{p_r-1}\times C_{p_1^{\alpha_1-1}}\times
\cdots \times C_{p_r^{\alpha_r-1}}.
\end{equation}

En particular $U_n$ es un grupo c{\'\i}clico $\iff n=2,4, p^{\alpha},
2p^{\alpha}$ con $p$ un n\'umero primo impar, $\alpha\geq 1$.
\end{teorema}

\begin{proof}
La ciclicidad se sigue del hecho de que si $p_i$ es impar,
entonces $p_i-1$ es par. $\fin$
\end{proof}

Recordemos que si $\cic n{}$ es un campo ciclot\'omico, entonces
$n\not\equiv 2\bmod 4$. Entonces se sigue

\begin{corolario}\label{C1.2.1.4}
La extensi\'on $\cic n{}/{\ma Q}$ es c{\'\i}clica para $n=4$ y para
$n=p^{\alpha}$ con $p$ un primo impar y $\alpha\geq 1$. $\fin$
\end{corolario}

M\'as adelante estudiaremos con m\'as detalle la correspondencia
de Galois entre los subcampos de $\cic n{}$ y los subgrupos
de $U_n$.

Recordemos que dado un campo num\'erico $K/{\ma Q}$,
$[K:{\ma Q}]=n<\infty$, {\em entonces el anillo\index{anillo
de enteros} de enteros} ${\cal O}_K$
de $K$ se define por
\[
{\cal O}_K=\{\alpha\in K\mid \Irr(\alpha,x,{\ma Q})\in{\ma Z}[x]\}.
\]
Equivalentemente, ${\cal O}_K$ es la cerradura entera de 
${\ma Z}$ en $K$ y ${\cal O}_K$ es un ${\ma Z}$--m\'odulo
libre de rango $n=[K:{\ma Q}]$.

Una {\em base entera\index{base entera}} $\{\alpha_1,\ldots,
\alpha_n\}$ es una base de ${\cal O}_K$ como ${\ma Z}$--m\'odulo,
es decir,
\[
{\cal O}_K\cong {\ma Z}\alpha_1\oplus\cdots\oplus {\ma Z}
\alpha_n.
\]
Finalmente, el {\em discriminante\index{discriminante}}
de $K$ se define por 
\begin{gather*}
\delta_K:=\det\Big(\alpha_i^{\sigma}\Big)_{\substack{1\leq i\leq n\\
\sigma\in T}}^{2}\in {\ma Z}\\
\intertext{donde $T$ es el conjunto de encajes de $K$ en
${\ma C}$,}
T:=\{\sigma\colon K\longto {\ma C}\mid \sigma \text{\ es
monomorfismo de campos}\}.
\end{gather*}

M\'as precisamente, si $T=\{\sigma_1,\ldots, \sigma_n\}$, entonces
sea
\[
C=\left(
\begin{array}{ccc}
\alpha_1^{\sigma_1}&\cdots&\alpha_n^{\sigma_1}\\
\vdots&&\vdots\\
\alpha_1^{\sigma_n}&\cdots&\alpha_n^{\sigma_n}
\end{array}
\right)
\qquad \text{y}\qquad \delta_K=\det C^2.
\]

Como de costumbre escribimos $n=r_1+2r_2$ donde
$r_1$ es el n\'umero de elementos $\sigma\in T$ tales que
$\sigma(K)\subseteq {\ma R}$, los cuales se llaman 
{\em encajes reales\index{encajes reales}}, y $2r_2$
es el n\'umero de elementos $\sigma \in T$ tales que
$\sigma(K)\nsubseteq {\ma R}$ los cuales se llaman
{\em encajes complejos\index{encajes complejos}} y son un
n\'umero par pues si $\sigma(K)\nsubseteq {\ma R}$ entonces
$\overline{\sigma(K)}\nsubseteq {\ma R}$, donde
$\overline{\sigma(K)}$ denota conjugaci\'on compleja.

\begin{teorema}\label{T1.2.1.4} Para cualquier campo
num\'erico, el signo del discriminante $\delta_K$ es $(-1)^{r_2}$.
\end{teorema}

\begin{proof}
Sea $C=\big(\alpha_i^{\sigma_j}\big)_{
1\leq i, j\leq n}$ con la notaci\'on anterior. Tomando la
matriz conjugada de $C$ la cual consiste en conjugar cada 
elemento de $C$, se tiene $\det \overline{C}=\det(\overline{
\alpha}_i^{\sigma_j})=(-1)^{r_2}\det C$ pues si $\sigma_j$ es
real, entonces $\overline{\sigma_j(\alpha_i)}=\sigma_j(
\alpha_i)$ y la fila respectiva permanece sin cambios y en
el caso en que $\sigma_j$ es complejo se intercambian las
filas $\sigma_j(\alpha_i)$ con $\overline{\sigma_j(\alpha_i)}$
y por cada permutaci\'on de filas hay un cambio de signo.
Se sigue que
\begin{gather*}
0<|\det C|^2 = (\overline{\det C})(\det C)=(-1)^{r_2}(\det C)^2=
(-1)^{r_2}\det C^2=(-1)^{r_2}\delta_K. \tag*{$\fin$}
\end{gather*}

\end{proof}

Notemos que cuando $K/{\ma Q}$ es Galois, $T=\Gal(K/{\ma Q})$
y $\sigma(K)=K$ para todo $\sigma \in T$. En particular
$r_1=0$ si $K\nsubseteq {\ma R}$ y en cuyo caso
$r_2=\frac{n}{2}$, en donde $n=[K:{\ma Q}]$ o $r_2=0$ si
$K\subseteq {\ma R}$ y en cuyo caso $r_1=n$. En particular,
si $K=\cic n{}$, $K\subseteq {\ma R}\iff K={\ma Q}, n=0, 1$. 
Por tanto $r_1=0$ y $r_2=\frac{\varphi(n)}{2}$. Si 
$n=2^m p_1^{\alpha_1}\cdots p_r^{\alpha_r}$, entonces 
$\varphi(n)/2$ es par excepto cuando $r=0, m=2$ o $r=1, m=0$
y $p=p_1$ es primo impar congruente con $3$ m\'odulo $4$.
Es decir

\begin{proposicion}\label{P1.2.1.5} Si $K=\cic n{}$ entonces
$\delta_K$ es positivo excepto para $\cic 4{}$ y para
$\cic p{\alpha}$ con $p$ n\'umero primo tal que
$p\equiv 3\bmod 4$. $\fin$
\end{proposicion}

La definici\'on que hemos usado para el discriminante $\delta_K$
es como un n\'umero entero. Recordemos como definimos
el discriminante como un ideal. En general, consideremos
un dominio Dedekind $A$ y sea $K$ una extensi\'on finita
y separable de $E:=\coc A$.

Sea $B$ la cerradura entera de $A$ en $K$: $B:=\{
\alpha\in K\mid \Irr(\alpha, x,E)\in A[x]\}$. Entonces $B$
es un dominio Dedekind (\cite[Cap. 1, Theorem 6.1]{Jan73}).
Esto \'ultimo se cumple a\'un cuando $K/E$ no sea separable.
Puesto que $K/E$ es separable, la traza $\Tr=\Tr_{K/E}\colon
K\to E$ es suprayectiva. El mapeo
\[
\varphi\colon K\times K\to E\qquad \text{dado por}\qquad
\varphi(x,y):=\Tr(xy)
\]
es $E$--bilineal y no degenerado, esto es, si $\Tr(xy)=0$ para
toda $y\in K$ entonces $x=0$ y si $\Tr(xy)=0$ para toda $x\in
K$, entonces $y=0$.

Se define $B^{\ast}:=\{x\in K\mid \Tr(xy)\in A
\text{\ para toda\ }y\in B\}$. Se tiene que $B^{\ast}$
es el m\'odulo complementario $B'$ definido en la
Secci\'on \ref{Sec0.1}. Entonces $B\subseteq B^{\ast}$ y
$B^{\ast}$ es un $B$--m\'odulo fraccionario. El inverso es un
ideal de $B$ llamado el {\em diferente\index{diferente}} de
$B/A$: ${\eu D}_{K/E}={\eu D}_{B/A}:=(B^{\ast})^{-1}$ y la
norma $N_{K/E}({\eu D}_{B/A})$ se llama el
{\em discriminante\index{discriminante}} de $B$ sobre $A$.

En nuestro caso, si $E$ es un campo num\'erico cualquiera
y $K$ es un extensi\'on finita de $E$, tomaremos $A={\cal O}_E$
y se tiene $B={\cal O}_K$ y ${\eu D}_{{\cal O}_K/{\cal O}_E}:=
{\eu D}_{K/E}$ es el {\em diferente\index{diferente}} de $K/E$.
Usaremos para el discriminante la siguiente notaci\'on:
\[
{\eu d}_{K/E}:=N_{K/E}({\eu D}_{K/E}).
\]

En el caso particular de $E={\ma Q}$, pondremos ${\eu d}_K:=
{\eu d}_{K/{\ma Q}}$ y ${\eu d}_K=\langle \delta_k\rangle$.

En general, cuando tenemos un campo num\'erico $K$ y $K=
{\ma Q}(\alpha)$ con $\alpha\in{\cal O}_K$, entonces
$\{1,\alpha,\ldots,\alpha^{n-1}\}$ es una base de $K/{\ma Q}$
donde $n=[K:{\ma Q}]$ y se tiene ${\ma Z}[\alpha]\subseteq 
{\cal O}_K$. Es raro que tengamos ${\ma Z}[\alpha]={\cal O}_K$
para alg\'un $\alpha\in{\cal O}_K$. Veamos que este es el caso
cuando $K=\cic n{}$ para $n\in{\ma N}$.

\begin{proposicion}\label{P1.2.1.6}
Sean $p$ un n\'umero primo y $m\in{\ma N}$. Entonces ${\ma Z}[
\zeta_{p^m}]$ es el anillo de enteros de $\cic pm$, es decir,
${\cal O}_{\cic pm}={\ma Z}[\zeta_{p^m}]$.
\end{proposicion}

\begin{proof}
Puesto que $\zeta_{p^m}\in{\cal O}_{\cic pm}$ se tiene 
${\ma Z}[\zeta_{p^m}]\subseteq {\cal O}_{\cic pm}$. Ahora bien,
dado $\alpha\in {\cal O}_{\cic pm}$, puesto que $\{
1,\zeta_{p^m},\ldots, \zeta_{p^m}^{\varphi(p^m)-1}\}$ es base
de $\cic pm/{\ma Q}$, se tiene que
\begin{equation}\label{Eq1.2.3}
\alpha=\sum_{i=0}^{\varphi(p^m)-1}a_i\zeta_{p^m}^i\qquad
\text{para}\qquad a_i\in{\ma Q}.
\end{equation}
Nuestro objetivo es probar que $a_i\in{\ma Z}$ y el resultado
se seguir\'a.

Primero recordemos que 
\begin{gather*}
\begin{align*}
\psi_{p^m}(x)&=\psi_p(x^{p^m-1})=
\prod_{\substack{i=0\\ \mcd(i,p)=1}}^{p^m-1}(x-\zeta_{p^m}^i)\\
&=(x^{p^m-1})^{p-1}+(x^{p^m-1})^{p-2}+\cdots + x^{p^m-1}+1
\end{align*}
\intertext{y en particular}
\psi_{p^m}(1)=p=\prod_{\substack{i=0\\
\mcd(i,p)=1}}^{p^m-1}(1-\zeta_{p^m}^i).
\end{gather*}

Sean $i,j$ primos relativos a $p$. Entonces existe $t\in{\ma Z}$
tal que $it\equiv j\bmod p^m$. En particular se sigue que
\begin{align*}
\frac{\zeta_{p^m}^j-1}{\zeta_{p^m}^i-1}=
\frac{\zeta_{p^m}^{it}-1}{\zeta_{p^m}^i-1}&=
(\zeta_{p^m}^i)^{t-1}+ (\zeta_{p^m}^i)^{t-2}+\cdots+
(\zeta_{p^m}^i)+1\\
&\in {\ma Z}[\zeta_{p^m}]\subseteq {\cal O}:={\cal O}_{\cic pm}.
\end{align*}

An\'alogamente $\frac{\zeta_{p^m}^i-1}{\zeta_{p^m}^j-1}\in
{\ma Z}[\zeta_{p^m}]$. De esto se sigue que existe $u\in
{\cal O}^{\ast}$ tal que
\[
1-\zeta_{p^m}^i= u(1-\zeta_{p^m}^j)
\]
y a nivel de ideales de ${\cal O}$ se tiene $\langle
1-\zeta_{p^m}^i\rangle=\langle 1-\zeta_{p^m}^j\rangle$ para
cualesquiera $i,j$ primos relativos a $p$. Por lo tanto
si definimos ${\eu p}:=\langle 1-\zeta_{p^m}\rangle$ se sigue
que
\[
\langle \psi_{p^m}(1)\rangle =\langle p\rangle =
\prod_{\substack{i=0\\ \mcd(i,p)=1}}^{p^m-1}\langle 1-\zeta_{p^m}^i
\rangle = {\eu p}^{\varphi(p^m)}.
\]

Puesto que $\varphi(p^m)=[\cic pm:{\ma Q}]$, se sigue que
$p$ es totalmente ramificado en $\cic pm/{\ma Q}$ y
${\eu p}=\langle 1-\zeta_{p^m}\rangle$ es un ideal primo de
${\cal O}$.

Se define $v:=v_{{\eu p}}$ la valuaci\'on correspondiente a
${\eu p}$, es decir, si $\alpha\in K^{\ast}$, $\alpha=\frac{a}{b}$ con
$a,b\in{\cal O}$ se tiene que $\langle \alpha\rangle = {\eu p}^n
{\eu a}$ con ${\eu a}$ un ideal fraccionario de ${\cal O}$ primo
relativo a ${\eu p}$ y entonces $v_{{\eu p}}(\alpha):=n$.
Se tiene que $v_{{\eu p}}$ cumple las mismas propiedades
de $v_p$ (ver Definici\'on \ref{D1.2.1.1}).

Se tiene que $v(p)=\varphi(p^m)$, $v({\eu p})=1$ y $v(1-\zeta_{
p^m}^i)=1$ para toda $\mcd (i,p)=1$.

Ahora bien, puesto que $\cic pm={\ma Q}(1-\zeta_{p^m})$ se
tiene que $\{1,1-\zeta_{p^m},(1-\zeta_{p^m})^2,\ldots, 
(1-\zeta_{p^m})^{\varphi(p^m)-1}\}$ es una base de $\cic pm/
{\ma Q}$. Puesto que $\alpha\in {\cal O}$, se tiene que
$v_{{\eu p}}(\alpha)\geq 0$. Escribamos
\begin{equation}\label{Eq1.2.4}
\alpha=\sum_{\substack{i=0\\ \mcd(i,p)=1}}^{\varphi(p^m)-1}
b_i(1-\zeta_{p^m})^i
\end{equation}
con $b_i\in{\ma Q}$. Ahora bien, puesto que $b_i\in {\ma Q}$
se tiene que $v(b_i)\equiv 0\bmod \varphi(p^m)$ pues
$v(p)=\varphi(p^m)$. Adem\'as $v((1-\zeta_{p^m})^i)=i
v(1-\zeta_{p^m})=i$ por tanto si $b_i,b_j\neq 0$ y $i\neq j$,
se tiene que $v((1-\zeta_{p^m})^i)\neq v((1-\zeta_{p^m})^j)$
de donde se sigue que
\[
0\leq v(\alpha)=\min_{i, i\neq 0} \{v(b_i(1-\zeta_{p^m})^i)\}\leq
v(b_i)+i \text{\ para toda $i$ con\ } b_i\neq 0.
\]
Por tanto $v(b_i)\geq -i\in [[1-\varphi(p^m),0]]$ y $v(b_i)\equiv
0\bmod \varphi(p^m)$ lo cual implica que $v(b_i)\geq 0$ para
toda $0\leq i\leq \varphi(p^m)-1$. Esto nos dice en particular
que $b_i$ se puede escribir en la forma $b_i=\frac{c_i}{d_i}$
con $c_i,d_i\in{\ma Z}$ primos relativos y $p\nmid d_i$ ya que
si $p|d_i$ entonces $p\nmid c_i$ y $v(b_i)=-v(d_i)<0$ lo cual
es absurdo.

Desarrollando la ecuaci\'on (\ref{Eq1.2.4}) y regresando a
nuestra expresi\'on original (\ref{Eq1.2.3}) se tiene
\[
\alpha=\sum_{i=0}^{\varphi(p^m)-1}a_i\zeta_{p^m}^i
\]
con $v(a_i)\geq 0$. Nuevamente esto significa que si
$a_i=\frac{\gamma_i}{\beta_i}$ con $\gamma_i,\beta_i\in
{\ma Z}$, con $\mcd (\gamma_i,\beta_i)=1$, $p\nmid \beta_i$.
Nuestro objetivo es probar que si alg\'un n\'umero primo
$q$ divide a $\beta_i$, entonces $q$ necesariamente debe
ser $p$. Esto \'ultimo, junto con lo que hemos probado de que
$p\nmid \beta_i$ implican que $\beta_i=1$ y que $a_i\in{\ma Z}$
como deseamos.

Tenemos que $G=G_{p^m}=\Gal(\cic pm/{\ma Q})\cong
U_{p^m}=({\ma Z}/p^m{\ma Z})^{\ast}$ e identificamos cada
$\sigma\in G$ con $c\in U_{p^m}$ donde $\sigma(\zeta_{p^m})=
\zeta_{p^m}^c$. De esta forma tenemos para $\sigma\in G$ y
denotando $\zeta:=\zeta_{p^m}$
\begin{equation}\label{Eq1.2.5}
\alpha^{\sigma}=\sum_{i=0}^{\varphi(p^m)-1}a_i\zeta_{p^m}^{ci}.
\end{equation}

Formando el vector columna $\big(\alpha^{\sigma}\big)_{\sigma
\in G}=\left(\begin{array}{c}\alpha\\ \vdots \\ \alpha^{\sigma}\\
\vdots\end{array}\right)$ y usando la ecuaci\'on (\ref{Eq1.2.5})
obtenemos la igualdad (donde el t\'ermino general $\sigma\in G$
lo identificamos con $c\in U_{p^m}$):
\begin{align}\label{Eq1.2.6}
\left(\begin{array}{c}\alpha\\ \vdots \\ \alpha^{\sigma}\\
\vdots\end{array}\right) &=\left(
\begin{array}{ccccc}
1& \zeta &\zeta^2&\cdots &\zeta^{\varphi(p^m)-1}\\
\cdots&\cdots&\cdots&\cdots&\cdots \\
1& \zeta^c &\zeta^{2c}&\cdots &(\zeta^c)^{\varphi(p^m)-1}\\
\cdots&\cdots&\cdots&\cdots&\cdots
\end{array}\right)
\left(\begin{array}{c}a_0\\
\vdots \\ a_{\varphi(p^m)-1}\end{array}\right)\nonumber \\
&=A \left(\begin{array}{c}a_0\\ 
\vdots \\ a_{\varphi(p^m)-1}\end{array}\right)
\quad \text{donde}\quad A=\big(\zeta^{cj}\big)_{\substack{
c\in U_{p^m}\\ 0\leq j\leq \varphi(p^m)-1}}
\end{align}

$A$ es una matriz cuadrada $\varphi(p^m)\times \varphi(p^m)$
con coeficientes en ${\ma Z}[\zeta]$. M\'as generalmente
se tiene que si $B$ es la {\em matriz de Vandermonde\index{matriz
de Vandermonde}}
\[
B:= \left(\begin{array}{cccc}
1&x_1&\ldots&x_1^{r-1}\\
1&x_2&\ldots&x_2^{r-1}\\
\vdots&\vdots&\cdots&\vdots\\
1&x_r&\ldots&x_r^{r-1}
\end{array}\right)
\]
entonces $\det B =\prod_{i<j}(x_i-x_j)$ (damos una demostraci\'on
al final de la proposici\'on).

Podemos ordenar $U_{p^m}=\{c_1,\ldots, c_{\varphi(p^m)}\}$ y 
poniendo $x_i=\zeta^{c_i}$, $1\leq i\leq \varphi(p^m)$ y $r=\varphi(
p^m)$ se tendr\'a que $A= \left(\begin{array}{cccc}
1&x_1&\ldots&x_1^{r-1}\\
1&x_2&\ldots&x_2^{r-1}\\
\vdots&\vdots&\cdots&\vdots\\
1&x_r&\ldots&x_r^{r-1}
\end{array}\right)$ y si $\Adj A$ denota la matriz adjunta de
$A$, se tiene $A^{-1}=\frac{1}{\det A}(\Adj A)$.
Ahora bien, $\det A=\prod_{i<j}(\zeta^{c_i}-\zeta^{c_j})$
y puesto que $\zeta^{c_i}-\zeta^{c_j}=\zeta^{c_i}(1-\zeta^{c_j-
c_i})=\zeta^{c_i} u_{ij}(1-\zeta)^{c_i-c_j}$ donde $u_{ij}$ es
una unidad de ${\cal O}$ se sigue que
$\det A=u(1-\zeta)^t$ para algunos $u\in{\cal O}^{\ast}$ y
$t\in{\ma Z}$. En particular tenemos de (\ref{Eq1.2.6}) que
$\left(\begin{array}{c}a_0\\ \vdots \\ a_{\varphi(p^m)-1}
\end{array}\right)= A^{-1} \left(\begin{array}{c}
\alpha\\ \vdots \\ \alpha^{\sigma}\\ \vdots \end{array}\right)$
lo cual implica que 
\begin{equation}\label{Eq1.2.7}
a_i=\frac{({\text{entero algebraico}})}{u(1-\zeta)^t}.
\end{equation}
Puesto que $p=v(1-\zeta)^{\varphi(p^m)}$ con $v\in
{\cal O}^{\ast}$, multiplicando por cierta unidad $w\in {\cal O}^{
\ast}$ y $(1-\zeta)^s$ para alg\'un $s$, se tiene que
$a_i=\frac{(\text{entero algebraico})}{p^w}\in{\ma Q}$. Por tanto
el \'unico primo $q$ que puede dividir a $\beta_i$ es $p$
y el resultado se sigue. $\fin$

\end{proof}

Ahora probamos que $\det \left(\begin{array}{cccc}
1&x_1&\ldots&x_1^{r-1}\\
1&x_2&\ldots&x_2^{r-1}\\
\vdots&\vdots&\cdots&\vdots\\
1&x_r&\ldots&x_r^{r-1}
\end{array}\right)= \prod_{i<j}(x_i-x_j)$.
Consideremos variables arbitrarias $X, X_2, \ldots, X_r$ y
consideremos el polinomio $f(X):=\det \left(\begin{array}{cccc}
1&X&\ldots&X^{r-1}\\
1&X_2&\ldots&X_2^{r-1}\\
\vdots&\vdots&\cdots&\vdots\\
1&X_r&\ldots&X_r^{r-1}
\end{array}\right)\in F[X]$ donde $F$ es el campo de 
las funciones racionales en las variables $X_2,\ldots, X_r$.
Entonces $f(X)$ es de grado $r-1$. Puesto que $f(X_2)=\cdots
=f(X_{r-1})=0$, se sigue que $f(X)=D \prod_{j=2}^{r}(X-X_j)$
donde $D$ es el coeficiente l{\'\i}der de $f$. Entonces
$D= \det \left(\begin{array}{cccc}
1&X_2&\ldots&X_2^{r-2}\\
1&X_3&\ldots&X_3^{r-2}\\
\vdots&\vdots&\cdots&\vdots\\
1&X_r&\ldots&X_r^{r-2}
\end{array}\right)$. Por hip\'otesis de inducci\'on en $r$, se
sigue que $D=\prod_{2\leq i<j\leq r}(X_i-X_j)$.
Definiendo $x_1:=X_1,x_2:=X_2,\ldots, x_r:=X_r$ se sigue
el resultado.

Notemos que puesto que $\{1,\zeta,\ldots, \zeta^{\varphi(p^m)-1}
\}$ es una base entera de ${\cal O}={\cal O}_{\cic pm}$ sobre
${\ma Z}$, se sigue que 
\begin{equation}\label{Eq1.2.8}
\delta_{\cic pm}=(-1)^{\varphi(p^m)/2}\det A^2.
\end{equation}
Calculando $\det A^2$ se seguir\'a el discriminante.

\begin{proposicion}\label{P1.2.1.7}
Con las notaciones anteriores, se tiene
\[
\det A^2=\pm p^{p^{m-1}(mp-m-1)}.
\]
\end{proposicion}

\begin{proof}
Para $c\in {\ma Z}$ se tiene $1-\zeta^c=-\zeta^c(1-\zeta^{-c})$
por lo que tenemos
\[
\det A =\pm \prod_{\substack{0<k<j< p^m\\
p\nmid kj}}(\zeta^j-\zeta^k)=u_1 \prod_{\substack{0<k<j<
p^m\\ p\nmid kj}}(1-\zeta^{k-j})
\]
donde $u_1\in{\cal O}^{\ast}$ y $\det A^2= u_2  \prod\limits_{\substack{0<k,j<p^m\\
p\nmid kj, k\neq j}}(1-\zeta^{k-j})$ con $u_2\in {\cal O}^{\ast}$.

Puesto que $\det A^2\in {\ma Z}$ y $u_3(1-\zeta)^{\varphi(p^m)}=
p$ con $u_3\in{\cal O}^{\ast}$, se tiene que $\det A^2=\pm p^s$
para alg\'un $s\in {\ma N}\cup \{0\}$. Para calcular $s$ 
consideremos nuevamente $v$ la valuaci\'on asociada a
${\eu p}=\langle 1-\zeta\rangle$.

Ahora bien ${\eu p}^{\varphi(
p^m)}=\langle p\rangle$, es decir $v(p)=\varphi(p^m)$ y para
$1\leq n\leq m$, $1-\zeta_{p^n}=1-\zeta_{p^m}^{p^{m-n}}=
u_4 (1-\zeta_{p^m})^{p^{m-n}}$, $u_4\in{\cal O}^{\ast}$.
M\'as precisamente se tiene
\begin{gather*}
\langle 1-\zeta_{p^m}^{p^{m-n}}\rangle ^{\varphi(p^n)}=
\langle 1-\zeta_{p^n}\rangle ^{\varphi(p^n)}=\langle p\rangle =\langle 1-\zeta_{p^m}\rangle ^{\varphi(p^m)},
\intertext{por lo que}
\langle 1-\zeta_{p^m}^{p^{n-m}}\rangle =\langle 1-\zeta_{p^m}\rangle ^{\varphi(p^m)/
\varphi(p^n)}=\langle 1-\zeta_{p^m}\rangle ^{p^{m-n}}
\intertext{lo cual implica}
1-\zeta_{p^m}^{p^{m-n}}=u_4(1-\zeta_{p^m})^{p^{m-n}}
\end{gather*}
con $u_4\in\*{\mc O}$.

Por lo tanto $v(1-\zeta_{p^n})=p^{m-n}$. Agrupando los
t\'erminos $1-\zeta_{p^m}^{p^{m-n}}$ en la expresi\'on de
$\det A^2$ obtenemos
\[
\det A^2=\omega \prod_{n=1}^m (1-\zeta_{p^m})^{s_{m-n}}
\]
donde $\omega\in{\cal O}^{\ast}$ y $s_i=|\{(k,j)\mid k\equiv j
\bmod p^i, k\not\equiv j\bmod p^{i+1}, p\nmid kj, 0<k,j
<p^m\}|$.
Una vez calculado $s_i$, se tendr\'a que
\[
\varphi(p^m)s=v(\det A^2)=\sum_{n=1}^m s_{m-n}p^{n-m}.
\]

Sea $0<j<\varphi(p^m)$, $p\nmid j$ fijo y sea
\[
c_j^{(n)}=\{k\mid 0<k<p^m, p\nmid k, k\equiv j\bmod p^n\}.
\]
Consideremos el epimorfismo natural
\begin{eqnarray*}
\varphi_n\colon U_{p^m}&\longto& U_{p^n}\\
x\bmod p^m&\longmapsto&x\bmod p^n
\end{eqnarray*}
con $\ker \varphi_n=D_{p^m,p^n}=\{x\in U_{p^m}\mid
x\equiv 1\bmod p^n\}$. Entonces, $c_j^{(n)}=\varphi_n^{-1}(
\{j\})$ y por tanto $|c_j^{(n)}|=|\ker \varphi_n|=\frac{\varphi(
p^m)}{\varphi(p^n)}=p^{m-n}$. Por otro lado tenemos que
\begin{gather*}
s_i=\Big| \bigcup_{j\in U_{p^m}}(c_j^{(i)}-c_j^{(i+1)})\Big|, \quad
1\leq i\leq m-2.\\
\intertext{Se sigue que}
s_i=\varphi(p^m)\cdot (p^{m-i}-p^{m-i-1})=\varphi(p^m)
p^{m-i-1}(p-1), \quad 1\leq i\leq m-2.
\end{gather*}

Para $i=m-1$, $s_{m-1}=|\{(k,j)\mid 0<k,j<p^m, p\nmid kj, k\neq j,
k\equiv j\bmod p^{m-1}\}|= \varphi(p^m)(p-1)$.

Finalmente, para $i=0$, consideremos $s_0$. Se tiene 
\[
s_0=|\{(k,j)\mid k\not\equiv j\bmod p\}|.
\]
Fijando $j$, $0<j<p^m$, $p\nmid j$, tenemos que existen
$p^{m-1}$ elementos $k$ con $0<k<p^m$ tales que
$k\equiv j\bmod p$. Por lo tanto hay $\varphi(p^m)-p^{m-1}$
elementos $k$ tales $k\not\equiv j \bmod p$. Puesto que 
existen $\varphi(p^m)$ tales elementos $j$, se sigue que
\[
s_0=\varphi(p^m)(\varphi(p^m)-p^{m-1})=
\varphi(p^m)(p^{m-1}(p-2)).
\]

Por tanto
\begin{align*}
\varphi(p^m) s&=v(\det A^2)=\sum_{i=0}^{m-1} s_i p^i =\\
&=\varphi(p^m)\big[p^{m-1}(p-2) +\sum_{i=1}^{m-1} p^{m-i-1}
(p-1) p^i\big] =\\
&=\varphi(p^m)[(p-2)p^{m-1}+(m-1)(p-1)p^{m-1}]=\\
&=\varphi(p^m)p^{m-1}(mp-m-1)
\end{align*}
de donde se sigue el resultado. $\fin$

\end{proof}

\begin{corolario}\label{C1.2.1.8}
Para un n\'umero primo $p$ y $m\in{\ma N}$ se tiene
$\delta_{\cic pm}=(-1)^{\varphi(p^m)/2} p^{p^{m-1}(mp-m-1)}$.
$\fin$
\end{corolario}

Otra demostraci\'on del Corolario \ref{C1.2.1.8} se encuentra
en la Subsecci\'on \ref{S12.3.1}.

Puesto que los primos finitos ramificados en $K/{\ma Q}$ son los
primos que dividen a $\delta_K$, donde $K$ es cualquier
campo num\'erico (ver Teorema \ref{T0.1.1} y 
Corolario \ref{C1.3.6}), se tiene

\begin{corolario}\label{C1.2.1.9}{\ }

\l
\item Si $p$ es un n\'umero primo, $m\in{\ma N}$, entonces
el \'unico primo finito ramificado en $\cic pm /{\ma Q}$ es $p$
y es totalmente ramificado.

\item Para $n\in{\ma N}$, $n\geq 3$, $n\not\equiv 2\bmod 4$,
un n\'umero primo $p$ se ramifica en $\cic n{}$ si y solamente
si $p$ divide a $n$.
\end{list}
\end{corolario}

\begin{proof}{\ }

\l
\item Es el Corolario \ref{C1.2.1.8}.

\item Del Corolario \ref{C1.2.10}, $\cic n{}=\prod_{i=1}^r
\cic {p_i}{\alpha_i}$ donde la descomposici\'on en primos de
$n$ est\'a dada por $n=p_1^{\alpha_1}\cdots p_r^{\alpha_r}$.
Finalmente un n\'umero primo $p$ se ramifica en $\cic n{}/
{\ma Q}\iff$ se ramifica en alg\'un $\cic {p_i}{\alpha_i}/{\ma Q}
\iff p=p_i$ para alg\'un $1\leq i\leq r\iff p|n$. $\fin$

\end{list}

\end{proof}

Para estudiar el anillo ${\cal O}_{\cic n{}}$ con $n\in {\ma N}$
arbitrario, salvo que $n\not\equiv 2\bmod 4$, veamos que
bajo ciertas condiciones para dos extensiones $E/{\ma Q}$ y
$K/{\ma Q}$, se tiene que ${\cal O}_E {\cal O}_K={\cal O}_{EK}$.

Primero recordemos que en una torre de campos num\'ericos
$K\subseteq L\subseteq M$, el diferente es multiplicativo
(ver Secci\'on \ref{Sec0.1}), esto es,
\[
{\eu D}_{M/K}={\eu D}_{M/L}\cdot \con_{L/M}{\eu D}_{L/K}
\]
donde $\con_{L/M}$ denota a la conorma de $L$ a $M$, es
decir, si $\pK$ es un ideal primo de ${\cal O}_L$ y el ideal
extendido $\pK{\cal O}_M=\pL_1^{e_1}\cdots \pL_g^{e_g}$,
entonces $\con_{L/M}\pK:=\pL_1^{e_1}\cdots \pL_g^{e_g}$
y el mapeo se extiende para cualquier ideal fraccionario
${\eu a}=\pK_1^{\alpha_1}\cdots \pK_r^{\alpha_r}$,
$\pK_1,\ldots, \pK_r$ ideales primos de ${\cal O}_L$ y
$\alpha_1, \ldots, \alpha_r\in {\ma Z}$. Se define
$\con_{L/M}{\eu a}:= \big(\con_{L/M} \pK_1\big)^{\alpha_1}
\cdots \big(\con_{L/M} \pK_r\big)^{\alpha_r}$.

Notemos que si ${\eu a}$ es un ideal fraccionario de $L$,
entonces $N_{L/M}\con_{L/M}{\eu a} = {\eu a}^{[M:L]}$.

\begin{window}[0,r,$\xymatrix{F\ar@{-}[r]\ar@{-}[d]&FH\ar@{-}[d]\\
E\ar@{-}[r]&H}$,{}]
Por otro lado decimos que dos extensiones $F/E$ y $H/E$ son
{\em linealmente disjuntas\index{linealmente disjuntas}} sobre
$E$ si una base de $F/E$ es tambi\'en una base de $FH/H$.
Equivalentemente, para extensiones finitas, si $[F:E]=
[FH:H]$.
\end{window}

\begin{teorema}\label{T1.2.1.10}
Sean $K$ y $E$ dos campos num\'ericos. Supongamos que los
discriminantes de $K$ y $E$ son primos relativos y que 
$K$ y $E$ son linealmente disjuntos sobre ${\ma Q}$. Entonces
\[
{\cal O}_{KE}={\cal O}_K {\cal O}_E \qquad \text{y}\qquad
{\eu d}_{KE}={\eu d}_K^{[E:{\ma Q}]} {\eu d}_E^{[K:{\ma Q}]}.
\]
\end{teorema}

\begin{proof}
Los diferentes satisfacen
\begin{equation}\label{Eq1.2.9}
{\eu D}_{KE/{\ma Q}}={\eu D}_{KE/K}\cdot\con_{K/KE}
{\eu D}_{K/{\ma Q}}={\eu D}_{KE/E}\cdot \con_{E/KE}
{\eu D}_{E/{\ma Q}}.
\end{equation}

Ahora bien, se tiene por hip\'otesis que $\con_{K/KE}
{\eu D}_{K/{\ma Q}}$ y $\con_{E/KE} {\eu D}_{E/{\ma Q}}$
son primos relativos. Por otro lado tenemos
\[
{\eu D}_{KE/K}|\con_{E/KE}{\eu D}_{E/{\ma Q}}\quad
\text{y}\quad {\eu D}_{KE/E}|\con_{K/KE}{\eu D}_{K/{\ma Q}}
\]
de donde obtenemos que ${\eu D}_{KE/K}$ y ${\eu D}_{KE/E}$
son primos relativos. De la ecuaci\'on (\ref{Eq1.2.9}) se sigue
que
\begin{equation}\label{Eq1.2.10}
{\eu D}_{KE/E}=\con_{K/KE}{\eu D}_{K/{\ma Q}} \quad \text{y}
\quad {\eu D}_{KE/K}=\con_{E/KE}{\eu D}_{E/{\ma Q}}.
\end{equation}

Para obtener una base de ${\cal O}_{KE}$ sobre ${\ma Z}$,
usaremos las bases complementarias, es decir, las del diferente
inverso.
Sea $W$ una base de ${\cal O}_K/{\ma Q}$ y $V$ una base de
${\cal O}_E/{\ma Z}$. Sea $W'$ la base dual de $W$ con 
respecto a la traza, esto es, si $W=\{w_1,\ldots, w_t\}$
entonces $W'=\{w_1',\ldots, w_t'\}\subseteq K$ satisface que
$\Tr_{K/{\ma Q}}(w_i'w_j)=\delta_{ij}$ para $1\leq i,j\leq t$. En
otras palabras
\[
{\ma Z}w_1'+\cdots+{\ma Z}w_t'=\{x\in K\mid \Tr_{K/{\ma Q}}
(x{\cal O}_K)\subseteq {\ma Z}\}={\eu D}_{K/{\ma Q}}^{-1}
\]
y $W'$ genera a ${\eu D}_{K/{\ma Q}}^{-1}$ sobre ${\ma Z}$.
Se sigue que $W'$ genera ${\eu D}_{KE/K}^{-1}$ sobre
${\cal O}_K$. Entonces el doble dual, $(W')'=W$ genera
${\cal O}_{KE}$ sobre ${\cal O}_E$. Por lo tanto
${\cal O}_{KE}={\cal O}_E(W)={\cal O}_E{\cal O}_K$.

Finalmente
\begin{align*}
{\eu d}_{KE/{\ma Q}}&= N_{KE/{\ma Q}}({\eu D}_{KE/{\ma Q}})
= N_{KE/{\ma Q}}(\con_{K/KE}{\eu D}_{K/{\ma Q}} \con_{E/KE}
{\eu D}_{E/{\ma Q}})\\
&=N_{KE/{\ma Q}}({\eu D}_{K/{\ma Q}}) N_{KE/{\ma Q}}(
{\eu D}_{E/{\ma Q}})\\
&= N_{K/{\ma Q}}(N_{KE/K} {\eu D}_{K/{\ma Q}}) 
N_{E/{\ma Q}}(N_{KE/E} {\eu D}_{E/{\ma Q}})\\
&= N_{K/{\ma Q}}({\eu D}_{K/{\ma Q}}^{[KE:K]})
N_{E/{\ma Q}}({\eu D}_{E/{\ma Q}}^{[KE:E]})
= {\eu d}_{K/{\ma Q}}^{[E:{\ma Q}]} 
{\eu d}_{E/{\ma Q}}^{[K:{\ma Q}]}.
\tag*{$\fin$}
\end{align*}

\end{proof}

El Teorema \ref{T1.2.1.10} nos facilita de manera substancial
el c\'alculo de una base entera de ${\cal O}_{\cic n{}}$ y el
discriminante $\delta_{\cic n{}}$ donde $n\in{\ma N}$, $n>1$,
$n\not\equiv 2\bmod 4$.

\begin{teorema}\label{T1.2.1.11}
Se tiene que para $n\in{\ma N}$, ${\ma Z}[\zeta_n]$ es el anillo
de enteros de $\cic n{}$, es decir, $\{1,\zeta_n,\ldots, \zeta_n^{
\varphi(n)-1}\}$ es una base entera de ${\cal O}_{\cic n{}}$.
\end{teorema}

\begin{proof}
Esto es consecuencia del Teorema \ref{T1.2.1.10}: Se tiene que si
$p,q$ son dos primos distintos, $\alpha,\beta\in{\ma N}$,
entonces $\cic p{\alpha}$ y $\cic q{\beta}$ son linealmente
disjuntos y que los discriminantes son primos relativos
(Corolarios \ref{C1.2.10} y \ref{C1.2.1.8}). Por lo tanto
si $n=p_1^{\alpha_1}\cdots p_r^{\alpha_r}$, se tiene que
$\cic n{}=\prod_{i=1}^r\cic {p_1}{\alpha_i}$ y por lo tanto
\begin{gather*}
{\cal O}_{\cic n{}}=\prod_{i=1}^r {\cal O}_{\cic {p_i}{\alpha_i}}=
\prod_{i=1}^r {\ma Z}[\zeta_{p_i^{\alpha_i}}]={\ma Z}[\zeta_n].
\tag*{$\fin$}
\end{gather*}

\end{proof}

\begin{teorema}\label{T1.2.1.12} Si $n=p_1^{\alpha_1}\cdots
p_r^{\alpha_r}$ es la descomposici\'on en primos de $n$,
entonces el diferente de $\cic n{}/{\ma Q}$ est\'a dado por
\[
{\eu D}_{\cic n{}/{\ma Z}}= \prod_{i=1}^r \con_{\cic {p_i}{\alpha_i}/
\cic n{}} {\eu p}_i^{p_i^{\alpha_i-1}(p_i\alpha_i-\alpha_i-1)}
\]
donde ${\eu p}_i$ es ideal primo de $\cic {p_i}{\alpha_i}$ dado
por ${\eu p}_i=\langle 1-\zeta_{p_i^{\alpha_i}}\rangle$.
\end{teorema}

\begin{proof}
Puesto que ${\eu d}_{\cic {p_i}{\alpha_i}}=\langle p_i\rangle^{
p_i^{\alpha_i}(p_i\alpha_i-\alpha_i-1)}$ (Corolario \ref{C1.2.1.8})
y $p_i$ es totalmente ramificado en $\cic {p_i}{\alpha_i}/{\ma Q}$:
${\eu p}_i^{\varphi(p_i^{\alpha_i}})=\langle p_i\rangle$, entonces
el grado relativo $f({\eu p}_i|p_i)$ es $1$ de donde obtenemos
que $N_{\cic {p_i}{\alpha_i}/{\ma Q}}{\eu p}_i =\langle p_i\rangle$.

En particular los diferentes ${\eu D}_{\cic {p_i}{\alpha_i}/{\ma Q}}$
son primos relativos a pares y ${\eu D}_{\cic {p_i}{\alpha_i}/
{\ma Q}}={\eu p}_i^{p_i^{\alpha_i-1}(p_i\alpha_i-\alpha_i-1)}$.
La conclusi\'on se sigue de la multiplicatividad de los diferentes.
$\fin$
\end{proof}

\begin{corolario}\label{C1.2.1.13} Se tiene para $n\in{\ma N}$,
$n>1$, $n\not\equiv 2\bmod 4$ que
\[
\delta_{\cic n{}}=(-1)^{\varphi(n)/2}\frac{n^{\varphi(n)}}{\prod_{p|n}
p^{\varphi(n)/(p-1)}}.
\]
\end{corolario}

\begin{proof}
Una primera demostraci\'on es usando el Teorema \ref{T1.2.1.12}
y el hecho de que ${\eu d}_{\cic n{}}=N_{\cic n{}/{\ma Q}}
{\eu D}_{\cic n{}/{\ma Q}}$.

Una segunda demostraci\'on es usando el Teorema \ref{T1.2.1.10}.
En ese caso, la igualdad ${\eu d}_{KE}={\eu d}_K^{[E:{\ma Q}]}
{\eu d}_E^{[K:{\ma Q}]}$ implica $|\delta_{KE}|=|\delta_K|^{
[E:{\ma Q}]} |\delta_E|^{[K:{\ma Q}]}$, la cual a su vez se puede
poner en forma aditiva tomando logaritmos:
\[
\log |\delta_{KE}| =[E:{\ma Q}]\log |\delta_K| + [K:{\ma Q}]
|\delta_E|.
\]
Dividiendo entre $[KE:{\ma Q}]$, obtenemos
\begin{equation}\label{Eq1.2.11}
\frac{\log |\delta_{KE}|}{[KE:{\ma Q}]}=
\frac{\log |\delta_{K}|}{[K:{\ma Q}]}+
\frac{\log |\delta_{E}|}{[E:{\ma Q}]}.
\end{equation}

La ventaja de la expresi\'on (\ref{Eq1.2.11}) es que es f\'acilmente
generalizable a una composici\'on de un n\'umero finito de
campos. En nuestro caso tenemos $\cic n{}=\prod_{i=1}^r
\cic {p_i}{\alpha_i}$ de donde
\begin{align*}
\frac{\log |\delta_{\cic n{}}|}{\varphi(n)}&=
\frac{\log |\delta_{\cic n{}}|}{[\cic n{}:{\ma Q}]}=
\sum_{i=1}^r \frac{\log |\delta_{\cic {p_i}{\alpha_i}}|}
{[\cic {p_i}{\alpha_i}:{\ma Q}]}= \sum_{i=1}^r
\frac{\log |\delta_{\cic {p_i}{\alpha_i}}|}{\varphi(p_i^{\alpha_i}}\\
&=\sum_{i=1}^r \frac{p_i^{\alpha_i-1}(p_i\alpha_i-\alpha_i-1)
\log p_i}{p_i^{\alpha_i-1}(p_i-1)}=
\sum_{i=1}^r \Big(\alpha_i-\frac{1}{p_i-1}\Big) \log p_i\\
&=\sum_{i=1}^r \alpha_i \log p_i-\sum_{i=1}^r \frac{\log p_i}{
p_i-1}=\log n - \sum_{p|n}\frac{\log p}{p-1}.\\
\intertext{Por tanto}
\log |\delta_{\cic n{}}|&=\varphi(n) \log n-\sum_{p|n}
\frac{\varphi(n)}{p-1} \log p= \log n^{\varphi(n)} -
\sum_{p|n} \log p^{\varphi(n)/(p-1)} \\
&= \log \frac{n^{\varphi(n)}}{\prod_{p|n} p^{\varphi(n)/(p-1)}}
\end{align*}
de donde, usando el Teorema \ref{T1.2.1.4}, se sigue que
\begin{gather*}
\delta_{\cic n{}} =(-1)^{\varphi(n)/2}\frac{n^{\varphi(n)}}{
\prod_{p|n} p^{\varphi(n)/(p-1)}}. \tag*{$\fin$}
\end{gather*}

\end{proof}

Uno de los problemas centrales que se estudian en cualquier
campo num\'erico, es su grupo de unidades. En el caso en que
$n=p^{\alpha}$ es una potencia de un primo, obtuvimos que
$\langle p\rangle =\langle 1-\zeta_{p^{\alpha}}\rangle^{\varphi(
p^{\alpha})}$ (ver la demostraci\'on de la Proposici\'on 
\ref{P1.2.1.6}). En particular $1-\zeta_{p^{\alpha}}$ no puede
ser unidad en ${\ma Z}[\zeta_{p^{\alpha}}]$. Resulta ser que 
cuando $n$ no es potencia de un primo la historia es
diferente.

Recordemos que $x^n-1=\prod_{d|n} \psi_d(x)$. Sea
$n=p_1^{\alpha_1}\cdots p_r^{\alpha_r}$, $r\geq 2$ un
n\'umero natural que no es potencia de un n\'umero primo.
Entonces
\[
f_n(x)=\frac{x^n-1}{x-1}=x^{n-1}+\cdots+x+1= \prod_{
\substack{d|n\\ d\neq 1}}\psi_d(x)=\prod_{j=1}^{n-1}
(x-\zeta_n^j).
\]

Sea ${\cal A}:=\{d\in {\ma N}\mid d|n, d \text{\ no es potencia
de primo}\}$. Se tiene que $n\in{\cal A}$ y ${\cal A}\neq
\emptyset$. Por otro lado tenemos
\begin{gather*}
n=f_n(1)=\prod_{j=1}^{n-1} (1-\zeta_n^j)=\prod_{i=1}^r\Big(
\prod_{\beta_i=1}^{\alpha_i} \psi_{p_i^{\beta_i}}(1)\Big)\cdot
\prod_{d\in{\cal A}}\psi_d(1).\\
\intertext{Ahora bien $\psi_{p_i^{\beta_i}}(1)=p$ para $1\leq \beta_i
\leq \alpha_i$. Por lo tanto}
\prod_{i=1}^r\prod_{\beta_i=1}^{\alpha_i}\psi_{p_i^{\beta_i}}(1)=
\prod_{i=1}^rp_i^{\alpha_i}=n,\\
\intertext{de donde se sigue que}
1=\prod_{d\in{\cal A}}\psi_d(1)=\prod_{d\in{\cal A}}\prod_{\mcd(j,d)=1}
(1-\zeta_d^j).
\end{gather*}

Puesto que $d=n\in{\cal A}$, $1-\zeta_n$ aparece en el producto
$\prod_{d\in{\cal A}}\psi_d(1)$ y por lo tanto $1-\zeta_n$ es
unidad. En consecuencia, tenemos que $\pm 1=N_{\cic n{}/
{\ma Q}}(1-\zeta_n)=\prod_{\mcd(j,n)=1}(1-\zeta_n^j)$.

Notemos el siguiente hecho general.

\begin{teorema}\label{T1.2.1.14}
Sea $K/{\ma Q}$ una extensi\'on finita de Galois tal que
la restricci\'on a $K$ de la conjugaci\'on compleja no 
es trivial, es decir, $J|_K\neq \Id_K$, donde $J$ denota
la conjugaci\'on compleja. Entonces $N_{K/{\ma Q}}
K^{\ast}\subseteq {\ma Q}^+:= \{x\in{\ma Q}\mid x>0\}$.
\end{teorema}

\begin{proof}

Sea $G:=\Gal(K/{\ma Q})$ y sea $H:=\{1, J|_K\}$, $|H|=2$ y
$H<G$. Sea $X\subseteq G$ un conjunto de representantes
de las clases derechas de $G$ m\'odulo $H$. Entonces
$HX=G$ y $G=X\uplus JX$.

Sea $\xi\in K^{\ast}$. Se tiene $N_{K/{\ma Q}}(\xi)=\prod_{
\sigma\in G}\sigma \xi = \big(\prod_{\sigma\in X}\sigma \xi\big)
\big(\prod_{\sigma\in X} J(\sigma \xi)\big)=\alpha\overline{\alpha}=
|\alpha|^2>0$. $\fin$

\end{proof}

Resumiendo el desarrollo anterior, tenemos

\begin{teorema}\label{T1.2.1.15}
Si $n$ no es potencia de un n\'umero primo, entonces
$1-\zeta_n$ es unidad en ${\cal O}_{\cic n{}}={\ma Z}[\zeta_n]$
y $N_{\cic n{}/{\ma Q}}(1-\zeta_n)= \prod_{(j,n)=1}(1-\zeta_n^j)
=1$. $\fin$
\end{teorema}

%% file: Capitulo4.tex
\chapter{Teorema de Kronecker--Weber\index{teorema de Kronecker--Weber}}\label{Ch1}

En este cap{\'\i}tulo presentamos la demostraci\'on del Teorema de Kronecker--Weber
usando los grupos de ramificaci\'on (ver Secci\'on \ref{Sec0.3}).

\section{El teorema y su demostraci\'on}\label{Sec1.1}

El Teorema de 
Kronecker--Weber\index{teorema!Kronecker--Weber}\index{teorema
de Kronecker--Weber}\index{Kronecker--Weber!teorema de $\sim$} 
establece que toda extensi\'on abeliana
de ${\ma Q}$ est\'a contenida en alg\'un campo ciclot\'omico
$\cic n{}$. El teorema fue originalmente afirmado por Kronecker
en 1853 \cite{Kro1853}.  Sin embargo la prueba 
estaba incompleta. El mismo Kronecker reconoci\'o
que hab{\'\i}a problemas con el primo $p=2$.
En 1886, Weber casi complet\'o 
 la demostraci\'on \cite{Web1886}. Todav{\'\i}a la prueba ten{\'\i}a
 una laguna, la cual no fue notada sino hasta 95 a\~nos despu\'es
 por Olaf Neumann \cite{Neu81}.
 Finalmente en 1896, D. Hilbert dio una nueva demostraci\'on
 completa de este resultado la cual se bas\'o en los grupos
 de ramificaci\'on \cite{Hil1896}. Esta es la primera prueba completa
 correcta del resultado, aunque el mismo Hilbert no lo supo
 pues consider\'o que la prueba de Weber estaba completa.
 
 En la actualidad hay muchas demostraciones, varias de ellas
 elementales, del Teorema de Kronecker--Weber: reducci\'on al
 caso local, usando Teor{\'\i}a de Campos de Clase en la cual
 se usa que un primo $p$ se descompone totalmente en
 $\cic n{} \iff p\equiv 1\bmod n$ lo cual prueba que $\cic n{}$ es el
 campo de clase de rayos correspondiente al ``m\'odulus'' $(n) \infty$
 y toda extensi\'on abeliana de ${\ma Q}$ tiene que estar
 contenida en alguno de estos campos, etc., 
 ver Corolario \ref{CClaseC4.5.13}.
 
 Aqu{\'\i} presentamos una prueba basada en los grupos de ramificaci\'on.
 
 Antes que nada recordemos que el Teorema de 
 Minkowski\index{teorema!Minkowski}\index{teorema de Minkowski}
 establece que si $1<[K:{\ma Q}]<\infty$, entonces existe
 un n\'umero primo $p$ el cual es ramificado en $K/{\ma Q}$.
 Como veremos a continuaci\'on, el caso central es que si
 $K/{\ma Q}$ es una extensi\'on c{\'\i}clica de grado $p$ con
 $p>2$ un n\'umero primo y \'unicamente se ramifica $p$, entonces
 $K\subseteq \cic p2$, esto es, $K$ es la \'unica subextensi\'on
 de $\cic p2$ de grado $p$ sobre ${\ma Q}$. Similarmente necesitamos
 el an\'alogo para el caso $p=2$.
 
 Como primer paso, estudiamos el caso moderadamente ramificado.
 
 \begin{proposicion}\label{P4.1} Sea $K/{\ma Q}$ una extensi\'on
 abeliana tal que $p\in{\ma N}$ es moderadamente ramificado
 en $K$. Entonces existe una extensi\'on $L$ de ${\ma Q}$
 y un subcampo $F\subseteq \cic p{}$  tal que:
 \lasa
 
\item Todo primo no ramificado en $K$ es no ramificado en $L$.

\item $p$ no se ramifica en $L$.

\item $FK=FL$.

\end{list}
\[
\xymatrix{
&&FK=FL\ar@{-}[dl]&\cic p{}\ar@{-}[lld]|!{[l];[ldd]}\hole\\
K\ar@{-}[urr]\ar@{-}[d]&F\\
{\ma Q}\ar@{-}[rr]\ar@{-}[ur]&&L\ar@{-}[uu]
}
\]
\end{proposicion}

\begin{proof}
Sea $\pK$ un primo en $K$ sobre $p$. Puesto que $p$
es moderadamente ramificado, el primer grupo de ramificaci\'on
$G_1$ de $p$ es trivial. Puesto que $K/{\ma Q}$ es abeliana
se sigue que el grupo de inercia $I(\pK|p)$ est\'a contenido
de manera natural en ${\ma F}_p^{\ast}=({\ma Z}/p{\ma Z})^{\ast}$,
Proposici\'on \ref{P1.3.16},
es decir, el {\'\i}ndice de ramificaci\'on $e$ de $\pK/p$ divide a $p-1$.
En particular $p\neq 2$.

Sea $\zeta:=\zeta_p$ y sea $F$ el \'unico subcampo de
$\cic p{}$ de grado $e$ sobre ${\ma Q}$: $[F:{\ma Q}]=e$.
Ahora $p$ es moderadamente ramificado en $F/{\ma Q}$ puesto
que $p\nmid e$ y totalmente ramificado en $F/{\ma Q}$ puesto
que $p$ es totalmente ramificado en $\cic p{}/{\ma Q}$.

Sea ${\eu q}$ el \'unico primo de $F$ sobre $p$. Consideremos
$FK=E$. Sea $\pL$ un primo en $E$ sobre $\pK$ y sea
$I:=I(\pL|\pK)$ el grupo de inercia. Se tiene $\pL\cap 
{\ma Z}=(p)$ y $\pL\cap {\cal O}_F={\eu q}$. Sea
$L:=E^I$ el subcampo de $E$ fijo bajo $I$.
\[
\xymatrix{
{\eu q}\ar@{--}[rrrr]\ar@{--}[dddd]\ar@{--}[rd]&&&&\pL
\ar@{--}[dddd]\ar@{--}\ar@{--}[ld]\\
&F\ar@{-}[rr]\ar@{-}[dd]&&FK=E\ar@{-}[dd]\\
&&L\ar@{-}[dl]\ar@{-}[ru]_I\\
&{\ma Q}\ar@{-}[rr]\ar@{--}[ld]&&K\ar@{--}[dr]\\
p\ar@{--}[rrrr]&&&&\pK
}
\]

Veamos que se cumplen las tres propiedades requeridas en
la proposici\'on. Primero, sea $q\in {\ma Z}$ un primo no ramificado
en $K/{\ma Q}$. Entonces $q\neq p$ y puesto que $q$ no es
ramificado en $F/{\ma Q}$, se sigue que $q$ no es ramificado
en $E=FK/{\ma Q}$. Por tanto $q$ no se ramifica en $L\subseteq
E$. Esto es la primera propiedad de $L$.

Ahora, como $L$ es el campo de inercia de $p$ en $E/{\ma Q}$, $p$
no es ramificado en $L/{\ma Q}$ lo cual prueba la segunda
propiedad de $L$.

Queda por demostrar que $FK=FL=E$. Por supuesto $FL\subseteq E$.
Se tiene que $p$ es no ramificado en $L/{\ma Q}$ y totalmente
ramificado en $E/L$. Por tanto $F\cap L={\ma Q}$.
\[
\xymatrix{
\hbox{$\begin{array}{ccc} \pK\\ & K& {}\\ {}
\end{array}$}\ar@{-}[rr]\ar@{-}[dd]_e
&&\hbox{$\begin{array}{ccc} &\pL\\ E=
FK&{}\\{}\end{array}$}\ar@{-}[dd]\ar@{-}[ld]\\
&L\\ \hbox{$\begin{array}{ccc} \\ &{\ma Q}&{}\\ p\end{array}$}\ar@{-}[rr]_e
\ar@{-}[ru]&& \hbox{$\begin{array}{ccc}\\ &F\\ && {\eu q}\end{array}$}
}
\]

Ahora bien, por el Lema de Abhyankar (ver Teorema
\ref{T2.3.B}) se tiene que
\[
e(\pL|p)=\mcm[e(\pK|p),e({\eu q}|p)]=\mcm[e,e]=e,
\]
esto es, $[E:L]=[E:E^I]=|I|=e$. Por otro lado, como $p$ es totalmente
ramificado en $F/{\ma Q}$ y no ramificado en $L/{\ma Q}$, se
tiene que $F$ y $L$ son linealmente disjuntos sobre ${\ma Q}$
y en particular $[FL:L]=[F:{\ma Q}]=e$. Por tanto
$[E:L]=[FL:L]$. Se sigue que
$FL=E=FK$. $\fin$
\end{proof}

\begin{corolario}\label{C4.2}
Sea $K/{\ma Q}$ una extensi\'on abeliana tal que todo primo
ramificado en $K/{\ma Q}$ es moderadamente ramificado.
Entonces existe $m\in{\ma N}$ tal que $K\subseteq \cic m{}$.
M\'as a\'un, si $p_1,\ldots,p_r$ son los primos ramificados en
$K$, se puede seleccionar $m=p_1\cdots p_r$.
\end{corolario}

\begin{proof} Como hicimos notar en la demostraci\'on
de la Proposici\'on \ref{P4.1}, $2\notin\{p_1,\ldots,p_r\}$.
Sean $p_1,\ldots,p_r$ todos los primos ramificados en
$K/{\ma Q}$. Por la Proposici\'on \ref{P4.1}, existe
$F_1\subseteq \cic {p_1}{}$ y $L_1$ tales que 
$p_2,\ldots, p_r$ son los primos ramificados en $L_1/{\ma Q}$
y $F_1K=F_1L_1$.

Aplicando lo mismo para $L_1$, existen $F_2\subseteq \cic {p_2}{}$
y $L_2$ tales que $p_3,\ldots, p_r$ son los primos ramificados
en $L_2$ y $F_2L_1=F_2L_2$. Continuando este proceso
para $i=2,\ldots, r-1$ se tiene que existen $F_i\subseteq
\cic {p_i}{}$ y $L_i$ tales que $p_{i+1},\dots p_r$ son los primos
ramificados en $L_i$ y $F_iL_{i-1}= F_iL_i$.

Finalmente, existe $F_r\subseteq \cic {p_r}{}$ y $L_r$ tales
que $L_r/{\ma Q}$ es no ramificada y $F_rL_{r-1}=F_rL_r$.
Por el Teorema de Minkowski, se tiene que $L_r={\ma Q}$ y 
por lo tanto $F_rL_{r-1}=F_r$, lo cual implica que $L_{r-1}
\subseteq F_r$. Del hecho $F_{r-1}L_{r-2}=F_{r-1}L_{r-1}
\subseteq F_r F_{r-1}$ se sigue que $L_{r-2}\subseteq
F_rF_{r-1}$.
Continuando este proceso, se sigue que $L_1\subseteq
F_rF_{r-1}\cdots F_2$, por tanto
\begin{gather*}
K\subseteq KF_1=LF_1\subseteq
F_r F_{r-1}\cdots F_2 F_1\subseteq \cic {p_1 \cdots p_r}{}. 
\tag*{$\fin$}
\end{gather*}
\end{proof}

\section{Caso central: ramificaci\'on salvaje}\label{S1.2}

Nuestro siguiente objetivo es probar un resultado an\'alogo
al Corolario \ref{C4.2} en el caso de ramificaci\'on salvaje.
Primero hagamos un caso particular. Separamos en dos
casos: $p$ par y $p$ impar.

\begin{proposicion}\label{P4.3} Las \'unicas extensiones 
cuadr\'aticas de ${\ma Q}$ con discriminante una potencia
de $2$ son ${\ma Q}(\sqrt{-1})=
{\ma Q}(i)=\cic 4{}$, ${\ma Q}(\sqrt{2})$ y
${\ma Q}(\sqrt{-2})={\ma Q}(\sqrt{2} i)$ todas ellas contenidas
en $\cic 8{}$.
\end{proposicion}

\begin{proof}
Sea $K={\ma Q}(\sqrt{d})$ con $d$ libre de cuadrados. Puesto
que $\delta_K=d, 4d$ y por hip\'otesis $\delta_K=
\pm 2^{\alpha}$, se sigue que necesariamente que $d=\pm 1, \pm 2$.
El caso $d=1$ queda descartado y por lo tanto $K$ es uno de
los tipos mencionados en la proposici\'on.
El rec{\'\i}proco es igual.

Ahora bien, se tiene $\zeta_8=\sqrt{\sqrt{-1}}=\sqrt[4]{-1} = 
\cos \pi/4 + i\sen \pi/4= \frac{\sqrt{2}}{2}(1+i)$. Por lo tanto
$\zeta_8+\bar{\zeta}_8=\sqrt{2}$, $\zeta_8-\bar{\zeta}_8=i\sqrt{2}$
de donde se sigue que $\cic 4{}, {\ma Q}(\sqrt{2}),
{\ma Q}(\sqrt{-2})\subseteq \cic 8{}$. $\fin$

\end{proof}

\begin{proposicion}\label{P1.4.1*} Sea $K/{\ma Q}$ una extensi\'on c{\'\i}clica de
grado $p$ sobre ${\ma Q}$, con $p$ un n\'umero primo impar, tal que el 
\'unico primo ramificado es $p$. Entonces el diferente de la extensi\'on es igual a
${\eu D}_K=\pK^{2(p-1)}$ donde $\pK$ es el \'unico ideal primo de $K$ sobre $p$.
\end{proposicion}

\begin{proof}
Puesto que $[K:{\ma Q}]=p$ y $p$ es ramificado, necesariamente tenemos que
$e(\pK|p)=p$. Sea $\pi\in\pK\setminus
\pK^2$. Se tiene que $v_{\pK}(\pi)=1$ por lo que
$\pi\notin{\ma Q}$ puesto que para cualquier $\alpha\in{\ma Q}$ tenemos que
$v_{\pK}(\alpha)=e(\pK|p)v_p(\alpha) = pv_{\pK}(\alpha)\neq 1$.

Se sigue que $K={\ma Q}(\pi)$. Sea 
\[
f(x):=\Irr(\pi,x,{\ma Q})=x^p+a_{p-1}x^{p-1}+
\cdots+a_1 x+a_0.
\]
Entonces $f(x)\in{\ma Z}[x]$ puesto que $\pi\in\pK\subseteq \AE K$, es decir,
$\pi$ es entero.
Como $\pi^p+a_{p-1}\pi^{p-1}+\cdots+a_1\pi+a_0=0$ con $a_0\neq 0$, se tiene
para $0\leq i\leq p-1$ que 
\begin{gather*}
v_{\pK}(a_i\pi^i)=pv_{p}(a_i)+i.\\
\intertext{Si para alg\'un $i$ 
tuvi\'esemos que $p\nmid a_i$, entonces se tendr{\'\i}a que $v_{\pK}(a_i)=0$ y
$v_{\pK}(a_i\pi^i)=i$ y entonces}
\infty =v_{\pK}(0)=v_{\pK}(\pi^p+a_{p-1}
\pi^{p-1}+\cdots+a_1\pi+a_0)=\min_{
\substack{a_i\neq 0\\ p\nmid a_i}}\{i\}\neq \infty
\end{gather*}
lo cual es absurdo.

Por el Corolario \ref{C1.3.7'} y de \cite[Proposition 5.6.7]{Vil2006}
(ver tambi\'en \cite[Corollary 5.7.20]{Vil2006}) se tiene que 
\[
{\eu D}_K=\langle f'(\pi)\rangle =\pK^k\quad
\text{donde}\quad k=\sum_{i=0}^{\infty} (|G_i|-1)
\]
 donde $G_i$ es el $i$--\'esimo grupo
de ramificaci\'on correspondiente a $\pK$ sobre $p$. Puesto que $G_i\subseteq
G:=\Gal(K|{\ma Q})$ y $|G|=p$, se tiene que $|G_i|=p \text{\ \'o\ } 1$ y por lo
tanto $|G_i|-1=p-1 \text{\ \'o\ } 0$ de donde se sigue que $p-1|k$.

Ahora, se tiene que $f'(\pi)=p\pi^{p-1}+(p-1)a_{p-1}\pi^{p-2}+\cdots+a_2\pi+a_1$
por lo que para las $i$, $1\leq i\leq p$ tales para $a_i\neq 0$, donde ponemos
 $a_p:=1$, se tiene que
\[
v_{\pK}(ia_i\pi^{i-1})=v_{\pK}(i)+v_{\pK}(a_i)+i-1\equiv (i-1)\bmod p.
\]
Por lo tanto si para $i\neq j$ tenemos que $a_i\neq 0$ y $a_j\neq 0$ entonces
$v_{\pK}(ia_i\pi^{i-1})\neq v_{\pK}(ja_j\pi^{j-1})$. Se sigue que
\[
k=v_{\pK}({\eu D}_K)=v_{\pK}(f'(\pi))=\min_{\substack{1\leq i\leq p\\a_i\neq 0}}
\{v_{\pK}(ia_i\pi^{i-1})\}= v_{\pK}(i_0)+v_{\pK}(a_{i_0})+i_0-1.
\]

Si $i_0=p$ se tendr{\'\i}a $v_{\pK}(pa_p\pi^{p-1})=v_{\pK}(p\pi^{p-1})=2p-1\not\equiv
0\bmod p-1$. Por tanto necesariamente $1\leq i_0\leq p-1$ y $v_{\pK}(a_{
i_0}\pi^{i_0-1})=pv_p(a_{i_0})+i_0-1<2p-1=
v_{\pK}(pa_p\pi^{p-1})$. Por otro lado $p|a_{i_0}$ lo cual
implica que $v_p(a_{i_0})=t\geq 1$ de donde obtenemos que
\[
2p-1>v_{\pK}(a_{i_0}\pi^{i_0-1})=tp+i_0-1
\]
lo cual implica que $tp<2p-1+1-i_0=2p-i_0$. Se sigue que $t=1$ y que
$k=p+i_0-1<2p-1$.

Por ser ramificaci\'on salvaje, se tiene que $k>p-1$ por lo que $p-1<p+i_0-1$.
Si acaso tuvi\'esemos que $i_0<p-1$, entonces $p+i_0-1<2(p-1)$ y por ende
$p-1<k<2(p-1)$ lo cual contradice que $p-1|k$. Se sigue que $i_0=p-1$ y que
$k=2(p-1)$. $\fin$
\end{proof}

El siguiente resultado nos prueba que \'unicamente existe una extensi\'on
c{\'\i}clica de ${\ma Q}$ de grado $p$ ramificada \'unicamente en $p$
en donde $p$ es un primo impar.

\begin{proposicion}\label{P1.4.2*} Sea $p$ un n\'umero primo con $p>2$
y sea $K/{\ma Q}$ una extensi\'on c{\'\i}clica de grado $p$ cuyo \'unico
primo ramificado es $p$. Entonces $K\subseteq {\ma Q}(\zeta_{p^2})$.
Esto es, $K$ es el \'unico subcampo de ${\ma Q}(\zeta_{p^2})$ de
grado $p$ sobre ${\ma Q}$.
\end{proposicion}

\begin{proof}
Primero consideremos un campo $L$ tal que $L/{\ma Q}$ es una
extensi\'on de Galois abeliana y tal que $[L:{\ma Q}]=p^2$ donde
$p$ es el \'unico primo ramificado. Sea $G_0$ el grupo de inercia
respectivo y sea $L^{G_0}=E$. Entonces $p$ no se ramifica
en $E/{\ma Q}$ y por tanto $E/{\ma Q}$ no es ramificada. Se sigue
del Teorema de Minkowski que $E={\ma Q}$ y que $G_0=G:=\Gal
(L/{\ma Q})$. Ahora bien, si $G_1$ es el primer grupo de ramificaci\'on,
entonces por ser la extensi\'on $L/{\ma Q}$ salvajemente ramificada
se tiene que $G_1\neq \{e\}$. Consideremos $F:=L^{G_1}$. Entonces
$p$ es moderadamente ramificado en la extensi\'on $F/{\ma Q}$.
Sin embargo, en el caso de que $F\neq {\ma Q}$ se tendr{\'\i}a que
que $p|[F:{\ma Q}]$ y $p$ necesariamente ser{\'\i}a salvajemente
ramificado, lo cual es absurdo. Se sigue que $G_1=G_0=G$ y
$|G_1|=|G_0|=|G|=p^2$.

Sea $G_r$ el primer grupo de ramificaci\'on tal que $|G_r|<p^2$.
Puesto que $G_{r-1}/G_r$ es un grupo $p$--elemental abeliano
y puesto que $G_0=G_1=G$, se sigue que $r-1\geq 1$, esto es,
$r\geq 2$. Entonces $G_{r-1}/G_r=G/G_r\subseteq {\eu p}^{r-1}/
{\eu p}^r\cong {\cal O}_L/{\eu p}\cong {\ma F}_p\cong {\ma Z}/
p{\ma Z}$ ya que al ser $p$ totalmente ramificado su grado de
inercia es $1$.

Se sigue que $|G_{r-1}/G_r|=|G_{r-1}|/|G_r| = p^2/|G_r|\leq p$
y por lo tanto $|G_r|=p$.

Consideremos $H$ cualquier subgrupo de $G$ de orden $p$
y sea $E:=L^H$. Sea ${\eu P}:={\eu p}\cap {\cal O}_{L^H}$.
Entonces por la Proposici\'on \ref{P1.4.1*} se tiene que
${\eu D}_{L^H/{\ma Q}}={\eu P}^{2(p-1)}$. Por lo tanto
\begin{gather*}
{\eu D}_{L/{\ma Q}}={\eu D}_{L^H/L} \con_{L^H/L} {\eu P}^{
2(p-1)}= {\eu D}_{L/L^H} {\eu p}^{2p(p-1)}.\\
\xymatrix{
&L\ar@{-}[d]\ar@{--}[r]&{\eu p}\ar@{--}[dd]\\
{\ma Q}\ar@{-}[r]\ar@{--}[d]&L^H\ar@{--}[dr]\\
p\ar@{--}[rr]&&{\eu P}}
\end{gather*}

Se sigue que el diferente 
${\eu D}_{L/L^H}={\eu D}_{L/{\ma Q}}{\eu p}^{-2p(p-1)}$
es independiente de $H$. Ahora bien, si $H\neq G_r$,
los grupos de ramificaci\'on de $L/L^H$ son
\begin{gather*}
G_i\cap H=\left\{
\begin{array}{ccl}
H&{\text{si}}&0\leq i\leq r-1\\
1&{\text{si}}&i\geq r
\end{array}
\right.\\
\intertext{de donde}
{\eu D}_{L/L^H}={\eu p}^s, \quad s=\sum_{i=0}^{\infty}
(|G_i\cap H|-1)=r(p-1).
\end{gather*}

Por otro lado, si $H=G_r$, se tiene que 
\[
{\eu D}_{L/L^{G_r}}={\eu p}^t\quad {\text{con}}\quad
t=\sum_{i=0}^{\infty}(|G_i\cap G_r|-1)\geq (r+1)(p-1).
\]

Puesto que ${\eu D}_{L/L^H}$ es independiente de
$H$, necesariamente $H= G_r$, lo que prueba que
$G$ tiene un \'unico subgrupo de orden $p$, 
a saber $G_r$. Por lo tanto $G$ es c{\'\i}clico.

Ahora sean $K$ y $K'$ dos extensiones c{\'\i}clicas de grado
$p$ de ${\ma Q}$ tales que $p$ es el \'unico primo ramificado.
En caso de que $K\neq K'$ se tendr{\'\i}a que $K K'$ es una
extensi\'on de Galois de grado $p^2$ de ${\ma Q}$ donde $p$ es
el \'unico primo ramificado. Por lo anterior se seguir{\'\i}a
que $K K'$ es una extensi\'on c{\'\i}clica de ${\ma Q}$ pero
\[
\Gal(KK'/{\ma Q})\cong \Gal(K/{\ma Q})\times \Gal(K'/{\ma Q})
\cong C_p\times C_p.
\]
Por tanto, necesariamente $K=K'$ de donde se sigue que hay una
\'unica extensi\'on $K/{\ma Q}$ c{\'\i}clica de grado $p$ en donde
\'unicamente $p$ es ramificado.

Finalmente puesto que $\Gal({\ma Q}(\zeta_{p^2})/{\ma Q})\cong
U_{p^2}\cong C_{p-1}\times C_p$, si consideramos $F:=
{\ma Q}(\zeta_{p^2})^H$ donde $H$ es el subgrupo de $U_{p^2}$
de orden $p-1$, se sigue que $\Gal(F/{\ma Q})\cong C_p$ y
$p$ es el \'unico primo finito ramificado, es decir, $F\subseteq {\ma Q}
(\zeta_{p^2})$ es el \'unico campo satisfaciendo estas condiciones.
$\fin$
\end{proof}

\begin{teorema}\label{T1.4.3*} Sea $p$ un primo impar. Sea
$K/{\ma Q}$ una extensi\'on abeliana de grado $p^m$ donde
$p$ es el \'unico primo ramificado. Entonces $K\subseteq {\ma Q}
(\zeta_{p^{m+1}})$ y de hecho $K$ es el campo fijo ${\ma Q}
(\zeta_{p^{m+1}})^H$ donde $H$ es el subgrupo de $U_{p^{m+1}}$
de orden $p-1$: $U_{p^{m+1}}\cong C_{p-1}\times C_{p^m}$.
En particular $K/{\ma Q}$ es una extensi\'on c{\'\i}clica.
\end{teorema}

\begin{proof}
Sea $L={\ma Q}(\zeta_{p^{m+1}})^H$. Entonces $KL$ es una
extensi\'on abeliana de ${\ma Q}$ donde $p$ es el \'unico primo
ramificado. Sea ${\mc G}=\Gal(K/{\ma Q})\cong C_{p^{\alpha_1}}
\times\cdots\times C_{p^{\alpha_t}}$ con $\alpha_1+\cdots
+\alpha_t=m$. Se tiene que
\begin{gather*}
\Gal(KL/{\ma Q})\subseteq \Gal(K/{\ma Q})\times \Gal(L/{\ma Q})
\cong {\mc G}\times C_{p^m} \quad {\text{con}}\quad\\
\begin{array}{rcl}
\varphi\colon\Gal(KL/{\ma Q})&\longrightarrow&
\Gal(K/{\ma Q})\times \Gal(L/{\ma Q})\\
\theta&\longmapsto&(\theta|_K,\theta|_L).
\end{array}\\
\xymatrix{
&KL\ar@{-}[dl]\ar@{-}[rd]\\
K\ar@{-}[dr]\ar@{-}[ddr]&&L\ar@{-}[ld]\ar@{-}[ddl]\\
&K\cap L\ar@{-}[d]\\
&{\ma Q}}
\end{gather*}

Si $KL/{\ma Q}$ no fuese una extensi\'on c{\'\i}clica, entonces
$KL$ contendr{\'\i}a una subextensi\'on de orden $p^2$ 
de tipo $C_p\times C_p$ sobre ${\ma Q}$ lo cual contradir{\'\i}a
la Proposici\'on \ref{P1.4.2*}. Por lo tanto $KL/{\ma Q}$ es c{\'\i}clica.
Puesto que $\Gal(KL/{\ma Q})\subseteq  C_{p^{\alpha_1}}
\times\cdots\times C_{p^{\alpha_t}}\times C_{p^m}$, se sigue
que $\Gal(KL/{\ma Q})\cong C_{p^m}$.  Por lo tanto
$K=L=KL$. $\fin$
\end{proof}

Para el primo par $p=2$ tenemos resultados similares. Primero
consideremos una extensi\'on $[K:{\ma Q}]=2$ y $2$ es el \'unico
primo ramificado. Entonces $K={\ma Q}(\sqrt{2})$, ${\ma Q}(\sqrt{-2})$
o ${\ma Q}(\sqrt{-1})={\ma Q}(i)={\ma Q}(\zeta_4)$.

\begin{teorema}\label{T1.4.4*}
Si $K/{\ma Q}$ es una extensi\'on c{\'\i}clica de grado $2^m$, $m\geq 2$,
y donde $2$ es el \'unico primo finito ramificado,
entonces $K\subseteq \cic 2 {m+2}$. M\'as a\'un $K$ es uno de dos:
$K=\cic 2{m+2}\cap {\ma R}={\ma Q}(\zeta_{2^{m+2}}+
\zeta_{2^{m+2}}^{-1})=:K_m$ o $K={\ma Q}(\zeta_4(\zeta_{2^{m+2}}+
\zeta_{2^{m+2}}^{-1}))={\ma Q}(\zeta_{2^{m+2}}-\zeta_{2^{m+2}}^{-1})$.
\end{teorema}

\begin{proof} Supongamos que $K$ es un campo real, 
y que la extensi\'on $K/{\ma Q}$ no es necesariamente 
c{\'\i}clica, \'unicamente satisfaciendo que 
$K/{\ma Q}$ sea abeliana, que $[K:{\ma Q}]=2^m$ y
que $2$ es el \'unico primo finito ramificado. Entonces $K K_m$
es una extensi\'on real abeliana 
de ${\ma Q}$ donde $2$ es el \'unico primo
finito ramificado. Nuevamente tenemos que
\[
\Gal(K K_m/{\ma Q})\subseteq \Gal(K/{\ma Q})\times \Gal(K_m/
{\ma Q})\cong \Gal(K/{\ma Q})\times C_{2^m}.
\]
Si $KK_m/{\ma Q}$ no fuese c{\'\i}clica entonces $KK_m$ tendr{\'\i}a
una subextensi\'on de tipo $C_2\times C_2$ sobre ${\ma Q}$ pero
esta ser{\'\i}a real lo cual implicar{\'\i}a que ${\ma Q}$ tendr{\'\i}a
tres extensiones cuadr\'aticas reales donde $2$ es el \'unico primo
ramificado. Esto contradice que la \'unica extensi\'on cuadr\'atica
real de ${\ma Q}$, donde $2$ es el \'unico primo
finito ramificado, es ${\ma Q}(\sqrt{2})$. Se sigue que $K=
K_m$ y en particular $K$ es una extensi\'on c{\'\i}clica de ${\ma Q}$
y $K\subseteq \cic 2{m+2}$.

Ahora consideremos $K$ no real y sea $M:=K(i)$. Sea $M^+:=
M\cap {\ma R}$. Se tiene que si $K=M$, entonces
$M=K^+(i)=K$ y $\Gal(K/{\ma Q})\cong \Gal(K^+/{\ma Q})
\times C_2$ lo que implica que $K/{\ma Q}$ no es c\'iclico
debido a que $M^+\cap {\ma Q}(i)={\ma Q}$ y
$M^+{\ma Q}(i)=M^+(i)=M$ puesto que
$[M:M^+]=2$. Por lo tanto $K\neq M$. 

Se tiene que $K^+=K_{m-1}\neq
M^+$ pues, por el caso anterior, $K^+$ es una
extensi\'on de grado $2^{m-1}$ y $[M^+:{\ma Q}]=2^m>
[K_{m-1}:{\ma Q}]=2^{m-1}$,
de donde se sigue que $M^+=K_m$. Puesto que
$M=M^+(i)$, $M=\cic 2{m+2}$ y $\Gal(M/{\ma Q}) \cong
C_2\times C_{2^m}$. Hay tres subcampos de {\'\i}ndice
$2$, a saber, $\cic 2{m+1}$, $K_m$ y ${\ma Q}(\zeta_{2^{m+2}}-
\zeta_{2^{m+2}}^{-1})$. Puesto
que $K/{\ma Q}$ es c{\'i}clica y no real, se tiene que
$K\neq \cic 2{m+1}$ y $K\neq K_m$ por lo que
necesariamente $K={\ma Q}(\zeta_{2^{m+2}}-
\zeta_{2^{m+2}}^{-1})$ (ver Teorema \ref{T10.1.1}).
En particular $K\subseteq \cic 2{m+2}$. $\fin$
\end{proof}

\begin{observacion}\label{R1.4.4**}
El subcampo ${\ma Q}(\zeta_{2^{m+2}}-\zeta_{2^{m+2}}^{-1})$ descrito
en el Teorema \ref{T1.4.4*} es la extensi\'on
${\ma Q}(\zeta_4(\zeta_{2^{m+2}}+\zeta_{2^{m+2}}^{-1}))$:
\[
\xymatrix{K_m\ar@{-}[dd]_2\ar@{-}[rr]^2&&{\ma Q}(\zeta_{2^{m+2}}+
\zeta_{2^{m+2}}^{-1})\ar@{-}[ld]^2\ar@{-}[dd]^2\\
&{\ma Q}(\zeta_4(\zeta_{2^{m+2}}+\zeta_{2^{m+2}}^{-1}))\ar@{-}[ld]^2\\
K_{m-1}\ar@{-}[rr]^2\ar@{-}[dd]_2&& K_{m-1}(\zeta_4)\ar@{-}[dd]^2\\
\\
{\ma Q}\ar@{-}[rr]^2&&{\ma Q}(\zeta_4)}
\]
Ahora 
\begin{align*}
\zeta_4(\zeta_{2^{m+2}}+\zeta_{2^{m+2}}^{-1})&=
\zeta_4\zeta_{2^{m+2}}+\zeta_4 \zeta_{2^{m+2}}^{-1}
=\zeta_{2^{m+2}}^{2^m+1}+\zeta_{2^{m+2}}^{2^m-1}\\
&=\zeta_{2^{m+2}}^{2^m+1}+\zeta_{2^{m+2}}^{2^{m+1}-2^m-1}
=\zeta_{2^{m+2}}^a+\zeta_2 \zeta_{2^{m+2}}^{-a}=
\zeta_{2^{m+2}}^a-\zeta_{2^{m+2}}^{-a}
\end{align*}
 con $\mcd(a,2)=1$,
por lo que ${\ma Q}(\zeta_{2^{m+2}}^a-\zeta_{2^{m+2}}^{-a})
={\ma Q}(\zeta_{2^{m+2}}-\zeta_{2^{m+2}}^{-1})$.
\end{observacion}

Con lo anterior ya tenemos todos los elementos para probar:

\begin{teorema}[Kronecker--Weber\index{teorema!Kronecker--Weber}\index{teorema
de Kronecker--Weber}\index{Kronecker--Weber!teorema de $\sim$}]\label{T1.4.5*}
Sea $K/{\ma Q}$ una extensi\'on abeliana finita. Entonces existe $n\in{\ma N}$
tal que $K\subseteq \cic n{}$.
\end{teorema}

\begin{proof}
Puesto que $K/{\ma Q}$ es una extensi\'on abeliana finita, se tiene
que $\Gal(K/{\ma Q})\cong \oplus_{i=1}^r C_{n_i}$ donde
cada $n_i$ es potencia de un n\'umero primo. Sea 
$K_i:=K^{H_i}$
el campo fijo bajo $H_i:=\oplus_{\substack{j=1\\ j\neq i}}^r C_{n_j}$.
Entonces $K=K_1\cdots K_r$. Si probamos que cada $K_i
\subseteq \cic {m_i}{}$ entonces
\[
K=K_1\cdots K_r\subseteq \cic {m_1}{} \cdots \cic {m_r}{}
\subseteq \cic {m_1\cdots m_r}{}.
\]
Por tanto podemos suponer que $K/{\ma Q}$ es c{\'\i}clica
de grado $p^m$ con $p$ un n\'umero primo.

Por la Proposici\'on \ref{P4.1}, existe una extensi\'on
$L$ de ${\ma Q}$ y $F\subseteq \cic n{}$ para alguna $n$
tal que $FK=FL$ tal que el \'unico primo ramificado en $L/
{\ma Q}$ puede ser $p$.
De hecho $n$ puede ser tomado como $q_1\cdots q_t$
donde los primos ramificados de $K/{\ma Q}$ son
$q_1, \ldots, q_t$ y posiblemente $p$. 

Entonces $L\cap F={\ma Q}$ y $\Gal(LF/{\ma Q})\cong
\Gal(L/{\ma Q})\times \Gal(F/{\ma Q}) \cong \Gal(FK/{\ma Q})$.
\begin{gather*}
\xymatrix{
K\ar@{-}[rr]\ar@{-}[dd]&&FK=FL\ar@{-}[dd]\ar@{-}[dl]\\
&L\ar@{-}[dl]\\
{\ma Q}\ar@{-}[rr]&& F\ar@{-}[r]&\cic {q_1\cdots q_t}{}
}
\\
\xymatrix{
L\ar@{-}[r]\ar@{-}[d]&FL\ar@{-}[d]\\ {\ma Q}\ar@{-}[r]& F
}
\quad
\xymatrix{
&K\ar@{-}[r]\ar@{-}[d]&FK\ar@{-}[d]\\
&K\cap F\ar@{-}[r]\ar@{-}[dl]&F\\ {\ma Q}
}
\quad 
\begin{array}{c}
{\ }\\ \\ \\
[0pt][FK:F]| [K:{\ma Q}]=p^m.
\end{array}
\end{gather*}

Por tanto $[L:{\ma Q}]=[FL:F]=[FK:F]|[K:{\ma Q}]$ lo cual implica
que $L/{\ma Q}$ es un $p$--extensi\'on c\'iclica donde $p$ es
el \'unico posible primo ramificado. 
Por los Teoremas \ref{T1.4.3*} y \ref{T1.4.4*} se
tiene que $L\subseteq \cic p\ell$ para alg\'un $\ell$.
Por tanto $K\subseteq FK=FL\subseteq \cic {q_1\cdots
q_t}{} \cic p\ell = \cic {p^\ell q_1\cdots q_t}{}$. $\fin$
\end{proof}

\begin{observacion}\label{O1.4.6*}
El Teorema de Kronecker--Weber no se cumple para ning\'un campo
num\'erico $K$ diferente de ${\ma Q}$.
\end{observacion}

\begin{proposicion}\label{P1.4.7*}
Sea $K$ un campo num\'erico, $K\neq{\ma Q}$. Entonces existe
una extensi\'on abeliana finita $L/K$ tal que $L\not\subseteq
K{\ma Q}(\zeta_{\infty})=K(\zeta_{\infty})$ donde ${\ma Q}(\zeta_{
\infty})=\bigcup_{n=1}^{\infty}\cic n{}$, $\zeta_{\infty}=\bigcup_{n=1}^{
\infty}\langle \zeta_n\rangle$.
\end{proposicion}

\begin{proof}
Sea $p\in{\ma Q}$ tal que existen dos primos $\pK_1,\pK_2$ en $K$
tales que $\pK_1\neq \pK_2$ y $\pK_i|_{\ma Q}=p, i=1,2$. Hay una infinidad
de tales primos $p$: por ejemplo, si $\tilde{K}$ es la cerradura de
Galois de $K/{\ma Q}$, sea $p$ tal que $p$ se descompone
totalmente en $\tilde{K}/{\ma Q}$. La existencia de tal $p$ se sigue
del Teorema de Densisdad de \v{C}evotarev o, de manera m\'as 
elemental, el n\'umero de primos ramificados es finito (Teorema
\ref{T0.1.1}) y si $p$ es no ramificado, el grupo de descomposici\'on
de $p$ es c\'iclico (ver Secci\'on \ref{Sec0.3}). 

Si $K/{\ma Q}$ no es
c\'iclico cualquier primo no ramificado funciona. 
Si $K/{\ma Q}$
es c\'iclico, por el Teorema de Kronecker--Weber, $K\subseteq \cic n{}$
para alg\'un $n$. Sea $p\equiv 1\bmod n$. Entonces $p$ se descompone
totalmente en $\cic n{}/{\ma Q}$ y en particular en $K$ (ver Teorema
\ref{T8.2} m\'as adelante).

Por el Teorema de Aproximaci\'on de Artin \cite[Theorem 2.5.3]{Vil2006}
(o m\'as simplemente, por propiedades elementales de Dominios Dedekind),
existe $\alpha\in K$ tal que $v_{\pK_1}(\alpha)=1$ y $v_{\pK_2}(\alpha)=0$.

Sea $L=K(\sqrt{\alpha})$ y sea $\pL_1$ es un primo sobre $\pK_1$.
Sea $\beta\in L$ tal que $\beta^2=\alpha\in K$. Entonces
\[
2v_{\pL_1}(\beta)=v_{\pL_1}(\beta^2)=v_{\pL_1}(\alpha)=
e_{L/K}(\pL_1|\pK_1)v_{\pK_1}(\alpha)=e_{L/K}(\pL_1|\pK_1).
\]
Por tanto $e_{L/K}(\pL_1|\pK_1)=2$ y $\pK_1$ es ramificado en $L/K$.

Por otro lado $v_{\pK_2}(\alpha)=0$. Veamos que $\pK_2$ es no
ramificado en $L/K$. En caso contrario, consideremos $\pL_2$ un
primo en $L$ sobre $\pK_2$. Se tiene $2v_{\pL_2}(\beta)=
v_{\pL_2}(\beta^2)=v_{\pL_2}(\alpha)=0$, esto es, $v_{\pL_2}(\beta)=0$.
Por tanto $\beta\in\AE{{\pL_2}}$ es entero en $\pL_2$. Sea $G=\Gal(L/K)=
I(\pL_2|\pK_2)=\langle\sigma\rangle$ el cual es de orden $2$. 
Se tiene $\sigma(\beta)=-\beta$ pues $\sigma(\beta^2)=(\sigma(\beta))^2=
\sigma(\alpha)=\alpha$ por lo que 
$\sigma(\beta)=\pm \sqrt{\alpha}=\pm \beta$
y $\sigma(\beta)\neq \beta$. Por definici\'on de grupo de inercia, $v_{\pL_2}
(\sigma(\beta)-\beta)>0$ pero $v_{\pL_2}(-2\beta)=0$. Se sigue
que $\pK_2$ es no ramificado.

Veamos que $L\not\subseteq K(\zeta_{\infty})$. En caso contrario, se
tendr\'ia que $L\subseteq K(\zeta_n)$ para alg\'un $n$. Sea $M=L\cap
\cic n{}$.
\[
\xymatrix{
\cic n{}\ar@{-}[d]\ar@{-}[r]&K(\zeta_n)\ar@{-}[d]\\
M=L\cap \cic n{}\ar@{-}[d]\ar@{-}[r]&L\ar@{-}[d]\\
K\cap\cic n{}\ar@{-}[r]\ar@{-}[d]&K\\{\ma Q}
}
\]
Como $\pK_1$ es ramificado en $L/K$, se tiene que $\pK_1|_{
\ma Q}=p$ es ramificado en $M/{\ma Q}$. Por otro lado, $\pL_2|_{
\ma Q}$ es no ramificado en $L/K$ y tampoco en $K/{\ma Q}$.
Por tanto $\pL_2|_{\ma Q}=p$ es no ramificado en $M/{\ma Q}$.
Esta contradicci\'on prueba que $L\not\subseteq K(\zeta_{\infty})$.
$\fin$
\end{proof}

%% file: Capitulo5.tex
\chapter{Propiedades y aplicaciones de los
campos ciclot\'omicos}\label{ch4}

\section{Caso especial del teorema de 
Dirichlet}\label{S4.1}

Una de las aplicaciones interesantes de los campos
ciclot\'omicos es la demostraci\'on de un caso particular
del Teorema de Dirichlet\index{teorema de
Dirichlet}\index{Dirichlet!teorema sobre $\sim$} de
n\'umeros primos en progresiones aritm\'eticas.
El Teorema de Dirichlet establece que si $a,b\in{\ma Z}$
son primos relativos, entonces existen una infinidad
de n\'umeros primos en la progresi\'on aritm\'etica
$a,a+2b,\ldots, a+mb,\ldots$. En otras palabras,
la congruencia $p\equiv a\bmod b$ tiene una infinidad
de soluciones para $p$ un n\'umero primo.

En la actualidad el Teorema de Dirichlet es una aplicaci\'on
del Teorema de Densidad de \v{C}ebotarev\index{teorema de
densidad de \v{C}ebotarev}\index{Cebotarev@\v{C}ebotarev!teorema
de densidad de $\sim$} aplicado a los campos ciclot\'omicos.
La demostraci\'on original de Dirichlet se basa en la
no anulaci\'on de una serie $L$ al ser evaluada en $s=1$.
Esto lo veremos m\'as adelante.

El caso particular de $a=1$, $b=n$ arbitrario ha sido
objeto de estudio y existen muchas demostraciones de este
caso, incluyendo varias elementales, algunas de ellas de
hecho m\'as elementales que la que presentamos aqu{\'\i}.
Otro caso especial es $a=-1$, $b=n$. Existe una demostraci\'on
elemental de este caso debido a E. Landau \cite{Lan09}.
Desconocemos si los polinomios ciclot\'omicos son aplicables
a este caso o alg\'un otro, diferente al que desarrollamos a
continuaci\'on.

\begin{proposicion}\label{P7.1}
Sean $p$ un n\'umero primo y $n\in{\ma N}$ tal que $p\nmid n$.
Sea $a\in{\ma Z}$. Entonces
\[
p|\psi_n(a)\iff o(a\bmod p)=n,
\]
es decir, $a^n\equiv 1\bmod p$ y $a^k\not\equiv 1\bmod p$
para $1\leq k\leq n-1$.
\end{proposicion}

\begin{proof}
{\ }

\noindent
$\Longrightarrow)$ Primero supongamos que $p|\psi_n(a)$.
Puesto que $x^n-1=\prod_{d|n}\psi_d(x)$, se sigue que
$a^n-1=\prod_{d|n}\psi_d(a)\equiv 0\bmod p$ y por lo
tanto $a^n\equiv 1 \bmod p$. Sea $k:=o(a\bmod p)$. Entonces
$k|n$. Ahora bien, si $k<n$ tendr{\'\i}amos que $0\equiv
a^k-1=\prod_{d|k}\psi_d(a)\bmod p$. En particular existe
$d_0|k$ (y por tanto $d_0<n$) tal que $\psi_{d_0}(a)\equiv
0\bmod p$. De esto se sigue que el polinomio
$x^n-1\bmod p$ tiene una ra{\'\i}z m\'ultiple pues
$\psi_{d_0}(a)\equiv \psi_n(a)\equiv 0\bmod p$ y $d_0\neq n$.
Sin embargo esto no es posible pues el polinomio
$p(x)=(x^n-1)\bmod p\in {\ma F}_p[x]$ no tiene ra{\'\i}ces
m\'ultiples debido a que $p'(x)=nx^{n-1}$ tiene una 
\'unica ra{\'\i}z, a saber $\alpha=0\bmod p$ (pues $p\nmid n$)
y $0\bmod p$ no es ra{\'\i}z de $p(x)$. Esta contradicci\'on
prueba que $o(a\bmod p)=n$.

\noindent
$\Longleftarrow)$ Rec{\'\i}procamente, si $o(a\bmod p)=n$,
entonces $a^n\equiv 1\bmod p$ pero $a^k\not\equiv1
\bmod p$ para $1\leq k\leq n-1$. De esto se sigue que
\begin{gather*}
a^n-1=\prod_{d|n}\psi_d(a)\equiv 0\bmod p \quad \text{y}\\
\quad a^k-1=\prod_{d|k}\psi_d(a)\not\equiv 0\bmod p
\quad \text{para toda} \quad k<n.
\end{gather*}
Por lo tanto $\psi_n(a)\equiv 0\bmod p$. $\fin$

\end{proof}

Usamos la Proposici\'on \ref{P7.1} para probar:

\begin{proposicion}\label{P7.2}
Si $p\nmid n$ donde $p$ es un n\'umero primo y $n\in
{\ma N}$, entonces $p|\psi_n(a)$ para alg\'un $a\in
{\ma Z}\iff p\equiv 1\bmod n$.
\end{proposicion}

\begin{proof}
{\ }

\noindent
$\Longrightarrow)$ Si $p|\psi_n(a)$ para alg\'un
$a\in{\ma Z}$, entonces $o(a\bmod p)=n$. Notemos
que necesariamente $p\nmid a$. Por lo tanto
$n|o(U_p)=p-1$ y por lo tanto $p\equiv 1\bmod n$.

\noindent
$\Longleftarrow)$ Rec{\'\i}procamente, como $U_p$ es un
grupo c{\'\i}clico de orden $p-1$, si $p\equiv 1\bmod n$, 
entonces $n| p-1$ y existe un elemento $\overline{a}
\in U_p$ de orden $n$, es decir $o(a\bmod p)=n$,
$a\in{\ma Z}$. Por la Proposici\'on \ref{P7.1} se sigue que
$p|\psi_n(a)$. $\fin$

\end{proof}

El caso especial del Teorema de Dirichlet se sigue
de este resultado.

\begin{corolario}[Caso especial del Teorema de
Dirichlet]\label{C7.3}
Sea $n\in{\ma N}$. Entonces existe una infinidad de n\'umeros
primos $p$ tales que $p\equiv 1\bmod n$.
\end{corolario}

\begin{proof}
Supongamos que para un $n$ cualquiera hemos probado
la existencia de un primo $p\equiv 1\bmod n$. En este caso
$n|p-1$ y por lo tanto $n\leq p-1$. Sea $\alpha\in {\ma N}$ tal que
$p<n^{\alpha}$. Aplicando la hip\'otesis para $n^{\alpha}$,
tenemos que existe un n\'umero primo
$q$ tal que $q\equiv 1\bmod n^{\alpha}$ y en particular tenemos
$q\equiv 1\bmod n$. Entonces $n^{\alpha}\leq q-1$ y por lo
tanto $p<q$. De esta forma se sigue que existe una  infinidad
de n\'umeros primos tales que $p\equiv 1\bmod n$. Es decir,
basta hallar uno solo de ellos.

Se tiene que para $m\in{\ma N}$, $nm\equiv 0\bmod n$ y
$\psi_n(nm)\equiv \psi_n(0)\bmod n$. Puesto que $\psi_n(x)=
\prod_{d|n} (x^d-1)^{\mu(n/d)}$, se sigue que $\psi_n(0)=
\pm 1$. Entonces $\psi_n(nm)\equiv \pm 1\bmod n$. En 
particular si $q$ es primo tal que $q|n$ entonces $q\nmid
\psi_n(nm)$.

Ahora bien, puesto que $\psi_n(x)$ es un
un polinomio de grado $\varphi(n)\geq 1$ y con coeficiente
l{\'\i}der $1$, en particular positivo, entonces $\lim_{m\to
\infty}\psi_n(nm)=\infty$.

Sea $m\in{\ma N}$ tal que $\psi_n(nm)>1$ y sea $p$ es un
n\'umero primo tal que $p|\psi_n(nm)$. Entonces $p\nmid n$ y por
la Proposici\'on \ref{P7.2}, se tiene que $p\equiv 1\bmod n$
y esto termina la demostraci\'on. $\fin$

\end{proof}

\begin{observacion}\label{O7.4} La demostraci\'on del
Corolario \ref{C7.3} puede hacerse suponiendo que
$p_1,\ldots, p_r$ son todos los n\'umeros primos congruentes
con $1$ m\'odulo $n$ y considerar $mnp_1\cdots p_r$ y como
hicimos, hallamos un primo $p_{r+1}|\psi_n(mn p_1\cdots p_r)$
para alguna $m\geq 1$. Este $p_{r+1}$ es diferente a
$p_1,\ldots, p_r$ y congruente a $1\bmod n$. Esta es la
demostraci\'on de Euclides con $n=2$. Es decir, se puede
considerar el polinomio $\psi_2(x)=x+1$ para probar que hay
una infinidad de primos impares.

Por ejemplo si usamos $\psi_4(x)=x^2+1$, probaremos que 
existen una infinidad de primos de la forma $4n+1$. Como 
ejercicio dejamos al lector probar que existen una infinidad
de n\'umeros primos de la forma $8n+1$ usando el polinomio
ciclot\'omico $\psi_8(x)=x^4+1$.
\end{observacion}

\begin{observacion}\label{O7.4'} La Proposici\'on \ref{P7.2}
nos proporciona una manera c\'omoda de producir primos
congruentes con $1$ m\'odulo $n$, pues simplemente
evaluamos $\psi_n(a)$ con $a\in{\ma Z}$ y todos los
factores primos de $\psi_n(a)$ son de esta forma.

Por ejemplo, $\psi_7(4)=\frac{4^7-1}{4-1}=\frac{2^{14}-1}{3}=
\frac{16(1024)-1}{3}= \frac{16383}{3}=5461= 43\cdot 127$ y
$43$ y $127$ son ambos primos congruentes con $1$ m\'odulo
$7$.
\end{observacion}

\section{Descomposici\'on de primos en $\cic n{}/{\ma Q}$}\label{S4.2}

Para analizar la descomposici\'on de un n\'umero primo $p$
en $\cic n{}$, $p\nmid n$, veremos cual es el s{\'\i}mbolo
de Artin asociado a $p$. Para ello empezamos con

\begin{proposicion}\label{P8.1}
Sea $p$ un n\'umero primo en ${\ma Z}$ y sea ${\eu p}$ un
ideal primo en $\cic n{}$ que divide a $p$. Entonces las
$n$--ra{\'\i}ces de la unidad son todas distintas m\'odulo ${\eu p}$.
\end{proposicion}

\begin{proof}
Puesto que $1+x+\cdots+x^{n-1}=\frac{x^n-1}{x-1}=\prod_{j=1}^{
n-1}(x-\zeta_n^j)$ se tiene que evaluando en $x=1$ obtenemos
$n=\prod_{j=1}^{n-1}(1-\zeta_n^j)$.

Como $p\nmid n$ se sigue que $n\notin {\eu p}$ y por lo
tanto $1-\zeta_n^j\notin {\eu p}$, $j=1,\dots, n-1$.

Si $\zeta_n^i\equiv \zeta_n^j\bmod {\eu p}$ entonces $\zeta_n^i(
1-\zeta_n^{j-i})\in{\eu p}$. Puesto que $\zeta_n^i$ es una
ra{\'\i}z de unidad y en particular una unidad, se sigue que
$1-\zeta_n^{j-i} \in{\eu p}$ lo cual implica que $j-i=0$, es decir,
$i=j$. Por tanto, para $1\leq i, j\leq n-1$, $i\neq j$, entonces
$\zeta_n^i\not\equiv \zeta_n^j\bmod {\eu p}$. $\fin$
\end{proof}

En la situaci\'on de la Proposici\'on \ref{P8.1} se tiene el diagrama
\[
\xymatrix{
{\eu p}\ar@{-}[r]&{\ma Z}[\zeta_n]\ar@{-}[r]\ar@{-}[d]&\cic n{}
\ar@{-}[d]\\
p\ar@{-}[r]&{\ma Z}\ar@{-}[r]&{\ma Q}
}
\]

Como $p\nmid n$, entonces ${\eu p}$ no es ramificado en
$\cic n{}/{\ma Q}$. Sean ${\cal O}_{\cic n{}}/{\eu p}= {\ma Z}
[\zeta_n]/{\eu p}$ y ${\cal O}_{{\ma Q}}/\langle p \rangle=
{\ma F}_p$ los campos residuales y sea $f:=[{\ma Z}[\zeta_n]/
{\eu p}:{\ma F}_p]$. Entonces ${\ma Z}[\zeta_n]/{\eu p}\cong
{\ma F}_{p^f}$ el campo finito de $p^f$ elementos. El
automorfismo de Frobenius $\sigma_p\in \Gal(\cic n{}/{\ma Q})$
correspondiente a $p$ es (ver Secci\'on \ref{Sec0.3}) $\sigma_p
=\Big[\frac{\cic n{}/{\ma Q}}{{\eu p}}\Big]$ el cual satisface
\[
\sigma_p y \equiv y^p \bmod {\eu p}, \quad y\in{\ma Z}[\zeta_n].
\]
En particular, $\sigma_p\zeta_n=\zeta_n^p\bmod {\eu p}$.

Por otro lado, para $\sigma\in \Gal(\cic n{}/{\ma Q})$, se tiene
$\sigma\zeta_n=\zeta_n^a$ con $\mcd (a,n)=1$ y por la
Proposici\'on \ref{P8.1}, $\zeta_n^a\equiv \zeta_n^p\bmod
{\eu p}\iff a=p$. Se sigue que $\sigma_p$ est\'a dado por
$\sigma_p\zeta_n=\zeta_n^p$.

As{\'\i}, $\overline{\sigma_p}=\sigma_p\bmod {\eu p}$ es 
generador del grupo c{\'\i}clico $\Gal\big(({\ma Z}[\zeta_n]/{\eu p})/
({\ma Z}/p{\ma Z})\big)$ el cual es de orden $f$. Por tanto
\[
\sigma_p^k\zeta_n=\zeta_n^{p^k}=\zeta_n\iff \sigma_p^k=\Id
\iff f|k.
\]
Es decir $f$ es el m{\'\i}nimo tal que $\zeta_n^{p^f}=\zeta_n$,
esto \'ultimo equivale a $p^f\equiv 1\bmod n$. 
Por lo tanto $f=o(p\bmod n)$.

Puesto que $p$ no es ramificado en $\cic n{}/{\ma Q}$ y
$[\cic n{}:{\ma Q}]=\varphi(n)$, si $p{\ma Z}[\zeta_n]=
{\eu p}_1\cdots {\eu p}_g$, $gf=\varphi(n)$ con $f=o(p\bmod
n)$. Es decir, hemos probado el llamado {\em Teorema de
Reciprocidad Ciclot\'omica\index{teorema de reciprocidad
ciclot\'omica}}.

\begin{teorema}[Reciprocidad ciclot\'omica]\label{T8.2}
Sea $p$ un primo racional, $n\in{\ma N}$. Si $p\nmid n$, 
entonces si ${\eu p}$ es un primo de ${\ma Z}[\zeta_n]$
sobre $p$ se tiene
\[
f=[{\ma Z}[\zeta_n]/{\eu p}:{\ma F}_p]=o(p\bmod n)
\]
y $g=\varphi(n)/f$ donde $p{\ma Z}[\zeta_n]={\eu p}_1\cdots
{\eu p}_g$ es la factorizaci\'on en
primos distintos y $\gr {\eu p}_i=f=o(p\bmod n)$.

En particular $p$ se descompone totalmente en $\cic n{}
\iff f=1\iff p\equiv 1\bmod n$. $\fin$
\end{teorema}

\begin{observacion}\label{O8.3} Como mencionamos 
anteriormente resulta ser que
 el Teorema de Kronecker--Weber\index{teorema
de Kronecker--Weber}\index{Kronecker--Weber!teorema de $\sim$}
es una consecuencia de la Teor{\'\i}a de Campos de Clases. Ah{\'\i}
se estudian los llamados ``m\'odulus''\index{m\'odulus} los cuales
son productos formales de la forma ${\eu m}=\prod {\eu p}^{n_{\eu
 p}}$ donde ${\eu p}$ recorre el conjunto de primos de un campo
 $K$ el cual, de momento, podemos suponer num\'erico,
 $n_{\eu p}\in{\ma Z}$, $n_{\eu p}\geq 0$ para todo ${\eu p}$
 y $n_{\eu p}=0$ para casi todo ${\eu p}$. Si ${\eu p}$ es un 
primo infinito real, entonces $n_{\eu p}=0$ o $1$ y si ${\eu p}$
es un primo infinito imaginario, $n_{\eu p}=0$. Entonces
se construye el grupo de clases de ideales de ``rayos''
correspondientes al m\'odulus ${\eu m}$ el cual se define por
$I_{K,{\eu m}}=D_{\eu m}/i(K_{{\eu m},1})$ donde $D_{\eu m}$
son los ideales fraccionarios de ${\cal O}_K$ primos 
relativos a todos los primos finitos que aparecen en ${\eu m}$,
$K_{{\eu m},1}$ son los elementos $x$ de $K^{\ast}=K\setminus
\{0\}$ tales que $x\equiv 1\bmod {\eu p}^{n_{\eu p}}$
para cada primo finito que aparece en ${\eu m}$ y tal que $\sigma
x>0$ donde $\sigma$ representa el encaje correspondiente
al primo infinito real ${\eu p}$ que aparece en ${\eu m}$, es decir,
$n_{\eu p}=1$ y finalmente, si $x\in K_{{\eu m},1}$, $i(x)$ denota
al ideal fraccionario principal generado por $x$.

Ahora si $H$ es un ``{\em subgrupo de
congruencia\index{subgrupo de congruencia}}'', es decir, si
$i(K_{{\eu m},1})\subseteq H\subseteq D_{\eu m}$, entonces
el Teorema de Existencia nos dice que existe un campo de
clase $L/K$ asociado a $H$, es decir, existe una \'unica
extensi\'on abeliana $L$ de $K$ tal que los primos de
${\cal O}_K$ que se descomponen totalmente en $L$ son,
salvo un conjunto de densidad $0$, los elementos de $H$.
Este es el {\em campo de clase\index{campo de clase}}
asociado a $H$.

En particular, la m\'axima extensi\'on abeliana de $K$ es
la uni\'on de todos los campos de clase asociados a $K$.

El Teorema \ref{T8.2} junto con esta observaci\'on nos dice
que $\cic m{}$ es el campo de clase asociado al m\'odulus
$m\infty$ y por tanto la m\'axima extensi\'on abeliana de
${\ma Q}$ es precisamente $\bigcup_{m=1}^{\infty} \cic m{}$
de donde obtenemos nuevamente el Teorema de 
Kronecker--Weber (ver Corolario \ref{CClaseT.4.5.8} y
Teorema \ref{CClaseT4.7.4}).
\end{observacion}

\begin{ejemplo}\label{Ej8.4} Consideremos $n=4$,
$\cic 4{}={\ma Q}(i)={\ma Q}(\sqrt{-1})$. El \'unico primo
finito ramificado en $\cic 4{}/{\ma Q}$ es $p=2$. Si $p$ es
impar, consideremos su clase m\'odulo $4$, $p\equiv 1\bmod 4$
o $p\equiv 3\bmod 4\equiv -1\bmod 4$. Entonces si $p\equiv 1
\bmod 4$, $f=1$ y por lo tanto $g=\frac{\varphi(4)}{1}=2$, es
decir, $p{\ma Z}[i]=\pK\overline{\pK}$. Si $p\equiv -1\bmod 4$,
$p^2\equiv 1\bmod 4$ y por lo tanto $f=2$, $g=1$ y en este
caso $p{\ma Z}[i]=\pK=\langle p \rangle$ permanece primo.

Ahora ${\ma Z}[i]$ es un anillo euclidiano con la norma:
\[
N(a+bi)=a^2+b^2=( a+bi)(a-bi), \quad a,b\in{\ma Z}.
\]
En particular ${\ma Z}[i]$ es de ideales principales. Sea
$p\equiv 1\bmod 4$, $p{\ma Z}[i]=\pK\overline{\pK}=\langle
\alpha \rangle \langle \beta \rangle$. Se tiene $N\alpha = 
N\beta = p$. 

Si $\alpha=a+bi$, $p=a^2+b^2$, es decir $p$ es suma de
dos cuadrados; por ejemplo: $13=2^2+3^2$; $5=1^2+2^2$;
$17=1^2+4^2$; $29=2^2+5^2$; etc.

Si $p\equiv 3\bmod 4$, la congruencia $a^2+b^2\equiv 3\bmod 4$
no tiene soluci\'on, pues para cualesquiera $a,b\in{\ma Z}$, se tiene
$a^2+b^2\equiv 0,1, 2\bmod 4$ por lo que $p\equiv 3\bmod 4$ no
es suma de cuadrados.

Finalmente, para $p=2$, $p$ es ramificado por lo que $2{\ma Z}[i]=
\pK^2$. Si $\gamma =a +bi$ genera $\pK$, se tiene $N\gamma
=a^2+b^2=2$ lo cual implica que $\gamma=\pm 1\pm i$. Adem\'as
$2=1^2+1^2$.

De lo anterior, dejamos deducir al lector que $m\in{\ma N}$ es suma
de dos cuadrados: $m=a^2+b^2$, $a,b\in{\ma Z}$ si y solamente
si $m=2^{\alpha}p_1\cdots p_r t^2$ donde $\alpha\in{\ma N}\cup
\{0\}$, $t\in{\ma N}$, $r\geq 0$, $p_1,\ldots, p_r$ son primos
congruentes con $1$ m\'odulo $4$.
\end{ejemplo}

\begin{ejemplo}\label{Ej8.5}
Sea $n=23$. Se tiene que $U_{23}$ es un grupo c{\'\i}clico
de orden $22$. Empezando por $2$, directamente se calcula
que $\langle 2\bmod 23\rangle=\{1,2,4,8,16,9,18,13,3,6,12\}$
y $o(2\bmod 23)=11$. Entonces puesto que $5^2\equiv 2\bmod
23$, se sigue que $5\bmod 23$ genera $U_{23}$, es decir $5$
es ra{\'\i}z $23$--\'esima primitiva. De lo anterior y puesto que
$\varphi(22)=10$, se sigue que los generadores de $U_{23}$ ser\'an:
\[
\{5,5^3,5^5,5^7,5^9,5^{13},5^{15},5^{17},5^{19},5^{21}\}=
\{5,7,10,11,14,15,17,19,20,21\}.
\]

As{\'\i}, $x\in U_{23}$ satisface
\begin{gather*}
o(x)=1\iff x\equiv 1\bmod 23,\\
o(x)=2\iff x\equiv -1\bmod 23\equiv 22\bmod 23,\\
o(x)=11 \iff x\in\{2,3,4,6,8,9,12,13,16,18\}\bmod 23,\\
o(x)=22 \iff x\in\{5,7,10,11,14,15,17,19,20,21\}\bmod 23.
\end{gather*}
As{\'\i} $p$, un n\'umero primo en ${\ma Z}$, se descompone 
totalmente en $\cic {23}{}\iff p\equiv 1\bmod 23$; $p$ se 
descompone en $11$ factores $\iff p\equiv 22\bmod 23$; $p$ se
descompone en $2$ factores $\iff p\equiv x\bmod 23$ con
$x\in \{2,3,4,6,8,9,12,13,16,18\}$ y finalmente $p$ permanece
primo en ${\ma Z}[\zeta_{23}]\iff p\equiv y\bmod 23$ donde
$y\in\{5,7,10,11,14,15,17,19,20,21\}$.
\end{ejemplo}

Ahora que conocemos la m\'axima extensi\'on abeliana
de ${\ma Q}$, establecemos cual es su grupo de Galois.
Primero recordemos

\begin{teorema}\label{T8.6} Sean $p$ un primo impar.
Sea ${\ma Q}^{(p)}:=\bigcup_{n=0}^{\infty} \cic pn$. Entonces
\[
G^{(p)}:=\Gal({\ma Q}^{(p)}/{\ma Q}) \cong {\ma Z}/(p-1){\ma Z}
\times {\ma Z}_p\cong C_{p-1}\times {\ma Z}_p
\]
donde ${\ma Z}_p$ es el anillo de los enteros $p$--\'adicos.
\end{teorema}

\begin{proof}
Se tiene que $\Gal(\cic pn/{\ma Q})\cong U_{p^n}\cong
{\ma Z}/(p-1){\ma Z}\times {\ma Z}/p^{n-1} {\ma Z}$ y se tiene el
diagrama conmutativo
\[
\xymatrix{\Gal(\cic p{n+1}/{\ma Q})\ar[r]^{\rest_{n+1}}
\ar[d]^{\theta_{n+1}}&\Gal(\cic pn/{\ma Q})\ar[d]^{\theta_n}\\
U_{p^{n+1}}\ar[r]_{\pi_{n+1}}&U_{p^n}
}
\]
donde $\rest_{n+1}$ es el mapeo restricci\'on: $\rest_{n+1}(
\sigma)=\sigma|_{\cic pn}$ y $\pi_{n+1}$ es la proyecci\'on
natural $\pi_{n+1}(x\bmod p^{n+1})=x\bmod p^n$. Notemos
que $\pi_{n+1}|_{{\ma Z}/(p-1){\ma Z}}= \Id_{{\ma Z}/(p-1){\ma Z}}$.

Entonces, usando el Teorema \ref{T6.2}, se tiene
\begin{align*}
G^{(p)}:&=\Gal({\ma Q}^{(p)}/{\ma Q}) = \Gal\Big(\bigcup_{n=1}^{
\infty}\cic pn/{\ma Q}\Big) =
\Gal(\lim_{\substack{\longto\\ n}}
\cic pn/{\ma Q})\\
&\cong \lim_{\substack{\longleftarrow\\ n, \rest_{n+1}}} \Gal(\cic pn/
{\ma Q})\cong \lim_{\substack{\longleftarrow
\\ n, \pi_{n+1}}} U_{p^n}\\
&\cong \lim_{\substack{\longleftarrow
\\ n,\pi_{n+1}}} \big({\ma Z}/(p-1){\ma Z}
\times {\ma Z}/p^{n-1}{\ma Z}\big)\\
&\cong \big({\ma Z}/(p-1){\ma Z}\big)
\times \lim_{\substack{\longleftarrow
\\ \pi_{n+1}}}
\big({\ma Z}/p^{n-1}{\ma Z}\big)
 \cong \big({\ma Z}/(p-1){\ma Z}\big)\times {\ma Z}_p
\cong C_{p-1}\times {\ma Z}_p. \tag*{$\fin$}
\end{align*}
\end{proof}

El resultado an\'alogo para $p=2$ es:

\begin{teorema}\label{T8.7}
Sea ${\ma Q}^{(2)}:=\bigcup_{n=1}^{\infty}\cic 2n$. Entonces
\[
G^{(2)}:=\Gal({\ma Q}^{(2)}/{\ma Q})\cong C_2\times {\ma Z}_2.
\]
\end{teorema}

\begin{proof}
Es igual a la demostraci\'on del Teorema \ref{T8.6} considerando
que $U_{2^n}\cong C_2\times C_{2{n-2}}$ y $\Gal(\cic 2n/{\ma Q})
\cong U_{2^n}$ donde el primer factor $C_2$ de $U_{2^n}$
est\'a generado por la conjugaci\'on compleja. $\fin$
\end{proof}

Los Teoremas \ref{T8.6} y \ref{T8.7} nos dan

\begin{teorema}\label{T8.8} Sea ${\ma Q}^{\ab}$ la m\'axima
extensi\'on abeliana de ${\ma Q}$. Entonces
\begin{align*}
G_{\ab}:&=\Gal({\ma Q}^{\ab}/{\ma Q}) \cong \big(
{\ma Z}/2{\ma Z}\big)
\times\prod_{\substack{p\text{\ primo}\\ p>2}}\big( {\ma Z}/
(p-1){\ma Z}\big)\times \prod_{p\text{\ primo}}{\ma Z}_p\\
&\cong C_2\times \prod_{p\text{\ impar}} C_{p-1}\times
\hat{{\ma Z}}.
\end{align*}
\end{teorema}

\begin{proof}
Por el Teorema de Kronecker--Weber se tiene que ${\ma Q}^{\ab}=
\bigcup_{n=1}^{\infty}\cic n{}$. Para cada $n\in{\ma N}$, sea
$n=p_1^{\alpha_1}\cdots p_r^{\alpha_r}$ su descomposici\'on
en primos. Entonces $U_n\stackrel{\delta}{\cong} U_{p_1^{\alpha_1}}
\times\cdots \times U_{p_r^{\alpha_r}}$ donde
\begin{eqnarray*}
\delta\colon U_n &\longto & U_{p_1^{\alpha_1}}
\times\cdots \times U_{p_r^{\alpha_r}}\\
x\bmod n&\longmapsto & (x\bmod p_1^{\alpha_1},\cdots,
x\bmod p_r^{\alpha_r}).
\end{eqnarray*}

A nivel de grupos de Galois, $\cic n{}=\prod_{i=1}^r\cic {p_i}{\alpha_i}$
y $\Gal(\cic n{}/{\ma Q})\stackrel{\varphi}{\cong} \prod_{i=1}^r\Gal( \cic
{p_i}{\alpha_i}/{\ma Q})$ donde $\varphi$ est\'a dado por la
restricci\'on: $\sigma\in \Gal(\cic n{}/{\ma Q})$, $\varphi(\sigma)=
\big(\sigma|_{\cic {p_1}{\alpha_i}},\cdots, \sigma|_{\cic {p_r}{\alpha_r}}
\big)$ y se tiene el diagrama conmutativo
\[
\xymatrix{
\Gal(\cic n{}/{\ma Q})\ar[rr]^{\varphi}\ar[d]^{\theta_n}&&
\prod\limits_{i=1}^r \Gal(\cic {p_i}{\alpha_i}/{\ma Q})\ar[d]^{(
\theta_{p_1^{\alpha_1}},\cdots, \theta_{p_r^{\alpha_r}})}\\
U_n\ar[rr]^{\delta}&&\prod\limits_{i=1}^{r} U_{p_i^{\alpha_i}}
}
\]

Por tanto se tiene los isomorfismos
\begin{gather*}
\lim_{\substack{\longleftarrow\\ n}}\Gal(\cic n{}/{\ma Q})\cong
\lim_{\substack{\longleftarrow\\ n}} U_n,\\
\lim_{\substack{\longleftarrow\\ n}} U_n \cong
\prod_{p\text{\ primo}} \Big(\lim_{\substack{\longleftarrow\\ \alpha(p)}}
 U_{p^{\alpha(p)}}\Big) \quad \text{y}\\
 \lim_{\substack{\longleftarrow\\ n}} \Gal(\cic n{}/{\ma Q})\cong
 \prod_{p\text{\ primo}} \Big(\lim_{\substack{\longleftarrow\\ \alpha(p)}}
\Gal(\cic p{\alpha(p)}/{\ma Q})\Big).
\end{gather*}
Estos tres isomorfismos y los Teoremas \ref{T6.2},
\ref{T8.6} y \ref{T8.7}
nos dan
\begin{align*}
G_{\ab}&=\Gal({\ma Q}^{\ab}/{\ma Q})= \Gal\Big(\bigcup_{n=1}^{\infty}
\cic n{}/{\ma Q}\Big)\\
&\cong \Gal\Big(\big(\lim_{\substack{\longto\\ n}}\cic n{}
\big)/{\ma Q}\Big)\cong
\lim_{\substack{\longleftarrow\\ n}} \Gal(\cic n{}/{\ma Q})\\
&\cong \prod_{p{\text{\ primo}}}\Big[\lim_{\substack{\longleftarrow\\ 
\alpha(p)}} \Gal(\cic p{\alpha(p)}/{\ma Q})\Big]\cong
\prod_{p{\text{\ primo}}}\big[
\lim_{\substack{\longleftarrow\\ \alpha(p)}}U_{p^{\alpha(p)}}\big]\\
&\cong \lim_{\substack{\longleftarrow\\ n}} U_n\cong
\lim_{\substack{\longleftarrow\\ m}} U_{2^m}\times 
\prod_{\substack{p{\text{\ primo}}\\ p>2}}
\lim_{\substack{\longleftarrow\\ m}} U_{p^m}\\
&\cong C_2\times {\ma Z}_2\times \prod_{p>2}(C_{p-1}\times
{\ma Z}_p)\cong C_2\times \prod_{p>2}C_{p-1}\times \prod_{p\text{\ 
primo}}
{\ma Z}_p.
\end{align*}

Finalmente, se tiene $\prod_{p{\text{\ primo}}} {\ma Z}_p\cong
\hat{{\ma Z}}$. $\fin$

\end{proof}

\begin{observacion}[Anillos euclidianos y dominios
de ideales principales]\label{O8.9-1}

Dada una extensi\'on abeliana finita $K$ de ${\ma Q}$, este campo
es subcampo de un campo de un campo ciclot\'omico
$K\subseteq \cic n{}$. 
Se sabe que \'unicamente hay $30$ campos ciclot\'omicos
con n\'umero de clase $1$ \cite[Chapter 11]{Was97}. 
Estos son $\cic n{}$ con
\begin{gather*}
n\in\{1,3,4,5,7,8,9,11,12,13,15,16,17,19,20,21,24,
25,27,28,32,33,\\
35,36,40,44,45,48,60,84\}
\end{gather*}

Dado que conocemos cuales campos $\cic n{}$ tiene n\'umero
de clase $1$, ?`eso nos dice que campos $K$ tienen n\'umero de 
clase $h_K$ igual a $1$? La respuesta es no. Es una conjetura
de K. F. Gauss que existe una infinidad de campos cuadr\'aticos
reales con n\'umero de clase $h_K=1$. En la actualidad, este
sigue siendo un problema abierto. De hecho no se sabe si
existen una infinidad de campos num\'ericos (de cualquier
grado) con n\'umero de clase $1$. Por ejemplo $h_K=1$ con
${\ma Q}(\sqrt{d})$, $d>0$ libre de cuadrados y $d<100$ para
\begin{gather*}
d\in\{2,3,5,6,7,11,13,14,17,19,21,22,23,29,31,33,37,38,41,43,
46,47,\\
53,57.59,61,62,67,69,71,73,77,83,86,89,93,94,97\}
\end{gather*}
(ver \cite[Chapter 1, Sect. 6]{Neu99}).

Para campos cuadr\'aticos $K={\ma Q}(\sqrt{d})$ con $d<0$,
$d$ libre de cuadrados, se sabe que hay exactamente nueve
campos con $h_K=1$. Estos son $K_d={\ma Q}(\sqrt{d})$
con 
\begin{gather*}
d\in\{-1,-2,-3,-7,-11,-19,-43,-67,-169\}.
\end{gather*}
(ver \cite{Sta69}). K. Yamamura \cite{Yam94} prob\'o que
hay exactamente 172 campos abelianos imaginarios $K$
con n\'umero de clase $h_K=1$.

Ahora bien, si el anillo de enteros de un campo es un anillo 
euclidiano, necesariamente el campo tiene n\'umero de clase
$1$. Cuando se tiene un anillo de ideales principales, es
decir, con n\'umero de clase $1$, no necesariamente el
anillo es euclidiano. Es un problema dif\'icil dar contraejemplos
pues, aunque en muchos casos se puede probar que la norma
no es una funci\'on euclidiana, la cual es la candidata natural,
se debe probar que no existe ninguna funci\'on euclidiana.
T. Motzkin \cite{Mot49} da un ejemplo de un anillo de ideales
principales que no es euclidiano. Este anillo es el anillo
de enteros del campo ${\ma Q}(\sqrt{-19})$. Por otro lado,
D. Clark \cite{Cla94} prob\'o que el anillo de enteros de
${\ma Q}(\sqrt{69})$ es euclidiano pero no euclidiano con
la norma.

Entonces una pregunta natural es: ?`cuales de los anillos
de enteros ciclot\'omicos que son de ideales principales son
euclidianos? No se conoce la respuesta completa y la
literatura sobre este problema es extensa. Particularmente,
las contribuciones de H. W. Lenstra Jr. y en especial
sus art\'iculos \cite{Len79} y el art\'iculo panor\'amico de
F. Lemmermeyer \cite{Lem95} publicado en 1995 y con
diversas actualizaciones, dan un panorama general sobre
el estado del arte de este problema.

Los anillos de enteros $\AE{\cic n{}}={\ma Z}[\zeta_n]$ de
$\cic n{}$ que se sabe que son euclidianos (con la norma)
son
\begin{gather*}
n\in\{1,3,4,5,7,8,9,11,12,13,15,16,20,24\}.
\end{gather*}
Para $n\notin\{13,16,24\}$ una prueba puede hallarse en
el art\'iculo de Lenstra \cite{Len75}. En adici\'on, Lenstra
prob\'o que ${\ma Z}[\zeta_{32}]$ no es euclidiano con la norma
\cite{Len79}. Ya era conocido que los anillos
${\ma Z}[\zeta_n]$ con $n\in\{1,3,4,5,8,12\}$
eran euclidianos. Kummer conjetur\'o que $n=17$ y $19$
son tambi\'en euclidianos. La respuesta no se conoce.
El caso de que ${\ma Z}[\zeta_{13}]$ es euclidiano
fue probado por R. G. McKenzie en su tesis
doctoral \cite{McK88}.

\end{observacion}

\begin{ejemplo}\label{Ej8.9}

En este ejemplo estudiaremos
$\cic {23}{}$ el cual es el primer campo ciclot\'omico con 
n\'umero de clase mayor a $1$.

Primero consideremos ${\ma Q}(\sqrt{-23})$. Ahora bien,
$\cic {23}{}$ tiene un \'unico subcampo cuadr\'atico $K$:
$K=\cic {23}{}^H$ donde $H$ es el \'unico subgrupo de $U_{23}$
de orden $11$, de hecho $H\cong  \langle 2\bmod 23\rangle
\cong \langle \sigma\in\Gal(\cic {23}{}/{\ma Q})\mid \sigma
\zeta_{23}=\zeta_{23}^2\rangle$.

Ahora si ${\ma Q}(\sqrt{d})$ es este subcampo, entonces $23$
es el \'unico primo ramificado en ${\ma Q}(\sqrt{d})$ por lo que
$\delta_{{\ma Q}(\sqrt{d})}=\pm 23=d$. Puesto que $23\equiv
3\bmod 4$ y $-23\equiv 1 \bmod 4$ se sigue que $d=-23$ y
${\ma Q}(\sqrt{d})={\ma Q}(\sqrt{-23})$. En otras palabras
${\ma Q}(\sqrt{-23})\subseteq \cic {23}{}$.

La base entera de ${\ma Q}(\sqrt{-23})/{\ma Q}$ es $\big\{
1,\frac{1+\sqrt{-23}}{2}\big\}$. Si ponemos $\alpha:=\frac{1+
\sqrt{-23}}{2}$, entonces $\Irr(\alpha,x,{\ma Q})=x^2-x+6=f(x)$.
Se tiene $f(x)\bmod 2=x^2-x=x(x-1)$. Por el Teorema de
Kummer se tiene que $2{\cal O}_{{\ma Q}(\sqrt{-23})} =
\pK\overline{\pK}$ donde
\[
\pK=\Big\langle 2,\frac{-1+\sqrt{-23}}{2}\Big\rangle
=\langle 2,\overline{\alpha}\rangle=\langle 2,1-\alpha\rangle,
\qquad \overline{\pK}=\Big\langle 2,\frac{1+\sqrt{-23}}{2}\Big\rangle=
\langle 2,\alpha\rangle.
\]

Por el Ejemplo \ref{Ej8.5}, se tiene que $2{\ma Z}[\zeta_{23}]=
\pL\overline{\pL}$. Por lo tanto $\pK{\ma Z}[\zeta_{23}] =\pL$
y $\overline{\pK}{\ma Z}[\zeta_{23}] =\overline{\pL}$ donde los
grados relativos $f(\pL|\pK)$ y $f(\overline{\pL}|\overline{\pK})$
son igual a $11$. En particular la norma $N:=N_{\cic {23}{}/
{\ma Q(\sqrt{-23})}}$
satisface $N\pL=\pK^{11}$ y $N\overline{\pL}=\overline{\pK}^{11}$.

Se tiene $\langle 2\rangle{\ma Z}[\zeta_{23}]=\langle 2,
\alpha\rangle \langle 2,\overline{\alpha}\rangle=\langle 2,\alpha
\rangle \langle 2,1-\alpha\rangle$. Ahora bien
\[
N_{\cic {23}{}/{\ma Q}}\pL=\langle 2\rangle^{11}=
N_{{\ma Q}(\sqrt{-23})/{\ma Q}} N_{\cic {23}{}/{\ma Q}(\sqrt{-23})} \pL=
N_{{\ma Q}(\sqrt{-23})_/{\ma Q}}\pK^{11}.
\]

Veamos que $\pK$ no es principal. Sea $\beta =a+b\alpha=
a+b\Big(\frac{1+\sqrt{-23}}{2}\Big)\in{\ma Z}[\alpha]={\ma Z}\Big[
\frac{1+\sqrt{-23}}{2}\Big]$, $a,b\in{\ma Z}$. Entonces
$N\beta=(a+b)^2+5b^2-ab$, es decir, $N\beta=2$ no tiene
soluci\'on para $a,b\in{\ma Z}$.

En caso de que $\pK$ fuese principal, digamos $\pK=\langle
\beta\rangle$, entonces $N\pK=\langle 2\rangle =\langle N\beta
\rangle$ y por lo tanto $0<N\beta =\beta\overline{\beta}=2$ lo 
cual es imposible.
Por lo tanto $\pK$ no es principal y el n\'umero de clase de
${\ma Q}(\sqrt{-23})$ es mayor a $1$.

Ahora bien, si $\beta=2-\alpha=\frac{3-\sqrt{-23}}{2}$, entonces
$N\beta =8=2^3$. Por tanto 
\[
\pK^3=\langle\beta\rangle =\Big\langle 2-\frac{1+\sqrt{-23}}{2}
\Big\rangle = \Big\langle\frac{3-\sqrt{-23}}{2}\Big\rangle
\]
es principal pues $\beta\in\pK^3$, $\langle 2\rangle =\langle N\beta
\rangle\subseteq \langle 2\rangle$. Por otro lado 
se sigue que $\pL$ no es principal pues $N_{\cic {23}{}/{\ma Q}(
\sqrt{-23})}\pL=\pK^{11}=\pK^9\pK^2=\langle\beta\rangle
\pK^2$ y $\pK^2$ no puede ser principal pues si lo fuese
se tendr\'ia que $\pK=\pK^{-2}\pK^3$ ser\'ia principal.
Adem\'as $\pL^3$ si es principal pues $\pK^3
{\ma Z}[\zeta_{23}]=
\beta{\ma Z}[\zeta_{23}]=\pL^3=\langle \beta\rangle$. De esto
se sigue que $3$ divide al n\'umero de clase de $\cic {23}{}$ y
$\cic {23}{}$ no tiene n\'umero de clase $1$.

De hecho se tiene que si $h_E$\index{numero de clase de un
campo@n\'umero de clase de un
campo}\label{numero de clase de un campo} denota el n\'umero
de clase de un campo $E$, se tiene que $h_{\cic {23}{}}=h_{
{\ma Q}(\sqrt{23}}=3$. Esto fue calculado por Helmut Hasse
\cite{Has64}. Se tiene que si $K_H$ denota al campo de clase
de Hilbert de $K={\ma Q}(\sqrt{-23})$, entonces
\[
K_H={\ma Q}\big(\sqrt{-23},\sqrt[3]{(25+3\sqrt{69})/2}+
\sqrt[3]{(25-3\sqrt{69})/2}\big).
\]
M\'as a\'un, si $\alpha=\sqrt[3]{(25+3\sqrt{69})/2}+
\sqrt[3]{(25-3\sqrt{69})/2}$, entonces $\alpha^3-3\alpha-25=0$.
El discriminante de $X^3-3X-25$ es 
\[
D=-(4(-3)^3+27(-25)^2)=-3^6\cdot 23.
\]
\end{ejemplo}

Ahora bien, nos falta describir el tipo de descomposici\'on
de un primo $p$ en $\cic n{}$ cuando $p|n$. Sin embargo
esto se sigue de lo que ya sabemos hasta ahora. Como de
costumbre suponemos que $n\not\equiv 2\bmod 4$.

Sean $p|n$, $a\geq 1$, $a\in{\ma N}$ tal que $p^a|n$ y $p^{a+1}
\nmid n$, lo cual lo escribiremos $p^a\| n$. Entonces por la 
Proposici\'on \ref{P1.2.6'} se tiene
\[
\cic n{}=\cic pa \cic {{n/p^a}} {}, \qquad \cic pa \cap \cic {{n/p^a}}{}=
{\ma Q}.
\]
Entonces, puesto que $p$ es totalmente ramificado en
$\cic pa /{\ma Q}$ y no ramificado en $\cic {{n/p^a}}{}/{\ma Q}$, se tiene
que
\begin{align*}
e&=\varphi(p^a)=p^{a-1}(p-1),\\
f&=o(p\bmod n/p^a),\\
g&=\frac{\varphi(n/p^a)}{o(p\bmod n/p^a)}.
\end{align*}

Resumimos esto en el siguiente resultado.

\begin{teorema}\label{T8.10}
Sea $n\in{\ma N}$, $n>1$, $n\not\equiv 2\bmod 4$. Sea $p$
un primo y sea $a$ la potencia exacta de $p$ que divide a $n$:
$p^a|n$, $p^{a+1}\nmid n$ ($p^a\| n$), $a\geq 0$. Entonces si
\[
p{\ma Z}[\zeta_n]=\pK_1^e\cdots \pK_g^e, \qquad \big[
{\ma Z}[\zeta_n]/\pK_i:{\ma Z}/p{\ma Z}\big]=f
\]
con $\pK_1,\ldots,\pK_g$ primos distintos de grado $f$
cada uno, se tiene
\begin{align*}
e&=\varphi(p^a)=p^{a-1}(p-1),\\
f&=o(p\bmod n/p^a),\\
g&=\frac{\varphi(n/p^a)}{o(p\bmod n/p^a)}. \tag*{$\fin$}
\end{align*}
\end{teorema}

Ahora veamos los grupos de inercia y de descomposici\'on de 
un primo $\pL$ en $\cic n{}$ sobre $p$. Sea $p^a\| n$. Entonces
se tiene que si $I:=I(\pL|p)$ y $D:=D(\pL|p)$ son los grupos
de inercia y de descomposici\'on de $\pL$ sobre $p$ en $\cic n{}/
{\ma Q}$, entonces
\begin{align*}
|I|&=e=\varphi(p^a)=p^{a-1}(p-1),\\
|D|&=ef=p^{a-1}(p-1)\cdot o(p\bmod n/p^a).
\end{align*}

Ahora bien, si $\pK:=\pL\cap \cic {{n/p^a}}{}$, se tiene que
$\pL/\pK$ es totalmente ramificado y $\pK/p$ es no ramificado.
\[
\xymatrix{
\pL &\cic n{}\ar@{-}[d]^{\qquad \Big\} \pL/\pK \text{\ 
es totalmente ramificado}}\\
\pK&\cic {{n/p^a}}{}\ar@{-}[d]^{\qquad \Big\} \pK/p \text{\ es
no ramificado}}\\
p& {\ma Q}
}
\]

Por lo tanto $I=\Gal(\cic n{}/\cic {{n/p^a}}{})$. Por otro lado, puesto
que $\cic n{}=\cic pa\cdot \cic {{n/p^a}}{}$ y $\cic pa\cap \cic {{n/p^a}}
{}={\ma Q}$, se tiene 
\begin{gather*}
\Gal(\cic n{}/{\ma Q})\cong
\Gal(\cic n{}/\cic {{n/p^a}}{})\times \Gal(\cic n{}/\cic pa)\quad  \text{y}\\
\Gal(\cic n{}/\cic {{n/p^a}}{})\cong \Gal(\cic pa/{\ma Q});\\
\Gal(\cic n{}/\cic pa)\cong \Gal(\cic {{n/p^a}}/{\ma Q}).\\
\xymatrix{
&\cic n{}\ar@{-}[dl]|{\backslash} \ar@{-}[dr]|{//}\\
\cic pa\ar@{-}[dr]|{//}&&\cic {{n/p^a}}{}\ar@{-}[dl]|{\backslash}\\
&{\ma Q}=\cic pa\cap \cic {n/{p^a}}{}
}
\end{gather*}
En particular $I\cong \Gal(\cic pa/{\ma Q})\cong U_{p^a}$.

Otra descripci\'on de $I$ la podemos obtener mediante
la sucesi\'on exacta de grupos abelianos:
\[
1\longto \Gal(\cic n{}/\cic {{n/p^a}}{})\longto \Gal(\cic n{}/{\ma Q})
\stackrel{\rest}{\longto} \Gal(\cic {{n/p^a}}{}/{\ma Q})\longto 1.
\]

Con los isomorfismos $\Gal(\cic n{}/{\ma Q})\stackrel{h}{\cong} U_n$
y $\Gal(\cic {{n/p^a}}{}/{\ma Q})\stackrel{k}{\cong} U_{n/p^a}$ y
la proyecci\'on natural 
\begin{eqnarray*}
\pi\colon U_n &\longto & U_{n/p^a}\\
x\bmod n&\longmapsto & x\bmod n/p^a,
\end{eqnarray*}
y $\ker \pi=D_{n,n/p^a}=\{x\bmod n\mid x\equiv 1\bmod n/p^a\}$
obtenemos el siguiente diagrama conmutativo donde las filas son
exactas
\begin{scriptsize}
\[
\xymatrix{
1\ar[r]&\Gal(\cic n{}/\cic {{n/p^a}}{})\ar[d]^{h'}\ar[r]^{i}&
\Gal(\cic n{}/{\ma Q})\ar[d]^h\ar[r]^{\rest}&\Gal(\cic {{n/p^a}}{}/{\ma Q})
\ar[r]\ar[d]^k & 1\\
1\ar[r] &D_{n,n/p^a}\ar[r]^j & U_n\ar[r]^{\pi}& U_{n/p^a}\ar[r]&1
}
\]
\end{scriptsize}

\noindent
donde $i, j$ son las inyecciones naturales y $h'$ es 
el mapeo restricci\'on
de $h$ a $\Gal(\cic n{}/\cic {{n/p^a}}{})$. Entonces
\[
I=\Gal(\cic n{}/\cic {{n/p^a}}{})\cong D_{n,n/p^a}\cong U_{p^a}
\cong\Gal({\ma Q}(\zeta_{p^a})/{\ma Q}).
\]

Por otro lado, el grupo de descomposici\'on de $\pK/p$ en
$\cic {{n/p^a}}{}/{\ma Q}$ est\'a dado por el automorfismo de
Frobenius $\sigma_p(\zeta_{n/p^a})=\sigma^p_{n/p^a}$ el cual
corresponde a $p\bmod (n/p^a)$ y es de orden $f$. El grupo
de descomposici\'on $D(\pL| p)$ es isomorfo a $D(\pK|p)\times
I(\pL|p)$ y por tanto a $\langle p\bmod (n/p^a)\rangle \times
U_{p^a}$.

Precisemos un poco m\'as al grupo $D(\pL| p)$. Se tiene
que $D(\pL| p)=\langle \sigma,I(\pL| p)\rangle$ donde
$\sigma|_{{\ma Q}(\zeta_{n/p^a})}$ satisface $\bar{\sigma}(
\zeta_{n/p^a})=\zeta_{n/p^a}^p$. Digamos que $\sigma(\zeta_n)=
\zeta_n^b$ con $\mcd(b,n)=1$. Entonces se tiene
\[
\zeta_n^{bp^a}=\sigma(\zeta_n^{p^a})=\sigma(\zeta_{n/p^a})=
\zeta_{n/p^a}^b=\zeta_{n/p^a}^p,
\]
lo cual implica que $b\equiv p\bmod (n/p^a)$.

Puesto que $\mcd(p^a,n/p^a)=1$, por el Teorema Chino del Residuo,
se tiene que existe $b\in {\ma Z}$ tal que $b\equiv p\bmod (n/p^a)$
y $b\equiv \alpha \bmod p^a$ con $\mcd(\alpha,p)=1$, por ejemplo 
$\alpha=1$. Sea $d=\mcd(b,n)$. Veamos que $d=1$. Si $d\neq 1$,
existe un n\'umero primo que divide a $d$.
Si $p\mid d$ entonces $p\mid b$ y de la congruencia $b\equiv p
\bmod n/p^a$ se sigue que $p\mid n/p^a$ lo cual es absurdo. Por tanto
existe un primo $q\neq p$ tal que $q\mid d$. Puesto que $d\mid n$,
y $q\neq p$, se sigue que $q\mid (n/p^a)$ y adem\'as $q\mid b$. De
la congruencia $b\equiv p\bmod (n/p^a)$ se sigue que $q\mid p$ lo
cual es absurdo. En resumen $d=\mcd(b,n)=1$.

Seleccionamos $b$ de la discusi\'on anterior satisfaciendo $b\equiv
p\bmod (n/p^a)$ y $b\equiv 1\bmod p^a$. Entonces tenemos que
$o(p\bmod (n/p^a))=o(b\bmod (n/p^a))=f$ y puesto que $b^i\equiv 
1\bmod p^a$ para toda $i$, se sigue que $o(b\bmod n)=f$. Sea
$\sigma$ dado por esta $b$, es decir, $\sigma(\zeta_n)=\zeta_n^b$.
Entonces veamos que $\langle \sigma\rangle\cap D_{n,n/p^a}=\{1\}$
pues si $\sigma^j\in D_{n,n/p^a}$ entonces $b^j\equiv 1\bmod (n/p^a)$
lo cual implica que $f\mid j$ y por tanto $\sigma^j=\Id$. Hemos obtenido
que
\begin{align*}
D(\pL| p)&\cong \langle b\rangle \times D_{n,n/p^a}\cong
\langle\sigma\rangle\times\Gal({\ma Q}(\zeta_n)/{\ma Q}(\zeta_{n/p^a}))\\
&=\langle p\bmod n/p^a, \Gal({\ma Q}(\zeta_n)/{\ma Q}(\zeta_{n/p^a}))
\rangle.
\end{align*}

\begin{teorema}\label{T8.11} Con las notaciones anteriores se
tiene que para cualquier $n\in{\ma N}$ y para cualquier n\'umero
primo $p$ en ${\ma Z}$, si $p^a$ es la potencia exacta de $p$
que divide a $n$, tenemos
\begin{gather*}
I(\pL|p)\cong I(\pK|p)\cong D_{n,n/p^a}\cong U_{p^a},\\
D(\pL|p)\cong D(\pK|p)\times I(\pL|p)\cong \langle p\bmod n/p^a
\rangle \times U_{p^a}. \tag*{$\fin$}
\end{gather*}
\end{teorema}

Podemos explicitar con ``{\em grupos de congruencias}''
la descomposici\'on de primos no ramificados.

\begin{teorema}\label{T3.1.1} Sean $n\in {\ma N}$ y $L\subseteq
{\ma Q}(\zeta_n)$ que corresponde al subgrupo $H\subseteq U_n
=\*{\big({\ma Z}/n{\ma Z}\big)}\cong \Gal({\ma Q}(\zeta_n){\ma Q})$
con $L={\ma Q}(\zeta_n)^H$. Si $p$ es un n\'umero primo con 
$p\nmid n$, entonces:
\las
\item $p$ es no ramificado.

\item M\'as precisamente, si $f$ es m{\'\i}nimo natural tal que 
$p^f\bmod n\in H$, entonces $p\o_L=\pK_1\cdots \pK_r$ con
$r=[L:{\ma Q}]/f$ y cada $\pK_i$ es de grado $f$.

\item $p$ es totalmente descompuesto en $L\iff p\equiv 1 \bmod n \in H$.
\end{list}
\end{teorema}

\begin{proof}
(1). Es el Corolario \ref{C1.2.1.9}.

(2). Se tiene $[L:{\ma Q}]=e_pf_pg_p=f_pg_p$.
Veamos que $f=f_p$. Sea $\varphi_p$ el automorfismo de 
Frobenius de $p$ en ${\ma Q}(\zeta_n)/{\ma Q}$, $\varphi_p(\zeta_n)
=\zeta_n^p$ (ver la discusi\'on despu\'es de la Proposici\'on
\ref{P8.1}). Si $\varphi$ es el correspondiente al automorfismo de 
Frobenius de $p$ en $L/{\ma Q}$,
entonces $\varphi=\varphi_p\bmod H$, es decir $\varphi=\varphi_p|L$.
\[
\xymatrix{{\ma Q}(\zeta_n)\ar@{-}[d]^H\\ L\ar@{-}[d]\\ {\ma Q}}
\]
Entonces $f_p=o(\varphi)$ y $\varphi^s=\Id\iff \varphi_p^s|_L=\Id_L
\iff \varphi_p^s\in H\cong \Gal({\ma Q}(\zeta_n)/L)$. Por tanto
$f_p=o(\varphi)=\min\{s\mid \varphi_p^s\in H\} =\min\{s\mid
\varphi_p^s(\zeta_n\mapsto \zeta_n^{p^s})\in H\}=
\min\{s\mid p^s \bmod n\in H\}=f$. 

(3) Se tiene que $p$ es totalmente descompuesto
$\iff f=1\iff$ $p\bmod n \in H$. $\fin$
\end{proof}

\begin{corolario}\label{C3.1.2.CC}
El primo $p$ es totalmente descompuesto
en ${\ma Q}(\zeta_n)^+={\ma Q}(\zeta_n+\zeta_n^{-1})\iff
p\equiv \pm 1 \bmod n$.

Similarmente, $p$ es totalmente descompuesto en ${\ma Q}(\zeta_n)
\iff p\equiv 1\bmod n$.
\end{corolario}

\begin{proof} En el primer caso tenemos $H=\{1,J\}=\{1,-1\}$,
donde $J$ denota la conjugaci\'on compleja. Por tanto $p$ es totalmente
descompuesto en ${\ma Q}(\zeta_n)^+\iff p\bmod n\in H\iff p\equiv \pm 1
\bmod n$. El segundo caso es similar. $\fin$
\end{proof}

Para explicitar los generadores de los ideales primos de $\pL$
sobre un primo $p$ en $\cic n{}/{\ma Q}$, usamos el Teorema de
Kummer \ref{T8.12}.

En el caso $\cic n{}/{\ma Q}$, en la notaci\'on del
Teorema de Kummer, tenemos $A={\ma Z}$, $\coc A=
{\ma Q}$, $E=\cic n{}$ y $B={\cal O}_{\cic n{}}={\ma Z}[\zeta_n]$,
$\Irr(x,\zeta_n,{\ma Q})=\psi_n(x)$, $\psi_n(x)\bmod p=
\overline{\psi_n(x)}=\big(\overline{P_1(x)}\cdots \overline{P_g(x)}\big)^{e}$,
entonces $p{\ma Z}[\zeta_n]=(\pL_1\cdots \pL_g)^e$, $\pL_i=
\langle p, P_i(\zeta_n)\rangle$.

En el caso particular de $p\equiv 1\bmod n$, $p$ se descompone
totalmente, $\psi_n(x)=\prod_{\mcd(i,n)=1}(x-\zeta_n^i)$, 
$\overline{\psi_n(x)}=\overline{(x-a_1)}\cdots \overline{(x-a_{
\varphi(n)})}$, $P_i(x)=x-a_i\in{\ma Z}[x]$, $1\leq i \leq \varphi(n)$
para algunos $a_i\in{\ma Z}$ tales que $o(a_i\bmod p)=n$ y
$\pL_i=\langle p,\zeta_n-a_i\rangle$, $P_i(\zeta_n)=\zeta_n-a_i$.

\section{Subcampos de $\cic n{}$}\label{S4.3}

Por el Teorema de Kronecker--Weber, toda extensi\'on abeliana
finita de ${\ma Q}$ est\'a contenida en alg\'un $\cic n{}$. Por 
tanto conocer la aritm\'etica de los subcampos de $\cic n{}$ 
equivale a conocer la aritm\'etica de las extensiones abelianas
finitas de ${\ma Q}$.

Primero notemos que si $K/{\ma Q}$ es cualquier extensi\'on
de Galois y si $J$ es la conjugaci\'on compleja, entonces
para $j:=J|_K$, $j=\Id\iff K\subseteq {\ma R}$. En caso de que
$j\neq \Id$, $K^j=K\cap {\ma R}=:K^+$ y $[K:K^+]=2$. El
campo $K^+$ se llama el {\em campo real\index{campo
real de un campo}} de $K$.

Notemos que lo anterior no se cumple cuando $K/{\ma Q}$
no es de Galois. Por ejemplo, si $K:={\ma Q}(\zeta_n \sqrt[n]{2})$,
con $n\in{\ma N}$, $n\geq 3$, entonces $K^+=K\cap
{\ma R}={\ma Q}$ y $[K:K^+]=[K:{\ma Q}]=n$.

Sea $K:=\cic n{}$. Entonces $J|_{\cic n{}}=j$ se identifica con
$-1\in U_n$.

\begin{lema}\label{L10.1} Se tiene
\[
\cic n{}^+=\cic n{}\cap {\ma R}={\ma Q}(\zeta_n+\zeta_n^{-1})=
{\ma Q}(\zeta_n+\overline{\zeta}_n)=\cic n{}^J=
{\ma Q}\Big(2\cos \frac{2\pi}{n}\Big).
\]
Adem\'as $\Irr(\zeta_n,x, \cic n{}^+)=x^2-(\zeta_n+\zeta_n^{-1})
x+1$.
\end{lema}

\begin{proof}
Sea $\alpha=\zeta_n+\zeta_n^{-1}=\frac{\zeta_n^2+1}{\zeta_n}$.
Por tanto $\zeta_n^2-\zeta_n \alpha+1=0$. En particular,
$\cic n{}={\ma Q}(\alpha)(\zeta_n)$, $[\cic n{}:{\ma Q}(\alpha)]
\leq 2=[\cic n{}:\cic n{}^+]$. El resultado se sigue del
hecho de que $\alpha\in{\ma R}$. $\fin$
\end{proof}

Ahora bien $\Gal(\cic n{}/\cic n{}^+)=\{1,J\}\cong \{\pm 1\}$, donde
$J$ es la conjugaci\'on compleja la cual es identificada con
$-1\in U_n$ puesto que $J$ corresponde al mapeo $\zeta_n\to
\zeta_n^J=\overline{\zeta_n}=\zeta_n^{-1}$.
Se tiene $\Gal(\cic n{}^+/{\ma Q})\cong 
U_n/\{\pm 1\}$.

Para cualquier $n\geq 3$, $n\not\equiv 2\bmod 4$, se tiene 
que la extensi\'on $\cic n{}/\cic n{}^+$ 
es ramificada en todos los primos
infinitos pues todos los primos infinitos de $\cic n{}^+$ son
reales y todos los primos infinitos en $\cic n{}$ son
complejos. En el caso de los primos finitos, tenemos una
diferencia entre si $n$ es potencia de un n\'umero primo y 
cuando $n$ es dividido por al menos dos n\'umeros primos
distintos.

\begin{teorema}\label{T10.2}{\ }

\las
\item Si $n=p^m$ donde $p$ es primo y $m\in{\ma N}$, entonces
$\cic n{}/\cic n{}^+$ es ramificada en los primos de $\cic n{}^+$
encima de $p$ y en los primos infinitos y es no ramificada en
ning\'un otro primo.

\item Si $n$ es dividido por al menos dos primos distintos,
$n\not\equiv 2\bmod 4$, entonces $\cic n{}/\cic n{}^+$ es
ramificada \'unicamente en los primos infinitos, es decir, es no
ramificada en todos los primos finitos.
\end{list}
\end{teorema}

\begin{proof}
Puesto que si $n=p^m$, $p$ es totalmente ramificada en
$\cic pm/{\ma Q}$, solo falta probar que ning\'un primo finito
se ramifica en $\cic n{}/\cic n{}^+$ cuando $n$ es dividido
por al menos dos primos distintos.

Sean $p,q$ dos n\'umeros primos impares distintos que dividen
a $n$ o $p$ impar y $q=4$. Entonces $\zeta_p, \zeta_q\in
\cic n{}\setminus \cic n{}^+$. Por tanto $\cic n{}=
{\ma Q}(\zeta_n+\zeta_n^{-1},\zeta_p) = {\ma Q}(\zeta_n+
\zeta_n^{-1},\zeta_q)=\cic n{}^+(\zeta_p)=\cic n{}^+(\zeta_q)$.

Ahora bien, al adjuntar $\zeta_p$ a $\cic n{}^+$ el \'unico primo
finito posible a ramificarse en $\cic n{}=\cic n{}^+
(\zeta_p)/\cic n{}^+$ es $p$ pero
$p$ no es ramificado en $\cic n{}=\cic n{}^+(\zeta_q)/\cic n{}^+$ pues
el \'unico primo finito posible a ramificarse en esta \'ultima
extensi\'on es $q$ si $q$ es impar o $2$ si $q=4$. Por lo
tanto $\cic n{}/\cic n{}^+$ es no ramificado en los primos
finitos. $\fin$

\end{proof}

\subsection{Subcampos de $\cic 2m$}\label{S4.10.1}

Consideremos los subcampos de $\cic 2m$ con $m\geq 2$. 
En este caso, tenemos que $\Gal(\cic 2m/{\ma Q})\cong U_{2^m}\cong
C_2\times C_{2^{m-2}}\cong \Gal(\cic 4{}/{\ma Q})
\times \Gal(\cic 2m^+/{\ma Q})$. Recordemos que la sucesi\'on
exacta
\[
1\longto D_{2^m,4}\longto U_{2^m}\stackrel{\pi}{\longto}
U_4\longto 1
\]
donde $\pi$ es el epimorfismo natural, se escinde pues
$D_{2^m,4}$ es c{\'\i}clico de orden $2^{m-2}$ que es de orden
maximo en $U_{2^m}$. Adem\'as, aunque esto no es necesario
para nuestro estudio,
se tiene que $D_{2^m,4}$ es generado por $1+2^2=5$ y $U_4=
\{\pm 1\}$. El encaje natural $i\colon U_4\to U_{2^m}$,
$i(-1)=-1$,  satisface
$\pi\circ i=\Id_{U_4}$ y por tanto es el mapeo de escisi\'on.
En otras palabras, $U_{2^m}\cong U_4\times D_{2^m.4}\cong
C_2\times C_{2^{m-2}}$.

Sean $\langle a\rangle = U_4\cong C_2$, $\langle b\rangle
=D_{2^m,4}\cong C_{2^{m-2}}$. Notemos que $\cic 2m^{U_{
4}}=\cic 2m^{\{1,J\}}=\cic 2m^+$ y $\cic 2m^{D_{2^m,4}}=
\cic 2m^{\Gal(\cic 2m/\cic 4{})} = \cic 4{}={\ma Q}(i)={\ma Q}
(\sqrt{-1})$.
\[
\xymatrix{
\cic 2m^+\ar@{-}[rr]_{U_4=\langle a\rangle \cong \{1,J\}}
\ar@{-}[d]&&\cic 2m\ar@{-}[d]^{\Gal(\cic 2m/\cic 4{})\cong
D_{2^m.4}\cong C_{2^{m-2}}\cong\langle b\rangle}\\
{\ma Q}\ar@{-}[rr]_{U_4}&&\cic 4{}
}
\]

Los subgrupos de orden $2^{r}$ de $U_{2^m}$
con $1\leq r \leq m-1$ son:
\las
\item $r=m-1$: $U_{2^m}=\langle a,b\rangle$;

\item $1\leq r\leq m-2$: 
\begin{itemize}
\item $\langle a,b^{2^{m-1-r}}\rangle = U_4\times D_{2^m,
2^{m+1-r}}\cong C_2\times C_{2^{r-1}}$;
\item $\langle b^{2^{m-2-r}}\rangle = D_{2^m,2^{m-r}}\cong C_{2^r}$;
\item $\langle a b^{2^{m-2-r}}\rangle\cong C_{2^r}$.
\end{itemize}
\end{list}

Los respectivos campos fijos son (ver Teorema \ref{T1.4.4*}):
\las
\item ${\ma Q}(\zeta_{2^m})^{U_{2^m}}={\ma Q}$;
\item 
\begin{itemize}
\item ${\ma Q}(\zeta_{2^m})^{U_4\times D_{2^m,2^{m+1-r}}}\subseteq
{\ma R}$. Por lo tanto ${\ma Q}(\zeta_{2^m})^{\langle a,b^{2^{m-1-r}}\rangle}
={\ma Q}(\zeta_{2^t})^+$ para alg\'un $t$ y 
\[
[{\ma Q}(\zeta_{2^m}):{\ma Q}(\zeta_{2^t})^+]
=\frac{[{\ma Q}(\zeta_{2^m}):{\ma Q}]}{[{\ma Q}(\zeta_{2^t})^+:{\ma Q}]}=
\frac{2^{m-1}}{2^{t-2}}=2^{m-t+1}=2^r,
\]
lo cual es equivalente a $r=m-t+1$ si y solamente si $t=m-r+1$.
Por tanto
\[
\cic 2m^{U_4\times D_{2^m,2^{m-r+1}}}=\cic 2{m-r+1}^+=
{\ma Q}(\zeta_{2^{m-r+1}}+\zeta^{-1}_{2^{m-r+1}}).
\]
\item $\cic 2m^{\langle b^{2^{m-r-2}}\rangle}=\cic 2m^{D_{2^m,2^{m-r}}}
=\cic 2t$ para alg\'un $t$ y
\[
[\cic 2m:\cic 2t]=\frac{[\cic 2m:{\ma Q}]}{[\cic 2t:{\ma Q}]}=\frac{2^{m-1}}
{2^{t-1}}=2^{m-t}=2^r
\]
lo cual es equivalente a $m-t=r$ si y solamente si $t=m-r$. Por tanto
\[
\cic 2m^{D_{2^m,2^{m-r}}}=\cic 2{m-r}.
\]
\item Sea $E=\cic 2m^{\langle a b^{2^{m-r-2}}\rangle}$.
Se tiene $(ab^{2^{m-r-2}})^2=a^2 b^{2^{m-r-1}}=b^{2^{m-r-1}}$.
Por tanto $E\subseteq \cic 2m^{\langle b^{2^{m-r-1}}\rangle}=
\cic 2m^{D_{2^m,2^{m-r+1}}}=\cic 2{m-r+1}$. Por otro lado
$a\notin \langle ab^{2^{m-2-r}}\rangle$ por lo que $E
\nsubseteq {\ma R}$. Adem\'as $[\cic 2{m-r+1}:E]=2$.
\[
\xymatrix{
\cic 2m^+\ar@{-}[rr]^2\ar@{-}[dd]_{2^{r-1}}&&\cic 2m\ar@{-}[dd]^{2^{r-1}}\ar@{-}[dddl]_{2^r}
|!{[ddll];[dd]}\hole\\ \\
\cic 2{m-r+1}^+\ar@{-}[rr]^2\ar@{-}[dd]_2&&\cic 2{m-r+1}\ar@{-}[dd]^2\\
&E\ar@{-}[dl]_2\ar@{-}[ur]_2\\
\cic 2{m-r}^+\ar@{-}[rr]^2\ar@{-}[dd]_{2^{m-r-2}}&&\cic 2{m-r}\ar@{-}[dd]^{2^{m-r-2}}\\ \\
{\ma Q}\ar@{-}[rr]^2&&\cic 4{}
}
\]

Sea $F:={\ma Q}((\zeta_{2^{m-r+1}}+\zeta^{-1}_{2^{m-r+1}})\zeta_4)\supseteq
{\ma Q}(\zeta_{2^{m-r}}+\zeta^{-1}_{2^{m-r}})$. Adem\'as $[F:{\ma Q}(
\zeta_{2^{m-r}}+\zeta^{-1}_{2^{m-r}})]=2$ y $F\neq \cic 2{m-r+1}^+$ pues
$F\nsubseteq {\ma R}$ y $F\neq \cic 2{m-r}$. Se sigue que $F=E$.

Se puede verificar que $\zeta_4(\zeta_{2^{m-r+1}}+\zeta^{-1}_{2^{m-r+1}})=
\zeta_{2^{m-r+1}}^{\alpha}-\zeta^{-\alpha}_{2^{m-r+1}}$ con $\mcd(\alpha,2)=1$.
De hecho, $\alpha=2^{m-r-1}-1$. Por tanto
\[
E=\cic 2m^{\langle a b^{2^{m-r-2}}\rangle}=
{\ma Q}(\zeta_{2^{m-r+1}}-\zeta^{-1}_{2^{m-r+1}}).
\]
\end{itemize}
\end{list}

Podemos hacer m\'as expl{\'\i}cita la descripci\'on de estos 
subcampos. Sea $\alpha_r=\zeta_{2^{r}}+\zeta_{2^{r}}^{-1}$.
Entonces $\alpha_r^2=\zeta_{2^{r-1}}+\zeta_{2^{r-1}}^{-1}+2=
\alpha_{r-1}+2$, esto es, $\alpha_r=\sqrt{\alpha_{r-1}+2}$
donde escogemos el signo positivo pues $\alpha_r>0$. Entonces
tenemos 
\begin{align*}
\zeta_{2^{3}}+\zeta_{2^{3}}^{-1}&= \zeta_8+\zeta_8^{-1}=\sqrt{2};\\
\zeta_{2^{4}}+\zeta_{2^{4}}^{-1}&=\sqrt{2 + \sqrt{2}};\\
\intertext{y en general}
\zeta_{2^{m}}+\zeta_{2^{m}}^{-1}&=\underbrace{\sqrt{2+\sqrt{
2+\sqrt{2+\cdots+\sqrt{2}}}}}_{m-2}=\alpha_m.
\end{align*}

Por tanto tenemos (ver Obsevaci\'on \ref{R1.4.4**}):

\begin{teorema}\label{T10.1.1}
Para $1\leq r\leq m-2$, $\cic 2m$ tiene tres subcampos de
grado $2^r$ sobre ${\ma Q}$:
\begin{itemize}
\item $\cic 2{r+2}^+={\ma Q}(\alpha_{r+2})$;
\item $\cic 2{r+1}$;
\item ${\ma Q}(\zeta_4(\zeta_{2^{r+1}}+\zeta_{2^{r+1}}^{-1}))=
{\ma Q}(\zeta_{2^{r+2}}-\zeta_{2^{r+2}}^{-1})={\ma Q}(\sqrt{-
\alpha_{r+2}})$
\end{itemize}
donde $\alpha_{r+2}=\underbrace{\sqrt{2+\sqrt{
2+\sqrt{2+\cdots+\sqrt{2}}}}}_{r}$. M\'as a\'un tanto el primer como
 el tercer campo son extensiones c{\'\i}clicas de ${\ma Q}$ e 
 inclusive el segundo campo para el caso $r=1$. En el caso $r\geq
 2$, 
 \begin{align*}
 \Gal({\ma Q}(\alpha_{r+2})/{\ma Q})&=\Gal(\cic 2{r+2}^+/{\ma Q})
 \cong D_{2^{r+2},2^2}\cong C_{2^r};\\
 \Gal ({\ma Q}(\sqrt{-\alpha_{r+2}})/{\ma Q})&\cong C_{2^r};\\
 \Gal (\cic 2{r+1}/{\ma Q})&\cong \langle 1,J\rangle \times
 D_{2^{r+1},2^2} \cong \langle \pm 1 \rangle \times
D_{2^{r+1},2^2}.
 \end{align*}
 
 En el caso $r=1$, las extensiones son c{\'\i}clicas de grado $2$
 y estas son ${\ma Q}(\sqrt{2})$, ${\ma Q}(\sqrt{-2})$ y
 ${\ma Q}(i)={\ma Q}(\sqrt{-1}) =\cic 4{}$. $\fin$
 \[
 \xymatrix{
 \cic 2m^+\ar@{-}[rrrr]\ar@{-}[dd]&&&&\cic 2m\ar@{-}[dd]\\ \\
 {\ma Q}(\alpha_{r+2})=\cic 2{r+1}^+\ar@{-}[rrrr]\ar@{-}[ddd]_{
 \substack{\Gal(\cic 2{r+2}^+/{\ma Q})\cong\\ \cong
 D_{2^{r+2},2^2}\cong C_{2^r}}}&&&&\cic 2{r+2}
 \ar@{-}[dll]\ar@{-}[d]\\
 &&{\ma Q}(\zeta_4 \alpha_{r+2})\ar@{-}[ddll]_{\substack{
 \Gal({\ma Q}(
 \sqrt{-\alpha_{r+2}}/{\ma Q}))\cong\\
 \cong C_{2^r}}}&&\cic 2{r+1}\ar@{-}[dd]\ar@{-}[ddllll]^{
 \substack{\Gal(\cic 2{r+1}/{\ma Q})\cong U_{2^{r+1}}\cong\\
 \cong C_2\times C_{2^{r-1}}}}\\ \\ {\ma Q}\ar@{-}[rrrr]&&&&\cic 4{}
 }
 \]
\end{teorema}

\subsection{Subcampos de $\cic pm$, $p$ primo, $p>2$}\label{S4.10.2}

En el caso en que $p$ es impar
se tiene que $\Gal(\cic pr/{\ma Q})\cong 
U_{p^r}$ es un grupo 
c{\'\i}clico y que $U_{p^r}\cong C_{p-1}\times C^{p^{r-1}}$,
el cual es de orden $\varphi(p^r)=p^{r-1}(p-1)$.

Un elemento de orden $p^{r-1}$ en $U_{p^r}$ es $1+p$ el cual
corresponde a $\sigma \zeta_{p^r}=\zeta_{p^r}^{1+p} = 
\zeta_{p^r}\zeta_{p^{r-1}}$ (ver Subsecci\'on \ref
{S1.2.1}). Consideremos una ra{\'\i}z primitiva
m\'odulo $p$, digamos $a$, es decir, $U_p=\langle a\rangle$ y
sea $\theta(\zeta_p)=\zeta_p^a$. Pongamos $\Gal(\cic pr/{\ma Q})
=\langle \theta, \sigma \rangle$ donde extendemos a $\theta$ como
una extensi\'on de $\theta (\zeta_p)=\zeta_p^a$ a $\theta(\zeta_{
p^r})$ y $o(\theta)=p-1$. Por ejemplo, si $\theta (\zeta_{p^r})=
\zeta_{p^r}^{\alpha}$, entonces
 $\zeta_{p^r}^{p^{r-1}a}=\zeta_p^a=
\theta(\zeta_p) =\theta(\zeta_{p^r}^{p^{r-1}})=\zeta_{p^r}^{p^{r-1}
\alpha}$. Es decir $
p^{r-1}a\equiv p^{r-1}\alpha\bmod p^r$, $a\equiv \alpha \bmod p$.
Adem\'as $\theta^{p-1}(\zeta_{p^r})=\zeta_{p^r}^{\alpha^{p-1}}=
\zeta_{p^r}$, esto es $\alpha^{p-1}\equiv 1\bmod p^r$. Es decir
$\theta(\zeta_{p^r})=\zeta_{p^r}^\alpha$,
$o(\alpha \bmod p^r)=p-1$.

Sea $T_r:=\cic pr^{\langle \theta\rangle}$, $\cic p{}=\cic pr^{
\langle \sigma \rangle}$.
\begin{gather*}
\xymatrix{
T_r\ar@{-}[dd]_{\langle\sigma\rangle}^{p^{r-1}}\ar@{-}@/^2pc/[rrr]^{
\langle \theta\rangle}\ar@{-}[rr]
&&\cic pr^+\ar@{-}[dd]\ar@{-}[r]_2&\cic pr\ar@{-}[dd]^{\langle \sigma\rangle}_{p^{r-1}}\\ \\
{\ma Q}\ar@{-}@/_2pc/[rrr]_{\langle \theta\rangle}
\ar@{-}[rr]&&\cic p{}^+\ar@{-}[r]^2&
\cic p{}
}\\
\cic pr=T_r\cic p{} =T_r(\zeta_p);\quad T_r\cap \cic p{}={\ma Q};\\
\cic pr^+=T_r\cic p{}^+=T_r(\zeta_p+\zeta_p^{-1});\quad
T_r\subseteq {\ma R}.
\end{gather*}

Notemos que si $s|\varphi(p^r)$, entonces escribiendo $s=bp^t$
con $t\leq r-1$ y $b|p-1$. Por unicidad de las extensiones de
${\ma Q}$ contenidas en $\cic pr$ pues $U_{p^r}$ es
c{\'\i}clico, el \'unico campo de grado $s$ contenido en $\cic pr$ es
$M_b T_{t+1}$ donde $M_b$ es el \'unico subcampo de $\cic p{}$
de grado $b$ sobre ${\ma Q}$.
\[
\xymatrix{
T_r\ar@{-}[rrr]\ar@{-}[d]&&&\cic pr\ar@{-}[d]\\
T_{t+1}\ar@{-}[rr]\ar@{-}[d]_{p^t}&&M_bT_{t+1}\ar@{-}[r]\ar@{-}[d]&
\cic p{t+1}\ar@{-}[d]\\
{\ma Q}\ar@{-}[rr]_b\ar@{-}[rru]^{s=bp^t}&&M_b\ar@{-}[r]&\cic p{}
}
\]

A continuaci\'on describiremos de manera m\'as o menos expl{\'\i}cita
los campos $M_b$ y $T_{t+1}$ en general. Empezamos con $M_b$.
Sea $c:=\frac{p-1}{b}$. Se tiene $[\cic p{}:M_b]=c$, $[M_b:{\ma Q}]=
b$, $bc =p-1$.

Notemos que $M_b$ es el campo fijo bajo el subgrupo de $U_p$
de orden $c$, es decir, bajo $\langle \theta^b\rangle$ ya que
$o(\theta^b)=c$ y $\Gal(M_b/{\ma Q})\cong \frac{\Gal(\cic p{}/{\ma Q})}{
\Gal(\cic p{}/M_b)}=\frac{\langle\theta\rangle}{\langle \theta^b\rangle}
= \langle\theta\bmod \theta^b\rangle=\langle\overline{\theta}\rangle$.
Esto es, $M_b=\cic p{}^{\langle\theta^b\rangle}$.

Puesto que en general se tiene que $\theta\zeta_p=\zeta_p^a$, entonces
$\theta^j\zeta_p=\zeta_p^{a^j}$. Consideremos $
\mu:=\zeta_p+\theta^b\zeta_p+\cdots+\theta^{(c-1)b}\zeta_p$.
As{\'\i} $\theta^b\mu=\theta^b\zeta_p+\theta^{2b}\zeta_p+\cdots+
\theta^{cb}\zeta_p=\mu$, es decir, $\mu\in \cic p{}^{\langle \theta^b\rangle}=
M_b$ y por lo tanto $[{\ma Q}(\mu):{\ma Q}]\leq [M_b:{\ma Q}]=b$.

Ahora bien, $\mu=\sum_{i=0}^{c-1}\theta^{ib}\zeta_p=\sum_{i=0}^{
c-1}\zeta_p^{a^{ib}}$ y se tiene 
\[
\theta\mu = \sum_{i=0}^{c-1}(\theta
\zeta_p)^{a^{ib}}= \sum_{i=0}^{c-1}\zeta_p^{a^{ib+1}}=
\sum_{i=0}^{c-1}\big(\zeta_p^{a^{ib}}\big)^a
\]
y, en general,
$\theta^j\mu =  \sum_{i=0}^{c-1}\big(\zeta_p^{a^{ib}}\big)^{a^j}=
\sum_{i=0}^{c-1}\zeta_p^{a^{ib+j}}$. Se sigue que $\theta^j\mu\neq \mu$
para $1\leq j\leq b-1$. 

Notemos que si $s=ib+j=i'b+j'$, $0\leq i,i'\leq c-1$, $0\leq j,j'\leq b-1$
entonces $b|(j'-j)$ y $0\leq |j'-j|<|b|=b$ lo cual implica que $j'=j$ y
por tanto que $i=i'$. De esta forma obtenemos que $\{ib+j\mid
0\leq i\leq c-1, 0\leq j\leq b-1\}=\{0,1,\cdots, p-2\}$. Se sigue que
$\{a^{ib+j}\}_{\substack{0\leq i\leq c-1\\ 0\leq j\leq b-1}}=
\{a^0,a^1,\ldots, a^{p-2}\}=\{1,2,\ldots,p-1\}$ y en consecuencia
$\{\zeta_p^{a^{ib+j}}\mid 0\leq i\leq c-1,0\leq j\leq b-1\}=\{
\zeta_p,\ldots,\zeta_p^{p-1}\}$.

Esto es, tenemos que como $a$ es una ra{\'\i}z
primitiva m\'odulo $p$, se tiene que $\{a^{ib+j}\mid 0\leq i\leq c-1, 
0\leq j\leq b-1\}=\{a^0=a^{p-1}=1,\ldots, a^{p-2}\}$ son todos
distintos m\'odulo $p$ y $\{\zeta_p,\ldots, \zeta_p^{p-1}\}$ es
base de $\cic p{}/{\ma Q}$, por lo que $\mu,\theta\mu,\ldots, \theta^{b-1}
\mu$ son $b$ conjugados distintos de $\mu$ de donde obtenemos
$[{\ma Q}(\mu):{\ma Q}]\geq b$.
Se sigue que 
\[
M_b={\ma Q}(\mu),
\]
donde $\mu=\sum_{i=0}^{c-1}
\zeta_p^{a^{ib}}$.

Ahora bien, de esta misma forma tenemos que
$T_r={\ma Q}(\delta)$, donde 
\[
\delta=\zeta_{p^r}+\theta\zeta_{p^r}+\cdots+\theta^{p-2}\zeta_{p^r}=
\zeta_{p^r}+\zeta_{p^r}^{\alpha}+\cdots+\zeta_{p^r}^{\alpha^{p-2}}=
\sum_{i=0}^{p-2}\zeta_{p^r}^{\alpha^i}
\]
donde $o(\alpha\bmod p^r)=p-1$ y se tiene $\sigma\delta=
\sum_{i=0}^{p-2}(\sigma \zeta_{p^r})^{\alpha^i}=\sum_{i=0}^{p-2}
\zeta_{p^r}^{(1+p)\alpha^i}$ y $\theta\delta=\delta$.

Para $T_{t+1}$, $t\leq r-1$, se tiene $T_{t+1}={\ma Q}(\delta_{t+1})$
donde $\delta_{t+1}=\sum_{i=0}^{p-2}\zeta_{p^{t+1}}^{\alpha^i}$ con
$o(\alpha \bmod p^{t+1})=p-1$.

\begin{teorema}\label{T10.2.1} Sean $p>2$ un primo impar y $m\in{\ma N}$.
Los subcampos de $\cic pm$ est\'an dados por
\begin{gather*}
M_b\cdot T_r,\quad b|p-1,\quad 0\leq r\leq m-1,
\intertext{donde}
M_b={\ma Q}(\mu_b),\quad T_r={\ma Q}(\delta_r),\quad
\mu_b=\sum_{i=0}^{c-1}\zeta_p^{a^{ib}}, \quad \delta_r=
\sum_{i=0}^{p-2}\zeta_{p^r}^{\alpha^i},\\
o(a\bmod p)=p-1,\quad o(\alpha\bmod p^r)=p-1, \quad c=\frac{p-1}{b},\\
[M_b:{\ma Q}]=b,\quad [T_r:{\ma Q}]=p^r.
\intertext{Adem\'as}
\Gal(M_b/{\ma Q})=\langle\overline{\theta}\rangle\cong C_b,
\quad \Gal(T_r/{\ma Q})
=\langle\overline{\sigma}\rangle \cong C_{p^r},\\
\theta(\zeta_p)=\zeta_p^a,\quad \sigma(\zeta_{p^m})=\zeta_{p^m}^{1+p},
\quad \overline{\theta}=\theta\bmod\langle\theta^b \rangle,\quad
\overline{\sigma}=\sigma\bmod \langle\sigma^{p^r}\rangle.
\tag*{$\fin$}
\end{gather*}
\end{teorema}

En particular se tiene que $\cic pm^+$ es de grado $\frac{\varphi(p^m)}{2}=
\frac{p-1}{2} p^{m-1}$, esto es, $\cic pm^+=M_{(p-1)/2} T_{m-1}$.
Notemos que $M_{(p-1)/2}=\cic p{}^+={\ma Q}(\mu_{(p-1)/2})$, 
$\mu_{(p-1)/2}=\zeta_p+\zeta_p^{-1}$. En otras palabras tenemos
$\cic pm^+=\cic p{}^+ T_{m-1}$.

En general $\cic n{}^+={\ma Q}(\zeta_n+\zeta_n^{-1})$. De hecho,
si $\sigma\in\Gal(\cic n{}/{\ma Q})$ es tal que $\sigma\zeta_n=\zeta_n^a$, 
entonces $\sigma(\zeta_n+\zeta_n^{-1})=\zeta_n^a+\zeta_n^{-a}
=\zeta_n+\zeta_n^{-1}$ si y
solamente si $\zeta_n^a=\zeta_n$ o $\zeta_n^a=\zeta_n^{-1}$ pues
$\Gal(\cic n{}/\cic n{}^+)\cong \{\pm 1\}$.

Una pregunta que surge es: ?`Cu\'al es el campo $\cic n{}^{-}:=
{\ma Q}(\zeta_n-\zeta_n^{-1})$?

Notemos que si $\alpha:=\zeta_n-\zeta_n^{-1}=\frac{\zeta_n^2-1}{\zeta_n}$
entonces $\zeta_n^2-\alpha\zeta_n-1=0$, esto es, $\Irr(\zeta_n,x,\cic n{}^{-})
| x^2-\alpha x-1$ y por lo tanto $[\cic n{}:\cic n{}^{-}]\leq 2$. 
Adem\'as, $x^2-\alpha x-1=(x-\zeta_n)(x+\zeta_n^{-1})$.

Queremos determinar
cuando $\cic n{}^-=\cic n{}$. Tendremos que $\cic n{}\neq \cic {}n^-$ si y
solamente si existe $\sigma\in \Gal(\cic n{}/{\ma Q})\cong U_n$, $\sigma
\zeta_n=\zeta_n^a$ con $a\neq 1$ tal que $\sigma\alpha=\alpha$.
Cuando $n$ es par se tiene $4|n$ y si tomamos
$a=\frac{n}{2}-1$, entonces $\sigma\zeta_n
=\zeta_n^a=\zeta_n^{n/2}\zeta_n^{-1}=\zeta_2\zeta_n^{-1}=-\zeta_n^{-1}$.
Por tanto $\sigma\alpha=\alpha$. Ahora, $a\neq 1$ para $n>4$. En el caso 
$n=4$ se tiene $\zeta_n-\zeta_n^{-1}=\zeta_4-\zeta_4^{-1}=i-(-i)=2i$ y $\cic n{}=
\cic n{}^-$. Para $n>4$, $4|n$, $\cic n{}\neq \cic n{}^-$
y $\Gal(\cic n{}/\cic n{}^-)\cong \{1,\frac{n}{2}-1\}$. Notemos adem\'as
que $\cic n{}^+\neq \cic n{}^-$ pues $\frac{n}{2}-1\neq -1$.

Ahora sea $n$ impar. Supongamos que $\sigma$ satisface que
$\sigma(\zeta_n-\zeta_n^{-1})=\zeta_n^a-\zeta_n^{-a}=\zeta_n-\zeta_n^{-1}$.
Puesto que $\zeta_n^a-\zeta_n^{-a}=2i\sen\frac{2\pi a}{n}$ entonces
debemos tener $\sen\frac{2\pi a}{n}=\sen\frac{2\pi}{n}$ lo cual equivale
a que 
\[
\frac{2(a+1)}{n}=2m+1,\quad m\in{\ma Z}\quad\text{o}\quad
\frac{a-1}{n}=m\in{\ma Z}.
\]
En el caso $\frac{2(a+1)}{n}=2m+1$, $m\in {\ma Z}$, tenemos que,
puesto que $1\leq a\leq n-1$, entonces
$\frac{4}{n}\leq \frac{2(a+1)}{n}\leq 2$. Para $n=3$ se tiene $\cic 3{}$ es
una extensi\'on cuadr\'artica y $\zeta_3-\zeta_3^{-1}=2i\frac{\sqrt{3}}{2}=
\sqrt{-3}$ y por ende ${\ma Q}(\zeta_3-\zeta_3^{-1})=\cic 3{} ={\ma Q}(
\sqrt{-3})$.

Para $n\geq 5$, $0< \frac{4}{n}<1$ y por tanto $1\leq \frac{2(a+1)}{n}\leq 2$.
Por lo tanto $2m+1\in\{1,2\}$ y $2m+1$ es impar lo que implica $m=0$.
Se seguir\'ia que $\frac{2(a+1)}{n}=1$ y $a=\frac{n}{2}-1\notin{\ma Z}$.

En el caso $\frac{a-1}{n}=m\in{\ma Z}$, nuevamente, como $1\leq a\leq n-1$
obtenemos $0\leq \frac{a-1}{n}\leq 1-\frac{2}{n}<1$, esto es, $\frac{a-1}{n}=0$
lo cual implica $a=1$. 

En resumen, si $n$ es impar entonces $\zeta_n^a-
\zeta_n^{-a}=\zeta_n-\zeta_n^{-1}$ implica $a=1$ y en consecuencia
$\zeta_n-\zeta_n^{-1}$ no es fijado por ning\'un $\sigma\in \Gal(\cic n{}/
{\ma Q})\setminus\{\Id\}$ y por tanto ${\ma Q}(\zeta_n-\zeta_n^{-1})=
\cic n{}^-=\cic n{}$.

Hemos probado el siguiente resultado.

\begin{proposicion}\label{P10.2.-1}
Sean $n\not\equiv 2\bmod 4$ y $\cic n{}^-={\ma Q}(\zeta_n-\zeta_n^{-1})$. Entonces,
\lasa
\item Si $n|4$ y $n>4$, entonces $\cic n{}^-
\neq \cic n{}$ y $[\cic n{}:\cic n{}^-]=2$.

\item Si $n=4$ o si $n$ es impar, entonces $\cic n{}^-=\cic n{}$. $\fin$
\end{list}
\end{proposicion}

\begin{observacion}\label{O10.2.2}
La obtenci\'on de todos los subcampos de $\cic p{}$ 
con $p>2$, es posible
debido a que la base de potencias de $\zeta_p$ es de
hecho una base normal. Esto no sucede en general
\cite{Joh85}.
\end{observacion}

\begin{teorema}\label{T10.2.3}
Sea $\cic n{}$ un campo ciclot\'omico con $n\not\equiv 2\bmod 4$. 
Entonces $\{\zeta_n^{\sigma}\}_{\sigma\in \Gal(\cic n{}/{\ma Q})
\cong U_n}$ es una base normal de $\cic n{}/{\ma Q}$ si y solamente
si $n$ es libre de cuadrados.
\end{teorema}

\begin{proof}
(1)
Si $n$ no es libre de cuadrado, entonces existe un n\'umero primo $p$ con
$p^2|n$. En el caso $n=4$, $\zeta_4=i$, $\sigma \zeta_4=\zeta_4^3=
-i$ y $\{\zeta_4,\zeta_4^3\}=\{i,-i\}$ no es base de $\cic 4{}/{\ma Q}$.
En el caso general, $n\neq 4$, sea $\zeta_p=\zeta_n^{n/p}$, Ahora bien
\begin{gather}\label{Ec10.2.3'}
\sum_{j=0}^{p-1} \zeta_n^{1+jn/p}=\sum_{j=0}^{p-1}\zeta_n\zeta_p^j=
\zeta_n\psi_p(\zeta_p)=0.
\end{gather}

Ahora $\mcd (n,1+\frac{jn}{p})=1$ pues si $q$ es un n\'umero primo tal
que $q|n$ entonces $q| j\frac{n}{p}$ puesto que $q|\frac{n}{p}$. Por tanto
(\ref{Ec10.2.3'}) es una relaci\'on entre los 
elementos $\{\zeta_n^{\sigma}\}_{
\sigma\in U_n}$ por lo que $\{\zeta_n^{\sigma}\}_{
\sigma\in U_n}$ no es base de $\cic n{}/{\ma Q}$. 

En el caso especial $n=4m$ se tiene $\zeta_p=\zeta_2=-1$, $\psi_2(x)=1+x$,
$\psi_2(\zeta_2)=1+\zeta_2=0$.

(2) Sea ahora $n=p_1\cdots p_r$ libre de cuadrados, por lo que 
$p_i\neq 2$ para $1\leq i\leq r$. Sea $p$ un n\'umero primo
con $p\geq 3$. Entonces una base de $\cic p{}/{\ma Q}$ es $\{1,
\zeta_p,\ldots,\zeta_p^{p-2}\}$. Multiplicando por $\zeta_p$ se tiene
que $\{\zeta_p,\zeta_p^2,\ldots,\zeta_p^{p-1}\}=\{\zeta_p^{\sigma}\}_{
\sigma\in\Gal(\cic [{}/{\ma Q})}$ es una base normal.

Ahora probamos el resultado 
por inducci\'on en $r\geq 2$. Si $m=p_2\cdots p_r$,
suponemos que $\{
\zeta_m^{\sigma}\}_{\sigma\in U_m}$ es una base normal de $\cic m{}/
{\ma Q}$ y si $p=p_1$, $\{\zeta_p^{\mu}\}_{\mu\in U_p}$ es
una base normal de $\cic p{}/{\ma Q}$.

Ahora bien $\{\zeta_p^{\mu}\zeta_m^{\sigma}\}_{\mu\in U_p,
\sigma\in U_m}$ es una base de 
la extensi\'on $\cic m{}\cic p{}/{\ma Q}$ pues
$\mcd (m,p)=1$. Adem\'as $[\cic {pm}{}:\cic m{}]=[\cic p{}:{\ma Q}]
=\varphi(p)$
y $[\cic {pm}{}:\cic p{}]=[\cic m{}:{\ma Q}]=\varphi(m)$.
\[
\xymatrix{
\cic m{}\ar@{-}[rr]^{\varphi(p)}\ar@{-}[d]_{\varphi(m)}&&
\cic m{}\cic p{}\ar@{-}[d]^{\varphi(m)}\\
{\ma Q}\ar@{-}[rr]_{\varphi(p)}&&\cic p{}
}
\]
Se tiene que $\{\zeta_m^{\sigma}\}_{\sigma\in U_m}$ es
base de $\cic m{}\cic p{}/\cic p{}$ y de $\cic m{}/{\ma Q}$ y
$\{\zeta_p^{\mu}\}_{\mu\in U_p}$ es base de 
$\cic m{}\cic p{}/\cic m{}$ y de $\cic p{}/{\ma Q}$.
Adem\'as, del hecho de que $\mcd(m,p)=1$, se sigue que
$\cic m{}\cic p{}=\cic {mp}{}=\cic {pm}{}=\cic {p_1\cdots p_r}{}
=\cic n{}$.

Por otro lado, se tiene $\{\zeta_p^{\mu}\zeta_m^{\sigma}\}_{\mu
\in U_p, \sigma \in U_m}=\{\zeta_p^{\mu}\zeta_m^{\sigma}\}_{
(\mu,\sigma)\in U_p\times U_m\cong U_{pm}=U_n}$.
De esta forma $\{\zeta_p^{\mu}\zeta_m^{\sigma}\}_{(\mu,\sigma)}=
\{\zeta_p^a\zeta_m^b\}_{\substack{\mcd(a,p)=1\\ \mcd(b,m)=1}}=
\{\zeta_{pm}^{ma}\zeta_{pm}^{pb}\}_{a,b}=\{\zeta_{pm}^{ma+pb}\}_{
a,b}$. Se tiene $1\leq a \leq p-1$ y $1\leq b \leq m-1$ por lo que
$m\leq ma\leq mp-m$ y $p\leq pb\leq mp-p$. Por tanto $m+p
\leq ma+pb\leq 2mp-m-p$.

(i) Sea $\mcd(ma+pb,n)=d$.
Si $q$ es un n\'umero primo y $q|n$, entonces
$q\in\{p_1,p_2,\ldots,p_r\}$. Ahora bien, $p_1=p$ y $p|n$ y por tanto $p|pb$
pero $p\nmid ma$ pues $p\nmid m$ y $p\nmid a$ por lo que
$p\nmid d=\mcd(m,p)$. Para $2\leq i\leq r$, tendremos $p_i|n$,
$p_i|ma$ pero $p_i\nmid pb$ pues $p_i\nmid p$ y $p_i\nmid b$
debido a que $\mcd (b,m)=1$. Por tanto $p_i\nmid d$. Se
sigue que $d=1$.

(ii) Si $ma_1+pb_1=ma_2+pb_2$ con $a_1,a_2\in\{1,\ldots, p-1\}$ y
$b_1,b_2\in\{1,\ldots,m-1\}$ con $\mcd(b_i,m)=1$, $1\leq i\leq 2$,
entonces $m(a_1-a_2)=p(b_2-b_1)$. Puesto que $\mcd(m,p)=1$ y
$p|m(a_1-a_2)$, obtenemos que $p|(a_1-a_2)$. Adem\'as $2-p\leq
a_1-a_2\leq p-2$ por lo que $\frac{a_1-a_2}{p}\in
\big[-\big(1-\frac{2}{p}\big),\big(1-\frac{2}{p}\big)\big]\cap {\ma Z}
=\{0\}$, esto es, $a_1=a_2$ y por tanto $b_1=b_2$.
Se sigue que $\{\zeta_p^{\mu}\zeta_m^{\sigma}\}_{\mu,\sigma}=
\{\zeta_n^{\theta}\}_{\theta\in\Gal(\cic n{}/{\ma Q})}$ es una base
normal de $\cic n{}/{\ma Q}$. $\fin$

\end{proof}

\subsection{Subcampos cuadr\'aticos}\label{S4.10.3}

Sea $p$ un primo impar. Entonces $[{\ma Q}(\zeta_p):{\ma Q}]=p-1$
el cual es un n\'umero par y $\Gal(\cic p{}/{\ma Q})\cong U_{p}\cong
C_{p-1}$. Por tanto existe un \'unico subcampo cuadr\'atico $K\subseteq
\cic p{}$, es decir, $[K:{\ma Q}]=2$. Se tiene que el \'unico primo finito
ramificado en $\cic p{}$ es $p$, lo cual implica que el \'unico primo
finito ramificado en $K/{\ma Q}$ es $p$. Escribiendo $K={\ma Q}(\sqrt{d})$
con $d$ libre de cuadrados, necesariamente $4\nmid \delta_K$ pues
$2$ no es ramificado. Por esto, tenemos $\delta_K=d=\pm p$ y puesto
que $4\nmid \delta_K$, se sigue que $d\equiv 1\bmod 4$. Entonces
\[
d=\begin{cases}
p&\text{si $p\equiv 1\bmod 4$}\\-p&\text{si $p\equiv 3\bmod 4$}
\end{cases} =(-1)^{(p-1)/2}p.
\]
As{\'\i}, el subcampo cuadr\'atico de $\cic p{}$ es ${\ma Q}
\Big(\sqrt{(-1)^{\frac{p-1}{2}}p}\Big)$.

Otra forma de probar que ${\ma Q}
\Big(\sqrt{(-1)^{\frac{p-1}{2}}p}\Big)\subseteq \cic p{}$ es la siguiente.
Consideremos $\psi_p(x)=\frac{x^p-1}{x-1}=x^{p-1}+\cdots+x+1=
\prod_{i=1}^{p-1}(1-\zeta_p^i)$, por lo que $p=\psi_p(1)=\prod_{i=1}^{p-1}
(1-\zeta_p^i)$. Para $\frac{p+1}{2}\leq j\leq p-1$, esto es, $j=p-i$ con
$1\leq i\leq \frac{p-1}{2}$, se tiene
\begin{gather*}
1-\zeta_p^j=1-\zeta_p^{p-i}=1-\zeta_p^{-i}=
\zeta_p^{-i}(\zeta_p^i-1)=-\zeta_p^{-i}(1-\zeta_p^i),\\
\intertext{por lo que}
\begin{align*}
p&=\prod_{i=1}^{p-1}(1-\zeta_p^i)=\prod_{i=1}^{(p-1)/2}(1-\zeta_p^i)\cdot
\prod_{i=1}^{(p-1)/2}(1-\zeta_p^{-i})=\\
&=(-1)^{(p-1)/2}\prod_{i=1}^{(p-1)/2}\zeta_p^{-i}\cdot \prod_{i=1}^{(p-1)/2}(1-\zeta_p^i)^2=(-1)^{(p-1)/2}\zeta_p^s\alpha^2,\quad\alpha\in\cic p{},\\
s&=-\Big(1+2+\cdots+\frac{p-1}{2}\Big)=-\frac{\big(\frac{p-1}{2}\big)
\big(\frac{p+1}{2}\big)}{2}=-\frac{(p^2-1)}{8}.
\end{align*}
\end{gather*}
Si $s$ es par, $s=2n$, entonces $\zeta_p^{2n}=(\zeta_p^n)^2$ y
$p=(-1)^{(p-1)/2}\beta^2$ con $\beta\in\cic p{}$. Si $s$ es impar, se
tiene que $p+s=2m$ es par y $\zeta_p^s=\zeta_p^{p+s}=(\zeta_p^m)^2$.
En cualquier caso, existe $\gamma\in\cic p{}$ tal que $(-1)^{(p-1)/2}p=
\gamma^2$ y por lo tanto $\gamma =\sqrt{(-1)^{\frac{p-1}{2}}p}\in \cic p{}$.

Ahora bien, por unicidad del subcampo cuadr\'atico, se tiene que
\[
{\ma Q}\big(\sqrt{-(-1)^{\frac{p-1}{2}}p}\big)=
{\ma Q}\big(\sqrt{(-1)^{\frac{p+1}{2}}p}\big)\nsubseteq \cic p{}.
\]
Adem\'as
\begin{align*}
{\ma Q}(\sqrt{p},\sqrt{-p})&={\ma Q}\Big(\sqrt{(-1)^{\frac{p-1}{2}}p},
\sqrt{(-1)^{\frac{p+1}{2}}p}\Big)= {\ma Q}\Big(\zeta_4,
\sqrt{(-1)^{\frac{p-1}{2}}p}\Big)\\
&\subseteq {\ma Q}(\zeta_4,\zeta_p)=
\cic {4p}{}.
\end{align*}

Para $p=2$, se tiene que 
\begin{gather*}
{\ma Q}(\sqrt{2}),{\ma Q}(\sqrt{-2})\subseteq \cic 8{}, \quad
{\ma Q}(\sqrt{2},\sqrt{-2})\nsubseteq \cic 4{}\\
\intertext{y de hecho}
\cic 8{}={\ma Q}(\zeta_4,\sqrt{2})={\ma Q}(\zeta_4,\sqrt{-2})=
{\ma Q}(\sqrt{2},\sqrt{-2})={\ma Q}(\zeta_4,\sqrt{2},\sqrt{-2}).
\end{gather*}

En general, si $K={\ma Q}(\sqrt{d})$ con $d$ libre de cuadrados, digamos
$d=\pm p_1\cdots p_r$ con $p_1,\ldots, p_r$ primos distintos, entonces
\[
{\ma Q}(\sqrt{d})\subseteq {\ma Q}(\zeta_4, \sqrt{p_1},\ldots, \sqrt{p_r})
\subseteq {\ma Q}(\zeta_8,\zeta_{p_1},\ldots, \zeta_{p_r})
\subseteq \cic {8p_1\cdots p_r}{}
\]
lo cual es una manera expl{\'\i}cita, en este caso, del Teorema de 
Kronecker--Weber.

\begin{observacion}[No existencia de extensiones de $\cic 4{}={\ma Q}(i)$]
\label{O11.0}
Veamos que no existe una extensi\'on 
c\'iclica $L/{\ma Q}$ de grado $4$  tal que
${\ma Q}(i)\subseteq L$. De otra forma, sea $L/{\ma Q}$ una extensi\'on
c\'iclica de grado $4$ y sea $G=\Gal(L/{\ma Q})=\langle\sigma\rangle$ y
$K={\ma Q}(i)\subseteq L$. Entonces $L=K(\sqrt{\beta})$ con alg\'un
$\beta\in{\ma Q}(i)$. Sea $\beta=a+i b$ con $a,b\in{\ma Q}$ y sea
$\alpha=\sqrt{\beta}$. Entonces $\Gal({\ma Q}(i)/{\ma Q})=\langle
\sigma\bmod \sigma^2\rangle$.
\[
\xymatrix{
L\ar@{-}[d]_{\langle \sigma^2\rangle}
\ar@/^2pc/@{-}[dd]^{\langle\sigma\rangle}\\
K\ar@{-}[d]_{\langle\sigma\bmod
\sigma^2\rangle}\\{\ma Q}
}
\quad
\begin{array}{l}
\\
\\
\sigma^2|_K=\Id_K
\end{array}
\]
Se tiene que $
\sigma (\beta)=\sigma(\alpha^2)=(\sigma \alpha)^2=a-bi\Rightarrow
\sigma(\alpha)=\pm \sqrt{a-bi}$, por tanto $\alpha\sigma(\alpha)=
\pm\sqrt{a^2+b^2}$. Necesariamente se debe tener $\sqrt{a^2+b^2}
\in {\ma Q}$ pues de lo contrario ${\ma Q}(\sqrt{a^2+b^2})$ ser\'ia
una subextensi\'on cuadr\'atica real de $L$ y $L$ tendr\'ia dos
subextensiones cuadr\'aticas contradiciendo el hecho de que
$L/{\ma Q}$ es una extensi\'on c\'iclica. Sea $c=\sqrt{
a^2+b^2}$, $\sqrt{a+bi}\sqrt{a-bi}=\sqrt{a^2+b^2}=c$.

Ahora $\alpha=\sqrt{a+bi}$, $\alpha^2=a+bi$ por lo que
$(\alpha^2-a)^2=(bi)^2=-b^2$, $\alpha^4-2a\alpha^2+a^2+b^2=0$.
Por tanto $\alpha^4-2a\alpha^2+c^2=0$ de donde se sigue
que $p(x)=\Irr(x,\alpha,{\ma Q})=x^4-2ax^2+c^2$. La ra\'ices 
de $p(x)$ son $\pm \alpha$, $\pm \gamma$ donde $\gamma=
\sqrt{a-bi}$. Se tiene $\alpha\gamma=c$ por lo que $\gamma=
\frac{c}{\alpha}$.

Se tiene que $o(\sigma)=4$, por lo que necesariamente
$\sigma(\alpha)=\pm\gamma$. Entonces
\begin{gather*}
a+b\sigma(i)=\sigma(a+bi)=\sigma(\alpha^2)=(\sigma(\alpha))^2
=\gamma^2=a-bi,
\intertext{por lo que $\sigma(i)=-i$. Se sigue que}
\sigma^2(\alpha)=\sigma(\sigma(\alpha))=\sigma(\pm \gamma)=
\pm \sigma(\gamma)=\pm \sigma\Big(\frac{c}{\alpha}\Big)=
\pm \frac{\sigma(c)}{\sigma(\alpha)}=\pm \frac{c}{\pm\gamma}=
\frac{c}{\gamma}=\alpha,
\end{gather*}
lo cual implica que $\sigma^2(\alpha)=\alpha$ por lo que
$o(\sigma)=2$. Esto contradice que $L/{\ma Q}$ es c\'iclica de
grado $4$.

Usando el Teorema de Kronecker--Weber, podemos dar
otro tipo de demostraci\'on de este hecho. Supongamos nuevamente
$L/{\ma Q}$ es una extensi\'on c\'iclica de grado $4$ de
${\ma Q}$ y tal que ${\ma Q}(i)\subseteq K$. Por el Teorema de
Kronecker--Weber ${\ma Q}(i)=\cic 4{}\subseteq L \subseteq \cic m{}$
para alg\'un $m$ con $4|m$. Nuevamente pongamos $\Gal(L/
{\ma Q})=\langle \sigma\rangle$. Entonces $\Gal(\cic 4{}/{\ma Q})
=\langle\sigma^2\rangle =\langle J\rangle$. 

Veamos que no existe $\sigma\in\Gal(\cic m{}/{\ma Q})$ tal que
$\sigma^2=J$ donde $J$ es conjugaci\'on compleja.
Supongamos que $\sigma^2=J$. Entonces $\sigma^2(i)=
\sigma^2(\zeta_4)=J(\zeta_4)=\bar{\zeta}_4=\zeta_4^3=-i$.
Sea $\sigma(\zeta_m)=
\zeta_m^{a_{\sigma}}$. Entonces
\[
\sigma(\zeta_4)=\sigma(\zeta_m^{m/4})=\zeta_m^{a_{\sigma}
m/4}=\zeta_4^{a_{\sigma}},
\]
por lo que $\sigma^2(\zeta_4)=\zeta_4^{a_{\sigma}^2}=\zeta_4^3$.
Por tanto $a_{\sigma}^2\equiv 3\bmod 4$. Sin embargo para 
cualquier $a\in{\ma Z}$ se tiene $a^2\equiv 0,1\bmod 4$. Esto
demuestra que no existe tal $L$.
\end{observacion}

La Observaci\'on \ref{O11.0} es un caso particular del siguiente
resultado.

\begin{teorema}[\cite{IsLuFa97}]\label{T11.-1}
Sean $k$ un campo de caracter\'istica diferente de $2$ y $K=
k(\sqrt{d})$ una extensi\'on cuadr\'atica. Entonces existe $L/k$ 
una extensi\'on c\'iclica de grado $4$ sobre $k$ con $k\subseteq
K\subseteq L\iff$ existe $\delta\in K$ tal que $\N_{K/k}\delta =-1$
($\iff d$ es representable como la suma de dos cuadrados en $k$).
\end{teorema}

\begin{proof}
Sea $\N_{K/k}\delta=-1$, digamos $\delta =a+b\sqrt{d}$ con 
$a,b\in k$. Sea $\Gal(K/k)=\{1,\sigma\}$. Entonces $\delta^{\sigma}
=a-b\sqrt{d}$ y $\N_{K/k}\delta=\delta\delta^{\sigma}=a^2-b^2d=-1$.
Por tanto $d=\big(\frac{a}{b}\big)^2+\big(\frac{1}{b}\big)^2$ y $d$
es representable como la suma de cuadrados. Rec\'iprocamente,
si $d=u^2+v^2$ entonces $u,v\neq 0$ puesto que $\sqrt{d}\notin k$
y si $t=\frac{v}{u}+\frac{\sqrt{t}}{u}$, $\N_{K/k} t=-1$.

Primero supongamos que existe $\delta$ tal que $\N_{K/k}\delta=
-1$. Sea $\lambda=1+\delta^2$. En caso de que $\lambda=0$ 
entonces $\delta^2=-1$ y $\delta=\zeta_4\in k$ y $L=k(\sqrt[4]{d})$
satisface lo requerido. Supongamos $\lambda\neq 0$. Entonces
\[
\lambda^{\sigma}=1+(\delta^{\sigma})^2=1+\Big(-\frac{1}{\delta}
\Big)^2=\frac{\delta^2+1}{\delta^2}=\frac{\lambda}{\delta^2}.
\]

Sean $\theta=\sqrt{\lambda}$, $\theta^2=\lambda$
y $L=K(\theta)$. Sea $\psi$ una extensi\'on
de $\sigma$ a $L$. Se tiene $(\theta^2)^{\psi}=(\theta^\psi)^2=\lambda^{
\sigma}=\lambda/\delta^2=(\theta/\delta)^2$. En particular
$\theta^{\psi}=\pm (\theta/\delta)\in L$. Por tanto $\psi$
es un automorfismo de $L$ sobre $k$, es decir, $\psi\in\Aut_k(L)$
y se tiene $\theta^{\psi^2}=-\theta$, $\theta^{\psi^3}=\mp \theta/\delta$
y $\psi^4=1$. Esto demuestra que $L/k$ es una extensi\'on normal
y adem\'as que $\Gal(L/k)\cong C_4$. Por otro lado, tenemos que
$K=k(\delta)$ y que $\delta=\pm \theta^{\psi}/\theta\in L$. En
particular $K\subseteq L$.

Rec\'iprocamente, supongamos que $k\subseteq K\subseteq L$ con
$K=k(\sqrt{d})$ y $L/k$ una extensi\'on c\'iclica de grado $4$. 
Sea $L=k(\theta)$ y $\Gal(L/k)=\langle\psi\rangle$. Sea $\delta:=
\frac{\theta-\theta^{\psi^2}}{\theta^{\psi}-\theta^{\psi^3}}\neq 0$ y
$\delta^{\psi}=-1/\delta$. Por tanto $\delta^{\psi^2}=\delta$ y en particular
$\delta\in K$. Finalmente, si $\sigma=\psi^2$ se tiene $\N_{K/k}\delta
=\delta\delta^{\sigma}=-1$.
$\fin$
\end{proof}

La Observaci\'on \ref{O11.0} se sigue del Teorema \ref{T11.-1} pues
para todo $a+bi\in{\ma Q}(i)$, se tiene $\N_{{\ma Q}(i)/{\ma Q}}(
a+bi)=a^2+b^2\neq -1$. Por otro lado el Teorema \ref{T11.-1} 
nos muestra que el problema del encaje en el problema inverso
de la Teor\'ia de Galois es altamente no trivial.

\begin{ejemplo}\label{Ej11.-2}
Sea $K={\ma Q}(\sqrt{65})$. Entonces $\delta=\sqrt{65}-8$
satisface que $\N_{K/{\ma Q}}(\delta)=-65+64=-1$, por lo que
existe $L/{\ma Q}$ c\'iclica de grado $4$ conteniendo $K$.
Un tal $L$ es $L:={\ma Q}(\sqrt{65+\sqrt{65}})$.
\end{ejemplo}

\section{Anillos de enteros}\label{S4.11}

Para $n\in{\ma N}$ se tiene que ${\ma Z}[\zeta_n]={\cal O}_{\cic n{}}$.
Una pregunta natural es si esto se cumple para todos los subcampos
que hemos estudiado, es decir, si $M_b T_r={\ma Q}(\mu_b,\delta_r)$, 
entonces, ?`es el anillo de enteros de $M_b T_r$ de la forma ${\ma Z}
[\alpha]$ para alg\'un $\alpha$?

\begin{teorema}\label{T11.1} El anillo de enteros de $
\cic n{}^+={\ma Q}(\zeta_n+
\zeta_n^{-1})$ es ${\ma Z}[\zeta_n+\zeta_n^{-1}]$.
\end{teorema}

\begin{proof}
Sea $\alpha =a_0+a_1(\zeta_n+\zeta_n^{-1})+\cdots+
a_m(\zeta_n+\zeta_n^{-1})^m$ con $a_i\in{\ma Q}$ un entero de
${\ma Q}(\zeta_n+\zeta_n^{-1})$. Puesto que $[\cic n{}^+:{\ma Q}]=
\varphi(n)/2$, se tiene que $m\leq (\varphi(n)/2)-1$. Queremos probar
que $a_i\in{\ma Z}$ con lo cual obtendremos que ${\cal O}_{\cic n{}^+}
\subseteq {\ma Z}[\zeta_n+\zeta_n^{-1}]$. Puesto que la otra contenci\'on
es inmediata, se seguir\'a la igualdad y el teorema.

Ahora bien, restando a ambos miembros de la igualdad anterior los
t\'erminos tales que $a_i\in{\ma Z}$, en caso de ser falsa nuestra
afirmaci\'on, podemos suponer que $a_m\notin{\ma Z}$.

Multiplicando a ambos miembros por $\zeta_n^m$, obtenemos
\[
\zeta_n^m\alpha=a_m+b_1 \zeta_n+\cdots+b_{2m-1}\zeta_n^{2m-1}
+a_m\zeta_n^{2m}
\]
el cual es un elemento de $\cic n{}$ y adem\'as es un entero algebraico,
es decir, $\zeta_n^m\alpha\in{\cal O}_{\cic n{}}={\ma Z}[\zeta_n]$. Adem\'as
se tiene que $2m\leq 2\big(\frac{\varphi(n)}{2}\big)-2 < \varphi(n)-1$ y puesto que
$\{1,\zeta_n,\cdots, \zeta_n^{\varphi(n)-1}\}$ es una base entera de
${\cal O}_{\cic n{}}$, se sigue que $a_m\in{\ma Z}$, probando
el resultado. $\fin$

\end{proof}

Para los campos $\cic 2m$, $m\in{\ma N}$ se sigue cumpliendo el
Teorema \ref{T11.1}. Es decir, tenemos

\begin{teorema}\label{T11.1'}
Para todo $m$ y todo subcampo $E\subseteq \cic 2m$ existe $\alpha
\in \AE E$ tal que $\AE E={\ma Z}[\alpha]$. M\'as a\'un, el anillo de enteros
de $E={\ma Q}(\xi)$ es ${\ma Z}[\xi]$ donde 
$\xi\in\{\zeta_{2^r},\zeta_{2^r}+\zeta_{2^r}^{-1},\menoszeta r\}$.
\end{teorema}

\begin{proof}
Del Teorema \ref{T10.1.1} tenemos que 
los subcampos de $\cic 2m$ son:
$\cic 2r,\cic 2r^+$ y ${\ma Q}(\menoszeta r)$, $0\leq r\leq m$.
El caso $\cic 2r$ es el Teorema \ref{T1.2.1.11}, el caso $\cic 2r^+$
es el Teorema \ref{T11.1}, restando el caso ${\ma Q}(\menoszeta r)$.

Sea $E={\ma Q}(\menoszeta r)$. Se tiene
\begin{gather*}
\xymatrix{
\cic 2r^+\ar@{-}[dd]\ar@{-}[rr]&&\cic 2r\ar@{-}[dd]\\
&E\ar@{-}[dl]\\ \cic 2{r-1}^+\ar@{-}[rr]&&\cic 2{r-1}
}\\
[E:{\ma Q}]=[\cic 2r^+:{\ma Q}]=\frac{[\cic 2r:{\ma Q}]}{2}=
\frac{\varphi(2^r)}{2}=2^{r-2}.
\end{gather*}
Entonces $\{1,(\menoszeta r),\ldots,(\menoszeta r)^{2^{r-2}-1}\}$ es
base de $E/{\ma Q}$.

Ahora bien, como $\zeta_{2^r}$ es un entero, $\menoszeta r=
\zeta_{2^r}-\zeta_{2^r}^{2^r-1}$ es entero y por tanto ${\ma Z}[
\menoszeta r]\subseteq \AE E$. Rec\'iprocamente, sea $\alpha\in
\AE E$, 
\[
\alpha=a_0+a_1(\menoszeta r)+\cdots+a_m(\menoszeta
r)^m, \quad m\leq 2^{r-2}-1, \quad a_i\in{\ma Q}.
\]
 En caso de que tuvi\'esemos
$\AE E\neq {\ma Z}[\menoszeta r]$, existir\'ia $\alpha\in \AE E$ con
alg\'un $a_j\notin {\ma Z}$. Pasando al lado izquierdo de la igualdad
los elementos con $a_j\in{\ma Z}$, podemos suponer que
$a_m\notin {\ma Z}$, $m\leq 2^{r-2}-1$.

Multiplicando por $\zeta_{2^r}^m$ a ambos lados de la ecuaci\'on
obtenemos
\[
\zeta_{2^r}^m\alpha=(-1)^m a_m+b_1\zeta_{2^r}+\cdots+b_{2m-1}
\zeta_{2^r}^{2m-1}+a_m\zeta_{2^r}^{2m}.
\]
Puesto que $\zeta_{2^r}^m\alpha$ es entero, $\zeta_{2^r}^m\alpha
\in \cic 2r$ y el anillo de enteros de $\cic 2r$ es ${\ma Z}[\zeta_{2^r}]$
se sigue que $a_m\in{\ma Z}$ puesto que $2m\leq 2^{r-1}-2<\varphi(
2^r)-1$ y $\{1,\zeta_{2^r},\ldots,\zeta_{2^r}^{\varphi(2^r)-1}\}$ es
base de $\cic 2r/{\ma Q}$. Esta contradicci\'on prueba que el anillo
de enteros de ${\ma Q}(\menoszeta r)$ es ${\ma Z}[\menoszeta r]$.
$\fin$
\end{proof}

El Teorema \ref{T11.1} prueba que, en la notaci\'on del Teorema
\ref{T10.2.1}, si $c=2$, $b=(p-1)/2$, $r=0$, entonces siempre se tiene
que el anillo de enteros de ${\ma Q}(\mu_{(p-1)/2})$ es ${\ma Z}[\mu_{
(p-1)/2}]$. Sin embargo, para $b<(p-1)/2$, inclusive con $r=0$, el 
resultado ya no sigue siendo cierto.

\begin{ejemplo}\label{Ej11.2}
Consideremos dos n\'umeros primos distintos, $p,q$ con $q>2$ tales que
si $f=o(p\bmod q)$, entonces $(q-1)/f>p$. Sea $E\subseteq \cic q{}$ el
subcampo de grado $\frac{q-1}{f}$ sobre ${\ma Q}$. Entonces
$p$ se descompone totalmente en la extensi\'on $E/{\ma Q}$ ya que
$p$ es no ramificado en $\cic q{}$ y el grupo de descomposici\'on de
$p$ es $\Gal(\cic q{}/E)$. Si el anillo de enteros ${\cal O}_E$ de $E$
fuese de la forma ${\ma Z}[\alpha]$ para alg\'un $\alpha$, entonces
si $f(x):=\Irr(\alpha, x,{\ma Q})\in {\ma Z}[x]$, tendr{\'\i}amos por el Teorema de
Kummer, Teorema \ref{T8.12}, que $f(x)$ se descompone en
en $\frac{q-1}{f}$ factores lineales m\'onicos distintos m\'odulo $p$.
Sin embargo, esto no es posible pues \'unicamente existen $p$ factores
lineales m\'onicos distintos m\'odulo $p$, a saber, $x, x+1,\ldots, x+
(p-1)$ y por hip\'otesis tenemos $\frac{q-1}{f}>p$, $f<\frac{q-1}{p}$.

Notemos que en este ejemplo, $b=\frac{q-1}{f}$, $c=f$.

Observamos que debemos tener necesariamente que $c=f>2$, y
$\frac{q-1}{f}>p$, esto es, $2<f<\frac{q-1}{p}$. Un ejemplo en que se
cumple lo anterior es $p=2$, $q=31$. Se tiene que $2^5=32\equiv 1
\bmod 31$, $o(2\bmod 31)=5$, es decir, $f=5$ y $\frac{q-1}{f}=\frac{31-1}{5}
=6>p=2$.

Otro ejemplo es $p=3$, $q=1093$. Entonces $3^7=2187\equiv
1\bmod 1093$ y $o(3\bmod 1093)=7=f$, $\frac{1093-1}{f}=\frac{1092}{7}
=156>p=3$. En otras palabras el anillo de enteros
del campo de grado $156$ sobre ${\ma Q}$
contenido en $\cic {1093}{}$ no es de la forma ${\ma Z}[\alpha]$
para ning\'un $\alpha$.
\end{ejemplo}

\section{Teorema de reciprocidad cuadr\'atica\index{teorema
de reciprocidad cuadr\'atica}\index{teorema!$\sim$ reciprocidad
cuadr\'atica}}\label{S4.5}

Recordemos el {\em s{\'\i}mbolo de Legendre\index{simbolo de Legendre@s{\'\i}mbolo
de Legrendre}\index{Legendre!s\'imbolo de $\sim$}}. 
Para dos n\'umeros $n, p\in {\ma Z}$, con
$p$ primo y $p\nmid n$, sea
 $\d\left(\frac{n}{p}\right)$ el s{\'\i}mbolo
 de Legendre, el cual est\'a definido por:
 \begin{align*}
 \d\left(\frac{n}{p}\right) &= \left\{
 \begin{array}{rl}
 1&\text{si $n$ es un residuo cuadr\'atico $\bmod\ p$}\\
 -1&\text{si $n$ es un residuo no cuadr\'atico $\bmod\ p$}
\end{array}
\right. =\\
&= \left\{
\begin{array}{rl}
1&\text{si $\exists\ x\in\ma Z$ tal que $n\equiv
x^2\bmod p$}\\
-1 &\text{si $\not\exists\ x\in \ma Z$ tal que $n\equiv x^2 \bmod p$}
\end{array}
\right..
\end{align*}

Notemos que 
\[
\left(\frac{n}{p}\right)\equiv n^{\frac{p-1}{2}}\bmod p
\]
pues si $\langle \xi\rangle= {\ma F}_p^{\ast}$, entonces existe $a\in{\ma Z}$
tal que $n\bmod p=\xi^a$ y $(n\bmod p)^{\frac{p-1}{2}}=
\xi^{\frac{a}{2}(p-1)}$ y se tiene que $(p-1)|\frac{a}{2}(p-1) \iff
a$ es par $\iff n$ es un cuadrado m\'odulo $p$. Por otro lado
$(p-1)|\frac{a}{2}(p-1) \iff \xi^{\frac{a}{2}(p-1)}=1$.

\begin{teorema}[Teorema de reciprocidad cuadr\'atica]\label{T4.5.1}
Sean $p$ y $q$ dos primos racionales impares distintos. Entonces
\[
\d\left(\frac{p}{q}\right)\cdot 
\d\left(\frac{q}{p}\right) = (-1)^{
\frac{p-1}{2}\frac{q-1}{2}}.
\]
\end{teorema}

\begin{proof} 
Primero, sean $p\equiv 1\bmod 4$ y $q$ arbitrario. Entonces
\begin{gather*}
\d\left(\frac{p}{q}\right)=1\iff \text{existe $x\in{\ma Z}$
tal que $x^2\equiv p \bmod q$}.\\
\intertext{Escribamos $p=1+4m$ para alg\'un
$m\in{\ma Z}$. Entonces}
x^2\equiv p\bmod q\iff x^2-1\equiv 4m\bmod q.\\
\intertext{Ahora bien, puesto que $q$ es impar,
$2$ es invertible m\'odulo $q$. Por tanto}
x^2-1\equiv 4m \bmod q \iff \frac{x+1}{2}\cdot \frac{x-1}{2}\equiv m \bmod q.\\
\intertext{Sea $y=\frac{x+1}{2}$. Por tanto $y-1=\frac{x-1}{2}$. Entonces}
\frac{x+1}{2}\cdot \frac{x-1}{2}\equiv m \bmod q\iff y(y-1)\equiv m \bmod q
\end{gather*}
lo cual es equivalente a que $y^2-y-m=y^2-y+\frac{1-p}{4}\equiv 0\bmod q$.

Ahora sea $\alpha=\frac{1+\sqrt{p}}{2}$. Entonces ${\mathcal O}_{{\ma Q}
(\sqrt{p})}={\ma Z}[\alpha]$. Tenemos $(2\alpha-1)^2=p$. Por tanto
$4\alpha^2-4\alpha+1=p$ lo cual equivale a $\alpha^2-\alpha+\frac{1-p}{4}=0$.
Es decir
\[
f(x):=\Irr(\alpha,x,{\ma Q})=x^2-x+\frac{1-p}{4}.
\]
Entonces $q$ se descompone en ${\ma Q}(\sqrt{p})$ si y solamente si,
por el Teorema de Kummer, 
$f(x)\bmod q$ se factoriza como el producto de dos factores lineales
distintos. Esto \'ultimo es equivalente a que $f(x)$ tenga ra{\'\i}ces
m\'odulo $q$, esto es, si y solamente si $y^2-y+\frac{1-p}{4}\equiv 0\bmod q$
tiene soluci\'on.

En resumen, tenemos que $\d\left(\frac{p}{q}\right)=1$ si y solamente si
$q$ se decompone en ${\ma Q}(\sqrt{p})$.

Por otro lado tenemos que $q$ se descompone en ${\ma Q}(\sqrt{p})$
si y solamente si ${\ma Q}(\sqrt{p})$ est\'a contenido en el campo de
descomposici\'on de $q$ en ${\ma Q}(\zeta_p)$ y puesto que este \'ultimo
es el campo fijo bajo el automorfismo de Frobenius correspondiente a $q$,
es decir a ${\ma Q}(\zeta_p)^{\langle\sigma_q\rangle}$, se tiene
\[
q \text{\ se descompone en\ } {\ma Q}(\sqrt{p})\iff {\ma Q}(\sqrt{p})\subseteq
{\ma Q}(\zeta_p)^{\langle\sigma_q\rangle}\iff \sigma_q \text{\ fija a\ }
{\ma Q}(\sqrt{p}).
\]

Ahora bien, $\sigma_q$ est\'a definido por $\sigma_q(\zeta_p)=\zeta_p^q$
y como la extensi\'on ${\ma Q}(\zeta_p)/{\ma Q}$ es una extensi\'on c{\'\i}clica,
se tiene la existencia de un \'unico campo de grado $h$ sobre ${\ma Q}$
para cada divisor $h$ de $p-1=[{\ma Q}(\zeta_p):{\ma Q}]$. Si definimos
$f:=o(\sigma_q)=o(q\bmod p)$, se tiene que
\[
 {\ma Q}(\sqrt{p})\subseteq
{\ma Q}(\zeta_p)^{\langle\sigma_q\rangle} \iff f\mid \frac{\varphi(p)}{2}=
\frac{p-1}{2}\iff q^{\frac{p-1}{2}}\equiv 1\bmod p.
\]

Sea $a$ un generador m\'odulo $p$, esto es, $o(a\bmod p)=p-1$. Sea 
$a^t\equiv q\bmod p$. Entonces $a^{t\frac{p-1}{2}}\equiv q^{\frac{p-1}{2}}
\equiv 1\bmod p$ lo cual equivale a que $p-1\mid t\frac{p-1}{2}$, esto es,
$t=2s$ es par y por tanto si $b=a^s$ se tiene $b^2=a^{2s}=a^t\equiv q
\bmod p$, es decir la ecuaci\'on $z^2\equiv q\bmod p$ es soluble. Esto
\'ultimo equivale a que $\d\left(\frac{q}{p}\right)=1$.

En resumen, hemos probado que si $p\equiv 1\bmod 4$, entonces para 
cualquier primo impar $q$, 
\[
\d\left(\frac{p}{q}\right)=1\iff \d\left(\frac{q}{p}\right)=1.
\]
Por tanto $\d\left(\frac{p}{q}\right)=\d\left(\frac{q}{p}\right)$ y el
resultado se sigue en este caso. Por simetr{\'\i}a lo mismo tenemos para
$q\equiv 1\bmod 4$ y $p$ cualquier primo impar.

Ahora sean $p$ y $q$ dos primos diferentes congruentes con $3$ m\'odulo $4$.
Primero notemos que $x^2\equiv -1 \bmod q$ no es soluble pues si lo fuese
$x^4\equiv 1\bmod q$ y por tanto $o(x\bmod q)=4$ por lo que $4\mid q-1=
o(U_q)$, lo que contradice que $q\equiv 3\bmod 4$.

Por tanto tenemos que $x^2\equiv a\bmod q$ es soluble si y solamente si
$y^2\equiv -a\bmod q$ no es soluble pues si ambas lo fuesen se tendr\'ia
que $(x/y)^2\equiv (a/-a)\equiv -1 \bmod q$.

De esta forma obtenemos que $\d\left(\frac{-p}{q}\right)=\d-\left(\frac{p}{q}\right)$.
Ahora bien, en general tenemos $\d\left(\frac{p}{-q}\right)=\d\left(\frac{p}{q}\right)$
pues $p\equiv x^2\bmod q\iff q|p-x^2\iff -q|p-x^2\iff p\equiv x^2\bmod (-q)$.
Se tiene $-p\equiv 1\bmod 4$ de donde, por la primera parte, obtenemos
\[
-\d\left(\frac{p}{q}\right)=\d\left(\frac{-p}{q}\right)=\d\left(\frac{q}{-p}\right)=
\d\left(\frac{q}{p}\right),
\]
es decir $\d\left(\frac{p}{q}\right)\cdot \d\left(\frac{q}{p}\right)=-1=
(-1)^{\frac{p-1}{2}\frac{q-1}{2}}$. $\fin$
\end{proof}

\begin{proposicion}\label{P4.5.2}
Sea $q$ una potencia de un n\'umero primo.
Se tiene que $\sqrt{-1}\in \F\iff q$ es potencia de $2$ o $q
\equiv 1\bmod 4$.
\end{proposicion}

\begin{proof}
Si $q$ es una potencia de $2$, entonces
$\sqrt{-1}=-1=1=1^2$, es decir $-1$ es un cuadrado
en $\F$.

Ahora consideremos $q$ impar. 

Una demostraci\'on inmediata es: $\*\F$ es el grupo de las
$(q-1)$ ra\'ices de la unidad y $\sqrt{-1}$ es una ra\'iz cuarta
primitiva de la unidad. Por lo tanto $\sqrt{-1}\in \F$ si y solo
si $4|q-1$.

Damos otra demostraci\'on, mucho m\'as larga pero instructiva.
Primero
consideremos el caso $q\equiv 3\bmod 4$. Veamos que
$\sqrt{-1}\notin \F$. Para este fin, tenemos,
por el Teorema \ref{T8.2}, que $p$ es inerte
en ${\ma Q}(\sqrt{-1})=
{\ma Q}(i)={\ma Q}(\zeta_4)$, donde $q=p^u$ con
$p\equiv 3\bmod 4$ y $u$ es impar pues
$o(p\bmod 4)=2$. Entonces 
$[{\ma F}_p(\sqrt{-1}):{\ma F}_p]=2$. Por otro lado
$[\F:{\ma F}_p]=u$ y $\mcd(2,u)=1$. Se sigue que
${\ma F}_p(\sqrt{-1})\nsubseteq \F$ y por tanto
$\sqrt{-1}\notin \F$.

Ahora veamos que si $q\equiv 1
\bmod 4$, entonces $\sqrt{-1}\in \F$. 
Sea $q=p^u$. Hay dos casos. Si $p\equiv
1\bmod 4$, entonces $u$ es arbitrario.
Se tiene, por el Teorema \ref{T8.2} que
$p$ se descompone en ${\ma Q}(i)/{\ma Q}$ por
lo que $[{\ma F}_p(\sqrt{-1}):{\ma F}_p]=1$
y $\sqrt{-1}\in{\ma F}_p\subseteq \F$.

Ahora, si $p\equiv 3\bmod 4$, $u$ es par.
El primo $p$ es inerte en ${\ma Q}(i)/{\ma Q}$ y
$[{\ma F}_p(\sqrt{-1}):{\ma F}_p]=2$ por lo que
${\ma F}_p(\sqrt{-1})={\ma F}_{p^2} \subseteq
{\ma F}_{p^u}=\F$. $\fin$
\end{proof}

%% file: Capitulo6.tex
\chapter{Caracteres de Dirichlet}\label{Ch12}

\section{Teor{\'\i}a de caracteres}\label{S12.1}

Primeramente recordemos los resultados b\'asicos de la teor{\'\i}a de
caracteres.

\begin{definicion}\label{D12.1}
Sea $G$ un grupo cualquiera. El {\em grupo de caracteres\index{grupo
de caracteres}} de $G$ se define como el grupo de homomorfismos de
$G$ en el grupo multiplicativo de los complejos: 
\[
\Hom(G,{\ma C}^{\ast})=\{f\colon G\to {\ma C}^{\ast}\mid
f \text{\ es homomorfismo de grupos}\}.
\]
\end{definicion}

Cuando $G$ es finito, si $n=|G|$, entonces si $f$ es un caracter y
$g\in G$, entonces $1=f(e)=f(g^n)=f(g)^n$, es decir, $f(g)\in W_n$
y podemos definir $\Hom(G,W_n)$ como el grupo de caracteres de $G$.

Si $f$ es un caracter y $G$ es un grupo topol\'ogico tal que para todo
$g\in G$, $f(g)$ es de orden finito, es decir, para cada $g\in G$, existe
$n_g\in {\ma N}$ tal que $f(g)^{n_g}=1$, entonces $f(g)\in W_{n_g}$.
Sea $R:={\ma Q}/{\ma Z}=\cup_{n=1}^{\infty}W_n\subseteq {\ma C}^{\ast}$
con la topolog{\'\i}a discreta. Entonces definimos
\[
\hat{G}:=\Hom(G,R)\subseteq \Hom(G,{\ma C}^{\ast})
\]
donde $\Hom(G,R):=\{f\colon G\to R\mid f \text{\ es continuo}\}$.

En general estaremos interesados en $\hat{G}$ y en general si $G$ es
un grupo finito, $\hat{G}=\Hom(G,{\ma C}^{\ast})$. En este caso $G$
tiene la topolog{\'\i}a discreta y todos los homomorfismos son
autom\'aticamente continuos.

\begin{ejemplo}\label{Ej12.1}
Sea $f\in\hat{{\ma Z}}=\Hom({\ma Z},R)$. Entonces $f$ est\'a completamente determinado
por $f(1)\in R$. Sea
\begin{eqnarray*}
\varphi\colon\hat{{\ma Z}}&\longto&R\\ f&\longmapsto&f(1).
\end{eqnarray*}
Entonces $\varphi$ es un isomorfismo de grupos y $\hat{{\ma Z}}\cong R$.
Observemos que en este ejemplo, $\hat{{\ma Z}}$ no denota al
anillo de Pr\"ufer (Ejemplo \ref{Ej5.10'}).
\end{ejemplo}

\begin{ejemplo}\label{Ej12.2}
Sea ${\ma Z}_p$ el anillo de los enteros $p$--\'adicos (Ejemplos
\ref{Ej5.10'} (3)) y sea ${\ma Q}_p$ el campo de los n\'umeros $p$--\'adicos,
${\ma Q}_p=\coc {\ma Z}_p$. Sea $f\in\hat{{\ma Z}}_p$. Si $g\in
{\ma Z}_p$, $f(ng)=nf(g)=0\in R$ para alguna $n\in{\ma N}$.
Ahora, si $n=p^ms$ con $\mcd(s,p)=1$, $s$ es invertible en ${\ma Z}_p$
y existe $t\in{\ma Z}_p$ tal que $st=1$. Como ${\ma Z}_p$ es la
cerradura de ${\ma Z}$ en la topolog{\'\i}a $p$--\'adica, existe
una sucesi\'on $\{t_i\}_{i=1}^{\infty}\subseteq {\ma Z}$ tal que
$\lim_{i\to\infty}t_i=t$ en la topolog{\'\i}a $p$--\'adica. Puesto que 
$f$ es continuo, 
\begin{align*}
f(p^mg)&=f(p^mst g)=\lim_{i\to\infty}f(st_ip^mg)=\\
&=\lim_{i\to\infty}t_if(sp^mg)=\lim_{i\to\infty}t_if(ng)=\lim_{i\to\infty}
0=0.
\end{align*}

Por tanto, $f(p^mg)=p^mf(g)=0$, es decir, $f\big(\hat{{\ma Z}}_p
\big)=R(p)=({\ma Q}/{\ma Z})(p)\cong {\ma Q}_p/{\ma Z}_p=:
R_p$.

Ahora bien,  por continuidad,
$f$ est\'a totalmente determinado por $f(1)$ pues si $x\in
{\ma Z}_p$, existe una sucesi\'on $\{x_n\}_{n=1}^{\infty}
\subseteq {\ma Z}$ tal que $x_n \xrightarrow[n\to\infty]{}
x$ en la topolog{\'\i}a $p$--\'adica y $f(x)=\lim_{n\to\infty}
(f(x_n))=\lim_{n\to\infty} x_nf(1)=xf(1)$.

Por tanto
\begin{eqnarray*}
\hat{{\ma Z}}_p&\longto& R_p={\ma Q}_p/{\ma Z}_p\\
f&\longmapsto&f(1)
\end{eqnarray*}
es un isomorfismo de grupos y $\hat{{\ma Z}}_p\cong R_p$.
\end{ejemplo}

\begin{proposicion}\label{P12.3}
Si $G$ es un grupo c{\'\i}clico finito, se tiene $\hat{G}\cong G$.
\end{proposicion}

\begin{proof}
Sea $G\cong C_m$
para alguna $m$ y sea $a$ un generador de $G$. Entonces 
cualquier $f\in\hat{G}$ est\'a completamente determinado por
$f(a)$ y $f(a)^m=f(a^m)=f(e)=1\in{\ma C}^{\ast}$, es decir,
$f(a)\in W_m$. Por lo tanto, $\hat{G}=\Hom(G,W_m)$. Finalmente,
tenemos que $\varphi\colon\hat{G}\to W_m$, $a\mapsto f(a)$
es un isomorfismo de grupos. $\fin$
\end{proof}

\begin{observacion}\label{O12.4} El isomorfismo de la Proposicion
\ref{P12.3} no es can\'onico pues depende del generador $a$
seleccionado.
\end{observacion}

\begin{teorema}\label{T12.5} Si $G$ es un grupo finito, entonces
$\hat{G}\cong G/G'$ donde $G'$ denota al subgrupo conmutador.
En particular, si $G$ es un grupo abeliano finito, entonces $\hat{G}
\cong G$.
\end{teorema}

\begin{proof}
Si $f\in\hat{G}$, puesto que, ya sea $R$ o ${\ma C}^{\ast}$ son
abelianos, se tiene que $G'\subseteq \ker f$ y por lo tanto $f$
se factoriza de manera \'unica 
\[
\xymatrix{
G\ar[rr]^f\ar[dr]_{\pi}&&R\\ & G/G'\ar[ur]_{\tilde{f}}}
\]
donde
 $\pi\colon G\to G/G'$ es la proyecci\'on natural y $f=\tilde{f}\circ \pi$.
 
Puesto que $G/G'$ es un grupo abeliano finito y se tiene $\hat{G}=
\widehat{G/G'}$, para terminar la demostraci\'on, basta probar que
$G\cong \hat{G}$ para un grupo abeliano finito.

Sea $G\cong C_{n_1}\times \cdots \times C_{n_r}$ un grupo
abeliano finito cualquiera. En general, si $G=H_1\times H_2$, la
funci\'on
\begin{eqnarray*}
\hat{G}&\longto&\hat{H}_1\times \hat{H}_2\\
f&\longmapsto& \big(f|_{h_1}, f|_{H_2}\big)
\end{eqnarray*}
es un isomorfismo de grupos. Por lo tanto, por la Proposici\'on \ref{P12.3}
\begin{gather*}
\hat{G}\cong \hat{C}_{n_1}\times\cdots\times \hat{C}_{n_r}
\cong C_{n_1}\times \cdots\times C_{n_r}\cong G. \tag*{$\fin$}
\end{gather*}
\end{proof}

\begin{observacion}\label{O12.6}
El Teorema \ref{T12.5} no es v\'alido para grupos abelianos 
infinitos, por ejemplo, $\hat{{\ma Z}}\cong R\not\cong {\ma Z}$
(ver Ejemplo \ref{Ej12.1}).
\end{observacion}

\begin{definicion}\label{D12.7}
Sean $G$, $H$ grupos topol\'ogicos cualesquiera. Un 
{\em apareamiento\index{apareamiento}} es una funci\'on bilineal
continua $\varphi\colon G\times H\to {\ma C}^{\ast}$, es decir,
\begin{gather*}
\varphi(xy,z)=\varphi(x,z)\varphi(y,z),\\
\varphi(u,vw)=\varphi(u,v)\varphi(u,w),
\end{gather*}
para cualesquiera $x,y,u\in G$ y $z,v,w\in H$.

Un apareamiento se llama {\em no degenerado\index{apareamiento!no
degenerado}} si
\[
\varphi(x,y)=1\ \forall\ y\in H\Longrightarrow x\in G' \quad \text{y}
\quad \varphi(u,v)=1\ \forall\ u\in G \Longrightarrow v\in H'.
\]
\end{definicion}

\begin{teorema}\label{T12.8}
Si $G$ y $H$ son grupos finitos y existe un apareamiento no
degenerado $\varphi\colon G\times H\to{\ma C}^{\ast}$, entonces
\begin{gather*}
\begin{array}{rcl}
\psi\colon G&\longto &\hat{H}\\g&\longmapsto &\varphi(g,\_)
\end{array}\quad,\quad
\begin{array}{rcl}
\theta\colon H&\longto& \hat{G}\\
h&\longmapsto&\varphi(\_,h)
\end{array}
\end{gather*}
son epimorfismos de grupos y $\ker \psi=G'$, $\ker \theta=H'$.
En particular, $G/G'\cong \hat{H}$; $H/H'\cong \hat{G}$ y $G/G'\cong H/H'$. Finalmente si $G$ y $H$ son grupos abelianos, entonces $G\cong H$.
\end{teorema}

\begin{proof}
Si $g\in G$, sea $\psi_g\colon H\to {\ma C}^{\ast}$ dada por $\psi_g(h)=
\varphi(g,h)$. Entonces $\psi_g$ es un homomorfismo y por tanto
$\psi_g\in\hat{H}$.

Adem\'as, si $g,g_1\in G$, se tiene $\psi_{gg_1}(h)=\varphi(gg_1,h)=
\varphi(g,h)\varphi(g_1,h)=\psi_g(h)\psi_{g_1}(h)$ para toda $h\in H$
por lo que  $\psi_{gg_1}=\psi_g\psi_{g_1}$ y 
\[
\begin{array}{rcl}
\psi\colon G&\to &\hat{H}\\g&\mapsto&\psi_g=\varphi(g,\_)\end{array}
\]
es un homomorfismo de grupos. Si $g\in G'$, puesto que ${\ma C}^{\ast}$
es abeliano, $\psi_g=1$ y por lo tanto $g\in\ker \psi$. Se sigue que
$G'\subseteq \ker\psi$. Ahora bien, si $g\in\ker \psi$, por ser $\varphi$ no
degenerada se sigue que $g\in G'$ y  por lo tanto $G'=\ker \psi$.

Bajo $\psi$ tenemos que $G/G'\subseteq \hat{H}\cong H/H'$. De manera
an\'aloga obtenemos $H/H'\subseteq \hat{G}\cong G/G'$.
El resultado se sigue. $\fin$
\end{proof}

Sea $G$ un grupo abeliano finito y consideremos el mapeo
\[
\varphi=\langle\ ,\ \rangle\colon G\times\hat{G}\to {\ma C}^{\ast},\quad
\langle\ ,\ \rangle (g,\psi)=\langle\chi, g\rangle:=\chi(g).
\]
Entonces 
\begin{gather*}
\langle\chi_1\chi_2,g\rangle =(\chi_1\chi_2)(g)=\chi_1(g)
\chi_2(g)=\langle\chi,g\rangle\langle\chi_2,g\rangle\\
\intertext{y}
\langle \chi,gh\rangle=\chi(gh)=\chi(g)\chi(h)=\langle\chi,g\rangle
\langle\chi,h\rangle,
\end{gather*}
es decir, $\varphi$ es un apareamiento.

Si $\langle g,\chi\rangle=1$ para toda $\chi\in\hat{G}$, entonces
$\chi(g)=1$ para toda $\chi\in\hat{G}$. Sea $H=\langle g\rangle$.
Entonces dado $\chi\in\hat{G}$, $\chi$ se puede factorizar de manera
\'unica:
\[
\xymatrix{
G\ar[rr]^{\chi}\ar[rd]_{\pi}&&{\ma C}^{\ast}\\ &G/H\ar[ru]_{\tilde{\chi}}
}
\qquad \tilde{\chi}\circ \pi=\chi
\]
y por tanto $\quad \begin{array}{rcl}\hat{G}&\to&\widehat{G/H}\\
\chi&\mapsto&\tilde{\chi}\end{array}\quad$ es un monomorfismo.
Por lo tanto $|G|=|\hat{G}|\leq \big|\big(\widehat{G/H}\big)\big|=
|G/H|\leq |G|$, lo cual implica que $|H|=1$ y $g=e$.

Rec{\'\i}procamente, si $\langle g,\chi\rangle=1$ para toda $g\in G$,
entonces $\chi(g)=1$ para toda $g\in G$ lo cual, por definici\'on,
nos dice que $\chi=1$.

Esto prueba que $\varphi$ es un mapeo no degenerado. En particular,
cuando $G$ es un grupo abeliano finito, $G\cong \hat{\hat{G}}$ de
manera can\'onica. Es decir, $\quad\begin{array}{rcl}
\theta\colon G&\to&\hat{\hat{G}}\\ g&\mapsto&\hat{g}\end{array}\quad$
definido por $\hat{g}\colon \hat{G}\to{\ma C}^{\ast}$, $\hat{g}(\chi)=\chi(g)$,
es un isomorfismo de grupos.

\begin{definicion}\label{D12.9}
Sean $G$ un grupo abeliano finito y $\langle\ ,\ \rangle\colon G\times \hat{G}
\to {\ma C}^{\ast}$ el mapeo bilineal $\langle g,\chi\rangle=\chi(g)$. Sea 
$H<G$. Definimos el {\em ortogonal\index{grupo ortogonal}} de $H$ por
\[
H^{\perp}:=\langle\chi\in\hat{G}\mid \chi(h)=1\ \forall\ h\in H\rangle=
\langle \chi\in\hat{G}\mid H\subseteq \ker \chi\rangle.
\]
\end{definicion}

\begin{proposicion}\label{P12.10} Sea $G$ un grupo
abeliano finito y $H$ un subgrupo de $G$. 
Se tiene $H^{\perp}\cong
(\widehat{G/H})$.
\end{proposicion}

\begin{proof}
Si $\chi\in H^{\perp}$, $\chi$ se factoriza de manera \'unica
\[
\xymatrix{
G\ar[rr]^{\chi}\ar[dr]_{\pi}&&{\ma C}^{\ast}\\&G/H\ar[ur]_{\tilde{\chi}}
}
\]
con $\tilde{\chi}\in (\widehat{G/H})$ y viceversa, si $\tilde{\chi}\in
(\widehat{G/H})$, entonces $\tilde{\chi}$ se puede levantar a un
elemento $\chi\in\hat{G}$: $\chi(g):=\tilde{\chi}(gH)$. Por lo tanto
$\chi\to\tilde{\chi}$ es un isomorfismo entre $H^{\perp}$ y $(\widehat{
G/H})$. $\fin$
\end{proof}

Ahora bien, si consideramos el mapeo restricci\'on
\begin{eqnarray*}
\rest\colon\hat{G}&\longto&\hat{H}\\
\chi&\longmapsto& \chi|_H
\end{eqnarray*}
se tiene que $\ker(\rest)=H^{\perp}$ y por lo tanto se tiene
$\hat{G}/H^{\perp}\subseteq \hat{H}$.

\begin{proposicion}\label{P12.11}
Se tiene $\hat{H}\cong \hat{G}/H^{\perp}$, donde
$G$ es un grupo abeliano finito y $H$ es un subgrupo de $G$.
\end{proposicion}

\begin{proof}
Hemos obtenido $\hat{G}/H^{\perp}\subseteq \hat{H}$ y
\[
\big|\hat{G}/H^{\perp}\big|=\frac{|\hat{G}|}{|H^{\perp}|}=
\frac{|G|}{\big|\big(\widehat{G/H}\big)\big|}=\frac{|G|}{|G/H|}=
|H|=|\hat{H}|
\]
de donde se sigue la igualdad. $\fin$.
\end{proof}

Si ${\mc H}<\hat{G}$, se define ${\mc H}^{\perp}=\{g\in G\mid
\chi(g)=1 \ \forall\ \chi\in{\mc H}\}=\bigcap_{\chi\in{\mc H}}\ker \chi$.

Identificamos $G$ con $\hat{\hat{G}}$ como antes. Con esta
identificaci\'on, tenemos:

\begin{proposicion}\label{P12.12}
$\big(H^{\perp}\big)^{\perp}=H^{\perp \perp}=H$.
\end{proposicion}

\begin{proof}
Por la Proposici\'on \ref{P12.10} $\big(H^{\perp}\big)^{\perp}\cong
\widehat{\big(\hat{G}/H^{\perp}\big)}$ y por lo tanto
\[
|H^{\perp \perp}|=\big|\widehat{\big(\hat{G}/H^{\perp}\big)\big|}
=\big|\hat{G}/H^{\perp}\big| =\frac{|\hat{G}|}{|H^{\perp}|}=
\frac{|G|}{\big|(\widehat{G/H})\big|}=\frac{|G|}{|G/H|}=|H|.
\]
Por otro lado, si $h\in H$, entonces $h\colon\chi\to\chi(h)$ satisface
que $h(H^{\perp})=1$ por definici\'on y por tanto $h\in(H^{\perp})^{\perp}$.
Es decir, $H\subseteq (H^{\perp})^{\perp}$ y como tiene la misma
cardinalidad se sigue que $H=(H^{\perp})^{\perp}= H^{\perp \perp}$.
$\fin$
\end{proof}

Hay varios resultados elementales que son de gran utilidad y que
probaremos a continuaci\'on.

\begin{proposicion}\label{P12.13} Sea $G$ un grupo 
abeliano finito.
Sea $\chi\in\hat{G}$. Si $\chi\neq 1$, entonces $\sum_{g\in G}\chi
(g)=0$. Si $\chi=1$, entonces $\sum_{g\in G}\chi(g)=|G|$.
\end{proposicion}

\begin{proof}
Si $\chi$ es trivial, esto es, $\chi=1$, es inmediato que
$\sum_{g\in G}\chi(g)=|G|$. Ahora, si $\chi\neq 1$, existe un
elemento $h\in G$ tal que $\chi(h)\neq 1$. Sea $t=\sum_{g\in G}
\chi(g)\in{\ma C}$. Entonces
\[
\chi(h) t=\chi(h)\sum_{g\in G}\chi(g)=\sum_{g\in G}\chi(h)\chi(g)=
\sum_{g\in G}\chi(hg)=\sum_{g\in G}\chi(g)=t.
\]
Por tanto $t(\chi(h)-1)=0$. Puesto que $\chi(h)\neq 1$ se sigue
que $t=0$. $\fin$
\end{proof}

\begin{proposicion}\label{P12.14}
Sean $G$ un grupo abeliano finito y $H<G$. Entonces $G$ tiene
un subgrupo isomorfo a $G/H$.
\end{proposicion}

\begin{proof}
Una prueba directa de este resultado es simplemente el
teorema de estructura de los grupos abelianos finitamente 
generados. Ahora, usando la teor{\'\i}a de caracteres, 
tenemos
\begin{gather*}
G/H\cong \big(\widehat{G/H}\big)\cong H^{\perp}\subseteq
\hat{G}\cong G. \tag*{$\fin$}
\end{gather*}
\end{proof}

\begin{observacion}\label{O12.15}
La segunda demostraci\'on de la Proposici\'on \ref{P12.14},
nos indica que el mapeo entre redes
\begin{eqnarray*} 
{\cal A}=\{H\mid H<G\}&\longto&{\cal B}=\{G/H\mid H<G\}\\
H&\longmapsto&G/H
\end{eqnarray*}
entre los subgrupos de un grupo abeliano y sus grupos cocientes,
es biyectiva y que la red de subgrupos ${\cal A}$ de $G$ es
sim\'etrica.
\end{observacion}

\section{Caracteres de Dirichlet}\label{S12.2}

Aplicamos toda la teor{\'\i}a de caracteres al caso especial
del grupo $U_n\cong\Gal(\cic n{}/{\ma Q})$ el cual es un grupo
abeliano finito.

\begin{definicion}\label{D12.2.1}
Un {\em caracter de Dirichlet\index{caracter de
Dirichlet}\index{Dirichlet!caracter de $\sim$}} es un elemento de
$\hat{U}_n$ para alg\'un $n\in{\ma N}$. Expl{\'\i}citamente, un
caracter de Dirichlet $\chi$ es un homomorfismo de grupos
$\chi\colon U_n\to {\ma C}^{\ast}$.
\end{definicion}

Ahora bien, si $\chi\in \hat{U}_n$ y $n|m$, entonces $\chi$ se
puede considerar un elemento $\hat{U}_m$ de la siguiente
forma. Sea $\varphi_{m,n}\colon U_m\to U_n$ la proyecci\'on
natural: $\varphi_{m,n}(x\bmod m):=x\bmod n$ y sea 
$\tilde{\chi}:=\chi\circ \varphi_{m,n}$:
\[
\xymatrix{
U_m\ar[rr]^{\tilde{\chi}}\ar[dr]_{\varphi_{m,n}}&&{\ma C}^{\ast}\\
&U_n\ar[ru]_{\chi}
}
\]

En cierta forma $\chi$ y $\tilde{\chi}$ son el mismo mapeo y $\chi$
puede considerarse m\'odulo $n$ o m\'odulo $m$.
Rec{\'\i}procamente, si existe $f|n$ y $\chi'\colon U_f\to{\ma C}^{\ast}$
tal que $\chi=\chi'\circ \varphi_{n,f}$, entonces tambi\'en podemos
definir $\chi$ m\'odulo $f$ usando $\chi'$ en lugar de $\chi$.

\begin{observacion}\label{O12.2.1'}
Para $n,m\in{\ma N}$ tales que $m|n$, se tiene que
el homomorfismo
$\varphi_{n,m}\colon U_n\to U_m$ es suprayectivo. Esto se sigue
del epimorfismo natural de anillos ${\ma Z}/n{\ma Z}\to {\ma Z}/m{\ma Z}$
cuyos grupos de unidades son $U_n$ y $U_m$ respectivamente. Una demostraci\'on
directa ser{\'\i}a la siguiente: sea $a\in{\ma Z}$ tal que $\mcd (a,m)=1$.
Queremos hallar $c\in{\ma Z}$ tal que $\mcd (c,n)=1$ y
$\varphi_{n,m}(c\bmod n)=
a\bmod m$, esto es, debemos hallar $c\in{\ma Z}$ tal que
$\mcd (c,n)=1$ y $c\equiv a\bmod m$. 

Sea $h:=\prod_{p\in{\mc A}}p$ donde ${\mc A}=\{p\mid p \text{\ es primo,
$p|m$ y $p\nmid n$}\}$. Si ${\mc A}=\emptyset$, entonces $h=1$.
Ahora bien $\mcd(h,n)=1$. Por el Teorema Chino del Residuo existe
$c\in{\ma Z}$ tal que $c\equiv a\bmod n$ y $c\equiv 1\bmod h$ (por ejemplo,
sean $\alpha,\beta\in{\ma Z}$ con $\alpha h +\beta n=1$. Por tanto
$\alpha(a-1)h+\beta(a-1)n=a-1$, de donde $1+\alpha(a-1)h=a+\beta(1-a)
n=:c$ satisface $c\equiv 1\bmod h$ y $c\equiv a \bmod n$). Se tiene que 
como $c\equiv a \bmod n$, para un primo $p$ que divide a $n$, entonces
$p$ no divide a $a$ y por tanto $p$ no divide a $c$. Se sigue que 
$\mcd (c,m)=1$, que $c\equiv a\bmod n$ y, por ende, $\varphi_{
n,m}$ es suprayectivo.
\end{observacion}

\begin{definicion}\label{D12.2.2}
Dado un caracter de Dirichlet $\chi$, el m{\'\i}nimo n\'umero 
natural $f$ tal que $\chi$ puede ser definido m\'odulo $f$ se
llama el {\em conductor\index{conductor de un caracter}} de $\chi$
y se denota por ${\eu f}_{\chi}$.
\end{definicion}

M\'as precisamente, definimos para $a,b\in{\ma N}$ el siguiente
orden
\[
a\leq_{\ast} b\iff a\mid b.
\]
Se tiene que $\leq_{\ast}$ es un orden parcial. Se definen {\em los
conductores} del caracter $\chi$ definido m\'odulo $n$
como los elementos minimales bajo
el orden $\leq_{\ast}$ del conjunto $\{m\leq_{\ast} n\mid \chi \text{\ 
puede ser definido m\'odulo $m$}\}$. Veremos que hay un \'unico
elemento m{\'\i}nimo de este conjunto y \'este ser\'a el conductor
de $\chi$ (ver Observaci\'on \ref{O12.2.7} y Teorema
\ref{T12.2.8}).

\begin{observacion}\label{O12.2.3}
Dado $\chi$ un caracter de Dirichlet definido m\'odulo $n$ y si
$m|n$, entonces $\chi$ puede definirse m\'odulo $m$, es decir,
existe $\tilde{\chi}\colon U_m\to {\ma C}^{\ast}$ tal que
$\chi=\tilde{\chi}\circ \varphi_{n,m}$, si y solamente
si para cualesquiera $a,b\in{\ma Z}$ con $\mcd (a,n)=
\mcd (b,n)=1$ y tales que $a\equiv b\bmod m$, se tiene
$\chi(a\bmod n)=\chi(b\bmod n)$ (en este caso se
define $\tilde{\chi}(c\bmod m):=\chi(a\bmod n)$ donde
$\mcd (a,n)=1$ y $a\equiv c\bmod m$).

Tambi\'en se tiene que un caracter $\chi\colon U_n\lra \*{\ma C}$
puede definirse m\'odulo $f$ para $f|n$ si y solamente si se tiene
que $\chi(a+f)=\chi(a)$ para toda $a\in U_n$ y tal que $a+f\in U_n$
(puede ser que $a+f\notin U_n$. Por ejemplo si $n=6,f=3,a=5$,
$a+f=8\notin U_6$).

En efecto, si $\chi$ puede definirse m\'odulo $f$, entonces
$\chi(a+f)=\tilde{\chi}( \varphi_{n,f}((a+f) \bmod n))=\tilde{\chi}
((a+f)\bmod f)=
\tilde{\chi}(a\bmod f)=\tilde{\chi}(\varphi_{n,f}(a\bmod n))=\chi(a)$.

Rec\'iprocamente, supongamos que 
$\chi(a+f)=\chi(a)$ para toda $a\in U_n$ tal que
$a+f\in U_n$. Sea $\tilde{\chi}\colon U_f\lra \*{\ma C}$ dada por:
si $x\in U_f$, existe $a\in U_n$ tal que $x=\varphi_{n,f}(a)$
(Observaci\'on \ref{O12.2.1'}). Se define $\tilde{\chi}(x\bmod f)=
\tilde{\chi}(a\bmod f)=\chi(a\bmod n)$. Se tiene que $\tilde{\chi}$
est\'a bien definida pues si $\mcd(a,n)=\mcd(n,b)=1$ y $a\equiv
b\bmod f$, digamos, $a=b-cf$ para alg\'un $c\in{\ma Z}$, entonces
$\tilde{\chi}(b\bmod f)=\chi(b\bmod n)=\chi((a+cf)\bmod n)=\chi(
a\bmod n)=\tilde{\chi}(a\bmod f)$; $\tilde{\chi}$ es homomorfismo
y satisface $\tilde{\chi}(x\bmod f)=\tilde{\chi}(a\bmod f)=\chi(
a\bmod n)=\tilde{\chi}(\varphi_{n,f}(a\bmod n))=(\tilde{\chi}\circ
\varphi_{n,f})(a\bmod n)$. Esto es, $\chi=\tilde{\chi}\circ \varphi_{n,f}$.
\end{observacion}

\begin{ejemplo}\label{Ej12.2.4}
Consideremos $n=8$, $U_8=\{1,3,5,7\}\cong C_2\times C_2$.
Sea $\chi\colon U_8\to {\ma C}^{\ast}$ definido por $\chi(1)=
\chi(5)=1; \chi(3)=\chi(7)=-1$. Puesto que $\chi(1)=\chi(1+4)=\chi(5)$
y $\chi(3)=\chi(3+4)=\chi(7)$ se tiene que $\chi$ puede definirse
m\'odulo $4$. Se tiene que $U_4=\{\pm 1\}$ entonces
 $\chi'\colon U_4\to
{\ma C}^{\ast}$ dado por $\chi'(1)=1$ y $\chi(-1)=-1$ satisface que
\begin{align*}
\chi'(1)&= \chi\circ \varphi_{8,4}(1)=\chi(1)=\chi(5)=1;\\
\chi'(-1)&=\chi\circ \varphi_{8,4}(-1)=\chi(3)=\chi(7)=-1.
\end{align*}

Adem\'as $\chi$ no puede definirse m\'odulo $2$ pues $U_2=
\{1\}$ y $\chi(1)\neq \chi(-1)$, $1\equiv -1 \bmod 2$. Por lo 
tanto ${\eu f}_{\chi}=4$.

Ahora consideremos $\sigma\in U_8$ dado por $\sigma(1)=
\sigma(3)=1$ y $\sigma(5)=\sigma(7)=-1$. Notemos que 
$\sigma(1)\neq \sigma(5)$ y $1\equiv 5\bmod 4$ por lo que
$\sigma$ no puede definirse m\'odulo $4$ y por lo tanto ${\eu f}_{
\sigma}=8$.
\end{ejemplo}

\begin{ejemplo}\label{Ej12.2.5}
Consideremos $U_{10}=\{1,3,7,9\}$. Se tiene que $U_5=\{
1,2,3,4\}$ y $U_5\cong U_{10}$ con $\varphi_{10,5}\colon
U_{10}\to U_5$, $\varphi_{10,5}(x\bmod 10)=x\bmod 5$. 
Entonces $\varphi_{10,5}(1)=1, \varphi_{10,5}(3)=3,
\varphi_{10,5}(7)=2, \varphi_{10,5}(9)=\varphi_{10,5}(-1)=
-1=4$.

Si $\chi\in U_{10}$ entonces $\chi$ puede autom\'aticamente
definirse m\'odulo $5$ con $\chi'=\chi\circ \varphi_{10,5}$ y por
lo tanto ${\eu f}_{\chi}= 1$ o $5$. Adem\'as ${\eu f}_{\chi}=1
\iff \chi(x)=1\ \forall\ x\in U_5$.

Por ejemplo, si $\chi(1)=1, \chi(3)=\chi(7)=-1, \chi(-1)=\chi(9)=-1$,
entonces ${\eu f}_{\chi}=5$. Similarmente, si $\sigma(3)=\zeta_4,
\sigma(7)=\zeta_4^3=-\zeta_4, 
\sigma(1)=1=\zeta_4^0$ y $\sigma(9)=
{-1}$, entonces ${\eu f}_{\sigma}=5$.
\end{ejemplo}

\begin{ejemplo}\label{Ej12.2.6}
Si $p$ es un n\'umero primo impar y $\chi$
es un caracter de $U_p$, entonces
si $\chi\neq 1$ necesariamente ${\eu f}_{\chi}=p$. Si
$\sigma\in U_{p^m}$ con $m\in{\ma N}$ y en caso de que
$p$ sea par, tomamos $m\geq 2$, entonces ${\eu f}_{\sigma}
=p^t$ para alg\'un $0\leq t\leq m$.
\end{ejemplo}

\begin{observacion}\label{O12.2.7}
Si $\chi\colon U_n\to {\ma C}^{\ast}$ es un caracter de Dirichlet,
el conductor de $\chi$ se defini\'o como
un n\'umero minimal ${\eu f}_{\chi}$
tal que $\chi=\chi'\circ \varphi_{n,{\eu f}_{\chi}}$ donde 
$\varphi_{n,{\eu f}_{\chi}}$ donde $\varphi_{n,{\eu f}_{\chi}}$ es el
epimorfismo natural 
\begin{eqnarray*}
\varphi_{n,{\eu f}_{\chi}}\colon U_n &\longto& U_{{\eu f}_{\chi}}\\
x\bmod n&\longmapsto& x\bmod {\eu f}_{\chi}.
\end{eqnarray*}

Es decir el conductor ${\eu f}_{\chi}$ de $\chi$, en caso de existir,
necesariamente divide a $n$. 
Lo anterior evita situaciones como por ejemplo 
\[
\xymatrix{
U_3=\{1,2\}\ar@{->}[rr]^{\varphi}\ar@{->}[rd]_{\chi_1}
&&\{1,3\}=U_4\ar@{->}[ld]^{\chi_2}\\
&\*{\ma C}}
\qquad
\begin{array}{l}
\varphi(1)=1, \varphi(2)=3,\\
\chi_1(1)=1, \chi_1(2)=-1,\\
\chi_2(1)=1, \chi_2(3)=-1,\\
\chi_2\circ \varphi=\chi_1\text{\ pero}
\varphi\neq \varphi_{m,n}.
\end{array}
\]

Aqu{\'\i} mencionamos que ``en caso de
existir'', lo cual es obvio si tomamos al pie de la letra la definici\'on,
es decir, ${\eu f}_{\chi}$ es m{\'\i}nimo con ${\eu f}_{\chi}|n$. Sin
embargo, si pensamos la minimalidad en t\'erminos de divisibilidad
(por ejemplo $2$ y $3$ son minimales dividiendo a $24$), ya no 
es tan obvio que ${\eu f}_{\chi}$ exista.

Para precisar, digamos que tenemos $\chi$ un caracter m\'odulo
$21$ y que $\chi$ se puede definir m\'odulo $3$ y tambi\'en $7$
y que $\chi$ no es trivial.

Aunque por definici\'on, ${\eu f}_{\chi}=3$, $7$ tambi\'en es minimal
en el sentido de divisibilidad y por qu\'e no pensar que $\chi$ tiene
dos conductores: $3$ y $7$.

En seguida veremos que esto no es posible y nuestra definici\'on de
${\eu f}_{\chi}$ quedar\'a dada sin ambig\"uedad alguna.
\end{observacion}

\begin{teorema}[Existencia del conductor]\label{T12.2.8}
Sea $\chi\colon U_n\to {\ma C}^{\ast}$ un caracter de Dirichlet y 
sean $a|n$, $b|n$ tales que $\chi=\chi_a\circ \varphi_{n,a}$ y
$\chi=\chi_b\circ \varphi_{n,b}$:
\[
\xymatrix{
U_n\ar[rr]^{\chi}\ar[dr]_{\varphi_{n,a}}&&{\ma C}^{\ast}\\
&U_a\ar[ur]_{\chi_a}
}
\qquad
\xymatrix{
U_n\ar[rr]^{\chi}\ar[dr]_{\varphi_{n,b}}&&{\ma C}^{\ast}\\
&U_b\ar[ur]_{\chi_b}
}
\]

Sea $c:=\mcd (a,b)$. Entonces $\chi$ puede definirse m\'odulo
$c$. En particular, si $a$ y $b$ son dos conductores de $\chi$,
en el sentido de que no existen $x|a$, $y|b$ tales que $\chi$
pueda ser definido m\'odulo $x$ o m\'odulo $y$, entonces $a=b$.
\end{teorema}

\begin{proof}
Sea $d$ el producto de todos los n\'umeros primos $p$ que 
dividen a $n$ pero que no dividen a $b$. Entonces $c=\mcd
(da,b)$ pues si $\beta=\mcd (da,b)$, como $c|a$ y $c|b$, entonces
$c|da$ y $c|b$ por lo que $c|\beta$. Por otro lado existen
$t,s,r,l\in{\ma Z}$ tales que $td+sb=1$; $ra+lb=c$, por lo que
al multiplicar ambas igualdades obtenemos
\[
c=(tr)(ad)+(tdl+ras+slb)b
\]
lo cual implica que $\beta|c$ y que $\beta=c$.

Para ver que $\chi$ puede ser definido m\'odulo $c$, debemos
probar que si $\mcd (x,n)=1$, $\mcd (y,n)=1$ y $x\equiv y\bmod
c$, entonces $\chi(x)=\chi(y)$. Por el Teorema Chino del
Residuo, existe $\alpha\in{\ma Z}$ tal que $\alpha\equiv x\bmod da$ y
$\alpha\equiv y\bmod b$. De hecho, si $x=y+lc$ y sea
$\varepsilon da+\delta b=c$, por lo que $l\varepsilon da+l\delta b=
lc =x-y$. Entonces $\alpha:=x-l\varepsilon da = y+l\delta b$.

Veremos primero que $\mcd (\alpha,n)=1$. Supongamos que $p$
es un n\'umero primo tal que $p|\alpha$ y $p|n$. Se tiene que
$p\nmid b$ pues en caso de que $p|b$ y debido a que
$\alpha=y+mb$, entonces se tendr{\'\i}a que $p|y$ y por tanto
$p|\mcd (y,n)=1$ lo cual es absurdo. As{\'\i}, $p\nmid b$.

Ahora bien, puesto que $p|n$, entonces $p|d$ por lo que
$p|da$ y $p|\alpha$ lo cual implica que $p|x$ pero en ese caso
tendr{\'\i}amos que $p|\mcd (x,n)=1$ lo cual nuevamente es
absurdo.

En resumen, tenemos que $\mcd (\alpha,n)=1$. Se tiene
\begin{align*}
\chi(\alpha)&=\chi_a\circ \varphi_{n,a}(\alpha)=\chi_{a}\circ
\varphi_{n,a}(x)=\chi(x),\\
\chi(\alpha)&=\chi_b\circ \varphi_{n,b}(\alpha)=\chi_b\circ
\varphi_{n,b}(y)=\chi(y).
\end{align*}

Por tanto $\chi(\alpha)=\chi(x)=\chi(y)$ y $\chi$ puede ser definido
m\'odulo $c$.
\[
\xymatrix{
U_n\ar[rr]^{\chi}\ar[dr]_{\varphi_{n,c}}&&{\ma C}^{\ast}\\
&U_c\ar[ur]_{\chi_c}
}
\]

Finalmente, si $a$ y $b$ son dos conductores de $\chi$, entonces
$\chi$ se puede definir m\'odulo $c:=\mcd (a,b)$ y $c|a$, $c|b$. Por
minimalidad de $a$ y $b$ se sigue que $a=c=b$. $\fin$
\end{proof}

\begin{observacion}\label{O12.2.9} Si $n=1$, $U_{1}=\{1\}$, 
$\chi\colon U_1\to {\ma C}^{\ast}$, $1\mapsto 1$ es el \'unico
caracter m\'odulo $1$ y $\chi$ es el {\em caracter 
trivial\index{caracter trivial}}. Para todo $m\in{\ma N}$,
$
\xymatrix{
U_m\ar[rr]^{\chi}\ar[dr]_{\varphi_{m,1}}&&{\ma C}^{\ast}\\
&U_1\ar[ur]_{\tilde{\chi}}
}$
$\varphi_{m,1}(a)=1$ y $\tilde{\chi}(x)=1$ para todo
$x\in U_m$. Es decir, para todo $m$, $\tilde{\chi}$ es el 
caracter trivial y es el \'unico caracter de conductor $1$.
\end{observacion}

\begin{observacion}\label{O12.2.10} Sea $n\in{\ma N}$ impar.
Entonces
\begin{eqnarray*}
\varphi_{2n,n}\colon U_{2n}&\longto & U_n,\\
a\bmod 2n&\longmapsto & a\bmod n.
\end{eqnarray*}
es un isomorfismo. Por lo tanto si $\chi$ es cualquier
caracter m\'odulo $2n$, $\chi$ puede definirse m\'odulo $n$.
$
\xymatrix{
U_{2n}\ar[rr]^{\chi}\ar[dr]_{\varphi_{2n,n}}^{\cong}&&{\ma C}^{\ast}\\
&U_n\ar[ur]_{\tilde{\chi}}
}
$
con $\tilde{\chi}:=\chi\circ \varphi^{-1}_{2n,n}$.

En particular no pueden existir caracteres con conductor $2n$
con $n$ impar y en especial no hay caracteres de conductor $2$.
Es decir, si $m\equiv 2\bmod 4$, $m$ no es el conductor de 
ning\'un caracter de Dirichlet.
\end{observacion}

\begin{proposicion}\label{P12.2.10}
Sean $\chi,\varphi$ dos caracteres de Dirichlet de conductores
${\eu f}_{\chi}$ y ${\eu f}_{\varphi}$. Supongamos que existe
$n\in{\ma N}$ tal que ${\eu f}_{\chi}|n$ y ${\eu f}_{\varphi}|n$
y tales que $\chi$ y $\varphi$ son iguales m\'odulo $n$,
es decir, $\chi,\varphi\colon U_n\to {\ma C}^{\ast}$ satisfacen
$\chi(a\bmod n)=\varphi(a\bmod n)$ para toda $\mcd (a,n)=1$.
Entonces $\chi=\varphi$, es decir, ${\eu f}_{\chi}={\eu f}_{\varphi}=
f$ y $\chi=\varphi\bmod f$.
\end{proposicion}

\begin{proof}
Consideremos
\[
\xymatrix{
U_n\ar[rr]^{\varphi_{n,{\eu f}_{\chi}=\pi}}\ar[dr]_{\chi}
&&U_{{\eu f}_{\chi}}\ar[dl]^{\tilde{\chi}}\\
&{\ma C}^{\ast}
}\qquad \tilde{\chi}\circ \pi =\chi=\varphi.
\]
Por tanto $\tilde{\varphi}\circ \pi =\chi=\varphi$, es decir $\varphi$
se puede definir m\'odulo ${\eu f}_{\chi}$ y en particular
${\eu f}_{\varphi}|{\eu f}_{\chi}$.Por simetr{\'\i}a $
{\eu f}_{\chi}|{\eu f}_{\varphi}$ y ${\eu f}_{\varphi}={\eu f}_{\chi}$.
Por lo tanto $\varphi$ y $\chi$ son el mismo caracter
m\'odulo $f={\eu f}_{\varphi}={\eu f}_{\chi}$. $\fin$
\end{proof}

Con lo visto hasta ahora, es claro que el elemento $J\sim -1$ en 
el grupo de Galois $\Gal(\cic n{}/{\ma Q})\cong U_n$ juega
un papel importante en nuestra teor{\'\i}a. Baste decir que
$K^J=K\iff K\subseteq {\ma R}$.

\begin{definicion}\label{D12.2.11}
Sea $\chi\colon U_n\to {\ma C}^{\ast}$ un caracter de Dirichlet
m\'odulo $n$. Entonces $\chi$ se llama {\em par\index{caracter
par}} si $\chi(-1)=1$ e {\em impar\index{caracter impar}} si $\chi(-1)
=-1$ (notemos que $\chi(-1)^2=\chi\big((-1)^2\big)=\chi(1)=1$
por lo que necesariamente $\chi(-1)\in\{\pm 1\}$).
\end{definicion}

Puesto que un caracter de Dirichlet puede ser definido m\'odulo muchos
$n\in{\ma N}$, es conveniente fijar, en algunas ocasiones, el
n\'umero $n$.

\begin{definicion}\label{D12.2.12} Si un caracter $\chi$ est\'a
definido m\'odulo su conductor ${\eu f}_{\chi}$, $\chi$ se
llama {\em primitivo\index{caracter primitivo}}.
\end{definicion}

Muchas veces es conveniente definir $\chi(a)$ para $a\in{\ma Z}$
para un caracter de Dirichlet.

\begin{definicion}\label{D12.2.13} Dado un caracter de Dirichlet
$\chi$, definimos $\chi(a)=0$ para $a\in{\ma Z}$ con $\mcd (a,{\eu f}_{
\chi})\neq 1$.
\end{definicion}

De esta forma podemos considerar $\chi\colon {\ma Z}\to{\ma C}^{
\ast}$. Notemos que es importante fijar ${\eu f}_{\chi}$ y no solo
considerar $\chi$ m\'odulo $n$ pues de esta forma $\chi\colon
{\ma Z}\to {\ma C}$ no estar{\'\i}a un{\'\i}vocamente determinado.

Por ejemplo, si $\chi$ est\'a dado por
 $\chi\colon U_3\to {\ma C}$, $\chi(1)=1,\chi(2)=-1$.
Si consideramos $\chi$ m\'odulo $12$, $\chi\colon U_{12}\to{\ma C}$,
entonces puesto que $\mcd (2,12)=2\neq 1$, si defini\'esemos
$\chi(a)=0$ para $\mcd (a,12)\neq 1$, necesariamente
tendr{\'\i}amos $\chi(2)=0$.

\begin{proposicion}\label{P12.2.14} Si consideramos 
al caracter $\chi$ como en 
la Definici\'on {\rm \ref{D12.2.13}}, $\chi\colon{\ma Z}\to{\ma C}$
satisface
\las
\item $\chi(a)=0$ si $\mcd (a,{\eu f}_{\chi})\neq 1$.
\item $\chi(1)\neq 0$.
\item $\chi(ab)=\chi(a)\chi(b)$ para toda $a,b\in{\ma Z}$.
\item $\chi(a)=\chi(b)$ para $a\equiv b\bmod {\eu f}_{\chi}$.
\end{list}
\end{proposicion}

\begin{proof}
(1), (2) y (4) son por definici\'on para $\chi$ definido m\'odulo
${\eu f}_{\chi}$.

Ahora bien, si $a,b\in{\ma Z}$, entonces si $\mcd (ab,{\eu f}_{\chi})
=1$, por definici\'on tenemos que $\chi(ab)=\chi(a)\chi(b)$. Si
$\mcd (ab,{\eu f}_{\chi})\neq 1$, entonces $\mcd (a,{\eu f}_{\chi})
\neq 1$ o $\mcd (b,{\eu f}_{\chi})\neq 1$ y por lo tanto $\chi(a)=0$
o $\chi(b)=0$ de donde se sigue $0=\chi(ab)=\chi(a)\chi(b)$. 
Esto demuestra (3). $\fin$
\end{proof}

\begin{observacion}\label{O12.2.15} Notemos que $\chi\colon {\ma Z}
\to {\ma C}$ no es homomorfismo pues $\chi(a+b)
\neq \chi(a)+\chi(b)$.
Por ejemplo, si $\chi$ est\'a definido de la siguiente forma:
$\chi\colon U_3\to {\ma C}^{\ast}$, $\chi(1)=1$ y $\chi(2)=-1$,
entonces $\chi(1+1)=\chi(2)=-1\neq 2=1+1=\chi(1)+\chi(1)$.
\end{observacion}

\begin{observacion}\label{O12.2.16} A menos que se diga lo
contrario, el m\'odulo de definici\'on de un caracter de Dirichlet
ser\'a su conductor, sin embargo cuando hablemos de los 
caracteres m\'odulo $n$ entenderemos todos los caracteres $\chi$
tales que ${\eu f}_{\chi}|n$.

Tambi\'en notemos que al definir $\chi(a)=0$ cuando $\mcd (a,{\eu f}_{
\chi})\neq 1$, es definir $\chi(a)=0$ tan poco como es posible y
adem\'as $\chi\colon {\ma Z}\to {\ma C}$ es una funci\'on
peri\'odica de per{\'\i}odo ${\eu f}_{\chi}$: $\chi(a+{\eu f}_{\chi})=
\chi(a)$ para toda $a\in{\ma Z}$.
\end{observacion}

\begin{observacion}\label{O12.2.17} Notemos que si $\chi$ no
est\'a definido m\'odulo su conductor, entonces $\chi$ no tiene
per{\'\i}odo ${\eu f}_{\chi}$. Por ejemplo, si consideramos
$\chi\colon U_6\to {\ma C}^{\ast}$ dada por $\chi(1)=1$ y
$\chi(5)=-1$ el cual est\'a definido m\'odulo $6$ aunque su
conductor ${\eu f}_{\chi}$ es igual a $3$, no tiene per{\'\i}odo
${\eu f}_{\chi}=3$ pues $1=\chi(1)\neq \chi(1+3)=\chi(4)=0$.

Comparar este fen\'omeno con la Observaci\'on \ref{O12.2.3}.
\end{observacion}

\begin{definicion}\label{D12.28} El {\em caracter trivial\index{caracter
trivial}} es el \'unico caracter de conductor $1$, es decir, este es
$\chi\colon{\ma Z}\to{\ma C}$ dado por $\chi(a)=1$ para toda $a
\in{\ma Z}$ (incluyendo $\chi(0)=1$).
\end{definicion}

\begin{definicion}\label{D12.28'}
Sean $\chi, \psi$ dos caracteres definidos m\'odulo
$n$ y $m$ respectivamente. Sea $t$ tal que $n,m|t$,
esto es, $\mcm[n,m]|t$. Sea $\gamma\colon U_t\lra
\*{\ma C}$ dado por $\gamma(a)=\chi(a)\psi(a)$.
Entonces $\gamma$ es el producto $\chi$ y $\psi$.
Por la Proposici\'on \ref{P12.2.10}, $\gamma$ est\'a
bien definido, es decir, $\gamma$ est\'a univocamente
definida.
\end{definicion}

\begin{ejemplo}\label{Ej12.2.20}
Sean $\chi\colon U_8\to{\ma C}^{\ast}$, $U_8=\{1,3,5,7\}$ dado por
$\chi(1)=\chi(3)=1$; $\chi(5)=\chi(7)=-1$, el cual tiene conductor
${\eu f}_{\chi}=8$ y $\sigma\colon U_{12}\to{\ma C}^{\ast}$, $U_{12}
=\{1,5,7,11\}$ dado por $\sigma(1)=\sigma(11)=1$; $\sigma(5)=
\sigma(7)=-1$, el cual tiene conductor ${\eu f}_{\sigma}=12$.

Entonces $\mcm[8,12]=24$. Definimos $\tilde{\chi}$, $\tilde{\sigma}$ 
m\'odulo $24$ donde se tiene
$U_{24}=\{1,5,7,11,13,17,19,23\}$. Tenemos:
\begin{eqnarray*}
\pi_{24,8}\colon U_{24}&\longto&U_8,\\
x \bmod 24&\longmapsto& x\bmod 8,\\
1,17&\longmapsto& 1,\\
5,13&\longmapsto& 5,\\
7,23&\longmapsto&7,\\
11,19&\longmapsto& 3.
\end{eqnarray*}
Por lo tanto
\begin{gather*}
\tilde{\chi}=\chi\circ \pi_{24,8}\colon U_{24}\longto {\ma C}^{\ast},\\
\tilde{\chi}(1)=\tilde{\chi}(17)=\chi(1)=1,\\
\tilde{\chi}(5)=\tilde{\chi}(13)=\chi(5)=-1,\\
\tilde{\chi}(7)=\tilde{\chi}(23)=\chi(7)=-1,\\
\tilde{\chi}(11)=\tilde{\chi}(19)=\chi(3)=1
\end{gather*}
y
\begin{eqnarray*}
\pi_{24,12}\colon U_{24}&\longto&U_{12},\\
x \bmod 24&\longmapsto& x\bmod 12,\\
1,13&\longmapsto& 1,\\
5,17&\longmapsto& 5,\\
7,19&\longmapsto&7,\\
11,23&\longmapsto& 11.
\end{eqnarray*}
Por lo tanto
\begin{gather*}
\tilde{\sigma}=\sigma\circ \pi_{24,12}\colon U_{24}
\longto {\ma C}^{\ast},\\
\tilde{\sigma}(1)=\tilde{\sigma}(13)=\sigma(1)=1,\\
\tilde{\sigma}(5)=\tilde{\sigma}(17)=\sigma(5)=-1,\\
\tilde{\sigma}(7)=\tilde{\sigma}(19)=\sigma(7)=-1,\\
\tilde{\sigma}(11)=\tilde{\sigma}(23)=\sigma(11)=1
\end{gather*}

Se sigue que
\begin{gather*}
\gamma=\tilde{\sigma}\tilde{\chi}\colon U_{24}
\longto {\ma C}^{\ast},\\
\gamma(1)=\tilde{\sigma}(1)\tilde{\chi}(1)=1,\\
\gamma(5)=\tilde{\sigma}(5)\tilde{\chi}(5)=1,\\
\gamma(7)=\tilde{\sigma}(7)\tilde{\chi}(7)=1,\\
\gamma(11)=\tilde{\sigma}(11)\tilde{\chi}(11)=1,\\
\gamma(13)=\tilde{\sigma}(13)\tilde{\chi}(13)=-1,\\
\gamma(17)=\tilde{\sigma}(17)\tilde{\chi}(17)=-1,\\
\gamma(19)=\tilde{\sigma}(19)\tilde{\chi}(19)=-1,\\
\gamma(23)=\tilde{\sigma}(23)\tilde{\chi}(23)=-1.
\end{gather*}
En este caso $\gamma$ es primitivo, ${\eu f}_{\gamma}=24$, y
$\gamma$ es el producto de $\chi$ y $\sigma$: 
$\gamma=\chi\sigma$.
\end{ejemplo}

\begin{ejemplo}
Sean ahora $\theta\colon U_3\to {\ma C}^{\ast}$ dada por
$\theta(1)=1$, $\theta(2)=-1$ y $\sigma\colon U_{12}\to {\ma C}^{\ast}
$ como en el Ejemplo \ref{Ej12.2.20}. Entonces ${\eu f}_{\theta}=
3$, ${\eu f}_{\sigma}=12$ y $\mcm[3,12]=12$. Sea $\tilde{\theta}=\theta
\circ \pi_{12,3}$ y $\tilde{\sigma}=\sigma$. Entonces
$\gamma=\tilde{\sigma}\tilde{\theta}\colon U_{12}\to {\ma C}^{\ast}$
est\'a dada por
\begin{gather*}
\tilde{\sigma}\tilde{\theta}(1)=\tilde{\sigma}(1)\tilde{\theta}(1)=1,\\
\tilde{\sigma}\tilde{\theta}(5)=\tilde{\sigma}(5)\tilde{\theta}(5)=1,\\
\tilde{\sigma}\tilde{\theta}(7)=\tilde{\sigma}(7)\tilde{\theta}(7)=-1,\\
\tilde{\sigma}\tilde{\theta}(11)=\tilde{\sigma}(11)\tilde{\theta}(11)=-1.
\end{gather*}

Ahora bien $\gamma$ puede ser definido m\'odulo $4$: Sea
$\xi\colon U_4\to{\ma C}^{\ast}$, dada por $\xi(1)=1$, $\xi(3)=-1$.
Entonces $\xi\circ
\pi_{12,4}=\gamma$.
Por lo tanto $\gamma$ no es primitivo.
 Se sigue que $\xi=\theta\sigma$ y ${\eu f}_{
\xi}=4$.
\end{ejemplo}

\begin{definicion}\label{D12.2.22}
Sea $\chi\colon U_n\to {\ma C}^{\ast}$ cualquier caracter de
Dirichlet definido m\'odulo $n$. Definimos el 
{\em conjugado\index{caracter conjugado}} de $\chi$ por
\[
\overline{\chi}\colon U_n\longto {\ma C}^{\ast}, \quad 
\overline{\chi}(a):=\overline{\chi(a)}=\chi(a)^{-1},
\]
pues $\chi(U_n)$ es un subgrupo de las $n$--ra\'ices de $1$ de 
$\*{\ma C}$.
Entonces $\chi\overline{\chi}$ es el caracter trivial y ${\eu f}_{\chi
\overline{\chi}}=1$.
\end{definicion}

\begin{observacion}\label{O12.2.22'}
Se tiene ${\eu f}_{\overline{\chi}}={\eu f}_{\chi}$.
\end{observacion}

El siguiente resultado es muy \'util para el c\'alculo de conductores.

\begin{teorema}\label{T12.2.23} Sean $\chi$, $\psi$ dos caracteres
de Dirichlet cuyos conductores son primos relativos $\mcd ({\eu f}_{
\chi},{\eu f}_{\psi})=1$. Entonces ${\eu f}_{\chi\psi}={\eu f}_{\chi}
{\eu f}_{\psi}$.
\end{teorema}

\begin{proof}
Sean $n={\eu f}_{\chi}$, $m={\eu f}_{\psi}$. Entonces $\chi\colon
U_n\to{\ma C}^{\ast}$, $\psi\colon U_m\to {\ma C}^{\ast}$. Se define
$\gamma\colon U_{[n,m]}=U_{nm}\to {\ma C}^{\ast}$ por
$\gamma(a):=\tilde{\chi}(a)\tilde{\psi}(a)$ donde $\tilde{\chi}=\chi\circ
\pi_{nm,n}$, $\tilde{\psi}=\psi\circ\pi_{nm,m}$. Definimos $\varphi
=\gamma\chi^{-1}\colon U_{[[n,m],n]}=U_{[nm,n]}=U_{nm}\to
{\ma C}^{\ast}$.
Se tiene que 
\begin{align*}
\varphi(a\bmod nm)&=\gamma(a\bmod nm)\tilde{\chi}^{
-1}(a\bmod nm)=\\
&=\tilde{\chi}(a\bmod nm)\tilde{\psi}(a\bmod nm)
\tilde{\chi}^{-1}(a\bmod nm)=\tilde{\psi}(a\bmod nm).
\end{align*}
Es decir, $\varphi=\tilde{\psi}$. Por tanto ${\eu f}_{\varphi}=
{\eu f}_{\tilde{\psi}}={\eu f}_{\psi}=m$.

Ahora bien, ${\eu f}_{\chi\psi}|[{\eu f}_{\chi},{\eu f}_{\psi}]=nm$, por lo
tanto
\[
m={\eu f}_{\psi}={\eu f}_{\varphi}={\eu f}_{(\chi\psi)\chi^{-1}}|
[{\eu f}_{\chi\psi},{\eu f}_{\chi^{-1}}]=[{\eu f}_{\chi\psi},{\eu f}_{\chi}]=
[{\eu f}_{\chi\psi},n]=\frac{{\eu f}_{\chi\psi} n}{\mcd ({\eu f}_{\chi\psi},
n)}.
\]
Sea $r:=\frac{{\eu f}_{\chi\psi} n}{\mcd 
({\eu f}_{\chi\psi},n)}$. Entonces
$m|r$. Tambi\'en tenemos que
\[
m\big|\frac{{\eu f}_{\chi\psi} n}{\mcd ({\eu f}_{\chi\psi},n)}
={\eu f}_{\chi\psi}
\frac{n}{\mcd ({\eu f}_{\chi\psi},n)}= {\eu f}_{\chi\psi} n_1\quad
\text{con}\quad n_1|n.
\]
Nuevamente, del hecho de que $\mcd(m,n)=1$, se sigue que
$m|{\eu f}_{\chi\psi}$.

Similarmente, se tiene que $n|{\eu f}_{\chi\psi}$ y puesto que
$\mcd (m,n)=1$, $mn=[{\eu f}_{\chi},{\eu f}_{\psi}]|{\eu f}_{\chi\psi}$.
Se sigue que ${\eu f}_{\chi\psi}=nm={\eu f}_{\chi}{\eu f}_{\psi}$.
$\fin$
\end{proof}

\begin{definicion}\label{D12.2.24}
Un caracter $\chi\colon G\to{\ma C}^{\ast}$ se dice que es un
{\em caracter de Galois\index{caracter de Galois}} si $G$ es el
grupo de Galois de una extensi\'on finita de campos $L/K$:
$G=\Gal(L/K)$.
\end{definicion}

Tenemos que los caracteres de Dirichlet pueden ser vistos como
caracteres de Galois pues si $\chi$ es un caracter de Dirichlet,
$\chi\colon U_n\to{\ma C}^{\ast}$, entonces $U_n\cong \Gal(\cic n{}/
{\ma Q})$.

Veamos como podemos usar los caracteres de Dirichlet para 
estudiar la aritm\'etica de la extensi\'on $\cic n{}/{\ma Q}$.

\begin{ejemplo}\label{Ej12.2.25}
Se tiene que $U_8\cong \Gal(\cic 8{}/{\ma Q})$. Sea $\chi\colon
U_8\to{\ma C}^{\ast}$ el caracter dado por $\chi(1)=\chi(7)=1$
y $\chi(3)=\chi(5)=-1$.

Se tiene que $\ker \chi=\{1,7\bmod 8\}$. Notemos que 
\[
\ker \chi\cong \Gal(\cic 8{}/{\ma Q}(\sqrt{2}))=\Gal(\cic 8{}/\cic 8{}^+)
\cong \{\pm 1\bmod 8\}.
\]
Esto es, ${\ma Q}(\sqrt{2})=\cic 8{}^+=\cic 8{}^{\ker \chi}$.
Por tanto $\chi$ es un caracter de $\frac{\Gal(\cic 8{}/{\ma Q})}{
\Gal(\cic 8{}/{\ma Q}(\sqrt{2}))}\cong \Gal({\ma Q}(\sqrt{2})/{\ma Q})$.

M\'as precisamente, si $G=\Gal(\cic 8{}/{\ma Q})$, $\chi\colon G\to
{\ma C}^{\ast}$ y $H=\ker \chi<G$, $\chi$ se factoriza:
\[
\xymatrix{
G\ar[rr]^{\chi}\ar[rd]_{\pi}&&{\ma C}^{\ast}\\&G/H\ar[ru]_{\tilde{\chi}}
}
\qquad \tilde{\chi}\colon G/H\to{\ma C}^{\ast},\qquad \tilde{\chi}\circ
\pi=\chi
\]
y $G/H\cong \Gal(\cic 8{}^H/{\ma Q})=\Gal({\ma Q}(\sqrt{2})/{\ma Q})$.

Similarmente, si $\varphi\colon U_8\to{\ma C}^{\ast}$ est\'a dado por
$\varphi(1)=\varphi(5)=1$ y $\varphi(3)=\varphi(7)=-1$, entonces
$\ker \varphi=\{1,5\bmod 8\mid 8\}$ y $\ker\varphi\cong
\Gal(\cic 8{}/\cic 4{})=
\{x\bmod 8\mid x\equiv 1\bmod 4\}=D_{8,4}$
y por tanto $\varphi$ es un caracter de 
\[
\frac{\Gal(\cic 8{}/{\ma Q})}{
\Gal(\cic 8{}/\cic 4{})}\cong \Gal(\cic 4{}/{\ma Q})\cong U_4.
\]
Es decir, $\cic 4{}=\cic 8{}^{\ker \varphi}$.

Notemos que ${\eu f}_{\varphi}=4$ y que ${\eu f}_{\chi}=8$.
\end{ejemplo}

Sea $\chi$ un caracter de Dirichlet definido m\'odulo $n$, $\chi\colon
U_n\to{\ma C}^{\ast}$. Entonces, si $H=\ker \chi$, $\chi$ es un
caracter de Galois de $\Gal(\cic n{}^H/{\ma Q})$:
\begin{gather*}
\xymatrix{
U_n\cong\Gal(\cic n{}/{\ma Q})\ar[rr]^{\chi}\ar[dr]_{\pi}&&
{\ma C}^{\ast}\\&U_n/\ker \chi\ar[ru]_{\tilde{\chi}}
}\\
 U_n/H\cong U_n/(\ker \chi)=\frac{\Gal(\cic n{}/{\ma Q})}{
\Gal(\cic n{}/\cic n{}^H)}\cong \Gal(\cic n{}^H/{\ma Q}).
\end{gather*}

\begin{definicion}\label{D12.2.26}
Sea $\chi\colon U_n\to {\ma C}^{\ast}$ un caracter de Dirichlet.
Se define el {\em campo que pertenece\index{campo perteneciente
a un caracter de Dirichlet}} a $\chi$ o que esta {\em
asociado\index{campo asociado a un caracter de Dirichlet}}
a $\chi$ por
$K:=\cic n{}^{\ker \chi}\subseteq \cic n{}$.
\end{definicion}

\begin{teorema}\label{T12.2.27} Sea $\chi\colon U_n\to{\ma C}^{\ast}$
un caracter de  Dirichlet con conductor $f:={\eu f}_{\chi}|n$. Sea
$\tilde{\chi}$ el caracter primitivo asociado a $\chi$: $\tilde{\chi}\circ
\pi_{n,f}=\chi$. Entonces los campos asociados a $\tilde{\chi}$ y a
$\chi$ son el mismo. En particular, el campo asociado a $\chi$
depende \'unicamente de $\chi$ y no de m\'odulo en donde est\'a
definido.
\end{teorema}

\begin{proof}
El resultado se sigue inmediatamente de la Observaci\'on
\ref{O12.2.28}. Presentamos aqu{\'\i} la prueba por completes.

Puesto que $\tilde{\chi}\circ \pi_{n,f}=\chi$, se tiene que 
$\ker \chi=\pi_{n,f}^{-1}(\ker \tilde{\chi})$. Sean $K$ y $L$ los
campos asociados a $\chi$ y $\tilde{\chi}$ respectivamente:
$K:=\cic n{}^{\ker \chi}$, $L:=\cic n{}^{\ker \tilde{\chi}}$. Primero
notemos que
\begin{align*}
[K:{\ma Q}]&=\Big|\frac{U_n}{\ker \chi}\Big|=\frac{|U_n|}{|\pi_{n,f}^{-1}(
\ker \tilde{\chi})|}=\frac{|U_n|}{|\ker \pi_{n,f}||\ker \tilde{\chi}|}\\
&=\frac{|U_f|}{|\ker \tilde{\chi}|}=[L:{\ma Q}].
\end{align*}

Por otro lado, si $\alpha\in L$, $\sigma\alpha=\alpha$ para toda
$\sigma\in\ker \tilde{\chi}$. Sea $\theta\in\ker \chi$. Entonces
$\pi_{n,f}\theta\in \ker \tilde{\chi}$, por lo tanto $(\pi_{n,f}\theta)(\alpha)
=\alpha=\pi_{n,f}(\theta(\alpha))=\tilde{\theta}(\alpha)$. Por lo 
tanto $L\subseteq K$ de donde se sigue que $L=K$.

La \'ultima parte se sigue inmediatamente pues todas los caracteres
de Dirichlet est\'an un{\'\i}vocamente determinados por el caracter
primitivo asociado. $\fin$
\end{proof}

\begin{observacion}\label{O12.2.28} De hecho lo establecido
en el Teorema \ref{T12.2.27} sucede
en general. Si $G:=\Gal(F/E)$, y $\pi\colon G\to G_1$ es un
epimorfismo, $G_1\cong G/H$. Si $M=F^H$, se tiene para
$R\subseteq G/H$ que $M^R=F^{\pi^{-1}(R)}$. $
\xymatrix{
F\ar@{-}[d]_H\\M\ar@{-}[d]\\E}$
\end{observacion}

\begin{observacion}\label{O12.2.29} Notemos que si $\chi$
es un caracter de Dirichlet y $X:=\langle\chi\rangle$ es el grupo
generado por $\chi$, entonces $\bigcap_{\varphi\in X}\ker \varphi=
\bigcap_{i=0}^{t-1}\ker \chi^i$ donde $o(\chi)=t$. Se tiene que
$\ker \chi^i\supseteq \ker \chi$ por lo que $\bigcap_{\varphi\in X}
\ker \varphi=\ker \chi$. En otras palabras $\cic n{}^{\cap_{\varphi
\in X}\ker \varphi}=\cic n{}^{\ker \chi}$.

M\'as generalmente, si $Y$ es cualquier conjunto de caracteres
de Dirichlet y si $X:=\langle Y\rangle$ es el grupo
generado por los elementos de $Y$, entonces si $n:=
\mcm \big[{\eu f}_{\chi}\mid \chi\in Y\big]$ y todos los elementos
de $Y$ los consideramos m\'odulo $n$, se sigue que $X$ es un grupo
de caracteres m\'odulo $n$, es decir, $X\subseteq \hat{U}_n$.
En este caso, con el argumento anterior, se tiene que
$\cic n{}^{\cap_{\chi\in Y}\ker \chi}=\cic n{}^{\cap_{\chi\in X}\ker \chi}$.
\end{observacion}

\begin{definicion}\label{D12.2.30} Sea $X$ cualquier grupo de
caracteres de Dirichlet definidos m\'odulo $n$. Entonces se define
el {\em campo que pertenece\index{campo perteneciente a un
grupo de caracteres de Dirichlet}} a $X$ o el {\em campo
asociado\index{campo asociado a un grupo de caracteres de
Dirichlet}} a $X$ por $K:=\cic n{}^H$ donde $H=\bigcap_{\chi\in X}
\ker \chi$.
\end{definicion}

Como antes, si $Y$ es un conjunto arbitrario de caracteres de
Dirichlet, entonces $H=\bigcap_{\chi\in Y}\ker \chi=\bigcap_{\chi\in
\langle Y\rangle=X}\ker \chi$.

Como veremos a continuaci\'on, si $X$ es un grupo de caracteres
de Dirichlet y $K$ es su campo asociado, entonces
\[
X=\widehat{\Gal(K/{\ma Q})}\cong \Gal(K/{\ma Q})
\]
y en particular $|X|=[K:{\ma Q}]$.

\begin{ejemplo}\label{Ej12.2.31} Sea $\chi$ un caracter de orden
$2$ y definido m\'odulo $n$: $\chi\in\hat{U}_n$, $\chi^2=1$, $\chi
\neq 1$. Entonces $\chi(U_n)=\{\pm 1\}$ y $\ker \chi = H=
\{a\in U_n\mid \chi(a)=1\}$. Se tiene que $\Gal(\cic n{}/\cic n{}^H)\cong
H$ y por tanto $\Gal(\cic n{}^H/{\ma Q})\cong \frac{\Gal(\cic n{}/
{\ma Q})}{\Gal(\cic n{}/\cic n{}^H)}\cong U_n/H\cong \langle\tilde{\chi}
\rangle$ donde $\tilde{\chi}$ es la factorizaci\'on de $\chi$:
\[
\xymatrix{
U_n\ar[rr]^{\chi}\ar[dr]_{\pi}&&{\ma C}^{\ast}\\
&U_n/H\cong U_n/(\ker \chi)\ar[ur]_{\tilde{\chi}}
}
\]
Por tanto $K$, el campo cuadr\'atico que pertenece a $\chi$, es una
extensi\'on cuadr\'atica de ${\ma Q}$ contenida en $\cic n{}$.
En particular, si $n$ es un primo impar, $n=p$, entonces
$K={\ma Q}\big(\sqrt{(-1)^{(p-1)/2}p}\big)$  (ver Secci\'on
\ref{S4.10.3}).
\end{ejemplo}

\begin{ejemplo}\label{Ej12.2.3}
Consideremos $\cic n{}^+$ el subcampo real de $\cic n{}$. Entonces
$\cic n{}^+=\cic n{}^{\{1,J\}}$ donde $J$ denota conjugaci\'on 
compleja. Entonces $\cic n{}^+$ pertenece al conjunto de caracteres
$X\subseteq \hat{U}_n\cong \widehat{\Gal(\cic n{}/{\ma Q})}$
tales que $\bigcap_{\chi\in X}\ker \chi=\{1, J\}=\{\pm 1\}$, esto es a
$X=\{\chi\mid \chi(-1)=1\}$.
En otras palabras, bajo el mapeo 
\begin{eqnarray*}
U_n\times \hat{U}_n&\longto&{\ma C}^{\ast}\\
(a,\sigma)&\longmapsto&\sigma(a)
\end{eqnarray*}
se tiene $X=\{\pm 1\}^{\perp}$ (ver Definici\'on \ref{D12.9}).

Ahora bien, si $\chi$ es cualquier caracter definido m\'odulo
$n$, y si $K$ es el campo asociado a $\chi$, se tiene que
\[
K\subseteq {\ma R}\iff \text{el mapeo\ }\zeta_n\stackrel{\sigma}{\longto}\zeta_n^{-1},\text{\ satisface\ }
\sigma\in\ker \chi\iff \chi(-1)=1.
\]
Es decir, $K$ es real si y solamente si $\chi$ es par.
\end{ejemplo}

\begin{ejemplo}\label{Ej12.2.32}
Consideremos en $\cic 8{}$ el subcampo cuadr\'atico real,
esto es,
${\ma Q}(\sqrt{2})=\cic 8{}^+$. Entonces ${\ma Q}(\sqrt{2})$ es
el campo que pertenece al caracter dado por $\chi\colon
U_8\to{\ma C}^{\ast}$, $\chi(1)=\chi(7)=1$ y $\chi(3)=\chi(5)=
-1$ pues $7\equiv -1\bmod 8$ y $\chi(-1)=1$ y ${\ma Q}(\sqrt{2})$
es un campo real.

Notemos que $\f{\chi}=8$ pues el \'unico caracter
de conductor $4$ necesariamente tiene como su campo asociado
a $\cic 4{}$.
\end{ejemplo}

\begin{ejemplo}\label{Ej12.2.33}
En $\cic {12}{}$, tenemos que $\cic {12}{}={\ma Q}(\zeta_3,\zeta_4)
={\ma Q}(\sqrt{3},\sqrt{-3})$. Hay tres subcampos cuadr\'aticos
de $\cic {12}{}$, a saber: $\cic 4{}={\ma Q}(\sqrt{-1})$, $\cic 3{}=
{\ma Q}(\sqrt{-3})$ y $\cic {12}{}^+={\ma Q}(\sqrt{3})$.

Entonces al considerar el caracter $\chi\colon U_{12}\to{\ma 
C}^{\ast}$, $\chi(1)=\chi(11)=1$ y $\chi(5)=\chi(7)=-1$, tenemos
que $\ker \chi=\{1,11\}$ y $o(\chi)=2$. Puesto que $\chi(11)=
\chi(-1)=1$ el campo que pertenece a $\chi$ es $K={\ma Q}
(\sqrt{3})=\cic {12}{}^+={\ma Q}(\zeta_{12}+\zeta_{12}^{-1})$.

Otra forma de concluir lo mismo, es notar que hay tres caracteres
cuadr\'aticos m\'odulo $12$, pero dos de ellos tienen conductores
$3$ (el asociado a $\cic 3{}$) y $4$ (el asociado a $\cic 4{}$). 
Puesto que el conductor de $\chi$ es $12$ se concluye que el
campo asociado es necesariamente ${\ma Q}(\sqrt{3})$.
\end{ejemplo}

Nuestro siguiente objetivo es mostrar que si $X$ es un grupo de
caracteres entonces $X\cong\widehat{\Gal(K/{\ma Q})}$ donde
$K$ es el campo asociado a $X$. Para ello consideraremos un
mapeo bilineal no degenerado.

Sea pues $X$ un grupo de caracteres de Dirichlet y sea $K$ el
campo asociado a $X$, es decir, $K=\cic n{}^H$ donde 
$n=\mcm \{\f{\chi}\mid \chi\in X\}$, $H=\bigcap_{\chi\in X}\ker \chi$.

Sea $\begin{array}{rcl}
\varphi\colon \Gal(K/{\ma Q})\times X&\longto&{\ma C}^{\ast}\\
(\sigma,\chi)&\longmapsto&\chi(\sigma)
\end{array}$ el mapeo natural. Veamos que $\varphi$ est\'a bien
definido.

Se tiene que $\sigma\in\Gal(K/{\ma Q})\cong U_n/H$. Ahora bien
$\chi\in X\subseteq \hat{U}_n$, es decir, $\chi\colon U_n\to
{\ma C}^{\ast}$. Puesto que $H\subseteq \ker \chi$, es decir,
$\chi(H)=1$, se tiene que $\chi$ se factoriza de manera
\'unica: $\tilde{\chi}\circ \pi=\chi$.
\[
\xymatrix{
U_n\ar[rr]^{\chi}\ar[rd]_{\pi}&&{\ma C}^{\ast}\\
&U_n/H\ar[ru]_{\tilde{\chi}}
}
\qquad\qquad _{U_n}\left\{\begin{array}{c}
\xymatrix{
\cic n{}\ar@{-}[d]^{H}\\ K\ar@{-}[d]^{\Big\} U_n/H}\\{\ma Q}
}\end{array}\right.
\]
y por tanto $\chi(\sigma)$ est\'a bien definido.

\begin{teorema}\label{T12.2.34} Se tiene que $\varphi$ es un
mapeo bilineal no degenerado, esto es, si $\chi(\sigma)=1$ para
toda $\sigma\in G:=\Gal(K/{\ma Q})$, entonces $\chi=1$ y si
$\chi(\sigma)=1$ para toda $\chi\in X$, entonces $\sigma=
\Id_K$.
\end{teorema}

\begin{proof}
Es inmediato que $\varphi$ es bilineal. Ahora bien, si $\chi(\sigma)
=1$ para toda $\sigma\in G\cong U_n/H$, entonces $\chi\colon
U_n\to {\ma C}^{\ast}$ es trivial pues $\chi(H)=1$ y al factorizar
$\chi$ a trav\'es de $H$, $\tilde{\chi}\colon U_n/H\to{\ma C}^{\ast}$,
$\tilde{\chi}(\sigma)=1$ para toda $\sigma\in U_n/H$, por lo que
$\chi=1$.

Rec{\'\i}procamente, si $\chi(\sigma)=1$ para toda $\chi\in X$,
entonces se tiene que
$\sigma\in\bigcap_{\chi\in X}\ker \chi=H$, por lo tanto
$\overline{\sigma}=1$ en $U_n/H=G$. $\fin$
\end{proof}

\begin{corolario}\label{C12.2.35} Se tiene que $X\cong \hat{G}=
\widehat{\Gal(K/{\ma Q})}$ y en particular $|X|=[K:{\ma Q}]$.
\end{corolario}

\begin{proof}
Esto no es m\'as que una aplicaci\'on del Teorema \ref{T12.8}. $\fin$
\end{proof}

El mismo mapeo $\varphi$ nos sirve, como vimos en general, para
dar un isomorfismo de redes entre los subcampos de $K$ y los
subgrupos de $X$. M\'as precisamente, sea $F\subseteq K$.

$\begin{array}{c}K\\ \big| \\F \\ \big| \\ {\ma Q}\end{array}$\qquad
Sea $Y=\{\chi\in X\mid \chi(g)=1\ \forall\ g\in\Gal(K/F)\}$.
Entonces,
$Y=\Gal(K/F)^{\perp}\cong \Big(\widehat{G/\Gal(K/F)}\Big)\cong
\widehat{\Gal(F/{\ma Q})}$ (Definici\'on \ref{D12.9}, Proposici\'on
\ref{P12.10}).

Rec{\'\i}procamente, dado $Y$ un subgrupo de $X$, sea $F$ el 
campo fijo bajo el ortogonal de $Y$: $F:=K^{Y^{\perp}}$, donde
recordemos que $Y^{\perp}=\{g\in G\mid \chi(g)=1\ \forall\ \chi\in Y\}$.
Por tanto $\Gal(K/F)=\Gal(K/K^{Y^{\perp}})\cong Y^{\perp}$. Por la
Proposici\'on \ref{P12.12} tenemos
\[
Y=Y^{\perp\perp}=\Gal(K/F)^{\perp}=\widehat{\Gal(F/{\ma Q})}.
\]
Esta correspondencia es biyectiva como lo prueba el siguiente
resultado.

\begin{teorema}\label{T12.2.36}
Existe una biyecci\'on entre los subgrupos de $X$ y los subcampos
de $K$ dada por
\begin{eqnarray*}
\Gal(K/F)^{\perp}&\longleftarrow&F\\
Y&\longto&K^{Y^{\perp}}
\end{eqnarray*}
\end{teorema}

\begin{proof}
Sean ${\cal A}=\{Y\mid Y<X\}$ y ${\cal B}=\{F\mid F \text{\ es
subcampo de\ } K\}$. Sean
\begin{alignat*}{6}
{\cal A}&\stackrel{\theta}{\longto}&&\  {\cal B} \qquad \qquad &
{\cal B}\ &&\stackrel{\delta}{\longto}&\ {\cal A}\\
Y&\longmapsto&& \ K^{Y^{\perp}}\qquad\qquad & F\ 
&&\longmapsto &\ 
\Gal(K/F)^{\perp}
\end{alignat*}

Se tiene que 
\begin{gather*}
(\theta\circ\delta)(F)=\theta\big(\Gal(K/F)^{\perp}\big)
=K^{(\Gal(K/F)^{\perp})^{\perp}}=K^{\Gal(K/F)}=F\\
\intertext{y}
(\delta\circ\theta)(Y)=\delta\big(K^{Y^{\perp}}\big)=\Gal(K/
K^{\perp})^{\perp}=(Y^{\perp})^{\perp}=Y
\end{gather*}
lo cual demuestra que $\theta$ y $\delta$ son biyecciones, cada
una de ellas inversa de la otra. $\fin$
\end{proof}

\begin{observacion}\label{O12.2.37}
Se tiene que el isomorfismo
\[
Y=\widehat{\Gal(F/{\ma Q})}\cong \Gal(F/{\ma Q})
\]
para $Y<X$, se expresa a trav\'es del mapeo bilineal
\begin{eqnarray*}
\mu\colon \Gal(F/{\ma Q})\times Y&\longto &{\ma C}^{\ast}\\
(g,\sigma)&\longmapsto & \mu(g,\sigma):=\sigma(g).
\end{eqnarray*}
\end{observacion}

\begin{proposicion}\label{P12.2.38}
Sean $X_1$, $X_2$ dos grupos de caracteres de Dirichlet y sean
$K_1$, $K_2$ sus campos asociados. Entonces
\las
\item $X_1\subseteq X_2\iff K_1\subseteq K_2$,
\item $X_1\cap X_2$ corresponde al campo $K_1\cap K_2$,
\item El grupo generado por $X_1$ y $X_2$: $\langle X_1,X_2\rangle
=X_1\cdot X_2$, corresponde al campo generado por $K_1$ y $K_2$:
$K_1 K_2$.
\end{list}
\end{proposicion}

\begin{proof}
Sean $X:=\langle X_1,X_2\rangle$ y $F$ es campo correspondiente
a $X$. M\'as precisamente, sean
\begin{gather*}
n:=\mcm \{\f{\chi}\mid\chi\in X\},\qquad n_i:=\mcm \{\f{\chi}\mid
\chi\in X_i\}, \quad i=1,2.\\
\intertext{Entonces $n_i|n$, $i=1,2$. Sean}
H_1=\bigcap_{\chi\in X_1}\ker \chi,\quad H_2=\bigcap_{\chi\in X_2}
\ker \chi, \quad H=\bigcap_{\chi\in X}\ker \chi,\\
K_1=\cic n{}^{H_1},\quad K_2=\cic n{}^{H_2}.
\end{gather*}
\[
\xymatrix{
&\cic n{}\\ \cic {n_1}{}\ar@{-}[ru]&F\ar@{-}[u]_{H}&\cic {n_2}{}
\ar@{-}[ul]\\
&K_1K_2\ar@{-}[u]\\ K_1\ar@{-}[uu]
\ar@{-}[ru]\ar@/^/@{-}[uuur]_{H_1}
&&K_2\ar@{-}[uu]\ar@{-}[lu]\ar@/_/@{-}[luuu]^{H_2}\\
&{\ma Q}\ar@{-}[lu]\ar@{-}[ru]
}
\]

Se tiene que si $X_1\subseteq X_2$ entonces $H_1\supseteq H_2$
y por tanto $\cic n{}^{H_1}\subseteq \cic n{}^{H_2}$ lo cual nos dice
que $K_1\subseteq K_2$.

Rec{\'\i}procamente, si $K_1\subseteq K_2$, entonces $K_1=
\cic n{}^{H_1}\subseteq \cic n{}^{H_2}=K_2$ y por tanto
$H_1=\Gal(\cic n{}/\cic n{}^{H_1})\supseteq \Gal(\cic n{}/\cic n{}^{H_2}
)=H_2$. Ahora queremos ver que $H_1\supseteq H_2$ si y 
solamente si $X_1\subseteq X_2$ lo cual implicar\'a (1).

Del mapeo bilineal $\varphi\colon \Gal(F/{\ma Q})\times
X\to {\ma C}^{\ast}$, $\varphi(\sigma,\chi)=\chi(\sigma)$ obtenemos
\begin{align*}
X_i^{\perp}&=\{g\in\Gal(F/{\ma Q})\mid \chi(g)= 1\ 
\forall\ \chi\in X_i\}\\
&=\{g\in \Gal(F/{\ma Q})\mid g\in\ker \chi \ \forall\ \chi\in X_i\}\\
&=\bigcap_{\chi\in X_i}\ker \chi = H_i, \quad i=1, 2,
\end{align*}
es decir, $H_i=X_i^{\perp}$. Se sigue que si $H_1\supseteq H_2$, 
entonces $X_1^{\perp}\supseteq X_2^{\perp}$ lo cual implica que
$(X_1^{\perp})^{\perp}=X_1\subseteq X_2=(X_2^{\perp})^{\perp}$
que a su vez demuestra que $K_1\subseteq K_2$ implica $X_1
\subseteq X_2$.

Con esto terminamos la demostraci\'on de (1).

Para ver (2), sea $F$ el campo asociado a $X_1\cap X_2$. Puesto
que $X_1\cap X_2\subseteq X_i$, $i=1,2$, se sigue de (1) que
$F\subseteq K_1\cap K_2$. Ahora bien, si $W$ es el grupo de
caracteres de Dirichlet asociado a $K_1\cap K_2$, por (1)
nuevamente, se tiene $W\subseteq X_i$, $i=1,2$ de donde
obtenemos que $W\subseteq X_1\cap X_2$ y por (1) se
sigue que $K_1\cap K_2\subseteq F$. Esto es (2).

Para probar (3), sabemos por Teor{\'\i}a de Galois que se tiene
 $K_1K_2=
\cic n{}^{H_1}\cic n{}^{H_2}=\cic n{}^{H_1\cap H_2}$ y $H_1
\cap H_2=X_1^{\perp}\cap X_2^{\perp}$. Veamos que
$(X_1\cup X_2)^{\perp}=\langle X_1,X_2\rangle^{\perp}$. Notemos
que $(X_1\cup X_2)^{\perp}=\{\sigma\in G\mid \chi(\sigma)=1
\text{\ para toda\ }\sigma\in X_1\cup X_2\}$, por lo que $\chi(\sigma)
=1$ para toda $\sigma\in\langle X_1,X_2\rangle$ lo cual prueba
que $(X_1\cup X_2)^{\perp}\subseteq \langle X_1,X_2\rangle^{
\perp}$.

Ahora si $\sigma\in\langle X_1,X_2\rangle^{\perp}$ entonces
$\chi(\sigma)=1$ para todo $\chi \in\langle X_1,X_2\rangle$
y por tanto $\chi(\sigma)=1$ para todo $\chi\in X_1\cup X_2$.
Se sigue que $\langle X_1,X_2\rangle^{\perp}\subseteq
(X_1\cup X_2)^{\perp}$ y obtenemos la igualdad.

Por otro lado, puesto que $X_i\subseteq X_1\cup X_2$ si
sigue que $X_i^{\perp}\supseteq (X_1\cup X_2)^{\perp}$. Se sigue
que $X_1^{\perp}\cap X_2^{\perp}\supseteq (X_1\cup X_2)^{\perp}$.
Rec{\'\i}procamente, si $\sigma\in X_1^{\perp}\cap X_2^{\perp}$
entonces $\chi_1(\sigma)=1$ y $\chi_2(\sigma)=1$ para 
cualesquiera $\chi_i\in X_i$, $i=1,2$. De esta manera obtenemos
que $\chi(\sigma)=1$ para todo $\chi\in X_1\cup X_2$ y 
en particular $\sigma\in (X_1\cup X_2)^{\perp}$. Esto prueba que
$X_1^{\perp}\cap X_2^{\perp}=(X_1\cup X_2)^{\perp}$.

Por lo anterior, tenemos que $H_1\cap H_2=X_1^{\perp}\cap
X_2^{\perp}=(X_1\cup X_2)^{\perp}=\langle X_1,X_2\rangle^{\perp}=
X^{\perp}$ y se sigue que $K_1K_2$ corresponde a $X=\langle
X_1,X_2\rangle$. $\fin$
\end{proof}

\begin{proposicion}\label{P12.2.39}
Sea $K/{\ma Q}$ una extensi\'on abeliana finita y sea $X$ el grupo
de caracteres de Dirichlet asociado a $K$. Sea $n=\mcd\{
\f{\chi}\mid \chi\in X\}$. Entonces $n$ es el m{\'\i}nimo natural
tal que $K\subseteq \cic n{}$.
\end{proposicion}

\begin{proof}
Si $K\subseteq \cic m{}$ entonces todo caracter $\chi\in X$
pueda ser definido m\'odulo $m$ y por tanto $\f{\chi}|m$. Se
sigue que $n|m$. $\fin$
\end{proof}

\begin{definicion}\label{D12.2.39'}
Sea $K/{\ma Q}$ una extensi\'on abeliana finita. El m\'inimo
$n\in{\ma N}$ tal que $K\subseteq \cic n{}$ se llama el
{\em conductor de $K$\index{conductor de una extensi\'on
abeliana finita de ${\ma Q}$}}.
\end{definicion}

\section{Aritm\'etica de $\cic n{}$ usando caracteres}\label{S12.3}

Los caracteres de Dirichlet resultan ser una herramienta poderosa
para el estudio de las extensiones abelianas de ${\ma Q}$.

Sea $n\in{\ma N}$, $n=p_1^{\alpha_1}\cdots p_r^{\alpha_r}$
su descomposici\'on en primos. Entonces $U_n
\stackrel{\phi}{\cong} U_{p_1^{
\alpha_1}}\times \cdots\times  U_{p_r^{\alpha_r}}$. Sea $\chi\in
\hat{U}_n$, $\chi\colon U_n\to {\ma C}^{\ast}$. Sea $
g_{p_i}\colon U_{p_i^{\alpha_i}} \to U_{p_1^{\alpha_1}}\times 
\cdots \times
 U_{p_r^{\alpha_r}}$ dado por $g_{p_i}(x)=(1,\ldots,1,
 \equis\limits_{\substack{\uparrow\\ i}},1\cdots, 1)$.

El isomorfismo $\phi\colon U_n\to U_{p_1^{
\alpha_1}}\times \cdots\times  U_{p_r^{\alpha_r}}$ est\'a dado por 
el Teorema Chino del Residuo: $\phi(x\bmod n)=(x\bmod
p_1^{\alpha_1},\ldots, x\bmod p_r^{\alpha_r})$. Se tiene el
diagrama conmutativo
\begin{gather*}
\xymatrix{
U_n\ar[r]^{\hspace{-1cm}\phi}\ar[rd]_{\chi}& U_{p_1^{
\alpha_1}}\times \cdots\times  U_{p_r^{\alpha_r}}&
U_{p_i^{\alpha_i}}\ar[l]_{\hspace{1cm}g_{p_i}}\ar[dl]^{\chi_{p_i}}\\
&{\ma C}^{\ast}
}\\
\intertext{donde $\chi_{p_i}\colon U_{p_i^{\alpha_i}}\to {\ma C}^{\ast}$
est\'a dado por:}
\chi_{p_i}:= \chi\circ \phi^{-1}\circ g_{p_i}.
\end{gather*}

Si $a\in{\ma Z}$ es primo relativo a $n$, entonces
\begin{gather*}
\chi_{p_i}(a)=\chi(\phi(g_{p_i}(a\bmod p_i^{\alpha_i})))=
\chi(\phi^{-1}(1,\ldots,1,a,1\ldots,1))=\chi(b_i)\\
\intertext{donde $b_i\in {\ma Z}$ satisface:}
\begin{align*}
b_i&\equiv 1\bmod p_j^{\alpha_j}, \quad j=1,\ldots, r, \quad j\neq i,\\
b_i&\equiv a_i \bmod p_i^{\alpha_i}.
\end{align*}
\end{gather*}

Notemos que $b_1\cdots b_i\cdots b_r\equiv 1\cdots a\cdots 1\equiv
a\bmod p_i^{\alpha_i}$ para toda $1\leq i\leq r$, lo cual implica que
$b_1\cdots b_r\equiv a\bmod n$. Por lo tanto
\[
\chi(a)=\chi(b_1\cdots b_r)=\chi(b_1)\cdots\chi(b_r)=\chi_{p_1}(a)
\cdots \chi_{p_r}(a)=(\chi_{p_1}\cdots\chi_{p_r})(a).
\]
Se sigue que $\chi=\chi_{p_1}\chi_{p_2}\cdots \chi_{p_r}$.

Por el Teorema \ref{T12.2.23} se tiene que ${\eu f}_\chi=
\prod_{i=1}^r{\eu f}_{\chi_{p_i}}$. En particular $p\mid {\eu f}_{\chi}$
si y solamente si $\chi_{p_i}$ es no trivial.

\begin{definicion}\label{D12.3.1} Con la notaci\'on anterior, si
$X$ es un grupo de caracteres m\'odulo $n$ y si $p$ es un
n\'umero primo, entonces se define
\[
X_p:=\{\chi_p\mid \chi\in X\}.
\]
Notemos que $X_p=\{1\}$ si $p\notin\{p_1,\ldots, p_r\}$.
\end{definicion}

\begin{ejemplo}\label{Ej12.3.2}
Sea $\chi\colon U_{12}\to{\ma C}^{\ast}$ el caracter cuadr\'atico
par, es decir $\chi(1)=\chi(11)=1$, $\chi(5)=\chi(7)=-1$.

Se tiene que $12=2^2\cdot 3$. Entonces $\chi=\chi_2\chi_3$, 
donde $\chi_2\colon U_4\to{\ma C}^{\ast}$, $\chi_3\colon U_3
\to {\ma C}^{\ast}$. Calculemos $\chi_2$ y $\chi_3$. Se tiene
\begin{gather*}
\begin{array}{ccccc}
U_4&\stackrel{g_2}{\longto}&U_4\times U_3&\stackrel{\phi^{-1}}{
\longto}&U_{12}\\
1&\longmapsto&(1,1)&\longmapsto& 1\\
3&\longmapsto& (3,1)&\longmapsto& 7
\end{array}
\intertext{pues $7\equiv 3\bmod 4$, $7\equiv 1\bmod 3$,}
\begin{array}{ccccc}
U_3&\stackrel{g_3}{\longto}&U_4\times U_3&\stackrel{\phi^{-1}}{
\longto}&U_{12}\\
1&\longmapsto&(1,1)&\longmapsto& 1\\
2&\longmapsto& (1,2)&\longmapsto& 5
\end{array}
\end{gather*}
pues $5\equiv 1\bmod 4$, $5\equiv 2\bmod 3$.

Por tanto
\begin{alignat*}{2}
\chi_2=\chi\circ \phi^{-1}\circ g_2,& \quad\chi_2(1)=1,&\quad
\chi_2(3)=\chi(7)=-1,\\
\chi_3=\chi\circ \phi^{-1}\circ g_3,& \quad\chi_3(1)=1,&\quad
\chi_3(3)=\chi(5)=-1.
\end{alignat*}

Entonces $X:=\langle \chi\rangle$, $X_2:=\langle \chi_2\rangle$,
$X_3:=\langle \chi_3\rangle$ y si $p$ es cualquier primo,
$p\neq 2,3$, $X_p=\{1\}$.
\end{ejemplo}

El resultado m\'as importante, es la relaci\'on entre $X_p$ y
el {\'\i}ndice de ramificaci\'on de $p$ en $K/{\ma Q}$, donde
$K$ es el campo asociado a $X$.

\begin{teorema}\label{T12.3.3} Sea $X$ un grupo de caracteres
de Dirichlet y sea $K$ el campo asociado a $X$. Sea $p$
un n\'umero primo y sea $e$ el {\'\i}ndice de ramificaci\'on de $p$
en $K/{\ma Q}$. Entonces $e=|X_p|$.
\end{teorema}

\begin{proof}
Sea $n:=\mcm\{\f{\chi}\mid\chi\in X\}$. Entonces $K\subseteq \cic n{}$.
Escribamos $n=p^a m$ con $\mcd (m,p)=1$. Definimos
$L:=K(\zeta_m)=K\cic m{}\subseteq \cic n{}$. Si $Y$ es el grupo
de caracteres m\'odulo $n$ asociado al campo $L$, entonces por
el Teorema \ref{T12.2.36} se tiene que $L=\cic n{}^{Y^{\perp}}$.

Ahora bien, puesto que $L=K\cic m{}$, por la Proposici\'on
\ref {P12.2.38} se tiene que el grupo de caracteres $Y$ est\'a
generado por $X$ y por los caracteres de $\cic n{}$ de conductor
un divisor de $m$, los cuales son precisamente $\hat{U}_m$.

Si $\varphi\in Y$, entonces $\varphi=\chi\psi$ con $\chi\in X$ y
$\psi\in \hat{U}_m$. Escribamos $\chi=\chi_p\chi'$ con 
$\chi'=\prod_{q|m}\chi_q\in\hat{U}_m$. Por tanto
$\varphi=\chi_p (\chi'\varphi)\in X_p\times \hat{U}_m$. En particular
$Y\subseteq X_p\times \hat{U}_m$.

Rec{\'\i}procamente, si $\xi\varphi\in X_p\times\hat{U}_m$, puesto
que $\xi\in X_p$, existe $\chi\in X$ tal que $\chi_p=\xi$, es decir,
$\chi=\xi \cdot \prod_{q|m}\chi_q=\xi\chi'$. Por lo tanto
\[
\xi\varphi=\chi_p\varphi=\chi_p\chi'((\chi')^{-1}\varphi)\in
\langle X,\hat{U}_m\rangle = Y.
\]
Se sigue que $Y=X_p\times\hat{U}_m$.

Nuevamente por la Proposici\'on \ref{P12.2.38}, $L$ se escribe como
$L=\cic m{} F$ donde $F\subseteq \cic n{}$ es el
campo perteneciente a $X_p$. Notemos que $F\subseteq \cic pa$
pues $X_p\subseteq \hat{U}_{p^a}$.

Tenemos el siguiente diagrama donde la ramificaci\'on indicada
se refiere a $p$:
\[
\xymatrix{
&&\cic n{}\ar@{-}[d]\ar@{-}@/_/[ddll]\ar@{-}@/^/[rrdd]\\
&&{L=K(\zeta_m)=F(\zeta_m)}=K\cic m{}\ar@{-}[ddl]_{\substack{\hbox{\rm\tiny
no}\\\hbox{\rm\tiny ramificado}}}\ar@{-}[d]^{\substack{\hbox{\rm\tiny{no
}}\\\hbox{\rm\tiny
ramificado}}}\ar@{-}[drr]\\
{\cic pa}\ar@{-}[dr]&&{K}\ar@{-}[dd]^{e=e_p(K/{\ma Q})}&&
{\cic m{}}\ar@{-}[ddll]^{\substack{\hbox{\rm
\tiny no}\\ \hbox{\rm\tiny ramificado}}}\\
&{F}\ar@{-}[dr]_{\substack{\hbox{\rm\tiny totalmente}\\ \hbox{\rm\tiny
ramificado}}}\\&&{{\ma Q}}}
\]

Entonces el {\'\i}ndice de ramificaci\'on $e$ est\'a dado por
\[
e=e_p(K/{\ma Q})=e_p(L/{\ma Q})=e_p(F/{\ma Q})=[F:{\ma Q}]
=|X_p|. \tag*{$\fin$}
\]

\end{proof}

Como consecuencia de este resultado, tenemos:

\begin{corolario}\label{C12.3.4}
Sea $\chi$ un caracter de Dirichlet y sea $K$ el campo asociado
a $X$. Entonces el n\'umero primo $p$ se ramifica en $K/{\ma Q}$
si y solamente si $\chi(p)=0$, es decir, si y solamente si
$p|\f{\chi}$.

M\'as generalmente, si $X$ es un grupo de caracteres
de Dirichlet y $L$ es el campo asociado a $X$, entonces $p$ se
ramifica en $L/{\ma Q}$ si y s\'olo si existe $\chi\in X$ tal que 
$\chi(p)=0$, es decir, si y solamente si $p\mid {\eu f}_{\chi}$ para
alg\'un $\chi\in X$.
\end{corolario}

\begin{proof}
Se tiene que $p$ se ramifica en $L/{\ma Q}$ si y s\'olo si $X_p
\neq \{1\}$, lo cual equivale a que existe $\chi\in X$ tal que
$\chi_p\neq 1$. Por tanto $p$ se ramifica en $L/{\ma Q}$ si y
s\'olo si existe $\chi\in X$ tal que $p|\f{\chi} \iff$ existe
$\chi\in X$ tal que $\chi(p)=0$. $\fin$
\end{proof}

El Teorema \ref{T12.3.3} se puede refinar. Se tiene

\begin{teorema}\label{T12.3.5}
Sea $X$ un grupo de caracteres de Dirichlet y sea $K$ su
campo asociado. Sea $p$ un n\'umero primo. Sean $Y=
\{\chi\in X\mid \chi(p)\neq 0\}$ y $Z=\{\chi\in X\mid \chi(p)=1\}$. 
Entonces con las notaciones usuales, tenemos
\[
e=[X:Y],\qquad f=[Y:Z],\qquad g=[Z:1]=|Z|.
\]
M\'as a\'un, $X/Y\cong \widehat{I(\pK|p)}$, $X/Z\cong
\widehat{D(\pK|p)}$, donde $I(\pK|p)$ y $D(\pK|p)$ denotan a los
grupos de inercia y descomposici\'on respectivamente de los
primos en $K$ que dividen a $p$.

Finalmente, el grupo de Galois de los campos residuales 
satisface $Y/Z\cong \widehat{\Gal({\ma F}_{p^f}/{\ma F}_p)}$.

En particular, $q\in{\ma Q}$ se descompone totalmente en $K/{\ma Q}
\iff Z=X\iff \chi(q)=1\ \forall\ \chi\in X$.
\end{teorema}

\begin{proof}
Sea $L$ el subcampo de $K$ que corresponde a $Y$. Por el
Corolario \ref{C12.3.4} se tiene que $L$ es el m\'aximo
subcampo de $K$ en donde $p$ es no ramificado. Por tanto,
$L$ es el campo fijo del grupo de inercia $I:=I(\pK|p)$.
\begin{window}[0,l,\xymatrix{
K\ar@{-}[d]^{\Big\} I=\Gal(K/L)}\\L\ar@{-}[d]^{p \text{\ no
ramificado}}\\ {\ma Q}},{}]
Se tiene que $L=K^{Y^{\perp}}$, $Y=\Gal(K/L)^{\perp}$. Por
tanto $I\cong Y^{\perp}\cong \big(\widehat{X/Y}\big)$. As{\'\i}
$X/Y\cong \widehat{\Gal(K/L)}=\hat{I}$. En particular
$e=|I|=\big|\hat{I}\big|=|X/Y|=[X:Y]$ y $X/Y\cong \widehat{
\Gal(K/L)}=\widehat{I(p)}$. Ahora bien, $Y\cong
\widehat{\Gal(L/{\ma Q})}$. Sea $n:=\mcm\{\f{\chi}\mid \chi\in Y\}$.
Puesto que $p$ es no ramificado en $L$, $p\nmid \f{\chi}$ para
toda $\chi\in Y$ y por tanto $p\nmid n$. Se tiene $L\subseteq
\cic n{}$. El automorfismo de Frobenius de $p$ en $\cic n{}$ es 
el automorfismo $\sigma_p\colon \zeta_n\to\zeta_n^p$. Por lo
tanto el automorfismo de Frobenius de $p$  en $L$ es
\end{window}
\[
\sigma_p \bmod \Gal(\cic n{}/L)=\overline{\sigma_p}=
\sigma_p\bmod H
\]
donde $H:=\Gal(\cic n{}/L)$.

Con la identificaci\'on $\Gal(\cic n{}/{\ma Q})\cong U_n$,
 tenemos que $\overline{\sigma_p}=p\bmod H$ donde consideramos
$H\subseteq U_n$.

Si $\chi\in Y$, entonces $\chi(\Gal(\cic n{}/L))=1$, es decir,
$\chi(H)=1$ o, lo que es lo mismo, $H\subseteq \ker \chi$. Se sigue
que $\chi(\overline{\sigma_p})=\chi(\sigma_p)$ y por lo tanto
$\chi(\overline{\sigma_p})=1\iff \chi(p)=1$.

De lo anterior obtenemos que
\begin{gather*}
\langle \overline{\sigma_p}\rangle^{\perp}=\{\chi\in Y\mid
\chi(p)=1\}=Z\\
\intertext{en el mapeo bilineal}
\Gal(L/{\ma Q})\times Y\longto {\ma C}^{\ast}.
\end{gather*}

Ahora bien, $\langle\overline{\sigma_p}\rangle$ es un grupo 
c{\'\i}clico de orden $f$ generado por $\overline{\sigma_p}$.
Se sigue que
\[
\frac{\widehat{\Gal(L/{\ma Q})}}{\langle\overline{\sigma_p}
\rangle^{\perp}}=\frac{Y}{\langle\overline{\sigma_p}
\rangle^{\perp}}\cong \frac{Y}{Z}\cong \widehat{
\langle\overline{\sigma_p}\rangle}\cong
\langle\overline{\sigma_p}\rangle,
\]
y $Y/Z\cong\widehat{\langle \overline{\sigma_p}\rangle}
=\widehat{D_{L/{\ma Q}}(p)}$.
Por lo tanto $[Y:Z]=f=o(\overline{\sigma_p})$.

Se tiene el diagrama
\[
\left.
\begin{array}{lcl}
&K\\ &\bigg| &\Bigg\} {\ }_{X/Y\to e\to \text{\ grupo de inercia}}\\& L\\
_{\langle\overline{\sigma_p}\rangle}&\bigg|&\Bigg\}
{\ }_{Y/Z\to f \text{\ inercia}}\\
&E=L^{\langle\overline{\sigma_p}\rangle}&\to {\ }_{\Gal(L/E)^{\perp}
=\langle \overline{\sigma_p}\rangle^{\perp}=Z}\\
&\bigg|& \Bigg\} {\ }_{Z\to g \text{\ descomposici\'on}}\\
&{\ma Q}\end{array}\right\} {\ }_X
\begin{array}{c}
\left.\phantom{\begin{array}{c}\xymatrix{K\ar[d]\\K\ar[d]\\K} \end{array}}\right\}
 {\ }_{\substack{X/Z\to ef\to\\ \to \text{\ grupo de descomposici\'on}}}
{\ }\\{\ }\\{\ }\end{array}
\]

El campo fijo del automorfismo de Frobenius $E$ corresponde
al campo de descomposici\'on de $p$. Por tanto $E$ corresponde
a $Z$ y $g=[E:{\ma Q}]=|Z|$ o, simplemente, 
\[
efg=[K:{\ma Q}]=|X|=[X:Y][Y:Z][Z:1]=ef[Z:1],
\]
por lo tanto $g=[Z:1]=|Z|$. Se sigue que $X/Z\cong
\widehat{D(\pK|p)}$. $\fin$
\end{proof}

\begin{corolario}\label{C12.3.6'}
Sea $K/{\ma Q}$ una extensi\'on abeliana finita y sea $X$ el grupo
de caracteres de Dirichlet asociado a $K$. Sea $p$ un primo
racional. Entonces $p$ se descompone totalmente en $K/{\ma Q}$
si y solamente si $Z=\{\chi\in X\mid \chi(p)=1\}=X$.
Esto es, $p$ se descompone totalmente en $K/{\ma Q}\iff
\chi(p)=1$ para toda $\chi\in X$.

En particular si $K/{\ma Q}$ es una extensi\'on c\'iclica y $X$
es generado por $\chi$, entonces $p$ se descompone totalmente
en $K/{\ma Q}$ si y solamente si $\chi(p)=1$.
$\fin$
\end{corolario}

\subsection{F\'ormula del conductor--discriminante}\label{S12.3.1}

Nuestro objetivo en esta secci\'on es probar que
\[
|\delta_K|=\prod_{\chi\in X}\f{\chi}
\]
donde $K/{\ma Q}$ es una extensi\'on abeliana finita  y $X$ es
el grupo de caracteres de Dirichlet asociado a $K$.

Primero consideremos un
subcampo $F\subseteq \cic pn$, $p$ con $p$ un n\'umero primo y
$n\in{\ma N}$. Para cualquier extensi\'on finita $K/{\ma Q}$.
denotamos por
$\epsilon_K:=|\delta_K|$ al valor absoluto del discriminante.
Se tiene $\delta_K=(-1)^{r_2}\epsilon_K$.

Sea $p^a:=\mcd\{\f{\varphi}\mid \varphi\in X\}$ donde $X$ es
el grupo de caracteres de Dirichlet asociado a $F$. Entonces
$F\subseteq \cic pa$, $F\nsubseteq \cic p{a-1}$ y $X\cong
\widehat{\Gal(F/{\ma Q})}$.

Sea $\pL=\langle 1-\zeta_{p^a}\rangle$ el \'unico ideal en
${\mathcal O}_{\cic pa}={\ma Z}[\zeta_{p^a}]$
sobre $p$ y sea $\pK:=\pL\cap {\cal  O}_F$.

Empezamos analizando los grupos de ramificaci\'on de 
$\cic pa$. Sea $G:=\Gal(\cic pa/{\ma Q})$. Se tiene
$(p){\cal O}_{\cic pa}={\eu P}^{\varphi(p^a)}=
(\zeta_{p^a}-1)^{\varphi(p^a)}$. Sea $\sigma\in G$, $\sigma\neq 1$,
dado por $\sigma(\zeta_{p^a})=\zeta_{p^a}^{a_{\sigma}}$,
$a_{\sigma}\in{\ma Z}$, $1\leq a_{\sigma}\leq p^a-1$,
$\mcd(a_{\sigma},p)=1$. Sea $a_{\sigma}-1=p^{\alpha_{\sigma}}
\ell_{\sigma}$ con $\mcd(\ell_{\sigma},p)=1$, $0\leq \alpha_{\sigma}
\leq a-1$. Entonces
\begin{align*}
i_G(\sigma):& = v_{{\eu P}}(\sigma(\zeta_{p^a})-\zeta_{p^a})=
v_{{\eu P}}(\zeta_{p^a}^{a_{\sigma}}-\zeta_{p^a})=
v_{\eu P}(\zeta_{p^a}(\zeta_{p^a}^{a_{\sigma}-1}-1))\\
&=v_{\eu P}(\zeta_{p^a}^{p^{\alpha_{\sigma}}\ell_{\sigma}}-1)=
v_{\eu P}(\zeta_{p^{a-\alpha_{\sigma}}}^{\ell_{\sigma}}-1)
= v_{\eu P}(\zeta_{p^{a-\alpha_{\sigma}}}-1)\\
&=v_{\eu P}((\zeta_{p^a}-1)^{p^{\alpha_{\sigma}}})=p^{\alpha_{\sigma}}.
\end{align*}
Es decir,
\begin{gather}\label{Ec12.3.1.1'}
i_G(\sigma)=p^{\alpha_{\sigma}}.
\end{gather}

Se tiene que $\sigma\in G_u \iff
v_{\eu P}(\sigma(\zeta_{p^a})-\zeta_{p^a})=p^{\alpha_{\sigma}}
\geq u+1 \iff u\leq p^{\alpha_{\sigma}}-1$. Se sigue que
\begin{gather}\label{Ec12.3.1.2'}
G_u=\{\sigma\in G\mid \sigma(\zeta_{p^a})=\zeta_{p^a}^{
a_{\sigma}}, v_p(a_{\sigma}-1)=\alpha_{\sigma}, u\leq
p^{\alpha_{\sigma}}-1\}.
\end{gather}

De (\ref{Ec12.3.1.2'}) y
recordando que $D_{p^a,p^m}=\Gal(\cic pa/\cic pm)$, 
$1\leq m\leq a$, $D_{p^a,p^0}=D_{p^n,1}=G$, se tiene que
\begin{alignat*}{2}
G_{-1}=G_0&=G,\\
G_u&\cong D_{p^a,p},&&\quad 1\leq u\leq p-1,\\
G_u&\cong D_{p^a,p^2}, &&\quad p\leq u\leq p^2-1,\\
\vdots &\qquad \vdots&&\quad \vdots\\
G_u&\cong D_{p^a,p^{a-1}}, &&\quad p^{a-2}\leq p^{a-1}-1,\\
G_u&=\{1\},&&\quad p^{a-1}\leq u.
\end{alignat*}

Como consecuencia del Teorema \ref{T1.3.11} se tiene
que se ${\eu D}_{\cic pa/{\ma Q}}={\eu P}^s$, entonces
\begin{align}
s&=\sum_{j=0}^{\infty}(|G_j|-1)=(|G|-1)+\sum_{j=1}^{a-1}
(p^j-p^{j-1})(|D_{p^a,p^j}|-1)\nonumber\\
&=[\cic pa:{\ma Q}]-1+\sum_{j=1}^{a-1}[\cic pj:{\ma Q}](
[\cic pa:\cic pj]-1)\nonumber\\
&=a[\cic pa:{\ma Q}]-\sum_{j=0}^{a-1}
[\cic pj:{\ma Q}].\label{Ec6.3.2'}
\end{align}

Notemos que 
\begin{align*}
s&=a\varphi(p^a)-\sum_{j=0}^{a-1}\varphi(p^j)= a(p^{a-1}(p-1))-
\sum_{j=1}^{a-1}(p^j-p^{j-1})-1\\
&=ap^a-ap^{a-1}-p^{a-1}+1-1=p^{a-1}(ap-a-1),
\end{align*}
lo cual nos da una nueva demostraci\'on del Corolario \ref{C1.2.1.8}.

Para el caso en que ${\ma Q}\subseteq F\subseteq \cic pa$,
sea $F_j:=F\cap \cic pj$, $0\leq j\leq a$. Sea $H:=\Gal(\cic pa/F)$ y 
sea ${\eu D}_{\cic pa/F}={\eu P}^t$. Se tiene que
$G_j\cong D_{p^a,p^{r_j}}$ para alguna $r_j$. Por 
la Proposici\'on \ref{P1.3.11'} se tiene que $H_j=
G_j\cap H=\Gal(\cic pa/F\cic p{r_j})$. Por tanto, similarmente
a (\ref{Ec6.3.2'}), como $[\cic pi:{\ma Q}]$ es el n\'umero
se sumandos $(|H_i|-1)$, se tiene
\begin{align}
t&=\sum_{\sigma\in H\setminus\{1\}}i_H(\sigma)=
\sum_{\sigma \in H\setminus\{1\}}i_G(\sigma)=
\sum_{j=0}^{\infty}(|H_j|-1)=\sum_{j=0}^{\infty}(|G_j\cap H|-1)\nonumber\\
&=([\cic pa:F]-1)+ \sum_{i=1}^{a-1}[\cic pi:{\ma Q}](
[\cic pa : F\cic pi]-1).\label{Ec6.3.2''}
\end{align}

Puesto que 
\begin{gather}\label{Ec12.3.1.3'}
{\eu D}_{{\ma Q}(\zeta_{p^a})/{\ma Q}} = 
{\eu D}_{{\ma Q}(\zeta_{p^a})/F}\cdot \con_{F/{\ma Q}(\zeta_{p^a})}
{\eu D}_{F/{\ma Q}},
\end{gather}
y ${\eu p}$ es totalmente ramificado en $\cic pa/F$, se tiene que
si ${\eu D}_{F/{\ma Q}}={\eu p}^r$, entonces 
de (\ref{Ec12.3.1.3'}) se obtiene que $\pL^s=\pL^t \pL^{r[\cic pa:F]}$,
esto es, por (\ref{Ec6.3.2'}) y (\ref{Ec6.3.2''}), se sigue que
\begin{align*}
r&=\frac{s-t}{[\cic pa:F]}=\frac{1}{[\cic pa:F]}\Big\{
([\cic pa:{\ma Q}]-[\cic pa:F])\\
&\qquad +\sum_{j=1}^{a-1}
[\cic pj:{\ma Q}]\big([\cic pa:\cic pj]-[\cic pa:F\cic pj]\big)\Big\}\\
&=([F:{\ma Q}]-1)+\frac{1}{[\cic pa :F]}\\
&\qquad \Big\{
\sum_{j=1}^{a-1}\big([\cic pa:{\ma Q}]-[\cic pa:F\cic pj]\big)\cdot
[\cic pj:{\ma Q}]\big)\Big\}.
\end{align*}
Se sigue que
\begin{align*}
r&=[F:{\ma Q}]-1+\sum_{j=1}^{a-1}\Big([F:{\ma Q}]-
\frac{[\cic pj:{\ma Q}]}{[F\cic pj:F]}\Big)\\
&=
a[F:{\ma Q}]-\sum_{j=0}^{a-1}\frac{[\cic pj:{\ma Q}]}
{[F\cic pj:F]}.
\end{align*}
Ahora bien, se tiene
\begin{gather*}
\xymatrix{
\cic pj\ar@{-}[rr]\ar@{-}[d]&&{F\cic pj}\ar@{-}[d]\\
F_j=\cic pj \cap F\ar@{-}[rr]&&F}
\qquad [F\cic pj:F]=]\cic pj:F_j],\\
\intertext{por lo que}
\frac{[\cic pj:{\ma Q}]}{[F\cic pj:F]}=\frac{[\cic pj:{\ma Q}]}
{[\cic pj :F_j]}=[F_j:{\ma Q}].
\end{gather*}

Se sigue que
\begin{equation}\label{Ec12.3.1.4'}
r=a[F:{\ma Q}]-\sum_{j=0}^{a-1}[F_j:{\ma Q}].
\end{equation}

Hemos probado

\begin{proposition}\label{P12.3.1.5'}
Sea $Q\subseteq F\subseteq \cic pa$, $p\geq 2$ un n\'umero
primo y $a\geq 1$. Entonces ${\eu D}_{F:{\ma Q}}=
{\eu p}^r$, donde $r=a[F:{\ma Q}]-\sum_{j=0}^{a-1}[F_j:{\ma Q}]$.
$\fin$
\end{proposition}

\begin{corolario}\label{C12.3.1.6'} Con las condiciones anteriores,
$|{\eu d}_{F/{\ma Q}}|=p^r$. $\fin$
\end{corolario}

Ahora bien, sea $X$ el grupo de caracteres de Dirichlet asociado
a $F$. Se tiene que $F_a=F$ y $F_0={\ma Q}$. Un caracter
$\chi$ tiene conductor $p^j$ si y solamente si $\chi$ es un
caracter asociado a $\cic pj$ pero no asociado a $\cic p{j-1}$.
Por lo tanto $X$ contiene precisamente $[F_j:{\ma Q}]-
[F_{j-1}:{\ma Q}]$ caracteres de conductor $p^j$, $1\leq j\leq a$.
Se sigue que $\prod_{\chi\in X}{\eu f}_{\chi}=p^{\alpha}$ donde
\begin{equation}\label{Ec12.3.1.7'}
\alpha=\sum_{j=1}^a j([F_j:{\ma Q}]-[F_{j-1}:{\ma Q}])=
n[F:{\ma Q}]-\sum_{j=0}^{a-1}[F_j:{\ma Q}].
\end{equation}

De (\ref{Ec12.3.1.4'}) y (\ref{Ec12.3.1.7'}) se sigue

\begin{proposicion}\label{P12.3.1.7}
Si $F\subseteq {\ma Q}(\zeta_{p^a})$ con $p$ un
n\'umero primo, y $X$ es el grupo de caracteres
de Dirichlet asociado a
$F$, entonces 
\[
\epsilon_F=\prod_{\chi\in 
X}{\eu f}_{\chi}. \tag*{$\fin$}
\]
\end{proposicion}

\begin{teorema}[F\'ormula del
conductor--discriminante]\label{T12.3.1.8}
Sea $K/{\ma Q}$ una extensi\'on abeliana finita. Entonces
\[
\delta_K=(-1)^{r_2}\prod_{\chi\in X}{\eu f}_{\chi}
\]
 donde $X$
denota el grupo de caracteres de Dirichlet asociado a $K$.
\end{teorema}

\begin{proof}
Basta probar que $\epsilon_K=\prod_{\chi\in X}\f{\chi}$. Fijemos 
un n\'umero primo $p$. Con las notaciones del Teorema \ref{T12.3.3},
esto es, $L=K(\zeta_m)$, $n=p^{a}m$ y $\mcd(m,p)=1$,
tenemos que puesto que ning\'un primo encima de $p$ es
ramificado ni en $L/K$ ni en $L/F$ y puesto que
\begin{gather}\label{Ec12.3.1.9}
{\eu D}_{L/{\ma Q}}={\eu D}_{L/K}\cdot
\con_{K/L}{\eu D}_{K/{\ma Q}}=
{\eu D}_{L/F}\cdot \con_{F/L}{\eu D}_{F/{\ma Q}},
\end{gather}
se tiene que si para cualquier extensi\'on $M/{\ma Q}$,
${\eu d}_{M/{\ma Q}}(p)$ denota la potencia exacta de
$\langle p\rangle$ que divide a ${\eu d}_{M/{\ma Q}}$
(y similarmente para $\epsilon_M(p)$), entonces
de (\ref{Ec12.3.1.9})
\begin{gather*}
{\eu d}_{L/{\ma Q}}(p)=(N_{L/{\ma Q}}{\eu D}_{L/{\ma Q}})(p)=
{\eu d}_{K/{\ma Q}}^{[L:K]}(p) = {\eu d}_{F/{\ma Q}}^{[L:F]}(p).\\
\intertext{Por lo tanto}
\epsilon_K(p)=\epsilon_L^{(1/[L:K])}(p)=\epsilon_F^{([L:F]/[L:K])}(p).
\end{gather*}

Ahora bien $[L:F]=[{\ma Q}(\zeta_m):{\ma Q}]=\varphi(m)$ y
$[L:K]=[Y:X]=\frac{|Y|}{|X|}= \frac{|X_p||\hat{\mathcal U}_m|}{|X|}=
\frac{|X_p|\varphi(m)}{|X|}$. Se sigue que $\frac{[L:F]}{[L:K]}=\frac{
|X|}{|X_p|}$. Puesto que $F\subseteq {\ma Q}(\zeta_{p^a})$, 
obtenemos de la Proposici\'on \ref{P12.3.1.7} que
 $\epsilon_F=\prod_{\varphi\in X_p}
{\eu f}_{\varphi}$. Por tanto
\[
\epsilon_K(p)=\epsilon_F^{(|X|/|X_p|)}(p)=\Big(\prod_{\varphi\in X_p}
{\eu f}_{\varphi}\Big)^{|X|/|X_p|}.
\]
Del epimorfismo natural $\pi\colon X\to X_p$, $\chi\mapsto
\chi_p$, obtenemos que $|\ker \pi|=\frac{|X|}{|X_p|}$. Del
Teorema \ref{T12.3.3}, $|X_p|=e$ y $|X|=[K:{\ma Q}]=efg$ (con
las notaciones usuales). Entonces, cada
$\chi_p$ aparece para exactamente
$fg=\frac{|X|}{|X_p|}$ 
elementos diferentes de $X$, esto es, $|\pi^{-1}(\chi_p)|=fg$.
Por lo tanto
\[
\epsilon_K(p)=\big(\prod_{\varphi\in X_p}{\eu f}_{\varphi}\big)^{fg}
=\prod_{\chi\in X}{\eu f}_{\chi_p}.
\]
Puesto que $\epsilon_K=\prod_p \epsilon_K(p)$, tenemos
\[
\epsilon_K=\prod_p \epsilon_K(p)=\prod_p\prod_{\chi\in X}
{\eu f}_{\chi_p}=\prod_{\chi\in X}\prod_p {\eu f}_{\chi_p}.
\]

Finalmente, puesto que para cualesquiera dos primos diferentes
$p$ y $q$, ${\eu f}_{\chi_p}$
y ${\eu f}_{\chi_q}$ son primos relativos y
 $\chi=\prod_p \chi_p$, se sigue del Teorema \ref{T12.2.23} que
  ${\eu f}_{\chi}=\prod_p {\eu f}_{\chi_p}$, as{\'\i} que:
\begin{gather*}
\epsilon_K=\prod_{p}\epsilon_K(p)=
\prod_{\chi\in X}\prod_p {\eu f}_{\chi_p}=
\prod_{\chi \in X}{\eu f}_{\chi}. \tag*{$\fin$}
\end{gather*}
\end{proof}

\begin{ejemplo}\label{Ej12.3.1.8'}
Sea $K={\ma Q}(\zeta_p+\zeta_p^{-1})=\cic p{}^+$ para $p$ un
n\'umero primo impar. Se tiene $X=\{\chi\in\widehat{U_p}\mid
\chi(-1)=1\}$. Todos los caracteres de $X$, con excepci\'on del
trivial, tiene conductor $p$ y hay $|X|-1=\frac{p-1}{2}-1=\frac{
p-3}{2}$ de estos caracteres. Puesto que $r_2=0$ por ser $K$
un campo real, se tiene $\delta_K=p^{\frac{p-3}{2}}$.

Otra forma de obtener este resultado es que $K/{\ma Q}$ \'unicamente
se ramifica en $p$, que $p$ es total y moderadamente ramificado, por
lo que $\delta_K=p^{e-1}$ y $e=[K:{\ma Q}]=\frac{p-1}{2}$.

Similarmente, obtenemos que $\delta_{\cic p{}}=(-1)^{\frac{p-1}{2}}
p^{p-2}$.
\end{ejemplo}

\section{Construcci\'on de extensiones abelianas}\label{S12.4}

El Teorema \ref{T12.3.3} es muy \'util para construir extensiones
abelianas con caracter{\'\i}sticas especiales. Primero veamos como
funciona con un ejemplo espec{\'\i}fico.

\begin{ejemplo}\label{Ej12.4.1}
Consideremos el caracter cuadr\'atico $\chi\colon U_{12}\to
{\ma C}^{\ast}$ dado por $\chi(1)=\chi(11)=1$ y $\chi(5)=\chi(7)=-1$.
Entonces el caracter es real, $\f{\chi}=12$ y el campo asociado
tiene que ser $\cic {12}{}^+={\ma Q}(\zeta_{12}+\zeta_{12}^{-1})=
{\ma Q}(\sqrt{3})$.

Sea $\chi=\chi_2\chi_3$. Entonces $\f{\chi_2}=4$ y $\f{\chi_3}=3$.
Por tanto los campos asociados a $\chi_2$ y $\chi_3$ son
$\cic 4{}$ y $\cic 3{}$ respectivamente. Sea $Y=\langle\chi_2\rangle
\times \langle\chi_3\rangle$. Se tiene que el campo asociado
a $Y$ es $\cic 4{}\cic 3{}=\cic {12}{}$. En particular $Y=\widehat{
U_{12}}$. Se tiene
\[
\xymatrix{
&\cic {12}{}={\ma Q}(\sqrt{3},\sqrt{-3})\ar@{-}[dl]\ar@{-}[dr]\ar@{-}[d]\\
{\ma Q}(i)=\cic 4{}\ar@{-}[dr]&{\ma Q}(\sqrt{3})\ar@{-}[d]
&\cic 3{}={\ma Q}(\sqrt{-3})\ar@{-}[dl]\\ &{\ma Q}
}
\]

Notemos que la ramificaci\'on est\'a dada por:
\l
\item En $\cic 4{}/{\ma Q}$: $2$ y el primo infinito $\infty$.
\item En ${\ma Q}(\sqrt{3})/{\ma Q}$: $2$ y $3$.
\item En $\cic 3{}/{\ma Q}$: $3$ e $\infty$.
\item En $\cic {12}{}/\cic 4{}$: $3$.
\item En $\cic {12}{}/{\ma Q}(\sqrt{3})$: $\infty$.
\item En $\cic {12}{}/\cic 3{}$: $2$.
\end{list}

En particular $\cic {12}{}/{\ma Q}(\sqrt{3})$ es no ramificada
en ning\'un primo finito.
\end{ejemplo}

\subsection{Campos de g\'eneros\index{campos de g\'eneros}}\label{S12.4.0}

\begin{teorema}[Leopoldt \cite{Leo53}\index{teorema!del
g\'enero de Leopoldt}\index{teorema de Leopoldt  del
g\'enero}\index{Leopoldt!teorema del
g\'enero}]\label{T12.4.2}
Sea $K/{\ma Q}$ una extensi\'on abeliana finita. Se $L$ la
m\'axima extensi\'on abeliana de ${\ma Q}$ que es no ramificada
en ning\'un primo finito. Entonces el grupo de caracteres de
Dirichlet $Y$ correspondiente a $L$ es 
\[
Y=\prod_{p} X_p
\]
donde $X$ es el grupo de caracteres de Dirichlet correspondiente
a $K$.
\end{teorema}

\begin{proof}
Primero notemos que para todo n\'umero primo $p$, $Y_p=X_p$
y por tanto $|Y_p|=|X_p|$. Por el Teorema \ref{T12.3.3} se
sigue que $e_p(L/K)=1$ y $L/K$ es no ramificada en ning\'un
primo finito.

Ahora, sea $E/K$ una extensi\'on no ramificada en ning\'un
primo finito y $E/{\ma Q}$ abeliana. Sea $Z$ el grupo de
caracteres de Dirichlet asociado a $E$. Entonces $|X_p|=
|Z_p|$ y $Z\supseteq X$. Por lo tanto $Z_p=X_p$ y se sigue que
$Z\subseteq \prod_p Z_p=\prod_p X_p=Y$ y por lo tanto
$E\subseteq L$. $\fin$
\end{proof}

Supongamos que $L/K$ es ramificada en los primos infinitos.
Entonces $K\subseteq {\ma R}$ y $L\nsubseteq {\ma R}$. Sea
$Y^+:=\{\chi\in Y\mid \chi(-1)=1\}$. Puesto que para toda
$\chi\in X$, $\chi(-1)=1$, $X\subseteq Y^+$. Por otro lado $Y^+=
\ker \theta$ donde $\theta\colon Y\to \{\pm\}$, $\theta(\chi)=\chi(-1)$
y puesto que $L\nsubseteq {\ma R}$, existe $\chi\in Y$ con
$\chi(-1)=-1$, es decir, $\theta$ es una funci\'on suprayectiva. Se 
sigue que $Y/Y^+\cong \{\pm \}$ y en particular $|Y^+|=|Y|/2$.

\begin{corolario}\label{C12.4.3} Con las notaciones del Teorema
{\rm \ref{T12.4.2}} se tiene que si $L^+$ es el campo correspondiente
a $Y^+$, entonces 
\las
\item Si $K$ y $L$ son ambos reales o ambos imaginarios, $L$
es la m\'axima extensi\'on de $K$  y abeliana sobre ${\ma Q}$, 
no ramificada en todo primo incluyendo los primos infinitos.
\item Si $K$ es real y $L$ es imaginario, $L^+/K$ es la m\'axima
extensi\'on no ramificada en ning\'un primo incluyendo los
primos infinitos y abeliana sobre ${\ma Q}$. Ahora, $L/K$ es la
m\'axima extensi\'on no ramificada en ning\'un primo finito y
abeliana sobre ${\ma Q}$. $\fin$
\end{list}
\end{corolario}

\begin{ejemplo}\label{Ej12.4.4}
En el Ejemplo \ref{Ej12.4.1}, $K={\ma Q}(\sqrt{3})$ con $X=\{\chi\}$,
$\chi\colon U_{12}\to {\ma C}^{\ast}$, $\chi(1)=\chi(11)=1$, $\chi(5)=
\chi(7)=-1$, $Y=\langle\chi_2\rangle\times\langle\chi_3\rangle$ y
$Y=\widehat{U_{12}}$, $Y^+=X$. Por tanto toda extensi\'on
$L/K$, $L/{\ma Q}$ abeliana, es ramificada en alg\'un primo.
\end{ejemplo}

\begin{ejemplo}\label{Ej12.4.5}
Sea $p$ un primo impar y sea $K={\ma Q}\big(\sqrt{(-1)^{(p-
1)/2}p}\big)\subseteq \cic p{}$ la subextensi\'on cuadr\'atica de
$\cic p{}$. Sea $\chi\colon U_p\to {\ma C}^{\ast}$ el caracter
asociado a $K$, $X=\langle\chi\rangle$. Entonces $o(\chi)=2=
[K:{\ma Q}]$, $\chi^2=1$ y $\chi\neq 1$. El conductor de $\chi$
es $\f{\chi}=p$. Se tiene $\chi(U_p)=\{\pm 1\}$.
Ahora bien, $K=\cic p{}^{\ker \chi}$ con
\[
\ker \chi=\{\sigma\in\Gal(\cic p{}/{\ma Q})\mid \chi(\sigma)=1\}.
\]

Sea $q$ cualquier n\'umero primo, $q\neq p$. Se tiene que
$q$ se descompone en $K/{\ma Q}$ si y s\'olo si $|Z|=2$
en la notaci\'on del Teorema \ref{T12.3.5}, es decir, 
$Z=\{\varphi\in X\mid \varphi(q)=1\}$. Esto es, $q$ se descompone
en $K/{\ma Q}\iff \chi(q)=1$.

Por otro lado, si $q\equiv a^2\bmod p$ para alg\'un $a\in{\ma Z}$, 
entonces $\chi(q)=\chi(a)^2=1$ y puesto que
\[
|\ker \chi|=\frac{|U_p|}{2}=\frac{p-1}{2}=|\{t\in U_p\mid t\equiv a^2
\bmod p\}|
\]
se sigue que $\chi(q)=1\iff q\equiv a^2\bmod p\iff \xbinom{q}{p}=1$
donde $\xbinom{q}{p}$ es el s{\'\i}mbolo de Legendre. Por 
tanto $\chi(q)=\xbinom{q}{p}$.

En resumen, el campo 
${\ma Q}\big(\sqrt{(-1)^{(p-1)/2} p}\big)$ corresponde
al s{\'\i}mbolo de Legendre: $\chi(q)=\xbinom{q}{p}$.
\end{ejemplo}

\begin{ejemplo}\label{Ej12.4.6}
Con respecto al Ejemplo \ref{Ej12.4.5}, nos preguntamos ahora
cual es el caracter cuadr\'atico de Dirichlet correspondiente al campo
cuadr\'atico ${\ma Q}\big(\sqrt{(-1)^{(p+1)/2}p}\big)$ donde
$p$ es un n\'umero primo impar.

Recordemos que si $p\equiv 1\bmod 4$, entonces
\begin{gather*}
\sqrt{(-1)^{(p-1)/2}p}=\sqrt{p} \quad \text{y} 
\quad \sqrt{(-1)^{(p+1)/2}p}=\sqrt{-p}\\
\intertext{y si $p\equiv 3\bmod 4$ entonces}
\sqrt{(-1)^{(p-1)/2}p}= \sqrt{-p}\quad \text{y} \quad
\sqrt{(-1)^{(p+1)/2}p}=\sqrt{p}.
\end{gather*}
Tenemos el siguiente diagrama
\begin{gather*}
\xymatrix{
\cic p{}\ar@{-}[rr]^2\ar@{-}[dd]_{\frac{p-1}{2}}&&\cic {4p}{}
\ar@{-}[dd]^{\frac{p-1}{2}}\\ \\
{\ma Q}\big(\sqrt{(-1)^{(p-1)/2}p}\big)\ar@{-}[rr]^2\ar@{-}[dd]_2&&
{\begin{array}{c}{\ma Q}\big(\zeta_4,\sqrt{(-1)^{(p-1)/2}p}\big)=\\ = {\ma Q}(\sqrt{p},
\sqrt{-p})\end{array}}\ar@{-}[dd]^{2}\ar@{-}[dl]_2\\
&{\ma Q}\big(\sqrt{(-1)^{(p+1)/2}p}\big)\ar@{-}[dl]_2\\
{\ma Q}\ar@{-}[rr]_2&& \cic 4{}={\ma Q}(\sqrt{-1})
}
\end{gather*}

Sean 
\begin{alignat*}{2}
\chi\colon U_p&\longto {\ma C}^{\ast}, \quad &\chi(q)&=\xbinom{q}{p},\\
\varphi\colon U_4&\longto {\ma C}^{\ast}, \quad &
\varphi(-1)&=-1.
\end{alignat*}
Por tanto $\varphi(q)=\begin{cases}
1& \text{si $q\equiv 1\bmod 4$}\\-1&\text{si $q\equiv -1\bmod 4
\equiv 3\bmod 4$}\end{cases}= (-1)^{(q-1)/2}$. 

Se sigue que ${\ma Q}\big(\sqrt{(-1)^{(p+1)/2}p}\big)$ corresponde
a $\varphi\chi$ el cual est\'a definido por $(\varphi\chi)(q)=
(-1)^{(q-1)/2}\xbinom{q}{p}$.

Finalmente, notemos que $\f{\varphi\chi}=\f{\varphi}\f{\chi}=4p$ y
que $\epsilon_{{\ma Q}\big(\sqrt{(-1)^{(p+1)/2}p}\big)}=
\big|\delta_{{\ma Q}\big(\sqrt{(-1)^{(p+1)/2}p}\big)}\big|=4p$.
\end{ejemplo}

\begin{observacion}\label{O12.4.6'}
El primo $q\neq 2$ se descompone en ${\ma Q}(i)/{\ma Q}\iff
q\equiv 1\bmod 4$.
\end{observacion}

\begin{ejemplo}\label{Ej12.4.7}
Sea $K_1={\ma Q}(\sqrt{10})$. Entonces $K_1={\ma Q}(
\sqrt{10})\subseteq {\ma Q}(\sqrt{2},\sqrt{5})\subseteq {\ma Q}(
\zeta_8,\zeta_5)=\cic {40}{}$ pues $5\equiv 1\bmod 4$.
Adem\'as $\Gal(\cic {40}{}/{\ma Q})\cong U_{40}\cong U_8\times
U_5\cong (C_2\times C_2)\times C_4=G$.
\[
\xymatrix{
\cic 8{}\ar@{-}[d]\\ {\ma Q}\ar@{-}[r]&{\ma Q}(\sqrt{5})=\cic 5{}^+
\ar@{-}[r]&\cic 5{}
}
\]
$G$ tiene $7$ subgrupos de orden $2$ y por tanto $7$ grupos
cociente de {\'\i}ndice $2$. Se sigue que $\cic {40}{}$ tiene
$7$ subcampos cuadr\'aticos:
\[
\cic 4{},{\ma Q}(\sqrt{2}),{\ma Q}(\sqrt{-2}), {\ma Q}(\sqrt{5}),
{\ma Q}(\sqrt{-5}), {\ma Q}(\sqrt{10})\rm{\ y\ } {\ma Q}(\sqrt{-10}).
\]
Puesto que ${\ma Q}(\sqrt{10})\nsubseteq \cic 4{},\cic 8{},
\cic 5{},\cic {10}{}, \cic {20}{}$ se sigue que el caracter $\chi$
asociado a ${\ma Q}(\sqrt{10})$ tiene conductor $\f{\chi}=40$
(o simplemente porque $\delta_{{\ma Q}(\sqrt{10})}=40=\f{\chi}$).

Se sigue que $\chi=\chi_2\chi_5$, $\f{\chi_2}=8$, $\f{\chi_5}=5$.
Adem\'as $\chi(-1)=1$ por lo que $\chi_2(-1)=\chi_5(-1)=\pm 1$.
En caso de que $\chi_2(-1)=\chi_5(-1)=-1$ se tendr{\'\i}a que
${\ma Q}(\sqrt{-5})\subseteq \cic 5{}$ pues $\chi_5(-1)=-1$
significa que el subcampo cuadr\'atico de $\cic 5{}$ ser{\'\i}a 
complejo. Se sigue que $\chi_2(-1)=\chi_5(-1)=1$ y
$\chi_2^2=\chi_5^2=1$.
Por lo tanto ${\ma Q}(\sqrt{2})$ es el campo asociado a $\chi_2$
y ${\ma Q}(\sqrt{5})$ es el campo asociado a $\chi_5$.

Si $Y=\langle \chi_2\rangle\oplus \langle\chi_5\rangle$ entonces
el campo asociado a $Y$ es $L={\ma Q}(\sqrt{2},\sqrt{5})$ y
${\ma Q}(\sqrt{2},\sqrt{5})$ es la m\'axima extensi\'on abeliana
de ${\ma Q}$ no ramificada sobre ${\ma Q}(\sqrt{10})$ pues
$\infty$ es no ramificado debido a que ${\ma Q}(\sqrt{2},\sqrt{5})$
es un campo real.

Ahora bien, usando el campo de clase de Hilbert (ver
Teorema \ref{T9.3}), esto es,
si $H_K$ es la m\'axima extensi\'on abeliana de $K$ no ramificada
en ning\'un primo incluyendo el $\infty$, se tiene que $I_K\cong
\Gal(H_K/K)$ donde $I_K$ es el grupo de clases de $K$. Puesto
que $L\subseteq H_K$, $2=[L:K]|[H_K:K]=|I_K|=h_K$, es decir,
$2|h_K$ y $K$, m\'as precisamente, ${\cal O}_K={\ma Z}[
\sqrt{10}]$, no es de ideales principales.
\end{ejemplo}

\begin{ejemplo}\label{Ej12.4.8}
Sea ahora $K={\ma Q}(\sqrt{-5})$. Se tiene $5\equiv 1\bmod 4$
y por lo tanto ${\ma Q}(\sqrt{-5})\nsubseteq \cic 5{}$ y 
${\ma Q}(\sqrt{-5})\subseteq {\ma Q}(\sqrt{-1},\sqrt{5})\subseteq
{\ma Q}(\zeta_4,\zeta_5)=\cic {20}{}$ y se tiene $\Gal(\cic {20}{}/
{\ma Q})\cong U_{20}\cong U_4\times U_5\cong C_2\times C_4$.
En particular $\cic {20}{}$ tiene tres subcampos cuadr\'aticos,
a saber, $\cic 4{}={\ma Q}(\sqrt{-1})$, ${\ma Q}(\sqrt{5})$ y
${\ma Q}(\sqrt{-5})$. Puesto que $\delta_K=-20$, $\f{\chi}=20$
donde $\chi$ es el caracter asociado a $K$.

Sea $\chi=\chi_2\chi_5$ con $\f{\chi_2}=4$ y $\f{\chi_5}=5$.
Puesto que $\chi(-1)=-1$, $\chi_2(-1)\neq \chi_5(-1)$. M\'as 
precisamente, \'unicamente hay un caracter cuadr\'atico de 
conductor $4$, el correspondiente a $\cic 4{}$, y un \'unico
caracter cuadr\'atico m\'odulo $5$, el correspondiente
a ${\ma Q}(\sqrt{5})$. Por tanto se tiene $\chi_2(-1)=-1$ y
$\chi_5(-1)=1$. As{\'\i}, el campo asociado a $Y=\langle \chi_2
\rangle\oplus\langle\chi_5\rangle$ es ${\ma Q}(\sqrt{-1},\sqrt{5})$
el cual es la m\'axima extensi\'on abeliana de ${\ma Q}$ no
ramificada en ning\'un primo sobre $K$, incluyendo al primo infinito
pues $K$ es un campo complejo.

Similarmente a como en el Ejemplo \ref{Ej12.4.7}, se tiene que
$2|h_K$ y ${\cal O}={\ma Z}\big[\frac{1+\sqrt{-5}}{2}\big]$
no es de ideales principales.
\end{ejemplo}

\begin{ejemplo}\label{Ej12.4.9}
Sea ahora $K={\ma Q}(\sqrt{30})$. Se tiene $30=2\cdot 3\cdot 5$,
$\delta_K=4\cdot 30=120$ por lo que el conductor de $\chi$, el
caracter asociado a $K$, es igual a $\f{\chi}=120$.
Ahora bien
\begin{gather*}
{\ma Q}(\sqrt{30})\subseteq {\ma Q}(\sqrt{2},\sqrt{3},\sqrt{5})
\subseteq {\ma Q}(\zeta_8,\zeta_{12},\zeta_5)={\ma Q}(
\zeta_8,\zeta_3,\zeta_5)=\cic {120}{}.\\
\intertext{Adem\'as}
\Gal(\cic {120}{}/{\ma Q})\cong U_{120}\cong U_8\times U_3\times
U_5\cong (C_2\times C_2)\times (C_2)\times (C_4).
\end{gather*}

En particular $\cic {120}{}$ tiene $\frac{2^4-1}{2-1}=15$ subcampos
cuadr\'aticos. Sea $\chi=\chi_2\chi_3\chi_5$, $\f{\chi_2}=8$,
$\f{\chi_3}=3$, $\f{\chi_5}=5$ y $\chi_2^2=\chi_3^2=\chi_5^2=1$.
Puesto que \'unicamente existen un caracter cuadr\'atico de 
conductores $3$ y $5$ respectivamente, $\chi_3$ corresponde
a $\cic 3{}={\ma Q}(\sqrt{-3})$, y $\chi_5$ corresponde a ${\ma Q}(
\sqrt{5})$. Adem\'as $\chi_3(-1)=-1$ y $\chi_5(-1)=1$. Puesto
que $\chi(-1)=1$, se tiene que $\chi_2(-1)=-1$ y por tanto
$\chi_2$ corresponde a ${\ma Q}(\sqrt{-2})$.

Sea $Y=\langle\chi_2\rangle\oplus\langle\chi_3\rangle\oplus\langle
\chi_5\rangle$. Entonces el campo $L$ asociado a $Y$ es
$L={\ma Q}(\sqrt{-2},\sqrt{-3},\sqrt{5})$.

Notemos $K$ es real y $L$ es imaginario, por lo que los primos
infinitos de $K$ son ramificados en $L$. Sea $Y^+=\{\varphi\in Y\mid
\varphi(-1)=1\}$, $Y^+$ corresponde a $L^+=L\cap {\ma R}$.
 Entonces $Y^+=\langle\chi_2\chi_3\rangle\oplus
\langle\chi_5\rangle$ y $\chi_2\chi_3$ corresponde a
${\ma Q}(\sqrt{-2}\sqrt{-3})={\ma Q}(\sqrt{6})$. Por lo tanto
$L^+={\ma Q}(\sqrt{6},\sqrt{5})$ y $L^+$ es la m\'axima
extensi\'on abeliana de ${\ma Q}$ no ramificada en ning\'un primo
de $K$ incluyendo al infinito.

Por otro lado $L={\ma Q}(\sqrt{-2},\sqrt{-3},\sqrt{5})$ es la
m\'axima extensi\'on abeliana de ${\ma Q}$ no ramificada en
ning\'un primo finito de $K$.
\[
\xymatrix{
&L={\ma Q}(\sqrt{-2},\sqrt{-3},\sqrt{5})\ar@{-}[dl]_{\substack{
\hbox{\rm\tiny{no}}\\ \hbox{\rm\tiny{ramificada}}}}
\ar@{-}[d]_{\hbox{\rm\tiny{no}}}^{\hbox{\rm\tiny{ramificada}}}
\ar@{-}[dr]^{\substack{\hbox{\rm\tiny{ramificada solo}}\\
\hbox{\rm\tiny{en $\infty$}}}}\\
{\ma Q}(\sqrt{-2},\sqrt{-15})\ar@{-}[dr]_{\substack{\hbox{\rm\tiny{
ramificada}}\\ \hbox{\rm\tiny{solo en $\infty$}}}}
&{\ma Q}(\sqrt{-10},\sqrt{-3})\ar@{-}[d]_{
\hbox{\rm\tiny{ramificada}}}^{\hbox{\rm\tiny{solo en $\infty$}}}&
{\ma Q}(\sqrt{5},\sqrt{6})\ar@{-}[dl]^{\substack{
\hbox{\rm\tiny{no ramificada}}\\ \hbox{\rm\tiny{en ning\'un primo}}}}\\
&{\ma Q}(\sqrt{30})
}
\]
\end{ejemplo}

Los Ejemplos \ref{Ej12.4.7}, \ref{Ej12.4.8} y \ref{Ej12.4.9} nos dan
la gu{\'\i}a del caso general que a continuaci\'on estudiamos.

\begin{ejemplo}\label{Ej12.4.10}
Sean $d\in{\ma Z}$ libre de cuadrados, $K={\ma Q}(\sqrt{d})$.
Escribamos $d=(-1)^{\varepsilon}2^{\delta} p_1\cdots p_s q_1
\cdots q_t$ donde $\varepsilon,\delta\in\{0,1\}$, $p_1,\ldots,p_s$
primos distintos congruentes con $1$ m\'odulo $4$ y $q_1,
\ldots, q_t$ primos distintos congruentes con $3$ m\'odulo $4$.

Sea $\chi$ el caracter cuadr\'atico asociado a $K$, 
$\chi=\chi_d$. Se tiene que
\[
\f{\chi}=\big|\delta_K\big|=\begin{cases}
|d|&\text{si $d\equiv 1\bmod 4$}\\ 4|d|&
\text{si $d\equiv 2,3 \bmod 4$}
\end{cases}.
\]

Por el Ejemplo \ref{Ej12.4.5} se tiene que $\chi_{p_i}$, $\chi_{q_j}$
corresponden al s{\'\i}mbolo de Legendre: $\chi_{p_i}(\ell)=
\xbinom{\ell}{p_i}$, $\chi_{q_j}(\ell)=\xbinom{\ell}{q_j}$,
$1\leq i\leq s$, $1\leq j\leq t$. Adem\'as
$\chi_{p_i}(-1)=1$, $i=1,2,\ldots, s$ y $\chi_{q_j}(-1)=-1$,
$j=1,2,\dots, t$. Por otro lado $\chi_{p_i}$ corresponde al campo
${\ma Q}\big(\sqrt{(-1)^{(p_i-1)/2}p_i}\big)={\ma Q}(\sqrt{p_i})$,
$1\leq i\leq s$ y $\chi_{q_j}$ corresponde al campo
${\ma Q}\big(\sqrt{(-1)^{(q_j-1)/2}q_j}\big)=
{\ma Q}(\sqrt{-q_j})$, $1\leq j
\leq t$.

M\'as a\'un $\chi(-1)=(-1)^{\varepsilon}$. El problema m\'as
complicado es ver que es $\chi_2$. Se tiene
\begin{gather*}
d\equiv (-1)^{\varepsilon}2^{\delta}(-1)^t\bmod 4\equiv (-1)^{
\varepsilon+t}2^{\delta}\bmod 4.\\
\intertext{Por tanto}
\begin{align*}
d\equiv 1\bmod 4 &\iff \delta =0 \text{\ y $\varepsilon+t$ es par},\\
d\equiv 2\bmod 4&\iff \delta =1,\\
d\equiv 3\bmod 4&\iff \delta =0 \text{\ y $\varepsilon+t$ es impar}.
\end{align*}
\end{gather*}

Por la f\'ormula del conductor--discriminante, Teorema \ref{T12.3.1.8},
tenemos que si $d\equiv 1\bmod 4$, $\f{\chi_2}=1$. Si $d\equiv 2\bmod 
4$, entonces $\f{\chi_2}=8$ y $\chi_2$ puede corresponder
a ${\ma Q}(\sqrt{2})$ o a ${\ma Q}(\sqrt{-2})$. Si $d\equiv 3
\bmod 4$, entonces $\f{\chi_2}=4$ y $\chi_2$ corresponde a
$\cic 4{}={\ma Q}(\sqrt{-1})$ y $\chi_2(-1)=-1$.

Veamos todos los casos:
\l
\item	{\underline{Si $d\equiv 1\bmod 4$,}}
\lasa
\item $K$ es real, $d>0$. En este caso $\varepsilon=0$, $\delta=0$,
$t$ es par. Se tiene que 
\begin{gather*}
Y=\Big(\bigoplus_{i=1}^s \langle \chi_{p_i}\rangle\Big)\bigoplus
\Big(\bigoplus_{j=1}^t \langle \chi_{q_j}\rangle\Big).\\
\intertext{El campo $L$ correspondiente a $Y$ ser\'a:}
L={\ma Q}(\sqrt{p_1},\ldots,\sqrt{p_s},\sqrt{-q_1},\ldots,\sqrt{-q_t}).
\intertext{Si $t>0$, $L$ es imaginario y entonces}
Y^+=\Big(\bigoplus_{i=1}^s \langle \chi_{p_i}\rangle\Big)\bigoplus
\Big(\bigoplus_{j=2}^t \langle \chi_{q_1}\chi_{q_j}\rangle\Big).\\
\intertext{y el campo asociado a $Y^+$ es}
{\ma Q}(\sqrt{p_1},
\ldots,\sqrt{p_s},\sqrt{q_1q_2},\ldots,\sqrt{q_1q_t}).
\end{gather*}

Notemos que $2^{s+t-2}|h_K$.

Si $t=0$, entonces $Y=Y^+$,
$L={\ma Q}(\sqrt{p_1},\ldots,\sqrt{p_s})$ y $2^{s-1}|h_K$.

\item Si $K$ es imaginario, $d<0$, $\varepsilon=1$, $t$ es impar y
$L={\ma Q}(\sqrt{p_1},\ldots,\sqrt{p_s},\sqrt{-q_1},\ldots,\sqrt{-q_t})$

Notemos que $2^{s+t-1}|h_K$.
\end{list}

\item {\underline{Si $d\equiv 2\bmod 4$.}} Se tiene que
$\f{\chi_2}=8$ en este caso y
\[
Y=\langle\chi_2\rangle\bigoplus
\Big(\bigoplus_{i=1}^s \langle \chi_{p_i}\rangle\Big)\bigoplus
\Big(\bigoplus_{j=1}^t \langle \chi_{q_j}\rangle\Big).
\]
Adem\'as $(-1)^{\varepsilon}=\chi(-1)=\chi_2(-1)(-1)^t$. Por tanto
$\chi_2(-1)=(-1)^{t+\varepsilon}$.

\lasa
\item Si $K$ es real, $d>0$, $\varepsilon=0$, $\chi_2(-1)=(-1)^t$.
Si $t$ es par, $\chi_2(-1)=1$ y $\chi_2$ corresponde a ${\ma Q}(
\sqrt{2})$. Si $t$ es impar, $\chi_2$ corresponde a ${\ma Q}(\sqrt{
-2})$. Se tiene
\[
L={\ma Q}(\sqrt{(-1)^t 2},\sqrt{p_1},\ldots,\sqrt{p_s},
\sqrt{-q_1},\ldots,\sqrt{-q_t}).
\]
Si $t=0$, $L^+=L={\ma Q}(\sqrt{2},\sqrt{p_1},\ldots,\sqrt{p_s})$
y $2^s|h_K$.

Si $t>0$, $t$ par, 
\begin{gather*}
L={\ma Q}(\sqrt{2},\sqrt{p_1},\ldots,\sqrt{p_s},
\sqrt{-q_1},\ldots,\sqrt{-q_t}),\\
L^+={\ma Q}(\sqrt{2}, \sqrt{p_1},\ldots,\sqrt{p_s},
\sqrt{q_1q_2},\ldots,\sqrt{q_1q_t})
\quad\text{y}\quad 2^{s+t-1}|h_K.\\
\intertext{Si $t$ es impar,}
L={\ma Q}(\sqrt{-2},\sqrt{p_1},\ldots,\sqrt{p_s},
\sqrt{-q_1},\ldots,\sqrt{-q_t}),\\
L^+={\ma Q}(\sqrt{p_1},\ldots,\sqrt{p_s},
\sqrt{2q_1},\ldots,\sqrt{2q_t})
\quad\text{y}\quad 2^{s+t-1}|h_K.
\end{gather*}

\item Si $K$ es imaginario, $d<0$, $\varepsilon=1$, 
$\chi_2(-1)=(-1)^{t+1}$ por lo que $\chi_2$ corresponde a
${\ma Q}(\sqrt{(-1)^{t+1}2})$. Se tiene
\[
L={\ma Q}(\sqrt{(-1)^{t+1}2},\sqrt{p_1},\ldots,\sqrt{p_s},
\sqrt{-q_1},\ldots,\sqrt{-q_t})
\]
y $2^{s+t}|h_K$.
\end{list}

\item {\underline{$d\equiv 3\bmod 4$.}} En este caso 
$\f{\chi}=4$ y por tanto $\chi(-1)=-1$ y $\chi_2$ corresponde
a $\cic 4{}={\ma Q}(\sqrt{-1})$.

\lasa
\item si $K$ es real, $d>0$, $\varepsilon=0$, $\chi_2(-1)=-1$ y
$\chi(-1)=(-1)^{\varepsilon}=1=\chi_2(-1)\prod_{j=1}^{t}
\chi_{q_j}(-1)=(-1)^{1+t}$. Por lo tanto $t$ es impar. Entonces
\begin{gather*}
L={\ma Q}(\sqrt{-1},\sqrt{p_1},\ldots,\sqrt{p_s},
\sqrt{-q_1},\ldots,\sqrt{-q_t}),\\
L^+={\ma Q}(\sqrt{p_1},\ldots,\sqrt{p_s},
\sqrt{q_1},\ldots,\sqrt{q_t})
\quad\text{y}\quad 2^{s+t-1}|h_K.
\end{gather*}

\item Si $K$ es imaginario, $d<0$, $\varepsilon=1$, $\chi(-1)=-1$.
Por lo tanto $\chi(-1)=(-1)^{\varepsilon}=-1=\chi_2(-1)\prod_{j=1}^{t}
\chi_{q_j}(-1)=(-1)^{1+t}$. Por lo tanto $t$ es par y
\[
L={\ma Q}(\sqrt{-1},\sqrt{p_1},\ldots,\sqrt{p_s},
\sqrt{-q_1},\ldots,\sqrt{-q_t})
\quad\text{y}\quad 2^{s+t}|h_K.
\]
\end{list}
\end{list}

Resumiendo, si $K$ es real $2^{m-2}|h_K$ 
y $2^{m-1}|h_K$ en el caso
en que $K$ es imaginario, donde $m$ es el n\'umero
de primos que dividen a $d$.
\end{ejemplo}

Ahora estudiamos el {\em s\'imbolo de 
Jacobi\index{Jacobi!s\'imbolo}\index{simbolo de Jacobi@s\'imbolo de Jacobi}}.
En realidad ya hicimos gran parte de lo que a continuaci\'on
exponemos en el Ejemplo \ref{Ej12.4.10}.
Sea $d\in{\ma Z}$, $d\neq 0$ libre de cuadrados. Consideremos
$K_d={\ma Q}(\sqrt{d})$ y sea $\sgn d:=\begin{cases}
1&\text{si $d>0$},\\ -1&\text{si $d<0$}\end{cases}$.
Entonces $d=(\sgn d)2^{\varepsilon}p_1\cdots p_s$ con
$\varepsilon\in\{0,1\}$. Se tiene $\delta_{{\ma Q}(\sqrt{d})}=
\begin{cases} d&\text{si $d\equiv 1\bmod 4$}, \\
4d&\text{si $d\equiv 2,3\bmod 4$}\end{cases}$.

Si $d\equiv 1\bmod 4$, los primos finitos ramificados en $K_d/
{\ma Q}$ son $p_1,\cdots, p_s$ y el caracter $\chi_d$ asociado
a $K_d$ es $\chi_d=\chi_{p_1}\cdots\chi_{p_s}$ donde cada
$\chi_{p_i}$ es el caracter cuadr\'atico de conductor $p_i$ y por
tanto $\chi_{p_i} =\artinp{}{p_i}$ y el campo asociado a $\chi_{p_i}$
es ${\ma Q}\Big(\sqrt{(-1)^{\frac{p_i-1}{2}}p_i}\Big)$.

Si $d\equiv 3\bmod 4$, $-d=(-1)^{\frac{d-1}{2}}d\equiv 1\bmod 4$
y se tiene que ${\ma Q}(\sqrt{d},i)={\ma Q}(\sqrt{d},\sqrt{-d})$.
Entonces $\chi_{-d}=\chi_{p_1}\cdots \chi_{p_s}$ y $\chi_d=
\varphi\chi_{p_1}\cdots \chi_{p_s}$ donde $\varphi$ es el
\'unico caracter de conductor $4$: $\varphi(q)=(-1)^{\frac{q-1}{2}}$
para $q$ primo impar.

Ahora consideremos $d\equiv 2\bmod 4$, esto es, $\varepsilon =1$
y $\delta_{{\ma Q}(\sqrt{d})}=4d$. Se sigue que $\chi_d=\psi
\chi_{p_1}\cdots \chi_{p_s}$ donde $\psi$ es un caracter de 
conductor $8$. Si $\psi$ es el caracter real, $\psi$ corresponde
a ${\ma Q}(\sqrt{2})$ y si $\psi$ es el caracter imaginario, $\psi$
corresponde a ${\ma Q}(\sqrt{-2})$.

Ordenamos los factores de $d$ de tal forma que $p_1,\ldots,p_t$ son
congruentes con $3$ m\'odulo $4$ y $p_{t+1},\ldots p_s$ son
congruentes con $1$ m\'odulo $4$. Se sigue que $\chi_d(-1)=
\psi(-1)(-1)^t=\sgn d$.

Si $\sgn d=1$, esto es, $d>0$, entonces $\psi(-1)=(-1)^t$
y si $\sgn d=-1$, es decir, $d<0$, entonces $\psi(-1)=(-1)^{t+1}$.
Por otro lado se tiene $\frac{|d|}{2}=\frac{(\sgn d)d}{2}=
p_1\cdots p_s\equiv p_1\cdots p_t\equiv (-1)^t\bmod 4$.

En particular $p_1\cdots p_s\equiv 1\bmod 4\iff t$ es par $\iff
\frac{|d|}{2}\equiv 1\bmod 4$. De esta forma obtenemos que si
$\psi_r$ es el caracter cuadr\'atico real de conductor $8$ y si
$\psi_i$ es el caracter cuadr\'atico imaginario de conductor $8$,
se tiene
\[
\chi_d=\begin{cases}
\chi_{p_1}\cdots \chi_{p_s}&\text{si $d\equiv 1\bmod 4$},\\
\varphi \chi_{p_1}\cdots \chi_{p_s} &\text{si $d\equiv 3\bmod 4$},\\
\psi_r\chi_{p_1}\cdots \chi_{p_s}&\text{si $d\equiv 2\bmod 4$ y 
$\frac{|d|}{2}\equiv 1\bmod 4$},\\
\psi_i\chi_{p_1}\cdots \chi_{p_s}&\text{si $d\equiv 2\bmod 4$ y 
$\frac{|d|}{2}\equiv 3\bmod 4$}.
\end{cases}
\]

Ahora bien, si $d\equiv 2\bmod 4$ y $\frac{|d|}{2}\equiv 1\bmod 4$,
entonces $|d|\equiv 2\bmod 8$ y si $d\equiv 2\bmod 4$ y
$\frac{|d|}{2}\equiv 3\bmod 4$, entonces $|d|\equiv 6\bmod 8$.

Sean $\delta(d)=\frac{1+(-1)^d}{2}$ y $\beta(d)=\frac{1-(-1)^{d^2
\frac{d-1}{2}}}{2}$. Se tiene que
\begin{gather*}
\delta(d)=\begin{cases} 0&\text{si $d\equiv 1\bmod 2$},\\
1&\text{si $d\equiv 0\bmod 2$}\end{cases},\qquad
\beta(d)=\begin{cases} 1&\text{si $d\equiv 3\bmod 4$},\\
0&\text{si $d\equiv 1,2\bmod 4$}\end{cases},\\
\frac{1-(-1)^{\frac{|d|+2}{4}}}{2}=\begin{cases} 1&\text{si $
|d|\equiv 2\bmod 8$},\\
0&\text{si $|d|\equiv 6\bmod 8$}\end{cases} \quad\text{y}\\
\frac{1-(-1)^{\frac{|d|+6}{4}}}{2}=\begin{cases} 0&\text{si $
|d|\equiv 2\bmod 8$},\\
1&\text{si $|d|\equiv 6\bmod 8$}.\end{cases}
\end{gather*}
Se sigue que si
\begin{gather}\label{E*.1}
\xi=\varphi^{\beta(d)}\psi_r^{\delta(d)\frac{1-(-1)^{\frac{|d|+2}{4}}}{2}}
\psi_i^{\delta(d)\frac{1-(-1)^{\frac{|d|+6}{4}}}{2}},
\end{gather}
entonces
\[
\xi=\begin{cases} 1&\text{si $d\equiv 1\bmod 4$},\\
\varphi&\text{si $d\equiv 3\bmod 4$},\\
\psi_r&\text{si $|d|\equiv 2\bmod 8$},\\
\psi_i&\text{si $|d|\equiv 6\bmod 8$}.
\end{cases}
\]

Definamos $\artinp{}{2}:=\xi$. Se sigue que $\chi_d=
\xi \chi_{p_1}\cdots \chi_{p_s}=\artinp{}{2}\artinp{}{p_1}\cdots
\artinp{}{p_s}$.

\begin{definicion}\label{D2.3-1}
Se define el {\em s\'imbolo de 
Jacobi\index{Jacobi!s\'imbolo}\index{simbolo de Jacobi@s\'imbolo de Jabobi}}
para $d\in{\ma Z}$, $d\neq 0$ y $d$ libre de cuadrados como
el caracter $\chi_d=\artinp{}{d}$ 
asociado al campo $K_d:={\ma Q}(\sqrt{d})$ y $\artinp{}{d}$
est\'a dado por
\[
\artinp{}{d}=\artinp{}{2}\artinp{}{p_1}\cdots \artinp{}{p_s}
\]
donde $\artinp{}{2}=\xi$ est\'a dado por (\ref{E*.1}).
\end{definicion}

\begin{observacion}\label{O2.3-2}
Del Corolario \ref{C12.3.6'}, se tiene que si $p$ es un primo
racional, entonces $p$ se descompone en $K_d/{\ma Q}$
si y solamente si $\artinp{p}{d}=1$.
\end{observacion}

Los Ejemplos \ref{Ej12.4.7}, \ref{Ej12.4.8}, \ref{Ej12.4.9} y \ref{Ej12.4.10}
son casos particulares de lo que se conoce como
{\em campos de g\'eneros\index{campos de g\'eneros}}.
El concepto de campo de g\'eneros se remonta a Gauss
\cite{Gau1801} en el contexto de formas cuadr\'aticas 
binarias. Para cualquier extensi\'on finita $K/{\ma Q}$, 
el campo de g\'eneros se define como la m\'axima extensi\'on $\g K $
no ramificada de $K$ tal que $\g K $ es la composici\'on de
 $K$ y de una extensi\'on abeliana
 $k^{\ast}$ de ${\ma Q}$: $\g K =Kk^{\ast}$. Esta definici\'on
 se debe a Fr\"olich \cite{Fro83}.
Si $K_H$ denota el campo de clase de Hilbert
de $K$ (Teorema \ref{T9.3}), 
se tiene $K\subseteq \g K \subseteq K_H$. Originalmente
la definici\'on de campos de g\'eneros fue dada para una
extensi\'on cuadr\'atica de ${\ma Q}$. De hecho Gauss
prob\'o que si $t$ es el n\'umero de primos
positivos diferentes que dividen al discriminante
$\delta_K$ de un campo cuadr\'atico num\'erico
 $K$, entonces el $2$--rango del grupo de clases de $K$
es $2^{t-2}$ si $\delta_K>0$ y existe un primo $p\equiv
3\bmod 4$ que divide a $\delta_K$ y $2^{t-1}$ en cualquier
otro caso (ver Ejemplos \ref{Ej12.4.10} y \ref{E2.3}).

H. Leopoldt \cite{Leo53} (ver Teorema \ref{T12.4.2})
determin\'o el campo de g\'eneros $\g K $ de una
extensi\'on abeliana $K$ de ${\ma Q}$ usando caracteres de Dirichlet,
generalizando de esta manera el trabajo de 
H. Hasse \cite{Has51} el cual introdujo
la teor{\'\i}a del g\'enero para campos cuadr\'aticos num\'ericos.

M. Ishida determin\'o el campo de g\'eneros $\g K $ de cualquier
extensi\'on abeliana finita de ${\ma Q}$ \cite{Ish76}. 
X. Zhang \cite{Xia85} dio una expresi\'on simple
de $\g K $ para cualquier extensi\'on abeliana $K$ de ${\ma Q}$ 
usando la teor{\'\i}a de ramificaci\'on de Hilbert.

En esta subsecci\'on presentamos brevemente la teor{\'\i}a
del g\'enero.

Sea $K$ un campo num\'erico, es decir, una extensi\'on 
finita de ${\ma Q}$.
Sea $K_H$ el campo de clase de Hilbert de $K$,
esto es, $K_H$ es la m\'axima extensi\'on abeliana
no ramificada de $K$. Entonces el campo de
g\'eneros $\g K $ de $K$ es la m\'axima extensi\'on de $K$
contenida en $K_H$ tal que sea la composici\'on de $K$
y de una extensi\'on abeliana $k^{\ast}$ de ${\ma Q}$. 
Equivalentemente, $\g K =K k^{\ast} \subseteq K_H$ con 
$k^{\ast}$ siendo la m\'axima extensi\'on abeliana
de ${\ma Q}$ contenida en $K_H$.

A continuaci\'on presentamos la teor{\'\i}a del g\'enero en
el caso abeliano para campos num\'ericos 
\cite{Leo53}. La teor{\'\i}a que presentamos est\'a basada
en el Teorema \ref{T12.4.2}. En este caso $\g K $ es la m\'axima extensi\'on
de $K$ contenida en $K_H$ tal que $\g K /{\ma Q}$ es 
abeliana. As{\'\i} pues, en esta secci\'on consideramos
$K/{\ma Q}$ una extensi\'on abeliana.
Por el Teorema de Kronecker--Weber (Teorema \ref{T1.4.5*}) existe
$n\in{\ma N}$ tal que $K\subseteq \cic n{}$, donde $\zeta_n$ denota
una ra{\'\i}z $n$--\'esima primitiva de uno.
Sea $X$ el grupo de caracteres de Dirichlet asociado a
$K$. 

El siguiente ejemplo es un teorema debido a Gauss. No es m\'as que
el Ejemplo \ref{Ej12.4.10} visto de manera m\'as estructural.

\begin{ejemplo}[Teorema del G\'enero de Gauss]\label{E2.3}
Sea $K={\ma Q}(\sqrt{d})$ una extensi\'on cuadr\'atica de ${\ma Q}$,
donde $d\in{\ma Z}$ es libre de cuadrados. Sea $m$ el n\'umero
de factores primos diferentes de $\delta_K$, el discriminante de
$K$. Si $p_1,\ldots, p_m$ son estos factores, seleccionamos,
$p_1=2$ si $2\mid \delta_K$.

Denotemos por ${\ma P}'_{\ma Q}$ al conjunto 
de los n\'umeros primos en ${\ma Q}$.
Sea $\chi$ el caracter cuadr\'atico asociado a $K$. Entonces
$\chi_{p_i}\neq 1$, $1\leq i\leq m$ y $\chi_q=1$ para todo
$q\in {\ma P}'_{\ma Q}\setminus \{p_1,\ldots,p_m\}$. Para $p_i\neq 2$,
$\chi_{p_i}$ es \'unico y $\chi_{p_i}(-1)=(-1)^{(p_i-1)/2}$.
En este caso el campo asociado a $\chi_{p_i}$ es ${\ma Q}\big(
\sqrt{(-1)^{(p_i-1)/2}p_i}\big)$. Si $p_1=2$, entonces hay
tres caracteres cuadr\'aticos $\chi_{p_1}=\chi_2$; dos de
ellos tienen conductor $8$, uno es real y el otro imaginario, y
el otro tiene conductor $4$. Si $\chi_2$ es real, $\chi(-1)=1$ y el
campo asociado es ${\ma Q}(\sqrt{2})$. Si $\chi_2$ es imaginario
de conductor $8$, $\chi(-1)=-1$ y el campo asociado es
${\ma Q}(\sqrt{-2})$. Finalmente, si $\chi_2$ es de conductor $4$, 
$\chi(-1)=-1$ y el campo asociado a $\chi_2$ es ${\ma Q}
(\zeta_4)={\ma Q}(i)={\ma Q}(\sqrt{-1})$. Se sigue 
que la m\'axima extensi\'on abeliana de ${\ma Q}$ no ramificada
en ning\'un primo finito es
$J={\ma Q}\big(\sqrt{\varepsilon}, \sqrt{(-1)^{
(p_i-1)/2}p_i}\mid 2\leq i\leq m\big)$ donde $\varepsilon=
(-1)^{(p_1-1)/2}p_1$ si $p_1\neq 2$ y $\varepsilon =2, -2$ o $-1$ si 
$p_1=2$.

Obtenemos que $[J:{\ma Q}]=2^m$ y que $[J:K]=2^{m-1}$. Se tiene que
$\g K =J$ excepto cuando $K$ es real y $J$ es imaginario y este \'ultimo caso
ocurre cuando $\delta_K>0$ ($d>0$) y existe
$p_i\equiv 3\bmod 4$. En este caso, $[J^+:K]=2^{m-2}$.
Por ejemplo, para la
 extensi\'on cuadr\'atica $K={\ma Q}(\sqrt{-14})$ sobre ${\ma Q}$,
tenemos $\g K ={\ma Q}(\sqrt{2},\sqrt{-7})$ y para $K=
{\ma Q}(\sqrt{79})$ obtenemos $J={\ma Q}(\sqrt{-79},i)$ y
$\g K =J^+=J\cap {\ma R}={\ma Q}(\sqrt{79})=K$.

Ahora si $I_K$ es el grupo de clase de $K$, $I_K
\cong \Gal(K_H/K)$ y $E$ es el campo fijo de ${I}_K^2$, entonces
$\Gal(E/K)\cong {I}_K/{I}_K^2$. Puesto que $\g K $ 
es la m\'axima extensi\'on abeliana de ${\ma Q}$ contenida en $K_H$, $\g K $ es el
subcampo fijado de $K_H$ bajo el grupo derivado $G'$ de $G$.
Puede ser verificado que (ver \cite{Jan73}) $G'={I}_K^2$ as{\'\i} que $\g K =E$ y
se sigue que el $2$--rango de ${I}_K$ es $2^{m-1}$ a menos que
$d>0$ y existe un primo $p\equiv 3\bmod 4$ que divide a $d$ y en
este caso el $2$--rango de ${I}_K$ es $2^{m-2}$.
\end{ejemplo}

\begin{ejemplo}\label{E2.4}
Si $p$ es un primo impar, $K$ es una extensi\'on c{\'\i}clica de
${\ma Q}$ de grado $p$ y si $m$ es el n\'umero de primos ramificados
en $K$, se sigue que $\g K $ es una $p$--extensi\'on elemental
abeliana de ${\ma Q}$ de grado $p^m$ y $[\g K :K]=
p^{m-1}$. En particular $p^{m-1}\mid |{I}_K|$.
\end{ejemplo}

Ahora sea $K$ una extensi\'on abeliana de ${\ma Q}$ con
grupo de caracteres de Dirichlet $X$. Consideremos para cada
 $p\in{\ma P}'_{\ma Q}$, $X_p$. Sea
$J$ el campo asociado a $\prod_{p\in{\ma P}'_{\ma Q}} X_p$. 
Sea $p^{m_p}:=\mcm\{{\eu f}_{\chi_p}\mid \chi\in X\}$ donde
${\eu f}_{\chi_p}$ denota al conductor de $\chi_p$. Entonces
el campo $K_p$ asociado a $X_p$ est\'a contenido en $\cic p{m_p}$
pero no en $\cic p{m_p-1}$. Si $p$ es impar, $K_p$ es el \'unico
subcampo de $\cic p{m_p}$ de grado $|X_p|$ sobre ${\ma Q}$ y
$K_p/{\ma Q}$ es una extensi\'on c{\'\i}clica. Si $p=2$, $K_2$ es uno de los
siguientes campos. Si $|X_2|=\varphi(2^{m_2})=2^{m_2-1}$, $K_2=
\cic 2{m_2}$. Si $|X_2|=\frac{\varphi(2^{m_2})}{2}=2^{m_2-2}$,
$K_2=\cic 2{m_2}^+=
{\ma Q}\big(\zeta_{2^{m_2}}+\zeta_{2^{m_2}}^{-1}\big)=\cic 2{m_2}
\cap {\ma R}$ si $\chi(-1)=1$ para todo $\chi\in X$ y $K_2=
{\ma Q}\big(\zeta_{2^{m_2}}-\zeta_{2^{m_2}}^{-1}\big)$ si existe
$\chi\in X$ con $\chi(-1)=-1$. 

Por lo tanto, si $K$ y $J$ son ambos reales o ambos imaginarios,
$\g K =J=\prod_{p\in{\ma P}'_{\ma Q}}K_p$. Si $K$ es real y
$J$ es imaginario, $\g K =J^+=J\cap {\ma R}$.

\begin{ejemplo}[ver Ejemplo \ref{Ej11.-2}]\label{E2.5}
Sean $K={\ma Q}\Big(\sqrt{65+\sqrt{65}}\Big)$ y $\alpha
=\sqrt{65+\sqrt{65}}$. Se tiene $(\alpha^2-65)^2=65$, esto es,
$\alpha^4-2\cdot 65\alpha^2+65^2-65=\alpha^4-130\alpha^2
+4160=0$. Sea $p(x):=(x^2-65)^2-65=x^4-130x^2+4160$.
Entonces $p(x)$ es Eiseinstein para $p=13$ y por tanto $p(x)$
es irreducible y $p(x)=\Irr(\alpha,x,{\ma Q})$.

Las ra\'ices de $p(x)$ son $\pm \alpha$ y $\pm \beta$ donde
$\beta=\sqrt{65-\sqrt{65}}$. Se tiene que $\alpha\beta=
\sqrt{65^2-65}=\sqrt{65\cdot 64}=8\sqrt{65}$. Adem\'as
$\alpha^2=65+\sqrt{65}=65+\frac{\alpha\beta}{8}$, es decir,
$\beta=\frac{8(\alpha^2-65)}{\alpha}\in{\ma Q}(\alpha)$. Por
tanto las cuatro ra\'ices de $p(x)$ est\'an en $K={\ma Q}(
\alpha)$ por lo que $K/{\ma Q}$ es una extensi\'on de
Galois de grado $4$.

Sea $\sigma\in G=\Gal(K/{\ma Q})$ tal que $\sigma(\alpha)=
\beta$. Entonces $\sigma(\alpha)=\beta
=\frac{8(\alpha^2-65)}{\alpha}=8\big(\alpha-\frac{65}{
\alpha}\big)$ de donde se sigue $\sigma(\beta)=
8\big(\sigma(\alpha) -\frac{65}{\sigma(\alpha)}\big)=
8\big(\beta-\frac{65}{\beta}\big)$. 

Ahora bien, $\beta^2=65-\sqrt{65}=65-\frac{\alpha\beta}{8}$
por lo que $\alpha=\frac{8(65-\beta^2)}{\beta}=8\big(\frac{65}
{\beta}-\beta\big)=-\sigma(\beta)$. Esto es, $\sigma(\beta)=
-\alpha=\sigma^2(\alpha)$. Por tanto $\sigma$ es de orden $4$
y $G=\langle\sigma\rangle\cong C_4$ es un grupo c\'iclico.

Puesto que $\sqrt{65}=\frac{\alpha\beta}{8}\in K={\ma Q}(\alpha)$,
se sigue que $K_0={\ma Q}(\sqrt{65})$ es el \'unico subcampo 
cuadr\'atico de $K/{\ma Q}$. Ahora bien, como $65\equiv 1\bmod 4$,
se tiene $\delta_{{\ma Q}(\sqrt{65})}=65$. Por tanto
$5$ y $13$ son los 
\'unicos primos ramificados en $K_0/{\ma Q}$ y se tiene $e_5(
K_0|{\ma Q})=e_{13}(K_0|{\ma Q})=2$.

Del hecho de que $K/{\ma Q}$ es c\'iclica de grado $4=2^2$,
se sigue que $5$ y $13$ son totalmente ramificados, es decir, $e_5(
K|{\ma Q})=e_{13}(K|{\ma Q})=4$. Adem\'as $5$y $13$ son
moderadamente ramificados. Los exponentes de $5$ y
$13$ en el diferente son ${\eu D}_{K/{\ma Q},5}=\pK_5^3$ y
${\eu D}_{K/{\ma Q},13}=\pK_{13}^3$, donde $\pK_5$ y $\pK_{13}$
son los \'unicos en $K$ sobre $5$ y $13$ respectivamente.

Sean $\mu$ el caracter asociado a $K_0/{\ma Q}$ y $\lambda$
el caracter asociado a $K/{\ma Q}$. Entonces $\mu=\lambda^2$,
$\mu=\mu_5\mu_{13}$ y $\mu_5$ y $\mu_{13}$ son los
\'unicos caracteres cuadr\'aticos de conductores $5$ y $13$
respectivamente. Por tanto sus campos asociados son
${\ma Q}(\sqrt{5})$ y ${\ma Q}(\sqrt{13})$ ya que $5\equiv 1
\bmod 4$ y $13\equiv 1\bmod 4$. Se sigue que el campo de
g\'eneros de $K_0$ es $\g{(K_0)}={\ma Q}(\sqrt{5},\sqrt{13})$.
Adem\'as tenemos $\mu_5(-1)=\mu_{13}(-1)=1$.

Falta por averiguar que primos son ramificados en $K/K_0$
adem\'as de $5$ y $13$. El discriminante de $\{1,\alpha,
\alpha^2,\alpha^3\}$ es (ver Secci\'on \ref{Sec0.1})
\[
{\eu d}_{K/{\ma Q}}(1,\alpha,\alpha^2,\alpha^3)=4^{4(4-1)/2}
\N_{L/K}(p'(\alpha)) \quad\text{y}\quad {\eu d}_{K/{\ma Q}}|
{\eu d}_{K/{\ma Q}}(1,\alpha,\alpha^2,\alpha^3).
\]

Se tiene $p'(x)=4x^3-260x$, $p'(\alpha)=4\alpha(\sqrt{65})$,
$\N_{K/{\ma Q}}(\alpha)=65\cdot 64$ (el t\'ermino constante
de $p(x)$) y $\N_{K/{\ma Q}}(\sqrt{65})=65^2$. Por tanto obtenemos
${\eu d}_{K/{\ma Q}}(1,\alpha,\alpha^2,\alpha^3)=2^{14}\cdot
5^3\cdot 13^3$. Se sigue que el \'unico primo, adem\'as $5$
y $13$, posiblemente ramificado en $K/K_0$, es $2$.

Ahora $\AE {{{\ma Q}(\sqrt{65})}}={\ma Z}\big[\frac{\sqrt{65}+1}
{2}\big]$ pues $65\equiv 1\bmod 4$.

Sea $\phi=\frac{\sqrt{65}+1}{2}$. Entonces $(2\phi-1)^2=
4\phi^2-4\phi+1=65$, esto es, $\Irr(\phi,x,{\ma Q})=
x^2-x-16$. Puesto que $x^2-x-16\bmod 2=x(x-1)$,
por el Teorema de Kummer tenemos que $2$ se descompone
en $K_0/{\ma Q}$ y $(2)={\mc P}_2{\mc P}'_2=
\big(2,\frac{\sqrt{65}+1}{2}\big)\big(2, \frac{\sqrt{65}-1}{2}\big)$,
${\mc P}_2\neq {\mc P}'_2$. 

Puesto que $\phi\in{\mc P}_2$,
se tiene $\phi-1=\frac{\sqrt{65}-1}{2}\notin {\mc P}_2$.
Por tanto 
\begin{align*}
0&=v_{{\mc P}_2}\big(\frac{\sqrt{65}-1}{2}\big)=
v_{{\mc P}_2}(\sqrt{65}-1)-v_{{\mc P}_2}(2)\\
&=v_{{\mc P}_2}(\sqrt{65}-1)-e_{K_0|{\ma Q}}({\mc P}_2|2)
v_2(2)=v_{{\mc P}_2}(\sqrt{65}-1)-1.
\end{align*}
As\'i, obtenemos que $v_{{\mc P}_2}(\sqrt{65}-1)=1$.
An\'alogamente $v_{{\mc P}'_2}(\sqrt{65}+1)=1$.

Ahora $64=65-1=(\sqrt{65}+1)(\sqrt{65}-1)$ por lo que
$6=v_2(64)=v_{{\mc P}_2}(64)=v_{{\mc P}_2}(\sqrt{65}+1)+
v_{{\mc P}_2}(\sqrt{65}-1)=v_{{\mc P}_2}(\sqrt{65}+1)+1$.
Por tanto $v_{{\mc P}_2}(\sqrt{65}+1)=5$. An\'alogamente
se tiene $v_{{\mc P}'_2}(\sqrt{65}-1)=5$.

Ahora bien $K=K_0(\alpha)$, $\alpha=\sqrt{65+\sqrt{65}}$.
Sea $\delta=\alpha^2=65+\sqrt{65}=\sqrt{65}(\sqrt{65}+1)$. Por
tanto $\Irr(\alpha,x,K_0)=x^2-\delta$. Se tiene
\[
v_{{\mc P}_2}(\delta)=v_{{\mc P}_2}(\sqrt{65})+
v_{{\mc P}_2}(\sqrt{65}+1)=0+5=5\quad\text{y}\quad
\mcd(2,5)=1,
\]
por tanto ${\mc P}_2$ es ramificado en $K/K_0$. De hecho,
si $\pK_2\in{\ma P}_K$ es un primo sobre ${\mc P}_2$, se 
tiene
\[
v_{\pK_2}(\alpha^2)=2v_{\pK_2}(\alpha)=v_{\pK_2}(\delta)
=e_{K|K_0}(\pK_2|{\mc P}_2)v_{{\mc P}_2}(\delta)=
e_{K|K_0}(\pK_2|{\mc P}_2)\cdot 5,
\]
de donde se sigue que $2|e_{K|K_0}(\pK_2|{\mc P}_2)$ y
por tanto $e_{K|K_0}(\pK_2|{\mc P}_2)=2$ y ${\mc P}_2$
es ramificado en $K/K_0$.

En resumen, los primos ramificados en $K/{\ma Q}$ son 
$2$, $5$ y $13$ y se tiene $e_2(K|{\ma Q})=2$, $e_5(K|{\ma Q})
=4$ y $e_{13}(K|{\ma Q})=4$.

Se tiene que $2$ es salvajemente ramificado en $K/{\ma Q}$ por
lo tanto ${\eu D}_{K/{\ma Q}}=(\pK_2\pK'_2)^{e^*} \pK_5^3
\pK_{13}^3$ con $e^*\geq 2$. Adem\'as
$\N_{K/{\ma Q}}\pK_2=(2)=\N_{K_0/{\ma Q}}
\N_{K/K_0}\pK_2=\N_{K_0/{\ma Q}}{\mc P}_2^f=(2)^f=(2)$.
Por tanto
\[
{\eu d}_{K/{\ma Q}}=(-1)^{r_2}(2\cdot 2)^{e^*}\cdot (5)^3
\cdot (13)^3=2^{2e^*}\cdot 5^3\cdot 13^3,\quad e^*\geq 2.
\]

Si $\chi$ es el caracter asociado a $K$, se tiene $\chi=
\chi_2\chi_5\chi_{13}$ con $\f{\chi_5}=5$, $\f{\chi_{13}}=13$
y $o(\chi_2)=2$, $o(\chi_5)=o(\chi_{13})=4$. De esta forma tenemos que
los campos asociados a $\chi_5$ y $\chi_{13}$ son extensiones
c\'iclicas de ${\ma Q}$ de grado $4$ contenidos en
$\cic 5{}$ y $\cic {13}{}$ respectivamente. Se sigue que
el campo asociado a $\chi_5$ es $\cic 5{}$. Sea $K_{13}$
el campo asociado a $\chi_{13}$. Se tiene $[K_{13}:{\ma Q}]
=4$. Por tanto $[\cic {13}{}:K_{13}]=3$ y en particular $K
\not\subseteq \cic {13}{}^+$ pues $[\cic {13}{}:\cic {13}{}^+]=2$.

Por tanto $\chi_5(-1)=\chi_{13}(-1)=-1$ y $\chi(-1)=1$.
Se sigue que $\chi_2(-1)=1$ por lo que $\chi_2$ es un caracter
positivo par. Necesariamente $\chi_2$ es el caracter cuadr\'atico
real de conductor $\f {\chi_2}=8=2^3$ con campo asociado
${\ma Q}(\sqrt{2})$.

El grupo de caracteres del campo de g\'eneros extendido de $K$ es
$Y=\langle \chi_2,\chi_5,\chi_{13}\rangle\cong C_2\times C_4
\times C_4$. El campo correspondiente a $Y$ es $K_Y=
{\ma Q}(\sqrt{2},\zeta_5,\beta)$ con $K_{13}={\ma Q}(\beta)$.
Se tiene $K_Y\not\subseteq {\ma R}$ y $K\subseteq {\ma R}$.
Por tanto $\g K=K_Y^+=K_Y\cap {\ma R}$. El grupo de
caracteres de $\g K$ es
\[
Y^+=\langle \chi\in Y\mid \chi(-1)=1\rangle =\langle \chi_2,
\chi_5^2,\chi_5\chi_{13}\rangle\cong C_2\times C_2\times C_4.
\]
En efecto este conjunto es $Y^+$ pues $[Y:Y^+]=2$ y para
toda $\chi\in Y^+$ se satisface $\chi(-1)=1$.

Se tiene $K\subseteq \g K\subseteq {\ma R}$ y $[\g K:K]=
\frac{[\g K:{\ma Q}]}{[K:{\ma Q}]}=\frac{16}{4}=4$ y 
$\Gal (\g K/K)\cong \frac{Y}{Y^+}\cong C_2\times C_2$. De
hecho $\g K={\ma Q}\big(\sqrt{2},\sqrt{5},\sqrt{65+\sqrt{65}}\big)$
puesto que $K\subseteq \g K\subseteq {\ma R}$, $\g K
\subseteq K_Y$ y $[K_Y:\g K]=2$.

Finalmente, por la f\'ormula del conductor--discriminante, tenemos
\begin{gather*}
{\eu d}_{K/{\ma Q}} =\f{\chi}\f{\chi^2}\f{\chi^3}=(8\cdot 5\cdot 13)
(5\cdot 13)(8\cdot 5\cdot 13)=2^6\cdot 5^3\cdot 13^3,\quad e^*=3,\\
{\eu D}_{K/{\ma Q}}=(\pK_2\pK'_2)^3 \cdot \pK_5^3\cdot\pK_{13}^3.
\end{gather*}

Usando los resultados de la Subsecci\'on \ref{S4.10.2}, se puede
calcular $K_{13}$: $K_{13}={\ma Q}\big(\sqrt{13}+\sqrt{2} i
\sqrt{\sqrt{13}+3}\big)$.

\end{ejemplo}

\subsection{Grupos abelianos como grupos de Galois y
 de clases}\label{S12.4.1}
 
 Sea $G$ un grupo abeliano finito dado. Escribamos $G=C_{n_1}
 \times \cdots \times C_{n_r}$. El primer resultado es que $G$ es
 realizable como grupo de Galois sobre cualquier campo num\'erico.
 
 \begin{teorema}\label{T12.4.1.1}
 Sea $E$ cualquier campo num\'erico y sea $G$ un grupo abeliano
 finito arbitrario. Entonces existen una infinidad de campos $F$
 tales que $F/E$ es una extensi\'on de Galois y tales que $G\cong
 \Gal(F/E)$.
 \end{teorema}
 
 \begin{proof}
 Sean $q_1,\ldots, q_s$ los primos ramificados en $E/{\ma Q}$.
 Puesto que hay una infinidad de primos congruentes a $1$
 m\'odulo $n$ con $n\in{\ma N}$ dado (Corolario \ref{C7.3}),
 seleccionamos primos $p_1,\ldots, p_r$ tales que $p_1<p_2
 <\cdots <p_r$, $q_1,\ldots,q_s<p_1$ y $p_i\equiv 1\bmod n_i$,
 $1\leq i\leq r$. 
 
 \begin{window}[3,l,\xymatrix{
 M\ar@{-}[r]\ar@{-}[d]_G&F=ME\ar@{-}[d]^G\\
 {\ma Q}\ar@{-}[r]&E
 },{}]
 Para cada $1\leq i\leq r$, sea $M_i$ el \'unico subcampo de
 $\cic {p_i}{}$ tal que $[M_i:{\ma Q}]=n_i$ el cual existe pues
 $n_i|p_i-1=[\cic {p_i}{}:{\ma Q}]$. Se tiene que $M_i/{\ma Q}$
 es c{\'\i}clica de grado $n_i$ y $p_i$ es el \'unico primo finito
 ramificado en $M_i/{\ma Q}$. Puesto que los primos ramificados
 en cada $M_i$ son distintos entre si y distintos a los primos
 ramificados en $E/{\ma Q}$, se tiene que si $M:=M_i\cdots M_r$,
 entonces $\Gal(M/{\ma Q})\cong \prod_{i=1}^r\Gal(M_i/{\ma Q})
 \cong \prod_{i=1}^r C_{n_i}\cong G$ y $M\cap E={\ma Q}$. Sea
 $F:=ME$. Por Teor{\'\i}a de Galois tenemos que $F/E$ es
 una extensi\'on de Galois y $\Gal(F/E)\cong \Gal(M/M\cap E)=
 \Gal(M/{\ma Q})\cong G$.
 \end{window}

Puesto que tenemos una infinidad de selecciones para los n\'umeros
primos $p_1,\ldots, p_r$, hay una infinidad de campos $M$ no
isomorfos que satisfacen lo anterior. $\fin$
 \end{proof}
 
 Podemos mejorar en ciertos aspectos el Teorema \ref{T12.4.1.1}.
 
 \begin{teorema}\label{T12.4.1.2}
 Sea $G$ un grupo abeliano finito. Entonces existen campos 
 num\'ericos $L$ y $K$ tales que:
 \lasa
 \item $\Gal(L/K)\cong G$,
 \item $L/K$ es no ramificada en ning\'un primo, incluyendo a
 los primos infinitos.
 \item Se puede seleccionar $L/{\ma Q}$ una extensi\'on abeliana
 y $K/{\ma Q}$ una extensi\'on c{\'\i}clica.
 \end{list}
 \end{teorema}
 
 \begin{proof}
 Sea $G\cong C_{n_1}\times \cdots\times C_{n_r}$. Nuevamente
 seleccionamos primos distintos $p_1,\ldots, p_r$ tales que $
 p_i\equiv 1\bmod 2n_i$, $1\leq i\leq r$ (el papel que juega el $2$
 en estas congruencias ser\'a relacionado m\'as adelante
 con el comportamiento de los primos infinitos en las extensiones).
 
 Se tiene que $\Gal (\cic {p_i}{}/{\ma Q})\cong U_{p_i}\cong
 C_{p_i-1}$. Sea $\psi_i$ un caracter de conductor $p_i$ que 
 genere al grupo c{\'\i}clico $\widehat{U_{p_i}}$, es decir,
 $\psi_i$ es de orden $p_i-1$. Consideremos el caracter $\chi_i:=
 \psi_i^{(p_i-1)/n_i}$. Puesto que $2\big|\frac{p_i-1}{n_i}$ se sigue
 que $\chi_i(-1)=1$ para toda $1\leq i\leq r$.
 
 Sea $p_{r+1}$ un primo impar distinto a $p_1,\ldots,p_r$ y tal que
 $p_{r+1}\equiv 1\bmod n_1\cdots n_r$. Ahora sea $\chi_{r+1}$
 cualquier caracter de conductor $p_{r+1}$ y tal que $\chi_{r+1}(-1)=
 -1$ y tal que $n_1\cdots n_r|o(\chi_{r+1})$. Tal caracter existe,
 por ejemplo, podemos tomar $\chi_{r+1}$ un generador de 
 $\widehat{U_{p_{r+1}}}$. La condici\'on $\chi_{r+1}(-1)=-1$ la
 usaremos para construir $K$ imaginario.
 
 Sea $\chi:=\chi_1\chi_2\cdots\chi_r\chi_{r+1}$. Sea $K$ el campo
 asociado a $\chi$, o m\'as precisamente, al grupo $X=\langle\chi
 \rangle$. Se sigue que $K/{\ma Q}$ es una extensi\'on
 c{\'\i}clica pues $\Gal(K/{\ma Q})\cong \widehat{\langle\chi\rangle}
 \cong\langle\chi\rangle$.
 
 Ahora bien, $\chi(-1)=\chi_1(-1)\cdots\chi_r(-1)\chi_{r+1}(-1)=-1$
 lo cual implica que $K$ es imaginario y en particular toda 
 extensi\'on de $K$ es no ramificada en los primos infinitos.
 
 Sea $Y:=\langle \chi_1,\ldots,\chi_r,\chi_{r+1}\rangle$ y sea $L$
 el campo asociado a $Y$. Entonces $L/{\ma Q}$ es una extensi\'on
 abeliana y se sigue que $L$ y $K$ satisfacen (c).
 
 Puesto que $\f{\chi_i}=p_i$ se sigue que $Y_{p_i}=X_{p_i}=
 \langle\chi_i\rangle$, $1\leq i\leq r+1$ y $Y_p=X_p=\{1\}$ para
 todo primo $p\neq p_1,\ldots,p_r,p_{r+1}$. Se sigue del
 Teorema \ref{T12.3.3} que $L/K$ es no ramificada en ning\'un
 primo pues
 \[
 e_{p_i}(L/{\ma Q})=|Y_{p_i}|=|\langle\chi_{p_i}\rangle|,\quad
 e_{p_i}(K/{\ma Q})=|X_{p_i}|=|\langle\chi\rangle_{p_i}|=
| \langle\chi_{p_i}\rangle|,
 \]
 $1\leq i\leq r+1$.
 
 Ahora bien,  puesto que $\chi=\chi_1\cdots\chi_r\chi_{r+1}$, 
 $\chi_{r+1}=\chi_1^{-1}\cdots\chi_r^{-1}\chi$ y por lo tanto
 \begin{equation}\label{Ec12.4.1.3}
 Y=\langle\chi_1,\ldots,\chi_r,\chi_{r+1}\rangle=\langle\chi_1,
 \ldots,\chi_r,\chi\rangle
 \end{equation}
 de donde se obtiene
 \begin{equation}\label{Ec12.4.1.4}
 \Gal(L/K)\cong \widehat{\Gal(L/K)}\cong\frac{\widehat{\Gal
 (L/{\ma Q})}}{\Gal(L/K)^{\perp}}\cong \frac{Y}{X}=\frac{Y}{
 \langle\chi\rangle}.
 \end{equation}
 
 Consideremos ahora los mapeos naturales
\[
\langle\chi_1,\ldots,\chi_r\rangle\stackrel{i}{\longto} X\stackrel{\pi}{
\longto} Y/\langle\chi\rangle
\]
y sea $\varphi=\pi\circ i$ y de (\ref{Ec12.4.1.3}) se sigue que 
$\varphi$ es suprayectiva.

Sea ahora $\chi_1^{\alpha_1}\cdots\chi_r^{\alpha_r}\in\ker\varphi$,
esto es, existe $\alpha$ tal que $\chi_1^{\alpha_1}
\cdots\chi_r^{\alpha_r}=\chi^{\alpha}=
\chi_1^{\alpha}\cdots\chi_r^{\alpha}
\chi_{r+1}^{\alpha}$. Puesto que todos los caracteres tienen
conductores primos relativos a pares, se sigue que
$\chi_{r+1}^{\alpha}=1$ y que $\alpha_i\equiv \alpha\bmod n_i$,
$n_i=o(\chi_i)$, para $1\leq i\leq r$.

Puesto que $n_1\cdots n_r|o(\chi_{r+1})$, se tiene que $\alpha
\equiv 0\bmod n_i$ para $1\leq i\leq r$ y por tanto
$\alpha_i\equiv 0\bmod n_i$ para toda $i$. Se sigue que
$\chi_1^{\alpha_1}\cdots \chi_r^{\alpha_r}=1$ y que $\ker \varphi
=\{1\}$, es decir $\varphi$ es un isomorfismo entre
$\langle\chi_1,\ldots,\chi_r\rangle$ y $Y/\langle\chi\rangle$.
Por tanto, de (\ref{Ec12.4.1.4}) se sigue que
\begin{align*}
\Gal(L/{\ma Q})&\cong \frac{Y}{\langle\chi\rangle}\cong
\langle\chi_1,\ldots,\chi_r\rangle\cong\bigoplus_{i=1}^r
\langle\chi_i\rangle\cong\\
&\cong\bigoplus_{i=1}^r\langle\psi_i^{(p_i-1)/n_i}\rangle\cong
\bigoplus_{i=1}^r C_{n_i}\cong G. \tag*{$\fin$}
\end{align*}
\end{proof}

Como corolarios al Teorema \ref{T12.4.1.2} y al Teorema de clase
de Hilbert, podemos enunciar:

\begin{corolario}\label{C12.4.1.5}
Sea $G$ un grupo abeliano finito. Entonces existe una extensi\'on
c{\'\i}clica del campo de los n\'umeros racionales $K/{\ma Q}$ tal
que el grupo de clases de ideales de $K$ contiene un subgrupo
isomorfo a $G$. En otras palabras, $G\subseteq I_K$.
\end{corolario}

\begin{proof}{\ }

\begin{window}[0,l, $\left.\begin{array}{c}
\xymatrix{
H_K\ar@{-}[d]_{H}\\ L\ar@{-}[d]_{G}\\K
}\end{array}\right\}{}_{I_K}$
,{}]
Sea $H_K$ el campo de clase de Hilbert de $K$, esto es, $\Gal(
H_K/K)\cong I_K$ y $H_K$ es la m\'axima extensi\'on abeliana
de $K$ no ramificada en ning\'un primo incluyendo a los primos
infinitos.
Sea $L/K$ como en el Teorema \ref{T12.4.1.2}. Entonces $L\subseteq
H_K$ y si $H:=\Gal(H_K/L)$, entonces $\Gal(L/K)\cong I_K/H
\cong G$. As{\'\i}, $I_K$ tiene un grupo cociente 
isomorfo a $G$. Por la
Proposici\'on \ref{P12.14} se sigue que $I_K$ tiene un subgrupo
isomorfo a $G$. $\fin$
\end{window}
\end{proof}

\begin{proposicion}\label{P12.4.1.6}
Sea $L/K$ una extensi\'on de campos num\'ericos tal que no
existe ninguna subextensi\'on $E/K$ con $E\subseteq L$, $E/K$
no ramificada en ning\'un primo incluyendo a los primos infinitos y
tal que $\Gal(E/K)$ es un grupo abeliano. Entonces $h_K|h_L$,
donde en general $h_{\ast}$ denota al n\'umero de clase
del campo $\ast$.

De hecho, existe un subgrupo $M$ de $I_L$ tal que $I_L/M\cong I_K$.
\end{proposicion}

\begin{proof}{\ }

\begin{window}[1,l,\xymatrix{
L\ar@{-}[r]\ar@{-}[d]&L H_K\ar@{-}[r]\ar@{-}[d]&H_L\\
K\ar@{-}[r]&H_K
},{}]
Sea $H_K$ la m\'axima extensi\'on abeliana de $K$ no ramificada.
Entonces $\Gal(H_K/K)\cong I_K$. Ahora, por hip\'otesis se sigue
que $H_K\cap L=K$ pues $H_K\cap L$ es una extensi\'on 
abeliana de $K$ no ramificada y contenida en $L$. En particular
$[LH_K:L]=[H_K:K]$.
Ahora bien $LH_K$ es una extensi\'on abeliana de $L$ no ramificada
por lo que $LH_K\subseteq H_L$. En particular
\end{window}
\[
[LH_K:L]=[H_K:K]=|I_K|=h_K|h_L=[H_L:L]. \tag*{$\fin$}
\]
\end{proof}

\begin{definicion}\label{D12.4.17} Un campo num\'erico $K$,
se llama {\em totalmente real\index{campo!totalmente 
real}\index{totalmente real!campo $\sim$}}
si todos sus encajes en ${\ma C}$ caen en ${\ma R}$.
Ahora $K$ se llama {\em totalmente imaginario\index{totalmente
imaginario}} si ninguno de sus encajes est\'a contenido en
${\ma R}$.

Un campo $K$ se llama de tipo $\MC$\index{campo!de tipo $\MC$}
($\MC$ significa {\em multiplicaci\'on compleja\index{multiplicaci\'on
compleja}}) si es un campo totalmente imaginario que es una
extensi\'on cuadr\'atica de un campo totalmente real $K^+$.
\end{definicion}

\begin{ejemplo}\label{Ej12.4.1.7}
Para $n\geq 3$, $\cic n{}$ es un campo de tipo $\MC$ con
$\cic n{}^+=\cic n{}\cap {\ma R}={\ma Q}(\zeta_n+\zeta_n^{-1})$.
\end{ejemplo}

\begin{corolario}\label{C12.4.1.7}
Si $K$ es un campo de tipo $MC$ entonces $h_{K^+}|h_K$.
En particular $h_{\cic n{}^+}|h_{\cic n{}}$.
\end{corolario}

\begin{proof}
La extensi\'on $K/K^+$ satisface las hip\'otesis de la Proposici\'on
\ref{P12.4.1.6} pues $K/K^+$ es ramificada en los primos
infinitos. $\fin$
\end{proof}

\begin{corolario}\label{C12.4.1.8}
Si $n|m$ entonces $h_{\cic n{}}|h_{\cic m{}}$.
\end{corolario}

\begin{proof}
Sea $E$ un campo tal que $\cic n{}\subsetneqq E\subseteq \cic m{}$.
Sea $U_n=\prod_{p|n}U_{p^{\alpha_p}}$ donde $n=\prod_{p|n} p^{
\alpha_p}$. Por tanto $X_p=\widehat{U_{p^{\alpha_p}}}$. Sea $Y$
el grupo de caracteres de Dirichlet asociado al campo $E$. Entonces
$\widehat{U_n}\subsetneqq Y$, por lo que existe un primo $p$ tal que
$Y_p\supsetneqq X_p$ y en particular, $p$ es ramificado en $E/\cic n{}$.
El resultado se sigue de la Proposici\'on \ref{P12.4.1.6}. $\fin$
\end{proof}

%% file: Capitulo7.tex
\chapter{Series $L$ de Dirichlet}\label{Ch13}

\section{Teorema de Dirichlet}\label{S13.1}

En esta secci\'on nos proponemos probar que si $a,b\in{\ma Z}$
son primos relativos, entonces existen una infinidad de primos
de la forma $a+nb$ con $n\in{\ma N}$.

Claramente la condici\'on de que $a$ y $b$ sean primos relativos
es necesaria pues si existe $d>1$ tal que $d|a$ y $d|b$ entonces
$d|a+nb$ y en este caso a lo m\'as podr{\'\i}a haber un
primo de la forma $a+nb$, a saber, $d$.

El Corolario \ref{C7.3} nos provee de una demostraci\'on
elemental para el caso $a=1$ y $b=n\in{\ma N}$ arbitrario.
E. Landau \cite{Lan09} da una 
demostraci\'on elemental para el caso $a=-1$ y $b=n$ arbitrario.
Mediante t\'ecnicas completamente elementales, se pueden probar
numerosos casos particulares del Teorema de Dirichlet:
primos de la forma $3n+1,3n+2,4n+1,4n+3, 8n+1,8n+3,8n+5,
8n+7,6n+5$, etc. El caso $1+tn$ tiene muchas demostraciones
elementales, algunas m\'as elementales que la que dimos en
el Corolario \ref{C7.3}. En el Teorema \ref{T13.1} indicamos
otra demostraci\'on, mucho m\'as b\'asica que la del Corolario
\ref{C7.3} pero no daremos todos los detalles.

\begin{teorema}\label{T13.1}
Sea $n\in{\ma N}$ cualquiera. Entonces hay una infinidad de
primos de la forma $p\equiv 1\bmod n$.
\end{teorema}

\begin{proof}
Sea $S$ un conjunto y sea $f\colon S\to S$ una
funci\'on. Sea $T_n:=\{s\in S\mid f^{(n)}(s)=s\}$, $n\in{\ma N}$.
Para $s\in T_n$ sea $d$ el m{\'\i}nimo n\'umero natural tal
que $f^{(d)}(s)=s$. Entonces $d\leq n$. Escribamos
$n=qd+r$ con $0\leq r\leq d-1$. Entonces $ f^{(r)}(s)=s$
Puesto que $r\leq d-1$ se sigue que
$r=0$ y $d|n$.

Definimos el conjunto $P_d:=\{s\in T_n\mid o(s)=d\}$ y
para $s\in P_n$, se tiene que $f^{(0)}(s),\ldots,
f^{(n-1)}(s)$ son todos distintos y en particular $n||P_n|$.

Adem\'as $T_n=\cup_{d|n}P_d$ y $P_{d_1}\cap P_{d_2}
=\emptyset$ para $d_1\neq d_2$. Se sigue que
$|T_n|=\sum_{d|n}|P_d|$. Por la f\'ormula de inversi\'on
de M\"obius\index{M\"obius!f\'ormula de 
inversi\'on de $\sim$},
obtenemos que $|P_n|=\sum_{d|n}\mu\big(
\frac{n}{d}\big)|T_d|$, donde $\mu(m)=
\begin{cases}
1&\text{si $m=1$}\\(-1)^r&\text{si $m=p_1\cdots p_r$}\\0&\text{
en otro caso}
\end{cases}$.

En particular $n|\sum_{d|n}\mu\big(
\frac{n}{d}\big)|T_d|$. 

Sean $a$ y $n$ enteros mayores
que uno y sea $n=p_1^{\alpha_1}\cdots p_r^{\alpha_r}$
la descomposici\'on de $n$ en primos. Sea $q$ un divisor 
com\'un de
\[
\frac{a^n-1}{a^{n/p_1}-1},\cdots, \frac{a^n-1}{a^{n/p_r}-1}.
\]
Para cualesquiera enteros $\alpha_0,\ldots, \alpha_{n-1}$ con
$0\leq \alpha_i\leq a-1$, $1\leq i\leq n$ definimos 
\[
(\alpha_0,\ldots, \alpha_{n-1})_a:=\sum_{i=0}^{n-1}\alpha_i
a^i=\alpha_0+\alpha_1a+\cdots+\alpha_{n-1}a^{n-1}.
\]
Sea $S:=\{\alpha= (\alpha_0,\ldots, \alpha_{n-1})_a\mid
q\nmid \alpha\}$. Si $q=1$ entonces $S=\emptyset$. Definimos
la funci\'on $f\colon S\to S$ dada por 
\[
f(\alpha)=(\alpha_{n-1},\alpha_0,\ldots, \alpha_{n-2})_a=
\alpha_{n-1}+\alpha_0 a+\cdots+\alpha_{n-2}a^{n-1}.
\]

Entonces obtenemos $T_n=S$, y por tanto
\[
|T_n|=|S|=a^n-1-\frac{a^n-1}{q}=\frac{(a^n-1)(q-1)}{q}
\]
lo cual implica que $n|\frac{(a^n-1)(q-1)}{q}$.

Si $n=1$ no hay nada que probar. Supongamos $n>1$ y
sean $\{q_1,\ldots,q_s\}$ primos de la forma $1+nt$. Sea
$n=p_1^{\alpha_1}\cdots p_r^{\alpha_r}$. Consideremos
los polinomios
\[
f_1(x)=\frac{x^n-1}{x^{n/p_1}-1},\ldots, f_r(x)=\frac{x^n-1}{x^{n/p_r}-1}.
\]

Todos estos polinomios tienen como ra{\'\i}z com\'un a $\zeta_n
$ por lo que existe un polinomio $g(x)$ no 
constante con coeficientes enteros y coeficiente l{\'\i}der $1$
tal que $g(x)$ divide a $f_i(x)$ para $i=1,2,\ldots, r$.

Se tiene que $f_i(0)=1$ lo cual implica que $g(0)=\pm 1$.
En particular, si $a$ es cualquier entero, entonces $g(a)\equiv
\pm 1\bmod a$, esto es, $a$ y $g(a)$ son primos relativos.
Adem\'as $\lim_{x\to\infty}g(x)=\infty$. Sea
 $t_0\in{\ma N}$ con $g(a)>1$ para todo $a>t_0$. Definimos
$a:=nt_0q_1\cdots q_s$. Si $s=0$, entonces $a:=nt_0$. 
Si $q$ es cualquier n\'umero primo que divide
a $g(a)>1$ entonces $q\nmid a$
por lo que $q\neq q_i$ para toda $i=1,\ldots,s$
y puesto que $n|(a^n-1)(q-1)$ y $n|a$, se tiene que
$n$ y $a^n-1$ son primos relativos lo cual implica que
$n|q-1$, es decir $q\equiv 1\bmod n$ y este es un nuevo
primo de la forma $1+nt$. $\fin$
\end{proof}

Ya que hemos mencionado que algunos casos particulares del 
Teorema de Dirichlet tienen demostraciones elementales, resulta
que no existe, hasta ahora, una demostraci\'on elemental del
caso general.
En esta secci\'on presentaremos la demostraci\'on del caso general.

\begin{definicion}\label{D13.2}
Sea $\chi\colon{\ma Z}\to{\ma C}$ un caracter de Dirichlet m\'odulo
$k$, es decir, entendemos $\chi(a)=0$ si $\mcd(a,k)>1$. Se
define la {\em $L$--serie de Dirichlet\index{Dirichlet!$L$--serie}}
por
\[
L(s,\chi)=\sum_{n=1}^{\infty}\frac{\chi(n)}{n^s}.
\]
\end{definicion}

A grandes rasgos, la demostraci\'on del Teorema de Dirichlet es
como sigue. Sean $\mcd (a,b)=1$ y consideremos todos los
caracteres de Dirichlet m\'odulo $a$. Entonces la serie $L(s,\chi)$
es una funci\'on anal{\'\i}tica para $s>1$ y se tiene
\[
L(s,\chi)=\prod_p\Big(1-\frac{\chi(p)}{p^s}\Big)^{-1}
\]
donde el producto es sobre todos los n\'umeros primos.

Tomando logaritmos y derivando, se obtiene
\[
-\frac{L'(s,\chi)}{L(s,\chi)}=\sum_p\frac{\chi(p)\ln p}{p^s-\chi(p)}.
\]

Sea $\Lambda\colon{\ma N}\to{\ma C}$ dada por $\Lambda(n)=
\begin{cases} \ln p&\text{si $n=p^c, c\geq 1$}\\0&\text{en otro caso}
\end{cases}$. Entonces
\[
-\frac{L'(s,\chi)}{L(s,\chi)}=\sum_{n=1}^{\infty}\frac{\chi(n)\Lambda(n)}
{n^s}.
\]

\begin{definicion}\label{D13.3}
La funci\'on $\Lambda(n)$ se llama la {\em funci\'on de
von--Mangoldt\index{von--Mangoldt!funci\'on de $\sim$}}.
\end{definicion}

Ahora bien, multiplicando por $\overline{\chi(b)}$ y sumando sobre
todos los caracteres $\delta$ se obtiene
\[
\sum_{n\equiv b\bmod a}\frac{\Lambda(n)}{n^s}=\frac{1}{\varphi(a)}
\sum_{\delta}\overline{\delta(b)}\frac{(-L'(s,\delta))}{L(s,\delta)}.
\]

Ahora cuando $s\to 1^+$ el lado izquierdo es aproximadamente
$\sum\limits_{p\equiv b\bmod a}\frac{\ln p}{p}$ y el lado derecho se va
a infinito lo cual prueba lo que queremos.

Ahora bien, si $\chi=1$, se puede probar que $\lim\limits_{s\to 1^+}
-\frac{L'(s,\chi_0)}{L(s,\chi_0)}=\infty$. Por tanto, para ver que el lado
derecho se va infinito basta probar que $\frac{L'(s,\chi)}{L(s,\chi)}$
permanece acotado cuando $s\to 1^+$ y $\chi\neq 1$, es decir,
el problema central es ver que $L(1,\chi)\neq 0$ para $\chi\neq 1$.
A continuaci\'on presentamos este desarrollo.

\begin{teorema}\label{T13.4}
Para todo caracter de Dirichlet $\chi$ m\'odulo $k$, $L(s,\chi)$ es
anal{\'\i}tica para $s\in{\ma C}$, $\partereal s>1$. Adem\'as, para
esta regi\'on, $L(s,\chi)$ tiene una representaci\'on como producto
de Euler
\[
L(x,\chi)=\prod_p \Big(1-\frac{\chi(p)}{p^s}\Big)^{-1}.
\]
\end{teorema}

\begin{proof}
Se tiene $|\chi(n)|\leq 1$ por lo que $\big|\frac{\chi(n)}{n^s}\big|
\leq \frac{1}{|n^s|}=\frac{1}{n^x}$ donde $x=\partereal s$
lo cual converge para $x>1$. Por tanto $L(s,\chi)$ es absolutamente
convergente para $x>1$.

Ahora para $\varepsilon >0$ y $x\geq 1+\varepsilon$,
$\big(\sum_{n=1}^m\frac{\chi(n)}{n^s}\big)'=-\sum_{n=1}^m
\frac{\chi(n)\ln n}{n^s}$ y $\big|\sum_{n=1}^m\frac{\chi(n)\ln n}{n^s}\big|
\leq \sum_{n=1}^m\frac{\ln n}{n^{1+\varepsilon}}$ la cual converge.

Por el criterio $M$ de Weierstrass se tiene que  $-\sum_{n=1}^{\infty}
\frac{\chi(n)\ln n}{n^s}$ converge absoluta y uniformemente para
$x\geq 1+\varepsilon$, $\varepsilon>0$. Puesto que 
$L(s,\chi)=\sum_{n=1}^{\infty}\frac{\chi(n)}{n^s}$ converge
absoluta y uniformemente para $x\geq 1+\varepsilon$, se sigue
que se puede diferenciar t\'ermino a t\'ermino y
\[
L'(s,\chi)=-\sum_{n=1}^{\infty}\frac{\chi(n)\ln n}{n^s}\quad\text{para}
\quad x>1.
\]

Ahora, por la misma raz\'on que antes, $\sum_{n=1}^{\infty}\frac{
\mu(n)\chi(n)}{n^s}$ converge absolutamente para $x>1$ y
uniformemente para $x\geq 1+\varepsilon$, $\varepsilon>0$.
Aqu{\'\i}, $\mu(n)$ es la funci\'on mu de M\"obius.

Se tiene que $\mu$ satisface $\sum_{d|n}\mu(d)=\begin{cases}
1&\text{si $n=1$}\\ 0&\text{si $n>1$}\end{cases}$
(Lema \ref{L1.2.13}).
Por tanto
\begin{align*}
\sum_{m=1}^{\infty}\frac{\chi(m)}{m^s}\cdot \sum_{n=1}^{\infty}
\frac{\chi(n)\mu(n)}{n^s}&=\sum_{t=1}^{\infty}\sum_{mn=t}
\frac{\chi(m)\chi(n)\mu(n)}{m^sn^s}\\
&=\sum_{t=1}^{\infty}
\frac{\chi(t)}{t^s}\sum_{n|t}\mu(n)=\frac{\chi(1)}{1^s}=1
\end{align*}
es decir, 
\begin{equation}\label{Ec13.4.1}
L(s,\chi)\cdot \sum_{n=1}^{\infty}\frac{\chi(n)\mu(n)}{n^s}=1
\end{equation}
y en particular $L(x,\chi)\neq 0$ para $x:=\partereal s>1$.

Finalmente, para $m>1$, sea $S$ el conjunto de los n\'umeros
naturales $n$ no divisibles por ning\'un primo $p>m$. Entonces
si $p_1,\ldots,p_r$ son todos los primos menores o iguales a $m$
se tiene
\begin{gather*}
\prod_{p\leq m}\Big(1-\frac{\chi(p)}{p^s}\Big)=
\prod_{i=1}^r\Big(1-\frac{\chi(p_i)}{p_i^s}\Big)=
\sum_{\substack{1\leq i_1<i_2<\cdots<i_t\leq r\\0\leq t\leq r}}
\frac{(-1)^t\chi(p_{i_1}\cdots p_{i_t})}{(p_{i_1}\cdots p_{i_t})^s}
\intertext{y}
\begin{align*}
\sum_{n\in S}\frac{\chi(n)\mu(n)}{n^s}&=\sum_{\alpha_1=0}^{\infty}
\cdots \sum_{\alpha_r=0}^{\infty}\frac{\chi(p_1^{\alpha_1}\cdots
p_r^{\alpha_r})\mu(p_1^{\alpha_1}\cdots
p_r^{\alpha_r})}{(p_1^{\alpha_1}\cdots
p_r^{\alpha_r})^s}\\
&= \sum_{\substack{1\leq i_1<i_2<\cdots<i_t\leq r\\0\leq t\leq r}}
\frac{\chi(p_{i_1}\cdots p_{i_t})}{(p_1^{\alpha_1}\cdots
p_r^{\alpha_r})^s}\mu(p_{i_1}\cdots p_{i_t})\\
&=
 \sum_{\substack{1\leq i_1<i_2<\cdots<i_t\leq r\\0\leq t\leq r}}
 \frac{(-1)^t \chi(p_{i_1}\cdots p_{i_t})}{(p_{i_1}\cdots p_{i_t})^s}.
\end{align*}
\intertext{Es decir}
\prod_{p\leq m}\Big(1-\frac{\chi(p)}{p^s}\Big)=
\sum_{n\in S}\frac{\chi(n)\mu(n)}{n^s}=
\sum_{n=1}^m \frac{\chi(n)\mu(n)}{n^s}+
\sum_{\substack{n'\in S\\ n'>m}}\frac{\chi(n')\mu(n')}{(n')^s}.
\end{gather*}

Por (\ref{Ec13.4.1}) se sigue que $\lim\limits_{m\to \infty}
\sum\limits_{n=1}^m \frac{\chi(n)\mu(n)}{n^s}=L(s,\chi)^{-1}$,
$\partereal s>1$. Ahora
\begin{gather*}
\sum_{\substack{n'\in S\\ n'>m}} \Big|\frac{\chi(n')\mu(n')}{(n')^s}
\Big|\leq \sum_{n=m+1}^{\infty}\frac{1}{n^s}\xrightarrow[m\to\infty]{}0
\quad\text{para}\quad \partereal s>1.\\
\intertext{Por tanto $\prod_p\big(1-\frac{\chi(p)}{p^s}\big)=\frac{1}{
L(s,\chi)}$ y}
L(s,\chi)=\prod_p\Big(1-\frac{\chi(p)}{p^s}\Big)^{-1},\quad
\partereal s>1. \tag*{$\fin$}
\end{gather*}
\end{proof}

\begin{definicion}\label{D13.5}
Sea $K$ un campo num\'erico, $[K:{\ma Q}]<\infty$.
La {\em funci\'on zeta de Dedekind}\index{funci\'on zeta de
Dedekind\index{Dedekind!funci\'on zeta de $\sim$}}
de $K$ se define por
\[
\zeta_K(s)=\sum_{{\eu a}}\frac{1}{(N{\eu a})^s}
\]
donde ${\eu a}$ recorre los ideales no cero de ${\cal O}_K$
y $N{\eu a}:=|{\cal O}_K/{\eu a}|$.

En particular si $K={\ma Q}$, obtenemos la {\em funci\'on
zeta de Riemann\index{funci\'on zeta de
Riemann}\index{Riemann!funci\'on zeta de $\sim$}}:
\[
\zeta(s)=\zeta_{{\ma Q}}(s)=\sum_{n=1}^{\infty}\frac{1}{n^s}.
\]
\end{definicion}

Como en el Teorema \ref{T13.4} obtenemos

\begin{teorema}\label{T13.6}
Se tiene que $\zeta_K(s)$ converge absolutamente para
$\partereal s>1$, uniformemente para $\partereal s\geq 1+\varepsilon
$, $\varepsilon>0$ y se tiene
\[
\zeta_K(s)=\prod_{\pK}\Big(1-\frac{1}{(N\pK)^s}\Big)^{-1},
\quad \partereal s>1
\]
donde $\pK$ recorre todos los ideales primos no cero de 
${\cal O}_K$. $\fin$
\end{teorema}

\begin{proposicion}\label{P13.7}{\ }

\lasa
\item Para $\partereal s>1$ tenemos
\[
-\frac{L'(s,\chi)}{L(s,\chi)}=\sum_{n=1}^{\infty}
\frac{\chi(n)\Lambda(n)}{n^s}
\]
donde $\Lambda(n)$ es la funci\'on de von--Mangoldt.
\item Si $\chi_0=1$ m\'odulo $k$, es decir, $\chi(n)=1$
si $\mcd (n,k)=1$ y $\chi_0(n)=0$ en otro caso, entonces
\[
\lim_{s\to 1}\Big|\frac{L'(s,\chi_0)}{L(s,\chi_0)}\Big|=\infty
\quad\text{y}\quad \lim_{\substack{s\in{\ma R}\\s\to 1^+}}
\frac{L'(s,\chi_0)}{L(s,\chi_0)}=-\infty.
\]
\end{list}
\end{proposicion}

\begin{proof}
Se tiene que $\big|\sum\limits_{n=1}^{\infty}\frac{\chi(n)\Lambda(n)}
{n^s}\big|\leq \sum\limits_{n=1}^{\infty}\frac{\ln n}{n^x}<\infty$ para
$x=\partereal s>1$.

Sea $n=p_1^{\alpha_1}\cdots p_r^{\alpha_r}$ y
\[
\sum_{d|n}\Lambda(d)=\sum_{i=1}^r\sum_{\beta_i=0}^{\alpha_i}
\Lambda(p_i^{\beta_i})=\sum_{i=1}^r\sum_{\beta_i=1}^{\alpha_i}
\ln p_i=\sum_{i=1}^r\ln p_i^{\alpha_i}=\ln n.
\]

Para $\partereal s=x>1$, se tiene:
\begin{align*}
L(s,\chi)\cdot \sum_{n=1}^{\infty}\frac{\chi(n)\Lambda(n)}{n^s}&=
\sum_{m=1}^{\infty}\frac{\chi(m)}{m^s}\sum_{n=1}^{\infty}
\frac{\chi(n)\Lambda(n)}{n^s}\\
&=\sum_{t=1}^{\infty}\frac{1}{t^s}\Big(\sum_{mn=t}\chi(mn)\Lambda
(n)\Big)=\sum_{t=1}^{\infty}\frac{\chi(t)}{t^s}\big(\sum_{n|t}
\Lambda(n)\big)\\
&=\sum_{t=1}^{\infty}\frac{\chi(t)\ln t}{t^s}=-L'(s,\chi)
\end{align*}
de donde se obtiene (a).

Ahora bien, si $\chi_0=1$ m\'odulo $k$, entonces
\begin{align*}
\frac{L'(s,\chi_0)}{L(s,\chi_0)}&=-\sum_{n=1}^{\infty}
\frac{\chi_0(n)\Lambda(n)}{n^s}=-\sum_{\substack{
n=1\\ \mcd(n,k)=1}}^{\infty}\frac{\Lambda(n)}{n^s}\\
&= -\sum_{n=1}^{\infty}\frac{\Lambda(n)}{n^s}+
\sum_{\mcd(n,k)>1}\frac{\Lambda(n)}{n^s}=
-\sum_{n=1}^{\infty}\frac{\Lambda(n)}{n^s}+
\sum_{p|k}\sum_{m=1}^{\infty}\frac{\Lambda(p^m)}{(p^m)^s}\\
&=-\sum_{n=1}^{\infty}\frac{\Lambda(n)}{n^s}+
\sum_{p|k}\sum_{m=1}^{\infty}\frac{\ln p}{(p^m)^s}\\
&=
-\sum_{n=1}^{\infty}\frac{\Lambda(n)}{n^s}+\sum_{p|k}
\ln p\Big(\frac{1}{1-p^s}-1\Big)\\
&=-\sum_{n=1}^{\infty}\frac{\Lambda(n)}{n^s}+
\sum_{p|k}\frac{p^s\ln p}{1-p^s}.
\end{align*}

Como $\{p\mid p|k\}$ es un conjunto finito y $\lim\limits_{s\to 1}
\frac{p^s\ln p}{1-p^s}=\frac{p\ln p}{1-p}<\infty$, para ver que
$\lim\limits_{s\to 1}\big|\frac{L'(s,\chi_0)}{L(s,\chi_0)}\big|=
\infty$ basta ver que $\lim\limits_{s\to 1}\big|\sum_{n=1}^{\infty}\frac{\Lambda(n)}
{n^s}\big|=\infty$.

Se tiene que $\sum\limits_{p\text{\ primo}}\frac{\ln p}{p}\geq\sum
\limits_{p
\text{\ primo}}\frac{1}{p}$. Ahora
\[
\sum_{n=1}^{\infty}\frac{\Lambda(n)}{n}=\sum_{p \text{\ primo}}
\sum_{c=1}^{\infty}\frac{\Lambda(p^c)}{p^c}\geq \sum_{p
\text{\ primo}}\frac{\ln p}{p}\geq \sum_{p \text{\ primo}}\frac{1}{p}.
\]

Veamos que $\sum\limits_{p \text{\ primo}}\frac{1}{p}=\infty$.
Sean $p_1,\ldots,p_r$ los primeros $r$ primos. Se tiene
$p_1=2=2^1=2^{2^0}=2^{2^{1-1}}$, $p_2=3<4=2^2=2^{2^{2-1}}$.
Supongamos que para $r>1$, $p_r<2^{2^{r-1}}$. Puesto que
 alg\'un n\'umero primo $p$ distinto de $p_1,\ldots,p_r$ divide
a $p_1\cdots p_r+1$, se tiene que $p_{r+1}\leq p_1\cdots p_r+1<
2\cdot 2^{2^2}\cdot 2^{2^3}\cdots 2^{2^{r-1}}+1<2^{2^r}$.

Si acaso $\sum\limits_{p \text{\ primo}}\frac{1}{p}=
\sum\limits_{i=1}^{\infty}\frac{1}{p_i}$ convergiese, existir{\'\i}a $n_0\in
{\ma N}$ tal que $\sum\limits_{n=n_0+1}^{\infty}\frac{1}{p_i}<\frac{1}
{2}$.

Para $x\in{\ma N}$ definamos $Q_{n_0}(x)$ el n\'umero de
n\'umeros naturales menores o iguales a $x$ y que no son divisibles
por ning\'un $p_n$ con $n\geq n_0+1$. Ahora, dado $p$ primo,
 el n\'umero de enteros positivos $m\leq x$ y divisibles por $p$
 es menor o igual a $x/p$. Por lo tanto
 \[
 x-Q_{n_0}<\frac{x}{p_{n_0+1}}+\frac{x}{p_{n_0+2}}+\cdots
 \leq x\sum_{n=n_0+1}^{\infty}\frac{1}{p_n}<\frac{x}{2}.
 \]
 Por lo tanto $x/2<Q_{n_0}(x)$.
 
 Sea $m<x$ y $m$ no es
 divisible por ning\'un primo $p_n$
 con $n\geq n_0+1$. Escribamos $m=s^2t$ donde $t$ es
 libre de cuadrados, es decir, $t=2^{a_1}3^{a_2}\cdots p_{n_0}^{
 a_{n_0}}$ con $a_i\in\{0,1\}$, esto es, $t$ tiene a lo m\'as $2^{n_0}$
 selecciones y hay a lo m\'as $\sqrt{m}<\sqrt{x}$ selecciones para
 $s$. Por lo tanto
$\frac{x}{2}<Q_{n_0}(x)<2^{n_0}\sqrt{x}$ pero esto es imposible
pues $\lim\limits_{x\to\infty}\frac{x/2}{2^{n_0}\sqrt{x}}=\infty$.
Se sigue que $\sum\limits_{p \text{\ primo}}\frac{1}{p}=\infty$.

Entonces $\sum\limits_{n=1}^{\infty}\frac{\Lambda(n)}{n}$ diverge.
Dado $M\in {\ma R}$, $M>0$, existe $m$ tal que $\sum\limits_{n
=1}^m\frac{\Lambda(n)}{n}>M$. Se sigue que existe $\varepsilon=
\varepsilon(M)>0$ tal que para $1<x<1+\varepsilon$, 
$\sum\limits_{n=1}^m \frac{\Lambda(n)}{n^x}>M$, por tanto
$\sum\limits_{n=1}^{\infty}\frac{\Lambda(n)}{n^x}>M$. De esto
se sigue (b). $\fin$
\end{proof}

Damos tres resultados m\'as antes de probar el Teorema de
Dirichlet.

\begin{lema}\label{L13.8}
Si $n\geq m\geq 1$ y $x\neq 1$, entonces si $\varphi$ es la
funci\'on fi de Euler, se tiene
\[
\big|\sum_{i=m}^{t}\chi(i)\big|<\frac{\varphi(k)}{2}.
\]
\end{lema}

\begin{proof}
Puesto que $\chi\neq 1$, se tiene que $\sum\limits_{\substack{
a=1\\ \mcd(a,k)=1}}^k\chi(a)=0$ donde $\chi$ est\'a definido m\'odulo
$k$. Puesto que $\chi(a)=0$ para $\mcd(a,k)>1$, se tiene que
$\sum_{a=1}^k\chi(a)=0$ y $\varphi(k)$ de estos t\'erminos satisfacen
$|\chi(a)|=1$ (estos son los que $\mcd (a,k)=1$).

Entonces si $n-m=tk+r$ con $0\leq r\leq k-1$, se tiene
\[
\sum_{i=m}^n\chi(i)=\sum_{i=m}^{m+tk}\chi(i)+\sum_{i=m+tk+1}^{
m+tk+r}\chi(i)=\sum_{i=m+tk+1}^{m+tk+r}\chi(i)
\]
es decir, podemos suponer que $n-m\leq r\leq k-1$.

Ahora bien, si en la suma $\sum_{i=m}^n\chi(i)$ hay a lo m\'as
$\varphi(k)/2$ t\'erminos con $|\chi(i)|=1$, se tiene $\big|
\sum_{i=m}^n\chi(i)\big|\leq \sum_{i=m}^n|\chi(i)|<
\varphi(k)/2$.

Si en la suma hay m\'as de $\varphi(k)/2$ t\'erminos con $|\chi(i)|=1$,
entonces puesto que $n\leq m+k-1$, se tiene
\[
\Big|\sum_{i=m}^n\chi(i)\Big|=\Big|\sum_{i=m}^{m+k-1}\chi(i)-
\sum_{n+1}^{m+k-1}\chi(i)\Big|=\Big|\sum_{i=n+1}^{m+k-1}\chi(i)\Big|
\leq \sum_{i=n+1}^{m+k-1}|\chi(i)|<\frac{\varphi(k)}{2}
\]
pues en esta \'ultima suma hay menos de $\varphi(k)/2$ t\'erminos
con $|\chi(i)|=1$. $\fin$

\end{proof}

\begin{lema}\label{L13.9} Si $\chi$ es cualquier caracter, $\chi_0=1$
y $s>1$, entonces tenemos 
\[
|L(s,\chi_0)|^3|L(s,\chi)|^4|L(s,\chi^2)|^2\geq 1 \quad (s\in{\ma R}).
\]
\end{lema}

\begin{proof}
Si $x,y\in{\ma R}$ y $0<x<1$, se tiene
\begin{equation}\label{Ec13.1.1}
(1-x)^3|1-xe^{iy}|^4|1-xe^{2iy}|^2<1
\end{equation}
Para verificar (\ref{Ec13.1.1}), hagamos $a:=\cos y$, $2a^2-1=
\cos 2y$. Entonces $|1-xe^{iy}|^4=(1+x^2-2xa)^2$, $|1-xe^{2iy}|^2
=(1+x)^2-4xa^2$. Si definimos $\ell(a)=(1+x^2-2xa)^2
((1+x)^2-4xa^2)$, $0<x<1$, $-1\leq a\leq 1$, $\ell'(a_0)=0$ para
$|a_0|<1 \iff a_0=\frac{1+x^2}{2(1+x)^2}$ y $\ell''(a)>0$ por lo que
$a_0$ es un m{\'\i}nimo local para $\ell(a)$, $|a|\leq 1$. Ahora bien
$\ell(-1)>\ell(a_0), \ell(1)$ por lo que el m\'aximo en $-1\leq a\leq 1$
es $\ell(-1)=(1+x)^4(1-x)^2$. Por lo tanto
\[
(1-x)^3|1-xe^{iy}|^4|1-xe^{2iy}|^2\leq (1-x)^3(1+x)^4(1-x)^2
=(1-x^2)^4(1-x)<1.
\]

Ahora, si $p$ es un n\'umero primo que no divide a $k$,
tenemos $\chi(p)=e^{iy}$ para alg\'un $y$. Sea $x:=1/p^s$.
Aplicamos (\ref{Ec13.1.1}), se tiene
\[
\Big|\Big(1-\frac{\chi_0(p)}{p^s}\Big)\Big|^3\Big|\Big(1-\frac{\chi(p)}{p^s}\Big)
\Big|^4\Big|\Big(1-\frac{\chi^2(p)}{p^s}\Big)\Big|^2\leq 1.
\]

Tomando la expresi\'on sobre todos los n\'umeros primos y
usando la representaci\'on en producto de Euler, obtenemos la
igualdad deseada. $\fin$
\end{proof}

\begin{lema}\label{L13.10}
Si $\chi\neq 1$ tenemos $|L'(s,\chi)|<\varphi(k)$ para $s\in{\ma R}$,
$s\geq 1$.
\end{lema}

\begin{proof}
De la demostraci\'on del Teorema \ref{T13.4}, tenemos
\[
L'(s,\chi_0)=-\sum_{n=1}^{\infty}\frac{\chi(n)\ln n}{n^s}
\quad \text{para}\quad s>1.
\]

Para $t\geq 3$, $t\in{\ma R}$, $f(t)=\frac{\ln t}{t^s}$ es
una funci\'on decreciente. Por tanto, por el Lema \ref{L13.8} se
tiene
\[
\Big|\sum_{i=m}^{n}\frac{\chi(i)\ln i}{i^s}\Big|\leq \frac{\varphi(k)}{2}
\frac{\ln m}{m^s}\leq \frac{\varphi(k)}{2}\frac{\ln m}{m}.
\]
En particular la serie $\sum_{n=1}^{\infty}\frac{\chi(n)\ln n}{n^s}$
converge absolutamente para $s\geq 1$. Se sigue que para $s\geq
1$
\begin{align*}
|L'(s,\chi)|&=\Big|\sum_{n=1}^{\infty}\frac{\chi(n)\ln n}{n^s}\Big|
\leq \frac{\chi(2)\ln 2}{2^s}+\Big|\sum_{n=3}^{\infty}
\frac{\chi(n)\ln n}{n^s}\Big|\\
&\leq \frac{\chi(2)\ln 2}{2}+\frac{\varphi(k)}{2}\frac{\ln 3}{3}<
\frac{1}{2}+\frac{\varphi(k)}{2}\leq \varphi(k). \tag*{$\fin$}
\end{align*}
\end{proof}

Un resultado fundamental que necesitamos es:

\begin{teorema}\label{T13.11}
Si $\chi\neq 1$, $L(1,\chi)\neq 0$ y $\frac{L'(s,\chi)}{L'(s,\chi_0)}$
est\'a acotada para $s>1$ y $\chi_0=1$.
\end{teorema}

\begin{proof}
Puesto que $\big|\sum_{i=m}^n\chi(i)\big|\leq \frac{\varphi(k)}{2}$,
si $\chi\neq 1$ y $s\in{\ma R}$, $s>1$, entonces se sigue que para
toda $n$, $\big|\sum_{i=1}^{n}\frac{\chi(i)}{i^s}\big|\leq
\sum_{i=1}^n\chi(i)\leq \frac{\varphi(k)}{2}$ y por tanto $|L(s,\chi)|
< \varphi(k)$.

Si $\chi$ es no real, esto es, $\chi({\ma Z})\nsubseteq {\ma R}$,
entonces $\chi^2\neq \chi_0$ pues en caso de que $\chi^2=\chi_0$,
$\chi({\ma Z})\subseteq \{0,1,-1\}$. En particular $|L(s,\chi^2)|<
\varphi(k)$. M\'as a\'un, tomando $1<s<2$, se tiene
\[
L(x,\chi_0)=\sum_{\substack{n=1\\ \mcd (n,k)=1}}^{\infty}
\frac{1}{n^s}\leq \sum_{n=1}^{\infty}\frac{1}{n^s}<1+\int_1^{\infty}
\frac{dx}{x^s}=\frac{s}{s-1}<\frac{2}{s-1}.
\]

Por el Lema \ref{L13.9}, se tiene
\[
|L(s,\chi)|\geq \frac{1}{|L(s,\chi_0)|^{3/4}}\frac{1}{|L(x,\chi^2)|^{2/4}}>
\frac{(s-1)^{3/4}}{2^{3/4}}\frac{1}{\sqrt{\varphi(k)}}>\frac{(s-1)^{3/4}}
{2\sqrt{\varphi(k)}}.
\]

En caso de que $L(1,\chi)=0$, se tendr{\'\i}a para $s>1$
\[
|L(s,\chi)|=|L(s,\chi)-L(1,\chi)|=\Big|\int_{1}^s L'(x,\chi)dx\Big|<
\varphi(k)(s-1).
\]
En particular para $1<s<2$ se seguir{\'\i}a que
$(s-1)^{1/4}>\frac{1}{2\varphi(k)^{3/2}}$ lo cual no se cumple:
por ejemplo para $s=1+\frac{1}{16\varphi(k)^{3/2}}\in (1,2)$.

Se sigue que $L(1,\chi)\neq 0$ para $\chi$ no real.

Ahora consideremos $\chi$ real, $\chi\neq
\chi_0$. Sea $f\colon{\ma N}
\to{\ma R}$ dada por $f(n)=\sum_{d|n}\chi(d)$. Sean $p_1,\ldots,
p_r$ los primos que dividen a $n$ y que no dividen a $k$ y sea
$\alpha_i$ tal que $p_i^{\alpha_i}|n$ y $p_i^{\alpha_i+1}\nmid n$,
$1\leq i\leq r$. Entonces $\chi(p_i)=(-1)^{\varepsilon_i}$,
$\varepsilon_i\in\{0,1\}$, $\chi(p_i^{\alpha_i})=(-1)^{\alpha_i
\varepsilon_i}$.

Escogemos $p_1,\ldots,p_t$ con $\varepsilon_i=1$ y $p_{t+1},
\ldots, p_r$ con $\varepsilon_i=0$. Entonces
\begin{align*}
f(n)&=\sum_{i=1}^r\sum_{\beta_i=0}^{\alpha_i}\chi(p_1^{\beta_1}
\cdots p_r^{\beta_r})=\sum_{i=1}^t\sum_{\substack{\beta_i=0\\
\text{alg\'un\ }\beta_i>0}}^{\alpha_i}\chi(p_1^{\beta_1}\cdots p_t^{\beta_t})+
\sum_{j=t+1}^r\sum_{\beta_j=0}^{\alpha_j}1\\
&=\sum_{i=1}^t\sum_{\beta_i=0}^{\alpha_i}(-1)^{\beta_1+\cdots
+\beta_t}+\big[(\alpha_{t+1}+1)\cdots(\alpha_r+1)-1\big]\\
&\geq
\sum_{i=1}^t\sum_{\beta_i=0}^{\alpha_i}(-1)^{\beta_1+\cdots+\beta_t}.
\end{align*}

Podemos hacer inducci\'on en $t$ para ver que la suma es mayor
o igual a $0$. Si $t=0$, la suma es $1=:h_0>0$.
Si $t>0$ y suponemos que $h_{t-1}:=\sum_{i=1}^{t-1}\sum_{\beta_i}^{
\alpha_i}(-1)^{\beta_1+\cdots+\beta_{t-1}}\geq 0$, entonces
$h_t=\sum_{\beta_t=0}^{\alpha_t}(-1)^{\beta_t} h_{t-1}=
\begin{cases} h_{t-1}&\text{si $\alpha_t$ es par}\\
0&\text{si $\alpha_t$ es impar}\end{cases}\geq 0$.

Ahora si $n=c^2$ es un cuadrado, $h_t=h_{t-1}>0$ y por tanto se
tiene $f(n)\geq 0$ para toda $n\geq 1$ y $f(n)\geq 1$ si $n=c^2$
es un cuadrado. Sean
\begin{equation}\label{Ec13.1.-1}
m= (4\varphi(k))^6\quad\text{y}\quad z=
\sum_{n=1}^m 2(m-n)f(n).
\end{equation}
Entonces
\[
z=\sum_{n=1}^m 2(m-n)\sum_{d|n}\chi(d)\igual_{\substack{\uparrow\\
\frac{n}{d}=v}}\sum_{n=1}^m\sum_{d|n} 2(m-vd)\chi(d)
=\sum_{vd\leq m}2(m-vd)\chi(d).
\]

Se tiene que, puesto que $f(n)\geq 0$, $n\geq 1$, $f(n)\geq 1$
si $n=c^2$
\begin{align}\label{Ec13.1.3}
z&\geq \sum_{v=1}^{\sqrt{m}}2(m-v^2)\geq \sum_{v=1}^{\sqrt{m}/2}
2(m-v^2)\geq \sum_{v=1}^{\sqrt{m}/2}2\big(m-\frac{m}{4}\big)
\nonumber\\
&= \frac{\sqrt{m}}{2} 2\big(\frac{3m}{4}\big)=\frac{3}{4}m^{3/2}
=\frac{3}{4}(4\varphi(k))^9.
\end{align}

Separamos $z$ en dos sumandos:
\begin{align}\label{Ec13.1.4}
z_1&=\sum_{d=1}^{m^{1/3}}\sum_{m^{3/2}<v\leq m/d}
2(m-vd)\chi(d),\nonumber\\
z_2&= \sum_{v=1}^{m^{2/3}}\sum_{0<d\leq m/v}
2(m-vd)\chi(d),\nonumber \\
z&=z_1+z_2.
\end{align}

Sean $z(n)$ una funci\'on valuada en ${\ma C}$,
$c\in{\ma N}$ y $t\geq c$. Sea $w(t):=\sum_{n=c}^t z(n)$ y 
definimos $w(c-1)=0$. Para $d\geq c$ sean $\mu_d:=
\max\limits_{c\leq t\leq d}|r(t)|$, 
y $\varepsilon_c\geq \varepsilon_{c+1}
\geq \cdots\geq \varepsilon_d\geq 0$. Entonces
\[
\sum_{n=c}^d \varepsilon_n z_n=\sum_{n=c}^d \varepsilon_n(
w(n)-w(n-1))=\sum_{n=c}^{d-1}
w(n)(\varepsilon_n-\varepsilon_{n-1})+w(d)\varepsilon_d.
\]
En particular
\begin{equation}\label{Ec13.1.2}
\Big|\sum_{n=c}^d \varepsilon_n z(n)\Big|\leq \mu_d\Big(
\sum_{n=c}^{d-1}(\varepsilon_n-\varepsilon_{n-1})+\varepsilon_d
\Big)=v\varepsilon_c.
\end{equation}

Aplicamos lo anterior a $\sum_{n=c}^d\chi(n)$. En este caso tenemos
$\big|\sum_{n=c}^d\chi(n)\big|\leq \frac{\varphi(k)}{2}$ y por tanto,
con $\varepsilon_n=\frac{1}{n^s}$ y $s\in{\ma R}$, $s>1$, 
se obtiene
\[
\Big|\sum_{n=c}^d\frac{\chi(n)}{n^s}\Big|\leq \frac{\varphi(k)}{2}
\frac{1}{c^s}\leq \frac{\varphi(k)}{2c}.
\]
Ahora bien, aplicando (\ref{Ec13.1.2}) a $z_1$, obtenemos
\begin{equation}\label{Ec13.1.5}
z_1\leq \Big|\sum_{d=1}^{m^{1/3}}\sum_{m^{2/3}<v\leq m/d}
2(m-dv)\chi(v)\Big|\leq \sum_{d=1}^{m^{1/3}}2m\frac{\varphi(k)}{2}
=m^{4/3}\varphi(k).
\end{equation}

Sea ahora $\theta=\frac{m}{v}-\big[\frac{m}{v}\big]$, donde $[x]$
denota la parte entera de $x\in{\ma R}$. Entonces $0\leq \theta
<1$ y se tiene
\begin{align*}
\sum_d (2m-2dv)&= 2m\sum_d 1-v\sum_d 2d =
2m\Big[\frac{m}{v}\Big]-v\Big[\frac{m}{v}\Big]\Big(\Big[\frac{m}{v}
\Big]+1\Big)\\
&= 2m\Big(\frac{m}{v}-\theta\Big)-v\Big(\Big(\frac{m}{v}-\theta\Big)^2
+\frac{m}{v}-\theta\Big)\\
&=\frac{2m^2}{v}-2m\theta-v\Big(\frac{m^2}{v^2}-2\theta\frac{m}{v}
+\theta^2+\frac{m}{v}-\theta\Big)\\
&=\frac{m^2}{v}-m+v(\theta-\theta^2).
\end{align*}

Puesto que $0\leq \theta<1$, $0\leq \theta-\theta^2\leq \theta<1$,
y por tanto
\begin{align*}
z_2&=m^2\sum_{v=1}^{m^{2/3}} \frac{\chi(v)}{v}-m\sum_{v=1}^{m^{
2/3}}\chi(v)+ \sum_{v=1}^{m^{2/3}}\chi(v) v(\theta-\theta^2)\\
&\leq m^2\Big(L(1,\chi)-\sum_{v=m^{2/3}+1}^{\infty}
\frac{\chi(v)}{v^s}\Big)+m\frac{\varphi(k)}{2}+m^{2/3}\sum_{v=1}^{
m^{2/3}}1.
\end{align*}
Ahora si aplicamos la desigualdad
\[
\Big|\sum_{n=c}^d\frac{\chi(n)}{n^s}\Big|\leq \frac{\varphi(k)}{2}
\frac{1}{c^s}\leq \frac{\varphi(k)}{2c}
\]
y tomando $c=m^{2/3}+1$, $v\to\infty$ se obtiene
\begin{align}\label{Ec13.1.6}
z_2&< m^2L(1,\chi)+m^2\frac{\varphi(k)}{2}\frac{1}{m^{2/3}}+m^{4/3}
\frac{\varphi(k)}{2}+m^{4/3}\varphi(k)\nonumber\\
&= m^2L(1,\chi)+m^{4/3}\varphi(k)\Big(\frac{1}{2}+\frac{1}{2}+1\Big)
=m^2L(1,\chi)+2m^{4/3}\varphi(k).
\end{align}

De (\ref{Ec13.1.-1}), (\ref{Ec13.1.3}), (\ref{Ec13.1.4}), (\ref{Ec13.1.5})
y (\ref{Ec13.1.6}) obtenemos
\begin{align*}
\frac{3}{4}(4\varphi(k))^9&\leq z=z_1+z_2\leq m^{4/3}\varphi(k) +
m^2L(1,\chi)+2 m^{4/3}\varphi(k)\\
&=m^2L(1,\chi)+3m^{4/3}\varphi(k)=m^2L(1,\chi)+3(4\varphi(k))^8
\varphi(k)\\
&=m^2L(1,\chi)+\frac{3}{4}(4\varphi(k))^9.
\end{align*}
En particular $m^2L(1,\chi)>0$ y por tanto $L(1,\chi)>0$. Se concluye
$L(1,\chi)\neq 0$.

Esto prueba la primera parte del teorema. La segunda parte
se sigue pues ya que $L(1,\chi)\neq 0$, $\frac{1}{L(s,\chi)}$
est\'a acotada para $s\geq 1$, $s\in{\ma R}$. Por el Lema
\ref{L13.10} se tiene que $L'(x,\chi)$ est\'a acotado para
$s\geq 1$, $s\in{\ma R}$. $\fin$.
\end{proof}

La \'ultima parte que se necesita para el Teorema de Dirichlet es:

\begin{teorema}\label{T13.12}
Sean $\mcd (c,k)=1$, $n>0$. Entonces para $s\in{\ma R}$, $s>1$
tenemos:
\[
-\frac{1}{\varphi(k)}\sum_{\chi\bmod k}\overline{\chi(c)}\frac{L'(s,\chi)}
{L(s,\chi)}=\sum_{n\equiv c\bmod k}\frac{\Lambda(n)}{n^s}.
\]
\end{teorema}

\begin{proof}
De la Proposici\'on \ref{P13.7} se tiene que
\begin{gather*}
-\frac{L'(s,\chi)}{L(s,\chi)}=\sum_{n=1}^{\infty}\frac{\chi(n)\Lambda(n)}
{n^s}.\\
\intertext{Se sigue que}
\begin{align*}
-\sum_{\chi\bmod k}\overline{\chi(a)}\frac{L'(s,\chi)}{L(s,\chi)}&=
-\sum_{\chi\bmod k}\frac{1}{\chi(c)}\sum_{n=1}^{\infty}
\frac{\chi(n)\Lambda(n)}{n^s}\\
&=\sum_{n=1}^{\infty}\frac{\Lambda(n)}
{n^s}\sum_{\chi\bmod k}\frac{1}{\chi(c)}\chi(n).
\end{align*}
\intertext{Ahora}
\begin{align*}
\sum_{\chi\bmod k}\frac{1}{\chi(c)}\chi(n)&=\sum_{\chi\bmod k}
\chi(c^{-1}n)=\sum_{\chi\in\widehat{U_k}}(\widehat{c^{-1}n})(\chi)\\
&=\begin{cases} 0&\text{si $\widehat{c^{-1}n}\neq 1$},\\
\varphi(k)&\text{si $\widehat{c^{-1}n}=1$}\end{cases}=
\begin{cases} 0&\text{si $c\not\equiv n\bmod k$},\\
\varphi(k)&\text{si $c\equiv n\bmod k$}\end{cases}.
\end{align*}
\intertext{Por lo tanto}
-\sum_{\chi\bmod k}\overline{\chi(a)}\frac{L'(s,\chi)}{L(s,\chi)}=
\varphi(k)\sum_{n\equiv c\bmod k}\frac{\Lambda(n)}{n^s}. 
\tag*{$\fin$}
\end{gather*}
\end{proof}

\begin{teorema}[Teorema de Dirichlet\index{Dirichlet!teorema
de $\sim$}\index{teorema de Dirichlet}]\label{T13.13}
Si $\mcd (a,b)=1$, $a,b\in{\ma N}$, entonces existen una infinidad
de n\'umeros primos $p$ tales $p\equiv b\bmod a$.
\end{teorema}

\begin{proof}
Consideremos $\widehat{U_a}$. Por el Teorema \ref{T13.12}
con $a=k$ y $b=c$ se tiene
\begin{equation}\label{Ec13.1.7}
-\frac{1}{\varphi(k)}\sum_{\chi\in\widehat{U_a}}\overline{\chi(b)}
\frac{L'(s,\chi)}{L(s,\chi)}=\sum_{n\equiv b\bmod a}\frac{\Lambda(n)}
{n^s}.
\end{equation}
Cuando $s\to 1^+$ el lado izquierdo de (\ref{Ec13.1.7}) se va a 
$\infty$ pues por la Proposici\'on \ref{P13.7}, $-\lim\limits_{s\to 1^+}
\frac{L'(s,\chi_0)}{L(s,\chi_0)}=\infty$ y por el Teorema \ref{T13.11}
los dem\'as $(\varphi(a)-1)$ t\'erminos est\'an acotados. Por lo tanto
el lado derecho de (\ref{Ec13.1.7}) satisface:
\[
\sum_{p\equiv b\bmod a}\frac{\ln p}{p^s}+\sum_{\substack{
p^m\equiv b\bmod a\\ m>1}}\frac{\ln p}{p^{ms}}\xrightarrow[s\to 
1^+]{} \infty.
\]
Por otro lado la segunda suma permanece acotada pues
\begin{align*}
\sum_{n=1}^{\infty}\frac{2\ln n}{n^s}&> \sum_{n=2}^{\infty}
\frac{\ln n}{n(n-1)}\geq \sum_p\frac{\ln p}{p(p-1)}\geq
\sum_{p,m>1}\frac{\ln p}{p^m}>\\
&>\sum_{p,m; m>1}\frac{\ln p}{p^{ms}}\geq \sum_{\substack{p,m; m>1\\
p^m\equiv b \bmod a}}\frac{\ln p}{p^{ms}}, \quad s>1.
\end{align*}

Por lo tanto $\sum\limits_{p\equiv b\bmod a}\frac{\ln p}
{p^s}\xrightarrow[
s\to 1^+]{}\infty$ de donde se sigue el resultado. $\fin$

\end{proof}

\begin{observacion}\label{base de Dirichlet}
Todo el desarrollo anterior lo hemos tomado de
 \cite[3.3]{FinRos2007}.
 
Si suponemos conocido que la funci\'on zeta de Dedekind
tiene un polo simple en $s=1$, la demostraci\'on del
Teorema \ref{T13.13} es casi inmediata. Primero se prueba
que $X$ es un grupo de caracteres de Dirichlet y $K$ es el
campo asociado a $X$, entonces
\[
\zeta_K(s)=\prod_{\chi\in X}L(s,\chi)
\]
lo cual se prueba comprobando que los factores de Euler son
los mismos. De ah{\'\i} se sigue que $L(1,\chi)\neq 1$ para
$\chi\neq 1$. A partir de aqu{\'\i} seguimos con los Teoremas
\ref{T13.12} y \ref {T13.13}.

Tambi\'en hacemos notar lo deseable de tener una demostraci\'on
no anal\'itica del Teorema de Dirichlet.
\end{observacion}

%% file: Capitulo8.tex
\chapter{Campos de funciones}\label{Ch5}

\section{Generalidades}\label{S5.1}

En este cap{\'\i}tulo presentamos un breve resumen,
usualmente sin demostraciones, de lo que son los campos
de funciones. Un desarrollo m\'as detallado de estos campos
puede consultarse en \cite{Che51,Deu73,Sti2009,Vil2003,Vil2006}. 
Aqu{\'\i} \'unicamente presentaremos algunos casos especiales
de campos de funciones que son los que aplicaremos a lo que nos
interesa que son los campos de funciones ciclot\'omicos (m\'odulos
de Carlitz) y m\'as generalmente, los campos de funciones
congruentes.

\begin{definicion}\label{D5.1.1} Sea $k$ un campo arbitrario.
Un {\em campo de funciones\index{campo de funciones}}
$K$ sobre $k$ es una extensi\'on finitamente generada de $k$
con grado de trascendencia $1$.
\end{definicion}

Nosotros nos restringiremos al caso en que $k$ es un campo
perfecto, de hecho nuestro principal inter\'es es cuando $k$ es
finito. En este caso, un campo de funciones $K/k$ es un campo
$K$ de la forma $K=k(x,y)$ donde $x$ es transcendente sobre
$k$ y $y$ es algebraico y separable sobre $k(x)$.

As{\'\i} mismo, supondremos que $k$ es algebraicamente cerrado
en $K$. Esto es, si $k'=\{\alpha\in K\mid \alpha \text{\ es algebraico
sobre k\ }\}$, entonces $k'=k$. En este caso $k$ se llama 
{\em el campo de constantes\index{campo de constantes}} de $K$.

\begin{definicion}\label{D5.1.1(2)}
Un campo de funciones $K$ sobre $k$ se llama {\em campo
de funciones congruentes\index{campos de funciones
congruentes}} o {\em campo global de funciones\index{campo
global de funciones}} si $k$ es un campo finito.
\end{definicion}

\begin{definicion}\label{D5.1.2}
Sea $v\colon K^{\ast}\to {\ma Z}$ una valuaci\'on discreta
de $K^{\ast}$ que es trivial en $k^{\ast}$, esto es, $v(\alpha)=0$
para $\alpha\in k^{\ast}$.

Ponemos $v(0)=\infty$ con el sobre
entendido de que $x<\infty$ para toda $x\in{\ma Z}$. Sea
${\cal O}_v=\{x\in K\mid v(x)\geq 0\}$ y $\pK_v=\{x\in K\mid
v(x)>0\}$ el anillo de valuaci\'on de $v$ y el ideal m\'aximo de
${\cal O}_v$ respectivamente. Entonces $v$ tambi\'en lo llamaremos
``{\em lugar\index{lugar}}'' de $K^{\ast}$. Al campo $k(v):=
{\cal O}_v/\pK_v$ se le llama el {\em campo
residual\index{campo residual}} de $v$.
\end{definicion}

Se tiene que $[k(v):k]<\infty$ y $[k(v):k]$ se le llame el
{\em grado\index{grado de un lugar}} de $v$.

Rec{\'\i}procamente, si ${\cal O}\subseteq K$ es un subanillo de $K$
que es de valuaci\'on con $k\subseteq {\cal O}$ y $\coc {\cal O}=K$,
${\cal O}$ da lugar a una valuaci\'on $v$ (que no describiremos
aqu{\'\i}) y si $\pK$ es el ideal maximal de ${\cal O}$, se denota
$k(\pK)=k(v)$ y $d(\pK)=d_K(\pK)=f_{\pK}=[k(\pK):k]$ al grado de $\pK$.

Cada uno de los tres objetos: el anillo de valuaci\'on ${\cal O}$,
la valuaci\'on $v$ y el ideal m\'aximo $\pK$ de ${\cal O}$, determinan
los otros dos. Todos ellos de manera indistinta los nombraremos
``{\em lugar}''.

Se tiene que si $v$ es una valuaci\'on de $K$, entonces $v|_{k(x)}$
es una valuaci\'on de $k(x)$ donde $x\in K\setminus k$, es decir,
$x$ es transcendente sobre $k$ y por tanto $[K:k(x)]$ es finita.

Rec{\'\i}procamente, dada una valuaci\'on $v$ en $k(x)$, $v$ se puede
extender a una valuaci\'on en $K$ y el n\'umero de tales extensiones
es menor o igual a $[K:k(x)]$.

\section{Valuaciones en $k(x)$}\label{S5.2}

Sea $f(x)\in k[x]$ un polinomio m\'onico e irreducible de $k[x]$.
Entonces si $\alpha(x)\in k(x)^{\ast}$, $\alpha(x)$ se puede
escribir de manera \'unica como $\alpha(x)=f(x)^s\frac{a(x)}{b(x)}$
donde $a(x),b(x)\in k[x]$ son polinomios primos relativos
a $f(x)$ y $s\in
{\ma Z}$. Entonces definimos $v_f(\alpha(x)):=s$. Se tiene que
$v_f$ es una valuaci\'on en $K$, asociada a, o correspondiente a, $f$.

Ahora sea $\beta(x)\in k(x)^{\ast}$, $\beta(x)=\frac{h(x)}{g(x)}$ con
$h(x),g(x)\in k[x]$. Definimos el {\em grado\index{grado de una
funci\'on racional}} de $\beta$ por: $\deg \beta=\deg h-\deg g$.
Sea $v_{\infty}\colon k(x)^{\ast}\to {\ma Z}$, $v_{\infty}(\beta(x))=
-\deg\beta$. Entonces $v_{\infty}$ es otra valuaci\'on diferente a
todas las valuaciones $v_f$, $f\in k[x]$ m\'onico e irreducible.

\begin{teorema}\label{T5.2.1}
El conjunto $\{v_{\infty},v_f\mid f(x)\in k[x] \text{\ es m\'onico e
irreducible}\}$ comprende a todas las valuaciones de $K$ que son
triviales en $k$. $\fin$
\end{teorema}

Usualmente denotaremos por $\pK_{\infty},\pK_f$ al lugar $v_{\infty}$,
$v_f$ respectivamente.

\begin{notacion}\label{N5.2.2} Si $K$ es un campo de funciones,
entonces ${\ma P}_K:=\{\pK\mid \pK \text{\ es lugar de\ }K\}$.
Si $\pK\in{\eu P}_K$ la valuaci\'on respectiva ser\'a denotada
por $v_{\pK}$.
\end{notacion}

A continuaci\'on definimos el substituto del dominio de Dedekind
${\cal O}_K$ en un campo num\'erico $K$.

\begin{definicion}\label{D5.2.3}
Sea $K$ un campo de funciones sobre $k$. El grupo abeliano libre
generado por ${\ma P}_K$ se llama {\em grupo de
divisores\index{grupo de divisores}} de $K$ y denota por $D_K$.
\end{definicion}

Los elementos $D_K$ se llaman {\em divisores\index{divisor}}.
Escribiremos $D_K$ multiplicativamente.

Un elemento ${\eu a}\in D_K$ es una expresi\'on formal
\[
{\eu a}=\prod_{i=1}^r \pK_i^{\alpha_i}
\]
con $\pK_1,\ldots,\pK_r\in {\ma P}_K$ y $\alpha_1,\ldots, \alpha_r
\in {\ma Z}$.

Equivalentemente, ${\eu a}=\prod_{\pK\in{\ma P}_K} \pK^{
v_{\pK}({\eu a})}$ donde $v_{\pK}({\eu a})\in{\ma Z}$ para toda
$\pK\in{\ma P}_K$ y $v_{\pK}({\eu a})=0$ para casi toda $\pK$, es
decir, $v_{\pK}({\eu a})$ es cero para todo $\pK\in {\ma P}_K$
salvo un n\'umero finito.

El elemento identidad de $D_K$ es el divisor ${\eu N}\in D_K$
tal que $v_{\pK}({\eu N})=0$ para todo $\pK\in{\ma P}_K$.

Los elementos de ${\ma P}_K$, es decir, los lugares, tambi\'en
reciben el nombre de {\em divisores
primos\index{primo!divisor}\index{divisor primo}}.

Se define el {\em grado\index{grado de un divisor}} ${\eu a}\in D_K$
por
\[
d_K({\eu a})=\sum_{\pK\in{\ma P}_K} d_K(\pK)v_{\pK}({\eu a}).
\]

Dado $x\in K^{\ast}$ se tiene que existe un n\'umero finito de
lugares $\pK$ de $K$ tales que $v_{\pK}(x)\neq 0$. Notemos que
si $x\in k^{\ast}$, $v_{\pK}(x)=0$ para toda $\pK\in{\ma P}_K$.
Se define el {\em divisor principal\index{divisor principal}} de
$x\in K^{\ast}$ por:
\[
(x)_K:=\prod_{\pK\in{\ma P}_K} \pK^{v_{\pK}(x)}.
\]

Se tiene que si $x\in K\setminus k$ existe al menos un lugar $\pK$
de $K$ tal que $v_{\pK}(x)>0$ y otro $\pK'$ tal que $v_{\pK'}(x)<0$.
Podemos tomar $\pK$ como una extensi\'on de $\pK_x$ en $k(x)$
a $K$ y $\pK'$ como una extensi\'on de $\pK_{\infty}$ a $K$.

Tambi\'en se tiene que si ${\eu Z}_x:=\prod_{v_{\pK}(x)>0}
\pK^{v_{\pK}(x)}$ y ${\eu N}_x:=\prod_{v_{\pK}(x)<0} \pK^{-v_{\pK}
(x)}$, entonces ${\eu Z}_x$ se llama el {\em divisor de
ceros\index{divisor de ceros de un elemento}} de $x$ y
${\eu N}_x$ se llama el
{\em divisor de polos\index{divisor de polos de un elemento}} de $x$.

\begin{teorema}\label{T5.2.4} Se tiene que si $x\in k^{\ast}$, 
${\eu Z}_x={\eu N}_x={\eu N}$ y si $x\in K\setminus k$,
$d({\eu Z}_x)=d({\eu N}_x)=[K:k(x)]$. $\fin$
\end{teorema}

\begin{corolario}\label{C5.2.5}
Para $x\in K^{\ast}$, $d_K((x)_K)=0$.
\end{corolario}

\begin{proof}
Notemos que $(x)_K=\frac{{\eu Z}_x}{{\eu N}_x}$ y que $d_K(
(x)_K)=d_K({\eu Z}_x)-d_K({\eu N}_x)=0$. $\fin$
\end{proof}

Sea $d_K\colon D_K\to{\ma Z}$ la funci\'on grado. Entonces
$d_K$ es un homomorfismo de grupos. Puesto que $d_K\neq 0$,
$\im d_K=m{\ma Z}$ para alg\'un $m\in{\ma N}$. En particular
$\im d_K\cong {\ma Z}$. Sea $\ker d_K=D_{K,0}:=
\{{\eu a}\in D_K\mid d_K({\eu a})=0\}$ el cual se llama el grupo de
divisores de grado $0$ de $K$. Adem\'as $P_K:=\{(x)_K\mid
x\in K^{\ast}\}\subseteq D_{K,0}$.

Se definen los siguientes grupos:
\begin{align*}
I_K:&=D_K/P_K= \text{\ grupo de clases de divisores de\ } k,\\
I_{K,0}:&=D_{K,0}/P_K= \text{\ grupo de clases de divisores
de grado $0$ de\ } K.
\end{align*}

Notemos que si $C\in I_K$ y ${\eu a}, {\eu b}\in D_K$ son tales
que ${\eu a},{\eu b}\in C$, entonces ${\eu a}={\eu b}(x)_K$ para
alg\'un $x\in K$. En particular $d_K({\eu a})=d_K({\eu b})$. Por
tanto podemos definir el {\em grado\index{grado de una clase}} de
$C$ por: $\tilde{d}_K(C):=d_K({\eu a})$ para ${\eu a}\in C$.
Se tiene que $\tilde{d}_K\colon I_K\to {\ma Z}$ tambi\'en es un
homomorfismo de grupos con $\im \tilde{d}_K=\im d_K=m{\ma Z}$
y $\ker \tilde{d}_K=I_{K,0}$. En adelante pondremos $\tilde{d}_K=
d_K$.

Se tienen las siguientes sucesiones exactas de grupos abelianos:
\begin{gather*}
1\longto D_{K,0}\longto D_K\stackrel{d_K}{\longto}\im d_K\cong
{\ma Z}\longto 0,\\
1\longto I_{K,0}\longto I_K\stackrel{d_K}{\longto}\im d_K\cong
{\ma Z}\longto 0,\\
1\longto P_K\longto D_{K,0}\longto I_{K,0}\longto 0,\\
1\longto P_K\longto D_K\longto I_K\longto 0,\\
\begin{array}{ccccccccc}
1&\longto& k^{\ast}&\longto &K^{\ast}&\longto& P_K&\longto& 1.\\
&&&&x&\longto&(x)_K
\end{array}
\end{gather*}

Puesto que ${\ma Z}$ es un grupo abeliano libre, en particular
proyectivo, se tiene
\begin{gather*}
I_K\cong I_{K,0}\oplus\im d_K\cong I_{K,0}\oplus {\ma Z},\\
D_K\cong D_{K,0}\oplus\im d_K\cong D_{K,0}\oplus {\ma Z}.
\end{gather*}

Los isomorfismos anteriores los podemos explicitar de la siguiente
forma. Sea ${\eu a}_1\in D_K$ tal que $d_K({\eu a}_1)=m$ donde
$\im d_K=m{\ma Z}$ y $m>0$. Sea ${\eu a}\in D_K$ arbitrario.
Entonces $d_K({\eu a})=tm$. Sea ${\eu a}_0:={\eu a}{\eu a}_1^{-t}$.
Entonces $d({\eu a}_0)=0$ y $\begin{array}{rrcl}
\varphi\colon&D_K&\longto &D_{K,0}\oplus {\ma Z}\\
&{\eu a}&\longmapsto &({\eu a}_0,t)\end{array}$
es el isomorfismo buscado. Similarmente para $I_K$.

En general se puede tener $m>1$. Sin embargo cuando $k$ es finito se tiene
que $m=1$.

\section{Reparticiones y diferenciales}\label{S5.3}

Sea $\pK\in{\ma P}_K$ y sea $K_{\pK}$ la completaci\'on de $K$
con respecto a la topolog{\'\i}a dada por la m\'etrica: $\|x\|_{\pK}:=
e^{-v_{\pK}(x)}$ donde entendemos $e^{-\infty}=0$, es decir,
$\|0\|_{\pK}=0$.

\begin{definicion}\label{D5.3.1}
Una {\em repartici\'on\index{repartici\'on}} o {\em
ad\`ele\index{ad\`ele}} de $K$ es una funci\'on 
\[
\varphi\colon
{\ma P}_K\longto \bigcup\limits_{\pK\in{\ma P}_K}K_{\pK}
\]
 tal que
\las
\item $\varphi(\pK)\in K_{\pK}\ \forall\ \pK\in{\ma P}_K$,
\item $v_{\pK}(\varphi(\pK))\geq 0$ para casi toda $\pK\in {\ma P}_K$.
\end{list}

Equivalentemente, una repartici\'on es una sucesi\'on $\xi=
\{\xi_{\pK}\}_{\pK\in{\ma P}_K}$ tal que $\xi_{\pK}\in {\cal O}_{K_{\pK}}
$ para casi toda $\pK\in{\ma P}_K$. Se define
\[
v_{\pK}(\xi):=v_{\pK}(\xi_{\pK}).
\]
\end{definicion}

\begin{notacion}\label{N5.3.2}
${\eu X}_K$ o $\Lambda_K$ denota el espacio de todas las
reparticiones de $K$.
\end{notacion}

${\eu X}_K$ tiene una estructura de $K$--\'algebra con las
operaciones
\begin{align*}
(\theta\xi)_{\pK}:&=\theta_{\pK}\xi_{\pK},\\
(\theta+\xi)_{\pK}:&=\theta_{\pK}+\xi_{\pK},\\
(x\xi)_{\pK}:&=x\xi_{\pK},
\end{align*}
para cualesquiera $\theta,\xi\in{\eu X}_K$, $x\in K$ y
$\pK\in{\ma P}_K$ y se tiene que $\varphi\colon K\to {\eu X}_K$
dada por $\varphi(x):=\xi_x$ donde $(\xi_x)_{\pK}=x$ para toda
$\pK\in {\ma P}_K$ es un monomorfismo de anillos. De esta
forma podemos considerar $K\subseteq {\eu X}_K$.

Sea ${\eu a}$ un divisor ${\eu a}\in D_K$ y $\xi$ una repartici\'on
de $K$. Entonces decimos que ${\eu a}$ {\em divide} a $\xi$
y ponemos ${\eu a}|\xi$ si $v_{\pK}(\xi)=v_{\pK}(\xi_{\pK})\geq
v_{\pK}({\eu a})$ para toda $\pK\in{\ma P}_K$.

Definimos ${\eu X}_K({\eu a})=\Lambda_K({\eu a})=\{\xi\in {\eu X}_K
\mid {\eu a}|\xi\}$. Entonces ${\eu X}_K({\eu a})$ es un
$k$--subespacio vectorial de ${\eu X}_K$.

\begin{definicion}\label{D5.3.3} Una {\em diferencial\index{diferencial
de Weil}} (de Weil) de $K$ es una funci\'on $k$--lineal 
($k$--funcional) $\omega\colon {\eu X}_K\to k$ tal que existe
un divisor tal que ${\eu X}_K({\eu a})+K\subseteq \ker\omega$.

En este caso decimos que ${\eu a}^{-1}$ {\em divide} a $\omega$
y ponemos ${\eu a}^{-1}|\omega$.
\end{definicion}

Por otro lado el espacio de diferenciales $\dif_K$ forma un
$K$--espacio vectorial con las operaciones de suma de funciones
y donde para $x\in K$, $\omega\in \dif_K$, $(x\omega)(\xi):=
\omega(x\xi)$, $\xi\in {\eu X}_K$.

De hecho $\dim_K\dif_K=1$, es decir, si $\omega_0$ es cualquier
diferencial no cero, entonces toda diferencial $\omega\in \dif_K$
puede escribirse de manera \'unica como $\omega=x\omega_0$,
para alg\'un $x\in K$.

Ahora bien, si $\omega$ es una diferencial no cero existe un \'unico
divisor $(\omega)_K\in \dif_K$ tal que para un divisor arbitrario
${\eu a}\in D_K$ se tiene
\begin{gather*}
{\eu a}|\omega\iff {\eu a}|(\omega)_K,\\
\intertext{esto es,}
{\eu X}_K({\eu a}^{-1})+K\subseteq \ker\omega \iff
v_{\pK}((\omega)_K)\geq v_{\pK}({\eu a}) \ \forall\ \pK\in{\ma P}_K.
\end{gather*}

El divisor $(\omega)_K$ se construye de la siguiente forma.
Dado $\omega\neq 0$, si ${\eu a}|\omega$ entonces se tiene que
$d_K({\eu a})$ est\'a acotado superiormente. Entonces $(\omega)_K$
es el divisor de m\'aximo grado que divide a $\omega$. De hecho
tenemos para $\omega$ una diferencial no cero de $K$ que
\[
\ker \omega ={\eu X}_K((\omega)_K^{-1})+K.
\]

Ahora bien, si $\dif_K({\eu a}):=\{\omega\mid \omega =0
\text{\ o\ } \omega\neq 0 \text{\ y\ } {\eu a}|\omega\}$, entonces
$\dif_K({\eu a})$ es un $k$--espacio vectorial y se tiene que
$\dif_K({\eu a})$ y $\frac{{\eu X}_K}{{\eu X}_K({\eu a}^{-1})+K}$
son isomorfos. De hecho el mapeo $k$--bilineal
\begin{eqnarray*}
\varphi\colon \dif_K({\eu a})\times \frac{{\eu X}_K}{
{\eu X}_K({\eu a}^{-1})+K}&\longto& k\\
(\omega,\overline{\xi})&\longmapsto& \omega(\xi)
\end{eqnarray*}
es no degenerado y se tiene que $\dif_K({\eu a})$ es de
dimensi\'on finita.

Para ver la dimensi\'on de estos espacios, consideremos:

\begin{definicion}\label{D5.3.4}
Sea ${\eu a}\in D_K$. Se define $L_K({\eu a})=\{x\in K\mid
x=0 \text{\ o\ }x\neq 0 \text{\ y\ } {\eu a}|(x)_K\}= \{x\in K\mid
v_{\pK}(x)\geq v_{\pK}({\eu a})\ \forall\ \pK\in{\ma P}_K\}$.
\end{definicion}

Se tiene que $L_K({\eu a})$ es  $k$--espacio vectorial.

\begin{teorema}\label{T5.3.5} Para todo ${\eu a}\in D_K$, $L_K(
{\eu a})$ es un $k$--espacio vectorial de dimensi\'on finita y
denotamos por $\ell_K({\eu a})$ la dimensi\'on de $L_K({\eu a})$:
$\ell_K({\eu a})=\dim_k L_K({\eu a})$. $\fin$
\end{teorema}

Notemos que si $d_K({\eu a})>0$, entonces $L_K({\eu a})=\{0\}$
y $\ell_K({\eu a})=0$.

Se tiene la sucesi\'on exacta de $k$--espacios vectoriales donde
${\eu a}|{\eu b}$
\[
0\longto \frac{L_K({\eu a})}{L_K({\eu b})}\longto \frac{{\eu X}_K({\eu
a})}{{\eu X}_K({\eu b})}\longto \frac{{\eu X}_K({\eu a})+K}
{{\eu X}_K({\eu b})+K}\longto 0.
\]

Un resultado central en la teor{\'\i}a de las funciones algebraicas
es el Teorema de Riemann--Roch, el cual enunciamos a
continuaci\'on.

\begin{teorema}[Riemann--Roch\index{Riemann--Roch!teorema
de $\sim$}]\label{T5.3.6}
Sea $K/k$ cualquier campo de funciones. Existe un entero no
negativo $g_K\geq 0$ que depende \'unicamente de $K$, llamado
el {\em g\'enero\index{g\'enero}} de $K$ tal que
\l
\item (Riemann) Para cualquier divisor ${\eu a}\in D_K$ se tiene
\[
\ell_K({\eu a})+d_K({\eu a})\geq 1-g_K.
\]
Es decir $\delta_K({\eu a}^{-1}):= \ell_K({\eu a})+d_K({\eu a})+
g_K-1\geq 0$.

\item Se tiene 
\begin{align*}
\delta_K({\eu a})&=\dim_k \dif_K({\eu a})=
\dim_k \frac{{\eu X}_K}{{\eu X}({\eu a}^{-1})+K}\\
&=\ell_K({\eu a}^{-1})
+d_K({\eu a}^{-1})+g_K-1.
\end{align*}

\item Se tiene $\delta_K({\eu a})=\ell_K({\eu a}(\omega)_K^{-1})$
donde $\omega$ es cualquier diferencial no cero.

\item (Riemann--Roch) Se tiene que para cualquier divisor ${\eu a}$
y cualquier diferencial no cero $\omega$,
\[
\ell_K({\eu a}^{-1})=d_K({\eu a})-g_K+1+
\ell_K({\eu a}(\omega)_K^{-1}). \tag*{$\fin$}
\]
\end{list}
\end{teorema}

Como consecuencia del Teorema de Riemann--Roch tenemos

\begin{corolario}\label{C5.3.7}
{\ }

\l
\item $\delta_K({\eu N})=g_K$.
\item Si $\omega\in\dif_K$, $\omega\neq 0$, $d_K((\omega)_K)
=2g_K-2$.
\item Si $d_K({\eu a})>2g_K-2$, $\ell_K({\eu a}^{-1})=d_K({\eu a})
-g_K+1$.
\item Si $\pK\in{\ma P}_K$ y $n\in{\ma N}$, $n>2g_K-1$, existe
$x\in K$ tal que ${\eu N}_x=\pK^n$, esto es, $(x)_K=\frac{{\eu a}}
{\pK^n}$ con ${\eu a}$ un divisor entero y primo relativo a $\pK$,
esto es, $v_{\pK}({\eu a})=0$. $\fin$
\end{list}
\end{corolario}

Finalmente, presentamos una aplicaci\'on importante que es
consecuencia del Teorema de Riemann--Roch.

\begin{teorema}[Teorema de aproximaci\'on fuerte\index{teorema!de
aproximaci\'on fuerte}\index{aproximaci\'on
fuerte!teorema de $\sim$}]\label{T5.3.8}
Sea $S$ un subconjunto propio del conjunto de lugares de un campo
de funciones $K/k$, $S\subsetneqq {\ma P}_K$.
Sean $\{\pK_1,\ldots,\pK_r\}\subseteq S$ un conjunto finito
de lugares de $K$, $a_1\in K_{\pK_1},\ldots,a_r\in K_{\pK_r}$ elementos
arbitrarios de las respectivas completaciones de $K$ en $\pK_i$,
$1\leq i\leq r$ y $n_1,\ldots,n_r$ enteros arbitrarios. Entonces, existe
$x\in\*K$ tal que 
\begin{gather*}
v_{\pK_i}(x-a_i)=n_i\quad \text{para toda}\quad 1\leq i\leq r,\\
v_{\pK}(x)\geq 0\quad \text{para todo}\quad \pK\in S\setminus\{
\pK_1,\ldots,\pK_r\}.
\end{gather*}
\end{teorema}

\begin{proof} Puesto que $K$ es denso en $K_{\pK_i}$, $1\leq i\leq r$,
existe, para cada $i$, $b_i\in K$ tal que $v_{\pK_i}(a_i-b_i)>n_i$. Entonces
$v_{\pK_i}(x-b_i)=v_{\pK_i}(x-a_i+a_i-b_i)=\min\{v_{\pK_i}(x-a_i),
v_{\pK_i}(a_i-b_i)\}=v_{\pK_i}(x-a_i)=n_i$, $1\leq i\leq r$. Por tanto
podemos suponer,
sin p\'erdida de generalidad, que $a_i\in K$, $1\leq i\leq r$.

Sea $\pK_0$ un lugar tal que $\pK_0\notin S$ y sea ${\eu a}:=\pK_0^{m}\prod_{
i=1}^r\pK_i^{-(n_i+1)}$ con $m$ cualquier entero tal que $m>2g_K-2
+\sum_{i=1}^r(n_i+1)\deg {\pK_i}$.
Entonces, por el Teorema de Riemann--Roch,
se tiene que
\[
\delta_K({\eu a})=\dim_k\frac{{\eu X}}{{\eu X}({\eu a}^{-1})+K}=\ell_K({\eu a}^{-1})+
d_K({\eu a}^{-1})+g_K-1=\dim_k\dif_K({\eu a}).
\]

Puesto que $d_K({\eu a})\geq m-\sum_{i=1}^r(n_i+1)\deg \pK_i
>2g_K-2=d_K(\omega)$
para toda diferencial diferente de cero, se sige que
${\eu a}\nmid \omega$
para cualquier diferencial no cero $\omega$ de $K$,
lo cual implica que $\dif_K({\eu a})=\{0\}$. Esto es,
\begin{gather}\label{Ec9.5.8(1)}
{\eu X}={\eu X}({\eu a}^{-1})+K.
\end{gather}

Sea $\xi=\{\xi_{\pK}\}_{\pK\in{\ma P}_K}$ la repartici\'on definida
por
\[
\xi_{\pK_i}=-a_i, \quad 1\leq i\leq r \quad\text{y}\quad 
\xi_{\pK}=0 \quad\text{para}
\quad \pK\notin\{\pK_1,\ldots,\pK_r\}.
\]
De (\ref{Ec9.5.8(1)}) obtenemos que existe $y\in K$ tal que $y+\xi\in{\eu X}_K(
{\eu a}^{-1})$. Esto nos dice que
\begin{gather}
v_{\pK_i}(y+\xi_{\pK_i})=v_{\pK_i}(y-a_i)\geq v_{\pK_i}({\eu a}^{-1})=n_i+1,
\quad 1\leq i\leq r\quad\text{y}\nonumber\\
 v_{\pK}(y)\geq v_{\pK}({\eu a}^{-1})=0 \quad\text{para todo}
\quad \pK\notin\{\pK_0,\pK_1,\ldots,\pK_r\}.\label{Ec9.5.8(2)}
\end{gather}

Sean $c_i\in K$ con $v_{\pK_i}(c_i)=n_i$, $1\leq i\leq r$. De 
(\ref{Ec9.5.8(2)}) obtenemos que existe $z\in K$ tal que
\begin{gather*}
v_{\pK_i}(z-c_i)\geq n_i+1>n_i,
\quad 1\leq i\leq r\quad\text{y}\\
 v_{\pK}(z)\geq 0 \quad\text{para todo}
\quad \pK\notin\{\pK_0,\pK_1,\ldots,\pK_r\}.
\end{gather*}

En particular, para $1\leq i\leq r$, 
se tiene que $v_{\pK_i}(z)=v_{\pK_i}(z-c_i+c_i)=
\min\{v_{\pK_i}(z-c_i),v_{\pK_i}(c_i)\}=n_i$. Finalmente, sea $x:=
y+z$. Entonces, 
\begin{gather*}
v_{\pK_i}(x-a_i)=v_{\pK_i}(y+z-a_i)=\min\{
v_{\pK_i}(y-a_i),v_{\pK_i}(z)\}=n_i\quad \text{para}\\
1\leq i\leq r \quad \text{y}\quad v_{\pK}(x)=
v_{\pK}(y+z)\geq 0\quad\text{para todo}\quad
\pK\notin\{\pK_0,\pK_1,\ldots,\pK_r\}.\tag*{$\fin$}
\end{gather*}
\end{proof}

\begin{observacion}\label{O5.3.9}
Puesto que todo elemento $x\in\*K$ es de grado $0$, la condici\'on
de que $S$ sea un subconjunto propio del conjunto de todos los
lugares, no puede suprimirse.
\end{observacion}

\section{Extensiones de Galois}\label{S5.4}

\begin{definicion}\label{D5.4.1} Sean $K/k$ y $L/\ell$ dos
campos de funciones. Decimos que $L$ es una {\em
extensi\'on\index{extensi\'on de campos de funciones}} de $K$ si
$K\subseteq L$ y $\ell\cap K=k$.
\end{definicion}

Si $\pK\in{\ma P}_K$ y $\pL\in{\ma P}_L$ tal que $\pL$ es una
extensi\'on de $\pK$. Se define el {\em grado relativo\index{grado
relativo}} por
\[
d_{L/K}(\pL|\pK)=[\ell(\pL):k(\pK)].
\]

Notemos que puesto que se tiene el diagrama
\[
\xymatrix{
k(\pK)\ar@{-}[r]\ar@{-}[d]&\ell(\pL)\ar@{-}[d]\\ k\ar@{-}[r]& \ell
}
\]
y $d_K(\pK)=[k(\pK):k]$, $d_L(\pL)=[\ell(\pL):\ell]$, se sigue
\[
d_{L/K}(\pL|\pK)d_K(\pK)=d_L(\pL)[\ell:k]
\] 
(finito o infinito).

Puesto que $d_K(\pK)$ y $d_L(\pL)$ son finitos, se tiene que
$d_{L/K}(\pL|\pK)<\infty \iff [\ell:k]<\infty$.

\begin{proposicion}\label{P5.4.2}
Sea $L/\ell$ una extensi\'on de $K/k$. Sea $\pK\in{\ma P}_K$.
Entonces el n\'umero de extensiones de la valuaci\'on $v_{\pK}$
a $L$ es finito.
\end{proposicion}

\begin{proof}
Por el Teorema de Riemann--Roch, existe $x\in K$ tal que 
${\eu N}_{x,K}=\pK^n$ para alg\'un $n\in{\ma N}$. Ahora
$\pL\in{\ma P}_L$ extiende a $\pK$ si y s\'olo si $v_{\pL}(x)<0$
lo cual es equivalente a que $\pL|{\eu N}_{x,L}$. Este \'ultimo
n\'umero es finito. $\fin$
\end{proof}

\begin{definicion}\label{D5.4.3}
Sea $\pL\in{\ma P}_L$ y sea $\pK$ la restricci\'on $\pL$ a $K$.
Esto es, $v_{\pL}|_K$ es equivalente a $v_{\pK}$. Ahora
bien $v_{\pL}\colon L^{\ast}\to{\ma Z}$ es suprayectiva pero
$v_{\pL}|_{K^{\ast}}\colon K^{\ast}\to {\ma Z}$ no necesariamente
lo es. Se define {\em {\'\i}ndice de ramificaci\'on\index{indice@{\'\i}ndice
de ramificaci\'on}} de $\pL$ sobre $\pK$ como el n\'umero natural
$e=e_{L/K}(\pL|\pK)$ tal que $v_{\pL}(\alpha)=ev_{\pK}(\alpha)$
para $\alpha\in K$.
\end{definicion}

\begin{definicion}\label{D5.4.4}
Si $L/\ell$ es una extensi\'on de $K/k$ y si $\pK\in{\ma P}_K$,
entonces si $\pL_1,\ldots,\pL_h$ son todos las extensiones de $K$
a $L$ se define la {\em conorma\index{conorma}} de $\pK$ a $L$
por
\[
\con_{K/L}\pK=\pL_1^{e_{L/K}(\pL_1|\pK)}\cdots
\pL_h^{e_{L/K}(\pL_h|\pK)}.
\]
Si ${\eu a}\in D_K$ es un divisor, ${\eu a}=\pK_1^{\alpha_1}\cdots
\pK_r^{\alpha_r}$ se define la conorma de ${\eu a}$ por
\[
\con_{K/L}({\eu a})=\prod_{i=1}^r \con_{K/L} (\pK_i)^{\alpha_i}.
\]
\end{definicion}

\begin{teorema}\label{T5.4.5}
Para cualquier extensi\'on $L/\ell$ de $K/k$, finita o infinita, se
tiene
\[
[L:K]=\sum_{i=1}^h e_{L/K}(\pL_i|\pK) d_{L/K}(\pL_i|\pK).
\tag*{$\fin$}
\]
\end{teorema}

Cuando $L/K$ es una extensi\'on de Galois se tiene que 
$d_{L/K}(\pL_i|\pK)=d_{L/K}(\pL_j|\pK)$ y $e_{L/K}(\pL_i|\pK)
=e_{L/K}(\pL_j|\pK)$ para todo $1\leq i,j\leq h$. Sean
\begin{gather*}
f:=d_{L/K}(\pL_i|\pK),\quad e=e_{L/K}(\pL_i|\pK), \quad
1\leq i\leq h.\\
\intertext{Entonces el Teorema \ref{T5.4.5} toma la forma}
[L:K]=efh.
\end{gather*}

\begin{definicion}\label{D5.4.6}
Sea $L/K$ una extensi\'on de Galois con grupo de Galois $G=
\Gal(L/K)$. Sea $\pK\in{\ma P}_K$ y $\pL\in{\ma P}_L$ una
extensi\'on de $\pK$. Se define
\l
\item $D=D_{L/K}(\pL|\pK):=\{\sigma\in G\mid \sigma\pL=\pL\}=$ 
grupo de descomposici\'on de $\pL/\pK$.

\item $I=I_{L/K}(\pL|\pK):=\{\sigma\in G\mid \sigma x\equiv
x\bmod \pL\ \forall\ x\in {\cal O}_\pL\}=$ 
grupo de inercia de $\pL/\pK$.
\end{list}
\end{definicion}

Como en el caso de campos num\'ericos, se tiene que 
\begin{gather*}
I\subseteq D, \quad |I|=e=e_{L/K}(\pL|\pK), 
\quad |D|=ef \quad \text{donde}\quad f=d_{L/K}(\pL|\pK),\\
D/I\cong \Gal(\ell(\pL/k(\pK)),\\
D\cong \Gal(L_{\pL}|K_{\pK}),
\end{gather*}
donde $L_{\pL}$ y $K_{\pK}$ son las completaciones de $L$ y $K$
en $\pL$ y $\pK$ respectivamente.

\section{Diferente, discriminante y ramificaci\'on}\label{S5.5}

Sea $L/K$ una extensi\'on separable de campos de funciones,
$\pL$ un lugar de $L$ y $\pK:=\pL|_K$. Sean $L_{\pL}$ y $K_{\pK}$
las respectivas completaciones. Se tiene que $[L_{\pL}:K_{\pK}]=
e_{L/K}(\pL|\pK) d_{L/K}(\pL|\pK)$. Sea $\tilde{L}$ la cerradura
de Galois de $L/K$ y sea $\pL$ un lugar en $\tilde{L}$ sobre
$\pL$.

Sea ${\cal O}_{\hat{\pL}}=\{x\in L_{\pL}\mid v_{\pL}(x)\geq 0\}$ la
completaci\'on de ${\cal O}_{\pL}$ y $\hat{\pL}=\{x\in L_{\pL}\mid
v_{\pL}(x)> 0\}$ la completaci\'on de $\pK$. Si $\Tr=\Tr_{L_{\pL}/
K_{\pK}}$ denota la traza de $L_{\pL}$ a $K_{\pK}$ se tiene:

\begin{teorema}\label{T5.5.1}
Existe $m\in{\ma Z}$, $m\geq 0$ tal que si $x\in L_{\pL}$ satisface
$v_{\pL}(x)\geq -m$, entonces $v_{\pK}(\Tr x)\geq 0$ y existe
$x_0\in L_{\pL}$ tal que $v_{\pL}(x_0)<-m$ y $v_{\pK}(\Tr x_0)
<0$. $\fin$
\end{teorema}

\begin{definicion}\label{D5.5.2}
El m\'aximo entero no negativo que satisface las condiciones
del Teorema \ref{T5.5.1} se denota por $m(\pL)$ y es llamado
el {\em exponente diferencial\index{exponente diferencial}} de $\pL$
con respecto a $K$.
\end{definicion}

\begin{teorema}\label{T5.5.3} Se tiene que $m(\pL)\geq e-1=
e_{L/K}(\pL|\pK)-1$. Adem\'as, puesto que $k$ es perfecto, $m(\pL)>
e-1$ si y s\'olo si la caracter{\'\i}stica de $k$ divide a $e$.

En particular $m(\pL)=0$ si $\pL$ no es ramificado. $\fin$
\end{teorema}

Como en el caso num\'erico, se define que $\pL$ es
{\em moderadamente ramificado\index{ramificaci\'on moderada}} si
$p\nmid e$ y {\em salvajemente ramificado\index{ramificaci\'on
salvaje}} si $p|e$, donde $p$ es la caracter{\'\i}stica de $k$.

\begin{teorema}\label{T5.5.4} Se tiene que $m(\pL)=0$ salvo un
n\'umero finito de lugares $\pL$. $\fin$
\end{teorema}

\begin{definicion}\label{D5.5.5}
El divisor ${\eu D}_{L/K}:=\prod\limits_{\pL\in{\ma P}_L}\pL^{m(\pL)}$
se llama el {\em diferente de la extensi\'on $L/K$\index{diferente
en campos de funciones}} y se tiene que $\pL|{\eu D}_{L/K} \iff
\pL$ es ramificado.

El {\em discriminante\index{discriminante en campos de funciones}}
${\eu d}_{L/K}$ de la extensi\'on $L/K$ se define como la norma
del diferente: ${\eu d}_{L/K}:=N_{L/K}({\eu D}_{L/K})$.
\end{definicion}

Sea $L/K$ una extensi\'on separable de campos de funciones.
Sea $\pK\in{\ma P}_K$ dado. Por el Teorema de Riemann--Roch
existe $x\in K$ tal que ${\eu N}_x=\pK^n$ para alg\'un $n\geq 1$.
Entonces $k[x]$ es un dominio Dedekind y sean ${\cal O}_K$ y
${\cal O}_L$ las cerraduras enteras de $k[x]$ en $K$ y $L$
respectivamente. Entonces ${\cal O}_K$ y ${\cal O}_L$ son
dominios Dedekind. Se puede definir el diferente de los dominios
${\cal O}_L$ y ${\cal O}_K$ de la manera usual; esto es,
\[
{\eu D}_{{\cal O}_L/{\cal O}_K}^{-1}:=\{x\in L\mid \Tr_{L/K}(xy)
\in{\cal O}_K\ \forall\ y\in {\cal O}_L\}.
\]

Entonces se tiene que si identificamos los ideales primos de 
${\cal O}_L$ con los lugares de $L$ que no est\'an sobre $\pK$,
se tiene
\[
{\eu D}_{{\cal O}_L/{\cal O}_K}=\prod_{\substack{\pL\in{\ma P}_L\\
\pL\nmid\pK}} \pL^{m(\pL)},
\]
es decir, los exponentes $m(\pL)$ son los mismos que los de la
Definici\'on \ref{D5.5.2} y en ${\eu D}_{{\cal O}_L/{\cal O}_K}$ solo
faltan los primos que dividen a $\pK$.

Es m\'as f\'acil estudiar ${\eu D}_{{\cal O}_L/{\cal O}_K}$ que
${\eu D}_{L/K}$ y para completar la informaci\'on sobre ${\eu D}_{
L/K}$ podemos tomar otro lugar $\pK'\neq \pK$ y repetir el proceso
para obtener los exponentes $m(\pL)$ de los lugares $\pL$ de $L$
sobre $\pK$.

Para dominios Dedekind, tenemos la siguiente forma para calcular
diferentes.

\begin{teorema}\label{T5.5.6}
Sea $A$ un dominio Dedekind, $K:=\coc A$ el campo de cocientes
y $L/K$ una extensi\'on finita y separable de $K$. Sea $B$ la
cerradura entera de $A$ en $L$. Entonces $B$ es un dominio
Dedekind y
\l
\item Si existe $\alpha\in B$ tal que $B=A[\alpha]$ entonces
${\eu D}_{B/A}=\langle f'(\alpha)\rangle$ donde $f(x)=\Irr(\alpha,x,
K)$.
\item En general se tiene
\[
{\eu D}_{B/A}=\langle f'(\alpha)\mid \alpha\in B, L=K(\alpha) 
\text{\ y\ } f(x)=\Irr(\alpha,x,K)\rangle. \tag*{$\fin$}
\]
\end{list}
\end{teorema}

\section{Formula de Riemann--Hurwitz}\label{S5.6}

Sea $L/\ell$ una extensi\'on finita de $K/k$.

\begin{definicion}\label{D5.6.1}
Sea $\xi\in{\eu X}_K$. Se define la {\em cotraza\index{cotraza 
de una repartici\'on}} de $\xi$, denotada por
$\cotr_{K/L}\xi$, como la repartici\'on $\theta\in {\eu X}_L$ dada 
como: si $\pL\in{\ma P}_L$, $\pL|_K=\pK$ y puesto que $K_{\pK}
\subseteq L_{\pL}$ se define $\theta_{\pL}:=\xi_{\pK}$.

Si $\Omega\in \dif_L$ es una diferencial de $L$, se define
la {\em traza\index{traza de una diferencial}} de $\Omega$, y se
denota por $\tr_{L/K}\Omega$ como la diferencial $\omega$ 
definida por:
\[
\omega\colon{\eu X}_K\to k, \quad \omega(\xi):=\Omega(\cotr_{K/L}
\xi).
\]
\end{definicion}

Se tiene que si $\xi\in{\eu X}_K$, $\Omega\in \dif_L$, entonces
$\cotr_{K/L}\xi\in {\eu X}_L$ y que $\tr_{L/K}\Omega\in\dif_K$.

\begin{definicion}\label{D5.6.2}
Si $\theta\in{\eu X}_L$ definimos la {\em traza\index{traza de una
repartici\'on}} de $\theta$, denotada por $\tr_{L/K}\theta$ como $\xi$
donde para $\pK\in{\ma P}_K$:
\[
\xi_{\pK}=\sum_{i=1}^h \tr_{L_{\pL_i}/K_{\pK}}\theta_{\pL_i}
\]
donde $\pL_1,\ldots,\pL_h$ son los lugares sobre $\pK$. Se tiene
que $\xi=\tr_{L/K}\theta \in{\eu X}_K$.
\end{definicion}

En el caso anterior, a la operaci\'on cotraza de reparticiones le
asociamos la operaci\'on traza de diferenciales.
Rec{\'\i}procamente a la operaci\'on traza de reparticiones queremos
asociarle una operaci\'on cotraza de diferenciales. En este punto
tenemos el problema que \'unicamente obtenemos $k$--linealidad
y no $\ell$--linealidad. Por lo pronto nos restringimos al caso
``{\em geom\'etrico\index{extensi\'on geom\'etrica}}'', esto es, 
cuando $\ell=k$.

\begin{definicion}\label{D5.6.3} Sea $L/K$ una extensi\'on finita
y geom\'etrica de campos de funciones, es decir, $k=\ell$. Sea
$\omega\in\dif_K$. Definimos la {\em cotraza\index{cotraza de
un diferencial}} de $\omega$, denotada por $\cotr_{K/L}\omega$,
por $\Omega$ donde para $\xi\in{\eu X}_L$, $\Omega(\xi)=
\omega(\tr_{L/K}\xi)$.
\end{definicion}

Se tiene:

\begin{teorema}\label{T5.6.4}
Si $L/K$ es geom\'etrica y finita, entonces para $\omega\in\dif_K$,
$\cotr_{K/L}\omega\in\dif_L$. M\'as a\'un si $L/K$ es separable
entonces $\cotr_{K/L}\omega\neq 0$ para $\omega\neq 0$ y se
tiene que
\[
(\cotr_{K/L}\omega)_L={\eu D}_{L/K}\con_{K/L}(\omega)_K.
\tag*{$\fin$}
\]
\end{teorema}

El Teorema \ref{T5.6.4} junto con el Corolario \ref{C5.3.7} ({\sc{ii}})
obtenemos la f\'ormula de Riemann--Hurwitz:

\begin{teorema}[Riemann--Hurwitz\index{Riemann-Hurwitz!f\'ormula
de $\sim$}\index{f\'ormula de Riemann--Hurwitz}]\label{T5.6.5}
Sea $L/K$ una extensi\'on finita, separable y geom\'etrica. Entonces
\begin{equation}\label{Ec5.6.5.1}
2g_L-2=[L:K](2g_K-2)+d_L({\eu D}_{L/K}).
\end{equation}
\end{teorema}

\begin{proof}
Se tiene que ${\eu a}$ es un divisor de $K$ entonces $d_L(
\con_{K/L}({\eu a}))=[L:K]d_K({\eu a})$. Por tanto de $(\Omega)_L=
\con_{K/L}(\omega)_K {\eu D}_{L/K}$ obtenemos que
\[
d_L((\Omega)_L)=d_L(\con_{K/L}(\omega)_K)+d_L({\eu D}_{L/K})
\]
la cual es precisamente (\ref{Ec5.6.5.1}). $\fin$
\end{proof}

Como en el caso num\'erico, tenemos los grupos de ramificaci\'on
y su relaci\'on con el diferente.

\begin{definicion}\label{D5.6.6} Sea $\pK\in{\ma P}_K$ y $\pL\in
{\ma P}_L$ sobre $\pK$. Sean $G_{-2}:=G$, $G_{-1}:=D_{L/K}(\pL|
\pK)$, $G_0:=I_{L/K}(\pL|\pK)$ y para $i\geq -i$, $i\in{\ma Z}$ se
define el {\em $i$--\'esimo grupo de ramificaci\'on\index{grupo
de ramificaci\'on}} $G_i$ por:
\[
G_i:=\{\sigma\in G_{-1}\mid v_{\pL}(\sigma a-a)\geq i+1
\text{\ para toda \ } a\in{\cal O}_L\}.
\]
\end{definicion}

Se tiene que $G_i$ es un subgrupo normal de $G_{-1}=D(\pL|\pK)$,
que $G_{i+1}\subseteq G_i$ para toda $i\in{\ma Z}$, $i\geq -1$ y
que existe $i_0$ tal que $G_{i_0}=\{\id\}$.

Si $\sigma\in G_{-1}$, $\sigma\neq \Id$, existe $i$ tal que $\sigma\in
G_i\setminus G_{i+1}$. Se pone 
$i_{G_{-1}}(\sigma)=i$. Si $\sigma=\Id$, 
definimos $i_{G_{-1}}(\sigma)=i_{G_{-1}}(\Id)=\infty$.
Notemos que 
\[
i_{G_{-1}}(\sigma)\geq j+1\iff \sigma\in G_j\quad \text{y que}
\quad \sum_{\substack{\sigma\neq \Id\\ \sigma\in G_{-1}}} i_{G_{-1}}
(\sigma)=\sum_{i=0}^{\infty}\big(|G_i|-1\big).
\]

Se tiene:

\begin{teorema}\label{T5.6.7} Si $m(\pL)$ es el exponente diferencial
de $\pL$, entonces 
\[
m(\pL)=\sum_{\substack{\sigma\in G_{-1}\\ \sigma\neq \Id}}
i_{G_{-1}}(\sigma)=\sum_{i=0}^{\infty}\big(|G_i|-1\big).
\tag*{$\fin$}
\]
\end{teorema}

%% file: Capitulo9.tex
\chapter{Campos globales de funciones y
campos de funciones ciclot\'omicos}\label{Ch6}

\section{Campos de funciones congruentes}\label{S6.1}

Se tiene que hay muchas similitudes entre los campos num\'ericos
y los campos de funciones. Cuando en estos \'ultimos el campo
de constantes es finito, los campos residuales tambi\'en son finitos
y nos permiten avanzar en esta analog{\'\i}a. Sin embargo, de entrada,
hay diferencias fundamentales. Si $K$ es un campo num\'erico
y $\pK$ es un ideal primo en ${\cal O}_K$, se tiene que ${\cal O}_K/{
\pK}$ es finito y en particular de caracter{\'\i}stica finita siendo que
el campo $K$ es de caracter{\'\i}stica $0$, esto es $K$ y ${\cal O}_K/
{\pK}$ tienen caracter{\'\i}sticas diferentes.

Por otro lado, si $K/k$ es un campo de funciones y $\pK\in{\ma P}_K$,
entonces $k(\pK)$ es una extensi\'on finita de $k$ y por tanto
$k$, $k(\pK)$ y $K$ tienen la misma caracter{\'\i}stica.

Cuando estudiamos ${\ma Q}$, tenemos que ${\ma Q}$ es el 
campo de cocientes de ${\ma Z}$. El an\'alogo a ${\ma Q}$ 
ser{\'\i}a un campo de funciones racionales $k(x)$ y el an\'alogo
a ${\ma Z}$ ser{\'\i}a $k[x]$ pues $k(x)$ es el campo de cocientes
de $k[x]$ siendo adem\'as que tanto ${\ma Z}$ como $k[x]$ son
anillos euclidianos. Sin embargo, a pesar de ser, en cierto sentido,
bastante parecidos, hay diferencias esenciales.

Por ejemplo, para $a,b,c,d\in k$ con $ad-bc\neq 0$, se tiene que
$k\big(\frac{ax+b}{cx+d}\big)=k(x)$, y por tanto el anillo de 
polinomios de $k(y)=k\big(\frac{ax+b}{cx+d}\big)$ donde
$y=\frac{ax+b}{cx+d}$, es $k[y]$ el cual, a pesar de ser isomorfo
a $k[x]$, no es igual siendo que los campos de cocientes si son
iguales. Esto no sucede en ${\ma Z}$, es decir, si $R$ es un
subanillo de ${\ma Q}$ isomorfo a ${\ma Z}$ como anillo, entonces
$R={\ma Z}$. Otra diferencia es de que si $F$ es un campo que
contiene a ${\ma Q}$ y es isomorfo a ${\ma Q}$, entonces $F=
{\ma Q}$. Esto tampoco sucede con los campos de funciones
racionales. De hecho si $n\in{\ma N}$, entonces $k(x^n)\cong k(x)
\cong k(x^{1/n})$ como campos pero $k(x^n)\subseteq k(x)
\subseteq k(x^{1/n})$ y $[k(x^{1/n}):k(x)]=n$ y $[k(x):k(x^n)]=n$.

\begin{definicion}\label{D6.1.1}
Un campo de funciones $K/k$ se llama {\em congruente\index{campo
de funciones congruente}} o {\em global\index{campo
global de funciones}} si $k$ es finito, $|k|=q$, $k\cong {\ma
F}_q$.
\end{definicion}

\begin{teorema}\label{T6.1.2}
Si $K/k$ es un campo de funciones congruente, $\ell$ es una 
extensi\'on finita de $k$ y $L:=K\ell$, entonces el campo de 
constantes de $L$ es $\ell$.

Adem\'as $[L:K]=[\ell:k]$.
\end{teorema}

\begin{proof}
Escribamos $\ell=k(\xi)$ con $[\ell:k]=f$. Entonces $\ell\cong
{\ma F}_{q^f}$ e $\Irr(\xi, x,k)|x^{q^f}-x=\prod_{\alpha\in\ell}
(x-\alpha)$. Ahora bien, $L=K\ell =Kk(\xi)=K(\xi)$ con
$\Irr(\xi, x,K)|\Irr(\xi,x,k)$. En particular $\Irr(\xi,x,K)\in
K[x]\cap \ell[x]=k[x]$. Por lo tanto
\begin{gather*}
\Irr(\xi,x,K)=\Irr(\xi,x,k)\quad\text{y}\\
 [L:K]=\gr \Irr(\xi,x,K)=
\gr\Irr(\xi,x,k)=[\ell:k].
\end{gather*}

Sea $\ell'$ el campo de constantes de $L$, $\ell\subseteq \ell'$. Por
tanto $L=K\ell'$. Imitando el paso que acabamos de realizar,
tendr{\'\i}amos que $[\ell':k]=[L:K]=[\ell:k]=f$ lo cual implica
que $\ell'=\ell$. $\fin$
\end{proof}

\begin{teorema}\label{T6.1.3} 
Sean $k$ un campo finito y $K/k$ un campo de funciones
congruente. Sean $\ell/k$ una extensi\'on finita y $L=K\ell$.
Entonces $L/K$ es una extensi\'on no ramificada.
\end{teorema}

\begin{proof}
Sean $\P$ un llugar de $K$ y $\pK$ un lugar de $L$ sobre
$\P$. Sea $I=I(\pK|\P)$ el grupo de inercia. Se tiene que
$\sigma \gamma\equiv \gamma \bmod \pK$ para toda $\gamma
\in\ell$, esto es, $v_{\pK}(\sigma\gamma-\gamma)>0$.

Por otro lado, puesto que $\sigma\gamma-\gamma\in\ell$, 
necesariamente $v_{\pK}(\sigma\gamma-\gamma)=\infty$,
es decir, $\sigma\gamma=\gamma$ para toda $\gamma\in
\ell$. Se sigue que $\sigma=\Id$ y $I=\{\Id\}$. Por tanto
$L/K$ es no ramificada.
$\fin$
\end{proof}

\begin{teorema}\label{T6.1.3(1)} Sea $L/K$ como antes. Sea
${\eu A}\in D_K$ un divisor de $K$. Entonces $d_L({\eu A})=
d_K({\eu A})$, esto es,
\[
d_L(\con_{K/L}{\eu A})=d_K({\eu A}).
\]
\end{teorema}

\begin{proof}
Sea $\P$ un lugar de $K$. Sea $\con_{K/L}\P=\pK_1\cdots
\pK_h$ y
\begin{align*}
d_L(\con_{K/L}\P)&=\sum_{i=1}^hd_L(\pK_i)=\sum_{i=1}^h
\frac{d_{L/K}(\pK_i|\P)d_K(\P)}{[\ell:k]}\\
&=\frac{d_K(\P)}{[\ell:k]}\sum_{i=1}^hd_{L/K}(\pK_i|\P)=
\frac{d_K(\P)}{[\ell:k]}[L:K].
\end{align*}

Por el Teorema \ref{T6.1.2} tenemos $[L:K]=[\ell:k]$ de donde
$d_L(\P)=d_K(\P)$ para todo lugar $\P\in{\ma P}_K$. Por tanto
$d_K({\eu A})=d_L({\eu A})$ para todo divisor ${\eu A}\in D_K$.
$\fin$
\end{proof}

\begin{teorema}\label{T6.1.3(2)} Sean $L/K$ y ${\eu A}$ como
en el Teorema {\rm{\ref{T6.1.3(1)}}}. Entonces $L_L({\eu A})=\ell
L_K({\eu A})$ y $\ell_L({\eu A})=\ell_K({\eu A})$. Adem\'as
$g_L=g_K$.
\end{teorema}

\begin{proof}
Se tiene
\begin{align*}
L_K({\eu A})&=\{x\in K\mid v_{\P}(x)\geq v_{\P}({\eu A})
\text{\ para toda\ } \P\in {\ma P}_K\} \quad \text{y}\\
L_L({\eu A})&=\{y\in K\mid v_{\pK}(y)\geq v_{\pK}({\eu A})
\text{\ para toda\ } \pK\in {\ma P}_L\} .
\end{align*}

Si $y\in L$ con $y=\sum_{i=1}^n a_ix_i$ con $a_i\in \ell$,
$x_i\in L_K({\eu A})$, entonces
\begin{align*}
v_{\pK}(y)&\geq \min_{1\leq i\leq n} v_{\pK}(a_ix_i)\geq
\min_{1\leq i\leq n}v_{\pK}(x_i) =\min_{1\leq i\leq n}v_{\P}
(x_i)\\
&\geq v_{\P}({\eu A})=v_{\pK}(\con_{K/L}{\eu A}) 
\end{align*}
donde $\pK|_K=\P$. Por tanto $\ell L_K({\eu A})\subseteq 
L_L({\eu A})$.

Rec\'iprocamente, consideremos $y\in L_L({\eu A})$. Entonces
$v_{\pK}(y)\geq v_{\pK}({\eu A})=v_{\P}({\eu A})$ para toda 
$\pK\in{\ma P}_L$. Sean $[\ell:k]=f$ y $\ell=k(\xi)$. Entonces existen
$a_i\in K$ tales que $y=a_0+a_1 \xi+\cdots+a_{f-1}\xi^{f-1}$.
Sean $y^{ (1)} = y$, $y^{ (2)},
\ldots, y^{ (f)}$ los conjugados de $y$. Tenemos que $y^{
(i)} = a_{ 0} + a_{ 1} \xi^{ (i)} + \cdots + a_{ f
- 1} \left( \xi^{ (i)} \right)^{ f - 1}$.
Resolviendo este sistema de ecuaciones para
$a_{ i}$, obtenemos
\[
a_{ i} = \frac{ \det \left[
\begin{array}{cccccc} 1 & \xi^{ (1)} & \cdots &
y^{ (1) } & \cdots & \left( \xi^{ (1) } \right)^{
f - 1} \\ \cdot & \cdot & \cdots & \cdot & \cdots
& \cdot \\ \cdot & \cdot & \cdots & \cdot & \cdots
& \cdot \\ \cdot & \cdot & \cdots & \cdot & \cdots
& \cdot \\ 1 & \xi^{ (f)} & \cdots & y^{ (f)} &
\cdots & \left( \xi^{ (f)} \right)^{ f - 1}
\end{array} \right] }{ \det \left[
\begin{array}{cccccc} 1 & \xi^{ (1)} & \cdots &
\left( \xi^{ (1)} \right)^{ i} & \cdots & \left(
\xi^{ (1)} \right)^{ f - 1} \\ \cdot & \cdot &
\cdots & \cdot & \cdots & \cdot \\ \cdot & \cdot &
\cdots & \cdot & \cdots & \cdot \\ \cdot & \cdot &
\cdots & \cdot & \cdots & \cdot \\ 1 & \xi^{ (f)}
& \cdots & \left( \xi^{ (f)} \right)^{ i} & \cdots
& \left( \xi^{ (f)} \right)^{ f - 1}
\end{array}
\right] } = \frac{b_{ i}}{c}
\]
con $c \in k^{ \ast}$,  $b_{ i} = \d \sum_{ j=1 }^{ f} t_{
j} y^{ (j)},$ y $ t_{ j} \in \ell$.

Puesto que $y^{(j)}$ es conjugado de $y$ y ${\eu A}\in D_K$,
entonces $y\in L_L({\eu A})$ implica $y^{(j)}\in L_L({\eu A})$. 
Por tanto
\[
v_{\P}(a_i)=v_{\P}\big(\frac{1}{c} \sum_{j=1}^f t_j y^{(j)}\big)\geq
\min_{1\leq j\leq f}v_{\pK}(y^{(j)})\geq v_{\pK}({\eu A})=
v_{\P}({\eu A}),
\]
por lo que $a_i\in L_K({\eu A})$. Se sigue que $\ell L_K({\eu A})
=L_L({\eu A})$. En particular $\ell_L(
{\eu A})=\dim_{\ell} L_L({\eu A})
=\dim_k L_K({\eu A})=\ell_K({\eu A})$.

Consideremos ahora un divisor ${\eu A}\in D_K$ tal que
\begin{align*}
d_K({\eu A})&=d_L({\eu A})>\max\{2g_K-2,2g_L-2\}.
\intertext{Por el Teorema de Riemann--Roch}
\ell_L({\eu A}^{-1})&=d_L({\eu A})-g_L+1,\\
\ell_K({\eu A}^{-1})&=d_K({\eu A})-g_K+1.
\end{align*}
Puesto que $\ell_L({\eu A}^{-1})=\ell_K({\eu A}^{-1})$ y que
$d_L({\eu A})=d_K({\eu A})$, se sigue que $g_L=g_K$.
$\fin$
\end{proof}

\begin{teorema}\label{T6.1.3(3)} Con las mismas hip\'otesis
que antes, sea $\pK$ un lugar de $L$ y sea $\P=\pK|_K$.
Entonces los campos residuales satisfacen $\ell(\pK)=
k(\P)\ell$.
\end{teorema}

\begin{proof}
Claramente $k(\P)\ell \subseteq \ell(\pK)$. Sea ahora $z\in
{\mc O}_{\pK}$. Podemos escribir $z$ como $\sum_{i=0}^{
f-1}a_i\xi^i$, $a_i\in K$ donde $\ell=k(\xi)$, $[\ell:k]=f$. 
Falta probar que $a_i\in{\mc O}_{\P}$. Esto se sigue con el
mismo argumento dado en la demostraci\'on del Teorema
\ref{T6.1.3(2)} ya que si $z^{(i)}$ es un conjugado de $z$, 
entonces existe $\sigma\in\Gal(\ell/k)$ tal que $z^{(i)}=z^{
\sigma}$ y $a^{(i)}=\sum_{i=1}^{f-1}t_i z^{(i)}$ con $t_i\in
\ell$ y $z^{(i)}\in {\mc O}_{\pK^{(i)}}$. Se sigue que $a_i
\in {\mc O}_{\P}$.
$\fin$
\end{proof}

Sean $K$, $L$, $k$ y $\ell$ como antes. Sea $\pK\in{\ma P}_K$ y
sean $\pL_1,\ldots,\pL_h$ los lugares de $L$ sobre $\pK$.
Ahora bien, $\ell/k$ siempre es una extensi\'on de Galois, de
hecho c{\'\i}clica: $\Gal(\ell/k)\cong {\ma Z}/f{\ma Z}\cong C_f$. Por
tanto $L/K$ es una extensi\'on de Galois y $d_{L/K}(\pL_i|\pK)=d$,
$1\leq i\leq h$.

Por los Teoremas \ref{T5.4.5} y \ref{T6.1.3} se sigue que
\[
[L:K]=[\ell:k]=f=dh.
\]
Ahora bien, si $r=d_K(\pK)=[k(\pK):k]$, entonces $k(\pK)\cong
{\ma F}_{q^r}$ y si $s=d_L(\pL_i)=[\ell(\pK_i):\ell]$, entonces
$\ell(\pL_i)\cong {\ma F}_{q^{fs}}$.

Se tiene que $k(\pK)\ell \cong {\ma F}_{q^r}{\ma F}_{q^f}={\ma F}_{
q^{\lcm[r,f]}}={\ma F}_{q^{fs}}\cong \ell(\pL_i)$. Por tanto $fs=\lcm[r,f]=\frac{
rf}{\mcd(r,f)}$ y $s=\frac{r}{\mcd(r,f)}$. Esto es, 
\[
d_L(\pL_i)=s=\frac{d_K(\pK)}{\mcd(d_K(\pK),f)}.
\]
Por otro lado tenemos $d_{L/K}(\pL_i|\pK)d_K(\pK)=d_L(\pL_i)
[\ell:k]$, es decir, $dr=sf$, lo cual equivale a $d=\frac{sf}{r}=\frac{r}{\mcd(r,f)}
\frac{f}{r}=\frac{f}{\mcd(r,f)}$ y $h=\frac{f}{d}=\mcd(r,f)$. En resumen, tenemos:

\begin{teorema}\label{T6.1.4}
Si $\pK$ es un lugar de $K$, $\pL_1,\ldots, \pL_h$ son los lugares de
$L=K\ell$ sobre $\pK$ y $[\ell:k]=f$, entonces
\las
\item $d_{L/K}(\pL_i|\pK)=\frac{f}{\mcd(d_K(\pK),f)}$,
\item $h=\mcd(d_K(\pK),f)$,
\item $d_L(\pL_i)=\frac{d_K(\pK)}{\mcd(d_K(\pK),f)}$. $\fin$
\end{list}
\end{teorema}

\section{Campos ciclot\'omicos}\label{S6.2}

Un campo ciclot\'omico num\'erico es de la forma $\cic n{}$, $\zeta_n
=\exp\big(\frac{2\pi i}{n}\big)$ y $\zeta_n^n=1$. Ahora bien,
el Teorema de Kronecker--Weber dice que toda extensi\'on abeliana
de ${\ma Q}$ est\'a contenida en una extensi\'on ciclot\'omica.
Equivalentemente, se tiene que si ${\ma Q}^{\ab}$ es la m\'axima
extensi\'on abeliana de ${\ma Q}$, entonces ${\ma Q}^{\ab}=
\bigcup_{n=1}^{\infty}\cic n{}$. Notemos que $\zeta_n$ es un
elemento de torsi\'on de ${\ma Z}$ actuando en $\overline{{\ma Q}}^{
\ast}=\overline{{\ma Q}}\setminus \{0\}$, por exponenciaci\'on
donde $\overline{{\ma Q}}$ es una cerradura algebraica de ${\ma Q}$.
Esto es, si $\alpha\in\overline{{\ma Q}}^{\ast}$ y $n\in{\ma Z}$, la
acci\'on est\'a definida por: $n\circ \alpha:=\alpha^n$. Entonces
\begin{gather*}
{\ma Q}^{\ab}={\ma Q}(\tor \overline{{\ma Q}}^{\ast}),
\quad\text{donde} \\
\tor \overline{{\ma Q}}^{\ast}=\text{torsi\'on de\ } \overline{{\ma
Q}}^{\ast}=\{\alpha\in\overline{{\ma Q}}^{\ast}\mid \text{existe $n\in
{\ma N}$ con $\alpha^n=1$}\}.
\end{gather*}

Pretendemos hacer un an\'alogo a todo lo anterior para campos de
funciones congruentes.

Sea $k={\ma F}_q$ y $K$ un campo de funciones racionales sobre
$k$: $K={\ma F}_q(T)$. Sea $R_T:={\ma F}_q[T]$ el anillo de
polinomios sobre ${\ma F}_q$. Se tiene que $K$ es el campo de
cocientes de $R_T$.

Sea $\overline{K}$ una cerradura algebraica de $K$ y $A$ el anillo
de endomorfismos de $\overline{K}$ sobre ${\ma F}_q$:
\begin{align*}
A=\End_{{\ma F}_q}(\overline{K})&=\{\varphi\colon\overline{K}\to
\overline{K}\mid \varphi(a+b)=\varphi(a)+\varphi(b),\\
&\hspace{1cm} \varphi(\alpha a)=
\alpha\varphi(a)\ \forall\ \alpha\in{\ma F}_q, \ \forall\ 
a,b\in\overline{K}\}.
\end{align*}

Entonces $A$ es un anillo y un ${\ma F}_q$--m\'odulo, es decir, en
este caso, ${\ma F}_q$--espacio vectorial, donde la multiplicaci\'on
de $A$ es la composici\'on. El anillo $A$ tiene dos elementos
sobresalientes:

\begin{definicion}\label{D6.2.1}
{\ }

\l
\item Sea $\varphi$ el automorfismo de Frobenius de $\overline{K}$
sobre ${\ma F}_q$, es decir, $\varphi\colon\overline{K}\to
\overline{K}$, $u\mapsto u^q$.

\item Sea $\mu_T$ la multiplicaci\'on por $T\colon \mu_T\colon
\overline{K}\to\overline{K}$, $u\mapsto Tu$.
\end{list}
\end{definicion}

Sea $\xi\colon R_T\to A$ la substituci\'on de $T$ por $\varphi+\mu_T$,
es decir, si $f(T)\in R_T$ es un polinomio, $\xi(f(T))=f(\varphi+\mu_T)
\in A$ es el endomorfismo dado por: si 
\begin{gather*}
f(T)=a_dT^d+\cdots+a_1T+
a_0, 
\intertext{entonces}
f(\varphi+\mu_T)(u)=a_d(\varphi+\mu_t)^d(u)+\cdots + a_1
(\varphi+\mu_T)(u)+a_0 u
\end{gather*} 
para $u\in\overline{K}$. Es decir
$\xi\colon R_T\to A$ est\'a dado por $\xi(T)=\varphi+\mu_T$.

Entonces $\xi$ es un homomorfismo de anillos y bajo $\xi$, $
\overline{K}$ es un $R_T$--m\'odulo, lo cual es el an\'alogo a que
$\overline{{\ma Q}}^{\ast}$ es un ${\ma Z}$--m\'odulo.

Notemos que para $u\in\overline{K}$,
\begin{align*}
(\varphi\circ \mu_t)(u)&=\varphi(Tu)=T^q u^q,\\
(\mu_T^q \circ \varphi)(u)&=\mu_T^q(u^q)=T^q u^q
\end{align*}
es decir, $\varphi\circ\mu_T=\mu_T^q\circ \varphi$ y en particular
$\varphi\circ \mu_T\neq \mu_T \circ \varphi$.

Con el fin de hacer la analog{\'\i}a con $\overline{{\ma Q}}^{\ast}$
y ${\ma Z}$ hacemos la siguiente notaci\'on:

\begin{notacion}\label{N6.2.2}
Si $u\in\overline{K}$ y $M\in R_T$ denotamos $u^M=M(\varphi+
\mu_T)(u)$. Esto es, $u^M=M\circ u=\xi(M)(u)=M(\varphi+\mu_T)(u)$.
\end{notacion}

Se tiene que para $u\in \overline{K}$, $M,N\in R_T$, entonces
\[
u^{M+N}=u^M+u^N\quad\text{y}\quad (u^M)^N=u^{MN}=u^{NM}
=(u^N)^M.
\]

\begin{teorema}\label{T6.2.3}
Sea $M=a_dT^d+\cdots+a_1T+a_0$ con $a_d\neq 0$. Entonces
$u^M=\sum\limits_{i=0}^d \carlitzbinom{M}{i}u^{q^i}$ donde
$\carlitzbinom Mi$ es un polinomio de $R_T$ de grado $(d-i)q^i$
y $\carlitzbinom M0=M$, $\carlitzbinom Md=a_d$.
\end{teorema}

\begin{proof}
Se tiene $u^T=(\varphi+\mu_T)(u)=u^q+Tu$. Vamos a probar por
inducci\'on en el grado de $M$ que $u^M=\sum\limits_{i=0}^d
\carlitzbinom Mi u^{q^i}$ para algunos $\carlitzbinom Mi\in R_T$,
es decir, $\gr_u u^M=q^{\gr_T M}$. Esto se cumple para $d=0$ y
$d=1$ donde $d=\gr_T M$.

Se tiene
\begin{align}
u^{T^{i+1}}&=(u^{T^i})^T=T(u^{T^i})=T\Big(\sum_{j=0}^i
\carlitzbinom {T^i}j u^{q^j}\Big)\nonumber \\
&=(\varphi+\mu_T)\Big(\sum_{j=0}^i \carlitzbinom {T^i}j u^{q^j}\Big)=
\sum_{j=0}^i \carlitzbinom {T^i}j^q u^{q^{j+1}}+\sum_{j=0}^i
T\carlitzbinom {T^i}j u^{q^j}\label{Ec6.2.0}
\end{align}
de donde se sigue lo afirmado. M\'as a\'un, de la expresi\'on
anterior se obtiene
\begin{gather}
u^{T^{i+1}}=\sum_{j=0}^{i+1}\carlitzbinom {T^{i+1}}j u^{q^j}=
\sum_{j=1}^{i+1} \carlitzbinom {T^i}{j-1}^q u^{q^j}+\sum_{j=0}^i
T \carlitzbinom {T^i}j u^{q^j}\nonumber \\
\intertext{por lo que}
\carlitzbinom {T^{i+1}}j=\carlitzbinom {T^i}{j-1}^q+T\carlitzbinom
{T^i}j, \quad 0\leq j\leq i+1\label{Ec6.2.1}
\end{gather}
donde definimos $\carlitzbinom {T^i}{\ell}=0$ si $\ell<0$ o $\ell
>i$.

Por lo lado tenemos que si $M,N\in R_T$, $\gr_T M, \gr_T N
\leq d$, $\alpha,\beta\in {\ma F}_q$, entonces
\begin{gather}
\begin{align*}
u^{\alpha M+\beta N}&=\sum_{i=0}^d \carlitzbinom{\alpha M+\beta N}i
u^{q^i}=\alpha u^M+\beta u^N\\
&=\alpha \sum_{i=0}^d \carlitzbinom Mi u^{q^i}+\beta \sum_{i=0}^d
\carlitzbinom Ni u^{q^i}
\end{align*}
\intertext{de donde}
\carlitzbinom{\alpha M+\beta N}i=\alpha \carlitzbinom Mi+\beta
\carlitzbinom Ni.\label{Ec6.2.2}
\end{gather}
En particular, si $M=a_d T^d+\cdots + a_0$, entonces
\[
\carlitzbinom Mi=\sum_{j=0}^d a_j\carlitzbinom {T^j}i.
\]

Ahora $\carlitzbinom {T^0}i=\carlitzbinom 1i=\begin{cases}
1&\text{si $i=0$}\\ 0&\text{si $i>0$}\end{cases}$ de donde
$\gr_T\carlitzbinom {T^0}{1}=0$ y $\carlitzbinom {T^{j+1}}i=
\carlitzbinom {T^j}{i-1}^q+T\carlitzbinom {T^j}i$.

Por inducci\'on en $j$, suponemos $\gr_T\carlitzbinom{T^j}i=
(j-i)q^i$, $0\leq i\leq j$, entonces 
\begin{align*}
\gr_T \carlitzbinom {T^j}{i-1}^q&=(j-(i-1))q^{i-1}q=(j-i+1)q^i,\\
\gr_T \carlitzbinom {T^j}i&=(j-i)q^i
\end{align*}
y por tanto, de (\ref{Ec6.2.1}) se tiene $\deg_T\carlitzbinom {T^{i+1}}i
(j-i+1)q^i=(j+1-i)q^i$ de donde obtenemos el resultado para $M=T^j$.
El caso general se sigue de (\ref{Ec6.2.2}).

Similarmente, por inducci\'on en $i$ suponemos $\carlitzbinom {T^i}0
=T^i$, entonces de (\ref{Ec6.2.1}) se tiene $\carlitzbinom {T^{i+1}}0=
\carlitzbinom {T^i}{-1}^q+T\carlitzbinom {T^i}0=T^{i+1}$.

Ahora 
\[
\carlitzbinom M0=\sum_{i=0}^d a_i\carlitzbinom {T^i}0=
\sum_{i=0}^d a_iT^i=M\quad\text{y}\quad
\carlitzbinom Md=\sum_{i=0}^da_i\carlitzbinom {T^i}d=a_d
\carlitzbinom {T^d}d=a_d. \tag*{$\fin$}
\]
\end{proof}

Resulta ser que la acci\'on de $R_T$ sobre $\overline{K}$:
$M\circ u:=u^M$ es la an\'aloga a la acci\'on de
${\ma Z}$ sobre $\overline{{\ma Q}}^{\ast}$: $n\circ \xi:=\xi^n$.
El campo ciclot\'omico num\'erico corresponde a
\[
\{\xi\in\overline{\ma Q}^{\ast}\mid \xi^n=1\}=\big\{\zeta_n^a\big\}_{
a=0}^{n-1},
\]
donde $\zeta_n=\exp(2\pi i/n)$.

Por analog{\'\i}a, el campo de funciones ciclot\'omico debe
corresponde a $\{u\in\overline{K}\mid u^M=0\}$.

\begin{definicion}\label{D6.2.4}
Sea $\Lambda_M$ los elementos de $\overline{K}$ que
corresponden a la $M$--torsi\'on de la acci\'on de $R_T$, es decir,
\[
\Lambda_M=\{u\in\overline{K}\mid u^M=0\}.
\]
$\Lambda_M$ recibe el nombre de
{\em m\'odulo de Carlitz\index{Carlitz!m\'odulo de 
$\sim$}\index{m\'odulo de Carlitz}} o {\em m\'odulo de 
Carlitz--Hayes\index{m\'odulo de Carlitz--Hayes}} 
de $M$.

M\'as precisamente, $\Lambda_M$ son los elementos
de torsi\'on del m\'odulo de Carlitz, el cual precisaremos
m\'as adelante.
\end{definicion}

\begin{proposicion}\label{P6.2.5}
Se tiene que $\Lambda_M$ es un $R_T$--subm\'odulo de
$\overline{K}$.
\end{proposicion}

\begin{proof}
Si $u\in \Lambda_M$ y $N\in R_T$, entonces $(u^N)^M=u^{NM}=
u^{MN}=(u^M)^N=0^N=0$, por lo que $u^N\in\Lambda_M$. $\fin$
\end{proof}

\begin{observacion}\label{O6.2.6}
Si $\alpha\in{\ma F}_q^{\ast}$, $\Lambda_{\alpha M}=\Lambda_M$
pues $u^{\alpha M}=(u^{\alpha})^M=(\alpha u)^M=\alpha u^M=0
\iff u^M=0$.

Debido a esto, siempre podemos considerar, sin p\'erdida de
generalidad, polinomios m\'onicos.
\end{observacion}

Los siguientes resultados nos muestran que $\Lambda_M$ es el
equivalente a $W_n=\{\xi\in{\ma C}\mid \xi^n=1\}\cong C_n\cong
{\ma Z}/n{\ma Z}$ el grupo c{\'\i}clico de $n$ elementos. En nuestro
caso, para ser an\'alogo, necesitamos que $\Lambda_M$ sea un
$R_T$--m\'odulo c{\'\i}clico isomorfo a $R_T/M$, es decir,
nuevamente $R_T$ es substituto de ${\ma Z}$ y $M$ de $n$.

\begin{proposicion}\label{P6.2.7} Se tiene que $u^M$ es un
polinomio separable en $u$ de grado $q^d$, donde $M\in R_T$
es de grado $d$, de donde $\Lambda_M$ es un conjunto finito
de $q^d$ elementos. M\'as a\'un, $\Lambda_M$ es un 
${\ma F}_q$--espacio vectorial de dimensi\'on $d$.
\end{proposicion}

\begin{proof} Sea $M=a_dT^d+\cdots+a_1T+a_0$. Entonces
\[
u^M=\sum_{i=0}^d\carlitzbinom Mi u^{q^i},\quad (u^M)'=
\frac{d}{du}(u^M)=\carlitzbinom M0=M\neq 0.
\]
As{\'\i}, $(u^M)'$ no tiene ra{\'\i}ces y $u^M$ es separable en $u$.
Adem\'as $\deg_u u^M=q^d$, por lo que $|\Lambda_M|=q^d$.
Claramente $\Lambda_M$ es un ${\ma F}_q$--espacio vectorial
y por tanto de dimensi\'on $d$. $\fin$
\end{proof}

Para continuar analizando la analog{\'\i}a entre $W_n$ y
$\Lambda_M$, notemos que si $n=p_1^{\alpha_1}\cdots
p_r^{\alpha_r}$ es la descomposici\'on en primos, $W_n\cong
\prod_{i=1}^r W_{p_i^{\alpha_i}}$. El an\'alogo para $\Lambda_M$ es:

\begin{proposicion}\label{P6.2.8}
 Si $M=\prod_{i=1}^r P_i^{\alpha_i}$ es la 
descomposici\'on de $M$ como producto de irreducibles, entonces
\[
\Lambda_M\cong \prod_{i=1}^r\Lambda_{P_i^{\alpha_i}}
\]
como $R_T$--m\'odulos.
\end{proposicion}

\begin{proof}
No es m\'as que un resultado general de m\'odulos sobre dominios
de ideales principales. $\fin$
\end{proof}

Para probar que $\Lambda_M$ es $R_T$--c{\'\i}clico, el paso 
esencial es cuando $M=P^n$, $P$ irreducible.

\begin{proposicion}\label{P6.2.9}
Si $M=P^n$, entonces $\Lambda_{P^n}\cong R_T/P^n$ como
$R_T$--m\'odulos y por lo tanto $\Lambda_{P^n}$ es 
$R_T$--c{\'\i}clico.
\end{proposicion}

\begin{proof}
Lo hacemos por inducci\'on en $n$. Para $n=1$, sea $\xi\in
\Lambda_P\setminus\{0\}$ y sea $\theta\colon R_T\colon R_T
\to \Lambda_P$ dada por $N\mapsto \xi^N$. Entonces $P\in\ker
\theta$ y el ideal $\langle P\rangle$ es maximal. 
Puesto que $\theta(1)=
\xi\neq 0$, se sigue que $\ker\theta=\langle P\rangle$ y $R_T/
\langle P\rangle$ es subm\'odulo de $\Lambda_P$. Puesto
que $|R_T/\langle P\rangle|=|\Lambda_P|=q^d$ donde $d=\gr P$,
se sigue que $R_T/\langle P\rangle \cong \Lambda_P$ y
cualquier $\xi\in \Lambda_P\setminus\{0\}$ es generador.

Supongamos que $\Lambda_{P^n}$ es c{\'\i}clico con generador
$\lambda$. Sea $\mu\colon \Lambda_{P^{n+1}}\to \Lambda_{P^n}$
dada por $\xi\mapsto \xi^P$. Se tiene que $\ker\mu=\Lambda_P$
y $|\Lambda_{P^{n+1}}/\Lambda_{P}|=|\Lambda_{P^n}|=q^{nd}$
por lo que $\Lambda_{P^{n+1}}/\Lambda_{P}\cong \Lambda_{P^n}$,
es decir 
\[
0\to \Lambda_P\to \Lambda_{P^{n+1}}\stackrel{\mu}
{\to}\Lambda_{P^n}\to 0
\]
es una sucesi\'on $R_T$--exacta.

Sea $\xi\in\Lambda_{P^{n+1}}$ tal que $\mu(\xi)=\xi^P=\lambda$
genera $\Lambda_{P^n}$ como $R_T$--m\'odulos. Sea $A$ el
$R_T$--m\'odulo generado por $\xi$, $A=R_T\circ \xi =\xi^{R_T}$.
Se tiene que $A\subseteq \Lambda_{P^{n+1}}$ y $A\cong R_T/
\an(\xi)$ donde $\an(\xi):=\{N\in R_T\mid \xi^N=0\}$.

Ahora bien $P^{n-1}\notin \an (\xi)$ pues $\lambda^{p^{n-1}}=
\mu(\xi^{P^{n-1}})\neq 0$. Sea $\alpha\in \Lambda_{P^{n+1}}$
cualquiera. Entonces $\mu(\alpha)=\alpha^P\in \Lambda_{P^n}=
R_T\circ \lambda$. Por tanto existe $B\in R_T$ tal que 
$\alpha^P=\lambda^B=\mu(\xi^B)=\xi^{PB}$. Entonces
$(\alpha-\xi^B)^P=0$, es decir, $\alpha-\xi^B\in\Lambda_P=\ker \mu$.

Se tiene que $\xi^P$ genera $\Lambda_{P^n}$ por lo que $(\xi^P)^{
P^{n-1}}=\xi^{P^n}\neq 0$ y $\xi^{P^n}\in \Lambda_P$. Por el caso
$n=1$, $\xi^{P^n}$ genera $\Lambda_P$. En particular, existe
$C\in R_T$ tal que $\xi^{P^n C}=\alpha-\xi^B$ o $\xi^{B+P^nC}=
\alpha$, es decir, $\xi$ genera $\Lambda_{P^{n+1}}$ como
$R_T$--m\'odulo. Finalmente, $\langle P^{n+1}\rangle \subseteq
\an(\xi)\subsetneqq \langle P^n\rangle$.

Sea $\an(\xi)=\langle Q\rangle$. Entonces $P^n|Q$, $Q\neq P^n$,
por lo que $Q=P^n Q_1|P^{n+1}$, es decir, $Q_1|P$ y $Q_1$
no es unidad. Se sigue que $Q_1=P$, $Q=P^{n+1}$, $\an(\xi)=
\langle P^{n+1}\rangle$ y $\Lambda_{P^{n+1}}=R_T\circ \xi\cong
R_T/\an(\xi)=R_T/\langle P^{n+1}\rangle$. $\fin$
\end{proof}

\begin{teorema}\label{T6.2.10} Para todo $M\in R_T$, $M\neq 0$,
$\Lambda_M$ es un $R_T$--m\'odulo c{\'\i}clico y $\Lambda_M
\cong R_T/\langle M\rangle$  como $R_T$--m\'odulos.
\end{teorema}

\begin{proof} Es inmediato de las Proposiciones \ref{P6.2.8} y
\ref{P6.2.9}. $\fin$
\end{proof}

Como en grupos, se tiene que:

\begin{proposicion}\label{P6.2.11} Sean $\lambda$ un generador de
$\Lambda_M$ y $A\in R_T$. Entonces $\lambda^A$ es generador
de $\Lambda_M$ si y s\'olo si $A$ y $M$ son primos relativos.
\end{proposicion}

\begin{proof}
Se tiene que $\lambda^A$ genera a $\Lambda_M$ si y s\'olo si $
\lambda\in R_T\circ \lambda^A=\{\lambda^{AB}\mid B\in R_R\}$,
es decir, si y s\'olo si existe $C\in R_T$ tal que $\lambda=
\lambda^{AC}$ lo cual equivale a que $\lambda^{1-AC}=0$. Puesto
que $\an(\lambda)=\langle M\rangle$, $\lambda^{1-AC}=0$ si y 
s\'olo si $M|1-AC$. $\fin$
\end{proof}

\begin{definicion}\label{D6.2.12} Sea $M\in R_T\setminus\{0\}$ y sea
$\Lambda_M$ el m\'odulo de Carlitz--Hayes de $M$. Se define el
{\em campo de funciones ciclot\'omico determinado por
$M$\index{campo de funciones ciclot\'omico}} al campo $K(
\Lambda_M)$.
\end{definicion}

\begin{teorema} Se tiene que si $\lambda_M$ es un generador de
$\Lambda_M$ como $R_T$--m\'odulo, entonces $K(\Lambda_M)=
K(\lambda_M)$ y $K(\Lambda_M)/K$ es una extensi\'on de Galois.
\end{teorema}

\begin{proof}
Se tiene $\lambda_M^{R_T}=\Lambda_M=\{\lambda_M^A\mid A\in
R_T\}$ por lo que $K(\Lambda_M)=K(\lambda_M)$ pues cada
elemento $\xi\in \Lambda_M$ es de la forma
\[
\xi=\lambda_M^A=A(\mu_T+\varphi)(\lambda_M)\in
K\big(\lambda_M^q,\{T^s\lambda_M\}\big)=K(\lambda_M).
\]
En particular, puesto que $\Lambda_M$ es el conjunto de ra{\'\i}ces
del polinomio $u^M\in R_T[u]\subseteq K[u]$ y $u^M$ es
separable (Proposici\'on \ref{P6.2.7}), $K(\Lambda_M)/K$ es normal
y separable, es decir, Galois. $\fin$
\end{proof}

\begin{definicion}\label{D6.2.14} Para $M\in R_T\setminus\{0\}$, 
$G_M$ denota al grupo de Galois de $K(\Lambda_M)/K$:
$G_M:=\Gal(K(\Lambda_M)/K)$.
\end{definicion}

Las acciones de $G_M$ y de $R_T$ sobre $\Lambda_M$ conmutan:
si $\sigma\in G_M$ y $N\in R_T$, entonces existe $A\in R_T$,
$\mcd(A,M)=1$ tal que $\sigma(\lambda_M)=\lambda_M^A\in
\Lambda_M$ y si $\lambda\in\Lambda_M$, existe $N\in R_T$
tal que $\lambda=\lambda_M^N$, por lo que
\begin{align*}
\sigma(\lambda)&=\sigma(\lambda_M^N)=
\sigma\Big(\sum_{i=0}^{\deg N}\carlitzbinom Ni \lambda_M^{q^i}\Big)=
\sum_{i=0}^{\deg N}\carlitzbinom Ni \sigma(\lambda_M)^{q^i}
=(\sigma(\lambda_M))^N\\
&=(\lambda_M^A)^N=\sum_{i=0}^{\deg N}\carlitzbinom Ni 
(\lambda_M^A)^{q^i}=((\lambda_M)^A)^N=(\lambda_M)^{AN}=
(\lambda_M)^{NA}.
\end{align*}
Esto es, $\sigma(\lambda_M^N)=\lambda_M^{NA}=
\lambda_M^{AN}=(\lambda_M^A)^N=\sigma(\lambda_M)^N$.

En el caso num\'erico, tenemos que $\Gal(\cic n{}/{\ma Q})\cong
({\ma Z}/n{\ma Z})^{\ast}$. Por tanto, en el caso de campos de
funciones podr{\'\i}amos esperar que $G_M\cong (R_T/\langle
M\rangle)^{\ast}$ el grupo de unidades del anillo $R_T/\langle
M\rangle$. Ahora bien, notemos que
\[
(R_T/\langle M\rangle)^{\ast}=\{A\bmod M\mid (A,M)=1\}
\]
y que el n\'umero de generadores de $\Lambda_M$
 es precisamente $\big|(R_T/\langle
M\rangle)^{\ast}\big|$ (Proposici\'on \ref{P6.2.11}). Ahora $\big|(
{\ma Z}/n{\ma Z})^{\ast}\big|$ est\'a dada por la funci\'on fi de Euler.
Aqu{\'\i} tenemos el an\'alogo.

\begin{definicion}\label{D6.2.15} Se define la funci\'on ``Fi'' de Euler
por: 
\[
\Phi(M):=\big|(R_T/\langle M\rangle)^{\ast}\big|
\]
para $M\in R_T\setminus\{0\}$.
\end{definicion}

Puesto que $\langle M\rangle=\langle \alpha M\rangle$
para $\alpha\in {\ma F}_q^{\ast}$, $\Phi(\alpha M)=\Phi(M)$ para
$\alpha\in {\ma F}_q^{\ast}$.

Por el Teorema Chino del Residuo, tenemos que si $M$ y $N$ son
primos relativos, entonces $R_T/\langle MN\rangle\cong R_T/\langle
M\rangle\times R_T/\langle N\rangle$ y $\big(R_T/\langle
MN\rangle\big)^{\ast}\cong \big(R_T/\langle
M\rangle\big)^{\ast}\times \big(R_T/\langle N\rangle\big)^{\ast}$.
En particular se tiene $\Phi(MN)=\Phi(M)\Phi(N)$ para $M$ y $N$
primos relativos.

Por otro lado, si $P$ es irreducible, entonces $R_T/\langle P\rangle$
es el campo de $q^d$ elementos donde $d=\gr P$. Por tanto
$\big|\big(R_T/\langle P \rangle\big)^{\ast}\big|=q^d-1=\Phi(P)$. Para
$n\geq 2$, se tiene que si $\big(R_T/\langle P^n \rangle\big)^{\ast}
\stackrel{\theta}{\longto}\big(R_T/\langle P^{n-1} \rangle\big)^{\ast}$
es el mapeo natural, $\theta(A\bmod P^n)=A\bmod P^{n-1}$, 
entonces $\ker \theta=\{ 1+P^{n-1}B\bmod P^n\mid B\in R_T\}\cong
R_T/\langle P\rangle$. En particular
\[
0\longto R_T/\langle P\rangle \stackrel{\mu}{\longto}
\big(R_T/\langle P^n \rangle\big)^{\ast}\stackrel{\theta}{\longto}
\big(R_T/\langle P^{n-1} \rangle\big)^{\ast}\longto 1
\]
es exacta, donde $\mu(B\bmod P)=1+P^{n-1}B\bmod P^n$.

En efecto $\mu\colon R_T/\langle P\rangle\lra \ker \theta$
es un isomorfismo: $\mu((B+C)\bmod P)=
\big(1+P^{n-1}(B+C)\big)\bmod P^n$ y se tiene que $\big((1+P^{n-1}B)
(1+P^{n-1}C)\big)\bmod P^n=\big(1+P^{n-1}B+P^{n-1}C+P^{2n-2}BC
\big)\bmod
P^n\equiv \big(1+P^{n-1}(B+C)\big)\bmod P^n$ pues $2n-2\geq n$
para $n\geq 2$. Adem\'as $\mu$ es inyectiva pues is $\mu(B\bmod P)
=\mu(C\bmod P)$, entonces $1+P^{n-1}B\equiv 1+P^{n-1}C\bmod P^n$.
Por tanto $P|B-C$ y $B\bmod P= C\bmod P$. El isomorfismo se sigue
de que $|R_T/\langle P\rangle|=q^d=|\ker \theta|$.

De inmediato obtenemos dos resultados:
\l
\item $\Phi(P^n)=\Phi(P^{n-1})\big|R_T/\langle P\rangle\big|=
\Phi(P^{n-1}) q^d$,
\item $R_T/\langle P\rangle\cong D_{P^n,P^{n-1}}=\{\overline{A}\in
\big(R_T/\langle P^n \rangle\big)^{\ast}\mid A\equiv 1\bmod P^{n-1}
\}=\ker \theta$.
\end{list}

Si por inducci\'on suponemos $\Phi(P^n)=q^{nd}-q^{(n-1)d}$, entonces
$\Phi(P^{n+1})=\Phi(P^n)q^d=(q^{nd}-q^{(n-1)d})q^d=q^{(n+1)d}-
q^{nd}$.

De manera an\'aloga, si 
\begin{gather*}
\pi\colon \big(R_T/\langle P^n
\rangle\big)^{\ast}\longto \big(R_T/\langle P \rangle\big)^{\ast}\\
\intertext{est\'a dada por $\pi(A\bmod P^n)=A\bmod P$, se tiene}
\ker \pi =D_{P^n,P}=\{\overline{A}\in \big(R_T/\langle P^n
\rangle\big)^{\ast}\mid A\equiv 1\bmod P\}
\rangle\\
\intertext{y en general se tiene, para $1\leq m<n$, que la sucesi\'on}
1\longto D_{P^n,P^m}\longto\big(R_T/\langle
 P^n \rangle\big)^{\ast}
\stackrel{\varphi}{\longto}\big(R_T/\langle P^m \rangle\big)^{\ast}
\longto 1\\
\intertext{es exacta donde $\varphi(A\bmod P^n)=A\bmod P^m$, y}
D_{P^n,P^m}=\{A\bmod 
P^n\in \big(R_T/\langle P^n \rangle\big)^{\ast}
\mid A\equiv 1\bmod P^m\}.
\end{gather*}

Se tiene una biyecci\'on $\mu\colon
R_T/\langle P^{n-m}\rangle \lra D_{P^n,P^m}$
dada por $B\bmod P^{n-m}\longmapsto 1+P^m B\bmod P^n$,
donde se toma como conjunto de representantes de 
$R_T/\langle P^{n-m}\rangle$ a ${\mc A}=\{B\in R_T\mid \deg B<
\deg P^{n-m}\}$. Hacemos notar que $\mu$ no es un homomorfismo
en general. Ahora bien, $\mu$ si es homomorfismo 
para $m=n-1$ y en particular para $n=2$.
Resumiendo, tenemos:

\begin{proposicion}\label{P6.1.16}
La funci\'on $\Phi$ de Euler satisface:
\las
\item $\Phi(MN)=\Phi(M)\Phi(N)$ para $M,N$ primos relativos.

\item $\Phi(P^n)=\big|R_T/\langle P^{n-1}\rangle\big|\Phi(P)=
q^{nd}-q^{(n-1)d}$ donde $P$ es irreducible, $n\geq 1$ y $d=\gr P$.

\item Para $1\leq m<n$ y $P$ irreducible, la sucesi\'on
\[
0 \lra D_{P^n,P^m}\longto\big(R_T/\langle
 P^n \rangle\big)^{\ast}
\stackrel{\varphi}{\longto}\big(R_T/\langle P^m \rangle\big)^{\ast}
\longto 1
\]
es exacta y $\mu\colon R_T/\langle P^{n-m}\rangle\lra D_{P^n,P^m}$,
$B\bmod P^{n-m}\longmapsto 1+P^m B\bmod P^n$ es una biyecci\'on
donde tomamos como conjunto de representantes de $R_T/\langle
P^{n-m}\rangle$ al conjunto ${\mc A}=\{B\in R_T\mid \deg B<
\deg P^{n-m}\}$.

\item $\Phi(M)$ es el n\'umero de generadores de $\Lambda_M$.
$\fin$
\end{list}
\end{proposicion}

\begin{proposicion}\label{P6.2.17} Sea $M\in R_T\setminus\{0\}$.
Entonces $G_M\subseteq \big(R_T/\langle M\rangle\big)^{\ast}$.
\end{proposicion}

\begin{proof}
Se tiene que si $\lambda=\lam M{}$ es un generador de $\Lam M{}$,
entonces $\cicl M{}=K(\lambda_M)$. Por tanto $\sigma\in G_M$
est\'a determinado por su acci\'on en $\lambda$ y puesto que
$\sigma \lambda$ es un conjugado de $\lambda$, $\sigma \lambda
\in \Lam M{}$. Si $\sigma \lambda=\beta$, entonces $\lambda=
\sigma^{-1}\beta$ por lo que necesariamente $\sigma\lambda$ es
generados $\Lam M{}$. Se sigue que $\sigma \lambda=\lambda^A$
para $\mcd (A,M)=1$. Denotemos $\sigma$ por $\sigma_A$, es
decir, $\sigma_A\lambda=\lambda^A$, $\mcd (A,M)=1$.

Sea $\varphi\colon G_M\to \big(R_T/\langle M\rangle\big)^{\ast}$,
$\sigma_A\mapsto A\bmod M$. Claramente $\varphi$ es un 
monomorfismo de grupos. $\fin$
\end{proof}

\begin{observacion}\label{O6.2.17'} Para $\sigma \in G_M$,
con $\sigma\stackrel{\varphi}{\longmapsto} A$, 
se tiene que para $\xi\in\Lambda_M$,
la acci\'on de $\sigma$ es $\sigma(\xi)=\xi^A$ pero para
$\xi\notin \Lambda_M$ no necesariamente $\sigma (\xi)=
\xi^A$. Por ejemplo, si $\alpha\in\*\F$, entonces $\alpha^A
\neq \alpha$ y $\sigma(\alpha)=\alpha$.
\end{observacion}

\begin{corolario}\label{C6.2.18} Para $M\in R_T\setminus\{0\}$,
$[\cicl M{}:K]\leq \Phi(M)$. $\fin$
\end{corolario}

\begin{definicion}\label{D6.2.19} Sea $S\in R_T$ un polinomio
m\'onico. Definimos el {\em polinomio
$S$--ciclot\'omico\index{polinomio ciclot\'omico}} por
\[
\Psi_S(u):=\prod_{\substack{\mcd (B,S)=1\\ \gr B<\gr S}} (u-
\lambda_S^B)
\]
donde $\lam S{}$ es un generador de $\Lam S{}$. Se tiene
$\Psi_S(u)\in \cicl S{}[u]$.
\end{definicion}

Notemos que $\Psi_S(u)$ es el an\'alogo al polinomio ciclot\'omico 
usual 
\[
\psi_n(x)=\prod_{\substack{\mcd(m,n)=1\\ 0\leq m\leq n}}
(x-\zeta_n^m).
\]

Result\'o ser que $\psi_n(x)\in{\ma Q}[x]$ es irreducible de grado
$\varphi(n)=\big({\ma Z}/n{\ma Z}\big)^{\ast}$. Veremos que $\psi_S(u)
\in K[u]$ es irreducible de grado $\Phi(S)$. Ahora bien, $\gr_u\Psi_S(
u)=\Phi(S)$ por la Proposici\'on \ref{P6.2.11}.

\begin{proposicion}\label{P6.2.20} Se tiene que $\Psi_S(u)\in K[u]$.
\end{proposicion}

\begin{proof}
Sea $\sigma_A\in G_S$. Entonces $\sigma_A(\lam S{})=\lam S{}^A$
y si $\mcd (B,S)=1$, $\sigma_A(\lam S{}^B)=\lam S{}^{BA}$ con
$\mcd (AB,S)=1$. Por lo tanto 
\[
\sigma_A(\Psi_S(u))=\prod\limits_{
\substack{\mcd (B,S)=1\\ \gr B<\gr S}}(u-\lam S{}^{AB}).
\]

Tomando $AB\bmod S$ y puesto que multiplicaci\'on por $A$ es
un automorfismo de $\big(R_T/\langle S\rangle\big)^{\ast}$, se sigue
que $\sigma_A(\Psi_S(u))=\Psi_S(u)$, $\sigma_A\in G_M$ de 
donde $\Psi_S(u)\in K[u]$. $\fin$
\end{proof}

Como en el caso num\'erico tenemos:

\begin{proposicion}\label{P6.2.21} {\ }

\las
\item Si $M,N\in R_T$ son polinomios m\'onicos con
 $M\neq N$, entonces $\mcd \big(
\Psi_M(u),\Psi_N(u)\big)=1$.

\item $u^M=\prod\limits_{\substack{N|M\\ N\text{\ m\'onico}}}\Psi_N(u)$,
$M\in R_T$ m\'onico.

\item $\Psi_M(u)=\prod\limits_{\substack{N|M\\ N \text{\ m\'onico}}} 
\big(u^N\big)^{\mu(M/N)}$ donde 
\[
\mu(D)=\begin{cases}
1& \text{si $D=1$},\\
(-1)^s& \text{si $D=P_1\cdots P_s$ con 
$P_1,\ldots,P_s$ m\'onicos e}\\ & \text{irreducibles distintos}\\
0& \text{en otro caso}
\end{cases}
\]
y $M$ m\'onico.
\end{list}
\end{proposicion}

\begin{proof}
{\ }

\las
\item Sea $D:=\mcd \big(\Psi_M(u),\Psi_N(u)\big)$. Si $D\neq 1$,
sea $\lambda\in \overline{K}$ ra{\'\i}z de $D$. Por tanto $\lambda$
es ra{\'\i}z de $\Psi_M(u)$ y $\Psi_N(u)$, es decir, $\lambda=
\lam M{}^A=\lam N{}^B$ para $\mcd (A,M)=1$, $\mcd (B,N)=1$
con $\deg A<\deg M$ y $\deg B<\deg M$. Por tanto se tiene 
$\lambda=\lambda_{MN}^{AN}=\lambda_{NM}^{BM}$ lo cual 
implica que $AN=BM$. Puesto que $B$ y $N$ son primos
relativos lo mismo que $A$ y $M$, se sigue que $B|A$ y $A|B$
con lo que se concluye que $A=B$ y por tanto $M=N$ lo cual
es contrario a lo supuesto.

Por tanto $D=\mcd \big(\Psi_M(u),\Psi_N(u)\big)=1$.

\item Claramente, si $N|M$, entonces $\Psi_N(u)|u^M$ y como
para $N_1\neq N_2$, $\mcd \big(\Psi_{N_1}(u),\Psi_{N_2}(u)\big)=1$,
se sigue que $\prod\limits_{\substack{N|M\\ N\text{\ m\'onico}}}
\Psi_N(u)$ divide a $u^M$.

Ahora $\gr \Big(\prod\limits_{\substack{N|M\\ N\text{\ m\'onico}}}
\Psi_N(u)\Big)=\sum\limits_{\substack{N|M\\ N\text{\ m\'onico}}}
\Phi(N)$.

Si $M=P^n$ con $P$ irreducible y $\gr P=d$, se tiene 
\begin{align*}
\sum_{N|M}\Phi(N)&=\sum_{i=0}^n\Phi(P^i)=
\sum_{i=1}^n\big(q^{id}-q^{(i-1)d}\big)+1\\
&= q^{nd}-1+1=q^{nd}=q^{\gr P^n}.\\
\intertext{En general si $M=P_1^{\alpha_1}\cdots P_r^{\alpha_r}$}
\sum_{\substack{N|M\\ N\text{\ m\'onico}}}\Phi(N)&=
\sum_{\substack{\beta_i=0\\ i=1,\ldots,r}}\Phi(P_1^{\beta_1})\cdots
\Phi(P_r^{\beta_r})=\prod_{i=1}^r \sum_{\beta_i=0}^{\alpha_i}
\Phi(P_i^{\beta_i})\\
&=\prod_{i=1}^r q^{\gr P_i^{\beta_i}}=q^{\gr M}.
\end{align*}

Por lo tanto $\gr\Big(\prod\limits_{\substack{N|M\\ N\text{\ m\'onico}}}
\Psi_N(u)\Big)=q^{\gr M}=\gr u^M$, lo cual implica $u^M=\prod\limits_{N|M}
\Psi_N(u)$.

\item La demostraci\'on es totalmente an\'aloga a la del Corolario
\ref{C1.2.14} (2). $\fin$
\end{list}
\end{proof}

Ahora veamos que la ramificaci\'on en campos de funciones es
totalmente paralela a la de los campos ciclot\'omicos num\'ericos.

\begin{definicion}\label{D6.2.23}
Al polo $\pK_{\infty}$ de $T$ en $K$ lo llamamos {\em primo
infinito\index{primo infinito}}.
\end{definicion}

\begin{proposicion}\label{P6.2.24} Sean $P\in R_T$ m\'onico e
irreducible de grado $d$ y $M=P^n$, $n\in{\ma N}$. Entonces
\l
\item Si ${\eu q}$ es cualquier otro divisor primo distinto a
$\pK_{\infty}$ y $\pK$, donde $\pK$ es el primo asociado a $P$,
entonces ${\eu q}$ no es ramificado.
\item El {\'\i}ndice de ramificaci\'on de ${\eu p}$ en $K(\Lam M{})/K$
es
\[
e(\pK)=\Phi(M)=q^{dn}-q^{d(n-1)}=[\cicl M{}:K].
\]
\end{list}
\end{proposicion}

\begin{proof}{\ }
\begin{window}[1,l,\xymatrix{
{\cal O}_M\ar@{-}[r]\ar@{-}[d] & \cicl M{}\ar@{-}[d]\\
R_T\ar@{-}[r] & K},{}]
Sea ${\cal O}_M$ la cerradura entera de $R_T$ en $\cicl M{}$.
Puesto que $R_T$ es un dominio Dedekind, ${\cal O}_M$ tambi\'en
lo es (ver Teorema \ref{DrinfeldT1.2.4}).
Se tiene que los primos $K$ ramificados en $\cicl M{}/K$, 
diferente a $\pK_{\infty}$, son aquellos que aparecen en el
discriminante ${\eu d}_{{\cal O}_M/R_T}$.
Sea $\lambda$ generador de $\Lam M{}$. Entonces $R_T[\lambda]
\subseteq {\cal O}_M$. Sean $g(u):=\Irr(\lambda,u,K)\in K[u]$ y
$f(u):=u^M$. Puesto que $f(\lambda)=0$, existe $h(u)\in K[u]$ tal que
$f(u)=h(u)g(u)$. Por tanto
\end{window}
\begin{gather*}
M=f'(u)=h'(u)g(u)+h(u)g'(u)\\
\intertext{de donde}
M=f'(\lambda)=h'(\lambda)g(\lambda)+h(\lambda)g'(\lambda)=
h(\lambda)g'(\lambda).\\
\intertext{En especial $\big(g'(\lambda)\big)_{{\cal O}_M}|(M)_{{\cal
O}_M}= P^n{\cal O}_M$. Se tiene (Teorema \ref{T5.5.6})}
{\eu D}_{{\cal O}_M/R_T}=\langle F'(\alpha)\mid \alpha\in{\cal O}_M,
\cicl M{}=K(\alpha), F(u)=\Irr(\alpha,u,K)\rangle.
\end{gather*}

En particular
\begin{equation}\label{Ec6.2.25}
{\eu D}_{{\cal O}_M/R_T}|\big(g'(\lambda))_{\cicl M{}}=P^n=
(\pK_1\cdots \pK_h)^{en}\quad\text{donde}\quad
P{\cal O}_M=\langle \pK_1\cdots \pK_h\rangle^e.
\end{equation}

Por tanto los \'unicos posibles primos ramificados en $\cicl M{}/K$ son
$\pK$ y $\pK_{\infty}$. Ahora
\begin{align*}
u^{P^n}&=(u^{P^{n-1}})^P=\sum_{i=0}^d \carlitzbinom Pi 
(u^{P^{n-1}})^{q^i}\\
&=u^{P^{n-1}}\Big(\sum_{i=0}^d\carlitzbinom Pi (u^{P^{n-1}})^{
q^i-1}\Big)=u^{P^{n-1}}t(u)
\end{align*}
con $t(u)=\frac{u^{P^n}}{u^{P^{n-1}}}=\sum\limits_{i=0}^d 
\carlitzbinom Pi (u^{P^{n-1}})^{q^i-1}$.

Se tiene $t(\alpha)=0\iff \alpha\in \Lam Pn\setminus \Lam P{n-1}\iff
\alpha$ es generador de $\Lam Pn$. Por tanto $t(u)=\Psi_{P^n}(u)$ y
\begin{align*}
t(u)&=\prod_{\mcd (A,M)=1}(u-\lambda^A)=\carlitzbinom P0+
\sum_{i=1}^d \carlitzbinom Pi (u^{P^{n-1}})^{q^i-1}\\
&=P+\sum_{i=1}^d \carlitzbinom Pi (u^{P^{n-1}})^{q^i-1}.
\end{align*}
Con $u=0$, se tiene
\begin{equation}\label{Ec6.2.26}
t(0)=\pm \prod_{\mcd (A,M)=1}\lambda^A=P.
\end{equation}

Ahora bien $u^A=um_A(u)$ con $m_A(u)\in R_T(u)$. En particular
$\lambda^A=\lambda F(\lambda)$ y $\lambda|\lambda^A$. Para
$(A,M)=1$, $\lambda^A$ es generador y por simetr{\'\i}a se sigue
que $\lambda^A|\lambda$, es decir $\lambda=\beta_A \lambda^A$
con $\beta_A\in {\cal O}_M^{\ast}$. Se sigue de (\ref{Ec6.2.26})
que $\pm P=\beta_0\lambda^{\Phi(M)}$ para $\beta_0\in
{\cal O}_M^{\ast}$. De (\ref{Ec6.2.25}) obtenemos $\langle
P\rangle_{{\cal O}_M}
=\langle\pL_1\cdots\pL_h
\rangle^e=(\lambda)^{\Phi(M)}$ y en particular
$v_{\pL_i}(\lambda)\geq 1$. Entonces $e=v_{\pL_i}((\pL_1\cdots
\pL_h)^e)=v_{\pL_i}(\lambda^{\Phi(M)})\geq \Phi(M)$, esto es,
\[
e\geq \Phi(M)=\big|\big(R_T/\langle M\rangle\big)^{\ast}\big|\geq
[\cicl M{}:K]\geq e
\]
de donde $e=\Phi(M)=\Phi(P^n)=[\cicl Pn:K]=q^{dn}-q^{d(n-1)}$.
$\fin$
\end{proof}

\begin{corolario}\label{C6.2.26'} Con las notaciones anteriores,
se tiene que $v_{\pL}(\lambda)=1$ donde $\pL$ es el
\'unico lugar de $\cicl Pn$ sobre $\pK$. $\fin$
\end{corolario}

\begin{observacion}\label{O6.2.27} De paso hemos obtenido que el
polinomio ciclot\'omico $\Psi_{P^n}(u)$ satisface
\[
\Psi_{P^n}(u)=\frac{u^{P^n}}{u^{P^{n-1}}}=\frac{\prod\limits_{N|P^n}
\Psi_N(u)}{\prod\limits_{N|P^{n-1}}\Psi_N(u)}=
\prod_{\mcd (A,M)=1}(u-\lambda^A).
\]
El caso general es consecuencia de la Proposici\'on \ref{P6.2.24}.
\end{observacion}

\begin{teorema}\label{T6.2.28} Si $M\in R_T\setminus\{0\}$ un
polinomio m\'onico. Entonces
\las
\item $t(u)=\Irr(\lambda,u,K)=\Psi_M(u)$ donde $\lambda$ es 
generador de $\Lam M{}$. En particular $\Psi_M(u)$
es irreducible.
\item $G_M=\Gal(\cicl M{}/K)\cong \big(R_T/\langle M\rangle\big)^{
\ast}$.
\item $[\cicl M{}:K]=\Phi(M)$.
\item Si $M=P^n$ donde $P$ es irreducible, entonces $\pK$ es 
totalmente ramificado en $\cicl Pn/K$ donde $\langle P\rangle_K=
\frac{\pK}{\pK_{\infty}^{\gr P}}$.
\end{list}
\end{teorema}

\begin{proof}
Si $M=P^n$, por la Proposici\'on \ref{P6.2.24}, se tiene
\[
[\cicl Pn:K]=\Phi(P^n)=\big|\big(R_T/\langle P^n\rangle\big)^{\ast}\big|
=\big|G_{P^n}\big|.
\]
Por otro lado, por la Proposici\'on \ref{P6.2.17}, se tiene $G_{P^n}
\subseteq \big(R_T/\langle P^n\rangle\big)^{\ast}$. Por tanto
$G_{P^n}\cong \big(R_T/\langle P^n\rangle\big)^{\ast}$ y $P$ es
totalmente ramificado pues $e=\Phi(P^n)=[\cicl Pn:K]$. Esto es (4).

En general, si $M=P_1^{\alpha_1}\cdots P_r^{\alpha_r}$ con 
$P_1,\ldots, P_r$ polinomios irreducibles distintos, $\Lam M{}\cong
\bigoplus\limits_{i=1}^r \Lam {P_i}{\alpha_i}$. Si probamos que
$[\cicl M{}:K]=\Phi(M)$ entonces, puesto que $G_M\subseteq
\big(R_T/\langle M\rangle\big)^{\ast}$, se seguir\'a la igualdad
y tambi\'en (2) y (3). Finalmente (1) se seguir\'a del hecho de que $t(\lambda)
=0$, $\gr_u t(u)=\Phi(\Psi_M(u))=\gr \Irr(\lambda,u,K)$ y de que
$\Irr(\lambda,u,K)| t(u)$.

En resumen, solo falta probar que $\Phi(M)=[\cicl M{}:K]$. Ahora
bien, puesto que $\pK_i$ es totalmente ramificado en $\cicl {P_i}{
\alpha_i}/K$ y no ramificado en $\prod\limits_{j\neq i}\cicl {P_j}{
\alpha_j}/K$, se sigue que
\[
[\cicl M{}:K]=\prod_{i=1}^r[\cicl {P_i}{\alpha_i}:K]=
\prod_{i=1}^r\Phi(P_1^{\alpha_i})=\Phi(M). \tag*{$\fin$}
\]
\end{proof}

\begin{corolario}\label{C6.2.29} El campo de constantes de 
$\cicl M{}$ es ${\ma F}_q$ y $\cicl M{}/K$ es una extensi\'on
geom\'etrica.
\end{corolario}

\begin{proof}
Sea $E_i:=K(\Lambda_{M/P_i^{\alpha_i}})$, $1\leq i\leq r$ donde
$M=P_1^{\alpha_1}\cdots P_r^{\alpha_r}$. Se tiene que $\Gal(\cicl
M{}/E_i)\cong \Gal (\cicl {P_i}{\alpha_i}/K)$. Sea $L$ la m\'axima
extensi\'on no ramificada de $K$ contenida en $\cicl M{}$. Ahora
bien, $\cicl M{}/E_i$ es totalmente ramificada en los primos que
est\'an sobre $\pK_i$ y $E_iL/E_i$ es no ramificada, por lo que
$E_iL=E_i$. Entonces $L\subseteq E_i$ para $1\leq i\leq r$.

Se sigue que $K\subseteq L\subseteq \cap_{i=1}^r E_i=K$. En 
particular $L=K$ y cada extensi\'on $K\subsetneqq F\subseteq
\cicl M{}$ es ramificada. Sea ${\ma F}_{q^s}$ el campo de 
constantes de $\cicl M{}$. Entonces
\[
K={\ma F}_q(T)\subseteq {\ma F}_{q^s}(T)\subseteq \cicl M{}.
\]
Se tiene que ${\ma F}_{q^s}(T)/{\ma F}_q(T)$ es no ramificada
(Teorema \ref{T6.1.3}). Por tanto ${\ma F}_{q^s}(T)={\ma F}_q(T)$
y ${\ma F}_{q^s}={\ma F}_q$, es decir, $s=1$. $\fin$
\end{proof}

\subsection{Estructura del grupo $\Gal(\cicl M{}/k)$}\label{S9.2.0}

Sea $M=\poly$ con $P_1,\ldots,P_r\in R_T^+$ distintos. Se tiene que
\begin{gather*}
\G M\cong \G {P_1^{\alpha_1}}\times \cdots\times \G {P_r^{\alpha_r}}
\intertext{y que}
\G {P^n}\cong D_{P^n,P}\times \G P,
\end{gather*}
para $P\in R_T^+$ y $n\in {\ma N}$, donde $D_{P^n,P}=
\{A\bmod P^n\mid A\equiv 1\bmod P\}$. Por tanto, 
para conocer $G_M=\Gal(
\cicl M{}/k)$, basta describir $D_{P^n,P}$ y $\G P$
para $P\in R_T^+$ y $n\in{\ma N}$.

Primero notemos que $\Gal(\cicl Pn/\cicl P{n-1})$ es un grupo
$p$--elemental abeliano pues
\[
\Gal(\cicl Pn/\cicl P{n-1})\cong D_{P^n,P^{n-1}}\cong R_T/\langle P\rangle.
\]
Se tiene que el grupo $R_T/\langle P\rangle$ es $p$--elemental abeliano
pues para cualquier $f\in R_T$, $p\cdot f=0$. Por tanto $R_T/\langle
P\rangle \cong C_p^{vd}$ como grupos, donde $q=p^v$ y $d=\deg P$.

Ahora bien, $\G P$ es el grupo multiplicativo de un campo, por lo
que es un grupo c\'iclico de orden $q^d-1$. Se sigue que que
\[
\G {P^n}\cong D_{P^n,P}\times C_{q^d-1}
\]
como grupos pues $D_{P^n,P}$ es el $p$--subgrupo de Sylow de
$\G {P^n}$.

En otro orden de ideas, se tiene que $[\cicl Pn:K]
=\Phi(P^n)=q^{(n-1)d}(q^d-1)$,
por lo que $\cicl Pn=K\iff q^{(n-1)d}=1$ y $q^d-1=1 \iff n=1, d=1$ y $q=2
\iff P^n\in\{T,T+1\}\subseteq {\ma F}_2(T)$.

De lo anterior, obtenemos que $\cicl M{}=K \iff q=2$ y $M\in\{T,T+1,T(T+1)\}$.
En otras palabras, $T$ y $T+1$ y $q=2$ juegan el papel de $\zeta_2$
en el caso ciclot\'omico num\'erico. De esta forma, siempre que se
considere $\cicl M{}$ con $q=2$, se supondr\'a que si $T|M$ o $(T+1)|
M$, entonces $T^2|M$ o $(T+1)^2|M$ respectivamente.

\begin{ejemplo}[Campos ciclot\'omicos con grupo de Galois 
c\'iclico]\label{Ej9.2.0.1}
Sea $M\in R_T$ con $M$ m\'onico no constante. Suponemos que, cuando
$q=2$, entonces si $N\in\{T,T+1\}$ y $N|M$, entonces $N^2|M$.
Sea $M=\poly$ con $d_i=\deg P_i$, $1\leq i\leq r$ tal que 
\begin{align*}
G_M&=\Gal(\cicl M{}/K)=\G M\cong 
\G {P_1^{\alpha_1}}\times \cdots\times \G {P_r^{\alpha_r}}\\
&\cong \Big(D_{P^{\alpha_1},P} \times \G {P_1}\Big)\times
\cdots \times \Big(D_{P^{\alpha_r},P} \times \G {P_r}\Big)\\
&\cong \big(D_{P^{\alpha_1},P}\times C_{q^{d_1}-1}\big)\times\cdots
\times \big(D_{P^{\alpha_r},P}\times C_{q^{d_r}-1}\big)
\end{align*}
es un grupo c\'iclico. Notemos que si $r\geq 2$,
puesto que $q-1|q^{d_i}-1$, necesariamente se debe tener $q=2$.

Primero consideremos el caso $r=1$, esto es, $M=P^n$ con $P\in
R_T^+$, $\deg P=d$ y $n\geq 1$. Cuando $q=2$ y $d=1$ se 
supone $n\geq 2$. En general, cuando $n\geq 2$,
se tiene que $G_M\cong D_{P^{n},P}\times
C_{q^d-1}$, por lo que $G_M$ es c\'iclico si y solamente si
$D_{P^{n},P}$ es c\'iclico. Ahora bien, $G_{P^2}$ es un grupo
cociente de $G_{P^n}$ por lo que $G_{P^2}$ debe ser c\'iclico.
Por tanto $D_{P^{2},P}\cong R_T/\langle P\rangle\cong C_p^{vd}$
debe ser c\'iclico, lo cual implica que $v=1,d=1$. En resumen,
si $G_{P^n}$ es c\'iclico para $n\geq 2$, necesariamente $q=p$ y
$d=1$.

Si $p\geq 3$, entonces $G_{P^3}$ no es c\'iciclo pues $D_{P^3,P}$
es de exponente $p$, es decir, es $p$--elemental abeliano pues
para cualquier $A\in R_T$ se tiene que $(1+AP\bmod P^3)^p=
1+A^pP^p\bmod P^3\equiv 1\bmod P^3$. Adem\'as $|D_{P^3,P}|=
p^{(3-1)1}=p^2$ lo cual implica que $D_{P^3,P}$ no es c\'iclico
y por tanto $G_{P^n}$ no es c\'iclico para $n\geq 3$ y $p\geq 3$.

Cuando $p=2$ y $d=1$, se tiene que $G_{P^2}\cong D_{P^2,P}
\cong R_T/\langle P\rangle\cong C_2$. Para $n=3$, $G_{P^3}
\cong D_{P^3,P}$ es c\'iclico pues tiene cardinalidad $4$ y $(1
+P)^2=1+P^2\not\equiv1\bmod P^3$, es decir, $D_{P^3,P}=\langle (1+P)
\bmod P^3\rangle$. Esto es, como grupos $D_{P^3,P}\cong C_4$.
Finalmente, si $n\geq 4$ el grupo $D_{P^n,P}$
no es c\'iclico, pues si lo fuese, $D_{P^,P}$
ser\'ia c\'iclico; sin embargo,
para toda $B\in R_T$ tenemos que $(1+PB)^4=
1+P^4B^4\equiv 1\bmod P^4$, es decir, $D_{P^4,P}$ es de 
exponente $4$ y de orden $8$.

En resumen, para $p=2$, $G_{P^n}$ es c\'iclico para $n\geq 2$
solamente en el caso $q=p$, $d=1$ y $n\in\{2,3\}$, esto es
$P^n\in\{T^2,T^3,(T+1)^2,(T+1)^3\}$.

Ahora bien, para el caso $r\geq 2$, debemos tener $q=2$. Sea
$M=\poly$ con $G_M\cong G_{P_1^{\alpha_1}}\times
\cdots\times G_{P_r^{\alpha_r}}$ un grupo c\'iclico. En particular
cada $G_{P_i^{\alpha_i}}$ es c\'iclico. Por el caso $r=1$,
suponemos $P_1=T, P_2=T+1$, $\alpha_i\in\{0,2,3\}$, 
$i=1,2$ y $\alpha_j=
1$, $3\leq j\leq r$. Puesto que $G_{T^2}\times G_{(T+1)^2}\cong
C_2\times C_2$ no es c\'iclico, debemos tener $\alpha_1\alpha_2
=0$. Ahora bien $G_{P_j}\cong C_{2^{d_j}-1}$ con
$d_j=\deg P_j\geq 2$, $3\leq j\leq r$.

Entonces $G_M$ es c\'iciclo con estas condiciones, si y solamente
si los n\'umeros $2^{d_j}-1$ son primos relativos a pares, $3\leq j
\leq r$. Finalmente, se tiene que $\mcd(2^c-1,2^d-1)=1\iff
\mcd(c,d)=1$ donde $c,d\geq 2$. En efecto, si $x$ es un divisor
mayor a $1$ de $c$ y de $d$, entonces $2^x-1|2^c-1$ y $2^x-1|
2^d-1$ y $2^x-1\geq 3$. Rec\'iprocamente, sea $\mcd(c,d)=1$.
Suponiendo que $2^c-1$ y $2^d-1$ no son primos relativos, podemos
seleccionar el m\'inimo n\'umero 
natural $c\geq 2$ tal que existe un n\'umero natural $d\in {\ma N}$,
$d>c$ primo relativo a $c$ y donde $2^c-1$ y $2^d-1$ no son primos
relativos. Sea $x>1$ tal que $x|2^c-1$ y $x|2^d-1$. Entonces $x|
(2^d-1)-(2^c-1)=2^c(2^{d-c}-1)$. Como $x$ es impar, $x|2^{d-c}-1$.
Sea $s\geq 2$ tal que $c<d<sc$. Entonces $x|2^{d-(s-1)c}-1$ y
$x|2^c-1$ lo cual es absurdo pues $0<d-(s-1)c<c$.

En resumen, hemos obtenido que $G_M$ es c\'iclico con $M=\poly$ 
si y solamente si tenemos uno de los siguientes casos.
\las
\item $M=P\in R_T^+$, $q$ arbitrario y donde suponemos que si $q=2$,
entonces $P\notin\{T,T+1\}$.

\item $M=P^2$ con $P\in R_T^+$, $q=p$ arbitrario y $\deg P=1$.

\item $M=P^3$ con $P\in R_T^+$, $q=2$, $P\in\{T,T+1\}$.

\item $r\geq 2$, $q=2$, $P_1=T$, $P_2=T+1$, con $\alpha_i\in\{0,2,3\}$, $i=
1,2$, $\alpha_1\alpha_2=0$ y 
$\alpha_j=1$, $\deg P_j=d_j\geq 2$, $3\leq j\leq r$ con $\{d_3,\ldots,
d_r\}$ primos relativos a pares.

\end{list}

Este ejemplo es el correspondiente del caso de los campos de funciones
ciclot\'omicos, al caso de los campos ciclot\'omicos num\'ericos
(ver Corolario \ref{C1.2.1.4}).
\end{ejemplo}

\subsection{Ramificaci\'on de $\p$ en $\cicl M{}/K$}\label{S9.2.1.Nueva}

Ahora estudiamos la descomposici\'on del primo infinito
$\p$ en una extensi\'on ciclot\'omica $\cicl M{}/K$.

Usaremos el pol\'igono de Newton para $M$
(ver \cite[Subsection 12.4.1]{Vil2006}). Sea $\pL$ un
primo de $K(\Lambda_M)$ sobre $\p$. Sea $e=e_{K(
\Lambda_M)/K}(\pL|\p)$. Se tiene $u^M=\sum\limits_{i=0}^d
\carlitzbinom{M}{i} u^{q^i}$ con $d=\deg M$ y $\carlitzbinom
{M}{i}\in R_T$ de grado $q^i(d-i)$, $\carlitzbinom{M}{0}
=M, \carlitzbinom{M}{d}=a_d$ donde $a_d$ es el coeficiente
l\'ider de $M$. Supondremos $a_d=1$. Analicemos las ra\'ices
de $u^M$ en $\bar{K}_{\infty}$, donde $K_{\infty}$ es la
completaci\'on de $K$ en $\p$ y $\bar{K}_{\infty}$ es una
cerradura algebraica de $K_{\infty}$.

Primero probaremos el siguiente resultado.

\begin{teorema}\label{T6.2.30(1)}
Sea $u^M\in K[u]\subseteq K_{\infty}[u]$. Para cada
$1\leq i\leq d$, existen exactamente $q^i-q^{i-1}$
ra\'ices $\tilde{\lambda}$ de $u^M/u$ en $\bar{K}_{\infty}$
tales que $v_{\pL}(\tilde{\lambda})
=e\big(d-i-\frac{1}{q-1}\big)$ (aqu\'i
no estamos usando las valuaciones normalizadas sino
que nos basamos en las valuaciones en $K
(\Lambda_M)$. Eventualmente
veremos que $e=q-1$).
\end{teorema}

\begin{proof} Se tiene $v_{\p}\Big(\carlitzbinom{M}{i}\Big)
=-q^i(d-i)$ y $v_{\pL}\Big(\carlitzbinom{M}{i}\Big)=-eq^i(d
-i)$. Consideremos el pol\'igono de Newton de $\frac{u^{
M}}{u}=\sum\limits_{i=0}^d\carlitzbinom{M}{i}u^{q^i-1}$.
Escribimos $\frac{u^M}{u}=\sum\limits_{j=0}^{q^d-1}
f_j(T)u^j$ donde $f_j(T)\neq 0\iff$ existe $i$ tal que
$j=q^i-1$ y en este caso se tiene
\begin{gather*}
f_j(T)=f_{q^i-1}(T)=\carlitzbinom{M}{i}, \quad
v_{\p}(f_{q^i-1}(T))=-\deg_T f_{q^i-1}(T)=-q^i(d-i), \\
\text{y}\quad
v_{\pL}(f_{q^i-1}(T))=e(-q^i(d-i))=-eq^i(d-i).
\end{gather*}

Los v\'ertices a considerar para el pol\'igono de
Newton son
\[
\{\beta_i\}_{i=0}^d=\{(q^i-1,\underbracket[0pt]{-eq^i(d-i))}_{
\substack{\uigual\\ v_{\pL}(f_{q^i-1}(T))}}\}_{0\leq i\leq d}.
\]

La pendiente entre $\beta_i$ y $\beta_{i+1}$, $0\leq i\leq d-1$, es:
\begin{align*}
\xi_i&=\frac{v_{\pL}(f_{q^{i+1}-1}(T))-v_{\pL}(f_{q^i-1}(T))}
{(q^{i+1}-1)-(q^i-1)}=\frac{v_{\pL}\Big(\carlitzbinom{M}{i+1}
\Big)-v_{\pL}\Big(\carlitzbinom{M}{i}\Big)}{q^i(q-1)}\\
&=\frac{-eq^{i+1}(d-(i+1))-(-eq^i(d-i))}{q^i(q-1)}=
\frac{(q^i-q^{i+1})e(d-i)+eq^{i+1}}{q^i(q-1)}\\
&=-e(d-i)+\frac{eq}{q-1}=ei-ed+\frac{eq}{q-1}<
e(i+1)-ed+\frac{eq}{q-1}=\xi_{i+1},
\end{align*}
esto es, $\xi_i<\xi_{i+1}$, $i=0,1,\ldots, d-1$.

Por tanto $\beta_0,\beta_1,\ldots,\beta_d$ son los v\'ertices
del pol\'igono de Newton pues $\xi_0<\xi_1<\cdots<\xi_{d-1}$.
Se sigue, para $0\leq i\leq d-1$, que hay $(q^{i+1}-1)-(q^i-1)$
ra\'ices de valuaci\'on $-\xi_i=e\big(d-i-\frac{q}{q-1}\big)$.
Equivalentemente, hay $q^{i+1}-q^i$ ra\'ices en $\bar{K}_{
\infty}$ de valuaci\'on (no normalizada)  $-\xi_i=e\big(d-(i+1)-
\frac{1}{q-1}\big)$.

En resumen hay $q^i-q^{i-1}$ ra\'ices $\tilde{\lambda}$ con
valuaci\'on $v_{\pL}(\tilde{\lambda})=e\big(d-i-\frac{1}{q-1}\big)$,
$1\leq i\leq d$. $\fin$
\end{proof}

El Teorema \ref{T6.2.30(1)} nos da en el caso $i=1$ que
hay $q-1$ ra\'ices de valuaci\'on $v_{\pK}(\tilde{\lambda})
=e\big(d-1-\frac{1}{q-1}\big)$. La importancia de estas
ra\'ices nos lo da el siguiente resultado.

\begin{proposicion}\label{P6.2.30(2)}
Sea $\tilde{\lambda}$ cualquiera de las $q-1$ ra\'ices con
valuaci\'on $v_{\pL}(\tilde{\lambda})=e\big(d-1-\frac{1}{q-1}\big)$.
Entonces $\lambda=\tilde{\lambda}$ es generador de
$\Lambda_M$.
\end{proposicion}

\begin{proof}
Sea $N\in R_T$ tal que $\deg N=s<d=\deg M$. Entonces
tenemos que
\begin{gather*}
\lambda^N=\sum_{i=0}^s\carlitzbinom{N}{i}\lambda^{q^i},\quad
\deg_T\carlitzbinom{N}{i}=q^i(s-i),\\
\begin{align*}
v_{\pL}\Big(\carlitzbinom{N}{i}\lambda^{q^i}
\Big)&=e\big(-q^i(s-i)+q^i\big(
d-1-\frac{1}{q-1}\big)\big)\\
&=eq^i\big(-s+i+d-1-\frac{1}{q-1}\big)\\
&=eq^i\big((d-s)+(i-1)-\frac{1}{q-1}\big),\\
v_{\pL}\Big(\carlitzbinom{N}{0}\lambda^{q^0}\Big)&
=v_{\pL}(N\lambda)=e\big(-s+\big(d-1-\frac{1}{q-1}\big)\big)\\
&=e\big(d-s-1-\frac{1}{q-1}\big).
\end{align*}
\end{gather*}

Si $i\geq 1$, puesto que $d-s\geq 1$, $d-s+i-1-\frac{1}{q-1}
> 0$, se tiene
\begin{gather*}
\begin{align*}
eq^i\big((d-s)+(i-1)-\frac{1}{q-1}\big)&>eq^0\big(d-s+0-1-
\frac{1}{q-1}\big)\\
&=e\big(d-s-1-\frac{1}{q-1}\big),
\end{align*}
\intertext{por lo tanto}
v_{\pL}\Big(\carlitzbinom{N}{i}\lambda^{q^i}\Big)>
v_{\pL}\Big(\carlitzbinom{N}{0}\lambda^{q^0}\Big)\quad
\text{para toda $i\geq 1$, lo cual implica}\\
v_{\pL}({\lambda}^N)=v_{\pL}\Big(\carlitzbinom{N}{0}
\lambda^{q^0}\Big)=v_{\pL}(N\lambda)=e\big(d-s-1-\frac{1}{q-1}
\big)\neq \infty
\end{gather*}
de donde se sigue que ${\lambda}^N\neq 0$ para toda
$N\in R_T$ tal que $\deg N<\deg M$. Por lo tanto 
$\lambda$ es generador de $\Lambda_M$.
$\fin$
\end{proof}

\begin{proposicion}\label{P6.2.30(3)} Las $q-1$ ra\'ices
$\tilde{\lambda}$ tales que $v_{\pL}(\tilde{\lambda})=
e\big(d-1-\frac{1}{q-1}\big)$, satisfacen que
$\tilde{\lambda}^{q-1}\in K_{\infty}$, donde $K_{\infty}$
denota la completaci\'on de $K$ en $\p$.
\end{proposicion}

\begin{proof}
Tenemos 
\[
\frac{u^M}{u}=\sum_{i=0}^d\carlitzbinom{M}{i}u^{q^i-1}=
\sum_{i=0}^d\carlitzbinom{M}{i}\big(u^{q-1}\big)^{\frac{q^i-1}{q-1}}
=f(u^{q-1})
\]
con $f(u)=\sum\limits_{i=0}^d\carlitzbinom{M}{i}
u^{\frac{q^i-1}{q-1}}$.

Se tiene que $\mu$ es ra\'iz de $u^M/u$ si y solamente si
$f(\mu^{q-1})=0$, esto es, $\mu^{q-1}$ es ra\'iz de $f$.
Puesto que si $\alpha\in\*\F$ se tiene
\[
(\alpha\mu)^{q-1}=\alpha^{q-1}\mu^{q-1}=\mu^{q-1},
\]
entonces $(\alpha\mu)^{q-1}$ es ra\'iz de $f$.
Por otro lado, si $\mu$ es ra\'iz de $u^M/u$, $\alpha
\mu$ tambi\'en es ra\'iz de $u^M/u$ para $\alpha\in\*\F$
pues $\carlitzbinom{M}{i}(\alpha\mu)^{q^i}=\alpha
\carlitzbinom{M}{i}\mu^{q^i}$, de donde $(\alpha\mu)^M/
(\alpha\mu)=(\alpha\mu^M)/(\alpha\mu)=\mu^M/\mu=0$.

Lo anterior prueba que el mapeo $\varphi\colon
\{\text{ra\'ices de $u^M/u$}\}\lra \{\text{ra\'ices de $f$}\}$,
$\mu\longmapsto \mu^{q-1}$ es suprayectivo y de hecho
es $(q-1)$ a $1$. Puesto que $u^M/u$ tiene $(q-1)$ ra\'ices
$\tilde{\lambda}$ tales que $v_{\pL}(\tilde{\lambda})=
e\big(d-1-\frac{1}{q-1}\big)$, y $v_{\pL}(\alpha\tilde{
\lambda})=v_{\pL}(\tilde{\lambda})$ para toda 
$\alpha\in \*\F$, $f$ tiene una ra\'iz $\mu$
de valuaci\'on $v_{\pL}(\mu)=v_{\pL}(\tilde{\lambda}^{q-1})=
(q-1)v_{\pL}(\tilde{\lambda})=e((q-1)(d-1)-1)$. Del
pol\'igono de Newton, se sigue que $\mu\in K_{\infty}$.
$\fin$
\end{proof}

En esta situaci\'on, tenemos que $K(\Lambda_M)=
K(\lambda_M)$ donde $\lambda_M$ es un generador de
$\Lambda_M$ como antes, esto es, $\lambda_M=\tilde{
\lambda}$. Sea $K(\Lambda_M)^+:=K(\lambda_M^{q-1})$.
Se tiene que $\Lambda_M\setminus \{0\}$ es el conjunto
de las ra\'ices de $u^M/u$ y $\Lambda_M^{q-1}\setminus
\{0\}$ consta de las ra\'ices de $f$ donde $f(u^{q-1})=u^M/u$. 
Como una ra\'iz de $f$ est\'a en $K_{\infty}$ la completaci\'on
del campo $K(\Lambda_M^{q-1})$ en $\pK:=\pL\cap
K(\Lambda_M)^+$ es $K(\Lambda_M)^+_{\pK}=K_{\infty}$.
Por tanto
\[
e_{K(\Lambda_M)^+/K}(\pK|\p)=f_{K(\Lambda_M)^+/K}(\pK|\p)
=1.
\]

Por otro lado $[\cicl M{}:\cicl M{}^+]\leq q-1$ y puesto que
\begin{gather*}
v_{\pL}(\lambda_M) = e\big(q-1-
\frac{1}{q-1}\big)=e(q-1)-\frac{e}{q-1}\in{\ma Z}
\intertext{lo cual implica que}
q-1\mid e=e_{K(\Lambda_M)/\cicl M{}^+}(\pL|\p)
\intertext{por lo tanto}
e_{K(\Lambda_M)/\cicl M{}^+}(\pL|\pK)=q-1
\quad\text{y}\quad [\cicl M{}:\cicl M{}^+]=q-1.
\end{gather*}

Para terminar, veamos cu\'al es el grupo
de inercia de $\p$ en $\cicl M{}/K$. Sea 
$\alpha\in\*\F$ y $\sigma_{\alpha}\in G_M=
\Gal(\cicl M{}/K)$ entonces se tiene
 $\sigma_{\alpha}\lambda_M=
\lambda_M^{\alpha}=\alpha\lambda_M$.

Se sigue que $\sigma_{\alpha}(\lambda_M^{q-1})=
(\sigma_{\alpha}\lambda_M)^{q-1}=(\alpha 
\lambda_M)^{q-1}=\alpha^{q-1}\lambda_M^{q-1}
=\lambda_M^{q-1}$. Por tanto $\cicl M{}^+
\subseteq \cicl M{}^J$ donde $J:=\{\sigma_{\alpha}\mid
\alpha\in\*\F\}\subseteq G_M$.

Por otro lado, tenemos que $[\cicl M{}:
\cicl M{}^J]=|J|=q-1=|\*\F|=[\cicl M{}:\cicl M{}^+]$
lo cual implica que $\cicl M{}^J=\cicl M{}^+$.

En resumen hemos probado el siguiente teorema.

\begin{teorema}\label{T6.2.30} Sea $M\in R_T$ un polinomio no
constante. Sean $G_M=\Gal(\cicl M{}/K)\cong \G M$ y
$J\cong \*\F\subseteq G_M$. Sea $\cicl M{}^+:=
\cicl M{}^J=K(\lambda_M^{q-1})$ donde $\lambda_M$
es un generador de $\Lambda_M$.

Entonces $\p$ se descompone totalmente en $\cicl M{}^+/K$
y si $\pL$ es un primo en $\cicl M{}$ sobre $\p$, $\pK=\pL\cap
\cicl M{}^+$, entonces $\pL/\pK$ es totalmente ramificado.
Esto es,
\begin{gather*}
e_{\cicl M{}^+/K}(\pK|\p)=f_{\cicl M{}^+/K}(\pK|\p)=1,\\
e_{\cicl M{}^/\cicl M{}^+}(\pL|\pK)=q-1,\quad 
f_{\cicl M{}^/\cicl M{}^+}(\pL|\pK)=1.
\end{gather*}

En particular $e_{\infty}=e_{\cicl M{}^/K}(\pL|\p)=q-1$, 
$f_{\infty}=f_{\cicl M{}^/K}(\pL|\pK)=1$.

Adem\'as $J=\{\sigma_{\alpha}\in G_M\mid \alpha
\in\*\F\}\cong \*\F$ es a la vez 
el grupo de descomposici\'on y el grupo de inercia de
$\pL/\p$:
\[
D_{\cicl M{}^/K}(\pL|\p)= I_{\cicl M{}^/K}(\pL|\p)\cong
\*\F.
\]

Por \'ultimo, $\cicl M{}^+=K(\lambda_M^{q-1})$.
\end{teorema}

\begin{proof}
El resultado lo hemos probado si $\lambda_M$ es un generador
de $\Lambda_M$ dado por la Proposici\'on \ref{P6.2.30(2)}.
Ahora $\Gal(\cicl M{}^+/K)\cong G_M/\*\F$. Los conjugados
de $\lambda_M^{q-1}$ son $\bar{\sigma}(\lambda_M^{q-1})$,
donde $\bar{\sigma}=\sigma\bmod \*\F$. Sea $\sigma(\lambda_M)
=\lambda_M^A$ con $A\in R_T$ y $\mcd(A,M)=1$. Se tiene
$\sigma(\lambda_M^{q-1})=\sigma(\lambda_M)^{q-1}$.
Por tanto, los conjugados de $\lambda_M^{q-1}$ son
$\{(\lambda_M^A)^{q-1}\}_{A\in \G M}$. Por tanto
$\cicl M{}^+=K(\lambda_M^{q-1})=K((\lambda_M^A)^{q-1})$.
Se sigue que $\cicl M{}^+=K(\lambda^{q-1})$ para cualquier
generador $\lambda$ de $\Lambda_M$. $\fin$
\end{proof}

\begin{definicion}\label{D6.2.30(5)} El campo
$\cicl M{}^+=K(\lambda_M^{q-1})$ se llama el {\em
subcampo real\index{subcampo real}} o el {\em m\'aximo
subcampo real} de $\cicl M{}$.
\end{definicion}

\begin{observacion}\label{O6.2.30(6)}
Lo anterior es an\'alogo al caso de campos num\'ericos. Ah\'i
tenemos
\[
\xymatrix{
\cic n{}\ar@{-}[d]\ar@{-}@/_1pc/[d]_{\{1,J\}}
\ar@{-}@/^1pc/[d]^{\substack{\text{$\infty$ es ramificado con}\\
\text{\'indice de ramificaci\'on $2$}}}\\
{\ma Q}(\zeta_n+\zeta_n^{-1})=\cic n{}^+=\cic n{}\cap {\ma R}
\ar@{-}[d]\ar@{-}@/^1pc/[d]^{\substack{\text{$\infty$ es totalmente}\\
\text{descompuesto}}}\\
{\ma Q}}\qquad
\xymatrix{
\cicl M{}\ar@{-}[d]\ar@{-}@/_1pc/[d]_{\*\F}\ar@{-}@/^1pc/[d]^{
\substack{\text{$\infty$ es totalmente}\\ \text{ramificado}}}\\
\cicl M{}^+=K(\lambda_M^{q-1})\ar@{-}[d]\ar@{-}@/^1pc/[d]^{\substack{\text{$\infty$ es totalmente}
\\ \text{descompuesto}}}\\K}
\]
$\{1,J\}=\{1,-1\}\longleftrightarrow \*\F=\F\setminus\{0\}$.

$\*\F$ juega el papel de $\{\pm 1\}$ (debido a esto, en los campos de funciones
se llaman ``{\em pares}'' a los enteros $a\in{\ma N}$ tales que
$(q-1)\mid a$).

\end{observacion}

\begin{observacion}\label{O6.2.31}
Notemos que $\cicl M{}=K$ si y solamente si $q=2$ 
y $M=T$, $M=T+1$ o $M=T(T+1)$. Es decir, ${\ma F}_2(T)
(\Lambda_T)$, ${\ma F}_2(T)(\Lambda_{T+1})$ y ${\ma F}_2(T)
(\Lambda_{T(T+1)})$ juegan el papel de $\cic 2{}={\ma Q}$
en el caso de los campos num\'ericos.
Debemos tener en cuenta siempre esta excepci\'on
para todo el desarrollo de los campos de funciones
ciclot\'omicos.

En particular, para $q=2$, $\p$ no es ramificado en ning\'un
$\cicl M{}/K$ y $\cicl {MT}{}=\cicl M{}$ para 
para todo $M\in R_T$ con $T\nmid M$. Similarmente
$\cicl {M(T+1)}{}=\cicl M{}$ para todo $M\in R_T$ tal que
$(T+1)\nmid M$ y finalmente tambi\'en se tiene que 
$\cicl {MT(T+1)}{}=\cicl M{}$ para toda $M\in R_T$ tal que
$\mcd (T(T+1),M)=1$.
\end{observacion}

\section{Ramificaci\'on en $\cicl M{}/K$}\label{S6.3}

Aqu{\'\i} estamos considerando $\cicl M{}/K$ donde $M\in R_T
\setminus\{0\}$ es un polinomio m\'onico.

\begin{proposicion}\label{P6.3.1} Se tiene que los primos ramificados
en $\cicl M{}/K$ son $\pK_{\infty}$ y los divisores de $M$ con la
excepci\'on de que si $q=2$, entonces $\p$ es no ramificado y
teniendo en cuenta la Observaci\'on {\rm{\ref{O6.2.31}}}.
\end{proposicion}

\begin{proof}
Se sigue de que $\cicl M{}=\prod_{P|M}\cicl p{\alpha}$, la Proposici\'on
\ref{P6.2.24} y el Teorema \ref{T6.2.30}. $\fin$
\end{proof}

Dado un lugar $\pK$ en $K$ con $\pK\neq \p$ 
y $\pL\in \cicl M{}$ sobre $\pK$, si 
$D$ e $I$ son los grupos de descomposici\'on e inercia, entonces
$\Gal(k(\pL)/k(\pK))\cong D/I$ (ver despu\'es de la Definici\'on
\ref{D5.4.6}). Si $\pK$ es no ramificado, $I=\{1\}$ y $D\cong \Gal(
k(\pL)/k(\pK))$. Puesto que $k(\pL)$ y $k(\pK)$ son campos finitos,
$\Gal(k(\pL)/k(\pK))$ es un grupo c{\'\i}clico generado por el 
automorfismo de Frobenius: 
\[
\sigma_{\pK}\colon k(\pL)\longto k(\pL),\quad \sigma_{\pK}(x)=
x^{|k(\pK)|}=x^{N(\pK)}
\]
donde denotamos $N(\pK)=|k(\pK)|=\big|{\cal O}_{\pK}/\pK\big|$.
As{\'\i} tenemos:

\begin{proposicion}\label{P6.3.2}
El automorfismo de Frobenius, el cual ser\'a denotado por
$\frobeniusbinom {\cicl M{}/K}
{\pL}$, est\'a caracterizado por la propiedad
\[
\frobeniusbinom{\cicl M{}/K}{\pL}(x)\equiv x^{N(\pK)}\bmod \pL 
\ \forall\ x\in
{\cal O}_{\pL}. \tag*{$\fin$}
\]
\end{proposicion}

Como $\cicl M{}/K$ es abeliana, $\frobeniusbinom{\cicl M{}/K}{\pL}$
es independiente de $\pL$ y s\'olo depende de $\pK$ y lo denotamos
$\varphi_P=\xbinom{\cicl M{}/K}{P}$ y se llama el {\em s{\'\i}mbolo de 
Artin\index{simbolo de Artin@s\'imbolo de Artin}\index{Artin!s{\'\i}mbolo de $\sim$}}.
Aqu\'i $P\in R_T^+$.

\begin{teorema}\label{T6.3.3}
Sea $P$ un polinomio irreducible que no divide a $M$. Entonces
el mapeo:
\begin{eqnarray*}
\varphi_P\colon \Lam M{}&\longto & \Lam M{}\\
\lambda&\longmapsto & \lambda^P
\end{eqnarray*}
corresponde al s{\'\i}mbolo de Artin $\xbinom{\cicl M{}/K}{P}
$.
\end{teorema}

\begin{proof}
Sea $\big(R_T\big)_P=\big\{\frac{f}{g}\mid f,g\in R_T, P\nmid g\big\}$ y sea
$( P)_K:=\frac{\pK}{\pK_{\infty}^{\gr P}}$, es decir, $P$ es un
elemento primo de $\pK$. Entonces $k(\pK)=\big(R_T
\big)_P/P\big(R_T\big)_P\cong R_T/\langle P\rangle\cong {\ma F}_{
q^d}$ donde $d=\gr P$.

Sea $\pL$ un lugar en $\cicl M{}$ sobre $\pK$. Entonces $N(\pK)=
|{\ma F}_{q^d}|=q^d$, $\cicl M{}\subseteq {\cal O}_{\pL}$. Entonces
\[
\xbinom{\cicl M{}/K}{P}(\lambda)\equiv \lambda^{q^d}
\bmod \pL.
\]

Ahora $u^P=u\Psi_P(u)=u\big(u^{q^d-1}+\beta_{q^d-2}u^{q^d-2}+
\cdots+\beta_1u+\beta_0\big)$.  Puesto que $\Psi_P(u)=\prod_{
\mcd (A,P)=1}(u-\lambda^A)$, $\lambda$ generador de $\Lam P{}$ y
$\Psi_P(0)=\pm \prod_{\mcd (A,P)=1}\lambda^A=P$
por (\ref{Ec6.2.26}), se tiene que $P|\beta_i$, $0\leq i\leq q^d-2$, lo
cual implica que $\lambda^P\equiv\lambda^{q^d}\bmod \pL$.

De la expresi\'on $u^M=\prod\limits_{A\bmod M}(u-\lambda^A)$, 
tomando derivadas con respecto a $u$, se tiene que $M=\sum\limits_{
A\bmod M}\Big(\prod\limits_{\substack{B\neq A\\ B\bmod M}} (u-
\lambda^B)\Big)$ que es constante en $u$. Sea $u=\lambda^C$.
Entonces $M=\prod\limits_{C\neq B}(\lambda^C-\lambda^B)$ y 
puesto que $P\nmid M$ se sigue que $\lambda^C\not\equiv
\lambda^B\bmod \pL$ para $C\not\equiv \bmod M$. En particular
$\lambda^P\equiv\lambda^Q\bmod \pL$ implica $\lambda^P=\lambda^Q$.

Finalmente $\lambda^P\equiv \xbinom{\cicl M{}/K}{P}\lambda\equiv
\lambda^{q^d}\bmod \pL$, de donde se sigue que $\varphi_P=
\xbinom{\cicl M{}/K}{P}$. $\fin$
\end{proof}

Con la notaci\'on usual de $e_P=$ ramificaci\'on de $P$, $f_P=$
grado de inercia y $h_P=$ n\'umero de primos encima de $P$,
tenemos:

\begin{proposicion}\label{P6.3.4}
Sea $M\in R_T\setminus\{0\}$ y sea $P$ un polinomio irreducible
que no divide a $M$. En $\cicl M{}/K$ tenemos 
\[
e_P=1,\quad f_P=o(P\bmod M),\quad h_P=\Phi(M)/f_P.
\]
\end{proposicion}

\begin{proof}
Sea $\lambda$ generador de $\Lam M{}$, $\cicl M{}=K(\lambda)$.
Sea $\pL$ un divisor primo en $\cicl M{}$ dividiendo a $\pK$ donde
$(P)_K=\frac{\pK}{\pK_{\infty}^{\gr P}}$. Entonces
\begin{align*}
{\cal O}_{\pL}&=\{\xi\in\cicl M{}\mid v_{\pL}(\xi)\geq 0\} \quad \text{y}\\
f_P&=\big[{\cal O}_{\pL}/\pL:(R_T)_P/P(R_T)_P\big]=
\big[({\cal O}_{M})_{\pL}/\pL({\cal O}_{M})_{\pL}:R_T/\langle 
P\rangle\big]\\
&=\big[{\cal O}_M/\pL{\cal O}_M:R_T/
\langle P\rangle\big]
\end{align*}
donde ${\cal O}_M$ es la cerradura entera de $R_T$ en $\cicl M{}$.

Sea $d=\gr P$. Puesto que $P\nmid M$, $\pK$ no es ramificado
en $\cicl M{}/K$ y el s{\'\i}mbolo de Artin $\varphi_P=\xbinom{\cicl
M{}/K}{P}$ en $P$ est\'a dado por $\varphi_P(\lambda)=\lambda^P$.
Entonces $e_P=1$ y $h_P=[\cicl M{}:K]/f_P=\Phi(M)/f_P$.

Finalmente, $f_P$ es el orden de $\varphi_P$, es decir $f_P$ es el
m{\'\i}nimo n\'umero natural tal que $\varphi_P^{f_P}=\Id\in G_M=
\Gal(\cicl M{}/K)$. Se tiene $\varphi_P^f=1\iff \varphi_P^f(\lambda)=
\lambda^{P^f}=\lambda\iff \lambda^{P^f-1}=0\iff M|P^f -1$.

Se sigue que $f_P=o(P\bmod M)$ es el m{\'\i}nimo n\'umero natural
tal que $M|P^{f_P}-1$. $\fin$
\end{proof}

El resultado general sobre ramificaci\'on est\'a dado por:

\begin{teorema}\label{T6.3.5} Sea $M=P_1^{\alpha_1}\cdots P_r^{
\alpha_r}\in R_T$ donde $P_1,\ldots,P_r$ son polinomios
irreducibles y sea $\cicl M{}/K$. Si $P\in R_T$ es distinto a $P_1,
\ldots, P_r$ y $\pK_{\infty}$, entonces
\begin{gather*}
e_P=1,\quad f_P=o(P\bmod M)\quad \text{y}\quad h_P=\Phi(M)/f_P.
\intertext{Si $P=P_i$ para alg\'un $1\leq i\leq r$, se tiene}
e_P=\Phi(P_i^{\alpha_i}),\quad f_P=o\big(P_i\bmod (M/P_i^{\alpha_i})
\big)\quad \text{y}\\
h_p=\frac{\Phi(M)}{\Phi(P_i^{\alpha_i}) f_{P_i}}=\frac{\Phi(M/P_i^{
\alpha_i})}{o\big(P_i\bmod (M/P_i^{\alpha_i})\big)}.
\end{gather*}

Finalmente, para $\pK_{\infty}$ se tiene:
\[
e_{\infty}=q-1,\quad f_{\infty}=1 \quad\text{y}\quad 
h_{\infty}=\Phi(M)/(q-1).
\]
\end{teorema}

\begin{proof}
El resultado se sigue de las Proposiciones \ref{P6.2.24}, \ref{P6.3.4}
y el Teorema \ref{T6.2.30}. $\fin$
\end{proof}

Tambi\'en tenemos el resultado an\'alogo al caso num\'erico que
establece ${\cal O}_{\cic n{}}={\ma Z}[\zeta_n]$ (Teorema 
\ref{T1.2.1.11}).

\begin{proposicion}\label{P6.3.6}
Sea $M=P^n$, $P\in R_T$ irreducible. Entonces ${\cal O}_M=R_T[
\lam M{}]$ donde ${\cal O}_M$ es la cerradura entera de $R_T$ en
$\cicl M{}$ y $\lam M{}$ es un generador de $\Lam M{}$.
\end{proposicion}

\begin{proof}
Sea $\lambda =\lam M{}$. Se tiene $R_T[\lambda]\subseteq 
{\cal O}_M$. Ahora sea $\alpha\in {\cal O}_M$. Puesto que
$\{1,\lambda,\ldots, \lambda^{\Phi(M)-1}\}$ es base de $\cicl M{}/K$,
existen $a_0,a_1,\ldots, a_r\in K$ tales que $\alpha=a_0+a_1 \lambda
+\cdots+a_r\lambda^r$ con $r=\Phi(M)-1$. Para probar el 
resultado, queremos probar que $a_i\in R_T$ para $i=0,1,\ldots, r$. 
Por el Corolario \ref{C6.2.26'} se tiene
que $v_{\pL}(\lambda)=1$ donde $\pL$ es el \'unico divisor en
$\cicl M{}$ sobre $\pK$, $(P)_K=\frac{\pK}{\pK_{\infty}^{\gr P}}$.

Se tiene que si $a_i\neq 0$, $v_{\pL}(a_i\lambda^i)=i+\Phi(P^n)
v_{\pK}(a_i)\equiv i\bmod \Phi(M)$ y por tanto $i\neq j$, $a_i\neq
0\neq a_J$, $v_{\pL}(a_i\lambda^i)\neq v_{\pL}(a_j\lambda^j)$.
Por tanto
\[
0\leq v_{\pL}(\alpha)=\min_{a_i\neq 0}\big\{v_{\pL}(a_i\lambda^i)
\big\}=\min_{a_i\neq 0}\big\{i+\Phi(M)v_{\pL}(a_i)\big\}.
\]

En particular $v_{\pL}(a_i)\geq 0$ para toda $i=0,1,\ldots, r$. Para
cualquier $\sigma_A\in G_{P^n}=\Gal(\cicl Pn/K)$, $\sigma_A(
\lambda)=\lambda^A$, se tiene
\[
\alpha_A:=\sigma_A \alpha=a_0+a_1\lambda^A+\cdots+
a_r(\lambda^A)^r
\]
$A\bmod P^n\in \big(R_T/\langle P^N \rangle\big)^{\ast}$. Si
$\big\{\overline{A}_1,\ldots, \overline{A}_{\Phi(M)}\big\}$ es un
conjunto de representantes de $\big(R_T/\langle P^n\rangle\big)^{
\ast}$, poniendo $\alpha_i:=\alpha^{A_i}$ y $\lambda_i=\lambda^{
A_i}$ se obtiene
\[
\left(\begin{array}{c}
\alpha_1\\ \vdots \\ \alpha_{\Phi(M)}\end{array}\right)=
\left(\begin{array}{ccccc}1& \lambda_1&\lambda_1^2&\cdots&
\lambda_1^r\\ \vdots&\vdots&\vdots&&\vdots\\
1&\lambda_{r+1}&\lambda_{r+1}^2&\cdots & \lambda_{r+1}^r
\end{array}\right)
\left(\begin{array}{c} a_0\\ \vdots\\ a_r\end{array}\right).
\]

El determinante de la matriz $\big[\lambda_i^j\big]_{\substack{0\leq
j\leq r\\ 1\leq i\leq r+1}}$ es un determinante de Vandermonde por lo
que $\det\big[\lambda_i^j\big]=\prod\limits_{1\leq t\leq \ell\leq r+1}
(\lambda_{\ell}-\lambda_t):=d$. Por tanto
\[
a_i=\frac{\det \left[\begin{array}{cccccccc}
1&\lambda_1&\cdots&\lambda_1^{i-1}&\alpha_1&\lambda_1^{i+1}&
\cdots&\lambda_1^r \\ \vdots&\vdots&
&\vdots&\vdots&\vdots&&\vdots \\
1&\lambda_{r+1}&\cdots&\lambda_{r+1}^{i-1}&\alpha_1&
\lambda_{r+1}^{i+1}&\cdots&\lambda_{r+1}^r\end{array}
\right]}{\det
\left[\begin{array}{ccccc}1& \lambda_1&\lambda_1^2&\cdots&
\lambda_1^r\\ \vdots&\vdots&\vdots&&\vdots\\
1&\lambda_{r+1}&\lambda_{r+1}^2&\cdots & \lambda_{r+1}^r
\end{array}\right]}=
\frac{b_i}{d}
\]
con $b_i\in {\cal O}_M$. Por la demostraci\'on de la Proposici\'on
\ref{P6.2.24} $\lambda =\beta_A\lambda^A$ para $A\bmod P^n\in
\big(R_T/\langle P^n\rangle\big)^{\ast}$ y $P=\beta_0\lambda^{\Phi
(P^n)}$, $\beta_n-\beta_0\in {\cal O}_M$.

Entonces para cualquier divisor primo ${\eu q}$ en $\cicl M{}$ que
no divide ni a $\pK$ ni a $\pK_{\infty}$, se tiene $v_{\pL}(\lambda)=
v_{\pL}(\lambda^A)=0$. Se sigue que el soporte del divisor de polos
de $a_i$ puede consistir solo de $\pK$ y $\pK_{\infty}$. Sin
embargo, puesto que $v_{\pK}(a_i)\geq 0$, se sigue que $a_i\in
R_T$ y que ${\cal O}_M=R_T[\lambda]$. $\fin$
\end{proof}

El resultado general, se sigue de la Proposici\'on \ref{P6.3.6} y del
an\'alogo del Teorema \ref{T1.2.1.10} (ver la demostraci\'on del Teorema
\ref{T1.2.1.11}).

\begin{teorema}\label{T6.3.7} Para cualquier $M\in R_T\setminus
\{0\}$, si ${\cal O}_M$ es la cerradura entera de $R_T$ en
 $\cicl M{}$ y $\lambda$ es un generador
de $\Lam M{}$, se tiene ${\cal O}_M=R_T[\lambda]$. $\fin$
\end{teorema}

El an\'alogo al Teorema \ref{T6.3.7} se sigue cumpliendo para
$\cicl M{}^+$.

\begin{teorema}\label{T6.3.7(1)} Para cualquier $M\in
R_T\setminus\{0\}$, se tiene ${\mc O}_{\cicl M{}^+}=
R_T[\lambda_M^{q-1}]$, donde $\lambda_M$ es cualquier
generador de $\Lambda_M$ tal que $\cicl M{}^+=
K(\lambda_M^{q-1})$.
\end{teorema}

\begin{proof} Se tiene $R_T[\lambda_M^{q-1}]\subseteq
{\mc O}_{\cicl M{}^+}$. 

Sea $\alpha\in\cicl M{}^+=K(\lambda_M^{q-1})$ entero sobre
$R_T$. Se tiene que $\{1,\beta,\ldots \beta^r\}$ es una base de la
extensi\'on de $\cicl M{}^+/K$, donde $\beta=\lambda_M^{q-1}$
y $r=\frac{\Phi(M)}{q-1}-1$. Por tanto
existen elementos $a_0,a_1,\ldots,a_r\in K$ tales que
$\alpha=a_0+a_1\beta+\cdots+a_r\beta^r$. Puesto que
$\beta^r=\big(\lambda_M^{q-1}\big)^{
\big(\frac{\Phi(M)}{q-1}-1\big)}=\lambda_M^{\Phi(M)-(q-1)}$
y $(q-1)r=\Phi(M)-(q-1)\leq \Phi(M)-1$, del hecho de que
$\alpha\in{\mc O}_{\cicl M{}^+}\subseteq {\mc O}_{\cicl M{}}$
se tiene que $a_0,\ldots,a_r\in R_T$ y por tanto
${\mc O}_{\cicl M{}^+}\subseteq R_T[\lambda_M^{q-1}]$.  $\fin$
\end{proof}

\begin{observacion}\label{O6.3.7(4)}
El Teorema \ref{T6.3.7(1)} se generaliza para cualquier subcampo
$\cicl M{}^+\subseteq F\subseteq \cicl M{}$. En efecto, puesto que
$\cicl M{}/\cicl M{}^+$ es una extensi\'on c\'iclica de Kummer de
grado $q-1$, se tiene que $F$ es de la forma $F=K(\lambda_M^m)$ con
$m\mid q-1$, $[\cicl M{}:F]=m$. Puesto que $\lambda_M^m$ es
entero, $R_T[\lambda_M^m]\subseteq {\mc O}_F$ y $\{1,\gamma,
\ldots, \gamma^r\}$ es una base de la extensi\'on $F/K$ donde
$\gamma=\lambda_M^m$, $r=\frac{\Phi(M)}{m}-1$. 

Sea $\alpha\in{\mc O}_F$. Entonces $\alpha\in{\mc O}_{\cicl M{}}=
R_T[\lambda_M]$ y $\alpha=\sum_{i=0}^r a_i\gamma^i$ con
$a_i\in K$. Ahora $\gamma^i=\lambda_M^{mi}$, $0\leq i\leq r$. Por
lo tanto $0\leq mi\leq mr=\Phi(M)-m\leq \Phi(M)-1$. Se sigue que
$a_0,\ldots,a_r\in R_T$. Por tanto $\alpha\in R_T[\lambda_M^m]$
y ${\mc O}_F\subseteq R_T[\lambda_M^m]$. Se sigue la
igualdad ${\mc O}_F=R_T[\lambda_M^m]$.
\end{observacion}

\begin{observacion}\label{O6.3.7(2)} Dado un campo $K\subseteq
E\subseteq \cicl M{}$ no necesariamente existe $\beta\in
{\mc O}_E$ tal que ${\mc O}_E=R_T[\beta]$.
\end{observacion}

\begin{ejemplo}\label{E6.3.7(3)} 
Consideremos dos polinomios $P, Q\in R_T$ polinomios m\'onicos
irreducibles distintos, de grados $d_0$ y $d$ respectivamente.
Digamos
\begin{align*}
P(T)&=T^{d_o}+a_{d_0-1}T^{d_0-1}+\cdots+a_1T+a_0,\\
Q(T)&=T^d+b_{d-1}T^{d-1}+\cdots+b_1T+b_0.
\end{align*}

Sea $f=o(P\bmod Q)$ tal que $\frac{q^d-1}{f}>q^{d_0}$. Sea $E$
el \'unico subcampo $K\subseteq E\subseteq \cicl Q{}$ tal que
$[E:K]=\frac{q^d-1}{f}$.
\[
\xymatrix{
\cicl Q{}\ar@{-}[d]_f\\ E\ar@{-}[d]_{\frac{q^d-1}{f}}\\ K}
\]

Entonces el grupo de descomposici\'on de $P$ en
la extensi\'on $\cicl Q{}/K$ es
$D=\Gal(\cicl Q{}/E)$ pues $e_P=1$, $f_P=f$, 
$h_P=\frac{[\cicl Q{}:K]}{e_Pf_P}
=\frac{q^d-1}{f}$. Entonces $P$ se descompone totalmente en $E/K$.

Supongamos que ${\mc O}_E=R_T[\delta]$ para alg\'un $\delta\in E$.
Sean $g(u)=\Irr(u,\delta, K)$ y $\overline{g(u)}=g(u)\bmod P\in \big(R_T/
P\big)[u]={\ma F}_{q^{d_0}}[u]$. Consideremos la factorizaci\'on
de $\overline{g(u)}$:
\begin{gather*}
\overline{g(u)}=\overline{G_1(u)}^{\alpha_1}\cdots \overline{G_h(u)}^{\alpha_h}.
\intertext{Por el Teorema de Kummer se tiene}
P{\mc O}_E=\pK_1^{\alpha_1}\cdots \pK_h^{\alpha_h}.
\intertext{Puesto que $P$ se descompone totalmente en $E/K$,
se sigue que}
P{\mc O}_E=\pK_1\cdots \pK_{[E:K]}.
\end{gather*}

El grado relativo $d_{E/K}(\pK_i|P)=1$, por lo que $\pK_i$ es de grado
$d_0$ sobre $\F$ o, lo que es lo mismo, $\pK_i$ es de grado $1$ sobre
${\ma F}_{q^{d_0}}$. En otras palabras, $G_1,\ldots, G_h$ son todos
polinomios lineales distintos en ${\ma F}_{q^{d_0}}[u]$ y $h=[E:K]=\frac{q^d-1}{f}
>q^{d_0}$ lo cual es absurdo pues \'unicamente existen $q^{d_0}$
polinomios lineales en ${\ma F}_{q^{d_0}}[T]$, a saber $\{T-\beta\}$,
$\beta\in {\ma F}_{q^{d_0}}$. Entonces ${\mc O}_E$ no puede ser
de la forma $R_T[\delta]$ con $\delta \in E$.

Para tener un ejemplo concreto, necesitamos hallar $q,d,d_0,Q,P,f$
con $\frac{q^d-1}{f}>q^{d_0}$.

Sean $q=2$, $P(T)=T^2+T+1$. Se tiene que $P$ es el \'unico polinomio
cuadr\'atico irreducible en ${\ma F}_2[T]$. Se tiene $d_0=2$. Sea $Q(T)=
T^4+T+1$. Entonces $Q$ es irreducible pues $Q(0)=Q(1)=1\neq 0$ por
lo que $Q(T)$ no tiene factores lineales y $P(T)^2=(T^2+T+1)^2=
T^4+T^2+1\neq Q(T)$ y como $P$ es el \'unico cuadr\'atico irreducible,
$Q$ no tiene factores cuadr\'aticos.

Ahora bien, se tiene $P(T)^3=T^6+T^5+T^3+T+1\equiv 1\bmod Q(T)$.
Se sigue que $f=3$, $d=4$, $d_0=2$ y $q=2$. En este caso tenemos
\[
\frac{q^d-1}{f}=\frac{2^4-1}{3}=\frac{16-1}{3}=5=2^2+1>4=q^{d_0}.
\]
Por tanto si $E$ es el subcampo tal que ${\ma F}_2(T)\subseteq E
\subseteq {\ma F}_2(T)(\Lambda_{T^4+T+1})$ con $[E:{\ma F}_2(T)]=
\frac{q^d-1}{f}=5$, ${\mc O}_E\neq R_T[\delta]$ para toda 
$\delta \in E$.
\end{ejemplo}

Se tiene $\big[\cicl M{}:\cicl M{}^+\big]=\big|{\ma F}_q^{\ast}\big|=
q-1$ y $\pK_{\infty}$ se descompone totalmente en $\Phi(M)/(q-1)$
divisores en $\cicl M{}^+/K$.

Como en el caso num\'erico (ver Teorema \ref{T10.2}),
el comportamiento de la extensi\'on
$\cicl M{}/\cicl M{}^{+}$ es diferente cuando $M$ es potencia
de un polinomio irreducible que cuando hay al menos dos primos
distintos dividiendo a $M$ como veremos a continuaci\'on.

Sea $M=\prod_{i=1}^r P_i^{\alpha_i}$ como producto de
polinomios irreducibles.
Sean $\p$ el primo infinito en $K$ y ${\eu P}_{\infty}$ un primo en
$\cicl M{}$ dividiendo a $\p$. Se tiene $e({\eu P}_{\infty}|\p)=q-1$
y $I({\eu P}_{\infty}|\p)={\ma F}_q^{\ast}$. Por lo tanto
$\p$ es no ramificado en $\cicl M{}^{+}/K$ y ${\eu P}_{\infty}$
es totalmente ramificado en $\cicl M{}/\cicl M{}^{+}$. Por otro
lado, si ${\eu Q}_{\infty}^{(i)}$ denota un primo en $\cicl {P_i}{\alpha_i}$
dividiendo a $\p$, se tiene $e({\eu Q}_{\infty}^{(i)}|\p)=q-1$.

\begin{lema}\label{L6.4.30(1)} Se tiene $\cicl M{}=\cicl M{}^{+}(
\Lambda_{P_i})$.
\end{lema}

\begin{proof} Sea $F:=\cicl M{}^{+}(\Lambda_{P_i})$, $F\subseteq
\cicl M{}$. Se tiene $[\cicl M{}:\cicl M{}^{+}]=q-1=|{\ma F}_q^{\ast}|$.
Con las notaciones naturales, tenemos:
\[
e_{F/K}(\p)=e_{F/\cicl M{}^{+}}(\p)
e_{\cicl M{}^{+}/K}(\p)=e_{F/\cicl M{}^{+}}(\p).
\]

Por otro lado, puesto que $K\subseteq \cicl {P_i}{}
\subseteq F$, se tiene $e_{F/K}(\p)\geq e_{\cicl {P_i}{}/K}(\p)
=q-1$. Por lo tanto
\begin{gather*}
\begin{align*}
e_{F/\cicl M{}^{+}}(\p)&\geq q-1=[\cicl M{}:\cicl M{}^{+}]
\\&\geq [F:\cicl M{}^{+}]\geq e_{F/\cicl M{}^{+}}(\p),
\end{align*}
\intertext{de donde se sigue que}
e_{F/\cicl M{}^{+}}(\p)=[F:\cicl M{}^{+}]=q-1=[\cicl M{}:
\cicl M{}^{+}].
\end{gather*}
Por lo tanto $F=\cicl M{}$. $\fin$
\end{proof}

\begin{corolario}\label{C6.4.30(2)}
En $\cicl M{}=\cicl M{}^{+}(\Lambda_{P_i})$, el \'unico posible
primo finito ramificado es $P_i$. $\fin$
\end{corolario}

\begin{corolario}\label{C6.4.30(3)}
Si $r\geq 2$, no hay ning\'un primo finito ramificado
en $\cicl M{}/\cicl M{}^+$. $\fin$
\end{corolario}

\begin{observacion}\label{O6.4.30(4)}
Si $r=1$, $M=P_1^{\alpha_1}$, $P_1$ es ramificado en
la extensi\'on
$\cicl M{}/\cicl M{}^+$, excepto en el caso $q=2$, $\alpha_1=1$
y $P_1=T$ o $P_1=T+1$.
\end{observacion}

\begin{observacion}\label{O6.4.30(5)}
En cualquier caso, excepto $q=2$, $\p$ es ramificado
en $\cicl M{}/\cicl M{}^+$.
\end{observacion}

\begin{definicion}\label{D6.4.30(6)}
Para $L\subseteq \cicl M{}$, se define el {\em subcampo
real\index{subcampo real}} de $L$ como $L^+:=\cicl M{}^+
\cap L$.
\end{definicion}

Tambi\'en se tiene que el comportamiento del generador $\lambda_M$
de $\Lambda_M$ (el cual equivale a $1-\zeta_n$ en el caso num\'erico),
es diferente cuando $M$ es potencia de un primo o cuando hay al
menos dos primos que dividen a $M$.

\begin{proposicion}\label{P6.4.30(7)} Sea $\lambda_M$ un generador
de $\Lambda_M$ como $R_T$--m\'odulo, $M\in R_T$ un polinomio
m\'onico no constante. Entonces
\las
\item Si $M$ es potencia de un primo, $\lambda_M$ es un elemento
primo de ${\mc P}$ donde ${\mc P}$ es el primo en $\cicl M{}$
dividiendo a $P$ y donde $M=P^n$.
\item Si $M$ es divisible por al menos dos polinomios irreducibles
distintos, entonces $\lambda_M$ es unidad de ${\mc O}_{\cicl M{}}$. 
\end{list}
\end{proposicion}

\begin{proof} (1) est\'a probado en la demostraci\'on de la Proposici\'on
\ref{P6.2.24} y Corolario \ref{C6.2.26'}
(ver en particular la Ecuaci\'on (\ref{Ec6.2.26})).

Para probar (2) se puede dar una prueba similar al caso
num\'erico. Alternativamente sea $M=P_1^{\alpha_1}\cdots P_r^{\alpha_r}$
con $r\geq 2$. Lo haremos por inducci\'on en $r$. Si $r=2$, $M=P_1^{\alpha_1}
P_2^{\alpha_2}$ y sin p\'erdida de generalidad podemos tomar los
generadores tales que $\lambda_{P_1^{\alpha_1}}=\lambda_M^{P_2^{\alpha_2}}$
y $\lambda_{P_2^{\alpha_2}}=\lambda_M^{P_1^{\alpha_1}}$.

Puesto que $\lambda_M^{P_2^{\alpha_2}}=\sum\limits_{j=0}^d \carlitzbinom{
P_2^{\alpha_2}}{j} \lambda_M^{q^j}=\lambda_M\cdot \xi$ con
$\xi\in{\mc O}_{\cicl M{}}$, se sigue que 
\[
\lambda_{P_1^{\alpha_1}}=\lambda_M\frac{\lambda_M^{P_2^{\alpha_2}}}{
\lambda_M}=\lambda_M\cdot \xi,
\]
por lo que $\lambda_M\mid \lambda_{P_1^{\alpha_1}}$ en ${\mc O}_{\cicl M{}}$.

Sea $N_M:=N_{\cicl M{}/K}$. Entonces se tiene $N_M\lambda_M\mid
N_{P_1^{\alpha_1}}(\lambda_{P_1^{\alpha_1}})^{a_1}$ y $N_M\lambda_M\mid
N_{P_2^{\alpha_2}}(\lambda_{P_2^{\alpha_2}})^{a_2}$ con
$a_i=[\cicl M{}:\cicl {P_i}{\alpha_i}]$, $i=1,2$.

Ahora bien $N_{P_i^{\alpha_i}}(\lambda_{P_i^{\alpha_i}})=\prod_{\mcd(A,P_i)=1}
\lambda_{P_i^{\alpha_i}}^A=\pm P_i$, $i=1,2$ por lo que $N_M\lambda_M\mid
P_i^{a_i}$, $i=1,2$. Se sigue que $N_M\lambda_M\mid\mcd(P_1^{a_1},
P_2^{a_2})=1$. Obtenemos que $\lambda_M$ es unidad.

Para $r\geq 3$, sean $M_1:=P_1^{\alpha_1}$ y $M_2=P_2^{\alpha_2}
\cdots P_r^{\alpha_r}$. Por el caso $r=2$ se tiene que $\lambda_{M_2}$
es unidad y $N_M\lambda_M\mid N_{M_2}(\lambda_{M_2})^a$ es unidad.
Se sigue que $\lambda_M$ es unidad. $\fin$
\end{proof}

\begin{proposicion}\label{P12*.2.2.G} Sea $K_{\infty}$ la completaci\'on
de $K$ en $\p$. Entonces 
\[
\bigcup_{M\in R_T} K(\Lambda_M)^+
\subseteq K_{\infty}.
\]
\end{proposicion}

\begin{proof} Se tiene que si $\pL$ es un divisor
primo en $K(\Lambda_M)^+$ sobre
$\p$, $K(\Lambda_M)^+\subseteq 
K(\Lambda_M)^+_{\pL}=K_{\infty}$
pues el grado relativo es igual $1$, es decir,
$d_{K(\Lambda_M)^+/K}
(\pL|\p)=1$. $\fin$
\end{proof}

\begin{proposicion}\label{P12*.2.2.G'} Se tiene
\[
\bigcup_{M\in R_T} K(\Lambda_M)
\subseteq K_{\infty}\big(\sqrt[q-1]{-1/T}\big).
\]
y esta es la m\'axima extensi\'on abeliana de $K$ dentro de 
$K_{\infty}\big(\sqrt[q-1]{-1/T}\big)$.
\end{proposicion}

\begin{proof} Ver Teorema \ref{P3.4}.
Sea $L=\bigcup_{M\in R_T} K(\Lambda_M)$. Entonces $L^+=
L^{G_0}\subseteq K_{\infty}$ donde $G_0$ es el grupo de inercia
de $\p$, $G_0\cong \*\F$.
Por la Proposici\'on \ref{P12*.2.2.G}, $L^+\subseteq K_{\infty}$. Ahora
bien, para toda $M\in R_T$, $M$ no constante,
se tiene que $L=L^+(\lambda_M)$,
donde $\lambda_M$ es un generador $\Lambda_M$ pues, por un lado
$L^+(\lambda_M)\subseteq L$ y por otro, el \'indice de ramificaci\'on
de $\p$ en $L^+(\lambda_M)/L^+$ es $q-1$ debido a que $K(\Lambda_M)
\subseteq L^+(\lambda_M)$ y $[L:L^+]=q-1$.

Por tanto, si tomamos $M=T$ se tiene que 
$\lambda=\lambda_T$ satisface
\[
\lambda^T=\lambda^q+T\lambda=0,
\]
por lo que $\lambda_T=
\sqrt[q-1]{-T}$. Se sigue que 
$L=L^+(\lambda_T)=L^+(1/\lambda_T)
=L^+\big(\sqrt[q-1]{-1/T}\big)\subseteq
K_{\infty}\big(\sqrt[q-1]{-1/T}\big)$. Se sigue que $K(\Lambda_M)
\subseteq K_{\infty}\big(\sqrt[q-1]{-1/T}\big)$ y $L\subseteq
K_{\infty}\big(\sqrt[q-1]{-1/T}\big)$.

Tambi\'en tenemos que la m\'axima extensi\'on abeliana de $K=
\F(T)$ dentro de $K_{\infty}\big(\sqrt[q-1]{-1/T}\big)$ es $L=
\bigcup_{M\in R_T}K(\Lambda_M)$. 

Alternativamente, si $S=1/T$, consideramos $\cicl S{}/K$ en
el cual $\p$ tiene \'indice de ramificaci\'on $q-1 ([\cicl S{}:K]
=q-1)$ y si $\lambda=\lambda_S$, entonces
$\lambda^S=\lambda_S^q+S\lambda_S=0$, por lo que
$\lambda_S=\sqrt[q-1]{-S}=\sqrt[q-1]{-1/T}$.
$\fin$
\end{proof}

Tambi\'en tenemos el resultado an\'alogo al Corolario \ref{C7.3}.

\begin{definicion}\label{D6.3.8} Sea $P\in R_T$ un polinomio m\'onico
e irreducible y sea $A\in R_T$. Decimos que
\[
\mu (A\bmod P)=M\in R_T
\]
si $M$ es m\'onico y de grado m{\'\i}nimo satisfaciendo $A^M
\equiv 0\bmod P$.
\end{definicion}

\begin{observacion}\label{O6.3.9} Sea $N\in R_T$ tal que
$A^N\equiv 0\bmod P$ y sea $N=QM+R$ con $Q, R\in R_T$ y
$R=0$ o $\gr R<\gr M$.

Entonces $A^N=(A^M)^Q+A^R$ por lo que $A^R\equiv 0\bmod P$.
Por tanto $R=0$ y $M$ divide a $N$. En particular el polinomio
$M$ dado en la Definici\'on \ref{D6.3.8} es \'unico.
\end{observacion}

Por otro, puesto que $R_T/\langle P\rangle$ es finito, el conjunto
$\{A^M\bmod P\mid M\in R_T\}$ es tambi\'en finito y por tanto 
existen dos elementos distintos $M_1,M_2\in R_T$ tales que
$A^{M_1}\equiv A^{M_2}\bmod P$. Por tanto $A^{M_1-M_2}\equiv
0\bmod P$ y $M_1-M_2\neq 0$.

\begin{proposicion}\label{P6.3.10} Sea $P\in R_T$ polinomio
irreducible y $M\in R_T$ m\'onico no divisible por $P$. Si $A\in R_T$,
entonces 
\[
P|\Psi_M(A) \iff \mu(A\bmod P)=M.
\]
\end{proposicion}

\begin{proof}
Primero supongamos que $P|\Psi_M(A)$. Puesto que 
\[
u^M=\prod_{
D|M}\Psi_D(u)\quad \text{se sigue que}\quad
A^M=\prod_{D|M}\Psi_D(A)\equiv
0\bmod P.
\]

Sea $\mu(A\bmod P)=N$. Entonces $N|M$ y por tanto $A^N=\prod_{
D|N}\Psi_D(A)\equiv 0\bmod P$. Por lo tanto existe $D_0|N$ tal que
$P|\Psi_{D_0}(A)$. Si $D_0\neq M$, entonces 
\[
A^M=\Psi_M(A)\Psi_{D_0}(A)\prod_{\substack{D\neq D_0, M\\
D|M}}\Psi_D(A)\equiv 0\bmod P^2.
\]

En particular $u^M\bmod P$ tiene una ra{\'\i}z m\'ultiple pero $(u^M)'=
M\not\equiv 0 \bmod P$, es decir, $u^M\bmod P$ es separable. Esta 
contradicci\'on prueba que $D_0=M$ y $\mu(A\bmod P)=M$.

Rec{\'\i}procamente, sea $\mu(A\bmod P)=M$. Entonces 
\[
A^M=
\prod_{D|M} \Psi_D(A)\equiv 0\bmod P.
\]
 Por tanto $P|\Psi_D(A)$
para alg\'un $D|M$. Si $D\neq M$, $A^D=\prod_{D'|D}\Psi_{D'}(A)
\equiv 0\bmod P$ lo que contradice el hecho de que $\mu(A\bmod P)
=M$. Por tanto $D=M$ y $P|\Psi_M(A)$. $\fin$
\end{proof}

\begin{proposicion}\label{P6.3.11} Sea $P\in R_T$ un polinomio
irreducible y $M\in R_T$ un polinomio m\'onico tal que $P\nmid M$.
Entonces $P$ divide a $\Psi_M(A)$ para alg\'un $A\in R_T$ si y
solamente si $P\equiv 1\bmod M$.
\end{proposicion}

\begin{proof}
Si $P|\Psi_M(A)$ para alg\'un $A\in R_T$, entonces por la 
Proposici\'on \ref{P6.3.10} se tiene que $\mu(A\bmod P)=M$. 
Ahora bien 
\[
u\Psi_P(u)=u^P=\sum_{i=0}^d \carlitzbinom Pi u^{q^i}\equiv
u^{q^d}\bmod P
\]
donde $d=\gr P$, pues $\Psi_P(u)=\prod\limits_{\mcd (C,P)=1}
(u-\lambda^C)$ y $\langle P\rangle =\langle \lambda\rangle^{\Phi(P)}$
(ver la prueba de la Proposici\'on \ref{P6.2.24}). Por tanto
$A^P\equiv A^{q^d}\bmod P$. 

Ahora bien, $\Phi(P)=q^d-1=\big|\big(R_T/\langle P\rangle\big)^{\ast}
\big|$ por lo que si $P\nmid A$ entonces $A^{q^d-1}\equiv 1\bmod P$
as{\'\i} que $A^{q^d}\equiv A\bmod P$. Si $P|A$, $A^{q^d}\equiv 0
\equiv A\bmod P$.

En cualquier caso tenemos $A^{q^d}\equiv A\bmod P$ y por tanto
$A^P\equiv A\bmod P$, lo cual equivale a $A^P-A=A^{P-1}\equiv 0
\bmod P$. Puesto que $\mu(A\bmod P)=M$ se sigue que
$M$ divide a $P-1$ lo cual implica que $P\equiv 1\bmod M$.

Rec{\'\i}procamente, supongamos ahora que $P\equiv 1\bmod M$.
Se tiene que $d=\gr (P-1)=\gr P$ y $u^{P-1}=\sum\limits_{i=0}^d
\carlitzbinom {P-1}i u^{q^i}$. Por lo tanto $\big(u^{P-1}\big)'\bmod P
\equiv (P-1)\bmod P\equiv -1\bmod P\not\equiv 0$ lo cual implica
que el polinomio $u^{P-1}\bmod P\in \big(R_T/\langle P\rangle\big)
[u]$ es separable.

Puesto que $\gr_u u^{P-1}=q^d=\big|R_T/\langle P\rangle\big|$ y
$A^{P-1}\equiv 0\bmod P$ para todo $A\in R_T$, se sigue que
\[
u^{P-1}\bmod P=\prod_{D|P-1}\Psi_D(u)\bmod P =\prod_{\substack{
A\bmod P\\ A\in R_T}} (u-A)\bmod P.
\]
Por lo tanto existe $A\in R_T$ tal que $\Psi_M(A)\equiv 0\bmod P$.
Se sigue que $P$ divide a $\Psi_M(A)$ y que $\mu(A\bmod P)
=M$. $\fin$
\end{proof}

\begin{corolario}[Caso particular al Teorema de
Dirichlet]\label{C6.3.12}
Sea $M\in R_T$ un polinomio m\'onico no constante. Entonces
existen una infinidad de polinomios irreducibles $P\in R_T$ tales que
$P\equiv 1\bmod M$.
\end{corolario}

\begin{proof}
Sea $\{P_1,\ldots, P_r\}$ un conjunto finito de cardinalidad $r\geq 0$
de polinomios irreducibles que satisfacen $P_i\equiv 1\bmod M$.
Sea $N:=MP_1\cdots P_r$ y sea $Q\in R_T$ arbitrario. Entonces
$\Psi_M(NQ)\equiv \Psi_M(0)\bmod N$. Ahora bien, puesto que
\begin{gather*}
\Psi_M(u)=\prod\limits_{\substack{D|M\\ D \text{\ m\'onico}}} \big(
u^D\big)^{\mu(M/D)}, \quad \Psi_P(u)=\frac{u^P}{u}\\
\intertext{(Proposici\'on \ref{P6.2.21}) y} 
\frac{u^M}{u}|_{u=0}=\carlitzbinom M0=M,\\
\intertext{entonces}
\Psi_M(0)=
\begin{cases} 
R&M=R^n\text{\ para alg\'un
$R$ irreducible, $n\geq 1$}\\ 1& \text{en otro caso}
\end{cases}.
\end{gather*}

Primero supongamos que $M$ no es potencia de un polinomio
irreducible, es decir, $\Psi_M(0)=1$. Entonces $\Psi_M(NQ)\equiv
1\bmod N$. Por lo tanto $\Psi_M(NQ)\equiv 1\bmod M$ y $\Psi_M(
NQ)\equiv 1\bmod P_i$, $1\leq i\leq r$. En particular, $P_i\nmid
\Psi_M(NQ)$.

Sea $P$ cualquier polinomio irreducible que divide a $\Psi_M(NQ)$.
Entonces $P\equiv 1\bmod M$ por la Proposici\'on \ref{P6.3.11} y
$P\neq P_i$.

Ahora si $M$ es potencia de un irreducible, $M=R^n$, se tiene
$P_i\neq R$, $1\leq i\leq r$ y $\Psi_M(NQ)\equiv R\bmod M$;
$\Psi_M(NQ)\equiv R\bmod P_i$.
Si $P|\Psi_M(NQ)$ con $P$ un polinomio irreducible, entonces
si $P=P_i$ para alg\'un $i$ con $1\leq i\leq r$. De 
aqu{\'\i} se seguir{\'\i}a que $P_i|R$ lo
cual implicar{\'\i}a que $P_i=R$ que es absurdo. Se sigue que
$P\equiv 1\bmod M$ y $P\neq P_i$ para todo $i$, $1\leq i\leq r$ (si
$r=0$ la condici\'on es vac{\'\i}a). $\fin$
\end{proof}

\begin{observacion}\label{O6.3.13} El Corolario \ref{C6.3.12} es un
caso particular del Teorema de Dirichlet sobre la infinitud de
polinomios irreducibles en progresiones geom\'etricas (ver Teorema
\ref{T6.3.14}) y el cual es consecuencia del Teorema de Densidad
de \v{C}ebotavev\index{Cebotarev@\v{C}ebotarev!teorema 
de densidad de $\sim$}
el cual no probaremos aqu{\'\i}.
\end{observacion}

\begin{teorema}[Dirichlet\index{teorema de
Dirichlet}\index{Dirichlet!teorema de $\sim$}]\label{T6.3.14}
Sean $M,N\in R_T$ cualesquiera dos polinomios m\'onicos no
constantes primos relativos. Entonces existe un infinidad de
polinomios irreducibles $P\in R_T$ tales que $P\equiv N\bmod M$.
$\fin$
\end{teorema}

\begin{observacion}\label{O6.3.14'}
La descripci\'on de los subcampos de un campo de funciones
ciclot\'omico $\cicl M{}$, $M\in R_T$, no es expl\'icita como
en el caso de campos num\'ericos (ver Secci\'on \ref{S4.3}).
Tenemos que la primera parte de esta descripci\'on es 
totalmente an\'aloga al caso num\'erico, esto es, si $M=
P_1^{\alpha_1}\cdots P_r^{\alpha_r}$ con $P_1,\ldots,
P_r\in R_T^+$ primos distintos, $\alpha_i\geq 1$ y $r\geq 1$,
entonces $\cicl M{}=\prod_{i=1}^r\cicl {P_i}{\alpha_i}$ y
por ende, basta describir los subcampos de $\cicl Pn$ con
$P\in R_T^+$ y $n\in{\ma N}$ para describir en general
los subcampos de $\cicl M{}$ para $M\in R_T$ arbitrario
no constante. Sin embargo el caso $\cicl Pn$ no es 
expl\'icito como en el caso num\'erico, debido fundamentalmente
a dos razones.

(1) $\cicl Pn/\cicl P{n-1}$ es una extensi\'on elemental abeliana
de grado $q^d$ con $d=\deg P$ y no c\'iclica de grado primo
como en el caso num\'erico. 

(2) Los subcampos de $\cic p{}$ con $p$ primo pudieron darse
de manera expl\'icita gracias a que $\{\zeta_p^{\sigma}\}_{\sigma
\in U_p}$ es una base normal de $\cic p{}/{\ma Q}$. En el caso
de campos de funciones, se tiene que si $\lambda$ es un
generador de $\Lambda_P$, en general $\{\lambda^A\}_{A\in
(R_T/\langle P\rangle)^*}$ no es una base normal de $\cicl P{}/
K$. Por ejemplo, si $q=2$, $p(T)=T^4+T+1\in R_T^+$, 
entonces $\big({\ma F}_2/\langle P(T)\rangle\big)\cong
{\ma F}_{2^4}={\ma F}_{16}$, $\big({\ma F}_2/
\langle P(T)\rangle\big)^*\cong {\ma F}_{16}^*\cong C_{15}$.
Si $A(T)=T$, entonces $o(A\bmod p(T))=15$ puesto que
$A^i\not\equiv 1\bmod p(T)$, $i\in\{1,3,5\}$ y $\lambda+
\lambda^A=\lambda+\lambda^T=\lambda^{1+T}=
\lambda^{A^4}$ por lo que $\{\lambda^{A^i}\}_{0\leq i\leq 14}$
no es una base normal de $\cicl {p(T)}{}/K$.

En el caso (1),  el estudio de las subextensiones
de $\cicl Pn/\cicl P{}$, se puede consultar el
Cap\'itulo \ref{Ch11} en donde se demuestra el Teorema
de Kronecker--Weber. Para (2) no tenemos ecuaciones
expl\'icitas. Las bases normales de $\cicl P{}/K$ se 
obtienen en \cite[Theorem 4]{Cha91}.
\end{observacion}

Terminamos esta secci\'on con una observaci\'on sobre los
{\em n\'umeros de clase\index{numero de clase@n\'umero de clase}} de
campos ciclot\'omicos. La definici\'on de n\'umero de clase 
la daremos m\'as adelante (Definici\'on \ref{D10.1.1.1}).

\begin{observacion}\label{O6.3.15} 
Si $M\in R_T$ y $h_M$ denota al n\'umero de clase del
campo ciclot\'omico $\cicl M{}$, entonces tenemos que si
$N,M\in R_T$ son tales que $N|M$, entonces
$h_N|h_M$. Este es un caso particular
del Teorema \ref{T10.1.2.15}.
\end{observacion}

\section{Caracteres de Dirichlet\index{caracteres de
Dirichlet}\index{Dirichlet!caracteres de $\sim$}}\label{S6.4}

\begin{definicion}\label{D6.4.1}
Sea $M\in R_T\setminus\{0\}$ un polinomio m\'onico. Un
{\em caracter de Dirichlet m\'odulo $M$\index{caracter de Dirichlet
m\'odulo un polinomio}} es un homomorfismo
\[
\chi\colon \units M{}\longto {\ma C}^{\ast}.
\]
\end{definicion}

\begin{observacion}\label{O6.4.2} Si $M$ divide a $N$ en $R_T$, 
tenemos el epimorfismo can\'onico
\begin{eqnarray*}
\varphi_{N,M}\colon \units N{}&\twoheadrightarrow& \units M{}\\
A\bmod N&\mapsto& A\bmod M.
\end{eqnarray*}

Entonces para todo caracter de Dirichlet m\'odulo $M$, $\chi\colon
\units M{}\to {\ma C}^{\ast}$, $\varphi_{N,M}$ induce un 
caracter de Dirichlet m\'odulo $N$: $\chi\circ \varphi_{N,M}\colon
\units N{}\longto {\ma C}^{\ast}$,
\[
\xymatrix{
\units N{}\ar[rr]^{\chi\circ\varphi_{N,M}}\ar[dr]_{\varphi_{N,M}}&&
{\ma C}^{\ast}\\
&\units M{}\ar[ru]_{\chi}
}
\]

Rec{\'\i}procamente, si $\chi$ es un caracter de Dirichlet m\'odulo
$M$, decimos que podemos definir $\chi$ m\'odulo $F$ para
$F|M$ si existe $\xi\colon \units F{}\to{\ma C}^{\ast}$ tal que
$\xi\circ \varphi_{M,F}=\chi$.
\[
\xymatrix{
\units M{}\ar[rr]^{\chi}\ar[ddr]_{\varphi_{M,F}}&& {\ma C}^{\ast}\\
&\text{\Huge{$\circlearrowleft$}}
\\ &\units F{}\ar[uur]_{\xi}
}
\]
\end{observacion}

\begin{observacion}\label{O6.4.2'} Si $\chi$ es un caracter
de Dirichlet definido m\'odulo $M\in R_T$, entonces si 
$F|M$, $F\in R_T$, se tiene que $\chi$ se puede definir m\'odulo
$F$ si y solamente si $\chi(A\bmod M)=\chi(B\bmod M)$ para
cualesquiera $A,B\in R_T$ primos relativos a $M$ y tales que
$A\equiv B\bmod F$.
\end{observacion}

\begin{teorema}[Existencia del conductor]\label{T6.4.3}
Sea $\chi$ un caracter de Dirichlet definido m\'odulo $M$. Entonces
existe un polinomio m\'onico \'unico $F$ en $R_T$ de grado
m{\'\i}nimo que divide a $M$ tal que $\chi$ puede ser definido 
m\'odulo $F$.
\end{teorema}

\begin{proof}
Sea $\chi\colon \units M{}\to {\ma C}^{\ast}$. Sean $A,B\in R_T$
m\'onicos que dividen a $M$ y tales que $\chi$ puede ser definido
m\'odulo $A$ y tambi\'en m\'odulo $B$, es decir, existen $\chi_A
\colon \units A{}\to{\ma C}^{\ast}$ y $\chi_B\colon \units B{}\to
{\ma C}^{\ast}$ tales que $\chi=\chi_A\circ\varphi_{M,A}$, 
$\chi=\chi_B\circ\varphi_{B,M}$. 

Consideremos $C:=\mcd (A,B)$
el m\'aximo com\'un divisor de $A$ y $B$. Sea $D$ el producto
de todos los polinomios m\'onicos irreducibles que dividen a $M$
pero que no dividen a $B$. Entonces $C=\mcd (DA,B)$. 

Para
ver que podemos definir $\chi$ m\'odulo $C$, consideremos
$U, V\in R_T$ primos relativos a $M$ tales que $U\equiv V\bmod C$.
Por el Teorema Chino del Residuo, existe $S\in R_T$ tal que 
$S\equiv U\bmod DA$ y $S\equiv V\bmod B$.

Veamos que $S$ y $M$ son primos relativos. En caso contrario
existir{\'\i}a $P\in R_T$ irreducible que divide a $S$ y a $M$. Sea
$S=V+QB$, entonces si $P|B$, entonces $P|V$ pero en este
caso se seguir{\'\i}a que $P|\mcd (V,M)=1$ lo cual es absurdo, 
esto es, $P\nmid B$. Ahora bien, puesto que $P|M$ y $P\nmid B$
entonces $P$ es un factor de $D$ y por lo tanto $P|DA$. Pero
$P|S$ por lo cual $P|U$ y por lo tanto $P|\mcd (U,M)=1$ lo cual
es absurdo. En resumen, $\mcd (S,M)=1$. Entonces
\begin{gather*}
\chi(S)=\chi_A\circ \varphi_{M,A}(S)=\chi_{A}\circ \varphi_{M,A}(U)
=\chi(U)
\intertext{y}
\chi(S)=\chi_B\circ \varphi_{M,B}(S)=\chi_B\circ \varphi_{M,B}(V)=
\chi(V).
\end{gather*}

Por lo tanto $\chi(S)=\chi(U)=\chi(V)$ de donde $\chi$ puede ser
definido m\'odulo $C$:
\[
\xymatrix{
\units M{}\ar[rr]^{\chi}\ar[rd]_{\varphi_{M,C}}&&{\ma C}^{\ast}\\
&\units C{}\ar[ru]_{\chi_C}
}
\]

Finalmente si $\chi$ puede ser definido m\'odulo $F_1$ y 
m\'odulo $F_2$ con $F_1$ y $F_2$ de grado m{\'\i}nimo y $F_1$,
$F_2$ m\'onicos, entonces $\chi$ puede ser definido m\'odulo $C
=\mcd (F_1,F_2)$. Puesto que $C|F_1$ y $C|F_2$ se sigue que
$C=F_1$ y $C=F_2$, es decir, $F_1=F_2$. $\fin$
\end{proof}

\begin{definicion}\label{D6.4.4}
El polinomio dado en el Teorema \ref{T6.4.3} se llama 
{\em conductor de\index{conductor de un caracter de Dirichlet}} $\chi$
y se denota por $F_{\chi}$. En otras palabras, si $\chi$ es un 
caracter de Dirichlet definido m\'odulo $M$, entonces $F_{\chi}$
es el \'unico polinomio de grado m{\'\i}nimo que divide a $M$ y
tal que $\chi$ puede definirse m\'odulo $F_{\chi}$.
\end{definicion}

\begin{observacion}\label{O6.4.5} Sean $q=2$ y $M\in R_T
\setminus\{0\}$ m\'onico tal que $\mcd (M,T)=\mcd (M,T+1)
=1$. Entonces no
existe ning\'un caracter de Dirichlet $\theta$ tal que $F_{\theta}=
TM$ ni $F_{\theta}=(T+1)M$.

En efecto notemos que $\Phi(MT)=\Phi(M)\Phi(T)=\Phi(M)$ y
$\Phi(M(T+1))=\Phi(M)\Phi(T+1)=\Phi(M)$ ya que $\Phi(T)=
\Phi(T+1)=1$ pues estamos en el caso $q=2$. Entonces
\[
\units {TM}{}\cong \units {(T+1)M}{}\cong \units M{}.
\]
En particular, con $M=1$, vemos que para $q=2$ no hay 
caracter de conductor $T$, $T+1$ o $T(T+1)$.
\end{observacion}

\begin{ejemplo}\label{Ej6.4.6} Sea $q=2$ y sea $\chi\colon \units T3
\to {\ma C}^{\ast}$ dado por
\begin{eqnarray*}
1&\longmapsto&1,\\
T+1&\longmapsto&-1,\\
T^2+1&\longmapsto&1,\\
T^2+T+1&\longmapsto&-1.
\end{eqnarray*}

Puesto que $\chi(T^2+A)=\chi(A)$ para toda $A\in \units T3$ entonces
$\chi$ se puede definir m\'odulo $T^2$ pues si $\xi\colon \units T2\to
{\ma C}^{\ast}$, $\xi(1)=1$, $\xi(1+T)=-1$ y $\varphi_{T^3,T^2}
\colon\units T3\to \units T2$, entonces
\begin{gather*}
\varphi_{T^3,T^2}(1)=\varphi_{T^3,T^2}(T^2+1)=1 \quad\text{y}\\
\varphi_{T^3,T^2}(T+1)=\varphi_{T^3,T^2}(T^2+T+1)=T+1
\end{gather*}
se sigue que $\xi\circ\varphi_{T^3,T^2}=\chi$. Por la Observaci\'on
\ref{O6.4.5}, se tiene que $F_{\chi}=T^2$.
\end{ejemplo}

\begin{ejemplo}\label{Ej6.4.7} Sea $q=2$ y sea $\chi\colon 
\big(R_T/\langle T^2(T+1)\rangle
\big)^{\ast}\to{\ma C}^{\ast}$ dado por $\chi(1)
=1$ y $\chi(T^2+T+1)=-1$. Entonces si $\xi\colon\units T2\to{\ma C}^{
\ast}$, $\xi(1)=1$, $\xi(1+T)=-1$ satisface $\xi\circ\varphi_{T^2(T+1),
T^2}=\chi$. Por lo tanto $F_{\chi}=T^2$.

La misma conclusi\'on puede obtenerse usando dos veces la
Observaci\'on \ref{O6.4.5}: $\big(R_T/\langle T^2(T+1)\rangle
\big)^{\ast}\cong \units T2$.
\end{ejemplo}

\begin{ejemplo}\label{Ej6.4.8}
Sea $q=2$, $N=T^2+T+1$, $M=NT$, $\omega=e^{2\pi i/3}$ y
$\theta\colon \units M{}\to {\ma C}^{\ast}$ dado por
\begin{eqnarray*}
1&\longmapsto& 1,\\
T^2+1&\longmapsto& \omega,\\
T+1&\longmapsto& \omega^2.
\end{eqnarray*}
Puesto que $\units {NT}{}\cong \units N{}$, se pueden definir 
$\tilde{\theta}\colon \units N{}\to{\ma C}^{\ast}$, $\tilde{\theta}(1)=1$,
$\tilde{\theta}(T)=\theta(T^2+1)=\omega$ y $\tilde{\theta}(T+1)=
\omega^2$. Por lo tanto $F_{\theta}=T^2+T+1$.
\end{ejemplo}

\begin{observacion}\label{O6.4.9} Dado un caracter de Dirichlet 
podemos considerar $\chi$ como un mapeo $\chi\colon R_T\to
{\ma C}$ definiendo $\chi(Q)=0$ si $\mcd (Q,F_{\chi})\neq 1$.
En caso de no especificarse, siempre consideraremos a un
caracter $\chi$ definido m\'odulo su conductor $F_{\chi}$.
\end{observacion}

\begin{definicion}\label{D6.4.10} Un caracter de Dirichlet $\chi$
 definido
m\'odulo su conductor se llama {\em primitivo\index{caracter
primitivo}}. En este caso se hace $\chi(Q)=0$ tan poco como sea posible.
\end{definicion}

Tambi\'en notemos que cuando $\chi$ est\'a definido m\'odulo
su conductor, tenemos que $\chi(A+F_{\chi})=\chi(A)$, esto es,
$\chi$ es peri\'odico de per{\'\i}odo $F_{\chi}$.

\begin{notacion}\label{N6.4.11}
Siempre que mencionemos los caracteres de $\units M{}$, $M\in
R_T$ o {\em caracteres m\'odulo $M$\index{caracteres m\'odulo un
polinomio}} incluiremos todos los caracteres cuyos conductores
dividan a $M$. El {\em caracter trivial\index{caracter trivial}} 
$\varepsilon$ satisface $\varepsilon(Q)=1$ para todo $Q\in R_T$.
Si $G$ es cualquier grupo, $\hat{G}$ denota al conjunto de sus
caracteres: $\hat{G}:=\Hom(G,{\ma C}^{\ast})=\{\chi\colon G\to
{\ma C}^{\ast}\mid \chi \text{\ es homomorfismo de grupos}\}$.
\end{notacion}

\begin{definicion}\label{D6.4.12} Diremos que un caracter es
{\em par\index{caracter par}} si $\theta(a)=1$ para todo $a\in{\ma
F}_q^{\ast}$.
\end{definicion}

\begin{proposicion}\label{P6.4.13} Sea $X:=\big\{\theta\in\widehat{
\units N{}}\mid \theta\text{\ es par}\big\}$. Entonces $X$ es 
subgrupo de $\widehat{\units N{}}$ de orden $\frac{\Phi(N)}{q-1}$.
\end{proposicion}

\begin{proof}
Se tiene la sucesi\'on exacta
\begin{gather*}
0\longto {\ma F}_q^{\ast}\longto \units N{}.\\
\intertext{Tomando duales, se obtiene la sucesi\'on exacta}
\begin{array}{rclcc}
\widehat{\units N{}}&\stackrel{{\eu X}}{\longto}&{\ma F}_q^{\ast}&
\longto &0\\
\theta&\longmapsto&\theta|_{{\ma F}_q^{\ast}}
\end{array}
\end{gather*}

Se tiene que $X=\ker {\eu X}$ y por tanto $\widehat{\units N{}}/X\cong
{\ma F}_q^{\ast}$ de donde se sigue que 
\[
|X|=\frac{\big|\widehat{\units N{}}\big|}{\big|{\ma F}_q^{\ast}\big|}=
\frac{\big|\units N{}\big|}{\big|{\ma F}_q^{\ast}\big|}=\frac{\Phi(N)}{
q-1}. \tag*{$\fin$}
\]
\end{proof}

\begin{definicion}\label{D6.4.14} Sean $\chi$, $\phi$ dos caracteres
de Dirichlet de conductores $F_{\chi}$ y $F_{\phi}$ respectivamente.
Se define el producto de $\chi$ y $\phi$ como sigue. Sea $Q:=
\big[F_{\chi},F_{\phi}\big]$ y se define $\gamma\colon \units Q{}
\to {\ma C}^{\ast}$ por $\gamma(A\bmod Q)=\chi(A\bmod Q)\phi(
A\bmod Q)$. Entonces el {\em producto\index{producto de 
caracteres}} $\chi\phi$ se define como el caracter primitivo asociado
a $\gamma$. En particular $F_{\chi\phi}|\big[F_{\chi}, F_{\phi}\big]$.
\end{definicion}

\begin{teorema}\label{T6.4.19} Si $\mcd (F_{\chi},F_{\phi})=1$
entonces $F_{\chi\phi}=F_{\chi}F_{\phi}$.
\end{teorema}

\begin{proof}
 Sean $N=F_{\chi}$, $M=F_{\phi}$, $S:=[N,M]=NM$. Definimos
 $\gamma\colon
 \units S{}\to{\ma C}^{\ast}$ dada por $\gamma(A\bmod S)=\gamma(
 A\bmod S)\phi(A\bmod S)$.
 
 Ahora bien, puesto que $[S,N]=S$, donde $[S,N]$
 denota al m\'inimo com\'un m\'ultiplo de $S$ y
$N$, se puede definir $\theta\colon
 \units S{}\to{\ma C}^{\ast}$ por $\theta(A\bmod S)=\gamma(A\bmod
 S)\chi^{-1}(A\bmod S)$. Se obtiene $\theta=\phi \bmod S$ lo que
 implica que $F_{\theta}=F_{\phi}=M$. Por lo tanto $M=F_{\theta}=
 F_{\gamma \chi^{-1}}|\big[F_{\gamma},F_{\chi^{-1}}\big]$.
 
 Es decir, $M|\big[F_{\gamma},F_{\chi^{-1}}\big]=\big[F_{\gamma},
 F_{\theta}\big]=\frac{F_{\gamma}N}{\mcd (F_{\gamma},N)}=F_{
 \gamma} N_1$ donde $N_1=\frac{N}{\mcd(F_{\gamma},N)}$.
Puesto que $\mcd (N,M)=1$ se sigue que $\mcd (N_1,M)=1$ y
$M|F_{\gamma}$. An\'alogamente $N|F_{\gamma}$ y puesto que
$\mcd (N,M)=1$, se tiene que $NM|F_{\gamma}$. Por otro lado
$F_{\gamma}=F_{\chi\phi}|\big[F_{\chi},F_{\phi}\big]=[N,M]=NM$.
Se sigue que $F_{\gamma}=NM=F_{\chi}F_{\phi}$.  $\fin$
\end{proof}

\begin{proposicion}\label{P6.4.20} Sean $\chi$, $\sigma$ dos
caracteres de Dirichlet de conductores $F_{\chi}$ y $F_{\sigma}$
respectivamente. Supongamos que existe $N$ tal que $F_{\chi}|N$,
$F_{\sigma}|N$ y $\chi, \sigma\colon \units N{}\to {\ma C}^{\ast}$ son
iguales m\'odulo $N$, es decir, $\chi(A\bmod N)=\sigma(A\bmod N)$
para todo $A$ primo relativo a $N$. Entonces $F_{\chi}=
F_{\sigma}$ y $\chi=\sigma \bmod F_{\chi}$, esto es, $\chi=\sigma$.
\end{proposicion}

\begin{proof}
Consideremos
\[
\xymatrix{
\units N{}\ar[r]^{\qquad \chi}\ar[d]_{\varphi_{N,F_{\chi}}}&
{\ma C}^{\ast}\\
\units {F_{\chi}}{}\ar[ru]_{\tilde{\chi}}
}
\]
Tenemos $\tilde{\chi}\circ \varphi_{N,F_{\chi}}=\chi=\sigma$, es
decir, $\sigma$ se puede definir m\'odulo $F_{\chi}$ lo cual implica
que $F_{\sigma}|F_{\chi}$. Por simetr{\'\i}a, se tiene que  $F_{\chi}|
F_{\sigma}$ y por tanto $F_{\chi}=F_{\sigma}$. Sean $\chi'\equiv
\chi\bmod F_{\chi}$, $\sigma'\equiv \sigma\bmod F_{\sigma}$. Sea
$A\in R_T\setminus\{0\}$ tal que $\mcd (A,F_{\chi})=1$. Existe
$A'\in R_T$ tal que $\mcd (A,N)=1$ y $A\bmod F_{\chi}=A'\bmod
F_{\chi}$. Luego $\chi'(A\bmod F_{\chi})=\chi'(A'\bmod F_{\chi})=
\chi(A'\bmod N)=\sigma(A'\bmod N)=\sigma(A'\bmod F_{\chi})=
\sigma'(A\bmod F_{\chi})$. $\fin$
\end{proof}

\begin{observacion}\label{O6.4.15} En general no se tiene que 
$(\chi\phi)({A})=\chi({A})\phi({A})$.
\end{observacion}

\begin{ejemplo}\label{Ej6.4.16} Sea $q=2$ y sea $\chi$ m\'odulo
$T^2(T^2+1)$ dado por 
\[
\chi(1)=1,\quad \chi(T^2+T+1)=1,\quad \chi(T^3+T^2+1)=-1,
\quad \chi(T^3+T+1)=-1.
\]

Por la Observaci\'on \ref{O6.4.5} se tiene que el conductor de $\chi$,
$F_{\chi}\in \{1,T^2,T^2+1,T^2(T^2+1)\}$. Como $\chi(T^3+T^2+1)=
-1$, $F_{\chi}\neq 1$.

Ahora $T^3+T^2+1\bmod T^2=1$ pero $\chi(T^3+T^2+1)=-1\neq 1$ 
y $T^3+T+1\bmod (T^2+1)=1$ pero $\chi(T^3+T+1)=-1\neq 1$, por
lo que $F_{\chi}\neq T^2,T^2+1$. Se sigue que $F_{\chi}=T^2(T^2
+1)$.

Ahora sea $\varphi\bmod T^2$ dada por $\varphi(1)=1$, $\varphi(1+
T)=-1$. Por tanto $F_{\varphi}=T^2$.

Consideremos el producto $\chi\varphi$. Se tiene que $\big[F_{\chi},
F_{\varphi}\big]=\big[T^2(T^2+1),T^2\big]=T^2(T^2+1)$ y definimos
$\gamma\colon \big(R_T/T^2(T^2+1)\big)^{\ast}\to {\ma C}^{\ast}$
por $\gamma(A)=\chi(A)\varphi(A)$. Entonces 
\begin{align*}
\gamma(1)&=\chi(1)\varphi(1)=1\cdot 1=1,\\
\gamma(T^2+T+1)&=\chi(T^2+T+1)\varphi(T^2+T+1)=(1)(-1)=-1,\\
\gamma(T^3+T^2+1)&=\chi(T^3+T^2+1)\varphi(T^3+T^2+1)
=(-1)(1)=-1,\\
\gamma(T^3+T+1)&=\chi(T^3+T+1)\varphi(T^3+T+1)
=(-1)(-1)=1.
\end{align*}

Sea $\xi\colon\units {T^2+1}{}\to{\ma C}^{\ast}$ dado por $\xi(1)=1$
y $\xi(T)=1$. Entonces $\xi\circ \varphi_{T^2(T^2+1),T^2+1}=
\gamma$. Por tanto $F_{\gamma}=T^2+1$ y $\xi=\chi\phi$. Notemos 
que $\xi(T)=-1\neq 0=\varphi(T)=\chi(T)\varphi(T)$.
\end{ejemplo}

\begin{ejemplo}\label{Ej6.4.16bis} Sean $q=2$, $\zeta=\zeta_6=
e^{2\pi i/6}$, $\omega=\zeta_3=e^{2\pi i/3}$, $M=T^2 N$ con
$N=T^2+T+1$. Definimos
\[
\begin{array}{rcl}
\chi\colon\units M{}&\longto&{\ma C}^{\ast}\\
1&\longmapsto&1\\
T+1&\longmapsto&\zeta\\
T^2+1&\longmapsto&\zeta^2\\
T^3+T^2+T+1&\longmapsto&-1\\
T^3+T^2+1&\longmapsto&-\zeta\\
T^3+T+1&\longmapsto&-\zeta^2
\end{array}
\qquad\qquad
\begin{array}{rcl}
\sigma\colon\units {{T^2}}{}&\longto&{\ma C}^{\ast}\\
1&\longmapsto&1\\
T+1&\longmapsto&-1
\end{array}
\]
Se tiene $F_{\chi}=M$, $F_{\sigma}=T^2$ y por tanto $\big[
F_{\chi},F_{\sigma}\big]=M$. Obtenemos $\psi\bmod M$:
\[
\begin{array}{rcccr}
\units M{}&\longto&\units {{T^2}}{}&\longto &{\ma C}^{\ast}\\
1 &\longmapsto&1&\longmapsto &1\\
T+1 &\longmapsto&T+1&\longmapsto &-1\\
T^2+1 &\longmapsto&1&\longmapsto &1\\
T^3+T^2+T+1 &\longmapsto&T+1&\longmapsto &-1\\
T^3+T^2+1 &\longmapsto&1&\longmapsto &1\\
T^3+T+1 &\longmapsto&T+1&\longmapsto &-1.
\end{array}
\]

Ahora sea $\gamma\colon \units M{}\to {\ma C}^{\ast}$ dada por
\[
\gamma(A)=\chi(A)\sigma(A)=
\begin{cases}
1&\text{si $A=1$}\\
\omega^2&\text{si $A=T+1$}\\
\omega&\text{si $A=T^2+1$}\\
1&\text{si $A=T^3+T^2+T+1$}\\
\omega^2&\text{si $A=T^3+T^2+1$}\\
\omega&\text{si $A=T^3+T+1$}
\end{cases}
\]
Podemos definir $\gamma$ m\'odulo $N$:
\[
\begin{array}{rcccc}
\units M{}&\longto&\units N{}&\longto &{\ma C}^{\ast}\\
1 &\longmapsto&1&\longmapsto &1\\
T+1 &\longmapsto&T+1&\longmapsto &\omega^2\\
T^2+1 &\longmapsto&T&\longmapsto &\omega\\
T^3+T^2+T+1 &\longmapsto&1&\longmapsto &1\\
T^3+T^2+1 &\longmapsto&T+1&\longmapsto &\omega^2\\
T^3+T+1 &\longmapsto&T&\longmapsto &\omega.
\end{array}
\]

El producto de caracteres $\chi$ y $\sigma$ es:
\[
\begin{array}{rcl}
\chi\sigma\colon\units N{}&\longto&{\ma C}^{\ast}\\
1&\longmapsto&1\\
T&\longmapsto&\omega\\
T+1&\longmapsto& \omega^2
\end{array}
\]
y $F_{\chi\sigma}=N$. Notemos que $(\chi\sigma)(T)=\omega\neq
0=\chi(T)\sigma(T)$.
\end{ejemplo}

\begin{definicion}\label{D6.4.17}
Si $\chi$ es un caracter de Dirichlet, definimos el 
{\em conjugado\index{caracter conjugado}} $\overline{\chi}$
de $\chi$ por 
$\overline{\chi}(A)=\overline{\chi(A)}$. Notemos que $\overline{\chi}
(A)=\chi(A)^{-1}$ para toda $\mcd (A,F_{\chi})=1$. Por tanto
$\chi\overline{\chi}$ es el caracter trivial y $F_{\overline{\chi}}=
F_{\chi}$.
\end{definicion}

\begin{observacion}\label{O6.4.18} Tenemos $G_M=\Gal(\cicl M{}/K)
\cong \units M{}$. Entonces un caracter de Dirichlet m\'odulo $M$
es un caracter de $G_M$, por lo que el caracter de Dirichlet
puede ser considerado un {\em caracter de
Galois\index{Galois!caracter de $\sim$}}.
\end{observacion}

\begin{definicion}\label{D6.4.21} Sea $\chi$ un caracter de Dirichlet
m\'odulo $M$, esto es, $\chi\in \widehat{G_M}\cong \widehat{\units
M{}}$. Entonces $\ker \chi\subseteq G_M$. Sea $K_{\chi}:=
\cicl M{}^{\ker \chi}$. El campo $K_{\chi}$ se llama el {\em campo
perteneciente a $\chi$\index{campo perteneciente a un caracter}}
o que $K_{\chi}$ {\em est\'a asociado\index{campo asociado a un
caracter}} a $\chi$.
\end{definicion}

\begin{ejemplo}\label{Ej6.4.29} Sea $\chi$ el caracter del Ejemplo
\ref{Ej6.4.6}. Entonces 
\begin{gather*}
\chi\colon \units T3\cong G_{T^3}=\Gal(
\cicl T3/K)\to{\ma C}^{\ast} \quad \text{y}\\
\ker \chi=\{1\bmod T^3,(T^2+1)
\bmod T^3\}.
\end{gather*}
Por lo tanto $\chi$ es un caracter de $\units T3/\ker
\chi\cong \units T2\cong \Gal(\cicl T2/K)$ y puede ser
considerado un caracter de $\Gal(\cicl T2/K)$.
\end{ejemplo}

\begin{ejemplo}\label{Ej6.4.30} Sea $\chi$ el caracter del 
Ejemplo \ref{Ej6.4.7}. Entonces $\units {T^2(T+1)}{}\cong
\units T2$ y puesto que cualquier caracter m\'odulo $T^2(T+1)$
o m\'odulo $T^2$ es el mismo caracter, se sigue que
$\cicl {T^2(T+1)}{}=\cicl T2$.
\end{ejemplo}

\begin{ejemplo}\label{Ej6.4.31}
Sean $q=3$, $M=T^2+1$, $\zeta=\zeta_8=e^{2\pi i/8}$. Se tiene
que $\Phi(M)=8$, $\units M{}=\{1,T+1,-T,-T+1,-1,-T-1,T,T-1\}$ y
$\Lam M{}=\{u\in\overline{K}\mid u^M=0\}=\{u\in\overline{K}\mid
\big((\varphi+\mu_T)^2+\Id\big)(u)=0\}=\{u\in\overline{K}\mid
u^9+(Tu)^3+Tu^3+T^2u+u=0\}$.
Entonces $u(u^8+T^3u^2+Tu^2+1)=0$. Por lo tanto
\[
\Psi_M(u)=u^8+T^3u^2+Tu^2+T^2+1=u^8+(T^3+T)u^2+(T^2+1).
\]

Sea $\lambda$ una ra{\'\i}z de $\Psi_M(u)$ y sean $\sigma_1=\Id$,
$\sigma_{-1}\colon\lambda\mapsto -\lambda$, $\cicl M{}^+=\{u\in
\cicl M{}\mid \sigma_{-1}(u)=u\}$ Puesto que
\begin{gather*}
\cicl M{}=\{A_0+A_1\lambda+\cdots+A_7\lambda^7\mid A_i\in K\}
\intertext{tenemos}
\cicl M{}^+=\{A_0+A_2\lambda^2+A_4\lambda^4+A_6\lambda^6\mid
A_i\in K\}.
\end{gather*}

Sea
\[
\begin{array}{rcl}
\theta\colon\units M{}&\longto&{\ma C}^{\ast}\\
1&\longmapsto& 1\\
T+1&\longmapsto& \zeta\\
-T&\longmapsto& \zeta^2=i\\
-T+1&\longmapsto& \zeta^3\\
-1&\longmapsto& \zeta^4=-1\\
-T-1&\longmapsto& \zeta^5\\
T&\longmapsto& \zeta^6\\
T-1&\longmapsto& \zeta^7
\end{array}
\qquad\text{luego}\qquad
\begin{array}{rcl}
\theta^2\colon\units M{}&\longto&{\ma C}^{\ast}\\
1&\longmapsto& 1\\
T+1&\longmapsto& i\\
-T&\longmapsto&-1\\
-T+1&\longmapsto& -i\\
-1&\longmapsto& 1\\
-T-1&\longmapsto&i\\
T&\longmapsto&-1\\
T-1&\longmapsto& -i.
\end{array}
\]
Por lo tanto $\ker \theta^2=\{1,-1\}\cong\{\sigma_1,\sigma_{-1}\}=
J<G=G_M=\Gal(\cicl M{}/K)$ luego $\cicl M{}^+=\cicl M{}^{\ker 
\theta^2}$. Entonces $\cicl M{}^+$ es el campo perteneciente
a $\theta^2$.
\end{ejemplo}

\begin{observacion}\label{O6.4.23} Sea $\chi$ un caracter de 
Dirichlet definido m\'odulo $M$ y sea $N\in R_T\setminus \{0\}$ 
un m\'ultiplo de $M$. Sea $\tilde{\chi}$ el caracter $\chi$ definido
m\'odulo $N$, es decir,
\[
\xymatrix{
\units N{}\ar[rr]^{ \tilde{\chi}}\ar[rd]_{\varphi_{N,M}}&&{\ma C}^{\ast}\\
&\units M{}\ar[ru]_{\chi}
}
\qquad \tilde{\chi}=\chi\circ \varphi_{N,M}.
\]

Sean $K_1=\cicl M{}^{\ker \chi}$ y $K_2=\cicl M{}^{\ker \tilde{\chi}}$.
Entonces
\[
\ker \varphi_{N,M}=D_{N,M}=\{A\bmod N\mid A\equiv 1\bmod M\},
\]
$\units N{}/\ker \varphi_{N,M}\cong \units M{}$. Puesto que $G_N\cong
\units N{}$ y $G_M\cong \units M{}$, $H=\Gal(\cicl N{},\cicl M{})$
\[
\qquad\qquad \left.\begin{array}{c}
\xymatrix{
\cicl N{}\ar@{-}[d]_{H \Bigg\{\qquad}\\ 
\cicl M{}\ar@{-}[d]^{G_M}\\K
}\end{array}\right\}{\ }_{G_N}
\quad
\begin{array}{c}
\text{Entonces $\cicl M{}=\cicl N{}^H$ y}\\
\ker \varphi_{N,M}\cong \Gal(\cicl N{}/\cicl M{})\cong H.
\end{array}
\]

Ahora bien, $\ker\tilde{\chi}=\varphi_{N,M}^{-1}(\ker \chi)$ y ya
que $\varphi_{N,M}^{-1}(\ker\chi)\supseteq \varphi_{N,M}^{-1}(
\{1\})=\ker \varphi_{N,M}$ se sigue que 
\[
K_2=\cicl N{}^{\ker \tilde{
\chi}}\subseteq \cicl N{}^{\ker \varphi_{N,M}}=\cicl M{}.
\]

Por tanto $K_2\subseteq \cicl M{}^{\ker \chi}=K_1$. Por otro lado
\begin{align*}
\big|\ker \tilde{\chi}\big|&=\big|\varphi_{N,M}^{-1}(\ker \chi)\big|=
\big|\ker \varphi_{N,M}\big|\big|\ker \chi\big|\\
&=\big[\cicl N{}:\cicl M{}\big]\big[\cicl M{}:K_1\big]=\big[\cicl N{}:K_1
\big]
\end{align*}
y $\big|\ker \tilde{\chi}\big|=\big[\cicl N{}:K_2\big]$. Se sigue que
$K_1=K_2$.

Lo anterior implica que, dado cualquier caracter de Dirichlet $\chi$ definido
m\'odulo $M$, sin importar su conductor, el campo $K_{\chi, M}
=\cicl M{}^{\ker \chi}$ depende \'unicamente de $\chi$ y no de $M$.
\end{observacion}

\begin{definicion}\label{D6.4.24} Sea $X$ cualquier grupo finito
de caracteres de Dirichlet. Sea $M$ el m{\'\i}nimo com\'un m\'ultiplo
de $\{F_{\chi}\mid\chi\in X\}$. Entonces $X$ es un subgrupo de
$\widehat{G_M}$. Sean $H:=\bigcap\limits_{\chi\in X}\ker\chi$ y
$K_X:=\cicl M{}^H$. Entonces $K_X$ se llama {\em el campo
que pertenece a $X$\index{campo que pertenece a un grupo de
caracteres}} o el {\em campo asociado a $X$\index{campo
asociado a un grupo de caracteres}}.
\end{definicion}

\begin{observacion}\label{O6.4.12(1)}
Sea $K\subseteq L\subseteq \cicl M{}$ y sea $X$ el grupo de
caracteres de Dirichlet asociado a $L$. Entonces $L=
\cicl M{}^H\subseteq 
\cicl M{}^+=\cicl M{}^{\*\F}$, donde $H:=\bigcap\limits_{\chi\in X}
\ker\chi$, si y solamente si $\*\F\subseteq H \iff \chi(a)=1$
para toda $\chi\in X$ y para toda $a\in\*\F$.

En general, si $X$ es el grupo de caracteres de Dirichlet
asociado a un campo $L\subseteq \cicl M{}$, 
entonces el grupo de caracteres de 
Dirichlet asociado a $L^+=L\cap \cicl M{}^+$ es 
$X^+=X\cap \{\chi\in\widehat{\G M}\mid \chi(a)=1
\text{\ para toda $a\in\*\F$}\}$, es decir
\[
X^+:=\{\chi\in X\mid \chi(a)=1\text{\ para toda\ } a\in \*\F\}.
\]

En particular se tiene que un caracter 
$\chi$ es par si y solamente
si $\pK_{\infty}$ se descompone totalmente 
en $K_{\chi}/K$.
\end{observacion}

Se tiene que si $X$ es c{\'\i}clico, $X=\langle\chi\rangle$, entonces
$K_X=K_{\langle\chi\rangle}$.

\begin{observacion}\label{O6.4.25} Con las notaciones anteriores,
tenemos que $H$ es subgrupo de $G_M$ y que $G_M/H\cong
\Gal(K_X/K)$. Por la Proposici\'on \ref{P12.10}, $H^{\perp}
\cong \widehat{\big(G_M/H\big)}\cong \widehat{\Gal(K_X/K)}$.
Puesto que $G_M$ es abeliano, $H^{\perp}\cong \Gal(K_X/K)$.
\end{observacion}

Tambi\'en, si $\chi\in X<\widehat{G_M}$, puesto que $\ker \chi
\supseteq H$, podemos considerar el mapeo inducido $\tilde{\chi}
\colon G_M/H\to {\ma C}^{\ast}$. Por lo tanto $X\subseteq \widehat{
G_M/H}\cong H^{\perp}$. Ahora $X^{\perp}<G_M$ y si $\alpha
\in X^{\perp}$, entonces $\chi(\alpha)=1$ para toda $\chi\in X$.
Por tanto $\alpha\in H$ y $X^{\perp}\subseteq H$ de tal forma que
$H^{\perp}\subseteq X^{\perp\perp}=X$. Se sigue que
\[
X=H^{\perp}\cong \widehat{\Gal(K_X/K)}\cong\Gal(K_X/K).
\]

Sea $X$ un grupo finito de caracteres de Dirichlet. Puesto que $X
\cong \widehat{\Gal(K_X/K)}$, podemos considerar el 
apareamiento natural
\begin{eqnarray*}
\Psi\colon \Gal(K_X/K)\times X&\longto & {\ma C}^{\ast}\\
(g,\chi)&\longmapsto& \chi(g).
\end{eqnarray*}

Bajo $\Psi$ tenemos que si $L$ es un subcampo de $K_X$, definimos
\[
Y_L=\Gal(K_X/L)^{\perp}\cong\Bigg(\widehat{\frac{\Gal(
K_X/K)}{\Gal(K_X/L)}}\Bigg)\cong \widehat{\Gal(L/K)}.
\]

Rec{\'\i}procamente, si $Y\subseteq X$ es un subgrupo de $X$,
sea $L_Y=K_X^{Y^{\perp}}$. Entonces $L_Y$ es el campo fijo
de $\{g\in\Gal(K_X/K)\mid \chi(g)=1\ \forall\ \chi\in Y\}$. Se tiene
$Y^{\perp}=\Gal(K_X/L_Y)$, as{\'\i} que $Y=Y^{\perp\perp}=
\Gal(K_X/L_Y)=Y_{L_Y}$.

Por otro lado, $L_{Y_L}=K_X^{Y^{\perp}}=K_X^{(\Gal(K_X/L)^{
\perp})^{\perp}}=K_X^{\Gal(K_X/L)}=L$. En otras palabras hemos
probado:

\begin{teorema}\label{T6.4.26} Existe una correspondencia biyectiva
entre ${\cal A}=\{Y\mid Y<X\}$ y ${\cal B}=\{L\mid L\subseteq K_X
\}$ dada por
\begin{eqnarray*}
{\cal A}&\longleftrightarrow&{\cal B}\\
Y&\longto&L_Y=K_X^{Y^{\perp}}\\
\widehat{\Gal(L/K)}\cong \Gal(K_X/L)^{\perp}=Y_L&\longleftarrow&L
\end{eqnarray*}

En particular obtenemos una correspondencia uno a uno entre todos
los subgrupos de caracteres de Dirichlet y subcampos de campos 
de funciones ciclot\'omicos. $\fin$
\end{teorema}

\begin{observacion}\label{O6.4.27}
Puesto que $\Gal(L/K)$ es un grupo finito, tenemos que 
$\Gal(L/K)\cong \widehat{\Gal(L/K)}\cong Y_L$ Esto se puede 
expresar por medio del pareo natural no degenerado
\begin{eqnarray*}
\Gal(L/K)\times Y_L\colon&\longto&{\ma C}^{\ast}\\
(g,\chi)&\longmapsto&\chi(g).
\end{eqnarray*}
\end{observacion}

Similar a la Proposici\'on \ref{P12.2.38}, tenemos

\begin{proposicion}\label{P6.4.28}
Sean $X_1$, $X_2$ dos grupos de caracteres de Dirichlet y sean
$K_i=K_{X_i}$, $i=1,2$, los campos pertenecientes a cada $X_i$.
Entonces
\las
\item $X_1\subseteq X_2$ si y solamente si $K_1\subseteq K_2$.
\item $K_{\langle X_1,X_2\rangle}=K_1K_2$. 
\item $K_{X_1\cap X_2}=K_1\cap K_2$. $\fin$
\end{list}
\end{proposicion}

\section[Caracteres de Dirichlet y aritm\'etica de campos
ciclot\'omicos]{Caracteres de Dirichlet y aritm\'etica de campos de 
funciones ciclot\'omicos}\label{S6.5}

En esta secci\'on veremos que los caracteres de Dirichlet pueden
ser aplicados para estudiar algunas propiedades aritm\'eticas de
campos de funciones ciclot\'omicos.

Sea $M\in R_T\setminus\{0\}$ m\'onico y sea $M=\prod_{i=1}^r
P_i^{\alpha_i}$ su descomposici\'on como producto de polinomios
irreducibles. Entonces
\begin{equation}\label{Ec6.5.1}
\units M{}\stackrel{\varphi}{\cong}\prod_{i=1}^r\units {P_i}{\alpha_i}
\end{equation}
con isomorfismo $\varphi$.

Si $\chi$ es un caracter de Dirichlet m\'odulo $M$, entonces 
correspondiente a (\ref{Ec6.5.1}) se tiene $\chi=\prod_{i=1}^r\chi_{
P_i}$ en donde $\chi_{P_i}$ es un caracter m\'odulo $P_i^{\alpha_i}$.
En otras palabras
\[
\chi(A\bmod M)=\prod_{i=1}^r \chi_{P_i}(A\bmod P_i^{\alpha_i}).
\]
Se tiene $\chi_{P_i}=\chi\circ \varphi^{-1}\circ g_{P_i}$ donde $g_{P_i}
\colon \units {P_i}{\alpha_i}\to \prod_{j=1}^r \units {P_j}{\alpha_j}$
est\'a dado por $g_{P_i}(A)=(1,\ldots,1,A,1,\ldots,1)$.

\begin{ejemplo}\label{Ej6.5.2}
Sea $\chi$ y $\varphi$ como en el Ejemplo \ref{Ej6.4.16}. Entonces
$\chi$ est\'a definido m\'odulo $T^2(T^2+1)$ y $\varphi$ definido
m\'odulo $T^2$. Sea $\phi:=\chi\varphi$, donde $\phi$ est\'a definido
m\'odulo $T^2+1$. Tenemos $\chi=(\chi\varphi)\varphi^{-1}=\phi
\varphi^{-1}$ y $\phi$ est\'a definido $T^2+1$ de tal forma que
$\chi_{T^2}=\varphi^{-1}=\varphi$ y $\chi_{T^2+1}=\phi$.
\end{ejemplo}

\begin{definicion}\label{D6.5.3} Sea $X$ un grupo finito de caracteres.
Entonces para un polinomio m\'onico e irreducible $P\in R_T$ 
definimos: $X_P=\{\chi_P\mid \chi\in X\}$.
\end{definicion}

\begin{ejemplo}\label{Ej6.5.4}
Sea $q=2$. Consideremos
 $M=T^2N$ con $N=T^2+T+1$ y $\theta\colon\units
M{}\to {\ma C}^{\ast}$ donde 
\begin{gather*}
\theta(1)=1,\quad \theta(T+1)=\zeta,\quad \theta(T^2+1)=\zeta^2,
\quad \theta(T^3+T^2+T+1)=-1,\\
\theta(T^3+T^2+1)=-\zeta,\quad \theta(T^3+T+1)=-\zeta^2\quad
\text{donde}\quad \zeta=\zeta_6=e^{2\pi i/6}.
\end{gather*}

Tenemos $\theta=\theta_{P_1}\theta_{P_2}$ donde $P_1=T$ y
$P_2=N$,
\begin{gather*}
\theta_{P_1}\colon \units T2\longto {\ma C}^{\ast},\quad
\theta_{P_2}\colon \units N{}\longto {\ma C}^{\ast},\\
\varphi^{-1}\circ g_{P_1}\colon \units T2\longto \units M{}
\quad\text{con}\\
1\mapsto (1,1)\mapsto 1,\quad T+1\mapsto (T+1,1)=\\
=(T^3+T^2+T+1,T^3+T^2+T+1)\mapsto T^3+T^2+T+1\\
\intertext{pues $T^3+T^2+T+1\equiv T+1\bmod T^2$,
$T^3+T^2+T+1\equiv 1\bmod N$ y}
\varphi^{-1}\circ g_{P_2}\colon \units N{}\to \units M{}\quad \text{con}\\
1\mapsto (1,1)\mapsto 1,\quad T\mapsto (1,T)=(T^2+1,T^2+1)
\mapsto T^2+1,\\
T+1\mapsto (1,T+1)=(T^3+T^2+1,T^3+T^2+1)\mapsto T^3+T^2+1
\quad\text{pues}\\
T^2+1\equiv 1\bmod T^2,\quad T^2+1\equiv T\bmod N,\quad
T^3+T^2+1\equiv 1\bmod T^2 \quad \text{y}\\
T^3+T^2+1\equiv T+1\bmod N.\\
\intertext{Luego $ \theta_{P_1}=\theta\circ\varphi^{-1}g_{P_1}$,
$ \theta_{P_2}=\theta\circ\varphi^{-1}g_{P_2}$,}
\theta_{P_1}(1)=1,\quad \theta_{P_1}(T+1)=\theta(T^3+T^2+T+1)=-1,\\
\theta_{P_2}(1)=1,\quad \theta_{P_2}(T)=\theta(T^2+1)=\omega,
\quad \theta_{P_2}(T+1)=\theta(T^3+T^2+1)=\omega^2,
\end{gather*}
donde $\omega=\zeta^2=e^{2\pi i/3}$.

Si $X=\langle \theta\rangle$, entonces $X_{P_1}=\langle\theta_{P_1}
\rangle$, $X_{P_2}=\langle\theta_{P_2}\rangle$, y $X_P=\{1\}$
si $P\notin\{T,N\}$.
\end{ejemplo}

\begin{teorema}\label{T6.5.5} Sean $X$ un grupo finito de caracteres
de Dirichlet y $K_X$ su campo asociado. Sea $P\in R_T\setminus
\{0\}$ un polinomio irreducible y sea $(P)_K=\frac{\pK}{\pK_{\infty}^{
\gr P}}$.
Sea $\pL$ un divisor primo de $K_X$ sobre $\pK$ y sea $e:=
e(\pL|\pK)$. Entonces $e=\big|X_P\big|$.
\end{teorema}

\begin{proof}
Sea $M$ el m{\'\i}nimo com\'un m\'ultiplo de $\{F_{\chi}\mid \chi\in X\}$.
Entonces $K_X\subseteq \cicl M{}$. Sea $M=P^a A$ donde $A\in
R_T$ y $P$ no divide a $A$. Sea $L=K_X(\Lam A{})=K_X\cicl A{}$.
Consideremos el siguiente diagrama donde la ramificaci\'on se
refiere a $P$.
\[
\xymatrix{
&&\cicl M{}\ar@{-}[d]\ar@{-}@/_/[ddll]\ar@{-}@/^/[rrdd]\\
&&{L=K_X(\Lam A{})=K_X\cicl A{}}
\ar@{-}[ddl]_{\substack{\hbox{\rm\tiny
no}\\\hbox{\rm\tiny ramificado}}}\ar@{-}[d]^{\substack{\hbox{\rm\tiny{no
}}\\\hbox{\rm\tiny
ramificado}}}\ar@{-}[drr]\\
{\cicl Pa}\ar@{-}[dr]&&{K_X}\ar@{-}[dd]^{e=e_p(K_X/K)}&&
\cicl A{}\ar@{-}[ddll]^{\substack{\hbox{\rm
\tiny no}\\ \hbox{\rm\tiny ramificado}}}\\
&K_{X_P}\ar@{-}[dr]_{\substack{\hbox{\rm\tiny totalmente}\\ 
\hbox{\rm\tiny ramificado}}}\\&&K}
\]

Entonces por la Proposici\'on \ref{P6.4.28} se tiene
$L=K_X\cicl A{}=K_XK_{\widehat{G_A}}=K_{\widehat{\langle X,G_A
\rangle}}$.

As{\'\i} que $L$ es el campo perteneciente al grupo generado
por $X$ y $\widehat{G_A}$, esto es, el grupo de caracteres de $L$
est\'a generado por $X$ y por los caracteres de Dirichlet de $G_M$
cuyo conductor es primo relativo a $P$.

Como en la demostraci\'on del Teorema \ref{T12.3.3} se tiene que
$\langle X,\widehat{G_A}\rangle\cong X_P\times\widehat{G_A}$.

Ahora bien $K_{X_P}\subseteq \cicl Pa$ y $L=K_{X_P}\cicl A{}$. Se
tiene que $\pK$ es no ramificado en $\cicl A{}/K$ y por tanto el
{\'\i}ndice de ramificaci\'on de $\pK$ en $K_X/K$ es el mismo que
el de $L/K$. Por otro lado, puesto que $L/K_{X_P}$ no es 
ramificado en los divisores primos que est\'an sobre $\pK$
y por el Teorema \ref{T6.2.28} (4), $\pK$ es totalmente ramificado
en $K_{X_P}/K$, concluimos que $e=\big[K_{X_P}:K\big]=\big|X_P
\big|$. $\fin$
\end{proof}

Con respecto al primo infinito, tenemos

\begin{proposicion}\label{P6.5.5+1}
Sea $\chi$ un caracter de Dirichlet y sea $L$ el campo perteneciente
a $\chi$. Sea $e_{\infty}$ el \'indice de ramificaci\'on de $\p$ en $L/K$.
Entonces $e_{\infty}=|\chi(\*\F)|$.

M\'as generalmente, sean $X$ un grupo de caracteres de Dirichlet,
$L$ el campo perteneciente a $X$ y $H:=\bigcap_{\chi\in X}
\ker \chi$. Si $e_{\infty}$ es el \'indice de ramificaci\'on de 
$\p$ en $L/K$, entonces
\begin{align*}
e_{\infty}&=\Big|\frac{\*\F}{H\cap\*\F}\Big|
=\frac{q-1}{\big|\{\xi\in\*\F\mid\chi(\xi)=1\text{\ para toda $\chi\in X$}\}\big|}\\
&=\frac{|\*\F|}{\big|\bigcap_{\chi\in X}\ker \chi_{|_{\*\F}}\big|}=\big|
\big\{\big(\chi(a)\big)_{\chi\in X}
\subseteq \big(\*\F\big)^{|X|}\mid a\in\*\F\big\}\big|.
\end{align*}
\end{proposicion}

\begin{proof}
Primero consideremos un caracter $\chi$. Entonces $\Gal(L/K)=\cicl N{}^{\ker \chi}$
donde $\chi$ est\'a definido m\'odulo $N$. Se tiene que $L^+=L\cap \cicl N{}=
L^{\ker\chi}\cap L^{\*\F}=L^{\ker \chi\cdot \*\F}$ de donde $e_{\infty}=
[L:L^+]=\Big|\frac{\ker \chi\cdot \*\F}{\ker \chi}\Big|=\Big|\frac{\*\F}{\ker\chi\cap\*\F}\Big|
=\Big|\frac{\*\F}{\ker\chi_{|_{\*\F}}}\Big|=|\chi(\*\F)|$.

Ahora consideremos $X$ un grupo de caracteres de Dirichlet m\'odulo $N$ y sea
$L$ el campo perteneciente a $X$. Entonces $e_{\infty}=[L:L^+]=\frac{|H\*\F|}
{|H|}=\frac{|\*\F|}{|H\cap \*\F|}$.

Se tiene que $H\cap\*\F=\Big(\bigcap_{\chi\in X}\ker\chi\Big)\cap \*\F=
\bigcap_{\chi\in X}\big(\ker\chi\cap\*\F)=\bigcap_{\chi\in X} \ker\chi_{|_{\*\F}}=
\{a\in\*\F\mid\chi(a)=1\text{\ para toda $\chi\in X$}\}$.

Sea $S=(\*\F)^{|X|}=\{(b_{\chi})_{\chi\in X}\mid b_{\chi}\in\*\F\}$. Sea 
$\varphi\colon\*\F\lra S$ dada por $\varphi(a)=\big(\chi(a)\big)_{\chi\in X}$.
Entonces $\varphi$ es un homomorfismo y $\ker \varphi=H\cap \*\F$ por lo
que $|\im\varphi|=\frac{|\*\F|}{|H\cap\*\F|}=e_{\infty}$.
$\fin$
\end{proof}

\begin{ejemplo}\label{Ej6.5.6} En el Ejemplo \ref{Ej6.5.4}, el
{\'\i}ndice de ramificaci\'on de $P_1=T$ es $\big|X_{P_1}\big|=2$
pues $X_{P_1}=\langle \theta_{P_1}\rangle$ y $\theta_{P_1}$ es
de orden $2$ y el {\'\i}ndice de ramificaci\'on de $P_2=T^2+T+1$
es $\big|X_{P_2}\big|=3$ pues $X_{P_2}=\langle\theta_{P_2}\rangle$
y $\theta_{P_2}$ es de orden $3$. Finalmente para otro primo
distinto a $P_1$, $P_2$ y $\pK_{\infty}$ se tiene que el {\'\i}ndice
de ramificaci\'on es $1$, es decir, es no ramificado.
\end{ejemplo}

\begin{ejemplo}\label{Ej6.5.7} Sea $q=2$ y consideremos el
caracter $\chi$ dado en el Ejemplo \ref{Ej6.4.16}. El conductor de
$\chi$ es $T^2(T^2+1)$. Por el Ejemplo \ref{Ej6.5.2} tenemos que
$\chi_{T^2}=\varphi$ y $\chi_{T^2+1}=\phi$. Notemos que $\Phi(T^2)
=\Phi(T^2+1)=q^{dn}-q^{d(n-1)}=2^{1(2)}-2^{1(2-1)}=2^2-2=4-2=2$.
De aqu{\'\i} que $\big[\cicl T2:K\big]=\big[\cicl {T^2+1}{}:K\big]=2$.
Tenemos
\begin{gather*}
u^{T^2}=\sum_{i=0}^2\carlitzbinom {T^2}i u^{q^i}=T^2 u+
\carlitzbinom {T^2}1 u^q+uq^2.\\
\intertext{Ahora bien,}
\carlitzbinom {T^2}1=T\carlitzbinom T1+\carlitzbinom T0^q=T+T^q
=T+T^2\\
\intertext{donde $\carlitzbinom T1=a_1=1$. Por lo tanto}
u^{T^2}=T^2+(T+T^2)u^2+u^4.\\
\intertext{Tambi\'en tenemos}
\Psi_{T^2}(u)=\frac{u^{T^2}}{u^T}=\frac{T^2u+(T+T^2)u^2+u^4}
{Tu+u^2}=u^2+Tu+T.
\end{gather*}

Por lo tanto cada ra{\'\i}z $\alpha$ de $\Psi_{T^2}(u)$ es de la
forma: $\big(\frac{\alpha}{T}\big)^2+\big(\frac{\alpha}{T}\big)=
-\frac{1}{T}=\frac{1}{T}$. Por lo tanto $\cicl T2=K(\beta)$ donde
$\beta$ es una ra{\'\i}z de la extensi\'on de Artin-Schreier que 
satisface $\beta^2-\beta=\frac{1}{T}$.

Similarmente, $\cicl {T^2+1}{}=K(\gamma)$ donde $\gamma^2+
\gamma=\frac{1}{T+1}$. Se sigue que $K_{\chi}=K(\varepsilon)$ con
$\varepsilon^2-\varepsilon=\frac{1}{T(T+1)}$ y tenemos el 
siguiente diagrama
\[
\xymatrix{
&\cicl {T^2(T+1)}{}=K(\beta,\gamma)\ar@{-}[dl]\ar@{-}[d]
\ar@{-}[dr]\\
K(\beta)\ar@{-}[dr]_{\chi_{T^2}}&K(\varepsilon)\ar@{-}[d]_{\chi}
&K(\gamma)\ar@{-}[dl]^{\chi_{T^2+1}}\\ &K
}
\]

En $K(\varepsilon)/K$, $T$ y $T+1$ son los primos ramificados.
En $K(\beta)/K$, $T$ es el \'unico primo ramificado y en $K(\gamma)
/K$, $T+1$ es el \'unico primo ramificado.
\end{ejemplo}

\begin{corolario}\label{C6.5.8}
Sea $\chi$ un caracter de Dirichlet. Entonces $P$ se ramifica en
$K_{\chi}/K$ si y s\'olo si $\chi(P)=0$, o equivalentemente, $P$
divide a $F_{\chi}$. Si $X$ es cualquier grupo finito de caracteres
de Dirichlet, entonces $P$ es no ramificado en $K_X/K$ si y s\'olo
si $\chi(P)\neq 0$ para toda $\chi\in X$.
\end{corolario}

\begin{proof}
Se tienen las equivalencias: $P$ es ramificado en $K_X/K\iff
X_P\neq\{1\}\iff$ existe $\chi\in X$ tal que $\chi_P\neq 1\iff$
existe $\chi\in X$ con $P|F_{\chi}\iff$ existe $\chi\in X$ con
$\chi(P)=0$. $\fin$
\end{proof}

Como en el caso num\'erico, Teorema \ref{T12.3.5}, los grupos
de inercia y de descomposici\'on est\'an relacionados con los
caracteres de Dirichlet de la siguiente manera.

\begin{teorema}\label{T6.5.9}
Sea $X$ un grupo finito de caracteres de Dirichlet y sea $K_X$ su
campo asociado. Sean $P\in R_T$ y $Y=\{\chi\in X\mid \chi(P)\neq
0\}$, $Z=\{\chi\in X\mid \chi(P)=1\}$. Entonces, si $(P)_K=\frac{\pK}
{\pK_{\infty}^{\gr P}}$, consideremos $\pL$ un divisor primo en 
$K_X$ sobre $\pK$. Entonces
\[
X/Y\cong\widehat{I(\pL|\pK)}\cong I(\pL|\pK)\quad\text{y}\quad
X/Z\cong D(\pL|\pK).
\]

En particular, $e=e(\pL|\pK)=[X:Y]$, $f=d(\pL|\pK)=[Y:Z]$ y $h=
[Z:1]=|Z|$ donde $h$ es el n\'umero de divisores primos en $K_X$
sobre $\pK$. Finalmente, el grupo $Y/Z$ es c{\'\i}clico de orden $f$.
\end{teorema}

\begin{proof}
An\'aloga a la del Teorema \ref{T12.3.5}. $\fin$
\end{proof}

Resolveremos el problema inverso de la Teor{\'\i}a de Galois para el
caso particular de un grupo abeliano. Primero necesitamos el 
siguiente resultado.

\begin{proposicion}\label{P6.5.10} Sea $P\in R_T$ un polinomio
irreducible m\'onico de grado $d$ y sea $n=p^t$. Entonces
$\units Pn$ contiene un subgrupo c{\'\i}clico de orden $p^t a$ para
cualquier $a$ que divide a $q^d-1$.
\end{proposicion}

\begin{proof}
Se tiene que $\big|\units Pn\big|=\Phi(P^n)=q^{dn}-q^{d(n-1)}=
q^{d(n-1)}(q^d-1)$.

Por tanto $\units Pn$ es isomorfo a una suma directa $H\oplus A$
donde $|H|=q^{d(n-1)}$ y $|A|=q^d-1$. Notemos que $A$ es el
\'unico subgrupo de orden $q^d-1$. Definimos
\begin{eqnarray*}
\theta\colon \units Pn & \longto & \units P{}\\
B \bmod P^n & \longmapsto & B\bmod P.
\end{eqnarray*}

Entonces $\theta$ es un epimorfismo y $\units Pn/\ker \theta \cong
\units P{}$.

Puesto que $\big|\units P{}\big|=\Phi(P)=q^d-1$, se sigue que
\[
A\cong\units P{}\quad\text{y}\quad H\cong \ker \theta\cong\{B\bmod
P^n\mid B\equiv 1\bmod P\}.
\]
Ahora $R_T/\langle P\rangle$ y ${\ma F}_{q^d}$ son isomorfos de
tal forma que $A$ es el grupo multiplicativo de los elementos 
diferentes de cero de un campo y por lo tanto $A$ es un grupo
c{\'\i}clico.

Sea $B=1+P$. Queremos determinar el orden de $B$ m\'odulo $P^n$
en $R_T/\langle P^n\rangle$. Ahora como $B\in \ker \theta$, $B\in
H$, y $o(B)=p^s$ para alg\'un $s\geq 0$. Entonces
\[
B^{p^s}=1+P^{p^s}\equiv 1\bmod P^n\iff p^s\geq n=p^t\iff s\geq t.
\]
Por tanto $o(B)=p^t$. $\fin$
\end{proof}

\begin{teorema}\label{T6.5.11} Sea $G$ un grupo abeliano finito. 
Entonces existen campos de funciones congruentes $E$ y $F$ tales
que
\l
\item $\Gal(F/E)\cong G$.
\item $F/E$ es no ramificada en todos los divisores primos.
\item $F/K$ es abeliana y $E/K$ es c{\'\i}clica.
\item El campo de constantes, tanto de $E$ como de $F$ es
${\ma F}_q$.
\end{list}
\end{teorema}

\begin{proof}
Sea $G\cong {\ma Z}/m_1{\ma Z}\times \cdots\times {\ma Z}/m_r
{\ma Z}$. Definimos
 $m_i=p^{t_i}a_i$ con $\mcd (a_i,p)=1$, $t_i\geq 0$
para $1\leq i\leq r$. Sea $d_i'=o(p\bmod a_i)$, es decir, $p^{d_i'}
\equiv 1\bmod a_i$, con $d_i'$ m{\'\i}nimo. Escojamos n\'umeros
naturales $d_1<d_2<\cdots <d_r$ donde cada $d_i'$ divide a
$d_i$. Por ejemplo, podemos tomar $d_1=d_1'$, $d_i=2d_{i-1}d_i'$,
$i=2,\ldots,r$. Sea $P_i\in R_T$ un polinomio m\'onico irreducible
de grado $d_i$. Tal polinomio $P_i$ existe pues si ${\ma F}_{q^{d_i}}
={\ma F}_q(\alpha_i)$ para alg\'un $\alpha_i$, entonces $P_i=\Irr(
\alpha_i,T,{\ma F}_q)$ es de grado $d_i$ y ${\ma F}(\alpha_i)\cong
R_T/\langle P_i\rangle$.

Por la Proposici\'on \ref{P6.5.10}, $\units {P_i}{p^{t_i}}$ contiene un
elemento de orden $p^{t_i}a_i=m_i$. Ahora, puesto que $\widehat{
\units {P_i}{p^{t_i}}}\cong \units {P_i}{p^{t_i}}$, existe un caracter
$\chi\bmod P_i^{t_i}$ de orden $m_i$. Es decir $\chi$ satisface
$o(\chi)=m_i$ y $F_{\chi_i}=p_i^{s_i}$ con $s_i\leq t_i$.

Sea $P_{r+1}$ otro polinomio m\'onico irreducible de grado $d_{r+1}
>d_r$ tal que $a_1\cdots a_r| q^{d_{r+1}}-1$. Tal $d_{r+1}$ existe
pues $\mcd (a_1\cdots a_r, q)=1$.

Sea $\chi_{r+1}$ un caracter de Dirichlet definido m\'odulo $P_{r+1}^{
p^t}$ para $t=t_1+\cdots +t_r$ y orden $m_{r+1}=p^t(q^{d_{r+1}}-1)$
(Proposici\'on \ref{P6.5.10}). Entonces
\[
m_1\cdots m_r=a_1\cdots a_r p^{t_1+\cdots t_r}|m_{r+1}.
\]
Sea $\chi:=\chi_1\cdots\chi_r\chi_{r+1}$ y $E:=K_X$ el campo 
asociado a $X:=\langle\chi\rangle$. Sea $Y:=\langle \chi_1,\ldots,
\chi_r,\chi_{r+1}\rangle$ y $F:=K_Y$ el campo correspondiente a $Y$.
Tenemos
\[
K\subseteq E=K_X\subseteq K_Y=F\subseteq \cicl M{}
\]
donde $M=P_1^{\alpha_1}\cdots P_r^{\alpha_r}P_{r+1}^{\alpha_{r+1}}
$ con $\alpha_i=p^{t_i}$, $1\leq i\leq r$ y $\alpha_{r+1}=p^t$. 
En particular el campo de constantes de $E$ y $F$ es ${\ma F}_q$
(Corolario \ref{C6.2.29}). Esto prueba ({\sc{iv}}) y tambi\'en tenemos
que $F/K$ es una extensi\'on abeliana.

Por otro lado se tiene $\Gal(E/K)\cong X\cong \langle\chi\rangle$
es c{\'\i}clico y obtenemos ({\sc iii}). Ahora bien, $Y=\langle \chi_1,
\ldots,\chi_r,\chi_{r+1}\rangle=\langle\chi_1,\ldots,\chi_r,\chi\rangle$
y $o(\chi)=o(\chi_{r+1})=m_{r+1}$ y puesto que $m_1\cdots m_r$ 
divide a $m_{r+1}$, $\chi$ es de orden maximal en  $Y$.

Por tanto $Y/X=Y/\langle\chi\rangle\cong \langle\chi_1,\ldots,\chi_r
\rangle$ y
\begin{align*}
Y/X&\cong \frac{\widehat{\Gal(K_Y/K)}}{\widehat{\Gal(K_X/K)}}
\cong \frac{\widehat{\Gal(K_Y/K)}}{\widehat{\Big(\frac{
\Gal(K_Y/K)}{\Gal(K_Y/K_X)}\Big)}}\cong \widehat{\Gal(K_Y/K_X)}\\
&\cong \widehat{\Gal(F/E)}\cong \Gal(F/E).
\end{align*}
Se obtiene que $\Gal(F/E)\cong\langle\chi_1,\ldots,\chi_r\rangle \cong
{\ma Z}/m_1{\ma Z}\times \cdots \times {\ma Z}/m_r{\ma Z}\cong
G$ lo cual prueba ({\sc i}).

Por el Teorema \ref{T6.3.5} se tiene que los primos ramificados en
$F/K$ son $\pK_1,\ldots,\pK_r, \pK_{r+1}$ y $\pK_{\infty}$, donde
$(P_i)_K=\frac{\pK_i}{\pK^{\gr P_i}_{\infty}}$.

Ahora bien el {\'\i}ndice de ramificaci\'on de $\pK_{\infty}$ en
$E/K$ es $q-1$ lo mismo que en $F/K$ puesto que $E$ es el
campo perteneciente a $\chi$, $q-1|o(\chi)$ y $\units {P_{r+1}}{p^t}$
contiene un \'unico subgrupo de orden $(q-1)$ (Proposici\'on 
\ref{P6.5.10}) y este grupo es el grupo de inercia $\pK_{\infty}$
(Teorema \ref{T6.2.30}). Por lo tanto $\pK_{\infty}$ es no ramificado
en $F/E$. Finalmente tenemos $Y_{P_i}=\langle \chi_i\rangle=X_{
P_i}$. Por el Teorema \ref{T6.5.5} se tiene que {\'\i}ndice de 
ramificaci\'on en cada divisor en $F$ que divide a $\pK_i$ es
$\frac{\big|Y_{P_i}\big|}{\big|X_{P_i}\big|}=1$. Por tanto $F/E$ es
no ramificada en cada divisor primo y esto prueba ({\sc ii}) y
el teorema. $\fin$
\end{proof}

Sea $E=K(\sqrt[t]{\gamma D})/K$ una extensi\'on de Kummer, donde
$t|q-1$, $D\in R_T$ es un polinomio m\'onico de grado positivo
y libre de $t$--potencias. ?`Se cumple que $E\subseteq \cicl D{}$?

La respuesta es que $E\subseteq \cicl D{}$ si y solamente si
$\gamma \equiv (-1)^{\deg D}\bmod (\F)^t$. Para ver esto, probamos
el siguiente resultado.

\begin{proposicion}\label{P5.1.1} Para un polinomio $P\in R_T^+$
de grado $d$, se tiene que
$K(\sqrt[t]{(-1)^d P})\subseteq K(\Lambda_P)$,
donde $t$ es cualquier divisor de $q-1$.
\end{proposicion}

\begin{proof}
Sea $\Phi_P(u)=\frac{u^P}{u}$ el $P$--\'esimo polinomio ciclot\'omico.
Tenemos
\[
\Phi_P(u)=\prod_{\substack{A\neq 0,A\in R_T\\
\deg A<\deg P}} (u-\lambda^A)=\sum_{i=0}^d \xbinom{P}{i} 
u^{q^i-1},
\]
 donde $\lambda\in \Lambda_P\setminus \{0\}$,
esto es, $\lambda$ es un generador como $R_T$--m\'odulo de
$\Lambda_P$. Entonces 
\[
\Phi_P(0)= (-1)^{q^d-1}\prod_{\substack{A\neq 0,A\in R_T\\
\deg A<\deg P}}\lambda^A=P.
\]

Ahora, todo polinomio $A\in R_T$, $A\neq 0$ puede ser un{\'\i}vocamente
escrito como producto de un elemento $\alpha\in {\ma F}_q^{\ast}$ y
un polinomio m\'onico $A_1$: $A=\alpha A_1$.
Ahora, $\lambda^A=\lambda^{\alpha
A_1}=\alpha \lambda^A$. Notemos que hay exactamente $q-1$
polinomios $A\in R_T$, $A\neq 0$ tal que $A_1$ aparece,
una para cada uno de los $q-1$ elementos de ${\ma F}_q^{\ast}$.
Por lo tanto
\begin{align*}
P&=(-1)^{q^d-1}\prod_{\substack{A\neq 0,A\in R_T\\
\deg A<\deg P}} \lambda^A=(-1)^{q^d-1}
\prod_{\substack{A_1\text{\ m\'onico}\\
\deg A_1<\deg P, \alpha\in{\ma F}_q^{\ast}}}\alpha\lambda^{A_1} \\
&=(-1)^{q^d-1}\Big(\prod_{\alpha\in{\ma F}_q^{\ast}}\alpha
\Big)^{\frac{q^d-1}{q-1}}\Big(\prod_{\substack{A_1\text{\ m\'onico}
\\ \deg A_1<\deg P}}
\lambda^{A_1}\Big)^{q-1}.
\end{align*}

Notemos que $\prod_{\alpha\in{\ma F}_q^{\ast}}\alpha=-1$ lo cual
se sigue de que
$\frac{x^q-x}{x}=x^{q-1}-1=\prod\limits_{\alpha\in{\ma F}_q^{\ast}}
(x-\alpha)$ por lo que $\prod\limits_{\alpha\in{\ma F}_q^{\ast}}
\alpha=-1$,
y que
$\xi:=\prod_{A_1\text{\ m\'onico}}\lambda^{A_1}\in K(\Lambda_P)$.
Entonces 
\[
(-1)^{q^d-1}(-1)^{(q^d-1)/(q-1)}\xi^{q-1}=(-1)^d\xi^{q-1} 
=P,
\]
con $\xi\in K(\Lambda_P)$. Se sigue que $\xi
=\sqrt[q-1]{(-1)^dP}\in K(\Lambda_P)$. En particular
$\sqrt[t]{(-1)^d P}=\xi^{(q-1)/t}\in K(\Lambda_P)$. $\fin$
\end{proof}

\begin{corolario}\label{C5.1.1'} 
Sea $D\in R_T$ un polinomio m\'onico. Entonces
\[
K(\sqrt[t]{(-1)^{\deg D} D})\subseteq K(\Lambda_D),
\]
donde $t$ es cualquier divisor de $q-1$. $\fin$
\end{corolario}

\subsection{S\'imbolos de Legendre,\index{Legendre!s\'imbolo de $\sim$}
\index{simbolo de Legendre@s\'imbolo de Legendre}
de Jacobi\index{simbolo de Jacobi@s\'imbolo
de Jacobi}\index{Jacobi!s\'imbolo de $\sim$} y ley de 
reciprocidad\index{ley de reciprocidad}\index{reciprocidad!ley
de $\sim$}}\label{SubLegendre}

Consideremos $P\in R_T^+$ con $\deg P=d$ y $t|(q-1)$.
Se tiene que $\cicl P{}/K$ es una extensi\'on c{\'\i}clica
de grado $q^d-1$. Por tanto existe un \'unico subcampo de $\cicl P{}$
de grado $q-1$ sobre $K$: $[L:K]=q-1$. Ahora bien, puesto que
${\ma F}_q^{\ast}\subseteq K$, se tiene que las $q-1$ ra{\'\i}ces
de unidad est\'an en $K$ y por tanto $L/K$ es una extensi\'on de
Kummer. Sea $L=K(\alpha)$ donde $\alpha^{q-1}=\beta\in K$. 
Ahora puesto que los primos ramificados en $L/K$ son $\pK$
y posiblemente $\pK_{\infty}$, donde $(P)_K
=\frac{\pK}{\pK_{\infty}^{\gr P}}$ y $\pK$
es totalmente ramificado como consecuencia de los
Teoremas \ref{T6.2.28} y \ref{T6.2.30}.
Se tiene $\beta\in R_T$ y $\beta=\gamma P^i$ con $\gamma\in
{\ma F}_q^{\ast}$ y $\mcd (i,q-1)=1$.

Por Teor{\'\i}a de Kummer, podemos tomar $i=1$ y se tiene $\alpha^{
q-1}=\gamma_1 P$, para alg\'un $\gamma_1
\in \F$. Si $\zeta$ es un generador de ${\ma F}_q^{\ast}$,
equivalentemente, $\zeta$ es una $q-1$ ra{\'\i}z primitiva de $1$, 
se tiene que si $G=\Gal(L/K)$, entonces $\sigma\alpha=\zeta\alpha$
para alg\'un generador $\sigma$ de $G$: $G=\langle \sigma\rangle$.

As{\'\i}, tenemos $L=K\big(\sqrt[q-1]{(-1)^dP}\big)$ con $d=\gr P$.
Ahora queremos es describir el caracter de Dirichlet
asociado a $L=K\big(\sqrt[t]{(-1)^d P}\big)\subseteq \cicl P{}$,
m\'as precisamente, el grupo de caracteres de Dirichlet 
asociado al campo $L$.

Sea $\theta\colon \units P{}\to {\ma C}^{\ast}$ este caracter, 
es decir, $\cicl P{}^{\ker \theta}=L$. Se tiene $F_{\theta}=P$.
Sea $X=\langle\theta\rangle$ y $|X|=t=o(\theta)=
[L:K]$. Entonces $\theta^t=1$, $\theta\neq \Id$ y $\theta\Big(
\G P\Big)=W_t=\{\zeta\in\*{\ma C}\mid \zeta^t=1\}$.

Se tiene $L=\cicl P{}^{\ker \theta}$, $\ker \theta=\{\sigma
\in G_P=\Gal(\cicl P{}/K)\mid \theta(\sigma)=1\}$.

Sea $Q$ un polinomio irreducible tal que $Q\neq P$. Se tiene
que $Q$ se descompone totalmente en $L/K$ si y s\'olo si
$Z=\{\chi\in X\mid \chi(Q)=1\}$ satisface que $|Z|=t$, es decir,
si s\'olo si $Z=X$ (ver Teorema \ref{T6.5.9}). Por tanto, $Q$
se descompone totalmente en $L/K$ si y s\'olo si $\theta(Q)=1$.

Consideremos $R_T/\langle P\rangle\cong {\ma F}_{q^d}$.
Los elementos de $\*\F\subseteq \*{{\ma F}_{q^d}}$ satisfacen
$\*\F=\{\alpha\in \*{{\ma F}_{q^d}}\mid
\alpha^{q-1}=1\}$. En particular, si $\alpha^t=1$, entonces
$\alpha\in\*\F$.

Se tiene que $\alpha=Q^{\frac{q^d-1}{t}}\bmod P\in
\G P$ y $\alpha^t=1$ por lo que $Q^{\frac{q^d-1}{t}}\bmod
P\in\*\F$.

\begin{definicion}\label{D6.5'.13}
Dado $P\in R_T^+$ de grado $d$, se define el {\em s\'imbolo
de Legendre\index{simbolo de Legendre@s\'imbolo de Legendre}
\index{Legendre!s\'imbolo de $\sim$}} 
o {\em s\'imbolo residual de las $t$ potencias
\index{simbolo residual de potencias@s\'imbolo
residual de potencias}} $\xbinom {}P_t$ por 
\[
\xbinom {}P_t\colon \G P\lra \*\F,\quad
\xbinom{M}{P}_t:=M^{\frac{q^d-1}{t}}\bmod P
\]
para $M\in R_T$ con $P\nmid M$. Si $P|M$ definimos
$\xbinom MP_t=0$.
\end{definicion}

Sea $M\in R_T$, $P\nmid M$. Si $M=B^t\bmod P$ para
alguna $B\in R_T$, se tiene $\theta(M)=\theta(B^t)=
\theta(B)^t=1$ por lo que $M\in\ker \theta$. Puesto que
\[
|\ker \theta|=\frac{\Big|\G P\Big|}{t}=\frac{q^d-1}{t}=
\Big|\Big(\G P\Big)^t\Big|,
\]
se sigue que $\ker \theta =\Big|\Big(\G P\Big)^t\Big|$. Esto es,
$\theta(M)=1$ si y s\'olo si existe $B\in R_T$ tal que $M=
B^t\bmod P$ si y s\'olo si $\xbinom MP_t=1$.

Por tanto $\ker \theta=\ker \xbinom {}P_t$. Se sigue que
$L=\cicl P{}^{\ker \artin {}P_t}$. Hemos obtenido:

\begin{proposicion}\label{P6.5'.14}
El caracter asociado a $L=K\big(\sqrt[t]{(-1)^d P}\big)$ es el
s\'imbolo de Legendre $\xbinom {}P_t$. $\fin$
\end{proposicion}

\begin{observacion}\label{O6.5'.15}
Recordemos que hemos hecho un abuso del lenguaje. Cuando
decimos ``caracter asociado $\chi$'', en realidad estamos 
diciendo ``el grupo de carcteres asociado al campo $L$ es el
grupo generado por $\chi$''.
\end{observacion}

En ese sentido, es el que estamos usando que el s\'imbolo de
Legendre $\xbinom {}P_t$ es el caracter asociado a
$L=K\big(\sqrt[t]{(-1)^d P}\big)$. En este mismo sentido,
tenemos que si $\beta\in{\ma N}$ es tal que $\mcd(t,\beta)
=1$, entonces 
\[
L=K\Big(\big(\sqrt[t]{(-1)^dP}\big)^{\beta}\Big)=L\Big(
\sqrt[t]{(-1)^{d\beta}P^\beta}\Big)=L\Big(\sqrt[t]{(-1)^{\deg P^{
\beta}}P^{\beta}}\Big)
\]
y $\Big\langle\xbinom {}P_t\Big\rangle=\Big\langle\xbinom
{}P_t^{\beta}\Big\rangle$ por lo que podemos decir que
$\xbinom {}P_t^{\beta}$ es el caracter asociado a $L=K\Big(
\sqrt[t]{(-1)^{\deg P^{\beta}}P^{\beta}}\Big)$.

Ahora, sea $\gamma\in{\ma N}$ tal que $\mcd(t,\gamma)=s$.
Entonces si $\gamma=sa$, $\mcd\big(\frac{t}{s},a\big)=1$ y
$\sqrt[t]{(-1)^{\deg P^{\gamma}} P^{\gamma}}=
\sqrt[t]{(-1)^{\deg P^{as}} P^{as}}=\sqrt[t/s]{
(-1)^{\deg P^{a}} P^{a}}$. Por tanto se tiene que 
\[
K\Big(\sqrt[t]{(-1)^{\deg 
P^{\gamma}} P^{\gamma}}\Big)=K\Big(\sqrt[t/s]{(-1)^{\deg P^{a}} 
P^{a}}\Big)=K(\Big(\sqrt[t/s]{(-1)^{\deg P} P}\Big).
\]

Sea $D=P_1^{\alpha_1}\cdots P_r^{\alpha_r}$ con $P_1,\ldots,
P_r\in R_T^+$ primos distintos, $r\geq 1$ y $1\leq\alpha_i\leq
t-1$, $1\leq i\leq r$. Sea $\alpha_i=a_ic_i$ con $a_i=\mcd(
\alpha_i,t)$. Entonces $\mcd(c_i,t/a_i)=\mcd\big(\frac{\alpha_i}
{a_i},\frac{t}{a_i}\big)=1$. Sea $L=K\Big(\sqrt[t]{(-1)^{\deg D}D}\Big)$
y sea $\chi_D$ el caracter de Dirichlet asociado a $L$. 
Se tiene
\begin{gather}\label{Ec9.5'.1}
\sqrt[t]{(-1)^{\deg D}D}=\prod_{i=1}^r \sqrt[t]{(-1)^{\deg P_i^{\alpha_i}}
P^{\alpha_i}}=\prod_{i=1}^r \sqrt[t/a_i]{(-1)^{\deg P_i^{c_i}}P^{c_i}}.
\end{gather}

Los primos finitos ramificados son $P_1,\ldots,P_r$ pues, por un
lado $L\subseteq \prod_{i=1}^r K\big(\sqrt[t]{(-1)^{\deg P_i^{\alpha_i}}
P^{\alpha_i}}\big)\subseteq \prod_{i=1}^r \cicl {P_i}{}$ de lo cual
obtenemos que si $Q\in R_T^+$ es tal que $Q\notin\{P_1,\ldots,P_r\}$,
entonces $Q$ es no ramificado en $L/K$. Por otro lado si $\pL$
es un primo en $L$ sobre $P_i$, se tiene
\begin{align*}
v_{\pL}(D)&=v_{\pL}\Big(\big(\sqrt[t]{(-1)^{\deg D}D}\big)^t\Big)=
tv_{\pL}\big(\sqrt[t]{(-1)^{\deg D}D}\big)\\
&=e_{L|K}(\pL|P_i)
v_{P_i}(D)=e_{L|K}(\pL|P_i)\alpha_i.
\end{align*}
Por tanto $t|e_{L|K}(\pL|P_i)\alpha_i$ y $1\leq \alpha_i\leq t-1$
de donde se sigue que $e_{L|K}(\pL|P_i)>1$. As\'i, $P_1,\ldots,
P_r$ son ramificados. De la igualdad anterior es f\'acil deducir
que el \'indice de ramificaci\'on de $P_i$ en $L/K$ es $e_{L|K}(
\pL|P_i)=e_{P_i}(L|K)=\frac{t}{a_i}=\frac{t}{\mcd(\alpha_i,t)}$
(ver el Teorema \ref{TRam1}). De esta forma, como consecuencia
de la igualdad (\ref{Ec9.5'.1}), obtenemos que $\chi_D=\d\prod_{i=1}^r
\xbinom{}{P_i}^{c_i}_{t/a_i}$.

Observemos en general que, para $P\in R_T^+$, $n\in{\ma N}$ y
$m|n$, se tiene
\[
\xbinom{}{P}_{n/m}=\xbinom{}{P}^m_n,
\]
pues $\xbinom MP_{n/m}=M^{\frac{q^d-1}{n/m}}\bmod P=\Big(
M^{\frac{q^d-1}{n}}\bmod P\Big)^m=\xbinom MP_n^m$, donde
$d=\deg P$.
En particular $\xbinom{}{P}_t=\xbinom{}{P}^{(q-1)/t}_{q-1}$.  Ponemos
$\xbinom{}{P}:=\xbinom{}{P}_{q-1}$. Se sigue que $\chi_D=\d\prod_{
i=1}^r\xbinom{}{P_i}^{\alpha_i}$.

\begin{definicion}\label{D9.5'.15}
Sean $P_1,\ldots,P_r\in R_T^+$ polinomios irreducibles distinto,
$r\geq 1$. Sea $D=P_1^{\alpha_1}\cdots P_r^{\alpha_r}$ con $1\leq
\alpha_i\leq t-1$, $1\leq i\leq r$ y $t|(q-1)$. Se define el {\em
s\'imbolo de Jacobi\index{Jacobi!s\'imbolo de $\sim$}
\index{simbolo de Jacobi@s\'imbolo de Jacobi}} por
\[
\xbinom {}D_t=\prod_{i=1}^r\xbinom{}{P_i}^{\alpha_i}.
\]

Se tiene que $\xbinom {}D_t$ es el caracter de Dirichlet asociado al
campo $L=K\big(\sqrt[t]{(-1)^{\deg D}D}\big)$.
\end{definicion}

\begin{teorema}\label{T9.5'.16}
Sea $D\in R_T$ como antes y $Q\in R_T^+$ con $Q\nmid M$. Entonces
el grado de inercia de $Q$ en $L=K\big(\sqrt[t]{(-1)^{\deg D}D}\big)$ es el
orden de $\xbinom QD_t$.

En particular $Q$ se descompone
totalmente en $L/K$ si y solamente si $\xbinom QP_{t}=1$.
\end{teorema}

\begin{proof}
Sea $X=\Big\langle\xbinom {}D_t\Big\rangle$ el grupo de caracteres de
Dirichlet asociado al campo $L$. Entonces, por el Teorema \ref{T6.5.9},
el grado de inercia de $Q$ en $L$ es $\big|X/Z\big|$ donde $Z=\{
\chi\in X\mid \chi(Q)=1\}$. Se tiene
\begin{gather*}
\xbinom QD_t^s=1\iff s\big|o\Big(\xbinom QD_t\Big)\quad\text{y}\\
\xbinom QD_t^s=1\iff \xbinom QD_t^s\in Z.
\end{gather*}
Se sigue que $X/Z\cong \Big\langle\xbinom QD_t\Big\rangle\cong
D_{L|K}({\eu Q}|Q)$ donde ${\eu Q}$ es un primo en $L$ sobre $Q$ y
$D_{L|K}({\eu Q}|Q)$ es el grupo de descomposici\'on.
$\fin$
\end{proof}

Terminamos esta secci\'on con la ley de reciprocidad.

\begin{teorema}[Ley de Reciprocidad\index{reciprocidad!ley de $\sim$}\index{ley
de reciprocidad}]\label{T9.5'.17}
Sean $P,Q\in R_T^+$ con $P\neq Q$. Entonces, si $t|(q-1)$, se tiene
\[
\xbinom PQ_t\xbinom QP_t=(-1)^{\frac{q-1}{t}(\deg P)(\deg Q)}.
\]
\end{teorema}

\begin{proof}
La demostraci\'on que presentamos se debe a L. Carlitz y puede encontrarse
en \cite[Theorem 3.3]{Ros2002}.

Primero notemos que basta probar el resultado para $t=q-1$ pues en general
$\xbinom {}P_t=\xbinom {}P_{q-1}^{(q-1)/t}$. 

Consideremos ${\ma F}_{q^m}$
el campo finito con $m=\deg P\cdot \deg Q$. Si $a$ es ra\'iz de $P$, esto es,
$P(a)=0$, puesto que $P(T)\in R_T$ y por tanto $(P(T))^q=P(T^q)$, se sigue
que $P(a^q)=0$. Ahora bien, debido a que $\F(a)={\ma F}_{q^{\deg P}}$, se
tiene que $a^{q^u}\neq a$ para $u<\deg P$ y por tanto $\{a,a^q,\ldots, a^{
\deg P-1}\}$ son todas las ra\'ices de $P$ . Se sigue que $P(T)=(T-a)
(T-a^q)\cdots(T-a^{q^{\deg P-1}})$.
Similarmente, si $b$ es una ra\'iz de $Q(T)$, se tiene $Q(T)=(T-b)
(T-q^q)\cdots(T-b^{q^{\deg Q-1}})$.

Para $f(T)\in {\ma F}_{q^m}[T]$ se tiene $f(T)\equiv f(a)\bmod (T-a)$. Ahora
bien, 
\begin{gather*}
\xbinom QP_{q-1}=Q^{\frac{q^{\deg P}-1}{q-1}}\bmod P=
Q^{1+q+\cdots+q^{\deg P-1}}\bmod P\quad \text{y} \\
\begin{align*}
Q(T)^{1+q+\cdots+q^{\deg P-1}}&
\equiv Q(T)Q(T^q)\cdots Q(T^{q^{\deg P-1}})\\
&\equiv Q(a)Q(a^q)\cdots
Q(a^{q^{\deg P-1}})\bmod (T-a).
\end{align*}
\end{gather*}

Aplicando la igualdad anterior a $T-a^{q^i}$, $1\leq i\leq \deg P-1$, se obtiene
\begin{align*}
\xbinom QP\equiv Q(T)^{1+q+\cdots+q^{\deg P-1}}\bmod P&\equiv Q(a)Q(a^q)\cdots
Q(a^{q^{\deg P-1}})\bmod P\\
&\equiv \prod_{i=0}^{\deg P-1}\prod_{j=0}^{\deg Q-1} (a^{q^i}-b^{q^j})\bmod P.
\end{align*}
Puesto que $\xbinom QP$ y $\prod_{i=0}^{\deg P-1}\prod_{j=0}^{\deg Q-1} 
(a^{q^i}-b^{q^j})\bmod P$ son elementos de ${\ma F}_{q^m}$, se sigue 
la igualdad. Por tanto
\begin{align*}
\xbinom QP&=\prod_{i=0}^{\deg P-1}\prod_{j=0}^{\deg Q-1} (a^{q^i}-b^{q^j})\\
&=(-1)^{\deg P \deg Q} \prod_{j=0}^{\deg Q-1}\prod_{i=0}^{\deg P-1} (b^{q^j}-a^{q^i})
=(-1)^{\deg P\deg Q}\xbinom PQ. \tag*{$\fin$}
\end{align*}
\end{proof}

\section{F\'ormula del conductor--discriminante}\label{S6.6}

Primero calculemos el diferente de una extensi\'on ciclot\'omica
$\cicl M{}/K$.

\begin{proposicion}\label{P6.6.1} Sea $P\in R_T$ un polinomio
m\'onico e irreducible de grado $d$ y sea $n\in{\ma N}$. Si $M=
P^n$, entonces el diferente de $\cicl Pn/K$ ${\eu D}_{P^n}$
est\'a dado por:
\begin{gather*}
{\eu D}_{P^n}=\pL^s\prod_{{\eu Q}|\pK_{\infty}}{\eu Q}^{q-2}\\
\intertext{donde $\pL$ es \'unico divisor primo sobre $\pK$,}
s=n\Phi(P^n)-q^{d(n-1)}=nq^{dn}-(n+1)q^{d(n-1)}\quad\text{y}\\
2g_{P^n}-2 = (dqn-dn-q)\frac{\Phi(P^n)}{q-1}-dq^{d(n-1)},
\end{gather*}
donde $g_{P^n}$ denota al g\'enero de $\cicl Pn$.
\end{proposicion}

\begin{proof}
Por los Teoremas \ref{T6.2.28} y \ref{T6.2.30} se tiene que cualquier
divisor primo diferente a $\pK$ y $\pK_{\infty}$ es no ramificado en
$\cicl Pn/K$, $\pK$ es totalmente ramificado y $\pK_{\infty}$ es
moderadamente ramificado con \'indice
de ramificaci\'on $q-1$. Se sigue que ${\eu D}_{P^n}=
\pL^s\prod\limits_{{\eu Q}|\pK_{\infty}}{\eu Q}^{q-2}$.

Solo falta determinar $s$. Se tiene que ${\cal O}_{P^n} =R_T[\lambda]
$ donde $\lambda$ es ra{\'\i}z de $\Psi_{P^n}(u)$ y $s=v_{\pL}(
\Psi'_{P^n}(\lambda))$. Ahora bien, puesto que $u^{P^n}=u^{P^{n-1}}
\Psi_{P^n}(u)$ se sigue que 
\begin{align*}
P^n&=(u^{P^n})'=(u^{P^{n-1}})'\Psi_{P^n}(u)+u^{P^{n-1}}\Psi_{P^n}'
(u)=\\
&=P^{n-1}\Psi_{P^n}(u)+u^{P^{n-1}}\Psi_{P^n}'(u).
\end{align*}
Por lo tanto $P^n=\lambda^{P^{n-1}}\Psi_{P^n}'(\lambda)$ y $
(\Psi_{P^n}'(\lambda))=\big(\frac{P^n}{\lambda^{P^{n-1}}}\big)$.

Puesto que $\lambda^{P^{n-1}}\in\Lam P{}$ y $\Psi_P(u)=
\prod\limits_{\mcd(A,P)=1}(u-\lambda_P^A)$, se tiene
\[
\Psi_P(0)=P=\pm \prod_{\mcd (S,P)=1}\lambda_P^S=\alpha
\lambda_P^{\Phi(P)}
\]
donde $\alpha$ es una unidad de ${\cal O}_{P^n}$. En particular
tenemos que $\big\langle\big(\lambda^{P^{n-1}}\big)^{\Phi(P)}\big
\rangle=
\langle P\rangle$.
Si ${\eu q}$ es el \'unico divisor primo de $\cicl P{}$ que divide a $\pK$,
entonces $v_{{\eu q}}(\lambda^{P^{n-1}})=\frac
{v_{{\eu q}}(P)}{\Phi(P)}=
\frac{e({\eu q}|\pK) v_{\pK}(P)}{\Phi(P)}=1$ pues ${\eu q}|\pK$ es 
totalmente ramificado en $\cicl P{}/K$. Entonces
\[
v_{\pL}(\lambda^{P^{n-1}})=e(\pL|{\eu q}) v_{{\eu q}}(
\lambda^{P^{n-1}})=\frac{\Phi(P^n)}{\Phi(P)}.
\]

Se sigue que
\begin{align*}
s&=v_{\pL}\big(\Psi_{P^n}'(\lambda)\big)=v_{\pL}\big(
\frac{P^n}{\lambda^{P^{n-1}}}\big)=v_{\pL}(P^n)-
v_{\pL}(\lambda^{P^{n-1}})\\
&=ne(\pL|\pK)-
\frac{\Phi(P^n)}{\Phi(P)}=n\Phi(P^n)-q^{d(n-1)}.
\end{align*}
Por la
f\'ormula del g\'enero de Riemann--Hurwitz (Teorema \ref{T5.6.5})
obtenemos
\begin{align*}
2g_{\cicl Pn}-2&=\big[\cicl Pn:K\big](2g_K-2)+d_{\cicl Pn}({\eu D}
_{P^n})\\
&=\Phi(P^n)(0-2)+d(n\Phi(P^n)-q^{d(n-1)})+\frac{\Phi(P^n)}{q-1}
(q-2)\\
&=\frac{\Phi(P^n)}{(q-1)}\big(-2(q-1)+dn(q-1)+(q-2)\big)-dq^{d(n-1)}\\
&=(dnq-dn-2q+2-q-2) \frac{\Phi(P^n)}{q-1}-dq^{d(n-1)}\\
&=(dqn-dn-q)\frac{\Phi(P^n)}{q-1}-dq^{d(n-1)}. \tag*{$\fin$}
\end{align*}
\end{proof}

El resultado general es:

\begin{teorema}[F\'ormula del g\'enero y del diferente]\label{T6.6.2}
Sea $M\in R_T$ un polinomio m\'onico no constante de la forma
$M=P_1^{\alpha_1}\cdots P_r^{\alpha_r}$ donde $P_1,\ldots, P_r$
son polinomios irreducibles distintos. Sea $d_i=\gr P_i$. Entonces
\[
{\eu D}_M=\prod_{i=1}^r\big(\prod_{\pL|\pK_i}\pL\big)^{s_i}
\prod_{{\eu Q}|\pK_{\infty}}{\eu Q}^{q-1}
\]
donde $(P_i)_K=\frac{\pK_i}{\pK_{\infty}^{d_i}}$, $s_i=\alpha_i
\Phi(P_i^{\alpha_i})-q^{d_i(\alpha_i-1)}$ y 
\[
2g_M-2=-2\Phi(M)+\sum_{i=1}^r d_is_i\frac{\Phi(M)}{\Phi(P_i^{\alpha_i
})}+(q-2)\frac{\Phi(M)}{q-1}.
\]
\end{teorema}

\begin{proof}
Para cada $i\in\{1,\ldots, r\}$, $\pK_i$ es completamente ramificado
en $\cicl {P_i}{\alpha_i}/K$ y no ramificado en $\cicl M{}/\cicl {P_i}
{\alpha_i}$. Para cada divisor primo ${\eu q}$ en $\cicl {P_i}{\alpha_i}$
que est\'a sobre $\pK_i$, hay $\frac{\Phi(M)/\Phi(P_i^{\alpha_i})}{f_i}$
divisores primos, cada uno de ellos de grado relativo $f_i$. Por  tanto
la contribuci\'on a ${\eu D}_M$ de $\pK_i$ es $\big(\prod\limits_{\pL|
\pK_i}\pL\big)^{s_i}$ donde $s_i$ es como en la Proposici\'on 
\ref{P6.6.1}. Tenemos 
\[
\gr_{\cicl M{}}\big(\prod_{\pL|\pK_i}\pL\big)=d_i\frac{\Phi(M)/\Phi(
P_i^{\alpha_i})}{f_i}=d_i\frac{\Phi(M)}{\Phi(P_i^{\alpha_i})}.
\]
Por lo
tanto ${\eu D}_M=\prod\limits_{i=1}^r \big(\prod\limits_{\pL|\pK_i}\pL
\big)^{s_i}\prod\limits_{{\eu Q}|\pK_{\infty}}{\eu Q}^{q-2}$ y
\begin{align*}
2g_M-2&= (2g_K-2)\big[\cicl M{}:K\big]+\gr_{\cicl M{}}{\eu D}_M\\
&=-2\Phi(M)+\sum_{i=1}^r s_id_i\frac{\Phi(M)}{\Phi(P_i^{\alpha_i})}+
(q-2)\frac{\Phi(M)}{q-1}. \tag*{$\fin$}
\end{align*}
\end{proof}

Veamos el siguiente ejemplo.

\begin{ejemplo}\label{Ej6.6.3}
Sean $q=3$, $M=T^2(T+1)$, $\zeta=\zeta_6=e^{2\pi i/6}$ y
\begin{eqnarray*}
\theta_T\colon \units T2&\longto& {\ma C}^{\ast}\\
1&\longmapsto &1\\
T-1&\longmapsto &\zeta\\
T+1&\longmapsto &\zeta^2\\
-1&\longmapsto &-1\\
-T+1&\longmapsto &-\zeta \\
-T-1&\longmapsto &-\zeta^2.
\end{eqnarray*}
Sea $\widetilde{\theta_T}=\theta_T\circ \varphi_{M,T^2}$, y sea
$\theta_{T+1}$ dado por
\begin{eqnarray*}
\theta_{T+1}\colon \units {T+1}{}&\longto &{\ma C}^{\ast}\\
1&\longmapsto &1\\
-1&\longmapsto &-1
\end{eqnarray*}
y sea $\widetilde{\theta_{T+1}}=\theta_{T+1}\circ \varphi_{M,T+1}$.

Sea $X=\langle \widetilde{\theta_T},\widetilde{\theta_{T+1}}\rangle$.
El campo perteneciente a $X$ es $L=\cicl M{}$. Tenemos $X_T=
\langle \theta_T\rangle$, $X_{T+1}=\langle\theta_{T+1}\rangle$.
As{\'\i}, $e_T=6$, $e_{T+1}=2$. Notemos que los grupos $Y$ y $Z$
(ver Teorema \ref{T6.5.9}) de $T$ y $T+1$ satisfacen:
\begin{gather*}
Y(T+1)=\{\tau\in X\mid \tau(T+1)\neq 0\}=\langle\widetilde{\theta_T}
\rangle \cong X_T\quad \text{y}\\
Z(T+1)=\{\tau\in X\mid \tau(T+1)=1\}=\langle\widetilde{\theta_{T+1}}^3
\rangle,\\
\intertext{luego $e_{T+1}=2$, $f_{T+1}=3$ y $h_{T+1}=2$. Ahora bien}
Y(T)=\{\tau\in X\mid \tau(T)\neq 0\}=\langle\widetilde{\theta_{T+1}}
\rangle\cong X_{T+1}\quad\text{y}\\
Z(T)=\{\tau\in X\mid \tau(T)=1\}=\{1\},
\end{gather*}
luego $e_T=6$, $f_T=2$ y $h_T=1$. Para $\pK_{\infty}$ tenemos
$e_{\infty}=2$, $f_{\infty}=1$, $h_{\infty}=6$. Por tanto
$s_1=2\Phi(T^2)-3=9$, $s_2=1\Phi(T+1)-3^0=1$, de donde el
diferente de $\cicl M{}/K$ es $
{\eu D}_M=\pL_T^9\pL_{T+1,1}\pL_{T+1,2}{\eu Q}_1{\eu Q}_2
{\eu Q}_3{\eu Q}_4{\eu Q}_5{\eu Q}_6$. Se sigue que el diferente
de ${\cal O}_M$ sobre $R_T$ es 
\begin{gather*}
{\eu D}_{{\cal O}_M/R_T}=\pL_T^p\pL_{T+1,1}\pL_{T+1,2}
\quad\text{y}\\
 {\eu d}_{{\cal O}_M/R_T}= N_{\cicl M{}/K}
{\eu D}_{{\cal O}_M/R_T}=T^{18}(T+1)^6.
\end{gather*}

Por otro lado tenemos
\begin{gather*}
\begin{tabular}{||c|c||}\hline\hline
\quad Caracter\quad{\ } & \quad Conductor\quad {\ } \\ \hline
$1$ & $1$ \\ \hline
$\widetilde{\theta_T}$ & $T^2$ \\ \hline
$\widetilde{\theta_T}^2$ & $T^2$ \\ \hline
$\widetilde{\theta_T}^3$ & $T$ \\ \hline
$\widetilde{\theta_T}^4$ & $T^2$ \\ \hline
$\widetilde{\theta_T}^5$ & $T^2$ \\ \hline
$\widetilde{\theta_{T+1}}$ & $T+1$ \\ \hline
$\widetilde{\theta_T}\widetilde{\theta_{T+1}}$ & $T^2(T+1)$ \\ \hline
$\widetilde{\theta_T}^2\widetilde{\theta_{T+1}}$ & $T^2(T+1)$ \\ \hline
$\widetilde{\theta_T}^3\widetilde{\theta_{T+1}}$ & $T(T+1)$ \\ \hline
$\widetilde{\theta_T}^4\widetilde{\theta_{T+1}}$ & $T^2(T+1)$ \\ \hline
$\widetilde{\theta_T}^5\widetilde{\theta_{T+1}}$ & $T^2(T+1)$ \\ \hline \hline
\end{tabular}\\
\intertext{As{\'\i}}
\begin{align*}
\prod_{\chi\in X}F_{\chi}&=1\cdot T^2 \cdot T^2\cdot T \cdot T^2
\cdot T^2\cdot (T+1)\cdot T^2(T+1)\cdot T^2(T+1)\cdot \\
&\hspace{1cm}\cdot T(T+1)
\cdot T^2(T+1)\cdot T^2(T+1)=\\
&= T^{18}(T+1)^6={\eu d}_{{\cal O}_M/R_T}.
\end{align*}
\intertext{Es decir}
\prod_{\chi\in X}F_{\chi}={\eu d}_{{\cal O}_M/R_T}.
\end{gather*}
\end{ejemplo}

\begin{teorema}[F\'ormula del conductor--discriminante\index{f\'ormula
del conductor--discriminante}]\label{T6.6.4}
Sea $L$ un subcampo de $\cicl M{}$ donde $M\in R_T$ es un 
polinomio m\'onico no constante. Sea ${\eu d}_{L/K}$ el
discriminante de ${\cal O}_L/R_T$ donde ${\cal O}_L$ es la 
cerradura entera de $R_T$ en $L$.
Sea $X_L$ el grupo de caracteres de Dirichlet asociado a $L$.
Entonces
\[
{\eu d}_{L/K}=\prod_{\chi\in X_L}F_{\chi}.
\]
\end{teorema}

\begin{proof}
La prueba se puede hacer de manera similar a la del caso 
num\'erico (Teorema \ref{T12.3.1.8}) y una demostraci\'on en el caso
de campos de funciones puede ser consultado en \cite{RzeVil2009}. Aqu{\'\i}
presentamos una demostraci\'on diferente.

Primero supongamos que $M=P^n$ donde $P\in R_T$ es un 
polinomio irreducible. Sean
\[
L_i:=L\cap \cicl Pi,\quad i=0,1,2,\ldots, n.
\]
Entonces $L_0=L$ y $L_0=K$. Tenemos que un caracter $\chi$
tiene conductor $P^j$ si y s\'olo si $\chi$ es un caracter asociado
al campo $\cicl Pj$ pero no a $\cicl P{j-1}$. Se sigue que $X_L$
contiene exactamente $[L_j:K]-[L_{j-1}:K]$ caracteres de conductor
$P^j$, $1\leq j\leq n$. Por lo tanto
\begin{gather}
\prod_{\chi\in X_L} F_{\chi}=P^{\alpha}\nonumber\\
\intertext{donde}
\alpha=\sum_{j=1}^nj\big([L_j:K]-[L_{j-1}:K]\big)=n[L_n:K]-
\sum_{j=0}^{n-1} [L_j:K].\nonumber\\
\intertext{Por lo tanto}
\prod_{\chi\in X_L}F_{\chi}=P^{\alpha}\quad\text{con}\quad
\alpha=n[L_n:K]-\sum_{j=0}^{n-1}[L_j:K]. \label{Ec6.6.5}
\end{gather}

Ahora bien, si probamos que ${\eu D}_{L/K}=\pK_L^{\alpha}$
donde $\pK_L:=\pK_n\cap {\cal O}_L$ y $\pK_n$ es el \'unico 
divisor primo de $\cicl Pn$ sobre $\pK$, puesto que el
grado relativo de $\pK_L$ sobre $\pK$ es $1$, se seguir\'a que
\[
{\eu d}_{L/K}=N_{L/K}{\eu D}_{L/K}=P^{\alpha}.
\]

Sea ${\eu D}_{L/K}=\pK_L^{\gamma}$. Tenemos ${\cal O}_{P^n}=
R_T[\lambda_{P^n}]$ y ${\cal O}_{P^n}={\cal O}_L[\lambda_{P^n}]$.
Sea $f(u):=\Irr(u,\lambda_{P^n}, L)$. Entonces $f(u)$ divide al
polinomio ciclot\'omico
\begin{gather*}
\Psi_{P^n}(u):=\prod_{\substack{(A,P)=1\\ \deg A
<\deg P^n}}(u-\lambda_{P^n}^A)=\Irr(u,\lambda_{P^n},K).\\
\intertext{Entonces}
\Psi_{P^n}(u)=\prod_{\sigma_A\in G}(u-\lambda_{P^n}^A),
\quad f(u)=\prod_{\sigma_A\in H}(u-\lambda_{P^n}^A)
\end{gather*}
donde para cualquier $A\in R_T$, primo relativo a $P$, definimos
$\sigma_A(
\lambda_{P^n})=\lambda_{P^n}^A$ y donde 
$G:=G_{P^n}=\Gal(\cicl Pn/K)$ y $H:=\Gal(k(\Lambda_{
P^n})/L)$.

Tenemos $\gr \Psi_{P^n}(u)=\Phi(P^n)=[K(\Lambda_{P^n}):K]=
q^{(n-1)d}(q^d-1)=|G|$ y $\gr f(u) =[K(\Lambda_{P^n}):K]=|H|$.
Escribimos $\Psi_{P^n}(u)=f(u)g(u)$. Por lo tanto
\[
g(u)=\prod_{\sigma_A\in G\setminus H}(u-\lambda_{P^n}^A).
\]

Se sigue que ${\eu D}_{K(\Lambda_{P^n})/K}=(\Psi'_{P^n}
(\lambda_{
P^n})) = {\eu p}_n^{\beta}$ donde
 $\beta=nq^{dn}-(n+1)q^{d(n-1)}$ (Proposici\'on \ref{P6.6.1}) y
${\eu D}_{K(\Lambda_{P^n})/K}=(f'(\lambda_{P^n}))= {\eu p}_n^{
\delta}$. Notemos que 
\begin{gather*}
\Psi_{P^n}'(u) = f'(u)g(u)+f(u)g'(u)\quad {\text{y}}\quad
\Psi'_{P^n}(\lambda_{P^n})=f'(\lambda_{P^n})g(\lambda_{P^n}).
\end{gather*}
Puesto que
 ${\eu D}_{K(\Lambda_{P^n})/K}={\eu D}_{K(\Lambda_{P^n})/L}
\con_{L/K(\Lambda_{P^n})} {\eu D}_{L/K}$ y ${\eu p}_L$ es
totalmente ramificado en $K(\Lambda_{P^n})/L$ obtenemos que
\begin{align*}
\gamma&=\frac{\beta-\delta}{[K(\Lambda_{P^n}):L]}=\frac{1}
{[K(\Lambda_{P^n}):L]}(v_{{\eu p}_n}(\Psi'_{P^n}(\lambda_{P^n}))
-v_{{\eu p}_n}(f'(\lambda_{P^n})))\\
&=\frac{1}{[K(\Lambda_{P^n}):L]}
 v_{{\eu p}_n}\Big(\frac{\Psi'_{P^n}(\lambda_{P^n})}
{f'(\lambda_{P^n})}\Big)
=\frac{1}{[K(\Lambda_{P^n}):L]}v_{{\eu p}_n}(g(\lambda_{P^n})).
\end{align*}

Esto es
\begin{equation}\label{Eq3}
\gamma=\frac{1}{[K(\Lambda_{P^n}):L]}v_{{\eu p}_n}(g(\lambda_{P^n})).
\end{equation}

Se tiene $g(\lam Pn)=\prod\limits_{\sigma_A\in G\setminus H}
\big(\lam Pn-\lam P{n}^A\big)$. Definimos la filtraci\'on
${\mathcal D}_i:=\{\sigma_A\in G\mid v_P(A-1)\geq i\}$,
$i=0,1,\ldots, n-1$. Tenemos ${\mathcal D}_{i+1}\subseteq
{\mathcal D}_i$.  Si $A\in R_T$ es tal que $v_P(A-1)=t$, entonces
$A=1+P^tR$ con $R$ y $P$ son primos relativos y 
\[
\sigma_A(\lam Pi)=\lam Pi^A=\lam 
Pi^{1+P^tR}=\lam Pi + (\lam Pi^{P^t})^R.
\]
Por lo tanto
 $\sigma_A(\lam Pi)= \lam Pi$ si y s\'olo si $t\geq i$, esto es,
\[
{\mathcal D}_i=\Gal(K(\Lam Pn)/K(\Lam Pi)).
\]

Consideramos
 ${\mathcal C}_i:=\{\sigma_A\in G\setminus H\mid v_P(A-1)=i\}=
(G\setminus H)\cap ({\mathcal D}_i\setminus {\mathcal D}_{i+1})=
{\mathcal D}_i\setminus (H\cup {\mathcal D}_{i+1})$.

Por lo tanto
\begin{align*}
|{\mathcal C}_i|&=|{\mathcal D}_i|-|{\mathcal D}_i \cap H|-
|{\mathcal D}_i \cap {\mathcal D}_{i+1}|+|{\mathcal D}_i\cap {\mathcal D}_{i+1} \cap H|\\
&=|{\mathcal D}_i|-|{\mathcal D}_i \cap H|-
|{\mathcal D}_{i+1}|+|{\mathcal D}_{i+1} \cap H|.
\end{align*}
Ahora,
\begin{align*}
{\mathcal D}_{i}\cap H&=\Gal(K(\Lam Pn)/LK(\Lam Pi))\quad\text{y}\\
{\mathcal D}_{i+1}\cap H&=\Gal(K(\Lam Pn)/LK(\Lam P{i+1})).
\end{align*}
Por lo tanto
\begin{multline*}
|{\mathcal C}_i|=\ext n{}i-\ext nLi\\
-\ext n{}{i+1}+\ext nL{i+1}, \quad 0\leq i\leq n-1.
\end{multline*}

Por otro lado, tenemos que $\sigma_A\in {\mathcal C}_i$ si y 
s\'olo si $A=1+P^iR$ con $R$ primo relativo
a $P$. Por tanto, si $\sigma_A \in 
{\mathcal C}_i$, entonces 
\[
\lam Pn-\lam Pn^A=\lam Pn
-\lam Pn-(\lam Pn^{P^i})^R=
-\lam P{n-i}^R, \quad 0\leq i\leq n-1.
\] 
As{\'\i} $\sigma_A\in {\mathcal C}_i$
si y s\'olo si $\sigma_A\in G\setminus H$ 
y $v_{{\eu p}_n}(\lam Pn-
\lam Pn^A)= v_{{\eu p}_n}(\lam P{n-i}^R)=
\frac{\Phi(P^n)}{\Phi(P^{n-i})}=q^{id}$.

Tenemos $G\setminus H=\bigcup\limits_{i=0}^{n-1}{\mathcal C}_i$ y
${\mathcal C}_i\cap {\mathcal C}_j=
\emptyset$ para $i\neq j$. Por lo tanto
\begin{align*}
v_{{\eu p}_n}(g(\lam Pn))&= \sum_{i=0}^{n-1}|{\mathcal C}_i| q^{id}=
\sum_{i=0}^{n-1}\ext n{}iq^{id}\\
&\hspace{1cm} -\sum_{i=0}^{n-1}\ext n{}{i+1} q^{id}\\
&\hspace{1cm} -\sum_{i=0}^{n-1}\ext nLi q^{id}\\
&\hspace{1cm}+\sum_{i=0}^{n-1}
\ext nL{i+1} q^{id}\\
&= [K(\Lam Pn):K]q^0-\ext n{}nq^{(n-1)d}\\
&\hspace{1cm}+\sum_{i=1}^{n-1}\ext n{}i q^{(i-1)d}
(q^d-1)\\
&\hspace{1cm} -[K(\Lam Pn):K]q^0+\ext n{}n q^{(n-1)d}\\
&\hspace{1cm}-\sum_{i=1}^{n-1}
\ext nLi q^{(i-1)d}(q^d-1).
\end{align*}

Por tanto
\begin{align*}
v_{{\eu p}_n}(g(\lam Pn))
&= [K(\Lam Pn):K]-[K(\Lam Pn):K]+\sum_{i=1}^{n-1}[K(\Lam Pn):L]\\
&\hspace{1cm} -
\sum_{i=1}^{n-1} \ext nLi [K(\Lam Pi):K]\\
&= n[K(\Lam Pn):K]-[K(\Lam Pn):L]\\
&\hspace{1cm} -\sum_{i=1}^{n-1}\ext nLi [K(\Lam Pi):K].
\end{align*}

Se sigue que
\[
\gamma=n\frac{[K(\Lam Pn):K]}{[K(\Lam Pn):L]}-1-\sum_{i=1}^{n-1}t_i,
\]
donde
\begin{align*}
t_i&=\frac{\ext nLi [K(\Lam Pi):K]}{[K(\Lam Pn):L]}\\
&=\frac{[K(\Lam Pn):K]}{[LK(\Lam Pi):K(\Lam Pi)][K(\Lam Pn):L]}\\
&= \frac{[L:K]}{[LK(\Lam Pi):K(\Lam Pi)]}
=\frac{[L:K]}{[L:L\cap K(\Lam Pi)]}\\
&=[L\cap K(\Lam Pi):K]=[L_i:K].
\end{align*}
\[
\xymatrix{
{L}\ar@{-}[r]\ar@{-}[d]&{LK(\Lam Pi)}\ar@{-}[d]\\
{L\cap K(\Lam Pi)}\ar@{-}[r]&{K(\Lam Pi)}}
\]

Por lo tanto
\begin{equation}\label{Eq4}
\gamma=n[L:K]-\sum_{i=0}^{n-1}[L_i:K].
\end{equation}

El caso $M=P^n$ se sigue de (\ref{Ec6.6.5}) y (\ref{Eq4}).

Para $M$ arbitrario, sea $P\in R_T$ un polinomio m\'onico e 
irreducible. Entonces escribimos $M=P^a A$ donde $A\in R_T$
y $P\nmid A$. Sea $E=L(\Lam A{})=LK(\Lam A{})$. Se tiene
el diagrama (ver Teorema \ref{T6.5.5})
\[
\xymatrix{
&&{K(\Lambda_{P^a A})}\ar@{-}[d]\ar@{-}[ddll]\ar@{-}[rrdd]\\
&&{E=LK(\Lambda_{A})}\ar@{-}[ddl]_{\hbox{\rm\tiny
no ramificado}}\ar@{-}[d]^{\hbox{\rm\tiny
no ramificado}}\ar@{-}[drr]\\
{K(\Lambda_{P^a})}\ar@{-}[dr]&&{L=K_{X_L}}
\ar@{-}[dd]^{e=e_P(L/K)}&&
{K(\Lambda_A)}\ar@{-}[ddll]^{\hbox{\rm
\tiny no ramificado}}\\
&{K_{X_P}}\ar@{-}[dr]_{\hbox{\rm\tiny totalmente
ramificado\phantom{xx}}}\\&&{K}}
\]

Tenemos
\[
 {\eu D}_{E/K}= {\eu D}_{E/L}
 \con_{L/E}{\eu D}_{L/K} = {\eu D}_{E/K_{X_P}}\con_{
 K_{X_P}/E}{\eu D}_{K_{X_P}/K}.
 \]
Denotemos para cualquier extensi\'on $E/F$ y $P\in R_T$,
${\eu d}_{E/F}(P)= P^s$ donde $P^s|{\eu d}_{E/F}$ y $P^{s+1}\nmid
{\eu d}_{E/F}$. Entonces
\begin{align*}
 {\eu d}_{E/F}(P)&=\big(N_{E/K}({\eu D}_{E/L})\big)(P)\cdot 
 {\eu d}_{L/K}^{[E:L]}(P)\\
& =\big(N_{E/K}({\eu D}_{E/K_{X_P}})\big)(P)\cdot
 {\eu d}_{K_{X_P}/K}^{[E:K_{X_P}]}(P).
\end{align*}

Puesto que $P$ es no ramificado en
 $E/L$ y en $E/K_{X_P}$ se tiene
 \begin{gather*}
 \big(N_{E/L}({\eu D}_{E/L})\big)(P)=
 \big(N_{E/K}({\eu D}_{E/K_{X_P}})\big)(P)=1.\\
\intertext{Por tanto}
{\eu d}_{L/K}(P)=\Big({\eu d}_{K_{X_P}/K}(P) 
\Big)^{[E:K_{X_P}]/[E:L]}.
\end{gather*}

Se tiene $[E:K_{X_P}]=[K(\Lambda_A):K]=\Phi(A)$ y $[E:L]=
  [Y:X_L]$ donde $Y$ es el grupo de caracteres de 
Dirichlet asociado a $E$. Se tiene
$Y=X_P\times \widehat{G_A}$. Se sigue que
\begin{align*}
[E:L]&=\frac{|Y|}{|X_L|}=\frac{|X_P|\big|
\widehat{G_A}\big|}{|X|}=\frac{\big|X_P\big|\Phi(A)}{|X|}\\
\intertext{y}
\frac{[E:K_{X_P}]}{[E:L]}&= \frac{\Phi(A)}{\Big(\frac{\big|
X_P\big|\Phi(A)}{|X|}\Big)}
=\frac{|X|}{\big|X_P\big|}=\frac{[L:K]}{\big|X_P\big|}=\frac{edh}{e}=dh
\end{align*}
donde $d$ es el grado relativo de los divisores primos de $L$ sobre
$\pK$ y $h$ el n\'umero de estos divisores primos.

Puesto que $K_{X_P}\subseteq \cicl Pa$, de la primera parte
de esta demostraci\'on obtenemos que
 ${\eu d}_{K_{X_P}/K}(P) =\prod\limits_{\varphi\in X_P}F_{\varphi}$.
 Por lo tanto
\[
{\eu d}_{K_{X}/K}(P)=\Big(\prod_{\varphi\in X_P}
F_{\varphi}\Big)^{([E:K_{X_P}]/
 [E:L])}=\Big(\prod_{\varphi\in X_P}F_{\varphi}\Big)^{dh}=
 \prod_{\varphi\in X_P}F_{\varphi}^{dh}.
 \]

Del epimorfismo natural $\pi\colon X_L\to X_P$, 
$\chi\mapsto \chi_P$ obtenemos
$|\ker \pi|=\frac{\big|X_L\big|}{\big|X_P\big|}=dh$. Por tanto,
para cada $\varphi\in X_P$,
$ |\pi^{-1}(\varphi)|=dh$. Esto es, para cada $\varphi\in X_P$ 
hay precisamente $dh$ elementos
$\chi\in X$ tales que $\pi(\chi)=\chi_P=\varphi$. Se sigue que
\[
{\eu d}_{L/K}(P)=\Big(\prod_{\varphi\in X_P}F_{\varphi}\Big)^{dh}=
\prod_{\chi\in X}F_{\chi_P}.
\]

Finalmente, tenemos que $F_{\chi}=\prod\limits_P F_{\chi_P}$
para $\chi\in X$ (Teorema \ref{T6.4.19}) y
$ {\eu d}_{L/K}=\prod\limits_{P} {\eu d}_{L/K}(P)$. Por lo tanto
\[
{\eu d}_{L/K}=\prod_{\chi\in X_L}F_{\chi}. \tag*{$\fin$}
\]
 \end{proof}

%% file: Capitulo10.tex
\chapter{Grupos de clase y ramificaci\'on de campos
de funciones congruentes}\label{ChRam}

En este cap\'itulo concentramos varios resultados sobre los grupos de
clase de campos de funciones congruentes as\'i como la aritm\'etica
de extensiones de este tipo de campos.

En este cap\'itulo consideraremos un campo de funciones $\K$
con campo de constantes $\F$, con $q=p^u$ para alguna $u\in
{\ma N}$. $K$ denotar\'a el campo global de
funciones racionales $\F(T)$. Para cualquier campo de 
funciones congruentes $\K$, $\K_m$ denotar\'a la extensi\'on
de constantes de grado $m$: $\K_m=\K{\ma F}_{q^m}$.

\section{Grupos de clase}\label{SRam.1}

En esta parte hacemos un compendio de algunos de los resultados
relevantes de grupos de clases de campos de funciones. El \'enfasis
por supuesto es para los campos de funciones congruentes, pero
presentaremos algunos resultados generales.

Empezamos por recordar algunos resultados. Las definiciones
b\'asicas est\'an dadas en el Cap\'itulo \ref{Ch5}, Secci\'on 
\ref{S5.2}.

\subsection{N\'umero de clase $h_{\K}$}\label{SRam.1.1}

\begin{definicion}\label{D10.1.1.1} Dado un campo $\K/k$
arbitrario de campos de funciones, el {\em n\'umero de
clase\index{numero de clase@n\'umero de clase}} $h_{\K}$ de $\K$ es la 
cardinalidad de su grupo de clases de divisores de grado $0$:
$h_{\K}=|I_{\K,0}|$.
\end{definicion}

\begin{observacion}\label{O10.1.1.2} A diferencia de los campos
num\'ericos, no se define $h_{\K}=|I_{\K}|$ pues este siempre es
infinito ya que $I_{\K}\cong I_{\K,0}\oplus {\ma Z}$.
\end{observacion}

\begin{teorema}\label{T10.1.1.3} Sea $k$ un campo arbitrario, y sea
$\K$ un campo de funciones sobre $k$. Si el g\'enero de $\K$ es
$g_{\K}=0$, entonces $h_{\K}=1$.
\end{teorema}

\begin{proof}
Sea ${\eu A}$ un divisor de grado $0$. Como $g_{\K}=0$, se tiene
$0>2g_{\K}-2=-2$, por un corolario al Teorema de Riemann--Roch, Corolario 
\ref{C5.3.7}, se tiene $\ell({\eu A}^{-1})=d({\eu A})-g_{\K}+1=0-0+1=1$.

Puesto que ${\eu A}^{-1}$ es de grado $0$ y de dimensi\'on $1$,
${\eu A}$ debe ser principal, esto es, ${\eu A}\in P_{\K}$ y $I_{\K,0}
=\{1\}$. $\fin$
\end{proof}

\begin{observacion}\label{O10.1.1.4} Notemos que se ha probado
que $h_{\K}=1$ con la \'unica suposici\'on de que $g_{\K}=0$, es
decir, no se ha pedido que $\K$ sea de funciones racionales. En
el caso de que un campo de funciones congruentes satisfaga $g_{\K}
=0$, $\K$ necesariamente es de funciones racionales. Esto es 
consecuencia de la Hip\'otesis de Riemann (Ecuaci\'on
(\ref{E10.1.2.1.1})).

Pueden existir campos de funciones $\K/k$ tales que $g_{\K}=0$ sin que
$\K$ sea de funciones racionales. Por ejemplo, si $k={\ma R}$ y $\K=
{\ma R}(x,y)$ con $x^2+y^2=1$, entonces $g_{\K}=0$ y $\K$ no es de
funciones racionales (ver \cite[despu\'es del Corollary 4.1.8]{Vil2006}).
\end{observacion}

Cuando $g_{\K}>0$ en general tenemos que $h_{\K}=\infty$. Sin
embargo cuando $\K$ es un campo de funciones congruentes, se
tiene $h_{\K}<\infty$. A continuaci\'on presentamos una demostraci\'on
de este hecho. Pero antes, hacemos notar que existen campos
de funciones $\K$ con campo de constantes infinito, esto es, no
congruentes, con $g_{\K}>0$ y $h_{\K}<\infty$ (ver 
Ejemplo \ref{Ej10.1.1.6}).

\begin{teorema}\label{T10.1.1.5} Sea $k$ un campo finito y $\K$ un 
campo de funciones con campo de constantes $k$. Entonces el
n\'umero de clase $h_{\K}$ es finito.
\end{teorema}

\begin{proof}
Sea $C$ cualquier clase y sea $N(C)=n$ la dimensi\'on de $C$. Sean
${\eu A}_1,\ldots,{\eu A}_m$ todos los divisores enteros pertenecientes
a la clase $C$. Sea ${\eu A}\in C$ arbitrario. Sea $A=\{x\in\*k\mid
(x)_{\K}=\frac{{\eu A}_i}{{\eu A}}\text{\ para alguna\ } 1\leq i\leq m\}$.
Sea $k=\F$. Notemos que $(x)_{\K}=(y)_{\K} \iff$ existe $\alpha\in
\*\F$ tal que $y=\alpha x$.

As\'i $|A|=(q-1)m$. Por otro lado $\dim_{\F}(A\cup\{0\})=n$, esto es,
$|A\cup \{0\}|=q^n$, $|A|=q^n-1$. Por tanto $m=\frac{q^n-1}{q-1}$. En
resumen, el n\'umero total de divisores enteros en una clase es
$\frac{q^{N(C)}-1}{q-1}$.

Sea ahora $x\in\K\setminus k$, $[\K:k(x)]=r<\infty$. Si $\pL$ es un
lugar de $\K$ de grado $t$, entonces $\pK:=\pL|_{k(x)}$ es de 
grado $\leq t$.

Como el n\'umero de polinomios en $k[x]$ de grado menor o igual
a $t$ es finito, se sigue que el n\'umero de divisores en $\K$ de grado
$t$ es finito.

Sea $t\geq g_{\K}$ y sea $C$ cualquier clase de grado $t$. Por el
Teorema de Riemann--Roch, $N(C)=d(C)-g_{\K}+1+N(W C^{-1})\geq
t-g_{\K}+1\geq 1>0$. Por tanto $C$ contiene un divisor entero ${\eu A}$.

Sea $\Lambda_t=\{C\mid C \text{\ es una clase de grado\ } t\}$. 
Entonces $\Lambda_t$ es finito pues cada elemento de $\Lambda_t$
contiene un divisor entero de grado $t$ y el n\'umero de estos
divisores es finito. Sea $\varphi\colon I_{\K,0}\lra
\Lambda_t$, $\varphi(C_0)=C_0{\eu A}\in \Lambda_t$. Es f\'acil
de verificar que $\varphi$ est\'a bien definida y que es biyectiva por
tanto $h_{\K}=|I_{\K,0}|=|\Lambda_t|<\infty$.
$\fin$
\end{proof}

Una pregunta natural es que si $|k|=\infty$ y $g_{\K}>0$, entonces
?`$h_{\K}=\infty$? La respuesta es no.

\begin{ejemplo}\label{Ej10.1.1.6}
Sea $\K/k$ un campo de funciones el\'ipticas con $k$ arbitrario.
Existe un lugar $\pK_0$ de grado $1$ y se tiene $g_{\K}=1$
(ver \cite[Section 4.2]{Vil2006}).
Se tiene \cite[Proposition 9.6.9]{Vil2006} que existe una biyecci\'on
entre los lugares de grado $1$ de $\K$ y el grupo $I_{K,0}=
D_{K,0}/P_K$.

Sea $k={\ma Q}$ y sea $\K={\ma Q}(x,y)$ con $y^2=x^3-x$.
Se tiene que $\K$ es un campo de funciones el\'ipticas con $\pK_0$
el lugar de grado uno el cual satisface
${\eu N}_x=\pK_0^2$ y ${\eu N}_y=
\pK_0^3$, donde ${\eu N}_z$ denota al divisor de polos de $z\in\K$.
Adem\'as $g_{\K}=1$. Sea $M=\{\pK\in {\ma P}_{\K}\mid 
d_{\K}(\pK)=1\}$. Entonces $h_{\K}=|M|$.

\[
\begin{minipage}{2cm}
\xymatrix{K\ar@{-}[d]_2\\ {\ma Q}}
\end{minipage}
\qquad\qquad
\begin{minipage}{7cm}
Sea $\pK$ un lugar de grado $1$ en $\K$. Entonces ${\pK}|_{{\ma Q}
(x)}$ es de grado $1$. Se tiene que $\pK_0|_{{\ma Q}(x)}={\mc P}_{
\infty}$ es el polo de $x$ en ${\ma Q}(x)$.
\end{minipage}
\]

Sea $\pK\neq \pK_0$. Entonces, si $\varphi_{\pK}$ es el lugar asociado
a $\pK$, se tiene $\varphi_{\pK}:K\lra \o_{\pK}/{\pK}\cup\{\infty\}$ y
$\o_{\pK}/\pK\cong k$ pues $\deg_K\pK=1$.

Como $\K={\ma Q}(x,y)$, entonces $\varphi_{\pK}(x)\in k\cup\{\infty\}$,
$\varphi_{\pK}(y)\in k\cup\{\infty\}$. Puesto que $v_{\pK}(x)\geq 0$, se 
tiene que $\varphi_{\pK}(x)\in k$ (recordemos que $\varphi_{\pK}(z)=0
\iff v_{\pK}(z)>0$ y que $\varphi_{\pK}(z)=\infty\iff
v_{\pK}(z)<0$). Similarmente, debido a que ${\eu N}_y=\pK_0^3$,
obtenemos que $\varphi_{\pK}(y)=\mu\in k$.

Ahora sea $\pK|_{{\ma Q}(x)}\sim x-\alpha$, $\alpha\in{\ma Q}$ ya que
$\pK|_{{\ma Q}(x)}$ es de grado $1$ y $\pK|_{{\ma Q}(x)}\neq {\mc P}_{
\infty}$. Se tiene que $\varphi_{\pK}(x-\alpha)=\varphi_{\pK}(\alpha)=
\varphi_{\pK}(x)-\alpha$, esto es, $\varphi_{\pK}(x)=\alpha$. Puesto
que $y^2=x^3-x$, se sigue que $\mu^2=\alpha^3-\alpha$, esto es,
$\alpha$ es ra\'iz de $f_{\mu}(X):=X^3-X-\mu^2$. 

Por otro lado, si $\mu=0$ ($\iff v_{\pK}(y)>0$), entonces $\alpha\in\{
0,1,-1\}$. Si $\mu\neq 0$, se verifica que $f_{\mu}$ es irreducible,
viendo directamente que $f_{\mu}$ no tiene ra\'ices racionales.
En resumen, para $\mu\neq 0$, $f_{\mu}$ es irreducible y por tanto
no existe $\alpha\in {\ma Q}\setminus\{0\}$ con $\varphi_{\pK}(y)=
\mu\in k$ para $\alpha\neq 0,1,-1$ que son las ra\'ices de $f_0(X)
=X^3-X$. Se sigue $|M|\leq 7$ (con esto basta para la finitud de
$I_{K,0}$), a saber, $\pK_{\alpha}={\mc P}_{\alpha}^2$
o ${\mc P}_{\alpha,1}{\mc P}_{\alpha,2}$, donde $\deg {\mc P}_{
\alpha}=\deg {\mc P}_{\alpha,1}=\deg{\mc P}_{\alpha,2}=1$ y
$\p=\pL_{\infty}^2$. De hecho, si $\pL_{\alpha}|{\mc P}_{\alpha}$,
$\alpha=0,1,-1$, entonces 
\begin{align*}
v_{\pL_{\alpha}}(y^2)&=2 v_{\pL_{\alpha}}(y)=v_{\pL_{\alpha}}(x^3-x)=
e(\pL_{\alpha}|{\mc P}_{\alpha})\cdot v_{{\mc P}_{\alpha}}(x^3-x)\\
&=e(\pL_{\alpha}|{\mc P}_{\infty})\cdot 1=e(\pL_{\alpha}|{\mc P}_{\alpha}).
\end{align*}
Por tanto $v_{\pL_{\alpha}}(y)=1$ y $e(\pL_{\alpha}|{\mc P}_{\alpha})=2$
y por la f\'ormula del g\'enero se obtiene que $\deg{\eu D}_{K|{{\ma Q}(x)}}
=4$. Por tanto $\deg{\eu D}_{K|{{\ma Q}(x)}}=\pL_{\infty}\pL_0
\pL_1\pL_{-1}$. Se sigue que 
\[
|M|=4=|I_{\K,0}|=h_{\K}=4<\infty.
\]
\end{ejemplo}

\subsection{La funci\'on zeta}\label{SRam1.2}

Sea $k=\F$ un campo finito y sea $\K$ un campo de funciones sobre 
$k$. Los resultados y definiciones que presentamos aqu\'i pueden
consultarse en \cite[Ch. 6 y 7]{Vil2006}.

\begin{teorema}[F.K. Schmidt]\label{T10.1.2.1} Sea
\[
\rho_{\K}=\min\{n\in{\ma N}\mid \text{existe ${\eu A}
\in D_{\K}$ con $\deg_{\K}{\eu A}=n$}\}.
\]
Entonces
$\rho_{\K}=1$.
\end{teorema}

\begin{proof}
Ver \cite[Theorem 6.3.8]{Vil2006}. 
$\fin$
\end{proof}

En particular tenemos que $\deg\colon D_{\K}\lra {\ma Z}$,
${\eu A}\longmapsto \deg_{\K}{\eu A}$, es suprayectiva.

\begin{definicion}\label{D10.1.2.2} Se define la 
{\em funci\'on zeta\index{funci\'on zeta}} de
$\K$, $\zeta_{\K}(s)$ por
\[
\zeta_{\K}(s)= \sum_{{\eu A}\text{\ entero}}\frac{1}{(\N {\eu A})^s}=
\sum_{{\eu A}\text{\ entero}} q^{-d_{\K}({\eu A}) s},
\]
esto es, $\N{\eu A}=q^{d_{\K}({\eu A})}$.
\end{definicion}

Se tiene que $\zeta_{\K}(s)$ converge absolutamente y uniformemente
en conjuntos compactos $\{s\in{\ma C}\mid \partereal s>1\}$.

Sea $u=q^{-s}$ y $Z_{\K}(u)=\zeta_{\K}(s)$. 

\begin{teorema}\label{T10.1.2.3}
Se tiene
\[
Z_{\K}(u)=\frac{P_{\K}(u)}{(1-u)(1-qu)}
\]
donde $P_{\K}(u)\in {\ma Z}[u]$ es un polinomio de grado $2g$. 
Adem\'as $P_{\K}(1)=h=h_{\K}$ es el n\'umero de clase de $\K$.
En particular $Z_{\K}(u)$ tiene un polo simple para $u=1$.

Por otro lado tenemos la f\'ormula del producto
\[
\zeta_{\K}(s)=\prod_{\pK\in{\ma P}_{\K}}(1-(\N\pK)^{-s})^{-1}\quad\text{para}
\quad \partereal s>1.
\]
\end{teorema}

\begin{proof}
Ver \cite[Theorems 6.3.5 y 6.3.7 y Corollary 6.3.6]{Vil2006}.
$\fin$
\end{proof}

\begin{definicion}\label{D10.1.2.4}
Sea $n\in{\ma N}$. Se define $A_n=A_n({\K})$ como el n\'umero de
divisores enteros de grado $n$ y sea $N_n=N_n({\K})$ el n\'umero
de divisores primos de grado $n$.
\end{definicion}

Se tiene
\begin{gather*}
Z_{\K}(u)=\sum_{n=0}^{\infty} A_n u^n \quad\text{y}\\
A_n=h_{\K}\big(\frac{q^{n-g_{\K}+1}-1}{q-1}\big)\quad \text{para}
\quad n>2g_{\K}-2
\end{gather*}
(ver \cite[Theorem 6.2.6]{Vil2006}).

\begin{definicion}\label{D10.1.2.5} Un {\em caracter de orden finito\index{caracter
de orden finito}} $\chi$ es un homomorfismo de grupos $\chi\colon I_{\K}\lra
\*{\ma C}$ tal que existe $n\in{\ma N}$ tal que $\chi^n=1$. En particular
$\chi(I_{\K})\subseteq W_n=\{\xi\in{\ma C}\mid\xi^n=1\}$.
\end{definicion}

Para un caracter de orden finito $\chi$ se define la {\em serie $L$ asociada 
a $\K$ y a $\chi$\index{serie $L$}} por
\[
L(s,\chi,\K)=\sum_{{\eu A}\text{\ entero}}\chi({\eu A})\frac{1}{\N({\eu A})^s},
\quad s\in {\ma C},\quad \partereal s>1.
\]

La serie converge absolutamente y uniformemente por compactos de 
$\{s\in{\ma C}\mid \partereal s>1\}$. Adem\'as tenemos la f\'ormula
del producto:
\[
L(s,\chi,\K)=\prod_{\pK\in{\ma P}_{\K}}\big(1-\frac{\chi(\pK)}{(\N \pK)^s}\big)^{-1},
\quad \partereal s>1.
\]

\begin{teorema}[Ecuaciones funcionales]\label{T10.1.2.6}
Sea $Z_{\K}(u)=\frac{P_{\K}(u)}{(1-u)(1-qu)}$, $P_{\K}(u)=
a_0+a_1u+\cdots+a_{2g}u^{2g}$ con $a_0=1$, $a_{2g}=q^g$.
Entonces $a_i=A_i-(q+1)A_{i-1}+qA_{i-2}$ donde definimos
$A_{-1}=A_{-2}=0$ y se tiene $A_0=1$ y en particular
$a_1=A_1-(q+1)=N_1-(q+1)$.

Adem\'as $a_{2g-i}=a_iq^{g-i}$, $0\leq i\leq 2g$ y
\[
q^{s(g-1)}\zeta_{\K}(s)=q^{(1-s)(g-1)} \zeta_{\K}(1-s)
\quad \text{para toda}\quad s\in{\ma C}.
\]

Para un caracter $\chi$ de orden finito, se tiene
\[
q^{s(g-1)}L(s,\chi,\K)=\chi(W) q^{(1-s)(g-1)} L(1-s,\bar{\chi},\K)
\]
donde $W=W_{\K}$ es la clase can\'onica de $\K$.
\end{teorema}

\begin{proof}
Ver \cite[Theorems 6.4.1, 6.4.3 y 6.4.6]{Vil2006}.
$\fin$
\end{proof}

Finalmente tenemos:

\begin{teorema}\label{T10.1.2.7} Sea $\K/k$ un campo de funciones,
$k=\F$, $l={\ma F}_{q^r}$, $L=\K l$. Sea $\chi_j$ el caracter
de $\K$ que satisface $\chi_j(C)=e^{2\pi i j/r}$ para cualquier clase
$C$ de grado $1$. Entonces
\begin{gather*}
\zeta_L(s)=\prod_{i=1}^rL(s,\chi_j,\K)\quad \text{y}\\
Z_L(u^r)=\prod_{j=1}^rZ_{\K}(\zeta_r^j u)
\end{gather*}
con $u=q^{-s}$ y $\zeta_r=e^{2\pi i /r}$.
\end{teorema}

\begin{proof}
\cite[Theorems 6.4.7 y 7.1.6]{Vil2006}.
$\fin$
\end{proof}

\subsection{Hip\'otesis de Riemann\index{hip\'otesis de
Riemann}\index{Riemann!hip\'otesis de $\sim$}}\label{SRam1.3}

Se tiene que $Z_{\K}(u)=\frac{P_{\K}(u)}{(1-u)(1-qu)}$,
$P_{\K}(u)=\sum_{i=0}^{2g}a_iu^i$ con $a_{2g-i}= a_i
q^{g-i}$, $0\leq i\leq 2g$. Adem\'as $a_1=A_1-(q+1)$,
$A_1=N_1$. Sean $w_1^{-1},\ldots, w_{2g}^{-1}$ las
ra\'ices de $P_{\K}(u)$. Entonces
\[
P_{\K}(u)=\prod_{i=1}^{2g}(1-w_iu).
\]

\begin{proposicion}\label{P10.1.2.7(1)} Se tiene que las siguientes
condiciones son equivalentes:
\las
\item Los ceros de $\zeta_{\K}(s)$ est\'a en $\partereal s=\frac{1}{2}$.
\item Los ceros de $Z_{\K}(u)$ est\'an en el c\'irculo $|u|=q^{-1/2}$.
\item $|w_i|=\sqrt{q}$, $1\leq i\leq 2g$.
\end{list}
\end{proposicion}

\begin{proof}
Ver \cite[Theorem 7.1.8]{Vil2006}.
$\fin$
\end{proof}

La {\em hip\'otesis de Riemann} establece que las condiciones
de la Proposicion \ref{P10.1.2.8} se cumplen para cualquier
campo de funciones congruente $\K$.

La hip\'otesis de Riemann es equivalente a
\begin{gather}\label{E10.1.2.1.1}
|N_1-(q+1)|\leq 2g\sqrt{q}.
\end{gather}

La hip\'otesis de Riemann fue probada por A. Weil en 1940--1941.

Como consecuencia de la hip\'otesis de Riemann, tenemos que
si $g_{\K}=0$, entonces $\K$ es un campo de funciones racionales
pues $|N_1-(q+1)|=0$, esto es, $N_1=q+1>1$ y por tanto $\K$
es de funciones racionales (ver \cite[Theorem 4.1.7]{Vil2006}).

Deducimos otra expresi\'on para la funci\'on zeta que es conveniente
para la estimaci\'on de $h_{\K}$. Por la f\'ormula del producto, 
tenemos
\[
\zeta_{\K}(s)=\prod_{\pK\in{\ma P}_{\K}}\big(1-\frac{1}{(\N\pK)^s}
\big)^{-1}=\prod_{m=1}^{\infty}\big(1-\frac{1}{q^{ms}}\big)^{N_m}
\quad \text{para}\quad \partereal s>1,
\]
donde $N_m$ es el n\'umero de divisores primos de $\K$
de grado $m$.

Se sigue que
\begin{align*}
\ln\zeta_{\K}(s)&=\sum_{m=1}^\infty -N_m[\ln (1-q^{-ms})]
\igual_{\substack{\uparrow\\ \ln(1-x)\\=-\sum\limits_{n=1}^{\infty}
\frac{x^n}{n}}} \sum_{m=1}^{\infty}-N_m\Big\{-\sum_{t=1}^{
\infty}\frac{q^{-mts}}{t}\Big\}\\
&=\sum_{m=1}^{\infty}\sum_{t=1}^{
\infty}\frac{N_m}{t}q^{-mts},\\
\frac{\zeta'_{\K}(s)}{\zeta_{\K}(s)}&=\big(\ln\zeta_{\K}(s)\big)'
\igual_{\substack{\uparrow\\ (q^{as})'=\\ (a\ln q)q^{as}}}
\sum_{m=1}^{\infty}\sum_{t=1}^{\infty}\frac{N_m(-mt\ln q)q^{-mts}}{t}\\
&=-\ln q\sum_{m=1}^{\infty}\sum_{t=1}^{\infty} mN_mq^{-mts}
=-\ln q\sum_{n=1}^{\infty} c_n q^{-ns}\\
\intertext{donde}
c_n&=\sum_{mt=n}mN_m=\sum_{m|n} mN_m.
\end{align*}

Haciendo $u=q^{-s}$, se tiene 
\begin{gather*}
\zeta_{\K}(s)=\sum_{m=1}^{\infty}
\sum_{t=1}^{\infty}\frac{N_m}{t} u^{mt}=\sum_{r=1}^{\infty}
d_r\frac{u^r}{r}
\intertext{donde}
\frac{dr}{r}=\sum_{mt=r}\frac{N_m}{t}=\sum_{m|r}\frac{N_m}{r/m}=
\frac{1}{r}\sum_{m|r}mN_m,
\end{gather*}
esto es, $d_r=\sum_{m|r}mN_m$.

Por tanto
\[
Z_{\K}(u)=\frac{P_{\K}(u)}{(1-u)(1-qu)}=\exp\Big(\sum_{r=1}^{\infty}
d_r\frac{u^r}{r}\Big).
\]

Ahora bien, consideremos $\K_r=\K{\ma F}_{q^r}$ la extensi\'on
de constantes. Sea $\pK$ es un lugar de grado $1$ en $\K_r$.
Sea $\P=\pK\cap \K$ el lugar respectivo en $\K$. Se tiene
\[
d_{\K_r/\K}(\pK|\P)d_{\K}(\P)=[{\ma F}_{q^r}:\F]d_{\K_r}(\pK)= r\cdot 1=r.
\]
Por tanto $d_{\K}(\P)|r$. Esto es, $d_{\K}(\P)=m|r$.

Rec\'iprocamente, si $d_{\K}(\P)=m|r$, entonces por el Teorema
\ref{T6.1.4}, si $\pK$ es un lugar en $\K_r$ sobre $\P$, 
\[
d_{\K_r}(\pK)=\frac{d_{\K}(\P)}{\mcd(d_{\K}(\P),r)}=\frac{m}
{\mcd(m,r)}=\frac{m}{m}=1.
\]

Se sigue que $d_r=N_1(\K_r)$ 
es el n\'umero de lugares de grado
$1$ en $\K_r/{\ma F}_{q^r}$. Entonces
\[
Z_{\K}(u)=\exp\Big(\sum_{r=1}^{\infty} N_1(\K_r)
\frac{u^r}{r}\Big).
\]

Ahora bien si $w_1,\ldots, w_{2g}$ son los inversos de las
ra\'ices de $P_{\K}(u)$, puesto que $P_{\K}(u)\in{\ma Z}[u]$,
si $w\in\{w_1,\ldots,w_{2g}\}$, entonces $\bar{w}=\{w_1,
\ldots,w_{2g}\}$.

Se tiene $q^g=\prod_{i=1}^{2g}w_i$ y $N_1(\K)-(q+1)=\sum_{
i=1}^{2g}w_i$. Adem\'as $P_{\K}(w_i^{-1})=0\iff P_{\K}\big(
\frac{w_i}{q}\big)=0$ (\cite[Proposition 7.1.7]{Vil2006}).

As\'i, en caso de que $\frac{1}{w_i}=\frac{w_i}{q}\iff w_i=\pm
\sqrt{q}$. Adem\'as, el n\'umero de veces que $\sqrt{q}$
aparece es par as\'i como el n\'umero de veces que $-
\sqrt{q}$ lo hace. Por tanto podemos escribir
\[
P_{\K}(u)=\prod_{i=1}^g (1-w_i u)(1-\bar{w}_i u)
\]
y por otro lado $|w_i|=\sqrt{q}$ para toda $1\leq i\leq g$.

\begin{proposicion}\label{P10.1.2.8} Se tiene
\[
(\sqrt{q}-1)^{2g_{\K}}\leq h_{\K}\leq (\sqrt{q}+1)^{2g_{\K}}.
\]
\end{proposicion}

\begin{proof}
Se tiene $h_{\K}=h=P_{\K}(1)=|P_{\K}(1)|=\prod_{i=1}^{2g}
|1-w_i|$. Puesto que $|w_i|=\sqrt{q}$, entonces
$\sqrt{q}-1\leq |1-w_i|\leq \sqrt{q}+1$. El resultado es
ahora inmediato. $\fin$
\end{proof}

\begin{observacion}\label{O10.1.2.9}
Sea $k$ un campo arbitrario y sea $E/F$ una extensi\'on finita
de campos de funciones sobre $k$. Entonces los mapeos
conorma $\con_{F/E}\colon I_{F,0}\lra I_{E,0}$ y 
$\con_{F/E}\colon I_F\lra I_E$ (ambos mapeos tienen mismo
n\'ucleo) no tiene por que ser inyectivos. El Ejemplo
\ref{E10.1.2.10} da una familia general de ejemplos.
M\'as adelante lo daremos de manera m\'as expl\'icita
para campos de funciones congruentes (Ejemplo \ref{E10.1.2.12}).
\end{observacion}

\begin{ejemplo}\label{E10.1.2.10} Sea $k$ un campo
algebraicamente cerrado de caracter\'istica $p>0$. Sea $\ell$
un primo tal que $\ell\neq p$. Sea $F$ un campo de funciones
sobre $k$ y sea $E=F(\sqrt[\ell]{u})$ una extensi\'on de
Kummer de grado $\ell$. 
Entonces $\con_{F/E}\colon I_{F,0}\lra I_{E,0}$
satisface que $\ker\con_{F/E}=\frac{D_{F,0}\cap P_E}{P_F}
\subseteq \frac{D_{E,0}^G\cap P_E}{P_F}= \frac{P_E^G}
{P_F}$. Puesto que los n\'ucleos de los dos mapeos conorma:
$\con_{F/E}\colon I_{F,0}\lra I_{E,0}$ y 
$\con_{F/E}\colon I_{F}\lra I_{E}$ son el mismo, se tiene
$\ker \con_{F/E}=\frac{D_{\K}\cap P_E}{P_{\K}}$.

Sean $R=\{x\in \*E\mid \text{tal que\ }x^{\ell}\in F\}$ y $\varphi
\colon R\lra P_E$ dada por $\varphi(x)=(x)_E$. Sea $G=\Gal(
E/F)=\langle \sigma\rangle$. Se tiene para $x\in R$ que
$\sigma(x)=\zeta_{\ell}^i x$ para alguna $0\leq i\leq \ell-1$
y $(\sigma x)_E=(x)_E$. Por tanto $\varphi(R)\subseteq P_E^G$.

Ahora sea $(z)_E\in P_E^G$, por tanto $(\sigma z)_E=(z)_E$,
es decir, existe $\xi\in\* k$ tal que $\sigma z=\xi z$. Se sigue
que $\sigma^2(z)=\sigma(\xi z)=\sigma \xi \sigma z=\xi \xi z=
\xi^2 z$ y en general $\sigma^t z=\xi^t z$. En particular
$z=\sigma^{\ell}(z)=\xi^{\ell}z$ y por tanto $\xi=\zeta_{\ell}^j$
para alguna $0\leq j\leq \ell-1$. Se sigue que $\sigma z^{\ell}=
(\sigma z)^{\ell}=(\xi z)^{\ell}=\xi^{\ell}z^{\ell}=z^{\ell}$, esto es,
$z^{\ell}\in F$. De esta forma tenemos que $\varphi(R)=P_E^G$.

Por tanto $\tilde{\varphi}\colon R\lra P_E^G/P_F$, $\tilde{\varphi}
(x)=(x)_E\bmod P_F$ es suprayectivo y $\ker \tilde{\varphi}=
\*F$. Por lo tanto $\frac{R}{\*F}\cong \frac{P_E^G}{P_F}$.

Ahora si $x\in R\setminus F$, $x^{\ell}=y\in F$ y $E=F(x)=
F(v)$ con $v^{\ell}=u$. Entonces existen $j\in\{0,1,\ldots, \ell-1\}$
y $c\in F$ tales que $x=v^j c$. Por tanto $\frac{R}{F}\cong
\frac{{\ma Z}}{\ell {\ma Z}}=\langle\bar{v}\rangle=\{1,\bar{v},
\ldots \bar{v}^{\ell-1}\}$ con $v=\sqrt[\ell]{u}$.

Cuando $E/F$ es no ramificada se tiene que $D_{E}^G=
\langle \con_{F/E}\P\mid \P\in{\ma P}_F\rangle$ pues si
$v_{\pK}({\eu A})=\alpha$, entonces para toda $\tau\in G$
se tiene $v_{\pK^{\tau}}({\eu A})=v_{\pK}({\eu A}^{\tau^{-1}})=
v_{\pK}({\eu A})$.

Por otro lado, como $E/F$ es no ramificado,
$\con_{F/E}\P=\pK_1\cdots \pK_{\ell}$ y $D_{E}^G=D_F$
y $\ker\con_{F/E}=\frac{D_{F}\cap P_E}{P_F}=
\frac{D_{E}^G\cap P_E}{P_F}=\frac{R}{F}\cong
\frac{{\ma Z}}{\ell{\ma Z}}$. Por tanto $\con_{F/E}$
no es inyectiva.
\end{ejemplo}

De hecho, el mapeo $\con_{F/E}\colon I_{F,0}\lra I_{E,0}$ puede
ser trivial.

\begin{ejemplo}\label{E10.1.2.10(1)}
Se tiene que si $\K$ es un campo de g\'enero $0$ que no es
separablemente generado (ver \cite[Section 8.2]{Vil2006}), entonces
$\K$ contiene subcampos de g\'enero arbitrariamente grande
(\cite{LanTat52}). Se sigue que la caracter\'istica de $\K$ es $2$ y que
$\K=k(x,y)$ donde $y^2=ax^2+b$ con $[k(a^{1/2},b^{1/2}):k]=4$.
El campo $L=k(z,w)$ con $z=x^2$ y $w=x^{2n+1}+y$ con $n\geq 1$
tiene g\'enero $g_L=-\big[\frac{-n}{2}\big]-1$. Seleccionamos
$k$ una cerradura separable de $k_0:={\ma F}_2(U,V)$ con
$U, V$ variables independientes sobre ${\ma F}_2$. El campo
$k$ satisface nuestras hip\'otesis.

Ahora bien, se puede probar que $L$ tiene una infinidad de
lugares de grado $1$. Sea $\{\pK_m\}_{m=0}^{\infty}$ una 
suceci\'on infinita de lugares de grado $1$ de $L$. Entonces
$\{C_m:=\overline{\pK_m\pK_0^{-1}}\}_{m=1}^{\infty}$, es una
colecci\'on infinita de clases de grado $0$ distintas. En particular
$|I_{L,0}|=h_{L}=\infty$. Ahora el mapeo $\con{L/\K}\colon
I_{L,0}\lra I_{\K,0}$ es trivial pues $h_{\K}=1$.
\end{ejemplo}

En el caso de extensiones de constantes, en general se
tiene

\begin{proposicion}\label{P10.1.2.11}
Sea $\K/k$ un campo de funciones arbitrario y sea $L=\K\ell$ 
donde $\ell/k$ es separablemente generado. Entonces
$\con_{\K/L}\colon I_{\K,0}\lra I_{L,0}$ y $\con_{\K/L}\colon
I_{\K}\lra I_L$, son inyectivas.
\end{proposicion}

\begin{proof}
El resultado es inmediato de \cite[Corollary 8.5.10]{Vil2006}. La
idea es que para ${\eu A}\in D_{\K}$, $\ell_{\K}({\eu A})=\ell_L({\eu A})$.
Tomando $d_{\K}({\eu A})>\max \{2g_{\K}-2, 2g_L-2\}$,
\begin{align*}
\ell_{\K}({\eu A}^{-1})&=d_{\K}({\eu A})-g_{\K}+1,\\
\ell_L({\eu A}^{-1})&=d_L({\eu A})-g_L+1.
\end{align*}

Puesto que en este caso $d_{\K}({\eu A})=d_L({\eu A})$ y 
$\ell_{\K}({\eu A}^{-1})=\ell_L({\eu A}^{-1})$, se sigue que
$g_{\K}=g_L$.

Finalmente si ${\eu A}\in\ker\con_{\K/L}$, $\con_{\K/L}\colon
I_{\K}\stackrel{\varphi}{\lra} I_L$, ${\eu A}=(\alpha)$ en $L$, por 
lo que $\ell_L({\eu A}^{-1})=1$ y $d_L({\eu A})=0$. Se sigue
que $\ell_{\K}({\eu A}^{-1})=1$ y $d_{\K}({\eu A})=0$. Por lo tanto
${\eu A}$ es principal en $\K$ y $\varphi$ es inyectiva.
$\fin$
\end{proof}

Volveremos sobre la Proposici\'on \ref{P10.1.2.11} m\'as
adelante en el caso de campos de funciones congruentes.

Presentamos a continuaci\'on un ejemplo concreto sobre
la no inyectividad del mapeo conorma en el caso de 
campos de funciones congruentes (ver Ejemplo 
\ref{E10.1.2.10}).

\begin{ejemplo}\label{E10.1.2.12}
Sean $K=\F(T)$, $\ell$ un primo  con $\ell|q-1$, $P_1(T),
P_2(T)\in R_T=\F[T]$ dos polinomios irreducibles m\'onicos
distintos de grado $\ell$. Sean $\K=K(\sqrt[\ell]{P_1P_2})$
y $L=\K(\sqrt[\ell]{P_2})$. Se tiene que el
primo infinito de $K$ no se ramifica
ni en $\K/K$ ni en $L/\K$ (ver Teorema \ref{TRam1} m\'as
adelante). En $\K/K$ los primos ramificados son precisamente
$P_1(T)$ y $P_2(T)$. Adem\'as en $L/\K$ el \'unico posible primo
ramificado es $P_2(T)$ pero $L=\K(\sqrt[\ell]{P_1})=\K(
\sqrt[\ell]{P_2})$ por lo que $L/\K$ es una extensi\'on c\'iclica
no ramificada de grado $\ell$.

Sean $\P_1$ y $\P_2$ los divisores primos de $K$ correspondientes
a $P_1$ y a $P_2$ respectivamente y sea $\P_{\infty}$ el primo
infinito de $K$. Entonces $(P_i)=\frac{\P_i}{\P_{\infty}^{\ell}}$, 
$i=1,2$. Adem\'as $(-1)^{\deg P_1P_2}=(-1)^{2\ell}=1$. Se sigue
que $\P_{\infty}$ se descompone totalmente en $K(\sqrt[\ell]{P_1})
=K(-\sqrt[\ell]{P_1})=K(\sqrt[\ell]{(-1)^{\ell}P_1})$,
en $K(\sqrt[\ell]{P_2})
=K(-\sqrt[\ell]{P_2})=K(\sqrt[\ell]{(-1)^{\ell}P_2})$ y en 
$K(\sqrt[\ell]{P_1P_2})$ (ver Corolario \ref{C5.1.1'}).
En particular $\P_{\infty}$ se descompone totalmente en $L/K$.

Sean $\pK_i=\P_i^{\ell}$ en $\K/K$, $i=1,2$. Sea ${\eu A}=
\frac{\pK_1}{(\con_{K/\K}\P_{\infty})^{\ell}}$ en $\K$. Se tiene $\con_{\K/L}
({\eu A})=(\sqrt[\ell]{P_1})_L$ el cual es un divisor principal.

Veamos que ${\eu A}$ no es principal $\K$. En caso de serlo,
digamos ${\eu A}=(\beta)_{\K}$ en $\K$, entonces $(\beta)_L=
(\sqrt[\ell]{P_1})_L$ por lo que $\beta=\xi\sqrt[\ell]{P_1}$
con $\xi$ una constante de $L$. Ahora bien, puesto que $\P_{\infty}$ se 
descompone totalmente en $L/K$, el campo de constantes de
$L$ es $\F$, esto es, $\xi\in\F$. Se sigue que $\sqrt[\ell]{P_1}\in \K$ lo cual es
absurdo pues si este fuese el caso, entonces $\sqrt[\ell]{P_2}=
\frac{\sqrt[\ell]{P_1P_2}}{\sqrt[\ell]{P_1}}\in\K$ y $\K=K(\sqrt[\ell]
{P_1}, \sqrt[\ell]{P_2})=L$.

As\'i $\bar{{\eu A}}\in \ker \con_{\K/L}$ y $\bar{{\eu A}}\neq 1$
por lo que $\con_{\K/L}$ no es inyectiva.
\end{ejemplo}

Volvemos al caso de extensiones de constantes, ahora 
expl\'icitamente para el caso de campos de funciones congruentes.

\begin{teorema}\label{T10.1.2.14}
Sea $\K$ un campo global de funciones con campo de constantes
$k=\F$. Sea $L/\K$ una extensi\'on de constantes, esto es
$L=\K l$. Entonces
$\con_{\K/L}\colon I_{\K,0}\lra I_{L,0}$ es inyectiva.
En particular $h_{\K}|h_L$.
\end{teorema}

\begin{proof}
Damos dos demostraciones. Para la primera, consideremos
$G=\Gal(L/\K)\cong\Gal(l/k)$ el cual es un grupo c\'iclico.
Ahora si existe $\K\subsetneqq E\subsetneqq L$, entonces $L/E$
y $E/\K$ son ambas extensiones de constantes y por 
inducci\'on en $[L:\K]$ podemos suponer que $\con_{\K/E}$ y
$\con_{E/L}$ son inyectivas. De esta forma el caso general
se reduce al caso $[L:\K]=\ell$, con $\ell$ primo.

Sea ${\eu A}\in I_{\K,0}$ tal que $\varphi(\bar{\eu A})=\con_{
\K/L}(\bar{\eu A})=(1)$, donde $\bar{\eu A}={\eu A}\bmod P_{\K}$.
Esto es, ${\eu A}$ es principal en $L$: ${\eu A}=(\alpha)$.
Suponemos que $\alpha\notin \K$ pues de otra forma 
${\eu A}\in P_{\K}$ y $\bar{\eu A}=1$. Sea $G=\langle\sigma\rangle$.

Sea $\sigma\alpha=\xi\alpha$ con $\xi\in\ell$ pues ${\eu A}^{
\sigma}=(\sigma\alpha)=(\alpha)={\eu A}$. Tomando normas
obtenemos
\[
\N_{L/\K}(\alpha)=\N_{L/\K}(\sigma\alpha)=\N_{L/\K}(\xi\alpha)=
\N_{L/\K} (\xi) \N_{L/\K} (\alpha)
\]
de donde se sigue que $\N_{L/\K}(\xi)=1=\N_{l/k}(\xi)$.

Por el Teorema 90 de Hilbert (ver Teorema \ref{CClaseT1.5.7}),
existe $\delta\in l$ tal que
$\xi=\delta^{\sigma-1}$, esto es, $\sigma \delta =\xi \delta$.
Por tanto
\begin{gather*}
\sigma(\delta^{-1}\alpha)=(\sigma\delta)^{-1}(\sigma\alpha)=
(\xi\delta)^{-1}(\xi\alpha)=\delta^{-1}\xi^{-1}\xi\alpha=\delta^{-1}\alpha
\intertext{por lo cual}
\beta=\delta^{-1}\alpha\in\K\quad\text{y}\quad (\alpha)=(\delta^{-1}
\alpha)=(\beta)\in P_{\K}.
\end{gather*}

De esta forma obtenemos que $\con_{\K/L}$ es inyectiva.

Una segunda demostraci\'on es usando otra vez el Teorema
90 de Hilbert. Como hemos visto si $\varphi$ denota la 
conorma, entonces $\ker \varphi\cong \frac{P_L^G}{P_{\K}}$.
Se tiene la sucesi\'on exacta
\begin{gather*}
1\lra \* l\lra\* L\lra P_L\lra 1.
\intertext{Se sigue que}
1\lra (\* l)^G=\* k\lra (\* L)^G=\* \K\lra P_L^G\lra\\
\lra H^1(G, \* l)
\lra H^1(G,\*L)=\{1\}.
\end{gather*}

De esta sucesi\'on exacta y del Teorema 90 de Hilbert aplicado
a $\* l$, se sigue que $\ker \con_{\K/L}=\frac{P_L^G}{P_{\K}}\cong
H^1(G,\* l)=\{1\}$ de donde obtenemos que $\con_{\K/L}$ es
inyectiva. $\fin$
\end{proof}

\begin{observacion}\label{O10.1.2.14(1)}
En la literatura se encuentran otras demostraciones del
Teorema \ref{T10.1.2.14}. Una de ellas es consecuencia
es un caso particular de la Proposici\'on \ref{P10.1.2.11}.
Otra m\'as se puede dar de la siguiente forma. Por los
Teoremas \ref{T6.1.3(1)} y \ref{T6.1.3(2)} se tiene que
$g_{\K}=g_L$ y para todo ${\eu A}\in D_{\K}$, $\ell_L(
{\eu A})=\ell_{\K}({\eu A})$. 

En particular, si $W_{\K}$ es la clase can\'onica de $\K$,
$W_{\K}$ es la \'unica clase $\K$ de grado $2g_{\K}-2=
2g_L-2$ tal que $\ell_{\K}(W_{\K}^{-1})=d_{\K}(W_{\K})-
g_{\K}+1+\ell_{\K}(W_{\K}^{-1}W_{\K})=2g_{\K}-2-g_{\K}
+1+1=g_{\K}$.

Por tanto $\con_{\K/L}(W_{\K})$ es la \'unica clase de grado
$2g_L-2$ tal que $\ell\big((\con_{\K/L}W_{\K})^{-1}\big)
=g_L=g_{\K}$. Se sigue que $\con_{\K/L}W_{\K}=W_L$.

Por tanto $\ell_L({\eu A}W_L^{-1})=\ell_{\K}({\eu A}W_{\K}^{-1})$
para toda ${\eu A}\in D_{\K}$.

Sea ${\eu A}\in\ker \con_{\K/L}$, $\con_{\K/L}\colon I_{\K}\lra
I_L$. Entonces $d_{\K}({\eu A})=d_L({\eu A})=0$ y ${\eu A}$
se hace principal en $L$. Por tanto, la clase de ${\eu A}$
en $L$ es $P_L$ y por tanto $\ell_L({\eu A})=1$.

De esta forma tenemos que ${\eu A}$ es una clase de grado
$0$ tal que $\ell_{\K}({\eu A})=\ell_L({\eu A})=1$ por lo que
${\eu A}\in P_{\K}$, esto es, ${\eu A}$ es principal en $\K$ y
$\con_{\K/L}$ es inyectiva.
\end{observacion}

A pesar de que el mapeo conorma no es inyectivo en general,
M. Madan \cite{Mad69, Mad70} prob\'o:

\begin{teorema}\label{T10.1.2.15} Sea $\K$ un campo de funciones 
congruente y $L/\K$ una extensi\'on finita de Galois. Entonces
$h_{\K}|h_L$.
\end{teorema}

\begin{proof}
Si probamos el resultado para $L/\K$ geom\'etrica, entonces
en general si $l$ es el campo de constantes de $L$, se tiene
que $\K\subsetneqq E=\K l\subseteq L$. En este caso $E/\K$
es una extensi\'on de contantes por lo que $h_{\K}|h_E$, y
$L/E$ es una extensi\'on geom\'etrica. Por tanto podemos
suponer que $L/\K$ es una extensi\'on geom\'etrica.

Ahora bien, podemos suponer que $G=\Gal(L/\K)$ es un
grupo simple pues si $\K\subsetneqq M\subsetneqq L$
con $M/\K$ tambi\'en es una extensi\'on de Galois, por
inducci\'on en $[L:\K]$ podemos suponer
que $h_{\K}|h_M$ y $h_M|h_L$ lo
cual implica que $h_{\K}|h_L$.

As\'i, suponemos que $L/\K$ es una extensi\'on geom\'etrica
simple. Sea $N=\ker \con_{\K/L}$. Si $N=\{\Id\}$, entonces
el resultado es inmediato. Suponemos que $N$ es no trivial.
Sea ${\eu A}$ un representante de una clase $\bar{\eu A}\in
N$.

Sea ${\eu A}=(\alpha)$ principal en $L$. Como ${\eu A}^{\sigma}
={\eu A}$ para toda $\sigma\in G$, se tiene $\sigma \alpha=
X_{\eu A}(\sigma)\alpha$  con $X_{\eu A}(\sigma)\in \* k$, donde
$k$ es el campo de constantes tanto de $\K$ como de $L$.

Se tiene que $X_{\eu A}\colon G\lra \* k$ es un homomorfismo
de grupos pues
\begin{align*}
X_{\eu A}(\sigma \mu)(\alpha)&=(\sigma\mu)(\alpha)=\sigma(\mu(
\alpha))=\sigma(X_{\eu A}(\mu)\alpha)=X_{\eu A}(\mu)\sigma
(\alpha)\\
&=X_{\eu A}(\mu)X_{\eu A}(\sigma)\alpha=
X_{\eu A}(\sigma)X_{\eu A}(\mu)\alpha,
\end{align*}
por lo que $X_{\eu A}(\sigma\mu)=X_{\eu A}(\sigma)X_{\eu A}
(\mu)$.

Veamos que $X_{\eu A}$ depende \'unicamente de la clase
$\bar{\eu A}$. Si ${\eu B}={\eu A}(\beta)$ con $\beta\in\K$. Entonces
${\eu B}=(\beta\alpha)$ en $L$ por lo que
\begin{align*}
\sigma(\beta\alpha)&=\sigma(\beta)\sigma(\alpha)=\beta
\sigma(\alpha)=\beta X_{\eu A}(\sigma)\alpha\\
\sigma(\beta\alpha)&=X_{\eu B}(\sigma)(\beta\alpha)=
\beta X_{\eu B}(\sigma)\alpha,
\end{align*}
esto es, $X_{\eu B}=X_{\eu A}$.

Ahora bien, si $X_{\eu A}(\sigma)=1$ para toda $\sigma\in G$, 
entonces $\alpha\in\K$ y por tanto ${\eu A}$ es principal en
$\K$, es decir, $\bar{\eu A}=1$. Se sigue que el mapeo
$\varphi\colon \bar{\eu A}\longmapsto X_{\eu A}$ es una
inyecci\'on de $N$ en $\Hom(G,\* k)$ por lo que $\Hom(G,
\* k)\neq \{1\}$.

Sea $\mu\in\Hom(G,\*k)$, $\mu\neq 1$. Entonces $\ker\mu
\neq G$ y $\ker \mu\normal G$. Se sigue que $\ker\mu=\{1\}$
pues $G$ es simple. Se sigue que $G$ es un grupo c\'iclico
de orden primo $\ell$.

Ya que $G\cong \mu(G)\subseteq \*k$, se tiene $\ell||\*k|=q-1$
donde suponemos $k\cong \F$. Tambi\'en se sigue que $N\cong
C_{\ell}$ es el grupo c\'iclico de orden $\ell$ y que $L/\K$
es una extensi\'on de Kummer de grado $\ell$.

En esta situaci\'on se puede concluir la demostraci\'on de
forma anal\'itica o de forma algebraica. Sea $G=\Gal(L/\K)$
y sea $\chi\in\hat{G}=\Hom(G,\*{\ma C})$ un generador de
$\hat{G}\cong G\cong C_{\ell}$. Para $\P\in {\ma P}_{\K}$
no ramificado, sea $\chi(\P)=\artin{L/\K}{\P}$ es el s\'imbolo
de Artin. Si $\P$ es ramificado se define $\chi(\P)=0$.
Tambi\'en entendemos $\chi^0=\Id$, esto es, $\chi^0(\P)=1$
para todo $\P\in{\ma P}_{\K}$.

Consideramos la serie $L$ dada por
\[
L(s,\chi)=\prod_{\P\in{\ma P}_{\K}}\Big(1-\frac{\chi(\P)}
{(\N\P)^s}\Big)^{-1}.
\]
Si $L_L(s)$ y $L_{\K}(s)$ son los numeradores de las funciones
zeta de $L$ y de $\K$ respectivamente, se tiene
\begin{gather*}
L_L(s)=L_{\K}(s)\prod_{i=1}^{\ell-1}L(s,\chi^i),
\intertext{ver \cite[(6)]{Mad69}. Por tanto}
h_L=h_{\K}\prod_{i=1}^{\ell-1}L(0,\chi^i).
\end{gather*}

Como $\prod_{i=1}^{\ell-1}L(0,\chi^i)$ es, por un lado un n\'umero
racional, a saber $h_L/h_{\K}$, y por otro, es un n\'umero algebraico,
si sigue que es un n\'umero entero de donde obtenemos
$h_{\K}|h_L$.

Tambi\'en podemos concluir de manera algebraica. En \cite{Mad70}
Madan prueba usando cohomolog\'ia que si $l/k$ es la extensi\'on
de constantes de grado $\ell$, entonces si $G\cong\bar{G}\cong
\Gal(L l/L)\cong\Gal(\K l/\K)\cong\Gal(l/k)$, se tiene que
\[
\big|I_{\bar{\K},0}^{\bar{G}}\big|=|I_{\K,0}|=\Big|\Big(\frac{I_{\K l,0}}
{\bar{N}}\Big)^{\bar{G}}\Big|=h_{\K}
\]
donde $\bar{N}=\ker\con_{\K l/Ll}\colon I_{\K l,0}\lra I_{Ll,0}$
el cual es no trivial.

Finalmente, puesto que $\Big(\frac{I_{\K l,0}}
{\bar{N}}\Big)^{\bar{G}}\subseteq I_{L,0}$, se sigue que $h_{\K}|h_L$.
$\fin$
\end{proof}

\begin{observacion}\label{O10.1.2.15'} Y. Aubry dio otra demostraci\'on
del Teorema \ref{T10.1.2.15} en \cite[Proposition 2.3]{Aub96} usando
m\'etodos geom\'etricos.
\end{observacion}

Como consecuencia obtenemos la Observaci\'on \ref{O6.3.15}.

\begin{corolario}\label{C10.1.2.16} Si $M,N\in\F[T]$ son dos polinomios
tales que $N$ divide a $M$, entonces si $h_M$ y $h_N$ denotan
los n\'umeros de clase de los campos de funciones ciclot\'omicos
$\cicl M{}$ y $\cicl N{}$ respectivamente, se tiene que
$h_N$ divide a $h_M$.
\end{corolario}

\begin{proof} Es consecuencia inmediata del Teorema \ref{T10.1.2.15}
y del hecho de que $\cicl M{}/\cicl N{}$ es una extensi\'on de Galois.
$\fin$
\end{proof}

\section{Dominios Dedekind y torsi\'on}\label{SDed.2}

Regresaremos a este tema en el Cap\'itulo \ref{DrinfeldCh15}, 
Secci\'on \ref{DrinfeldC1}.
\begin{definicion}\label{D10.2.1.1}
Un dominio entero (conmutativo con unidad) $D$ se llama un
{\em dominio o anillo Dedekind\index{anillo Dedekind}\index{dominio
Dedekind}} si $D$ no es un campo y adem\'as satisface las
siguientes tres condiciones:
\las
\item Todo ideal primo ${\mathcal P}$ no cero es maximal, esto es,
$D$ es de dimensi\'on uno.
\item $D$ es noetheriano.
\item $D$ es enteramente cerrado, esto es, si $\K=\coc D$ es el
campo de cocientes de $D$ y si $x\in \K$ satisface una relaci\'on
$x^n+a_{n-1}x^{n-1}+\cdots+a_1x+a_0=0$ con $a_i\in D$, $0
\leq i\leq n-1$, entonces $x\in D$.
\end{list}
\end{definicion}

\begin{ejemplos}\label{Ej10.2.1.2}
${\ma Z}$, $k[x]$ con $k$ un campo arbitrario, $D$ un dominio de ideales
principales (DIP), ${\mathcal O}_{\K}$ el anillo de enteros de cualquier
campo num\'erico $\K$, son todos dominios Dedekind.
\end{ejemplos}

Dado un dominio Dedekind y $\K=\coc D$ su campo de cocientes, un
$D$--m\'odulo $M\neq 0$ con $M\subseteq \K$ se llama {\em ideal
fraccionario\index{ideal fraccionario}} si $M$ es finitamente generado
como $D$--m\'odulo. Equivalentemente, existe $d\in D$, $d\neq 0$ tal
que $dM\subseteq D$. Como ejemplos de ideales fraccionarios
tenemos a los ideales no ceros usuales
de $D$.

\begin{teorema}\label{T10.2.1.3} Si $D$ es un dominio de Dedekind, todo ideal
fraccionario ${\eu A}$ se escribe de manera \'unica como 
\[
{\eu A}={\mathcal P}_1^{\alpha_1}\cdots {\mathcal P}_r^{\alpha_r}
\]
con ${\mathcal P}_i$ ideales primos no cero de $D$, y $\alpha_i\in
{\ma Z}$, donde para un ideal primo no cero ${\mathcal P}$ de $D$,
${\mathcal P}^{-1}=\{x\in \K\mid x{\mathcal P}\subseteq D\}$. 
\end{teorema}

\begin{proof} \cite[Theorem 5.7.4]{Vil2006}. $\fin$
\end{proof}

\begin{teorema}\label{T10.2.1.4}
Sea $A$ un dominio Dedekind, $\K=\coc A$. Sea $L$ una extensi\'on
finita de $\K$ con $[L:\K]=n$. Sea $B=\{\alpha\in L\mid \Irr(\alpha,x,\K)\in A[x]\}$
la cerradura entera de $A$ en $L$. Entonces $B$ es un dominio Dedekind.
\end{teorema}

\begin{proof} \cite[Theorem 5.7.7]{Vil2006}. $\fin$
\end{proof}

\begin{teorema}\label{T10.2.1.5} Si $A$ es un dominio Dedekind, $\K=\coc A$
y $A\subseteq B\subsetneqq \K$, con $B$ anillo. Entonces $B$ es dominio
Dedekind.
$\fin$
\end{teorema}

\begin{proposicion}\label{P10.2.1.6}
Todo ideal no cero de un dominio Dedekind puede ser generado por a lo
m\'as dos elementos.
\end{proposicion}

\begin{proof} \cite[Exercise 5.10.34]{Vil2006}. $\fin$
\end{proof}

Sea $\K/k$ un campo arbitrario de funciones y sean ${\mathcal P}_1,\ldots,
{\mathcal P}_r$ un n\'umero finito de lugares de $\K$ ($r\geq 1$). Sea 
${\mathcal O}=\cap_{{\mathcal P}\notin\{{\mathcal P}_1,\ldots,{\mathcal P}_r\}}
{\mathcal O}_{\mathcal P}=\{x\in \K\mid v_{\mathcal P}(x)\geq 0 \text{\ para
todo lugar\ }{\mathcal P}\notin\{{\mathcal P}_1,\ldots,{\mathcal P}_r\}\}$.

\begin{teorema}\label{T10.2.1.7} 
El anillo ${\mathcal O}$ es un dominio Dedekind.
\end{teorema}

\begin{proof} 
Por el Teorema de Riemann--Roch, existe $x_i\in \K$ tal que el divisor de
polos de $x_i$ es $\eta_{x_i}={\mathcal P}_i^{\alpha_i}$, $\alpha_i\geq 1$. Sea
$y:=\sum_{i=1}^r x_i$. Entonces
$\eta_y=\prod_{i=1}^r {\mathcal P}_i^{\alpha_i}$. Veamos que ${\mathcal O}$
es la cerradura entera de $k[y]\subseteq k(y)$ en $\K$.
\[
\xymatrix{
R={\mathcal O}_{\K}\ar@{-}[d]\ar@{-}[r]&\K\ar@{-}[d]\\ k[y]\ar@{-}[r]&k(y)}
\]

Sea $R:=\{\alpha\in \K\mid \Irr(\alpha,T,k(y))\in k[y][T]\}$ la cerradura entera
de $k[y]$ en $\K$. Entonces, si $\alpha\in R$, se tiene $\alpha^n+a_{n-1}
\alpha^{n-1}+\cdots+a_1\alpha+a_0=0$ con $a_i\in k[y]$. Notemos que la
conorma $\con_{k(y)/\K}{\mathcal P}_{\infty}
={\mathcal P}_1^{\alpha_1}\cdots {\mathcal P}_r^{
\alpha_r}$ puesto que $\eta_y={\mathcal P}_{\infty}$ en $k(y)$. Si ${\mathcal P}
\notin\{{\mathcal P}_1,\ldots,{\mathcal P}_r\}$, $v_{\mathcal P}(a_i)\geq 0$ puesto
que ${\mathcal P}|_{k(y)}\neq {\mathcal P}_{\infty}$. Si
tuvi\'esemos $v_{\mathcal P}(\alpha)
<0$, entonces $v_{\mathcal P}(\alpha^n)=nv_{\mathcal P}(\alpha)<
iv_{\mathcal P}(\alpha)
\leq iv_{\mathcal P}(\alpha)+v_{\mathcal P}(a_i)=v_{\mathcal P}(a_i\alpha^i)$,
$0\leq i\leq n-1$.
Por lo tanto 
\[
\infty =v_{\mathcal P}(0)=v_{\mathcal P}(\alpha^n+a_{n-1}\alpha^{n-1}+\cdots
+a_1\alpha+a_0)=nv_{\mathcal P}(\alpha)<0,
\]
lo cual es absurdo. Se sigue que $v_{\mathcal P}(\alpha)
\geq 0$ y $\alpha\in{\mathcal O}$.
De esta forma obtenemos que $R\subseteq {\mathcal O}$. Probaremos la
otra contenci\'on, sin embargo, es suficiente esta primera contenci\'on pues $R$
es dominio Dedekind y $R\subseteq {\mathcal O}\subseteq \coc R=\K$ (ver
Teorema \ref{T10.2.1.5}).

Sea $\alpha\in {\mathcal O}$. Por tanto $v_{\mathcal P}(\alpha)\geq 0$ para
toda ${\mathcal P}\notin\{{\mathcal P}_1,\ldots, {\mathcal P}_r\}$. Sea $f(T)=
\Irr(\alpha,T,k(y))=T^m+b_{m-1}T^{m-1}+\cdots+b_1T+b_0\in k(y)[T]$. Entonces
$\alpha^m+b_{m-1}\alpha^{m-1}+\cdots+b_1\alpha+b_0=0$ con $b_0\neq 0$.

Sea $\tilde{\K}$ la cerradura normal de $\K/k(y)$. Sean $\alpha=\alpha^{(1)},
\alpha^{(2)},\ldots, \alpha^{(m)}$ los conjugados de $\alpha$, esto es, 
$f(T)=\prod_{i=1}^m(T-\alpha^{(i)})$ ya sea separable o no. Entonces $b_i$
es una funci\'on sim\'etrica en $\alpha^{(1)},\ldots,\alpha^{(m)}$, $b_i=
\sum \alpha^{(j_1)}\cdots \alpha^{(j_{m-i})}$ y $v_{\eu p}(\alpha^{(j)})\geq 0$
donde ${\eu p}$ es un primo en $\tilde{\K}$ sobre ${\mathcal P}$. Por lo tanto
$v_{\eu p}(b_i)\geq 0$, de donde $v_{\eu q}(b_i)\geq 0$ para todo 
primo de $k(y)$ con ${\eu q}\neq {\mathcal P}_{\infty}$. Se sigue
que $b_i\in k[y]$ y por tanto $\alpha\in R$ y por tanto $\o\subseteq R$.

Se sigue que ${\mathcal O}=R$. Finalmente, puesto que $k[y]$ es un DIP,
$k[y]$ es dominio Dedekind y por el Teorema \ref{T10.2.1.4} se sigue que
$R={\mathcal O}$ es dominio Dedekind.
$\fin$
\end{proof}

\subsection{N\'umero de clase de ${\mc O}_S$}\label{SRamDed2.1}

Dado un dominio Dedekind $D$, el {\em n\'umero 
de clase de $D$\index{numero de clase de un dominio Dedekind@n\'umero de
clase de un dominio Dedekind}} se define 
como $|Cl_D|=h_D=\big|\frac{D_D}{P_D}\big|$,
donde $D_D$ es el grupo de ideales fraccionarios 
de $D$ y $P_D$ son los ideales 
fraccionarios principales.

En general $h_D$ no necesariamente es finito. De hecho Luther Claborn \cite{Cla66}
prob\'o que todo grupo abeliano es el grupo de clase de un dominio Dedekind.

Sabemos que el grupo de clases de grado $0$ 
de un campo de funciones globales
$I_{\K,0}$ es finito y $|I_{\K,0}|=h_{\K}$ (Teorema \ref{T10.1.1.5}). 
Sea $S=\{\pK_1,\ldots,\pK_s\}$ un conjunto
finito y no vac{\'\i}o, esto es, $s\geq 1$, de lugares de $\K$. Sea 
${\mc O}_S:=\{x\in \K\mid v_{\pK}(x)
\geq 0\text{\ para todo $\pK\notin S$}\}
=\bigcap_{\pK\notin S}\o_{\pK}$. Entonces ${\mc O}_S$
es la cerradura entera de $\F[x]$,
donde $x$ es un elemento de $\K$ con $\eta_x
=\pK_1^{\alpha_1}\cdots\pK_s^{\alpha_s}$
para algunos $\alpha_i>0$, $i=1,\ldots, s$. Entonces 
${\mc O}_S$ es un dominio Dedekind (Teorema \ref{T10.2.1.4}).

Sean $Cl_S=Cl_{{\mc O}_S}$ el grupo de clases de 
${\mc O}_S$ y $h_{S}=|Cl_{_S}|$, 
finito o infinito. Probaremos que
$h_S$ es finito. Sea $D_{\K}$ el grupo de divisores de 
$\K$, $D_{\K,0}$ el subgrupo de $D_{\K}$
de los divisores de grado $0$, $P_{\K}$ el subgrupo de $D_{\K,0}$ de los divisores
principales, $I_{\K}=D_{\K}/P_{\K}$ el 
grupo de clase de divisores y $I_{\K,0}=D_{\K,0}/P_{\K}$
el subgrupo de clase de divisores de grado $0$. 

Sean $D_S=D_{\K}/\langle S\rangle=
D_{\K}/(\oplus_{i=1}^s{\ma Z}\pK_i)$ 
el subgrupo de $D_{\K}$ generado por los divisores $\pK
\notin S$. Para $x\in \* \K$ se define el {\em $S$--divisor de $x$\index{S@$S$--divisor}}
por $(x)_S=\prod_{\pK\notin S}\pK^{v_{\pK}(x)}$ y $P_S=\{(x)_S\mid x\in\* \K\}$ es
el {\em grupo de los $S$--principales\index{grupo de los $S$--principales}}, $P_S
\subseteq D_S$. Sea $C_S=D_S/P_S$. 
De hecho $Cl_S\cong C_S$ por la correspondencia
de divisores. Esto lo veremos en detalle m\'as
adelante. Sea $D_{\K}(S)=\langle \pK\in S\rangle$ y 
$P_{\K}(S)=P_{\K}\cap D(S)$ es el
subgrupo de divisores principales con soporte en $S$.

Notemos que en $\K$, el divisor $(x)_S$ no necesariamente es
principal, de hecho, ni de grado $0$. Por ejemplo, si $\pK\notin S$
y ${\eu q}\in S$. Entonces $\frac{{\pK}^{\deg {\eu q}}}{{\eu q}^{\deg \pK}}$
es de grado $0$. Existe $m\in{\ma N}$ tal que $\Big(\frac{{\pK}^{\deg 
{\eu q}}}{{\eu q}^{{\deg \pK}}}\Big)^m=(x)$ es principal. Sin embargo $(x)_S=
\pK^{m\deg {{\eu q}}}$ no es de grado $0$ y en particular no es
principal. Ahora bien, $(x)_S$ es el ideal fraccionario generado
por $x$ considerado como elemento de ${\mc O}_S$.

Sea $\deg\colon D_{\K}\to {\ma Z}$ el mapeo grado. 
Se tiene que $\deg(D_{\K})=\mu {\ma Z}$
para alg\'un $\mu\in {\ma N}$. Puesto que $\F$ es finito, $\mu=1$ 
por el teorema de F. K. Schmidt, Teorema \ref{T10.1.2.1}.
En particular $\deg(D_{\K})={\ma Z}$. Sea $\deg (D(S))=d{\ma Z}$, $d\in{\ma N}$.
Sea $E_{\K}(S)$ el grupo de las {\em $S$--unidades\index{S@$S$--unidades}}, es decir,
$E_{\K}(S)=\{x\in \*\K\mid v_{\pK}(x)=0\text{\ para toda\ } \pK\notin S\}$.

En resumen, tenemos las siguientes definiciones, donde
$\K$ es el campo fijado de antemano:

\las
\item[$\bullet$] $D_S=\langle \pK\mid\pK\notin S\rangle$.

\item[$\bullet$] Para $x\in\*\K$, $(x)_S:=\prod_{\pK\notin{\ma P}_{\pK}}
\pK^{v_{\pK}(x)}$.

\item[$\bullet$] $P_S=\{(x)_S\mid x\in\* \K\}$.

\item[$\bullet$] $C_S=D_S/P_S$.

\item[$\bullet$] $D_{\K}(S)=D(S)=\langle \pK\in S\rangle$.

\item[$\bullet$] $P_{\K}(S)=P(S)=P_{\K}\cap D(S)$.

\item[$\bullet$] $E_{\K}(S)=E(S)=\{x\in \*\K\mid v_{\pK}(x)=0
\text{\ para toda\ } \pK\notin S\}$.

\end{list}

\begin{teorema}\label{TRamDed2.1.3-1}
Las unidades de ${\mc O}_S$ es $E(S)$, esto es,
$E(S)=\*{{\mc O}_S}$ y $Cl_S\cong C_S$.
\end{teorema}

\begin{proof}
Si $x\in{\mc O}_S$, $v_{\pK}(x)\geq 0$ para toda $\pK\notin S$.
Si $x$ es unidad de ${\mc O}_S$, existe 
$y\in {\mc O}_S$ tal que $xy=1$ por lo que $v_{\pK}(x)+v_{\pK}(
y)=0$ para toda $\pK\notin S$. Por tanto $v_{\pK}(x)=0$ para
toda $\pK\notin S$. De esta forma tenemos que $\*{{\mc O}_S}
\subseteq E(S)$.

El rec\'iproco es claro ya que si $x\in E(S)$, entonces $v_{\pK}(x)
=0$ para toda $\pK\notin S$ lo cual implica que $v_{\pK}(x^{-1})=
-v_{\pK}(x)=0$ para toda $\pK\notin S$.

Ahora bien, $Cl_S\cong I_S/Pr_S$ donde $I_S$ es el conjunto
de los ideales fraccionarios de ${\mc O}_S$ y $Pr_S$ son los
ideales fraccionarios $(x)\in I_S$ donde $x\in\coc {\mc O}_S$.

Tenemos que $I_S$ es el grupo libre generado por los lugares
$\pK\in{\ma P}_{\K}\setminus S$, esto es, con la identificaci\'on
${\ma P}_{\K}\setminus S\ni\pK\longleftrightarrow \pK\cap
{\mc O}_S$, se tiene $I_S=D_S$. El mapeo $I_S\stackrel{\varphi}
{\lra} D_S/P_S$ es un epimorfismo de grupos. Adem\'as
$\ker\varphi =I_S\cap P_S$ es el conjunto de ideales fraccionarios
$M$ de ${\mc O}_S$ tales que existe $x\in\*\K$, $M=(x)_S$.

De esta forma, $x\in \*\K=\coc{\mc O}_S$ y el divisor principal
con respecto a $I_S$ es precisamente $M$ pues con respecto
a ${\mc O}_S$, el ideal de un elemento $y\in\*\K$ es
$(y)=\prod_{\P\notin S}\P^{\alpha_{\P}}$ con $\alpha_{\P}=
v_{\P}(y)$. Se sigue que $Cl_S\cong C_S$.
$\fin$
\end{proof}

La relaci\'on entre $h_S$ y $h_{\K}$ se obtiene del 
siguiente teorema de F. K. Schmidt.

\begin{teorema}[F. K. Schmidt]\label{TRamDed2.1.3}
Se tiene que $d=\mcd\{\deg \pK\mid \pK\in S\}$ 
donde $\deg (D(S))=d{\ma Z}$
y las siguientes sucesiones de grupos
son exactas:
\las
\item $1\longrightarrow \*\F\longrightarrow E_{\K}
(S)\longrightarrow P_{\K}(S)\longrightarrow 1$,
\item $1\lra D_{\K}(S)_0/P_{\K}(S)\lra I_{\K,0}
\lra C_S\lra {\ma Z}/d{\ma Z}\lra 0$.
\end{list}
\end{teorema}

\begin{proof} 
Para ${\eu A}\in D_{\K}(S)$, ${\eu A}=
\pK_1^{\alpha_1}\cdots \pK_s^{\alpha_s}$ con $s=|S|$,
$\deg {\eu A}=\sum_{i=1}^s\alpha_i\deg \pK_i$. Por tanto
\[
\deg (D_{\K}(S))=\langle\deg\pK\mid \pK\in 
S\rangle=\langle \mcd\{\deg\pK\mid\pK\in S\}\rangle=(d),
\]
por lo que $d=\mcd\{\deg\pK\mid\pK\in S\}$.

Sea $\varphi\colon E_{\K}(S)\lra P_{\K}(S)$ dado por 
$\varphi(x)=(x)_S\in P_{\K}(S)$. Entonces
$\varphi$ es un mapeo suprayectivo pues si 
${\eu A}\in P_K(S)$, entonces ${\eu A}\in P_K(S)
\cap D(S)$. Por tanto ${\eu A}=(y)_K=\prod_{\eu p}
{\eu p}^{v_{\eu p}(y)}$ y $v_{\eu p}(y)=0$ para toda
$y\in S$. Por tanto $(y)_K=(y)_S$ y $\varphi(y)=
(y)_S={\eu A}$, por lo tanto $\varphi$
es suprayectiva. Ahora
\begin{gather*}
\begin{align*}
\ker \varphi&=\{x\in E_{\K}(S)\mid v_{\pK}(x)=0 \text{\ para toda\ }\pK\in S\}\\
&=\{x\in\*\K\mid v_{\pK}(
x)=0 \text{\ para todo lugar\ } \pK\in \P \K\}=\*\F,
\end{align*}
\intertext{por lo que}
1\lra\*\F\lra E_{\K}(S)\lra P_{\K}(S)\lra 1
\end{gather*}
es exacta.

Sea $\tau\colon D_{\K}\lra D_S$, dada por
$\tau({\eu A})=\prod_{\pK\notin S}\pK^{v_{\pK}({\eu A})}$.
Por definici\'on de $D_S$, $\tau$ es suprayectiva y $\ker \tau=\{{\eu A}\in D_{\K}\mid 
v_{\pK}=0 \text{\ para toda\ } \pK\notin S\}=D_{\K}(S)$. Por lo tanto la sucesi\'on
$1\lra D_{\K}(S)\lra D_{\K}\stackrel{\tau}{\lra} D_S\lra 1$ es exacta. Ahora bien
$\tau(P_{\K})=\{\tau((x)_{\K})\mid x\in\* \K\}
=\big\{\prod_{\pK\notin S} \pK^{v_{\pK}(x)}\mid
x\in\* \K\big\}=\{(x)_S\mid x\in \* \K\}= P_S$.

En particular, $\tau$ induce un epimorfismo
$\tilde{\tau}\colon I_{\K}=D_{\K}/P_{\K}\lra C_S=D_S/P_S$,
${\eu A}\bmod P_{\K}\longmapsto \tau({\eu A})\bmod P_S$ y
\begin{gather*}
\begin{align*}
\ker\tilde{\tau}&=\{{\eu A}\bmod P_{\K}\mid\tau({\eu A})\in P_S\}\\
&=\{{\eu A}\bmod P_{\K}\mid
\tau({\eu A})=(x)_S \text{\ para alg\'un\ } x\in\* \K\}\\
&=\{{\eu A}\bmod P_{\K}\mid {\eu A}\in D_{\K}(S)\}=
\frac{D_{\K}(S) P_{\K}}{P_{\K}}\\
&\cong
\frac{D_{\K}(S)}{P_{\K}\cap D_{\K}(S)}=\frac{D_{\K}(S)}{P_{\K}(S)},
\end{align*}
\intertext{por lo que}
1\lra \frac{D_{\K}(S)}{P_{\K}(S)}\lra I_{\K}\stackrel{\tilde{\tau}}{\lra}C_S\lra 1
\intertext{es exacta. Se sigue que}
1\lra \frac{D_{\K}(S)_0}{P_{\K}(S)}\lra I_{\K,0}\stackrel{\tilde{\tau}}{\lra} C_S
\end{gather*}
es exacta.

Falta ver que $C_S/\tilde{\tau}(I_{\K,0})\cong d{\ma Z}$. Para ver esto \'ultimo, 
consideremos $\tau\colon D_{\K}\lra D_S$ nuevamente y $\tau(P_{\K})=P_S$,
$\tau(D_{\K}(S))\subseteq P_S$. Por tanto obtenemos el mapeo
$D_{\K}/(P_{\K} D_{\K}(S))\stackrel{\bar{
\tau}}{\lra} C_S=D_S/P_S$. Como $\tau$ es un epimorfismo, $\bar{\tau}$ tambi\'en
lo es y $\tau^{-1}(P_S)=P_{\K} D_{\K}(S)$, por lo que $D_{\K}/(P_{\K} D_{\K}(S))\cong C_S
\cong D_S/P_S$.

Estamos interesados en el con\'ucleo de 
\[
I_{\K,0}=\frac{D_{\K,0}}{P_{\K}}\stackrel{
\theta}{\lra} C_S\cong \frac{D_{\K}}{P_{\K} D_{\K}(S)}.
\]
Aqu{\'\i}, para ${\eu A}\in D_{\K,0}$,
$\theta({\eu A}\bmod P_{\K})={\eu A}\bmod P_{\K} D_{\K}(S)$, 
por lo que $\im \theta=
\frac{D_{\K,0} P_{\K} D_{\K}(S)}{P_{\K} D_{\K}(S)}$, 
de donde 
\begin{align*}
\coker \theta&=\frac{D_{\K}/(P_{\K} D_{\K}(
S))}{(D_{\K,0}P_{\K}D_{\K}(S))/(P_{\K}D_{\K}(S))}\cong 
D_{\K}/(D_{\K,0}P_{\K}D_{\K}(S))\\
&\underbracket[0pt]{=}_{\substack{\uparrow\\ 
P_{\K}\subseteq D_{\K,0}}} D_{\K}/(D_{\K,0}D_{\K}(S)).
\end{align*}

Ahora $\deg\big(D_{\K}/(D_{\K,0}D_{\K}(S))\big)=
\frac{\deg (D_{\K})}{\deg (D_{\K}(S))}{\ma Z}\cong
{\ma Z}/d{\ma Z}$. Por tanto la sucesi\'on
\[
1\lra \frac{D_{\K}(S)_0}{P_{\K}(S)}\lra I_{\K,0}
\lra C_S\lra \frac{\ma Z}{d{\ma Z}}\lra 0
\]
es exacta.
$\fin$
\end{proof}

\begin{observacion}\label{ORamDed1.2.4}
No es necesario suponer que $\K$ sea un campo 
de funciones congruente en el Teorema
\ref{TRamDed2.1.3}. El \'unico cambio a considerar 
en general para un campo de funciones $\K/k$ arbitrario es
que $\deg (D_{K})=\mu {\ma Z}$ para $\mu\in{\ma N}$. Entonces
$\mu\mid d$ y se obtiene la sucesi\'on exacta
\[
1\lra \frac{D_{\K}(S)_0}{P_{\K}(S)}\lra I_{\K,0}\lra C_S\lra 
\frac{\ma Z}{(d/\mu){\ma Z}}\cong
\frac{\mu{\ma Z}}{d{\ma Z}}\lra 0.
\]
\end{observacion}

\begin{corolario}\label{CRamDed1.2.5} Se tiene que 
$E_{\K}(S)/\*\F$ es un grupo abeliano libre
finitamente generado de rango $s-1=|S|-1$.
\end{corolario}

\begin{proof} 
Se tiene $E_{\K}(S)/\*\F\cong P_{\K}(S)\subseteq 
D_{\K}(S)_0$ y este \'ultimo
es un grupo libre en $|S|-1$ generadores. Puesto que $D_{\K}(S)_0/P_{\K}(S)$ es finito
se sigue que $P_{\K}(S)_0$ es libre en $|S|-1$ generadores.
$\fin$
\end{proof}

\begin{corolario}[Teorema de las Unidades de
Dirichlet\index{teorema de las unidades de
Dirichlet}\index{Dirichlet!teorema de las unidades de
$\sim$}]\label{CRamDed1.2.5(0)}
Se tiene $E_K(S)\cong \*\F\times {\ma Z}^{s-1}$. $\fin$
\end{corolario}

\begin{observacion}\label{ORamDed1.2.5(1)} Con respecto al Corolario
\ref{CRamDed1.2.5}, de hecho, si fijamos $\pK_0\in S$, entonces para
$\pK\in S$, $\pK\neq \pK_0$, $\deg \pK_0^{\deg\pK}=\deg\pK^{\deg \pK_0}$. Puesto
que $I_{\K,0}$ es finito, existe $m_{\pK}\in{\ma N}$ tal que
$\Big(\frac{\pK^{\deg\pK_0}}{\pK_0^{\deg\pK}}\Big)^{m_{\pK}}=(x_{\pK})$
es principal, $x_{\pK}\in \*\K$. Entonces $x_{\pK}\in E_{\K}(S)$ y adem\'as
veremos que 
$\langle x_{\pK}\mid\pK\in S\setminus\{\pK_0\}\rangle$ es un grupo libre
de rango $|S|-1$ que, en adici\'on, es de \'indice finito en $E_{\K}(S)$.

Primero veamos que $\langle x_{\pK}\mid\pK\in S\setminus\{\pK_0\}\rangle$
es ${\ma Z}$--libre. Sea $S'=S\setminus\{\pK_0\}$. Si se tiene una relaci\'on
$a=\prod_{\pK\in S'} x_{\pK}^{c_{\pK}}=1$ con $c_{\pK}\in{\ma Z}$. Entonces
$v_{\pK}(a)=c_{\pK}(\deg \pK_0) m_{\pK}=0$ para toda $\pK\in S'$. Por tanto
$c_{\pK}=0$ para toda $\pK\in S'$.

Ahora veamos que $\langle x_{\pK}\mid\pK\in S'\rangle$
es de \'indice finito en $E_{\K}(S)$. Sea $x\in E_{\K}(S)$. Entonces
$v_{\eu q}(x)=0$ para toda ${\eu q}\notin S$. Para $\pK\in S$ sea
$d_{\pK}=v_{\pK}(x)\in {\ma Z}$. Sea $y=x^f$ donde $f=\prod_{\pK\in S'}
m_{\pK} (\deg \pK_0)$. Entonces para $\pK\in S'$:
\begin{align*}
v_{\pK}\big(x^f x_{\pK}^{(-f d_{\pK}/m_{\pK}(\deg \pK_0))}\big)&=fv_{\pK}(x)-
\frac{fd_{\pK}}{m_{\pK}(\deg \pK_0)}v_{\pK}(x_{\pK})\\
&=fd_{\pK}-
\frac{fd_{\pK}}{m_{\pK}(\deg\pK_0)}m_{\pK}(\deg\pK_0)=0.
\end{align*}

Se sigue que $v_{\eu q}(x)=0$ para todo lugar ${\eu q}\neq \pK_0$.
Como $(x)$ es de grado $0$ se tiene que
$v_{\pK_0}(x)=0$ por lo que $x^f
z^{-1}\in \*k=\*\F$, donde $z=\prod_{\pK\in S'}
x_{\pK}^{(f d_{\pK}/m_{\pK}(\deg \pK_0))}$. Por lo tanto $x^f\in
\langle x_{\pK}\mid \pK\in S'\rangle$ y $E_{\K}(S)/\langle x_{\pK}\mid
\pK\in S'\rangle$ es un grupo finito.
\end{observacion} 

\begin{corolario}\label{CRamDed1.2.6} Se tiene $|C_S|=h_S<\infty$. 
$\fin$
\end{corolario}

\begin{corolario}\label{CRamDed1.2.7} Se tiene que
\[
h_{\K} d = \gamma h_S,
\]
donde $\gamma=\big|\frac{D_{\K}(S)_0}{P_{\K}(S)}\big|$. 
$\fin$
\end{corolario}

\subsection{Sobre el grupo $I_A/Q_A$}\label{S1.3.5}

Muchas veces tenemos que los ideales en un dominio Dedekind act\'uan en
algunos grupos sin que se pueda extender f\'acilmente la acci\'on al grupo
de los ideales fraccionarios. Sin embargo, muchas veces esta acci\'on s\'i
puede hacerse sobre $C_S$, el grupo de clases de ideales fraccionarios
del dominio Dedekind.

Sea $A$ un dominio Dedekind con $h_A$ finito.

\begin{definicion}\label{D1.3.5.1}
Se definen
\begin{gather*}
I_A=\{\text{ideales no cero de } A\},\\
Q_A=\{\text{ideales principales no cero de } A\},\\
D_A=\{\text{grupo de ideales fraccionarios de } A\},\\
P_A=\{\text{grupo de ideales fraccionarios principales de } A\},\\
\pic (A)=C_A=\frac{D_A}{P_A}.
\end{gather*}
\end{definicion}

Se tiene que $D_A$ y $P_A$ son grupos pero $I_A$ y $Q_A$ no son
grupos sino simplemente monoides con identidad $A=(1)$.

Definimos una relaci\'on en $I_A$: ${\eu A}\sim {\eu B}$
si y solamente si existen $x,y\in A\setminus \{0\}$ con $x{\eu A}=
y{\eu B}$. 

Esta es una relaci\'on de equivalencia. Denotemos la clase
de ${\eu A}\in I_A$ por $[{\eu A}]$ y definimos el conjunto cociente
\[
I_A/\sim:=\{[{\eu A}]\mid {\eu A}\in I_A\}.
\]

En $I_A/\sim$ se define la operaci\'on natural
\[
[{\eu A}][{\eu B}]=[{\eu A}{\eu B}].
\]

Se verifica f\'acilmente que la operaci\'on no depende de los representantes.

Esta operaci\'on define una estructura de grupo en $I_A/\sim$ y
la clase de $A$ es la identidad de este grupo. Tambi\'en es f\'acil ver
que $[A]=\{(x)\mid x\in A, x\neq 0\}$, esto es, la clase de $A$ consiste
de los ideales principales $Q_A$. Por tanto podemos denotar
\[
I_A/\sim = I_A/Q_A.
\]
Como $A$ es un anillo conmutativo, $I_A/Q_A$ es un grupo abeliano.

Sea ${\eu A}\in I_A$. En particular ${\eu A}\in D_A$ y $\bar{\eu A}\in
D_A/P_A$ es de orden finito, digamos $\bar{\eu A}^m=1$, esto es,
${\eu A}^m=(x)$ es principal. Sin embargo $x\in {\eu A}^m\subseteq A$,
por lo que $[{\eu A}][{\eu A}^{m-1}]=[{\eu A}^m]=[(x)]=[A]$. Por
tanto $[{\eu A}]^{-1}=[{\eu A}^{m-1}]$.

El resultado siguiente se verifica de manera rutinaria.

\begin{teorema}\label{T1.3.5.2} El mapeo natural $\varphi\colon
I_A/Q_A\lra D_A/P_A$,
$[{\eu A}]\lra \bar{\eu A}$ define un isomorfismo de grupos
abalianos
\begin{gather*}
I_A/Q_A\cong D_A/P_A=\pic (A). \tag*{$\fin$}
\end{gather*}
\end{teorema}

\section{Aritm\'etica de extensiones de campos de funciones
globales}\label{SRamDed.2}

Varios de los resultados que presentamos en esta secci\'on son
aplicables a campos de funciones m\'as 
generales pero nos concentraremos
en extensiones finitas $L/\K$ de campos de 
funciones congruentes. Seguimos con la notaci\'on
$K=\F(T)$.

Para un campo $E$ de caracter\'istica $p>0$, usamos la
notaci\'on $\wp(x)=x^p-x$ para $x\in E$ y $\wp(E)=\{
\wp(x)\mid x\in E\}$.

\subsection{Extensiones radicales}\label{S5.2.B-1}

\begin{teorema}\label{TRam1}
Sea $\K =K(\sqrt[n]{\gamma D})$, donde 
$D\in R_T$ es un polinomio m\'onico y $\gamma \in
\ff$. Sea $D=P_1^{\alpha_1}\cdots P_s^{\alpha_s}$ 
la descomposici\'on de $D$ como producto de
polinomios irreducibles. Sea $\P_i$ el divisor
primo correspondiente a $P_i\in R_T^+$.
Supondremos que $D$
est\'a libre de $n$--potencias, es decir, 
$0<\alpha_i<n$ para $1\leq i\leq s$ 
y tambi\'en supondremos que $p\nmid n$.

Si $P\in R_T^+$ y $P\nmid D$, entonces $P$
es no ramificado,
los primos finitos ramificados son
$\P_1,\ldots,\P_s$ y ellos son moderadamente
ramificados. 

Sea
${\eu p}_i$ un primo en $\K$ encima de $\P_i$. Entonces
\[
e({\eu p}_i|\P_i)=\frac{n}{d_i},
\]
donde $d_i=\mcd(\alpha_i,n)$.

Sea ahora $d=\mcd(\deg D,n)=\mcd(n-\deg D,n)$.
Entonces el \'indice
de ramificaci\'on del primo infinito en $\K/K$ es
\[
e_{\infty}=e_{\infty}(\K |K) = \frac{n}{d}.
\]
Esto es, $\p$ es ramificado en $\K/K \iff
n\nmid \deg D$.

Finalmente los grados de inercia de los divisores primos
de $\K$ sobre $\p$ son los grados de los polinomios
$\Irr(\zeta_d^j\mu, \bar{X},\F)$, $0\leq j\leq d-1$,
donde $\mu$ es una
ra\'iz $d$--\'esima fija de $\gamma$ en alguna cerradura
algebraica de $\F$ y $\zeta_d$ denota una ra\'iz 
$d$--\'esima primitiva de la unidad.
\end{teorema}

\begin{proof} Cuando $P\nmid D$, $P\in R_T^+$, entonces
$P$ es no ramificado (ver \cite[Example 5.8.9]{Vil2006}). 
Presentamos otra demostraci\'on de este hecho usando
la teor\'ia de campos ciclot\'omicos. Primero supongamos
que $\K/K$ es una extensi\'on de Kummer. Sea
$\K_1:=K(\sqrt[n]{(-1)^{\deg D} D})\subseteq \cicl M{}$.
Entonces $P$ es no ramificado en $\K_1/K$ (Proposici\'on
\ref{P6.3.4}). Puesto que $\K\K_1=\K_1\F(\sqrt[n]
{(-1)^{\deg D}\gamma})$, $P$ es no ramificado en $\K\K_1/K$
y en particular $P$ no es ramificado en $\K/K$.

Otra demostraci\'on de lo anterior se puede hacer de
manera general, ver la Proposici\'on \ref{P6.4-1.Ram7}.

Cuando $\K/K$ no es de Kummer, sea $K_0$ la extensi\'on de
constantes de $K$ que contiene las $n$--\'esimas ra\'ices
de la unidad, $K_0=K(\zeta_n)$. 
Entonces $\K K_0=K_0(\sqrt[n]{\gamma D})$
es una extensi\'on de Kummer. Por lo anterior
$P$ no es ramificado en $\K K_0/K_0$. Puesto que
$P$ tampoco es ramificado en $K_0/K$ se sigue que $P$
no es ramificado en $\K K_0/K$ y en particular $P$
es no ramificado en $\K/K$.

Para $P_i$, $1\leq i\leq s$, tenemos
\begin{gather}\label{Eq5.1.B}
e({\eu p}_i|\P_i) v_{\P_i}(D)=e({\eu p}_i|\P_i) \alpha_i=
v_{{\eu p}_i}(D)=v_{{\eu p}_i}((\sqrt[n]{\gamma D})^n)=nv_{{\eu p}_i}
(\sqrt[n]{\gamma D}).
\end{gather}
Obtenemos de (\ref{Eq5.1.B}) que
$\frac{n}{d_i}\mid e({\eu p}_i|\P_i)$. 
Por otro lado, si escribimos $\K =K(y)$,
donde $y^n=\gamma D$, establezcamos
$z=y^{n/d_i}$. Por tanto $z^{d_i}=y^n =
\gamma D= \gamma (P_i^{\alpha_i/d_i})^{d_i}
(D/P_i^{\alpha_i})$. Por lo tanto
$K(z)=K(\sqrt[d_i]{\gamma D/P_i^{\alpha_i}})$. 
Por la primera parte, tenemos que $\P_i$ es no ramificado en
$K(z)/K$ .
Se sigue que $e({\eu p}_i|\P_i)=n/d_i$.
Por tanto $e_{\P_i}=n/d_i$.

Ahora veamos $\p$. Se tiene
\begin{align*}
D(T)&=T^s+a_{s-1}T^{s-1}+\cdots+a_1T+a_0\\
&=T^s\big(1+a_{s-1}(1/T)+\cdots+a_1(1/T)^{s-1}+a_0(1/T)^s\big).
\end{align*}

Sea $X:=1/T$ y sea $D_1(X)=1+a_{s-1}X+\cdots+a_1X^{s-1}
+a_0X^s$. Notemos que $\mcd(X, D_1(X))=1$ y que
$D_1(X)=a_j D_2(X)$ donde $a_j\neq 0$ y donde $a_i=0$
para $0\leq i\leq j-1$, es decir, $a_j$ es el coeficiente l\'ider
de $D_1(X)$ y $D_2(X)$ es m\'onico.

Sea $a_j=\varepsilon$. Por tanto, escribiendo $s=an+b$ con
$a\in{\ma N}\cup \{0\}$ y $0\leq b< n$, se tiene
\begin{gather*}
\gamma D(T)=\frac{\gamma\varepsilon D_2(X)}{X^s}=
\frac{\gamma\varepsilon D_2(X)\cdot X^{n-b}}{X^{(a+1)n}}.
\intertext{Por tanto}
\sqrt[n]{\gamma D(T)}=\frac 1{X^{a+1}}\sqrt[n]{\gamma\varepsilon D_2(X)X^{n-b}}
\quad\text{y}\quad \mcd(D_2(X),X)=1,\quad \gamma\varepsilon\in \*\F.
\intertext{Por el caso finito, se tiene, usando que $\mcd(n-s,s)=
\mcd(n,s)=\mcd(n,b)=d$ y que $s=\deg D$,}
e_{\infty}=\frac n{\mcd(n-b,n)}=\frac n{\mcd(n-\deg D,n)}=\frac n
{\mcd(\deg D,n)}=\frac nd.
\end{gather*}

Ahora veamos como es la descomposici\'on de $\p$
en $\K/K$. Consideremos la subextensi\'on $E=
K(\sqrt[d]{\gamma D})$. Entonces aplicando lo que acabamos
de probar sobre la ramificaci\'on $\p$ ahora en la 
extensi\'on $E/K$, se tiene que $\p$ no es ramificado.
Entonces si $\pL$ es un primo en $\K$ sobre $\p$ y
${\mc P}=\pL\cap E$, se tiene que al ser $\p$ totalmente
ramificado en $\K/E$, entonces $f_{\K/K}(\pL|\p)=
f_{E/K}({\mc P}|\p)=f$. 

Haremos el c\'alculo de $f$ en las completaciones. Para 
esto consideremos 
\begin{gather*}
D(T)=T^s+a_{s-1}T^{s-1}+\cdots +a_1 T+a_0,
\intertext{donde $s=\deg D$, $d\mid s$ y $a_0,a_1,\ldots,a_{s-1}\in \F$.
Entonces}
D(T)=T^s\big(1+a_{s-1}(1/T)+\cdots+
a_1(1/T)^{s-1}+a_0(1/T)^s\big)=T^s D_1(1/T).
\end{gather*}
Puesto
que $d\mid s$, se tiene que $E=K(\sqrt[d]{\gamma D_1(1/T)})$.

Las completaci\'on de $E$ en ${\mc P}$ corresponde a
$K_{\infty}(\delta)$, donde $\delta$ es una ra\'iz de 
$\sqrt[d]{\gamma D_1(1/T)}$. Ahora bien si tomamos
$\gamma D_1(1/T)\bmod \p\equiv \gamma\bmod (1/T)$. Entonces
el polinomio $X^d-\gamma D(1/T)\in K[X]$ m\'odulo $\p$
es $\bar{X}^d-\gamma \in \F[\bar{X}]$. Si $\mu$ es una ra\'iz
$d$--\'esima fija de $\gamma$ en una cerradura algebraica
de $\F$, se tiene que
\[
\bar{X}^d-\gamma=\prod_{i=0}^{d-1} (\bar{X}-\zeta_d^i \mu)=
\prod_{j=1}^r\Irr(\zeta_d^{i_j}\mu,\bar{X},\F)
\]
para algunas $0\leq i_1<i_2<\cdots<i_r\leq d-1$. 
Usando el Lema de Hensel, se tiene que existen $r$ polinomios
m\'onicos irreducibles distintos $F_1(X),\ldots, F_r(X)\in
K_{\infty}[X]$ tales que $\deg_X F_j(X)=\deg_{\bar{X}}
\Irr(\zeta_d^{i_j}\mu,\bar{X},\F)$ y 
\[
X^d-\gamma D_1(1/T)=\prod_{j=1}^r F_j(X).
\]

Se sigue que las completaciones de $E$ en los diversos
divisores de $\p$ son de grado $\deg_X F_j(X)=
\deg_{\bar{X}}\Irr(\zeta_d^{i_j}\mu,\bar{X},\F)$, $1\leq j
\leq r$, de donde se sigue lo afirmado. $\fin$
\end{proof} 

\begin{observacion}\label{O10.1.Ram2}
Con las notaciones del Teorema \ref{TRam1}, cuando
$n\mid q-1$, entonces la extensi\'on radical $\K=K(
\sqrt[n]{\gamma D})$ es una extensi\'on c\'iclica de Kummer.
Cuando $K$ no contiene a las ra\'ices $n$--\'esimas de la
unidad, podemos hacer una extensi\'on de constantes
adecuada y en la nueva extensi\'on obtenemos una
extensi\'on c\'iclica de Kummer.

Ahora bien, si $\K/K$ es una extensi\'on c\'iclica de grado
$n$ pero $n\nmid q-1$ (y donde adem\'as suponemos que
$p\nmid n$), entonces al agregar constantes adecuadas,
la nueva extensi\'on obtenida es una extensi\'on de Kummer.
En el siguiente resultado, presentamos un m\'etodo para
hallar el elemento cuya $n$--\'esima ra\'iz es generadora
de la extensi\'on.
\end{observacion}

\begin{proposicion}\label{P10.1.Ram3}
Sea $n\in{\ma N}$ tal que $n\nmid q-1$ y tal que $p\nmid n$.
Sea $\K/K$ una extensi\'on c\'iclica de grado $n$ con grupo
de Galois $G=\Gal(\K/K)=\langle\sigma\rangle$. Sea $m=
o(q\bmod n)$, es decir, $m$ es el m\'inimo 
n\'umero natural tal que
$n\mid q^m-1$. Sea $\zeta$ una ra\'iz $n$--\'esima primitiva
de la unidad. Se tiene ${\ma F}_q(\zeta)={\ma F}_{q^m}$.
Sea $\K(\zeta)=\K_m$ y $K(\zeta)=K_m$ 
las extensiones de constantes de grado $m$. Sea $\K=
K(\delta)$, con $\delta \in \K$ y tal que 
\[
\nu=\sum_{i=0}^{n-1}\zeta^i \sigma^i(\delta)\neq 0.
\]

Entonces $D=\nu^n\in \K_m$ y $\K_m=K_m(\sqrt[n]
{D})=K_m(\nu)$.
\end{proposicion}

\begin{proof}
Se tiene que ${\mc G}=\Gal(\K_m/K_m)=\langle\tau\rangle 
=\{1,\tau,\ldots, \tau^{n-1}\}$ y $\sigma =\tau|_{\K}$.
Por el teorema de la independencia de automorfismos
de Artin, $\{1,\tau,\ldots, \tau^{n-1}\}$ es 
linealmente independiente
sobre $\K(\zeta)$. Por tanto, existe $\mu\in \K(\zeta)$
 tal que $\sum_{i=0}^{n-1}\zeta^i \tau^i(\mu)\neq 0$.
 
Veamos que podemos tomar $\mu\in \K$.
Supongamos que $\sum_{i=0}^{n-1}\zeta^i
\tau^i(\xi)=0$ para toda $\xi\in \K$. Se tiene que
$\K(\zeta)/K$ tiene como base $\{1,\zeta,\ldots,
\zeta^{m-1}\}$, por lo que existen $\xi_1,\ldots,
\xi_{m-1}\in \K$ tales que $\mu=\sum_{j=0}^{m-1}
\zeta^j \xi_j$. Entonces
\begin{align*}
\sum_{i=0}^{n-1}\zeta^i\tau^i(\mu)&=
\sum_{i=0}^{n-1}\zeta^i\tau^i\Big(\sum_{j=0}^{m-1}
\zeta^j \xi_j\Big)=\sum_{i=0}^{n-1}\zeta^i\sum_{j=0}^{m-1}
\tau^i\big(\zeta^j\xi_j)\\
&= \sum_{i=0}^{n-1}\zeta^i \sum_{j=0}^{m-1}
\zeta^j\tau^i(\xi_j)=\sum_{i=0}^{n-1}\sum_{j=0}^{m-1}
\zeta^i\zeta^j\sigma^i(\xi_j)\\
&=\sum_{j=0}^{m-1}\sum_{i=0}^{n-1}\zeta^j\zeta^i
\sigma^i(\xi_j)=\sum_{j=0}^{m-1}\zeta^j\Big(
\sum_{i=0}^{n-1}\zeta^i\sigma^i(\xi_j)\Big)=
\sum_{j=0}^{m-1}\zeta^j (0)=0
\end{align*}
lo cual es absurdo. Por tanto existe $\delta\in\K$ tal que
$\nu=\sum_{i=0}^{n-1}\zeta^i \sigma^i(\delta)\neq 0$.
Se tiene que $\nu\in \K_m$. 

Por tanto
\begin{align*}
\tau(\nu)&=\tau\Big(\sum_{i=0}^{n-1}
\zeta^i \sigma^i(\delta)\Big)=\sum_{i=0}^{n-1}
\zeta^i\sigma^{i+1}(\delta)=\zeta^{-1}\sum_{i=0}^{n-1}
\zeta^{i+1}\sigma^{i+1}(\delta)\\
&=\zeta^{-1}\sum_{i=0}^{n-1}\zeta^i\sigma^i(\delta)=
\zeta^{-1}\nu\neq \nu.
\end{align*}
En particular $\nu\notin \K$. El conjunto de 
conjugados de $\nu$ es $\big\{\tau^j(\nu)
\big\}_{j=0}^{n-1}=\big\{\zeta^{-j}\nu\big\}_{j=0}^{n-1}$
el cual es de cardinalidad $n$ por lo que $\K_m=
K_m(\nu)$ y adem\'as $\tau(\nu^n)=
\tau(\nu)^n=(\zeta^{-1}\nu)^n=\nu^n$,
lo cual implica que $\nu^n=D\in K_m$. 
Por tanto $\K_m=K_m(\nu)=K_m(\sqrt[n]{D})$.
$\fin$
\end{proof}

\begin{proposicion}\label{P3.7.V.Ram} Sea $l$ un primo distinto
de $p$ y tal que $l^n\mid q-1$.
Sea $\K=K(\sqrt[l^n]{\gamma D})$
donde $\gamma\in {\ma F}_q^{\ast}$ y $D\in R_T$ es un
polinomio m\'onico libre de $l^n$--potencias. 

Entonces, si $e_{\infty}$, $f_{\infty}$ y $h_{\infty}$ denotan
el \'indice de ramificaci\'on, el grado de inercia y el
\'indice de descomposici\'on de $\p$ respectivamente
en $\K/K$, entonces
\begin{gather*}
e_{\infty}=l^{n-t},\quad f_{\infty}=l^m,\quad\text{y}\quad
h_{\infty}=l^{t-m},
\end{gather*}
donde $\deg D=l^{t^{\prime}} a$ con $\mcd (a,l)=1$, 
$t=\min\{n,t^{\prime}\}$ y
${\ma F}_q(\sqrt[l^t]{(-1)^{\deg D}\gamma})=
{\ma F}_{q^{l^m}}$.
\end{proposicion}

\begin{proof}
El c\'alculo del \'indice de ramificaci\'on $e_{\infty}$
est\'a dado en el Teorema \ref{TRam1}.

Se tiene que $\p$ es no ramificado en $K(\sqrt[l^t]{
(-1)^{\deg D} D})$. Por tanto, del
Teorema \ref{T6.3.5}, $\p$ obtenemos que se descompone
totalmente en $K(\sqrt[l^t]{(-1)^{\deg D}D})
 \subseteq \cicl D{}$. Por otro
$\p$ es totalmente inerte en $K{\ma F}_{q^{l^m}}/K$
por ser $\p$ de grado uno (Teorema \ref{T6.1.4}).

Se sigue que el grado de inercia de
$\p$ en $K(\sqrt[l^t]{D})
{\ma F}_{q^{l^m}}/K$ es $l^m$. De esta forma obtenemos
que el campo de inercia de $\p$ en
$K(\sqrt[l^t]{D}){\ma F}_{q^{l^m}}/K$ es $K(\sqrt[l^t]{(-1)^{\deg D}D})$.
Por lo tanto $\p$ es totalmente descompuesto en 
la extensi\'on $K(\sqrt[l^t]{
\gamma D}){\ma F}_{q^{l^m}}/K(\sqrt[l^t]{\gamma D})$:
\[
\xymatrix{
K(\sqrt[l^t]{\gamma D})\ar@{-}[d]_{l^t}\ar@{-}[r]&K(\sqrt[l^t]{\gamma D})
{\ma F}_{q^{l^m}}=
K(\sqrt[l^t]{\gamma D})K(\sqrt[l^t]{(-1)^{\deg D}D})\ar@{-}[d]\\
K\ar@{-}[r]^{l^t}_{\substack{\p\text{\ totalmente}\\ \text{descompuesto}}}&
K(\sqrt[l^t]{(-1)^{\deg D}D})
}
\]
Por lo tanto $f_{\infty}=l^m$. Se sigue el resultado. $\fin$
\end{proof}

\begin{observacion}\label{O3.8.V.Ram}
La obtenci\'on de $f_{\infty}$ en la Proposici\'on \ref{P3.7.V.Ram}
pudo haber sido realizada tambi\'en usando el Lema de
Hensel como lo hicimos en el Teorema \ref{TRam1}.
\end{observacion}

\subsection{Extensiones c\'iclicas de grado primo}\label{SRam2.1}

En esta parte presentaremos el tipo de descomposici\'on de primos
de un campo $\K$ en una extensi\'on c\'iclica $L/\K$ de grado primo
$l$ en ambos casos: $l=p$ y $l\neq p$.

Enunciamos un caso particular del Proposici\'on \ref{P3.7.V.Ram}.

\begin{proposicion}\label{P5.1.2.Ram}
Sea $l$ un n\'umero primo tal que $l\mid q-1$. Sea $\K=K(\sqrt[l]{\gamma D})$
un extensi\'on de Kummer c\'iclica de grado $l$ donde $\gamma\in\*\F$ y
$D\in R_T$ es un polinomio m\'onico libre de $l$--potencias. Entonces
el comportamiento de $\p$ en $\K/K$ es el siguiente:
\l
\item Si $l\nmid \deg D$, $\p$ es ramificado.
\item Si $l\mid \deg D$ y $\gamma\in ({\ma F}_q^{\ast})^l$, $\p$ se descompone.
\item Si $l\mid \deg D$ y $\gamma\notin ({\ma F}_q^{\ast})^l$, $\p$ es inerte.
\end{list}
\end{proposicion}

\begin{proof}
Si $l\nmid \deg D$, entonces por el Teorema \ref{TRam1}, el
\'indice de ramificaci\'on de $\p$ en $\K/K$ es $l/1=l$ y por tanto
$\p$ es ramificado.

Ahora, si $l\mid \deg D$, en la notaci\'on de 
la Proposici\'on \ref{P3.7.V.Ram}, $t=1$.
Sea ${\ma F}:=\F(\sqrt[l]{(-1)^{\deg D} \gamma})=
\F(\sqrt[l]{\gamma})$. Entonces ${\ma F}=\F\iff
\sqrt[l]{\gamma}\in \*\F \iff \p$ se descompone. $\fin$
\end{proof}

Podemos usar la Proposici\'on \ref{P5.1.2.Ram} para estudiar
el comportamiento de cualquier primo de $K$ en una extensi\'on de
Kummer de grado $l$. Sea $\K$ como en la Proposici\'on
\ref{P5.1.2.Ram}. Sea $P\in R_T^+$ de grado $d$. 
Si $P\mid D$, entonces $P$ es ramificado en $\K/K$. 
Supongamos que $P\nmid D$. Consideremos
la extensi\'on de constantes $K_d$ de $K$. Entonces
$P$ se descompone totalmente en $K_d/K$ (Teorema \ref{T6.1.4}).

Sea $\pK$ un primo en $K_d$ sobre $P$. $\pK$
es un primo de grado $1$ en $K_d={\ma F}_{q^d}(T)$.
Se tiene $\K_d=K_d(\sqrt[l]{\gamma D})$. Hacemos de $\pK$
el primo infinito en $K_d$. Para esto, sea $X:=\frac{1}{T-\alpha}$
donde $T-\alpha\in {\ma F}_{q^d}[T]$ es el polinomio asociado
a $\pK$. En particular, puesto que $P\nmid D$, $D(\alpha)\neq 0$.
Escribimos
\begin{align*}
D(T)&=(T-\alpha)^s+a_{s-1}(T-\alpha)^{s-1}+\cdots+a_1(T-\alpha)+a_0\\
&= (T-\alpha)^s\Big(\frac{a_0}{(T-\alpha)^s}+\frac{a_1}{(T-\alpha)^{s-1}}+
\cdots +\frac{a_{s-1}}{(T-\alpha)}+1\Big)\\
&=X^{-s}(a_0 X^s+a_1 X^{s-1}+\cdots+a_{s-1}X+1),
\end{align*}
con $s=\deg D$. Sean $s=tl+r$ con $1\leq r\leq l$ y
$D_1(X):=a_0 X^s+a_1 X^{s-1}+\cdots + a_{s-1}X +1$. 

Entonces
\[
K_d(\sqrt[l]{\gamma D})=K_d(\sqrt[l]{\gamma X^{-l(t+1)}X^{l-r}D_1(X)})=
K_d(\sqrt[l]{\gamma X^{l-r}D_1(X)}).
\]
Puesto que $\pK$ corresponde al primo infinito (es decir a $1/X$),
el primo infinito no es ramificado por lo que $l\mid \deg{X^{l-r}
D_1(X)}=l-r+s=l-r+\deg D$.

Tenemos $\K_d=K_d(\sqrt[l]{\mu D_2(X)})$ 
con $\mu =\gamma a_0
=\gamma D(\alpha)$ y $D_2(X)=X^{l-r}D_1(X)$.
Entonces $\pK$ se descompone en $\K_d/K_d$ si y solamente
si $\mu\in (\*{{\ma F}_{q^d}})^l$. Se sigue que $P$ se
descompone en $\K/K$ si y solamente si $\gamma D(\alpha)\in
(\*{{\ma F}_{q^d}})^l$.

En resumen, tenemos 

\begin{proposicion}\label{P5.1.2.Ram2}
Sea $P$ un primo no ramificado en la extensi\'on de Kummer
c\'iclica $\K/K$ de grado $l$ 
dada por $\K=K(\sqrt[l]{\gamma D})$ donde 
$\gamma\in \*\F$ y $D\in R_T$ es un polinomio
de grado $d$. Sean $\alpha
\in\bar{\ma F}_q$ una
ra\'iz de $P$ y $\mu=\gamma D(\alpha)
\in \*{{\ma F}_{q^d}}$. Entonces
$P$ se descompone en $\K/K$ si y solamente
si $\mu \in (\*{{\ma F}_{q^d}})^l$. $\fin$
\end{proposicion}

\begin{observacion}\label{O5.1.2.Ram2-1}
Una demostraci\'on alternativa de la Proposici\'on \ref{P5.1.2.Ram2}
es usando el teorema de Kummer \cite[Theorem 5.8.2]{Vil2006}.
M\'as precisamente, sea $P\in R_T^+$ de grado $d$
tal que $P\nmid D$.
Consideremosel polinomio $G(X)=X^l-\gamma D\bmod P=
X^l-\gamma D(\alpha)\in {\ma F}_{q^d}[X]$, donde
$\alpha$ es una ra\'iz de $P(T)$.
Entonces $G(X)$ se descompone en $l$
factores si y solamente si $G(X)$ tiene una ra\'iz en ${\ma F}_{q^d}$,
esto es, existe $\beta\in {\ma F}_{q^d}$ tal que $\beta^l-
\gamma D(\alpha)=0$.
\end{observacion}

Resumiendo, el comportamiento de los primos finitos en una
extensi\'on de Kummer de grado primo $l$ en un campo
de funciones congruente es el siguiente:

\begin{teorema}\label{T5.1.2.Ram2-2}
Sean $l$ un primo dividiendo a $q-1$, $D\in R_T$, $\gamma
\in\*\F$, $\K=K(\sqrt[l]{\gamma D(T)})$, $P\in R_T^+$ 
de grado $d$ y $\alpha\in{\ma F}_{q^d}$
una ra\'iz de $P(T)$. Entonces
\lasa
\item Si $P|D$ se tiene que $P$ es ramificado en $\K/K$.

\item Si $P\nmid D$ y $\gamma D(\alpha)\in (\*{{\ma F}_{q^d}})^l$,
entonces $P$ se descompone en $\K/K$.

\item Si $P\nmid D$ y $\gamma D(\alpha\notin (\*{{\ma F}_{
q^d}})^l$, entonces $P$ es inerte en $\K/K$. $\fin$
\end{list}
\end{teorema}

Ahora consideremos $l=p$. Estudiaremos la ramificaci\'on
de los primos en extensiones de Artin--Schreier
de $K$.

Sea $\K/K$ una extensi\'on c\'iclica de grado
$p$. Entonces $\K=K(y)$ donde $y^p-y=\alpha\in K$.

Usando fracciones parciales y modificando $y$ con
transformaciones del tipo $y\leftrightarrow jy+c$
con $j\in\{1,2,\ldots, p-1\}$ y $c\in K$, y por tanto
$\alpha\leftrightarrow j\alpha +\wp(c)$, donde $\wp
(c)=c^p-c$, obtenemos que
la ecuaci\'on que genera a $\K$ puede ser normalizada
como
\begin{equation}\label{E2.2.Ram}
y^p-y=\alpha=\sum_{i=1}^r\frac{Q_i}{P_i^{e_i}} + f(T)=
\frac{Q}{P_1^{e_1}\cdots P_r^{e_r}}+f(T),
\end{equation}
donde $P_i\in R_T^+$, $Q_i\in R_T$, 
$\gcd(P_i,Q_i)=1$, $e_i>0$, $p\nmid e_i$, $\deg Q_i<
\deg P_i^{e_i}$, $1\leq i\leq r$,
$\deg Q<\sum_{i=1}^r\deg P_i^{e_i}$, $f(T)\in R_T$,
con $p\nmid \deg f$ cuando $f(T)\not\in {\ma F}_q$
y $f(T)\notin \wp({\ma F}_q)$ cuando $f(T)
\in {\ma F}_q^{\ast}$ (ver \cite[Example 5.8.8]{Vil2006}).

\begin{proposition}\label{P2.4.Ram3}\label{P2.4}
Sea $\K=K(y)$ dada por la Ecuaci\'on {\rm (\ref{E2.2.Ram})}.
Entonces los primos finitos ramificados son precisamente
$P_1,\ldots,P_r$. 

Con respecto al primo $\p$, tenemos
\l
\item descompuesto si $f(T)=0$.
\item inerte si $f(T)\in {\ma F}_q$ y $f(T)\not\in \wp({\ma F}_q)$.
\item ramificado si $f(T)\not\in {\ma F}_q$ (por lo 
tanto $p\nmid\deg f$).
\end{list}
\end{proposition}

\begin{proof} Sea $P\in R_T^+$, tal que $P\notin \{P_1,\ldots,
P_r\}$. Sea $\pK$ un primo en $\K$ sobre $P$. Al completar
se tiene $\K_{\pK}=K_P(\alpha)$. Por hip\'otesis, $\alpha$
es entero en $\pK$. Entonces $F(Y)=\Irr(\alpha,Y,K_P)=
Y^p-Y-\alpha$ y el diferente local divide a $F'(\alpha)=-1$
(Corolario \ref{C1.3.7'}) por lo que $P$ es no ramificado.

Otra demostraci\'on de lo anterior la damos en la
Proposici\'on \ref{P6.3-1.Ram5}.

Ahora, para $1\leq i\leq r$, si $\pK_i$ es un primo en $\K$
sobre $P_i$, se tiene de la Ecuaci\'on (\ref{E2.2.Ram}) que
si $y$ es una ra\'iz, es decir $y^p-y=\alpha$,
entonces
\begin{gather}\label{EcRam1}
v_{\pK_i}(y^p-y)= pv_{\pK_i}(y)=v_{\pK_i}(\alpha)=
e(\pK_i|P_i)v_{P_i}(\alpha)=e(\pK_i|P_i)e_i.
\end{gather}

Puesto que $\mcd(p,e_i)=1$, se sigue que $p\mid e(\pK_i|
P_i)$ por lo que $P_i$ es ramificado en $\K/K$, $1\leq i
\leq r$.

Con respecto a $\p$, procedemos de la siguiente forma.
Primero consideremos el caso
$f(T)=0$. Entonces $v_{\p}(\alpha)=
\deg(P_1^{e_1}\cdots P_r^{e_r})-\deg Q>0$.
Por el mismo argumento que al principio, se sigue
que $\p$ es no ramificado. 

Ahora bien $y^p-y=\prod_{i=0}^{p-1}
(y-i)$. Sea ${\eu P}_{\infty}\mid\p$. Entonces 
\[
v_{{\eu P}_{\infty}}(y^p-y)=\sum_{i=0}^{p-1}v_{{\eu P}_{\infty}}
(y-i)=e({\eu P}_{\infty}|\p)v_{\p}(\alpha)=v_{\p}(\alpha)>0.
\]
Por lo tanto, existe $0\leq i\leq p-1$ tal que
$v_{\p}(y-i)>0$. Sin p\'erdida de generalidad, podemos
suponer que $i=0$. Sea
$\sigma\in \Gal(K/k)\setminus \{\Id\}$. Suponemos que
${\eu P}_{\infty}^{\sigma}={\eu P}_{\infty}$. 
Tenemos que $y^{\sigma}=
y-j$, $j\neq 0$. Por lo tanto, por un lado
\[
v_{{\eu P}_{\infty}}(y-j)=v_{{\eu P}_{\infty}}(y^{\sigma})=
v_{\sigma({\eu P}_{\infty})}(y)=v_{{\eu P}_{\infty}}(y)>0.
\]
Por otro lado, puesto que
$v_{{\eu P}_{\infty}}(y)>0=v_{{\eu P}_{\infty}}(j)$,
se sigue que
\[
v_{{\eu P}_{\infty}}(y-j)=\min\{
v_{{\eu P}_{\infty}}(y),v_{{\eu P}_{\infty}}(j)\}=0.
\]
Esta contradicci\'on prueba que
${\eu P}_{\infty}^{\sigma}\neq {\eu P}_{\infty}$
por lo que $\p$ se descompone en $K/k$.

Ahora consideramos el caso $f(T)\neq 0$.
Si $f(T)\not\in \F$, entonces si $\pL$ es un
primo en $\K$ sobre $\p$, puesto que 
$v_{\p}(f(T))<0$, se sigue que
\[
v_{\pL}(y^p-y)=pv_{\pL}(y)=e(\pL|\p)
v_{\p}(y)=e(\pL|\p)(-\deg f(T)).
\]
Puesto que $\mcd(p,\deg f(T))=1$, se sigue que $p\mid
e(\pL|\p)$ y por tanto $\p$ es ramificado en $\K/K$.

El \'ultimo caso es cuando
$f(T)\in\F$, $f(T)\not\in \wp(\F)$. Sea
$b\in {\ma F}_{q^p}$ con $b^p-b=a=f(T)$. Puesto que $\deg \p=1$,
$\p$ es inerte en la extensi\'on de constantes $K(b)/K$
(Teorema \ref{T6.1.4}). Supongamos que $\p$
se descompone en $\K=K(y)/K$. Tenemos el siguiente diagrama 
\[
\xymatrix{
K(y)\ar@{-}[r]^{\p}_{\text{inerte}}\ar@{-}[d]_{\substack{\p\\
\text{se descompone}}}
&K(y,b)\ar@{-}[d]\\ K\ar@{-}[r]^{\p}_{\text{inerte}}&K(b)
}
\]

El grupo de descomposici\'on de $\p$ en $K(y,b)/K$ es
$\Gal(K(y,b)/K(y))$. Por tanto $\p$ es inerte  en cada campo
de grado $p$ sobre $K$ diferente a $\K=K(y)$.
Puesto estos campos de grado $p$ son
$K(y+ib), K(b)$, $0\leq i\leq p-1$, se tiene que
en $K(y+b)/K$
\[
(y+b)^p-(y+b)=(y^p-y)+(b^p-b)=\alpha-a =\frac{Q}{P_1^{
e_1}\cdots P_r^{e_r}}
\]
con $\deg (\alpha-a)<0$. Por tanto, por la primera parte,
$\p$ se descompone en $K(y+b)/K$ y en
$K(y)/K$ lo cual es imposible. Por tanto
$\p$ es inerte en $K(y)/K$. $\fin$
\end{proof}

La ramificaci\'on de un primo en extensiones radicales o
extensiones de Artin--Schreier para campos de funciones
arbitrarios puede darse en general.

\begin{proposicion}\label{P6.3-1.Ram5}
Sea $k$ un campo arbitrario de caracter\'istica $p$ y sea
$\K$ un campo de funciones con campo de constantes
$k$. Sea $L=\K(y)$ donde 
\[
y^p-y=\alpha\in \K.
\]
Sea $\P$ un primo en $\K$. Se tiene
\las
\item Si $v_{\P}(\alpha)\geq 0$, entonces
$\P$ es no ramificado.

\item Si $v_{\P}(\alpha)=-\lambda<0$ con $\mcd(\lambda,p)
=1$, entonces $\P$ es ramificado.
\end{list}
\end{proposicion}

\begin{proof}
Supongamos primero que $v_{\P}(\alpha)\geq 0$. Supongamos
que $\P$ es ramificado. Sea $\pK$ el primo en $L$ sobre $\P$.
Entonces $v_{\pK}(y)\geq 0$ pues de lo contrario
$0>v_{\pK}(y^p-y)=pv_{\pK}(y)=e(\pK|\P)v_{\P}(\alpha)>0$
lo cual es absurdo. Por tanto $y\in{\mc O}_{\pK}$ es entero
en $\pK$. Ahora bien $G=\Gal(L/\K)=I(\pK|\P)=\langle\sigma\rangle$
es un grupo c\'iclico de orden $p$. Por la definici\'on del
grupo de inercia tenemos $v_{\pK}(\sigma y-y)>0$ pero
$\sigma(y)=y+1$ y por tanto $v_{\pK}(\sigma y-y)=v_{\pK}
(1)=0$. Esta contradicci\'on muestra que $\P$ en no ramificado
en $L/\K$.

Ahora si $v_{\P}(\alpha)=-\lambda<0$ con $\mcd(\lambda,p)=$,
entonces si $e$ denota el \'indice de ramificaci\'on de $\P$, entonces
si $\pK$ es un primo en $L$ sobre $\P$, entonces
\[
v_{\pK}(y^p-y)=pv_{\pK}(y)=v_{\pK}(\alpha)=ev_{\P}(\alpha)
=-e \lambda.
\]
Puesto que $\mcd(\lambda,p)=1$, entonces $p\mid e$ y $\P$
es ramificado.
$\fin$
\end{proof}

\begin{observacion}\label{O6.3-2.Ram6}
Si el campo $k$ es perfecto, toda extensi\'on c\'iclica $L/\K$
de grado $p$ y cualquier $\P$ primo en $\K$, puede ser
llevado a la forma de la Proposici\'on \ref{P6.3-1.Ram5}.
Cuando $k$ no es perfecto, esto no necesariamente se
cumple.
\end{observacion}

An\'alogo al caso de Artin--Schreier, tenemos el caso
de las extensiones radicales.

\begin{proposicion}\label{P6.4-1.Ram7}
Sea $k$ un campo perfecto cualquiera y sea $\K/k$ un
campo de funciones con campo de constantes $k$.
Sea $n\in{\ma N}$ primo relativo a la caracter\'istica
de $k$. Sea $L=\K(\sqrt[n]{\alpha})$. Sea $\P$
un primo en $\K$. Entonces
\las
\item si $v_{\P}(\alpha)=0$, entonces $\P$ es no
ramificado en $L/\K$,

\item si $v_{\P}(\alpha)=m$ con $1\leq m\leq n-1$, entonces
$\P$ es ramificada en $L$ con \'indice de ramificaci\'on
$\frac{n}{\mcd(n,m)}$.
\end{list}

\end{proposicion}

\begin{proof}
Notemos que podemos suponer que la extensi\'on $L/\K$
es de Kummer pues en caso de no serlo, consideremos $\K'=
\K(\zeta_n)$ y $L'=L(\zeta_n)$. Entonces $L'=\K'(\sqrt[n]{
\alpha})$ es una extensi\'on de Kummer y puesto que 
$k$ es perfecto, $L'/L$ y $\K'/\K$ son no ramificadas por lo
que $\P$ es ramificado en $L/\K$ si y solo si es ramificado
en $L'/\K'$ y con el mismo \'indice de ramificaci\'on.

Sea $y^n=\alpha$, $L=\K(y)$.

Consideremos el caso $v_{\P}(\alpha)=0$. Entonces si $\pK$
es un primo en $L$ sobre $\K$, tenemos que
\[
v_{\pK}(y^n)=nv_{\pK}(y)=v_{\pK}(\alpha)=
e(\pK|\P)v_{\P}(\alpha)=0.
\]
Por tanto $y$ es entero en $\pK$. Si $\Gal(L/\K)=\langle
\sigma\rangle$, el grupo de inercia es de la forma
$I(\pK|\P)=\langle\sigma^j\rangle$ con $1\leq j\leq n-1$
y se tiene $v_{\pK}(\sigma^j (y)-y)>0$. Por otro lado
$\sigma(y)=\zeta_n y$ y $\sigma^j(y)=\zeta_n^j y$ por lo
que $v_{\pK}(\sigma^j(y)-y)=v_{\pK}(\zeta_n^j y-y)=
v_{\pK}((\zeta_n^j-1)y)=v_{\pK}(y)=0$. Esta contradicci\'on
muestra que $\P$ es no ramificado.

El resto de la demostraci\'on es completamente similar
a la del Teorema \ref{TRam1}.
$\fin$
\end{proof}

Regresamos al caso de campos de funciones congruentes.
Usamos la Proposici\'on \ref{P2.4.Ram3} 
para estudiar la descomposici\'on
de los primos no ramificados en $\K/K$.

Sea $\K=K(y)$ dado por (\ref{E2.2.Ram})
donde suponemos que $\K/K$ es geom\'etrica. Sea $\P$
un divisor primo tal que $\P\notin\{\P_1,\ldots,\P_r,\p\}$. 
Entonces $\P$ es ya sea inerte o se descompone $\K/K$. 
El siguiente resultado establece el tipo de
descomposici\'on de $\P$.

\begin{proposicion}\label{P6.3.Ram4} Sean
$\P$ y $\K/K$ como antes. Escribamos la
Ecuaci\'on {\rm{(\ref{E2.2.Ram})}} como
\[
y^p-y=\alpha=u(T),
\]
con $u(T)=\frac{g(T)}{h(T)}\in K$ tal que
$\mcd(g(T),h(T))=1$. Sea
$P(T)\in R_T^+$ el polinomio irreducible asociado a
$\P$, digamos de grado 
$\deg P(T)=m$. Sea $\nu\in {\ma F}_{q^m}$ 
una ra\'iz de $P(T)$.
Entonces $\P$ se descompone en $\K/K$ si y solamente
si $u(\nu)\in \wp({\ma F}_{q^m})$.
\end{proposicion}

\begin{proof} Se tiene $[{\ma F}_{q^m}:\F]=m$. Sean
$\nu=\nu_1,\ldots,\nu_m$ las ra\'ices de $P$ 
en ${\ma F}_{q^m}$, $P(T)=\prod_{i=1}^m
(T-\nu_i)$. Tenemos que $\P$ se descompone
totalmente en  $K_m/K$. Aqu\'i $K_m$ 
denota la extensi\'on de constantes de $K$ de grado $m$.
Puesto que estamos suponiendo que $\K/K$ es geom\'etrica,
se tiene que $K_m\cap \K=K$ y por tanto $\P$
se descompone totalmente en $\K/K$ si y solamente si
${\eu Q}$ se descompone en $\K_m/K_m$ donde $\K_m=
\K K_m$ y ${\eu Q}$ es un primo de $K_m$ sobre $\P$.
\[
\xymatrix{\K\ar@{-}[rr]\ar@{-}[d]&&\K_m=\K K_m\ar@{-}[d]\\
K\ar@{-}[rr]_{\substack{\text{$\P$ se}
\\ \text{ descompone}}}&&K_m}
\]

Digamos que ${\eu Q}$ es el primo asociado a
$T-\nu\in K_m$. Puesto que $v_{\P}(u(T))\geq 0$
se sigue que $v_{{\eu Q}}(u(T))\geq0$ de tal forma que
$T-\nu\nmid h(T)$ y $h(\nu)\neq 0$. 
M\'as a\'un $g(\nu)=0\iff v_{\P}(u(T))>0$.
Tenemos que $\deg_{K_m} {\eu Q}=1$. 
Hacemos de ${\eu Q}$ el primo infinito en
$K_m$, esto es, sea $T'=\frac{1}{T-\nu}$,
$(T')_{K_m}=\frac{{\eu Q}'_0}{{\eu Q}'_{\infty}}=
\frac{{\eu Q}_{\infty}}{{\eu Q}}$, donde $(T)_{k_m}=\frac{
{\eu Q}_0}{{\eu Q}_{\infty}}$. Tenemos $T=\frac{1}{T'}+\nu$.

Escribamos $u_1(T'):=u(T)=u\big(\frac{1}{T'}+
\gamma\big)=\frac{g\big(\frac{1}{T'}+\nu\big)}
{h\big(\frac{1}{T'}+\nu\big)}=\frac{g\big(\frac{1}{T'}(1+T'\nu)\big)}
{h\big(\frac{1}{T'}(1+T'\nu)\big)}$.

Sea $g(T)=a_sT^s+a_{s-1}T^{s-1}+\cdots+a_1T+
a_0$, $a_s\neq 0$, $a_i\in \F$, 
$0\leq i\leq s$; $h(T)=b_tT^t+b_{t-1}T^{t-1}+\cdots+
b_1T+b_0$, $b_t\neq 0$, $b_j\in \F$, 
$0\leq j\leq t$. Entonces
\begin{gather*}
g\Big(\frac{1}{T'}(1+T'\nu)\Big)=
\frac{1}{(T')^s}\Big(a_s+\cdots+g(\nu)(T')^s\Big)=
\frac{1}{(T')^s}g_1(T');\\
h\Big(\frac{1}{T'}(1+T'\nu)\Big)=
\frac{1}{(T')^t}\Big(b_t+\cdots+h(\nu)(T')^t\Big)=
\frac{1}{(T')^t}h_1(T').
\intertext{Se sigue que}
\deg_{T'}g_1(T')\leq s;\quad\deg_{T'}h_1(T')=t.
\intertext{Por lo tanto}
\deg_{T'}u_1(T')=\deg_{T'}g_1(T')
-s\leq 0,\quad\text{y}\quad v_{{\eu Q}}(u_1(T'))=s-\deg_{
T'}(g_1(T'))\geq 0.
\end{gather*}

La forma reducida de $u_1(T')$ es
\[
u_1(T')=\sum_{j=1}^{r'}\frac{Q'_j(T')}{(P'_j(T'))^{e'_j}}+u(\nu).
\]
Por tanto ${\eu Q}$ se descompone en 
$\K_m/K_m\iff u(\nu)\in 
\wp({\ma F}_{q^m})$. Esto prueba la proposici\'on. $\fin$
\end{proof}

\begin{observacion}\label{O6.3.Ram5}
Sea $y^p-y=\alpha=u(T)$ como en la Proposici\'on
\ref{P6.3.Ram4}. Sea $P(T)
\in R_T^+$ el polinomio irreducible asociado al primo ${\mc P}$
con $\deg P(T)=m$ y donde $\mu\in{\ma F}_{q^m}$ es una ra\'iz de
$P(T)$. Supongamos que la ecuaci\'on $y^p-y=u(T)$
est\'a en forma normal para ${\mc P}$, es decir,
$v_{\mc P}(u(T))\geq 0$ \'o $v_{\mc P}(u(T))<0$ y $\mcd(v_{\mc P}
(u(T)),p)=1$. 

Supongamos que ${\mc P}$ se descompone totalmente en $\K/K$.
Entonces $u(\mu)\in \wp({\ma F}_{q^m})$. Veamos otra forma
equivalente de esta afirmaci\'on.

Se tiene $X^p-X-u(T)\in K[X]$, $K=\F(T)$. Entonces $X^p-X-u(T)
\bmod P(T)=X^p-X-u(\mu)\in R_T[X]/\langle P(T)\rangle\cong 
{\ma F}_{q^m}[X]$. Esto se debe a que en $R_T/\langle P(T)
\rangle\cong \F(\mu)$ se tiene que $(T-\mu)|(u(T)-u(\mu))$ por
lo que $u(T)-u(\mu)\equiv 0\bmod P(T)$.

Entonces $u(\mu)\in\wp({\ma F}_{q^m})\iff X^p-X-u(\mu)=
X^p-X-\wp(\xi)$ para alg\'un $\xi\in {\ma F}_{q^m}$.
Ahora bien, $X^p-X-\wp(\xi)=\prod_{i=0}^{p-1} (X-\xi-i)\iff
X^p-X-u(T)\bmod P(T)$ se descompone totalmente.

En resumen, ${\mc P}$ se descompone totalmente en $\K/K
\iff X^p-X-u(T)\bmod P(T)\in {\ma F}_{q^m}[X]$ se descompone
en $p$ factores lineales.

Notemos que esta \'ultima afirmaci\'on es consecuencia inmediata
del Teorema de Kummer (\cite[Theorem 5.8.1]{Vil2006}).
\end{observacion}

Terminamos esta secci\'on resumiendo el comportamiento
de un primo en una extensi\'on de Artin-Schreier $L/\K$ de
campos de funciones con campo de constantes un campo
perfecto. Parte de este resultado puede consultarse en 
\cite{WuSh2010}.

\begin{teorema}\label{T6.3.Ram6-1}
Sea $\K/k$ un campo de funciones con $k$ un campo perfecto
de caracter\'istica $p$.
Sea $L=\K(y)$ donde $y^p-y=\alpha\in \K\setminus \wp(\K)$.
Entonces $L/\K$ es una extensi\'on c\'iclica de grado $p$.
Sea $\pK$ un divisor primo de $\K$ y normalizamos $\alpha$
de tal manera que $v_{\pK}(\alpha)\geq 0$ \'o $v_{\pK}(
\alpha)<0$ y $\mcd(v_{\pK}(\alpha),p)=1$. Entonces
\las
\item $L/\K$ es una extensi\'on de constantes si y solamente
si $\alpha$ puede seleccionarse en $k$.

\item Si $v_{\pK}(\alpha)<0$ entonces $\pK$ es ramificado en $L/\K$.

\item Si $v_{\pK}\geq 0$ y $X^p-X-\alpha\bmod \pK$ se descompone,
entonces $\pK$ se descompone en $L/\K$.

\item Si $v_{\pK}\geq 0$ y $X^p-X-\alpha\bmod \pK$ es 
irreducible, entonces $\pK$ es inerte en $L/\K$.
\end{list}

En particular, si $v_{\pK}(\alpha)>0$, $\pK$ se descompone
en $L/\K$.
\end{teorema}

\begin{proof}{\ }

(1): Si $\alpha\in k$ entonces $L=l\K$ donde $l=k(y)$. Rec\'iprocamente,
sea $L/\K$ es una extensi\'on de constantes. Sea $L=l\K$ con $l/k$
de grado $p$ (ver \cite[Theorem 8.4.4]{Vil2006}). Entonces existen $z\in l$ y $\beta\in k$
tales que $l=k(z)$ con $z^p-z=\beta$, esto es, podemos seleccionar
$\alpha=\beta\in k$.

(2): Es el contenido de la Proposici\'on \ref{P6.3-1.Ram5}.

(3) y (4): Se siguen de la Proposici\'on \ref{P6.3-1.Ram5} y del teorema
de Kummer (\cite[Theorem 5.8.2]{Vil2006}).

Cuando $v_{\pK}(\alpha)>0$, se tiene que $X^p-X-\alpha\bmod \pK=
X^p-X=\prod_{i=0}^{p-1} (X-i)$ se descompone y la afirmaci\'on se
sigue de(3).
$\fin$
\end{proof}

\section{Resultados generales sobre 
ramificaci\'on}\label{SRam3}

En esta secci\'on, algunos de los resultados los enunciamos
para campos de funciones en general, en particular para
ejemplificar que algunas conclusiones para campos de
funciones congruentes se cumplen gracias a que el campo
de constantes $k$ es un campo finito $\F$ que en particular
es un campo perfecto. Algunas conclusiones no se cumplen
para campos $k$ m\'as generales, particularmente cuando
$k$ no es perfecto.

Empezamos con el {\em Lema de Abhyankar}, el cual
juega un papel central en la teor\'ia de ramificaci\'on.

\begin{teorema}[Lema de Abhyankar]\label{T2.3.B}
Sea $F/E$ una extensi\'on finita y separable de campos
de funciones. Supongamos que
$F=E_1E_2$ con $E\subseteq E_i\subseteq F$. 
Sean ${\P}$ un primo de $E$ y $\pL$ un primo de $F$ dividiendo
a $\P$. Sean ${\eu p}_i={\eu P}
\cap E_i$ para $i=1,2$. Si al menos una de las
extensiones $E_i/E$ es moderadamente ramificada en
${\eu p}_i$, entonces
\[
e_{F|E}({\eu P}|{\P})=\mcm[e_{E_1|E}({\eu p}_1|{\P}), e_{E_2|E}
({\eu p}_2|{\P})].
\]
\end{teorema}

\begin{proof} \cite[Theorem 12.4.4]{Vil2006}. $\fin$
\end{proof}

\begin{observacion}\label{ORam6}
Con las notaciones del Teorema \ref{T2.3.B}, puesto que
$e_{E_i|E}({\eu p}_i|{\P})|e_{F|E}({\eu P}|{\P})$, para $i=1,2$,
se sigue que 
\[
\mcm[e_{E_1|E}({\eu p}_1|{\P}), e_{E_2|E}
({\eu p}_2|{\P})]|e_{F|E}({\eu P}|{\P})
\]
en completa generalidad.

En el Lema de Abhyankar, si las dos extensiones son
salvajemente ramificadas, entonces la conclusi\'on no
se cumple. Por ejemplo, consideremos $F={\ma F}_9
(T)$ y sean $E_1=F(y)$ con $y^3-y=T$ y $E_2=F(z)$
con $z^3-z=\mu T$, donde $\mu\in {\ma F}_9\setminus
{\ma F}_3$. Sea $E=E_1E_2=F(y,z)$.
Veamos que $E_1\neq E_2$. Si
tuvi\'esemos $E_1=E_2$ entonces existir\'ian
$j\in\{1,2\}$ y $c\in F$ tal que $\mu T=jT+\wp(c)$.
Como $v_{\pK_F}(T)=-1$, donde $\pK_F$ es 
el primo infinito de $F$, se sigue que $v_{\pK_F}
(\wp(c))=v_{\pK_F}((\mu-j)T)=-1$ y por tanto 
$v_{\pK_F}(c)\geq 0$ pues si $v_{\pK_F}(c)
<0$, entonces $v_{\pK_F}(\wp(c))=-3$.
Esto es absurdo pues $(\mu-j)T=\wp(c)$.

Todas las extensiones intermedias de grado $3$
entre $E/F$ son $E_1,E_2,F(y+z)=F(w)$ y $F(y+2z)=
F(u)$ y como $w^3-w=(1+\mu)T$ y $u^3-u=(1+2\mu)
T$, $\pK_F$ se ramifica totalmente en $E/F$. Por
tanto $e(E|F)=9\neq \mcm[e(E_1|F),e(E_2|F)]=
\mcm[3,3]=3$.
\end{observacion}

\begin{observacion}\label{O.Ram11}
El Lema de Abhyankar no es aplicable 
para extensiones inseparables.
\end{observacion}

\begin{ejemplo}\label{E.Ram12}
Consideremos $k$ un campo imperfecto de caracter\'istica $p>0$
y sea $a\in k\setminus k^p$, esto es, $b=a^{1/p}\notin k$. Sean
$F=k(x)$ y $L=F(y)=k(x,y)$
donde $y^p-y=ax^p$. Entonces $L/k(x)=F$ es
una extensi\'on separable, no ramificada y el campo de
constantes de $L$ es $k$ (ver m\'as adelante). 
\[
\xymatrix{
&&\underbracket[0pt]{M}_{\substack{\uigual\\
L(b)=k(b,x,y)=F(b,y)=E(y)}}
\ar@{-}[ddll]_p\ar@{-}[ddrr]^p\ar@{-}[d]\\
&&\pL\\ \underbracket[0pt]{L}_{
{\substack{ \uigual\\ k(x,y)=F(y)}}}\ar@{-}[ddrr]_p
\ar@{-}[r]&{\eu q}\ar@{--}[ur]_{e=p}
\ar@{--}[dr]^{e=1}&&\p\ar@{--}[lu]^{e=p}
\ar@{--}[dl]_{e=1}\ar@{-}[r]
&\underbracket[0pt]{E}_{\substack{\uigual\\ 
F(b)=k_1(x)}}\ar@{-}[ddll]^p\\
&& {\mc P}_{\infty}\ar@{-}[d]\\
&&\underbracket[0pt]{F}_{\substack{\uigual\\ k(x)}}
}
\]
Sea $E=F(a^{1/p})=F(b)=k_1(x)$
donde $k_1=k(b)=k(a^{1/p})$. Sea ${\eu q}$ un primo en $L$ sobre
${\mc P}_{\infty}$, el primo infinito de $k(x)$. Se tiene que
${\eu q}/{\mc P}_{\infty}$ es puramente inseparable (ver m\'as 
adelante).

Ahora consideremos $y^p-y=ax^p$ sobre $E=F(b)$. Sea
$z=y-a^{1/p}x=y-bx$. Entonces
\[
z^p-z=(y^p-b^px^p)-(y-bx)=(y^p-y)-b^px^p+bx=bx.
\]
Por tanto $M=F(b)(y)=E(y)=E(z)$ est\'a definida por $z^p-z=bx$
y por tanto el primo infinito $\p$ de $E=k_1(x)$ es ramificado
en $M/E$.

Puesto que ${\mc P}_{\infty}$ no es ramificado en $L/F$
y $\p$ es ramificado en $M/E$, se sigue que $e_{M/F}(
\pL|{\mc P}_{\infty})=p$. Como consecuencia se tiene
que $e_{E/F}(\p|{\mc P}_{\infty})=1$ y $e_{M/L}(\pL|{\eu q})
=p$. En particular
\begin{gather*}
\mcm[e_{L/F}({\eu q}|{\mc P}_{\infty}), e_{E/F}(
\p|{\mc P}_{\infty})]=\mcm[1,1]=1<e_{M/F}(\pK|\p)=p.
\end{gather*}

Aclaremos algunas afirmaciones que nos sirvieron
para el ejemplo. Primero, $L/F$ es separable pues
$L=F(y)$ y $y^p-y=ax^p$ tiene ra\'ices distintas:
$(y^p-y)'=-1\neq 0$. Veamos que el campo
de constantes de $L$ es $k$. Si el campo de
constantes de $k$ fuese $k'\neq k$, entonces
$[k':k]\geq p$, $L=k'(x)$ y existir\'ia $c\in k'\setminus k$
tal que $L=F(c)$ con $c^p-c=d\in k$. Entonces,
puesto que $L=F(c)=F(y)$ se tendr\'ia 
\begin{gather}\label{EqRam11}
ax^p=jd+\xi^p-\xi,
\end{gather}
con $\xi\in F$. Por tanto $-p=v_{{\mc P}_{\infty}}(ax^p)=
-pv_{{\mc P}_{\infty}}(\xi)$. Se ve que la Ecuaci\'on
(\ref{EqRam11}) no tiene soluci\'on para $\xi$, lo cual
se sigue directamente al escribir $\xi=\frac{l(x)}{g(x)}$
y usando que $a\notin k^p$.

Alternativamente, consideremos los divisores:
\begin{align*}
(y^p-y)_L&=(y)_L (y-1)_L\cdots (y-(p-1))_L=
\frac{{\eu q}_{0,1}}{\eu q}\frac{{\eu q}_{0,2}}{\eu q}
\cdots \frac{{\eu q}_{0,p}}{\eu q}\\
&=(ax^p)_L=
\Big(\frac{{\mc P}_0}{{\mc P}_{\infty}}\Big)^p
\end{align*}
y los divisores $\{{\eu q}_{0,i}\}_{i=1}^p$ son
primos relativos a pares.

Por tanto $\con_{F/L}({\mc P}_0)={\eu q}_{0,1}\cdots
{\eu q}_{0,p}$, es decir ${\mc P}_0$ se descompone
totalmente en $L/F$ por lo que
\[
{\mc O}_{{\eu q}_{0,i}}/{\eu q}_{0,i}\cong {\mc O}_{
{\mc P}_0}/{\mc P}_0\cong k.
\]
Por tanto el campo de constantes de $L$ es $k$.

Ahora veamos que $L/F$ es no ramificada. Si $\pK$
es cualquier lugar de $L$ diferente a ${\mc P}_{\infty}$,
entonces $v_{\pK}(ax^p)\geq 0$ por lo que $\pK$
en no ramificado en $L/F$ (Proposici\'on \ref{P6.3-1.Ram5}).
Analicemos ${\mc P}_{\infty}$ y ${\eu q}\mid {\mc P}_{\infty}$.
Se tiene de $y^p-y=ax^p$ que $v_{\eu q}(y^p-y)=v_{
\eu q}(ax^p)=e_{L/F}({\eu q}|{\mc P}_{\infty}) p v_{{\mc P}_{
\infty}}(x)=-pe_{L/F}({\eu q}|{\mc P}_{\infty})$,
de donde $v_{\eu q}(y)=-e_{L/F}({\eu q}|{\mc P}_{\infty})
=v_{\eu q}(x)$. 

Se sigue que $\frac{y}{x}\in {\mc O}_{\eu q}\setminus
{\eu q}$. Sea $\frac{y}{x}\equiv c\bmod{\eu q}$,
$\big(\frac{y}{x}\big)^p\equiv c^p\bmod {\eu q}$. 
Por otro lado, se tiene $\frac{y}{x^p}\in{\eu q}$. Por
lo tanto 
\[
a=\frac{y^p-y}{x^p}=\Big(\frac{y}{x}\Big)^p-\Big(\frac{y}
{x^p}\Big)\equiv \Big(\frac{y}{x}\Big)^p\equiv c^p\bmod {\eu q},
\]
es decir, $c^p\equiv a\bmod {\eu q}$. Por lo tanto
$c=a^{1/p}\in {\mc O}_{\eu q}/{\eu q}$. Concluimos que
${\mc O}_{\eu q}/{\eu q}=k(a^{1/p})$. Se sigue
${\mc O}_{\eu q}/{\eu q}$ es puramente inseparable
sobre $k={\mc O}_{{\mc P}_0}/{\mc P}_0$. Esto es,
${\eu q}/{\mc P}_{\infty}$ es puramente inseparable
y en particular no ramificada. 

Por tanto $L/F$ es una extensi\'on no ramificada.
\end{ejemplo}

\begin{observacion}\label{O.Ram13}
El Ejemplo \ref{E.Ram12} sirve para ejemplificar
varias cosas. Entre otras:
\las
\item El Lema de Abhyankar no se cumple para
extensiones inseparables.

\item La extensiones $LE/E$ y $LE/L$ son ramificadas
pero $L/F$ y $E/F$ no lo son (ver Corolario \ref{C.10.4.2.Ram}
(2) m\'as adelante).

\item Existen extensiones de constantes ramificadas
($M/L$).

\item Existen extensiones separables geom\'etricas
no ramificadas ($L/F$). 

\item Existen extensiones separables con extensiones
de campos residuales inseparables ($L/F$ y ${\eu q}/
{\mc P}_{\infty}$).
\end{list}

Todos estos fen\'omenos se deben ya sea a la inseparabilidad
o que el campo de constantes no es perfecto.
\end{observacion}

\begin{teorema}\label{T.Ram10}
Sea $L/\K$ una extensi\'on finita y separable de campos
de funciones congruentes. Sea $\pK$ un primo en $\K$.
Sea $\tilde{L}$ la cerradura de Galois de $L/\K$. Entonces
\las
\item $\pK$ es ramificado en $L/\K$ si y solamente si $\pK$
es ramificado en $\tilde{L}/\K$.
\item $\pK$ es totalmente descompuesto en $L/\K$ si y
solamente si $\pK$ es totalmente descompuesto en
$\tilde{L}/\K$.
\end{list}
\end{teorema}

\begin{proof} (1). (Ver Teorema \ref{T1.5.3}).
Si $\pK$ es ramificado en $L/\K$, es inmediato que $\pK$
es ramificado en $\tilde{L}/\K$.

Rec\'iprocamente, supongamos que
$\pK$ es no ramificado en $L/\K$. Entonces $e_{
L/\K}({\eu q}|\pK)=1$ para todo primo
${\eu q}$ de $L$ sobre $\pK$.
Por tanto $e_{L^{\sigma}/\K}({\eu q}^{\sigma}|\pK)=1$
para todo primo ${\eu q}^{\sigma}$ de $L^{\sigma}$
sobre $\pK$.

Sea $R=\{\sigma\colon L\lra 
\bar{\K}\mid \sigma(a)=a\text{\ 
para $a\in \K$}\}$  donde $\bar{\K}$
denota una cerradura algebraica de $\K$.
 Entonces $\tilde{L}=
\prod_{\sigma\in R}L^{\sigma}$.

Sea $\pL$ un primo en $\tilde{L}$ dividiendo a $\pK$.
Sea $I=I_{\tilde{L}/\K}(\pL|\pK)$ el grupo de inercia de
$\pL/\pK$. Sea $F:=\tilde{L}^I$. Entonces $I=
\Gal(\tilde{L}/F)$. Sea $\sigma\in R$ y sea
$H_{\sigma}:=\Gal(\tilde{L}/L^{\sigma})$. Entonces
\[
I_{\tilde{L}/\K}(\pL|\pK)\cap H_{\sigma}=I_{\tilde{L}/
L^{\sigma}}(\pL|\pL\cap L^{\sigma}).
\]
Sea $\pL\cap L^{\sigma}={\eu q}_1$ el cual es un primo
de $L^{\sigma}$ sobre $\pK$. En particular
tenemos $I_{\tilde{L}/L^{\sigma}}(\pL|{\eu q}_1)
\subseteq I_{\tilde{L}/\K}(\pL|\pK)$.

Adem\'as, debido que $e_{L^{\sigma}/\K}(
{\eu q}_1|\pK)=1$ 
y por tanto 
\begin{gather*}
e_{\tilde{L}/\K}(\pL|\pK)=e_{\tilde{L}
/L^{\sigma}}(\pL|{\eu q}_1) e_{L^{\sigma}/
\K}({\eu q}_1|\pK)=e_{\tilde{L}
/L^{\sigma}}(\pL|{\eu q}_1),\\
e_{\tilde{L}/\K}(\pL|\pK)=|I_{\tilde{L}/\K}(\pL|\pK)|=
|I_{\tilde{L}/L^{\sigma}}(\pL|{\eu q}_1)|=e_{\tilde{L}
/L^{\sigma}}(\pL|{\eu q}_1)
\end{gather*}
se sigue que $I_{\tilde{L}/\K}(\pL|\pK)=
I_{\tilde{L}/L^{\sigma}}(\pL|{\eu q}_1)$.

Puesto que $I_{\tilde{L}/\K}(\pL|\pK)\cap H_{\sigma}=I_{\tilde{L}/
L^{\sigma}}(\pL|{\eu q}_1)=I_{\tilde{L}/\K}(\pL|\pK)$,
se sigue que $I_{\tilde{L}/\K}(\pL|\pK)\subseteq H_{\sigma}$.

Por lo tanto $\tilde{L}^I=F\supseteq \tilde{L}^{H_{\sigma}}=
L^{\sigma}$. Se sigue que $\tilde{L}=\prod_{\sigma\in R}
L^{\sigma}\subseteq F\subseteq \tilde{L}$. Entonces
tenemos $F=\tilde{L}$ por lo que $I_{\tilde{L}/\K}(\pL|\pK)=
\{1\}$ y por tanto $\pK$ no es ramificado en $\tilde{L}$.

(2). (Ver Teorema \ref{T1.5.4}) Si $\pK$ es totalmente descompuesto en $\tilde{L}/\K$,
entonces $\pK$ es totalmente descompuesto en $L/\K$.

Rec\'iprocamente, si $\pK$ es totalmente descompuesto
en $L/\K$ se tiene que $\pK$ es no ramificado en
$\tilde{L}/\K$. Si en el inciso anterior cambiamos
los grupos de inercia por los grupos de descomposici\'on
y los \'indices de ramificaci\'on por los grados de inercia,
obtenemos que si ahora $F=\tilde{L}^D$, donde 
$D$ es el grupo de descomposici\'on de $\pL/\pK$,
entonces llegamos a $F=\tilde{L}$ y por tanto
$\pK$ es totalmente descompuesto en $\tilde{L}/\K$.
$\fin$
\end{proof}

Con las mismas t\'ecnicas de la demostraci\'on del Teorema
\ref{T.Ram10}, se puede demostrar el siguiente resultado.

\begin{proposicion}\label{P.Ram10(1)}
Sean $L/E$ y $F/E$ dos extensiones separables finitas
donde el primo $\pK$ de $E$ es totalmente descompuesto.
Entonces $\pK$ es totalmente descompuesto en $LF/E$.
\end{proposicion}

\begin{proof}
Por ejemplo podemos tomar $\tilde{L}$ y $\tilde{F}$ las
cerraduras de Galois de $L/E$ y de $F/E$ respectivamente.
Entonces $\pK$ se descompone totalmente en
$\tilde{L}\tilde{F}/E$. $\fin$
\end{proof}

Dos resultados importantes sobre extensiones
abelianas finitas, son los siguientes.

\begin{proposicion}\label{Palestine3.1}\label{P3.1} Sean $L/K$ 
una extensi\'on abeliana finita 
de campos de funciones globales,
$P\in R_T^+$ y $d:=\deg P$. Supongamos que $P$
es moderadamente ramificada en
$L/K$. Si $e$ denota el \'indice de ramificaci\'on de
$P$ en $L/K$, entonces $e\mid q^d-1$,
donde el campo de constantes de $K$ es $\F$.
\end{proposicion}

\begin{proof}
Primero consideramos en general una extensi\'on finita
de Galois $L/K$. Sean $G_{-1}=D$ el grupo de descomposici\'on
de $P$, $G_0=I$ el grupo de inercia y
$G_i$, $i\geq 1$ los grupos de ramificaci\'on. Sea
${\eu P}$ un divisor primo de $L$ que divide a
$P$. Entonces si ${\mathcal O}_{{\eu P}}$
denota el anillo de valuaci\'on de ${\eu P}$, tenemos
\[
U^{(i)}=1+{\eu P}^i\subseteq {\mathcal O}_{\eu P}^{\ast}=
{\mathcal O}_{\eu P}\setminus {\eu P},
i\geq 1,  U^{(0)}={\mathcal O}_{\eu P}^{\ast}.
\]

Sea $l({\eu P}):= {\mathcal O}_{\eu P}/{\eu P}$ el campo
residual de ${\eu P}$. Los siguientes son monomorfismos:
\begin{eqnarray*}
G_i/G_{i+1}&\stackrel{\varphi_i}{\hooklongrightarrow} &U^{(i)}/
U^{(i+1)}\cong 
\begin{cases}
l({\eu P})^{\ast}, i=0\\
{\eu P}^i/{\eu P}^{i+1}\cong l({\eu P}), i\geq 1.
\end{cases}\\
\overline{\sigma}&\longmapsto& \sigma\pi/\pi
\end{eqnarray*}
donde $\pi$ denota un elemento primo de ${\eu P}$.

Probaremos que si $G_{-1}/G_1=D/G_1$ es abeliano, entonces
\[
\varphi=\varphi_0\colon G_0/G_1\longrightarrow 
U^{(0)}/U^{(1)}\cong
\big({\mathcal O}_{\eu P}/{\eu P}\big)^{\ast}
\]
satisface que $\im \varphi\subseteq {\mathcal O}_P/(P)\cong
R_T/(P)\cong {\ma F}_{q^d}$. En particular,
en nuestro caso, se seguir\'a que
$\big|G_0/G_1\big|\mid \big|{\ma F}_{q^d}^{\ast}\big|=q^d-1$.

Para probar la afirmaci\'on anterior, notemos que
\[
\Gal(({\mathcal O}_{\eu P}/
{\eu P})/({\mathcal O}_P/(P)))\cong D/I=G_{-1}/G_0
\]
(ver \cite[Corollary 5.2.12]{Vil2006}).

Sea $\sigma\in G_0$ y $\varphi(\bar{\sigma})=\varphi(\sigma
\bmod G_1)=[\alpha]= 
\alpha \bmod  {\eu P}\in \big({\mathcal O}_{\eu P}/
{\eu P}\big)^{\ast}$.
Por lo tanto $\sigma\pi\equiv \alpha\pi\bmod {\eu P}^2$.

Sea $\theta\in G_{-1}=D$ arbitrario y sea
$\pi_1:=\theta^{-1} \pi$. Entonces $\pi_1$ 
es un elemento primo de ${\eu P}$. Puesto que
$\varphi$ es independiente del elemento primo,
se sigue que
$\sigma \pi_1\equiv \alpha \pi_1\bmod {\eu P}^2$,
esto es, $\sigma\theta^{-1}\pi\equiv \alpha
\theta^{-1}\pi\bmod {\eu P}^2$.
Puesto que $G_{-1}/G_1$ es un grupo
abeliano, se tiene que
\[
\sigma\pi=(\theta\sigma\theta^{-1})(\pi)
\equiv \theta(\alpha)\pi\bmod{\eu P}^2.
\]
Por tanto $\sigma\pi\equiv \theta(\alpha)
\pi\bmod{\eu P}^2$ y
$\sigma\pi\equiv \alpha\pi\bmod {\eu P}^2$. 
Se sigue que
$\theta(\alpha)\equiv \alpha\bmod {\eu P}$
para toda $\theta\in G_{-1}$.

Si escribimos $\tilde{\theta}=
\theta\bmod G_0$, entonces
$\tilde{\theta}[\alpha]
=[\alpha]$, esto es, $[\alpha]$ es un elemento fijo
bajo la acci\'on del grupo
$G_{-1}/G_0\cong \Gal(({\mathcal O}_{\eu P}/
{\eu P})/({\mathcal O}_P/(P)))$. Obtenemos que
$[\alpha]\in {\mathcal O}_P/(P)$.
Por lo tanto $\im \varphi\subseteq 
\big({\mathcal O}_P/(P)\big)^{\ast}$
y $\big|G_0/G_1\big|\mid \big|
\big({\mathcal O}_P/(P)\big)^{\ast}\big|=
q^d-1$.

Finalmente, puesto que $L/K$ es abeliano
y $P$ es moderadamente ramificado, se tiene
$G_1=\{1\}$, y se sigue que
$e = |G_0|=|G_0/G_1|\mid q^d-1$. $\fin$
\end{proof}

\begin{observacion}\label{OPalestine2.Ram}
Cuando $L/K$ es una extensi\'on de Galois moderadamente
ramificada de campos de funciones congruentes, 
entonces el grupo de inercia $|G_0|$
es un grupo c\'iclico por ser el grupo de
Galois de los campos residuales de $L/K$ y estos
son campos finitos.

En particular, en cualquier extensi\'on moderadamente
ramificada $L/K$ que no sea
c\'iclica, no hay primos totalmente ramificados.
\end{observacion}

\begin{proposicion}\label{PLista10.4}
Sea $K$ un campo de funciones arbitrario de caracter\'istica $p>0$.
Sean $L$ y $E$  dos extensiones de Galois de $K$
de grado $l$ distintas, donde $l\neq p$ es un
n\'umero primo, y tales que $L \cap E = K$.

Sea ${\eu P}_K$ un divisor primo
de $K$. Sean ${\eu P}_L$ y ${\eu P}_E$
lugares de $L$ y $E$ respectivamente, tales
que ${\eu P}_K = {\eu P}_L^l$, 
${\eu P}_K = {\eu P}_E^l$ en $L/K$ y $E/K$
respectivamente, es decir, ${\eu P}_K$ es ramificado tanto en
$L/K$ como en $E/K$.

Sea $F = LE$ y sea
${\eu P}_F$ un lugar de $F$ tal que
${\eu P}_F \mid {\eu P}_K$.

Entonces existe un \'unico campo $M$, $K 
\varsubsetneqq M \varsubsetneqq F$, es decir, $[M:K]=l$, tal
que ${\eu P}_K$ no es
ramificado en $M/K$.
\end{proposicion}

\begin{proof}
Por la demostraci\'on de la
Proposici\'on \ref{Palestine3.1}, se tiene que el
grupo de inercia $I({\eu P}_F | {\eu P}_K)$ es c{\'\i}clico.
Se tiene que $\Gal(F/K)\cong C_l\times C_l$. Puesto
que ${\eu P}_K$ es ramificado tanto en $E/K$ como en 
$L/K$, se tiene que $I({\eu P}_F | {\eu P}_K)\neq
\{\Id\}$ y como $I({\eu P}_F | {\eu P}_K)$ es c{\'\i}clico,
necesariamente $I=I({\eu P}_F | {\eu P}_K)\cong C_l$.
Esto tambi\'en puede obtenerse como consecuencia
del Lema de Abhyankar \ref{T2.3.B}.

Sea $M:=F^I$. Se tiene que $\Gal(M/K)\cong
I({\eu P}_F | {\eu P}_K)\supseteq I({\eu P}_F | {\eu P}_M)$
y todos ellos son grupos de orden $l$. Se sigue
que $I({\eu P}_F | {\eu P}_K)=I({\eu P}_F | {\eu P}_M)$.
\[
\xymatrix{
L\ar@{-}[dd]_l\ar@{-}[rr]&&F=LE\ar@{-}[dd]\\
&M\ar@{-}[ru]^I\ar@{-}[dl]\\
K\ar@{-}[rr]_l&&E
}
\]
Puesto  que $l=e({\eu P}_F|{\eu P}_K)=e({\eu P}_F|{\eu P}_M)
e({\eu P}_M|{\eu P}_K)=l\cdot e({\eu P}_F|{\eu P}_K)$, se sigue que $e({\eu P}_M|{\eu P}_K)=1$,
esto es, ${\eu P}_K$ es no ramificado en $M/K$.

Finalmente, si $N$ fuese otro campo de grado $l$ sobre $K$,
esto es, $N\neq M$, tal que ${\eu P}_K$ no sea ramificado en $N/K$,
entonces $\Gal(F/N)\neq I({\eu P}_F | {\eu P}_K)$,
lo cual implica que
\[
I({\eu P}_F | {\eu P}_N)\subseteq \Gal(F/N)\cap
I({\eu P}_F | {\eu P}_K)\neq I({\eu P}_F | {\eu P}_K),
\]
por lo que $I({\eu P}_F | {\eu P}_N)=\{\Id\}$
y por tanto ${\eu P}_K$ no ser\'ia ramificado en $F/K$
de donde obtenemos que $M$ es \'unico.
$\fin$
\end{proof}

\begin{proposition}\label{Palestine3.3}\label{P11.3.2}
Sea $L/K$ una extensi\'on abeliana de campos
de funciones globales, donde $\F$ es
el campo de constantes de $K$ y donde a lo m\'as
un divisor primo ${\eu p}_0$ de grado $1$ es
ramificado y adem\'as la extensi\'on es moderadamente
ramificada. Entonces $L/K$ es una extensi\'on de
constantes.
\end{proposition}

\begin{proof} Por la Proposici\'on \ref{Palestine3.1} 
se tiene $e:=e_{L/k}({\eu p}_0)|q-1$.
Sea $I$ el grupo de inercia de ${\eu p}_0$.
Entonces $|I|=e$ y ${\eu p}_0$ es no ramificado
en $E:=L^I/K$. Por tanto $E/K$ es una 
extensi\'on no ramificada. Se sigue que $E/K$
es una extensi\'on de constantes (ver la 
Observaci\'on \ref{Palestine5}). La extensi\'on
$L/E$ es geom\'etrica. Supongamos que $L\neq E$.

Sea $[E:K]=m$. Entonces ${\eu P}_0$ es un divisor de
primo en $E$ que divide a ${\eu p}_0$ entonces el grado
relativo $d_{E/K}({\eu P}_0|{\eu p}_0)$ es igual a $m$,
el n\'umero de divisores primos en $E/K$ es $1$ y
el grado de ${\eu P}_0$ es $1$ 
con respecto a ${\ma F}_{q^m}$
(Teorema \ref{T6.1.4}).
Se sigue que ${\eu P}_0$ es el \'unico divisor primo
ramificado en $L/E$ y es de grado $1$
y totalmente ramificado.
Adem\'as $[L:E]=e\mid q-1=|{\ma F}_q^{\ast}|$.

Las $(q-1)$--\'esimas ra\'ices de unidad pertenecen a $\F\subseteq K$.
Por tanto $K$ contiene las $e$--\'esimas ra\'ices de unidad y $L/E$
es una extensi\'on de Kummer, digamos 
$L=E(y)$ con $y^e=\alpha \in E=k{\ma F}_{q^m}
={\ma F}_{q^m}(T)$. Escribimos $\alpha$ en
su forma normal como prescrita por Hasse
\cite{Has35}, \cite[Theorem 5.8.10]{Vil2006}: 
$(\alpha)_E=\frac{{\eu P}_0^a
{\eu a}}{{\eu b}}$, $0<a<e$ y $\mcd({\eu P}_0,
{\eu a}{\eu b})=1$. Ahora bien, puesto
que $\deg (\alpha)_E=0$, se sigue que
$\deg_E {\eu a}$ o $\deg_E {\eu b}$ no es
un m\'ultiplo de $e$. Esto contradice que
${\eu P}_0$ es el \'unico primo ramificado en
$L/E$. M\'as precisamente, existe otro primo
${\eu q}\neq {\eu P}_0$ de $E$ tal que $v_{\eu q}(\alpha)=b$
no es m\'ultiplo de $e$ y por tanto, si ${\eu Q}$ es
un primo en $L$ sobre ${\eu q}$,
\[
v_{\eu Q}(y^e)=ev_{\eu Q}(y)=v_{\eu Q}(\alpha)=
e({\eu Q}|{\eu q})v_{\eu q}(\alpha)= e({\eu Q}|{\eu q})b,
\]
de donde obtenemos que $e({\eu Q}|{\eu q})>1$.
Por lo tanto $L/K$ es una extensi\'on de
constantes. $\fin$
\end{proof}

\begin{observacion}\label{Palestine5}
En la Proposici\'on \ref{Palestine3.3} hemos usado
que si $L/K$ es no ramificada, entonces $L/K$
es una extensi\'on de constantes. Esto se cumple
pues $K$ es de g\'enero $0$. M\'as precisamente,
si $K$ es cualquier campo de funciones (con campo
de constantes arbitrario) de g\'enero $0$ (aunque
no sea de funciones racionales), entonces si
$L/K$ es una extensi\'on geom\'etrica separable no trivial,
entonces $L/K$ es ramificada. De hecho, si $L/K$
es una extensi\'on geom\'etrica separable no trivial y $g_L$
denota el g\'enero de $L$, usando que $g_K=0$, se
obtiene de la f\'ormula de Riemann--Hurwitz que
\[
2g_L-2=[L:K](2g_K-2)+d,
\]
donde $d$ es el grado del diferente. Si $L/K$ no fuese
ramificada, se tendr\'ia $d=0$ y por tanto
$2g_L-2=[L:K](0-2)+0<-2$ lo que implicar\'ia $g_L<0$.

Por otro lado, en nuestro caso de campos de funciones
congruentes, las extensiones de constantes son no
ramificadas, pero en general para campos de funciones
esto no necesariamente se cumple. Por ejemplo,
consideremos $k$ un campo de caracter\'istica $p$
el cual no es perfecto. Sea $a\in k\setminus k^p$.
Sea $K=k(x)$ y sea $L=K(y)=k(x,y)$ donde
$y^p=ax^p$. Entonces $[L:K]=p$ es una extensi\'on
inseparable y $L=k'(x)$ donde $k'=k(a^{1/p})$. Se sigue que
$L/K$ es una extensi\'on de constantes. Consideremos
el primo $\pK$ asociado a $x^p-a$ en $K$. Entonces
puesto que $(x-b)^p=x^p-b^p=x^p-a$, donde
$b=a^{1/p}$, obtenemos que $\pK$
es ramificado en $L/K$.

La diferencia es que en nuestro caso, los campos finitos
son perfectos y por tanto las extensiones de constantes
son separables.
\end{observacion}

\begin{observacion}\label{Palestine6}
Todas las hip\'otesis de la Proposici\'on \ref{Palestine3.3}
son necesarias.
\las
\item Sea $q=3$. El polinomio $P(T)=T^2-T-1$ es irreducible
en ${\ma F}_3[T]$. Sean $K={\ma F}_3(T)$ y $\K=K(\sqrt[2]{P(T)})$.
Entonces, $\K/K$ es una extensi\'on c\'iclica geom\'etrica donde
\'unicamente hay un primo ramificado, a saber, $P(T)$. Por
tanto es necesario pedir que el primo ramificado sea de grado
$1$.

\item Sea $K={\ma F}_p(T)$ y sea $\K=K(y)$ donde $y^p-y=T^{p+1}$
es una extensi\'on geom\'etrica c\'iclica de grado $p$ donde
\'unicamente el primo infinito es ramificado y el cual es de grado
$1$. Por tanto es necesario pedir que la extensi\'on sea moderadamente
ramificada.

\end{list}
\end{observacion}

\begin{observacion}\label{O10.4.1.Ram0}
Cuando el campo de constantes $k$ es un campo algebraicamente
cerrado de caracter\'istica $p\geq 0$ y $L/K=k(x)$ es una
extensi\'on c\'iclica de grado $\ell^n$ con $\ell$ un primo
diferente a $p$, entonces en la extensi\'on $L/K$ hay al menos
dos primos totalmente ramificados.

En efecto, si $L=K(y)$ con $y^{\ell^n}=f(x)\in k[x]$ y tal que $f(x)$ 
m\'onico y libre de $\ell^n$ potencias. Sea $f(x)=(x-a_1)^{\alpha_1}
\cdots (x-a_r)^{\alpha_r}$ con $1\leq \alpha_i\leq \ell^n-1$, $a_i
\neq a_j$ para $i\neq j$ y $a_i\in k$. Sea $\alpha_i=\ell^{m_i}
c_i$ con $\mcd(\ell,c_i)=1$. Por el Teorema \ref{TRam1} se tiene
que el \'indice de ramificaci\'on de primo asociado a
$x-a_i$, ${\mc P}_i$, es igual a $\frac{\ell^n}{\mcd(\alpha_i,
\ell^n)}$. Si $m_i\geq 1$ para toda $i$, entonces si $m_0=
\min_{1\leq i\leq r}\{m_i\}\geq 1$. Entonces $L=K(z)$ con
$z^{\ell^{n-m_0}}=(x-a_1)^{\ell^{m_1-m_0}c_1}\cdots
(x-a_r)^{\ell^{m_r-m_0}c_r}$ lo que contradice que $L/K$
es una extensi\'on de grado $\ell^n$.

Sea $m_1=0$. Entonces $\mcd(\alpha_1,\ell)=1$. Entonces
${\mc P}_1$ es totalmente ramificado. El divisor de $f(x)$
es $(f(x))_K={\mc P}_1^{\alpha_1}\cdots {\mc P}_r^{\alpha_r}
\p^{\alpha_{r+1}}$ con $\sum_{i=1}^{r+1}\alpha_i=0$. Si
$\alpha_{r+1}=\ell^{m_{r+1}}c_{r+1}$ con $\mcd(\ell,c_{r+1})=1$,
se tiene $\sum_{i=1}^{r+1}\ell^{m_i}c_i=0$. Si $m_j\geq 1$ para
$2\leq j\leq r+1$, entonces, debido a que $m_1=0$, tendremos
que $\ell|c_1=-\sum_{i=2}^{r+1}\ell^{m_i} c_i$ lo cual es absurdo.

Por tanto, existe alg\'un $j$ con $2\leq j\leq r+1$ tal que
$m_j=0$. Poniendo $\p={\mc P}_{r+1}$, alguno de ${\mc P}_2,
\ldots {\mc P}_r, {\mc P}_{r+1}$ es totalmente ramificado en $L/K$.
\end{observacion}

\begin{proposicion}\label{P10.4.1.Ram}
Sea $E/F$ una extensi\'on finita de Galois de campos de
funciones congruentes. Sea $K/F$ una extensi\'on
arbitraria finita y sea $L=EK$. Sea $\pK_L$ un divisor
primo en $L$ y sean $\pK_E, \pK_F$ y $\pK_K$
sus restricciones a $E, F$ y $K$ respectivamente.
Sea $\rest\colon \Gal(L/K)\lra \Gal(E/F)$ el monomorfismo
de restricci\'on. Sean
$D$ e $I$ los grupos de descomposici\'on
y de inercia respectivamente. Entonces
\las
\item $I(\pK_L|\pK_K)|_E \subseteq I(\pK_E|\pK_F)$,
\item $D(\pK_L|\pK_K)|_E \subseteq D(\pK_E|\pK_F)$.
\end{list}
\end{proposicion}

\begin{proof}
\[
\xymatrix{
E\ar@{-}[r]\ar@{-}[d]&L=EK\ar@{-}[d]\\
E\cap K\ar@{-}[r]\ar@{-}[d]&K\ar@{-}[dl]\\
F}
\]

Primero consideremos $\sigma\in I(\pK_L|\pK_K)$. Entonces
$\sigma x\equiv x\bmod \pK_L$ para toda $x\in {\mc O}_{\pK_L}$.
Si $y\in {\mc O}_{\pK_E}$, entonces $y\in {\mc O}_{\pK_L}$
y por tanto $\sigma y\equiv y \bmod ({\pK_L}\cap E)$.
Entonces (1) es consecuencia de lo anterior y de que
${\pK_L}\cap E={\pK_E}$.

Ahora consideremos $\tau\in D(\pK_L|\pK_K)$. Entonces
$\pK_L^{\tau}=\pK_L$. Por lo tanto $\pK_E^{\tau}=(\pK_L
\cap E)^{\tau}=\pK_L^{\tau}\cap E^{\tau}=\pK_L\cap E=
\pK_E$ por lo que $\tau|_E\in D(\pK_E|\pK_F)$,
probando (2). $\fin$
\end{proof}

La Proposici\'on \ref{P10.4.1.Ram}
tiene varias consecuencias acerca de la
descomposici\'on de primos al trasladar las extensiones.

\begin{corolario}\label{C.10.4.2.Ram}
Con las hip\'otesis de la Proposici\'on {\rm{\ref{P10.4.1.Ram}}},
se tiene:
\las
\item Si $\pK_F$ se descompone totalmente en $E/F$,
entonces $\pK_K$ se descompone totalmente en $L/K$.

\item Si $\pK_F$ es no ramificado en $E/F$, entonces 
$\pK_K$ es no ramificado en $L/K$. Equivalentemente
si $\pK_K$ es ramificado en $L/K$, entonces $\pK_F$
es ramificado en $E/F$.
\end{list}
\end{corolario}

\begin{proof}
(1) se sigue del hecho de que $D(\pK_E|\pK_F)=\{1\}$.

(2) se sigue del hecho de que $I(\pK_L|\pK_K)\neq \{1\}$.
$\fin$
\end{proof}

\begin{observacion}\label{O10.4.3.Ram}
\las
\item Si $\pK_F$ es ramificado en $E/F$, $\pK_K$
no necesariamente es ramificado en $L/K$. Por ejemplo, sean
$F=\F(T)$, $E=F(y)$ con $y^p-y=\frac{1}{T}$ y
$K=F(z)$ con $z^p-z=\frac{1}{T}+\frac{1}{T+1}$. Entonces
si $\pK_F$ es el primo asociado a $T$, $\pK_F$
es ramificado en $E/F$ y en $K/F$
(Proposici\'on \ref{P2.4.Ram3}).
Ahora $L=EK=\F(T,y,z)=K(w)$ donde 
$w=z-y$ y $w^p-w=
\frac{1}{T+1}$, por lo que $\pK_K$ no es ramificado
en $L/K$.

\item Si $\pK_K$ es inerte en $L/K$, $\pK_F$ no 
necesariamente es inerte en $E/F$. Por ejemplo,
sea $q=p=3$, $\sqrt{-1}\notin {\ma F}_3$
(Proposici\'on \ref{P4.5.2}). Adem\'as
$P(T)=T^3-T-1$ es irreducible en ${\ma F}_3(T)$.
Sean $F={\ma F}_3(T)$, $E=F(\sqrt{P})$ y $K=F(
\sqrt{-P})$. Entonces $\pK_F=\p$ se ramifica en $E/F$
(Proposici\'on \ref{P5.1.2.Ram}). Por otro lado,
puesto que $L=EK={\ma F}_3(T,\sqrt{-P},\sqrt{P})=
K(\sqrt{-1})$ y como $\pK_F$ se ramifica en
$K/F$, $\pK_K$ es de grado $1$ y $L/K$ es una
extensi\'on de constantes de grado
$2$, por lo que $\pK_K$ es inerte en $L/K$.

Este mismo ejemplo muestra que ramificaci\'on en
$E/F$ puede transformarse en inercia en $L/K$.

\item Ramificaci\'on en $E/F$ puede transformarse
en cualquier cosa (ramificaci\'on, inercia o descomposici\'on)
en $L/K$. Por ejemplo sean $F={\ma F}_3(T)$, $E=F(y)$
con $y^3-y=T$ y $K=F(z)$ con $z^3-z=\frac{T^2+T-2}{T+1}$.
Entonces $\pK_F=\p$ es ramificado en $E/F$ y como
$L=K(w)$ con $w=y-z$ y $w^3-w=\frac{2}{T+1}$, $\pK_K$
se descompone en $L/K$.

Por otro lado si $E=F(u)$, $K=F(v)$ con $F={\ma F}_3(T)$,
$u^3-u=T$, $v^3-v=T^4$, entonces $\pK_F=\p$ es
totalmente ramificado en $L/F$ y en particular $\pK_K$
se ramifica en $L/K$.
\end{list}
\end{observacion}

\begin{observacion}\label{O1.5.10}
Si $E/F$ es una extensi\'on finita y separable, $\pK$ es un primo
de $F$ tal que $f_{E/F}({\mc P}|\pK)=1$ para todo lugar ${\mc P}$ de
$E$ sobre $\pK$, entonces no necesariamente $f_{E/F}
(\tilde {\mc P}|\pK)=1$ para $\tilde{\mc P}$ en $\tilde E$ sobre $\pK$.
\end{observacion}

\begin{ejemplos}\label{E1.5.11}{\ }
\las
\item Sean $F=\F(T)$ y $l$ un n\'umero primo no dividiendo a 
$q-1$. Sea $E=F(\sqrt[l]{P})$ para alg\'un $P\in R_T$ tal que
$l\nmid \deg P$. El primo infinito $\p$ de $F$
es ramificado en $E/F$ (Teorema \ref{TRam1}) por lo que 
el grado de inercia de $\p$ es 1 en $E$. Se tiene $\tilde E=
E(\zeta_l)={\ma F}_{q^r}(\sqrt[l]{P})$ con $\F(\zeta_l)=
{\ma F}_{q^r}$. Ahora bien $\zeta_l^{q^r-1}=1$ por lo que
$l|q^r-1$. Por tanto $r=o(q\bmod l)>1$. Se sigue que
$f_{{\ma F}_{q^r}(T)/\F(T)}({\mc P}_{\infty}|\p)=r>1$ donde
${\mc P}_{\infty}$ es el primo infinito de ${\ma F}_{q^r}(T)$.

\item Sean $F={\ma Q}$ y $E={\ma Q}(\sqrt[3]{2})$. Entonces
$\tilde E={\ma Q}(\sqrt[3]{2},\zeta_3)$. Se tiene que $2,3$ e
$\infty$ son los primos ramificados en $\tilde E/F$. Ahora bien,
$5$ es inerte en $\cic 3{}/{\ma Q}$ pues $5\equiv 2\bmod 3$ y
$o(5\bmod 3)=2$ (Teorema \ref{T8.2}). Puesto que el grupo
de inercia de $5$ en $\tilde E/F$ es c\'iclico, necesariamente
es de orden $2$. Se sigue que el grado de inercia de $5$ en
$E/F$ es 1 pero es $2$ en $\tilde E/F$.
\end{list}
\end{ejemplos}

\begin{proposicion}\label{P10.4.4.Ram}
Sea $L/K$ un extensi\'on c\'iclica de campos de
funciones congruentes de grado una potencia
de un primo, $[L:K]=l^n$. Sea
$L_0=K\subseteq L_1\subseteq L_2\subseteq \cdots
\subseteq L_{n-1}\subseteq L_n=L$ con $[L_i:L_{i-1}]=
l$, $1\leq i\leq n$. Sea $\pK$ un primo de $K$
tal que $\pK$ se ramifica en $L_i/L_{i-1}$. Entonces
$\pK$ es totalmente ramificado en $L/L_{i-1}$.
\end{proposicion}

\begin{proof}
Sea $G_0$ el grupo de inercia de $L/K$. Ya que $
G=\Gal(L/K)=\langle\sigma\rangle$, $o(\sigma)=l^n$, los
\'unicos subgrupos de $G$ son $\langle \sigma^{l^j}\rangle$,
$0\leq j\leq n$. Por tanto $G_0=\langle\sigma^{l^j}\rangle$
para alguna $j$. Ahora bien $L^{\langle\sigma^{l^j}\rangle}=
L_{n-j}$. Por tanto el campo de inercia de $\pK$ es
$L_{n-j}\subseteq L_{i-1}$ y $\pK$ es totalmente ramificado
en $L/L_{n-j}$. $\fin$
\end{proof}

Coloquialmente, podemos decir que en las extensiones del
tipo de la Proposici\'on \ref{P10.4.4.Ram}, cuando
un primo empieza a ramificarse, continua ramific\'andose
hasta $L$.

\begin{observacion}\label{O10.4.4.Ram+1}
En el caso de campos globales,
tanto num\'ericos como de funciones, los campos residuales
son campos finitos. Esto implica en particular que la inercia
de un primo en una extensi\'on de Galois, corresponde al
grupo de Galois de campos finitos y por lo tanto es c\'iclico.
M\'as precisamente, si $D$ denota al grupo de descomposici\'on
e $I$ al grupo de inercia, $D/I$ debe ser c\'iclico. En particular
no podemos tener primos totalmente inertes en extensiones
no c\'iclicas y en extensiones del tipo $C_p\times \cdots\times
C_p$, con $p$ un n\'umero primo, 
el grado de inercia de cualquier primo debe ser $1$ \'o $p$.

Esto no se cumple en general. Por ejemplo, si consideramos
$\K={\ma Q}(x)$ y $L=\K(\sqrt{x+3},\sqrt{x+5})$, entonces
$\Gal(L/\K)\cong C_2\times C_2$. Entonces, si $\pK$ es
es primo correspondiente a $x$, se tiene que $\pK$ es no
ramificado en $L/\K$ y tanto $T^2-(x+3)\bmod \pK=T^2-3\in
({\ma Q}[x]/\langle x\rangle)[T]\cong {\ma Q}[T]$ como
$T^2-(x+5)\bmod \pK=T^2-5$ son irreducibles. Por el
teorema de Kummer \cite[Theorem 5.8.2]{Vil2006},
$\pK$ es inerte, tanto en $L_1=\K(\sqrt{x+3})$
como en $L_2=\K(\sqrt{x+5})$, lo que implica
que $\pK$ es totalmente inerte en $L/\K$. De hecho, el
campo residual del primo en $L$ sobre $\pK$ es
${\ma Q}(\sqrt{3},\sqrt{5})$.
\end{observacion}

%% file: Capitulo11.tex
\chapter{Extensiones radicales de campos de funciones}\label{Ch9*}

\section{Introducci\'on}\label{S9*.1}

El contenido de este cap{\'\i}tulo est\'a basado en \cite{SanVil2013-1,
SanVil2013-2}.

Sea $L/\K$ una extensi\'on arbitraria de campos. Decimos que $L/\K$
es una {\em extensi\'on radical\index{extensi\'on radical}} si existe
$\alpha \in L$ tal que $L=\K(\alpha)$ y existe $n\in{\ma N}$ tal que
$\alpha^n=a\in \K$. En otras palabras $\alpha$ es una ra{\'\i}z
del polinomio $x^n-a\in \K[x]$. El elemento $\alpha$ usualmente
se representa como $\alpha=\sqrt[n]{a}$.

M\'as generalmente, una extensi\'on radical $L/\K$ es una extensi\'on
generada por ra{\'\i}ces de polinomios $x^{n_i}-a_i\in \K[x]$, esto es,
\[
L=\K\big(\sqrt[n_1]{a_1},\ldots, \sqrt[n_m]{a_m}\big)
\quad{\text{con}}\quad n_i\in {\ma N}, \quad a_i\in \K,\quad
1\leq i\leq m.
\]

En la teor{\'\i}a de las extensiones radicales se tienen de manera
natural dos grupos que permiten el estudio de tales extensiones.

\begin{definicion}\label{D9*.1.1} Sea $L/\K$ una extensi\'on
arbitraria de campos. Se definen
\l
\item {\em el grupo de torsi\'on\index{grupo de torsi\'on}} de $L/\K$ por
\[
T(L/\K):=\{\alpha\in L^{\ast}\mid {\text{existe\ }} n\in{\ma N} {\text{\ tal que\ }}
\alpha^n\in \K\},
\]
\item {\em el grupo de cogalois\index{grupo cogalois}} de $L/\K$ por
\[
\cog(L/\K):=\frac{T(L/\K)}{\K^{\ast}}.
\]
\end{list}
\end{definicion}

Para el estudio de las extensiones radicales, Greither y Harrison
\cite{GreHar86} desarrollaron una teor{\'\i}a que, en cierta forma,
es dual a la teor{\'\i}a de Galois, la cual ha sido generalizada en
innumerables direcciones.

\begin{definicion}[Ver \cite{GreHar86} y \cite{BaRzVi91}]
\label{D9*.1.2} Una extensi\'on finita de campos
$L/\K$ se llama
\l
\item {\em coseparable\index{extensi\'on coseparable}} si $L=
\K(T(L/\K))$,

\item {\em conormal\index{extensi\'on conormal}} si
$|\cog(L/\K)|\leq [L:\K]$,

\item {\em cogalois\index{extensi\'on cogalois}} si es conormal y
coseparable.
\end{list}
\end{definicion}

Notemos que coseparable y radical significan lo mismo.

\begin{definicion}\label{D9*.1.3} Una extensi\'on {\em $n$ de
Kummer\index{extensi\'on de Kummer}} es una extensi\'on $L/\K$
tal que 
\[
L=\K(\sqrt[n]{\Delta}),
\]
donde $n\in{\ma N}$ es primo relativo a la caracter{\'\i}stica de $\K$,
$\mu_n\subseteq \K$ donde $\mu_n$ es el grupo de las $n$--\'esimas
ra{\'\i}ces de la unidad y $\Delta$ es un subgrupo de $\K^{\ast}$
que contiene a $\K^{n\ast}$ y $\K(\sqrt[n]{\Delta})$
es el campo generado por todas las ra{\'\i}ces $\sqrt[n]{a}$
con $a\in \Delta$.
\end{definicion}

Este es el origen de la {\em Teor{\'\i}a de Kummer\index{teor{\'\i}a
de Kummer}}, la cual es una parte importante en el estudio de los
campos de clase. Se tienen los siguientes resultados.

\begin{teorema}\label{T9*.1.4} Sea $\K$ cualquier campo que contiene al grupo $\mu_n$
de las $n$--ra{\'\i}ces de la unidad donde $n$ es primo relativo
a la caracter{\'\i}stica de $\K$. Entonces
\l
\item Toda extensi\'on $n$ de Kummer $L/\K$ es una extensi\'on de 
Galois con grupo de Galois $\Gal(L/\K)$ abeliano
de exponente $n$.

\item Si $L/\K$ es una extensi\'on abeliana de exponente $n$,
entonces existe un subgrupo $\K^{n\ast}\subseteq \Delta \subseteq
\K^{\ast}$ tal que $L=\K(\sqrt[n]{\Delta})$.
\end{list}
\end{teorema}

\begin{proof} Ver Teoremas \ref{CClaseT1.1.4} y
\ref{CClaseT1.6.2}, as\'i como la
Subsecci\'on \ref{CClaseS1.6}. $\fin$
\end{proof}

La analog{\'\i}a existente entre campos num\'ericos y campos
de funciones congruentes, y de manera m\'as precisa, de campos
ciclot\'omicos con campos de funciones ciclot\'omicos, nos llevan
a la pregunta natural si existe lo an\'alogo de la torsi\'on usual
con la torsi\'on modular definida por la acci\'on de Carlitz-Hayes.

Daremos una nueva definici\'on de extensi\'on radical usando la
acci\'on de Carlitz--Hayes. As{\'\i}, una extensi\'on $L/\K$ ser\'a
llamada {\em radical\index{extensi\'on radical}\index{radical!acci\'on
de Carlitz--Hayes}} si $L$ puede
ser generada por algunos elementos $u$ con $u^{M_u}\in \K$
sobre $\K$, donde $M_u$ son polinomios en $R_T$. Entre
estas extensiones, estamos en especial interesados en las llamadas
extensiones {\em radicales ciclot\'omicas\index{radicales ciclot\'omicas}}.
Una extensi\'on se llamar\'a radical ciclot\'omica si es radical,
separable y {\em pura\index{extensi\'on pura}}. Estas extensiones
pueden ser vistas como generalizaciones naturales de las
extensiones de Carlitz--Kummer. Estas extensiones tiene
propiedades an\'alogas a las extensiones cogalois definidas
en \cite{GreHar86}, ver las Secciones
\ref{prop_extensiones_radicales} y \ref{rad_ciclotomicas_props}. 
Una extensi\'on radical ciclot\'omica $L/\K$ satisface
que $L=\K(T(L/\K))$. Notemos la analog{\'\i}a
con la definici\'on previa.

En este cap{\'\i}tulo estudiaremos la torsi\'on dada por la
acci\'on de Carlitz--Hayes. Por tanto entendemos por
``radical'' en el sentido de esta acci\'on. Estudiaremos
la estructura de campos de funciones congruentes generadas
por torsi\'on. En la Secci\'on \ref{radicales_ciclotomicas} 
definimos el concepto de extensi\'on radical ciclot\'omica como
un an\'alogo natural de las extensiones cogalois cl\'asicas.
Daremos ejemplos tanto de extensiones radicales como de
no radicales ciclot\'omicas as\'i como de extensiones puras y 
de extensiones no puras y mostraremos que, como en el
caso cl\'asico, la extensi\'on $K(\Lambda_{P^n})/K(\Lambda_P)$
es pura en donde $P\in R_T$ es un polinomio irreducible y
$n\in {\ma N}$. En las Secciones \ref{prop_extensiones_radicales} y
\ref{rad_ciclotomicas_props} daremos algunas propiedades
tanto de extensiones radicales como de extensiones radicales
ciclot\'omicas y probaremos que, como en el caso cl\'asico,
para extensiones de Galois, el grupo de cogalois es isomorfo
al grupo de los homomorfismos cruzados.

En la Secci\'on \ref{teoremas_de_estructura} 
obtendremos los resultados principales del cap{\'\i}tulo:
caracterizaremos las extensiones radicales ciclot\'omicas
finitas. En particular, veremos que extensiones
radicales finitas son $p$--extensiones donde $p$
es la caracter{\'\i}stica del campo base. Esto se probar\'a
en los Teoremas 
\ref{tdim_p_galois} y \ref{tdim_p} 
y en el Corolario \ref{cogalois_grado_p_pureza}.
En la Secci\'on \ref{examples}, daremos ejemplos y aplicaciones
de estos resultados. Finalmente, en la Secci\'on
\ref{estimacion_para_cogalois} hallaremos una cota 
superior para la cardinalidad del grupo de cogalois
de una extensi\'on radical ciclot\'omica finita.

Durante este cap{\'\i}tulo, usaremos la siguiente notaci\'on.

$p$ denota a un n\'umero primo.

$q = p^{\nu}$, $\nu \in \mathbb{N}$.

$K= {\mathbb{F}}_{q}(T)$ denota al campo de funciones racionales.

$R_{T}={\mathbb{F}}_{q}[T]$.

$\mu(\K)$ denota al conjunto de las ra{\'\i}ces de
Carlitz contenidas en un campo $\K$.

$\overline{K}$ denota a una cerradura algebraica de $K$.

$\car(L)$ denota a la caracter{\'\i}stica de un campo $L$.

Si $E/L$ es una extensi\'on de campos tal que
$K\subseteq L\subseteq
E\subseteq \overline{K}$, denotamos por
$T(E/L)$ al conjunto $\{u\in
E\mid \text{existe $M\in R_{T}$ tal que $u^{M}\in L$}\}$.

$C_{m}$ denota al grupo c{\'\i}clico de orden $m$.

\section{Extensiones de Kummer de campos de funciones}\label{S9*.2}

En esta secci\'on presentaremos una generalizaci\'on de extensiones de
Kummer por medio de la acci\'on de Carlitz--Hayes. En lo que resta
en este cap{\'\i}tulo $p$ siempre denotar\'a un n\'umero primo y $q=p^{\nu}$
donde $\nu\in{\ma N}$.
Denotaremos $k={\ma F}_q$, $K=k(T)$ y $R_T=k[T]$. Llamaremos
$\Lambda_M$, $M\in R_T\setminus\{0\}$, las {\em $M$--ra{\'\i}ces
de Carlitz\index{M--raices de Carlitz@$M$--ra{\'\i}ces de Carlitz}} y si $\lambda_M$ es
generador de $\Lambda_M$, $\lambda_M$ se llamar\'a
{\em ra{\'\i}z primitiva de Carlitz
\index{raiz primitiva de Carlitz@ra{\'\i}z primitiva de Carlitz}\index{Carlitz!ra{\'\i}z
primitiva de $\sim$}}.

Notemos que si $a\in\overline{K}$, entonces el conjunto de todas
las ra{\'\i}ces del polinomio $z^M-a\in\overline{K}[z]$ es el conjunto
$\{\alpha+\lambda\mid\lambda\in\Lambda_{M}\}$ donde 
$\alpha$ es cualquier ra{\'\i}z fija $z^{M}-a$ en $\overline{k}$.

Necesitaremos varios resultados de teor{\'\i}a de m\'odulos en
esta secci\'on.

La Proposici\'on \ref{P9*.2.2} y el Teorema
\ref{T9*.2.4} son an\'alogos a ({\sc i}) y ({\sc ii}) del Teorema
\ref {T9*.1.4} con la salvedad que consideraremos \'unicamente
extensiones finitas.

\subsection{Algo sobre la teor{\'\i}a de m\'odulos}\label{S9*.1.2.1}

En esta subsecci\'on, a menos que se indique lo contrario, todos
los m\'odulos y homomorfismos considerados son $R_T$--m\'odulos
y $R_T$--homomorfismos respectivamente.

Sea $A$ un m\'odulo, $a\in A$. Se define el homomorfismo
\[
\varphi_{a}\colon R_{T}\rightarrow A, \quad\text{definido por}\quad
\varphi_a(M):=Ma.
\]

\begin{definicion}\label{D9*.1.2.1}
Decimos que $A$ en un {\em m\'odulo c{\'\i}clico\index{m\'odulo
c{\'\i}clico}\index{m\'odulo!c{\'\i}clico}} si existe $a\in A$ tal que
$\varphi_a$ es un epimorfismo.
\end{definicion}

Notemos que la Definici\'on \ref{D9*.1.2.1} es equivalente a decir
que existe $a\in A$ tal que $A=(a)= R_T a$.

Para $a\in A$ consideremos el n\'ucleo del homomorfismo
$\varphi_{a}$, $\ker(\varphi_{a})$. Si $\ker(\varphi_{a})\neq\{0\}$
existe un polinomio no cero $M$, al cual lo podemos suponer
sin p\'erdida de generalidad m\'onico, tal que
$\ker(\varphi_{a})=(M)$.

\begin{definicion}\label{D9*.1.2.2}
Sea $a\in A$. Decimos que $a$ tiene {\em orden infinito} si el n\'ucleo
de $\varphi_{a}$ es cero. Decimos que $a$ tiene {\em orden
finito} si existe un polinomio m\'onico $M\in R_T\setminus \{0\}$
tal que $(M)= \ker(\varphi_{a})=(M)$. 
Si $A$ es un m\'odulo, {\em un exponente\index{exponente de un
m\'odulo}} de $A$ 
es un elemento no cero $M\in R_{T}$, tal que $Ma=0$
para todo $a\in A$.
\end{definicion}

\begin{observacion}\label{O9*.2.1.1}
Si $A$ es un m\'odulo finito, existe $a\in A$ tal que 
$\varphi_{a}$ es un epimorfismo de tal forma que existe
$M\in R_{T}\setminus\{0\}$ tal que
\l
\item $\ker(\varphi_{a})=(M)$ y
\item $R_{T}/(M)\cong A.$
\end{list}

Como antes, podemos reemplazar a $M$ por un polinomio m\'onico
y entonces diremos que $A$ tiene {\em orden $M$\index{orden
de un m\'odulo}}.
\end{observacion}

La demostraci\'on del siguiente lema es directa y no la presentamos.

\begin{lema}\label{L9*.2.1.2}
Sea $A$ un m\'odulo c{\'\i}clico de orden $N$ con $N\neq 0$. Sea
$N_{1}$ un divisor m\'onico de $N$. Entonces existe un subm\'odulo
de $A$ de orden $N_{1}$. $\fin$
\end{lema}

\begin{observacion}\label{O9*.2.1.3}
Con las condiciones del Lema \ref{L9*.2.1.2}, se sigue de la
Observaci\'on \ref{O9*.2.1.1} que $Na=0$ para toda $a\in A$.

Por otro lado, si $B$ es un subm\'odulo c{\'\i}clico de $A$ de
orden $N_1$, entonces nuevamente de la Observaci\'on
\ref{O9*.2.1.1} obtenemos que existe $b\in
B$, tal que $\varphi_{b}$ es un epimorfismo y
$\ker(\varphi_{b})=(N_{1})$. Puesto que $b\in A$, existe $N_{2}\in
R_{T}$, tal que $b=N_{2}a$.

Ahora, puesto que $Nb = N(N_{2}a)=N_{2}(Na)=0$ tenemos que
$N\in\ker(\varphi_{b})= (N_{1})$. Por tanto $N_{1}$ es un divisor de
$N$. Puesto que todos los m\'odulos c{\'\i}clicos de orden $M$ son isomorfos
a $R_{T}/(M)$, esto es, esencialmente \'unicos, se sigue que para cada
divisor m\'onico $M$ de $N$, existe un \'unico subm\'odulo
c{\'\i}clico de orden $M$ de $A$.
\end{observacion}

Sea $A$ un m\'odulo c{\'\i}clico con exponente $M$. Denotamos
por $C_{M}$ al m\'odulo $R_{T}/(M)$ el cual es c{\'\i}clico de orden $M$.
Esto en analog{\'\i}a a la notaci\'on $C_{m}$ de los grupos c{\'\i}clicos.

\begin{definicion}\label{D9*.2.1.4}
Se denota por $\hat{A}$ o por $\Hom_{R_{T}}(A,C_{M})$ al grupo de
homomorfismos de $A$ en $C_{M}$, donde $A$ es de exponente $M$.
Este m\'odulo
se llama el {\em m\'odulo dual\index{m\'odulo dual}} de $A$.
\end{definicion}

Supongamos que $f\colon A\rightarrow B$ es un homomorfismo, y que 
tanto $A$ como $B$ tienen exponente $M$. Entonces se tiene un
homomorfismo $\widehat{f}\colon \widehat{B}\rightarrow \widehat{A}$
definido por $\widehat{f}(\psi) = \psi\circ f$. Notemos que
$\widehat{(\ )}$ es un funtor contravariante, es decir
\[
\widehat{(\ )}:R_{T}\text{--m\'odulos de exponente $M$}\longrightarrow R_{T}
\text{--m\'odulos}
\]
es tal que si $f:A\rightarrow B$ y $g:B\rightarrow C$ son
homomorfismos, entonces

(1) $\widehat{g\circ f} = \widehat{f}\circ\widehat{g}$ y

(2) $\widehat{1}= 1$.

\begin{lema}\label{L9*.2.1.5}
Si $A$ es un m\'odulo finito de exponente $M$, tal que $A =
B\times C$, entonces $\widehat{A}$ es isomorfo a $\widehat{B}\times
\widehat{C}$.
\end{lema}

\begin{proof}
Las proyecciones naturales $\pi_{1}:B\times C\rightarrow B$ y
$\pi_{2}:B\times C\rightarrow C$, inducen los homomorfismos
$\widehat{\pi}_{1}:\widehat{B}\rightarrow \widehat{B\times C}$ y
$\widehat{\pi}_{2}:\widehat{C}\rightarrow \widehat{B\times C}$, por
lo que podemos definir $\theta\colon \widehat{B}\times \widehat{C}\rightarrow
\widehat{B\times C}$ dada por $\theta(\psi_{1},\psi_{2}) =
\widehat{\pi_{1}}(\psi_{1}) + \widehat{\pi_{2}}(\psi_{2})$, donde
$(\psi_{1},\psi_{2})\in \widehat{B}\times \widehat{C}$.

Se tiene que $\theta$ es un homomorfismo. Por otra parte si $\psi\in
\widehat{B\times C}$ entonces, puesto que $\psi$ es un homomorfismo,
$\psi(x,y) = \psi(x,0) + \psi(0,y)$ para todo $(x,y)\in B\times C$.
Ahora se define $\psi_{1}:B\rightarrow A_{M}$ por $\psi_{1}(x) =
\psi(x,0)$ y $\psi_{2}:C\rightarrow A_{M}$ dado por $\psi_{2}(y) =
\psi(0,y)$. Entonces $\psi_{1}$ y $\psi_{2}$ son homomorfismos. 
De esta forma se induce una funci\'on
\[
\delta\colon\widehat{B\times C}\rightarrow \widehat{B}\times
\widehat{C}
\]
dada por $\delta(\psi)=(\psi_{1},\psi_{2})$ el cual es un
homomorfismo de m\'odulos y cuya inversa es $\theta$. El resultado
se sigue. $\fin$
\end{proof}

\begin{proposicion} \label{P9*.2.1.6}
Un m\'odulo finito $A$ es isomorfo a su dual. 
Esto es
\[
A\cong\widehat{A}= \Hom_{R_{T}}(A,C_{M}),
\]
donde $A$ tiene exponente $M$.
\end{proposicion}

\begin{proof}
Tenemos que se puede escribir
$A\cong \oplus_{P} A_{P}$ (ver Teorema \ref{DrinfeldT1.3.3}). 
La suma anterior es sobre todos los
polinomios m\'onicos irreducibles $P$ y $A_{P}$ denota los elementos de $A$ que
tiene orden una potencia de $P$.

Ahora por el Teorema 4.9, Cap{\'\i}tulo 5 de \cite{HilWu82} se tiene que
$A_{P}$ se puede escribir como $A_{P}\cong C_{P^{\alpha_{1}}}\oplus
\cdots \oplus C_{P^{\alpha_{k}}}$, donde $\alpha_{1}\geq \cdots \geq
\alpha_{k}\geq 1$ y cada $C_{P^{\alpha_{i}}}$ es un m\'odulo
c{\'\i}clico cuyo generador tiene orden $P^{\alpha_{i}}$. De esta forma, cada
$C_{P^{\alpha_{i}}}$ tiene orden $P^{\alpha_{i}}$. N\'otese que cada
$A_{P}$ y cada $C_{P^{\alpha_{i}}}$ tiene exponente $M$.

Por la observaci\'on anterior y el Lema
\ref{L9*.2.1.5}, podemos suponer que $A$ es
c{\'\i}clico generado por $a$ de orden $P^{\alpha}$, con $\alpha\in
{\mathbb{N}}$ y $P\in R_{T}$ irreducible. Por lo tanto la funci\'on
$\varphi_{a}$ es un epimorfismo y $(P^{\alpha})= \ker(\varphi_{a})$.

Puesto que $M$ es de exponente de $A$, se tiene que $P^{\alpha}| M$. Ahora
del Lema \ref{L9*.2.1.2}, junto con la Observaci\'on
\ref{O9*.2.1.1}, se tiene $C_{M}$ tiene un \'unico subm\'odulo
c{\'\i}clico de orden $P^{\alpha}$, que denotamos por $C_{P^{\alpha}}$.
El homomorfismo $\varphi_{a}\colon R_{T}\rightarrow A$ induce un
isomorfismo, que seguiremos denotando por
$\varphi_{a}\colon R_{T}/(P^{\alpha})\rightarrow A$.

Al inverso del isomorfismo $\varphi_{a}$, lo denotaremos por $\psi$.
Sea $y = \psi(a)$, entonces $y$ es un generador de $C_{P^{\alpha}}$.
Al componer $\psi$ con la inclusi\'on natural
$C_{P^{\alpha}}\hookrightarrow C_{M}$, se obtiene un elemento de
$\widehat{A}$, que seguiremos denotando por $\psi$.

Ahora sea $\varphi\in \widehat{A}$. Notemos que
$\im(\varphi)\subseteq C_{M}$ es un subm\'odulo c{\'\i}clico de orden
$N$. As{\'\i}, existe $w\in \im(\varphi)$, $w=\varphi(a_{w})$ con
$a_{w}\in A$, que genera a $\im(\varphi)$ y su orden es $N$.

Por otra parte se tiene que $P^{\alpha}\in (N)$. Por lo tanto
$P^{\alpha}=ND$, para alg\'un $D\in R_{T}$. Por lo que
$N=P^{\gamma}$ para alg\'un $\gamma\leq \alpha$, es decir, $w$ tiene
orden $P^{\gamma}$ y como $a_w$ genera a $\im(\varphi)$, se tiene que
$\im(\varphi)\subseteq C_{P^{\alpha}}$.

Por otro lado $\varphi$ est\'a determinado completamente por su
acci\'on en $a$, donde $a\in A$ es un generador de $A$. Por lo tanto
$\varphi(a) = Ny$. Ahora si $\psi_{N}=N\psi$ entonces
$\psi_{N}(a)=N\psi(a)=Ny=\varphi(a)$, es decir, $\varphi=\psi_{N}\in
(\psi)$.

De esta forma se tiene $\widehat{A}=(\psi)$ el cual es de 
orden $q^{\deg (P^{\alpha})}$. Por lo tanto $A\cong \widehat{A}$. $\fin$
\end{proof}

\begin{definicion}\label{D9*.2.1.7} Sean $A$ y $B$ m\'odulos.
Una {\em funci\'on bilineal\index{funci\'on bilineal}}
de $A\times B$ en un m\'odulo $C$ es
una funci\'on $A\times B\rightarrow C$, denotada por $(a,b)\mapsto
\langle a,b\rangle$, que tiene la propiedad siguiente: para cada $a\in A$, la
funci\'on $b\mapsto \langle a,b \rangle$ es un homomorfismo y, para cada $b\in B$,
la funci\'on $a\mapsto \langle a,b\rangle$ es un homomorfismo. Un
elemento $a\in A$ se
dice {\em ortogonal\index{elemento ortogonal}}
a $S\subseteq B$ si $\langle a,b\rangle = 0$ para cada $b\in S$.

De modo an\'alogo tenemos la definici\'on de que $b\in B$ sea
ortogonal a $S\subseteq A$, esto es, si $\langle a,b\rangle=0$
para toda $a\in A$. El {\em n\'ucleo izquierdo
\index{nucleo izquierdo@n\'ucleo izquierdo}} de la
funci\'on bilineal es el subm\'odulo de $A$, que denotamos por
$N_{I}$, ortogonal a $B$.

El {\em n\'ucleo derecho\index{nucleo derecho@n\'ucleo derecho}}
de la funci\'on bilineal es el subm\'odulo
de $B$, que denotamos por $N_{D}$, ortogonal a $A$.
\end{definicion}

Un elemento $b\in B$ da lugar a un elemento de $\Hom_{R_{T}}(A,C)$, dado por
$a\mapsto \langle a,b\rangle$, que denotamos por $\psi_{b}$. Entonces $\psi_{b}$
se anula en $N_{I}$, es decir, $\psi_{b}(a) = 0$ para cada $a\in
N_{I}$. As{\'\i}, $\psi_{b}$ induce un homomorfismo $A/N_{I}\rightarrow
C$, dado por $a+ N_{I}\mapsto\psi_{b}(a)$.

Por otro lado, si $b\equiv b^{\prime} \modulo N_{D}$ entonces
$\psi_{b} = \psi_{b^{\prime}}$, esto da lugar, en primer t\'ermino, a un
homomorfismo $\psi:B/N_{D}\rightarrow \Hom_{R_{T}}(A/N_{I},C)$ dado
por $\psi(b+N_{D}) = \psi_{b}$, y, en segundo t\'ermino, a la sucesi\'on
exacta de m\'odulos
\begin{equation}\label{ec1_3}
0\rightarrow B/N_{D}\rightarrow \Hom_{R_{T}}(A/N_{I},C).
\end{equation}

De modo similar se obtiene
\begin{equation}\label{ec2_3}
0\rightarrow A/N_{I}\rightarrow \Hom_{R_{T}}(B/N_{D},C).
\end{equation}

\begin{proposicion}\label{P9*2.1.8}
Sea $A\times A^{\prime}\rightarrow C$ una funci\'on bilineal de
m\'odulos, con $C$ un m\'odulo c{\'\i}clico finito de orden $M$. Sean $B$ y
$B^{\prime}$ los n\'ucleos izquierdos y derecho, respectivamente.
Supongamos que $A^{\prime}/B^{\prime}$ es finito. Entonces $A/B$ es
finito y $A^{\prime}/B^{\prime}$ es isomorfo al m\'odulo dual de
$A/B$.
\end{proposicion}

\begin{proof}
De las sucesiones exactas (\ref{ec1_3}) y (\ref{ec2_3}), se deduce
que las sucesiones siguientes son exactas
\begin{equation}
0\rightarrow A^{\prime}/B^{\prime}\rightarrow
\Hom_{R_{T}}(A/B,C)\label{ec1_bilineal}
\end{equation}
\noindent y
\begin{equation}
0\rightarrow A/B\rightarrow
\Hom_{R_{T}}(A^{\prime}/B^{\prime},C)\label{ec2_bilineal}.
\end{equation}

De (\ref{ec2_bilineal}) deducimos que $A/B$ puede ser visto como un
subm\'odulo de $\Hom_{R_{T}}(A^{\prime}/B^{\prime}, C)$, de aqu{\'\i}
la finitud de $A/B$. Por otro lado se tienen las desigualdades, que
se infieren de las sucesiones (\ref{ec1_bilineal}) y (\ref{ec2_bilineal}) y de la
Proposici\'on \ref{P9*.2.1.6}:
\begin{gather*}
\card (A/B) \leq \card (\widehat{A^{\prime}/B^{\prime}}) = \card
(A^{\prime}/B^{\prime})\\
\intertext{y}
\card(A^{\prime}/B^{\prime}) \leq \card(\widehat{A/B}) =
\card(A/B).
\end{gather*}
La segunda igualdad se debe a la Proposici\'on
\ref{P9*.2.1.6}. De esto se deduce la suprayectividad de
la sucesi\'on exacta (\ref{ec1_bilineal}), y de esto se
sigue el resultado. $\fin$
\end{proof}

\section{Teor{\'\i}a de Kummer}\label{S9*.3}

En esta secci\'on se dar\'a una generalizaci\'on de las extensiones
de Kummer, un poco diferente a las dadas por Chi y Li en
\cite{CheLi2001} y por Schultheis en \cite{Sch90}. En lo que
sigue supondremos que las extensiones a considerar son
subextensiones de $\overline{K}/K$. Siguiendo a \cite{Lan93}, sean
$M\in R_{T}$ un polinomio no constante y 
$\varphi\colon \K\rightarrow \K$
definido por $\varphi(u) = u^{M}$, donde $\K=K(\Lambda_{M})$.
Entonces $\varphi$ es un $R_{T}$--homomorfismo. Por otra parte
consideremos un $R_{T}$--subm\'odulo $B$ de $\K$, bajo la acci\'on de
Carlitz Hayes, que contenga a $\K^{M} = \varphi(\K)$.

Denotamos por $\K_{B}$ la composici\'on de todos los campos
$\K(\sqrt[M]{a})$ con $a\in B$. Esto \'ultimo quiere decir que
adjuntamos a $\K$ una ra{\'\i}z arbitraria $\alpha$ de la ecuaci\'on
$z^{M} - a = 0$, donde $\alpha\in\overline{K}$. Puesto que las
$M$-ra{\'\i}ces de Carlitz est\'an en $\K$, tal campo no depende de la
elecci\'on de la ra{\'\i}z $\alpha$, y por lo tanto $\K_{B}$ es de
Galois sobre $\K$.

\begin{definicion}\label{D9*.2.1}
Diremos que una extensi\'on abeliana $L/\K$, con grupo
de Galois $G$, es una
extensi\'on {\em $R_{T}$--abeliana\index{extensi\'on $R_T$--abeliana}}
si $G$ tiene estructura de
$R_{T}$-m\'odulo; una extensi\'on $R_{T}$-abeliana $L/\K$ se dice que
tiene {\em exponente\index{exponente de una extensi\'on $R_T$--abeliana}}
$M\in R_{T}$ si $M\cdot\sigma = 1$ para cada
$\sigma\in G$ (ver \cite{CheLi2001}).
\end{definicion}

\begin{proposicion}\label{P9*.2.2}
\l
\item Sea $B$ un $R_{T}$-m\'odulo de $\K$
que contiene a $\K^{M}$ y sea
$\K_{B}$ la composici\'on de todos los campos $\K(\sqrt[M]{a})$, para
cada $a\in B$. Entonces $\K_{B}/\K$ es Galois y abeliana.

\item Supongamos
que $\K_{B}/\K$ es una extensi\'on $R_{T}$-abeliana y de exponente
$M$. Entonces existe una funci\'on bilineal:
\[
G\times B\rightarrow \Lambda_{M}\quad\text{dada por}\quad
(\sigma,a)\mapsto \langle\sigma,a\rangle
\]
donde $\langle\sigma,a\rangle = \sigma(\alpha) - \alpha$ y $\alpha$
satisface $\alpha^{M} = a$. El n\'ucleo izquierdo es $1$ y el n\'ucleo
derecho es $\K^{M}$.

La extensi\'on
$\K_{B}/\K$ es finita si y s\'olo
si $(B:\K^{M})$ es finito. Si esto ocurre, entonces
\[
B/\K^{M} \cong \widehat{G}.
\]

En particular se tiene que
\[
[\K_{B}:\K] = (B:\K^{M}).
\]
\end{list}
\end{proposicion}

\begin{proof}
(I) Sea $b\in B$ y sea $\beta$ una $M$--ra{\'\i}z de $b$. El polinomio $z^M-b$
se descompone en factores lineales en $\K_B$ para todo $b\in B$.
Entonces $\K_B/\K$ es una extensi\'on de Galois.
Sean $G = \Gal(\K_{B}/\K)$, $\sigma\in G$, $b\in B$ y $\beta$ una
ra{\'\i}z del polinomio $z^{M} - b$. Entonces $\sigma(\beta) = \beta
+ \lambda^{M_{\sigma}}$, para alg\'un $M_{\sigma}\in R_{T}$, donde
$\lambda$ es un generador de $\Lambda_{M}$, por lo que se tiene un
monomorfismo $G\rightarrow \Lambda_{M}$, $\sigma\mapsto
\lambda^{M_{\sigma}}$ de donde se sigue que $G$ es
un grupo abeliano.

(II) Definimos $G\times B\rightarrow \Lambda_{M}$ por
$(\sigma,b)\mapsto \langle\sigma,b\rangle$, donde 
$\langle\sigma,b\rangle = \sigma(\beta) -
\beta$ y $\beta^{M}=b$. Esta definici\'on es independiente
de la elecci\'on de la $M$--ra{\'\i}z de $b$.
Se tiene que $\langle\sigma,a + b\rangle
= \langle\sigma,a\rangle
+ \langle\sigma,b\rangle$ para cada $a,b\in B$ y puesto que $(\sigma(\beta) -
\beta)\in \Lambda_{M}$, se sigue que $\langle\sigma\cdot\tau,b\rangle=
\langle\sigma,b\rangle + \langle\tau,b\rangle$.

Sea $\sigma\in G$ y supongamos que $\langle\sigma,b\rangle = 0$ para cada $b\in
B$. Por lo tanto, si $\beta$ satisface que $\beta^{M}= b$ se tiene
que $\sigma(\beta) = \beta$, y como esto vale para cada generador
se tiene que $\sigma = 1$, es decir el n\'ucleo izquierdo es $1$.

Por otro lado si $b\in B$ satisface que $\langle\sigma,b\rangle = 0$ para todo
$\sigma\in G$ entonces $\sigma(\beta) = \beta$ para toda
$\sigma\in G$. Por lo tanto, $\beta\in \K$ y $b = \beta^{M}\in
\K^{M}$. De aqu{\'\i} se sigue que el n\'ucleo derecho es $\K^{M}$.

Ahora supongamos que $B/\K^{M}$ es finito. 
Entonces $G/1 = G$ es finito. En particular
$\K_{B}/\K$ es finito. 
Ahora bien, si $\K_{B}/\K$ es finito, puesto que el
n\'ucleo derecho es $\K^M$, de la Proposici\'on
\ref{P9*2.1.8} se obtiene que la siguiente sucesi\'on es exacta
\[
0\rightarrow B/\K^{M}\rightarrow \Hom_{R_{T}}(G/1,\Lambda_{M}).
\]
De esta sucesi\'on y de que
$\Hom_{R_{T}}(G/1,\Lambda_{M})$ es finito, se sigue que
$(B:\K^{M})$ es finito.

Finalmente, puesto que por la
Proposici\'on \ref{P9*.2.1.6} $B/\K^M$ es isomorfo al m\'odulo
dual de $G$, se tiene que $B/\K^M\cong G$, as{\'\i} que
$[\K_B:\K]=(B:\K^M)$. $\fin$
\end{proof}

Antes de mostrar la proposici\'on siguiente necesitamos algunas
definiciones, dadas en \cite{CheLi2001}. 

\begin{definicion}\label{D9*.2.3}
Una extensi\'on $R_{T}$--abeliana $L/\K$ se dice que es {\em
$R_{T}$--c{\'\i}clica} si $\Gal(L/\K)$ es un $R_{T}$-m\'odulo c{\'\i}clico.
En este caso si $\Gal(L/\K)\cong R_{T}/(M)$, con $M$ un polinomio
m\'onico, diremos que $L/\K$ es una extensi\'on c{\'\i}clica de orden
$M$.
\end{definicion}

En el siguiente teorema denotamos por ${\mathfrak{M}}$ el conjunto
de $R_{T}$-subm\'odulos de $\K$, que contienen a $\K^{M}$ y
${\mathfrak{F}}$ denota el conjunto de extensiones $R_{T}$-abelianas
de $\K$ de exponente $M$.

\begin{teorema}\label{T9*.2.4}
Con las notaciones de la Proposici\'on
{\rm \ref{P9*.2.2}}, la funci\'on
$\varphi:{\mathfrak{M}}\rightarrow {\mathfrak{F}}$ dada por
$\varphi(B)=\K_{B}$ es inyectiva. Adem\'as si $L/\K$ es una
extensi\'on $R_{T}$-abeliana, finita, de exponente $M$ entonces
existe un $R_{T}$-subm\'odulo $B$, de $\K$, que contiene a $\K^{M}$,
tal que $L = \K_{B}$.
\end{teorema}

\begin{proof}
Para mostrar la inyectividad de la funci\'on anterior bastar\'a probar
que si $\K_{B_{1}}\subseteq \K_{B_{2}}$ entonces $B_{1} \subseteq
B_{2}$, puesto que de la igualdad $\varphi(B_{1})=\varphi(B_{2})$,
se deducen las contenciones $\K_{B_{1}}\subseteq \K_{B_{2}}$ y
$\K_{B_{2}}\subseteq \K_{B_{1}}$.

Sea $b\in B_{1}$. Se tiene que
$\K(\sqrt[M]{b})\subseteq \K_{B_{2}}$
por lo que $\K(\sqrt[M]{b})$ est\'a
contenido en una subextensi\'on finitamente generada de $\K_{B_{2}}$,
es decir, existen un n\'umero finito de elementos
$b_{i}\in B_{2}$ de modo que
$\K(\sqrt[M]{b})\subseteq \K(b_{1}, \ldots , b_{m})$. As{\'\i} podemos
suponer que $B_{2}/\K^{M}$ es finitamente generada y
por tanto es una extensi\'on finita.

Sea $\beta$ tal que $\beta^M=b$.
Sea $B_{3}$ el subm\'odulo de $\K$ generado por $B_{2}$ y $b$. Veamos
$\K_{B_{2}} = \K_{B_{3}}$. Tenemos que $\K_{B_{2}} \subseteq
\K_{B_{3}}$. Para mostrar la otra contenci\'on, sea $\alpha$ una
ra{\'\i}z $M$-\'esima de $c\in B_{3}$. Si $c\in B_{2}$ entonces
$\K(\alpha)\subseteq \K_{B_{2}}$. Si $c$ es de la forma $b^{N} + \sum
b^{N_{i}}_{i}$, con $b_{i}\in B_{2}$, entonces $\alpha^{M} = b^{N} +
\sum b^{N_{i}}_{i} = \beta^{MN} + \sum \beta^{MN_{i}}_{i}$, con
$\beta_i^M=b_i$, $i=1,\ldots, s$, es
decir, $\alpha = \beta^{N} + \sum \beta^{N_{i}}_{i} + \lambda^{A}$,
donde $\lambda$ es un generador de $\Lambda_M$.
Por lo tanto $\K(\alpha) \subseteq \K_{B_{2}}$. De aqu{\'\i}
se sigue que $\K_{B_{3}}\subseteq \K_{B_{2}}$.

Entonces, por la Proposici\'on \ref{P9*.2.2} (II) se tiene
$(B_{2}:\K^{M})=(B_{3}:\K^{M})$, de esta manera $b\in
B_{2}$, por lo que $B_{1}\subseteq B_{2}$.

Por otro lado, sea $\K^{\prime}$ una extensi\'on $R_{T}$--abeliana de
$\K$ de exponente $M$, finita. Sea $G=\Gal(\K^{\prime}/
\K)$. Entonces, por los Teoremas 4.7 y 4.9, Cap{\'i}tulo 4
de \cite{HilWu82}, $G$ es suma
directa, finita, de $R_{T}$-subm\'odulos de exponente $M$. Aplicando
Teor{\'\i}a de Galois podemos suponer que la extensi\'on es
c{\'\i}clica de exponente $M$. Ahora por la Proposici\'on 2.6 de
\cite{CheLi2001}, se tiene que toda extensi\'on c{\'\i}clica $
\K^{\prime}/\K$ de
exponente $M$, se obtiene adjuntando una $M$-ra{\'\i}z de un elemento
de $\K$.

As{\'\i} $\K^{\prime}$ es la adjunci\'on de $M$-ra{\'\i}ces, es decir,
existen $\{b_{j}\}\subseteq \K$ y $\{\alpha_{j}\}\subseteq
\K^{\prime}$ tales que $\alpha^{M}_{j} = b_{j}$ y $\K^{\prime} =
\K(\{\alpha_{j}\})$. Sea $B$ el subm\'odulo de $\K$ generado por
$\{b_{j}\}$ y $\K^{M}$. Entonces 
$\K^{\prime}\subseteq \K_{B}$. Por otro
lado consideremos una ra{\'\i}z $M$-\'esima de $c\in B$, digamos
$\alpha$. As{\'\i} $\alpha^{M} = c$. Se tiene que $
c=\sum_{j=1}^s b_j^{N_j} + a^M$, $a\in \K$. Entonces
$\alpha=\sum_{j=1}^s \alpha_j^{N_j}+a$ por lo que
$\K(\alpha)\subseteq 
\K^{\prime}$. Se sigue que
$\K_{B} \subseteq \K^{\prime}$ y $\varphi(B)
=\K^{\prime}$. Esto
termina la demostraci\'on. $\fin$
\end{proof}

\begin{proposicion}\label{P9*2.5}
Sea $L/\K$ una extensi\'on $R_{T}$--abeliana, finita, supongamos que
$\Lambda_{N}\subseteq \K$, con $N\in R_{T}$ no constante. Sea
\[
W = \{\overline{\alpha}= \alpha+\K^{N}\in \K/\K^{N}\mid \sqrt[N]{\alpha}\in L\}.
\]
Entonces $W\cong\Hom(G,\Lambda_{N})$, donde $G=\Gal(L/\K)$.
\end{proposicion}

\begin{proof}
Dado $\overline{a}\in W$, se define una funci\'on
$\varphi_{\overline{a}}:G\rightarrow\Lambda_{N}$ definida as{\'\i}
$\varphi_{\overline{a}}(\sigma) = \sigma(\alpha)-\alpha$, donde
$\alpha$ es una ra{\'\i}z $N$-\'esima de $a$;
$\varphi_{\overline{a}}$ es independiente de la ra{\'\i}z usada.
Notemos que
\begin{align*}
\varphi_{\overline{a}}(\sigma\circ\tau) &= \sigma(\tau(\alpha))-\alpha
= \sigma(\tau(\alpha)-\alpha+\alpha)-\alpha\\
&= \sigma(\tau(\alpha)-\alpha)+\sigma(\alpha)-\alpha
= \tau(\alpha)-\alpha+\sigma(\alpha)-\alpha.
\end{align*}
Por tanto $\varphi_{\overline{a}}$ es un
homomorfismo de grupos abelianos. Por lo tanto es posible definir
$f:W\rightarrow \Hom(G,\Lambda_{N})$ dado por $f(\overline{a})
= \varphi_{\overline{a}}$.

Se tiene que $f$ es un homomorfismo de grupos abelianos. Ahora si
$f(\overline{a})= \varphi_{\overline{a}} = 0$ entonces
$\sigma(\alpha)-\alpha = 0$, para cada $\sigma\in G$. De esta manera
se tiene que $\alpha\in \K$. Puesto que $a=\alpha^N$,
entonces $a\in \K$. De esta modo
$\overline{a}=0$, por lo tanto $f$ es inyectiva.

Ahora sea $\varphi:G\rightarrow\Lambda_{N}$ un homomorfismo de
grupos abelianos. Entonces
\[
\varphi(\sigma\circ\tau) = \varphi(\sigma)+\varphi(\tau) 
= \varphi(\sigma)+\sigma(\varphi(\tau)),
\]
es decir, $\varphi$ es un homomorfismo cruzado, por lo
tanto por el Teorema 90 de Hilbert aditivo, existe un $\alpha\in L$ tal
que $\varphi(\sigma) = \sigma(\alpha)-\alpha$. As{\'\i} tenemos
$(\sigma(\alpha)-\alpha)^{N}= \sigma(\alpha^{N})-\alpha^{N} = 0$, 
por lo que $a=\alpha^{N}\in
\K$, lo cual prueba la suprayectividad de $f$. $\fin$
\end{proof}

\section{Extensiones radicales
ciclot\'omicas.}\label{radicales_ciclotomicas}

Los siguientes resultados ser\'an \'utiles en esta secci\'on.

\begin{proposicion}\label{L10.10.1^{-1}}
Sea $\K/K$ una extensi\'on finita. Entonces $\mu(\K)
=\Lambda_M$ para alg\'un $M\in R_T$.
\end{proposicion}

\begin{proof}
Damos dos demostraciones. Para la primera
se tiene $\mu(\K)=\{u\in \K\mid \text{existe $M\in R_T$ 
tal que $u^M=0$}\}$. Entonces $\mu(\K)$ es un $R_T$--m\'odulo
pues si $z\in \mu(\K)$, consideremos $N\in R_T$ tal que $z^N=0$.
Sea $N^{\prime}\in R_T$ arbitrario. Se tiene $(z^{N^{\prime}})^N=z^{N^{\prime}
N}=z^{N N^{\prime}}=(z^N)^{N^{\prime}}=0^{N^{\prime}}=0$.

Sea $P\in R_T$ un polinomio m\'onico e irreducible y sea $\mu(\K)(P)=\{
u\in \K\mid \text{existe $n\in{\ma N}$ tal que $u^{P^n}=0$}\}$ la 
$P$--torsi\'on de $\mu(\K)$. Entonces $\mu(\K)(P)$
es un $R_T$--subm\'odulo de $\mu(\K)$ y se tiene
\[
\mu(\K)=\bigoplus_{\substack{P\in R_T\\ \text{$P$ m\'onico e
irreducible}}} \mu(\K)(P).
\]

Fijemos $P\in R_T$ m\'onico e irreducible tal que $\mu(\K)(P)\neq 0$
y consideremos $n\in{\ma N}$ el m{\'\i}nimo 
tal que $u^{P^n}=0$ para todo $u\in\mu(\K)(P)$ , esto es,
existe $z\in \mu(\K)(P)$ tal que $z^{P^{n-1}}\neq 0$.
Sea $u\in \mu(\K)(P)$. Entonces $u^{P^s}=0$ con $s\leq n$.
Por tanto $u^{P^n}=(u^{P^s})^{P^{n-s}}=0^{P^{n-s}}=0$. Se sigue
que $\mu(\K)(P)\subseteq \Lambda_{P^n}$. 

Por otro lado, existe $\lambda\in \mu(\K)(P)$ con $\lambda^{P^n}
=0$ y $\lambda^{P^{n-1}}\neq 0$. En particular $\lambda\in \Lambda_{P^n}
\setminus \Lambda_{P^{n-1}}$ por lo que $\lambda$ es un generador
de $\Lambda_{P^n}$ y se tiene que por ser $\mu(\K)(P)$
es un $R_T$--m\'odulo $\{\lambda^A\}_{A\in R_T}=\Lambda_{P^n}
\subseteq \mu(\K)(P)$. Se sigue que $\mu(\K)(P)=
\Lambda_{P^n}$.

Por tanto
\begin{align*}
\mu(\K)&=\bigoplus_{\substack{P\in R_T\\ \text{$P$ m\'onico e
irreducible}}} \mu(\K)(P)
=\bigoplus_{\substack{P\in R_T\\ \text{$P$ m\'onico e
irreducible}\\ \mu(\K(P))\neq 0}} \mu(\K)(P)\\
&=\bigoplus_{j=1}^{r}\Lambda_{P_j^{\alpha_j}}
=\Lambda_{P_1^{\alpha_1}\cdots P_r^{\alpha_r}}=\Lambda_M
\end{align*}
donde $M=P_1^{\alpha_1}\cdots P_r^{\alpha_r}$.
Esto termina la primera demostraci\'on.

Para la segunda demostraci\'on, se tiene que $\mu(L)$ es
el subm\'odulo de torsi\'on del $R_T$--m\'odulo $L$. Si
$x\in \mu(L)$, entonces existe $N\in R_T$, $N\neq 0$ con
$x^N=0$. Se sigue que tomamos tal $N$ de grado m\'inimo,
entonces $(N)=\an(x)$ es el anulador de $x$. Por tanto
$\Lambda_N\subseteq \mu(L)$.

En particular $\mu(L)$ es finito pues de lo contrario existir\'ian
una infinidad de $N_i\in R_T$ distintos con $\Lambda_{N_i}
\subseteq \mu(L)$ y $\cicl {N_i}{}\subseteq L$ de donde 
obtendr\'iamos $\infty>[L:K]\geq [\cicl {N_i}{}:K]\xrightarrow
[i\to\infty]{} \infty$, lo cual es absurdo.

Sea $\mu(L)=\{x_1,\ldots,x_m\}$ y sea ${\eu a}_i$ el anulador
de $x_i$. Esto es, ${\eu a}_i=\{A\in R_T\mid x_i^A=0\}$.
Entonces, si ${\eu a}_i=\langle A_i\rangle$, $x_i^{A_i}=0$,
$x_i\in \Lambda_{A_i}$ y $\Lambda_{A_i}\subseteq \mu(L)$.

Sea ${\eu b}=\bigcap_{i=1}^m{\eu a}_i=\langle M\rangle$.
Entonces ${\eu b}\subseteq {\eu a}_i$, $1\leq i\leq m$ y por
tanto $A_i|M$, $1\leq i\leq m$. Sean $C_i\in R_T$ tales que
$M=A_iC_i$, $1\leq i\leq m$. Si existiese $P\in R_T^+$ con
$P|C_i$ para toda $i$, entonces si escribimos $C_i=B_i P$
y obtendr\'iamos $M=A_iB_iP$ para toda $i$. Esto es,
si $M_1=\frac{M}{P}=A_iB_i\in {\eu a}_i$ para toda
$i$. En consecuencia, se tendr\'ia que $\langle M_1\rangle
\subseteq \bigcap_{i=1}^m {\eu a}_i={\eu b}=\langle M\rangle$.
Por tanto se obtendr\'ia que $M|M_1$ y si $M_1=MN_1$ con
$N_1\in R_T$, se seguir\'ia el absurdo $MN_1=\frac{M}{P}$.
Por tanto $\{C_1,\ldots,C_m\}$ son primos relativos.

Ahora para $x\in\mu(L)$, $x_i^{A_i}=0$ para alguna $i$ y
$x\in \Lambda_M$. Se sigue que $\mu(L)\subseteq \Lambda_M$.

Sea ahora $x\in \Lambda_M$. Puesto que $M=A_iC_i$, se tiene
$x^M=0=(x^{C_i})^{A_i}$, lo cual implica que $x^{C_i}\in 
\Lambda_{A_i}\subseteq \mu(L)$. Como $\{C_1,\ldots, C_m\}$
son primos relativos, existen $B_1,\ldots,B_m\in R_T$ tales que
$1=\sum_{i=1}^m B_iC_i$.

Por tanto $x=x^1=x^{\sum_{i=1}^m B_iC_i}=\sum_{i=1}^m
(x^{C_i})^{B_i}\in \mu(L)$. De esta forma obtenemos que
$\Lambda_M\subseteq \mu(L)$ y en consecuencia 
$\Lambda_M=\mu(L)$. Esto termina la segunda demostraci\'on.
$\fin$
\end{proof}

\begin{proposicion}\label{finitud_RC}
Sean $q > 2$, $M\in R_{T}$ no constante. Consideremos la extensi\'on
$K(\Lambda_{M})/K$. Entonces $\mu(K(\Lambda_{M})) = \Lambda_{M}$.
\end{proposicion}

\begin{proof} 
Sea $\mu(\cicl M{})=\Lambda_{N}$. Puesto que $\Lambda_{M}\subseteq
\mu(\cicl M{})$, si para un polinomio irreducible $P\in R_{T}$ y
$\alpha\in {\mathbb{Z}}$, $\alpha\geq 0$, tenemos que
$P^{\alpha}\mid M$, entonces $P^{\alpha}\mid N$. Si
$P^{\alpha+1}\nmid M$, no podemos tener que $P^{\alpha+1}\mid N$
puesto que en caso contrario el \'indice de ramificaci\'on de $P$ in
$K(\Lambda_{M})/K$ debe ser divido por $\Phi(P^{\alpha +
1})=[K(\Lambda_{P^{\alpha + 1}}):K]$, pero la ramificaci\'on de $P$ en
$K(\Lambda_{M})/K$ es $\Phi(P^{\alpha})$. As\'i $N=M$.
$\fin$
\end{proof}

En lo que sigue, a menos que se especifique otra cosa, las
extensiones de campos consideradas $L/\K$ satisfacen que $K\subseteq
\K\subseteq L\subseteq \overline{K}$. Por otro lado a las extensiones
anteriores se les da estructura de $R_{T}$-m\'odulo, usando la
acci\'on de Carlitz Hayes definida anteriormente. El primer objeto a considerar, asociado a la
extensi\'on $L/\K$, es el siguiente:
\[
T(L/\K) = \{u\in L\mid \text{existe un $M\in R_{T}
\setminus \{0\}$ tal que $u^{M} \in \K$}\}.
\]

N\'otese que $T(L/\K)\subseteq L$ es un subgrupo del grupo {\it
aditivo} $L$. Por otro lado $T(L/\K)$ es un $R_{T}$ - m\'odulo y el
$R_{T}$-m\'odulo $T(L/\K)/\K$ es de $R_{T}$-torsi\'on. A este \'ultimo
$R_{T}$ - m\'odulo lo denotamos por $\cog (L/\K)$. Se tiene que
$\cog(L/\K)$ es an\'alogo al grupo $T(L/\K)/\K^{*}$, en el caso de
una extensi\'on $L/\K$ de campos y $T(L/\K)$ denota el grupo de
torsi\'on usual, ver \cite{BaRzVi91} p.2.

\begin{definicion}\label{definicion1_c2}
Diremos que una extensi\'on $L/\K$ es {\it radical} si existe un
subconjunto $A\subseteq T(L/\K)$ tal que $L= \K(A)$.
Decimos que $L/\K$
es {\em pura\index{extensi\'on pura}} si para cada polinomio m\'onico irreducible $M\in
R_{T}$ y cada $u\in L$ tal que $u^{M} = 0$ se tiene que $u\in \K$.
Finalmente diremos que $L/\K$ es una extensi\'on {\it radical
ciclot\'omica} si:

(1) es radical,

(2) separable y

(3) pura.

Al m\'odulo $\cog(L/\K)= T(L/\K)/\K$ lo llamaremos {\em m\'odulo
de cogalois de la extensi\'on}.
\end{definicion}

A continuaci\'on probamos un resultado debido a Schultheis \cite{Sch90}.

\begin{teorema}\label{ShultheisProp2.3}
Sean $\K$ una extensi\'on finita de $K(\Lambda_M)$ y $z\in
\K\setminus\K^M$ y $F(u)= u^M-z$. Aqu\'i $\K^M=\{x^M\mid
x\in\K\}$. Sean
$F_1(u),\ldots,F_s(u)$ los distintos factores irreducibles de $F(u)$ en
$\K[u]$ y sea $\alpha\in \overline{K}$ cualquier ra{\'\i}z de $F_1(u)$.
Entonces el campo de descomposici\'on de $F(u)$ sobre $\K$
es $\K(\alpha)$. Adem\'as la extensi\'on $\K(\alpha)/
\K$ es elemental abeliana y en particular $[\K(\alpha):\K]=
p^t$ para alg\'un $t\in {\ma N}\cup\{0\}$.
\end{teorema}

\begin{proof} Puesto que las ra{\'\i}ces de $F(u)$ son los elementos del conjunto
$\{\alpha+\lambda\mid \lambda\in \Lambda_M\}$, se sigue que $\K
(\alpha)$ es el campo de descomposici\'on de $F(u)$. Si $G:=\Gal(
\K(\alpha)/\K)$, definimos $\varphi\colon G\to
\Lambda_M$ por $\varphi(\sigma)=\lambda_{\sigma}\in\Lambda_M$
donde $\sigma(\alpha)=\alpha+\lambda_{\sigma}$. Claramente $\varphi$
es un monomorfismo de grupos y puesto que $\Lambda_M$ es un 
$p$--grupo elemental abeliano, se sigue $G$ lo es y en particular
$|G|=[\K(\alpha):\K]=p^t$ para alg\'un $t\in{\ma N}\cup
\{0\}$. $\fin$
\end{proof}

\begin{ejemplo}\label{ejemplo1}
La extensi\'on $K(\Lambda_{M})/K$, con $M\in R_{T}$, es radical ya
que existe $W = \Lambda_{M}\subseteq T(K(\Lambda_{M})/K)$ tal que
$K(\Lambda_{M}) = K(W)$, es separable, pero no pura, ya que por la
Proposici\'on \ref{finitud_RC}, se tiene que la \'unicas ra\'ices de
Carlitz que est\'an en $K(\Lambda_{M})$ son $\Lambda_{M}$ y si $Q$ es
un factor irreducible de $M$, $\lambda_{Q}\in \Lambda_{M}$, pero no
est\'a en $K$. Por lo tanto $K(\Lambda_{M})/K$ no es una extensi\'on
radical ciclot\'omica.
\end{ejemplo}

El siguiente ejemplo muestra la existencia de extensiones radicales
ciclot\'omicas.

\begin{ejemplo}\label{ejemplo4}
Sea $p$ un primo impar, $q=p$ y $M = T$. Considere la extensi\'on
$K(\Lambda_{M})/K$ cuyo grado es $q - 1=p-1$. Se ha visto que
$K(\Lambda_{M})/K$ no es pura, ver Ejemplo \ref{ejemplo1}. Ahora
consideremos el polinomio $F(X) = X^{T} - 1 = X^{p} + X T - 1$.

Se afirma que $1\in K(\Lambda_{M})\setminus K(\Lambda_{M})^{M}$, ya
que si ocurre lo contrario, existe un $u\in K(\Lambda_{M})$ tal que
$u^{M} = 1$. Sea $\alpha$ un generador de $\Lambda_{M}$. Notemos
que $[K(\alpha):K]=p-1$, por lo que $\{1,\alpha,\alpha^{2},\ldots,
\alpha^{p-2}\}$ es una base de $K(\Lambda_{M})$ sobre $K$.

Por lo tanto $u$ se puede escribir como $u = a_{0} + a_{1} \alpha +
\cdots a_{p-2} \alpha^{p - 2}$ con $a_{0},a_{1}\cdots a_{p - 2}\in
K$. Por lo tanto
\begin{align}\label{ec5_1}
u^{T} &=  a^{T}_{0} + (a_{1} \alpha)^{T} + \cdots (a_{p-2} \alpha^{p
- 2})^{T}\nonumber\\
&= (a^{p}_{0} + a_{0} T) + (a^{p}_{1} \alpha^{p} + a_{1} \alpha T) +
\cdots + (a^{p}_{p-2} \alpha^{p (p-2)} + a_{p-2} \alpha^{p-2} T)
\end{align}

Como $\alpha^{T} = \alpha^{p} + \alpha T = 0$ entonces $\alpha^{p} =
- \alpha T$. Por lo tanto, puesto que $u^{T}=1$, de
(\ref{ec5_1}) se obtiene
\begin{align*}
1 &= (a^{p}_{0} + a_{0} T) + (a^{p}_{1} \alpha^{p} + a_{1} \alpha T)+
\cdots + (a^{p}_{p-2} \alpha^{p (p-2)} + a_{p-2} \alpha^{p-2} T)\\
&= (a^{p}_{0} + a_{0} T)+(-a^{p}_{1}\alpha T+a_{1}\alpha T) +\cdots
+ (-a^{p}_{p-2}\alpha^{p-2}T^{p-2}+ a_{p-2} \alpha^{p-2} T)
\end{align*}
es decir
\[
0=(a^{p}_{0} + a_{0} T-1)+ c_{1}\alpha+ c_{2}\alpha^{2}+\cdots + c_{p-2}\alpha^{p-2}
\]
donde $c_{i}= (-1)^{i}a^{p}_{i}T^{i}+ a_{i}T$,
$i=1,\ldots, p-2$, pertenecen a $K$.

Por lo tanto llegamos a la ecuaci\'on
\begin{equation}\label{Eq3.2}
0 = a^{p}_{0} + a_{0} T  - 1
\end{equation}
ya que $\{1,\alpha,\alpha^{2}\ldots, \alpha^{q - 2}\}$ es
base de $K(\Lambda_{M})$ sobre $K$. En particular $a_0\neq 0$.
Sea $a_{0} =
\frac{f(T)}{g(T)}$, con $(f(T),g(T)) = 1$, de aqu\'i derivamos la
ecuaci\'on $f^{p}(T) + f(T)g^{p - 1}(T) T = g^{p}(T)$. Se sigue que 
$f(T), g(T)\in{\ma F}_q^{\ast}$ y por lo tanto $a_0\in{\ma F}_q^{\ast}$
lo cual contradice (\ref{Eq3.2}).

Sea $L$ el campo de descomposici\'on de $F(X)$, sobre
$K(\Lambda_{M})$, entonces la extensi\'on $L/K(\Lambda_{M})$ es
separable. Del Teorema \ref{ShultheisProp2.3} obtenemos 
que $[L:K(\Lambda_{M})] = p^{t}$, con $t \geq 1$. Si $\beta$ es una
ra\'iz de $F(X)$ se tiene que $L = K(\Lambda_{M})(\beta)$, as\'i la
extensi\'on $L/K(\Lambda_{M})$ es radical. Notemos que como el
polinomio irreducible de $\beta$ divide a $F(X) = X^{p} + X T - 1$, entonces
tal irreducible es $F(X)$. En particular de esto se deduce que $t =
1$.

Para mostrar que la extensi\'on $L/K(\Lambda_{M})$ es radical
ciclot\'omica, bastar\'a mostrar que es pura, puesto que hemos
mostrado que es radical y separable. Para este fin consideremos
polinomios m\'onicos irreducibles $N$, con grado de $N$ $> 1$ y sea
$u\in L$ tal que $u^{N} = 0$. Se afirma que $u = 0$ pues en caso
contrario $u\neq 0$ es un generador de $\Lambda_{N}$, debido a que
$N$ es irreducible y a la Proposici\'on  \ref{P6.2.9}. 
As\'i se puede considerar el diagrama
\[
\xymatrix{& L\ar@{-}[dl]\ar@{-}[dr]  &\\
K(u)\ar@{-}[dr] & & K(\Lambda_{M})\ar@{-}[dl] \\
& K&}
\]

Ahora del Teorema \ref{T6.2.28} se tiene $[K(u):K] =
\Phi(N) =p^{\deg(N)}-1\geq p(p-1) = [L:K]$, pero esto contradice
que $[K(u):K]\mid [L:K]$. Por lo tanto $u=0\in K(\Lambda_{M})$. Esto
muestra la propiedad (3) de la Definici\'on \ref{definicion1_c2},
para los polinomios de grado mayor que 1.

Resta mostrar la propiedad (3) de la Definici\'on
\ref{definicion1_c2}, para los polinomios de grado 1. Para ello se
consideran los polinomios $T, T+1, \ldots, T+(p-1)$. Bastar\'a
considerar, por ejemplo, $N = T+1$. Sea $u\in L$ tal que $u^{T+1} =
0$ y supongamos que $u\notin K(\Lambda_{M})$, en particular $u\neq
0$. De esta manera $\Irr(u,X,K(\Lambda_{M}))\mid (X^{p-1} +
T+1)$, pero esto contradice nuestra suposici\'on de que $\deg(\Irr
(u,X,K(\Lambda_{M})))= p$. Por lo tanto $u\in
K(\Lambda_{M})$. 
\end{ejemplo}

Para el siguiente ejemplo necesitamos la siguiente proposici\'on

\begin{proposicion}\label{pureza_carlitz}
Sea $q > 2$, $P\in R_{T}$ m\'onico e irreducible y $n\in
{\mathbb{N}}$. Entonces la extensi\'on
$K(\Lambda_{P^{n}})/K(\Lambda_{P})$ es pura.
\end{proposicion}

\begin{proof} 
Si $\lambda_{Q}\in K(\lambda_{P^{n}})$, entonces $Q$ es ramificado
en $K(\lambda_{P^{n}})/K$ lo cual implica que $Q=P$, por la
Teorema \ref{T6.2.28}. Por
lo tanto $K(\lambda_{P^{n}})/K(\lambda_{P})$ es pura.
$\fin$
\end{proof}

\begin{ejemplo}\label{ejemplo5}

La extensi\'on $K(\Lambda_{P^{n}})/K(\Lambda_{P})$ es radical
ciclot\'omica, ya que, ciertamente, es radical, separable ya que el
polinomio, con coeficientes en $R_{T}$, $U^{P^{n}}$, es separable y
por la Proposici\'on \ref{pureza_carlitz} la extensi\'on es pura. 
\end{ejemplo}

\section{Algunas propiedades de las extensiones
radicales.}\label{prop_extensiones_radicales}

Las extensiones radicales $L/\K$ estudiadas aqu{\'\i},
tienen propiedades an\'alogas a las extensiones radicales usuales
consideradas en \cite{GreHar86} y en \cite{Alb2003}.

\begin{definicion}\label{modulotorsion_2bis}
Si $G$ es un m\'odulo de torsi\'on se pondr\'a
\[
{\mathcal{O}}_{G}=\{\ord(g)\mid g\in G\}.
\]
\end{definicion}

\begin{definicion}\label{modulotorsion2_2bis}
Un m\'odulo $G$ se dice {\em acotado} si $G$ es un m\'odulo de
torsi\'on y los grados de los elementos de
${\mathcal{O}}_{G}\subseteq R_{T}$ forman un conjunto acotado, o de
modo equivalente, ${\mathcal{O}}_{G}$ es finito.
\end{definicion}

Sea $A$ un $R_{T}$-m\'odulo de torsi\'on. Consideremos
${\mathcal{O}}_{A}$. Supongamos que $A$ es
un $R_{T}$-m\'odulo acotado. Al m\'inimo com\'un m\'ultiplo de los
elementos de ${\mathcal{O}}_{A}$, lo llamaremos el {\it
$R_{T}$-exponente de $A$} o, si el contexto lo permite, el exponente
de $A$, y lo denotamos por $\ex(A)$.

Ahora sea $E/F$ una extensi\'on radical, no necesariamente finita.
Existe un subconjunto $A\subseteq T(E/F)$ tal que $E=F(A)$. Podemos reemplazar
$A$ por el subm\'odulo de $E$ generado por $A$ y $F$, que seguiremos
denotando por $A$.

Ahora $A/F$ es un $R_{T}$-m\'odulo de torsi\'on, por lo que tiene sentido
considerar ${\mathcal{O}}_{A/F}$. Diremos que una extensi\'on de
$R_{T}$-torsi\'on, $E/F$, es una {\it extensi\'on acotada} si $A/F$
es un $R_{T}$-m\'odulo acotado, en este caso si $N= \ex(A/F)$,
diremos que $E/F$ es una extensi\'on {\it $N$ acotada}.

En este contexto se tiene la siguiente proposici\'on.

\begin{proposicion}\label{galois_torsion_raices}
Sea $E/F$ una extensi\'on radical acotada, no necesariamente finita,
y sea $N = \ex(A/F)$. Entonces $E/F$ es de Galois si y s\'olo si
$\lambda_{M}\in E$ para todo $M\in {\mathcal{O}}_{A/F}$.
\end{proposicion}

\begin{proof} 
Sea $\alpha\in E$ cuyo orden es $M\in R_{T}$. As\'i tenemos que
$\alpha^{M}=a\in F$. Consideremos el polinomio
$f(X)=X^{M}-a=\prod_{N} (X-(\alpha+\lambda^{N}_{M}))\in F[X]$. Por
lo tanto los conjugados de $\alpha$ son
\[
\{\alpha+\xi_{1},\ldots , \alpha+\xi_{s}\}
\]
para algunos $\xi_{i}\in \Lambda_{M}$.

Supongamos que la extensi\'on $E/F$ es Galois.
Sea $B$ el $R_{T}$-m\'odulo generado por $\{\xi_{1},\ldots
,\xi_{s}\}$. Entonces $B\subseteq E$ y existe un $M^{\prime}\in
R_{T}$, que divide a $M$, tal que $B=\Lambda_{M^{\prime}}$. Si
$M^{\prime}\neq M$, entonces $\alpha^{M^{\prime}}=a^{\prime}\in F$
lo cual es una contradicci\'on. Por lo tanto $M^{\prime}=M$ y
$\lambda_{M}\in E$.

Ahora supongamos que $\lambda_{M}\in E$ para todo
$M\in{\mathcal{O}}_{A/F}$. Sea $u\in A$ y $M = \ord(u)$.
Puesto que todo conjugado de $u$, sobre $F$, es de la forma
$u+\lambda^{N}_{M}\in E$, se sigue que la extensi\'on $E/F$ es
normal, y como $u$ es separable sobre $F$, entonces
$E/F$ es una extensi\'on de Galois.
$\fin$
\end{proof}

En algunas extensiones radicales $L/\K$, es posible encontrar un
elemento primitivo expl\'icito y que pertenezca
$\cog(L/\K)$, como lo muestra la siguiente proposici\'on.

\begin{proposicion}\label{primitivo_explicito}
Sea $L/\K$ es una extensi\'on tal que $L=\K(\alpha,\beta)$ y existen
$M,N\in R_{T}$ con $\alpha^{M}=a$, $\beta^{N}=b$, $a,b\in \K$, $M$ y
$N$ primos relativos. Entonces $L=\K(\alpha+\beta)$, es decir,
$\alpha+\beta$ es un elemento primitivo.
\end{proposicion}

\begin{proof} 
Puesto que $\alpha+\beta\in \K(\alpha,\beta)$ se tiene
$\K(\alpha+\beta)\subseteq \K(\alpha,\beta)$. Por otro lado
$(\alpha+\beta)^{M}= \alpha^{M}+\beta^{M}= a+\beta^{M}\in
\K(\alpha+\beta)$ y $(\alpha+\beta)^{N}= \alpha^{N}+\beta^{N}=
\alpha^{N}+b\in \K(\alpha+\beta)$. Por lo tanto se tiene que
$\beta^{M},\alpha^{N}\in \K(\alpha+\beta)$.

Ahora puesto que existen $S_{1},S_{2}\in R_{T}$ tales que
$1=MS_{1}+NS_{2}$ se tiene que
\begin{gather*}
\alpha=\alpha^{1}=\alpha^{MS_{1}+NS_{2}}=a^{S_{1}}+(\alpha^{N})^{S_{2}}\in \K(\alpha+\beta)\\
\intertext{y}
\beta=\beta^{1}=\beta^{MS_{1}+NS_{2}}= (\beta^{N})^{S_{1}}+ b^{S_{2}}\in \K(\alpha+\beta).
\end{gather*}

As\'i pues, $\K(\alpha,\beta)= \K(\alpha+\beta)$. M\'as a\'un,
$(\alpha+\beta)^{MN}= (\alpha^{M})^{N}+(\beta^{N})^{M}\in \K$.
$\fin$
\end{proof}

En particular, con las hip\'otesis de la Proposici\'on
\ref{primitivo_explicito}, se tiene que
\[
[L:\K]\leq \mid \cog(L/\K)\mid.
\]

Notemos que el argumento de la Proposici\'on
\ref{primitivo_explicito} se puede generalizar a extensiones de la
forma $L/\K$, con $L=\K(\alpha_{1},\ldots , \alpha_{s})$ de modo que
existen $M_{i}\in R_{T}$ con $\alpha^{M_{i}}_{i}= a_{i}\in \K$ y los
polinomios $M_{i}$ primos relativos a pares.

\section{Algunas propiedades de las extensiones radicales ciclot\'omicas}\label{rad_ciclotomicas_props}

Las extensiones radicales ciclot\'omicas tiene algunas propiedades
an\'alogas a las propiedades de las extensiones cogalois cl\'asicas.
Necesitamos primero un lema.

\begin{lema}\label{pureza}
Sea $\K\subseteq L\subseteq L^{\prime}$ una torre de campos. Entonces
$L^{\prime}/\K$ es pura si y s\'olo si $L^{\prime}/L$ y $L/\K$ son
puras.
\end{lema}

\begin{proof} 
Supongamos que $L^{\prime}/\K$ es pura. Sean $\lambda\in L^{\prime}$
y $P\in R_{T}$, m\'onico e irreducible, tal que $\lambda^{P}_{P} =
0$. Entonces $\lambda_{P}\in \K\subseteq L$, puesto que
$L^{\prime}/\K$ es pura. Por lo tanto $L^{\prime}/L$ es pura. De modo
completamente an\'alogo se prueba que $L/\K$ es pura.

Por otro lado supongamos que $L^{\prime}/L$ y $L/\K$ son puras. Sean
$\lambda_{P}\in L^{\prime}$ y $P\in R_{T}$, m\'onico e irreducible,
tal que $\lambda^{P}_{P} = 0$. Puesto que $L^{\prime}/L$ es
pura, $\lambda_{P}\in L$ y como $L/\K$ es pura, entonces
$\lambda_{P}\in \K$.
$\fin$
\end{proof}

\begin{proposicion}\label{prop_cog}
Sea $\K\subseteq L\subseteq L^{\prime}$ una torre de campos.
Se tienen las siguientes propiedades.
\begin{itemize}
\item[{\rm (1)}] Existe una sucesi\'on exacta de $R_{T}$-m\'odulos
\[
0\rightarrow \cog(L/\K)\rightarrow \cog(L^{\prime}/\K)\rightarrow \cog(L^{\prime}/L).
\]
\item[{\rm (2)}] Si la extensi\'on $L^{\prime}/\K$ es radical ciclot\'omica, entonces la
extensi\'on $L^{\prime}/L$ es radical ciclot\'omica.
\item[{\rm (3)}] Si la extensi\'on $L^{\prime}/\K$ es radical, y las extensiones $L^{\prime}/L$ y $L/\K$ son radicales ciclot\'omicas, entonces
$L^{\prime}/\K$ es radical ciclot\'omica.
\end{itemize}
\end{proposicion}

\begin{proof} 
(1) El homomorfismo can\'onico
\[
\cog(L^{\prime}/\K)\rightarrow \cog(L^{\prime}/L),\,\, x+\K\mapsto x+L
\]
es un $R_{T}$-homomorfismo con n\'ucleo $\cog(L/\K)$.
Esto prueba que la sucesi\'on de $R_{T}$-m\'odulos
\[
0\rightarrow \cog(L/\K)\rightarrow \cog(L^{\prime}/\K)\rightarrow \cog(L^{\prime}/L)
\]
es exacta.

(2) Como $L^{\prime}/\K$ es separable, entonces $L^{\prime}/\K$ es
separable y, por el Lema \ref{pureza}, $L^{\prime}/\K$ es pura.
Finalmente puesto que $T(L^{\prime}/\K)\subseteq T(L^{\prime}/L)$
se tiene que $L^{\prime}/L$ es radical.

(3) Como $L^{\prime}/L$ y $L/\K$ son extensiones
radicales ciclot\'omicas
entonces ambas son separables y puras. Por lo tanto, por el Lema
\ref{pureza}, la extensi\'on $L^{\prime}/\K$ es pura, adem\'as
separable. Se sigue que $L^{\prime}/\K$ es 
una extensi\'on radical ciclot\'omica.
$\fin$
\end{proof}

Veremos que para algunas extensiones
$L/\K$ se tiene que el $R_{T}$-m\'odulo
$\cog(L/\K)$ es finito. Para empezar considere $L/\K$ una
extensi\'on de Galois de campos de funciones, con grupo de Galois $G
= \Gal(L/\K)$. Notemos que $\mu(L)$ es un $G$-m\'odulo, mediante la
acci\'on siguiente: dado $\sigma\in G$ y $u\in \mu(L)$ pongamos
$\sigma\cdot u = \sigma(u)$. Puesto que la acci\'on de Carlitz Hayes
conmuta con $\sigma$, $\sigma\cdot u$ est\'a bien definida.

\begin{definicion}\label{definicion_morfismo_cruzado}
Una funci\'on $f:G\rightarrow \mu(L)$ se dice que es un {\it
homomorfismo cruzado de $G$ con coeficientes en $\mu(L)$} si para
cada $\sigma,\tau\in G$ se tiene que $f(\sigma\circ\tau) = f(\sigma)
+ \sigma\cdot f(\tau)$.

Al conjunto de homomorfismos cruzados los denotamos por
\begin{gather*}
Z^{1}(G,\mu(L)),\\
\intertext{y $B^1(G,\mu(L))$ denota al subconjunto de 
$Z^1(G,\mu(L))$ dado por}
\{\chi\in Z^1(G,\mu(L))\mid \text{existe $u\in \mu(L)$ tal que\ }
\chi=f_u\},\\
\intertext{donde $f_u$ es la funci\'on definida por}
f_u(\sigma):=\sigma u -u \quad \text{para cada}\quad \sigma \in G.
\end{gather*}
\end{definicion}

\begin{teorema}\label{finitud_TC/L}
Sea $L/\K$ una extensi\'on finita de Galois, $G$ su grupo de Galois.
Entonces la funci\'on $\phi:\cog(L/\K)\rightarrow
Z^{1}(G,\mu(L))$, dada por $\phi(u + \K) = f_{u}$ donde
$f_{u}(\sigma) = \sigma(u) - u$, es un isomorfismo de grupos.
\end{teorema}

\begin{proof} 
Se define $\theta:T(L/\K)\rightarrow Z^{1}(G,\mu(L))$ mediante
$\theta(u) = f_{u}$. Obs\'ervese que $f_{u}(\sigma\circ\tau) =
\sigma(\tau(u)) - u$, adem\'as $f_{u}(\sigma) = \sigma(u) - u$ y
$f_{u}(\tau) = \tau(u) - u$. Aplicando a esta \'ultima ecuaci\'on
$\sigma$ se obtiene $\sigma(f(\tau)) = \sigma(\tau(u)) - \sigma(u)$.
Al sumar esta ecuaci\'on con la primera se obtiene que $f_{u}$ es un
homomorfismo cruzado. Notemos de paso que si $\sigma\in G$ entonces
$f_{u}(\sigma) = \sigma(u) - u$ est\'a en $\mu(L)$, puesto que
existe un $N\in R_{T}$ tal que $u^{N}\in \K$ por lo tanto $(\sigma(u)
- u)^{N} = (\sigma(u))^{N} - u^{N} = \sigma(u^{N}) - u^{N} = 0$.

Adem\'as $\theta(u + v) = f_{u + v}$ y $f_{u + v}(\sigma) = \sigma(u
+ v) - (u + v) = \sigma(u) + \sigma(v) - u - v = \sigma(u) - u +
\sigma(v) - v$, es decir, $\theta(u + v) = \theta(u) + \theta(v)$.
Por lo tanto $\theta$ es un homomorfismo. Por otra parte, sea $u\in
\ker(\theta)$. As\'i $\theta(u) = f_{u} = 0$, es decir,
$f_{u}(\sigma) = \sigma(u) - u = 0$, y como $L/\K$ es de Galois,
entonces $u\in \K$.

Rec\'iprocamente si $u\in \K$, ciertamente $\theta(u) = 0$. As\'i
$\ker(\theta) = \K$ y por lo tanto tenemos un monomorfismo de grupos
abelianos
\[
\phi:\cog(L/\K)\rightarrow Z^{1}(G,\mu(L)).
\]

Por otro lado $Z^{1}(G,\mu(L))\subseteq Z^{1}(G,L)$ y por el Teorema
90 de Hilbert aditivo, se tiene que
\[
Z^{1}(G,L) = B^{1}(G,L)
= \{f\in Z^{1}(G,L)\mid \text{existe un $u\in L$ tal que $f =
f_{u}$}\}.
\]

Entonces, dado $f\in Z^{1}(G,\mu(L))$ existe un $u\in L$ tal que $f =
f_{u}$, por lo que para cada $\sigma\in G$, $f(\sigma) =
f_{u}(\sigma) = \sigma(u) - u \in \mu(L)$. Para probar que
$\phi$ es suprayectiva, es necesario probar que $u\in
\cog(L/\K)$.

Ahora, $u$ es algebraico
sobre $\K$, y se puede considerar la cerradura de Galois $\K^{\prime}$
de $\K(u)/\K$. Se tiene que $\K\subseteq \K(u)\subseteq
\K^{\prime}\subseteq L$.

Sea $H = \Gal(L/\K^{\prime})$. Entonces $H\lhd G$ y
$\textrm{card}(G/H)$ es finita. Los conjugados de $u$ son
$\{\overline{\sigma}(u)\mid \overline{\sigma}\in \overline{G} =
G/H\}$, as\'i
\[
\overline{\sigma}(u) = \sigma(u) = u + z_{\sigma}
\text{ con $z_{\sigma}\in \mu(L)$}.
\]

Ahora bien, puesto que \'unicamente hay un n\'umero finito de elementos
$\overline{\sigma}\in \overline{G} = \{\sigma_{1}H,\ldots,
\sigma_{s}H\}$ existen $N_{\sigma_{1}},\ldots, N_{\sigma_{s}}\in
R_{T}$ tales que $(z_{\sigma_{i}})^{N_{\sigma_{i}}} = 0$, sea $N =
N_{\sigma_{1}}\cdots N_{\sigma_{s}}$ entonces
\[
\overline{\sigma}(u^{N})= (u + z_{\sigma})^{N} = u^{N} + z_{\sigma}^{N} = u^{N} +
(z^{N_{\sigma}}_{\sigma})^{P_{\sigma}}= u^{N},
\]
donde $N=N_{\sigma}P_{\sigma}$.

Como la extensi\'on $\K^{\prime}/\K$ es de Galois, esto implica que
$u^{N}\in \K$,  por lo que $u\in \cog(L/\K)$ y $\phi$ es suprayectiva.
$\fin$
\end{proof}

Del Teorema \ref{finitud_TC/L}, obtenemos el siguiente resultado.

\begin{proposicion}\label{dummit}
Sean $E/F$ una extensi\'on finita de Galois con grupo
de Galois $\Gamma = \Gal(E/F)$ y
$\Delta$ un subgrupo normal de $\Gamma$. Entonces la sucesi\'on can\'onica
de grupos abelianos
\[
0\rightarrow Z^{1}(\Gamma/\Delta, \mu(E/F)^{\Delta})\stackrel{\theta_{1}}
{\rightarrow} Z^{1}(\Gamma,\mu(E/F))\stackrel{\theta_{2}}{\rightarrow} Z^{1}(\Delta,\mu(E/F))
\]
es exacta, donde $\mu(E/F)^{\Delta}= \{\zeta\in \mu(E/F)\mid
\sigma(\zeta) = \zeta\text{\ para toda\ } \sigma\in \Delta\}$.
\end{proposicion}

\begin{proof} 
Supongamos que $\theta_{1}(f) = 0$. Entonces si $\overline{\sigma}\in
\Gamma/\Delta$, se tiene que $f(\overline{\sigma}) =
\theta_{1}(f)(\sigma) = 0$. De este modo $\theta_{1}$ es inyectiva.
Por otro lado $\im(\theta_{1})\subseteq \ker(\theta_{2})$ ya
que si $f = \theta_{1}(f^{\prime})$, con $f^{\prime}\in
Z^{1}(\Gamma/\Delta, \mu(E/F)^{\Delta})$, entonces
$\theta_{2}(f)(\sigma) = \theta_{1}(f^{\prime})(\sigma) =
f^{\prime}(\overline{\sigma}) = 0$.

Ahora si $f\in \ker(\theta_{2})$ entonces para cada $\sigma\in
\Delta$, se tiene que $f(\sigma) = 0$. Por lo se puede definir
$f^{\prime}:\Gamma/\Delta\rightarrow \mu(E/F)$ mediante
$f^{\prime}(\overline{\sigma}) = f(\sigma)$. Por la condici\'on
impuesta a $f$, $f^{\prime}$ esta bien definida y es un morfismo
cruzado. Finalmente si $\tau\in \Delta$ entonces
$\tau(f^{\prime}(\overline{\sigma})) = \tau(f(\tau^{-1}\circ\sigma))
= \tau(f(\sigma)) = f^{\prime}(\overline{\sigma})$, es decir,
$f^{\prime}\in Z^{1}(\Gamma/\Delta, \mu(E/F)^{\Delta})$ y $f =
\theta_{1}(f^{\prime})$.
$\fin$
\end{proof}

\begin{corolario}\label{cor_finitud_TC/L}
Sea $L/\K$ una extensi\'on de Galois finita. Si la cardinalidad de
$\mu(L)$ es finita entonces el $R_{T}$-m\'odulo $\cog(L/\K)$ es
finito.
\end{corolario}

\begin{proof} 
Se sigue de la Proposici\'on \ref{finitud_TC/L}.
$\fin$
\end{proof}

\section{Algunos teoremas de
estructura de extensiones radicales
ciclot\'omicas}\label{teoremas_de_estructura}

\begin{proposicion}\label{cogalois_disjunto}
Sea $L/\K$ una extensi\'on de campos, tal que $[L:\K]=\ell$ con $\ell$
un primo diferente a $p = \car(K)$. Entonces $L/\K$ no es
radical ciclot\'omica.
\end{proposicion}

\begin{proof} 
Supongamos que $L/\K$ es radical ciclot\'omica, por lo tanto
$\cog(L/\K)$ es no trivial. Sea $\overline{\alpha}\in
\cog(L/\K)$ distinto de $0$, esto significa que $\alpha\notin
\K$. As\'i, existe un $M\in R_{T}$ tal que $\alpha^{M}\in \K$.
Podemos
suponer que $M$ es m\'onico y que es el polinomio de grado m\'inimo
con tal propiedad, es decir, el orden de $\overline{\alpha}$ es $M$.
Por lo que es posible suponer que existe un polinomio irreducible $Q$,
reemplazando a $\alpha$ si es necesario, tal que $\alpha^{Q} = a\in
\K$.

Sea $f(X) = \Irr(\alpha,X,\K)\in \K[X]$,
puesto que $\alpha^{Q}- a =
0$ entonces $f(X)\mid
X^{Q}- a$. Por lo tanto $f(X) = \prod(X - (\alpha +
\lambda^{B}_{Q}))$, para ciertos $B\in R_{T}$. Observemos que
$\deg(f(X))= \ell$, pues $L=\K(\alpha)$ y por lo tanto
$\sum(\alpha + \lambda^{B}_{Q})=
\ell\alpha + \lambda^{\sum B}_{Q}\in \K$. Por otro lado, puesto que
$\ell\neq p$ entonces $\ell\neq 0$ en $\K$. As\'i pues $D = \sum B$
es diferente de cero pues, en caso contrario, tendr{\'\i}amos que
 $\alpha\in \K$. Por tanto
podemos suponer que el grado de $D$ es menor que el grado de
$Q$.

Por otra parte $\lambda^{D}_{Q}\notin \K$, pero $\lambda^{D}_{Q}\in
L$ lo que contradice que la extensi\'on
$L/\K$ es pura. Por lo tanto $L/\K$ no es 
una extensi\'on radical ciclot\'omica.
$\fin$
\end{proof}

\begin{corolario}\label{cogalois_disjuntocor1}
Sea $L/\K$ una extensi\'on de Galois, tal que $[L:\K]= p^{s}n$, con
$p\nmid n$, $n>1$ y $p=\car(\K)$. Entonces $L/\K$ no es
una extensi\'on radical ciclot\'omica.
\end{corolario}

\begin{proof} 
Por el teorema de Cauchy el grupo $G = \Gal(L/\K)$ tiene un elemento
de orden $\ell$, digamos $g$, donde $\ell$ es un primo que divide a
$n$. Considere el subgrupo $H = (g)$ de $G$. Si $L/\K$ fuese radical
ciclot\'omica entonces, por la Proposici\'on \ref{prop_cog}, la
extensi\'on $L/L^{\prime}$, donde $L^{\prime} = L^{H}$, ser{\'\i}a radical
ciclot\'omica. Pero $[L:L^{\prime}] = \ell$ y por la Proposici\'on
\ref{cogalois_disjunto} tal extensi\'on no es radical ciclot\'omica.
Por lo tanto $L/\K$ no es radical ciclot\'omica.
$\fin$
\end{proof}

\begin{corolario}\label{estructura1}
Si $L/\K$ es Galois y radical ciclot\'omica, entonces $[L:\K]$ es de
la forma $p^{s}$, con $s\in {\mathbb{N}}$ y $p=\car(\K)$.
\end{corolario}

\begin{proof} 
Si ocurre lo contrario, se tendr\'a que $[L:\K]= p^{n} m$, con $n$
un entero mayor o igual a $0$, $p\nmid m$ y $m > 1$. Sin embargo por el
Corolario \ref{cogalois_disjuntocor1} $L/\K$ no ser\'ia radical
ciclot\'omica, lo cual es una contradicci\'on.
$\fin$
\end{proof}

\begin{lema}\label{pureza_p}
Sea $L/\K$ una extensi\'on tal que $[L:\K] = p^{s}$ con $s\in
{\mathbb{N}}$ y $p=\car(\K)$. Entonces $L/\K$ es pura.
\end{lema}

\begin{proof} 
Supongamos que $L/\K$ no es pura, as\'i existe un $a=\lambda_{P}\in
L$, con $P\in R_{T}$ irreducible tal que $a^{P}= 0$ pero $a\notin \K$.
Considere el diagrama siguiente
\[
\xymatrix{& L\ar@{-}[d]\\
K(\lambda_{P})\ar@{-}[d]\ar@{-}[r]  & K(\lambda_{P})\K = \K(\lambda_{P}) \ar@{-}[d]\\
K \ar@{-}[r] & \K  }
\]

Sea $\widetilde{\K}=\K\cap K(\lambda_{P})$.
Entonces, por Teor\'ia de
Galois, se tiene que $\K(\lambda_{P})/\K$ es Galois, con grupo de
Galois $G$ isomorfo a $\Gal(K(\lambda_{P})/\widetilde{\K})$. Por otro
lado
\[
\mid G\mid\mid [L:\K]=p^{s}\quad \text{y}\quad \mid G\mid\mid (q^{d}-1)
\]
donde $d=\deg(P)$. Por lo tanto $\mid G\mid = 1$, es
decir, $\lambda_{P}\in \K$.
$\fin$
\end{proof}

\begin{ejemplo}\label{schu_cogalois}
Una extensi\'on de {\it Carlitz--Kummer}, ver \cite{Sch90}, es
una extensi\'on $L/\K$ tal que

\medskip
(1) $\K$ es una extensi\'on finita de $K(\Lambda_{M})$, para alg\'un
$M\in R_{T}$.

(2) $L$ es campo de descomposici\'on del polinomio $f(X)=X^{M}-z\in
\K[u]$, sobre $\K$, donde 
$z\in \K\setminus \K^{M}$.

Por el Teorema \ref{ShultheisProp2.3}, se tiene
que $[L:\K]=p^{t}$, donde $p=\car(\K)$. Ahora el Lema
\ref{pureza_p} muestra que las extensiones de Carlitz Kummer son
extensiones radicales ciclot\'omicas.
\end{ejemplo}

Por otro lado, como consecuencia del Lema \ref{pureza_p},
se tiene el siguiente teorema.

\begin{teorema}\label{tdim_p_galois}
Una extensi\'on, de Galois, $L/\K$ es radical ciclot\'omica si y
s\'olo si es radical, separable y $[L:\K]=p^{s}$ con $s\in
{\mathbb{N}}$ y $p=\car(\K)$. $\fin$ 
\end{teorema}

En este contexto se tiene el siguiente teorema.

\begin{teorema}\label{tdim_p}
Si $L/\K$ es radical ciclot\'omica, entonces $[L:\K]=p^{n}$ para
alguna $n\geq 0$, donde $p=\car(\K)$.
\end{teorema}

\begin{proof} 
Sea $L/\K$ una extensi\'on
radical ciclot\'omica. Entonces $L=\K(\alpha_{1},\ldots,
\alpha_{t})$ de tal modo que $\alpha^{M_{i}}_{i}=a_{i}\in \K$ donde
$M_{i}\in R_{T}$. Entonces
\[
[L:\K]=[L:\K(\alpha_1,\ldots,\alpha_{t-1})]\cdots
[\K(\alpha_1,\alpha_2):\K(\alpha_1)]
[\K(\alpha_1):\K].
\]
Puesto que cada $\K(\alpha_1,\ldots,\alpha_i)/
\K(\alpha_1,\ldots,\alpha_{i-1})$ es una extensi\'on finita
radical ciclot\'omica, es suficiente considerar el caso 
$L=\K(\alpha)$.

Supongamos que $L=\K(\alpha)$ con $\alpha^M\in
\K$ para alg\'un $M\in R_T$. Sea $M=P_1^{e_1}
\cdots P_s^{e_s}$ su factorizaci\'on como producto de polinomios
irreducibles distintos. Sea $\beta_i:=\alpha^{M/P_j^{e_j}}$ para
$1\leq j\leq s$. Se tiene que $L=\K(\beta_1,\ldots \beta_s)$.
Por el mismo argumento anterior, podemos suponer
$L=\K(\alpha)$ con $\alpha^{P^e}\in \K$ para
alg\'un polinomio irreducible $P\in R_T$.

Ahora sea $L=\K(\alpha)$ tal que $\alpha^{P^e}\in \K$
para alg\'un polinomio irreducible $P\in R_T$. Sean $\gamma_i=
\alpha^{P^{e-i}}$ para $1\leq i\leq e$. Entonces 
$\K\subseteq \K(\gamma_1)
\subseteq \K(\gamma_2)\subseteq \cdots \subseteq
\K(\gamma_e)=L$. Por lo tanto
\[
[L:\K]=[L:\K(\gamma_{e-1})]\cdots
[\K(\gamma_2):\K(\gamma_1)]
[\K(\gamma_1):\K].
\]
Puesto que cada $K(\gamma_i)/\K(\gamma_{
i-1})$ es una extensi\'on radical ciclot\'omica, es suficiente
considerar el caso $L=K(\alpha)$ con $\alpha^P\in\K$
para alg\'un polinomio irreducible $P\in R_T$.

Supongamos $\lambda_{P}\in L$. Entonces $L/\K$ es de Galois
por ser el campo de descomposici\'on de $X^P-\alpha^P\in \K[X]$. Por el
Corolario \ref{estructura1}, $L/\K$ es una $p$-extensi\'on.

Ahora supongamos que $\lambda_P\notin L$. Consideremos
el diagrama
\[
\xymatrix{& & L=\K(\alpha)\ar@{-}[r]^{a} & L(\lambda_{P})=\K(\lambda_{P},\alpha)\\ 
& \K\ar@{-}[r]^{d} & \K(\alpha)\cap \K(\lambda_{P})\ar@{-}[r]^{a}\ar@{-}[u]^{b}&  \K(\lambda_{P})\ar@{-}[u]^{b}\\ K\ar@{-}[r]&\K\cap K(\lambda_{P})\ar@{-}[r]^{d}\ar@{-}[u]^{c} & \K(\alpha)\cap K(\lambda_{P})\ar@{-}[r]^{a}\ar@{-}[u]^{c}&
K(\lambda_{P})\ar@{-}[u]^{c}}
\]

Puesto que $\K(\lambda_{P},\alpha)/\K(\lambda_{P})$ es Galois
adem\'as, por el Teorema \ref{ShultheisProp2.3} se
tiene que $N=\Gal(L(\lambda_{P})/\K(\lambda_{P}))$ puede considerarse
como un subgrupo de $\Lambda_{P}$, es decir, $N$ es un $p$-grupo
elemental abeliano y $\mid N\mid=b=p^{n}$.

Puesto que
\[
[L:\K]=[L:\K(\alpha)\cap \K(\lambda_{P})][\K(\alpha)
\cap \K(\lambda_{P}):\K]=bd=p^{n}d
\]
basta mostrar que $d=1$.

Se tiene que $L(\lambda_P)/\K(\alpha)\cap \K(\lambda_P)=M$
es una extensi\'on de Galois pues $L=\K(\lambda_P,\alpha)$
y si $\sigma\colon L(\lambda_P)\lra\overline{L(\lambda_P)}$ con
$\sigma|_M=\Id_M$, entonces $\sigma(\lambda_P)=\lambda_P^Q
\in L(\lambda_P)$ y $\sigma\alpha=\alpha+\lambda_P^A\in
L(\lambda_P)$.

Sean $H=\Gal(L(\lambda_{P})/(\K(\alpha)\cap \K(\lambda_{P})))$,
$G=\Gal(L(\lambda_{P})/\K)$ y
$N=\Gal(L(\lambda_{P})/\K(\lambda_{P}))$. N\'otese que $N$ es un
subgrupo normal de $G$.

Se tiene que
\[
G/N\cong \Gal(\K(\lambda_{P})/\K)<
\Gal(K(\lambda_{P})/K)\cong
C_{q^{d}-1}.
\]

As\'i pues $G/N$ es un grupo c\'iclico de orden $q^{d}-1$, en
particular, primo relativo a $p$. Adem\'as se tiene que $\mid
G/N\mid=ad$.

Por el Teorema de Hall, ver \cite{Hal69}, Teorema 9.3.1, como $G$ es
soluble, existe un subgrupo $R$ de $G$ con $R$ c\'iclico de orden $ad$,
tal que $G=NR$ (de hecho $G$ es el producto semidirecto
$G\cong N\rtimes R$ ya que $(\mid R\mid,\mid N \mid)=1$).

Por el mismo Teorema de Hall, todo subgrupo de orden un divisor de
$\mid R\mid =ad$ esta contenido en un conjugado $R^{\prime}$ de $R$
y se tiene que $G=NR^{\prime}\cong N\rtimes R^{\prime}$.

Sea $S=\Gal(L(\lambda_{P})/\K(\alpha))\cong C_{a}$. Por lo tanto
podemos suponer $S\subseteq R$ y $\mid R/S\mid=d$. Notemos que
$(d,p)=1$.

Sea $E=L(\lambda_{P})^{R}$. Observemos que $L(\lambda_{P})^{S}=
\K(\alpha)=L$. Por lo tanto $\K\subseteq E\subseteq L$,
$[L:E]=[R:S]=d=\mid R/S\mid$. Ahora, como $L/\K$ es
una extensi\'on radical
ciclot\'omica lo es tambi\'en $L/E$. Por lo tanto $d=1$.
$\fin$
\end{proof}

\begin{corolario}\label{tdim_p_colateral}
Con las notaciones del Teorema {\rm \ref{tdim_p}} tenemos
$\K(\alpha)\cap
\K(\lambda_{P})=\K$,
$[L:\K]=[L(\lambda_{P}):\K(\lambda_{P})]$.
Adem\'as
\[
\Irr(\alpha,X,\K)=
\Irr(\alpha,X,\K(\lambda_{P}))=F_{1}(X)=\prod(X-(\alpha+\lambda^{A}_{P})).
\]
\end{corolario}

\begin{proof} 
Se sigue de la demostraci\'on del Teorema \ref{tdim_p}.
$\fin$
\end{proof}

\begin{corolario}\label{cogalois_grado_p_pureza}
Una extensi\'on finita $L/\K$ es radical ciclot\'omica si y s\'olo si es
separable, radical y $[L:\K]=p^{m}$ para alg\'un $m\in {\mathbb{N}}$.
\end{corolario}

\begin{proof} 
Se sigue del Teorema \ref{tdim_p} y del Lema \ref{pureza_p}.
$\fin$
\end{proof}

\section{Ejemplos y aplicaciones}\label{examples}

En esta secci\'on veremos algunas aplicaciones de los resultados
anteriores. En primer lugar tenemos la siguiente consecuencia del
Teorema \ref{finitud_TC/L}.

\begin{corolario}\label{finitud_redes}
Si $E/L$ es una extensi\'on finita y de Galois, con grupo de Galois $\Gamma$,
entonces la funci\'on:
\[
\phi:\{H\mid L\leq H \leq T(E/L)\}\rightarrow \{U\mid U\leq Z^{1}(\Gamma,\mu(E))\},
\]
dada por $\phi(H) = \{f_{\alpha}\in
Z^{1}(\Gamma,\mu(E))\mid \alpha\in H\}$, es un isomorfismo de redes.
\end{corolario}

\begin{proof} 
Se sigue del isomorfismo dado en el Teorema
\ref{finitud_TC/L}.
$\fin$
\end{proof}

Ahora, sea $E/L$ una extensi\'on de Galois con grupo de Galois
$\Gamma$. Definimos
\[
f:\Gal(E/L)\times \cog(E/L)\rightarrow \mu(E)
\]
dado por $f(\sigma,\overline{u}) = \sigma(u) - u$. Puesto
que $\cog(E/L)\rightarrow Z^{1}(\Gamma, \mu(E))$ es un
isomorfismo, se tiene la funci\'on evaluaci\'on
\[
\langle\ ,\ \rangle:\Gamma\times Z^{1}(\Gamma, \mu(E))\rightarrow \mu(E)
\]
dado por $\langle\sigma,h\rangle = h(\sigma)$.

Ahora para cada $\Delta\leq \Gamma$, $U\leq Z^{1}(\Gamma, \mu(E))$ y
$\chi\in Z^{1}(\Gamma, \mu(E))$ definimos:
\begin{gather*}
\Delta^{\bot} = \{h\in Z^{1}(\Gamma, \mu(E))\mid \langle\sigma,h\rangle = 0 \text{ para cada $\sigma\in \Delta$}\},\\
U^{\bot} = \{\sigma\in \Gamma\mid \langle\sigma,h\rangle = 0 \text{ para cada $h\in U$}\},\\
\chi^{\bot} = \{\sigma\in \Gamma\mid \langle\sigma,\chi\rangle = 0\}.
\end{gather*}
As\'i $\Delta^{\bot}\leq Z^{1}(\Gamma,\mu(E))$ y
$U^{\bot}\leq \Gamma$.

\begin{proposicion}\label{redes_casineat}
Sea $E/L$ una extensi\'on finita y de Galois con
grupo de Galois $\Gamma$.
Sea $L^{\prime}$ una extensi\'on intermedia de $E/L$. Entonces
$L^{\prime}/L$ es radical si y solamente si existe un subgrupo
$U\leq Z^{1}(\Gamma,\mu(E))$ tal que $\Gal(E/L^{\prime}) =
U^{\bot}$.
\end{proposicion}

\begin{proof} 
Si $L^{\prime}/L$ es una extensi\'on radical, existe $\widetilde{G}\subseteq T(E/L)$
tal que $L^{\prime} = L(\widetilde{G})$. Podemos reemplazar
$\widetilde{G}$ por el subgrupo {\it aditivo} generado por
$\widetilde{G}$ y $L$, que denotamos por $G$. As\'i $L\leq G\leq
T(E/L)$ y $L^{\prime} = L(G)$. Sea
\[
U = \phi(G)=\{f_{\alpha}\mid \alpha\in G\}\leq Z^{1}(\Gamma,\mu(E))
\]
donde $\phi$ es la funci\'on dada en el Corolario
\ref{finitud_redes}. Entonces
\begin{align*}
U^{\bot} &= \{\sigma\in \Gamma\mid \langle\sigma,f_{\alpha}\rangle = 0 \text{ para cada $f_{\alpha}\in U$}\}\\
&=\{\sigma\in \Gamma\mid f_{\alpha}(\sigma) = 0\text{ para cada $f_{\alpha}\in U$}\}\\
&=\{\sigma\in \Gamma\mid  \sigma(\alpha)=\alpha \text{ para cada $f_{\alpha}\in U$}\}\\
&=\{\sigma\in \Gamma\mid \sigma(x) = x \text{ para cada $x\in L(G)$}\}\\
&=\Gal(E/L(G)) = \Gal(E/L^{\prime}).
\end{align*}

Rec\'iprocamente, en caso de que exista
un subgrupo $U\leq Z^{1}(\Gamma,\mu(E))$
tal que $\Gal(E/L^{\prime}) = U^{\bot}$, entonces veamos que
\[
\Gal(E/L^{\prime}) = U^{\bot} = \Gal(E/L(G)),
\]
con $G=\{\alpha\in
E\mid f_{\alpha}\in U\} = \phi^{-1}(U)$ donde $\phi$ es la funci\'on
dada en el Corolario \ref{finitud_redes}.

Para mostrar las igualdades anteriores s\'olo debemos mostrar
$U^{\bot} = \Gal(E/L(G))$. Para este fin consideremos $\tau\in
U^{\bot} = \{\sigma\in\Gamma\mid h(\sigma)=0\text{\ para toda\ } h\in U\}$.
Ahora si $\alpha\in G$ entonces $f_{\alpha}\in U$, en particular,
$f_{\alpha}(\tau) = 0 = \tau(\alpha) - \alpha$. Por lo tanto para
todo $\alpha\in G$, $\tau(\alpha) = \alpha$ y, de este modo, $\tau$
fija a $L(G)$ as\'i que $\tau\in\Gal(E/L(G))$.

Ahora si $\tau\in \Gal(E/L(G))$, sea $h\in U$. Entonces existe un
$\alpha\in G$ tal que $h = f_{\alpha}$, por la definici\'on de $G$ y
el hecho de que $\phi$ es biyectiva. Se sigue que $h(\tau) =
f_{\alpha}(\tau) = 0$ por lo que $\tau\in U^{\bot}$. Ahora por
Teor\'ia de Galois, se tiene que $L^{\prime} = L(G)$
$\fin$
\end{proof}

El siguiente resultado es una aplicaci\'on de la Proposici\'on
\ref{redes_casineat}, ver \cite{BarVel93}. El s\'imbolo
$\sqrt[N]{\alpha}$ denota una ra\'iz del polinomio $u^{N}-\alpha$.

\begin{proposicion}\label{redes2}
Sean $\K/F$ una extensi\'on finita y separable y $E$ la cerradura
normal de $\K/F$. Supongamos que existe una extensi\'on finita $L/F$
tal que

{\rm (1)} $E(\lambda_{N})\cap L = F$ donde $N\in R_{T}$ es un polinomio no
constante.

{\rm (2)} $\K L = L(\sqrt[N]{\alpha})$ para alg\'un $\alpha\in L$ distinto
de $0$.

Entonces $\K=F(\sqrt[N]{\alpha})$.
\end{proposicion}

\begin{proof} 
Consideremos el diagrama siguiente
\[
\xymatrix{E(\lambda_{N})\ar@{-}[d]\ar@{-}[r]& E(\lambda_{N})L\ar@{-}[d]\\
\K\ar@{-}[d]\ar@{-}[r]  & \K L \ar@{-}[d]\\
F \ar@{-}[r] & L  }
\]

Puesto que la extensi\'on $E(\lambda_{N})/F$ es de Galois
se tiene que $E(\lambda_{N})L/L$ es una extensi\'on de
Galois y de la hip\'otesis (1) se tiene
\[
G=\Gal(E(\lambda_{N})/F)\cong \Gal(E(\lambda_{N})L/L)=G_{1}.
\]

Por la hip\'otesis (2) se tiene $\K L = L(\sqrt[N]{\alpha})$. Sea
$\beta = \sqrt[N]{\alpha}$ y consideremos $\sigma\in
\Gal(E(\lambda_{N})L/\K L)$. Entonces
\begin{equation}\label{ec4}
(\sigma(\beta)- \beta)^{N} = \sigma(\beta^{N})-\beta^{N} = \sigma(\alpha)-\alpha = 0.
\end{equation}

Por lo tanto definimos $\chi:G_{1}\rightarrow \mu(E(\lambda_{N})L)$
por $\chi(\sigma) = \sigma(\beta)-\beta$.

Entonces $\Gal(E(\lambda_{N})L/\K L) = \ker(\chi)$, ya que si
$\sigma\in \Gal(E(\lambda_{N})L/\K L)$ se tiene que $\chi(\sigma) =
\sigma(\beta) - \beta = 0$ y rec{\'\i}procamente. Adem\'as, de (\ref{ec4}) se
tiene que $\chi$ toma valores en $\Lambda_{N}$. Puesto que $G$ y
$G_{1}$ son isomorfos, $\chi$ puede ser definido en $G$.

Por lo anterior $\chi$ puede considerarse como un elemento de
$Z^{1}(G,E(\lambda_{N}))$ y $\ker(\chi)$ es igual a
$\Gal(E(\lambda_{N})/\K)$, puesto que $G$ y $G_{1}$ es isomorfo. Por
la Proposici\'on \ref{redes_casineat} $\K/F$ es una extensi\'on radical.
$\fin$
\end{proof}

Sea $E/F$ una extensi\'on finita de Galois, con grupo de Galois $G$.
Sea $L/F$ otra extensi\'on tal que $L\cap E = F$, considere la
composici\'on $EL$. La funci\'on de restricci\'on
\[
\Gal(EL/L)\rightarrow \Gal(E/F),\quad \sigma\mapsto \sigma\mid_{E}
\]
es un isomorfismo de grupos. Denotamos por
$S(L_{1}/L_{2})$ al subconjunto de extensiones de $L_{1}$ contenidas
en $L_{2}$. Entonces las funciones
\begin{gather*}
\varepsilon:S(E/F)\rightarrow S(EL/L),\quad \K^{\prime}/F\mapsto L\K^{\prime}/L
\intertext{y}
\lambda:S(EL/L)\rightarrow S(E/F),\quad \K_{1}/L\mapsto (\K_{1}\cap E)/F
\end{gather*}
son isomorfismo de redes, inversas una de la otra.

Denotamos por $ST(E/F)$ al conjunto de todas las subextensiones
$\K^{\prime}/F$ de $E/F$ que son radicales. Entonces para todo
$\K^{\prime}/F\in ST(E/F)$ existe un $R_{T}$-m\'odulo $G$, no
necesariamente \'unico, tal que $F\subseteq G\subseteq T(E/F)$ y
$\K^{\prime} = F(G)$. Definimos $G_{1} = G + L$. Entonces
$L\K^{\prime} = L(G_{1})$, ya que $L\K^{\prime} = L(G)$, $L\subseteq
G_{1}\subseteq T(EL/L)$ y $G_{1}$ es un $R_{T}$-m\'odulo. Por tanto
$\varepsilon(\K^{\prime}/L)\in ST(EL/L)$. De este modo la
restricci\'on de $\varepsilon$ a las extensiones radicales da lugar
a una funci\'on inyectiva
\begin{gather*}
\rho:ST(E/F)\rightarrow ST(EL/L)
\intertext{definida por}
F(G)/F\mapsto F(G)L/L = L(G + L)/L
\end{gather*}
donde $G$ es un $R_{T}$-m\'odulo tal que $F\subseteq G
\subseteq T(E/F)$.

\begin{proposicion}\label{redes_torsion_galois}
Sea $E/F$ una extensi\'on finita de Galois con grupo de Galois
$\Gamma$ y sea $L/F$ una extensi\'on arbitraria, con $L\subseteq
\overline{K}$, tal que $E\cap L = F$. Si $\mu(EL) = \mu(E)$,
entonces se tiene:
\begin{description}
\item[{\rm (1)}] $(G + L) \cap E = G$ para todo $R_{T}$-m\'odulo $G$ con
$F\subseteq G\subseteq T(E/F)$.
\item[{\rm (2)}] $G_{1} = (G_{1}\cap E) + L$ para todo $R_{T}$-m\'odulo
$G_{1}$, con $L\subseteq G_{1} \subseteq T(EL/L)$.
\item[{\rm (3)}] La funci\'on
\begin{gather*}
\rho:ST(E/F)\rightarrow ST(EL/L)\\
F(G)/F\mapsto L(G +L)/L, \,\,\, F\leq G\leq T(E/F)\\
\intertext{es biyectiva, y la funci\'on}
ST(EL/L)\rightarrow ST(E/F),\\
L(G_{1})/L\mapsto F(G_{1}\cap E)/F, \,\,\, L\leq G_{1} \leq T(EL/L)
\end{gather*}
es su inversa.
\end{description}
Aqu\'i, la notaci\'on $F\leq G$ indica que $F$ es un
subm\'odulo del $R_{T}$-m\'odulo $G$.
\end{proposicion}

\begin{proof} 
(1) Sea $w\in (G + L)\cap E$ as\'i $w = x + y$ donde $x\in G$ e
$y\in L$, por lo que $y = w - x\in E$. As\'i $y\in F$ ya que $E\cap
L = F$. Por lo tanto $w \in G$.

Rec\'iprocamente si $x\in G$ ciertamente $x\in (G + L)\cap E$.

(2) Denotamos por $\Gamma_{1}$ al grupo de Galois de $EL/L$. Hemos
visto anteriormente que se tiene un isomorfismo de grupos
\begin{equation}\label{ec3}
\Gamma_{1}\rightarrow \Gamma,\quad \sigma_{1}\rightarrow
\sigma_{1}\!\!\mid_{E}
\end{equation}

Como $\mu(EL) = \mu(E)$, el isomorfismo anterior induce un
isomorfismo de grupos
\[
\upsilon:Z^{1}(\Gamma,\mu(E))\rightarrow Z^{1}(\Gamma_{1},\mu(EL))
\]
dado como sigue: sea $h\in Z^{1}(\Gamma,\mu(E))$. Si
$\sigma_{1}\in \Gamma_{1}$ se tiene que $\sigma_{1}\!\!\mid_{E}\in
\Gamma$, y definimos $\upsilon(h)(\sigma_{1}) :=
h(\sigma_{1}\mid_{E})$. Ahora
\[
\upsilon(h)(\sigma_{1}\circ
\sigma_{2}) = h(\sigma_{1}\circ \sigma_{2}\mid_{E}) =
h(\sigma_{1}\mid_{E}\circ \sigma_{2}\mid_{E}).
\]

As\'i $\upsilon(h)$ es un homomorfismo cruzado. Por construcci\'on
$\upsilon$ es un homomorfismo de grupos, y por (\ref{ec3}) se tiene
el que $\upsilon$ es un isomorfismo de grupos.

Sea $G_{1}$ con $L\leq G_{1} \leq T(EL/L)$. Ahora si $w\in
(G_{1}\cap E) + L$ entonces $w = x + y$ con $x\in (G_{1}\cap E)$ e
$y\in L$, as\'i $w\in G_{1}$. Ahora sea $a_{1}\in G_{1}$. Entonces
$f_{a_{1}}\in Z^{1}(\Gamma_{1},\mu(EL)$. Tenemos que existe $f\in
Z^{1}(\Gamma,\mu(E))$ tal que $f_{a_{1}} = \upsilon(f)$. De la
Proposici\'on \ref{finitud_TC/L}, existe $a\in T(E/F)$ tal que
$f_{a_{1}} = \upsilon(f = f_{a})$.

Se tiene que $f_{a_{1}}(\sigma_{1}) = f_{a}(\sigma_{1}\mid_{E})$
para todo $\sigma_{1}\in \Gamma_{1}$. De aqu\'i se sigue que
$\sigma_{1}(a_{1}) - a_{1} = \sigma_{1}(a) - a$, es decir,
$\sigma_{1}(a_{1} - a) = a_{1} - a$ para cada $\sigma_{1}\in
\Gamma_{1}$. As\'i pues $a_{1} - a\in L$. Por lo tanto $a_{1} = a +
b$ donde $b\in L$. Puesto que $a\in G_{1}$ se sigue que $a\in
(G_{1}\cap E) + L$.

(3) Por la observaci\'on hecha previamente a esta proposici\'on se
tiene que $\rho$ es inyectiva, por lo que basta mostrar que $\rho$
es suprayectiva. Para ello, sea $\K_{1}/L\in ST(EL/L)$. Entonces
$\K_{1} = L(G_{1})$ para alg\'un $G_{1}$ con $L\leq G_{1}\leq
T(EL/L)$. Por lo tanto si ponemos $G = G_{1}\cap E$, obtenemos que
$F(G)/F\in ST(E/F)$ y que
\[
\rho(F(G)/F) = L(F(G))/L = L(F(G_{1}\cap E))/L = L(L + (G_{1} \cap E))/L.
\]
Por (2), se tiene que $L(L + (G_{1} \cap E))=L(L+G) =
L(G_{1}) = \K_{1}$.
$\fin$
\end{proof}

Por otro lado, el rec{\'\i}proco del Teorema \ref{tdim_p} no siempre es
v\'alido como lo muestra el siguiente lema.
\begin{lema}\label{auxiliar_contra_cogalois}
Sea $L/\K$ una extensi\'on Galois tal que $[L:\K]=p^{2}$,
$\mu(L)=\mu(\K)$ y $G = \Gal(L/\K)\cong C_{p^{2}}$. Entonces $L/\K$ no
es una extensi\'on radical.
\end{lema}

\begin{proof} 
Supongamos que la extensi\'on $L/\K$ es radical. Considere el grupo
de cohomolog\'ia $H^{1}(G,\mu(L))$. 
Puesto que $\mu(L)=\mu(\K)$ se tiene que
$B^{1}(G,\mu(L))=\{1\}$ (ver Lema \ref{auxiliar_coho_cogalois}).
Por tanto, $H^{1}(G,\mu(L))=Z^{1}(G,\mu(L))/B^{1}(G,\mu(L))\cong
\Hom(G,\mu(L))$. En consecuencia, por la Proposici\'on
\ref{finitud_TC/L}, se tendr\'a
\[
\cog(L/\K)\cong \Hom(G,\mu(\K)).
\]
Consideremos un elemento de orden $p$, digamos
$\tau$, en $G$. Sea $H=(\tau)$ y $L^{\prime}= L^{H}$. Notemos que,
al ser $H$ normal en $G$, se tiene que $L^{\prime}/\K$ es una
extensi\'on normal y, por lo tanto, de Galois. Adem\'as
$G^{\prime}=\Gal(L^{\prime}/\K)$ es isomorfo a $C_{p}$. N\'otese que
$\mu(L^{\prime})=\mu(\K)$. As\'i
\[
\cog(L^{\prime}/\K)\cong \Hom(G^{\prime},\mu(\K)).
\]
Notemos que la cardinalidad de $\cog(L/\K)\cong 
\Hom(G,\mu(\K))$ es $\mid \mu(\K)\mid$.
Para ver esto sea $a\in G\cong
C_{p^{2}}$ un generador. Un homomorfismo $\psi:G\rightarrow \mu(\K)$
queda completamente determinado por su acci\'on en $a$. Por lo tanto
hay $\mid \mu(\K)\mid$ homomorfismos de $G$ en $\mu(\K)$. Del mismo
modo podemos mostrar que la cardinalidad de
$\cog(L^{\prime}/\K)\cong \Hom(G^{\prime},\mu(\K))$ es $\mid
\mu(\K)\mid$.

Por otro lado tenemos que $\cog(L^{\prime}/\K)\subseteq
\cog(L/\K)$, ver Proposici\'on \ref{prop_cog}. Entonces,
como ambos grupos tienen la misma
cardinalidad, se tiene $\cog(L^{\prime}/\K)= \cog(L/\K)$.
Se sigue que $L=L^{\prime}$, ya que si $\alpha_{1},
\ldots , \alpha_{s}$ generan a $L$ sobre $\K$, entonces por lo
mostrado se tendr\'a que $\alpha_{1}, \ldots , \alpha_{s}\in
L^{\prime}$ y de aqu\'i la afirmaci\'on. Se tendr\'ia que
$[L:\K]=p^{2}=[L^{\prime}:\K]=p$ lo cual es una contradicci\'on.
$\fin$
\end{proof}

El siguiente ejemplo muestra que la propiedad de ser extensi\'on
radical no es hereditaria.

\begin{ejemplo}\label{noredes_cogalois}
Sea $M=P^{n}$, $n\in {\mathbb{N}}$ y $P\in R_{T}$ irreducible, se
considera la extensi\'on $K(\Lambda_{M})/K(\lambda_{P})$. Sea $t\in
{\mathbb{N}}$ de tal modo que $p^{t-1}< n \leq p^{t}$ y $n_{0}$ la
parte entera de $\frac{n}{p^{t-1}}$.

Del Corolario 1 de \cite{LamVil2001}, se obtiene
\[
H_{M}\cong ({\mathbb{Z}}/p^{t}{\mathbb{Z}})^{\alpha}\times
 {\mathbb{Z}}/p^{n_{1}}{\mathbb{Z}}\times \cdots \times {\mathbb{Z}}/p^{n_{s}}{\mathbb{Z}}
\]
con $t> n_{1}\geq \cdots \geq n_{s}\geq 0$. Aqu{\'\i} $H_M$ es 
el grupo de Galois de la extensi\'on $k(\Lambda_M)/k(\Lambda_P)$.

Sea $n=5$ y $p=3$. Si $t= 2$ se cumple que $p^{t-1}< n \leq p^{t}$.
Adem\'as $n_{0}=1$. El valor de $\alpha$ esta dado por el Corolario
1 de \cite{LamVil2001}.

Se puede escoger un subgrupo de $H_{M}$ de la forma
\[
H_M=({\mathbb{Z}}/p^{t}{\mathbb{Z}})^{\alpha-1}\times 
{\mathbb{Z}}/p^{n_{1}}{\mathbb{Z}}\times \cdots \times
 {\mathbb{Z}}/p^{n_{s}}{\mathbb{Z}}.
\]
Sea $L^{\prime}=L^{H}$, as\'i
$\Gal(L^{\prime}/K(\lambda_{P}))\cong C_{p^{2}}$. Se tiene tambi\'en
que $\mu(K(\lambda_{P}))=\mu(L^{\prime})$, esto es posible,
escogiendo adecuadamente $q= p^{\nu}$.

Por el Lema \ref{auxiliar_contra_cogalois},
$L^{\prime}/K(\lambda_{P})$ no es radical. Por lo tanto
$K(\lambda_{P^{5}})/K(\lambda_{P})$ es una extensi\'on Galois
radical ciclot\'omica, pero no cumple la propiedad de que si $L$ es
un campo tal que $K(\lambda_{P})\subseteq L \subseteq
K(\lambda_{P^{5}})$, entonces $L/K(\lambda_{P})$ es 
una extensi\'on radical.
\end{ejemplo}

\begin{ejemplo}\label{ejemplo_schultheis}
En este ejemplo tendremos $q=p\geq 3$.
Considere la extensi\'on $L/K(\Lambda_{T})$, donde $L$ es el campo
de descomposici\'on del polinomio $f(X) = X^{T} - 1$, con
coeficientes en $K(\Lambda_{T})$. El grado de esta extensi\'on es
$[L:K(\Lambda_{T})] = p$, ver Ejemplo \ref{ejemplo4}. Trataremos de
determinar la estructura de $\cog(L/K(\Lambda_{T}))$.

Supongamos que $\overline{\beta}\in \cog(L/K(\Lambda_{T}))$
tiene orden $Q^{r}$, con $Q$ m\'onico irreducible, $r\geq 1$ y
$Q\neq T$. Por la Proposici\'on \ref{finitud_RC} se tiene que
$\lambda_{Q^{r}}\in L$. Puesto que
$\lambda_{Q}=\lambda^{Q^{r-1}}_{Q^{r}}\in L$, por pureza se tiene
que $\lambda_{Q}\in K(\Lambda_{T})$, pero esto implica que $Q=T$ lo
cual es una contradicci\'on. Por lo tanto
\[
\cog(L/K(\Lambda_{T}))\cong \cog(L/K(\Lambda_{T}))_{T}
\]
donde $\cog(L/K(\Lambda_{T}))_{T}$ es el conjunto de
elementos de $\cog(L/K(\Lambda_{T}))$ cuyo orden es una
potencia de $T$.

Necesitaremos un lema para obtener la cardinalidad de
$\cog(L/K(\Lambda_{T}))$. Para empezar sea $z\in K$, $z\neq 0$,
y $N\in R_{T}$ un polinomio no constante. Consideremos $g(X) = X^{N} -
z\in K(\Lambda_{N})[X]$. El campo de descomposici\'on de $g(X)$,
sobre $K$, es de la forma $\K=K(\alpha,\lambda_{N})$ donde $\alpha$
es una ra\'iz arbitraria de $g(X)$ y $\lambda_{N}$ un generador de
$\Lambda_{N}$. Como el polinomio $g(X)$ es separable, la extensi\'on
$\K/K$ es de Galois.

Sea $G = \Gal(L/K)$ entonces dado $\sigma\in G$ se tiene que
$\sigma(\alpha) = \alpha + \lambda^{M_{\sigma}}$ y $\sigma(\lambda)
= \lambda^{N_{\sigma}}$, donde $M_{\sigma}$ y $N_{\sigma}$ se
determinan salvo un m\'ultiplo de $N$, y $N_{\sigma}$ es primo
relativo a $N$.

Por otro lado considere $G(N)$ el subgrupo de $GL_{2}(R_{T}/(N))$ de
todas las matrices de la forma
\begin{displaymath}
\left(\begin{array}{ccc} 1 & 0  \\
\overline{B} & \overline{A}
\end{array}\right)
\end{displaymath}
donde $\overline{B}\in R_{T}/(N)$ y $\overline{A}\in
(R_{T}/(N))^{*}$. De esta descripci\'on se sigue que $
\card(G(N)) = q^{\deg(N)} \Phi(N)$. Sea
$\theta:G\rightarrow G(N)$ definida por:
\begin{displaymath}
\theta(\sigma) = \left(\begin{array}{ccc} 1 & 0  \\
\overline{M_{\sigma}} & \overline{N_{\sigma}}
\end{array}\right).
\end{displaymath}

Tenemos el siguiente lema.

\begin{lema}\label{auxiliar_noabeliano_galois}
Sea $L/K$ la extensi\'on anteriormente descrita y $\theta$ la
funci\'on anteriormente definida. Entonces $\theta$ es un
monomorfismo de grupos. Por otra parte si $N=P$, $P$ m\'onico e
irreducible, $z\in R_{T}$ como antes y la ecuaci\'on $g(X)=0$ no
tiene soluciones en $R_{T}$, entonces $\theta$ es un isomorfismo de
grupos.
\end{lema}

\begin{proof} 
Sean $\sigma,\tau\in G$. Se tiene que $\sigma(\tau(\alpha)) =
\sigma(\alpha + \lambda^{M_{\tau}}) = \alpha + \lambda^{M_{\sigma}}
+ \lambda^{M_{\tau}N_{\sigma}}$ adem\'as $\sigma(\tau(\lambda)) =
\sigma(\lambda^{N_{\tau}}) = \lambda^{N_{\sigma} N_{\tau}}$, por lo
tanto
\[
\theta(\sigma\cdot\tau) = \left(\begin{array}{ccc} 1 & 0  \\
\overline{M_{\sigma} + M_{\tau}N_{\sigma}} & \overline{N_{\sigma} N_{\tau}}
\end{array}\right)= \left(\begin{array}{ccc} 1 & 0  \\
\overline{M_{\sigma}} & \overline{N_{\sigma}}
\end{array}\right) \left(\begin{array}{ccc} 1 & 0  \\
\overline{M_{\tau}} & \overline{N_{\tau}}
\end{array}\right) =\theta(\sigma) \theta(\tau).
\]

Por lo tanto $\theta$ es un homomorfismo de grupos. Si
$\theta(\sigma)$ es la matriz identidad se tiene que $M_{\sigma}$ es
un m\'ultiplo de $N$ y que $N_{\sigma} = 1  + N Q$. As\'i $\sigma =
e$, es decir, $\theta$ es un monomorfismo de grupos.

Si $N=P$, donde $P$ es un polinomio
m\'onico e irreducible, $z\in R_{T}$ como antes y la
ecuaci\'on $g(X)=0$ no tiene soluciones en $R_{T}$, entonces por el
Teorema 1.7 (3) de \cite{Hsu97}, se tiene que $\Gal(L/K(\lambda_{P}))$
tiene cardinalidad $q^{\deg(P)}$, es decir, el monomorfismo
anterior es un isomorfismo.
$\fin$
\end{proof}

Regresando a nuestro ejemplo,
se mostrar\'a que $\mu(L)=\Lambda_{T}$. Para empezar, ciertamente
$\Lambda_{T}=\mu(K(\Lambda_{T}))\subseteq\mu(L)$. Por otro lado sea
$u\in \mu(L)$ no nulo. Existe un $N\in R_{T}$ tal que
$u^{N}=0$. Por lo tanto $u$ es de la forma $\lambda^{M}_{N}$.
Podemos suponer que $(M,N)=1$. Por tanto, por la Proposici\'on
\ref{P6.2.11}, podemos afirmar que $\lambda_{N}\in L$.
Sea $N=P^{\alpha_{1}}_{1}\cdots P^{\alpha_{s}}_{s}$. Entonces
$\lambda_{P_{i}}=\lambda^{P^{\alpha_{1}}_{1}\cdots
P^{\alpha_{i}-1}_{i}\cdots P^{\alpha_{s}}_{s}}_{N}\in L$. Puesto
que $L/K(\Lambda_T)$ es pura,
tendremos que $\lambda_{P_{i}}\in K(\Lambda_{T})$. As\'i
$P_{i}=T$. Por lo tanto $N=T^{n}$ con $n\geq 1$ y $n\in
{\mathbb{N}}$.

Supongamos que $n\geq 2$ y considere el diagrama
\[
\xymatrix{& L\ar@{-}[dl]\ar@{-}[dr]  &\\
K(u)\ar@{-}[dr] & & K(\Lambda_{T})\ar@{-}[dl] \\
& K& }
\]

Del diagrama anterior obtenemos $[L:K(u)]\Phi(T^{n})=p(p-1)$, puesto
que $\Phi(T^{n})=p^{n-1}(p-1)$, se sigue que $[L:K(u)]p^{n-1}=p$. Si
$n\geq 3$ entonces $n-2\geq 1$, as\'i $[L:K(u)]p^{n-2}=1$ lo cual es
una contradicci\'on. Solo resta considerar el caso $n=2$, que
implica que $L=K(u)$, pero del Lema \ref{auxiliar_noabeliano_galois}
se tiene que $\Gal(L/K)$ es un grupo no abeliano, lo cual contradice que el
grupo $\Gal(K(\Lambda_{T^{2}})/K)$ es abeliano. Por lo tanto $n=1$ y
$u=\lambda^{M}_{T}\in K(\Lambda_{T})$.

Por el Lema \ref{auxiliar_coho_cogalois} se tiene
que $B^{1}(G,\mu(L))=\{0\}$. Se sigue
$H^{1}(G,\mu(L))=Z^{1}(G,\mu(L))/B^{1}(G,\mu(L))\cong 
\Hom(G,\mu(L))$.  De esta manera, utilizando la
demostraci\'on del Lema \ref{auxiliar_contra_cogalois}, se tiene que
\[
\mid\cog(L/K(\Lambda_{T}))\mid=[L:K(\Lambda_{T})]=p.
\]
\end{ejemplo}

\begin{ejemplo}\label{ejemplo_entre_ciclotomicos}
Consideremos la extensi\'on $K(\Lambda_{P^{n}})/K(\Lambda_{P})$.
Calcularemos la cardinalidad del m\'odulo
$\cog(K(\Lambda_{P^{n}})/K(\Lambda_{P}))$ en el siguiente caso:
Sea $P(T) = T$, $q = p
> 2$ y $n = 2$. Sea $H_{T^2}=\{\overline{N}\in R_T/(T^2)\mid
(N,T^2)=1\text{\ y\ } N\equiv 1\bmod T\}$.
Entonces $\card(H_{T^{2}}) = q^{d (n-1)} = p$,
con $d = \deg(P(T)) = 1$. En particular el grupo
$H_{T^{2}}$ es c\'iclico. Se tiene que
\begin{gather*}
H^{1}(H_{T^{2}},\Lambda_{T^{2}}) \cong \ker(N_{H_{T^{2}}})/ D\Lambda_{T^{2}}\\
\intertext{donde definimos $N_{H_{T^{2}}}:\Lambda_{T^{2}}\rightarrow
\Lambda_{T^{2}}$ y $D:\Lambda_{T^{2}}\rightarrow \Lambda_{T^{2}}$
como}
N_{H_{T^{2}}}(x) = x + \sigma\cdot x + \cdots + \sigma^{p-1}\cdot x,\\
D(x) = \sigma\cdot x - x,
\end{gather*}
donde $\sigma = 1 + T + (T^{2})$ es un generador de $H_{T^{2}}$ y
$x\in \Lambda_{T^{2}}$. Por otro lado si $
x=\lambda_{T^2}^M$ se tiene que
\begin{gather*}
N_{H_{T^{2}}}(x) = \lambda_{T^2}^M + \lambda_{T^2}^{M(1+T)}+\cdots
\lambda_{T^2}^{M(1+(p-1)T)}=\lambda_{T^2}^{pM+(1+2+\cdots+p-1)MT}=0.
\end{gather*}

Notemos que $1 + 2 + \cdots + (p-1) = 0$ ya que tal suma es igual
a $\frac{p(p-1)}{2}=0$ en ${\mathbb{F}}_{p}$. De esta
manera se tiene que $\ker(N_{H_{T^{2}}}) = \Lambda_{T^{2}}$.
Observemos tambi\'en que $D(x) = \lambda_{T^2}^{M(1+T)}
-\lambda_{T^2}^M=\lambda_{T}^M$. Por tanto
\[
D \Lambda_{T^{2}} = \Lambda_{T}.
\]
Se sigue que $H^{1}(H_{T^{2}},\Lambda_{T^{2}}) =
\Lambda_{T^{2}}/\Lambda_{T}$. Por otra parte del Lema
\ref{auxiliar_coho_cogalois} se tiene que $
\card(B^{1}(H_{T^{2}},\Lambda_{T^{2}})) = 
\card(\Lambda_{T^{2}}/\Lambda_{T})$ y recordando que
\begin{gather*}
H^{1}(H_{T^{2}},\Lambda_{T^{2}}) = Z^{1}(H_{T^{2}},\Lambda_{T^{2}})/B^{1}(H_{T^{2}},\Lambda_{T^{2}})\\
\intertext{se obtiene, por la Proposici\'on \ref{finitud_TC/L},}
\mid(\cog(\K(\Lambda_{T^{2}})/\K(\Lambda_{T})))\mid = \mid(Z^{1}(H_{T^{2}},\Lambda_{T^{2}}))\mid = [\K(\Lambda_{T^{2}}):\K(\Lambda_{T})]^{2}.
\end{gather*}
\end{ejemplo}

El siguiente lema, muestra que ciertas extensiones tienen
propiedades an\'alogas a las enunciadas en el Lema 1.3 de
\cite{GreHar86}, concretamente los pasos 1 y 2. Sin embargo veremos
despu\'es que, en general, estas propiedades no se cumplen.

\begin{lema}\label{clave_cogalois}
Considere la extensi\'on $L/K(\lambda_{P})$, donde $L$ es campo de
descomposici\'on del polinomio $X^{P}-a$, donde $P\in R_{T}$ es
irreducible y $a\in K(\lambda_{P})\setminus K(\lambda_{P})^{P}$. El
m\'odulo $\cog(L/K(\lambda_{P}))$ no tiene elementos de orden
$Q$, donde $Q$ es un polinomio irreducible, distinto de $P$.
Adem\'as si $\nu_{{\mathfrak{p}}}(a)\geq q^{d}$, donde $d =
\deg(P)$, se tiene que $\cog(L/\K)$ no tiene elementos de
orden $P^{2}$.
\end{lema}

\begin{proof} 
Supongamos que $\cog(L/K(\lambda_{P}))$ tiene un elemento de
orden $Q$. Entonces como $L/K(\lambda_{P})$ es Galois,
se tiene que $\lambda_{Q}\in L$ y como $L/K(\lambda_{P})$ es radical
ciclot\'omica, se tendr\'a que $\lambda_{Q}\in \mu(\K)= \Lambda_{P}$,
por la Proposici\'on \ref{finitud_RC}, por lo tanto $Q = P$.

Ahora supongamos que $\cog(L/K(\lambda_{P}))$ tiene un elemento
de orden $P^{2}$, es decir, existe $\widehat{\beta}\in 
\cog(L/\K)$ tal que $\beta^{P^{2}} = b \in \K$. Entonces, como
$L/K(\lambda_{P})$ es radical, se sigue de la Proposici\'on
\ref{galois_torsion_raices} que $\lambda_{P^{2}}\in L$. Consideremos
el siguiente diagrama
\[
\xymatrix{{\mathcal{O}}_{L}\ar@{-}[d]\ar@{-}[r] & L\ar@{-}[d] \\
{\mathcal{O}}_{K(\lambda_{P^{2}})}\ar@{-}[d]\ar@{-}[r]& K(\lambda_{P^{2}})\ar@{-}[d]\\
{\mathcal{O}}_{K(\lambda_P)}\ar@{-}[d]\ar@{-}[r]& K(\lambda_P)\ar@{-}[d]\\
R_{T}\ar@{-}[r] & K }
\]

El \'indice de ramificaci\'on del primo $P$ en la extensi\'on
$K(\lambda_{P^{2}})/K$, es $\Phi(P^{2})$ por lo que el \'indice de
ramificaci\'on de $P$ en la extensi\'on $L/K$ es $\widetilde{d}
\Phi(P^{2})$, donde $\widetilde{d} = e_{L/\K(\lambda_{P^{2}})}
(P)$. Del
Teorema 3.9. de \cite{Sch90} tenemos que
el \'indice de ramificaci\'on es
$\Phi(P)$. En otras palabras, $d \Phi(P^{2}) = \Phi(P)$, lo cual es
absurdo, de donde se sigue la afirmaci\'on.
$\fin$
\end{proof}

\begin{ejemplo}\label{ejemplo6_1}
Sean $P,Q\in R_{T}$, irreducibles y distintos. Considere la
extensi\'on $L=K(\Lambda_{P^{2}Q^{2}})/K$. Notemos que
$L=K(\lambda_{P^{2}},\lambda_{Q^{2}})$. Sea $\sigma=1+PQ\in
G=\Gal(L/K)$. Se tiene que $\sigma\neq 1$ ya que
$\lambda^{1+PQ}_{P^{2}Q^{2}}=\lambda_{P^{2}Q^{2}}+\lambda_{PQ}\neq
\lambda_{P^{2}Q^{2}}$.

Se tiene que
$\sigma(\lambda_{PQ})=\lambda^{1+PQ}_{PQ}=\lambda_{PQ}$. Por lo
tanto si $\K$ es el campo fijo de $(\sigma)$, se tiene que
$\lambda_{PQ}\in \K$.

Por otro lado se tiene que $\sigma^{p}=(1+PQ)^{p}=1+P^{p}Q^{p}\equiv
1 \bmod\, P^{2}Q^{2}$, es decir, el orden de $\sigma$ es $p$. Por lo
tanto $[L:\K]=p$.

Puesto que
$\sigma(\lambda_{P^{2}})= \lambda_{P^{2}}+\lambda^{Q}_{P}\neq
\lambda_{P^{2}}$, se tiene que $\lambda_{P^2}\notin
\K$. De modo an\'alogo se puede mostrar que
$\beta=\lambda_{Q^{2}}\notin \K$.

Ahora, como $[L:\K]=p$, se tiene 
que $L=\K(\alpha)=\K(\beta)$, adem\'as
$\alpha^{P}=\lambda_{P}$ y $\beta^{Q}=\lambda_{Q}$. Por lo tanto el
m\'odulo $\cog(L/\K)$, tiene elementos de orden $P$ y de orden
$Q$.
\end{ejemplo}

\begin{ejemplo}\label{ejemplo7_1} Sea $q=p^{\nu}$ con $p\geq 3$.
Sea $L=K(\Lambda_{P^{3}})$ y $\sigma=1+P\in \Gal(L/K(\Lambda_{P}))$.
Se tiene que $\sigma^{p}=(1+P)^{p}\equiv 1\bmod\, P^{3}$.
Adem\'as $\sigma\neq 1$ ya
que $\sigma(\lambda_{P^{3}})= \lambda_{P^{3}}+\lambda_{P^{2}}\neq
\lambda_{P^{3}}$.

Sea $\K= L^{(\sigma)}$, se tiene que $[L:\K]=p$. Por otra parte
$\sigma(\lambda_{P^{2}})= \lambda_{P^{2}}+\lambda_{P}\neq
\lambda_{P^{2}}$. Por lo tanto $\alpha=\lambda_{P^{2}}\notin \K$.
As\'i $L=\K(\alpha)$ y $\alpha^{P}= a\in \K$.

Ahora bien, $\lambda_{P^{3}}\in L$ tiene orden $P^{2}$ ya que
$\lambda^{P^{2}}_{P^{3}}\in \K$ y $\lambda^{P}_{P^{3}}\notin \K$. Por
lo tanto el m\'odulo $\cog(L/\K)$ tiene elementos de orden
$P^{2}$.
\end{ejemplo}

Los Ejemplos \ref{ejemplo6_1} y \ref{ejemplo7_1} 
muestran que no tenemos los an\'alogos
del Lema 1.3. de \cite{GreHar86}, a saber si $L/\K$ es
cogalois, en el sentido cl\'asico, y es tal que $[L:\K]=p$, con
$L=\K(\alpha)$ con $\alpha^{p}=a\in \K$ y $L/\K$ es separable y pura,
entonces

(a) El grupo $\cog(L/\K)$ no tiene elementos de orden $q\neq p$,
$q$ un n\'umero primo.

(b) El grupo $\cog(L/\K)$ no tiene elementos de orden $p^{2}$.

\section{Una cota para $\mid\cog(L/\K)\mid$}\label{estimacion_para_cogalois}

En esta secci\'on establecemos una cota superior para la
cardinalidad del m\'odulo $\cog(L/\K)$. En lo que sigue sea
$q=p^{\nu}$ y sea $L/\K$ una extensi\'on
radical.

\begin{observacion}\label{observacion_elemental_abeliano}
Si la extensi\'on $L/\K$ es de Galois y radical ciclot\'omica y
tal que $\mu(\K)=\mu(L)$,
entonces al ser radical $L$ es de la forma $\K(\rho_{1},\ldots
,\rho_{t})$, con $\rho_i^{M_{i}}=a_{i}\in \K$, para algunos $M_{i}\in
R_{T}$. Por otra parte las ra\'ices de $X^{M_{i}}-a_{i}$ son
$\{\rho_{i}+\lambda^{A}_{M_{i}}\}_{A\in R_{T}}$. Por lo tanto
$\Gal(\K(\rho_{i})/\K)\subseteq \Lambda_{M_{i}}$. De esta manera
$\Gal(\K(\rho_{i})/\K)$ es un $p$--grupo
elemental abeliano.

Puesto que se tiene una inyecci\'on
\[
\Gal(L/\K)\hookrightarrow \prod^{t}_{i=1}
\Gal(\K(\rho_{i})/\K)
\]
se sigue que $\Gal(L/\K)$ es un $p$--grupo elemental abeliano.

Si $L/K$ es una extensi\'on finita y $\mu(L)=\Lambda_M$,
definimos
\[
\deg(\mu(L))=\deg(M).
\]
\end{observacion}

\begin{lema}\label{auxiliar_coho_cogalois}
Sea $L/\K$ una extensi\'on finita de Galois radical ciclot\'omica.
Entonces existe un isomorfismo
\[
B^{1}(G,\mu(L))\cong \mu(L)/\mu(\K)
\]
como $R_T$--m\'odulos.
\end{lema}

\begin{proof} 
Se define $\psi\colon \mu(L)\rightarrow B^{1}(G,\mu(L))$ como sigue:
$\psi(u)=f_{u}$. Observemos que
$\psi(u+v)=\psi(u)+\psi(v)$ ya que
\begin{align*}
(f_{u}+f_{v})(\sigma) &= f_{u}(\sigma)+f_{v}(\sigma)= \sigma(u)-u+\sigma(v)-v\\
&=\sigma(u+v)-(u-v) =(f_{u+v})(\sigma).
\end{align*}
y si $M\in R_{T}$ se tendr\'a que $\psi(u^{M})=f_{u^M}=f^{M}_{u}$,
ya que
\begin{align*}
f_{u^{M}}(\sigma) &= \sigma(u^{M})-u^{M} =(\sigma(u))^{M}-u^{M}\\
&=(\sigma(u)-u)^{M} =(f_{u}(\sigma))^{M}.
\end{align*}

Por tanto $\psi$ es un homomorfismo de $R_{T}$-m\'odulos,
suprayectivo por la definici\'on de $B^{1}(G,\mu(L))$. Puesto que
$L/\K$ es una extensi\'on de Galois
tenemos $\ker(\psi)=\mu(\K)$.
$\fin$
\end{proof}

\begin{proposicion}\label{acotacion1}
Sea $L/\K$ una extensi\'on Galois y radical ciclot\'omica. Supongamos
que $\mu(L)=\mu(\K)$. Entonces
\[
\mid\cog(L/\K)\mid=q^{m\deg(\mu(L))},
\]
donde $[L:\K]=p^m$.
\end{proposicion}

\begin{proof} 
Por la Observaci\'on \ref{observacion_elemental_abeliano} se tiene que
$\Gal(L/\K)\cong C^{m}_{p}$, para alg\'un $m\in {\mathbb{N}}$. Puesto
que $B^{1}(G,\mu(L))=\{0\}$ y $H^{1}(G,\mu(L))\cong 
\Hom(G,\mu(L))$, entonces de la Proposici\'on \ref{finitud_TC/L}, se
tiene que
\[
\cog(L/\K)\cong Z^{1}(G,\mu(L))/B^{1}(G,\mu(L))
\cong H^{1}(G,\mu(L))\cong \Hom(G,\mu(L)).
\]
Adem\'as $\mu(L) \cong C^{\nu\deg(\mu(L))}_{p}$. Por
tanto, si denotamos por ${\eu L}_p({\ma F}_p^m,{\ma F}_p^{\nu
\deg(\mu(L))})$ al conjunto de las transformaciones lineales de 
${\ma F}_p^m$ a ${\ma F}_p^{\nu\deg(\mu(L))}$ y 
al conjunto de las
matrices $m\times \nu\deg(\mu(L))$ con coeficientes en ${\ma F}_p$
lo denotamos por
${\eu M}_{m\times \nu\deg(\mu(L))}({\ma F}_p)$,
se tiene
\begin{align*}
\Hom(G,\mu(L)) &= \Hom(C^{m}_{p},C^{\nu\deg(\mu(L))}_{p})
={\mathfrak{L}}_{p}({\mathbb{F}}^{m}_{p},
{\mathbb{F}}^{\nu\deg(\mu(L))}_{p})\\
&={\mathfrak{M}}_{m\times \nu\deg(\mu(L))}({\mathbb{F}}_{p}).
\end{align*}

Por tanto $\mid \Hom(G,\mu(L))\mid=p^{ms
\deg(\mu(L))}= q^{m\deg(\mu(L))}$.
$\fin$
\end{proof}

\begin{ejemplo}\label{ejemplo_schultheis_12}
Del Ejemplo \ref{ejemplo_schultheis}, se sigue que la extensi\'on
$L/K(\Lambda_{T})$, donde $L$ es el campo de descomposici\'on del
polinomio $f(X) = X^{T} - 1$, cumple que
$\mid\cog(L/K(\Lambda_{T}))\mid=[L:K(\Lambda_{T})]=q=q^{m\deg(\mu(L))}$,
en concordancia con la Proposici\'on \ref{acotacion1}.
\end{ejemplo}

\begin{proposicion}\label{acotacion2}
Sea $L/\K$ una extensi\'on Galois y radical ciclot\'omica y
supongamos que $L=\K(\mu(L))$. Entonces $\mid\cog(L/\K)\mid\leq
q^{m\deg(\mu(L))}$ para alguna $m\in{\ma N}$.
\end{proposicion}

\begin{proof} Por el Corolario \ref{cogalois_grado_p_pureza} se tiene que
$[L:\K]=p^{m}$ para alguna $m\in
{\mathbb{N}}$. Ahora la
demostraci\'on es por inducci\'on sobre $m$. Sea $L/\K$ una
extensi\'on Galois, radical ciclot\'omica, tal que $L=\K(\mu(L))$ y
$[L:\K]=p$. Por lo tanto $L/\K$ es c\'iclica de grado $p$. Sean
$M=P^{\alpha_{1}}_{1}\cdots P^{\alpha_{r}}_{r}$ y
$N=P^{\beta_{1}}_{1}\cdots P^{\beta_{r}}_{r}$, con $1\leq \beta_{i}
\leq \alpha_{i}$, donde $i=1,\ldots , r$, tales que $\mu(L)=\Lambda_M$
y $\mu(\K)=\Lambda_N$.

Sea $G:=\Gal(L/\K)
=(\sigma)$. Se tiene $\sigma(\lambda_{M})=\lambda^{A}_{M}$ ya
que la acci\'on de Carlitz-Hayes conmuta con $\sigma$.
Notemos que $\sigma(\lambda_{M})\neq\lambda_{M}$, en caso contrario
esto implicar\'ia que $\lambda_{M}\in \K$, es decir, 
$L=\K$ lo cual es
una contradicci\'on. Por lo tanto $M\nmid (A-1)$. Sea $M=ND$. Por lo tanto
$\lambda_{N}=\lambda^{D}_{N}$. Se tiene que
\[
\lambda_{N}=\sigma(\lambda_{N})=\sigma(\lambda^{D}_{M})=
(\sigma(\lambda_{M}))^{D}=\lambda^{AD}_{M}=\lambda^{A}_{N}.
\]
Se sigue que $\lambda^{A-1}_{N}=0$, es decir, $N\mid (A-1)$.

Por otro lado
\begin{align*}
\Tr_{G}(\lambda_{M}) &=
\lambda_{M}+\lambda^{A}_{M}+\lambda^{A^{2}}_{M}+\cdots
+\lambda^{A^{p-1}}_{M}\\
&=\lambda^{1+A+A^{2}+\cdots +
A^{p-1}}_{M}=\lambda^{\frac{A^{p}-1}{A-1}}_{M}
=\lambda^{(A-1)^{p-1}}_{M}
\end{align*}
donde la \'ultima igualdad se debe a que
$\frac{A^{p}-1}{A-1}=\frac{(A-1)^{p}}{A-1}$.

Por lo tanto $\Tr_{G}(\lambda_{M})\in \K\cap
\Lambda_{M}=\Lambda_{N}$. De aqu\'i se obtiene que existe $C\in
R_{T}$ tal que $\lambda^{(A-1)^{p-1}}_{M}=\lambda^{C}_{N}$. Como
$\sigma^{p}=1$ se tiene que
$\sigma^{p}(\lambda_{M})=\lambda^{A^{p}}_{M}=\lambda_{M}$, es decir,
$\lambda^{A^{p}-1}_{M}=0$. Puesto que $A^{p}-1= (A-1)^{p}$, tenemos
que $M\mid (A-1)^{p}$.

Podemos escribir $A-1= P^{\gamma_{1}}_{1}\cdots
P^{\gamma_{r}}_{r}Q$ con $(Q,P_{1}\cdots P_{r})=1$. Ahora si
$\beta_{i}< \alpha_{i}$ se tiene que $\lambda_{NP_{i}}\in L
\setminus \K$ por lo tanto
$\sigma(\lambda_{NP_{i}})=
\lambda^{A}_{NP_{i}}\neq
\lambda_{NP_{i}}$. Por lo tanto $NP_i\nmid (A-1)$.

De aqu\'i se sigue lo siguiente:

(i) Puesto que $M\nmid (A-1)$ se tiene que
$\gamma_{i_{0}}<\alpha_{i_{0}}$ para alg\'un $i_{0}\in \{1,\ldots
r\}$.

(ii) Puesto que $N\mid (A-1)$ se tiene que
$\beta_{i}\leq\gamma_{i}$. Ya que $NP_{i}\nmid (A-1)$, entonces
$\beta_{i}+1>\gamma_{i}$. Se sigue que $\beta_{i}=\gamma_{i}$.

(iii) Como $\lambda^{(A-1)^{p-1}}_{M}=\lambda^{C}_{N}$ se tiene que
$\alpha_{i}-(p-1)\gamma_{i}\leq\beta_{i}$, para $1\leq i\leq r$.

(iv) De $M\mid (A-1)^{p}$ se sigue que $\alpha_{i}\leq p\gamma_{i}$,
para $1\leq i\leq r$.

Ahora $\Tr_{G}(\lambda^{B}_{M})=(\Tr
(\lambda_{M}))^{B}=\lambda^{B(A-1)^{p-1}}_{M}$ para cualquier
$B\in R_T$.

Sea $B=P^{\delta_{1}}_{1}\cdots P^{\delta_{r}}_{r}R$ con
$(R,P_{1}\cdots P_{r})=1$. Se tiene
\begin{align*}
\lambda^{B}_{M}\in \ker \Tr_{G} &\Leftrightarrow \delta_{i}+(p-1)\gamma_{i}\geq\alpha_{i}\, \text{para cada $i$}\\
&\Leftrightarrow \delta_{i}\geq 0 \,\text{y}\, \delta_{i}+(p-1)\gamma_{i}\geq\alpha_{i}\, \text{para cada $i$}\\
&\Leftrightarrow \delta_{i}\geq
\max\{0,\alpha_{i}-(p-1)\gamma_{i}\} \text{ para cada $i$.}
\end{align*}

Por lo tanto $\ker \Tr_{G} = (\lambda^{B}_{M})$ con
$B=P^{\delta_{1}}_{1}\cdots P^{\delta_{r}}_{r}$ y $\delta_{i}=
\max\{0,\alpha_{i}-(p-1)\gamma_{i}\}$ con $1\leq i\leq r$. As\'i
$(\lambda^{B}_{M})=(\lambda_{M^{\prime}})$, con
$M^{\prime}=P^{\mu_{1}}_{1}\cdots P^{\mu_{r}}_{r}$ donde $\mu_{i}=
\alpha_{i}-\delta_{i}$, $1\leq i\leq r$.

Adem\'as $I_{G}(\lambda_{M})=((\sigma-1)\lambda_{M})=
(\lambda^{A-1}_{M})$, donde $I_{G}:\mu(L)\rightarrow\mu(L)$ es el
homomorfismo definido por $I_{G}(u)=\sigma(u)-u$. Por otro lado
$I_{G}(\lambda_{M})=(\lambda_{M^{\prime\prime}})$, con
$M^{\prime\prime}=P^{\varphi_{1}}_{1}\cdots P^{\varphi_{r}}_{r}$,
donde $\varphi_{i}=\max\{\alpha_{i}-\gamma_{i},0\}$, $1\leq
i\leq r$.

De (ii) obtenemos que $\varphi_{i}= \alpha_{i}-\beta_{i}$ si
$\beta_{i}<\alpha_{i}$. Si
$\alpha_{i}=\beta_{i}$ de (ii) se obtiene que
$\alpha_{i}-\gamma_{i}\leq  0$. Por lo tanto
$\varphi_{i}=\alpha_{i}-\beta_{i}$.

As\'i, se tiene que
\begin{gather*}
\mid H^{1}(G,\mu(L))\mid = \frac{\mid
(\lambda_{M^{\prime}})\mid}{\mid(\lambda_{M^{\prime\prime}})\mid}
=\mid(\lambda_{M^{\prime\prime\prime}})\mid\\
\intertext{con $M^{\prime\prime\prime}=P^{\varepsilon_{1}}_{1}\cdots
P^{\varepsilon_{r}}_{r}$ donde}
\varepsilon_{i}=\mu_{i}-\varphi_{i}=\alpha_{i}-\delta_{i}
-(\alpha_{i}-\beta_{i})=\beta_i-\delta_i, \quad 1\leq i\leq r.
\end{gather*}

Obviamente $\varepsilon_i\leq \beta_{i}$, por tanto
\[
\mid H^{1}(G,\mu(L))\mid=q^{\deg M^{\prime\prime\prime}}\leq q^{\deg N}.
\]
Combinando esta desigualdad y el Lema
\ref{auxiliar_coho_cogalois}, obtenemos
\begin{align*}
\mid\cog(L/\K)\mid &= \mid H^{1}(G,\mu(L))\mid\mid
B^{1}(G,\mu(L))\mid = \mid H^{1}(G,\mu(L))\mid\frac{\mid \mu(L)
\mid}{\mid
\mu(\K)\mid}\\
&=q^{\deg M^{\prime\prime\prime}}q^{\deg M-\deg N} \leq
q^{\deg M}= q^{\deg (\mu(L))}.
\end{align*}

Ahora sea $[L:\K]=p^{m}$ para alg\'un $m\in {\mathbb{N}}$,
$m\geq 2$.
Sea $H$ un subgrupo de $G$ de orden $p^{m-1}$. Sea $E=L^{H}$. Entonces
$\K\subseteq E\subseteq L$. Tenemos $[E:\K]=p$,
$[L:E]=p^{m-1}$ y $L=E(\mu(L))$.

Si $E/\K$ no fuera radical ciclot\'omica, entonces
$\cog(E/\K)=\{0\}$, ya que en caso contrario existe
$\overline{\alpha}\in \cog(E/\K)$ no cero. En
particular $\alpha\notin \K$. De esta manera
$E=\K(\alpha)$, pero esto implica que 
$E/\K$ es radical ciclot\'omica,
lo cual es absurdo.

Si $E/\K$ es radical ciclot\'omica se tienen dos casos a considerar

(i) $\mu(E)\neq\mu(\K)$ y

(ii) $\mu(E)=\mu(\K).$

En el caso (i), por lo demostrado para el caso $[E:\K]=p$ se tiene
que
\[
\mid\cog(E/\K)\mid\leq q^{\deg(\mu(E))}.
\]

En el caso (ii), por la Proposici\'on \ref{acotacion1} se tiene que
$\mid\cog(E/\K)\mid=q^{\deg(\mu(E))}$. As\'i, en cualquier
caso,
\[
\mid\cog(E/\K)\mid\leq q^{\deg(\mu(E))}\leq q^{\deg(\mu(L))}.
\]

Por lo tanto, puesto que $L=E(\mu(L))$ y $[L:E]=p^{m-1}$, por
inducci\'on se tiene que $\mid\cog(L/E)\mid\leq
q^{(m-1)\deg(\mu(L))}$. Por lo tanto de la sucesi\'on exacta
\[
0\rightarrow \cog(E/\K)\rightarrow \cog(L/\K)\rightarrow\cog(L/E)
\]
se tiene que $\mid \cog(L/\K)\mid\leq\mid
\cog(E/\K)\mid\mid \cog(L/E)\mid\leq q^{m\deg(\mu(L))}$.
$\fin$
\end{proof}

De la demostraci\'on de la Proposici\'on \ref{acotacion2}, para el
caso $m=1$, se tiene el siguiente teorema.

\begin{teorema}
Sea $L/\K$ c\'iclica de grado $p$, $L=\K(\alpha)$ tal que $\alpha\in
\cog(L/\K)$. Entonces $L/\K$ es radical ciclot\'omica,
$\mid\cog(L/\K)\mid=p^{\nu t}$, donde $q=p^{\nu}$ y
\[ 
t = \left\lbrace
\begin{array}{c l}
\deg(\mu(L))-\deg(\mu(\K))+\deg(\frac{B-1}{C}) & \text{si $\mu(L)\neq\mu(\K)$},\\
\deg(\mu(L)) & \text{si $\mu(L)=\mu(\K)$}.
\end{array}
\right.
\]
donde $\sigma(\lambda_{M})=\lambda^{B}_{M}$, $B-1=\mcd(A-1,B)$ y $C$ es
de grado m\'inimo tal que $C\mid (B-1)$ y $M\mid C(B-1)^{p-1}$. $\fin$
\end{teorema}

\begin{proposicion}\label{acotacion3}
Sea $L/\K$ una extensi\'on Galois radical ciclot\'omica. Entonces
\[
\mid\cog(L/\K)\mid\leq q^{m\deg(\mu(L))}.
\]
donde $[L:\K]=p^{m}$.
\end{proposicion}

\begin{proof} 
Sea $E=\K(\mu(L))$ con $\K\subseteq E\subseteq L$, entonces
\[
\mid \cog(L/\K)\mid\leq\mid \cog(E/\K)\mid\mid \cog(L/E)\mid\leq q^{m\deg (\mu(L))}
\]
por la Proposiciones \ref{acotacion1} y \ref{acotacion2}.
$\fin$
\end{proof}

El ejemplo siguiente muestra que la desigualdad de la Proposici\'on
\ref{acotacion3} puede ser estricta.

\begin{ejemplo} \label{ejemplo_no_se_alcanza_cota}
Sea $L=K(\Lambda_{P^{2p-1}})$, con $P\in R_{T}$ irreducible, y
$\sigma=1+P^{2}\in \Gal(K(\Lambda_{P^{2p-1}})/K)$. Se tiene
\[
\sigma(\lambda_{2p-1})=\lambda_{P^{2p-1}}+\lambda_{P^{2p-3}}
\neq \lambda_{P^{2p-1}}
\]
de este modo $\sigma\neq 1$.

Por otro lado se tiene que $\sigma^{p}= (1+P^{2})^{p}=1+P^{2p}$, de
esta manera
\[
\sigma^{p}(\lambda_{P^{2p-1}})=\lambda_{P^{2p-1}}+\lambda^{P^{2p}}_{P^{2p-1}}=\lambda_{P^{2p-1}}.
\]
Por lo tanto $\sigma^{p}=1$, as\'i el orden de $\sigma$ es
$p$.

Sea $E = L^{(\sigma)}$. Entonces $[L:E]=p$ y $L/E$ es radical
ciclot\'omica. Se tiene que $\sigma(\lambda^{M}_{P^{2p-1}})=
\lambda^{M}_{P^{2p-1}}+\lambda^{M}_{P^{2p-3}}=\lambda^M_{P^{2p-1}}$ si
y s\'olo si el exponente en que aparece $P$ en la descomposici\'on
de $M$ es mayor o igual a $2p-3$. En este caso
$\lambda^{P^{2p-3}}_{P^{2p-1}}= \lambda_{P^{2}}\in E$. Notemos que
$\lambda_{P^{3}}\notin E$. Por lo tanto $\mu(E)=\Lambda_{P^{2}}$.
Adem\'as $\mu(L)= \Lambda_{P^{2p-1}}$.

Ahora sea $N_{\mu(L)}$ el mapeo traza de $L$ a $E$, esto es, $N_{\mu(L)}
=\sum_{i=0}^{p-1}\sigma^i$. Tenemos
 $N_{\mu(L)}(\lambda^{M}_{P^{2p-1}})=
\lambda^{M(\frac{(1+P^{2})^{p}-1 }{(1+P^{2})-1 } )
}_{P^{2p-1}}=\lambda^{MP^{2p-2}}_{P^{2p-1}}=\lambda^{M}_{P}=0$ si y
s\'olo si $P$ divide a $M$. Por lo tanto $\ker
N_{\mu(L)}=(\lambda^{P}_{P^{2p-1}})=\Lambda_{P^{2p-2}}$.

Sea $G:=\Gal(L/E)=\langle \sigma\rangle$ y $I_G
:=\langle\sigma-1\rangle$. Entonces
$I_{G}(\mu(L))=(\sigma(\lambda_{P^{2p-1}})-\lambda_{P^{2p-1}})=
(\lambda_{P^{2p-3}})=\Lambda_{P^{2p-3}}$. Por lo tanto
\[
\mid\cog(L/E)\mid=\mid H^{1}(G,\mu(L))\mid\frac{\mid \mu(L) \mid }{\mid 
\mu(E)  \mid}=\frac{\mid\Lambda_{ P^{2p-2}} \mid }{\mid\Lambda_{P^{2p-3}} 
\mid}\frac{\mid\Lambda_{P^{2p-1}}  \mid}{\mid\Lambda_{P^{2}}   \mid}= q^{d(2p-2)}
\]
donde $d=\deg(P)$. Puesto que $m=1$, se tiene
\[
\mid\cog(L/E)\mid= q^{d(2p-2)}< q^{d(2p-1)}=q^{m\deg(\mu(L))}.
\]
\end{ejemplo}

\begin{teorema}\label{cota_superior_cog}
Sea $L/\K$ una extensi\'on radical ciclot\'omica. Entonces si
$\widetilde{L}$ es la cerradura de Galois de $L$, se tendr\'a
\[
\mid\cog(L/\K)\mid\leq q^{m\deg(\mu(\widetilde{L}))}
\]
donde $[\widetilde{L}:\K]= p^{m}.$
\end{teorema}

\begin{proof} 
Sea $G=\Gal(\widetilde{L}/\K)=HN$ con 
$H$ un subgrupo normal de $G$,
$N$ el $p$ subgrupo de Sylow. Sea $F=\widetilde{L}^{H}$. Se puede
suponer que $F=L$ cambiando $H$ por un conjugado. De aqu\'i
tenemos el diagrama
\[
\xymatrix{L \ar@{-}[r]^{H}\ar@{-}[d]
& \widetilde{L} \ar@{-}[d]^{N} \\
\K \ar@{-}[r] & E }
\]

Sea $\alpha\in \cog(L/\K)$, distinto de cero, as\'i existe un
$N\in R_{T}$ tal que $\alpha^{N}=a\in \K$. Puesto que $\alpha\in
\widetilde{L}$, $\alpha^{N}\in \K\subseteq E$, es decir, $\alpha\in
\cog(\widetilde{L}/E)$. Si $\alpha=0$ en
$\cog(\widetilde{L}/E)$ tendr{\'\i}amos $\alpha\in E\cap L=\K$
por lo que $\alpha=0$ en $\cog(L/\K)$ lo cual es una contradicci\'on.

Por lo tanto $\cog(L/\K)\subseteq \cog(\widetilde{L}/E)$.

De esta manera $\mid\cog(L/\K)\mid\leq
\mid\cog(\widetilde{L}/E)\mid\leq
q^{m\deg(\mu(\widetilde{L}))}$.
$\fin$
\end{proof} 

%% file: Capitulo12.tex
\chapter{$p$--extensiones abelianas 
en caracter{\'\i}stica $p$}\label{Ch9'}

\section{Introducci\'on}\label{S1.pea}

Recordemos la definici\'on de extensiones de Kummer
(ver Definici\'on \ref{D9*.1.3}).

\begin{definicion}\label{D.9`.1} Una {\em extensi\'on de Kummer\index{extensi\'on de
Kummer}\index{Kummer!extensi\'on de $\sim$}} finita es una extensi\'on
de Galois $L/K$ donde el grupo de Galois de la extensi\'on
es $G=\Gal(L/K)$ un grupo c{\'\i}clico finito
$C_n$ de orden $n$ y tal que $K$ contiene al conjunto 
$\mu_n=\langle \zeta_n\rangle$ de las $n$--ra{\'\i}ces de unidad y donde
la caracter{\'\i}stica de $K$ no divide a $n$.
\end{definicion}

Las extensiones de Kummer est\'an caracterizadas por
\[
\langle\zeta_n\rangle\subseteq K,\quad L=K(\sqrt[n]{\alpha}),\quad \alpha\in L
\quad \text{y}\quad f(x):=\Irr(\alpha,x,K)=x^n -\alpha^n.
\]

Las ra{\'\i}ces de $f(x)$ son $\big\{\zeta_n^i\alpha\big\}_{i=0}^{n-1}$.
En el caso en que la caracter{\'\i}stica $p$ de $K$ divide a $n$, 
el grupo de las $p$--ra{\'\i}ces de la unidad es el grupo trivial
$\mu_p=\{1\}$ ya que en caracter{\'\i}stica $p$ tenemos que $\xi^p
=1=1^p$ lo cual implica que $(\xi^p-1^p)=(\xi-1)^p=0$ por lo
que $\xi=1$.

Los primeros en considerar las extensiones c{\'\i}clicas de grado
$p$ en caracter{\'\i}stica $p$ fueron E. Artin y O. Schreier \cite{ArtSch27}.
Ellos probaron que una ecuaci\'on $x^p-x-a$ o bien
proporciona una 
extensi\'on c{\'\i}clica de grado $p$ o bien
todas sus ra{\'\i}ces est\'an
en el campo base pues si $\alpha$ es una ra{\'\i}z cualquiera de
$x^p-x-a$ entonces $\{\alpha,\alpha+1,\ldots,\alpha+(p-1)\}$
son todas las ra{\'\i}ces del polinomio. 
Artin y Schreier probaron el rec{\'\i}proco, esto es, toda extensi\'on
c{\'\i}clica de grado $p$ en caracter{\'\i}stica $p$ est\'a generada
por un polinomio irreducible de la forma $f(x)=x^p-x-a\in K[x]$.
En ese mismo trabajo, Artin y Schreier estudiaron extensiones 
c{\'\i}clicas de grado $p^2$.

Usando las t\'ecnicas de Artin--Schreier, A. Albert \cite{Alb34} encontr\'o
una descripci\'on recursiva de generaci\'on de todas las extensiones
c{\'\i}clicas de grado $p^n$ en caracter{\'\i}stica $p$ y prob\'o
que cualquier extensi\'on c{\'\i}clica de grado $p$ es una 
subextensi\'on de una extensi\'on c{\'\i}clica de grado $p^n$. En particular,
si podemos generar una extensi\'on c{\'\i}clica de grado $p$
en caracter{\'\i}stica $p$, entonces existen extensiones c{\'\i}clicas
de grado $p^n$ para toda $n\in{\ma N}$.

En \cite{Wit36-1}, E. Witt encontr\'o condiciones necesarias y 
suficientes para que un elemento $\theta\in K$, sea tal que la extensi\'on
$K=k(\theta)$ sea c{\'\i}clica de grado $p^f$, $f\geq 2$,
sobre $k$ que contiene a una extensi\'on c{\'\i}clica $R/k$ dada
de grado $p^{f-1}$. Usando este resultado, H.L. Schmid \cite{Sch36-0}
dio una caracterizaci\'on de todas las extensiones c{\'\i}clicas
de grado $p^n$, $K/k$, por medio de elementos $\beta_1,\beta_2,
\ldots,\beta_n\in k$. Schmid dio, de manera recursiva, las ecuaciones
que generan estas extensiones. E. Witt \cite{Wit36-2} encontr\'o
una forma vectorial de describir las extensiones c{\'\i}clicas de grado
$p^n$ halladas por Schmid. Esta forma vectorial es lo que conocemos
como {\em vectores de Witt\index{Witt!vectores de $\sim$}\index{vectores
de Witt}}.

Inmediatamente despu\'es de los resultados de Witt, Schmid \cite{Sch36}
interpret\'o los resultados aritm\'eticos que hab{\'\i}a obtenido 
anteriormente en \cite{Sch36-0} y los puso en t\'erminos de los vectores
de Witt.

En este cap\'itulo estudiamos las $p$--extensiones abelianas
$E/F$ donde $F$ es un campo de caracter\'istica $p>0$. El
\'enfasis ser\'a cuando $F=K=\F(T)$ es el campo de funciones
racionales sobre el campo finito $\F$ y tambi\'en cuando
$F=\K$ donde $\K$ es un campo de funciones congruente.

En primer lugar estudiamos
las $p$--extensiones elementales abelianas la cual ha sido
abordado intensivamente por numerosos autores. Usualmente,
la forma de estudiar este tipo de extensiones 
es considerar todas las subextensiones de grado
$p$ y aplicar tanto el criterio de Hasse \cite{Has35}
(ver Proposici\'on \ref{P2.4.Ram3}) como hechos
conocidos sobre el comportamiento de los primos en esas
subextensiones. (ver Cap\'itulo \ref{ChRam}).
Una vez hecho este estudio se regresa a la
extensi\'on total. Garcia y Stichtenoth \cite{GarSti91},
cambiando el punto
de vista m\'as usual, consideraron las 
extensiones de este tipo que pueden ser
dadas por un \'unico polinomio aditivo, a saber extensiones
dadas por la ecuaci\'on 
$y^q-y=\alpha$. Veremos el comportamiento
de los primos en estas extensiones dada por 
esta ecuaci\'on y de hecho lo haremos para 
polinomios aditivos cuyas ra\'ices se encuentren
en el campo base. 

Consideramos un polinomio
aditivo que tiene sus ra\'ices en el campo base
y probamos, entre otras cosas que 
dado un polinomio aditivo $f(X)$ cuyas ra\'ices se encuentren en
el campo base entonces toda $p$--extensi\'on 
elemental abeliana se puede
describir por medio de una ecuaci\'on del tipo $f(X)=u$.
Describiremos la descomposici\'on de primos en
este tipo de extensiones, generalizando los resultados
dados en el Cap\'itulo \ref{ChRam}.

Cuando el campo base es un campo global de funciones racionales,
daremos una cota inferior sobre el \'indice de ramificaci\'on de
los primos ramificados sin tener que referirnos a sus subextensiones
de grado $p$. Tambi\'en caracterizaremos los primos totalmente
descompuestos. En el caso de extensiones c\'iclicas de grado $p$
dadas por ecuaciones de Artin--Schreier,
la relaci\'on entre los generadores es bien conocida
(ver por ejemplo \cite[Secci\'on 5.8]{Vil2006}). En este cap\'itulo
damos el resultado correspondiente para $p$--extensiones elementales
abelianas obtenidas a partir de polinomios aditivos.

El siguiente objetivo de este cap{\'\i}tulo es presentar la teor{\'\i}a
elemental de los vectores de Witt como generadores de extensiones
c{\'\i}clicas de grado $p^n$ y la aritm\'etica de estos campos obtenida
por Schmid. En particular daremos una breve introducci\'on del concepto
de {\em conductor\index{conductor}} de un campo.

Al final generalizamos los resultados obtenidos
para $p$--extensiones elementales abelianas y
$p$--extensiones c\'iclicas
a extensiones que llamamos {\em multic\'iclicas},
las cuales est\'an dadas
por una ecuaci\'on de Witt de la forma $\vec y^q\Witt - \vec y=\vec
\alpha$. Mucho del formalismo de las $p$--extensiones 
elementales abelianas
puede ser usado en esta situaci\'on y de hecho
los resultados pueden seguir
siendo generalizados a otros {\em polinomios vectoriales aditivos}
cuyas ra\'ices se encuentren en el anillo 
de Witt del campo base. Aqu\'i
nos restringimos \'unicamente a la ecuaci\'on $\vec y^q\Witt -
\vec y=\vec \alpha$.

En este cap{\'\i}tulo todos los campos considerados ser\'an de
caracter{\'\i}stica $p$.

\section{Extensiones de Artin--Schreier}\label{S9'.2}

Empezamos por recordar la teor{\'\i}a de las extensiones c{\'\i}clicas
de grado $p$.

Otra versi\'on del siguiente resultado puede consultarse en
el Teorema \ref{CClaseT1.1.3}.

\begin{teorema}[E. Artin y O. Schreier \cite{ArtSch27}]\label{T9'.2.1}
Sea $f(x)=x^p-x-a\in k[x]$ un polinomio irreducible. Sea $K
=k(\alpha)$, donde $\alpha$ es una ra{\'\i}z de $f(x)$. Entonces
$K/k$ es una extensi\'on c{\'\i}clica de grado $p$.

Rec{\'\i}procamente, si $K/k$ es una extensi\'on c{\'\i}clica de grado
$p$, existe $\alpha\in K$ tal que $\Irr(\alpha,x,k)=x^p-x-a$ para 
alg\'un $a\in k$.
\end{teorema}

\begin{proof}
Sea $\alpha\in \bar{k}$, $\bar{k}$ una cerradura algebraica de $k$,
cualquier ra{\'\i}z de $f(x)$. Entonces las ra{\'\i}ces de $f(x)$ son
$\{\alpha+i\}_{i=0}^{p-1}$. En particular $K/k$ es una extensi\'on 
normal y separable y existe $\sigma \in G:=\Gal(K/k)$ con $\sigma
\alpha=\alpha +1$. Entonces $o(\sigma)=p$ y $G=\langle \sigma
\rangle$ es c{\'\i}clico de orden $p$.

Rec{\'\i}procamente, sea $K/k$ una extensi\'on c{\'\i}clica de
grado $p$, digamos $K=k(\beta)$. Sea $\sigma$ un generador de
$G=\Gal(K/k)$, $o(\sigma)=p$. Sean $\sigma^i\beta=\beta_i$, $i=
0,\ldots,p-1$, los conjugados de $\beta$.
Se tiene $\sigma \beta_i=\beta_{i+1}$, $i=0,\ldots,p-1$. 
Notemos que $\beta_0=\beta_p=\beta$.
Consideremos el determinante de Vandermonde:
\[
\det\big(\beta_i^j\big)_{0\leq i,j\leq p-1}=\prod_{i<j}(\beta_i-\beta_j)\neq 0.
\]
En particular, la matriz $\big(\beta_i^j\big)_{0\leq i,j\leq p-1}$ es no singular y 
existe $0\leq j\leq p-1$ tal que $\gamma:=\sum_{i=0}^{p-1}
\beta_i^j\neq 0$. Sea $\alpha:=-\frac{1}{\gamma}\sum_{i=0}^{p-1}
i\beta_i^j$. Se tiene, usando $\beta_p=\beta_0$, que
\begin{align*}
\sigma\alpha&=-\frac{1}{\sigma\gamma}\sum_{i=0}^{p-1}i\sigma(
\beta_i^j)=-\frac{1}{\gamma}\sum_{i=0}^{p-1}i(\beta_{i+1}^j)=
-\frac{1}{\gamma}\Big(\sum_{i=0}^{p-1}(i+1)\beta^j_{i+1}-\sum_{
i=0}^{p-1}\beta_i^j\Big)\\
&=-\frac{1}{\gamma}\sum_{i=1}^{p}i\beta_i^j+\frac{1}{\gamma}
\sum_{i=0}^{p-1}\beta_i^j =-\frac{1}{\gamma}\sum_{i=0}^{p-1}
i\beta_i + 1=\alpha+1.
\end{align*}

Entonces $\big\{\sigma^i\alpha\big\}_{i=0}^{p-1}=
\big\{\alpha+i\big\}_{i=0}^{p-1}$ son conjugados. Notemos que 
\[
\sigma(\alpha^p-\alpha)=(\sigma \alpha)^p-(\sigma\alpha)=
(\alpha+1)^p-(\alpha+1)=\alpha^p+1-\alpha-1=\alpha^p-\alpha,
\]
de donde se sigue que $\alpha^p-\alpha\in k$. 
De esta forma obtenemos que
$\Irr(\alpha,x,k)=x^p-x-a$ donde $a:=\alpha^p-\alpha$.
$\fin$
\end{proof}

\begin{corolario}\label{C.9'.2.2} Si $K/k$ es una extensi\'on c{\'\i}clica
de grado $p$ con $K=k(\alpha)=k(\beta)$ y $\Irr(\alpha,x,k)=x^p-x-a$,
$\Irr(\beta,x,k)=x^p-x-b$, entonces existen $j\in\{1,2,\ldots,p-1\}$ y $c
\in k$ tales $\alpha=j\beta+c$ y $a=jb+c^p-c$ y rec{\'\i}procamente.
\end{corolario}

\begin{proof} Sea $\sigma\in \Gal(K/k)$ tal que $\sigma\alpha=\alpha+1$. Puesto
que $\sigma\beta$ es conjugado a $\beta$, existe $i\in\{1,2,\ldots, p-1\}$
tal que $\sigma\beta=\beta+i$. Sea $j\in\{1,2,\ldots,p-1\}$ tal que
$ji\equiv 1\bmod p$. Entonces 
\begin{gather*}
\sigma(\alpha-j\beta)=\sigma\alpha-j\sigma \beta = (\alpha+1)-
j(\beta+i)=\alpha-j\beta+1-ji=\alpha-j\beta,\\
\intertext{de donde se sigue que $\alpha-j\beta\in k$. Ahora}
\begin{align*}
a&=\alpha^p-\alpha=(j\beta+c)^p-(j\beta+c)=j^p\beta^p+c^p-j\beta-c\\
&=j\beta^p-j\beta+c^p-c=j(\beta^p-\beta)+(c^p-c)= jb+c^p-c.
\end{align*}
\end{gather*}

El rec{\'\i}proco es similar.
$\fin$
\end{proof}

Para otra versi\'on del Corolario \ref{C.9'.2.2}, ver el 
Teorema \ref{CClaseT1.1.5}.

\begin{observacion}\label{O12.2.2(1)}
Se tiene que si $\alpha\in\bar{k}$ es una ra\'iz
de $x^p-x-a$, entonces
\[
x^p-x-a=\prod_{i=0}^{p-1}(x-(\alpha+i)).
\]
\end{observacion}

\begin{notacion}\label{N.9'.2.3}
Sea $a\in k$, una campo de caracter{\'\i}stica $p$. Se denota
$\wp(a):=a^p-a$.
\end{notacion}

\section{Generalidades sobre $p$--extensiones 
elementales abelianas}\label{S4.pea}

\subsection{El polinomio $X^q-X$}\label{S3.pea}

Empezamos con un resultado fundamentales sobre 
extensiones dadas por polinomios $X^q-X=u$.
Aqu\'i hacemos un breve estudio de las
extensiones $K/k$ cuando $k$ es un campo
de funciones con campo de constantes $k_0$, 
$K=k(y)$ con $y^q-y=u\in k$ y 
$\F\subseteq k$ y por tanto $\F\subseteq k_0$.

\begin{proposicion}\label{P3.1.pea} Supongamos que 
$K/k$ es  una $p$--extensi\'on elemental
abeliana de grado $q=p^n$ y tal que ${\ma F}_{q}
\subseteq k_0$. Entonces
existe $y\in K$ tal que $K=k(y)$ y cuyo polinomio m\'inimo es 
$\Irr(X,y,k)=X^{q}-X-a$ para alg\'un $a\in k$.
 
Rec\'iprocamente, si ${\ma F}_{q}\subseteq k_0$ y $\varphi(X)=
X^{q}-X-a\in k[X]$ es irreducible, entonces $K=k(y)$ con $\varphi
(y)=0$ es una $p$--extensi\'on elemental abeliana de grado $q$.
Los campos intermedios $k\subseteq E_{\mu}\subseteq K$
de grado $p$ sobre $k$, est\'an
dados por $E_{\mu}=k(y_{\mu})$ con $\mu\in{\ma F}_{q}^{\ast}$ y
\[
y_{\mu}:=(\mu y)^{p^{n-1}}+(\mu y)^{p^{n-2}}+\cdots+(\mu y)^p+(\mu y),
\]
$y_{\mu}^p-y_{\mu}=\mu a$, esto es, $k(y)=k(\wp^{-1}(U))$ con
$U=\{\mu a\mid \mu \in {\ma F}_{p^n}\}$.
\end{proposicion}

\begin{proof}
Primero supongamos que $K/k$ es una $p$--extensi\'on
abeliana de grado $q$ y tal que $\F\subseteq k_0$. Sean
$y_1,\ldots,y_n\in K$ tales que $K=k(y_1,\ldots,y_n)$ con
cada $y_i$ satisfaciendo una ecuaci\'on de Artin--Schreier
$y_i^p-y_i=\gamma_i\in k$, $1\leq i\leq n$. Sea $\{\mu_1,
\ldots,\mu_n\}$ una base de $\F$ sobre ${\ma F}_p$
y definimos
\[
y:=\mu_1 y_1+\cdots+\mu_n y_n\in K.
\]

Sea $G=\Gal(K/k)$ con $G=\langle \sigma_1,\ldots,\sigma_n
\rangle$ tal que $\sigma_i(y_j)=y_j+\delta_{ij}$ donde $1\leq i,j
\leq n$ y $\delta_{ij}$ es la delta de Kronecker:
$\delta_{ij}=0$ si $i\neq j$ y $\delta_{i,j}=1$ si $i=j$.

Para cualquier $\sigma\in G$, existen $a_1,\ldots,
a_n\in \{0,1,\ldots,p-1\}$ tales que $\sigma =\sigma_1^{a_1}\cdots
\sigma_n^{a_n}$. Entonces
\begin{align*}
\sigma(y)&=\sigma_1^{a_1}\cdots\sigma_{n-1}^{a_{n-1}}\Big(
\sigma_n^{a_n}\big(\sum_{i=1}^n \mu_i y_i\big)\Big)\\
&=\sigma_1^{a_1}\cdots\sigma_{n-1}^{a_{n-1}}\Big(
\mu_1 y_1+\cdots +\mu_{n-1} y_{n-1} + \mu_n(y_n+a_n)\Big)\\
&=\sigma_1^{a_1}\cdots\sigma_{n-2}^{a_{n-2}}\Big(
\sigma_{n-1}^{a_{n-1}}\big(
\mu_1 y_1+\cdots \mu_{n-1} y_{n-1} + \mu_ny_n+\mu_na_n)\Big)\\
&=\sigma_1^{a_1}\cdots\sigma_{n-2}^{a_{n-2}}\Big(
\mu_1 y_1+\cdots +\mu_{n-2} y_{n-2}+\mu_{n-1}y_{n-1}+
\mu_n y_n\\
&\hspace{3cm}+\mu_{n-1}a_{n-1}+\mu_na_n\Big) \\
&=\ldots=y+\sum_{i=1}^n \mu_i a_i.
\end{align*}

Por tanto $\sigma(y)=y\iff \sum_{i=1}^n \mu_i a_i=0\iff
a_1=\ldots=a_n=0$ por ser $\{\mu_1,\ldots,\mu_n\}$ base
de $\F/{\ma F}_p$. Se sigue que $K=k(y)$. Ahora bien,
puesto que $\mu_i\in \F$ para $1\leq i\leq n$, se tiene que
$\mu_i^q=\mu_i$. Por otro lado
\[
y_i^q-y_i=\sum_{j=1}^n (y_i^{p^j}-y_i^{p^{j-1}})=
\sum_{j=1}^n (y_i^p-y_i)^{p^{j-1}}=\sum_{j=1}^n
\gamma_i^{p^{j-1}}
\]
 y por tanto
\begin{align*}
y^q-y=\sum_{i=1}^n \mu_i(y_i^q-y_i)=\sum_{i=1}^n
\sum_{j=1}^n \mu_i\gamma_i^{p^{j-1}}\in k.
\end{align*}

Rec\'iprocamente, supongamos ahora que $\varphi(X)=
X^q-X-a\in k[X]$ es irreducible. Sea $y\in \bar{k}$ una
ra\'iz de $\varphi(X)$, donde $\bar{k}$ es una cerradura
algebraica de $k$. Sea $K=k(y)$.

Para toda $\mu\in\F$ se tiene $\mu^q=\mu$ y por tanto
$y+\mu$ es ra\'iz de $\varphi(X)$ pues 
\[
\varphi(y+\mu)=
(y+\mu)^q-(y+\mu)-a= y^q-y+\mu^q-\mu-a=y^q-y-a=
\varphi(y)=0.
\]

Por tanto $\{y+\mu\}_{\mu\in\F}$ es el conjunto de ra\'ices
de $\varphi(X)$. Se sigue que $K/k$ es una extensi\'on 
normal y separable y por ser $\varphi(X)$ irreducible, $[K:k]=
q=p^n$. Como $y$ y $y+\mu$ son conjugados para toda
$\mu\in \F$, se sigue que $G=\Gal(K/k)=\{\sigma_{\mu}\mid
\mu\in\F\}$ donde $\sigma_{\mu}(y)=y+\mu$. El mapeo
\[
G\lra \F, \quad \sigma_{\mu}\longmapsto \mu
\]
es un isomorfismo de grupos (aditivos) y por tanto,
$G\cong \F\cong C_p^n$, donde $C_p$ denota el 
grupo c\'iclico de $p$--elementos.

Finalmente, para cada $\mu\in \*F$, se define 
\begin{gather*}
y_{\mu}:=(\mu y)^{p^{n-1}}+(\mu y)^{p^{n-2}}+\cdots+(\mu y)^p+(\mu y).
\intertext{Entonces tenemos}
\begin{align*}
y_{\mu}^p-y_{\mu}&=
(\mu y)^{p^{n}}+(\mu y)^{p^{n-1}}+\cdots+(\mu y)^{p^2}+(\mu y)^p\\
&\hspace{1cm}- (\mu y)^{p^{n-1}}-(\mu y)^{p^{n-2}}
-\cdots-(\mu y)^p-(\mu y)\\
&=\mu^qy^q-\mu y=\mu(y^q-y)=\mu a.
\end{align*}
\end{gather*}

Ahora bien, puesto que $\{1,y,\ldots,y^{q-1}\}$ es base de $K/k$,
y $\mu\neq 0$, se tiene que $y_{\mu}\notin k$ y adem\'as, puesto que,
$k\subseteq k(y_{\mu})\subseteq k(y)$,
$k(y_{\mu})/k$ es un $p$--extensi\'on y $[k(y_{\mu}):k]\leq p$,
se sigue que $[k(y_{\mu}):k]=p$.

Se tiene que $\F\lra K$ dado por $\mu\longmapsto y_{\mu}$
es un homomorfismo de grupos aditivos ($y_0=0$).
Si $k(y_{\mu})=k(y_{\delta})$ con $\mu,\delta\in
\*\F$, entonces por el Corolario \ref{C.9'.2.2}
tenemos que $y_{\mu}=jy_{\delta} + c$ para
algunos $j\in \*{{\ma F}_p}$ y $c\in k$. 

Se tiene que
$c=y_{\mu}-jy_{\delta}=y_{\mu-j\delta}$. Si $\mu-j\delta\neq 0$
se tendr\'ia $c=y_{\mu-j\delta}\notin k$ lo cual es absurdo. Por
tanto $\mu=j\delta$ para alg\'un $j\in\{0,1,\ldots, p-1\}$.

En consecuencia de esta forma hemos obtenido 
$\frac{p^n-1}{p-1}$ extensiones distintas
de grado $p$ sobre $k$ contenidas en $K$ dadas por
$\{E_{\mu}=k(y_{\mu})\}_{\mu\in\*\F}$. Como el grupo
$C_p^n$ tiene exactamente $\frac{p^n-1}{p-1}$ subgrupos
de \'indice $p$, se sigue que $\{E_{\mu}\}_{\mu\in
\*\F}$ son todas las subextensiones de $K$ de
grado $p$ sobre $k$.
$\fin$
\end{proof}

\subsection{Resultados generales}\label{S2.pea}

\begin{definicion}\label{D.12.3.2.pea}
Sean $p$ un n\'umero primo, $n\in{\ma N}$ y $q=p^n$. Sea $k$
un campo arbitrario de caracter\'istica $p$. Cuando $k$ sea un
campo de funciones, $k_0$ denotar\'a al campo de constantes de
$k$, al cual supondremos perfecto. 

Un polinomio $f(X)\in k[X]$ se llama {\em aditivo\index{polinomio
aditivo}} si $f(x+y)=f(x)+f(y)$ para cualesquiera
$x,y\in \bar{k}$, una cerradura algebraica fija de $k$.
\end{definicion}

Para otra versi\'on de la Proposici\'on \ref{P12.3.3.pea}, ver la Proposici\'on
\ref{DrinfeldP1.1.6}.

\begin{proposicion}\label{P12.3.3.pea} Sea $f$ un polinomio
aditivo. Entonces $f$ est\'a dado por
\begin{equation}\label{Eq2.1.pea}
f(x)=\sum_{i=0}^n a_ix^{p^i},
\end{equation}
con $a_i\in k$. 
\end{proposicion}

\begin{proof}
Sea $f(x)=\sum_{i=0}^n a_i x^i$. Entonces $f(0)=f(0+0)=f(0)+f(0)$, por tanto $f(0)=
a_0=0$. Tomando la derivada formal de $f$, se tiene
$f^{\prime}(x)=\sum_{i=1}^n ia_i x^{i-1}$. 

Sea $g(x)=f(x+\alpha)-
f(x)-f(\alpha)$ para $\alpha\in F=\bar{k}$.
Entonces $g(\beta)=0$ para toda $\beta\in F$. Puesto
que $F$ es infinito, se sigue que $g(x)=0$, esto es, $f(x+\alpha)=f(x)+f(\alpha)$.
Derivando esta expresi\'on, obtenemos
\[
f^{\prime}(\alpha)=\frac{d}{dx}(f(x+\alpha))|_{x=0}=\frac{d}{dx}\big(
f(x)+f(\alpha)\big)|_{x=0}=f^{\prime}(0).
\]
Es decir, $f^{\prime}(\alpha)=f^{\prime}(0)$ para toda $\alpha\in F$. Puesto que $F$
es infinito, $f^{\prime}(x)=f^{\prime}(0)=c\in F$ es constante. 

Se sigue que $f'(x)=\sum_{i=1}^n ia_i x^{-1}=f'(0)=c$ por lo que $ia_i=0$
para $i=2,\ldots, n$ en $F$. En particular para $i$ tal que $i\not\equiv 0\bmod p$ se
tiene $a_i=0$.

Escribimos
\[
f(x)=\sum_{l=1}^m
a_{lp}x^{lp}=\sum_{j=0}^{m^{\prime}}a_{p^j}x^{p^j}+\sum_{\substack{
\text{$s$ no es po-}\\
\text{tencia de $p$}}}a_s x^s=f_0(x)+f_1(x).
\]
Puesto que $f(x)$ es aditivo, tenemos que
$f_1(x)=f(x)-f_0(x)$ es aditivo.

Sea $f_1(x)=f_2(x)^p$ con $f_2(x)\in F^{1/p}[x]=F[x]$ ($F=\bar{k}$). 
Veamos que $f_2(x)$ es aditivo. Se tiene para $\alpha,\beta\in\bar{k}$
\begin{align*}
f_2(\alpha+\beta)&=f_1(\alpha+\beta)^{1/p}= f_1(\alpha)^{1/p}
+f_1(\beta)^{1/p}=f_2(\alpha)+f_2(\beta).
\end{align*}
Por inducci\'on en el grado, 
podemos suponer que $f_2(x)=\sum_{i=1}^{m'} c_ix^{p^i}$. Se sigue que
$f_1(x)=f_2(x)^p=\sum_{i=1}^{m'} c_i^p x^{p^{i+1}}$ con $c_i^p\in F$. 
Finalmente
obtenemos que $f(x)=f_0(x)+f_1(x)=\sum_{j=0}^m a_j x^{p^j}$. 
$\fin$
\end{proof}

\begin{observacion}\label{O12.3.4.pea} 
Ver la Observaci\'on \ref{DrinfeldO1.1.5}.
Se puede definir que un polinomio $f(x)=\sum_{i=0}^m a_i x^i\in
F[x]$ es aditivo si para $\alpha,\beta\in F$, $f(\alpha+\beta)=
f(\alpha)+f(\beta)$. Si $F$ es infinito, la demostraci\'on de 
la Proposici\'on \ref{P12.3.3.pea} es aplicable a $F$ y entonces
$f$ est\'a dado por $f(x)=\sum_{j=0}^n a_j x^{p^j}$. 

En caso de que $F$ sea finito, el resultado no se cumple.
Por ejemplo, si $F={\ma F}_p$, $p\geq 3$ y $f(x)=g(x)+
(x^p-x)^n$ para $n\in {\ma N}$ y $g(x)=\sum_{i=0}^m a_ix^{p^i}$, 
es un polinomio aditivos pues
$f(\alpha)=g(\alpha)$ para toda $\alpha\in F$, pero, por ejemplo, para $g(x)=0, n=2$
se tiene 
$f(x)=(x^p-x)^2=x^{2p}-2x^{p+1}+x^2$.

La raz\'on de lo anterior, es que $F$ es finito, digamos $F=\F=\{u\in\bar{\ma F}_p\mid
u^q-u=0\}$, entonces $(\alpha^q-\alpha)^n=0$ para toda $n\in{\ma N}$ y toda $\alpha
\in\F$.
\end{observacion}

Siempre supondremos que el polinomio aditivo
$f(X)$ dado por (\ref{Eq2.1.pea}) es m\'onico y separable,
es decir, $a_n=1$ y $a_0\neq 0$. M\'as a\'un, supondremos que
las ra\'ices de $f(X)$ se encuentran en el campo base, esto es,
\begin{equation}\label{Eq2.2.pea}
{\mc G}_f:=\{\xi\in \bar{k}\mid f(\xi)=0\}\subseteq k.
\end{equation}
Como caso especial consideraremos el polinomio aditivo $f(X)=
X^{p^n}-X=X^q-X$.

En general cuando consideremos polinomios 
de la forma $F(X)=f(X)-u\in k[X]$,
supondremos que $F(X)$ es irreducible. 
Sea $K=k(y)$ con $f(y)=u\in k$, es decir, 
$F(y)=0$. Entonces $f(y+\xi)=f(y)+f(\xi)=
f(y)=u$ para todo $\xi\in{\mc G}_f$ por lo que el conjunto de
ra\'ices de $F(X)$ es 
\begin{gather*}
y+{\mc G}_f=\{y+\xi\mid\xi\in{\mc G}_f\}.
\end{gather*}

Se tiene que todo elemento $\sigma$
de $G=\Gal(K/k)$ est\'a totalmente determinado por
$\sigma(y)$ y puesto que $y$ y $\sigma (y)$ son conjugados,
existe $\xi_{\sigma}\in {\mc G}_f$ tal que $\sigma(y)=y+\xi_{\sigma}$.
Por (\ref{Eq2.2.pea}), la extensi\'on 
$K=k(y)$ de $k$ es una extensi\'on de Galois de grado $p^n$.

\begin{proposicion}\label{P2.1.pea} En general, para cualquier
polinomio aditivo $f(X)\in k[X]$ de grado $p^n$, se tiene que
${\mc G}_f$ es un grupo aditivo ${\mc G}_f\subseteq (\bar{k},+)$
isomorfo a $C_p^n=\big({\ma Z}/p{\ma Z}\big)^n$.
Esto es, ${\mc G}_f$ es un $\ma F_p$--espacio vectorial de 
dimensi\'on $n$.
\end{proposicion}

\begin{proof} Para $\alpha,\beta\in{\mc G}_f$ se tiene
\begin{gather*}
f(\alpha+\beta)=f(\alpha)+f(\beta)=0+0=0,\\
f(0)=0,\quad f(-\alpha)=f((p-1)\alpha)=
\underbrace{f(\alpha)+\cdots +f(\alpha)}_{p-1}=
0+\cdots+0=0.
\end{gather*}

Por tanto ${\mc G}_f\subseteq (\bar{k},+)$. Finalmente, 
$p\beta =0$ para todo $\beta\in \bar{k}$ y $|{\mc G}_f|=p^n$,
por lo que ${\mc G}_f\cong C_p^n$. $\fin$
\end{proof}

En general, si $V$ es un $p$--subgrupo finito de
$\bar{k}$, entonces denotamos:
\begin{gather*}
f_V(X)=\prod_{\delta\in V}(X-\delta)
\end{gather*}
el cual es un polinomio aditivo (Proposici\'on
\ref{DrinfeldP2.2.15}). En particular un polinomio
$f(X)\in k[X]$ aditivo (\ref{Eq2.1.pea}) satisface $f(X)=f_{{\mc G}_f}(X)$.

\begin{proposicion}\label{P2.2.pea} Con
 las notaciones anteriores, tenemos que $\theta\colon G\lra {\mc G}_f$
dado por $\theta(\sigma)=\xi_{\sigma}$, donde $\sigma y=y+\xi_{\sigma}$,
es un monomorfismo de grupos y por tanto podemos considerar $G\subseteq
{\mc G}_f$. Cuando $F(X)$ es irreducible, se tiene la igualdad $G={\mc G}_f$.
\end{proposicion}

\begin{proof}
Para $\sigma,\tau\in G$, se tiene
\begin{gather*}
y+\xi_{\sigma\tau}=\sigma\tau(y)=
\sigma(y+\xi_{\tau})=\sigma(y)+\xi_{\tau}=y+\xi_{\sigma}+\xi_{\tau},\\
\intertext{por lo que}
\theta(\sigma\tau)=\xi_{\sigma\tau}=\xi_{\sigma}+\xi_{\tau}=
\theta(\sigma)+\theta(\tau).
\end{gather*}
Si $\theta(\sigma)=0$, entonces $\sigma(y)=y$ por lo que 
$\sigma=\Id$. Se sigue que $\theta$ es un monomorfismo de grupos.

En el caso de que $F(X)$ sea irreducible, $|G|=[K:k]=\deg F(X)=
p^n=|{\mc G}_f|$. $\fin$
\end{proof}

\begin{observacion}\label{O2.3.pea} 
En el caso en que ${\mc G}_f\nsubseteq k$,
se tiene que el campo de descomposici\'on de $f(X)$ es $k({\mc G}_f,y)$ y
por tanto $K=k(y)$ no es una extensi\'on normal de $k$. De hecho
$\Gal(k({\mc G}_f,y)/k)$ es el  producto semidirecto 
$\Gal(k({\mc G}_f,y)/k)\cong 
\Gal(k({\mc G}_f)/k)\ltimes\Gal(k({\mc G}_f,y)/k({\mc G}_f))$.
\end{observacion}

En general, si $\beta_1,\ldots,\beta_m\in k$, ${\mc L}_{{\ma F}_p}\{\beta_1,\ldots,\beta_m\}$
denota al ${\ma F}_p$--espacio vectorial generado por $\beta_1,
\ldots,\beta_m$.

Denotamos por $\wp$ al operador de Artin--Schreier, es
decir, $\wp(c)=c^p-c$ y por $\wp_a$ a $\wp_a(c)=c^p-a^{p-1}c=
a^p\wp\big(\frac{c}{a}\big)$. En el caso de un campo de funciones
$k/k_0$, con $k_0$ un campo finito,
$R_T$ denota al anillo de polinomios
$k_0[T]$ y  $R_T^+$ denota al conjunto de 
polinomios m\'onicos de $k_0[T]$.

De momento consideraremos el caso especial del polinomio $X^q-X$
con las notaciones y convenciones dadas al principio de la secci\'on. Sea
$k/k_0$ un campo de funciones, 
$G=\Gal(K/k)\cong C_p^n$, $K=k(y)$, $y^q-y=u\in k$. 

\begin{teorema}\label{T4.3.pea}
Con las notaciones anteriores, dado un primo $\P$ de $k$, se tiene
que existe $y\in K$ con $K=k(y)$, $y^q-y=u$ con $v_{\P}(u)\geq 0$
o $v_{\P}(u)=-\lambda p^m$ donde $\lambda>0$, $\mcd(\lambda, p)=1$
y $0\leq m < n$.
Si $v_{\P}(u)\geq 0$, $\P$ es no ramificado en $K/k$.
Si $v_{\P}(u)=-\lambda p^m$, entonces $p^{n-m}\mid e_{\P}$ donde
$e_{\P}$ denota el \'indice de ramificaci\'on de $\P$ en $K/k$.

En particular, si $m=0$, $\P$ es totalmente ramificado.
\end{teorema}

\begin{proof}
Ver la demostraci\'on del Teorema \ref {T5.1.pea}. 
Ah\'i presentamos el mismo resultado en un caso
m\'as general. $\fin$
\end{proof}

\begin{observacion}\label{O4.1.pea} 
El n\'umero $m$ 
que satisface el Teorema \ref{T4.3.pea} no es \'unico.
\end{observacion}

\begin{ejemplo}\label{Ej4.2.pea}
Sea $k=k_0(T)$ con ${\ma F}_{p^2}\subseteq k_0$ y
sea $K=k(y)$ con $y^{p^2}-y
=u=T^{\lambda p}$ con $\lambda\in{\ma N}$, $\mcd(\lambda,p)=1$. 
Entonces $v_{\p}(u)=-\lambda p$ y en este caso $m=1$, $n=2$.

Sea $z:=y^p-T^{\lambda}$. Entonces
\begin{align*}
z^{p^2}-z&=(y^p)^{p^2}-(T^{\lambda})^{p^2}-y^p+T^{\lambda}=
(y^{p^2}-y)^p-T^{\lambda p^2}+T^{\lambda}\\
&=(T^{\lambda p})^p-
T^{\lambda p^2}+T^{\lambda}=T^{\lambda}=\nu,
\end{align*}
y en este caso $v_{\p}(\nu)=-\lambda$, $m=0$ y $n=2$.

Notemos que necesariamente $k(z)=k(y)=K$ pues $e_{\p}=
p^2=[K:k]$. M\'as generalmente, esto es consecuencia inmediata
del Teorema \ref{T7.1.pea}.
\end{ejemplo}

\begin{definicion}\label{D4.4.pea}
Cuando un primo $\P$ satisface las condiciones del Teorema \ref{T4.3.pea}
con respecto a la ecuaci\'on $y^q-y=u$, decimos que la ecuaci\'on
est\'a en una {\em forma normal con respecto a $\P$}. 
Por la Observaci\'on \ref{O4.1.pea} tenemos que
la forma normal de la ecuaci\'on no es \'unica.
\end{definicion}

\begin{observacion}\label{O4.5.pea}
No necesariamente se tiene 
$e_{\P}=p^{n-m}$, ver el Ejemplo \ref{Ej4.2.pea}.
En ese ejemplo tenemos con $m=0$, $p^2\mid e_{\P}$ y por tanto
$e_{\P}=p^2$. Para $m=1$, $p^{2-1}=p\mid e_{\P}$ pero $e_{\P}
\neq p$. M\'as a\'un, aunque $m$ sea m\'inimo con las propiedades
anteriores e inclusive aunque $m$ sea \'unico en alg\'un caso,
tampoco se sigue que necesariamente $e_{\P}=p^{n-m}$.
Sin embargo, si $m=0$, esto es, $v_{\P}(u)=-\lambda$,
entonces $p^n\mid e_{\P}$, por lo que $e_{\P}=p^n$ y
el primo es totalmente ramificado.
\end{observacion}

\begin{ejemplo}\label{Ej4.6.pea}
Sea $k=k_0(T)$ con ${\ma F}_{p^2}
\subseteq k_0$. Sea
$K=k(y_1,y_2)$ con $y_1^p-y_1=T$, $y_2^p-y_2=T^2$ para
$p>2$ y $y_2^p-y_2=T^3$ si $p=2$. Entonces si $\mu\in {\ma F}_{
p^2}\setminus \ma F_p$, consideramos $y=y_1+\mu y_2$, con 
lo cual tenemos $K=k(y)$ y
\begin{align*}
y^{p^2}-y&=y_1^{p^2}+\mu^{p^2}y_2^{p^2}-y_1-\mu y_2=
(y_1^{p^2}-y_1)+\mu(y_2^{p^2}-y_2)\\
&=[(y_1^p-y_1)^p+(y_1^p-y_1)]+\mu[(y_2^p-y_2)^p+(y_2^p-y_2)]\\
&=T^p+T+\mu T^{2p}+\mu T^2=T^p(1+\mu T^p)+(T+\mu T^2)=\gamma,
\end{align*}
$v_{\p}(\gamma)=-2p$.

Las otras extensiones intermedias $K/k$ de grado $p$
est\'an dadas por  $k(y_1+\xi y_2)$, $1\leq \xi
\leq p-1$ y satisfacen
\[
(y_1+\xi y_2)^p-(y_1+\xi y_2)=(y_1^p-y_1)+\xi (y_2^p-y_2)=T+\xi T^2,
\]
$v_{\p}(T+\xi T^2)=-2$. En particular $\p$ se ramifica totalmente, 
$e_{\p}=p^2$, $n=2$, $m=1$, $n-m=1<2$.

Notemos que no podemos obtener $m=0$ con ning\'un cambio de
variable pues si esto fuese posible tendr\'iamos $K=k(z)$ con $z^{p^2}-
z=\nu\in k$ y $v_{\p}(\nu)=-\lambda$ con $\mcd(\lambda, p)=1$, 
entonces en todas las extensiones intermedias se tendr\'ia
(ver la Proposici\'on \ref{P3.1.pea})
\[
z_{\mu}^p-z_{\mu}=\mu\nu,\quad \nu\in{\ma F}_{p^2}^{\ast}
\quad\text{y}\quad v_{\p}(\mu\nu)=-\lambda
\]
y por tanto el exponente del diferente de ${\eu p}_{\infty}$, donde
${\eu p}_{\infty}$ es un primo sobre $\p$, ser\'ia $(\lambda+1)(p-1)$
y $\lambda+1$ ser\'ia el \'unico n\'umero de ramificaci\'on para todas
estas extensiones intermedias. Sin embargo en las subextensiones
$k(y_1)/k$ y $k(y_2)/k$ los n\'umeros de ramificaci\'on son diferentes,
a saber, $1+1=2$ y $2+1=3$ (o $3+1=4$ en el caso $p=2$),
ver \cite[Theorem 5.8.11]{Vil2006}.

As\'i, $m=1$ es \'unico, $p^{n-m}=p^{2-1}=p\mid e_{\p}$ pero 
$e_{\p}= p^2\neq p$.
\end{ejemplo}

\begin{ejemplo}\label{Ej4.7.pea}
Sean $y_1^p-y_1=T$ y $y_2^p-y_2=\frac{1}{T}$ y $K=
k(y_1,y_2)$. Entonces $\p$
no es totalmente ramificado pues es no ramificado en $k(y_2)/
k$ y en este caso $e_{\p}=p$, $n=2$, $m=1$ y $p^{n-m}=p=e_{\p}$.
\end{ejemplo}

\section{Polinomios aditivos}\label{S5.pea}

Con las notaciones y convenciones de la Secci\'on \ref{S2.pea}, consideremos
$f(X)\in k[X]$ un polinomio aditivo m\'onico y separable de grado $p^n$:
\begin{equation}\label{E2.pea}
f(X)=X^{p^n}+a_{n-1}X^{p^{n-1}}+\cdots+a_2X^{p^2}+a_1 X^p +a_0 X
\in k[X],\quad a_0\neq 0.
\end{equation}

Por la Proposici\'on \ref{P2.1.pea} tenemos que ${\mc G}_f$ es un grupo aditivo
isomorfo a $C_p^n$ con ${\mc G}_f\subseteq \bar{k}$, es decir, ${\mc G}_f
\subseteq (\bar{k},+)$.

Supondremos como de costumbre
que ${\mc G}_f\subseteq k$. Sea $K=k(y)$ con $f(y)=u\in k$.
Estamos suponiendo que $F(X)=f(X)-u\in k[X]$ es irreducible. Por tanto
de la Proposici\'on \ref{P2.2.pea} tenemos que $G=\Gal(K/k)\cong {\mc G}_f$.

As\'i, de la Proposici\'on \ref{P2.2.pea}, con $F(X)$ irreducible, obtenemos que
el conjunto $\{\varepsilon_1,\ldots,\varepsilon_n\}\subseteq {\mc G}_f$ 
es base de ${\mc G}_f$ si y solamente
si $G=\langle \sigma_{\varepsilon_1},\ldots,\sigma_{\varepsilon_n}\rangle$.

Se tiene que $G$ tiene $\frac{p^n-1}{p-1}$ subgrupos de \'indice $p$, esto es,
$K/k$ tiene $\frac{p^n-1}{p-1}$ subextensiones de grado $p$ sobre $k$.
Veamos un poco m\'as de cerca estas subextensiones.

Si denotamos por $z$ a los elementos tales que $k\subseteq k(z)\subseteq
K$ con $[k(z):k]=p$, entonces $k(z)$ es el campo fijo bajo un subgrupo $H$
de $G$ de \'indice $p$: $k(z)=K^H$. En este caso, si $G=H\oplus
{\ma F}_p \sigma_z$,
entonces $\Gal(k(z)/k)\cong \langle\sigma_z\rangle$ con $\sigma_z
(z)=z+1$.

De esta forma tenemos que si $E=k(z_1,\ldots,z_n)$, entonces
$E=K$ si y solamente si $G=\langle\sigma_{z_1},\ldots,\sigma_{z_n}
\rangle$. Ahora sea $K=k(z_1,\ldots,z_n)$ y denotemos $\Gal(k(z_i)/k)=
\langle\sigma_i\rangle$ y $G\cong \langle\sigma_1,\ldots,\sigma_n\rangle$.
Sea $z:=\alpha_1 z_1+\cdots+\alpha_n z_n$ con $\alpha_1,\ldots,
\alpha_n\in{\ma F}_p$ no todos cero. Entonces si $z_i^p-z_i=
\wp(z_i)=\gamma_i$, se tiene
\[
z^p-z=\wp(z)=\wp\Big(\sum_{i=1}^n\alpha_iz_i\Big)=\sum_{i=1}^n
\alpha_i\wp(z_i)=\sum_{i=1}^n\alpha_i\gamma_i.
\]

Notemos que $\{\gamma_1,\ldots,\gamma_n\}\subseteq k$ es linealmente
independiente sobre ${\ma F}_p$ pues si $\sum_{i=1}^n\alpha_i\gamma_i
=0$ con alg\'un $\alpha_{i_0}\neq 0$, entonces $\gamma_{i_0}=
\sum_{i\neq i_0}\alpha_{i_0}^{-1}\alpha_i\gamma_i$ y por tanto
\begin{align*}
\wp(z_{i_0})&=z_{i_0}^p-z_{i_0}=\gamma_{i_0}=\sum_{i\neq i_0}
\alpha_{i_0}^{-1}\alpha_i\gamma_i\\
&=\sum_{i\neq i_0}\alpha_{i_0}^{-1}\alpha_i
(z_i^p-z_i)=\wp\big(\sum_{i\neq i_0}\alpha_{i_0}^{-1}\alpha_iz_i\big),
\end{align*}
lo cual implica que
\begin{gather*}
\wp\big(z_{i_0}-\sum_{i\neq i_0}\alpha_{i_0}^{-1}\alpha_iz_i\big)=0.
\intertext{Por tanto $z_{i_0}-\sum_{i\neq i_0}\alpha_{i_0}^{-1}
\alpha_iz_i=\beta
\in{\ma F}_p$. Se sigue que}
z_{i_0}\in k(z_1,\ldots,z_{i_0-1},
z_{i_0+1},\ldots,z_n)\quad \text{y que}\quad [K:k]\leq p^{n-1}
\end{gather*}
lo cual es absurdo.
De esta forma tenemos que $\{\gamma_1,\ldots,\gamma_n\}\subseteq
k$ es un conjunto linealmente independiente sobre ${\ma F}_p$.

Regresando a la expresi\'on $z=\sum_{i=1}^n\alpha_i z_i$, se tiene
$\wp(z)=\sum_{i=1}^n\alpha_i\gamma_i=\gamma$. Si tuvi\'esemos que
$\gamma\in \wp(k)$, digamos $\gamma=\wp(A)$ con $A\in k$, entonces
$\wp(z-A)=\wp\big(\sum_{i=1}^n \alpha_iz_i-A\big)=0$. Por tanto
$\sum_{i=1}^n\alpha_iz_i-A=\beta\in {\ma F}_p$. Puesto que 
$\alpha_{i_0}\neq 0$, tendr\'iamos
\[
z_{i_0}=-\sum_{i\neq i_0}\alpha_{i_0}^{-1}\alpha_iz_i+\alpha_{i_0}^{-1}
\beta+\alpha_{i_0}A\in k(z_1,\ldots,z_{i_0-1},z_{i_0+1},\ldots,z_n),
\]
lo cual es absurdo. Por tanto $\gamma\notin \wp(k)$ y $[k(z):k]=p$.

Con este procedimiento obtenemos $p^n-1$ extensiones de grado
$p$. Ahora bien, si $k(z)=k(w)$ con $z=\sum_{i=1}^n\alpha_iz_i$ y
$w=\sum_{i=1}^n\beta_iz_i$, $\alpha_i,\beta_i\in{\ma F}_p$, se tiene
que $z=jw+c$ con $j\in{\ma F}_p^{\ast}$ y $c\in k$. Se sigue que
$c=0$ y que $z=jw$. De esta forma obtenemos $\frac{p^n-1}{p-1}$
extensiones de grado $p$ distintas y por tanto son todas.
Resumiendo, tenemos

\begin{proposicion}\label{P5.0.pea}
Si $K=k(z_1,\ldots,z_n)/k$ es una $p$--extensi\'on elemental
abeliana de grado $p^n$ y $[k(z_i):k]=p$, $1\leq i\leq n$, entonces
las subextensiones de grado $p$ sobre $k$ est\'an dadas por
$k(z)$ donde $z=\alpha_1 z_1+\cdots+\alpha_n z_n$ con $\alpha_1,
\ldots,\alpha_n\in{\ma F}_p$ no todos cero. $\fin$
\end{proposicion}

Ahora consideremos $k/k_0$ un campo de funciones, como de costumbre,
con $k_0$ perfecto y $f(X)\in k_0[X]$ y ${\mc G}_f\subseteq k_0$. Sea
$f(X)$ dado por (\ref{E2.pea}). Sea $K=k(y)$ con $f(y)=u\in k$. Sea $\P$
un lugar de $k$. Se tiene el mismo resultado que en el Teorema \ref{T4.3.pea},
a saber:

\begin{teorema}\label{T5.1.pea} Se puede seleccionar $u\in k$ tal que
o bien $v_{\P}(u)\geq 0$ o bien $v_{\P}(u)=-\lambda p^m$ con 
$\lambda\in {\ma N}$, $\mcd(\lambda,p)=1$ y $0\leq m< n$. En el primer
caso $\P$ es no ramificado en $K/k$ y en el segundo $\P$
es ramificado y $p^{n-m}\mid e_{\P}$.
\end{teorema}

\begin{proof}
M\'as adelante (Teorema \ref{T5.3.pea})
veremos c\'omo obtener todas las subextensiones
de grado $p$ y una vez hecho esto, se puede usar esta informaci\'on
para determinar el tipo de descomposici\'on de $\P$.
Aqu\'i presentamos otra demostraci\'on en
el esp\'iritu de la de Hasse para extensiones de Artin--Schreier \cite{Has35}.

(1) Si $v_{\P}(u)\geq 0$, se tiene $u\in \o_{\P}$ con $f(y)=u$. Si definimos
$h(X)=\Irr(X,y,k)$, entonces $h(X)\mid f(X)-u$. Por tanto $f(X)-u=h(X)l(X)$ y
$f'(X)=h'(X)l(X)+h(X)l'(X)$. Se sigue que $f'(y)=h'(y)l(y)+0$.
De esta forma obtenemos $h'(y)\mid f'(y)$.

Ahora, si $y\in\o_{\eu p}$ donde $\eu p$ es un primo de $K$ sobre
$\P$, $f'(y)=a_0\neq 0$, $v_{\P}(a_0)=0$ pues $a_0\in k_0^{\ast}$.
Por tanto el diferente local satisface $\eu D_{\o_{\eu p}/\o_{\P}}\mid
\langle f'(y)\rangle=\{1\}$. Por tanto $\eu p\nmid \eu D_{K/k}$ y
$\eu p$ es no ramificado.

Otra demostraci\'on de lo anterior puede darse de manera
an\'aloga a la Proposici\'on \ref{P6.3-1.Ram5}.

(2) Si $v_{\P}(u)<0$, sea $v_{\P}(u)=-\lambda p^m$. Si $m<n$, $u$
satisface las condiciones del teorema. Si $m\geq n$, ponemos
$\lambda p^m=\lambda_1 p^n$,
$u\in k\subseteq k_{\P}$ y $\pi$ un elemento primo para $\P$,
$v_{\P}(\pi)=1$, $\pi\in k$. 
Escribamos $u$ en la forma
\begin{equation}\label{Eq4.1.pea}
u=\frac{b_{-\lambda p^m}}{\pi^{\lambda p^m}}+
\frac{b_{-\lambda p^m+1}}{\pi^{\lambda p^m-1}}+\cdots+
\frac{b_{-1}}{\pi}+b_0+b_1\pi+\cdots \in k_{\P}\cong k(\P)((\pi)).
\end{equation}

Existe $c\in k(\P)$ tal que $c^{p^n}=
b_{-\lambda_1 p^n}$. Sea $C\in\o_{\P}$ donde $k(\P)=\o_{\P}/\P$,
tal que $c=C\bmod \P\in k(\P)$. Sea $z=y-C\pi^{-\lambda_1}$. 
Entonces $k(z)=k(y)=K$ y
\begin{align*}
f(z)&=f(y)-f(C\pi^{-\lambda_1})\\
&=u-(C^{p^n}\pi^{-\lambda_1 p^n}+a_{n-1}
C^{p^{n-1}}\pi^{-\lambda_1 p^{n-1}}+\cdots\\
&\hspace{1cm} +a_1C^p\pi^{-\lambda_1 p}+
a_0C \pi^{-\lambda_1})\\
&=\frac{b_{-\lambda p^m}}{\pi^{\lambda p^m}}+
\frac{b_{-\lambda p^m+1}}{\pi^{\lambda p^m-1}}+\cdots+
\frac{b_{-1}}{\pi}+b_0+b_1\pi+\cdots \\
&\hspace{1cm}-\Big(\frac{b_{-\lambda p^m}}{\pi^{\lambda p^m}}+
\frac{b_{-\lambda p^m+1}}{\pi^{\lambda p^m-1}}+\cdots+
\frac{b_{-\lambda_1}}{\pi^{\lambda_1}}\Big)\\
&=\sum_{i\geq -\lambda p^m+1}\gamma_i\pi^i,\quad \gamma_i\in k(\P).
\end{align*}

Por tanto $v_{\P}(u-f(C\pi^{-\lambda p^m}))\geq 
-\lambda p^m+1>-\lambda p^m$.
En resumen, si $v_{\P}(u)=-\lambda p^m$ con $m\geq n$, existe $\delta \in k$
tal que si $z:=y-\delta$, entonces $k(z)=k(y)$ y 
$v_{\P}(u-f(\delta))>-\lambda p^m$.

Con este procedimiento podemos llevar la ecuaci\'on a
$f(y)=u$ con $v_{\P}(u)\geq 0$ o $v_{\P}(u)=
-\lambda p^m$, $\mcd(\lambda,p)=1$,
$\lambda>0$ y $0\leq m< n$. Con esto terminamos el proceso.

Si $u$ es de esta \'ultima forma,
veamos que para todo $\delta\in k$, $v_{\P}(u-f(\delta))
\leq v_{\P}(u)$, es decir, el valor de $v_{\P}(u)$ es el m\'aximo
posible con sustituciones del tipo $z=y-\delta$ con $\delta\in k$.
 Se tiene $f(\delta)=\delta^{p^n}+a_{n-1}\delta^{p^{n-1}}
+\cdots+a_1\delta^p+a_0 \delta$.

Si $v_{\P}(\delta)\geq 0$ entonces $v_{\P}(f(\delta))\geq 0$ y puesto
que $v_{\P}(u)<0$, se sigue
\[
v_{\P}(u-f(\delta))=\min\{v_{\P}(u),v_{\P}(f(\delta))\}=v_{\P}(u).
\]

Si $v_{\P}(\delta)<0$, $v(\delta^{p^n})=p^n v_{\P}(\delta)<p^iv_{\P}(\delta)
\leq v_{\P}(a_i)+p^i v_{\P}(\delta)=v_{\P}(a_i\delta^{p^i})$. Por tanto
$v_{\P}(f(\delta))=p^nv_{\P}(\delta)=v_{\P}(\delta^{p^n})$.
Ahora bien $v_{\P}(f(\delta))\equiv 0\bmod p^n$ y $v_P(u)\not\equiv
0\bmod p^n$. Por tanto $v_{\P}(f(\delta))\neq v_{\P}(u)$ y
$v_{\P}(u-f(\delta))=\min\{v_{\P}(u),v_{\P}(f(\delta))\}\leq v_{\P}(u)$.

Sea $\eu p$ un lugar de $K$ sobre $\P$. Entonces $v_{\eu p}(y)<0$
y por tanto $v_{\eu p}(f(y))=p^nv_{\eu p}(y)=v_{\eu p}(y^{p^n})$.
Por tanto $v_{\eu p}(f(y))=p^nv_{\eu p}(y)=v_{\eu p}(u)=e_{\P}
v_{\P}(u)=-e_{\P}\lambda p^m$ por lo que se sigue que $p^{n-m}
\mid e_{\P}$ debido a que $\mcd(\lambda,p)=1$. $\fin$
\end{proof}

Cuando $u$ se escribe con respecto a un primo $\P$ como
en el Teorema \ref{T5.1.pea}, decimos que $u$ est\'a en una
{\em forma normal con respecto a $\P$}. Una forma normal no
es \'unica en general.

En el caso en que $k=k_0(T)$ es un campo de funciones racionales
y ${\mc G}_f\subseteq k_0$ se tiene:

\begin{teorema}\label{T5.2.pea}
Sean $k=k_0(T)$ un campo de funciones racionales, $f(X)\in k_0[X]$
un polinomio aditivo dado por {\rm{(\ref{E2.pea})}} y $K=k(y)$ una 
$p$--extensi\'on elemental abeliana donde $f(y)=u\in k$ y
$F(X)=f(X)-u\in k[X]$ es irreducible de grado $p^n$. Entonces se
puede seleccionar $u$ satisfaciendo
\begin{equation}\label{Eq5.2.pea}
u=\sum_{i=1}^r \frac{Q_i(T)}{P_i(T)^{\alpha_i}}+ R(T)
\end{equation}
donde $P_1,\ldots, P_r$ son polinomios m\'onicos irreducibles
distintos, $Q_1,\ldots, Q_r\in k_0[T]$, $\mcd(Q_i,P_i)=1$,
$\deg Q_i<\deg P_i^{\alpha_i}$, $\alpha_i=\lambda_i p^{m_i}>0$ con
$0\leq m_i< n$ y $\mcd(\lambda_i,p)=1$ para $1\leq i\leq r$ y
$R(T)$ es un polinomio tal que si $R(T)\notin k_0$ entonces
$\deg R(T)=\lambda_0 p^m>0$ con $\mcd(\lambda_0, p)=1$,
$0\leq m<n$ y si $R(T)\in k_0$ entonces o bien $R(T)=0$ o
$R(T)\notin f(k_0)=\{f(\delta)\mid \delta\in k_0\}$.

Adem\'as $P_1,\ldots,P_r$ son exactamente los primos finitos de $k$
ramificados en $K$ y $\p$ es ramificado si y solamente si $R(T)
\notin k_0$.
\end{teorema}

\begin{proof}
Sea $f(y)=u=\frac{g(T)}{h(T)}$ con $\mcd(g(T),h(T))=1$. Desarrollando
en fracciones parciales obtenemos
\[
u=\sum_{i=1}^r\sum_{j=1}^{\beta_i}\frac{Q_j^{(i)}(T)}{P_i(T)^j} +R(T),
\]
donde $\deg Q_j^{(i)}<\deg P_i^j$ para cualesquiera 
$1\leq j\leq \beta_i$ y $1\leq i\leq r$ y $S(T)\in k_0[T]$.

Si $\beta_1=\lambda p^n>0$, entonces seleccionamos $C\in k_0[T]$
tal que 
\[
C(T)^{p^n}\equiv Q^{(1)}_{\beta_1}(T)\bmod P_1(T)
\]
 lo cual
es posible hacer pues $k_0[T]/(P_1)$ es un campo perfecto.
Usando la sustituci\'on $z=y-\frac{C}{P_1(T)^{\lambda}}$ 
nos da que $K=k(z)$
y $f(z)=f(y)-f\big(\frac{C}{P_1(T)^{\lambda}}\big)
=u-f\big(\frac{C}{P_1(T)^{\lambda}}\big)=w$. Se sigue
que las valuaciones de $w$ para un primo arbitrario $\P\neq \p$ de $k_0$
satisfacen:
\[
v_{\P}(w) \begin{cases} \geq 0&\text{si $v_{\P}(u)\geq 0$}\\
=-\beta_j& \text{si $\P$ es el primo asociado a $P_j(T)$ para $2\leq j\leq r$}\\
>-\beta_1&\text{si $\P$ es el primo asociado a $P_1(T)$}.
\end{cases}
\]
Repitiendo el proceso obtenemos para $\beta_j$, $2\leq j\leq r$,
que $K=k(y)$ con $f(y)=u$
y $u$ tiene la forma (\ref{Eq5.2.pea}) excepto por $R(T)$. 

Ahora si $R(T)=b_dT^d+\cdots +b_0$ satisface que $d=\lambda p^n$
hacemos la sustituci\'on $y=z-cT^{\lambda}$ donde $c^{p^n}=b_d$. 
Siguiendo este proceso llegamos a que $R(T)\in k_0$ o $\deg R(T)=
\lambda p^m$ con $0\leq m< n$. Finalmente si $R(T)\in k_0$ y
$R(T)=f(\delta)$ para alg\'un $\delta\in k_0$, tomamos $z=y-\delta$.

La ramificaci\'on es consecuencia inmediata del Teorema \ref{T5.1.pea}.
$\fin$
\end{proof}

\begin{definicion}\label{D5.2'.pea}
Cuando la ecuaci\'on $f(y)=u$ que define a la
extensi\'on $K = k(y)$ satisface las condiciones del
Teorema \ref{T5.2.pea}
decimos que la ecuaci\'on  est\'a en una
{\em forma reducida}.
\end{definicion}

Notemos  que la forma reducida en general no es \'unica.

A continuaci\'on generalizamos la Proposici\'on \ref{P3.1.pea} para
polinomios aditivos en general.

Sea $K=k(y)$ con $f(y)=u\in k$, $f(X)\in k[X]$ un polinomio m\'onico
aditivo separable cuyas ra\'ices est\'an en $k$ y $F(X)=f(X)-u
\in k[X]$ irreducible, $G=\Gal(K/k)\cong {\mc G}_f\cong C_p^n$.

\begin{teorema}\label{T5.3.pea} Las subextensiones de grado $p$
sobre $k$ de $K/k$ est\'an dadas por $k(z_{{\mc H}})$ con
\[
z_{{\mc H}}^p-z_{{\mc H}}=\frac{u}{f_{{\mc H}}(\varepsilon_{{\mc H}})^p},
\]
donde ${\mc H}<{\mc G}_f$ es un subespacio de ${\mc G}_f$ de codimensi\'on
$1$, $f_{{\mc H}}(X)=\prod_{\delta\in{\mc H}}(X-\delta)$ y ${\mc G}_f=
{\mc H}+{\ma F}_p\varepsilon_{{\mc H}}$, $\varepsilon_{{\mc H}}\in {\mc G}_f$. M\'as
a\'un, $k(z_{{\mc H}})$ es el campo fijo bajo ${\mc H}$: $k(z_{{\mc H}})=K^{{\mc H}}$
y ${\mc G}al(k(z_{{\mc H}})/k)=\langle \sigma_{\varepsilon_{{\mc H}}}\rangle$ donde
$\sigma_{\varepsilon_{{\mc H}}}(y)=y+\varepsilon_{{\mc H}}$.
Se tiene $\sigma_{\varepsilon_{{\mc H}}}(z_{{\mc H}})=z_{{\mc H}}+1$.
\end{teorema}

\begin{proof} Este resultado debe compararse con el Teorema
\ref{T7.2.pea} el cual es m\'as general pero menos expl\'icito.
Se tiene que $f(X)=\prod_{\alpha\in{\ma F}_p}
f_{{\mc H}}(X-\alpha \varepsilon_{{\mc H}})$. Ahora bien, $f_{{\mc H}}(X)$ es un
polinomio aditivo (ver Proposici\'on \ref{DrinfeldP2.2.15}) 
y $f_{{\mc H}}(X-\alpha\varepsilon_{{\mc H}})=
f_{{\mc H}}(X)-\alpha f_{{\mc H}}(\varepsilon_{{\mc H}})$.

Denotemos $Y:=f_{{\mc H}}(X)$. Entonces
\begin{align*}
f(X)&=\prod_{\alpha=0}^{p-1}(Y-\alpha f_{{\mc H}}(\varepsilon_{{\mc H}}))
=f_{{\mc H}}(\varepsilon_{{\mc H}})^p\cdot\prod_{\alpha=0}^{p-1}\Big(\frac{Y}{
f_{{\mc H}}(\varepsilon_{{\mc H}})}-\alpha\Big)\\
&=f_{{\mc H}}(\varepsilon_{{\mc H}})^p
\Big(\Big(\frac{Y}{f_{{\mc H}}(\varepsilon_{{\mc H}})}\Big)^p-\frac{Y}{f_{{\mc H}}(
\varepsilon_{{\mc H}})}\Big).
\end{align*}

As\'i, 
\begin{gather*}
f(X)=f_{{\mc H}}(X)^p-f_{{\mc H}}(\varepsilon_{{\mc H}})^{p-1}f_{{\mc H}}(X)=
f_{{\mc H}}(\varepsilon_{{\mc H}})^p\Big(\Big(\frac{f_{{\mc H}}(X)}{f_{{\mc H}}(\varepsilon_{{\mc H}})}
\Big)^p-\Big(\frac{f_{{\mc H}}(X)}{f_{{\mc H}}(\varepsilon_{{\mc H}})}\Big)\Big).
\intertext{Esto es,}
f(X)=f_{{\mc H}}(\varepsilon_{{\mc H}})^p\wp\Big(\frac{f_{{\mc H}}(X)}{f_{{\mc H}}(\varepsilon_{
{\mc H}})}\Big)=\wp_{f_{{\mc H}}(\varepsilon_{{\mc H}})}(f_{{\mc H}}(X)).
\end{gather*}

De esta forma obtenemos $f(X)=f_{{\mc H}}(\varepsilon_{{\mc H}})^p(z^p-z)$ 
donde $z=\frac{Y}{f_{{\mc H}}(\varepsilon_{{\mc H}})}=\frac{f_{{\mc H}}(X)}{
f_{{\mc H}}(\varepsilon_{{\mc H}})}$. Sea 
\begin{equation}\label{Eq5.3.pea}
z_{{\mc H}}:=\frac{f_{{\mc H}}(y)}{f_{{\mc H}}(
\varepsilon_{{\mc H}})}.
\end{equation}
 Entonces 
\[
z_{{\mc H}}^p-z_{{\mc H}} =\frac{f(y)}{f_{{\mc H}}(\varepsilon_{{\mc H}})^p}=
\frac{u}{f_{{\mc H}}(\varepsilon_{{\mc H}})^p} \qquad\text{o}\qquad
\wp_{f_{{\mc H}}(\varepsilon_{{\mc H}})}\big(
f_{{\mc H}}(\varepsilon_{{\mc H}})z_{{\mc H}}\big)=u.
\]

M\'as a\'un, si $\{\varepsilon_1,\ldots,\varepsilon_{n-1}\}$ es una
$\ma F_p$--base de ${\mc H}$ y si ponemos 
$\varepsilon_n:=\varepsilon_{{\mc H}}$,
entonces $\{\varepsilon_1,\ldots,\varepsilon_{n-1},\varepsilon_n\}$
es una base de ${\mc G}_f/\ma F_p$. Se tiene que si $\sigma_{
\varepsilon_i}(y)=y+\varepsilon_i$, entonces obtenemos que
$G=\langle
\sigma_{\varepsilon_1},\ldots,\sigma_{\varepsilon_{n-1}},
\sigma_{\varepsilon_n}\rangle$ y se tiene para $1\leq i\leq n-1$
\begin{gather*}
\sigma_{\varepsilon_i}(z_{{\mc H}})=\sigma_{\varepsilon_i}\Big(
\frac{f_{{\mc H}}(y)}{f_{{\mc H}}(\varepsilon_n)}\Big)=\frac{f_{{\mc H}}(y+
\varepsilon_i)}{f_{{\mc H}}(\varepsilon_n)}=\frac{f_{{\mc H}}(y)+
f_{{\mc H}}(\varepsilon_i)}{f_{{\mc H}}(\varepsilon_n)}=\frac{
f_{{\mc H}}(y)+0}{f_{{\mc H}}(\varepsilon_n)}=z_{{\mc H}},\\
\intertext{y}
\sigma_{\varepsilon_n}(z_{{\mc H}})=\frac{f_{{\mc H}}(y+\varepsilon_n)}{
f_{{\mc H}}(\varepsilon_n)}=z_{{\mc H}}+1,
\end{gather*}
as\'i que $k(z_{{\mc H}})/k$ es una subextensi\'on de $K/k$ de grado
$p$, el campo $k(z_{{\mc H}})$ es el campo 
fijo bajo ${\mc H}$ y $\Gal(k(z_{{\mc H}})/k)\cong\langle \sigma_{
\varepsilon_{{\mc H}}}\rangle$ donde $\sigma_{\varepsilon_{{\mc H}}}(y)=
y+\varepsilon_{{\mc H}}$. $\fin$
\end{proof}

\begin{teorema}\label{T5.4.pea} Sea $f(X)\in k[X]$ un polinomio m\'onico
aditivo separable de grado $p^n$
con ${\mc G}_f\subseteq k$. Sea $K/k$ una $p$--extensi\'on elemental
abeliana de grado $p^n$. Entonces existen $y\in K$ y $u\in k$ tales
que $K=k(y)$ y $f(y)=u$.
\end{teorema}

\begin{proof}
Sean $y_1,\ldots, y_n\in K$ tales que $K=k(y_1,\ldots, y_n)$ y
$y_i^p-y_i=\gamma_i\in k$. Sea $G=\Gal(K/k)$, $G=\langle
\sigma_1,\ldots, \sigma_n\rangle$ con $\sigma_i(y_j)=
y_j+\delta_{ij}$ con $\delta_{ij}$ la delta de Kronecker. Sea
$\{\mu_1,\ldots,\mu_n\}$ una base de ${\mc G}_f$ sobre $\ma F_p$.
Sea $y:=\sum_{i=1}^n\mu_i y_i$ y sea $f(X)$ dado como en 
(\ref{E2.pea}).

Se tiene $y_i^p=y_i+\gamma_i$, $y_i^{p^2}=y_i^p+\gamma_i^p=
y_i+\gamma_i+\gamma_i^p$ y en general $y_i^{p^m}=y_i+
l_m(\gamma_i)$ donde $l_m(\gamma_i)=\gamma_i+\gamma_i^p
+\cdots+\gamma_i^{p^{m-1}}$, $m\in \ma N$. Entonces
\[
f(\mu_iy_i)=\sum_{j=0}^n a_j(\mu_i y_i)^{p^j}=
\sum_{j=0}^n a_j(\mu_i^{p^j}y_i+\mu_i^{p^j} l_j(\gamma_i))=
y_i f(\mu_i)+h_i,
\]
con $h_i=\sum_{j=0}^n a_j\mu_i^{p^j}l_j(\gamma_i)$ y
$y_if(\mu_i)=0$ pues $\mu_i\in {\mc G}_f$. As\'i,
\[
f(y)=f\big(\sum_{i=1}^n\mu_iy_i\big)=\sum_{i=1}^nf(\mu_iy_i)=
\sum_{i=1}^nh_i=u\in k.
\]

Si $\sigma\in G$, existen $\nu_1,\ldots,\nu_n\in \ma F_p$ tales que
$\sigma=\sigma_1^{\nu_1}\cdots\sigma_n^{\nu_n}$ y
\[
\sigma(y)=\sigma\big(\sum_{i=1}^n \mu_iy_i\big)=\sum_{i=1}^n
\mu_i(y_i+\nu_i)=y+\sum_{i=1}^n\nu_i\mu_i.
\]
Por tanto $\sigma=\Id\iff \sigma(y)=y\iff \sum_{i=1}^n \nu_i\mu_i=0
\iff \nu_1=\cdots=\nu_n=0$, confirmando que $\{\mu_1,\ldots,\mu_n\}$
es base de ${\mc G}_f/{\ma F}_p$ y que $G=\Gal(K/k)=\langle \sigma_1,
\ldots,\sigma_n\rangle$. $\fin$
\end{proof}

\section[Descomposici\'on en $p$-extensiones]{Descomposici\'on 
de primos en $p$--extensiones elementales abelianas
de campos de funciones}\label{S6.pea}

Vamos a considerar $k=k_0(T)$ un campo 
de funciones donde supondremos que
$k_0$ es un campo finito con $\F\subseteq k_0$.
El objetivo de esta secci\'on es describir la
descomposici\'on de los primos no ramificados en una $p$--extensi\'on
elemental abeliana $K$ de $k$ dada por $K=k(y)$ donde
$f(y)=u$, $f(X)\in k_0[X]$ dado por (\ref{E2.pea}) y ${\mc G}_f\subseteq k_0$.
Notemos que el grupo de descomposici\'on
de cualquier primo no ramificado, es c\'iclico y por tanto es de orden
$1$ o $p$. Supondremos que la extensi\'on $K/k$ es geom\'etrica.

Un resultado fundamental que usaremos
en esta secci\'on es la descomposici\'on
de $\p$ en extensiones de Artin--Schreier
(ver Proposici\'on \ref{P2.4.Ram3}).

\begin{proposicion}\label{P6.1.pea}
Sea $K/k$ una extensi\'on c\'iclica de grado $p$
tal que $K=k(w)$ con $w\in K$ de la forma
\begin{gather}\label{Eq6.1.pea}
w^p-w=u=\sum_{i=1}^r\frac{Q_i}{P_i^{e_i}} + R(T)=
\frac{Q}{P_1^{e_1}\cdots P_r^{e_r}}+R(T),
\end{gather}
donde $P_i\in R_T^+$, $Q_i\in R_T$, 
$\mcd(P_i,Q_i)=1$, $e_i>0$, $p\nmid e_i$, $\deg Q_i<
\deg P_i^{e_i}$, $1\leq i\leq r$,
$\deg Q<\sum_{i=1}^r\deg P_i^{e_i}$, $R(T)\in R_T$,
con $p\nmid \deg R(T)$ cuando $R(T)\not\in k_0$
y $R(T)\notin \wp(k_0)$ cuando $R(T)
\in k_0^{\ast}$. Entonces
el divisor primo $\p$ es
\l
\item descompuesto si $R(T)=0$.
\item inerte si $R(T)\in k_0$ y $R(T)\not\in \wp(k_0)$.
\item ramificado si $R(T)\not\in k_0$
(por lo que $p\nmid\deg R(T)$). 
$\fin$
\end{list}
\end{proposicion}

\begin{ejemplo}\label{Ej5.2(1).pea}
Sean $k=\mathbb{F}_{27}(T)$ y $k(z)/k$ la $3$--extensi\'on 
elemental abeliana  de grado $27$,
definida por la siguiente ecuaci\'on: 
\begin{gather*}
z^{27}-z=\frac{1}{(T+1)^{54}}+\frac{1}{T+1}+T^{9}+T^{3}+T+\omega+1=u(T),
\end{gather*}
donde $\omega\in \mathbb{F}_{27} $, $\omega^{3}=
\omega+2$ y $\mathbb{F}_{27}=\mathbb{F}_{3}(\omega)$. 
Primero observemos que si $y=z-\frac{1}{\left(T+1\right)^{2}}$, 
entonces $k(y)=k(z)$ y $y^{27}-y=(z-\frac{1}{\left(T+
1\right)^{2}})^{27}-(z-\frac{1}{\left(T+1\right)^{2}})=u(T)-
\frac{1}{\left(T+1\right)^{54}}+\frac{1}{\left(T+1\right)^{2}}$. 
Esto es 
\[
y^{27}-y=\frac{1}{\left(T+1\right)^{2}}+\frac{1}{T+1}
+T^{9}+T^{3}+T+\omega+1=r(T).
\]
Por el Teorema \ref{T4.3.pea} 
los lugares $\mathcal{P}_{1}$,  asociado al polinomio irreducible 
$p_{1}(T)=T+1$ y $\mathcal{P}_{\infty}$, son ramificados, 
ahora se procede al c\'alculo del \'indice de ramificaci\'on 
para estos lugares, para este fin se utiliza la Proposici\'on 
\ref{P3.1.pea} en donde es posible obtener las $13=\frac{3^{3}-1}{3-1}$ 
sub--extensiones de Artin-Schreier donde las tres 
principales est\'an determinadas por las ecuaciones:  
\begin{align*}
y_{1}^{3}-y_{1}&=r(T),\\
y_{1}&=y^{9}+y^{3}+y.\\
y_{2}^{3}-y_{2}&=\omega r(T),\\
y_{2}&=(\omega y)^{9}+(\omega y)^{3}+\omega y.\\
y_{3}^{3}-y_{3}&=\omega^{2} r(T),\\
y_{3}&=(\omega^{2} y)^{9}+(\omega^{2} y)^{3}+\omega^{2} y.
\end{align*}
El siguiente diagrama representa las 
3 sub--extensiones principales y algunas 
3--extensiones elementales abelianas de grado 9.
\begin{small}
\[
\xymatrix{ 
& & k(y)   & & \\
&k(y_{1},y_{2}) \ar@{-}[ur]   &   k(y_{1},y_{3})\ar@{-}[u] & 
 k(y_{2},y_{3})=k(y)^{I}\ar@{-}[ul]  &\\
k(y_{1}) \ar@{-}[ur] \ar@{-}[urr]|!{[ur];[rr]}\hole & &  
 k(y_{2})=k(y)^{D}\ar@{-}[ul] \ar@{-}[ur] &  & 
  k(y_{3})\ar@{-}[ul] \ar@{-}[ull] |!{[ll];[ul]}\hole\\
 &  & k=\mathbb{F}_{27}(T) \ar@{-}[urr]^{3} 
 \ar@{-}[u]^{3} \ar@{-}[ull]_{3} &  & 
}
\]
\end{small}

Observemos que el lugar $\mathcal{P}_{1}$ 
es totalmente ramificado, esto es 
$e(\wp_{1}|\mathcal{P}_{1})=27$ para un 
lugar $\wp_{1}$ sobre $\mathcal{P}_{1}$ ya que 
$v_{\mathcal{P}_{1}}(r(T))=-2$ el cual es primo 
relativo con $3$. Por otro lado, si $z_{1,1}=y_{1}-T^{3}$, 
entonces $k(z_{1,1})=k(y_{1})$ y $z^{3}_{1,1}-z_{1,1}=
\frac{1}{\left(T+1\right)^{2}}+\frac{1}{T+1}+2T^{3}+T+
\omega+1$ ahora sea $z_{1,2}=z_{1,1}+T-\omega^{2}$, 
se tiene que $k(z_{1,2})=k(z_{1,1})=k(y_{1})$ y 
\[
z^{3}_{1,2}-z_{1,2}=\frac{1}{\left(T+1\right)^{2}}
+\frac{1}{T+1}=r_{1}(T).
\]
Por la Proposici\'on \ref{P6.1.pea}, 
$\mathcal{P}_{\infty}$ es descompuesto en $k(y_{1})/k$.

En $k(y_{2})/k$, sea $z_{2,1}=y_{2}+(2\omega+2)T^{3}$, 
entonces $k(z_{2,1})=k(y_{2})$
 y $z^{3}_{2,1}-z_{2,1}=(2\omega+1)T^{3}+\omega T+
 \omega^{2}+\omega+\frac{\omega}{\left(T+1\right)^{2}}+
 \frac{\omega}{T+1}$. De nuevo, si a la extensi\'on anterior 
 se le realiza la sustituci\'on: $z_{2,2}=z_{2,1}-2\omega T$, 
 entonces $k(z_{2,2})=k(z_{2,1})=k(y_{2})$ y 
\[
z^{3}_{2,2}-z_{2,2}=
 \frac{\omega}{\left(T+1\right)^{2}}+\frac{\omega}{T+1}+
 \omega^{2}+\omega=r_{2}(T),
\]
como $\omega^{2}+\omega\in \mathbb{F}_{27}-
\mathcal{A}(\mathbb{F}_{27})$, por la 
Proposici\'on \ref{P6.1.pea} se tiene que el primo al infinito es 
inerte en $k(y_{2})/k$. Por \'ultimo, de manera an\'aloga y 
con las siguientes sustituciones:
 $z_{3,1}=y_{3}-(2\omega^{2}+\omega+2)T^{3}$ 
y $z_{3,2}=z_{3,1}-(2\omega^{2}+2)T$, se obtiene $k(z_{3,2})
=k(z_{3,1})=k(y_{3})$, donde $$z^{3}_{3,2}-z_{3,2}=
\frac{\omega^{2}}{\left(T+1\right)^{2}}+\frac{\omega^{2}}{T+1}
+2T+\omega^{2}+\omega+2=r_{3}(T).$$ Por lo tanto, 
$\mathcal{P}_{\infty}$ es ramificado en $k(y_{3})/k$.
Concluimos que el campo de 
descomposici\'on del primo al infinito $\mathcal{P}_{\infty}$ 
es $k(y)^{D}=k(y_{2})$ y por tanto el campo 
de inercia es $k(y)^{I}=k(y_{2},y_{3})$. Luego la descomposici\'on 
para el primo $\mathcal{P}_{\infty}$ es: 
$\left(e(\wp_{\infty}|\mathcal{P}_{\infty}),
f(\wp_{\infty}|\mathcal{P}_{\infty}),
h(\wp_{\infty}|\mathcal{P}_{\infty})\right)=(3,3,3)$, 
para un lugar ${\eu p}_{\infty}$ sobre $\mathcal{P}_{\infty}$ en $k(y)$.
\end{ejemplo}

\bigskip

Regresando a nuestro estudio,
observamos que con la hip\'otesis de ser de grado $p$, separable,
y con ra\'ices en el campo base, esencialmente
hay \'unicamente un polinomio aditivo de grado $p$. Si $X^p+aX$
es aditivo, sus ra\'ices son $i\alpha$ con $0\leq i\leq p-1$ y 
$\alpha=\sqrt[p-1]{-a}$ una ra\'iz no cero fija de $X^p+aX$.
Estamos suponiendo que $\alpha\in {\ma F}_q^{\ast}$. Entonces $\alpha^p
=-\alpha a$ y
\[
X^p+aX=\alpha^p\big(Z^p-Z\big)\quad\text{con}\quad Z=\frac{X}{\alpha},
\quad \alpha\in {\ma F}_q^{\ast},
\]
es decir, $X^p+aX=\alpha^p\wp\big(\frac{X}{\alpha}\big)$ con 
$\alpha^{p-1}=-a$.

Volvamos a la extensi\'on $K=k(y)$, $\Gal(K/k)\cong {\mc G}_f$. Se tiene que
la descomposici\'on de $\p$ en $K$, en caso de que sea no ramificado,
es

\begin{proposicion}\label{P6.2.pea}
Sea $K=k(y)/k$ con $f(y)=u$ y $u$ dado en forma reducida 
{\rm{(\ref{Eq5.2.pea})}}. Supongamos que $\p$ es no ramificado en $K/k$.
Si $R(T)=0$, entonces $\p$ se descompone totalmente en $K/k$.
Rec\'iprocamente, si $\p$ se descompone totalmente,
entonces existe una forma reducida $f(y)=u$ en donde $R(T)=0$.
\end{proposicion}

\begin{proof}
Si $R(T)=0$, entonces por el Teorema \ref{T5.3.pea} y la Proposici\'on
\ref{P6.1.pea}, se sigue que $\p$ es descompuesto en todas las
subextensiones de grado $p$ y por tanto $\p$ es totalmente
descompuesto en $K/k$.

Rec\'iprocamente, supongamos que $\p$ es totalmente descompuesto.
Sea $K=k(y_1,\ldots,y_n)$ con $k(y_i)/k$ extensiones c\'iclicas de grado
$p$ dadas por ecuaciones de
Artin--Schreier en forma reducida. Como $\p$ es descompuesto
en todas ellas, por el Teorema \ref{T5.4.pea}, se tiene
$K=k(y_0)$ con $y_0=\sum_{
i=1}^n\mu_iy_i$, $\{\mu_1,\ldots,\mu_n\}$ una base de ${\mc G}_f$ sobre
${\ma F}_p$, y en la forma reducida de $f(y_0)=u_0$,
el polinomio correspondiente al comportamiento de $\p$ es $0$. $\fin$
\end{proof}

En el caso especial del polinomio aditivo $f(X)=X^q-X$, vamos a probar
algo m\'as que la Proposici\'on \ref{P6.2.pea}. Para esto probamos

\begin{lema}\label{L6.2'.pea} Sea $S\in{\ma F}_{q^m}$. 
Entonces se tiene que
$\mu S\in \im\wp$ para todo $\mu\in\F$ si y solamente si existe
$\lambda\in {\ma F}_{q^m}$ tal que $S=\lambda^q-\lambda$.
\end{lema}

\begin{proof}
Consideremos el homomorfismo $g\colon{\ma F}_{q^m}\lra {\ma F}_{q^m}$,
$g(\lambda)=\lambda^q-\lambda$. Entonces $\ker g=\F$ y $|\im g|=
\frac{q^m}{q}$. El homomorfismo de Artin--Schreier, $\wp\colon {\ma F}_{q^m}\to
{\ma F}_{q^m}$, $\wp(\lambda)=\lambda^p-\lambda$, satisface que
$\ker \wp={\ma F}_{p}$ y $|\im\wp |=\frac{q^m}{p}$. Finalmente consideremos
$h\colon {\ma F}_{q^m}\lra {\ma F}_{q^m}$ dado por $h(\lambda)=\lambda
+\lambda^p+\cdots+\lambda^{p^{n-1}}$, $q=p^n$. Entonces 
$h(\lambda)^p-h(\lambda)=\lambda^q-\lambda$.
Se tiene que 
\[
(\wp\circ h)=(h\circ \wp)=g.
\]
 Ahora bien, si $S=\lambda^q-\lambda$
para alg\'un $\lambda\in {\ma F}_{q^m}$,
entonces $\mu S=(\mu\lambda )^q-(\mu \lambda)$ para todo $\mu\in \F$ y
se tiene 
\[
(\mu\lambda )^q-(\mu \lambda)=g(\mu\lambda)=\wp(h(\mu\lambda))
\in \im \wp,
\]
para todo $\mu\in\F$.

Rec\'iprocamente, supongamos que $A:=\big\{\mu S\big\}_{\mu\in \F}
\subseteq \im\wp$. Si
$S=0$ no hay nada que probar. Consideremos $S\neq 0$. 
En caso de que $S\notin\im g$
tendr\'iamos $\mu S\notin \im g$ para todo
$\mu\in{\ma F}_q^{\ast}$ pues en caso contrario $\mu S=\lambda^q-\lambda$ para alg\'un
$\lambda\in {\ma F}_{q^m}$ y para todo $\mu\in\F$, $\mu^q=\mu$ y
\[
S=\frac{\lambda^q}{\mu}-\frac{\lambda}{\mu}=\Big(\frac{\lambda}{\mu}\Big)^q
-\Big(\frac{\lambda}{\mu}\Big)\in \im g,
\]
lo cual es absurdo.

Bajo esta suposici\'on el subgrupo aditivo $A=\big\{\mu S\big\}_{\mu\in \F}$ de
${\ma F}_{q^m}$, tiene cardinalidad $q$
y cumple que $A\cap \im g=\{0\}$.
Ya que $A\subseteq \im \wp$ y puesto que
$g=\wp\circ h$, de donde $\im g\subseteq \im \wp$, 
se sigue que $A+\im g\subseteq
\im\wp$. Por otro lado tenemos
\[
|A+\im g|=|A||\im g|=q\cdot \frac{q^m}{q}=q^m>\frac{q^m}{p}=|\im \wp|,
\]
lo cual es una contradicci\'on. Por tanto $A\cap \im g\neq \{0\}$. 
Luego existe $\mu\in{\ma F}_q^{\ast}$
con $\mu S\in \im g$ lo cual implica que $S\in\im g$. Esto termina la
demostraci\'on. $\fin$
\end{proof}

Sea $f(X)$ un polinomio aditivo dado por (\ref{E2.pea}).
Se debe tener el mismo resultado
que en el Lema \ref{L6.2'.pea}, pero no tenemos la demostraci\'on.
Sea ${\mc G}_f={\mc L}_{{\ma F}_p}\{\varepsilon_1,\ldots,
\varepsilon_n\}$ y consideremos
los siguientes $n$ subespacios de ${\mc G}_f$ de codimensi\'on $1$: 
\[
{\mc H}_i:=
{\mc L}_{{\ma F}_p}\{\varepsilon_1,\ldots,\varepsilon_{i-1},
\varepsilon_{i+1},\ldots,\varepsilon_n\},\quad 1\leq i\leq n.
\]
Sean $f_{{\mc H}_i}(X)=f_i(X)$, $a_i=f_i(\varepsilon_i)\neq 0$. 

 Lo que necesitamos es tener:
\begin{gather}\label{Eq*.pea}
\bigcap_{i=1}^n\im\wp_{a_i}=\im f.
\end{gather}

Se tiene:

\begin{proposicion}\label{P6.2''.pea} Sea $K=k(y)/k$ 
con $f(y)=u$ y $u$ dado en la forma reducida 
{\rm{(\ref{Eq5.2.pea})}}. Supondremos que
$\p$ es no ramificado en $K/k$.
Si $f(x)=x^q-x$ o si se satisface {\rm{(\ref{Eq*.pea})}}, 
entonces $\p$ se descompone totalmente si y
solamente si $R(T)=0$.
\end{proposicion}

\begin{proof}
Primero hacemos el caso particular $f(X)=X^q-X$.
Se tiene que $\p$ se descompone totalmente en $K/k$
si y solamente si se descompone en todas las extensiones
intermedias de grado $p$: $y_{\mu}^p-y_{\mu}=\mu u$ donde
$y_{\mu}=h(\mu y)=(\mu y)^{p^{n-1}}+\cdots +(\mu y)^p + 
(\mu y)$, $\mu\in {\ma F}_q^{\ast}$. Esto \'ultimo es equivalente a que
$\mu R(T)\in\im\wp({\ma F}_{q^m})$ para toda $\mu\in\F$.
Por el Lema \ref{L6.2'.pea} esto \'ultimo equivale 
a que $R(T)=\lambda^q-\lambda$
para alg\'un $\lambda\in{\ma F}_{q^m}$. Puesto que $u$
est\'a en forma reducida, se sigue que $R(T)=0$.

Para el caso general, por el Teorema \ref{T5.3.pea}, $\p$ se descompone
totalmente en $K/k$ si y solamente si $\p$ se descompone en todas
las subextensiones de grado $p$ dadas por 
\begin{equation}\label{infinito.pea}
z_{{\mc H}}^p-z_{{\mc H}}=\frac{u}{f_{{\mc H}}(\varepsilon_{{\mc H}})^p},
\end{equation}
para todo hiperplano ${\mc H}$ de ${\mc G}_f$. El t\'ermino de
(\ref{infinito.pea}) que determina el comportamiento de $\p$
es $\frac{R(T)}{f_{{\mc H}}(\varepsilon_{{\mc H}})^p}$. Por tanto
$\p$ se descompone totalmente en (\ref{infinito.pea}) si y solamente
si $R(T)\in\im \wp_{f_{{\mc H}}(\varepsilon_{{\mc H}})}$. 
El resultado es ahora consecuencia
inmediata de la ecuaci\'on (\ref{Eq*.pea}). $\fin$
\end{proof}

A continuaci\'on presentamos un resultado 
an\'alogo a la Proposici\'on \ref{P6.1.pea}
para cierto tipo de $p$-extensiones 
elementales abelianas. Obs\'ervese que en 
el caso en el que el primo infinito es ramificado, 
no se precisa el comportamiento del primo 
en la extensi\'on. 
 
\begin{corolario}\label{C6.2'''.pea}
Sean $k=k_0(T)$ un campo de funciones 
racionales y $K = k(y)$ una $p$-extensi\'on 
elemental abeliana de $k$ dada por una forma reducida 
\begin{gather*}
y^q-y=u=\sum_{i=1}^r\frac{Q_i(T)}{P_i(T)^{\alpha_i}} + R(T)
\end{gather*}
satisfaciendo las condiciones del Teorema {\rm{\ref{T5.2.pea}}}.

Entonces el divisor primo $\p$ es
\l
\item totalmente descompuesto si $R(T)=0$.
\item no ramificado con grado de inercia $p$ si 
$R(T)\in k_0$ y $R(T)\not\in\{\lambda^q - \lambda
 | \lambda \in k_0 \}$.
\item ramificado si $R(T)\not\in k_0$. $\fin$
\end{list}
\end{corolario}

Sea $K=k(y)$ dado por (\ref{Eq5.2.pea}) y sean $\P_i$ los divisores
primos asociados a
$P_i(T)$, $1\leq i\leq r$. Sea $\P$ un divisor primo tal que $\P\notin\{\P_1,
\ldots,\P_r,\p\}$. Entonces $\P$ es no ramificado y por tanto
su grupo de descomposici\'on es c\'iclico pues corresponde
al grupo de Galois de una extensi\'on de campos finitos,
a saber, de los campos residuales. Por tanto $\P$
o bien es totalmente descompuesto en $K/k$
o tiene grado de inercia $p$. Los siguiente resultados 
nos dan el tipo de descomposici\'on de $\P$.

Primero estudiaremos el caso Artin--Schreier.

\begin{proposicion}\label{P6.3.pea} 
Sean $\P$ y $K/k$ como antes. Escribamos
la ecuaci\'on {\rm{(\ref{Eq6.1.pea})}} como
\[
y^p-y=u(T),
\]
con $u(T)=\frac{g(T)}{h(T)}\in k=k_0(T)$ tal que
$\mcd(g(T),h(T))=1$. Sea 
$P(T)\in R_T^+$ el polinomio irreducible asociado
a $\P$, digamos que $\deg P(T)=m$. 
Sea $\nu\in k'$ una ra\'iz de $P$, donde
$k'$ es el campo de descomposici\'on del polinomio $P(T)$.
Entonces $\P$ se descompone
totalmente en $K/k$ si y solamente si $u(\nu)\in \wp(k')$.
\end{proposicion}

\begin{proof}. Se tiene $[k':k]=m$. Sean
$\nu=\nu_1,\ldots,\nu_m$ las ra\'ices de $P$ en $k'$, $P(T)=\prod_{i=1}^m
(T-\nu_i)$. Se tiene que $\P$ se descompone totalmente en $k'(T):=k_m/k=
k_0(T)$. Aqu\'i $k_m$ denota la extensi\'on de constantes de $k$ de
grado $m$. Puesto que estamos suponiendo que la extensi\'on $K/k$ es geom\'etrica,
se tiene que $k_m\cap K=k$ y por tanto $\P$ se descompone totalmente en $K/k$
si y solamente si ${\mc Q}$ se descompone totalmente en $K_m/k_m$ donde $K_m=
Kk_m$ y ${\mc Q}$ es un primo en $k_m$ sobre $\P$.
\[
\xymatrix{K\ar@{-}[rr]\ar@{-}[d]&&K_m\ar@{-}[d]\\
k\ar@{-}[rr]_{\substack{\text{$\P$ es totalmente}
\\ \text{descompuesto}}}&&k_m}
\]

Digamos que ${\mc Q}$ es el primo asociado a $T-\nu\in k_m$. 
Puesto que $v_{\P}(u(T))\geq 0$ se tiene que $v_{{\mc Q}}(u(T))\geq0$ por lo que
$T-\nu\nmid h(T)$ y $h(\nu)\neq 0$. Adem\'as $g(\nu)=0\iff v_{\P}(u(T))>0$.
Se tiene que $\deg_{k_m}
{\mc Q}=1$. Hacemos a ${\mc Q}$ el primo infinito en $k_m$, esto es, sea $T'=\frac{1}{T-\nu}$,
$(T')_{k_m}=\frac{{\mc Q}'_0}{{\mc Q}'_{\infty}}=\frac{{\mc Q}_{\infty}}{{\mc Q}}$, donde $(T)_{k_m}=\frac{
{\mc Q}_0}{{\mc Q}_{\infty}}$. Se tiene $T=\frac{1}{T'}+\nu$.

Escribimos $u_1(T'):=u(T)=u\big(\frac{1}{T'}
+\nu\big)=\frac{g\big(\frac{1}{T'}+\nu\big)}
{h\big(\frac{1}{T'}+\nu\big)}=\frac{g\big(\frac{1}{T'}(1+T'\nu)\big)}
{h\big(\frac{1}{T'}(1+T'\nu)\big)}$.

Sean $g(T)=a_sT^s+a_{s-1}T^{s-1}+\cdots+
a_1T+a_0$, $a_s\neq 0$, $a_i\in k_0$, 
$0\leq i\leq s$; $h(T)=b_tT^t+b_{t-1}T^{t-1}+\cdots+
b_1T+b_0$, $b_t\neq 0$, $b_j\in k_0$, 
$0\leq j\leq t$. Entonces
\begin{gather*}
g\Big(\frac{1}{T'}(1+T'\nu)\Big)=\frac{1}{(T')^s}\Big(a_s+\cdots+g(\nu)(T')^s\Big)=
\frac{1}{(T')^s}g_1(T');\\
h\Big(\frac{1}{T'}(1+T'\nu)\Big)=\frac{1}{(T')^t}\Big(b_t+\cdots+h(\nu)(T')^t\Big)=
\frac{1}{(T')^t}h_1(T').
\intertext{Se sigue que}
\deg_{T'}g_1(T')\leq s;\quad\deg_{T'}h_1(T')=t.
\intertext{Por tanto}
\deg_{T'}u_1(T')=\deg_{T'}g_1(T')-s\leq 0,\quad\text{y}\quad v_{{\mc Q}}(u_1(T'))=s-\deg_{
T'}(g_1(T'))\geq 0.
\end{gather*}

La forma reducida de $u_1(T')$ es
\[
u_1(T')=\sum_{j=1}^{r'}\frac{Q'_j(T')}{(P'_j(T'))^{e'_j}}+u(\nu).
\]
Por tanto ${\mc Q}$ se descompone en $K_m/k_m\iff u(\nu)\in 
\wp(k')$. Esto demuestra la proposici\'on. $\fin$
\end{proof}

El caso general es consecuencia de la Proposici\'on \ref{P6.3.pea}. De hecho, 
sea $K=k(y)$ donde $y$ est\'a dado por la forma reducida
(\ref{Eq5.2.pea}). Suponemos $\P$
es no ramificado en $K/k$. Por el Teorema \ref{T5.3.pea}, usando
la notaci\'on de ah\'i, se tiene que todas las
subextensiones de $K$ de grado $p$ sobre $k$ est\'an dadas
por (\ref{infinito.pea}).

Entonces como consecuencia de estas expresiones y de las
Proposiciones \ref{P6.2''.pea} y \ref{P6.3.pea}, obtenemos el 
resultado principal de esta secci\'on.

\begin{teorema}\label{T6.4.pea}
Con las notaciones anteriores, sea $\P$ un divisor primo de grado $m$ no
ramificado en $K/k$. Entonces se tiene que $\P$ se descompone
totalmente en $K/k\iff$ para todo hiperplano ${\mc H}$ se cumple
\begin{gather*}
u(\nu) \frac{1}{f_{{\mc H}}(\varepsilon_{{\mc H}})^p}\in \wp (k'),
\intertext{y $\P$ tiene grado de inercia $p$ en $K/k \iff$ existe un hiperplano ${\mc H}$ tal que}
u(\nu)\frac{1}{f_{{\mc H}}(\varepsilon_{{\mc H}})^p}\notin \wp (k').
\end{gather*}

Equivalentemente, si {\rm{(\ref{Eq*.pea})}} se satisface, entonces
se tiene que $\P$ se descompone totalmente en $K/k
\iff u(\nu)\in f(k')$. $\fin$
\end{teorema}

\section{Generaci\'on de $p$--extensiones 
elementales abelianas}\label{S7.pea}

Sea $k$ un campo arbitrario de caracter\'istica $p>0$ y sea
$f(X)\in k[X]$ un polinomio aditivo dado por (\ref{E2.pea}) y donde
${\mc G}_f\subseteq k$. Sea
$u\in k$ tal que $F(X)=f(X)-u\in k[X]$ es irreducible.
Sea $K=k(y)$ donde $f(y)=u$,  $\Gal(K/k)
\cong {\mc G}_f$. Se tiene que $K/k$ tiene $\frac{(p^n-1)(p^n-p)\cdots (p^n-p^{m-1})}
{(p^m-1)(p^m-p)\cdots (p^m-p^{m-1})}$ subextensiones $k\subseteq E
\subseteq K$ tales que $[E:k]=p^m$.

Queremos ver qu\'e relaci\'on satisfacen $y$ y $z$ si 
$K=k(y)=k(z)$ y $f(z)=\chi\in k$. 
En el caso de extensiones de Artin--Schreier la relaci\'on
la conocemos por el Corolario \ref{C.9'.2.2}. El siguiente
resultado generaliza el Corolario \ref{C.9'.2.2}. Este
es el resultado principal de esta secci\'on.

\begin{teorema}\label{T7.1.pea} 
Con las notaciones anteriores, los siguientes
enunciados son equivalentes:
\l
\item $k(y)=k(z)$,

\item existen $A_{n-1},A_{n-2},\ldots, A_1,A_0
\in {\mc G}_f$ que satisfacen que 
\begin{gather*}
\elemental A{\beta}n=0\\
\text{con}\quad \beta \in  {\mc G}_f \iff \beta=0
\end{gather*}
y $D\in k$ tales que
\begin{equation}\label{Eq7.1.pea}
z=\elemental Ayn+D.
\end{equation}
\end{list}
\end{teorema}

El Teorema \ref{T7.1.pea} es consecuencia inmediata
del siguiente teorema m\'as general.

\begin{teorema}\label{T7.2.pea} Sea $K=k(y)$. Entonces
las siguientes condiciones son equivalentes:
\l
\item $E=k(z)$ con $k\subseteq E\subseteq K$,
$[E:k]=p^m$ y donde $g(z)=\chi\in k$ para alg\'un
$\chi\in k$ y para alg\'un polinomio aditivo $g(X)$ tal que
$g(X)\mid f(X)$, esto es, $g=f_V$ para un subgrupo aditivo
de ${\mc G}_f$ de dimensi\'on $m$ sobre ${\ma F}_p$,

\item
existen $A_{n-1}, A_{n-2}, \ldots,A_1,A_0\in {\mc G}_f$, 
$C\in k$  y un ${\ma F}_p$--subespacio vectorial
${\mc H}$ de ${\mc G}_f$ de dimensi\'on $n-m$ tales que 
\las
\item $z$ satisface
\begin{multline}\label{Eq7.2.pea}
z=\\
\elemental Ayn +C.
\end{multline}

\item
para $\beta\in {\mc G}_f$,
\begin{gather*}
\elemental A{\beta}n=0\\ \iff \beta\in {\mc H}.
\end{gather*}
\end{list}
\end{list}

La relaci\'on entre $E$ y  {\rm{(\ref{Eq7.2.pea})}} est\'a dada de la
siguiente forma. Si
${\mc H}={\mc L}_{{\ma F}_p}\{\mu_{m+1},\cdots,\mu_{n}\}$ donde $\{\mu_1,\ldots,
\mu_n\}$ es la base de ${\mc G}_f$ tal que si $G=\langle \sigma_1,
\ldots,\sigma_n\rangle$, $\sigma_i (y)=y+\mu_i$, $1\leq i\leq n$, 
entonces $E$ es el campo fijo por el subgrupo 
$H:=\langle \sigma_{m+1},
\ldots,\sigma_{n}\rangle$ de $G$. Esto es, $H$ corresponde
a ${\mc H}$ bajo el isomorfismo dado en la 
Proposici\'on {\rm{\ref{P2.2.pea}}} y si ${\mc G}_f={\mc H}\oplus V$, $g(X)=
f_V(X)=\prod_{\delta\in V}(X-\delta)\mid f(X)$.
\end{teorema}

\begin{proof}
Sea $G=\Gal(K/k)=\langle \sigma_1,\ldots, \sigma_n\rangle$ con
$\sigma_i(y)=y+\mu_i$, donde 
$\mu_i\in{\mc G}_f$ y $\{\mu_1,\ldots,\mu_n\}$ es una
$\ma F_p$--base de ${\mc G}_f$.
M\'as precisamente, se tiene que si $\sigma_1,\ldots,\sigma_n\in G$
con $\sigma_i(y)=y+\mu_i$, entonces $G=\langle\sigma_1,\ldots,
\sigma_n\rangle\iff \{\mu_1,\ldots,\mu_n\}$ es base de ${\mc G}_f/{\ma F}_p$.
Notemos que para $0\leq \alpha_i\leq
p-1$, $1\leq i\leq n$, $\sigma=\sigma_1^{\alpha_1}\cdots\sigma_n^{
\alpha_n}$ se satisface $\sigma(y)=y+\sum_{i=1}^n\alpha_i\mu_i$.

Primero consideremos un subcampo $E$ de $K$ de grado $p^m$
sobre $k$. Podemos seleccionar un conjunto de $n$ generadores
de $G=\Gal(K/k)=\langle\sigma_1,\ldots,\sigma_n\rangle$ de tal
forma que $E=K^{\langle\sigma_{m+1},\ldots,\sigma_n\rangle}$ es
el campo fijo bajo $H=\langle\sigma_{m+1},\ldots,\sigma_n\rangle$. Se tiene
que $\Gal(K/E)= \langle\sigma_{m+1},\ldots,\sigma_n\rangle$.

Sea $\theta\colon G\to{\mc G}_f$ el isomorfismo dado por $\sigma_i
\mapsto \mu_i$, $1\leq i\leq n$. Sea
${\mc H}=\theta(H)<{\mc G}_f$ y sea $V$ cualquier secci\'on de la sucesi\'on
exacta
\[
0\lra {\mc H}\stackrel{i}{\lra}{\mc G}_f\stackrel{\pi}{\lra}{\mc G}_f/{\mc H}\lra 0,
\]
es decir, $V=\varphi({\mc G}_f/{\mc H})<{\mc G}_f$ donde $\varphi\colon {\mc G}_f/{\mc H}
\lra{\mc G}_f$ satisface $\pi\circ\varphi =\Id_{{\mc G}_f/{\mc H}}$. Se tiene ${\mc G}_f\cong
{\mc H}\oplus V$ como ${\ma F}_p$--espacios vectoriales. 

Por el Teorema \ref{T5.4.pea}, tenemos que existen $z\in E$ y $\chi\in k$ 
tales que $E=k(z)$ con $f_V(z)=\chi\in k$ y $f_V(X)=\prod_{\delta\in V}
(X-\delta)\mid \prod_{\delta\in{\mc G}_f}(X-\delta)=f(X)$.

Se tiene que $\Gal(E/k)\cong\langle\bar{\sigma}_1,\ldots,
\bar{\sigma}_m
\rangle$ donde $\bar{\sigma}_i=\sigma_i|_E$ o $\bar{\sigma}_i=
\sigma_i\bmod \Gal(K/E)$. 

Ahora bien, sea $\sigma_i(z)=z+\gamma_i$, $1\leq i\leq m$, donde 
$\{\gamma_1,\ldots,\gamma_m\}$ es una base de $V$ y
$\sigma_j(z)=z$ para $m+1\leq j\leq n$. Por notaci\'on
ponemos $\gamma_j=0$ para $m+1\leq j\leq n$.

Sean $A_{n-1},A_{n-2},\ldots, A_1, A_0\in {\mc G}_f$ arbitrarios y sea
\begin{equation}\label{Eq7.4.pea}
w:=\elemental Ayn.
\end{equation}
Esto es, si denotamos $l(X)=\elementalito AXn$, entonces $w=l(y)$.

Probemos que existen $A_{n-1},A_{n-2},\ldots, A_1, A_0
\in {\mc G}_f$ y $D\in k$ tales que 
\begin{equation}\label{Eq7.3.pea}
z=w+D.
\end{equation}

Se tiene
\begin{equation}\label{Eq7.5.pea}
\sigma_i(w)=\sigma_i(l(y))=l(\sigma_i(y))=l(y+\mu_i)=l(y)+l(\mu_i)=w+l(\mu_i).
\end{equation}

Se sigue que
\begin{gather}
\sigma_i(w)=w+\gamma_i,\quad 1\leq i\leq n\iff l(\mu_i)=\gamma_i,
\quad 1\leq i\leq n \label{Eq7.6.pea}\\
\iff M\left[\begin{array}{c}A_0\\A_1\\ \vdots\\ A_{n-2}\\
A_{n-1}\end{array}\right]=\left[\begin{array}{c}\gamma_1\\ \gamma_2\\ \vdots\\ 
\gamma_{n-1}\\ \gamma_n\end{array}\right]=
\left[\begin{array}{c}\gamma_1\\ \vdots\\ 
\gamma_{m}\\0 \\ \vdots \\ 0\end{array}\right], \nonumber
\end{gather}
donde $M$ es la matriz
\[
M=\matriz{\mu_1}\mu.
\]
A $M$ se le conoce como {\em matriz de Moore\index{matriz de
Moore}}.

Veamos que $M$ es no singular. Sea
\[
B(X):=\matriz X\mu=\left[\begin{array}{c}F(X)\\F(\mu_2)\\ \vdots\\ 
F(\mu_{n-1})\\F(\mu_n)\end{array}\right],
\]
donde $F(Z):=[Z\quad Z^p\ \cdots\  Z^{p^{n-2}}\quad Z^{p^{n-1}}]$ 
con $Z\in\{X,\mu_2,\ldots,\mu_n\}$ denota
las filas de $B(X)$. Se tiene que
$B(\mu_1)=M$ y $\det B(X)$ es un polinomio aditivo en $k[X]$ de
grado $p^{n-1}$.

Sean $(i_2,\ldots,i_n)\in \ma F_p^{n-1}$ y $\xi=i_2\mu_2+
\cdots+i_n\mu_n$. Entonces:
\[
B(\xi)=\left[\begin{array}{c}F(i_2\mu_2+\cdots+i_n\mu_n)\\F(\mu_2)\\ \vdots\\ 
F(\mu_{n-1})\\F(\mu_n)\end{array}\right]=
\left[\begin{array}{c}i_2F(\mu_2)+\cdots+i_nF(\mu_n)\\F(\mu_2)\\ \vdots\\ 
F(\mu_{n-1})\\F(\mu_n)\end{array}\right].
\]
Por tanto $\det B(\xi)=0$ para todo $\xi\in\{i_2\mu_2+\cdots +
i_n\mu_n\mid i_2,\ldots,i_n\in\ma F_p\}=C$. Puesto que $\{\mu_2,\ldots,
\mu_n\}$ es un conjunto linealmente independiente sobre ${\ma F}_p$,
se tiene que $|C|=p^{n-1}=\deg B(X)$. De esta forma, tenemos que
$C$ es el conjunto de ra\'ices de $\det B(X)$. En particular, puesto
que $\mu_1\notin C$, $\det B(\mu_1)=\det M\neq 0$ y $M$ es no singular.

Por tanto (\ref{Eq7.6.pea}) tiene una soluci\'on \'unica:
\begin{equation}\label{Eq7.7.pea}
\left[\begin{array}{c}A_0\\A_1\\ \vdots\\ A_{n-2}\\
A_{n-1}\end{array}\right]=M^{-1}\left[\begin{array}{c}\gamma_1\\ \vdots\\ 
\gamma_{m}\\0 \\ \vdots \\ 0\end{array}\right].
\end{equation}

Sea $\beta=\sum_{i=1}^n c_i\mu_i\in {\mc G}_f$ con $c_i\in{\ma F}_p$,
$1\leq i\leq n$.
Por tanto 
\[
l(\beta)=l(\sum_{i=1}^nc_i\mu_i)=\sum_{i=1}^nl(c_i\mu_i)=
\sum_{i=1}^nc_il(\mu_i)=\sum_{i=1}^nc_i\gamma_i=\sum_{i=1}^mc_i\gamma_i.
\]
Se sigue que $l(\beta)=0\iff c_1=\ldots=c_m=0\iff \beta\in{\mc L}_{{\ma F}_p}\{\mu_{m+1},
\ldots,\mu_n\}={\mc H}$.

Finalmente, se tiene que $\sigma_i(z-w)=z-w$ para todo $1\leq i\leq n$, por lo que
$z-w=D\in k$ y $z$ est\'a en la forma (\ref{Eq7.2.pea}).

Para probar el rec\'iproco, sea $z$ dado por (\ref{Eq7.2.pea}), $z=l(y)+D$. Entonces
\[
\sigma_i(z)=\sigma_i(l(y)+D)=l(\sigma_i(y))+D=l(y+\mu_i)+D=
l(y)+l(\mu_i)+D=z+l(\mu_i)
\]
y se tiene $l(\mu_i)=0\iff i\geq m+1$. Por tanto $k(z)\subseteq K^{\langle
\sigma_{m+1},\ldots,\sigma_n\rangle}$. Ahora bien, para cualesquiera $c_1,\ldots,
c_m\in {\ma F}_p$, no todos cero, $\sigma_1^{c_1}\cdots\sigma_m^{c_m}(z)
=z+l(\beta)$ con $\beta=\sum_{i=1}^m c_i\mu_i\neq 0$, $l(\beta)\neq 0$, lo
que implica que $[k(z):k]\geq p^m$.
Se sigue que $[k(z):k]=p^m$ y que $k(z)=K^{\langle \sigma_{m+1},\ldots,
\sigma_n\rangle}$.

Sean $\xi_i:=l(\mu_i)$, $1\leq i\leq m$, $V={\mc L}_{{\ma F}_p}\{\xi_1,\ldots,\xi_m\}$ y $f_V(X)\mid
f(X)$. Entonces 
\begin{gather*}
f_V(z)=f_V(l(y)+D)=f_V\big(\sum_{i=0}^{n-1}A_iy^{p^i}+D\big)=
\sum_{i=0}^{n-1} f_V(A_iy^{p^i})+f_V(D)=\chi.
\end{gather*}

Veamos que $\chi=f_V(z)\in k$.

Se tiene que si $\sigma:=\sigma_1^{c_1}\cdots\sigma_n^{c_n}\in G=
\Gal(K/k)$ entonces para
 $\mu:=\sum_{i=1}^n c_i\mu_i$, $\sigma(y)=y+\mu$, se tiene
\begin{align*}
\sigma(f_V(z))&=\sigma(f_V(l(y)+D))=f_V(l(y+\mu)+D)=f_V(l(y)+l(\mu)+D)\\
&=f_V(l(y)+D)+f_V(l(\mu))=f_V(z)+f_V(l(\mu)).
\end{align*}
Finalmente,
\[
l(\mu)=\sum_{i=1}^n c_il(\mu_i)=\sum_{i=1}^m c_i\xi_i\in V,
\]
por lo que $f_V(l(\mu))=0$ y $\sigma(f_V(z))=f_V(z)$ para todo $\sigma\in
\Gal(K/k)$. Se sigue que $\chi=f_V(z)\in k$ y esto termina la demostraci\'on.
$\fin$
\end{proof}

\section{Existencia de $p$--extensiones c\'iclicas}

Como mencionamos en la introducci\'on, Artin
y Schreier profundizaron su m\'etodo y encontraron las extensiones
c{\'\i}clicas de grado $p^2$. Su t\'ecnica fue el primer paso para
lo que vendr{\'\i}a en los siguientes a\~nos: los resultados de Albert,
Schmid y Witt.

Veamos como podemos construir extensiones c{\'\i}clicas de
grado $p^n$ en caracter{\'\i}stica $p$. Sea $G$ un grupo c{\'\i}clico
de orden $p^n$. Sea $K/k$ una extensi\'on c{\'\i}clica de grado
$p^{n-1}$. Sean $\langle\varphi\rangle=
\Gal(K/k)$, $o(\varphi)=p^{n-1}$, $\chi\in {\ma F}_p^{\ast}=
\{1,\ldots,p-1\}$ y $L=K(\theta)$ una extensi\'on c{\'\i}clica
de grado $p$ tal que $L/k$ es una extensi\'on
c{\'\i}clica de grado $p^n$. Sea
$\langle \sigma\rangle=\Gal(L/k)$, $o(\sigma)=p^n$ tal que
$\varphi=\sigma\bmod \Gal(L/K)$, es decir, $\varphi=\sigma_{\mid K}$.
Sea $\sigma^{p^{n-1}}=\psi$, $\langle\psi\rangle=\Gal (L/K)$.

$
\d \left.
\begin{array}{cl}
L\\
\biggl|\biggr. & \langle \psi\rangle\\
K\\
\biggl| & \langle\varphi\rangle\biggr.\\
k
\end{array}\right\}\langle\sigma\rangle
$
\quad Ahora, $L/K$ es una extensi\'on c{\'\i}clica de grado $p$, es decir,
una extensi\'on de Artin--Schreier. As{\'\i}, podemos escoger
$\theta$ tal que $\wp\theta=\theta^p-\theta=\gamma\in K$ y
tal que $\psi\theta=\theta+\chi$ o, equivalentemente,
$(\psi-1)\theta=\chi$. Se tiene que
\begin{gather*}
(\psi-1)(\sigma-1)\theta=(\sigma-1)(\psi-1)\theta=(\sigma-1)\chi=0,\\
\intertext{es decir, $\delta:=(\sigma-1)\theta \in K$. Adem\'as}
(\varphi-1)\gamma=(\varphi-1)(\wp\theta)=(\sigma-1)(\wp\theta)=
\wp((\sigma-1)\theta)=\wp \delta,\\
\intertext{y se tiene}
\begin{align*}
\Tr_{K/k}\delta&=\sum_{i=0}^{p^{n-1}-1}\varphi^i\delta=
\frac{\varphi^{p^{n-1}}-1}{\varphi-1}\delta=
\frac{\sigma^{p^{n-1}}-1}{\sigma-1}(\sigma-1)\theta\\
&=(\sigma^{p^{n-1}}-1) \theta =(\psi-1)\theta =\chi.
\end{align*}
\end{gather*}

En resumen, si $L/k$ es una extensi\'on c{\'\i}clica de grado $p^n$
que contiene a $K$, entonces existen $\theta\in L$, $\gamma\in K$,
$\chi\in\{1,2,\ldots,p-1\}$ y $\delta\in K$ tales que si $\sigma_{\mid K}
=\varphi$, $\sigma_{p^{n-1}}=\psi$, entonces
\lasa
\item $\wp\theta=\gamma$,
\item $(\psi-1)\theta=\chi$,
\item $(\sigma-1)\theta=\delta$,
\item $(\varphi-1)\gamma=\wp \delta$,
\item $\Tr_{K/k}\delta =\chi$.
\end{list}

Rec{\'\i}procamente, sea $K/k$ una extensi\'on
c{\'\i}clica de grado $p^{n-1}$.
Como $K/k$ es una extensi\'on separable, existe
$\delta\in K$ con $\Tr_{K/k}\delta =\chi$, $\chi\in {\ma F}_p^{\ast}$.
 En particular,
\[
\Tr_{K/k}\wp(\delta)=\wp(\Tr_{K/k}\delta)=\wp(\chi)=\chi^p-\chi=\chi-\chi=0.
\]
Por el Teorema 90 de Hilbert, existe $\gamma\in K$ tal que
$(\varphi-1)\gamma=\wp\delta$, donde $\langle\varphi\rangle=\Gal(K/k)$.
Si $\delta'$ es cualquier otro elemento tal que $\Tr_{K/k}\delta^{\prime}=\chi$, 
entonces $\Tr_{K/k}(\delta^{\prime}-\delta)=0$. Nuevamente por
el Teorema 90 de Hilbert, existe $\alpha\in K$ tal que
$\delta^{\prime}-\delta=(\varphi-1)\alpha$. Se tiene
$\delta^{\prime}=\delta+
(\varphi-1)\alpha$. Al sustituir $\delta$ por $\delta'=\delta + (\varphi-1)
\alpha$, tenemos 
\begin{align*}
(\varphi-1)(\gamma+\wp\alpha)&=(\varphi-1)\gamma+\wp((\varphi-1)\alpha)
=\wp\delta+\wp((\varphi-1)\alpha)\\
&=\wp(\delta+(\varphi-1)\alpha)=\wp\delta'.
\end{align*}
Es decir, la sustituci\'on $\delta\leftrightarrow \delta+(\varphi-1)\alpha$
corresponde a la sustituci\'on $\gamma\leftrightarrow \gamma +
\wp\alpha$, $\alpha\in K$.

Veamos que $\gamma\neq \wp\beta$ para $\beta\in K$. En caso contrario,
si $\gamma=\wp \beta$ para alg\'un
$\beta\in K$, cambiando $\gamma$ por $\gamma'=\gamma
-\wp \beta=0$, se tendr{\'\i}a $(\varphi-1)\gamma'=0=\wp\delta'$, esto es,
$\delta'\in{\ma F}_p$ y $\Tr_{K/k}\delta'=0=\chi$ lo cual es absurdo.
Esto prueba que $\gamma\notin \wp(K)$.

Sea $\theta$ una soluci\'on de la ecuaci\'on $x^p-x-\gamma\in K[x]$,
es decir, $\wp\theta=\gamma$, $\theta\notin K$. Sea $L=K(\theta)$.
Se tiene $[L:k]=p^n$. Ahora, $\wp(\theta+\delta)=\wp\theta+\wp\delta=
\gamma+(\varphi-1)\gamma = \varphi\gamma$. Sea
\begin{gather*}
\sigma\colon L\to L\quad \text{definida por}\quad
 \sigma\theta=\theta +\delta \quad \text{y}\quad \sigma_{\mid K}
=\varphi.\\
\intertext{Se tiene que}
\begin{align*}
\sigma^{p^{n-1}}\theta-\theta &=(\sigma^{p^{n-1}}-1)\theta =
\Big(\sum_{i=0}^{p^{n-1}-1}\sigma^i\Big)(\sigma-1)\theta\\
&=\Big(\sum_{i=0}^{p^{n-1}-1}\sigma^i\Big)\delta=
\sum_{i=0}^{p^{n-1}-1}\varphi^i\delta = \Tr_{K/k}\delta=\chi.
\end{align*}
\end{gather*}
Esto es, $\sigma^{p^{n-1}}\theta=\theta+\chi\neq \theta$ por lo que
$\sigma$ tiene orden $p^n$. As{\'\i}, $L/k$ es c{\'\i}clica generada por
$\sigma$. En resumen, tenemos

\begin{teorema}[Witt \cite{Wit36-1}]\label{T9'.2.4}
Sea $K/k$ una extensi\'on c{\'\i}clica de grado $p^{n-1}$, $n\geq 2$.
Entonces para construir cualquier
extensi\'on c{\'\i}clica $L/k$ de grado $p^n$ que contenga
a $K$, se eligen de manera arbitraria los siguiente:
\l
\item Un generador $\varphi$ de $\Gal(K/k)$.
\item Un elemento $\chi\neq 0$ en ${\ma F}_p$, es decir, $\chi\in
\{1,2,\ldots,p-1\}$.
\item Una soluci\'on $\delta\in K$ de la ecuaci\'on $\Tr_{K/k}
\delta=\chi$.
\item Una soluci\'on $\gamma\in K$ de
la ecuaci\'on $(\varphi-1)\gamma =\wp\delta$.
\end{list}
La extensi\'on $L$ se obtiene como $L=K(\theta)$ donde
$\wp\theta=\gamma$. Cualquier otra extensi\'on de este tipo se
obtiene sustituyendo $\gamma$ por $\gamma +c$ con $c\in k$. $\fin$.
\end{teorema}

Este es el resultado clave usado por Schmid para generar
extensiones c{\'\i}clicas de grado $p^n$ en caracter{\'\i}stica $p$.

\subsection{La construcci\'on de Schmid}\label{S9'.3}

Sea $k$ un campo arbitrario de caracter{\'\i}stica $p$. Sea $K_i$ una
extensi\'on c{\'\i}clica de $k$ de grado $p^i$ con $\Gal(K_i/k)=\langle
\varphi_i\rangle$, $i=1,2,\ldots,n$. Sea $T_i:=\Tr_{K_i/k}$. Suponemos
$k\subseteqq K_1\subseteqq\ldots\subseteqq K_n=L$. Seleccionamos
$\chi=1$ en el resultado de Witt (Teorema \ref{T9'.2.4}).
Sea $c_i\in K_i$ un elemento tal que
$T_i c_i=1$ y como de costumbre sea $\wp x=x^p-x$. Sea $\Delta_i$
el operador $\varphi_i-1$. La extensi\'on $K_i/K_{i-1}$ est\'a dada
por $K_i=K_{i-1}(v_i)$, $i=2,3,\ldots $ con $\wp v_i=z_{i-1}\in K_{i-1}$
y $\Delta_{i-1}z_{i-1} =\wp c_{i-1}$. Se tiene que $\varphi_i(v_i)=
v_i+c_{i-1}$.

Consideremos el elemento $\alpha=-v_1^{p-1}\in K_1$. Sea $X_i:=
v_i+i$, $1\leq i\leq p$. Sean $\sigma_i$ las funciones sim\'etricas
elementales en $X_1,\ldots,X_p$:
\begin{gather*}
\sigma_0=1,\quad \sigma_1=
\sum_{i=1}^pX_i,\quad \sigma_2=\sum_{i<j}X_iX_j\quad,
\ldots,\quad \sigma_p=
X_1\cdots X_p \\
\intertext{y sean}
\rho_m:=X_1^m+\cdots+X_p^m\quad \text{para}\quad
m\geq 1 \quad \text{y}\quad \rho_0= p=0.\\
\intertext{Se tiene}
Y^p-Y-\beta_1=\prod_{i=1}^p(Y-(v_i+i)=\prod_{i=1}^p (Y-X_i)=
\sum_{i=0}^p (-1)^i\sigma_i Y^{p-i},
\end{gather*}
por lo que $(-1)^p\prod_{i=1}^pX_i=-\beta_1$, esto es,
$\sigma_1=\cdots=\sigma_{p-2}=0$,
$\sigma_{p-1}=-1$ y $\sigma_p=\beta_1$.

Por las identidades de Newton
\begin{gather*}
\rho_{p-1}-\rho_{p-2}\sigma_1+\cdots+(-1)^{p-2}\rho_1\sigma_{p-2}
+(-1)^{p-1}\sigma_{p-1}(p-1)=0\\
\intertext{se obtiene $\rho_{p-1}=(-1)^{p-1}\sigma_{p-1}=(-1)^p=-1$.
Por otro lado}
\rho_{p-1}=\sum_{i=1}^p(v_i+i)^{p-1}=
\Tr_{K_1/k}v_1^{p-1}=-1.
\end{gather*}
Por lo tanto $\Tr_{K_1/k}(-v_1^{p-1})=1$.
Este es el $\delta$ correspondiente al Teorema \ref{T9'.2.4} ({\sc iii}).

En el caso general consideramos $c_n:=(-1)^n\prod_{i=1}^n
v_i^{p-1}$, se tiene
\begin{align*}
\Tr_{K_n/k}c_n&=\Tr_{K_{n-1}/k}\Tr_{K_n/K_{n-1}}\big\{
(-1)^n \prod_{i=1}^n v_i^{p-1}\big\}\\
&=\Tr_{K_{n-1}/k}
\big((-1)^n\prod_{i=1}^{n-1}v_i^{p-1} \Tr_{K_n/K_{n-1}}
v_n^{p-1}\big)\\
&=\Tr_{K_{n-1}/k}\Big((-1)^n\prod_{i=1}^{n-1}v_i^{p-1} (-1)\Big)=
\Tr_{K_{n-1}/k}c_{n-1}.
\end{align*}

Por lo tanto, se sigue por inducci\'on que $\Tr_{K_n/k}c_n=1$ y
estos elementos sirven para la construcci\'on. Schmid 
construy\'o en general las extensiones y prob\'o que la extensi\'on
c{\'\i}clica $K_n/k$ en general est\'a dada por
\begin{gather}\label{Ec9'.2.1}
\begin{array}{lllll}
K_1&=k(v_1),& \wp v_1=\beta_1, & \Delta_1 v_1=1,\\
K_2&=K_1(v_2),& \wp v_2=z_1+\beta_2,& \Delta_2 v_2=c_1,
& \Delta_1z_1=\wp c_1\\
K_3&=K_2(v_3),& \wp v_3=z_2+\beta_3, & \Delta_3 v_3=c_2,
& \Delta_2z_2=\wp c_2\\
\hspace{0.3cm}\vdots &\hspace{.8cm}\vdots&\hspace{1cm} \vdots&\hspace{1cm}\vdots
&\hspace{1cm}\vdots\\
K_n&=K_{n-1}(v_n), & \wp v_n=z_{n-1}+\beta_n, & \Delta_n v_n=c_{n-1},
& \Delta_{n-1}z_{n-1}=\wp c_{n-1}
\end{array}
\end{gather}
donde $K_n$ est\'a determinada por los elementos $\beta_1,\ldots,\beta_n\in k$
arbitrarios, $\beta_1\notin \wp(k)$ y donde $\Delta_i=\varphi_i-1$, $\varphi_i(v_i)
=v_i+c_{i-1}$, $\langle\varphi_i\rangle=\Gal(K_i/k)$.

Las ecuaciones que encontr\'o Schmid para generar las extensiones c{\'\i}clicas
fueron reconocidas por Witt en forma vectorial y esto dio lugar a los {\em vectores
de Witt\index{vectores de Witt}\index{Witt!vectores de $\sim$}}.

\subsection{Vectores de Witt\index{vectores de
Witt}\index{Witt!vectores de $\sim$}}\label{S9'.4}

Sea $p$ un n\'umero primo fijo. Para un vector $\vec x=(x_1,x_2,\ldots )$
con una cantidad
a lo m\'as numerable de componentes $x_n$, en 
caracter{\'\i}stica $0$, se definen las
{\em componentes fantasmas\index{componentes fantasma}}
de $\vec x$ por
\begin{equation}\label{Eq9'.2.2}
x^{(t)}=x_1^{p^{t-1}}+px_2^{p^{t-2}}+\cdots+p^{t-1}x_t=\sum_{i=1}^t p^{i-1}x_i^{p^{t-i}},
\quad t=1,2, \ldots
\end{equation}

Rec{\'\i}procamente, $x_t$ puede ser calculado recursivamente como un polinomio
en $x^{(1)}, x^{(2)},\ldots, x^{(t)}$ a partir de (\ref{Eq9'.2.2}). Esta correspondencia puede
ser expresada como
\[
\vec x=(x_1,x_2,x_3,\ldots \mid x^{(1)},x^{(2)},x^{(3)},\ldots).
\]

{\em La suma\index{suma de Witt}\index{Witt!suma de $\sim$}} $\Witt +$,
{\em la diferencia\index{diferencia de Witt}\index{Witt!diferencia de $\sim$}} $\Witt -$ y
{\em el producto\index{producto de Witt}\index{Witt!productor de $\sim$}} $\Witt \times$
de Witt se definen por
\begin{equation}\label{Eq9'.2.3}
\vec x \Witt {\masmenos_{\times}} \vec y:=\big(?,?,\ldots\mid x^{(1)} {\masmenos_{\times}} y^{(1)}, 
x^{(2)} {\masmenos_{\times}} y^{(2)},\ldots\big).
\end{equation}

Esto es, las componentes fantasma se operan t\'ermino a t\'ermino y las
componentes usuales se calculan a partir de los resultados que se obtengan
en las componentes fantasma.

Notemos que $(\vec x)^{(n)}$ y $x^{(n)}$ denotan lo mismo, a saber, la
$n$--componente fantasma del vector $\vec x$. Similarmente 
$(\vec x)_n$ y $x_n$ denotan lo mismo, la $n$--componente
de $\vec x$.

Lo anterior puede precisarse de la siguiente forma. Consideremos tres familias
$\big\{x_i,y_j,z_{\ell}\big\}_{i,j,\ell=1}^N$ donde $N\in{\ma N}\cup \{\infty\}$, de
variables independientes sobre ${\ma Q}$
y consideramos el anillo $R={\ma Q}[x_i,y_j,z_{\ell}]_{i,j,\ell}$.
Sea $R^N$ el producto $\underbrace{R\times \ldots \times R\times \ldots}_{N}$. 
Por abuso del lenguaje,
denotamos por $R^N$ al anillo que como conjunto base tiene al mismo conjunto
$R^N$ y cuyas operaciones son t\'ermino a t\'ermino (esto corresponde a las
componentes fantasma) y sea $R_N$ el anillo que como conjunto sigue siendo
$R^N$ pero con las {\em operaciones de Witt\index{Witt!operaciones
de $\sim$}\index{operaciones de Witt}}: sea $\varphi\colon R_N\to R^N$ dado por
$\varphi(a_1,a_2,\ldots, a_N)=\big(a^{(1)},a^{(2)},\ldots,a^{(N)}\big)$ donde
\[
a^{(m)}:= a_1^{p^{m-1}}+pa_2^{p^{m-2}}+\cdots+p^{m-1}a_m,\quad m=1,2,\ldots, N.
\]
Se tiene que $\varphi$ es un mapeo biyectivo y el inverso
$\psi\colon R^N\to R_N$ 
est\'a dado por $\psi\big(a^{(1)},a^{(2)},\ldots,a^{(N)}\big)
=(a_1,a_2,\ldots,a_N)$ donde
\[
a_m:=\frac{1}{p^{m-1}}\big(a^{(m)}-a_1^{p^{m-1}}-pa_2^{p^{m-2}}-\cdots-p^{m-2}
a_{m-1}^p\big),\quad m=1,\ldots, N.
\]

Entonces las operaciones $\Witt +,\Witt -,\Witt \times$ sobre $R^N$ se definen por
\begin{equation}\label{Eq9'.3.4}
\vec a\Witt {\masmenos_{\times}}\vec b:= \big(\vec a^{\varphi} {\masmenos_{\times}}
\vec b^{\varphi}\big)^{\varphi^{-1}} =\big(\vec a^{\varphi} {\masmenos_{\times}}
\vec b^{\varphi}\big)^{\psi}.
\end{equation}

En otras palabras, dados dos vectores en $R_N$, los trasladamos a $R^N$ y 
ah{\'\i} los operamos de la manera usual, es decir, componente por componente
y al resultado lo volvemos a $R_N$. Como $R^N$ es conmutativo con unidad,
$R_N$ tambi\'en es conmutativo con unidad.
 Por ejemplo si $N=2$, entonces dados $\vec x, \vec y$
\begin{gather*}
\vec x=(x_1,x_2\mid x^{(1)},x^{(2)})=(x_1,x_2\mid x_1, x_1^p+px_2),\\
\vec y=(y_1,y_2\mid y^{(1)},y^{(2)})=(y_1,y_2\mid y_1, y_1^p+py_2),\\
\intertext{se tiene que}
\begin{align*}
\vec z&=\vec x\Witt +\vec y=(z_1,z_2\mid z^{(1)},z^{(2)})=(z_1,z_2\mid z_1, z_1^p+pz_2)\\
&=(?,?\mid x_1+y_1, x_1^p+px_2+ y_1^p+py_2).
\end{align*}
\end{gather*}
Esto es $z_1=x_1+y_1$, $z_1^p+pz_2=x_1^p+px_2+y_1^p+py_2$. Por lo tanto
\begin{align*}
z_2&=\frac{1}{p}\big(x_1^p+px_2+y_1^p+py_2-(x_1+y_1)^p\big)\\
&=\frac{1}{p}\big(px_2+py_2-\sum_{i=2}^{p-1}\binom{p}{i}x_1^iy_1^{p-i}\big)=
x_2+y_2-\sum_{i=2}^{p-1} \frac{1}{p}\binom{p}{i} x_1^iy_1^{p-i}.
\end{align*}
Por lo tanto
\[
\vec z=\vec x \Witt + \vec y= \big(x_1+y_1, x_2+y_2-\sum_{i=2}^{p-1}
\frac{1}{p}\binom{p}{i}x_1^iy_1^{p-i}\mid x_1+y_1, x_1^p+px_2+y_1^p+py_2\big).
\]

A continuaci\'on introducimos las siguientes operaciones en los vectores
de Witt que son de gran utilidad para obtener informaci\'on  de la naturaleza
del anillo $R_N$.

\begin{definicion}\label{D9'.3.1} Se define el {\em operador de corrimiento\index{operador
de corrimiento}} $V\colon R_N\to R_N$ por
\begin{gather*}
V(x_1,\ldots,x_n,\ldots)=(0,x_1,\ldots, x_n,\ldots),\\
V^i(x_1,\ldots,x_n,\ldots)=
(0,\ldots,0,\underbracket[0pt]{x_1}_{\substack{
\uparrow\\ i+1}},\ldots, x_n,\ldots), \quad i\in{\ma N},
\end{gather*}
y se define la {\em funci\'on componente\index{funci\'on componente}} $\{\ \}\colon
R\to R_N$,
\begin{gather*}
\{u\}:=(u,0,\ldots, 0,\ldots)=(u,0,\ldots,0,\ldots \mid u,u^p,u^{p^2},\ldots),\\
\{u\}^{(n)}=u^{p^{n-1}},\quad n\geq 1, \quad \text{y}\quad
 u_1=u,\quad u_n=0, \quad n\geq 2.
\end{gather*}
\end{definicion}

Notemos que $V^i(\{u\})=(0,\ldots, 0,
\underbracket[0pt]{u}_{\substack{\uparrow\\ i+1}},0,\ldots)$.

\begin{observacion}\label{O9'.3.2} Se tiene 
\begin{equation}\label{Eq9'.3.(c)}
(V\vec x)^{(n)}=px^{(n-1)}, \quad n=1,2,\ldots,  \quad \text{donde}\quad x^{(0)}=0.
\end{equation}
En efecto $(Vx)^{(n)}=(0,x_1,\ldots,x_m,\ldots)^{(n)}=0^{p^{n-1}}+px_1^{p^{n-2}}+
\cdots+p^{n-1} x_{n-1}$ y $x^{(n-1)}=x_1^{p^{n-2}}+px_2^{p^{n-3}}+\cdots+
p^{n-2}x_{n-1}$, de donde se sigue la igualdad.
\end{observacion}

En otras palabras tenemos
\begin{equation}\label{Eq9'.3.6}
V\vec x=(0,x_1,x_2,\ldots\mid 0,px^{(1)},px^{(2)},\ldots).
\end{equation}
Para $s=0,1,2,\ldots$ se tiene
\begin{equation}\label{Eq9'.3.7}
V^s\vec x=(\underbrace{0,\ldots,0}_s,x_1,x_2,\ldots\mid 0,\ldots,0, p^sx^{(1)},
p^sx^{(2)},\ldots)
\end{equation}
donde $V^0=\Id$, es decir, $V^0\vec x=\vec x$ para toda $\vec x\in R^N$.
En particular
\[
(V^s \vec x)^{(n)}=
\begin{cases}
0&\text{si $s\leq n$}\\
p^sx^{(n-s)}&\text{si $n\geq s+1$}
\end{cases}
=p^sx^{(n-s)}
\]
donde $x^{(1-s)}=\cdots=x^{(-1)}=x^{(0)}=0$.
Esta \'ultima igualdad puede ser verificada por inducci\'on en $s$:
\begin{align*}
(V^s\vec x)^{(n)}&=(V(V^{s-1}\vec x))^{(n)}=p(V^{s-1}\vec x)^{(n-1)}\\ 
&=p(p^{s-1}x^{(n-1-(s-1))})=p^s x^{(n-s)},\\
x^{(1-s)}&=x^{(2-s)}=\cdots=x^{(-1)}=x^{(0)}=0.
\end{align*}
Aplicado lo anterior a $\{u\}$ se tiene la igualdad
\[
(V^s\{u\})^{(n)}=p^s\{u\}^{(n-s)}=
\begin{cases}
p^s u^{p^{n-s-1}}&\text{si $n\geq s+1$}\\
0&\text{si $n\leq s$}
\end{cases}.
\]

\begin{proposicion}\label{P9'.3.3}
Para $\vec x,\vec y\in R_N$, $u\in R$ se tiene
\begin{gather}
V(\vec x\Witt +\vec y)=V\vec x\Witt +V\vec y,\label{Eq9'.3.8}\\
\vec x=(x_1,x_2,\ldots)=\sum_{j=0}^{\substack{r\\ \bullet}} V^j(\{x_{j+1}\})\Witt +
V^{r+1}(x_{r+2},x_{r+3},\ldots), \label{Eq9'.3.9}\\
\{u\}(x_1,x_2,\ldots,x_n,\ldots)=(ux_1,u^px_2,\ldots,
u^{p^{n-1}}x_{n-1},\ldots).\label{Eq9'.3.10}
\end{gather}
\end{proposicion}

\begin{proof} Puesto que $\vec x=\vec y$ si y s\'olo si $x^{(n)}=y^{(n)}$
para toda $n\in {\ma N}$, basta verificar que las $n$--componentes
fantasma coinciden.

\l
\item $(V(\vec x\Witt +\vec y))^{(n)}=p(x+y)^{(n-1)}=p(x^{(n-1)}+ y^{(n-1)})=
px^{(n-1)}+py^{(n-1)}= (V\vec x)^{(n)}+(V\vec y)^{(n)}$ de donde 
$V(\vec x \Witt + \vec y)=V\vec x\Witt +V\vec y$.

\item Se tiene
\begin{gather*}
\Big(\sum\limits_{j=0}^{\substack{r\\ \bullet}} V^j(\{x_{j+1}\})
+V^{r+1}(x_{r+2},x_{r+3},\ldots)\Big)^{(n)}\\
=\sum\limits_{j=0}^r (V^j(\{x_{j+1}\}))^{(n)}
+(V^{r+1}(x_{r+2},x_{r+3},\ldots))^{(n)}\\
=\sum\limits_{j=0}^r (V^j(\{x_{j+1}\}))^{(n)}
+p^{r+1}(x_{r+2},x_{r+3},\ldots)^{(n-(r+1))}=A.
\end{gather*}

Para $n=1,2,\ldots, r+1$, $0\leq j\leq r$, 
\[
(V^j(\{x_{j+1}\}))^{(n)}=(p^j\{x_{j+1}\})^{(n-j)}=
\begin{cases}
p^j x_{j+1}^{p^{n-j-1}}, & n\geq j+1\\
0, & n\leq j
\end{cases},
\]
y $p^{r+1}(x_{r+2},x_{r+3},\ldots)^{(n-(r+1))}=0$.

Por tanto, para $n=1,2,\ldots, r+1$, se tiene
$A=\sum_{j=0}^{n-1}p^j x_{j+1}^{
p^{n-j-1}}=\sum_{j=1}^n p^{j-1}x_j^{n-j}=x^{(n)}$.

Ahora bien, para $n\geq r+2$, $(V^j(\{x_{j+1}\}))^{(n)}=
p^jx_{j+1}^{p^{n-j-1}}$, $j=0,1,\ldots, r$ y 
\begin{align*}
p^{r+1}(x_{r+2},x_{r+3},\ldots)^{(n-(r+1))}&=p^{r+1}
(x_{r+2},x_{r+3},\ldots)^{(n-r-1)}\\
&= p^{r+1}\big(
\sum_{i=1}^{n-r-1}p^{i-1}x_{r+1+i}^{p^{n-r-1-i}}\big)\\
&=\sum_{i=1}^{n-r-1}p^{r+i}x_{r+1+i}^{p^{n-r-1-i}}\qquad i\leftrightarrow r+1+i\\
&=\sum_{i=r+2}^n p^{i-1}x_i^{p^{n-i}}.
\end{align*}

Por tanto
\begin{align*}
A&=\sum_{j=0}^rp^j x_{j+1}^{p^{n-j-1}}+\sum_{j=r+2}^n p^{j-1}x_j^{p^{n-j}}\\
&=\sum_{j=1}^{r+1}p^{j-1}x_j^{p^{n-j}}+\sum_{j=r+2}^n p^{j-1}x_j^{p^{n-j}}\\
&=\sum_{j=1}^n p^{j-1}x_j^{p^{n-j}}=x^{(n)}.
\end{align*}

Se sigue (\ref{Eq9'.3.9}).

\item Se tiene que $\big(\{u\}(x_1,x_2,\ldots,x_n\ldots)\big)^{(n)}=\{u\}^{(n)}
(x_1,x_2,\ldots)^{(n)}=u^{p^{n-1}}x^{(n)}$.
Por otro lado 
\begin{align*}
(ux_1,u^px_2,\ldots,u^{p^{n-1}}x_{n-1},\ldots)^{(n)}&=
\sum_{j=1}^np^{j-1}(u^{p^{j-1}} x_j)^{p^{n-j}}\\
&= u^{p^{n-1}}\sum_{j=1}^n
p^{j-1}x_j^{p^{n-j}}=u^{p^{n-1}}x^{(n)}.
\end{align*}
Se sigue (\ref{Eq9'.3.10}). $\fin$
\end{list}
\end{proof}

\begin{notacion}\label{N9'.3.4}
Para $m\in {\ma N}\cup\{0\}$, se denota
\begin{equation}\label{Eq9'.3.11}
\vec 0:=(0,0,\ldots,0,\ldots),\quad \vec 1:=(1,0,\ldots,0,\ldots),
\quad \vec m=m\vec 1:=\underbrace{\vec 1\Witt +\vec 1\Witt +\cdots\Witt +
\vec 1}_{m\text{\ veces}}.
\end{equation}
\end{notacion}

Sea $p$ un n\'umero primo, $\vec x=(x_1,\ldots,x_n,\ldots)$. Se define
\begin{equation}\label{Eq9'.3.12}
F(\vec x)=\vec x^p:=(x_1^p,\ldots, x_n^p,\ldots).
\end{equation}

Observemos que $\vec x^p$ {\underline{no}} es la $p$--potencia de
la multiplicaci\'on de Witt, es decir $\vec x^p\neq \underbrace{\vec x \Witt \times
\vec x \Witt \times \cdots \Witt \times \vec x}_{p}$.

\begin{definicion}\label{D9'.3.5}
Al homomorfismo $F$ se la llama el {\em automorfismo de
Frobenius\index{automorfismo de Frobenius para vectores de
Witt}\index{Frobenius!automorfismo de $\sim$ para vectores de Witt}}
en $R_N$.
\end{definicion}

Observemos que
\begin{gather}
x^{(n)}=\sum_{j=1}^n p^{j-1}x_j^{p^{n-j}}=\sum_{j=1}^{n-1}
p^{j-1}\big(x_j^p\big)^{p^{n-1-j}}+p^{n-1}x_n=\big(x^p\big)^{(n-1)}+p^{n-1}x_n.\nonumber\\
\intertext{Esto es}
x^{(n)}=\big(x^p\big)^{(n-1)}+p^{n-1}x_n,\quad n=1,2,\ldots, \quad
\big(x^p\big)^{(0)}=0.\label{Eq9'.3.13}
\end{gather}

Observemos que $\big(x^p\big)^{(n-1)}$ denota la $(n-1)$--componente fantasma
de $x^p$ y no la $p(n-1)$ potencia $x$. De (\ref{Eq9'.3.4}), obtenemos que
si $I={\ma Z}[x_i,y_j,z_k]$, entonces
\[
p^{n-1}(x+y)_n\equiv (x+y)^{(n)}\bmod I=(x^{(n)}+y^{(n)})\equiv p^{n-1}x_n
+p^{n-1}y_n \bmod I.
\]
Por tanto existe $f\in I$, tal que
\begin{equation}\label{Eq9'.3.14}
(x+y)_n=x_n+y_n+f(x_1,y_1,\ldots,x_{n-1},y_{n-1}).
\end{equation}

\subsection{Aritm\'etica de los vectores de Witt}\label{S9'.5}

Sea ${\eu F}$ un dominio entero de caracter{\'\i}stica $0$, de tal forma
que ${\ma Z}\subseteqq {\eu F}$. Sea $p\in{\ma Z}$ un n\'umero primo.

\begin{lema}\label{L9'.5.1}
Sean $\vec x,\vec y$ dos vectores cuyas componentes regulares est\'an en ${\eu F}$.
Entonces para $r>0$, se tienen que las congruencias
\[
x_n\equiv y_n\bmod p^r {\eu F}\quad \text{y}\quad x^{(n)}\equiv
y^{(n)}\bmod p^{r+n-1} {\eu F}
\]
son equivalentes.
\end{lema}

\begin{proof} Procedemos por inducci\'on en $n$. Si para $n-1$ se tiene la equivalencia,
entonces si $x_n\equiv y_n \bmod p^r {\eu F}$, entonces $(x^p)_n=
x^p_n\equiv y^p_n=(y^p)_n\bmod p^{r+1} {\eu F}$ por lo que
$(x^p)^{(n-1)}\equiv (y^p)^{(n-1)}\bmod p^{r+n-1} {\eu F}$.

Ahora, por (\ref{Eq9'.3.13}), $x^{(n)}=(x^p)^{(n-1)}+p^{n-1}x_n$, entonces
\[
(x^{(n)}-y^{(n)})-(p^{n-1}x_n-p^{n-1}y_n)=(x^p)^{(n-1)}-(y^p)^{(n-1)}
\bmod p^{r+n-1} {\eu F}.
\]

Por lo tanto $x^{(n)}\equiv y^{(n)}\bmod p^{r+n-1}{\eu F} \iff p^{n-1}x_n
-p^{n-1}y_n\equiv 0 \bmod p^{r+n-1} {\eu F} \iff x_n\equiv y_n \bmod p^r {\eu F}$.
$\fin$
\end{proof}

Como consecuencia inmediata del Lema \ref{L9'.5.1} tenemos:

\begin{teorema}\label{T9'.5.2} Se tiene que $(\vec x \Witt {\masmenos\limits_{\times}} \vec y)_n
\in {\ma Z}[x_1,y_1,\ldots,x_{n-1},y_{n-1},x_n,y_n]$, $n=1,2,\ldots$.
\end{teorema}

\begin{proof} Por definici\'on $x^{(n)}, y^{(n)}\in {\eu F}={\ma Z}[x_1,y_1,\ldots,x_n,y_n]$.
Ahora bien, puesto que por (\ref{Eq9'.3.13}) se tiene
 $x^{(n)}\equiv (x^p)^{(n-1)}\bmod
p^{n-1}{\eu F}$ y $y^{(n)}\equiv (y^p)^{(n-1)}\bmod p^{n-1}{\eu F}$, se sigue que
$(\vec x \Witt {\masmenos\limits_{\times}} \vec y)^{(n)}=x^{(n)} {\masmenos\limits_{\times}} y^{(n)}\equiv
(x^p)^{(n-1)} \masmenos\limits_{\times} (y^p)^{(n-1)}\bmod p^{n-1}{\eu F}$.

Por inducci\'on, si para $j<n$ el teorema ya est\'a demostrado, entonces
$(\vec x \Witt {\masmenos\limits_{\times}}\vec y)_j^p\equiv (x^p\masmenos\limits_{\times} y^p)_j
\bmod p {\eu F}$. Por el Lema \ref{L9'.5.1} para $n-1$, 
\[
\big((\vec x\Witt {\masmenos_{\times}} \vec y)^p\big)^{(n-1)}\equiv (x^p \masmenos_{\times}
y^p)^{(n-1)}\bmod p^{n-1} {\eu F}.
\]
As{\'\i} $p^{n-1}(\vec x \Witt {\masmenos\limits_{\times}}\vec y)_n=
(\vec x \Witt {\masmenos\limits_{\times}}\vec y)^{(n)}-\big((\vec x\Witt {\masmenos\limits_{\times}}
\vec y)^p\big)^{(n-1)} \equiv 0\bmod p^{n-1} {\eu F}$, de donde
se sigue que $\vec x\Witt {\masmenos\limits_{\times}}\vec y \in {\eu F}$.
$\fin$
\end{proof}

Cuando estudiemos la acci\'on de Galois, necesitaremos el siguiente resultado.

\begin{teorema}\label{T9'.5.3} Se tiene que, por componentes, 
\[
\vec p \vec x\equiv
V\vec x^p \bmod p{\ma Z}[x_1,x_2,\ldots],
\]
 es decir $(\vec p \vec x)_n\equiv
(V\vec x^p)_n \bmod p{\ma Z}[x_1,x_2,\ldots]$.
\end{teorema}

\begin{proof} Por (\ref{Eq9'.3.13}), (\ref{Eq9'.3.(c)}) y (\ref{Eq9'.2.3}) se tiene que las componentes
fantasma satisfacen
\begin{align*}
(\vec p\vec x)^{(n)}&=px^{(n)}\equiv p(\vec x^p)^{(n-1)}+p^nx_n\equiv
p(\vec x^p)^{(n-1)}\\
&\equiv (V \vec x^p)^{(n)} \bmod p^n{\ma Z}[x_1,\ldots,x_n,\ldots ].\\
\intertext{Por el Lema \ref{L9'.5.1}, con $r=1$, se obtiene}
(\vec p\vec x)_n&=
(V\vec x^p)_n\bmod p{\ma Z}[x_1,\ldots,x_n,\ldots ]. \tag*{$\fin$}
\end{align*}
\end{proof}

\subsection{Vectores de Witt en caracter{\'\i}stica $p$}\label{S9'.6}

Hasta ahora hemos considerado las operaciones de Witt en caracter{\'\i}stica
$0$ pues al pasar de las componentes fantasma a las componentes de Witt,
se est\'a dividiendo entre una potencia de $p$. Esto es, con la notaci\'on 
de (\ref{Eq9'.3.4})
\[
\vec a\Witt {\masmenos_{\times}}\vec b= \big(\vec a^{\varphi} {\masmenos_{\times}}
\vec b^{\varphi}\big)^{\varphi^{-1}} =\big(\vec a^{\varphi} {\masmenos_{\times}}
\vec b^{\varphi}\big)^{\psi}
\]
lo cual hace de $R_N$ un anillo pues todas las reglas se cumplen en $R^N$
y son transformadas a $R_N$ bajo $\varphi^{-1}$. Por ejemplo
\begin{align*}
(\vec a\Witt + \vec b)\Witt + \vec c&=(\vec a^{\varphi} +\vec b^{\varphi})^{
\varphi^{-1}}\Witt +\vec c=\Big(\big((\vec a^{\varphi}  +\vec b^{\varphi})^{
\varphi^{-1}}\big)^{\varphi} +\vec c^{\varphi}\Big)^{\varphi^{-1}}\\
&=\big((\vec a^{\varphi}  +\vec b^{\varphi}) + \vec c^{\varphi}\big)^{
\varphi^{-1}} =\big(\vec a^{\varphi} +(\vec b^{\varphi} + \vec c^{\varphi})\big)^{
\varphi^{-1}}\\
&=\Big(\vec a^{\varphi}  +\big((\vec b^{\varphi} + \vec c^{\varphi})^{
\varphi^{-1}}\big)^{\varphi}\Big)^{\varphi^{-1}}\\
&=\Big(\vec a^{\varphi} +\big(\vec b\Witt + \vec c\big)^{\varphi}\Big)^{
\varphi^{-1}}=\vec a\Witt +(\vec b\Witt + \vec c).
\end{align*}

Por el Teorema \ref{T9'.5.3} las operaciones de Witt pueden hacerse m\'odulo
$p$ y de esta forma obtenemos

\begin{teorema}\label{T9'.6.1} Sea $k$ un campo de caracter{\'\i}stica $p$
y sea $W_N(k)$ el {\em anillo de Witt\index{anillo de Witt}\index{Witt!anillo de $\sim$}}
\[
W_N(k):=\{(x_1,\ldots,x_n,\ldots )\mid x_i\in k\}, \quad N\in {\ma N}\cup
\{\infty\}.
\]
Entonces $W_N(k)$ es un anillo conmutativo con unidad y se 
tiene para $\vec x,\vec y\in W_N(k)$
\begin{gather}
(\vec x \Witt {\masmenos_{\times}}\vec y)^p= \vec x^p
 \Witt {\masmenos_{\times}}
\vec y^p. \label{Eq9'.6.12}\\
\vec p \Witt \times \vec x=\vec p\vec x=V\vec x^p=
(V\vec x)^p.\label{Eq9'.6.13}\\
(V^i\vec x) \Witt \times (V^j \vec y)=
V^{i+j}(\vec x^{p^j}\Witt \times \vec y^{p^i}).\label{Eq9'.6.14}
\end{gather}
\end{teorema}

\begin{proof} Todas las propiedades de anillo se cumplen formalmente por lo que
$W_N(k)$ es un anillo. Ahora bien, se tiene que (\ref{Eq9'.6.12}) se
sigue del Lema \ref{L9'.5.1} y de (\ref{Eq9'.3.13}):
\begin{multline*}
\big((\vec x\Witt + \vec y)^p\Witt -
\vec x^p\Witt - \vec y^p\big)^{(n)} =
\big((\vec x\Witt+\vec y)^p\big)^{(n)} -
(\vec x^p)^{(n)} - (\vec y^p)^{(n)}\\
=(\vec x \Witt + \vec y)^{(n+1)}
- p^n(\vec x\Witt +\vec y)_n
- x^{(n+1)}+ p^nx_n- y^{(n+1)}+ p^n y_n.
\end{multline*}

Puesto que $(\vec x\Witt + \vec y)^{(n+1)}=x^{(n+1)}+y^{(n+1)}$, se sigue
que $\big((\vec x\Witt + \vec y)^p\Witt -\vec x^p\Witt -\vec y^p\big)^{(n)}
\equiv 0\bmod p^n= 0\bmod p^{1+(n-1)}$. Por el Lema \ref{L9'.5.1} se
sigue $(\vec x\Witt +\vec y)^p=\vec x^p\Witt + \vec y^p$. Similarmente
$(\vec x\Witt - \vec y)^p=\vec x^p\Witt - \vec y^p$ y
$(\vec x\Witt {\times}\vec y)^p=\vec x^p\Witt {\times} \vec y^p$.

Se tiene que (\ref{Eq9'.6.13}) es el Teorema \ref{T9'.5.3}. Adem\'as tenemos
por (\ref{Eq9'.3.(c)}) y (\ref{Eq9'.3.13}) que
\begin{gather*}
\big((V\vec x)^p\big)^{(n)}= (V\vec x)^{(n+1)}-p^n x_{n+1}=px^{(n)}\\
\intertext{y}
(V\vec x^p)^{(n)}=p(\vec x^p)^{(n-1)}=p(x^{(n)}-p^{n-1}x_n)\\
\intertext{de donde se obtiene}
\big((V\vec x)^p\big)^{(n)}\equiv (V\vec x^p)^{(n)}\bmod p^n.\\
\intertext{Por el Lema \ref{L9'.5.1} obtenemos que}
(V\vec x)_n^p\equiv (V\vec x^p)_n\bmod p.\\
\intertext{Por lo tanto}
F\circ V=V\circ F=\vec p.
\end{gather*}

Para probar (\ref{Eq9'.6.14}) notemos primero que de (\ref{Eq9'.3.13})
obtenemos
\begin{align*}
(\vec x^p)^{(n)}&=x^{(n+1)}-p^nx_n\equiv x^{(n+1)}\bmod p^n,\\
(\vec x^{p^2})^{(n)}&=(\vec x^p)^{(n+1)}-p^n(\vec x^p)_n\\
&=x^{(n+2)}-p^{n+1}x_{n+1}-p^n(\vec x^p)_n\equiv x^{(n+2)}\bmod p^n.\\
\intertext{En general obtenemos}
(\vec x^{p^j})^{(n)}&\equiv x^{(n+j)}\bmod p^n.
\end{align*}

Ahora para $i=0,j=0$ se tiene 
\[
(V^i\vec x)\Witt \times (V^j \vec y)= (V^0\vec x)\Witt \times (V^0 \vec y)=
\vec x\Witt \times \vec y= V^{0+0}(\vec x^{p^0} \Witt \times \vec y^{p^0}),
\]
por lo que se cumple (\ref{Eq9'.6.14}) para $i=j=0$.

Para $i=0, j\geq 1$ se tiene
\begin{gather*}
(V^i\vec x)\Witt \times (V^j \vec y) = V^0(\vec x)\Witt \times V^j(\vec y)=
\vec x\Witt \times V^j\vec y\\
\intertext{y}
V^{i+j}(\vec x^{p^j}\vec y^{p^i}) = V^j(\vec x^{p^j} \Witt \times \vec y).
\end{gather*}
Se sigue que
\begin{align*}
(\vec x\Witt \times V^j\vec y)^{(n)}&=x^{(n)}(V^j\vec y)^{(n)}=p^jx^{(n)}y^{(n-j)},\\
\big(V^j(\vec x^{p^j}\Witt \times \vec y)\big)^{(n)}&=p^j\big((\vec x^{p^j})
\Witt \times \vec y\big)^{(n-j)}=p^j(\vec x^{p^j})^{(n-j)} y^{(n-j)}\\
&\equiv
p^jx^{(n)}y^{(n-j)}\bmod p^n.
\end{align*}
Por lo tanto $(\vec x V^j \vec y)_n\equiv (V^j(\vec x^{p^j}\Witt \times \vec y))_n
\bmod p$ lo cual implica que $(V^i \vec x)\Witt \times (V^j\vec y)=
V^{i+j}(\vec x^{p^j}\Witt \times \vec y^{p^i})$ para $i=0, j\geq 1$.

Para $i\geq 1,j=0$, $(V^i\vec x)\Witt \times (V^j\vec y)=(V^j \vec y)
\Witt \times (V^i \vec x)=V^{j+i}(\vec y^{p^i}\Witt \times \vec x^{p^j})=
V^{i+j}(\vec x^{p^j}\Witt \times \vec y^{p^i})$.

Para el caso $i\geq 1,j\geq 1$, de la relaci\'on $V\vec x^p=\vec p\vec x$, se tiene
\begin{gather*}
(V^j\vec x^p)=V(V^{j-1}\vec x^p)=\vec p(V^{j-1}\vec x^p)=\vec p^2
(V^{j-2}\vec x^p)=\cdots=\vec p^{j-1} V\vec x^p= \vec p^j\vec x,\\
\intertext{de donde}
\begin{align*}
(V^i\vec x)\Witt \times (V^j\vec y)&= (V^{i-1}(V\vec x))\Witt \times
(V^j\vec y)=V^{i+j-1}((V\vec x)^{p^j}\Witt \times \vec y^{p^{i-1}})\\
&=V^{i+j-1}((V(\vec x^{p^j}))\Witt \times \vec y^{p^{i-1}}) = V^{i+j-1}((
V\vec x^{p^j})\Witt \times (V^0\vec y^{p^{i-1}}))\\
&=V^{i+j-1}(V\vec x^{p^j}\Witt \times (\vec y^{p^{i-1}})^p)=
V^{i+j}(\vec x^{p^j}\Witt \times \vec y^{p^i}). \tag*{$\fin$}
\end{align*}
\end{gather*}
\end{proof}

Como veremos posteriormente, el siguiente resultado nos permite obtener
el an\'alogo al Teorema \ref{T9'.2.1} para $n\in {\ma N}$ (el Teorema
\ref{T9'.2.1} es el caso $n=1$).

\begin{teorema}\label{T9'.6.2} En el anillo $W_N(k)$, donde $k$ 
es un campo de caracter{\'\i}stica $p$, se tiene que el vector
$\vec a=(a_1,\ldots,a_n,\ldots)\in W_N(k)$ es invertible si y s\'olo si
$a_1 \neq 0$.
\end{teorema}

\begin{proof} Se tiene que $\{a_1^{-1}\}\Witt \times \vec a=(1,y_1,\ldots)$. Por tanto
$\vec 1\Witt - \vec a\Witt \times \{a_1^{-1}\}=V\vec y$ para alg\'un $\vec y\in W_N(k)$. Sea 
$(V\vec y)^i$ la $i$--\'esima potencia de $V\vec y$ en $W_N(k)$, es decir,
$(V\vec y)^i= \underbrace{V\vec y\Witt \times V\vec y\Witt \times \cdots \Witt
\times V\vec y}_{i}$. Entonces
\begin{align*}
\vec a\Witt \times \{a_1^{-1}\}\Witt \times \sum_{j=0}^{\substack{\infty\\ \bullet}}
(V\vec y)^j &= (1\Witt - V\vec y)\Witt {\times} \sum_{j=0}^{\substack{\infty\\ \bullet}}
(V\vec y)^j\\
&= \sum_{j=0}^{\substack{\infty\\ \bullet}}
(V\vec y)^j\Witt -\sum_{j=1}^{\substack{\infty\\ \bullet}}
(V\vec y)^j=(V\vec y)^0=\vec 1
\end{align*}
y por tanto $\vec a$ es invertible y de hecho se tiene
$\vec a^{-1}=\{a_1^{-1}\}\Witt \times \sum\limits_{j=0}^{\substack{\infty\\ \bullet}}
(V\vec y)^j$.

El rec{\'\i}proco es inmediato pues existe $\vec b\in W_N(k)$ tal que
$\vec a\Witt \times \vec b=(a_1 b_1,\ldots)=\vec 1=(1,0,\ldots)$, esto
es, $a_1b_1=1$ y en particular $a_1\neq 0$. $\fin$
\end{proof}

De la demostraci\'on del Teorema \ref{T9'.6.1}, obtenemos
\begin{equation}\label{Eq9'.6.15}
\vec p^j\Witt \times \vec 1=V^j(\vec 1^p)=(\underbrace{0,\ldots,0}_j,1,0,\ldots),
\quad j\in {\ma N}.
\end{equation}

\begin{ejemplo}\label{E9'.6.3}
Consideremos $N=n\in{\ma N}$ y $W_n({\ma F}_p)$. Notemos que
$\vec 1\in W_n({\ma F}_p)$ y se tiene de (\ref{Eq9'.6.15}) que
$\vec p^n=\vec p^n \Witt \times \vec 1 = \vec 0$ pero $\vec p^{n-1}
=\vec p^{n-1} \Witt \times \vec 1\neq \vec 0$. Adem\'as se tiene que $|
W_n({\ma F}_p)|=p^n$, lo cual implica que, como grupo con la adici\'on
de Witt, $W_n({\ma F}_p)$ es c{\'\i}clico de orden $p^n$.
M\'as a\'un, $\varphi\colon {\ma Z}\to W_n({\ma F}_p)$, $1\mapsto
\vec 1$, es un epimorfismo de anillos con $\ker \varphi (p^n)$,
se sigue que $W_n({\ma F}_p)\cong {\ma Z}/p^n{\ma Z}$ como anillos.

En particular tenemos que $W_n({\ma F}_p)$ es de caracter{\'\i}stica
$p^n$. Por otro lado, las unidades de $W_n({\ma F}_p)$ son precisamente
$\{\vec i\mid 1\leq i\leq p^n-1, \mcd (i,p)=1\}$.
\end{ejemplo}

\begin{ejemplo}\label{E9'.6.4}
En el caso $N=\infty$ se tiene que $W_N({\ma F}_p)$ es un
anillo de caracter{\'\i}stica $0$ pues para toda $n$, $\vec p^n=
\vec p^n\Witt \times \vec 1=(0,\ldots,0,\underbracket[0pt]{1}_{
\substack{\uparrow\\ n+1}},0,\cdots)
\neq 0$ y si $\mcd (i,p)=1$, $i\in{\ma N}$, $\vec i$ es unidad pues
$i=\sum_{j=0}^{m}a_jp^j$, $a_j\in\{0,1,\ldots, p-1\}$ y $a_1\neq 0$.

De hecho se tiene que $W_{\infty}({\ma F}_p)\cong {\ma Z}_p$,
${\ma Z}_p$ el anillo de los enteros $p$--\'adicos.
\end{ejemplo}

\begin{proposicion}\label{P9'.6.5} Sea $k$ un campo de caracter{\'\i}stica
$p$. Entonces $W_N(k)$ es de caracter{\'\i}stica $p^n$ si $N=n\in{\ma N}$
y $0$ si $N=\infty$.
\end{proposicion}

\begin{proof} Igual que la de los Ejemplos \ref{E9'.6.3} y \ref{E9'.6.4}.
$\fin$
\end{proof}

\subsection{Extensiones c{\'\i}clicas de grado $p^n$ en caracter{\'\i}stica
$p$}\label{S9'.7}

Usando los vectores de Witt, se tiene una teor{\'\i}a para
$p$--extensiones c{\'\i}clicas finitas paralela a la
Teor{\'\i}a de Artin--Schreier\index{Artin--Schreier!teor{\'\i}a
de $\sim$}\index{teor{\'\i}a de Artin--Schreier}
 para extensiones c{\'\i}clicas de grado
$p$ en caracter{\'\i}stica $p$ y siendo esta \'ultima
a su vez una teor{\'\i}a
paralela a la Teor{\'\i}a de Kummer\index{Kummer!teor{\'\i}a de $\sim$}
\index{teor{\'\i}a de Kummer}. Esta teor{\'\i}a recibe con
frecuencia el nombre de {\em Teor{\'\i}a Aditiva de
Kummer\index{Kummer!teor{\'\i}a aditiva de $\sim$}\index{teor{\'\i}a
aditiva de Kummer}}.

En el caso de una extensi\'on c{\'\i}clica $L/K$ en caracter{\'\i}stica
$p$ de grado $p^n$ con $n=1$, Artin y Schreier probaron que toda
tal extensi\'on est\'a dada por un ecuaci\'on del tipo $\wp y=x$, donde
$\wp y:=y^p-y$ y $x\in K, x\notin \wp(K)=\{a^p-a\mid a\in K\}$
(Teorema \ref{T9'.2.1}). La demostraci\'on se basa en que si un
elemento $z\in L$ satisface $\Tr_{L/K} z=0$, entonces existe
$w\in L$ tal que $(\sigma-1)w=z$ donde $\langle \sigma\rangle
=\Gal (L/K)$, aunque la demostraci\'on original de Artin--Schreier,
que es la que presentamos nosotros, no lo hizo de esta forma.
Esto mismo se cumple pr\'acticamente palabra por
palabra  para extensiones c{\'\i}clicas de grado $p^n$ usando
el lenguaje de los vectores de Witt.

Sea $K$ un campo arbitrario de caracter{\'\i}stica $p$ y consideremos
$W_n(K)=\{(x_1,\ldots,x_n)\mid x_i\in K\}$ el anillo de los vectores
de Witt de longitud $n$ con coeficientes en $K$. Sea $L/K$
una extensi\'on finita de Galois con grupo de Galois $G=\Gal(
L/K)$.

\begin{definicion}\label{D9'.7.1} Si $\vec y\in W_n(L)$,
$\vec y=(y_1,\ldots, y_n)$, se define para $\sigma \in G$,
\[
\sigma\vec y:=(\sigma y_1,\ldots,\sigma y_n)= \vec y^{\sigma}
\]
y la {\em traza\index{traza de vectores de Witt}\index{vectores de
Witt!traza}} $\Tr_{L/K}\colon W_n(L)\to W_n(K)$ se define
por
\[
\Tr_{L/K}\vec y=\sum_{\sigma\in G}^{\bullet}\sigma \vec y=
(\Tr_{L/K} y_1,?,\ldots,?)\in W_n(K).
\]
\end{definicion}

Si $y_1\in L$ es tal que $\Tr_{L/K}y_1\neq 0$, $\Tr_{L/K} \vec y$
es invertible (Teorema \ref{T9'.6.2}). Adem\'as tenemos que
$\sigma(\vec y\Witt +\vec z)=\sigma \vec y\Witt +\sigma \vec z$
y $\sigma(\vec y\Witt \times \vec z)=\sigma \vec y\Witt \times
\sigma \vec z$ pues si $\vec a=(a_1,\ldots, a_n)$ entonces
$\sigma \vec a=(\sigma a_1,\ldots,\sigma a_n)$ con 
$(\sigma \vec a)^{(t)}=\sum_{i=1}^t p^{i-1}(\sigma a_i)^{p^{t-i}}=
\sigma\big(\sum_{i=1}^t p^{i-1}a_i^{p^{t-i}}\big)=\sigma a^{(t)}$.

El siguiente resultado nos prueba que el primer grupo de 
cohomolog{\'\i}a $H^1(W_n(L), G)$ es igual a $\{0\}$. M\'as precisamente

\begin{teorema}\label{T9'.7.2} Sea $\varphi\colon G\to W_n(L)$
con $\varphi(\sigma)=\vec a_{\sigma}$. Si se tiene $\vec a_{\sigma}
\Witt + \sigma \vec a_{\tau}=\vec a_{\sigma\tau}$ para cualesquiera
$\sigma,\tau\in G$, entonces existe $\vec b\in W_n(L)$ tal que
$\vec a_{\sigma}=(1\Witt - \sigma)\vec b$ para toda $\sigma\in G$.
\end{teorema}

\begin{proof} Sea $\vec c=(c_1,\ldots,c_n)\in W_n(L)$ tal que $\Tr_{
L/K}c_1\neq 0$. Tal $\vec c$ existe pues $L/K$ es separable.
Ahora por el Teorema \ref{T9'.6.2}, se tiene que
$\Tr_{L/K}\vec c\in W_n(K)$ es invertible. Sea
\[
\vec b:=(\Tr_{L/K} \vec c)^{-1}\Witt \times \sum_{\tau\in G}^{\bullet} \vec a_{\tau} \Witt \times\tau \vec c.
\]
Entonces, para $\sigma\in G$, se tiene
\begin{align*}
(1\Witt - \sigma)\vec b&=\vec b\Witt - (\Tr_{L/K}\vec c)^{-1}\Witt \times
\sum_{\tau\in G}^{\bullet} \sigma \vec a_{\tau}\Witt \times (\sigma\tau) \vec c\\
&=\vec b\Witt - (\Tr_{L/K}\vec c)^{-1}\Witt \times\Big(
\sum_{\tau\in G}^{\bullet} (\vec a_{\sigma\tau}\Witt - \vec a_{\sigma})\Witt \times (\sigma\tau) \vec c\Big)\\
&= \vec b\Witt - (\Tr_{L/K}\vec c)^{-1}\Witt\times \Big(\sum_{\tau\in G}^{\bullet} 
\vec a_{\sigma\tau} \Witt \times (\sigma\tau) \vec c\Big)\\
&\hspace{2cm}\Witt +(\Tr_{L/K} \vec c)^{-1}
\Witt\times \vec a_{\sigma}\Witt \times\Big(\sum_{\tau\in G}^{\bullet} (\sigma\tau) \vec c\Big)\\
&= \vec b\Witt -(\Tr_{L/K} \vec c)^{-1}\Witt \times (\Tr_{L/K} \vec c)\Witt \times \vec b\\
&\hspace{2cm} \Witt +
(\Tr_{L/K} \vec c)^{-1}\Witt \times \vec a_{\sigma} \Witt \times (\Tr_{L/K} \vec c)\\
&=\vec b\Witt - \vec b\Witt + \vec a_{\sigma}= \vec a_{\sigma}
\end{align*}
para toda $\sigma\in G$. $\fin$
\end{proof}

\begin{definicion}\label{D9'.7.3} Para $\vec y\in W_n(L)$, se define
\[
\wp \vec y:=\vec y^p\Witt - \vec y= (y_1^p,\ldots,y_n^p)\Witt -
(y_1,\ldots,y_n).
\]
\end{definicion}

Se tiene que $\wp(\vec y\Witt + \vec z)=\wp \vec y\Witt + \wp \vec z$
para cualesquiera $\vec y, \vec z\in W_n(L)$ (Teorema \ref{T9'.6.1})).

\begin{proposicion}\label{P9'.7.4} Se tiene que $\wp \vec x=\vec 0\iff
\vec x\in W_n({\ma F}_p)$.
\end{proposicion}

\begin{proof} Se tienen las equivalencias
\begin{gather*}
\wp \vec x=\vec 0\iff \vec x^p=\vec x\iff (x_1^p,\ldots, x_n^p)=(x_1,\ldots,x_n)\\
\iff x_i^p=x_i, 1\leq i\leq n\iff x_i\in{\ma F}_p, 1\leq i\leq n\iff \vec x\in
W_n({\ma F}_p). \tag*{$\fin$}
\end{gather*}
\end{proof}

\begin{definicion}\label{D9'.7.5-1}
Sea $K$ un campo de caracter\'istica $p>0$. Sea $W_n(K)$ el anillo
de los vectores de Witt de longitud $n$. Se define la relaci\'on $\vec
\alpha\sim\vec \beta$ si y solamente si existe $\vec c\in W_n(K)$ tal que
$\vec\alpha\Witt -\vec\beta=\wp(\vec c)$ o, equivalentemente, $\vec 
\alpha=\vec \beta\Witt+\wp(\vec c)$.
\end{definicion}

Se tiene que $\sim$ es una relaci\'on de equivalencia y si $\vec \alpha
\sim\vec \beta$ y $\vec \alpha_1\sim\vec \beta_1$ entonces $\vec
\alpha \Witt+ \vec \alpha_1\sim \vec \beta\Witt +\vec \beta_1$. 
Finalmente, si $\vec \alpha\sim\vec\beta$ entonces $\Witt -\vec \alpha
\sim \Witt -\vec \beta$ pues $\vec 0=\wp(\vec c\Witt -\vec c)=
\wp(\vec c)\Witt+\wp(\Witt -\vec c)$, es decir, $\wp(\Witt - \vec c)=
\Witt -\wp(\vec c)$.

Notemos que $\vec p^n=\vec 0$ (Ejemplo \ref{E9'.6.4}), por lo que
$\vec{p^n-1}=\Witt - \vec 1$.

\begin{definicion}\label{D9'.7.5} Un vector $\vec x\in W_n(K)$ se llama
{\em descompuesto\index{vector de Witt descompuesto}\index{Witt!vector
descompuesto} o descomponible\index{vector de Witt
descomponible}\index{Witt!vector descomponible}} si
$\vec x\sim \vec 0$, es decir, 
si existe $\vec y\in W_n(K)$ tal que $\wp \vec y=\vec x$.
Los vectores descomponibles son los elementos de la clase de
equivalencia del vector $\vec 0$.

De esta forma $\vec x\sim \vec 0$
significa que $\vec x\in \wp(W_n(K))$.
\end{definicion}

\begin{proposicion}\label{P.9'.7.6} Se tiene que $(0,x_2,\ldots, x_n)
\in W_n(K)$ es descompuesto si y s\'olo si $(x_2,\ldots, x_n)\in
W_{n-1}(K)$ es descompuesto.
\end{proposicion}

\begin{proof} Si $(0,x_2,\ldots,x_n)=\wp \vec y=(\wp y_1,\ldots)$ se tiene que
$\wp y_1=0$ y por tanto $y_1\in{\ma F}_p$. Por tanto
$(y_1,0,\ldots,0)\in W_n({\ma F}_p)$ por lo que $\wp(y_1,0,\ldots,0)
=\vec 0$.

Por otro lado, tenemos de (\ref{Eq9'.3.9})
\[
(y_1,y_2,\ldots,y_n)=\{y_1\}\Witt + V(y_2,\ldots,y_n)=
(y_1,0,\ldots, 0)\Witt + (0,y_2,\ldots, y_n).
\]
Por tanto $\wp \vec y = \wp (\{y_1\})\Witt + \wp((0,y_2,\ldots,y_n))=
\vec 0\Witt + \wp ((0,y_2,\ldots,y_n))=(0,x_2,\ldots,x_n)$. Se sigue que
\begin{gather*}
\wp ((y_2,\ldots,y_n))=(x_2,\ldots,x_n)\iff \wp((0,y_2,\ldots,y_n)=
(0,x_2,\ldots,x_n). \tag*{$\fin$}
\end{gather*}
\end{proof}

\begin{observacion}\label{O9'.7.7-1}
Si $\alpha\in K$, entonces se tiene que
$(0,\ldots,0,\wp(\alpha))\sim\vec 0$ pues
$(\underbrace{0,\ldots,0}_{n-1},\wp(\alpha))\sim \vec 0\iff
(\underbrace{0,\ldots,0}_{n-2},\wp(\alpha))\sim \vec 0\iff
\cdots\iff \wp(\alpha)\sim \vec 0=0$.

Sin embargo para $\alpha\in K$, $(\wp(\alpha),0,\ldots,0)$ no
necesariamente es descomponible.
\end{observacion}

\begin{ejemplo}\label{E9'.7.7-2}
Consideremos $n=2$ y sea 
\begin{align*}
(z_1,z_2)=\vec z=\wp(\vec \beta)&=
\vec\beta^p \Witt -\vec\beta=(\beta_1^p,\beta_2^p)\Witt -(\beta_1,
\beta_2)\\
&=(\beta_1^p,\beta_2^p\mid \beta_1^p,\beta_1^{p^2}+p\beta_2^p)
\Witt - (\beta_1,\beta_2\mid \beta_1, \beta_1^p+p\beta_2)\\
&=(\beta_1^p-\beta_1,z_2\mid \beta_1^p-\beta_1,\beta_1^{p^2}+p
\beta_2^p-\beta_1^p-p\beta_2)\\
&=(\beta_1^p-\beta_1,z_2\mid \beta_1^p-\beta_1,\beta_1^{p^2}
-\beta_1^p+p(\beta_2^p-\beta_2))\\
&=(\beta_1^p-\beta_1,z_2\mid \beta_1^p-\beta_1, z^{(2)}).
\end{align*}

Por tanto $z^{(2)}=z_1^p+pz_2=(\beta_1^p-\beta_1)^p+pz_2=
\beta_1^{p^2}-\beta_1^p+p(\beta_2^p-\beta_2)$. Se sigue que
$z_2=\frac{(\beta_1^{p^2}-\beta_1^p)-(\beta_1^p-\beta_1)^p}{p}
+\beta_2^p-\beta_2$.

En particular, si queremos hallar $\vec\beta$ tal que 
$\wp(\vec \beta)=(\wp(\alpha),0)$, entonces $\wp(\beta_1)=
\wp(\alpha)$ y $z_2=0$, esto es, 
\[
\beta_2^p-\beta_2=\frac{(\beta_1^p-\beta_1)^p-
(\beta_1^{p^2}-\beta_1^p)}{p}=\sum_{i=1}^{p-1}\frac{1}{p}
\binom pi (-1)^i\beta_1^{p(p-i)+i}.
\]

Entonces $\beta_1=\alpha$ (o, en general, $\beta_1=\alpha+
\chi$ con $\chi\in\F$) y $\wp(\beta_2)=\sum_{i=1}^{p-1}\frac{1}{p}
\binom pi (-1)^i\beta_1^{p(p-i)+i}$.

Por ejemplo, sean $K=\F(T)$, $\alpha=\frac{1}{T}$, $\beta_1=
\frac{1}{T}$ y $p=3$. En este caso, tenemos
\[
\beta_2^3-\beta_2=-\frac{1}{3}\binom 31\beta^{3(2)+1}+\frac{1}{3}
\binom 32\beta_1^{3(1)+2}=\frac{1}{T^5}-\frac{1}{T^7}=\frac{
T^2-1}{T^7}.
\]
Por tanto $v_T\big(\frac{T^2-1}{T^7}\big)=-7$ el cual es primo 
relativo a $3$ de donde obtenemos que $\beta_2\notin K$ y por
tanto no existe $\vec\beta\in W_2({\ma F}_2(T))$ tal que 
$\wp(\vec \beta)=\big(\wp\big(\frac{1}{T}\big),0\big)$.
\end{ejemplo}

\begin{teorema}\label{T9'.7.7} Sea $K$ un campo arbitrario de
caracter{\'\i}stica $p$. Entonces dado $\vec x\in W_n(K)$, existe
$\vec y\in W_n(\bar{K})$ tal que $\wp \vec y=\vec x$ donde
$\bar{K}$ denota una cerradura algebraica de $K$.
\end{teorema}

\begin{proof} Sea $\vec x=(x_1,x_2,\ldots,x_n)\in W_n(K)$. Ahora bien,
como $x_1\in K$ existe $y_1\in \bar{K}$ tal que $y_1^p-y_1=x_1$,
es decir, $\wp y_1=x_1$. Se tiene
\[
(\wp y_1,x_2,\ldots,x_n)=\wp((y_1,0,\ldots,0))\Witt + (0,x'_2,\ldots,x'_n).
\]
Por inducci\'on en $n$, existe $(y_2,\ldots,y_n)$ tal que
$\wp((y_2,\ldots,y_n))=(x'_2,\ldots,x'_n)$. Por lo tanto
\begin{align*}
\wp\big(\{y_1\}\Witt + (0,y_2,\ldots, y_n)\big)&=
\wp((y_1,0,\ldots,0))\Witt +\wp((0,y_2,\ldots,y_n))\\
&=\wp((y_1,\ldots,y_n))\\
&=(\wp y_1,x_2,\ldots, x_n)\Witt -(0,x'_2,\ldots, x'_n)\\
&\hspace{1cm} \Witt + \wp
((0,y_2,\ldots,y_n)) \\
&= (\wp y_1,x_2,\ldots,x_n)= (x_1,\ldots, x_n).
\tag*{$\fin$}
\end{align*}
\end{proof}

\begin{notacion}\label{N9'.7.8}
El campo de descomposici\'on de la ecuaci\'on $\wp \vec y=\vec x$,
donde $\vec x=(x_1,\ldots,x_n)\in W_n(K)$, se denota por 
$K(\vec y)=K(y_1,\ldots,y_n)$ donde $\vec y=(y_1,\ldots,y_n)
\in W_n(\bar{K})$. Este campo tambi\'en se denota por
$K(\vec y)=K(\wp^{-1}\vec x)$.
\end{notacion}

\begin{proposicion}\label{P9'.7.9} Sea $\vec y_0$
una soluci\'on de la ecuaci\'on
$\wp \vec y=\vec x$. Entonces todas las soluciones son $\vec y_0
\Witt +\vec m$ con $m\in\{0,1,\ldots,p^n-1\}$.
\end{proposicion}

\begin{proof} Si $\vec y$, $\vec y'\in W_n(\bar{K})$ son tales que
$\wp \vec y=\wp \vec y'$, entonces $\wp(\vec y\Witt - \vec y')
=\vec 0$ de donde se sigue que $\vec y\Witt - \vec y'\in
W_n({\ma F}_p)$. $\fin$
\end{proof}

Ahora consideremos $K$ cualquier campo de caracter{\'\i}stica
$p$, $n\in{\ma N}$ y $\vec x=(x_1,\ldots, x_n)\in W_n(K)$.
Consideremos $\vec y=(y_1,\ldots,y_n)$ soluci\'on a la
ecuaci\'on $\wp \vec y=\vec x$ y $L=K(\wp^{-1} \vec x)=
K(\vec y)=K(y_1,\ldots,y_n)$. Puesto que todas las soluciones de la
ecuaci\'on $\wp \vec y=\vec x$ son $\{\vec y_0\Witt +\vec m\}_{0\leq
m\leq p^n-1}$ con $\vec y_0$ una soluci\'on fija, $L/K$ es una
extensi\'on normal y separable pues $\vec m\in W_n({\ma F}_p)
\subseteqq W_n(K)$. Por lo tanto $L/K$ es una extensi\'on 
de Galois. 

M\'as a\'un, consideremos la funci\'on $\varphi\colon G:=\Gal(
L/K)\to W_n({\ma F}_p)\cong {\ma Z}/p^n{\ma Z}=C_{p^n}$
dada de la siguiente forma:
si $\sigma \in G$, $\sigma \vec y_0=\vec y_0\Witt + \vec m_{\sigma}$
para alg\'un $m_{\sigma}\in\{0,1,\ldots,p^n-1\}$, entonces
$\varphi(\sigma)=\vec m_{\sigma}$. Se tiene que $\varphi$
es un monomorfismo de grupos
y en particular $G\cong C_{p^t}$ para alg\'un $t\leq n$. Puesto
que los subgrupos de $C_{p^n}$ est\'an generados por $\vec p^i$,
$0\leq i\leq n$, se tiene que $G=\langle \sigma_m\rangle$ para
alg\'un $0\leq m\leq n$ donde $\sigma_m(\vec y_0)=\vec y_0
\Witt + \vec p^m$, $o(G)= p^{n-m}$.

Ahora bien, por el (\ref{Eq9'.6.13}) se tiene que 
\begin{gather*}
\vec p=\vec p\Witt \times \vec 1=V(\vec 1^p)=V(\vec 1),\\
\vec p^2=\vec p^2 \Witt \times \vec 1=\vec p\Witt \times (\vec p\Witt \times
\vec 1)= \vec p\Witt \times V(\vec 1) = V(V(\vec 1)^p)=V^2\vec 1,\\
\intertext{y en general}
\vec p^m=V^m(\vec 1), \quad 0\leq m\leq n, \quad \vec p^n=\vec 0.
\end{gather*}

Se tiene 
\begin{align*}
\sigma_m(\vec y_0)&=(\sigma_m y_1,\ldots,\sigma_m y_n)=
\vec y_0\Witt +\vec p^m=(y_1,\ldots,y_n)+V^m(\vec 1)\\
&=(y_1,\ldots,y_m,y_{m+1},\ldots,y_n)\Witt + (\underbrace{0,0,\ldots,0}_{m},
1,0,\ldots,0) \\
&= (y_1,\ldots,y_m,y'_{m+1},\ldots, y'_n),
\end{align*}
por lo que $\sigma_my_j=y_j$ para $1\leq j\leq m$ (para $m=0$ no
hay tal que $y_j$). Esto es, $y_1,\ldots,y_m\in K$. En particular
$\wp y_1=x_1\in \wp(K)$.

En la otra direcci\'on, si $x_1\notin \wp(K)$, como $\wp \vec y=\vec x$, se sigue
que $\wp y_1=x_1$ y por lo tanto $y_1\notin K$. 
Observemos que existe $\sigma\in G$
tal que $\sigma y_1=y_1+1$. Digamos que $\sigma(\vec y_0)=\vec y_0
\Witt + \vec m$, por tanto $\vec m=(1,\alpha_2,\ldots,\alpha_n)$ es invertible
en $W_n({\ma F}_p)$. Sea $t\in {\ma N}$ tal que $\vec t\Witt \times \vec m=\vec 1$,
$\sigma^t(\vec y_0)=\vec y_0\Witt + \vec t\Witt \times \vec m= \vec y_0\Witt +
\vec 1$ y $o(\sigma)= p^n$. En particular $G\cong C_{p^n}$ si y s\'olo si
$x_1\notin \wp(K)$.

Hemos obtenido

\begin{teorema}\label{T9'.7.10}
Sea $K$ un campo de caracter{\'\i}stica $p$, $n\in {\ma N}$ y
$\vec x\in W_n(K)$. Entonces la ecuaci\'on $\wp \vec y=\vec x$ define
una extensi\'on de Galois c{\'\i}clica de $K$: $L=K(\vec y)=K(
y_1,\ldots,y_n)= K(\wp^{-1}\vec x)$. Adem\'as $\Gal(L/K)\cong
C_{p^{n-m}}$ donde $y_1,\ldots, y_m\in K$, $y_{m+1}\notin K$.
Finalmente $L/K$ es c{\'\i}clica de grado $p^n$ si y s\'olo si
$x_1\notin\wp(K)$ donde $\vec x=(x_1,\ldots,x_n)$. En este
\'ultimo caso, $G=\Gal(L/K)$ est\'a generado por $\sigma
(\vec y)=\vec y\Witt + \vec 1$.

Rec{\'\i}procamente, si $L/K$ es una extensi\'on c{\'\i}clica de 
grado $p^n$, existe $\vec x\in W_n(K)$ tal que $L=K(\wp^{-1} \vec x)$, 
esto es, toda extensi\'on c{\'\i}clica de grado $p^n$ se obtiene por
medio de una ecuaci\'on del tipo $\wp \vec y=\vec x$.
\end{teorema}

\begin{proof} \'Unicamente falta demostrar que si $L/K$ es c{\'\i}clica
de grado $p^n$, entonces existe $\vec y\in W_n(\bar{K})$ tal que
$L=K(\vec y)$ y donde $\wp \vec y=\vec x$.

Sea $G=\Gal(L/K)=\langle \sigma\rangle$, $o(\sigma)=p^n$. 
Se tiene que $\vec 1\in W_n(L)$ satisface $\Tr_{L/K} \vec 1=
\sum\limits_{\sigma\in G}^{\bullet} \sigma \vec 1= \vec p^n\Witt \times
\vec 1 = \vec p^n=\vec 0$. Por el Teorema \ref{T9'.7.2}, existe
$\vec y\in W_n(L)$ tal que $(\sigma\Witt - \vec 1)\vec y=\vec 1$,
esto es, $\sigma \vec y= \vec y\Witt + \vec 1$. Sea $\wp
\vec y=\vec x$. Entonces
\[
\sigma(\wp \vec y)=\wp(\sigma \vec y)=\wp (\vec y\Witt +\vec 1)=
\wp(\vec y)\Witt +\wp (\vec 1)=\wp (\vec y),
\]
es decir, $\wp \vec y=\vec x\in W_n(K)$. Puesto que $\sigma
\vec y=\vec y\Witt +\vec 1$, $K(\vec y)\subseteqq L$ y $K(\vec y)/K$
es una extensi\'on c{\'\i}clica de orden $p^n$ de donde se
sigue que $L=K(\vec y)=K(\wp^{-1} \vec x)$. $\fin$
\end{proof}

\begin{corolario}\label{C9'.7.11} Sea $L=K(\vec y_1)=
K(\vec y_2)$, $\vec y_1, \vec y_2\in W_n(L)$ una extensi\'on
c{\'\i}clica de grado $p^n$ con $\wp \vec y_i=\vec x_i\in W_n(K)$,
$i=1,2$. Entonces existen $\vec j\in W_n({\ma F}_p)$ invertible,
es decir, $\mcd(j,p)=1$,
y $\vec z\in W_n(K)$ tales que 
\[
\vec y_1=\vec j\Witt \times \vec y_2
\Witt + \vec z\quad \text{y}\quad \vec x_1=\vec 
j\Witt \times \vec x_2\Witt + \wp \vec z
\]
y rec{\'\i}procamente.
\end{corolario}

\begin{proof} Sea $G:=\Gal (L/K)=\langle \sigma\rangle$ tal que $\sigma
\vec y_1=\vec y_1\Witt +\vec 1$. Como $\sigma$ est\'a completamente
determinado por $\sigma \vec y_2$, se tiene que $\sigma 
\vec y_2=\vec y_2\Witt + \vec i$ con $\vec i$ invertible en
$W_n({\ma F}_p)$ o, equivalentemente, $\mcd (i,p)=1$.
Sea $\vec j\in W_n({\ma F}_p)$ tal que $\vec j\Witt \times
\vec i=\vec 1$. Entonces $\sigma(\vec j\Witt \times \vec y_2)=
\vec j\Witt \times \vec y_2\Witt + \vec j\Witt \times \vec i=
\vec j\Witt \times \vec y_2 \Witt +\vec 1$. Por tanto
$\sigma(\vec y_1\Witt -\vec j\Witt \times \vec y_2)=
\vec y_1\Witt - \vec j\Witt \times \vec y_2$ lo cual implica
que $\vec y_1\Witt -\vec j\Witt \times \vec y_2=\vec z
\in W_n(K)$.

El rec{\'\i}proco es claro. $\fin$
\end{proof}

Consideramos $L=K(y_1,\ldots,y_n)/K$ una extensi\'on c{\'\i}clica
de grado $p^n$ con $\wp \vec y=\vec x$. Puesto que $y_n\notin
K(y_1,\ldots,y_{n-1})$ (pues de lo contrario tendr{\'\i}amos
$[L:K]\leq p^{n-1}$), se sigue que $L=K(y_n)$. Ahora bien si
$K_{n-1}=K(y_1,\ldots,y_{n-1})$, se tiene que
\begin{align*}
\wp \vec y&= \wp\Big(\sum_{i=0}^{\substack{n-1\\ \bullet}}
V^i(\{y_{i+1}\})\Big)=\wp\Big(\sum_{i=0}^{\substack{n-2\\ \bullet}}
V^i(\{y_{i+1}\})\Big) \Witt + \wp\big(V^{n-1}(\{y_n\})\big)\\
&= \wp((y_1,\ldots, y_{n-1},0))\Witt + \wp((0,\ldots, 0, y_n))=\vec x,
\end{align*}
$\wp((0,\ldots,0, y_n))=(0,\ldots,0, y_n^p-y_n)$. Por tanto, tomando
la componente fantasma $n$, si sigue que $y_n^p-y_n=
z_{n-1}+x_n$ con $z_{n-1}\in K_{n-1}$.

Por el Teorema \ref{T9'.2.4} se tiene que si $\langle \sigma\rangle
=\Gal(L/K)$, $\varphi_{\mid K_{n-1}}$, entonces $\wp y_n=
y_n^p-y_n=z_{n-1}+x_n$, $(\sigma -1)y_n = \delta$ con $
(\varphi-1)z_{n-1}=\wp \delta$. Con esto recuperamos el
resultado de Schmid (\ref{Ec9'.2.1}) para la generaci\'on de una
extensi\'on c{\'\i}clica de grado $p^n$.

\begin{teorema}\label{T9'.7.12} Sea $L/K$ una extensi\'on c{\'\i}clica
de grado $p^n$ y sean $K_0=K\subseteqq K_1\subseteqq \ldots
\subseteqq K_{n-1}\subseteqq K_n=L$ tales que $[K_i:K]=p^i$ y sea
$\langle \sigma_i\rangle=\Gal(K_i/K)$, $K_i=K_{i-1}(y_i)$. Se tiene
$\sigma_i=\sigma_{n\mid K_i}$. Entonces $L/K$ satisface
\begin{gather}\label{Ec9'.7.13}
\begin{array}{rlllr}
K_1&=k(y_1),& \wp y_1=x_1, &(\sigma_1-1) y_1=1,\\
K_2&=K_1(y_2),& \wp y_2=z_1+x_2,& (\sigma_2-1) y_2=c_1,
& (\sigma_1-1)z_1=\wp c_1\\
K_3&=K_2(y_3),& \wp y_3=z_2+x_3, & (\sigma_3-1) y_3=c_2,
& (\sigma_2-1)z_2=\wp c_2\\
\hspace{0.3cm}\vdots &\hspace{.8cm}\vdots&\hspace{1cm} \vdots&\hspace{1cm}\vdots
&\hspace{1cm}\vdots\\
L=K_n&=K_{n-1}(y_n), & \wp y_n=z_{n-1}+x_n, & (\sigma_n-1) y_n=c_{n-1},
&(\sigma_{n-1}-1)z_{n-1}\\
&&&&=\wp c_{n-1}
\end{array}
\end{gather}
donde $z_i, c_i\in K_i$, $1\leq i\leq n-1$, $y_i\in K_i$ ($y_i\notin K_{i-1}$)
y $x_1,\ldots,x_n\in K$ con $x_1\notin\wp (K)$. Toda extensi\'on c{\'\i}clica
de grado $p^n$ est\'a determinado por elementos
arbitrarios $x_1,\ldots,x_n\in K$ con $x_1
\notin \wp(K)$. $\fin$
\end{teorema}

\section{Sobre el conductor}\label{9'.8}

Primero se define el {\em conductor\index{conductor}} para campos locales.
Un {\em campo local\index{campo local}} es un campo completo bajo una
valuaci\'on no arquimediana con campo residual finito.

Sea $K$ un campo completo con 
respecto a una valuaci\'on discreta $v$ y el cual tiene campo residual finito.
Sea $L/K$ una extensi\'on abeliana finita. Sea $n$ el m{\'\i}nimo
$n\in{\ma N}\cup\{0\}$ tal que $U_K^{(n)}\subseteqq N_{L/K} L^{\ast}$ donde
${\eu p}$ es el ideal primo del anillo de enteros ${\cal O}_K=\{x\in K\mid
v(x)\geq 0\}$, es decir, ${\eu p}=\{x\in K\mid v(x)>0\}$ y $U^{(n)}:= 1+
{\eu p}^n$.

Por teor{\'\i}a de campos de clase, Teorema \ref{T17.3.2.7}, Subsecci\'on
\ref{STeoremadeExistencia} y 
Proposici\'on \ref{CCUnidades} se tiene
que $N_{L/K} L^{\ast}$ es abierto en $K^{\ast}$ y $1\in N_{L/K} L^{\ast}$
por lo que tal $n$ existe por ser $\big\{U^{(n)}_K\big\}_{n=0}^{\infty}$
un sistema fundamental de vecindades abiertas de $1$.

\begin{definicion}\label{D9'.8.2}
Se define el {\em conductor\index{conductor}}
de $L/K$ por 
\[
{\eu f}={\eu f}_K={\eu f}(L/K):={\eu p}^n.
\]
\end{definicion}

Se tiene que una extensi\'on abeliana de campos locales $L/K$
es no ramificada si y s\'olo si su conductor es ${\eu f}=1$ (Teorema
\ref{CClaseT3.2.24}).

\begin{teorema}\label{T9'.8.2'}
Sea $K$ un campo local. El mapeo 
\[
L\mapsto {\cal N}_L:=
N_{L/K} L^{\ast}
\]
establece una correspondencia uno a uno entre las extensiones abelianas
finitas $L/K$ y los subgrupos abiertos de {\'\i}ndice finito en $K^{\ast}$.
M\'as a\'un:
\begin{gather*}
L_1\subseteq L_2\iff {\cal N}_{L_1}\supseteq {\cal N}_{L_2},
\quad {\cal N}_{L_1L_2}= {\cal N}_{L_1}\cap {\cal N}_{L_2},
\quad {\cal N}_{L_1\cap L_2}={\cal N}_{L_1} {\cal N}_{L_2}.
\end{gather*}
\end{teorema}

\begin{proof}
Teorema \ref{CClaseT3.2.29}.
$\fin$
\end{proof}

Ahora, en {\em campos globales\index{campo global}}, es decir, extensiones
finitas de ${\ma Q}$ o de ${\ma F}_p(t)$ se usan los conceptos de
{\em ad\`eles\index{ad\`eles}} y de {\em id\`eles\index{id\`eles}}. Sea $K$ un
campo global. Un {\em ad\`ele} o {\em id\`ele aditivo} o {\em repartici\'on}
es una familia $\alpha=(\alpha_{\eu p})_{\eu p}$ con $\alpha_{\eu p}\in 
K_{\eu p}$, la completaci\'on de $K$ en $\eu p$ y adem\'as $\alpha_{\eu p}$
es entero para casi todo $\eu p$, esto es, $v_{\eu p}(\alpha_{\eu p})\geq 0$
para casi todo $\eu p$ y donde $\eu p$ recorre todos los divisores primos
de $K$, incluyendo los primos infinitos.

Se denota ${\ma A}_K:= \prod'_{\eu p}K_{\eu p}$ (producto
restringido\index{producto restringido})
el {\em anillo de los ad\`eles\index{anillo de ad\`eles}} con suma
y multiplicaci\'on por componentes. El grupo de {\em id\`eles\index{id\`eles}}
de $K$ se define como el grupo de las unidades de ${\ma A}_K^{\ast}$,
es decir, $ J_K={\ma A}_K^{\ast}$. Esto es,
un id\`ele es una familia $\alpha=(\alpha_{\eu p})_{\eu p}$, $\alpha_{\eu p}
\in K_{\eu p}^{\ast}$ y $\alpha_{\eu p}$ es una unidad en el anillo de 
enteros ${\cal O}_{\eu p}$ de $K_{\eu p}$ para casi toda $\eu p$. Se
escribe $ J_K=\prod'_{\eu p} K_{\eu p}^{\ast}$. De manera natural
$K^{\ast}\subseteqq  J_K$ v{\'\i}a el mapeo diagonal,
$K^{\ast}$ son los {\em id\`eles
principales\index{id\`eles!principales}\index{id\`eles principales}}.
Sea $C_K:= J_K/K^{\ast}$. El grupo $C_K$ recibe el nombre de el
{\em grupo de clases de id\`eles\index{id\`eles!grupo de clases de
$\sim$}\index{grupo de clases de id\`eles}}.

\begin{proposicion}\label{P9'.8.3} Si $K$ es un campo num\'erico,
$I_K$ denota el grupo de clases de $K$,
y $ J_K^{S_{\infty}}=\prod_{{\eu p}\mid \infty}K_{\eu p}^{\ast}\times
\prod_{{\eu p}\nmid \infty}U_{\eu p}$, entonces $I_K\cong  J_K/
 J_K^{S_{\infty}} K^{\ast}\cong C_K/\big(( J_K^{S_{\infty}}K^{\ast})/
K^{\ast}\big)$.

En el caso de $K$ un campo de funciones, si $I_{K,0}$ denota
el grupo de clases de $K$ de grado $0$ y $C_{K,0}$ denota
al grupo de id\`eles de grado $0$, $I_{K,0}\cong \frac{C_{K,0}}
{\tilde{C}_K}\cong \frac{C_{K,0}}{\bar{U}}$, donde $\bar{U}=
U\*K/\*K$, $U=\prod_{\pK\in{\ma P}_K}$ y $I_K\cong 
\frac{C_K}{\tilde{C}_K}$.
\end{proposicion}

\begin{proof}
Subsecci\'on \ref{CClaseS4.7}. $\fin$
\end{proof}

Un {\em m\'odulus\index{m\'odulus}} es un ideal entero ${\eu m}
=\prod_{{\eu p}\nmid \infty} {\eu p}^{n_{\eu p}}$ de ${\cal O}_K$, el
anillo de enteros de $K$ al cual lo consideramos como ${\eu m}=
\prod_{{\eu p}} {\eu p}^{n_{\eu p}}$ con $n_{\eu p}=0$ para
${\eu p}\mid \infty$. Se define para cada lugar ${\eu p}$ de $K$:
\[
U_{\eu p}^{(0)}= U_{\eu p}\quad \text{y}\quad U_{\eu p}^{(n_{\eu p})}=
\begin{cases}
1+{\eu p}^{n_{\eu p}}&\text{si ${\eu p}\nmid \infty$},\\
{\ma R}_+^{\ast}\subseteqq K_{\eu p}^{\ast}&\text{si ${\eu p}$ es real},\\
{\ma C}^{\ast}=K_{\eu p}^{\ast}&\text{si ${\eu p}$ es complejo},
\end{cases}
\]
para $n_{\eu p}>0$. Se define $\alpha_{\eu p} \equiv 1
 \bmod {\eu p}^{n_{\eu p}}
\iff \alpha_{\eu p}\in U_{\eu p}^{(n_{\eu p})}$. Esta definici\'on 
corresponde a la congruencia usual si ${\eu p}$ es finito, $\alpha_{\eu p}
>0$ si ${\eu p}$ es real y es una condici\'on vac{\'\i}a para ${\eu p}$
complejo.

\begin{definicion}[Definici\'on \ref{CClaseD4.4.5}]\label{D9'.8.4}
Si $K$ es un campo num\'erico,
se define $C_K^{\eu m}:= J_K^{\eu p} K^{\ast}/K^{\ast}$
donde $ J_K^{\eu m}:=\prod_{\eu p}U_{\eu p}^{(n_{\eu p})}$.
El grupo $C_K^{\eu m}$ recibe el nombre de {\em grupo de
congruencia\index{grupo de congruencia}\index{congruencia!grupo de
$\sim$} m\'odulo ${\eu m}$}. Al grupo $C_K/C_K^{\eu m}$ se le
llama el {\em grupo de rayos\index{grupo de rayos} m\'odulo ${\eu m}$}.
\end{definicion}

Si $L/K$ es una extensi\'on de Galois, hay una norma $N_{L/K}\colon
 J_L\to  J_K$ definida como sigue (Teorema \ref{CClaseTCCG}).
Sea ${\eu p}$ un lugar de $K$ y sea $L_{\eu p}:=\prod_{{\eu P}\mid
{\eu p}} L_{\eu P}$. Cada $\alpha_{\eu p}\in L_{\eu p}^{\ast}$ define
un automorfismo $\alpha_{\eu p}\colon L_{\eu p}\to L_{\eu p}$,
$x\mapsto \alpha_{\eu p}x$ del $K_{\eu p}$--espacio vectorial
$L_{\eu p}$. Se define la norma de $\alpha_{\eu p}$ por:
$N_{L_{\eu p}/K_{\eu p}}(\alpha_{\eu p}) =\deg (\alpha_{\eu p})$.
Se induce un homomorfismo $N_{L/K}\colon  J_L\to
 J_K$ dado por: si $\alpha=(\alpha_{\eu P})\in  J_L$,
las componentes locales de $N_{L/K}(\alpha)$ est\'an dadas por
(Teorema \ref{CClaseT4.2.1}).
\[
N_{L/K}(\alpha)_{\eu p}=\prod_{{\eu P}\mid {\eu p}} N_{
L_{\eu P}/K_{\eu p}}(\alpha_{\eu P}).
\]

Ahora $N_{L/K}$ manda id\`eles principales en id\`eles principales
y por tanto la norma induce otra norma $N_{L/K}\colon C_L\to
C_K$ (Teorema \ref{CClaseT4.2.1}). Se tiene

\begin{teorema}\label{T9'.8.5}
Sea $K$ un campo global. El mapeo 
\[
L\mapsto {\cal N}_L:=
N_{L/K} C_L
\]
es una correspondencia uno a uno entre las extensiones abelianas
finitas $L/K$ y los subgrupos cerrados de {\'\i}ndice finito en $C_K$.
M\'as a\'un:
\begin{gather*}
L_1\subseteqq L_2\iff {\cal N}_{L_1}\supseteq {\cal N}_{L_2},
\quad {\cal N}_{L_1L_2}= {\cal N}_{L_1}\cap {\cal N}_{L_2},
\quad {\cal N}_{L_1\cap L_2}={\cal N}_{L_1} {\cal N}_{L_2}.
\end{gather*}
El campo $L/K$ que corresponde al subgrupo ${\cal N}$ de $C_K$
se llama el {\em campo de clase\index{campo de clase} de ${\cal N}$}.
Se tiene
\begin{gather*}
\Gal(L/K)\cong C_K/{\cal N}.
\end{gather*}
\end{teorema}

\begin{proof}
Teoremas \ref{CClaseT4.2.1} y \ref{CClaseTC.1}.
$\fin$
\end{proof}

Ahora bien, entre los grupos cerrados de $C_K$ de {\'\i}ndice finito
est\'an los grupos de congruencias $C_K^{\eu m}$ donde ${\eu m}=
\prod_{\eu p}{\eu p}^{n_{\eu p}}$.

\begin{definicion}[{Definici\'on \ref{CClaseD4.5.3}}]\label{D9'.8.6}
Si $K$ es un campo num\'erico,
el campo de clase $K^{\eu m}/K$ que corresponde a $C_K^{\eu m}$
se llama el {\em campo de clase de rayos\index{campo de
clase de rayos} m\'odulo ${\eu m}$}.
\end{definicion}

Se tiene $\Gal(K^{\eu m}/K)\cong C_K/C_K^{\eu m}$ y 
si ${\eu m}\mid {\eu m}'$ entonces $K^{\eu m}\subseteqq
K^{{\eu m}'}$ (Observaci\'on \ref{CClaseO4.5.4}). Ahora bien si
${\cal N}$ es cualquier grupo cerrado de {\'\i}ndice finito en $C_K$,
${\cal N}$ contiene a un subgrupo de congruencia $C_K^{\eu m}$
y por tanto toda extensi\'on abeliana $L/K$ est\'a contenida
en un campo de clase $K^{\eu m}/K$.

\begin{definicion}[{Definici\'on \ref{CClaseD4.5.5}}]\label{D9'.8.7}
Sea $L/K$ una extensi\'on abeliana finita. Sea ${\cal N}_L=
N_{L/K} C_L$. Se define el {\em conductor\index{conductor}
${\eu f}(L/K)={\eu f}$} de $L/K$ como el m\'aximo com\'un
divisor de los m\'odulus ${\eu m}$ tales que $L\subseteqq K^{
\eu m}$: ${\eu f}:=\mcd\{{\eu m}\mid {\eu m}\text{\ es un m\'odulus
y\ }L\subseteqq K^{\eu m}\}$. Esto es, ${\eu f}=\mcd\{{\eu m}\mid
{\eu m}\text{\ es un m\'odulus y\ } C_K^{\eu m}\subseteqq
{\cal N}_L\}$.
\end{definicion}

En otras palabras $K^{\eu f}/K$ es el m{\'\i}nimo campo de clase
de rayos que contiene a $L/K$. Una observaci\'on interesante
es que ${\eu m}$ no necesariamente es el conductor de 
$K^{\eu m}/K$, es decir, puede existir ${\eu f}\mid {\eu m}$,
${\eu f}\neq {\eu m}$ y $K^{\eu f}=K^{\eu m}$ (o $C_K^{\eu m}=
C_K^{\eu f}$).

La relaci\'on entre los conductores locales y los conductores
globales es:

\begin{teorema}\label{T9'.8.8}
Si ${\eu f}$ es el conductor de una extensi\'on abeliana finita
$L/K$ de campos globales y ${\eu f}_{\eu p}$ es el conductor
de la extensi\'on local $L_{\eu p}/K_{\eu p}$ ($L_{\eu p}$ una
extensi\'on de $K_{\eu p}$, es decir, si ${\eu P}_1,\ldots,{\eu P}_r
\mid {\eu p}$, escogemos cualquier ${\eu P}_i\mid {\eu p}$
y ponemos $L_{\eu p}:=L_{{\eu P}_i}$), entonces si definimos
${\eu f}_{\eu p}=1$ para ${\eu p}$ infinito, se tiene
\begin{gather*}
{\eu f}=\prod_{\eu p}{\eu f}_{\eu p}.
\end{gather*}
\end{teorema}

\begin{proof}
Teorema \ref{CClaseT4.5.10}. 
$\fin$
\end{proof}

\begin{corolario}\label{C9'.8.9}
Dada una extensi\'on abeliana finita de campos
globales $L/K$, se tiene que
${\eu p}$ se ramifica en $L$ si y s\'olo si ${\eu p}\mid{\eu f}$.
\end{corolario}

\begin{proof}
Corolario \ref{CClaseC4.5.11}.
$\fin$
\end{proof}

\begin{ejemplo}[Corolario \ref{CClaseT.4.5.8}]\label{E9'.8.10}
Los campos de clase de rayos de ${\ma Q}$ son 
precisamente los campos ciclot\'omicos pues los m\'odulus
est\'an dados por ${\eu m}=(m)$, $m\in{\ma N}$ y ${\ma Q}^{
\eu m}={\ma Q}(\zeta_m)$ y en particular toda extensi\'on
abeliana de ${\ma Q}$ est\'a contenida en un campo
ciclot\'omico. Se sigue de esto una prueba del Teorema de
Kronecker--Weber\index{Kronecker--Weber!teorema de
$\sim$}\index{teorema de Kronecker--Weber}
\index{Kronecker--Weber!teorema de $\sim$} para campos
num\'ericos: la m\'axima extensi\'on abeliana de ${\ma Q}$
es ${\ma Q}^{ab}=\cup_{n=1}^{\infty}{\ma Q}(\zeta_n)$.
\end{ejemplo}

\begin{observacion}\label{O9'.8.11} Vemos que dado un campo
num\'erico $K$, los campos $K^{\eu m}$ son los an\'alogos
a los campos ${\ma Q}(\zeta_m)$, pues cada extensi\'on abeliana
de $K$ est\'a contenida en alg\'un $K^{\eu m}$ y cada extensi\'on
abeliana de ${\ma Q}$ est\'a contenida en alg\'un ${\ma Q}(
\zeta_m)$. Ahora bien, la gran diferencia es que sabemos la
existencia de los campos $K^{\eu m}$ pero no como est\'an
generados, a diferencia de los campos ciclot\'omicos ${\ma Q}
(\zeta _m)$ que expl{\'\i}citamente est\'an dados por las ra{\'\i}ces
de la ecuaci\'on $x^m-1$.
\end{observacion}

\begin{observacion}\label{O9'.8.12} Dado un campo num\'erico
$K$, el campo de clase $K^1$ corresponde a la m\'axima 
extensi\'on abeliana de $K$ no ramificada en ning\'un primo
finito. Este campo es usualmente llamado el {\em campo
de clase de Hilbert extendido\index{Hilbert!campo de clase extendido de
$\sim$}\index{campo de clase extendido de Hilbert}}. Los
primos infinitos pueden o no ser ramificados en $K^1/K$.
Notemos que ${\ma Q}^1={\ma Q}$.
\end{observacion}

\begin{definicion}\label{D9'.8.13} El subcampo $K\subseteqq
K_H\subseteqq K^1$ tal que los primos infinitos de $K$ son
no ramificados o, equivalentemente, se descomponen totalmente,
se llama {\em el campo de clase de Hilbert\index{Hilbert!campo
de clase de Hilbert}\index{campo de clase de Hilbert}}.
\end{definicion}

Recordemos el Teorema \ref{T9.3}.

\begin{teorema}[Campo de clase de Hilbert\index{campo de clase
de Hilbert}\index{Hilbert!campo de clase de $\sim$}]\label{T9'.8.14}
Se tiene que el grupo de Galois de $K_H/K$ satisface
$\Gal(K_H/K)\cong I_K$, el grupo de clases de $K$.
\end{teorema}

\begin{proof}
Corolario \ref{CClaseC4.5.13}.
$\fin$
\end{proof}

\subsection{Representaciones, caracteres y conductores}\label{S9'.8.1}

\begin{definicion}\label{D9'.8.15} Una {\em representaci\'on\index{representaci\'on
de un grupo}} de un grupo finito $G$ es una acci\'on de $G$ en un
${\ma C}$--espacio vectorial de dimensi\'on finita $V$. Equivalentemente,
una representaci\'on es un homomorfismo de grupos
\[
\rho\colon G\to \GL(V)=\Aut_{\ma C}(V).
\]
\end{definicion}

Una acci\'on la podemos entender como: $\varphi\colon G\times V\to V$
con $\varphi(\sigma,v):=\sigma\circ v:= \rho(\sigma)(v)$. Tambi\'en es com\'un
usar la notaci\'on $(V, \rho)$ para indicar la representaci\'on de $G$ en $V$.

La {\em representaci\'on trivial\index{representaci\'on
trivial}} es $(\rho,{\ma C})$ con $\rho(\sigma)=1$ para
toda $\sigma\in G$. El {\em grado\index{grado de una
representaci\'on}\index{representaci\'on!grado de una $\sim$}} de la 
representaci\'on es la dimensi\'on de $V$.

Una representaci\'on $\rho$ se llama {\em irreducible}\index{representaci\'on
irreducible} si $V$ no admite ning\'un subespacio propio $0\varsubsetneqq W
\varsubsetneqq V$ que sea $G$--invariante, es decir, $\sigma\circ
W\subseteqq W$ para toda $\sigma\in G$.

\begin{proposicion}\label{P9'.8.16} Si $G$ es abeliano, toda representaci\'on
irreducible de $G$ es de grado $1$, es decir es un caracter
\begin{gather*}
\rho\colon G\to \GL_1({\ma C})\cong {\ma C}^{\ast}. \tag*{$\fin$}
\end{gather*}
\end{proposicion}

Se tiene que toda representaci\'on $(\rho,V)$ se factoriza a trav\'es de una
suma directa $V=V_1\oplus\cdots\oplus V_s$ de representaciones
irreducibles. M\'as precisamente, $\rho_i\colon G\to V_i$ es una 
representaci\'on irreducible y $\rho=\rho_1\oplus\cdots\oplus \rho_s$.
Expl{\'\i}citamente si $\{v_{ij}\}_{j=1}^{m_i}$ es una base de $V_i$
y consideramos la base de $V$ dada por $\{v_{ij}\}_{1\leq i\leq s}^{1\leq j\leq m_i}$
y si la matriz de $\rho_i(\sigma)$ con respecto a la base
$\{v_{ij}\}_{j=1}^{m_i}$ es la matriz $(m_i\times m_i)$ $A_{\sigma, i}$,
entonces $\rho(\sigma)$ es la matriz
\[
A_{\sigma}=\left(
\begin{array}{ccc}
\text{\fbox{${A_{\sigma,1}}$}}&&0\\
&\ddots&\\
0&& \text{\fbox{$A_{\sigma, s}$}}
\end{array}
\right)
\]
con respecto a la base $\{v_{ij}\}_{1\leq i\leq s}^{1\leq j\leq m_i}$.

Dos representaciones $(\rho,V)$ y $(\rho',V')$ se llaman
{\em equivalentes\index{representaciones equivalentes}} si existe
un isomorfismo $\varphi\colon V\to V'$ de $G$--espacios vectoriales,
esto es, $\varphi$ es un isomorfismo de espacios vectoriales
tal que $\varphi(\sigma\circ v)=\sigma\circ \varphi(v)$ para toda
$\sigma\in G$ y toda $v\in V$.

Si en la suma de $G$--espacios $V=V_1\oplus\cdots\oplus V_s$ una
representaci\'on $(\rho_{\alpha}, V_{\alpha})$ tiene $r_{\alpha}$
representaciones equivalentes entre las representaciones 
$(\rho_1, V_1),\ldots, (\rho_s, V_s)$, se usa la notacion:
\[
\rho\sim \sum_{\alpha}r_{\alpha}\rho_{\alpha}
\]
y $r_{\alpha}$ se llama la {\em multiplicidad} de $\rho_{\alpha}$
en $\rho$.

\begin{definicion}\label{D9'.8.17} Dada una representaci\'on
$(\rho, V)$ el {\em caracter\index{caracter de una representaci\'on}}
de $\rho$ se define por
\[
\chi_{\rho}\colon G\to {\ma C},\quad \chi_{\rho}(\sigma)=
\text{traza de\ }\rho(\sigma).
\]
\end{definicion}

Se tiene que si $\rho\sim \sum_{\alpha}r_{\alpha}
\rho_{\alpha}$, entonces $\chi_{\rho}=\sum_{\alpha}
r_{\alpha}\chi_{\rho_{\alpha}}$.

Un caracter $\chi$ se llama {\em irreducible\index{caracter irreducible}}
si $\chi$ es el caracter de una representaci\'on irreducible. Una 
{\em funci\'on central\index{funci\'on central}} o {\em funci\'on
de clase\index{funci\'on de clase}} es una funci\'on $f\colon G\to
{\ma C}$ tal que $f(\sigma \tau \sigma^{-1})= f(\tau)$ para
cualesquiera $\sigma,\tau\in G$.

\begin{teorema}\label{T9'.8.18} Se tiene que dos representaciones son
equivalentes si y s\'olo si  sus caracteres son iguales. $\fin$
\end{teorema}

Se tiene que toda funci\'on central $\varphi$ se puede escribir
un{\'\i}vocamente como una combinaci\'on lineal 
\[
\varphi=\sum_{\chi}c_{\chi}\chi, \quad c_{\chi}\in {\ma C},
\]
donde los $\chi$'s son caracteres irreducibles.

\begin{teorema}\label{T9'.8.19} Se tiene que $\varphi$ es el
caracter de una representaci\'on si y s\'olo si $c_{\chi}\in
{\ma N}\cup\{0\}$ para toda $\chi$. $\fin$
\end{teorema}

Se tiene que si $(\rho, V)$ es una representaci\'on con
caracter $\chi$, se tiene $\dim V^G=\frac{1}{|G|}
\sum_{\sigma\in G}\chi(\sigma)$.

\subsection{Conductores de Artin}\label{S9'.8.2}

Sea $L/K$ una extensi\'on de Galois de campos globales. Sea 
$G=\Gal(L/K)$. Dado un caracter irreducible $\chi$ de $G$,
se define el ideal ${\eu f}(\chi)$ por:
\begin{gather*}
{\eu f}(\chi):=\prod_{{\eu p}\nmid \infty}{\eu p}^{f_{\eu p}(\chi)}
=\prod_{{\eu p}\nmid \infty} {\eu f}_{\eu p}(\chi)\\
\intertext{con}
f_{\eu p}(\chi)=\sum_{i\geq 0} \frac{g_i}{g_0}\codim V^{G_i}
\end{gather*}
donde $V$ es una representaci\'on con caracter $\chi$,
$G_i$ es el $i$--\'esimo grupo de ramificaci\'on de $L_{\eu P}/
K_{\eu p}$, y $g_i$ es el orden de $G_i$ donde ${\eu P}$
es cualquier divisor en $L$ dividiendo a ${\eu p}$.

\begin{definicion}\label{D9'.8.20} Al ideal ${\eu f}(\chi)$ se le
llama el {\em conductor de Artin\index{Artin!conductor de $\sim$}\index{conductor
de Artin}} del caracter $\chi$.

Para campos locales, ${\eu f}_{\eu p}(\chi)={\eu p}^{f(\chi)}$ se define
como el {\em conductor local de Artin\index{Artin!conductor local
de $\sim$}\index{conductor local de Artin}} del caracter $\chi$.
\end{definicion}

En realidad el conductor local de Artin se define de manera m\'as general
como sigue. Sea $L/K$ una extensi\'on de campos locales, con grupo
de Galois $G=\Gal(L/K)$. Sea $f$ el grado de inercia de $L/K$. Se
define $i_G(\sigma)=v_L(\sigma x-x)$ donde $x$ es cualquier
elemento tal que ${\cal O}_L={\cal O}_K[x]$ y $v_L$ es la valuaci\'on
de $L$. Sea 
\[
a_G(\sigma)=
\begin{cases}
-f i_G(\sigma)&\text{para $\sigma\neq 1$},\\
f\sum_{\tau\neq 1} i_G(\tau)&\text{para $\sigma =1$}.
\end{cases}
\]
Se tiene que $a_G$ es una funci\'on central sobre $G$ y se puede
escribir
\[
a_G=\sum_{\chi}f(\chi)\chi,\qquad f(\chi)\in {\ma C},
\]
donde $\chi$ var{\'\i}a sobre los caracteres irreducibles de $G$.
Se tiene que $f(\chi)$ es un entero no negativo y por lo tanto
podemos formar el ideal $f_{\eu p}(\chi)={\eu p}^{f(\chi)}$, que 
ser\'a la ${\eu p}$--componente del conductor de Artin global.

La relaci\'on entre los conductores de Artin y los conductores
local y global antes definidos, se obtiene de los siguientes
resultados.

\begin{teorema}\label{T9'.8.21}
Sea $L/K$ una extensi\'on de Galois de campos locales y sea
$\chi$ un caracter de $\Gal(L/K)$ de grado $1$. Sea $L_{\chi}=
L^{\ker \chi}$ el campo fijo del n\'ucleo de $\chi$. Sea ${\eu f}$
el conductor de $L_{\chi}/K$. Entonces 
\begin{gather*}
{\eu f}={\eu f}_{\eu p}(\chi).
\end{gather*}
\end{teorema}

\begin{proof}
\cite[Proposition 11.6, p. 532]{Neu99}.
$\fin$
\end{proof}

En el caso global tenemos:

\begin{teorema}\label{T9'.8.22}
Sea $L/K$ una extensi\'on de campos globales, $\chi$ un caracter
de $\Gal(L/K)$ de grado $1$. Sea $L_{\chi}$ el campo fijo de 
$\ker \chi$ y sea ${\eu f}_{\chi}$ el conductor global de la
extensi\'on $L_{\chi}/K$. Entonces
\begin{gather*}
{\eu f}_{\chi}={\eu f}(\chi).
\end{gather*}
\end{teorema}

\begin{proof}
\cite[Proposition 11.10, p. 535]{Neu99}.
$\fin$
\end{proof}

Como consecuencia de lo anterior, puesto que una extensi\'on abeliana
$L/K$, ya sea de campos locales o globales, se tiene que los
conductores de Artin y los conductores usuales son lo mismo.
As{\'\i}, en el caso abeliano, los conductores se pueden considerar, en
caso de as{\'\i} convenir, como conductores de Artin.
Por otro lado se tiene la siguiente f\'ormula hallada por H.
Hasse y E. Artin.

\begin{teorema}[F\'ormula del conductor--discriminante\index{f\'ormula
del conductor--discriminante}]\label{T9'.8.23}
Para cualquier extensi\'on finita de Galois $L/K$ de campos
globales, se tiene
\begin{gather*}
{\eu d}_{L/K}=\prod_{\chi}{\eu f}(\chi)^{\chi(1)}
\end{gather*}
donde $\chi$ recorre el conjunto de todos los caracteres irreducibles
de $\Gal(L/K)$ y ${\eu d}_{L/K}$ denota el discriminante de la
extensi\'on $L/K$.
\end{teorema}

\begin{proof}
\cite[11.9, p. 534]{Neu99}.
$\fin$
\end{proof}

Notemos que $\chi(1)$ es precisamente la dimensi\'on
de la representaci\'on asociada a $\chi$.

En lo que resta de este cap{\'\i}tulo el campo  de funciones 
racionales sobre ${\ma F}_q$ ser\'a denotado por $K$: 
$K={\ma F}_q(T)$. Como es usual $R_T$ denota el
anillo de polinomios en $T$ sobre ${\ma F}_q$:
$R_T={\ma F}_q[T]$.

Sea ahora $\chi\colon \big(R_T/(M)\big)^{\ast}\to {\ma C}^{\ast}$
un caracter de Dirichlet. Se tiene que si ${\eu f}_{\chi}$ es el
conductor de $\chi$ (como caracter de Dirichlet) y si ${\eu f}'_{\chi}$
es el conductor de Artin de $\chi$, entonces ${\eu f}_{\chi}
={\eu f}'_{\chi}$.
Por otro lado el conductor de Dirichlet de un caracter $\chi$ es
$P^{\alpha}$ con $P\in R_T$ es m\'onico e irreducible, si y s\'olo
si $\chi\colon \big(R_T/(P^{\alpha})\big)^{\ast}\to {\ma C}^{\ast}$
pero no puede definirse m\'odulo $P^{\alpha-1}$:
$\chi\colon\big(R_T/(P^{\alpha-1})\big)^{\ast}\to {\ma C}^{\ast}$.

\subsection{Conductor local de $K(\Lam P {\alpha})$}\label{S9'.8.3}

Se tiene que $K(\Lam P{\alpha})_{\eu B}\cong {\ma F}_q((\lam P{\alpha}))$
pues $P$ es totalmente ramificado y $v_{\eu B}(\lam P{\alpha})=1$,
donde $\lambda =\lam P{\alpha}$ es generador de $\Lam P{\alpha}$ y
${\eu B}$ es el primo en $K(\Lam P{\alpha})$ sobre $P$. Ahora se tiene
que $N_{K(\Lam P{\alpha})_{\eu B}/K_P}\big(K(\Lam P{\alpha})^{\ast}_{\eu B}\big)
=(P)\times U_P^{(\alpha)}$ (Teorema \ref{CClaseT3.2.5.32}).

Como consecuencia, se tiene que el conductor local de $K(\Lam P{\alpha})
/K$ en $P$ es $P^{\alpha}$ y $1$ para cualquier otro $Q\neq P$, $Q
\in R_T$ irreducible.

\begin{lema}\label{L9'.8.24} Supongamos que
$\K\subseteqq K(\Lam P{\beta})$ para alg\'un $\beta$ y sea ${\eu f}_\K=
P^{\gamma}$. Entonces $\gamma\leq \alpha\iff \K\subseteqq
K(\Lam P{\alpha})$.
\end{lema}

\begin{proof} 

\noindent
$\Longrightarrow$)  Supongamos que $\K\nsubseteqq K(\Lam P{\alpha})$
y sea $L=\K K(\Lam P{\alpha})\varsupsetneqq K(\Lam P{\alpha})$. Entonces
por el Teorema \ref{T9'.8.2'} se tiene que ${\cal N}_L\varsubsetneqq
{\cal N}_{K(\Lam P{\alpha})}=(P)\times U^{(\alpha)}$. Sea ${\eu f}_L=P^{\delta}$.
Entonces $U^{(\delta)}\subseteqq {\cal N}_L$ y $U^{(\delta -1)}\nsubseteqq
{\cal N}_L$. Adem\'as $P\in {\cal N}_L$ lo cual implica que
\[
(P)\times U^{(\delta)}\subseteqq {\cal N}_L\varsubsetneqq {\cal N}_{
K(\Lam P{\alpha})} = (P)\times U^{(\alpha)},
\]
de donde se sigue que $\delta > \alpha$. Ahora bien se tiene que
\[
{\eu f}_L={\eu f}_{\K K(\Lam P{\alpha})}={\eu f}^{\max\{\gamma,\alpha\}}=
P^{\gamma}=P^{\delta}
\]
por lo que obtenemos $\delta=\gamma>\alpha$. Este absurdo prueba que
$\K\subseteqq K(\Lam P{\alpha})$.

\noindent
$\Longleftarrow$) Se tiene que $U^{(\alpha)}\subseteqq {\cal N}_{K(\Lam P{\alpha})}
\subseteqq {\cal N}_\K$, de donde se sigue $\alpha\geq \gamma$. $\fin$
\end{proof}

Como consecuencia tenemos

\begin{proposicion}\label{P9'.8.25} Sea $\K\subseteqq K(\Lam P{\beta})$
para alg\'un $\beta\in {\ma N}$. Entonces ${\eu f}_\K=P^{\alpha}$ si y s\'olo
si $\K\subseteqq K(\Lam P{\alpha})$ y $\K\nsubseteqq K(\Lam P{\alpha-1})$.
\end{proposicion}

\begin{proof} 
Si $\K\subseteqq K(\Lam P{\alpha})$ y $\K\nsubseteqq K(\Lam P{\alpha-1})$
entonces si ${\eu f}_\K=P^{\gamma}$ se sigue del Lema \ref{L9'.8.24} que
$\gamma\leq \alpha$ y $\gamma \not\leq \alpha-1$ por lo que $\gamma=\alpha$.

Rec{\'\i}procamente, si ${\eu f}_\K=P^{\alpha}$ entonces nuevamente
por el Lema \ref{L9'.8.24}, y puesto que $\alpha\leq \alpha$ se sigue que
$\K\subseteqq K(\Lam P{\alpha})$. Ahora bien, si tuvi\'esemos
$\K\subseteqq K(\Lam P{\alpha-1})$ entonces se seguir{\'\i}a que
$\alpha\leq \alpha-1$ lo que prueba 
que $\K\nsubseteqq K(\Lam P{\alpha-1})$.
$\fin$
\end{proof}

\subsection{El conductor de acuerdo a Schmid}\label{S9'.8.4}

Nuestro objetivo en esta subsecci\'on es enunciar el c\'alculo de
Schmid para el conductor en una extensi\'on c{\'\i}clica determinada
por un vector de Witt. Primero volvemos a presentar un caso 
particular de la Proposici\'on \ref{P2.4.Ram3}.

\begin{proposicion}\label{P9'.8.26} Sea $\K/K$ una extensi\'on 
c{\'\i}clica de grado $p$ tal que $\K\subseteqq K(\Lam P{\beta})$ para
alg\'un $\beta\in{\ma N}$. Entonces existe $y\in \K$ tal que
$\K=K(y)$ con $\wp y=y^p-y=h(T)\in K$ con $h(T)=\frac{g(T)}{
P(T)^{\lambda}}$ con $g(T)\in R_T$, $\mcd (P(T), g(T))=1$, $\lambda>0$ y
$\mcd(\lambda,p)=1$.
\end{proposicion}

\begin{proof} Por el Teorema \ref{T9'.2.1} se tiene que existe $y\in \K$
tal que $\K=K(y)$ y $y^p-y=h(T)\in K$. Se  tiene
$\mu(X)=\Irr(y,X,K)=X^p-X-h(T)$ con $h(T)\in K$ y $h(T)\notin \wp(K)=
\{a^p-a\mid a\in K\}$. Sea
$h(T)=\frac{g(T)}{f(T)}$ con $g(T),f(T)\in R_T$,
$\mcd (g(T),f(T))=1$, $f(T)=\prod_{i=1}^r P_i^{\alpha_i}$, donde
$P_1, \ldots, P_r$ son polinomios irreducibles distintos.
Descomponiendo a $h(T)$ en fracciones parciales, se tiene
\[
\frac{g(T)}{f(T)}=s(T)+\sum_{i=1}^{r} \frac{t_i(T)}{P_i(T)^{\alpha_i}},
\quad \deg t_i(T)<\deg P_i(T)^{\alpha_i},
\quad t_i(T), s(T)\in R_T.
\]

Se tiene que para cualquier divisor primo ${\eu p}\notin\{{\eu p}_1,\ldots,
{\eu p}_r,{\eu p}_{\infty}\}$, donde ${\eu p}_i$ denota al divisor primo
correspondiente a $P_i$, se tiene $v_{\eu p}(y^p-y)=v_{\eu p}(h)\geq 0$.
Es decir, $y$ es entero con respecto a ${\eu p}$. Se sigue que
\[
\mu(X)=\prod_{i=0}^{p-1}(X-y-i)\quad \text{y}\quad \mu'(X)=
\sum_{i=0}^{p-1}\prod_{j\neq i}(X-y-i).
\]
Por lo tanto $\mu'(y)=\sum_{i=0}^{p-1}\prod_{j\neq i}(y-y-j)=
\prod_{j=1}^{p-1}(-j)$ es una unidad y por el Teorema \ref{T5.5.6},
$\eu p$ es no ramificado.

Ahora si para alg\'un $1\leq i\leq r$, $p\mid \alpha_i$ entonces si
$\alpha_i=\lambda_i p$, podemos escribir
\[
\frac{g(T)}{f(T)}=\frac{t_0(T)}{P_i(T)^{\lambda_i p}}+s_i(T)\quad\text{con}
\quad v_{{\eu p}_{i}}(s_i(T))\geq 0, \quad \deg t_0(T)<\deg P_i(T)^{\lambda_i p}.
\]
Ahora $K(T)/(P_i(T))$ es un campo finito, por lo tanto perfecto. Existe
$m(T)\in K(T)$ tal que $m(T)^p\equiv t_0(T)\bmod P_i(T)$. Sea $n(T):=
-\frac{m(T)}{P_i(T)^{\lambda_i}}$. Sea $u=y+n(T)$. Entonces $\K=
K(u)=K(y)$ y 
\[
u^p-u=h(T)+n(T)^p-n(T)=h_0(T)
\]
con $v_{{\eu p}_i}(h_0(T))>-\lambda_i p$, $v_{{\eu p}_j}(h_0(T))=
v_{{\eu p}_j}(h(T))$ para $j\neq i$ y $v_{\eu p}(h(T))\geq 0$ para
${\eu p}\notin\{{\eu p}_1,\ldots,{\eu p}_r, {\eu p}_{\infty}\}$. Continuando
con este proceso, se obtiene $\K=K(w)$ con $w^p-w= \ell(T)$ con
$(\ell(T))_K=\frac{{\eu C}}{{\eu p}_1^{\lambda_1} \cdots {\eu p}_m^{\lambda_m}}
{\eu p}_{\infty}^s$ con ${\eu C}$ un divisor entero primo relativo a
${\eu p}_1,\ldots,{\eu p}_m$, $m\leq r$, $\lambda_i>0$, 
$\mcd(\lambda_i,p)=1$ donde reenumeramos a los
elementos de $\{{\eu p}_1,\ldots,{\eu p}_r\}$ del conjunto inicial
en caso de ser necesario.

Para $1\leq i\leq m$ se tiene $v_{{\eu p}_i}(w^p-w)=v_{{\eu B}_i}(\ell(T))=
e({\eu P}_i\mid{\eu p})v_{{\eu p}_i}(\ell(T))=-e({\eu P}_i\mid{\eu p})\lambda_i<0$,
donde ${\eu P}_i$ es un divisor en $\K$ sobre ${\eu p}_i$.
Por tanto $v_{{\eu P}_i}(w)<0$ y $v_{{\eu P}_i}(w^p-w)=v_{{\eu P}_i}
(w^p)=pv_{{\eu P}_i}(w)=-e({\eu P}_i\mid{\eu p})\lambda_i$.
Puesto que $\mcd(p,\lambda_i)=1$, se tiene que $p\mid
e({\eu P}_i\mid{\eu p})$ y por tanto ${\eu p}_i$ es ramificado en $\K/K$.
Por otro lado, como $\K\subseteqq K(\Lam P{\alpha})$, el \'unico
primo finito ramificado es ${\eu p}$, el divisor asociado a $P$ y por
tanto se tiene $\K=K(y)$ con
\begin{gather*}
y^p-y=\frac{g(T)}{P(T)^{\lambda}}\quad\text{con}\quad g(T)\in R_T,
\quad \mcd(g(T),P(T))=1,\\
\lambda>0\quad\text{y}\quad
\mcd(\lambda,p)=1. \tag*{$\fin$}
\end{gather*}
\end{proof}

\begin{corolario}\label{C9'.8.27} Si $\K/K$ es una extensi\'on
c{\'\i}clica de grado $p^n$ con $\K\subseteqq K(\Lam P{\alpha})$ para
alg\'un $\alpha\in{\ma N}$, entonces existe $\vec y$ tal que
$\K=K(\vec y)$ con $\vec y^p\Witt- \vec y =\vec {\beta}\in W_n(K)$ con
$\beta_i(T)=\frac{g_i(T)}{P(T)^{\lambda_i}}$ con $g_i(T)\in R_T$,
$\lambda_i\geq 0$ y si $\lambda_i>0$ entonces $\mcd(g_i(T), P(T))=1$ y
$\mcd(\lambda_i,p)=1$. Finalmente, $\lambda_1>0$.
\end{corolario}

\begin{proof} Procedemos por inducci\'on en $n$. El caso $n=1$ es la Proposici\'on
\ref{P9'.8.26}. Para el caso $n+1$, $\K_{n+1}=
\K_n(y_{n+1})$, $y_{n+1}^p-y_{n+1}=
z_n+\beta_{n+1}$ con $z_n\in \K_n$ y
$v_{{\eu P}}(z_n)\geq 0$ para todo
divisor ${\eu P}$ que no divide a $P$. Por tanto, por el proceso de
la demostraci\'on de la
Proposici\'on \ref{P9'.8.26} se tiene que $\beta_{n+1}$ tiene la forma
requerida. Notemos que puede ser que $\lambda_{n+1}=0$ pues
la codificaci\'on de la ramificaci\'on de los primos de $\K_n$ sobre $P$
bien se pudiera presentar en $z_n$. $\fin$
\end{proof}

\begin{observacion}\label{O9'.8.28} 
Hemos desarrollado el caso
particular en que la extensi\'on $\K/K$ es c{\'\i}clica de grado $p^n$ y
$\K\subseteqq K(\Lam P{\alpha})$ para alg\'un $\alpha\in{\ma N}$.
En este caso particular se tiene que $\lambda_1>0$ pues $P$
es totalmente ramificado en $\K/K$. Sin embargo la Proposici\'on
\ref{P9'.8.26} y el Corolario \ref{C9'.8.27} pueden ser generalizados
de manera natural a una extensi\'on arbitrara $\K/K$ c{\'\i}clica de
grado $p^n$. En este caso tenemos el {\em criterio de Schmid\index{criterio
de Schmid}\index{Schmid!criterio de $\sim$}}
\l
\item Para $n=1$, $\K=K(y)$ con $y^p-y=h(T)\in K$ y tal que
\begin{gather*}
(h(T))_K=\frac{\eu C}{{\eu p}_1^{\lambda_1}\cdots {\eu p}_r^{\lambda_r}}
\quad \text{con}\quad \mcd ({\eu C},{\eu p}_i)=1, \quad 1\leq i\leq r, \\
\lambda_i>0 \quad \text{y}\quad \mcd (\lambda_i,p)=1.
\end{gather*}
Los primos ramificados son precisamente ${\eu p}_1,\ldots, {\eu p}_r$.

\item Para $n$ arbitraria, $\K=K(\vec y)$, $\vec y^p\Witt - \vec y=\vec \beta
\in W_n(K)$ tal que 
\begin{gather*}
(\beta_i)_K=\frac{{\eu C}_i}{{\eu p}_1^{\lambda_{1,i}}
\ldots {\eu p}_r^{\lambda_{r,i}}}\quad \text{con}\quad \lambda_{j,i}\geq 0
\quad \text{y si}\quad \lambda_{j,i}>0,\\
\mcd ({\eu C}_i, {\eu p}_j)=1
\quad \text{y}\quad \mcd(\lambda_{j,i},p)=1.
\end{gather*}
El {\'\i}ndice de ramificaci\'on de cada ${\eu p}_j$ es $p^{n-i+1}$
donde $i$ es el primer {\'\i}ndice $i$ tal que $\lambda_{j,i}>0$. En
otras palabras, ${\eu p}_j$ es no ramificado en $\K_{i-1}:=
K(y_1,\ldots,y_{i-1})/K$ y totalmente ramificado en $\K=\K_n/\K_i$.
\end{list}

Este desarrollo se puede encontrar en
\cite[Section 2, p\'agina 162]{Sch36},
\cite[Section 3, p\'agina 115]{Sch36-0}.
\end{observacion}

Aqu\'i presentamos una demostraci\'on de esta reducci\'on de Schmid.

\begin{teorema}[Reducci\'on de Schmid,
\cite{Sch36,Sch36-0}]\label{T9'.8.28-1}
Sea $k$ un campo perfecto de caracter\'istica $p>0$ y sea $K$
un campo de funciones sobre $k$. Sean $\vec \alpha=(\alpha_1,
\ldots,\alpha_n)\in W_n(K)$ y $\pK$ un lugar de $K$. Entonces
existe $\vec\gamma\in W_n(K)$ tal que $\vec \beta:=\vec \alpha
\Witt+\wp(\vec\gamma)=(\beta_1,\ldots,\beta_n)$ satisface lo siguiente:
para cada $1\leq i\leq n$, $\beta_i$ cumple $v_{\pK}(\beta_i)
\geq 0$ \'o $v_{\pK}(\beta_i)<0$ y $p\nmid v_{\pK}(\beta_i)$.
En otras palabras, si $\beta_i\neq 0$, el divisor de $\beta_i$,
$(\beta_i)_K$ es de la forma
\begin{gather*}
(\beta_i)_K=\frac{{\eu a}_i}{\pK^{\lambda_i}}\quad\text{con}\quad
v_{\pK}({\eu a}_i)\geq 0\quad\text{y}\quad \lambda_i=0 \quad
\text{\'o}\\ 
v_{\pK}({\eu a}_i)=0, \quad
 \lambda_i>0\quad\text{y}\quad \mcd(p,\lambda_i)=1.
\end{gather*}
\end{teorema}

\begin{proof}
Lo haremos por inducci\'on en $n$. Si $n=1$, este es el resultado
de H. Hasse \cite{Has35} (ver \cite[Theorem 5.8.10]{Vil2006}).

Suponemos que, para $m<n$, las primeras $m$ componentes 
del vector $\vec\alpha\Witt+\wp(\vec\gamma)=
(\beta_1,\ldots,\beta_m,\xi_{m+1},\ldots,\xi_n)$ 
est\'an dados de la forma descrita, es decir, $v_{\pK}
(\beta_i)\geq 0$ \'o $v_{\pK}(\beta_i)<0$ y $\mcd(p,v_{\pK}(\beta_i))
=1$ para $1\leq i\leq m$. 
El paso de inducci\'on se sigue inmediatamente si $\xi_{m+1}$
tambi\'en est\'a en la forma prescrita.

Por tanto suponemos que $\xi_{m+1}$ no est\'a en la forma 
prescrita, esto es, $v_{\pK}(\xi_{m+1})=-p\mu$ con $\mu\in{\ma N}$.
Por el procedimiento de Hasse, consideremos $\nu
\in K$ tal que $\xi_{m+1}+\wp(\nu)=:\beta_{m+1}$ 
 y $v_{\pK}(\beta_{m+1})\geq 0$
\'o $v_{\pK}(\beta_{m+1})<0$ y $\mcd(p,v_{\pK}(\beta_{m+1}))=1$.

Ahora
\begin{align*}
\vec\alpha\Witt+\wp(\vec\gamma)\Witt+&V^m(\wp(\{\nu\})
=\vec\alpha\Witt+\wp(\vec\gamma)\Witt+\wp(V^m\{\nu\})\\
&=(\beta_1,\ldots,\beta_m,\xi_{m+1},\ldots,\xi_n)\Witt+
\wp((0,\ldots,0,\nu,0,\ldots,0))\\
&=(\beta_1,\ldots,\beta_m,\xi_{m+1},\ldots,\xi_n)\Witt+
(0,\ldots,0,\wp(\nu),x_{m+2},\ldots,x_n)\\
&=(\beta_1,\ldots,\beta_m,\xi_{m+1}+\wp(\nu),x'_{m+2},\ldots,
x'_m)\\
&=(\beta_1,\ldots,\beta_m,\beta_{m+1},x'_{m+2},\ldots,x'_n)
\end{align*}
con $\beta_i$, $1\leq i\leq m+1$ en la forma prescrita. Aqu\'i hemos
usado que $\wp\circ V^m=V^m\circ \wp$ (Teorema \ref{T9'.6.1}).
El resultado se sigue.
$\fin$
\end{proof}

\begin{corolario}\label{C9'.8.28-2}
Sea $\vec \beta$ el vector de Witt obtenido en el Teorema
{\rm{\ref{T9'.8.28-1}}} y sea $(\beta_i)_K=\frac{{\eu a}_i}{\pK^{\lambda_i}}$
con $v_{\pK}({\eu a}_i)\geq 0$ y $\lambda_i=0$ \'o $v_{\pK}({\eu a}_i)=0$,
$\lambda_i>0$ y $\mcd(p,\lambda_i)=1$. Sean $\wp(\vec \theta)=
\vec\theta^p\Witt - \vec\theta=\vec \beta$ con $\vec\theta\in W_n(\bar{K})$
y $L=K(\vec \theta)=K(\theta_1,\ldots,\theta_n)$. Sean $L_i=
K(\theta_1,\ldots,\theta_i)$ los subcampos intermedios entre $K$ y $L$
donde ponemos $L_0=K$ y tenemos $L_n=L$. Sea $1\leq \mu\leq n$
el primer \'indice tal que $\lambda_{\mu}>0$ en caso de existir. Entonces
el \'indice de ramificaci\'on de $\pK$ en $L/K$ es $p^{n-\mu+1}$. En
particular $\pK$ es totalmente ramificado si y solamente si $\mu=1$.
En el caso de que $\lambda_i=0$ para todo $1\leq i\leq n$, se tiene que
$\pK$ es no ramificado.
\end{corolario}

\begin{proof}
La extensi\'on $L_i/K$ est\'a dada por $L_i=K(\theta_1,\ldots,\theta_i)$ donde
$\wp((\theta_1,\ldots,\theta_i))=(\theta_1,\ldots,\theta_i)^p\Witt -(\theta_1,
\ldots,\theta_i)=(\beta_1,\ldots,\beta_i)$. Por el Teorema \ref{T9'.7.12} tenemos
que $\wp(\theta_i)=z_{i-1}+\beta_i$. Por inducci\'on se sigue que, si $
\lambda_1=\cdots=\lambda_{i-1}=0$, entonces $v_{\pK}(z_{i-1})\geq 0$
y $\pK$ es no ramificado en $L_{i-1}/K$. Para $i$, si $i<\mu$
tenemos que $\lambda_i=0$ y por tanto $v_{\pK}(z_{i-1}+\beta_i)
\geq 0$ de donde se sigue que $\pK$ es no ramificado en $L_i/L_{i-1}$
por el Teorema de Hasse \cite{Has35} o \cite[Theorem 5.8.10]{Vil2006}.

As\'i, si $\lambda_1=\cdots=\lambda_i=0$, entonces $\pK$ es no 
ramificado en $L_i/K$. Ahora $\lambda_1=\cdots=\lambda_{\mu-1}=0$,
$\lambda_{\mu}>0$, por lo anterior, se tiene que $\pK$ es no ramificado
en $L_{\mu-1}/K$. Adem\'as $v_{\pK}(z_{\mu-1}+\beta_{\mu})=
v_{\pK}(\beta_{\mu})=\lambda_{\mu}<0$. Por el Teorema de Hasse,
$\pK$ es ramificado en $L_{\mu}/L_{\mu-1}$. Puesto que $L/L_{\mu-1}$
es una extensi\'on c\'iclica de orden una potencia de un primo, se sigue
que $\pK$ es totalmente ramificado en $L/L_{\mu-1}$ y por tanto
el \'indice de ramificaci\'on de $\pK$ en $L/K$ es $[L:L_{\mu-1}]=p^{n-\mu+1}$.
$\fin$
\end{proof}

H. L. Schmid introdujo los siguientes invariantes. Sea $\K/K$ una extensi\'on
c{\'\i}clica de grado $p^n$ con $\K\subseteqq K(\Lam P{\alpha})$ para
alg\'un $\alpha\in{\ma N}$. Sea $\K=K(\vec y)$ tal que $\vec y^p\Witt -
\vec y= \vec \beta\in W_n(K)$, $(\beta_i)=\frac{{\eu C}_i}{{\eu p}^{\lambda_i}}$
con $\lambda_i\geq 0$ y si $\lambda_i>0$ entonces $\mcd ({\eu C}_i,
{\eu p})=1$ y $\mcd(\lambda_i,p)=1$ donde ${\eu p}$ es el divisor
asociado a $P$.

Sea $M_n:=\max\limits_{1\leq i\leq n}\{p^{n-i}\lambda_i\}$. Notemos que
$M_i=\max\{pM_{i-1}, \lambda_i\}$, $M_1<M_2<\cdots <M_n$ y que el valor
m\'aximo se alcanza en un \'unico $p^{n-i}\lambda_i$ pues si
$p^{n-i}\lambda_i=p^{n-j}\lambda_j$ con $j>i$, entonces $p^{j-i}
\lambda_i=\lambda_j$ pero esto contradice que $\mcd(p,\lambda_j)=1$.

\begin{teorema}\label{T9'.8.29} 
Con las condiciones anteriores, se tiene que
el conductor de $\K/K$ es
\[
{\eu f}_\K=P^{M_n+1}.
\]
\end{teorema}

\begin{proof} \cite[p. 163]{Sch36}. $\fin$
\end{proof}

\begin{corolario}\label{C9'.8.30} Sea $\K/K$ una extensi\'on c\'iclica de
grado $p^n$ con $\K\subseteqq K(\Lam P{\alpha})$ para alg\'un $\alpha
\in{\ma N}$. Entonces $M_n+1\leq \alpha$.
\end{corolario}

\begin{proof} Se sigue inmediatamente del 
Lema \ref{L9'.8.24} y del Teorema
\ref{T9'.8.29}. $\fin$
\end{proof}

\subsection{Ramificaci\'on en $p$--extensiones c\'iclicas en
caracter\'istica $p$}\label{S8(1).pea}

En esta secci\'on mencionamos brevemente el tipo de ramificaci\'on
de los primos en una $p$--extensi\'on c\'iclica $\K/\F(T)=K$ dada
por una ecuaci\'on de Witt: $\K=\F(T)(\vec y)$:
\[
\vec y^p \Witt - \vec y=\vec\alpha\in W_n(\F(T)).
\]

Este caso es similar al de extensiones de Artin--Schreier. 
Aqu{\'\i} consideramos
$\K=K(\vec{y})$ donde $\vec{y}^p\Witt - \vec{y}=\vec {\beta}$, y
la operaci\'on es la diferencia de Witt. La extensi\'on es una
$p$--extensi\'on finita de grado menor o igual a $p^n$ donde
$\vec y$ es de longitud $n$. Sean $P_1,\ldots P_r$ los
divisores primos finitos que son ramificados en $\K/K$.

\begin{teorema}\label{T5.3.1.pea}\label{T5.3.1}
 Sea $\K/K$ una extensi\'on c{\'\i}clica
de grado $p^n$ donde se tiene que
$P_1,\ldots,P_r\in R_T^+$ y posiblemente $\p$, 
son los divisores primos ramificados. Entonces
$\K=K(\vec y)$ donde
\[
\vec y^p\Witt -\vec y=\vec \beta={\vec\delta}_1\Witt + \cdots \Witt + {\vec\delta}_r
\Witt + \vec\mu,
\]
con $\beta_1^p-\beta_1\notin \wp(K)$,
$\delta_{ij}=\frac{Q_{ij}}{P_i^{e_{ij}}}$, $e_{ij}\geq 0$, $Q_{ij}\in R_T$ y
\l
\item si $e_{ij}=0$ entonces $Q_{ij}=0$;
\item si $e_{ij}>0$ entonces $p\nmid e_{ij}$, $\mcd(Q_{ij},P_i)=1$ y
$\deg (Q_{ij})<\deg (P_i^{e_{ij}})$, 
\end{list}
y ${\mu}_j=f_j(T)\in R_T$ con
\l
\setcounter{bean}{2}
\item $p\nmid \deg f_j$ cuando $f_j\not\in {\ma F}_q$ y
\item $\mu_j\notin\wp({\ma F}_q):=\{a^p-a\mid a\in{\ma F}_q\}$ cuando 
$\mu_j\in {\ma F}_q^{\ast}$.
\end{list}
\end{teorema}

\begin{proof}
Consideremos $\K/K$ una extensi\'on c{\'\i}clica de grado $p^n$ definida
por $\K:=K(\vec y)$, $\vec y^p\Witt -\vec y=\vec \beta$
con $\vec y \in W_n(\K)$ un vector de Witt de longitud $n$ en $\K$
y $\vec\beta\in W_n(K)$ un vector de Witt de longitud $n$ en $K$.

Sea $\vec \beta=(\beta_1,\ldots,\beta_n)$ tal que 
\begin{gather}\label{Eq2.cap13}
\beta_j=\sum_{i=1}^r \frac{Q_{ij}}{P_i^{e_{ij}}}+f_j(T), \text{\ donde\ }
P_1,\ldots,P_r\in R_T^+, \big\{Q_{ij}\big\}_{1\leq i\leq r}^{1\leq j\leq n}
\subseteq R_T,\nonumber\\
f_j(T)\in R_T, e_{ij}\in{\ma N}\cup \{0\} \text{\ para toda\ }
1\leq i\leq r \text{\ y\ } 1\leq j\leq n.
\end{gather}

Sea $\varphi$ definido como en (\ref{Eq9'.2.2}).
Aplicando $\varphi$ a $\vec \beta$ obtenemos el vector
$\vWitt {\beta}n$ y de la definici\'on de $\beta^{(j)}$, se tiene que
que estos elementos son de la forma
\begin{gather*}
\beta^{(j)}=\sum_{i=1}^r \frac{Q'_{ij}}{P_i^{e'_{ij}}} +f'_j(T)\text{\ para todo\ }
1\leq j\leq n.\\
\intertext{Escribimos}
\vec \beta=\vec {\gamma_1}+\cdots+\vec {\gamma_r}+\vec {\xi},\\
\vWitt {\beta}n=\vWitt {\gamma_1}n+\cdots+\vWitt {\gamma_r}n+
\vWitt {\xi}n\\
\intertext{con}
\gamma_i^{(j)}=\frac{Q'_{ij}}{P_i^{e'_{ij}}},\quad 1\leq i\leq r,\quad
1\leq j\leq n\quad \text{y}\quad \xi^{(j)}=f'_j(T).
\end{gather*}

Cuando aplicamos $\varphi^{-1}$, obtenemos
\[
(\beta_1,\ldots,\beta_n)=\vWitt {\beta}n^{\varphi^{-1}} =
({\vec {\gamma}_1})^{\varphi^{-1}}\Witt +\cdots\Witt + ({\vec {\gamma}_r})^{
\varphi^{-1}}\Witt + ({\vec {\xi}})^{\varphi^{-1}}
\]
y cada vector $({\vec {\gamma_i}})^{\varphi^{-1}}$ es de la forma
$\Big(\frac{Q''_{i1}}{P_i^{e''_{i1}}},\cdots, \frac{Q''_{in}}{P_i^{e''_{in}}}\Big)$
y el vector $(\vec \xi)^{\varphi^{-1}}$ es de la forma $\big(f''_1(T),\ldots,
f''_n(T)\big)$. En otras palabras
\[
\vec \beta={\vec {\delta}}_1\Witt +\cdots\Witt + {\vec {\delta}}_r
\Witt + \vec \mu
\]
en donde las componentes de cada ${\vec \delta}_i$ tienen polos
a lo m\'as en $P_i$ y $\vec \mu$ tiene componentes con polos
a lo m\'as en $\p$. Sea ${\eu p}_i$ el divisor correspondiente
a $P_i$.

Ahora cada ${\vec\delta}$ y $\vec\mu$ pueden ser normalizados de manera que
cada componente $({\vec\delta}_i)_j:=\delta_{ij}$ tiene divisor
\begin{gather}
\big(\delta_{ij}\big)_K=\frac{{\eu a}_{ij}}{{\eu p}_i^{\lambda_i}} \text{\ con\ }
\lambda_i\geq 0;
\text{\ si \ }\lambda_i=0, \text{\ entonces\ }v_{{\eu p}_i}({\eu a}_{ij})\geq 0;\nonumber\\
\text{\ si \ } \lambda_i>0, \text{\ entonces\ } \mcd(p,\lambda_i)=1 \text{\ y\ }
v_{{\eu p}_{ij}}({\eu a}_{ij})=0\label{E5.3.1(1)}
\end{gather}
y similarmente para $\vec \mu$ con respecto a $\p$ (ver \cite[p\'agina 62]{Sch36}
y Observaci\'on \ref{O9'.8.28}).
De hecho, la normalizaci\'on puede ser obtenida mediante el cambio de variable
${\vec y}_i\mapsto {\vec y}_i
\Witt +{\vec \alpha}_i$ la cual corresponde a la substituci\'on
${\vec \delta}_i\mapsto {\vec \delta}_i\Witt + {\vec\alpha}_i^p\Witt -
{\vec\alpha}_i$ y por lo tanto las componentes de cada ${\vec \alpha}_i$
no tiene otro polo que no sea ${\eu p}_i$. M\'as directamente, en el
nivel $j$, aplicando la t\'ecnica de Hasse, una sustituci\'on $y_j\to
y_j+\xi_j$ nos lleva a la normalizaci\'on de $\delta_{ij}$ en la forma normal
dada en (\ref{E5.3.1(1)}). $\fin$
\end{proof}

Ahora estudiamos el comportamiento de $\p$ en $K/k$.

\begin{proposition}\label{P5.3.2.pea}\label{P5.3.2}
Sea $K/k$ dado como el en Teorema {\rm{\ref{T5.3.1.pea}}}. Sea
$\mu_1=\cdots=\mu_s=0$, $\mu_{s+1}\in {\ma F}_q^{\ast}$, 
$\mu_{s+1}\not\in \wp({\ma F}_q)$ y finalmente sea $t+1$ el prime 
{\'\i}ndice con $f_{t+1}\not\in{\ma F}_q$ (y por lo tanto $p\nmid \deg f_{t+1}$).
Entonces el {\'\i}ndice de ramificaci\'on de 
$\p$ es $p^{n-t}$, el grado de inercia de
$\p$ es $p^{t-s}$ y el n\'umero de descomposici\'on de
$\p$ es $p^s$. M\'as precisamente, si $\Gal(K/k)=\langle \sigma\rangle
\cong C_{p^n}$, entonces el grupo de inercia de $\p$ es
${\eu I}=\langle \sigma^{p^t} \rangle$ y el grupo de descomposici\'on de
$\p$ es ${\eu D}=\langle \sigma^{p^s}\rangle$.
\end{proposition}

\begin{proof}
Puesto que la extensi\'on $K/k$ es una extensi\'on de Galois
de orden una potencia de un n\'umero primo, el campo de
inercia es el primer nivel en donde $\p$ se ramifica. El 
{\'\i}ndice de este primer nivel es $t+1$ (ver \cite{Sch36}).
Por otro lado, por la misma raz\'on, el campo de descomposici\'on
es el primer nivel donde $\p$ es inerte y este est\'a dado por $s+1$
(Proposici\'on \ref{P2.4.Ram3}). $\fin$
\end{proof}

\section{Extensiones multic\'iclicas}\label{S8.pea}

Sea $q=p^n$ y consideremos un campo $k$ tal que $\F\subseteq k$. 
Sea $\vec x^q:=\big(x^q_1,\ldots, x^q_m
\big)$. Entonces $\vec x^q\Witt - \vec x=0\iff
\vec x\in W_m(\F)\subseteq W_m(k)$. El anillo $W_m(\F)$ es lo que se
conoce como un {\em anillo de Galois}. Como grupo, veamos que
$W_m(\F)$ es un $W_m({\ma F}_p)$--m\'odulo libre de rango $n$
y en particular, $W_m(\F)\cong\Big({\ma Z}/p^m{\ma Z}\Big)^n$
como grupo.

\begin{proposicion}\label{P8.1.pea} Se tiene que $W_m(\F)$ es un
$W_m({\ma F}_p)$--m\'odulo libre de rango $n$, donde $q=p^n$.
M\'as precisamente, sea $\{\mu_1,\ldots,\mu_n\}$ una base de
$\F$ sobre ${\ma F}_p$ y sean $\vec \mu_i:=\{\mu_i\}=
\big(\mu_i,0,\ldots,0\big)$, $1\leq i\leq n$. Entonces $\{\vec \mu_1,
\ldots,\vec \mu_n\}$ es una $W_m({\ma F}_p)$--base de $W_m(\F)$.
Esto es
\[
W_m(\F)=\bigoplus_{i=1}^n W_m({\ma F}_p) \cdot \vec \mu_i.
\]
\end{proposicion}

\begin{proof}
Sean $\vec \alpha_1,\ldots, \vec\alpha_n\in W_m({\ma F}_p)$,
con
\[
\vec \alpha_i=\big(\alpha_{i1},\ldots,\alpha_{im}\mid 
\alpha_i^{(1)},\ldots,\alpha_i^{(m)}\big)\quad
\vec \mu_i=\big(\mu_i,0\ldots,0\mid
\mu_i,\mu_i^p,\ldots,\mu_i^{p^{m-1}}\big).
\]
Entonces se tiene que
$\vec \alpha_i\Witt\cdot \vec\mu_i=\big(\alpha_{i1}\mu_i,\ldots,?\mid
\alpha_{i1}\mu_i,\ldots, ?\big)$. En caso de que $\Witt \sum_{i=1}^n
\vec\alpha_i\vec\mu_i=
\big(\ldots\mid \sum_{i=1}^n \alpha_{i1}\mu_i,\ldots\big)=\vec 0
=\big(0,\ldots,0\mid 0,\ldots, 0\big)$, obtenemos que $\sum_{i=1}^n
\alpha_{i1}\mu_i=0$ lo cual implica que $\alpha_{i1}=0$ para todo
$1\leq i\leq n$.

As\'i tenemos que $\vec \alpha_i=(0,\alpha_{i2},\ldots,\alpha_{im}\mid
0,p\alpha_{i2},\ldots)$. Por tanto la segunda entrada
fantasma de $\Witt\sum_{i=1}^n
\vec\alpha_i\Witt\cdot \vec\mu_i=\vec 0$ queda
$\vec 0^{(2)}=\sum_{i=1}^n p\alpha_{i2}\mu_i^p$. Por lo tanto
\[
0=\vec 0_2=\frac{1}{p}\big(\vec 0^{(2)}-\vec 0_1^p\big)=
\sum_{i=1}^n\alpha_{i2}\mu_i^p=\sum_{i=1}^n\alpha_{i2}^p\mu_i^p=
\Big(\sum_{i=1}^n\alpha_{i2}\mu_i\Big)^p=0,
\]
de donde $\sum_{i=1}^n\alpha_{i2}\mu_i=0$ lo cual implica que
$\alpha_{i2}=0$ para todo $1\leq i\leq n$. Continuando este proceso
obtenemos $\vec \alpha_i=\vec 0$ de donde se sigue el resultado. $\fin$
\end{proof}

\begin{observacion}\label{O8.1'.pea}
En general 
se tiene que $\{\vec \xi_1,\ldots,\vec 
\xi_n\}\subseteq W_m(\F)$ es una base sobre $W_m({\ma F}_p) \iff \{\xi_{11},
\ldots,\xi_{n1}\}$ es base de $\F$ sobre ${\ma F}_p$. Esto se puede 
demostrar siguiendo la demostraci\'on de la Proposici\'on \ref{P8.1.pea} y
observando que si $\vec\alpha_1,\ldots,\vec\alpha_n\in W_m({\ma F}_p)$,
entonces $\Witt \sum_{i=1}^n \vec\alpha_i \Witt\cdot \vec\xi_i=\big(
\sum_{i=1}^n \alpha_{i1}\xi_{i1},\ldots \big)$, etc. Dejamos los detalles de la
demostraci\'on al lector interesado.
\end{observacion}

Consideremos la ecuaci\'on $\vec y^q\Witt - \vec y=\vec \alpha$ donde
$\vec \alpha\in W_m(k)$. Sea $\vec y_0\in W_m(\bar{k})$ una soluci\'on
de $\vec y^q\Witt - \vec y=\vec\alpha$, donde $\bar{k}$
denota una cerradura algebraica de $k$.
Notemos que si $K=k(\vec y)=k\va ym\big)$
es c\'iclica de grado $p^m$ sobre
$k$ entonces $y_m\notin k\va y{m-1}\big)$ 
pues de lo contrario se tendr\'ia que
$[K:k]=[k\va y{m-1}\big):k]\leq 
p^{m-1}$. As\'i $K=k(y_m)$.

El conjunto de ra\'ices
de $\vec y^q\Witt -\vec y=\vec \alpha$ es el conjunto $\{\vec y_0\Witt +
\vec \mu\}_{\vec \mu\in W_m(\F)}$. Sea $K=k(\vec y_0)$. Entonces, puesto
que $\F\subseteq k$, se tiene que $W_m(\F)\subseteq W_m(k)$ y por tanto
$K/k$ es normal. Puesto que $|W_m(\F)|=q^m=p^{nm}$, todas las 
ra\'ices de $\vec y^q\Witt -\vec y=\vec \alpha$ son diferentes y por tanto
$K/k$ es separable y de Galois.

Sean $G:=\Gal(K/k)$ y $\sigma\in G$. Entonces $\vec y_0$ y $\sigma
(\vec y_0)$ son conjugados por lo que existe $\vec \xi\in W_m(\F)$ tal que
$\sigma (\vec y_0)=\vec y_0\Witt + \vec\xi$. Ponemos $\sigma_{\vec \xi}:=
\sigma$.

\begin{proposicion}\label{P8.3.pea}
Con las notaciones anteriores, se tiene que $\varphi\colon G\lra W_m(\F)$,
dado por $\varphi(\sigma_{\vec \xi})=\vec \xi$ es un monomorfismo de grupos
lo que implica que la extensi\'on es abeliana y $G\subseteq W_m(\F)\cong
\Big({\ma Z}/p^m{\ma Z}\Big)^n$. As\'i, se tiene que
$G\cong{\ma Z}/p^{a_1}{\ma Z}\times\cdots\times {\ma Z}/p^{a_n}{\ma Z}$ con
$m\geq a_1\geq a_2\geq \cdots\geq a_n\geq 0$. $\fin$
\end{proposicion}

Rec\'iprocamente, sea $K=k(z)$ una $p$--extensi\'on abeliana finita de 
exponente $m$ y rango $n$, es decir 
$G\cong{\ma Z}/p^{a_1}{\ma Z}\times\cdots\times {\ma Z}/p^{a_n}{\ma Z}$ con
$m\geq a_1\geq a_2\geq \cdots\geq a_n\geq 1$ y suponemos que 
$\F\subseteq k$ donde $q=p^n$. Sea $K=k(z_1,\ldots,z_n)$ donde
$\Gal(k(z_i)/k)\cong {\ma Z}/p^{a_i}{\ma Z}$, $1\leq i\leq n$. Entonces
existe $\vec y_i$ tal que
$k(z_i)=k(\vec y_i)$ con $\vec y_i^p\Witt -\vec y_i=\vec \alpha_i\in
W_m(k)$, $\vec y_i\in W_m(K)$ donde ponemos $\vec y_i=\big(
\underbrace{0,\ldots,0}_{m-a_i},y_{i,m-a_i+1},\ldots, y_{i,m}\big)$, esto
es, completamos con ceros las entradas para hacer los vectores
de longitud $m$.

Sea $G=\langle\sigma_1,\ldots,\sigma_n\rangle$ con 
$\sigma_j(\vec y_i)=\begin{cases}\vec y_i\Witt +\vec 1&\text{si $i=j$}\\
\vec y_i&\text{si $i\neq j$}\end{cases}$, y $o(\sigma_i)=p^{a_i}$.
Definimos
\[
\vec y:=\vec \xi_1\Witt \cdot \vec y_1\Witt +\cdots \Witt +\vec\xi_n
\Witt \cdot \vec y_n,
\]
donde $\{\vec\xi_1,\ldots,\vec\xi_n\}$ es una base de $W_m(\F)$
sobre $W_m({\ma F}_p)$.
Entonces $\vec y_i\in W_m(K)$, $\vec \xi_i\in W_m(\F)\subseteq W_m(k)$.
Se sigue que $k(\vec y)\subseteq K$. Veamos que $k(\vec y)=K$.

Sea $\sigma\in G$, digamos $\sigma=\sigma_1^{b_1}\cdots \sigma_n^{b_n}$,
$0\leq b_i\leq a_i-1$, $1\leq i\leq n$. Entonces
\begin{align*}
\sigma\vec y&=\sigma\Big(\Witt \sum_{i=1}^n\big(\vec \xi_i\Witt \cdot
\vec y_i\big)\Big)=\Witt \sum_{i=1}^n\sigma\big(\vec\xi_i\Witt \cdot \vec y_i\big)\\
&\underbracket[0pt]{=}_{\substack{\uparrow\\ \vec \xi_i\in W_m(\F)}}
\Witt \sum_{i=1}^n\vec \xi_i\Witt \cdot \sigma(\vec y_i)=\Witt \sum_{i
=1}^n\vec \xi_i\Witt\cdot (\vec y_i\Witt +\vec b_i)\\
&=\Witt \sum_{i=1}^n\vec \xi_i\Witt \cdot \vec y_i\Witt +
\Witt \sum_{i=1}^n \vec b_i\Witt \cdot \vec\xi_i =\vec y\Witt + \vec \xi,
\end{align*}
donde $\vec \xi:=
\Witt\sum_{i=1}^n\vec b_i\Witt\cdot \vec \xi_i\in W_m(\F)$. Por tanto
$\sigma\vec y=\vec y\iff \vec \xi=\vec 0\iff \vec b_1=\cdots=\vec b_n=
\vec 0\iff b_1=\cdots=b_n=0\iff \sigma=\Id$. Se sigue que  $K=k(\vec y)$.

En resumen tenemos:

\begin{teorema}\label{T8.4.pea}
Sea $k$ un campo de caracter\'istica $p>0$ tal que $\F\subseteq k$. Sea
$\vec\alpha\in W_m(k)$. Si $K=k(\vec y_0)$ donde $\vec y_0$ es una
ra\'iz de $\vec y^q\Witt -\vec y=\vec \alpha\in W_m(k)$, entonces
$K/k$ es una
$p$--extensi\'on abeliana de exponente $p^h$ con $h\leq m$ y rango
$l$ con $l\leq n$, donde $q=p^n$. Adem\'as si $G:=\Gal(K/k)$, se
tiene que $G$ es isomorfo de manera natural con un subgrupo del
anillo de Galois $W_m(\F)$.

Rec\'iprocamente, si $K=k(\vec y_0)$ es una $p$--extensi\'on de Galois
de exponente $p^h$ y rango $l$, entonces $\vec y_0$ es ra\'iz de
alguna ecuaci\'on de la forma $\vec y^q\Witt -\vec y=\vec \alpha$ para
alg\'un $\vec \alpha\in W_m(k)$. $\fin$
\end{teorema}

Consideremos ahora $K/k$ una $p$--extensi\'on abeliana finita dada por
$\vec y^q\Witt -\vec y=\vec\alpha\in W_m(k)$ y donde suponemos que
$\F\subseteq k$. Digamos que $G=\Gal(K/k)\cong \prod_{i=1}^n{\ma Z}/
p^{a_i}{\ma Z}$ con $m\geq a_1\geq a_2\geq \cdots\geq a_n\geq 0$. Sea
$\vec \xi\in W_m(\F)$ y sea
\begin{gather}\label{Eq**.pea}
\vec y_{\vec\xi}:=(\vec \xi^{p^{n-1}}\Witt \cdot \vec y^{p^{n-1}})\Witt +
(\vec \xi^{p^{n-2}}\Witt \cdot \vec y^{p^{n-2}})\Witt +\cdots \Witt +
(\vec \xi^{p}\Witt \cdot \vec y^{p})\Witt +(\vec\xi\Witt\cdot \vec y).
\end{gather}
Entonces 
\[
\vec y_{\vec\xi}^p\Witt -\vec y_{\vec\xi}=\vec \xi\Witt\cdot \vec \alpha,
\]
esto es, $k(\vec y_{\vec\xi})/k$ es una extensi\'on c\'iclica de orden $p^h$
con $h\leq m$.

De ahora en adelante supondremos que $G:=\Gal(K/k)=
\langle\sigma_1,\ldots,\sigma_n\rangle \cong \Big(
{\ma Z}/p^m{\ma Z}\Big)^n$ con $o(\sigma_i)=p^m$ para todo
$1\leq i\leq n$. El grupo $G$ tiene $\frac{q^m-q^{m-1}}{p^m-p^{m-1}}$
subgrupos c\'iclicos de orden $p^m$ distintos. En particular debemos
tener que $\alpha_1\neq 0$. Queremos ver que
entre las extensiones $k(\vec y_{\vec\xi})$ est\'an todas las subextensiones
c\'iclicas de orden $p^m$. Con el isomorfismo $G\cong W_m(\F)$, 
consideremos $\sigma_i:=\sigma_{\vec \xi_i}$, $1\leq i\leq n$ donde
$\{\xi_1,\ldots,\xi_n\}$ es una base de $\F$ sobre ${\ma F}_p$ y donde
recordemos que $\vec\xi_i=\{\xi_i\}=(\xi_i,0,\ldots, 0)$.
De hecho notemos que $\{\xi_1,\ldots,\xi_n\}$ 
es una base de $\F$ sobre ${\ma F}_p$ si y solamente si $G=
\langle\sigma_1,\ldots,\sigma_n\rangle$.

Es f\'acil ver que para $\sigma_{\vec\delta}\in G$ se tiene que
$\sigma_{\vec\delta}(\vec y_{\vec\xi}) =
\vec y_{\vec\xi}\Witt +\Witt \sum_{i=0}^{n-1}
(\vec \xi\Witt \cdot \vec \delta)^{p^i}$. Entonces
$\sigma_{\vec\delta}(\vec y_{\vec\xi})=\vec y_{\vec\xi}\iff
\Witt \sum_{i=0}^{n-1} (\vec \xi\Witt \cdot \vec \delta)^{p^i}=\vec 0$.

En general sea $g_{\vec \xi}(\vec\delta):=\Witt \sum_{i=0}^{n-1}
(\vec \xi\Witt \cdot \vec \delta)^{p^i}$. Entonces
$g_{\vec \xi}(\vec \delta)^p\Witt - g_{\vec \xi}(\vec \delta)=\vec 0$,
esto es, $g_{\vec \xi}(\vec \delta)\in W_m({\ma F}_p)$. El
mapeo $g_{\vec\xi}\colon W_m(\F)\lra W_m({\ma F}_p)$ no necesariamente
es suprayectivo y se tiene
\[
\frac{W_m(\F)}{\ker g_{\vec \xi}}\cong \im g_{\vec\xi}\subseteq W_m(
{\ma F}_p),\quad |W_m({\ma F}_p)|=p^m,
\]
por lo que $|\ker g_{\vec \xi}|\geq \frac{|W_m(\F)|}{|W_m({\ma F}_p)|}=
\frac{q^m}{p^m}$.

Ahora bien $k(\vec y_{\vec\xi})$ es el campo fijo de $K$ bajo $\ker g_{\vec \xi}$.

\begin{proposicion}\label{P8.5.pea} Con las notaciones anteriores,
se tiene que $[k(\vec y_{\vec\xi}):k]=p^m\iff \vec
\xi$ es invertible en $W_m(\F)$
lo cual es equivalente a que si $\vec\xi=(\xi_1,\ldots,\xi_m)$, entonces
$\xi_1\neq 0$.
\end{proposicion}

\begin{proof}
Tenemos que en la expresi\'on $\Witt \sum_{i=0}^{n-1}
(\vec \xi\Witt \cdot \vec \delta)^{p^i}$ la primera componente es
$\sum_{i=0}^{n-1}(\xi_1\delta_1)^{p^i}$.
Supongamos que $\vec \xi$ es invertible, es decir, $\xi_1\neq 0$.
Consideremos el mapeo $\psi\colon\F\lra{\ma F}_p$ dado por
$\psi(\delta)=\sum_{i=0}^{n-1}(\delta\xi_1)^{p^i}$el cual es no cero pues
el polinomio $p(x)=(\xi_1x)^{p^{n-1}}+\cdots+(\xi_1x)^p+
(\xi_1 x)=0$ tiene $p^{n-1}$ ra\'ices. Por tanto
si consideramos la extensi\'on $k(y_{\xi_1})/k$
dada por $y_{\xi_1}^p-y_{\xi_1}=\alpha_1\neq 0$,
el grupo que fija a la extensi\'on no es todo $\F$ y en particular
$[k(y_{\xi_1}):k]=p$. Por tanto $[k(\vec y_{\vec \xi}):k]=p^m$.

Rec\'iprocamente, en caso de que $[k(\vec y_{\vec \xi}):k]=p^m$,
necesariamente $[k(y_{\xi_1}):k]=p$ y el argumento es
reversible por lo que $\vec\xi$ es invertible.  $\fin$
\end{proof} 

\begin{corolario}\label{C8.5'.pea}
Las subextensiones c\'iclicas de grado $p^m$ est\'an dadas 
por $k(\vec y_{\vec\xi})$ donde $\vec y_{\vec\xi}$ est\'a dada por
{\rm{(\ref{Eq**.pea})}}, $\vec\xi$ es invertible y se tiene
\begin{gather*}
\vec y^p_{\vec\xi}\Witt -\vec y_{\vec \xi}=\vec\xi\Witt\cdot\vec\alpha.
\tag*{$\fin$}
\end{gather*}
\end{corolario}

En particular, tomando una base $\{\mu_1,\ldots,\mu_n\}$ de $\F$
sobre ${\ma F}_p$, se tiene que $K=k(\vec y)=k(\vec y_{\vec \mu_1},
\ldots, \vec y_{\vec \mu_n})$.

\begin{proposicion}\label{P8.5''.pea}
Sea $K/k$ una extensi\'on de Galois con grupo de Galois isomorfo a
$W_m(\F)$ con $\F\subseteq k$. 
Supongamos $K=k(\vec z_1,\ldots,\vec z_n)$
con $\vec z_i\in W_m(K)$, 
$\Gal(k(\vec z_i)/k)\cong {\ma Z}/p^m{\ma Z}$, $1\leq i\leq n$.
Entonces todas las subextensiones $k\subseteq k(\vec z)\subseteq K$
tales que $\Gal(k(\vec z)/k)\cong {\ma Z}/p^m{\ma Z}$ est\'an dadas por
\[
\vec z=\Witt \sum_{i=1}^n\vec \alpha_i\Witt\cdot \vec z_i,
\]
con $\vec \alpha_i\in W_m({\ma F}_p)$, $1\leq i\leq n$
y alg\'un $\vec \alpha_{i_0}$ invertible.
\end{proposicion}

\begin{proof}
Sea $G:=\Gal(K/k)=\langle\sigma_1,\ldots,\sigma_n\rangle$ de tal forma
que $\sigma_i\vec z_j=\vec z_j\Witt +\vec \delta_{ij}$ con $\vec \delta_{ij}=
\begin{cases} \vec 1&\text{si $i=j$}\\ \vec 0&\text {si $i\neq j$}\end{cases}$.

Sean $\vec \alpha_1,\ldots,\vec \alpha_n\in W_m({\ma F}_p)$ y sea
$\vec z =\Witt \sum_{i=1}^n\vec \alpha_i\Witt\cdot \vec z_i$. Sea
$\wp(\vec z_i)=\vec z_i^p\Witt -\vec z_i=\vec \gamma_i$ con 
$\vec \gamma_i=\big(\gamma_{i1},\ldots,\gamma_{im}\big)$ y 
$\gamma_{i1}\notin \wp(k)$. Entonces
\[
\wp(\vec z)=\Witt \sum_{i=1}^n\vec \alpha_i\Witt\cdot \wp(\vec z_i)
=\Witt \sum_{i=1}^n\vec \alpha_i\Witt\cdot \vec \gamma_i=:\vec \gamma,
\]
con $\gamma_1=\sum_{i=1}^n\alpha_{i1}\gamma_{i1}$. Se tiene que
$[k(\vec z):k]=p^m\iff \gamma_1\notin \wp(k)$.

Ahora supongamos que $\wp(\vec z)=\Witt \sum_{i=1}^n\vec 
\alpha_i\Witt\cdot \vec\gamma_i=\vec \gamma=\wp(\vec A)$ para alg\'un
$\vec A\in W_m(k)$. Entonces $\wp(\vec z\Witt - \vec A)=\vec 0$, esto
es, $\vec z\Witt - \vec A\in W_m({\ma F}_p)$, $\vec z=\vec \beta
\Witt + \vec A$ con $\vec \beta\in W_m({\ma F}_p)$.
En este caso, si existiese $\vec \alpha_{i_0}$ invertible, entonces
\begin{gather*}
\vec z=\Witt \sum_{i=1}^n\vec \alpha_i\Witt\cdot \vec z_i=\vec \beta
\Witt + \vec A
\intertext{y por ende}
\vec z_{i_0}=\Witt - \Witt \sum_{\substack{i=1\\ i\neq i_0}}^n
\vec \alpha_{i_0}^{-1}\Witt \cdot \vec \alpha_i\Witt\cdot \vec z_i\Witt +
\vec \alpha_{i_0}^{-1}\Witt \cdot \vec \beta\Witt + \vec \alpha_{i_0}^{-1}
\Witt \cdot \vec A.
\intertext{Ahora bien, 
puesto que $\vec \beta\in W_m({\ma F}_p)\subseteq W_m(k)$ y
$\vec A\in W_m(k)$, se sigue que $\vec z_{i_0}\in k(\vec z_1,\ldots,
\vec z_{i_0-1},\vec z_{i_0+1},\ldots, \vec z_n)$, que}
K=k(\vec z_1,\ldots,\vec z_n)=k(\vec z_1,\ldots,
\vec z_{i_0-1},\vec z_{i_0+1},\ldots, \vec z_n)
\end{gather*}
y que $[K:k]\leq p^{m(n-1)}<p^{mn}$, lo cual es absurdo.

En resumen, si alg\'un $\vec \alpha_{i_0}$ es invertible, $\vec \gamma=
\wp(\vec z)\notin \wp(W_m(k))$.

Con este procedimiento obtenemos $t$ extensiones $k(\vec z)$ con 
$[k(\vec z):k]=p^m$ y $t=\big|\{(\vec\alpha_1,\ldots,\vec\alpha_n)\mid
\vec\alpha_i\in W_m({\ma F}_p)\text{\ y alg\'un $\vec \alpha_i$ invertible}
\}\big|$.

Se tiene que $\vec\alpha_1,\ldots,\vec\alpha_n\in W_m({\ma F}_p)$ son
no invertibles si y solamente si $\alpha_{11}=\alpha_{21}=\cdots=\alpha_{n1}
=0$ donde $\vec \alpha_i=(\alpha_{i1},\ldots,\alpha_{im})$, $1\leq i\leq n$.
Entonces $t=\big|W_m({\ma F}_p)^n\big|-\big|W_{m-1}({\ma F}_p)^n\big|=
p^{nm}-p^{n(m-1)}=q^m-q^{m-1}$.

Ahora bien, dos de estas extensiones 
$k(\vec z)$, $k(\vec w)$ satisfacen que
$k(\vec z)=k(\vec w)\iff \vec z=\vec j\Witt \cdot \vec w\Witt + \vec c$ con
$\vec j\in W_m({\ma F}_p)$ invertible y $\vec c\in W_m(k)$
(Corolario \ref{C9'.7.11}). Puesto que
$\vec z$ y $\vec w$ son combinaciones ``lineales'' de $\vec z_1,\ldots,
\vec z_n$ sobre $W_m({\ma F}_p)$, $\vec c=\vec 0$ y $\vec j\in W_m(
{\ma F}_p)^{\ast}$. Finalmente
$\big|W_m({\ma F}_p)^{\ast}\big|=\big|W_m({\ma F}_p)\big|-
\big|W_{m-1}({\ma F}_p)^\big|=p^m-p^{m-1}$.

De esta forma hemos obtenido $\frac{q^m-q^{m-1}}{p^m-p^{m-1}}$
extensiones c\'iclicas distintas $k(\vec z)/k$ de grado $p^m$ y por tanto todas.
$\fin$
\end{proof}

Para estudiar la generaci\'on de las extensiones
de este tipo, tenemos el mismo resultado que
el del Teorema \ref{T7.1.pea}. Sea
$K=k(\vec y)$ tal que $\vec y^q\Witt -\vec y=\vec \alpha$ y
donde $\Gal(K/k)\cong W_m(\F)$. Sea
$L=k(\vec z)$ tal que $\vec z^q\Witt-\vec z=\vec \beta$.

Para $\vec A_{n-1},\vec A_{n-2},\ldots, \vec A_1, \vec A_0 \in W_m(\F)$
se define $\mc R (\vec X)\in W_m(\F)[\vec X]$ por 
\[
\mc R(\vec X):=
\vec A_{n-1}\Witt \cdot \vec X^{p^{n-1}}\Witt +
\vec A_{n-2}\Witt \cdot \vec X^{p^{n-2}}\Witt +\cdots\Witt +
\vec A_{1}\Witt \cdot \vec X^{p}\Witt + \vec A_0\Witt \cdot \vec X.
\]

\begin{teorema}\label{T8.6.pea} Con las notaciones anteriores,
$k(\vec y)=k(\vec z)$ si y solamente si (existen $\vec A_{n-1},\vec A_{n-2},\ldots, 
\vec A_1, \vec A_0
\in W_m(\F)$ que satisfacen que $\mc R(\vec \beta)=0$ con $\vec \beta \in 
W_m(\F) \iff \vec \beta=0$) y $\vec D\in W_m(k)$ con 
\begin{equation}\label{Eq8.6.pea}
\vec z=\mc R(\vec y)\Witt +\vec D.
\end{equation}
\end{teorema}

\begin{proof} La demostraci\'on es paralela a la demostraci\'on del
Teorema \ref{T7.2.pea} usando el formalismo de las operaciones de Witt.
Daremos \'unicamente unos pocos detalles. 

Supongamos primero que $k(\vec y)=k(\vec z)$.
Sea $\sigma_i\in G$ dado
por $\sigma_i(\vec y)=\vec y\Witt + \vec \mu_i$, $1\leq i\leq n$. Sea
$\vec w=\mc R(\vec y)$. Entonces $\sigma \in G$ est\'a dado
por $\sigma=\sigma_1^{b_1}\cdots\sigma_n^{b_n}$ con $b_i\in {\ma Z}$,
$0\leq b_i\leq p^m-1$, $1\leq i\leq n$. Se tiene
$\sigma(\vec w)=\vec w\Witt + \mc R\big(\Witt \sum_{i=1}^n \vec b_i
\Witt \cdot \vec \mu_i\big)$. En particular $\sigma_i(\vec w)=\vec w\Witt +
\mc R(\vec \mu_i)$.

Si $\sigma_i(\vec z)=\vec z\Witt +\vec \xi_i$, $1\leq i\leq n$, entonces
$\{\vec\xi_1,\ldots,\vec\xi_n\}$ es una base de $W_m(\F)$ sobre
$W_m({\ma F}_p)$.
Se quieren hallar $\vec A_{0}, \ldots,\vec A_{n-1}\in W_m(\F)$ tales que
$\mc R(\vec\mu_i)=\vec \xi_i$, $1\leq i\leq n$.

Se tiene
\begin{gather*}
\mc R(\vec\mu_i)=\vec \xi_i, 1\leq i\leq n 
\iff \vec M\Witt \cdot\left[\begin{array}{c}\vec A_0\\ \vec A_1\\ \vdots\\ \vec A_{n-2}\\
\vec A_{n-1}\end{array}\right]=\left[\begin{array}{c}\vec \xi_1\\ \vec \xi_2\\ \vdots\\ 
\vec \xi_{n-1}\\ \vec\xi_n\end{array}\right]
\end{gather*}
donde $\vec M$ es la matriz
\[
\vec M=\vmatriz\mu.
\]
Ahora bien, es claro que $\det \vec M=(\det M,\ldots )$ donde 
\[
M=\matriz{\mu_1}{\mu}.
\]
Se tiene que $\det M\in {\ma F}_q^{\ast}$ y por tanto $\det \vec M$ es una unidad
de $W_m(\F)$ de donde se sigue que $\vec M$ es invertible y por
tanto los $\vec A_i\in W_m(\F)$ existen y son \'unicos satisfaciendo
que $\sigma(\vec w)=\vec w\Witt +\vec \xi_i$. El resto de la demostraci\'on
es igual a la del Teorema \ref{T7.2.pea}. $\fin$
\end{proof}

Ahora estudiaremos el caso de campos de funciones racionales. Sea
$k=k_0(T)$ un campo de funciones racionales donde $k_0$ es
un campo finito tal que $\F\subseteq k_0$.
Tenemos el resultado an\'alogo a \cite[Theorem 5.5]{MalRzeVil2012}.

\begin{teorema}\label{T8.7.pea} Sea $K/k$ una extensi\'on tal que
$\Gal(K/k)\cong W_m(\F)$ y donde se tiene que
$P_1,\ldots,P_r\in R_T^+$ y posiblemente $\p$, 
son los primos ramificados.
Entonces $K=k(\vec y)$ est\'a dada por
\[
\vec y^q\Witt -\vec y=\vec \beta={\vec\delta}_1\Witt 
+ \cdots \Witt + {\vec\delta}_r
\Witt + \vec\gamma,
\]
con $y_1^q-y_1=\beta_1$ es irreducible,
$\delta_{ij}=\frac{Q_{ij}}{P_i^{e_{ij}}}$, $e_{ij}\geq 0$, $Q_{ij}\in R_T$
y si $e_{ij}>0$, entonces $e_{ij}=
\lambda_{ij}p^{m_{ij}}$, $\mcd (\lambda_{ij},p)=1$,
$0\leq m_{ij}< n$, $\mcd(Q_{ij},P_i)=1$ y 
$\deg (Q_{ij})<\deg (P_i^{e_{ij}})$, y $\gamma_j=f_j(T)\in R_T$ con
$\deg f_j=\nu_j p^{m_j}$ con $\mcd(p,\nu_j)=1$, $0\leq m_j<n$
cuando $f_j\not\in k_0$.
\end{teorema}

\begin{proof}
Para la primera reducci\'on de separar los polinomios primos 
en el denominador se procede como
en el caso de Teorema \ref{T5.3.1.pea}. Una vez obtenida esta
simplificaci\'on se procede como Schmid (Observaci\'on \ref{O9'.8.28}),
como en la demostraci\'on del Teorema \ref{T5.2.pea}
usando el Corolario \ref{C8.5'.pea}. $\fin$
\end{proof}

\begin{observacion}\label{O8.8'.pea}
La descomposici\'on de $\p$ en una $p$--extensi\'on c\'iclica,
est\'a dada por la Proposici\'on \ref{P5.3.2} y en particular $\p$
se descompone totalmente en $K/k$ si y solamente si
$\vec\gamma=\vec 0$.
Otra forma equivalente es decir que $\p$
se descompone totalmente en $K/k$ si y solamente si existe $\vec \theta
\in W_m(k)$ tal que $\vec\gamma =\vec 
\theta^p\Witt - \vec\theta=\wp(\vec\theta)$.
\end{observacion}

De la Observaci\'on \ref{O8.8'.pea} se obtiene:

\begin{teorema}\label{P8.9.pea}
Sea $K/k$ como en el Teorema {\rm{\ref{T8.7.pea}}}. 
Si $\vec\gamma=\vec 0$, entonces $\p$ se descompone totalmente.

Rec\'iprocamente, si $\p$ se descompone totalmente,
entonces existe una descomposici\'on como en el Teorema
{\rm{\ref{T8.7.pea}}} con $\vec\gamma=\vec 0$.
\end{teorema}

\begin{proof}
La prueba es similar a la de la Proposici\'on \ref{P6.2.pea}. $\fin$
\end{proof}

%% file: Capitulo13.tex
\chapter{El teorema de Kronecker--Weber en 
caracter{\'\i}stica $p$\index{Kronecker--Weber!teorema de $\sim$ en caracter\'istica $p$}}\label{Ch11}

\section{Introducci\'on}\label{S11.1}

El teorema cl\'asico de Kronecker--Weber establece que toda
extensi\'on finita ${\ma Q}$ est\'a contenida en un campo ciclot\'omico,
ver Cap\'itulo \ref{Ch1}.
Equivalentemente, la m\'axima extensi\'on abeliana de ${\ma Q}$
es la uni\'on de todos los campos ciclot\'omicos. En 1974
D. Hayes \cite{Hay74}, prob\'o el resultado an\'alogo para
campos de funciones racionales congruentes. 
Tenemos que la uni\'on de todos los campos de funciones
ciclot\'omicos no es la m\'axima extensi\'on abeliana
del campo de funciones racionales congruente
$K={\ma F}_q(T)$ puesto que todas estas extensiones 
son geom\'etricas y el primo infinito es moderadamente
ramificado. Hayes prob\'o que la m\'axima extensi\'on 
abeliana de $K$ es la composici\'on de la uni\'on de todos
los campos de funciones ciclot\'omicos con la uni\'on de todas
las extensiones de constantes y con la uni\'on  de los
subcampos del campo de funciones ciclot\'omico
para el primo infinito donde el primo infinito es total y
salvajemente ramificado. La demostraci\'on de Hayes
usa teor{\'\i}a de campos de clase.

La demostraci\'on del caso cl\'asico la dimos 
en el Cap{\'\i}tulo \ref{Ch1} la cual usa grupos de ramificaci\'on.
La herramienta clave en la demostraci\'on es que 
hay una \'unica extensi\'on c{\'\i}clica de ${\ma Q}$
de grado $p$ ($p$ impar), y $p$ es el \'unico primo
ramificado. En el caso de campos de funciones racionales la
situaci\'on es bastante diferente. Existen muchas
extensiones c{\'\i}clicas de $K$ de grado $p$
donde \'unicamente un divisor primo es ramificado.

En este
cap{\'\i}tulo se presenta una demostraci\'on
del an\'alogo al Teorema de Kronecker--Weber\index{Kronecker--Weber!teorema
de $\sim$}\index{teorema de Kronecker--Weber}
\index{Kronecker--Weber!teorema de $\sim$}
para campos de funciones racionales congruentes usando
argumentos de conteo en el caso de ramificaci\'on salvaje.
Primero, como en el caso cl\'asico, mostramos que
cualquier extensi\'on c{\'\i}clica de $K$ est\'a contenida
en la composici\'on de un campo de funciones
ciclot\'omicos y una extensi\'on de constantes.
El siguiente paso, el fundamental, es mostrar que
toda extensi\'on c{\'\i}clica de grado una potencia
de $p$ donde \'unicamente hay un primo ramificado
y \'este es completamente ramificado, est\'a contenida
en un campo de funciones ciclot\'omica.
Una vez que esto est\'e probado, el resto de la prueba
se sigue f\'acilmente. Usamos la aritm\'etica
de vectores de Witt desarrollada por Schmid en \cite{Sch36}
(ver Cap{\'\i}tulo \ref{Ch9'}).

\section{El Teorema de Kronecker--Weber para campos de 
funciones\index{Kronecker--Weber!teorema de $\sim$ para
campos de funciones}}\label{S11.2}

Para esta secci\'on establecemos nuevamente la notaci\'on
que usaremos.
Sea $K_T:=\bigcup_{M\in R_T} K(\Lambda_M)$, ${\ma F}_{\infty}:=
\bigcup_{m\in{\ma N}}{\ma F}_{q^m}$.
Denotamos por ${\eu p}_{\infty}$
el divisor de polos de $T$ en $K$. 
Denotamos por $L_n$ al m\'aximo subcampo
de $K(\Lambda_{1/T^n})$ donde $\p$
es total y salvajemente ramificado, $n\in{\ma N}$. 
Sea $L_{\infty}:=\bigcup_{n\in{\ma N}}L_n$.

El principal objetivo de este cap{\'\i}tulo es probar el siguiente
resultado.

\begin{teorema}[Kronecker--Weber, {\cite{Hay74}},
{\cite[Theorem12.8.31]{Vil2006}}]\label{T11.2.1}
La m\'axima extensi\'on abeliana 
 $A$ de $K$ es $A=K_T
{\ma F}_{\infty} L_{\infty}$.
\end{teorema}

Para probar el Teorema \ref{T11.2.1} es suficiente probar que
toda extensi\'on abeliana finita de $K$ est\'a contenida en
$K(\Lambda_N) {\ma F}_{q^m} L_n$ para algunos $N\in R_T$, y
$m,n\in{\ma N}$.

Sea $L/K$ una extensi\'on abeliana finita. Sea $G:=\Gal(L/K)\cong
C_{n_1}\times \cdots\times C_{n_l}\times C_{p^{a_1}}\times
\cdots \times C_{p^{a_h}}$ donde $\mcd (n_i,p)=1$, $1\leq
i\leq l$ y $a_j\in{\ma N}$, $1\leq j\leq h$. 
Sea $S_i\subseteq
L$ tal que $\Gal(S_i/K)\cong C_{n_i}$, $1\leq i\leq l$ y sea
$R_j\subseteq L$ tal que $\Gal(R_j/K)\cong C_{p^{a_j}}$,
$1\leq j\leq h$. Para probar el Teorema \ref{T11.2.1} es suficiente
mostrar que cada $S_i$ y cada $R_j$ est\'an contenidos
en $K(\Lambda_N) {\ma F}_{q^m} L_n$ para algunas $N\in R_T,
m,n\in{\ma N}$.

En resumen, podemos suponer que
$L/K$ es una extensi\'on c{\'\i}clica de grado $h$ donde
ya sea $\mcd(h,p)=1$ o $h=p^n$ para alg\'un
$n\in{\ma N}$.

\subsection{Extensiones geom\'etricas moderadamente
ramificadas}\label{S11.3}

En esta subsecci\'on probaremos el Teorema \ref{T11.2.1}
para el caso particular de una extensi\'on moderadamente
ramificada.
Sea $L/K$ una extensi\'on abeliana.
Sea $P\in R_T$, $d:=\deg P$.

Por la Proposici\'on \ref{P3.1}, se tiene que si 
$P$ es moderadamente ramificado en 
$L/K$, entonces si $e$ denota el {\'\i}ndice de ramificaci\'on de
$P$ en $L$, tenemos $e\mid q^d-1$.

Ahora consideremos una extensi\'on abeliana finita 
moderadamente ramificada $L/K$
donde $P_1,\ldots,P_r$ donde son los primos finitos ramificados.
Sean $P\in\{P_1,\ldots,P_r\}$ y $e$ 
el {\'\i}ndice de ramificaci\'on de
$P$ in $L$. Entonces, como consecuencia de la 
Proposici\'on \ref{P3.1}, tenemos que
$e\mid q^d-1$. Ahora bien, $P$ es totalmente
ramificado en $K(\Lambda_P)/K$ con {\'\i}ndice de ramificaci\'on
$q^d-1$. En esta extensi\'on $\p$ tiene {\'\i}ndice de 
ramificaci\'on igual a $q-1$.

Sea $K\subseteq E\subseteq K(\Lambda_P)$ con $[E:K]=e$.
Pongamos $\tilde{{\eu P}}$ un divisor primo en $LE$ que divide a $P$.
Sean ${\eu q}:= \tilde{\eu P}|_E$ y ${\eu P}:=\tilde{\eu P}|_L$.
\[
\xymatrix{
{\eu P}\ar@/^1pc/@{-}[rrrr]\ar@/_1pc/@{-}[dddd]\ar@{--}[dr]
&&&&\tilde{\eu P}
\ar@/^1pc/@{-}[dddd]|!{[dddl];[dddr]}\hole\ar@{--}[dl]\\
&L\ar@{--}[rr]\ar@{-}[dd]&&LE\ar@{--}[dd]\\
&&M\ar@{-}[dl]\ar@{-}[ur]_H\\
&K\ar@{-}[rr]_e&&E\ar@{-}[rr] \ar@{--}[dr]&&K(\Lambda_P)\\
P\ar@/_1pc/@{-}[rrrr]\ar@{--}[ru]&&&&{\eu q}
}
\]

Tenemos que
$e=e_{L/K}({\eu P}|P)=e_{E/K}({\eu q}|P)$. Como consecuencia
del Lema de Abhyankar \cite[Theorem 12.4.4]{Vil2006}, 
se obtiene que
\[
e_{LE/K}(\tilde{\eu P}|P)=\mcm [e_{L/K}({\eu P}|P), e_{E/K}(
{\eu q}|P)]=\mcm[e,e]=e.
\]
Sea $H\subseteq \Gal(LE/K)$ el grupo de inercia de $\tilde{
\eu P}/P$. Pongamos $M:=(LE)^H$. Entonces $P$ 
es no ramificado en la extensi\'on $M/K$.
Queremos probar que $L\subseteq MK(\Lambda_P)$. 
De hecho se tiene que
$[LE:M]=e$ y $E\cap M=K$ puesto que $P$ 
es totalmente ramificado en
$E/K$ y no ramificado en $M/K$. Se sigue que
$[ME:K]=[M:K][E:K]$. Luego
\[
[LE:K]=[LE:M][M:K]=e\frac{[ME:K]}{[E:K]}=
e\frac{[ME:K]}{e}=[ME:K].
\]
Puesto que $ME\subseteq LE$ se sigue que
$LE=ME\subseteq MK(\Lambda_P)$. Por tanto
$L\subseteq MK(\Lambda_P)$.

En $M/K$ los primos finitos ramificados son
$\{P_2,\cdots, P_r\}$. En caso de que $r-1\geq 1$
podemos aplicar el argumento anterior a
$M/K$ obteniendo de esta forma una extensi\'on
$M_2/K$ de tal manera que a lo m\'as $r-2$ primos finitos
de $K$ son ramificados en $M_2 K$ y se tiene que
$M\subseteq M_2K(\Lambda_{P_2})$,
por lo que $L\subseteq MK(\Lambda_{P_1})\subseteq
M_2K(\Lambda_{P_1})K(\Lambda_{P_2})$. 

Llevando a cabo el proceso anterior a lo m\'as
$r$ veces, obtenemos
\begin{equation}\label{E3.2}
L\subseteq M_0K(\Lambda_{P_1P_2\cdots P_r})
\end{equation}
en donde en la extensi\'on $M_0/K$ 
el \'unico posible primo ramificado es $\p$.

Notemos la similitud con la Proposici\'on \ref{P4.1}.

Por la Proposici\'on \ref{P11.3.2}, se sigue que $M_0/K$
es una extensi\'on de constantes.

Como corolario a (\ref{E3.2}) y a la Proposici\'on \ref{P11.3.2}
obtenemos el Teorema \ref{T11.2.1} para el caso moderadamente
ramificado.

\begin{proposicion}\label{C3.5}
Si $L/K$ es una extensi\'on finita abeliana moderadamente ramificada
donde los divisores primos finitos ramificados son
$P_1,\ldots,P_r$, entonces
\[
L\subseteq {\ma F}_{q^m} K(\Lambda_{P_1\cdots P_r}).
\]
para alguna $m\in{\ma N}$. $\fin$
\end{proposicion}

\subsection{Extensiones salvajemente ramificadas}\label{S11.4}

\subsubsection{Reducciones}\label{S11.4'}

Como consecuencia de la Proposici\'on \ref{C3.5}, el Teorema
\ref{T11.2.1} se seguir\'a si probamos el caso particular de una
extensi\'on c{\'\i}clica $\K/K$ de grado $p^n$ para
alguna $n\in{\ma N}$. Ahora, este tipo de extensiones est\'an
dadas por medio de un vector de Witt:
\[
\K=K(\vec y)=K(y_1,\ldots,y_n) \quad\text{con}\quad
 \vec y^p\Witt - \vec y =
\vec \beta=(\beta_1,\ldots,\beta_n)\in W_n(K).
\]

Primero consideremos una extensi\'on de Artin--Schreier.
Sea $\K:=K(y)$ donde $y^p-y=\alpha\in K$. 
La ecuaci\'on puede ser normalizada como sigue:
\begin{equation}\label{Eq1}
y^p-y=\alpha=\sum_{i=1}^r\frac{Q_i}{P_i^{e_i}} + f(T),
\end{equation}
donde $P_i\in R_T^+$, $Q_i\in R_T$, 
$\mcd(P_i,Q_i)=1$, $e_i>0$, $p\nmid e_i$, $\deg Q_i<
\deg P_i^{e_i}$, $1\leq i\leq r$, $f(T)\in R_T$,
con $p\nmid \deg f$ cuando $f(T)\not\in {\ma F}_q$.

Tenemos que los primos finitos ramificados en $\K/K$ son
precisamente $P_1,\ldots,P_r$. Con respecto a $\p$ 
la Proposici\'on \ref{P2.4} establece que $\p$ es
\l
\item descompuesto si $f(T)=0$.
\item inerte si $f(T)\in {\ma F}_q$ y $f(T)\not\in \wp({\ma F}_q):=
\{a^p-a\mid a\in {\ma F}_q\}$.
\item ramificado si $f(T)\not\in {\ma F}_q$ (por tanto $p\nmid\deg f$).
\end{list}

El Teorema \ref{T5.3.1} nos prueba que en una extensi\'on c{\'\i}clica
de grado $p^n$ dada por un vector de Witt, podemos separar cada
primo ramificado.

El comportamiento de $\p$ en $\K/K$ est\'a establecido en 
la Proposici\'on \ref{P5.3.2}.

Consideremos el campo $\K=K(\vec y)$ como en el
Teorema \ref{T5.3.1}, donde \'unicamente un divisor
prime $P\in R_T^+$ se ramifica, con
\begin{gather}\label{EqNew1}
\begin{array}{l}
\text{$\beta_i=
\frac{Q_i}{P^{\lambda_i}}$, $Q_i\in R_T$
tal que $\lambda_i\geq 0$},\\
\text{si $\lambda_i=0$ entonces $Q_i=0$,}\\
\text{si $\lambda_i>0$ entonces $\mcm(\lambda_i,p)=1$,
$\mcd(Q_i,P)=1$ y $\deg Q_i<\deg P^{\lambda_i}$},\\
\lambda_1>0.
\end{array}
\end{gather}

Un caso particular del Teorema \ref{T5.3.1} adecuado para
nuestro estudio se da en la siguiente proposici\'on.

\begin{proposicion}\label{P2.6} Supongamos que toda extensi\'on
$\K_1/K$ que cumpla con las condiciones de {\rm (\ref{EqNew1})} 
satisface que $\K_1
\subseteq K(\Lambda_{P^{\alpha}})$ para alguna
$\alpha\in{\ma N}$. Sea $\K/K$ la extensi\'on
definida por $\K=K(\vec{y})$ donde
 $\wp(\vec{y})=  \vec y^p\Witt - \vec y=
\vec \beta$ con $\vec \beta=(\beta_1,\ldots,\beta_n)$,
$\beta_i$ dado en forma normal:
$\beta_i\in{\ma F}_q$
o $\beta_i=\frac{Q_i}{P^{\lambda_i}}$, 
$Q_i\in R_T$ y $\lambda_i>0$,
$\mcd(\lambda_i,p)=1$, $\mcd(Q_i,P)=1$ y $\deg Q_i\leq
\deg P^{\lambda_i}$. Entonces $\K\subseteq {\ma F}_{q^{p^n}}
K(\Lambda_{P^{\alpha}})$ para alguna $\alpha\in{\ma F}_q$.
\end{proposicion}

\begin{proof} Del Teorema \ref{T5.3.1} se tiene que podemos
descomponer el vector $\vec \beta$ como
$\vec \beta=\vec \varepsilon\Witt + \vec \gamma$ con
$\varepsilon_i\in{\ma F}_q$ para toda $1\leq i\leq n$ y
$\gamma_i=0$
o $\gamma_i=\frac{Q_i}{P^{\lambda_i}}$,
$Q_i\in R_T$ y $\lambda_i>0$,
$\mcd(\lambda_i,p)=1$, $\mcd(Q_i,P)=1$ y $\deg Q_i<
\deg P^{\lambda_i}$. 

Sea
$\gamma_1=\cdots=\gamma_r=0$, y
$\gamma_{r+1}\notin {\ma F}_q$.
Tenemos $\K\subseteq K(\vec \varepsilon)K(\vec \gamma)$.
Ahora $K(\vec \varepsilon)\subseteq {\ma F}_{q^{p^n}}$ y
$K(\vec \gamma)=K(0,\ldots,0,\gamma_{r+1},\ldots,\gamma_n)$.

Para cualquier vector de Witt $\vec x=(x_1,\ldots,x_n)$ tenemos
la descomposici\'on dada por la Proposici\'on \ref{P9'.3.3}
\begin{align*}
\vec x=&(x_1,0,0,\ldots,0)\Witt +(0,x_2,0,\ldots,0)\Witt +\cdots
\Witt + (0,\ldots, 0,x_j,0,\ldots,0)\\
&\Witt +(0,\ldots, 0, x_{j+1}, \ldots, x_n)
\end{align*}
para cada $0\leq j\leq n-1$. Se sigue que
$K(\vec \gamma)=K(\gamma_{r+1},\ldots,\gamma_n)$.
Puesto que este campo satisface las condiciones de (\ref{EqNew1}),
tenemos $K(\vec \gamma)\subseteq K(\Lambda_{P^{\alpha}})$
para alguna $\alpha\in{\ma N}$. El resultado se sigue. $\fin$
\end{proof}

\begin{observacion}\label{R2.6'}
El primo $\p$ puede ser manejado de la misma manera que cualquier
$P\in R_T^+$. Las condiciones
(\ref{EqNew1}) para $\p$ son las siguientes.
Sea $\K=K(\vec \mu)$ con
$\mu_j=f_j(T)\in R_T$, adem\'as $f_j(0)=0$ para toda $j$ y
ya sea $f_j(T)=0$ o $f_j(T)\neq 0$ y $p\nmid \deg f_j(T)$.
La condici\'on $f_j(0)=0$ significa que el primo infinito
para $T^{\prime}=1/T$ es o descompuesto o ramificado
en cada nivel, esto es, el grado de inercia es
$1$ en $\K/K$.
En este caso, con el cambio de variable $T^{\prime}=1/T$
la hip\'otesis en la Proposici\'on \ref{P2.6} debe decir que
cualquier campo que cumpla estas condiciones
satisface que
$\K\subseteq K(\Lambda_{T^{\prime m}})=
K(\Lambda_{T^{-m}})$ para alguna $m\in{\ma N}$.
Sin embargo, puesto que el grado de la extensi\'on
$\K/K$ es una potencia de $p$ necesariamente
tenemos que $\K$ est\'a contenida en
$K(\Lambda_{T^{-m}})^{
{\ma F}_q^{\ast}}=L_{m-1}$.
\end{observacion}

Con las notaciones del Teorema \ref{T5.3.1}
obtenemos que si $\vec z_i^p\Witt - \vec z_i
=\vec \delta_i$, $1\leq i\leq r$ y si $\vec v^p\Witt -\vec v=\vec \mu$,
entonces $\K=K(\vec y)\subseteq K(\vec z_1,\ldots,\vec z_r,\vec v)=
K(\vec z_1)\ldots K(\vec z_r) K(\vec v)$. Por lo tanto, si el
Teorema \ref{T11.2.1} se cumple para cada $K(\vec z_i)$,
$1\leq i\leq r$ y para $K(\vec v)$, entonces se cumple para $\K$.

Del Teorema \ref{T5.3.1}, de la Proposici\'on \ref{P2.6} y
de la observaci\'on despu\'es de esta proposici\'on, obtenemos que
para probar el Teorema \ref{T11.2.1}, es suficiente mostrar que
cualquier extensi\'on $\K/K$ que cumpla las condiciones
de (\ref{EqNew1}) satisface que ya sea
$\K\subseteq K(\Lambda_{P^{\alpha}})$
para alguna $\alpha\in{\ma N}$ o $\K\subseteq L_m$
para alguna $m\in{\ma N}$.
Es suficiente estudiar el caso $P\in R_T^+$.

De las Proposiciones \ref{P2.4} y \ref{P5.3.2} obtenemos

\begin{proposicion}\label{P2.5'}
Si $\K$ es un campo definido por una ecuaci\'on del tipo
dado en {\rm (\ref{EqNew1})}, entonces
$\K/K$ es una extensi\'on c{\'\i}clica de grado
$p^n$, $P$ es el \'unico primo ramificado, es
totalmente ramificado y $\p$ es totalmente descompuesto.

Similarmente, si $\K=K(\vec v)$
donde $v_i=f_i(T)\in R_T$, $f_i(0)=0$
para toda $1\leq i\leq n$ y $f_1(T)\notin {\ma F}_q$,
$p\nmid \deg f_1(T)$, entonces $\p$ es el \'unico
primo ramificado en $\K/K$, es totalmente
ramificado y el divisor de ceros de $T$, el cual es
ahora el primo al infinito en $R_{1/T}$, es totalmente
descompuesto. $\fin$
\end{proposicion}

Hemos reducido la demostraci\'on del Teorema \ref{T11.2.1}
a probar que cualquier extensi\'on del tipo dado en
la Proposici\'on \ref{P2.5'} est\'a contenido ya sea 
en $K(\Lambda_{P^{\alpha}})$ para alguna 
$\alpha\in{\ma N}$ o en $L_m$ para alguna
$m\in{\ma N}$. El segundo caso es consecuencia del
primero con el cambio de variable $T^{\prime}=1/T$.

Sea $n,\alpha\in{\ma N}$. Denotemos por $v_n(\alpha)$ 
al n\'umero de grupos c{\'\i}clicos de orden $p^n$ contenidos
en $\G {P^{\alpha}}\cong
\Gal(K(\Lambda_{P^{\alpha}})/K)$. Tenemos que
$v_n(\alpha)$ es el n\'umero de extensiones c{\'\i}clicas
$\K/K$ de grado $p^n$ y $\K
\subseteq K(\Lambda_{P^{\alpha}})$.
Toda extensi\'on de este tipo satisface que su
conductor ${\eu F}_\K$ divide a $P^{\alpha}$.

Sea ahora $t_n(\alpha)$ el n\'umero de extensiones
de campos $\K/K$ de grado $p^n$ tal que
$P$ es el \'unico primo ramificado, es totalmente
ramificado, $\p$ es totalmente descompuesto y
su conductor ${\eu F}_\K$ es un divisor de $P^{\alpha}$.
Puesto que toda extensi\'on c{\'\i}clica $\K/K$
de grado $p^n$ tal que $K\subseteq 
\K\subseteq
K(\Lambda_{P^{\alpha}})$ satisface estas
condiciones, tenemos que $v_n(\alpha)\leq t_n(\alpha)$.
Si probamos que $t_n(\alpha)\leq v_n(\alpha)$ entonces
toda extensi\'on satisfaciendo la ecuaci\'on (\ref{EqNew1})
est\'a contenida en una extensi\'on ciclot\'omica y por ende
se sigue el Teorema \ref{T11.2.1}.

En resumen, para probar el Teorema \ref{T11.2.1}, 
es suficiente probar 
\begin{gather}\label{EqNew2}
t_n(\alpha)\leq v_n(\alpha)\quad\text{para toda}\quad
n, \alpha\in{\ma N}.
\end{gather}

\subsubsection{Demostraci\'on de (\ref{EqNew2})}\label{S11.4N}

Probaremos ahora por inducci\'on en $n$ la
ecuaci\'on (\ref{EqNew2})
y como consecuencia obtenemos el Teorema \ref{T11.2.1}.
Primero calcularemos $v_n(\alpha)$ para todas
$n,\alpha\in{\ma N}$.

\begin{proposicion}\label{P11.4.1}
El n\'umero $v_n(\alpha)$ de grupos c{\'\i}clicos de
orden $p^n$ contenidos en $\G {P^{\alpha}}$ es
\begin{gather*}
v_n(\alpha)=\frac{q^{d(\alpha-\mayorchico{\alpha}{p^n})}-q^{d(\alpha-\mayorchico{\alpha}
{p^{n-1}})}}{p^{n-1}(p-1)}=
\frac{q^{d(\alpha-\mayorchico{\alpha}
{p^{n-1}})}\big(q^{d(\mayorchico{\alpha}{p^{n-1}}
-\mayorchico{\alpha}{p^{n}})}-1\big)}{p^{n-1}(p-1)},
\end{gather*}
donde $\lceil x\rceil$ denota la {\em funci\'on techo\index{funci\'on techo}},
esto es, $\lceil x\rceil$ denota al entero m\'as peque\~no
tal que es mayor o igual a $x$.
\end{proposicion}

\begin{proof}
Sea $P\in R_T^+$ y $\alpha\in{\ma N}$ con $\deg P=d$.
Primero calculamos cuantas extensiones c{\'\i}clicas
de grado $p^n$ est\'an contenidas en
$\lama {\alpha}$. Puesto que $\p$ es moderadamente
ramificada en $\lama \alpha$, si $\K/K$ es
una extensi\'on c{\'\i}clica de grado $p^n$,
$\p$ se descompone totalmente en $\K/K$
(Teorema \ref{T6.3.5}). Tenemos
$\Gal(\lama \alpha/K)\cong \G{P^{\alpha}}$ y
la sucesi\'on exacta
\begin{gather}\label{E4.1}
0\longrightarrow D_{P,P^{\alpha}}\longrightarrow
\G {P^{\alpha}}\stackrel{\varphi}{\longrightarrow} \G P\longrightarrow 0,
\end{gather}
donde
\begin{eqnarray*}
\varphi\colon \G {P^{\alpha}}&\longrightarrow&\G P\\
A\bmod P^{\alpha}&\longmapsto& A\bmod P
\end{eqnarray*}
y $D_{P,P^{\alpha}}=\ker \varphi=\{N\bmod P^{\alpha}\mid
N\equiv 1\bmod P\}$. Podemos considerar, sin peligro 
alguno, que $D_{P,P^{\alpha}}=\{1+hP\mid h\in R_T, \deg h<\deg P^{\alpha}=
d\alpha\}$.

Tenemos que $\G {P^{\alpha}}\cong \G P \times D_{P,P^{\alpha}}$
y que $\G P\cong C_{q^d-1}$.
En primer lugar calculamos cuantos elementos de
orden $p^n$ existen en $\G {P^{\alpha}}$. Estos elementos
pertenecen a $D_{P,P^{\alpha}}$. Sea $A=1+hP\in D_{P,P^{\alpha}}$
de orden $p^n$. Escribimos $h=gP^{\gamma}$ con
$g\in R_T$, $\mcd (g,P)=1$ y $\gamma\geq 0$.
Tenemos $A=1+gP^{1+\gamma}$. Puesto que $A$ es
de orden $p^n$, se sigue que
\begin{gather}
A^{p^n}=1+g^{p^n}P^{p^n(1+\gamma)}\equiv 1\bmod P^{\alpha}\label{E4.2}\\
\intertext{y}
A^{p^{n-1}}=1+g^{p^{n-1}}P^{p^{n-1}(1+\gamma)}
\not\equiv 1\bmod P^{\alpha}.\label{E4.3}
\end{gather}

De (\ref{E4.2}) y (\ref{E4.3}) se obtiene
\begin{gather}\label{E4.4}
p^{n-1}(1+\gamma)<\alpha\leq p^n(1+\gamma)\\
\intertext{y (\ref{E4.4}) es equivalente a}
\mayor{\alpha}{p^n}-1\leq \gamma <\mayor{\alpha}
{p^{n-1}}-1.
\end{gather}
Notemos que para la existencia de al menos un elemento
 de orden $p^n$ es necesario que $\alpha>p^{n-1}$.

Ahora para cada $\gamma$ que satisface (\ref{E4.4}) tenemos
$\mcd (g,P)=1$ y $\deg g + d(1+\gamma) <d \alpha$,
esto es, $\deg g<d(\alpha-\gamma-1)$. Luego, existen
$\Phi(P^{\alpha-\gamma-1})$ tales $g$'s, donde
para cualquier
 $N\in R_T$, $\Phi(N):=\big|\G N\big|$.
 
 Por lo tanto el n\'umero de elementos de orden $p^n$ 
en $D_{P,P^{\alpha}}$ es 
\begin{gather}\label{E4.5}
\sum_{\gamma=\mayorchico{\alpha}{p^n}-1}^{
\mayorchico{\alpha}{p^{n-1}}-2}\Phi(P^{\alpha-\gamma-1})=
\sum_{\gamma'=\alpha-\mayorchico{\alpha}{p^{n-1}}+1}^{
\alpha-\mayorchico{\alpha}{p^n}}\Phi(P^{\gamma'}).
\end{gather}

Notemos que para cualquier $1\leq r\leq s$ tenemos
\begin{align*}
\sum_{i=r}^{s}\Phi(P^i)&=\sum_{i=r}^s q^{d(i-1)}(q^d-1)=
(q^d-1)q^{d(r-1)}\sum_{j=0}^{s-r}q^{dj}\\
&=(q^d-1)q^{d(r-1)}\frac{q^{d(s-r+1)}-1}{q^d-1}=q^{ds}-q^{d(r-1)}.
\end{align*}
Luego, (\ref{E4.5}) es igual a
\[
q^{d(\alpha-\mayorchico{\alpha}{p^n})}-q^{d(\alpha-\mayorchico{\alpha}
{p^{n-1}})}=q^{d(\alpha-\mayorchico{\alpha}
{p^{n-1}})}\big(q^{d(\mayorchico{\alpha}{p^{n-1}}
-\mayorchico{\alpha}{p^{n}})}-1\big).
\]

Puesto que cada grupo c{\'\i}clico de orden
$p^n$ tiene $\varphi(p^n)=
p^{n-1}(p-1)$ generadores, se sigue el resultado. $\fin$
\end{proof}

Notemos que si $\K/K$ es cualquier campo contenido en
$K(\Lambda_{P^{\alpha}})$ entonces es conductor
${\eu F}_\K$ de $\K$ es un divisor de $P^{\alpha}$
(Lema \ref{L9'.8.24}).

Ahora calculamos el n\'umero de extensiones
c{\'\i}clicas $\K/K$ de grado $p$ tales que
$P$ es el \'unico primo ramificado (es completamente ramificado),
$\p$ se descompone en $\K/K$ y
${\eu F}_\K\mid P^{\alpha}$.

\begin{proposition}\label{P11.4.3}
Toda extensi\'on c{\'\i}clica
$\K/K$ de grado $p$ tal que $P$ es el \'unico primo ramificado,
$\p$ se descompone en $\K/K$ y ${\eu F}_\K\mid P^{\alpha}$
est\'a contenido en $\lama \alpha$.
\end{proposition}

\begin{proof}
De la Teor{\'\i}a de Artin--Schreier 
(ver (\ref{Eq1})) y Proposici\'on \ref{P2.4}
se tiene que $\K$ satisface $\K=K(y)$ con la ecuaci\'on
de Artin--Schreier de $y$ normalizada como
prescrita por Hasse (\cite{Has35} y 
Observaci\'on \ref{O9'.8.28}). Por lo tanto
\[
y^p-y=\frac{Q}{P^{\lambda}},
\]
donde $P\in R_T^+$, $Q\in R_T$, 
$\mcd(P,Q)=1$, $\lambda >0$, $p\nmid \lambda$, $\deg Q<
\deg P^{\lambda}$. Ahora el conductor ${\eu F}_\K$ satisface
 ${\eu F}_\K=P^{\lambda+1}$ as{\'\i} que
$\lambda\leq \alpha-1$.

Tenemos que si $\K=K(z)$ con $z^p-z=a$ entonces existen
$j\in {\ma F}_p^{\ast}$ y $c\in K$ tales que
$z=jy+c$ y $a=j\frac{Q}{P^{\lambda}}+\wp(c)$ donde
$\wp (c):=c^p-c$. Si $a$ est\'a dado en forma normal
entonces $c=\frac{h}{P^{\gamma}}$ con $p\gamma\leq \lambda$
(de hecho, $p\gamma<\lambda$ puesto que $\mcd (\lambda,p)=1$)
y $\deg h<\deg P^{\gamma}$ o $h=0$. Sea $\gamma_0:=
\integerchico{\alpha-1}{p}$. Entonces cualquier tal $c$ puede ser
escrito como $c=\frac{hP^{\gamma_0-\gamma}}{P^{\gamma_0}}$.
Por lo tanto $c\in {\mathcal G}:=\Big\{\frac{h}{P^{\gamma_0}}\mid
h\in R_T, \deg h<\deg P^{\gamma_0}=d\gamma_0
\text{\ o\ } h=0\Big\}$.

Si $c\in{\mathcal G}$ y $j\in\{1,2,\ldots,p-1\}$ tenemos
\begin{align*}
a&=j\frac{Q}{P^{\lambda}}+\wp(c)=j\frac{Q}{P^{\lambda}}+
\frac{h^p}{P^{p\gamma_0}}+\frac{h}{P^{\gamma_0}}\\
&=\frac{jQ+P^{\lambda-p\gamma_0}h+P^{\lambda-\gamma_0}h}
{P^{\lambda}}=\frac{Q_1}{P^{\lambda}},
\end{align*}
con $\deg Q_1<\deg P^{\lambda}$.
Puesto que $\lambda-p\gamma_0>0$ y $\lambda-\gamma_0>0$,
tenemos $\mcd(Q_1,P)=1$. Por lo tanto $a$ est\'a
en forma normal.

Se sigue que el mismo campo tiene $|{\ma F}_p^{\ast}||\wp({\mathcal
G})|$ representaciones diferentes, todas ellas en forma normal.
Ahora ${\mathcal G}$ y $\wp({\mathcal G})$ son grupos
aditivos y $\wp\colon {\mathcal G}\to \wp({\mathcal G})$ es un
epimorfismo de grupos con n\'ucleo 
\[
\ker \wp = {\mathcal G}\cap
\{c\mid \wp(c)=c^p-c=0\}={\mathcal G}\cap {\ma F}_p=\{0\}.
\]
Tenemos $|\wp({\mathcal G})|=|{\mathcal G}|=|R_T/(P^{\gamma_0})|=
q^{d\gamma_0}$.

De la discusi\'on anterior se desprende que el n\'umero de 
extensiones c{\'\i}clicas diferentes $\K/K$ de
grado $p$ tales que el conductor
de $\K$ es ${\eu F}_\K= P^{\lambda+1}$ es igual a
\begin{equation}\label{E4.8}
\frac{\Phi(P^{\lambda})}{|{\ma F}_p^{\ast}||\wp({\mathcal G})|}=
\frac{q^{d(\lambda-1)}(q^d-1)}{(p-1)q^{d\gamma_0}}=
\frac{q^{d(\lambda-\integerchico{\lambda}{p}-1)}(q^d-1)}{p-1}=
\frac{1}{p-1}\Phi\big(P^{\lambda-\integerchico{\lambda}{p}}\big).
\end{equation}

Por lo tanto, el n\'umero de extensiones c{\'\i}clicas
$\K/K$ de grado $p$ tal que el conductor de
$\K$ es ${\eu F}_\K$ es un divisor de $P^{\alpha}$ est\'a dada por
$\frac{w(\alpha)}{p-1}$ donde
\begin{equation}\label{E3.10}
w(\alpha)=\sum_{\substack{\lambda=1\\ \mcd(\lambda,p)=1}}^{\alpha-1}
\Phi\big(P^{\lambda-\integerchico{\lambda}{p}}\big).
\end{equation}

Para calcular $w(\alpha)$ escribamos $\alpha-1=pt_0+r_0$
con $t_0\geq 0$ y
$0\leq r_0\leq p-1$. Ahora $\{\lambda\mid 1\leq \lambda \leq 
\alpha-1, \mcm(\lambda,p)=1\}={\mathcal A}\cup {\mathcal B}$ donde
\begin{gather*}
{\mathcal A}
=\{pt+r\mid 0\leq t\leq t_0-1, 1\leq r\leq p-1\}
\quad \text{y}\quad {\mathcal B}=
\{pt_0+r\mid 1\leq r\leq r_0\}.
\intertext{Entonces}
w(\alpha)=\sum_{\lambda\in{\mathcal A}}
\Phi\big(P^{\lambda-\integerchico{\lambda}{p}}\big)+\sum_{\lambda\in
{\mathcal B}}\Phi\big(P^{\lambda-\integerchico{\lambda}{p}}\big)
\end{gather*}
donde entendemos que si un conjunto,
${\mathcal A}$ o ${\mathcal B}$ es vac{\'\i}o,
la suma respectiva es $0$.

Entonces
\begin{align}\label{E4.8'}
w(\alpha)&=\sum_{\substack{0\leq t\leq t_0-1\\ 1\leq r\leq p-1}}
q^{d(pt+r-t-1)}(q^d-1)+\sum_{r=1}^{r_0}(q^{d(pt_0+r-t_0-1})(q^d-1)\nonumber\\
&=(q^d-1)\Big(\sum_{t=0}^{t_0-1}q^{d(p-1)t}\Big)\Big(\sum_{
r=1}^{p-1}q^{d(r-1)}\Big)+
(q^d-1)q^{d(p-1)t_0}\sum_{r=1}^{r_0}q^{d(r-1)}\\
&=(q^d-1)\frac{q^{d(p-1)t_0}-1}{q^{d(p-1)}-1}\frac{q^{d(p-1)}-1}{q^d-1}
+(q^d-1)q^{d(p-1)t_0}\frac{q^{dr_0}-1}{q^d-1}\nonumber\\
&=q^{d((p-1)t_0+r_0}-1=q^{d(pt_0+r_0-t_0}-1=q^{d(\alpha-1-\integerchico
{\alpha-1}{p})}-1.\nonumber
\end{align}

Entonces, el n\'umero de extensiones c{\'\i}clicas
$\K/K$ de grado $p$ tales que
 $P$ es el \'unico primo ramificado, ${\eu F}_\K\mid P^{\alpha}$ y
$\p$ se descompone, es
\begin{equation}\label{E4.7}
\frac{w(\alpha)}{p-1}=\frac{q^{d(\alpha-1-\integerchico
{\alpha-1}{p})}-1}{p-1}.
\end{equation}

Para finalizar la demostraci\'on de la Proposici\'on 
\ref{P11.4.3} necesitamos el siguiente

\begin{lemma}\label{L4.2}
Para cualquier $\alpha\in{\ma Z}$ y $s\in {\ma N}$ se tiene
\l
\item $\integer{\integerchico{\alpha}{p^s}}{p}=\integer{\alpha}{p^{s+1}}$.

\item $\mayor{\alpha}{p^s}=\integer{\alpha-1}{p^s}+1$.
\end{list}
\end{lemma}

\begin{proof}
Para ({\sc i}), notemos que el caso $s=0$ es claro. Pongamos $\alpha=tp^{s+1}+r$
con $0\leq r\leq p^{s+1}-1$. Sea $r=lp^s+r'$ con $0\leq r'\leq p^s-1$. Notemos
que $0\leq l\leq p-1$. Por lo tanto $\alpha=tp^{s+1}+lp^s+r'$, $0\leq r'\leq p^s-1$
y $0\leq l\leq p-1$. Luego $\integer{\alpha}{p^s}=tp+l$, y
$\d\frac{\integer{\alpha}{p^s}}{p}=t+\frac{l}{p}$, $0\leq l\leq p-1$. As{\'\i}
$\integer{\integerchico{\alpha}{p^s}}{p}=t=\integer{\alpha}{p^{s+1}}$.

Para ({\sc ii}), escribamos $\alpha=p^s t+r$ con
$0\leq r\leq p^s-1$. Si $p^s\mid \alpha$ entonces $r=0$ y $\mayor{\alpha}{p^s}
=t$, $\integer{\alpha-1}{p^s}=\integer{p^s t-1}{p^s}=\d\Big[t-\frac{1}{p^s}\Big]=
t-1=\mayor{\alpha}{p^s}-1$.

Si $p^s\nmid \alpha$, entonces $1\leq r\leq p^s-1$ y $\alpha-1=p^st+(r-1)$
con $0\leq r-1\leq p^s-2$. Por tanto $\mayor{\alpha}{p^n}=\d\Big\lceil q+\frac{r}{p^s}
\Big\rceil=t+1$ y $\integer{\alpha-1}{p^s}=\d\Big[t+\frac{r-1}{p^s}\Big]=t=
\mayor{\alpha}{p^s}-1$. $\fin$
\end{proof}

Del Lema \ref{L4.2} ({\sc i}) obtenemos que (\ref{E4.7}) es igual a
\begin{equation}\label{E3.11}
\frac{w(\alpha)}{p-1}=\frac{q^{d(\alpha-1-(\mayorchico
{\alpha}{p}-1))}-1}{p-1}=\frac{q^{d(\alpha-\mayorchico
{\alpha}{p})}-1}{p-1}=v_1(\alpha).
\end{equation}

Como consecuencia de (\ref{E3.11}), tenemos
la Proposici\'on \ref{P11.4.3}. $\fin$
\end{proof}

La Proposici\'on \ref{P11.4.3} prueba (\ref{EqNew2})
para $n=1$ y toda $\alpha \in{\ma N}$.

Ahora consideremos una extensi\'on c{\'\i}clica
$\K_n/K$ de grado $p^n$ tal que
$P$ es el \'unico primo ramificado, es completamente ramificado, $\p$
se descompone completamente en $\K_n/K$ y
${\eu F}_\K\mid P^{\alpha}$. Queremos probar que
$\K_n\subseteq \lama \alpha$. Esto ser\'a probado por
inducci\'on en $n$. El caso
$n=1$ es la Proposici\'on \ref{P11.4.3}. Suponemos que toda
extensi\'on c{\'\i}clica $\K_{n-1}$ de grado $p^{n-1}$, $n\geq 2$
tal que $P$ es el \'unico primo ramificado,
$\p$ se descompone completamente en
$\K_{n-1}/K$ y tal que ${\eu F}_{\K_{n-1}}
\mid P^{\delta}$
est\'a contenida en $\lama \delta$ donde $\delta\in{\ma N}$.

Sea $\K_n$ cualquier extensi\'on c{\'\i}clica de grado $p^n$ tal que
$P$ es el \'unico primo ramificado y es totalmente ramificado,
 $\p$ se descompone totalmente en $\K_n/K$
y ${\eu F}_{\K_n}\mid P^{\alpha}$.
Sea $\K_{n-1}$ el subcampo de $\K_n$ de grado $p^{n-1}$.
Ahora consideremos $\K_n/K$ generado por el 
vector de Witt $\vec{\beta}=
(\beta_1,\ldots,\beta_{n-1},\beta_n)$, esto es, $\wp(\vec{y})=
\vec{y}^p\Witt - \vec{y}=\vec{\beta}$,
y suponemos que $\vec{\beta}$ est\'a en su forma normal
descrita por Schmid
(ver Observaci\'on \ref{O9'.8.28}, \cite{Sch36-0,Sch36}).
Entonces $\K_{n-1}/K$ est\'a dado por el vector de Witt
$\vec{\beta'}=(\beta_1,\ldots,\beta_{n-1})$.

Sea $\vec{\lambda}:=(\lambda_1,\ldots,\lambda_{n-1},\lambda_n)$
el vector de los par\'ametros de Schmid, esto es, donde cada
$\beta_i$ est\'a dado por
\begin{gather*}
\beta_i=\frac{Q_i}{P^{\lambda_i}}, \text{\ donde\ } Q_i=0 \text{\ 
(esto es, $\beta_i=0$) y $\lambda_i:=0$ o}\\
 \mcd(Q_i,P)=1, \deg Q_i<\deg P^{\lambda_i}, \lambda_i>0 \text{\ y\ }
\mcd(\lambda_i,p)=1.
\end{gather*}
Ahora, puesto que $P$ es totalmente ramificado, se tiene $\lambda_1>0$.

Ahora calculamos cuantas extensiones $\K_n/\K_{n-1}$
diferentes pueden ser construidas por medio de
$\beta_n$. 

\begin{lema}\label{LN1}
Para un campo fijo $\K_{n-1}$, el n\'umero de campos diferentes
$\K_n$ es menor o igual a
\begin{equation}\label{EqNew5}
\frac{1+w(\alpha)}{p}=
\frac{1}{p}q^{d(\alpha-\mayorchico {\alpha}p)}.
\end{equation}
\end{lema}

\begin{proof}
Para $\beta_n\neq 0$, cada ecuaci\'on en forma normal
est\'a dada por
\begin{gather}\label{E4.11}
y_n^p-y_n=z_{n-1}+\beta_n,
\end{gather}
donde $z_{n-1}$ es el elemento en $\K_{n-1}$ obtenido por
la generaci\'on de Witt de $\K_{n-1}$ del vector $\vec{\beta'}$ (ver
(\ref{Ec9'.2.1}) y Teorema \ref{T9'.7.12}). De hecho 
$z_{n-1}$ est\'a dado, formalmente, por
\[
z_{n-1}=\sum_{i=1}^{n-1}\frac{1}{p^{n-i}}\big[
y_i^{p^{n-i}}+\beta_i^{p^{n-i}}-(y_i+\beta_i+
z_{i-1})^{p^{n-i}}\big],
\]
con $z_0=0$.

Como en el caso $n=1$ tenemos que existen
$\Phi(P^{\lambda_n})$ extensiones para los diferentes
$\beta_n$ con
$\lambda_n>0$. El n\'umero de elementos $\beta_n$ diferentes 
que nos dan el mismo campo $\K_n$ con el cambio
$y_n\to y_n+c$, $c\in {\mathcal G}_{\lambda_n}
:=\big\{\frac{h}{P^{\gamma_0}}\mid
h\in R_T, \deg h<\deg P^{\gamma_0}=d\gamma_0 \text{\ o\ }
h=0\big\}$ donde $\gamma_n=\integer{\lambda_n}{p}$, obtenemos
$\beta_n\to \beta_n+\wp(c)$ est\'a tambi\'en en forma normal.
Por tanto, el n\'umero de elementos $\beta_n$ que nos
dan el mismo campo $\K_n$ con este cambio de
variable es $
q^{d(\integerchico{\lambda_n}{p})}$. Por lo tanto obtenemos
a lo m\'as $\Phi\big(P^{\lambda_n-\integerchico{\lambda_n}{p}}\big)$ 
posibles campos $\K_n$ para cada $\lambda_n>0$ (ver (\ref{E4.7})).
M\'as precisamente, si para cada $\beta_n$ con $\lambda_n>0$ definimos
$\overline{\beta_n}:=\{\beta_n+\wp(c)\mid c\in {\mathcal G}_{\lambda_n}\}$,
entonces cada elemento de $\overline{\beta_n}$ 
nos da el mismo campo $\K_n$.

Sea $v_P$ la valuaci\'on en $P$ y
\begin{gather*}
{\mathcal A}_{\lambda_n}:=\{\overline{\beta_n}\mid v_P(\beta_n)=-
\lambda_n\},\\
{\mathcal A}:=\bigcup_{\substack{\lambda_n=1\\
\mcd(\lambda_n,p)=1}}^{\alpha-1} {\mathcal A}_{\lambda_n}.
\end{gather*}

Entonces cada campo $\K_n$ est\'a dado por
$\beta_n=0$ o $\overline{\beta_n}\in
{\mathcal A}$.  De (\ref{E4.8'})
Tenemos que el n\'umero de campos $\K_n$
conteniendo un campo fijo $\K_{n-1}$ que
obtuvimos en (\ref{E4.11}) es menor o igual a
\begin{gather}\label{E410'}
1+|{\mathcal A}|=1+w(\alpha)=q^{d\big(\alpha-1-\integerchico
{\alpha-1}{p}\big)}=q^{d\big(\alpha-1-\mayorchico{\alpha}{p}+1\big)}=
q^{d\big(\alpha-\mayorchico{\alpha}{p}\big)}.
\end{gather}

Ahora con la substituci\'on $y_n\to y_n+jy_1$, $j=0,1,\ldots,p-1$,
en (\ref{E4.11}) obtenemos
\[
(y_n+jy_1)^p-(y_n+jy_1)=y_n^p-y_n +j(y_1^p-y_1)=
z_{n-1}+\beta_n+j\beta_1.
\]

Por lo tanto, cada una de las extensiones obtenidas en
(\ref{E4.11}) se repite al menos $p$ veces, esto es, para
cada $\beta_n$, obtenemos la misma extensi\'on con
$\beta_n,\beta_n+\beta_1,\ldots, \beta_n+(p-1)\beta_1$.
Probaremos que diferentes $\beta_n+j\beta_1$ corresponden
a elementos diferentes de $\{0\}\cup {\mathcal A}$.

Fijemos $\beta_n$. Modificamos cada  $\beta_n+j\beta_1$ en su
forma normal: $\beta_n+j\beta_1 +\wp(c_{\beta_n,j})$ para alguna
$c_{\beta_n,j}\in K$. De hecho $\beta_n+j\beta_1$ est\'a siempre
en forma normal con la posible excepci\'on de $\lambda_n=
\lambda_1$ y a\'un en este caso esto sucede para a lo m\'as
un {\'\i}ndice $j\in
\{0,1,\ldots,p-1\}$: si $\lambda_n\neq \lambda_1$, 
\[
v_P(\beta_n+j\beta_1)=
\begin{cases}
-\lambda_n&\text{si $j=0$}\\
-\max\{-\lambda_n,-\lambda_1\}&\text{si $j\neq 0$}
\end{cases}.
\]
Cuando $\lambda_n= \lambda_1$ y si $v_P(\lambda_n+j
\lambda_1)=u>-\lambda_n=-\lambda_1$ y $p|u$,
entonces para  $i\neq j$, $v_P(\beta_n+i\beta_i)=v_P(
\beta_n+j\beta_1+(i-j)\beta_1)=-\lambda_n=-\lambda_1$.
En otras palabras $c_{\beta_n,j}=0$ con muy pocas
excepciones.

Cada $\mu=\beta_n+j\beta_1+\wp(c_{\beta_n,j})$, $j=0,1,\ldots,p-1$
satisface que o bien $\mu=0$ o $\overline{\mu}\in{\mathcal A}$. 
Veremos que todos estos elementos nos dan elementos
diferentes de $\{0\}\cup{\mathcal A}$.

Si $\beta_n=0$, entonces para $j\neq 0$, $v_P(j\beta_1)=-\lambda_1$, 
de tal forma que $\overline{j\beta_1}\in{\mathcal A}$. Ahora si
$\overline{j\beta_1}=\overline{i\beta_1}$,
entonces
\[
j\beta_1=\beta_n^{\prime}+\wp(c_1)\quad\text{y}\quad
i\beta_1=\beta_n^{\prime}+\wp(c_2)
\]
para alguna $\beta_n^{\prime}\neq 0$ y algunas
$c_1,c_2\in{\mathcal G}_{\lambda_1}$. Se sigue que
$(j-i)\beta_1=\wp(c_2-c_1)\in\wp(K)$. Esto no es
posible por la elecci\'on de $\beta_1$ a menos que $j=i$.

Sea $\beta_n\neq 0$. El caso $\beta_n+j\beta_1=0$
para alguna $j\in\{0,1,\ldots,p-1\}$ ha sido ya considerada en
en el primer caso. Por tanto consideramos el caso
$\beta_n+j\beta_1
+\wp(c_{\beta_n,j})\neq 0$ para toda $j$. Si para algunas
 $i,j\in\{0,1,\ldots,p-1\}$
tenemos $\overline{\beta_n+j\beta_1+\wp(c_{\beta_n,j})}=\overline{
\beta_n+i\beta_1+\wp(c_{\beta_n,i})}$ entonces existen $\beta_n^{\prime}$
y $c_1,c_2\in K$ tales que 
\[
\beta_n+j\beta_1+\wp(c_{\beta_n,j})=\beta_n^{\prime}+\wp(c_1)\quad
\text{y}\quad
\beta_n+i\beta_1+\wp(c_{\beta_n,i})=\beta_n^{\prime}+\wp(c_2).
\]
Se sigue que $(j-i)\beta_1=\wp(c_1-c_2+
c_{\beta_n,i}-c_{\beta_n,j})\in\wp(K)$ de tal forma que
$i=j$.

Por tanto cada campo $\K_n$ es representado por al menos
$p$ elementos diferentes de $\{0\}\cup {\mathcal A}$. 
El resultado se sigue. $\fin$
\end{proof}

Ahora bien, de acuerdo con Schmid (Teorema \ref{T9'.8.29}), el conductor
de $\K_n$ es $P^{M_n+1}$ donde
$M_n=\max\{pM_{n-1},\lambda_n\}$ y 
$P^{M_{n-1}+1}$ es el conductor de $\K_{n-1}$.
Puesto que ${\eu F}_{\K_n}\mid
P^{\alpha}$, se tiene $M_n\leq \alpha-1$. Por lo tanto
$pM_{n-1}\leq \alpha-1$ y $\lambda_n\leq \alpha-1$. Luego
${\eu F}_{\K_{n-1}}\mid P^{\delta}$ con $\delta=\integer{\alpha-1}{p}+1$.

\begin{proposicion}\label{P4.4}
Se tiene
\begin{gather*}
\frac{v_n(\alpha)}{v_{n-1}(\delta)}=
\frac{q^{d(\alpha-\mayorchico{\alpha}{p})}}{p},
\end{gather*}
donde $\delta=\integer{\alpha-1}{p}+1$.
\end{proposicion}

\begin{proof}
De la Proposici\'on \ref{P11.4.1} obtenemos 
\begin{align*}
v_n(\alpha)&=\frac{q^{d(\alpha-\mayorchico{\alpha}
{p^{n-1}})}\big(q^{d(\mayorchico{\alpha}{p^{n-1}}
-\mayorchico{\alpha}{p^{n}})}-1\big)}{p^{n-1}(p-1)}\\
&=\frac{q^{d(\alpha-\mayorchico{\alpha}
{p^{n-1}})}}{p^{n-1}(p-1)}\big(q^{d(\mayorchico{\alpha}{p^{n-1}}
-\mayorchico{\alpha}{p^{n}})}-1\big),\\
\intertext{y}
v_{n-1}(\delta)&=\frac{q^{d(\delta-\mayorchico{\delta}
{p^{n-2}})}\big(q^{d(\mayorchico{\delta}{p^{n-2}}
-\mayorchico{\delta}{p^{n-1}})}-1\big)}{p^{n-2}(p-1)}\\
&=\frac{q^{d(\delta-\mayorchico{\delta}
{p^{n-2}})}}{p^{n-2}(p-1)}\big(q^{d(\mayorchico{\delta}{p^{n-2}}
-\mayorchico{\delta}{p^{n-1}})}-1\big).
\end{align*}

Ahora del Lema \ref{L4.2} obtenemos
\begin{align*}
\mayor{\delta}{p^{n-2}}-\mayor{\delta}{p^{n-1}}&=\Big(
\integer{\delta-1}{p^{n-2}}+1\Big)-\Big(\integer{\delta-1}{p^{n-1}}+1\Big)\\
&=\integer{\delta-1}{p^{n-2}}-\integer{\delta-1}{p^{n-1}}=
\integer{\integerchico{\alpha-1}{p}}{p^{n-2}}-
\integer{\integerchico{\alpha-1}{p}}{p^{n-1}}\\
&=\integer{\alpha-1}{p^{n-1}}-\integer{\alpha-1}{p^n}=
\Big(\mayor{\alpha}{p^{n-1}}-1\Big)-\Big(\mayor{\alpha}{p^{n}}-1\Big)\\
&=\mayor{\alpha}{p^{n-1}}-\mayor{\alpha}{p^n},\\
\delta-\mayor{\delta}{p^{n-2}}&=\Big(\integer{\alpha-1}{p}+1\Big)
-\Big(\integer{\delta-1}{p^{n-2}}+1\Big)\\
&=\integer{\alpha-1}{p}-\integer{\delta-1}{p^{n-2}}=\integer{\alpha-1}{p}
-\integer{\integerchico{\alpha-1}{p}}{p^{n-2}}\\
&=\integer{\alpha-1}{p}-\integer{\alpha-1}{p^{n-1}}.
\end{align*}

Por tanto
\begin{gather*}
v_{n-1}(\delta)=\frac{q^{d(\integerchico{\alpha-1}{p}-\integerchico{
\alpha-1}{p^{n-1}})}}{p^{n-2}(p-1)}\Big(q^{d(\mayorchico{\alpha}{p^{n-1}}
-\mayorchico{\alpha}{p^n})}-1\Big).
\end{gather*}
As{\'\i}, por el Lema \ref{L4.2},
\begin{align*}
\frac{v_n(\alpha)}{v_{n-1}(\delta)}&=
\frac{\frac{q^{d(\alpha-\mayorchico{\alpha}
{p^{n-1}})}}{p^{n-1}(p-1)}\big(q^{d(\mayorchico{\alpha}{p^{n-1}}
-\mayorchico{\alpha}{p^{n}})}-1\big)}
{\frac{q^{d(\integerchico{\alpha-1}{p}-\integerchico{
\alpha-1}{p^{n-1}})}}{p^{n-2}(p-1)}\Big(q^{d(\mayorchico{\alpha}{p^{n-1}}
-\mayorchico{\alpha}{p^n})}-1\Big)}\\
&=\frac{1}{p}q^{d(\alpha-\mayorchico{\alpha}{p^{n-1}}-\integerchico{
\alpha-1}{p}+\integerchico{\alpha-1}{p^{n-1}})}\\
&=\frac{1}{p}q^{d(\alpha-\mayorchico{\alpha}{p^{n-1}}-(\mayorchico{
\alpha}{p}-1)+(\mayorchico{\alpha}{p^{n-1}}-1))}\\
&=\frac{1}{p}q^{d(\alpha-\mayorchico{\alpha}{p})}.
\end{align*}

Esto prueba el resultado. $\fin$
\end{proof}

De aqu{\'\i}, de la Proposici\'on \ref{P4.4}, Lema
\ref{LN1} (\ref{EqNew5}) y puesto que por hip\'otesis de inducci\'on
tenemos
$t_{n-1}(\delta)=v_{n-1}(\delta)$,
obtenemos
\[
t_n(\alpha)\leq t_{n-1}(\delta)\big(\frac{1}{p}q^{d\big(\alpha-\mayorchico
{\alpha}p\big)}\big)=v_{n-1}(\delta)\big(\frac{1}{p}q^{d\big(\alpha-\mayorchico
{\alpha}p\big)}\big)=v_n(\alpha).
\]
Esto prueba (\ref{EqNew2}) y el Teorema \ref{T11.2.1}.

\subsubsection{Prueba alternativa de (\ref{EqNew2})}\label{S6}

Mantenemos la misma notaci\'on de las subsecciones previas.
Sea $\K/K$ una extensi\'on satisfaciendo las
condiciones de (\ref{EqNew1}) y con el conductor un divisor de
$P^{\alpha}$. Tenemos que 
${\eu F}_\K=P^{M_n+1}$ donde
\begin{gather*}
M_n=\max\{p^{n-1}\lambda_1,p^{n-2}\lambda_2,\ldots,
p\lambda_{n-1},\lambda_n\}\\
\intertext{ver Teorema \ref{T9'.8.29}. Por lo tanto}
{\eu F}_\K\mid P^{\alpha}\iff M_n+1\leq \alpha \iff
p^{n-i}\lambda_i\leq \alpha-1, \quad i=1,\ldots, n.
\end{gather*}

Entonces $\lambda_i\leq \integer{\alpha-1}{p^{n-i}}$. Estas
condiciones proporcionan todas las extensiones c{\'\i}clicas
de grado $p^n$ donde $P\in R_T^+$ es el \'unico primo
ramificado, es totalmente ramificado, $\p$ se decompone
totalmente y su conductor divide a $P^{\alpha}$. Ahora
estimamos el n\'umero de diferentes formas para generar
$\K$.

Sea $\K=K(\vec y)$. Primero notemos que con el cambio
de variable $y_i$ por $y_i+c_i$ para cada $i$,
$c_i\in K$ obtenemos el mismo campo. Para estas 
nuevas formas de generar $\K$ que a su vez
cumpla (\ref{EqNew1}), debemos tener:
\l
\item Si $\lambda_i=0$, $c_i=0$.
\item Si $\lambda_i>0$, entonces $c_i
\in \Big\{\frac{h}{P^{\gamma_i}}\mid
h\in R_T, \deg h<\deg P^{\gamma_i}=d\gamma_i
\text{\ o\ } h=0\Big\}$,
donde $\gamma_i=\integer{\lambda_i}{p}$.
Por lo tanto tenemos a lo m\'as
$\Phi\big(P^{\lambda_i-\integerchico{\lambda_i}{p}}\big)$ 
extensiones para este $\lambda_i$ 
(ver (\ref{E4.8})). Puesto que $1\leq \lambda_i
\leq\integer{\alpha-1}{p^{n-i}}$ y $\mcd(\lambda_i,p)=1$,
si definimos $\delta_i:=\integer{\alpha-1}{p^{n-i}}+1$, 
de (\ref{E3.10}) y de (\ref{E4.8'}) obtenemos que debemos
tener a lo m\'as
\begin{gather}\label{EqNew3}
w(\delta_i)=\sum_{\substack{\lambda_i=1\\ 
\mcm(\lambda_i,p)=1}}^{\delta_i-1}
\Phi\big(P^{\lambda_i-\integerchico{\lambda_i}{p}}\big)=
q^{d\big(\delta_i-1-\integerchico{\delta_i-1}{p}\big)}-1
\end{gather}
expresiones diferentes para todos los posibles
$\lambda_i>0$.

Ahora del Lema \ref{L4.2} tenemos
\[
\delta_i-1-\integer{\delta_i-1}{p}=\integer{\alpha-1}{p^{n-i}}-
\integer{\integer{\alpha-1}{p^{n-i}}}{p}=\integer{\alpha-1}
{p^{n-i}}-\integer{\alpha-1}{p^{n-i+1}}.
\]
Por tanto
\begin{gather}\label{EqNew4}
w(\delta_i)=q^{d\big(\integerchico{\alpha-1}{p^{n-i}}-\integerchico{
\alpha-1}{p^{n-i+1}}\big)}-1.
\end{gather}
\end{list}

Cuando $\lambda_i=0$ es permitido tenemos a lo m\'as $w(\delta_i)+1$
extensiones con par\'ametro $\lambda_i$. Por tanto, puesto que $\lambda_1
>0$ y $\lambda_i\geq 0$ para $i=2,\ldots, n$, tenemos que el
n\'umero de extensiones satisfaciendo
(\ref{EqNew1}) y con conductor un divisor de 
$P^{\alpha}$ es a lo m\'as
\begin{gather*}
s_n(\alpha):=w(\delta_1)\cdot \prod_{i=2}^n\big(w(\delta_i)+1\big).\\
\intertext{De (\ref{EqNew3}) y de (\ref{EqNew4}), obtenemos}
s_n(\alpha)=\Big(q^{d\big(\integerchico{\alpha-1}{p^{n-1}}
-\integerchico{\alpha-1}{p^n}\big)}-1\Big)\cdot \prod_{i=2}^n
q^{d\big(\integerchico{\alpha-1}{p^{n-i}}
-\integerchico{\alpha-1}{p^{n-i+1}}\big)}.\\
\intertext{Por lo tanto $\prod_{i=2}^n(w(\delta_i)+1)=q^{d\mu}$ donde}
\begin{align*}
\mu &= \sum_{i=2}^n \Big(\integer{\alpha-1}{p^{n-i}}
-\integer{\alpha-1}{p^{n-i+1}}\Big)=
\sum_{i=2}^n \integer{\alpha-1}{p^{n-i}}
-\sum_{j=1}^{n-1}\integer{\alpha-1}{p^{n-j}}\\
&=\integer{\alpha-1}{p^{n-n}}-\integer{\alpha-1}{p^{n-1}}=
\alpha-1-\integer{\alpha-1}{p^{n-1}}.
\end{align*}
\end{gather*}

Se sigue que
\begin{align*}
s_n(\alpha)&=\Big(q^{d\big(\integerchico{\alpha-1}{p^{n-1}}
-\integerchico{\alpha-1}{p^n}\big)}-1\Big)\cdot q^{d\big(\alpha
-1-\integerchico{\alpha-1}{p^{n-1}}\big)}\\
&= q^{d\big(\integerchico{\alpha-1}{p^{n-1}}-
\integerchico{\alpha-1}{p^{n}}+\alpha-1-
\integerchico{\alpha-1}{p^{n-1}}\big)}-
q^{d\big(\alpha-1-\integerchico{\alpha-1}{p^{n-1}}\big)}\\
&=q^{d\big(\alpha-1-\integerchico{\alpha-1}{p^{n}}\big)}
-q^{d\big(\alpha-1-\integerchico{\alpha-1}{p^{n-1}}\big)}.
\end{align*}

Del Lema \ref{L4.2} ({\sc ii}) obtenemos
\begin{align*}
\alpha-1-\integer{\alpha-1}{p^n}=\alpha-\mayor{\alpha}{p^n}
\quad\text{y}\quad \alpha-1-\integer{\alpha-1}{p^{n-1}}=
\alpha-\mayor{\alpha}{p^{n-1}}.\\
\intertext{Por tanto}
s_n(\alpha)=q^{\big(\alpha-\mayorchico{\alpha}{p^n}\big)}
-q^{\big(\alpha-\mayorchico{\alpha}{p^{n-1}}\big)}=
p^{n-1}(p-1)v_n(\alpha).
\end{align*}

Finalmente, el cambio de variable $\vec y\to \vec j \Witt \times \vec y$ con
$\vec j\in W_n({\ma F}_p)^{\ast}\cong \big({\ma Z}/p^n{\ma Z}
\big)^{\ast}$ da el mismo campo y tenemos
$\vec \beta\to \vec j\Witt \times \vec \beta$. Por tanto
\[
t_n(\alpha)\leq \frac{s_n(\alpha)}{\varphi(p^n)}=\frac{s_n(\alpha)}
{p^n(p-1)}=v_n(\alpha).
\]
Esto prueba (\ref{EqNew2}) y el Teorema \ref{T11.2.1}.

%% file: Capitulo14.tex
\chapter{Campos de g\'eneros}\label{Ch12*}

\section{Introducci\'on}\label{S12*.1}

Dado un campo num\'erico $K$, esto es, $[K:{\ma Q}]<\infty$, el {\em campo
de clase de Hilbert} $K_H$ de $K$ se define como la m\'axima extensi\'on
abeliana de $K$ no ramificada. Esto significa que ning\'un primo, finito o
infinito, de $K$ se ramifica en $K_H$. Por teor\'ia de campos de clase, 
sabemos que $K_H/K$ es una extensi\'on finita y que $\Gal(K_H/K)
\cong I_K$, el grupo de clases de $K$. Adem\'as los primos totalmente 
descompuestos de $K$ en $K_H$ son precisamente los ideales primos
principales (Proposici\'on \ref{CClaseP4.7.10}).
 
Similarmente, se puede definir el {\em campo de clase de Hilbert extendido}
$K_H^+$ como la m\'axima extensi\'on abeliana de $K$ tal que los primos
finitos no sean ramificados. Se tiene que $K_H^+/K$ es un extensi\'on finita
y $\Gal(K_H^+/K)\cong I_K^+:=D_K/P_K^+$ donde $P_K^+=\{(\alpha)\mid
\alpha\in K$ y $(\alpha)$ est\'a generado por un elemento totalmente
positivo$\}$ (un elemento $\beta$ se llama {\em totalmente
positivo\index{elemento totalmente positivo}} si para todo
encaje real, $K\stackrel{\varphi}{\hooklongrightarrow} {\ma R}$, se tiene
$\varphi(\beta)>0$).

M\'as a\'un, $K_H^+/K_H$ es una extensi\'on $2$--elemental 
abeliana finita y los
primos que se descomponen totalmente en $K_H^+/K$ son los ideales primos
principales generados por un elemento totalmente 
positivo (Proposici\'on \ref{CClaseP4.7.10'}).

Dada un extensi\'on $E/F$ finita de campos num\'ericos, el {\em campo de
g\'eneros\index{campo de g\'eneros}} $\g E$ relativo a $F$ es el campo
$\g E=E\* k$ donde $\*k$ es la m\'axima extensi\'on abeliana de $F$ de tal
suerte que $E\*k\subseteq E_H$.
\[
\xymatrix{
E\ar@{-}[r]\ar@{-}[d]&E\*k=\g E\ar@{-}[r]\ar@{-}[d]&E_H\\ F\ar@{-}[r]&\*k}
\]
Equivalentemente, $\*k$ es la m\'axima extensi\'on abeliana de $F$
contenida en $E_H$. Esta es la definici\'on dada por Fr\"ohlich \cite{Fro83}.

En cierta forma $\g E$ es la parte m\'as f\'acil de estudiar del campo
de clase de Hilbert $E_H$ de $E$.

Cuando $\K$ es un campo de funciones congruente, para estudiar el
campo de g\'eneros de $\K$, digamos relativo a $K=\F(T)$, donde
$\F$ es el campo de constantes de $\K$, el problema al que nos
enfrentamos es, ?`qu\'e es $\K_H$?

Si definimos $\K_H$ como la m\'axima extensi\'on abeliana 
no ramificada de $K$, tendremos $\K_H=\nr\K\supseteq \K\bar{\ma F}_q$
donde $\bar{\ma F}_q=\bigcup_{n=1}^{\infty}{\ma F}_{q^n}$ es a la vez la
m\'axima extensi\'on abeliana de $\F$ como su cerradura algebraica.
De hecho se tiene que
\[
\Gal(\nr\K/\K)\cong \Gal(\nr\K/\K\bar{\ma F}_q)
\times \Gal(\K\bar{\ma F}_q/\K)\cong
I_{\K,0}\times \hat{\ma Z}
\]
donde $\hat{\ma Z}=\lim_{\leftarrow n}{\ma Z}/n
{\ma Z}\cong \prod_{p\text{\ primo}}{\ma Z}_p$ es el anillo de 
Pr\"ufer y $I_{\K,0}$ es el grupo de clases de divisores de grado
$0$ de $\K$ (ver Cap\'tulo \ref{CClaseC17}, Secci\'on \ref{CClaseC4}).

Este grupo es demasiado grande
 para su estudio. Notemos que todos los
primos de $\K$ son eventualmente inertes en $\K\bar{\ma F}_q/\K$. Esto es,
si $\pK$ es un primo de $\K$, existe $n_0(\pK)\in{\ma N}$ tal que
$\pK$ se descompone en $\K{\ma F}_{q^{n_0}}$ y $\pK$ es
totalmente inerte en $\K\bar{\ma F}_q/\K{\ma F}_{q^{n_0}}$ (Teorema
\ref{T6.1.4}). En otras palabras, si exigimos como condici\'on extra
que alg\'un conjunto de primos de $\K$ se descomponga totalmente
en la definici\'on de $\K_H$, entonces $\K_H/\K$ ser\'a finito.
Esta es definici\'on de campos de clase de Hilbert dada por 
M. Rosen (\cite{Ros87}). Ver Cap\'itulo \ref{CClaseC17}, Secci\'on
\ref{CClaseS4.9-1}.

En este cap\'itulo estudiaremos el campo de g\'eneros $\g\K$ de una
extensi\'on $\K/K$, $K=\F(T)$ donde el campo de constantes de
$\K$ es una extensi\'on finita de $\F$.

\section{Algo de historia y algunos antecedentes}\label{S12*.1-1}

Para el caso num\'erico, ver Secci\'on \ref{S12.4}. 

El concepto del {\em campo de g\'eneros\index{campo de g\'eneros}}
se remonta a F. Gauss \cite{Gau1801} en el contexto de formas
cuadr\'aticas binarias. Si $\K/{\ma Q}$ es un campo num\'erico,
$\g\K=\K\*k$ donde $\*k$ es la m\'axima extensi\'on abeliana de 
${\ma Q}$ y tal que $\K\subseteq \g\K\subseteq \K_H$. Originalmente,
la definici\'on de campos de g\'eneros fue dada para una extensi\'on
cuadr\'atica de ${\ma Q}$. Tenemos que, para una extensi\'on
cuadr\'atica de un campo num\'erico $\K$, el grupo de Galois de
$\g\K/\K$ es isomorfo al m\'aximo subgrupo de exponente $2$ del
grupo de clase de ideales (ver Secci\'on \ref{S12.4}). 
D. Hilbert fue el primero en traducir el trabajo de Gauss al concepto
de campos de g\'eneros que es usado actualmente. H. Hasse \cite{Has51}
fue quien introdujo la teor\'ia de g\'eneros para campos cuadr\'aticos
num\'ericos describiendo el trabajo de Gauss por medio de la
teor\'ia de campos de clase. H. W. Leopoldt \cite{Leo53} determin\'o
el campo de g\'eneros $\g\K$ de una extensi\'on abeliana $\K$ de
${\ma Q}$ generalizando el trabajo de Hasse.

El desarrollo de Leopoldt fue por medio del uso de los
caracteres de Dirichlet adem\'as de relacionar estos caracteres
con la aritm\'etica de $\K$ (Secci\'on \ref{S12.3}). Fue A. 
Fr\"ohlich \cite{Fro59-1,
Fro59-2} quien introdujo la noci\'on que hemos descrito y que sirve
para el estudio del campo de g\'eneros para extensiones no 
necesariamente abelianas de ${\ma Q}$.

En este cap\'itulo estamos interesados en campos de funciones
congruentes $\K$. Como mencionamos anteriormente, no hay
una noci\'on en general de campo de clase de Hilbert (CCH).
Usaremos la noci\'on de campo de clase de Hilbert introducida
por M. Rosen, en la cual se fija un conjunto no vac\'io $S_{\infty}$
de primos de $\K$.

Usando la definici\'on de Rosen de CCH es posible dar un concepto
de campos de g\'eneros de $\K$ similar al caso de campos 
num\'ericos. En la literatura ha habido diferentes definiciones de
campos de campos de g\'eneros de acuerdo a las diversas definiciones
de CCH dada. R. Clement \cite{Cle92} encontr\'o un campo de
g\'eneros extendido de una extensi\'on c\'iclica $\K$ de $K=\F(T)$
de grado primo dividiendo a $q-1$. Clement us\'o un concepto de
CCH similar al usado por Hasse de campos cuadr\'aticos de n\'umeros
$\K$. Espec\'ificamente, esta es la extensi\'on abeliana finita de
$\K$ tal que los ideales primos del anillo de enteros ${\mc O}_{\K}$
de $\K$ que se descomponen totalmente son precisamente los
ideales principales generados por un elemento de norma positiva,
lo cual coincide, en este caso, con los ideales principales generados
por elementos totalmente positivos.

S. Bae y J.K. Koo \cite{BaeKoo96} generalizaron los resultados de
Clement con los m\'etodos desarrollado por Fr\"ohlich \cite{Fro83}.
Ellos definieron el campo de g\'eneros para cualquier campo de
funciones global y desarrollaron el an\'alogo a la teor\'ia cl\'asica
de campos de g\'eneros. B. Angl\`es y J.-F. Jaulent \cite{AngJau2000}
usaron $S$--grupos de clase extendido para establecer los 
resultados para la teor\'ia de extensiones finitas de campos 
globales, donde $S$ es un conjunto finito no vac\'io de lugares.

Para esta cap\'itulo usaremos las siguientes notaciones.
$K=\F(T)$ ser\'a el campo de funciones racionales.
Sea $R_T^+$ el conjunto de polinomios m\'onicos e irreducibles
en $R_T$. 
Para cualquier campo de funciones $\K/{\ma F}_q$, 
$\K_m:= \K{\ma F}_{q^m}$ denota a la extensi\'on de constantes.
Para cualquier $m\in{\ma N}$,
$C_m$ denota un grupo c{\'\i}clico de orden $m$.

Para cualquier extensi\'on finita $\K/K$ usaremos el s\'imbolo
$\S \K$ para denotar ya sea un primo o todos los primo en $\K$
sobre $\p$, el divisor de polos de $T$ en $K$. 
Recordemos que los
primos que se ramifican en $K(\Lambda_N)/K$ son ${\eu p}_{\infty}$
($q\neq 2$)
y los polinomios $P\in R_T^+$ tales que $P\mid N$, con la excepci\'on
cuando $q=2$ y $N\in\{T, T+1, T(T+1)\}$ 
pues en este caso $K(\Lambda_N)=K$.

Establecemos $L_n$ como el m\'aximo subcampo de
$K(\Lambda_{1/T^n})$ donde ${\eu p}_{
\infty}$ es total y salvajemente ramificado,
$n\in{\ma N}$. Para cualquier campo $F$,
$_nF$ denota la composici\'on $FL_n$.

La definici\'on de Rosen para un campo de clase de
Hilbert relativo de un campo de funciones congruente $\K$,
es la siguiente.

\begin{definicion}[\cite{Ros87}]\label{D3.1}
Sea $\K$ un campo de funciones con campo de constantes ${\ma F}_q$.
Sea $S$ cualquier conjunto finito no vac{\'\i}o de
divisores primos de $\K$. El
{\em campo de clase de Hilbert de $\K$ relativo a $S$}, $\K_{H,S}$, es el
la m\'axima extensi\'on abeliana no ramificada de
$\K$ donde cada elemento de $S$ se descompone totalmente.
\end{definicion}

A partir de ahora, para cualquier extensi\'on finita 
$\K$ de $K$ consideraremos  $S$ 
como el conjunto de divisores primos que dividen a
${\eu p}_{\infty}$, el divisor de polos de
 $T$ en $K$ y escribiremos $\K_H$ en lugar de 
 $\K_{H,S}$.

\begin{definicion}\label{D3.2}
Sea $\K$  una extensi\'on finita de $K$. El {\em
campo de g\'eneros\index{campo de g\'eneros}} 
$\K_{\eu {ge}}$ de $\K$ es la m\'axima extensi\'on de $\K$
contenida en $\K_H$ que sea la composici\'on de
$\K$ y una extensi\'on abeliana de $K$. Equivalentemente,
$\K_{\eu {ge}}=\K K^{\ast}$ donde $K^{\ast}$ es la m\'axima extensi\'on abeliana
de $K$ contenida en $\K_H$.
\end{definicion}

\begin{definicion}\label{D12*.2.2'} El grado $[\g\K:\K]$ se llama
el {\em n\'umero de g\'enero\index{numero de genero@n\'umero de g\'enero}} de $\K$ y
al grupo $\Gal(\g\K/\K)$ se le llama {\em grupo de g\'eneros\index{grupo
de g\'eneros}} de $\K$.
\end{definicion}

\section{Campo de g\'eneros. Resultados generales}\label{S6.7}

\begin{teorema}\label{T12*.2.2.A} Sea $\K$ un campo de funciones
global con campo de constantes $\F$. Sea $S$ cualquier conjunto
finito no vac\'io de lugares de $\K$. Sea ${\mc O}_S=\{x\in\K\mid
v_{\pK}(x)\geq 0\text{\ para toda\ } \pK\notin S\}$. Sea $Cl({\mc O}_S)=
Cl_S$
el grupo de clases de ideales del anillo Dedekind ${\mc O}_S$. 
Entonces el campo de clase de Hilbert 
$\K_{H,S}$ de $\K$ con respecto 
a $S$, satisface
\[
\Gal(\K_{H,S}/\K)\cong Cl_S
\]
y el campo de constantes de $\K_{H,S}$ es ${\ma F}_{q^d}$
donde $d=\mcd(\deg \pK\mid \pK\in S)$. En particular, la extensi\'on
$\K_{H,S}/\K$ es finita.
\end{teorema}

\begin{proof} Por el Teorema \ref{T6.1.4}, si $\K_r$ es la extensi\'on de
constantes de grado $r$, se tiene que los elementos de $S$
se descomponen totalmente en $\K_r/\K\iff r|\deg \pK$ para toda
$\pK\in S$. De ah\'i se sigue que $\K_d$ es la m\'axima 
extensi\'on de constantes de $\K$ donde los elementos de $S$
se descomponen totalmente. Por tanto ${\ma F}_{q^d}$ es el
campo de constantes  de $\K_{H,S}$. 
Ver Corolario \ref{CClaseC4.9.11}.

Ahora $\Gal(\K_{H,S}/\K)\cong Cl_S$ bajo el
mapeo de reciprocidad de Artin (Teorema \ref{CClaseT4.9.6}).
La finitud de la extensi\'on $\K_{H,S}/\K$ es
el contenido del Corolario \ref{CRamDed1.2.6}.
$\fin$
\end{proof}

\begin{teorema}\label{T12*.2.2.B} Sea $K_0$ un campo de funciones
congruente y sea $\K/K_0$ una extensi\'on abeliana finita con grupo
de Galois $G=\Gal(\K/K_0)$. Sea $S$ cualquier conjunto finito no vac\'io
de lugares de $\K$. Sea $\K_{H,S}$ el respectivo campo de clase
de Hilbert de $\K$. Entonces $\K_{H,S}/K_0$ es una extensi\'on
de Galois y ${\mc G}:=\Gal(\K_{H,S}/K_0)$ es el producto 
semidirecto de $H=\Gal(\K_{H,S}/\K)$ y $G$ con $H\normal {\mc G}$:
${\mc G}=H\rtimes G$:
\[
\Gal(\K_{H,S}/K_0)\cong \Gal(\K_{H,S}/\K)\rtimes \Gal(\K/K_0).
\]
\end{teorema}

\begin{proof} Sea $\sigma$ un $K_0$--monomorfismo de $\K_{H,S}$
en una cerradura normal de $\K_{H,S}/K_0$. Entonces 
$\sigma(\K)=\K$ y $\K_{H,S}\cdot \sigma(\K_{H,S})$ es una extensi\'on 
abeliana de $\K$ no ramificada y donde todos los primos de $S$
se descomponen totalmente (ver Proposici\'on \ref{P.Ram10(1)}). 
Por maximalidad se sigue que
$\K_{H,S}\cdot \sigma(\K_{H,S})=\K_{H,S}$, 
de donde $\sigma(\K_{H,S})=
\K_{H,S}$ y por tanto $\K_{H,S}/K_0$ es de Galois.
\[
\xymatrix{
&\K_{H,S}\ar@{-}[dl]_H\ar@{-}[ddl]^{\mc G}\\
\K\ar@{-}[d]_G\\ K_0}
\]

Se tiene que $H\normal {\mc G}$ y tenemos la sucesi\'on exacta
\[
1\lra H\xrightarrow{\ i\ }{\mc G}\xrightarrow{\ \pi\ }G\lra 1
\]
donde $G$ act\'ua en $H$ bajo conjugaci\'on.

M\'as precisamente, si $g\in G$ y $h\in H$, consideremos $x\in
{\mc G}$ tal que $g=\pi(x)$. Entonces $g\circ h:=xhx^{-1}=
h^{x^{-1}}$. Por ser $H$ abeliano, $g\circ h$ no depende de $x$.

Ahora bien, por ser $\K_{H,S}$ el campo de clase de Hilbert,
$H\cong Cl_S$ bajo el mapeo de reciprocidad de Artin.
Expl\'icitamente estamos considerando el mapeo
$Cl_S\xrightarrow{\artin{\K_{H,S}/\K}{}} H$
donde para un ideal 
\[
{\eu a}=\prod_{i=1}^r \pK_i^{\alpha_i},\quad
\xbinom{\K_{H,S}/\K}{\eu a}:=\prod_{i=1}^r\xbinom{\K_{H,S}/\K}{\pK_i}^{\alpha_i}
\]
y $\artin{\K_{H,S}/\K}{\pK_i}$ es el automorfismo de Frobenius correspondiente
a $\pK_i$ (ver Cap\'itulo \ref{CClaseC17}, Secci\'on \ref{CClaseC4}).

Identificamos $Cl_S\subseteq {\mc G}$ mediante el
mapeo $\artin{\K_{H,S}/\K}{}$. Esto es $Cl_S\cong
H\hooklongrightarrow {\mc G}$. La acci\'on de $G$ en 
$Cl_S$ es conjugaci\'on, esto es, si $\bar{\eu a}\in
Cl_S\subseteq {\mc G}$ y $\sigma \in G$, entonces
$\sigma\circ \bar{\eu a}=\sigma \bar{\eu a}\sigma^{-1}$, lo cual
tiene sentido pues ${\eu a}$ es un ideal obtenido en $\K$ y
$G=\Gal(\K/K_0)$.

Sea $R\subseteq {\mc G}$ un conjunto de representantes de $G$
en ${\mc G}$, esto es, para cada $\sigma \in G$, seleccionamos 
\'unicamente un elemento 
$r_{\sigma}\in{\mc G}$ tal que $\pi(r_{\sigma})=
\sigma$. Entonces si $y\in{\mc G}$, $\pi(y)=\sigma=\pi(r_{\sigma})$
se tiene $r_{\sigma}^{-1}y\in\ker\pi=H\cong Cl_S$
y por tanto $y=r_{\sigma}\bar{\eu a}=\sigma \bar{\eu a}$ para
alg\'un $\bar{\eu a}\in Cl_S$.

Se tiene para $\sigma,\tau\in G$ y para $\bar{\eu a}, \bar{\eu b}\in
Cl_S$,
\[
(\sigma\bar{\eu a})(\tau\bar{\eu b})=\sigma\tau\tau^{-1}
\bar{\eu a}\tau\bar{\eu b}=(\sigma\tau)(\bar{\eu a}^{\tau}
\bar{\eu b}).
\]
Por tanto $y\longmapsto (\bar{\eu a},\sigma)$ es un
isomorfismo, ${\mc G}\cong Cl_S\rtimes G\cong
\Gal(\K_{H,S}/\K)\rtimes \Gal(\K/K_0)$. $\fin$
\end{proof}

El Teorema \ref{T12*.2.2.B} tiene varias consecuencias importantes.
Ver el Teorema \ref{CClaseT4.10.1}.

\begin{teorema}\label{T12*.2.2.C} Con las hip\'otesis 
del Teorema {\rm{\ref{T12*.2.2.B}}},
se tiene que el grupo de g\'eneros $\Gal(\g\K/\K)$ satisface
\begin{gather*}
\Gal(\g\K/K_0)\cong {\mc G}/I_GCl_S\quad
\text{y}\quad\Gal(\g\K/\K)
\cong Cl_S/I_G Cl_S,\\
\intertext{donde}
I_G{\mc G}=
\langle (\sigma-1)x\mid \sigma\in G, x\in{\mc G}\rangle=
\langle x^{\sigma^{-1}}x^{-1}=(\sigma\circ x) x^{-1}\mid
\sigma\in G, x\in{\mc G}\rangle.
\end{gather*}

En particular, ${\mc G}'=I_GCl_S$.
\end{teorema}

\begin{proof} Se tiene que $\g\K$ es la m\'axima 
extensi\'on abeliana de
$K_0$ contenida en $\K_{H,S}$. Por tanto $\Gal(\g\K/K_0)\cong
{\mc G}/{\mc G}'$ y $\Gal(\g\K/\K)\cong Cl_S/{\mc G}'$.

Veamos que ${\mc G}'=I_GCl_S$. Se tiene para
$\sigma\in G$, $\bar{\eu a}\in Cl_S$, 
$(\sigma-1)\bar{\eu a}=(\sigma\circ\bar{\eu a})\bar{\eu a}^{-1}
=\sigma\bar{\eu a}\sigma^{-1}\bar{\eu a}^{-1}\in {\mc G}'$ por
lo que $I_G Cl_S\subseteq {\mc G}'$.

Para el rec\'iproco, primero notemos que como ${\mc G}/Cl_S
\cong G$ es abeliano, se tiene ${\mc G}'\subseteq Cl_S$. Ahora,
sean $x,y\in{\mc G}$, $x=\sigma\bar{\eu a}$, $y=\tau\bar{\eu b}$
con $\sigma,\tau\in G$, $\bar{\eu a},\bar{\eu b}\in Cl_S$.

Se tiene que $x^{-1}=\bar{\eu a}^{-1}\sigma^{-1}=\sigma^{-1}
\sigma\bar{\eu a}^{-1}\sigma^{-1}=\sigma^{-1}\bar{\eu a}^{-\sigma^{-1}}$
y $y^{-1}=\tau^{-1}\bar{\eu b}^{-\tau^{-1}}$. Por tanto
\[
xyx^{-1}y^{-1}=(\sigma\tau\bar{\eu a}^{\tau}\bar{\eu b})
(\sigma^{-1}\tau^{-1}\bar{\eu a}^{-\sigma^{-1}\tau^{-1}}\bar{\eu b}^{
-\tau^{-1}}).
\]

Tambi\'en notemos que si $\bar{\eu c}\in Cl_S$ y
$\delta\in G$, $\bar{\eu c}\delta=\delta\delta^{-1}\bar{\eu c}\delta
=\delta\bar{\eu c}^{\delta}$. Se sigue que
\begin{align*}
xyx^{-1}y^{-1}&=\sigma\tau\sigma^{-1}\tau^{-1}\bar{\eu a}^{
\tau\sigma^{-1}\tau^{-1}}\bar{\eu b}^{\sigma^{-1}\tau^{-1}}
\bar{\eu a}^{-\sigma^{-1}\tau^{-1}}\bar{\eu b}^{-\tau^{-1}}\\
&=\bar{\eu a}^{\sigma^{-1}(1-\tau^{-1})}\bar{\eu b}^{\tau^{-1}(
\sigma^{-1}-1)}\in I_G Cl_S.
\end{align*}.

Por lo tanto ${\mc G}'\subseteq I_G Cl_S$ y
${\mc G}'=I_G Cl_S$. $\fin$
\end{proof}

\begin{corolario}\label{C12*.2.2.D} Con las hip\'otesis del
Teorema {\rm{\ref{T12*.2.2.B}}},
si en adici\'on $G$ es c\'iclico, generado por $\sigma$, entonces
\begin{gather*}
\Gal(\g\K/\K)\cong \frac{Cl_S}{Cl_S^{1-\sigma}}.
\tag*{\fin}
\end{gather*}
\end{corolario}

\begin{corolario}\label{C12*.2.2.E} Sea $\K/K_0$ una extensi\'on
c\'iclica de grado primo $l$. Suponemos que
el orden del n\'umero de $S$--clases
de $K_0$ es $1$, en especial $K_0=K=\F(T)$.
Entonces $\g\K/K_0$ y $\g\K/\K$ son
ambas extensiones $l$--elementales abelianas.
\end{corolario}

\begin{proof} Sea $G=\Gal(\K/K_0)$.
Se tiene $\Gal(\g\K/\K)\cong \frac{Cl_S}{
Cl_S^{1-\sigma}}$. Con las hip\'otesis que tenemos sobre
$K_0$, la norma de las $S$--clases de $\K$ sobre las
$S$--clases de $K_0$ es trivial por lo que
$\frac{Cl_S}{Cl_S^{1-\sigma}}\cong H^{-1}(G,Cl_S)$
donde $G=\Gal(\K/K_0)$. Puesto que $G^l=\{1\}$, se tiene
que $l H^{-1}(G,Cl_S)=\{0\}$. En particular, $\Gal(\g\K/\K)$
es un grupo $l$--elemental abeliano.

Finalmente, $\Gal(\g\K/K_0)\cong \Gal(\g\K/\K)\times \Gal(\K/K_0)$
tambi\'en es un grupo $l$--elemental abeliano. $\fin$
\end{proof}

\begin{observacion}\label{O14.3.5'}
Se tiene, con la misma demostraci\'on, que si el orden del
n\'umero de $S$-clases de $K_0$ es $1$ y si $\K/K_0$ es
c\'iclica de grado $m$, entonces $\g\K/K_0$ y $\g\K/\K$ son
ambas extensiones abelianas de exponente $m$.
\end{observacion}

Consideremos $K=\F(T)$, $S_0=\{\p\}$ y $\K$ una extensi\'on
c\'iclica de grado primo $l$ de $K$. Sea $\AE {\K}$ la 
la cerradura entera de $R_T=\F[T]$ en $\K$ y sea $S=\{
\pK\in{\ma P}_{\K}\mid \pK|\p\}$. Sea $Cl(\AE {\K})=Cl_S$
el grupo de clases de ideales del dominio Dedekind $\AE\K=
\{x\in\K\mid v_{\eu q}(x)\geq 0\text{\ para todo\ }{\eu q}\notin S\}$.

Por el Corolario \ref{C12*.2.2.E}, $\g\K/K$ y $\g\K/\K$ son 
$l$--extensiones elementales abelianas. Ahora bien, si $G=
\Gal(\K/K)=\langle\sigma\rangle$, se tiene que si $\Phi_l(\sigma)
=1+\sigma+\cdots+\sigma^{l-1}$, entonces para $\bar{\eu a}\in
Cl_S$, $\Phi_l(\sigma)\bar{\eu a}=\bar{\eu a}^{\Phi_l(\sigma)}\in
Cl_{S_0}$ en $K$, pero $Cl_{S_0}=\{1\}$ por lo que
$\Phi_l(\sigma)\bar{\eu a}=\bar{\eu a}^{\Phi_l(\sigma)}=
\N_{\K/K}(\bar{\eu a})=1$.

\begin{lema}\label{L12*.2.2.F}
Se tiene
\begin{gather*}
(\sigma -1)^{l-1}=lh(\sigma)+\Phi_l(\sigma),\\
l=(\sigma-1)^{l-1}f(\sigma)+\Phi_l(\sigma)g(\sigma),
\end{gather*}
para algunos $h(X), f(X), g(X)\in{\ma Z}[X]$.
\end{lema}

\begin{proof}
Sea $\zeta_l$ una $l$--ra\'iz primitiva de la unidad y consideremos
${\ma Z}[\zeta_l]=\AE {\cic l{}}\cong {\ma Z}[X]/\langle\Phi_l(X)\rangle$,
(Proposici\'on \ref{P1.2.1.6}), donde el \'ultimo isomorfismo est\'a dado
por $\zeta_l\mapsto X\bmod \Phi_l(X)$, y $\langle\Phi_l(1)\rangle=
\langle l\rangle=\prod_{i=1}^{l-1}\langle 1-\zeta_l^i\rangle=\langle
1-\zeta_l\rangle^{l-1}$. 

En particular $l (1-\zeta_l)^{1-l}=\eta$ es una unidad en ${\ma Z}[
\zeta_l]$. Sea $\eta=f(X)\bmod \Phi_l(X)\in {\ma Z}[X]/\langle\Phi_l
(X)\rangle \cong {\ma Z}[\zeta_l]$ para $f(X)\in {\ma Z}[X]$. Por
tanto 
\[
l=(1-X)^{l-1}f(X)\bmod \Phi_l(X).
\] 
Se sigue que existe
$g(X)\in{\ma Z}[X]$ tal que $l=(1-X)^{l-1}f(X)+\Phi_l(X) g(X)$, de
donde obtenemos $l=(1-\sigma)^{l-1}f(\sigma)+\Phi_l(\sigma)g(\sigma)$.

Por otro lado, se tiene que
\[
(X-1)^l=\sum_{i=0}^l\binom li (-1)^iX^{l-i}=X^l+\sum_{i=1}^{l-1}
\binom li (-1)^i X^{l-i}+(-1)^l,
\]
y $l|\binom li$ para $1\leq i\leq l-1$ por lo que $(X-1)^l=X^l+(-1)^l
+l\varphi(X)$ para alg\'un $\varphi(X)\in {\ma Z}[X]$. Adem\'as,
$0=(1-1)^l=1+(-1)^l+l\varphi(1)$ lo cual es igual a $l\varphi(1)$ para
$l\neq 2$. Por tanto, para $l\neq 2$, $\varphi(1)=0$.

Consideremos $l\neq 2$. Entonces $\varphi(X)=(X-1)h(X)$ con
$h(X)\in {\ma Z}[X]$. Por tanto
\[
(X-1)^l=X^l-1+l(X-1)h(X),
\]
por lo que $\frac{(X-1)^l}{X-1}=(X-1)^{l-1}=\frac{X^l-1}{X-1}+lh(X)=
\Phi_l(X)+lh(X)$.

Para $l=2$, $(X-1)^2=X^2+1-2X$, $\frac{(X-1)^2}{X-1}=X-1=
(X+1)-2=\Phi_2(X)+2(-1)$. Por tanto, para todo $l$,
\[
(X-1)^{l-1}=\Phi_l(X)+lh(X)
\]
con $h(X)\in {\ma Z}[X]$. As\'i,
$(\sigma-1)^{l-1}=lh(\sigma)+\Phi_l(\sigma)$ con $h(X)\in{\ma Z}[X]$.
$\fin$
\end{proof}

\begin{corolario}\label{C12*.2.2.G}
Con las condiciones anteriores, tenemos que
\[
lCl_S=
Cl^l_S=Cl_S^{\sigma-1)^{l-1}}.
\]
\end{corolario}

\begin{proof}
Puesto que $\Phi_l(\sigma)\bar{\eu a}=1$ para $\bar{\eu a}\in Cl_S$,
se sigue del Lemma \ref{L12*.2.2.F} que, para $\bar{\eu a}\in Cl_S$,
$\bar{\eu a}^{(\sigma-1)^{l-1}}=\bar{\eu a}^{h(\sigma)l+\Phi_l(\sigma)}=
(\bar{\eu a}^{h(\sigma)})^l$
y $\bar{\eu a}^l=\bar{\eu a}^{f(\sigma)(\sigma-1)^{l-1}+\Phi(\sigma)
g(\sigma)}=(\bar{\eu a}^{f(\sigma)})^{(\sigma-1)^{l-1}}$ de donde
se sigue el resultado.
$\fin$
\end{proof}

Se tiene la cadena
\[
Cl_S\supseteq Cl_S^{(\sigma-1)}\supseteq Cl_S^{(\sigma-1)^2}
\supseteq \ldots\supseteq Cl_S^{(\sigma-1)^{l-1}}=Cl_S^l.
\]
En particular se tiene que $\frac{Cl_S^{(\sigma-1)^j}}{Cl_S^{(\sigma-1)^{j+1}}}$,
es un $l$--grupo elemental abeliano para $j=0,1,\ldots,l-2$. Definimos
\[
\lambda_i:=\dim_{{\ma F}_l} \frac{Cl_S^{(\sigma-1)^{i-1}}}{Cl_S^{
(\sigma-1)^{i}}},\quad i=1,\ldots, l-1.
\]

\begin{lema}\label{L12*.2.2.H}
Se tiene $l^{\lambda_1}=[\g\K:\K]$. Adem\'as $\lambda_i\geq
\lambda_{i+1}$ para $i=1,2,\ldots,l-2$ y el $l$--rango de $Cl_S$
es igual a $\lambda_1+\cdots+\lambda_{l-1}$.
\end{lema}

\begin{proof}
Puesto que $\Gal(\K/\K)\cong \frac{Cl_S}{Cl_S^{(\sigma-1)}}$, se
sigue que $[\g\K:\K]=l^{\lambda_1}$.

Ahora sea $\theta_i\colon \frac{Cl_S^{(\sigma-1)^{i-1}}}{Cl_S^{
(\sigma-1)^{i}}}\lra \frac{Cl_S^{(\sigma-1)^i}}{Cl_S^{
(\sigma-1)^{i+1}}}$ dada por $\theta_i(\bar{\eu a})=\bar{\eu a}^{
(\sigma-1)}$, $i=1,2.\ldots,l-2$. Se tiene que $\theta_i$ es una
funci\'on bien definida y es un epimorfismo, de donde se sigue
que $\lambda_i\geq \lambda_{i+1}$, $1\leq i\leq l-2$.

Finalmente, se tiene
\begin{gather*}
\lambda_1+\cdots+\lambda_{l-1}=\sum_{i=1}^{l-1}\dim_{{\ma F}_l}
\frac{Cl_S^{(\sigma-1)^{i-1}}}{Cl_S^{(\sigma-1)^{i}}}=\dim_{
{\ma F}_l}\frac{Cl_S}{Cl_S^{l}}=l\text{--rango de $Cl_S$}.
\tag*{$\fin$}
\end{gather*}
\end{proof}

\begin{corolario}\label{C12*.2.2.I}
Se tiene que $l\nmid h_S=|Cl_S|\iff \K=\g\K$.
\end{corolario}

\begin{proof}
Si $\K=\g\K$ entonces $\lambda_1=0$. Por tanto $0\leq\lambda_i
\leq \lambda_1=0$ para $1\leq i\leq l-1$, por tanto $\lambda_1+
\cdots+\lambda_{l-1}=0$.

Rec\'iprocamente, si $\lambda_1+\cdots+\lambda_{l-1}=0$
se sigue que $\lambda_1=0$.
$\fin$
\end{proof}

\begin{observacion}\label{O12*.2.2(1).E}
El resultado del Corolario \ref{C12*.2.2.E} no se cumple cuando
el n\'umero de clases de ideales no es $1$. Por ejemplo, si denotamos
$Cl_{K_0,S}$ y $Cl_{\K,S}$ los respectivos grupos de clases de ideales
de $K_0$ y de $\K$, supongamos que $Cl_{S,K_0}\cong Q_{l}\times
C$ con $Q_l$ el $l$--subgrupo de Sylow de $Cl_{S,K_0}$ y $C$ de
orden primo relativo a $l$. Suponemos que $C\neq \{1\}$.

Se tiene que el n\'ucleo de la conorma $\con_{K_0/\K}\colon
Cl_{K_0,S}\lra Cl_{\K,S}$ es de orden $1$ o $l$ 
(Ejemplo \ref{E10.1.2.10}). En particular
la conorma es inyectiva en $C$ y en especial $C\subseteq Cl_{\K,S}$.
Ahora la norma en $C$ satisface $\N(C)=C^l=C$, esto es
$\ker\N\cap C=\{1\}$. De la sucesi\'on exacta
\[
1\lra \ker\N\lra Cl_{\K,S}\lra \N(Cl_{K_0,S})\lra 1,
\]
y del hecho de que $\frac{\ker \N}{Cl_{\K,S}^{(1-\sigma)}}\cong
H^{-1}(G,Cl_{\K,S})$ es un $l$--grupo elemental abeliano, se sigue
que $C$ tiene una copia isomorfa en $\frac{Cl_{\K,S}}{Cl_{\K,S}^{(1-
\sigma)}}\cong \Gal(\g{\K}/\K)$ por lo que $\Gal(\g{\K}/\K)$
no es un $l$--grupo.
\end{observacion}

Como hemos visto en el Cap\'itulo \ref{Ch11}, 
el Teorema de Kronecker--Weber
establece que si $\K/K$ es una extensi\'on abeliana finita, entonces
existen $n,m\in{\ma N}$, $N\in R_T^+$ tales 
que $\K\subseteq K(\Lambda_N)
{\ma F}_{q^m}L_n:={_nK}(\Lambda_N)_m$. 
Usando este resultado
podemos establecer m\'as sobre $\K$.

\begin{teorema}\label{P3.4} Sea $\K/K$ una extensi\'on abeliana
finita. Entonces
\las
\item Si $\p$ se descompone totalmente en $\K$, entonces $\K\subseteq 
\cicl N{}^+$ para alg\'un $N\in R_T$.

\item Si $d_{\K}(\p)=1$ entonces $\K\subseteq \cicl N{} L_n
={_n\cicl N{}}$ para alg\'un $N\in R_T$ y alg\'un $n\in{\ma N}\cup
\{0\}$.

\item Si el primo 
$\p$ es moderadamente ramificado en $\K$, entonces
$\K\subseteq \cicl N{}{\ma F}_{q^m}=
\cicl N{}_m$ para alg\'un $N\in R_T$ 
y alg\'un $m\in {\ma N}$.

\item Si $\p$ es moderadamente ramificado en $\K$ y $d_{\K}
(\p)=1$, entonces $\K\subseteq \cicl N{}$ para alg\'un $N\in R_T$.
\end{list}
\end{teorema}

\begin{proof}
Para ver un resultado relacionado con
(1), ver la Proposici\'on \ref{P12*.2.2.G}.

Sea $D:=D_{\K/K}(\pL|\p)$ el grupo de descomposici\'on de
$\pL|\p$, donde $\pL$ es cualquier primo de $\K$ sobre $\p$.
Se tiene que $D$ es de orden
 $|D|=e_{\K/K}(\pL|\p)f_{\K/K}(\pL|\p)=
e_{\infty}(\K|K)f_{\infty}(\K|K)=ef$. Sea $L:=\K^D$
el campo fijo de $D$, $K\subseteq L=\K^D\subseteq \K$.

Sea $I=I_{\K/K}(\pL|\p)$ el grupo de inercia
de $\pL|\p$ y sea $R=\K^I$.
Entonces $\K^I/K$ es la m\'axima subextensi\'on de $\K/K$
no ramificada en $\p$ y $\K^D/K$ es la m\'axima subextensi\'on
de $\K/K$ donde $\p$ se descompone totalmente.

Ahora, sea $V=V_{\K/K}^{(1)}(\pL|\p)=V_{\infty}^{(1)}(\K/K)$ 
el primer grupo de ramificaci\'on
de $\p$ en $\K/K$. Sea $S=K^V$. Entonces $\K^V/K$ es la
m\'axima subextensi\'on de $\K/K$ con $\p$ moderadamente
ramificado. 

Se tiene $|I|=e$, $|D|=ef$, $|V|=p^s$ con $e=e_0p^s$ y
$\mcd(e_0,p)=1$.
\[
\xymatrix{
\K\ar@/_5pc/@{-}[ddd]_D\ar@/^2pc/@{-}[dd]^I
\ar@/_1pc/@{-}[d]_S\ar@{-}[d]^{p^s}\ar@/_3pc/@{-}[dd]^e\\ 
\K^V=S \ar@{-}[d]_{e_0}\\ \K^I=R\ar@{-}[d]^f\\ 
\K^D=L\ar@{-}[d]_{h=\frac{[\K:K]}{ef}}\\
K\ar@/_5pc/@{-}[uuuu]_G
}
\]

Por el Teorema de Kronecker--Weber, existen $n\in{\ma N}\cup
\{0\}$, $m\in{\ma N}$ y $N\in R_T$ tales que $\K\subseteq
{_n\cicl N{}}_m=\cicl N{}M$ donde $M=L_nK_m$ y $K_m=
{\ma F}_{q^m}[T]$. Se tiene que $e_{\infty}(M|K)=q^n$,
$f_{\infty}(M|K)=m$, $[M:K]=q^n m=e_{\infty}(M|K)f_{\infty}
(M|K)$, $h_{\infty}(M|K)=1$ y $e_{\infty}(\cicl N{}|K)=q-1$,
$f_{\infty}(\cicl N{}|K)=1$.

Sea $E:=\K M\cap \cicl N{}$. Por la correspondencia de
Galois tenemos $\K M=EM$. 
\[
\xymatrix{
\cicl N{}\ar@{-}[rr]\ar@{-}[d]&& \cicl N{}M={_n\cicl N{}}_m
\ar@{-}[d]\\
E\ar@{-}[rr]\ar@{-}[dd]&& \K M=EM\ar@{-}[dl]\ar@{-}[ddd]\\ 
&\K\ar@{-}[dl]\\ \cicl N{}\cap \K\ar@{-}[d]\\ K\ar@{-}[rr]&&M
}
\]

Usamos la siguiente notaci\'on:
$e_{\infty}(U|W)=e^0_{\infty}(U|W)w_{\infty}(U|W)$ donde
$w_{\infty}(U|W)=|V^{(1)}_{\infty}(U|W)|$ y $\mcd\big(e^0_{
\infty}(U|W), p\big)=1$. Se tiene $f_{\infty}(E|K)=1$,
$e_{\infty}(E|K)|(q-1)$, $f_{\infty}(M|K)=f_{\infty}(K_m|K)
=m$, $e_{\infty}(M|K)=e_{\infty}(L_n|K)=q^n$, 
$e_{\infty}(EM|K)=q^ne_{\infty}(E|K)=w_{\infty}(EM|K)
e_{\infty}^0(EM|K)$ y $f_{\infty}(EM|K)=m$.

Se sigue que $(EM)^V=E_m$ pues $w_{\infty}(E_m|K)=1$
y $w_{\infty}(EM|E_m)=w_{\infty}(EM|K)=q^n=[EM:E_m]$.
La \'ultima igualdad se sigue de que $[EM:E_m]=\frac{
[EM:K]}{[E_m:K]}=\frac{[E:K]q^nm}{[E:K]m}$. 

Adem\'as
$(EM)^I=E_m^+$ y $(EM)^D=E^+$.

Finalmente, $f_{\infty}(E|K)=1$ y $f_{\infty}(EM|{_nE})=m=
[EM:{_nE}]$.

\las
\item Supongamos que $\p$ se descompone totalmente en
$\K/K$. Entonces $\K\subseteq (\K M)^D=(EM)^D=E^+
\subseteq \cicl N{}^+$.

\item Supongamos que $\p$ satisface $d_{\K}(\p)=1$. Entonces
$f_{\infty}(\K|K)=1$. Se sigue que $\K\subseteq (\K M)^{\Gal(EM/
{_nE})}={_nE}$.

\item Ahora sea $\p$ moderadamente ramificado en $\K/K$.
Entonces $\K\subseteq(\K M)^V=(EM)^V=E_m\subseteq \cicl N{}_m$.

\item Finalmente, si $\p$ es moderadamente ramificado y $d_{\K}(\p)
=1$, se tiene $\K\subseteq (\K M)^V\cap (\K M)^{\Gal(EM/{_nE})}
\subseteq \cicl N{}_m\cap {_n\cicl N{}}=\cicl N{}$.
$\fin$
\end{list}
\end{proof}

\begin{lema}\label{L12*.2.2.H-1} Sea $\K$ cualquier campo de funciones
congruente y sea $\g\K$ el campo de g\'eneros de $\K/K$. Sea 
$\K\subseteq L$ tal que $L/\K$ es abeliana y es no ramificada en
los primos finitos de $\K$ y tal que $\g\K\subseteq L$. Entonces
$\g\K=L^D$ donde $D$ es el grupo generado por todos
los grupos de descomposici\'on $D(\pL|\pK)$ para $\pK\in \S\K$
y $\pL|\pK$, $\pL\in{\ma P}_L$.
\end{lema}

\begin{proof} $
\xymatrix{
&&&\K_{H,S}\ar@{-}[dl]\\
&&L\ar@{-}[dl]\\
&\g\K\ar@{-}[dl]\\
\K\ar@{-}[d]\\ K}$

Se tiene que $L^D/\K$ satisface que $L/\K$ es no ramificada y todos los
primos de $\S\K$ se descomponen totalmente. Puesto que $\K\subseteq
L^D\subseteq L$, $L^D/\K$ es abeliana y por tanto $L^D\subseteq \K_{H,S}$.
Se sigue que $L^D\subseteq \g\K$.

Finalmente, $L^D\subseteq \g\K\subseteq L$, y si $\g\K\neq L^D$,
entonces alg\'un $\p\in\S\K$ no se descompondr\'ia totalmente en 
$\g\K/K$. Se sigue que $L^D=\g\K$. $\fin$
\end{proof}

\begin{proposicion}\label{P12*.2.2.I} Sean $M, \K$ dos campos de 
funciones congruentes, $K\subseteq M\subseteq \K$. Entonces,
si $\g M$ y $\g\K$ denotan los campos de g\'eneros de $M/K$ y de
$\K/K$ respectivamente. Si $\K/M$ es separable y finita,
se tiene $\g M\subseteq \g\K$.
\end{proposicion}

\begin{proof} Se tiene 
\[
\xymatrix{
\g M\ar@{-}[r]\ar@{-}[d]&\K\g M\ar@{-}[d]\\M\ar@{-}[r]&\K=\K M}
\]

Como $\g M/M$ es no ramificada y $\S M$ 
se descompone totalmente en $\g M/M$, entonces
$\K\g M/\K M=\K$ es no ramificada y $\S {\K}$ se
descompone totalmente en $\K\g M/\K$ (ver 
Corolario \ref{C.10.4.2.Ram}). 
De ah\'i se sigue que $\g M\subseteq \K \g M\subseteq \g\K$
(de hecho $\g M=M\* k$ con $\*k/K$
abeliana, $\K\g M=\K M\*k=\K\* k$). $\fin$
\end{proof}

\begin{proposicion}\label{P12*.2.2.J} Sean $M/K$ y 
$E/K$ dos extensiones
finitas y separables, $\g M$ y $\g E$ los campos de g\'eneros 
respectivos. Entonces
$\g M\g E\subseteq \g{(ME)}$.
\end{proposicion}

\begin{proof} Se sigue de la Proposici\'on \ref{P12*.2.2.I} 
pues $M,E\subseteq
ME$. $\fin$
\end{proof}

\begin{observacion}\label{O12*.2.2.K} En general veremos que no
necesariamente $\g M\g E=\g {(ME)}$, esto es, se puede tener
$\g M\g E\subsetneqq \g{(ME)}$ (ver Observaci\'on \ref{O3.6.B}).
\end{observacion}

\begin{observacion}\label{O12*.2.2.L}
Sea $K=\F(T)$ y sea $K_m$ la extensi\'on de constantes de 
$K$ de grado $m$. Sea $\K$ cualquier extensi\'on finita de 
$K$ tal que ${\ma F}_{q^m}\subseteq \K$. Sean $\g\K$ y
$\g {\K^{\prime}}$ los campos g\'eneros de $\K/K$ y de
$\K/K_m$ respectivamente. Puesto que los divisores
primos de $K$ se descompone totalmente en $\g \K$
y en $\g {\K^{\prime}}$ y $\g\K/\K$ y $\g {\K^{\prime}}/\K$
son abelianas y no ramificadas, se sigue que 
$\g\K=\g{\K^{\prime}}$.
\end{observacion}

\begin{proposicion}\label{P12*.2.2.M}
Sea $\K/K$ una extensi\'on finita y separable. Entonces
$\g{(\g\K)}=\g\K$.
\end{proposicion}

\begin{proof}
Por definici\'on, tenemos que $\g\K\subseteq \g{(\g\K)}$. 
Ahora bien, se tiene $\g{(\g\K)} = \g\K E$
donde $E/K$ es una extensi\'on abeliana, tal que 
$\g\K E/\g\K$ es no ramificada y $\S{\g\K}$ se descompone 
totalmente en $\g\K E$.
 
Sea $\g\K=\K F$. Entonces $\g{(\g\K)}=\g\K E=\K FE$.
Ahora $F/K$ y $E/K$ son extensiones abelianas por lo que
$FE/K$ es abeliana. Ahora, $\g\K=\K F/\K$ es no ramificada
y $\g{(\g\K)}=\g\K E/\g\K$ es no ramificada. Por tanto
$\g{(\g\K)}=\K FE/\K$ es no ramificada. Similarmente,
tenemos que $\S\K$ se descompone totalmente en $\g\K/\K$
y $\S{\g\K}$ se descompone totalmente en $\g{(\g\K)}/
\g\K$ por lo que $\S\K$ se descompone totalmente en
$\g{(\g\K)}$. Se sigue que $\g{(\g\K)}\subseteq \g\K$
de donde se obtiene la igualdad $\g{(\g\K)}=\g\K$.
$\fin$
\end{proof}

\begin{observacion}\label{O12*.2.2.N}
Cuando $\K/k$ es abeliana, la demostraci\'on de la Proposici\'on
\ref{P12*.2.2.M} es m\'as directa. En este caso $\g\K$ es la
m\'axima extensi\'on  abeliana de $K$ tal que $\g\K/\K$ es no
ramificada y $\S\K$ se descompone totalmente en $\g\K/\K$.

Por definici\'on $\g\K\subseteq \g{(\g\K)}$. Ahora bien, puesto que
$\g\K/K$ es abeliana, $\g{(\g\K)}$ es la m\'axima extensi\'on
abeliana de $K$ con $\g{(\g\K)}/\g\K$ no ramificada y tal que
$\S{\g\K}$ se descompone totalmente en $\g{(\g\K)}$. Se
sigue que $\g{(\g\K)}/\K$ es no ramificada y $\S\K$ se
descompone totalmente en $\g{(\g\K)}$. Por tanto 
$\g{(\g\K)}\subseteq \g\K$ y $\g{(\g\K)}=\g\K$.
\end{observacion}

Con respecto al campo de constantes de $\g\K$
de una extensi\'on abeliana de $\K/K$, se tiene:

\begin{lema}\label{L4.1.1}
Si $\K/K$ es una extensi\'on abeliana finita y si el grado
de cualquier divisor primo en $\S \K$ es
$t$, entonces el campo
de constantes de $\K_{{\eu {ge}}}$ es ${\ma F}_{q^t}$.
\end{lema}

\begin{proof}
Consideremos la extensi\'on de constantes
$\K_r:=\K{\ma F}_{q^r}$ de $\K$. Entonces el
n\'umero de primos en $\K_r$ sobre cualquier primo en
 $\S \K$ es $h=\mcd
(d_\K(\S \K),r)=\mcd(t,r)$ (Teorema \ref{T6.1.4}). 
Por lo tanto $\S \K$ se descompone completamente en
$\K_r/\K$ si y s\'olo si $h=r$ y esto es equivalente a $r\mid
d_\K(\S \K)=t$. Se sigue que la m\'axima extensi\'on de
constantes de $\K$ en donde
$\S \K$ se descompone totalmente es
$\K_t=\K{\ma F}_{q^t}$. Luego el campo de constantes de
$\K_{{\eu {ge}}}$ es ${\ma F}_{q^t}$.
$\fin$
\end{proof}

\section{Campos de g\'eneros para $\K$ subcampo de un 
campo ciclot\'omico}\label{S12*.3}

En el resto de este cap{\'\i}tulo desarrollamos un an\'alogo a la teor{\'\i}a
de Leopoldt del g\'enero para campos de funciones congruentes.
Se da una descripci\'on del campo de g\'eneros $\K_{\eu {ge}}$ 
de una extensi\'on abeliana finita de un campo de funciones
racionales congruente por medio de su grupo de caracteres
de Dirichlet en campos de funciones ciclot\'omicos.
Aqu{\'\i} consideramos el campo de clases de Hilbert
$\K_H$ de un campo de funciones $\K$ usando la construcci\'on de
Rosen para $S_{\infty}=
\{{\eu p}_{\infty}\}$, donde ${\eu p}_{\infty}$ es el divisor de polos
de $T$ en el campo de funciones racionales $K={\ma F}_q(T)$.

M\'as precisamente, sea $\K$ una extensi\'on abeliana finita
de $K$. Entonces, si $\K$ est\'a contenido en una extensi\'on
ciclot\'omica, veremos que $\K_{\eu {ge}}$ tambi\'en est\'a
contenido en una extensi\'on ciclot\'omica y encontraremos
el grupo de caracteres de Dirichlet asociado a $\K_{\eu {ge}}$.

Cuando $\K/K$ es una extensi\'on abeliana, 
$\K_{\eu {ge}}$ es la m\'axima
extensi\'on de $K$ contenida en $\K_H$. 
El principal objetivo
en esta secci\'on es encontrar $\K_{\eu {ge}}$ 
donde $\K$ es un subcampo
de un campo de funciones ciclot\'omico. 
En lo que sigue, $\K$
siempre denotar\'a una extensi\'on finita geom\'etrica de $K$.
Primero notamos que tenemos el an\'alogo al resultado de
Leopoldt (Teorema \ref{T12.4.2}).

\begin{proposicion}\label{P3.3}
Si $\K\subseteq K(\Lambda_N)$ y el grupo de caracteres asociado 
a $\K$ es $X$, entonces la m\'axima extensi\'on abeliana
$J$ de $\K$ no ramificada en ning\'un primo finito
$P\in R_T^+$, contenida en una extensi\'on ciclot\'omica, es el
campo asociado a $Y=\prod_{P\in R_T^+} X_P=\prod_{
P\mid N} X_P$.
\end{proposicion}

\begin{proof}
An\'aloga a la demostraci\'on del Teorema \ref{T12.4.2}. $\fin$
\end{proof}

\begin{observacion}\label{O12*.3.A} Podemos considerar al
campo asociado a $Y$ de la Proposici\'on \ref{P3.3} como el
{\em campo de g\'eneros extendido\index{campo de g\'eneros
extendido}} dentro del ciclot\'omico y lo denotaremos por
$\g \K^{(+)}$.
\end{observacion}

Notemos que $\g \K^{(+)}$ no es lo mismo que $\g \K^+=
\g\K\cap K(\Lambda_N)^+$, donde $\K\subseteq K(\Lambda_N)$.
Se tiene $\g\K^+\subseteq \g\K\subseteq \g\K^{(+)}$.

Ahora, $K(\Lambda_N)^+$ es el campo fijo de $\*\F\subseteq \G N
\cong \Gal(K(\Lambda_N)/K)$. El grupo de caracteres asociado a
$K(\Lambda_N)^+$ es 
\[
\{\chi\in X_N\mid \chi(a)=1\text{\ para toda\ } a\in\*\F\}
\]
 y donde $X_N$ es el grupo de caracteres asociado a
$K(\Lambda_N)$, esto es, $X_N\cong \G N$.

\begin{teorema}\label{T12*.3.B} Sea $K\subseteq \K\subseteq
K(\Lambda_N)$ para alg\'un $N\in R_T^+$. Entonces $\g\K
\subseteq K(\Lambda_N)$ y $\g\K=\K(\g{\K^{(+)}})^+$. M\'as
precisamente, si el grupo de caracteres de Dirichlet de $\K$
es $X$ y $L$ es el campo asociado a $Y=\prod_{P\in R_T^+}
X_P$, entonces 
\[
\g\K=\K L^+.
\] 
En particular se tiene que si
$\K\subseteq \cicl N{}$ entonces $\g\K\subseteq \cicl N{}$.
\end{teorema}

\begin{proof} Sea $F/\K$ tal que $F/K$ una extensi\'on abeliana no
ramificada y tal que los elementos de $\S\K$ se descomponen
totalmente $F/\K$. En particular $\p$ es moderadamente
ramificada.

Por la Proposici\'on \ref{P3.4}, $F\subseteq K(\Lambda_M)_m$
para algunos $M\in R_T^+$, $m\in {\ma N}$.

Sea ${\mc I}$ el grupo de inercia de $\S\K$ en $K(\Lambda_M)/K$
y sea $B=K(\Lambda_M)^{\mc I}$.
\[
\xymatrix{K(\Lambda_M)\ar@{-}[dd]_{\mc I}\ar@{-}[rr]&&K(\Lambda_M)_m
\ar@{-}[dd]\\ &F\ar@{-}[ddl] |!{[dl];[dr]}\hole\\  B\ar@{-}[d]\ar@{-}[rr]&&B_m\ar@{-}[d]\\
\K\ar@{-}[rr]\ar@{-}[d]&&\K_m\ar@{-}[d]\\
K\ar@{-}[rr]&&K_m}
\]

Puesto que los elementos de $\S B$ tienen grado $1$, estos son
totalmente inertes en $B_m/B$. Adem\'as los elementos de $\S B$ son
totalmente ramificados en $K(\Lambda_M)/B$. Ahora bien, los elementos
de $\S\K$ son totalmente descompuestos en $B/\K$ de donde se
obtiene que $B$ es el campo de descomposici\'on de $\S\K$ en
$K(\Lambda_M)_m/\K$. Se sigue que 
$F\subseteq B\subseteq K(\Lambda_M)$.

Sea $Z$ el grupo de caracteres de Dirichlet asociado a $F$. Puesto
que $F/\K$ es no ramificada, se sigue que $X\subseteq Z\subseteq Y$,
esto es, $F\subseteq L$ pues $L$ es la m\'axima extensi\'on abeliana
contenida en alg\'un campo de funciones ciclot\'omico tal que
$L/\K$ es no ramifica en los primos finitos. En particular, podemos
tomar $M=N$. Por lo tanto $\g\K=L^{\mc D}$ donde
${\mc D}$ es el grupo de descomposici\'on de $\S\K$ en $L/\K$.

Ahora bien, $\S\K$ se descompone totalmente en $L^+\K/\K$ pues
$\p$ se descompone totalmente en $L^+/K$. Puesto que $L/\K$ es
no ramificada $L^+\K\subseteq L$ por lo que $L^+\K/\K$ es no
ramificada, de donde $L^+\K\subseteq \g\K$ y se tiene $L^+\K
\subseteq \g\K\subseteq L$.

Para finalizar, veamos que $\S{L^+\K}$ es totalmente ramificada
en la extensi\'on 
$L/L^+\K$, lo cual se sigue del hecho que $L^+\subseteq L^+\K
\subseteq L$ y de que $\S{L^+}$ es totalmente ramificada en $L/L^+$.
Puesto que $L^+\K\subseteq \g\K\subset L$ y 
$\g\K/L^+\K$ es no ramificada, se sigue que $\g\K=L^+\K$. $\fin$
\end{proof}

\section{El campo de g\'eneros para extensiones 
abelianas de $K$}\label{S2.A}

Sea $K=\F(T)$.
Sea $\p$ el polo $T$ en $K$. Sea $\K/K$
una extensi\'on abeliana finita. Del Teorema de
Kronecker--Weber tenemos que existen
$n,m\in{\ma N}$ y $N\in R_T$ tales que
\[
\K\subseteq \ _nK(\Lambda_N)_m
:=L_n K(\Lambda_N){\ma F}_{q^m},
\]
donde $L_n$ denota el subcampo de
$K(\Lambda_{1/T^{n+1}})$ de
grado $q^n$ y ${\ma F}_{q^m}(T):=K_m$ 
es la extensi\'on de constantes de $K$ de grado
$m$. Tenemos que $\p$ es total y salvajemente
ramificado en $L_n/K$. Adem\'as $\p$ 
es totalmente inerte en $K_m/K$.

Para cualquier extensi\'on de
campos $E/F$, $e_{\infty}(E|F)$, $f_{\infty}(E|F)$ y $h_{\infty}(E|F)$ 
denotan el \'indice de ramificaci\'on, el grado de
inercia y el n\'umero de descomposici\'on de 
$\S F$ en $E$ respectivamente. Para
$P\in R_T^+$, $e_P(E|F)$ denota
el \'indice de ramificaci\'on de cualquier primo 
en $F$ sobre $P$ en $E/F$.

Sea $M:=L_nK_m$. Entonces
\begin{gather}\label{Eq1.A}
e_{\infty}(M|K)=q^n, \quad f_{\infty}(M|K)=m  \quad \text{y} 
\quad h_{\infty}(M|K)=1.
\end{gather}
Tenemos que $M\cap K(\Lambda_N)=K$. 

El resultado principal sobre el campo de g\'eneros
para extensiones abelianas de $K$ es el siguiente.

\begin{teorema}\label{T2.1.A} Con las notaciones anteriores, 
sea $\K/K$ una extensi\'on abeliana finita. Sea
$E:=\K M\cap K(\Lambda_N)$. Entonces 
\[
\g \K=\g E^{H_1} \K=(\g E \K)^H,
\]
donde
$H_1=H|_{\g E}$ y $H$ es el grupo de descomposici\'on
de cualquier primo de $\S \K$ en $\g E \K/\K$. Tambi\'en
tenemos que $|H_1|=d\mid q-1$. El campo de constantes de $\g \K$
es ${\ma F}_{q^t}$ donde $t$ es el grado de $\S \K$ en $\K$.
Finalmente, $\g E \K/\g \K$ y $E\K/E^{H_1}\K$
son extensiones de constantes de grado $d$.
\end{teorema}

\begin{proof} Que el campo de constantes 
$\g \K$ es ${\ma F}_{q^t}$
es el contenido del Lema \ref{L4.1.1}.

Puesto que $K(\Lambda_N)\cap M=K$ y 
$E=\K M\cap K(\Lambda_N)$,
de la correspondencia de Galois entre $K(\Lambda_N)/K$
y $K(\Lambda_N)M/M$, tenemos que
$E$ corresponde a $\K M$. Por tanto
$EM$ corresponde a $E$. Se sigue que
\begin{gather*}
\K M=EM.
\end{gather*}
\[
\xymatrix{
K(\Lambda_N)\ar@{-}[d]\ar@{-}[rrr]&&&K(\Lambda_N)M\ar@{-}[d]\\
E\ar@{-}[rrr]\ar@{-}[dd]&&&\K M=EM\ar@{-}[ddd]\\
&&\K\ar@{-}[ru]\ar@{-}[ddl]\\
E\cap \K\ar@{-}[rru]\ar@{-}[d]\\
K\ar@{-}[r]&\K\cap M\ar@{-}[rr]&&M
}
\]

Ahora $E\cap \K\subseteq \g E\cap 
\K\subseteq K(\Lambda_N)\cap \K
= (K(\Lambda_N)\cap \K M)\cap K(\Lambda_N)\cap \K=
E\cap K(\Lambda_N)\cap \K=E\cap \K$. Esto es
\begin{gather*}
E\cap \K=\g E\cap \K=K(\Lambda_N)\cap \K.
\end{gather*}

Tenemos
$[E:K]=[EM:M]= [\K M:M]=[\K:\K\cap M]$. Por tanto
\begin{gather}\label{Eq4.A}
[\K:K]=[E:K][\K\cap M:K].
\end{gather}

Veamos que $E\K/\K$ es no ramificada. Primero notamos
que $E\subseteq E\K\subseteq E\K M=E\cdot EM=EM$. En
la extensi\'on $M/K$, $\p$ es el \'unico primo ramificado.
Por tanto, en $\K M/E$ los \'unicos
posibles primos ramificados son aquellos en $\S E$.
Tambi\'en se tiene que en la extensi\'on $\K M/\K$ 
los \'unicos posibles primos ramificados son los elementos
de $\S\K$ y puesto que  $\K\subseteq E\K\subseteq
\K M=EM$, los \'unicos posibles primos ramificados en
$E\K/\K$ son aquellos en $\S \K$.
\[
\xymatrix{
E\ar@{-}[rr]\ar@{-}[dd]&&E\K\ar@{-}[dl]\ar@{-}[r]&EM=\K M
\ar@{-}[dll]\ar@{-}[ddd]\\
&\K\ar@{-}[dl]\\
E\cap \K\ar@{-}[d]\\
K\ar@{-}[rrr]&&&M}
\]

De la Ecuaci\'on (\ref{Eq1.A}) y del Proposici\'on 
\ref{P10.4.1.Ram}, obtenemos que
\begin{gather*}
\xymatrix{
\K\ar@{-}[r]\ar@{-}[d]&E\K\ar@{-}[d]\\
M\cap \K\ar@{-}[r]&\bullet \ar@{-}[r]&M}\\
e_{\infty}(E\K|\K)\mid e_{\infty}(M|\K\cap M)\quad \text{y}\quad
e_{\infty}(M|\K\cap M)\mid e_{\infty}(M|K)=q^n.
\end{gather*}

Por otro lado, se tiene
\begin{gather*}
e_{\infty}(E\K|\K)\mid e_{\infty}(E|E\cap \K) \quad \text{y} \quad
e_{\infty}(E|E\cap \K)\mid e_{\infty}(K(\Lambda_N)|K)=q-1.
\intertext{Por tanto}
e_{\infty}(E\K|\K)\mid \mcd(q^n,q-1)=1
\end{gather*}
de donde obtenemos que $E\K/\K$ es no ramificada.

Ahora bien, nuevamente por la Proposici\'on 
\ref{P10.4.1.Ram} se tiene
\[
\underbracket[0pt]{e_{\infty}(E\K|\K)}_{
\substack{\uigual\\ 1}} f_{\infty}(E\K|\K)\mid e_{\infty}(E|E\cap \K)
\underbracket[0pt]{f_{\infty}(E|E\cap \K)}_{\substack{\uigual\\ 1}},
\]
y $e_{\infty}(E\K|\K)=1$, $f_{\infty}(E|E\cap \K)=1$. Se obtiene que
$f_{\infty}(E\K|\K)\mid e_{\infty}(E|E\cap \K)$ y $e_{\infty}(
E|E\cap \K)\mid q-1$, de donde $f_{\infty}(E\K|\K)\mid q-1$.

Sea $d=f_{\infty}(E\K|\K)$. Tenemos que $E\K/\K$ es no ramificada
y el grado de inercia de $\S\K$ en $E\K/\K$ es $d\mid q-1$.
Puesto que $\g E/E$ es no ramificada y adem\'as
los elementos de  $\S E$ se descomponen totalmente en
$\g E/E$, lo mismo sucede en $\g E\K/E\K$. De esta forma
obtenemos que $\g E \K/\K$ es una extensi\'on no ramificada
y el grado de inercia de $\S \K$ es $d$.

Sea $H$ el grupo de descomposici\'on de cualquier primo
en $\S\K$ en $\g E\K/\K$ y sea 
$H_1:=H|_{\g E}$. Puesto que $\g E\cap \K=E\cap \K$, de
la correspondencia de Galois, obtenemos,
bajo el mapeo de restricci\'on, que
$H\cong H_1$, $|H|=|H_1|$
y $\g E^{H_1} \K=(\g E \K)^H$. Adem\'as, $H_1\subseteq I_{\infty}(
K(\Lambda_N)|K)\cong C_{q-1}$, donde $I_{\infty}$ denota
al grupo de inercia de $\p$. 
Por tanto $H$ es un grupo c\'iclico y $H\cong
H_1\cong C_d$.

Ya que $\S \K$ se descompone totalmente en 
$\g E^{H_1}\K/\K$, se sigue que
\begin{gather*}
\g E^{H_1} \K\subseteq \g \K.
\end{gather*}

\begin{turn}{270}
\xymatrix{
&K(\Lambda_N)\ar@{-}[rrrr]\ar@{-}[d]&&&&K(\Lambda_N)M\ar@{-}[d]\\
&C\ar@{--}[dd]\ar@{-}[rrrr]&&&&CM=
\g \K M\ar@{--}[dd]\ar@{-}[dl]\ar@/^2pc/@{-}[dddd]^{
\text{no ramificada}}\\
&&&&\g \K\ar@{-}[ddl]
|!{[ld];[ddd]}\hole
\ar@/^2pc/@{-}[ddddl]^{\text{no ramificada}}
|!{[dl];[dr]}\hole
|!{[dl];[ddd]}\hole
|!{[dddll];[ddd]}\hole\\
&\g E=\g E^{H_1} E\ar@{-}[ddr]
\ar@{-}[rr]\ar@{-}[dl]_{H_1=H|_{\g E}}&&
\g E \K\ar@{-}[d]_H\ar@{-}[rr]
|!{[ur];[d]}\hole
\ar@{-}[ddr]&&\g E M\ar@{-}[dd]\\
\g E^{H_1}\ar@{-}[rrr]
|!{[ru];[rrd]}\hole
\ar@{-}[rrrruu]
|!{[ru];[rrd]}\hole
|!{[ru];[rru]}\hole
\ar@{-}[ddr]&&&\g E^{H_1}\K=(\g E \K)^H\ar@{-}[dd]\\
&&E\ar@{-}[rr]
|!{[ru];[rd]}\hole
\ar@{-}[dl]&&E\K\ar@{-}[r]\ar@{-}[dl]&E\K M=EM=
\K M\ar@{-}[lld]\ar@{-}[dd]\\
&E\cap \K\ar@{-}[rr]\ar@{-}[d]&&\K\\
&K=K(\Lambda_N)\cap M\ar@{-}[rrrr]^{e_{\infty}=q^n,\quad f_{\infty}=m}&&&&M
}
\end{turn}

Sea $E_1:=E \g E^{H_1}\subseteq \g E$. Ahora bien $H_1\subseteq
I_{\infty}(E|E\cap \K)$, de tal forma que 
$\S {\g E^{H_1}}$ es totalmente ramificado en 
$\g E/\g E^{H_1}$. Por lo tanto $\S {E_1}$ es totalmente ramificado en 
$\g E/E_1$. Por otro lado $\S E$ se descompone totalmente
en $\g E/E$. Se sigue que $\S {E_1}$ 
se descompone totalmente en $\g E/E_1$. Esto es, $\S {E_1}$
se ramifica y se descompone totalmente en $\g E/E_1$. 
Por tanto
\begin{gather*}
\g E=E_1=E\g E^{H_1}.
\intertext{Se sigue que}
(\g E \K)^H=\g E^{H_1} \K\subseteq \g \K \quad \text{y} \quad 
E\g E^{H_1}=\g E.
\end{gather*}

Para probar la otra contenci\'on, definimos
$C:=\g \K M\cap K(\Lambda_N)$.
Tenemos
\begin{gather*}
E\subseteq EM=\K M\subseteq \g \K M,\quad E\subseteq K(\Lambda_N).
\intertext{Por tanto}
E\subseteq \g \K M\cap K(\Lambda_N)=C,
\quad \text{esto es}\quad E\subseteq C.
\end{gather*}

Adem\'as, $\g E^{H_1}\subseteq \g E^{H_1} \K
\subseteq \g \K\subseteq \g \K M$ y $\g E^{H_1}\subseteq \g E
\subseteq K(\Lambda_N)$. Se sigue que $\g E^{H_1}
\subseteq \g \K M\cap K(\Lambda_N)=C$.
Obtenemos $\g E^{H_1}\subseteq C$. Por tanto
\begin{gather}\label{Eq6.A}
\g E =E \g E^{H_1} \subseteq C.
\end{gather}

Puesto que $C=\g \K M\cap K(\Lambda_N)$, 
de la correspondencia de Galois obtenemos que
$CM=\g \K M$. Ahora, puesto que 
$\g \K/\K$ es no ramificada y $\S \K$
se descompone totalmente, se sigue que
\begin{gather}
CM=\g \K M/\K M=EM\quad \text{es no ramificada y}\nonumber\\
\S {\K M} \text{\ se descompone totalmente}. \label{Eq7.A}
\end{gather}

Ahora probaremos que $C/E$ es no ramificada. De (\ref{Eq7.A})
se sigue que $CM/\K M$ es no ramificada.
Ahora, en  $\K M=EM/E$, los \'unicos primos ramificados
son aquellos en $\S E$ y ellos tienen \'indice de ramificaci\'on
igual a $q^n$. 
Se sigue que los \'unicos primos ramificados en $CM/E$
son aquellos en $\S E$. Por tanto los \'unicos posibles primos
ramificados en $C/E$ son aquellos en $\S E$. Ahora bien
\begin{gather*}
e_{\infty}(C|E)\mid e_{\infty}(CM|E)=q^n
\quad\text{y}\quad e_{\infty}(C|E)\mid
e_{\infty}(K(\Lambda_N)|K)=q-1
\intertext{por lo que}
e_{\infty}(C|E)\mid \mcd(q^n,q-1)=1.
\end{gather*}
Por tanto $C/E$ es una extensi\'on no ramificada.

Por otro lado, siendo $\S E$ no ramificada en $C/E$, 
$\S E$ se descompone totalmente en
$C/E$ ya que $C\subseteq K(\Lambda_N)$. Se sigue que
$C\subseteq \g E$. De esto y de la Ecuaci\'on
(\ref{Eq6.A}), obtenemos
\begin{gather*}
C=\g E \quad\text{y}\quad \g E M=CM=\g \K M.
\end{gather*}

Tenemos que $\g E \K\subseteq \g E \g \K$. 
Puesto que $\g \K/\K$ no es ramificada y
$\S \K$ se descompone totalmente en $\g\K$,
lo mismo sucede en la extensi\'on $\g E\g \K/\g E \K$. 
En particular $h_{\infty}(\g E\g \K|\g E \K)=[\g E\g \K|\g E \K]$.

Ahora, en la extensi\'on $\g EM/\g E$, 
los \'unicos primos ramificados son aquellos en
$\S {\g E}$ y tenemos que $e_{\infty} (\g EM|\g E)=q^n$
y $f_{\infty}(\g EM|\g E)=m$. Por otro lado, se tiene que
$e_{\infty} (\g E|K)\mid q-1$, el cual es relativamente
primo a  $q$, $f_{\infty}(\g E|K)=1$,
$e_{\infty} (M|K)=q^n$ y $f_{\infty} (M|K)=m$.
\[
\xymatrix{
\g E\ar@{-}[rrr]\ar@{-}[d]&&&\g E M\ar@{-}[d]\\
\g E\cap M=K\ar@{-}[rrr]_{e_{\infty}=q^n, f_{\infty}=m}&&& M
}
\]

Sean $F_1$ y $F_2$ dos campos tales que 
$K\subseteq F_1\subseteq F_2
\subseteq M$. Sean $R_i=\g E F_i$, 
$i=1,2$. Puesto que $f_{\infty}(\g E|K)=1$
y $e_{\infty}(\g E|K)\mid q-1$, se sigue de la
correspondencia de Galois entre
$M/K$ y $\g E M/\g E$ que
$e_{\infty}(R_i|\g E)=e_{\infty}(F_i|K)$ y que
$f_{\infty}(R_i|\g E)=f_{\infty} (F_i|K)$, $i=1,2$. 
Por tanto $e_{\infty}(F_2|F_1)=e_{\infty}(R_2|R_1)$
y $f_{\infty}(F_2|F_1)=f_{\infty}(R_2|R_1)$.

Ya que $h_{\infty}(M|K)=1$, tenemos
$h_{\infty}(R_2|R_1)=1$. En particular
\begin{gather}
R_1\neq R_2\iff F_1\neq F_2\iff e_{\infty}
(F_2|F_1)>1\text{\ o\ }f_{\infty}(F_2|F_1)>1\nonumber\\
\iff e_{\infty}(R_2|R_1)>1
\text{\ o\ }f_{\infty}(R_2|R_1)>1. \label{Eq4.5'}
\end{gather}

Puesto que 
\[
\g E\subseteq \g E \K\subseteq \g E \g \K \subseteq \g \K M=\g E M,
\]
$\S {\g E \K}$ es no ramificada en 
$\g E\g \K/\g E \K$ y $\S{\g E \K}$
se descompone totalmente pues esto mismo es lo que sucede
en $\g\K/\K$ (Proposici\'on \ref{P10.4.1.Ram}).
Por lo tanto, se tiene que
$e_{\infty}(\g E\g \K|\g E \K)=1$ y
$f_{\infty}(\g E\g \K|\g E \K)=1$. De (\ref{Eq4.5'}), se sigue que
\[
\g E\g \K=\g E \K.
\]
Por tanto $\g \K\subseteq \g E\g \K=\g E \K$. 
Puesto que $\g E \K/\K$ es no ramificada,
si $D$ es el grupo de descomposici\'on de
$\S \K$ in $\g E \K/\K$, obtenemos que
$\g \K=(\g E \K)^D$. Finalmente, tenemos
\[
f_{\infty}(\g E \K|\K)=f_{\infty}(\g E \K|E\K)f_{\infty}(E\K|\K)=1\cdot d=d.
\]
Por tanto $D=H$ y $(\g E \K)^D=(\g E \K)^H=\g E^{H_1} \K=\g \K$. 

Hemos obtenido
\[
\g\K=\g E^{H_1}\K.
\]

Finalmente, falta probar que $\g E \K/\g \K$ 
y que $E\K/\K$ son extensiones de constantes.
Puesto que $\g \K M=\g E M$ y $\g E\g \K=\g E \K$, tenemos
\begin{gather*}
\g \K=(\g E \K)^H\subseteq \g E \K\subseteq \g 
E\g \K\subseteq \g E\g \K M= \g E M.
\intertext{Sean $F_1=\g \K\cap M$ y $F_2=\g E \K\cap M$. Tenemos
$d=[\g E\K:\g \K]=f_{\infty}(\g E \K|\g \K)=[F_2:F_1]=
e_{\infty}(F_2|F_1) f_{\infty}(F_2|F_1)
h_{\infty}(F_2|F_1)$. Puesto que $e_{\infty}(F_2|F_1)\mid q^n$ 
y $h_{\infty}(F_2|F_1)=1$,
se sigue que}
e_{\infty} (F_2|F_1)=e_{\infty}(\g E\K|\g \K)=1\quad\text{y}\quad
f_{\infty} (F_2|F_1)=f_{\infty}(\g E\K|\g \K)=d.
\end{gather*}
Esto es, $K\subseteq F_1\subseteq 
F_2\subseteq M$ y $e_{\infty}(F_2|F_1)=1$.

Sean $a$ y $b$ tales que $F_2\subseteq F_1 K_bL_a$.
Consideremos $A_i=F_i K_b\cap L_a$, $i=1,2$. 
Puesto que $e_{\infty}(F_2|F_1)=1$, notemos que
tenemos $e_{\infty}(A_2|A_1)=1$ puesto que
$F_iK_b=A_iK_b/A_i$, $i=1,2$, son extensiones de constantes,
esto es
\begin{align*}
e_{\infty}(A_2|A_1)&=\frac{e_{\infty}(A_1K_b|A_1)e_{\infty}(
A_2K_b|A_1K_b)}{e_{\infty}(A_2K_b|A_2)}=
\frac{1\cdot e_{\infty}(A_2K_b|A_1K_b)}{1}\\
&= e_{\infty}(A_2K_b|A_1K_b)=e_{\infty}(F_2K_b|F_1K_b)\mid
e_{\infty}(F_2|F_1)=1.
\end{align*}
\begin{gather*}
\xymatrix{
L_a\ar@{-}[rrr]\ar@{-}[d]&&& L_a K_b\ar@{-}[d]\\
A_2\ar@{-}[rrr]\ar@{-}[dd]&&& F_2 K_b=A_2 K_b\ar@{-}[dd]\\
&&F_2\ar@{-}[ur]\ar@{-}[ddl]\\
A_1\ar@{-}[rrr]
|!{[rru];[rd]}\hole
\ar@{-}[dd]&&& F_1 K_b=A_1 K_b\ar@{-}[dd]\\
&F_1\ar@{-}[rru]\ar@{-}[dl]\\
K\ar@{-}[rrr]&&&K_b
}
\\
e_{\infty}(F_2K_b|F_1K_b)=e_{\infty}(F_2|F_1)=e_{\infty}(A_2|A_1)=1.
\end{gather*}

Puesto que $\p$ es totalmente ramificado en $L_a/K$,
se sigue que $A_1=A_2$. Por lo tanto
$F_2 K_b=F_1K_b$ y $F_2/F_1$
es una extensi\'on de constantes. 

Se tiene $F_1=\g\K\cap M$. Consideremos
$\g\K\subseteq \g E\K\subseteq \g\K M=\g E M$:
\[
\xymatrix{
\g\K\ar@{-}[r]\ar@{-}[d]&\g E \K\ar@{-}[r]\ar@{-}[d]&
\K M=\g EM\ar@{-}[d]\\
F_1\ar@{-}[r]&F_2\ar@{-}[r] &M}
\]
Por tanto $\g\K\subseteq F_2\g\K=\g E \K$.
Se sigue que $\g E \K/\g \K$ es una extensi\'on de
constantes de grado $[\g E\K:\g \K]=|H|=d$. 

La demostraci\'on
de que $E\K/E^{H_1}\K$ es una extensi\'on de constantes
es completamente similar.

Esto termina la demostraci\'on del teorema. $\fin$
\end{proof}

Para el caso particular de una $p$--extensi\'on abeliana
finita $\K/K$, se tiene, por un lado, que $d\mid q-1$ y, por otro
lado, que $d\mid [E\K:\K]$. Puesto que
$\K$ es una $p$--extensi\'on, obtenemos de 
(\ref{Eq4.A}), que $E/K$ es tambi\'en una
$p$--extensi\'on. Finalmente, puesto que
$\Gal(E\K/K)\to \Gal(E/K)\times 
\Gal(\K/K)$, $\sigma\mapsto (\sigma|_{E},\sigma|_{\K})$
es inyectiva, se sigue que $E\K/K$ es tambi\'en
una $p$--extensi\'on. Por tanto $d\mid p^a$ para alguna
$a$. Se sigue que $d=1$. Hemos probado

\begin{teorema}\label{T2.2.A} Con las notaciones 
del Teorema {\rm{\ref{T2.1.A}}},
sea $\K/K$ una $p$--extensi\'on abeliana finita. Sea
$E:=\K M\cap K(\Lambda_N)$. Entonces
\[
\g \K=\g E \K
\]
y $\g \K/K$ es una $p$--extensi\'on.
\end{teorema}

\begin{proof} La \'ultima afirmaci\'on se sigue del hecho de que $\g E$
es tambi\'en una $p$--extensi\'on. $\fin$
\end{proof}

\begin{ejemplo}\label{Ej2.3.A}
Sea $P\in R_T^+$ y sea $F_P$ cualquier campo tal que
$K\subsetneqq F_P\subseteq \cicl P{}$. Entonces $F_P/K$
es una extensi\'on c\'iclica. Puesto que el \'unico primo
finito ramificado en $F_P/K$ es $P$, el grupo de caracteres
de Dirichlet asociado a $F_P/K$ es el generado por un
caracter $\chi$ de conductor $P$. En particular si $X$
es el grupo de caracteres de Dirichlet asociado a $F_P$
es $X=\langle \chi\rangle =X_P$. Por el Teorema \ref{T12*.3.B}
obtenemos que $\g {(F_{P})}=F_P$.

M\'as generalmente, sean $P_1,\ldots,P_s\in R_T^+$
primos distintos. Sea $F:=\prod_{i=1}^s F_{P_i}$ donde
$K\subseteq F_{P_i}\subseteq \cicl {P_i}{}$, $1\leq i\leq s$.
Entonces si $\chi_{P_i}=\chi_i$ es el caracter asociado a $F_{P_i}$,
$1\leq i\leq s$. Entonces si $X$ es el grupo de caracteres
asociado a $F$, se tiene $X=\langle \chi_i\mid 1\leq i\leq s\rangle
=\prod_{i=1}^s X_{P_i}$. Por el Teorema \ref{T12*.3.B} se
sigue que
\[
\g F= \prod_{i=1}^s \g {(F_{P_i})}=\prod_{i=1}^s F_{P_i}=F.
\]
\end{ejemplo}

\begin{ejemplo}\label{Ej2.4.A}
Sea $\K/K$ una extensi\'on abeliana finita, la cual es moderadamente
ramificada y tal que $\g \K=\K$. Sean $P_1,\ldots, P_s$ los
primos finitos de $K$ ramificados en $\K$. Sea $e_i$ el
\'indice de ramificaci\'on de $P_i$ en $\K/K$, $1\leq i\leq s$.
Se tiene que $e_i|q-1$. Para cada $1\leq i\leq s$, sea 
$F_i$ el \'unico subcampo de $\cicl {P_i}{}$ de grado $e_i$ sobre $K$.
Sea $F:=\prod_{i=1}^s F_i$.

Por el Ejemplo \ref{Ej2.3.A} se tiene que $\g F=F$. Sea $n\in{\ma N}$ tal que
$\K\subseteq \cicl D{} K_n$ donde $D=\prod_{i=1}^s P_i$. Sea $E:=
\K_n\cap \cicl D{}$. Entonces los primos ramificados en $E/K$ son
precisamente $P_1,\ldots,P_s$ con \'indices de ramificaci\'on 
$e_1,\ldots,e_s$ respectivamente. Por el Teorema \ref{T2.1.A}
tenemos que $\K=\g \K=E^{H_1}\K$. Se sigue que $E^{H_1}
\subseteq \K$.

Nuevamente por el Teorema \ref{T2.1.A}, $E\K/\K$ es una extensi\'on
de constantes. Adem\'as, $E_n=\K_n$. Finalmente, veamos que
$F\K =F K_u$ para alg\'un $u\in{\ma N}$. De hecho tenemos que
$F\subseteq F\K\subseteq F_n$, de donde se sigue nuestra
afirmaci\'on.
\end{ejemplo}

\subsection{Conductor de constantes}\label{S14.5.1}

Dada una extensi\'on abeliana finita $\K$ de $K$, por el
Teorema de Kronecker--Weber, existen $n,m\in{\ma N}$ y 
$N\in R_T$ tales que $\K\subseteq {_nK(\Lambda_N)_m}$.
Los m\'inimos $n$ y $N$ que satisfacen esta contenci\'on
est\'an dados por teor\'ia de campos de clase y son simplemente
los conductores locales de la extensi\'on $\K/K$: $n$
para $\p$ y $N$ para los primos finitos.

En esta parte determinaremos el m\'inimo $m$, el cual
como veremos, est\'a
relacionado con el n\'umero $d$ dado en el Teorema \ref{T2.1.A}.
El n\'umero $m$ se llamar\'a el {\em conductor de
constantes\index{conductor de constantes}} de la
extensi\'on abeliana $\K/K$.

Primero sean $n\in{\ma N}$, $N\in R_T$ y $m\in{\ma N}$ 
m\'inimo (que
en principio depende de $n$ y $N$) tales que 
$\K\subseteq {_nK(\Lambda_N)_m}$. Consideremos el siguiente
cuadro de extensiones de Galois
\[
\xymatrix{
{_n K(\Lambda_N)}\ar@{-}[rr]\ar@{-}[dd]&&U= {_nK(\Lambda_N)}\K
\ar@{-}[r]\ar@{-}[dl]\ar@{--}[dd] & {_nK(\Lambda_N)_m}\ar@{-}[dd]\\
&\K\ar@{-}[dl] \\ K\ar@{-}[rr] && K_{m'}\ar@{-}[r] &K_m
}
\]
Esto es, sean $U={_nK(\Lambda_N)} \K$ y $K_{m'}=U\cap K_m$.
Por la correspondencia de Galois, tenemos que $U={_nK(\Lambda_N)}
K_{m'}={_nK(\Lambda_N)_{m'}}={_nK(\Lambda_N)} \K\supseteq \K$.

Puesto que $m$ es m\'inimo, se tiene que $m'=m$. En otras palabras,
$m$ est\'a determinado por la igualdad 
\begin{gather}\label{conductor de constantes}
{_nK(\Lambda_N)}\K={_nK(\Lambda_N)_m}.
\end{gather}

Veamos ahora que $m$ es independiente de $n$ y de $N$. Sean
$n_i\in{\ma N}$, $N_i\in R_T$ y $m_i\in{\ma N}$ m\'inimo
tales que $\K\subseteq {_{n_i} K(\Lambda_{N_i})_{m_i}}$,
$i=1,2$.

Sean $n_0:=\max \{n_1,n_2\}$, $N_0=\mcm [N_1,N_2]$ y $m_0\in
{\ma N}$ m\'inimo tal que $\K\subseteq {_{n_0}K(\Lambda_{N_0})_{m_0}}$.
Entonces por (\ref{conductor de constantes}) obtenemos
\begin{align*}
{_{n_0}K(\Lambda_{N_0})}\K&= L_{n_0}({_{n_i}K(\Lambda_{N_i})}
K(\Lambda_{N_0}))\K=L_{n_0}\big({_{n_i}K(\Lambda_{N_i})}\K\big)
K(\Lambda_{N_0})\\
&=L_{n_0} ({_{n_i}K(\Lambda_{N_i})_{m_i}}) K(
\Lambda_{N_0})={_{n_0}K(\Lambda_{N_0})_{m_i}},\\
{_{n_0}K(\Lambda_{N_0})}\K&={_{n_0}K(\Lambda_{N_0})_{m_0}}.
\end{align*}
Por tanto $m_1=m_2=m_0$.

De esta forma consideramos $\K\subseteq {_n\cicl N{}_m}$ con
$m$ m\'inimo.
Sea $F:= \K\cap {_n\cicl N{}}$ y consideremos el cuadro de Galois
(ver (\ref{conductor de constantes}))
\[
\xymatrix{
{_n\cicl N{}}\ar@{-}[rr]^{m\phantom{xxxxx}}\ar@{-}[d] && {_n\cicl N{}_m}
 ={_n\cicl N{}} \K
\ar@{-}[d]\\ F\ar@{-}[rr]^m\ar@{-}[d]&&\K \ar@{-}[dll]\\ K
}
\]

Sea $t$ el grado de $\S \K$ en $\K$. Esto es, $t=f_{\infty}(\K|K)$. Se tiene que
\begin{gather*}
e_{\infty}({_n\cicl N{}_m}|{_n\cicl N{}})=1,\quad 
f_{\infty}({_n\cicl N{}_m}|{_n\cicl N{}})=m.
\intertext{En particular}
\{1\}=I_{\infty}({_n\cicl N{}_m}|{_n\cicl N{}})\subseteq I_{\infty}(\K|F),\\
C_m\cong D_{\infty}({_n\cicl N{}_m}|{_n\cicl N{}})\subseteq
D_{\infty}(\K|F).
\end{gather*}

Puesto que $[\K:F]=m$ y $m\leq |D_{\infty}(\K|F)|\leq [\K:F]=m$, 
se sigue que $|D_{\infty}(\K|F)|=m$ y que $D_{\infty}(\K|F)\cong C_m$. En
particular $h_{\infty}(\K|F)=1$ y $h_{\infty}({_n\cicl N{}_m}|{_n\cicl N{}})=1$.

Por otro lado, se tiene
\begin{gather*}
t=f_{\infty}(\K|K)=f_{\infty}(\K|F)f_{\infty}(F|K) =f_{\infty}(\K|F)\cdot 1=f_{\infty}(\K|F),
\end{gather*}
esto es, $f_{\infty}(\K|F)=t$. Adem\'as
\begin{gather*}
e_{\infty}(\K|F) f_{\infty}(\K|F)h_{\infty}(\K|F)=e_{\infty}(\K|F)\cdot t\cdot 1=m,
\end{gather*}
por lo cual $e_{\infty}(\K|F)=\frac{m}{t}$.  Se sigue que
\begin{gather}\label{Eq14.5.7}
m=[\K:F]=f_{\infty}(\K|F)e_{\infty}(\K|F)=t e_{\infty}(\K|F)
=t\frac{e_{\infty}(\K|k)}{e_{\infty}(F|k)}.
\end{gather}

Ahora investigaremos la relaci\'on entre $m$ y $d=f_{\infty}(\g E \K|\g \K)$ dada
en el Teorema \ref{T2.1.A}. Recordemos que $M=L_nK_m$, que
$E=\K M\cap \cicl N{}$ y que $EM=\K M$.
Tenemos 
\[
\g E\subseteq \g E \K \subseteq \g E \K  L_n\subseteq \g E \K  M =\g E EM=\g E M.
\]

Sean $A:=\g E\K \cap M$ y $B:=\g E \K  L_n \cap M$. 
De la correspondencia de Galois, tenemos que
$\g E \K =\g E A$ y que $\g E \K  L_n=\g E B$. 
\[
\xymatrix{
\g E\ar@{-}[d]\ar@{-}[r]&\g E \K \ar@{-}[d]\ar@{-}[r]&\g E \K L_n\ar@{-}[d]
\ar@{-}[r]&\g E M\ar@{-}[d]\\
k\ar@{-}[r]&A\ar@{-}[r]&B\ar@{-}[r]&M
}
\]
Tenemos $L_n\subseteq
\g E \K  L_n\cap M=B\subseteq M=L_nk_m$.  Por lo 
tanto $B/L_n$ es una extension de constantes.
Digamos $B=L_n k_{m'}$ con $m'|m$. De la
correspondencia de Galois, obtenemos
\[
\K \subseteq \g E  \K  L_n=\g E B=\g E L_n k_{m'}\subseteq 
\cicl N{} L_n k_{m'}={_n\cicl n{}}_{m'}.
\]
Puesto que $m$ es el m\'inimo, 
se sigue que $m'=m$, que $B=M$ y que $\g E \K L_n=\g EM$.

Ahora bien, $\g E(A L_n)=(\g E A)L_n=(\g E \K )L_n=\g E M$. 
De la correspondencia de Galois se obtiene que $A L_n=M$. 
Consideramos el siguiente cuadro de Galois:
\[
\xymatrix{
L_n\ar@{-}[d]\ar@{-}[r]&AL_n=M=L_nk_m\ar@{-}[d]\\
A\cap L_n\ar@{-}[r]&A
}
\]

Tenemos que $f_{\infty}(AL_n|L_n)=f_{\infty}(M|L_n)=m$
y que $e_{\infty}(AL_n|L_n)= e_{\infty}(M|L_n)=1$. Por tanto
\begin{align*}
\{1\}&=I_{\infty}(AL_n|L_n)\subseteq I_{\infty}(A|A\cap L_n)\quad\text{y}\\
C_m&\cong D_{\infty}(AL_n|L_n)\subseteq D_{\infty}(A|A\cap L_n).
\end{align*}

Puesto que $[A:A\cap L_n]=[M:L_n]=m$ se sigue que
 $D_{\infty}(A|A\cap L_n)
\cong C_m$, que $e_{\infty}(A|A\cap L_n)=1$ 
y que $f_{\infty}(A|A\cap L_n)=m$. Por tanto
$f_{\infty}(\g E\K |k)=f_{\infty}(\g E\K |\g \K)
f_{\infty}(\g \K |\K)f_{\infty}(\K |k)=d\cdot 1\cdot t= dt=td$. As\'i,
\begin{gather*}
f_{\infty}(\g E M|\g E \K)=\frac{f_{\infty}(\g EM|k)}
{f_{\infty}(\g E \K |k)}=\frac{m}{td}.
\intertext{Finalmente}
\begin{align*}
\frac{m}{td}&=f_{\infty}(\g EM|\g E \K)|[\g EM:\g E\K ]=[M:A]\\
&=[L_n:A\cap L_n]|[L_n:k]=q^n.
\end{align*}
\intertext{Se sigue que}
m=td p^s
\end{gather*}
para alguna $s\in{\ma N}\cup \{0\}$.

En adici\'on, se tiene
$f_{\infty}(\K _m|\K)=\frac{m}{t}=e_{\infty}(\K |F)$. Notemos
que $td=f_{\infty}(\K|K)f_{\infty}(E\K|\K)=f_{\infty}(E\K|K)$.

Hemos obtenido

\begin{teorema}[Conductor de constantes 1]
\label{T2.1.AA} Sea $\K$ una extensi\'on
abeliana finita de $k$. Sean
$n,m\in{\ma N}$ y $N\in R_T$ tales que
$\K \subseteq {_n\cicl N{}_m}$ y tal que $m$ es
el m\'inimo con esta propiedad. Entonces $m$
es independiente de $n$ y $N$. Sea
$t=f_{\infty}( \K |k)$ el grado de los primos
infinitos de $\K $. Sean $M=L_nK_m$, $E=\K M\cap \cicl N{}$,
$F=\K \cap {_n\cicl N{}}_m$ y 
$d=f_{\infty}(E\K |\K)=f_{\infty}(\g E \K |\g \K)$.
Entonces ${_n\cicl N{}} \K ={_n\cicl N{}}_m$ y
\begin{gather*}
m=[\K :F]=t e_{\infty}(\K |F)=td p^s=f_{\infty}(E\K|K) p^s
\end{gather*}
para alguna $s\geq 0$. En particular
$$
e_{\infty}(\K |F)=dp^s=f_{\infty}(\K _m|\K).
\eqno{\fin}
$$
\end{teorema}

\begin{observacion}\label{R2.2AA}
Cuando $p\nmid \frac{m}{t}$, en particular cuando $\K |k$
es moderadamente ramificada en
$\p$, tenemos $s=0$ y 
$m=td$. En el caso general, se puede tener $s\geq 1$.
\end{observacion}

\begin{ejemplo}\label{Ej2.3AA}
Sea $p$ un primo y sea $q=p$. Definimos $X:=1/T$. Tenemos
que  $L_1:=\cicl X2^{{\ma F}_q^{\ast}}$ y que $[L_1:K]
=q=p$. Por tanto la extensi\'on $L_1/K$ es una extensi\'on de 
Artin--Schreier. Para nuestros fines, no es necesario dar
$L_1$ expl\'icitamente pero, por conveniencia del lector,
daremos un generador de $L_1$.

Sea $\lambda$ un generador de $\Lambda_{X^2}$. Entonces
$\lambda^{p-1}$ is a generator of $\cicl X2^+=L_1$ 
(ver Teorema \ref{T6.2.30}).
Ahora $\lambda$ es una ra\'iz del polinomio ciclot\'omico
$\Psi_{X^2}(u)$. Se tiene que
$\Psi_{X^2}(u)=\Psi_X(u^X)$ donde $u^P$
denotes la acci\'on de Carlitz. Puesto que $\Psi_X(u)=u^P/u=u^{q-1}+X$,
se obtiene que $\Psi_{X^2}(\lambda)=(\lambda^q+X\lambda)^{q-1}+X$.
Sea $\mu:=\lambda^{q-1}$ y $\xi:=\mu+X$. Entonces
\begin{gather*}
\xi^q-X\xi^{q-1}+X=0.
\intertext{Finalmente, si $\delta:=1/\xi$, entonces $L_1=K(\delta)$ con}
\delta^q-\delta=-1/X=-T,\quad \delta=\frac{T}{T\lambda^{q-1}+1}.
\end{gather*}

Consideremos a una soluci\'on $\alpha$ de
la ecuaci\'on $y^p-y=1$. Entonces ${\ma F}_p(\alpha)=
{\ma F}_{p^p}$. Sea $K_p={\ma F}_p(\alpha)(T)={\ma F}_{p^p}(T)$
y $L_1K_p=K(\alpha, \delta)$. Las $p+1$ extensiones $\K /K$ de
grado $p$ sobre $K$ tales que $K\subseteq \K \subseteq L_1K_p$
son $\{K(\alpha+i\delta)\}_{i=0}^{p-1}$ y $L_1$. Sea $\K :=K(
\alpha+\delta)$. Entonces $\K \neq K_p$ y $\K \neq L_1$. Se sigue que
$\K =K(z)$ con $z^p-z=1-T$.

Sea $N\in R_T$ un polinomio no cero arbitrario. 
Entonces $\K \subseteq L_1K_p\subseteq
{_1\cicl N{}}_p$ y $\K \nsubseteq {_1\cicl N{}}_1$. Por lo tanto 
$m=p$. Sea $M=L_1K_p$. 
Tenemos $f_{\infty}(\K |K)=1$, $e_{\infty}(\K |K)=p$. Tambi\'en se tiene
$E:=\K M\cap \cicl N{}=M\cap \cicl N{}=K$. Se sigue que $\g E=K$ y
$\g \K =\g E \K =\K $. Se obtiene que 
$E\K =\K $ y que $f_{\infty}(E\K |\K )=d=1$. Por lo tanto
$td=1\neq m=p$. En este ejemplo tenemos $s=1$.
\end{ejemplo}

Calcularemos $m$ de otra forma. Primero, usando una demostraci\'on
totalmente paralela a la del Teorema \ref{T2.1.A} obtenemos

\begin{teorema}\label{T2.4AA} Sea $\K /K$
una extensi\'on abeliana finita. Sea
\[
R:=\K_m\cap{_n\cicl N{}}.
\]

Entonces
\[
\g \K =\g {R^{{\mathcal H}_1}} \K = (\g R \K )^{\mathcal H},
\]
donde ${\mathcal H}$ es el grupo de descomposici\'on de cualquier
primo de $\S \K $
en $\g R \K /\K $, ${\mathcal H}_1:= {\mathcal H}|_{\g R}$ 
y ${\mathcal H}_2:={\mathcal H}_1|_R$.

Sea $d^{\ast}:= f_{\infty}(R\K /\K )$. Tenemos ${\mathcal H}\cong 
{\mathcal H}_1\cong {\mathcal H}_2\cong
C_{d^{\ast}}$ y $d^{\ast}|q-1$. Tambi\'en se tiene que
las extensiones $\g R \K /\g \K $ y $R\K 
/R^{{\mathcal H}_2}\K $ son extensiones de constantes de grado
$d^{\ast}$. Finalmente, el campo de constantes de
$\g \K $ es ${\ma F}_{q^t}$, 
donde $t$ es el grado de $\S \K $ en $\K $. $\fin$
\end{teorema}

Consideremos nuevamente a
$F= \K \cap {_n\cicl N{}}$ y los siguientes cuadros de Galois:
\begin{gather*}
\xymatrix{
{_n\cicl N{}}\ar@{-}[rr]\ar@{-}[d]&&{_n\cicl N{}_m}\ar@{-}[d]\\
R\ar@{-}[rr]\ar@{-}[dd]\ar@{-}[r]&& \K_m=R_m\ar@{-}[dd]\ar@{-}[dl]\\
& \K \ar@{-}[dl]\\ K\ar@{-}[rr]&&K_m
}
\end{gather*}

\begin{gather*}
\xymatrix{
{_n\cicl N{}} \ar@{-}[rr]\ar@{-}[d]&&{_n\cicl N{}} \K ={_n \cicl N{}_m}
\ar@{-}[d]\\ 
C\ar@{-}[rr]\ar@{-}[d]&&R_m=\K_m\ar@{-}[d]\\
R= \K_m\cap{_n\cicl N{}}\ar@{-}[rr]\ar@{-}[d]&&R \K \ar@{-}[d]\\
F= \K \cap {_n\cicl N{}}\ar@{-}[rr]&& \K 
}
\end{gather*}

Ya que $R= \K_m\cap {_n\cicl N{}}$, se sigue que $ \K_m=R_m$.
Ahora bien, $ \K ,R\subseteq R\K  \subseteq  \K_m=R_m$.

Definimos $C:= \K_m\cap {_n\cicl N{}}$. Entonces $C=R$ y,
de la correspondencia de Galois, obtenemos 
que $R \K =R_m= \K_m$.

Se sigue que el campo de constantes de
$R \K $ es ${\ma F}_{q^m}$.
El campo de constantes de $R\g \K $ tambi\'en es 
${\ma F}_{q^m}$.

El campo de constantes de $\g \K $ es ${\ma F}_{q^t}$.
Por otro lado, tenemos que 
$R\g \K /\g{R^{{\mathcal H}_1}} \K =\g \K $ 
es una extensi\'on de constantes de grado
$d^{\ast}=|{\mathcal H}_1|$. Por tanto, el campo de constantes de
$R\g \K $ es ${\ma F}_{q^{td^{\ast}}}$. Se sigue que $td^{\ast}=m$.

Notemos que $td^{\ast}=f_{\infty}(\K|K)f_{\infty}(R\K|\K)=f_{\infty}(R\K|K)$.
Hemos obtenido

\begin{teorema}[Conductor of constantes 2]
\label{T2.5.AA} Sea $ \K $ una extensi\'on abeliana finita de
$K$. Sean $n,m\in{\ma N}$ y $N\in R_T$ tales que
$\K \subseteq {_n\cicl N{}_m}$ y tal que $m$ es
el m\'inimo con esta propiedad. Sea
$t=f_{\infty}( \K |K)=f_{\infty}(\K |F)$ el grado de los 
primos infinitos de $\K $. Sean $R=
\K_m\cap {_n\cicl N{}}$ y $d^{\ast}=f_{\infty}(R\K |\K )$.
Entonces
\begin{gather*}
m=t e_{\infty}( \K |F)=td^{\ast}=f_{\infty}(R\K|K).
\end{gather*}
En particular
$$
d^{\ast}=f_{\infty}(R\K |\K )=e_{\infty}( \K |F).
\eqno{\fin}
$$
\end{teorema}

\begin{observacion}\label{R2.6.AA} Cuando una extensi\'on
abeliana $\K/K$ es moderadamente ramificada, entonces
el conductor de constantes $m$ de la extensi\'on
satisface $m=f_{\infty}(\K/K) d$. En particular, cuando
$d=1$, se tiene $m=f_{\infty}(\K/K)$.
\end{observacion}

\section[Kummer y $p$--extensiones]{Descripci\'on expl\'icita
de campos de g\'eneros de
extensiones de Kummer y de $p$--extensiones contenidas
en un campo de funciones ciclot\'omico}\label{S4.A}

Para una extensi\'on abeliana finita $\K/K$, la descripci\'on
de $\g \K$ depende de la descripci\'on de $\g E$
(Theorem \ref{T2.1.A}). En esta secci\'on presentamos
descripciones expl\'icitas de algunos subcampos ciclot\'omicos
$E$ con el fin de hallar $\g E$. Aqu\'i $\K$ denota
un campo $K\subseteq \K\subseteq K(\Lambda_N)$
para alguna $N\in R_T$ y $K={\ma F}_q(T)$.

\begin{observacion}\label{R4.0.A}
Sea $K\subseteq \K\subseteq K(\Lambda_N)$ y sea $X$ el
grupo de caracteres de Dirichlet asociado a $\K$. Sea $L$
el campo asociado a $\prod_{P\in R_T^+}X_P$ entonces
\[
\g \K=L^D,
\]
donde $D$ es el grupo de descomposici\'on de cualquier primo
${\eu p}\in\S \K$ en $L/\K$.
\end{observacion}

\begin{proposicion}\label{P4.1.A}
Con las notaciones anteriores, sea $X$ el grupo de 
caracteres de Dirichlet correspondiente a $\K$.
Fijemos $P\in R_T^+$. Sea $Y$ un grupo de caracteres de
Dirichlet tal que $Y=Y_P$, esto es, par cualquier
$\chi\in Y$, el conductor de $\chi$ es una potencia de $P$:
${\mc F}_{\chi}=P^{\alpha_{\chi}}$ para alguna
$\alpha_{\chi}\in {\ma N} \cup\{0\}$. 
Sea $L$ el campo asociado a $\langle X,Y\rangle$,
esto es, si $F$ es el campo asociado a $Y$, entonces $L=\K F$. 
Supongamos que $\K F/\K$ es no ramificada en $P$. 
Entonces $Y\subseteq X_P$.
\end{proposicion}

\begin{proof} Tenemos que 
\[
|\langle X,Y\rangle_P|= e_P(\K F|K)=e_P(\K F|\K)
e_P(\K|K)=e_P(\K|K)=|X_P|.
\]
Puesto que $X_P\subseteq \langle
X,Y\rangle_P$, se sigue que $X_P=\langle X,Y\rangle_P$.
Debido a que $Y_P\subseteq \langle X,Y\rangle_P$,
el resultado se sigue. $\fin$
\end{proof}

\begin{corolario}\label{C4.2.A}
Si $|Y|=|X_P|$, entonces $Y=X_P$. $\fin$
\end{corolario}

Aplicaremos la Proposici\'on \ref{P4.1.A} 
tanto a extensiones de Kummer de $K$ como a $p$--extensiones
abelianas finitas de $K$.

\subsection{Extensiones de Kummer}\label{S4.1.A}

Sea $\K=K(\sqrt[d]{\gamma D})$ una extensi\'on de Kummer
con $\K\subseteq K(\Lambda_D)$. Podemos suponer
sin p\'erdida de generalidad que $\gamma=(-1)^{\deg D}$ y
que $D\in R_T$ es un polinomio m\'onico (ver 
Proposici\'on \ref{P5.1.1}).
Suponemos de $D$ est\'a libre de $d$--potencias y que
$d\mid q-1$. Digamos que
\[
D= P_1^{\alpha_1}\cdots P_r^{\alpha_r}, \quad r\geq 1,\quad
1\leq \alpha_i\leq d-1,\quad 1\leq i\leq r.
\]
Pongamos $d_i:=\mcd(\alpha_i,d)$. Entonces
$\mcd\big(\frac{\alpha_i}{d_i},\frac{d}{d_i}\big)=1$.

Sea ${\eu p}_i$ un primo en $\K$ dividiendo a
$P_i\in R_T^+$. Sea $\beta:=\sqrt[d]{\gamma D}$,
$\beta^d=\gamma D=
\gamma P_1^{\alpha_1}\cdots P_r^{\alpha_r}$. Se tiene que
$e_i:=e_{P_i}(\K|K) = d/d_i$ (ver 
Teorema \ref{TRam1}).

Sea $F_i=K\Big(\sqrt[d/d_i]{(-1)^{\deg P_i^{
\alpha_i/d_i}}P_i^{\alpha_i/d_i}}\Big)$.
Escribamos $\gamma_i=(-1)^{\deg P_i^{\alpha_i/d_i}}$. 
Sea $X$ el grupo de caracteres de Dirichlet asociados 
a $\K$. De hecho, $X$ es un grupo c\'iclico de orden $d$ y
sea $X=\langle \chi\rangle$. Sea $Y$ el grupo de caracteres
de Dirichlet asociado a $F_i$. Entonces $Y=Y_{P_i}$ y
$|Y_{P_i}|=e_{P_i}(F_i|K)=d/d_i$ 
ya que $\mcd \big(\frac{d}{d_i},\frac{\alpha_i}{d_i}\big)=1$,
y $|X_{P_i}|=e_{P_i}(\K|K)=d/d_i=|Y_{P_i}|$. 

Veremos que $\K F_i/\K$ es no ramificado en $P_i$. Tenemos
\begin{align*}
\K F_i&=K\Big(\sqrt[d]{\gamma D}, \sqrt[d/d_i]{\gamma_i
P_i^{\alpha_i/d_i}}\Big)=K\Big(\sqrt[d]{\gamma D},\sqrt[d]{
\gamma_i^{d_i}P_i^{\alpha_i}}\Big)\\
&= \K\Big(\sqrt[d]{(-1)^{\deg P_i^{\alpha_i}}P_i^{\alpha_i}}\Big)=
\K\Big(\sqrt[d]{\frac{\gamma D}{\gamma_i^{d_i} P_i^{\alpha_i}}}\Big),
\end{align*}
y $P_i\nmid \frac {D}{P_i^{\alpha_i}}$. Por tanto $P_i$ es no
ramificado en $\K F_i/\K$ (Teorema \ref{TRam1}).
Por tanto $Y_{P_i}=X_{P_i}=Y$.

Se sigue que el campo asociado al grupo
$\prod_PX_P$ es $K(\xi_1,\ldots
\xi_r)$ donde $\xi=\sqrt[d/d_i]{\gamma_i P_i^{\alpha_i/d}}$.

Hemos probado

\begin{teorema}\label{T4.3.A-1}
Sea $X$ el grupo de caracteres de Dirichlet asociados
a $\K=K\big(\sqrt[d]{\gamma D}\big)$ con
$d\mid q-1$, $D\in R_T$ libre de $d$--potencias,
$D=P_1^{\alpha_1}\cdots P_r^{\alpha_r}$, $r\geq 1$, $1\leq \alpha_i
\leq d-1$, $1\leq i\leq r$, $\gamma =(-1)^{\deg D}$. 
Sea $d_i=\mcd(\alpha_i,d)$,
$1\leq i\leq r$. Entonces el campo asociado a $\prod_PX_P=
\prod_{i=1}^rX_{P_i}$ es $L=K(\xi_1,\ldots,\xi_r)$ donde $\xi_i
=\sqrt[d/d_i]{\gamma_i P_i^{\alpha_i/d_i}}$ y $\gamma_i=(-1)^{
\deg P_i^{\alpha_i/d_i}}$. Esto es, 
$$
L=K\big(\sqrt[d]{(-1)^{\deg P_1^{\alpha_1}}P_1^{\alpha_1},\ldots,
(-1)^{\deg P_r^{\alpha_r}}P_r^{\alpha_r}}\big). \eqno{\fin}
$$
\end{teorema}

\subsection{$p$--extensiones}\label{S4.2.A}

Ahora consideramos $\K=K(\vec y)$ donde 
\[
\vec y^{p^u}\Witt -\vec y=
\vec \delta_1\Witt +\cdots \Witt + \vec \delta_r
\]
con
$\vec \delta_i=(\delta_{i,1},\ldots,\delta_{i,v})$ 
para algunas $v\in{\ma N}$, $\delta_{i,j}=
\frac{Q_{i,j}}{P_i^{e_{i,j}}}$, $e_{i,j}\geq 0$, $Q_{i,j}\in R_T$.
Aqu\'i estamos suponiendo que
${\ma F}_{p^u}\subseteq {\ma F}_q$.

Sea $X$ el grupo de caracteres asociado a $\K$
donde suponemos que $\K\subseteq K(\Lambda_N)$
para alguna $N\in R_T$. De acuerdo a 
Schmid \cite{Sch36}, el \'indice de ramificaci\'on de
$P_i$ en $\K/K$ es determinado por el primer \'indice 
$j$ tal que podemos escribir $\delta_{i,j}=
\frac{Q_{i,j}}{P_i^{e_{i,j}}}$ con $\mcd(Q_{i,j},P_i)=1$, $e_{i,j}>0$
y $\mcd(e_{i,j},p)=1$ (Observaci\'on \ref{O9'.8.28}).

En otras palabras, el \'indice de ramificaci\'on de
$P_i$ en $\K/K$ depende \'unicamente de
$\vec \delta_i$ y no de $\vec \delta_1,\ldots, \vec\delta_{i-1},
\vec \delta_{i+1},\ldots, \vec \delta_r$.
Por lo tanto, si $Y$ es el grupo de caracteres de
Dirichlet asociado a 
\[
F_i=K(\vec y_i)\quad \text{con}\quad \vec y_i^{p^u}
\Witt -\vec y_i=\vec\delta_i,\quad 1\leq i\leq r,
\]
tenemos $|X_{P_i}|=|Y|=|Y_{P_i}|$.
M\'as a\'un la extensi\'on $\K F_i=K(\vec y,\vec y_i)=K(\vec y, 
\vec y\Witt - \vec y_i)=\K(\vec y\Witt - \vec y_i)$ 
es no ramificada en $P_i$ en $\K$. 
Se sigue que el campo asociado a
$\prod_P X_P=\prod_{i=1}^r X_{P_i}$
es $K(\vec y_1,\ldots, \vec y_r)$. Aqu\'i, el grupo de 
descomposici\'on $D$ es trivial.

Entonces, tenemos

\begin{teorema}\label{T4.3.A}
Con las notaciones anteriores, si
$\K=K(\vec y)$, entonces el campo asociado
a $\prod_P X_P=\prod_{i=1}^r X_{P_i}$
es 
\begin{gather*}
K(\vec y_1,\ldots, \vec y_r). \tag*{$\fin$}
\end{gather*}
\end{teorema}

\subsubsection{Campos de g\'eneros de
extensiones de Kummer}\label{Su5.1}

Aqu{\'\i} supondremos que $q\geq 3$. Primeramente queremos saber
cuando un campo $K(\sqrt[l]{P})$, donde
 $l\mid q-1$ y $P\in R_T^+$, est\'a contenido en
 $K(\Lambda_P)$. Sea $d=\deg P$.
 Por la Proposici\'on \ref{P5.1.1} se tiene que 
 $K(\sqrt[l]{(-1)^d P})\subseteq K(\Lambda_P)$.
 El grupo de Galois $\Gal(K(\Lambda_P)/K)
\cong \G P\cong {\ma F}_{q^d}^{\ast}$ es un grupo c{\'\i}clico de orden
$q^d-1$. Por lo tanto
existe una \'unica extensi\'on de la forma $K(\sqrt[l]{\alpha P})$,
$\alpha \in {\ma F}_q^{\ast}$, contenido en $K(\Lambda_P)$.
Notemos que si $\alpha\notin \big({\ma F}_q^{\ast}\big)^l$,
$K(\sqrt[l]{P})\neq K(\sqrt[l]{\alpha P})$ puesto que de otra forma
$\sqrt[l]{\alpha}\in K$ y as{\'\i} $\alpha\in \big({\ma F}_q^{\ast}
\big)^l$.

Usaremos que para cualquier $\alpha\in {\ma F}_q^{\ast}$, $1
\leq e\leq l-1$, $K(\sqrt[l]{\alpha P^e})=K(\sqrt[l]{\alpha^f P})$
donde $fe\equiv 1\bmod l$.
Puesto que tenemos $l$ clases $\bmod ({\ma F}_q^{\ast})^l$ en
${\ma F}_q^{\ast}$, los $l$ campos distintos $K(\sqrt[l]{\alpha P})$,
$\alpha\in {\ma F}_q^{\ast}$ est\'an dados por las clases
$\bmod ({\ma F}_q^{\ast})^l$. Por lo tanto $K(\sqrt[l]{\alpha^f P})
\subseteq K(\Lambda_P)$ si y s\'olo si $\alpha^f\equiv (-1)^d \bmod
({\ma F}_q^{\ast})^l$.

Sea $\K:=K(\sqrt[l]{\gamma D})$ donde $D\in
R_T$ es un polinomio m\'onico sin factores que sean $l$ potencias,
$\gamma\in {\ma F}_q^{\ast}$ y
$D=P_1^{e_1}\cdots P_r^{e_r}$ donde $P_i\in R_T^+$, $1\leq e_i
\leq l-1$, $1\leq i\leq r$. M\'as a\'un arreglamos el producto de tal forma
que $l\mid \deg P_i$ para $1\leq i\leq s$ y $l\nmid \deg
P_j$ para $s+1\leq j\leq r$, $0\leq s\leq r$.
En general, siempre se tiene que $E=K(\sqrt[l]{(-1)^{\deg D}D})$,
donde $E=\K{\ma F}_{q^l}\cap \cicl D{}$ y
${\ma F}_q^{\ast} \subseteq ({\ma F}_{q^l}^{\ast})^l$.

Usaremos la descomposici\'on de $\p$ dado por la 
Proposici\'on \ref{P5.1.2.Ram}.

\begin{observacion}\label{R4.2} Sean $\langle \sigma\rangle=
\Gal(\cicl N{}_m/\cicl N{})\cong\Gal(K_m/K)$. Entonces
cuando $\p$ es moderadamente ramificado en $\K/K$ y
$\K\subseteq \cicl N{}_m$ y el grupo de descomposici\'on
es $\eu D\cong \langle\sigma^t\rangle$, se tiene
\[
[\g\K:\K]=\frac{[E_{\eu{ge},m}:\K]}{|{\eu D}|}=\frac{[E_{\eu{ge},m}
:\K][\K_m:\K]}{m/t}=[E_{\eu{ge}}:E]t,
\]
donde $t$ es el grado de cualquier primo en $\S {\K}$.
\end{observacion}

En nuestro caso, la Observaci\'on \ref{R4.2}, tenemos que
$[\K_{\eu {ge}}:\K]=[E_{\eu {ge}}:E]t$ donde
\[
t=\deg\S \K=
\begin{cases}
1&\text{si $\p$ no es inerte en $\K/K$}\\
l&\text{si $\p$ es inerte en $\K/K.$}
\end{cases}
\]

Cuando $\K=E$, esto es, cuando $\K\subseteq \cicl D{}$, si $\chi$ es el
caracter de orden $l$ asociado a $\K$, $\chi=\chi_{P_1}\cdots
\chi_{P_r}$, consideramos $Y=\langle \chi_{P_i}\mid 1\leq i\leq r\rangle$.
El campo asociado a $Y$ es $F=K(\raizm 1,\ldots, \raizm r)$ y $\K_{\eu {ge}}=F$
si $l\nmid \deg D$ o si $l\mid \deg P_i$ para toda $i$ (esto es, $s=r$). 
Esto es as{\'\i} puesto que en el primer caso $\p$ es
ya ramificado en $\K$ y en segundo $\p$ es no ramificado en $F/K$
(Proposici\'on \ref{P5.1.2.Ram}).

Cuando $l\mid \deg D$ y $l\nmid \deg P_r$, $\p$ se ramifica en $F/K$ y es
no ramificado en $E/K$. En este caso $[F:E_{\eu {ge}}]=l$. 
Sean $a_{s+1},\ldots,
a_{r-1}\in{\ma Z}$ tales que $l\mid \deg(P_iP_r^{a_i})$, esto es,
$\deg P_i+a_i\deg P_r\equiv 0\bmod l$, $s+1\leq i\leq r-1$. Sea
\begin{align*}
F_1:&=K\big(\raiz 1,\ldots, \raiz s, \sqrt[l]{P_{s+1}P_r^{a_{s+1}}},\ldots,
\sqrt[l]{P_{r-1}P_r^{a_{r-1}}}\big)\\
&\subseteq F \subseteq K(\Lambda_{P_1P_2\cdots P_r}).
\end{align*}
Entonces $\S E$ se descompone en
$F_1/E$, $\K\subseteq F_1\subseteq E_{\eu {ge}}$ y $[F:F_1]=l$.
Se sigue que $E_{\eu {ge}}=F_1$.

En el caso general, del Teorema \ref{T2.1.A} obtenemos
$\K_{_{\eu {ge}}}=E_{_{\eu {ge}}} \K$. Por lo tanto

\begin{teorema}\label{T5.1.6}
Sea $D=P_1^{e_1}\cdots P_r^{e_r}\in R_T$ un polinomio
m\'onico que no tiene $l$ potencias, donde
$P_i\in R_T^+$, $1\leq e_i\leq l-1$, $1\leq i\leq r$. Sea $0\leq s\leq r$
tal que $l\mid \deg P_i$ para $1\leq i\leq s$ y $l\nmid \deg P_j$
para $s+1\leq j\leq r$. Sea $\K:=K(\sqrt[l]{\gamma D})$ donde $\gamma\in
{\ma F}_q^{\ast}$. 
Entonces $\K_{_{\eu {ge}}}$ est\'a dado por:
\l
\item
$K\big(\sqrt[l]{\gamma D},\raizm 1,\ldots, \raizm r\big)$ si $l\nmid \deg D$
o si $l\mid \deg P_i$ para toda $1\leq i\leq r$,

\item
$K\big(\sqrt[l]{\gamma D},\raiz 1,\ldots, \raiz s, \sqrt[l]{P_{s+1}P_r^{a_{s+1}}},\ldots,
\sqrt[l]{P_{r-1}P_r^{a_{r-1}}}\big)$, donde el exponente $a_j$ satisface
 $\deg P_j+a_j\deg P_r\equiv
0\bmod l$, $s+1\leq j\leq r-1$, si $l\mid \deg D$ y $l\nmid \deg P_r$. $\fin$
\end{list}
\end{teorema}

\subsubsection{Extensiones de Kummer en general}\label{S14.8.1.1}

El objetivo de esta subsecci\'on es dar una expresi\'on expl\'icita
de $\g {\K}$ donde ${\K}/K$ es una extensi\'on finita de Kummer del
campo de funciones racionales $K=\F(T)$. Una dificultad
t\'ecnica que se nos presenta es que, en general, para dos
extensiones finitas, se tiene que
${\K}_1,{\K}_2$ de $K$, $\g{({\K}_1)}\g{({\K}_2)}
\subsetneqq \g{({\K}_1{\K}_2)}$
(ver la Observaci\'on \ref{O12*.2.2.K}).
Sin embargo probaremos que cuando los grados de ${\K}_i/K$
son primos relativos, se cumple la igualdad. Con eso en mano
veremos que el caso de una extensi\'on de Kummer general
es consecuencia directa del caso de la igualdad y del
caso de una extensi\'on que es de grado una potencia
de un n\'umero primo.

En el estudio de campos de g\'eneros, el principal
obst\'aculo es la aparici\'on de inercia en la composici\'on
de campos y la contenci\'on de campos de g\'eneros
mencionado anteriormente. Esta es la raz\'on que el caso de una
extensi\'on abeliana de grada una potencia de $p$, donde
$p$ es la caracter\'istica, es mucho m\'as directo que
el caso Kummer.

Antes de continuar, para futura referencia recordamos el
llamado Lema de Abhyankar.\index{Abhyankar!Lema de
$\sim$}\index{Lema de Abhyankar} Se puede consultar una
demostraci\'on en \cite[Theorem 12.4.4]{Vil2006}.

\begin{teorema}[Lema de Abhyankar]\label{T14.8.1.1.1}
Sea $L/{\K}$ una extensi\'on separable de campos de funciones
congruentes. Supongamos que
$L={\K}_1{\K}_2$ con ${\K}\subseteq {\K}_i\subseteq L$, $1\leq i\leq 2$. Sea
$\pK$ un divisor primo de ${\K}$ y $\pL$ in divisor primo de $L$
que divide a $\pK$. Sean $\pL_i:=\pL\cap {\K}_i$, $i=1,2$. 
Si al menos una de las dos extensiones
${\K}_i/{\K}$ es moderadamente ramificada en $\pK$, entonces
\[
e_{L/{\K}}(\pL|\pK)=\lcm[e_{{\K}_1/{\K}}(\pL_1|\pK),e_{{\K}_2/{\K}}(\pL_2|\pK)],
\]
donde $e_{L/{\K}}(\pL|\pK)$ denota el \'indice de ramificaci\'on. $\fin$
\end{teorema}

\begin{observacion}\label{O14.8.1.1.2} El Lema de 
Abhyankar es v\'alida para extensiones de campos de funciones y
no \'unicamente para campos de funciones congruentes.
\end{observacion}

\subsubsection{Extensiones de Kummer de grado una potencia
de un primo}\label{S14.8.1.1.3}

Sean $l$ un n\'umero primo tal que $l^n$ divide a $q-1$
y ${\K}/K$ una extensi\'on de Kummer de exponente $l^n$.
De la Teor\'ia de Kummer, se tiene que ${\K}$ es el
compuesto ${\K}={\K}_1\cdots {\K}_s$ de extensiones de
Kummer c\'iclicas linealmente disjuntas. De manera m\'as
precisa, tenemos que $\K$ puede escribirse como
\begin{gather*}
{\K}=K\big(\sqrt[l^{n_1}]{\gamma_1 D_1},
\cdots, \sqrt[l^{n_s}]{\gamma_s D_s}\big)=
{\K}_1\cdots {\K}_s\\
\intertext{donde}
{\K}_{\varepsilon}=K\big(\sqrt[l^{n_{\varepsilon}}]{\gamma_{\varepsilon} D_{\varepsilon}}\big), \quad 1\leq \varepsilon \leq s,\\
\Gal({\K}/K)\cong \Gal({\K}_1/K)\times\cdots\times
\Gal({\K}_s/K)\cong C_{l^{n_1}}\times\cdots\times C_{l^{n_s}}
\end{gather*}
con $n=n_1\geq n_2\geq \cdots \geq n_s$,
$\gamma_{\varepsilon}\in \*\F$, $D_{\varepsilon}\in R_T$, $1\leq \varepsilon
\leq s$,  y $l^n|q-1$.
Entonces ${\K}_{\varepsilon}=K(\sqrt[l^{n_{\varepsilon}}]
{\gamma D_{\varepsilon}})$ es una extensi\'on c\'iclica de $K$ de
grado $l^{n_{\varepsilon}}$.

Sean $P_1,\ldots,P_r$ los primos finitos de $K$ ramificados en
$\K$, donde los primos $P_1,\ldots,P_r\in R_T^+$ son distintos. 
Sin p\'erdida de generalidad, se puede suponer que
\[
D_{\varepsilon}=\polyn {\varepsilon}\quad\text{con}\quad 0\leq
\alpha_{j,\varepsilon}\leq l^{n_{\varepsilon}}-1,\quad 1\leq j\leq r,
\quad 1\leq \varepsilon\leq s.
\]
De hecho, $\alpha_{j,\varepsilon}=0$ si y solamente si $P_{j}$ es 
no ramificado en ${\K}_{\varepsilon}/K$.

Sea $\alpha_{j,\varepsilon}=b_{j,\varepsilon} l^{a_{j,\varepsilon}}$
con $\mcd(b_{j,\varepsilon},l)=1$ cuando $\alpha_{j,
\varepsilon}\neq 0$ y sea $\deg P_{j}=c_{j} l^{d_{j}}$ con
$\mcd(c_{j},l)=1$, $1\leq j\leq r$.

Para $x\in{\ma Z}$,
$v_l(x)$ denota la valuaci\'on de $x$ en $l$. Esto es, $v_l(x)=\gamma$
si $l^{\gamma}|x$ y $l^{\gamma+1}\nmid x$. Escribimos $v_l(0)=\infty$.

\subsubsection{El caso ciclot\'omico}\label{S14.8.1.1.4}

Como primer paso, suponemos que ${\K}$ est\'a contenido
en un campo de funciones ciclot\'omicos, m\'as
precisamente, en $\cicl {D_1\cdots D_s}{}$.
Esto equivale a que $\gamma_{\varepsilon}\equiv
(-1)^{\deg D_{\varepsilon}}\mod ({\*\F})^{l^{n_{\varepsilon}}}$ para $1\leq
\varepsilon\leq s$ (ver Proposici\'on \ref{P5.1.1}).
Cuando ${\K}_{\varepsilon}$ est\'a contenido en un campo ciclot\'omico
de funciones podemos suponer, sin p\'erdida de generalidad, que
${\K}_{\varepsilon}=K(\sqrt[l^{n_{\varepsilon}}]{
D_{\varepsilon}^*})$, $1\leq \varepsilon\leq s$,
donde $D_{\varepsilon}^*=(-1)^{\deg  D_{\varepsilon}} D_{\varepsilon}$.
Notemos que si
$l^{n_{\varepsilon}}|\deg D_{\varepsilon}$, entonces 
$\Ku{l^{n_{\varepsilon}}}
{\*{D_{\varepsilon}}}=\Ku{l^{n_{\varepsilon}}}{D_{\varepsilon}}$.

Primero consideramos $F=\Ku{l^n}{\*D}$ con $D=\poly$,
una extensi\'on c\'iclica de Kummer de $K$. Sea $X=
\langle\chi\rangle$ el grupo de caracteres de Dirichlet
asociado a $F$. Notemos que para cualquier $\nu\in {\ma N}$ 
primo relativo a $l$, el campo asociado a
$\chi^{\nu}$ is $F$ pues $X=\langle\chi^{\nu}\rangle$.
Este hecho corresponde a que $F=\Ku{l^{n}}{\*{(D^{\nu})}}$. 

Cuando $D=P\in R_T^+$ se tiene que el caracter asociado a
$F$ es $\xbinom{}P_{l^{n}}$, el s\'imbolo de Legendre (Subsecci\'on
\ref{SubLegendre}). Por tanto, si $\chi_D$ es el caracter asociado a
$\Ku {l^n}{\*D}$, se tiene que $\chi_{D}=
\prod_{j=1}^r\chi_{P_{j}}^{\alpha_{j}}$.

Del Teorema \ref{TRam1} se tiene que si
$F=\Ku{m}{\gamma D}$ es una extensi\'on geom\'etrica y
separable de $K$ con $\gamma\in\*\F$ y $D=\poly\in R_T$, entonces
\begin{gather*}
e_{F/K}(P_{j})=\frac m{\mcd(\alpha_j,m)},\quad 1\leq j
\leq r\quad \text{y}\\
e_{\infty}(F|K):=e_{F/K}(\p)=\frac m{\mcd(\deg D,m)}.
\end{gather*}

De la Proposici\'on \ref{P3.3}, tenemos que si $X$ es
el grupo de caracteres asociado a la extensi\'on
ciclot\'omica c\'iclica de Kummer  $F=K(\sqrt[l^n]{\*D})$,
$Y=\prod_{P\in R_T^+} X_P$ es el grupo de caracteres
a ${\mc M}$, la m\'axima extensi\'on ciclot\'omica de $F$
no ramificada en los primos finitos. 

Sean $P=P_{j}$, $X=X_P=\langle\chi_P\rangle$ y $F_0$
es el campo asociado a $X_P$. Entonces $F_0$ est\'a
contenido en un campo de funciones ciclot\'omicos, 
$P$ es el \'unico primo finito ramificado en
$F_0/K$ y $P$ es moderadamente ramificado en $F_0/K$. 
Esto implica que $F_0
\subseteq \cicl P{}$ y que $\Gal(\cicl P{}/K)\cong C_{q^{d_P}-1}$ con
$d_P:=\deg P$. Por tanto $F_0$ es el \'unico subcampo de $\cicl P{}$
grado $o(\chi_P)=l^{\beta_P}$ sobre $K$. Puesto que
$F_0/K$ es una extensi\'on de Kummer, se sigue que
$F_0=\Ku{l^{\beta_P}}{\*P}$.

Expl\'icitamente tenemos que la m\'axima extensi\'on
ciclot\'omica de $F=\Ku {l^n}{\*D}$ no ramificada en los primos
finitos es
\begin{gather}\label{Ec14.8.1.1.5}{\mc M}:=K(\sqrt[l^n]{(P_1^{\alpha_1})^*},\ldots, \sqrt[l^n]
{(P_r^{\alpha_r})^*}).
\end{gather}

\begin{observacion}\label{O14.8.1.1.6}
Sea $\alpha=l^a b$ con $\mcd(b,l)=1$ y
$a < n$. Entonces $\Ku {l^n}{(P^{\alpha})^*}=
\Ku {l^{n-a}}{\*P}$. En particular, si $\alpha_{\varepsilon}
=l^{a_{\varepsilon}}b_{\varepsilon}$ con $\mcd(l,b_{\varepsilon})=1$,
$1\leq \varepsilon\leq r$, entonces
\[
{\mc M}=K(\sqrt[l^{n-a_1}]{P_1^*},\ldots, \sqrt[l^{n-a_r}]{P_r^*})=
F_1\cdots F_r,
\]
con $F_{\varepsilon}=\Ku {l^{n-a_{\varepsilon}}}{P_{
\varepsilon}^*}$, $1\leq \varepsilon\leq r$.
\end{observacion}

Otra demostraci\'on de (\ref{Ec14.8.1.1.5}) es usando el Lema
de Abhyankar, Teorema \ref{T14.8.1.1.1}. Por un lado, tenemos
\[
[{\mc M}:K]=\prod_{P\in R_T^+}|X_P|=\prod_{j=1}^r
|X_{P_{j}}|= \prod_{j=1}^r e_{{\mc M}/K}(P_{j})=\prod_{j=1}^r 
l^{n-a_{j}}.
\]
Por otro lado, si $F_{j}=\Ku{l^{n-a_{j}}}{\*{(P_{j})}}$, 
por el Lema de Abyankar se tiene que $F F_{j}/{\K}$ es
no ramificada en ning\'un primo finito por lo que
$F F_1\cdots F_r/F$ es no ramificada en los primos finitos y
$F\subseteq F_1\cdots F_r$. Obtenemos que $F_1\cdots
F_r\subseteq \gex F$ y que $[F_1\cdots F_r:K]=[{\mc M}:K]$.
Por tanto ${\mc M}=F_1\cdots F_r$.

Sean  $\alpha_{j}=b_{j} l^{a_{j}}$ con $\mcd(l,b_{j})
=1$ y $\deg P_{j}=l^{d_{j}}c_{j}$ con $\mcd(l,c_{j})=1$, 
$1\leq j\leq r$. Entonces
\begin{gather}
e_{P_j}(F|K)=l^{n-a_{j}},\nonumber\\
e_{\infty}(F|K)=\frac {l^n}{\gcd(l^n,\deg D)}:=l^t,\label{Ec3.1}\\
e_{\infty}(F_{j}|K)=\frac{l^{n-a_{j}}}
{\gcd(l^{n-a_{j}},\deg P_{j})}=
\frac{l^{n-a_{j}}}{l^{\min\{n-a_{j},d_{j}\}}}=
l^{n-a_{j}-\min\{n-a_{j},d_{j}\}}.\nonumber
\end{gather}

Del Lema de Abhyankar obtenemos que
\[
e_{\infty}({\mc M}|K)=\lcm_{1\leq j\leq r}[l^{n-a_{j}-\min\{
n-a_{j},d_{j}\}}]:=l^m.
\]

Por tanto $[{\mc M}:\g F]=l^{m-t}$. Para hallar $\g F$,
debemos encontrar un subcampo 
$F\subseteq L$ de ${\mc M}$ tal que $L/F$ sea no ramificada
y que $[{\mc M}:L]=l^{m-t}$. Si $L$ es tal campo, tenemos que
$L\subseteq \g F$ pues $\p$ es no ramificado en $L/F$ y $L$, 
siendo ciclot\'omico, satisface que $\p$ se descompone 
totalmente en $L/F$. Puesto que
$[{\mc M}:L]=[{\mc M}:\g F]$, se sigue que $L=\g F$.

Regresando a una extensi\'on de Kummer general ${\K}$ of $K$,
donde ${\K}=K\big(\sqrt[l^{n_1}]{\*{D_1}},
\cdots, \sqrt[l^{n_s}]{\*{D_s}}\big)={\K}_1\cdots {\K}_s$.
Sean $P_1,\ldots,P_r$ los primos finitos ramificados en $\K/K$.
Por el Lema de Abhyankar y por (\ref{Ec3.1}) se tiene, para
$P_j\in R_T^+$, que
\begin{gather}
e_{P_j}({\K}|K)=\lcm_{1\leq\varepsilon\leq s}
[e_{P_j}({\K}_{\varepsilon}|K)]=l^{\beta_j}\nonumber
\intertext{con}
\beta_j=\max_{1\leq \varepsilon\leq s}\{n_{\varepsilon}-
v_l(\alpha_{j,\varepsilon})\}=\max_{\substack{1\leq \varepsilon\leq s\\
b_{j,\varepsilon}\neq 0}} \{n_{\varepsilon}-a_{j,\varepsilon}\},\label{Ec3.5}
\intertext{y}
\begin{align*}
l^t:&=e_{\infty}({\K}|K)=\lcm_{1\leq \varepsilon\leq s}
[e_{\infty}({\K}_{\varepsilon}|K)]=
\lcm_{1\leq \varepsilon\leq s}\Big[\frac{l^{n_{\varepsilon}}}{\gcd(
l^{n_{\varepsilon}},\deg D_{\varepsilon})}\Big]\\
&=\lcm_{1\leq \varepsilon\leq s}\big[l^{n_{\varepsilon}-\min\{
n_{\varepsilon},v_l(\deg D_{\varepsilon})\}}\big],
\end{align*}
\end{gather}
esto es, 
\begin{gather}\label{Ec3.6}
e_{\infty}(\K|K)=l^t\quad\text{con}\quad
t=\max\limits_{1\leq\varepsilon\leq s}\big\{n_{\varepsilon}-
\min\{n_{\varepsilon},v_l(\deg D_{\varepsilon}\}\big\}.
\end{gather}

Sea $X_{\varepsilon}=\langle\chi_{\varepsilon}\rangle$ el grupo
de caracteres de Dirichlet correspondiente a ${\K}_{\varepsilon}$,
$1\leq \varepsilon\leq s$. Sea $\chi_{\varepsilon}=\prod_{P\in 
R_T^+} \chi_{\varepsilon,P}$ el producto $\chi_{\epsilon}$ 
en t\'ermino de sus $P$--componentes. Tenemos que
$e_P({\K}_{\varepsilon}|K)=o(\chi_{\varepsilon,P})$.

De esta forma, obtenemos que $X=\langle
\chi_1,\ldots,\chi_s\rangle$ es el grupo de caracteres de
Dirichlet asociado a ${\K}$, $X=X_1\cdots X_s$. Se tiene que
\[
X_P=(X_1)_P\cdots (X_s)_P=\langle\chi_{1,P}\rangle\cdots
\langle\chi_{s,P}\rangle=\langle\chi_{\gamma_P,P}\rangle
\]
con $o(\chi_{\gamma_P,P})=\max_{1\leq\varepsilon\leq s}
\{o(\chi_{\varepsilon,P})\}=e_P({\K}|K)$ para $P\in R_T^+$.

Sea ${\mc M}$ el campo asociado a $Y:=\prod_{P\in R_T^+}X_P$.
${\mc M}$ es la m\'axima extensi\'on ciclot\'omica de $\K$ no 
ramificada en ning\'un primo finito. De la Ecuaci\'on 
(\ref{Ec3.5}) obtenemos
\[
e_{P_{j}}({\K}|K)=l^{\beta_j}.
\]

Por tanto ${\mc M}=F_1\cdots F_r$ con
$F_{j}=\Ku{l^{\beta_j}}{\*{P_{j}}}$
y de la Ecuaci\'on (\ref{Ec3.1}) obtenemos que
\begin{align}
l^m:&=e_{\infty}({\mc M}|K)=\max_{1\leq j\leq r}\{e_{\infty}(F_j|K)\}=
\max_{1\leq j\leq r}\Big\{\frac{l^{\beta_j}}
{\gcd(l^{\beta_j},\deg P_j)}\Big\}\nonumber\\
&=\max_{1\leq j\leq r}\big\{
l^{\beta_j-\min\{\beta_j,d_j\}}\big\},\label{Ec3.7}
\end{align}
de tal forma que
$m=\max\limits_{1\leq j\leq r}\big\{\beta_j-\min\{\beta_j,d_j\}\big\}$.

El procedimiento para obtener $\g {\K}$
es el siguiente. Ordenamos $P_1,\ldots,P_r$ de tal forma que
$n=\beta_1\geq \beta_2\geq \cdots\geq \beta_r \geq 1$, 
esto es, ordenamos $P_1,\ldots, P_r$ en orden decreciente con respecto
a sus \'indices de ramificaci\'on en ${\K}/K$. Hay al menos un
$F_i$ tal que $e_{\infty}(F_i|K)=l^m$. Seleccionamos $i$ como el
m\'aximo \'indice con esta propiedad. Probaremos que existen ciertas
potencias $z_j$, $1\leq j\leq i-1$ tales que 
$\p$ es no ramificado en $\Ku {l^{\beta_j}}{(P_jP_i^{z_j})^*}/K$.
Escribamos $Q_j =P_jP_i^{z_j}$.

Para $j>i$ tenemos dos casos, $Q_j=P_jP_i^{y_jl^{\varepsilon_j}}$ 
o $Q_j=P_j^{y_j}P_i^{l^{\varepsilon_j}}$ para algunas $y_j,\varepsilon_j
\in {\ma Z}$ tales que el campo
$F_j=\Ku {l^{\gamma_j}}{Q_j}$, para alguna $\gamma_j$,
satisface que el \'indice de ramificac\'on de $P_j$
en $F_j/K$ es $l^{\beta_j}$, en $P_i$  es menor o igual a $l^{\beta_i}$
y en $\p$ es no ramificado. El resto se seguir\'a tomando la
composici\'on de todos estos campos y uno m\'as de la forma
$\Ku{l^{\xi_i}}{P_i^*}$ para alg\'un $\xi_i$.

El resultado para extensiones ciclot\'omicos c\'iclicas de grado una
potencia de un primo es el siguiente.

\begin{teorema}\label{T14.8.1.1.7} 
Sea ${\K}/K$ una $l$--extensi\'on ciclot\'omica finita de Kummer,
${\K}={\K}_1\cdots {\K}_s$,
${\K}_{\varepsilon}=\Ku{l^{n_{\varepsilon}}}{D_{\varepsilon}^*}$,
$D_{\varepsilon}\in R_T$, $1\leq \varepsilon\leq s$ y
$\Gal({\K}/K)\cong C_{l^{n_1}}\times \cdots \times C_{l^{n_s}}$
con $n=n_1\geq n_2\geq \cdots \geq n_s$ y $l^n|q-1$.
Entonces ${\K}$ es un subcampo del campo de funciones ciclot\'omicos
${\K}\subseteq \cicl {D_1\cdots D_s}{}$.
Sean $P_1,\ldots,P_r$ los primos finitos en $K$ ramificados en
${\K}/K$ con $P_1,\ldots,P_r\in R_T^+$ polinomios m\'onicos
distintos. Sean
\begin{gather*}
e_{P_j}({\K}|K)=l^{\beta_j},\quad
1\leq \beta_j\leq n, \quad 1\leq j\leq r,\quad\text{y}\\
e_{\infty}({\K}|K)=l^t, \quad 0\leq t\leq n
\end{gather*}
dadas por las Ecuaciones
{\rm{(\ref{Ec3.5})}} y {\rm{(\ref{Ec3.6})}} y sea
$\deg P_j=c_jl^{d_j}$ con $\gcd(c_j,l)=1$, $1\leq j\leq r$.

Ordenamos $P_1,\ldots, P_r$ de tal forma que $n=\beta_1\geq
\beta_2\geq \ldots \geq \beta_r$.

La m\'axima extensi\'on ciclot\'omica no ramificada
en los primos finitos ${\mc M}$ de $\K$ est\'a dada por
${\mc M}=\gex {\K}=K\Big(\sqrt[l^{\beta_1}]{P_1^*},\ldots,
\sqrt[l^{\beta_r}]{P_r^*}\Big)$. Sea $l^m=e_{\infty}({\mc M}|K)$
dado por la Ecuaci\'on {\rm{(\ref{Ec3.7})}}.

Seleccionamos el \'indice $i$ tal que $m=\beta_i-\min\{\beta_i,d_i\}$
y tal que para $j>i$ se tiene
$m>\beta_j-\min\{\beta_j,d_j\}$. Esto es, $i$ es el m\'aximo
\'indice donde se obtiene
$l^m$ como el \'indice de ramificaci\'on de $\p$.

En el caso $m=t$ se tiene ${\mc M}=\g {\K}=
\prod_{j=1}^r\Ku{l^{\beta_j}}{P_j^*}$. 

En el caso $m>t\geq 0$, tenemos $\min\{
\beta_i,d_i\}=d_i$ y $m=\beta_i-d_i$. Sean $a,b\in{\ma Z}$ tales que
$a\deg P_i+b l^{n+d_i}=l^{d_i}=\gcd(l^{n+d_i},\deg P_i)$. Definimos $z_j=-a
\frac{\deg P_j}{l^{d_i}}=-ac_jl^{d_j-d_i}\in {\ma Z}$ para $1\leq j\leq i-1$.
Para $j>i$, consideremos $y_j\in{\ma Z}$ con $y_j\equiv -c_i^{-1}c_j\bmod l^n
=-ac_j \bmod l^n$. Sean
\[
E_j=
\begin{cases}
\Ku{l^{\beta_j}}{P_jP_i^{z_j}}&\text{si $j<i$},\\
\Ku{l^{d_i+t}}{\*{P_i}}&\text{si $j=i$},\\
\Ku{l^{\beta_j}}{P_jP_i^{y_jl^{d_j-d_i}}}&\text{si $j>i$ y $d_j\geq d_i$},\\
\Ku{l^{\beta_j+d_i-d_j}}{P_j^{l^{d_i-d_j}}{P_i^{y_j}}}&\text{si $j>i$ y $d_i> d_j$}.
\end{cases}
\]

Entonces $\g {\K}=E_1\cdots E_{i-1}E_iE_{i+1}\cdots E_r$.
\end{teorema}

\begin{proof}
Cuando $m=t$ se tiene que $\g {\K}={\mc M}=\prod_{j=1}^r F_j$ donde
$F_j=\Ku{l^{\beta_j}}{\*{P_j}}$.

Supongamos que $m>t\geq 0$. Entonces $d_i<\beta_i$ y $\beta_i-d_i=m$. 

Para $j<i$ tenemos $\beta_j\geq \beta_i$ y
$\beta_j-d_j\leq \beta_j-\min\{\beta_j,d_j\}\leq m=\beta_i-d_i$. Por tanto $d_j\geq
\beta_j-m= \beta_j-(\beta_i-d_i)=\beta_j-\beta_i+d_i\geq d_i$. 
En particular $d_i|\deg P_j$,
$1\leq j\leq i-1$. Puesto que $\mcd(l^{n+d_i},\deg P_i)=l^{d_i}$,
existen $a,b\in{\ma Z}$, tales que $a\deg P_i+bl^{n+d_i}=l^{d_i}$. 
Multipicando por $\deg P_j$ y dividiendo entre $l^{d_i}$, $j<i$, obtenemos 
\[
a\frac{\deg P_j}{l^{d_i}}\deg P_i+b\deg P_j l^{n}=\deg P_j.
\]
Esto es, $\deg P_j+z_j\deg P_i=b\deg P_j l^{n}$, donde
$z_j= -a\frac{\deg P_j}{l^{d_i}}=-ac_j l^{d_j-d_i}$. 
Notemos que $\mcd(a,l)=1$.
Se tiene $l^n|\deg(P_jP_i^{z_j})$.
Sea $E_j:=\Ku{l^{\beta_j}}{P_jP_i^{z_j}}$. Se sigue que 
\begin{gather*}
e_{\infty}(E_j|K)=1, \quad e_{P_j}(E_j|K)=l^{\beta_j}=
e_{P_j}({\K}|K)\quad \text{y}\\
e_{P_i}(E_j|K)=l^{\beta_j-v_l(z_j)}|l^{\beta_i}=e_{P_i}({\K}|K),
\end{gather*}
debido a que $v_l(z_j)=d_j-d_i$ y a que $\beta_j-v_l(z_j)=
\beta_j-d_j+d_i\leq m+d_i=\beta_i$. En particular $E_j\subseteq \g {\K}$.

Ahora consideramos el caso $j>i$. Sea $y_j\in{\ma Z}$ tal que
$y_j\equiv -c_i^{-1}c_j\bmod l^n$. Puesto que $ac_i\equiv 1\bmod
l^n$ tenemos $c_i^{-1}\equiv a \bmod l^n$.
Esto es posible puesto que $\mcd(c_ic_j,l)=1$.
Notemos que $\mcd(y_j,l)=1$. 

Primero consideremos el caso $d_j\geq d_i$. Sean
\begin{gather*}
Q_j=P_jP_i^{y_j l^{d_j-d_i}}\quad \text{y}\quad E_j:=\Ku{l^{\beta_j}}{Q_j}.
\intertext{Tenemos}
\deg Q_j=\deg P_j+y_j l^{d_j-d_i} \deg P_i=c_jl^{d_j}+y_jl^{d_j-d_i} c_i
l^{d_i}=l^{d_j}(c_j+y_jc_i).
\end{gather*}
Se sigue que $l^n|\deg Q_j$ y que $e_{\infty}(E_j|K)=1$. 
Por otro lado 
\[
e_{P_j}(E_j|K)=l^{\beta_j}=e_{P_j}({\K}|K)\quad \text{y}\quad
e_{P_i}(E_j|K)=l^{\beta_j-d_j+d_i-v_l(y_j)}.
\]
Puesto que $\gcd(y_j,l)=1$, tenemos
$\beta_j-d_j+d_i\leq m+d_i=\beta_i$ y $e_{P_i}(E_j|K)|e_{P_i}({\K}|K)$.
Se sigue que $E_j\subseteq \g {\K}$.

Ahora consideremos el caso $d_i\geq d_j$. Definimos
$Q_j=P_j^{l^{d_i-d_j}}P_i^{y_j}$
y sea $E_j=\Ku{l^{\beta_j+d_i-d_j}}{Q_j}$. Tenemos
\[
\deg Q_j=l^{d_i-d_j}\deg P_j+y_j\deg P_i=l^{d_i-d_j}c_jl^{d_j}+y_jc_il^{d_i}=
l^{d_i}(c_j+y_jc_i).
\]
Por tanto $l^n|\deg Q_j$ y $e_{\infty}(E_j|K)=1$. Por otro lado tenemos
\begin{gather*}
e_{P_j}(E_j|K)=l^{\beta_j+d_i-d_j-d_i+d_j}=l^{\beta_j}=e_{P_j}(E_j|K) 
\intertext{y}
e_{P_i}(E_j|K)=l^{\beta_j+d_i-d_j}|l^{\beta_i}=e_{P_i}({\K}|K).
\end{gather*}
Se sigue que $E_j\subseteq \g {\K}$.

Por tanto $L:=E_1\cdots E_{i-1}E_{i+1}\cdots E_r\subseteq \g {\K}$.
Del Lema de Abhyankar obtenemos que
\begin{gather*}
e_{\infty}(L|K)=1,\\\
e_{P_j}(L|K)=e_{P_j}({\K}|K)=l^{\beta_j}, j\neq i
\intertext{y}
e_{P_i}(L|K)|e_{P_i}({\K}|K)=l^{\beta_i}.
\end{gather*}
De hecho, podemos dar un argumento directo para probar que
realmente $e_{P_i}(L|K)= e_{P_i}({\K}|K)=l^{\beta_i}$. 
Ver la observaci\'on \ref{O14.8.1.1.8}.

Para cualquier $1\leq j\leq r$, denotamos por $I_j$ al grupo de inercia de
$P_j$ en ${\mc M}'/K$ donde ${\mc M}'$ es cualquier
subcampo de ${\mc M}$ que contiene a $E_j$. Para cualquier tal campo 
${\mc M}'$, se tiene $|I_j|=l^{\beta_j}$.

Sean ${\mc J}:=\{j>i\mid d_i>d_j\}$ y ${\mc I}:=\{1,2,\ldots, i-1,i+1,\ldots,
r\}\setminus {\mc J}$. En otras palabras, si $j\in{\mc I}$
entonces $E_j=\Ku{l^{\beta_j}}{
P_jP_i^{x_j}}$ para alg\'un $x_j\in {\ma Z}$. 

Definimos ${\mc I}=\{m_1,\ldots,m_u\}$.
Sea $F:=E_{m_1}\cap E_{m_2}$. Puesto que se tiene que
$P_{m_1}$ es totalmente ramificado en $E_{m_1}$ 
y no ramificado en $E_{m_2}$,
se sigue que $F=K$ y por tanto 
$[E_{m_1}E_{m_2}:K]=[E_{m_1}:K][E_{m-2}:K]$.
\[
\xymatrix{
E_{m_1}\ar@{-}[r]\ar@{-}[d]&E_{m_1}E_{m_2}\ar@{-}[d]\\
K=E_{m_1}\cap E_{m_2}\ar@{-}[r]&E_{m_2}
}
\]
Adem\'as, $\Gal(E_{m_1}E_{m_2}/K)\cong \Gal(E_{m_1}/K)\times
\Gal(E_{m_2}/K)\cong I_{m_1}\times I_{m_2}$ y $(E_{m_1}E_{m_2})^{
I_{m_1}I_{m_2}}=K$. Por inducci\'on se obtiene para $1\leq v\leq u$:
\las
\item $\big(E_{m_1}\cdots E_{m_{v-1}}\big)\cap E_{m_v}=K$,
\item $\big[E_{m_1}\cdots E_{m_{v}}:K\big]=[E_{m_1}:K]\cdots [E_{m_v}:K]$,
\item $\big(E_{m_1}\cdots E_{m_{v}}\big)^{I_{m_1}\cdots I_{m_v}}=K$,
\item $I_{m_1}\cdots I_{m_v}\cong I_{m_1}\times\cdots\times I_{m_v}$.
\end{list}

Esto es, para cualquier $\mu\in {\mc I}$ tenemos que
$\Big(\prod_{j\in{\mc I}\setminus \{\mu\}}E_j\Big)\cap 
E_{\mu}=K$ puesto que para cualquier subcampo no trivial de
$A:=\prod_{j\in{\mc I}\setminus \{\mu\}}E_j$ al menos uno de $P_j$ con $j\in
{\mc I}\setminus \{\mu\}$ es ramificado en este subcampo y $P_j$
es no ramificado en $E_{\mu}$. En particular se tiene
\begin{gather}\label{Ec3.2}
\big[\prod_{j\in{\mc I}}E_j:K\big]=\prod_{j\in{\mc I}}[E_j:K].
\end{gather}

Tambi\'en tenemos que
\begin{gather}\label{Ec3.3}
\Big(\prod_{j\in{\mc I}}E_j\Big)\bigcap \Big(\prod_{j\in{\mc J}}E_j\Big)=K
\end{gather}
puesto que en cualquier subcampo no trivial de
$\prod_{j\in{\mc I}}E_j$ al menos uno de los primos
$P_{\mu}$ con $\mu\in {\mc I}$ es ramificado y $P_{\mu}$ es no ramificado
en $\prod_{j\in{\mc J}}E_j$. En otras palabras,
\[
[L:K]=\Big[\prod_{j\neq i}E_j:K\Big]=\Big[\prod_{j\in{\mc I}}E_j:K\Big]
\Big[\prod_{j\in{\mc J}}E_j:K\Big].
\]

Para calcular $\Big[\prod_{j\in{\mc J}}E_j:K\Big]$ ordenamos ${\mc J}$ como
sigue. Escribamos ${\mc J}=\{j_1,\ldots,j_s\}$ con $d_i-d_{j_1}\leq d_i-
d_{j_2}\leq\cdots \leq d_i-d_{j_s}$. Tenemos que $E_j^{I_j}=\Ku{l^{d_i-
d_j}}{P_j^{y_j}}=\Ku{l^{d_i-d_j}}{P_j}$.

En primer lugar consideramos
$E_{j_1}E_{j_2}$. Tenemos que $E_{j_1}\cap E_{j_2}=
C_1$ donde denotamos $C_u:=\Ku{l^{d_i-d_{j_u}}}{P_i}$, 
$j_u\in{\mc J}$, $1\leq u\leq s$. De hecho,
si $\Lambda:=E_{j_1}\cap 
E_{j_2}$, entonces $P_{j_1}$ es no ramificado en $E_{j_2}$ y
$P_{j_2}$ es no ramificado en $E_{j_1}$, as\'i que,
el \'unico ramificado primo en
$\Lambda/K$ es $P_i$. M\'as a\'un, $C_1\subseteq E_{j_1}$
y $C_2\subseteq E_{j_2}$ y $C_1=C_1\cap C_2\subseteq
E_{j_1}\cap E_{j_2}$. Ahora, $P_{j_1}$ es totalmente ramificado en
$E_{j_1}/C_1$. En particular, si $C_1\subsetneqq C'\subseteq F_{j_1}$, $P_{
j_1}$ es ramificado en $C'/C_1$ y puesto que $P_{j_1}$ es no ramificado en
($E_{j_1}\cap E_{j_2})/C_1$, se sigue que $E_{j_1}\cap E_{j_2}=C_1$.

Consideremos el siguiente diagrama
\[
\xymatrix{
& E_{j_1}\ar@{-}[r]\ar@{-}[d]&E_{j_1}E_{j_2}\ar@{-}[d]\\
&C_1\ar@{-}[r]\ar@{-}[dl]&E_{j_2}\\ K
}
\]
Tenemos $[E_{j_1}E_{j_2}:K]=[E_{j_1}:C_1][E_{j_2}:C_1][C_1:K]=
\frac{[E_{j_1}:K][E_{j_2}:K]}{[C_1:K]}$.

Ahora consideremos $E_{j_1}E_{j_2}E_{j_3}$. Con un argumento
similar al del caso previo, se tiene
$E_{j_1}E_{j_2}\cap E_{j_3}=C_2$. Consideremos el siguiente diagrama
\[
\xymatrix{
& E_{j_1}E_{j_2}\ar@{-}[r]\ar@{-}[d]&E_{j_1}E_{j_2}E_{j_3}\ar@{-}[d]\\
&C_2\ar@{-}[r]\ar@{-}[dl]&E_{j_3}\\ K
}
\]
Por tanto $[E_{j_1}E_{j_2}E_{j_3}:K]=[E_{j_1}E_{j_2}:C_2][E_{j_3}:C_2][C_2:K]
=[E_{j_1}E_{j_2}:K][E_{j_3}:C_2]=
\frac{[E_{j_1}:K][E_{j_2}:K][E_{j_3}:K]}{[C_1:K][C_2:K]}$.

Por inducci\'on obtenemos, para $1\leq v\leq s$,
\las
\item $\big(E_{j_1}\cdots E_{j_{v-1}}\big)\cap E_{j_v}=C_{v-1}$,
\item $\big[E_{j_1}\cdots E_{j_{v}}:K\big]=[E_{j_1}:C_1]\cdots 
[E_{j_{v-1}}:E_{j_v}]
[E_{j_v}:C_v][C_v:K]=\big(\prod_{\mu=1}^v l^{\beta_{j_{\mu}}}\big)
l^{d_i-d_{j_v}}$,
\item $\big(E_{j_1}\cdots E_{j_{v}}\big)^{I_{j_1}\cdots I_{j_v}}=K$,
\item $I_{j_1}\cdots I_{j_v}\cong I_{j_1}\times\cdots\times I_{j_v}$.
\end{list}
Esto es,
\begin{align}
\Big[\prod_{j\in{\mc J}}E_j:K\Big]&=[E_{j_1}\cdots E_{j_s}:K]=
\Big(\prod_{u=1}^s [E_{j_u}:C_u]\Big)[C_s:K]\nonumber\\
&=\big(\prod_{j\in{\mc J}}l^{\beta_j}\big)l^{d_i-d_s}.\label{Ec3.4}
\end{align}

De las Ecuaciones (\ref{Ec3.2}), (\ref{Ec3.3}) y (\ref{Ec3.4}), obtenemos
\begin{gather*}
[L:K]=\Big(\prod_{\substack{j=1\\ j\neq i}}^r l^{\beta_j}\Big) l^{d_i-d_s}
\quad\text{y}\quad L\cap E_i=\Ku{l^{d_i-d_{j_s}}}{P_i}=C_s,
\intertext{donde $E_i:=\Ku{l^{d_i+t}}{\*{P_i}}$. Ahora, tenemos que}
e_{P_i}(E_i|K)=l^{d_i+t}|l^{\beta_i}=l^{d_i+m}=e_{P_i}({\K}|K)
\intertext{y}
e_{\infty}(F_i|K)=\frac{l^{d_i+t}}{\gcd(l^{d_i},l^{d_i+t})}=l^{d_i+t-d_i}
=l^t=e_{\infty}({\K}|K).
\end{gather*}

Se sigue que $LE_i\subseteq \g {\K}$ y que
\begin{align*}
[LE_i:K]&=[L:C_s][E_i:C_s]
[C_s:K]=[L:C_s][E_i:K]\\
&=\frac{\prod_{j=1}^rl^{\beta_j}}{l^{\beta_i-(d_i+t)}}
=\frac{[{\mc M}:K]}{l^{m-t}}=\frac{[{\mc M}:K]}{[{\mc M}:\g {\K}]}=[\g {\K}:K].
\end{align*}
Por tanto $\g {\K}=LE_i=E_1\cdots E_{i-1}E_iE_{i+1}\cdots E_r$. $\fin$
\end{proof}

\begin{observacion}\label{O14.8.1.1.8}
En la notaci\'on del Teorema \ref{T14.8.1.1.7}, tenemos que en el
caso $m>t$, ${\K}\subseteq \g {\K}=E_1\cdots E_{i-1} E_i E_{i+1}\cdots
E_r$. Por tanto existe $j\neq i$ tal que $e_{P_i}(E_j|K)=l^{\beta_i}$.
Sin embargo uno se pregunta cual es la raz\'on de esto.
Aqu\'i damos una prueba directa. Para $j<i$
tenemos que $E_j=\Ku{l^{\beta_j}}{P_jP_i^{z_j}}$ con $z_j=-acl^{d_j-d_i}$.
Ahora, $d_j\geq d_i$ y $\beta_j\geq \beta_i$ ($j<i$) y $v_l(z_j)=
d_j-d_i$ puesto que $\mcd(ac_j,1)=1$. Entonces $e_{P_i}(E_j|K)=
l^{\beta_j-v_l(z_j)}=l^{\beta_j-d_j+d_i}$. As\'i, requerimos que
para alg\'un $j<i$ se tiene $\beta_j-d_j+d_i=\beta_i$ o, equivalentemente,
$\beta_j-d_j=\beta_i-d_i$.

De la definici\'on del \'indice $i$ tenemos que
$\beta_j-d_j\leq m=\beta_i-d_i$
y que para $j>i$ se tiene $\beta_j-d_j<m$.

Supongamos que $\beta_j-d_j<m$ para toda $j\neq i$. 
De la Ecuaci\'on (\ref{Ec3.5}) obtenemos que
$\beta_j=\max_{1\leq \varepsilon\leq s}\{n_{\varepsilon}-
a_{j,\varepsilon}\}$ donde, por conveniencia seleccionamos $a_{j,\varepsilon}=
n$ en el caso $\alpha_{j,\varepsilon}=b_{j,\varepsilon}l^{a_{j,\varepsilon}}=0$,
esto es, cuando $b_{j,\varepsilon}=0$ 
debido a que, de esta forma, $n_{\varepsilon}
-a_{j,\varepsilon}\leq 0$ y el m\'aximo no puede ser obtenido en $\varepsilon$
puesto que $1\leq \beta_j\leq n$.

Sea $1\leq \mu\leq s$ tal que $\beta_i=\max_{1\leq \varepsilon\leq s}
\{n_{\varepsilon}-a_{i,\varepsilon}\}=n_{\mu}-a_{i,\mu}$ de forma
que $m=\beta_i-d_i=n_{\mu}-a_{i,\mu}-d_i$ y
\begin{gather}\label{Ec3.8}
a_{i,\mu}+d_i=n_{\mu}-m.
\end{gather}
Definimos $\deg D_{\mu}=c_0l^{d_0}$ con $\mcd 
(c_0,l)=1$, esto es, $v_l(\deg D_{\mu})
=d_0$. Puesto que $D_{\mu}=\prod_{j=1}^rP_j^{\alpha_{j,\mu}}$ tenemos
\begin{gather}\label{Ec3.9}
\deg D_{\mu}=\sum_{j=1}^r \alpha_{j,\mu}\deg P_j=\sum_{j=1}^r b_{j,\mu}
l^{a_{j,\mu}}c_jl^{d_j}=\sum_{j=1}^r b_{j,\mu}c_jl^{a_{j,\mu}+d_j}.
\end{gather}

Fijemos $j\neq i$ y sea $\beta_j=\max_{1\leq\varepsilon\leq s}\{n_{\varepsilon}-
a_{j,\varepsilon}\}\geq n_{\mu}-a_{j,\mu}$. Por tanto $a_{j,\mu}\geq n_{\mu}
-\beta_j$ y $a_{j,\mu}+d_j\geq n_{\mu}-\beta_j+d_j=n_{\mu}-(\beta_j-d_j)
>n_{\mu}-m$. De la Ecuaci\'on (\ref{Ec3.8}) obtenemos
\begin{gather}\label{Ec3.10}
a_{j,\mu}+d_j>n_{\mu}-m\quad \text{para}\quad j\neq i\quad\text{y}\quad
a_{i,\mu}+d_i=n_{\mu}-m.
\end{gather}

De las Ecuaciones (\ref{Ec3.9}) y (\ref{Ec3.10}) se sigue que
\[
d_0=a_{i,\mu}+d_i=n_{\mu}-m.
\]
Por otro lado tenemos que
\begin{align*}
t&=\max_{1\leq \varepsilon\leq s}\{n_{\varepsilon}-\min\{n_{\varepsilon},
v_l(\deg D_{\varepsilon})\}\geq n_{\mu}-\min\{n_{\mu},d_0\}\\
&=n_{\mu}-\min\{n_{\mu},n_{\mu}-m\}=n_{\mu}-(n_{\mu}-m)=m,
\end{align*}
esto es, $t\geq m$, lo cual contradice nuestra suposici\'on
de que: $t<m$. Por tanto,
existe $j<i$ tal que $\beta_j-d_j=m$ y $e_{P_i}(E_j|K)=l^{\beta_i}
=e_{P_i}({\K}|K)$.
\end{observacion}

\subsubsection{El caso general de grado una
potencia de un primo}\label{S14.8.1.1.9}

Ahora consideramos ${\K}/K$ una extensi\'on de Kummer de
orden una potencia de  $l$. Si el grupo de Galois
$\Gal({\K}/K)$ es de exponete $l^n$, se tiene que $\Gal({\K}/K)\cong
C_{l^{n_1}}\times\cdots\times C_{l^{n_s}}$ con $n=n_1\geq\cdots
\geq n_s$ y $l^n|q-1$. As\'i, tenemos que ${\K}$ es de la forma
${\K}=K\big(\sqrt[l^{n_1}]{\gamma_1 D_1},\ldots, 
\sqrt[l^{n_s}]{\gamma_s D_s}\big)$ con $D_{\varepsilon}\in
R_T$, $\gamma_{\varepsilon}\in
\*\F$, $1\leq\varepsilon\leq s$ and ${\K}_{\varepsilon}=\Ku{l^{n_{\varepsilon}}}
{\gamma_{\varepsilon} D_{\varepsilon}}$.

Sea $E={\K}K_{l^n}\cap \cicl {D_1\cdots D_s}{}=
K\big(\sqrt[l^{n_1}]{D_1^*},\ldots, \sqrt[l^{n_s}]{D_s^*}\big)$.
Por el Teorema \ref{T2.1.A} tenemos que $\g {\K}=\g E^H {\K}$ donde $H$
es el grupo de descomposici\'on de $\p$ tanto en
${\K}\g E/\g {\K}$ como en ${\K}\g E/{\K}$.
Tambi\'en se tiene que ${\K}E/{\K}E^{H_1}$ 
es una extensi\'on de constantes donde
$H_1=H|_E$ es el grupo de descomposici\'on de $\p$ en ${\K}E/{\K}$.

Tenemos que $E{\K}=E\big(\sqrt[l^{n_1}]{\gamma_1 D_1},\ldots, 
\sqrt[l^{n_s}]{\gamma_s D_s}\big)=E\big(\sqrt[l^{n_1}]{\varepsilon_1},\ldots, 
\sqrt[l^{n_s}]{\varepsilon_s}\big)$ donde $\varepsilon_{j}=(-1)^{\deg D_{j}}
\gamma_{j}$, $1\leq j\leq s$ puesto que 
$\sqrt[l^{n_j}]{\varepsilon_j}=\frac{\sqrt[l^{n_j}]{\gamma_j D_j}}{\sqrt[l^{n_j}]
{D_j^*}}$ y $\sqrt[l^{n_j}]{D_j^*}\in E$. En particular, $E{\K}/E$ es una
extensi\'on de constantes y el grupo de inercia de $\p$ en $E{\K}/K$ es
$f=l^v$ donde ${\ma F}_{q^{l^v}}=\F\big(\sqrt[l^{n_1}]{\varepsilon_1},\ldots, 
\sqrt[l^{n_s}]{\varepsilon_s}\big)$ puesto que, $E$ 
siendo ciclot\'omico, satisface que el grado de inercia de
$\p$ en $E/K$ es $1$.

Se sigue que
\[
|H|=\frac{[\F\big(\sqrt[l^{n_1}]{\varepsilon_1},\ldots, 
\sqrt[l^{n_s}]{\varepsilon_s}\big):\F]}{\deg_{\K} \p}:=l^u.
\]

Puesto que $H_1=H|_E\subseteq I_{\infty}(E|K)$ donde $I_{\infty}(E|K)$
denota al grupo de inercia de $\p$ en $E/K$ tambi\'en denota
al grupo de inercia de $\p$ en $\g E/K$. En particular $\g E^+=
\g E^{I_{\infty}(E|K)}$. Denotamos por ${\mc H}_1$ el
grupo de inercia de $\p$ en $\g E/K$. Del Teorema \ref{T14.8.1.1.7}
tenemos que $\p$ es no ramificado en $L=E_1\cdots E_{i-1}E_{i+1}\cdots
E_r$ y totalmente ramificado en $\g E/L\Ku{l^{d_i}}{P_i^*}$.
Se sigue que $\g E^+=L \Ku{l^{d_i}}{P_i^*}$ y que
\[
I_{\infty}(\g E|K)=\Gal\big(\g E/L\Ku{l^{d_i}}{P_i^*}\big).
\]
El grupo ${\mc H}_1$ es el subgrupo de orden
$l^u$ of $I_{\infty}(\g E|K)$.

Se sigue que $\g E^{{\mc H}_1}=L\Ku{l^{d_i+t-u}}{P_i^*}$ y por tanto 
\[
\g {\K}=\g E^{{\mc H}_1} {\K}=E_1\cdots E_{i-1} \Ku{l^{d_i+t-u}}{P_i^*}
E_{i+1}\cdots E_r {\K}.
\]

De esta forma hemos probado el resultado principal sobre
extensiones de Kummer de grado una potencia de un n\'umero 
primo sobre $K$.

\begin{teorema}\label{T14.8.1.1.10}
Sea ${\K}=K\big(\sqrt[l^{n_1}]{\gamma_1 D_1},\ldots, 
\sqrt[l^{n_s}]{\gamma_s D_s}\big)$ 
una extensi\'on de Kummer de $K$ de grado una potencia
de un n\'umero primo con $D_{\varepsilon}\in
R_T$, $\gamma_{\varepsilon}\in
\*\F$, $1\leq\varepsilon\leq s$ y ${\K}_{\varepsilon}=\Ku{l^{n_{\varepsilon}}}
{\gamma_{\varepsilon} D_{\varepsilon}}$. Sea $E=
K\big(\sqrt[l^{n_1}]{D_1^*},\ldots, \sqrt[l^{n_s}]{D_s^*}\big)$.

Con las notaciones del Teorema {\rm{\ref{T14.8.1.1.7}}}
tenemos que
\[
\g {\K}=L\Ku{l^{d_i+t-u}}{P_i^*}{\K}=E_1\cdots E_{i-1}E_{i+1}\cdots E_r
\Ku{l^{d_i+t-s}}{P_i^*}{\K}
\]
donde $l^u=\frac{[\F(\sqrt[l^{n_1}]{\varepsilon_1},\ldots, 
\sqrt[l^{n_s}]{\varepsilon_s}):\F]}{\deg_{\K} \p}$
y $\varepsilon_j =(-1)^{\deg D_j}\gamma$, $1\leq j\leq s$. $\fin$
\end{teorema}

\subsubsection{Extensiones de Kummer generales}\label{S14.8.1.1.11}

Sea${\K}/K$ una extensi\'on de Kummer. Sea $G:=\Gal({\K}/K)$
el grupo de Galois de ${\K}/K$. Si $S_1,\ldots, S_h$ son
los diversos subrupos de Sylow de $G$ entonces $G\cong S_1\times
\cdots\times S_h$ y cada $S_j$ es un grupo de orden
una potencia de un primo, digamos $|S_j|=l_j^{m_j}$, $1\leq j\leq h$.

Podemos escribir ${\K}={\K}_1\cdots {\K}_h$ con $\Gal({\K}_j/K)
\cong S_j$, $1\leq j\leq h$.

Del Teorema \ref{T14.8.1.1.10},
conocemos cada $\g {({\K}_j)}$, $1\leq j\leq h$.
El conocimiento de $\g {\K}$ se seguir\'ia inmediatamente
si tuvi\'esemos $\g {\K}=\prod_{j=1}^h\g{({\K}_j)}$. 
Desafortunadamente, en general, para cualesquiera dos
campos $L_1$ and $L_2$, tenemos que
$\g{(L_1)}\g{(L_2)}\subsetneqq \g{(L_1L_2)}$,
Observaci\'on \ref{O12*.2.2.K}.

Ahora bien, en nuestro caso, tenemos que
$[{\K}_i:K]=l_i^{m_i}$ y $\mcd([{\K}_i:K],
[{\K}_j:K])=1$ para toda $i\neq j$. El conocimiento de $\g {\K}$
es una consecuencia inmediata del siguiente resultado.

\begin{proposicion}\label{P14.8.1.1.12}
Sean $L_i/K$, $i=1,2$ dos extensiones abelianas finitas tales que
$\mcd([L_1:K],[L_2:K])=1$. Entonces $\g{(L_1)}\g{(L_2)}=\g{(L_1
L_2)}$.
\end{proposicion}

Para probar la Proposici\'on \ref{P14.8.1.1.12}, primero 
probamos la siguiente proposici\'on.

\begin{proposicion}\label{P14.8.1.1.13}
Sea $L/K$ una extensi\'on abeliana finita y sea $l$
un n\'umero primo. Entonces $l|[L:K]\iff l|[\g L:K]$.
\end{proposicion}

\begin{proof}
$\Rightarrow)$ Es claro puesto que $L\subseteq \g L$.

\noindent
$\Leftarrow)$ Supongamos que $l|[\g L:K]$. Primero supongamos
que $L\subseteq \cicl N{}$. Sea $X$ el grupo de caracteres
de Dirichlet asociado a $L$ y sea $Y=\prod_{P\in R_T^+} X_P$.
Sea $\gex L$ el campo asociado a $Y$. Tenemos que
$\g L\subseteq \gex L$. De hecho $\g L=\gex L^+L$. Entonces $l|[\gex L:K]$.
Por tanto, existe $\chi\in Y$ de orden $l$. Sea $\chi
=\prod_{P\in R_T^+}\chi_P$ con $\chi_P\in X_P$ y $\chi
\neq 1$, $\chi^l=1$.

En particular existe $P\in R_T^+$ con $o(\chi_P)=l$. Sea
$\varphi\in X$ con
$\varphi_P=\chi_P$. Ahora, si $l\nmid [L:K]$, $\mcd([L:K],l)=1$. Sea
$m=o(\varphi)$. Entonces $l\nmid m$ y $\varphi^m=\prod_{
P\in R_T^+}\varphi_P^m=1$. Por tanto $\varphi_P^m=\chi_P^m=1$
y $\chi_P^l=1$. Se sigue que $\chi_P=1$. Esto contradice que
$o(\chi_P)=l$. Se sigue que $l|[L:K]$.

Para el caso genral, sea $\g L=\g{E^H}L$. Si $l|[\g L:K]$ entonces $l|
[\g{E^H}:K]$ o $l|[L:K]$. Si $l|[\g{E^H}:K]$ entonces $l|[\g E:K]$ as\'i,
por el caso ciclot\'omico, obtenemos que $l|[E:K]|[L:K]$. Por
tanto $l|[L:K]$. $\fin$
\end{proof}

\begin{corolario}\label{C14.8.1.1.14} Cuando $L$ es
ciclot\'omico, tenemos $l|[L:K]\iff l|[\gex L:K]$. $\fin$
\end{corolario}

En el caso ciclot\'omico, el hecho de que $\g{(L_1)}\g{(L_2)}
\subsetneqq \g{(L_1L_2)}$ es consecuencia del hecho de que,
en general, tenemos $L_1^+L_2^+\subsetneqq (L_1L_2)^+$.
Probaremos que, en nuestro caso, tenemos igualdad.

\begin{lema}\label{L14.8.1.1.15}
Sean $L_i$, $i=1,2$ dos campos ciclot\'omicos.
Entonces, si $\gcd ([L_1:K],[L_2:K])=1$, 
tenemos $(L_1L_2)^+=L_1^+L_2^+$.
\end{lema}

\begin{proof}
Puesto que $L_i^+\subseteq L_i$, $i=1,2$, 
se sigue que $\mcd([L_1^+:K],[L_2^+:K])=1$.

\begin{tiny}
\[
\xymatrix{
L_1\ar@{-}[dd]\ar@{-}[rrrr]&&&&L_1L_2\ar@{-}[dddd]\ar@{-}[dl]\\
&&&(L_1L_2)^+\ar@{-}[dl]\\
L_1^+\ar@{-}[rr]\ar@{-}[dd]&&L_1^+L_2^+\ar@{-}[dd]\\ \\
K=L_1^+\cap L_2^+\ar@{-}[rr]&&L_2^+\ar@{-}[rr]&&L_2
}
\]
\end{tiny}

Tenemos que $e_i:=e_{\infty}(L_i|K)=[L_i:L_i^+]$, $i=1,2$.
Por lo tanto $\gcd(e_1,e_2)=1$. Se sigue que $e_{\infty}(L_1L_2|
L_1^+L_2^+)=e_1e_2$ y $[L_1L_2:L_1^+L_2^+]=e_1e_2$.

Esto es, $L_1L_2/L_1^+L_2^+$ es totalmente ramificado en $\p$.
Por otro lado, $(L_1L_2)^+$ es el m\'aximo subcampo de
$L_1L_2$ donde $\p$ se descompone totalmente en $(L_1L_2)^+/K$
y es totalmente ramificada en $L_1L_2/(L_1L_2)^+$.
Se sigue que $e_{\infty}(
L_1L_2|(L_1L_2)^+)|e_1e_2=e_{\infty}(L_1L_2|L_1^+L_2^+)$. 
As\'i $(L_1L_2)^+\subseteq L_1^+L_2^+\subseteq
(L_1L_2)^+$ y $(L_1L_2)^+=L_1^+L_2^+$. $\fin$
\end{proof}

\begin{observacion}\label{O14.8.1.1.16}
En el caso de que $L_1$ y $L_2$ sean ciclot\'omicos, es posible que
$\mcd([L_1^+:K],[L_2^+:K])=1$ pero que $\mcd([L_1:K],
[L_2:K])\neq 1$.

Por ejemplo, sean $P,Q\in R_T^+$ distintos y de grado $1$. Entonces
$[\cicl P{}:K]=q-1=[\cicl Q{}:K]$ y $\cicl P{}^+=\cicl Q{}^+=K$.
Por tanto, si $L_1=\cicl P{}$ y $L_2=\cicl Q{}$, se tiene que
$\mcd([L_1^+:K],[L_2^+:K])=\mcd(1,1)=1$ pero $\gcd([L_1:K],
[L_2:K])=\gcd(q-1,q-1)=q-1>1$ para $q>2$.
\end{observacion}

\begin{corolario}\label{C14.8.1.1.17}
Si $\mcd([L_1:K],[L_2:K])=1$ con $L_i\subseteq \cicl {N_i}{}$,
$i=1,2$, entonces $\g{(L_1L_2)}=\g{(L_1)}\g{(L_2)}$.
\end{corolario}

\begin{proof}
Del Corolario \ref{C14.8.1.1.14} obtenemos
$\mcd([\gex{(L_1)}:K],[\gex{(L_2)}:K])=1$. Ahora bien, tenemos que
$\g{(L_1)}=\gex{(L_1)}^+ L_1$ y $\g{(L_2)}=\gex{(L_2)}^+ L_2$.
Por tanto
\begin{align*}
\g{(L_1)}\g{(L_2)}&=\gex{(L_1)}^+ L_1\gex{(L_2)}^+ L_2=
\gex{(L_1)}^+\gex{(L_2)}^+ L_1L_2\\
&=(\gex{(L_1)}\gex{(L_2)})^+ (L_1L_2)=(\gex{(L_1L_2)})^+(L_1L_2)
=\g{(L_1L_2)}. \tag*{$\fin$}
\end{align*}
\end{proof}

Los detalles del siguiente resultado se dejan al cuidado del lector.

\begin{proposicion}\label{P14.8.1.1.18}
Sean $A,B$ y $C$ tres campos de funciones globales tales que
$B/A$ y $C/A$ son extensiones de Galois finitas que satisfacen
$\mcd([B:A],[C:A])=1$. Entonces $B\cap C=A$. Sea $D=BC$.

Si $\pK_A$ es un divisor primo de $A$ y $\pK_B,\pK_C,\pK_D$
satisfacen $\pK_B\cap A=\pK_C\cap A=\pK_D\cap A=\pK_A$, entonces
\begin{gather*}
e_{\pK_A}(B|A)=e_{\pK_C}(D|C), \quad
f_{\pK_A}(B|A)=f_{\pK_C}(D|C), \quad
h_{\pK_A}(B|A)=h_{\pK_C}(D|C),\\
e_{\pK_A}(C|A)=e_{\pK_C}(D|B),\quad
f_{\pK_A}(C|A)=f_{\pK_C}(D|B),\quad
h_{\pK_A}(C|A)=h_{\pK_C}(D|B),
\end{gather*}
donde $e,f$ y $h$ denotan el \'indice de ramificaci\'on,
el grado de inercia y el grado de descomposici\'on respectivamente.
$\fin$
\end{proposicion}

\subsubsection*{Demostraci\'on de la Proposici\'on \ref{P14.8.1.1.12}}

Sean $L_i/K$ dos extensiones abelianas finitas, $i=1,2$ tales que
$\mcd([L_1:K],[L_2:K])=1$. Sean $E_i=L_iM\cap \cicl N{}$, $i=1,2$
donde $L_i\subseteq {_n\cicl N{}_m}$, $i=1,2$ y $M=L_nK_m$ como
de costumbre. Sean
$L=L_1L_2$ y $E=E_1E_2$. Entonces $E=LM\cap \cicl N{}$.

Ahora $\g{(L_i)}=L_i \g{(E_i)}^{H_i}$, $i=1,2$ con
$|H_i|=f_{\infty}(L_iE_i|L_i)$, $i=1,2$ y
$\g L=L \g E^H$ donde $|H|=f_{\infty}(LE|L)$.
\[
\xymatrix{
\cicl N{}\ar@{-}[d]\\
E\ar@{-}[rr]\ar@{-}[dd]&&LM=EM\ar@{-}[dl]\ar@{-}[dd]\\
&L\ar@{-}[dl]\\
K=M\cap \cicl N{}\ar@{-}[rr]&&M
}
\qquad
\xymatrix{
\\
\\
L\ar@{-}[r]\ar@{-}[d]&LM\ar@{-}[d]\\
L\cap M\ar@{-}[r]&M
}
\]
Tenemos que $[E:K]=[LM:M]$ y $[LM:M]=[L:L\cap M]|[L:K]$.
Por tanto $[E:K]|[L:K]$. An\'alogamente, se tiene que
$[E_i:K]|[L_i:K]$, $i=1,2$.

Se sigue que $\mcd([E_1:K],[E_2:K])=1$. De la Proposici\'on
\ref{P14.8.1.1.18} obtenemos que
$|H|=|H_1||H_2|$. Entonces, del Corolario \ref{C14.8.1.1.17},
obtenemos que
\begin{gather*}
[\g E:\g E^H]=|H|=|H_1||H_2|=[\g{(E_1)}:\g{(E_1)}^{H_1}]
[\g{(E_2)}:\g{(E_2)}^{H_2}]
\intertext{y tenemos, con $a=|H_1|$, $b=|H_2|$ y $ab=
|H_1||H_2|=|H|$, que}
\xymatrix{
\g{(E_1)}\ar@{-}[r]\ar@{-}[d]^a&\bullet\ar@{-}[d]^a
\ar@{-}[r]^{b\phantom{xxxxx}}&\g{(E_1)}\g{(E_2)}
=\g{(E_1E_2)}\ar@{-}[d]^a\ar@{-}[dl]_{ab}\\
\g{(E_1)}^{H_1}\ar@{-}[r]\ar@{-}[d] &\g{(E_1)}^{H_1}
\g{(E_2)}^{H_2}\ar@{-}[d]\ar@{-}[r]^b&\bullet
\ar@{-}[d]&\hspace{-1cm} \mcd(a,b)=1\\
K\ar@{-}[r]&\g{(E_2)}^{H_2}\ar@{-}[r]^b&\g{(E_2)}
}
\intertext{Se sigue que}
\begin{align*}
[\g{(E_1E_2)}:\g{(E_1)}^{H_1}\g{(E_2)}^{H_2}]&=
[\g{(E_1)}:\g{(E_1)}^{H_1}][\g{(E_2)}:\g{(E_2)}^{H_2}]\\
&=|H_1||H_2|=|H|
\end{align*}
\end{gather*}
puesto que $\g{(E_1E_2)}=\g{(E_1)}\g{(E_2)}$.
Por otro lado
$[\g{(E_1E_2)}:\g{(E_1E_2)}^{H}]=|H|$. Ya que
$\g{(E_1)}^{H_1}\g{(E_2)}^{H_2}\subseteq \g{(E)}^{H}$,
obtenemos que $\g{(E_1)}^{H_1}\g{(E_2)}^{H_2}=
\g{(E_1E_2)}^{H}$.
Por tanto
\begin{align*}
\g{(L_1L_2)}&=\g{(E_1E_2)}^{H}(L_1L_2)=\g{(E_1)}^{H_1}
\g{(E_2)}^{H_2}L_1L_2\\
&=\g{(E_1)}^{H_1}L_1\g{(E_2)}^{H_2}
L_2=\g{(L_1)}\g{(L_2)}.
\tag*{$\fin$}
\end{align*}

Como consecuencia, obtenemos el resultado principal
de esta subsecci\'on.

\begin{teorema}\label{T14.8.1.1.19}
Sea ${\K}/K$ una extensi\'on de Kummer finita de orden
$n=l_1^{m_1}\cdots l_s^{m_s}$ con $l_1,\ldots,l_s$
n\'umeros primos distintos. Sea
${\K}={\K}_1\cdots {\K}_s$ con $[{\K}_j:K]=
l_j^{m_j}$, $1\leq j\leq s$. Entonces 
\[
\g {\K}=\prod_{j=1}^s\g{({\K}_j)}
\]
donde cada $\g{({\K}_j)}$, $1\leq j\leq s$ se calcula en los
Teoremas {\rm{\ref{T14.8.1.1.7}}} y {\rm{\ref{T14.8.1.1.10}}}.
\end{teorema}

\begin{proof}
Como consecuencia de la
Proposici\'on \ref{P14.8.1.1.12} obtenemos que
$\g {\K}=\g{({\K}_1)}\cdots \g{({\K}_s)}$. 
La descripci\'on expl\'icita de cada $\g{({\K}_i)}$ es precisamente
el contenido de los Teoremas
\ref{T14.8.1.1.7} and \ref{T14.8.1.1.10}. $\fin$
\end{proof}

\subsection[Descripci\'on expl\'icita de $p$--extensiones abelianas]
{Descripci\'on expl\'icita de campos de g\'eneros
de $p$--extensiones abelianas finitas de $K$}\label{S3.A}

En esta secci\'on usaremos las notaciones de vectores de Witt
que fueron introducidas en el Cap\'itulo \ref{Ch9'}. En particular,
para un anillo conmutativo con unidad $R$, $W_v(R)$ denota
al anillo de vectores de Witt de longitud $v$ con coeficientes
en $R$.

Sea $\K/K$ una $p$--extensi\'on
abeliana finita. Recordemos que $K={\ma F}_q(T)$, digamos
$q=p^{\nu}$. Supondremos que 
${\ma F}_{p^u}\subseteq \F$, esto es, $u\mid \nu$. 

Entonces tenemos que
\[
\Gal(\K/K)\cong \big({\ma Z}/p^{\alpha_1}{\ma Z}\big)\times
\cdots\times \big({\ma Z}/p^{\alpha_u}{\ma Z}\big) \quad
\text{con}\quad 1\leq \alpha_1\leq \cdots\leq \alpha_u=v.
\]
Existen $\vec w_1, \ldots, \vec w_u\in W_v(\bar{K})$
tales que $\vec w_i^p\Witt - \vec w_i=
\vec\xi_i\in W_v(K)$, con
$\K=K(\vec w_1,\cdots,\vec w_v)$. Tambi\'en tenemos que
existe $\vec y_0\in W_v(\bar{K})$ tal que 
$\K=K(\vec y_0)$ con
\[
\vec y_0^{p^u}\Witt - \vec y_0=\vec\xi_0 \quad\text{para alguna}
\quad \vec \xi_0\in W_v(K)
\]
(ver Teorema \ref{T8.4.pea}). 
Aqu\'i $\bar{K}$ denota una cerradura algebraica de $K$.

Sean $P_1,\ldots,P_r\in R_T^+$ los primos finitos de $K$
ramificados en $\K$.
Del Teorema \ref{T8.7.pea} se sigue que podemos
descomponer $\vec\xi_0$ como
\begin{gather}\label{Eq3.0.A}
\vec \xi_0={\vec\delta}_1
\Witt + \cdots \Witt + {\vec\delta}_r \Witt + \vec\gamma,
\end{gather}
$\delta_{i,j}=\frac{Q_{i,j}}{P_i^{e_{i,j}}}$, $e_{i,j}\geq 0$, $Q_{i,j}\in R_T$
y si $e_{i,j}>0$, entonces $e_{i,j}=
\lambda_{i,j}p^{m_{i,j}}$, $\mcd (\lambda_{i,j},p)=1$,
$0\leq m_{i,j}< n$, $\mcd(Q_{i,j},P_i)=1$ y
$\deg (Q_{i,j})<\deg (P_i^{e_{i,j}})$, y $\gamma_j=f_j(T)\in R_T$ con
$\deg f_j=\nu_j p^{m_j}$ y $\mcd(q,\nu_j)=1$, $0\leq m_j<n$
cuando $f_j\not\in \F$.

Si el \'indice de ramificaci\'on de
$P_i$ es $p^{a_i}<p^v$, podemos escribir
$\vec \delta_i=(\delta_{i,1},\ldots,\delta_{i,v})=(0,\ldots,0,\delta_{
i,(v-a_i+1)},\ldots,\delta_{i,v})$. En particular $\p$ 
se descompone totalmente en $K(\vec y_i)/K$, donde
$\vec y_i^{p^u}\Witt - \vec y_i=\vec \delta_i$
(ver Observaci\'on \ref{O9'.8.28}).

Sea $\vec z^{p^u}\Witt - \vec z=\vec \gamma$. 
En $K(\vec z)/K$ el \'unico posible primo ramificado es
$\p$. Notemos que si 
\[
\vec y=\vec y_1\Witt +\cdots\Witt + \vec y_r, \quad \text{entonces}\quad 
\vec y^{p^u}\Witt -\vec y=\vec\xi_0\Witt - \vec \gamma=\vec\delta_1\Witt +
\cdots\Witt + \vec\delta_r
\]
y que $\p$ se descompone totalmente en $K(\vec y)/K$.

El primer resultado principal de esta secci\'on es

\begin{teorema}\label{T3.1.A}
Con las notaciones de la Secci\'on {\rm{\ref{S2.A}}}, sea
$E=\K M\cap K(\Lambda_N)$. Entonces
$E=K(\vec y)$, 
\[
\g E=K(\vec y_1,\ldots,
\vec y_r)\quad \text{y}\quad \g \K=K(\vec y_1,\ldots,\vec y_r, \vec z).
\]
\end{teorema}

\begin{proof} De la correspondencia de Galois entre $\cicl N{}/K$
y $M/K$, tenemos
$EM=\K M$. Probar que $E=K(\vec y)$ es equivalente
a mostrar que $K(\vec y)M=\K M$ 
puesto que $K(\vec y)\subseteq K(\Lambda_N)$.

Ahora, $K(\vec z)\subseteq M$ puesto que 
$M=L_n{\ma F}_{q^m}(T)$ codifica toda la inercia
y toda la ramificaci\'on, la cual totalmente salvaje, de
$\p$. Tenemos
\begin{gather*}
K(\vec y)M=K(\vec y)K(\vec z)M\supseteq 
K(\vec y\Witt +\vec z)M=\K M.
\intertext{Rec\'iprocamente,}
\K M=\K K(\vec z)M=K(\vec y_0)K(\vec z)
M\supseteq K(\vec y_0\Witt - \vec z)M=
K(\vec y)M.
\intertext{Por tanto}
\K M=K(\vec y)M\quad \text{y} \quad E=K(\vec y).
\end{gather*}

Del Theorem \ref{T4.3.A}
obtenemos que $\g E=K(\vec y_1,\ldots,\vec y_r)$.
Finalmente
\begin{align*}
\g \K&=\g E \K=K(\vec y_1,\ldots, 
\vec y_r)K(\vec y_0)=K(\vec y_1,\ldots,\vec y_r)
K(\vec y_0\Witt -\vec y_1\Witt -\cdots\Witt -\vec y_r)\\
&=K(\vec y_1,
\ldots,\vec y_r)K(\vec z)=K(\vec y_1,\ldots,\vec y_r,\vec z).
\end{align*}
Esto termina la demostraci\'on. $\fin$
\end{proof}

\begin{observacion}\label{R3.2.A} Notemos que la
demostraci\'on del Teorema \ref{T3.1.A}
funciona a\'un en el caso de que
$\vec \delta_i$ y $\vec \gamma$ no est\'an en 
forma reducida descrita arriba. Lo \'unico que 
necesitamos es que en cada extensi\'on
$\vec y_i^{p^u}
\Witt -\vec y_i=\vec \delta_i$, $1\leq i\leq r$
y $\vec z^{p^u}\Witt -\vec z=
\vec \gamma$ hay a lo m\'as primo ramific\'andose.
\end{observacion}

Del Teorema \ref{T2.2.A}, los casos de 
Artin--Schreier, de extensiones Witt, y de $p$--extensiones
elementales abelianas son consecuencia 
inmediata del Teorema \ref{T3.1.A}.

\begin{corolario}\label{C3.3.A} Sea $K=\F(T)$.
\l
\item (Extensiones de Artin--Schreier) Sea $\K=K(y)$ con 
\[
y^p-y=\alpha=\sum_{i=1}^r\frac{Q_i}{P_i^{e_i}} + f(T),
\]
donde $P_i\in R_T^+$, $Q_i\in R_T$, 
$\mcd(P_i,Q_i)=1$, $e_i>0$, $p\nmid e_i$, $\deg Q_i<
\deg P_i^{e_i}$, $1\leq i\leq r$, $f(T)\in R_T$,
con $p\nmid \deg f$ cuando $f(T)\not\in \F$.

Entonces 
\[
\K_{{\eu {ge}}}=K(y_1,\ldots,y_r,\beta),
\]
donde
$y_i^p-y_i=\frac{Q_i}{P^{e_i}}$, $1\leq i\leq r$ y
$\beta^p-\beta=f(T)$.

\item (Extensiones de Witt) Sea $\K=K(\vec y)$ donde
\[
\vec y^p\Witt -\vec y=\vec \beta={\vec\delta}_1\Witt 
+ \cdots \Witt + {\vec\delta}_r
\Witt + \vec\mu,
\]
con $\delta_{i,j}=\frac{Q_{i,j}}{P_i^{e_{i,j}}}$, 
$e_{i,j}\geq 0$, $Q_{i,j}\in R_T$
y si $e_{i,j}>0$, entonces $p\nmid e_{i,j}$, $\mcd(Q_{i,j},P_i)=1$ y 
$\deg (Q_{i,j})<\deg (P_i^{e_{i,j}})$, y $\mu_j=f_j(T)\in R_T$ con
$p\nmid \deg f_j$ cuando $f_j\not\in \F$.

Entonces
\[
\K_{{\eu {ge}}}=K({\vec y}_1,\ldots,{\vec y}_r,\vec z),
\]
donde ${\vec y}_i^p\Witt -{\vec y}_i ={\vec \delta}_i$, $1\leq i\leq r$
y ${\vec z}^p\Witt -\vec z=\vec\mu$. 

\item ($p$--extensiones elementales abelianas)
Supongamos que ${\ma F}_{p^u}\subseteq \F$. 
Sea $\K=K(y)$ con
\[
y^{p^u}-y=\alpha=\sum_{i=1}^r\frac{Q_i}{P_i^{e_i}} + f(T),
\]
donde $P_i\in R_T^+$, $Q_i\in R_T$ y $f(T)\in \F[T]$.

Entonces 
\[
\g \K=K(y_1,\ldots, y_r, z),
\]
donde $y_i^{p^u}-y_i=\frac{Q_i}{P_i^{e_i}}$,
$1\leq i\leq r$ y $z^{p^u}-z=f(T)$. $\fin$

\end{list}
\end{corolario}

\begin{ejemplo}\label{E5.12}
Sea $K={\ma F}_3(T)$ y $\K=K(\vec y)$ donde
 $\vec y^3\Witt - \vec y=\vec {\beta}=
\big(\frac{1}{T}+1,\frac{1}{T+1}+T\big)$. Entonces la descomposici\'on
prescrita en el Teorema \ref{T5.3.1} es:
\[
\vec \beta=\Big(\frac{1}{T},\frac{T+1}{T^2}\Big)\Witt +
\Big(0,\frac{1}{T+1}\Big)\Witt + \big(1,T\big).
\]
Por lo tanto, si $\vec y_1^3\Witt -\vec y_1=\vec\delta_1
=\big(\frac{1}{T},\frac{T+1}{T^2}\big)$,
$\vec y_2^3\Witt -\vec y_2=\vec \delta_2=\big(0,\frac{1}{T+1}\big)$ y
$\vec z^3\Witt -\vec z=\vec \mu=\big(1,T\big)$, entonces
$\K_{\eu {ge}}=K(\vec y_1,\vec y_2,\vec z)$.
\end{ejemplo}

\section{Extensiones abelianas finitas generales de $K$}\label{S5.A}

Hasta ahora hemos dado la descripci\'on expl\'icita del campo
de g\'eneros de $p$--extensiones abelianas $\K$ de $K=\F(T)$
y donde hemos considerado el caso tal que
${\ma F}_{p^u}\subseteq \F$, $\K=K(\vec y)$ y $\vec y$
est\'a dada por una ecuaci\'on de la forma
$\vec y^{p^u} \Witt - \vec y=\vec \beta
\in W_m(K)$. Cuando ${\ma F}_{p^u}\nsubseteq \F$ 
el campo $\K$ no puede ser dado por este tipo de ecuaciones. 

En esta secci\'on daremos la descripci\'on expl\'icita de
$\g\K$ donde $\K/K$ es una extensi\'on finita de grado 
$t$ con $\mcd(t,q-1)=1$. El caso $t\mid q-1$ est\'a 
parcialmente considerado en la Subsecci\'on \ref{S4.1.A}.

Seguimos usando las notaciones de la Secci\'on \ref{S12*.3}.

\begin{observacion}\label{R5.1.A} Para cualquier
extensi\'on abeliana $\K/K$ de grado $t$ con $\mcd(t,q-1)=1$,
tenemos que si $E=\K M\cap K(\Lambda_N)$, 
entonces $[E:K]\mid t$ (ver Ecuaci\'on (\ref{Eq4.A})). Si $X$ 
es el conjunto de caracteres de Dirichlet de $E$, tenemos que
$\mcd(|X|,q-1)=\mcd([E:K],q-1)=1$. Puesto que para $\chi\in X$
y cualquier $P\in R_T^+$, tenemos que $\chi_P^{|X|}=1$, 
obtenemos que $\mcd([\g E:K],q-1)=1$. En particular
$H=\{1\}$. Por lo tanto 
\[
\g \K=\g E \K.
\]
\end{observacion}

En general si $\K_1$ y $\K_2$ son dos extensiones finitas de
$K$ tenemos
\[
\g{(\K_1)}\g{(\K_2)}\subseteq \g{(\K_1\K_2)},
\]
pero podr\'iamos tener
$\g{(\K_1)}\g{(\K_2)}\subsetneqq \g{(\K_1\K_2)}$,
ver Observaciones \ref{O12*.2.2.K} y \ref{O3.6.B}. 
Probaremos que para
cualquier extensi\'on abeliana de $K$ de grado primo
relativo a $q-1$ tendremos la igualdad.
En particular si $\K_1$ y $\K_2$ son dos $p$--extensiones
abelianas finitas de $K$, tendremos igualdad.

Para un subcampo $\K\subseteq K(\Lambda_N)$
para alg\'un $N\in R_T$, denotemos por
$\g {\K'}$ a la m\'axima extensi\'on abeliana de $K$
contenida en $K(\Lambda_N)$, no ramificada en los
primos finitos. Tenemos (ver Observaci\'on \ref{R4.0.A})
\begin{gather}\label{Eq3.1.A}
\g \K=(\g {\K'})^D,
\end{gather}
donde $D$ es el grupo descomposici\'on de cualquier
elemento de $\S \K$ en $\g {\K'}/\K$.

Consideramos $\K_i\subseteq K(\Lambda_N)$, $i=1,2$ 
y sean $X_i$ los grupos de los caracteres de Dirichlet
asociados a $\K_i$. Por lo tanto 
$Y=X_1X_2=\langle X_1,X_2\rangle$ es el grupo de los
caracteres de Dirichlet asociado a 
$L=\K_1\K_2$. Sea $P\in R_T^+$. Es f\'acil verificar que
\begin{gather*}
\langle X_1,X_2\rangle_P=\langle (X_1)_P,(X_2)_P\rangle,
\intertext{por lo que obtenemos que}
\prod_{P\in R_T^+}Y_P=\prod_{P\in R_T^+}
\langle X_1,X_2\rangle_P=\Big(\prod_{P\in R_T^+}(X_1)_P\Big) 
\Big(\prod_{P\in R_T^+} (X_2)_P\Big).
\intertext{Se sigue que}
\g {(\K_1)}' \g {(\K_2)}'=\g{(\K_1\K_2)}'.
\end{gather*}

Hemos probado

\begin{proposicion}\label{P3.4.A} Para $\K_i\subseteq K(\Lambda_N)$, 
$i=1,2$, tenemos
$$
\g {(\K_1)}' \g {(\K_2)}'=\g{(\K_1\K_2)}'. \eqno{\fin}
$$
\end{proposicion}

\begin{corolario}\label{C3.5.A} Sean $\K_i\subseteq K(\Lambda_N)$, 
$i=1,2$ tales que $\K_1/K$ y $\K_2/K$ 
son extensiones abelianas finitas de grado primo relativo
a $q-1$. Entonces
$\g {(\K_1)} \g {(\K_2)}=\g{(\K_1\K_2)}$.
\end{corolario}

\begin{proof} Puesto que los grupos de descomposici\'on de
$\p$ en $\K_1/K$, en $\K_2/K$ y en
$\K_1\K_2/K$ son el grupo trivial, se sigue de
la Ecuaci\'on (\ref{Eq3.1.A}) que $\g {(\K_i)}=
\g{(\K_i)}'$, $i=1,2$ y $\g{(\K_1\K_2)}=\g{(\K_1\K_2)}'$. 
El resultado se sigue de la Proposici\'on \ref{P3.4.A}. $\fin$
\end{proof}

\begin{corolario}\label{C3.6.A}
Sean $\K_i/K$, $i=1,2$ dos extensiones abelianas finitas de
grado primo relativo a $q-1$. Entonces
\[
\g {(\K_1)} \g {(\K_2)}=\g{(\K_1\K_2)}.
\]
\end{corolario}

\begin{proof}
Sea $\F={\ma F}_{p^{\nu}}$,
$\K_i\subseteq L_nK(\Lambda_N) 
{\ma F}_{p^{\nu m}}(T)$, $i=1,2$, y sea
$M:=L_n{\ma F}_{p^{\nu m}}(T)$. Establecemos
$E_i:=\K_iM\cap K(\Lambda_N)$, $i=1,2$ y $E:=\K_1\K_2M
\cap K(\Lambda_N)$. Usando la correspondencia de Galois
puede ser probado que $E=E_1E_2$. 

Del Corolario \ref{C3.5.A} tenemos
$\g E=\g {(E_1)} \g {(E_2)}$. Por lo tanto
\begin{align*}
\g {(\K_1)} \g {(\K_2)}&=\g{(E_1)}\K_1\cdot \g{(E_2)}\K_2=
\g {(E_1)}\g{(E_2)} \cdot \K_1\K_2\\
&=\g E \cdot \K_1\K_2=
\g{(\K_1\K_2)}.
\end{align*}
De esta forma $\g {(\K_1)} \g {(\K_2)}=\g{(\K_1\K_2)}$. $\fin$
\end{proof}

\begin{corolario}\label{C3.7.A} Sean
$\K_i/K$, $i=1,2$ dos $p$--extensiones abelianas finitas.
Entonces
$$
\g {(\K_1)} \g {(\K_2)}=\g{(\K_1\K_2)}. \eqno{\fin}
$$
\end{corolario}

Como consecuencia de lo anterior, obtenemos una
descripci\'on del campo de g\'eneros de una $p$--extensi\'on
abeliana finita de $K$.

\begin{corolario}\label{C3.8.A} Sea $\K/K$ una $p$--extensi\'on
abeliana con grupo de Galois $\Gal(\K/K)=G
\cong G_1\times\cdots\times G_s$ con $G_i\cong
{\ma Z}/p^{\alpha_i}{\ma Z}$, $1\leq i\leq s$. Sea $\K$
la descomposici\'on $\K=\K_1\cdots \K_s$ tal que
$\Gal(\K_i/K)\cong G_i$. Sean $P_1,\ldots, P_r$
los primos finitos ramificados en $\K/K$. Sea
$\K_i=K(\vec w_i)$ dado por la ecuaci\'on
\begin{gather*}
\vec w_i^p\Witt - \vec w_i=\vec \xi_i,\quad 1\leq i\leq s.
\intertext{Escribamos cada
$\vec\xi_i$ como en Ecuaci\'on {\rm{(\ref{Eq3.0.A})}} esto es,}
\vec \xi_i={\vec\delta}_{i,1}
\Witt + \cdots \Witt + {\vec\delta}_{i,r} \Witt + \vec\gamma_i,
\intertext{tal que todas las componentes de 
$\vec\delta_{i,j}$ se escriben de tal manera que el
grado del numerador es menor al grado del denominador,
el soporte del denominador es a lo m\'as $\{P_j\}$ 
y las componentes de
$\vec\gamma_i$ son polinomios. Sea}
\vec w_{i,j}^p\Witt -\vec w_{i,j}=\vec\delta_{i,j},\quad
1\leq i\leq s,\quad 1\leq j\leq r
\intertext{y}
\vec z_i^p\Witt -\vec z_i=\vec \gamma_i, \quad 1\leq i\leq s.
\intertext{Entonces}
\g \K=K\big(\vec w_{i,j},\vec z_i\mid 1\leq i\leq s, 1\leq j\leq r\big).
\end{gather*}
\end{corolario}

\begin{proof} El resultado es consecuencia de la Observaci\'on
\ref{R3.2.A}, del Corolario \ref{C3.3.A} (b) y 
del Corolario \ref{C3.7.A}. $\fin$
\end{proof}

Como siguiente paso, consideramos una extensi\'on c\'iclica
$\K/K$ de grado $t$ tal que $\mcd(t,p(q-1))=1$. 

En este caso,
tenemos que $E=\K M\cap K(\Lambda_N)$ satisface que
$[E:K]$ es primo relativo a $q-1$. De aqu\'i tendremos que
$\g {E'}=\g E$ y $\g \K=\g E \K$. 
Por tanto, debemos describir $\g E$.

\begin{proposicion}\label{P3.9.A}Sea
 $E\subseteq K(\Lambda_N)$ una extensi\'on c\'iclica de
 $K$ de grado $t$ primo relativo a $p(q-1)$. Sean
 $P_1,\ldots,P_r\in R_T^+$ los primos finitos de $K$
 ramific\'andose en $E$. Entonces
 \[
\g E=\prod_{j=1}^r F_j,
\]
donde $K\subseteq F_j\subseteq K(\Lambda_{P_j})$
es el campo de grado $a_j$ sobre $K$ con
$o(\chi_{P_j})=a_j$, y $\chi$ es el caracter asociado a $E$.
\end{proposicion}

\begin{proof} Se sigue del hecho de que $X=\langle \chi\rangle$ 
es el grupo de caracteres de Dirichlet asociado a $E$,
$\g E$ es el campo correspondiente a
$\prod_{j=1}^r X_{P_j}$, $X_{P_j}=
\langle\chi_{P_j}\rangle$ y $F_j$ es el campo
asociado a $\chi_{P_j}$. $\fin$
\end{proof}

Tenemos nuestro resultado final de esta secci\'on.

\begin{teorema}\label{T3.10.A} Sea
$\K/K$ una extensi\'on abeliana de grado
$t$ con $\mcd (t,q-1)=1$. Sean $P_1,\ldots, P_r\in R_T^+$ 
los primos en $K$ ramific\'andose en $\K$.
Sea $E=\K M\cap K(\Lambda_N)=E_0E_1\cdots E_s$
donde $E_i/K$ es una extensi\'on c\'iclica de grado
$t_i$, $\mcd(t_i,p(q-1))=1$,
$1\leq i\leq s$ y $E_0/K$ es una $p$--extensi\'on 
abeliana. Entonces
\[
\g \K=\g E \K\quad\text{donde}\quad
\g E=\g {(E_0)}\g {(E_1)}\cdots \g {(E_s)},
\]
$\g {(E_0)}$ est\'a dado por el Corolario {\rm{\ref{C3.8.A}}} 
y $\g {(E_i)}=\prod_{j=1}^s F_{i,j}$ est\'a dado por la 
Proposici\'on {\rm{\ref{P3.9.A}}}, $1\leq i\leq s$. 

M\'as a\'un, si $[F_{i,j}:K]=b_{i,j}$, entonces si
$F_j:=\prod_{i=1}^s F_{i,j}$ es el subcampo de
$K(\Lambda_{P_j})$ de grado $b_j:=\mcm[b_{i,j},
1\leq i\leq s]$ sobre $K$. Finalmente tenemos
$$
\g \K=\g{(E_0)}\Big(\prod_{j=1}^r F_j\Big) \K. \eqno{\fin}
$$
\end{teorema}

\section{Campos de g\'eneros de extensiones
no abelianas de $K$}\label{S1.B}

Como mencionamos antes, 
A. Fr\"ohlich \cite{Fro59-1, Fro59-2} introdujo
el concepto de campos de g\'eneros para extensiones
no necesariamente abelianas de campos num\'ericos.
Fr\"ohlich defini\'o el campo de g\'eneros $\g\K$ de un
campo num\'erico arbitrario $\K/{\ma Q}$ como
$\g \K := \K K^{\ast}$ donde $K^{\ast}$ 
es el m\'aximo campo num\'erico abeliano tal que $\K K^{\ast}/
\K $ es no ramificada. Se tiene que $K^{\ast}/{\ma Q}$
es el m\'axima extensi\'on abeliana de ${\ma Q}$
contenida en $\g\K$.  

En esta parte consideramos un campo de funciones
congruente $\K$ con campo de constantes
conteniendo a $\F$ y tal que
$\K$ es una extensi\'on finita y separable de $K=\F(T)$.
Usamos, como antes, la definici\'on de Campo de Clase
de Hilbert de $\K$ para hallar $\g\K$ como $\g\K=\K K^{\ast}$,
donde $K^{\ast}$ es la m\'axima extensi\'on abeliana
de $K$ tal que $\K \*K$ est\'a contenido en el campo de
clase de Hilbert de $\K$. Tenemos que $\*K$ es
el compuesto de dos campos
$K_1^{\ast}$ y $K_2^{\ast}$ con $\mcd([\*{K_1}:K],p)=1$
y $\*{K_2}/K$ es una $p$--extensi\'on abeliana finita.

Probaremos que $K_1^{\ast}$ est\'a contenido en el
compuesto de una extensi\'on de constantes con
campos $F_{\P}$ tales que
$K\subseteq  F_{\P} \subseteq K(\Lambda_P)$, $\P$
es totalmente ramificado en
$F_{\P}/K$ y donde tenemos que $\P$ recorre 
el conjunto del primos finitos ramificados de $\K/K$.
Aqu\'i $P$ es el polinomio m\'onico irreducible en $T$
asociado a $\P$ y $K(\Lambda_P)$ es el
campo de funciones $P$--\'esimo ciclot\'omico.
El campo $K_2^{\ast}$ codifica la ramificaci\'on
salvaje de la extensi\'on $K^{\ast}/K$.
La mayor dificultad al considerar $K_1^{\ast}$
es la descomposici\'on de los primos infinitos.

\subsection{Algunas notaciones}\label{S2.B}

En esta parte, $\K$ denota una extensi\'on 
finita y separable de $K=\F(T)$. ${\ma P}_K^*$ es el conjunto 
de divisores primos finitos de $K$

Primero notemos que es posible tener
$\g \K =\K  E$ para varios subcampos
$E\subsetneqq K^{\ast}$. Estamos interesados en 
$K^{\ast}$ mismo.
En particular tenemos $K^{\ast}=\g {K^{\ast}}$ 
y puesto que $K^{\ast}/K$ es abeliana, la descripci\'on
de $K^{\ast}$ la dimos en el Teorema \ref{T2.1.A}.

Sean $\P_1,\ldots, \P_s,\P_{s+1},\ldots, \P_t$ 
los primos finitos de $K$ ramificados en $\K$. Sea
$P_i\in R_T^+$ tal que el divisor $(P_i)_K$ es 
$(P_i)_K=\frac{\P_i}{\p^{\deg P_i}}$ para $1\leq i\leq t$.
Para un primo $\P\in{\ma P}_K$, si $\con_{K/\K}\P=\producto$,
denotamos 
\begin{gather}\label{Eq0.1.B}
e_{\P}=\mcd(e_1,\ldots,e_r)= p^{u_{\P}} e_{\P}^{(0)}, \quad u_{\P}\geq
0,\quad \mcd(p,e^{(0)}_{\P})=1.
\end{gather}
Supondremos que
$p\nmid e_{\P_i}$ para $1\leq i\leq s$ y $p\mid e_{\P_j}$ para
$s+1\leq j\leq t$. Esto es, $u_{\P_i}=0$ para
$1\leq i\leq s$ y $u_{\P_j}
\geq 1$ para $s+1\leq j\leq t$.

Una de las herramientas principales en esta parte es
el Lema de Abhyankar (Teorema \ref{T2.3.B}).

Enunciamos un caso particular del Teorema \ref{T2.1.A}
que necesitaremos en la Subsecci\'on \ref{S4.B}.

\begin{teorema}\label{T2.6.B}
Sea $\K/{\ma F}_q$ una extensi\'on abeliana finita de
$K$ donde $\p$ es moderadamente ramificada. Sean $N\in R_T$
y $m\in{\ma N}$ tales que $\K \subseteq K(\Lambda_N)
{\ma F}_{q^m}$.
Sea $E_{{\eu {ge}}}$ el campo de g\'eneros de
$E:=K(\Lambda_N)\cap
\K{\ma F}_{q^m}$. 
Sea $H_1$ el subgrupo que corresponde al
grupo de descomposici\'on de los primos infinitos de $\K$
en $\g E/\K$ bajo la correspondencia de Galois. 
Entonces el campo de g\'eneros de $\K$
est\'a dado por
\begin{gather*}
\K_{{\eu {ge}}}= \g {E^{H_1}}\K .
\tag*{\fin}
\end{gather*}
\end{teorema}

\subsection{Caso general}\label{S4.B}

Primero probaremos algunos resultados generales.

\begin{proposicion}\label{P3.1.B} 
Sea $E/K$ es una extensi\'on abeliana finita moderadamente
ramificada. Para un primo finito ${\mathcal P}\in {\ma P}_K$, sea
$\con_{K/\K}{\mathcal P}=\producto$ y sea
$e^{\ast}_{\P}$ el \'indice de ramificaci\'on de $\P$ en $E/K$. 
Entonces $\K E/\K $ es no ramificada en los primos que
dividen a $\P$ si y solamente si
$e^{\ast}_{\P}\mid e_{\P}$, donde $e_{\P}$ esta dado por
{\rm{(\ref{Eq0.1.B})}}.
\end{proposicion}

\begin{proof} Sea ${\eu P}$ un primo en $\K E$ sobre
 $\P={\eu P}\cap K$. Por tanto
${\eu P}\cap \K ={\eu p}_i$ para alg\'un $i$.
Del Lema de Abhyankar, obtenemos
\begin{align*}
e({\eu P}|\P)&=\mcm[e({\eu p}_i|\P),e({\eu P}\cap E|\P)]
= \mcm [e_i,e^{\ast}_{\P}]\\
&=e({\eu P}|{\eu p}_i)e({\eu p}_i|\P)
=e({\eu P}|{\eu p}_i) e_i,
\end{align*}
donde $e_i=e(\pK_i|\P)$.

Por lo tanto ${\eu P}$ es no ramificado en $\K E/\K 
\iff e({\eu P}|{\eu p}_i)=1 \iff \mcm[e_i,e^{\ast}_{\P}]=
e_i\iff e_{\P}^{\ast}\mid e_i$. De esto 
se sigue el resultado. $\fin$
\end{proof}

Consideremos la conorma de $K$ en $\K$
del primo infinito de $K$:
\begin{gather}\label{Eq3.2'.B}
\con_{K/\K} \p=\pinfty. 
\end{gather}
Sean $t_i$ el grado de ${\eu p}_{i,\infty}$
(con respecto a $\F$) y
\begin{gather}\label{Eq3.2''.B}
t_0:=\mcd (t_1,\ldots,t_{r_{\infty}}).
\end{gather} 

Tenemos el resultado an\'alogo al Lema \ref{L4.1.1}.

\begin{proposition}\label{P3.2.B} El campo de constantes
de $\g \K $ es ${\ma F}_{q^{t_0}}$.
\end{proposition}

\begin{proof} Por conveniencia, supondremos que el
campo de constantes de $\K$ es $\F$. Ver la 
Observaci\'on \ref{O.3.2(1).B} para la demostraci\'on
del caso general.

Consideremos la extensi\'on de constantes
$\K{\ma F}_{q^m}/\K $ con $m\geq 1$. 
Tenemos que ${\eu p}_{i,\infty}$ se descompone
en $\mcd(t_i,m)$ factores (Teorema \ref{T6.1.4}).
Por lo tanto ${\eu p}_{i,\infty}$ se descompone totalmente
en $\K{\ma F}_m \iff \mcd(t_i,m)=m\iff m\mid t_i$.
Se sigue que los primos infinitos de $\K$ se decomponen
totalmente en $\K{\ma F}_m \iff m\mid t_0$. De
esta forma obtenemos que ${\ma F}_{q^{t_0}}$ 
es el campo de constantes $\g \K $. $\fin$
\end{proof}

\begin{observacion}\label{O.3.2(1).B}
Cuando el campo de constantes de $\K$ es
${\ma F}_{q^{\nu}}$, entonces $\nu\mid t_0$ y
$\nu\mid t_i$ y $\nu\mid m$. Poniendo, 
$t_0=\nu t_0'$, $t_i=\nu t_i'$ y $m=\nu m'$, todo
el argumento de la Proposici\'on \ref{P3.2.B} es 
aplicable para $t_0',t_i'$ y $m'$ con la conclusi\'on
de que el campo de contantes de $\K$ es
${\ma F}_{q^{\nu t_0'}}={\ma F}_{q^{t_0}}$.
\end{observacion}

\begin{proposicion}\label{P3.3.B} Sea $\P$ 
un divisor primo de $K$ de grado
$d$ y tal que satisface $\P\neq \p$ y
$p\nmid e_{\P}$. Sea $\g \K =\K K^{\ast}$ y sea
$e_{\P}^{\ast}$ el \'indice de ramificaci\'on de
$\P$ in $K^{\ast}/K$.  Entonces $\mcd\big(e_{\P},
\frac{q^d-1}{q-1}\big)\mid e_{\P}^{\ast}$ 
y $e_{\P}^{\ast}\mid e_{\P}$.
\end{proposicion}

\begin{proof} De la Proposici\'on \ref{P3.1.B} tenemos
que $e^{\ast}_{\P} \mid 
\mcd (e_1,\ldots,e_r)=e_{\P}$. M\'as a\'un $\P$ 
es totalmente ramificado en
$K(\Lambda_P)/K$, donde $(P)_K=\frac{\P}{\p^{\deg P}}$ y
$\p$ se descompone totalmente en
$K(\Lambda_P)^{\ff}/K$. El grado de la extensi\'on
$K(\Lambda_P)^{\ff}/K$ es $(q^d-1)/(q-1)$. 
Sea $S$ el subcampo  $K\subseteq
S\subseteq K(\Lambda_P)$ de grado
$\mcd\big(e_{\P},\frac{q^d-1}{q-1}\big)$.
Entonces por la Proposici\'on \ref{P3.1.B} 
tenemos que $S$ satisface que 
$\K S/\K $ es no ramificado y los primos infinitos de $\K$
se descomponen totalmente en
$\K S/\K $ puesto que $\p$ se descompone totalmente
en  $S/K$. Por lo tanto 
$\K S\subseteq \K K^{\ast}$, $S\subseteq K^{\ast}$ y
$\mcd\big(e_{\P},\frac{q^d-1}{q-1}\big)\mid e_{\P}^{\ast}$. $\fin$
\end{proof}

Sea $G:=\Gal(K^{\ast}/K)$ y sea
$G_p$ el $p$--subgrupo de Sylow de $G$.
Entonces $G=G_0\times G_p$ con $p\nmid |G_0|$.
Por tanto tenemos la descomposici\'on
\begin{gather}\label{Eq3.1''.B}
K^{\ast}=K_1^{\ast}K_2^{\ast},\quad 
K_1^{\ast}\cap K_2^{\ast}=K, \quad
G_0=\Gal(K_1^{\ast}/K),\quad G_p=\Gal(K_2^{\ast}/K).
\end{gather}

La extensi\'on $K_1^{\ast}/K$ es moderadamente
ramificada y $K_2^{\ast}/K$ es una $p$--extensi\'on por
lo que es salvajemente ramificado a menos que sea
una extensi\'on de constantes.

Ahora estudiamos el campo $K_1^{\ast}$.
Para hallar la descripci\'on expl\'icita de
$K_1^{\ast}$ procedemos como sigue. Sea
\begin{gather}\label{Eq3.1'.B}
F_0:=\prod_{\P\in{\ma P}_K^{\ast}}F_{\P}=\prod_{i=1}^t F_{\P_i}
\end{gather}
donde $K\subseteq F_{\P}\subseteq K(\Lambda_P)$ 
es el \'unico subcampo de la extensi\'on de
$K(\Lambda_P)/K$ de grado 
\begin{gather}\label{Eq3.3'.B}
c_{\P}:=\mcd (e_{\P},q^{d_{\P}}-1)=\mcd (e^{(0)}_{\P},q^{d_{\P}}-1).
\end{gather}
Por tanto $F_0$ satisface que 
$\K F_0/\K $ es no ramificada en todos los primos
finitos (Proposici\'on \ref{P3.1.B}). 

Sea $R:=K(\Lambda_{P_1\cdots P_t})$ y $R^+:=K(
\Lambda_{P_1\cdots P_t})^+$. Entonces $F_0\subseteq R$.

\begin{teorema}\label{T3.4'.B} Con las 
notaciones anteriores, tenemos
\[
K_1^{\ast}\subseteq F_0{\ma F}_{q^{u_1}}\quad \text{y}\quad
\K (F_0\cap R^+){\ma F}_{q^{t_0^{\prime}}} \subseteq \K  K_1^{\ast}
\subseteq \K F_0 {\ma F}_{q^{u_1}},
\]
para alguna $u_1\in{\ma N}$ y donde
$F_0$ est\'a dada por {\rm (\ref{Eq3.1'.B})}, $t_0$
est\'a dada por la Proposici\'on {\rm{\ref{P3.2.B}}} 
y $t_0=t_0^{\prime}
p^v$ con $\mcd(t_0^{\prime}, p)=1$.
\end{teorema}

\begin{proof} Probaremos que 
$K_1^{\ast}\subseteq F_0{\ma F}_{q^{u_1}}$ para alguna
$u_1\in{\ma N}$. Para cualquier primo
$\P\in {\ma P}_K^{\ast}$ obtenemos de la Proposici\'on
\ref{P3.1.B} que si el \'indice de ramificaci\'on de
$\P$ en $K_1^{\ast}/K$ es $b_{\P}$, entonces 
$b_{\P}|e_{\P}$ y puesto que $K_1^{\ast}/K$
es una extensi\'on abeliana finita moderadamente 
ramificada, tenemos $b_{\P}\mid q^{d_{\P}}-1$
(Proposici\'on \ref{Palestine3.1}). Por lo tanto
$b_{\P}\mid c_{\P}=\mcd(e_{\P},q^{d_{\P}}-1)=[F_{\P}:K]$. 
Sea $F_{\P}^{\prime}$ el subcampo de 
$F_{\P}$ de grado $b_{\P}$ sobre $K$.

Podemos suponer que los primos finitos ramificados
en $K_1^{\ast}/K$ son todos los $\P_1,\ldots,
\P_t$ puesto que, si alguno de los $b_{\P_i}$ son 
iguales a $1$, el argumento dado continuaci\'on
funciona a\'un en este caso.

Empezamos con $\P_1$. 
Del Lema de Abhyankar, tenemos que el \'indice
ramificaci\'on de $\P_1$ en $K_1^{\ast}F_{\P_1}^{\prime}$
sobre $K$ es $b_{\P_1}$. Sea $I_{\P_1}$ 
el grupo de inercia de $\P_1$ en
$K_1^{\ast} F^{\prime}_{\P_1}$ el cual es de orden $b_{\P_1}$.
Sea $E_1$ el campo fijo de $K_1^{\ast}F^{\prime}_{P_1}$ bajo
$I_{\P_1}$. Puesto que $\P_1$ 
es totalmente ramificado en $F^{\prime}_{\P_1}/K$
y no ramificado en $E_1/K$ tenemos que
$E_1\cap F^{\prime}_{\P_1}=K$ y
\begin{gather*}
[E_1F^{\prime}_{\P_1}:K]=[E_1:K][F^{\prime}_{\P_1}:K]=
\frac{[K_1^{\ast}F^{\prime}_{\P_1}:K]}{|I_{\P_1}|}|I_{\P_1}|=
[K_1^{\ast}F^{\prime}_{\P_1}:K].
\end{gather*}
\[
\xymatrix{K_1^{\ast}\ar@{-}[rr]
\ar@{-}[dd]&&K_1^{\ast}F^{\prime}_{\P_1}=
E_1F^{\prime}_{\P_1}\ar@{-}[dd]
\ar@{-}[dl]_{I_{\P_1}}\\ &E_1=(K_1^{\ast}F^{\prime}_{\P_1})^{I_{\P_1}}
\ar@{-}[dl]\\
K\ar@{-}[rr]^{b_{\P_1}}&&F^{\prime}_{\P_1}}
\]

Por lo tanto $K_1^{\ast}F^{\prime}_{\P_1}=E_1F^{\prime}_{\P_1}$.
M\'as a\'un, debido a que
$\P_2, \ldots, \P_t$ son no ramificados en $F^{\prime}_{\P_1}$ 
sus \'indices de ramificaci\'on son
$b_{\P_2},\ldots, b_{\P_t}$ 
en $E_1F^{\prime}_{\P_1}/F^{\prime}_{\P_1}$.
Se obtiene de esta forma que $\P_2, \ldots, \P_t$ 
tienen \'indices de ramificaci\'on
$b_{\P_2},\ldots, b_{\P_t}$ en $E_1/K$.

Tomemos ahora $E_1$ en lugar de
$K_1^{\ast}$ y $F^{\prime}_{\P_2}$ en lugar de
$F^{\prime}_{\P_1}$. Obtenemos $E_2$ tal que 
$E_1F_{\P_2}^{\prime}=E_2F_{\P_2}^{\prime}$ y $\P_3,
\ldots,\P_t$ son ahora los \'unicos primos finitos de $K$
ramificados en $E_2$ con \'indices de ramificaci\'on 
$b_{\P_3},\ldots b_{\P_t}$ respectivamente. Notemos que
\[
K_1^{\ast}F^{\prime}_{\P_1}F^{\prime}_{\P_2}=
E_1F^{\prime}_{\P_1}F^{\prime}_{\P_2}=F^{\prime}_{\P_1}E_1
F^{\prime}_{\P_2}=F^{\prime}_{\P_1}E_2F^{\prime}_{\P_2}=
E_2F^{\prime}_{\P_1}F^{\prime}_{\P_2}.
\]

En el paso general, tenemos 
$E_{i-1}F_{\P_i}^{\prime}=E_iF_{\P_i}^{\prime}$, los
\'indices de ramificaci\'on de $\P_{i+1},\ldots, \P_t$
en $E_i/K$ son $b_{\P_{i+1}}, \ldots, b_{\P_t}$
y $K_1^{\ast}F^{\prime}_{\P_1}\ldots F^{\prime}_{\P_i}=
E_i F^{\prime}_{\P_1}\ldots F^{\prime}_{\P_i}$.

Continuando con este proceso, finalmente
obtenemos un campo $E_t$ el cual satisface 
$E_{t-1}F_{\P_t}^{\prime}=E_tF_{\P_t}^{\prime}$,
ning\'un primo finito es ramificado en $E_t/K$ y
$K_1^{\ast} F_0^{\prime} = E_t F_0^{\prime}$ donde
\begin{gather}\label{Ec3.13'}
F_0^{\prime}=\prod_{i=1}^t F^{\prime}_{\P_i}.
\end{gather}

Puesto que el \'unico posible primo ramificado
en $E_t/K$ es $\p$ y es moderadamente ramificado,
de la Proposici\'on \ref{Palestine3.3} obtenemos
que $E_t/K$ es una extensi\'on de constantes,
digamos $E_t={\ma F}_{q^{u_1}}(T)=K_{u_1}$.

Puesto que los campos $\big\{F_{\P_i}^{\prime}\big\}_{i=1}^t$ son
linealmente disjuntos a pares y
$F_0^{\prime}/K$ es una extensi\'on geom\'etrica, tenemos que
\begin{gather*}
[F_0^{\prime}:K]=\prod_{i=1}^t [F_{\P_i}^{\prime}:K]=\prod_{i=1}^t
b_{\P_i}, \quad E_t\cap F_0^{\prime}=K \quad\text{y}\\
[K_1^{\ast}F_0^{\prime}:K]=[E_t:K][F_0^{\prime}:K].
\end{gather*}
En particular, ${\ma F}_{q^{u_1}}$ es el campo de constantes de
$K_1^{\ast}F_0^{\prime}$.

Por lo tanto $K_1^{\ast}\subseteq K_1^{\ast}F_0^{\prime}
= E_t F_0^{\prime} \subseteq  F_0{\ma F}_{q^{u_1}}$ y
$\K K_1^{\ast}\subseteq \K F_0{\ma F}_{q^{u_1}}$.
Finalmente, puesto que la extensi\'on
$\K (F_0\cap R^+){\ma F}_{q^{t_0}}/\K$ es
no ramificada y los primos infinitos son totalmente
descompuestos, se sigue que 
$\K (F_0\cap R^+){\ma F}_{q^{t_0^{\prime}}}
\subseteq \K K_1^{\ast}$. $\fin$
\end{proof}

\begin{observacion}\label{O3.4''.B}
En la demostraci\'on del Teorema
\ref{T3.4'.B} de hecho hemos obtenido que
$K_1^{\ast}\subseteq E_tF_0^{\prime}$ 
y que ${\ma F}_{q^{u_1}}$ es el 
campo de constantes de $K_1^{\ast}F_0^{\prime}$.
\end{observacion}

Mencionamos que es posible usar las t\'ecnicas de la 
demostraci\'on del Teorema \ref{T2.1.A} para demostrar
el Teorema \ref{T3.4'.B}. Por otro lado, n\'otese
la similitud del Teorema \ref{T3.4'.B} con el Corolario
\ref{C4.2}.

Para estudiar $K_2^{\ast}$ primero probamos
(ver Observaci\'on \ref{O12*.2.2.K}):

\begin{lema}\label{L3.5.B}
Tenemos $\g {K^{\ast}} =\g {({K^{\ast}_1})}
\g {({K^{\ast}_2})}= K^{\ast}$. M\'as a\'un
$\g {({K^{\ast}_1})}
=K^{\ast}_1$ y $\g {({K^{\ast}_2})}=K^{\ast}_2$.
\end{lema}

\begin{proof} Tenemos que $K^{\ast}=K_1^{\ast}K_2^{\ast}$ 
y ya hemos hecho notar que $\g {K^{\ast}}=K^{\ast}$.
Puesto que $\g {({K^{\ast}_1})}/K_1^{\ast}$ 
es no ramificada y que los primos infinitos se
descomponen totalmente, lo mismo sucede en
la extensi\'on $K^{\ast}
\g {({K^{\ast}_1})}/K^{\ast}$ de tal forma
$\g {({K^{\ast}_1})}\subseteq
\g {K^{\ast}}$. Similarmente 
$\g {({K^{\ast}_2})}\subseteq \g {K^{\ast}}$.
De aqu\'i $\g {({K^{\ast}_1})}\g {({K^{\ast}_2})}
\subseteq \g {K^{\ast}}$.

Ahora, puesto que $\g {({K^{\ast}_1})}
\supseteq K_1^{\ast}$ y
$\g {({K^{\ast}_2})}\supseteq K_2^{\ast}$, obtenemos
\[
\g {K^{\ast}}=K^{\ast} =K_1^{\ast} K_2^{\ast}\subseteq 
\g {({K^{\ast}_1})}\g {({K^{\ast}_2})}\subseteq \g {K^{\ast}}.
\]

Sea ahora $[K_1^{\ast}:K]=a$ y $[K_2^{\ast}:K]=p^v$ 
donde $p\nmid a$.
Si $K_1^{\ast}\subsetneqq \g {({K^{\ast}_1})}$, consideremos
$S:=\g {({K^{\ast}_1})}\cap K_2^{\ast}$. 
De la correspondencia de Galois, obtenemos que
$S\neq K$.

\[
\xymatrix{
\*{K_1}\ar@{-}[r]\ar@{-}[d]&\g{(\*{K_1})}\ar@{-}[r]\ar@{-}[d]
&\g{\*K}=\*{K_1}\*{K_2}\ar@{-}[d]\\
K=\*{K_1}\cap \*{K_2}\ar@{-}[r]&S=\g{(\*{K_1})}\cap \*{K_2}
\ar@{-}[r]&\*{K_2}
}
\]

Sea $[S:K]=p^b$ con $b\geq 1$. 
Tenemos que $S/K$ es no ramificado puesto que
de otra forma existir\'ia un primo en $K$
con \'indice ramificaci\'on $p^{c}$ con
$c\geq 1$ en $S$. Puesto que $p\nmid a$, 
se sigue que existe un primo ramificado en
$\g {({K^{\ast}_1})}/K_1^{\ast}$ con \'indice de 
ramificaci\'on $p^c$. 
Esta contradicci\'on muestra que
$S/K$ es no ramificada. Por tanto $S/K$
es una extensi\'on de constantes.
Se sigue que $\p$ tiene grado de inercia
$p^b$ en $S/K$ pero esto implicar\'ia que el 
grado de inercia de los primos infinitos en
$\g {({K^{\ast}_1})}/K_1^{\ast}$ es $p^b$
lo cual es imposible. Por tanto
$\g {({K^{\ast}_1})}=K_1^{\ast}$.
Similarmente $\g {({K^{\ast}_2})}=K^{\ast}_2$. $\fin$
\end{proof}

\begin{observacion}[Ver Observaci\'on \ref{O12*.2.2.K}]\label{O3.6.B}
En general se tiene que
si $L=L_1L_2$, entonces $\g {(L_1)}\g {(L_2)} \subseteqq
\g L$ pero no necesariamente $\g L=
\g {(L_1)}\g {(L_2)}$.
\end{observacion}

\begin{ejemplo}\label{Ej12*.2.2.M}
Sean $q>2$ y sean $P,Q,R,S\in R_T$ 
cuatro polinomios m\'onicos en $K$. Sea
$L_1:=K(\Lambda_{PQ})^+$
y $L_2:=K(\Lambda_{RS})^+$. 
Entonces, por el Teorema \ref{T12*.3.B}, se tiene que
$L_1\subseteq K(\Lambda_{PQ})$
y $K(\Lambda_{PQ})/L_1$ es totalmente
ramificada en $\S{L_1}$ de donde se sigue que
$L_1=\g {(L_1)}$. Similarmente
$L_2=\g {(L_2)}$. 

Sea $L:=L_1L_2$. Veamos que $\g L=K(\Lambda_{
PQRS})^+$. 

Primero, como consecuencia del
Teorema \ref{T12*.3.B}, $\g L \subseteq 
K(\Lambda_{PQRS})$. Puesto que
$L\subseteq K(\Lambda_{PQRS})^+$,
se sigue que $\g L\subseteq 
K(\Lambda_{PQRS})^+$.

Ahora bien, por el Corolario \ref{C6.4.30(3)}
se tiene que $K(\Lambda_{PQ})/L_1$ es no
ramificada en todos los primos finitos. De aqu\'i
obtenemos que $P$ y $Q$ son no ramificados en
$K(\Lambda_{PQRS})/L_1$. Similarmente
$R$ y $S$ son no ramificados en 
$K(\Lambda_{PQRS})/L_2$.

Por tanto los \'unicos primos ramificados en
$K(\Lambda_{PQRS})/L$ son los elementos
de $\S L$. Ahora bien
\[
L\subseteq K(\Lambda_{PQRS})^+
\subseteq K(\Lambda_{PQRS}),
\]
y $\S L$ se descompone totalmente en
$K(\Lambda_{PQRS})^+/L$.
Por tanto obtenemos que
 $K(\Lambda_{PQRS})^+
\subseteq \g L$ y se sigue que
$\g L=K(\Lambda_{PQRS})^+$.

Finalmente, 
\[
[K(\Lambda_{PQRS}):L]
=(q-1)^2, \quad [K(\Lambda_{PQRS}):
K(\Lambda_{PQRS})^+]=q-1,
\]
de donde $[\g L:L]=q-1>1$. 
Por tanto 
\[
\g L=\g {(L_1 L_2)}\neq \g {(L_1)}\g {(L_2)}= L.
\]
\end{ejemplo}

En la Subsecci\'on \ref{SS5.A.1} se estudiar\'a con m\'as
detalle el fen\'omeno de la Observaci\'on \ref{O3.6.B}.

Regresando a nuestro desarrollo, sea 
$\Gal(K_2^{\ast}/K)\cong C_{p^{n_1}}\times \cdots
\times C_{p^{n_{\nu}}}$ y si por cada $1\leq i\leq \nu$, $E_i$
es un subcampo $K\subseteq E_i\subseteq K_2^{\ast}$
tal que $\Gal(E_i/K)\cong C_{p^{n_i}}$ entonces, del
Corolario \ref{C3.3.A}, obtenemos que
$\g {(E_i)}$ es la composici\'on de $p$--extensiones
c\'iclicas de $K$ en cada una de las cuales, o bien \'unicamente
un primo es ramificado o bien es una extensi\'on de
constantes.

Por tanto $\g {({K^{\ast}_2})}=K_2^{\ast}$ es
la composici\'on de este tipo de $p$--extensiones
c\'iclicas.

Finalmente, puesto que $u_{\P_j}\geq 1$
para $s+1\leq j\leq t$ (ver (\ref{Eq0.1.B})),
obtenemos el siguiente resultado.

\begin{teorema}\label{T3.6.B}
El campo $K_2^{\ast}$ es de la forma 
$K_2^{\ast} = J_{s+1} J_{s+2}\cdots J_t J_{\infty}$ donde
$\P_j$ es el \'unico primo ramificado en $J_j/K$,
$[J_j:K]=p^{v_j}$ con $0\leq v_j\leq u_{\P_j}$ 
para $s+1\leq j\leq t$, y
$J_{\infty}$ es una $p$--extensi\'on abeliana finita la cual
es, o bien una extensi\'on de constantes o bien,
$\p$ es el \'unico primo ramificado. $\fin$
\end{teorema}

\subsection{El campo de g\'eneros en
un caso especial}\label{S3.B}

Sea $\K/K$ una extensi\'on finita y separable tal que
para todo $\P\in{\ma P}_K$, $p\nmid e_{\P}=
\mcd(e_1,\ldots,e_r)$ donde $\con_{K/\K}\P=\producto$. 
Esto es, suponemos que $t=s$ 
y tambi\'en suponemos que $p\nmid e_{\p}$.

Tenemos que $K_1^{\ast}$ est\'a dado por
(\ref{Eq3.1''.B}) y en este caso tenemos que
$K_2^{\ast}/K$ es una extensi\'on no ramificada.
Se sigue que $K_2^{\ast}/K$ es una extensi\'on
de constantes.

Con el fin de encontrar expresiones m\'as 
expl\'icitas de $K_1^{\ast}$ procedemos de
la siguiente forma. Primero consideramos 
el comportamiento de $\p$. Sea 
$\con_{K/\K}\p$ dado por (\ref{Eq3.2'.B}).

Tenemos que
$e_{\infty}(F_{\P}|K)\mid
\mcd(c_{\P},q-1)$ para $\P\in {\ma P}_K^{\ast}$. Del
Lema de Abhyankar obtenemos que si 
\[
c_{\infty}:=e_{\infty}(F_0|K),
\]
es el \'indice de ramificaci\'on de $\p$ in $F_0/K$
entonces
\begin{align*}
c_{\infty}&\mid \mcm\big[\mcd(e_{\P_1},q-1),\ldots,\mcd(e_{\P_s},q-1)\big]\\
&=\mcd\big(\mcm[e_{\P_1},\ldots, e_{\P_s}],q-1\big).
\end{align*}

Para obtener una f\'ormula para
$c_{\infty}$, procedemos como sigue.
Tenemos de (\ref{Eq3.3'.B}) que
\[
c_{\P}=[F_{\P}:K]=\mcd (e_{\P},q^{d_{\P}}-1),
\]
donde $d_{\P}=d_K(\P)$. Sea $H:=\Gal(R/F_0)$, donde $R=K(
\Lambda_{P_1\cdots P_s})$. Sea $S:=F_0R^+$, con $R^+=K(
\Lambda_{P_1\cdots P_s})^+$.
Por tanto
\begin{align*}
e_{\infty}(S|K)&=e_{\infty}(S|R^+)e_{\infty}(R^+|F_0\cap R^+)
e_{\infty}(F_0\cap R^+|K)\\
&=[S:R^+]\cdot 1\cdot 1=[S:R^+]=
[F_0:F_0\cap R^+];\\
e_{\infty}(S|K)&=e_{\infty}(S|F_0)e_{\infty}(F_0|F_0\cap R^+)
e_{\infty}(F_0\cap R^+|K)\\
&=1\cdot e_{\infty}(F_0|F_0\cap R^+)\cdot 1=
e_{\infty}(F_0|F_0\cap R^+).
\end{align*}
Se sigue que
\begin{align}
c_{\infty}&=e_{\infty}(F_0|K)=e_{\infty}(F_0|F_0\cap R^+)=e_{\infty}
(S|K)\nonumber\\
& =[F_0:F_0\cap R^+]=[S:R^+]. \label{Eq3.2.B}
\end{align}

Seleccionamos al m\'aximo campo $F$ que satisface
$F_0\cap R^+\subseteq F
\subseteq F_0$ y tal que los primos infinitos de $\K$
se descomponen totalmente en $\K F$. Notemos que 
tal campo $F$ existe puesto que si 
$F_1,F_2$ son dos campos tales que
$F_0\cap R^+\subseteq F_i\subseteq F_0$ 
y tales que los primos infinitos de
$\K$ se descomponen totalmente en $\K F_i/\K $, $i=1,2$,
entonces $F_1F_2$ satisface las mismas propiedades.

\begin{observacion}\label{O4.-1.B}
Con las notaciones del Teorema \ref{T3.4'.B},
observamos que, puesto que $\K F/\K $ 
es no ramificada y como $\p$ se descompone totalmente
en $\K F/\K $, se sigue que $F\subseteq K_1^{\ast}$  
de tal forma que $K_1^{\ast} F=K_1^{\ast}
\subseteq K_1^{\ast}F_0^{\prime}$. 

Ahora bien, se tiene que $F_0\cap R^+\subseteq \K$
y por tanto $\K(F_0\cap R^+)/\K$ es no ramificada. Por
la Proposici\'on \ref{P3.1.B} tenemos que si el \'indice
de ramificaci\'on de un primo $\P$ de $K$
en $F_0\cap R^+/K$ es $a_{\P}$ entonces $a_{\P}\mid
b_{\P}=[F'_{\P}:K]$. Si $X$ es el grupo de caracteres
de Dirichlet asociado a $F_0\cap R^+$, entonces
$|X_{\P}|=a_{\P}\mid b_{\P}$ y como el grupo de
caracteres de Dirichlet $Y$ asociado a
$F_0'=\prod_{\P} F_{\P}'$ (ver (\ref{Ec3.13'})) 
satisface $Y=\prod_{\P} Y_{\P}$, se sigue
que $X\subseteq Y$ y por tanto $F_0\cap
R^+\subseteq F_0^{\prime}\subseteq F_0$. Se obtiene que 
$F\subseteq F_0^{\prime}$ .
En general es posible que $F_0^{\prime}\neq F$, 
ver Ejemplo \ref{Ex5.1.B}.
\end{observacion}

Como siguiente paso, determinamos $F$ para
una extensi\'on abeliana $\K/K$.

\begin{proposicion}\label{P4.0.B} 
Sea $\K/K$ una extensi\'on abeliana finita
moderadamente ramificada. Con la notaci\'on del 
Teorema {\rm{\ref{T2.6.B}}} tenemos
\begin{gather*}
F\subseteq \g E \subseteq F_0,
\intertext{m\'as precisamente}
F=\g {E^{H_1}}\quad\text{y}\quad \g \K =\K F.
\end{gather*}
\end{proposicion}

\begin{proof} En este caso se tiene que
$s=t$ y que $N=P_1\cdots P_t$. Puesto que para
cualquier primo $\P$ en $K$, el \'indice ramificado
en $\K/K$ es el mismo que el \'indice de
ramificaci\'on en $E/K$ y 
$F_0=\prod_{\P\in{\ma P}_K^{\ast}} F_{\P}$,
tenemos $\g E\subseteq F_0$.

El primo infinito se descompone totalmente en
$F_0\cap R^+/K$. Por tanto los primos infinitos 
se descomponen totalmente en $E(F_0\cap R^+)/E$. 
Puesto que la extensi\'on $E(F_0\cap R^+)/E$
es no ramificada, se tiene que
$F_0\cap R^+\subseteq \g E$. 

Por el Lema de Abhyankar (Teorema \ref{T2.3.B}) vemos
que la extensi\'on $\K (F_0\cap R^+)/\K$ es no ramificada
y los primos infinitos se descomponen totalmente. 
Por tanto $F\subseteq \g E$.

Finalmente, nuevamente por el Lema de Abhyankar,
$\K \g E/\K $ es no ramificada y la inercia
de los primos infinitos corresponde a $H_1$, esto es,
$\g {E^{H_1}}$ es la m\'axima extensi\'on tal que
$F_0 \cap R^+\subseteq \g {E^{H_1}}\subseteq
F_0$ y que en $\K \g {E^{H_1}}/\K$
los primos infinitos se descomponen totalmente.
Por tanto $F=\g {E^{H_1}}$. Del Teorema \ref{T2.6.B},
se sigue que $\g \K =\K F$. $\fin$
\end{proof}

\begin{observacion}\label{O5.4.B} Sea $\K/K$ una
extensi\'on abeliana finita moderadamente ramificada.
Sean $\P_1,\ldots,\P_s$ los primos finitos ramificados.
Entonces $F_0=\prod_{i=1}^s F_{\P_i}$ 
con $K\subseteq F_{\P_i}
\subseteq K(\Lambda_{P_i})$. Tenemos 
$[F_{\P_i}:K]=c_{\P_i}=\mcd
(e_{\P_i},q^{\deg P_i}-1)$. Puesto que $\K/K$ 
es abeliana y moderadamente ramificada,
se tiene que $e_{\P_i}\mid q^{\deg P_i}-1$ 
(Proposicion \ref{Palestine3.1}). Por tanto, 
$c_{\P_i} = e_{\P_i}$. 
\end{observacion}

Ahora sea
\begin{gather*}
c^{\prime}_{\infty}:=[F:F_0\cap R^+]=e_{\infty}(F|K).
\intertext{Puesto que} 
c_{\infty}=[F_0:F_0\cap R^+] =[F_0:F][F:F_0\cap R^+]=
[F_0:F]c^{\prime}_{\infty}
\end{gather*}
tenemos $c^{\prime}_{\infty}\mid c_{\infty}$.
\[
\xymatrix{&R\ar@{-}[d]_{H\cap \ff}
\ar@/^5pc/@{-}[ddd]^{q-1}\ar@/^3pc/@{-}[dd]^{H_2}\\
F_0\ar@{-}[r]\ar@{-}[d]_{H_2/(H\cap \ff)}\ar@/^1pc/@{-}[ur]^H
\ar@/_6pc/@{-}[dd]_{c_{\infty}}&
S=F_0 R^+\ar@{-}[d]_{H_2/(H\cap \ff)}\\
F\ar@{-}[r]\ar@{-}[d]_{c^{\prime}_{\infty}}&N=FR^+\ar@{-}[d]_{
c^{\prime}_{\infty}}\\
F_0\cap R^+
\ar@{-}[r]&R^+}
\]

El grado $c^{\prime}_{\infty}$ debe satisfacer lo siguiente. Por
el Lema de Abhyankar, tenemos que si
${\eu P}$ es un primo en $\K F$ que divide
$\p$ y si ${\eu P}\cap \K ={\eu p}_{i,\infty}$, entonces
\begin{align*}
e({\eu P}|\p)&=\mcm[e_{i,\infty},c^{\prime}_{\infty}]=\frac{
e_{i,\infty}c^{\prime}_{\infty}}{\mcd(e_{i,\infty},c^{\prime}_{\infty})}\\
&=e({\eu P}|{\eu p}_{i,\infty})e({\eu p}_{i,\infty}|\p)=
e({\eu P}|{\eu p}_{i,\infty})e_{i,\infty}.
\end{align*}
Se sigue que
\begin{gather}\label{Eq3.3.B}
e({\eu P}|{\eu p_{i,\infty}})=\frac{c^{\prime}_{\infty}}
{\mcd(e_{i,\infty},c^{\prime}_{\infty})}.
\end{gather}
Por lo tanto
\begin{gather*}
e({\eu P}|{\eu p}_{i,\infty})=1\iff
\mcd (e_{i,\infty},c^{\prime}_{\infty})=c^{\prime}_{\infty} 
\iff c^{\prime}_{\infty}\mid e_{i,\infty}.
\end{gather*}
De esta forma se tiene que, 
$\K F/\K $ es no ramificada si y solamente si
$c^{\prime}_{\infty}\mid e_{\P_{\infty}}=\mcd (e_{1,\infty},\ldots,
e_{r_{\infty},\infty})$.

Por tanto $c^{\prime}_{\infty}$ debe ser maximal en
el sentido de que $c^{\prime}_{\infty}\mid c_{\infty}$,
$c^{\prime}_{\infty}\mid 
e_{\infty}$ donde $e_{\infty}=e_{\p}$,  
$c_{\infty}$ est\'a dado por (\ref{Eq3.2.B})
y los primos infinitos de $\K$ se descomponen
totalmente en $\K F$. Por lo tanto
\begin{gather}\label{Eq3.30.B}
c^{\prime}_{\infty}\mid\mcd (c_{\infty}, e_{\infty}).
\end{gather}
Esto es, $F$ es el campo
\begin{gather}\label{Eq3.31.B}
F_0\cap R^+\subseteq F\subseteq F_0\quad\text{tal que}\quad
[F:F_0\cap R^+]=c^{\prime}_{\infty}.
\end{gather}

Sea $H_2$ el subgrupo de $\ff$ de orden $\frac{q-1}
{c_{\infty}^{\prime}}$ y sea $N:=R^{H_2}$. Notemos que 
$|H\cap \ff|=[R:F_0R^+]=
\frac{q-1}{c_{\infty}}$. Por tanto $|H\cap \ff|\mid |H_2|$ y de
(\ref{Eq3.30.B}) obtenemos
\[
[S:N]=[F_0:F]=\frac{c_{\infty}}{c_{\infty}^{\prime}}.
\]

Con las notaciones anteriores, se tiene el siguiente resultado.

\begin{teorema}\label{T3.4.B} Sea $\K/K$ una extensi\'on
abeliana finita y separable tal que cada primo
$\P\in {\ma P}_K$ satisface que si
$\con_{K/\K}\P=\producto$, entonces 
$p\nmid e_{\P}= \mcd(e_1,\ldots,e_r)$.
De esta forma tenemos
\[
\K F{\ma F}_{q^{t_0}} \subseteq \g \K  \subseteq \K F_0 {\ma F}_{q^{u}},
\]
donde $F_0$ est\'a dada por {\rm (\ref{Eq3.1'.B})}, $F$ est\'a
dado por {\rm (\ref{Eq3.31.B})}, $t_0$ 
est\'a dado por {\rm{(\ref{Eq3.2''.B})}}
y $u\in{\ma N}$.

M\'as a\'un, $\K F_0{\ma F}_{q^{u}}/\K $ es no ramificada
en cada primo finito y el \'indice de ramificaci\'on del
primo infinito ${\eu p}_{i,\infty}$ es
$\frac{c_{\infty}}{\mcd(e_{i,\infty},c_{\infty})}$, 
$1\leq i\leq r_{\infty}$
donde $c_{\infty}$ est\'a dado por {\rm{(\ref{Eq3.2.B})}}.
\end{teorema}

\begin{proof} Puesto que $\K F/\K $ es no ramificada y los
primos infinitos se descomponen totalmente, obtenemos
que $F\subseteq K_1^{\ast}$. Por lo tanto, por la
Proposici\'on \ref{P3.2.B} tenemos que
 $\K F\finite {t_0}\subseteq \K  K^{\ast}=\g \K$.

Puesto que $p\nmid e_{\P}$ para toda
$\P\in{\ma P}_K$ se sigue del Teorema \ref{T3.6.B} y de la
Proposici\'on  \ref{Palestine3.3} que $K_2^{\ast}/K$
es una extensi\'on de constantes, de tal forma que
$K_2^{\ast}={\ma F}_{q^{u_2}}(T)$. 
M\'as a\'un, puesto que $K_2^* \subseteq \g \K $ 
tenemos $u_2|t_0$. 

Del Teorema \ref{T3.4'.B} tenemos que $K_1^{\ast}\subseteq F_0
{\ma F}_{q^{u_1}}$ para alguna $u_1\in{\ma N}$
y $\K K_1^{\ast}\subseteq \K F_0{\ma F}_{q^{u_1}}$.
De aqu\'i $\g \K =\K K^{\ast}=\K K_1^{\ast}K_2^{\ast}\subseteq
\K F_0{\ma F}_{q^{u}}$, donde $u=\mcm [u_1,u_2]$.

El \'indice de ramificaci\'on de $\p$ en
$F_0/K$ es $c_{\infty}$ donde $c_{\infty}$ est\'a
dado por (\ref{Eq3.2.B}).
Aplicando (\ref{Eq3.3.B}) a $F_0$ y $c_{\infty}$
obtenemos que el \'indice de ramificaci\'on de ${\eu p}_{i,\infty}$
en $\K F_0{\ma F}_{q^{u}}/\K $ es $\frac{c_{\infty}}
{\mcd(e_{i,\infty},c_{\infty})}$.
$\fin$
\end{proof}

\begin{observacion}\label{O3.5.B} Notemos que
$F\subseteq \g \K \cap F_0$. Puesto que
$\g \K \cap F_0 \subseteq \g \K $, tenemos que
los primos infinitos de $\K$ se descomponen totalmente
en $(\g \K \cap F_0 )\K $. Adem\'as,
$F_0 \cap R^+ \subseteq \g \K \cap F_0 \subseteq F_0$. 
Se sigue de la maximalidad de $F$ que
\begin{gather*}
\g \K \cap F_0=F.
\intertext{De esta forma obtenemos}
\K F\finite {t_0} \cap F_0\subseteq \g\K\cap F_0
\subseteq F\subseteq \K F \finite {t_0} \cap F_0,
\intertext{esto es, $\K F \finite {t_0}\cap F_0=F$.
Por otro lado se tiene que $\K,F,\finite {t_0}\subseteq
\g\K$ y $F\subseteq F_0$, por lo tanto
$\K F\finite {t_0}\subseteq \g\K\cap \K F_0\finite {t_0}$. Ahora}
\xymatrix{
F_0\ar@{-}[r]\ar@{-}[d]&\K F_0 \finite{t_0}
\ar@{-}[d]\\ F\ar@{-}[r]& \K F \finite {t_0}
}
\end{gather*}

En la extensi\'on $\K F_0\finite{t_0}/\K F\finite {t_0}$
los primos infinitos no tienen ninguna descomposici\'on
y en $\g\K/\K F\finite {t_0}$ los primos infinitos de
descomponen totalmente, por lo que
\[
\g \K \cap \K F_0\finite {t_0}=\K F\finite {t_0}.
\]

Notemos que si hubi\'esemos tenido en la demostraci\'on
que $(\g \K )_u\cap F_0=F$, entonces, por
la correspondencia de Galois se tendr\'ia que
\begin{gather*}
(\g \K )_u=((\g \K )_u\cap F_0) K_u=F\K{\ma F}_{q^u}.
\intertext{Por tanto}
F\K{\ma F}_{q^{t_0}} 
\subseteq \g \K \subseteq (\g \K )_u=
F\K{\ma F}_{q^u}=\K F{\ma F}_{q^{t_0}}{\ma F}_{q^u}.
\end{gather*}
Se seguir\'ia que $\g \K/\K F {\ma F}_{q^{t_0}}$ ser\'ia 
una extensi\'on de
constantes y puesto que el campo de constantes de
$\g \K $ es ${\ma F}_{q^{t_0}}$, tendr\'iamos la igualdad
\[
\g \K  =\K F {\ma F}_{q^{t_0}}.
\]

En caso de que $F=F_0$ entonces 
$\K F{\ma F}_{q^{t_0}}\subseteq
\g \K \subseteq \K F{\ma F}_{q^u}$ por
lo que $\g \K/\K F{\ma F}_{q^{t_0}}$
es una extensi\'on de constantes y entonces
$\g \K = \K F{\ma F}_{q^{t_0}}$.

Finalmente, si $u=t_0$, entonces 
$\g \K =\g \K  \cap \K F_0{\ma F}_{q^{t_0}}=
\K F{\ma F}_{q^{t_0}}$ (ver el diagrama abajo). Por tanto 
$\g \K =\K F{\ma F}_{q^{t_0}}$.

Tambi\'en, cuando $\K/K$ es una extensi\'on abeliana
finita moderadamente ramificada, tenemos
$\g \K = \K F$ (ver Proposici\'on \ref{P4.0.B}).

In resumen, es muy posible que siempre
tengamos la igualdad $\g \K =\K F{\ma F}_{q^{t_0}}$.
\begin{scriptsize}
\[
\xymatrix{F_0\ar@{-}[d]\ar@{-}[r]\ar@/_3pc/@{-}[ddd]^{c_{\infty}}
&\K F_0\ar@{-}[dd]\ar@{-}[r]&(\K F_0)_ {t_0}\ar@{-}[r]
\ar@{-}[dd]&\g \K  F_0\ar@{-}[d]\ar@{-}[r]&
(\g \K  F_0)_u=(\K F_0)_u\ar@{-}[d]\\
({\g \K})_u\cap F_0\ar@{-}[d]\ar@{-}[rrr]
|!{[r];[dd]}\hole|!{[rr];[dd]}\hole&&& ({\g \K})_v\ar@{-}[r]^{u/v}
\ar@{-}[d]^{v/t_0}&({\g \K})_{u}\\
F\ar@{-}[r]\ar@{-}[d]_{c^{\prime}_{\infty}}
&\K F\ar@{-}[r]\ar@{-}[d]&(\K F)_{t_0}\ar@{-}[r]
&\g \K \ar@{-}[ur]_{u/t_0}\ar@{-}[dll]\\
F_0\cap R^+\ar@{-}[d]&\K \ar@{-}[dl]\\
K}
\]
\end{scriptsize}
\end{observacion}

Resumimos la discusi\'on anterior en la siguiente
proposici\'on.

\begin{proposicion}\label{P3.6.B} Con las notaciones
del Teorema {\rm{\ref{T3.4.B}}}, en caso de que $\K/K$ satisfaga
al menos una de las siguientes condiciones
\las
\item $(\g \K )_u\cap F_0=F$,
\item $F=F_0$,
\item $u=t_0$,
\item $\K/K$ es una extensi\'on abeliana finita moderadamente
ramificada,
\end{list}
se tiene que $\g\K=\K F{\ma F}_{q^{t_0}}$. $\fin$
\end{proposicion}

\subsection{Sobre la extensi\'on $\ggge L/\gge L1\gge L2$}\label{SS5.A.1}

A continuaci\'on estudiamos con detalle
la contenci\'on $\gge L1\gge L2\subseteq \ggge L$ dada en la Observaci\'on
\ref{O3.6.B}. Sea, como de costumbre, $K=\F(T)$. Notemos que
el Ejemplo \ref{Ej12*.2.2.M} est\'a basado en la parte ramificada
del primo infinito $\p$ en una extensi\'on ciclot\'omica. Esto no 
es casualidad. Esto se debe, entre otras razones, 
a la siguiente observaci\'on.

\begin{observacion}\label{O.SSA.1.1}
Si $\cicl N{}^+\subseteq E\subseteq \cicl N{}$ entonces 
$\g E=E$. En efecto, $\cicl N{}^+\subseteq E\subseteq
\g E\subseteq \cicl N{}$ y $\cicl N{}/\cicl N{}^+$ es totalmente
ramificada en $\p$ de donde se sigue la igualdad.
\end{observacion}

Tambi\'en es \'util el siguiente resultado.

\begin{proposicion}\label{P.SSA.1.2}
Si $P\in R_T^+$, $\gamma\in\*\F$, $L=K(\sqrt[l]{\gamma P})$
y $E=LM\cap \cicl P{}=K\big(\sqrt[l]{(-1)^{\deg P}P}\big)$, entonces
$\g E=E$ y $\g L=L$.
\end{proposicion}

\begin{proof}
Que $\g E=E$ se sigue del hecho de que el caracter de 
Dirichlet $\chi$ asociado a $E$ satisface que $\chi=\chi_P$.
Ahora bien, $LE=L\big(\sqrt[l]{(-1)^{\deg P} \gamma}\big)$.
Si $(-1)^{\deg P} \gamma\in (\*\F)^l$ entonces $LE=L$ y el
grupo de descomposici\'on $H$ de los primos infinitos en
$LE/L$ es trivial. Por tanto $\g L=L\g E^H=LE=L$.
Si $(-1)^{\deg P} \gamma\notin (\*\F)^l$, $H\cong C_l$,
$\g E^H=E^H=K$ y $\g L=L\g E^H=LK=L$.
$\fin$
\end{proof}

Sean $L_i/K$,
$i=1,2$, dos extensiones abelianas finitas y sean $E_i=L_iM
\cap \cicl N{}$, $i=1,2$, donde $L_i\subseteq {_n\cicl N{}}_m$
y $M=L_n{\ma F}_{q^m}(T)=L_nK_m$. Sean $L=L_1L_2$ y $E=LM\cap
\cicl N{}$. Entonces $E=E_1E_2$. En general, $\gge E1\gge E2
\subseteq \g E$ y $\gge L1\gge L2\subseteq \g L$. 

\begin{proposicion}\label{P.SSA.1.3}
Si $\g L=\gge L1\gge L2$ entonces $\g E=\gge E1\gge E2$.
\end{proposicion}

\begin{proof}
De la demostraci\'on del Teorema \ref{T2.1.A} obtenemos
que $\g EM=\g LM$ y que $\gge Ei M=\gge Li M$ para $i=1,2$.
Ahora bien,
\begin{align*}
\g EM&=\g LM=\gge L1\gge L2 M=\gge L1 M \gge L2 M\\
&=\gge E1 M \gge E2 M= \gge E1 \gge E2 M.
\end{align*}

Por la correspondencia de Galois, obtenemos 
que $\g E=\gge E1 \gge E2$. 
$\fin$
\end{proof}

\begin{observacion}\label{O.SSA.1.4}
El rec\'iproco de la Proposicion \ref{P.SSA.1.3} no se cumple
en general.
\end{observacion}

\begin{ejemplo}\label{E.SSA.1.5}
Sean $P_1, P_2\in R_T^+$ dos primos distintos tales que
$\deg P_1=a$ con $1\leq a\leq l-1$ y $\deg P_2=l-a$
donde $l$ es un n\'umero primo tal que $l|q-1$. 
Sean $L_1=K(\sqrt[l]{P_1})$ y $L_2=K(\sqrt[l]{\gamma P_2})$
con $\gamma\in \*\F$ y $\gamma\notin (\*\F)^l$.

Se tiene que $P_i$ es ramificado en $L_i/K$, $i=1,2$
y en particular el campo de constantes de cada $L_i$
es $\F$. Ahora bien, el campo de constantes de $L=
L_1L_2$ es $\F$
 pues el grupo de inercia de 
$P_1$ en $L/K$ es $\Gal(L/L_2)$ de donde obtenemos que
$P_1$ es ramificado en todas las subextensiones de $L$
de grado $l$ con excepci\'on de $L_2$, pero en $L_2/K$,
$P_2$ es ramificado. 

Ahora bien, $\p$ es ramificado tanto en $L_1/K$ como en 
$L_2/K$. Se tiene $E_1=K\big(\sqrt[l]{(-1)^{\deg P_1}P_1}\big)
=K\big(\sqrt[l]{(-1)^{a}P_1}\big)$ y $E_2= 
K\big(\sqrt[l]{(-1)^{\deg P_2}P_2}\big)
=K\big(\sqrt[l]{(-1)^{l-a}P_2}\big)=
K\big(\sqrt[l]{(-1)^a P_2}\big)$. De la 
Proposici\'on \ref{P.SSA.1.2} se tiene
$\gge Ei=E_i$ y $\gge Li=L_i$, $i=1,2$.

Por otro lado, tenemos 
$E=E_1E_2=K\big(\sqrt[l]{(-1)^a P_1},\sqrt[l]{(-1)^{l-a} P_2}\big)$
y $X=\langle\chi,\tau\rangle=\langle \chi_{P_1},\tau_{P_2}\rangle$
es el grupo de caracteres asociado a $E$ y puesto que $\prod_{
P\in R_T^+}X_P=X$ se tiene $\g E=E$ de donde obtenemos
\[
\ggge E=\g E=E_1E_2=\gge E1\gge E2.
\]

Por otro lado $L=L_1L_2=K(\sqrt[l]{P_1},\sqrt[l]{\gamma P_2})$. Puesto
que $\sqrt[l]{\gamma}\notin \*\F$, se tiene que $\p$ es inerte en $K(\sqrt[l]{
\gamma P_1P_2})$ pues $\deg (P_1P_2)=a+(l-a)=l$, $\sqrt[l]{\gamma}\notin
(\*\F)^l$ y $\sqrt[l]{(-1)^{\deg P_1P_2}\gamma}=\sqrt[l]{(-1)^l \gamma}
=-\sqrt[l]{\gamma}$ (ver Proposici\'on \ref{P5.1.2.Ram}).
Se sigue que el grado del primo infinito en $L$ es $d_{\p}(L)=l$ pues
$\sqrt[l]{\gamma P_1P_2}\in L$. Por tanto $\g L\supseteq L {\ma F}_{
q^l}$. De hecho $\g L=\g E^H L=E^H L$ y $EL=L\big(\sqrt[l]{(-1)^a}\big)
\subseteq L{\ma F}_{q^l}$.

En resumen, se tiene que $\g L=L{\ma F}_{q^l}=
L_1L_2{\ma F}_{q^l}=\gge L1\gge L2 
{\ma F}_{q^l}$, por lo que $[\g L:\gge L1\gge L2]=l$ y 
por tanto $\g L\neq \gge L1\gge L2$.
\end{ejemplo}

\begin{proposicion}\label{P.SSA.1.6} 
Para cualesquier $E_i\subseteq \cicl N{}$,
$i=1,2$, se tiene $[\ggge E:\gge E1\gge E2]|(q-1)$.
\end{proposicion}

\begin{proof}
Para cualquier $F\subseteq \cicl N{}$ con grupo de caracteres de
Dirichlet asociado $X$, si $Y=\prod_{P\in R_T^+}X_P$ y
$Y^+=\{\chi\in Y\mid \chi(a)=1 \text{\ para toda\ } a\in\*\F\}$,
el campo $J$ asociado a $Y$ es la m\'axima 
extensi\'on abeliana, dentro
de $\cicl N{}$, no ramificada en los primos finitos y $J^+$,
el campo asociado a $Y^+$, satisface que $\g F=FJ^+$
(ver Teorema \ref{T12*.3.B}).

Con esta notaci\'on, sean $\g E=EJ^+$, $\gge Ei=E_iJ_i^+$,
$i=1,2$. Entonces
\[
\gge E1\gge E2=E_1J_1^+E_2J_2^+=E_1E_2J_1^+J_2^+=
EJ_1^+J_2^+.
\]

Notemos que $J=J_1J_2$ (ver Proposici\'on \ref{P3.4.A}).
De esta forma tenemos que $[\g E:\gge E1\gge E2]=
[EJ^+:EJ_1^+J_2^+]$. Adem\'as se tiene que
$[EJ^+:EJ_1^+J_2^+]|[J^+:J_1^+J_2^+]$. Se sigue
que para probar la proposici\'on,
basta probar que $[J^+:J_1^+J_2^+]|(q-1)$ donde
$J_1$ y $J_2$ son campos arbitrarios tales que
$J_i\subseteq \cicl N{}$ y $J=J_1J_2$.

Se tiene que $\cicl N{}^+=\cicl N{}^{\*\F}$. Sean
$S=\Gal(\cicl N{}/J)$, $S_i=\Gal(\cicl N{}/J_i)$, $i=1,2$. Entonces
$J=\cicl N{}^S$ y $J_i=\cicl N{}^{S_i}$, $i=1,2$.
\[
\xymatrix{&\cicl N{}\ar@{-}[dl]_S\ar@{-}[dr]^{\*\F}\\
J\ar@{-}[rd]&& \cicl N{}^+\ar@{-}[dl]\\ &J^+}
\]
Entonces $J^+=J\cap \cicl N{}^+=\cicl N{}^S\cap 
\cicl N{}^{\*\F}=\cicl N{}^{S\*\F}$ y $J_i^+=\cicl N{}^{
S_i\*\F}$, $i=1,2$.

Por otro lado $J=\cicl N{}^S=J_1J_2=\cicl N{}^{S_i}
\cicl N{}^{S_2}=\cicl N{}^{S_1\cap S_2}$. Por tanto
$S=S_1\cap S_2$.

As\'i, $J_1^+J_2^+=\cicl N{}^{S_1\*\F}\cicl N{}^{S_2\*\F}=
\cicl N{}^{S_1\*\F\cap S_2\*\F}$. Por tanto 
\[
\xymatrix{
&\cicl N{}\ar@{-}[dl]\ar@/_1pc/@{-}[dl]_{S\*\F}
\ar@/^2pc/@{-}[ddl]^{S_1\*\F\cap S_2\*\F}\\
J^+\ar@{-}[d]\\ J_1^+J_2^+}
\quad
\Gal(J^+/J_1^+J_2^+)\cong\frac{S_1\*\F\cap S_2\*\F}{S\*\F}
=\frac{S_1\*\F\cap S_2\*\F}{(S_1\cap S_2)\*\F}.
\]
Se sigue que $[J^+:J_1^+J_2^+]=\Big|\frac{S_1
\*\F\cap S_2\*\F}{(S_1\cap S_2)\*\F}\Big|$.

Para cualesquiera dos grupos finitos $H_1,H_2$, tenemos
$|H_1H_2|=\frac{|H_1||H_2|}{|H_1\cap H_2|}$. Por tanto
\begin{align*}
|S_1\*\F\cap S_2\*\F|&=\frac{|S_1\*\F||S_2\*\F|}{|S_1S_2\*\F|}=
\frac{\frac{|S_1||\*\F||S_2||\*\F|}{|S_1\cap \*\F||S_2\cap\*\F|}}{
\frac{|S_1S_2||\*\F|}{|S_1 S_2\cap\*\F|}}=\frac{\frac{|S_1||S_2|
|\*\F|^2}{|S_1\cap\*\F||S_2\cap\*\F|}}
{\frac{|S_1||S_2|}{|S_1\cap S_2|} \frac{|\*\F|}{
|S_1S_2\cap \*\F|}}\\
&= \frac{|S_1\cap S_2||S_1S_2\cap\*\F||\*\F|}
{|S_1\cap\*\F||S_2\cap\*\F|}.
\end{align*}

Por otro lado $|(S_1\cap S_2)\*\F|=|S\*\F|=
\frac{|S||\*\F|}{|S\cap \*\F|}$, por tanto
\begin{align*}
[S_1\*\F\cap S_2\*\F:(S_1\cap S_2)\*\F]&=\frac
{|S_1\cap S_2||S_1S_2\cap\*\F||\*\F|}{|S_1\cap
\*\F||S_2\cap \*\F|}\cdot \frac{|S\cap \*\F|}{|S||\*\F|}\\
&= \frac{|S_1S_2\cap \*\F||S\cap\*\F|}{|S_1\cap\*\F|
|S_2\cap \*\F|}.
\end{align*}

Ahora bien, $S\cap\*\F\subseteq S_1\cap\*\F$. Sea
$[S_2\cap \*\F:S\cap \*\F]=\alpha\in{\ma N}$. Por tanto
$[S_1\*\F\cap S_2\*\F:(S_1\cap S_2)\*\F]=\frac{1}{\alpha}
\frac{|S_1S_2\cap\*\F|}{|S_1\cap \*\F|}=\frac{1}{\alpha}
\frac{\frac{|S_1S_2||\*\F|}{|S_1S_2\*\F|}}{\frac{|S_1||\*\F|}
{|S_1\*\F|}}=\frac{1}{\alpha}\frac{|S_1S_2||S_1\*\F|}{
|S_1S_2\*\F||S_1|}$.

Se tiene $S_1S_2\subseteq S_1S_2\*\F$ y sea $[S_1S_2\*\F:
S_1S_2]=\beta\in{\ma N}$. Se sigue que 
\begin{align*}
[J^+:J_1^+J_2^+]&=
[S_1\*\F\cap S_2\*\F:(S_1\cap S_2)\*\F]=
\frac{1}{\alpha\beta} \frac{|S_1\*\F|}{|S_1|}\\
&=\frac{1}{\alpha\beta}
\frac{|S_1||\*\F|}{|S_1\cap \*\F||S_1|}
=\frac{1}{\alpha\beta}\cdot \frac{|\*\F|}{|S_1\cap \*\F|}=
\frac{q-1}{\alpha\beta |S_1\cap \*\F|}.
\end{align*}

Puesto que $\alpha\beta|S_1\cap\*\F|\in{\ma N}$, se sigue
que $[J^+:J_1^+J_2^+]|(q-1)$ y que $[\g E:\gge E1\gge E2]|(q-1)$.
$\fin$
\end{proof}

Para el caso general, primero probamos el siguiente lema.

\begin{lema}\label{L.SSA.1.6'}
Sean $E_i\subseteq \cicl {N_i}{}$ para $i=1,2$ y $E=E_1E_2$.
Entonces $[\g E^+:\gge E1^+\gge E2^+]=\varepsilon | q-1$.
\end{lema}

\begin{proof} Sea $F\subseteq \cicl N{}$ cualquiera y sea $X$
el grupo de caracteres de Dirichlet asociado a $F$. Entonces, si
$Y=\prod_{P\in R_T^+}X_P$, el campo asociado a $Y$ lo denotamos
por $\g {\tilde F}$. Se tiene que $\g F\subseteq \g {\tilde F}$ y de
hecho $\g F$ se obtiene de $\g {\tilde F}$ al remover la ramificaci\'on
de $\p$ que no aparece en $F$. M\'as precisamente, si $e_{\infty}(
\g {\tilde F}|F)=h$ denota el \'indice de ramificaci\'on de los primos
infinitos en $\g {\tilde F}/F$ y si $H$ es el subgrupo de orden $h$
de $I_{\infty}(\g {\tilde F}|K)$, donde
\'este \'ultimo denota al grupo de inercia, el cual
es naturalmente isomorfo a un subgrupo del grupo c\'iclico $\*\F
\cong I_{\infty}(\cicl N{}|\cicl N{}^+)$, entonces $\g F=\g{\tilde F}^H$.
Puesto que $\g {\tilde F}^+\subseteq \g F\subseteq \g {\tilde F}$,
se sigue que $\g F^+=\g {\tilde F}^+$. M\'as a\'un, tenemos que
$\g F=F\g {\tilde F}^+$.

Sean ahora $F_1, F_2$ subcampos de alg\'un campo ciclot\'omico
y sea $F=F_1F_2$. Entonces $\g {\tilde F}=\gge {\tilde F}1 \gge 
{\tilde F}2$ (Proposici\'on \ref{P3.4.A}).

Sean $L=\g {\tilde F}$ y $L_i=\gge {\tilde F}i$, $i=1,2$. Entonces
$L=L_1L_2$  y de la demostraci\'on de la Proposici\'on \ref{P.SSA.1.6},
obtenemos que $[L^+:L_1^+L_2^+]|q-1$. Ahora bien, $L^+=
\g {\tilde F}^+=\g F^+$ y $L_i^+ =\gge {\tilde F}i^+=\gge Fi^+$ para
$i=1,2$. Por tanto para cualesquiera $F_1,F_2$ subcampos de un
campo de funciones ciclot\'omico, se tiene $[\g F^+:\gge F1^+\gge F2^+]|
q-1$.
$\fin$
\end{proof}

\begin{proposicion}\label{P.SSA.1.7}
Para $L_i/K$ dos extensiones abelianas finitas, $i=1,2$, se tiene
$[\g L:\gge L1\gge L2]|(q-1)^2$ donde $L=L_1L_2$.
\end{proposicion}

\begin{proof}
En general, sea $L/K$ una extensi\'on abeliana finita con
$L\subseteq {_n\cicl N{}}_m$ y sea $F=LM\cap \cicl N{}$.
Sean ${\mc H}=D_{\infty}(FL/L)$, el grupo de descomposici\'on
de los primos infinitos en $FL/L$ y sea ${\mc H}_1:={\mc H}|_F$.

\[
\xymatrix{
F\ar@{-}[r]\ar@{-}[d]&LF\ar@{-}[d]\\
F^{{\mc H}_1}\ar@{-}[d]\ar@{-}[r]
&LF^{{\mc H}_1}=(LF)^{\mc H}\ar@{-}[d]\\
L\cap F\ar@{-}[r]&L
}
\]

Se tiene que $F/F^{\mc H}$ es totalmente ramificada en los
primos infinitos. Esto es, $F^{\mc H}$ es la m\'axima extensi\'on
de $L\cap F$ tal que los primos infinitos se descomponen
totalmente. Esto mismo es lo que sucede en 
$LF^{{\mc H}_1}=(LF)^{\mc H}/L$
y en $L\g F^{{\mc H}_1}/L$ pues ${\mc H}=D_{\infty}(L\g F/L)=
D_{\infty}(LF/L)$.

Regresando a nuestro caso, $L=L_1L_2$, $E=LM\cap \cicl N{}$, 
$E_i=L_iM\cap \cicl N{}$, $i=1,2$ y $E=E_1E_2$. Como $\g LM=
\g EM$, se tiene $\gge Ei=\gge Li\cap \cicl N{}$, $i=1,2$ y
$\g E=\g L M\cap \cicl N{}$. Sean $H=D_{\infty}(LE/L)$ y $H_i
=D_{\infty}(L_iE_i/L_i)$, $i=1,2$.

Se tiene
\begin{gather*}
\xymatrix{
\g E\ar@{-}[d]\ar@{-}[r]& L\g E\ar@{-}[d]\\
\g E^H\ar@{-}[r]\ar@{-}[d]&L\g E^H\ar@{-}[d]^{
\substack{\p\text{\ totalmente}\\ \text{descompuesto}}}\\
L\cap \g E\ar@{-}[r]&L}
\qquad
\xymatrix{
E_i\ar@{-}[d]\ar@{-}[r]& L_i E_i\ar@{-}[d]\\
E_i^{H_i}\ar@{-}[r]\ar@{-}[d]&L_i E_i^{H_i}\ar@{-}[d]^{
\substack{\p\text{\ totalmente}\\ \text{descompuesto}}}& i=1, 2.\\
L_i\cap E_i\ar@{-}[r]&L_i}
\intertext{Se sigue que}
\xymatrix{
L_1 \gge E1^{H_1}\ar@{-}[d]_{\substack{\p\text{\ totalmente}\\ 
\text{descompuesto}}}\ar@{-}[r]& L_1L_2\gge E1^{H_1}=
L\gge E1^{H_1}\ar@{-}[d]^{\substack{\p\text{\ totalmente}\\ 
\text{descompuesto}}}\ar@{-}[r]&L\gge E1\ar@{-}[r]&L\g E
\ar@/^1pc/@{-}[dll]\\ L_1\ar@{-}[r]&L=L_1L_2
}
\end{gather*}
Esto es, $L\subseteq L\gge E1^{H_1}\subseteq L\g E$ y
$L\gge E1^{H_1}/L$ es totalmente descompuesto en los
primos infinitos.

Se sigue que $L\gge E1^{H_1}\cap \g E\subseteq \g E^H$. Por 
otro lado $\gge E1^{H_1}\subseteq L\gge E1^{H_1}$ y $
\gge E1^{H_1}\subseteq \gge E1\subseteq \g E$. Por tanto
$\gge E1^{H_1}\subseteq L \gge E1^{H_1}\cap \g E\subseteq
\g E^H$, esto es, $\gge E1^{H_1}\subseteq \g E^H$.
\[
\xymatrix{
\g E\ar@{-}[r]\ar@{-}[dd]&L\g E\ar@{-}[d]\\
&L \gge E1\ar@{-}[d]\\
L\gge E1^{H_1}\cap E\ar@{-}[r]\ar@{-}[d]&L\gge E1^{H_1}\ar@{-}[d]\\
L\cap \g E\ar@{-}[r]& L}
\]

Similarmente, $\gge E2^{H_2}\subseteq  \g E^H$ de donde se
sigue que $\gge E1^{H_1}\gge E2^{H_2}\subseteq \g E^H$.

Ahora bien, $\g E/\g E^H$ y $\gge Ei/\gge Ei^{H_i}$, $i=1,2$
son totalmente ramificados en los primos infinitos por lo que
$\g E^+=\g E\cap \cicl N{}^+\subseteq \g E^H\subseteq \g E$
y $\gge Ei^+=\gge Ei\cap \cicl N{}^+\subseteq \gge Ei^{H_i}
\subseteq \gge Ei$, $i=1,2$. Entonces
\[
\xymatrix{
\cicl N{}^+\ar@{-}[r]\ar@{-}[d]&\g E^H\cicl N{}^+\ar@{-}[r]
\ar@{-}[d]&\g E \cicl N{}^+\ar@{-}[r]\ar@{-}[d]&\cicl N{}\\
\g E^+=\g E\cap \cicl N{}^+\ar@{-}[r]&\g E^H\ar@{-}[r]&\g E
}
\]
Se tiene que $\delta=[\g E^H:\g E^+]=[\g E^H\cicl N{}^+:\cicl N{}^+] |
[\cicl N{}:\cicl N{}^+]=q-1$ y $[\g E^+:\gge E1^+\gge E2^+]=
\varepsilon|(q-1)$ (Lema \ref{L.SSA.1.6'}).

Se sigue que
\begin{align*}
[\g L:\gge L1\gge L2]&=[L\g E^H:L_1\gge E1^{H_1}L_2\gge E2^{H_2}]\\
&=[L\g E^H:L\gge E1^{H_1}\gge E2^{H_2}] \big| [\g E^H:\gge E1^{H_1}
\gge E2^{H_2}]
\intertext{y}
[\g E^H:\gge E1^{H_1}\gge E2^{H_2}]&=
\frac{[\g E^H:\g E^+][\g E^+:\gge E1^+\gge E2^+]}
{[\gge E1^{H_1}\gge E2^{H_2}:\gge E1^+\gge E2^+]}\\
&=\frac{\delta\varepsilon}
{[\gge E1^{H_1}\gge E2^{H_2}:\gge E1^+\gge E2^+]}\big|(q-1)^2
\end{align*}
por lo que $[\g L:\gge L1\gge L2]\big|(q-1)^2$.
$\fin$
\end{proof}

\section{Aplicaciones y ejemplos}\label{S5.B}

\begin{ejemplo}\label{Ex5.1.B}  Consideremos 
$q=3$ y $P=T^3+2T+1$. Tenemos que 
$P$ es irreducible en ${\ma F}_3(T)$. Sea $\K =K(\sqrt{P})$. 
En nuestra construcci\'on, si $\P$ 
es el primo correspondiente a $P$, tenemos 
$F_0=F_{\P}=K(\sqrt{(-1)^{\deg P}P})=K(\sqrt{-P})$. Ahora
$\p$ no ramificado en $\K $ y en $F_0=F_{\P}$. Por tanto
$t_0=1$, esto es, el campo de constantes de $\g\K$ es
${\ma F}_3$. Puesto que $[R^+:K]=13$ ya que
$[\cicl P:K]=q^{\gr P-1}=3^3-1=27-1=26$,
y $[F_0:K]=[F_{\P}:K]=2$, tenemos que $F_0\cap R^+
=K$. Ahora bien $\K F_0=\K (\sqrt{-1})$. 
Puesto que $\sqrt{-1}\notin
{\ma F}_3$ (Proposici\'on \ref{P4.5.2})
se sigue que $\K (\sqrt{-1})=\K{\ma F}_9$ 
y los primos infinitos son inertes en
$\K F_0/\K $. Por tanto $F=K$ y $\g \K =\K$.
Aqu\'i tenemos $F_0^{\prime}=F_0=K(\sqrt{-P})\neq K=F$.
\end{ejemplo}

\subsection{Extensiones c\'iclicas de grado primo
que no divide a $q(q-1)$}\label{S5.1.B}

Sea $l$ un primo que no divide a $q(q-1)$ y sea
$\K/K$ una extensi\'on c\'iclica de grado $l$. Sean
$\P_1,\ldots, \P_t$ los primos de $K$ ramificados en $\K$.
Por la Proposici\'on \ref{Palestine3.1}, 
$l\mid (q^{\deg P_i}-1)$ para $1\leq i\leq t$. 
Puesto que $l\nmid q-1$, en particular se tiene que
$\p$ es no ramificado. En este caso tenemos
$K\subseteq F_{\P i}\subseteq K(\Lambda_{{P_i}})$
para $1\leq i\leq t$ donde $F_{\P i}$ es
el \'unico subcampo de
$K(\Lambda_{{P_i}})$ de grado $c_{\P_i}=\mcd(
e_{\P_i},q^{d_{\P_i}}-1)=l$. Entonces
\[
F_0=\prod_{i=1}^t F_{\P i}\subseteq 
K(\Lambda_{P_1\cdots P_t})^+.
\]
Por lo tanto tenemos $c_{\infty}=e_{\p}=1$, $F_0\cap R^+=F_0$.
Se sigue que $F=F_0$ y que $c_{\infty}^{\prime}=1$. 
M\'as a\'un, si $t_0$ es el grado de los primos infinitos
sobre $\p$ en $\K$, entonces
$t_0=1$ o $l$. De hecho $t_0=1$ si y solamente
si $\p$ se descompone en $\K/K$. Esto equivale a
$\K \subseteq K(\Lambda_{P_1\cdots P_t})^+$. 
Tenemos $t_0=l$ si y solamente si $\p$ es inerte
en $\K/K$ si y solamente si
$\K \nsubseteq K(\Lambda_{P_1\cdots P_t})^+$.

De la Proposici\'on \ref{P3.6.B}
tenemos $\g \K =\K F{\ma F}_{q^{t_0}}$, puesto que en este
caso tenemos $F=F_0$ y $u=t_0$. 

Primero consideramos
$\K \subseteq K(\Lambda_{P_1\cdots P_t})^+$. 
Entonces $\g \K =
\K F{\ma F}_{q^{t_0}}=\K F=F$ y $[\g \K  {}:\K ]=l^{t-1}$.

Ahora consideramos 
$\K \nsubseteq K(\Lambda_{P_1\cdots P_t})^+$.
Entonces $\K \nsubseteq F$ y en particular 
$K\subseteq \K \cap F
\subsetneqq \K $ de tal forma
$\K \cap F=K$ y $[\K F:\K ]=[F:K]=l^t$. 
Probaremos que ${\ma F}_{q^l}\subseteq \K F$. 
Primeramente, tenemos que
\[
[\K F:K]=[{\ma F}_{q^l} F:K]=l^{t+1}.
\]
Ahora, si
$K_l:={\ma F}_{q^l}(T)$, entonces $K_l
\cap \K =K$. Ahora $\p$ es inerte en $\K/K$ 
y en $K_l/K$. El grupo de descomposici\'on
${\mathcal D}$ de $\p$ en $K_l=\K K_l$ es
es un grupo c\'iclico de orden 
$l$. Consideremos $L:=(K_l)^{\mathcal D}$. 
El primo $\p$ se descompone totalmente en
$L/K$ y $\P_1,\ldots,\P_t$ son los primos ramificados en
$L/K$. Se sigue que $L\subseteq F$. 
Puesto que $L\neq \K $ obtenemos que
$\K L=K_l$ y $\K L=K_l=\K{\ma F}_{q^l}\subseteq 
\K F$. Entonces ${\ma F}_{q^l}\subseteq \K F$.
Por lo tanto 
\[
\g \K =\K F{\ma F}_{q^l}=\K F\quad \text{y}\quad
[\g \K  {}:\K ]= [\K F:\K ]=l^t.
\]

\subsection{Extensiones radicales}\label{S5.2.B}

Sea $\K =K(\sqrt[n]{\gamma D})$, donde 
$D\in R_T$ es un polinomio m\'onico y $\gamma \in
\ff$. Sea $D=P_1^{\alpha_1}\cdots P_s^{\alpha_s}$ 
la descomposici\'on de $D$ como producto de
polinomios irreducibles. Supondremos que $D$
est\'a libre de $n$--potencias, es decir, 
$0<\alpha_i<n$ para $1\leq i\leq s$ 
y tambi\'en supondremos que $p\nmid n$.

Los primos finitos ramificados son
$\P_1,\ldots,\P_s$ y ellos son moderadamente
ramificados. El \'indice de ramificaci\'on de cada 
$\P_i$ en $\K/K$ es igual a $e_{\P_i}=n/d_i$.
Similarmente, obtenemos $e_{\infty}=e_{\infty}(\K |K) = 
n/d$, donde $d=\mcd(\deg D,n)$ (ver
Teorema \ref{TRam1}).

De esta forma obtenemos 
$F_0=\prod_{i=1}^s F_{\P_i}$ con 
$c_{\P_i}=\mcd(e_{\P_i},q^{\deg P_i}-1)=
\mcd\big(\frac{n}{d_i},q^{\deg P_i}-1\big)$. Entonces 
\begin{gather*}
e_{\infty}(F_{\P_i}|K)\mid \mcd(c_{\P_i},q-1)
=\mcd(e_{\P_i},q-1)
=\mcd\big(\frac{n}{d_i},q-1\big).\\
\intertext{Por lo tanto}
\begin{align*}
c_{\infty}&\mid \mcd(\mcm[e_{\P_i},
\ldots,e_{\P_s}],q-1)\\
&=\mcd(\mcm\big[
\frac{n}{d_1},
\ldots,\frac{n}{d_s}\big],q-1)=\mcd(\frac{n}{d_0},q-1),
\end{align*}
\intertext{donde $d_0=\mcd[d_1,\ldots,d_s]$. Tambi\'en tenemos}
c^{\prime}_{\infty}\mid \mcd(c_{\infty},e_{\infty})|\mcd\big(\frac{n}{d_0},
\frac{n}{d},q-1\big).
\end{gather*}

Del Teorema \ref{T3.4.B} obtenemos 
\[
\K F{\ma F}_{q^{t_0}} =
K(\sqrt[n]{\gamma D})F{\ma F}_{q^{t_0}}\subseteq \g{\K} =
\g {K(\sqrt[n]{\gamma D})}\subseteq \K F_0{\ma F}_{q^u},
\]
donde $t_0, u\in{\ma N}$.

Para hallar $t_0$, repetimos lo hecho en el
Teorema \ref{TRam1}. Consideramos los subcampos
$E=K(\sqrt[d]{\gamma
D})\subseteq K(\sqrt[n]{\gamma D})$. Puesto que 
$\p$ es no ramificada en
$E/K$ y completamente ramificado en
$\K/E$, tenemos que el grado de inercia
de $\p$ en $\K/K$ es igual al grado de inercia
de $\p$ en $E/K$. Notemos que 
\begin{align*}
D(T)&=T^l+a_{l-1}T^{l-1}+\cdots+a_1
T+a_0\\
&=T^l\big(1+a_{l-1}(\frac{1}{T})+\cdots +a_1(\frac{1}{T})^{l-1}
+a_0 (\frac{1}{T})^l\big)=T^l D_1(\frac{1}{T})
\end{align*}
con $D_1(0)=1$ y $d|l$. Por lo tanto $E=K(\sqrt[d]{\gamma D_1(1/T)})$
con $D_1(1/T)\in {\ma F}_q[1/T]$ y $D_1(1/T)
\equiv  1 \bmod  (1/T) $. De esta forma tendremos que
$X^d-\gamma D_1(1/T) \bmod\ \p$
se reduce a $\bar{X}^d-\gamma\in {\ma F}_q [\bar{X}]$.

Sea $\mu\in\bar{\ma F}_q$ 
una ra\'iz $d$--\'esima fija de $\gamma$. Si
$\zeta_d$ denota una ra\'iz $d$--\'esima primitiva
de la unidad, tenemos que la factorizaci\'on de
$\bar{X}^d-\gamma$ en ${\ma F}_q[\bar{X}]$ 
es de la forma
\[
\bar{X}^d-\gamma=\prod_{j=1}^r\Irr(\zeta_d^{i_j}\mu,\bar{X},{\ma F}_q)
\]
para algunos $0\leq i_1<i_2<\cdots <i_r\leq d-1$. 
Del Lema de Hensel, obtenemos que
$X^d-\gamma D_1\big(\frac{1}{T}\big)=\prod_{j=1}^r
F_j(X)$ con $F_j(X)\in K_{\infty}[X]$ 
polinomios irreducibles distintos. En particular
$\con_{K/\K}\p={\eu p}_{\infty, 1}\cdots {\eu p}_{\infty,r}$
con $\deg_\K  {\eu p}_{\infty, j}=\deg F_j(X)$, 
$1\leq j\leq r$. Por lo tanto
\begin{align*}
t_0&=\mcd_{1\leq j\leq r}\{\deg F_j(X)\}=\mcd_{1\leq j\leq r}
\{[{\ma F}_q(\zeta_d^{i_j}\mu):{\ma F}_q]\}\\
&=\mcd_{0\leq i\leq d-1}
\{[{\ma F}_q(\zeta_d^i \mu):{\ma F}_q]\}.
\end{align*}
En resumen si escribimos $\sqrt[d]{
\gamma}=\mu$,
\begin{gather}\label{Eq5.1*.B}
t_0=\mcd_{0\leq j\leq r}\big\{\big[{\ma F}_q
(\zeta_d^{i_j}\sqrt[d]{\gamma}):{\ma F}_q\big]\big\}
=\mcd_{0\leq i\leq d-1}\big\{\big[{\ma F}_q(\zeta_d^{i}
\sqrt[d]{\gamma}):{\ma F}_q\big]\big\}.
\end{gather}

Se sigue de (\ref{Eq5.1.B}) que si $\mcd(\alpha_i,n)=1$ 
para alguna $i$, entonces $\K/K$ es una
extensi\'on geom\'etrica.

\begin{ejemplo}\label{Ex5.3.B} Consideremos 
$q=3$, $P_1 = T$, $P_2=T^2 - T -1$
y $D=P_1^2P_2$. Tenemos que
$P_1$ y $P_2$ son irreducibles en 
${\ma F}_3(T)$. Sea $\K =K(\sqrt[10]{-D})$. 
Se tiene $d_1=\mcd(2,10)=2$, $e_{\P_1}=\frac{10}
{2}=5$, $e_{\P_2}=\frac{10}{1}=10$, $c_{\P_1}=
\mcd(5,2)=1$ y $c_{\P_2}=\mcd(10,2)=2$.
En nuestra
construcci\'on, si $\P_i$ es el primo correspondiente a
$P_i$, tenemos $F_{\P_1}=K$ y $F_{\P_2}=K(\sqrt{P_2})$. 
Por lo tanto $F_0=K(\sqrt{P_2})$.
Por un lado, $\p$ se descompone en $K(\sqrt{P_2})/K$, 
por lo que $K(\sqrt{P_2}) \subseteq
K(\Lambda_{P1P2})^+$. Por lo tanto 
$F =F_0=K(\sqrt{P_2})$. 

Puesto que $d=\mcd(10,4)=2$, obtenemos
$t_0 = 2$. Puesto que $u_2$ es una potencia de
$3$ y $u_2$ divide a $t_0$, tenemos
$u_2 = 1$. De la demostraci\'on del Teorema 
\ref{T3.4'.B}, obtenemos que en este caso
$E_1 = K_1^*$ y ${\ma F}_{q^{u_1}} =
E_2 \subseteq K_1^*$. Por lo tanto 
$u_1$ divide a $t_0$. Obtenemos que
$u = u_1 \in \{1,2 \}$. Se sigue del
Teorema \ref{T3.4.B} que $\g \K  = \K K(\sqrt{P_2}) 
{\ma F}_9$.
\end{ejemplo}

\subsection{Extensiones radicales de grado una
potencia de primo que divide a $q-1$}\label{S5.3.B}

Como un caso particular a la Subsecci\'on 
\ref{S5.2.B}, consideremos $l$ un n\'umero primo tal que
$l^n\mid q-1$. Sea $D\in R_T$ 
un polinomio m\'onico que es libre de $l^n$--potencias.
Sea $D=P_1^{\alpha_1}\cdots P_s^{\alpha_s}$ 
con $P_1,\ldots,
P_s\in R_T^+$ y $v_l(\alpha_i)=a_i<n$. 
Sea $\gamma \in \ff$ y
$\K =K(\sqrt[l^n]{\gamma D})$. Entonces por 
el Teorema \ref{TRam1}, se tiene
$e_{\P_i}=l^{n-a_i}$, $1\leq i\leq s$.
Puesto que $\K/K$ es una extensi\'on c\'iclica de
grado $l^n$, $\K/K$ es una extensi\'on geom\'etrica si
y solamente si $a_i=0$ para alguna $1\leq i\leq s$.

Ahora, tenemos $F_{\P_i}\subseteq K(
\Lambda_{P_i})$ y $c_{\P_i}=
\mcd(e_{\P_i},q^{\deg P_i}-1)= e_{\P_i}=l^{n-a_i}$. 
Por lo tanto
$F_{\P_i}=K(\sqrt[l^{n-a_i}]{(-1)^{\deg P_i}P_i})$ y 
$F_0=\prod_{i=1}^s F_{\P_i}$.

Se tiene $e_{\p}=e_{\infty}=l^{n-d}$, donde 
$d=\min\{n,d^{\prime}\}$
y $v_l(\deg D)=d^{\prime}$. M\'as a\'un, el grado
de inercia de $\p$ es 
$f_{\infty}=l^m$, donde ${\ma F}_{q^{l^m}}
=\F(\sqrt[l^d]
{(-1)^{\deg D}\gamma})$ (ver Proposici\'on \ref{P3.7.V.Ram}).
Por lo tanto $t_0=l^m$ y el campo de
constantes de $\g \K $ es ${\ma F}_{q^{l^m}}$.

Ahora, con respecto a $\p$ tenemos
$e_{\infty}(F_{\P_i}|K)=l^{n-a_i-d_i}$, donde
$d_i=\min\{n-a_i,d_i^{\prime}\}$, 
$v_l(\deg P_i)=d^{\prime}_i$.
Del Lema de Abhyankar obtenemos
\begin{align*}
e_{\infty}(F_0|K)&=\mcm[e_{\infty}(F_{\P_i}|K)\mid 1\leq i\leq s]\\
&=\mcm[l^{n-a_i-d_i}\mid 1\leq i\leq s]=l^{n-\delta},
\end{align*}
donde $\delta=\min\limits_{1\leq i\leq s}\{a_i+d_i\}=
\min\limits_{1\leq i\leq s}\{a_i+\min\{n-a_i,d_i^{\prime}\}\}=
\min\limits_{1\leq i\leq s}\{n,v_l(\deg P_i^{\alpha_i})\}$. 
Esto es
\[
c_{\infty}=[F_0:F_0\cap R^+]=e_{\infty}(F_0|K)=l^{n-\delta}.
\]
Notemos que $d\geq \delta$. 
Entonces $c^{\prime}_{\infty}\mid\mcd(c_{\infty},
e_{\infty})=\mcd(l^{n-\delta},l^{n-d})=l^{n-d}$.
Se tiene que $F$ es el subcampo
$F_0\cap R^+\subseteq F
\subseteq F_0$ tal que 
$[F:F_0\cap R^+]=c_{\infty}^{\prime}\mid
l^{n-d}$.
\[
\xymatrix{
F_0\ar@{-}[d]^{\frac{c_{\infty}}{c_{\infty}^{\prime}}}
\ar@/_2pc/@{-}[dd]_{l^{n-\delta}=c_{\infty}}\\
F\ar@{-}[d]^{c^{\prime}_{\infty}\mid l^{n-d}}\\ F_0\cap R^+}
\]

\begin{ejemplo}\label{Ex5.3'.B}
Sean $K={\ma F}_5(T)$ y 
$\K :=K\big(\sqrt[3]{T(T^2+T+1)}\big)\cdot {\ma F}_{5^2}=
K\big(T,\sqrt[3]{D(T)}\big)\cdot {\ma F}_{5^2}$, 
donde $D(T)=T(T^2+T+1)$. Sea
$K_0:=K(\sqrt[3]{D(T)})$. Notemos que 
si $\zeta_3$ denota una ra\'iz tercera primitiva de la
unidad, entonces $\zeta_3\notin {\ma F}$ y
$\K =K_0(\zeta_3)
=K_0\cdot {\ma F}_{5^2}$
es la cerradura de Galois de $K_0/K$.
Se tiene que $T$ y $T^2+T+1$ 
son irreducibles en ${\ma F}_5(T)$
puesto que $\zeta_3\notin {\ma F}_5$
y $T^2+T+1=\frac{T^3-1}{T-1}=
(T-\zeta_3)(T-\zeta_3^2)$. 
De hecho ${\ma F}_5(\zeta_3)={\ma F}_{5^2}
={\ma F}_{25}$. 

Sean ${\mathcal P}_T$ y ${\mathcal P}_{
T^2+T+1}$ los divisores primos en
$K={\ma F}_5(T)$ correspondientes a
$T$ y a $T^2+T+1$ respectivamente. 
El primo infinito $\p$ es no ramificado ni en
$K_0/K$ ni en $\K/K$ debido a que $\deg
D(T)=3$. Sea $t_0(K_0)$ y
$t_0(\K )$ dado por (\ref{Eq3.2''.B}) con
respecto a los campos
$K_0$ y $\K $ respectivamente. Puesto
que $\gamma=1$, de (\ref{Eq5.1*.B}) obtenemos
\begin{gather*}
t_0(K_0)=\mcd_{0\leq i\leq 2}\big\{\big[
{\ma F}_5(\zeta_3^i):{\ma F}_5\big]\big\}=
\mcd\{1,2,2\}=1
\intertext{y}
\begin{align*}
t_0(\K )&=\mcd_{0\leq i\leq 2}\big\{\big[
{\ma F}_{5^2}(\zeta_3^i):{\ma F}_5\big]\big\}
=\mcd_{0\leq i\leq 2}\{[{\ma F}_{5^2}:{\ma F}_5]\}\\
&= \mcd\{2,2,2\}=2.
\end{align*}
\end{gather*}
Puesto que $t_0(K_0)\neq t_0(\K )$, 
se sigue que $\g{(K_0)}\neq \g \K $.

Sea $F_0=F_{{\mathcal P}_T} 
F_{{\mathcal P}_{T^2+T+1}}$. Tenemos
\begin{gather*}
[F_{{\mathcal P}_T}:K]=c_{{\mathcal P}_T}
=\mcd(e_{{\mathcal P}_T},
5^{d_{{\mathcal P}_T}}-1)=(3,4)=1.
\intertext{Por tanto $F_{{\mathcal P}_T}=K$.
Ahora bien}
\begin{align*}
[F_{{\mathcal P}_{T^2+T+1}}:K]&=c_{{\mathcal P}_{T^2+T+1}}=
\mcd(e_{{\mathcal P}_{T^2+T+1}}, 
q^{\deg {\mathcal P}_{T^2+T+1}}-1)\\
&= \mcd(3,24)=3.
\end{align*}
\end{gather*}

Puesto que $q-1=e_{\p}(K(\Lambda_{T^2+T+1})|K)=4$
se sigue que 
\[
F_{{\mathcal P}_{T^2+T+1}}
\subseteq K(\Lambda_{T^2
+T+1})^+\quad\text{y}\quad F_0=F_{{\mathcal P}_{T^2+T+1}}=F_0\cap
K(\Lambda_{T(T^2+T+1)})^+. 
\]
Por lo tanto
$F_0$ es el \'unico subcampo de
$K(\Lambda_{T^2+T+1})$ de grado $3$ 
sobre $K$ y $F_0=F$.

Tenemos que $F_{{\mathcal P}_{T^2+T+1}}\cdot {\ma F}_{5^2}/
{\ma F}_{5^2}(T)$ es una extensi\'on de
Kummer de grado $3$ donde los primos finitos
ramificados son los primos en ${\ma F}_{5^2}(T)$ 
que dividen a $T^2+T+1=(T-\zeta_3)
(T-\zeta_3^2)$. Puesto que $\p$ se descompone
totalmente en $F_{{\mathcal P}_{T^2+T+1}}
/K$, ${\eu p}_{\infty}$, el primo infinito en
${\ma F}_{5^2}(T)$, se descompone totalmente
en $F_{{\mathcal P}_{T^2+T+1}}
\cdot {\ma F}_{5^2}/{\ma F}_{5^2}(T)$. Por lo tanto 
$F_{{\mathcal P}_{T^2+T+1}} \cdot
{\ma F}_{5^2}={\ma F}_{5^2}(T)
\big(\sqrt[3]{(-1)^{\deg Q}Q(T)}\big) = {\ma F}_{5^2}(T)
\big(\sqrt[3]{Q(T)}\big)$ con 
$\deg Q(T)=3$ y $T-\zeta_3, T-\zeta_3^2$ 
son los \'unicos polinomios irreducibles que
dividen a $Q(T)$.
\[
\xymatrix{
F_{{\mathcal P}_{T^2+T+1}}\ar@{-}[r]
\ar@{-}[d]_3 & F_{{\mathcal P}_{T^2+T+1}}
\cdot {\ma F}_{5^2}={\ma F}_{25}(T)\big(\sqrt[3]{Q(T)}\big)
\ar@{-}[d]^3\\
K\ar@{-}[r]&{\ma F}_{25}(T)
}
\]

Se sigue que 
\begin{align*}
F_{{\mathcal P}_{T^2+T+1}}
\cdot {\ma F}_{5^2}&={\ma F}_{25}(T)\Big(\sqrt[3]{(T-\zeta_3)(T-\zeta_3^2)^2}\Big)\\
&={\ma F}_{25}(T)\Big(\sqrt[3]{(T-\zeta_3)^2(T-\zeta_3^2)}\Big).
\end{align*}
Ahora, puesto que $F_0=F$, 
de la Observaci\'on \ref{O3.5.B} obtenemos
\begin{gather*}
\g \K =\K F{\ma F}_{5^2}={\ma F}_{25}\Big(T,\sqrt[3]{T(T-\zeta_3)(T-\zeta_3^2)},
\sqrt[3]{(T-\zeta_3)(T-\zeta_3^2)^2}\Big).
\end{gather*}

Sea $\g {\K ^{\prime}}$ el campo de g\'eneros de
$\K/{\ma F}_{25}(T)$. Podemos aplicar el Teorema
\ref{T5.1.6}. Con las notaciones de ah\'i,
tenemos $r=3, P_1=T, P_2=T-\zeta_3, P_3=
T-\zeta_3^2, \gamma =1, \alpha=(-1)^{\deg 
D}\gamma =-1\in({\ma F}_{25}^{
\ast})^3, a_1=a_2=2$. Por lo tanto
\[
\g {\K ^{\prime}}={\ma F}_{25}\Big(T, \sqrt[3]{T
(T-\zeta_3^2)^2},\sqrt[3]
{(T-\zeta_3)(T-\zeta_3^2)^2}\Big).
\]

Tambi\'en tenemos que 
\begin{multline*}
{\ma F}_{25}\Big(T,\sqrt[3]{T(T-\zeta_3)(T-\zeta_3^2)},
\sqrt[3]{(T-\zeta_3)(T-\zeta_3^2)^2}\Big) \\
={\ma F}_{25}\Big(T, \sqrt[3]{T(T-\zeta_3^2)^2},\sqrt[3]
{(T-\zeta_3)(T-\zeta_3^2)^2}\Big).
\end{multline*}
Por lo tanto $\g \K =\g {\K ^{\prime}}$ (ver
Observaci\'on \ref{O12*.2.2.L}).

Finalmente, de la Observaci\'on 
\ref{O3.5.B} tenemos que
$\g {(K_0)}=\K F=K\big(T,\sqrt[3]{D(T)}\big)
F_{{\mathcal P}_{T^2+T+1}}$.
\end{ejemplo}

\subsection{Sobre el problema del encaje en el problema
inverso de la teor\'ia de Galois}\label{S14.9.2}

Consideremos $K=\F(T)$ y sea $l$ un n\'umero primo tal que
$l|q-1$ y $l^2\nmid q-1$. Sea $\K:=K(\sqrt[l]{\gamma D})$ con
$\gamma\in\*\F$, $D=P_1^{\alpha_1}\cdots P_r^{\alpha_r}$
con $P_1,\ldots,P_r\in R_T^+$ distintos y $1\leq \alpha_i\leq
l-1$, $1\leq i\leq r$. Se tiene que $\p$ es ramificado en $\K/K$
si y solamente si $l\nmid \deg D$. 

Sea $E:=\K K_m\cap \cicl N{}$ donde ponemos
$K_m={\ma F}_{q^m}(T)$, $\K\subseteq \cicl N{}_m$ y $n=0$
pues $\p$ es moderadamente ramificado en $\K/K$.
Entonces $E=K(\sqrt[l]{D^*})$ donde $D^*=(-1)^{\deg D} D$. 
Notemos que existe $K\subseteq \K\subseteq L$ con $\Gal(L/K)
\cong C_{l^2}$ si y solamente si existe $K\subseteq E\subseteq
F\subseteq \cicl {N_1}{}$ para alg\'un $N_1$ con $\Gal(F/K)\cong
C_{l^2}$. En efecto, si $L$ existe, entonces $F:=LK_m\cap 
\cicl {N-1}{}$. Rec\'iprocamente, si $F$ existe, sea $L:=F\K$.
\[
\xymatrix{
F\ar@{-}_l[d]\ar@{-}[r]&F\K=L\ar@{-}^l[d]\\
E\ar@{-}_l[d]\ar@{-}[r]&\K \ar@{-}[dl]^l\\ K
}
\]

Si $\p$ se ramifica en $E$, entonces $\p$ es totalmente
ramificado en $F/K$ pues si $I_{\infty}(F|K)$ es el grupo de inercia
de $\p$ en $F/K$, se tiene $I_{\infty}(F|K)\subseteq \Gal(F/K)\cong
C_{l^2}$ y como $\p$ es ramificado en $E/K$, $I_{\infty}(F|K)
\not\cong \Gal(F/E)$ el cual es el \'unico subgrupo de orden $l$ de
$\Gal(F/K)$. Se sigue que 
$I_{\infty}(F|K)\cong C_{l^2}$, $I_{\infty}(F|K)=
\Gal(F/K)\cong C_{l^2}$ y $e_{\infty}(F|K)=l^2$. Por otro lado,
se tiene $e_{\infty}(\cicl N{}|K)=q-1$ lo cual 
implica que $e_{\infty}(F|K)
=l^2|q-1=e_{\infty}(\cicl N{}|K)$ lo cual contradice las hip\'otesis.
Se sigue que $\p$ no puede ser ramificado en $E/K$ y por tanto
$l|\deg D$.

Ahora supongamos $l|\deg D$. Sea $d_i=\deg P_i$,
$1\leq i\leq r$. Sea $X=\langle
\mu\rangle$ el grupo de caracteres de Dirichlet asociado a $E=
K(\sqrt[l]{D^*})$. Entonces, en caso de existir $F$, $P_i$ ser\'a
totalmente ramificado en $F/K$ y por tanto 
$e_{P_i}(F|K)=l^2$. En general
tenemos que una extensi\'on abeliana moderadamente ramificada,
el \'indice de ramificaci\'on en $P$ divide a $q^d-1$ donde $d=\deg P$.

Necesitamos que $l^2|q^{d_i}-1$, $d_i=\deg P_i$, $1\leq i\leq r$. 
Ahora, $l|q-1$ por tanto $q\equiv 1\bmod l$ y $q^{d_i}-1=(q-1)(
q^{d_i-1}+q^{d_i-2}+\cdots+q+1)$. Por otro lado se tiene
 $q^{d_i-1}+q^{d_i-2}+\cdots+q+1
\equiv \underbrace{1+1+\cdots+1+1}_{d_i}\equiv d_i\bmod l$. En
particular, $l^2|q^{d_i}-1$ si y solamente si $l|q^{d_i-1}+q^{d_i-2}
+\cdots+q+1$ si y solamente si $l|d_i$. Por tanto $d_i=l d'_i$ y 
$l|\deg P_i$.

Rec\'iprocamente, si $E=K(\sqrt[l]{D^*})$ con $ld_i=\deg P_i$, $1\leq i
\leq r$ y $l^2|q^{d_i}-1$, entonces existe un caracter de conductor 
$P_i$ y orden $l^2$ ya que $\Gal(\cicl Pi/K)\cong C_{q^{d_i}-1}$.
Sea $\chi_{P_i}^l=\mu_{P_i}$ y sean $Q_1,\ldots,Q_s\in R_T^+$,
distintos tales que 
$\{P_1,\ldots,P_r\}\cap \{Q_1,\ldots,Q_s\}=\emptyset$.
Sean $\chi_{Q_j}$ caracteres de conductor $Q_j$ con
$\chi_{Q_j}^l=1$, $1\leq j\leq s$. 
Sea $\chi=\chi_{P_1}\cdots\chi_{P_r}\chi_{Q_1}\cdots \chi_{Q_s}$,
$\chi^l=\mu$. Si $F$ es el campo asociado a $\langle\chi\rangle$,
$E\subseteq F$ y $\Gal(F/K)\cong\langle\chi\rangle\cong C_{l^2}$.
Se tiene

\begin{teorema}\label{T14.9.2.1}
Con las hip\'otesis anteriores, sea $\K=K(\sqrt[l]{\gamma D})$. Entonces
existe $K\subseteq\K\subseteq L$ con $\Gal(L/K)\cong C_{l^2}$
si y solamente si $l|\deg P_i$, $1\leq i\leq r$. $\fin$
\end{teorema}

\begin{ejemplo}\label{T14.9.2.2}
No existe una extensi\'on $F/K$ c\'iclica de grado $4$ con $q=3$
donde $T$ sea totalmente ramificado.
En efecto, en caso de existir, sea $\K$ la subextensi\'on de grado $2$:
$K\subseteq \K\subseteq F$, $[\K:K]=2$. Entonces $\K/K$ es Kummer
y $\K=K(\sqrt[2]{\gamma D})$, $T|D$ y $T^2\nmid D$. Esto no
es posible pues $2\nmid \deg T=1$.

Se sigue que para $q=3$, el problema del encaje 
no tiene soluci\'on para $\K=
K(\sqrt{TD})$ donde $D\in R_T$ y $T\nmid D$. 
Es decir, no existe una extensi\'on c\'iclica $L$ de $K$ de
grado $4$ conteniendo a $\K$. Esto se cumple para $q$
impar arbitrario.
\end{ejemplo}

%% file: Capitulo15.tex
\chapter{M\'odulos de Drinfeld}\label{DrinfeldCh15}

\section{Introducci\'on}\label{DrinfeldS15.1}

Los {\em m\'odulos de Drinfeld} aparecieron como tales en el trabajo de
V. Drinfeld en 1974 (\cite{Dri74}), aunque ya L. Carlitz hab\'ia
descubierto en 1935 (\cite{Car35}) el primer m\'odulo de Drinfeld,
a saber, el ahora conocido como el m\'odulo de Carlitz. En el
trabajo de Drinfeld, los que ahora conocemos como m\'odulos
de Drinfeld fueron llamados {\em m\'odulos el\'ipticos} por el
mismo Drinfeld por su similitud con las curvas el\'ipticas.

El objetivo de este cap\'itulo es presentar una introducci\'on a los
m\'odulos de Drinfeld, incluyendo su aplicaci\'on a la teor\'ia de
campos de clase. 

En la Secci\'on \ref{DrinfeldC1} presentamos las propiedades b\'asicas de los
polinomios aditivos. Aqu\'i hacemos menci\'on de que la diferencia
entre la teor\'ia expl\'icita de campos de clase que se puede hacer
en el caso de los campos de funciones, con respecto a la de los
campos num\'ericos, es debido a que hay muchos polinomios aditivos
en caracter\'istica positiva, a diferencia de los pocos existentes en
caracter\'istica $0$.

En este mismo cap\'itulo, presentamos algunos resultados que necesitaremos
sobre dominios Dedekind y la torsi\'on de m\'odulos sobre estos
dominios Dedekind. Finalizamos con el estudio del n\'umero de clase
del tipo de dominios Dedekind que se presentan en el estudio
de los m\'odulos de Drinfeld.

La Secci\'on \ref{DrinfeldC2} presenta las propiedades b\'asicas de los m\'odulos
de Drinfeld. La parte central de esta secci\'on es probar la existencia
de m\'odulos de Drinfeld sobre el campo $\Ci=\ma C_p$, el cual
es el an\'alogo al campo de los n\'umeros complejos en caracter\'istica
$p>0$. La construcci\'on de estos m\'odulos se debe a Drinfeld mismo
y usa teor\'ia anal\'itica: redes, funciones exponenciales, teorema de
uniformizaci\'on anal\'itica, etc. Finalizamos el cap\'itulo presentando
el c\'alculo de la red para el m\'odulo de Carlitz, calculado por
\'el mismo y damos el concepto de morfismo entre m\'odulos de
Drinfeld.

En la Secci\'on \ref{DrinfeldC3} presentamos una introducci\'on a la teor\'ia 
de campos de clase expl\'icita, desarrollada primero por Drinfeld mismo
y despu\'es explicitada por D. Hayes. Seguimos muy de cerca el desarrollo
de Hayes de \cite{Hay92}. Se estudia el campo de clase de Hilbert y
el campo de clase de Hilbert extendido de un campo de funciones
congruente arbitrario. Aqu\'i se presenta la diferencia fundamental
entre los campos de funciones y los campos num\'ericos: en nuestro
caso, Hayes hace una descripci\'on expl\'icita de la m\'axima extensi\'on
abeliana de un campo de funciones global $\K $. En el caso num\'erico,
esta descripci\'on expl\'icita \'unicamente se conoce para el campo de
los n\'umeros racionales $\ma Q$ y para extensiones cuadr\'aticas
imaginarias de $\ma Q$.

\section{Dominios Dedekind y el m\'odulo de Carlitz}\label{DrinfeldC1}

El Cap\'itulo \ref{Ch6} presenta un estudio detallado del m\'odulo de
Carlitz, esto es, los campos de funciones ciclot\'omicos.

\subsection{Propiedades b\'asicas y polinomios aditivos}

Sean $A=R_T=\F[T]$\label{Drinfeldanillopolinomios}, 
$K=\F(T)$\label{Drinfeldcampofuncionesracionales} el campo de funciones racionales. Sea
$\bar{K}$ una cerradura algebraica de $K$. Sea
\begin{gather*}
\begin{align*}
\End_{\F}(\bar{K})=&\{\varphi\colon \bar{K}\to \bar{K}\mid \varphi(a+b)=\varphi(a)+
\varphi(b), \varphi(\alpha a)=\alpha \varphi(a)\\
&\text{\ para toda $\alpha\in
\F$ y para cualesquiera $a,b\in \bar{K}$}\},
\end{align*}
\intertext{la $\F$--\'algebra
de endomorfismos de $\bar{K}$ sobre $\F$\index{grupo de endomorfismos}.
Sean $\tau,\mu_T\in\End_{\F}(\bar{K})$, dados por}
\mu_T\colon \bar{K}\lra \bar{K},\quad \mu_T(u)=Tu\quad\text{y}\quad
\tau\colon \bar{K}\lra \bar{K},\quad \tau(u)=u^q.
\end{gather*}\label{Drinfeldendomorfismos}

Entonces el {\em m\'odulo de Carlitz\index{m\'odulo de Carlitz}\index{Carlitz!m\'odulo
de $\sim$}}\label{DrinfeldmoduloCarlitz} 
es el homomorfismo de $\F$--\'algebras $C\colon A\to \End_{\F}(\bar{K})$
dado por $C(M)(u):=M(\tau+\mu_T)(u)$ para $M\in A$. Esto es, si $M\in R_T$ con
$M=a_dT^d+\cdots + a_1T+a_0$, entonces 
\[
C(M)(u)=a_d(\tau+\mu_T)^d(u)+\cdots+a_1(\tau+\mu_T)(u)+a_0(u),
\]
donde $(\tau+\mu_T)^i=\underbrace{(\tau+\mu_T)\circ\cdots\circ (\tau+\mu_T)}_{
i\text{\ veces}}$ es la composici\'on y $a_0(u)=a_0u$.
Denotamos $C_M=C(M)$. En otros contextos se utiliza la notaci\'on
$C_M(u)=u^M=M(\tau+\mu_T)(u)$.

En resumen, se tiene $C(\alpha)=\alpha\colon \bar{K}\to\bar{K}$, $u\mapsto \alpha u$,
esto es, $C_{\alpha}(u)=\alpha u$,
para $\alpha\in \F$ y $u\in \bar{K}$, y $C(T)=\tau+\mu_T\colon \bar{K}\to\bar{K}$, $u\mapsto
u^q+Tu$, es decir, $C_T(u)=u^q+Tu=Tu+u^q$ y $C(MN)=C(M)\circ C(N)$ para $M,N
\in R_T$. Esto \'ultimo significa 
\[
C_{MN}(u)=C_M(C_N(u))=C_N(C_M(u))=C_{NM}(u).
\]

Tambi\'en notemos que
\begin{gather*}
(\tau\circ \mu_T)(u)=\tau(Tu)=T^qu^q\quad\text{y}\\
(\mu_T^q\circ \tau)(u)=\mu_T^q(u^q)=\underbrace{\mu_T\circ\cdots\circ \mu_T}_{
q\text{\ veces}}(u^q)=\underbrace{\mu_T\circ\cdots\circ \mu_T}_{
q-1\text{\ veces}}(Tu^q)=\ldots = T^qu^q.
\end{gather*}
Es decir $\tau\circ \mu_T=\mu_T^q\circ \tau$.

Sea ahora $a\in\F$, $C_a(u)=au$, y se tiene $\tau\circ a\colon \bar{K}\to\bar{K}$,
$(\tau\circ a)(u)=\tau(au)=a^qu^q=(a^q\circ \tau)(u)$, esto es, \fbox{$\tau\circ
a=a^q \circ \tau$}.

Ahora bien, puesto que $C_M(u+w)=C_M(u)+C_M(w)$ para $u,w\in \bar{K}$ y
$C_M(au)=aC_M(u)$ para $a\in\F$, entonces $C_M$ es un {\em polinomio
aditivo\index{polinomio aditivo}}.

\begin{definicion}\label{DrinfeldD1.1.1} Sea $F$ un campo cualquiera y sea $f(x)\in F[x]$
un polinomio. $f(x)$ se llama {\em aditivo} si para todas $\alpha,\beta\in F$,
se tiene $f(\alpha+\beta)=f(\alpha)+f(\beta)$.
\end{definicion}

\begin{ejemplo}\label{DrinfeldEj1.1.2}
Si $f$ y $g$ son aditivos, entonces $f+g$, $\alpha f$ con $\alpha\in F$, $f
\circ g$ y $f(x)=ax$ con $a\in F$, son aditivos.

Si $\car F=p>0$, entonces $f(x)=x^{p^i}$ es aditivo. En particular, $\tau_p^i(x):=
x^{p^i}$ y los polinomios generados por $\{\tau_p^i\}_{ieq 0}$ son polinomios
aditivos.
\end{ejemplo}

\begin{definicion}\label{DrinfeldD1.1.3}
Sea $F$ un campo de caracter{\'\i}stica $p>0$. Se define
\[
F\langle\tau_p\rangle=\Big\{\sum_{i=0}^n a_i\tau_p^i=\sum_{i=0}^n a_ix^{p^i}\mid
a_i\in F\Big\}\subseteq F[x].
\]
\end{definicion}

Se tiene que $\aditivos Fp$ consiste de polinomios aditivos.

\begin{observacion}\label{DrinfeldO1.1.4}
Notemos que con la multiplicaci\'on, $F\langle\tau_p\rangle$ no es un subanillo
de $F[x]$: $\tau_p\cdot \tau_p=x^px^p=x^{2p}\notin F\langle\tau_p\rangle$ para
$p\geq 3$ (si $p=2$, $\tau_2\cdot\tau_2\cdot\tau_2=x^6\notin F\langle \tau_2\rangle$).

Por otro lado, $F\langle \tau_p\rangle$ es un anillo con la {\underline{composici\'on}}:
$\tau_p\circ \tau_p=\tau_p^2$, $\tau_p^2(x)=x^{p^2}$.

Si $F\neq {\ma F}_p$, $F\langle\tau_p\rangle$ es no conmutativo pues si $\alpha\in F$,
entonces $\tau_p \alpha=\alpha^p\tau_p$, esto es, $(\tau_p\circ\alpha)(x)=(\alpha x)^p
=\alpha^px^p=\alpha^p\tau_p(x)=(\alpha^p\circ \tau_p)(x)$ y existe $\alpha\in F$
con $\alpha\notin {\ma F}_p$, es decir, $\alpha^p\neq \alpha$.
\end{observacion}

\begin{observacion}\label{DrinfeldO1.1.5} Ver tambi\'en la Observaci\'on
\ref{O12.3.4.pea}. En caracter{\'\i}stica $p>0$, todos los polinomios en 
$F\langle\tau_p\rangle$ son aditivos pero pueden existir polinomios aditivos que no
est\'an en $\aditivos Fp$. Por ejemplo, si $F={\ma F}_p$, $p\geq 3$ y $f(x)=g(x)+
(x^p-x)^n$ para $n\in {\ma N}$ y $g(x)\in\aditivos Fp$, es un polinomio aditivo pues
$f(\alpha)=g(\alpha)$ para toda $\alpha\in F$ pero, por ejemplo, para $g(x)=0, n=2$
$f(x)=(x^p-x)^2=x^{2p}-2x^{p+1}+x^2\notin \aditivos Fp$.
\end{observacion}

La raz\'on de lo anterior, es que $F$ es finito, digamos $F=\F=\{u\in\bar{\ma F}_p\mid
u^q-u=0\}$, y por lo tanto $(\alpha^q-\alpha)^n=0$ para toda $n\in{\ma N}$ y toda $\alpha
\in\F$.

\begin{proposicion}\label{DrinfeldP1.1.6}
Sean $F$ un campo infinito y $f(x)\in F[x]$ un polinomio aditivo. Entonces, si $\car F=0$,
se tiene que $f(x)=ax$ para alg\'un $a\in F$. Si $\car F=p>0$, entonces
$f(x)\in \aditivos Fp$.
\end{proposicion}

\begin{proof} Sea $f(x)=\sum_{i=0}^n a_ix^i$. 
Entonces $f(0)=f(0+0)=f(0)+f(0)$, por tanto $f(0)=
a_0=0$. Se tiene $f^{\prime}(x)=\sum_{i=1}^n ia_i x^{i-1}$. Sea $g(x)=f(x+\alpha)-
f(x)-f(\alpha)$ para $\alpha\in F$. Entonces $g(\beta)=0$ para toda $\beta\in F$. Puesto
que $F$ es infinito, se sigue que $g(x)=0$, esto es, $f(x+\alpha)=f(x)+f(\alpha)$.
Derivando esta expresi\'on, obtenemos
\[
f^{\prime}(\alpha)=\frac{d}{dx}(f(x+\alpha))|_{x=0}=\frac{d}{dx}\big(
f(x)+f(\alpha)\big)|_{x=0}=f^{\prime}(0).
\]
Es decir, $f^{\prime}(\alpha)=f^{\prime}(0)$ para toda $\alpha\in F$. Puesto que $F$
es infinito, $f^{\prime}(x)=f^{\prime}(0)=c\in F$. Por tanto, puesto que $f(x)=\sum_{
i=1}^n a_ix^i$, entonces $f^{\prime}(x)=\sum_{i=1}^n i a_i x^{i-1}=c
=a_1=f^{\prime}(0)$. Se
sigue que $ia_i=0$, $i=2,3,\ldots, n$, en $F$.

Si $\car F=0$, $a_i=0$, $i=2,\ldots,n$, por tanto $f(x)=cx$. Si $\car F=p>0$, $ia_i=0$,
por lo que $a_i=0$ para $i\not\equiv 0\bmod p$. Se sigue que 
\[
f(x)=\sum_{l=1}^m
a_{lp}x^{lp}=\sum_{j=0}^{m^{\prime}}a_{p^j}x^{p^j}+\sum_{\substack{
\text{$s$ no es po-}\\
\text{tencia de $p$}}}a_s x^s=f_0(x)+f_1(x).
\]
Puesto que $f(x)$ es aditivo, entonces
$f_1(x)=f(x)-f_0(x)$ es aditivo.

Sea $f_1(x)=f_2(x)^p$ con $f_2(x)\in F^{1/p}[x]$. Veamos que $f_2(x)$ es aditivo
en $F^{1/p}=\{\alpha^{1/p}\mid \alpha\in F\}$, el cual es un campo que
contiene a $F$. Sean $\alpha^{1/p},\beta^{1/p}\in
F^{1/p}$, $\alpha,\beta\in F$. Entonces
\begin{align*}
f_2(\alpha^{1/p}+\beta^{1/p})&=f_1(\alpha^{1/p}+\beta^{1/p})^{1/p}= f_1(\alpha^{1/p})^{1/p}
+f_1(\beta^{1/p})^{1/p}\\
&=f_2(\alpha^{1/p})+f_2(\beta^{1/p}).
\end{align*}
Esto es, $f_2$ es aditivo por lo que $f_2^p$ es aditivo. Por inducci\'on en el grado, 
podemos suponer que $f_2(x)=\sum_{i=1}^m c_ix^{p^i}$ pues $\deg f_2
<\deg f$. Se sigue que
$f_1(x)=f_2(x)^p=\sum_{i=1}^m c_i^p x^{p^{i+1}}$ con $c_i^p\in F$. Finalmente
obtenemos que $f(x)=f_0(x)+f_1(x)=\sum_{j=0}^m a_j x^{p^j}$. 
$\fin$ \end{proof}

Otra versi\'on de la Proposici\'on \ref{DrinfeldP1.1.6} la dimos en la Proposici\'on
\ref{P12.3.3.pea}.

\begin{observacion}\label{DrinfeldO1.1.6(1)} 
Si $F=\F$ y $f(x)\in F[x]$ es aditivo, entonces dividiendo entre $x^q-x$, obtenemos
$f(x)=h(x)+(x^q-x)l(x)$ con $\deg h(x)<q$. 

Se tiene $h(\alpha)=f(\alpha)$ para
toda $\alpha\in F$. Sea $g(x)=h(x+\alpha)-h(x)-h(\alpha)$ con $\deg g<q$.
Adem\'as $g(\beta)=0$ para toda $\beta\in F$, lo cual implica que $g(x)\equiv 0$,
esto es, $h(x+\alpha)=h(x)+h(\alpha)$.

Se puede repetir la demostraci\'on de la Proposici\'on \ref{DrinfeldP1.1.6} usando que
$F^p=F$ y obtenemos que $h(x)=\sum_{i=0}^ma_i x^{p^i}$ donde $m<r$, con
$q=p^r$. En otras palabras
\[
f(x)=\sum_{i=0}^ma_ix^{p^i}+(x^q-x)l(x)
\]
con $m<r$ y $l(x)\in F[x]$.
\end{observacion}

De ahora en adelante supondremos que $F$ es infinito y de caracter{\'\i}stica $p>0$.
Sea ${\mathcal P}(F)=\{f(x)\in F[x]\mid \text{$f$ es aditivo}\}$\label{Drinfeldaditivos}. Entonces
$\theta\colon {\mathcal P}(F)\longrightarrow \aditivos Fp$, $\theta(f(x))=g(\tau_p)(x)$, donde
$f(x)=\sum_{i=0}^m a_ix^{p^i}$, $g(\tau_p)=\sum_{i=0}^ma_i\tau_p^i$, es una
biyecci\'on. 

Recordemos que ${\mathcal P}(F)$ no es cerrado bajo multiplicaci\'on y que $\aditivos 
Fp$ es un anillo bajo la composici\'on.

\begin{definicion}\label{DrinfeldD1.1.7} El anillo $\aditivos Fp$ se llama {\em el anillo de polinomios
torcidos sobre $F$\index{polinomios torcidos}}. Se tiene $\tau_p\alpha=\alpha^p \tau_p$
para $\alpha\in F$.
\end{definicion}

Ahora bien, si queremos que $f(\tau_p)(\alpha u)=\alpha f(\tau_p)(u)$ para toda
$\alpha\in \F$ con $q=p^r$, entonces veamos que $f(\tau_p)=g(\tau)$ con $\tau=\tau_p^r$.
Primero, $\tau(\alpha u)=(\alpha u)^q=\alpha^q u^q \underbracket[0pt]{=}_{\substack{
\uparrow\\ \alpha\in\F}}\alpha u^q=\alpha\tau(u)$, por tanto $g(\tau)\in\End_{\F}(B)$
donde $B$ es cualquier $\F$--\'algebra y
\begin{align*}
\End_{\F}(B)=&
\{\varphi\colon B\to B\mid \varphi(u+v)=\varphi(u)+\varphi(v), \varphi(\alpha u)=
\alpha \varphi(u) \\
&\text{\ para toda $\alpha\in \F$ y para todas $u,v\in B$}\}.
\end{align*}

Sea $f(\tau_p)=\sum_{i=0}^d a_i\tau_p^i$, $f(\tau_p)(u)=\sum_{i=0}^d a_i u^{p^i}$
para $u\in F$ y donde adem\'as suponemos que $\F\subseteq F$. Entonces tenemos
$f(\tau_p)(\alpha u)=\sum_{i=0}^d a_i \alpha^{p^i} u^{p^i}=\alpha \sum_{i=0}^d a_i u^{
p^i} = \alpha f(\tau_p)(u)\iff$ para toda $i$ tal que $a_i\neq 0$, $\alpha^{p^i}=\alpha$ para toda $\alpha
\in\F \iff r\mid i$ para toda $i$ tal que $a_i\neq 0$, por tanto, $f(\tau_p)=
\sum_{j=0}^{d^{\prime}} a_j \tau_p^{rj} = \sum_{j=0}^{d^{\prime}} a_j \tau^j= g(\tau)
\in \aditivos {F}{} $ donde $\tau=\tau_q=\tau_p^r$.

De ahora en adelante, a menos que se diga lo contrario, se supondr\'a $\F\subseteq
F$ y $g(\tau)\in \aditivos F{}$ con $\tau=\tau_p^r$, es decir, $\tau(u)=u^q$.

\begin{definicion}\label{DrinfeldD1.1.8} Decimos que $f(\tau)$ es {\em divisible por la
derecha (resp. por la izquierda)\index{divisible por la 
derecha}\index{divisible por la izquierda}}
por $g(\tau)$ si existe $h(\tau)\in \aditivo F$ tal que $f(\tau)=h(\tau)g(\tau)$ (resp. $f(\tau)=
g(\tau)h(\tau)$).
\end{definicion}

\begin{observacion}\label{DrinfeldO1.1.9}
Si $\serie fda$, $a_d\neq 0$, se define $\deg_{\tau} f$ o, simplemente,
$\deg f$, por $d$. Se tiene que $\deg (f\circ g)=\deg f+\deg g$ pues si
$\serie gsb$, $b_s\neq 0$, $f(\tau)\circ g(\tau)=\Big(\sum_{i=0}^d a_d\tau^d\Big)
\Big(\sum_{j=0}^s b_j\tau^j\Big)=\sum_{i,j}a_i \tau^i b_j \tau^j=\sum_{i,j}a_ib_j^{q^i}
\tau^{i+j}=\sum_{l=0}^{d+s} c_l \tau^l$ con $c_{d+s}=a_d b_s^{q^d}\neq 0$.

En particular $\aditivo F$ es un dominio entero no conmutativo.

Notemos que $\deg_x f=q^{\deg_{\tau} f}$.
\end{observacion}

\begin{proposicion}[Algoritmo derecho de la divisi\'on\index{algoritmo
derecho de la divisi\'on}]\label{DrinfeldP1.1.10}
Consideremos $f(\tau), g(\tau)\in \aditivo F$ con $g(\tau)\neq 0$. Entonces existen
$q(\tau), r(\tau)\in \aditivo F$ con $\deg r(\tau)<\deg g(\tau)$ tales que
\[
f(\tau)=q(\tau)g(\tau)+r(\tau).
\]

Adem\'as, $q(\tau)$ y $f(\tau)$ son \'unicos.
\end{proposicion}

\begin{proof} Sean $\serie fda$ y $\serie gsb$. Si $f(\tau)=0$, 
sean $q(\tau)=r(\tau)=0$.

Sea $f(\tau)\neq 0$. Si $\deg g(\tau)=s>d=\deg f(\tau)$, sean $q(\tau)=0$ y
$r(\tau)=f(\tau)$. Sea ahora $\deg g=s\leq d=\deg f$. Entonces
\begin{align*}
f(\tau)-a_d(b_s^{-1})^{q^{d-s}}\tau^{d-s}g(\tau)=&f(\tau)-\sum_{j=0}^{s-1}
a_db_s^{-q^{d-s}}b_j^{q^{d-s}}\tau^{d-s+j}\\
&-\underbrace{a_db_s^{-q^{d-s}}b_s^{q^{d-s}}\tau^{d-s+s}}_{=a_d\tau^d}=f_1(\tau)
\end{align*}
con $\deg f_1<d$.

Repitiendo el proceso, ahora con $f_1(\tau)$ y por hip\'otesis de
inducci\'on en $\deg f$, finalmente obtenemos $f(\tau)=q(\tau)g(\tau)+r(\tau)$ con
las condiciones requeridas, esto es, con $\deg r(\tau)<\deg g(\tau)$.

Supongamos ahora que  $f(\tau)=q_1(\tau)g(\tau)+r_1(\tau)=q(\tau)g(\tau)+r(\tau)$.
Entonces $(q-q_1)g=r_1-r$ con $\deg (r_1-r)<\deg g$. Por lo tanto $q-q_1=0$,
de donde $r-r_1=0$. 
$\fin$ \end{proof}

\begin{corolario}\label{DrinfeldC1.1.11}
Todo ideal izquierdo de $\aditivo F$ es principal.
\end{corolario}

\begin{proof} Sea $I\neq 0$ un ideal izquierdo de 
$\aditivo F$. Sea $g\in I$ tal que $g\neq 0$
y $\deg g\leq \deg g_1$ para toda $g_1\in I$, $g_1\neq 0$. 
Se sigue que $Rg\subseteq I$ donde
$R=\aditivo F$.

Sea $f\in I$. Por el algoritmo derecho de la divisi\'on, se tiene que $f=qg+r$ para
$q,r\in R$ y $\deg r<\deg g$. Puesto que $f\in I$ y $qg\in I$, entonces $r\in I$
por lo que $r=0$ y $f=qg\in Rg$. Se sigue que $I=Rg$. 
$\fin$ \end{proof}

\begin{observacion}\label{DrinfeldO1.1.12} En general no existe el algoritmo izquierdo de
la divisi\'on. Sin embargo, cuando $F$ es perfecto, esto es, $F^p=F$, si existe el
algoritmo izquierdo de la divisi\'on.
\end{observacion}

\begin{proposicion}[Algoritmo izquierdo de la divisi\'on\index{algoritmo izquierdo
de la divisi\'on}]\label{DrinfeldP1.1.13}
Consideremos $f(\tau), g(\tau)\in \aditivo F=R$, 
$g(\tau)\neq 0$. Si $F$ es perfecto, existen
$q(\tau),r(\tau)\in R$ con $\deg r(\tau)<\deg g(\tau)$ \'unicos tales que $f(\tau)=
g(\tau)q(\tau)+r(\tau)$.
\end{proposicion}

\begin{proof} Sean $f(\tau)\neq 0$, $\serie fda$, $a_d\neq 0$ y $\serie gsb$, $b_s\neq 0$.
Sea $s\leq d$ el cual es el \'unico caso en que se tiene que hacer algo nuevo, pues
los casos $f(\tau)=0$ y $s>d$ son lo mismo que en el caso del algoritmo derecho.

Sea
\begin{align*}
f_1(\tau)&=f(\tau)-g(\tau)b_s^{-1/q^s} a_d^{1/q^s}\tau^{d-s}\\
&=\sum_{i=0}^d a_i\tau^i-
\sum_{j=0}^{s}b_j\tau^j b^{-1/q^s}a_d^{1/q^s} \tau^{d-s}\\
&=\sum_{i=0}^d a_i\tau^i-\sum_{j=0}^s
b_jb_s^{-q^{j-s}}a_d^{q^{j-s}}\tau^{j+d-s}\\
&=\sum_{i=0}^d a_i\tau^i-\sum_{j=0}^{s-1}b_j
b_s^{-q^{j-s}}a_d^{q^{j-s}}a_d^{q^{j-s}}\tau^{j+d-s}-a_d\tau^d
\end{align*}
donde $b_s^{-1/q^s}, a_d^{1/q^s}\in F$. Por tanto $\deg f_1(\tau)<d$.
El resto se sigue como en el caso del algoritmo derecho.
$\fin$ 
\end{proof}

\begin{corolario}\label{DrinfeldC1.1.14} Si $F$ es perfecto, todo ideal 
derecho de $\aditivo F$ es principal. $\fin$
\end{corolario}

\begin{observacion}\label{DrinfeldO1.1.15}
Si $f(\tau)$ es divisible por la derecha por $g(\tau)$, entonces $g(x)\mid
f(x)$ en el sentido ordinario.

En efecto, consideremos $\serie hmc$. Para cualquier
$l(\tau)\in\aditivo F$, denotemos $l(\tau)|_x=l(x)$. M\'as
precisamente, si $\serie lr\alpha$, $l(\tau)|_x=l(x)=\sum_{i=0}^r
\alpha_i x^{q^i}$. 

Sea $g(\tau)=\sum_{j=0}^t\beta_j\tau^j$. Se
tiene $\tau g(\tau)=\sum_{j=0}^t \beta_j^q\tau^{j+1}$.

Por tanto $\tau g(\tau)|_x=\sum_{j=0}^t\beta_j^q x^{q^{j+1}}=
\Big(\sum_{j=0}^t \beta_j x^{q^j}\Big)^q=\big(g(\tau)|_x\big)^q=
g(x)^q$.

Para $i\geq 2$, $\tau^i g(\tau)|_x=\tau(\tau^{i-1}g(\tau))|_x=
\big(\tau^{i-1}g(\tau)\big)^q|_x=\Big(\tau^{i-1}g(\tau)|_x\Big)^q$.

Por inducci\'on en $i$, suponemos $(\tau^{i-1}g(\tau))|_x=
g(x)^{q^{i-1}}$. Se sigue que $\tau^ig(\tau)|x=g(x)^{q^i}$. Por tanto
\begin{align*}
h(\tau)g(\tau)|_x&=\Big[\Big(\sum_{i=0}^m c_i\tau^i\Big)g(\tau)\Big]_x=
\Big[\sum_{i=0}^t c_i\tau^i g(\tau)\Big]_x=\sum_{i=0}^t c_i \big[\tau^i
g(\tau)\big]_x\\
&=\sum_{i=0}^t c_i(g(x))^{q^i}=\Big(\sum_{i=0}^t c_i g(x)^{q^i-1}\Big) g(x).
\end{align*}

Por tanto, si $f(\tau)=h(\tau)g(\tau)$, entonces, $f(\tau)|_x=f(x)=
p(x) g(x)$ con $p(x)=\sum_{i=0}^t c_i(g(x))^{q^i-1}$, de donde
obtenemos que $g(x)\mid f(x)$.
\end{observacion}

Como consecuencia de la Observaci\'on \ref{DrinfeldO1.1.15}, se sigue que
si $f(\tau)=q(\tau) g(\tau)+r(\tau)$ con $\deg_{\tau} r(\tau)<\deg_{\tau}
g(\tau)$, entonces $f(x)=q_1(x)g(x)+r(x)$ con $\deg_x r(x)<\deg_x
g(x)$ pues como $q(x)$ y $g(x)$ son polinomios aditivos, si $m(\tau)=
q(\tau)g(\tau)$, $g(\tau)$ divide a $m(\tau)$ por la derecho, por lo que
$g(x)\mid m(x)$, esto es $m(x)=q_1(x)g(x)$. Ahora bien
\begin{gather*}
f(\tau)=q(\tau)g(\tau)+r(\tau)=m(\tau)+r(\tau)\quad\text{y}\quad
r(\tau)=f(\tau)-m(\tau)\\
\intertext{por lo tanto}
r(x)=f(x)-m(x)=f(x)-q_1(x)g(x),\quad f(x)=q_1(x)g(x)+r(x)\\
\intertext{y}
\deg_x r(x)=q^{\deg_{\tau}r(\tau)}<q^{\deg_{\tau}g(\tau)}=\deg_x g(x).
\end{gather*}

Es decir, si $r(\tau)$ es el residuo del algoritmo derecho de la divisi\'on
de $f(\tau)$ por $g(\tau)$, entonces $r(x)$ es el residuo del algoritmo
de la divisi\'on ordinario de $f(x)$ por $g(x)$.

\subsection{Dominios Dedekind y torsi\'on}\label{DrinfeldS1.2}

Recordemos la definici\'on y las propiedades elementales de los
Dominios Dedekind, los cuales ya hemos usado repetidamente
a lo largo de este trabajo.

\begin{definicion}\label{Drinfeld1.2.1}
Un dominio entero (conmutativo con unidad) $D$ se llama un
{\em dominio o anillo Dedekind\index{anillo Dedekind}\index{dominio
Dedekind}} si $D$ no es un campo y adem\'as satisface las
siguientes tres condiciones:
\las
\item Todo ideal primo ${\mathcal P}$ no cero es maximal, esto es,
$D$ es de dimensi\'on uno.
\item $D$ es noetheriano.
\item $D$ es enteramente cerrado, esto es, si $\LL =\coc D$ es el
campo de cocientes de $D$ y si $x\in \LL $ satisface una relaci\'on
$x^n+a_{n-1}x^{n-1}+\cdots+a_1x+a_0=0$ con $a_i\in D$, $0
\leq i\leq n-1$, entonces $x\in D$.
\end{list}
\end{definicion}

\begin{ejemplos}\label{DrinfeldEj1.2.2}
${\ma Z}$, $\ell [x]$ con $\ell $ un campo arbitrario, $D$ un dominio de ideales
principales (DIP), ${\mathcal O}_\LL $ el anillo de enteros de cualquier
campo num\'erico $\LL $ son todos dominios Dedekind.
\end{ejemplos}

Dado un dominio Dedekind y $\LL =\coc D$ su campo de cocientes, un
$D$--m\'odulo $M\neq 0$ con $M\subseteq \LL $ se llama {\em ideal
fraccionario\index{ideal fraccionario}} si $M$ es finitamente generado
como $D$--m\'odulo. Equivalentemente, existe $d\in D$, $d\neq 0$ tal
que $dM\subseteq D$. Como ejemplos, los ideales no ceros usuales
de $D$ son ideales fraccionarios.

\begin{teorema}\label{DrinfeldT1.2.3} Si $D$ es un dominio de Dedekind, todo ideal
fraccionario ${\eu A}$ se escribe de manera \'unica como 
\[
{\eu A}={\mathcal P}_1^{\alpha_1}\cdots {\mathcal P}_r^{\alpha_r}
\]
con ${\mathcal P}_i$ ideales primos no cero de $D$, y $\alpha_i\in
{\ma Z}$, donde para un ideal primo no cero ${\mathcal P}$ de $D$,
${\mathcal P}^{-1}=\{x\in \LL \mid x{\mathcal P}\subseteq D\}$. 
\end{teorema}

\begin{proof} \cite[Theorem 5.7.4]{Vil2006}. 
$\fin$
\end{proof}

\begin{observacion}\label{DrinfeldO1.2.3(1)}
Si $x\in \* D$ el ideal fraccionario principal se define como
$(x)=\prod_{\substack{{\mc P}\text{\ primo}\\ {\mc P}\neq 0}} 
{\mc P}^{v_{{\mc P}}(x)}$
y donde $x\in {\mc P}^{v_{{\mc P}}(x)}\setminus {\mc P}^{v_{{\mc P}}(x)+1}$ si
$v_{{\mc P}}(x)\geq 0$ y $x^{-1}\in {\mc P}^{-v_{{\mc P}}(x)}\setminus {\mc P}^{
-v_{{\mc P}}(x)+1}$ si $v_{{\mc P}}(x)<0$.
\end{observacion}

\begin{teorema}\label{DrinfeldT1.2.4}
Sea $A$ un dominio Dedekind, $\LL =\coc A$. Sea $L$ una extensi\'on
finita de $\LL $ con $[L:\LL ]=n$. Sea $B=\{\alpha\in L\mid \Irr(\alpha,x,\LL )\in A[x]\}$
la cerradura entera de $A$ en $L$. Entonces $B$ es un dominio Dedekind.
\end{teorema}

\begin{proof} \cite[Theorem 5.7.7]{Vil2006}. 
$\fin$ \end{proof}

\begin{teorema}\label{DrinfeldT1.2.5} Si $A$ es un dominio Dedekind, $\LL =\coc A$
y $A\subseteq B\subsetneqq \LL $, con $B$ anillo. Entonces $B$ es dominio
Dedekind.
$\fin$
\end{teorema}

\begin{proposicion}\label{DrinfeldP1.2.6}
Todo ideal no cero de un dominio Dedekind puede ser generado por a lo
m\'as dos elementos.
\end{proposicion}

\begin{proof} \cite[Exercise 5.10.34]{Vil2006}.
$\fin$ \end{proof}

Sea $\LL /\ell $ un campo arbitrario de funciones y sean ${\mathcal P}_1,\ldots,
{\mathcal P}_r$ un n\'umero finito de lugares de $\LL $ ($r\geq 1$). Sea 
\begin{align*}
{\mathcal O}&=\bigcap_{{\mathcal P}\notin\{{\mathcal P}_1,\ldots,{\mathcal P}_r\}}
{\mathcal O}_{\mathcal P}\\
&=\{x\in \LL \mid v_{\mathcal P}(x)\geq 0 \text{\ para
todo lugar\ }{\mathcal P}\notin\{{\mathcal P}_1,\ldots,{\mathcal P}_r\}\}.
\end{align*}

Tambi\'en se denota ${\mc O}={\mc O}_S$ donde $S=\{{\mc P}_1,
\ldots, {\mc P}_r\}$.

\begin{teorema}\label{DrinfeldT1.2.7} El anillo ${\mathcal O}$ es un dominio Dedekind.
\end{teorema}

\begin{proof} Por el Teorema de Riemann--Roch, existe $x_i\in \LL $ tal que el divisor de
polos de $x_i$, $\eta_{x_i}={\mathcal P}_i^{\alpha_i}$, $\alpha_i\geq 1$. Sea
$y:=\sum_{i=1}^r x_i$. Entonces
$\eta_y=\prod_{i=1}^r {\mathcal P}_i^{\alpha_i}$. Veamos que ${\mathcal O}$
es la cerradura entera de $\ell [y]\subseteq \ell (y)$ en $\LL $.
\[
\xymatrix{
R={\mathcal O}_\LL \ar@{-}[d]\ar@{-}[r]&\LL \ar@{-}[d]\\ \ell [y]\ar@{-}[r]&\ell (y)}
\]

Sea $R:=\{\alpha\in \LL \mid \Irr(\alpha,T,\ell (y))\in \ell [y][T]\}$ la cerradura entera
de $\ell [y]$ en $\LL $. Entonces, si $\alpha\in R$, se tiene $\alpha^n+a_{n-1}
\alpha^{n-1}+\cdots+a_1\alpha+a_0=0$ con $a_i\in \ell [y]$. Notemos que la
conorma de ${\mc P}_{\infty}$ satisface que
$\con_{\ell (y)/\LL }{\mathcal P}_{\infty}
={\mathcal P}_1^{\alpha_1}\cdots {\mathcal P}_r^{
\alpha_r}$ puesto que $\eta_y={\mathcal P}_{\infty}$ en $\ell (y)$. Si ${\mathcal P}
\notin\{{\mathcal P}_1,\ldots,{\mathcal P}_r\}$, $v_{\mathcal P}(a_i)\geq 0$ puesto
que ${\mathcal P}|_{\ell (y)}\neq {\mathcal P}_{\infty}$. 

Veamos que $v_{\mc P}(\alpha)\geq 0$.
Si $v_{\mathcal P}(\alpha)
<0$, $v_{\mathcal P}(\alpha^n)=nv_{\mathcal P}(\alpha)<iv_{\mathcal P}(\alpha)
\leq iv_{\mathcal P}(\alpha)+v_{\mathcal P}(a_i)=v_{\mathcal P}(a_i\alpha^i)$,
$0\leq i\leq n-1$.
Por lo tanto 
\[
\infty =v_{\mathcal P}(0)=v_{\mathcal P}(\alpha^n+a_{n-1}\alpha^{n-1}+\cdots
+a_1\alpha+a_0)=nv_{\mathcal P}(\alpha)<0,
\]
lo cual es absurdo. Se sigue que $v_{\mathcal P}(\alpha)
\geq 0$ y $\alpha\in{\mathcal O}$.
De esta forma obtenemos que $R\subseteq {\mathcal O}$. Probaremos la
otra contenci\'on, sin embargo, es suficiente esta primera contenci\'on
para el resultado pues $R$
es dominio Dedekind y $R\subseteq {\mathcal O}\subseteq \coc R=\LL $ (ver
Teorema \ref{DrinfeldT1.2.5}).

Sea $\alpha\in {\mathcal O}$. Por tanto $v_{\mathcal P}(\alpha)\geq 0$ para
toda ${\mathcal P}\notin\{{\mathcal P}_1,\ldots, {\mathcal P}_r\}$. Sea $f(T)=
\Irr(\alpha,T,\ell (y))=T^m+b_{m-1}T^{m-1}+\cdots+b_1T+b_0\in \ell (y)[T]$. Entonces
$\alpha^m+b_{m-1}\alpha^{m-1}+\cdots+b_1\alpha+b_0=0$ con $b_0\neq 0$.

Sea $\tilde{\LL }$ la cerradura normal de $\LL /\ell (y)$. Sean $\alpha=\alpha^{(1)},
\alpha^{(2)},\ldots, \alpha^{(m)}$ los conjugados de $\alpha$
con multiplicidades, esto es, 
$f(T)=\prod_{i=1}^m(T-\alpha^{(i)})$ ya sea separable o no. Entonces $b_i$
es una funci\'on sim\'etrica en $\alpha^{(1)},\ldots,\alpha^{(m)}$, $b_i=
\sum \alpha^{(j_1)}\cdots \alpha^{(j_{m-i})}$ y $v_{\eu p}(\alpha^{(j)})\geq 0$
donde ${\eu p}$ es un primo en $\tilde{\LL }$ sobre ${\mathcal P}$. Por lo tanto
$v_{\eu p}(b_i)\geq 0$, de donde $v_{\eu q}(b_i)\geq 0$ para todo 
primo de $\ell (y)$ con ${\eu q}\neq {\mathcal P}_{\infty}$. Se sigue
que $b_i\in \ell [y]$ y por tanto $\alpha\in R$ y por tanto $\o\subseteq R$.

Se sigue que \fbox{${\mathcal O}=R$}. Finalmente, puesto que $\ell [y]$ es un DIP,
$\ell [y]$ es dominio Dedekind y por el Teorema \ref{DrinfeldT1.2.4} se sigue que
$R={\mathcal O}$ es dominio Dedekind.
$\fin$ \end{proof}

\subsection[Estructura de m\'odulos sobre un
dominio Dedekind]{Estructura general de m\'odulos finitamente generados sobre
un dominio Dedekind}\label{DrinfeldS1.3}

\begin{teorema}\label{DrinfeldT1.3.2} Sean $R$ es un dominio Dedekind y $M_1,M_2$
dos $R$--m\'odulos libres de torsi\'on que se pueden escribir como
\begin{gather*}
M_1\cong I_1\oplus\cdots\oplus I_m; \quad M_2\cong J_1\oplus\cdots
\oplus J_n,
\intertext{donde $I_i, J_j$ son ideales fraccionarios de $R$. Entonces $M_1\cong
M_2 \iff m=n$ y existe $a\in \LL =\coc R$, $a\neq 0$, tal que}
I_1\cdots I_m\cong a J_1\cdots J_n.
\end{gather*}
\end{teorema}

\begin{proof} \cite[Theorem 1.39, p\'agina 30]{Nar2004}. 
$\fin$ \end{proof}

\begin{teorema}\label{DrinfeldT1.3.3} Si $R$ es un dominio Dedekind y si $M\neq 0$ es un
$R$--m\'odulo finitamente generado, entonces $M\cong \tor M\oplus P$ donde
$\tor M$ es el m\'odulo de torsi\'on de $M$ y $P$ es libre de torsi\'on.

M\'as a\'un, $\tor M\cong \oplus_{i=1}^m R/{\eu p}_i^{n_i}$ para algunos ideales
primos ${\eu p}_i\neq 0$, $n_i\in{\ma N}$. Las parejas $({\eu p}_i, n_i)$ est\'an
un\'ivocamente determinadas salvo permutaciones.
\end{teorema}

\begin{proof} \cite[Theorem 3.7, p\'agina 5]{Kor2010}. $\fin$ \end{proof}

\begin{teorema}\label{DrinfeldT1.3.4} Sea $R$ un dominio Dedekind y sea $M$ un
$R$--m\'odulo finitamente generado. Entonces las siguientes condiciones son
equivalentes.
\las
\item $M$ es libre de torsi\'on.
\item $M$ es plano.
\item $M$ es proyectivo.
\end{list}
\end{teorema}

\begin{proof} \cite[Theorem 3.6, p\'agina 4]{Kor2010}. $\fin$ 
\end{proof}

\begin{observacion}\label{DrinfeldO1.3.5} Se tiene que si $\LL =\coc R$, entonces
\[
\tor M=\ker \Big(M\to M\otimes_R \LL \Big), m\mapsto m\otimes 1.
\]
\end{observacion}

\begin{teorema}\label{DrinfeldT1.3.6} Sea $R$ un dominio Dedekind y sea $M$
un $R$--m\'odulo finitamente generado. Entonces $M$ es proyectivo si
y solamente si $M$ es libre de torsi\'on y esto \'ultimo se cumple si y
solamente si $M\cong R^{n-1}\oplus I$ donde $n$ es el rango de $R$, esto
es, $n=\dim_\LL  M\otimes_R \LL $ e $I$ es un ideal no cero de $R$.
\end{teorema}

\begin{proof} \cite[Theorem 7.2, p\'agina 12]{May}. $\fin$ \end{proof}

\begin{teorema}\label{DrinfeldT1.3.1}
Sea $M$ un $R$--m\'odulo finitamente generado donde $R$ es un dominio
Dedekind. Sea $\tor M=\{m\in M\mid \text{existe\ } r\in R, r\neq 0, rm=0\}$
el {\em subm\'odulo de torsi\'on\index{submodulo de torsion@subm\'odulo de torsi\'on}} de $M$. Entonces
\[
M\cong R^m\oplus I\oplus \tor M
\]
donde $m\in {\ma Z}$, $m\geq 0$, e $I$ es un ideal de $R$.
\end{teorema}

\begin{proof} \cite[Theorem 1.32, p\'agina 24]{Nar2004}. 
$\fin$ \end{proof}

En resumen, si $M$ es un $R$--m\'odulo finitamente generado, donde
$R$ es un dominio Dedekind, se tiene que
\[
M\cong R^{n-1}\oplus I \oplus \tor(M)
\]
donde $n=\dim_\LL  M\otimes_R \LL $, $\LL =\coc R$, $I$ es un ideal de $R$
y $\tor(M)=\{m\in M\mid \text{existe\ } r\in R, r\neq 0 \text{\ con\ }rm=0\}
=\bigoplus_{i=1}^m R/{\eu p}_i^{n_i}$ con ${\eu p}_1,\ldots, {\eu p}_m$
ideales primos no cero 
de $R$, $n_1,\ldots, n_m\in{\ma N}$ y 
el conjunto $\{(n_i,{\eu p}_i)\}_{i=1}^m$ es
\'unico salvo permutaci\'on.

\subsubsection{Torsi\'on de m\'odulos sobre dominios Dedekind}\label{DrinfeldS1.3.1}

Sean $R$ un dominio Dedekind y $M$ un $R$--m\'odulo. Sean $I$ un
ideal no cero de $R$ y sea 
\[
M[I]:=\{m\in M\mid am=0 \text{\ para toda\ } a\in I\}=\cap_{a\in I} M[a],\label{Drinfeldtorsion}
\]
donde $M[a]:=\{m\in M\mid am=0\}$. Se tiene que tanto $M[a]$ como $M
[I]$ son subm\'odulos de $M$. Notemos que $M[a]=M[(a)]$.

\begin{proposicion}\label{DrinfeldP1.3.7} Sean $I,J$ dos ideales no cero de $R$
tales que $I+J=R$, es decir, $I,J$ son primos relativos. Entonces
$M[IJ]=M[I]\oplus M[J]$.
\end{proposicion}

\begin{proof} Sean $a\in I$, $b\in J$ tales que $1=a+b$. Entonces para toda
$m\in M$ se tiene $m=1\cdot m=am+bm$. Sea $m\in M[IJ]$, entonces
$am\in M[J]$ y $bm\in M[I]$ pues para toda $x\in I$ y para toda $y\in J$
se tiene que $xb, ya\in IJ$ y $yam=0$, $xbm=0$ con la hip\'otesis de que
$m\in M[IJ]$. 

Ahora bien,  $M[I]$ y $M[J]$ son subm\'odulos de $M[IJ]$,
ya que si $s\in M[I]$, para cualesquiera $\alpha, \beta\in IJ$, $\alpha\beta
s=\beta\alpha s=\beta 0=0$. Por tanto se tiene que $M[IJ]=M[I]+M[J]$.

Finalmente, si $m\in M[I]\cap M[J]$, $am=0$, $bm=0$ por lo que
$m=am+bm=0+0=0$. Se sigue que $M[I]\cap M[J]=\{0\}$ y
\fbox{$M[IJ]=M[I]\oplus M[J]$}.
$\fin$ \end{proof}

\begin{corolario}\label{DrinfeldC1.3.8} Sea $I$ un ideal no cero de $R$, $I=
P_1^{e_1}\cdots P_r^{e_r}$ con $P_i$ ideales primos no cero distintos.
Entonces $M[I]=\oplus_{i=1}^r M[P_i^{e_i}]$. $\fin$
\end{corolario}

\begin{definicion}\label{DrinfeldD1.3.9} Sea $P$ un ideal primo de $R$, $P\neq 0$.
Se define  la {\em componente $P$--primaria
$M(P)$ de $M$\index{componente primaria}} por
\[
M(P):=\bigcup_{e=1}^{\infty} M[P^e]=\{m\in M\mid \text{existe $e$ tal que
$xm=0$ para toda\ } x\in P^e\}.
\]
\end{definicion}

\begin{proposicion}\label{DrinfeldP1.3.10} Si $M$ es de torsi\'on, entonces
$M=\bigoplus\limits_{\substack{P\text{\ primo}\\ P\neq 0}}M(P)$.
\end{proposicion}

\begin{proof} Se tiene para todo ideal primo $P\neq 0$ 
que $M(P)\subseteq \tor M=M$.
Ahora sea $m\in M$. Existe $a\in R$ con $a\neq 0$, $am=0$. Sea $(a)=
P_1^{e_1}\cdots P_r^{e_r}$. Entonces $m\in M[a]=M[(a)]=\bigoplus_{i=1}^r
M[P_i^{e_i}]\subseteq \sum_{i=1}^r M(P_i)\subseteq \sum_{\substack{P
\text{\ primo}\\ P\neq 0}} M(P)$.

Si $0=\sum_{i=1}^t m_i$ con $m_i\in M(P_i)$, $P_1,\ldots, P_t$ distintos,
entonces $m_i\in M[P_i^{c_i}]$ para algunos $c_i$. Por tanto
$0=\sum_{i=1}^t m_i\in \bigoplus_{i=1}^t M[P_i^{c_i}]$. Se sigue que $m_i=0$,
$1\leq i\leq t$, por lo que $\sum_{\substack{P\text{\ primo}\\ P\neq 0}} M(P)=
\bigoplus_{\substack{P\text{\ primo}\\ P\neq 0}} M(P)=M$.
$\fin$ \end{proof}

\begin{proposicion}\label{DrinfeldP1.3.11} Sea $0\longrightarrow M_1
\stackrel{f}{\longrightarrow} M_2 \stackrel{g}{\longrightarrow} M_3
\longrightarrow 0$ una sucesi\'on exacta de $R$--m\'odulos de 
torsi\'on y $R$ dominio Dedekind. Sea $P\neq 0$ un ideal primo
de $R$. Entonces $0\longrightarrow M_1(P)
\stackrel{f_1}{\longrightarrow} M_2(P) \stackrel{g_1}{\longrightarrow} M_3(P)
\longrightarrow 0$ es exacta, donde $f_1=f|_{M_1(P)}$, $g_1=g|_{M_2(P)}$.
\end{proposicion}

\begin{proof} Sea $m\in M_1(P)$, entonces $m\in M_1[P^c]$ para alg\'un $c$. Se
sigue que $xm=0$ para toda $x\in P^c$. De esta forma $xf_1(m)=f_1(xm)=f(0)=0$.
Por lo tanto $f_1(m)\in M_2[P^c]\subseteq M_2(P)$. En particular obtenemos
que $f_1(M_1(P))\subseteq M_2(P)$. 

Sea $n\in M_2(P)$. Entonces $n\in M_2[P^d]$
para alg\'un $d$. Si $x\in P^d$, $x g_1(n)=g_1(xn)=g_1(0)=0$. 
Por tanto $g_1(n)\in M_3[P^d]\subseteq
M_3(P)$ y $g_1(M_2(P))\subseteq M_3(P)$.

Ahora bien, $f|_{M_1(P)}$ es inyectiva y $g\circ f|_{M_1(P)}=0|_{M_1(P)}=0$
por lo que $f_1$ es uno a uno e 
$\im f_1\subseteq  \ker g_1$. Si $x\in \ker g_1$, $x\in \ker g=\im f$
por lo que existe $y\in M_1$ tal que $f(y)=x$. Adem\'as, $x\in \ker g_1\subseteq
M_2(P)$, por lo que $x\in M_2[P^d]$ para alg\'un $d$. Por tanto $\alpha x=0$
para toda $\alpha\in P^d$. De esta forma obtenemos que $f(\alpha y)=
\alpha f(y)=\alpha x=0$. Se sigue que $\alpha y=0$ puesto que la funci\'on $f$
es inyectiva. As{\'\i} $y\in M_1[P^d]$. Finalmente obtenemos que $\im f_1=
\ker g_1$.

Falta ver que $g_1\colon M_2(P)\to M_3(P)$ es suprayectiva. Sea $n\in M_3(P)$.
Por tanto $n\in M_3[P^c]$ alg\'un $c$. Sea $m\in M_2$ tal que $g(m)=n$. Sea
$m=\sum_{i=1}^t m_i$, $m_i\in M_2[P_i^{e_i}]$ para algunos ideales primos
$P_i\neq 0$ y algunos $e_i$. De esta forma se tiene $g(m)=\sum_{i=1}^t g(m_i)=n\in
M_3[P^c]$. Puesto que los distintos subm\'odulos $M_3(P)$ forman suma directa,
se tiene que $g(m_i)=0$ para $P_i\neq P$. Por tanto $P=P_{i_0}$ para alg\'un
$i_0$ y $g(m_{i_0})=n$. Se sigue que $g_1\colon M_2(P)\to M_3(P)$ es
suprayectiva.
$\fin$ \end{proof}

\begin{teorema}\label{DrinfeldT1.3.12} Sea $M$ un $R$--m\'odulo arbitrario y sea $P$
un ideal primo no cero de $R$ y donde $R$ es un dominio Dedekind. Sea $\pi
\in P\setminus P^2$, esto es, $v_P(\pi)=1$. Entonces $M[P^e]=M[\pi^e](P)$.
\end{teorema}

\begin{proof} Se tiene $(\pi^e)=P^e I$ con $I, P$ primos relativos puesto que $v_P(\pi^e)
=e$. Entonces $M[\pi^e]=M[P^e I]=M[P^e]\oplus M[I]$. Finalmente, $M[P^e](P)=
M[P^e]$ y $M[I](P)=0$, por lo que $M[\pi^e](P)=M[P^e](P)\oplus M[I](P)=
M[P^e]\oplus 0=M[P^e]$.
$\fin$ \end{proof}

\begin{teorema}\label{DrinfeldT1.3.13} Sea $R$ un dominio Dedekind y 
sea $M$ un $R$--m\'odulo divisible, es decir,
si $x\in M$ y $r\neq 0$ con $r\in R$, entonces existe $y\in M$ tal que $ry=x$
(esto es algo as{\'\i} como $y=\frac{x}{r}$).
Sea $P$ un ideal primo no cero de $R$. Entonces para $e\in {\ma N}$, $e>1$,
la sucesi\'on
\[
0\longrightarrow M[P]\stackrel{i}{\longrightarrow} M[P^e]\stackrel{\varphi}{
\longrightarrow} M[P^{e-1}]\longrightarrow 0
\]
es exacta, donde $i$ es la inclusi\'on y donde $\varphi(x)=\pi x$ con $\pi\in
P\setminus P^2$, esto es, $v_P(\pi)=1$.

En caso de que $M$ no sea divisible, entonces 
\[
0\longrightarrow M[P]
\stackrel{i}{\longrightarrow} M[P^e]\stackrel{\varphi}{
\longrightarrow} M[P^{e-1}]
\]
es exacta, es decir, la \'ultima suprayectividad
puede fallar.
\end{teorema}

\begin{proof} Primero veamos que $0\longrightarrow M[\pi]\stackrel{\tilde{i}}{\longrightarrow} 
M[\pi^e]\stackrel{\tilde{\varphi}}{\longrightarrow} M[\pi^{e-1}]\longrightarrow 0$,
es exacta en donde $\tilde{\varphi}(x)=\pi x$.

Si $x\in M[\pi^e]$, $\pi^{e-1}(\pi x)=\pi^e x=0$. Por lo tanto $\pi x=\tilde{\varphi}(x)
\in M[\pi^{e-1}]$. Ahora $\ker \tilde{\varphi}=\{x\in M[\pi^e]\mid \pi x=\tilde{\varphi}(x)
=0\}=M[\pi]$. Por tanto $\ker \tilde{\varphi}=\im \tilde{i}$.

Falta ver que $\tilde{\varphi}$ es suprayectiva. Sea $y\in M[\pi^{e-1}]\subseteq M$.
Puesto que $M$ es divisible y $\pi\neq 0$, existe $x\in M$ tal que $\pi x=y=\varphi(x)$.
Se tiene $\pi^e x=\pi^{e-1}(\pi x)=\pi^{e-1} y=0$ por lo que $x\in M[\pi^e]$. Se sigue
que 
\begin{gather*}
0\longrightarrow M[\pi]\stackrel{\tilde{i}}{\longrightarrow} 
M[\pi^e]\stackrel{\tilde{\varphi}}{\longrightarrow} M[\pi^{e-1}]\longrightarrow 0,
\intertext{es exacta. Por tanto}
0\longrightarrow M[\pi](P)\stackrel{\tilde{i}}{\longrightarrow} 
M[\pi^e](P)\stackrel{\tilde{\varphi}}{\longrightarrow} M[\pi^{e-1}](P)\longrightarrow 0,
\intertext{es exacta. Finalmente $M[\pi](P)=M[P]$, $M[\pi^e](P)=M[P^e]$ y $M[\pi^{e-1}](P)=
M[P^{e-1}]$ por lo que}
0\longrightarrow M[P]\stackrel{i}
{\longrightarrow} M[P^e]\stackrel{\varphi}{
\longrightarrow} M[P^{e-1}]\longrightarrow 0
\end{gather*}
es exacta.
$\fin$ \end{proof}

\begin{corolario}\label{DrinfeldC1.3.14} Para $e\geq 1$, $|M[P^e]|=|M[P]|^e$.
\end{corolario}

\begin{proof} Lo hacemos por inducci\'on en $e$. Se tiene del Teorema \ref{DrinfeldT1.3.13} que
$|M[P^e]| =|M[P^{e-1}]||M[P]|
\underbracket[0pt]{=}_{\substack{\uparrow\\
\text{inducci\'on}\\ \text{en $e$}}} |M[P]|^{e-1}|M[P]|=|M[P]|^e$.
$\fin$ \end{proof}

Sea ahora $\K $ un campo de funciones global, es decir, con campo de constantes
un campo finito $\F$, $q=p^u$. Sea $\p$ un lugar fijo de $\K $ de grado $d_{\infty}
=[\K (\p):\F] \geq 1$. Sea $A=\{x\in \K \mid v_{\eu p}(x)\geq 0 
\text{\ para toda ${\eu p}\neq \p$}\}$. Sea $F$ cualquier campo tal que $\F\subseteq
F$.

\begin{definicion}\label{DrinfeldD1.3.15}
El campo $F$ se llama {\em $A$--campo\index{A@$A$--campo}} si existe $\delta\colon
A\lra F$ un homomorfismo de anillos. Sea ${\eu q}=\ker \delta$. Entonces ${\eu q}$
es un ideal primo de $A$.
\end{definicion}

\begin{definicion}\label{DrinfeldD1.3.16} Un {\em m\'odulo de Drinfeld\index{m\'odulo de
Drinfeld} sobre $F$}, donde $F$ es un $A$--campo, 
es un homomorfismo de anillos $\rho\colon A\to \aditivo F$
con $\tau u=u^q$, tal que $D\circ \rho=\delta$ y $\rho(a)\neq \delta(a) \tau^0$
para alguna $a\in A$ y donde $D\big(\sum_{i=1}^r\alpha_i \tau^i\big)=\alpha_0$.
Equivalentemente, existe $a\in A$ tal que $\rho(a)\notin F$.

En otras palabras, $\rho$ es un m\'odulo de Drinfeld si $\rho$ es un homomorfismo
de anillos tal que para toda $a\in A$ se tiene que $\rho(a)=\delta(a)\tau^0+
\sum_{i=1}^r\alpha_i \tau^i$ y $\rho(a)\neq \delta(a)\tau^0$ para alguna $a\in A$.
Se denota tambi\'en $\rho(a)=\rho_a$ y $\rho$ denotar\'a siempre un $A$--m\'odulo
de Drinfeld.
\end{definicion}

\begin{observacion}\label{DrinfeldO1.3.17} No estamos garantizando que para
cualesquiera $A$ y $F$ como antes existan m\'odulos de Drinfeld.
De hecho, no siempre existen.
\end{observacion}

\begin{ejemplo}\label{DrinfeldEj1.3.18}
Sea $A=R_T=\F[T]$ y sea $F$ cualquier campo de contenga a $A$ y por
tanto a $K=\F(T)=\coc A\subseteq F$. Sea $\delta\colon A\lra F$ cualquier $\F$--homomorfismo.
Sean $r\in{\ma N}$ y $\rho_T:=\delta(T)+\sum_{i=1}^r \alpha_i \tau^i$ arbitrario
con $\alpha_r\neq 0$. Entonces $\rho$ se puede extender de manera
\'unica a un homomorfismo $\rho\colon A\to \aditivo F$ por $\rho_{M(T)}:=
M(\rho_T)$, es decir, $\rho(M(T))= M(\rho(T))$ para $M\in R_T$. En particular
para $A=R_T$ los m\'odulos de Drinfeld existen en abundancia.

La raz\'on de lo anterior se debe a que $A$ es una $\F$--\'algebra libre.
\end{ejemplo}

Un caso particular del Ejemplo \ref{DrinfeldEj1.3.18} es el m\'odulo de Carlitz
$C\colon A\to \aditivo {\bar{K}}$ donde $\bar{K}$ es la cerradura algebraica
de $K$ y $C_T=T+\tau$.

\begin{definicion}\label{DrinfeldD1.3.19}
Se denota $\Drin_A(F)$\label{DrinfeldDrinf} al conjunto de los $A$--m\'odulos de Drinfeld sobre
$F$ una vez que el mapeo $\delta\colon A\lra F$ ha sido dado.
\end{definicion}

Casi siempre $\delta$ ser\'a la inclusi\'on natural o un mapeo de reducci\'on
m\'odulo un ideal primo no cero de $A$.

Dada $B$ cualquier $F$--\'algebra, $A$ act\'ua sobre $B$ v{\'\i}a $\delta$
si definimos 
\begin{gather*}
a\circ v:=\delta(a) v
\intertext{para $v\in B$ y $a\in A$. Es decir,
$B$ se hace un $A$--m\'odulo. Por otro lado, si consideramos $\rho$,
$B$ se hace tambi\'en $A$--m\'odulo con la acci\'on}
a\ast v:=\rho_a(v),
\intertext{esto es, si $\rho_a=\sum_{i=0}^r\alpha_i \tau^i$, entonces}
a\ast v=\rho_a(v)=\sum_{i=0}^r \alpha_i v^{q^i}=\delta(a)v+\sum_{i=1}^r\alpha_i
v^{q^i}.
\end{gather*}

Por definici\'on, existe $a$ tal que $\rho_a\neq \delta(a)$ por lo que 
$a\ast v\neq a\circ v$. A veces se escribe $B_{\rho}$ para denotar la
estructura de $A$--m\'odulo de $B$ bajo la acci\'on de $\rho$.

\begin{definicion}\label{DrinfeldD1.3.20}
Al ideal ${\eu q}=\ker \delta$, $\delta\colon A\to F$, se le llama
{\em caracter{\'\i}stica\index{caracter{\'\i}stica de un m\'odulo de
Drinfeld} de $\rho$}. Si ${\eu q}=(0)$ se dice que $\rho$ tiene
{\em caracter{\'\i}stica gen\'erica\index{caracter{\'\i}stica gen\'erica}}
o {\em caracter{\'\i}stica infinita\index{caracter{\'\i}stica infinita}}.
Si ${\eu q}\neq 0$ se dice que $\rho$ tiene {\em caracter{\'\i}stica
finita\index{caracter{\'\i}stica finita}}. 

Se denota ${\eu q}=\car (\rho)$.
\end{definicion}

\begin{proposicion}\label{DrinfeldP1.3.21} Cualquier m\'odulo de Drinfeld
$\rho\colon A\to \aditivo F$ es un mapeo inyectivo.
\end{proposicion}

\begin{proof} Sea ${\eu p}=\ker \rho$ un ideal primo de $A$. Si ${\eu p}
\neq 0$, entonces $A/{\eu p}$ es un campo finito pues ${\eu p}$
es maximal y $\tilde{\rho}\colon A/{\eu p}\to \aditivo F$ es un
homomorfismo de anillos.

Sea $a\in A$ tal que $\rho_a\notin F$, por tanto $\tilde{\rho}(\bar{a})
\notin F$. En particular $\deg_{\tau}\tilde{\rho}(\bar{a})=r\geq 1$ y
$\deg_{\tau}\tilde{\rho}(\bar{a}^n)=\deg_{\tau}\tilde{\rho}(
\bar{a})^n=nr\geq 1$, lo que contradice que
existe $n\in{\ma N}$ con $\bar{a}^n=1$ y por ende que $\tilde{\rho}
(\bar{a}^n)=\tilde{\rho}(1)=1$ que es de grado $0$. Se sigue
que ${\eu p}=0$ y que $\rho$ es inyectivo.
$\fin$ \end{proof}

La Proposici\'on \ref{DrinfeldP1.3.21} demuestra que los $A$--m\'odulos de
Drinfeld son encajes no triviales de $A$ en $\aditivo F$.

Ahora bien, >por qu\'e \'unicamente se toma un \'unico primo al infinito
y no m\'as de uno?

Sea $S=\{{\eu p}_1,\ldots, {\eu p}_s\}$ primos de $\K $
con $s\geq 2$. Sea ${\eu A}=\frac{{\eu p}_1^{\deg {\eu p}_2}}
{{\eu p}_2^{\deg {\eu p}_1}}$ que es de grado $0$. Puesto que el grupo
de clases $I_{\K ,0}$ es finito, existe $n\in{\ma N}$ tal que
${\eu A}^n=(x)_\K =\frac{{\eu p}_1^a}{{\eu p}_2^b}$ es principal. Entonces
$x, x^{-1}\in {\mathcal O}_S=\cap_{{\eu p}\notin S}{\mathcal O}_{\eu p}$.
Esto es, $x$ es una unidad de ${\mathcal O}_S$, $x\in {\mathcal O}_S^{\ast}$.

Ahora si $\rho\colon {\mathcal O}_S\to \aditivo F$ es un m\'odulo de
Drinfeld, entonces $\rho(x)\in \aditivo F^{\ast}=F^{\ast}$ y $\rho(x)=
\delta(x)\tau^0\in F^{\ast}$ pues 
$\rho({\mathcal O}_S^{\ast})\subseteq F^{\ast}$.

Sea $\F(x)\subseteq \K $ y $x\in \K \setminus \F$. Por tanto $[\K :\F(x)]<\infty$.
Por tanto, si $y\in {\mathcal O}_S\subseteq \K $, existen $n\in{\ma N}$ y
$f_0(x),\ldots, f_{n-1}(x)\in \F(x)$ tales que 
\[
y^n+f_{n-1}(x)y^{n-1}+\cdots+ f_1(x)y+f_0(x)=0.
\]

Sea $y\in {\mathcal O}_S$ tal que
$\rho_y=\rho(y)=\sum_{i=0}^r\alpha_i\tau^i$ con $r\geq 1$, 
$\alpha_r\neq 0$ y $\deg_{\tau}\rho(y)^n=nr$. 
Tal $y$ existe si $\rho$ es un m\'odulo de Drinfeld.
Se tiene que si $f_i(x)\neq 0$, $\rho_{f_i(x)}\in F$ puesto que $\rho_x\in F$.
Por tanto $\deg_{\tau}\rho_{f_i(x)y^i}=ir$ para $f_i(x)\neq 0$.
Se sigue que
\[
\deg_{\tau}\rho_{(y^n+f_{n-1}(x)y^{n-1}+\cdots+f_1(x)y+f_0(x))}=nr
\neq \deg_{\tau}\rho_0=\deg_{\tau}0,
\]
lo cual es un absurdo. Esto implica que $\deg_{\tau}\rho_y=0$ y $\rho(
{\mathcal O}_S)\subseteq F$ lo que contradice la definici\'on de m\'odulo
de Drinfeld. Esta es la raz\'on por la cual \'unicamente se toma un
\'unico primo al infinito.

Sea $v_{\infty}$ la valuaci\'on asociada a $\p$. Sea $d_{\infty}=\deg_\K \p$.

\begin{definicion}\label{DrinfeldD1.3.22} Para $x\in \K $ se define $\deg(x):=
-d_{\infty}v_{\infty}(x)$\index{grado de un elemento} 
y $|x|_{\infty}:=q^{\deg (x)}$\index{valor absoluto de un elemento}.
\end{definicion}

>Que significa $\deg(x)$? 

Se tiene que si $x\in A$, $x\neq 0$, $(x)_\K =
\frac{\pK_1^{\alpha_1}\cdots\pK_s^{\alpha_s}}{\p^t}$ con el grado de $(x)_\K $
igual a $0$ y donde $t=-v_{\infty}(x)$ con $\alpha_1,\ldots, \alpha_s\geq 0$.

Ahora se tiene que para todo $\pK$ lugar de $\K $, $d_\K (\pK)
:=[\o_{\pK}/\pK:\F]$ y
$|\o_{\pK}/\pK|=\N\pK=q^{d_\K (\pK)}$.

Para $\pK\neq \p$, sea $\iota\colon
A\to \o_{\pK}$ el encaje natural y $\pK\cap A=\pK^{\prime}$ satisface que
$\iota(\pK^{\prime})\subseteq \pK$, por lo que $\iota$ induce una inyecci\'on
$\tilde{\iota}\colon A/\pK^{\prime}\hookrightarrow \o_{\pK}/\pK$, $\tilde{\iota}(
a\bmod \pK^{\prime})=a\bmod \pK$.

Sea $A_{\pK^{\prime}}=\big\{\frac{\alpha}{\beta}\mid \alpha,\beta\in A, v_{
\pK^{\prime}}(\beta)=0\big\}$. El ideal m\'aximo  del anillo localizado $A_{
\pK^{\prime}}$ es $\pK^{\prime}A_{\pK^{\prime}}=\big\{\frac{\alpha}{\beta}
\mid \alpha,\beta\in A, v_{\pK^{\prime}}(\alpha)>0,
v_{\pK^{\prime}}(\beta)=0\big\}$.

Se tiene que el encaje $\theta\colon A\hookrightarrow A_{\pK^{\prime}}$,
$\theta(a)=\frac{a}{1}$, 
satisface $\theta(\pK^{\prime})\subseteq \pK^{\prime} A_{\pK^{\prime}}$,
por lo que $\tilde{\theta}:A/\pK^{\prime}\hookrightarrow A_{\pK^{\prime}}/
\pK^{\prime}A_{\pK^{\prime}}$ es una inyecci\'on.

Para $\frac{\alpha}{\beta}\in A_{\pK^{\prime}}$ queremos hallar $\gamma
\in A$ tal que $\frac{\gamma}{1}-\frac{\alpha}{\beta}\in\pK^{\prime}A_{\pK^{\prime}}$,
esto es, $\frac{\gamma}{1}-\frac{\alpha}{\beta}=\frac{\gamma\beta-\alpha}
{\beta}\in \pK^{\prime}A_{\pK^{\prime}}$, es decir $\gamma \beta-\alpha
\in \pK^{\prime}$ o, equivalentemente, $v_{\pK^{\prime}}(\gamma\beta
-\alpha)\geq 1$ pues en este caso se tiene $\tilde{\theta}
(\gamma)=\frac{\gamma}{1}\equiv \frac{\alpha}{\beta}\bmod
\pK^{\prime} A_{\pK^{\prime}}$. 
Ahora bien, $\bar{\beta}\in A/\pK^{\prime}$, $\bar{\beta}
\neq 0$.

Sea $\gamma\equiv \alpha\beta^{-1}\bmod \pK^{\prime}$ con $\gamma
\in A$ el cual existe pues $\overline{\alpha\beta^{-1}}=\bar{\alpha}
\overline{\beta^{-1}}$ existe. Por tanto $\gamma \beta-\alpha\in
\pK^{\prime}$. Se sigue que $\tilde{\theta}$ es suprayectiva y
$A/\pK^{\prime}\cong A_{\pK^{\prime}}/\pK^{\prime}A_{\pK^{\prime}}$.

Finalmente veamos que $A_{\pK^{\prime}}/\pK^{\prime}A_{\pK^{\prime}}
\stackrel{\tilde{\varphi}}{\hookrightarrow} \o_{\pK}/\pK$ es suprayectiva.

Sea $\xi\in \o_{\pK}$. Si $\xi\in \pK$, $\tilde{\varphi}(0)=
\bar{\xi}=0$ en $\o_{\pK}/\pK$. Supongamos $\xi\notin\pK$. Se quiere
un elemento $\frac{\alpha}{\beta}\in A_{\pK^{\prime}}$, $v_{\pK^{\prime}}
(\beta)=0$, tal que $\frac{\alpha}{\beta}-\xi\in\pK$. Basta que
$\alpha-\xi\beta\in\pK$ con $\beta$ unidad en $\o_{\pK}$ puesto que
$v_{\pK}(\beta)=0=v_{\pK}(\beta^{-1})$. Se tiene $\coc A=\K $ y
$\xi\in\o_{\pK}\subseteq \K $. Por tanto $\xi=\frac{a}{b}$ con $a,b\in A$.
Puesto que $v_{\pK}(\xi)\geq 0$, $v_{\pK}(b)=0$ y $a\xi-b=0$.
Se sigue que $\tilde{\varphi}\big(\overline{\frac{a}{b}}\big)=\bar{\xi}$.
Por tanto $\tilde{\varphi}$ es un isomorfismo y
\[
A/\pK^{\prime}\cong A_{\pK^{\prime}}/\pK^{\prime}A_{\pK^{\prime}}
\cong \o_{\pK}/\pK.
\]

En particular $|A/\pK^{\prime}|=|\o_{\pK}/\pK|=q^{d_\K (\pK)}$.

En general $A/(\pK^{\prime})^m\cong A_{\pK^{\prime}}/(\pK^{\prime}
A_{\pK^{\prime}})^m$ y $\pK^{\prime}A_{\pK^{\prime}}$ es un ideal principal
pues $A_{\pK^{\prime}}$ es un dominio de valuaci\'on discreta.
Por tanto, si $\pK^{\prime}A_{\pK^{\prime}}=(\pi)$, entonces
$A_{\pK^{\prime}}\stackrel{\mu}{\longrightarrow}\big(\pK^{\prime}A_{\pK^{\prime}}
\big)^{m-1}$, $a\mapsto a\pi^{m-1}$, es un epimorfismo de
$A_{\pK^{\prime}}$--m\'odulos. Por tanto 
\[
\begin{array}{ccccc}
A_{\pK^{\prime}}&\stackrel{\mu}{\longrightarrow}&\big(\pK^{\prime}
A_{\pK^{\prime}}\big)^{m-1}&\stackrel{p}{\longrightarrow}&\big(
A_{\pK^{\prime}}\pK^{\prime}\big)^{m-1}/\big(\pK^{\prime}
A_{\pK^{\prime}}\big)^m\\
a&\longrightarrow&a\pi^{m-1}&\longrightarrow&a\pi^{m-1}\bmod
\big(\pK^{\prime}A_{\pK^{\prime}}\big)^m
\end{array}
\]
es suprayectiva y $\ker p\circ \mu=(\pi)=\pK^{\prime}A_{\pK^{\prime}}$.
Se sigue que $\frac{A_{\pK^{\prime}}}{\pK^{\prime}A_{\pK^{\prime}}}\cong
\frac{(\pK^{\prime}A_{\pK^{\prime}})^{m-1}}{(
\pK^{\prime}A_{\pK^{\prime}})^m}$ como $A_{\pK^{\prime}}$--m\'odulos.
De esta forma se obtiene la sucesi\'on exacta
\begin{gather*}
0\longrightarrow\frac{A_{\pK^{\prime}}}{\pK^{\prime}A_{\pK^{\prime}}}\cong
\frac{\big(\pK^{\prime}A_{\pK^{\prime}}\big)^{m-1}}{\big(
\pK^{\prime}A_{\pK^{\prime}}\big)^m}
\stackrel{i}{\longrightarrow}\frac{A_{\pK^{\prime}}}{\big(\pK^{\prime}
A_{\pK^{\prime}}\big)^m}\\
\longrightarrow\frac{\frac{
A_{\pK^{\prime}}}{\big(\pK^{\prime}A_{\pK^{\prime}}\big)^m}}
{\frac{(\pK^{\prime}A_{\pK^{\prime}})^{m-1}}{\big(\pK^{\prime}
A_{\pK^{\prime}}\big)^m}}\cong \frac{
A_{\pK^{\prime}}}{\big(\pK^{\prime}A_{\pK^{\prime}}\big)^{m-1}}
\longrightarrow 0.
\end{gather*}

En particular $\Big|\frac{A_{\pK^{\prime}}}{\big(\pK^{\prime}
A_{\pK^{\prime}}\big)^m}\Big|=\Big|\frac{A_{\pK^{\prime}}}{\big(\pK^{\prime}
A_{\pK^{\prime}}\big)^{m-1}}\Big|\cdot \big|\frac{A_{\pK^{\prime}}}{\pK^{\prime}
A_{\pK^{\prime}}}\big|\underbracket[0pt]{=}_{\substack{\uparrow\\ \text{inducci\'on}}}
\big|\frac{A_{\pK^{\prime}}}{\pK^{\prime}
A_{\pK^{\prime}}}\big|^m$.

De esta forma obtenemos que 
\[
\Big|\frac{A}{\pK^m}\Big|=\Big|\frac{A_{\pK^{\prime}}}{\big(\pK^{\prime}
A_{\pK^{\prime}}\big)^m}\Big|= \Big|\frac{A_{\pK^{\prime}}}{\pK^{\prime}
A_{\pK^{\prime}}}\Big|^m=\Big|\frac{\o_{\pK}}{\pK}\Big|^m=
\Big|\frac{\o_{\pK}}{\pK^m}\Big|=q^{d_\K (\pK^m)}.
\]

Finalmente, si $x\in A$, $(x)_\K =\frac{\pK_1^{\alpha_1}\cdots \pK_s^{
\alpha_s}}{\p^t}$, $xA=(\pK^{\prime}_1)^{\alpha_1}\cdots 
(\pK^{\prime}_s)^{\alpha_s}$. Por el Teorema Chino del Residuo obtenemos
\begin{align*}
A/xA&\cong \big(A/(\pK^{\prime}_1)^{\alpha_1}\big)\times\cdots\times
\big(A/(\pK^{\prime}_s)^{\alpha_s}\big)\qquad\qquad\text{y}\\
\big|A/xA\big|&=\big|A/(\pK^{\prime}_1)^{\alpha_1}\big|\cdots
\big|A/(\pK^{\prime}_s)^{\alpha_s}\big|\\
&=q^{\sum_{i=1}^s \alpha_i d_\K (\pK_i)}=
q^{t d_\K (\p)}=q^{-v_{\infty}(x)d_{\infty}}=q^{\deg (x)}.
\end{align*}

En resumen, hemos obtenido

\begin{proposicion}\label{DrinfeldP1.3.23}
Si $x\in A$, $x\neq 0$, entonces si $\deg x=-v_{\infty}(x)d_{\infty}$,
se tiene $q^{\deg(x)}=\big|A/xA\big|$. $\fin$
\end{proposicion}

\begin{definicion}\label{DrinfeldD1.3.24}
Sea $\rho\colon A\to \aditivo F$ un m\'odulo de Drinfeld. Sea
$\phi\colon A\to {\ma Z}$ definido por $\phi(a)=-\deg_{\tau}
\rho_a=-\deg \rho_a$ para $a\neq 0$. 
\end{definicion}

Notemos que $\deg_u \rho_a(u)=q^{\deg_{\tau}\rho_a}$.

Entonces $\phi$ es una valuaci\'on de $\K $ pues:
\begin{align*}
\phi(ab)&=-\deg_{\tau}\rho_{ab}=-\deg_{\tau}\rho_a\rho_b=-\deg_{\tau}
\rho_a-\deg_{\tau}\rho_b=\phi(a)+\phi(b);\\
\phi(a+b)&=-\deg_{\tau}\rho_{a+b}=-\deg_{\tau}(\rho_a+\rho_b)\\
&\geq
\min\{-\deg_{\tau}\rho_a,-\deg_{\tau}\rho_b\}=\min\{\phi(a),\phi(b)\}.
\end{align*}

Por tanto $\phi$ es una valuaci\'on en $A$ que se extiende a $\K =\coc
A$: $\phi\big(\frac{a}{b}\big):=\phi(a)-\phi(b)=-\deg_{\tau}\rho_a+
\deg_{\tau}\rho_b$. Puesto que existe $a\in A$ tal que $\deg_{\tau}
\rho_a\geq 1$, entonces $\phi(a)\leq -1$ y, por otro lado, para todo lugar
$\pK\neq \p$, $v_{\pK}(a)\geq 0$, se sigue que $\phi$ es equivalente
a $v_{\infty}$. Por tanto existe $r_{\rho}\in {\ma Q}^+$ tal que
$\phi(a)=r_{\rho}d_{\infty}v_{\infty}(a)$, es decir, $\phi(a)=-\deg_{\tau}
\rho_a=-r_{\rho}\deg (a)$, de donde
\centerline{\fbox{$\deg\rho_a=r_{\rho}\deg(a)$}.}

\begin{definicion}\label{DrinfeldD1.3.25}
El racional $r_{\rho}$ se llama el {\em rango de $\rho$\index{rango
de un m\'odulo de Drinfeld}}.
\end{definicion}

Veremos que $r_{\rho}\in {\ma N}$.

\begin{ejemplo}\label{DrinfeldEj1.3.26} Sea $A=R_T$ y $C\colon A\to
\aditivo {\bar{K}}$ el m\'odulo de Carlitz, esto es, $C_T=T+\tau$.
Se tiene $d_{\infty}=1$, $v_{\infty}(T)=-1$ y $\deg C_T=1=-r_C
\cdot 1\cdot (-1)=r_C$, esto es $r_C=1$ y el rango del m\'odulo
de Carlitz es $1$.
\end{ejemplo}

\begin{definicion}\label{DrinfeldD1.3.27} Dado un m\'odulo de Drinfeld
$\rho\colon A\to \aditivo F$, se define la {\em altura\index{altura
de un m\'odulo de Drinfeld} de $\rho$}, $h_{\rho}$ como sigue:
\las
\item Si $\car(\rho)=0$, entonces $h_{\rho}=0$,
\item Si $\car(\rho)=\pL\neq (0)$, sea $v_{\pL}$ es la valuaci\'on
asociada a $\pL$. Para $a\in A$, $a\neq 0$, sea $\rho_a=
\sum_{i=0}^n\alpha_i\tau^i$, $\alpha_0=\delta(a)$. Sea $i_0$ tal que
$\alpha_{i_0}\neq 0$ y $\alpha_j=0$ para $j=0,1,\ldots, i_0-1$. Se define
$j_{\rho}(a)=\ord(\rho_a)=i_0$\index{orden de un elemento}.
En particular 
\[
\ord(\rho_a)>0\iff \delta(a)=0 \iff a\in\pL. 
\]
Adem\'as $j_{\rho}$ define
una valuaci\'on no trivial en $A$ la cual es equivalente a $v_{\pL}$
pues $\ord(\rho_a)>0\iff v_{\pL}(a)>0$. Se sigue que existe $h_{\rho}
\in {\ma Q}^+$ tal que si $h_{\rho}$ es la altura, entonces
\[
j_{\rho}(a)=\ord(\rho_a)=h_{\rho}v_{\pL}(a)\deg_\K \pL.
\]
\end{list}
\end{definicion}

Se demostrar\'a que $h_{\rho}\in{\ma N}\cup \{0\}$.

\begin{ejemplo}\label{DrinfeldEj1.3.27'} Si $C$ es el m\'odulo de Carlitz,
$\delta$ es inyectiva y por tanto $h_C=0$.
\end{ejemplo}

\begin{ejemplo}\label{DrinfeldEj1.3.28} Sean $A=R_T$ y $F$ cualquier campo
que contiene a $\F$ y $\rho\colon A\to \aditivo F$ un m\'odulo de Drinfeld
de rango $r$ y altura $h$. Sea $\rho_T=\alpha_0+\sum_{i=1}^r
\alpha_i \tau^i$, $\alpha_1,\ldots,\alpha_r\in F$, $\alpha_r\neq 0$ pues
$\deg\rho_T=-d_{\infty}r_{\rho}v_{\infty}(T)=-(1)r_{\rho}(-1)=r_{\rho}=r$.

Adem\'as $\delta(T)=\alpha_0$, por lo tanto $\delta(f(T))=f(\delta(T))=
f(\alpha_0)$ para $f(T)\in R_T$. Por tanto $\rho_{f(T)}=f(\alpha_0)+
\cdots+ \beta\tau^{r\deg_T f}$ y
\[
\car(\rho)=\begin{cases}
(0)&\text{si $\alpha_0$ es trascendente sobre $\F$}\\
\big(\Irr(\alpha_0,T,\F)\big)&\text{si $\alpha_0$ es algebraico sobre $\F$}.
\end{cases}
\]
\end{ejemplo}

\begin{ejemplo}\label{DrinfeldEj1.3.30}
Ahora si $A=R_T$ y $\pL$ es cualquier lugar de $A$, $F=A/\pL$ y
$\delta\colon A\to F=A/\pL$ es la proyecci\'on can\'onica, $\car (\rho)=
\pL$.

En general, sea $\car(\rho)=\pL\neq(0)$. En el caso de que $\alpha_0=0$ se
tiene $\delta(T)=\alpha_0$ por lo que $(T)=\pL$ y $\ord(\rho_T)=
i_0= h_{\rho}v_{\pL}(T)\deg_K(T)=h_{\rho}\cdot 1\cdot 1=h_{\rho}=h$.

Por tanto si $F=A/(T)\cong\F$, $\delta\colon A\lra F$ es la proyecci\'on
can\'onica y $\rho\colon A\lra\aditivo F$ dado por $\rho_T=\tau^h
+\tau^r$ es un $A$--m\'odulo de Drinfeld de altura $h$ y rango $r$.
Es decir siempre existen m\'odulos de Drinfeld de altura $h$ y 
rango $r$.

Finalmente, si $\alpha_0\neq 0$, $\alpha_0$ es algebraico sobre $\F$
y $\pL=(\Irr(\alpha_0,T,\F))=(f(T))$, se tiene
\[
\ord\rho_{f(T)}=h_{\rho}v_{\pL}(f(T))\deg_K(f(T))=h_{\rho}\cdot 1
\cdot \deg_T f(T)=h\deg_T(f(T)),
\]
es decir $h=\frac{\ord(\rho_{f(T)})}{\deg_T(f(T))}$.
\end{ejemplo}

\subsubsection{El rango y la altura de un m\'odulo de Drinfeld}\label{DrinfeldS1.3.2''}

Probaremos que $r_{\rho}\in{\ma N}$ y que $h_{\rho}\in
{\ma N}\cup\{0\}$.

Sea $\rho\colon A\to \aditivo F$ un $A$--m\'odulo de Drinfeld. Sea
$\bar{F}$ una cerradura algebraica de $F$. Entonces $A$ act\'ua
en $\bar{F}$ v\'ia $\rho$. Sean $I$ un ideal no cero de $A$ y $a\in A$.

\begin{definicion}\label{DrinfeldD1.3.2''.1}
Se define $\rho[I]=\bigcap_{a\in I}\rho[a]=\{u\in\bar{F}\mid
\rho_a(u)=0\text{\ para toda\ } a\in I\}$, donde $\rho[a]=\rho[(a)]
=\{u\in \bar{F}\mid \rho_a(u)=0\}$.
\end{definicion}

Notemos que $\rho[I]$ es el m\'odulo $M[I]$ definido para dominios
Dedekind y $M$ un $R$--m\'odulo. En este caso $\rho=M$
es un m\'odulo de Drinfeld.

Ahora bien, puesto que $\rho_a(u)$ es un polinomio en $u$ y $\bar{F}$
es algebraicamente cerrado, $\bar{F}$ es $A$--admisible, es decir,
si $a\neq 0$, $a\in A$ y $v\in\bar{F}$, $\rho_a(u)=\sum_{i=0}^s
\alpha_i u^{q^i}=v$ es soluble para $u\in \bar{F}$.

Por el Teorema \ref{DrinfeldT1.3.13}, se tiene que para toda $e\geq 2$,
la sucesi\'on
\[
0\longrightarrow \rho[\pK]\longrightarrow \rho[\pK^e]
\longrightarrow\rho[\pK^{e-1}]\longrightarrow 0
\]
es exacta con $\pK$ un ideal primo no cero de $A$. Sea $\pi\in
\pK\setminus \pK^2$. Sea $a\neq 0$, $A/(a)$ es finito y para
todo ideal $I\neq 0$, $A/I\cong \bigoplus_{i=1}^s A/\pK_i^{e_i}$,
donde $I=\prod_{i=1}^s \pK_i^{e_i}$ y $\big|A/\pK_i^{e_i}\big|=
\big|A/\pK_i\big|^{e_i}=q^{(\deg_\K \pK_i)e_i}<\infty$ y finalmente
$\big|A/(a)\big|=q^{\deg(a)}$ (Proposici\'on \ref{DrinfeldP1.3.23}).

Sea $\pK\neq \car(\rho)=\pL$. Sea $b\in \pK$, $b\neq 0$,
\[
\rho[b]=\{u\in\bar{F}\mid \rho_b(u)=\sum_{i=0}^t\alpha_iu^{q^i}
=0\},
\]
de donde se sigue que $|\rho[b]|<\infty$ y puesto que
 $\rho[\pK]\subseteq \rho[b]$
obtenemos que $\rho[\pK]$ es finito. Por otro lado $\rho[\pK]$
es un $A$--m\'odulo y $\rho_a(u)=0$ para toda $a\in \pK$ y
$u\in \rho[\pK]$ por lo que $\rho[\pK]$ es un $A/\pK$--m\'odulo
(espacio vectorial), $\big|A/\pK\big|=q^{\deg_K \pK}$.
Sea $d:=\dim_{A/\pK}\rho[\pK]$. Entonces
\[
|\rho[\pK]|=|A/\pK|^d=q^{d\deg_K\pK}\quad\text{y}\quad
|\rho[\pK^e]|=|\rho[\pK]|^e=q^{ed\deg_\K \pK}.
\]

Se tiene que $Cl_A=C_S$, donde $S=\{\p\}$,
es finito (Corolario \ref{CRamDed1.2.6}),
donde $Cl_A$ denota al grupo de clase de $A$. Sea $m\in{\ma N}$
tal que $\pK^m=(a)$. Por tanto se tiene
\[
\rho[\pK^m]=\rho[(a)]=\rho[a]=\{u\in \bar{F}\mid \rho_a(u)=0\}.
\]

Sea $\rho_a=\delta(a)\tau^0+\cdots+\alpha_n\tau^n$. Ahora
$a\notin \pL=\ker \delta$ pues en caso contrario, tendr{\'\i}amos que
$\pK^m=(a)\subseteq \pL$ lo cual es absurdo. Por tanto $\delta(a)
\neq 0$. Se sigue que
\[
\rho_a(u)=\delta(a)u+\alpha_1u^q+\cdots+\alpha_nu^{q^n}
\]
es un polinomio separable en $u$ ya que $\rho_a(u)^{\prime}=
\delta(a)\neq 0$. Por tanto $|\rho[a]|=q^n=q^{\deg_{\tau}\rho_a}=
q^{r_{\rho}\deg a}$ con $r_{\rho}$ el rango de $\rho$. De esta
forma obtenemos
\[
q^{r_{\rho}\deg a}=|\rho[a]|=|\rho[\pK^m]|=q^{md\deg_\K  \pK}.
\]

Puesto que $\pK^m=(a)$, $\deg a=m\deg \pK$ y se tiene $r_{\rho}
\deg a=r_{\rho}m\deg_\K \pK=dm\deg_\K \pK$. Se sigue que 
\fbox{$r_{\rho}=d\in{\ma N}$}.

Adem\'as $\rho[\pK]\cong (A/\pK)^{r_{\rho}}$ como
$A$--m\'odulos y como $A/\pK$--m\'odulos (espacio vectoriales).

\begin{teorema}\label{DrinfeldT1.3.29} Sea $\rho$ un $A$--m\'odulo de 
Drinfeld. Entonces $r_{\rho}\in{\ma N}$ y $h_{\rho}\in {\ma N}\cup
\{0\}$. M\'as a\'un, si $\pK$ es un ideal primo no cero de $A$,
se tiene: si $\pK\neq \pL=\car(\rho)$, $\rho[\pK^e]\cong
(A/\pK^e)^{r_{\rho}}$ para $e\geq 1$.

Si $\pK=\pL$, $\rho[\pL^e]=(A/\pL^e)^{r_{\rho}-h_{\rho}}$.
\end{teorema}

\begin{proof} Se tiene, para $\pK\neq \pL$, $\rho[\pK]=(A/\pK)^{r_{\rho}}$. 
Suponemos por inducci\'on $\rho[\pK^{e-1}]=(A/\pK^{e-1})^{r_{\rho}}$.
Consideremos la sucesi\'on exacta
\[
0\lra\rho[\pK]\lra\rho[\pK^e]\lra\rho[\pK^{e-1}]\lra 0.
\]

Se tiene $\rho[\pK^e]({\eu q})=
\begin{cases}
\rho[\pK^e]&\text{si ${\eu q}=\pK$}\\
0&\text{si ${\eu q}\neq \pK$}
\end{cases}$.
Se sigue del Teorema \ref{DrinfeldT1.3.3} que 
\[
\rho[\pK^e]=\bigoplus_{i=1}^s A/\pK^{f_i}
\]
 para algunos $s, f_i$.

Ahora bien, $\pK^e$ anula a $\rho[\pK^e]$, esto es, si $a\in \pK^e$,
$\rho_a(u)=0$ para toda $u\in \rho[\pK^e]$. Por tanto necesariamente
se tiene que $f_i\leq e$. En caso de que alg\'un $f_i$ fuese menor a $e$,
puesto que $\rho[\pK^{e-1}]\cong (A/\pK^{e-1})^{r_{\rho}}$, esto es $f_i
\geq e-1$ y por tanto se tendr\'ia
\begin{gather*}
\rho[\pK^e]=(A/\pK^e)^{r_1}\oplus (A/\pK^{e-1})^{r_2},
\intertext{con $r_2\geq 1$. Se seguir{\'\i}a que}
|\rho[\pK^{e-1}]|=|A/\pK^{e-1}|^{r_1+r_2}=|A/\pK^{e-1}|^{r_{\rho}},
\end{gather*}
lo cual implica que $r_1+r_2=r_{\rho}$. Por tanto obtendr\'iamos
\[
|\rho[\pK^e]|=|A/\pK^e|^{r_1}\cdot |A/\pK^{e-1}|^{r_2}<|A/\pK^e|^{r_{\rho}},
\]
lo cual es absurdo. Se sigue que

\centerline{\fbox{$\rho[\pK^e]\cong (A/\pK^e)^{r_{\rho}}$}.}

Ahora si $\pK=\pL=\ker \delta$ en el caso en que $\pL\neq 0$, todo el
procedimiento anterior es igual con la \'unica excepci\'on de que si
$\pL^m=(b)$, entonces 
\begin{align*}
\rho[\pL^m]=\rho[(b)]=\rho[b]=\{u\in\bar{F}\mid \rho_b(u)&=
\alpha_hu^{q^h}+\cdots+\alpha_nu^{q^n}\\
&=\big(\beta_hu+\cdots+\beta_nu^{q^{n-h}}
\big)^{q^h}=0\},
\end{align*}
por lo que $|\rho[\pL^m]|=|\rho(b)|=q^{n-h_{\rho}}$.

Si $d^{\prime}$ es la dimensi\'on de $\rho[\pL]$ sobre $A/\pL$, entonces
$d^{\prime}=r_{\rho}-h_{\rho}$ y $\rho[\pL^e]\cong(A/\pL)^{r_{\rho}-h_{\rho}}$.
En particular $h_{\rho}\in {\ma N}\cup\{0\}$. 
$\fin$ \end{proof}

\begin{corolario}\label{DrinfeldC1.3.30}
Si ${\eu A}$ es un ideal de $A$, ${\eu A}\neq 0$ con ${\eu A}$ primo relativo
a $\pL$, entonces $\rho[{\eu A}]\cong (A/{\eu A})^{r_{\rho}}$.
\end{corolario}

\begin{proof}
Sea ${\eu A}=\pK_1^{e_1}\cdots\pK_t^{e_t}$, $\pK_i\neq \pL$ para toda
$i$ y
\begin{gather*}
\rho[{\eu A}]\cong \bigoplus_{i=1}^t\rho[\pK_i^{e_i}]\cong \bigoplus_{i=1}^t\big(
A/\pK_i^{e_i}\big)^{r_{\rho}}=\Big(\bigoplus_{i=1}^t A/\pK_i^{e_i}\Big)^{r_{\rho}}=
\big(A/{\eu A}\big)^{r_{\rho}}. \tag*{$\fin$}
\end{gather*}
\end{proof}

Ahora sean $I$ un ideal no cero de
$A$ y $\rho$ un $A$--m\'odulo de Drinfeld. Sea $J$ el ideal
izquierdo de $\aditivo F$ generado por $\{\rho_a\mid a\in A\}$, es decir
\[
J=\big\{\sum_{i=1}^n r_{a_i}\rho_{a_i}\mid r_{a_i}\in R=\aditivo F, a_i\in I\big\}.
\]\

Entonces $J$ es principal, digamos generado por $\rho_I\in J$, $J=R\rho_I$.
Notemos que en general $\rho_I$ no necesariamente est\'a en la imagen de $I$ bajo 
$\rho$, es decir, $\rho_I\notin \{\rho_a\mid a\in I\}$ (incluso, puede ser que
$\rho_I\notin \{\rho_a\mid a\in A\}$). Seleccionamos como $\rho_I$
al generador m\'onico, esto es, $\rho_I$ es m\'onico.

Si $I=(a)$ es principal, $a\neq 0$, entonces veamos 
que $\rho_I=c\rho_a$ para alg\'un $c\in\*F$. Si $\xi\in J$, entonces
$\xi=\sum_{i=1}^n r_{i}\rho_{b_i}$ para algunos $r_i\in\aditivo F$ y
$b_i\in A$, $1\leq i\leq n$. Sea $b_i=d_i a$ para alg\'un $d_i\in A$.
Entonces $\xi=\big(\sum_{i=1}^n r_i \rho_{d_i}\big)\rho_a$
con $\sum_{i=1}^n r_i\rho_{d_i}\in \aditivo F$, de
donde se sigue que $\rho_a$ es generador de $J$.

\begin{proposicion}\label{DrinfeldP1.3.31} Para $I\neq 0$
se tiene $\rho[I]=\ker \rho_I$, esto es, $\rho_I\colon\bar{F}\to\bar{F}$,
$\ker \rho_I=\{\xi\in\bar{F}\mid \rho_I(\xi)=0\}$ y $\deg_{\tau} \rho_I=r_{\rho}
\deg I$, donde $|A/I|=q^{\deg I}$.
\end{proposicion}

\begin{proof} Si $u\in \rho[I]$, $\rho_a(u)=0$ para toda $a\in I$. Puesto que $\rho_I\in
J$, $\rho_I=\sum_{a\in I}\xi_a\rho_a$ para algunas $\xi_a\in \aditivo F$. Por tanto
\[
\rho_I(u)=\sum_{a\in I}\xi_a(\rho_a(u))=\sum_{a\in I}\xi_a(0)=\sum_{a\in I}0=0,
\]
de donde $u\in \ker\rho_I$.

Rec{\'\i}procamente, si $u\in \ker \rho_I$ y $a\in I$, entonces $\rho_a\in J$ por lo que
$\rho_a=\xi\rho_I$ para alguna $\xi\in\aditivo F$. Se sigue que $\rho_a(u)=
\xi(\rho_I(u))=\xi(0)=0$. Por tanto $u\in \rho[I]$ y 
\fbox{$\rho[I]=\ker \rho_I$}.

Ahora si consideramos $\bar{F}$ una cerradura algebraica de $F$, $\rho[I]
\cong (A/I)^{r_{\rho}}$ para $I$ primo relativo a la caracter{\'\i}stica,
se tendr\'ia
\[
|\rho[I]|=q^{r_{\rho}\deg I}=|\ker \rho_I|=|\{u\in\bar{F}\mid \rho_I(u)=
\sum_{i=0}^{\deg\rho_I}\alpha_iu^{q^i}=0\}|=q^{\deg \rho_I},
\]
por lo tanto \fbox{$r_{\rho}\deg I=\deg \rho_I$}.

Ahora supongamos que la caracter{\'\i}stica de $A$ es $\pL\neq 0$ y que 
${\eu A}$ es un ideal no cero de $A$ divisible por $\pL$. Presentamos
la demostraci\'on en este caso tambi\'en aqu{\'\i} aunque usaremos definiciones
y resultados que veremos m\'as adelante. Por la teor{\'\i}a general
de dominios de Dedekind, existen elementos $a,b\in A$ y un ideal ${\eu B}$ 
primo relativo a $\pL$ tal que $a{\eu A}=b{\eu B}$. De hecho tenemos

\begin{lema}\label{DrinfeldL1.3.32}
Sea $R$ un dominio Dedekind y sean ${\eu q}, {\eu a}$ ideales no cero de $R$
tales que ${\eu q}\mid {\eu a}$. Entonces existen $a,b\in R$ tales que $a{\eu a}=
b{\eu b}$ con $\mcd({\eu b},{\eu q})=1$.
\end{lema}

\begin{proof} Sea $b\in A$ tal que $v_{\eu q}(b)=v_{\eu q}({\eu a})$ y $v_{\eu t}(b)=0$
para todo ideal ${\eu t}\mid {\eu c}$ donde ${\eu a}={\eu q}^n{\eu c}$,
$v_{\eu q}({\eu a})=n$ y $\mcd({\eu c},{\eu q})=1$. 
Por tanto $(b)={\eu q}^n{\eu d}$
con $\mcd({\eu d},{\eu q})=1$ y $\mcd({\eu d},{\eu c})=1$.

Sea  $a\in A$ con $v_{\eu q}(a)=0$ y $v_{\eu t}(a)=v_{\eu t}({\eu d})$ para todo
lugar ${\eu t}\mid {\eu d}$. Entonces $(a)={\eu d}{\eu d}^{\prime}$, $\mcd({\eu d}{
\eu d}^{\prime},{\eu q})=1$. Se tiene $a{\eu a}=(a){\eu q}^n{\eu c}={\eu d}{\eu d}^{
\prime}{\eu q}^n{\eu c}=(b){\eu d}^{\prime}{\eu c}$. Si ${\eu b}={\eu d}^{\prime}{\eu c}$,
$a{\eu a}=b{\eu b}$ y $\mcd({\eu b},{\eu q})=1$.
$\fin$ \end{proof}

M\'as adelante veremos (Proposici\'on \ref{DrinfeldP3.1.27} (1))
que $\rho_{{\eu A}{\eu B}}=({\eu A}\star\rho)_{\eu B}\rho_{\eu A}$
donde $\star$ es una cierta acci\'on de
ideales de $A$ sobre los m\'odulos de Drinfeld (ver Ecuaci\'on
(\ref{DrinfeldEq3.1.15'}) antes de la Observaci\'on \ref{DrinfeldO3.1.16}).

Usando lo anterior, se tiene $\rho_{{\eu A}a}=({\eu A}\star \rho)_{(a)}\rho_{\eu A}=
({\eu B}\star \rho)_{(b)}\rho_{\eu B}=\rho_{{\eu B}b}$. 
Puesto que $\rho_{(a)}$ difiere de 
$\rho_a$ por un m\'ultiplo no cero de $\K $ y 
puesto que m\'odulos is\'ogenos (lo cual
quiere decir que existe un homomorfismos entre ellos), tienen el mismo rango,
aplicando esto a $\rho_{{\eu A}a}=\rho_{{\eu B}b}$, obtenemos que
\[
r\deg a+\deg_{\tau}\rho_{\eu A}=r\deg b+\deg_{\tau} 
\rho_{\eu B}\underbracket[0pt]{=}_{
\substack{\uparrow\\ \mcd({\eu B},\pL)=1}} r\deg b+ r\deg {\eu B}.
\]

Adem\'as, puesto que $a{\eu A}=b{\eu B}$, $\deg a+\deg 
{\eu A}=\deg b+\deg_{\tau}{\eu B}$.
Por tanto
\begin{align*}
\deg_{\tau}\rho_{\eu A}&=r\deg b+r\deg {\eu B}-r\deg a=r\deg
{\eu B}-r\deg\big(\frac{a}{b}\big)\\
&=r\deg {\eu B}-r\deg \big(\frac{\eu B}{\eu A}\big)=r\deg{\eu A}.
\end{align*}
Finalmente obtenemos que 
\fbox{$\deg_{\tau}\rho_{\eu A}=r\deg {\eu A}$}.
$\fin$ \end{proof}

\section{Construcci\'on y propiedades de m\'odulos de Drinfeld}\label{DrinfeldC2}

\subsection{Funciones exponenciales}\label{DrinfeldS2.1}

La funci\'on exponencial usual, $e\colon{\ma C}\lra {\ma C}$ dada por 
$e^z=\sum_{n=0}^{\infty}\frac{z^n}{n!}$ es una funci\'on entera, cada elemento
de $\*{\ma C}$ se toma una infinidad (numerable) de veces y $0\notin \im e$.
Adem\'as $e^{z+w}=e^ze^w$ para cualesquiera $z,w\in {\ma C}$, es decir,
$e$ es una funci\'on multiplicativa.

Ahora, si $F$ es un campo de caracter{\'\i}stica $p>0$ y $f\colon F\to F$
satisface que $f(x+y)=f(x)f(y)$ para cualesquiera $x,y\in F$, esto es, $f$ es
multiplicativa, entonces $f(x)^p=f(px)=f(0)=f(0+0)=f(0)^2$ de donde se sigue
que $f(0)^2=f(0)$ lo cual implica que $f(0)=0$ o $1$. Puesto que $f(x)^p=f(0)$,
se sigue que $f(x)=0$ o $1$ para toda $x\in F$. De esta forma obtenemos
que $f$ es id\'enticamente igual a $0$ o a $1$. Esto es, no hay funciones 
exponenciales como tales en caracter{\'\i}stica positiva excepto las funciones
triviales $f(x)=0$ o $1$ para toda $x$.

Sin embargo, como vimos en la Secci\'on \ref{DrinfeldC1}, hay numerosas funciones
aditivas $f\colon F\lra F$, por ejemplo, los polinomios aditivos.

Ahora bien, necesitamos un campo similar a ${\ma C}$ pero en caracter{\'\i}stica
$p>0$. Se tiene que ${\ma C}$ es completo, lo cual se refiere a una de sus propiedades
topol\'ogicas m\'as importantes, y es algebraicamente cerrado, la cual es una
de sus propiedades centrales desde el punto de vista algebraico. 

Sea $A$ en general,
m\'as precisamente, sea $\K $ un campo de funciones congruente, $\p$ un lugar
fijo de $\K $, $d_{\infty}=\deg \p$ y $A=\{x\in \K \mid v_{\pK}(x)\geq 0  \text{\ para
todo lugar\ } \pK\neq \p\}$. Sea $\K_{\infty}=\K_{\p}$ la completaci\'on de $\K $
en $\p$. Se tiene que $\K_{\infty}\cong {\ma F}_{q^{d_{\infty}}}((\pi))$ donde $\pi
\in \K $ es un elemento primo en $\p$, es decir, $v_{\infty}(\pi)=v_{\p}(\pi)=1$
(\cite[Theorem 2.5.20]{Vil2006}).

Entonces $\Ki$ es completo pero no es algebraicamente cerrado. 
Por ejemplo, $\sqrt[2]{\pi}$ es algebraico pero $\sqrt[2]{\pi}\notin \K_{\infty}$.
Si $p\neq 2$, $\sqrt[2]{\pi}$ es adem\'as separable. Si $p=2$, se puede
tomar $\sqrt[3]{\pi}$. En general ${\ma F}_{q^{2d_{\infty}}}((\pi))$ es una
extensi\'on algebraica propia de $\K_{\infty}$.

Sea $\bar{\K }_{\infty}$ una
cerradura algebraica de $\Ki$. Entonces $\bar{\K }_{\infty}$ es 
algebraicamente cerrado
pero no es completo. 
Por ejemplo $y=\sum_{m=1}^{\infty} \pi^{m/2}$ o $\sum_{m=1}^{\infty}
\pi^{m/3}$ no es algebraico sobre $\K_{\infty}$ por lo que $y\notin
\bar{\K }_{\infty}$ pero si $y_n=\sum_{m=1}^n \pi^{m/2}$, se
tiene $|y-y_n|=c^{-(n+1)/2}$ con alg\'un $0<c<1$, por lo que 
$|y-y_n|\xrightarrow[n\to\infty]{} 0$.

Consideremos una extensi\'on de 
$\vi$ a $\bar{\K }_{\infty}$. Sea
${\ma C}_p=\Ci$ la completaci\'on de $\bar{\K }_{\infty}$ 
con respecto a $\vi$. Entonces
$\Ci$ es completo y algebraicamente cerrado
(ver la Proposici\'on \ref{DrinfeldP2.1.2}). 
$\Ci$ es el an\'alogo a ${\ma C}$ en
caracter{\'\i}stica $p$.

Notemos que si $x\in \K $ es tal que $\eta_x=\p^a$, $a\geq 1$ 
y si ${\mathcal P}_{\infty}$
es el primo infinito en $\F(x)$, $[\K :\F(x)]<\infty$ por lo que 
$[\Ki:\F(x)_{{\mathcal P}_{\infty}}]
<\infty$. Se sigue que $\bar{\K }_{\infty}=\overline{\F(x)}_{{\mathcal P}_{\infty}}$. Esto es,
$\Ci$ se puede construir a partir del primo infinito en un campo de funciones racionales
$\F(T)$.

Para probar las afirmaciones sobre $\Ci$, tenemos en general que si $\K $ es un campo
completo con respecto a una valuaci\'on (no necesariamente discreta), $v\colon \K \lra
{\ma R}\cup\{0\}$, $|x|_v=\alpha^{v(x)}$ con $0< \alpha <1$ fijo. Entonces
\begin{gather*}
\o_\K =
\{\alpha \in \K \mid |x|_v\leq 1\}=\{x\in \K \mid v(x)\geq 0\}, 
\intertext{es el anillo de enteros o anillo de valuaci\'on de $\K $ y}
\pK_\K =\{x\in \K \mid |x|_v<1\}=\{x\in \K \mid v(x)>0\}
\end{gather*}
es el ideal m\'aximo de $\o_\K $. Se tiene que
$\o_\K $ es un anillo local con ideal m\'aximo $\pK_\K $ y $\K (\pK)=\o_\K /\pK_\K $ es el campo
residual de $\K $.

\begin{proposicion}\label{DrinfeldP2.1.1} Si $\K $ es completo con respecto a una valuaci\'on $v$
y $L/\K $ es una extensi\'on finita, entonces existe una \'unica extensi\'on $w$ de $v$
a $L$, $L$ es completo con respecto a la valuaci\'on $w$ y $w$ est\'a dada de la siguiente
forma: si $y\in \*L$, $\N_{L/\K }y\in \K $, entonces, 
\[
w(y)=\frac{1}{[L:\K ]}v(x)\quad \text{y}\quad |y|_w=\sqrt[n]{\big|\N_{L/\K }(y)\big|},
\]
con $n=[L:\K ]$.
\end{proposicion}

\begin{proof} La existencia de $w$ est\'a dada por el Lema de Chevalley
(\cite[Theorem 2.4.3]{Vil2006}). Veamos la unicidad.
Sean $\alpha\in \*L$ y $\beta=\alpha^n/\N\alpha$ donde $\N=\N_{L/\K }$. Entonces
$\N(\beta)=\frac{\N\alpha^n}{(\N\alpha)^n}=\frac{(\N\alpha)^n}{(N\alpha)^n}=1$.

Sea $\gamma\in L$ con $|\gamma|_w<1$. Sea $\{\xi_1,\ldots,\xi_n\}$ una base de
$L/\K $. Dado $t\in{\ma N}$, escribimos $\gamma^t=x_1^{(t)}\xi_1+\cdots+x_n^{(t)}\xi_n$
con $x_i^{(t)}\in \K $.

Puesto que $|\gamma|_w<1$ se sigue que $\gamma^t\xrightarrow[t\to \infty]{}0$
y por ende $x_i^{(t)}\xrightarrow[t\to\infty]{}0$ para toda $1\leq i\leq n$.

Ahora bien, $\N(\gamma^t)=(\N\gamma)^t$ es un polinomio homog\'eneo en 
$x_1^{(t)},\ldots,x_n^{(t)}$ por lo que $\N(\gamma^t)=(\N\gamma)^t\xrightarrow[
x\to\infty]{}0$ lo cual implica que $|\N\gamma|_v<1$.

De manera similar, si $|\gamma|_w>1$ se obtiene que $|\N\gamma|_v>1$. De
ah{\'\i} se sigue que $|\gamma|_w=1$ si $|\N\gamma|_v=1$.
Puesto que $|\N\beta|_v=1$, se sigue que $1=|\beta|_w=
\frac{|\alpha|^n_w}{|\N\alpha|_v}$. Por tanto 
\fbox{$|\alpha|_w=\sqrt[n]{\big|\N\alpha\big|}$}
por lo que $w$ es \'unica y se sigue la f\'ormula propuesta.
$\fin$
\end{proof}

\begin{proposicion}\label{DrinfeldP2.1.2}
Sea $\K $ un campo completo con valuaci\'on $v$. Sea $\bar{\K }$ una cerradura 
algebraica de $\K $ junto con la \'unica extensi\'on de $v$ a $\bar{\K }$, la cual 
tambi\'en llamaremos $v$. Sea $\hat{\bar{\K }}$ la completaci\'on de $\bar{\K }$
con respecto a $v$. Entonces $\hat{\bar{\K }}$ es completo y algebraicamente
cerrado. En particular $\Ci={\ma C}_p$ es completo y algebraicamente cerrado.
\end{proposicion}

\begin{proof} Sea $f(x)=\sum_{i=0}^n\alpha_i x^i\in \hat{\bar{\K }}[x]$ con $\alpha_n=1$.
Se probar\'a que $f(x)$ tiene una ra{\'\i}z en $\hat{\bar{\K }}$. Sea $L=\hat{\bar{\K }}
(\beta)$ con $\beta$ una ra{\'\i}z de $f(x)$ y nuevamente denotamos por $v$
a la extensi\'on de $v$ a $L$. Suponemos $\beta\neq 0$. 
Sea $M_1>0$ y sea $f_1(x)=\sum_{i=0}^n
\tilde{\alpha}_ix^i\in \bar{\K }[x]$ tal que $\delta<\min_{0\leq i\leq n}\{v(\alpha_i-
\tilde{\alpha}_i)\}$, es decir, $f_1(x)\in\bar{\K }[x]$ es un polinomio ``cercano'' a
$f(x)$. Entonces
\[
f_1(\beta)=f_1(\beta)-f(\beta)=\sum_{i=0}^n (\tilde{\alpha}_i-\alpha_i)\beta^i.
\]

Sean $\{\beta_1,\ldots,\beta_n\}$ las ra{\'\i}ces de $f_1$ en $\bar{\K }$, el cual
es algebraicamente cerrado. De esta forma
\[
f_1(x)=\prod_{j=1}^n(x-\beta_j)\quad \text{y}\quad f_1(\beta)=\prod_{j=1}^n
(\beta-\beta_j)=\sum_{i=0}^n(\tilde{\alpha}_i-\alpha_i)\beta^i.
\]

Sea $\xi_{\beta}=\min_{0\leq i\leq n}\{v(\beta^i)\}=\min_{0\leq i\leq n}\{iv(\beta)\}$.
Entonces
\begin{align*}
\sum_{j=1}^n v(\beta-\beta_j)&=v\big(\prod_{j=1}^j(\beta-\beta_j)\big)=
v\big(\sum_{i=0}^n(\tilde{\alpha}_i-\alpha_i)\beta^i\big)\\
&\geq \min_{0\leq i\leq n}\{v(\tilde{\alpha}_i-\alpha_i)+iv(\beta)\}>
M_1 + \xi_{\beta}.
\end{align*}

Para alg\'un {\'\i}ndice $1\leq j(1)\leq n$ se tiene $v(\beta-\beta_{j(1)})>\frac{
M_1 +{\xi_{\beta}}}{n}$. Para $m\in{\ma N}$, sean $M_m=m$ y $f_m(x)$
en lugar de $M_1$ y de $f_1(x)$ respectivamente en el desarrollo anterior.
Entonces para $m\in{\ma N}$, con las notaciones anteriores, se tiene:
\[
v(\beta-\beta_{j(m)})>\frac{m+\xi_{\beta}}{n}\xrightarrow[m\to\infty]{} \infty
\Lra |\beta-\beta_{j(m)}|_v\xrightarrow[m\to\infty]{}0,
\]
y $\{\beta_{j(m)}\}\subseteq \bar{\K }$. Por lo tanto $\beta\in \hat{\bar{\K }}$ y 
por tanto $\hat{\bar{\K }}$ es algebraicamente cerrado.
$\fin$
\end{proof}

\begin{observacion}\label{Drinfeld2.1.3} Como $|\ |_v$ es no arquimediano, $\sum_{n=0}^{
\infty}a_n$ converge en $\K \iff \lim_{n\to\infty}a_n=0$ pues $\K $ es completo.
\end{observacion}

\subsection{Redes y funciones exponenciales}\label{DrinfeldS2.2}

Si $A=R_T$ y $F$ es cualquier $\F$--\'algebra, existen muchos m\'odulos de 
Drinfeld $\rho \colon A\to\aditivo F$, donde $\delta\colon A\lra F$ ha sido dado. Esto se
debe a que $A$ es una $\F$--\'algebra libre generada por $T$, por lo que
$\rho_T=\delta(T)+\sum_{i=1}^r\alpha_i\tau^i$, con $\alpha_r\neq 0$, $r\geq 1$,
es arbitrario, entonces $\rho(M(T))=M(\rho(T))$ para $M\in R_T$.

Cuando $A\neq R_T$ no sabemos que existan m\'odulos de Drinfeld $\rho\colon
A\lra \aditivo F$; de hecho para algunos campos $F$ no existen m\'odulos de
Drinfeld sobre $A$. 

\begin{ejemplo}\label{DrinfeldEj2.2.0}
Sean $q=p=3$, $\K ={\ma F}_q(T)$, $\p$ el lugar asociado a $T^2+1$ y $A=
\{x\in \K \mid v_{\eu p}(x)\geq \text{\ para todo lugar ${\eu p}\neq \p$}\}$.
Entonces
\[
A=\Big\{\frac{G(T)}{(T^2+1)^n}\mid G(T)\in{\ma F}_q[T], n\in{\ma N}, 
\deg G(T)\leq 2n\Big\}.
\]

Puesto que $\p$ es de grado $2$ y $h_\K =1$, se tiene $h_A=2$
(ver Corolario \ref{CRamDed1.2.7}).

Sea $\xi=\frac{1}{T^2+1}$ y consideremos $R_{\xi}={\ma F}_q[\xi]$, 
${\ma F}_q(\xi)=\coc R_{\xi}$. Entonces $A$ es la cerradura entera
de $R_{\xi}$ en $\K $. Usando el algoritmo de la divisi\'on, se sigue que
si $x\in A$, entonces $x=\frac{G(T)}{(T^2+1)^n}$, $\deg G(T)\leq 2n$ y
\[
G(T)=\alpha_0+\alpha_1 (T^2+1)+\cdots +\alpha_n(T^2+1)^n
=\alpha_0+\alpha_1\xi^{-1}+\cdots +\alpha_n\xi^{-n},
\]
con $\alpha_i\in {\ma F}_q[T]$ de grado menor o igual a $1$. Adem\'as,
puesto que $\deg G(T)\leq n$, se tiene que $\alpha_n\in{\ma F}_q$.

Por tanto tenemos
\begin{align*}
x=\frac{G(T)}{(T^2+1)^n}=\xi^n G(T)&=\alpha_n+\alpha_{n-1}\xi+\cdots
+\alpha_1\xi^{n-1}+\alpha_0\xi^n\\
&=\beta_0+\beta_1\xi+\cdots+\beta_{n-1}\xi^{n-1}+\beta_n\xi^n,
\end{align*}
con $\beta_i=\alpha_{n-1}=a_i+b_iT\in{\ma F}_q[T]$, $0\leq i\leq n$
y $\beta_0=a_0$.

De esta forma obtenemos
\begin{align}\label{DrinfeldEq2}
x=\xi^nG(T)&=\sum_{i=0}^n a_i\xi^i+T\sum_{i=1}^n b_i\xi^i\nonumber\\
&=\sum_{i=0}^n a_i\xi^i+(T\xi)\sum_{i=0}^{n-1}b_{i+1}\xi^i=F(\xi)+(T\xi)H(\xi)
\end{align}
con $F(\xi), G(\xi)\in R_{\xi}$, $\deg F(\xi)\leq n$, $\deg H(\xi)\leq n-1$.

Notemos que el grado de $F(\xi)$ en $T$ es par y el grado de $(T\xi)H(\xi)$
en $T$ es impar, de donde se sigue que $x=0\iff F(\xi)=H(\xi)=0$. En particular
$\{1,T\xi\}$ es una base entera de $A/R_{\xi}$. Por otro lado,
puesto que $\xi=\frac{1}{T^2+1}$, se sigue que $(\xi T)^2=-\xi^2+\xi$.
Por tanto 
\[
\ell(Z):=\Irr(Z,T\xi,{\ma F}_q(\xi))=Z^2+\xi^2-\xi.
\]

Sea ${\ma F}_q[X,Y]\xrightarrow[\phantom{xxxx}]{\phi} A$ dada por $\phi(f(X,Y))=
f(\xi,T\xi)$. Entonces por (\ref{DrinfeldEq2}), se sigue que $\phi$ es un epimorfismo
de anillos. Adem\'as $\phi(Y^2+X^2-X)=0$, esto es, $\langle Y^2+X^2-X\rangle
\subseteq \ker \phi$ y $\phi$ induce el epimorfismo $\tilde{\phi}\colon
{\ma F}_q[X,Y]/\langle Y^2+X^2-X\rangle\longrightarrow A$ dada por
$\tilde{\phi}(f(X,Y)\bmod \langle Y^2+X^2-X\rangle) = f(\xi,T\xi)$.
De (\ref{DrinfeldEq2}) se sigue que $\tilde{\phi}$ es un isomorfismo de anillos.

Podemos aplicar el Teorema de Kummer sobre descomposici\'on de
ideales primos en la extensi\'on $A/R_{\xi}$ pues $A=R_{\xi}[T\xi]$.
En particular tenemos
\begin{gather*}
\ell(Z)\bmod \xi=Z^2;\qquad \ell(Z)\bmod (\xi-1)=Z^2;
\intertext{de donde}
(\xi)={\eu p}_{\xi}^2\quad \text{con}\quad {\eu p}_{\xi}=(\xi,T\xi)\quad\text{y}\\
(\xi-1)={\eu p}_{\xi-1}^2\quad \text{con}\quad {\eu p}_{\xi-1}=
(\xi-1,T\xi),
\intertext{${\eu p}_{\xi}$ y ${\eu p}_{\xi-1}$
ideales primos de $A$. Adem\'as $(T\xi)^2=\xi(1-\xi)$, por lo que}
(T\xi)={\eu p}_{\xi} {\eu p}_{\xi-1}.
\end{gather*}
As{\'\i}, $\xi,\xi-1$ y $T\xi$ son elementos irreducibles no primos de $A$
y $(T\xi)^2=\xi(1-\xi)$.

Veamos que no existe ning\'un m\'odulo de Drinfeld de rango $1$ sobre
$\K $.

Sea $\rho\colon A\lra \aditivo \K $ un $A$--m\'odulo de Drinfeld
de rango $1$ donde $\K =\coc A$.
De hecho, puesto que $\deg \xi=\deg T\xi=2$,
 $\rho$ est\'a determinado por
\begin{gather*}
\rho_{\xi}=\xi + \gamma_1\tau+\gamma_2 \tau^2;\qquad
\rho_{T\xi}=T\xi+\epsilon_1\tau+\epsilon_2 \tau^2
\intertext{y puesto que $(T\xi)^2=\xi(1-\xi)$, por la relaci\'on}
\rho_{(T\xi)^2}=\rho_{T\xi} \rho_{T\xi}=\rho_{\xi}(1-\rho_{\xi})=\rho_{
\xi(1-\xi)}.
\end{gather*}

Se tiene que 
\begin{align*} \rho_{(T\xi)^2}&=\{\text{t\'erminos de grado menor}\}
+ \epsilon_2 \tau^2 \epsilon_2\tau^2\\
&=\{\text{t\'erminos de grado
menor}\} + \epsilon_2^{10}\tau^4
\intertext{y} 
\rho_{\xi(1-\xi)}&=\{\text{t\'erminos de grado menor}\} -
 \gamma_2\tau^2\gamma_2\tau^2\\
 &= \{\text{t\'erminos de grado menor}\}-\gamma_2^{10}\tau^4.
 \end{align*}
 
 En particular $\epsilon_2^{10}=-\gamma_2^{10}$, es decir
 $\sqrt[10]{-1}=\frac{\epsilon_2}{\gamma_2}\in \K $ lo cual es absurdo
 pues $q=3$ y $\sqrt[2]{-1}\notin \K $.
 
 En particular no existen m\'odulos de Drinfeld de rango $1$ sobre $\K $.
\end{ejemplo}

Necesitamos un procedimiento especial para hallar m\'odulos
de Drinfeld sobre $\Ci$. Para hacer esto necesitaremos varios conceptos de
tipo anal{\'\i}tico: funciones exponenciales, funciones enteras, productos infinitos,
redes, etc. Empezamos por generalizar los m\'odulos de Drinfeld.

Sean $A$ y $F$ como siempre. Se usa $\delta\colon \K \lra F$ para la extensi\'on
de $\delta$ a $\K =\coc A$ cuando $\delta$ es un monomorfismo. Esto es,
$\delta\big(\frac{a}{b}\big)=\frac{\delta(a)}{\delta(b)}$. Sea $F\torcido=
\big\{\sum_{i=0}^{\infty}\alpha_i\tau^i\mid\alpha_i\in F\big\}$ con la multiplicaci\'on
de series de Taylor sujeta a $\tau\alpha=\alpha^q\tau$ con $\alpha\in F$.

\begin{definicion}\label{DrinfeldD2.2.1} Al anillo $F\torcido$ se le llama {\em
las series de potencias torcidas\index{series de potencias torcidas}}.
\end{definicion}

Como antes, se define la aumentaci\'on $D\colon F\torcido
\lra F$ por $D(f(\tau))=f(0)$.

\begin{definicion}\label{DrinfeldD2.2.2}  Sea $\delta\colon A\lra F$ un monomorfismo
y $\delta\colon \K \lra F$ la extensi\'on. Por un
 {\em $\K $--m\'odulo formal\index{m\'odulo formal}} entenderemos
a un homomorfismo de anillos $\rho\colon \K \lra F\torcido$ 
tal que $D\circ \rho=
\delta$ y tal que $\rho(\alpha)\notin F$ para alg\'un $\alpha\in \K $.
\end{definicion}

\begin{observacion}\label{DrinfeldO2.2.3} Sea $\rho\in\Drin_A(F)$ tal que $\car(\rho)=0$.
Entonces $\rho_a(0)\in \*F$ para $a\neq 0$, donde $\rho_a(0)$ es la constante
del polinomio torcido $\rho_a$. Entonces, como elemento de $F\torcido\supseteq
\aditivo F$, $\rho_a$ es invertible. La demostraci\'on es igual que en series formales.
Como $\rho_a$ es invertible para toda $a\in A\setminus\{0\}$, $\rho$ se extiende
a $\rho\colon \K \lra F\torcido$ por $\rho_{a/b}=\rho_a\rho_b^{-1} \big(=\frac{\rho_a}
{\rho_b}\big)$.
\end{observacion}

De ahora en adelante consideraremos $\delta\colon A\lra \Ci$ el mapeo inclusi\'on.
El objetivo es probar el {\em Teorema de uniformizaci\'on anal{\'\i}tica} para $A$--m\'odulos
de Drinfeld sobre $\Ci$, esto es, si $\rho\in\Drin_A(\Ci)$, existe una \'unica red 
$\Gamma$ en $\Ci$ tal que $\rho=\rho^{\Gamma}$ (m\'as adelante veremos el
significado de $\rho^{\Gamma}$, Teorema \ref{DrinfeldT2.2.19}).
Por supuesto, primero veremos que existen 
m\'odulos de Drinfeld $\Drin_A(\Ci)$ a granel.

\subsubsection{Redes}\label{DrinfeldS2.Redes}

\begin{definicion}\label{DrinfeldD2.2.4} Una {\em red\index{red}} $\Gamma$ es un $A$--subm\'odulo
discreto finitamente generado de $\Ci$. Es decir, $\Gamma$ es discreto en la
topolog{\'\i}a de $\Ci$ y la acci\'on de $A$ en $\Gamma$ es multiplicaci\'on en $\Ci$
y, por supuesto, $a\circ \alpha\in\Gamma$ para toda $a\in A$ y para toda $\alpha\in\Gamma$.
\end{definicion}

\begin{definicion}\label{DrinfeldD2.2.5} Si $\Gamma$ es una red, la dimensi\'on sobre
$\Ki$ del $\Ki$ espacio vectorial $\Ki\otimes_A \Gamma$ se llama el
{\em rango\index{rango de una red} de $\Gamma$} y se denota por
$r_{\Gamma}:=\dim_{\Ki}\Ki \Gamma$.
\end{definicion}

\begin{ejemplo}\label{DrinfeldEj2.2.5(1)}
Sea $x\in \K_{\infty}\setminus A$ y sea $\Gamma=A+Ax$. El n\'umero m\'inimo
de generadores de $\Gamma/A$ es $2$ y $\K_{\infty}\Gamma =\Gamma$ por
lo que $\dim_{\K_{\infty}}\Gamma=1$. Se tiene que $\Gamma$ no es discreto.

De hecho, puesto que $x\in \K_{\infty}$, existe $\{x_n\}_{n=1}^{\infty}
\subseteq \K $ tal que $\lim_{n\to \infty}x_n=x$. Sea $x_n=\frac{a_n}{b_n}$
con $a_n,b_n\in A$. Entonces 
\[
x-x_n=x-\frac{a_n}{b_n}=\frac{1}{b_n}(b_n x-a_n)\xrightarrow[n\to\infty]{} 0.
\]

Por tanto $M$ no es una red.
\end{ejemplo}

Se tiene

\begin{proposicion}\label{DrinfeldP2.2.5(2)}
Sea $M$ un $A$--m\'odulo proyectivo contenido en un $\K $ espacio
vectorial $V$ de dimensi\'on finita. Entonces $M$ es discreto si y solamente
si $\theta\colon \K \otimes_A M\lra V$, $\alpha\otimes_A m\longmapsto 
\alpha m$ es inyectivo.
\end{proposicion}

\begin{proof} \cite[Proposition 4.6.3, p\'agina 74]{Gos96}. $\fin$
\end{proof}

\begin{observacion}\label{DrinfeldO2.2.5(3)}
Dada una $A$--red $M$, se define el {\em rango de $M$} como su
rango de $A$--subm\'odulo finitamente generado y libre de
torsi\'on $\Ci$, esto es, $M$ como $A$--m\'odulo proyectivo finitamente
generado.

Del hecho de que $M$ es discreto, se tiene que el rango de $M$ coincide con
$\dim_{\K_{\infty}} (\K_{\infty}\otimes_A M)=\dim_{\K_{\infty}} \K_{\infty} M$.

En el Ejemplo \ref{DrinfeldEj2.2.5(1)} $M$ tiene rango $2$ como $A$--m\'odulo
sin embargo $\dim_{\K_{\infty}}\K_{\infty}M=1\neq 2$.
\end{observacion}

\begin{proposicion}\label{DrinfeldP2.2.5(5)} Si $\xi\in\Ci$, $\xi\neq 0$ y $I$
es un ideal fraccionario de $A$, $I\xi=\Gamma$ es discreto.
\end{proposicion}

\begin{proof} Si $x\in A$, $x\neq 0$, $v_{\infty}(x)=v_{{\eu p}_{\infty}}(x)
\leq 0$ por lo que $|x|_{\infty}=q^{-v_{\infty}(x)}\geq 1$.
Sea $\alpha\in A$, $\alpha\neq 0$ tal que $\alpha I=J
\subseteq A$ y $|x|_{\infty}\geq 1$ para toda $x\in J$ por lo que
$|\alpha y|_{\infty}\geq 1$ para toda $y\in I$. Se sigue que $|y|_{
\infty}\geq\frac{1}{|\alpha|_{\infty}}$ para toda $y\in I$. De esta
forma se obtiene que $|\xi y|_{\infty}=|\xi|_{\infty}|y|_{\infty} \geq
\frac{|\xi|_{\infty}}{|\alpha|_{\infty}}$ para toda $y\in I$. Por tanto
\[
B\big(0,r\big)\cap \Gamma=\{0\},\quad r=\frac{|\xi|_{\infty}}{|\alpha|_{\infty}}
\]
de donde obtenemos que $\Gamma$ es discreto. $\fin$
\end{proof}

\begin{corolario}\label{DrinfeldC2.2.5(6)} Si $I$ es un ideal de $A$, $J$ es 
un ideal fraccionario de $A$, entonces $A$, $I$ y $J$ son $A$--m\'odulos
discretos. $\fin$
\end{corolario}

\begin{proposicion}\label{DrinfeldP2.2.5(7)} Si $\{\xi_1,\ldots,\xi_r\}\subseteq \Ci$
es un conjunto linealmente independiente sobre $\K_{\infty}$ y 
$J_1,\ldots,J_r$ son ideales fraccionarios de $A$, entonces el $A$--m\'odulo
\[
\Gamma=J_1\xi_1+\cdots+J_r\xi_r
\]
es discreto.
\end{proposicion}

\begin{proof}
En caso contrario, para toda $n\in{\ma N}$, se tiene $B\big(0,\frac{1}{n}\big)
\cap \Gamma\neq \{0\}$. Sea $x_n\in \Gamma$, $x_n\neq 0$ tal que 
$|x_n|_{\infty}<\frac{1}{n}$. 

Sea $x_n=a_n^{(1)}\xi_1+\cdots+a_n^{(r)}\xi_r\in \Gamma$, $a_n^{(i)}\in J_i$,
$1\leq i\leq r$.

Se tiene que $\{a_n^{(j)}\}_{n\in{\ma N}}$ converge y $a_n^{(j)}\neq 0$
para una infinidad de \'indices $n$ y alg\'un $j$, digamos $a_n^{(j)}
\xrightarrow[n\to\infty]{} a_0^{(j)}\neq 0$ pues $|a_n^{(j)}|_{\infty}\geq 1$.

Se sigue que
\[
x_n\xrightarrow[n\to\infty]{} 0=a_0^{(1)}\xi_1+\cdots+a_0^{(r)}\xi_r
\]
con alg\'un $a_0^{(j)}\neq 0$ lo que contradice la independencia
lineal de $\{\xi_1,\ldots,\xi_r\}$ sobre $\K_{\infty}$. Por tanto
$\Gamma=J_1\xi+\cdots+J_r\xi_r$ es discreto. $\fin$
\end{proof}

\begin{ejemplo}\label{DrinfeldEj2.2.6} Se tiene $[\Ci:\Ki]=\infty$. Sea $r\geq 1$ arbitrario
y sea $\{\alpha_1,\ldots,\alpha_r\}\subseteq \Ci$ un conjunto linealmente independiente
sobre $\Ki$. Sean $a_1,\ldots,a_r\in A\setminus\{0\}$. Sea
\[
\Gamma = A\frac{\alpha_1}
{a_1}+\cdots+A\frac{\alpha_r}{a_r}=\Big\{\sum_{i=1}^r 
x_i\frac{\alpha_i}{a_i}\mid x_i\in A\Big\}.
\]
Entonces $\K_{\infty}\Gamma=\K_{\infty}\alpha_1\oplus\cdots\oplus
\K_{\infty}\alpha_r$ por lo que $\Gamma$ 
es una red de rango $r_{\Gamma}=r$.
\end{ejemplo}

\begin{observacion}\label{DrinfeldO2.2.6(1)}
Todas las redes $\Gamma$ de rango $r$ est\'an dadas por $\Gamma=
J_1\xi_1+\cdots+J_r\xi_r$ con $\{\xi_1,\ldots,\xi_r\}\subseteq \Ci$ un conjunto
linealmente independiente sobre $\K_{\infty}$ y $J_1,\ldots, J_r$ son
ideales fraccionarios de $A$.
\end{observacion}

\begin{observacion}\label{DrinfeldO2.2.7} El Ejemplo \ref{DrinfeldEj2.2.6} muestra la gran diferencia
con ${\ma C}$ pues las \'unicas redes tienen rango $1$ o $2$ sobre ${\ma R}$ que es
la completaci\'on arquimediana de ${\ma Q}$.
\end{observacion}

\subsubsection{M\'odulos de Drinfeld y redes}\label{DrinfeldS2.1.ModDrinfRed}

Para los preliminares consideramos $F$ un campo completo con respecto a una
valuaci\'on. Una serie $\sum_{j=0}^{\infty} a_j$ con $a_j\in F$ converge $\iff
\lim_{j\to\infty} v(a_j)=\infty\iff \lim_{j\to\infty}|a_j|=0$. De hecho tenemos
\[
\Big|\sum_{j=0}^ma_j-\sum_{j=0}^n a_j\Big|=\Big|\sum_{n+1}^ma_j\Big|\leq
\max_{n+1\leq j\leq m}|a_j|\xrightarrow[m\to\infty]{}0.
\]

En consecuencia si $f(x)=\sum_{i=0}^{\infty} a_ix^i$ 
es una serie de potencias, sea
$\alpha\in F$ y $f(\alpha)=\sum_{i=0}^{\infty}a_i\alpha^i$. Entonces
\begin{gather*}
v(a_i\alpha^i)=v(a_i)+iv(\alpha)\xrightarrow[i\to\infty]{}\infty\iff
 \text{\ para toda $M>0$
existe $i_0$}\\
\text{ tal que para toda\ } i\geq i_0,
v(a_i)+iv(\alpha)>M\\
\iff v(\alpha)>\frac{M}{i}-\frac{v(a_i)}{i} \text{\ para toda\ } i\geq i_0.
\end{gather*}

Sea ${\mathcal P}(f):=\limsup\limits_{i\to\infty}\Big( \dps 
-\frac{v(a_i)}{i}$\Big).\label{Drinfeldradioconvergencia}

\begin{proposicion}\label{DrinfeldP2.2.8} Sea $\alpha\in \K $. Entonces $f(x)$ converge
en $\alpha$ si $v(\alpha)>{\mathcal P}(f)$ y diverge si $v(\alpha)<{\mathcal P}
(f)$. Si $v(\alpha)=
{\mathcal P}(f)$ no podemos afirmar nada.
\end{proposicion}

\begin{proof} Si $v(\alpha)>{\mathcal P}
(f)$ entonces consideremos $\epsilon= \frac{v(\alpha)-{\mathcal P}(f)}{2}>0$.
Entonces existe $n_0$ tal que para toda $n\geq n_0$, 
$-\frac{v(a_n)}{n}<v(\alpha)-\epsilon$, esto es,
$v(\alpha)>\epsilon-\frac{v(a_n)}{n}$. Por
tanto $v(a_n\alpha^n)=nv(\alpha)+v(a_n)>n\epsilon\xrightarrow[n\to\infty]{}\infty$,
por lo tanto la serie converge.

Ahora si $v(\alpha)<{\mathcal P}
(f)$, sea $\delta=\frac{{\mathcal P}(f)-v(\alpha)}{2}>0$. Existe 
una subsucesi\'on $\{a_{n_k}\}_{k=1}^{\infty}$ tal que
existe $k_0$
tal que para toda $k\geq k_0$, $-\frac{v(a_{n_k})}{n_k}>v(\alpha)+\delta$.
Por tanto $v(a_{n_k})<n_k\big(-\delta-v(\alpha)\big)=
n_k\big(-\frac{{\mathcal P}(f)}{2}+\frac{v(\alpha)}{2}-v(\alpha)\big)=
n_k\big(-\frac{{\mathcal P}
(f)}{2}-\frac{v(\alpha)}{2}\big)< n_k\big(-\frac{v(\alpha)}{2}-\frac{
v(\alpha)}{2}\big)=-n_kv(\alpha)$.

Por tanto $v(a_{n_k})+n_kv(\alpha)=v(a_{n_k}\alpha^{n_k})<0$ 
de donde $|a_{n_k}\alpha^{n_k}|>1$
y en particular $a_{n_k}\alpha^{n_k}$ no 
converge a $0$ cuando $n_k\lra\infty$ y la serie diverge.
$\fin$
\end{proof}

\begin{ejemplo}\label{DrinfeldEj2.2.9} Sea $f(x)=\sum_{n=0}^{\infty}x^n$, $a_n=1$, por lo
que $v(a_n)=0$ y ${\mathcal P}(f)=0$. Por tanto converge para $v(\alpha)>0$, es decir,
cuando $|\alpha|_v<1$ y diverge para $v(\alpha)<0$, esto es, para $|\alpha|_v>1$.

$\sum_{n=0}^{\infty}\alpha^n$ diverge para $|\alpha|=1$ pues $|\alpha^n|=1
\lra 1\neq 0$.

Si $f(x)=\sum_{i=0}^{\infty} \pi^{-i}x^{i^2}$ donde $\pi$ es un elemento primo,
entonces 
\begin{gather*}
a_i=\begin{cases} 0&\text{si $i\neq j^2$ para toda $j\in {\ma N}$},\\
\pi^{-j}&\text{si $i=j^2$ para alg\'un $j\in {\ma N}$}.
\end{cases}
\intertext{Se sigue que}
\frac{v(a_i)}{i}=\begin{cases} \infty&
\text{si $i\neq j^2$ para toda $j\in {\ma N}$},\\
-\frac{1}{j}&\text{si $i=j^2$ para alg\'un $j\in {\ma N}$}.
\end{cases}
\end{gather*}

Por tanto $\limsup\limits_{i\to\infty}\big(-\frac{v(a_i)}{i}\big)=0$. Para
$\alpha=1$, $f(1)=\sum_{i=0}^{\infty}\pi^{-i}$ diverge pues $v(\pi^{-i})
=-i\xrightarrow[i\to\infty]{} -\infty$.
\end{ejemplo}

Como en el caso polinomial, podemos trabajar con pol{\'\i}gonos de Newton para
series formales $F[[x]]$ y usando estos pol{\'\i}gonos se puede probar que si $F$ es
completo y algebraicamente cerrado, entonces:

\begin{proposicion}\label{DrinfeldP2.2.10} Sean $\{m_i\}_i$ la sucesi\'on de pendientes
del pol{\'\i}gono de Newton de $f(x)=\sum_{j=0}^{\infty}a_jx^j$. 
Entonces la sucesi\'on
$\{m_i\}_i$ es creciente y se tiene $\lim_{i\to\infty}m_i =-{\mathcal P}(f)$.

\las
\item Para $f(x)=\sum_{j=0}^{\infty} a_jx^j\in L[[x]]$ y $\hat{\bar{L}}=F$, entonces las
ra{\'\i}ces de $f(x)$ son algebraicas sobre $L$.

\item Hay un n\'umero finito de ceros de $f(x)\in F[[x]]$ en la bola $\{x\mid v(x)\geq t\}$
donde $t>{\mathcal P}(f)$. En particular $f(x)$ tiene una cantidad a lo m\'as numerable de
ceros. 

\item El conjunto de ceros de $f(x)$ es un conjunto discreto.

\item Para $f(x)$ de la forma $f(x)=\sum_{j=0}^{\infty}a_jx^{p^j}$
con $a_0\neq 0$ se tiene que $f'(x)=a_0$ y las ra\'ices de
$f(x)$ son simples. 
$\fin$
\end{list}
\end{proposicion}

\begin{definicion}\label{DrinfeldD2.2.11} Se dice que {\em $f(x)$ es entera\index{funci\'on entera}}
si ${\mathcal P}(f)=-\infty$, es decir, si $f(x)$ converge para toda $x$.
\end{definicion}

Se tiene el siguiente resultado (comparar con el Teorema de Liouville cl\'asico):

\begin{lema}\label{DrinfeldL2.2.12} Si $f(x)$ es una funci\'on entera sin ceros, entonces
$f(x)$ es constante. Como consecuencia, si $f$ es entera y no constante, $f$ es
suprayectiva.
\end{lema}

\begin{proof} Para la \'ultima parte, si $f(x)$ es entera y no constante, entonces para
cualquier $c$, consideremos $g(x)=f(x)-c$.
$\fin$
\end{proof}

Como consecuencia del Lema \ref{DrinfeldL2.2.12} obtenemos el an\'alogo al teorema cl\'asico
de factorizaci\'on de Weierstrass\index{teorema de factorizaci\'on de Weierstrass}.

\begin{teorema}\label{DrinfeldT2.2.13} Sea $f\colon \Ci\lra\Ci$ una funci\'on entera
dada por $f(u)=
\sum_{n=0}^{\infty} a_n u^n\in \Ci[[u]]$. Sean $\{\lambda\}_{\lambda\in I}$ los
diferentes ceros de $f$ en $\Ci$ con $\lambda\neq 0$,
cada $\lambda$ con multiplicidad $m_{\lambda}$.
Entonces $I$ es a lo m\'as numerable y si $|I|=\infty$, $\{\lambda\}_{\lambda\in I}=
\{\lambda\}_{i=1}^{\infty}$, $\lim_{t\to\infty}\vi(\lambda_t)=-\infty$ y si $n$ es la
multiplicidad del cero de $f$ en $u=0$, se tiene la {\em expansi\'on de Euler de
$f$\index{expansi\'on de Euler}}:
\[
f(u)=cu^n\prod_{t=1}^{\infty}\big(1-\frac{u}{\lambda_t}\big)^{m_t}
\]
para alguna constante $c\in \Ci$ y donde $m_t:=m_{\lambda_t}$.

Rec{\'\i}procamente, todas estas funciones son enteras.
\end{teorema}

\begin{proof} \cite[Chapter 2, p\'aginas 38--42]{Gos96}. $\fin$
\end{proof}

La definici\'on de 
convergencia para productos infinitos es la siguiente: 
en general, en un campo completo, se
dice que $\prod_{n=0}^{\infty}a_n$ converge a $a$ si $\lim_{n\to\infty}\prod_{
i=0}^na_i=a$ y $a \neq 0$. En particular, si $\prod_{n=0}^{\infty} a_n$ converge,
entonces $\lim_{n\to\infty}a_n=1$. 

As{\'\i}, si $a_n=1+v_n$, $\prod_{n=0}^{\infty}a_n=
\prod_{n=0}^{\infty}(1+v_n)$ converge entonces $\lim_{n\to\infty}v_n=0$. En 
${\ma C}$, $\prod_{n=0}^{\infty}(1+v_n)$ converge $\iff \sum_{n=0}^{\infty}\ln
(1+v_n)$ converge. En caracter{\'\i}stica $p>0$, $\prod_{n=0}^{\infty}(1+v_n)$
converge $\iff \lim_{n\to\infty}v_n=0$.

Ahora sea $f$ una funci\'on entera con ceros simples en $X\subseteq \Ci$ y
$f(x)=x\prod_{\gamma\in X\setminus \{0\}}\big(1-\frac{x}{\gamma}\big)$. Si $f(x)$
es aditiva, esto es, $f(x+y)=f(x)+f(y)$, entonces $X$ es un subgrupo aditivo de
$\Ci$. 

Rec{\'\i}procamente supongamos que $X$ es un subgrupo de $\K $ tal que
$X\cap B(a,r)$ es finito para toda $a\in\Ci$ y para toda $r>0$, donde $B(a,r)=
\{x\in\Ci\mid|x-a|<r\}$. Sea 
\[
e_X(x):=\exponencial x{\gamma}X.\label{Drinfeldfuncionexponencial}
\]
Entonces
$\lim_{\gamma\in X\setminus\{0\}}x/\gamma=0$ para toda $x\in\Ci$ puesto
que dado $M>0$ \'unicamente un n\'umero finito de $\gamma\in X$ satisfacen
que $|\gamma|\leq M$, por tanto $|\gamma|>M$ para casi toda $\gamma$,
por lo que $|1/\gamma|<1/M$ para casi toda $\gamma$. Se sigue que $e_X$
es una funci\'on entera.

\begin{proposicion}\label{DrinfeldP2.2.15} 
Si $Y\subseteq \Ci$ es un $\F$--espacio vectorial
de dimensi\'on finita, entonces $f_Y(u)=c\prod_{\gamma\in Y}
(u-\gamma)$ es un $\F$--polinomio lineal aditivo.
\end{proposicion}

\begin{proof} Probaremos el resultado por 
inducci\'on en $\dim_{\ma F} Y=n$ con ${\ma F}=\F$.
Para $n=0$, $Y=\{0\}$, $f_Y(u)=cu$ satisface lo pedido.

Sea $n\geq 1$ y sea $W\subseteq Y$ un subespacio de dimensi\'on $n-1$.
Para $\gamma_0\in Y\setminus W$, $Y=W+{\ma F}\gamma_0$, por lo tanto
$f_Y(u)=c\prod_{\substack{w\in W\\ \mu\in{\ma F}}}(u-(w+\mu \gamma_0))$.

Sea $f_W(u)=\prod_{w\in W}(u-w)$ el cual es un polinomio ${\ma F}$--lineal y
\begin{align*}
f_Y(u)&=c\prod_{\substack{w\in W,\\ \mu\in {\ma F}}}(u-(w+\mu\gamma_0))=
c\prod_{\mu\in{\ma F}}\prod_{w\in W}((u-\mu\gamma_0)-w)\\
&=c\prod_{\mu\in{\ma F}}f_W(u-\mu\gamma_0)=c\prod_{\mu\in{\ma F}}
(f_W(u)-\mu f_W(\gamma_0))\\
&\igual_{\substack{\uparrow\\ z=\frac{f_W(u)}{f_W(\gamma_0)}}}
cf_W(\gamma_0)^{|{\ma F}|}\cdot \prod_{\mu\in{\ma F}}(z-\mu)=
cf_W(\gamma_0)^{|{\ma F}|}(z^{|{\ma F}|}-z)\\
&= cf_W(u)^{|{\ma F}|}-cf_W(\gamma_0)^{|{\ma F}|-1}f_W(u).
\end{align*}

Se sigue que $f_Y(u)=cf_Y(u)^{|{\ma F}|}-cf_W(\gamma_0)^{|{\ma F}|-1}
f_W(u)$ el cual es ${\ma F}$--lineal, por lo que
\[
f_Y(\alpha x+\beta y)=\alpha f_Y(x)+\beta f_Y(y)\quad \text{para cualesquiera}\quad
x,y\in \Ci,\quad \alpha,\beta\in{\ma F}.
\]

Se sigue que $f_Y$ es ${\ma F}$--lineal. $\fin$
\end{proof}

\begin{proposicion}\label{DrinfeldP2.2.14} Sea $X$ un subgrupo discreto como
antes. Se tiene que $e_X(x+y)=e_X(x)+e_X(y)$
para cualesquiera $x,y\in \Ci$, esto es, $e_X$ es una funci\'on aditiva.
M\'as a\'un, $e_X(u)$ es $\F$--lineal.
\end{proposicion}

\begin{proof} Sea $N\in{\ma N}$ y sea $X_N:=\{\lambda\in X\mid |\lambda|_{\infty}\leq
N\}$. Entonces $X_N$ es un subgrupo finito de $X$ pues si $\lambda,\mu\in X_N$
entonces $|\lambda+\mu|_{\infty}\leq \max\{|\lambda|_{\infty},|\mu|_{\infty}\}\leq
N$.

Sea $p_N(u):= \exponencial u{\gamma}{X_N}\in \Ci[u]$. Para toda $u\in\Ci$, se
tiene $\lim_{N\to \infty} P_N(u)=e_X(u)$. 

Puesto que, $\lim_{\gamma\to\infty}\frac{u}{\gamma}=0$, $e_X(u)$ existe
donde $e_X(u)=\lim_{n\to\infty} P_N(u)$.
Ahora, para probar que $e_X(x+y)=e_X(x)+e_X(y)$,
basta probar que $P_N(x)$ es aditivo pues en
este caso se tiene se tiene para toda $N\in{\ma N}$
que 
\begin{align*}
e_X(x+y)&=\lim_{N\to\infty}P_N(x+y)=
\lim_{N\to\infty}(P_N(x)+P_N(y))\\
&=\lim_{N\to\infty}P_N(x)+\lim_{N\to\infty}P_N(y)
=e_X(x)+e_X(y).
\end{align*}

El resultado se sigue de la Proposici\'on \ref{DrinfeldP2.2.15}. $\fin$
\end{proof}

\begin{definicion}\label{DrinfeldD2.2.16}
La funci\'on $e_X$ se llama la {\em funci\'on exponencial\index{funci\'on exponencial}}
de $X$ y se tiene 
\[
e_X(x)=\exponencial x{\gamma}X.
\]
\end{definicion}

Si $X$ es una $A$--red, entonces si ponemos $\Gamma=X$, la funci\'on exponencial
de $\Gamma$, $e_{\Gamma}(u)=\exponencial u{\gamma}{\Gamma}$ la cual es
una funci\'on entera con ceros exactamente en $\Gamma$
los cuales son simples y $e_{\Gamma}$ es
$\F$--lineal.

Para definir el m\'odulo de Drinfeld $\rho^{\Gamma}$ asociado a una $A$--red $\Gamma$,
primero consideremos dos redes $\Gamma$, $\Gamma^{\prime}$ tales que
$\Gamma\subseteq \Gamma^{\prime}$ de tal forma que $\Gamma$ es de 
{\'\i}ndice finito, esto es, $\big|\Gamma^{\prime}/\Gamma\big|<\infty$.

Consideremos $e_{\Gamma}(u)$. Se tiene que $e_{\Gamma}(u+\gamma)=e_{\Gamma}
(u)$ para toda $\gamma\in \Gamma$. Adem\'as, si $\gamma$ es cualquier per{\'\i}odo
de $e_{\Gamma}$, es decir, $e_{\Gamma}(u+\gamma)=e_{\Gamma}(u)$, entonces,
$\gamma\in \Gamma$ pues $e_{\Gamma}(u+\gamma)=e_{\Gamma}(u)+e_{\Gamma}
(\gamma)=e_{\Gamma}(u)$ por lo que $e_{\Gamma}(\gamma)=0$ y $\gamma\in
\Gamma$.

Sea $\varphi\colon \Gamma^{\prime}\lra e_{\Gamma}(\Gamma^{\prime})=
\{e_{\Gamma}(\gamma)\mid \gamma\in \Gamma^{\prime}\}$, $\xi\longmapsto
e_{\Gamma}(\xi)$.

Entonces $\varphi$ es un epimorfismo de $\F$--espacios vectoriales y $\ker\varphi=
\Gamma$ por lo que $\Gamma^{\prime}/\Gamma\cong e_{\Gamma}(\Gamma^{\prime})$
como $\F$--espacios vectoriales.

\begin{definicion}\label{DrinfeldD2.2.17} Sean $\Gamma,\Gamma^{\prime}$ dos redes
tales que $\Gamma\subseteq \Gamma^{\prime}$ y es de {\'\i}ndice finito. Se
define
\[
P(\Gamma^{\prime}/\Gamma, u):=\exponencial u{\lambda}
{e_{\Gamma}(\Gamma^{\prime})}.\label{DrinfeldpolinomioDrinfeld}
\]

Entonces $P(\Gamma^{\prime}/\Gamma, u)$ es un polinomio $\F$--lineal de
grado $\big|\Gamma^{\prime}/\Gamma\big|=\big|e_{\Gamma}(\Gamma')\big|$.
\end{definicion}

\begin{proposicion}\label{DrinfeldP2.2.18} Sean $\Gamma\subseteq \Gamma^{\prime}
\subseteq \Gamma^{\prime\prime}$ tres redes con $[\Gamma^{\prime\prime}:\Gamma]
<\infty$. Entonces:
\lasa
\item $e_{\Gamma^{\prime}}(u)=P(\Gamma^{\prime\prime}/\Gamma,e_{\Gamma}
(u))$ para $u\in \Ci$.
\item $P(\Gamma^{\prime\prime}/\Gamma,u)=P(\Gamma^{\prime\prime}/\Gamma^{\prime},
P(\Gamma^{\prime}/\Gamma,u))$.
\end{list}
\end{proposicion}

\begin{proof} El lado izquierdo de (a) tiene como ra{\'\i}ces a los elementos $\lambda\in 
\Gamma^{\prime}$ y del lado derecho de (a), las ra\'ices
son los elementos $\mu$ tales que $e_{
\Gamma}(u)\in e_{\Gamma}(\Gamma^{\prime})$, es decir, los elementos $u\in
\Gamma^{\prime}$. Por lo tanto ambas funciones son enteras con los mismos ceros.
Se sigue que existe una constante $c$ tal que $e_{\Gamma^{\prime}}(u)=
cP(\Gamma^{\prime}/\Gamma,e_{\Gamma}(u))$. Se obtiene que $c=1$ 
observando que $\frac{e_{\Gamma^{\prime}}(u)}{u}$ y $\frac{P(\Gamma^{\prime}/
\Gamma,e_{\Gamma}(u))}{u}$ coinciden al ser
evaluadas en $u=0$: ambas son igual a $1$.

De hecho, $f_{\Gamma'}(u)=\frac{e_{\Gamma'}(u)}{u}=\prod_{\substack{
\lambda\in e_{\Gamma}(\Gamma')\\ \lambda\neq 0}}\big(1-\frac{u}{\lambda}
\big)$. Por tanto
$f_{\Gamma'}(0)=\prod_{\substack{
\lambda\in e_{\Gamma}(\Gamma')\\ \lambda\neq 0}}1=1$. Por otro
lado, sea 
\[
g(u)=\frac{P(\Gamma'/\Gamma, e_{\Gamma}(u))}{u}=
\frac{e_{\Gamma}(u)}{u}\cdot \prod_{\substack{
\lambda\in e_{\Gamma}(\Gamma')\\ \lambda\neq 0}}\big(1-\frac{
e_{\Gamma}(u)}{\lambda}\big), 
\]
por lo que 
$g(0)=f_{\Gamma}(0)\cdot \prod_{\substack{
\lambda\in e_{\Gamma}(\Gamma')\\ \lambda\neq 0}}1=1$.

Para probar (b), vemos que $P(\Gamma^{\prime\prime}/
\Gamma,u)$ es un polinomio de grado
$\big|\Gamma^{\prime\prime}/\Gamma\big|$ cuyos ceros son precisamente 
los elementos de
$e_{\Gamma}(\Gamma^{\prime\prime})$. Ahora bien,
se tiene que $P(\Gamma^{\prime\prime}/
\Gamma^{\prime},P(\Gamma^{\prime}/\Gamma,u))$ es un polinomio de grado
$\big|\Gamma^{\prime\prime}/\Gamma^{\prime}\big|\cdot\big|\Gamma^{\prime}/
\Gamma\big|=\big|\Gamma^{\prime\prime}/\Gamma\big|$. Ambos polinomios
son m\'onicos.

Ahora veamos que $e_{\Gamma}(\Gamma^{\prime\prime})$ son los ceros de $f(u)=
P(\Gamma^{\prime\prime}/\Gamma^{\prime},P(\Gamma^{\prime}/\Gamma,u))$.
Sea $\lambda\in e_{\Gamma}(\Gamma^{\prime\prime})$, digamos $\lambda=
e_{\Gamma}(\beta)$, $\beta\in \Gamma^{\prime\prime}$. Entonces
\begin{align*}
f(\lambda)&=P(\Gamma^{\prime\prime}/\Gamma^{\prime},P(\Gamma^{\prime}/
\Gamma,\lambda))=P(\Gamma^{\prime\prime}/\Gamma^{\prime},P(\Gamma^{\prime}/
\Gamma,e_{\Gamma}(\beta))\\
&\underbracket[0pt]{=}_{\substack{\uparrow\\(a)}}P(\Gamma^{\prime\prime}/
\Gamma^{\prime},e_{\Gamma^{\prime}}(\beta))\underbracket[0pt]{=}_{\substack{
\uparrow\\ (a)}}e_{\Gamma^{\prime\prime}}(\beta)=0. \tag*{$\fin$}
\end{align*}
\end{proof}

Ahora veremos que dada una red $\Gamma$ de rango $r\geq 1$ en $\Ci$, 
obtenemos un m\'odulo de Drinfeld de rango $r$ asociado a $\Gamma$.

Sea $\Gamma$ una red y sea $a\in A\setminus\{0\}$. Entonces $\Gamma\subseteq
a^{-1}\Gamma$.

\begin{teorema}\label{DrinfeldT2.2.19} Sea $\Gamma$ una red de rango $r$. Para $a\in
A\setminus\{0\}$, definimos
\[
\rho_a^{\Gamma}\colon \Ci\lra\Ci \e\text{por}\e \rho_a^{\Gamma}(u):=
a P(a^{-1}\Gamma/\Gamma,u).
\]

Entonces $\rho_a^{\Gamma}\in\aditivo {\Ci}$. Sea $\rho^{\Gamma}\colon
A\to \aditivo {\Ci}$ definido por $\rho^{\Gamma}(a):=\rho_a^{\Gamma}$ si $a\neq 0$
y $\rho^{\Gamma}(0)=0$. Entonces $\rho^{\Gamma}$ es un $A$--m\'odulo
de Drinfeld de rango $r$ sobre $\Ci$.
\end{teorema}

\begin{proof} Como $P(a^{-1}\Gamma/\Gamma,u)$ es $\F$--lineal, se sigue que
$\rho_a^{\Gamma}\in\aditivo {\Ci}$. Ahora
\[
\rho_a^{\Gamma}(u)=aP(a^{-1}\Gamma/\Gamma,u)=a\exponencial u{\lambda}{
e_{\Gamma}(a^{-1}\Gamma)}.
\]

Por tanto $D(\rho_a^{\Gamma})=a=\delta(a)$ donde $\delta\colon A\lra \Ci$ es el
encaje natural. En particular $D\circ\rho^{\Gamma}=\delta$.

Ahora 
\begin{align*}
e_{a^{-1}\Gamma}(u)&=\exponencial u{\lambda}{a^{-1}\Gamma}
=u\prod_{\mu\in\Gamma\setminus\{0\}}\big(1-\frac{u}{a^{-1}\mu}\big)\\
&=a^{-1}\exponencial {au}{\mu}{\Gamma}=a^{-1}e_{\Gamma}(au).
\end{align*}
Por tanto
\begin{equation}\label{DrinfeldEq2.2.20}
\text{\fbox{$e_{a^{-1}\Gamma}(u)=a^{-1}e_{\Gamma}(au)$}}.
\end{equation}

Se sigue que
$e_{\Gamma}(au)=ae_{a^{-1}\Gamma}(u)=aP(a^{-1}\Gamma/\Gamma,
e_{\Gamma}(u))=\rho_a^{\Gamma}(e_{\Gamma}(u))$
(Proposici\'on \ref{DrinfeldP2.2.18} (a)).
Por tanto tenemos multiplicaci\'on compleja
\begin{equation}\label{DrinfeldEq2.2.21}
\text{\fbox{$e_{\Gamma}(au)=\rho_a^{\Gamma}(e_{\Gamma}(u))$}}.
\end{equation}

Entonces
\begin{gather*}
\rho_{a+b}^{\Gamma}(e_{\Gamma}(u))=e_{\Gamma}((a+b)u)=e_{\Gamma}
(au)+e_{\Gamma}(bu)=\rho_a^{\Gamma}(e_{\Gamma}(u))+\rho_b^{\Gamma}
(e_{\Gamma}(u));\\
\rho_{ab}^{\Gamma}(e_{\Gamma}(u))=e_{\Gamma}((ab)u)=\rho_a^{
\Gamma}(e_{\Gamma}(bu))=\rho_a^{\Gamma}\rho_b^{\Gamma}(e_{\Gamma}(u)).
\end{gather*}

Puesto que el mapeo exponencial es suprayectivo, se sigue que
\begin{gather*}
\rho_{a+b}^{\Gamma}=\rho_a^{\Gamma}+\rho_b^{\Gamma}\e\text{y}\e
\rho_{ab}^{\Gamma}=\rho_a^{\Gamma}\rho_b^{\Gamma}\e\text{para cualesquiera}\e
a,b\in A.
\end{gather*}

Veamos que el rango de $\rho^{\Gamma}$ es $r$. De esto se seguir\'a que $\rho$
es un m\'odulo de Drinfeld, es decir, $\rho_a^{\Gamma}\neq a$ para alg\'un $a\in A$
pues $r\geq 1$.

Se tiene $\rho_a^{\Gamma}(u)=aP(a^{-1}\Gamma/\Gamma, u)$, por lo tanto
$\deg_u \rho_a^{\Gamma}(u)=\big|a^{-1}\Gamma/\Gamma\big|$.

Ahora bien, puesto que $\Gamma$ es de rango $r$, $\Gamma=
I_1\xi_1\oplus\cdots\oplus I_r\xi_r$
con $I_1,\ldots,I_r$ ideales fraccionarios (Observaci\'on \ref{DrinfeldO2.2.6(1)})
la cual es la suma de $r$ ideales fraccionarios. M\'as precisamente, $\Gamma$
es un $A$--m\'odulo finitamente generado libre de torsi\'on y puesto que $A$ es
dominio Dedekind, $\Gamma\cong A^{r-1}\oplus I$ como $A$--m\'odulos con $I$
un ideal fraccionario (Teorema \ref{DrinfeldT1.3.1}).

Ahora pongamos $\Gamma=\oplus_{i=1}^r J_i$ con $J_i$ un ideal fraccionario.
Afirmamos que si $J$ es cualquier ideal fraccionario, entonces si $a\in A\setminus\{0\}$,
entonces $a^{-1}J/J\cong a^{-1}A/A\cong A/aA$.

El isomorfismo $a^{-1}A/A\cong A/aA$ se obtiene al multiplicar por $a$. M\'as generalmente,
si ${\eu L}$ es ideal fraccionario y $x\in\*\K $ donde $\K =\coc A$, entonces
$a^{-1}(x{\eu L}/)x{\eu L}\cong a^{-1}{\eu L}/{\eu L}$ lo cual se obtiene de
\[
\varphi\colon a^{-1}{\eu L}\lra a^{-1}x{\eu L}\lra a^{-1}x{\eu L}/x{\eu L}, 
\quad a^{-1}l\longmapsto
a^{-1}xl\longmapsto \overline{a^{-1}xl}.
\]
Entonces $\varphi$ es epimorfismo de
$A$--m\'odulos y $\ker \varphi={\eu L}$.

Para probar que $a^{-1}J/J\cong a^{-1}A/A$ primero veamos que existe $x\in
\*\K $ tal que $xJ\subseteq A$ y $\mcd(xJ,aA)=1$. Para ver esto, sean
\[
aA={\mathcal P}_1^{\alpha_1}\cdots {\mathcal P}_s^{\alpha_s},\e \alpha_i\geq 0\e
\text{y}\e J={\mathcal P}_1^{\beta_1}\cdots {\mathcal P}_s^{\beta_s} J_1
\]
con $J_1$ un ideal fraccionario
primo relativo a ${\mathcal P}_1,\ldots,{\mathcal P}_s$ y $\beta_i\in
{\ma Z}$, $1\leq i\leq s$. 
Supongamos $\beta_1<0,\ldots,\beta_t<0, \beta_{t+1}\geq 0,\ldots,
\beta_s\geq 0$. Sean $x_i\in {\mathcal P}_i^{|\beta_1|}\setminus {\mathcal P}_i^{
|\beta_i|+1}$, $1\leq i\leq t$ y $y_j\in {\mathcal P}_j^{\beta_j}\setminus
{\mathcal P}_j^{\beta_j+1}$, $t+1\leq j\leq s$. Por el teorema chino del residuo,
existen $x,y\in A$ tales que $x\equiv x_i\bmod {\mathcal P}_i^{|\beta_i|+1}$,
$1\leq i\leq t$; $y\equiv y_j\bmod {\mathcal P}_j^{\beta_j+1}$, $t+1\leq j\leq s$.

Entonces $xy^{-1}J=J_2$ es un ideal fraccionario primo relativo a $aA$. Ahora sea 
\[
J_2={\eu q}_1^{\gamma_1}\cdots {\eu q}_n^{\gamma_n} \pK_1^{\varepsilon_1}
\cdots \pK_m^{\varepsilon_m}\e\text{con}\e \gamma_i>0, \varepsilon_j<0.
\]
Sean $z_l\in \pK_l^{|\varepsilon_l|}\setminus \pK_l^{|\varepsilon_l|+1}$,
$1\leq l\leq m$ y sea
$z\in A$ tal que 
\[
z\equiv z_l\bmod \pK_l^{|\varepsilon_l|+1},\quad 1\leq l\leq m\quad
\text{y}\quad z\equiv 1\bmod {\mathcal P}_i, \quad 1\leq i\leq s.
\]
En particular $z\notin {\mathcal P}_i$.

Entonces $zJ_2=J_3\subseteq A$ y es primo relativo a $aA$.

En resumen, sea $w\in \*\K $ tal que $I=wJ\subseteq A$ y $\mcd(I,aA)=1$.
Entonces $a^{-1}I/I\cong a^{-1}(wJ)/(wJ)\cong a^{-1}J/J$.

Finalmente, sea
\[
\xymatrix{
a^{-1}I\ar@{^{(}->}[r]^i\ar@/_2pc/^{\psi}[rr]&a^{-1}A\ar@{>>}[r]^{\pi\ \ }
&a^{-1}A/A}, \quad a^{-1}x \stackrel{\psi}{\longmapsto} a^{-1}x\bmod A.
\]
Veamos que $\psi$ es un epimorfismo. Sean $x\in A$ y $\overline{a^{-1}x}\in
a^{-1}A/A$. Se quiere encontrar $y\in I$ tal que $\overline{a^{-1}y}=
\overline{a^{-1}x}$, esto es, $a^{-1}x-a^{-1}y\in A$, lo cual implicar\'a que
$x\in y+aA$. Puesto que $\mcd(I,aA)=1$, se tiene que $I+aA=A$ de donde
se sigue que existe $y\in I$ tal que $x=y+aA$ y por tanto obtenemos
que $\psi$ es suprayectiva con n\'ucleo $I$. Por tanto 
\fbox{$a^{-1}I/I\cong a^{-1}A/A$}.

Se sigue que
\begin{align*}
\big|a^{-1}\Gamma/\Gamma\big|&=\Big|\bigoplus_{i=1}^r a^{-1}J_i/J_i\Big|=
\prod_{i=1}^r\big|a^{-1}J_i/J_i\big|=\prod_{i=1}^r\big|a^{-1}A/A\big|\\
&=\big|A/aA\big|^r=q^{r\deg a}=\deg_{u}\rho_a^{\Gamma}(u)
=q^{\deg_{\tau}\rho_a^{\Gamma}}.
\end{align*}

Por tanto \fbox{$\rho^{\Gamma}$ es de rango de $r$}. $\fin$
\end{proof}

\begin{observacion}\label{DrinfeldO2.2.19(1)}
Dada una red $\Gamma$, el $A$--m\'odulo de Drinfeld $\rho^{\Gamma}$
admite {\em multiplicaci\'on compleja\index{multiplicaci\'on compleja}}:
\[
e_{\Gamma}(au)=\rho_a^{\Gamma}(e_{\Gamma}(u)),
\]
es decir, $e_{\Gamma}a=\rho_a^{\Gamma}e_{\Gamma}$ o, equivalentemente,
$e_{\Gamma}ae_{\Gamma}^{-1}=\rho_a^{\Gamma}$. Por tanto $\rho^{\Gamma}$
es \'unico y can\'onicamente definido.
\end{observacion}

El Teorema \ref{DrinfeldT2.2.19} muestra que existen $A$--m\'odulos de Drinfeld $\Drin_A(\Ci)$
sobre $\Ci$ para cualquier $A$ y para cualquier $r\geq 1$. Estos 
son los m\'odulos $\rho^{\Gamma}$
con $\Gamma$ una $A$--red. Este resultado, como hemos mencionado antes, no es
un hecho obvio y de hecho pueden no existir m\'odulos de Drinfeld $A\lra\aditivo F$
para algunos $F$ (ver Ejemplo \ref{DrinfeldEj2.2.0}).

El rec{\'\i}proco  del resultado anterior es cierto y ya lo mencionamos antes. Primero
demostramos el siguiente resultado el cual es central para el 
Teorema de Uniformizaci\'on Anal\'itica y para muchos otros resultados.

\begin{teorema}\label{DrinfeldNT1}
Sea $\rho$ un $A$--m\'odulo de Drinfeld sobre $\Ci$. Entonces, si $\rho\colon
A\lra \aditivo F$ con un campo $F$ conteniendo a $\K =\coc A$, existe una
\'unica serie $\xi_{\rho}\in F\torcido$ tal que
\begin{gather}
D\xi_{\rho}=1\quad\text{y}\quad \xi_{\rho}a\tau^0=\rho_a\xi_{\rho}\quad\text{para
toda $a\in A$},\nonumber
\intertext{es decir}
\xi_{\rho}a\xi_{\rho}^{-1}=\rho_a\quad\text{para toda $a\in A$}.\label{DrinfeldEq2.2.24}
\end{gather}
\end{teorema}

\begin{proof} Ya vimos que $\rho$ se puede extender a un $\K $--m\'odulo formal
 $\rho\colon \K \lra F\torcido$.
Para probar (\ref{DrinfeldEq2.2.24}), primero consideremos el
caso en que $\rho$ trivial, es decir, no es
un m\'odulo de Drinfeld, esto es, $\rho(a)=a$ para toda $a\in \K $. Entonces definimos
$\xi_{\rho}=\tau^0=1$. Si $\xi'_{\rho}$ fuese otra serie tal que $\xi'_{\rho}
a=\rho_a \xi'_{\rho}=a\xi'_{\rho}$. Expandiendo se ve inmediatamente que
$\xi'_{\rho}\in \K $ y puesto que $D\xi'_{\rho}=1$, se sigue que $\xi'_{\rho}
=\xi_{\rho}$.

Ahora consideremos $\rho$ no trivial, es decir, un verdadero
m\'odulo de Drinfeld. Entonces existe $\alpha\in \K $ tal que $\rho_{\alpha}\neq
\alpha\tau^0$.

Afirmamos que $\alpha$ es trascendente sobre $\F$ pues en caso de ser algebraico
se tendr\'ia una relaci\'on
\begin{gather*}
\alpha^n+a_{n-1}\alpha^{n-1}+\cdots+a_1\alpha+a_0=0,\e a_i\in\F\\
\intertext{y por tanto}
0=\rho_0=\rho_{\alpha}^n+a_{n-1}\rho_{\alpha}^{n-1}+\cdots+a_1\rho_{\alpha}
+a\rho_1\e\text{pues}\e \rho_{\alpha}^0=\rho_1.
\end{gather*}
Por tanto se tendr{\'\i}a que $\rho_{\alpha}$ 
es algebraico sobre $\F$ de donde se seguir\'ia que $\rho_{
\alpha}=\alpha\tau^0$ lo cual es contradice la naturaleza de $\alpha$.

Ahora veamos que existe una \'unica serie de potencias $\lambda_{\rho_{\alpha}}=
\sum_{i=0}^{\infty}c_i\tau^i\in F\torcido$ con $c_0=1$ y $\lambda_{\rho_{\alpha}}
\alpha\tau^0=\rho_{\alpha}\lambda_{\rho_{\alpha}}$ con la multiplicaci\'on de
$F\torcido$.

Para demostrar la afirmaci\'on anterior, escribamos $\rho_{\alpha}=
\sum_{j=0}^{\infty}a_j\tau^j$ con $a_0=\alpha$ 
y la suma es finita. Entonces
\begin{gather*}
\begin{align*}
\lambda_{\rho_{\alpha}}\alpha\tau^0&=\Big(\sum_{i=0}^{\infty} c_i\tau^i\Big)
\alpha\tau^0
=\sum_{i=0}^{\infty}c_i\alpha^{q^i}\tau^i;\\
\rho_{\alpha}\lambda_{\rho_{\alpha}}&=\Big(\sum_{j=0}^{\infty}a_j\tau^j\Big)
\Big(\sum_{i=0}^{\infty}c_i\tau^i\Big)\\
&=\sum_{i=0}^{\infty}\Big(\sum_{j=0}^i a_j\tau^j
c_{i-j}\tau^{i-j}\Big)=\sum_{i=0}^{\infty}\Big(\sum_{j=0}^i a_jc_{i-j}^{q^i}\Big)\tau^i.
\end{align*}
\intertext{Por tanto}
\lambda_{\rho_{\alpha}}\alpha\tau^0=\rho_{\alpha}\lambda_{\rho_{\alpha}}
\iff c_i \alpha^{q^i}=\sum_{j=0}^i a_j c_{i-j}^{q^j},\e i=0,1,\ldots\\
\iff \big(\alpha^{q^i}-\underbracket[0pt]{\alpha}_{\substack{\uigual\\ a_0}}\big)c_i=
\sum_{j=1}^ia_j c_{i-j}^{q^j}, \e i=1,2,\ldots
\end{gather*}

Puesto que $\alpha$ es transcendente sobre $\F$, se tiene que
$\alpha^{q^i}-\alpha\neq 0$ para toda $i\geq 1$. Por lo tanto
\[
c_i=\frac{1}{\alpha^{q^i}-\alpha}\sum_{j=1}^i a_jc_{i-j}^{q^j}=
\frac{1}{\alpha^{q^i}-\alpha}\sum_{j=0}^{i-1}a_{i-j}c_j^{q^{i-j}},
\]
por lo que $c_0=1$ y los $c_i$'s, $i\geq 1$ son \'unicos.

Como consecuencia notemos que $\lambda_{\rho_a}F\lambda_{
\rho_a}^{-1}$ es el centralizador de $\rho_a$ en $F\torcido$. De hecho, en
$F\torcido$, el centralizador de $\alpha\tau^0$ es $F\tau^0=F$ pues si
$\big(\sum_{i=0}^{\infty}d_i\tau^i\big)(\alpha\tau^0)=(\alpha\tau^0)\big(
\sum_{i=0}^{\infty}d_i\tau^i\big)$ entonces $d_i\alpha^{q^i}=d_i\alpha$
para toda $i=1,2,\ldots$. Obtenemos que $(\alpha^{q^i}-\alpha)d_i=0$
para toda $i\geq 1$. Se sigue que $d_i=0$ para toda $i\geq 1$.

Ahora, $\lambda_{\rho_{\alpha}}(\alpha\tau^0)\lambda_{\rho_{\alpha}}^{-1}=
\rho_{\alpha}$. Por tanto si $\mu\rho_{\alpha}=\rho_{\alpha}\mu$, entonces
$\mu \lambda_{\rho_{\alpha}}\alpha\lambda_{\rho_{\alpha}}^{-1}=
\lambda_{\rho_{\alpha}}\alpha\lambda_{\rho_{\alpha}}^{-1}\mu$. De esta
forma obtenemos
\[
(\lambda_{\rho_{\alpha}}^{-1}\mu\lambda_{\rho_{\alpha}})(\alpha\tau^0)=
(\alpha \tau^0)(\lambda_{\rho_{\alpha}}^{-1}\mu \lambda_{\rho_{\alpha}})
\Lra \lambda_{\rho_{\alpha}}^{-1}\mu\lambda_{\rho_{\alpha}}\in F\Lra
\mu\in \lambda_{\rho_{\alpha}}F\lambda_{\rho_{\alpha}}^{-1}.
\]

Seguimos con la construcci\'on de $\xi_{\rho}$ que satisfaga $\xi_{\rho}
\alpha\tau^0=\rho_{\alpha}\xi_{\rho}$. Recordemos que estamos en el
caso $\rho_{\alpha}\neq\alpha\tau^0$. Sea $\xi_{\rho}:=\lambda_{\rho_{\alpha}}$.

Sea $x\in \K $, entonces $\rho_x\rho_{\alpha}=\rho_{\alpha}\rho_x$. Se sigue
que $\rho_x\in \lambda_{\rho_{\alpha}}F\lambda_{\rho_{\alpha}}^{-1}$, es
decir, existe $t\in F$ tal que $\rho_x=\lambda_{\rho_{\alpha}} t\lambda_{
\rho_{\alpha}}^{-1}$. Igualando al coeficiente constante, tenemos $t=x$,
por lo que $\xi_{\rho}x\tau^0=\lambda_{\rho_{\alpha}}x\tau^0=\rho_x
\lambda_{\rho_{\alpha}}=\rho_x \xi_{\rho}$. Es decir, tal $\xi_{\rho}$
existe.

Veamos que $\xi_{\rho}$ es \'unico. Sea 
$\xi_{\rho}^{\prime}$ satisface las mismas condiciones
que $\xi_{\rho}$. De la unicidad de $\lambda_{\rho_{\alpha}}$ se sigue que
$\xi_{\rho}^{\prime}=\xi_{\rho}=\lambda_{\rho_{\alpha}}$. $\fin$
\end{proof}

\begin{observacion}\label{DrinfeldNO3}
La Ecuaci\'on (\ref{DrinfeldEq2.2.24}) es fundamental para la demostraci\'on
del  Teorema \ref{DrinfeldT2.2.22} y para muchos otros resultados.
Adem\'as se tiene que $\xi_{\rho}$ es
la funci\'on exponencial de $\rho$ y los ceros de $\xi_{\rho}$
forman la red relacionada con el m\'odulo de Drinfeld.
\end{observacion}

\begin{proposicion}\label{DrinfeldNP2} La serie $\xi_{\rho}$ dada
por el Teorema {\rm{\ref{DrinfeldNT1}}} es una funci\'on entera.
\end{proposicion}

\begin{proof}
Sea $\xi_{\rho}=\sum_{i=0}^{\infty}c_i\tau^i$ con $c_0=1$. Se tiene que
\begin{gather*}
\xi_{\rho}\text{\ es una funci\'on entera\ } \iff \sum_{i=0}^{\infty}c_ix^{q^i}=
\sum_{i=0}^{\infty}d_{q^i} x^{q^i}=f \\
\text{es una funci\'on entera\ }\iff {\mathcal P}(f)
=-\lim_{i\to\infty} \frac{\vi(d_{q^i})}{q^i}=-\infty\\
\iff \lim_{i\to\infty}\frac{\vi(c_i)}
{q^i}=\lim_{i\to\infty}\vi(c_i^{q^{-i}})=\infty\iff \lim_{i\to\infty}|c_i|^{1/q^i}=0.
\end{gather*}

Sea $a\in A$ con $\deg a>0$ por lo que $a$ es transcendente sobre $\F$, con
el mismo argumento que usamos para $\alpha$. Sea
$\rho_a=a\tau^0+\sum_{i=1}^t a_i\tau^i$ y sea $n\geq t$.
Puesto que $\xi_{\rho}a\tau^0=\rho_a\xi_{\rho}$, usando las recursiones halladas
anteriormente, se tiene
\[
(a^{q^n}-a)c_n=\sum_{j=1}^t a_jc_{n-j}^{q^j}.
\]

Entonces, puesto que $\vi(a)<0$, $\vi(a^{q^n}-a)=q^n\vi(a)$ y 
\begin{gather*}
q^n\vi(a)+\vi(c_n)\geq \min_{1\leq j\leq t}\big\{\vi(a_j)+q^{j}\vi(c_{n-j})\big\},
\intertext{lo cual equivale a}
\frac{\vi(c_n)}{q^n}\geq \min_{1\leq j\leq t}\big\{\frac{\vi(a_j)}{q^n}+q^{j-n}\vi(c_{n-j})
\big\}-\vi(a).
\end{gather*}

Sea $\theta$ tal que $\theta<\vi(a)<0$. Para $n\gg 0$, digamos $n\geq n_0$,
se tiene $\min_{1\leq j\leq t}\big\{\frac{\vi(a_j)}{q^n}\big\}<\vi(a)-\theta$. Por tanto
$\frac{\vi(c_n)}{q^n}\geq \min_{1\leq j\leq t}\big\{\frac{\vi(c_{n-j})}{q^{n-j}}\big\}-\theta$.

Por tanto, por recursi\'on, tenemos que $\frac{\vi(c_n)}{q^n}\xrightarrow[n\to\infty]{}
\infty$, de donde se sigue que $\xi_{\rho}$ es una funci\'on entera. $\fin$
\end{proof}

\begin{teorema}[Teorema de uniformizaci\'on
anal{\'\i}tica\index{teorema de uniformizaci\'on anal{\'\i}tica}]\label{DrinfeldT2.2.22}
Sea $\rho$ un $A$--m\'odulo de Drinfeld sobre $\Ci$. Entonces existe una \'unica red
$\Gamma$ tal que $\rho=\rho^{\Gamma}$.
\end{teorema}

\begin{proof} 
Sea $\rho$ un $A$--m\'odulo de Drinfeld sobre $\Ci$ de rango $d$.
Para la demostraci\'on basta considerar un subcampo $F$ de $\Ci$ que
sea completo y que
contenga a $\Ki$. El argumento lo haremos sobre $F$.

Consideramos la funci\'on entera $\xi_{\rho}$ dada por el Teorema \ref{DrinfeldNT1}.

Sea $M_{\rho}=\{x\in\Ci\mid \xi_{\rho}(x)=0\}$ el cual est\'a
contenido en la cerradura separable de 
$F\subseteq \Ci$. 
Adem\'as los ceros de $\xi_{\rho}$ son simples (Proposici\'on \ref{DrinfeldP2.2.10} (4)).
Veamos que $M_{\rho}$ es una red. Por la Proposici\'on \ref{DrinfeldP2.2.10} (3),
se tiene que $M_{\rho}$ es un $A$--m\'odulo discreto. Falta ver que
$M_{\rho}$ es finitamente generado.

Sea $V$ el $\K_{\infty}$ espacio vectorial $V=\K_{\infty}\otimes_A M_{\rho}=
\K_{\infty} M_{\rho}$ y queremos ver que $\dim_{\K_{\infty}} \K_{\infty} M_{\rho}
<\infty$.

En caso contrario, sea $\{m_1,\ldots, m_s,\ldots \}$ un infinidad de 
elementos linealmente independientes sobre $\Ki$ y sea
\[
V_i=\Ki m_1\oplus \cdots\oplus \Ki m_i, \e M_i:=M_{\rho}\cap V_i.
\]

Se verifica que $M_i$ es una $A$--red pues si $x,y\in M_i$, $\alpha x+\beta y\in
M_i$, $\alpha,\beta\in A$.

Veamos que $M_i$ es finitamente generado sobre $A$. Se tiene $\Ki M_i\subseteq
V_i$ y si $\{w_1,\ldots, w_s\}$, $s\leq i$, es una base de $\Ki M_i=W$ sobre $\Ki$,
sea
\[
L:=Aw_1+\cdots+Aw_s\subseteq M_i\subseteq W.
\]

Puesto que $M_i$ es
discreto en $V$, existe una vecindad $U$ de $0$ en $V$ con $M_i\cap U=\{0\}$.
Existe una vecindad $U_0$ de $0$ tal que $U_0=-U_0$ y $U_0+U_0=U$, de 
hecho se puede suponer que $U=B(0,\epsilon)$ y $U+U=U$ pues si $x,y\in U$,
$|x|_{\infty}<\epsilon$, $|y|_{\infty}<\epsilon$, $|x+y|_{\infty}\leq \max\{|x|_{\infty},
|y|_{\infty}\}<\epsilon$. Con esto se tiene que para $x,y\in M_i$, $(U_0+x)\cap
(U_0+y)\neq \emptyset\iff x=y$, pues si $z\in (U_0+x)\cap (U_0+y)$ entonces
$x=z-a$, $y=z-b$ con $a, b\in U_0$. Por tanto $x-y=b-a\in U_0-U_0\in M_i
\cap U=\{0\}$, por lo que $x=y$.

Se sigue que $M_i/L$ es un subgrupos discreto tanto de 
$V_i/L$ como del grupo 
$W/L=(\Ki w_1+\cdots+\Ki w_n)/(Aw_1+\cdots+Aw_n)$, 
el cual es compacto pues $\K_{\infty}/A$ es compacto al ser $\K_{\infty}$
localmente compacto y $\K $ denso en $\K_{\infty}$, $\coc A=\K $.

Por lo tanto $M_i/L$ es
finito y en particular finitamente generado. Adem\'as 
\begin{gather*}
\dim_{\Ki}(\Ki\otimes_A L)=
\dim_{\Ki}(\Ki\otimes_A M_i)=\dim_{\Ki}(W)=s. 
\intertext{Finalmente}
a^{-1}M_i/M_i\cong \big(A/(a)\big)^i\e\text{y}\e a^{-1}M_i/M_i\subseteq a^{-1}M_{\rho}/
M_{\rho}\cong \big(A/(a)\big)^d,
\end{gather*}
lo cual es una contradicci\'on para $i>d$.

En resumen, $M_{\rho}$ es una red de rango $d$.

Sea $\rho^{M_{\rho}}$ el $A$--m\'odulo de Drinfeld asociado a $M_{\rho}$
(Teorema \ref {DrinfeldT2.2.19}) y sea $e_{M_{\rho}}$ su funci\'on exponencial.
Puesto que $\xi_{\rho}$ y $e_{M_{\rho}}$ tienen los mismos ceros,
los cuales son simples y $D\xi_{\rho}=De_{M_{\rho}}=1$, se sigue
que $\xi_{\rho}=e_{M_{\rho}}$. Por tanto
\begin{gather*}
\rho_a=\xi_{\rho}a\xi_{\rho}^{-1}=e_{M_{\rho}}a e_{M_{\rho}}^{-1}=
\rho_a^{M_{\rho}}\quad\text{y}\quad \rho=\rho^{M_{\rho}}. \tag*{$\fin$}
\end{gather*}
\end{proof}

Resaltamos el siguiente hecho de la demostraci\'on del Teorema \ref{DrinfeldT2.2.22}.

\begin{corolario}\label{DrinfeldC2.2.23} Dado un $A$--m\'odulo de Drinfeld $\rho\in\Drin_A(F)$,
$\Ki\subseteq F\subseteq \Ci$ con $F$ un campo completo, existe una \'unica serie
formal $\xi_{\rho}\in F\torcido$ con $D\xi_{\rho}=1$ y $\xi_{\rho}a=\rho_a\xi_{\rho}$
para toda $a\in A$. $\fin$
\end{corolario}

\begin{observacion}\label{DrinfeldO2.2.25} Sea $\rho\colon A\lra \aditivo {\Ci}$ un m\'odulo
de Drinfeld de rango $r$. Sea $\Gamma$ la red de rango $r$ correspondiente a $\rho$,
es decir, $\rho=\rho^{\Gamma}$. Ahora bien, $\rho_a^{\Gamma}=aP(a^{-1}\Gamma/
\Gamma, x)$ donde 
\begin{gather*}
P(\Gamma^{\prime}/\Gamma,x)=x\prod_{\mu\in\Gamma^{\prime}/\Gamma\setminus\{0\}}
\Big(1-\frac{x}{e_{\Gamma}(\mu)}\Big), \e e_{\Gamma}(u)=\exponencial u{\gamma}
{\Gamma}.
\intertext{Por tanto}
aP(a^{-1}\Gamma/\Gamma,x)=a\Big[x\cdot\prod_{\mu\in a^{-1}\Gamma/
\Gamma\setminus \{0\}}\Big(1-\frac{x}{e_{\Gamma}(\mu)}\Big)\Big]=\rho_a(x).
\end{gather*}

Ahora bien, se tiene $\rho_a(x)=0\iff ax\prod_{\mu\in a^{-1}\Gamma/\Gamma\setminus\{0\}}
\big(1-\frac{x}{e_{\Gamma}(u)}\big)=0
\iff x=0$ o $x=e_{\Gamma}(\mu)$, $\mu\in a^{-1}\Gamma/\Gamma$. Por tanto
\begin{gather*}
\text{\fbox{$\rho[a]=e_{\Gamma}\big(a^{-1}\Gamma/\Gamma)$}}, \e 
\text{\fbox{$\rho_a^{\Gamma}(x)=\rho_a(x)=ax \prod_{\mu\in 
a^{-1}\Gamma/\Gamma\setminus
\{0\}}\Big(1-\frac{x}{e_{\Gamma}(\mu)}\Big)$}}.
\intertext{M\'as a\'un, se tiene}
e_{\Gamma}(au)=\rho_a^{\Gamma}(e_{\Gamma}(u))=\rho_a(e_{\Gamma}(u)).
\end{gather*}
\end{observacion}

\begin{proposicion}\label{DrinfeldP2.2.23'}
Si $\rho$ es un m\'odulo de Drinfeld, $a\in A$, $a\neq 0$, la $a$--torsi\'on de
$\rho$ est\'a dada por $\rho[a]=e_{\Gamma}(a^{-1}\Gamma/\Gamma)$. $\fin$
\end{proposicion}

\begin{ejemplo}\label{DrinfeldEj2.2.26} Si $C$ es el 
m\'odulo de Carlitz y $\Gamma_C$ es la
red correspondiente a $C$, se tiene
\begin{gather*}
C_T(e_{\Gamma_C}(u))=\rho_T(e_{\Gamma_C}(u))=Te_{\Gamma_C}(u)+
(e_{\Gamma_C}(u))^q=e_{\Gamma_C}(Tu).
\intertext{Es decir,}
\text{\fbox{$e_{\Gamma_C}(Tu)=Te_{\Gamma_C}(u)+e_{\Gamma_C}(u)^q$}}.
\intertext{En general,}
\text{\fbox{$e_{\Gamma_C}(Mu)=\rho_M(e_{\Gamma_C}(u))=C_M(e_{\Gamma_C}
(u))$}},\e M\in R_T.
\end{gather*}

Como  $\Gamma_C$ es una red de rango $1$ y $A$ es de ideales principales,
$\Gamma=\Gamma_C=A\xi_C=R_T\xi_C$
para alg\'un $\xi=\xi_C\in \Ci$.

Sea $M\in R_T$, $C[M]=e_{\Gamma_C}(M^{-1}\Gamma/\Gamma)
=e_{\Gamma_C}(
M^{-1}(A\xi)/A\xi)=\big\{e_{\Gamma_C}\big(\frac{h}{M}\xi\big)$, $h\in R_T\big\}$.

La $M$--torsi\'on de $C$ es 
\[
\Big\{e_{\Gamma_C}\Big(\frac{h}{M}\xi\Big), h\in R_R\Big\}=\Lambda_M.
\]
\end{ejemplo}

\begin{notacion}\label{DrinfeldEj2.2.27} $e_{\Gamma_C}=e_C$ y 
\fbox{$C_M(e_C(u))=e_C(Mu)$}, 
$M\in R_T$ y $C[M]=\Lambda_M=\big\{e_C\big(\frac{h}{M}\xi\big)\mid h\in R_T\big\}$.
Un generador de $\Lambda_M$ es \fbox{$\frac{1}{M}\xi=\lambda_M$}.
\end{notacion}

En general, para cualquier $\Gamma$:
\begin{gather}
\text{\fbox{$e_{\Gamma}(au)=\rho_a^{\Gamma}(e_{\Gamma}(u))$}}\e\text{o}\e
\text{\fbox{$e_{\Gamma}(u)=\rho_a^{\Gamma}\big(e_{\Gamma}\big(\frac{u}{a}
\big)\big)$}}. \label{DrinfeldEq2.2.28}\\
\text{\fbox{$e_{a^{-1}\Gamma}(u)=a^{-1}e_{\Gamma}(au)$}}. \label{DrinfeldEq2.2.29}\\
\text{\fbox{$e_{\Gamma}(x)=\exponencial x{\gamma}{\Gamma}$}}. \label{DrinfeldEq2.2.30}\\
\begin{align}
\rho_a^{\Gamma}(u)&=aP(a^{-1}\Gamma/\Gamma,u))
=a \exponencial u{\lambda}{
e_{\Gamma}(a^{-1}\Gamma)}\nonumber\\
&=au\prod_{\lambda\in(a^{-}\Gamma/\Gamma)\setminus
\{0\}}\Big(1-\frac{u}{e_{\Gamma}(\lambda)}\Big).
\end{align}\\
\text{\fbox{$\rho[a]=e_{\Gamma}(a^{-1}\Gamma/\Gamma)\e$ con $\e \rho=\rho^{
\Gamma}$}}. \label{DrinfeldEq2.2.32}
\end{gather}

\subsection{C\'alculo de $\xi=\tilde{\pi}$ y de la 
exponencial de Carlitz}\label{DrinfeldS2.3}

Se tiene $e_C(Tu)=C_T(e_C(u))=Te_C(u)+(e_C(u))^q$. Pongamos la expansi\'on
de $e_C(u)$: $e_C(u)=\sum_{i=0}^{\infty}\frac{u^{q^i}}{D_i}$ con radio de convergencia
infinito. Ahora
\begin{align*}
e_C(Tu)&=\sum_{i=0}^{\infty}\frac{T^{q^i}u^{q^i}}{D_i}=Te_C(u)+(e_C(u))^q\\
&=\sum_{i=0}^{\infty}\frac{Tu^{q^i}}{D_i}+\sum_{i=0}^{\infty}\frac{u^{q^{i+1}}}{D_i^q}=
\sum_{i=0}^{\infty}\frac{Tu^{q^i}}{D_i}+\sum_{i=1}^{\infty}\frac{u^{q^i}}{D_{i-1}^q}.
\end{align*}

Por tanto $\frac{T^{q^i}}{D_i}=\frac{T}{D_i}+\frac{1}{D_{i-1}^q}$, $i=1,2,\ldots$
y como $\frac{u^{q^0}}{D_0}=1$, $D_0=1$. Se sigue
\[
T^{q^i}D_{i-1}^q=TD^q_{i-1}+D_i\e\text{y por tanto}\e D_i=(T^{q^i}-T)D_{i-1}^q, i\geq 1.
\]

\begin{definicion}\label{DrinfeldD2.3.1} Para $i\geq 0$, se define $[i]:=T^{q^i}-T$, $[0]=0$.
\end{definicion}

Obtenemos que
\begin{gather*}
\begin{align*}
D_i&= [i]D_{i-1}^q=[i][i-1]^qD_{i-2}^{q^2}=[i][i-1]^q[i-2]^{q^2}D_{i-3}^{q^3}=\ldots \\
&=[i][i-1]^q[i-2]^{q^2}\ldots[1]^{q^{i-1}}D_0^{q^i}=
[i][i-1]^q[i-2]^{q^2}\ldots[1]^{q^{i-1}}.
\end{align*}
\intertext{Esto es}
\text{\fbox{$D_i=\prod_{j=1}^i [j]^{q^{i-j}}$}}.
\end{gather*}

En particular, $\deg_TD_i=\sum_{j=1}^iq^{i-j}\deg_T[j]=
\sum_{j=1}^iq^{i-j}q^j=iq^i$. Por tanto
\[
\text{\fbox{$\deg_T D_i=iq^i$}},\e \vi(D_i)=-iq^i\e\text{y}\e \vi\Big(
\frac{u^{q^i}}{D_i}\Big)=q^i\vi(u)+iq^i.
\]

Para $i\gg 0$, $\vi\big(\frac{u^{q^i}}{D_i}\big)=q^i(\vi(u)+i)>0$ y $\lim_{i\to\infty}\vi\big(
\frac{u^{q^i}}{D_i}\big)=\infty$, $\sum_{i=0}^{\infty}\frac{u^{q^i}}{D_i}$ es entera y
\[
\text{\fbox{$e_C(u)=\sum_{i=0}^{\infty}\frac{u^{q^i}}{D_i}
=\sum_{i=0}^{\infty}\frac{u^{q^i}}{
\prod_{j=1}^i[j]^{q^{i-j}}}$}}.
\]

Ahora $e_C(u)$ no es inyectivo pero podemos definir una inversa de $e_C(u)$
en una vecindad de $0$.

Sea $\log_C(u)=\sum_{i=0}^{\infty}\frac{(-1)^iu^{q^i}}{L_i}$ la inversa de $e_C(u)$
en una vecindad de $0$, $L_0=1$. 
Ahora se tiene $e_C(\alpha+\beta)=e_C(\alpha)+e_C(\beta)$.
Si $a=e_C(\alpha)$, $b=e_C(\beta)$, entonces $a+b=e_C(
\alpha)+e_C(\beta)=e_C(\alpha+\beta)$ y $
\log_C(e_C(\alpha+\beta))=\alpha+\beta=\log_C(a+b)=\log_C(a)+\log_C(b),$
por lo que 
\[
\text{\fbox{$\log_C(a+b)=\log_C(a)+\log_C(b)\e\text{para cualesquiera}\e a,b$}}.
\]

As{\'\i}, $e_C(Tu)=Te_C(u)+e_C(u)^q$. Se sigue que $\log_C(e_C(Tu))=Tu=
\log_C(Te_C(u))+\log_C(e_C(u)^q)$. Sea $v=e_C(u)$, $u=\log_C(v)$. Por tanto
\[
\text{\fbox{$T\log_C(v)=\log_C(Tv+v^q)$}}.
\]

Obtenemos 
\[
T\sum_{i=0}^{\infty}\frac{(-1)^iv^{q^i}}{L_i}=T\log_C(v)=\log_C(Tv+v^q)=
\sum_{i=0}^{\infty}\frac{(-1)^i(Tv+v^q)^{q^i}}{L_i}.
\]
Por lo tanto
\begin{align*}
\sum_{i=0}^{\infty}\frac{(-1)^iTv^{q^i}}{L_i}&=\sum_{i=0}^{\infty}\frac{(-1)^iT^{q^i}v^{q^i}+
(-1)^iv^{q^{i+1}}}{L_i}\\
&=\sum_{i=0}^{\infty}\frac{(-1)^iT^{q^i}v^{q^i}}{L_i}+\sum_{i=0}^{\infty}\frac{(-1)^i
v^{q^{i+1}}}{L_i}\\
&= \sum_{i=0}^{\infty}\frac{(-1)^iT^{q^i}v^{q^i}}{L_i}+\sum_{i=1}^{\infty}\frac{(-1)^{i-1}
v^{q^i}}{L_{i-1}}.
\end{align*}

De esta forma $\frac{T}{L_i}=\frac{T^{q^i}}{L_i}-\frac{1}{L_{i-1}}$, por lo que 
$L_{i-1} T=L_{i-1}T^{q^i}-L_i$, $L_i=(T^{q^i}-T)L_{i-1}$, $i\geq 1$. Se sigue que
\[
L_i=[i]L_{i-1}=[i][i-1]L_{i-2}=\cdots=[i][i-1]\cdots[1],\e L_0=1,\e\text{y}\e
L_i=\prod_{j=1}^i[j].\label{DrinfeldLi}
\]

Se tiene $\deg_T L_i=\sum_{j=1}^i\deg_T[j]=\sum_{j=1}^{i}{q^j}=q\frac{q^i-1}{q-1}$
y la serie $\log_C(u)=\sum_{i=0}^{\infty}\frac{(-1)^iu^{q^i}}{L_i}$ converge en $u
\iff \big|\frac{u^{q^i}}{L_i}\big|\xrightarrow[i\to\infty]{}0\iff \vi\big(\frac{u^{q^i}}{L_i}\big)
\xrightarrow[i\to\infty]{}\infty$. Ahora bien
\[
\vi\big(\frac{u^{q^i}}{L_i}\big)=q^i\vi(u)-q\frac{q^i-1}{q-1}=q^i\Big(\vi(u)-\frac{q}{q-1}
\big(1-\frac{1}{q^i}\big)\Big).
\]
Se tiene que $\vi(u)-\frac{q}{q-1}\big(1-\frac{1}{q^i}\big)\xrightarrow[i\to\infty]{}
\vi(u)-\frac{q}{q-1}$. Por lo tanto
\[
\lim_{i\to\infty}\vi\big(\frac{u^{q^i}}{L_i}\big)=\infty\iff \vi(u)>\frac{q}{q-1}\iff
\deg_T (u)<-\frac{q}{q-1}.
\]

Nuestro siguiente objetivo es encontrar una expresi\'on de $\tilde{\pi}$ y hallar
su grado, el cual debe ser $-\frac{q}{q-1}$ o $\vi(u)=\frac{q}{q-1}$ pues el conjunto
de ceros de $e_C(u)$ es $A\tilde{\pi}$ y su inversa alrededor de $0$ est\'a
definida para elementos $u$ con $\vi(u)>\frac{q}{q-1}$.

Sea 
\begin{align*}
e_C(x\log_C(u))&=C_x(e_C(\log_C(u)))=C_x(u)=
\sum_{j=0}^{\infty}\carlitzbinom xj u^{q^j}\\
&=\sum_{i=0}^{\infty}\frac{x\log_C(u)^{q^i}}{D_i}=\sum_{i=0}^{\infty}
\frac{x^{q^i}}{D_i}\Big(
\sum_{j=0}^{\infty}\frac{(-1)^ju^{q^j}}{L_j}\Big)^{q^i}\\
&=\sum_{i=0}^{\infty}\frac{x^{q^i}}{D_i}\Big(\sum_{j=0}^{\infty}\frac{(-1)^{j q^i}u^{q^{i+j}}}
{L_j^{q^i}}\Big)=\sum_{i=0}^{\infty}\sum_{j=0}^{\infty}\frac{x^{q^i}}{D_i}\frac{(-1)^{jq^i}
u^{q^{i+j}}}{L_j^{q^i}}\\
&\underbracket[0pt]{=}_{\substack{\uparrow\\ t=i+j}}\sum_{t=0}^{\infty}\Big(
\sum_{j=0}^t \frac{x^{q^{t-j}}}{D_{t-j}}\frac{(-1)^{jq^{t-j}}}{L_j^{q^{t-j}}}\Big)u^{q^t}.
\end{align*}

Por tanto
\begin{gather*}
\carlitzbinom xt=\sum_{j=0}^t\frac{x^{q^{t-j}}}{D_{t-j}}\frac{(-1)^{jq^{t-j}}}{L_j^{q^{t-j}}}
\underbracket[0pt]{=}_{\substack{\uparrow\\ (-1)^{q^s}=-1}}\sum_{j=0}^t \frac{(-1)^{t-j}}{
L_{t-j}^{q^j}}\frac{x^{q^j}}{D_j}\e\text{por lo tanto}\\
\text{\fbox{$\carlitzbinom xj =\sum_{i=0}^j(-1)^{j-i}\frac{x^{q^i}}{D_i L_{j-i}^{q^i}}$}}.
\end{gather*}

Por otro lado, puesto que para $M\in R_T$, $C_M=\sum_{i=0}^d C_{M,i} \tau^i$
es de grado $d=\deg M$ y $C_{M,i}=\carlitzbinom Mi$, (aqu{\'\i} $\carlitzbinom Mi$
significa lo mismo en esta nueva notaci\'on que en la notaci\'on usual de campos
de funciones ciclot\'omicas, ver Teorema \ref{T6.2.3}), se sigue
\[
\sum_{i=0}^d\carlitzbinom Mi u^{q^i}=C_M(u)=e_C(M\log_C(u))=\sum_{i=0}^{\infty}
\carlitzbinom Mi u^{q^i}.
\]

Para $j>\deg M$, $\carlitzbinom Mj=0$, esto es, todo polinomio de grado menor
a $j$ es cero de $\carlitzbinom xj$ y hay $q^j$ de estos polinomios. Definimos
\begin{gather*}
\text{\fbox{$\ex_t(x)=\prod_{\substack{M\in R_T\\ \deg M<t}}(x-M)
\underbracket[0pt]{=}_{\substack{
\uparrow\\ (x-M)=-M\big(1-\frac{x}{M}\big)}} A_t x \prod_{\substack{M\in R_T
\setminus\{0\}\\\deg M<t}}\big(1-\frac{x}{M}\big)$}},
\end{gather*}
con \fbox{$A_t=\prod_{\substack{M\in R_T\setminus\{0\}\\ \deg M<t}}(-M)$}.

Se tiene que $\ex_t(x)$ es de grado $q^t$ en $T$ y se anula en 
todo $M\in R_T$ con $\deg M
<t$. Se sigue que $\carlitzbinom xt=\d\frac{\ex_t(x)}{B_t}$ para alg\'un $B_t\in \K =
\coc A$, $A=R_T$.

Calculemos $B_t$ y $A_t$.

Se tiene $C_{T^t}(u)=\sum_{j=0}^t\carlitzbinom {T^t}j u^{q^j}$ y $\carlitzbinom  {T^t}t=1$.
Por lo tanto $\carlitzbinom {T^t}t=1=\frac{\ex_t(T^t)}{B_t}$. Se sigue que
\[
B_t=\ex_t(T^t)=\prod_{\deg M<t}(T^t-M)=\prod_{\substack{\text{$N$ m\'onico}\\
\deg N=t}}N.
\]
El coeficiente de $x$ en $\carlitzbinom xt=\sum_{i=0}^t (-1)^{t-i}\frac{x^{q^i}}{D_i
L_{t-i}^{q^i}}$ aparece con $i=0$ y por lo tanto es igual a $\frac{(-1)^t}{D_0 L_t}=
\frac{(-1)^t}{L_t}$. El coeficiente de $x$ en $\frac{\ex_t(x)}{B_t}$ es la constante de
$\frac{\ex_t(x)}{x B_t}$ la cual es $\frac{A_t}{B_t}$.

En resumen, \fbox{$\frac{(-1)^t}{L_t}=\frac{A_t}{B_t}$} y $A_t=\prod_{\substack{M\in R_T
\setminus\{0\}\\ \deg M<t}}(-M)$.

Ahora 
\begin{gather*}
\begin{align*}
A_t&=\prod_{\substack{M\in R_T\setminus\{0\}\\ \deg M<t}}(-M)=
\prod_{j=0}^{t-1}\prod_{\deg M=j}(-M)\\
&\underbracket[0pt]{=}_{\substack{\uparrow\\
M=a_jT^j+\cdots+a_0\\a_j\big(T^j+\cdots+\frac{a_0}{a_j}\big)}}\prod_{j=0}^{t-1}
\Big(\prod_{\alpha\in\*\F}\alpha\Big)\Big\{\prod_{\substack{\text{$M$ m\'onico}\\
\deg M=j}}(-M)\Big\}^{q-1}\\
&\underbracket[0pt]{=}_{\substack{\uparrow\\ \prod_{\alpha\in \*\F}\alpha=-1}}
\prod_{j=0}^{t-1} (-1)B_j^{q-1}=(-1)^t(B_0\cdots B_{t-1})^{q-1}.
\end{align*}
\intertext{Por lo tanto}
\text{\fbox{$A_t=(-1)^t (B_0\cdots B_{t-1})^{q-1}$}}.
\end{gather*}

Se sigue que 
\[
\frac{(-1)^t}{L_t}=\frac{A_t}{B_t}=\frac{(-1)^t (B_0\cdots B_{t-1})^{q-1}}
{B_t},
\]
es decir $B_t=L_t(B_0\cdots B_t)^{q-1}$. Para $t+1$ tendremos
\begin{align*}
B_{t+1}&=L_{t+1}(B_0\cdots B_{t-1}B_t)^{q-1}=L_{t+1}(B_0\cdots B_{t-1})^{q-1}
B_t^{q-1}\\
&=L_{t+1}\frac{B_t}{L_t}B_t^{q-1}=\frac{L_{t+1}}{L_t} B_t^q=[t+1]B_t^q.
\end{align*}

Por tanto \fbox{$B_{t+1}=[t+1] B_t^q$}. 
Por inducci\'on se sigue que $B_t=D_t$ de donde
tendremos 
\[
\text{\fbox{$\carlitzbinom xt=\frac{\ex_t(x)}{D_t}$}}, \quad t\in{\ma N}\cup\{0\}.
\]

As{\'\i}
\begin{gather*}
\carlitzbinom xt =\frac{\ex_t(x)}{D_t}=\frac{A_t}{B_t} x\prod_{\substack{
M\in R_T\setminus\{0\}\\ \deg M<t}}\Big(1-\frac{x}{M}\Big)=\frac{(-1)^t}{L_t}
x\prod_{\substack{M\in R_T\setminus\{0\}\\ \deg M<t}}\Big(1-\frac{x}{M}\Big).
\intertext{En resumen}
\text{\fbox{$\carlitzbinom xt =\frac{(-1)^t}{L_t}x
\prod_{\substack{M\in R_T\setminus\{0\}\\
\deg M<t}}\Big(1-\frac{x}{M}\Big)$}}.
\end{gather*}

Se sigue que
\begin{align*}
x\prod_{\substack{M\in R_T\setminus\{0\}\\
\deg M<t}}\Big(1-\frac{x}{M}\Big)&=(-1)^tL_t\carlitzbinom xt=(-1)^tL_t
\sum_{i=0}^t(-1)^{t-i}\frac{x^{q^i}}{D_iL_{t-i}^{q^i}}\\
&=\sum_{i=0}^t(-1)^i\frac{L_t}{D_iL_{t-i}^{q^i}} x^{q^i}=\frac{\ex_t(u)}{A_t}.
\end{align*}

Recordemos que estamos buscando $\tilde{\pi}$. Obtenemos
\begin{gather*}
\begin{align*}
\sum_{i=0}^{\infty}\frac{\tilde{\pi}^{q^i}u^{q^i}}{D_i}&=
e_C(\tilde{\pi}u)=\tilde{\pi}\exponencial xM{R_T}\\
&=\tilde{\pi}\lim_{t\to\infty}\frac{\ex_t(u)}{A_t}=\tilde{\pi}
\lim_{t\to\infty}\Big(\sum_{i=0}^t (-1)^i\frac{L_t}{D_iL_{t-i}^{q^i}}u^{q^i}\Big)\\
&=\tilde{\pi}\sum_{i=0}^{\infty}\Big\{\lim_{t\to\infty}(-1)^i\frac{L_t}{L_{t-i}^{q^i}}
\Big\}\frac{u^{q^i}}{D_i}.
\end{align*}
\intertext{Por tanto}
\text{\fbox{$\tilde{\pi}^{q^i-1}=\lim_{t\to\infty}(-1)^i \frac{L_t}{L_{t-i}^{q^i}}$}}.
\end{gather*}

Ahora bien, $\deg \frac{L_t}{L_{t-i}^{q^i}}=q\frac{q^t-1}{q-1}-q^i q\frac{q^{t-i}-1}{q-1}=
\frac{q}{q-1}(q^i-1)$. Por tanto 
\[
\text{\fbox{$\vi (\tilde{\pi})=\frac{q}{q-1}$}}.
\]

Continuamos con nuestro desarrollo para hallar $\tilde{\pi}$. Se tiene
$[i+1]-[i]=T^{q^{i+1}}-T-T^{q^i}+T=(T^q-T)^{q^i}=[1]^{q^i}$. Sea

\[
\alpha_i:= \prod_{j=2}^i\big(1-\frac{[j-1]}{[j]}\big)=\prod_{j=2}^i\frac{[j]-[j-1]}{[j]}=
\prod_{j=2}^i\frac{[1]^{q^{j-1}}}{[j]}=\frac{[1]^{(q^i-1)/(q-1)}}{L_i}.
\]

Puesto que $\sum_{j=2}^{\infty}\frac{[j-1]}{[j]}$ es convergente, $\lim_{i\to\infty}
\alpha_i=\alpha$ existe. Entonces $\alpha\in \Ci$ y $|\alpha_i|_{\infty}=1$, por lo que
$|\alpha|_{\infty}=1$, $\deg \alpha_i=0$. Ahora $\deg(\alpha_{i+1}-\alpha_i)=
-q^i(q-1)$. Sea $\delta_i=\alpha_i-\alpha$, $\deg \delta_i=-q^i$.

Carlitz \cite{Car35} dedujo que 
\begin{gather*}
\lim_{d\to\infty}\sum_{i=0}^d (-1)^i\frac{L_d}{D_iL_{d-i}^{q^i}}
u^{q^i}=\sum_{i=0}^{\infty}\frac{(-1)^i}{D_i}u^{q^i}\alpha^{q^i-1}x_i,
\intertext{donde $x_i= [1]^{(q^i-1)/(q-1)}$. En particular}
\frac{\tilde{\pi}^{q-1}}{D_1}=\lim_{t\to\infty}\Big(-\frac{L_t}{D_1L_{t-1}^q}\Big)=
(-1)\alpha^{q-1}x_1=(-1)\alpha^{q-1}
\end{gather*}
 y por tanto $\tilde{\pi}=\sqrt[q-1]{D_1}
\alpha=\sqrt[q-1]{-[1]}\alpha$.

Sea $\xi_0$ una ra{\'\i}z fija de $\sqrt[q-1]{-[1]}$. Entonces
\begin{gather}
\text{\fbox{$\tilde{\pi}=\xi_0\prod_{i=1}^{\infty}\Big(1-\frac{[i-1]}{[i]}\Big)$}},\label{DrinfeldEc2.9'}\\
\tilde{\pi} =T(-T)^{\frac{1}{q-1}}\prod_{i=1}^{\infty}(1-T^{1-q^i})^{-1}\in (-T)^{\frac{1}{q-1}}
\Ki,\nonumber
\end{gather}
donde $\Ki=\F((1/T))$. Todo este desarrollo puede consultarse en \cite{Car35}.

En el caso cl\'asico $\tilde{\pi}=\pm 2\pi i$ y $\pm$ es la elecci\'on que hagamos
de $\sqrt[2]{-1}=i$.

\begin{corolario}\label{DrinfeldC2.3.-1}
Se tiene
\begin{gather*}
\frac{1}{\tilde{\pi}}e_C(\tilde{\pi}x)=x\prod_{\substack{M\in R_T\\ M\neq 0}}
\Big(1-\frac{x}{M}\Big). \tag*{$\fin$}
\end{gather*}
\end{corolario}

\subsection{Morfismos (homomorfismos) entre m\'odulos de Drinfeld}\label{DrinfeldS2.4}

\begin{definicion}\label{Drinfeld2.4.1} Sean $\rho,\rho^{\prime}\in \Drin_A(F)$ dos m\'odulos
de Drinfeld sobre $F$. Una {\em isogen{\'\i}a\index{isogen{\'\i}as entre m\'odulos de
Drinfeld}} o {\em morfismo\index{morfismos entre m\'odulos de Drinfeld}} de $\rho$
a $\rho^{\prime}$ es un polinomio torcido $f\in \aditivo F$ tal que $f\rho_a=
\rho_a^{\prime}f$ para toda $a\in A$.
\end{definicion}

El producto de $2$ isogen{\'\i}as es una isogen{\'\i}a. De esta forma $\Drin_A(F)$ es una
categor{\'\i}a cuyos morfismos son las isogen{\'\i}as y el conjunto de isogen{\'\i}as
entre $\rho$ y $\rho^{\prime}$ se denota por 
$\Isog(\rho,\rho^{\prime})$\label{Drinfeldisogenias}. Si $0\neq
f\in\Isog(\rho,\rho^{\prime})$, se tiene
\[
\deg_{\tau}(f\rho_a)=\deg_{\tau}(f)+\deg_{\tau}\rho_a=\deg_{\tau}\rho^{\prime}_a
+deg_{\tau}(f)=\deg_{\tau}(\rho^{\prime}_a f).
\]
Esto es, \fbox{$\deg \rho_a=\deg \rho^{\prime}_a$ para toda $a\in A$}. Por tanto
$\rho$ y $\rho^{\prime}$ tienen el mismo rango, la misma altura y 
pertenecen a la misma categor{\'\i}a.

Si $\phi\in \Isog(\rho,\rho^{\prime})$, $\psi\in\Isog(\rho^{\prime},\rho^{\prime\prime})$,
$\psi\phi\in\Isog(\rho,\rho^{\prime\prime})$, es decir 
\[
\Isog(\rho,\rho^{\prime})\times
\Isog(\rho^{\prime},\rho^{\prime\prime})\lra \Isog(\rho,\rho^{\prime\prime}),\e (\phi,
\psi)\longmapsto \psi\phi
\]
es un mapeo biaditivo y $\psi\phi$ es la composici\'on de
$\phi$ y $\psi$. Es decir, si $\phi=f$ y $\psi=g$, $f\rho_a=\rho^{\prime}_a f$ y 
$g\rho^{\prime}_a=\rho^{\prime\prime}_a g$ para toda $a\in A$, entonces
$gf\rho_a=g\rho^{\prime}_a f=\rho^{\prime\prime}_a gf$ es la composici\'on.

Si $f$ es un isomorfismo y $g=f^{-1}$, se tiene $gf\rho_a=\rho_a=\rho_a gf$ para
toda $a\in A$, por lo tanto $gf=1$ y por ende \fbox{$f,g\in\*{\aditivo F}=\*F$}.

\begin{ejemplo}\label{DrinfeldEj3.4.2}
Sea $A=R_T=\F[T]$, $K=\coc A=\F(T)$ y sean $\rho:=C\colon A\lra \aditivo K$,
$T\longmapsto C_T=T+\tau$; $\rho^{\prime}=
C_T^{\prime}\colon A\lra \aditivo K$, $T\longmapsto
C^{\prime}_T=T-\tau$. Entonces
\begin{gather*}
\text{$\rho$ y $\rho^{\prime}$ son isomorfos $\iff$ existe $\alpha\in\*{\F(T)}=\*\F$
tal que $\alpha C_T=C_T^{\prime}\alpha$}\\
\iff \alpha(T+\tau)=\alpha T+\alpha\tau=(T-\tau)\alpha=
T\alpha-\alpha^q\tau\iff \alpha^{q-1}=-1.
\end{gather*}

Si $p\neq 2$, $-1\neq 1$ y $\alpha^{q-1}=1$ para toda $\alpha\in\*\F$. Por tanto $C$
y $C^{\prime}$ no son isomorfos sobre $K$ para toda $p\neq 2$. Para $p=2$,
$-1=1$ y $C=C^{\prime}$.

Ahora sea $L$ cualquier campo conteniendo $K(\sqrt[q-1]{-1})$ y $\alpha =\sqrt[q-1]
{-1}$, por lo tanto $C^{\prime}_T=T-\tau$ es isomorfo a $C_T=T+\tau$ sobre $L$,
es decir cuando consideramos $C,C^{\prime}\colon A\lra\aditivo L$.
\end{ejemplo}

\section{Teor{\'\i}a de campos de clase}\label{DrinfeldC3}

\subsection{Antecedentes}\label{DrinfeldS3.1}

Aqu{\'\i} presentamos la teor{\'\i}a de campos de clase 
expl{\'\i}cita desarrollada por D. Hayes, la cual
utiliza m\'odulos de Drinfeld de rango $1$ y 
nos da extensiones abelianas expl{\'\i}citas.

En los a\~nos 1930's, L. Carlitz desarroll\'o la teor\'ia de campos de clase
sobre $K={\ma F}_q(T)$ usando el m\'odulo de Carlitz sobre $A=\F[T]$.
Esta teor\'ia fue de alguna forma olvidada. Sin embargo D. Hayes,
alumno de Carlitz, not\'o la gran similitud entre la teor\'ia de grupos
formales desarrollada por L. Lubin y J. Tate para hacer teor\'ia de
campos de clase locales con la acci\'on de Carlitz: grupos formales:
$F(X)=\pi X+X^q$; m\'odulo de Carlitz: $C_T=T+ \tau^q$, y con esto
desarroll\'o la teor\'ia empezada por Carlitz para hacer teor\'ia de campos
de clase global expl\'icita: \cite{Hay74,Hay79,Hay92}. En este
cap\'itulo presentamos las ideas b\'asicas desarrolladas por Hayes.

Sea $\K $ un campo de funciones sobre 
$\F$ y sea $\p$ un lugar fijo, $\deg_\K \p=d_{\infty}$.
Sea $\Fi={\ma F}_{q^{d_{\infty}}}$ el campo residual de 
$\K $ en $\p$ y tambi\'en de la
completaci\'on $\Ki=\K_{\p}$. Sea $A=\{x\in \K 
\mid v_{\pK}(x)\geq 0 \text{\ para toda\ }
\pK\neq \p\}$. Para $x\in\*\K $, $\deg x=-d_{\p}v_{\p}(x)
=-\di\vi(x)$. Si $x\in A$, $\big|A/xA\big|=
q^{\deg x}$ y en general ponemos $\N x=q^{\deg x}$ 
y para ${\eu A}$ ideal de $A$,
\[
\rho[{\eu A}]=\{u\in\bar{\K }\mid\rho_a(u)=0\text{\ para toda\ } 
a\in{\eu A}\}=\{u\in \bar{\K }
\mid \rho_{\eu A}(u)=0\}
\]
donde $R\rho_{\eu A}$ es el 
ideal izquierdo generado por
$\{\rho_a\}_{a\in{\eu A}}$ en $R=\aditivo \K $.

\begin{definicion}\label{DrinfeldD3.1.1} El {\em grupo de Picard\index{grupo de Picard}} de $A$,
$\pic A$ se define por $\pic A=Cl_A=\frac{D_A}{P_A}$ donde $D_A$ es el grupo de
divisores fraccionarios de $A$ y $P_A$ es el subgrupo de los ideales principales
$P_A$ y $h_A=|\pic A|=\di h_\K $.
\end{definicion}

Esta definici\'on ya la hab\'iamos dado (despu\'es del 
Corolario \ref{CRamDed1.2.7}).

Sea $\rho$ un m\'odulo de Drinfeld de rango $r$, $\rho\colon A\lra\aditivo F$ con
$\K \subseteq F$, es decir $\delta\colon \lra F$ es el encaje natural.

Sea $\rho[a]\cong (A/(a))^r$ y sea $F_{\rho,a}:=F(\rho[a])$, $a\in A$,
$a\neq 0$.

\begin{teorema}\label{DrinfeldT3.1.2} Se tiene que $F_{\rho,a}/F$ es una extensi\'on
de Galois. Sea $G_{\rho,a}=\Gal(F_{\rho,a}/F)$. Entonces existe un monomorfismo
natural $G_{\rho,a}\lra \GL_r(A/(a))$ donde $\GL_r(M)=\Aut_{A/aA}(A/(a))^r$ que
es el grupo de matrices invertibles $r\times r$ con entradas en $A/(a)$.
\end{teorema}

\begin{proof} Se tiene $\rho[a]=\{\xi\in\bar{F}\mid \rho_a(\xi)=0\}$ y $\rho_a=a+
\sum_{i=1}^{r\deg a}\alpha_i\tau^i$, $r\geq 1$, $\alpha_r\neq 0$. Por tanto
$\rho_a(\xi)=a\xi +\sum_{i=1}^{r\deg a}\alpha_i\xi^{q^i}$ es un polinomio
separable de grado $q^{r\deg a}$. Por tanto $F_{\rho,a}/F$ es una extensi\'on
de Galois pues al ser $F_{\rho,a}=F(\rho[a])$, la extensi\'on es normal,
esto es, $\rho_a(u)\in F[u]$ se descompone en $F_{\rho,a}$. 

Sea $\sigma\in G_{\rho,a}$. Entonces si $\xi\in\rho[a]$,
$\rho_a(\xi)=0$ y $\sigma\rho_a(\xi)=\rho_a(\sigma \xi)=0$, por lo tanto
$\sigma\xi\in\rho[a]$. Adem\'as,
si $b\in A$,
\[
\rho_b=b+\sum_{i=1}^{r\deg b} \beta_i \tau\in \aditivo F, \quad
\rho_b(\sigma\xi)=b(\sigma\xi)+\sum_{i=1}^{r\deg b}\beta_i
(\sigma\xi)^{q^i}
\] 
con $\beta_1,\cdots,\beta_{r\deg b}\in \bar{F}$ por lo que
$\rho_b(\sigma\xi)=\sigma(\rho_b\xi)$. Por lo tanto $\sigma\in \Aut(\rho[a])$
y $\sigma\in \GL_r(A/aA)$ (se tiene que $\sigma\in\Aut(\rho[a])$ pues $\sigma
\circ\sigma^{-1}=\Id_{\rho[a]}$). Por lo tanto $G_{\rho,a}\subseteq \GL_r(A/(a))$.
$\fin$
\end{proof}

\begin{corolario}\label{Drinfeld3.1.3}
Si $r=1$, $G_{\rho,a}$ es abeliano y $G_{\rho,a}\subseteq \*{(A/(a))}$.
\end{corolario}

\begin{proof} Se tiene $G_{\rho,a}\subseteq \GL_1(A/(a))\cong \*{(A/(a))}$. $\fin$
\end{proof}

\begin{observacion}\label{Drinfeld3.1.4}
En general el problema de encontrar la imagen de $G_{\rho,a}$ en $\GL_r
(A/aA)$ es complicado.
\end{observacion}

\subsubsection{Campos de clase para el caso $h_A=1$}\label{DrinfeldS3.1.2}

Para este caso tenemos que $h_A=1=\di h_\K $
(ver Corolario \ref{CRamDed1.2.7}) por lo que $\di=1$ y $h_\K =1$.
\'Unicamente hay $8$ campos $\K $ de g\'enero $g_\K >0$ y $h_\K =1$ y 
\'unicamente $4$ de ellos tienen un primo de grado $1$ y adem\'as este primo
es \'unico en cada caso. Por tanto hay \'unicamente $5$ casos $A$ con
$h_A=1$: los $4$ anteriormente y el caso $g_\K =0$.

Sea $A$ con $h_A=1$ y sea $\rho\in\Drin_A(\K )$ de rango $1$. Se probar\'a
que $\rho$ si se pueden definir sobre $\K $ en este caso. Sea ${\eu A}$ cualquier
ideal de $A$. Entonces ${\eu A}$ es principal y sea $\alpha_{\eu A}$ un 
generador de ${\eu A}$. 

Ahora, el ideal generado por ${\eu A}$ en $R=\aditivo \K $
es principal. Sea $\rho_{\eu A}\in \aditivo \K $ tal que
 $R{\eu A}=R\rho_{\eu A}$. Podemos seleccionar $\rho_{\eu A}$
m\'onico pues $\rho_{\eu A}\in \aditivo \K $ y si $\xi$ es el coeficiente l{\'\i}der,
$R\rho_{\eu A}=R\xi^{-1}\rho_{\alpha_{\eu A}}$, $\xi\in \K $. 
Notemos que $R\rho_{\eu A}=R\rho_{\alpha_{\eu A}}$.
Sea $\beta_{\eu A}\in \K $ tal que
$\rho_{\alpha_{\eu A}}=\beta_{\eu A}\rho_{\eu A}$, es decir, $\beta_{\eu A}$ es
el coeficiente l{\'\i}der de $\rho_{\alpha_{\eu A}}$.

Ahora bien, puesto que $r_{\rho}=1$, tenemos $\rho[{\eu A}]\cong A/{\eu A}$
(Corolario \ref{DrinfeldC1.3.30}).
Sea $\lambda_{\eu A}$ un generador de $\rho[{\eu A}]=\rho[\alpha_{\eu A}]=
\rho[\rho_{\eu A}]$ como $A$--m\'odulo y sea $G_{\rho,{\eu A}}=\Gal(
\K_{\rho,{\eu A}}/\K )$ con $\K_{\rho,{\eu A}}=\K (\rho[{\eu A}])$.

\begin{definicion}\label{DrinfeldD3.1.5} Se define el {\em polinomio 
ciclot\'omico\index{polinomio
ciclot\'omico} con respecto a ${\eu A}$} por:
\[
\psi_{\eu A}(u)=\prod_{\bar{\alpha}\in \*{(A/{\eu A})}}(u-\rho_{\alpha}(\lambda_{\eu A})).
\]
\end{definicion}

Se tiene que $\deg_u\psi_{\eu A}(u)=
\big|\*{(A/{\eu A})}\big|$. Adem\'as, $\psi_{\eu A}
(u)\in \K (\rho[{\eu A}])[u]=\K_{\rho,{\eu A}}[u]$. 
Si $\rho_{\alpha}(\lambda_{\eu A})=\rho_{\beta}(\lambda_{\eu A})$,
entonces $\rho_{\alpha-\beta}(\lambda_{\eu A})=0$ y $\alpha-\beta\in {\eu A}$,
esto es, $\bar{\alpha}=\bar{\beta}$, es decir, $\psi_{\eu A}(u)$ es separable.

\begin{proposicion}\label{DrinfeldP3.1.6} Si $h_A=1$, entonces para todo ideal ${\eu A}$ no
cero de $A$, se tiene $\psi_{\eu A}(u)\in \K [u]$.
\end{proposicion}

\begin{proof} Si $\sigma\in G_{\rho,{\eu A}}$, se tiene que $\sigma\lambda_{\eu A}$ es 
generador de $\rho[{\eu A}]$ pues $\sigma^{-1}(\sigma \lambda_{\eu A})=
\lambda_{\eu A}$ y para $\bar{\alpha}\in(\*{A/{\eu A})}$, $\rho_{\alpha}(\lambda_{
\eu A})$ es generador de $\rho[{\eu A}]$ pues existe $\beta\in A$ con
$\alpha\beta\equiv 1\bmod {\eu A}$, por lo que $\rho_{\beta}\rho_{\alpha}(
\lambda_{\eu A})=\rho_1(\lambda_{\eu A})=\lambda_{\eu A}$.

Por tanto $\sigma\lambda_{\eu A}=\rho_{\alpha}(\lambda_{\eu A})$ para alg\'un
$\bar{\alpha}\in \*{(A/{\eu A})}$. Digamos $\sigma\lambda_{\eu A}=\rho_{\alpha_{
\sigma}}(\lambda_{\eu A})$, $\bar{\alpha}_{\sigma}\in\*{(A/{\eu A})}$. Se tiene
$\sigma\rho_{\beta}(\lambda_{\eu A})=\rho_{\beta}(\sigma\lambda_{\eu A})=
\rho_{\beta}(\rho_{\alpha_{\sigma}}(\lambda_{\eu A}))=\rho_{\beta\alpha_{\sigma}}(
\lambda_{\eu A})$ y si $\rho_{\beta\sigma_{\sigma}}
(\lambda_{\eu A})=\rho_{\gamma\alpha_{\sigma}}(\lambda_{\eu A})$,
entonces $\beta\alpha_{\sigma}-\gamma\alpha_{\sigma}\in {\eu A}$ y puesto que
$\alpha_{\sigma}\in(\*{A/{\eu A})}$, se sigue que $\bar{\beta}=\bar{\gamma}$ lo
cual implica que $\rho_{\beta}(\lambda_{\eu A})
=\rho_{\gamma}(\lambda_{\eu A})$. En otras palabras,
\begin{align*}
\sigma(\psi_{\eu A}(u))&=\sigma\Big(\prod_{\bar{\beta}\in\*{(A/{\eu A})}}(u-
\rho_{\beta}(\lambda_{\eu A}))\Big)=\prod_{\bar{\beta}\in\*{(A/{\eu A})}}(u-
\rho_{\beta\alpha_{\sigma}}(\lambda_{\eu A}))\\
&= \prod_{\bar{\beta}\in\*{(A/{\eu A})}}(u-\rho_{\beta}(\lambda_{\eu A}))=\psi_{
\eu A}(u),
\end{align*}
de donde se sigue que $\psi_{\eu A}\in \K [u]$. $\fin$
\end{proof}

En la siguiente proposici\'on supondremos que $\rho\colon A\lra\aditivo \K $, esto es,
que $\rho$ est\'a definido sobre $\K $, lo cual es cierto pues $h_A=1$ pero lo probaremos
m\'as adelante.

\begin{proposicion}\label{DrinfeldP3.1.7} Con las notaciones 
anteriores, $h_A=1$ y $\rho$
de rango $1$, consideremos ${\eu A}=
\pK^m$ con $\pK$ un ideal primo no cero de
$A$. Entonces:
\las
\item $\pK$ es totalmente ramificado en $\K (\rho[\pK^m])/\K $ por lo que su
{\'\i}ndice de ramificaci\'on es $e_{\pK}(\K_{\rho,\pK^m}|\K )=[\K_{\rho,\pK^m}:\K ]$.

\item Si ${\eu q}$ es cualquier otro ideal primo 
no cero distinto a $\p$ y a $\pK$, entonces
${\eu q}$ es no ramificado en $\K_{\rho,\pK^m}/\K $.
\end{list}
\end{proposicion}

\begin{proof} Sea $\alpha$ generador de $\eu p$. Por tanto $\alpha^m$ es generador de
$\pK^m$: $D(\rho_{\pK^m})=\rho_{\pK^m}(0)=\alpha^m$ donde $\rho_{\pK^m}(0)$
denota la constante de $\rho_{\pK^m}$, o lo que es lo mismo, a $\rho^{\prime}_{
\pK^m}(u)$.

Sea $\lambda=\lambda_{\pK^m}$ un generador de $\rho[\pK^m]$. Sea $f(u)=
\rho_{\alpha^m}(u)$ y $\rho_{\alpha^m}(\lambda)=0$. Por tanto $g(u)=\Irr(\lambda,
u,\K )\mid f(u)$. Sea $f(u)=g(u)h(u)$, por lo tanto 
\[
\alpha^m=f^{\prime}(u)=
g^{\prime}(u)h(u)+g(u)h^{\prime}(u). \
\]
Sea $\o$ la cerradura entera de $A$ en 
$\K_{\rho,\pK^m}=\K (\rho[{\eu p}^m])$ 
y $\alpha^m=f^{\prime}(\lambda)=g^{\prime}(\lambda)h(\lambda)$.
\[
\xymatrix{\o\ar@{-}[r]\ar@{-}[d]&\K_{\rho,\pK^m}\ar@{-}[d]\\ A\ar@{-}[r]&\K }
\]

Por lo tanto $(g^{\prime}(\lambda))_{\o}\mid (\alpha^m)_{\o}=\pK^m \o$. Puesto
que 
\[
{\eu D}_{\o/A}=\mcd\{F^{\prime}(\xi)\mid \xi \text{\ es entero,\ } \K_{\rho,\pK^m}=
\K (\xi), F(u)=\Irr(\xi,u,\K )\},
\]
se obtiene que ${\eu D}_{\o/A}
\mid (\pL_1\cdots\pL_h)^{me}$
donde $\pK\o=(\pL_1\cdots \pL_h)^e$.

Como consecuencia, obtenemos que los \'unicos posibles primos ramificados en
$\K_{\rho,\pK^m}/\K $ son $\pK$ y $\p$. Veamos a continuaci\'on el valor de $e:=
e_{\K_{\rho,\pK^m}/\K }(\pL_i|\pK)$. Sea $d:=\deg \pK$. Se tiene
\begin{gather*}
\rho_{\pK^m}(u)=\rho_{\pK}(\rho_{\pK^{m-1}}(u))=\sum_{i=0}^d c_i\big(\rho_{\pK^{m-1}}
(u)\big)^{q^i}\e\text{con}\e \rho_{\pK}=\sum_{i=0}^dc_i\tau^i,\e c_0=\alpha.
\intertext{Por tanto}
\rho_{\pK^m}(u)=\rho_{\pK^{m-1}}(u)\Big(\sum_{i=0}^d c_i\big(\rho_{\pK^{m-1}}(u)
\big)^{q^i-1}\Big)=\rho_{\pK^{m-1}}(u)t(u),\\
\intertext{con}
t(u)=\frac{\rho_{\pK^m(u)}}{\rho_{\pK^{m-1}}(u)}=\sum_{i=0}^d c_i\big(\rho_{\pK^{m-1}}
(u)\big)^{q^i-1}\in \K [u].
\end{gather*}

Adem\'as $t(\xi)=0\iff \xi\in\rho[\pK^m]\setminus\rho[\pK^{m-1}]\iff \xi$ es
generador de $\rho[\pK^m]$. Por tanto
\[
t(u)=\prod_{\bar{\delta}\in\*{(A/\pK^m)}}(u-\rho_{\delta}(\lambda))\e\text{y}\e
t(0)=\pm\prod_{\delta}\rho_{\delta}(\lambda)=\frac{\rho_{\pK^m}(0)}
{\rho_{\pK^{m-1}(0)}}=\alpha.
\]

Se sigue que \fbox{$\alpha=\pm\prod_{\delta}\rho_{\delta}(\lambda)$}.

Ahora $\rho_a(u)=uH(u)$ con alg\'un $H(u)\in \K [u]$, por lo tanto $\rho_a(\lambda)=
\lambda H(\lambda)$ y $\lambda\mid \rho_a(\lambda)$. Si $\bar{\delta}\in \*{\big(
A/\pK^m\big)}$ entonces $\rho_{\delta}(\lambda)$ es generador $\rho[\pK^m]$ por
lo tanto $\pm \alpha=\beta_0\lambda^s$ con $s:=\big|\*{\big(A/\pK^m\big)}\big|$ y
$\beta_0$ una unidad en $\o$. Adem\'as $v_{\pL_i}(\lambda)\geq 1$
pues $\pm \alpha=\lambda\cdot\prod_{\delta\neq 1} \rho_{\delta}(\lambda)$
y por lo tanto si se tuviese $v_{\pL_i}(\lambda)\leq 0$, 
entonces se tendr\'ia $v_{\pL_i}(\alpha)
\leq 0$. Se sigue que
\begin{gather*}
e=v_{\pL_i}(\pK)=v_{\pL_i}(\alpha)=v_{\pL_i}\big(\prod_{\delta}\rho_{\delta}(\lambda)\big)=
v_{\pL_i}(\beta_0\lambda^s)=sv_{\pL_i}(\lambda)\geq s.
\intertext{Por lo tanto}
e\geq s=\big|\*{\big(A/\pK^m\big)}\big|\geq \big|\Gal(\K_{\rho,\pK^m}/\K )\big|=
[\K_{\rho,\pK^m}:\K ]\geq e.
\end{gather*}
Se sigue que $e=[\K_{\rho,\pK^m}:\K ]$ y $\pK$ es totalmente ramificado en
$\K_{\rho,\pK^m}/\K $. 
$\fin$
\end{proof}

\begin{corolario}\label{DrinfeldC3.1.8} Si $h_A=1$ y $\rho$ es de rango $1$,
entonces
\las
\item $\psi_{\pK^m}(u)=\d\frac{\rho_{\pK^m}(u)}{\rho_{\pK^{m-1}}(u)}=\d\prod_{\lambda
\in \rho[\pK^m]\setminus\rho[\pK^{m-1}]}(u-\lambda)=\Irr(\lambda,u,\K )=t(u)$.

\item Si ${\eu A}$ es cualquier ideal no cero de $A$, $\K_{\rho,{\eu A}}=
\K (\rho[{\eu A}])=\K (\lambda_{\eu A})$. Entonces 
$G_{\rho,{\eu A}}=\Gal(\K_{\rho,{\eu A}}/\K )\cong
\*{\big(A/{\eu A}\big)}=\*{\big(A/(a)\big)}$ donde ${\eu A}=(a)$.

\item $\big|G_{\rho,{\eu A}}\big|=[\K_{\rho,{\eu A}}:\K ]=\Phi({\eu A}):=\big|\*{
\big(A/{\eu A}\big)}\big|$.
\end{list}
\end{corolario}

\begin{proof} Igual que en el caso ciclot\'omico, es 
decir, para el m\'odulo de Carlitz
(Teorema \ref{T6.2.28}). $\fin$
\end{proof}

\begin{teorema}\label{DrinfeldT3.1.9} Sea $A$ con n\'umero de clase $1$ y sea $\rho\in
\Drin_A(\K )$ un m\'odulo de Drinfeld de rango $1$ sobre $\K $. Sea $\pK$ un ideal
primo no cero de $A$ primo relativo a ${\eu A}$. Entonces si $\varphi_{\pK}$ denota
al {\em s{\'\i}mbolo de Artin\index{simbolo de Artin@s\'imbolo de Artin}}
(automorfismo de Frobenius\index{automorfismo de Frobenius}),
$\varphi_{\pK}=\artinp {\K_{\rho,{\eu A}}/\K }{\pK}$, se tiene $\lambda^{\varphi_{\pK}}=
\varphi_{\pK}(\lambda)=\rho_{\pK}(\lambda)$ para $\lambda\in \rho[{\eu A}]$.
\end{teorema}

\begin{proof} Sea $\lambda=\lambda_{\eu A}$ 
un generador de $\rho[{\eu A}]$ y sea $\pL$ un
divisor primo en $\K_{\rho,{\eu A}}$ encima de $\pK$. Entonces $\frac{\rho_{\pK}(u)}{u}=
\prod_{a\in\*{(A/\pK)}}(u-\rho_a(\lambda_{\pK}))$ y de ah{\'\i} se sigue, como en el caso
ciclot\'omico, que $\frac{\rho_{\pK}(u)}{u}$ es Eisenstein (ver 
\cite[Proposition 12.3.18]{Vil2006} o la demostraci\'on del 
Teorema \ref{T6.3.3}).

Si $\frac{\rho_{\pK}(u)}{u}=\prod_{a\in\*{(A/(a))}}(u-\rho_a(\lambda_{\pK}))=
u^{q^d-1}+\beta_{q^d-2}u^{q^d-2}+\cdots+\beta_1 u +\beta_0$, entonces $\pK$ divide
a $\beta_i$, $0\leq i\leq q^d-2$. Se sigue que
\[
\rho_{\pK}(\lambda_{\eu A})=\lambda_{\eu A}\cdot \prod_{a\in\*{(A/{\eu A})}}
(\lambda_{\eu A}-\rho_a(\lambda_{\pK}))\equiv \lambda_{\eu A}^{q^d}\bmod \pL.
\]

Adem\'as, $\varphi_{\pK}(\lambda_{\eu A})=\lambda_{\eu A}^{
\varphi_{\pK}}\equiv\lambda_{\eu A}^{q^d}\bmod \pL$. Para ver
la igualdad $\lambda_{\eu A}^{\varphi_{\pK}}=\rho_{\pK}(\lambda_{\eu A})$ se procede
como en el caso ciclot\'omico, (ver Teorema \ref{T6.3.3}).  $\fin$
\end{proof}

\begin{observacion}\label{DrinfeldO3.1.10} Puesto que el automorfismo de Frobenius en el
caso $h_A=1$ y en el caso ciclot\'omico (Carlitz) est\'a dado de manera totalmente
an\'aloga, toda la descomposici\'on de los primos finitos es igual en ambos casos. M\'as
dif{\'\i}cil es probar que la descomposici\'on de $\p$ es totalmente similar al
primo infinito ${\mathcal P}_{\infty}$ en el caso del m\'odulo de Carlitz (notemos que
$\di=1$). Ver despu\'es de la demostraci\'on del Teorema \ref{DrinfeldT3.4.7} y el
grupo $I_{\p}$ descrito ah\'i (\cite[Proposition 4.15]{Hay85}).
\end{observacion}

\subsubsection{Campos de clase para el caso general}\label{DrinfeldS.3.1.2}

Sea $\K $ cualquier campo de funciones congruente, $\p$ un lugar fijo de $\K $, $\di=
\deg \p$ y $A=\{x\in \K \mid v_{\pK}(x)\geq 0, \pK\neq \p\}$. Sea $\rho\in\Drin_A(\Ci)$
un $A$--m\'odulo de Drinfeld sobre $\Ci$ de rango $1$.

\begin{definicion}\label{DrinfeldD3.1.3} Un {\em campo de clase\index{campo de clase}}
de $A$ significa una extensi\'on abeliana finita $L$ de $\K $ tal que $\p$ se descompone
totalmente en $L/\K $. Un {\em campo de clase restringido\index{campo de clase
restringido}} de $A$ significa una extensi\'on abeliana finita de $\K $.
\end{definicion}

\begin{definicion}\label{DrinfeldD3.1.4} Sea $\rho$ un $A$--m\'odulo de Drinfeld sobre $\Ci$ tal que
$\delta(a)=a$ para toda $a\in A$ donde $\rho$ es de cualquier
rango. Sea $\K \subseteq E\subseteq \Ci$ un subcampo que
contiene a $\K $. Se dice que $\rho$ {\em est\'a definido sobre $E$\index{campo de definici\'on
de un m\'odulo de Drinfeld}} o que {\em $E$ es un campo de definici\'on para $\rho$} si
$\rho$ es isomorfo sobre $\Ci$ a un $A$--m\'odulo de Drinfeld $\rho^{\prime}$ tal que $\rho_a^{
\prime}\in\aditivo E$ para toda $a\in A$.
\end{definicion}

\begin{observacion}\label{DrinfeldO3.1.12'}
Se tiene que $\K {\ma F}_{\infty}\subseteq \K_{\rho}$ pues todo $\rho\in\Drin_A(
\Ci)$ est\'a definido sobre $\K $ y por la Observaci\'on \ref{DrinfeldO3.1.20''}, 
${\ma F}_{\infty}\subseteq \K_{\rho}$.
\end{observacion}

\begin{ejemplo}\label{DrinfeldEj3.1.5} Si $\rho$ es un $A$--m\'odulo de Drinfeld de rango uno,
entonces $\Ki$ es un campo de definici\'on para $\rho$. En efecto existe una red $\Gamma$
de rango $1$ en $\Ci$ tal que $\rho=\rho^{\Gamma}$. Puesto que $r_{\Gamma}=1$,
entonces $\Gamma\cong I$ donde $I$ es un ideal no cero de $A$ y por tanto existe $\xi
\in \Ci$ tal que $\Gamma^{\prime} =\xi \Gamma\subseteq \Ki\subseteq \Ci$. La construcci\'on
de $\rho^{\Gamma^{\prime}}$ es de hecho sobre $\Ki$. Ahora $\rho^{\Gamma^{\prime}}=
\xi \rho \xi^{-1}$.

M\'as precisamente, $\Gamma=\xi_1I_1$ con $\xi_1\in\Ci$, $\xi_i\neq 0$ e
$I_1$ ideal de $A$. Si $\xi:=\xi_1^{-1}$, entonces $\Gamma':=\xi\Gamma
=I_1$ y la construcci\'on de $\rho^{\Gamma'}$ se puede realizar sobre
$\K_{\infty}$.
\end{ejemplo}

\begin{teorema}\label{DrinfeldT3.1.6} Sea $\rho$ un $A$--m\'odulo de Drinfeld sobre $\Ci$
de cualquier rango. Entonces existe un campo de definici\'on $\K_{\rho}$ de $\rho$ el
cual es finitamente generado sobre $\K $ y que est\'a contenido en cualquier campo de
definici\'on para $\rho$, es decir, $\K_{\rho}$ es el m{\'\i}nimo campo de definici\'on para
$\rho$.
\end{teorema}

\begin{proof}
Para $a\in A$, sea $\rho_a=a+\sum_{i=1}^{r_{\rho} \deg a}c_i\tau^i$, $c_i\in\Ci$. Para
$\xi\in\*{\Ci}$ tenemos
\begin{gather*}
(\xi\rho\xi^{-1})_a=\xi\big(a+\sum_{i=1}^{r_{\rho}\deg a}c_i\tau^i\big)\xi^{-1}=a+
\sum_{i=1}^{r_{\rho}\deg a}\xi c_i\xi^{-q^i}\tau^i=a+\sum_{i=1}^{r_{\rho} \deg a}
\xi^{1-q^i}c_i\tau^i.
\end{gather*}

Como notaci\'on ponemos:
\begin{gather}\label{DrinfeldEq3.1.6'}
\rho_a=a+\sum_{i=1}^{r_{\rho} \deg a} c_i(\rho,a)\tau^i,\\
c_i(\xi\rho\xi^{-1},a)=\xi^{1-q^i}c_i(\rho,a).\nonumber
\end{gather}

Sea $a\in A$ no constante fija y sea $S=\{c_i\mid c_i\neq 0\}$. 
Notemos que $a$ es trascendente sobre $\F$ y $\rho_a\neq a$. Sea 
\begin{gather*}
g:=\mcd\{q^i-1\mid i\in S\}\e\text{y sea}\e g=\sum_{j\in S}\alpha_j (q^j-1)
\end{gather*}
para algunas $\alpha_j\in{\ma Z}$. Para $i\in S$, sea
\begin{gather}\label{DrinfeldEj3.1.6''}
I_i=I_i(\rho,a):=c_i\big(\prod_{j\in S} c_j^{\alpha_j}\big)^{(1-q^i)/g}\in \Ci.
\end{gather}

Se tiene $I_i(\rho,a)=c_i(\rho,a)\big(\prod_{j\in S}c_i(\rho,a)^{\alpha_j}\big)^{(1-q^i)/g}$.
El conjunto $\{I_i(\rho,a)\mid i\in S\}$ est\'an en el campo de definici\'on de $\rho$.

Sea $\theta(c_j)=\xi^{1-q^j}c_j$. Entonces
\begin{align*}
\theta(I_i)&=\xi^{1-q^i}c_i\Big[\prod_{j\in S} (\xi^{1-q^j}c_j)^{\alpha_j}
\Big]^{(1-q^i)/g}\\
&=\xi^{1-q^i}c_i\Big[\prod_{j\in S}\xi^{(1-q^j)\alpha_j}\Big]^{(1-q^i)/g} \Big(
\prod_{j\in S}c_j^{\alpha_j}\Big)^{(1-q^i)/g}\\
&=\xi^{1-q^i}c_i\big(\xi^{-g}\big)^{(1-q^i)/g}\big(\prod_{j\in S}
c_j^{\alpha_J}\big)^{(1-q^i)/g}=I_i.
\end{align*}

Los $A$--m\'odulos de Drinfeld isomorfos a $\rho$, est\'an dados por 
$\xi\rho\xi^{-1}$ con $\xi\in\*\Ci$, es decir, son $\xi\rho_a\xi^{-1}$. Por tanto se tiene
que $\{I_j\mid j\in S\}$ pertenecen a cualquier campo de definici\'on de $\rho$ pues
$I_i(\xi\rho\xi^{-1},a)=I_i(\rho,a)$, por lo que $I_i(\rho,a)$ est\'a en cualquier campo
de definici\'on de $\rho$.

Sea $\K_{\rho}:=\K (I_j\mid j\in S)$. Sea $\xi\in\Ci$ seleccionado tal que $\xi^g=
\prod_{i\in S}c_i^{\alpha_i}$. Entonces $I_i=c_i(\xi^g)^{(1-q^i)/g}=\xi^{1-q^i} c_i$.
Por tanto
\[
(\xi\rho\xi^{-1})_a=a+\sum_{i=1}^{r_{\rho}\deg a}\xi^{1-q^i}c_i\tau^i=
a+\sum_{i=1}^{r_{\rho}\deg a}I_i\tau^i
\]
tiene coeficientes en $\K_{\rho}$.

Falta ver que $\K_{\rho}$ es un campo de definici\'on de $\rho$. Se tiene que
$\xi\rho_a\xi^{-1}\in\aditivo {\K_{\rho}}$ con $\xi^g=\prod_{i\in S} c_i^{\alpha_i}$.
Sea $x\in A$. Se quiere ver que $\xi\rho_x\xi^{-1}\in \aditivo {\K_{\rho}}$. 

Ahora si $\rho'=\xi\rho\xi^{-1}$, $\rho'_a\in \aditivo {\K_{\rho}}$. Se quiere 
probar que $\rho'_x\in \aditivo {\K_{\rho}}$.

Sea $\rho'_a=\sum_{i=0}^rb_i\tau^i$ con $b_0=a$ y $c_i\in \K_{\rho}$
y sea $\rho'_x=\sum_{j=0}^s d_j\tau^j$. Queremos probar que $d_j
\in \K_{\rho}$.

Se tiene $\rho'_a\rho'_x=\rho'_x\rho'_a$ por lo que
\begin{align*}
\mu&=\big(\sum_{i=0}^r b_i\tau^i\big)\big(\sum_{j=0}^s d_j \tau^j\big)=
\sum_{i=0}^r\sum_{j=0}^s b_i d_j^{q^i}\tau^{i+j}=\sum_{l=0}^{r+s}
\big(\sum_{t=0}^l b_t d_{l-t}^{q^t}\big)\tau^l,\\
\delta&=\big(\sum_{j=0}^s d_j\tau^j\big)\big(\sum_{i=0}^r b_i\tau^i\big)=
\sum_{j=0}^s\sum_{i=0}^r d_j b_i^{q^j}\tau^{j+i}=\sum_{l=0}^{r+s}
\big(\sum_{t=0}^l b_t^{q^{l-t}} d_{l-t}\big)\tau^l,
\end{align*}
$\mu=\delta$ lo que implica que $\sum_{t=0}^l b_t d_{l-t}^{q^t}=
\sum_{t=0}^l b_t^{q^{l-t}} d_{l-t}$. Por lo tanto
\begin{gather*}
b_0^{q^l} d_l-b_0d_l=\sum_{t=1}^l \big(b_t d_{l-t}^{q^t}-
b_t^{q^{l-t}}d_{l-t}\big).
\intertext{Puesto que $b_0=a$, obtenemos}
(a^{q^l}-a)d_l=\sum_{t=1}^l \big(b_t d_{l-t}^{q^t}-
b_t^{q^{l-t}}d_{l-t}\big).
\end{gather*}
Se sigue que $d_j\in \K_{\rho}$ y $\K_{\rho}$ es un campo de definici\'on 
$\rho'$. $\fin$
\end{proof}

\begin{observacion}\label{DrinfeldO3.1.15} Se ver\'a que si $\rho$ es de rango $1$,
$\K_{\rho}$ es el campo de clase de Hilbert de $\K $, es decir $\K_{\rho}/\K $
es la m\'axima extensi\'on abeliana no ramificada en donde $\p$ se
descompone totalmente. En particular para $A$--m\'odulos de Drinfeld de
rango $1$ sobre $\Ci$, $\K_{\rho}$ es independiente de $\rho$.
\end{observacion}

Nuevamente consideramos 
$\pic (A)=\frac{D_A}{P_A}$ el grupo de Picard de $A$, $D_A$ los ideales
fraccionarios de $A$ y $P_A$ los ideales fraccionarios principales de $A$.
Tambi\'en se tiene $\pic(A)\cong I_A/Q_A$ donde $I_A$ es el conjunto de los
ideales no cero de $A$ y $Q_A$ es el conjunto de ideales principales
no cero de $A$ (Teorema \ref{T1.3.5.2}).

Sea $h_A=\big|\pic(A)\big|=\di h_\K $ (Corolario
\ref{CRamDed1.2.7}) y sea $\Drin_A(\Ci)$ el conjunto de m\'odulos
de Drinfeld sobre $\Ci$. Sea ${\eu A}$ un ideal de $A$ y $\rho\in \Drin_A(\Ci)$.
Sea $I_{\eu A}$ el ideal izquierdo generado por $\{\rho_a\}_{a\in{\eu A}}$, es
decir, $I_{\eu A}=R\cdot \{\rho_a\}_{a\in{\eu A}}$, $R=\aditivo {\Ci}$. 

Sea $\rho_{\eu A}$ un generador de $I_{\eu A}$. 
Sea $\rho_{\eu A}=f_1(\tau)\rho_{a_1}
+\cdots+f_m(\tau)\rho_{a_m}$, $a_1,\ldots, a_m\in{\eu A}$, $f_i(\tau)\in R$. 
Notemos que para toda $a\in A$, $I_{\eu A}\rho_a\subseteq
I_{\eu A}$ pues si $\mu=\sum_{i=1}^n r_i\rho_{a_i^{\prime}}$, 
$a_i^{\prime}\in{\eu A}$,
$\mu \rho_a=\sum_{i=1}^nr_i \rho_{a_i^{\prime}}\rho_a=\sum_{i=1}^n
r_i\rho_{a_i^{\prime}a}\in I_{\eu A}$ ya que $a_i^{\prime}a\in {\eu A}$.

Dado $a\in A$, $\rho_{\eu A}\rho_a\in I_{\eu A}$ por lo que existe un \'unico
$\rho_a^{\prime}\in R$ tal que $\rho_{\eu A}\rho_a
=\rho_a^{\prime}\rho_{\eu A}$.

Sea $\rho^{\prime}\colon A\lra \aditivo {\Ci}$, $a\longmapsto \rho_a^{\prime}$.
Este es un $A$--m\'odulo de Drinfeld sobre $\Ci$ y se denota
\begin{gather}\label{DrinfeldEq3.1.15'}
\rho^{\prime}:={\eu A}\star \rho.
\end{gather}

\begin{observacion}\label{DrinfeldO3.1.16} Todo el desarrollo anterior es v\'alido
para cualquier m\'odulo de Drinfeld $\rho\in \Drin_A(F)$ con $F$ arbitrario.
A continuaci\'on probamos que $\rho^{\prime}\in \Drin_A(F)$ y de hecho
que si $\delta^{\prime}$ es el mapeo estructural $\delta^{\prime}\colon
A\lra F$ de $\rho^{\prime}$, entonce $\delta^{\prime}=\delta$ donde $\delta$
es el mapeo estructural de $\rho$.
\end{observacion}

Primero notemos que $\rho$ y $\rho^{\prime}$ son is\'ogenos con isogen{\'\i}a
$\rho_{\eu A}$, es decir, $\rho_{\eu A}\rho_a=\rho_a^{\prime}\rho_{\eu A}$.
El hecho de que $\rho^{\prime}$ es un m\'odulo de Drinfeld se verifica de manera
rutinaria:
\begin{gather*}
\left.
\begin{array}{rcl}
\rho_{\eu A}\rho_{ab}&=&\rho_{ab}^{\prime}\rho_{\eu A},\\
\rho_{\eu A}\rho_{ab}&=&\rho_{\eu A}\rho_a\rho_b=\rho_a^{\prime}\rho_{\eu A}
\rho_b=\rho_a^{\prime}\rho_b^{\prime}\rho_{\eu A}, \e
\end{array}\right\}
\Lra  \rho_{ab}^{\prime}=\rho_a^{\prime}\rho_b^{\prime}.\\
\left.
\begin{array}{rcl}
\rho_{\eu A}\rho_{a+b}&=&\rho_{a+b}^{\prime}\rho_{\eu A},\\
\rho_{\eu A}\rho_{a+b}&=&\rho_{\eu A}(\rho_a+\rho_b)=\rho_{\eu A}\rho_a+
\rho_{\eu A}\rho_b\\
&=&\rho_a^{\prime}\rho_{\eu A}+\rho_b^{\prime}\rho_{\eu A}=
(\rho_a^{\prime}+\rho_b^{\prime})\rho_{\eu A},
\end{array}\right\}\e
\Lra \rho_{a+b}^{\prime}=\rho_a^{\prime}+\rho_b^{\prime}.
\end{gather*}

Sea $\rho_{\eu A}=d_{i_0}\tau^{i_0}+\sum_{i=i_0+1}^m d_i\tau^i$, $d_{i_0}\neq 0$.
El t\'ermino $i_0$ en cada lado de la ecuaci\'on $\rho_{\eu A}\rho_a=\rho_a^{\prime}
\rho_{\eu A}$ es $d_{i_0}\delta(a)^{q^{i_0}}=\beta_0d_{i_0}$ donde $\rho_a^{\prime}=
\beta_0+\sum_{i=1}^n \beta_i\tau^i$. Por lo tanto $\beta_0=\delta(a)^{q^{i_0}}$,
es decir, $D(\rho_a^{\prime})=\delta(a)^{q^{i_0}}$.

Ahora $\rho[I]=\ker \rho_I$ (Proposici\'on \ref{DrinfeldP1.3.31})
 y $|\rho[{\eu A}]|=|\ker \rho_{\eu A}|=q^{\deg {\eu A}
-i_0}$. Si $h_{\rho}=0$, $i_0=0$, $\delta(a)^{q^{i_0}}=\delta(a)$ y $\rho^{\prime}$
es un m\'odulo de Drinfeld (este es nuestro caso $\delta\colon A\hookrightarrow \Ci$,
sin embargo podemos probar el caso general). Si $\car \rho={\eu q}\neq 0$,
$\delta(A)\cong A/{\eu q}$ subcampo de $F$. Cuando se prob\'o que $h_{\rho}
\in {\ma N}\cup \{0\}$, se prob\'o que $i_0=h_{\rho}m\deg_\K  {\eu q}$ con
${\eu A}^m=(y)$ principal (ver la demostraci\'on del Teorema \ref{DrinfeldT1.3.29}).

Si $x\in \delta(A)=A/{\eu q}$, $x^{q^{\deg {\eu q}}}=x$ por lo que 
$x^{q^{i_0}}=x^{q^{\deg {\eu q} h_{\rho} m}}=x$. 
Por lo tanto $\delta(a)^{q^{i_0}}=\delta(a)$ y $\rho^{\prime}$
es un m\'odulo de Drinfeld.

Ahora, para $\xi\in\Ci$, 
\begin{gather*}
\rho_a=a+\sum_{i=1}^{r_{\rho}\deg a}c_i(\rho,a)\tau^i\e\text{y}\\
\xi\rho_a\xi^{-1}=a+\sum_{i=1}^{r_{\rho}\deg a}\xi^{1-q^i}c_i(\rho,a)\tau^i=
a+\sum_{i=1}^{r_{\rho}\deg a}c_i(\xi\rho\xi^{-1}, a)\tau^i,
\intertext{es decir}
c_i(\xi\rho\xi^{-1},a)=\xi^{1-q^i}c_i(\rho,a).
\end{gather*}

Ahora en general 
\begin{align*}
\rho_{\eu A}\rho_a&=\Big(\sum_{i=0}^md_i\tau^i\Big)\Big(\sum_{j=0}^{
r_{\rho} \deg a}c_j(\rho,a)\tau^j\Big)\\
&=\sum_{i=0}^m\sum_{j=0}^{r_{\rho}\deg a} d_ic_j(\rho,a)^{q^i}\tau^{i+j}=
\sum_{l=0}^{m+r_{\rho}\deg a}\Big(\sum_{t=0}^l d_t c_{l-t}(\rho,a)^{q^t}\Big)\tau^l,\\
\rho_a^{\prime}\rho_{\eu A}&=\Big(\sum_{j=0}^{r_{\rho}\deg a}c_j(\rho^{\prime},a)
\tau^j\Big)\Big(\sum_{i=0}^m d_i\tau^i\Big)=\sum_{j=0}^{r_{\rho^{\prime}}\deg a}
\sum_{i=0}^mc_j(\rho^{\prime},a)d_i^{q^j}\tau^{j+i}\\
&=\sum_{l=0}^{m+r_{\rho^{\prime}}\deg a}\Big(\sum_{t=0}^l d_t^{q^{l-t}} c_{l-t}
(\rho^{\prime},a)\Big)\tau^l.
\end{align*}

Por lo tanto $r_{\rho}=r_{\rho^{\prime}}$ y
\[
\sum_{t=0}^l d_tc_{l-t}(\rho,a)^{q^t}=\sum_{t=0}^l d_t^{q^{l-t}}c_{l-t}(\rho^{\prime},a),
\e 0\leq l\leq m+r_{\rho}\deg a.
\]

Sea $\rho^{\prime\prime}$ tal que $\rho_{\eu A}(\xi\rho_a\xi^{-1})=\rho^{\prime\prime}_a
\rho_{\eu A}$. Entonces
\begin{align*}
\sum_{t=0}^ld_tc_{l-t}(\xi\rho\xi^{-1},a)^{q^t}&=\sum_{t=0}^l d_t \big((\xi^{1-q^{l-t}})^{q^t}
c_{l-t}(\rho,a)\big)^{q^t}\\
&= \sum_{t=0}^l d_t(\xi^{q^t-q^l})c_{l-t}(\rho,a)^{q^t}\\
&=\frac{1}{\xi^{q^l}}
\Big(\sum_{t=0}^l d_t\big(\xi c_{l-t}(\rho,a)\big)^{q^t}\Big)\\
&= \sum_{t=0}^l d_t^{q^{l-t}}c_{l-t}(\rho^{\prime\prime},a).
\end{align*}

Es decir, $\frac{1}{\xi^{q^l}}\sum_{t=0}^l d_t\big(\xi c_{l-t}(\rho,a)\big)^{q^t}=
\sum_{t=0}^l d_t^{q^{l-t}} c_{l-t}(\rho^{\prime\prime},a)$. ?`Cual es $\rho^{\prime\prime}$?

Tambi\'en se tiene $\rho_{\eu A}\rho_a=\rho_a^{\prime}\rho_{\eu A}$ lo cual implica
que 
\begin{gather*}
\xi\rho_a\xi^{-1}=\xi(\rho_{\eu A}^{-1}\rho_a^{\prime}\rho_{\eu A})\xi^{-1},
\intertext{por lo que} 
\rho_{\eu A}(\xi\rho_a\xi^{-1})=(\rho_{\eu A}\xi\rho_{\eu A}^{-1})\rho_a^{\prime}(
\rho_{\eu A}\xi^{-1}\rho_{\eu A}^{-1})\rho_{\eu A}=(\mu\rho_a^{\prime}\mu^{-1})\rho_{
\eu A}=\rho_a^{\prime\prime}\rho_{\eu A},
\end{gather*}
donde $\mu=\rho_{\eu A}\xi\rho_{\eu A}^{-1}$. El problema es que en general
$\mu\in E\torcido$. Por ejemplo, $\rho_{\eu A}=1-\tau$, $\rho_{\eu A}^{-1}=
\sum_{i=0}^{\infty}\tau^i$, $\rho_{\eu A}\xi\rho_{\eu A}^{-1}=(1-\tau)\xi
\sum_{i=0}^{\infty}\tau^i=(\xi-\xi^q\tau)\big(\sum_{i=0}^{\infty}\tau^i\big)=
\sum_{i=0}^{\infty}\xi\tau^i-\sum_{i=0}^{\infty}\xi^q\tau^{i+1}=\xi+
(\xi-\xi^q)\sum_{i=1}^{\infty}\tau^i$.

Ahora  si ${\eu A}=(a)$ es principal con $a\neq 0$, sea $\alpha=
c_{r_{\rho} \deg a}(\rho,a)$ el coeficiente l{\'\i}der de $\rho_a$. 
Se tiene que $\rho_a$ es generador de $I_{\eu A}$. Entonces
$\rho_{\eu A}=\alpha^{-1}\rho_a$, es decir, el coeficiente l{\'\i}der de
$\rho_{\eu A}$ es $1$. Se tiene $({\eu A}\star \rho)_b=\rho^{\prime}_b$,
\[
\rho_b^{\prime}=\rho_{\eu A}\rho_b\rho_{\eu A}^{-1}=\alpha^{-1}
\rho_a\rho_b\rho_a^{-1}\alpha=\alpha^{-1}\rho_b\rho_a\rho_a^{-1}\alpha=
\alpha^{-1}\rho_b\alpha,
\]
por lo tanto $\rho$ y ${\eu A}\star \rho$ son isomorfos. Entonces, si consideramos
las clases de isomorfismo de $A$--m\'odulos de Drinfeld $\rho$ sobre $F$,
el conjunto $I_A$ de ideales no cero de $A$ act\'ua en los m\'odulos de
Drinfeld por medio de ${\eu A}\star \rho$. Si ${\eu A}$ es principal, $\rho$
y ${\eu A}\star \rho$ son isomorfos, por lo tanto $\pic(A)\cong\frac{I_A}{Q_A}$
act\'ua sobre las clases de isomorfismos de $A$--m\'odulos de Drinfeld.

Adem\'as, ${\eu A}\star \rho$ y $\rho$ tienen el mismo rango y $D({\eu A}
\star\rho)=D(\rho)=\delta$, por lo que $\pic(A)$ act\'ua en las clases de
isomorfismos de m\'odulos de Drinfeld de rango fijo $r$.

\begin{definicion}\label{DrinfeldD3.1.17}
El {\em campo de clase de Hilbert\index{campo de clase de Hilbert}} $H_A
\subseteq \Ci$ es la m\'axima extensi\'on abeliana no
ramificada de $\K =\coc A$ donde $\p$
se descompone totalmente.
\end{definicion}

\begin{observacion}\label{DrinfeldO3.1.18}
Veremos que $H_A=\K_{\rho}$ donde $\rho$ es un $A$--m\'odulo de Drinfeld
de rango $1$. En particular $\K_{\rho}$ es independiente de $\rho$
(ver Proposici\'on \ref{DrinfeldP3.3.1}).
\end{observacion}

Para poder mostrar que $\K_{\rho}=H_A$ se usar\'a una {\em funci\'on
signo\index{funci\'on signo}\index{signo}}
la cual dar\'a lugar a un grupo $\pic^+(A)$ el cual
de hecho es una extensi\'on de $\pic(A)$. Se tiene que $\pic^+(A)$ corresponde
a una extensi\'on $H_A^+$ de $\K $, donde todos los primos finitos de $\K $
son no ramificados. La idea de la funci\'on signo es controlar el coeficiente
l\'ider de $\rho_a$ lo cual resulta m\'as eficiente que controlar $\K_{\rho}$.
De esta forma evitamos tratar con clases de isomorfismos de $A$--m\'odulos
de Drinfeld de rango $1$.

El campo $H_A^+$ es el an\'alogo al caso num\'erico del campo
de clase de Hilbert extendido (una vez fijado el signo).
Ver tambi\'en la Secci\'on \ref{CClaseS4.10}. En el caso
num\'erico, el campo de clase extendido corresponde a que los ideales
primos no cero totalmente descompuestos, son los principales
generados por elementos totalmente positivos. En el caso de $H_A^+$
se tendr\'a el an\'alogo con elementos positivos con respecto al signo.

En el resto de este cap\'itulo consideraremos \'unicamente $A$--m\'odulos de Drinfeld
sobre $\Ci$ de rango $1$, aunque en algunos casos, los resultados
son v\'alidos para rango $r$ arbitrario.

\begin{definicion}\label{DrinfeldD3.1.19} Una {\em funci\'on signo\index{funci\'on
signo}\index{signo}} $\sgn\colon \*{\Ki}\lra\*{\Fi}$ es una homomorfismo tal que $\sgn|_{
\*{\Fi}}=\Id_{\*{\Fi}}$ y $1+\hat{\pK}_{\infty}=U_{\Ki}^{(1)}=U_{\infty}^{(1)}\subseteq
\ker\sgn$. Por tanto $\sgn$ est\'a totalmente determinado por
$\sgn(\pi_{\infty})$ puesto que como grupos $\*{\Ki}=U_{\Ki}^{(1)}\times
\*{\Fi}\times (\pi_{\infty})$. Usamos la convenci\'on $\sgn(0)=0$.
\end{definicion}

Para $\sigma\in\Gal(\Fi/\F)$, la composici\'on $\sigma\circ \sgn$ se llama {\em torcer
 la funci\'on signo\index{signo torcido}} o un {\em torcimiento de la funci\'on
 signo}. La funci\'on $\sigma\circ \sgn$ se llama {\em signo 
 torcido\index{signo torcido}}.
 
 Ahora $\big|\*{\Fi}\big|=q^{\di}-1$ por lo que hay $q^{\di}-1$ funciones signo 
 dependiendo de nuestra elecci\'on $\sgn(\pi_{\infty})\in\*{\Fi}$.
 
 M\'as a\'un, si $\sgn$
 y $\sgn^{\prime}$ son dos funciones signo, sea $f\colon\*{\Ki}\lra \*{\Fi}$, $f(x)=
 \frac{\sgn(x)}{\sgn^{\prime}(x)}$. Se tiene que
  $f(\xi)=1$ para toda $\xi\in U_{\p}=U$ y para
 toda $\xi\in\Fi$. Por tanto $f$ est\'a totalmente determinado por $f(\pi_{\infty})$.
Sea $f(\pi_{\infty})=\xi\in\*\Fi\cong \ma Z/(q^{\di}-1) \ma Z$. Entonces el diagrama
\[
\xymatrix{
\*\Ki\ar[rr]^f\ar[dr]_{\vi}&&\*\Fi\\ & \ma Z\cong \*{\Ki}/(U_{\infty}^{(1)}
\times \*\Fi)\ar[ru]_{\phi}}
\]
es conmutativo pues $(\phi\circ \vi)(\pi_{\infty})=\phi(1)=\xi=f(\pi_{\infty})$.

Puesto que para $x\in\*\Ki$ se tiene $\deg x=-\di \vi(x)$, se sigue
que
\[
f(x)=\phi(\vi(x))=\xi^{\vi(x)}=\xi^{(-\deg x)/\di}=\xi_0^{(\deg x)/\di},
\]
donde $\xi_0=\xi^{-1}\in\*\Fi$.

Finalmente obtenemos que
 $\sgn(x)=\sgn^{\prime}(x)\xi_0^{(\deg x)/\di}$, alg\'un $\xi_0\in \*{\F}$.

\begin{definicion}\label{DrinfeldD3.1.20} Un $A$--m\'odulo de Drinfeld sobre $\Ci$
de cualquier rango,
$\rho\in\Drin_A(\Ci)$ se llama {\em normalizado\index{m\'odulo de Drinfeld
normalizado}} si el coeficiente l{\'\i}der $\mu_{\rho}(x)$ de $\rho_x$ pertenece
a $\Fi$ para toda $x\in A$. Si para alguna funci\'on signo $\sgn$, el
mapeo $x\longmapsto \mu_{\rho}(x)$ es un signo torcido, entonces $\rho$
se llama {\em signo normalizado\index{m\'odulo de Drinfeld signo normalizado}}.
\end{definicion}

Para $x\in A$, sea $\mu_{\rho}(x)$ el coeficiente l{\'\i}der de $\rho_x$. Para
$x,y\in A$, se tiene
\begin{align*}
\rho_{xy}&=\rho_x\rho_y=\Big(\sum_{i=0}^{r_{\rho}\deg x} c_i(\rho,x)\tau^i\Big)
\Big(\sum_{j=0}^{r_{\rho}\deg y} c_j(\rho,y)\tau^j\Big)\\
&=\sum_{i=0}^{r_{\rho}\deg x}\sum_{j=0}^{r_{\rho}\deg y}c_i(\rho,x)c_j(\rho,y)^{q^i}
\tau^{i+j},
\end{align*}
por tanto $\mu_{\rho}(xy)=c_{r_{\rho}\deg x}(\rho,x)c_{r_{\rho}\deg y}(
\rho,y)^{q^{r_{\rho}\deg x}}=\mu_{\rho}(x)\mu_{\rho}(y)^{q^{r_{\rho}\deg x}}$.

As\'i
\begin{gather}\label{DrinfeldEq3.1.20'}
\mu_{\rho}(xy)=\mu_{\rho}(x)\mu_{\rho}(y)^{q^{r_{\rho}\deg x}}=\mu_{\rho}(y)
\mu_{\rho}(x)^{q^{r_{\rho}\deg y}}=\mu_{\rho}(yx).
\end{gather}

Sea ahora $n_0\in{\ma N}$ tal que para toda $m\geq n_0$, existen elementos
de $A$ con $\vi(a)=-m$. Tales elementos $A$ y tal $n_0$ existen pues de hecho
por el Teorema de Riemann--Roch, para $n>\max\{2g_K-1,0\}$ existe $x_n\in \K $
tal que $\eta_{x_n}=\p^n$ donde $\eta_y$ denota el divisor de polos de $y$.
Esto es, $\vi(x_n)=-n$ y para toda $\pK\neq \p$, $v_{\pK}(x)\geq 0$, esto
es $x_n\in A$.

Para extender la definici\'on de $\mu_{\rho}$ a todo $\*{\Ki}$ procedemos de
la siguiente forma. Sea primero $x\in\*{\Ki}$ con $\vi(x)=-m<-n_0$. Notemos
que $A/\p\cong\p^{-m}/\p^{-m+1}\cong \tilde{\pK}_{\infty}^{-m}/
\tilde{\pK}_{\infty}^{-m+1}$ donde $\tilde{\pK}_{\infty}$ 
es la completaci\'on de $\p$ en $\infty$ o, equivalentemente, es la 
extensi\'on de $\p$ a $\Ki$. En particular $x\bmod \p^{-m+1}\neq 0$ en
$\tilde{\pK}_{\infty}^{-m}/\tilde{\pK}_{\infty}^{-m+1}$.

Existe $a\in A$ tal que $a\bmod \p^{-m+1}=a\bmod \tilde{\pK}_{\infty}^{-m+1}=x\bmod
\tilde{\pK}_{\infty}^{-m+1}$. Se define $\mu_{\rho}(x):=\mu_{\rho}(a)$. Lo primero que
debemos verificar es que esta definici\'on es independiente de $a$. Sea 
$b\in A$ tal que $\bar{a}=\bar{b}\bmod \tilde{\pK}_{\infty}^{-m+1}$. En particular $\deg
a=\deg b$ puesto que $\deg a=-\di\vi(a)$ y $\vi(a)=-m=\vi(b)$ y $a-b\in \p^{-m+1}$,
esto es $\vi(a-b)\geq -m+1$ y por tanto
$\deg (a-b)=-\vi(a-b)\di\leq (m-1)\di<m\di=\deg a=\deg b$.

Por tanto $\deg \rho_{a-b}<\deg \rho_a$ y $\rho_a=\rho_{b+(a-b)}=\rho_b+
\rho_{a-b}$ y $\mu_{\rho}(a)=\mu_{\rho}(b)=\mu_{\rho}(x)$.
Se sigue que $\mu_{\rho}(x)$ est\'a 
bien definido para $\vi(x)\leq n_0$, $x\in \*{\Ki}$.

Ahora sea $\*{\Ki}$ arbitrario. Para definir $\mu_{\rho}(x)$ consideremos $\alpha
\in \*{\Ki}$ tal que $\vi(x\alpha)\leq -n_0$ y tal que $\vi(\alpha)\leq -n_0$. Se tiene que
$\mu_{\rho}(x\alpha)$ y $\mu_{\rho}(\alpha)$ est\'an definidas. Puesto que queremos
$\mu_{\rho}(x\alpha)=\mu_{\rho}(x)\mu_{\rho}^{q^{r_{\rho}\deg x}}
(\alpha)$, m\'as precisamente $\mu_{\rho}(\alpha)$, definimos
\[
\mu_{\rho}(x)=\mu_{\rho}(x\alpha)\mu_{\rho}(\alpha)^{-q^{r_{\rho}\deg x}}.
\]

Veamos que $\mu_{\rho}(x)$ est\'a bien definido. Sea $\beta\in \*{\Ki}$ tal que
$\vi(x\beta)\leq -n_0$ y $\vi(\beta)\leq -n_0$. Entonces
\[
\mu_{\rho}(x)=\mu_{\rho}(x\beta)\mu_{\rho}(\beta)^{-q^{r_{\rho}\deg x}}.
\]

Consideremos $x\alpha\beta$. Sean $a,b,c\in A$ tales que $\overline{x\alpha}=
\bar{c}$, $\bar{\alpha}=\bar{a}$, $\bar{\beta}=\bar{b}$, por lo que $\overline{x\alpha
\beta}=\overline{bc}$. Se tiene
\begin{gather*}
\deg c=\deg(x\alpha)=\deg x+\deg \alpha=\deg x+\deg a.
\intertext{Ahora, pongamos $r_{\rho}=r$ y}
\begin{align*}
\mu_{\rho}(x)&=\mu_{\rho}(x\alpha)\mu_{\rho}(\alpha)^{-q^{r\deg x}}=\mu_{\rho}(c)
\mu_{\rho}(a)^{-q^{r\deg x}}\e\text{y}\\
\mu_{\rho}(x)&=\mu_{\rho}(x\alpha\beta)\mu_{\rho}(\alpha\beta)^{-q^{r\deg x}}=
\mu_{\rho}(bc)\mu_{\rho}(ab)^{-q^{r\deg x}}\\
&=\mu_{\rho}(c)\mu_{\rho}(b)^{q^{r\deg c}}\big[\mu_{\rho}(a)\mu_{\rho}(b)^{q^{r\deg a}}
\big]^{-q^{r\deg x}}\\
&=\mu_{\rho}(c)\mu_{\rho}(b)^{q^{r\deg c}}\mu_{\rho}(a)^{-q^{r\deg x}}
\mu_{\rho}(b)^{-q^{r(\deg a+\deg x)}}\\
&=\mu_{\rho}(c)\mu_{\rho}(a)^{-q^{r\deg x}}\mu_{\rho}(b)^{q^{r\deg c}}
\mu_{\rho}(b)^{-q^{r\deg c}}\\
&=\mu_{\rho}(c)\mu_{\rho}(a)^{-q^{r\deg x}}.
\end{align*}
\end{gather*}

Es decir, 
\begin{gather*}
\mu_{\rho}(x\alpha)\mu_{\rho}(\alpha)^{-q^{r\deg x}}=
\mu_{\rho}(x\alpha\beta) \mu_{\rho}(\alpha\beta)^{-q^{r\deg x}}.
\intertext{Similarmente}
\mu_{\rho}(x\beta)\mu_{\rho}(\beta)^{-q^{r\deg x}}
=\mu_{\rho}(x\alpha\beta)\mu_{\rho}(\alpha\beta)^{-q^{r\deg x}}.
\intertext{Por lo tanto}
\text{\fbox{$\mu_{\rho}(x)\text{\ est\'a bien definida}$}}.
\end{gather*}

\begin{proposicion}\label{DrinfeldP3.1.20'}
Se tiene que $\mu_{\rho}(x)=1$ para toda $x\in U^{(1)}_{\p}$.
\end{proposicion}

\begin{proof}
Sea $x=1+\alpha \pi$ con $v_{\infty}(\alpha)\geq 0$. Se tiene 
$\deg x=-d_{\infty}v_{\infty}(x)=0$. Sea $y\in A$ con $v_{\infty}(y)=
-m<-n_0$. Por tanto tendremos $yx=y+\alpha y \pi$ y $yx\equiv
y\bmod \tilde{\p}^{-m+1}$ pues $yx-y=\alpha y \pi \in \tilde{\p}^{-m+1}$.
Por tanto $\mu_{\rho}(xy)=\mu_{\rho}(y)$ y 
\begin{gather*}
\mu_{\rho}(x)=\mu_{\rho}(xy)\mu_{\rho}^{q^{-r_{\rho}\deg x}}(y)=
\mu_{\rho}(y)\mu_{\rho}(y)^{-1}=1. \tag*{$\fin$}
\end{gather*}
\end{proof}

Ahora fijemos una funci\'on signo $\sgn$ y estudiaremos el objeto $(\K ,\p,\sgn)$.
Un elemento $x\in A$ (o $x\in \K $ o $x\in \Ki$) se llama {\em positivo\index{elemento
positivo}} si $\sgn(x)=1$.

Un resultado central es:

\begin{teorema}\label{DrinfeldT3.1.21}
Todo m\'odulo de Drinfeld $\rho\in \Drin_A(\Ci)$ sobre $\Ci$ de rango $1$ es
isomorfo sobre $\Ci$ a un m\'odulo de Drinfeld $\rho^{\prime}$ que es signo
normalizado, es decir, $\mu_{\rho^{\prime}}(x)=\sigma\circ \sgn (x)$ para todo
$x\in A$ y alg\'un $\sigma\in \Gal(\Fi/\F)$.
\end{teorema}

\begin{proof} Sea $\pi$ un elemento primo\label{Drinfeldelementoprimo}
 en $\tilde{\pK}_{\infty}$, es decir,
$\vi(\pi)=1$, que sea positivo para la funci\'on signo $\sgn$. De hecho,
si $\pi^{\prime}$ es cualquier primo $\tilde{\pK}_{\infty}$ y $\sgn \pi^{\prime}
=\gamma\in\*{\Fi}$, se toma $\pi:=\gamma^{-1} \pi^{\prime}$.

Sea $\xi\in \Ci$ seleccionado tal que $\xi^{q^{\di}-1}=\mu_{\rho}(\pi^{-1})$ y
sea $\rho^{\prime}=\xi\rho\xi^{-1}$ el cual es isomorfo a $\rho$ sobre $\Ci$.
Entonces $c_i(\xi\rho\xi^{-1},a)=c_i(\rho^{\prime},a)=\xi^{1-q^i}c_i(\rho,a)$,
$r_{\rho}=1$.

En general por la definici\'on extendida de $\mu_{\rho}$ a todo $\*{\Ki}$ se tiene
que $\mu_{\rho^{\prime}}(y)=\xi^{1-q^{\deg y}}\mu_{\rho}(y)$ para toda
$y\in \*{\Ki}$. Por lo tanto
\[
\mu_{\rho^{\prime}}(\pi^{-1})=\xi^{1-q^{\deg \pi^{-1}}}\mu_{\rho}(\pi^{-1})=
\xi^{1-q^{\di}}\mu_{\rho}(\pi^{-1})=1.
\]
As{\'\i} $\mu_{\rho'}(\pi^{-1})=1$.

Sea $x\in\*{\Ki}$, $x=c\mu\pi^n$, $c\in \*{\Fi}$, $\mu\in U^{(1)}$, $n\in{\ma Z}$.
Entonces $\sgn(x)=\sgn(c) \sgn(\mu) \sgn(\pi)^n
=\sgn c\cdot 1\cdot 1=\sgn c$.
Por tanto $\sgn(x)=c\in\Fi$.

Se tiene que si $m\leq 0$ y si
$\mu_{\rho'}(\mu)=1$, entonces $\mu_{\rho'}(\mu \pi^m)=1$.
En particular, para $a\in A$, 
por la Proposici\'on \ref{DrinfeldP3.1.20'}, tendremos que
 $\mu_{\rho^{\prime}}(a)=
\mu_{\rho^{\prime}}(c\mu\pi^m)=\mu_{\rho^{\prime}}(c)=\mu_{\rho^{\prime}}(
\sgn (a))$.

Para $\xi\in {\ma F}_{\infty}$, $\xi\neq 0$, $A/\p\cong {\ma F}_{\infty}$.
Sea $a\in A$ tal que $a\equiv \xi\bmod\p$. Sea $y\in A$ con $v_{\infty}
(y)=-m<n_0$, $\xi y\equiv ay \bmod \p^{-m+1}$ por lo tanto,
$\mu_{\rho}(\xi y)=\mu_{\rho}(ay)$. Se sigue 
\begin{align*}
\mu_{\rho}(\xi)&=\mu_{\rho}(\xi y)\mu_{\rho}^{-q^{r_{\rho}\deg \xi}}(y)=
\mu_{\rho}(ay)\mu_{\rho}(y)^{-1}\\
&=\mu_{\rho}(y)\mu_{\rho}^{q^{
r_{\rho} \deg y}} (a) \mu_{\rho}(y)^{-1}=\mu_{\rho}(a)^{q^{r_{\rho} \deg y}}
\in {\ma F}_{\infty}.
\end{align*}
Esto es, \fbox{$\mu_{\rho}({\ma F}_{\infty})\subseteq {\ma F}_{\infty}$}.

La restricci\'on $\mu_{\rho^{\prime}}|_{\Fi}$ es un automorfismo $\iota_{\rho^{\prime}}
\colon \Fi\lra\Fi$ que deja fijo a $\F$, es decir $\iota_{\rho^{\prime}}\in \Gal(\Fi/\F)$ y
$\mu_{\rho^{\prime}}(a)=\iota_{\rho^{\prime}}(\sgn (a))$ y $\rho^{\prime}$ es
isomorfo a $\rho$ y es signo--normalizado. $\fin$
\end{proof}

\begin{observacion}\label{DrinfeldO3.1.20''}
Restringiendo $\mu_{\rho}$ a ${\ma F}_{\infty}$ como subcampo de 
$\K_{\infty}\cong {\ma F}_{\infty}((\pi))$, obtenemos un automorfismo
$\nu\colon {\ma F}_{\infty}\lra {\ma F}_{\infty}$ que fija $\F$. Por
tanto $\mu_{\rho}({\ma F}_{\infty})={\ma F}_{\infty}\subseteq \K_{\rho}$
pues corresponde a los coeficientes l\'ider de $\rho$ extendido a
$\K_{\infty}$.

De hecho, por el Teorema \ref{DrinfeldT3.1.21}, podemos tomar $\rho$
signo normalizado. Si $\rho'=\xi\rho \xi^{-1}$ son isomorfos,
se tiene $\mu_{\rho'}(x)=\xi^{1-q^{d_{\infty}}}\mu_{\rho}(x)$
y por tanto $\nu\colon {\ma F}_{\infty}\lra \K_{\rho}$ da
${\ma F}_{\infty}\subseteq \K_{\rho}$.
\end{observacion}

\begin{observacion}\label{DrinfeldO3.1.22} Si $\rho$ es signo normalizado y $\pi$ es un
elemento de $\tilde{\pK}_{\infty}$ con $\sgn(\pi)=1$, entonces $\mu_{\rho}(\pi^{-1})=1$.
\end{observacion}

\begin{proof} Se tiene $\mu_{\rho}(\pi^{-1})=\sigma(\sgn(\pi^{-1}))=\sigma((\sgn \pi)^{-1})=
\sigma(1)=1$.
$\fin$
\end{proof}

\begin{definicion}\label{DrinfeldD3.1.23}
Un {\em $A$--m\'odulo de Hayes\index{m\'odulo de Hayes}} es un $A$--m\'odulo
de Drinfeld sobre $\Ci$ de rango $1$ el cual es signo--normalizado.
\end{definicion}

Denotamos por ${\mathcal H}$ al conjunto de
los $A$--m\'odulos de Hayes.\label{DrinfeldHayes}

\begin{ejemplo}\label{DrinfeldEj3.1.24}
El m\'odulo de Carlitz es un m\'odulo de Hayes.
\end{ejemplo}

\begin{proposicion}\label{DrinfeldP3.1.25}
Si $\rho$ y $\rho^{\prime}=\xi\rho\xi^{-1}$ son dos $A$--m\'odulos de Drinfeld de rango
$1$ sobre $\Ci$ que son signo normalizados, entonces $\xi\in\*{\Fi}$ y $\mu_{\rho}=
\mu_{\rho^{\prime}}$.
\end{proposicion}

\begin{proof} Puesto que $\mu_{\rho}(\pi^{-1})=\mu_{\rho^{\prime}}(\pi^{-1})=1$ por ser signo
normalizados, se tiene que 
\[
1=\mu_{\rho^{\prime}}(\pi^{-1})=\xi^{1-q^{\di}}\mu_{\rho}(\pi^{-1})=\xi^{1-q^{\di}}\]
por lo que $\xi^{q^{\di}-1}=1$ y $\xi\in\*{\Fi}$.

Por tanto para cualquier $a\in A$, $\mu_{\rho^{\prime}}(a)=\xi^{1-q^{\deg a}}
\mu_{\rho}(a)=\mu_{\rho}(a)$. $\fin$
\end{proof}

\begin{corolario}\label{DrinfeldC3.1.26}
Cada clase de isomorfismo de $A$--m\'odulos de Drinfeld de rango uno sobre
$\Ci$ tiene exactamente $(q^{\di}-1)/(q-1)$ $A$--m\'odulos 
de Drinfeld con signo normalizado.
\end{corolario}

\begin{proof} Dado $\rho\in \Drin_A(\Ci)$ de rango $1$, $\rho$ es isomorfo a un
$A$--m\'odulo de Hayes $\rho^{\prime}$. Ahora bien, todo $A$--m\'odulo
$\rho^{\prime\prime}$ isomorfo a $\rho^{\prime}$ est\'a dado por
$\rho^{\prime\prime}=\xi\rho^{\prime}\xi^{-1}$, $\xi\in \*{\Ci}$. Si $\rho^{
\prime\prime}$ es tambi\'en signo normalizado, $\xi\in \*{\Fi}$ 
y finalmente $\Aut(\rho')$,
el grupo de automorfismos de $\rho'$, es isomorfo a $\*{\F}$. Se sigue que hay
$\frac{|\*{\Fi}|}{|\*{\F}|}=\frac{q^{\di}-1}{q-1}$ $A$--m\'odulos de Drinfeld 
signo normalizados isomorfos a $\rho$. $\fin$
\end{proof}

Si $\rho\in{\mathcal H}$ es un m\'odulo de Hayes y ${\eu A}$ es un ideal no cero de $A$,
sea $\rho^{\prime}={\eu A}\star \rho$, es decir, $\rho_{\eu A}\rho_a=\rho^{\prime}_a
\rho_{\eu A}$ para toda $a\in A$.

Puesto que $\rho$ es signo normalizado, 
para $\alpha\in{\eu A}$, $\mu_{\rho}(\alpha)\in
\Fi$. Por tanto el coeficiente l{\'\i}der de $\rho_{\eu A}
\rho_a$ est\'a en $\Fi$ pues $\rho_{\eu A}$ es m\'onico y de ah{\'\i} se 
sigue que $\rho^{\prime}$ es signo normalizado.

Se tiene que $\Aut_{\Ki}(\Ci)$ act\'ua de 
manera natural en $\aditivo {\Ci}$: si $\sigma\in
\Aut_{\Ki}(\Ci)$, $f\in \aditivo {\Ci}$, 
$f=\sum_{i=0}^m \alpha_i\tau^i$ entonces $\sigma f=
\sum_{i=0}^m(\sigma\alpha_i)\tau^i$. Si 
$\rho\in\Drin_A(\Ci)$, $\sigma\in\Aut_{\Ki}(\Ci)$, entonces
$\sigma\rho\colon A\lra \aditivo{\Ci}$ est\'a dado por $(\sigma\rho)_a
=\sigma(\rho_a)$. Se tiene
$\sigma\rho\in \Drin_A(\Ci)$ pues
\begin{gather*}
(\sigma\rho)_{ab}=\sigma(\rho_{ab})
=\sigma(\rho_a\rho_b)=\sigma(\rho_a)\sigma(\rho_b)=
(\sigma\rho)_a(\sigma\rho)_b,\\
(\sigma\rho)_{a+b}=\sigma(\rho_{a+b})=\sigma(\rho_a+\rho_b)=\sigma(\rho_a)+\sigma(\rho_b)=
(\sigma\rho)_a+(\sigma\rho)_b,\\
D\circ (\sigma\rho)=\sigma\delta=\delta.
\end{gather*}

Adem\'as si ${\eu A}$ es un ideal no cero de $A$, $(\sigma\rho)_{\eu A}$,
el cual es m\'onico, es tal que
$(\sigma\rho)_{\eu A}(\sigma\rho)_a=\rho^{\prime\prime}_a(\sigma\rho)_{\eu A}$ para toda
$a\in A$. Por tanto si $\rho_{\eu A}\rho_a=\rho_a^{\prime}\rho_{\eu A}$, es decir,
$\rho^{\prime}={\eu A}\star\rho$, $\sigma(\rho_{\eu A}\rho_a)=\sigma\rho_{\eu A}(\sigma\rho)_a
=(\sigma\rho^{\prime})_a\sigma\rho_{\eu A}$.

Adem\'as $\sigma\rho_{\eu A}$ es generador de $\sigma(R\rho_{\eu A})
=\sigma(R)\sigma\rho_{\eu A}
=R\sigma\rho_{\eu A}$. Por lo tanto $\sigma\rho_{\eu A}=(\sigma\rho)_{\eu A}$ y
$(\sigma\rho)_{\eu A}(\sigma\rho)_a=(\sigma\rho^{\prime})_a(\sigma\rho)_a$. Se
sigue que $\rho^{\prime\prime}_a = (\sigma\rho^{\prime})_a$. Esto es,
$\rho^{\prime\prime}={\eu A}\star\sigma\rho=\sigma({\eu A}\star\rho)=\sigma\rho^{\prime}$.
Hemos obtenido
\begin{gather}\label{DrinfeldEq3.1.26'}
\text{\fbox{$\sigma({\eu A}\star \rho)={\eu A}\star\sigma\rho$}}.
\end{gather}

Antes de continuar, veamos algunas propiedades de $\rho_{\eu A}$ y $\star$.

\begin{proposicion}\label{DrinfeldP3.1.27}
Sea $\rho\in \Drin_A(F)$ un $A$--m\'odulo de Drinfeld arbitrario y sean ${\eu A}$ y
${\eu B}$ dos ideales no cero de $A$. Entonces
\las
\item $\rho_{{\eu A}{\eu B}}=({\eu B}\star \rho)_{\eu A}\rho_{\eu B}=({\eu A}\star\rho)_{\eu B}
\rho_{\eu A}$.

\item ${\eu A}\star({\eu B}\star \rho)=({\eu A}{\eu B})\star\rho$.

\item $D(\rho_{{\eu A}{\eu B}})=D(({\eu B}\star \rho)_{\eu A}) D(\rho_{\eu B})$.
\end{list}
\end{proposicion}

\begin{proof} (1) $R\rho_{{\eu A}{\eu B}}=\langle 
R\rho_{ab}\mid a\in {\eu A}, b\in {\eu B}\rangle$.
Sea $x\in R\rho_{{\eu A}{\eu B}}$, entonces
\begin{align*}
x&\underbracket[0pt]{=}_{\substack{\uparrow\\ r_t\in R, a_t\in
{\eu A},\\ b_t\in{\eu B}}}
\sum_{t}r_t\rho_{a_tb_t}=\sum_t r_t\rho_{a_t}\rho_{b_t}
=\sum_t \underbrace{r_t \rho_{b_t}}_{
\substack{\inmenosnoventa\\ R\rho_b}}\rho_{a_t}\\
&=\sum_t r^{\prime}_t\rho_{\eu B}\rho_{a_t}
\underbracket[0pt]{=}_{\substack{\uparrow\\
{\eu B}\star \rho=\rho^{\prime}}} \sum_t\underbrace{r_t^{\prime}\rho^{\prime}_{a_t}}_{
\substack{\inmenosnoventa\\ R\rho^{\prime}_{\eu A}}}\rho_{\eu B}=\sum_t r_t^{\prime\prime}\rho^{\prime}_{\eu A}
\rho_{\eu B}=\Big(\sum_{t}r_t^{\prime\prime}\Big)\rho_{\eu A}^{\prime}\rho_{\eu B}.
\end{align*}

Por lo tanto $R\rho_{{\eu A}{\eu B}}\subseteq R\rho^{\prime}_{\eu A}\rho_{\eu B}$.

Rec{\'\i}procamente, sea $y\in R\rho^{\prime}_{\eu A}\rho_{\eu B}$, entonces 
\begin{align*}
y&=\underbrace{s\rho^{\prime}_{\eu A}}_{\substack{\inmenosnoventa\\ 
R\rho_{\eu A}^{\prime}}}\rho_{\eu B}=
\sum_l s_l\underbracket[0pt]{\rho_{a_l}^{\prime}}_{\substack{\uparrow\\
a_l\in{\eu A}}}\rho_{\eu B}
\underbracket[0pt]{=}_{\substack{\uparrow\\ {\eu B}
\star \rho=\rho^{\prime}}}\sum_l\underbrace{
s_l\rho_{\eu B}}_{\substack{\inmenosnoventa\\ R\rho_{\eu B}}}\rho_{a_l}\\
&= \sum_{l,i}s_{l,i}\underbrace{\rho_{b_i}}_{\substack{\inmenosnoventa\\ 
{\eu B}}}\rho_{a_l}=\sum_{l,i}
s_{l,i}\rho_{b_i a_l}=\sum_{l,i}\underbrace{s_{l,i}
\rho_{a_lb_i}}_{\substack{\inmenosnoventa\\ R\rho_{{\eu A}{\eu B}}}}
=s^{\prime}\rho_{{\eu A}{\eu B}}.
\end{align*}

Se sigue que \fbox{$R\rho^{\prime}_{\eu A}\rho_{\eu B}\subseteq 
R\rho_{{\eu A}{\eu B}}$}.

Obtenemos que $R\rho^{\prime}_{\eu A}\rho_{\eu B}=R\rho_{{\eu A}{\eu B}}$ 
y puesto que
$\rho_{\eu A}^{\prime}\rho_{\eu B}$ es m\'onico, se sigue que 
\[
\text{\fbox{$\rho_{{\eu A}{\eu B}}=\rho^{\prime}_{\eu A}\rho_{\eu B}
=\big({\eu B}\star\rho\big)_{\eu A}\rho_{\eu B}$}}.
\]

Puesto que ${\eu A}{\eu B}={\eu B}{\eu A}$, $\rho_{{\eu A}{\eu B}}=\rho_{{\eu B}{\eu A}}=
\big({\eu A}\star \rho\big)_{\eu B}\rho_{\eu A}$.

\noindent
(2) Sean $\left.\begin{array}{l} {\eu B}\star \rho=\rho^{\prime}\\ {\eu A}\star \rho^{\prime}=\rho^{\prime\prime}\\
({\eu A}{\eu B})\star \rho=\rho^{\prime\prime\prime}\end{array}\right\}$ es decir
$\left\{\begin{array}{l} \rho_{\eu B}\rho_a=\rho_a^{\prime}\rho_{\eu B} \text{\ para toda $a\in A$}\\
\rho_{\eu A}^{\prime}\rho_a^{\prime}=\rho_a^{\prime\prime}\rho_{\eu A}^{\prime} \text{\ para toda $a\in A$}\\
\rho_{{\eu A}{\eu B}}\rho_a=\rho_a^{\prime\prime\prime}\rho_{{\eu A}{\eu B}} \text{\ para toda $a\in A$}.
\end{array}\right.$

De $\rho_{\eu B}\rho_a=\rho_a^{\prime}\rho_{\eu B}$ para toda $a\in A$ obtenemos
\[
\rho_{{\eu A}{\eu B}}\rho_a=({\eu B}\star\rho)_{\eu A}\rho_{\eu B}
\rho_a=\rho_{\eu A}^{\prime}
\rho_{\eu B}\rho_a=\rho_{\eu A}^{\prime}\rho_a^{\prime}\rho_{\eu B}
=\rho_a^{\prime\prime}
\rho_{\eu A}^{\prime}\rho_{\eu B}=\rho_a^{\prime\prime}
\rho_{{\eu A}{\eu B}}.
\]

Se sigue que $\rho_{{\eu A}{\eu B}}\rho_a=\rho_a^{\prime\prime}
\rho_{{\eu A}{\eu B}}$ por lo
que $\rho^{\prime\prime}=\rho^{\prime\prime\prime}$. Por tanto
\[
\rho^{\prime\prime}=({\eu A}\star \rho^{\prime})={\eu A}\star ({\eu B}\star \rho)=
({\eu A}{\eu B})\star \rho=\rho^{\prime\prime\prime}.
\]

Finalmente obtenemos \fbox{${\eu A}\star({\eu B}\star\rho)=({\eu A}{\eu B})\star\rho$}.

\noindent
(3) Se tiene $D(\rho_{{\eu A}{\eu B}})=D\big[({\eu B}\star\rho)_{\eu A}\rho_{\eu B}\big]=
D(({\eu B}\star \rho)_{\eu A})D(\rho_{\eu B})$. $\fin$
\end{proof}

Volviendo a nuestra exposici\'on, puesto que 
todo $\rho\in\Drin_A(\Ci)$ corresponde a una
red $\Gamma$ con $\rho=\rho^{\Gamma}$, veamos 
los morfismos de redes. Se tiene

\begin{teorema}\label{DrinfeldT3.1.28}
Se tiene que $\xi\in\Ci$ es un isomorfismo entre $\rho^{\Gamma}$ y $\rho^{\Gamma^{\prime}}$, es decir,
$\rho_a^{\Gamma^{\prime}}=\xi \rho_a^{\Gamma}\xi^{-1}$, si y solamente si 
$\Gamma^{\prime}=\xi \Gamma$.
\end{teorema}

\begin{proof} Recordemos  que ten{\'\i}amos para $a$ y $\Gamma$ lo siguiente
(ver Ecuaci\'on (\ref{DrinfeldEq2.2.20}))
\begin{gather}
e_{a^{-1}\Gamma}(u)=a^{-1}e_{\Gamma}(au)\e\text{o}\e ae_{a^{-1}\Gamma}(u)
=e_{\Gamma}(au)\nonumber\\
\rho_a^{\Gamma}=au\prod_{\lambda\in \frac{a^{-1}
\Gamma}{\Gamma}\setminus\{0\}}\Big(1-\frac{u}
{e_{\Gamma}(\lambda)}\Big)\label{DrinfeldEq3.1.28'}\\
e_{\Gamma}(au)=\rho_a^{\Gamma}(e_{\Gamma}(u)).\nonumber
\end{gather}
\
Apliquemos (\ref{DrinfeldEq3.1.28'}). Sea $\rho_a^{\Gamma^{\prime}}
=\xi\rho_a^{\Gamma}\xi^{-1}$. Por
tanto
\begin{align*}
\rho_a^{\Gamma^{\prime}}(u)&=(\xi\rho_a^{\Gamma}\xi^{-1})(u)
=\xi\rho_a^{\Gamma}(\xi^{-1}u)
=\xi\Big[a\xi^{-1}u\prod_{\lambda\in\frac{a^{-1}
\Gamma}{\Gamma}\setminus\{0\}}\Big(
1-\frac{\xi^{-1}u}{e_{\Gamma}(\lambda)}\Big)\Big]\\
&= au\prod_{\lambda\in\frac{a^{-1}
\Gamma}{\Gamma}\setminus\{0\}}\Big(1-\frac{u}{\xi e_{\Gamma}(\lambda)}\Big).
\end{align*}

Usando (\ref{DrinfeldEq3.1.28'}) con $\Gamma_1:=\xi \Gamma$, esto es, $\Gamma=\xi^{-1}
\Gamma_1$, $a=\xi$, $u=\lambda$:
\begin{gather*}
e_{\xi^{-1}\Gamma_1}(\lambda)=\xi^{-1}e_{\Gamma_1}(\xi \lambda)\Lra
\xi e_{\Gamma}(\lambda)=e_{\xi\Gamma}(\xi\lambda),
\intertext{por lo que}
\begin{align*}
(\xi\rho_a^{\Gamma}\xi^{-1})(u)&=au\prod_{\lambda\in\frac{a^{-1}
\Gamma}{\Gamma}\setminus\{0\}}\Big(1-\frac{u}{\xi e_{\Gamma}(\lambda)}\Big)\\
&\underbracket[0pt]{=}_{\substack{\uparrow\\
\lambda\in\frac{a^{-1}\Gamma}{\Gamma}\iff\xi \lambda\in
\frac{(a^{-1}\xi \Gamma)}{\xi \Gamma}}}
au\prod_{\xi\lambda\in\frac{a^{-1}(\xi
\Gamma)}{(\xi\Gamma)}\setminus\{0\}}\Big(1-\frac{u}{\xi e_{\xi\Gamma}(\xi\lambda)}\Big)\\
&=\rho_a^{\xi\Gamma}(u)=\rho_a^{\Gamma^{\prime}}(u).
\end{align*}
\end{gather*}

Por tanto \fbox{$\xi\Gamma=\Gamma^{\prime}$}.

El rec{\'\i}proco es an\'alogo. $\fin$
\end{proof}

\begin{definicion}\label{DrinfeldD3.1.29}
Dadas dos redes $\Gamma_1$, $\Gamma_2$, si $c\in \Ci$ es tal que $c\Gamma_1
\subseteq \Gamma_2$, entonces $f\colon \Gamma_1\lra \Gamma_2$, $f(x)=cx$ es
un $A$--homomorfismo. Se define $\Hom(\Gamma_1,\Gamma_2)=\{c\in\Ci\mid
c\Gamma_1\subseteq \Gamma_2\}$.
\end{definicion}

Sea $\Red_A(\Ci)$ el conjunto de las $A$--redes de $\Ci$.\label{Drinfeldredes}
Entonces el teorema de
uniformizaci\'on anal{\'\i}tica, Teorema \ref{DrinfeldT2.2.22},
establece que el mapeo $\Red_A(\Ci)\lra \Drin_A(\Ci)$,
$\Gamma\longmapsto \rho^{\Gamma}$, es biyectiva.

\begin{teorema}\label{DrinfeldT3.1.30} Sea $\Gamma,\Gamma^{\prime}\in\Red_A(\Ci)$
dos redes del mismo rango y sea $c\neq 0$, $c\in\Hom(\Gamma,\Gamma^{\prime})$.
Sea $f_c(x)=cP(c^{-1}\Gamma/\Gamma,x)$ (esencialmente igual a
$\rho_c^{\Gamma}(x)$).

Entonces $f_c\in\Hom(\rho^{\Gamma},\rho^{\Gamma^{\prime}})$ y $c\mapsto
f_c$ es un isomorfismo de grupos abelianos y de $\F$ espacios vectoriales de
$\Hom(\Gamma,\Gamma^{\prime})$ con $\Hom(\rho^{\Gamma},\rho^{\Gamma^{
\prime}})$.
\end{teorema}

\begin{proof} Se tiene $e_{c^{-1}\Gamma^{\prime}}(u)=P(c^{-1}\Gamma^{\prime}/
\Gamma,e_{\Gamma}(u))=c^{-1}e_{\Gamma^{\prime}}(cu)$. En particular
\[
e_{\Gamma^{\prime}}(cu)=cP(c^{-1}\Gamma^{\prime}/\Gamma,e_{\Gamma}
(u))=f_c(e_{\Gamma}(u)).
\]

Se sigue que $e_{\Gamma^{\prime}}\circ c=f_c\circ e_{\Gamma}$. Por tanto
\[
e_{\Gamma^{\prime}}(ca)=e_{\Gamma^{\prime}}(ac)=\rho_a^{\Gamma^{\prime}}(
e_{\Gamma^{\prime}}(c))=\rho_a^{\Gamma^{\prime}} f_ce_{\Gamma}.
\]

Puesto que $e_{\Gamma}(au)=\rho_a^{\Gamma}(e_{\Gamma}(u))$ y
$e_{\Gamma^{\prime}} c=f_c e_{\Gamma}$, se sigue que
\begin{gather*}
f_c\rho_a^{\Gamma}e_{\Gamma}=f_c\rho_a^{\Gamma}(e_{\Gamma}(1))=
f_c e_{\Gamma}(a)=e_{\Gamma^{\prime}}(ca).
\intertext{Por lo tanto}
\rho_a^{\Gamma^{\prime}}f_c e_{\Gamma}=f_c\rho_a^{\Gamma}e_{\Gamma}
\Lra \text{\fbox{$f_c\rho_a^{\Gamma}=\rho_a^{\Gamma^{\prime}}f_c$}}.
\end{gather*}

Se sigue que $f_c\in\Hom(\rho^{\Gamma},\rho^{\Gamma^{\prime}})$.

Se verifica que $c\mapsto f_c$ es $\F$--lineal. Por otro lado, puesto que
$Df_c=c$, el mapeo es inyectivo.

Ahora sea $f\in\Hom(\rho^{\Gamma},\rho^{\Gamma^{\prime}})$. Si $f=0$, sea
$c=0$ y $f_c=f_0=0=f$. Sea $f\neq 0$. Se tiene $f\rho_a^{\Gamma}=
\rho_a^{\Gamma^{\prime}}f$. Multiplicando por $e_{\Gamma}$, se obtiene
$f\rho_a^{\Gamma}e_{\Gamma}=fe_{\Gamma}a=\rho_a^{\Gamma^{\prime}}
fe_{\Gamma}$. Sea $c:=Df$. Puesto que $f\rho_a^{\Gamma}=\rho_a^{
\Gamma^{\prime}}f$ se sigue que $c\neq 0$. Por tanto
\begin{gather*}
e_{\Gamma^{\prime}}ca=e_{\Gamma^{\prime}}ac=\rho_a^{\Gamma^{\prime}}
(e_{\Gamma}(c))=\rho_a^{\Gamma^{\prime}}e_{\Gamma}c.
\intertext{Se sigue}
(fe_{\Gamma}-e_{\Gamma^{\prime}}c)a=\rho_a^{\Gamma^{\prime}}(
fe_{\Gamma}-e_{\Gamma}c).
\end{gather*}

Puesto que $c=Df$, el coeficiente de $\tau^0$ en $fe_{\Gamma}-e_{\Gamma}c$
es $0$. Por el argumento usado en la \'ultima
parte de la demostraci\'on del Teorema \ref{DrinfeldNT1} 
se sigue $fe_{\Gamma}=e_{\Gamma^{\prime}}c$ y
$f(e_{\Gamma}(u))=e_{\Gamma^{\prime}}(cu)$. 

M\'as precisamente,
sea $g=(fe_{\Gamma}-e_{\Gamma'}c)\in \Ci\torcido$ y $ga=\rho_a^{\Gamma'}g$.
Entonces, para $a\in A$, $\serie g{\infty}b$, $b_0=0$, 
$\rho_a^{\Gamma'}=\sum_{j=0}^{\infty} d_j\tau^j$
\begin{align*}
g(\tau)a&=\sum_{i=0}^{
\infty}b_i\tau^i a=\sum_{i=0}^{\infty} a^{q^i}b_i\tau^i=\big(\sum_{j=0}^m d_j
\tau^j\big)\big(\sum_{i=0}^{\infty} b_i\tau^i\big)\\
&\underbracket[0pt]{\igual}_{\substack{\uparrow\\
d_0=a}} a\sum_{i=0}^{\infty}b_i\tau^i+\sum_{j=1}^m\sum_{i=0}^{\infty}
d_jb_i^{q^j} \tau^{j+i}.
\end{align*}
 Igualando el t\'ermino $l$:
\[
a^{q^i}b_l=ab_l+\Big(\sum_{i+j=l}d_jb_i^{q^j}\Big)=ab_l +\sum_{j=1}^l
d_jb_{l-j}^{q^j}.
\]
Si $b_0=b_1=\cdots=b_{l-1}=0$, $a^{q^i}b_l=ab_l+0$, lo cual implica
que $(a^{q^i}-a)b_l=0$. Por tanto $b_l=0$. Se sigue que $g\equiv 0$
y que
\[
fe_{\Gamma}=e_{\Gamma'}c.
\]

Si $\gamma\in \Gamma$,
$\gamma$ es ra{\'\i}z de $f(e_{\Gamma}(u))$ por lo que $0=e_{
\Gamma^{\prime}}(c\gamma)$. Se sigue que $c\gamma\in\Gamma^{\prime}$
de tal forma que $c\Gamma\subseteq \Gamma^{\prime}$, esto es, $c\in
\Hom(\Gamma,\Gamma^{\prime})$. Para finalizar hay que probar que $f=f_c$.
La demostraci\'on se debe hacer usando que $D(f)=D(f_c)=c$ y que para
toda $a\in A$, 
\[
(f-f_c)\rho_a^{\Gamma}=\rho_a^{\Gamma^{\prime}}(f-f_c)
\]
lo cual se sigue de que $f\rho_a^{\Gamma}=\rho_a^{\Gamma'} f$
para toda $a\in A$ y de que $f_c\rho_a^{\Gamma}=f_c^{\Gamma'} f$
para toda $a\in A$. Entonces, definiendo $g(\tau)=f(\tau)-
\serie{f_c}{\infty}b$ se tiene que $b_0=0$.  Repetimos el argumento anterior
para probar que $g(\tau)\equiv 0$ obteniendo el resultado.
$\fin$
\end{proof}

En particular $\Gamma\cong \Gamma^{\prime}$ si y solamente si
existe $\xi$ tal
que $\Gamma^{\prime}=\xi \Gamma$ si y solamente si
$\rho^{\Gamma^{\prime}}\cong
\rho^{\Gamma}$ si y solamente si
$\xi\rho_a^{\Gamma}\xi^{-1}=\rho_a^{\Gamma^{\prime}}$.

\begin{ejemplo}\label{DrinfeldEj3.1.30'} Sea $C$ 
el m\'odulo de Carlitz. Entonces
$C=\rho^{\Gamma}$ donde $\Gamma=\xi A$ con $\xi=\tilde{\pi}$ dado
en (\ref{DrinfeldEc2.9'}). Sean $\Gamma':=\xi^{-1}\Gamma=A$ y $\rho^{\Gamma'}=
\rho^{A}$ corresponde a $\xi^{-1}C\xi$, esto es $\rho_T^A=
\xi^{-1} C_T \xi=T+\xi^{-1}\tau\xi=T+\xi^{q-1}\tau=T-\alpha \tau$ con
$\alpha=\prod_{i=1}^{\infty}\big(1-\frac{[i-1]}{[i]}\big)$ como consecuencia
de (\ref{DrinfeldEc2.9'}).
\end{ejemplo}

Sea $\Red_{A,r}(\Ci)$\label{Drinfeldredes1} el 
conjunto de redes de rango $r$ en $\Ci$ y sea
${\mathcal R}_{A,r}$ el conjunto de clases de isomorfismos en $\Red_{A,r}
(\Ci)$\label{Drinfeldredes2}. Se tiene que si ${\eu A}$ es un ideal fraccionario de $A$,
$\Gamma\lra {\eu A}^{-1}\Gamma$ define una acci\'on de $D_A$, los ideales
fraccionarios de $A$ en $\Red_{A,r}(\Ci)$ y si  ${\eu A}=(a)\in P_A$ es principal,
$\Gamma\lra a^{-1}\Gamma\cong \Gamma$, por lo tanto $\pic(A)=\frac{D_A}
{P_A}$ define esta acci\'on en ${\mathcal R}_{A,r}$. Ahora, puesto que cada
red de rango $1$ es isomorfa a un ideal fraccionario de $A$, se sigue:

\begin{proposicion}\label{DrinfeldP3.1.31} ${\mathcal R}_{A,1}$ es un espacio homog\'eneo
principal para $\pic (A)$, es decir, ${\mathcal R}_{A,1}$ tiene una sola \'orbita
en $\pic (A)$ y el estabilizador es $\{1\}$. En otras palabras, 
\[
{\mathcal R}_{A,1}=\pic(A)\cdot \bar{\Gamma}.
\]

En particular $|{\mathcal R}_{A,1}|=|\pic(A)|=h_A$ es finito. $\fin$
\end{proposicion}

Puesto que $\Gamma\lra \rho^{\Gamma}$ es una biyecci\'on, se sigue:

\begin{teorema}\label{DrinfeldT3.1.32} Hay 
exactamente $h_A$ clases de isomorfismos
de $A$--m\'odulos de Drinfeld de rango $1$ sobre $\Ci$ y ${\mathcal D}_1$,
el conjunto de clases de isomorfismos de m\'odulos de Drinfeld de rango $1$
sobre $\Ci$, es un espacio homog\'eneo principal con la acci\'on $\star$: 
$({\eu A}\star\rho)$, esto es, ${\mathcal D}_1=\pic(A)\star \bar{\rho}$. $\fin$
\end{teorema}

Sea ahora $\rho\in{\mathcal H}$ un $A$--m\'odulo de Hayes, es decir de rango
$1$ y con signo normalizado. Sea ${\eu A}$ un ideal no cero de $A$ tal que
${\eu A}\star \rho=\rho$. En especial, si ${\eu A}$ estabiliza a $\rho$,
${\eu A}$ estabiliza a la clase de isomorfismo de $\rho$ m\'as no al rev\'es.

En particular, si $\bar{\eu A}\in\pic(A)$, $\bar{\rho}$ la
clase de isomorfismos de $\rho$, entonces
$\bar{\eu A}\circ \bar{\rho}=\bar{\rho}$, es un elemento en
${\mathcal D}_1$. Puesto que la acci\'on de $\pic(A)$ en ${\mathcal R}_{A,1}$
es transitiva, el estabilizador de $\bar{\rho}$ es el conjunto de ideales 
principales no cero de $A$. En particular ${\eu A}=(x)$ es principal.
Ahora sea ${\eu A}\star \rho=\rho^{\prime}$.
Entonces se tiene $\rho_{\eu A}\rho_a=\rho_a^{\prime}
\rho_{\eu A}$ para toda $a\in A$ y 
\begin{gather*}
\rho_{\eu A}=\mu_{\rho}(x)^{-1}\rho_x,\e
\rho_a^{\prime}=({\eu A}\star\rho)_a.\\
\intertext{Por tanto}
\rho_{\eu A}\rho_a=\mu_{\rho}(x)^{-1}\rho_x\rho_a=\mu_{\rho}(x)^{-1}
\rho_{xa}=\mu_{\rho}(x)^{-1}\rho_a\rho_x
=\rho_a^{\prime}\rho_{\eu A}=\rho_a^{\prime}\mu_{\rho}(x)^{-1}\rho_x\\
\Lra \mu_{\rho}(x)^{-1}\rho_a=\rho_a^{\prime}\mu_{\rho}(x)^{-1}\Lra
\rho_a^{\prime}=\mu_{\rho}(x)^{-1}\rho_a\mu_{\rho}(x).
\intertext{Se sigue que}
({\eu A}\star \rho)_a=\mu_{\rho}(x)^{-1}\rho_a\mu_{\rho}(x),
\end{gather*}
donde ${\eu A}=xA$.

Sea $\xi\in\Aut(\rho)$,
$c_i(\xi\rho\xi^{-1},a)=\xi^{1-q^i}c_i(\rho,a)=c_i(\rho,a)$ para toda
$i$ y para toda $a$. Se sigue que $\xi^{1-q}=1$, es decir, $\xi^q=
\xi$ lo cual implica que $\xi\in\*{\F}$. Esto es, $\Aut(\rho)=\F$.
Entonces ${\eu A}\star\rho=\rho$ implica que $\mu_{\rho}(x)\in
\Aut(\rho)=\*{\F}$. 
Volviendo a nuestro estudio, 
$\mu_{\rho}(x)\in\*{\F}$ y por tanto el estabilizador de $\rho$ es
$\{xA\mid x\in A,\mu_{\rho}(x)\in \*{\F}\}=\{x A\mid x\in A, \sgn(x)=1\}$
pues $xA=(\mu_{\rho}(x)^{-1}x)A$ ya que $\mu_{\rho}(x)\in \*{\F}
\subseteq A$.

Si ${\eu A}\in D_A$ estabiliza a $\rho$, entonces
${\eu A}=xA$ es principal pues
$\pic(A)=\frac{D_A}{P_A}$ act\'ua transitivamente sobre ${\mathcal D}_1$
y como $({\eu A}\star\rho)_a=\rho_a^{\prime}=\mu_{\rho}(x)^{-1}
\rho_a\mu_{\rho}(x)$, entonces $\mu_{\rho}(x)\in\*{\F}$ y existe
$y\in A$ con $\sgn y=1$, $y{\eu A}\in A$. Por tanto la acci\'on de los 
ideales de $A$ se extiende a $D_A$ y el estabilizador de $\rho$ es
\[
P_A^{+}=\{xA\mid x\in\*\K ,\sgn x=1\}\subseteq P_A.
\]

\begin{definicion}\label{DrinfeldD3.1.33} El grupo $\pic^+(A)=\frac{D_A}{P_A^+}$
se llama el {\em grupo de clases extendido\index{grupo de clase 
extendido} de $A$ relativo a $\sgn$}.
\end{definicion}

Se tiene la sucesi\'on exacta
\[
1\lra \frac{P_A}{P^+_A}\lra \frac{D_A}{P_A^+}\lra \frac{D_A}{P_A}\lra 1
\]
y se tiene $\Big|\frac{P_A}{P_A^+}\Big|=\frac{q^{\di}-1}{q-1}=$ 
n\'umero de diferentes funciones signo. De hecho $\frac{P_A}{P_A^+}
\underbracket[0pt]{\cong}_{\substack{\uparrow\\ \sgn}}\frac{
\*{\Fi}}{\*{\F}}$, y
\[
h_A^+=|\pic^+(A)|=\frac{q^{\di}-1}{q-1}|\pic(A)|=\frac{q^{\di}-1}{q-1}h_A.
\]

Ahora bien, puesto que el conjunto de 
los m\'odulos de Hayes ${\mathcal H}$ es
esencialmente igual al conjunto 
$\pic^+(A)\cdot \rho$, se sigue que

\begin{teorema}\label{DrinfeldD3.1.34} El conjunto de $A$--m\'odulos de Hayes es
un espacio homog\'eneo para $\pic^+(A)$ con la acci\'on de $\star$ y
$|{\mathcal H}|=h_A^+=|\pic^+(A)|$. $\fin$
\end{teorema}

\subsection{El campo de clase de Hilbert extendido $H_A^+$
con respecto a $\sgn$}\label{DrinfeldS3.2}

\begin{definicion}\label{D15.4.38'}
Sea $\rho\in{\mathcal H}$ y sea $y\in A$, $y\notin\F$. Se
define $H_A^+$ como
el campo generado sobre $\K $ por los coeficientes de $\rho_y$.
\end{definicion}

Como vimos para $\K_{\rho}$, $H_A^+$ es independiente de $y$
(ver la demostraci\'on del Teorema \ref{DrinfeldT3.1.6}). Todos
los $A$--m\'odulos ${\eu A}\star \rho$, ${\eu A}$ ideal no cero de $A$,
est\'an definidos sobre $H_A^+$,
lo cual se sigue de que si ${\eu A}\star \rho=\rho'$,
entonces $\rho'_a=\rho_{\eu A}\rho_a \rho_{\eu A}^{-1}$.
Adem\'as este conjunto $\{{\eu A}\star
\rho\}_{{\eu A}\neq 0}$ se compone de todos los $A$--m\'odulos
signo normalizados. Se sigue que $H_A^+/\K $ es independiente de la
elecci\'on de $\rho$, pero si depende de la funci\'on signo $\sgn$ seleccionada.

\begin{definicion}\label{DrinfeldD3.2.1}
El campo $H_A^+$ se llama {\em campo normalizador\index{campo
normalizador}} para $A$--m\'odulos de Drinfeld sobre $(\K ,\p,\sgn)$.
\end{definicion}

Para $\sigma\in\Aut_\K (\Ci)$, $\sigma\rho$ son signos normalizados, por 
lo que est\'an definidos sobre $H_A^+$. Se sigue que $H_A^+$
contiene a todos los conjugados de su conjunto 
finito de generadores (ver la Ecuaci\'on (\ref{DrinfeldEq3.1.26'})). Por
tanto $H_A^+/\K $ es una extensi\'on finita y normal.

\begin{teorema}\label{DrinfeldT3.2.2}
La extensi\'on $H_A^+/\K $ es una extensi\'on de Galois y 
se tiene que $\Gal(H_A^+/\K )\subseteq
\pic^+(A)$. En particular $H_A^+/\K $ es una extensi\'on abeliana.
\end{teorema}

\begin{proof} Se tiene que $\K_{\rho}\subseteq H_A^+$ y $\K_{\rho}/\K $
es una subextensi\'on algebraica de $\Ki/\K $ puesto que $\K_{\rho}
\subseteq \Ki$. Veamos que cada 
elemento de $\Ki$ que es algebraico sobre $\K $
es separable, por lo que se seguir\'a que 
$\K_{\rho}/\K $ es una extensi\'on separable.

En efecto si $\alpha\in\Ki$ es algebraico sobre $\K $, sea $L=\K (\alpha)$.
Podemos suponer que $\alpha$ es normal, esto es considerando
la cerradura normal $\tilde{L}$ de $L/\K $ y suponiendo $\tilde{L}=L$.
Si $L/\K $ no fuese separable, se tendr\'ia una subextensi\'on
$\K \subseteq E \subsetneqq L$ con $L/E$ una extensi\'on puramente
inseparable. Sobre los primos infinitos de $E$, existe para cada
uno un \'unico primo en $L$. Por otro lado, puesto que $\K (\alpha)
\subseteq \Ki$, se tiene $\K (\alpha)_{\infty}=\Ki$ y $\p$ se descompone
totalmente en $L/\K $. Por tanto $E=L$ lo cual es absurdo y de
donde se sigue que $L/\K $ es una extensi\'on separable.

Sea $\xi\in\Ci$ un isomorfismo de $\rho$ con $\rho^{\prime}$, esto es,
$\rho^{\prime}=\xi\rho\xi^{-1}$ tal que $\rho^{\prime}$ est\'a definido
sobre $\K_{\rho}$. Sea $x\in A$ no constante y positivo. Entonces,
con $r=1$ y puesto que $x$ es positivo, $\mu_{\rho}(x)=\sigma(
\sgn(x))=\sigma(1)=1$ (ver Observaci\'on \ref{DrinfeldO3.1.22}) se tiene
\[
c_{r\deg x}(\rho^{\prime},x)=\xi^{1-q^{r\deg x}}\mu_{\rho}(x)=\mu_{\rho^{\prime}}
(x)\in \K_{\rho},
\]
por tanto $\K_{\rho}(\xi)/\K_{\rho}$ es separable pues $\mcd(1-q^{r\deg x},
p)=1$. Puesto que $\K_{\rho}\incientoochenta c_i(\rho^{\prime},y)=\xi^{1-q^i}
c_i(\rho,y)\in H_A^+$, se tiene $H_A^+\subseteq \K_{\rho}(\xi)$ y por
tanto $H_A^+/\K $ es Galois.

Se tiene para $\sigma\in\Aut_\K (\Ci)$, ${\eu A}\star \sigma\rho=\sigma(
{\eu A}\star \rho)$ (ecuaci\'on (\ref{DrinfeldEq3.1.26'})), 
por lo que la acci\'on de $\Gal(H^+_A/\K )$ conmuta con
la acci\'on $\pic^+(A)$. Sea $\theta\colon \Gal(H_A^+/\K )\lra \pic^+(A)$ dada
como sigue. Si $\sigma\in\Gal(H^+_A/\K )$, $\theta(\sigma):=\bar
{\eu A}_{\sigma}$ donde ${\eu A}_{\sigma}$ satisface $\sigma\rho=
{\eu A}_{\sigma}\star \rho$. 

Notemos que $\theta$ est\'a bien definida pues todo $\sigma\rho$ es de la
forma ${\eu A}\star \rho$ para alg\'un ${\eu A}\in D_A$
y el estabilizador de $\rho$ es $P_A^+$. Adem\'as es claro que $\theta$
es un homomorfismo de grupos. 

El homomorfismo
$\theta$ es inyectivo puesto que si $\sigma
\neq \Id$, $\sigma\rho\neq \rho$. $\fin$
\end{proof}

Para estudiar la ramificaci\'on en $H_A^+/\K $
en un primo infinito $\pK$, necesitamos considerar el grupo
de inercia el cual est\'a relacionado con el mapeo de reducci\'on 
m\'odulo $\pK$
para un lugar $\pK$. Sea $B^+$ la cerradura entera de $A$ en $H_A^+$.
\[
\xymatrix{
B^+\ar@{-}[r]\ar@{-}[d]&H_A^+\ar@{-}[d]\\A\ar@{-}[r]&\K }
\]

En general, sea $\rho\in\Drin_A(F)$ de rango $r$. Supongamos que $F$
es un campo con una valuaci\'on discreta $v$ y que todos los coeficientes
de $\rho_a$ son enteros con respecto a $v$, 
esto es, $v(\alpha)\geq 0$ para $\alpha$ coeficiente de
$\rho_a$. Sea $\o_v$ el anillo de valuaci\'on de $v$
con ideal m\'aximo $\pK$ y campo residual $F(\pK)=\o_v/\pK$. Tomamos
los coeficientes de $\rho_a\bmod \pK$ y denotamos esta reducci\'on por
$\rho^{(\pK)}$. En general $\rho^{(\pK)}$ no es un $A$--m\'odulo de 
Drinfeld si todos los t\'erminos no constantes son congruentes a $0
\bmod \pK$, aunque de cualquier forma $\rho^{(\pK)}\colon
A\lra \aditivo{F(\pK)}$
sigue siendo un homomorfismo de $\F$--\'algebras y $\delta^{(\pK)}\colon
A\lra F(\pK)$ es tal que $\delta^{(\pK)}=\delta \bmod \pK$.

\begin{definicion}\label{DrinfeldD3.2.3}
Se dice que $\rho$ tiene {\em reducci\'on estable\index{reducci\'on
estable} en $\pK$} si existe un $A$--m\'odulo de Drinfeld $\rho^{\prime}
\in \Drin_A(F)$ isomorfo a $\rho$ tal que todos los coeficientes de 
$\rho_a^{\prime}$ son enteros en $v$ para toda $a\in A$ y $\rho^{\prime
(\pK)}$ es un $A$--m\'odulo de Drinfeld $\rho^{\prime (\pK)}\in
\Drin_A(F(\pK))$.

Decimos que $\rho$ tiene {\em buena reducci\'on\index{buena 
reducci\'on} en $\pK$} si adem\'as $\rho^{\prime (\pK)}$ tiene rango $r$,
el mismo que el de $\rho$.
\end{definicion}

\begin{observacion}\label{DrinfeldO3.2.4} Si $\rho$ tiene rango $1$, entonces toda
reducci\'on estable es buena reducci\'on.
\end{observacion}

El resultado clave es, a\'un si $\rho$ no tiene reducci\'on estable (resp. buena
reducci\'on), existe una extensi\'on $F^{\prime}$ de $F$ tal que $\rho$
tiene reducci\'on estable (resp. buena reducci\'on) en $F^{\prime}$.

\begin{definicion}\label{DrinfeldD3.2.5} Decimos que $\rho$ tiene {\em reducci\'on
estable potencial\index{reducci\'on estable potencial}} (resp. {\em
buena reducci\'on potencial\index{buena reducci\'on potencial}}) si
existe $F\subseteq F^{\prime}$ tal que $\rho$ tiene reducci\'on estable
(resp. buena reducci\'on) sobre $F'$.
\end{definicion}

\begin{ejemplo}\label{DrinfeldEj3.2.6}
Si $A=R_T=\F[T]$, entonces para toda $r>1$, el m\'odulo de Drinfeld
$\rho_T=T+\tau+a_2\tau^2+\cdots+a_{r-1}\tau^{r-1}+T\tau^r$,
$\rho\colon A\lra \aditivo \K $, tiene
reducci\'on estable pero no buena reducci\'on en $v_T$ sobre
$\coc A=\K $ pues 
\[
\rho_T\bmod T=T+\tau+\sum_{i=2}^r a_i\tau^i
\bmod T = \tau+\sum_{i=2}^r a_i\tau^i \bmod T \neq 0\bmod T. 
\]
Si $\rho'$ es isomorfo a $\rho$, $\rho'=\xi\rho \xi^{-1}$ con
$\xi\in \aditivo \K $. Entonces
$c_r(\rho',T)=\xi T\xi^{-q^r}=\xi^{1-q^r} T$ y se tendr\'a $v_T(\xi^{1-q^r})
=(1-q^{r})v_T(\xi)=-1$ si y solamente si $v_T(\xi)=\pm 1$, 
$1-q^{r}=\mp 1$ si y solamente si $q^{r}=1\pm 
1=2$ o $0$ lo cual es imposible pues $r>1$.

Similarmente
el $A$--m\'odulo de Drinfeld $\varphi_T=T+T\tau+T\tau^2+\cdots
+T\tau^r$ no tiene reducci\'on estable en $v_T$ (por la misma raz\'on 
anterior).
\end{ejemplo}

\begin{teorema}\label{DrinfeldT3.2.7} Toda $A$--m\'odulo de Drinfeld $\rho$
sobre un campo con valuaci\'on discreta $v$ tiene reducci\'on estable
potencial. En particular, si $\rho$ es de rango $1$, $\rho$ tiene buena
reducci\'on potencial.
\end{teorema}

\begin{proof} Sea $a\in A$, $\rho_a=\sum_i a_i\tau^i$ y definimos
\[
\gamma_a:=\min_{i>0}\frac{v(a_i)}{q^i-1}.
\]

Sean $x_1,\ldots,x_s$ un conjunto de generadores de $A$ como 
$\F$--\'algebra y sea $\gamma:=\min_{1\leq j\leq s} \gamma_{x_j}$. Si
$\gamma=0$, entonces hay un elemento $\rho_a^{(\pK)}$ de grado
$>0$ y no hay nada que hacer. Supongamos $\gamma\neq 0$, es
decir, $\rho_a^{(\pK)}=\delta^{(\pK)}(a)$ para toda $a\in A$.

Sean $i,j$ tales que 
$\gamma=\frac{c_i(\rho,x_j)}{q^i-1}$. Sea $F^{\prime}$ una extensi\'on
de $F$ con ramificaci\'on de grado $q^i-1$ en $v$. Sea $x\in F^{\prime}$ 
con $v^{\prime}(x)=c_i(\rho,x_j)=(q^i-1)\gamma$. Sea $\rho^{\prime}=x
\rho x^{-1}$ y $c_l(\rho^{\prime},a)=x^{1-q^l}c_l(\rho,a)$. Por tanto
\begin{align*}
v^{\prime}(c_i(\rho^{\prime},x_j))&=(1-q^i)v^{\prime}(x)+v^{\prime}(
c_i(\rho,x_j))\\
&=(1-q^i)(q^i-1)\gamma+e(v^{\prime}|v)v(c_i(\rho,x_j))\\
&=(1-q^i)(q^i-1)\gamma +(q^i-1)(q^i-1)\gamma=0.
\end{align*}

Por lo tanto $\rho^{\prime}_{x_j}$ es una $A$--m\'odulo de Drinfeld.
$\fin$
\end{proof}

Regresamos a nuestra discusi\'on donde $\rho$ es signo normalizado
de rango $1$, el cual est\'a definido en $B^+$ y puede ser reducido para 
cada ideal ${\eu q}$ de $B^+$. Sea $\pi_{\eu q}\colon B^+\lra B^+/{\eu q}$
el mapeo de reducci\'on y sea $\pK={\eu q}\cap A$. Recordemos que $B^+$
es la cerradura entera de $A$ en $H_A^+$.

\begin{proposicion}\label{DrinfeldP3.2.8} El mapeo de 
reducci\'on $\rho\longmapsto
\pi_{\eu q}\circ \rho$, $\pi_{\eu q}\circ 
\rho\colon A \lra \aditivo {H_A^+({\eu q})}$
es inyectivo sobre ${\mathcal H}$ el conjunto de $A$--m\'odulos de Hayes.
\end{proposicion}

\begin{proof} Supongamos que $\rho,\rho^{\prime}\in{\mathcal H}$ y que ambos
se reducen m\'odulo $\eu q$
 al mismo $\phi\in\Drin_A(B^+/{\eu q})$, esto es, $\phi=
\rho^{(\pK)}=\rho^{\prime (\pK)}$.

Sea $\eu A$ ideal de $A$ tal que $\rho^{\prime}={\eu A}\star \rho$. Veremos
que se puede suponer que ${\eu A}$ y ${\eu q}$ son primos relativos.
Por el teorema de aproximaci\'on de Artin (ver \cite[Theorem 2.5.3]{Vil2006}),
podemos hallar
$z\in \K $ tal que $z\equiv 1\bmod\p$ y $v_{\eu q}(z)=-v_{\eu q}({\eu A})$.
Notemos que $\sgn(z)=1$ por lo que $z\in P_A^+$.

Sea $z\cdot {\eu A}=\frac{\eu B}{\eu C}$ con ${\eu B}$ y ${\eu C}$ primos
relativos. El ideal ${\eu B}{\eu C}^{h^+_A-1}$ es primo relativo a ${\eu q}$
pues $v_{\eu q}({\eu B})=v_{\eu q}({\eu C})=0$ ya que $v_{\eu q}(z{\eu A})
=0$. Se tiene ${\eu C}^{h_A^+}= (z^{\prime})\in P_A^+$, 
y ${\eu B}{\eu C}^{h_A^+-1}=\frac{\eu B}{\eu C}{\eu C}^{h^+_A}=
z{\eu A}\cdot z^{\prime}=zz^{\prime}{\eu A}=z^{\prime\prime}{\eu A}$.

Sea $\frac{\eu B}{\eu C}\eu C^{h^+_A-1}=\eu D=z'' {\eu A}$. Puesto que 
$\eu D\equiv \eu A\bmod P_A^+$, $\eu D$ satisface $\eu D\star \rho=
z^{\prime\prime}\eu A\star \rho$. Ahora bien, puesto que
$\sgn z^{\prime\prime}=1$, se tiene
\[
\eu A\star\rho=\rho^{\prime\prime}, \quad z^{\prime\prime}(\eu A\star\rho)=
\mu_{\rho}(z^{\prime\prime})^{-1}\rho_a^{\prime}\mu_{\rho}(z^{\prime\prime})=
\rho_a^{\prime\prime},
\]
por lo que $z^{\prime\prime}\eu A\star \rho=
\eu A\star \rho$.

As\'i, podemos suponer $\rho^{\prime}=\eu A\star\rho$ y $\mcd(\eu q,
\eu A)=1$. Reduciendo la ecuaci\'on $\rho_{\eu A}\rho_x=\rho_x^{\prime}
\rho_{\eu A}$ m\'odulo $\eu q$, obtendremos
\[
\rho_{\eu A}\rho_x\underbracket[0pt]{\equiv}_{\substack{\uparrow\\
\rho^{(\eu q)}=\rho^{\prime (\eu q)}\equiv \rho_x\rho_{\eu A}}}
\rho_x^{\prime}\rho_{\eu A}
\bmod \eu q\e\text{para toda}\e x\in A.
\]

Por tanto se tiene $\pi_{\eu q}(\rho_{\eu A})\in\End(\pi_{\eu q}\circ \rho)=A$
(ver \cite[Corollary 5.14]{Hay92}).
Se sigue que existe $a\in A$ tal que $\rho_{\eu A}\equiv \rho_a\bmod \eu q$.
Puesto que $\mu_{\rho}(\rho_{\eu A})=1$
se tiene que $\mu_{\rho}(a)=1$ y $a$ es positivo.

Si probamos que $\eu A$ es principal entonces, por ser $a$ positivo,
\[
\eu A\star\rho=\mu_{\rho}(a)^{-1}\rho\mu_{\rho}(a)=\rho.
\]

Sea $\eu B=\eu A+aA$. Puesto que $\rho_{\eu A}\equiv \rho_a\mod \eu q$,
$\rho^{(\eu q)}[\eu B]=\rho^{(\eu q)}[\eu A]=\rho^{(\eu q)}[a]$ en
$\overline{B^+/\eu q}$ una cerradura algebraica de $B^+/\eu q$,
se obtiene que $\big|A/\eu A\big|=\big|A/\eu B\big|=\big|A/(a)\big|$ y por
tanto $\eu A=\eu B=aA$. $\fin$
\end{proof}

\begin{proposicion}\label{DrinfeldP3.2.9} La extensi\'on $H_A^+/\K $ es no ramificada
en cada lugar finito $\pK$ de $A$.
\end{proposicion}

\begin{proof} Sea $\sigma\in I$ el grupo de inercia de $\pL/\pK$ donde $\pL$ es un
lugar de $H_A^+$ sobre $\pK$. Entonces $\sigma\rho\equiv \rho\bmod \pL$
para $\rho\in{\mathcal H}$, pues $H_A^+=\K (\text{coeficientes de $\rho_y$})$.
Como el mapeo de reducci\'on es inyectivo en $\mathcal H$, se sigue que 
$\sigma\rho=\rho$ y puesto que $H_A^+$ es generado por los coeficientes
de $\rho_y$ sobre $\K $, se tiene que $\sigma|_{H_A^+}=\Id$. Por lo tanto
$\sigma=\Id$, $I=\{\id\}$ y $\pK$ es no ramificado. $\fin$
\end{proof}

Sea $\pK$ un ideal primo no cero de $A$ y sea $\sigma_{\pK}=\artinp{
H_A^+/\K }{\pK}$ el s\'imbolo de Artin de $\pK$. Si $\eu A=\prod_{i=1}^s
\pK_i^{\alpha_i}$ es un ideal no cero de $A$, se define el s\'imbolo de Artin
$\sigma_{\eu A}$ por 
\[
\sigma_{\eu A}=\prod_{i=1}^s\sigma_{\pK_i}^{\alpha_i}.
\]

Uno de los resultados principales en la teor\'ia de $A$--m\'odulos de
Drinfeld de rango $1$ sobre $\Ci$ es:

\begin{teorema}\label{DrinfeldT3.2.10}
Si $\rho\in \mathcal H$ es un $A$--m\'odulo de Hayes, se tiene
\[
\sigma_{\eu A}\rho=\eu A\star \rho,
\]
es decir, el mapeo de Artin coincide con la acci\'on $\eu A\star \rho$.

En particular 
\[
\Gal(H_A^+/\K )\cong \pic^+(A)\quad \text{y}\quad [H_A^+:\K ]=
\frac{q^{\di}-1}{q-1}h_A=\frac{q^{\di}-1}{q-1}\di h_\K .
\]

En el caso particular del m\'odulo de Carlitz se tiene que, $\rho_T=C_T$, $\di=1$, 
$h_\K =h_{\F(T)}=1$ por lo que $H_A^+=\K =\F(T)$.
\end{teorema}

\begin{proof} Si $\eu A$, $\eu B$ son dos ideales no cero de $A$, tenemos
$\eu A\star(\eu B\star\rho)=(\eu A\eu B)\star \rho$. Por tanto, basta probar
el resultado para $\eu A=\pK$ un ideal primo no cero de $A$.

Sea $\pL$ un divisor primo de $B^+$ sobre $\pK$ y consideremos
el automorfismo de Frobenius $\sigma_{\pL}$ en $\pL$. Se tiene 
$\sigma_{\pL} x\equiv x^{\N\pK}\bmod\pL$ para toda $x\in B^+$.

Sea $\rho^{\prime}=\pK\star\rho$ por lo que $\rho_{\pK}\rho_y=
\rho^{\prime}_y\rho_{\pK}$ para toda $y\in A$. Sea $\varphi=\rho
\bmod \pL$ la reducci\'on de $\rho$ m\'odulo $\pL$. Se tiene
$r_{\varphi}=1$ y $\car (\varphi)=\pL$ y $1\leq h_{\varphi}\leq 
r_{\varphi}$ implica $h_{\varphi}=1$. En particular 
\begin{equation}\label{Ec15.4.16'}
\varphi_{\pK}=\tau^{\deg \pK}. 
\end{equation}
Reduciendo la ecuaci\'on $\rho_{\pK}\rho_y=
\rho_y^{\prime}\rho_{\eu p}$ m\'odulo $\pL$, obtenemos
\begin{gather}\label{DrinfeldEq3.2.10'}
\tau^{\deg \pK}\rho_y\equiv\rho^{\prime}_y\tau^{\deg \pK}\bmod \pL.
\end{gather}

Sean $\rho_y=\sum_{i=0}^{\deg y} a_i\tau^i$, $\rho_y^{\prime}
=\sum_{j=0}^{\deg y}b_j\tau^j$. De la Ecuaci\'on (\ref{DrinfeldEq3.2.10'})
obtenemos 
\begin{equation}\label{DrinfeldEc3.7'}
a_i^{\N(\pK)}\equiv b_i\bmod \pL.
\end{equation}
 Por tanto
\begin{align*}
\rho'_y\rho_{\pK}&= \rho_{\pK}\rho_y=\sum_{i=0}^{\deg y} (\sigma_{\pK} a_i)\tau^i
\underbracket[0pt]{\equiv}_{\substack{\uparrow\\\text{$\sigma_{\pK}$
el Frobenius}}}\sum_{i=0}^{\deg y}a_i^{\N(\pK)}\tau^i\\
&\underbracket[0pt]{\equiv}_{\substack{\uparrow\\ \text{(\ref{DrinfeldEc3.7'})}}}
\sum_{i=0}^{\deg y} b_i\tau^i=(\pK\star \rho)_y \bmod
\pL.
\end{align*}

Puesto que la reducci\'on m\'odulo $\pL$ es inyectiva en $\mathcal H$,
se sigue que $\sigma_{\pK}\rho=\pK\star\rho$. 

El mapeo natural $\varphi\colon\Gal(H_A^+/\K )\lra \pic^+(A)$ es inyectivo.
Dado $\pK$ ideal primo no cero de $A$, se tiene $\sigma_{\pK}\rho=
\pK\star\rho$, por lo que $\varphi(\sigma_{\pK})=\bar{\pK}\in \pic^+(A)$.
Para todo ${\eu A}$ ideal no cero de $A$, 
$\varphi(\sigma_{\eu A})=\bar{\eu A}$.
Por tanto $\varphi$ es suprayectiva. $\fin$
\end{proof}

\begin{corolario}\label{DrinfeldC3.2.11}
Para $x\in\*\K $, se define $\sigma_x$ como el s\'imbolo de Artin
$\sigma_{xA}$ correspondiente al ideal principal $xA$. Por tanto
$\sigma_x\rho=\mu_{\rho}(x)^{-1}\rho\mu_{\rho}(x)$ para $\rho
\in\mathcal H$.
\end{corolario}

\begin{proof} Si $x\in A$, entonces
$\sigma_x\rho=(Ax\star \rho)=\mu_{\rho}(x)^{-1}\rho\mu_{\rho}(x)$. 

Sea $xy^{-1}\in\*\K $ con $x,y\in A$. Entonces
\begin{align*}
\sigma_y\sigma_{xy^{-1}}\rho &=\sigma_x\rho
=\mu_{\rho}(x)^{-1}\rho\mu_{\rho}(x)\\
&= \mu_{\rho}(y^{-1})(\mu_{\rho}(xy^{-1}))^{-1}
\rho \mu_{\rho}(xy^{-1})\mu_{\rho}(y)\\
&=\sigma_y(\mu_{\rho}(xy^{-1})^{-1}
\rho \mu_{\rho}(xy^{-1}))
\end{align*}
lo cual implica que $\sigma_{xy^{-1}}\rho=\mu_{\rho}(xy^{-1})^{-1}
\rho\mu_{\rho}(xy^{-1})$. $\fin$
\end{proof}

A continuaci\'on probamos el ``{\em teorema del ideal 
principal\index{teorema del ideal principal para $H_A^+$}}'' para $H_A^+$.

\begin{teorema}\label{DrinfeldT3.2.12} Sea $w$ una valuaci\'on que corresponde
a un primo no ramificado sobre $\pK$ en una extensi\'on finita $L$ de $\K $.
Si $\rho$ es un m\'odulo de Drinfeld de rango uno definido sobre el anillo 
de valuaci\'on $\o_w$, entonces $\rho_{\pK^{e-1}}(x)$ divide a $\rho_{\pK^e}
(x)$ en $\o_w$ y el cociente $\frac{\rho_{\pK^e}(x)}{\rho_{\pK^{e-1}}(x)}$
es Eisenstein en $\pK$, es decir, si $\frac{\rho_{\pK^e}(x)}{\rho_{\pK^{e-1}}
(x)}=\sum_{i=0}^t a_ix^i$, $a_i\in \pK$, $0\leq i\leq t-1$, $a_0\notin\pK^2$ y
$a_t\notin \pK$.
\end{teorema}

\begin{proof} Sea $e=1$. Puesto que la 
reducci\'on m\'odulo el ideal maximal $\pL$ sobre $\pK$,
tiene caracter\'istica $\pK$, todos, excepto $\mu_{\rho^{(\pK)}}(x)$ pertenecen
al ideal maximal pues el rango y la altura son $1$, es decir, ya hab\'iamos 
probado que $\phi_{\pK}=\tau^{\deg \pK}$ donde 
$\phi=\rho\bmod \pL$ (ver Ecuaci\'on (\ref{Ec15.4.16'})). Ahora
hay que probar que $w(D(\rho_{\pK}))\leq 1$. 

Sea $a\in A$ tal que $w(a)=1$
y sea $(a)=\pK \eu A$ y $\pL/\pK$ es no ramificado. Se tiene $\rho_a
\underbracket[0pt]{=}_{\substack{\uparrow\\ \rho_{aA}=\mu(\rho_a)^{-1}\rho_a
\\=\mu_{\rho}(a)^{-1}\rho_a}} \mu_{\rho}(a)(\pK\star\rho)_{\eu A}\rho_{\pK}$
por lo que $1=w(D(\rho_a))=w(a)\geq w(D(\rho_{\pK}))$.

Ahora para $e>1$, 
\begin{align*}
\rho_{\pK^e}(x)&=[(\pK^{e-1}\star\rho)_{\pK}\rho_{\pK^{e-1}}](x)=\frac{
(\pK^{e-1}\star \rho)_{\pK}(\rho_{\pK^{e-1}}(x))}{\rho_{\pK^{e-1}}(x)}
\rho_{\pK^{e-1}}(x)\\
&=f(\rho_{\pK^{e-1}}(x))\rho_{\pK^{e-1}}(x),
\end{align*}
donde $f(t)=(\pK^{e-1}\star \rho)_{\pK}(t)/t$.

Por el caso $e=1$, si $\rho^{\prime}=\pK^{e-1}\star \rho$, $\rho^{\prime}_{\pK}=
(\pK^{e-1}\star \rho)_{\pK}$ es Eisenstein, por lo que $f(t)$ es Eisenstein
y $\rho_{\pK^{e-1}}(x)\equiv x^{q^{(e-1)\deg\pK}}\bmod \pL$ y $f$ es el
cociente $\rho_{\pK^e}(x)/\rho_{\pK^{e-1}}(x)$. $\fin$
\end{proof}

Por el Teorema \ref{DrinfeldT3.3.4}, el siguiente resultado se aplica en general.

\begin{teorema}[Ideal principal para $H_A^+$]\label{DrinfeldT3.2.13}
Sea $\eu A$ un ideal no cero de $A$, entonces $\con_{A/B^+}
\eu A=\eu AB^+=D(\rho_{\eu A}) B^+$ es principal, donde $D(\rho_{\eu A})$
es el t\'ermino constante de $\rho_{\eu A}$, donde $\rho$ est\'a definido
sobre $B^+$, es decir, $\rho_x\in\aditivo {B^+}$ para toda $x\in A$.
\end{teorema}

\begin{proof} Recordemos que $H_A^+$ es independiente de $\rho$. Ahora bien,
se tiene que $D(\rho_{\eu A\eu B})=D((\eu B\star\rho)_{\eu A})D(\rho_{\eu B})
=D(\rho^{\prime}_{\eu A})D(\rho_{\eu B})$ con $\rho^{\prime}=\eu B
\star \rho$, Proposici\'on \ref{DrinfeldP3.1.27} (3). 
Si probamos el teorema para $\eu A$ y $\eu B$, entonces, puesto que
${\eu A}B^+$ es independiente de $\rho$, se tiene $D(\rho_{
\eu A})B^+=D(\rho'_{\eu A})B^+={\eu A}B^+$. Por lo tanto
\begin{align*}
\con_{A/B^+}(\eu A\eu B)&= (\eu A\eu B)B^+={\eu A}(\eu B B^+)=
\eu A(D(\rho_{\eu B})B^+)\\
&=D(\rho_{\eu B})(\eu AB^+)=D(\rho_{\eu B})D(\rho^{\prime}_{\eu A})B^+
=D(\rho_{\eu A\eu B})B^+.
\end{align*}

As{\'\i}, basta probar el resultado para $\eu A=\pK$ un ideal primo no cero
de $A$. Se tiene que $\rho_{\pK}$ es Eisenstein. Es decir, todos los
coeficientes de $\rho_{\pK}$, excepto el l{\'\i}der, pertenecen a cualquier
ideal $\pL$ sobre $\pK$ y $v_{\pL}(D(\rho_{\pK}))=1$.

El resultado se seguir\'a probando que todo ideal $\eu q$ de $B^+$
que no divide a $\pK$, satisface que $v_{\eu q}(D(\rho_{\pK}))=0$ pues
en este caso $(D(\rho_{\pK}))_{B^+}=\pK B^+$.

Sea $e\geq 1$ tal que $\pK^e=yA$ es principal y sea $\eu t=\pK^{e-1}$.
Entonces, por la Proposici\'on \ref{DrinfeldP3.1.27} (3), se tiene que 
$D(\rho_{\pK^e})=
D(\rho_{{\eu t}\pK})=D((\pK\star \rho)_{\eu t}) D(\rho_{\pK})$. Se sigue que
\[
v_{\eu q}(D(\rho_{\pK^e}))=
v_{\eu q}(D(\pK\star \rho)_{\eu t})+v_{\eu q}(D(\rho_{\pK}))=v_{\eu q}(
\mu_{\rho}^{-1}(y)y)=0,
\]
puesto que $\mu_{\rho}^{-1}(y)$ es una unidad en $B^+$ y puesto que
$yA=\pK^e$.

Ahora $\rho$ y $\pK\star\rho$ est\'an definidos sobre $B^+$ por lo que ambos
valuaciones $v_{\eu q}(D(\pK\star\rho)_{\eu t})$ y $v_{\eu q}(D(\rho_{\pK}))$ son
no negativos y por ende $v_{\eu q}(D(\rho_{\pK}))=0$. $\fin$
\end{proof}

\subsection{El campo de clase de Hilbert $H_A$}\label{DrinfeldS3.3}

Regresamos a $\K_{\rho}$, el m\'inimo campo de definici\'on para un $A$--m\'odulo
de Drinfeld $\rho\in\Drin_A(\Ci)$, $\rho$ de rango uno. Suponemos que $\rho$
es signo--normalizado. Se tiene que $\K \Fi\subseteq \K_{\rho}\subseteq H_A^+$
(Observaci\'on \ref{DrinfeldO3.1.12'}).

Sea $\xi\in\Ci$ tal que $\rho^{\prime}=\xi\rho\xi^{-1}$ est\'a definido sobre
$\K_{\rho}$. Puesto que el grupo de automorfismos de $\rho$ es $\*\F$
(Proposici\'on \ref{DrinfeldP3.1.25}), el
m\'aximo com\'un denominador de
los elementos invariantes dados en la obtenci\'on de $\K_{\rho}$
es $q-1$. Esto es, recordemos que $\K_{\rho}=\K (I_j\mid j\in S)$, $I_i:=
c_i\big(\prod_{j\in S}c_j^{\alpha_j}\big)^{(1-q^i)/g}$, donde $g=\mcd\big\{
q^i-1\mid i\in S\big\}$.

De hecho tenemos que si $n,m\in{\ma N}$, $\mcd(q^n-1,q^m-1)=q^d-1$
con $d=\mcd(n,m)$ pues como $d\mid n$ y $d\mid m$, entonces $q^d-1
\mid q^n-1$ y $q^d-1\mid q^m-1$ por lo que $q^d-1\mid h=\mcd\big(
q^n-1,q^m-1\big)$.

Puesto que $\*{\ma F_{q^n}}$ y $\*{\ma F_{q^m}}$ son c\'iclicos, existe un
\'unico subgrupo $T\subseteq \*{\ma F_{q^n}}\cap \*{\ma F_{q^m}}=
\*{\ma F_{q^d}}$ de orden $h$, por lo que $h\mid q^d-1$ y $h=q^d-1$.
Por lo tanto 
\[
g=\mcd\big\{q^i-1\mid i\in S\big\}=q^{i_0}-1
\]
para alg\'un
$i_0$. Si $i_0>1$, se tendr\'ia que para toda $\xi\in\*{\ma F_{q^{i_0}}}$,
$\xi\rho\xi^{-1}=\rho$ pues $c_j(\xi\rho\xi^{-1},a)=\xi^{1-q^i}c_j(\rho,a)$
y si $j\in S$, $i_0\leq j$, por lo que $\xi^{1-q^j}=1$ y por tanto se
seguir\'ia que $\Aut(\rho)=\*{\ma F_{q^{i_0}}}$ lo cual es absurdo. Por tanto
$g=q-1$.

Ahora bien, $\xi^g=\xi^{q-1}=\prod_{i\in S}c_i^{\alpha_i}$ (ver la demostraci\'on
del Teorema \ref{DrinfeldT3.1.6}), por tanto $\xi_0:=\xi^{q-1}\in H_A^+$ pues
$H_A^+$ es el campo generado por los coeficientes de $\rho$. Puesto que
$\rho$ es signo normalizado, $\mu_{\rho}(\pi^{-1})=1$ y
\begin{gather*}
\mu_{\rho^{\prime}}(\pi^{-1})=\xi^{1-q^{\di}}\mu_{\rho}(\pi^{-1})=
\xi_0^{(1-q^{\di})/(q-1)}\mu_{\rho}(\pi^{-1})=\xi_0^{(1-q^{d_{\infty}})/q-1},
\intertext{lo cual implica que}
\xi_0^{(q^{\di}-1)/(q-1)}=\mu_{\rho'}(\pi^{-1})^{-1}\in \K_{\rho}.
\end{gather*}

Adem\'as, hab\'iamos visto que $H_A^+=\K_{\rho}(\xi_0)$ puesto que los
coeficientes de $\rho_x^{\prime}=\xi\rho_x\xi^{-1}$ generan $H_A^+$ para
cualquier $x\in A$ no constante (ver Teoremas \ref{DrinfeldT3.2.2}
y \ref{DrinfeldT3.2.10}). Por lo tanto $[H_A^+:\K_{\rho}]\leq
\frac{q^{\di}-1}{q-1}$.

Ahora consideremos la sucesi\'on exacta
\[
1\lra \frac{P_A}{P_A^+}\lra \pic^+(A)\stackrel{\theta}{\lra}\pic (A)\lra 1,
\]
donde $\theta$ es el mapeo natural. Identificamos $\pic^+(A)$ con
$\Gal(H_A^+/\K )$. Se tiene

\begin{proposicion}\label{DrinfeldP3.3.1}
El campo $\K_{\rho}$ es el subcampo de $H_A^+$ fijado por el subgrupo
de $\pic^+(A)$ generado por $\sigma_x$, $x\in\*\K $. De hecho
$\sigma_x=\sigma_{xA}$, $\sigma_{\eu A}=\eu A\star \rho$ es el
s\'imbolo de Artin.

Adem\'as, la extensi\'on $H_A^+=\K_{\rho}(\xi_0)/\K_{\rho}$ es una 
extensi\'on de Kummer c\'iclica de grado $(q^{\di}-1)/(q-1)$ y para 
cualquier $x\in \*\K $, tenemos 
\[
\xi_0^{\sigma_x}=\mu_{\rho}(x)^{q-1}\xi_0.
\]

En particular, $\K_{\rho}$ es independiente de la elecci\'on de $\rho$.
\end{proposicion}

\begin{proof} Puesto que $\sigma_x\rho=\mu_{\rho}(x)^{-1}\rho\mu_{\rho}(x)$,
se sigue de $c_i(\sigma_x\rho,a)=\mu_{\rho}^{q^i-1}c_i(\rho,a)$ que 
$\sigma_x$ fija a los invariantes $I_j$ y por lo tanto a $\K_{\rho}=
\K (I_j\mid j\in S)$.

Consideremos $\sigma_x$ como monomorfismo de $H_A^+(\xi)$ en
$\Ci$. Se tiene 
\begin{align*}
\rho^{\prime}&=\sigma_x\rho^{\prime}=\xi^{\sigma_x}\sigma_x\rho
\xi^{-\sigma_x}=\xi^{\sigma_x}\mu_{\rho}(x)^{-1}\rho\mu_{\rho}(x)
\xi^{-\sigma_x}\\
&= (\xi^{\sigma_x-1}\mu_{\rho}(x)^{-1})\rho^{\prime}(\xi^{\sigma_x-1}
\mu_{\rho}(x)^{-1})^{-1},
\end{align*}
donde $\rho_x=\xi\rho_x^{\prime}\xi^{-1}$ para toda $x\in A$.

Esto es, $\xi^{\sigma_x-1}\mu_{\rho}(x)^{-1}$ es un automorfismo de
$\rho^{\prime}$  por tanto un elemento de $\*\F$. Puesto que $\xi_0=
\xi^{q-1}$ se sigue que 
\begin{align*}
\xi_0^{\sigma_x}&= (\xi^{q-1})^{\sigma_x}=(\xi^{\sigma_x})^{q-1}=
\underbrace{(\xi^{\sigma_x-1}\mu_{\rho}(x)^{-1})^{q-1}}_{\substack{\uigual\\ 1}}
\xi^{q-1}\mu_{\rho}(x)^{(q-1)}\\
&=\xi_0\mu_{\rho}(x)^{q-1},
\end{align*}
esto es, $\xi_0^{\sigma_x}=\mu_{\rho}(x)^{q-1}\xi_0$.

Seleccionando $\mu_{\rho}(x)$ como un generador de
$\*\Fi$, esto es, $o(\mu_{\rho}(x)
\bmod (q^{\di}-1))=q^{\di}-1$, se sigue que $[\K_{\rho}(\xi_0):\K_{\rho}]
\geq \frac{q^{\di}-1}{q-1}$ y puesto que se ten\'ia la otra desigualdad,
obtenemos $[H_A^+:\K_{\rho}]=\frac{q^{\di}-1}{q-1}$.

Por otro lado, puesto que  $\frac{P_A}{P_A^+}
\cong\frac{\*\Fi}{\*\F}$ que es de
orden $\frac{q^{\di}-1}{q-1}$, se sigue que
$\K_{\rho}$ es el campo fijo del subgrupo
$\{\sigma_x\mid x\in\*\K \}$ pues $\sigma_xI_i=I_i$. Finalmente, la
extensi\'on es c\'iclica de Kummer pues $o(\sigma_x)=\frac{q^{\di}-1}{q-1}\mid
q^{\di}-1$ y $\*\Fi\subseteq \K_{\rho}$. $\fin$
\end{proof}

\begin{definicion}\label{DrinfeldD3.3.2} El campo com\'un de definici\'on de las 
$A$--m\'odulos de Drinfeld de rango $1$ sobre $\Ci$ se denota por
$H_A$ y se llama {\em campo de clase de Hilbert\index{campo de clase
de Hilbert} de $A$} (ver Definici\'on \ref{DrinfeldD3.1.17}).
\end{definicion}

\begin{observacion}\label{DrinfeldO3.3.2'}
Notemos que $H_A^+$ depende de la funci\'on $\sgn$ y $H_A$
depende \'unicamente de $\p$.
\end{observacion}

\begin{teorema}\label{DrinfeldT3.3.3}
Se tiene que $\p$ se descompone totalmente en $H_A/\K $ y todo divisor
primo $\pK$ de $\K $ es no ramificado en $H_A/\K $. La extensi\'on 
$H_A/\K $ es de 
grado $h_A$ con grupo de Galois isomorfo a $\pic(A)$ bajo el mapeo de
Artin. Si $\rho$ es un $A$--m\'odulo de Drinfeld de rango $1$ definido sobre
$H_A$, entonces 
\[
\sigma_{\eu A}\rho=\eu A\star \rho
\]
para todo ideal no cero $\eu A$ en $A$ y donde $\sigma_{\eu A}$ es el 
mapeo de Artin.

Finalmente, en $H_A^+/H_A$, $\p$ se ramifica totalmente.
\end{teorema}

\begin{proof} Puesto que $\Ki$ es un campo de definici\'on de $\rho$, $\K_{\rho}
=H_A\subseteq \Ki$. Ahora, $\Ki\cong \ma F_{q^{\di}}((\pi))=\Fi((\pi))$ y si
$\hat{\o}_{\infty}=\Fi[[\pi]]$, $\hat{\pK}_{\infty}=\pi\hat{\o}_{\infty}$, se tiene
$\hat{\o}_{\infty}/\hat{\pK}_{\infty}\cong \o_{\infty}/\p$ por lo que el grado relativo 
de $\p$ en $\Ki/\K $ es $1$ y $\p$ es no ramificado pues $\pi$ es un elemento
primo tanto para $\K $ como para $\Ki$. Se sigue que en $\Ki/\K $, $\p$ es
totalmente descompuesto. Como $H_A\subseteq \Ki$, se tiene que $\p$
es totalmente descompuesto en $H_A/\K $. Si $\pK$ es un primo finito,
$\pK$ es no ramificado en $H_A^+/\K $ por lo que es no ramificado en
$H_A/\K $. El mapeo de Artin fue obtenido previamente.

Se tiene $[H_A^+:H_A]=\frac{q^{\di}-1}{q-1}=\big|\frac{P_A}{P_A^+}\big|$ y
$P_A/P_A^+$ deja fijo a $H_A$. Por tanto
\[
\left.
\begin{array}{rcl}
\Gal(H_A^+/\K )&\lra&\Gal(H_A/\K )=G\\
\sigma&\lra& \sigma|_{H_A}\\
\pic^+(A)&\lra &G\\
P_A/P_A^+&\lra & \{1\}
\end{array}\right\}
\e\Lra\e 
\frac{\pic^+(A)}{(P_A/P^+_A)}\cong \pic(A)\subseteq G
\]
y $|\pic(A)|=|G|=h_A$ por lo que \fbox{$\Gal(H_A/\K )\cong \pic(A)$}.

Ahora, por lo que sabemos de teor\'ia de campos de clase, la m\'axima
extensi\'on abeliana no ramificada de $\K $ donde $\p$ se descompone
totalmente satisface $\Gal(L/\K )\cong Cl_A=\pic(A)$ (Teorema
\ref{CClaseT4.9.6}) y puesto que
$H_A\subseteq L$ y $\Gal(H_A/\K )\cong\pic(A)$ se sigue que $L=H_A$.
Se tiene que $H_A^+/H_A$ es totalmente ramificada en $\p$
(ver \cite{Hay92}, despu\'es del Theorem 15.6).
$\fin$
\end{proof}

\begin{corolario}\label{DrinfeldC3.3.3'}
Se tiene
\begin{gather*}
\Gal(H_A^+/H_A)\cong \*{\Fi}/\*{\F}. \tag*{$\fin$}
\end{gather*}
\end{corolario}

El siguiente resultado lo enunciamos sin demostrarlo.

\begin{teorema}\label{DrinfeldT3.3.4} Sea $B$ la cerradura entera de $A$ en
$H_A$. Entonces todo $A$--m\'odulo de Drinfeld $\rho$ de rango $1$
es isomorfo a un $A$--m\'odulo de Drinfeld $\rho^{\prime}$ el cual est\'a
definido en $B$ y donde $\mu_{\rho^{\prime}}(a)$ es una unidad en
$B$ para toda $a\in A\setminus\{0\}$.
\end{teorema}

\begin{proof} \cite[Theorem 15.8]{Hay92}. $\fin$
\end{proof}

Se tiene tambi\'en

\begin{teorema}[Ideal principal]\label{DrinfeldT3.3.5} Sea $\rho$ un $A$--m\'odulo
de Drinfeld el cual est\'a definido en $B$. Si $\eu A$ es cualquier ideal no 
cero en $A$, entonces $\eu A B=D(\rho_{\eu A}) B$ es el ideal principal
de $B$ generado por $D(\rho_{\eu A})$.
\end{teorema}

\begin{proof} Similar al Teorema \ref{DrinfeldT3.2.13}. $\fin$
\end{proof}

\subsection{Teor\'ia de campos de clase expl\'icita y campos de clase de 
rayos}\label{DrinfeldS3.4}

Construiremos la m\'axima extensi\'on abeliana de $\K $. Esta construcci\'on,
despu\'es de haber construido $H_A$ y $H_A^+$, es totalmente an\'aloga
al caso de los campos de funciones ciclot\'omicos. Fijamos una funci\'on
signo $\sgn$. Sea $\eu m$ un ideal no cero de $A$, $\eu m\neq A$. Sea
$\K_{\eu m}:=\K (\rho[\eu m])$\label{DrinfeldKm}. Como en el caso de campos de funciones
ciclot\'omicos $\K_{\eu m}/\K $ es una extensi\'on abeliana no ramificada 
fuera de $\eu m$ y $\p$, donde $\rho$ es un $A$--m\'odulo de Drinfeld de
rango $1$ sobre $\Ci$. 

Se tiene que $\K_{\eu m}$ es un campo de clase de rayos
extendidos de conductor $\eu m$. Como de costumbre, definimos $\K^+_{\eu m}$
como la m\'axima extensi\'on abeliana en $\K_{\eu m}$ en la cual $\p$
se descompone totalmente. Resulta ser que $\K^+_{\eu m}$ es el campo de
clase de rayos m\'odulo $\eu m$. De esta forma se obtiene una descripci\'on
expl\'icita de la m\'axima extensi\'on abeliana de $\K $ en la cual $\p$
se descompone totalmente. Se obtienen todos los campos de clase
variando $\p$. Las t\'ecnicas para estudiar $\K_{\eu m}$ son similares a los de
$H_A^+$.

Sea $D_{\eu m,A}$\label{Drinfeldidealesm}
 el grupo de ideales fraccionarios de $A$ generado por
los ideales primos no cero de $A$ tales que $\pK\nmid \eu m$ y sea
\[
P_{\eu m,A}^+:=\{xA\mid x\in\*\K , x\text{\ es positivo\ }, 
x\equiv 1\bmod \eu m\}.\label{Drinfeldprincipalesm}
\]

\begin{definicion}\label{DrinfeldD3.4.1}
El grupo cociente $\pic_{\eu m}^+(A)=
\frac{D_{\eu m,A}}{P_{\eu m,A}^+}$\label{Drinfeldpicm}
se llama el {\em grupo de clases de 
rayos extendido\index{grupo de clase de rayos
extendido} m\'odulo $\eu m$ relativo a $\sgn$}.
\end{definicion}

Recordemos que en el lenguaje de id\`eles
(ver Secci\'on \ref{CClaseC4} y la Subsecci\'on
\ref{CClaseS4.9}), dado un m\'odulus $\eu m$, 
$J_\K^{\eu m}=\{\alpha\in J_\K \mid \alpha\equiv 1\bmod \eu m\}=\prod_{
\pK\in\PP \K } U_{\pK}^{(n_{\pK})}$, $\eu m=\prod_{\pK\in \PP \K }\pK^{n_{\pK}}$ y
$C_\K^{\eu m}=\JKm\*\K /\*\K \subseteq C_\K $ y $C_\K /\CKm$ es el grupo de
clases de rayos m\'odulo $\eu m$.

Ahora consideremos el conjunto de m\'odulos de Hayes $\mathcal H$. Se
tiene $\rho[\eu m]\cong A/{\eu m}$ como $A$--m\'odulos 
puesto que $r_{\rho}=1$ (ver Corolario \ref{DrinfeldC1.3.30}).
Sea $\Phi(\eu m):=\big|\*{\big(A/{\eu m}\big)}\big|$. Se tiene que $\rho[\eu m]$
tiene $\Phi(\eu m)$ generadores como $A$--m\'odulos.

Sea $X_{\eu m}:=\{(\rho,\lambda)\mid \rho\in \mathcal H, \lambda \text{\ generador
de\ }\rho[\eu m]\}$. Se define la acci\'on de $D_{\eu m, A}$ sobre $X_{\eu m}$
por 
\[
\eu A\star(\rho,\lambda)=(\eu A\star \rho,\rho_{\eu A}(\lambda)).
\]

Veamos que $\rho_{\eu A}(\lambda)$
es generador de $({\eu A}\star \rho)[\eu m]$. Sea $\rho'={\eu A}\star
\rho$. Para $\mu\in\rho[{\eu m}]$ y $m\in{\eu m}$, tenemos $\rho_m
(\mu)=0$ por lo que $\rho'_m(\rho_{\eu A}(\mu))=\rho_{\eu A}(\rho_m(u))
=\rho_{\eu A}(0)=0$.

Por tanto $\varphi\colon \rho[{\eu m}]\lra \rho'[{\eu m}]$ dada por
$\varphi(\mu)=\rho_{\eu A}(\mu)$ est\'a bien definida y $\varphi$
es un homomorfismo de $A$--m\'odulos. Veamos que $\varphi$
es inyectiva. Si $\varphi(\lambda)=\varphi(\lambda')$, se tiene
$\rho_{\eu A}(\lambda)=\rho_{\eu A}(\lambda')$ lo cual implica que
$\rho_{\eu A}(\lambda-\lambda')=0$.

Por otro lado, $\rho_m(\lambda-\lambda')=0$ para toda $m\in{\eu m}$.
Puesto que $\mcd({\eu A},{\eu m})=1$, se tiene que ${\eu A}+{\eu m}
=(1)$ y existen $a\in {\eu A}$ y $m\in{\eu m}$ tales que $a+m=1$.
Ahora $\rho_a\in R\rho_{\eu A}$, es decir, existe $\xi\in R=\aditivo{\Ki}$
tal que $\rho_a=\xi\rho_{\eu A}$. Se sigue que 
\begin{gather*}
\rho_a(\lambda-\lambda')=\xi(\rho_{\eu A}(\lambda-\lambda'))=\xi(0)=0
\quad\text{y}\quad \rho_m(\lambda-\lambda')=0.
\intertext{Por tanto}
\lambda-\lambda'=\rho_1(\lambda-\lambda')=
\rho_{a+m}(\lambda-\lambda')=\rho_a(\lambda-\lambda')
+\rho_m(\lambda-\lambda')=0+0=0,
\end{gather*}
por tanto $\lambda=\lambda'$ y $\varphi$ es inyectiva.

Puesto que $\rho[{\eu m}]\cong\rho'[{\eu m}]\cong A/{\eu m}$, tenemos que
$|\rho[{\eu m}]|=|\rho'[{\eu m}]|$ de donde se sigue que $\varphi$ es
suprayectiva y $\lambda$ es generador de $\rho[{\eu m}]\iff \rho_{\eu A}
(\lambda)$ es generador de $\rho'[{\eu m}]$.

Probemos que el estabilizador de cualquier 
punto $(\rho,\lambda)$ es $P_{\eu m,A}^+$.
Tenemos que $\rho_{\eu A}(\lambda)$ es generador de
$({\eu A}\star\rho)[{\eu m}]$
y si $\eu A\star(\rho,\lambda)=(\rho,\lambda)$
entonces $\eu A\star \rho=\rho$ y $\rho_{\eu A}(\lambda)=\lambda$
y como $\rho$ es de Hayes para que
${\eu A}\star {\rho}={\rho}$, ${\eu A}$ debe ser principal
(ver Teorema \ref{DrinfeldT3.1.32} y posterior).

Sea $\eu A=(x)$, y $(\eu A\star\rho)_a=\mu_{\rho}(x)^{-1}\rho_a\mu_{\rho}
(x)$. Adem\'as $\rho_{\eu A}(\lambda)=\rho_{xA}(\lambda)=
\mu_{\rho}(x)^{-1}\rho_x(\lambda)=\lambda$ lo cual implica 
$\rho_x(\lambda)=\mu_{\rho}(x)\lambda$. Por tanto obtenemos
que $\rho_{x-\mu_{\rho}(x)}(\lambda)=0$, esto es,
$x-\mu_{\rho}(x)\in{\eu m}$. As\'i, $\mu_{\rho}(x)^{-1} x-1\in{\eu m}$
y $\mu_{\rho}(x)^{-1}x\equiv 1\bmod {\eu m}$.

Como $\eu A\star\rho=\rho$, $\mu_{\rho}(x)\in \Aut (\rho)=\*\F$ y el
estabilizador de $\rho$ es 
\[
\{xA\mid x\in A, \mu_{\rho}(x)\in\*\F\}
\underbracket[0pt]{=}_{\substack{\uparrow\\ xA=(\mu_{\rho}(x)^{-1}
x) A}}\{xA\mid x\in A, \sgn(x)=1\}
\]
y adem\'as $\mu_{\rho}(x)^{-1}x\equiv 1\bmod {\eu m}$. Por lo
tanto $\eu A\star (\rho,\lambda)=(\rho,\lambda)\iff \eu A=xA, \sgn x=1$,
$x \equiv 1\bmod \eu m$, por lo que el estabilizador es $P_{\eu m,A}^+$.

Por un lado tenemos que $\rho[{\eu m}]\cong (A/{\eu m})$ por lo que
el n\'umero de generadores de $\rho[{\eu m}]$ es $\Phi({\eu m})=
\big|((A/{\eu m})^*)\big|$ y en particular
$|X_{\eu m}|=|\mathcal H|\Phi(\eu m)=|\pic^+(A)|\Phi(\eu m)$.

Por otro lado tenemos la sucesi\'on exacta
\[
1\lra \frac{D_{\eu m,A}\cap P_A^+}{P_{\eu m,A}^+}\lra \pic_{\eu m}^+(A)
\xrightarrow[]{\text{natural}} \pic^+(A)\lra 1.
\]

Ahora si $\rho_{\eu A}(\lambda)$ es generador de $\rho[{\eu m}]$
en particular ${\eu A}$ debe ser estabilizador de $\rho$, esto es
${\eu A}\in P_A^+$ (Teorema \ref{DrinfeldD3.1.34}) y por otro lado
${\eu A}$ es primo relativo a ${\eu m}$. De esto se deduce que
\begin{gather*}
\Big|\frac{D_{\eu m,A}\cap P_A^+}{P_{\eu m,A}^+}\Big| =
\Phi({\eu m}).
\intertext{Por tanto}
|\pic_{\eu m}^+(A)|=|\pic^+(A)|\Big|\frac{D_{\eu m,A}\cap P_A^+}{P_{\eu m,A}^+}
\Big|=|\pic^+(A)|\Phi(\eu m)=|X_{\eu m}|.
\end{gather*}

Obtenemos que $\pic_{\eu m}^+(A)$ act\'ua 
transitivamente sobre $X_{\eu m}$ y se sigue:

\begin{teorema}\label{DrinfeldT3.4.2}
El conjunto $X_{\eu m}$ es un espacio homog\'eneo para $\pic_{\eu m,A}^+(A)$
con la acci\'on $\star$, es decir, $X_{\eu m}=\pic_{\eu m}^+(A)\star(
\rho_0,\lambda_0)$. $\fin$
\end{teorema}

A continuaci\'on copiamos toda la teor\'ia de campos de funciones ciclot\'omicos
pero en lugar de tomar como campo base a $K=\F(T)$ o $\K $, tomamos 
$H_A^+$.

\begin{definicion}\label{DrinfeldD3.4.3} 
Sea $\K (\eu m):=H_A^+(\rho[\eu m])$. 
\end{definicion}

Como en el caso $H_A^+$, se probar\'a que $\K (\eu m)/\K $ es una extensi\'on
de Galois y no ramificada fuera de $\p$ y de los primos que dividen a $\eu m$.

\begin{proposicion}\label{DrinfeldP3.4.4}
Sea $L/K$ una extensi\'on finita y sea $\rho\in\Drin_A(\Ci)$ de rango $1$ el cual
est\'a definido sobre un anillo de valuaci\'on $\o_{\eu q}$ en $L$ donde
$\eu q$ es no ramificado en $L/\K $. Sea $\pK=\eu q\cap A$. Sea $\eu A=
\pK^e$, $\eu B=\pK^{e-1}$. Entonces $\rho_{\eu B}(t)$ divide a $\rho_{\eu A}
(t)$ en $\o_{\eu q}[t]$ y el cociente es Eisenstein en $\pK$.
\end{proposicion}

\begin{proof} Para $e=1$ la demostraci\'on es como en el caso para $H_A^+$ con igual
argumento (ver Teorema \ref{DrinfeldT3.2.12}).

Para $e>1$, sea $f(t)=(\eu B\star \rho)_{\pK}(t)/t$ y se tiene $\rho_{\eu A}(t)=
f(\rho_{\eu B}(t))\rho_{\pK}(t)$. Para el caso $e=1$, $f(t)$ es Eisenstein en
$\pK$ y $\rho_{\pK}(t)\equiv t^{\N(\pK)}\bmod \eu q$. $\fin$
\end{proof}

Por la Proposici\'on \ref{DrinfeldP3.4.4}, se sigue que $\K (\eu m)/\K $ tiene el mismo
tipo de ramificaci\'on como en el caso ciclot\'omico, es decir

\begin{proposicion}\label{DrinfeldP3.4.5} Sea $\eu m=\pK^e$ donde $\pK$ es un ideal
primo no cero de $A$. Entonces $\K (\pK^e)=H_A^+(\rho[\pK^e])/H_A^+$ es
totalmente ramificado en $\eu q$ donde $\eu q$ es un divisor primo de 
$H_A^+$ sobre $\pK$ y el \'indice de ramificaci\'on es $\Phi(\pK^e)$. M\'as
a\'un, la extensi\'on $\K (\pK^e)/H_A^+$ es no ramificada en todo primo finito
$\pK_1\neq \pK$. Finalmente, se tiene $[\K (\pK^e):H_A^+]=\Phi(\pK^e)$.
\end{proposicion}

\begin{proof} El polinomio $f(u)=\rho_{\pK^e}(u)/\rho_{\pK^{e-1}}(u)$ es Eisenstein
en $\pK$ y $f(u)=\prod (u-\lambda)$ donde $\lambda$ recorre los 
generadores del $A$--m\'odulo c\'iclico $\rho[\pK^e]\cong A/\pK^e$. Se 
tiene $\deg f(u)=\Phi(\pK^e)$. El resto de la demostraci\'on es como en el
caso ciclot\'omico (ver Teorema \ref{T6.2.28}). $\fin$
\end{proof}

\begin{corolario}\label{DrinfeldC3.4.6} Para cualquier ideal no cero $\eu m$ de $A$,
$\K (\eu m)/H_A^+$ es una extensi\'on de Galois con grupo
de Galois isomorfo a $\*{\big(
A/\eu m\big)}$. Los primos finitos ramificados son los ideales primos $\pK$
que dividen a $\eu m$ con \'indice de ramificaci\'on $\Phi(\pK^e)$ donde
$\pK^e$ es la potencia exacta de $\pK$ que dividen a $\eu m$.
\end{corolario}

\begin{proof} Igual que en el caso ciclot\'omico (ver Teorema \ref{T6.2.28}). $\fin$
\end{proof}

Ahora recordemos que $\K_{\eu m}=\K (\rho[\eu m])$
con $\rho[{\eu m}]=\{u\in \bar{\K }_{\infty}\mid \rho_m(u)=0
\text{\ para toda } m\in{\eu m}\}$. Entonces $\K_{\eu m}/\K $
es una extensi\'on normal debido a que $\sigma(\K_{\eu m})=\K (\sigma\rho[
\eu m])=\K_{\eu m}$ para todo $\sigma\in \Aut_\K (\Ci)$.

Puesto que $H_A^+$ est\'a generado sobre $\K $
 por los coeficientes de $\rho_a$ con
$\rho\in\mathcal H$ y $a$ un elemento no constante de $A$, se sigue que
$\K (\rho[{\eu m}])=\K_{\eu m}=\K ({\eu m})=H_A^+(\rho[{\eu m}])$.

De hecho se tiene:

\begin{teorema}\label{DrinfeldT3.4.6'} Si $\lambda$ es un generador de $\rho[{\eu m}]$,
entonces $H_A^+(\rho[{\eu m}])=\K (\lambda)$. En particular $\K ({\eu m})
=\K_{\eu m}$.
\end{teorema}

\begin{proof} Ver \cite[Theorem 3.6.7, p\'agina 77]{Tha2004} o
\cite[Theorem 7.5.15, p\'agina 207]{Gos96}. $\fin$
\end{proof}

\begin{teorema}\label{DrinfeldT3.4.7} Sea $\eu A$ un ideal no cero de $A$ que es primo
relativo a $\eu m$ y sea $\lambda\in\rho[\eu m]$. Entonces, si $\sigma_{\eu A}$
es el automorfismo de Artin, se tiene 
\[
\lambda^{\sigma_{\eu A}}:=
\sigma_{\eu A}\lambda=\rho_{\eu A}(\lambda).
\]
\end{teorema}

\begin{proof} Si $\eu A$ y $\eu B$ son dos ideales no cero de $A$, ambos primos relativos
a $\eu m$ que satisfacen $\lambda^{\sigma_{\eu A}}=\rho_{\eu A}(\lambda)$
y $\lambda^{\sigma_{\eu B}}=\rho_{\eu B}(\lambda)$, entonces 
\begin{align*}
\sigma_{\eu{AB}}&=\sigma_{\eu A}\sigma_{\eu B}(\lambda)=\sigma_{\eu B}
\sigma_{\eu A}(\lambda)=\sigma_{\eu B}(\rho_{\eu A}(\lambda))=
(\sigma_{\eu B}\rho_{\eu A})(\sigma_{\eu B}(\lambda))\\
&=(\sigma_{\eu B}\rho_{\eu A})(\rho_{\eu B}(\lambda))
\underbracket[0pt]{\igual}_{\substack{\uparrow\\ \sigma_{\eu B}\rho=
({\eu B}\star \rho)\\ \text{Teorema \ref{DrinfeldT3.3.3}}}}
(\eu B\star\rho)_{\eu A}(\rho_{\eu B}(\lambda))=\rho_{\eu {AB}}(\lambda).
\end{align*}

Por tanto podemos suponer que $\eu A=\pK$ es primo con $\pK\nmid {\eu m}$. 
Se $\eu q$ en $\K_{\eu m}$
es un primo sobre $\pK$, $\sigma_{\pK}$ satisface $\sigma_{\pK}\lambda
=\lambda^{\N \pK}\bmod \eu q$ (ver Teorema \ref{DrinfeldT3.2.10}).

Puesto que $\rho\bmod \eu q=\phi$ satisface $\phi_{\pK}=\tau^{\deg\pK}$,
se sigue como en la demostraci\'on en el caso $H_A^+/\K $ que
$\phi_{\pK}\lambda=\lambda^{\N\pK}\bmod \eu q$ y por tanto $\sigma_{\pK}
(\lambda)=\rho_{\pK}(\lambda)$. $\fin$
\end{proof}

Como consecuencia de este resultado, vemos que $\K_{\eu m}$ es independiente
de $\rho\in {\mc H}$, pues si $\lambda$ es generador de $\rho[{\eu m}]$,
$\lambda^{\sigma_{\eu A}}=\rho_{\eu A}(\lambda)$ es generador
de $({\eu A}\star\rho)[{\eu m}]$ con ${\eu A}\star\rho\in{\mc H}$. 
Adem\'as $\pic_{\eu m}^+(A)$
act\'ua como automorfismos v\'ia
\[
\sigma_{\eu A}(\lambda)=\rho_{\eu A}(\lambda).
\]

Puesto que $|\pic_{\eu m}^+(A)|=|\pic^+(A)|\Phi(\eu m)=[H_A^+:\K ]
[H_A^+(\rho[\eu m]):H_A^+]=[H_A^+(\rho[\eu m]):\K ]=[\K_{\eu m}:\K ]$,
se sigue que 
\begin{gather*}
\pic_{\eu m}^+(A)\cong \Gal(\K_{\eu m}/\K ).\\
\xymatrix{
\K (\eu m)=\K_{\eu m}\ar@{-}[d]\ar@{-}@/^2pc/[d]^{
\frac{D_{\eu m,A}\cap P_A^+}{P_{\eu m,A}^+}\cong \*{(A/\eu m)}}
\ar@{-}@/_4pc/[dd]_{\pic^+_{\eu m}(A)}\\
H_A^+\ar@{-}[d]\ar@{-}@/^2pc/[d]^{\pic^+(A)}\\\K }
\end{gather*}

Ahora bien, los elementos positivos de $A$ generan $A/{\eu m}$ as\'i que
el mapeo $Aa\longmapsto \sigma_a:=\sigma_{aA}$ donde $a\in A$ es un
elemento de $A$ primo relativo a $\eu m$, induce un isomorfismo entre
$\*{\big(A/\eu m\big)}$ y $\Gal(\K_{\eu m}/H_A^+)$
(recordemos que $|\Gal(\K_{\eu m}/H_A^+)|=|(A/{\eu m})^*|=
\Phi({\eu m})$). Para un elemento $\lambda
\in\rho[\eu m]$ y si $x\in A$ es tal que $x\equiv 1\bmod \eu m$ en $\*\K $,
usando el s\'imbolo de Artin se obtiene
\[
\sigma_x(\lambda)=\rho_{xA}(\lambda)=\mu_{\rho}(x)^{-1}\lambda
\]
y $\mu_{\rho}(x)\in\*\Fi$. Por tanto $\Gal(\K_{\eu m}/\K )$ contiene un
subgrupo $I_{\p}$ isomorfo a $\*{\Fi}$. 
M\'as precisamente
\[
I_{\p}\cong\{xA\bmod P_{{\eu m},A}^+\mid x\equiv 1\bmod {\eu m}\}
\cong \*{\Fi}.
\]
De hecho Hayes prob\'o que este grupo es tanto
el grupo de descomposici\'on como el grupo de inercia
de $\p$ en $\K_{\eu m}/\K $ (ver
\cite[Proposition 4.15, p\'agina 229]{Hay85}).

\begin{definicion}\label{DrinfeldD3.4.8} El campo fijo de $\K_{\eu m}$ bajo
$I_{\p}$, $\K_{\eu m}^+:=\K_{\eu m}^{I_{\p}}$ se llama el {\em campo de clase
de rayos\index{campo de clase de rayos de conductor $\eu m$} de conductor
$\eu m$}. Estos campos corresponden a los campos de clase a
$\JKSm=\Big(\prod_{\pK\in S}\*{\K_{\pK}}\times\prod_{\pK\in \PP \K \setminus S}
U_{\pK}^{(m_{\pK})}\Big)$ donde $S=\{\p\}$.
\end{definicion}

Por tanto $\K_{\eu m}^+$ es el campo de clase de rayos de $\K $ de conductor
$\eu m$ en donde $\p$ se descompone totalmente. Estos 
campos, $\K_{\eu m}$ y $\K_{\eu m}^+$ son totalmente
similares a los campos ciclot\'omicos, donde $\K_{\eu m}$ juega
el papel del campo ciclot\'omico usual $\ma Q(\zeta_m)$ y $\K_{\eu m}^+$
el subcampo real de $\ma Q(\zeta_m)$, $\ma Q(\zeta_m)^+$ y lo mismo en
campos de funciones ciclot\'omicos $\K_{\eu m}\longleftrightarrow K(\Lambda_M)$
y $\K_{\eu m}^+\longleftrightarrow K(\lambda_M)^+$.

Notemos que el
\'indice de ramificaci\'on en $\K_{\eu m}/\K_{\eu m}^+$ es $q^{\di}-1
=|\*{\Fi}|$.
Se tiene $e_{\infty}(H_A^+|H_A)=\frac{q^{\di}-1}{q-1}$, 
$e_{\infty}(\K_{\eu m}/\K_{\eu m}^+)=q^{\di}-1$.
\begin{gather*}
\xymatrix{
&\K_{\eu m}\ar@{-}[dd]\ar@{-}[ld]^{
\substack{I_{\p}\\ \cong \*{{\ma F}_{\infty}}}}\ar@{-}[dr]^{(A/{\eu m})^{\ast}}
\ar@{-}@/^-6pc/[ddd]_{\pic_{\eu m}^+(A)}\\
\K_{\eu m}^+\ar@{-}[dr]&&H_A^+\ar@{-}[dl]^{\substack{\*{\Fi}/\*{\F}\cong\\ P_A/P_A^+}}
\ar@{-}@/^4pc/[ddl]^{\pic^+(A)}\\
&H_A\ar@{-}[d]_{\pic(A)}\\& \K }
\end{gather*}

Sea ahora $\Ki^+:=\cup_{\eu m}\K_{\eu m}^+$ donde $\eu m$ recorre el conjunto
de todos los ideales propios no cero de $A$. De esta forma $\Ki^+$ es la 
m\'axima extensi\'on abeliana de $\K $ tal que $\p$ se descompone totalmente.

Si $\Ki=\cup_{\eu m}\K_{\eu m}$, entonces $\Ki^+$ es el campo fijo de 
$\Ki$ bajo $I_{\p}$: $\Ki^+=\Ki^{I_{\p}}$. Notemos que ahora $\Ki$ denota algo
diferente a la completaci\'on de $\K $ en $\p$. Por otro lado $\Ki^+$ est\'a
contenido en la completaci\'on de $\K $ en $\p$ puesto que $\p$ se
descompone totalmente en $\Ki^+$.

De manera totalmente an\'aloga al caso de campos de funciones ciclot\'omicas,
se tiene que $\Gal(\Ki/\K )\cong V_{\p}$ donde $V_{\p}=\{1\}\times\prod_{\pK
\neq \p}U_{\pK}\subseteq J_\K $, es decir, $V_{\p}$ es el subgrupo  de los 
id\`eles cuya componente en $\p$ es $1$ y en las dem\'as son elementos
de $U_{\pK}=\*\o_{\pK}$. Esto es, $\Gal(\Ki/\K )\cong J_\K /(\*\K \times \Ki^{(1)}\times
\ma Z)$ y $\Ki$ corresponde a $\*\K \times \Ki^{(1)}\times \ma Z$ y $\Ki^{(1)}=
\ker (\sgn)\cap U_{\p}$.

Ahora tomamos $\p^{\prime}\neq \p$, $\p^{\prime}$ un divisor primo y
consideramos $\K_{\infty^{\prime}}$. Entonces
la intersecci\'on de los dos subgrupos de 
id\`eles correspondientes es $\*\K $ lo cual implica que $\Ki \K_{\infty^{\prime}}$ es la
m\'axima extensi\'on abeliana de $\K $ (ver Teorema \ref{CClaseTCCG}).

\begin{teorema}\label{DrinfeldT3.4.9} Sean $\p,\p^{\prime}$ dos divisores primos 
distintos de $\K $. Si $\Ki$ y $\K_{\infty^{\prime}}$ son como antes, entonces $\Ki
\K_{\infty^{\prime}}$ es la m\'axima extensi\'on abeliana de $\K $. $\fin$
\end{teorema}

%% file: Capitulo16.tex
\chapter{Teor{\'\i}a de Iwasawa\index{Iwasawa!teor{\'\i}a de
$\sim$}\index{teor{\'\i}a de Iwasawa}}\label{Ch7}

Este cap{\'\i}tulo presenta el inicio de la Teor{\'\i}a de Iwasawa
que fue desarrollada por Kenkichi Iwasawa, especialmente en
\cite{Iwa59, Iwa73, Iwa81}. El desarrollo que aqu{\'\i} presentamos
est\'a basado fuertemente en \cite[Cap{\'\i}tulo 13]{Was97}.

\section{Campos ciclot\'omicos infinitos}\label{S7.2}

Sea $p$ un n\'umero primo en ${\ma Z}$ y sea $\cic p{\infty}:
=\bigcup\limits_{n=0}^{\infty}\cic pn$. Sea $G:=\Gal(\cic p{\infty}/
{\ma Q})$. Se tiene que $\sigma\in G$ est\'a determinado por su
acci\'on en $\zeta_{p^n}$, $n\geq 1$. Sea $\sigma\zeta_{p^n}=
\zeta_{p^n}^{a_n}$ con $a_n\bmod p^n\in \big({\ma Z}/p^n{\ma Z}
\big)^{\ast}=U_{p^n}$. Se tiene que $a_n\equiv a_{n-1}\bmod p^{n
-1}$ por lo que obtenemos un elemento $\big(a_n\big)_{n\in{\ma N}}
\in{\ma Z}_p^{\ast}=\lim\limits_{\longleftarrow}
\big({\ma Z}/p^n{\ma Z}\big)^{
\ast}=\Gal\limits_{\longleftarrow}(\cic pn/{\ma Q})$.
Rec{\'\i}procamente, si $a\in{\ma Z}_p^{\ast}$, $\sigma \zeta_{p^n}=
\zeta_{p^n}^{a}$ da un automorfismo de $\cic p{\infty}$. 

Si $p>2$, tenemos ${\ma Z}_p^{\ast}=\lim\limits_{\longleftarrow}
\big({\ma Z}/p^n{\ma Z}\big)^{\ast}=\lim\limits_{\longleftarrow}
\big({\ma Z}/(p-1){\ma Z}\times {\ma Z}/p^{n-1}{\ma Z}\big)\cong
{\ma Z}/(p-1){\ma Z}\times {\ma Z}_p$.
Si $p=2$, ${\ma Z}_2^{\ast}=\lim\limits_{\longleftarrow}
\big({\ma Z}/2^n{\ma Z}\big)^{\ast}=\lim\limits_{\longleftarrow}
\big({\ma Z}/2{\ma Z}\times {\ma Z}/2^{n-2}{\ma Z}\big)\cong
{\ma Z}/2{\ma Z}\times {\ma Z}_2$.
\[
\xymatrix{
\cic p{}\ar@{-}[r]^{{\ma Z}_p}\ar@{-}[d]_{U_p\cong}^{{\ma Z}/(p-1)
{\ma Z}}&\cic p{\infty}\ar@{-}[d]^{C_{p-1}}\\
{\ma Q}\ar@{-}[r]_{{\ma Z}_p}&{\ma Q}_{\infty}
}
\xymatrix{
{\ }\ar@{}[d]_{p>2}\\ {\ }
}
\qquad
\xymatrix{
{\ }\ar@{-}[d]\\ {\ }
}
\qquad\qquad
\xymatrix{
\cic 4{}\ar@{-}[r]^{{\ma Z}_2}\ar@{-}[d]_{C_2}
&\cic 2{\infty}\ar@{-}[d]^{C_{2}}\\
{\ma Q}\ar@{-}[r]_{{\ma Z}_2}&{\ma Q}_{\infty}
}
\]

Estamos particularmente interesados en $\cic p{\infty}$
y en ${\ma Q}_{\infty}=\cic p{\infty}^{C_{p-1}}$, $p\geq 2$.

\section{Ramificaci\'on en extensiones algebraicas}\label{S7.3}

Aqu{\'\i} presentamos una teor{\'\i}a de ramificaci\'on para extensiones
algebraicas arbitrarias.

Sea $k/{\ma Q}$ una extensi\'on algebraica no necesariamente 
finita. Sea ${\cal O}_k:=\{\alpha\in k\mid \Irr(\alpha,x,{\ma Q})\in
{\ma Z}[x]\}$. ${\cal O}_k$ es el {\em anillo de enteros\index{anillo
de enteros}} de $k$. Se tiene
\[
{\cal O}_k=\lim_{\substack{\longto\\ [E:{\ma Q}]<\infty\\ E\subseteq k}}
{\cal O}_E=\bigcup_{\substack{[E:{\ma Q}]<\infty\\ E\subseteq k}}
{\cal O}_E.
\]

\begin{observacion}\label{O7.3.1} En general ${\cal O}_k$ no es
noetheriano y por lo tanto no es dominio Dedekind. M\'as adelante
daremos un ejemplo.
\end{observacion}

Lo que si tenemos es:

\begin{proposicion}\label{P7.3.2} Si $\pK$ es un ideal primo no cero
de ${\cal O}_k$, entonces $\pK$ es maximal.
\end{proposicion}

\begin{proof}
Se tiene que ${\cal O}_k/\pK$ es un dominio entero y 
$\pK\cap {\ma Z}$ es un ideal primo de ${\ma Z}$.
Sea $\alpha \in \pK$,
$\alpha\neq 0$. Entonces $[{\ma Q}(\alpha):{\ma Q}]<\infty$ y
$0\neq N_{{\ma Q}(\alpha)/{\ma Q}}\alpha\in \pK\cap {\ma Z}$,
es decir, $\pK\cap {\ma Z}\neq \{0\}$. Por lo tanto $\pK\cap {\ma Z}=
\langle p\rangle$ con $p$ un n\'umero primo. Se tiene
\[
{\ma Z}/p{\ma Z}={\ma Z}/(\pK\cap {\ma Z})\cong (\pK+{\ma Z})/\pK
\subseteq {\cal O}_k/\pK.
\]

Por otro lado si $\overline{\beta}\in {\cal O}_k/\pK$, $\beta\in
{\cal O}_k$, $\beta$ entero sobre ${\ma Z}$, por lo tanto $\overline{
\beta}$ es algebraico sobre ${\ma F}_p$ lo cual implica que
${\cal O}_k/\pK\subseteq \overline{{\ma F}_p}$. Es decir, ${\ma F}_p
\cong {\ma Z}/p{\ma Z}\subseteq {\cal O}_k/\pK\subseteq 
\overline{{\ma F}_p}=\bigcup\limits_{n=1}^{\infty}{\ma F}_{p^n}$.

Ahora veamos que ${\cal O}_k/\pK$ es un campo. Si
$\overline{\alpha}\in {\cal O}_k/\pK$, $\overline{\alpha}\neq 0$,
$\overline{\alpha}\in\overline{{\ma F}_p}$ y $\overline{\alpha}$
satisface una ecuaci\'on
\begin{gather*}
\overline{\alpha}^n+a_{n-1}\overline{\alpha}^{n-1}+\cdots+
a_1\overline{\alpha}+a_0=0\\
\intertext{con $a_i\in {\ma F}_p$ y $a_0\neq 0$. Por tanto}
-a_0^{-1}\big(\overline{\alpha}^{n-1}+\cdots+a_2\overline{\alpha}
+a_1\big)\overline{\alpha}=1
\end{gather*}
es decir $\overline{\alpha}^{-1}=-a_0^{-1}\big(\overline{\alpha}^{n-1}
+\cdots+a_2\overline{\alpha}+a_1\big)\in{\cal O}_k/\pK$ y 
$\overline{\alpha}$ es invertible. Por lo tanto ${\cal O}_K/\pK$ es
un campo y $\pK$ es maximal. $\fin$
\end{proof}

\begin{corolario}\label{C7.3.3} ${\cal O}_k/\pK$ es una extensi\'on
abeliana de ${\ma F}_p$.
\end{corolario}

\begin{proof}
Se sigue de que ${\ma F}_p\subseteq {\cal O}_k/\pK\subseteq 
\overline{{\ma F}_p}$ y de que 
\[
\Gal(\overline{{\ma F}_p}/{\ma F}_p)=\Gal\big(\lim_{\longto}
{\ma F}_{p^n}/{\ma F}_p\big)\cong \lim_{\longleftarrow}\Gal(
{\ma F}_{p^n}/{\ma F}_p)\cong \lim_{\longleftarrow}\big({\ma Z}/
n{\ma Z}\big)\cong \hat{{\ma Z}}
\]
el cual es abeliano. $\fin$
\end{proof}

M\'as generalmente, tenemos

\begin{corolario}\label{C7.3.4}
Si $L/K$ es una extensi\'on algebraica cualquiera
 de campos num\'ericos,
entonces si $\pL$ es un ideal primo no cero de ${\cal O}_L$, $\pL
\cap {\cal O}_K=\pK$ es un ideal primo no cero de $K$ y ${\cal O}_L
/\pL$ es una extensi\'on de Galois abeliana de ${\cal O}_K/\pK$.
\end{corolario}

\begin{proof}
Se sigue de que ${\ma F}_p\subseteq {\cal O}_k/\pK\subseteq
{\cal O}_L/\pL\subseteq \overline{{\ma F}_p}$. $\fin$
\end{proof}

Rec{\'\i}procamente, tenemos:

\begin{proposicion}\label{P7.3.5} Sea $L/K$ una extensi\'on 
algebraica de campos. Si $\pK$ es un ideal no cero de ${\cal O}_K$,
existe un ideal primo $\pL$ de ${\cal O}_L$ tal que $\pL\cap{\cal O}_K
=\pK$.
\end{proposicion}

\begin{proof}
Primero veamos que $\pK{\cal O}_L\neq {\cal O}_L$. Se tiene
que la localizaci\'on $\big({\cal O}_{L}\big)_{\pK}=\big\{
\frac{a}{b}\mid a\in{\cal O}_L, b\in {\cal O}_K\setminus \pK\big\}$ es
un anillo entero sobre $\big({\cal O}_K\big)_{\pK}$ y $\big({\cal O}_K
\big)_{\pK}$ es un anillo local con ideal m\'aximo $\pK
\big({\cal O}_K\big)_{\pK}$. Si $\pK{\cal O}_L={\cal O}_L$ se 
tendr{\'\i}a que $\pK\big({\cal O}_L\big)_{\pK}=
\big({\cal O}_L\big)_{\pK}$ y en este caso tendr{\'\i}amos
\[
1=a_1b_1+\cdots+a_nb_n\quad \text{con}\quad a_i\in\pK,
b_i\in \big({\cal O}_L\big)_{\pK}.
\]

Sea $B=\big({\cal O}_K\big)_{\pK}[b_1,\ldots,b_n]$. Se tiene 
$\pK B=B$ y $B$ es un $A:=\big({\cal O}_K\big)_{\pK}$ m\'odulo
finitamente generado pues cada $b_i$ es entero sobre $A$. Por
el Lema de Nakayama\index{Nakayama!lema de $\sim$}
(Teorema \ref{T7.3.6} m\'as adelante), se sigue que $B=0$ lo
cual es absurdo. Por tanto $\pK{\cal O}_L\neq {\cal O}_L$.

Se tiene el diagrama
\[
\xymatrix{
{\cal O}_L\ar[r]&\big({\cal O}_L\big)_{\pK}\\
{\cal O}_K\ar[r]\ar[u]&\big({\cal O}_K\big)_{\pK}\ar[u]
}
\]
$\pK\big({\cal O}_L\big)_{\pK}\subseteq {\eu m}\big(
{\cal O}_L\big)_{\pK}$ con ${\eu m}$ es un ideal maximal. En 
particular ${\eu m}\cap \big({\cal O}_K\big)_{\pK}=\pK
\big({\cal O}_K\big)_{\pK}$. Por tanto ${\eu m}':={\eu m}\cap
{\cal O}_L$ es un ideal primo de ${\cal O}_L$ y se tiene
\[
{\eu m}'\cap {\cal O}_K={\eu m}\big({\cal O}_L\big)_{\pK}\cap
\big({\cal O}_K\big)_{\pK}\cap {\cal O}_K=\pK\big({\cal O}_L\big)_{\pK}
\cap {\cal O}_K=\pK. \tag*{$\fin$}
\]
\end{proof}

\begin{teorema}[Lema de Nakayama]\label{T7.3.6}
Sea $A$ un anillo conmutativo con unidad y sea ${\eu a}\subseteq
\bigcap\limits_{{\eu m}\text{\ maximal}}{\eu m}$. Si $M$ es
un $A$--m\'odulo finitamente generado tal que ${\eu a}M=M$
entonces $M=0$.
\end{teorema}

\begin{proof}
Supongamos que $M\neq 0$ y sea $m$ el m{\'\i}nimo n\'umero
de generadores de $M$ como $A$--m\'odulo. Digamos que $M=
A\omega_1+\cdots+A\omega_m$. Se tiene que $\omega_1=
a_1\omega_1+\cdots+a_m\omega_m$ para algunos $a_i\in{\eu a}$.
Por tanto $(1-a_1)\omega_1=a_2\omega_2+\cdots+a_m\omega_m$.

Ahora $1-a_1\in A^{\ast}$ pues en caso contrario existir{\'\i}a ${\eu m}
$ maximal tal que $1-a_1\in{\eu m}$ pero $a_1\in{\eu a}\subseteq
{\eu m}$ implicar{\'\i}a $1\in{\eu m}$. Por tanto 
Si $m=1$ entonces $(1-a_1)\omega_1=0$ por lo que $\omega_1
=0$, lo cual es absurdo. Si $m\geq 2$, entonces
$\omega_1\in\langle
\omega_2,\ldots,\omega_m\rangle$ y en este caso $M$ est\'a
generado por $m-1$ elementos lo cual es absurdo. Por tanto $M=0$.
$\fin$
\end{proof}

Como en el caso finito, tenemos la transitividad sobre los ideales
primos que se encuentran sobre uno dado del campo base, es
decir:

\begin{proposicion}\label{P7.3.7} Sea $L/K$ una extensi\'on de
Galois de campos num\'ericos. Sean $\pL$ y $\pL'$ dos ideales
primos de ${\cal O}_{L}$ sobre el ideal primo $\pK$ de ${\cal O}_K$.
Entonces existe $\sigma \in \Gal(L/K)$ tal que $\sigma\pL=\pL'$.
\end{proposicion}

\begin{proof}
Este resultado lo conocemos cuando la extensi\'on $L/K$ es finita.
En general sean
\[
K=F_0\subseteq F_1\subseteq \cdots \subseteq F_n\subseteq
\cdots \subseteq L, \quad L=\bigcup_{n=1}^{\infty} F_n
\]
tales que $F_n/K$ son extensiones finitas de Galois,
lo cual es posible hacerlo pues $L$ es un conjunto numerable. Sean
$\pK_n:=\pL\cap {\cal O}_{F_n}$, $\pK'_n:=\pL'\cap{\cal O}_{F_n}$.
Existe $\tau_n\in \Gal(F_n/K)$ tal que $\tau_n(\pK_n)=\pK'_n$. Sean
$\sigma_n\in \Gal(L/K)$ tales que $\sigma_n|_{F_n}=\tau_n$. 
Puesto que $\Gal(L/K)$ es compacto, la sucesi\'on $\big\{\sigma_n
\big\}_{n=1}^{\infty}$ tiene un punto de acumulaci\'on $\sigma$.
Sea $\sigma_{n_i}\xrightarrow[i\to\infty]{}\sigma$.
Por facilidad suponemos $\sigma_n\xrightarrow[n]{}\sigma$. Sea
$m$ arbitrario. Se tiene que $\Gal(L/F_m)$ es una vecindad
abierta de $\Id$ y $\sigma^{-1}\sigma_n\in\Gal(L/F_m)$ para $n\gg
m$. Por tanto para $n$ suficientemente grande tenemos
\[
\sigma^{-1}\sigma_n \pK_m=\pK_m.\quad\text{Por tanto}\quad
\sigma\pK_m=\sigma_n\pK_m=\sigma_n\big(\pK_n\cap {\cal O}_{F_m}
\big)=\pK'_n\cap{\cal O}_{F_m}=\pK'_m.
\]
Puesto que $\pL=\bigcup\limits_{m=1}^{\infty}\pK_m$ y $\pL'=
\bigcup\limits_{m=1}^{\infty}\pK'_m$ se sigue que $\sigma\pL=
\pL'$. $\fin$
\end{proof}

\begin{ejemplo}\label{Ej7.3.8} Este es un ejemplo de un anillo
${\cal O}_K$ que no es dominio Dedekind. Sea $K=\cic p{\infty}$
y sea $\pK=\bigcup\limits_{n=1}^{\infty}\pK_n$ donde $\pK_n=
\langle 1-\zeta_{p^n}\rangle$ es el ideal primo sobre $p$ de $\cic pn$.
Se tiene que 
 $\langle 1-\zeta_{p^n}\rangle$ es el ideal primo sobre $p$
de $\cic pn$ y $\langle 1-\zeta_{p^n}\rangle^{\varphi(p^n)}=
\langle p\rangle$. En particular $\pK=\langle 1-\zeta_p, 1-\zeta_{p^2},
\ldots, 1-\zeta_{p^n},\ldots\rangle$ y puesto que $\langle 1-
\zeta_{p^n}\rangle^p=\langle 1-\zeta_{p^{n-1}}\rangle$ se sigue que
$\pK^p=\pK$. Esto implica lo que queremos.

Otra forma de obtener lo mismo es notando que $\pK$ no es
finitamente generado pues si $\pK=\langle a_1,\ldots, a_m\rangle$,
entonces cada $a_i\in \pK_{t_i}$ para alg\'un $t_i$. Sea $t:=\max\{
t_i\mid 1\leq i\leq m\}$. Por tanto $a_1,\ldots,a_m\in \pK_t=\langle
1-\zeta_{p^t}\rangle$ pero esto es absurdo pues $1-\zeta_{p^{t+1}}
\notin\pK_t$.
\end{ejemplo}

El Ejemplo \ref{Ej7.3.8} nos indica, entre otras cosas, que no 
podemos estudiar ramificaci\'on v{\'\i}a la descomposici\'on de
primos. En su lugar lo hacemos usando los grupos de inercia.
Esto lo hacemos como en el caso finito.

\begin{definicion}\label{D7.3.9} {\ }
\begin{window}[0,l,\xymatrix{
\pL\ar@{-}[r]\ar@{-}[d]&K\ar@{-}[d]\\ \pK\ar@{-}[r]& k},{}]
Sea $K/k$ una extensi\'on de
Galois. Sea $\pK$ un primo de $k$ y sea $\pL$ un primo de $L$
sobre $\pK$. Se define el {\em grupo de descomposici\'on\index{grupo
de descomposici\'on}} de $\pL$ sobre $\pK$ por $D(\pL|\pK)=D=
\{\sigma\in \Gal(K/k)\mid \sigma \pL=\pL\}$.
El {\em grupo de inercia\index{grupo de inercia}} de $\pL$ sobre 
$\pK$ se define por $
I(\pL|\pK)=I=\{\sigma\in D\mid \sigma\alpha\equiv \alpha\bmod \pL
\ \forall\ \alpha\in{\cal O}_K\}$.
\end{window}
\end{definicion}

\begin{proposicion}\label{7.3.10} Se tiene que $D$ y $I$ son 
subgrupos cerrados de $G=\Gal(K/k)$ y por tanto son compactos.
\end{proposicion}

\begin{proof}
Sea $k=F_0\subseteq F_1\subseteq \ldots F_n\subseteq \ldots
\subseteq K$ tal que $K=\bigcup\limits_{n=1}^{\infty}F_n$ y $F_n/k$
de Galois y finita. Sean $\pL_n:=\pL\cap {\cal O}_{F_n}$, 
$D_n=\{\sigma \in G\mid \sigma \pL_n=\pL_n\}$. Se tiene que
$\pL=\bigcup\limits_{n=1}^{\infty}\pL_n$, $D\subseteq D_n$ y $
D=\bigcap\limits_{n=1}^{\infty} D_n$. Ahora bien, se tiene $\Gal(K/F_n)
\subseteq D_n$ y $\Gal(K/F_n)$ es abierto por lo que $D_n$ es
abierto. Se sigue que $D_n$ es cerrado
por la compacidad de $G$. Se sigue que $D$ es cerrado.

\begin{window}[2,l,\xymatrix{K\ar@{-}[d]_N\\ \widetilde{k(\alpha)}},{}]
Veamos que $I$ es cerrado. Sea $\sigma \notin I$. Entonces existe
$\alpha\in{\cal O}_K$ tal que $\sigma \alpha-\alpha\notin \pL$. 
Consideremos $\widetilde{k(\alpha)}/k$ la cerradura de Galois de
$k(\alpha)/k$. Se tiene que $\widetilde{k(\alpha)}/k$ es finita y
$N:=\Gal(K/\widetilde{k(\alpha)})$ es un subgrupo abierto de $G$. 
Ahora bien $\sigma N$ es una vecindad abierta de $\sigma$. Sea
$\psi\in N$. Entonces $\psi(\alpha)=\alpha$ por lo que $\sigma\psi(
\alpha)=\sigma(\alpha)$ y $(\sigma\psi)(\alpha)-\alpha=\sigma\alpha
-\alpha\notin \pL$. Por tanto $\sigma N\cap I=\emptyset$ de donde
se sigue que el complemento de $I$ es abierto y por tanto $I$
es cerrado. $\fin$
\end{window}
\end{proof}

Sea $\varphi\colon D\lra\Gal\big({\cal O}_K/\pL:{\cal O}_k/\pK\big)=:
G(\pL)$ el mapeo natural, es decir, si $\sigma\in D$, $\sigma
{\cal O}_K={\cal O}_K$, $\sigma\pL=\pL$ y $\sigma|_{{\cal O}_K}=
\Id_{{\cal O}_k}$ por lo que $\tilde{\sigma}\colon {\cal O}_K/\pL
\to{\cal O}_{K}/\pL$, $\sigma(x\bmod \pL)=\sigma x\bmod \pL$ est\'a
bien definida y $\tilde{\sigma}\in G(\pL)$. Entonces $\varphi(\sigma)=
\tilde{\sigma}$. Se tiene:

\begin{teorema}\label{T7.3.11}
La sucesi\'on $1\longto I\longto D\stackrel{\varphi}{\longto}G(\pL)
\longto 1$ es exacta.
\end{teorema}

\begin{proof} Sabemos que en el caso finito, $\varphi$ es
suprayectiva. Ahora $G(\pL)=\lim\limits_{\longleftarrow}\Gal\big(
({\cal O}_{F_n}/\pL_n)/({\cal O}_k/\pK)\big)$. Sea $\tilde{D}_n:=
D(\pL_n|\pK)\subseteq \Gal(F_n/k)$. Se tiene $D=\lim\limits_{
\longleftarrow}\tilde{D}_n$ bajo los mapeos:
\begin{align*}
\varphi\colon D&\longto \prod_{n=1}^{\infty}\tilde{D}_n\subseteq
\prod_{n=1}^{\infty}\Gal(F_n/k)\\
\sigma& \longmapsto \prod_{n=1}^{\infty}\sigma|_{F_n},\quad
\sigma|_{F_n}\in \Gal(F_n/k)
\end{align*}
$\sigma|_{F_n}(\pL_n)=\sigma\big(\pL\cap {\cal O}_{F_n}\big)=
\pL\cap {\cal O}_{F_n}=\pL_n$. Entonces tenemos el diagrama
conmutativo
\[
\xymatrix{ 
\tilde{D}_n\ar[rr]^{\hspace{-1.2cm}\varphi_n}&&
\Gal\big(({\cal O}_{F_n}/\pL_n)/({\cal O}_k/
\pK)\big)\\ \tilde{D}_{n+1}\ar[u]\ar[rr]^{\hspace{-1.75cm}
\varphi_{n+1}}&&\Gal\big(({\cal O}_{
F_{n+1}}/\pL_{n+1})/({\cal O}_k/\pK)\big)\ar[u]
}
\]
donde las flechas verticales son los mapeos de restricci\'on.

Ahora bien puesto que los mapeos $\varphi_n$ son suprayectivos,
pasando al l{\'\i}mite, se sigue que $\varphi$ es suprayectiva y por 
definici\'on tenemos que $\ker \varphi=I$. $\fin$
\end{proof}

\begin{corolario}\label{C7.13.12} Se tiene $D/I\cong\Gal
\big(({\cal O}_K/\pL)/({\cal O}_k/\pK)\big)$. $\fin$
\end{corolario}
\[
\xymatrix{\overline{{\ma Q}}\ar@{-}[r]\ar@{-}[d]&
{\eu B}\ar@{-}[d]\\ K\ar@{-}[d]\ar@{-}[r]
&\pL\ar@{-}[d]\\ k\ar@{-}[r]& {\pK}}
\]

Ahora supongamos que $K/k$ es algebraica pero no necesariamente
de Galois. Fijemos una cerradura algebraica $\overline{\ma Q}$
de ${\ma Q}$. Entonces $\overline{\ma Q}/K$ y $\overline{\ma Q}/k$
son extensiones de Galois. Sean $\pL$ y $\pK$ como antes y sea
${\eu B}$ un primo ${\cal O}_{\overline{\ma Q}}$ tal que ${\eu B}|_K
=\pL$. Entonces
\begin{align*}
I({\eu B}|\pK)&\subseteq \Gal(\overline{\ma Q}/k),\\
I({\eu B}|\pL)&\subseteq \Gal(\overline{\ma Q}/K)\subseteq
\Gal(\overline{\ma Q}/k),\\
I({\eu B}|\pL)&=I({\eu B}|\pK)\cap\Gal(\overline{\ma Q}/K).
\end{align*}

\begin{definicion}\label{D7.13.13} Se define el {\em {\'\i}ndice de
ramificaci\'on\index{indice@{\'\i}ndice de ramificaci\'on}} $e(\pL|\pK)$ por:
$e(\pL|\pK)=[I({\eu B}|\pK):I({\eu B}|\pL)]$ el cual puede ser
infinito.
\end{definicion}

\begin{observacion}\label{O7.13.14}
$e(\pL|\pK)$ no depende de ${\eu B}$ pues si ${\eu B}'$ es otro
primo de $\overline{\ma Q}$ sobre $\pL$, entonces existe $\sigma\in
\Gal(\overline{\ma Q}/K)$ tal que ${\eu B}'=\sigma{\eu B}$ y por 
tanto $I({\eu B}'|\pK)=\sigma I({\eu B}|\pK)\sigma^{-1}$ y 
$I({\eu B}'|\pL)=\sigma I({\eu B}|\pL)\sigma^{-1}$ y por lo tanto
\[
[I({\eu B}'|\pK):I({\eu B}'|\pL)]=[\sigma I({\eu B}|\pK)\sigma^{-1}
:\sigma I({\eu B}|\pL)\sigma^{-1}]=[I({\eu B}|\pK):I({\eu B}|\pL)].
\]
\end{observacion}

En el caso en que $K/k$ sea una extensi\'on de Galois, la restricci\'on
\begin{eqnarray*}
 \Gal\big(\overline{\ma Q}/k\big)&\longto&\Gal(K/k)\\
\sigma&\longmapsto&\sigma|_K
\end{eqnarray*}
 tiene como n\'ucleo
a $\Gal\big(\overline{\ma Q}/K\big)$ y $I({\eu B}|\pK)\longto I(\pL|\pK)$
es suprayectiva con n\'ucleo $I({\eu B}|\pL)$, es decir
\[
\frac{I({\eu B}|\pK)}{I({\eu B}|\pL)}\cong I(\pL|\pK)\quad \text{y}\quad
e(\pL|\pK)=|I(\pL|\pK)|.
\]

En el caso de lugares arquimedianos o infinitos, procedemos de la
siguiente forma:

\begin{definicion}\label{D7.13.15} Un {\em lugar
arquimediano\index{lugar arquimediano}} de $k$ o un {\em lugar
infinito\index{lugar infinito}} de $k$ es, ya sea un encaje
real $\phi\colon k\lra {\ma R}$ o bien un par $(\psi,\overline{\psi})$
de encajes complejos $\psi,\overline{\psi}\colon k\lra{\ma C}$, 
$\psi\neq \overline{\psi}$.
\end{definicion}

Por el Lema de Zorn, tenemos que cualquier encaje $\phi$ o $(\psi,
\overline{\psi})$ se puede extender a un encaje: $\overline{\ma Q}
\lra{\ma C}$ y en particular se puede extender a un encaje de $K$:
$K\lra {\ma C}$.

Sea $K/k$ una extensi\'on de Galois y sean $\phi_1$, $\phi_2$ dos
extensiones a $K$ de un encaje $\phi$ de $k$. Entonces $\phi_2^{-1}
\circ \phi_1\in\Gal(K/k)$ y por tanto $\phi_1=\phi_2\sigma$ para
alg\'un $\sigma\in\Gal(K/k)$. Similarmente, si $\big(\psi_1,
\overline{\psi}_1\big)$ y $\big(\psi_2,\overline{\psi}_2\big)$ son dos
extensiones de $\phi$, $\psi_1=\psi_2\sigma$ y $\big(\psi_1,
\overline{\psi}_1\big)=\big(\psi_2,\overline{\psi}_2\big)\circ \sigma =
\big(\psi_2\circ \sigma, \overline{\psi}_2\circ\sigma\big)$.

De manera an\'aloga si $\big(\psi_1,\overline{\psi}_1\big)$ y $\big(
\psi_2,\overline{\psi}_2\big)$ se pueden extender a $\big(\psi,
\overline{\psi}\big)$, entonces $\psi_2^{-1}\psi_1\in\Gal(K/k)$ y
$\psi_1=\psi_2\sigma$, $\overline{\psi}_1=\overline{\psi}_2\sigma$.

Si $w$ es un lugar arquimediano de $K$, $w|_k=v$ es
un lugar arquimediano de $k$ y se define el grupo de inercia y el 
grupo de descomposici\'on de $w$ sobre $v$ por:

\begin{definicion}\label{D7.3.16} Sea $K/k$ una extensi\'on de 
Galois. Entonces se define
\[
I(w|v)=D(w|v)=\{\sigma\in \Gal(K/k)\mid w\sigma=w\}.
\]
\end{definicion}

Si $w$ es real, $w\sigma=w$ implica $\sigma=w^{-1}w=\Id$.

Si $w$ es complejo, $w=\big(\psi,\overline{\psi}\big)$ y $\big(\psi,
\overline{\psi}\big)\sigma=\big(\psi,\overline{\psi}\big)$ entonces
$\psi\sigma=\psi$ o $\overline{\psi}\sigma=\psi$, es decir, $\sigma
=\Id$ o $\sigma=\overline{\psi}^{-1}\circ \psi$. Si $v=\big(\phi,\overline{
\phi}\big)$ es complejo, $\Id=\sigma|_k=\overline{\psi}\circ\psi|_k=
\overline{\phi}^{-1}\circ \phi\neq \Id$ no es posible por lo que si $v$
es complejo, $\sigma=\Id$. Si $v$ es real, $\overline{\psi}^{-1}
\circ \psi|_k=\phi^{-1}\circ \phi=\Id$ por lo que existe $\overline{\psi}^{
-1}\circ \psi=\psi^{-1}\circ \overline{\psi}=\sigma\neq \Id$.

En resumen,

\begin{teorema}\label{T7.3.16} Se tiene $|I(w,v)|=1$ o $2$ y $|I(w,v)|=
2\iff v$ es real y $w$ es complejo. En este \'ultimo caso, si
$w=\big(\psi,\overline{\psi}\big)$, $I(w,v)=\{\Id, \overline{\psi}^{-1}
\circ \psi=\psi^{-1}\circ\overline{\psi}\}$. $\fin$
\end{teorema}

\section{Pro--$p$--grupos}\label{S7.4}

\begin{definicion}\label{D7.4.1} Sea $p$ un n\'umero primo. Un grupo
profinito $H$ se llama {\em pro--$p$--grupo\index{pro--$p$--grupo}}
si $H$ es el l{\'\i}mite proyectivo de $p$--grupos finitos:
\[
H=\lim_{\substack{\longleftarrow\\ N}}H/N, \quad  [H:N]=p^n,\quad
n\in{\ma N}, \quad N\lhd H.
\]
\end{definicion}

\begin{ejemplo}\label{Ej7.4.2} ${\ma Z}_p$ es un pro--$p$--grupo
(abeliano): ${\ma Z}_p=\lim\limits_{\substack{\longleftarrow\\ n}}
{\ma Z}/p^n{\ma Z}$.
\end{ejemplo}

\begin{definicion}\label{D7.4.3} Sea $I$ un conjunto y sea $L(I)$ el
grupo discreto libre generado por elementos $x_i$, $i\in I$. Sea
$X$ la familia de subgrupos normales $M$ de $L(I)$ tales que:
\l
\item $L(I)/M$ es un $p$--grupo finito.
\item $M$ contiene a casi todo $x_i$, es decir, a todo $x_i$ salvo
un n\'umero finito.
\end{list}

Sea $F(I):=\lim\limits_{\substack{\longleftarrow\\M\in X}} 
L(I)/M$. Entonces el pro--$p$--grupo
$F(I)$ se llama el {\em pro--$p$--grupo libre generado
por los $x_i$\index{pro--$p$--grupo libre}}.
\end{definicion}

\begin{ejemplo}\label{Ej7.4.4} Sea $I$ finito, $|I|=n\in {\ma N}\cup \{
0\}$. Ponemos $F(I)=F(n)$. Se tiene $F(0)=\{1\}$, $F(1)\cong
{\ma Z}_p$ pues $L(1)={\ma Z}$ y $F(1)=\lim\limits_{\substack{
\longleftarrow\\ n}}{\ma Z}/p^n{\ma Z}\cong{\ma Z}_p$.
\end{ejemplo}

\begin{proposicion}\label{P7.4.5} Se tiene $F(n)/F'(n)\cong
{\ma Z}_p^n$ y m\'as generalmente, para $I$ arbitrario, $F(I)/F'(I)
\cong\bigoplus\limits_{i\in I}{\ma Z}_p$.
\end{proposicion}

\begin{proof} (Esquema). Se tiene 
\begin{gather*}
F(I)=\lim\limits_{\longleftarrow}
L(I)/H,\quad F'(I)=\lim\limits_{\longleftarrow}L'(I)/(H\cap L'(I))\\
\intertext{y $L'(I)/(H\cap L'(I)) \cong HL'(I)/H$. Por tanto}
\frac{L(I)/H}{L'(I)/(H\cap L'(I))}\cong \frac{L(I)/H}{HL'(I)/H}=
L(I)/HL'(I)=\Big(\bigoplus_{i\in I}{\ma Z}\Big)/H\quad\text{y}\\
\lim_{\longleftarrow}\Big(\bigoplus_{i\in I}{\ma Z}\Big)/H\cong
\bigoplus_{i\in I}{\ma Z}_p. \tag*{$\fin$}
\end{gather*}
\end{proof}

\section{Extensiones ${\ma Z}_p$}\label{S7.5}\label{S7.4.1}

Consideremos una cadena de campos
\[
K=K_0\subseteq K_1\subseteq K_2\subseteq \ldots\subseteq K_n
\subseteq \ldots\subseteq K_{\infty}=\bigcup_{n=0}^{\infty}K_n
\]
tal que $\Gal(K_n/K)\cong {\ma Z}/p^n{\ma Z}$. Entonces
$\Gal(K_{\infty}/K)=\lim\limits_{\longleftarrow}{\ma Z}/p^n{\ma Z}\cong
{\ma Z}_p$.

\begin{definicion}\label{D7.4.1.1} En el desarrollo anterior, $K_{\infty}/
K$ es una {\em extensi\'on ${\ma Z}_p$\index{extensi\'on
${\ma Z}_p$}}
o una {\em extensi\'on $\Gamma$\index{extensi\'on $\Gamma$}}
donde ponemos $\Gamma:={\ma Z}_p$.
\end{definicion}

\begin{ejemplo}[Extensi\'on $p$--ciclot\'omica de
${\ma Q}$\index{extensi\'on $p$--ciclot\'omica de ${\ma Q}$}]
\label{Ej7.4.1.2}
(Ver Teoremas \ref{T8.6}, \ref{T8.7} y \ref{T8.8}
y la Secci\'on \ref{S4.3}). Sea $K={\ma Q}$.
Se tiene que $\cic pn/{\ma Q}$ tiene como grupo de Galois a 
$\Gal(\cic pn/{\ma Q})\cong U_{p^n}\cong \big({\ma Z}/p^n{\ma
Z}\big)^{\ast}$.

Primero consideremos $p>2$. Entonces $U_{p^n}$ es c{\'\i}clico y se
tiene $\big({\ma Z}/p^n{\ma Z}\big)^{\ast}\cong {\ma Z}/(p-1){\ma Z}
\times \big({\ma Z}/p^{n-1}{\ma Z}\big)$ (Teorema \ref{T1.2.1.3}).

En particular $U_{p^n}$ contiene un \'unico subgrupo de 
orden $p-1$, a saber, ${\ma Z}/(p-1){\ma Z}$ en la descomposici\'on
anterior. Sea ${\ma Q}_{n-1}:=\cic pn^{{\ma Z}/(p-1){\ma Z}}$. 
Entonces $\Gal({\ma Q}_{n-1}/{\ma Q})\cong{\ma Z}/(p^{n-1}{\ma Z})$.

Notemos que $2$ divide a $p-1$ por lo que ${\ma Q}_{n-1}
\subseteq {\ma Q}(\zeta_{p^n}+\zeta_{p^n}^{-1})=\cic pn^+
\subseteq {\ma R}$.

Sea ${\ma Q}_{\infty}:=\bigcup\limits_{n=1}^{\infty}{\ma Q}_n$.
Entonces $\Gal({\ma Q}_{\infty}/{\ma Q})\cong {\ma Z}_p$.

Ahora si $p=2$, tenemos para $n\geq 2$, que
\[
\Gal(\cic 2n/{\ma Q})\cong U_{2^n}\cong\langle \pm 1\rangle\times
\langle 1+2^2\rangle\cong ({\ma Z}/2{\ma Z})\times ({\ma Z}/2^{n-2}
{\ma Z}).
\]

\begin{window}[0,l,\xymatrix{
\cic 4{}\ar@{-}[rr]\ar@{-}[dd]&&\cic 8{}=\cic 23\ar@{-}[dd]\\
&E_1\ar@{-}[dl]\\{\ma Q}\ar@{-}[rr]&&\cic 8{}^+},{}]
Si $n=0$, $\zeta_{2^0}=\zeta_1=1$; si $n=1$, $\zeta_{2^1}=\zeta_2
=-1$; si $n=2$, $\zeta_{2^2}=\zeta_4=i$. Para $n=3$, $U_8\cong
C_2\times C_2$ y hay $3$ subgrupos de orden $2$. Por lo tanto
hay $3$ subextensiones cuadr\'aticas, a saber, $\cic 4{}$, $\cic 8{}^+=
{\ma Q}(\zeta_8+\zeta_8^{-1})={\ma Q}(\sqrt{2})$ y $E_1={\ma Q}(
\zeta_4(\zeta_8+\zeta_8^{-1}))={\ma Q}(\sqrt{-2})$.
\end{window}

Para $n\geq 4$, $U_2^n\cong C_2\times C_{2^{n-2}}$ y por lo tanto 
$U_{2^n}$ tiene $2$ subgrupos c{\'\i}clicos de orden $2^{n-2}$ y por
lo tanto $2$ grupos cociente c{\'\i}clicos de orden $2^{n-2}$.
\begin{scriptsize}
\[
\xymatrix{
\cic 4{}\ar@{--}[r]\ar@{-}[dd]&\cic 2{n-1}\ar@{-}[rr]\ar@{-}[dd]
&&\cic 2n\ar@{-}[rr]\ar@{-}[dd]&&\cic 2{n+1}\ar@{-}[dd]\\
&& E_{n-2}\ar@{-}[dl]&&E_{n-1}\ar@{-}[dl]\\
{\ma Q}\ar@{--}[r]&\cic 2{n-1}^+\ar@{-}[rr]&&\cic 2n^+\ar@{-}[rr]
&&\cic 2{n+1}^+
}
\]
\end{scriptsize}

Hay $3$ subcampos de $\cic 2n$ de grado $2^{n-2}$ sobre
${\ma Q}$, $2$ de ellos son c{\'\i}clicos sobre ${\ma Q}$ y el otro
no (para $n\geq 4$). Los $3$ subcampos son $\cic 2{n-1}$,
$\cic 2n^+$ y $E_{n-2}={\ma Q}\big(\zeta_4(\zeta_{2^{n-1}}+
\zeta_{2^{n-1}}^{-1})\big)$. Puesto que $\cic 2{n-1}/{\ma Q}$ no
es c{\'\i}clico, $\cic 2n^+$ y $E_{n-2}$ son las extensiones 
c{\'\i}clicas de ${\ma Q}$ de grado $2^{n-2}$.

Se tiene que $E_{n-2}\nsubseteqq E_{n-1}$ pues de lo contrario
tendr{\'\i}amos que $\cic 2n=E_{n-2}\cic 2{n-1}^+\subseteq E_{n-1}$
lo cual es absurdo.

Ahora bien $\cic 2n^+$ corresponde al subgrupo $\{\pm 1\}$ de
$U_{2^n}$, es decir, $\cic 2n^+=\cic 2n^{\{\pm 1\}}$ y $\Gal\big(
\cic 2n^+/{\ma Q}\big)\cong U_n/\{\pm 1\}\cong \langle 1+2^2\rangle
\cong C_{2^{n-2}}$.

En particular $\bigcup\limits_{n=1}^{\infty}\cic 2n^+=\cic 2{\infty}^+$
satisface:
\[
\Gal\big(\cic 2{\infty}^+/{\ma Q}\big)\cong \lim_{\substack{\longleftarrow
\\ n}}C_{2^{n-2}}\cong {\ma Z}_2.
\]
Ponemos ${\ma Q}_{\infty}:=\cic 2{\infty}^+$.
\end{ejemplo}

\begin{observacion}\label{O7.4.1.3}
Tenemos que para cualquier $p\geq 2$, ${\ma Q}_{\infty}$ es la
\'unica extensi\'on de ${\ma Q}$ tal que $\Gal({\ma Q}_{\infty}/{\ma Q})
\cong {\ma Z}_p$.

Para ver esto, sea $K/{\ma Q}$ tal que $\Gal(K/{\ma Q})\cong
{\ma Z}_p$, $p\geq 2$. Si $K\subseteq \bigcup\limits_{n=1}^{\infty}
\cic pn=\cic p{\infty}$, el resultado se sigue del hecho de que 
\[
\Gal
\big(\cic p{\infty}/{\ma Q}\big)\cong \begin{cases}
C_2\times {\ma Z}_2&\text{si $p=2$}\\
C_{p-1}\times {\ma Z}_p&\text{si $p\geq 3$}\end{cases}
\]
y el \'unico subgrupo finito de $\Gal\big(\cic p{\infty}/{\ma Q}\big)$ es
$C_2$ si $p=2$ y $C_{p-1}$ si $p\geq 3$.

Ahora bien, en general, sea $p^n{\ma Z}_p<{\ma Z}_p\cong
\Gal(K/{\ma Q})$ y sea $K_n:=K^{p^n{\ma Z}_p}$, $\Gal(K_n/{\ma Q})
\cong {\ma Z}_p/p^n{\ma Z}_p\cong {\ma Z}/p^n{\ma Z}$.

Por el Teorema de Kronecker--Weber, de hecho por el Teorema
\ref{T1.4.3*}, $K\subseteq \cic p{n+1}$ y por tanto $K=\bigcup\limits_{
n=1}^{\infty}K_n\subseteq \bigcup\limits_{n=1}^{\infty}\cic p{n+1}=\cic
p{\infty}$ y se sigue la unicidad.
\end{observacion}

Veamos con m\'as detalle $\cic 2n^+={\ma Q}(\zeta_{2^n}+\zeta_{
2^n}^{-1})$ (ver Teorema \ref{T10.1.1}).
Sea $\alpha_n:=\zeta_{2^n}+\zeta_{2^n}^{-1}$. Entonces
\begin{align*}
\alpha_0&= \zeta_1+\zeta_1^{-1}=1+1=2,\\
\alpha_1&= \zeta_2+\zeta_2^{-1}=-1-1=-2,\\
\alpha_2&= \zeta_4+\zeta_4^{-1}=i+\frac{1}{i}=0,\\
\alpha_3&= \zeta_8+\zeta_8^{-1}.
\end{align*}
Notemos que para $n\geq 3$, $\alpha_n^2=\big(\zeta_{2^n}+
\zeta_{2^n}^{-1}\big)^2=\zeta_{2^{n-1}}+\zeta_{2^{n-1}}^{-1}+2=
\alpha_{n-1}+2$. Adem\'as $\alpha_n\geq 0$, por lo que 
$\alpha_n=\sqrt{\alpha_{n-1}+2}$, esto es, $\alpha_3=\sqrt{2}$,
$\alpha_4=\sqrt{2+\sqrt{2}}$ y en general $\alpha_n=\sqrt{
\underbrace{2+\sqrt{2+\sqrt{2+\sqrt{2+\cdots+\sqrt{2}}}}}_{n-2}}$.

\begin{definicion}\label{D7.4.1.4}
Para todo primo $p\geq 2$, la extensi\'on ${\ma Q}_{\infty}\subseteq
\cic p{\infty}$ que satisface $\Gal({\ma Q}_{\infty}/{\ma Q})\cong
{\ma Z}_p$ se le llama {\em la ${\ma Z}_p$--extensi\'on
ciclot\'omica\index{extensi\'on ${\ma Z}_p$--ciclot\'omica}} de ${\ma Q}
$. Notemos que ${\ma Q}_{\infty}\subseteq {\ma R}$, y por tanto,
${\ma Q}_{\infty}\subseteq \cic p{\infty}^+=\cic p{\infty}\cap {\ma R}$.
\end{definicion}

\[
\xymatrix{
&K\ar@{-}[r]\ar@{-}[d]&K_{\infty}=K{\ma Q}_{\infty}\ar@{-}[d]\\
&K\cap {\ma Q}_{\infty}\ar@{-}[r]\ar@{-}[dl]&{\ma Q}_{\infty}\\ {\ma Q},}
\]
Ahora sea $K$ cualquier campo num\'erico finito, es decir, $[K:{\ma Q}
]<\infty$. Sea $K_{\infty}:=K{\ma Q}_{\infty}$. Se tiene que $\Gal(K_{
\infty}/K)\cong \Gal({\ma Q}_{\infty}/K\cap {\ma Q}_{\infty})\subseteq
{\ma Z}_p$. Por tanto $\Gal({\ma Q}_{\infty}/K\cap {\ma Q}_{\infty})
\cong p^n{\ma Z}_p$ para alg\'un $n\geq 0$ y en particular $\Gal(
K_{\infty}/K)\cong p^n{\ma Z}_p\cong {\ma Z}_p$, es decir $K_{\infty}
/K$ es una extensi\'on ${\ma Z}_p$.

\begin{definicion}\label{D7.4.1.5}
Esta extensi\'on $K_{\infty}/K$ se llama la {\em ${\ma Z}_p$
extensi\'on ciclot\'omica\index{extensi\'on ${\ma Z}_p$--ciclot\'omica
de un campo num\'erico}} de $K$.
\end{definicion}

\begin{observacion}\label{O7.4.1.6} Sea $p\geq 2$ un n\'umero
primo cualquiera. Entonces $\cic {2p}{}=\cic p{}$ si $p>2$ y $\Gal\big(
\cic {2p}{}/{\ma Q}\big)\cong C_2$ si $p=2$ o $C_{p-1}$ si $p>2$.
Entonces
\[
\xymatrix{
\cic {2p}{}\ar@{-}[r]^{{\ma Z}_p}\ar@{-}[d]_{C_{p-1} \text{\ o\ }C_2}
&\cic p{\infty}\ar@{-}[d]^{C_{p-1}\text{\ o\ }C_2}\\
{\ma Q}\ar@{-}[r]_{{\ma Z}_p}&{\ma Q}_{\infty}
}
\qquad \Gal\big(\cic p{\infty}/\cic {2p}{}\big)\cong {\ma Z}_p.
\]
\end{observacion}

\begin{proposicion}\label{P7.4.1.7} Sea $K$ una extensi\'on finita
de ${\ma Q}$ y $K_{\infty}$ la ${\ma Z}_p$--extensi\'on ciclot\'omica
de $K$. Entonces $\zeta_{p^{\infty}}\in K_{\infty}\iff \zeta_{2p}\in K$
y en este caso $K_{\infty}=K(\zeta_{p^{\infty}})$.
\end{proposicion}

\begin{proof}
Primero, si $\zeta_{2p}\in K$, entonces $K_{\infty}=K{\ma Q}_{\infty}
\supseteq \cic {2p}{}{\ma Q}_{\infty}=\cic p{\infty}$ por lo que $\zeta_{
p^{\infty}}\in K$.

Rec{\'\i}procamente, si $\zeta_{p^{\infty}}\in K_{\infty}=K{\ma Q}_{
\infty}$ se obtiene $K_{\infty}=K(\zeta_{p^{\infty}}){\ma Q}_{
\infty}=K(\zeta_{p^{\infty}})$. Se tiene que
\begin{gather*}
\xymatrix{
&K\ar@{-}[r]^{\hspace{-2cm}{\ma Z}_p}
\ar@{-}[d]&K_{\infty}=K{\ma Q}_{\infty}=K
{\ma Q}(\zeta_{p^{\infty}})\ar@{-}[d]\\
&K\cap {\ma Q}(\zeta_{p^{\infty}})\ar@{-}[r]^{{\ma Z}_p}
\ar@{-}[dl]&{\ma Q}(\zeta_{p^{\infty}})\\{\ma Q}}\\
 \Gal\Big(\cic p{\infty}/(K\cap \cic p{\infty})\Big)\cong
\Gal(K{\ma Q}_{\infty}/K)\cong {\ma Z}_p.
\end{gather*}

De esta forma se sigue que $\Gal(\cic p{\infty}/{\ma Q})\cong C_t\times
{\ma Z}_p$ donde $t=2$ si $p=2$ y $t=p-1$ para $p\geq 3$. Los
subgrupos cerrados de $C_t\times {\ma Z}_p$ isomorfos a ${\ma Z}_p
$ son \'unicamente los grupos $p^n{\ma Z}_p$ y por tanto $K\cap
\cic p{\infty}$ corresponde a alg\'un $p^n{\ma Z}_p$. Se sigue que
$K\cap \cic p{\infty}={\ma Q}(\zeta_{\infty})^{p^n{\ma Z}_p}=
\cic {2p}{n+1}$ lo cual implica que $\zeta_{2p}\in K$.

Como vimos, en este caso $K_{\infty}=K{\ma Q}_{\infty} =K
{\ma Q}_{\infty}(\zeta_{2p})=K\cic p{\infty}=K(\zeta_{p^{\infty}})$.
$\fin$
\end{proof}

\section{Ramificaci\'on y descomposici\'on de primos en ${\ma Q}_{
\infty}/{\ma Q}$}\label{S7.6}

Sea $\ell$ un primo de ${\ma Q}$ que sea ramificado en ${\ma Q}_{
\infty}/{\ma Q}$. Entonces si ${\eu l}$ es un primo de ${\ma Q}_{\infty}
$ sobre $\ell$, $I:=I({\eu l}|\ell)$ es un subgrupo cerrado de $\Gal(
{\ma Q}_{\infty}/{\ma Q})\cong {\ma Z}_p$ y por lo tanto $I=p^n
{\ma Z}_p$ para alg\'un $n\geq 0$ o $I=0$. Puesto que $\ell$
es ramificado, $I\neq 0$ y por tanto 
\[
I\cong p^n{\ma Z}_p\cong{\ma Z}_p, {\ma Q}^I={\ma Q}_n\quad
\text{pues}\quad {\ma Z}_p/I={\ma Z}_p/p^n{\ma Z}_p\cong
{\ma Z}/p^n{\ma Z}\cong \Gal({\ma Q}_n/{\ma Q}).
\]

En particular $\ell$ es ramificado en ${\ma Q}_{n+1}/{\ma Q}_n$ pero
el \'unico primo ramificado en ${\ma Q}_{n+1}/{\ma Q}_n$ es $p$ de
donde se sigue que $\ell=p$. Por otro lado, $p$ es totalmente
ramificado y por lo tanto s\'olo hay un primo $\pK$ en ${\ma Q}_{
\infty}$ sobre $p$ y $I(\pK|p)={\ma Z}_p=\Gal({\ma Q}_{\infty}
/{\ma Q})$.

Ahora consideremos $\ell\neq p$. Sea $f_n$ el grado de inercia de
$\ell$ en $\cic pn/{\ma Q}$. Por facilidad estudiaremos el caso $p>2$.
El caso $p=2$ es similar aunque t\'ecnicamente m\'as complicado.

Tenemos (Teorema \ref{T8.2} o Teorema \ref{T8.10}) que $f_n=
o(\ell\bmod p^n)$, es decir $\ell^{f_n}\equiv 1\bmod p^n$, con $f_n$
m{\'\i}nimo con esta propiedad. Consideremos el diagrama
\[
\xymatrix{
\cic p{}\ar@{-}[r]^{\tilde{f}_n}\ar@{-}[d]_{f_1}&\cic p{n+1}
\ar@{-}[ld]_{f_{n+1}}\ar@{-}[d]^{f_1'}\\
{\ma Q}\ar@{-}[r]_{\tilde{f}_{n}'}&{\ma Q}_n
}
\]

Puesto que $\mcd (\tilde{f}_n,f_1)=1$, entonces $\tilde{f}_n=
\tilde{f}_n'$. Por otro lado se tiene que $f_{n+1}=f_1\tilde{f}_n$.
Adem\'as, $\tilde{f}_n|\big[\cic p{n+1}:\cic p{}\big]=p^n$. Se tiene
que $f_n\neq 1$ para alg\'un $n$ pues de lo contrario, si $f_n=1$
para toda $n$, entonces $p^n|\ell-1$ para todo
$n$, lo cual es absurdo pues $\ell
\neq 1$.

Sea $n_0\in {\ma N}$ tal que 
$\ell^{f_{n_0}}\equiv 1\bmod p^n$, $f_{n_0}>1$. Si $\big\{f_n\big\}_{
n=1}^{\infty}$ fuese acotada, digamos $f_n\leq M$, entonces $\ell^{
f_n}-1\leq \ell^M-1$ y existir{\'\i}a $m\in{\ma N}$ tal que $p^m\nmid
\ell^{f_m}-1$ lo cual es absurdo.
 Por tanto $f_n\xrightarrow[n\to\infty]{}\infty$.
 
 Sea ${\eu l}$ un primo en $\cic p{\infty}$ sobre $\ell$ y sea $D=D(
 {\eu l}|\ell)$ el cual es cerrado en ${\ma Z}_p$ lo cual implica que 
 $D\cong p^n {\ma Z}_p$ para alguna $n\in{\ma N}\cup \{0\}$ o $D
 =0$. Esto \'ultimo no puede ocurrir pues por lo anterior $\ell$ no 
 puede ser totalmente descompuesto. As{\'\i} $D\cong p^n{\ma Z}_p$,
 $n\geq 0$ y en particular ${\ma Q}_n={\ma Q}^D$ y $\ell$ es 
 totalmente inerte ${\ma Q}_{\infty}/{\ma Q}_n$ y $\ell$ es 
 totalmente descompuesto en ${\ma Q}_n/{\ma Q}$.
 
 Si $\infty$ es el primo infinito, $I=D=\{0\}$.
 
 \begin{observacion}\label{O7.6.1} La \'ultima parte es 
 de hecho otra demostraci\'on
 de que ${\ma Q}_{\infty}=\cic 2{\infty}^+$
 es la \'unica extensi\'on ${\ma Z}_2$ de
 ${\ma Q}$ pues al no ser $\infty$ ramificado entonces se sigue
 que ${\ma Q}_{\infty}\subseteq {\ma R}$ y por tanto ${\ma Q}_{\infty}
 =\cic 2{\infty}^+$.
 \end{observacion}
 
 Resumimos la discusi\'on anterior en el siguiente resultado:
 
 \begin{teorema}\label{T7.6.2} Para $p\geq 2$, se tiene que
 $p$ es el \'unico primo (finito o infinito) que es ramificado en
 ${\ma Q}_{\infty}/{\ma Q}$ siendo adem\'as totalmente ramificado.
 
 Para $\ell\neq p$, existe $n=n(\ell)$ tal que $\ell$ es totalmente
 descompuesto en ${\ma Q}_{n(\ell)}/{\ma Q}$ y totalmente inerte en
 ${\ma Q}_{\infty}/{\ma Q}_{n(\ell)}$. No hay primos finitos totalmente
 descompuestos en ${\ma Q}_{\infty}/{\ma Q}$. Finalmente el primo
 $\infty$ es totalmente descompuesto en ${\ma Q}_{\infty}/{\ma Q}$.
 $\fin$
 \end{teorema}
 
 Sea ahora cualquier campo num\'erico finito y sea $K_{\infty}=K
 {\ma Q}_{\infty}$ la extensi\'on ${\ma Z}_p$--ciclot\'omica de $K$.
 Sea $K\cap {\ma Q}_{\infty}=E$. Entonces ${\ma Q}\subseteq E
 \subseteq {\ma Q}_{\infty}$ y $\Gal({\ma Q}_{\infty}/E)$ es un
 subgrupo cerrado de $\Gal({\ma Q}_{\infty}/{\ma Q})\cong {\ma Z}_p$.
 En particular $\Gal({\ma Q}_{\infty}/E)$ es isomorfo a $p^n{\ma Z}_p$
 para alguna $n\geq 0$ y por tanto $E={\ma Q}_n$.
 \[
 \xymatrix{
&K\ar@{-}[r]\ar@{-}[d]&K_m\ar@{-}[r]\ar@{-}[d]&K_{\infty}\ar@{-}[d]\\
&E={\ma Q}_n\ar@{-}[dl]\ar@{-}[r]&{\ma Q}_{n+m}
 \ar@{-}[r]&{\ma Q}_{\infty}\\
 {\ma Q}
 }
 \]
 
 Puesto que $K_{\infty}/K$ es una extensi\'on abeliana e infinita, por
 el Teorema de clase de Hilbert (Teorema \ref{T9'.8.14}
y Corolario \ref{CClaseC4.5.13}),
 $K_{\infty}/K$ tiene necesariamente
 que ser ramificada. De hecho, si $\pK$ es un primo en $K$ 
 sobre $p$, $\pK$ es ramificado en $K_{\infty}/K$.
 
 Por otro lado, si $\pK$ es cualquier primo de $K$ ramificado en
$K_{\infty}/K$, entonces $\pK|_{{\ma Q}_n}$ es ramificado en
${\ma Q}_{\infty}$ de donde, por el Teorema \ref{T7.6.2}, se sigue
que $\pK|_{{\ma Q}_n}=p$. Esto es, todos los $p$--primos de $K$
son ramificados en $K_{\infty}/K$ y estos son los \'unicos.

Ahora, sea ${\eu l}$ primo de $K$ tal que ${\eu l}|_{\ma Q}=\ell
\neq p$. Entonces el grado de inercia de $\ell$ en $K_{\infty}/{\ma 
Q}_n$ es $\infty$ de donde se sigue que si ${\eu l}_{\infty}$
es un primo de 
$K_{\infty}$ sobre ${\eu l}$, entonces $D({\eu l}_{\infty}|{\eu l})\cong
p^m{\ma Z}_p$ para alguna $m\geq 0$. En otras palabras, ${\eu l}$
es totalmente descompuesto en $K_m/K$ y totalmente inerte en
$K_{\infty}/K_m$. Los primos infinitos son totalmente descompuestos
en $K_{\infty}/K$.
\[
 \xymatrix{
 &K_0=K\ar@{-}[r]\ar@{-}[d]&K_m\ar@{-}[r]\ar@{-}[d]&K_{\infty}\ar@{-}[d]\\
 &{\ma Q}_n\ar@{-}[dl]\ar@{-}[r]&{\ma Q}_{n+m}
 \ar@{-}[r]&{\ma Q}_{\infty}\\
 {\ma Q}
 }
 \]

Resumiendo tenemos:

\begin{teorema}\label{T7.6.3} Sea $K$ un campo num\'erico
finito y sea $K_{\infty}/K$ la extensi\'on ${\ma Z}_p$--ciclot\'omica
de $K$. Entonces los primos ramificados en $K_{\infty}/K$ son
precisamente los primos sobre $p$. Si ${\eu l}$ no es un $p$--primo,
entonces existe $m\in {\ma N}\cup \{0\}$ tal que ${\eu l}$ es
totalmente descompuesto en $K_m/K$ y totalmente inerte en
$K_{\infty}/K_m$. 

Finalmente los primos infinitos de $K$ son totalmente descompuestos
en $K_{\infty}/K$. $\fin$
\end{teorema}

A continuaci\'on enunciamos un teorema que corresponde a la 
teor{\'\i}a local de campos de clase.

\begin{teorema}\label{T7.6.4}
Sea ${\ma Q}_p$ el campo de los n\'umeros $p$--\'adicos. Sea
$K/k$ una extensi\'on abeliana finita donde $[k:{\ma Q}_p]<\infty$.
Entonces existe un mapeo, llamado el {\em mapeo de reciprocidad
de Artin\index{Artin!mapeo de reciprocidad de $\sim$}}
\begin{eqnarray*}
k^{\ast}&\longto&\Gal(K/k)\\
a&\longmapsto&(a,K/k)
\end{eqnarray*}
que induce un isomorfismo de grupos $k^{\ast}/N_{K/k}K^{\ast}\cong
\Gal(K/k)$, donde $N_{K/k}$ es el mapeo norma.

Sea $I$ el subgrupo de inercia de $\Gal(K/k)$. Entonces si $U_k$ y
$U_K$ son los grupos de unidades de $k$ y $K$ respectivamente,
se tiene $U_k/N_{K/k}U_K\cong I$.

Si $K/k$ es no ramificada, entonces $\Gal(K/k)$ es una extensi\'on
c{\'\i}clica generada por el automorfismo de Frobenius $F$, $\langle
F\rangle =\Gal(\overline{K}/\overline{k})$, $\overline{K}$ y $\overline{
k}$ los campos residuales y se tiene $(a,K/k)=F^{v_{\pi}(a)}$ donde
$v_{\pi}$ es la valuaci\'on de $k$, es decir, $\pi$ es un elemento 
primo, esto es, un elemento de valuaci\'on $1$.
\end{teorema}

\begin{proof}
Ver Teorema \ref{CClaseT3.2.2} (TCCL) y Teorema \ref{CCLT17.6.10}.
$\fin$
\end{proof}

Hasta ahora hemos considerado \'unicamente extensiones 
${\ma Z}_p$--ciclot\'omicas $K_{\infty}/K$, es decir, donde $K_{\infty}
=K{\ma Q}_{\infty}$. Recordemos que ${\ma Q}$ solo tiene una
\'unica extensi\'on ${\ma Z}_p$ y esta es la ciclot\'omica.

Cuando $K\neq {\ma Q}$, $K$ puede tener otras extensiones ${\ma Z
}_p$, de hecho, si $K$ no es totalmente real, lo cual quiere decir que
$r_1\neq [K:{\ma Q}]$, $K$ tiene m\'as 
extensiones ${\ma Z}_p$. Algunos de los resultados que hemos 
obtenido para extensiones ${\ma Z}_p$--ciclot\'omicas siguen siendo
ciertas para toda extensi\'on ${\ma Z}_p$ pero otros no.

\begin{proposicion}\label{P7.6.5} Sea $K_{\infty}/K$ una extensi\'on
${\ma Z}_p$ con $K$ un campo num\'erico finito. Entonces para
cada $n\geq 0$ hay un \'unico subcampo $K_n$ de grado $p^n$
sobre $K$ y $K_n$ y $K_{\infty}$ son los \'unicos subcampos entre
$K$ y $K_{\infty}$.
\end{proposicion}

\begin{proof}
Sea $S$ un subgrupo cerrado de ${\ma Z}_p$ y sea $x\in S$ tal que
$v_p(x)$ sea m{\'\i}nimo. Entonces $x{\ma Z}\subseteq S$ y como
$S$ es cerrado, se sigue que $x{\ma Z}_p\subseteq S$. Ahora para
$\xi\in S$, $v_p(\xi)\geq v_p(x)$, esto es, $\xi x^{-1}\in {\ma Z}_p$ y
por tanto $S\subseteq x{\ma Z}_p$, de donde se sigue que $S=x
{\ma Z}_p$. Ahora bien, claramente se tiene que $x{\ma Z}_p=0$
o $p^n{\ma Z}_p$ para alg\'un $n\geq 0$. Se sigue que los
subcampos de la extensi\'on $K_{\infty}/K$ son $K_{\infty}=K^{\{0\}}$
y $K_n=K_{\infty}^{p^n{\ma Z}_p}$, $K_0=K$. $\fin$
\end{proof}

\begin{teorema}\label{T7.6.6} Sea $K_{\infty}/K$ una extensi\'on 
${\ma Z}_p$ con $[K:{\ma Q}]<\infty$. Sea $\pK$ un divisor primo,
posiblemente infinito, de $K$ tal que  $\pK|_{\ma Q}\neq p$. Entonces
$\pK$ es no ramificado en $K_{\infty}/K$, esto es, $K_{\infty}/K$ es
no ramificada fuera de $p$.
\end{teorema}

\begin{proof}
Sea $I\subseteq \Gal(K_{\infty}/K)\cong{\ma Z}_p$ el grupo de
inercia de $\pK$. Puesto que $I$ es cerrado, entonces $I=\{0\}$ o
$I=p^u{\ma Z}_p$ para alg\'un $u\geq 0$. Si $I=\{0\}$, entonces
$\pK$ es no ramificado. Supongamos que $I=p^u{\ma Z}$, $u\geq
0$. En particular $I$ es infinito. Para $\pK$ un primo infinito,
tenemos $|I|=1$ o $2$, por lo que $\pK|_{\ma Q}\neq \infty$. Ahora
para cada $n$, sea $\pK_n$ un primo en $K_n$ tal que $\pK_n|_{
K_{n-1}}=\pK_{n-1}$ y $\pK_n|_K=\pK_n|_{K_0}=\pK$. Sea $\tilde{K
}_n$ la completaci\'on de $K_n$ en $\pK_n$ y $\tilde{K}_{\infty}=
\bigcup\limits_{n=1}^{\infty}\tilde{K}_n$. Entonces $I\subseteq \Gal(
\tilde{K}_{\infty}/\tilde{K})$. Sea $U$ el grupo de unidades de $
\tilde{K}$.

Si $I_n$ es el grupo de inercia de $\pK$ en $K_n/K$, entonces se
tiene que $I=\lim\limits_{\substack{\longleftarrow\\ n}}I_n$. Pero por
teor{\'\i}a de campos de clase (Teorema \ref{T7.6.4}) se tiene que
$U_n\stackrel{\varphi_n}{\longrightarrow}
I_n$ es suprayectiva. Para cualquiera $n\leq m$, se tiene
\[
\xymatrix{
U\ar[rr]^{\varphi_n}\ar[dr]_{\varphi_m}^{\hspace{.5cm}
\text{\Large{$\circlearrowleft$}}}&& I_n\\
&I_m\ar[ur]_{\pi_{m,n}}
}\qquad \pi_{m,n}\circ \varphi_m=\varphi_n,
\]
entonces existe un epimorfismo continuo $U\to I=\lim\limits_{\substack{
\longleftarrow\\ n}}I_n\cong p^u{\ma Z}_p$.

\begin{window}[0,l,\xymatrix{
{\cal O}_{\tilde{K}}\ar@{-}[r]\ar@{-}[d] & \tilde{K}\ar@{-}[d]^r\\
{\ma Z}_{\ell}\ar@{-}[r]&{\ma Q}_{\ell}},{}]
Por otro lado $U\subseteq {\cal O}_{\tilde{K}}$, $[\tilde{K}:{\ma Q}_{
\ell}]<\infty$ donde $\ell=\pK|_{\ma Q}\neq p$. Sabemos que 
$\tilde{K}^{\ast}
\cong \langle\pi\rangle\times U_{\tilde{K}}\cong \langle\pi
\rangle\times\mu_{q-1}\times U_{\tilde{K}}^{(1)}$ donde $\pi$ es un
elemento primo, $\big|{\cal O}_{\tilde{K}}/\tilde{\pK}\big|$, $U_{\tilde{K
}}^{(1)}=1+\tilde{\pK}$ y se tiene $
U_{\tilde{K}}\supseteq U_{\tilde{K}}^{(1)}\supseteq \ldots\supseteq
U_{\tilde{K}}^{(m)}\supseteq \ldots, \quad U_{\tilde{K}}^{(m)}= 1+
\tilde{\pK}^m$,
$U_{\tilde{K}}/U_{\tilde{K}}^{(1)}\cong\big( {\cal O}_{\tilde{K}}
/\tilde{\pK}
\big)^{\ast}=K(\tilde{\pK})^{\ast}$ (campo residual) y 
$U_{\tilde{K}}^{(m)}/U_{\tilde{K}}^{(m+1)}\cong K(\tilde{\pK})$.
\end{window}

Si $e=e(\pK|\ell)$ en $\tilde{K}/{\ma Q}_{\ell}$, se tiene para $n>
e/(\ell-1)$ los isomorfismos y homeomorfismos, uno inverso del otro
$\tilde{\pK}_K^n \biy\limits^{\exp}_{\log} U_K^{(n)}$.

As{\'\i} se tiene que ${\cal O}_{\tilde{K}}/\tilde{\pK}_K$ es finito y $
{\cal O}_{\tilde{K}} \cong {\ma Z}_{\ell}^r$ por lo que  $\tilde{\pK}_K^m
\cong {\ma Z}_{\ell}^r$. Por tanto $U_K^{(n)}\cong {\ma Z}_{\ell}^r$ y
$U_{\tilde{K}}/U_{\tilde{K}}^{(n)}$ finito implica que 
\[
U_{\tilde{K}}\cong(\text{grupo finito})\times {\ma Z}_{\ell}^r,
\quad r=[\tilde{K}:{\ma Q}],
\]
y el grupo finito es la torsi\'on de $U_{\tilde{K}}$.

Se tiene que $U_{\tilde{K}}\to p^n{\ma Z}_p\cong {\ma Z}_p$ el 
cual no tiene torsi\'on, por lo que existe un epimorfismo ${\ma Z}_{
\ell}^r\to p^n{\ma Z}_p$ y por tanto
\[
{\ma Z}_{\ell}^r\longto p^n{\ma Z}_p\stackrel{\pi}{\longto}
 p^n{\ma Z}_p/p^{n+1}
{\ma Z}_{p+1}\cong {\ma Z}/p{\ma Z}.
\]
Sin embargo ${\ma Z}_{\ell}^r$ no tiene subgrupos de {\'\i}ndice
$p$ pues $p\neq \ell$. Si sigue que $I=\{0\}$ y $\pK$ no es 
ramificado en $K_{\infty}/K$. $\fin$.
\end{proof}

\begin{proposicion}\label{P7.6.7} Sea $[K:{\ma Q}]<\infty$ y sea $
K_{\infty}/K$ una extensi\'on ${\ma Z}_p$. Existe al menos un primo
ramificado en $K_{\infty}/K$ y existe $n\geq 0$ tal que todo primo 
ramificado en $K_{\infty}/K$ es totalmente ramificado en $K_{\infty}/
K_n$.
\end{proposicion}

\begin{proof}
La m\'axima extensi\'on abeliana no ramificada corresponde, por
el Teorema de Clase de Hilbert (Teorema \ref{T9'.8.14}),
al grupo de ideales de $K$ y por
tanto es finita. Siendo $K_{\infty}/K$ abeliana e infinita,
necesariamente es ramificada.

Por el Teorema \ref{T7.6.6} s\'olo hay un n\'umero finito de ideales
primos ramificados, a saber, a lo m\'as los primos sobre $p$. Sean
$\pK_1,\ldots, \pK_r$ los primos ramificados, $r\geq 1$, y sean $I_1,
\ldots, I_r$ los grupos de inercia respectivos. Se tiene $I_j\neq \{0\}$,
$1\leq i\leq r$, por lo que $I_j\cong p^{n_j}{\ma Z}_p$ para alg\'un 
$n_j\geq 0$. En particular tenemos que $\bigcap\limits_{j=1}^r I_j=
p^m{\ma Z}_p$ donde $m=\max\{n_j\mid j=1,\ldots, r\}$. Por tanto,
puesto que $K_m=K_{\infty}^{p^m{\ma Z}_p}$, se tiene que $\Gal(
K_{\infty}/K_m)\subseteq I_j$ y cada $\pK_j$ es totalmente 
ramificado en $K_{\infty}/K_m$. $\fin$
\end{proof}

Podemos dar otra demostraci\'on de la unicidad de ${\ma Q}_{\infty}$.

\begin{corolario}\label{C7.6.8} ${\ma Q}_{\infty}/{\ma Q}$ es la
\'unica extensi\'on ${\ma Z}_p$ de ${\ma Q}$.
\end{corolario}

\begin{proof}
Sea $K_{\infty}/{\ma Q}$ una extensi\'on ${\ma Z}_p$ de ${\ma Q}$.
Tenemos que en $K_{\infty}/{\ma Q}$ \'unicamente $p$ se ramifica y
puesto que $K_1/{\ma Q}$ es ramificada, se tiene que $p$ es
totalmente ramificado en $K_{\infty}/{\ma Q}$.

Ahora bien, usando  caracteres de Dirichlet deducimos que
$K_n\subseteq \cic p{n_m}$ para alg\'un $m_n\in{\ma N}$ y adem\'as
$\infty$ no es ramificado, lo cual implica que $K_n\subseteq {\ma Q}
\big(\zeta_{p^{n_m}}+\zeta_{p^{n_m}}^{-1}\big)$ de donde se
sigue que $K_{\infty}=\bigcup\limits_{n=1}^{\infty}K_n\subseteq 
\cic p{\infty}^+$. Por lo tanto $K_n={\ma Q}_n$ y $K_{\infty}=
{\ma Q}_{\infty}$. $\fin$
\end{proof}

Se tiene que todo campo num\'erico finito $K$ tiene al menos una
extensi\'on ${\ma Z}_p$, a saber, $K_{\infty}:=K{\ma Q}_{\infty}$. 
Nos preguntamos ahora: dado un campo num\'erico $K$, $[K:{\ma Q}
]<\infty$, ?`cu\'antas extensiones ${\ma Z}_p$ tiene K?
Usamos la teor{\'\i}a de campos de clase para responder a esta
pregunta.

Sea $K$ dado, $[K:{\ma Q}]=d=r_1+2r_2<\infty$ donde $r_1$ denota
el n\'umero de encajes reales de $K\hookrightarrow {\ma R}$ y $r_2$
es el n\'umero de pares de encajes complejos no reales de $K$,
$K\hookrightarrow{\ma C}$.

Sea $U_K$ el grupo de unidades de $K$ (m\'as precisamente, de
${\cal O}_K$) y sea $E_1=\{u\in U_K\mid u\equiv 1\bmod \pK$ para
toda $\pK$ de $K$ tal que $\pK|_{\ma Q}= p\}$.

Sea $U_{K_{\pK}}=U_{\pK}$ el grupo de unidades del campo local 
$K_{\pK}$, $\pi$ un elemento primo y $U_{1,\pK}=\{x\in U_{\pK}\mid
x\equiv 1\bmod \pK\}$. Sea $q=N\pK$, esto es, $q=\big|
{\cal O}_{K_{\pK}}/\pK\big|$. Entonces ${\cal O}_{K_{\pK}}/\pK\cong
{\ma F}_q$. Sea $x\in U_{\pK}$, entonces se tiene
\[
x=a_0+a_1\pi+\cdots, \quad\text{con}\quad a_0\neq 0
\quad \text{y}\quad x^q\equiv a_0^q\bmod \pi\equiv a_0\bmod \pi
\equiv x\bmod \pi,
\]
es decir, $x^{N\pK-1}\equiv 1\bmod \pi$.

En particular, si $x\in U_K\subseteq {\cal O}_{K_{\pK}}$, entonces
$x^{N\pK-1}\in U_{1,\pK}$ y si $t:=\prod\limits_{\pK|p}(N\pK-1)\in
{\ma N}$, entonces $x^t\in E_1$, es decir, $U_K^t\subseteq E_1$.
En particular, $E_1$ es de {\'\i}ndice finito en $U_K$. Adem\'as, 
puesto que $U_K\cong {\ma Z}^{r_1+r_2-1}\times W(K)$ (Teorema
de las Unidades de Dirichlet) como grupos, donde $W(K)$ denota
a las ra{\'\i}ces de unidad en $K$, se tiene que el rango sobre ${\ma Z}
$ de $E_1$ es $r_1+r_2-1$.

M\'as precisamente, $E_1\cong {\ma Z}^{r_1+r_2-1}\times (\text{finito})
$ como grupos.

Por otro lado $U_{1,\pK}$ es un ${\ma Z}_p$--m\'odulo con acci\'on
$s\circ u:=u^s$ con $u\in U_{1,\pK}$, $s\in {\ma Z}_p$. Por lo anterior,
$U_{1,\pK}$ es de {\'\i}ndice finito en $U_{\pK}$ y por otro lado, para
$n$ suficientemente grande
\[
U_{n,\pK}\cong U_{\pK}^{(n)}\isomo\limits_{\substack{\to\\ \log}}^{\substack{
\exp\\ \leftarrow}} \pK^n{\cal O}_{K_{\pK}}.
\]

Se sigue, como grupos, $U_{n,\pK}\cong \pK^n\cong \pK\cong
{\cal O}_K$ 
(${\cal O}_{K_{\pK}}/\pK\cong {\ma F}_q$ es el campo residual) y
${\cal O}_{K_{\pK}}\cong {\ma Z}_p^{[K_{\pK}:{\ma Q}_p]}$, $
U_{n-1,\pK}/U_{n,\pK}\cong{\ma F}_q$.

De lo anterior obtenemos que 
\[
U_{1,\pK}\cong (\text{finito})\times \pK\cong (\text{finito})\times
{\ma Z}_p^{[K_{\pK}:{\ma Q}_p]}.
\]

En resumen tenemos
\begin{eqnarray*}
{\ma Z}^{r_1+r_2-1}\times(\text{finito})\cong E_1&\stackrel{\varphi}
{\longto}&\prod_{\pK|p}U_{1,\pK}\cong
(\text{finito})\times{\ma Z}_p^{\sum\limits_{\pK|p}[K_{\pK}:
{\ma Q}_p]}\\
&&\hspace{1.2cm}\cong (\text{finito})\times 
{\ma Z}_p^{[K:{\ma Q}]}\cong
(\text{finito})\times {\ma Z}_p^d\\
\varepsilon&\longmapsto& (\varepsilon,\ldots,\varepsilon)
\end{eqnarray*}

Usando la topolog{\'\i}a producto $\prod\limits_{\pK|p} U_{1,\pK}$ se
tiene que la cerradura $\overline{E}_1$ de $E_1$ (m\'as 
precisamente, de $\varphi(E_1)$) en $\prod\limits_{\pK|p}U_{1,\pK}$
es un ${\ma Z}_p$--m\'odulo, es decir, $\overline{E}_1\cong (
{\text{finito}})\times {\ma Z}_p^s$ para alg\'un $s$, esto es,
$\ran_{{\ma Z}_p} \overline{E}_1=s$. Puesto que
$\ran_{\ma Z} E_1=r_1+r_2-1$, se tiene $s\leq r_1+r_2-1$.

Existe la siguiente conjetura:

\begin{conjetura}[Leopoldt\index{Leopoldt!conjetura de
$\sim$}]\label{Co7.6.9}
Se tiene $\ran_{{\ma Z}_p}\overline{E}_1=s=r_1+r_2-1$.
\end{conjetura}

\begin{observacion}\label{O7.6.10} Se sabe que la conjetura de 
Leopoldt es cierta si la extensi\'on $K/{\ma Q}$ es abeliana.
\end{observacion}

Como consecuencia de la teor{\'\i}a de campos de clase, se tiene:

\begin{teorema}\label{T7.6.11} Sea $0\leq \delta \leq r_1+r_2-1$ tal
que $\ran_{{\ma Z}_p}\overline{E}_1=r_1+r_2-1-\delta$. Entonces
hay exactamente $r_2+1+\delta=d-\ran_{{\ma Z}_p}\overline{E}$
($d=[K:{\ma Q}]$) ${\ma Z}_p$--extensiones independientes de $K$.
En otras palabras, si $\tilde{K}$ denota la composici\'on de todas las
extensiones ${\ma Z}_p$ de $K$, se tiene que $\Gal(\tilde{K}/K)\cong
{\ma Z}_p^{r_2+1+\delta}$, ($r_2+1\leq r_2+1+\delta\leq r_1+2r_2=
d$).
\end{teorema}

Antes de presentar un esquema de demostraci\'on de este resultado,
recordemos un teorema de campos de clase. Sea $k$ un campo
num\'erico finito, $[k:{\ma Q}]<\infty$. Sea $\pK$ un primo finito o
infinito de $k$. Sea $k_{\pK}$ la completaci\'on de $k$ en $\pK$, 
y $U_{\pK}$ las unidades locales de $k_{\pK}$. Si $\pK$ es
arquimediano, es decir infinito, ponemos $U_{\pK}=k_{\pK}^{\ast}$ 
($={\ma C}^{\ast}$ o ${\ma R}^{\ast}$). Se define el 
{\em grupo de id\`eles\index{id\`eles}}
de $k$ (ver Subsecci\'on \ref{S17.6.3N}) por:
\[
J_k:=\big\{\big(\ldots,x_{\pK},\ldots\big)\in\prod_{\pK}k_{\pK}^{\ast}\mid
x_{\pK}\in U_{\pK}\text{\ para casi todo\ } \pK\big\}.
\]

A $J_k$ se le considera con la siguiente topolog{\'\i}a. Si $U=
\prod\limits_{\pK}U_{\pK}$, a $U$ se le dota de la topolog{\'\i}a
producto y definimos a $U$ como abierto en $J_k$. M\'as
precisamente, $J_k$ es un grupo topol\'ogico donde una base de
vecindades de la unidad $1$ ($=(\ldots,1,\ldots)$) est\'a dado por:
para $S\subseteq {\ma P}_k=\{\pK\mid \pK\text{\ es lugar de\ }k\}$,
$S$ finito y 
\[
\prod_{\pK\in S}V_{\pK}\times \prod_{\pK\notin S}U_{\pK}\subseteq
J_k,\quad V_{\pK}\subseteq k_{\pK}^{\ast}
\]
es una vecindad b\'asica de $1\in k_{\pK}^{\ast}$.

Se tiene que $\overline{V}_{\pK}$ es compacto y $\prod\limits_{\pK
\notin S}U_{\pK}$ es compacto, por lo que cada vecindad de $1$ en
$J_k$ contiene una vecindad
\[
\prod_{\pK\in S}V_{\pK}\times \prod_{\pK\notin S}U_{\pK}
\]
cuya cerradura es compacta. En otras palabras $J_k$, con la
topolog{\'\i}a definida as{\'\i}, es un grupo topol\'ogico localmente
compacto.

Sea 
\begin{eqnarray*}
k^{\ast}&\stackrel{\theta}{\longto}&J_k\\
x&\longmapsto&(\ldots, x,\ldots)=\big(x\big)_{\pK\in{\ma P}_K}
\end{eqnarray*}
Entonces $\theta$ es inyectiva, la imagen es un conjunto discreto.
$k^{\ast}$, o m\'as precisamente, $\theta(k^{\ast})$ (Proposici\'on
\ref{P7.6.12}) recibe el nombre de id\`eles principales.

\begin{proposicion}\label{P7.6.12} Se tiene que $k^{\ast}$ es
discreto en $J_k$. En particular $k^{\ast}$ es cerrado en $J_k$.
\end{proposicion}

\begin{proof}
Es suficiente probar que $1\in J_k$ tiene una vecindad que no 
contiene ning\'un otro id\`ele principal. Sea $S$ el conjunto de primos
arquimedianos y sea $V=\{\alpha\in J_k\mid |\alpha_{\pK}-1|<1
\text{\ para $\pK\in S$ y $|\alpha_{\pK}|=1$, para\ } \pK\notin S\}$.
Entonces si $B_{\pK}(1,1)$ es la bola abierta de radio $1$ y centro
$1$ en $k_{\pK}$, se tiene $V=\prod\limits_{\pK\in S}B_{\pK}(1,1)
\times \prod\limits_{\pK\notin S}U_{\pK}$.

Si hubiese un id\`ele principal $x\in V$ con $x\neq 1$, se tendr{\'\i}a:
\begin{align*}
1&=\prod_{\pK}|x-1|_{\pK}=\prod_{\pK\in S}|x-1|_{\pK}\cdot
\prod_{\pK\notin S}|x-1|_{\pK}<\prod_{\pK\notin S}|x-1|_{\pK}\\
&\leq \prod_{\pK\in S}\max \{|x|_{\pK}, 1\}=1\quad \text{pues}\quad
x\in k^{\ast}, x\in V.
\end{align*}
lo cual es absurdo. $\fin$
\end{proof}

\begin{definicion}\label{D7.6.13} Se define el {\em grupo de clases
de id\`eles\index{id\`eles!grupo de $\sim$}} por $C_k:=J_k/k^{\ast}$.
\end{definicion}

Sea $L/k$ una extensi\'on finita. Si $\pL$ es un primo de $L$ y $
\pL|_k=\pK$, sea $N_{\pL/\pK}\colon L_{\pL}\to k_{\pK}$ la norma de
campos locales. Sea $x=(\ldots, x_{\pL},\ldots)\in J_L$. Se define la
norma de $J_L$ a $J_k$ por $N_{L/k}\colon J_L\to J_k$, $N_{L/k}=
(\ldots,y_{\pK}, \ldots)\in J_k$ donde $y_{\pK}=\prod\limits_{\pL|\pK}
N_{\pL/\pK} x_{\pL}$.

Si $x$ es principal, esto es, $x``="(\ldots,x,\ldots)$, entonces
$y_{\pK}=\prod\limits_{\pL|\pK}N_{\pL/\pK}x_{\pL}=\prod\limits_{\pL|
\pK}N_{\pL/\pK}x=N_{L/k} x$ para toda $\pK$ y $N_{L/k}$ es la
norma usual de $L$ en $k$. Por tanto $N_{L/k}(x)$ es principal. Por
tanto la norma $J_L\to J_k$ induce una norma de $C_L$ en $C_k$:
$N_{L/k}\colon C_L\to C_k$.

A continuaci\'on presentamos los teoremas fundamentales de la
teor{\'\i}a de campos de clase sin demostraci\'on (ver los Teoremas
\ref{T9'.8.2'}, \ref{T9'.8.5}, \ref{CClaseT3.2.29}, \ref{CClaseT4.2.1} y
\ref{CClaseTC.1}, la Proposici\'on \ref{P9'.8.3} 
y la Subsecci\'on \ref{CClaseS4.7}).

\begin{teorema}\label{T7.6.14} Si $L/k$ es una extensi\'on abeliana
finita, entonces existe un isomorfismo
\[
J_k/\big(k^{\ast} N_{L/k}J_L\big)\cong C_k/\big(N_{L/k}C_L\big)
\cong\Gal(L/k).
\]
Adem\'as, el primo finito o infinito $\pK$ es no ramificado en $L/k$
si y s\'olo si $U_{\pK}\subseteq k^{\ast}N_{L/k}J_L$ donde $U_{\pK}
\hookrightarrow J_k$ est\'a dado por $u\mapsto (\ldots,1,\ldots,1,u,1,
\ldots,1,\ldots)$. $\fin$
\end{teorema}

\begin{teorema}\label{T7.6.15}
Si $H$ es un subgrupo abierto de $C_k$ de {\'\i}ndice finito, entonces
existe una \'unica extensi\'on abeliana $L/k$ tal que $N_{L/k}C_L=H$.
Equivalentemente, si $H$ es un abierto de {\'\i}ndice finito en $J_k$
y $k^{\ast}\subseteq H$, entonces existe una \'unica extensi\'on
abeliana $L/k$ tal que  $k^{\ast}N_{L/k}J_L=H$. $\fin$
\end{teorema}

\begin{teorema}\label{T7.6.16} Si $L_1$ y $L_2$ son dos extensiones
abelianas de $k$, entonces $L_1\subseteq L_2\iff k^{\ast}N_{L_1/k}
J_{L_1}\supseteq k^{\ast}N_{L_2/k}J_{L_2}$. $\fin$
\end{teorema}

Los Teoremas \ref{T7.6.14}, \ref{T7.6.15} y \ref{T7.6.16} se pueden
enunciar para las extensiones infinitas de campos num\'ericos. Sea
$D_k$ la componente conexa de la identidad del grupo de clases de
id\`eles $C_k$.

\begin{teorema}\label{T7.6.17}
\lasa
\item Si $L/k$ es una extensi\'on abeliana entonces existe un
subgrupo cerrado $H$ de $C_k$ con $D_k\subseteq H\subseteq
C_k$ tal que $C_k/H\cong \Gal(L/k)$.

Adem\'as el primo $\pK$ es no ramificado en $L/k \iff k^{\ast}U_{\pK}
/k^{\ast}\subseteq H$.

\item Dado un subgrupo cerrado $H$ de $C_k$ tal que $D_k
\subseteq H\subseteq C_k$ o, equivalentemente, $C_k/H$ es
totalmente disconexo, existe una \'unica extensi\'on abeliana $L/k$
tal que $C_k/H\cong \Gal(L/k)$. $\fin$
\end{list}
\end{teorema}

\begin{ejemplo}\label{Ej7.6.18}
Sea $k$ un campo num\'erico y sea $L$ el campo de clase de
Hilbert de $k$, es decir, $L$ es la m\'axima extensi\'on abeliana no
ramificada en todo primo finito o infinito. Puesto que $L/k$ es
no ramificada en todas partes, se tiene $U:=\prod\limits_{\pK\in
{\ma P}_k}U_{\pK}\subseteq k^{\ast}N_{L/k}J_L$. Puesto que
$L/k$ es maximal con respecto a esta propiedad, $k^{\ast} U$ es el
grupo que corresponde a $L$ de donde 
\[
J_k/k^{\ast}U\cong \Gal(L/k).
\]

Sea $D_k$ el grupo de ideales de $k$ y sea $\varphi\colon J_k\to
D_k$ dada por: $\varphi\big((x)_{\pK}\big)=\prod\limits_{\pK \text{\ 
finito}} \pK^{v_{\pK}(x_{\pK})}$. Se tiene que $\varphi$ est\'a bien
definida pues $v_{\pK}(x_{\pK})=0$ para casi toda $\pK$. Se tiene
que $\ker \varphi=U$ y finalmente, ya que $\varphi^{-1}(
\text{ideales principales})=\varphi^{-1}(P_k)$ son los id\`eles principales, 
se sigue que $J_k/U k^{\ast}\cong D_k/P_k=I_k=$ grupo de
clases de $k$. En otras palabras
\[
\Gal(L/k)\cong I_k= \text{\ grupo de clases de ideales de\ }k.
\]
\end{ejemplo}

Ahora regresamos al Teorema \ref{T7.6.11}.

\medskip

\noindent{\bf Teorema \ref{T7.6.11}.} {\em
Si $[K:{\ma Q}]<\infty$ y $\tilde{K}$ es la composici\'on de todas las
${\ma Z}_p$--extensiones de $K$, entonces $\Gal(\tilde{K}/K)\cong
{\ma Z}_p^{r_2+1+\delta}$ donde $\ran_{{\ma Z}_p}\overline{E}_1
=r_1+r_2-1-\delta$, $\delta \geq 0$.}

\begin{proof}
Sea $F$ la m\'axima extensi\'on abeliana de $K$ que es no 
ramificada fuera de $p$. Entonces $\tilde{K}\subseteq F$. Sea
$J=K_K$. Del Teorema \ref{T7.6.17} obtenemos la existencia de un
grupo cerrado $H$ tal que $K^{\ast}\subseteq H\subseteq J$ tal que
$J/H\cong \Gal(F/K)$.

Sea $U_{\tilde{\ell}}$ el grupo de unidades locales de primo $\tilde{\ell}
$ de $K$ si $\tilde{\ell}$ es finito y ponemos $U_{\tilde{\ell}}=
K^{\ast}_{\tilde{\ell}}$ si $\tilde{\ell}$ es arquimediano.

Sean $U':=\prod\limits_{\pK|p}U_{\pK}$, $U'':=\prod\limits_{\tilde{\ell}
\nmid p}U_{\tilde{\ell}}$ y $U:=U'\times U''$. Poniendo $1$ en el resto
de las componentes, $U'$, $U''$ son subgrupos de $J$ y $U$ es abierto. Puesto
que $F/K$ es no ramificado fuera $p$, se tiene que $U''\subseteq H$.

Puesto que $F$ es maximal y $H$ es cerrado, entonces $H=
\overline{K^{\ast}U''}$. Sean $J':=J/H$ y $J'':=K^{\ast}U/H=K^{\ast}
U'U''/H=U'H/H\cong U'/(U'\cap H)$. Ahora sea $U_1:=\prod\limits_{
\pK|p}U_{\pK,1}$, donde $U_{\pK,1}$ es el grupo de las unidades
locales congruentes con $1$ m\'odulo $\pK$. Como vimos 
anteriormente, se tiene que $U_{\pK}/U_{\pK,1}\cong {\ma F}_q^{\ast}
$, ${\ma F}_q={\cal O}_{K_{\pK}}/\pK$ es el campo residual. Por lo
tanto $U'/U_1$ es finito. De hecho tenemos $U_{\pK}\cong
\mu_{q-1}\times U_{\pK,1}$ pues $U_{\pK,1}$ contiene a lo m\'as
$p^s$ ra{\'\i}ces de $1$, $q=p^n$, $p=\car{\ma F}_q$. Por lo tanto
$U'=\prod\limits_{\pK|p}U_{\pK}\cong (\text{finito})\times \prod\limits_{
\pK|p}U_{\pK,1}=(\text{finito})\times U_1=T\times U_1$ con $T$
finito y de orden primo relativo a $p$. As{\'\i} $U_1\cong(p\text{--grupo
finito})\times {\ma Z}_p^{[K:{\ma Q}]}$. Por otro lado si $U_{\pK,n}
:=\{x\in U_{\pK}\mid x\equiv 1\bmod \pK^n\}$, se tiene $U_{\pK,n}/
U_{\pK, n+1}\cong {\ma F}_q$.

Por lo tanto $\big|U_{\pK,1}/U_{\pK, n}\big|=q^{n-1}$ y para $n$
suficientemente grande, $U_{\pK,n}\cong \pK^n\cong {\cal O}_{K_{
\pK}}\cong {\ma Z}_p^{[K_{\pK}:{\ma Q}]}$.
Se sigue que $U'=\prod\limits_{\pK|p}U_{\pK}=T\times \prod\limits_{
\pK|p}U_{\pK,1}=T\times U_1$, $\mcd (|T|,p)=1$. As{\'\i} se tiene que
si $C=J''=U'/(U'\cap H)$ entonces $C=(\text{finito})\times \frac{U_1
(U'\cap H)}{U'\cap H}$ pues $U'\cong T\times U_1$.

Si definimos $A_1:=TM/M$ donde $M:=U\cap H$ y $B_1:=U_1M/M$,
entonces $A_1,B_1\subseteq C$ y puesto que $U'=TU_1$,
se sigue que $A_1+B_1=C$. Por otro lado $A_1$ es un cociente
de $T$ y por tanto finito y de orden primo relativo a $p$ y $B$ es un
cociente de $U_1$ de donde tenemos $B_1\cong (p\text{--grupo 
finito})\times {\ma Z}_p^{\beta}$.

En particular, $A_1\cap B_1=\{0\}$ por lo que $C\cong A_1\oplus
B_1$ y por tanto $C\cong (\text{finito})\times U_1M/M$. As{\'\i},
tenemos
\[
J''=(\text{finito})\times \frac{U_1(U'\cap H)}{U'\cap H}\cong
(\text{finito})\times \frac{U_1}{U_1\cap U\cap H}=(\text{finito})\times
\frac{U_1}{U_1\cap H}.
\]

Sea $\psi\colon E_1\to U_1\subseteq J$ el mapeo diagonal en $U_1$
y con $1$ en las otras entradas, es decir,
\[
\psi(\varepsilon)=(\ldots,1,\ldots,\varepsilon,\ldots,\ldots)\quad\text{con
$\varepsilon$ si $\pK|p$ y $1$ si $\pK\nmid p$}.
\]
En otras palabras $\psi(\varepsilon)=\big(x_{\pK}\big)_{\pK}$, $x_{\pK}
=\begin{cases} \varepsilon&\text{si $\pK|p$}\\ 1&\text{si $\pK\nmid p$}
\end{cases}$.

Antes de continuar con la demostraci\'on, probemos:

\begin{lema}\label{L7.6.19}
Se tiene $U_1\cap H=U_1\cap \overline{K^{\ast} U''}=\overline{\psi(
E_1)}$.
\end{lema}

\begin{proof}[Lema {\rm \ref{L7.6.19}}] Sea $\varepsilon\in E_1$, $\psi(
\epsilon)\in U_1$ y puesto que $\frac{\psi(\varepsilon)}{\varepsilon}$
tiene componente $1$ en todas las entradas tales que $\pK\nmid p$,
se tiene $\psi(\varepsilon)=\varepsilon\cdot\frac{\psi(\varepsilon)}{
\varepsilon}\in K^{\ast}U''$.
Por tanto $\overline{\psi(E_1)}\subseteq U_1\cap \overline{K^{\ast}
U''}$.

Rec{\'\i}procamente, sea $U_{\pK,n}=1+\pK^n=\{x\in U_{\pK}\mid
x\equiv 1\bmod \pK^n\}$ y sea $U_n:=\prod\limits_{\pK|p}U_{\pK,n}$.
Poniendo $1$ en todas las componentes tales que $\pK\nmid p$ se
tiene que 
\[
K^{\ast}U''U_n=\bigcup_{x\in K^{\ast}}\big(x\cdot U''\cdot U_n\big)
\]
es abierto y $\overline{K^{\ast}U''}=\bigcap\limits_{n=1}^{\infty}
K^{\ast}U''U_n$ pues si $x\in K^{\ast}U''$, entonces $x\big(U''\times
U_n\big)\cap K^{\ast}U''\neq \emptyset$, $x y=\xi t$ para $y\in
U''\times U_n$, $\xi\in K^{\ast}$, $t\in U''$ por lo que $x=y^{-1}\xi t=
\xi t y^{-1}\in K^{\ast}U'U_n$ de donde $\overline{K^{\ast}U''}
\subseteq\bigcap\limits_{n=1}^{\infty}K^{\ast}U''U_n$ y 
rec{\'\i}procamente si $x\notin \overline{K^{\ast}U''}$ existe una
vecindad $V$ de $1$ tal que $xV\cap K^{\ast}U''=\emptyset$.
Entonces existe $n\geq 1$ tal que $V\supseteq U_n\times 
\prod\limits_{\pK\in S}W_{\pK}\times \prod\limits_{\substack{\pK\notin
S\\ \pK\nmid p}}U_{\pK}$ donde $S$ es un conjunto finito de lugares
$\pK$ donde $\pK\nmid p$. 

Si $x\in K^{\ast}U''U_n$, $x=\xi ty$ con $\xi \in K^{\ast}$, $t\in U''$,
$y\in U_n$. Por tanto $xy^{-1}=\xi t\in xV\cap K^{\ast}U''$ lo cual
es absurdo. As{\'\i} pues, $x\notin K^{\ast}U''U_n$ y por tanto $
\overline{K^{\ast}U''}\supseteq \bigcap\limits_{n=1}^{\infty} K^{\ast}
U''U_n$ de donde obtenemos la igualdad anunciada: $\overline{K^{
\ast}U''}=\bigcap\limits_{n=1}^{\infty}K^{\ast}U''U_n$.

Similarmente tenemos $\overline{\psi(E_1)}=\bigcap\limits_{n=1}^{
\infty}\psi(E_1)U_n$, y $U_n$ es compacto. Para verificar que
$U_1\cap \overline{K^{\ast}U''}\subseteq \overline{\psi(E_1)}$, basta
probar que $U_1\cap K^{\ast}U''U_n\subseteq \psi(E_1)U_n$ para 
toda $n$.

Sea $x\in K^{\ast}$, $u''\in U''$, $x\in U_n\subseteq U_1$, $xu''u\in
U_1$ por lo que $xu''\in U$. Ahora bien, $u''$ tiene componente $1$
en todas las entradas tales que $\pK\nmid p$ y $x$ debe ser una
unidad principal en las componentes $\pK|p$.

Ahora $U_1$ tiene componentes $1$ para $\pK\nmid p$  y $u''$ es
una unidad ah{\'\i}, por lo que $x$ debe ser una unidad en esas
entradas. En resumen, $x$ es unidad local para toda $\pK$ por lo
que $x$ debe ser una unidad global adem\'as de que $x\equiv 1\bmod
\pK$ lo cual implica que $x\in E_1$.

As{\'\i} tenemos $xu''=x\in E_1$ para $\pK|p$. Para $\pK\nmid p$,
$xu''=1$ por tanto $xu''\in \psi(E_1)$ y por lo tanto $xu''u\in\psi(E_1)
U_n$. Esto termina la demostraci\'on del lema. $\fin$
\end{proof}

Regresando a la demostraci\'on del Teorema \ref{T7.6.11}, 
recordemos que $U_1\cong (\text{finito})\times {\ma Z}_p^{[K:{\ma Q}]}
$ por lo que $U_1/(U_1\cap H)= U_1/\big(\overline{\psi(E_1)}\big)
\cong (\text{finito})\times {\ma Z}_p^{\alpha}$ donde 
\[
\alpha=[K:{\ma Q}]-\ran_{{\ma Z}_p}\overline{\psi(E_1)}=r_1+2r_2-
(r_1+r_2-1-\delta)=r_2+1+\delta.
\]
Por tanto $J''\cong (\text{finito})\times \frac{U_1}{(U_1\cap H)}\cong
(\text{finito})\times {\ma Z}_p^{r_2+1+\delta}$, ($J''\cong \frac{
K^{\ast}U}{H}$).

Ten{\'\i}amos que $J'=J/H$, $H=\overline{K^{\ast}U''}$. Ahora
\[
J'/J''\cong J/K^{\ast}U=I_K
\]
donde $I_K$ es el grupo de clases de ideales
de $K$ el cual es finito.
Adem\'as, $J'/{\ma Z}_p^{r_2+1+\delta}\cong (\text{finito})=T_1$, 
digamos $|T_1|=m$. 

Se tiene que $m{\ma Z}_p^{r_2+1+\delta}\subseteq mJ'\subseteq
{\ma Z}_p^{r_2+1+\delta}$ por lo que $mJ'\cong {\ma Z}_p^{r_2+1+
\delta}$ como ${\ma Z}_p$--m\'odulos.

Sea $J'_m=\{x\in J'\mid mx=0\}$. Entonces $J'_m$ es cerrado en $J'$
y $J'/J'_m\cong m J'\cong {\ma Z}_p^{r_2+1+\delta}$.

Ahora bien $|J'_m|\leq m$ pues en caso contrario, esto es, si $|J'_m|
>m$, y puesto que $J'/{\ma Z}_p^{r_2+1+\delta}\cong T_1$, entonces
tendr{\'\i}amos dos elementos $x,y\in J'_m$, $x\neq y$, $x\bmod 
{\ma Z}_p^{r_2+1+\delta}=y\bmod {\ma Z}_p^{r_2+1+\delta}$ y $
m(x-y)=0$ lo cual no es posible pues $x-y\in {\ma Z}_p^{r_2+1+\delta}
$ es no cero y ${\ma Z}_p^{r_2+1+\delta}$ es libre de torsi\'on.

As{\'\i}, $|J'_m|\leq m<\infty$. Finalmente el campo fijo de $J'_m
\subseteq J'\cong\Gal(F/K)$ debe ser $\tilde{K}$ y por lo tanto
\[
\Gal(\tilde{K}/K)\cong J'/J'_m\cong {\ma Z}_p^{r_2+1+\delta}.
\tag*{$\fin$}
\]
\end{proof}

\begin{corolario}\label{C7.6.20} Sea $K_H$ el campo de clase de
Hilbert de un campo num\'erico finito $K$ y sea $F$ la m\'axima
extensi\'on abeliana de $K$ no ramificada fuera de $p$. Entonces
\[
\Gal(F/K_H)\cong\frac{\prod_{\pK|p}U_{\pK}}{\overline{E}}
\]
donde $\overline{E}$ es la cerradura del grupo de unidades $E$ de 
$K$ cuando en encajado diagonalmente en $\prod\limits_{\pK|p}
U_{\pK}$ y $U_{\pK}$ es el grupo de unidades locales de $K_{\pK}$.
\end{corolario}

\begin{proof}
Se tiene que $\Gal(F/K)\cong J'=J/H$ y el subgrupo cerrado $J''=
KU/H$ corresponde a $K_H$. Por lo tanto $\Gal(F/K_H)\cong J''
\cong U'/(U'\cap H)$. Como antes, se tiene que $U'\cap H=
\overline{\psi(E)}$ y el resultado se sigue. $\fin$
\end{proof}

\begin{corolario}\label{C7.6.21} ${\ma Q}$ s\'olo tiene una
extensi\'on ${\ma Z}_p$.
\end{corolario}

\begin{proof}
$[{\ma Q}:{\ma Q}]=1=r_1+2r_2$, $r_1=1$, $r_2=0$, $0\leq
\ran_{{\ma Z}_p} \overline{E}_1=r_1+r_2-1-\delta\leq r_1+r_2-1=
1+0-1=0$. Por lo tanto $\delta=0$ y $r_2+1+\delta=0+1+0=1$. 

De hecho, $E_1=\{x\in U_{\ma Q}\mid x\equiv 1\bmod p\}=\{1\}$, por
lo que $\overline{E}_1=\{1\}$ y por lo tanto $\ran_{{\ma Z}_p} E_1
=0$.

Una tercera demostraci\'on es simplemente que 
\[
1=0+1=r_2+1\leq \text{\ n\'umero de extensiones ${\ma Z}_p$ de
${\ma Q}$\ }\leq d=[{\ma Q}:{\ma Q}]=1. \tag*{$\fin$}
\]
\end{proof}

\begin{observacion}\label{O7.6.22} Se tiene que ${\ma Z}_p\times
{\ma Z}_p$ no puede ser el grupo de Galois de ninguna extensi\'on
$K/{\ma Q}$. 

El problema inverso de la Teor{\'\i}a de Galois pregunta si dado
$G$ un grupo finito, existe una extensi\'on de Galois $K$ de ${\ma Q}
$ tal que $\Gal(K/{\ma Q})\cong G$. Ya vimos que si $G$ es abeliano
la respuesta es si (Teorema \ref{T12.4.1.1}; ver tambi\'en Teorema
\ref{T12.4.1.2}).

Sin embargo, si $G$ es infinito vemos que la respuesta al problema
inverso a la Teor{\'\i}a de Galois es en general no, aunque el grupo
$G$ sea abeliano, por ejemplo $G={\ma Z}_p\times {\ma Z}_p$.
\end{observacion}

\section{Estructura de $\Lambda={\ma Z}_p
[[T]]$--m\'odulos}\label{S7.7}

En esta secci\'on consideraremos una ${\ma Z}_p$--extensi\'on
$K_{\infty}/K$ y ponemos $\Gal(K_{\infty}/K)\cong \Gamma =
{\ma Z}_p$, $\Gamma$ escrito multiplicativamente.

Por un lado tenemos que varios m\'odulos relacionados con la
extensi\'on $K_{\infty}/K$ tienen una acci\'on dada por $\Gamma$,
considerado $\Gamma$ como el grupo de Galois y, por otro lado,
tienen una acci\'on del grupo aditivo ${\ma Z}_p$, no considerado
como grupo de Galois. Juntando ambas acciones tendremos varios
${\ma Z}_p[\Gamma]$--m\'odulos. Esta es la raz\'on de denotar
$\Gal(K_{\infty}/K)$ por $\Gamma$ en lugar de ${\ma Z}_p$.

Sea $\gamma$ un generador topol\'ogico de $\Gamma$, es decir,
$\overline{\langle\gamma\rangle}=\Gamma$. Por ejemplo, bajo el
isomorfismo $\Gamma\stackrel{\cong}{\longto}{\ma Z}_p$,
$\gamma$ puede corresponder a $1\in{\ma Z}_p$.

Se tiene que $\Gamma^{p^n}\cong p^n{\ma Z}_p$. Sea $\Gamma_n
:=\Gamma/\Gamma^{p^n}\cong C_{p^n}$, el grupo c{\'\i}clico de $p^n
$ elementos. Entonces $\Gamma_n\cong\Gal(K_n/K)=\langle
\gamma\bmod \Gamma^{p^n}\rangle$.

Sea ${\cal O}$ el anillo de enteros de una extensi\'on finita de ${\ma Q
}_p$ y sea $\pK$ el ideal maximal de ${\cal O}$ y $\pi$ un elemento
primo de $\pK$, es decir, $\pK=\langle\pi\rangle$.

Consideremos el anillo grupo ${\cal O}[\Gamma_n]$. Si $m\geq n
\geq 0$, sea $\phi_{m,n}\colon {\cal O}[\Gamma_m]\to{\cal O}[
\Gamma_n]$ el mapeo inducido por el mapeo natural $\psi_{m,n}
\colon \Gamma_m\to\Gamma_n$. 

Se tiene que ${\cal O}[\Gamma_n]\cong\frac{{\cal O}[T]}
{\big((1+T)^{p^n}
-1\big)}$ con isomorfismo
\[
\gamma\bmod \Gamma^{p^n}\stackrel{\xi}{\longto}(1+T)\bmod\big(
(1+T)^{p^n}-1\big).
\]
Se tiene el siguiente diagrama conmutativo para $m\geq n\geq 0$:
\[
\xymatrix{
{\cal O}[\Gamma_m]\ar[r]^{\hspace{-1.2cm}
\xi}\ar[d]^{\phi_{m,n}}&{\cal O}[T]/\big(
(1+T)^{p^m}-1\big)\ar[d]_{\phi'_{m,n}}\\
{\cal O}[\Gamma_n]\ar[r]_{\hspace{-1.2cm}\xi}&{\cal O}[T]/\big(
(1+T)^{p^n}-1\big)
}
\]
donde $\phi'_{m,n}$ es el mapeo natural pues $\big(
(1+T)^{p^n}-1\big)|\big((1+T)^{p^m}-1\big)$.

Sea $\lim\limits_{\substack{\longleftarrow\\ n}}{\cal O}[\Gamma_n]
\cong {\cal O}[[\Gamma]]$ el anillo de grupo profinito. En general
se tiene que ${\cal O}[\Gamma]\subseteq {\cal O}[[\Gamma]]$ pues
si $\alpha\in{\cal O}[\Gamma]$, entonces $\alpha=\sum\limits_{
\sigma\in \Gamma}a_{\sigma}\sigma$ con $a_{\sigma}=0$ para
casi toda $\sigma\in \Gamma$. Sea $\alpha_n:=\sum\limits_{\sigma
in \Gamma} a_{\sigma}\phi_n(\sigma)$ donde $\phi_n\colon \Gamma
\to \Gamma_n$ es el mapeo asociado con $\lim\limits_{\substack{
\longleftarrow\\ n}} \Gamma_n=\Gamma$. Entonces $\alpha\in
{\cal O}[\Gamma_n]$ y claramente $\phi_{m,n}(\alpha_m)=\alpha_n$.
Sin embargo ${\cal O}[\Gamma]\subsetneqq {\cal O}[[\Gamma]]$
como veremos m\'as adelante.

Se tiene 
\[
{\cal O}[[\Gamma]]=\lim_{\substack{\longleftarrow\\ n}}{\cal O}[
\Gamma_n]\cong \lim_{\substack{\longleftarrow\\ n}} {\cal O}[T]/
\big((1+T^{p^n}-1\big).
\]
Veremos que ${\cal O}[[\Gamma]]\cong{\cal O}[[T]]$.

\begin{proposicion}[Algoritmo euclidiano\index{algoritmo euclidiano}]
\label{P7.7.1}
Sean  $f,g\in {\cal O}[[T]]$ y $f=a_0+a_1T+\cdots $ con $a_i\in\pK$,
$0\leq i\leq n-1$ y $a_n\in{\cal O}^{\ast}$. Entonces se tiene una
\'unica expresi\'on del tipo $g=qf+r$ con $q\in{\cal O}[[T]]$, $r\in
{\cal O}[T]$, $r$ un polinomio de grado menor o igual a $n-1$.
\end{proposicion}

\begin{proof}
Primero probaremos la unicidad de la expresi\'on. Basta probar que 
si $qf+r=0$ entonces $q=r=0$. Supongamos que $q$ o $r$ no es
cero. Podemos suponer que $\pi\nmid r$ o $\pi\nmid q$. Puesto que
el grado de $r$ es menor o igual a $n-1$ y $\pi|a_i$, $0\leq i\leq n-1$,
tomando m\'odulo $\pi$, se tiene que $\pi|r$ y por tanto $\pi|qf$. 
Puesto que $a_n\in{\cal O}^{\ast}$, $\pi\nmid f$ de donde se sigue
que $\pi|q$ lo cual contradice nuestra elecci\'on. Por lo tanto $q=r=0$
lo cual prueba la unicidad.

Para demostrar la existencia, sea $\tau=\tau_n\colon {\cal O}[[T]]\to
{\cal O}[[T]]$ dado por $\tau\big(\sum\limits_{i=0}^{\infty} b_iT^i\big)=
\sum\limits_{i=n}^{\infty}b_iT^{i-n}=b_n+b_{n+1}T+\cdots $.

Se tiene que $\tau$ es ${\cal O}$--lineal y satisface:
\lasa
\item $\tau\big(T^nh(T\big)=h(T)$ para toda $h(T)\in{\cal O}[[T]]$.
\item $\tau(h(T))=0\iff h(T)\in{\cal O}[T]$, con $\gr h(T)\leq n-1$.
\end{list}

Por hip\'otesis podemos escribir 
\[
f(T)=\pi P(T)+T^n U(T)
\]
con $P(T)\in{\cal O}[T]$, $\gr P\leq n-1$ y $U(T)=a_n+a_{n+1}T+
\cdots =\tau(f(T))\in{\cal O}[[T]]$.

Para tener $g=qf+r$, $\gr r\leq n-1$ es necesario y suficiente que
$\tau(g)=\tau(qf)$. Puesto que $qf=\pi P q +T^n qU$ necesitamos
resolver
\[
\tau(g)=\tau(\pi q f)+\tau(T^n q U)=\tau(\pi q P)+q U=\tau(\pi q P)
+q\tau(f).
\]

Ahora bien, $\tau(f)=U$ es invertible. Sea $z=qU=q\tau(f)$. La 
ecuaci\'on anterior equivale a
\[
\tau(g)=\tau\Big(\frac{\pi q}{U}P U\Big)+z=\tau\Big(z\cdot\frac
{\pi P}{\tau(f)}\Big)+z=\Big(I+\tau\cdot \frac{\pi P}{\tau (f)}\Big)(z).
\]

Notemos que $\tau\cdot \frac{\pi P}{\tau(f)}\colon{\cal O}[[T]]\to \pK
{\cal O}[[T]]$ pues $\frac{\pi P}{\tau(f)}\in\pK{\cal O}[[T]]$ por lo tanto
se puede invertir $\big(I+\tau\cdot\frac{\pi P}{\tau(f)}\big)$:
\begin{gather*}
\begin{align*}
U q&=z=\Big(I+\tau\cdot\frac{\pi P}{\tau(f)}\Big)^{-1}\tau(g)=
\sum_{j=0}^{
\infty}\Big(\tau\cdot\frac{\pi P}{\tau(f)}\Big)^j\tau(g)\\
&=
\sum_{j=0}^{\infty}(-1)^j\pi^j\Big(\tau\cdot\frac{P}{U}\Big)^j\cdot \tau(g),
\end{align*}
\intertext{por lo tanto}
q=\frac{1}{U(T)}\sum_{j=0}^{\infty}(-1)^j\pi^j\Big(\tau\cdot \frac{P}{U}
\Big)^j\cdot\tau(g). \tag*{$\fin$}
\end{gather*}
\end{proof}

\begin{definicion}\label{D7.7.1} Un polinomio $p(T)\in{\cal O}[T]$
se llama {\em distinguido\index{polinomio distinguido}} si $p(T)=
T^n+a_{n-1}T^{n-1}+\cdots+a_1T+a_0$ con $a_i\in\pK$ para
$0\leq i<n-1$. Notemos que $p(T)$ ``casi'' es de Eisenstein con
la salvedad de que posiblemente $\pi^2$ divide a $a_0$.
\end{definicion}

\begin{teorema}[Teorema $p$--\'adico de Preparaci\'on de
Weierstrass\index{Weierstrass!teorema de preparaci\'on
de $\sim$}\index{teorema de preparaci\'on de Weierstrass}]
\label{T7.7.3}
Sea $f(T)=\sum_{i=0}^{\infty}a_iT^i\in{\cal O}[[T]]$ tal que para 
alguna $n\geq 0$, se tiene que $a_i\in \pK$, $0\leq i\leq n-1$ pero
$a_n\notin \pK$, esto es, $a_n\in{\cal O}^{\ast}$. Entonces $f$ puede
ser escrito de manera \'unica en la forma 
\[
f(T)=p(T)U(T)
\]
con $U(T)\in{\cal O}[[T]]$ es unidad y $p(T)$ es un polinomio 
distinguido de grado $n$.

M\'as generalmente, si $f(T)\in {\cal O}[[T]]\setminus\{0\}$, $f$ puede
ser escrito de manera \'unica como $f(T)=\pi^uP(T)U(T)$ como
antes y $u\geq 0$.
\end{teorema}

\begin{proof}[Manin]
Notemos que la segunda parte es inmediata de la primera donde 
$\pi^u$ es la m\'axima potencia de $\pi$ que divide a todos los
coeficientes de $f$.

Sea $g=T^n$. Entonces, por el algoritmo euclidiano, se tiene que
$T^n=q(T)f(T)+r(T)$, $\gr r\leq n-1$.

Puesto que $q(T)f(T)\equiv q(T)\big(a_nT^n+T^{n+1}h(T)\big)\bmod
\pi$ por lo que $r(T)\equiv 0\bmod \pi$. Por lo tanto $p(T)=T^n-r(T)$
es un polinomio distinguido de grado $n$. Sea $q_0:=q(0)$.
Comparando coeficientes de $T^n$, se tiene que $1=q_0 a_n
\bmod \pi$ lo cual implica que $q_0\in{\cal O}^{\ast}$ as{\'\i} que 
$q(T)$ es una unidad.

Sea $U(T)=\frac{1}{q(T)}$ lo que nos da $f(T)=P(T)U(T)$.

Para la unicidad, puesto que todo polinomio distinguido de grado $n$
puede ser escrito como $p(T)=T^n-r(T)$, la ecuaci\'on $f(T)=P(T)U(T)
$ puede ser escrito de nuevo como $T^n=U(T)^{-1}f(T)+r(T)$ y
por la unicidad del algoritmo euclidiano, se sigue la unicidad de $
f(T)=P(T)U(T)$. $\fin$
\end{proof}

\begin{corolario}\label{C7.7.4} Sea $\overline{\ma Q}_p$ una 
cerradura algebraica de ${\ma Q}_p$ y sea ${\ma C}_p$ la 
completaci\'on de $\overline{\ma Q}_p$. Se tiene que como campos,
${\ma C}\cong {\ma C}_p$ m\'as no topol\'ogicamente. En particular
${\ma C}_p$ es algebraicamente cerrado. Entonces, existen
un n\'umero finito de elementos $x\in {\ma C}_p$ tales que $|x|<1$ y
$f(x)=0$.
\end{corolario}

\begin{proof}
Sea $f(x)=0$ y escribamos $f(T)=\pi^u P(T)U(T)$ como en el Teorema
de Preparaci\'on de Weierstrass. Puesto que $U(x)\neq 0$, es
invertible, $U(x)$ converge para $|x|<1$. Por lo tanto $P(x)=0$.
De ah{\'\i} se sigue el resultado. $\fin$
\end{proof}

\begin{lema}\label{L7.7.5} Supongamos $P(T)\in{\cal O}[T]$ un
polinomio distinguido y sea $g(T)\in{\cal O}[T]$ arbitrario. Si $g(T)/
P(T)\in {\cal O}[[T]]$, entonces $g(T)/p(T)\in{\cal O}[T]$.
\end{lema}

\begin{proof}
Supongamos $g(T)=f(T)P(T)$ con $f(T)\in{\cal O}[[T]]$. Sea $x\in
{\ma C}_p$ un cero de $P(T)$. Por tanto $0=P(x)=x^n+\pi \alpha$
lo cual implica que $|x|<1$ y por lo tanto $f(x)$ converge. Se sigue
que $g(x)=0$. Dividiendo entre $T-x$ y ampliando el anillo ${\cal O}$
en caso necesario, obtenemos que $P(T)|g(T)$ como polinomios,
por lo tanto en ${\cal O}[T]$. $\fin$
\end{proof}

\begin{teorema}[J.-P. Serre]\label{T7.7.6} Se tiene ${\cal O}[[\Gamma
]]\cong {\cal O}[[T]]$, donde el isomorfismo 
es inducido por $\gamma\colon \longmapsto 1+T$.
\end{teorema}

\begin{proof}
Es suficiente probar que ${\cal O}[[T]]\cong \lim\limits_{\substack{
\longleftarrow\\ n}}\frac{{\cal O}[T]}{\big((1+T)^{p^n}-1\big)}$.

Notemos que $p_n(T)=(1+T)^{p^n}-1$ es un polinomio distinguido.
De hecho, $\langle \pi, T\rangle \supseteq \langle p,T\rangle$ es 
un ideal maximal de ${\cal O}[T]$ y tambi\'en da el ideal maximal de
${\cal O}[[T]]$. Ahora $P_0(T)=(1+T)^{p^0}-1=T\in\langle p,T\rangle$.
Adem\'as
\begin{align*}
\frac{P_{n+1}(T)}{P_n(T)}&=\frac{(1+T)^{p^{n+1}}-1}{(1+T)^{p^{n}}-1}
\igual_{\substack{\uparrow\\ y=(1+T)^{p^{n}}}}\frac{y^p-1}{y-1}=
y^{p-1}+y^{p-2}+\cdots+y+1\\
&= (1+T)^{p^{n}(p-1)}+(1+T)^{p^{n}(p-2)}+\cdots+(1+T)^{p^{n}}+1
\in \langle p,T\rangle.
\end{align*}

Por inducci\'on en $n$ se sigue que $P_n(T)\in\langle p,T\rangle^{
n+1}$. Por el algoritmo euclidiano, para cada $f(T)\in{\cal O}[[T]]$,
$f(T)=q_n(T)P_n(T)+f_n(T)$, $\gr f_n(T)<p^n$.
Entonces $f(T)\to f_n(T)$ es un mapeo natural: ${\cal O}[[T]]\to
\frac{{\cal O}[T]}{\langle P_n(T)\rangle}$ y si $m\geq n\geq 0$ 
entonces
\[
f_m(T)-f_n(T)=p_n(T)\Big(q_n(T)-q_m(T)\frac{P_m(T)}{P_n(T)}\Big).
\]
Por lo tanto $P_n(T)|f_m(T)-f_n(T)$ como polinomios (Lema 
\ref{L7.7.5}).

Por lo tanto $f_m(T)\equiv f_n(T)\bmod P_n(T)$. Se sigue que 
$(f_0,f_1,\ldots, )\in \lim\limits_{\substack{\longleftarrow\\ n}}
\frac{{\cal O}[T]}{\langle P_n(T)\rangle}$.
As{\'\i} tenemos un mapeo:
\begin{eqnarray*}
{\cal O}[[T]]&\stackrel{\varphi}{\longto}&\lim\limits_{\substack{
\longleftarrow\\ n}}\frac{{\cal O}[T]}{\langle P_n(T)\rangle}\\
f&\longmapsto& (f_0,f_1,\ldots)
\end{eqnarray*}

Ahora si $f_n=0$ para toda $n$, $P_n|f$ para toda $n$. Por tanto
$f\in\bigcap\limits_{n=0}^{\infty}\langle p, T\rangle^{n+1}=(0)$ de
donde se sigue que $\varphi$ es $1$--$1$.

Sea ahora $(f_0,f_1,\ldots)\in\lim\limits_{\substack{\longleftarrow\\ n}}
\frac{{\cal O}[T]}{\langle P_n(T)\rangle}$. Por lo tanto, para $m\geq
n\geq 0$, $f_m\equiv f_n\bmod P_n$, es decir, $f_m\equiv f_n
\bmod \langle p,T\rangle^{n+1}$. Se sigue que los t\'erminos
constantes son congruentes m\'odulo $p^{n+1}$, los t\'erminos
lineales son congruentes m\'odulo $p^n$, etc.

As{\'\i}, los coeficientes $f_m$ forman una sucesi\'on de Cauchy.
Alternativamente $f=\lim\limits_{n\to\infty} f_n$ existe pues ${\cal O}
[[T]]$ es completo en la topolog{\'\i}a $\langle p,T\rangle$--\'adica.
Sea $f(T):=\lim\limits_{n\to \infty} f_n(T)\in{\cal O}[[T]]$. Veremos que
$\varphi(f)=(f_0,f_1,\ldots)$.

Si $m\geq n\geq 0$, $f_m-f_n=q_{m,n}P_n$ para alg\'un $q_{m.n}\in
{\cal O}[T]$. Entonces $q_{m,n}=\frac{f_m-f_n}{P_n}\xrightarrow[m\to
\infty]{}\frac{f-f_n}{P_n}$. Puesto que $q_{m,n}\in{\cal O}[T]$, el
l{\'\i}mite debe estar en ${\cal O}[[T]]$, es decir, sin denominadores.
Por lo tanto $f=(P_n)\big(\lim\limits_{m\to\infty}q_{m,n}\big)+f_n$, 
esto es, $\varphi(f)=(f_0,f_1,\ldots)$ y $\varphi$ es suprayectiva. $\fin$
\end{proof}

\begin{lema}\label{L7.7.7}
Se tiene que $\Lambda={\ma Z}_p[[T]]$ es un dominio de
factorizaci\'on \'unica cuyos irreducibles son $p$ y los polinomios
distinguidos irreducibles.

Las unidades de $\Lambda^{\ast}$ de $\Lambda$ son las series de potencias
con termino constante en ${\ma Z}_p^{\ast}$.
\end{lema}

\begin{proof}
Dado $f(T)\in{\ma Z}_p[[T]]$, se tiene que por el Teorema de
Preparaci\'on de Weierstrass, $f(T)$ se puede escribir de manera
\'unica como $f(T)=p^uP(T)U(T)$ con  $u\in{\ma Z}$, $u\geq 0$, 
$p(T)$ un polinomio distinguido y $U(T)\in\Lambda^{\ast}$. Por
el Lema \ref{L7.7.5} se tiene que si $f(T)\in{\ma Z}_p[T]$, 
entonces $U(T)\in{\ma Z}_p[T]$.

Ahora, en ${\ma Z}_p[T]$, $P(T)$ es producto \'unico de polinomios
m\'onicos e irreducibles: $P(T)=P_1(T)\cdots P_r(T)$. Ahora
$P(T)\bmod p=T^n$. Por lo tanto $P_i(T)\bmod p=T^{n_i}$, $1\leq
i\leq r$, $n=n_1+\cdots+n_r$ y cada $P_i(T)$ es un polinomio
distinguido. Por lo tanto ${\ma Z}_p[[T]]$ es un dominio de
factorizaci\'on \'unica (D.F.U.) y sus irreducibles son $p$ y los
polinomios distinguidos irreducibles. $\fin$
\end{proof}

\begin{observacion}\label{O7.7.8} Se tiene que $\Lambda={\ma Z}_p
[[T]]$ no es un dominio de ideales principales. En particular no
tenemos los teoremas estructurales de los dominios de ideales
principales en esta clase de dominios. Sin
embargo, $\Lambda={\ma Z}_p[[T]]$ es dominio de factorizaci\'on
\'unica y noetheriano pues ${\ma Z}_p$ es noetheriano.
\end{observacion}

\begin{lema}\label{L7.7.9} Sean $f,g\in \Lambda$ primos relativos.
Entonces el anillo $\Lambda/\langle f,g\rangle$ es finito.
\end{lema}

\begin{proof} Sea $h\in\langle f,g\rangle$ un elemento de grado
m{\'\i}nimo. Entonces $h=p^sH$ con $H=1$ o $H$ distinguido. 
Supongamos que $H\neq 1$. Puesto que $f$ y $g$ son primos
relativos, podemos suponer que $H\nmid f$. Sin embargo
\[
f=qH+r,\quad \gr r<\gr H=\gr h\quad\text{por lo que}\quad
p^sf=qh+p^sr
\]
y puesto $\gr(p^s r)=\gr r<\gr h$, $p^s r\in\langle f,g\rangle$ y
$p^r s\neq 0$ pues $H\nmid f$, esto es absurdo. Por lo tanto
$H=1$ y $h=p^s$. Sin p\'erdida de generalidad, podemos suponer
que $p\nmid f$ y que $f$ es distinguido 
pues en otro caso se puede usar $g$ o dividir por
una unidad ya que $f$ y $g$
son primos relativos. Se tiene que $\langle f,g\rangle \supseteq\langle
p^s,f\rangle$. 

Por el algoritmo euclidiano, dado $\alpha(T)\in\Lambda$, $\alpha
\bmod f\equiv r$ con $\gr r<\gr f$ y puesto que s\'olo hay un n\'umero
finito de tales polinomios m\'odulo $p^s$, entonces
\[
\big|\Lambda/\langle p^s,f\rangle\big|<\infty \quad\text{y}\quad
\big|\Lambda/\langle f,g\rangle\big|\leq \big|\Lambda/\langle p^s,f
\rangle\big|<\infty. \tag*{$\fin$}
\]
\end{proof}

\begin{lema}[Casi ``Teorema Chino del Residuo'']\label{L7.7.10}
Sean $f,g\in\Lambda$ primos relativos. Entonces
\l
\item el mapeo natural
\begin{eqnarray*}
\Lambda/\langle f,g\rangle&\stackrel{\tilde{\varphi}}{\longto}&
\Lambda/\langle f\rangle\oplus \Lambda/\langle g\rangle\\
h\bmod \langle f,g\rangle&\stackrel{\tilde{\varphi}}{\longmapsto}&
(h\bmod \langle f\rangle,h\bmod \langle g\rangle)
\end{eqnarray*}
es inyectivo y tiene con\'ucleo finito.

\item Existe un mapeo inyectivo $\Lambda/\langle f\rangle\oplus
\Lambda/\langle g\rangle\stackrel{\psi}{\longto}\Lambda/\langle 
f,g\rangle$ que tiene con\'ucleo finito.
\end{list}
\end{lema}

\begin{proof}
\l
\item Por ser $\Lambda$ un dominio de factorizaci\'on \'unica, si $h
\in\ker \varphi$, donde $\varphi\colon \Lambda\to 
\Lambda/\langle f\rangle\oplus \Lambda/\langle g\rangle$, se tiene
que $h\bmod \langle f\rangle=0$ y $h\bmod\langle g\rangle=0$,
de donde obtenemos que $f|h$ y $g|h$ y por lo tanto $fg|h$ lo cual
implica que $\tilde{\varphi}\colon \Lambda/\langle f,g\rangle
\longto\Lambda/\langle f\rangle\oplus \Lambda/\langle g\rangle$
es inyectivo.

Sea ahora $(a\bmod \langle f\rangle,b\bmod \langle g\rangle)\in
\Lambda/\langle f\rangle\oplus \Lambda/\langle g\rangle$. Si
$a-b\in\langle f,g\rangle$ entonces $a-b=fA+gB$ para algunos
$A,B\in \Lambda$. Sea $c:=a-fA=b+gB$. Entonces $c\equiv a\bmod
\langle f\rangle$ y $c\equiv b\bmod \langle f\rangle$ y en este caso
$(a\bmod \langle f\rangle,b\bmod \langle g\rangle)$ est\'a en la imagen
de $\tilde{\varphi}$.

Puesto que $\big|\Lambda/\langle f,g\rangle=n<\infty$, podemos 
seleccionar $r_1,\ldots,r_n\in\Lambda$ un conjunto de representantes
de $\Lambda/\langle f,g\rangle$. En particular, tendremos que
$\{(0\bmod \langle f\rangle,r_j\bmod \langle g\rangle\mid 1\leq i\leq n
\}$ es un conjunto de representantes del con\'ucleo de $\tilde{\varphi}$
(ejercicio) y por lo tanto el con\'ucleo es finito.

\item Sea $\Lambda/\langle f, g\rangle\cong M\subseteq
\Lambda/\langle f\rangle\oplus \Lambda/\langle g\rangle=:N$. Se tiene
que $|N/M|<\infty$.

Sea $P$ cualquier polinomio distinguido en $\Lambda$ que sea
primo relativo a $f\cdot g$. Para ver su existencia, 
notemos que, por el criterio de Eiseinstein, $T^n+p$
es distinguido e irreducible para toda $n$.
Puesto que hay una infinidad de ellos, existe uno que es primo 
relativo a $fg$.

Si $(x,y)\in N$, entonces $p^i(x,y)\cong p^j(x,y)\bmod M$ para 
algunos $i<j$. Puesto que $1-p^{j-i}\in\Lambda^{\ast}$, entonces
\begin{align*}
p^i(x,y)&=(1-p^{j-i})^{-1}(1-p^{j-i})p^i(x,y)=\\
&=(1-p^{j-i})^{-1}(p^i(x,y)-p^j(x,y))\in M.
\end{align*}

Por lo tanto $p^kN\subseteq M$ para alg\'un $k$ (recordemos que
$|N/M|<\infty$). Supongamos $p^k(x,y)=0$ en $N$. Entonces
$f|p^k x$, $g|p^k y$. Puesto que $\mcd (p,fg)=1$, entonces
$f|x$ y $g|y$ lo cual implica que $(x,y)=0$ en $N$. Se sigue que
si $\cdot p^k$ denota multiplicaci\'on por $p^k$, entonces
\[
N\stackrel{\cdot p^k}{\longto}M\stackrel{\sim}{\longto}
\Lambda/\langle f, g\rangle\quad \text{es inyectivo}.
\]

La imagen contiene al ideal $\langle p^k,fg\rangle$ y puesto que
$p^k$ y $fg$ son primos relativos, se sigue que $\Lambda/\langle p^k
fg\rangle$ es finito y por lo tanto $\coker (\cdot p^k)$ es finito. $\fin$
\end{list}
\end{proof}

\begin{proposicion}\label{P7.7.11} Los ideales primos de $\Lambda$
son $(0)$, $\langle p,T\rangle$, $\langle p\rangle$ y $\langle P(T)
\rangle$ donde $P(T)$ es un polinomio distinguido irreducible. M\'as
a\'un, $\Lambda$ es un anillo local con ideal m\'aximo $\langle p,T
\rangle$.
\end{proposicion}

\begin{proof} Se tienen los siguientes isomorfismos
\begin{align*}
\Lambda/(0)&\cong \Lambda,\\
\Lambda/\langle p\rangle&\cong {\ma F}_p[[T]],\\
\Lambda/\langle p,T\rangle&\cong {\ma F}_p,\\
\Lambda/\langle P(T)\rangle&\cong {\ma Z}_p[T]/\langle P(T)\rangle,
\end{align*}
este \'ultimo isomorfismo se puede obtener
\ usando el algoritmo euclidiano. Puesto que
todos los cocientes anteriores son dominios enteros, todos los
ideales son primos y adem\'as contenidos en $\langle p,T\rangle$.

Sea ahora $\pK$ un ideal primo no cero de $\Lambda$. Sea $h\in
\pK$ de grado m{\'\i}nimo. Sea $h=p^s H$ con $H=1$ o $H$ 
distinguido. Puesto que $\pK$ es primo se sigue que $p\in\pK$ o
$H\in \pK$. En el caso en que $1\neq H\in\pK$, entonces, puesto que
$h$ es de grado m{\'\i}nimo, $H$ necesariamente es irreducible.
En cualquier de estos dos casos, $\langle f\rangle \subseteq \pK$
donde $f=p$ o $f$ es irreducible y distinguido. Si $\langle f\rangle
=\pK$ el resultado se sigue. Si $\langle f\rangle \neq \pK$, 
consideremos $g\in\pK$, $g\notin\langle f\rangle$, esto es,
$f\nmid g$.

Puesto que  $f$ es irreducible, $f$ y $g$ son primos relativos y
$\langle f,g\rangle\subseteq \pK$ implica que $|\Lambda/\pK|\leq
|\Lambda/\langle f,g\rangle|<\infty$. En particular $\Lambda/\pK$
es un ${\ma Z}_p$--m\'odulo finito y $p^{n_0}\in \pK$ para alguna
$n_0\in{\ma N}$ y como $\pK$ es primo, $p\in\pK$.

Por otro lado, tenemos que $T^i\equiv T^j\bmod \pK$ para algunos
$i,j\in{\ma N}$, $i<j$ y ya que $1-T^{j-i}\in\Lambda^{\ast}$, entonces
\[
T^i=(1-T^{j-i})^{-1}(T^i(1-T^{j-i}))=(1-T^{j-i})^{-1}(T^i-T^j)\in\pK
\]
lo que a su vez implica que $T\in \pK$ y por lo tanto $\langle p,T
\rangle\subseteq \pK$. Finalmente, puesto que $\Lambda/\langle
p,T\rangle \cong {\ma F}_p$ se sigue que $\langle p,T\rangle$ es
maximal y por lo tanto $\langle p,T\rangle=\pK$. $\fin$
\end{proof}

En contraste a la finitud de $\Lambda/\langle f,g\rangle$ (Lema
\ref{L7.7.9}) se tiene:

\begin{proposicion}\label{P7.7.12} Sea $f\in\Lambda$, $f\notin
\Lambda^{\ast}$. Entonces $\Lambda/\langle f\rangle$ es 
infinito.
\end{proposicion}

\begin{proof}
Supongamos $f\neq 0$ pues de otra manera el resultado es
inmediato. Por el Teorema de Preparaci\'on de Weierstrass, basta
suponer que $f=p$ o $f=P(T)$ con $P(T)$ un polinomio distinguido.
Si $f=p$, $\Lambda/\langle p\rangle\cong {\ma F}_p[[T]]$ y si $f=
P(T)$ es un polinomio distinguido, $\Lambda/\langle f\rangle=
\Lambda/\langle P(T)\rangle \cong {\ma Z}_p[T]/\langle P(T)\rangle$
el cual es isomorfo, como grupo, a ${\ma Z}_p^{\gr_TP(T)}$.
El resultado se sigue. $\fin$
\end{proof}

\begin{definicion}\label{D7.7.13} Dos $\Lambda$--m\'odulos $M$ y
$M'$ se
llaman {\em seudo--isomorfos\index{m\'odulos!seudo--isomorfos}},
lo cual denotamos por $M\sim M'$, si existe un
$\Lambda$--homomorfismo $\varphi\colon M\to M'$ tal que tanto
$\ker \varphi$ como $\coker\varphi$ son conjuntos finitos.

Equivalentemente, $M$ y $M'$ son seudo isomorfos si existe una
sucesi\'on exacta de $\Lambda$--m\'odulos 
\[
0\longto A\longto M\longto M'\longto B\longto 0
\]
con $A$ y $B$ finitos.
\end{definicion}

\begin{observacion}\label{O7.7.14} Si $M\sim M'$, {\underline{no}}
necesariamente $M'\sim M$. En particular, $\sim$ no es una
relaci\'on de equivalencia.
\end{observacion}

\begin{ejemplo}\label{Ej7.7.15}
Se tiene que $\langle p,T\rangle \sim \Lambda$ pues la inyecci\'on
natural $\langle p,T\rangle \to \Lambda$ tiene con\'ucleo ${\ma F}_p$
el cual es finito.

Por otro lado, si $\varphi\colon\Lambda\to \langle p,T\rangle$ es 
cualquier $\Lambda$--homomorfismo, sea $\varphi(1)=f(T)$.
Entonces se tiene que $\coker \varphi=\langle p,T\rangle/\langle
f\rangle$. Puesto que $|\Lambda/\langle f\rangle|=\infty$ y 
$|\Lambda/\langle p,T\rangle|<\infty$, necesariamente $|\langle p,T
\rangle/\langle f\rangle|=\infty$.
\end{ejemplo}

Sin embargo, lo que si se tiene es:

\begin{teorema}\label{T7.7.16} Si $M$ y $M'$ son
$\Lambda$--m\'odulos de torsi\'on finitamente generados, entonces
$M\sim M'\iff M'\sim M$.
\end{teorema}

\begin{proof}
Ver \cite[p\'agina 272]{Was97}. $\fin$
\end{proof}

\begin{observacion}\label{O7.7.17} Si $f,g\in\Lambda$ son primos relativos,
entonces por el casi Teorema Chino del Residuo, $\Lambda/\langle fg\rangle
\sim \Lambda/\langle f\rangle\oplus \Lambda/\langle f\rangle$ y $\Lambda/
\langle f\rangle\oplus \Lambda/\langle g\rangle \sim 
\Lambda/\langle fg\rangle$.
\end{observacion}

\begin{teorema}[Estructura de $\Lambda$--m\'odulos finitamente generados]
\label{T7.7.18}
Sea $M$ un $\Lambda$--m\'odulo finitamente generado. Entonces
\[
M\sim \Lambda^r\oplus\Big(\bigoplus_{i=1}^s\Lambda/\langle p^{n_i}\rangle
\Big)\oplus\Big(\bigoplus_{j=1}^t\Lambda/\langle f_j(T)^{m_j}\rangle\Big)
\]
donde $r,s,t,n_i,m_j\in{\ma N}\cup \{0\}$ y cada $f_j(T)$ es un polinomio
irreducible distinguido.
\end{teorema}

\begin{proof}
El enunciado es similar al de la estructura de un m\'odulo sobre un dominio
de ideales principales, pero en nuestro caso tenemos \'unicamente un
seudo--isomorfismo en lugar de un isomorfismo.

Supongamos que $M$ est\'a generado por $u_1,\ldots,u_n$ con relaciones
$\lambda_1u_1+\cdots+\lambda_n u_n=0$, $\lambda_i\in \Lambda$.

Sea $R$ el $\Lambda$--subm\'odulo de $\Lambda^n$ formado por las
relaciones. Puesto que $\Lambda$ es un anillo noetheriano, $R$ es 
$\Lambda$--finitamente generado. As{\'\i} $M$ puede ser representado
por una matriz cuyas filas son de la forma $(\lambda_1,\ldots, \lambda_n)$
y donde $\sum\limits_{i=1}^n\lambda_i u_i=0$ es una relaci\'on. Se tiene la
sucesi\'on exacta $0\longto R\longto \Lambda^n\longto M\longto 0$ de
$\Lambda$--m\'odulos.
Por abuso del lenguaje, esta matriz tambi\'en ser\'a llamada $R$.

Las siguientes son las operaciones b\'asicas de filas y columnas que
corresponden a cambiar generadores de $R$ y $M$.

\noindent
{\underline{Operaci\'on A}}: Permuta de filas o de columnas.

\noindent
{\underline{Operaci\'on B}}: Adicionamos un m\'ultiplo de una fila o una
columna a otra fila o columna respectivamente. Como caso especial, si
$\lambda'=q\lambda+r$,
\[
\left(
\begin{array}{cccc}
\vdots&&\vdots\\ \lambda&\cdots&\lambda'&\cdots\\ \vdots&&\vdots
\end{array}\right) \longmapsto
\left(
\begin{array}{cccc}
\vdots&&\vdots\\ \lambda&\cdots&r&\cdots\\ \vdots&&\vdots
\end{array}\right)
\quad r=\lambda'-q\lambda.
\]

\noindent
{\underline{Operaci\'on C}}: Podemos multiplicar una fila o una columna por
un elemento de $\Lambda^{\ast}$.

Notemos que estas tres operaciones son las mismas que las usadas
para los dominios de ideales principales. Tendremos tres operaciones
adicionales que es donde intervienen los seudos--isomorfismos en 
lugar de los isomorfismos.

\noindent
{\underline{Operaci\'on 1}}: Si $R$ contiene una fila $(\lambda_1,
p\lambda_2,\ldots,p\lambda_n)$ con $p\nmid \lambda_1$ entonces
cambiamos $R$ por la matriz $R'$ cuya primer fila es $(\lambda_1,
\lambda_2,\ldots, \lambda_n)$ y las dem\'as filas de $R'$ son las
filas de $R$ con el primer elemento multiplicado por $p$:
\[
\left(
\begin{array}{cccc}
\lambda_1&p\lambda_2&\cdots&\cdots\\
\alpha_1&\alpha_2&\cdots&\cdots\\
\beta_1&\beta_2&\cdots&\cdots\\
\vdots&\vdots&\vdots&\vdots
\end{array}
\right)
\longto
\left(
\begin{array}{cccc}
\lambda_1&\lambda_2&\cdots&\cdots\\
p\alpha_1&\alpha_2&\cdots&\cdots\\
p\beta_1&\beta_2&\cdots&\cdots\\
\vdots&\vdots&\vdots&\vdots
\end{array}
\right)
\]

Como caso especial, si $\lambda_2=\cdots=\lambda_n=0$, entonces
podemos multiplicar $\alpha_1,\beta_1,\ldots$ por una potencia 
arbitraria de $p$.

Afirmamos que el m\'odulo $M'$ es seudo--isomorfo a $M$. En 
efecto, en $R$ tenemos la relaci\'on 
\[
\lambda_1u_1+p(\lambda_2 u_2+\cdots+\lambda_n u_n)=0
\]
Sea $\tilde{M}'=M\oplus \Lambda v$ con un nuevo generador $v$ y
sean las relaciones adicionales $(-u_1,pv)=0$, $(\lambda_2 u_2+
\cdots +\lambda_n u_n,\lambda_1 v)=0$. Sea $M'$ igual a $\tilde{M}'$
m\'odulo estas nuevas relaciones.

Sea $M\stackrel{\varphi}{\longto}M'$ el mapeo natural. Si $m\in\ker
\varphi$, $m$ est\'a en el m\'odulo de relaciones as{\'\i} que:
\[
(m,0)=a(-u_1,pv)+b(\lambda_2u_2+\cdots+\lambda_n u_n,\lambda_1
v), \quad a,b\in \Lambda
\]
de donde $ap=-b\lambda_1$. Puesto que $p\nmid \lambda_1$,
$p|b$ y $\lambda_1|a$. En la primera componente obtenemos
\begin{align*}
m&=-au_1+b(\lambda_2 u_2+\cdots+\lambda_nu_n)\\
&=\frac{-a}{\lambda_1}(\lambda_1 u)-\frac{ap}{\lambda_1}(\lambda_2
u_2+\cdots+\lambda_nu_n)\\
&=-\frac{a}{\lambda_1}(\lambda_1u_1+p\lambda_2 u_2+\cdots+
p\lambda_nu_n)=-\frac{a}{\lambda_1}(0)=0
\end{align*}
lo cual prueba que $\varphi$ es inyectiva.

Ahora bien las im\'agenes de $pv ``=" u_1$ y $\lambda_1v ``="
-(\lambda_2u_2+\cdots+\lambda_nu_n)$ en $M'$ est\'an en la
imagen de $M$ y el ideal $\langle p,\lambda_1\rangle$
aniquila a $M'/M$. Puesto que $p\nmid \lambda_1$, $\mcd (p,
\lambda_1)=1$, $\big|\Lambda/\langle p,\lambda_1\rangle\big|<
\infty$ y $M'$ es finitamente generado, se sigue que $|M'/M|<
\infty$ y por tanto $M\sim M'$.

El nuevo m\'odulo tiene generadores $v,u_2,\ldots, u_n$, ``$u_1=
pv$''. Cualquier relaci\'on $\alpha_1u_1+\cdots+\alpha_nu_n=0$
viene a ser $p\alpha_1v+\alpha_2u_2+\cdots+\alpha_n u_n=0$.
As{\'\i} que la primer columna es multiplicada por $p$. Finalmente
tambi\'en tenemos la relaci\'on $\lambda_1v+\lambda_2u_2+\cdots
+\lambda_nu_n=0$ as{\'\i} que la nueva matriz $R'$ tiene la forma
dada pues removemos la fila redundante $(p\lambda_1,p\lambda_2,
\ldots,p\lambda_n)$.

Continuamos con la demostraci\'on del teorema.

\noindent
{\underline{Operaci\'on 2}}: Si todos los elementos de la primer
columna de $R$ son divisibles por $p^k$ y si existe una fila $(
p^k\lambda_1,\ldots,p^k\lambda_n)$ con $p\nmid \lambda_1$,
podemos cambiar $R$ a una matriz $R'$ que es la misma que 
$R$ excepto que $(p^k\lambda_1,\ldots,p^k\lambda_n)$ es
reemplazado por $(\lambda_1,\ldots,\lambda_n)$. Es decir
\[
\left(
\begin{array}{cccc}
p^k\lambda_1&p^k\lambda_2&\cdots&\cdots\\
p^k\alpha_1&\alpha_2&\cdots&\cdots\\
p^k\beta_1&\beta_2&\cdots&\cdots\\
\vdots&\vdots&\vdots&\vdots
\end{array}
\right)
\longto
\left(
\begin{array}{cccc}
\lambda_1&\lambda_2&\cdots&\cdots\\
p^k\alpha_1&\alpha_2&\cdots&\cdots\\
p^k\beta_1&\beta_2&\cdots&\cdots\\
\vdots&\vdots&\vdots&\vdots
\end{array}
\right).
\]
Al hacer esto obtenemos un m\'odulo $M'$ tal que $M'=M''\oplus S$,
donde $M''$ est\'a dado por la matriz $R'$ y $S\cong \Lambda/\langle
p^k\rangle$ para alg\'un $k$. Por lo tanto, como veremos a 
continuaci\'on,  $S$ es de la forma de los 
sumandos directos del enunciados del teorema y $M\sim M'$.

Se tiene $M'=\big(M\oplus \Lambda v\big)/\big\langle (p^ku_1, -p^k v),
(\lambda_2u_2+\cdots+\lambda_nu_n, \lambda_1v)\big\rangle$.
Como en la Operaci\'on 1, puesto que $p\nmid \lambda_1$, podemos
encajar $M$ en $M'$ y tambi\'en el ideal $\langle p^k,\lambda_1
\rangle$ aniquila a $M'/M$. Por lo tanto $|M'/M|<\infty$ y $M\sim M'$.

Ahora bien, usando que $p^k(u_1-v)=0$ y que $p^k$ divide al 
primer coeficiente entre las relaciones que involucran a $u_1$, se
tiene que $M'=M''\oplus \Lambda(u_1-v)$ donde $M''$ est\'a generado
por $v, u_2,\ldots, u_n$ y tiene relaciones generadas por $(\lambda_1,
\ldots, \lambda_n)$ y $R$. Por tanto $M''$ tiene $R'$ por relaci\'on.
Notemos que $\Lambda(u_1-v)\cong \Lambda/\langle p^k\rangle$
que es de la forma de la suma directa buscada.

\noindent
{\underline{Operaci\'on 3}}: Si $R$ contiene una fila $(p^k\lambda_1,
\ldots, p^k \lambda_n)$ y para alg\'un $\lambda$ con $p\nmid \lambda
$, $(\lambda\lambda_1,\ldots, \lambda \lambda_n)$ es tambi\'en
una relaci\'on, no necesariamente expl{\'\i}citamente contenida en $R
$, entonces podemos cambiar $R$ a $R'$, donde $R'$ es la misma
que $R$ excepto que $(p^k\lambda_1,\ldots, p^k\lambda_n)$ se 
reemplaza por $(\lambda_1,\ldots,\lambda_n)$. Veamos que el
m\'odulo $M'$ obtenido es seudo isomorfo a $M$. 

Para probar nuestra afirmaci\'on, consideremos la 
proyecci\'on can\'onica $\varphi\colon M\to
M':=M/\Lambda(\lambda_1u_1+\cdots+\lambda_n u_n)$.
Se tiene que $\ker \varphi$ es aniquilado
por el ideal $\langle \lambda,p^k\rangle$. Puesto que $M$ es
finitamente generado, $\ker \varphi$ es finitamente generado y
$\Lambda/\langle \lambda,p^k\rangle$ es finito si y s\'olo si $\ker
\varphi$ es finito. Puesto que $\varphi$ es suprayectiva, $M\sim M'$
y $M'$ tiene la matriz de relaciones $R'$.

Las seis operaciones A, B, C, 1, 2 y 3 se llaman {\em
admisibles\index{operaciones admisibles}} y todas ellas conservan las
dimensiones de la matriz original.

Continuando con la demostraci\'on del teorema, sea $f\in\Lambda
\setminus\{0\}$. Entonces por el Teorema de Preparaci\'on de
Weierstrass, tenemos $f(T)=p^uP(T)U(T)$ con $P(T)$ un polinomio
distinguido y $U\in\Lambda^{\ast}$. Definimos
\[
\gr_{\omega} f:=
\begin{cases}
\infty&\text{si $\mu>0$}\\ \gr P&\text{si $\mu=0$}
\end{cases}.
\]
$\gr_{\omega} f$ se llama el {\em grado de Weierstrass de
$f$\index{grado de Weierstrass}}. Dada una matriz $R$ definimos
\[
\gr^{(k)}R:=\min \gr_{\omega}(a_{i,j}')\quad \text{para}\quad i,j\geq k
\]
donde $(a_{i,j}')$ recorre todo el conjunto de matrices obtenidas a
partir de $R$ por medio de operaciones admisibles que dejan las
primeras $(k-1)$ filas sin cambio. Permitimos que $a_{i,j}$ cambie
para $i\geq k$ y cualquier $j$. Tambi\'en permitimos operaciones
del tipo B que usa, pero que no cambia, las primeras $(k-1)$ filas.

Si la matriz $R$ tiene la forma
\[
\left(
\begin{array}{ccccccc}
\lambda_{1,1}&0&\cdots&0&0&\cdots&0\\
\vdots&\vdots&&\vdots&\vdots&&\vdots\\
0&0&\cdots&\lambda_{r-1,r-1}&0&\cdots&0\\
\ast&\ast&\cdots&\ast&\ast&\cdots&\ast\\
\ast&\ast&\cdots&\ast&\ast&\cdots&\ast
\end{array}
\right):=\left(
\begin{array}{cc}
D_{r-1}&0\\ M&N\end{array}\right)
\]
con $\lambda_{k,k}$ un polinomio distinguido y $\gr\lambda_{k,k}=
\gr_{\omega}\lambda_{k,k}=\gr^{(k)}(R)$ para $1\leq k\leq r-1$, 
decimos que $R$ est\'a en la forma $(r-1)$--normal. 

Antes de continuar la demostraci\'on del teorema, probemos

\begin{lema}\label{L7.7.19} Si la submatriz $N$ es no cero, entonces
$R$ puede ser transformada con operaciones admisibles en una
matriz $R'$ que est\'a en la forma $r$--normal y que tiene los primeros
$(r-1)$ elementos diagonales iguales a los de $R$.
\end{lema}

\begin{proof}[Lema {\rm{\ref{L7.7.19}}}] 
Usando la Operaci\'on 1, podemos suponer, en caso necesario,
que una potencia grande de $p$ divide cada $\lambda_{i,j}$ con
$i\geq r$ y $j\leq r-1$, es decir, los elementos de $M$.

Esto es, $p^t|M$, para $t$ suficientemente grande y tal que $p^t\nmid
N$. Usando la Operaci\'on 2, podemos suponer $p\nmid N$. As{\'\i}
mismo, podemos suponer que $N$ tiene una entrada $\lambda_{i,j}$
tal que $\gr_{\omega}\lambda_{i,j}=\gr^{(r)}R<\infty$.

Si $\lambda_{i,j}=P(T)U(T)$ con $P(T)$ polinomio distinguido y
$U(T)\in \Lambda^{\ast}$, multipliquemos la columna $j$ por $U^{
-1}$. De esta forma podemos suponer que $\lambda_{i,j}$ es
distinguido puesto que las primeras $r-1$ filas tienen $0$ en la
columna $j$, y por lo tanto esos elementos no cambian. Por la
Operaci\'on A, podemos suponer $\lambda_{ij}=\lambda_{rr}$
nuevamente por la raz\'on que tenemos $0$ en los primeros lugares.

Por el algoritmo de la divisi\'on, el cual es un caso especial de la
Operaci\'on B, podemos suponer que $\lambda_{rj}$ es un polinomio
con $\gr \lambda_{rj}<\gr \lambda_{rr}$, $j\neq r$ y $\gr \lambda_{rj}
<\gr \lambda_{jj}$, $j<r$.

Ahora bien $\lambda_{rr}$ tiene grado de Weierstrass minimal en
$N$ por lo que se debe tener que $p|\lambda_{rj}$ para $j>r$.
Por la Operaci\'on 1, podemos suponer que $p^t|\lambda_{rj}$,
$j>r$, $t$ suficientemente grande. Si $\lambda_{rj}\neq 0$ para
alg\'un $j>r$, por la Operaci\'on 1, podemos quitar la potencia de $p$
de alg\'un $\lambda_{rj}$ con $j>r$ y los ceros siguen sin cambios
alguno. Tenemos
\[
\gr_{\omega}\lambda_{rj}\leq \gr \lambda_{rj}<\gr \lambda_{rr}=
\gr_{\omega} \lambda_{jj}
\]
lo cual es imposible. Consecuentemente, $\lambda_{rj}=0$ para
$j>r$.

 Si $\lambda_{rj}\neq 0$ para $j<r$ por la Operaci\'on $1$ podemos
 obtener $p\nmid \lambda_{rj}$ para alguna $j$, pero en este
caso se tiene
\[
\gr_{\omega}\lambda_{rj}\leq \gr \lambda_{rj}<\gr \lambda_{jj}
=\gr_{\omega}\lambda_{jj}
\]
y ya que $\gr_{\omega}\lambda_{jj}=\gr^{(j)}(R)$ se obtiene una
contradicci\'on a la definici\'on de $\gr^{(j)}(R)$. De esto 
concluimos que $\lambda_{rj}=0$ para todo $j\neq r$ que es lo que
quer{\'\i}amos probar. $\fin$
\end{proof}

Continuamos con la demostraci\'on del teorema. Por el Lema 
\ref{L7.7.19} empezamos con la matriz con $r=1$ y podemos ir 
cambiando $R$ hasta obtener una matriz
\[
\left(
\begin{array}{ccc|cc}
\lambda_{11}&&\text{\huge{$0$}}\\
&\ddots&&&\text{\Huge{$0$}}\\
\text{\huge{$0$}}&&\lambda_{rr}\\
{\ }&{\ }&{\ }&{\ }&{\ }\\
\hline
{\ }&{\ }&{\ }&{\ }&{\ }\\
&\text{\huge{$M$}}&&&\hspace{.5cm}\text{\huge{$0$}}\hspace{.5cm}
\end{array}
\right)
\]
con cada $\lambda_{jj}$ un polinomio distinguido y $\gr \lambda_{jj}
=\gr^{(j)} R$ para $j\leq r$. Por el algoritmo de la divisi\'on podemos
suponer $\lambda_{ij}$ es un polinomio y que $\gr\lambda_{ij}<
\gr\lambda_{jj}$ para $i\neq j$.
Si tuvi\'esemos $\lambda_{ij}\neq 0$ para alg\'un $i\neq j$, puesto 
que $\gr_{\omega}\lambda_{jj}$ es minimal, necesariamente $p|
\lambda_{ij}$ y por ende tenemos una relaci\'on no cero $(\lambda_{
i1},\ldots,\lambda_{ir},0,\ldots, 0)$ que es divisible por $p$. Sea
$\lambda:=\lambda_{11}\cdots \lambda_{rr}$. Entonces $p$ no 
divide a $\lambda$ y puesto que $\lambda_{jj}$ es un polinomio
distinguido para $j=1,\ldots, r$, se tiene que 
\[
\Big(\lambda\frac{1}{p}\lambda_{i1},\ldots, \lambda\frac{1}{p}\lambda_{
ir},0,\ldots,0\Big)
\]
tambi\'en es una relaci\'on puesto que $\lambda_{jj}u_j=0$.

Por la Operaci\'on 3, podemos suponer que $P$ no divide a $
\lambda_{ij}$ para alg\'un $j$, as{\'\i} que $\gr_{\omega}\lambda_{ij}
\leq \gr\lambda_{ij}<\gr \lambda_{jj}=\gr^{(j)}R$ lo cual es imposible.
Por lo tanto $\lambda_{ij}=0$ para todo $i,j$ tales que $i\neq j$.
Esto significa que $M=0$.

En t\'erminos de $\Lambda$--m\'odulos, tenemos:
\[
M'=\Lambda/\langle \lambda_{11}\rangle \oplus
\Lambda/\langle \lambda_{22}\rangle\oplus \cdots \oplus
\Lambda/\langle \lambda_{rr}\rangle\oplus \Lambda^{n-r}.
\]

Volviendo a escribir los factores $\Lambda/\langle p^k\rangle$ que
quitamos durante la Operaci\'on 2 obtenemos el resultado con la
salvedad de que $\lambda_{ii}$ no necesariamente es una
potencia de un polinomio irreducible, pero puesto que si $f$ y $g$
son primos relativos tenemos
\[
\Lambda/\langle f,g\rangle \sim \Lambda/\langle f \rangle\oplus
\Lambda/\langle g\rangle\quad \text{y}\quad 
\Lambda/\langle f\rangle\oplus \Lambda/\langle g\rangle\sim
\Lambda/\langle f, g\rangle
\]
esto termina la demostraci\'on del teorema. $\fin$
\end{proof}

\section{Los invariantes de Iwasawa\index{Iwasawa!invariantes de
$\sim$}}\label{S7.8}

Nuestro objetivo en esta secci\'on es probar lo siguiente: Sea $K$
un campo num\'erico finito y sea $K_{\infty}/K_0=K$ una extensi\'on
${\ma Z}_p$. Si $p^{e_n}$ es la potencia exacta de $p$ que divide
al n\'umero de clase de $K_n$, $h(K_n)$, entonces existen enteros
no negativos $\lambda\geq 0$, $\mu\geq 0$ y un entero $\gamma$,
independientes de $n$ y un n\'umero natural $n_0$ tal que para 
$n\geq n_0$, 
\begin{equation}\label{Ec7.8.1}
e_n=\mu p^n+\lambda n+\gamma.
\end{equation}

\begin{definicion}\label{D7.8.2} Los n\'umeros $\mu,\lambda, \gamma$
dados en (\ref{Ec7.8.1}) se llaman los {\em los invariantes de
Iwasawa}.
\end{definicion}

Sea $K_{\infty}/K$ una extensi\'on ${\ma Z}_p$, $\Gamma=\Gal(K_{
\infty}/K)\cong {\ma Z}_p$ y sea $\gamma_0$ un generador
topol\'ogico de $\Gamma$.
Sea $L_n$ la m\'axima $p$--extensi\'on abeliana no ramificada de
$K_n$. Se tiene que $X_n:=\Gal(L_n/K_n)\cong A_n$, donde $A_n$
es el $p$--subgrupo de Sylow del grupo de clases $I_{K_n}$ del
campo $K_n$.

\[
\xymatrix{
L_0\ar@{-}[r]\ar@{-}[d]^{A_0}&L_n\ar@{--}[r]\ar@{-}[d]^{A_n}&
L\\
K=K_0\ar@{-}[r]&K_n\ar@{--}[r]&K_{\infty}}
\]
Sea $L:=\bigcup\limits_{n=1}^{\infty}L_n$. Notemos que $L_n
\subseteq L_{n+1}$ ya que $K_{n+1}L_n/K_{n+1}$
es una extensi\'on no ramificada.

\begin{window}[0,r,\xymatrix{
&L\\K_{\infty}\ar@{-}[ru]^X\ar@{-}[d]_{\Gamma}\\ K\ar@{-}[ruu]_G},{}]
Sea $X:=\Gal(L/K_{\infty})$.
Puesto que $L_n$ es maximal, se tiene que $L_n/K$ es una 
extensi\'on de Galois.
Sea $G:=\Gal(L/K)$. Se tiene la sucesi\'on exacta
$ 1\longto X\longto G\longto \Gamma\longto 1$.
\end{window}

\vspace{1cm}

\begin{definicion}\label{D7.8.3} La condici\'on (A) se define por:
$$
\text{Todos los primos ramificados de $K_{\infty}/K$ son 
totalmente ramificados.}\leqno{\text{(A)}}
$$
\end{definicion}

\begin{observacion}\label{O7.8.4} Como vimos anteriormente
(Proposici\'on \ref{P7.6.7}), si $K_{\infty}/K$ no satisface la
condici\'on (A), existen un n\'umero natural $m$ tal que $K_{\infty}/
K_m$ es totalmente ramificado en cada primo ramificado y por lo
tanto $K_{\infty}/K_m$ satisface la condici\'on (A).
\end{observacion}

\begin{window}[0,r,\xymatrix{
L_n\ar@{-}[d]\\K_n\ar@{-}[r]&K_{n+1}},{}]
Supondremos que $K_{\infty}/K$ satisface la condici\'on (A). 
Entonces $L_n/K_n$ es no ramificada y $K_{n+1}/K_n$ es 
totalmente ramificada por lo que $K_{n+1}\cap L_n=K_n$ lo cual
implica que $\Gal(L_n/K_n)\cong \Gal(L_nK_{n+1}/K_{n+1})
\cong \frac{\Gal(L_{n+1}/K_{n+1})}{\Gal(L_{n+1}/L_nK_{n+1})}$.
Por lo tanto $X_n\cong \frac{X_{n+1}}{\Gal(L_{n+1}/L_nK_{n+1})}$.
\end{window}

Se tiene el mapeo restricci\'on suprayectivo: $X_{n+1}\stackrel
{\rest}{\longto}X_n$. Se tiene el siguiente resultado de la Teor{\'\i}a
de Campos de Clase:

\begin{teorema}\label{T7.8.5} (Ver el Teorema {\rm{\ref{T17.6.135N}}}).
Si $L/K$ y $L'/K'$ son dos extensiones
abelianas de campos num\'ericos tales que $K\subseteq K'$ y $L
\subseteq L'$. Entonces tenemos el siguiente diagrama
conmutativo
\[
\xymatrix{
I_{K'}\ar[rr]^{\hspace{-.5cm}
(\ ,L'|K')}\ar[d]^{N_{K'/K}}&&\Gal(L'/K')\ar[d]^{\rest}\\
I_K\ar[rr]_{\hspace{-.5cm}(\ ,L|K)}&&\Gal(L/K)
}
\]
donde $I_K$ e $I_{K'}$ son los grupos ya sea de divisores o de 
id\`eles de $K$ y $K'$ respectivamente, $N_{K'/K}$ es la norma de
$I_{K'}$ a $I_K$; $[\ ,L|K]$ y $[\ ,L'|K']$ son lo mapeos de Artin y
$\rest\colon \Gal(L'/K')\to \Gal(L/L\cap K')\subseteq \Gal(L/K)$,
$\sigma\mapsto \sigma|_L$, es el mapeo restricci\'on.
\end{teorema}

\begin{proof}
Sea $\pK$ un ideal de $K'$ no ramificado en $L'/K'$ y sea
$\pL$ un ideal de $L'$ sobre $\pK$. Similarmente, sean $\tilde{\pK}$,
$\tilde{\pL}$ para $L$ y $K$ de tal forma que $\tilde{\pK}=\pK|_K$,
$\tilde{\pL}=\pL|_L$, $\tilde{\pK}$ no ramificado en $L/K$.
\[
\xymatrix{
\pK\ar@{--}[r]&K'\ar@{-}[r]\ar@{-}[d]&K'L\ar@{-}[r]\ar@{-}[d]&L'
\ar@{--}[r]&\pL\\
\tilde{\pK}\ar@{--}[r]&K\ar@{-}[r]&L\ar@{--}[r]&\tilde{\pL}
}
\]

Sea $f=\big[{\cal O}_{K'}/\pK :{\cal O}_K/\tilde{\pK}\big]$ el grado
relativo. Entonces $N_{K'/K}\pK=\tilde{\pK}^f$ y si $N$ es la norma
absoluta, $N\pK=\big(N\tilde{\pK}\big)^f$. Puesto que ${\cal O}_L
\subseteq {\cal O}_{L'}$, se tiene 
\begin{gather*}
\sigma_{\pK}^{L'/K'}:=(\pK, L'|K')\quad \text{y}\quad
\sigma_{\pK}^{L'/K'}|_L(x)\equiv x^{N\pK}\bmod \tilde{\pL}
\quad\text{para}\quad x\in{\cal O}_L.\\
\intertext{Ahora}
\begin{align*}
\sigma_{N_{K'/K}\pK}^{L/K}(x)&=\big(N_{K'/K}\pK, L|K\big)(x)=
\big(\tilde{\pK}, L/K\big)^f(x)\\
&\equiv x^{(N\tilde{\pK})^f}=x^{N\pK}
\bmod \tilde{\pL}\quad \text{para}\quad x\in{\cal O}_L.
\intertext{Por lo tanto}
\rest (\ ,L'|K')(\pL)&=\rest\circ (\pL,L'|K')=\big(N_{K'/K}\pL, L|K\big)\\
&=\big((\ ,L|K)\circ N_{K'/K}\big)(\pL). \tag*{$\fin$}
\end{align*}
\end{gather*}
\end{proof}

Aplicando el Teorema \ref{T7.8.5}, tenemos el diagrama conmutativo
\[
\xymatrix{
X_{n+1}\ar[rr]^{\rest}
\ar[d]_{\big(\ ,L_{n+1}|K_{n+1}\big)}^{\cong}&&X_n\ar[d]_{\cong}^{
\big(\ ,L_n|K_n\big)}\\
A_{n+1}\ar[rr]_{N_{K_{n+1}/K_{n}}}&&A_n
}
\]
corresponde a la norma: $A_{n+1}\to A_n$ sobre el grupo de clases.
Ahora bien, $\Gal(L_nK_{\infty}/K_{\infty})\cong \Gal(L_n/L_n\cap
K_{\infty}=K_n)\cong X_n$. Por tanto
\begin{align*}
\lim_{\substack{\longleftarrow\\ n}}X_n&\cong \lim_{\substack{
\longleftarrow\\ n}} \Gal(L_nK_{\infty}/K_{\infty})=\Gal\Big(\Big(
\bigcup_{n=1}^{\infty}L_nK_{\infty}\Big)/K_{\infty}\Big)\\
&=\Gal(L/K_{\infty})\cong X.
\end{align*}
Esto es, $\lim\limits_{\substack{\longleftarrow\\ n}}X_n=\lim\limits_{
\substack{\longleftarrow\\ n}}\Gal\big(L_n/K_n)\cong \Gal(L/K_{\infty})
\cong X$.

Sea $\gamma\in \Gamma_n=\Gamma/\Gamma^{p^n}$ y 
extendamos $\gamma$ a $\tilde{\gamma}\in \Gal(L_n/K)$. Sea $x\in
X_n$. Entonces $\gamma$ act\'ua en $x$ por: $x^{\gamma}:=\tilde{\gamma}
x\tilde{\gamma}^{-1}$.
\[
\xymatrix{
&L_n\ar@{-}[d]^{X_n}\\ K\ar@{-}[r]_{\Gamma_n}\ar@{-}[ru]&K_n
}\qquad
\begin{array}{c}
{\ }\\{\ }\\{\ }\\
1\longto X_n\longto \Gal(L_n/K)\longto \Gamma_n\longto 1.
\end{array}
\]

Puesto que $X_n=\Gal(L_n/K_n)$ es abeliano, $x^{\gamma}$ est\'a
bien definida pues si $\gamma'$ es tal que $\tilde{\gamma}
(\tilde{\gamma}')^{-1}\in X_n$, entonces 
\[
\big(\tilde{\gamma}^{-1}(\tilde{\gamma}')\big)x\big(\tilde{\gamma}^{-1}
\tilde{\gamma}'\big)^{-1}=x\quad \text{por lo que}\quad
\tilde{\gamma}'x\tilde{\gamma}'^{-1}=\tilde{\gamma}x\tilde{\gamma}^{
-1}.
\]

Esta acci\'on corresponde a la acci\'on de $\Gamma_n$ en $A_n$ 
ya que $\artinp{L|K}{\sigma {\pK}}=\sigma
\artinp{L|K}{\pK}\sigma^{
-1}$, donde $\artinp{\star}{\star}$ denota al mapeo de reciprocidad
de Artin.  Esto es, $X_n$ es un ${\ma Z}_p[\Gamma_n]$--m\'odulo. Si 
representamos un elemento $X\cong\lim\limits_{\substack
{\longleftarrow\\ n}}X_n$ como un vector $(x_0,x_1,\ldots)$ y 
dejando ${\ma Z}_p[\Gamma_n]$ actuar en la $n$--\'esima 
componente, entonces $X$ es un $\Lambda$--m\'odulo, $\Lambda
=\lim\limits_{\substack{\longleftarrow\\ n}}{\ma Z}_p[\Gamma_n]$:
si $\vec{x}\in X$ y $\vec{\lambda}\in\Lambda$, $\vec{\lambda}=
(\lambda_0,\lambda_1,\ldots)$, $\vec{x}^{\vec{\lambda}}:=
\big(\tilde{\lambda}_0x_0\tilde{\lambda}_0^{-1}, 
\tilde{\lambda}_1x_1\tilde{\lambda}_1^{-1},\ldots,
\tilde{\lambda}_nx_n\tilde{\lambda}_n^{-1}\ldots\big)$. Ahora
\[
\tilde{\lambda}_{n+1}x_{n+1}\tilde{\lambda}_{n+1}^{-1}
\xrightarrow[]{\rest}
\tilde{\lambda}_{n+1}\big|_{K_n}x_{n+1}\big|_{K_n}
\tilde{\lambda}_{n+1}^{-1}\big|_{K_n}=
\tilde{\lambda}_{n}x_{n}\tilde{\lambda}_{n}^{-1},
\]
por lo tanto
$\vec{x}^{\vec{\lambda}}\in X$.

Con el isomorfismo $\Lambda\cong {\ma Z}_p[[T]]$, $1+T\in
\Lambda$ act\'ua como $\gamma_0\in \Gamma$. Se tiene
$x^{\gamma}=\tilde{\gamma}x\tilde{\gamma}^{-1}$ para
$\gamma\in \Gamma$, $x\in X$ y $\tilde{\gamma}$ es una
extensi\'on de $\gamma$ a $G$.

Sean $\pK_1,\ldots, \pK_s$ los primos ramificados en $K_{\infty}/K$
y sean $\tilde{\pK}_1,\ldots, \tilde{\pK}_s$ primos de $L$ sobre cada
$\pK_i$. Sean $I_i\subseteq G$ los grupos de inercia de $\tilde{\pK}_i
$, $1\leq i\leq s$.

\begin{window}[0,l,\xymatrix{&L\ar@{-}[dl]_X\\ K_{\infty}\ar@{-}[d]_{
\Gamma}\\K},{}]
Puesto que $L/K_{\infty}$ es no ramificada, $I_i\cap X=\{1\}$. Puesto
que $K_{\infty}/K$ es totalmente ramificado en $\pK_i$, se tiene que
el mapeo $I_i\hookrightarrow G/X\cong \Gamma$ es suprayectivo
y por tanto un isomorfismo. De esta forma tenemos que $G=I_iX=
XI_i$, $1\leq i\leq s$. Sea $\sigma_i\in I_i$ el elemento que se 
mapea a $\gamma_0$. Entonces $\sigma_i$ es un generador
topol\'ogico de $I_i$. Puesto que $I_i\subseteq XI_1$, se tiene
que existe $a_i\in X$ tal que $\sigma_i=a_i\sigma_1$, y $a_1=1$.
\end{window}

\begin{proposicion}\label{P7.8.6} Si la extensi\'on $K_{\infty}/K$
satisface la condici\'on {\rm (A)}, entonces si $G'$ es la cerradura del
subgrupo conmutador de $G$, se tiene $G'=X^{\gamma_0-1}=TX$.
\end{proposicion}

\begin{proof}
Puesto que $\Gamma\cong I_1\subseteq G$, y $I_1$ se mapea sobre
$\Gamma\cong G/X$, se puede levantar $\gamma$ al 
correspondiente elemento en $I_1$ para as{\'\i} definir la acci\'on de
$\Gamma$ en $X$. Por simplicidad, identificamos $\Gamma$ y $
I_1$ de tal manera que $x^{\gamma}=\gamma x\gamma^{-1}$. Sea
$a=\alpha x$, $b=\beta y$ con $\alpha,\beta\in \Gamma$ y $x,y\in
X$ elementos arbitrarios de $G=\Gamma X$. Entonces
\begin{align*}
aba^{-1}b^{-1}&= \alpha x \beta y x^{-1}\alpha^{-1}y^{-1}\beta^{-1}=
x^{\alpha}\alpha \beta yx^{-1}\alpha^{-1}y^{-1}\beta^{-1}\\
&=x^{\alpha}\big(yx^{-1}\big)^{\alpha\beta}\underbrace{(\alpha
\beta)\alpha^{-1}}_{\in \Gamma}y^{-1}\beta^{-1}\igual_{\substack{
\uparrow\\ \Gamma\text{\ abeliano}}}x^{\alpha}\big(yx^{-1}\big)^{\alpha
\beta}\big(y^{-1}\big)^{\beta}\\
&\igual_{\substack{\uparrow\\ X\text{\ abeliano}}} x^{\alpha}
\big(x^{-1}\big)^{\alpha\beta}y^{\alpha\beta}\big(y^{-1}\big)^{\beta}=
\big(x^{\alpha}\big)^{1-\beta}\big(y^{\beta}\big)^{\alpha-1}\in G'.
\end{align*}

Tomemos $\beta=1$, $\alpha=\gamma_0$, entonces $\big(x^{\alpha}
\big)^0y^{\gamma_0-1}=y^{\gamma_0-1}\in G'$, esto es, $X^{
\gamma_0-1}\subseteq G'$.

Para $\beta\in\Gamma$ arbitrario, existe $c\in{\ma Z}_p$ tal que
$\beta=\gamma_0^c$. Por tanto
\[
1-\beta =1-\gamma_0^c=\big(1-(1+T)^c\big)=1-\sum_{n=0}^{\infty}
\binom{c}{n}T^n\in T\Lambda
\]
donde $\binom{c}{n}:=\frac{c(c-1)\cdots (c-n+1)}{n!}$, $n\in {\ma N}
\cup \{0\}$. Puesto que $\gamma_0-1=T$, $\big(x^{\alpha}\big)^{1-
\beta}\in X^{\gamma_0-1}$.

Similarmente obtenemos que $\big(y^{\beta}\big)^{1-\alpha}\in
X^{\gamma_0-1}$. Finalmente, puesto que $X^{\gamma_0-1}=TX$
es cerrado por ser la imagen de conjunto compacto $X$, se sigue
que $G'\subseteq X^{\gamma_0-1}$ y por ende $G'=X^{\gamma_0
-1}$. $\fin$
\end{proof}

\begin{proposicion}\label{P7.8.7} Sea $K_{\infty}/K$ que satisface
la condici\'on {\rm (A)}. Sea $Y_0$ el ${\ma Z}_p$ subm\'odulo de $X$
generado por $\{a_i\mid 2\leq i\leq s\}$ y por $X^{\gamma_0-1}
=TX$ donde tenemos $\sigma_i=a_i\sigma_1$, $\overline{\langle
\sigma_1\rangle}=I_i$. Sea $Y_n:=\gamma_n Y_0$, donde
\[
\gamma_n:=1+\gamma_0+\gamma_0^2+\cdots+\gamma_0^{p^n-1}
=\frac{\gamma_0^{p^n}-1}{\gamma_0-1}=\frac{(1+T)^{p^n}-1}{T}.
\]

Entonces $X_n\cong X/Y_n$ para $n\geq 0$.
\end{proposicion}

\begin{proof}
Primero sea $n=0$. Se tiene $K\subseteq L_0\subseteq L$ y 
$L_0$ es la m\'axima $p$--extensi\'on abeliana no ramificada de $K$.
Entonces $\Gal(L/L_0)$ es el subgrupo cerrado de $G$ generado
por $G'$, pues $L_0/K$ es abeliana, y por todos los subgrupos
de inercia $I_i$, $1\leq i\leq s$ ya que $L_0/K$ es no ramificada. Por
lo tanto $\Gal(L/L_0)$ es el subgrupo cerrado generado por $X^{
\gamma_0-1}$, $I_1$ y $a_2,\ldots, a_s$ de tal forma que 
\begin{align*}
X_0&=\Gal(L_0/K)=\frac{G}{\Gal(L/L_0)}=\frac{XI_1}{\Gal(L/L_0)}\\
&\cong \frac{I_1X}{\langle X^{\gamma_0-1}, I_1,a_2,\ldots,a_2\rangle}
\cong \frac{X}{\langle X^{\gamma_0-1},a_2,\ldots,a_2\rangle}=
\frac{X}{Y_0}.
\end{align*}

Consideremos $n\geq 1$, Reemplazando $K$ por $K_n$ y $
\gamma_0$ por $\gamma_0^{p^n}$ se tiene que  $\sigma_i^{p^n}$
toma el papel de $\sigma_i$. Se sigue que
\begin{gather*}
\xymatrix{
&L_n\ar@{-}[d]^{X_n}\ar@{-}[r]&L\ar@{-}[d]\\ 
K\ar@{-}[r]&K_n\ar@{-}[r]&K_{
\infty}}\\
\begin{align*}
\sigma_i^{k+1}&=(a_i\sigma_1)^{k+1}=a_i\sigma_1a_i\sigma_1^{-1}
\sigma_1^2a_i\sigma_1^{-2}\sigma_1^3\cdots\sigma_1^ka_i
\sigma_1^{-k}\sigma_1^{k+1}\\
&=a_i^{1+\sigma_1+\cdots+\sigma_1^k}
\sigma_i^{k+1}.
\end{align*}
\end{gather*}

Por tanto $\sigma_i^{p^n}\igual\limits_{\substack{\uparrow\\ \sigma_1
=\gamma_0}} (\gamma_n a_i)\sigma_1^{p^n}$, esto es, $a_i$ es
reemplazado por $\gamma_n a_i$.

Finalmente $X^{\gamma_0-1}$ es reemplazado por $X^{\gamma_0^{
p^n}-1}=\gamma_nX^{\gamma_0-1}$ y $Y_0$ es reemplazado por
$\gamma_n Y_0$ lo cual implica el resultado. $\fin$
\end{proof}

\begin{observacion}\label{O7.8.8} La Proposici\'on \ref{P7.8.7}
es crucial pues nos permite recabar informaci\'on para $X_n$ de la 
informaci\'on que se tenga de $X$.
\end{observacion}

Veamos un resultado b\'asico de \'algebra conmutativa aplicado 
a nuestro caso (ver Teorema \ref{T7.3.6}).

\begin{proposicion}[Lema de Nakayama\index{lema de
Nakayama}\index{Nakayama!lema de $\sim$}]\label{P7.8.9}
Sea $X$ un $\Lambda$--m\'odulo compacto, es decir, $X$ tiene una
topolog{\'\i}a compacta Hausdorff, $X$ es un $\Lambda$--m\'odulo
y la acci\'on de $\Lambda$ en $X$; $\Lambda\times X\to X$, es una
funci\'on continua.

Entonces $X$ es finitamente generado como $\Lambda$--m\'odulo
si y solamente si $X/\langle p,T\rangle X$ es finito.

Si $x_1,\ldots, x_n$ generan $X/\langle p,T\rangle X$ sobre ${\ma Z}
$, entonces $x_1,\ldots,x_n$ generan $X$ como 
$\Lambda$--m\'odulo. En particular 
\[
X/\langle p,T\rangle X=0\iff X=0.
\]
\end{proposicion}

\begin{proof} Sea $U$ una vecindad abierta de $0\in X$. Puesto que
$\langle p, T\rangle^n\longto 0$ en $\Lambda$, cada $z\in X$ tiene
una vecindad abierta $U_z$ tal que $\langle p, T\rangle^{n_z}U_x
\subseteq U$ para alg\'un $n_z\in {\ma N}$ suficientemente grande.
Puesto que $X$ es compacto, un n\'umero finito de las vecindades
$U_z$ cubren a $X$. En particular, para $n$ suficientemente grande,
$\langle p,T\rangle^n X\subseteq U$ y por tanto 
\[
\bigcap_{n=1}^{\infty}\big(\langle p,T\rangle^n X\big)=\{0\}.
\]

Ahora bien, $\Lambda/\langle p, T\rangle \cong{\ma F}_p$ es
finito y $X/\langle p, T\rangle X$ es un $\Lambda/\langle p, T
\rangle$--m\'odulo, de donde se sigue que $X//\langle p, T\rangle X$
es finitamente generado si y solamente si es finito.

Sean $x_1,\ldots, x_n$ generadores $X/\langle p, T\rangle X$ y sea
$Y=\Lambda x_1+\cdots +\Lambda x_n\subseteq X$. Puesto que 
$Y$ es imagen continua de $\Lambda^n$, $Y$ es un conjunto
compacto y por lo tanto $Y$ es cerrado en $X$. Se sigue que 
$X/Y$ es un $\Lambda$--m\'odulo compacto. Tenemos
\[
\langle p, T\rangle\big(X/Y\big)=\frac{T+\langle p, T\rangle X}{Y}=
\frac{X}{Y}
\]
de donde se sigue que $\langle p, T\rangle^n\big(X/Y\big)=X/Y$ 
para toda $n\geq 0$. Por lo tanto $X/Y=\bigcap\limits_{n=0}^{\infty}
\langle p, T\rangle^n\big(X/Y\big)=\{0\}$ y por lo tanto $X=Y$, esto es,
$\{x_1,\ldots, x_n\}$ genera a $X$. $\fin$
\end{proof}

\begin{proposicion}\label{P7.8.10} Si la extensi\'on $K_{\infty}/K$
satisface la condici\'on {\rm (A)} y $X=\Gal(L/K_{\infty})$, entonces $X$ es
un $\Lambda$--m\'odulo finitamente generado.
\end{proposicion}

\begin{proof} Se tiene $\gamma_1=1+\gamma_0+\cdots+\gamma_0^{
p-1}\in \langle p, T\rangle$, por lo que $Y_0/\langle p, T\rangle Y_0$
es cociente de $Y_0/\gamma_1 Y_0=Y_0/Y_1\subseteq X/Y_1=
X_1$ donde $Y_0=\langle X^{\gamma_0-1}, a_2,\ldots, a_s\rangle$,
$Y_n=\gamma_n Y_0$, $X_n\cong X/Y_n$.

Puesto que $X_1$ es finito, $Y_0$ es finitamente generado. Por otro
lado, $X/Y_0=X_0$ es finito lo cual implica que $X$ es finitamente
generado como $\Lambda$--m\'odulo. $\fin$
\end{proof}

Ahora analicemos el caso en que $K_{\infty}/K$ no necesariamente
satisface la condici\'on (A). Sea $e\geq 0$ tal que en $K_{\infty}/K_e$
todos los primos ramificados son totalmente ramificados. Los 
resultados anteriores se satisfacen para $K_{\infty}/K_e$. En 
particular $X$, el cual es el misma tanto para $K$ como para $K_e$,
es finitamente generado como $\Lambda$--m\'odulo. Para $n\geq e$
tenemos
\[
1=\gamma_0^{p^e}+\gamma_0^{2p^e}+\cdots +\gamma_0^{p^n-p^e}
=\frac{\gamma_n}{\gamma_e}=:\gamma_{n,e},
\]
de hecho, $\frac{\gamma_n}{\gamma_e}=\frac{
\gamma_0^{p^n}-1}{\gamma_0^{p^e}-1}=\frac{\big(\gamma_0^{p^e}
\big)^{p^{n-e}}-1}{\gamma_0^{p^e}-1}$.

Entonces $\gamma_{n,e}$ reemplaza $\gamma_n$ en la extensi\'on
$K_{\infty}/K_e$ puesto que $\gamma_0^{p^e}$ es un generador
topol\'ogico de $\Gal(K_{\infty}/K_e)$. Sea $Y_e$ el respectivo 
m\'odulo $Y_0$ para $K_e$ en lugar de $K$. Entonces $Y_n=
\gamma_{n,e}Y_e$ y $X_n=X/Y_n$ para $n\geq e$. Entonces
tenemos

\begin{proposicion}\label{P7.8.11} Sea $K_{\infty}/K$ una
extensi\'on ${\ma Z}_p$. Entonces se tiene que 
 $X:=\Gal(L/K_{\infty})$ es un
$\Lambda$--m\'odulo finitamente generado y existe $e\geq 0$ tal que
$X_n\cong X/\gamma_{n,e}Y_e$ para toda $n\geq e$. $\fin$
\end{proposicion}

Podemos aplicar el teorema de estructura de $\Lambda$--m\'odulos
finitamente generados para $X$ y para $Y_e$ y puesto que
$X/Y_e$ es finito, se tiene que 
\begin{equation}\label{Ec7.8.12}
Y_e \sim X\sim \Lambda^r\oplus \Big(\bigoplus_{i=1}^s \Lambda/
\langle p^{k_i}\rangle\Big)\oplus\Big(\bigoplus_{j=1}^t \Lambda/
\langle f_j(T)^{m_j}\rangle\Big).
\end{equation}

Calcularemos $V/\gamma_{n,e}V$ para cada sumando del seudo
isomorfismo dado en (\ref{Ec7.8.12}).
\las
\item Si $V=\Lambda$ se tiene que $\Lambda/\langle \gamma_{n,e}
\rangle$ es infinito. Puesto que $Y_e/\gamma_{n,e}Y_e\subseteq
X/Y_n$ es finito, se sigue que $\Lambda$ no puede aparecer como
sumando, esto es, $r=0$ en (\ref{Ec7.8.12}) y $X$ es de torsi\'on.

\item Si $V=\Lambda/\langle p^k\rangle$, entonces $V/\gamma_{n,e}
V=\Lambda/\langle p^k, \gamma_{n,e}\rangle$.

Se tiene que 
\begin{align*}
\gamma_{n,e}&=\frac{\gamma_n}{\gamma_e}=\frac{\big((1+T)^{p^n}-1
\big)/T}{\big((1+T)^{p^e}-1\big)/T}\\
&=1+(1+T)^{p^e}+(1+T)^{2p^e}+
\cdots+(1+T)^{p^n-p^e}
\end{align*}
es un polinomio distinguido.

Ahora, por el algoritmo de la divisi\'on, se tiene que cada elemento
de $\Lambda/\langle p^k,\gamma_{n,e}\rangle$ se representa 
un{\'\i}vocamente por un polinomio m\'odulo $p^k$ de grado menor
a $\gr \gamma_{n,e}=p^n-p^e$. Por lo tanto
\[
\big|V/\gamma_{n,e}V\big|=p^{k(p^n-p^e)}=p^{kp^n+c}
\]
donde $c$ es la constante $c=-kp^e$.

\item $V=\Lambda/\langle f(T)^m\rangle$ donde $f(T)$ es un 
polinomio distinguido e irreducible. Sea $g(T):=f(T)^m$. Entonces
$g$ tambi\'en es un polinomio distinguido, digamos de grado $d$.
Entonces
\[
T^d\equiv pQ(T)\bmod g
\]
para alg\'un $Q(T)\in{\ma Z}_p[T]$. Por lo tanto $T^k\equiv (p)S
\bmod g$, $S$ es un polinomio, $k\geq d$. As{\'\i}, si $p^n\geq d$,
se tiene
\[
(1+T)^{p^n}=1+(p)S_1+T^{p^n}\equiv 1+(p)S_2\bmod g
\]
con $S_1$, $S_2$ polinomios.

Por lo tanto $(1+T)^{p^{n+1}}\equiv 1+p^2 S_3\bmod g$, $S_3$ un
polinomio. En general se sigue que
\begin{align*}
P_{n+2}(T):&=(1+T)^{p^{n+2}}-1\\
&=\big((1+T)^{(p-1)p^{n+1}}+\cdots+
(1+T)^{p^{n+1}}+1\big)\\
&\hspace{1cm}\big((1+T)^{p^{n+1}}-1\big)\\
&\equiv (1+\cdots+1+p^2 S_4)P_{n+1}(T)\\
&\equiv p(1+pS_5)P_{n+1}
(T)\bmod g
\end{align*}
donde $S_4$, $S_5$ son polinomios.

Ahora bien, $1+pS_6\in\Lambda^{\ast}$, para $S_6$ polinomio, se
tiene que $\frac{P_{n+2}}{P_{n+1}}$ act\'ua como multiplicaci\'on 
$pu$, $u\in\Lambda^{\ast}$ sobre $V=\Lambda/\langle g(T)\rangle$
para $p^n\geq d$. Supongamos $n_0>e$, $p^{n_0}\geq d$ y $n
\geq n_0$. Entonces
\begin{gather*}
\frac{\gamma_{n+2,e}}{\gamma_{n+1,e}}=\frac{\gamma_{n+2}}
{\gamma_{n+1}}=\frac{P_{n+2}}{P_{n+1}}\quad \text{y}\quad
\gamma_{n+2,e}V=\frac{P_{n+2}}{P_{n+1}}\big(\gamma_{n+1,e}V
\big)=p\gamma_{n+1,e}V.\\
\intertext{Por tanto}
\big|V/\gamma_{n+2,e}V\big|=|V/pV|\big|pV/p\gamma_{n+1,e}V\big|
\quad\text{para}\quad n\geq n_0.
\end{gather*}

Puesto que $\mcd (g(T),p)=1$, se tiene que multiplicaci\'on por $p$
es inyectiva en donde obtenemos
\[
\big|pV/p\gamma_{n+1,e}V\big|=\big|V/\gamma_{n+1,e}V\big|
\quad\text{para}\quad n\geq n_0.
\]

Puesto que $V/pV=\Lambda/\langle p, g(T)\rangle =\Lambda/
\langle p, T^d\rangle$ se sigue que $|V/pV|=p^d$. Por inducci\'on
en $n$, obtenemos
\[
\big|V/\gamma_{n,e}V\big|=p^{d(n-n_0-1)}\big|V/\gamma_{n_0+1,e}
V\big|, \quad n\geq n_0+1.
\]

Si $V/\gamma_{n,e}V$ es finito para toda $n$, entonces 
$\big|V/\gamma_{n,e}V\big|=p^{dn+c}$ para $n\geq n_0+1$ y 
alguna constante $c$. Si $V/\gamma_{n,e}V$ es infinito, $V$
no puede ser sumando en el seudo isomorfismo (\ref{Ec7.8.12})
debido a que $\big|Y_e/\gamma_{n,e}Y_e\big|<\infty$. Este caso
puede suceder \'unicamente cuando $\mcd (\gamma_{n,e},f)\neq 1$.
\end{list}

Resumimos toda la discusi\'on anterior en:

\begin{teorema}\label{T7.8.13}
Sea $E=\Lambda^r\oplus\Big(\bigoplus\limits_{i=1}^s\Lambda/\langle
p^{k_i}\rangle\Big)\oplus\Big(\bigoplus\limits_{j=1}^t \Lambda/\langle
g_j(T)\Big)$ donde cada $g_j(T)$ es un polinomio distinguido, no
necesariamente irreducible. Sea $m:=\sum_{i=1}^s k_i$, $\ell=\sum_{
j=1}^t \gr g_j(T)$. Si $E/\gamma_{n,e}E$ es finito para toda $n$,
entonces $r=0$ y existen $n_0\in{\ma N}$  y $c\in{\ma Z}$ tales que
para $n>n_0$, $\big|E/\gamma_{n,e}E\big|=p^{mp^n+\ell n+c}$. 
$\fin$
\end{teorema}

\begin{teorema}\label{T7.8.14} Sea $E$ como en el Teorema 
{\rm \ref{T7.8.13}} con $r=0$. Entonces se tiene que  $m=0$ si y solamente
si el $p$--rango de $E/\gamma_{n,e}E$ permanece
acotado cuando $n\to\infty$.
\end{teorema}

\begin{proof}
Si $A$ es un grupo abeliano finito, tenemos que el $p$--rango de
$A$ es igual a $\dim_{{\ma F}_p}A/pA$. Ahora bien, recordemos que
$\gamma_{n,e}$ es distinguido de grado $p^n-p^e$, as{\'\i}, si $\gr
\gamma_{n,e}\geq \max\limits_{1\leq j\leq t}\gr g_j(T)$, se tiene que
\begin{align*}
E/\gamma_{n,e}E&=\Big(\bigoplus_{i=1}^s \Lambda/\langle p,
\gamma_{n,e}\rangle\Big)\oplus \Big(\bigoplus_{j=1}^t \Lambda/\langle
p,g_j(T),\gamma_{n,e}\rangle\Big)\\
&\cong \Big(\bigoplus_{i=1}^s \Lambda/\langle p, T^{p^n-p^e}\rangle
\Big)\oplus\Big(\bigoplus_{j=1}^t \Lambda/\langle p, T^{\gr g_j(T)}
\rangle\Big)\\
&\cong \big({\ma Z}/p{\ma Z}\big)^{s(p^n-p^e)+\ell}.
\end{align*}

Por tanto, el $p$--rango de $E$ est\'a acotado si y solamente si
$s=0$ si y solamente si $m=0$. $\fin$
\end{proof}

Regresando a la extensi\'on $K_{\infty}/K$, tenemos una sucesi\'on
exacta: \[
0\longto A\longto Y_e\longto E\longto B\longto 0
\]
 donde
$A$ y $B$ son finitos y $E$ es como antes. Conocemos
$\big|E/\gamma_{n,e}E\big|$ para $n>n_0$. De aqu{\'\i} ya 
podemos concluir que si $p^n$ es la potencia exacta que divide
a $|A_n|=|X_n|=\big|\frac{X}{\gamma_{n,e}Y_e}\big|$ se tiene 
$e_n=mp^n+\ell n+c_n$ con $c_n$ acotado.

La siguiente pieza es:

\begin{proposicion}\label{P7.8.15} Sean $Y$ y $E$ $\Lambda
$--m\'odulos seudo isomorfos de tal forma que $Y/\gamma_{n,e}Y$
es finito para $n\geq e$. Entonces existen una constante $c$ y alguna
$n_0$ tales que $\big|Y/\gamma_{n,e}Y\big|=p^c
\big|E/\gamma_{n,e}E\big|$ para toda $n\geq n_0$.
\end{proposicion}

\begin{proof}
Se tiene el siguiente diagrama conmutativo:
\[
\xymatrix{
0\ar[r]&\gamma_{n,e} Y\ar[r]\ar[d]^{\phi'_n}&Y\ar[r]\ar[d]^{\phi}&
Y/\gamma_{n,e}Y\ar[d]^{\phi''_n}\ar[r]&0\\
0\ar[r]&\gamma_{n,e}E\ar[r]&E\ar[r]&E/\gamma_{n,e}E\ar[r] &0
}
\]

Se tiene:
\l
\item Puesto que $\ker \phi'_n\subseteq \ker \phi$, $\big|\ker\phi'_n
\big|\leq |\ker \phi|$.

\item Se tiene $\coker \phi=E/\Phi(Y)$. Sea $E/\Phi(Y)=
\{\overline{u}_1,\ldots, \overline{u}_r\}$. Es decir, si $x\in E$, es tal 
que $x\equiv u_i\bmod \Phi(Y)$, se tiene que $\gamma_{n,e}x-
\gamma_{n,e}u_i\in\gamma_{n,e}\Phi(Y)=\Phi(\gamma_{n,e}Y)$ lo
cual implica que $\frac{\gamma_{n,e}E}{\Phi(\gamma_{n,e}Y)}=
\{\overline{\gamma_{n,e}u_1},\ldots, \overline{\gamma_{n,e}u_r}\}$ y
por lo tanto $\big|\coker \phi'_n\big|\leq |\coker \phi|$. 

\item Puesto que los representantes de $\coker \phi$ dan 
representantes de $\coker \phi''_n$ se sigue que 
$\big|\coker \phi''_n\big|\leq |\coker \phi|$

\item Aplicando el Lema de la Serpiente, se tiene la sucesi\'on 
exacta
\begin{align*}
0&\longto \ker \phi'_n\longto \ker \phi\longto\ker \phi''_n\longto\\
&\longto \coker \phi'_n\longto
\coker \phi\longto \coker \phi''_n\longto 0.
\end{align*}
\begin{align*}
\text{Por tanto\ }& |\ker \phi''_n|\leq |\ker \phi||\coker \phi'_n|
\quad \text{y por ({\sc {ii}})}\\
&|\ker \phi||\coker \phi'_n|\leq |\ker \phi||\coker \phi|.\\
\text{Esto es\ }& |\ker \phi''_n|\leq |\ker \phi||\coker \phi|.
\end{align*}
\end{list}

Ahora bien, si $m\geq n\geq 0$, tenemos
\lasa
\item $|\ker \phi'_n|\geq |\ker \phi'_m|$:

Se tiene $\gamma_{m,e}=\frac{\gamma_{m,e}}{\gamma_{n,e}}
\gamma_{n,e}$, por lo que $\gamma_{m,e}Y\subseteq \gamma_{
n,e}Y$ y por ende $\ker\phi'_m\subseteq \ker \phi'_n$.

\item $|\coker \phi'_n|\geq |\coker \phi'_m|$:

Sean $\gamma_{m,e}y\in 
\gamma_{m,e}E$ y $z\in\gamma_{n,e}E$ el representante del
elemento 
$\gamma_{n,e}y $ en $\coker \phi'_n=\frac{\gamma_{n,e}E}{
\gamma_{n,e}\phi'_n(Y)}$. Por tanto $\gamma_{n,e}y-z=\phi(
\gamma_{n,e}x)$ para alg\'un $x\in Y$. Multiplicando por $\frac{
\gamma_{m,e}}{\gamma_{n,e}}$, se tiene que
\[
\gamma_{m,e}y-\big(\frac{\gamma_{m,e}}{\gamma_{n,e}}\big)(z)=
\phi(\gamma_{m,e}x)=\phi'_n(\gamma_{m,e}x).
\]

Es decir, los representantes de $\coker\phi'_m$ se obtienen de
multiplicar los representantes de $\coker\phi'_n$ por $\frac{\gamma_{
m,e}}{\gamma_{n,e}}$.

\item $|\coker \phi''_n|\leq |\coker \phi''_m|$:

Se tiene que $\gamma_{m,e}E\subseteq \gamma_{n,e}E$. Por lo
tanto, del epimorfismo $E/\gamma_{m,e}E\to E/\gamma_{n,e}E$ se
tiene que $\coker \phi''_n\twoheadrightarrow \coker\phi''_n$.
\end{list}

Obtenemos $|\ker \phi'_m|\des\limits_{\substack{\uparrow\\
(a)}}|\ker \phi'_n|\des\limits_{\substack{\uparrow\\ (\text{\sc{i}})}}
|\ker \phi|$ para $m\geq n$,
\begin{gather*}
|\coker \phi'_m|
\des\limits_{\substack{\uparrow\\ (b)}}|\coker \phi'_n|
\des\limits_{\substack{\uparrow\\ (\text{\sc{ii}})}} |\coker \phi|
\quad \text{y}\\
|\coker \phi''_n|\des\limits_{\substack{\uparrow\\ (c)}}|\coker \phi''_m|
\des\limits_{\substack{\uparrow\\ (\text{\sc{iii}})}}|\coker \phi''|.
\end{gather*}

Por tanto existe $n_0\in{\ma N}$ tal que para $n\geq n_0$ los 
\'ordenes de $\ker \phi'_n$, $\coker \phi'_n$
 y $\coker \phi''_n$ son constantes.

Finalmente, por el Lema de la Serpiente, tenemos que 
\[
|\ker\phi'_n||\ker\phi''_n||\coker\phi|=|\ker\phi||\coker\phi'_n|
|\coker\phi''_n|
\]
por lo que $|\ker\phi''_n|$ es constante para $n\geq n_0$. As{\'\i},
para toda $n\geq n_0$, tenemos
\[
0\longto \ker \phi''_n\longto Y/\gamma_{n,e}Y
\stackrel{\phi''}{\longto} E/\gamma_{n,e}
E\longto \coker \phi''_n\longto 0
\]
lo cual implica que $\big|Y/\gamma_{n,e}Y\big|=
\frac{|\ker \phi''_n|}{|\coker \phi''_n|}\big|E/\gamma_{n,e}E\big|$. $\fin$
\end{proof}

\begin{teorema}[Iwasawa\index{Iwasawa!teorema de $\sim$}]
\label{T7,8.16} Sea $M$ un $\Lambda$--m\'odulos finitamente 
generado de torsi\'on y supongamos que para $n\in{\ma Z}$, $n\geq
0$, $M/\gamma_n M$ es un grupo finito de orden $p^{e_n}$. 
Entonces existen enteros no negativos $m$, $\ell$, y un entero $c$
tales que para $n$ suficientemente grande, $e_n=mp^n+\ell n+c$.
\end{teorema}

\begin{proof}
Sea $M\sim \Lambda^r\oplus\Big(\bigoplus\limits_{i=1}^s \Lambda/
\langle p^{k_i}\rangle\Big)\oplus \Big(\bigoplus\limits_{j=1}^t\Lambda/
\langle f_j^{m_j}(T)\rangle\Big)$. Puesto que $\big|M/\gamma_nM\big|
<\infty$, $r=0$ y por la Proposici\'on \ref{P7.8.15} se tiene $\big|
M/\gamma_nM\big|=p^{mp^n+\ell n+c}$ con $m=\sum_{i=1}^s k_i$,
$\ell=\sum_{j=1}^t m_j \gr f_j(T)$ y $c\in{\ma Z}$. $\fin$
\end{proof}

Como caso especial, obtenemos el resultado principal:

\begin{teorema}[Iwasawa]\label{T7.8.17} Sea $K_{\infty}/K$
una extensi\'on ${\ma Z}_p$. Sea $p^{e_n}$ la potencia exacta que
divide al n\'umero de clase de
 $K_n$. Entonces existen enteros $\lambda
\geq 0$, $\mu\geq 0$ y $\gamma$ independientes de $n$, y  un
natural $n_0$ tal que para $n\geq n_0$, $e_n=\mu p^n+\lambda n+
\gamma$. $\fin$
\end{teorema}

\begin{teorema}\label{T7.8.18} Sea $K_{\infty}/K$ una extensi\'on 
${\ma Z}_p$ en la cual exactamente un primo es ramificado y este es
totalmente ramificado. Entonces
\[
I_{K_n}(p)=A_n\cong X_n\cong\frac{X}{\big((1+T)^{p^n}-1)\big)}
\]
y $p\nmid h_0\iff p\nmid h_n$ para toda $n\geq 0$, donde $h_n=
\big|I_{K_n}(p)\big|$.
\end{teorema}

\begin{proof}
Se tiene que $K_{\infty}/K$ satisface la condici\'on (A) con $s=1$, es
decir, un \'unico primo ramificado. Entonces
\begin{gather*}
Y_0=\langle TX=X^{\gamma_0-1}, a_2,\ldots, a_s\rangle =TX,\\
Y_n=\gamma_n Y_0=\frac{(1+T)^{p^n}-1}{T}Y_0=\big((1+T)^{p^n}-1
\big)X.
\end{gather*}
Por tanto $X_n\cong X/Y_n=\frac{X}{\big((1+T)^{p^n}-1\big)X}$.

Si $p\nmid h_0$, entonces $X/TX=X/Y_0\cong A_0=\{1\}$, esto es,
$X/TX=0$ lo cual implica $X/\langle p, T\rangle X=0$. Por el Lema
de Nakayama, $X=0$ y $A_n=\{1\}$ para toda $n\geq 0$. $\fin$
\end{proof}

\begin{corolario}\label{C7.8.19} Para $K={\ma Q}$ y $p$ cualquier
primo, $\mu=\lambda=\gamma=0$.
\end{corolario}

\begin{proof}
Se tiene que ${\ma Q}_{\infty}/{\ma Q}$ es una extensi\'on ${\ma Z}_p
$ de ${\ma Q}$ y $h_0=\big|I_{\ma Q}\big|=1$. Puesto que
\'unicamente $p$ se ramifica, $A_n=\{1\}$ para toda $n$, esto es, $|A_n|
=p^{e_n}$, $e_n=0=\mu p^n\lambda n+\gamma$. $\fin$
\end{proof}

\begin{teorema}\label{T7.8.20} Se tiene que $\mu=0$ si y s\'olo si
si el $p$--rango de $A_n$ se mantiene acotado cuando $n\to\infty$.
\end{teorema}

\begin{proof}
Se tiene $Y_e\sim E=\Big(\bigoplus\limits_{i=1}^s \Lambda/
\langle p^{k_i}\rangle\Big)\oplus \Big(\bigoplus\limits_{j=1}^t\Lambda/
\langle f_j^{m_j}(T)\rangle\Big)$. Se tiene que $s=0$ si y solamente
si el $p$--rango de $E/\gamma_{n,e}E$ est\'a acotado.
Como antes, tenemos una sucesi\'on exacta
\[
0\longto C_n\longto Y_e/\gamma_{n,e}Y_e\longto E/\gamma_{n,e}E
\longto B_n\longto 0
\]
con $|C_n|$, $|B_n|$ acotados para toda $n$. Por lo tanto $\mu=0$
si y solamente si el $p$--rango de $Y_e/\gamma_{n,e}Y_e$ est\'a
acotado.  Se tiene $A_n\cong X_n=X/\gamma_{n,e}Y_e$ y $X/Y_e$
es finito. El resultado se sigue. $\fin$
\end{proof}

\begin{conjetura}[Iwasawa]\label{C7.8.21} Si $K$ es cualquier
campo num\'erico finito y $K_{\infty}/K$ es la extensi\'on 
${\ma Z}_p$--ciclot\'omica de $K$, entonces $\mu=0$.
\end{conjetura}

Iwasawa primero conjetur\'o esto para cualquier extensi\'on 
${\ma Z}_p$. Sin embargo \'el mismo encontr\'o contraejemplos cuando
$K_{\infty}/K$ no es la extensi\'on ciclot\'omica.
Washington y Ferrero probaron la Conjetura \ref{C7.8.21} cuando
$K/{\ma Q}$ es una extensi\'on abeliana.

Ahora bien, tenemos

\begin{teorema}\label{T7.8.22}
Se tiene que $\lambda=0$ si y solamente si el exponente de $A_n$
est\'a acotado.
\end{teorema}

\begin{proof}
{\ }

\noindent
$\Longrightarrow)$ Sea $X\sim E=
\Big(\bigoplus\limits_{i=1}^s \Lambda/
\langle p^{k_i}\rangle\Big)\oplus \Big(\bigoplus\limits_{j=1}^t
\Lambda/ \langle f_j(T)^{m_j}\rangle\Big)$. Entonces $\mu=
\sum_{i=1}^s k_i$, $\lambda=\sum_{j=1}^t m_j \gr f_j(T)$. Se tiene
\[
\lambda =0\iff X\sim \bigoplus\limits_{i=1}^s \Lambda/
\langle p^{k_i}\rangle.
\]

Sea $k_0=\max\limits_{1\leq i\leq s} k_i$, entonces $p^{k_0}X\sim
0$, es decir, $p^{k_0}X$ es finito, por lo que existe $a\in{\ma N}$ tal
que $p^a X=\{0\}$ lo cual implica que $p^a A_n\cong p^a X_n=0$,
es decir, el exponente de $A_n$ est\'a acotado.

\noindent
$\Longleftarrow)$ Sea $p^a A_n=0$ para alguna $a\in{\ma N}$ y 
toda $n$. Entonces $p^a X_n=p^a\big(X/\gamma_{n,e} Y_e\big)=0$
para $n\geq e$. Por lo tanto $p^a X\subseteq
\bigcap\limits_{n=e}^{\infty}
\gamma_{n,e}Y_e =\{0\}$, esto es, $t=0$ pues para toda
$f(T)\neq 0$ distinguido e irreducible se tiene $p^a\big(\Lambda/
\langle f(T)^m\rangle\big)\neq 0$. Por lo tanto $\lambda=0$. $\fin$
\end{proof}

\begin{conjetura}[R. Greenberg\index{conjetura de Ralph Greenberg}]
\label{C7.8.23}
Sea $K$ una campo num\'erico finito. Si $K$ es totalmente real,
es decir, $r_2=0$, entonces $\lambda=\mu=0$.
\end{conjetura}

\begin{proposicion}\label{P7.8.24}
Sea $K$ un campo num\'erico finito, $K_{\infty}/K$ una extensi\'on
${\ma Z}_p$ y $s$ el n\'umero de divisores primos de $K$ 
ramificados en $K_{\infty}/K$. Sea $L/K_{\infty}$ una $p$--extensi\'on
de $K_{\infty}$ no ramificada tal que $L/K$ es abeliana. Entonces
$\Gal(L/K)\cong {\ma Z}_p^t\oplus R$ donde $R$ es un $p$--grupo
finito y $t\leq s$.
\end{proposicion}

\begin{proof}
{\ }

\begin{window}[0,l,\xymatrix{&L\\ K_{\infty}\ar@{-}[ru]^N
\ar@{-}[d]_{\Gamma}\\ K\ar@{-}[uur]_G},{}]
Sean $G:=\Gal(L/K)$ y $N:=\Gal(L/K_{\infty})$. Sean $\pK_1,\ldots,
\pK_s$ los primos de $K$ ramificados en $K_{\infty}/K$ y sean
$T_1,\ldots, T_s$ los respectivos grupos de inercia de cada $\pK_i$
en $L/K$. Puesto que $L/K_{\infty}$ es no ramificada se tiene que
$T_i\cap N=\{0\}$ para $i=1,\ldots, s$. Por lo tanto $T_i\cong T_i N/N
$, esto es, $T_i$ es isomorfo a un subgrupo cerrado no trivial de $
\Gamma =G/N$. Se tiene que $T_i\cong \Gamma$ para $1\leq i
\leq s$.
\end{window}

Sean $T:=T_1\cdots T_s$, $E=L^T$ y $\Gal(L/E)=T$. Se tiene que 
$E/K$ es no ramificada y abeliana y por ende finita, esto es $G/T$
es un grupo finito. Ahora bien, $G/T$ es un $p$--grupo
debido a que tanto $G$ como $T$ son ${\ma Z}_p$--m\'odulos.

As{\'\i}, $T\cong {\ma Z}_p^t\oplus R_1$ con $R_1$ un $p$--grupo
finito y $t\leq s$ ya que $T=T_1\cdots T_s$ y $T_i\cong \Gamma$,
$1\leq i\leq s$. Por lo tanto $G\cong {\ma Z}_p^t\oplus R$ con $R$
un $p$--grupo finito. $\fin$
\end{proof}

\begin{corolario}\label{C7.8.25} Si $s=1$, es decir, \'unicamente hay
un primo ramificado, entonces la m\'axima $p$--extensi\'on no
ramificada $L$ de $K_{\infty}$ tal que $L/K$ es abeliana satisface 
que $L/K_{\infty}$ es finita. 
\end{corolario}

\begin{proof}
Se tiene en este caso que $G\cong {\ma Z}_p^t\oplus R$, $1\leq t\leq
1=s$ y $\Gal(K_{\infty}/K)\cong \Gamma$, por lo que $R=\Gal(L/
K_{\infty})$ es un $p$--grupo finito. $\fin$
\end{proof}

\begin{observacion}\label{O7.8.26} Ya hab{\'\i}amos probado que si
$F$ es la m\'axima extensi\'on abeliana $K$ no ramificada fuera de 
$p$, entonces $\Gal(F/K)\cong {\ma Z}_p^{r_2+1+\delta}\times (
\text{finito})$ (ver Teorema \ref{T7.6.11}).
\end{observacion}

\begin{definicion}\label{D7.8.27} Un campo num\'erico $K$ se llama
de tipo $\MC$ ({\em multiplicaci\'on compleja\index{multiplicaci\'on
compleja}}) si $K$ es totalmente imaginario que es una extensi\'on
cuadr\'atica de un campo totalmente real.
\end{definicion}

\begin{ejemplo}\label{Ej7.8.28} Si $K\nsubseteq {\ma R}$ y $K/
{\ma Q}$ es una extensi\'on abeliana, entonces para todo encaje
$\sigma\colon K\to {\ma C}$, $\sigma(K)=K\nsubseteq {\ma R}$
por lo que $K$ es totalmente imaginario. Sea $J\colon K\to K$,
$x\mapsto \overline{x}$ la conjugaci\'on compleja. Entonces $J\in
\Gal (K/{\ma Q})$ y $o(J)=2$ pues $J|_K\neq \Id_K$. Sea $K^+:=
K^{\{1,J\}}= K\cap {\ma R}$. Entonces $K$ es totalmente real pues
$K^+\subseteq {\ma R} $ y $K^+/{\ma Q}$ es de Galois. Finalmente
$[K:K^+]=|\{1,J\}| =2$. Por tanto $K$ es de tipo $\MC$.
\end{ejemplo}

\begin{ejemplo}\label{Ej7.8.29} Si $n\geq 3$, $\cic n{}$ es un campo
de tipo $\MC$.
\end{ejemplo}

\begin{proposicion}\label{P7.8.30} Si $K$ es de tipo $\MC$ y si $K^+$
es su subcampo real, esto es, $[K:K^+]=2$ y $K^+$ es totalmente
real, entonces la conjugaci\'on compleja $J$ induce un automorfismo
en $K$ el cual es independiente del encaje $K$ en ${\ma C}$. 
Adem\'as $K^+=K\cap {\ma R}$.
\end{proposicion}

\begin{proof}
Sean $\phi,\psi\colon K\to {\ma C}$ dos encaje de $K$. Notemos que 
$\phi(K)/\phi(K^+)$ es una extensi\'on cuadr\'atica y por tanto normal.
Adem\'as, puesto que $\phi(K^+)\subseteq {\ma R}$, entonces
$J(\phi(K^+))=\phi(K^+)$, es decir, $J\circ \phi(\alpha)=\phi(\alpha)$
para todo $\alpha\in K^+$. Sea $\overline{\phi}:=J\circ \phi$. Puesto
que $\phi(K)/\phi(K^+)$ es normal, se tiene $J(\phi(K))=\phi(K)$, es
decir, $\overline{\phi}(K)=\phi(K)$. Se sigue que $\phi^{-1}\circ 
\overline{\phi}$ est\'a definida en $K$. De esta forma tenemos
que $\psi^{-1}\circ \overline{\psi}$ y $\phi^{-1}\circ \overline{\phi}$
son automorfismos de $K$ que fijan a $K^+$ puesto que $K^+$
es totalmente real.

Ahora bien, puesto que $K$ es totalmente imaginario, ni $\phi^{-1}\circ
\overline{\phi}$ ni $\psi^{-1}\circ \overline{\psi}$ puede ser la identidad
y puesto que  $\Gal(K/K^+)$ es de orden $2$ se sigue que 
$\phi^{-1}\circ \overline{\phi}=\psi^{-1}\circ \overline{\psi}=J$. Por
lo tanto $J\colon K\to K$ es un automorfismo de $K$, $J\neq \Id$ y
$\Gal(K/K^+)=\{1,J\}$. Finalmente $K^+=K^{\{1,J\}}=K\cap {\ma R}$.
$\fin$
\end{proof}

\begin{observacion}\label{O7.8.31}
Dado un campo num\'erico finito, $[K:K\cap {\ma R}]$ puede ser
arbitrariamente grande. Por ejemplo, si $K={\ma Q}(\zeta_n\sqrt[n]{2})$,
se tiene que $K\cap {\ma R}={\ma Q}$ y $[K:{\ma Q}]=n$. Otro ejemplo
es  si $p\geq 3$ un n\'umero 
primo y sea $f(x)=x^{2p}+3$. Si $\omega_0,\ldots,\omega_{2p-1}$
son las ra{\'\i}ces de $f(x)$, $\omega_t=3^{1/2p} e^{((\pi+2t\pi)i)/2p}$.

\begin{window}[2,l,\xymatrix{
{\ma Q}(\omega_0)\ar@{-}[r]\ar@{-}[d]&{\ma Q}(\zeta_p,\omega_0)
\ar@{-}[d]^{\Big\} 2p}\\{\ma Q}\ar@{-}[r]&{\ma Q}(\zeta_p)=
{\ma Q}(\zeta_{2p})},{}]
Tomemos por ejemplo $\omega_0=3^{1/2p} e^{\pi i/2p}$ y sea $K:=
{\ma Q}(\omega_0)$. Entonces $[K:{\ma Q}]=2p$ y $K$ es 
totalmente imaginario pues si $\sigma\colon K\to {\ma C}$ es un
encaje, $\sigma(K)= {\ma Q}(\omega_i)$ para alguna $i$, $0\leq i\leq
2p-1$. Por otro lado, la cerradura de Galois de $K/{\ma Q}$ es $
\tilde{K}:={\ma Q}(\zeta_{2p},\omega_0)={\ma Q}(\zeta_p,\omega_0)$.
Se tiene que ${\ma Q}(\zeta_{2p},\omega_0)/{\ma Q}(\zeta_{2p})$
es una extensi\'on de Kummer c{\'\i}clica de grado $2p$. Sea $A:=
{\ma Q}(\omega_0)\cap {\ma R}$. Notemos que en general
si $[K:K\cap {\ma R}]=2$ entonces necesariamente $\Gal(K/K\cap
{\ma R})=\{1,J\}$, esto es, $K^J=K$.
\end{window}

Por tanto si tuvi\'esemos $[{\ma Q}(\omega_0):A]=2$, entonces
$A={\ma Q}(\omega_0^2)$ puesto que hay una correspondencia
entre los subcampos de ${\ma Q}(\omega_0)$ que contienen a ${\ma
Q}$ y los de ${\ma Q}(\omega_0,\zeta_{2p})$ que contienen a ${\ma
O}(\zeta_{2p})$ y esta \'ultima es c{\'\i}clica. 

En particular hay un \'unico subcampo de {\'\i}ndice $2$ en ${\ma Q}(
\omega_0)$ y por tanto este corresponde a ${\ma Q}(\omega_0^2)$.
Sin embargo $\omega_0^2=3^{1/p}e^{\pi i/p}\notin {\ma R}$. Todos
los subcampos de ${\ma Q}(\omega_0)$ son ${\ma Q}(\omega_0^j)$,
con $j=0,1,\ldots, 2p-1$ y $\omega_0^j=3^{j/2p}e^{\pi ji/2p}\in{\ma R}
\iff j=2p$. En particular 
\[
{\ma Q}(\omega_0)\cap {\ma R}={\ma Q}\quad \text{y}\quad
[{\ma Q}(\omega_0):{\ma Q}(\omega_0)\cap{\ma R}]=[{\ma Q}(
\omega_0):{\ma Q}]=2p.
\]
\end{observacion}

\begin{teorema}\label{T7.8.32}
Sea $K$ un campo num\'erico finito de tipo $\MC$ y sea $K^+:=
K\cap {\ma R}$ su subcampo real. Sean $h:=|I_K|$, $h^+=
|I_{K^+}|$. Entonces $h^+|h$.
\end{teorema}

\begin{proof}
{\ }

\begin{window}[0,r,\xymatrix{
H_{K^+}\ar@{-}[r]\ar@{-}[d]&H_{K^+}K\ar@{-}[d]\\K^+\ar@{-}[r]&K},{}]
Se tiene que $K/K^+$ es totalmente ramificado en todos los primos
arquimedianos de $K^+$. Sean $H_K$ y $H_{K^+}$ los campo de
clase de Hilbert de $K$ y $K^+$ respectivamente. Entonces
$K\cap H_{K^+}=K^+$ pues $H_{K^+}/K^+$ es no ramificada. Por
tanto $H_{K^+}K/K$ es una extensi\'on abeliana y no ramificada. Se
sigue que $H_{K^+}K\subseteq H_K$. Obtenemos que
\end{window}
\begin{align*}
h&=|I_K|=[H_K:K]=[H_K:H_{K^+}K][H_{K^+}K:K]\\
&=[H_K:H_{K^+}K][H_{K^+}:K^+]=[H_K:H_{K^+}K]|I_{K^+}|\\
&= [H_K:H_{K^+}K]h^+. \tag*{$\fin$}
\end{align*}
\end{proof}

\begin{teorema}\label{T7.8.33}
Sea $K$ un campo de tipo $\MC$ y sea $E$ el grupo de unidades
de $K$. Sea $E^+$ el grupo de unidades de $K^+$ y sea $W$ el 
grupo de ra{\'\i}ces de unidad en $K$. Entonces si definimos
$Q=[E:WE^+]$, se tiene que $Q=1$ o $2$.
\end{teorema}

\begin{proof}
Sea $\phi\colon E\to W$ definida por $\phi(\varepsilon)=\varepsilon/
\overline{\varepsilon}$. Puesto que para todo encaje $\sigma\colon K
\to{\ma C}$ se tiene que $\overline{\varepsilon^{\sigma}}=
\overline{\varepsilon}^{\sigma}$ se sigue que
\[
\big|\phi(\varepsilon)^{\sigma}\big|=\big|\big(\varepsilon/\overline{
\varepsilon}\big)^{\sigma}\big|=\big|\varepsilon^{\sigma}/\overline{
\varepsilon}^{\sigma}\big|=\big|\frac{\varepsilon^{\sigma}}{\overline{
\varepsilon^{\sigma}}}\big|=1
\]
lo cual implica que $\phi(\varepsilon)\in W$, Sea $\psi\colon E\to W/
W^2$ el mapeo inducido por $\phi$, es decir, el mapeo $\pi\circ \phi$
donde $\pi$ es la proyecci\'on natural $W\to W/W^2$
\begin{gather*}
E\stackrel{\phi}{\longto}W\stackrel{\pi}{\longto}W/W^2.\\
\intertext{Supongamos que $\varepsilon =\xi\varepsilon_1$ con
$\xi\in W$, $\varepsilon_1\in E^+$, entonces}
\phi(\varepsilon)=\phi(\xi\varepsilon_1)=\frac{\xi\varepsilon_1}{
\overline{\xi}\overline{\varepsilon}_1}=\frac{\xi}{\overline{\xi}}=
\frac{\xi}{\xi^{-1}}=\xi^2\in W^2.
\end{gather*}
En particular $\varepsilon\in\ker \psi$.

Rec{\'\i}procamente, si $\varepsilon\in\ker \psi$, entonces $\phi(
\varepsilon)=\varepsilon/\overline{\varepsilon}=\xi^2\in W^2$. 
Definimos $\varepsilon_1:=\xi^{-1}\varepsilon$. Entonces $\overline{
\varepsilon_1}=\overline{\xi^{-1}}\overline{\varepsilon}=\xi\frac{
\varepsilon}{\xi^2}=\xi^{-1}\varepsilon=\varepsilon_1$, esto es, 
$\varepsilon_1\in{\ma R}$. Por lo tanto $\ker\psi=WE^+$ y tenemos
la inyecci\'on inducida $E/WE^+\hookrightarrow W/W^2$ y $\big|
W/W^2\big|=2$ puesto que $W$ es c{\'\i}clico y $2||W|$. Se sigue
que $Q=1$ o $2$. De hecho, si $\phi(E)=W$, entonces $Q=2$ y si
$\phi(E)=W^2$, $Q=1$. $\fin$
\end{proof}

\begin{corolario}\label{C7.8.34} Sea $K=\cic n{}$ el campo 
ciclot\'omico de las $n$--ra{\'\i}ces de la unidad y $n\geq 3$. 
Entonces $Q=\begin{cases} 1&\text{si $n=p^m$}\\2&\text{en 
otro caso}\end{cases}$.
\end{corolario}

\begin{proof}
Ver \cite[Corolario 4.13, p\'agina 40]{Was97}. $\fin$
\end{proof}

\begin{teorema}\label{T7.8.35} Sea $C$ el grupo de clases de
ideales de $\cic n{}$, $n\geq 3$ y $C^+$ el grupo de clases de ideales
de $\cic n{}^+$. Entonces el mapeo natural $C^+\to C$ es inyectivo.
\end{teorema}

\begin{proof}
Ver \cite[Teorema 4.14, p\'ag. 40]{Was97}. $\fin$
\end{proof}

\begin{observacion}\label{O7.8.36} El Teorema \ref{T7.8.35}
no se cumple para campos de tipo $\MC$ en general. Por ejemplo,
si $K={\ma Q}(\sqrt{10},\sqrt{-2})$, entonces $K$ es de tipo $\MC$ y
$K^+={\ma Q}(\sqrt{10})$. El ideal $\pK:=\langle 2,\sqrt{10}\rangle$
en ${\ma Q}(\sqrt{10})$ no es principal pues si $\langle 2,\sqrt{10}
\rangle=\langle \alpha\rangle$ en ${\ma Q}(\sqrt{10})$, $N
\alpha=N\pK=2$ pero si $\alpha=a+b\sqrt{10}$, $N\alpha=a^2-10 b^2
\neq \pm 2$ ya que si fuese posible, $a^2\equiv \pm 2 \bmod 5\equiv
2,3\bmod 5$. Sin embargo los residuos m\'odulo $5$ son $0$, $1$,
y $4$. 

Por otro lado, sea $C^{+}=I_{{\ma Q}(\sqrt{10})}
\stackrel{\phi}{\longto}C=I_{{\ma Q}(\sqrt{10},\sqrt{-2})}$ y 
$\pK=\langle 2, \sqrt{10}\rangle$, $\phi(\pK)=\pL^2$ pues como $2$
es ramificado en ${\ma Q}(\sqrt{10})/{\ma Q}$, ${\ma Q}(\sqrt{-5})/
{\ma Q}$ y ${\ma Q}(\sqrt{-2})/{\ma Q}$, $2$ es totalmente
ramificado en $K/{\ma Q}$. Sea $\alpha:=-\sqrt{-2}\in K$. Entonces
$2=\alpha (\sqrt{-2})$, $\sqrt{10}=\alpha(-\sqrt{-5})$ por lo que
$\pK=\langle \alpha\rangle$ en $C$ y en particular $\overline{\pK}
\in \ker \phi$ y $\overline{\pK}\neq 1$.
\end{observacion}

Por otro lado se sabe en general que para $K$ de tipo $\MC$, $
C^+\stackrel{\phi}{\longto}C$ satisface que $|\ker\phi|=1$ o $2$.

Sea $K_{\infty}/K$ una extensi\'on ${\ma Z}_p$ tal que $K_n$ es de
tipo $\MC$ para toda $n$. Entonces $K_{\infty}^+/K^+$ es una
extensi\'on ${\ma Z}_p$. Notemos que $r_2=0$ para $K^+$ y por
tanto $K_{\infty}^+/K^+$ ser\'a la extensi\'on ${\ma Z}_p$--ciclot\'omica
en caso de que la Conjetura de Leopoldt (Conjetura \ref{Co7.6.9})
se cumpla y por lo tanto, en este caso, $K_{\infty}/K$ es la 
extensi\'on ${\ma Z}_p$--ciclot\'omica pues $K_{\infty}=KK_{\infty}^+
=KK^+{\ma Q}_{\infty}=K{\ma Q}_{\infty}$.

En general, si $K$ es de tipo $\MC$, se tiene que $J\in\Gal(K/K^+)$,
$J\neq \Id$ y en particular $J$ es un automorfismo de $K$. Se tiene
que $J$ act\'ua en varios grupos y m\'odulos asociados a $K$. Por
ejemplo, $J$ act\'ua en $I_K$ o en $A:=I_K(p)$.
Si $A$ es un grupo abeliano tal que $J$ act\'ua en $A$, definimos
\[
A^+:=\{a\in A\mid J(a)=a\}\quad \text{y}\quad A^-:=\{a\in A\mid J(a)
=-a\}.
\]

Si $2$ es una unidad en $A$, lo cual, en el caso finito, significa que
$2$ no divide al orden de $A$, se tiene que
\[
A=A^++A^-\quad \text{y}\quad A^+\cap A^-=\{0\}.
\]
En efecto, si $x\in A$, $x=\frac{x+\overline{x}}{2}+\frac{x-\overline{x}}
{2}=y+z$, donde $y=\frac{x+\overline{x}}{2}$, $z=\frac{x-\overline{x}}
{2}$. Ahora bien, $\overline{y}=\frac{\overline{x}+\overline{\overline{x}}}
{2}=\frac{x+\overline{x}}{2}=y$ y $\overline{z}=\frac{\overline{x}-
\overline{\overline{x}}}{2}=-\frac{x-\overline{x}}{2}=-z$, es decir, $y\in
A^+$ y $z\in A^-$ lo cual implica que $A^++A^-=A$.

Ahora bien, si $x\in A^+\cap A^-$, entonces $\overline{x}=x=-x$,
es decir, $2x=0$ lo cual implica $x=0$ pues $2$ es unidad en $A$.
Por lo tanto $A^+\cap A^-=\{0\}$. En consecuencia tenemos $A=
A^+\oplus A^-$.

Por ejemplo, si $p$ es un n\'umero primo, $p>2$ y $A$ es el 
$p$--subgrupo de Sylow del grupo de clases de un campo $K$ de
tipo $\MC$, tenemos la situaci\'on anterior.

Ahora consideremos $K_{\infty}/K$ una extensi\'on ${\ma Z}_p$ de
tipo $\MC$, $A_n=I_{K_n}(p)$ con $p>2$. Entonces $A_n=
A_n^+\oplus A_n^-$, $X_n=X_n^+\oplus X_n^-$.
\[
\xymatrix{
&&L_n\\K_0\ar@{-}[r]\ar@{-}[d]&K_n\ar@{-}[d]^{\{1,J\}}\ar@{-}[ru]^
{X_n}\\
K_0^+\ar@{-}[r]&K_n^+
}
\]

\begin{window}[0,r,\xymatrix{
K\ar@{-}[r]^{\Gamma}\ar@{-}[d]_J&K_{\infty}\ar@{-}[d]_J\\
K^+\ar@{-}[r]_{\Gamma}&K_{\infty}^+},{}]
Puesto que la acci\'on de $\gamma_0$ en $\Gamma=\Gal(K_{\infty}/
K)$ conmuta con la acci\'on de $J$, se tiene que $X=X^+\oplus X^-$
con $X^+:=\lim\limits_{\substack{\longleftarrow\\n}} X_n^+$ y 
$X^-:=\lim\limits_{\substack{\longleftarrow\\n}} X_n^-$.
Obtenemos como antes $A_n^{\pm}\cong X_n^{\pm}\cong
\frac{X^{\pm}}{\gamma_{n,e}Y_e^{\pm}}$. Adem\'as tenemos que $
h_n^+|h_n$ y si definimos $h_n^-:=\frac{h_n}{h_n^+}$, se tiene que
si $p^{e_n^{\pm}}$ es la potencia exacta de $p$ que divide a $h_n^{
\pm}$, entonces 
\end{window}
\begin{gather*}
e_n=e_n^++e_n^-\quad \text{y}\quad e_n^{\pm}=\mu^{\pm}+
\lambda^{\pm}n+\gamma^{\pm}\quad\text{con}\\
\mu=\mu^++\mu^-,\quad \lambda=\lambda^++\lambda^-\quad
\text{y}\quad \gamma=\gamma^++\gamma^-.
\end{gather*}

Como obtuvimos antes, se sigue que $\mu^{\pm}=0\iff$ el $p$--rango
de $A_n^{\pm}$ est\'a acotado.

En el caso $p=2$, si $x\in A^+\cap 
A^-$, entonces $\overline{x}=x=-x$, esto es $2x=0$, pero no podemos
concluir que $x=0$ y por ende que $A^+\cap A^-=\{0\}$. Por ejemplo,
se tiene que $h({\ma Q}(\sqrt{-5}))=2$, esto es, $A\cong {\ma Z}/2
{\ma Z}$. Se tiene que $\{I,J\}=\Gal({\ma Q}(\sqrt{-5})/{\ma Q})$ y
$A^+=A^-=A$.

Si $X$ es un clase de ideales de $K^+$ con $X^2=(1)$, entonces
$X\in A^+\cap A^-$ pues $\overline{X}=X=X^{-1}$.

Sean $p=2$ y $K$ un campo de clase $\MC$. Definimos $A=I_K
(2)$, $A^-=\{a\in A\mid \overline{a}=a^{-1}\}$ y $A^+=I_{K^+}
(2)$. Sea $N\colon A\to A^+$, $N=1+J$ la norma de $K$ en $K^+$
pues $\Gal(K/K^+)=\{1, J\}$.

Usando el mapeo de Reciprocidad de Artin, tenemos el diagrama
conmutativo (Teorema \ref{T7.8.5})
\[
\xymatrix{
I_{K}=D_K/P_K\ar[r]^{\cong}_{\text{Artin}}\ar[d]^{\text{Norma}}&
\Gal(H_K/K)\ar[d]^{\rest}\\
I_{K^+}=D_{K^+}/P_{K^+}\ar[r]^{\cong}_{\text{Artin}}&
\Gal(H_{K^+}/K^+)
}
\]
donde $H_K$ y $H_{K^+}$ denotan los campos de clase de Hilbert de
$K$ y $K^+$ respectivamente.

Puesto que $K\cap H_{K^+}=K^+$ debido a que los primos 
arquimedianos son ramificados el mapeo de restricci\'on $\rest$ es
suprayectivo de donde se sigue que el mapeo norma $N\colon
I_K=D_K/P_K\to I_{K^+}=D_{K^+}/P_{K^+}$ es sobre.
\[
\xymatrix{
&H_K\ar@{-}[d]\\ 
H_{K^+}\ar@{-}[r]\ar@{-}[d]|{=}&H_{K^+}K\ar@{-}[d]|{=}\\
K^+\ar@{-}[r]&K
}
\]

Por tanto $N\colon A\to A^+$ es suprayectivo y 
\begin{gather*}
\ker N=\{a\mid (1+J)a=a^{1+J}=a\overline{a}=1\}=\{a\in A\mid
\overline{a}=a^{-1}\}=A^-\\
\intertext{por lo que tenemos la sucesi\'on exacta}
0\longto A^-\longto A\longto A^+=I_{K^+}(2)\longto 0.
\end{gather*}

Como antes se sigue que $\mu=\mu^++\mu^-$, $\lambda=\lambda^+
+\lambda^-$, $\gamma=\gamma^++\gamma^-$, $\mu^+=0$ si
y solamente si el $2$--rango de $I_{K_n^+}(2)$ est\'a acotado
y $\mu^-=0$ si y solamente si el $2$--rango de $A_n^-$ est\'a
acotado.

Ahora bien, cuando $p>2$, se tiene que si $A=I_{K}(p)$, $
A^+=\{a\in A\mid a^J=a\}$, entonces $A^+\cong I_{K^+}(p)$.
M\'as generalmente, tenemos:

\begin{proposicion}\label{P7.8.37} Sea $F$ un campo num\'erico 
finito y sea $K$ una extensi\'on de Galois de $F$ de grado $d$. Sea
$\ell$ un n\'umero primo que no divide a $d$. El homomorfismo 
natural $I_F(\ell)\to I_K(\ell)$ es inyectivo, la norma
$N_{K/F}\colon I_K(\ell)\to I_F(\ell)$ es suprayectivo
y las siguientes condiciones son equivalentes
\l
\item $I_F(\ell)=I_K(\ell)$.
\item $\ell$--rango de $I_F(\ell)=\ell$--rango de $I_K
(\ell)$.
\item La norma $N_{K/F}\colon I_K(\ell)\to I_F(\ell)$
es un isomorfismo.
\end{list}
\end{proposicion}

\begin{proof}
Sea $N=N_{K/F}$. Sea ${\eu a}\in I_F(\ell)$ el cual es 
principal en $I_K$: ${\eu a}=\langle \alpha\rangle$, $\alpha\in
K$. Tomando la norma $N{\eu a}={\eu a}^d=\langle N\alpha\rangle=
\langle \beta\rangle$, $\beta\in F$ principal. Por otro lado, existe
$n\in{\ma N}$ tal que ${\eu a}^{\ell^n}=\langle \gamma\rangle$ con
$\gamma\in F$ pues ${\eu a}\in I_F(\ell)$. Puesto que
$\mcd (d,\ell^n)=1$, existen $a,b\in{\ma Z}$ tales que $1=ad+b\ell^n$
y ${\eu a}={\eu a}^{da}{\eu a}^{\ell^nb}=\langle \beta^a\gamma^b
\rangle$ con $\beta^a\gamma^b\in F$ por lo que ${\eu a}$ es
principal y $I_F(\ell)\to I_K(\ell)$ es 1-1.

Se tiene que $I_F(\ell)\supseteq NI_K(\ell)
\supseteq NI_F(\ell)=I_F(\ell)^d=I_F(\ell)$
esta igualdad debido a que $d$ y $\ell$ son primos relativos.
Se sigue que la norma es suprayectiva.

Probemos ahora la \'ultima parte.

\noindent
\underline{({\sc{i}}) $\Rightarrow$ ({\sc{ii}})}. Es inmediato

\noindent
\underline{({\sc{ii}}) $\Rightarrow$ ({\sc{iii}})}. Puesto que el
$\ell$--rango de $I_F(\ell)=\ell$--rango de $
I_K(\ell)$, tenemos que $\frac{I_K(\ell)}{\ell
I_K(\ell)}=\frac{I_F(\ell)}{\ell I_F(\ell)}$.

Se tiene para un $\ell$--grupo abeliano finito $A$ que $\big|\{x\in A
\mid x^{\ell}=1\}\big|=\big| A/A^{\ell}\big|$ de donde obtenemos que
si ${\eu a}\in I_K(\ell)$ con ${\eu a}\in \ker N$ y $o({\eu a})=
\ell^n$ para $n\geq 1$, entonces $o\big({\eu a}^{\ell^{n-1}}\big)=\ell$ y
${\eu a}^{\ell^{n-1}}\in\ker N$.  Puesto que $\big|\frac{I_K(\ell)}{
\ell I_K(\ell)}\big|=\big|\frac{I_F(\ell)}{\ell I_F(\ell)}
\big|$, se tiene que $\{x\in I_K(\ell)\mid x^{\ell}=1\}=\{
x\in I_F(\ell)\mid x^{\ell}=1\}$.

Puesto que ${\eu a}^{\ell^{n-1}}\in \{x\in I_K(\ell)\mid x^{\ell}=1
\}$, ${\eu a}^{\ell^{n-1}}\in I_F(\ell)$. Por otro lado $1=
N{\eu a}^{\ell^{n-1}}={\eu a}^{\ell^{n-1}d}$ lo cual contradice que
$o\big({\eu a}^{\ell^{n-1}}\big)=o\big({\eu a}^{\ell^{n-1}d}\big)=\ell$.
Por tanto $N$ es 1-1 y puesto que es sobre, $N$ es un isomorfismo.

\noindent
\underline{({\sc{iii}}) $\Rightarrow$ ({\sc{i}})}. Como la conorma
es 1-1, podemos poner $I_F(\ell)\subseteq I_K(\ell)$
y puesto que la norma es un isomorfismo, ambos grupos son
del mismo orden de donde se sigue la igualdad $I_F(\ell)=
I_K(\ell)$. $\fin$
\end{proof}

En  particular obtenemos

\begin{proposicion}\label{P7.8.38} Sea $K$ un campo de tipo $\MC$
y sea $p$ un n\'umero primo. Sea $I_K(p)$ el $p$--subgrupo
de Sylow del grupo de clases de ideales $I_K$ de $K$.
Entonces
\l
\item El mapeo natural $I_{K^+}(p)\to I_K(p)$ es 
1-1 si $p>2$ y tiene n\'ucleo de orden $1$ o $2$ si $p=2$.

\item La norma $N=N_{K/K^+}\colon I_{K}(p)\to I_{K^+}
(p)$ es suprayectiva.
\end{list}
\end{proposicion}

\begin{proof}
{\ }

\l
\item Se tiene $[K:K^+]=2$, por lo que si $p>2$, $I_{K^+}(p)
\to I_{K}(p)$ es inyectiva por la Proposici\'on \ref{P7.8.37}.

Sea $p=2$ y sea ${\eu a}$ un ideal en $K^+$ tal que ${\eu a}$ es
principal en $K$: ${\eu a}=\langle \alpha\rangle$, $\alpha\in K$.
Entonces $\overline{\eu a}=\langle\overline{\alpha}\rangle ={\eu a}=
\langle \alpha\rangle$, por lo que $\frac{\alpha}{\overline{\alpha}}$
es una unidad en ${\cal O}_K$.

Si $\sigma$ es un encaje de $K$ en ${\ma C}$, entonces, como $K$
es de tipo $\MC$, se tiene $\sigma\big(\frac{\alpha}{\overline{\alpha}}
\big)=\frac{\sigma\alpha}{\sigma(\overline{\alpha})}=\frac{\sigma
(\alpha)}{\overline{\sigma(\alpha)}}$, esto es, $\big|\frac{\sigma\alpha}
{\overline{\sigma\alpha}}\big|=1$ de donde se sigue que $\frac{
\alpha}{\overline{\alpha}}$ es una ra{\'\i}z de unidad.

Si en $K$ tenemos ${\eu a}=\langle \beta\rangle$, esto es, $\langle
\alpha\rangle=\langle\beta\rangle$, entonces existe $u\in E_K$,
unidades de $K$, tal que $\alpha=\beta u$. Por lo tanto $\frac{
\alpha}{\overline{\alpha}}=\frac{\beta}{\overline{\beta}}=\frac{u}{
\overline{u}}$. Sea $\varphi\colon E_K\to W_K=$ ra{\'\i}ces de unidad
en $K$, $u\mapsto \frac{u}{\overline{u}}$. Entonces se tiene
$\frac{\alpha}{\overline{\alpha}}=\frac{\beta}{\overline{\beta}}\bmod
\varphi(E_K)$.

Sea
\begin{eqnarray*}
\phi\colon \ker \big(I_{K^+}\to I_{K}\big)&\longto&
W/\varphi(E_K)\\
{\eu a}&\longmapsto& \frobeniusbinom{\alpha}{\overline{\alpha}}.
\end{eqnarray*}
Veamos que $\phi$ es inyectiva. De hecho, si $\phi({\eu a})=
\frobeniusbinom{\alpha}{\overline{\alpha}}=[1]$, 
entonces $\frac{\alpha}
{\overline{\alpha}}=\frac{u}{\overline{u}}$ con $u\in E_K$ y por lo
tanto $\frac{\alpha}{u}=\frac{\overline{\alpha}}{\overline{u}}=
\overline{\big(\frac{\alpha}{u}\big)}$, esto es $\frac{\alpha}{u}\in K^+$
y ${\eu a}=\langle\alpha\rangle=\langle \alpha/u\rangle$ por tanto
$\overline{{\eu a}}=(1)$.
As{\'\i} $\phi$ es 1-1. Se sigue que $\big|\ker \big(I_{K^+}\to
 I_{K}\big)\big|\leq \big|W_K/\varphi(E_K)\big|\leq [W_K:
 W_K^2]=2$ pues $W_K/\varphi(E_K)$ se puede encajar en
 $W_K/W_K^2$ el cual es de orden $2$ debido a que $W_K$ es
 c{\'\i}clico de orden par.
 
\item Si $H_K$ y $H_{K^+}$ son los campos de clase de Hilbert, 
entonces tenemos el diagrama conmutativo:
\[
\xymatrix{
I_K=D_K/P_K\ar[r]^{\cong}_{\text{Artin}}\ar[d]_{N}&
\Gal(H_K/K)\ar@{>>}[d]^{\rest}\\
I_{K^+}=D_{K^+}/P_{K^+}\ar[r]^{\cong}_{\text{Artin}}&
\Gal(H_{K^+}/K^+)
}
\]
de donde se sigue que $N$ es suprayectiva. $\fin$
\end{list}
\end{proof}

\begin{proposicion}\label{P7.8.39} Si $p$ es un n\'umero primo, $p>2
$, y $K$ es un campo num\'erico finito de tipo $\MC$, entonces si $
A=I_K(p)$, se tiene $A^+\cong I_{K^+}(p)$.
\end{proposicion}

\begin{proof}
Se tiene que $A^+=\{a\in A\mid \overline{a}=a\}=\{a\in A\mid 
a^{J-1}=1\}$. Sea $N\colon I_K(p)\to I_{K^+}(p)$,
$N=1+J$. Entonces $N$ es suprayectivo (Proposicion \ref{P7.8.38}) y
\begin{align*}
\ker N&=\{x\in I_K(p)\mid Nx=x^{J+1}=\overline{x}x=1\}=\\
&=\{x\in I_K(p)\mid \overline{x}=x^{-1}\}=A^-.
\end{align*}

Por tanto tenemos la sucesi\'on exacta
\begin{gather*}
0\longto A^-\longto A\longto I_{K^+}(p)\longto 0.\\
\intertext{Por otro lado, puesto que $A=A^+\oplus A^-$ se tiene}
A^+\cong A/A^-\cong I_{K^+}(p). \tag*{$\fin$}
\end{gather*}
\end{proof}

Resumiendo los resultados anteriores, tenemos:
para $A=I_K(p)$ donde $K$ es un campo num\'erico finito
de clase $\MC$, se tiene
\[
A^-=\{x\in A\mid \overline{x}=x^{-1}\}\quad \text{y}\quad
A^+=\begin{cases}
\{x\in A\mid \overline{x}=x\}\cong I_{K^+}(p)&\text{si $p>2$}\\
I_{K^+}(2)&\text{si $p=2$}
\end{cases}.
\]

Si $p>2$, $A\cong A^+\oplus A^-$ y si $p=2$, entonces  se tiene la
sucesi\'on exacta $1\longto A^-\longto A\longto A^+\longto 1$.

\begin{teorema}\label{T7.8.40} Sean $K$ un campo num\'erico
finito de $\MC$, $p$ un n\'umero primo impar. Supongamos que
$\zeta_p\in K$. Sea $A=I_K(p)$. Entonces $\ran_p A^+\leq
\ran_p A^-+1$.

Sea $W_K$ el grupo de las ra{\'\i}ces de unidad en $K$. Si $K(W_K^{
1/p})/K$ es ramificada, entonces $\ran_p A^+\leq \ran_p A^-$.
\end{teorema}

\begin{proof}
Se $L$ la m\'axima extensi\'on abeliana no ramificada de
$K$ de exponente $p$ y sea $G:=\Gal(L/K)$. Entonces, por teor{\'\i}a
de campos de clase, $G\cong I_K/I_K^p\cong {}_p
I_K$ y sea $G^+:=I_{K^+}/I_{K^+}^p\cong {}_p
I_{K^+}$.

Puesto que $\zeta_p\in K$, $L/K$ es una extensi\'on de Kummer. 
Sea $B$ el grupo $(K^{\ast})^p\subseteq B\subseteq K^{\ast}$ tal 
que $L=K(B^{1/p})$. Se tiene el mapeo de Kummer
\[
\xymatrix{
&L\ar@{-}[dl]_{{\cal G}}\ar@{-}[d]^G\\ K^+\ar@{-}[r]_{\{1,J\}}&K}
\qquad
\begin{array}{rcl}{\ }\\{\ }\\
G\times B/(K^{\ast})^p&\longto & {}_pW_K=\{\xi\in{\ma C}^{\ast}\mid
\xi^p=1\}\\
\big(\sigma,\overline{b}\big)&\longmapsto& \frac{\sigma\big(b^{1/p}
\big)}{b^{1/p}}.
\end{array}
\]

Se tiene que $L$ es Galois sobre $K^+$ y si ${\cal G}=\Gal(L/K^+)$,
entonces $1\to G\to {\cal G}\to \{1,J\}\to 1$ es exacta. Es decir,
$\{1,J\}$ act\'ua en $G$ por medio de conjugaci\'on. Puesto que
$p$ es impar, se tiene $G\cong G^+\times G^-$.

Sea $V:=B/(K^{\ast})^p$. Entonces $J$ act\'ua en $V$ y nuevamente
tenemos que $V\cong V^+\times V^-$. Sea $\sigma\in G^+$, esto es,
$J\circ \sigma =\sigma$, es decir, $\overline{\sigma(b^{1/p})}=
\sigma(b^{1/p})$. Ahora bien, consideremos $R:=\{[b]\in V\mid
\big(\sigma,[b]\big)=1 \ \forall\ \sigma\in G^+\}$. Entonces $\overline{b}
\in R\iff \frac{\sigma(b^{1/p})}{b^{1/p}}=1$ para toda $\sigma \in G^+
\iff \sigma(b^{1/p})=b^{1/p}\ \forall\ \sigma \in G^+\iff b^{1/p}=\sigma(
b^{1/p})=\overline{\sigma(b^{1/p})}=\overline{b^{1/p}}$.

Se sigue que $R=V^+$ y por tanto $G^+\times V/V^+\cong G^+\times
V^-\to {}_pW_K$ es un mapeo bilineal no degenerado. Se sigue que
$G^+\cong V^-$.

En particular $\ran_p G^+=\ran_pA^+=\ran_pV^-$. Se tiene que si
$b\in B$, puesto que $L=K(B^{1/p})$ y $L/K$ es no ramificada, 
entonces $K(b^{1/p})/K$ es no ramificada.

Sea ${\eu a}:=\langle b^{1/p}\rangle$ en $K(b^{1/p})$, esto es, 
${\eu a}^p=\langle b\rangle$. Sea $\langle b\rangle= \pK_1^{s_1}\cdots
\pK_r^{s_r}$ con $\pK_1,\ldots, \pK_r$ divisores primos de $K$ y
sean $\pK_i=\pL_{i1}\cdots \pL_{ir}$. Entonces como $K(b^{1/p})/K$
es no ramificada, $p|s_1,\ldots, p|s_r$ y $\langle b\rangle ={\eu B}^p$
para un ideal ${\eu B}$ de $K$.

El mapeo $\begin{array}{rlc}B&\longto&I_K\\ b&\longmapsto&
\overline{\eu B}\end{array}$ se puede levantar a un homomorfismo a
$V$:
$V=B/K^{\ast p}\stackrel{\varphi}{\longto} {}_pI_K=\{x\in I_K\mid x^p=1\}$ 
el cual induce a su vez un homomorfismo $\varphi^-
\colon V^-\longto {}_pI_K^-$ ya que $\varphi\circ J=
J\circ \varphi$.
Ahora bien probaremos que existe un mapeo inyectivo $\ker \varphi^-
\longto \big(E_K/E_K^p\big)^-$. Sea $b\in \ker \varphi^-$, ${\eu B}^p
=\langle b\rangle$, $\varphi(b)=\overline{\eu B}=\langle 1\rangle$, es
decir, ${\eu B}=\langle a\rangle$ por lo que $\langle a^p\rangle=
\langle b\rangle$. Entonces existe $u\in E_K$ tal que $b=a^pu$. 
Puesto que $b\in V^-$, $a^{-p}u^{-1}=b^{-1}=\overline{b} =
\overline{a^p}\overline{u}=\overline{a}^p\overline{u}$. Por tanto $
\overline{u}=u^{-1}$ por lo que $u\in E_K^-$. Entonces se tiene el
mapeo $\ker \varphi^- \longto  \big(E_K/E_K^p
\big)^-$, $b\mapsto u$.

Se tiene $[E_K:E_{K^+}W_K]=1$ o $2$, lo cual implica que $\frac{
E_K}{E_K^p}\cong \frac{(E_KW_K)}{(E_KW_K)^p}$ debido a que
$p>2$. De hecho si $C$ y $D$ son grupos abelianos $D<C$ y $[C:
D]=2$ y si $p$ es un n\'umero primo, $p>2$, entonces el mapeo
natural $D/D^p\to C/C^p$ es un isomorfismo; sea $C=C\uplus  \xi D$
y $D\stackrel{\theta}{\longto} C$ el encaje natural. Por tanto $\theta(
D^p)\subseteq C^p$ y $\theta$ induce $D/D^p\stackrel{\tilde{\theta}}
{\longto}C/C^p$. Para $x=\xi y$, $y\in D$ entonces $x\equiv xx^{-p}
\bmod C^p\equiv \xi y\xi^{-p}y^{-p}=\big(\xi^2\big)^{-(p-1)/2} y^{1-p}\in D
$ y en particular $\tilde{\theta}\big(\big(\xi^2\big)^{-(p-1)/2} y^{1-p}\big)
=x\bmod C^p$ y $\tilde{\theta}$ es un isomorfismo.

Volviendo a nuestra demostraci\'on, se tiene la sucesi\'on exacta
\[
1\longto \frac{W_K}{W_K^p}\longto \frac{E_{K^+}W_K}{(E_{K^+}W_K
)^p}\longto \text{subgrupo de\ }\Big(\frac{E_{K^+}}{E_{K^+}^p}\Big)
\longto 1.
\]
Puesto que $p>2$ y $\Big(\frac{E_{K^+}}{E_{K^+}^p}\Big)^+=
\frac{E_{K^+}}{E_{K^+}^p}$ tomando la parte $(-)$ en la sucesi\'on
exacta se tiene que 
\[
{\ma Z}/p{\ma Z}\cong \frac{W_K}{W_K^p}\cong \Big(\frac{W_K}{
W_K^p}\Big)^-=\Big(\frac{E_{K^+}W_K}{\big(E_{K^+}W_K\big)^p}
\Big)^-.
\]

Por lo tanto $G^+\cong V^-$, $V^-\stackrel{\varphi^-}{\longto} {}_p
I_K^{-1}$ y 
\[
\ran_p\big(\ker \varphi^{-}\big)\leq \ran_p
\big(E_K/E_K^p\big)^-=1.
\]
Resumiendo
\begin{align*}
\ran_p A^+&=\ran_p G^+=\ran_p V^-\leq\\
&\leq  \ran_p A^-+\ran_p\ker 
\varphi^-=\ran_p A^-+1.
\end{align*}

Finalmente, si $K\big(W_K^{1/p}\big)$ es ramificada, entonces $B$
no contiene a $W_K$ y por lo tanto
\begin{eqnarray*}
\ker\varphi^-&\stackrel{\lambda}{\longto}&\big(E_K/E_K^p\big)^-=
\big(W_K/W_K^p\big)^-={\ma Z}/p{\ma Z}\\
b&\longmapsto& u
\end{eqnarray*}
por lo que si $b\in\ker \varphi^-$, entonces $u=1$ y $b=1$, esto es,
$\lambda$ no puede ser suprayectiva de donde $\lambda=1$,
$\ker \varphi^-=\{1\}$ y $\ran_p A^+\leq \ran_p A^-$. $\fin$
\end{proof} 

\begin{corolario}\label{C7.8.41}
Sean $p>0$ y $K=\cic p{}$. Sea $h_p$ el n\'umero de clase de $K$.
Se tiene que si $p\nmid h_p^-$ entonces $p\nmid h_p^+$ y 
consecuentemente $p\nmid h_p$.
\end{corolario}

\begin{proof} Se tiene que $W_K=\langle \zeta_{2p}\rangle$ y $W_K^{
1/p}=\langle \zeta_{2p^2}\rangle$, $K(W_K^{1/p})=\cic {2p}2
=\cic p2$, esto es, $K(W_K^{1/p})/K$ es ramificada por lo que
$\ran_p A^+\leq \ran_pA$ donde $A=I_{K}$. Por hip\'otesis, 
tenemos que $\ran_p A^-=0$ pues $p\nmid h_p^-$ lo cual implica
que $\ran_p A^+=0$ y por ende $p\nmid h_p^+$. $\fin$
\end{proof}

\begin{observacion}\label{O7.8.42} Se tiene que si $p=37$, $h_{37}^-
=37$ por lo que $37| h_{37}^-$.

Similarmente, se tiene que $h_{69}^-=69$ aunque en este caso 
$69$ no es primo.
\end{observacion}

Para $p=2$, se tiene:

\begin{teorema}\label{T7.8.43} Sea $K$ un campo num\'erico finito
de tipo $\MC$ y sean $A_K$, $A_{K^+}$ los $2$--subgrupos de Sylow de
los grupos de clases de $K$ y $K^+$ respectivamente. Entonces
\[
\ran_2 A_{K^+}\leq 1+\ran_2 A_K^-.
\]
\end{teorema}

\begin{proof}
Sea $\varphi\colon A_{K^+}\to A_K$ el mapeo natural. Se tiene que
$|\ker \varphi|=1$ o $2$. Ahora si $x\in {}_2\varphi(A_{K^+})=\{
x\in\varphi(A_{K^+})\mid x^2=1\}$, entonces $x\in {}_2A_K^-=\{x\in
A_K^-\mid x^2=1\}$ pues si $x^2=1$, se tiene $x=x^{-1}$ y puesto
que $x\in\varphi(A_{K^+})$, $\overline{x}=x$ lo cual implica que
$\overline{x}=x=x^{-1}$, es decir, $x\in A_K^-$.

De esta forma obtenemos $\ran_2 {}_2\big(\varphi(A_{K^+})\big)
\leq \ran_2 {}_2 A_K^-$. Por lo tanto
 $|\ker \varphi|=\frac{\big|{}_2A_{K^+}
\big|}{\big|{}_2\varphi(A_{K^+})\big|}=1$ o $2$, por tanto $\ran_2
\varphi(A_{K^+})\geq \ran_2 A_{K^+}-1$. Finalmente obtenemos 
\[
\ran_2 A_{K^+}\leq \ran_2 \varphi(A_{K^+})+1\leq \ran_2 A_K^- 
+1. \tag*{$\fin$}
\]
\end{proof}

Sea ahora $K_{\infty}/K$ una extensi\'on ${\ma Z}_p$ tal que $K_n$
es de tipo $\MC$. Se tiene $\mu=\mu^++\mu^-$ y $\mu=0$ si y
solamente si $\ran_p A_n$ est\'a acotado.

\begin{corolario}\label{C7.8.44} Se tiene que si $\mu^-=0$, entonces
$\mu^+=0$ y por ende $\mu=0$. En consecuencia, para cualquier
n\'umero primo $p$, $\mu=0\iff \mu^-=0$.
\end{corolario}

\begin{proof} Si $\mu^-=0$, el $p$--rango de $A_n^-$ est\'a acotado
y por lo tanto el $p$--rango de $A_n^-$ est\'a acotado. Se sigue
que $\mu^+=0$ y $\mu=\mu^++\mu^-=0+0=0$. $\fin$
\end{proof}

\begin{proposicion}\label{P7.8.45} Sea $K_{\infty}/K$ una extensi\'on
${\ma Z}_p$ tal que $\mu=0$. Entonces si $L$ es la m\'axima 
$p$--extensi\'on abeliana de $K_{\infty}$ no ramificada, se tiene
\[
X:=\Gal(L/K_{\infty})\cong \lim_{\substack{\longleftarrow\\ n}} A_n
\cong {\ma Z}_p^{\lambda}\oplus (p\text{--grupo finito}).
\]
\end{proposicion}

\begin{proof}
 Se tiene que $X\sim E=\bigoplus\limits_{j=1}^s \Lambda/\langle
f_j(T)^{m_j}\rangle$, $\lambda=\sum\limits_{j=1}^s m_j\gr f_j(T)$.

Por el algoritmo de la divisi\'on, obtenemos que $E\cong {\ma Z}_p^{
\lambda}$ como ${\ma Z}_p$--m\'odulo. Puesto que $X$ es ${\ma
Z}_p$--m\'odulo finitamente generado pues $E$ lo es, por
el teorema sobre la estructura de m\'odulos sobre un dominio de
ideales principales (pues ${\ma Z}_p$ lo es) se sigue
que el ${\ma Z}_p$--rango de $X$ es $\lambda$. $\fin$
\end{proof}

Ahora consideremos una situaci\'on dual. Sea el mapeo natural
$\begin{array}{rcl} A_n&\stackrel{\phi_{n,m}}{\longto}& A_m\\
\overline{{\eu a}}&\longmapsto& \overline{{\eu a}}\end{array}$
para $m\geq n$. Se tiene 
\begin{gather*}
\Gal(L_n/K_n)=X_n\cong X/\gamma_{n,e}Y_e \xrightarrow[\text{Artin}]
{\cong} A_n\quad \text{para}\quad n\geq e.\\
\xymatrix{
&&L\ar@{-}[dd]^{X}\\ & L_n\ar@{-}[d]^{X_n}\\
K_0\ar@{-}[r]&K_n\ar@{-}[r]&K_{\infty}
}
\end{gather*}
\begin{eqnarray*}
A_n&\xrightarrow[]{\cong}&X_n=\Gal(L_n/K_n)\\
C&\longmapsto&\xbinom{L_n|K_n}{{\eu a}}
\end{eqnarray*}
donde $\xbinom{L_n|K_n}{{\eu a}}$ donde es el s{\'\i}mbolo de Artin
y ${\eu a}\in C$. As{\'\i}
\begin{eqnarray*}
A_n&\stackrel{\cong}{\longto}&X/\gamma_{n,e}Y_e\\
C&\longmapsto& x\bmod \gamma_{n,e} Y_e,
\end{eqnarray*}
con $x|_{L_n}=\xbinom {L_n|K_n}{{\eu a}}$ y ${\eu a}\in C$.

\begin{teorema}\label{T7.8.46}
Para $m\geq n\geq e$, los siguientes diagramas son conmutativos:
\begin{gather*}
\xymatrix{
A_n\ar[r]^{\hspace{-.5cm}
\sim}_{\hspace{-.5cm}t_n}\ar[d]_{\phi_{n,m}}&X/\gamma_{n,e}Y_e
\ar[d]_{\varepsilon_{n,m}}\\ 
A_m\ar[r]^{\hspace{-.5cm}\sim}_{
\hspace{-.5cm}t_m}&X/\gamma_{m,e}Y_e
}
\qquad\qquad
\xymatrix{
A_m\ar[r]^{\hspace{-.5cm}\sim}_{
\hspace{-.5cm}t_m}\ar[d]_{N_{m,n}}^{\rm{norma}}
&X/\gamma_{m,e}Y_e
\ar[d]^{\rest}\\ 
A_n\ar[r]^{\hspace{-.5cm}\sim}_{\hspace{-.5cm}
t_n}&X/\gamma_{n,e}Y_e
}\\
\intertext{donde}
\varepsilon_{n,m}\colon X/\gamma_{n,e}Y_e\longto
X/\gamma_{m,e}Y_e, \quad\varepsilon_{n,m}(x\bmod \gamma_{n,e}
Y_e)=\gamma_{m,n} x\bmod \gamma_{m,e}Y_e,\\
\gamma_{m,n}
=\frac{\gamma_m}{\gamma_n}=\frac{(1+T)^{p^m}-1}{(1+T)^{p^n}-1}
\end{gather*}
y $t_n$, $t_m$ son los mapeos de Artin.
\end{teorema}

\begin{proof}
El segundo diagrama ya lo conocemos de la teor{\'\i}a de campos de
clase (Teorema \ref{T7.8.5}).

Sea $x\in X$ y sea ${\eu A}$ un ideal de $K_m$ en $A_m$ tal que
$x|_{K_m}=\xbinom{L_m|K_m}{{\eu A}}$. Para el primer diagrama,
si $\rho$ var{\'\i}a sobre todos los elementos de $\Gal(K_m/K_n)$,
entonces 
\begin{gather*}
{\eu a}=N_{m,n}{\eu A}\prod_{\rho\in\Gal(K_m/K_n)} \rho({\eu A}),\\
\begin{align*}
\big(t_m\circ \phi_{n,m}\big)\big(\overline{{\eu a}}\big)&=t_m\big(
\overline{\eu a})=\xbinom{L_m|K_m}{\eu a}=\prod_{\rho\in
\Gal(K_m/K_n)}\xbinom{L_m|K_m}{\rho({\eu A})}\\
&=\prod_{\rho}\rho\circ \xbinom{L_m|K_m}{\eu A}=\prod_{\rho}
\rho\big(x|_{K_m}\big)=\gamma_{m,n}x|_{K_m}\\
&=\big(\varepsilon_{n,m}
\circ t_n\big)\big(\overline{\eu a}\big). \tag*{$\fin$}
\end{align*}
\end{gather*}
\end{proof}

Por lo tanto si $A=\lim\limits_{\substack{\longto\\n}}
A_n\cong\lim\limits_{\substack{\longto\\ \varepsilon_{n,m}}} X/
\gamma_{n,e}Y_e$, se tiene el diagrama conmutativo de
 sucesiones exactas
\[
\xymatrix{
0\ar[r]&Y_e/\gamma_{n,e}Y_e\ar[r]\ar[d]^{\varepsilon_{n,m}}&
X/\gamma_{n,e}Y_e\ar[r]\ar[d]^{\varepsilon_{n,m}}&X/Y_e\ar[r]
\ar[d]^{\varepsilon_{n,m}=\gamma_{m,n}}&0
\ar@{}[d]^{\hspace{1cm}\hbox{$n\leq m$.}}\\
0\ar[r]&Y_e/\gamma_{m,e}Y_e\ar[r]&X/\gamma_{m,e}Y_e\ar[r]&
X/Y_e\ar[r]&0
}
\]

Para $m \gg n$, $\varepsilon_{n,m}=\gamma_{m,n}\colon
X/Y_e\to X/Y_e$ es $0$ pues $\gamma_{m,e}X\subseteq Y_e$.
Es decir, tenemos $A=\lim\limits_{\substack{\longto\\ n}}
A_n\cong \lim\limits_{\substack{\longto\\n}}Y_e/\gamma_{n,e}Y_e$.

Similarmente si $Y_e\sim E=\Big(\bigoplus\limits_{i=1}^t \Lambda/
\langle p^{k_i}\rangle\Big)\oplus \Big(\bigoplus\limits_{j=1}^s
\Lambda/\langle g_j(T)\rangle\Big)$, $A\cong \lim\limits_{\substack{
\longto\\ n}}E/\gamma_{n,e}E$.

Adem\'as para $V:=\Lambda/\langle p^k\rangle$, 
\[
V/\gamma_{n,e}V
\cong\{p(T)\bmod \langle p^k, \gamma_{n,e}\rangle\mid p(T)\in
{\ma Z}_p[T]\}\cong \frac{\big({\ma Z}/p^k{\ma Z}\big)[T]}{T^{p^n-p^e}}
\cong C_{p^k}^{p^n-p^e}
\]
pues $\gamma_{n,e}$ es un polinomio
distinguido de grado $p^n-p^e$. Aqu\'i $C_{p^k}$ denota al grupo
c\'iclico de $p^k$ elementos.

Para $V=\Lambda/\langle g(T)\rangle$ con $g(T)$ polinomio
distinguido de grado $d$, se tiene $\frac{\gamma_{n+2,e}}{\gamma_{
n+1,e}}=\frac{P_{n+2}(T)}{P_{n+1}(T)}$ donde $P_n(T)=(1+T)^{p^n}
-1$ y $\frac{P_{n+2}(T)}{P_{n+1}(T)}$ act\'ua en
$\Lambda/\langle g(T)\rangle$, para $p^n\geq d$,
como $p\times\text{unidad}$.

De esta forma tenemos:
\begin{gather*}
0\longto V/pV\cong \Lambda/\langle p,g(T)\rangle\cong \Lambda/
\langle p, T^d\rangle \cong {\ma F}_p[T]/\langle T^d\rangle \longto\\
\longto
V/\gamma_{n+2,e}V\xrightarrow[]{\frac{\gamma_{n+2,e}}
{\gamma_{n+1,e}}} V/\gamma_{n+1,e}V
\longto 0,\\
V/\gamma_{n,e}V\cong C_{p^{a_{1,n}}}\oplus\cdots\oplus
C_{p^{a_{d,n}}}\quad\text{con $a_{i,n+1}=1+a_{i,n}$ para toda $i$}.
\end{gather*}
Por tanto $\lim\limits_{\substack{\longto\\ n}}E/\gamma_{n,e}E
\cong \big({\ma Q}_p/{\ma Z}_p\big)^{\sum\limits_{j=1}^s \gr g_j}
\oplus B$ con $p^aB=0$ donde $a\geq \max\{k_j\mid 1\leq i\leq t_i\}$.
Adem\'as $B=0\iff t=0\iff \mu=0$. Se tiene
\begin{align*}
{\ma Q}_p/{\ma Z}_p&=:R\cong W(p)=\{\xi\in{\ma C}^{\ast}\mid
\xi^{p^n}=1 \text{\ para alg\'un\ }n\in{\ma N}\}=\\
&=\bigcup_{n=1}^{\infty}\frac{\big(1/p^n\big){\ma Z}}{{\ma Z}}=
\bigcup_{n=1}^{\infty}\langle \zeta_{p^n}\rangle.
\end{align*}

Por tanto tenemos

\begin{teorema}\label{T7.8.47}
Si $K_{\infty}/K$ es una extensi\'on ${\ma Z}_p$ y si $A_n$ es el 
$p$--subgrupo de Sylow de $I_{K_n}$, se tiene que si $A:=
\lim\limits_{\substack{\longto\\ n}}A_n$, entonces $A\cong \big(
{\ma Q}_p/{\ma Z}_p\big)^{\lambda}\oplus A'$ con $A'$ de exponente
acotado, esto es, $p^a A=0$ para alg\'un $a\in{\ma N}$. Adem\'as
$A'=0\iff \mu=0$ y por tanto $A$ es un grupo divisible $\iff \mu=0$.
$\fin$
\end{teorema}

\begin{observacion}\label{O7.8.48} Lo anterior muestra que, puesto
que ${\ma Q}_p/{\ma Z}_p$ es divisible, la conjetura de Iwasawa de
que $\mu=0$ es equivalente a que $A$ sea divisible.
\end{observacion}

\begin{proposicion}\label{P7.8.49}
Sean $p$ un primo impar y $K$ un campo num\'erico finito de tipo
$\MC$. Sea $K_{\infty}/K$ la extensi\'on ${\ma Z}_p$--ciclot\'omica.
Entonces el mapeo $A_n^-\to A_{n+1}^-$ es inyectivo.
\end{proposicion}

\begin{proof}
Sea $I$ un ideal en $K_n$ que se hace principal $K_{n+1}$, $
\overline{I}\in A_n$. Digamos $I=\langle \alpha\rangle$, $\alpha\in K_{
n+1}$. Sea $\langle \sigma\rangle =\Gal(K_{n+1}/K_n)$. Entonces
$\langle \alpha^{\sigma-1}\rangle =\frac{I^{\sigma}}{I}=(1)$, esto es,
$\alpha^{\sigma-1}=\varepsilon \in E_{n+1}$, $E_{n+1}$ las
unidades de $K_{n+1}$.

Sea $N$ la norma de $K_{n+1}/K_n$. Por lo tanto
\begin{gather*}
N\varepsilon =N\alpha^{\sigma-1}=\frac{N\alpha^{\sigma}}{N\alpha}
=\frac{N\alpha}{N\alpha}=1.\\
\intertext{Por tanto $\tilde{I}\to \varepsilon$ induce el mapeo}
\ker(A_n\to A_{n+1})\longto\frac{\ker N|_{E_{n+1}}}{I_{\langle\sigma
\rangle} E_{n+1}}=\frac{\ker N|_{E_{n+1}}}{E_{n+1}^{\langle\sigma-1
\rangle}}\\
:=H^{-1}\big(\Gal\big(K_{n+1}/K_n\big), E_{n+1}\big).
\end{gather*}

Sea $I$ que representa a una clase en $A_{n}^-$. Entonces
$\langle\alpha^{1+J}\rangle = I^{1+J}=\langle \beta\rangle$, $\beta
\in K_n$, por lo que $\beta^{\sigma}=\beta$. Se sigue que $\alpha^{
1+J}=\beta\eta$ con $\eta\in E_{n+1}$.

Sea $\alpha_1=\frac{\alpha^2}{\eta}$, $\varepsilon_1=\alpha_1^{
\sigma-1}=\frac{\big(\alpha^{\sigma-1}\big)^2}{\eta^{\sigma-1}}=
\frac{\varepsilon^2}{\eta^{\sigma-1}}\in E_{n+1}$ y
\begin{align*}
\varepsilon_1^{1+J}&= \big(\alpha_1^{1+J}\big)^{\sigma-1}=
\Big(\Big(\frac{\alpha^2}{\eta}\Big)^{1+J}\Big)^{\sigma-1}=\Big(\frac{
\big(\alpha^{1+J}\big)^2}{\eta^{1+J}}\Big)^{\sigma-1}=
\Big(\frac{\big(\beta\eta\big)^2}{\eta^{1+J}}\Big)^{\sigma-1}\\
&= \big(\beta^2\big)^{\sigma-1}\big(\eta^2\eta^{-1-J}\big)^{\sigma-1}=
1\cdot \big(\eta^{1-J}\big)^{\sigma-1}=\big(\eta^{\sigma-1}\big)^{1-J}
\in E_{n+1}^-.\\
\intertext{Adem\'as}
E_{n+1}^-&=\{\varepsilon\in E_{n+1}\mid \varepsilon^{1+J}=1\}=
\{\varepsilon\in E_{n+1}\mid \overline{\varepsilon}=\varepsilon^{-1}\}\\
&=\{\varepsilon\in E_{n+1}\mid |\varepsilon|=1\}\quad\text{y}\quad
\big(\varepsilon^{\sigma}\big)^{1+J}=\big(\varepsilon^{1+J}\big)^{
\varepsilon}=1^{\sigma}=1
\end{align*}
por lo que $E_{n+1}^-=W_{n+1}$ las ra{\'\i}ces de unidad en
$K_{n+1}$.

Por lo tanto $\big(\varepsilon_1^m\big)^{1+J}=1$ para alguna $m$.
Tambi\'en $N\varepsilon_1=\big(N\alpha_1\big)^{\sigma-1}=1$. Se
tiene que $H^1\big(\Gal\big(K_{n+1}/K_n\big), W_{n+1}\big)=\{1\}$,
esto es, si $\varepsilon_1\in W_{n+1}$ y $N\varepsilon_1=1$,
entonces existe $\varepsilon_2\in W_{n+1}$ tal que $\varepsilon_1
=\varepsilon_2^{\sigma-1}$.

En efecto, $\varepsilon_1=\alpha_1^{\sigma-1}$ pero requerimos
$y\in W_{n+1}$ tal que $y^{\sigma-1}=\varepsilon_1$. Sean las
dos sucesiones exactas
\begin{gather*}
1\longto W_n\longto W_{n+1}\stackrel{\sigma-1}{\longto}
W_{n+1}^{\sigma-1}\longto 1\\
1\longto (W_{n+1}\cap \ker N)\longto W_{n+1}\stackrel{N}{\longto}
W_n\longto 1.
\end{gather*}

Para verificar la exactitud de la segunda sucesi\'on hay que ver
que $N$ es sobre. Si $\zeta_p\notin K_0$, $\zeta_p\notin K_m$ para
toda $m$ y en este caso $NW_n=W_n^p=W_n$ y por tanto $N$
es sobre. Ahora, si $\zeta_p \in K_0$, $K_{n+1}=K_n(\zeta_p^m)$
para alg\'un $m\geq n+1$ y $W_{n+1}=\langle \zeta_{p^m}\rangle
\times \langle \zeta_t\rangle$ con $\mcd (t,p)=1$ y $W_n=\langle
\zeta_{p^m}^p\rangle \times \langle \zeta_t\rangle$, 
$N\zeta_{p^m}=\zeta_{p^{m-1}}=\zeta_{p^m}^p$ y 
$N\zeta_t=\zeta_t^p$ y como $\mcd (t,p)=1$
se tiene que $\langle\zeta_t\rangle=\langle \zeta_t^p\rangle$ de
donde $NW_{n+1}=W_n$.

As{\'\i} $\big|W_{n+1}^{\sigma-1}\big|=\frac{|W_{n+1}|}{|W_n|}=
\big|W_{n+1}\cap \ker N\big|$ y puesto que $W_{n+1}^{\sigma-1}
\subseteq W_{n+1}\cap \ker N$ se sigue la igualdad.

Regresando a nuestra demostraci\'on, se tiene que $\alpha_1^{\sigma
-1}=\varepsilon_1=\varepsilon_2^{\sigma-1}$ para alguna $
\varepsilon_2\in W_{n+1}$ por lo tanto $\big(\frac{\alpha_1}{
\varepsilon_2}\big)^{\sigma}=\frac{\alpha_1}{\varepsilon_2}$ lo cual
implica que $\frac{\alpha_1}{\varepsilon_2}\in K_n$.

Sin embargo $\big\langle\frac{\alpha_1}{\varepsilon_2}\big\rangle=
\langle\alpha_1\rangle =\langle \alpha^2\rangle =I^2$ en $K_{n+1}$
lo cual implica $\big\langle \frac{\alpha_1}{\varepsilon_2}\big\rangle=
I^2$ en $K_n$ por lo que $I^2$ es principal. Puesto que $p\neq 2$
y $o(\overline{I})$ es una potencia de $p$, $\overline{I}\in A_n^-$,
se sigue que $I$ es principal en $K_n$ y por tanto $A_n^-\to A_{n+
1}^-$ es inyectiva como quer{\'\i}amos probar. $\fin$
\end{proof}

\begin{observacion}\label{O7.8.50} El mapeo $A_n^+\to A_{n+1}^+$
no necesariamente es inyectivo y si $p=2$, entonces el mapeo
$A_n^-\to A_{n+1}^-$ no necesariamente es inyectivo. Finalmente
tenemos que si $C_n=I_{K_n}$, $C_n \stackrel{\phi}{\longto}
C_{n+1}$ es la conorma, entonces se tiene 
$\begin{array}{ccccc}
C_n&\stackrel{\phi}{\longto}&C_{n+1}&\stackrel{N}{\longto}&C_n\\
x&\longmapsto&\phi(x)&\longmapsto&x^p
\end{array}$ por lo que $\ker \phi\subseteq A_n=C_n(p)$.
\end{observacion}

\begin{proposicion}\label{P7.8.51} Sean $p$ un primo impar y $K$
un campo n\'umero finito de tipo $\MC$. Sea $K_{\infty}/K$ la 
extensi\'on ${\ma Z}_p$--ciclot\'omica. Entonces $X^-:=\lim\limits_{
\substack{\longleftarrow\\ n}}A_n^-$ no contiene 
$\Lambda$--subm\'odulos finitos. Por tanto hay una inyecci\'on
$X^-\hookrightarrow \big(\bigoplus\limits_{i=1}^t \Lambda/\langle
p^{k_i}\rangle\big)\oplus \big(\bigoplus\limits_{j=1}^s\Lambda/\langle
g_j(T)\rangle\big)$ con con\'ucleo finito.
\end{proposicion}

\begin{proof}
Sea $F\subseteq X^-$ un $\Lambda$--subm\'odulo finito. Sea
$\gamma_0$ un generador topol\'ogico de $\Gal(K_{\infty}/K)$. 
Puesto que $F$ es finito, existe un n\'umero natural $n_0$ tal que 
para $n\geq n_0$, $\gamma^{p^n}_0$ act\'ua trivialmente en $F$. 
Supongamos que existe $0\neq x=(\ldots,x_m,x_{m+1},\ldots)\in F
\subseteq \lim\limits_{\substack{\longleftarrow\\ n}}A_n^-$ con
$X_{m+1}\to X_m$ con la norma y $x_m\neq 0$ para $m\geq m_0$.

Sea $m\geq m_0,n_0$. Se tiene, por la
Proposici\'on \ref{P7.8.49}, que
$x_m\neq 0$ cuando lo mandamos en $A_{m+1}^-$. Apliquemos
$N=1+\gamma_0^{p^m}+\gamma_0^{2p^m}+\cdots+\gamma_0^{(p-1)
p^m}$ a $x$. Puesto que $m\geq n_0$, este elemento act\'ua como
multiplicaci\'on por $p$ puesto que $\gamma_0^{ip^m}x=x$. Tambi\'en
se tiene que $N=N_{K_{m+1}/K_m}$ por lo que $px_{m+1}=x_m\neq
0$ en $A_{m+1}^-$. Por lo tanto $px\neq 0$. Se sigue que
$\begin{array}{ccc}
F&\stackrel{p}{\longto}&F\\ x&\longmapsto&px\end{array}$ es
1--1 en el $p$--grupo finito lo cual \'unicamente puede suceder si 
$F=0$. $\fin$
\end{proof}

\begin{corolario}\label{C7.8.52} Sea $p$ un primo impar y sea $K$
un campo num\'erico finito de tipo $\MC$. Sea $K_{\infty}/K$ la 
extensi\'on ${\ma Z}_p$--ciclot\'omica. Si $\mu^-=0$ se tiene que
$X^-=\lim\limits_{\substack{\longleftarrow\\n}}A_n^-\cong {\ma Z}_p^{
\lambda^-}$ como ${\ma Z}_p$--m\'odulos.
\end{corolario}

\begin{proof}
Se tiene que 
\[
X^-\hookrightarrow \bigoplus\limits_{j=1}^s \Lambda/
\langle g_j(T)\rangle\quad\text{y}\quad 
\Lambda/\langle g(T)\rangle \cong
{\ma Z}_p[T]/\langle T^{\gr g}\rangle\cong {\ma Z}_p^{\gr g}.
\]
Por la estructura de ${\ma Z}_p$ m\'odulos, usando que ${\ma Z}_p$
es de ideales principales, se sigue que $X^-\cong {\ma Z}_p^{\sum
\limits_{j=1}^s \gr g_j}={\ma Z}_p^{\lambda^-}$. $\fin$
\end{proof}

\begin{window}[1,l,\xymatrix{
&M_{\infty}\\K_{\infty}\ar@{-}[ru]^{{\eu X}_{\infty}}\ar@{-}
[d]_{\Gamma}\\ K_0\ar@{-}[ruu]_G},{}]
Ahora estudiaremos la m\'axima $p$--extensi\'on abeliana de $K_{
\infty}$ no ramificada fuera de $p$. Sea $F$ un campo totalmente
real y sea $p$ un n\'umero primo impar. Sea $K_0:=F(\zeta_p)$, 
$K_{\infty}/K$ la extensi\'on ${\ma Z}_p$--ciclot\'omica, $M_{\infty}$
la m\'axima $p$--extensi\'on abeliana de $K_{\infty}$ no ramificada
fuera de $p$, ${\eu X}_{\infty}=\Gal(M_{\infty}/K_{\infty})$ y $G=\Gal(
M_{\infty}/K_0)$. Se tiene la sucesi\'on exacta $1\longto {\eu X}_{
\infty}\longto
G\stackrel{\pi}{\longto}
 \Gamma\longto 1$.
 
 Entonces ${\eu X}_{\infty}$
es un $\Gamma$--m\'odulo por conjugaci\'on: si $x\in\Gamma$, sea
$g\in G$ tal que $\pi(g)=x$ y $\xi\in{\eu X}_{\infty}$. Entonces
$x\circ \xi:=g\xi g^{-1}$.

Sea $M_n$ la m\'axima $p$--extensi\'on abeliana de $K_n$ no 
ramificada fuera de $p$ y consideremos $\omega_n:=\gamma_0^{p^n
}-1=(1+T)^{p^n}-1$. Entonces $\Gal(M_n/K_{\infty})={\eu X}_{\infty}/
\omega_n{\eu X}_{\infty}$. La demostraci\'on de esto es la misma
que la de cuando $X_n\cong X/Y_n$, es decir, con el subgrupo
conmutador pero sin grupos de inercia (ver Proposiciones \ref{P7.8.7},
y \ref{P7.8.11}).
\end{window}

Como vimos anteriormente, si $r_2=r_2(K_0)$, entonces 
\begin{gather*}
\begin{align*}
{\ma Z}_p^{r_2(K_n)+1+\delta_n}\times (\text{$p$--grupo finito})&=
{\ma Z}_p^{r_2p^n+1+\delta_n}\times (\text{$p$--grupo finito})\\
&\cong\Gal(M_n/K)
\end{align*}
\intertext{y donde $\delta_n$ es el error en la Conjetura de Leopoldt,
esto es,}
{\eu X}_{\infty}/\omega_n{\eu X}_{\infty}\cong {\ma Z}_p^{r_2 p^n+
\delta_n}\times (\text{$p$--grupo finito}).
\end{gather*}

Ahora bien, tenemos que $\langle p,\omega_n\rangle \subseteq\langle
p,T\rangle$ y ${\eu X}_{\infty}/\langle p,\omega_n\rangle {\eu X}_{
\infty}\cong {\ma F}_p^{r_2p^n+\delta_n}\times (\text{$p$--grupo finito})
$.

Se sigue del Lema de Nakayama (Proposici\'on \ref{P7.8.9}) y puesto
que ${\eu X}_{\infty}$ es un $\Lambda$--m\'odulo finitamente
generado, que ${\eu X}_{\infty}\sim \Lambda^a\oplus (
\Lambda\text{--torsi\'on})$ para alguna $a\geq 0$.

\begin{lema}\label{L7.8.53} La sucesi\'on $\delta_n$ est\'a acotada.
\end{lema}

\begin{proof}
Supongamos que $\delta_n>0$ para alguna $n$. Sea 
$\{\varepsilon_1, \ldots,\varepsilon_r\}$ 
una base de las unidades de $K_n$ congruentes
con $1$: $E_1:=E_1(K_n)=\{\xi\mid \xi\in K_n, \xi\text{\ unidad y\ }
\xi\equiv 1\bmod \pK \ \forall\ \pK|p\}$.

Rearreglando, suponemos que $\varepsilon_{\delta_n+1}, \ldots,
\varepsilon_r$ son independientes y generan $\overline{E}_1$ sobre
${\ma Z}_p$. Entonces se tiene que existen $a_{ij}\in{\ma Z}_p$
tales que $\varepsilon_i=\prod\limits_{j=\delta_n+1}^r \varepsilon_j^{
a_{ij}}$, $1\leq i\leq \delta_n$.

Sea $a'_{ij}$ la $n$--\'esima parcial de $a_{ij}$. M\'as precisamente
$a_{ij}\equiv a'_{ij}\bmod p^n$. Sea $\eta_i:=\varepsilon_i\cdot
\prod\limits_{j=\delta_n+1}^r\varepsilon_j^{-a'_{ij}}=\prod\limits_{j=
\delta_n+1}^r \varepsilon_j^{(a_{ij}-a'_{ij})}$, $1\leq i\leq \delta_n$.

Se tiene que $\eta_i$ es una $p^n$--potencia en $\overline{E}_1
\subseteq \prod\limits_{\pK|p}U_{1,\pK}$ y $\eta_1,\ldots, \eta_{
\delta_n}$ generan un subgrupo $\big({\ma Z}/p^n{\ma Z}\big)^{
\delta_n}$ de $K_n^{\ast}/\big(K_n^{\ast}\big)^{p^n}$ ya que
$\eta_i=\varepsilon_i\varepsilon_{\delta_n+1}^{\alpha_0}\cdots
\varepsilon_r^{\alpha_r}$, $\alpha_j=a'_{ij}$.
Puesto que $\zeta_p\in K_0$, $\zeta_{p^n}\in K_n$ y por Teor{\'\i}a
de Kummer, $K_n\big(\big\{\eta_i^{1/p^n}\big\}^{\delta_n}_{i=1}\big)
/K_n$ es una extensi\'on abeliana con grupo de Galois isomorfo
a $\big({\ma Z}/p^n{\ma Z}\big)^{\delta_n}$. Tambi\'en por Teor{\'\i}a
de Kummer, puesto que $\eta_i$ es unidad, la extensi\'on es no
ramificada fuera de $p$. Ahora bien, $\eta_i$ es una $p^n$ potencia
en $U_{1,\pK}$ para toda $\pK|p$ por lo que los campos locales
respectivos satisfacen, $
K_n\big(\big\{\eta_i^{1/p^n}\big\}^{\delta_n}_{i=1}\big)_{\pK}=\big(
K_n\big)_{\pK}$, esto es $e_{\pK}=f_{\pK}=1$ y  por tanto cada $\pK
|p$ se descompone totalmente en $
K_n\big(\big\{\eta_i^{1/p^n}\big\}^{\delta_n}_{i=1}\big)/K_n$ y
en particular no son ramificadas. Todo lo anterior
prueba que $K\big(\big\{\eta_i^{1/p^n}\big\}_{i=1}^{\delta_n}\big)
/K_n$ es no ramificada y en particular $S:=
K_n\big(\big\{\eta_i^{1/p^n}\big\}^{\delta_n}_{i=1}\big)\subseteq L_n$
donde $L_n$ es la m\'axima $p$--extensi\'on abeliana de $K_n$ no
ramificada.

\[
\xymatrix{
&S\ar@{-}[r]&L_n\\ K_n\ar@{-}[ru]^{\big({\ma Z}/p^n{\ma Z}
\big)^{\delta_n}}\ar@/_/@{-}[rru]_{X_n}}
\]

En particular, $X_n$ tiene un cociente isomorfo a $\big({\ma Z}/
p^n{\ma Z}\big)^{\delta_n}$. Para $n\gg 0$ los t\'erminos $\Lambda/
\langle p^k\rangle$ de $X\sim \Big(\bigoplus\limits_{i=1}^t
\Lambda/\langle p^{k_i}\rangle \Big)\oplus \Big(\bigoplus\limits_{j=1
}^s\Lambda/\langle g_j(T)\rangle\Big)$ no pueden contribuir a 
${\ma Z}/p^n{\ma Z}$ pues son de exponente acotado. Los
t\'erminos $\bigoplus\limits_{j=1
}^s\Lambda/\langle g_j(T)\rangle$ a lo m\'as contribuyen con $\lambda
=\sum\limits_{j=1}^s\gr g_j$ de estos t\'erminos y en particular
$\delta_n\leq \lambda$. $\fin$
\end{proof}

\begin{lema}\label{L7.8.54} El Lema {\rm \ref{L7.8.53}} sigue siendo 
v\'alido en el caso en que $\zeta_p\notin K_0$.
\end{lema}

\begin{proof}
Sea $K'_0:=K_0(\zeta_p)$ y $\delta'_n=\delta\big(K_n(\zeta_p)\big)$.
Se tiene que $\delta_n\leq \delta'_n\leq \lambda_{K_0(\zeta_p)}$.
$\fin$
\end{proof}

Notemos que $\ran_{{\ma Z}_p}{\eu X}_{\infty}/\omega_n{\eu X}_{
\infty}=r_2 p^n+\delta_n$, $0\leq \delta_n\leq \lambda$, lo cual implica:

\begin{teorema}\label{T7.8.55} El m\'odulo ${\eu X}_{\infty}$ satisface
que ${\eu X}_{\infty}\sim \Lambda^{r_2}\oplus (\Lambda
\text{--torsi\'on})$.
\end{teorema}

\begin{proof}
La torsi\'on de ${\eu X}_{\infty}$ contribuye de la siguiente manera:
\lasa
\item[(1)] Si $V=\Lambda/\langle p^k\rangle$, entonces $V/
\omega_n V=\Lambda/\langle p^k,\omega_n\rangle \Lambda$ y
$\big|V/\omega_nV\big|=p^{kp^n}<\infty$.

\item[(2)] Si $V=\Lambda/\langle g(T)\rangle$ donde $g(T)$ es un
polinomio distinguido de grado $d$, entonces $\big|V/\omega_n V\big|
=p^{nd+c}$ donde $c$ es una constante.
\end{list}

De (1)  y (2) se sigue que la torsi\'on de ${\eu X}_{\infty}$ contribuye
con subm\'odulos finitos de ${\eu X}_{\infty}/\omega_n {\eu X}_{\infty}$
y se tiene que el ${\ma Z}_p$--rango de ${\eu X}_{\infty}/\omega_n
{\eu X}_{\infty}$ es igual a $r_2 p^n +\delta_n=r_2 p^n + o(1)$.

Finalmente, $\Lambda/\omega_n \Lambda\cong {\ma Z}_p^{p^n}$,
lo cual implica que ${\eu X}_{\infty}\sim \Lambda^{r_2}\oplus
(\Lambda\text{--torsi\'on})$. $\fin$
\end{proof}

Ahora estudiaremos una relaci\'on entre invariantes y los
$\ell$--subgrupos de Sylow de los niveles $K_n$ en una extensi\'on
${\ma Z}_p$.

Recordemos que si $K/F$ es una extensi\'on de Galois de grado $d$
y $\ell$ es un primo con $\ell\nmid d$, entonces si $B_F:={}_{\ell}
I_F$ y $B_K={}_{\ell}I_K$ son las $\ell$ partes 
primarias de los grupos de clases de ideales de $F$ y $K$
respectivamente, entonces la conorma $B_F\hookrightarrow B_K$
es 1--1 y la norma $B_K\stackrel{N}{\longto}B_F$ es suprayectiva
y $I_F(\ell)=I_K(\ell)$ si y solamente si $\ran_{\ell}
B_F=\ran_{\ell} B_K$ lo cual tambi\'en es equivalente a que la norma
$N_{K/F}\colon I_K(\ell)\to I_F(\ell)$ es un isomorfismo
(Proposici\'on \ref{P7.8.37}).

\begin{proposicion}\label{P7.8.56} Sea $p$ un n\'umero primo 
distinto a $\ell$ y sea $K_n$ una extensi\'on c{\'\i}clica de un campo
num\'erico $K_0$ de grado $p^n$. Sean $f=o(\ell \bmod p^n)$ y $
C_{\gamma}$ la parte $\ell$--primaria, es decir, el $\ell$--subgrupo
de Sylow del grupo de clases de ideales del subcampo de grado 
$p^{\gamma}$, $\gamma\leq n$.
Sea $D_n$ el kernel de la norma:
$0\longto D_n\longto {}_{\ell}C_n\stackrel{N}{\longto}{}_{\ell}C_{n-1}
\longto 0$
donde en general, si $A$ es un subgrupo abeliano ponemos ${}_{\ell}
A:=\{x\in A\mid \ell x=0\}$. Entonces $D_n\neq 0$ si y solamente si
$C_n\neq C_{n-1}$ y en este caso $\dim_{{\ma F}_{\ell}}D_n\geq f$.
$\fin$
\end{proposicion}
\[ 
{}_{p^n}\left\{\begin{array}{c}
\xymatrix{
K_n \ar@{-}[d]\\ K_{\gamma}\ar@{-}[d]^{\Big\} p^{\gamma}}\\K_0
}\end{array}\right.
\]

\begin{proof} Por la Proposici\'on \ref{P7.8.37} se tiene que $N$ es
suprayectiva y $N$ es un isomorfismo en ${}_{\ell}C_n\to {}_{\ell}C_{
n-1}\iff D_n=0\iff {}_{\ell}C_n={}_{\ell}C_{n-1}$. Sea $G=\Gal(K_n/K_0)
$. Entonces $D_n$ es un $G$--m\'odulo. Sea $\rho\colon G\to \Aut(
D_n)$ la representaci\'on de $G$ en $D_n$. Veamos que $\rho$ es
1--1, esto es, que la representaci\'on es {\em
fiel\index{representaci\'on fiel}}. Si $\rho$ no fuese 1--1, entonces $
\rho$ es trivial en el \'unico subgrupo de orden $p$, a saber, en
$S=\Gal(K_n/K_{n-1})$. Se tiene
\[
N_{n,n-1}D_n=\big\{\sum_{g\in S} gd\mid d\in D_n\big\}=\big\{\sum_{g
in S}d\mid d\in D_n\big\}={}_p D_n=0.
\]

Por otro lado, $|D_n|=\ell^s$ con $\ell\neq p$ lo cual implica que
$pD_n=D_n=0$. Supongamos pues que $D_n\neq 0$. Sea $D:=
D_n\otimes_{{\ma F}_{\ell}}\tilde{\ma F}_{\ell}$ la extensi\'on de 
escalares donde $\tilde{\ma F}_{\ell}$ es una cerradura algebraica
de ${\ma F}_{\ell}$. Entonces $D$ es un $G$--m\'odulo con acci\'on
$g\circ(d\otimes x):= gd\otimes x$.

Adem\'as, $\dim_{{\ma F}_{\ell}}D_n=\dim_{\tilde{\ma F}_{\ell}} D$. Sea
$D\cong \bigoplus\limits_{i=1}^t D_i$ con cada $D_i$--m\'odulo
irreducible. Se tiene que  $D$ es un 
$\tilde{\ma F}_{\ell}[G]$--m\'odulo.
Puesto que $p\neq \ell$, veremos que $\tilde{\ma F}_{\ell}[G]$ es
{\em semisimple\index{anillo semisimple}},
esto es, si $M$ es un $\tilde{\ma F}_{\ell}[G]$--m\'odulo
y $N<M$ subm\'odulo de $M$, entonces
necesariamente $M\cong N\oplus N'$ como 
$\tilde{\ma F}_{\ell}$--m\'odulos para alg\'un $N'$. 

En efecto, sea $N'$ un $\tilde{\ma F}_{\ell}$--espacio vectorial
tal que como espacios vectoriales, $M\cong N\oplus N'$. Sea
$\pi\colon M\to N$ la proyecci\'on lineal. En particular, $\pi(x)=x$ para
$x\in N$. Sea $\varphi:=\frac{1}{p^n}\Tr_G(\pi)$, esto es,
\[
\varphi(y)=\frac{1}{p^n}\sum_{g\in G} g\circ \pi(y).
\]
Puesto que $\pi\colon M\to N$, se define la acci\'on de $g$ en
$\Hom(M,N)$ por $(g\circ \varphi)(z):=g\varphi\big(g^{-1}z\big)$.
En particular $(g\circ \pi)(y)=g\pi\big(g^{-1}y\big)$.

Sea $\varphi\colon M\to N$ un $\tilde{\ma F}_{\ell}[G]$ homomorfismo.
Sea $j\colon N\to M$ el encaje. Entonces
\begin{align*}
(\varphi\circ j)(y)&=\varphi(y)=\frac{1}{p^n}\sum_{g\in G}(g\circ \pi)(y)=
\frac{1}{p^n}\sum_{g\in G}g\pi\big(g^{-1}y\big)\\
&\igual_{\substack{\uparrow\\ g^{-1}y\in N}}\frac{1}{p^n}\sum_{g\in G}
gg^{-1}y=\frac{1}{p^n}\sum_{g\in G}y=y,
\end{align*}
es decir $\varphi\circ j=\Id_N$ y por lo tanto la sucesi\'on
\[
\xymatrix{
0\ar[r]&N\ar[r]^j&M\ar[r]\ar@/^1pc/[l]^{\varphi}&M/N\ar[r]&0}
\]
se escinde y $M\cong N\oplus M/N$ como $G$--m\'odulos.

Ahora bien, como $\tilde{\ma F}_{\ell}$ es algebraicamente cerrado,
$\dim_{\tilde{\ma F}_{\ell}} D_i=1$ para todo $i$ debido a que
$\tilde{\ma F}_{\ell}$ es semisimple, conmutativo y por lo tanto
es un producto directo de {\em anillos simples\index{anillo
simple}}. Cada anillo simple
es un anillo de matrices sobre $\tilde{\ma F}_{\ell}$ y como es 
conmutativo esto sucede \'unicamente cuando es igual a $\tilde{\ma
F}_{\ell}$.

Ahora, $G$ es fiel en $D_n$ por lo $G$ es fiel en alg\'un $D_i$ pues
de lo contrario no ser{\'\i}a fiel en el \'unico subgrupo de orden $p$.
As{\'\i}, $G$ opera fielmente en alg\'un $D_i$. Si $G=\langle \sigma
\rangle$, $D_i\cong \tilde{\ma F}_{\ell}$ y
entonces $\varphi\colon G\to \Aut(
D_i)$ es 1--1. Sea $\varphi(\sigma)=a$, entonces $\varphi(\sigma^2)=
a^2, \ldots, 1=\varphi(1)=\varphi\big(\sigma^{p^n}\big)=a^{p^n}$ y 
$\varphi\big(\sigma^{p^{n-1}}\big)=a^{p^{n-1}}\neq 1$. Por lo tanto
$a=\zeta=\zeta_{p^n}$ es una ra{\'\i}z $p^n$--\'esima primitiva de 
$1$. M\'as precisamente
\[
\xymatrix{
G\ar[rr]^{\hspace{-.9cm}\varphi}
\ar[rd]_{\rho}&&\Aut(D\otimes\tilde{\ma F}_{\ell})=\Aut(D)
\\ &\Aut(D_n)& \text{$\begin{array}{rcl}
\varphi(\sigma)|_{D_i}\colon
D_i&\longto& D_i\\ x&\longmapsto& \zeta x
\end{array}$}}
\]

Sea $\{v_1,\ldots,v_s\}$ una
base de $D_n$ sobre ${\ma F}_{\ell}$, $s=t$.
Entonces $\{1\otimes v_1,\ldots,1\otimes v_s\}$ es 
una base de $D$ sobre
$\tilde{\ma Z}_{\ell}$. Sea $D_i=\langle \omega\rangle$ con
$\omega=\sum\limits_{i=1}^s\lambda_i\otimes v_i$, $\lambda_i\in
{\ma F}_{\ell}$, $\sigma\omega=\zeta \omega=\sum\limits_{i=1}^s
\zeta\lambda_i\otimes v_i$.

Por otro lado, $\sigma\omega =\sum\limits_{i=1}^s \lambda_i\otimes
\sigma v_i$. Sea
\begin{eqnarray*}
\rho\colon G&\longto& \Aut(D_n)\\
\sigma&\longmapsto &A=(a_{ij})\quad
 (\text{matriz sobre\ } {\ma F}_{\ell}).
\end{eqnarray*}
Por tanto
\begin{align*}
\sigma\omega&=\sum_{i=1}^s \Big\{\lambda_i\otimes \Big(
\sum_{j=1}^s a_{ij}v_j\Big)\Big\}=\sum_{i=1}^s\sum_{j=1}^s \lambda_i
a_{ij}\otimes v_j\\
&=\sum_{k=1}^s \Big(\sum_{j=1}^s \lambda_ja_{jk}\Big)\otimes v_k,
\end{align*}
lo cual implica que $\zeta\lambda_k=\sum\limits_{j=1}^s a_{jk}
\lambda_j$.

Se sigue que $\big(A-\zeta I)\left(\begin{array}{c}
\lambda_1\\ \vdots\\ \lambda_s\end{array}\right)=0$, esto es, $\zeta$
es ra{\'\i}z del polinomio caracter{\'\i}stico $f(T)\in {\ma F}_{\ell}[T]$
de $A$. Ahora bien, $s=\gr f(T)$ y $\gr f(T)\geq \gr \Irr(T,\zeta,{\ma 
F}_{\ell})$. Puesto que $f=\big[{\ma F}_{\ell}(\zeta):{\ma F}_{\ell}\big]$
se sigue que $f\leq s$. $\fin$
\end{proof}

\begin{corolario}\label{C7.8.57} Sea $D$ una representaci\'on finito
dimensional de un grupo c{\'\i}clico $G$ de orden $p^n$ sobre 
${\ma F}_{\ell}$ con $\ell\neq p$. Si $D$ es fiel, entonces $\dim_{
{\ma F}_{\ell}}D\geq f$ donde $f=o(\ell\bmod p^n)$. $\fin$
\end{corolario}

\begin{teorema}\label{T7.8.58} Sea $K_{\infty}/K_0$ una extensi\'on
${\ma Z}_p$ y sea $C_n$ el $\ell$--subgrupo de Sylow del grupo de
clase $I_{K_n}$ de $K_n$ donde $\ell$ es un primo diferente
a $p$. Si los $\ell$--rangos de $C_n$ est\'an acotados, es decir,
$\dim_{{\ma F}_{\ell}}\frac{C_n}{\ell C_n}\leq M$ para toda $n$ y 
alg\'un $M>0$, entonces los \'ordenes de $C_n$ est\'an acotados.
\end{teorema}

\begin{observacion}\label{O7.8.59} El enunciado del Teorema 
\ref{T7.8.58} se puede pensar en algo as{\'\i} como si $\mu_{\ell}=0$,
entonces $\mu_{\ell}=\lambda_{\ell}=0$.
\end{observacion}

\begin{proof}[Teorema {\rm \ref{T7.8.58}}]
En caso de que los \'ordenes no estuviesen acotados se tendr{\'\i}a
que $D_n\neq 0$ para $n$ suficientemente grande. Puesto
\[
f_n=o(\ell \bmod p^n)\xrightarrow[n\to \infty]{} \infty\quad \text{se
sigue que}\quad \dim_{{\ma F}_{\ell}} D_n\xrightarrow[n\to \infty]{}
\infty
\]
lo cual implica que $\dim_{{\ma F}_{\ell}} {}_{\ell}C_n=\dim_{{\ma F}_{
\ell}}\frac{C_n}{\ell C_n}\xrightarrow[n\to\infty]{}\infty$. $\fin$
\end{proof}

\begin{observacion}\label{O7.8.60} El Teorema \ref{T7.8.58} no es
\'unicamente aplicable a las extensiones ${\ma Z}_p$. Por ejemplo,
dado un campo num\'erico $K$ y $\ell$ un n\'umero primo, o bien
los $\ell$--rangos de $C_L$ tienen a infinito o $I_{L}(\ell)=
I_{K}(\ell)$ cuando $p\to \infty$.
\end{observacion}

\begin{teorema}\label{T7.8.61}
Sea $K_{\infty}/K_0$ una extensi\'on ${\ma Z}_p$ tal que $K_n$ es
campo de tipo $\MC$. Sea $\ell$ un n\'umero primo distinto de
$p$ y supongamos que $\zeta_{\ell}\in K_0$. Si $\big|C_n^-(\ell)\big|$
est\'a acotado, entonces $\big|C_n^+(\ell)\big|$ tambi\'en est\'a
acotado.
\end{teorema}

\begin{proof}
Por el Teorema \ref{T7.8.40} se tiene $\dim_{{\ma F}_{\ell}} {}_{\ell}
C_n^+\leq \dim_{{\ma F}_{\ell}} {}_{\ell} C_n^-+1$ por lo tanto los
$\ell$--rangos de $C_n^+(\ell)$ est\'an acotados. Por el Teorema
\ref{T7.8.58}, $\big|C_n^+(\ell)\big|$ est\'an acotados. $\fin$
\end{proof}

\section{Ejemplo de $\mu>0$}\label{S7.9}

En esta secci\'on presentamos los ejemplos construidos por
Iwasawa de extensiones ${\ma Z}_p$ donde $\mu>0$.
Para la parte de cohomolog\'ia de grupos, ver la Subsecci\'on
\ref{CClaseS1.5}.

\begin{proposicion}[Takagi y Chevalley]\label{P7.8.62}
Sea $L/K$ una extensi\'on c{\'\i}clica finita de campos num\'ericos con
grupos de Galois $G$. Sea $I_L$ el grupo de clases de ideales de
$L$ y sea $I_L^G:=\{c\in I_L\mid c^g=c\ \forall\ g\in G\}$.
Sean
\[
e_0(L|K):=\prod_{\pK\in{\ma P}_K}e(\pL|\pK),\quad e_{\infty}(L|K):=
\prod_{\substack{\pK_{\infty}, \rm{\ primos}\\ \rm{infinitos}}}
e(\pL_{\infty}|\pK_{\infty})
\]
y sean $E_L$ y $E_K$ los grupos de unidades de $L$ y $K$ 
respectivamente. Entonces
\[
\big|I_L^G\big|=\frac{h(K)e(L|K)}{[L:K][E_K: N_{L/K}L^{\ast}\cap
E_K]}
\]
donde $h(K)$ es el n\'umero de clase de $K$ y $e(L|K)
=e_0(L|K)e_{\infty}(L|K)$.
\end{proposicion}

\begin{proof}
Sean $P_L$ es grupo de ideales fraccionarios principales de $L$ y
$D_L$ el grupo de ideales fraccionarios no cero de $L$. Se tiene
la sucesi\'on exacta $0\longto P_L\longto D_L\longto I_L\longto 0$.
En cohomolog{\'\i}a, obtenemos la sucesi\'on exacta
\[
0\longto P_L^G\longto D_L^G\longto I_L^G\longto H^1(G, P_L)
\longto H^1(G, D_L)\longto \cdots
\]

Ahora bien $H^1(G,D_L)=\frac{\ker N_G|_{D_L}}{I_G D_L}$. Se tiene que
si $G=\langle \sigma\rangle$, $I_GD_L=\langle \sigma-1\rangle I_L$ y
\begin{gather*}
H^1(G,D_L)=\frac{\langle \pL/\pL^{\sigma}\mid \pL\in {\ma P}_L\rangle}
{\langle \pL/\pL^{\sigma}\mid \pL\in{\ma P}_L\rangle}=\{0\}.\\
\intertext{Por tanto obtenemos la sucesi\'on exacta}
0\longto P_L^G\longto D_L^G\longto I_L^G\longto
H^1(G,P_L)\longto 0.\\
\intertext{De donde}
0\longto D_L^G/P_L^G\longto I_L^G
\longto H^1(G, P_L)\longto 0.
\end{gather*}
Se sigue que $\big|I_L^G\big|=\big[D_L^G:P_L^G\big]
\big|H^1(G,P_L)\big|$. 

Ahora bien, $D_L^G\supseteq P_L^G\supseteq P_K$ por lo que
\begin{gather*}
\big[D_L^G:P_L^G\big]=\frac{\big[D_L^G:P_K\big]}{\big[P_L^G:P_K\big]}
=\frac{\big[D_L^G:D_K\big]\big[D_K:P_K\big]}{\big[P_L^G:P_K\big]}.\\
\intertext{Por otro lado,}
D_L^G=\big\langle \prod_{\pL|\pK}\pL\mid \pK\in{\ma P}_K\big\rangle,
\quad D_K=\langle \pK\mid \pK\in {\ma P}_K\rangle=
\big\langle \big(\prod_{\pL|\pK}\pL\big)^{e(\pL|\pK)}\mid \pK\in{\ma P
}_K\big\rangle
\end{gather*}
lo cual implica que $\big[D_L^G: D_K\big]=e_0(L|K)$ y  por ende
$\big[D_L^G:P_L^G\big]=e_0(L|K)\frac{h_K}{\big[P_L^G:P_K\big]}$.

Ahora bien, se tiene la sucesi\'on exacta
\begin{gather*}
\begin{array}{ccccccccc}
1&\longto&E_L&\longto&L^{\ast}&\longto&P_L&\longto&0\\
&&&&\alpha&\longmapsto&\langle\alpha\rangle
\end{array}\\
\intertext{con lo cual obtenemos en cohomolog{\'\i}a la
sucesi\'on exacta}
0\longto E_L^G=E_K\longto \big(L^{\ast}\big)^G=K^{\ast}
\longto P_L^G\longto H^1(G,E_L)\longto H^1(G, L^{\ast})
\end{gather*}
con $H^1(G,L^{\ast})=\{0\}$ por Teorema 90 de Hilbert.

Se sigue que
\[
0\longto K^{\ast}/E_K\cong P_K\longto P_L^G\longto H^1(G,E_L)
\longto 0
\]
es exacta, lo cual implica que $\big[P_L^G:P_K\big]=\big|H^1(G, E_L)
\big|$. 

Para un $G$--m\'odulo $A$ se define le {\em cociente de
Herbrand\index{cociente de Herbrand}} por $
\varphi(A):=\frac{h^0(A)}{h^1(A)}=\frac{|H^0(G,A)|}{|H^1(G,A)|}.$
Se tiene que si $A$ es finito, entonces $\varphi(A) =1$ y si $0\longto
A\longto B\longto C\longto 0$ es una sucesi\'on exacta de $G
$--m\'odulos, entonces $\varphi(B)=\varphi(A)\varphi(C)$.

Se tiene que $\varphi(E_L)=\frac{e_{\infty}(L|K)}{[L:K]}$ y
\begin{align*}
\varphi(E_L)^{-1}|H^0(G,E_L)|&=|H^1(G,E_L)|=\big[P_L^G:P_K\big]\\
&=|H^0(G,E_L)|\cdot \frac{[L:K]}{e_{\infty}(L|K)}.
\intertext{Ahora bien, $H^0(G,E_L)=\frac{E_L^G}
{N_{L/K}E_L}=\frac{E_K}{N_{L/K}E_L}$ por lo que}
\big[P_L^G:P_K\big]&=\big[E_K:N_{L/K}
E_L\big] \cdot \frac{[L:K]}{e_{\infty}(L|K)}.
\end{align*}

De la sucesi\'on exacta $1\longto E_L\longto L^{\ast}\longto
P_L\longto 1$ se tiene la sucesi\'on exacta
\[
0=H^1(G,K^{\ast})\longto H^1(G,P_L)\longto H^0(G,E_L)
\stackrel{\varphi}{\longto}H^0(G,L^{\ast})
\]
donde $\varphi\colon \frac{E_L^G}{N_{L/K}E_L}=\frac{E_K}{N_{
L/K} E_L}\longto \frac{(L^{\ast})^G}{N_{L/K} L^{\ast}}=\frac{
K^{\ast}}{N_{L/K}L^{\ast}}$.

Se sigue que $H^1(G,P_L)\cong \ker \varphi=\frac{E_K\cap N_{L/K}
L^{\ast}}{N_{L/K}E_L}$.
Por tanto
\begin{align*}
\big|I_L^G\big|&=\big[D_L^G:I_K\big]\big|H^1(G,P_L)\big|=
\frac{e_0(L|K)h_K}{\big[P_L^G:P_K\big]}\big|H^1(G, P_L)\big|\\
&= \frac{e_0(L|K)h_K}{\frac{|H^0(G,E_L)|[L:K]}{e_{\infty}(L|K)}}\cdot
\big[E_K\cap N_{L/K}L^{\ast}:N_{L/K}E_L\big]\\
&=\frac{e(L|K)h_K}{\big[E_K:N_{L/K}E_L\big][L:K]}\cdot
\big[E_K\cap N_{L/K}L^{\ast}:N_{L/K}E_L\big]\\
&=\frac{e(L|K)h_K}{[L:K] \big[E_K: E_K\cap N_{L/K} L^{\ast}\big]}.
\tag*{$\fin$}
\end{align*}
\end{proof}

\begin{lema}\label{L7.8.63}
Sea $\ell$ un entero, $\ell\geq 2$. Sea $K_d$ una extensi\'on de un
campo num\'erico $K$ de grado $d$. Sean ${\eu q}_1,\ldots, {\eu q}_t
$ ideales primos de $K$ que se descomponen totalmente en $K_d$.
Sea $K'$ una extensi\'on c{\'\i}clica de $K$ de grado $\ell$ en la
cual ${\eu q}_1,\ldots, {\eu q}_t$ son totalmente ramificados. Sea
$K'_d:=K'K_d$. Entonces $\frac{\big|I_{K'_d}\big|}
{\big|I_{K_d}\big|}$ es divisible por $\ell^{(t-[K:{\ma Q}])d-1}$.
\end{lema}

\begin{proof}
{\ }

\begin{window}[0,l,\xymatrix{
&K'_d\\ K_d\ar@{-}[ru]^{\ell}&&K'\ar@{-}[lu]_{d}\\ & K\ar@{-}[lu]^{d}
\ar@{-}[ru]_{\ell}},{}]
Se tiene que $K_d$ y $K'$ son linealmente disjuntos pues los ideales
primos ${\eu q}_i$ son totalmente ramificados en $K'$ y totalmente
descompuestos en $K_d$, $[K'_d:K_d]=\ell$. Aplicamos la
Proposici\'on de Takagi--Chevalley a la extensi\'on c{\'\i}clica $K'_d/
K_d$ y donde obtenemos, con $G=\Gal(K'_d/K_d)$, que:
\end{window}
\[
\big|I_{K'_d}^G\big|=\frac{e(K'_d|K_d)\big|I_{K_d}\big|}
{\big[K'_d:K_d\big]\big[E_{K_d}:E_{K_d}\cap N_{K'_d/K_d}(K'_d)^{
\ast}\big]}.
\]

Se tiene $E_{K_d}^{\ell}\subseteq E_{K_d}\cap N_{K'_d/K_d}(K'_d)^{
\ast}$ lo cual implica que
\begin{align*}
\frac{\big|I_{K'_d}\big|}
{\big|I_{K_d}\big|}&=\frac{\big|I_{K'_d}\big|}
{\big|I_{K'_d}^G\big|}\frac{\big|I_{K'_d}^G\big|}
{\big|I_{K_d}\big|}\\
&=\frac{e(K'_d|K_d)}{\big[K'_d:K_d\big]\big[E_{K_d}:E_{K_d}^{\ell}
\big]}\cdot [E_{K_d}\cap N_{K'_d/K_d}K'_d: E_{K_d}^{\ell}\big]
\cdot \frac{\big|I_{K'_d}\big|}
{\big|I_{K_d}^G\big|}.
\end{align*}

Se sigue que $\frac{e(K'_d|K_d)}{\ell\cdot [K'_d:K_d]
[E_{K_d}:E_{K_d}^{\ell}]}\Big| \frac{|I_{K'_d}|}
{|I_{K_d}|}$.

Ahora bien, sobre cada ${\eu q}_i$ hay $d$ primos en $K_d$ y
cada uno de ellos es totalmente ramificado en $K'_d/K_d$ lo cual
implica que $\ell^{td}|e(K'_d|K_d)$. Adem\'as tenemos $\big[K'_d:
K_d\big]=\ell$.

Por el teorema de la unidades de Dirichlet, tenemos que si $s=r_1+
r_2\leq r_1+2r_2=\big[K_d:{\ma Q}\big]=d[K:{\ma Q}]$, entonces
$E_{K_d}\cong {\ma Z}^{s-1}\times (\text{grupo finito})$.

Se tiene que $\big[E_{K_d}:E_{K_d}^{\ell}\big]= \ell^s$ o $\ell^{s-1}$
correspondiendo cada caso a si $W_{K_d}(\ell)\neq \{1\}$ o $
W_{K_d}(\ell)=\{1\}$ respectivamente. Por lo tanto $
\big[E_{K_d}:E_{K_d}^{\ell}\big]|\ell^{d[K:{\ma Q}]}$ lo cual implica que
$\frac{e(K'_d|K_d)}{\ell\cdot 
[E_{K_d}:E_{K_d}^{\ell}]}$ es dividido por $\ell^{td-d[K:{\ma Q}]-1}$ y
finalmente obtenemos que $\ell^{d(t-[K:{\ma Q}])-1}\Big|
\frac{|I_{K'_d}|} {|I_{K_d}|}$. $\fin$
\end{proof}

Todo lo anterior es la base en los ejemplos en que Iwasawa encontr\'o
extensiones ${\ma Z}_p$ tales que $\mu>0$.

\begin{teorema}\label{T7.8.64} Sea $K_{\infty}/K_0=K$ una 
extensi\'on ${\ma Z}_p$. Sean ${\eu q}_1,\ldots ,{\eu q}_t$ ideales
primos de $K$ que se descomponen totalmente en $K_{\infty}$.
Sea $K'$ una extensi\'on c{\'\i}clica de $K$ de grado $\ell$ en
donde ${\eu q}_1,\ldots ,{\eu q}_t$ son totalmente ramificados. 
Entonces si $K'_{\infty}=K_{\infty}K'$ y $K'_n=K_nK'$, se tiene
que si $\ell^{e'_n}=\big|I_{K'_n}(\ell)\big|$ entonces
$e'_n\geq (t-[K:{\ma Q}])p^n-1$.
\end{teorema}

\begin{proof}
Se tiene
\[
\xymatrix{K'_0\ar@{--}[r]\ar@{-}[d]&K'_n\ar@{--}[r]\ar@{-}[d]^{\ell}&
K'_{\infty}\ar@{-}[d]\\K_0\ar@{--}[r]&K_n\ar@{--}[r]&K_{\infty}
}
\]
$\big[K'_n:K_n\big]=\ell$ y el resultado es consecuencia del Lema
\ref{L7.8.63}. $\fin$
\end{proof}

Para dar un ejemplo donde $\mu>0$, necesitamos una extensi\'on
${\ma Z}_p$ donde exista una infinidad de primos totalmente
descompuestos en $K_{\infty}/K$. En es caso, podemos tomar $t$
arbitrariamente grande. Procedemos as{\'\i}: primero agregamos la
ra{\'\i}z $\zeta_{\ell}$: \quad
$\xymatrix{K_0(\zeta_{\ell})\ar@{-}[r]\ar@{-}[d]&
K_{\infty}(\zeta_{\ell})\ar@{-}[d]\\K_0\ar@{-}[r]&K_{\infty}}$.

En ese caso se verifica que si hay una infinidad de primos
totalmente descompuestos en $K_{\infty}(\zeta_{\ell})/K_0(\zeta_{\ell}
)$, entonces tambi\'en hay una infinidad de primos totalmente 
descompuestos en $K_{\infty}/K_0$, esto es, podemos suponer que
$\zeta_{\ell}\in K_0=K$.

Sea $\alpha\in K$ tal que $v_{{\eu q}_i}(\alpha)=1$ para $1\leq i\leq t$,
es decir
\[
\langle \alpha\rangle ={\eu q}_1\cdots {\eu q}_t{\eu a}
\]
con ${\eu a}$ primo relativo a ${\eu q}_i$, $1\leq i\leq t$. Sea $K':=
K_0\big(\alpha^{1/\ell}\big)$.

Entonces ${\eu q}_1,\ldots, {\eu q}_t$ son ramificados en $K'$ y por
lo tanto $\mu_{K'}\geq t-[K:{\ma Q}]>0$ cuando $\ell=p$ y
$e'_n\geq (t-[K:{\ma Q}])p^n-1$, $\ell^{e'_n}$ es la potencia de 
$\ell$ que divide a $h(K'_n)$.

\begin{teorema}\label{T7.8.65} Sea $K$ un campo num\'erico finito
de tipo $\MC$. Entonces
\l
\item Existe una extensi\'on ${\ma Z}_p$, $K_{\infty}/K$ que es 
Galois sobre $K^+$ tal que si $\Gamma=\Gal(K_{\infty}/K)$, 
entonces $\Gamma=\Gamma^-$.
\[
\xymatrix{K\ar@{-}[r]^{\Gamma}\ar@{-}[d]_{\{1,J\}}&K_{\infty}\\
K^+\ar@{-}[ru]_{G}}
\]
Se tiene la sucesi\'on $1\longto \Gamma\longto G\longto \{1,J\}
\longto 1$ donde $J$ act\'ua en $\Gamma$ por conjugaci\'on.

\item Sea ${\eu q}^+$ un primo de $K^+$ que no divide a $p$ y es
inerte en $K$, y sea ${\eu q}$ el respectivo primo en $K$. Entonces
${\eu q}$ se descompone completamente en $K_{\infty}$ y por el
teorema de densidad de Cevotarev, hay una infinidad de tales
primos ${\eu q}^+$ y ${\eu q}$. Por la observaci\'on anterior, existe
$K'/K$ tal que $K'_{\infty}/K_0$ tiene invariante de Iwasawa $\mu'
>0$.
\end{list}
\end{teorema}

\begin{proof}
Sea $M$ la m\'axima $p$--extensi\'on abeliana de $K$ no ramificada
fuera de $p$.
En particular $K_{\infty}\subseteq M$. Se tiene ${\cal G}:=\Gal(M/K)
\sim U_p/\overline{E}$, donde $U_p:=\prod\limits_{\pK|p}U_{\pK}$
y $\overline{E}$ es la cerradura de las unidades globales en $U_p$.
Se tiene $\xymatrix{K\ar@{-}[r]^{{\cal G}}\ar@{-}[d]_{\{1,J\}}&M\\
K^+\ar@/_1pc/@{-}[ru]^{{\eu G}}}\qquad 1\longto {\cal  G}\longto {\eu G}
\longto \{1,J\}\longto 1$.

Cada primo ``real'' $\pK^+$ sobre $p$ en $K^+$ permanece primo
en $K$ o se descompone en un producto de dos primos $\pK_1$,
$\pK_2$ en $K$. En cada caso tenemos que $\prod\limits_{\pK|\pK^+}
U_{\pK}\subseteq U_p$ y v{\'\i}a el mapeo exponencial obtenemos
que hay un subgrupo de {\'\i}ndice finito isomorfo a $\prod\limits_{
\pK|\pK^+}p^m{\cal O}_{\pK}$.

Sea $A:=\prod\limits_{\pK|p}p^m{\cal O}_{\pK}$. Se tiene la sucesi\'on
exacta
\[
0\longto A^-\longto A\xrightarrow[]{1+J}(1+J)A\longto 0
\]
donde $A^-=\ker (1+J)$ y $\ran_{{\ma Z}_p} A^-=\ran A/(1+J)A$. Se
tiene que  $A^{1+J}\subseteq p^m{\cal O}_{\pK^+}$, por tanto
$\ran_{{\ma Z}_p} A/(1+J)A\geq 1$ lo cual implica que $\ran_{{\ma Z
}_p} A^-\geq 1$. Por lo tanto $\ran_{{\ma Z}_p} U_p^-\geq 1$.

M\'as a\'un, $\overline{E}$ contiene un subgrupo de {\'\i}ndice finito
el cual est\'a fijo bajo $J$ ($[E:E^+]<\infty$) y por lo tanto 
$\overline{E^+}$ est\'a fijo bajo conjugaci\'on compleja. Por 
simplicidad suponemos $p>2$ y para
${\cal G}=\Gal(M/K)$ tenemos que ${\cal G}\cong {\cal G}^+\times
{\cal G}^-$ y $\ran_{{\ma Z}_p}{\cal G}^-\geq 1$. Por lo tanto
 existe un factor $\Gamma$ de ${\cal G}$ tal que $\Gamma=
 \Gamma^-\cong{\ma Z}_p$. 
 
 Se sigue que $\Gamma$ es el grupo de Galois de una extensi\'on
 ${\ma Z}_p$ $K_{\infty}$ de $K$ la cual es normal sobre $K^+$
 pues si  $\sigma\colon K_{\infty}\to \overline{\ma Q}$ es un
 encaje con $\sigma|_{K^+}=\Id_{K^+}$ entonces $\sigma(K)=
 K$ y $\sigma|_K\in \Gal(K/K^+)=\{1,J\}$. Se tiene $K_{\infty}^J
 =K_{\infty}$. Si $\sigma|_K=\Id$, entonces $\sigma\in\Gamma$ por 
lo que $\sigma(K_{\infty})= K_{\infty}$. Si $\sigma|_K=J$
entonces $J\circ \sigma=\Id$ por lo que 
$J\circ \sigma(K_{\infty})$ y $\sigma(K_{\infty})=J(K_{\infty})=
K_{\infty}$.

Sea ${\eu q}^+$ un primo de $K^+$ inerte en $K/K^+$, con ${\eu q}^+
|p$. Sea $D$ el grupo de descomposici\'on de ${\eu q}^+$ en $G=
\Gal(K_{\infty}/K^+)$. Entonces $D$ es topol\'ogicamente c{\'\i}clico
pues ${\eu q}$ es no ramificado en $K_{\infty}$. Ahora bien, se 
tiene que $G=\langle \gamma, J\mid J\gamma J=\gamma^{-1},
\gamma\in \Gamma\rangle$ ya que $\Gal(K/K^+)=\langle 1,\gamma
\rangle \cong \{1,J\}$ y $J\circ \gamma=J\gamma J^{-1}=
J\gamma J=\gamma^{-1}$. $\xymatrix{K\ar@{-}[r]^{\Gamma}
\ar@{-}[d]_{\{1,J\}}&K_{\infty},\\ K^+}\quad \Gamma^-=\Gamma$.

Puesto que ${\eu q}^+$ permanece primo en $K$, no se puede
tener que $D\subseteq \Gamma$ lo cual implica que $D$ contiene
un elemento $J\gamma$ de orden $2$ lo cual finalmente implica
que ${\eu q}$ se descompone totalmente en $K_{\infty}/K$. $\fin$
\end{proof}

\begin{observacion}\label{O7.8.66} Para $p=2$ se toma el grupo
factor $G/G^{1+J}$ para obtener la menos parte y el argumento es
esencialmente el mismo que en el caso $p$ impar.
\end{observacion}

Veamos otros resultados.

\begin{teorema}[L. Washington]\label{T7.8.67} Sea $k$ una 
extensi\'on abeliana de ${\ma Q}$ y sea $K/k$ la extensi\'on ${\ma Z
}_p$--ciclot\'omica de $k$. Sea $\ell\neq p$ un n\'umero primo y
sea $\ell^{e_n}$ la potencia exacta de $\ell$ que divide al n\'umero
de clase $h_n=h(K_n)$ de $K_n$. Entonces $e_n$ est\'a acotado
y m\'as precisamente, $e_n$ es constante para $n$ suficientemente
grande. $\fin$
\end{teorema}

\begin{observacion}\label{O7.8.68} El Teorema \ref{T7.8.67} no se
cumple si la extensi\'on $K/k$ no es la ciclot\'omica (ver Teorema
\ref{T7.8.64}).
\end{observacion}

\begin{teorema}[L. Washington]\label{T7.8.69}
Sea $k$ un campo num\'erico abeliano imaginario, esto es, $\Gal
(k/{\ma Q})$ es un grupo abeliano y $K\nsubseteqq{\ma R}$. Sea
$K/k$ la extensi\'on ${\ma Z}_p$--ciclot\'omica de $k$. Sean 
\[
h_n=
h(K_n)\quad \text{y}\quad H=\{\ell\mid \ell 
\text{\ es primo y $\ell|h_n^-$ para alg\'un $n$}\}.
\]
 Entonces $H$ es infinito. $\fin$
\end{teorema}

Finalmente veamos los invariantes de Iwasawa en campos de 
funciones.

Sea $\K$ un campo de funciones con campo de constantes el
campo finito de $q$ elementos ${\ma F}_q$, $q=\ell^m$ con $\ell$
un n\'umero primo. Sea $\K_n:=\K {\ma F}_{q^{p^n}}$. Entonces
$\Gal(\K_n/\K)\cong \Gal({\ma F}_{q^{p^n}}/{\ma F}_q)\cong{\ma Z}/p^n
{\ma Z}$. Sea $\K_{\infty}:=\bigcup\limits_{n=1}^{\infty}\K_n$.
Entonces $\Gal(\K_{\infty}/\K)\cong \Gamma\cong {\ma Z}_p$.

Consideramos la {\em funci\'on zeta\index{funci\'on zeta en un
campo de funciones}} de $\K$
(ver la Subsecci\'on \ref{SRam1.2}): 
$Z_{\K}(u)=Z_0(u) =\frac{P_0(u)}
{(1-u)(1-qu)}$ donde $P(u)\in {\ma Z}[u]$ es un polinomio de grado
$2g$ y $g$ es el g\'enero de $\K$. Se tiene que el n\'umero de
clase de $\K$ est\'a dado por $h_0=h(\K)=P_0(1)$.

Sea $Z_{\K_n}(v)$ la funci\'on zeta de $\K_n$. Entonces
\[
Z_{\K_n}(v)=Z_{\K_n}\big(u^{p^n}\big)=\prod\limits_{j=1}^{p^n}
Z_{\K}\big(\zeta_{p^n}^ju\big).
\]

Por otro lado tenemos que $P_0(u)=\prod\limits_{i=1}^{2g}
\big(1-\alpha_i^{-1} u\big)$ donde $\alpha_1,\ldots, \alpha_{2g}$
son las ra{\'\i}ces de $P_0(u)$. Por lo tanto $
P_n\big(u^{p^n}\big)=\prod\limits_{
i=1}^{2g}\big(1-\alpha_i^{-p^n}u^{p^n}\big)$ lo cual implica que 
\begin{align*}
\frac{h_n}{h_0}&=\frac{P_n(1)}{P_0(1)}=\frac{\prod\limits_{
i=1}^{2g}\big(1-\alpha_i^{-p^n}\big)}{\prod\limits_{
i=1}^{2g}\big(1-\alpha_i^{-1}\big)}=
\frac{\prod\limits_{i=1}^{2g}\prod\limits_{j=1}^{p^n}
\big(1-\zeta_{p^n}^j \alpha_i^{-1}\big)}{\prod\limits_{
i=1}^{2g}\big(1-\alpha_i^{-1}\big)}\\
&=\prod\limits_{i=1}^{2g}\prod\limits_{j=1}^{p^n-1}
\big(1-\zeta_{p^n}^j \alpha_i^{-1}\big)=:A_n.
\end{align*}

Se sigue que $v_p(A_n)=\lambda n+\gamma$ para $n$ 
suficientemente grande y $\lambda\leq 2g$. En particular $\mu=0$.
Es decir tenemos:

\begin{teorema}\label{T7.8.70} Para campos de funciones 
congruentes los invariantes de Iwasawa satisfacen que $\lambda\leq
2g$, donde $g$ es el g\'enero del campo y $\mu=0$. $\fin$
\end{teorema}

%% file: Capitulo17.tex
\chapter{Teor\'ia de campos de clase}\label{CClaseC17}

\section{Introducci\'on}

La teor{\'\i}a de campos de clase es esencialmente el estudio de
las extensiones abelianas, en general finitas, 
aunque tambi\'en se estudian las extensiones infinitas, de cuatro clases
de campos: campos globales (campos num\'ericos y campos de
funciones) y campos locales (de igual caracter{\'\i}stica y de caracter{\'\i}stica
distinta). En primera instancia, no es as{\'\i} como se definen los campos
de clase, sino m\'as bien de otras formas (campos en donde un cierto
conjunto de primos se descomponen totalmente o subgrupos de ciertos
grupos que son normas, etc.), lo cual lleva, al final del d{\'\i}a, a que es lo
mismo que extensiones abelianas finitas.

Hay otros objetos los cuales pueden ser estudiados con la teor{\'\i}a de
campos de clase que no est\'an comprendidas en las cuatro anteriores:
ciertos campos de funciones no congruentes. Tambi\'en hay una teor{\'\i}a
de campos de clase no abeliana: el programa de Langlands. No haremos
m\'as historia de la teor{\'\i}a de campos de clase pues nos reservamos
esto para la Subsecci\'on \ref{CClaseC2}.

El objetivo de este cap\'itulo es presentar
los resultados fundamentales de la teor{\'\i}a de campos de clase con la
profundidad necesaria para poder aplicarla a la teor{\'\i}a de n\'umeros algebraica,
principalmente para campos de funciones. Una 
opci\'on pudiera haber sido presentar la
teor{\'\i}a de campos de clase para campos num\'ericos desde el punto de
vista de ideales.
Un magn{\'\i}fico libro con este enfoque es el libro de Janusz \cite{Jan96}.
Sin embargo esta aproximaci\'on 
no nos dir{\'\i}a cu\'ales son las diferencias con respecto
a los campos de funciones que son nuestro inter\'es primario.

Para los lectores interesados en las demostraciones y desarrollos completos
de la teor{\'\i}a de clase, hay muchos textos excelentes, tanto para campos
globales como para locales y algunos de ellos para ambos. El texto 
fundamental es el trabajo de Artin--Tate \cite{ArtTat61}. Para campos
locales podemos indicar los libros de Iwasawa y de Serre \cite{Iwa86,Ser}
as{\'\i} como el art{\'\i}culo de Serre \cite{Ser67} comprendido en el
libro cuyos editores son Cassels y Fr\"olich
\cite{CasFro67}. Para campos globales podemos mencionar el art{\'\i}culo
de Tate \cite{Tat67} el cual tambi\'en se encuentra en el libro de Cassels
y Fr\"olich. Para el estudio de ambos casos,
el local y el global, debemos mencionar los libros
de Neukirch \cite{Neu69,Neu86,Neu99} y el libro de Chevalley \cite{Che54}.
Tratados m\'as cl\'asicos son los trabajos de
Hasse \cite{Has50,Has,Has67-1}.

En la Secci\'on \ref{CClaseC1} presentamos algunos
preliminares aislados que ayudan a comprender el trabajo subsiguiente.
Estos preliminares constan fundamentalmente de la teor{\'\i}a de Kummer,
el automorfismo de Frobenius y el s\'imbolo de Artin.
La Subsecci\'on \ref{CClaseC2} tiene como objeto presentar una panor\'amica
general de la historia de la teor{\'\i}a de campos de clase. Esta subsecci\'on
es muy \'util para el lector cuando est\'e leyendo los cap{\'\i}tulos posteriores,
pues puede consultar en que contexto los resultados fueron apareciendo
y su raz\'on de ser. Hay muchas referencias para esta parte pero nos basamos
fundamentalmente en el magn{\'\i}fico art{\'\i}culo de Conrad \cite{Con}. Otras
referencias estupendas son el art{\'\i}culo de Hasse \cite{Has67} en el
libro de Cassels y Fr\"olich y los art{\'\i}culos y libros de Roquette 
\cite{Roq2001,Roq2005,Roq2013}.

La Secci\'on \ref{CClaseS1.2} da un r\'apido repaso a la teor\'ia de campos locales.
En la Secci\'on \ref{CClaseS1.5} estudiamos la cohomolog\'ia de grupos
que es el enfoque que damos para la teor\'ia de campos de
clase. La Secci\'on \ref{CClaseC3} estudia la teor{\'\i}a
local de campos de clase con \'enfasis en los grupos de Lubin--Tate.
Decidimos presentar en la Subsecci\'on \ref{IsomorfismoNeukirch} 
la teor\'ia de Neukirch \cite{Neu86} sobre el
isomorfismo que lleva su nombre y el cual obtiene, con otro enfoque
distinto a la de cohomolog\'ia de grupos, el teorema de reciprocidad.
Este mismo isomorfismo se puede aplicar a campos globales.
Este subsecci\'on puede omitirse sin p\'erdida de continuidad.

La Secci\'on \ref{CClaseC4} es la teor{\'\i}a global de campos de clase.
El teorema de existencia es el tema de la Secci\'on
\ref{S17.6.11N}. La Secci\'on \ref{CClaseS4.4} estudia los grupos de 
congruencias que transparentan
un poco la correspondencia dada por el teorema de reciprocidad
de Artin.
En la Subecci\'on \ref{CClaseS4.5} estudiamos la teor\'ia de campos 
de clase en campos num\'ericos, primero con del enfoque de id\`eles
y posteriormente en la Subsecci\'on \ref{CClaseS4.7} con el enfoque de
ideales y divisores. La Subsecci\'on \ref{CClaseS4.9-1} 
comprende  de manera central
el estudio de los campos de clase de Hilbert y de Hilbert extendido
en campos de funciones.
Finalizamos con un estudio de los campos de g\'eneros desde el
punto de vista de campos de clase en la Secci\'on \ref{CClaseS4.10}.

\section{Preliminares y antecedentes}\label{CClaseC1}

En esta primera parte, presentamos varios temas que ya sea
son necesarios para el desarrollo de estas notas, o bien
porque son partes de las demostraciones, posiblemente
no presentadas aqu{\'\i}, de los teoremas  fundamentales
de la teor{\'\i}a.

\subsection{Extensiones de Kummer\index{extensiones!Kummer}
 y de Artin--Schreier\index{extensiones!Artin--Schreier}}\label{CClaseS1.1}

Este tipo de extensiones ya fueron parcialmente discutidas en
la Subsecci\'on \ref{SRam2.1}.

En esta secci\'on estudiaremos las llamadas {\em extensiones de
Kummer\index{extensiones!Kummer}} y las {\em extensiones
de Artin--Schreier\index{extensiones!Artin--Schreier}}. Estas
\'ultimas pueden ser consideradas las extensiones de Kummer
aditivas. Las extensiones de Kummer juegan un papel preponderante
para los teoremas de existencia de campos de clase para las
extensiones c{\'\i}clicas de orden primo relativo a la caracter{\'\i}stica
y las extensiones de Artin--Schreier son usadas con el mismo
fin en el caso de campos de funciones para extensiones c{\'\i}clicas
de grado igual a la caracter{\'\i}stica. Nosotros no presentaremos
todos los detalles de la demostraci\'on de los teoremas de existencia
desde este punto de vista,
pero si el lector quiere profundizar m\'as en los detalles de las demostraciones
de los teoremas de existencia, es necesario aplicar este tipo de
extensiones.

\begin{teorema}[Teorema de independencia
de Artin\index{Artin!independencia de caracteres de $\sim$}]\label{CClaseT1.1.1} 
Sean $G$ un grupo multiplicativo, $F$ un campo y $\sigma_1,
\sigma_2,\ldots,\sigma_n$, $n$  
homomorfismos de grupos
distintos de $G$ en $\* F$. Entonces $\sigma_1,
\sigma_2,\ldots,\sigma_n$ son independientes, es decir, si 
$a_1,\ldots,a_n\in F$ son tales que
\[
a_1\sigma_1(x)+a_2\sigma_2(x)+\cdots+a_n\sigma_n(x)=0
\]
para todo $x\in G$, entonces $a_1=a_2=\ldots=a_n=0$.
\end{teorema}
\begin{proof} \cite[Ch. IV, Section 4,Theorem 4.1]{Lan93}. $\fin$
\end{proof}

Una aplicaci\'on del Teorema \ref{CClaseT1.1.1} 
es para el caso en que $L$ es un campo
arbitrario y $\sigma_1,\ldots,\sigma_n$ son distintos automorfismos
de $L$. Entonces $\sigma_1,\ldots,\sigma_n$ son independientes.
En este caso se toma $G=\* L$.

\begin{teorema}\label{CClaseT1.1.2} Sea $L/K$ una extensi\'on c{\'\i}clica
de grado $n$ y sea $G=\Gal(L/K)=\langle \sigma \rangle$. Sea
$\alpha\in L$. Entonces
\lasa
\item $\Tr_{L/K}\alpha=0\iff$ existe $\beta\in L$ tal que $\alpha=\beta
-\sigma \beta$, donde $\Tr_{L/K}\index{traza}\label{CClasetraza}=\Tr$ 
denota la traza de $L$ a $K$.

\item $\N_{L/K}\alpha=1\iff$ existe $\beta\in L$ tal que $\alpha=
\beta/\sigma\beta$, donde $\N_{L/K}\index{norma}\label{CClasenorma}
=\N$ denota la norma
de $L$ a $K$.
\end{list}
\end{teorema}
\begin{proof}

\noindent
(a) $\Leftarrow)$ Si $\alpha=\beta-\sigma\beta$, entonces 
\[
\Tr_{L/K}\alpha
=\Tr_{L/K}\beta-\Tr_{L/K}(\sigma\beta)=\Tr_{L/K}\beta-\Tr_{L/K}\beta=0.
\]

\noindent
$\Rightarrow)$ Puesto que $L/K$ es separable, por el Teorema \ref{CClaseT1.1.1}
se tiene que existe $\gamma\in L$ tal que
$\Tr_{L/K}\gamma=a\neq 0$ con $a\in K$ (de hecho,
$L/K$ es separable $\iff \Tr\neq 0$). Por tanto $\Tr_{L/K}(a^{-1}\gamma)=
a^{-1}\Tr_{L/K}\gamma =1$.

Sea $\alpha\in L$ tal que 
$\Tr_{L/K}\alpha=0$. Entonces $\sigma^0\alpha=\alpha=-\sum_{j=1}^{n-1}
\sigma^j\alpha$. Sea $\beta=\sum_{i=0}^{n-2}\Big(\sum_{j=0}^i 
\sigma^j \alpha\Big)
\sigma^i \gamma_1$ con $\Tr_{L/K}\gamma_1=1$. Entonces 
$\beta-\sigma\beta=\alpha$.

\noindent
(b) $\Leftarrow)$ Si $\alpha=\beta/\sigma\beta$, $\N\alpha=\N\beta/
\N(\sigma\beta)=\N\beta/\N\beta=1$.

\noindent
$\Rightarrow)$ Sea ahora $\N_{L/K}\alpha=1$. Consideremos 
\begin{align*}
\xi:&=c+\alpha\sigma(c)+\alpha\sigma(\alpha)\sigma^2
(c)+\cdots+\alpha\sigma(\alpha)\cdots\sigma^{n-2}(\alpha)\sigma^{n-1}(c)\\
&=c+\sum_{j=1}^{n-1}\Big(\prod_{i=0}^{j-1} \sigma^i(\alpha)\Big)
\sigma^j(c)\\
\intertext{con $c\in L$. Entonces}
\alpha\sigma(\xi)&= \alpha\sigma(c)+\sum_{j=1}^{n-1}\Big(
\prod_{i=0}^{j-1}\alpha\sigma^{i+1}(\alpha)\Big)\sigma^{j+1}(c)\\
&=\underbracket[0pt]{\alpha\sigma(c)}_{\substack{\uparrow\\ j=1}}+
\sum_{j=2}^{n-1}\Big(\prod_{i=0}^{j-1}\sigma^i(\alpha)
\Big)\sigma^j(c)
+\underbracket[0pt]{\Big(\prod_{i=0}^{n-1}\sigma^i(\alpha)\Big)}_{\substack{\uigual\\1}}
\underbracket[0pt]{\sigma^n(c)}_{\substack{\uigual\\c}}\\
&=c+\sum_{j=1}^{n-1}\Big(\prod_{i=0}^{j-1}\sigma^i(\alpha)\Big)
\sigma^j(c)=\xi.
\end{align*}

Esto es, $\alpha\sigma(\xi)=\xi$. Por el Teorema \ref{CClaseT1.1.1}, existe $c$
tal que $\xi=\beta\neq 0$ y se tiene $\alpha\sigma(\beta)=\beta$ por lo que
$\alpha=\beta/\sigma(\beta)$. $\fin$
\end{proof}

Otra versi\'on del siguiente resultado fue dada en el Teorema 
\ref{T9'.2.1}.

\begin{teorema}[Extensiones de Artin--Schreier\index{extensiones!Artin--Schreier}]\label{CClaseT1.1.3}
Sea $K$ un campo de caracter{\'\i}stica $p>0$, $\ca K=p$. Entonces $L/K$
es una extensi\'on c{\'\i}clica de grado $p \iff$ existe $z\in L$ tal que
$L=K(z)$ e $\Irr(z,T,K)=T^p-T-a\in K[T]$ y, en particular, $z^p-z=a$.
\end{teorema}
\begin{proof}

\noindent
$\Rightarrow)$ Sea $G:=\Gal(L/K)=\langle\sigma\rangle$, $o(\sigma)=p$.
Se tiene $\Tr_{L/K}1=p=0$. Por tanto existe $z\in L$ tal que $\sigma z-z=1$,
esto es, $\sigma z=z+1$. Se sigue que $\sigma^i z=z+i$ y $\sigma^i z=z
\iff p|i$. Por lo tanto $\Irr(z,T,K)=\prod_{i=0}^{p-1}(T-(z+i))$ es de
grado $p$.

Se tiene $\sigma(z^p-z)=(\sigma z)^p-(\sigma z)=(z+1)^p-(z+1)=z^p-z$
de donde obtenemos que $z^p-z=a\in K$ y $z^p-z-a=0$. Por lo tanto
$\Irr(z,T,K)=T^p-T-a=\prod_{i=0}^{p-1}(T-(z+i))$.

\noindent
$\Leftarrow) $Si $L=K(z)$ e $\Irr(z,T,K)=T^p-T-a$, entonces para toda
$i\in{\ma Z}$ setiene $i^p\equiv i\bmod p$ y $(z+i)^p-(z+i)=
z^p+i^p-z-i=z^p-z=a$.
Se sigue que $z,z+1,\ldots,z+(p-1)$ son los ra{\'\i}ces de $\Irr(z,T,K)$. En
particular $z$ y $z+1$ son conjugados sobre $K$ y $L=K(z)$ es de
Galois sobre $K$. Sean $G=\Gal(L/K)$ y $\sigma\in G$ tal que
$\sigma z=z+1$, de donde 
$\sigma^i z=z+i$ y $o(\sigma)=p$ por lo que
se tiene $G=\langle\sigma\rangle$ es un grupo c{\'\i}clico de orden $p$. $\fin$
\end{proof}

\begin{teorema}[Extensiones de Kummer\index{extensi\'on 
de Kummer}]\label{CClaseT1.1.4}
Sea $K$ un campo de caracter{\'\i}stica $p\geq 0$ y sea
$n\in{\ma N}$ tal que $p\nmid n$ (en el caso $p=0$, $n$ es arbitrario).
Supongamos que $\zeta_n\in K$ donde $\zeta_n$ es una ra{\'\i}z $n$--\'esima
primitiva de la unidad. Entonces $L/K$ es una extensi\'on c{\'\i}clica de 
grado $n\iff$ existe $z\in L$ tal que $L=K(z)$ e $\Irr(z,T,K)=T^n-a\in K[T]$,
esto es, $L=K(\sqrt[n]{a})$.
\end{teorema}
\begin{proof}

\noindent
$\Rightarrow)$ Sea $G=\Gal(L/K)=\langle\sigma\rangle$, $o(\sigma)=n$.
Se tiene $\N\zeta_n=\zeta_n^n=1$. Por lo tanto existe $z\in L$ tal que
$\sigma z=\zeta_n z$ y $\sigma^i z=\zeta_n^i z$, por lo que
$\sigma^i z=z\iff n|i$. Por lo tanto $z, \zeta_n z,\ldots, \zeta_n^{n-1} z$
son los distintos conjugados de $z$. Se sigue que $\Irr(z,T,K)=\prod_{i=0}^{
n-1}(T-\zeta_n^i z)$.

Por otro lado $\sigma z^n=(\sigma z)^n=(\zeta_n z)^n=\zeta_n^n z^n=
z^n$, esto es, $z^n=a\in K$ y $z,\zeta_n z,\ldots, \zeta_n^{n-1} z$ son 
ra{\'\i}ces de $T^n-a\in K[T]$. Por lo tanto $\Irr(z,T,K)=T^n-a$ y
$z^n=a$.

\noindent
$\Leftarrow)$ Para $a\neq 0$, $T^n-a$ es un polinomio separable debido
a que $p|n$ y tiene distintas ra{\'\i}ces $z,\zeta_n z,\ldots, \zeta_n^{n-1} z$
donde $z\in \bar{K}$ es tal que $z^n=a$, donde $\bar{K}$ es una 
cerradura algebraica de $K$. Se sigue que $L=K(z)$ es una
extensi\'on de Galois. Puesto que se supone que $T^n-a$ es irreducible
y $z$ y $\zeta_n z$ son conjugados. Por lo tanto existe $\sigma\in G
=\Gal(L/K)$ tal que $\sigma z=\zeta_n z$. Por tanto $o(\sigma)=n
=o(G)=[L:K]$ de donde se sigue que $L/K$ es una extensi\'on c{\'\i}clica
de grado $n$. $\fin$
\end{proof}

Para otra versi\'on del Teorema \ref{CClaseT1.1.5} a continuaci\'on,
ver el Corolario \ref{C.9'.2.2}.

\begin{teorema}\label{CClaseT1.1.5} Sea $K$ tal que $\ca K=p>0$ y sean
$L_i=K(z_i)/K$, $i=1,2$ dos extensiones c{\'\i}clicas de grado $p$ dadas por
$z_i^p-z_i=a_i\in K$, $i=1,2$. Lo siguiente es equivalente
\las
\item $L_1=L_2$.
\item $z_1=jz_2+b$ para $1\leq j\leq p-1$ y $b\in K$.
\item $a_1=ja_2+(b^p-b)=ja_2+ \wp(b)$ para $1\leq j\leq p-1$ y $b\in K$. Aqu{\'\i}
usamos la notaci\'on $\wp(b)=b^p-b$\label{CClaseoperadorartinschreier}.
\end{list}
\end{teorema}

\begin{proof} {\underline{(1) $\Leftrightarrow$ (2)}}.
Si $z_1=jz_2+b$, entonces $z_2=iz_1-ib$ con $ij\equiv 1\bmod p$, por lo
que $L_1=L_2$. Rec{\'\i}procamente, si $L_1=L_2$, entonces si $G=
\Gal(L_1/K)=\Gal(L_2/K)=\langle\sigma\rangle$ donde seleccionamos
$\sigma$ tal que $\sigma z_1=z_1+1$. Ahora bien, puesto que $\sigma z_2$
es conjugado de $z_2$ sobre $K$, entonces existe $1\leq i\leq p-1$ tal que 
$\sigma z_2=z_2+i$. Sea $1\leq j\leq p-1$ tal que $ij\equiv 1\bmod p$. Entonces
\[
\sigma (jz_2)=j\sigma(z_2)=jz_2+ji=jz_2+1.
\]
Por tanto $\sigma(z_1-jz_2)=z_1-jz_2$, por lo que $z_1-jz_2=b\in K$.

\noindent
$\underline{(2)\Rightarrow (3)}:$ 
Ahora si $z_1=jz_2+b$ se tiene $z_1^p-z_1=a_1=(jz_2+b)^p-(jz_2+b)=
j(z_2^p-z_2)+\wp(b)= ja_2+\wp (b)$.

\noindent
$\underline{(3)\Rightarrow (2)}:$
Rec{\'\i}procamente, si $a_1=ja_2+\wp(b)$, $z_1^p-z_1=(jz_2+b)^p-(jz_2+b)$,
es decir, $(z_1-(jz_2+b))^p-(z-(jz_2+b))=0$, por lo tanto $\omega=z_1-jz_2-b$
es una ra{\'\i}z de $\omega^p-\omega=0$ por lo que 
$\omega\in {\ma F}_p$. Se sigue que $z_1=jz_2+b+\omega$ y $\wp(
b+\omega)=\wp(b)$, por lo que $a_1=ja_2+\wp(b+\omega)$.
$\fin$
\end{proof}

Similarmente se puede probar

\begin{teorema}\label{CClaseT1.1.6} Sea $K$ un campo de caracter{\'\i}stica
$p\geq 0$ tal que $\zeta_n\in K$ con $p\nmid n$ y $\zeta_n$ una ra{\'\i}z
$n$--\'esima primitiva de la unidad. Sean $L_i=K(z_i)$, $i=1,2$, dos
extensiones c{\'\i}clicas de $K$ de grado $n$ dadas por $z_i^n=a_i\in K$.
Entonces lo siguiente es equivalente:
\las
\item $L_1=L_2$.
\item $z_1=z_2^jc$ para alg\'un $1\leq j\leq n-1$ con $\mcd(j,n)=1$ y $c\in K$.
\item $a_1=a_2^j c^n$ para alg\'un $1\leq j\leq n-1$ tal que $\mcd(j,n)=1$
y $c\in K$. $\fin$
\end{list}
\end{teorema}

Con respecto a la ramificaci\'on en las extensiones de Kummer y de Artin--Schreier,
los siguientes dos resultados se deben a Hasse. En este caso se tiene $K/k$
un campo de funciones y supondremos $k$ perfecto 
(en general consideraremos el caso $k=\F$, el cual es perfecto). Se tiene:

\begin{teorema}[H. Hasse]\label{CClaseT1.1.7}
Sea $k$ un campo perfecto de caracter{\'\i}stica $p>0$ y sea $\pK$ un lugar
fijo de $K$. Si $L/k$ es una extensi\'on c{\'\i}clica de grado $p$, entonces
existe $y\in L$ tal que $L=K(y)$ con $y^p-y=a$ tal que, 
o bien $v_\pK(a)\geq 0$,
o bien $v_\pK(a)=-\lambda<0$ y $p\nmid \lambda$.

Si $v_\pK(a)\geq 0$, entonces $\pK$ es no ramificado. Si $v_\pK(a)=
-\lambda<0$ y $p\nmid \lambda$, entonces $\pK$ es ramificado y el diferente
local est\'a dado por ${\eu D}_\pL=\pL^{(\lambda+1)(p-1)}$ donde $\pL$ es el
lugar de $L$ encima de $\pK$, es decir, $\pK=\pL^p$.
\end{teorema}

\begin{proof} Ver Proposici\'on \ref{P6.3-1.Ram5}, 
Observaci\'on \ref{O6.3-2.Ram6}
y \cite[Theorems 5.8.10 y 5.8.11]{Vil2006}. $\fin$
\end{proof}

\begin{teorema}[H. Hasse]\label{CClaseT1.1.8}
Sea $k$ un campo de caracter{\'\i}stica $p\geq 0$. Sea $L/K$ una extensi\'on
c{\'\i}clica de grado $n$ con $p\nmid n$ y tal que $\zeta_n\in k$ donde
$\zeta_n$ es una ra{\'\i}z $n$--\'esima
primitiva de la unidad. Sea $\pK$ un lugar fijo de
$K$. Entonces $L=K(y)$ tal que
 $y^n=a$ con $0\leq v_\pK(a)\leq n-1$. Se tiene
que $\pK$ es no ramificado en $L/K\iff v_\pK(a)=0$.

Si $v_\pK(a)=m>0$ y $\pL$ es un divisor de $L$ encima de $\pK$, tenemos 
$e(\pL|\pK)=\frac{n}{\mcd(n,m)}$ y $v_\pL({\eu D}_\pL)=\frac{n}{\mcd(n.m)}-1$,
donde ${\eu D}_\pL$ denota al diferente local.
\end{teorema}

\begin{proof} Proposici\'on \ref{P6.4-1.Ram7} y
\cite[Theorem 5.8.12]{Vil2006}. $\fin$
\end{proof}

\subsection{Automorfismo de Frobenius\index{automorfismo de Frobenius}
y s{\'\i}mbolo de Artin\index{simbolo de Artin@s\'imbolo de Artin}\index{Artin!s{\'\i}mbolo de
$\sim$}}\label{CClaseS1.3}

Dada una extensi\'on finita de campos finitos ${\ma F}_{q^d}/\F$, se tiene
que el grupo de Galois $G=\Gal({\ma F}_{q^d}/\F)$ es un grupo 
c{\'\i}clico de orden $d$. El generador $\tau\colon {\ma F}_{q^d}
\to {\ma F}_{q^d}$ dado por $\tau(x)=x^q$ se le llama el {\em
automorfismo de Frobenius}.\index{automorfismo de Frobenius}
\index{Frobenius!automorfismo de $\sim$} Se tiene que $\tau(x)=x
\iff x^q=x\iff x\in\F$.

\begin{definicion}\label{CClase1.3.0}\label{CClaseD1.3.-1}
Un {\em campo global\index{campo global}} es, o bien una 
extensi\'on finita de ${\ma Q}$, o bien un campo de funciones
con campo de constantes un campo finito.

Los campos globales de funciones, reciben tambi\'en el nombre de
{\em campos de funciones congruentes\index{campos de
funciones congruentes}}.
\end{definicion}

Ahora bien, dada una extensi\'on de Galois de campos globales
$L/K$, sea $\pK$ un primo de $K$ y sea $\pL$ un primo de $L$ sobre
$\pK$. Sea $L(\pL)/K(\pK)$ la extensi\'on de campos residuales.
Si $D=D(\pL|\pK)$ es el grupo de descomposici\'on de $\pL/\pK$, el
mapeo $D\xhookrightarrow \varphi \Gal(L(\pL)/K(\pK))$,
$\sigma\mapsto \bar{\sigma}$ = clase de 
$\sigma$ m\'odulo $\pK$, es un epimorfismo de grupos con n\'ucleo
$I=I(\pL|\pK)$ el grupo de inercia de $\pL/\pK$. Es decir, tenemos
$D/I\cong G$. Si $\pL/\pK$ es no ramificado, entonces $I=\{1\}$ y
por tanto $D\cong G$. En particular, existe un \'unico $\sigma_\pK \in
D$ tal que $\sigma_\pK\xrightarrow{\varphi} \tau$ es el automorfismo
de Frobenius. El automorfismo $\sigma_\pK$ se llama el
{\em automorfismo de Frobenius de $\pL/\pK$}.

\begin{definicion}\label{CClaseD1.3.0}
Se define la {\em norma absoluta\index{norma absoluta}\label{CClasenormaabsoluta}}
del campo residual $L(\pL)$ de un campo ya sea local o global
como $\N(\pL):=|L(\pL)|$.
\end{definicion}

Sean $\N(\pL)=|L(\pL)|$ y $\N(\pK)=|K(\pK)|$, digamos
$\N(\pK)=q^f$ con $f=[K(\pK):\F]$ en el caso de campos de funciones,
entonces $\sigma_\pL(\pL)=\pL$, $\sigma_\pL|_\pK=\Id$ y
$\sigma_\pL(\o_\pL)=\o_\pL$. Es decir,
$\sigma_\pL$ est\'a caracterizados por
\[
\overline{\sigma_\pL}(x)\equiv x^{\N(\pK)}\bmod \pL \quad \text{
para toda}\quad x\in \o_\pL.
\]

Escribimos $\sigma_\pL=\frobenius{L/K}{\pL}$\label{CClasefrobenius} y $\frobenius{L/K}{\pL}(x)
\equiv x^{\N(\pK)}\bmod \pL$ para toda $x\in \o_\pL$. 

\subsubsection{Propiedades del automorfismo de Frobenius}\label{CClaseS1.3.1}

Aqu\'i se supone que $L/K$ es una extensi\'on de Galois finita de
campos globales.

\begin{proposicion}\label{CClaseP1.3.1.1}
Si $\sigma\in
{\mathcal G}=\Gal(L/K)$, se tiene 
\[
\frobenius{L/K}{\sigma \pL}=
\sigma\frobenius{L/K}{\pL}\sigma^{-1}.
\]
\end{proposicion}

\begin{proof} Sea $y\in \o_L$ y sea $x\in \o_L$ tal que $y=\sigma^{-1}
x$. Entonces
\[
\frobenius{L/K}{\pL}y=\frobenius{L/K}{\pL}\sigma^{-1}x\equiv
(\sigma^{-1} x)^q\bmod \pL.
\]
Aplicando $\sigma$ a esta igualdad, obtenemos $\sigma
\frobenius{L/K}{\pL}\sigma^{-1}x\equiv x^q\bmod \sigma\pL$.
Por la unicidad del automorfismo de Frobenius se sigue que
$\sigma\frobenius{L/K}{\pL}\sigma^{-1}=\frobenius{L/K}{
\sigma\pL}$. $\fin$
\end{proof}

\begin{proposicion}\label{CClaseP1.3.1.2}
Sea una torre $K\subseteq E\subseteq L$. Sea ${\eu q}=\pL\cap E$.
Entonces
\[
\frobenius{L/K}{\pL}^{f({\eu q}|\pK)}=\frobenius{L/E}{\pL},
\]
donde $f({\eu q}|\pK)$ es el grado de inercia de ${\eu q}$ sobre
$\pK$, esto es, $f({\eu q}|\pK)=[\o_L/{\eu q}:\o_K/\pK]$.
\end{proposicion}

\begin{proof}
Se tiene $\F\subseteq {\ma F}_{q^{f_0}}\subseteq {\ma F}_{q^f}$ donde
$f_0=f({\eu q}|\pK)$ y $f=f(\pL|\pK)$.

El automorfismo de Frobenius generando $\Gal({\ma F}_{q^f}/
{\ma F}_{q^{f_0}})$ corresponde al mapeo $y\xrightarrow{\tau}y^{q^{f_0}}$.
Si $\sigma$ es el automorfismo de Frobenius de $\Gal({\ma F}_{q^f}/
\F)$, entonces $\tau =\sigma^{f_0}$. Ahora bien $\tau=\frobenius{
L/E}{\pL}$ y $\sigma=\frobenius{L/K}{\pL}$. $\fin$
\end{proof}

\begin{proposicion}\label{CClaseP1.3.1.3} Sea $K\subseteq E
\subseteq L$ donde $E/K$ es tambi\'en una extensi\'on
de Galois. Sea ${\eu q}=\pL\cap E$. Entonces
se tiene que la restricci\'on satisface
\[
\frobenius{E/K}{{\eu q}}=\rest_E\frobenius{L/K}{\pL}=
\frobenius{L/K}{\pL}\big|_E.
\]
Esto es,
$\rest_E \frobenius{L/K}{\pL}=\frobenius{E/K}{\pL\cap E}$.
\end{proposicion}

\begin{proof} Se tiene $\rest_E\colon\Gal(L/K)\lra \Gal(E/K)$ est\'a dada
por $\sigma\longmapsto \sigma|_E$ y se tiene que $\ker \rest_E=
H=\Gal(L/E)$.

Sea $\sigma=\frobenius{L/K}{\pL}$. Entonces $\sigma x\equiv
x^q\bmod \pL$ para $x\in\o_L$. Por tanto, si $x\in\o_E$, se tiene
$\sigma x-x^q\in \o_E$ y $\sigma x\equiv x^q\bmod \pL$. Se
sigue que $\sigma x\equiv x^q\bmod (\pL\cap \o_E)=x^q\bmod
{\eu q}$. As\'i, $\sigma|_E=\frobenius{E/K}{{\eu q}}$. $\fin$
\end{proof}

\begin{corolario}\label{CClaseC1.3.1.4} Sean $K\subseteq E, F
\subseteq L$ tales que $E/K$ y $F/K$ son extensiones de Galois.
Sean ${\eu q}=\pL\cap E$ y ${\eu t}=\pL\cap F$. Supongamos que
$L=EF$. Entonces el mapeo $\sigma\lra (\sigma|_E, \sigma|_F)$
de $\Gal(L/K)$ en $\Gal(E/K)\times \Gal(F/K)$, el cual es inyectivo,
da lugar a
\[
\begin{array}{ccccc}
\frobenius{L/K}{\pL}&\lra&\frobenius{E/K}
{{\eu q}}&\times&\frobenius{F/K}{{\eu t}}\\
\uigual &&\uigual&&\uigual\\
\frobenius{EF/K}{\pL}&\lra&\frobenius{E/K}{\pL\cap E}&\times&
\frobenius{F/K}{\pL\cap F}.
\end{array}
\]
\end{corolario}

\begin{proof} Se sigue inmediatamente 
de la Proposici\'on \ref{CClaseP1.3.1.3}. $\fin$
\end{proof}

\begin{proposicion}\label{CClaseP1.3.1.5} Se tiene $\pK$ in $\o_K$ se
descompone totalmente en $L/K$ si y s\'olo si $\frobenius{L/K}
{\pL}=1$.
\end{proposicion}

\begin{proof} Puesto que $\pL|\pK$ es no ramificada, $\pK$ se descompone
totalmente si y s\'olo si $f(\pL|\pK)=1$ si y s\'olo si 
\begin{gather*}
\o_L/\pL=\o_K/\pK
\iff \Gal((\o_L/\pL)/(\o_K/\pK))=\{1\}\iff \frobenius{L/K}{\pL}=1. \tag*{$\fin$}
\end{gather*}
\end{proof}

\begin{corolario}\label{CClaseC1.3.1.6} Sean $E/K$ y $F/K$ extensiones
de Galois y $L=EF$. Entonces $\pK$ en $\o_K$ se descompone
totalmente en $L/K$ si y s\'olo si $\pK$ se descompone totalmente
tanto en $E/K$ como en $F/K$.
\end{corolario}

\begin{proof} De la Proposici\'on \ref{CClaseP1.3.1.5} se tiene que $\pK$ se descompone
totalmente en $L/K \iff \frobenius{L/K}{\pL}=1$ si y s\'olo si bajo el
mapeo $\frobenius{L/K}{\pL}\xhookrightarrow{\phantom{xxxx}}
\frobenius{E/K}{\pL\cap E}\times\frobenius {F/K}{\pL\cap F}$, se
tiene $\frobenius{E/K}{\pL\cap E}=\frobenius{F/K}{\pL\cap F}=1
\iff \pK$ se descompone totalmente tanto en $E/K$ como en $F/K$. $\fin$
\end{proof}

Sea $L/K$ una extensi\'on abeliana finita de campos globales. En este
caso, para $\sigma\in\Gal(L/K)$, $\frobenius{L/K}{\pL}=\sigma
\frobenius{L/K}{\pL}\sigma^{-1}=\frobenius{L/K}{\sigma \pL}$. Esto
es, si $L/K$ es abeliana, $\frobenius{L/K}{\pL}$ no depende de
$\pL$ sino \'unicamente de $\pK$. En este caso escribimos
\[
\frobenius{L/K}{\pL}\label{CClaseartin}=\artinp{L/K}{\pK}=(L/K,\pK)=
(\_, L/K,\pK)=(\_,L/K)=\psi_{L/K}(\pK),
\]
el cual se llama el {\em s{\'\i}mbolo de Artin}.\index{simbolo de Artin@s\'imbolo
de Artin}\label{CClasesimbolodeArtin}

\begin{observacion}\label{CClaseO1.3.1}
Se tiene que $(L/K,\pK)=1\iff (\pL|\pK$ es no ramificada y $L(\pL)=
K(\pK))\iff (\pK$ es totalmente descompuesto en $L$).

En general se tiene que $o((L/K,\pK))=d_{L/K}(\pL|\pK)=
[L(\pL):K(\pK)]$ el cual es el orden del automorfismo de
Frobenius.
\end{observacion}

\begin{teorema}\label{CClaseT1.3.1.7}
Sea $L/K$ una extensi\'on abeliana finita de campos globales.
Sea $E$ una extensi\'on finita de $K$ y sea $F/E$ una extensi\'on
abeliana finita tal que $L\subseteq F$ (por tanto $LE\subseteq F$).
Sea $\N_{E/K}\colon E\lra K$ la norma de $E$ en $K$. Sea 
$\theta=\rest\colon \Gal(F/E)={\mc G}\lra\Gal(L/K)=G$ la
restricci\'on. Sea ${\eu q}$ un primo en $E$ no ramificado en
$F$ y sea $\pK=\N_{E/F}{\eu q}$ un primo en $K$ no ramificado
en $L$. Sea ${\eu t}$ un primo en $F$ sobre ${\eu q}$ y $\pL$
un primo en $L$ sobre $\pK$. 

Entonces $\artinp{L/K}{\N_{E/K}{\eu q}}=\rest_L\artinp{F/E}{
{\eu q}}=\artinp{F/E}{{\eu q}}\Big|_L$.

En otras palabras, si $S$ es un conjunto finito de primos en
$K$ que contienen a todos los primos infinitos y a todos los
primos ramificados y si $S'$ es el conjunto de primos de $E$
que est\'an sobre $S$, entonces el diagrama
\[
\xymatrix{
D_E^{S_0'}\ar[rr]^{\psi_{F/E}}\ar[d]_{\N_{F/K}}&&{\mc G}
\ar[d]^{\rest_L}\\D_K^{S_0}\ar[rr]_{\psi_{L/K}}&&G
}
\]
es conmutativo, donde $D_E^{S_0'}$ es el grupo libre generado
por los divisores primos de $E$ que no est\'an en $S'$ y
an\'alogamente $D_K^{S_0}$. Aqu\'i $\psi_{F/E}$ y $\psi_{L/K}$
denotan los mapeos de Artin.
\end{teorema}

\begin{proof}
Sea $f=f({\eu q}|\pK)=f(E|F)=[\tilde{E}:\tilde{K}]=[\o_E/{\eu q}:
\o_K/\pK]$. Entonces $\N_{E/K}=\pK^f$.

Sea $\sigma=\psi_{F/E}({\eu q})=\artinp{F/E}{{\eu q}}=
\frobenius{F/E}{{\eu Q}}$ donde ${\eu Q}$ es un primo en $F$ sobre
${\eu q}$ y sea $\tau=\psi_{L/K}(\pK)=\artinp{L/K}{\pK}=
\frobenius{L/K}{\pL}$ con $\pL$ un primo en $L$ sobre $\pK$.
\[
\xymatrix{
&F\ar@{-}[d]\ar@/^2pc/@{-}[dd]^{\mc G}\\
L\ar@{-}@/_2pc/[dd]_G\ar@{-}[r]\ar@{-}[d]&LE\ar@{-}[d]\\
L\cap E\ar@{-}[d]\ar@{-}[r]&E\ar@{-}[dl]\\ K
}
\]
Entonces $\psi_{L/K}(\N_{E/K}{\eu q})=\psi_{L/K}(\pK^f)=
\psi_{L/K}(\pK)^f=\tau^f$ y $\rest_L\psi_{F/E}({\eu q})=
\artinp{F/E}{{\eu q}}\Big|_L=\sigma|_L$.

Debemos probar que $\sigma =\tau^f$. Ahora bien, $\o_L/\pL
\subseteq \o_F/{\eu Q}$. Para $x\in \o_L/\pL$ se tiene 
$\sigma x=x^{\N({\eu Q})}$. Puesto que $\N({\eu Q})=|\o_E/
{\eu Q}|$ y $[\o_E/{\eu Q}:\o_K/\pK]=f$, se sigue que
$|\o_E/{\eu Q}|=|\o_K/\pK|^f$, esto es, $\N({\eu Q})=\N(\pK)^f$.

Por tanto $\sigma x=x^{\N({\eu Q})}=x^{\N(\pK)^f}=\tau^f x$
puesto que $\tau x=x^{\N(\pK)}$ por lo que $\tau^f x=
x^{\N(\pK)^f}$. De esta forma obtenemos que $\sigma =\tau^f$. $\fin$
\end{proof}

Siempre que usamos el s{\'\i}mbolo de Artin o el automorfismo de Frobenius,
estaremos suponiendo que $\pL|\pK$ es no ramificado.

Recordemos que si $K=\F(T)$ y $K_M=K(\Lambda_M)$ es una
extensi\'on de campos de funciones ciclot\'omicas con $M\in R_T=\F[T]$, 
entonces si $\lambda_M$ es un generador del $R_T$--m\'odulo $
\Lambda_M$ y si $P\in R_T$ es un polinomio m\'onico e irreducible, 
con $P\nmid M$, entonces el s{\'\i}mbolo de Artin est\'a dado por
$(K_M/k,P)=(\lambda_M\mapsto \lambda_M^P)$ (Teorema \ref{T6.3.3}).

\subsection{Extensiones de Galois infinitas}\label{CClaseS1.4}

La base del siguiente desarrollo se hizo en el Cap\'itulo \ref{ch2}.
Recordemos los hechos fundamentales.

Sea $\Omega/k$ una extensi\'on de Galois y sea $G=\Gal(\Omega/k)$.
A $G$ se le da la {\em topolog{\'\i}a de Krull\index{topolog{\'\i}a de Krull}}
la cual se define como sigue. Para $\sigma\in G$, tomamos las clases
$\sigma\Gal(\Omega/K)$ como una base de vecindades de $\sigma$
donde $K/k$ recorre el conjunto de todas las extensiones finitas de 
Galois con $k\subseteq K\subseteq \Omega$. Se tiene que las operaciones
de grupo:
\[
\begin{array}{rclcccccrcl}
G\times G&\to &G&&&{\qquad}\text{y}\qquad&&&G&\to&G\\
(\sigma,\varphi)&\mapsto& \sigma\varphi&&&{\qquad}&&&
\sigma&\mapsto&\sigma^{-1}
\end{array}
\]
son funciones continuas en la topolog{\'\i}a de Krull.
De esta forma $G$ se hace un 
{\em grupo topol\'ogico\index{grupo topol\'ogico}}.
Cuando $G$ es finito, la topolog{\'\i}a de Krull es la topolog{\'\i}a
discreta, es decir todos los subconjuntos son abiertos y todos
son cerrados.

\begin{proposicion}\label{CClaseP1.4.1} Si $\Omega/k$ es una extensi\'on de
Galois, finita o infinita, $G=\Gal(\Omega/k)$ es Hausdorff y compacto.
\end{proposicion}
\begin{proof} Teorema \ref{T6.2}.
$\fin$
\end{proof}

\begin{teorema}[Correspondencia de Galois\index{Galois!correspondencia
de $\sim$}]\label{CClaseT1.4.2}
 Sea $\Omega/k$ una extensi\'on de Galois, finita o infinita.
Entonces $K\to \Gal(\Omega/K)$ da lugar a una correspondencia
biyectiva entre todas las subextensiones $K/k$ con $k\subseteq K\subseteq
\Omega$ y todos los subgrupos cerrados de $G=\Gal(\Omega/k)$.

Los subgrupos abiertos de $G$ corresponden a las subextensiones
finitas $k/k$ de $\Omega/k$, es decir $[K:k]<\infty$.
\end{teorema}
\begin{proof} Teorema \ref{T6.3}.
$\fin$
\end{proof}

\begin{observacion}\label{CClaseO1.4.2'} Dado cualquier subgrupo $H<G$,
el campo fijo de $H$, $\Omega^H$ es efectivamente un campo. De
hecho, $\Omega^H=\Omega^{\bar{H}}$ donde $\bar{H}$ es la
cerradura de $H$.

En efecto, el verificar que $\Omega^H$ es un campo es rutinario
pues en esta parte no importan las caracter{\'\i}sticas topol\'ogicas
de $H$. Ahora bien, puesto que $H\subseteq \bar{H}$ se tiene
$\Omega^{\bar{H}}\subseteq \Omega^H$. Por otro lado, si
$x\in \Omega^H$, entonces $\sigma x=x$ para toda $\sigma \in H$.
Dado $\psi\in\bar{H}$, entonces, puesto que estamos considerando
la topolog{\'\i}a de Krull, si $N<G$ es de {\'\i}ndice finito, $\psi N\cap
H\neq \emptyset$. 
Se tiene $[K(x):K]<\infty$ y si $\widetilde{K(x)}$ es la cerradura
de Galois de $K(x)/K$, $[\widetilde{K(x)}:K]<\infty$. Sea $N=
\Gal(\Omega/\widetilde{K(x)})\lhd G$ y $[G:N]=[\widetilde{K(x)}:K]<
\infty$. Por tanto existe $n\in N$, $h\in H$ tal que $\psi n=h$. 
Por tanto
\[
(\psi n) x\underbracket[0pt]{=}_{\substack{\uparrow\\ nx=x}}
\psi (x)= h(x)=x,
\]
por lo que $x\in \Omega^{\bar{H}}$ y por tanto \fbox{$\Omega^H=
\Omega^{\bar{H}}$}.
\end{observacion}

Recordemos lo siguiente discutido en el Cap\'itulo \ref{ch2}.
En el caso de un grupo profinito $G$, el conjunto de {\'\i}ndices lo tomamos
como $I:=\{N\mid N\lhd G, N \text{\ es abierto en \ } G\}$. Definimos
$N\leq N_1\iff N_1\subseteq N$ y $G_N:=G/N$. Entonces, para
$N\leq N_1$, $\phi_{N_1,N}\colon G/N_1\to G/N$, $gN_1\mapsto gN$ y
se tiene que $G\cong \lim\limits_{\substack{\leftarrow\\ N}}G/N$.

Para calcular el grupo de Galois de una extensi\'on $L/K$ de campos,
donde $L/K$ es una extensi\'on algebraica, normal y separable, 
procedemos de la siguiente forma. Para
 $\alpha\in L$ se tiene
$[K(\alpha):K]<\infty$ y si $\widetilde{K(\alpha)}/K$ es la cerradura
normal de $K(\alpha)/K$, entonces $\widetilde{K(\alpha)}/K$ es
una extensi\'on de Galois finita. Definimos $L_{\alpha}:=
\widetilde{K(\alpha)}$. Seleccionamos los $\beta\in L$ tales que
$L_{\beta}=K(\beta)$ y si $I=\{\alpha\in L\mid K(\alpha)/K \text{\ es
una extensi\'on de Galois (necesariamente es finita)}\}$, entonces
$L=\cup_{\beta\in I}L_{\beta}=\lim\limits_{\substack{\rightarrow\\ \beta}}
L_{\beta}$ donde definimos $\beta\leq \beta_1\iff L_{\beta}
\subseteq L_{\beta_1}$ y donde $\psi_{\beta,\beta_1}\colon L_{\beta}
\lra L_{\beta_1}$ es el encaje.

\begin{teorema}\label{CClaseT1.4.8}
Con las notaciones anteriores, se tiene que 
\begin{gather*}
\Gal(L/K)=\Gal\big(\lim\limits_{\substack{\rightarrow\\ \beta\in I}}
L_\beta/K\big)\cong \lim\limits_{\substack{\leftarrow\\ \beta\in I}}
\Gal(L_{\beta}/K)
\end{gather*}
\end{teorema}
\begin{proof} Ver Teorema \ref{T6.2}. $\fin$
\end{proof}

\begin{observacion}\label{CClaseO1.4.9} Notemos que si $L_\beta\subseteq
L_{\beta_1}$, entonces $\Gal(L_{\beta}/K)\cong
\frac{\Gal(L_{\beta_1}/K)}{\Gal(L_{\beta_1}/L_{\beta})}$.

Definimos $G_{\beta}:=\Gal(L_{\beta}/K)$ y $G_{\beta}\leq G_{\beta_1}
\iff G_{\beta}$ es un grupo cociente de $G_{\beta_1}$ bajo el mapeo
de restricci\'on y definimos 
\begin{eqnarray*}
\psi_{\beta_1,\beta}\colon\Gal(L_{\beta_1}/K)& \longrightarrow & \Gal(L_{\beta}/K)\\
\sigma &\longmapsto &\sigma|_{L_{\beta}}.
\end{eqnarray*}

El isomorfismo anunciado en el Teorema \ref{CClaseT1.4.8} est\'a dado por
\begin{eqnarray*}
\Gal(L/K)&\longrightarrow&\lim\limits_{\substack{
\leftarrow\\ \beta\in I}}\Gal(L_{\beta}/K)\\
\sigma&\longmapsto & \{\sigma|_{L_{\beta}}\}_{\beta\in I}.
\end{eqnarray*}
\end{observacion}

\begin{ejemplo}[Ejemplos \ref{Ej5.10'}]\label{CClaseE1.4.10}
Sean $p$ un n\'umero primo, $m,n\in{\ma N}$, $m\geq n$ y
\[
\phi_{m,n}\colon {\ma Z}/p^m{\ma Z}\longrightarrow {\ma Z}/p^n {\ma Z}
\]
dada por $\phi_{m,n}(x\bmod p^m)=x\bmod p^n$. Se define
${\ma Z}_p=\big\{\sum_{n=0}^{\infty} a_n p^n\mid a_n\in
\{0,1,\ldots,p-1\}\big\}$ y sea $\varphi_m\colon {\ma Z}_p\to
{\ma Z}/p^m {\ma Z}$, dada por 
\[
\varphi_m\big(\sum_{n=0}^{\infty}
a_n p^n\big) = \sum_{n=0}^{m-1} a_n p^n \bmod p^m.
\]

De la discusi\'on anterior se puede probar que 
${\ma Z}_p\cong \lim\limits_{\substack{\leftarrow\\ m}}
{\ma Z}/p^m {\ma Z}$.

Sea ${\ma Q}_n$ la subextensi\'on de grado $p^n$ de
${\ma Q}(\zeta_{p^{n+1}})/{\ma Q}$ cuando $p$ es impar y
si $p=2$, ${\ma Q}_n$ es la subextensi\'on c{\'\i}clica de
grado $2^n$ de ${\ma Q}(\zeta_{2^{n+2}})/{\ma Q}$ contenida en
${\ma R}$. Sea ${\ma Q}_{\infty}=\cup_{n=1}^{\infty}
{\ma Q}_n$, entonces ${\ma Q}_{\infty}/{\ma Q}$ es de Galois
y
\begin{align*}
\Gal({\ma Q}_{\infty}/{\ma Q})&= \Gal\big(\cup_{n=1}^{\infty}
{\ma Q}_n/{\ma Q}) \cong\Gal(\lim\limits_{\substack{\rightarrow\\ n}}
{\ma Q}_n/{\ma Q}\big)\\
&\cong \lim\limits_{\substack{\leftarrow\\ n}}\Gal({\ma Q}_n/{\ma Q})
\cong \lim\limits_{\substack{\leftarrow\\ n}}{\ma Z}/p^n{\ma Z}
\cong {\ma Z}_p.
\end{align*}

${\ma Q}_{\infty}/{\ma Q}$ recibe el nombre de la
{\em ${\ma Z}_p$--extensi\'on ciclot\'omica de 
${\ma Q}$\index{Z@${\ma Z}_p$--extensi\'on ciclot\'omica}}
(ver Cap\'itulo \ref{Ch7}).
\end{ejemplo}

\begin{ejemplo}\label{CClaseE1.4.11}
Sea $\F$ el campo finito de $q$ elementos, $q=p^r$. Si
$K/\F$ es una extensi\'on de grado $d$, $K\cong
{\ma F}_{q^d}$ y $\Gal(K/\F)\cong {\ma Z}/d{\ma Z}$. Por
tanto $\overline{\F}=\abe {\F}=\cup_{n=1}^{\infty}{\ma F}_{q^n}$ y
\begin{align*}
\Gal(\overline{\F}/\F)&= \Gal\big(\cup_{n=1}^{\infty} {\ma F}_{q^n}\big/
\F\big)\cong \Gal(\lim\limits_{\substack{\rightarrow\\ n}} {\ma F}_{q^n}/\F)\\
&\cong \lim\limits_{\substack{\leftarrow\\ n}}\Gal({\ma F}_{q^n}/\F)
\cong \lim\limits_{\substack{\leftarrow\\ n}} {\ma Z}/n{\ma Z}
\cong \hat{{\ma Z}}\cong \prod_p {\ma Z}_p
\end{align*}
donde $\hat{{\ma Z}}$ es el {\em anillo o grupo de 
Pr\"ufer\index{anillo de Pr\"ufer}\index{grupo de Pr\"ufer}} y
de hecho es la completaci\'on de ${\ma Z}$.
(ver comentario despu\'es del Teorema \ref{CClaseT4.8.2} y
\cite[Ch. X]{AtiMac69}).
\end{ejemplo}

\subsection{Teor{\'\i}a de Kummer\index{teor{\'\i}a de Kummer}
\index{Kummer!teor{\'\i}a de $\sim$}}\label{CClaseS1.6}

La Teor\'ia de Kummer juega un papel relevante en el {\em Teorema
de Existencia} en la teor\'ia de campos de clase global (ver Subsecci\'on
\ref{STeoremadeExistencia} y Teorema \ref{T17.6.160N}).

Sea $K$ un campo de caracter{\'\i}stica $p\geq 0$ y $n\in {\ma N}$ tal
que $p\nmid n$. Sea $W_n$ el {\em grupo de las $n$--ra{\'\i}ces de 
unidad\index{grupo de las $n$--ra\'ices de unidad}}.
Se supone que $W_n\subseteq K$. Una {\em extensi\'on de
Kummer\index{extensi\'on de Kummer}\index{Kummer!extensi\'on de $\sim$}}
de $K$ de la forma $K(\sqrt[n]{\Delta}):=K(\{\sqrt[n]{\alpha}\mid \alpha
\in \Delta\})$ donde $\Delta$ es un subgrupo de $\* K$ tal que 
$(\*K)^n \subseteq \Delta \subseteq \* K$.

Si $L/K$ es una extensi\'on de Kummer, $L/K$ es una extensi\'on abeliana
de exponente $n$, es decir, $\sigma^n=\Id_K$ para toda $\sigma \in
\Gal(L/K)$, aunque $L/K$ no es necesariamente finita.

\begin{proposition}\label{CClaseP1.6.1} Si $L/K$ es una extensi\'on abeliana
no necesariamente finita de exponente $n$, entonces $L=K(\sqrt[n]{\Delta})$
con $\Delta =(\* L)^n\cap \* K$, es decir, $\sqrt[n]{\Delta}=\* L\cap
\sqrt[n]{\* K}$.
\end{proposition}

\begin{proof} Por definici\'on se tiene que $K(\sqrt[n]{\Delta})\subseteq L$.
Ahora bien, $L$ es la composici\'on de sus subextensiones c{\'\i}clicas.
Sea $M/K$ una subextensi\'on c{\'\i}clica de $L/K$, por lo tanto
$\Gal(M/K)$ tiene orden un divisor de $n$, por lo que
$M=K(\sqrt[n]{a})$ con $a\in (\* L)^n\cap \*K$ (Teorema \ref{CClaseT1.1.4}),
es decir, con $a\in K$ tal que $\sqrt[n]{a}\in L$. Se sigue que
$M\subseteq K(\sqrt[n]{\Delta})$ y 
por tanto $L\subseteq K(\sqrt[n]{\Delta})$. $\fin$
\end{proof}

\begin{teorema}[Teor{\'\i}a de
Kummer\index{Kummer!teor{\'\i}a de $\sim$}\index{teor{\'\i}a de
Kummer}]\label{CClaseT1.6.2}
Las extensiones de Kummer de exponente $n$
est\'an en correspondencia biyectiva
con los subgrupos $\Delta$ de $\*K$ que contienen a $(\*K)^n$. Si
$L=K(\sqrt[n]{\Delta})$, entonces $\Delta=(\*L)^n\cap \*K$ y se
tiene el isomorfismo
\[
\widehat{\Gal(L/K)}=\Hom(\Gal(L/K),W_n)\cong \Delta/(\*K)^n.
\]
El isomorfismo proviene de que dado $a\bmod (\*K)^n\in \Delta/(\*K)^n$,
se le asocia el caracter $\chi_a\colon \Gal(L/K)\to W_n$ dado por
$\chi_a(\sigma)=\frac{\sigma(\sqrt[n]{a})}{\sqrt[n]{a}}$.

M\'as expl\'icitamente, la correspondencia de Kummer est\'a dada por
\begin{eqnarray*}
\*L &\longmapsto & \Delta=(\*L)^n\cap \*K,\\
\Delta &\longmapsto & L=K\big(\sqrt[n]{\Delta}\big)=K_{\Delta}.
\end{eqnarray*}
Adem\'as se tiene
\begin{eqnarray*}
\Delta_1\Delta_2&\longleftrightarrow & K_{\Delta_1}K_{\Delta_2};\\
\Delta_1\cap \Delta_2&\longleftrightarrow & K_{\Delta_1}\cap K_{\Delta_2}.
\end{eqnarray*}
\end{teorema}

\begin{observacion}\label{CClaseO1.6.3} Si $L/K$ es una extensi\'on 
de Kummer de exponente $n$ infinita,
entonces $\Gal(L/K)$ tiene la topolog{\'\i}a de Krull y $\Hom(
\Gal(L/K),W_n)$ es el grupo de todos los homomorfismos 
{\underline{continuos}} $\chi\colon \Gal(L/K)\lra W_n$ y donde
$W_n$ tiene la topolog{\'\i}a discreta.

Si $\Gal(L/K)$ es finito, $\Gal(L/K)$ tiene la topolog{\'\i}a discreta
y todo homomorfismo $\chi\colon\Gal(L/K)\lra W_n$ es autom\'aticamente
continuo.
\end{observacion}

\begin{proof} (Teorema \ref{CClaseT1.6.2}).
Sea $L/K$ una extensi\'on de Kummer de exponente $n$. Entonces
$L=K(\sqrt[n]{\Delta})$ con $\Delta=(\* L)^n\cap \* K$ (Proposici\'on
\ref{CClaseP1.6.1}). Se define el homomorfismo 
\[
\Delta \xrightarrow{\ \theta\ } \Hom(\Gal(L/K),W_n),\quad a\mapsto \chi_a
\]
definido por $\chi_a(\sigma)=\frac{\sigma(\sqrt[n]{a})}{\sqrt[n]{a}}$
con $\sigma\in \Gal(L/K)$. Ahora bien se tiene
\begin{gather*}
a\in\ker \theta\iff \chi_a(\sigma)=\frac{\sigma(\sqrt[n]{a})}{\sqrt[n]{a}}=1
\text{\ para toda\ } \sigma \in \Gal(L/K) \\
\iff \sqrt[n]{a}\in \* K \iff a\in (\* K)^n.
\end{gather*}
Esto es, $\ker \theta =(\* K)^n$. El homomorfismo $\tilde{\theta}$ inducido
por $\theta$
\[
\Delta/(\* K)^n\xrightarrow{\ \tilde{\theta}\ } \Hom(\Gal(L/K),W_n)
\]
es inyectivo.

Para demostrar que $\tilde{\theta}$ es suprayectivo, primero consideremos
el caso en que $L/K$ es una extensi\'on finita. Sea $\chi\in\Hom(\Gal(
L/K),W_n)$. Entonces $\chi\colon\Gal(L/K)\lra\* L$ es en particular
un homomorfismo cruzado. Por el Teorema \ref{CClaseT1.5.7} se tiene que 
existe $y\in\* L$ tal que $\chi(\sigma)=\frac{\sigma y}{y}$ para toda
$\sigma \in \Gal(L/K)$.

Ahora bien, $\sigma(y^n)=(\sigma y)^n=\chi(\sigma)^n y^n= y^n$
para toda $\sigma\in\Gal(L/K)$ lo cual implica que $y^n\in \* K$.
Por tanto $y^n=x\in (\* L)^n\cap \* K= \Delta$ y 
\[
\chi_x(\sigma)=\frac{\sigma(\sqrt[n]{x})}{\sqrt[n]{x}}=
\frac{\sigma(y)}{y}=\chi(\sigma)\quad \text{para toda}\quad
\sigma \in\Gal(L/K),
\]
de donde se sigue que $\chi=\chi_x=\tilde{\theta}(x)$.

Ahora consideremos una extensi\'on infinita $L/K$. Sea $\{\Delta_{\alpha}
/(\* K)^n\}_{\alpha\in{\mc A}}$ el conjunto de los subgrupos finitos
de $\Delta/(\* K)^n$. Sea $K_{\alpha}:=K(\sqrt[n]{\Delta_{\alpha}})$,
$\alpha\in{\mc A}$. Se tiene
\[
\Delta/(\* K)^n=\bigcup_{\alpha\in{\mc A}} \Delta_{\alpha}/(\* K)^n
\quad\text{y}\quad L=\bigcup_{\alpha\in {\mc A}} K_{\alpha}.
\]

Se sigue que $\{\Gal(L/K_{\alpha})\}_{\alpha\in{\mc A}}$ forman una
base de vecindades de $\Id\in\Gal(L/K)$.

Dado $\chi\colon \Gal(L/K)\lra W_n$ continuo, $\ker \chi=\chi^{-1}(
\{1\})$ es cerrado y de \'indice finito en $\Gal(L/K)$ y por tanto
es abierto. Se sigue que existe $\alpha\in{\mc A}$ tal que $\Gal(L/
K_{\alpha})\subseteq \ker \chi$.

Ahora bien, $\chi$ induce un homomorfismo $\bar{\chi}\colon
\Gal(K_{\alpha}/K)\lra W_n$ de tal forma que $\chi(\sigma)=
\bar{\chi}(\sigma|_{K_{\alpha}})$. Puesto que $\Gal(K_{\alpha}/K)$
es finito, por lo anterior, existe $a\in K_{\alpha}$ tal que $\bar{\chi}
=\chi_a\colon\Gal(K_{\alpha}/K)\lra W_n$. Por tanto
\[
\chi(\sigma)=\bar{\chi}(\sigma|_{K_{\alpha}})=\frac{\sigma(
\sqrt[n]{a})}{\sqrt[n]{a}}=\chi_a(\sigma),
\]
por lo que $\chi=\chi_a$ y $\tilde{\theta}$ es suprayectiva.

Sea ahora $(\* K)^n\subseteq \Delta\subseteq \* K$ y $L=
K(\sqrt[n]{\Delta})$. Entonces veamos que $\Delta = (\* L)^n
\cap \* K$. Sea $\Delta_1=(\* L)^n\cap \* K$. Entonces
$\Delta \subseteq \Delta_1$. Ahora bien, por lo anteriormente
probado, se tiene
\begin{gather*}
\Delta_1/\* K\cong \Hom(\Gal(L/K), W_n).
\intertext{Al subgrupo $\Delta/(\* K)^n\subseteq \Delta_1/(\* K)^n$ le
corresponde el subgrupo}
\Hom(\Gal(L/K)/H, W_n) \text{\ de\ }
\Hom(\Gal(L/K),W_n) {\text{\ donde\ }}\\
H=\{\sigma\in\Gal(L/K)\mid
\chi_a(\sigma)=1 \text{\ para toda\ } a\in \Delta\}.
\end{gather*}

Puesto que $\sigma(\sqrt[n]{a})=\chi_a(\sigma)\sqrt[n]{a}$, obtenemos
que el subgrupo $H$ fija a los elementos de $\sqrt[n]{\Delta}$.
Ya que $L=K(\sqrt[n]{\Delta})$, se tiene que $H$ fija a $L$ y
por tanto $H=\{1\}$. Se sigue que
$\Hom(\Gal(L/K)/H, W_n)=\Hom(\Gal(L/K),W_n)$ por lo que
$\Delta/(\* K)^n=\Delta_1/(\* K)^n$, de donde se sigue que $\Delta
=\Delta_1$.

De esta forma, el mapeo $\Delta\lra L=K(\sqrt[n]{\Delta})$ es la
correspondencia dada en el enunciado del teorema.
$\fin$
\end{proof}

La composici\'on de dos extensiones de Kummer de exponente $n$
es nuevamente una extensi\'on de Kummer de exponente $n$ y por lo
tanto todas las extensiones de Kummer de exponente $n$ est\'an
contenidos en la extensi\'on de Kummer de exponente $n$ m\'axima:
$\tilde{K}=K(\sqrt[n]{\*K})$ y se tiene
\[
\Hom(\Gal(\tilde{K}/K),W_n)\cong \frac{\*K}{(\*K)^n}.
\]

\subsection{Historia de la teor{\'\i}a de campos de clase}\label{CClaseC2}

En esta subsecci\'on pretendemos dar un panorama m\'as o menos
general de lo que trata la teor{\'\i}a de campos de clase.

Lo tratado en esta peque\~na subsecci\'on es un extracto de
\cite{Con}.

\subsubsection{?`Que es la teor{\'\i}a de campos de clase?}\label{CClaseS2.1}

La teor{\'\i}a de campos de clase es la descripci\'on de las extensiones
abelianas de {\em campos globales\index{campos globales}} 
(extensiones finitas de ${\ma Q}$ y campos
de funciones con campo de constantes un campo finito $\F$) y de 
{\em campos locales\index{campos locales}} 
(extensiones finitas del campo de los racionales
n\'umeros $p$--\'adicos ${\ma Q}_p$
y las series de Laurent ${\ma F}((x))$ donde ${\ma F}$ es un campo finito).
La raz\'on de llamar a uno de estos campos {\em un campo de clase\index{campo
de clase}} se refiere a que estos campos est\'an relacionados a grupos
de clases de ideales. Uno de los teoremas principales es que los campos de
clase son los mismos que las extensiones abelianas.

Para una extensi\'on de campos $L/K$ ponemos 
\[
\Spl(L/K)=\{
\pK\text{\ lugar de\ } K\mid \pK\text{\ se descompone totalmente en\ } L\}
\index{Spl}\label{CClaseSpl}.
\]

Podemos pensar que la teor{\'\i}a de campos de clase se origina con el
trabajo de Kronecker y m\'as espec{\'\i}ficamente con el {\em Teorema de
Kronecker--Weber\index{Kronecker--Weber!teorema de $\sim$}\index{teorema
de Kronecker--Weber}} (1853): toda extensi\'on abeliana finita de 
${\ma Q}$ est\'a contenida en alg\'un campo ciclot\'omico. La primera
demostraci\'on y completa y correcta del Teorema de Kronecker--Weber
la dio Hilbert en 1896 (sin \'el saber que las anteriores ten{\'\i}an alguna
laguna, ver Cap\'itulo \ref{Ch1}).

Se tiene que si $L/K$ es una extensi\'on de Galois de campos num\'ericos,
$\Spl(L/K)$ tiene densidad $1/[L:K]$ (Kronecker). Como consecuencia se
tiene (Bauer): sean $L_1,L_2$ dos extensiones de Galois de $K$
entonces $L_1\subseteq L_2\iff \Spl(L_2/K)\subseteq \Spl(L_1/K)$
(salvo un n\'umero finito de primos) y en particular $L_1=L_2
\iff \Spl(L_1/K)=\Spl(L_2/K)$ con la igualdad salvo un n\'umero finito
de primos.

El t\'ermino {\em campo de clase\index{campo de clase}} fue introducido
en 1891 por Weber. En 1897 Weber extendi\'o el concepto de grupo
de clases de ideales: para un campo num\'erico $K$ y un ideal no
cero ${\eu m}$ de $\o_K$, sea $D_K^{\eu m}$ el grupo de ideales 
fraccionarios en $K$ primos relativos a ${\eu m}$
y sea $P^+_{K,{\eu m}}$\label{CClasePKm+}
el grupo de ideales fraccionarios $\big(\alpha/\beta\big)$ con $\alpha,
\beta\in \o_K$ tales que $(\alpha)$ y $(\beta)$ son primos relativos a
${\eu m}$, $\alpha\equiv\beta\bmod{\eu m}$, en el sentido de
que $\alpha/\beta-1\in\mK$ y $\alpha/\beta$ es totalmente
positivo, es decir, si $\varphi\colon K\to {\ma R}$ es un encaje real de 
$K$, $\varphi\big(\alpha/\beta\big)\in{\ma R}^+$, es decir, $\varphi\big(
\alpha/\beta\big)>0$.

Se tiene que $[D_K^{\eu m}:P_{K,{\eu m}}^+]<\infty$ y todo grupo intermedio
$P_{K,{\eu m}}^+\subseteq H\subseteq D_K^{\eu m}$ se llama un {\em grupo
de ideales con m\'odulus ${\eu m}$\index{grupo de ideales con m\'odulus
${\eu m}$}} y el cociente $D_K^{\eu m}/H$ se llama {\em grupo de ideales
generalizado\index{grupo de ideales generalizado}}.

Si ${\eu m}=(1)$ y $P_K$ es el grupo de los ideales principales, se tiene
$P_{K,(1)}^+\subseteq P_K\subseteq D_K^{(1)}$ y $I_K:=
D_K^{(1)}/P_K$ es el grupo de clases de ideales usuales.

\begin{teorema}[Weber]\label{CClaseT2.1.1} 
Para cualquier ideal no cero de ${\eu m}$ de $\o_K$ y para
un grupo de ideales $H$ con m\'odulus ${\eu m}$, supongamos que hay
una extensi\'on de Galois $L/K$ tal que $\Spl(L/K)\subseteq H$ con un
n\'umero finito de excepciones. Entonces $[D_K^{\eu m}:H]\leq [L:K]$.

Si $\Spl(L/K)=H$ con un n\'umero finito de excepciones, entonces
$[D_K^{\eu m}:H]=[L:K]$ y hay una infinidad de primos en cada clase de
$D_K^{\eu m}/H$. $\fin$
\end{teorema}

\begin{definicion}[Weber]\label{CClaseD2.1.2} Para un ideal no cero ${\eu m}$ de
$\o_K$ y para un grupo de ideales $H$ con m\'odulus ${\eu m}$, el
{\em campo de clase\index{campo de clase}} sobre $K$ para $H$ es una
extensi\'on de Galois $L/K$ tal que para los primos $\pK\nmid {\eu m}$
en $K$, $\pK$ se descompone totalmente en $L\iff\pK\in H$ (si tal $L$
existe, entonces $L$ es \'unico).
\end{definicion}

David Hilbert propuso una ley cuadr\'atica de reciprocidad:
$\prod_v(a,b)_v=1$ para cualesquiera $a,b\in\*K$ y $v$ recorre los lugares
de $K$. La prueba de Hilbert de esta f\'ormula no funciona para extensiones
no ramificadas. La prueba de Hilbert del Teorema de Kronecker--Weber
funcion\'o en parte debido a que ${\ma Q}$ no tiene extensiones propias no
ramificadas. Debido a lo anterior, Hilbert se interes\'o en las extensiones
abelianas no ramificadas como un obst\'aculo en las demostraciones.

\begin{conjetura}[Hilbert 1898]\label{CClaseC2.1.3}
Para cualquier campo num\'erico $K$ hay una
\'unica extensi\'on $K_H/K$ tal que:
\l
\item $K_H/K$ es Galois y $\Gal(K_H/K)\cong I_K$.

\item $K_H/K$ es no ramificada y toda extensi\'on abeliana con
esta propiedad est\'a contenida en $K_H$.

\item Para cualquier primo $\pK$ de $K$, el grado $f_{\pK}$ de los campos
residuales es el orden de $\pK$ en $I_K$.

\item Todo ideal $K$ se hace principal en $K_H$.
\end{list}
\end{conjetura}

As{\'\i} $K_H$ es un campo de clase en el sentido de Weber: el campo de
clase para el grupo de los ideales principales fraccionarios de $K$.

Takagi empez\'o, aproximadamente en 1914, con una nueva definici\'on
de campos de clase, usando normas de ideales ($\N_{L/K}\pL=
\pK^{f(\pL|\pK)}$, $\pK=\pL\cap K$, $f(\pL|\pK)$ el grado de inercia)
en lugar de leyes de
descomposici\'on y tambi\'en incorpor\'o los primos infinitos dentro
de la definici\'on de los m\'odulus .

La liga entre los puntos de vista de Weber y  de Takagi es que cuando $L/K$
es Galois y $\pK$ es no ramificada en $L$, $\pK$ se descompone
totalmente en $L\iff \pK$ es la norma de alg\'un ideal de $L$.

As{\'\i}, de ahora en adelante, un {\em m\'odulus\index{m\'odulus}}
es ${\eu m}={\eu m}_f{\eu m}_{\infty}$,
${\eu m}_f$ la parte finita de ${\eu m}$ (es el de antes) y ${\eu m}_{\infty}$
es un producto formal de encajes reales de $K$. Un ideal fraccionario se
dice primo relativo a ${\eu m}$ si lo es a ${\eu m}_f$. Sea $D_K^{\eu m}$
el grupo de los ideales fraccionarios primos relativos a ${\eu m}$ y sea 
$P_{K,{\eu m}}^+$ el grupo de los ideales fraccionarios principales $\big(\alpha/
\beta\big)$ con $\alpha,\beta\in\o_K\setminus\{0\}$ tales que:
\l
\item $(\alpha)$ y $(\beta)$ son primos relativos a ${\eu m}$,
\item $\alpha\equiv \beta\bmod {\eu m}_f$, en el sentido $\alpha/\beta-1\in\mK$,
\item $v\big(\alpha/\beta\big)>0$ para toda $v|{\eu m}_{\infty}$.
\end{list}

Un subgrupo intermedio $P_{K,{\eu m}}^+\subseteq H\subseteq 
D_K^{\eu m}$ se llama
un {\em grupo de ideales con m\'odulus ${\eu m}$\index{grupo de 
ideales con m\'odulus ${\eu m}$}}. Para una extensi\'on finita $L/K$, sea
$\N_{\eu m}(L/K)=\{{\eu a}\text{\ en\ }
 K\mid {\eu a}=\N_{L/K}({\eu c}) \text{\ para ${\eu c}$
en $L$ y ${\eu a}$ y ${\eu c}$ primos relativos a ${\eu m}$}\}$. 
Sea $H^{\eu m}(L/K)=
P_{K,{\eu m}}^+\N_{\eu m}(L/K)\label{CClaseH^m}$ el cual se llama 
{\em subgrupo de normas\index{subgrupo
de normas}}.

Resulta ser que todo subgrupo de $D_K^{\eu m}/
P_{K,{\eu m}}^+$ es el grupo de normas
de alguna extensi\'on abeliana finita de $K$. 
Se tiene que $P_{K,{\eu m}}^+(\text{Weber})
=P_{K,{\eu m}\infty}^+$ donde $\infty$ denota
el producto de todos los lugares reales de $K$.

Se tiene que los primos que no dividen a ${\eu m}$ y se descomponen 
totalmente en $L$, est\'an en $\N_{\eu m}(L/K)\subseteq H^{\eu m}(L/K)$,
es decir, $\Spl(L/K)\subseteq H^{\eu m}(L/K)$ con excepci\'on de los primos
que dividen a ${\eu m}$.

As{\'\i}, como antes, $[D_K^{\eu m}:H^{\eu m}(L/K)]\leq [L:K]$ (primera
desigualdad\index{primera desigualdad}).

\begin{definicion}[Takagi]\label{CClaseD2.1.4} Una extensi\'on de Galois de
campos num\'ericos $L/K$ se llama {\em campo de clase\index{campo de clase}}
si $[D_K^{\eu m}:H^{\eu m}(L/K)]=[L:K]$ para alg\'un m\'odulus ${\eu m}$
(un tal m\'odulus 
${\eu m}$ se llama un {\em m\'odulus admisible\index{m\'odulus admisible}}
o {\em m\'odulus de definici\'on\index{m\'odulus de definici\'on}} para $L/K$).
\end{definicion}

\begin{teorema}[Takagi\index{teorema de Takagi}
1920]\label{CClaseT2.1.5} Sea $K$ un campo num\'erico.
\l
\item {\underline{\rm{Existencia:}}} \quad \begin{minipage}{8cm}
Para cada grupo de ideales $H$ hay un campo de clase sobre $K$.
\end{minipage}

\item[]

\item[]

\item {\underline{\rm{Isomorfismo:}}} \quad \begin{minipage}{7.7cm}
Si $H$ es un grupo de ideales con m\'odulus ${\eu m}$ y tiene campo
de clase $L/K$, entonces $\Gal(L/K)\cong D_K^{\eu m}/H$. En particular
$L/K$ es abeliana.
\end{minipage}

\item[]

\item[]

\item {\underline{\rm{Completitud:}}} \quad \begin{minipage}{7.5cm}
Toda extensi\'on abeliana finita de $K$ es un campo de clase.
En particular, campo de clase de $K$ y extensi\'on abeliana
finita de $K$, es lo mismo.
\end{minipage}

\item[]

\item[]

\item {\underline{\rm{Comparaci\'on:}}} \quad \begin{minipage}{7.4cm}
Si $H_1$ y $H_2$ son grupos de ideales con m\'odulus com\'un ${\eu m}$
y ellos tiene campos de clase $L_1$ y $L_2$, entonces, $L_1\subseteq
L_2\iff H_2\subseteq H_1$.
\end{minipage}

\item[]

\item[]

\item {\underline{\rm{Conductor:}}} \quad \begin{minipage}{8cm}
Para toda extensi\'on abeliana $L/K$, los lugares de $K$ 
que aparecen en el soporte del conductor
${\eu f}_{L/K}$ son los primos ramificados en $L/K$.
\end{minipage}

\item[]

\item[]

\item {\underline{\rm{Descomposici\'on:}}} \quad \begin{minipage}{7.1cm}
Si $H$ es un grupo de ideales con m\'odulus ${\eu m}$ y campo de clase
$L/K$, entonces cualquier primo $\pK\nmid {\eu m}$ es no ramificado en
$L$ y el grado de inercia $f_\pK$ es igual al orden de $\pK$ en 
$D_K^{\eu m}/H$. $\fin$
\end{minipage}

\end{list}
\end{teorema}

En su demostraci\'on, Takagi prob\'o la {\em segunda desigualdad\index{segunda
desigualdad}} para una extensi\'on abeliana: $[D_K^{\eu m}:H^{\eu m}(L/K)]
\geq [L:K]$ para alguna ${\eu m}$. La primera desigualdad vale para toda
extensi\'on de Galois y la segunda desigualdad es v\'alida \'unicamente para
extensiones abelianas, es decir, si $L/K$ es una extensi\'on de Galois
no abeliana, entonces 
\[
[D_K^{\eu m}:H^{\eu m}(L/K)]<[L:K]\quad \text{para todo m\'odulus}\quad {\eu m},
\]
ver Teorema \ref{T17.6.192N}.

?`Como funciona la teor{\'\i}a de campos de clase? Si queremos una
correspondencia tipo Galois, se tiene que si tomamos todos los campos de
clase de golpe, tenemos el siguiente problema de comparaci\'on:
los moduli admisibles para dos campos de clase pueden no ser el mismo
por lo que tenemos que pasar a un m\'odulus com\'un 
para poder compararlos.

Para poder tener una biyecci\'on tipo Galois necesitamos identificar
todos los grupos de ideales que tienen el mismo campo de clase.
?`Como hacerlo? 

Si $H$ y $H^{\prime}$ son grupos de ideales para $K$
definidos moduli ${\eu m}$ y ${\eu m}^{\prime}$, es decir, $P_{K,{\eu m}}^+
\subseteq H\subseteq D_K^{\eu m}$ y $P_{K,{\eu m}^{\prime}}^+\subseteq
H^{\prime}\subseteq D_K^{{\eu m}^{\prime}}$, llamamos a $H$ y $H^{\prime}$
{\em equivalentes\index{subgrupos equivalentes}} 
si existe un m\'odulus ${\eu m}^{\prime\prime}$
divisible tanto por ${\eu m}$ como por ${\eu m}^{\prime}$ tal que los
homomorfismos naturales $D_K^{{\eu m}^{\prime\prime}}\mapsto D_K^{\eu m}
/H$ y $D_K^{{\eu m}^{\prime\prime}}\mapsto D_K^{{\eu m}^{\prime}}/H^{\prime}$
tienen el mismo n\'ucleo, es decir $H\cap D_K^{{\eu m}^{\prime\prime}}
=H^{\prime}\cap D_K^{{\eu m}^{\prime\prime}}$. 

Los grupos de ideales
equivalentes, tienen el mismo campo de clase y entonces la correspondencia
entre campos de clase sobre $K$ y los grupos de ideales en $K$,
hasta equivalencia, es biyectiva. Esto hace las cosas complicadas.

Cuando
pasamos al lenguaje de {\em id\`eles\index{id\`eles}}\label{CClaseid'eles}, todos los grupos de
ideales equivalentes se fusionan en un \'unico subgrupo de id\`eles, 
haciendo la teor{\'\i}a de campos de clase un poco m\'as simple,
o mejor dicho, menos complicada.

Ahora bien, el teorema de descomposici\'on de Takagi muestra que
para un primo $\pK\nmid \mK$, $\pK$ se descompone totalmente $\iff
\pK\in H$ as{\'\i} que las nociones de campos de clase de Weber y de Takagi
coinciden.

Por otro lado, las condiciones sobre los primos 
para que se descompongan totalmente
se da por condiciones de congruencia. Por ejemplo, los primos que
se descomponen totalmente 
en ${\ma Q}(i)/{\ma Q}$ son los primos $p\equiv 1\bmod 4$,
y $2$ es el \'unico primo ramificado. Los primos 
que se descomponen totalmente
en ${\ma Q}(\sqrt{6})/{\ma Q}$ son los primos $p$ tales que $p\equiv
1,5,19,23\bmod 24$. Finalmente, los primos que se descomponen totalmente
en ${\ma Q}(\zeta_n)/{\ma Q}$, donde ${\ma Q}(\zeta_n)$ es el $n$--\'esimo
campo ciclot\'omico, son los primos $p$ tales que $p\equiv 1\bmod n$. 

En general los primos de que no dividen a $\mK$ y que est\'an en 
$\Spl(L/K)$ son aquellos en el subgrupo $H_\mK/P^+_{K,\mK}$ 
de $D_K^\mK/P^+_{K,\mK}$
y que pertenecen a un subgrupo que puede ser pensado como condiciones
generalizadas de congruencias (por esto, los grupos $H_\mK$ se llaman
{\em grupos de congruencia\index{grupos de congruencia}}).

Ahora bien, puesto que los campos de clase y extensiones 
abelianas son lo mismo, la descomposici\'on total en una extensi\'on
abeliana est\'a descrita por congruencias. Resulta ser que
el rec{\'\i}proco tambi\'en se cumple.

\begin{teorema}\label{CClaseT2.1.6} Sea $L/K$ una extensi\'on finita
de campos num\'ericos y supongamos que existe un m\'odulus
$\mK$ y un conjunto finito $S$ de primos que contienen a todos
los que dividen a $\mK$, de tal forma que la condici\'on de que
un primo $\pK\notin S$
es o no totalmente descompuesto en $L$ est\'a determinado
por la clase de $\pK$ en $D_K^\mK/P^+_{K,\mK}$. Entonces $L/K$
es una extensi\'on abeliana. $\fin$
\end{teorema}

\begin{corolario}\label{CClaseC2.1.7} Para un campo num\'erico
$L/{\ma Q}$ y $m\in {\ma N}$ las siguientes condiciones son
equivalentes:
\l
\item Para cualquier primo positivo $p\nmid m$, la 
descomposici\'on de $p$ est\'a determinada por una condici\'on
de congruencia en $p\bmod m$.

\item $L\subseteq {\ma Q}(\zeta_m)$. $\fin$
\end{list}
\end{corolario}

Takagi prob\'o que hay un isomorfismo $D_K^\mK/H_\mK\cong\Gal(
L/K)$ para todos los moduli $\mK$ que son $K$--admisibles.
Sin embargo no dio ning\'un isomorfismo; el isomorfismo fue 
obtenido de manera indirecta. Hoy sabemos que sus argumentos
pertenecen a la cohomolog{\'\i}a de grupos. Artin describi\'o este
isomorfismo por medio de la {\em Ley de Reciprocidad\index{ley
de reciprocidad}}.

\begin{definicion}\label{CClaseD2.1.8} Para una extensi\'on abeliana
$L/K$ y un $K$--m\'odulus divisible por todos los primos que se
ramifican en $L$, el {\em mapeo de Artin\index{automorfismo
de Artin}\index{mapeo de Artin}\index{Artin!automorfismo de $\sim$}}
$\psi_{L/K, \mK}\colon D_K^\mK\to\Gal(L/K)\label{CClaseartinconm}$ est\'a dado por
$\psi_{L/K,\mK}(\pK)=(\pK,L/K)=\artinp{L/K}{\pK}$.
\end{definicion}

\begin{teorema}[Ley de Reciprocidad de Artin\index{ley
de reciprocidad de Artin}\index{Artin!ley de reciprocidad de $\sim$},
Artin 1927]\label{CClaseT2.1.9}

El mapeo de Artin $\psi_{L/K,\mK}$ es suprayectivo y su
n\'ucleo contiene a $\N_\mK(L/K)$. Cuando $\mK$ es admisible,
el n\'ucleo de $\psi_{L/K,\mK}$ es $P^+_{K,\mK}\N_\mK(L/K)=
H_\mK(L/K)$, esto es, $D_K^\mK/H_\mK(L/K)\cong\Gal(L/K)$
mediante el mapeo de Artin. $\fin$
\end{teorema}

La parte m\'as dif{\'\i}cil en la ley de reciprocidad es probar que
el n\'ucleo de $\psi_{L/K,\mK}$, con $\mK$ admisible, contiene
a $P^+_{K,\mK}$. Es decir, probar que si $(\alpha)\in P^+_{K,\mK}$, entonces
$\psi_{L/K,\mK}((\alpha))=1$.

\subsubsection{Campos de clase v{\'\i}a id\`eles (C. Chevalley)}\label{CClaseS2.1.1}

La teor{\'\i}a de campos de clase locales (o teor{\'\i}a local de campos
de clase) fue establecida por sus propios m\'eritos y no como
reducci\'on de la teor{\'\i}a global principalmente por H. Hasse.

\begin{teorema}[H. Hasse, F.K. Schmidt, 1930]\label{CClaseT2.1.10}
Para una extensi\'on abeliana de campos locales $E/F$ (de
caracter{\'\i}stica $0$), el {\em mapeo local de Artin\index{mapeo
local de Artin}\index{Artin!mapeo local de $\sim$}} o
{\em s{\'\i}mbolo residual de la norma \index{simbolo residual de la norma@s\'imbolo residual de la
norma}}, $\psi_{E/F}\colon \*F\twoheadrightarrow\Gal(E/F)$
es un epimorfismo con n\'ucleo $\N_{E/F}\*E$ por lo que
$\*F/\N_{E/F}\*E\cong \Gal(E/F)$. $\fin$
\end{teorema}

De esta forma, asociando a $E$ el grupo $\N_{E/F}(\*E)$ 
obtenemos una correspondencia biyectiva que voltea el orden entre
las extensiones abelianas finitas de $F$ y los subgrupos abiertos de
{\'\i}ndice finito en $\* F$.

La imagen de las unidades de $F$, 
$U_F:=\*{\o_F}$ bajo el mapeo local de Artin
es el grupo de inercia $I(E/F)$, es decir,
$\psi_{E/F}(U_F)=I(E/F)$, as{\'\i} que
\begin{gather*}
e(E|F)=[U_F\N_{E/F}(\*E):\N_{E/F}(\*E)]=[U_F:\N_{E/F}U_E]
\intertext{(ver Corolario \ref{C17.2.15'}). En consecuencia,}
f(E|F)=\frac{[E:F]}{e(E|F)}=[\*F:U_F\N_{E/F}(\*E)]
\end{gather*}
es el orden de $\pi$ en $\*F/U_F\N_{E/F}(\*E)$ para cualquier
elemento primo $\pi$ de $F$.

Cuando la extensi\'on $E/F$ es no abeliana, se tiene
\[
[\*F:\N_{E/F}(\*E)]<[E:F].
\]

Una vez que la teor{\'\i}a local de campos de clase fue 
establecida, el siguiente paso fue obtener los teoremas
de la teor{\'\i}a global de campos de clase por medio de 
aqu\'ellos de la teor{\'\i}a local. El concepto que permite hacer
esto son los {\em id\`eles\index{id\`eles}},
los cuales tambi\'en permiten teor{\'\i}a
de campos de clase globales para extensiones 
abelianas infinitas.

\begin{definicion}[Chevalley, 1936]\label{CClase2.1.11}
El {\em grupo de id\`eles\index{grupo de id\`eles}} $J_K$ de un
campo num\'erico $K$ es el conjunto de sucesiones
$\vec x=(x_v)_v$ indexadas por el conjunto de lugares
$v$ de $K$, tales que $x_v\in \*{K_v}$ para toda $v$
y adem\'as $x_v\in \*{\o_v}=U_{K_v}=U_v$ para casi todos los
lugares $v$ (es decir, para todos, salvo un n\'umero finito)
y donde $\o_v$ denota al anillo de enteros de $K_v$, el
campo completado de $K$ en $v$,
y $U_{K_v}=\*{\o_v}$ es el grupo de unidades de $\o_v$.
\end{definicion}

Un elemento de $J_K$ se llama {\em id\`ele\index{id\`eles}}.
Este fue el nombre sugerido por Hasse a Chevalley, el cual
los hab{\'\i}a llamado originalmente como {\em elemento
ideal}. $J_K$ es un grupo bajo la multiplicaci\'on entrada
por entrada. Se tiene el encaje diagonal $\*K\hookrightarrow
J_K$ y la imagen se llama el {\em grupo de los id\`eles
principales\index{id\`eles principales}}. Similarmente se
tiene el encaje $\enc{\ }v\colon \*{K_v}\hookrightarrow J_K$
tal que $x_v\in \*{K_v}$ se mapea a $x_v$ en la entrada $v$
y con $1$ en las dem\'as componentes.

Para $\vec x\in J_K$ se tiene el ideal fraccionario (recordemos
que estamos en campos num\'ericos), 
\[
\imath(\vec x)={\eu a}_{\vec x}
=\prod_{v\nmid
\infty}\pK_v^{v(x_v)}
\] 
lo cual permite pasar de id\`eles a ideales.

Usando este paso de id\`eles a ideales se tiene que cualquier
grupo de clase generalizado de $K$ se puede realizar como
un grupo cociente de $J_K$ como sigue: sea $\mK$ un $K$--m\'odulus.
Sean $\vec x\in J_K$, $\alpha_0\in\*K$ tal que para todo $\pK$
en el soporte de $\mK$, se cumple $v_\pK(x_\pK/\alpha_0-1)
\geq v_\pK(\mK)$ cuando $\pK| \mK_f$ (la parte finita de $\mK$) y
$x_\pK/v_\pK(\alpha_0)>0$ cuando $\pK|\mK_{\infty}$ (la parte
infinita de $\mK$).
Lo anterior se puede lograr gracias al 
Teorema de Aproximaci\'on de Artin, Teorema \ref{CClaseT1.1.1}.

El id\`ele $\frac{\vec x}{\alpha_0}=\big(\ldots, \frac{x_{\eu q}}
{\alpha_0},\ldots\big)_{\eu q}$ tiene ideal correspondiente
$\imath\big(\frac{\vec x}{\alpha_0}\big)\in D_K^\mK$ y adem\'as
$\overline{\imath\big(\frac{\vec x}{\alpha_0}\big)}\in 
D_K^\mK/P^+_{K,\mK}$
est\'a bien definido. El n\'ucleo del mapeo $J_K\to D_K^\mK/P^+_{D,\mK}$
contiene a los ideales principales, por lo que $D_K^\mK/P^+_{K,\mK
}\cong J_K/\*KS_\mK$ para alg\'un $S_\mK\subseteq J_K$.

Resulta que dos grupos de ideales $H$ y $H^{\prime}$ son equivalentes
(es decir $H\cap D_K^{\mK^{\prime\prime}}=
H^{\prime}\cap D_K^{\mK^{\prime\prime}}$
para alg\'un m\'ultiplo $\mK^{\prime\prime}$ de $\mK$ y $\mK^{\prime}$)
corresponden al mismo grupo de id\`eles.

El mapeo de Artin es el siguiente:
\[
\psi_{L/K}\colon J_K\longrightarrow D_K^\mK/P^+_{K,\mK}\xrightarrow{
\psi_{L/K,\mK}} \Gal(L/K),
\]
y se tiene que $\psi_{L/K}$ es suprayectiva e independiente de la
elecci\'on del m\'odulus admisible $\mK$. Adem\'as $\psi_{L/K}
(\*K)=1$.

Si $L/K$ es una extensi\'on de Galois,
se define la norma de $J_L$ a $J_K$ como
$\N_{L/K}\colon J_K\to J_K$ definida por
$\N_{L/K}(\vec y)=\vec x$ donde 
$x_\pK=\prod_{\pL | \pK}\N_{L_\pL/K_\pK}
(y_\pL)$.

Entonces $\ker \psi_{L/K}=\*K\N_{L/K}(J_L)$ en el caso de que $L/K$
sea una extensi\'on abeliana.

Para que la correspondencia entre extensiones abelianas de $K$ y
subgrupos de id\`eles sea biyectiva, se necesita hacer de $J_K$ un
espacio topol\'ogico. La topolog{\'\i}a dada es la topolog{\'\i}a del 
producto restringido: una base de vecindades abiertas de $\vec 1\in
J_K$ est\'a formado por los conjuntos 
$\prod_\pK V_\pK$ donde $V_\pK$ es una 
vecindad abierta de $1\in \*{K_\pK}$ para toda $\pK$ y $V_\pK=
U_\pK=\*{\o_\pK}$ para casi todo $\pK$. Entonces $J_K$ es un grupo
topol\'ogico localmente compacto (la topolog{\'\i}a producto 
{\underline{no}} es localmente compacta).

\begin{teorema}\label{CClase2.1.12} Para una extensi\'on abeliana finita
de campos num\'ericos $L/K$, el mapeo de Artin $\psi_{L/K}$ es
un epimorfismo de $J_K$ sobre $\Gal(L/K)$ con n\'ucleo $\*K 
\N_{L/K}(J_L)$ por lo que $J_K/\*K \N_{L/K}(J_L)\cong\Gal(L/K)$.
La correspondencia que asocia a cada extensi\'on abeliana finita $L$ de $K$
con el subgrupo $\*K\N_{L/K}(J_L)$ es biyectiva entre las extensiones
abelianas finitas de $K$ y los subgrupos abiertos de {\'\i}ndice finito de
$J_K$ y que contienen a $\*K$. La correspondencia voltea contenciones.
$\fin$
\end{teorema}

Para un lugar $\pK$ de $K$, la composici\'on
\[
\*{L_\pL}\xrightarrow{\norma} \*{K_\pK}\xrightarrow{\enc{\ }{\pK}}J_K
\xrightarrow{\psi_{L/K}}\Gal(L/K)
\]
tiene como imagen el grupo de descomposici\'on $D(\pL|\pK)$ y
la imagen de $\*{\o_\pL}=U_\pL$ es el grupo de inercia $I(\pL|\pK)$.
Adem\'as $\Gal(\abe K/K)$ es el m\'aximo grupo cociente de $J_K/\*K$
totalmente disconexo.

M\'as a\'un, para $\vec x\in J_K$, $\psi_{L/K}(\vec x)=\prod_\pK
(x_\pK, L_\pL/K_\pK)$ donde $(x_\pK, L_\pL/K_\pK)$ es mapeo local
de Artin, $(x_\pK, L_\pL/K_\pK)\in \Gal(L_\pL/K_\pK)\cong D(\pL|\pK)$ y
para cada $\pK$ seleccionamos un \'unico $\pL|\pK$, cualquiera pero
\'unicamente uno y en general se denota $L_\pK$ en lugar de $L_\pL$.

\subsubsection{Campos de funciones}\label{CClaseS2.1.2}

Hasse prob\'o que los teoremas de teor{\'\i}a local de campos de clase
son los mismos en caracter{\'\i}stica $0$ que en
caracter{\'\i}stica $p>0$, excepto que necesitamos
ser expl{\'\i}citos acerca de usar {\underline{subgrupos abiertos}}
de {\'\i}ndice finito. 

En 1935 Witt prob\'o el {\em Teorema de Existencia\index{teorema
de existencia}} para extensiones
abelianas con grado divisible por $p$ lo cual complet\'o el trabajo de 
Schmidt para extensiones abelianas de grado no divisible por $p$.

El punto de vista de Chevalley por medio de id\`eles funciona en ambos
casos, campos num\'ericos y campos de funciones, sin embargo
tenemos una diferencia entre los dos casos para extensiones abelianas
infinitas. Para un campo de funciones $K$, como en el caso num\'erico,
el mapeo de Artin $J_K/\*K\to \Gal(\abe K/K)$ tiene imagen densa,
pero ahora el mapeo es inyectivo en lugar de suprayectivo (la
demostraci\'on de la suprayectividad en el caso de campos de n\'umeros
falla en el caso de campos de funciones pues no existen lugares
arquimedianos en estos \'ultimos).

La imagen del mapeo de Artin en el caso de campos de funciones
est\'a caracterizado como el conjunto de elementos de $\Gal(\abe K/K)$ los
cuales, en la cerradura algebraica del campo de constantes de $K$,
son potencias enteras del automorfismo de Frobenius.

M\'as precisamente, se tiene $\bar {\ma F}_q=\abe\F$ y el siguiente diagrama
\[
\xymatrix{
&\abe K\ar@{-}[ld]_{\mc G}\ar@{-}[d]^H\\
K\ar@{-}[r]^G\ar@{-}[d]&K\bar {\ma F}_q\ar@{-}[d]\\
\F\ar@{-}^G[r]&\bar {\ma F}_q
}
\]

Sean ${\mc G}=\Gal(\abe K/K)$, $H=\Gal(\abe K/K\bar {\ma F}_q)$ y
$G=\Gal(K\bar {\ma F}_q/K)\cong \Gal(\bar {\ma F}_q/\F)\cong
\hat{\ma Z}$ la completaci\'on de ${\ma Z}$. Entonces $G\cong
{\mc G}/H$. 

El mapeo $\psi\colon J_K/\*K\lra{\mc G}=\Gal(\abe K/K)$
satisface que $\im\psi=\{\sigma\in{\mc G}\mid \sigma|_{K
\bar {\ma F}_q}\in{\ma Z}\}$ donde $\sigma|_{K\bar {\ma F}_q}\in{\ma Z}$
significa que si $\tau$ es el Frobenius de $K$, esto es, $\tau x=x^q$,
entonces $\sigma|_{K\bar {\ma F}_q}=\tau^m$ para alguna $m\in{\ma Z}$.

Para finalizar, la teor{\'\i}a de campos de clase no provee campos de clase
de manera expl{\'\i}cita. En el caso de campos de n\'umeros, \'unicamente
tenemos los campos de clase expl{\'\i}citos para ${\ma Q}$ (Teorema
de Kronecker--Weber, campos ciclot\'omicos) y para los
campos cuadr\'aticos imaginarios. De hecho, en 1880 Kronecker
en una carta a Dedekind describi\'o su 
``{\em sue\~no de juventud\index{sue\~no de juventud de Kronecker}}''
(``{\em Jugendtraum}'' en alem\'an), como poder describir las extensiones
abelianas de un campo num\'erico por medio de extensiones generadas
por ra{\'\i}ces de algunas funciones transcendentes. El sue\~no de Kronecker
no se ha materializado todav\'ia.

Para campos de funciones, D. Hayes en 1974, basado en el trabajo de su
asesor, L. Carlitz, construy\'o una teor{\'\i}a de campos de clase 
expl{\'\i}cita sobre el campo de funciones racionales $\F(T)$.
Ver Cap\'itulo \ref{Ch6}.

Drinfeld en
el mismo a\~no (1974), usando ``{\em m\'odulos el{\'\i}pticos\index{m\'odulos
el{\'\i}pticos}}'', ahora conocidos como ``{\em m\'odulos de 
Drinfeld\index{m\'odulos de Drinfeld}}'', hizo expl{\'\i}cita la obtenci\'on
de los campos de clase sobre cualquier campo de funciones congruente.

Lo que hace diferente lo expl{\'\i}cito entre los campos num\'ericos y los
campos de funciones es que, en caracter{\'\i}stica $p>0$, hay muchas
funciones aditivas (ver el Cap\'itulo \ref{DrinfeldCh15} sobre los m\'odulos de
Drinfeld).

\section{Campos locales}\label{CClaseS1.2}

En esta secci\'on presentamos los resultados fundamentales
de los campos locales. Hacemos notar que varios de los 
resultados para campos locales siguen siendo v\'alidos 
para campos completos con respecto a una valuaci\'on.
Ver tambi\'en la Secci\'on \ref{S7.6}.

\begin{definicion}\label{CClaseD1.2.1}
Dados un campo $K$ y $v=v_K\colon \*K\to {\ma R}$ una 
valuaci\'on\index{valuaci\'on}, es decir, $v$ satisface
$v(xy)=v(x)+v(y)$ y $v(x+y)\geq \min\{v(x),v(y)\}$ 
para cualesquiera $x,y\in \*K$, donde escribimos $v(0)=
\infty$, se define el {\em valor absoluto
\index{valor absoluto}}\label{CClasevalorabsoluto}
por $|x|_v=|x|=c^{v(x)}$ con $0<c<1$ fijo
arbitrario. La topolog{\'\i}a dada a $K$ por $|\ |_v$ es
independiente del $c$ seleccionado.
\end{definicion}

\begin{definicion}\label{CClaseD1.2.3} Una valuaci\'on se llama
{\em discreta\index{valuaci\'on discreta}\index{discreta!valuaci\'on}}
si $v(\*K)\cong {\ma Z}$. Una valuaci\'on discreta $v$ se llama
{\em normalizada\index{valuaci\'on normalizada}} si $v(\*K)={\ma Z}$.
\end{definicion}

En este cap\'itulo \'unicamente se estudiar\'an valuaciones
discretas.

Consideraremos campos completos con respecto a $|\ |_v$. Sea
$\o_K=\{x\in K\mid v(x)\geq 0\}\index{anillo 
de valuaci\'on}\label{CClaseanillovaluacion}$
el {\em anillo de valuaci\'on\index{anillo
de valuaci\'on}}. Entonces $\o_K$ es un anillo local (de hecho,
es un anillo de valuaci\'on) con ideal m\'aximo\index{ideal m\'aximo}
\label{CClaseidealmaximo} $\pK_K=\pK=\{x\in K\mid v(x)>0\}$.
Se tiene que $\o_K=\bar{B}_{|\ |_v}(0,1)=\{x\in K\mid |x|_v\leq 1\}$
es la bola cerrada,
$\pK=B_{|\ |_v}(0,1)=\{x\in K\mid |x|_v<1\}$ es la bola abierta
y $\tilde{K}=K(\pK)=
\o_K/\pK$ es el {\em campo
residual\index{campo residual}\label{CClasecamporesidual}}.

\begin{definicion}\label{CClaseD1.2.3+1}
Las unidades de $K$ se definen por $U_K:=\*{{\mc O}_K}=
{\mc O}_K\setminus \pK=\{a\in\*K\mid v(x)=0\}=
\{x\in K\mid |x|_v=1\}$, es decir, $U_K$ es la frontera de la bola
unitaria con centro en $0$.
\end{definicion}

Si $v_K$ es discreta, entonces $v_K(\*K)=\beta{\ma Z}$ para alg\'un
$\beta\in{\ma R}^+$. Sea $\pi\in K$ con $v_K(\pi)=\beta$. Sea
$\xi\in\pK$, digamos $v_K(\xi)=\beta n$ para alg\'un $n\geq 1$.
De esta forma tenemos que $v_K(\xi \pi^{-n})=0$, esto es, $\xi
\pi^{-n}=u\in U_K$ es una unidad y $\xi=u\pi^{n}\in \pK$. Se sigue
que $\pK=\langle\pi\rangle$ es principal.
Rec\'iprocamente, si $\pK=\langle\pi\rangle$ es principal, sea
$v_K(\pi)=\alpha\in{\ma R}^+$. Para $\xi\in\*K$, $\xi=\frac ab$, 
$a, b\in\o_K$, con $b\neq 0$, se tiene $v_K(\xi)=v_K(a)-v_K(b)$.
Para $a\in\o_K$, si $a\in U_K$, $v_K(a)=0=0\cdot \alpha$; si
$a\in\pK_K$, $a=a_1\pi$ con $a_1\in \o_K$. Repitiendo el proceso
para $a_1$ concluimos que existe $n\in{\ma N}$ y $a_n\in U_K$ con
$a=a_n \pi^n$. Se sigue que $v_K(a)=n\alpha$. Por tanto
$v_K(\*K)=\alpha{\ma Z}$ y $v_K$ es discreta. En resumen, hemos
probado que 

\begin{proposicion}\label{CClaseP1.2.3+2}
Una valuaci\'on $v_K$ es discreta si y s\'olo si $\pK_K$ es un ideal
principal. $\fin$
\end{proposicion}

Si $\ca \tilde{K}=0$ entonces $\ca K=0$. Si $\ca \tilde{K}=p>0$, entonces
puede ser $\ca K=p$ o $\ca K=0$. 

\begin{definicion}\label{CClaseD1.2.3+3}
Un campo $K$ completo con respecto
a una valuaci\'on se llama {\em local\index{campo local}} si $\tilde{K}$ es
un campo finito.
\end{definicion}

Por razones t\'ecnicas para la teor\'ia de campos globales, a veces
se consideran tanto a ${\ma R}$ como a ${\ma C}$ como campos
locales.

 En el caso de un campo local, se seleccionar\'a al 
n\'umero $0<c<1$ que define el valor absoluto como $c=q^{-1}$ donde el
campo residual de $K$ es $\F$. 

\begin{definicion}\label{CClaseD1.2.3+4}
Sea $v$ normalizada. Un {\em elemento
primo\index{elemento primo}} o {\em elemento
uniformizante\index{elemento uniformizante}} $\pi_K=\pi$ de $K$
es cualquier elemento de valuaci\'on $1$: $v(\pi)=1$.
\end{definicion}

Sea $L/K$ una extensi\'on de campos, $w$ una valuaci\'on de $L$
y sea $v:=w|_K$. Entonces $\o_K=K\cap \o_L$, $\pK_K=K\cap
\pK_L$, $\tilde{K}=\o_K/\pK_K=\o_K/(\o_K\cap \pK_L)\cong
(\pK_L+\o_K)/\pK_L\subseteq \o_L/\pK_L=\tilde{L}$. Se
tiene que $v(\*K)\subseteq w(\*L)\subseteq {\ma R}^+$.
Adem\'as $U_K=\{x\in\o_K\mid v(x)=0\}$ y $\*K/U_K\cong
v(\*K)$.

\begin{definicion}\label{ramCL} El {\em \'indice de ramificaci\'on} $e$,
y el {\em grado de inercia} $f$, se definen por:
\begin{align*}
e&=e(\pK_L|\pK_K)=e(w|v)=e(L|K)=[w(\*L):v(\*K)],\\
f&=f(\pK_L|\pK_K)=f(w|v)=f(L|K)=[\tilde{L}:\tilde{K}].
\end{align*}
\end{definicion}

Dado cualquier campo $K$ con valuaci\'on $v$, sea $\bar{K}$ la
completaci\'on de $K$ con respecto a la topolog\'ia dada por $v$.
Entonces dado $x\in\bar{K}$, existe $\{x_n\}_{n=0}^{\infty}\subseteq
K$ con $\lim_{n\to\infty}x_n=x$. Entonces $v$ se extiende a
$\bar{K}$ como $\bar{v}(x)=\lim_{n\to\infty}v(x_n)=\lim_{n\to\infty}
\bar{v}(x_n)$. Es claro que $\bar{v}$ no depende la sucesi\'on
$\{x_n\}_{n=0}^{\infty}$ y que $\bar{v}$ es \'unica. Sea $w:=\bar{v}$.
Como $w$ no depende de la sucesi\'on, si $x=0$ podemos tomar
$x_n=0$ para toda $n$ y $w(0)=\infty=v(x_n)$.

Sea $x\neq 0$, $x_n\neq 0$ y $v(x_n)\xrightarrow[n\to\infty]{}w(x)
\neq 0$. Por otro lado $x-x_n\xrightarrow[n\to\infty] {} 0$, por lo
que $v(x-x_n)\xrightarrow[n\to\infty]{} \infty$. Por tanto
\[
w(x)=w(x-x_n+x_n)=\min\{w(x-x_n),w(x_n)\}=w(x_n)=v(x_n),
\]
para toda $n$ suficientemente grande. Por tanto $w(\*{\bar{K}})=
v(\*K)$ de donde se sigue que $e=e(\bar{K}|K)=1$. Tambi\'en, tenemos
el isomorfismo $\tilde{\bar{K}}=\tilde{K}$ de donde se obtiene
que $f=f(\bar{K}|K)=1$.

Un resultado para campos completos con respecto a una
valuaci\'on es el Lema de Hensel. La versi\'on que es \'util para
nuestros fines es la siguiente.

\begin{teorema}[Lema de Hensel\index{teorema!Lema de 
Hensel}\index{Lema de Hensel}\index{Hensel!Lema de $\sim$}]\label{Lema de Hensel}
Sea $K$ un campo completo con respecto a una valuaci\'on $v$
y anillo de enteros $A=\{x\mid v(x)\geq 0\}$ un anillo de valuaci\'on
discreta. Sean $f(x)\in A[x]$ un polinomio y $\bar{f}:=f\bmod
\pK$ donde $\pK$ es el ideal m\'aximo de $A$. Entonces 
cualquier ra\'iz simple $\lambda$ de $\bar{f}$ en $A/\pK$ se
levanta de manera \'unica a una ra\'iz de $f$ en $A$, es decir,
existe un \'unico $\alpha\in A$ tal que $f(\alpha)=0$ y donde
$\bar{\alpha}=\lambda$.
\end{teorema}

\begin{proof}
Sea $\beta\in A$ una tal ra\'iz de $f$, es decir, $\bar{\beta}=\lambda$.
Entonces $f(x)=(x-\beta)g(x)$ con $\bar{g}(\lambda)\neq 0$. Si $\gamma$
es otra ra\'iz de $f$ con $\bar{\gamma}=\lambda$, entonces $f(\gamma)=
(\gamma-\beta)g(\gamma)=0$. Ahora bien $g(\gamma)\bmod\pK=\bar{g}
(\lambda)\neq 0$, por lo que $g(\gamma)\in A\setminus \pK$, es decir,
$g(\gamma)$ es invertible en $A$. Por tanto $\gamma=\beta$, de donde
se sigue la unicidad del levantamiento de la ra\'iz.

Para la existencia, sea $\beta_1\in A$ tal que $\bar \beta_1=\lambda$.
Entonces $f(\beta_1)\equiv 0\bmod \pK$. Supongamos que hemos
hallado $\beta_n\in A$ tal que $\bar\beta_n=\lambda$ y $f(\beta_n)
\equiv 0\bmod \pK^n$. Hallaremos $\beta_{n+1}\in A$ tal que 
$\beta_{n+1}=\beta_n\bmod \pK^n$ y $f(\beta_{n+1})\equiv 0
\bmod \pK^{n+1}$.

Sea $h\in\pK^n$ por determinarse y sea $\beta_{n+1}=\beta_n+h$.
Desarrollando en serie de Taylor se tiene $f(\beta_{n+1})=
f(\beta_n)+hf'(\beta_n)
+h^2 y$ para alg\'un $y\in A$. Se tiene que $h^2y\in \pK^{2n}$. 
Necesitamos que $f(\beta_n)+hf'(\beta_n)\equiv 0\bmod \pK^{n+1}$.
Puesto que $\lambda$ es ra\'iz simple de $\bar{f}$, se tiene que
$\bar{f'}(\lambda)=\bar{f'}(\bar\beta_n)\neq 0$. Se sigue que
$f'(\beta_n)\in A\setminus \pK$, esto es, $f'(\beta_n)$ es
invertible y la ecuaci\'on 
\[
h=\frac{\xi-f(\beta_n)}{f'(\beta_n)}\in \pK^n,
\]
con $\xi\in \pK^{n+1}$ arbitrario ($\xi-f(\beta_n)\in\pK^n$) tiene soluci\'on y 
$f(\beta_n)+hf'(\beta_n)\equiv\xi\equiv 0\bmod \pK^{n+1}$.

Sea $\beta:=\lim_{n\to\infty}\beta_n$ el cual existe, $\beta\in A$ y
$f(\beta)=0$. $\fin$
\end{proof}

En el caso de un campo de funciones congruente, $\K /\F$, si $\pK$ es un
lugar de $\K $ y $\K_{\pK}$ es la completaci\'on con respecto a la valuaci\'on
$v_\pK$, $\K_\pK\cong {\ma F}_{q^d}((\pi))$ donde $[\tilde{\K}:\F]=d$,
esto es, $\tilde{\K}={\ma F}_{q^d}$ y $\pi\in \K$ es un elemento tal que $
v_\pK(\pi)=1$.

Se tiene que los campos residuales de $\K_{\pK}$ y de $\K$ en $\pK$ son 
isomorfos, es decir, $\o_{\tilde{\pK}}/\tilde{\pK}\cong \o_{\pK}/\pK$. El
isomorfismo se sigue del mapeo $\xi\colon \o_\pK\hookrightarrow \o_{
\tilde{\pK}}\twoheadrightarrow \o_{\tilde{\pK}}/\tilde{\pK}$ donde el primer
mapeo es la inyecci\'on natural y el segundo es la proyecci\'on natural;
$\xi$ es suprayectiva y $\ker \xi =\pK$ (ver \cite[Proposition 2.3.10]{Vil2006}).
Ponemos $\pK=\pK_K$.

\begin{proposicion}\label{CClaseP1.2.3} Si $K$ es un campo local, 
entonces como grupo multiplicativo se tiene 
\begin{gather*}
\* K\label{CClasecampolocal}=\langle\pi\rangle
\times U_K=\langle\pi\rangle\times \* {\F}\times U_K^{(1)}
\intertext{donde $\F$ es el campo residual $\tilde{K}$ de $K$,}
U_K^{(n)}= 1+\pK_K^n\index{n@$n$--unidades locales}
\label{CClasenunidadeslocales}=\{a\in K\mid a\equiv 1\bmod\pi^n\},
\quad n\geq 1,\\
U_K=U_K^{(0)}=\{a\in K\mid v_K(a)=0\}=\{a\in K\mid |a|_v=1\}\index{unidades
locales}\label{CClaseunidadeslocales},
\end{gather*}
 esto es,
$U_K=\* {\o_K}$ es la circunferencia unitaria y $q=|\tilde{K}|$. 
Notemos que $\unidades 0\neq 1+\pK^0=\o_K$.

Al grupo $U_K^{(n)}$ se le llama
el grupo de las {\em $n$--unidades principales\index{n@$n$--unidades
principales}}, $n\geq 1$  y se define
$|x|_v=q^{-v(x)}$, es decir, seleccionamos $c=q^{-1}$
en la Definici\'on {\rm{\ref{CClaseD1.2.1}}}. Adem\'as
$\* \F\cong W_{q-1}$, donde $W_t$ denota al grupo de las 
$t$--ra{\'\i}ces de uno, $\pK=\pi\label{CClaseelementoprimo} \o_K$, 
$\langle\pi\rangle=\{\pi^m\mid
m\in {\ma Z}\}\cong {\ma Z}$.

Finalmente tenemos 
\begin{gather*}
U_K/U_K^{(1)}\cong \* {\tilde{K}}=\* \F\ (\text{multiplicativo$)$
 y \ }\\
 U_K^{(n)}/U_K^{(n+1)}\cong \pK_K^n/\pK_K^{n+1}
\cong \tilde{K}\cong \F,  n\geq 1\ (\text{aditivo$)$}.
 \end{gather*}
\end{proposicion}
\begin{proof}
Sea $a\in\*K$ con $v_K(a)=n\in{\ma Z}$. Entonces $v_K(a\pi^{-n})=0$,
es decir, $a\pi^{-n}=u\in U_K$ y $a=\pi^n u$. Por tanto la
descomposici\'on $a=\pi^{v_K(a)}u$ es \'unica.

Veamos que la representaci\'on es \'unica. Digamos que $a=\pi^mu_1=
\pi^nu$, por lo que $v_K(a)=m=n$ y $u_1=u$. Se sigue que la funci\'on
$\*K\to\langle\pi\rangle \times U_K$, $a\mapsto (\pi^{v_K(a)},a\pi^{-v_K(a)})$
es un isomorfismo de grupos.

Ahora $U_K\stackrel{\varphi}{\lra}\*{\big(\o_K/\pK_K\big)}$, $u\mapsto u
\bmod \pK_K$ es un homomorfismo de grupos y $\ker\varphi=\{u\in U_K\mid
u\equiv 1\bmod \pK_K\}=\unidades 1$. Se tiene que $\varphi$ es un
epimorfismo pues si $\xi\bmod \pK\in\*{\big(\o_K/\pK_K\big)}$ se tiene
$v_K(\xi)= 0$, por lo que $\xi\in U_K$. 

La ecuaci\'on $X^{q-1}-1$ se descompone totalmente
en $\o_K$ por el Lema de Hensel, Teorema \ref{Lema de Hensel}. Las
ra\'ices de $X^{q-1}-1$ son los elementos de $\*\F$ y por tanto la
sucesi\'on
\[
1\lra \unidades 1\lra U_K\lra \*{\big(\o_K/\pK_K\big)}\cong \*\F\lra 1,
\]
se escinde y $U_K\cong \*\F\times \unidades 1$ como grupos.

Para $n\geq 1$ se tiene que $\unidades n\stackrel{\theta}{\lra}\o_K/\pK_K$
dada por $1+x\pi^n\lra x\bmod p_K$ es un epimorfismo y $\ker \theta
=\{1+x\pi^n\mid x\in\pK_K\}=\unidades {n+1}$. Por tanto 
$\unidades n/\unidades {n+1}\cong \o_K/\pK_K\cong \F$.

Finalmente, se verifica directamente que el mapeo 
\[
\unidades n/
\unidades {n+1}\stackrel{g}{\lra}\pK_K^n/\pK_K^{n+1}, \qquad (1+a\pi^n)
\bmod \unidades {n+1}\mapsto a\pi^n\bmod \pK_K^{n+1}
\]
es un isomorfismo de grupos.
$\fin$
\end{proof}

Notemos que, puesto que $U_K\supseteq \unidades 1\supseteq \cdots
\supseteq \unidades n$, se obtiene el siguiente corolario.

\begin{corolario}\label{C17.2.12'-1}
Para $n\geq 1$,
\begin{gather*}
\begin{align*}
\big|U_K/\unidades n\big|&=\prod_{i=0}^{n-1}\big|\unidades i/\unidades {i+1}
\big|=\big|U_K/\unidades 1\big|\cdot \big|\prod_{i=1}^{n-1}\big|\unidades i/
\unidades {i+1}\big|\\
&= \big|\*\F\big|\cdot\big|\prod_{i=1}^{n-1}\pK^i/\pK^{i+1}\big|=
\big|\*\F\big|\cdot \big|\prod_{i=1}^{n-1}\o_K/\pK\big|=(q-1)q^{n-1}.
\end{align*}
\intertext{Adem\'as, tenemos}
\o_K/\pK_K\stackrel{\cong}{\lra}\pK_K^n/\pK_K^{n+1},\quad x\longmapsto
x\pi_K^n
\intertext{y}
\big|\o_K/\pK_K\big|=\Big|\prod_{m=1}^n \pK_K^{m-1}/\pK_K^m\Big|=q^n.
\tag*{$\fin$}
\end{gather*}
\end{corolario}

\begin{ejemplo}\label{E17.2.12'}.
 Consideremos $K={\ma Q}_p$ el campo
de los n\'umeros $p$--\'adicos. El anillo de enteros de $K$ es el anillo de
los enteros $p$--\'adicos: $\AE K={\ma Z}_p$. 
Ahora bien $U_K=\*{\AE K}\cong
\*{\ma Z}_p$. Se tiene $\*{\ma Z}_p=
\big\{\xi\big(1+p \sum_{n=0}^{\infty}a_np^n\big)\mid \xi\in\*{\ma F}_p, a_n\in
{\ma F}_p\big\}\cong \*{\ma F}_p \times (1+p{\ma Z}_p)=\*{\ma F}_p\times
U_K^{(1)}$. Esto es $U_K^{(1)}\cong 1+p{\ma Z}_p$. Notemos que si
$p=2$, $U_K=\unidades 1\cong 1+2{\ma Z}_2$. De esta forma
tenemos que $\*{{\ma Z}_p}\cong \*{\ma F}_p \times (1+p{\ma Z}_p)$ y
\[
1+p{\ma Z}_p=
\begin{cases} {\ma Z}_p,& p\geq 3,\\
C_2\times (1+4{\ma Z}_2), &p=2.
\end{cases}
\]

M\'as en general, para $n\in{\ma N}$, se tiene $U_K^{(n)}=1+p^n{\ma Z}_p$.
\end{ejemplo}

\begin{definicion}\label{valorabsoluto}
Sea $K$ un campo cualquiera con un valor absoluto $|\ |_K$. Sea $V$
un espacio vectorial de dimensi\'on $n$, entonces una {\em norma\index{norma}}
$\norm$ de $V$ es funci\'on $\norm\colon V\lra {\ma R}^+\cup\{0\}$ tal que
\las
\item $\norm(v)\geq 0$ para toda $v\in V$ y $\norm(v)=0\iff v=0$.

\item $\norm(v_1+v_2)\leq \norm(v_1)+\norm(v_2)$ para cualesquiera 
$v_1,v_2\in V$.

\item $\norm(\alpha v)=|\alpha|_K\norm(v)$ para toda $\alpha\in K$ y toda
$v\in V$.
\end{list}
\end{definicion}

\begin{teorema}\label{unicidaddenorma}
Sean $K$ un campo valuado completo y $V$ un $K$-espacio vectorial
de dimensi\'on finita $n$. Entonces cualesquiera dos normas $\norm_1,
\norm_2$ de $V$ definen la misma topolog\'ia de $V$ y $V$ es un
espacio completo con la m\'etrica inducida.
\end{teorema}

\begin{proof}
Por inducci\'on en $n$. Para $n=1$, si $\{w_1\}$ es base de $V/K$ ($V
\cong K$), entonces si $\norm$ es una norma en $V$, entonces para $v=
\alpha_1w_1\in V$, $\norm(v)=|\alpha_1|_K\norm(w_1)=\lambda |
\alpha_1|_K$ con $\lambda=\norm(w_1)>0$. Por tanto $\norm=\lambda
|\ |_K$ y el resultado se sigue.

Sea $n>1$ y suponemos cierto el resultado para cualquier $K$-espacio
vectorial de dimensi\'on $m\leq n-1$. Sea $\{w_1,\ldots, w_n\}$ una
base de $V/K$. Sea $\norm_{\infty}$ la norma infinita de $V$, es decir,
\[
\norm_{\infty}(\alpha_1w_1+\cdots+\alpha_nw_n)=\max_{1\leq j\leq n}
\{|\alpha_j|_K\}.
\]
Se tiene que $\norm_{\infty}$ es una norma de $V$ y $V$ es completo
con respecto a $\norm_{\infty}$.

Sea $\norm$ una norma arbitraria de $V$. Sea $v=\alpha_1w_1+\cdots+
\alpha_nw_n$, entonces $\norm(v)=\norm(\alpha_1w_1+\cdots+\alpha_n
w_n)\leq \sum_{j=1}^n\norm(\alpha_jw_j)=\sum_{j=1}^n|\alpha_j|_K\norm(
w_j)\leq c\norm_{\infty}(v)$ donde $c=\sum_{j=1}^n\norm(w_j)>0$. Por
tanto $\norm(v)\leq c\norm_{\infty}(v)$ para toda $v\in V$.

Lo anterior prueba que la topolog\'ia definida por $\norm_{\infty}$ es m\'as
fina que la definida de $\norm$. Para probar que la topolog\'ia de $\norm$
es m\'as fina que la definida por $\norm_{\infty}$,
basta probar que dado $r>0$, existe $\delta=\delta(r)>0$
tal que si $v=\sum_{j=1}^n\alpha_jw_j$ es tal que $\norm(v)<\delta$, 
entonces $\norm_{\infty}(v)<r$.

Supongamos que no se cumple la condici\'on anterior, es decir, existe
$r>0$ tal que dado $\delta>0$, existe $v_{\delta}\in V$ con $\norm(
v_{\delta})<\delta$ pero $\norm_{\infty}(v_{\delta})\geq r$. Para 
$\delta_m=\frac 1m>0$, existe $v_{\delta_m}=:v_m$ con $\norm(
v_m)<\frac 1m$ pero $\norm_{\infty}(v_m)\geq r$. Consideremos $\{
v_m\}_{m=1}^{\infty}$, $\norm_{\infty}(v_m)=|\alpha_j^{(m)}|_K\geq r$
donde $v_m=\alpha_1^{(m)}w_1+\cdots+\alpha_n^{(m)}w_n$.
Puesto que $\{v_m\}_{m=1}^{\infty}$ es infinito, 
existe un \'indice $j_0$ y una subsucesi\'on $\{v_{m_t}\}_{
t=1}^{\infty}$ con $\norm_{\infty}(v_{m_t})=|\alpha_{j_0}^{(m_t)}|_K
\geq r$. Sin p\'erdida de generalidad, podemos suponer $j_0=n$
y $\{v_{m_t}\}_{t=1}^{\infty}=\{v_m\}_{m=1}^{\infty}$, es decir $\norm(
v_m)<\frac 1m$ pero $|\alpha_n^{(m)}|_K\geq r>0$.

Sea $u_m:=\frac {v_m}{\alpha_n^{(m)}}=\xi_1^{(m)}w_1+\cdots+
\xi_{n-1}^{(m)}w_{n-1}+w_n$ donde $\xi_j^{(m)}=\frac{\alpha_j^{(m)}}
{\alpha_n^{(m)}}$. Entonces $\norm(u_m)<\frac 1{mr}$, $r\leq |
\alpha_n^{(m)}|_K$, $\frac 1{|\alpha_n^{(m)}|_K}\leq \frac 1r$.

Adem\'as $u_m-u_t=\sum_{j=1}^{n-1}(\xi_j^{(m)}-\xi_j^{(t)})w_j$ y
$\norm(u_m-u_t)\leq \norm(u_m)+\norm(u_t)<\frac 1{mr}+\frac 1{tr}$.

Por hip\'otesis de inducci\'on, $\norm|_{V'}$ induce la misma 
topolog\'ia a la de $\norm_{\infty}|_{V'}$ donde $V'$ es el subespacio
de $V$ generado por $\{w_1,\ldots, w_{n-1}\}$. Ahora, dado $\epsilon
>0$, existe $t\in {\ma N}$ tal que para toda $m\geq t$, se tiene
$u_m-u_t\in B_{\norm}(0,\epsilon)$. Como $\norm_{\infty}$ induce
la misma topolog\'ia, existe $\delta >0$ tal que $u_m-u_t\in
B_{\norm_{\infty}}(0,\delta)\subseteq B_{\norm}(0,\epsilon)$ y
$\norm_{\infty}(u_m-u_t)=\max_{1\leq j\leq n-1}|\xi_j^{(m)}-\xi_j^{
(t)}|_K$ para toda $j$ con $1\leq j\leq n-1$. Por tanto $\{\xi_j^{(m)}
\}_{m=1}^{\infty}$ es una sucesi\'on de Cauchy y como $K$ es
completo, la sucesi\'on es convergente.

Sea $\xi_j^{(0)}=\lim\limits_{\substack{m\to\infty\\ |\ |_K}}\xi_j^{(m)}$
y sea $u_0:=\xi_1^{(0)}w_1+\cdots+\xi_{n-1}^{(0)}w_{n-1}+w_n$.
Entonces $\norm(u_0-u_m)=\norm\big(\sum_{j=1}^{n-1}(\xi_j^{
(0)}-\xi_j^{(m)})w_j\big)\leq \sum_{j=1}^{n-1}|\xi_j^{(0)}-\xi_j^{(
m)}|_K\norm(w_j)\xrightarrow[m\to\infty]{}0$. Por tanto $\norm_{
\infty}(u_0-u_m)\xrightarrow[m\to\infty]{}0$. Se sigue que
$u_m\xrightarrow[m\to\infty]{\norm_{\infty}} u_0$.

En particular, tenemos 
\[
\norm(u_0)=\norm(u_0-u_m+u_m)\leq \norm(u_o-u_m)+
\norm(u_m)\m 0,
\]
puesto que $\norm(u_m)<\frac 1{mr}$. Se sigue que $u_0=0=
\xi_1^{(0)}w_1+\cdots+\xi_{n-1}^{(0)}w_{n-1}+w_n$ lo cual
contradice que $\{w_1,\ldots,w_n\}$ es un conjunto linealmente
independiente. $\fin$
\end{proof}

\begin{teorema}\label{unicidaddevalorabsoluto}
Sea $K$ un campo valuado completo y sea $L/K$ una extensi\'on
finita de grado $n$. Entonces existe una \'unica extensi\'on $|\ |_L$
a $L$ del valor absoluto $|\ |_K$ de $K$ el cual est\'a dado 
por 
\[
|y|_L=|\N_{L/K}(y)|_K^{1/n},
\]
para $y\in L$ y $\N_{L/K}\colon L\to K$ es la norma. M\'as a\'un,
$L$ es un campo completo con el valor absoluto $|\ |_L$. 

Finalmente, la valuaci\'on normalizada $v_L$ est\'a dada por
\[
v_L(y):=\frac en v_K(\N_{L/K}(y)).
\]
\end{teorema}

\begin{proof}
La existencia de $v_L$ est\'a garantizada por el Lema de
Chevalley \cite[Corollary 2.4.5]{Vil2006}.

Ahora bien, si $\alpha\in\*L$, sea $\beta=\frac {\alpha^n}{\N(\alpha)}$,
donde $\N=\N_{L/K}$. Entonces 
\[
\N(\beta)=\frac{N(\alpha^n)}{(\N(\alpha))^n}=\frac{(\N(\alpha))^n}
{(\N(\alpha))^n}=1.
\]

Veamos que $|\gamma|_L<1$ implica $|\N(\gamma)|_K<1$. Sea
$\gamma^t=x_1^{(t)}w_1+\cdots+x_n^{(t)}w_n$ donde $\{w_1,
\ldots,w_n\}$ es una base $L/K$ y donde $x_j^{(t)}\in K$.

Puesto que $|\gamma|_L<1$, se tiene $\gamma^t\xrightarrow[t\to
\infty]{} 0$. Puesto que $|\ |_L$ y $\norm_{\infty}$ definen la 
misma topolog\'ia, se tiene que $\norm_{\infty}(\gamma^t)=
\max_{1\leq j\leq n}|x_j^{(t)}|\xrightarrow[t\to\infty]{} 0$. Se sigue
que $x_j^{(t)}\xrightarrow[t\to\infty]{}0$ para toda $1\leq j\leq n$.
Por tanto $|\N(\gamma)|_K<1$.

Ahora si $|\gamma|_L>1$, entonces $|1/\gamma|_L<1$ por lo que
$|\N(1/\gamma)|_K<1$ y por tanto $|\N(\gamma)|_K<1$.

En resumen, si $\N(\gamma)=1$, necesariamente $|\gamma|_L=1$.
Puesto que $\N(\beta)=1$, se tiene
\[
1=|\beta|_L=\frac{|\alpha|_L^n}{|\N(\alpha)|_K}.
\]
Por tanto $|\alpha|_L=\sqrt[n]{|\N(\alpha)|_K}$.

Finalmente $v_L(\beta)=0=v_L\big(\frac{\alpha^n}{\N(\alpha)}\big)=
nv_L(\alpha)-v_L(\N(\alpha))=nv_L(\alpha)-ev_K(\N(\alpha))$ por
lo que $v_L(\alpha)=\frac en v_K(\N(\alpha))$.
$\fin$
\end{proof}

\begin{corolario}\label{C17.3.2.7}
En el caso de un campo local, $n=ef$ (Teorema {\rm{\ref{T17.3.2.4}}})
por lo que $v_L(\alpha)=\frac 1f v_K(\N(\alpha))$. $\fin$
\end{corolario}

\begin{corolario}\label{C17.3.2.8}
Si $K$ es un campo valuado completo y $L/K$ es una extensi\'on
finita, entonces
\las
\item Si $a,b\in L$ son conjugados, entonces $|a|_L=|b|_L$ y 
$v_L(a)=v_L(b)$.
\item Si $L/K$ is Galois con grupo $G$, entonces para toda
$\sigma \in G$ se tiene $v_L\circ \sigma =v_L$.
\end{list}
\end{corolario}

\begin{proof}
\las
\item Si $F$ es la cerradura normal de $L/K$ y si $b=\sigma a$,
$|b|_L=|\N_{L/K}(\sigma a)|_K^{1/n}=|\N_{L/K}(a)|_K^{1/n}=|a|_L$.
\item Se tiene 
\begin{align*}
(v_L\circ \sigma)(\alpha)&=v_L(\sigma \alpha)=\frac en v_K(
\N_{L/K}(\sigma \alpha))\\
&=\frac en v_K(\N_{L/K}(\N_{L/K}(\alpha))=
v_L(\alpha)
\end{align*}
de donde se sigue que $v_L\circ \sigma=v_L$. $\fin$
\end{list}
\end{proof}

\subsection{Propiedades de las unidades de un 
campo local}\label{CClaseS1.2.1}

Sea $K$ un campo local con campo residual $\tilde K\cong \F$.

Notemos que $U_K^{(n)}=1+\pK_K^n$, $n\geq 0$, es un subgrupo abierto
de $\o_K$ pues
\begin{align*}
\unidades n&=\{x\in K^{\ast}\mid v_K(x-1)>n-1\}\\
&=\{x\in K^{\ast}\mid
|x-1|<q^{-(n-1)}\}=B(1,q^{-(n-1)}),
\end{align*}
y claramente $\bigcap_{n=0}^{\infty} \unidades n=\{1\}$
por lo que $\big\{\unidades n\big\}_{n\geq 0}$ forman
un sistema fundamental de vecindades abiertas de $1\in K^{\ast}$.

Por otro lado $\unidades n$, $n\geq 0$, es un subgrupo de 
$\*K$. Sea $R$ un conjunto de representantes de $\*K/\unidades n$
con $1\in R$,
$\unidades n=\*K\setminus \Big(\bigcup_{r\in R, r\neq 1}r\unidades n
\Big)$. Se sigue que $\unidades n$ es cerrado. Resumimos lo anterior
en

\begin{proposicion}\label{CCUnidades}
Los conjuntos $\unidades n$, $n\geq 0$, son a la vez abiertos y
cerrados en $\*K$. $\fin$
\end{proposicion}

\begin{proposicion}\label{CCCompactos}
El conjunto $\unidades n$, $n\geq 0$, es compacto.
\end{proposicion}

\begin{proof}
Sea $I$ una colecci\'on de subconjuntos abiertos de $\*K$ que
recubren a $\unidades n$. Si $\unidades n$ no pudiesen ser 
cubiertas por un n\'umero finito de conjuntos de $I$ y puesto que
$\unidades {n+1}$ es de \'indice finito en $\unidades n$, alguna clase
$u\unidades {n+1}\subseteq \unidades n$ no ser\'ia cubierta por
un n\'umero finito de elementos de $I$. Continuando con este
proceso, tenemos que existe una sucesi\'on
\[
u_1\unidades {n+1}\supseteq u_2\unidades {n+2}\supseteq \ldots,
\]
tal que ning\'un $u_j\unidades {n+j}$ puede ser cubierto por un
n\'umero finito de elementos de $I$.

Puesto que $\unidades n$ es un subconjunto cerrado de $\*K$
y $\*K$ es completo, se tiene que $\unidades n$ es completo,
de donde se sigue que $\bigcap_{j=1}^{\infty} u_j\unidades {n+j}
\neq \emptyset$. Sea $u_0\in \bigcap_{j=1}^{\infty} u_j\unidades {n+j}$.
Entonces $u_0\unidades {n+j}=u_j\unidades {n+j}$ para todo $j\geq 1$.
Los conjuntos
de $I$ son conjuntos abiertos y $u_0\unidades {n+j}$ es una
vecindad fundamental de $u_0$, por lo que existe $T\in I$ con 
$u_0\unidades {n+j}\subseteq T$. Esto \'ultimo
contradice que $u_j\unidades{n+j}=u_0\unidades{n+j}$ no puede
ser cubierto por un n\'umero finito de elementos de $I$. Por
tanto, $\unidades n$ es un conjunto compacto.
$\fin$
\end{proof}

Con respecto a los ideales $\pK_K^n$ tenemos:

\begin{proposicion}\label{idealesabiertos}
Los ideales $\pK_K^n$, $n\geq 0$ son a la vez conjuntos abiertos
y cerrados de $\*K$.
\end{proposicion}

\begin{proof} Se sigue inmediatamente de que $\unidades n=1+
\pK_K^n$ y de la Proposici\'on \ref{CCUnidades}. Una demostraci\'on
directa es la siguiente. Se tiene que para $n\geq 0$, 
$\pK_K^n=\{x\in K\mid v_K(x)\geq n\}=\{x\in K\mid
v_K(x)>n-1\}$, es decir, $\pK_K^n=v_K^{-1}\big((n-1,\infty]\big)=
v_K^{-1}\big([n,\infty]\big)$, $v_K\colon K\to {\ma R}\cup\{0\}=(-\infty,
\infty]$ es claramente una funci\'on continua (si $x_n\xrightarrow[n\to\infty]{}
x$, $v_K(x_n)\xrightarrow[n\to \infty]{}v_K(x)$). Puesto que
$(n-1,\infty]$ es abierto y $[n,\infty]$ es cerrado en
${\ma R}\cup\{0\}$ se sigue lo afirmado. $\fin$
\end{proof}

\begin{corolario}\label{CCdisconexo}
Sea $K$ un campo local. Entonces $K$ es un campo
totalmente disconexo y no discreto. El grupo multiplicativo $\*K$
de $K$, es localmente compacto, no compacto y totalmente 
disconexo.
\end{corolario}

\begin{proof}
Puesto que $\unidades n$, $n\geq 0$ es un sistema fundamental de
vecindades de $1\in \*K$, y puesto que $\unidades n$ es compacto, 
se sigue que $\*K$ es localmente compacto. Ahora bien $\*K$
no es compacto pues si $\pi=\pi_K$, $\*K=\bigcup_{n=0}^{\infty}
\pi^n U_K$ donde cada conjunto $\pi^n U_K$ es
abierto y la uni\'on es una
uni\'on disjunta a pares. El hecho de que $\*K$ es totalmente
disconexo se prueba de manera an\'aloga a que $K$ es
totalmente disconexo y cuya prueba damos enseguida.

Para ver que $K$ es totalmente disconexo, consideremos $x\neq 0$.
Entonces $v_K(x)=n\in{\ma Z}$. Entonces $x\notin\pK_K^m$ con 
$m\geq \max\{1,n+1\}$, esto es $0\in\pK_K^m$, $x\in K\setminus
\pK_K^m$ siendo $\pK_K^m$ y $K\setminus\pK_K^m$ conjuntos
abiertos disjuntos. Finalmente $K$ no es discreto pues $\pK_K^n$,
$n\geq 0$ es un sistema fundamental de vecindades de $\{0\}$ y
$\pK_K^n\neq \{0\}$, esto es, $\{0\}$ no es un conjunto abierto.
$\fin$
\end{proof}

El isomorfismo $\tilde{K}\cong \pK_K^n/\pK_K^{n+1}$ se sigue del
hecho de que $\pK=\pi\o_K$, por lo que la multiplicaci\'on por 
$\pi^n$ induce el isomorfismo. M\'as precisamente, sea 
$\varphi\colon \o_K\lra \o_K$ dada por $\varphi(x)=\pi^n x$. Entonces
$\varphi$ es un homomorfismo de grupos aditivos. Adem\'as,
$\im \varphi=\pi^n\o_K=\pK_K^n$ y $\varphi^{-1}(\pK_K^{n+1})=
\pi\o_K=\pK_K$, de donde se sigue que $\o_K/\pK_K\cong
\pK_K^n/\pK_K^{n+1}$.

Adem\'as el mapeo $\psi\colon\o_K\lra \o_K$ dado por $\psi(u)
=u-1$ da lugar al homomorfismo de grupos $\unidades n\lra
\pK_K^n/\pK_K^{n+1}$ con $\unidades n$ un grupo multiplicativo
y $\pK_K^n/\pK_K^{n+1}$ un grupo aditivo, el cual es suprayectivo
por definici\'on y $\ker \psi=1+\pK_K^{n+1}$. De hecho, lo anterior
se debe a que 
\[
\psi(uv)=uv-1=(u-1)(v-1)+(u-1)+(v-1)
\]
por lo que $uv-1\equiv (u-1)+(v-1)\bmod \pK_K^{2n}$ y
$2n\geq n+1$.

Como consecuencia se tiene:

\begin{proposicion}\label{CClaseP1.2.1.1} Sea $p$ la caracter\'istica
del campo residual del campo local $K$.
\las
\item Para todo $n\geq 1$, $(\unidades n)^p\subseteq \unidades {n+1}$.
\item Para $m\in{\ma N}$ con $\mcd(m,p)=1$, se tiene que para cada
$n\geq 1$, el mapeo $u\lra u^m$ es un automorfismo de $\unidades n$,
en particular $\big(\unidades n\big)^m=\unidades n$ para toda $n\geq 1$
y toda $m$ con $p\nmid m$.
\end{list}
\end{proposicion}

\begin{proof}
(1) Se tiene $\unidades n/\unidades {n+1}\stackrel{\theta}{\cong}\F$.
Por tanto, si $x\in \unidades n$, se tiene $\theta(\bar{x}^p)=p
\theta(\bar{x})=0$, de lo cual obtenemos que $x^p\in\unidades {n+1}$.

(2) Sea $f\colon \unidades n\lra \unidades{n}$ dado por 
$f(u)=u^m$, con $n\geq 1$ arbitrario. Entonces $f$ es un endomorfismo.
Sea $x\in\unidades n$ con $x^m=1$. Se tiene $x=1+u\pi^n$ con $u
\in \o_K$. Si $u\neq 0$, $u=\pi_K^su_0$ para algunos $u_0\in U_K$
y $s\geq 0$. Por tanto 
\[
x=1+u_0\pi_K^{n+s}\quad \text{y}\quad x^m=1+mu_0\pi_K^{
n+s}+\text{\ t\'erminos con valuaci\'on mayor}.
\]
Puesto que $v_K(m)
=0$, $mu_0\neq 0$ lo cual es absurdo. Por tanto $u=0$ y $x=1$. Por
tanto $f$ es inyectiva.

Ahora sea $x\in\unidades n$. Para toda $r\in{\ma N}$, existen $a_r,b_r
\in{\ma Z}$ tales que $a_rm+b_rp^r=1$, $x=x^1=x^{a_rm+b_rp^r}=
(x^{a_r})^m\cdot (x^{b_r})^{p^r}=\xi_r\mu_r$. Sea $x^{b_r}=1+c_r$ con
$c_r\in \pK_K$. Entonces $\mu_r=(x^{b_r})^{p^r}\equiv 1+p^rc_r\bmod
p^{r+1}\pK_K$ y $v_K(\mu_r-1)\geq r+1$ y $\lim_{r\to\infty}\mu_r=1$. Se
sigue que $y:=\lim_{r\to\infty}x^{a_r}$ existe y $x=y^m$, por lo que
$f$ es suprayectiva.
$\fin$
\end{proof}

\begin{proposicion}\label{CClaseP1.2.1.2'} Sean $v=v_K$ la valuaci\'on
de $K$ y $m\in{\ma N}$ tal que la caracter\'istica de $K$ no divida a $m$.
Entonces para $n>v_K(m)$, el mapeo $g$, $x\mapsto x^m$ da lugar a
un isomorfismo $\unidades n\lra \unidades {n+v_K(m)}$. En particular,
si $\ca K=p>0$, $v_K(m)=0$ y $\unidades n=\big(\unidades n\big)^m$
para toda $n\geq 1$ y para $p\nmid m$.
\end{proposicion}

\begin{proof} Sea $\pi$ un elemento primo de $K$, esto es, $v_K(\pi)=1$.
Sea $x=1+a\pi^n\in\unidades n$, $a\in\o_K$. Entonces
\[
x^m=1+ma\pi^n+\binom{m}{2}a^2\pi^{2n}+\cdots\equiv 1\bmod
\pK_K^{n+v_K(m)}.
\]
Por tanto $x^m\in \unidades {n+v_K(m)}$ para toda $n$.

Para probar que el mapeo es suprayectivo, debemos probar que
para $a\in\o_K$ arbitrario, existe $x\in\o_K$ tal que 
\[
1+a\pi^{n+v_K(m)}=(1+x\pi^n)^m,
\]
lo cual equivale a que $1+a\pi^{n+v_K(m)}=1+m\pi^n x+\pi^{2n}f(x)$
donde $f(x)\in\o_K[x]$. Si $K$ es de caracter\'istica prima relativa a
$m$, $m=u\pi^{v_K(m)}$ con $u\in U_K$. Si $\car K=p>0$ y $p|m$,
$m=0$ en $K$. As\'i
\begin{gather*}
1+a\pi^{n+v_K(m)}=
\begin{cases}
1+u\pi^{n+v_K(m)}x+\pi^{2n}f(x)&\text{si $p\nmid m$},\\
1+\pi^{2n}f(x)&\text{si $p\mid m$}.
\end{cases}
\intertext{Por lo tanto necesitamos resolver}
0=\begin{cases}
-a+ux + \pi^{n-v_K(m)}f(x)&\text{si $p\nmid m$},\\
-a+\pi^{n-v_K(m)}f(x)&\text{si $p\mid m$}.
\end{cases}
\end{gather*}

La primera ecuaci\'on tiene soluci\'on por el Lema de Hensel,
pero la segunda ecuaci\'on no tiene soluci\'on. 
M\'as precisamente, la primera ecuaci\'on es equivalente
a resolver la ecuaci\'on $-a+ux+\pi^{n-v_K(m)}f(x)=0$.
Puesto que $n>v_K(m)$, m\'odulo $\pK_K$ se tiene 
$\bar a+\bar u\bar x=0$ con $\bar u\neq 0$. Por tanto
$\overline{h(x)}=\overline{-a+ux+\pi^{n-v_K(m)}f(x)}=0$ tiene una
ra\'iz simple y por el Lema de Hensel, Teorema \ref{Lema de Hensel},
la ra\'iz $\bar x=\frac{\bar a}{\bar u}$ se levanta a $\o_K$ a
una ra\'iz de $h(x)$ y por tanto $g\colon \unidades n\lra
\unidades {n+v_K(m)}$, $x\mapsto x^m$, es suprayectiva.

Ahora, $\ker g=W_m=\{\xi\in\unidades n\mid \xi^m=1\}$. Sea
$\xi\in\ker g$, $\xi=1+u\pi^n$ con $u\in\o_K$. Entonces $\xi\equiv
1\bmod \pK_K^n$, por lo tanto $\xi^r\equiv 1\bmod \pK_K^n$ para
toda $r\geq 1$. Se tiene
\[
0=(\xi^m-1)=(\xi-1)(\xi^{m-1}+\cdots+\xi+1).
\]
Puesto que $n>v_K(m)$ y $\xi^{m-1}+
\cdots+\xi+1\equiv m\bmod \pK_K^n=
u\pi^{v_K(m)}\bmod \pK_K^n$ con 
$u\in U_K$, se sigue que $\xi^{m-1}+\cdots+
\xi+1\notin \pK_K^n$ y en particular $\xi^{m-1}+\cdots+\xi+1
\neq 0$. Por tanto $\xi=1$ y $g$ es inyectiva.
$\fin$
\end{proof}

\begin{corolario}\label{mpotencias}
Para $m\in{\ma N}$ tal que $\ca K\nmid m$, se tiene que 
$\big(\*K\big)^m$ es un subgrupo abierto de $\*K$.

Adem\'as, si $\ca K=0$, se tiene que $\bigcap_{m=1}^{\infty}
\big(\*K\big)^m=\{1\}$.
\end{corolario}

\begin{proof}
Para $x\in\*K$, $x^m\in\big(\*K\big)^m$, y para $n>v_K(m)$ con
$p\nmid m$, $p=\ca K$, $x^m\unidades {n+v_K(m)}=\big(
x\unidades n\big)^m\subseteq \big(\*K\big)^m$, la cual es una
vecindad abierta de $x^m$. Se sigue que $\big(\*K\big)^m$
es abierto en $\*K$.

Sea ahora $\ca K=0$ y sea $a\in\bigcap_{m=1}^{\infty}
\big(\*K\big)^m\subseteq \bigcap_{\substack{m=1\\ p\nmid m}}^{
\infty}\big(\*K\big)^m$. Si $v_K(a)\neq 0$, para cada $m\in
{\ma N}$ con $p\nmid m$, existe $y_m\in\*K$ con $a=y_m^m$,
por lo tanto $v_K(a)=mv_K(y_m)$ con $v_K(y_m)\neq 0$.
Por tanto, $m|v_K(a)$ para toda $m$ con $p\nmid m$, lo cual
es absurdo. Se sigue que $v_K(a)=0$ y $a\in U_K$. Por tanto
$a\in\bigcap_{m=1}^{\infty}\big(U_K)^m$, esto es, $a=u_m^m$
con $u_m\in U_K$ para toda $m$. 

Sea $n\in {\ma N}$ y sea
$[U_K:\unidades n]=m$. Entonces $u_m^m\in \unidades n$
de donde obtenemos que $a\in\bigcap_{m=1}^{\infty} \unidades n$
y $v_K(a-1)\geq n$ para toda $n\in{\ma N}$. Por tanto $v_K(a-1)
=\infty$ y $a=1$.
$\fin$
\end{proof}

\begin{observacion}\label{noppotencia}
Si $K$ es un campo local de caracter\'istica $p>0$, entonces 
$\big(\*K\big)^p$ no es abierto en $\*K$ pues si este fuese
el caso, existir\'ia $m$ tal que $\unidades m\subseteq
\big(\*K\big)^p$. En particular $1+\pi_K^m=x^p$
para alg\'un $x\in\*K$. De esta forma, tendr\'iamos que
$x^p-1=(x-1)^p=\pi_K^m$ lo cual no es posible puesto
que de lo contrario tendr\'iamos que $m=pv_K(x-1)$ y
entonces $p|m$. Puesto que $\unidades n\subseteq
\unidades m$ esto dir\'ia que $p|n$ para toda $n\geq m$
lo cual es absurdo.

Por otro lado, $\big(\*K\big)^p$ si es cerrado pues
$\*K\cong\langle\pi_K\rangle \times U_K$ como grupos
y $\big(\*K\big)^p=\langle\pi_K^p\rangle U_K^p$ y como
$U_K$ es compacto, $U_K^p$ es compacto y por tanto
cerrado.

Cuando $\ca K=0$, un resultado que probaremos usando el
cociente de Herbrand en cohomolog\'ia es que para $m\in
{\ma N}$, $\big(\*K\big)^m$ es de \'indice finito en $\*K$ y
de hecho $[\*K:\big(\*K\big)^m]=mq^{v_K(m)}|\mu_m(K)|=
m\cdot |m|_{\pK_K}^{-1}\cdot |\mu_m(K)|$, donde $|\tilde K|
=q$ y $\mu_m(K)=\{\xi\in\*K\mid \xi^m=1\}$.
\end{observacion}

\subsection{Representaci\'on, normas, valores absolutos
y extensiones de campos locales}\label{S17.3.2}

Empezamos con un resultado fundamental sobre extensiones
finitas de campos locales.

\begin{proposicion}\label{P17.3.2.1}
Sea $L/K$ una extensi\'on finita de campos locales con
anillos de valuaci\'on $\o_L$ y $\o_K$ respectivamente.
Sean $\pK_K=\langle\pi_K\rangle$ y $\pK_L=\langle
\pi_L\rangle$ con $\pi_K$ y $\pi_L$ elementos primos
de $K$ y de $L$ respectivamente. Sean $f=f(L|K)=
[\o_L/\pK_L:\o_K/\pK_K]=[\tilde L:\tilde K]$ y $e=e(L|K)$
el \'indice de ramificaci\'on. Sea $\{\omega_1,\ldots,
\omega_f\}$ una base de $\tilde L/\tilde K$. Sean 
$\alpha_i\in\o_L$ tales que $\bar{\alpha}_i=\omega_i$,
$1\leq i\leq f$. Sea $\eta_{ij}=\alpha_i\pi_L^j$, $1\leq i
\leq f$, $0\leq j\leq e-1$. Entonces $\{\eta_{ij}\}_{i,j}$ es
un sistema linealmente independiente sobre $\*K$ y
en particular $[L:K]\geq ef$. M\'as a\'un se tiene que:
\las
\item Si $y=\sum_{i=1}^f A_i\alpha_i$ con $A_i\in K$,
entonces se tiene que $v_L(y)=\min_{1\leq i\leq f}\{ev_K(A_i)\}$.

\item Si $z=\sum_{i=1}^f\sum_{j=0}^{e-1}c_{ij}\eta_{ij}$ con $c_{ij}\in K$, 
entonces 
\[
v_L(z)=\min_{\substack{1\leq i\leq f\\
0\leq j\leq e-1}}\{ev_K(c_{ij})+j\}.
\]
\end{list}
\end{proposicion}

\begin{proof}
\las
\item Sea $y=\sum_{i=1}^f A_i\alpha_i$ con $A_i\in K$. Si
$A_1=\cdots A_f=0$, el resultado se sigue inmediatamente.
Por tanto, rearreglando, podemos suponer que $v_K(A_1)=
\cdots =v_K(A_s)<v_K(A_i)$, $s+1\leq i\leq f$. Sea $B_i:=
\frac{A_i}{A_1}$, $1\leq i\leq f$. Entonces $B_1=1$ y $B_i
\in\o_K$ para toda $i$. Adem\'as se tiene que $\frac{y}{A_1}
=\sum_{i=1}^f B_i\alpha_i\in \o_L$. Si $y/A_1\in\pK_K$ entonces
se tendr\'ia que $\frac{\bar{y}}{\bar {A}_1}=0=\bar{\alpha_1}+
\sum_{i=1}^f\bar{B_i}\bar{\alpha_i}$ lo que contradice que
$\{\omega_1,\ldots, \omega_f\}$ es un conjunto linealmente
independiente sobre $\tilde K$. Se sigue que $y/A_1\notin
\pK_L$ por lo que $v_L\big(y/A_1\big)=0$ y $v_L(y)=v_L(
A_1)=ev_K(A_1)=e\min_{1\leq i\leq f}\{v_K(A_i)\}=
\min_{1\leq i\leq f}\{ev_K(A_i)\}$.
\item Sea $z=\sum_{ij}c_{ij}\eta_{ij}$ con $c_{ij}\in K$. Sea
$y_j=\sum_{i=1}^f c_{ij}\alpha_i$, $0\leq j\leq e-1$ y por 
tanto $z=\sum_{j=0}^{e-1}y_j\pi_L^j$. Si $y_j\neq 0$, 
$v_L(y_j)$ es m\'ultiplo de $e$ y $v_L(y_j\pi_L^j)\equiv
j\bmod e$. Por tanto los valores $\{v_L(y_j)\mid y_j\neq
0\}_{j=0}^{e-1}$ son todos distintos y por tanto
\begin{align*}
v_L(z)&=\min_{0\leq j\leq e-1}\{v_L(y_j\pi_L^j)\}=
\min_{0\leq j\leq e-1}\{v_L(y_j)+j\}\\
&=\min_{0\leq j\leq e-1}\{\min_{1\leq i\leq f}\{ev_K
(c_{ij})+j\}\}=\min_{\substack{0\leq j\leq e-1\\ 1\leq i\leq f}}
\{ev_K(c_{ij})+j\}.
\end{align*}
\end{list}

Ahora sea $\sum_{ij}x_{ij}\eta_{ij}=0$ para algunos $x_{ij}\in K$.
Por tanto $\infty=v_L(0)=\min_{ij}\{ev_K(x_{ij})+j\}$. Por tanto
$v_K(x_{ij})=\infty$ para todo $i, j$ y por tanto $x_{ij}=0$ para
todo $i,j$, probando que $\{\eta_{ij}\}_{ij}$ es un conjunto
linealmente independiente de $L$. $\fin$
\end{proof}

Ahora probaremos que $\{\eta_{ij}\}_{ij}$ es de hecho una
base de $L/K$. Esto lo haremos viendo que los campos locales,
de hecho los campos completos bajo una valuaci\'on discreta,
se pueden representar por medio de series de Laurent.

Sea $K$ un campo completo con respecto a una valuaci\'on
discreta $v_K$ normalizada, esto es, $v_K(\*K)={\ma Z}$ y
sea $\o_K/\pK_K=\tilde K$ el campo residual. Sea $R$ un
sistema completo de representantes de $\tilde K$, $R
\subseteq \o_K$ y tal que $0\in R$, $0$ es el representante
de $0+\pK_K=\pK_K$ en $R$. Para cada $n\in {\ma Z}$, sea
$\pi_n\in K$ tal que $v_K(\pi_n)=n$.  Consideremos cualquier
suma $\sum_{n=m}^{\infty}r_n\pi_n$ donde $m\in{\ma Z}$ y
$r_n\in R$. Se tiene que $v_K(r_n\pi_n)=v_K(r_n)+v_K(\pi_n)
\geq n$ ($= n$ si $r_n\neq 0$ y $\infty$ si $r_n=0$). Por tanto
$v_K(r_n\pi_n)\xrightarrow[n\to\infty]{}\infty$ y por tanto
$\sum_{n=m}^{\infty}r_n\pi_n$ converge en $K$.

Entenderemos $\sum_{n=m}^{\infty}r_n\pi_n=\sum_{-\infty}^{
\infty}r_n\pi_n$ con $r_n=0$ para toda $n<m$.

\begin{teorema}\label{T17.3.17.2} 
\las
\item Cada $x\in K$ se expresa de manera \'unica como
$x=\sum_{n=m}^{\infty} r_n\pi_n$, $r_n\in R$. Si $x\neq 0$
y $r_m\neq 0$, entonces $v_K(x)=m$.

\item Sean $x=\sum_{n=m}^{\infty}r_n\pi_n$, $y=\sum_{n=m}^{
\infty}s_n\pi_n$ con $r_n,s_n\in R$. Entonces para $i\in{\ma Z}$,
se tiene que $v_K(x-y)\geq i \iff r_n=s_n$ para toda $n<i$.
\end{list}
\end{teorema}

\begin{proof}
\las
\item La unicidad se sigue de (\ref{(2)}): Si $x=\sum r_n\pi_n=
\sum r'_n\pi_n$ fuesen dos expresiones de $x$ y si $r_m\neq
r'_m$ para alg\'un $m\in{\ma Z}$, podemos seleccionar el
m\'inimo tal $m$ y $v_K(0)=v_K(x-x)=m$ lo cual es absurdo.

La valuaci\'on de $x$ tambi\'en se sigue de (\ref{(2)}). Se tiene
que $0=\sum 0 \pi_n$ y $x=\sum r_n\pi_n$, $x\neq 0$, si $m\in
{\ma Z}$ es m\'inimo $m\in{\ma Z}$ con $r_m\neq 0$, Entonces
$v_K(x-0)=v_K(x)=m$.

Para la representaci\'on de $x\in K$, $x\neq 0$ y $v_K(x)=m
<\infty$, se tiene, para $n\in{\ma Z}$, $\pK_K^n=\{x\in K\mid
v_K(x)\geq n\}=\o_K\cdot \pi_n$. De que $\o_K=R+\pK_K$,
se sigue que
\[
\pK_K^n=R\pi_n+\pK_K^{n+1}=R\pi_n+R\pi_{n+1}+\cdots +
R\pi_t+\pK_K^{t+1}
\]
para toda $t\geq n$.
Puesto que $x\in \pK_K^m$, existen $r_m,r_{m+1},r_{m+2},
\ldots $ tales que $x\equiv \sum_{n=m}^tr_n\pi_n\bmod \pK_K^{
t+1}$ para $t\geq m$, de donde, $x=\lim_{t\to\infty}\sum_{n=
m}^{\infty}r_n\pi_n$.

As\'i pues, toda la demostraci\'on se seguir\'a de probar
(\ref{(2)}).

\item\label{(2)} Si $r_n=s_n$ para toda $n$, entonces $x=y$
y $v_K(x-y)=\infty$. Supongamos ahora que existe $m\in{\ma Z}$
con $r_m\neq s_m$ y $r_n=s_n$ para toda $n<m$. As\'i,
$x-y=\sum_{n=m}^{\infty}(r_n-s_n)\pi_n$. Puesto que $r_m,s_m
\in R$, $r_m\not\equiv s_m\bmod\pK$ por lo que $v_K(r_m-
s_m)=0$ y $v_K((r_m-s_m)\pi_m)=m$.

Para $n>m$, $v_K((r_n-s_n)\pi_n)\geq v_K(\pi_n)=n>m=v_K((
r_m-s_m)\pi_m)$. Puesto que $m=v_K((r_m-s_m)\pi_m)<v_K\big(
\sum_{n=m+1}^{\infty}(r_n-s_n)\pi_n\big)(\geq m+1)$, se sigue
que $v_K(x-y)=v_K((r_m-s_m)\pi_m)=m$.
Esto prueba (\ref{(2)}). $\fin$
\end{list}
\end{proof}

\begin{corolario}\label{C17.3.17.3}
Sea $R^{\infty}=\prod_{i=0}^{\infty}R=\{(r_i)_{i=0}^{\infty}\mid
r_i\in R\}=\prod_{n=0}^{\infty} R_n$ con $R_n=R$ para toda
$n$. Consideremos $R^{\infty}$ con la topolog\'ia producto
y cada $R_n$ tiene la topolog\'ia discreta. Entonces
\[
R^{\infty}\xrightarrow[]{\ \mu\ }\o_K,\qquad (r_n)_{n=
0}^{\infty}\longmapsto
\sum_{n=0}^{\infty}r_n\pi_n,
\]
es un homeomorfismo de $R^{\infty}$ en $\o_K$.
En particular $\o_K$ es un conjunto compacto.
Adem\'as, $\o_K\cong \pK_K^n$ para toda $n\geq 0$ y
$\pK_K^n$ es compacto.
\end{corolario}

\begin{proof}
La parte (1) del Teorema \ref{T17.3.17.2} prueba que $\mu$ es
biyectiva y la parte (2) prueba que es un homeomorfismo.
La compacidad de $\o_K$ se sigue del Teorema de Tychonoff.

El mapeo $x\longmapsto x\pi_K^n$ es un homeomorfismo de
$\o_K$ en $\pK_K^n$.
$\fin$
\end{proof}

Notemos que $\o_K=\bigcup_{n=0}^{\infty}\pi_K^n U_K$ y cada
$\pi_K^n U_K$ es abierto, por tanto $\o_K$ es abierto.
Por otro lado $\bigcap_{n=0}^{\infty}\pK_K^n=\{0\}$ pues si
$x\in\bigcap_{n=0}^{\infty}\pK_K^n$, $v_K(x)\geq n$ para toda
$n$ de donde se sigue que $v_K(x)=\infty$ y $x=0$.

Tomando los mapeos can\'onicos $\o_K/\pK_K^m\lra \o_K/\pK_K^n$
para $m\geq n$, se sigue que
\begin{gather*}
\o_K\cong \lim_{\substack{\leftarrow\\n}}\o_K/\pK_K^n.
\intertext{Similarmente}
\unidades 1\cong \lim_{\substack{\leftarrow\\ n}}\unidades 1/
\unidades n
\end{gather*}
con respecto a los mapeos can\'onicos $\unidades 1/\unidades m
\lra \unidades 1/\unidades n$ con $m\geq n$.

\begin{teorema}\label{T17.3.17.4}
Sea $K$ un campo local.
\las
\item Si $\ca K=0$, $[K:{\ma Q}_p]=ef=d$, $\unidades 1\cong W_K\oplus 
{\ma Z}_p^d$ con $W_K=\{\xi\in K\mid \xi^{p^m}=1 \text{\ para alg\'un\ }
m\}\cong {\ma Z}/p^a{\ma Z}$ para alg\'un $a\geq 0$.

\item Si $\ca K=p>0$, $\unidades 1\cong \prod_{i=1}^{\infty}{\ma Z}_p$.
\end{list}
\end{teorema}

\begin{proof}
\cite[Propositions 2.6, 2.7, 2.8]{Iwa86}.
$\fin$
\end{proof}

\begin{teorema}\label{T17.3.2.4} Con las notaciones de
la Proposici\'on {\rm{\ref{P17.3.2.1}}}, se tiene que $\{\eta_{ij}\}_{
\substack{1\leq i \leq f\\ 0\leq j\leq e-1}}$
es una base de $L/K$ y $[L:K]=ef$.

M\'as a\'un, $\{\eta_{ij}\}_{\substack{1\leq i \leq f\\ 0\leq j\leq e-1}}$
es una base
$\o_L$ sobre $\o_K$, esto es, $\o_L$ es un $\o_K$-m\'odulo
libre de rango ef.
\end{teorema}

\begin{proof}
Sea $R$ un conjunto de representantes de $\o_K/\pK_K$ en $\o_K$
con $0\in R$ y sea $R'=\{\sum_{i=1}^fr_i\alpha_i\mid r_1,\ldots, r_f
\in R\}$. Entonces $R'$ es un conjunto de representantes de
$\o_L/\pK_L$ en $\o_L$ y $0\in R'$. Sean $\pi_K$ y $\pi_L$ elementos
primos de $K$ y $L$ respectivamente.

Para $m\in{\ma Z}$, escribamos $m=te+j$, $t\in{\ma Z}$ y $0\leq j\leq
e-1$ y definimos $\pi_{L,m}:=\pi_K^t\pi_L^j$. Ahora como $v_L|_K=ev_K$,
se tiene $v_L(\pi_{L,m})=et+j=m$.

Dado $y\in L$, $y$ se escribe de manera \'unica en la forma 
\[
y=\sum_{m=-\infty}^{\infty}r'_m\pi_{L,m},
\]
con $r'_m\in R'$. Sea $r'_m=\sum_{i=1}^f r_{i,m}\alpha_i$
con $r_{i,m}\in R$. Se sigue que
\[
y=\sum_m\sum_{i=1}^f r_{i,m}\alpha_i\pi_{L,m}=
\sum_{i=1}^f\sum_{j=0}^{e-1}\beta_{ij}\alpha_i\pi_L^j,
\]
donde $\beta_{ij}=\sum_{i=-\infty}^{\infty}r_{i,te+j}\pi_K^i\in K$. Por tanto
$\{\eta_{ij}\}_{ij}$ genera $L/K$ y por la Proposici\'on \ref{P17.3.2.1},
este conjunto es linealmente independiente, por lo tanto es base
y $[L:K]=ef$.

Ahora bien $\{\eta_{ij}\}_{ij}\subseteq \o_L$. Sea $x\in \o_L$,
$x=\sum_{ij}c_{ij}\eta_{ij}$. Entonces, por la Proposici\'on
\ref{P17.3.2.1}, tenemos que $v_K(x)=\min_{ij}\{ev_K(c_{ij})
+j\}\geq 0$ por lo que
$ev_K(c_{ij})\geq -j$ para todas $i,j$. Si se tuviese $v_K(c_{ij})
\leq -1$, entonces $ev_K(c_{ij})=-le\geq -j$, $l\geq 1$ lo que
implicar\'ia que $j\geq e$ lo cual es absurdo. Por tanto $v_K(
c_{ij})\geq 0$ para todas $i, j$. Es decir $c_{ij}\in \o_K$ y el
resultado se sigue.
$\fin$
\end{proof}

\begin{teorema}\label{T17.3.2.6}
Sea $K$ un campo local y sea $L/K$ una extensi\'on finita de $K$.
Sean $\N=\N_{L/K}$ y $T:=\Tr_{L/K}$ la norma y la traza
de $L$ en $K$ respectivamente: $\N, T\colon L\lra K$. Entonces
$\N$ y $T$ son funciones continuas. Adem\'as $\N(\o_L)
\subseteq \o_K$ y $T(\o_L)\subseteq \o_K$.
\end{teorema}

\begin{proof}
Sea $\{\eta_{ij}\}_{ij}$ la base de $L/K$ dada por la Proposici\'on
{\rm{\ref{T17.3.2.4}}}. Sean $y\in L$ y $y\eta_{ij}=\sum_{s=1}^f
\sum_{t=0}^{e-1}x_{ijst}\eta_{st}$ para toda $i,j$ con $x_{ijst}\in K$.
Se tiene
\begin{gather*}
\begin{align*}
v_L(y)+j&=v_L(y)+v_L(\eta_{ij})=v_L(y\eta_{ij})=\min_{s,t}\{ev_K(x_{ijst})
+t\}\\
&\leq ev_K(x_{ijst})+t.
\end{align*}
\intertext{Por tanto}
v_K(x_{ijst})\geq \frac 1e(v_L(y)+j-t)\quad \text{para todas\ } i,j,s,t.
\end{gather*}
Se sigue que $y\lra 0\iff v_L(y)\lra \infty \Longrightarrow v_K(x_{ijst})
\lra \infty \iff x_{ijst}\lra 0$ para todas $1\leq i, s\leq f$ y para todas
$0\leq j, t\leq e-1$.

Se sigue que las $x_{ijst}$'s dependen continuamente de $y\in L$.
Puesto que $T$ y $\N$ son la traza y el determinante de la $n\times
n$ matriz $(x_{ijst})$ respectivamente, tenemos que $T$ y $\N$ son
funciones continuas.

Finalmente, si $y\in\o_L$, $v_L(y)\geq 0$ lo que implica que $v_K(
x_{ijst})\geq 0$ de donde obtenemos que $T(\o_L)\subseteq \o_K$
y $\N(\o_L)\subseteq \o_K$.
$\fin$
\end{proof}

Un resultado de especial relevancia para nosotros para la teor\'ia local
de campos de clase es que si $L/K$ es una extensi\'on finita de
campos locales, entonces $\N_{L/K}(\*L)$ es un conjunto cerrado
en $\*K$.

\begin{teorema}\label{T17.3.2.7}
Sea $L/K$ una extensi\'on finita de campos locales y se $\N$
la norma de $L$ en $K$. Entonces $\N(U_L)$ es un subgrupo
compacto y por tanto cerrado en $\*K$ y $\N(\*L)$ es un
subgrupo cerrado de $\*K$.
\end{teorema}

\begin{proof}
Puesto que la norma $\N$ es continua (Teorema \ref{T17.3.2.6})
y $U_L$ es compacto (Proposici\'on \ref{CCCompactos}) por lo
que $\N(U_L)$ es compacto en $\*K$ y $\N(U_L)\subseteq U_K$.
Puesto que $\*K$ es Hausdorff, $\N(U_L)$ es cerrado en $\*K$.

Ahora bien, $\N(U_L)=\N(\*L)\cap U_K\subseteq U_K$. Como 
$U_K$ es abierto, $\N(U_L)$ es abierto en $\N(\*L)$ y 
tambi\'en es cerrado en $\N(\*L)$.
Se sigue que $\N(\*L)/\N(U_L)$ es un espacio discreto: si $\xi\in
\N(\*L)$, $\xi\N(U_L)$ es abierto en $\N(\*L)$ y $\xi\N(U_L)=
\pi^{-1}(\bar{\xi})$, $\pi\colon \N(\*L)\lra \N(\*L)/\N(U_L)$ la
proyecci\'on natural de donde se sigue que $\bar{\xi}$ es abierto
en $\N(\*L)/\N(U_L)$. Por tanto $\N(\*L)/\N(U_L)$ es un
subgrupo discreto de $\*K/\N(U_L)$. Ahora bien, por ser
$\N(U_L)$ cerrado, se sigue que $\*K/\N(U_L)$ es Hausdorff
(se tiene en general que si $G$ es un grupo topol\'ogico
Hausdorff y $H$ es un subgrupo cerrado de $G$, entonces
$G/H$ es Hausdorff). Por el Lema \ref{discretocerrado} abajo,
se sigue que $\N(\*L)$ es un subgrupo cerrado de $\*K$.
$\fin$
\end{proof}

\begin{lema}\label{discretocerrado}
Sean $G$ un grupo topol\'ogico Hausdorff y $H<G$ un subgrupo
discreto de $G$. Entonces $H$ es cerrado.
\end{lema}

\begin{proof}
\noindent
{\bf Paso 1:} Sea $U$ una vecindad abierta tal que $U\cap H=\{e\}$,
entonces existe una vecindad abierta $V\subseteq U$ tal que
$VV^{-1}\subseteq U$ y $e\in V$.

En efecto, sea $\sigma:U\times U\lra G$, $\sigma(y_1,y_2)=y_1
y_1^{-1}$. Por continuidad, existe una vecindad $N\subseteq
U\times U$ de $(e,e)$ tal que $\sigma(N)=U$ y $N$ contiene una
vecindad abierta de la forma $V_1\times V_2$ con $V_1,V_2
\subseteq U$ y $e\in V_1\cap V_2$. Entonces $V:=V_1\cap
V_2$ es la vecindad buscada.

\noindent
{\bf Paso 2:} Sea $x\in G\setminus H$. Sea $U$ una vecindad abierta
de $e$ con $U\cap H=\{e\}$. El mapeo $L_x\colon G\to G$, $L_x(y)
=xy$ con $x\in G$ es un homeomorfismo con inverso $L_{x^{-1}}$.
Sea $V\subseteq U$ una vecindad de $e$ tal que $e\in V$ y $VV^{
-1}\subseteq U$. Sea $W=L_x(V)$. Se tiene que $W$ tiene a lo m\'as
un elemento de $H$. Si $W\cap H=\emptyset$, $W$ es la vecindad
buscada, Si $W\cap H=xV\cap H=\{h\}$, existen abiertos $U_x$ y
$U_h$ de $G$ con $x\in U_x$ y $h\in U_h$ y $U_x\cap U_h=\emptyset$.
Finalmente $xV\cap U_x$ es una vecindad abierta de $x$
disjunta de $H$.

Por tanto $H$ es cerrado en $G$.
$\fin$
\end{proof}

\subsection{Clasificaci\'on de campos locales. Extensiones no
ramificadas}\label{S17.3.3}

A continuaci\'on caracterizamos los campos locales.

\begin{proposicion}\label{CCCar0}
Sea $K$ un campo local de caracter\'istica $0$, esto es, ${\ma Q}
\subseteq K$. Sea $e:=v_K(p)$ donde $p=\ca (\o_K/\pK_K)=\ca
\tilde K$. Entonces $K$ es una extensi\'on finita de ${\ma Q}_p$,
los n\'umeros $p$-\'adicos y $[K:{\ma Q}_p]=ef$ donde $f=[\o_K/
\pK_K:{\ma F}_p]$. Finalmente, $\o_K$ es un ${\ma Z}_p$-m\'odulo
libre de rango $ef$: $\o_K\cong {\ma Z}_p^{ef}$ como 
${\ma Z}_p$-m\'odulo.
\end{proposicion}

\begin{proof}
Sean $v$ la valuaci\'on de $K$ y $w:=v|_{\ma Q}$. Se tiene que
$p=p\cdot 1_K\neq 0$ y $)<w(p)=v(p)=e<\infty$. Se tiene que $w$
es equivalente a la valuaci\'on $p$-\'adica $v_p$ de ${\ma Q}$. Sea
$L$ la cerradura de ${\ma Q}$ en $K$. Puesto que $K$ es completo,
$L$ es completo con $v|_L$ y es una completaci\'on de ${\ma Q}$.

Puesto que $w$ es equivalente a $v_p$, se sigue que $L\cong
{\ma Q}_p$. Ahora $v(p)=e(K|{\ma Q}_p)v_p(p)=e(K|{\ma Q}_p)=e$.
Se tiene que $\o_K/\pK_K={\ma F}_{p^f}$ con $f=[\o_K/\pK_K:
{\ma Z}_p/p{\ma Z}_p]=[\o_K/\pK_K:{\ma F}_p]$. Se sigue del
Teorema \ref{T17.3.2.4} que $[K:{\ma Q}_p]=ef$ y que $\o_K$ es
un ${\ma Z}_p$-m\'odulo libre de rango $ef$.
$\fin$
\end{proof}

\begin{lema}\label{L17.3.2.5}
Sea $K$ un campo local con $\tilde K\cong \F$, $q$ una potencia
de $p$. 
\las
\item Para cada $x\in\o_K$, el l\'imite
\begin{gather*}
\omega(x):=\lim_{n\to\infty}x^{q^n}
\intertext{existe en $\o_K$ y el mapeo $\omega$ satisface}
\omega(x)\equiv x\bmod \pK_K,\qquad \omega(x)^q=\omega(x),
\qquad \omega(xy)=\omega(x)\omega(y).
\end{gather*}
\item Sea $A=\{x\in K\mid x^{q-1}=1\}$, $R=A\cup\{0\}=
\{x\in K\mid x^q=x\}$. Entonces $R$ es un conjunto completo
de representantes de $\o_K/\pK_K$ en $\o_K$, conteniendo
a $0$. $A$ es el conjunto de todas las $(q-1)$-ra\'ices de la
unidad en $K$ y el homomorfismo can\'onico $\o_K\to \tilde K$
induce un isomorfismo de grupos multiplicativos: $A\cong
\*\F$.
\end{list}
\end{lema}

\begin{proof}
\las
\item
Por inducci\'on probaremos las congruencias $x^{q^n}\equiv
x^{q^{n-1}}\bmod \pK_K^n$ para $n\geq 1$. Para $n=1$,
$x^q\equiv x\bmod \pK_K$ se sigue del hecho de que 
$\tilde K$ es el campo finito de $q$ elementos. Supongamos
que la congruencia se cumple para $n\geq 1$ de tal manera
que $x^{q^n}=x^{q^{n-1}}+y$ con $y\in \pK_K^n$. El caso de
$\ca K=p$ se sigue inmediatamente tomando la $q$ potencia
de esta igualdad. Aqu\'i presentamos el caso general. Se tiene
\[
x^{q^{n+1}}=\sum_{i=0}^q \binom qi x^{iq^{n-1}}y^{q-i}.
\]
Para $0<i<q$, el entero $\binom qi=\frac qi \binom {q-1}{i-1}$ es
divisible por $p$ por lo que $\binom qi y^{q-i}\in \pK_K^{n+1}$.
Para $i=0$ lo anterior tambi\'en se cumple, por lo cual  obtenemos
lo deseado: $x^{q^{n+1}}\equiv x^{q^{n}}\bmod \pK_K^{n+1}$.

Por tanto $\{x^{q^{n}}\}_{n=1}^{\infty}$ es una sucesi\'on de
Cauchy por lo que es convergente a un elemento
$\omega (x)$ en $\o_K$ pues $\o_K$ es
cerrado y $K$ es completo. De las congruencias, obtenemos
que $x^{q^n}\equiv x\bmod \pK_K$ y por tanto que $\omega(x)
\equiv x\bmod \pK_K$. Tambi\'en obtenemos que
\begin{gather*}
\omega(x)^q=\lim_{n\to\infty}x^{q^{n+1}}=\omega(x),\qquad
\omega(xy)=\lim_{n\to\infty}(x^{q^n}y^{q^n})=\omega(x)
\omega(y).
\end{gather*}
\item Sea $V=\{\omega(x)\mid x\in\o_K\}$. Puesto que
$\omega(x)\equiv x\bmod \pK_K$, cada clase residual de
$\o_K\bmod \pK_K$ contiene al menos un elemento de $V$.
Por otro lado, puesto que $\omega(x)^q=\omega(x)$, $V$
es un subconjunto de $R$. Ahora bien, el n\'umero de 
elementos $x$ en $K$ que satisfacen $x^q-x=0$ es a lo
m\'as $q$ mientras que el n\'umero de elementos de 
$\tilde K$ es $q$. Se sigue que $V=R$ y $R$ es un
sistema completo de representantes de $\o_K/\pK_K$ en
$\o_K$. Se tiene $0=\omega(0)\in R$. Puesto que $\omega
(xy)=\omega(x)\omega(y)$, las propiedades enunciadas
para $A$ se siguen inmediatamente. $\fin$
\end{list}
\end{proof}

\begin{proposicion}\label{CCCp}
Sea $K$ un campo local de caracter\'istica $p>0$, es decir,
${\ma F}_p\subseteq K$. Entonces $(K,v_K)\cong (\F((T)), v_T)$
donde $v_T\big(\sum_{n=m}^{\infty}a_nT^n)=m$ si $a_m\neq
0$ para alg\'un $q=p^f$, $f\geq 1$. De hecho, $\F\cong
\o_K/\pK_K$ es el campo residual.
\end{proposicion}

\begin{proof}
Puesto que $K$ es un campo de caracter\'istica $p$, el
conjunto $R=\{x\in K\mid x^q=x\}$ en el Lema
\ref{L17.3.2.5} es un subcampo de $K$ de $q$ elementos. 
Por tanto $\F=R\subseteq K$. Sea $\pi=\pi_K$ un
elemento primo y representamos a $K$ como en el
Teorema \ref{T17.3.17.2} con $\pi_n=\pi^n$, $n\in{\ma Z}$.
Entonces el mapeo $\sum_na_n\pi^n\longmapsto \sum_n
a_nT^n$, $a_n\in\F$, define un $\F$-isomorfismo $(K,v_K)
\cong (\F((T)), v_T)$. $\fin$
\end{proof}

Resumimos en el siguiente resultado la caracterizaci\'on
de los campos locales.

\begin{teorema}\label{CamposLocales}
Sea $K$ un campo local. Entonces
\las
\item Si $\ca K=0$, $K/{\ma Q}_p$ es una extensi\'on finita
de grado $ef$, $e=v_K(p)=e(K|{\ma Q}_p)$ y $f=[\o_K/\pK_K:
{\ma F}_p]=f(K|{\ma Q}_p)$.
\item Si $\ca K=p>0$ y $\tilde K=\o_K/\pK_K\cong \F$, entonces
$(K,v_K)\cong (\F((T)), v_T)$. $\fin$
\end{list}
\end{teorema}

A continuaci\'on veremos que las extensiones no ramificadas
de campos locales est\'an en correspondencia con las
extensiones de los campos finitos.

\begin{teorema}\label{T17.3.3.1}
Sea $K$ un campo local con campo residual $\tilde K=\o_K/
\pK_K\cong \F$. Para cada $n\in{\ma N}$, existe una extensi\'on
$L/K$ no ramificada de grado $n$, $L$ es \'unico, es el
campo de descomposici\'on del polinomio $x^{q^n}-x$ sobre
$K$ y $L/K$ es una extensi\'on c\'iclica de grado $n$. Cada
elemento $\sigma\in\Gal(L/K)$ induce un automorfismo 
$\sigma'$ de $\tilde L/\tilde K$ y el mapeo $\sigma\to\sigma'$
induce un isomorfismo $\Gal(L/K)\cong \Gal(\tilde L/\tilde K)$.

M\'as a\'un, $L=K(A)$ donde $A=\{x\in \bar{K}\mid x^{
q^n}=x\}$.
\end{teorema}

\begin{proof}
El campo finito $\tilde K\cong \F$ tiene una extensi\'on de
grado $n$, la cual es \'unica y es c\'iclica. Por tanto existe
un polinomio m\'onico irreducible $g(x)\in \tilde K[x]$ de
grado $n$. Sea $f(x)\in\o_K[x]$ m\'onico tal que 
$\overline{f(x)}=g(x)=f(x)\bmod \pK_K$. Sea $\omega$
una ra\'iz de $f(x)$, $f(\omega)=0$. Sea $L:=K(\omega)$.

Sea $f(x)=x^n+\sum_{i=0}^{n-1}a_ix^i\in\o_K[x]$. Entonces
$\omega^n=-\sum_{i=0}^{n-1}a_i\omega^i$ con $a_i\in\o_K$.
Por tanto 
\begin{align*}
v_L(\omega^n)&=nv_L(\omega)=v_L\big(-\sum_{i=0}^{n-1}
a_i\omega^i\big)\geq \min_{0\leq i\leq n-1}\{v_L(a_i)+
iv_L(\omega)\} \\
&\geq\min_{0\leq i\leq n-1}\{iv_L(\omega)\}
=i_0v_L(\omega),
\end{align*}
para alg\'un $0\leq i_0\leq n-1$. Por tanto $(n-i_0)v_L(\omega)
\geq 0$ de donde se sigue que $v_L(\omega)\geq 0$ y
$\omega \in\o_L$. Sea $\alpha=\bar{\omega}\in \o_L/\pK_L$.
Por tanto $g(\alpha)=0=\overline{f(\omega)}$. Como $g(x)$
es irreducible, se sigue que $[\tilde K(\alpha):\tilde K]=\deg
g(x)=n$.

Por otro lado, $f(\omega)=0$ implica $[L:K]=[K(\omega):K]
\leq \deg f(x)=n$. Se sigue que 
\[
n=[\tilde K(\alpha):\tilde K]\leq [\tilde L:\tilde K]=f(L|K)\leq
[L:K]\leq n,
\]
de donde se obtiene $[L:K]=f(L|K)=n$. Se sigue que $L/K$
es una extensi\'on no ramificada de grado $n$.

Sea $E/K$ cualquier extensi\'on no ramificada de grado $n$.
Entonces $n=[E:K]=f(E|K)=[\tilde E:\tilde K]$ por lo que 
$\tilde E\cong {\ma F}_{q^n}$. Se tiene que $A=\{x\in\o_E\mid
x^{q^n}=x\}$ es un conjunto completo de representantes de
$\tilde E$ en $\o_E$ y $A$ tiene $q^n$ elementos, es decir,
$A$ es el conjunto de ra\'ices de $x^{q^n}-x\in\o_E[x]$.
Se sigue que $K(A)$ es el campo de descomposici\'on de
$x^{q^n}-x$ sobre $K$ y como este polinomio es separable,
la extensi\'on $F=K(A)/K$ es separable y por tanto $K(A)/K$
es Galois.

Se tiene que $\tilde K\subseteq \tilde F\subseteq \tilde E$ y puesto
que $A\subseteq F$, $\tilde E=\tilde F\cong {\ma F}_{q^n}$.
Adem\'as $\tilde F/\tilde K$ es una extensi\'on c\'iclica de
grado $n$. Cada $\sigma\in\Gal(F/K)$ induce un automorfismo
$\sigma'\in\Gal(\tilde F/\tilde K)$ y el mapeo $\sigma\mapsto
\sigma'$ define un homomorfismo $\Gal(F/K)\stackrel{
\varphi}{\lra}\Gal(\tilde F/\tilde K)$.

Si $\sigma'=1$, esto es, $\sigma\to 1$, $\sigma(A)=A$ pues
$A$ es el conjunto de ra\'ices de $x^{q^n}-x$ y puesto que
$A$ es un conjunto completo de representantes de $\tilde F$,
$\sigma'=1$ implica que $\sigma$ fija a cada elemento de $A$
y, puesto que $F=K(A)$, $\sigma =1$. Por lo tanto $\varphi$
es inyectiva y 
\[
[F:K]=|\Gal(F/K)|\leq |\Gal(\tilde F/\tilde K)|=[\tilde F:\tilde K]=
f(F|K)\leq [F:K],
\]
de donde se sigue que $[F:K]=f(F|K)=[\tilde F:\tilde K]=n$ y
$F/K$ es una extensi\'on no ramificada de grado $n$. Por
tanto $F=E$. Se sigue que $E=K(A)$ y por tanto $E$
es \'unico y $\Gal(E/K)\cong C_n$.
$\fin$
\end{proof}

Para finalizar esta secci\'on, presentamos algunos
resultados sobre la imagen de la norma de unidades
y $n$-unidades indicando algunas referencias para su
demostraci\'on m\'as adelante.

Sea $L/K$ una extensi\'on finita y separable de campos locales con
valuaciones $v_K$ y $v_L$ respectivamente. Se tiene la
sucesi\'on exacta $1\lra U_K\lra \* K\xrightarrow{v_K}
{\ma Z}\lra 0$ y el diagrama conmutativo
\[
\xymatrix{
1\ar[r]&U_K\ar[r]\ar[d]^i&\* K\ar[r]^{v_K}
\ar[d]_i&{\ma Z}\ar[r]\ar[d]^{\id}&0\\
1\ar[r]&U_L\ar[r]&\* L\ar[r]^{v_L}&{\ma Z}\ar[r]&0}
\]
donde las filas son exactas e $i$ denota al encaje.

\begin{proposicion}\label{CClaseP1.2.1.2}
Sea $\N_{L/K}\colon \* L\lra \* K$ la norma. Entonces se tiene
el siguiente diagrama conmutativo, donde las filas son exactas:
\[
\xymatrix{
1\ar[r]&U_L\ar@{^{(}->}[r]^{\iota}\ar[d]^{\N_{L/K}}_{\phi}&\* L\ar[r]^{v_L}
\ar[d]^{\N_{L/K}}_{\psi}&{\ma Z}\ar[r]\ar[d]^f&0\\
1\ar[r]&U_K\ar@{^{(}->}[r]^{\iota}&\* K\ar[r]^{v_K}&{\ma Z}\ar[r]&0}
\]
es decir, $fv_L(x)=v_K(\N_{L/K} x)$ donde $\N_{L/K}\pK_L=
\pK_K^f$, $\pK_L$ denota al ideal m\'aximo de $\o_L$ y donde
$f=[\o_L/\pK_L:\o_K/\pK_K]$.

En particular tenemos $v_K(\N_{L/K} y)=fv_L(y)$ para 
$y\in \*L$.
\end{proposicion}

\begin{proof}
Se tiene que $\N_{L/K}U_L\subseteq U_K$ y $\iota\circ\varphi=\psi\circ
\iota$ pues esto es simplemente $\iota\circ \N_{L/K}=\N_{L/K}\circ
\iota$. La igualdad $v_K\circ \N_{L/K}=f v_L$ es el contenido del
Teorema \ref{unicidaddevalorabsoluto} y del Corolario \ref{C17.3.2.7}.
$\fin$
\end{proof}

\begin{corolario}\label{C17.2.15'}
Si $L/K$ es una extensi\'on abeliana finita de campos locales,
se tiene que $e=[U_K:\N_{L/K}U_L]$.
\end{corolario}

\begin{proof}
En general, para una extensi\'on finita y separable de campos locales, tenemos,
aplicando el Lema de la Serpiente (ver Teorema
\ref{CClaseT1.5.2} m\'as adelante), al diagrama de la Proposici\'on
\ref{CClaseP1.2.1.2}, la sucesi\'on exacta
\begin{gather*}
\cdots\lra \ker f=\{0\}\lra\coker \phi\lra \coker \psi
\lra \coker f\lra 0,
\intertext{donde el mapeo $f$ es multiplicaci\'on por $f$. Por tanto}
0\lra \frac{U_K}{\N_{L/K}U_L}\lra \frac{\*K}{\N_{L/K}\*L}\lra
\frac{{\ma Z}}{f{\ma Z}}\lra 0.
\end{gather*}
Se sigue que $[\*K:\N_{L/K}\*L]=[U_K:\N_{L/K}U_L] f$.

En el caso particular en que la extensi\'on $L/K$ sea abeliana,
como veremos m\'as adelante (Teorema \ref{CClaseT3.2.2}), se tiene
$[\*K:\N_{L/K}\*L]=[L:K]=ef$ y por tanto $e=[U_K:\N_{L/K}U_L]$.
$\fin$
\end{proof}

\begin{proposicion}\label{CClaseP1.2.1.3}
Sea $L/K$ una extensi\'on finita no ramificada de campos locales
(y por tanto Galois). Entonces $\N:=\N_{L/K}\colon U_L^{(n)}\lra
\unidades n$ para toda $n$, esto es, $\N(U_L^{(n)})\subseteq
\unidades n$.
\end{proposicion}

\begin{proof}
Sea $x=1+y\in U_L^{(n)}$ con $y\in \pK_L^n$. Para $\sigma\in G=
\Gal(L/K)$ se tiene $\sigma x=1+\sigma y$ y por tanto
\[
\N x=\prod_{\sigma\in G}\sigma x=\prod_{\sigma\in G}
(1+\sigma y)\equiv \Big(1+\sum_{\sigma\in G}\sigma y\Big)
\bmod \pK_L^{2n}.
\]

Ahora bien, como $L/K$ es una extensi\'on no ramificada, $\pK_L^n
\cap U_K=\pK_K^n$ y $\sum_{\sigma\in G}\sigma y\in \pK_L^n\cap
U_K=\pK_K^n$ lo cual implica que $\N x\equiv 1\bmod \pK_K^n$. $\fin$
\end{proof}

Pasando a los cocientes en la Proposici\'on \ref{CClaseP1.2.1.3}, tenemos
\[
\N\colon U_L^{(n)}/U_L^{(n+1)}\lra \unidades n/\unidades {n+1}.
\]
Para $n=0$, $U_L^{(0)}/U_L^{(1)}=U_L/U_L^{(1)}\cong \*{\tilde{L}}$
y $\unidades 0/\unidades 1=U_K/\unidades 1\cong \*{\tilde{K}}$.

As\'i, $\N_0\colon \*{\tilde{L}}\lra \*{\tilde{K}}$ es la norma de campos
residuales. Para $n\geq 1$, $U_L^{(n)}/U_L^{(n+1)}\cong \pK_L^n/
\pK_L^{n+1}$ es un espacio vectorial de dimensi\'on $1$ sobre
$\tilde{L}$, esto es, $\pK_L^n/\pK_L^{n+1}\cong \tilde{L}$.

\begin{teorema}\label{CClaseT1.2.1.4} Sea $L/K$ una extensi\'on no
ramificada de campos locales (de Galois). Entonces
\las
\item $\N(U_L^{(n)})=\unidades n$ para toda $n\geq 1$.
\item $U_K/\N U_L\cong \* {\tilde{K}}/\N\*{\tilde{L}}$.
\item $\* K/\N \* L\cong {\ma Z}/f{\ma Z}\times \*{\tilde{K}}/\N
\*{\tilde{L}}$ donde $f=[L:K]=[\tilde{L}:\tilde{K}]$ (pues
$L/K$ es no ramificada).
\end{list}

En nuestro caso, $\tilde{K}$ y $\tilde{L}$ son campos finitos, por lo
que $\*{\tilde{K}}/\N \*{\tilde{L}}=\{1\}$ y $\*K/\N \* L\cong {\ma Z}/
f{\ma Z}$.
\end{teorema}

\begin{proof} (1) Ver Teorema \ref{CCLT17.6.3}.
Recordemos en general que si ${\eu D}_{L/K}$ denota al diferente
de la extensi\'on $L/K$, entonces 
\[
{\eu D}_{L/K}^{-1}=\{x\in L\mid \Tr(x\o_L)\subseteq \o_K\}.
\]
Puesto que $L/K$ es no ramificada, ${\eu D}_{L/K}^{-1}=(1)=\o_L$.
Se sigue que $\Tr\o_L=\o_K$.

De esta forma si $\pi\in\o_K$ con $v_K(\pi)=1$, se tiene que $v_L(
\pi)=1$ al ser $L/K$ no ramificada. Si $u\in U_L^{(n)}$, $u=1+\pi^n$
con $x\in\o_L$, 
\[
\N u=\prod_{\sigma\in G}(\sigma u)=\prod_{\sigma\in G}
\sigma(1+\pi^n x)=\prod_{\sigma\in G}(1+\pi^n\sigma x)=
1+\pi^n\Tr x+\pi^{n+1}\xi
\]
con $\xi\in \o_K$ pues $\N u\in \o_K$.

Veamos que $(\N U_L^{(n)}) \unidades{n+1}=\unidades n$.
Se tiene $(\N U_L^{(n)})\unidades {n+1}\subseteq \unidades n$.
Sea $a\in \unidades n$ arbitrario. Escribamos $a=1+\pi^n z$ con
$z\in\o_K$. Sea $x\in\o_L$ tal que $\Tr x=z$. Se tiene
\[
\N(1+\pi^n x)=1+\pi^n\Tr x+ \pi^{n+1}\xi=1+\pi^n z+\pi^{n+1}\xi
\]
con $\xi\in\o_K$ puesto que $\N(1+\pi^nx)\in\o_K$.

Queremos hallar $y\in\o_K$ tal que $(\N(1+\pi^n x))(1+\pi^{n+1}y)=
1+\pi^n z=a$. Se tiene
\begin{align*}
\N(1+\pi^n x)(1+\pi^{n+1}y)&=1+\pi^n z+\pi^{n+1}\xi +\pi^{n+1}y
+\pi^{2n+1}zy+\pi^{2n+2}y\xi\\
&= a+\pi^{n+1}(\xi+y(1+\pi^nz+\pi^{n+1}\xi)).
\end{align*}

Esto es, necesitamos resolver $\xi+y(1+\pi^nz+\pi^{n+1}\xi)=0$, es
decir, $yw=-\xi$ donde $w=1+\pi^nz+\pi^{n+1}\xi\in \unidades n
\subseteq U_K$. En este caso, $y=-\xi w^{-1}$ satisface lo deseado.

As\'i, dado $a\in\unidades n$, existen $x_1\in U_L^{(n)}$ y $z_1
\in \unidades {n+1}$ tales que $a=(\N x_1)z_1$. Ahora,
puesto que $z_1\in
\unidades {n+1}$, entonces existen $x_2\in U_L^{(n+1)}$ y
$z_2\in\unidades {n+2}$ tales que $z_1=(\N x_2)z_2$ por lo
cual $a=\N(x_1x_2)z_2$.

Continuando el proceso, tenemos que para toda $t\in{\ma N}$,
existen $x_1,\ldots, x_t\in\o_L$ y $z_t\in\unidades {n+t}$ tales
que $a=\N(x_1\cdots x_t)\cdot z_t$. Ahora $z_t\xrightarrow
[t\to 1]{}1$ por lo que $\lim_{t\to\infty}\N(x_1\cdots x_t)z_t=
\N\big(\lim_{t\to\infty}(x_1\cdots x_t)\big)=a$.

As\'i, $x_1\cdots x_t$ converge por ser $L$ un campo completo
y puesto que $x_i\in U_L^{(n+i)}\subseteq U_L^{(n)}$, el l\'imite
es un elemento de $U_L^{(n)}$. De esta forma, si $x_0=\lim_{
t\to \infty}(x_1\cdots x_t)$, $a=\N x_0$ con $x_0\in U_L^{(n)}$.
Se sigue que
\[
\N U_L^{(n)}=U_K^{(n)}.
\]

(2) Se tiene el siguiente diagrama conmutativo
\[
\xymatrix{
1\ar[r]&U_L^{(1)}\ar[r]\ar[d]^{\N_{L/K}}&U_L\ar[r]
\ar[d]^{\N_{L/K}}&\tilde{L}^{\ast}\ar[r]\ar[d]^{\N_{L/K}}&0\\
1\ar[r]&\unidades 1\ar[r]\ar[d]&U_K\ar[r]\ar[d]
&\tilde{K}^{\ast}\ar[r]\ar[d]&0\\ 
&0&U_K/\N U_L&\*{\tilde{K}}/
\N\*{\tilde{L}}}
\]
Por el Lema de la Serpiente (ver Teorema
\ref{CClaseT1.5.2} m\'as adelante) se obtiene (2).

(3) Sea $\pi\in\o_K$, $v_K(\pi)=v_L(\pi)=1$. Se tiene $K^{\ast}
\cong (\pi)\times U_K$, $\L\cong (\pi)\times U_L$. Por tanto
$\N \L=(\pi^f)\times \N U_L$. Se sigue
\begin{gather*}
\frac{K^{\ast}}{\N\L}\cong \frac{(\pi)}{(\pi^f)}\times \frac{U_K}{
\N U_L}\cong \frac{{\ma Z}}{f{\ma Z}}\times 
\frac{\*{\tilde{K}}}{\N \*{\tilde{L}}}. \tag*{\fin}
\end{gather*}
\end{proof}

\begin{corolario}\label{CClaseC1.2.1.5} Sea $L/K$ una
extensi\'on finita no ramificada de campos locales. 
Entonces las siguientes condiciones son equivalentes
\las
\item $[K^{\ast}:\N\L]=f=[\tilde{L}:\tilde{K}]$.
\item $U_K=\N U_L$.
\item $\*{\tilde{K}}=\N\*{\tilde{L}}$. $\fin$
\end{list}
\end{corolario}

\begin{observacion}\label{CClaseO1.2.1.6}
Si $L/K$ es una extensi\'on no ramificada no necesariamente Galois,
la Proposici\'on \ref{CClaseP1.2.1.3} y el Teorema \ref{CClaseT1.2.1.4} siguen
siendo v\'alidos.
\end{observacion}

Esto es aplicable para cuando $L/K$ es una extensi\'on de campos 
completos con respecto a valores absolutos no arquimedianos con
campos residuales no necesariamente finitos. Cuando $L$ y $K$
son campos locales, si la extensi\'on es finita no ramificada,
necesariamente $L/K$ es una extensi\'on c\'iclica (ver el Teorema 
{\rm{\ref{T17.3.3.1}}}).

\section{Cohomolog{\'\i}a de grupos finitos}\label{CClaseS1.5}

En esta primera parte damos una serie de propiedades generales de 
m\'odulos que detallaremos m\'as adelante. Los resultados b\'asicos
de cohomolog{\'\i}a se pueden consultar en
\cite{Ser64,Ser,Vil2006}.

\begin{definicion}\label{CoD17.5.1}
Sea $G$ un grupo finito con identidad $1$ y sea $A$ un grupo abeliano.
Entonces $A$ se llama {\em $G$-m\'odulo\index{$G$-m\'odulo}} 
(izquierdo) si $G$
act\'ua en $A$ ($G\times A\lra A$, $(\sigma,a)\longmapsto \sigma a$)
tal que para todos $\sigma, \tau\in G$
y para todos $a,b\in A$, se tiene
\las
\item $1\cdot a=a$
\item $\sigma(a+b)=\sigma a+\sigma b$
\item $(\sigma\tau)(a)=\sigma(\tau a)$
\end{list}
\end{definicion}

Sea $\zg=\{\suma \sigma G n\mid n_{\sigma}\in{\ma Z}\}$. Entonces
$\zg$ es un grupo abeliano libre en $|G|$ generadores. La multiplicacion
definida por
\[
\Big(\suma \sigma G n\Big)\Big(\suma \sigma G m\Big)=
\sum_{\tau\in G}\Big(\sum_{\varepsilon\delta=\tau}n_{\varepsilon}
m_{\delta}\Big)\tau,
\]
hace de $\zg$ un anillo y $\zg$ se llama el 
{\em anillo grupo $\zg$\index{anillo grupo}}. Notemos que los conceptos
de $G$-m\'odulo y $\zg$-m\'odulo coinciden:
\[
\Big(\suma \sigma G n\Big)a=\sum_{\sigma\in G} n_{\sigma}(\sigma a),
\qquad a\in A.
\]
Se tiene que $\zg$, como grupo aditivo, es un $G$-m\'odulo.

Sea $\N_G=\sum_{\sigma\in G}\sigma$, es decir, $n_{\sigma}=1$
para toda $\sigma \in G$.
Entonces ${\ma Z}\N_G=\{\sum_{\sigma\in G}n\sigma\mid n\in{\ma Z}\}$
es un ideal de $\zg$ y $I_G:=\{\suma \sigma G n\mid \sum_{\sigma\in
G}n_{\sigma}=0\}$ es tambi\'en un
ideal llamado el {\em ideal aumentaci\'on\index{ideal
aumentaci\'on}\index{aumentaci\'on!ideal}}.

Se tiene que $I_G=\ker \varphi$, donde $\varphi\colon \zg\lra {\ma Z}$,
$\varphi\big(\suma \sigma G n\big)=\sum_{\sigma\in G}n_{\sigma}$ y la
sucesi\'on $0\lra I_G\lra \zg \stackrel{\varphi}{\lra} {\ma Z}\ (\cong {\ma Z}\N_G)
\lra 0$ es exacta. $\N_G$ se llama {\em norma\index{norma}} o 
{\em traza\index{traza}} de $\zg$.

El mapeo $\mu\colon {\ma Z}\lra \zg$, $\mu(n):=n\N_G$ se llama 
{\em coaumentaci\'on\index{coaumentaci\'on}}. Sea $J_G:=\zg/{\ma Z}
\N_G$. Entonces se tienen las sucesiones exactas de $G$-m\'odulos
\begin{gather*}
0\lra I_G\lra \zg\stackrel{\varphi}{\lra}{\ma Z}\lra 0,\\
0\lra {\ma Z}\stackrel{\mu}{\lra} \zg\lra J_G\lra 0.
\end{gather*}

\begin{definicion}\label{CoD17.5.2}
Sea $A$ un $G$-m\'odulo. Se tienen los siguientes subm\'odulos:
\las
\item $A^G=\{a\in A\mid \sigma a \text{\ para toda\ }a\in A\}$ el {\em m\'odulo
fijo\index{modulo fijo@m\'odulo fijo}} de $A$.
\item $\N_GA=\{\N_G a\mid a\in A\}=\{\sum_{\sigma\in G}\sigma a\mid
a\in A\}$ el {\em grupo norma\index{grupo norma}} de $A$.
\item ${_G}A={_{\N_G}}A=\{a\in A\mid \N_Ga=\sum_{\sigma \in G}
\sigma a=0\}$ el n\'ucleo de la norma.
\item $I_G A=\{\sum_{\sigma\in G}n_{\sigma}(\sigma a_{\sigma}-
a_{\sigma})\mid a_{\sigma}\in A\}=\langle (\sigma -1)A\mid
\sigma\in G\rangle$.
\end{list}
\end{definicion}

Se verifica directamente que $\N_G A\subseteq A^G$ y $I_G A
\subseteq {_G}A$. Formamos los cocientes $A^G/\N_G A$ y ${_G}A/I_G A$
los cuales son los grupos de cohomolog\'ia
(de Tate) de $A$ en dimensiones $0$
y $-1$ respectivamente.

Sea $A$ un $G$-m\'odulo y $H$ un subgrupo de $G$. Entonces $A$
es un $H$-m\'odulo. Si adem\'as $H\normal G$, $A^H$ es un $G/H$-m\'odulo.

\begin{definicion}\label{Co17.5.3}
Sean $A$ y $B$ dos $G$-m\'odulos. Entendemos
por un $G$-homomorfismo $f\colon A\lra B$
a un homomorfismo de grupos abelianos tal que $f(\sigma a)=\sigma
f(a)$ para todas $a\in A$ y $\sigma \in G$.
\end{definicion}

Dados dos $G$-m\'odulos $A$ y $B$, $\Hom(A,B)$ denota el grupo de  los
${\ma Z}$-homomorfismos (es decir, homomorfismos de grupos abelianos)
de $A$ en $B$. El grupo $\Hom(A,B)$ se puede hacer un $G$-m\'odulo
as\'i: Si $f\in\Hom(A,B)$ y $\sigma \in G$, se define la acci\'on de $\sigma$
en $f$ por:
\begin{gather*}
\sigma \circ f=\sigma f=\sigma\circ f\circ \sigma^{-1},
\intertext{es decir,}
(\sigma f)(a)=\sigma f(\sigma^{-1}a),
\end{gather*}
para $a\in A$.
Se tiene que el grupo de $G$-homomorfismos $\Hom_G(A,B)$ de $A$ en
$B$ es un $G$-subm\'odulo de $\Hom(A,B)$ y de hecho se tiene
\[
(\Hom(A,B))^G=\Hom_G(A,B)\quad \text{(Proposici\'on \ref{CoP17.5.14})}.
\]

Para dos grupos abelianos $A$ y $B$, el producto tensorial de $A$
y $B$ sobre ${\ma Z}$ se denota simplemente $A\otimes B$.
Esto es, $A\otimes B=A\otimes_{\ma Z} B$. Se tiene que
\begin{gather*}
{\ma Z}\otimes A\cong A, \quad (x\otimes a)\mapsto xa,\\
A\otimes B\cong B\otimes A.
\end{gather*}

Las siguientes propiedades son de verificaci\'on directa.

\begin{proposicion}\label{CoP17.5.4}
Sea $\{X_i\mid i\in I\}$ una familia de $G$-m\'odulos y sea $Y$
otro $G$-m\'odulo. Entonces se tienen los siguientes isomorfismos
\las
\item $Y\otimes\Big(\bigoplus_{i\in I} X_i\Big)\cong \bigoplus_{i\in I}(Y\otimes X_i)$.
\item $\Hom\big(\bigoplus_{i\in I}X_i,Y\big)\cong \prod_{i\in I}\Hom_G(X_i,Y)$.
\item $\Hom_G\big(Y,\prod_{i\in I}X_i\big)\cong \prod_{i\in I}\Hom_G(Y,X_i)$.

Si adem\'as $Y$ es finitamente generado como grupo abeliano, se tiene
\item $Y\otimes \big(\prod_{i\in I} X_i\big)\cong \prod_{i\in I} (Y\otimes X_i)$.
\item $\Hom_G\big(Y,\bigoplus_{i\in I} X_i\big)\cong \bigoplus_{i\in I}\Hom_G
(Y,X_i)$. $\fin$
\end{list}
\end{proposicion}

Sean $A,B, C, D$ cuatro $G$-m\'odulos y $h\colon A\to C$ un $G$-homomorfismo.
Entonces se tiene homomorfismos
\begin{gather*}
 \*h\colon \Hom(C,B)\lra \Hom(A,B),\quad f\longmapsto \*h(f)=f\circ h
\intertext{y}
h'\colon A\otimes B\lra C\otimes B, \quad a\otimes b\longmapsto h(a)\otimes b.
\intertext{Sea $g\colon B\lra D$ un $G$-homomorfismo, se tienen los homomorfismos}
g_{\star}\colon \Hom(A,B)\lra\Hom (A,D), \quad f\longmapsto g_{\star}(f)=g\circ f 
\intertext{y}
g'\colon A\otimes B\lra A\otimes D, \quad a\otimes b\longmapsto a\otimes g(b).
\intertext{Tambi\'en se tiene}
h\otimes g\colon A\otimes B\lra C\otimes D, \quad (h\otimes g)
(a\otimes b)=h(a)\otimes g(b).
\intertext{Finalmente, si $H\colon C\lra A$ y $G\colon B
\lra D$ son $G$-homomorfismos, se tiene}
(H,G)\colon \Hom(A,B)\lra\Hom
(C,D), \quad f\longmapsto G\circ f\circ H.
\end{gather*}

\begin{definicion}\label{CoD17.5.5}
Un $G$-m\'odulo $X$ se llama {\em proyectivo\index{modulo proyectivo@m\'odulo
proyectivo}\index{proyectivo!m\'odulo $\sim$}} si satisface cualquiera de las
siguientes condiciones equivalentes
\las
\item Todo diagrama 
$\xymatrix{X\ar@{->}[dr]^{\varphi}\\
B\ar@{->}[r]_{\mu}&C\ar@{->}[r]&0}$ de $G$-modulos
con $\mu$ suprayectiva, se puede extender 
$\psi\colon X\lra B$ tal que
$\xymatrix{X\ar@{->}[dr]^{\varphi}\ar@{-->}[d]_{\psi}\\
B\ar@{->}[r]_{\mu}&C\ar@{->}[r]&0}$ es conmutativo, esto es,
$\mu\circ \psi=\varphi$.

\item $X$ es sumando directo de un $G$-m\'odulo libre.

\item Cualquier sucesi\'on exacta $0\lra A\lra B\lra X\lra 0$ se escinde
y $B\cong A\oplus X$.
\end{list}
\end{definicion}

\begin{proposicion}\label{CoP17.5.6}
Si $X$ es un $G$-m\'odulo y $0\lra A\stackrel{h}{\lra}B
\stackrel{g}{\lra}C\lra 0$ es una sucesi\'on exacta de $G$-m\'odulos,
entonces la sucesi\'on inducida
\[
0\lra\Hom_G(X,A)\stackrel{\*h}{\lra}\Hom_G(X,B)\stackrel{\*g}{\lra}
\Hom_G(X,C)
\]
es exacta. Si $X$ es proyectivo, entonces $\*g$ es suprayectiva.
\end{proposicion}

\begin{proof}
Se deja al cuidado del lector. $\fin$
\end{proof}

\begin{proposicion}\label{CoP17.5.7}
\las
\item Si $\cdots \lra X_{q+1}\stackrel{d_{q+1}}{\lra}X_q\stackrel{
d_q}{\lra}X_{q-1}\lra\cdots $ es una sucesi\'on exacta de ${\ma Z}$-m\'odulos
libres y si $D$ es cualquier ${\ma Z}$ m\'odulo, entonces la sucesi\'on
\[
\cdots\lra \Hom(X_{q-1},D)\stackrel{\*{d_q}}{\lra}\Hom(X_q,D)\xrightarrow[]
{\*{d_{q+1}}}\Hom(X_{q+1},D)\lra \cdots
\]
es exacta.

\item Si $0\lra A\lra B\lra C\lra 0$ es una sucesi\'on exacta de 
${\ma Z}$-m\'odulos libres y si $X$ es cualquier ${\ma Z}$-m\'odulo, entonces
la siguiente sucesi\'on es exacta:
\[
0\lra A\otimes X\lra B\otimes X\lra C\otimes X\lra 0.
\]

\item Si $0\lra A\lra B\lra C\lra 0$ es una sucesi\'on exacta de 
${\ma Z}$-m\'odulos y si $X$ es un ${\ma Z}$-m\'odulo libre, entonces
la siguiente sucesi\'on es exacta:
\[
0\lra A\otimes X\lra B\otimes X\lra C\otimes X\lra 0.
\]
\end{list}
\end{proposicion}

\begin{proof}
Se deja al cuidado del lector.
$\fin$
\end{proof}

Uno de los resultados centrales en cohomolog{\'\i}a de grupos es

\begin{teorema}[Lema de la serpiente\index{lema
de la serpiente}]\label{CClaseT1.5.2}
Sea
\[
\begin{CD}
&&A @>f>> B @ >g>>C @>>> 0\\
&&@VV{\alpha}V @VV{\beta}V @VV{\gamma}V\\
0 @>>> A' @>{f'}>> B' @>{g'}>> C'\\
\end{CD}
\]
un diagrama conmutativo de 
$G$--m\'odulos, en donde las filas son exactas. Entonces
existe un {\em homomorfismo de conexi\'on\index{homomorfismo
de conexi\'on}}
$\delta: \ker \gamma
\longrightarrow \coker \alpha$ tal que
\[
\ker \alpha \stackrel{ \tilde{ f } }{
\longrightarrow } \ker \beta \stackrel{
\tilde{ g } }{ \longrightarrow } \ker
\gamma \stackrel{ \delta }{ \longrightarrow }
\coker \alpha \stackrel{ \tilde{ f' } }{
\longrightarrow } \coker \beta \stackrel{
\tilde{ g' } }{ \longrightarrow } \coker
\gamma
\]
es una sucesi\'on exacta, donde 
$\tilde{f'}$ y $\tilde{g'}$ son los homomorfismos inducidos
de $f'$ y $g'$ respectivamente y $\tilde{f}$ y $\tilde{g}$ son las
restricci\'on de  $f$ y $g$ respectivamente.

Adem\'as, si $f$ es inyectiva, entonces $\tilde{f}$ es
inyectiva y si $g'$ es suprayectiva, $\tilde{g'}$ es
suprayectiva.

Formalmente, tenemos $\delta=(f')^{-1}\circ \beta\circ g^{-1}$.
\end{teorema}

\begin{proof} Los detalles pueden ser consultados en
\cite[Theorem A.1.16]{Vil2006}. 
El mapeo de conexi\'on $\delta$ est\'a dado de la siguiente
forma. Sea $z\in \ker \gamma$, entonces $z=g(y)$ para
alg\'n $y\in B$. Por tanto $0=\gamma(z)=(\gamma\circ g)(y)=
(g'\circ \beta)(y)$. Por tanto $\beta(y)\in\ker g'=\im f'$. Sea
$\beta(y)=f'(a)$ para alg\'un $a\in A'$. Sea $\delta(z):=a
+\im \alpha$. El mapeo de conexi\'on se puede
representar por $\delta:=(f')^{-1}\circ \beta\circ
g^{-1}$. 
\[
\xymatrix{
& B\ar@{->}[d]_{\beta}\ar@{<-}[r]_{g^{-1}}&C\ar@{->}@/_2pc/[dll]_{\delta}\\
A'\ar@{<-}[r]^{(f')^{-1}}&B'
}
\]
Se verifica que $\delta$ est\'a bien definida y que la
sucesi\'on es exacta.
 $\fin$
\end{proof}

\begin{definicion}\label{CoD17.5.8}
Una {\em resoluci\'on proyectiva\index{resolucion proyectiva@resoluci\'on
proyectiva}} $P$ de ${\ma Z}$ es una sucesi\'on exacta de $G$-m\'odulos de
la forma
\[
 \qquad \cdots\lra P_n\too n{} P_{n-1}{}\too{{n-1}}{}\cdots\lra P_1\too 1{} P_0
\too 0{} {\ma Z}\lra 0,\leqno{P:}
\]
donde ${\ma Z}$ se considera el $G$-m\'odulo trivial: $g\circ x=x$ para todas
$g\in G$ y $x\in{\ma Z}$ y los m\'odulos $P_n$ son proyectivos para toda
$n\in{\ma N}\cup\{0\}$. En particular $\partial_n\circ\partial_{n+1}=0$ para
toda $n\geq 0$.
\end{definicion}

Dada una resoluci\'on proyectiva de ${\ma Z}$ y dado un $G$-m\'odulo
$A$ se tienen las sucesiones
\begin{gather*}
0\lra\Hom_G(P_0,A)\too 1* \Hom_G(P_1,A)\too 2*\cdots \lra\\
\lra\Hom_G(P_{n-1},A)\too n*\Hom_G(P_n,A)\lra\cdots
\intertext{y}
\cdots\lra P_n\otimes_G A\too n+P_{n-1}\otimes_G A\lra\cdots
\lra P_1\otimes_G A\too 1+ P_0\otimes_G A\lra 0,
\end{gather*}
donde $\partial^*_n(\varphi)=\varphi\circ \partial_n$, $\partial^+_n(
x\otimes a)=\partial_n x\otimes a$ y donde $\partial^*_0
=0$, $\partial_0^+=0$.

\begin{definicion}\label{CoD17.5.9*}
Para $n=0,1,\ldots, $ el {\em $n$-\'esimo grupo de 
cohomolog\'ia\index{grupo de cohomolog\'ia}} del $G$-m\'odulo
$A$ con respecto a la resoluci\'on proyectiva $P:$ se define como el grupo:
\begin{gather*}
H^n(P,A):=\ker \partial_{n+1}^*/\im \partial_n^*
\intertext{y el {\em $n$-\'esimo grupo de homolog\'ia\index{grupo de homolog\'ia}}
como}
H_n(P,A):=\ker \partial_n^+/\im\partial_{n+1}^+
\end{gather*}
\end{definicion}

Uno de los resultados centrales para la cohomolog\'ia de grupos es que
los grupos de cohomolog\'ia y de homolog\'ia no dependen de la resoluci\'on
dada.

\begin{teorema}\label{unicidadcohomologia}
Si $P$ y $P'$ son dos resoluciones proyectivas de ${\ma Z}$, se tiene
$H^n(P,A)\cong H^n(P',A)$ $H_n(P,A)\cong H_n(P',A)$ para toda $n=
0,1,\ldots$ y para todo $G$-m\'odulo $A$.
\end{teorema}

\begin{proof}
Ver \cite[Theorem A.1.19]{Vil2006}.
$\fin$
\end{proof}

\begin{definicion}\label{CoD17.5.9}
Para un $G$-m\'odulo arbitrario $A$ y para $n=0,1,\ldots $ se definen los
{\em grupos de cohomolog\'ia\index{grupos de cohomolog\'ia}}
$\co nGA$ como $\co nPA$ y los {\em grupos de 
homolog\'ia\index{grupos de homolog\'ia}} $\co nGA$ como $\co nPA$
donde $P$ es cualquier resoluci\'on proyectiva de ${\ma Z}$.
\end{definicion}

Gracias al Teorema \ref{unicidadcohomologia}, para tener los grupos de
homolog\'ia y de cohomolog\'ia, basta hallar una resoluci\'on proyectiva
de ${\ma Z}$. Sea $G^{n+1}=\underbrace{G\times\cdots\times G}_{
n+1}$, $P_n:={\ma Z}[G^{n+1}]$ el anillo grupo. $G$ act\'ua en $P_n$
as\'i: $x\circ (g_0,\ldots,g_n)=(xg_0,\ldots,xg_n)$, $x\in G$, $(g_0,\ldots,
g_n)\in G^{n+1}$. Entonces $P_n$ es un ${\ma Z}$-m\'odulo libre con
base $\{(g_0,\ldots, g_n)\mid g_i\in G\}$. Adem\'as $P_n$ es un 
$\zg$-m\'odulo libre con base $\{(1,x_1,\ldots,x_n)\mid x_i\in G\}$.

Sea $\partial_n\colon P_n\lra P_{n-1}$ dada por
\[
\partial_n(g_0,g_1,\ldots,g_n):=\sum_{i=0}^n (-1)^i (g_0,g_1,\ldots,
\hat{g_i},\ldots, g_n),
\]
donde $\hat{g_i}$ significa que el elemento $g_i$ ha sido removido.
Ahora $\partial_0\colon P_0=\zg\lra{\ma Z}$ es la aumentaci\'on
$\partial_0(g)=1$ para toda $g\in G$. Se tiene que la sucesi\'on
\begin{gather*}
\cdots\lra P_n\too n{} P_{n-1}\lra\cdots\lra P_1\too 1{} P_0\too 0{}
{\ma Z}\lra 0
\end{gather*}
es $G$-exacta y se llama la {\em resoluci\'on 
can\'onica\index{resolucion canonica@resoluci\'on can\'onica}}.

Esto prueba la existencia tanto de los grupos de homolog\'ia como
los de cohomolog\'ia para cualquier $G$-m\'odulo $A$.

Sean $A$ y $B$ dos $G$-m\'odulos y sea $f\colon A\to B$ un
$G$-homomorfismo. Se definen $H^n(f)\colon \co nGA\to
\co nGB$ y $H_n(f)\colon \ho nGA\to \ho nGB$, $n\geq 0$, de
la siguiente forma: Sea
\[
\cdots\lra P_n\too n{} P_{n-1}\lra\cdots\lra P_1\too 1{} P_0\too 0{}
{\ma Z}\lra 0\leqno{P:}
\]
cualquier resoluci\'on proyectiva. Denotemos $P_i\otimes_G A=
P_i\otimes_{\zg} A$ y sea
\[
f_n\colon P_n\otimes_G A\lra P_n\otimes_G B,
\quad f_n(x\otimes_G a)=x\otimes_G f_n(a)=(1_{P_n}\otimes f)(x
\otimes_G a).
\]
Entonces $f$ induce naturalmente el mapeo $H_n(f)\colon \ho nGA
\lra\ho nGB$.

Ahora consideremos la sucesi\'on
\begin{gather*}
0\lra\Hom_G(P_0,A)\too 1*\Hom_G(P_1,A)\too 2*\cdots\lra\\
\lra\Hom_G(P_{n-1},A)\too n*\Hom_G(P_n,A)\lra\cdots
\end{gather*}
Sea $\*{f_n}\colon \Hom_G(P_n,A)\lra \Hom_G(P_n,B)$ dado por
$\*{f_n}(\varphi)=f_n\circ \varphi$. Se tiene que $\*{f_n}$ induce
naturalmente $H^n(f)\colon \co nGA\lra \co nGB$, $n=0,1,2,\ldots$.

Como consecuencia del lema de la serpiente, se puede demostrar

\begin{teorema}\label{CClaseT1.5.3}
Si $0\to A\xrightarrow{f} B\xrightarrow{g}C\to 0$ es una sucesi\'on exacta
de $G$--m\'odulos, entonces existen homomorfismos de grupos 
$\varepsilon_n\colon \ho {{n+1}}GC\lra \ho nGA$ y $\delta_n\colon
\co nGC\lra \co {{n+1}}GA$, $n=0,1,\ldots$ llamados de 
{\em conexi\'on\index{mapeos de conexion@mapeos de conexi\'on}}, tales
que las siguientes sucesiones infinitas de homolog\'ia y de cohomolog\'ia
\begin{gather*}
\cdots\lra \ho {{n+1}}GB\xrightarrow{H_{n+1}(g)}{}\ho {{n+1}}GC
\xrightarrow{\varepsilon_n}{}\ho nGA\xrightarrow{H_n(f)}{}\ho nGB
\lra\\
\cdots\xrightarrow{\varepsilon_0}{}\ho 0GA\xrightarrow{H_0(f)}{}
\ho 0GB\xrightarrow{H_0(g)}\ho 0GC\lra 0
\intertext{y}
0\lra \co 0GA\xrightarrow{H^0(f)}{}
\co 0GB\xrightarrow{H^0(g)}\ho 0GC \xrightarrow{\delta_0}{}\co1GA
\lra\cdots \\
\xrightarrow{\delta_{n-1}}{}\co nGA\xrightarrow{H^n(f)}{}
\co nGB \xrightarrow{H^n(g)}{}\co nGC\xrightarrow{\delta_n}{}
\co {{n+1}}GA\lra\cdots
\end{gather*}
son sucesioes exactas de grupos.

Los mapeos de conexi\'on est\'an dados por el Lema de la Serpiente, 
a saber: si $\{P_n\}_n$ es una resoluci\'on proyectiva, $\partial_n
\colon P_n\lra P_{n-1}$, $X_n=P_n\otimes_G X$, $X^n=\Hom_G(
P_n,X)$, $X\in\{A,B,C\}$, $f_n\colon A_n\lra B_n$, $f_n(x\otimes a)=
x\otimes f(a)$, $g_n \colon B_n
\lra C_n$, $g_n(y\otimes b)=y\otimes g(b)$; 
$f^n\colon A^n\lra B^n$, $f^n(\mu)=f\circ \mu$,
$g^n \colon B^n\lra C^n$, $g^n(\psi)=g\circ \psi$;
$\partial^+_{n,X}\colon X_n\lra X_{n-1}$, $\partial^+_{n,X}(\alpha
\otimes x)=\partial_{n,X}(\alpha)\otimes x$ y $\partial^*_{n,X}\colon
X^{n-1}\lra X_n$, $\partial^*_{n,X}(\mu)=\mu\circ \partial_{n,X}$,
entonces, formalmente se definen:
\[
\varepsilon_n=f_n^{-1}\circ \partial_{n+1,B}^+\circ g_{n+1}^{-1}
\qquad \text{y}\qquad
\delta_n=(f^{n+1})^{-1}\circ \partial^*_{n+1,B}\circ (g^n)^{-1}.
\]
\end{teorema}
\begin{proof} \cite[Theorem A.3.6]{Vil2006}. $\fin$
\end{proof}

\subsection{Homolog\'ia y cohomolog\'ia en bajas
dimensiones}\label{SBajasDimensiones}

Sea $A$ un $G$-m\'odulo y $P$ la resoluci\'on proyectiva can\'onica.
Se tiene que $P_0\otimes A\cong A$. Por tanto
\[
\ho 0GA=(P_0\otimes A)/\im(\partial_1\otimes 1).
\]
Ahora $(\partial_1\otimes 1)((g_1,g_2)\otimes a)=g_1a-g_2a$ por lo que
$\im(\partial_1\otimes 1)=\langle a-ga\mid g\in G, a\in A\rangle=D A
\subseteq A$. Por tanto
\[
\ho 0GA=A_G=A/D A.
\]

Sea $I_G=\big\{\suma \sigma G a\mid \sum_{\sigma\in G}a_{\sigma}=0\big\}$
el ideal aumentaci\'on de $\zg$. Se tiene $\zg/I_G\cong {\ma Z}$. Si
$\suma \sigma G a\in I_G$, $a_1=-\sum_{\sigma\neq 1}a_{\sigma}$. Por
tanto 
\begin{align*}
\suma \sigma G a&=a_1\cdot 1+\sum_{\sigma\neq 1}a_{\sigma}\sigma=
\big(-\sum_{\sigma\neq 1}a_{\sigma}\big)\cdot 1+\sum_{\sigma\neq 1}
a_{\sigma}\sigma\\
&=\sum_{\sigma\in G}a_{\sigma}(\sigma-1)\in\langle \sigma
-1\mid\sigma\in G\rangle.
\end{align*}

Rec\'iprocamente, $\sigma-1\in I_G$ para toda $\sigma\in G$. Por tanto
$DA=\langle a-\sigma a\mid \sigma\in G, a\in A\rangle =I_G A$ y
\[
\ho 0GA=A/I_G A.
\]

Sea ahora $P$ la resoluci\'on can\'onica de $A$. Se tiene que la
sucesi\'on
\[
0\lra\Hom_G({\ma Z},A)\too 0*\Hom_G(P_0,A)\too 1*\Hom_G(P_1.A)
\]
es exacta tanto en $\Hom_G({\ma Z},A)$ como en $\Hom_G(P_0,A)$ de
donde obtenemos que
\[
\co 0GA=\ker \partial^*_1=\im \partial^*_0=\Hom_G({\ma Z},A)\cong (\Hom(
{\ma Z},A))^G\cong A^G.
\]

Para calcular $\co 1GA$ consideremos $P$ la resoluci\'on can\'onica
\begin{gather*}
\cdots\lra P_n\too n{} P_{n-1}\lra\cdots\lra P_1\too 1{} P_0\too 0{}
{\ma Z}\lra 0,
\intertext{donde $P_n={\ma Z}[G^{n+1}]$. Sea $K_n:=\Hom_G(P_n,A)$.
Se tiene}
0\lra K_0\too 0*K_1\lra\cdots\too n*K_n\too {{n+1}}*\cdots
\end{gather*}

Entonces 
\begin{align*}
\co nGA&=\ker\partial^*_{n+1}/\im\partial^*_n=Z^n(G,A)/B^n(G,A)\\
&= \text {$n$-cociclos}/\text {$n$-cofronteras}.
\end{align*}
En particular se tiene
$\co 1GA=\ker\partial^*_2/\im \partial^*_1=Z^1(G,A)/B^1(G,A)$ donde
\begin{gather*}
Z^1(G,A)=\{f\colon G\lra A\mid f(gh)=gf(h)+f(g), g, h\in G\}\\
=\text{grupo de los homomorfismos cruzados de $G$ en $A=1$-cociclos} 
\intertext{y}
B^1(G,A)=\{f\colon G\lra A\mid \text{existe $a\in A$ con $f(g)=ga-a$ para
toda $g\in G$}\}\\
=\text{$1$-cofronteras}.
\intertext{En otras palabras}
\co 1GA=\frac{\{f(gh)=gf(h)+f(g)\}}{\{f(g)=ga-a\}}.
\end{gather*}

Por ejemlo, si $G$ act\'ua de manera trivial en $A$, entonces
$Z^1(G,A)=\Hom(G,A)$ y $B^1(G,A)=0$, por lo que $\co 1GA=\Hom(G,A)$.
En particular, si $G$ es un grupo finito y ${\ma Z}$ tiene acci\'on
$G$-trivial, entonces $\co 1G{\ma Z}=\Hom(G,{\ma Z})=\{0\}$.

Un ejemplo importante es el grupo
$H^2(G,A)$ el cual est\'a en correspondencia biyectiva
con las clases de equivalencia de extensiones de $G$ por $A$.
M\'as precisamente, un elemento $f\in Z^2(G,A)$ (llamado
{\em conjunto de factores\index{conjunto de factores}}) determina
un \'unico grupo $E$ tal que $A\lhd E$ y $E/A\cong G$, es decir,
$E$ est\'a definido por medio de la sucesi\'on exacta
\begin{gather*}
0\to A\hookrightarrow E\xrightarrow{\pi} G\to 1
\intertext{donde 
$G$ act\'ua en el grupo abeliano $A$ por:} 
\text{si $g=\pi(e)\in G$,
entonces $g\circ a:=eae^{-1}$}.
\end{gather*}

Dos tales campos $E, E^{\prime}$, se llaman 
{\em equivalentes\index{grupos equivalentes}}
si existe un isomorfismo $\varphi\colon E\to E^{\prime}$ tal que
el diagrama 
\[
\xymatrix{
&&{E}\ar[rd]\ar[dd]^{\varphi}\\
{0}\ar[r]&{A}\ar[ur]\ar[dr]&&
{G}\ar[r]&{0}\\&&{E'}\ar[ur]}
\]
es conmutativo.

\begin{observacion}\label{CClaseO1.5.1} Si $E$ es equivalente
a $E^{\prime}$, entonces $E$ y $E^{\prime}$ son isomorfos pero
puede ser que $E$ y $E^{\prime}$ sean isomorfos pero no
equivalentes.

As{\'\i}, si $A$ es un grupo abeliano y $G$ es otro grupo que 
act\'ua sobre $A$, $g\circ a$: $G\times A\to A$, entonces hay
exactamente $|H^2(G,A)|$ grupos ${\mathcal G}$, salvo equivalencia,
tales que
$1\to A\to {\mathcal G}\to G\to 1$ es exacta con la acci\'on
de $G$ en $A$ dada.

Por ejemplo, se tiene que $H^2({\ma Z}/p{\ma Z},{\ma Z}/p{\ma Z})
\cong {\ma Z}/p{\ma Z}$ donde $G={\ma Z}/p{\ma Z}$ act\'ua
en $A={\ma Z}/p{\ma Z}$ de manera trivial. Por tanto hay $p$
clases de equivalencia de grupos de orden $p^2$, los cuales son
necesariamente abelianos, pero \'unicamente hay dos grupos 
abelianos de orden $p^2$ no isomorfos entre s{\'\i}, a saber
${\ma Z}/p^2{\ma Z}$ y ${\ma Z}/p{\ma Z}\oplus {\ma Z}/p{\ma Z}$.
De hecho el elemento identidad de $H^2({\ma Z}/p{\ma Z},
{\ma Z}/p{\ma Z})$ corresponde a ${\ma Z}/p{\ma Z}\oplus
{\ma Z}/p{\ma Z}$ y los otros $p-1$ elementos corresponden
a diversos grupos ${\ma Z}/p^2{\ma Z}$ (isomorfos pero
inequivalentes).

\end{observacion}

\subsection{Cohomolog\'ia de Galois y grupos de cohomolog\'ia de 
Tate}\label{SCo17.5.2'}

Sea $G$ un grupo de Galois de alguna extensi\'on de campos $L/K$
y sea $A$ un $G$-m\'odulo asociado de alguna manera a la extensi\'on
$L/K$. Entonces los $G$-grupos de cohomolog\'ia y de homolog\'ia de
$A$ si llaman {\em grupos de cohomolog\'ia de Galois\index{grupos
de cohomolog\'ia de Galois}}.

Uno de los resultados centrales de {\em cohomolog{\'\i}a de
Galois\index{cohomolog\'ia de Galois}} es

\begin{teorema}[Teorema 90 de Hilbert]\index{Hilbert!teorema 90
de $\sim$}\index{teorema 90 de Hilbert}\label{CClaseT1.5.7}
Si $L/K$ es una extensi\'on finita de Galois con grupo $G=\Gal(L/K)$
entonces $H^1(G,\* L)=\{1\}$.
\end{teorema}
\begin{proof} Sea $f\in Z^1(G,\*L)$, $f\colon G\to \* L$, $f(\theta\sigma)=
f(\theta)\cdot \theta f(\sigma)$ para cualesquiera $\theta, \sigma\in G$.

Del Teorema \ref{CClaseT1.1.1} seleccionamos $x\in\* L$ tal que $y=
\sum_{\sigma\in G} f(\sigma)\sigma(x)\in \* L$, es decir, $\sum_{\sigma\in G}
f(\sigma)\sigma\neq 0$. Entonces $\theta(y)=f(\theta)^{-1}y$
para toda $\theta \in G$ y por tanto
$f(\theta)=\theta(y^{-1})y=\theta(y)^{-1}y=\frac{y}{\theta(y)}
\in B^1(G,\* L)$. Se sigue que $H^1(G,\* L)=\{1\}.$ $\fin$
\end{proof}

\begin{observacion}\label{CClaseO1.5.7'} En realidad, el Teorema \ref{CClaseT1.5.7}
es una generalizaci\'on del Teorema 90 de Hilbert, el cual
es este mismo resultado pero \'unicamente para $G$ un grupo c\'iclico
finito. El Teorema \ref{CClaseT1.5.7} se debe a E. Noether y por tanto este
teorema es de Hilbert--Noether.
\end{observacion}

\begin{observacion}\label{CClaseO1.5.8} El Teorema \ref{CClaseT1.1.2} (b) 
(extensiones de Kummer) es una aplicaci\'on del
Teorema 90 de Hilbert con $G=\langle \sigma\rangle$ un grupo c{\'\i}clico
y usando que en este caso $H^1(G,\* L)\cong H^{-1}(G,\* L)$.

Para las extensiones de Artin--Schreier (Teorema \ref{CClaseT1.1.2} (a)), se usa
que $H^1(G, L)=0$ lo cual se cumple para cualquier grupo $G=\Gal(L/K)$
y cualquier campo. Ver la Observaci\'on \ref{CClaseO1.5.9}.
\end{observacion}

\begin{observacion}\label{CClaseO1.5.9} Si $L/K$ es una extensi\'on
finita de Galois, entonces el Teorema de la Base Normal
(ver Teorema \ref{T1.4.3}) establece
que existe $z\in L$ tal que $\{\sigma z\}_{\sigma\in G}$, con
$G=\Gal(L/K)$, es una base de $L/K$. Esto nos dice que, como
$G$--m\'odulos, se tiene 
\[
L\cong K[G]=\Big\{\sum_{\sigma\in G}\alpha_{\sigma}
\sigma\mid \alpha_{\sigma}\in K\Big\}\cong K\otimes_{\ma Z} {\ma Z}[G]
\]
de donde se sigue que $H^n(G,L)=0$ para toda $n\in {\ma Z}$.
Ver Corolario \ref{CoC17.5.8}.
\end{observacion}

A continuaci\'on presentamos los llamados grupos de cohomolog\'ia de Tate,
que son una amalgama entre la homolog\'ia y la cohomolog\'ia, modificando
ligeramente los grupos $\ho 0GA$ y $\co 0GA$ y re-indexando
adecuadamente. El proceso tambi\'en se puede obtener por medio de la
resoluci\'on completa, ver Subsecci\'on \ref{SResolucionCompleta}.

Sean $\N=\sum_{\sigma\in G}\sigma\in\zg$ y $I_G$ el ideal aumentaci\'on.
Se tiene $\N((\sigma-1)a)=0$ para cualquier $\sigma \in G$ y por tanto
$I_G\subseteq \ker \N$. Por otro lado $\N\sigma a=\sigma \N a=\N a$,
de donde obtenemos que $\im \N\subseteq A^G$. Recordemos que
$\ho 0GA=A_G=A/I_GA$ y que $\co 0GA=A^G$. Por tanto, pasando a
cocientes, $\*\N$ define un homomorfismo $\*\N\colon \ho 0GA\lra\co 0GA$.
Sea $\hot 0GA=\ker \*\N=\ker \N/I_G A$ y 
$\cot 0GA=\coker\*\N=A^G/\N A$.

Se tienen las sucesi\'on exactas
\[
0\lra \hot 0GA\lra \ho 0GA\xrightarrow[]{\N^*_A}\co 0GA\lra\cot 0GA\lra 0.
\]

\begin{teorema}\label{CoT17.5.2.0}
Sea $G$ un grupo finito y sea $0\lra A\stackrel{f}{\lra} B
\stackrel{g}{\lra} C\lra 0$ una sucesi\'on exacta de $G$-m\'odulos.
Entonces el diagrama
\[
\xymatrix{
\ho 1GC\ar@{->}[r]^{\varepsilon_0}\ar@{->}[d]&\ho 0GA
\ar@{->}[r]^{H_0(f)}\ar@{->}[d]^{\*\N_A}&\ho 0GB\ar@{->}[r]^{
H_0(g)}\ar@{->}[d]^{\N^*_B}&\ho 0GC\ar@{->}[d]^{\*\N_C}
\ar@{->}[r]&0\ar@{->}[d]\\
0\ar@{->}[r]& \co 0GA\ar@{->}[r]_{H^0(f)}&\co 0GB\ar@{->}[r]_{H^0
(g)}&\co 0GC\ar@{->}[r]_{\delta_0}&\co 1GA
}
\]
es conmutativo y las filas son exactas, donde $\varepsilon_0$ y
$\delta_0$ son los homomorfismos de conexi\'on.
$\fin$
\end{teorema}

Por el Lema de la Serpiente (Teorema \ref{CClaseT1.5.2}), se obtiene:

\begin{corolario}\label{CoC17.5.2.1}
Existe un homomorfismo can\'onico ($\ker\*\N_C=\hat H_0(G,C)$
y $\coker \*\N_A=\hat H^0(G,A)$)
\begin{gather*}
\delta\colon \hot 0GC\lra \cot 0GA
\intertext{que hace que la sucesi\'on de grupos}
\hot 0GC\lra\hot 0GB\lra \hot 0GC\stackrel{\delta}{\lra}\\
\stackrel{\delta}{\lra}\cot 0GA
\lra \cot 0GB\lra \cot 0GC
\end{gather*}
sea exacta. Formalmente $\delta=f^{-1}\circ \N_B^*\circ g^{-1}$.

Adem\'as $\delta$ nos da una sucesi\'on exacta
\begin{gather*}
\cdots\lra \ho 1GC \stackrel{\varepsilon_0}{\lra}
\hot 0GC\lra\hot 0GB\lra \hot 0GC\stackrel{\delta}{\lra}\\
\stackrel{\delta}{\lra}\cot 0GA
\lra \cot 0GB\lra \cot 0GC\stackrel {\delta_0}{\lra}\co 1GA
\lra \cdots \tag*{$\fin$}
\end{gather*}
\end{corolario}

\begin{definicion}\label{grupos de Tate}
Sea $G$ un grupo finito. Se definen
los {\em grupos de cohomolog\'ia de Tate\index{grupos de
cohomolog\'ia de Tate}} de $A$, con exponentes en ${\ma Z}$,
se definen por:
\begin{align*}
\cot nGA&=\co nGA \text{\ para\ } n \geq 1.\\
\cot 0GA&= A^G/\N A.\\
\cot {{-1}}GA&=\ker\N_A/I_G A.\\
\cot {{-n}}GA&= \ho {{n-1}}GA \text{\ para\ } n\geq 2.
\end{align*}
\end{definicion}

Como consecuencia del Corolario \ref{CoC17.5.2.1} y del
Teorema \ref{CClaseT1.5.3}, se obtiene

\begin{teorema}\label{sucesion exacta}
Si $0\lra A\lra B\lra C\lra 0$ es una sucesi\'on exacta de $G$-m\'odulos,
entonces
\begin{gather*}
\cdots\lra \cot {{n-1}}GC \lra \cot nGA \lra \cot nGB\lra\\
\lra \cot nGC\lra
\cot {{n+1}}GA \lra\cdots
\end{gather*}
es una sucesi\'on exacta para toda $n\in {\ma Z}$. 
Los mapeos de conexi\'on est\'an dados como en el Teorema
{\rm{\ref{CClaseT1.5.3}}}. $\fin$
\end{teorema}

\begin{notacion}\label{grupos Tate}
Sean $G$ un grupo finito y $A$ un $G$-m\'odulo. Los grupos de cohomolog\'ia
de Tate de $G$ con coeficientes en $A$ ser\'an indistintamente denotados
por $\cot nGA$ y por $\co nGA$. En adelante, a menos que se diga lo
contrario, los grupos de cohomolog\'ia ser\'an los grupos de cohomolog\'ia
de Tate.
\end{notacion}

\subsubsection{Grupos c\'iclicos}

Sea $G=\langle\sigma\rangle$ un grupo finito de orden $n$. Sea
$\N=\sum_{i=0}^{n-1}\sigma^i$ y $D=\sigma -1$. Se tiene que
$ND=DN=\sigma^n-1=0$. Consideremos el ideal aumentaci\'on $I_G=
\langle g-1\mid g\in G\rangle=\langle\sigma-1\rangle =D\zg$.

\begin{proposicion}\label{CoPCo17.5.2.2}
Se tiene $\ker \N=I_G=\im D$ y $\ker D=\zg^G=\im \N$.
$\fin$
\end{proposicion}

Sean $T_i:=\zg$, $i=0,1,\ldots$ y $\partial_i\colon T_i\lra T_{i-1}$
dada por 
\[
\partial_i=\begin{cases} D&\text{si $i$ es impar,}\\
\N& \text{si $i$ es par},\end{cases},\quad i\geq 1.
\]
Sea
$\varepsilon\colon \zg\lra {\ma Z}$ el homomorfismo aumentaci\'on: $\varepsilon
\big(\sum_{i=0}^{n-1} a_i\sigma^i\big)=\sum_{i=0}^{n-1} a_i$.

Se tiene que la sucesi\'on 
\[
\cdots\lra T_i\too i{} T_{i-1}\lra\cdots\lra T_1\too 1{}T_0
\stackrel{\varepsilon}{\lra}{\ma Z}\lra 0,
\]
de $G$-m\'odulos es exacta.

Esta es una resoluci\'on libre de ${\ma Z}$ y por tanto, para $n=1,2,\ldots $,
se tiene
\begin{gather*}
\cot {{2n-1}}GA=\co {{2n-1}} GA=\frac{\ker \*\N}{\im \*D}=
\frac{\ker\N_A}{D A}=\cot {{-1}}GA,\\
\cot {{2n}}GA=\co {{2n}} GA=\frac{\ker\*D}{\im \*\N}=\frac{A^G}{\N A}
=\cot 0GA.
\intertext{Similarmente, para homolog\'ia, obtenemos para $n=1,2,\ldots $,}
\cot {{-2n}}GA=\ho {{2n-1}} GA=\frac{\ker D_*}{\im \N_*}=\frac{A^G}{\N A}
=\cot 0GA,\\
\cot {{-(2n+1)}}GA=\ho {{2n}} GA
=\frac{\ker \N_*}{\im D_*}=\frac{\ker \N_A}{D A}=\cot {{-1}}GA
=\cot 1GA.
\end{gather*}

\begin{teorema}\label{CoT17.5.2.3}
Sea $G$ un grupo c\'iclico finito. Entonces, para cualquier $G$-m\'odulo $A$,
tenemos
\begin{gather*}
\cot {{2n}}GA=\cot 0 GA=\frac{A^G}{\N A},\\
\cot {{2n+1}}GA=\cot {{-1}} GA=
\frac{\ker\N_A}{D A},
\end{gather*}
para $n\in{\ma Z}$. $\fin$
\end{teorema}

\begin{definicion}\label{cociente de Herbrand}
Sea $G$ es un grupo c{\'\i}clico finito y sea $A$ un $G$--m\'odulo
tal que $H^0(G,A)$ y $H^1(G,A)$ son finitos de \'ordenes $h_0(A)$ y
$h_1(A)$ respectivamente. Se define el {\em cociente de
Herbrand\index{cociente de Herbrand}\index{Herbrand!cociente de $\sim$}}
de $A$ por $h(A):=\frac{h_0(A)}{h_1(A)}$.
\end{definicion}

Se tiene 

\begin{teorema}\label{CClaseT1.5.5}
Sea $G$ un grupo c{\'\i}clico finito y sea $0\to A \xrightarrow{f}B\xrightarrow{g}
C\to 0$ una sucesi\'on exacta de $G$--m\'odulos. Se tiene que el
siguiente hex\'agono es exacto:
\[
\xymatrix{
&{{H}^0(G,A)} \ar[r]^{f_0}
&{{H}^0(G,B)}\ar[rd]^{g_0}\\
{{H}^1(G,C)}\ar[ru]^{\delta_1}
&&&{{H}^0(G,C)}\ar[dl]^{\delta_0}\\
&{{H}^1(G,B)}\ar[lu]^{g_1}
&{{H}^1(G,A)}\ar[l]^{f_1}}
\]
y si dos de $h(A)$, $h(B)$ y $h(C)$ est\'an definidos, el tercero tambi\'en
est\'a definido y se tiene $h(B)=h(A)h(C)$.
\end{teorema}

\begin{proof} Se sigue inmediatamente de la sucesi\'on larga de cohomolog{\'\i}a
(Teorema \ref{CClaseT1.5.3}) y del hecho de que $H^0(G,X)\cong H^2(G,X)$
(Teorema \ref{CoT17.5.2.3}). $\fin$
\end{proof}

\begin{proposicion}\label{CClaseP1.5.6}
Si $G$ es un grupo c{\'\i}clico finito y $A$ es un $G$--m\'odulo finito,
entonces $h(A)=1$, esto es, $|H^i(G,A)|$ es constante para toda 
$i\in{\ma Z}$.
\end{proposicion}

\begin{proof} \cite[Proposition A.4.6]{Vil2006}. $\fin$
\end{proof}

\begin{corolario}\label{CoC1.5.7}
Sea $f\colon A\lra B$ un $G$-homomorfismo tal que $\ker f$ y $\coker f$
son finitos. Entonces $h(A)$ est\'a definido $\iff h(B)$ lo est\'a y en este
caso $h(A)=h(B)$.

En particular, si $A<B$ es de \'indice finito, entonces $h(A)=h(B)$.
$\fin$
\end{corolario}

Como consecuencia de la Proposici\'on \ref{CoP17.5.4}, se tiene 
la aditividad de la cohomolog\'ia. Ver tambi\'en 
\cite[Propositions (3.7) y (3.8), Part I]{Neu69}.

\begin{proposicion}\label{CClaseP1.5.9'}
Sean $G$ un grupo finito y $\{A_i\}_{i\in{\mc I}}$ una familia de
$G$--m\'odulos. Entonces
\las
\item $H^n(G,\bigoplus_i A_i)\cong \bigoplus_i H^n(G,A_i)$,

\item $H^n(G,\prod_i A_i)\cong \prod_i H^n(G,A_i)$,
\end{list}
para toda $n\in{\ma Z}$. $\fin$
\end{proposicion}

\subsection{M\'odulos inducidos y co-inducidos\index{modulos
inducidos@m\'odulos inducidos}\index{modulos co-inducidos@m\'odulos
co-inducidos}}\label{inducidos}

\begin{definicion}\label{coinducido}
Sea $X$ un grupo abeliano y sea $G$ actuando en $X$ de manera trivial
$g\circ x=x$ para toda $g\in G$ y para toda $x\in X$. El $G$-m\'odulo
$A:=\Hom(\zg,X)$ se llama {\em $G$-m\'odulo coinducido\index{modulo
coinducido@m\'odulo coinducido}} por $X$.
\end{definicion} 

Se tiene para $\varphi\in A$, $g,g'\in G$, $(g\circ \varphi)(g')=g\varphi(
g^{-1}g')=\varphi(g^{-1}g')$.

\begin{definicion}\label{inducido}
Un $G$-m\'odulo $A$ se llama inducido\index{modulo inducido@m\'odulo
inducido} si tiene la forma $A=\zg\otimes D$ donde $D$ es un grupo 
abeliano (considerado como $G$-m\'odulo trivial).
\end{definicion}

\begin{proposicion}\label{CoP17.5.1.3}
Sea $G$ un grupo finito y sea $D$ un grupo abeliano con $G$-acci\'on
trivial. Entonces
\[
\zg\otimes D\cong \Hom(\zg, D)\cong \bigoplus_{\sigma \in G}\sigma D.
\]
Esto es, cuando $G$ es finito, inducido y co-inducido es lo mismo.
\end{proposicion}

\begin{proof}
Sea $A=\zg\otimes D$. Para $a\in A$, $a=\big(\suma \sigma G \alpha\big)\otimes
d=\sum_{\sigma\in G}\sigma(\alpha_{\sigma}\otimes d)=\sum_{\sigma\in G}\sigma
(a_{\sigma}d)\in\sum_{\sigma\in G}\sigma D$ y claramente la suma es
directa. De aqui obtenemos que $\zg\otimes 
D\cong \bigoplus_{\sigma\in G}\sigma D$.

Por otro lado, Sea $\Lambda\colon\Hom(\zg, D)\lra \bigoplus_{\sigma\in G}
\sigma D$, $f\longmapsto \sum_{\sigma\in G}\sigma f(\sigma)=\big(f(
\sigma)\big)_{\sigma\in G}$. Entonces $\Lambda$ es un isomorfismo.
$\fin$
\end{proof}

\begin{teorema}\label{CoT17.5.4}
Sea $A$ un $G$-m\'odulo co-inducido, entonces 
$\co qGA=0$ para toda $q\geq 1$.
\end{teorema}

\begin{proof}
Sea $B$ cualquier $G$-m\'odulo y sea $A$ un $G$-m\'odulo co-inducido, $A
=\sum_{\sigma\in G}\sigma D$ con $D$ un grupo abeliano, $G$-m\'odulo trivial.

Veamos que $\Hom_G(B,A)=\Hom(B,D)$.

Sea $f\in\Hom_G(B,A)=\big(\Hom(B,A)\big)^G$. Entonces $f(b)=\sum_{\sigma
\in G}\sigma f_{b,\sigma}$ con $f_{b,\sigma}\in D$. Se tiene que $\tau f(b)=
\sum_{\sigma\in G}\tau\sigma f_{b,\sigma}=\sum_{\theta}\theta f_{b,\tau^{-1}\theta}
=f(b)$. Por tanto $f_{b,\mu}=f_b$ constante para toda $\mu\in G$.

Por tanto $f(b)=\sum_{\sigma\in G}\sigma f_b=\big(\sum_{\sigma\in G}\sigma\big)
f_b$, $f_b\in D$. Definimos 
\[
\Lambda\colon\Hom_G(B,A)\lra\Hom(B,D)\quad \text{dada
por\ } \Lambda(f)=\tilde f, \tilde f(b)=f_b\in D.
\]
Entonces $\Lambda$ es un
isomorfismo.

Sea 
\[
\cdots\lra P_i\lra P_{i-1}\lra \cdots \lra P_1\lra P_0\lra {\ma Z}\lra 0
\leqno{(P)}
\]
una resoluci\'on y sea $K_i=\Hom_G(P_i,A)=\Hom(P_i,D)$. Puesto que
$P_i$ es ${\ma Z}$-libre, la sucesi\'on
\begin{gather*}
0\lra\Hom({\ma Z},D)\tooo 0{}\Hom(P_0,D)\tooo 1{}\Hom(P_1,D)\tooo 2{}
\cdots\\
\lra \Hom(P_{i-1},D)\tooo i{}\Hom(P_i,D)\lra \cdots
\end{gather*}
es exacta (Proposici\'on \ref{CoP17.5.7})
y $\co qGA=\frac{\ker d_{q+1}}{\im d_q}=0$ para $q\geq 1$.
$\fin$
\end{proof}

\begin{teorema}\label{CoT17.5.5}
Si $A$ es un $G$-m\'odulo inducido, entonces $\ho qGA=0$ para
toda $q\geq 1$.
\end{teorema}

\begin{proof}
Sea $B$ cualquier $G$-m\'odulo. Se tiene 
\begin{gather*}
B\otimes_G A\cong B\otimes_G
\big(\zg\otimes_{\ma Z}D\big)\cong\big(B\otimes_G\zg\big)\otimes_{\ma Z}
D\cong B\otimes_{\ma Z}D.
\intertext{Esto es}
B\otimes_G A\cong B\otimes_{\ma Z} D.
\end{gather*}

Consideremos la resoluci\'on
\[
\cdots\lra P_i\lra P_{i-1}\lra\cdots\lra P_1\lra P_0\lra {\ma Z}\lra 0.
\leqno{(P)}
\]
Obtenemos que
\begin{gather*}
0\lra {\ma Z}\otimes_G A\lra P_0\otimes_G A\lra P_1\otimes_G A
\lra \cdots\\
 \lra P_{i-1}\otimes D\lra P_i\otimes D\lra\cdots \\
\intertext{y por tanto}
 0\lra {\ma Z}\otimes D\lra P_0\otimes D\lra P_1\otimes D
 \lra \cdots\\
 \lra P_{i-1}\otimes D\lra P_i\otimes D\lra\cdots
 \end{gather*}
 Puesto que los m\'odulos $P_i$ son ${\ma Z}$-libres, la \'ultima sucesi\'on
 es exacta (Proposici\'on \ref{CoP17.5.7})
 de donde se sigue que $\ho qGA=0$ para toda $q\geq 1$.
$\fin$
\end{proof}

\begin{observacion}\label{CoO17.5.6}
Para la cohomolog\'ia usual, es decir, no de Tate, $\ho 0GA$ y $\co 0GA$
no tiene por que ser $0$ si $A$ es inducido y/o co-inducido. Por ejemplo
$\co 0GA=A^G\cong D$. Veremos que para la cohomolog\'ia de Tate,
estos dos grupos son $0$.
\end{observacion}

\begin{teorema}\label{CoT17.5.7}
Si $G$ es un grupo finito y $A$ es un $G$-m\'odulo inducido o co-inducido,
$\cot 0GA=\hot 0GA=0$.
\end{teorema}

\begin{proof}
Veamos que $\hot 0GA=0$ con $A$ un $G$-m\'odulo inducido. Sea $A=
\bigoplus_{\sigma\in D}\sigma D$. Sea $a\in A$, $a=\big(\sigma a_{\sigma}
\big)_{\sigma\in G}=\suma \sigma G a$. Si $a\in A^G$, 
\[
\tau a=\sum_{
\sigma\in G}\tau\sigma a_{\sigma}=\sum_{\theta\in G}\theta a_{\tau^{-1}
\theta}=\sum_{\sigma\in G}\sigma a_{\tau^{-1}\sigma}=\suma \sigma Ga
=a
\] 
para toda $\tau\in G$. Por tanto $a_{\sigma}=a_0\in D$ constante.
Se sigue que $a=\big(\sum_{\sigma\in G}\sigma\big)a_0=\N_G a_0$
de donde se sigue que $A^G\subseteq \N_G A$ y por tanto $\cot
0GA=\frac{A^G}{\N_G A}=\{0\}$.

Ahora veamos que $\hot 0GA=\frac{\ker\N_G}
{I_GA}=\{0\}$. Sea $a\in A$ como antes.
Entonces
\begin{align*}
\N_G a&=\sum_{\tau\in G}\tau\sum_{\sigma\in G}\sigma 
a_{\sigma}=\sum_{\tau\in G}
\sum_{\sigma\in G}\tau\sigma a_{\sigma}\\
&=\sum_{\theta\in G}\theta\big(
\sum_{\tau\in G}a_{\tau^{-1}\theta}\big)=\big(\sum_{\theta\in G}
\theta\big)\big(\sum_{\sigma\in G}a_{\sigma}\big).
\end{align*}

Por tanto $\N_G a=0\iff \sum_{\sigma\in G}a_{\sigma}=0$, esto
es, $\sum_{\sigma\in G}\sigma a_{\sigma}\in I_G A$ y por
tanto $\hat H_0(G,A)=\{0\}$.
$\fin$
\end{proof}

\begin{definicion}\label{CoD17.5.7'}
Un $G$-m\'odulo $A$ se llama {\em cohomol\'ogicamente
trivial\index{cohomologicamente trivial@cohomol\'ogicamente
trivial}} si $\cot qHA=\{0\}$ para toda $q\in {\ma Z}$ y para todo
subgrupo $H$ de $G$.
\end{definicion}

\begin{corolario}\label{CoC17.5.8}
Si $G$ es un grupo finito y $A$ es un $G$-m\'odulo inducido o 
co-inducido, entonces $A$ es cohomol\'ogicamente trivial.
\end{corolario}

\begin{proof}
Se sigue de los Teoremas \ref{CoT17.5.4}, \ref{CoT17.5.5} y
\ref{CoT17.5.7} y del hecho de que si $A$ es $G$-inducido o 
$G$-co-inducido, entonces $A$ es $H$-inducido o 
$H$-co-inducido para todo subgrupo $H$ de $G$.
$\fin$
\end{proof}

Sea $L/K$ una extensi\'on finita de Galois con grupo $G$. Entonces
tanto $L$ como $\*L$ son $G$-m\'odulos. Adem\'as, $L/K$ tiene una
base normal, esto es, existe $\alpha\in L$ tal que $\{\sigma \alpha\}_{
\sigma\in G}$ es una base de la extensi\'on $L/K$
(Teorema \ref{T1.4.3}). Como $G$-m\'odulo
se tiene que $L=\bigoplus_{\sigma\in G}K(\sigma \alpha)\cong
K\otimes_{\ma Z}\zg$. En particular $L$ es $G$-inducido y por
el Corolario \ref{CoC17.5.8} $L$ es cohomol\'ogicamente trivial.

El mapeo $f\colon G\lra I_G/I_G^2$ dado por $f(\sigma)=(\sigma-1)+I_G^2$
induce un isomorfismo
\[
G/G'\cong I_G/I_G^2,
\]
donde $G'$ es el subgrupo conmutador de $G$. El isomorfismo
inverso est\'a dado
por $h((\sigma-1)+I_G^2)=\sigma G'$.

\subsection{Resoluci\'on completa}\label{SResolucionCompleta}

El estudio de la homolog\'ia y la cohomolog\'ia hecha hasta aqu\'i, se
rompe en dos partes, una para cada una de ellas. John Tate desarroll\'o
un complejo en el cual la homolog\'ia y la cohomolog\'ia se funden
en un solo grupo a partir de una resoluci\'on completa, la cual presentamos
aqu\'i por conveniencia y para futura referencia.

Sea $C$ un ${\ma Z}$-m\'odulo. Definimos $\widehat C:=\Hom(C,{\ma Z})$.
Entonces $\widehat C$ es un ${\ma Z}$-m\'odulo. Sea $A$ cualquier 
${\ma Z}$-m\'odulo. Definimos
\begin{gather}
\Phi\colon C\otimes A\lra \Hom(\widehat C,A)\nonumber
\intertext{el homomorfismo de ${\ma Z}$-m\'odulos dado por}
\Phi(c\otimes a)(\xi)=\xi(c)\cdot a,\quad c\in C, a\in A, \xi\in \hat c,
\xi(c)\in{\ma Z}.\label{Iso}
\end{gather}

Si $C$ es un ${\ma Z}$-m\'odulo libre finitamente generado, esto es,
$C\cong {\ma Z}^m$ para alg\'un $m\in{\ma N}$, entonces $\widehat C
\cong C\cong {\ma Z}^m$. El isomorfismo est\'a dado de la siguiente
forma. Si $\{e_i\}_{i=1}^m$ es una base de $C$ como ${\ma Z}$-m\'odulo,
se define $\*{e_i}\in\widehat C$ por $\*{e_i}(e_j)=\langle \*{e_i},e_j\rangle=
\delta{ij}$, la delta de Kronecker. Entonces $\{\*{e_i}\}_{i=1}^m$ es una
${\ma Z}$-base de $\widehat C$. 
Se tiene un isomorfismo natural $\widehat{\widehat{C}}
\cong C$ identificando $e_i^{**}$ con $e_i$.

\begin{proposicion}[homotop\'ia de contracci\'on]\label{CoP17.5.10}
Dada una sucesi\'on exacta de ${\ma Z}$-m\'odulos libres
\begin{gather*}
\cdots\lra X_{n+1}\tooo {{n+1}}{}X_n\tooo n{}X_{n-1}\lra\cdots,
\intertext{existe un ${\ma Z}$-homomorfismo $s_n\colon X_n\lra X_{n+1}$
tal que}
d_{n+1}s_n+s_{n-1}d_n=\Id_{X_n}.
\end{gather*}
\end{proposicion}

\begin{proof}
Sea $U_n=\ker d_n=\im d_{n+1}$. Como $U_n<X_n$ y $X_n$ es un
${\ma Z}$-m\'odulo libre, $U_n$ es un ${\ma Z}$-m\'odulo libre. Se tiene
la sucesi\'on exacta
\begin{gather}\label{Ec17.5.2}
0\lra U_n\stackrel{i_n}{\hooklongrightarrow} X_n\tooo n{} U_{n-1}\lra 0.
\end{gather}
Como $U_{n-1}$ es un ${\ma Z}$-m\'odulo libre, la sucesi\'on (\ref{Ec17.5.2})
se escinde. Sean $\psi_n\colon U_{n-1}\lra X_n$, $\varphi\colon X_n\lra
U_n$ los mapeos de escisi\'on. Entonces $d_n\circ i_n=0$, $\varphi_n
\circ\psi_n=0$, $d_n\circ \psi_n=\Id_{U_{n-1}}$ y $\varphi_n\circ i_n =
\Id_{U_n}$. Sea $s_n:=\psi_{n+1}\circ \varphi_n\colon X_n\lra X_{n+1}$.
Escribiendo $X_n=i_n(U_n)\oplus \psi_n(U_{n-1})$, se verifica
el resultado.
$\fin$
\end{proof}

\begin{corolario}\label{CoC17.5.11}
Dada una sucesi\'on exacta de ${\ma Z}$-m\'odulos libres
\[
\cdots\lra X_{n+1}\tooo {{n+1}}{} X_n\tooo n{} X_{n-1}\lra\cdots,
\leqno{(X)}
\]
entonces la sucesi\'on
\[
\cdots\lra \widehat{X}_{n-1}\toooh {{n+1}} 
\widehat{X}_n\toooh n\widehat{X}_{n+1}\lra\cdots,
\leqno{(\widehat{X})}
\]
es exacta, donde $\widehat{d}_n\colon \widehat{X}_n\lra \widehat{X}_{n+1}$
est\'a dada por $\widehat{d}_n(f)=f\circ d_{n+1}$.
\end{corolario}

\begin{proof}
Sea $\widehat{s}_{n-1}\colon \widehat{X}_{n-1}\lra\widehat{X}_{n}$, definido
por $\widehat{s}_{n-1}(f)=f\circ s_{n-1}$ $\xymatrix{
X_n\ar@{->}[r]^f\ar@{<-}[d]_{s_{n-1}} & {\ma Z}\\
X_{n-1}\ar@{-->}[ru]_{f\circ s_{n-1}}}$. Entonces $\widehat{d}_n\circ
\widehat{s}_n+\widehat{s}_{n+1}\circ \widehat{d}_n=\Id_{\widehat{X}_n}$.
De aqu\'i, se verifica directamente que la sucesi\'on $(\widehat{X})$ es
exacta. $\fin$
\end{proof}

Ahora consideremos $\zg$-m\'odulos, $G$ un grupo finito. Sea $C$ un
$G$-m\'odulo izquierdo y $\widehat C=\Hom(C,{\ma Z})$ es un
$G$-m\'odulo izquierdo con la acci\'on usual:
\[
(\sigma f)(c)=\sigma(f(\sigma^{-1}c))=f(\sigma^{-1}c),
\]
pues $G$ act\'ua trivialmente en ${\ma Z}$.

Consideremos $\Phi$ dado en
(\ref{Iso}). Directamente se verifica que $\sigma(\Phi(c\otimes a))=\Phi(\sigma(
c\otimes a))$ y por tanto $\Phi$ es un $G$-homomorfismo.

Si $C$ es un ${\ma Z}$-m\'odulo libre de rango finito, $\Phi$ es
biyectiva. De hecho, sea $\{e_i\}_{i=1}^n$ una base de $C$ como
${\ma Z}$-m\'odulo. Si $\Phi(c\otimes a)=0$, entonces si $c=\sum_{
i=1}^n\alpha_ie_i$, evaluando en $\xi=\*{e_j}$ obtenemos
$\xi(c)\cdot a=\alpha_ja=0$ para toda $j$. Por tanto $c\otimes a=
\sum_{i=1}^n\alpha_ie_i\otimes a=\sum_{i=1}^ne_i\otimes \alpha_i a=0$
de donde se sigue que $\Phi$ es inyectiva.

Ahora, sea $f\in\Hom(\widehat C,A)$. Sea $f(\*{e_j})=a_j$ para $1\leq j
\leq n$. Entonces $\Phi\big(\sum_{i=1}^n e_i\otimes a_i\big)=f$ y
$\Phi$ es suprayectiva.

Con lo anterior, hemos probado

\begin{teorema}\label{CoT17.5.11}
Sea $C$ un $\zg$-m\'odulo libre con base finita. Entonces 
\[
C\otimes A\cong\Hom(\widehat C, A)
\]
como $G$-m\'odulos izquierdos y donde el isomorfismo
esta dado por $\Phi$ definido en {\rm{(\ref{Iso})}}. Tambi\'en
tenemos que $\widehat{\widehat{C}}\cong C$ como
$G$-m\'odulos. $\fin$
\end{teorema}

\begin{proposicion}\label{CoP17.5.12}
Si $C$ es un $\zg$-m\'odulo libre con una base finita, entonces $\widehat
C$ tambi\'en lo es y $\widehat C\cong C$ como $\zg$-m\'odulos.
\end{proposicion}

\begin{proof}
Primero supongamos $C\cong\zg$. Si probamos que $\widehat C\cong
\zg$, entonces para $C=\bigoplus_{i=1}^n\zg f_i$, se seguir\'a que $\widehat
C\cong\bigoplus_{i=1}^n \widehat{\zg f_i}\cong \bigoplus_{i=1}^n \zg
\*{f_i}$. Por tanto basta probar el caso $C\cong \zg$, esto es, debemos 
probar que $\Hom(\zg,{\ma Z})\cong \zg$.

Sea $\{e_i\}_{i=1}^n$ una base de $C=\zg$ como ${\ma Z}$-m\'odulo
(por ejemplo, $\{e_i\}_{i=1}^n=G$). Sea $\{\*{e_i}\}_{i=1}^n\subseteq
\widehat C$ la base dual, esto es, $\*{e_j}(e_i)=\delta_{ij}$ para $1\leq
i,j\leq n$. Sea $\varphi$ el ${\ma Z}$-isomorfismo $\varphi\colon C\lra
\widehat C$ dado por $\varphi(e_i)=\*{e_i}$, $1\leq i\leq n$. Es inmediato
que $\varphi$ es $\zg$-isomorfismo. $\fin$
\end{proof}

Hemos probado

\begin{proposicion}\label{CoP17.5.13}
Dada una sucesi\'on exacta de $\zg$-m\'odulos izquierdos y $\zg$-homomorfismos
\[
\cdots\lra X_{n+1}\tooo {{n+1}}{} X_n\tooo n{} X_{n-1}\lra\cdots,
\leqno{(X)}
\]
tales que cada $X_n$ es un $\zg$-m\'odulo libre con base finita, entonces la
sucesi\'on
\[
\cdots\lra \widehat{X}_{n-1}\toooh {{n+1}} \widehat{X}_n\toooh n\widehat{X}_{n+1}\lra\cdots,
\leqno{(\widehat{X})}
\]
es exacta, donde cada $\widehat {X}_n$ es un $\zg$-m\'odulo libre con base finita.
$\fin$
\end{proposicion}

Ahora consideremos dos $G$-m\'odulos izquierdos $C$ y $A$. En el producto tensorial,
$C$ se considera como un $G$-m\'odulo derecho definiendo
\[
c\sigma:=\sigma^{-1}c,\quad c\in C, \sigma\in G.
\]
Se sigue que, considerando $C\otimes A$ como un $G$-m\'odulo izquierdo con:
$\sigma(c\otimes a)=(\sigma a)\otimes(\sigma a)$, $c\in C$, $a\in A$ y $\sigma \in
G$, se tiene

\begin{proposicion}\label{CoP17.5.14}
$C\otimes_G A\cong (C\otimes A)_G=\frac{C\otimes A}{I_G(C\otimes A)}$
y $\Hom_G(C,A)=\Hom(C,A)^G$.
\end{proposicion}

\begin{proof}
Por definici\'on $C\otimes_G A$, donde $C$ es un $G$-m\'odulo derecho y 
$A$ es un $G$-m\'odulo izquierdo, es $F/R_G$ donde $F$ es el grupo
abeliano libre generado por los elementos $(c,a)$ de $C\times A$: $F=
\bigoplus_{\substack{c\in C\\a\in A}}{\ma Z} (c,a)$ y $R_G$ es el 
subgrupo generado por
\begin{gather*}
(c+c_1,a)-(c,a)-(c_1,a);\quad (c,a+a_1)-(c,a)-(c,a_1);\\
(c\sigma,a)-(c,\sigma a),\quad \sigma\in G.
\end{gather*}

Ahora bien, si consideramos $C\otimes A=C\otimes_{\ma Z} A$, $F$ es el
mismo y $R_{\ma Z}$ es el subgrupo generado por
\begin{gather*}
(c+c_1,a)-(c,a)-(c_1,a);\quad (c,a+a_1)-(c,a)-(c,a_1).
\intertext{Se tiene la sucesi\'on exacta}
0\lra R_G/R_{\ma Z}\lra F/R_{\ma Z}\lra F/R_G\lra 0.
\end{gather*}
El grupo $R_G/R_{\ma Z}$ visto en $F/R_{\ma Z}=C\otimes A$
es el grupo generado por $\langle (c\sigma\otimes a)-(c\otimes \sigma a)
\mid c\in C, a\in A, \sigma\in G\rangle$.

Vemos a $C$ como $G$-m\'odulo izquierdo: $c\sigma :=\sigma^{-1} c$.
De esta forma tenemos que $C\otimes_G A=(C\otimes A)/{\mc R}$ donde
${\mc R}=\langle (\sigma c,a)-(c,\sigma a)\rangle$. Se tiene que
\[
(c\sigma,a)-(c,\sigma a)=(\sigma^{-1}c,a)-(c,\sigma a)=(\sigma^{-1}-1)
(c,\sigma a)\in I_G(C\otimes A),
\]
de donde se sigue que $C\otimes_G A=\frac{C\otimes A}{I_G(C\otimes A)}$.

Sea $f\in \Hom_G(C,A)$. En particular $f\in \Hom(C,A)$. Ahora $G$ act\'ua
en $\Hom(C,A)$ por $(\sigma\circ f)(c)=\sigma f(\sigma^{-1} c)$. Ahora si
$f\in \Hom_G(C,A)$, $f(\sigma^{-1} c)=\sigma^{-1} f(c)$, entonces
\[
(\sigma\circ f)(c)=\sigma f(\sigma^{-1} c)=\sigma\sigma^{-1}f(c)=f(c),
\]
esto es, $(\sigma\circ f)=f$ de donde se sigue que $f\in\Hom(C,A)^G$.

Rec\'iprocamente, si $f\in\Hom(C,A)^G$, $(\sigma\circ f)=f$, por lo tanto,
$f(\sigma c)=(\sigma\circ f)(\sigma c)=(\sigma f)(\sigma^{-1}(\sigma c))=
(\sigma f)(c)$, por lo que $f\in\Hom_G(C,A)$.
$\fin$
\end{proof}

Sea $C$ un $\zg$-m\'odulo libre de rango finito. Por el Teorema 
\ref{CoT17.5.11} se sigue que $\widehat C\otimes A\cong\Hom(C.A)$.

\begin{proposicion}\label{CoP17.5.15}
Sea $C$ un $\zg$-m\'odulo libre de rango finito y sea $A$ un grupo abeliano.
Entonces la norma
$\N\colon (\widehat C\otimes A)\lra (\widehat C\otimes A)$ induce un 
isomorfismo
\[
\*{\N}\colon (\widehat C\otimes A)_G\lra (\widehat C\otimes A)^G.
\]
\end{proposicion}

\begin{proof}
Notemos que $\zg\otimes A\cong (\bigoplus_{\sigma \in G}\sigma {\ma Z})
\otimes A\cong  \bigoplus_{\sigma \in G} \sigma({\ma Z}\otimes A)\cong
\bigoplus_{\sigma\in G} \sigma A$. Se sigue que $X:=\widehat C\otimes A=
\bigoplus_{\sigma\in G}\sigma A$.

Sea $a\in X$, $a=(\sigma a_{\sigma})_{\sigma \in G}=\sum_{\sigma\in G}
\sigma a_{\sigma}$. Si $a\in X^G$, entonces 
\[
\tau a=\sum_{\sigma \in G}
\tau\sigma a_{\sigma}=\sum_{\theta\in G}\theta a_{\tau^{-1}\theta}=
\sum_{\sigma\in G}\sigma a_{\tau^{-1}\sigma}=\sum_{\sigma\in G}
\sigma a_{\sigma}=a,
\]
para toda $\tau\in G$. Se sigue que $a_{\sigma}=a_{\sigma'}$ para
cualesquiera $\sigma, \sigma'\in G$. Por tanto $a_0=a_{\sigma}$
es constante para toda $\sigma\in G$. Se sigue que $a=\big(\sum_{
\sigma\in G}\sigma \big)a_0=\N_G a_0$. Por tanto $X^G\subseteq \N_G
X$, lo cual implica que $\*{\N}$ es suprayectiva.

Ahora sea $a\in X$, $a=\sum_{\sigma\in G}\sigma a_{\sigma}$ con
$a_{\sigma}\in D$. Entonces
\[
\N_Ga=\sum_{\tau\in G}\sum_{\sigma\in G}\sigma a_{\sigma}=
\big(\sum_{\theta\in G}\theta\big)\big(\sum_{\sigma\in G} a_{\sigma}\big).
\]
Por tanto $\N_Ga=0\iff \sum_{\sigma\in G}a_{\sigma}=0$ y en este caso
se tiene $\sum_{\sigma\in G}\sigma a_{\sigma}=\sum_{\sigma\in G}
\sigma a_{\sigma}-\sum_{\sigma\in G}1 a_{\sigma}=\sum_{\sigma\in G}
(\sigma -1)a_{\sigma}\in I_G X$. Por tanto $\*\N$ es inyectiva y es un
isomorfismo. $\fin$
\end{proof}

Como consecuencia, tenemos

\begin{teorema}\label{CoT17.5.16}
Sean $C$ y $A$ dos $G$-m\'odulos izquierdos con $C$ un $\zg$-m\'odulo
libre de rango finito. Entonces existe un isomorfismo $\tau\colon
\widehat C\otimes_G A\lra \Hom_G(C,A)$ dado por
\[
[\tau(f\otimes_G a)](c)=\sum_{\sigma\in G} f(\sigma^{-1}c)\sigma a,
\]
con $f\in\widehat C$, $c\in C$ y $a\in A$.
\end{teorema}

\begin{proof}
Se hab\'ia probado que $\widehat C\otimes A\cong \Hom(C,A)$
(Teorema \ref{CoT17.5.11}) con isomorfismo $\Phi$ dado por
$\Phi\colon \widehat C\otimes A\lra \Hom(\widehat{\widehat{C}},A)
\cong \Hom(C,A)$, 
\[
\Phi(f\otimes a)(\xi)=\hat{\xi}(f) a=f(\xi)a.
\]

Por otro lado $\*\N\colon (\widehat C\otimes A)_G\xrightarrow[]{\cong}
\Hom(C,A)^G=\Hom_G(C,A)$ es un isomorfismo. Se tiene
\[
\xymatrix{
(\widehat C\otimes A)_G\ar@{->}[r]^{\*\N}_{\cong}\ar@{-->}[rd]_{\tau}
&(\widehat C\otimes A)^G\ar@{->}[d]^{\Phi}\\
&\Hom(C,A)^G=\Hom_G(C,A)}
\]
Entonces $\tau=\Phi\circ \*\N$ es el isomorfismo. Se tiene
\begin{align*}
\tau(f\otimes a)(c)&=\sum_{\sigma\in G}\Phi(\sigma f\otimes\sigma a)(c)=
\sum_{\sigma\in G}\hat c(\sigma f)\sigma a=\sum_{\sigma\in G}
(\sigma f)(c)\sigma a\\
&=\sum_{\sigma\in G}\sigma f(\sigma^{-1} c)\sigma a=
\sum_{\sigma\in G}f(\sigma^{-1} c)\sigma a.\tag*{$\fin$}
\end{align*}
\end{proof}

\begin{definicion}\label{CoD17.5.17}
Una {\em resoluci\'on completa\index{resolucion completa@resoluci\'on
completa}} $P$ para un grupo finito $G$ es una sucesi\'on exacta
\[
\cdots\lra P_n\tooo n{} P_{n-1}\lra\cdots\lra P_0\tooo 0{}P_{-1}\lra
\cdots\lra P_{-n}\tooo {{-n}}{}\cdots \leqno{(P)}
\]
de m\'odulos $P_n$ que son $\zg$-libres finitamente generados
junto con un elemento $e\in(P_{-1})^G$, $e\neq 0$, tal que la imagen
de $d_0$ est\'a generada por $e$.
\end{definicion}

Puesto que $\sigma e=e$ para toda $\sigma\in G$, se sigue de que
$\im d_0$ es un ${\ma Z}$-m\'odulo generado por $e$. Adem\'as 
como $P_{-1}$ es ${\ma Z}$-libre, se tiene $ne\neq 0$ para toda
$n\in{\ma Z}$, $n\neq 0$. Por tanto $d_0$ admite una factorizaci\'on
$d_0=\mu\circ \varepsilon$, donde $\varepsilon$ es un $G$-epimorfismo
y $\mu$ es un $G$-monomorfismo dado por $\mu(1)=e$.
\[
\xymatrix{
P_0\ar@{->}[r]^{d_0}\ar@{-->}[ddr]_{\varepsilon}&\langle e\rangle
\ar@{->}[dd]_{\lambda}^{\begin{array}{c}e\\ \downarrow \\1\end{array}}&
P_0\ar@{->}[rr]^{d_0}\ar@{->}[rd]_{\varepsilon}&&P_{-1}\\
&&&{\ma Z}\ar@{->}[ru]_{\mu}\ar@{->}[dr]\\
&{\ma Z}&0\ar@{->}[ru]&&0
}
\]

Consideremos las sucesiones exactas
\[
\cdots\lra P_n\tooo n{} P_{n-1}\tooo {{n-1}}{}\cdots\lra P_0\stackrel{
\varepsilon}{\lra}{\ma Z}\lra 0
 \leqno{(P_i)}
 \]
 \[
0\lra {\ma Z}\stackrel{\mu}{\lra} P_{-1}\tooo {{-1}}{}\cdots\tooo {{-n+2}}{}
P_{-n+1}\tooo {{-n+1}}{} P_{-n}\lra\cdots \leqno{(P_d)}
\]

La sucesi\'on $(P_i)$ da una resoluci\'on proyectiva de ${\ma Z}$ con
$\zg$-m\'odulos libres finitamente generados. Por el Corolario \ref{CoC17.5.11}
y el dual de $(P_d)$: $\widehat{P}_{-n}=\Hom(P_{-n},{\ma Z})$, se tiene que
\[
\cdots\lra \widehat{P}_{-n}\toooh {{-n+1}} \widehat{P}_{-n+1}
\toooh {{-n+2}}\cdots\toooh{{-1}} \widehat{P}_{-1}\stackrel{
\hat\mu}{\lra}{\ma Z}\lra 0
 \leqno{(\widehat{P}_d)}
\]
tambi\'en da una resoluci\'on proyectiva por medio de $\zg$-m\'odulos
libres finitamente generados. Rec\'iprocamente, dadas dos resoluciones
$(P_i)$ y $(P_i')$ de ${\ma Z}$ por $\zg$-m\'odulos finitamente
generados, podemos construir una resoluci\'on completa empalmando
$(P_i)$ y $(\widehat{P}_i')$ reenumerados adecuadamente.

Dada una resoluci\'on completa $(P)$ y un $G$-m\'odulo izquierdo $A$,
podemos considerar el {\em complejo\index{complejo}}
\[
\Hom_G(P,A)=\big\{\Hom_G(P_n,A)\big\}_{n=-\infty}^{\infty}.
\]

Para $n\geq 0$, dejamos el grupo $\Hom_G(P_n,A)$ como est\'a. Para
$n<0$ reemplazamos $\Hom_G(P_n,A)$ por el grupo isomorfo 
$\widehat{P}_n\otimes_G A$ usando el isomorfismo $\tau$
dado en el Teorema \ref{CoT17.5.16}.

Veamos el mapeo $\widehat{P}_{-1}\otimes_G A\xrightarrow{\ \ \alpha\ \ }{}
\Hom_G(P_0,A)$, inducido por $d_0\colon P_0\lra P_{-1}$. Se tiene
$d_0=\mu\circ\varepsilon$ y por tanto obtenemos el diagrama conmutativo
\[
\xymatrix{
\ar@{-->}@/^2pc/[rr]_{\alpha}
\widehat{P}_{-1}\otimes_G A\ar@{->}[r]^{\Phi\ \ }_{\cong\ \ }\ar@{->}[d]_{
\hat{\mu}\otimes_G 1}&\Hom_G(P_{-1},A)\ar@{->}[r]^{\widehat{d}_0}
\ar@{->}[d]_{\hat{\mu}}&\Hom_G(P_0,A)\\
{\ma Z}\otimes_G A\ar@{->}[r]_{\Phi}&\Hom_G({\ma Z},A)\ar@{->}[ru]_{
\hat{\varepsilon}}}
\]
Esto es, $\alpha$ se factoriza:
\[
\xymatrix{
\widehat{P}_{-1}\otimes_G A\ar@{->}[r]^{\alpha\ \ }\ar@{->}[d]&\Hom_G(P_0,A)\\
H_0(G,A)\ar@{->}[r]^{\*\N}&H^0(G,A)\ar@{->}[u]
}
\]
De esto obtenemos finalmente.

\begin{teorema}\label{CoT17.5.18}
Para cualquier $G$-m\'odulo izquierdo $A$, los grupos de cohomolog\'ia de
Tate $\hat H^n(G,A)$, $n\in{\ma Z}$ se obtiene
como $H^n\big(\Hom_G(P,A)\big)$ donde
$P$ es cualquier resoluci\'on completa de $G$. Si $f\colon A\lra B$ es un
homomorfismo de $G$-m\'odulos, entonces los $H^n(f)$ puede ser
calculados de $\Hom_G(P,A)\lra\Hom_G(P,B)$. Si $0\lra A\lra B\lra C
\lra 0$ es una sucesi\'on exacta de $G$-m\'odulos izquierdos, los mapeos
de conexi\'on $\hat H^n(G,C)\xrightarrow{\ \ \delta \ \ }{}\hat H^{n+1}(G,A)$
pueden ser calculados de la sucesi\'on exacta de resoluciones
\begin{gather*}
0\lra \Hom_G(P,A)\lra \Hom_G(P,B)\lra \Hom_G(P,C)\lra 0,
\end{gather*}
como en el Teorema {\rm{\ref{CClaseT1.5.3}}}. $\fin$
\end{teorema}

\subsection{Resoluci\'on can\'onica}\label{SCo17.5.6}

Para $q\geq 1$ consideremos las $q$-tuplas $(\sigma_1,\ldots,\sigma_q)$ con
los $\sigma_i$ recorriendo el grupo $G$. Una $q$-tupla de estas se llama una
{\em $q$-celdas\index{celdas}} con v\'ertices $\sigma_1,\ldots,\sigma_q$. Las
$q$-celdas dan generadores libres de nuestros $G$-m\'odulos
\[
P_q=P_{-q-1}=\bigoplus_{\substack{(\sigma_1,\ldots,\sigma_q)\in\\
\underbrace{_{G\times\cdots\times G}}_{q}}}\zg(\sigma_1,\ldots,\sigma_q)
\cong {\ma Z}[\overbrace{G\times\cdots\times G}^{q+1}].
\]

Para $q=0$, $P_0=P_{-1}=\zg$ con el $1\in\zg$ como el generador de la
{\em celda nula o vac\'ia} $(\cdot)$. Sea $\varepsilon\colon P_0\lra {\ma Z}$,
$\varepsilon\big(\suma {\sigma}Ga\big)=\sum_{\sigma\in G}a_{\sigma}$ y
$\mu\colon {\ma Z}\lra P_{-1}$, $\mu(n)=n\N_G=\sum_{\sigma\in G}n\sigma$
(llamado {\em coaumentaci\'on\index{coaumentaci\'on}}).

Los mapeos $d_q$, est\'an dados por:
\begin{align}
d_0(1)&=\N_G,\label{d_0}\\
d_1(\sigma)&=\sigma-1,\label{d_1}\\
d_q(\sigma_1,\ldots,\sigma_q)&=\sigma_1(\sigma_2,\ldots,\sigma_q)\nonumber\\
&\quad+
\sum_{i=1}^{q-1}(-1)^i (\sigma_1,\ldots,\sigma_{i-1},\sigma_i\sigma_{i+1},
\sigma_{i+2},\ldots, \sigma_q)\nonumber\\
&\quad +(-1)^q(\sigma_1,\ldots,\sigma_{q-1}) \quad 
\text{para\ } q>1,\label{d_q}\\
d_{-1}(1)&=\sum_{\sigma\in G}(\sigma^{-1}(\sigma)-(\sigma))\quad \text{para\ }
q=-1,\label{d_-1}\\
d_{-q-1}(\sigma_1,\ldots,\sigma_q)&=\sum_{\sigma\in G}\sigma^{-1}(\sigma,
\sigma_1,\ldots,\sigma_q)\nonumber\\
&\quad +\sum_{\sigma\in G}\sum_{i=1}^q (-1)^i (\sigma_1,\ldots,\sigma_{i-1},
\sigma_i\sigma,\sigma^{-1},\sigma_{i+1},\ldots, \sigma_q) \nonumber\\
&\quad + \sum_{\sigma\in G}(-1)^{q+1}
(\sigma_1,\ldots,\sigma_q,\sigma)\quad \text{para\ } q>0.\label{d_-q-1}
\end{align}

Se puede verificar directamente que los mapeos $d_q$, $q\in{\ma Z}$ definidos
en (\ref{d_0}, \ref{d_1}, \ref{d_q}, \ref{d_-1}, \ref{d_-q-1})
hacen de la resoluci\'on completa una
sucesi\'on exacta. Sin embargo, m\'as adelante indicaremos como obtenerlos.

Se tiene que $P_{-q-1}\cong \Hom(P_q,{\ma Z})=\widehat {P}_q$, $q\geq 0$.
Ahora, si $A$ es un $G$-m\'odulo izquierdo, sea $A_q=\Hom_G(P_q, A)$. Los
elementos $A_q$, es decir, los $G$-homomorfismos $f\colon P_q\lra A$ se 
llaman las {\em{$q$-cadenas\index{cadenas@$q$-cadenas}}}
y $f$ est\'a determinado
por los valores $\{f(\sigma_1,\ldots,\sigma_q)\}$, $\sigma_1,\ldots,\sigma_q\in G$.

De la sucesi\'on exacta
\begin{gather*}
\cdots \tooi{{-2}}{} P_{-2}\tooi{{-1}}{}P_{-1}\tooi 0{}P_0\tooi 1{}P_1\tooi 2{}P_2
\tooi 3{}\cdots
\intertext{se obtiene la sucesi\'on}
\cdots \too{{-2}}{} A_{-2}\too{{-1}}{}A_{-1}\too 0{}A_0\too 1{}A_1\too 2{}A_2
\too 3{}\cdots
\end{gather*}
donde $\partial_q(f)=f\circ d_q$, $\partial_q\colon A_{q-1}\lra A_q$.
Se tiene $\partial_{q+1}\circ \partial_q=0$, por lo que $\im\partial_q
\subseteq \ker\partial_{q+1}$. Se definen los {\em 
$q$-cociclos\index{cocilos@$q$-cocilos}} $Z_q$ y los {\em
$q$-cofronteras\index{cofronteras@$q$-cofronteras}} $B_q$ por
\[
Z_q:=\ker\partial_{q+1},\quad B_q:=\im \partial_q,\quad q\in{\ma Z}.
\]

\begin{definicion}\label{CoD17.5.19}
Se define el {\em $q$-\'esimo grupo de cohomolog\'ia\index{grupo de
cohomolog\'ia@grupo de cohomolog\'ia}} (de Tate) por
\[
H^q(G,A)=Z_q/B_q.
\]
$H^q(G,A)$ es el grupo de cohomolog\'ia de dimensi\'on $q\in{\ma Z}$ del
$G$-m\'odulo $A$ o el $q$-\'esimo grupo de cohomolog\'ia de $G$ con
coeficientes en $A$.
\end{definicion}

Se tiene que todo elemento de $H^q(G,A)$ es representado por un mapeo
$f\colon \underbrace{G\times\cdots\times G}_q\lra A$. Adem\'as $A_0=
A_{-1}=\Hom_G(\zg,A)\cong A$.

De (\ref{d_0}, \ref{d_1}, \ref{d_q}, \ref{d_-1}, \ref{d_-q-1}) obtenemos
\begin{align*}
\partial_0(f)&=\N_G f, f\in A_{-1}=A,\\
(\partial_1f)(\sigma)&=\sigma(f)-f, f\in A_0=A,\\
(\partial_qf)(\sigma_1,\ldots,\sigma_q)&=\sigma_1f(\sigma_2,\ldots,\sigma_q)\\
&\quad+
\sum_{i=1}^{q-1}(-1)^i f(\sigma_1,\ldots,\sigma_{i-1},\sigma_i\sigma_{i+1},
\sigma_{i+2},\ldots, \sigma_q)\\
&\quad +(-1)^qf(\sigma_1,\ldots,\sigma_{q-1}) \quad 
\text{para\ } f\in A_{q-1}, q\geq 1,\\
(\partial_{-1}f)&=\sum_{\sigma\in G}(\sigma^{-1}f(\sigma)-f(\sigma))\quad \text{para\ }
f\in A_{-2}\\
(\partial_{-q-1}f)(\sigma_1,\ldots,\sigma_q)&=\sum_{\sigma\in G}\sigma^{-1}f(\sigma,
\sigma_1,\ldots,\sigma_q)\\
&\quad +\sum_{\sigma\in G}\sum_{i=1}^q (-1)^i f(\sigma_1,\ldots,\sigma_{i-1},
\sigma_i\sigma,\sigma^{-1},\sigma_{i+1},\ldots,\sigma_q)\\
&\quad + \sum_{\sigma\in G}(-1)^{q+1}f
(\sigma_1,\ldots,\sigma_q,\sigma)\quad \text{para\ } q\geq 0.
\end{align*}

As\'i los $q$-cociclos son los mapeos $f\colon G\times \cdots\times G\lra A$
con $\partial_{q+1} f=0$ y las $q$-cofronteras son los mapeos para los
cuales existe $g\in A_{q-1}$ con $f=\partial_q g$.

Ahora verificaremos las f\'ormulas
(\ref{d_0}, \ref{d_1}, \ref{d_q}, \ref{d_-1}, \ref{d_-q-1}), particularmente
el caso de $d_m$ con $m\leq 0$.
Consideremos 
\[
P_{-q-1}=\bigoplus\limits_{\substack{(\sigma_1,\ldots,\sigma_q)\in\\
G\times\cdots\times G}}\zg(\sigma_1,\ldots,\sigma_q).
\]
Enumeremos el sistema de $\zg$-generadores de $P_q$ consistente
de las $q$-celdas como $\{\Lambda_i\}$. Esto es, $\{\Lambda_i\}=
\{(\sigma_1,\ldots,\sigma_q)\mid \sigma_j\in G, 1\leq j\leq q\}$. Definimos
$\{\*{\Lambda_i}\}$ la base dual de $\{\Lambda_i\}$ por
\[
\*{\Lambda_i}(\sigma \Lambda_j)=\begin{cases}
1&\text{si $\sigma=1$ y $i=j$}\\
0&\text{de otra forma.}
\end{cases}
\]
Se tiene $\*{\Lambda_i}\in\Hom(P_q,{\ma Z})$. Se tiene que
$\{\*{\Lambda_i}\}$ es un sistema de generadores de $\Hom(P_q,
{\ma Z})$. Esto hace que los $G$-m\'odulos $\Hom(P_q,{\ma Z})$
y $P_q$ sean can\'onicamente isomorfos. Ahora bien
\[
P_{-q-1}=P_q \isomo_{\substack{\uparrow\\\Lambda_i\leftrightarrow 
\*{\Lambda_i}}}\Hom(P_q,{\ma Z}),
\]
$q\geq 0$ y ${\ma Z}\cong\Hom({\ma Z},{\ma Z})$.

Se tiene
\begin{scriptsize}
\[
\xymatrix{
0\ar@{->}[r]&{\ma Z}\ar@{->}[r]^{\mu}\ar@{-->}@/_1pc/[dd]&P_{-1}
\ar@{->}[r]^{d_{-1}}\ar@{-->}@/_1pc/[dd]&P_{-2}
\ar@{->}[r]^{d_{-2}}\ar@{-->}@/^1pc/[dd]&P_{-3}
\ar@{->}[r]^{d_{-3}}\ar@{-->}@/^1pc/[dd]&\cdots\\
&\ucong&\ \ucong \ \psi_0&\ucong&\ucong\\
0\ar@{->}[r]&\Hom({\ma Z},{\ma Z})\ar@{->}[r]\ar@{-->}@/_1pc/[dd]
&\Hom(P_0,{\ma Z})\ar@{->}[r]\ar@{-->}@/_1pc/[dd]&\Hom(P_1,{\ma Z})
\ar@{->}[r]\ar@{-->}@/^1pc/[dd]&\Hom(P_2,{\ma Z})\ar@{->}[r]
\ar@{-->}@/^1pc/[dd]&\cdots\\
&\ucong&\ \ucong \ \varphi_0&\ucong&\ucong\\
0\ar@{<-}[r]&{\ma Z}\ar@{<-}[r]^{\varepsilon}&P_{0}
\ar@{<-}[r]^{d_{1}}&P_{1}
\ar@{<-}[r]^{d_{2}}&P_{2}
\ar@{<-}[r]^{d_{3}}&\cdots
}
\]
\end{scriptsize}

La situaci\'on es como sigue para $q\geq 1$:
\[
\xymatrix{
P_q\ar@{->}[r]^{d_q}\ar@{->}[d]_{*}^{\cong}&P_{q-1}\ar@{->}[d]_{\cong}^{*}\\
\widehat P_q\cong\Hom(P_q,{\ma Z})\ar@{<-}[r]^{d^*_q}
\ar@{-->}@/_2.5pc/[dd]&
\Hom(P_{q-1},{\ma Z})\cong \widehat P_{q-1}
\ar@{-->}@/^2.5pc/[dd]\\
\ucong\  x^*_i\leftrightarrow x_i&\ucong \ y^*_i\leftrightarrow y_i\\
P_q=P_{-q-1}\ar@{<-}[r]^{d_{-q}}&P_{-q}=P_{q-1}
}
\]

Por tanto $d_q=d^*_q$ y $d^*_q(f)=f\circ d_q\in \Hom(P_q,{\ma Z})$
para $f\in\Hom(P_{q-1},{\ma Z})$.

Sea $\{x_i\}$ el sistema de generadores libres de $P_q$ como el $\zg$-m\'odulo
dado por las $q$-celdas: $\{x_i\}_{i\in I}=\{(\sigma_1,\ldots,\sigma_q)\mid
\sigma_i\in G\}$. Entonces el correspondiente sistema de generadores de
$P_q$ como ${\ma Z}$-m\'odulo es $\{\sigma x_i\}_{\substack{\sigma\in G\\
i\in I}}$. Se tiene que $\{\sigma x_i\}_{\substack{\sigma\in G\\
i\in I}}=\{\sigma(\sigma_1,\ldots,\sigma_q)\mid \sigma, \sigma_i\in G\}$.

Sean $\{y_j\}_{i\in J}$ y $\{\sigma y_j\}_{\substack{\sigma\in G\\
j\in J}}$ las bases correspondientes de $P_{q-1}$.

Sean
\begin{align*}
d_qx_t&=\sum_{\sigma\in G, s}a_{t,\sigma,s}\sigma y_s\quad \text{con\ }
a_{t,\sigma,s}\in {\ma Z},\\
d^*_qy^*_l&=\sum_{\tau\in G, u}b_{l,\tau,u}\tau x^*_u \quad \text{con\ }
b_{l,\tau,u}\in{\ma Z}.
\intertext{Se tiene}
(d^*_qy^*_l)(\theta x_m)&=\sum_{\tau,u}b_{l,\tau,u}(\tau x^*_u)(\theta x_m)
=\sum_{\tau,u}b_{l,\tau,u}\tau(x^*(\tau^{-1}\theta(x_m)))=b_{l,\theta,m}.
\intertext{Por otro lado}
(d^*_qy^*_l)(\theta x_m)&= y^*_ld_q(\theta x_m)=y^*_l(\theta(d_q(x_m)))=
y^*_l(\theta\big(\sum_{\sigma,s}a_{m,\sigma,s}\sigma y_s\big))\\
&=y^*_l\big(\sum_{
\sigma,s}a_{m,\sigma,s}(\theta \sigma)(y_s)\big)=\sum_{\sigma,s}
a_{m,\sigma,s}y^*_l((\theta\sigma)(y_s))=a_{m,\theta^{-1},l}.
\end{align*}

Por tanto se tiene
\begin{gather*}
(d^*_qy^*_l)(\theta x_m)=b_{l,\theta,m}=a_{m,\theta^{-1},l}.
\intertext{Poniendo $l ``="(\sigma_1,\ldots.\sigma_{q-1})$ y $m``="(\sigma'_1,
\ldots,\sigma'_q)$, se tiene}
b_{(\sigma_1,\ldots,\sigma_{q-1}),\theta,(\sigma'_1,\ldots,\sigma'_q)}=
a_{(\sigma'_1,\ldots,\sigma'_q),\theta^{-1},(\sigma_1,\ldots,\sigma_{q-1})}.
\end{gather*}

Ahora bien, de acuerdo a la ecuaci\'on (\ref{d_q}), los coeficientes no cero
para las $(q+1)$-celdas son
\begin{align}
a_{(\mu_1,\ldots,\mu_{q+1}),\mu_1,(\mu_2,\ldots,\mu_{q+1})}&=1,\label{d_q(1)}\\
a_{(\mu_1,\ldots,\mu_{q+1}),1,(\mu_1,\ldots,\mu_{i-1},\mu_i\mu_{i+1},
\mu_{i+2},\ldots,\mu_{q+1})}&=(-1)^i,\quad 1\leq i\leq q,\label{d_q(2)}\\
a_{(\mu_1,\ldots,\mu_{q+1}),1,(\mu_1,\ldots,\mu_{q})}&=(-1)^{q+1}.\label{d_q(3)}
\end{align}

Para (\ref{d_q(1)}), poniendo $(\sigma_1,\ldots,\sigma_q)=(\mu_2,\ldots,
\mu_{q+1})$ y $\mu_1=\tau^{-1}$, se tiene que
los $b's$ correspondientes son $\big\{
b_{(\sigma_1,\ldots,\sigma_{q}),\tau^{-1},(\tau,\sigma_1,\ldots,\sigma_{q})},
\tau\in G\big\}$.
Para (\ref{d_q(2)}) ponemos $\theta=1$ y $(\sigma_1,\ldots,\sigma_q)=
(\mu_1,\ldots,\mu_{i-1},\mu_i\mu_{i+1},\mu_{i+2},\ldots,\mu_{q+1})$,
se obtiene que los $b's$ son 
\[
\big\{
b_{(\sigma_1,\ldots,\sigma_{q}),1,(\sigma_1,\ldots,\sigma_{i-1},\sigma_i
\tau,\tau^{-1},\sigma_{i+1},\ldots,\sigma_{q})},
\tau\in G\big\}.
\]

Para (\ref{d_q(3)}) se obtiene que los $b's$ correspondientes son
\[
\big\{b_{(\sigma_1,\ldots,\sigma_{q}),1,(\sigma_1,\ldots,\sigma_{q},\tau)},
\tau\in G\big\}.
\]
De esto, se obtiene (\ref{d_-q-1}). Las otras f\'ormulas se obtienen de manera
similar.

\subsubsection{Cambio de dimensi\'on}\label{SCambio de dimension}

Sean $I_G$ el ideal aumentaci\'on y $J_G=\zg/{\ma Z}\N_G$. Se tienen las
sucesiones exactas
\begin{gather*}
0\lra I_G\lra \zg\stackrel{\varepsilon}{\lra}{\ma Z}\lra 0,\\
0\lra {\ma Z}\stackrel{\Lambda}{\lra}\zg\lra J_G\lra 0.
\end{gather*}

Puesto que ${\ma Z},\zg, I_G$ y $J_G$ son todos ${\ma Z}$-m\'odulos libres,
para un $G$-m\'odulo $A$, se tienen las siguientes sucesiones $G$-exactas
(Proposici\'on \ref{CoP17.5.7}):
\begin{gather}\label{I_GJ_G}
\begin{array}{cccccccc}\
0&\lra &I_G\otimes A&\lra &\zg\otimes A&\lra& A&\lra 0,\\
0&\lra& A&\lra& \zg\otimes A&\lra& J_G\otimes A&\lra 0,
\end{array}
\end{gather}
puesto que $A\cong {\ma Z}\otimes A$.

Como $\zg$ es un $G$-m\'odulo inducido, $\zg$ es cohomol\'ogicamente
trivial y por el Teorema \ref{sucesion exacta} se tienen isomorfismos
\begin{gather*}
\mu\colon H^{q-1}(H,J_G\otimes A)\stackrel{\cong}{\lra} \co qHA,\\
\mu^{-1}\colon H^{q+1}(H,I_G\otimes A)\stackrel{\cong}{\lra} \co qHA,
\end{gather*}
para todo $q\in{\ma Z}$ y para todo subgrupo $H<G$, donde
$\mu$ y $\mu^{-1}$ son los mapeos de conexi\'on.

Continuamos el proceso, y para cualquier $m\in {\ma Z}$ definimos
\begin{gather*}
A^m:=\overbrace{J_G\otimes\cdots\otimes J_G}^m\otimes A \quad \text{para\ } m\geq 0,\\
A^m:=\underbrace{I_G\otimes\cdots\otimes I_G}_{-m}\otimes A \quad \text{para\ } m\leq 0.
\end{gather*}

Se tienen isomorfismos
\begin{gather}
\co {{q-m}}H{A^m}\xrightarrow[\mu,\mu^{-1}]{\cong} \co
{{q-(m-1)}}HA^{m-1} \xrightarrow[\mu,\mu^{-1}]{\cong}\cdots \xrightarrow[\mu,\mu^{-1}]{\cong}
\co qHA,\nonumber
\intertext{esto es,}
\mu^m\colon \co {{q-m}}H{A^m}\stackrel{\cong}{\lra}\co qHA,\quad m\in{\ma Z}.\label{cambiodim}
\end{gather}

Se usa el isomorfismo (\ref{cambiodim}) para deducir propiedades de grupos de cohomolog\'ia
de dimensi\'on $q$ a propiedades an\'alogas para cohomolog\'ia en dimensiones superiores
o inferiores. En particular, esta t\'ecnica nos permitir\'a reducir muchas definiciones y
demostraciones al caso $0$-dimensional que es m\'as f\'acil de trabajar. Este m\'etodo
se llama {\em cambio de dimensi\'on\index{cambio de dimensi\'on}}. Se tiene que
\[
\co qGA\cong \co 0G{A^q} =(A^q)^G/\N_G A^q.
\]

El siguiente resultado, es un ejemplo de como usar el cambio de dimensi\'on.

\begin{teorema}\label{mH=0}
Sea $G$ un grupo finito de $n$ elementos. Sea $A$ un $G$-m\'odulo. Entonces
$n\co qHA=0$ para toda $q\in{\ma Z}$.
\end{teorema}

\begin{proof}
Si $q=0$, $\co 0GA=A^G/\N_G A$, y si $a\in A^G$, $na=\sum_{\sigma\in G}
\sigma a\in\N_G A$. Por tanto $n\co 0GA=0$. Se sigue que $n\co qGA=
n\co 0G{A^q}=0$ para toda $q\in{\ma Z}$.
$\fin$
\end{proof}

\begin{definicion}\label{divisible}
Si $X$ es un grupo abeliano, $X$ se llama {\em un\'ivocamente o 
\'unicamente divisible\index{unicamente divisible@\'unicamente divisible}}
si para toda $n\in{\ma N}$ y para toda $a\in X$, la ecuaci\'on $nx=a$
tiene una \'unica soluci\'on $x\in X$ (formalmente $x=\frac an$). 
Se tiene que $X$ es un\'ivocamente divisible si para toda $n\in{\ma N}$
el mapeo $X\lra X$, $x\longmapsto nx$ es un isomorfismo de grupos,
\end{definicion}

\begin{ejemplo} 
El campo de los n\'umeros racionales ${\ma Q}$ es
un\'ivocamente divisible. M\'as generalmente, cualquier campo de
caracter\'istica $0$ es un\'ivocamente divisible.
\end{ejemplo}

\begin{corolario}\label{unicamente divisible}
Si $A$ es un $G$-m\'odulo que es un\'ivocamente divisible, entonces $A$
es cohomol\'ogicamente trivial.
\end{corolario}

\begin{proof}
Se tiene que $n\Id:A\lra A$, $a\longmapsto na$ es una biyecci\'on para toda
$n\in{\ma Z}$. Por tanto, si $n=|H|$, $\co qHA=n\co qHA=0$ para todo $q\in
{\ma Z}$ y para todo subgrupo $H<G$.
$\fin$
\end{proof}

Aplicando el Corolario \ref{unicamente divisible} a ${\ma Q}$ y a la suceci\'on
exacta de $G$-m\'odulos triviales $0\lra {\ma Z}\lra {\ma Q}\lra {\ma Q}/{\ma Z}
\lra 0$, obtenemos que $\co qG{{\ma Q}/{\ma Z}}\cong \co {{q+1}}G{\ma Z}$
para toda $q\in{\ma Z}$.

En particular, $\co 2G{\ma Z}\cong \co 1G{{\ma Q}/{\ma Z}}\cong \Hom(G,
{\ma Q}/{\ma Z})\cong \widehat G\cong G/G'=\abe G$, donde
$\chi(G)=\widehat G=\Hom(G,{\ma Q}/{\ma Z})$ es el {\em grupo de
caracteres\index{grupo de caracteres}} de $G$.

Por otro lado de la sucesi\'on $0\lra I_G\lra \zg\lra {\ma Z}\lra 0$, obtenemos
$\co {{-2}}G{\ma Z}\cong \co {{-1}}G{I_G}\cong I_G/I_G^2\cong G/G'=\abe G$.
Este \'ultimo isomorfismo est\'a dado por:
\[
(\sigma -1)+I_G^2\lra \sigma G'\quad {\text para\ }\sigma \in G.
\]
En particular $\co {{-2}}G{\ma Z}\cong \abe G\cong \co 2G{\ma Z}$.

Como hicimos notar antes $\co 1G{\ma Z}=\Hom(G,
{\ma Z})=0$ y $\co {{-1}}G{\ma Z}=
\frac{\ker \N_G}{I_G{\ma Z}}=0$. Finalmente $\co 0G{\ma Z}={\ma Z}/
n{\ma Z}=C_n$, donde $|G|=n$.

De hecho tenemos 
\[
\co {{-q}}G{\ma Z}\cong \chi\big(\co qG{\ma Z}\big)=
\Hom(\co qG{\ma Z},{\ma Q}/{\ma Z}), \quad \text{para\ } q>0,
\]
ver \cite[Corollary 4-4-7]{Wei69}.

\subsection{Cambio de grupo}\label{SCambiogrupo}

Sea $\varphi\colon G'\lra G$ un homomorfismo de grupos finitos (aqu\'i $G'$
denota a un grupo finito arbitrario y no al conmutador de $G$). Si $A$ es
un $G$-m\'odulo, $A$ puede hacerse un $G'$-m\'odulo de la siguiente forma.
Para $\tau'\in G'$ y para $a\in A$, se define $\tau'\circ a:=\varphi(\tau')\circ a$.
En caso de ser necesario, el $G'$-m\'odulo $A$ se denotar\'a por $\varphi^* A$.

Si $A$ y $B$ son dos $G$-m\'odulos y $f\colon A\lra B$ es un 
$G$-homomorfismo, entonces $f$ tambi\'en es un $G'$-homomorfismo
$f\colon \varphi^* A\lra \varphi^* B$,
\[
f(\tau' a)=f(\varphi(\tau') a)=\varphi(\tau')f(a)=\tau'f(a).
\]

M\'as generalmente, si $\varphi\colon G'\lra G$ es un homomorfismo
de grupos, si $A'$ es un $G'$-m\'odulo y si $g\colon A\lra A'$ es un mapeo
aditivo, esto es, $g$ es un homomorfismo de grupos abelianos o $g$
es un ${\ma Z}$-homomorfismo, entonces $\varphi$ y $g$ se llaman
{\em compatibles\index{compatible!homomorfismos}} si
\[
g(\varphi(\tau')a)=\tau'g(a)\quad \text{para toda\ } \tau'\in G'.
\]
En otras palabras, si $g$ es un $G'$-homomorfismo entre 
$\varphi^* A$ y $A'$.

Consideremos cualesquiera dos complejos $P'$ y $P$, esto es,
resoluciones completas para $G'$ y para $G$. Podemos suponer
que tanto $P'$ como $P$ son las resoluciones can\'onicas. Entonces
$P'$ y $P$ son resoluciones proyectivas, de hecho libres.
\[
\xymatrix{
P':\cdots\ar@{->}[r]^{\partial'_3}&P'_2 \ar@{->}[r]^{\partial'_2}\ar@{-->}[dd]_{\Lambda_2}&P'_1
\ar@{-->}[ddl]_{\Delta_1}\ar@{-->}[dd]_{\Lambda_1}\ar@{->}[r]^{\partial'_1}
&P'_0\ar@{->}[rr]^{\partial'_0}\ar@{->}[dr]^{\varepsilon'}\ar@{-->}[dd]^{
\Lambda_0}\ar@{-->}[ddl]_{\Delta_0}&&P'_{-1}\ar@{->}[r]^{\partial'_{-1}}
\ar@{-->}[dd]_{\Lambda_{-1}}&P'_{-2}\ar@{-->}[ddl]_{\Delta_{-2}}
\ar@{->}[r]^{\partial'_{-2}}\ar@{-->}[dd]^{\Lambda_{-2}}&\cdots
\\
&&&& {\ma Z}\ar@{->}[ru]^{\mu'}\ar@{->}[rd]_{\mu}
\\
P:\cdots\ar@{->}[r]^{\partial_3}&P_2\ar@{->}[r]_{\partial_2}&P_1\ar@{->}[r]_{\partial_1}
&P_0\ar@{->}[rr]_{\partial_0}\ar@{->}[ru]_{\varepsilon}&&P_{-1}
\ar@{->}[r]_{\partial_{-1}}&P_{-2}\ar@{->}[r]_{\partial_{-2}}&\cdots
}
\]

Se tiene que $\partial', \varepsilon',\mu'$ son $G'$-homomorfismos; $\partial,
\varepsilon,\mu$ son $G$-homomorfismos y por tanto, cuando los m\'odulos
$P_n$ son vistos como $G'$-m\'odulos, $\partial,\varepsilon,\mu$ son 
tambi\'en $G'$-homomorfismos.

\begin{teorema}\label{T17.5.6.1}
Para cualquier homomorfismo $\varphi\colon G'\lra G$ existen $G'$-homomorfismos
$\Lambda_n\colon P'_n\lra P_n$, para $n\geq 0$ tales que $\varepsilon\circ
\Lambda_0=\varepsilon'$ y $\partial_{n+1}\circ \Lambda_{n+1}=\Lambda_n\circ
\partial'_{n+1}$.
\end{teorema}

\begin{proof}
Se tiene que $P'_0$ es $G'$-libre y $\varepsilon$ es suprayectiva, por lo que
existe un $G'$-homomorfismo $\Lambda_0\colon P'_0\lra P_0$ tal que
$\varepsilon\circ \Lambda_0=\varepsilon'$, de hecho, si $\{\xi_i\}$ es 
$G'$-base de $X'_0$ y $\varepsilon'(\xi_i)=n_i\in {\ma Z}$, como $\varepsilon$
es sobre, existe $\delta_i\in P_0$ con $\varepsilon(\delta_i)=n_i$. Se define
$\Lambda_0(\xi_i):=\delta_i$.

Se tiene $\partial_0\Lambda_0\partial'_1=\mu\varepsilon\Lambda_0
\partial'_1=\mu\varepsilon'\partial'_1=0$.

Supongamos que $\Lambda_n$ ha sido definido para $n\geq 0$ y que
$\partial_n\Lambda_n\partial'_{n+1}=0$. Puesto que $P'_{n+1}$ es
$G'$-libre, $P_{n+1}$ es $G'$-inducido. Sea $P_{n+1}=\sum_{\sigma
\in G'} \sigma X$ y sea $\pi\colon P'_{n+1}\lra X$ la proyecci\'on, esto es,
si $a\in P'_{n+1}$, $a=\sum_{\sigma'\in G'}\sigma' b_{\sigma'}$ con $
b_{\sigma'}\in X$, entonces $\pi(a)=b_1$. Para $\tau\in G'$, se tiene
\begin{gather*}
\begin{align*}
(\tau'\circ \pi)(a)&=\tau'\pi((\tau')^{-1}a)=\tau'\pi\big(\sum_{\sigma'\in G'}
(\tau')^{-1}\sigma' b_{\sigma'}\big)\\
&=\tau'\pi\big(\sum_{\theta'\in G'}
\theta' b_{\tau'\theta'}\big)=\tau' b_{\tau'}.
\end{align*}
\intertext{Esto es, $\tau'\circ \pi$ es la proyecci\'on de $P'_{n+1}$ en $\tau' X$.
Por tanto}
(\sigma'\circ \pi)\circ(\theta'\circ \pi)(a)=(\sigma'\circ \pi)(\theta'b_{\theta'})=
\begin{cases} 0&\text{si $\sigma'\neq\theta'$}\\
\theta'b_{\theta'}&\text{si $\sigma'=\theta'$} \end{cases} =
\delta_{\sigma',\theta'}(\theta'\circ \pi)(a).
\intertext{Se sigue que}
(\sigma'\circ \pi)\circ(\theta'\circ \pi)=\delta_{\sigma',\theta'}(\theta'\circ \pi).
\end{gather*}

Por otro lado, se tiene que $\Id_{P'_{n+1}}(a)=a=\sum_{\sigma'\in G}
\sigma' b_{\sigma'}=\sum_{\sigma'\in G'}(\sigma'\circ \pi)(a)$. Por tanto
\[
\Id_{P'_{n+1}}=\sum_{\sigma'\in G'}\sigma'\circ \pi.
\]
De esta forma, existe $\pi'_{n+1}(=\pi)\in \Hom(P'_{n+1},P'_{n+1})$ con
$\Id_{P'_{n+1}}=\sum_{\sigma'\in G'}\sigma'\circ \pi'_{n+1}$. Sea
$\Lambda_{n+1}:=\sum_{\sigma'\in G'}\sigma'(D_n\Lambda_n\partial'_{
n+1}\pi'_{n+1})$ donde $D_n\colon P_n\lra P_{n+1}$ satisface
$D_{n-1}\circ \partial_n+\partial_{n+1}\circ D_n=\Id_{P_n}$ 
(Proposici\'on \ref{CoP17.5.10}).
\[
\xymatrix{
P'_{n+1}\ar@{->}[r]^{\pi'_{n+1}}\ar@{-->}[dr]_{\Lambda_{n+1}}&
P'_{n+1}\ar@{->}[r]^{\partial'_{n+1}}&P'_n\ar@{->}[d]^{\Lambda_n}\\
&P_{n+1}\ar@{<-}[r]_{\ D_n}&P_n
}
\]

Entonces $\Lambda_{n+1}$ es un $G'$-homomorfismo y se tiene
\begin{align*}
\partial_{n+1}\Lambda_{n+1}&=\sum_{\sigma'\in G'} \partial_{n+1}(
\sigma'(D_n\Lambda_n\partial'_{n+1}\pi'_{n+1}))\\
&=\sum_{\sigma'\in G'}
(\sigma'(\partial_{n+1} D_n\Lambda_n\partial'_{n+1}\pi'_{n+1}))\\
&\igual_{\substack{\uparrow\\ \partial_{n+1}D_n=\\1-D_{n-1}\partial_n}}
\sum_{\sigma'\in G'}\sigma'(\Lambda_n\partial'_{n+1}\pi'_{n+1})-
\sum_{\sigma'\in G'}\sigma'(D_{n-1}
\underbrace{\partial_n\Lambda_n\partial'_{n+1}}_{0}\pi'_{n+1})\\
&=\Lambda_n\partial'_{n+1}\sum_{\sigma'\in G'}\sigma'\pi'_{n+1}
-0=\Lambda_n\partial'_{n+1} \Id_{P'_{n+1}}=\Lambda_n\partial'_{n+1}.
\end{align*}

Esto es, $\partial_{n+1}\Lambda_{n+1}=\Lambda_n\partial'_{n+1}$.
Finalmente 
\begin{gather*}
\partial_{n+1}\Lambda_{n+1}\partial'_{n+2}=
\Lambda_n\partial'_{n+1}\partial'_{n+2}=0.
\tag*{$\fin$}
\end{gather*}
\end{proof}

\begin{observacion}\label{O17.5.6.2}
El proceso del Teorema \ref{T17.5.6.1} 
no puede ser continuado como est\'a en general para $n<0$.
Sin embargo, cuando $\varphi\colon G'\lra G$ es un monomorfismo, el
proceso puede ser continuado para $n<0$.
\end{observacion}

\begin{teorema}\label{T17.5.6.3}
Si $\varphi\colon G'\lra G$ es un monomorfismo de grupos, entonces
existen $G'$ homomorfismos $\Lambda_n\colon P'_n\lra P_n$ tales
que $\varepsilon\circ \Lambda_0=\varepsilon'$ y $\partial_{n+1}\circ
\Lambda_{n+1}=\Lambda_n\circ \partial'_{n+1}$ para toda $n\in{\ma Z}$.
\end{teorema}

\begin{proof}
Se tiene $\Lambda_0$ tal que $\varepsilon\circ \Lambda_0=
\varepsilon'$ y $\partial_0\Lambda_0\partial'_1=0$. Supongamos
inductivamente que para $n\leq 0$ se tiene $\Lambda_n$ con
$\partial_n\Lambda_n\partial'_{n+1}=0$. Ahora $G'$ se puede 
considerar como un subgrupo de $G$, por lo que $P_n$ es $G'$-libre
y en particular es $G'$-inducido. Por el argumento del Teorema
\ref{T17.5.6.1}, existe $\pi_n\in\Hom(P_n,P_n)$ tal que $S'(\pi_n):=
\sum_{\sigma\in G}\sigma \pi_n=\Id_{P_n}$.

Sea $\Lambda_{n-1}:=S'(\pi_{n-1}\partial_n\Lambda_n D'_{n-1})=
\sum_{\sigma\in G}\sigma(\pi_{n-1}\partial_n\Lambda_n D'_{n-1})$,
donde $D'_{n+1}\colon P'_n\lra P'_{n+1}$ satisface $D'_{n-1}\circ \partial'_n
+\partial'_{n+1}\circ D'_n=\Id_{P'_n}$.
Entonces $\partial_n\circ \Lambda_n=\Lambda_{n-1}\circ \partial'_n$
y $\partial_{n-1}\Lambda_{n-1}\partial'_n=0$. Adem\'as $\mu=
\Lambda_{-1}\circ \mu'$.
$\fin$
\end{proof}

\begin{teorema}\label{T17.5.6.4}
Sea $\varphi\colon G'\lra G$ un homomorfismo de grupos finitos.
Sean $A$ y $A'$, un $G$-m\'odulo y un $G'$-m\'odulo respectivamente,
y sea $g\colon A\lra A'$ compatible con $\varphi$, es decir, $
g(\varphi(\tau')a)=\tau'g(a)$ para toda $\tau'\in G'$ y para toda $a\in A$.

Existe una \'unica familia de homomorfismos de grupos
\[
(\varphi, g)_n\colon \co nGA\lra \co n{G'}{A'},
\]
 para $n\geq 0$ (y para
$n\in {\ma Z}$ si $\varphi$ es inyectiva) que satisface:
\las
\item $(1,1)_n=\Id_n$ para $n$.\label{T17.5.6.4(1)}
\item Si $\varphi'\colon G''\lra G'$ es un homomorfismo de grupos finitos,
$A''$ es un $G''$-m\'odulo y si $g'\colon A'\lra A''$ es compatible con
$\varphi'$, esto es, $g'(\varphi'(\tau'')a')=\tau''g'(a')$ para toda $\tau''
\in G''$, $a'\in A'$, entonces $(\varphi\circ\varphi', g'\circ g)_n=(\varphi',
g')_n\circ (\varphi,f)_n$ para toda $n\geq 0$ (a veces $n\in{\ma Z}$).
\label{T17.5.6.4(2)}
\item Si $g_1\colon A\lra A'$ es compatible con $\varphi$, entonces
$(\varphi,g+g_1)_n=(\varphi,g)_n+(\varphi,g_1)_n$.\label{T17.5.6.4(3)}
\end{list}
\end{teorema}

\begin{proof}
Sea $\Lambda\colon P'\lra P$ un $G'$-homomorfismo conmutando con
fronteras: 
\begin{gather*}
\xymatrix{P'_{n+1}\ar@{->}[rr]^{\partial'_{n+1}}\ar@{->}[d]_{
\Lambda_{n+1}}&&P'_n\ar@{->}[d]^{\Lambda_n}\\
P_{n+1}\ar@{->}[rr]_{\partial_{n+1}}&&P_n}
\\
\partial_{n+1}\circ \Lambda_{
n+1}=\Lambda_n\circ\partial'_{n+1}
\\
\xymatrix{
\Hom_G(P'_n,A')\ar@{->}[rr]^{(\partial'_{n+1})^*}\ar@{->}[d]_{(\varphi,g)_n}^{
\mu_n}&&\Hom_G(P'_{n+1},A')\ar@{->}[d]^{\mu_{n+1}}\\
\Hom_G(P_n,A)\ar@{->}[rr]_{\partial^*_{n+1}}&&\Hom_G(P_{n+1},A)}
\end{gather*}

Si $\psi\in\Hom_G(P'_n,A')$, 
\[
\xymatrix{P'_n\ar@{->}[r]^{\psi}\ar@{<-}[d]_{
\partial'_{n+1}}&A'\ar@{<--}[dl]^{\partial^*_{n+1}(\psi)=\psi\circ \partial'_{n+1}}
\\ P_{n+1}}\qquad \qquad
\xymatrix{P'_n\ar@{->}[rr]^{(\varphi,g)_n(\xi)}\ar@{->}[d]_{\Lambda_n}&&
A'\ar@{<-}[d]^g\\ P_n\ar@{->}[rr]_{\xi}&&A
}
\]

Sea $\xi\in\Hom_G(P_n,A)$, $(\varphi,g)_n(\xi)=g\circ\xi\circ\Lambda_n$, 
$\partial^*_{n+1}(\xi)=\xi\circ \partial_{n+1}$.

Se tiene $\xi\in\ker \*{\partial_{n+1}}$ por lo que $\xi\circ \partial_{n+1}=0$.
Se sigue que
\begin{align*}
\*{(\partial'_{n+1})}((\varphi,g)_n(\xi))&=(\varphi,g)_n(\xi)\circ \partial'_{n+1}=
g\circ \xi\circ \Lambda_n \circ \partial'_{n+1}\\
&=g\circ \underbrace{\xi\circ 
\partial_{n+1}}_{\substack{\uigual\\0}}\circ \Lambda_{n+1}=0
\end{align*}
Por tanto $(\varphi,g)_n(\xi)\in\ker \*{\partial_{n+1}}$.

Ahora si $\xi\in\im \*{\partial_n}$, $\xi=\*{\partial_n}(\nu)$ con $\nu\in
\Hom_G(P_{n-1},A)$ y
\begin{align*}
(\varphi,g)_n(\xi)&= g\circ \xi\circ \Lambda_n=g\circ (\*{\partial_n}(\nu))\circ
\Lambda_n=g\circ \nu\circ \partial_n\circ \Lambda_n\\
&=g\circ \nu\circ \Lambda_{n-1}
\circ \partial'_n=\*{(\partial'_n)}(g\circ \nu\circ \Lambda_{n-1})\in \im\*{(
\partial'_n)}.
\end{align*}

Para la unicidad, sean $\Lambda,\Lambda'\colon P'\lra P$ dos $G'$-homomorfismos
que conmutan con fronteras. Sea $\Phi:=\Lambda-\Lambda'$. Se construir\'an 
$G'$-homomorfismos $\Lambda_n\colon P'_n\lra P_{n+1}$ tales que 
\begin{gather}\label{Homo}
\Phi_n=\partial_{n+1}\Delta_n+\Delta_{n-1}\partial'_n.
\end{gather}

Sean, en cualquier caso, $n\geq 0$ o $n\in{\ma Z}$
\begin{gather*}
\Lambda_0=S'(D_0\Phi_0\pi'_0)=\sum_{\tau'\in G'}\tau'D_0 \Phi_0\pi'_0
\quad \text{y}\quad \Lambda_{-1}=0.\\
\intertext{Entonces}
\begin{align*}
\partial_1\Lambda_0&=S'(\underbrace{\partial_1D_0}_{\substack{\uigual\\
1-D_{-1}\partial_0}}\Phi_0 \pi'_0)=S'(\Phi_0\pi'_0)-S'(D_{-1}\partial_0\Phi_0
\pi'_0)\\
&=\Phi_0\underbrace{S'(\pi'_0)}_{\substack{\uigual\\1}}-S'(D_{-1}
\mu\underbrace{\varepsilon\Phi_0}_{\substack{\uigual\\ 0}}\pi'_0)=\Phi_0,
\end{align*}
\end{gather*}
por lo que (\ref{Homo}) es v\'alido para $n=0$.

Para $n>0$, procedemos por inducci\'on. Supongamos (\ref{Homo}) cierta
$n$ y sea 
\begin{gather*}
\Delta_{n+1}:=S'(D_{n+1}(\Phi_{n+1}-\Delta_n \partial'_{n+1})\pi'_{n+1}), \quad n>0.
\intertext{Ahora}
\begin{align*}
\partial_{n+1}(\Phi_{n+1}-\Delta_n\partial'_{n+1})&=
\partial_{n+1}\Phi_{n+1}-
\partial_{n+1}\Delta_n \partial'_{n+1}\\
&=\Phi_n\partial'_{n+1}-\partial_{n+1}\Delta_n\partial'_{n+1}
\underbracket[0pt]{=}_{\substack{\uparrow\\ \text{(\ref{Homo}) para $n$}}}\\
&=\partial_{n+1}\Delta_n \partial'_{n+1}+\Delta_{n-1}
\underbrace{\partial'_n\partial'_{n+1}}_{0}-\partial_{n+1}\Delta_n\partial'_{n+1}=0.
\end{align*}
\end{gather*}

Se sigue que
\begin{align*}
\partial_{n+2}\Delta_{n+1}&=S'(\partial_{n+2}D_{n+1}(\Phi_{n+1}-\Delta_n\partial'_{
n+1})\pi'_{n+1})\\
&\underbracket[0pt]{=}_{\substack{\uparrow\\ \partial_{n+2}D_{n+1}=\\1-D_n\partial_{n+1}}}
S'((\Phi_{n+1}-\Delta_n\partial'_{n+1})\pi'_{n+1})\\
&\qquad \qquad -S'(D_n\underbrace{
\partial_{n+1}(\Phi_{n+1}-\Delta_n\partial'_{n+1}}_{\substack{\uigual \\0}})\pi'_{n+1})\\
&=S'(\Phi_{n+1}\pi'_{n+1})-S'((\Delta_n \partial'_{n+1})\pi'_n)\\
&=\Phi_{n+1}S'(\pi'_{n+1})-
\Delta_n\partial'_{n+1}S'(\pi'_n)\\
&=\Phi_{n+1}\cdot 1-\Delta_n\partial'_{n+1}\cdot 1=\Phi_{n+1}-\Delta_n\partial'_{n+1}.
\end{align*}

Se sigue que (\ref{Homo}) se cumple para $n+1$.

Por tanto, si $\xi\in\ker \partial'_n$, entonces se tiene 
$\Lambda_n(\xi)-\Lambda'_n(\xi)=\Phi_n(\xi)=\partial_{n+1}\Delta_n(\xi)+\Delta_{n-1}
\partial'_n(\xi)=\partial_{n+1}\Delta_n(\xi)\in\im \partial_{n+1}$, de donde,
$\Lambda_n(\xi)\equiv\Lambda_n(\xi)\bmod (\im \partial_{n+1})$ y por lo tanto
$(\varphi,g)_n$ no depende de $\Lambda$.

(\ref{T17.5.6.4(1)}), (\ref{T17.5.6.4(2)}) y (\ref{T17.5.6.4(3)}) se verifican directamente.
$\fin$
\end{proof}

\begin{observacion}\label{O17.5.6.5}
Todo el desarrollo anterior, es con respecto a la cohomolog\'ia de Tate. Ahora, para
$H^0$ no de Tate, esto es, $\co 0GA=A^G$, tambi\'en se satisface
\[
(\varphi,g)_0\colon \co 0GA=A^G\lra \co 0{G'}{A'}=(A')^{G'}:
\]
Se tiene para $\tau'\in G'$, $a\in A$, $\tau'g(a)=g(\varphi(\tau'))=g(a)$, por lo
tanto $(\varphi,g)_0(A^G)\subseteq (A')^{G'}$.
\end{observacion}

\begin{ejemplo}[\text{Restricci\'on}\index{Restricci\'on}]\label{Ej17.5.6.6(1)}
Sean $G$ un grupo finito y $H< G$ un subgrupo. Entonces, si
$\varphi=i\colon H\lra G$ es la inyecci\'on natural, para $A$ un $G$-m\'odulo, $A$
es un $H$-m\'odulo con la misma acci\'on. Entonces
\[
(i,\Id_A)_n:=\res_n\colon \co nGA\lra \co nHA,
\]
se llama {\em restricci\'on\index{restricci\'on}} y como $i$ es inyectiva, el mapeo
$\res_n$ es v\'alido para toda $n\in{\ma Z}$.

A nivel de cadenas $P_H: P_{H,n}={\ma Z}[H^{n+1}]$, $n\geq 0$ y $P_G:
P_{G,n}={\ma Z}[G^{n+1}]$, dada una $q$-cadena de $G$: ($q\geq 0$),
\begin{gather*}
x\colon \underbrace{G\times G\times\cdots\times G}_{q}\lra A,\\
\res x=y\colon \underbrace{H\times H\times\cdots\times H}_{q}\lra A,\quad
y=x|_{H\times\cdots\times H}.
\end{gather*}
\end{ejemplo}

Tambi\'en se tiene:

\begin{proposicion}\label{P17.5.6.6(1)-1}
Sean $A$ y $B$ dos $G$-m\'odulos y $H<G$ un subgrupo de $G$, $G$
un grupo finito. Sea $f\colon A\lra B$ un $G$-homomorfismo. Entonces el
diagrama
\[
\xymatrix{
\co qGA\ar@{->}[rr]^{H^q{(f)}}\ar@{->}[d]_{\res_q}&&\co qGB\ar@{->}[d]^{
\res_q}\\ \co qHA\ar@{->}[rr]_{H^q(f)}&&\co qHB
}
\]
es conmutativo. Notemos que un $G$-homomorfismo $f\colon A\lra B$,
es tambi\'en un $H$-homomorfismo $f\colon A\lra B$.
$\fin$
\end{proposicion}

En bajas dimensi\'on, tenemos para $q=0$, no de Tate, $\res_0\colon H
\to G$, $h\mapsto h$, $\co 0GA\stackrel{\res}{\lra}\co 0HA$, $\res\colon
A^G\ \lra A^H$, $a\mapsto a$ pues $a\in A^G$ implica $a\in A^H$.

Ahora, para $q=1$, $\res\colon \co 1GA\lra \co 1HA$ est\'a dada de la
siguiente forma. Si $f\in Z^1(G,A)$, $f\colon G\lra A$, $f(\sigma\tau)=\sigma
f(\tau)+f(\sigma)$ para todas $\sigma, \tau\in G$. La restricci\'on $f|_H$
satisface la misma relaci\'on pero para $\sigma, \tau\in H$. Si $f\in
B^1(G,A)$, existe $a\in A$ tal que $f(\sigma)=\sigma a-a$ para toda $\sigma
\in G$, en particular, para toda $\sigma\in H$.

\begin{ejemplo}[Inflaci\'on\index{Inflaci\'on}]\label{Ej17.5.6.6(2)}
Sea $\varphi=\pi\colon G\lra G/H$ donde $H\normal G$ (aqu\'i $G$
corresponde a $G'$ de la definici\'on general y $G/H$ corresponde a $G$).
Si $A$ es un $G$-m\'odulo, $A^H$ es un $G/H$-m\'odulo. Sea $i:A^H
\lra A$ la inyecci\'on natural. A las funciones
\[
(\pi,i)_n\colon \co n{G/H}{A^H}\lra \co nGA,\quad n\geq 0,
\]
se les llama {\em inflaci\'on\index{inflaci\'on}}. Ahora bien, como $\pi$ no
es monomorfismo, inflaci\'on s\'olo est\'a definida para $n\geq 0$. 

A nivel de $q$-cadenas, inflaci\'on se define de la siguiente forma. 
Sean $q\geq 1$ y $H\normal G$ un subgrupo normal de $G$. Sea
$\sigma\colon \underbrace{G/H\times\cdots\times G/H}_{q}
\lra A^H$ una $q$-cadena. Se define $\tau\colon\underbrace{
G\times\cdots\times G}_q\lra A$ por
\[
\tau(g_1,\ldots,g_q)=\sigma(g_1 H,\ldots, g_q H).
\]
El mapeo $\tau$ es la inflaci\'on de $\sigma$: $\tau=\infla \sigma$.
Se tiene $\partial_{q+1}\circ \infla=\infla\circ \partial_{q+1}$ y por
tanto se obtiene el mapeo en cohomolog\'ia inflaci\'on.

Para ver la inflaci\'on en bajas dimensiones, procedemos as\'i.
Sea $\pi\colon G\lra G/H$ la proyecci\'on natural. Ahora $(A^H)^{G/H}
\subseteq A^G$. Por tanto 
\[
\co 0{G/H}{A^H}\stackrel{\infla}{\lra}
\co 0GA=A^G
\]
 es la inyecci\'on natural, en el caso de $q=0$,
tanto de Tate como no de Tate.

Para $q=1$, si $f\in \co 1{G/H}{A^H}$, tomando un representante
$f\in Z^1(G/H,A^H)$, $f\colon G/H\lra A^H$ satisface $f(\bar{\sigma}
\bar{\tau})=\bar{\sigma}f(\bar{\tau})+f(\bar{\sigma})$ para todas $
\bar{\sigma},\bar{\tau}\in G/H$. Entonces $\infla f=\tilde f$ satisface
$\tilde f\colon G\lra A$, $\tilde f(g):=f(gH)=f(\bar g)\in A^H\subseteq A$.
Se tiene $\tilde f(gk)=f(gk H)=f(gH kH)=gH\cdot f(kH)+f(gH)$ el cual
est\'a bien definido pues $f(kH)\in A^H$.

Independientemente del criterio general, se puede ver que restricci\'on
se podr\'ia definir por cambio de dimensi\'on usando si $0\lra A\lra B
\lra C\lra 0$ es una sucesi\'on exacta de $G$-m\'odulos, entonces con el
mapeo de conexi\'on $\delta$, tenemos el diagrama conmutativo
\begin{gather*}
\xymatrix{
\co qGC\ar@{->}[r]^{\delta}\ar@{->}[d]_{\res_q}&\co {{q+1}}GA
\ar@{->}[d]^{\res_{q+1}}\\
\co qHC\ar@{->}[r]_{\delta}&\co {{q+1}}HA
}
\intertext{para toda $q\in{\ma Z}$. Sea 
$\mu^q\colon \co 0G{A^q}\stackrel{\cong}{\lra}
\co qGA$, se tiene el diagrama conmutativo}
\xymatrix{
\co 0G{A^q}\ar@{->}[r]^{\mu^q}\ar@{->}[d]_{\res_0}&\co qGA\ar@{->}[d]^{
\res_q}\\ \co 0H{A^q}\ar@{->}[r]_{\mu^q}&\co qHA
}
\end{gather*}

El caso de inflaci\'on es diferente, pues la sucesi\'on exacta respectiva,
deber\'ia ser: $0\lra A^H\lra B^H\lra C^G\lra 0$ pero esta sucesi\'on no
es exacta en general. De hecho, la sucesi\'on que es exacta es
$0\lra A^H\lra B^H\lra C^H\lra \co 1HA\lra\cdots$.

Ahora bien, resulta ser que E. Weiss \cite{Wei59} define un mapeo inflaci\'on
para \'indices negativos y se llama {\em deflaci\'on}.
\end{ejemplo}

Similarmente al caso de restricci\'on, tenemos:

\begin{proposicion}\label{P17.5.6.6(3)-1}
Sean $A$ y $B$ dos $G$-m\'odulos y $H\normal G$ 
un subgrupo normal de $G$, $G$
un grupo finito. Sea $f\colon A\lra B$ un $G$-homomorfismo. Entonces el
diagrama
\[
\xymatrix{
\co q{G/H}{A^H}\ar@{->}[rr]^{H^q{(f)}}\ar@{->}[d]_{\infla_q}
&&\co q{G/H}{B^H}\ar@{->}[d]^{
\infla_q}\\ \co qGA\ar@{->}[rr]_{H^q(f)}&&\co qGB
}
\]
es conmutativo para toda $q\geq 0$.
Notemos que un $G$-homomorfismo $f\colon A\lra B$,
induce un $G/H$-homomorfismo $f\colon A^H\lra B^H$.
$\fin$
\end{proposicion}

\begin{ejemplo}[Corestricci\'on\index{corestricci\'on}]\label{Ej17.5.6.6(3)}
Sea $G$ un grupo y sea $H<G$, digamos $[G:H]=m$. Sea $G=
\cupdot_{i=1}^m \sigma_i H$ la descomposici\'on de $G$ en clases
izquierdas de $H$. Se define $\Tr_{H\to G}\colon A^H\lra A^G$, $\Tr_{
H\to G}(a)=\sum_{i=1}^m\sigma_i a$. Esto es, $\Tr_{H\to G}(a)=
\N_{G/H}(a)$, aunque entendiendo que $H$ no necesariamente es
normal en $G$.

Si $\{\sigma'_i\}_{i=1}^m$ es otro conjunto de representantes, esto es,
$\sigma'_i=\sigma_ih_i$ par algunos $h_i\in H$, entonces $\sigma'_i a=
\sigma_i h_i a\underbracket[0pt]{=}_{\substack{\uparrow\\ a\in A^H}} \sigma a$.
Por tanto, $\Tr_{H\to G}(a)$ no depende de los representantes
$\{\sigma_i\}_{i=1}^m$.

Si $\theta\in G$, $\theta\sigma_i H=\sigma_{\theta(i)}H$, esto es, 
$\theta\sigma_i=\sigma_{\theta(i)} h_i$ para alguna $h_i\in H$.
Adem\'as $\{\sigma_{\theta(i)}\}_{i=1}^m\subseteq \{\sigma_i\}_{i=1}^m$.
Si $\sigma_{\theta(i)}=\sigma_{\theta(j)}$, entonces $\sigma_{\theta(i)}=
\theta\sigma_ih_i^{-1}=\theta\sigma_jh_j^{-1}=\sigma_{\theta(j)}$.
Por tanto $\sigma_ih_i^{-1}=\sigma_jh_j^{-1}$ por lo que $\sigma_i\in
\sigma_j H$ de donde se sigue que $\sigma_i=\sigma_j$. Entonces
$\tilde \theta\colon \{\sigma_i\}\lra\{\sigma_{\theta(i)}\}$ es una biyecci\'on
y por tanto
\[
\theta\big(\sum_{i=1}^m\sigma_i a\big)=\sum_{i=1}^m\theta\sigma_i a=
\sum_{i=1}^m \sigma_{\theta(i)} a=\sum_{i=1}^m\sigma_i a,
\]
de donde se obtiene que $\Tr_{H\to G} a\in A^G$.

Se sigue que 
\las
\item $\Tr_{H\to G}\colon \Hom_H(A,B)\lra \Hom_G(A,B)$ es un
homomorfismo.
\item Si $f\in \Hom_G(A,B)$, entonces $\Tr_{H\to G}(f)=[G:H] f$.
\end{list}

Sea $P$ una resoluci\'on de $\zg$-m\'odulos libres de ${\ma Z}$. 
Entonces $P$ tambi\'en es una resoluci\'on de ${\ma Z}[H]$-m\'odulos.
Sea $f\in \Hom_H(P_n,A)$, entonces
\begin{gather*}
\*\partial\big(\Tr_{H\to G}(f)\big)=\Tr_{H\to G}(f)\circ \partial=
\Tr_{H\to G}(f\circ \partial)=\Tr_{H\to G}\big(\*\partial(f)\big).
\intertext{Es decir, $\*\partial\circ \Tr_{H\to G}=\Tr_{H\to G}\circ \*\partial$.}
\xymatrix{
\Hom_H(P_n,A)\ar@{->}[r]^{\*{\partial_{n+1}}}\ar@{->}[d]_{\Tr_{H
\to G}}&\Hom_H(P_{n+1},A)\ar@{->}[d]^{\Tr_{H\to G}}\\
\Hom_G(P_n,A)\ar@{->}[r]_{\*{\partial_{n+1}}}&\Hom_G(P_{n+1},A)
}
\end{gather*}

Si $f\in\ker \*{\partial_{n+1}}$, entonces $\big(\*{\partial_{n+1}}\circ
\Tr_{H\to G}\big)(f)=\big(\Tr_{H\to G}
\circ \*{\partial_{n+1}}\big)(f)=0$, es decir, $\Tr_{H\to G} f
\in \ker \*{\partial_{n+1}}$ y por tanto tenemos homomorfismos
\[
\co nHA\xrightarrow[]{\Tr_{H\to G}}\co nGA
\]
para $n\in {\ma Z}$.

Este grupo de mapeos se llama {\em corestricci\'on}. Expl\'icitamente,
se tiene para 
\begin{align*}
\cores_0 \colon & \co 0HA\lra \co 0GA\\
& a+\N_H A\lra \N_{G/H}a +\N_G A,\\
\cores_{-1}\colon & \co {{-1}}HA\lra \co {{-1}}GA\\
& a+I_H A\lra a+I_G A \quad \text{para\ } 
a\in {_{\N_H}}A\subseteq \in {_{\N_G}} A.
\end{align*}

Sea $0\lra A\stackrel{i}{\lra} A\stackrel{j}{\lra} C\lra 0$ una sucesi\'on
exacta de $G$-m\'odulos. Entonces el siguiente diagrama es
conmutativo
\[
\xymatrix{
\co {{-1}}HC \ar@{->}[r]^{\delta}\ar@{->}[d]_{\cores_{-1}}&\co 0HA
\ar@{->}[d]^{\cores_0}\\ \co {{-1}}GC\ar@{->}[r]^{\delta}&\co 0GA
}
\]

De aqu\'i tenemos que la corestricci\'on es la \'unica
familia de homomorfismos
$\cores_q\colon \co qHA\lra \co qGA$, $q\in{\ma Z}$, tal que:
\las
\item Si $q=0$, $\cores_0\colon \co 0HA\lra \co 0GA$, 
$a+\N_H A\lra \Tr_{H\to G}a +\N_G A =\N_{G/H}a+\N_G A$, $a\in A^H$.
\item Para cada sucesi\'on $G$-exacta $0\lra A\lra B\lra C\lra 0$,
el siguiente diagrama es conmutativo
\[
\xymatrix{
\co qHC \ar@{->}[r]^{\delta}\ar@{->}[d]_{\cores_q}&\co {{q+1}}HA
\ar@{->}[d]^{\cores_{q+1}}\\ \co qGC\ar@{->}[r]^{\delta}&\co {{q+1}}GA
}
\]
Formalmente, $\cores_{q+1}=\delta\circ \cores_q\circ \delta^{-1}$.

El mapeo $\cores_q$ est\'a determinado por el diagrama conmutativo
\[
\xymatrix{
\co 0H{A^q} \ar@{->}[r]^{\mu^q}\ar@{->}[d]_{\cores_{0}}&\co qHA
\ar@{->}[d]^{\cores_q}\\ \co 0GC\ar@{->}[r]^{\mu^q}&\co qGA
}
\]
$\cores_q=\mu^q\circ \cores_0\circ (\mu^q)^{-1}$.

\end{list}
\end{ejemplo}

\begin{ejemplo}[Conjugaci\'on]\label{Ej17.5.6.6(4)}
Sea $\varphi=\varphi_{\sigma}\colon G(=G')\lra G$ dado por $\varphi_{\sigma}
(g)=\sigma g\sigma^{-1}$. Sean $A$ un $G$-m\'odulo y $\varphi^*_{\sigma} A
(=A)$ es tambi\'en un $G$-m\'odulo con acci\'on $g A=\varphi_{\sigma}(g) a=
\sigma g\sigma^{-1} a$.

Sea $f_{\sigma}\colon A\lra A$ dada por $f_{\sigma}(a)=\sigma^{-1}a$. 
Entonces 
\[
f_{\sigma}(\varphi_{\sigma}(g)a)=f_{\sigma}(\sigma g\sigma^{-1} a)=
\sigma^{-1}(\sigma g\sigma^{-1}a)=g\sigma^{-1}a=gf_{\sigma}(a)=
g\sigma^{-1}a,
\]
esto es, $f_{\sigma}(\varphi_{\sigma}(g)a)=gf_{\sigma}(a)$. Por tanto
$(\varphi_{\sigma},f_{\sigma})$ es compatible y puesto que $\varphi_{\sigma}$
es inyectiva, se tienen homomorfismos $(\varphi_{\sigma},f_{\sigma})_n
\colon \co nGA\lra \co nGA$ para toda $n\in{\ma Z}$.

Para $n=0$, $(\varphi_{\sigma},f_{\sigma})_0\colon \Id_{\co 0GA}\colon
A^G/\N_G A\lra A^G/\N_G a$ y $f_{\sigma}(a)=\sigma^{-1}a=a$
para toda $a\in A^G$.
El diagrama
\[
\xymatrix{
\co 0G{A^q}\ar@{->}[r]^{\mu^q}\ar@{->}[d]_{(\varphi_{\sigma},f_{\sigma})_0}
&\co qGA\ar@{->}[d]^{(\varphi_{\sigma},f_{\sigma})_q}\\
\co 0G{A^q}\ar@{->}[r]^{\mu^q}&\co qGA
}
\]
es conmutativo y $(\varphi_{\sigma},f_{\sigma})_0=\Id$ implica que
$(\varphi_{\sigma},f_{\sigma})_q=\Id$ para toda $q\in{\ma Z}$.

\end{ejemplo}

\subsection{Lema de Shapiro, la sucesi\'on Inflaci\'on-restricci\'on,
mapeo de transferencia}\label{SShapiro}

Empezamos analizando la composici\'on corestricci\'on-restricci\'on:
$\cores\circ \res$.

\begin{teorema}\label{corestriccionrestriccion}
Sea $G$ un grupo finito y sea $H$ un subgrupo, Entonces 
la composici\'on
\[
\co qGA\xrightarrow[]{\res}\co qHA\xrightarrow[]{\cores}\co qGA
\]
es el endomorfismo $\cores\circ \res= n(=n\cdot \Id)$ donde $n=
[G:H]$.
\end{teorema}

\begin{proof}
Si $a\in A^G$, $\cores_0\circ \res_0(\bar a)=\cores_0(a+\N_H A)=
\sum_{\sigma\in G/H} \sigma a+\N_G A=na +\N_G A= n\bar a$.

Para $q$ general, tenemos
\[
\xymatrix{
\co 0G{A^q} \ar@{->}[rr]^{\cores_0\circ \res_0=n}\ar@{->}[d]_{\cong}^{\mu^q}
&&\co 0G{A^q}\ar@{->}[d]^{\mu^q}_{\cong}\\ \co qGA\ar@{->}[rr]^{
\cores_q\circ \res_q}&&\co qGA
}
\]
por lo que $\cores_q\circ \res_q= n$ para toda $n\in{\ma Z}$.
$\fin$
\end{proof}

\begin{observacion}\label{O17.5.6.10}
El Teorema \ref{corestriccionrestriccion} implica que si tomamos
$H=\{1\}$, entonces $n=|G|$ y en particular $n\co qGA=0$
para toda $q\in{\ma Z}$.
\end{observacion}

Sea $f\colon A\lra B$ es un $G$-homomorfismo y si $H<G$, los 
siguientes diagramas son conmutativos
\[
\xymatrix{
\co 0G{A} \ar@{->}[rr]^{H^0(f)}\ar@<1ex>[d]^{\res}
&&\co 0GB\ar@<1ex>[d]^{\res}\\ \co qHA\ar@<1ex>[u]^{
\cores}\ar@{->}[rr]^{H^q(f)}&&\co qHB\ar@<1ex>[u]^{
\cores}
}
\]

Es decir, tanto restricci\'on como restricci\'on conmutan con los
homomorfismos en cohomolog\'ia inducidos de un homomomorfismo
de $G$-m\'odulos.

Ahora bien, los grupos $\co qGA$ son grupos abelianos de torsi\'on pues
$n\co qGA=0$ donde $n=|G|$. Por tanto $\co qGA$ es suma directa
de sus $p$-subgrupos de Sylow:
\[
\co qGA=\bigoplus_p\co qGA_p,
\]
donde $\co qGA_p$ es el $p$-subgrupo de Sylow de $\co qGA$.

\begin{proposicion}\label{P17.5.6.10-1}
La restricci\'on $\res\colon \co qGA_p\lra \co q{G_p}A$ es inyectiva y la
corestricci\'on $\cores \colon \co q{G_p}A\lra \co qGA_p$ es biyectiva donde
$G_p$ es un $p$-subgrupo de Sylow de $G$.
\end{proposicion}

\begin{proof}
Se tiene $\cores \circ \res=[G:G_p]\Id$ y $\mcd(p,[G:G_p])=1$, entonces el
mapeo $\co qGA_p\xrightarrow[]{\cores\circ \res} \co qGA_p$ es biyectiva.
Por tanto, si $x\in \co qGA_p$ es tal que $\res x=0$, $\cores\circ \res x=0$,
por tanto $\res$ es inyectiva.

Ahora si $|G_p|=p^m$, $p^m\co q{G_p}A=0$ por lo que $\cores \big(\co
q{G_p}A\big)\subseteq \co qGA_p$:
\[
\xymatrix{
\co qGA_p\ar@{->}[r]^{\res}\ar@{->}@/_1pc/[rr]_{\text{biyectiva}}&
\co q{G_p}A\ar@{->}[r]^{\cores}&\co qGA_p.
}
\]
Por tanto, $\cores \big(\co q{G_p}A\big)=\co qGA_p$.
$\fin$
\end{proof}

\begin{corolario}\label{C17.5.6.10-2}
Si para cada primo $p$ se tiene que $\co q{G_p}A=0$ para alg\'un
$p$-subgrupo de Sylow $G_p$ de $G$, entonces $\co qGA=0$.
\end{corolario}

\begin{proof}
Se tiene $\res\colon\co qGA_p\lra \co q{G_p}A$ es inyectiva,  por lo tanto
$\co qGA_p=0$ para toda $p$. Se sigue que $\co qGA=0$.
$\fin$
\end{proof}

\begin{definicion}\label{D17.5.6.11}
Sea $G$ un grupo finito y sea $H<G$ un subgrupo de $G$. Un 
$G$-m\'odulo $A$ se llama {\em $G/H$-inducido\index{inducido!m\'odulo
relativamente $\sim$}} si $A=\bigoplus_{\sigma\in G/H}\sigma B$ donde
$B$ es un $H$-subm\'odulo de $A$ y $\sigma$ var\'ia en un conjunto de
representantes de clases izquierdas de $H$ en $G$. Cuando $H=\{1\}$,
$G/H=G$ y $G/H$-inducido significa $G$-inducido.
\end{definicion}

Sea $\varphi=i\colon H\lra G$ la inclusi\'on y sea $g=\pi \colon A\lra B$
donde $\pi(a)=b_1$ donde $a=\sum_{\sigma_{G/H}} \sigma b_{\sigma}$. Si $h
\in H$, 
\begin{align*}
g(\varphi(h)a)&=\pi(i(h)a)=\pi(h a)=\pi\big(\sum_{\sigma\in G/H} h\sigma(
hb_{\sigma})\big)\\
&= \pi\Big(h \cdot 1 (h(b_1))+\sum_{\substack{\sigma\in G/H\\ \sigma
\neq 1}}h\sigma(hb_{\sigma})\Big)=h(b_1),
\intertext{y}
h\circ g(a)&=h\pi (a)=h(b_1),
\end{align*}
por lo que $\varphi$ y $g$ son compatibles y tenemos homomorfismos:
\[
\co nGA\xrightarrow[]{(\varphi,g)_n} \co nHB,\quad n\in{\ma Z}.
\]

\begin{proposicion}\label{P17.5.6.12}
Sea $C$ un $G$-m\'odulo. Entonces
\[
\Hom_G(C,A)\isomo_{\psi}\Hom_H(C,B),
\]
donde $\psi(f)=\pi \circ f$, $\pi\colon A\lra B$ con $\pi(a)=b_1$
donde $a=\sum_{\sigma\in G/H}\sigma b_{\sigma}$. El inverso
de $\psi$ es $\psi^{-1}(g)=S_g=\sum_{\sigma\in G/H}\sigma g\sigma^{-1}$.
\end{proposicion}

\begin{proof}
Sea $f\in\Hom_G(C,A)$. Entonces $\pi\circ f\colon C\lra B$ es un
$H$-homomorfismo. Si $\pi\circ f=0$, entonces $\pi\circ f(c)$ para toda
$c\in C$. Si $f\neq 0$, existir\'ia $f(c)=a\neq 0$. Sea $a=\sum_{\sigma
\in G/H}\sigma b_{\sigma}$ y sea $b_{\theta}\neq 0$. Entonces
\begin{align*}
\theta^{-1}f(c)&=f(\theta^{-1}c)=\theta^{-1}a=\sum_{\sigma\in G/H}
\theta^{-1}\sigma b_{\sigma}=\theta^{-1}\theta b_{\theta}+\sum_{
\substack{\theta\in G/H\\ \sigma\neq \theta}}\theta^{-1}\sigma b_{\sigma}\\
&=b_{\theta}+\sum_{\substack{\sigma\in G/H\\ \sigma\neq \theta}}
\theta^{-1}\sigma b_{\sigma}.
\intertext{Por lo tanto}
(\pi\circ f)(\theta^{-1}c)&=\pi(f(\theta^{-1}c))=\pi\big(b_{\theta}+\sum_{
\substack{\sigma\in G/H\\ \sigma\neq \theta}}\theta^{-1}\sigma
b_{\sigma}\big)=b_{\theta}\neq 0.
\end{align*}
Por tanto $\psi$ es inyectiva.

Sea $g\in\Hom_H(C,B)$ y sea $S_g\in\Hom_G(C,A)$ dada por
\[
S_g(c)=\sum_{\sigma\in G/H}\sigma g(\sigma^{-1}c).
\]

Si $\sigma'=\sigma h$, esto es, $\sigma'H=\sigma H$ es la misma clase
izquierda, entonces
\[
\sigma'g((\sigma')^{-1}(c))=\sigma hg(h^{-1}\sigma^{-1}c)=\sigma h h^{-1}
g(\sigma^{-1}c)=\sigma g(\sigma^{-1}c),
\]
por lo que $S_g(c)$ est\'a bien definido. Se verifica directamente que 
$S_g(\theta c)=\theta S_g(c)$ para toda $\theta\in G$ y toda $c\in C$.
Por tanto $S_g\in\Hom_G(C,A)$ y $(\pi\circ S_g)(c)=S_g(c)_1=g(1^{-1}
c)=g(c)$ y por tanto $\pi\circ S_g=g$ y $\psi$ es suprayectiva.
$\fin$
\end{proof}

Como corolario, obtenemos el siguiente resultado.

\begin{teorema}[Lema de Shapiro\index{lema de 
Shapiro}\index{Shapiro!lema de $\sim$}]\label{CClaseT1.5.9''}
Sean $G$ un grupo finito y $H$ un subgrupo de $G$.
Sea $A$ un $G$--m\'odulo que es 
{\em $G/H$--inducido\index{m\'odulo inducido}}, esto es,
existe $B\subseteq A$ es un $H$--m\'odulo tal que 
\begin{gather*}
A\cong \bigoplus_{\sigma \in G/H}\sigma B,
\intertext{donde $\sigma$ corre sobre un conjunto de representantes
izquierdos de $H$ en $G$. Entonces}
H^n(G,A)\cong H^n(H,B)\quad \text{para toda $n\in{\ma Z}$}
\intertext{bajo el isomorfismo inducido por}
H^n(G,A)\xrightarrow{\res}H^n(H,A)\xrightarrow{\ \bar{\pi}\ }
H^n(H,B),
\end{gather*}
donde $\bar{\pi}$ es inducido por la proyecci\'on natural
$\pi\colon A\lra B$ y $\res$ es el mapeo de restricci\'on.
\end{teorema}

\begin{proof} Aplicar lo anterior con $C=P_n$, donde
$P$ es una $G$-resoluci\'on completa. $\fin$
\end{proof}

\begin{teorema}[inflaci\'on-restricci\'on]\label{inflacionrestriccion}
Sean $G$ un grupo finito, $A$ un $G$-m\'odulo y $H\normal G$ un
subgrupo normal de $G$. Entonces
\[
0\lra \co 1{G/H}{A^H}\stackrel{\infla}{\lra} \co 1GA\stackrel{\res}{\lra}
\co 1HA,
\]
es exacta.
\end{teorema}

\begin{proof}
Para ver la inyectividad de $\infla$, sea $f\colon G/H\lra A^H$ un
$1$-cocliclo tal que $\infla f$ es una $1$-cofrontera del $G$-m\'odulo $A$.
Entonces
\[
\infla f(\sigma)=f(\sigma H)=\sigma a-a \quad \text{para alguna $a\in A$
y para toda $\sigma\in G$.}
\]
Para $\tau\in H$, se tiene $\sigma\tau a-a=\sigma a-a$, por tanto $\tau a
=a$ para toda $\tau \in H$. Se sigue que $a\in A^H$. Por lo tanto
$f(\sigma H)=\sigma H a-a$ es una $1$-cofrontera. Se sigue que $\infla$
es inyectiva.

Ahora veamos la exactitud en $\co 1GA$. Sea $f\colon G/H\lra A^H$ un
$1$-cociclo de $A^H$. Para $\sigma\in H$, se tiene que
\[
\res\circ\infla f(\sigma)=\infla f(\sigma)=f(\sigma H)=f(H)=f(\bar 1).
\]
Ahora $f(\bar 1)=f(\bar 1\cdot\bar 1)=f(\bar 1)+f(\bar 1)$ de donde
obtenemos $f(\bar 1)=0$ y por tanto $\res\circ \infla=0$. En particular
$\im \infla \subseteq \ker \res$.

Rec\'iprocamente, sea $f\colon G\lra A$ un $1$-cociclo del $G$-m\'odulo
$A$ cuya restricci\'on a $H$ es una $1$-cofrontera del $H$-m\'odulo $A$:
\[
f(\tau)=\tau a-a,\quad a\in A\quad \text{para toda $\tau\in H$.}
\]

Sea $g\colon G\lra A$, $g(\sigma)=\sigma a-a$ para todo $\sigma \in G$,
$g$ es una $1$-cofrontera. Entonces $f-g$ es un $1$-cociclo $h(\sigma)=
f(\sigma)-g(\sigma)$ y $h$ y $f$ est\'an en la misma clase de
cohomolog\'ia y adem\'as $h(\tau)=0$ para $\tau\in H$. Entonces
\begin{gather*}
h(\sigma\tau)=h(\sigma)+\sigma h(\tau)=h(\sigma)\quad \text{para
toda $\tau \in H$}
\intertext{y adem\'as}
h(\tau\sigma)=h(\tau)+\tau h(\sigma)=\tau h(\sigma)\quad \text{para toda $\tau\in H$}.
\end{gather*}
Por tanto $\tau h(\sigma)=h(\sigma)$ para toda $\tau\in H$ y para
toda $\sigma \in G$.

Sea $F\colon G/H\lra A$ dado por $F(\sigma H)=h(\sigma)$. Entonces 
$f(\sigma H)\in A^H$ y $F$ es un cociclo con $\infla F=h$, lo cual prueba
que $\ker \res\subseteq \im \infla$. El resultado se sigue.
$\fin$
\end{proof}

Para $q\geq 2$ se tiene el an\'alogo siempre y cuando se cumpla que 
$\co iHA=0$ para $i=1,\ldots,q-1$. Esto es,

\begin{teorema}\label{inflares}
Sea $A$ un $G$-m\'odulo y sea $H\normal G$ un subgrupo normal de $G$.
Si $\co iHA=0$ para $i=1,\cdots,q-1$ y $q\geq 1$, entonces la sucesi\'on
\[
0\lra \co q{G/H}{A^H}\stackrel{\infla}{\lra}\co qGA\stackrel{\res}{\lra}\co qHA
\]
es exacta.
\end{teorema}

\begin{proof}
Se probar\'a por inducci\'on en $q$. El caso $q=1$ es el Teorema \ref{inflacionrestriccion}.
De la sucesi\'on exacta $0\lra {\ma Z}\lra\zg\lra J_G=\zg/{\ma Z}\lra 0$, se 
obtiene la sucesi\'on exacta
\begin{gather*}
0\lra \underbracket[0pt]{A}_{\substack{\ucong\\ A\otimes {\ma Z}}}\lra
\zg\otimes A\lra J_G\otimes A\lra 0.
\intertext{Puesto que $\co 1HA=0$, se sigue que la sucesi\'on}
0\lra A^H\lra \big(\zg\otimes A\big)^H\lra \big(J_G\otimes A\big)^H\lra 0
\end{gather*}
es exacta. Por tanto tenemos que el siguiente diagrama es
conmutativo.
\begin{scriptsize}
\begin{gather*}
\xymatrix{
0\ar@{->}[r]&\co {{q-1}}{G/H}{(J_G\otimes A)^H}\ar@{->}[r]^{\infla}\ar@{->}[d]_{
\delta=\mu}&\co {{q-1}}G{(J_G\otimes A)}\ar@{->}[r]^{\res}\ar@{->}[d]^{\delta
=\mu}&\co {{q-1}}H{(J_G\otimes A)}\ar@{->}[d]^{\delta=\mu}\\
0\ar@{->}[r]&\co q{G/H}{A^H}\ar@{->}[r]^{\infla}&\co qGA
\ar@{->}[r]^{\res}&\co qHA
}
\end{gather*}
\end{scriptsize}

Puesto que $\zg\otimes A$ es tanto $G$-m\'odulo
inducido como $H$-m\'odulo
inducido, los mapeos de conexi\'on $\delta=\mu$ son isomorfismos. Por
tanto $\co iH{(J_G\otimes A)}\cong \co {{i+1}}HA=0$ para $i=1,\cdots,
q-2$. Entonces por inducci\'on, la fila superior en el diagrama anterior
es exacta, por tanto lo es la fila inferior, probando el resultado.
$\fin$
\end{proof}

El caso que se nos presenta en teor\'ia de campos de clase es que si
$K\subseteq L\subseteq M$ es una torre de campos con $q=2$,
$M/K$ y $L/K$ extensiones de Galois, entonces si $G=\Gal(M/K)$,
$H=\Gal(M/L)$ y $G/H\cong \Gal(L/K)$ entonces $\co 1H{\*M}=\{1\}$
y por tanto 
\[
0\lra\co 2{\Gal(L/K)}{\*L}\stackrel{\infla}{\lra} \co 2{\Gal(
M/K)}{\*M}\stackrel{\res}{\lra}\co 2{\Gal(M/L)}{\*M}
\]
es exacta, y en particular se tiene que
podemos considerar la contenci\'on
$\co 2{\Gal(L/K)}{\*L}\subseteq
\co 2{\Gal(M/K)}{\*M}$.

Ahora veremos un mapeo importante en teor\'ia de campos de clase,
el llamado {\em mapeo de transferencia\index{mapeo de transferencia}}.
En alem\'an, el mapeo de transferencia se llama 
{\em Verlagerung\index{Verlagerung}}. Primero recrdemos que
\[
\co {{-2}}G{\ma Z}\cong G/G'=\abe G.
\]
Si $H<G$ tenemos que el mapeo de transferencia, $\Ver$, es el mapeo
restricci\'on en dimensi\'on $-2$:
$\Ver:=\res_{-2}\colon \co {{-2}}G{\ma Z}\to \co {{-2}}H{\ma Z}$,
$\Ver\colon \abe G\to \abe H$.

Sea $g\in G$. Sea $\{x_i\}_i$ un conjunto de representantes de las
clases derechas m\'odulo $H: G=H x_1\cupdot \ldots \cupdot H x_n$.
Ahora, para cada \'indice $i$, existe un \'indice $\sigma(i)$ y un elemento
$\xi\in H$ tal que $x_i g=\xi x_{\sigma(i)}$, Sea $\xi:=\prod_{i=1}^n$.
Entonces
\[
\Ver (gG')=\xi H', \quad \text{esto es,} \quad
\Ver(gG')=\Big(\prod_{i=1}^n x_i g x_{\sigma(i)}^{-1}\Big) H'.
\]

\subsection{Producto copa\index{producto copa}}\label{CClaseS1.5.1}

Uno de los resultados fundamentales en teor\'ia de campos de clase,
es el teorema de reciprocidad. El metodo cohomol\'ogico para obtener
el teorema de reciprocidad es el Teorema de Tate-Nakayama que es
un resultado del producto copa.

Para definir el producto copa de manera directa, se hace por medio
del producto tensorial a nivel de los valores de cadenas homog\'eneas.
M\'as precisamente, sean $G$ un grupo finito,
$A$ y $B$ dos $G$-m\'odulos
y consideramos $A\otimes B$ como $G$-m\'odulo: $\sigma(a\otimes b)=
\sigma a \otimes \sigma b$, $a\in A, b\in B$ y $\sigma \in G$. Consideremos
una resoluci\'on $P$ de $G$ y 
\begin{gather*}
\Hom_G(P_p,A)\times \Hom_G(P_q, B)\xrightarrow[]{\ \Cup\ }\Hom_G(
P_{p+q},A\otimes B)
\intertext{para $p,q\geq 0$, dado por la f\'ormula}
(a\Cup b)(g_0,\ldots,g_{p+q})=a(g_0,\ldots,g_p)\otimes b(g_p,\ldots,
g_{p+q}).
\end{gather*}
Para $p=q=0$, $(a\Cup b)(g_0)=a(g_0)\otimes b(g_0)=(a\otimes b)(g_0)$.

Desarrollamos el producto copa en todas las dimensiones $p,q\in {\ma Z}$.

Dados los dos $G$-m\'odulos $A$ y $B$, el mapeo $A\times B\lra A\otimes B$,
induce un mapeo bilineal can\'onico 
\begin{gather*}
A^G\times B^G\lra (A\otimes B)^G,
\intertext{que mapea $\N_G A\times \N_G B$ a $\N_G(A\otimes B)$. Por tanto, se
induce un mapeo bilineal}
\co 0GA\times \co 0GB\lra \co 0G{A\otimes B},\\ (\bar a,\bar b)\longmapsto
\overline{(a\otimes b)},\quad \bar a=a+\N_G A.
\end{gather*}

\begin{definicion}\label{CClaseD1.5.10}
El elemento $\overline{a\otimes b}\in \co 0G{A\otimes B}$ se llama el
{\em producto copa\index{producto copa}}, $\bar a\in \co 0GA$ y $\bar b
\in\co 0GB$ y se denotar\'a $\bar a\Cup \bar b=\overline{a\otimes b}$.
\end{definicion}

\begin{teorema}[Producto copa\index{producto copa}]\label{CClaseT1.5.13}
Existe una \'unica familia de mapeos bilineales $\Cup$,
llamada el {\em producto copa\index{producto copa}} tal que
para cualesquiera $p,q\in{\ma Z}$ se tiene
\[
\Cup\colon H^p(G,A)\times H^q(G,B)\longrightarrow 
H^{p+q}(G,A\otimes B)
\]
que satisface:
\las
\item\label{copa1} Si $p=q=0$, el producto copa est\'a dado por el producto tensorial:
$(\bar{a},\bar{b})\mapsto \bar{a}\Cup \bar{b}=
\overline{a\otimes b}$, $\bar{a}\in H^0(G,A)=A^G/\N A$, $\bar{b}
\in H^0(G,B)=B^G/\N B$.

\item\label{copa2} Si $0\to A\to B\to C \to 0$ y
$0\to A\otimes D\to B\otimes D\to C\otimes D\to 0$
son sucesiones $G$-exactas, entonces el siguiente diagrama conmuta:
\[
\begin{CD}
H^p(G,C)&\times& H^q(G,D)@>\Cup>>H^{p+q}(G,C\otimes D)\\
@V{\delta}VV@VV{1}V@VV{\delta}V\\
H^{p+1}(G,A)&\times& H^q(G,D)@>\Cup>>H^{p+q+1}(G,A\otimes D)\\
\end{CD}
\]
es decir, $\delta({c}\Cup {d})=\delta(
{c})\Cup {d}$, ${c}
\in H^p(G,C)$, ${d}\in H^q(G,D)$ con $\delta$
los mapeos de conexi\'on.

\item\label{copa3} Si $0\to A\to B \to C\to 0$ y
$0\to D\otimes A\to D\otimes B \to D\otimes C\to 0$
son $G$-exactas, entonces el siguiente diagrama conmuta:
\[
\begin{CD}
H^p(G,D)&\times& H^q(G,C)@>\Cup>>H^{p+q}(G,D
\otimes C)\\
@V{1}VV@VV{\delta}V@VV{(-1)^p \delta}V\\
H^{p}(G,D)&\times& H^{q+1}(G,A)@>\Cup>>H^{p+q+1}(G,D\otimes A)\\
\end{CD}
\]
es decir, $\delta({d}\Cup c)=(-1)^p\big(
{d}\Cup \delta c\big)$, $d
\in H^p(G,D)$, $c\in H^q(G,C)$.
\end{list}
\end{teorema}

\begin{observacion}\label{CClaseO1.5.14-1}
El factor $(-1)^p$ se debe a la anticonmutatividad del mapeo de
conexi\'on $\delta$.
\end{observacion}

\begin{observacion}\label{CClaseO1.5.14-2}
Lo que se est\'a haciendo es definir un producto copa $\co pGA\times
\co qGB\stackrel{\Cup}{\lra}\co {{p+q}}G{A\otimes B}$ por medio del
diagrama (\ref{diagramaproductocopa}) m\'as adelante y ver que satisface
las condiciones (\ref{copa1})--(\ref{copa3}) del Teorema 
\ref{CClaseT1.5.13}. Las condiciones
(\ref{copa1})--(\ref{copa3}) hace \'unica la posible definici\'on de $\Cup$ por los 
isomorfismos $\co 0G{A^p}\isomo\limits_{\mu^p}\co pGA$ y de $\co 0G{B^q}
\isomo\limits_{\mu^q}\co qGB$. Si cambiamos el isomorfismo $(-1)^{pq}
\mu^q$ por otro, obtendr\'iamos otro producto copa con otras
condiciones (\ref{copa1})--(\ref{copa3}) del Teorema \ref{CClaseT1.5.13}.
\end{observacion}

\begin{proof} (Teorema \ref{CClaseT1.5.13}) Como en el caso de los
mapeos de restricci\'on, se obtiene el producto copa en general del
caso particular $p=0, q=0$ usando el cambio de dimensi\'on.

Se tiene:
\begin{gather*}
A^p\otimes B=J_G\otimes\cdots\otimes J_G\otimes A\otimes B=(A\otimes B)^p
\intertext{y}
A\otimes B^q=A\otimes J_G\otimes\cdots\otimes J_G\otimes B=(
J_G\otimes\cdots\otimes J_G)\otimes A\otimes B=(A\otimes B)^q,
\end{gather*}
para $A$ y $B$ dos $G$-m\'odulos y para $p,q\geq 0$. Similarmente para
$p, q\leq 0$ con $I_G$ en lugar de $J_G$.

Usando el isomorfismo $\mu^m\colon \co {{q-m}}G{A^m}\stackrel{\cong}{\lra}
\co qHA$, $m\in{\ma Z}$, empezando con el caso $p=0$, $q=0$ y obtenemos
el producto copa del diagrama
\begin{gather}\label{diagramaproductocopa}
\xymatrix{
\co 0G{A^p}\times \co 0G{B^q}\ar@{->}[rr]^{\hspace{-38pt}\Cup}_{\hspace{-38pt}
(\bar x,\bar y)\to \overline{(
x\otimes y)}}\ar@<-6ex>[d]_{\cong}^{\mu^p}\ar@<7ex>[d]_{\cong}^{1}&&
\co 0G{(A\otimes B^q)^p}=\co 0G{A^p\otimes B^q}\ar@<-8ex>[d]^{\mu^p}\\
\co pGA\times \co 0G{B^q}\ar@{->}[rr]^{\hspace{-30pt}\Cup}_{\hspace{-30pt}
\varphi}\ar@<-6ex>[d]_{\cong}^{1}\ar@<7ex>[d]_{\cong}^{\mu^q}
&& \co pG{(A\otimes B)^q}=\co pG{A\otimes B^q}
\ar@<-8ex>[d]^{(-1)^{pq}\mu^q}\\
\co pGA\times \co qGB\ar@{->}[rr]^{\hspace{-30pt}\Cup}_{\hspace{-30pt}
\psi}&& \co {{p+q}}G{A\otimes B}{\phantom{=\co pG{A\otimes B^q}}}
}
\end{gather}
Los isomorfismos 
\begin{gather*}
\co pG{(A\otimes B)^q}\xrightarrow[]{(-1)^{pq}\mu^q}
\co {{p+q}}G{A\otimes B}
\intertext{y}
\co 0G{(A\otimes B^q)^p}\xrightarrow[]{\mu^p}
\co pG{(A\otimes B)^q}
\end{gather*}
 determinan el producto copa. Si cambiamos los
isomorfismos $\mu^p$ y $(-1)^{pq}\mu^q$ por algunos otros isomorfismos,
el primero dependiendo de $p$ y el segundo de $p$ y de $q$, obtendr\'iamos
otro producto copa. No sabemos de alg\'un otro producto copa que sea
\'util desde nuestro punto de vista: El Teorema de Reciprocidad en
teor\'ia de campos de clase y la definici\'on de estos isomorfismos como
$(-1)^{pq}$ y $\mu^p$ de debe a que a nivel de cocadenas queremos el
producto tensorial
\begin{gather*}
(a\Cup b)(g_0,\ldots, g_{p+q})=a(g_0,\ldots,g_p)\otimes b(g_p,\ldots, g_{p+q})
\intertext{y por lo tanto}
\partial (a\Cup b)=(\partial a\Cup b)+(-1)^p(a\Cup \partial b).
\end{gather*}

Regresando al diagrama (\ref{diagramaproductocopa}), se tiene que
el mapeo $\varphi$ se obtiene del diagrama superior pues $\mu^p$ y $
\mu^p\times 1$ son isomorfismos. El mapeo $\psi$ se obtiene del 
diagrama inferior. Veremos que $\psi=\Cup$ es el producto que satisface
(\ref{copa2}) y (\ref{copa3}) del teorema. Adem\'as si $\psi$ satisface 
las condiciones del
teorema, este mapeo es \'unici como consecuencia de (\ref{copa1}), 
(\ref{copa2}) y (\ref{copa3}).

Antes de probar las condiciones para $\psi$, veamos $\psi=\Cup$ expl\'icitamente
en t\'erminos de cociclos en los casos $(0,q)$ y $(p,0)$ y que corresponden
a la definici\'on directa en $m$-cadenas dadas al inicio de la subsecci\'on.

\begin{proposicion}\label{productocopapyq}
Sean $a_p$ y $b_q$ un $p$-cociclo de $A$ y un $q$-cociclo de $B$ y $\bar a_p$,
$\bar b_q$ sus clases de cohomolog\'ia, entonces
\begin{gather*}
\bar a_0\Cup \bar b_q=\overline{a_0\otimes b_q}\quad \text{y}\quad
\bar a_p\Cup \bar b_0=\overline{a_p\otimes b_0}.
\intertext{Notemos que si $b_q(\sigma_1,\ldots,\sigma_q)\in B$ es un $q$-cociclo,
entonces} 
a_0\otimes b_q(\sigma_1,\ldots,\sigma_q)\in A\otimes B
\end{gather*}
con $a_0\in A^G$, es tambi\'en un $q$-cociclo.
\end{proposicion}

\begin{proof}
Se puede verificar directamente que $\bar a_0\Cup \bar b_q$ y $\bar a_p
\Cup \bar b_0$ definidos en la proposici\'on satisfacen las condiciones
(\ref{copa1}), (\ref{copa2}) y (\ref{copa3}) del teorema para $(0,q)$ y
$(p,0)$ respectivamente, viendo el comportamiento de los cociclos bajo
los respectivos mapeos. Finalmente la parte inferior del diagrama 
(\ref{diagramaproductocopa}) para $p=0$ y la parte superior para $q=0$,
se sigue que el producto dado por (\ref{diagramaproductocopa}) coincide
con lo dado en la proposici\'on.
$\fin$
\end{proof}

Sean $0\lra B\lra C\lra 0$ y $0\lra A\otimes D\lra B\otimes D\lra C\otimes
D\lra 0$ sucesiones $G$-exactas y similarmente $0\lra A\lra B\lra C\lra 0$
y $0\lra D\otimes A\lra D\otimes B\lra D\otimes C\lra 0$. Se tienen las
siguientes  sucesiones $G$-exactas:
\begin{gather*}
0\lra A^q\lra B^q\lra C^q\lra 0\quad\text{y}\\
0\lra(A\otimes D)^q\lra (B\otimes D)^q\lra (C\otimes D)^q\lra 0
\intertext{y tambi\'en}
0\lra A^p\lra B^p\lra C^p\lra 0\quad\text{y}\\
0\lra (D\otimes A)^p\lra (D\otimes B)^p\lra (D\otimes C)^p\lra 0.
\end{gather*}

Se tienen diagramas:
\begin{tiny}
\begin{gather*}
 \xymatrix{
 \co pGC\times \co 0G{D^q}\ar@{->}[rr]^{\Cup}\ar@{->}[dd]_{(1,\mu^q)}
 \ar@{->}[rd]^{(\delta,1)}&& \co pG{(C\otimes D)^q}\ar@{->}[rd]^{\delta}
 \ar@{->}[dd]|!{[ld];[dr]}\hole^{\substack{(-1)^{pq}\mu^q\\ {\ }\\{\ }\\{\ }\\ {\ }}}\\
 &\co {{p+1}}GA\times \co 0G{D^q}\ar@{->}[rr]_{\hspace{-20pt}
 \Cup}\ar@{->}[dd]^{\substack{(1,\mu^q)\\ {\ }\\{\ }\\{\ }\\ {\ }}}
 &&\co {{p+1}}G{(A\otimes D)^q}\ar@{->}[dd]^{(-1)^{(p+1)q}\mu^q}\\
 \co pGC\times \co qGD\ar@{->}[rr]|!{[rd];[ur]}\hole_{\hspace{-20pt}
 \Cup}\ar@{->}[rd]_{(\delta,1)}
 &&\co {{p+q}}G{C\otimes D}\ar@{->}[rd]^{\delta}\\
 &\co {{p+1}}GA\times \co qGD\ar@{->}[rr]^{\Cup}&&\co {{p+q+1}}G{A\otimes D}
 }\\
\\
 \xymatrix{
 \co 0G{D^p}\times \co qGC\ar@{->}[rr]^{\Cup}\ar@{->}[dd]_{(\mu^p,1)}
 \ar@{->}[rd]^{(1,\delta)}&& \co qG{(D\otimes C)^p}\ar@{->}[rd]^{\delta}
 \ar@{->}[dd]|!{[ld];[dr]}\hole^{\substack{\mu^p\\ {\ }\\{\ }\\{\ }\\ {\ }}}\\
 &\co 0G{D^p}\times \co {{q+1}}GA\ar@{->}[rr]_{\hspace{-20pt}
 \Cup}\ar@{->}[dd]^{\substack{(\mu^p,1)\\ {\ }\\{\ }\\{\ }\\ {\ }}}
 &&\co {{q+1}}G{(D\otimes A)^p}\ar@{->}[dd]^{\mu^q}\\
 \co pGD\times \co qGC\ar@{->}[rr]|!{[rd];[ur]}\hole_{\hspace{-20pt}
 \Cup}\ar@{->}[rd]_{(1,\delta)}
 &&\co {{p+q}}G{D\otimes C}\ar@{->}[rd]^{(-1)^p\delta}\\
 &\co {{p}}GD\times \co {{q+1}}GA\ar@{->}[rr]^{\Cup}&&\co {{p+q+1}}G{D\otimes A}
 }
\end{gather*}
\end{tiny}

Las partes (\ref{copa2}) y (\ref{copa3}) del teorema, son las partes inferiores de los
diagramas anteriores.

Es inmediato que los lados izquierdos conmutan. Ahora, los lados derechos
se componen de $q$ (resp. $p$) cuadrados de cohomolog\'ia de cambio de
dimensi\'on por lo que son conmutativos. Las partes frontales y traseras de
los diagramas conmutan por la definici\'on del producto copa dado en el
diagrama (\ref{diagramaproductocopa}) y finalmente, los cuadrados superiores
conmutan por la Proposici\'on \ref{productocopapyq}. Puesto que los mapeos
verticales son biyectivos, la conmutatividad de los cuadros superiores 
implican la conmutatividad de los cuadros inferiores, lo cual prueba (\ref{copa2})
y (\ref{copa3}) del teorema. Esto finaliza la prueba.
$\fin$
\end{proof}

\begin{proposicion}\label{homoprodcopa}
Sean $f\colon A\lra A'$ y $g\colon B\lra B'$ dos $G$-ho\-mo\-mor\-fis\-mos y sea
$f\otimes g\colon A\otimes A\lra A'\otimes B'$ el homomorfismo inducido por 
$f$ y $g$. Si $\bar a\in \co pGA$ y $\bar b\in \co qGB$, entonces
\[
\bar f(\bar a) \Cup \bar g(\bar b)=\overline{(f\otimes g)}(\bar a\Cup \bar b)\in 
\co {{p+q}}G{A'\otimes B'}.
\]
\end{proposicion}

\begin{proof}
La afirmaci\'on es equivalente a que el diagrama
\[
\xymatrix{
\co pGA\times \co qGB\ar@{->}[d]_{(\bar f,\bar g)}\ar@{->}[r]^{\Cup}&
\co {{p+q}}G{A\otimes B}\ar@{->}[d]^{f\otimes g}\\
\co pG{A'}\times \co qG{B'}\ar@{->}[r]_{\Cup}&\co {{p+q}}G{A'\otimes B'}
}
\]
es conmutativo para todas $p,q\in{\ma Z}$.

Para $p=q=0$, la igualdad se sigue de lo siguiente:
\begin{gather*}
\Cup(\bar f,\bar g)(\bar a,\bar b)=\bar f(\bar a)\Cup \bar g(\bar b)=
\overline{f(a)\otimes g(b)},\\
(f\otimes g)(\bar a\Cup \bar b)=(f\otimes g)(\overline{a\otimes b})=
\overline{f(a)\otimes g(b)}.
\end{gather*}

Para el caso general, consideremos el diagrama
\begin{tiny}
\begin{gather*}
\xymatrix{
\co 0G{A^p}\times \co 0G{B^q}\ar@{->}[rr]^{\Cup}\ar@{->}[dd]_{(\bar f,\bar g)}
\ar@{->}[rd]^{(\mu^p,\mu^q)}&& \co 0G{A^p\otimes B^q}\ar@{->}[rd]^{\mu^{p+q}}
\ar@{->}[dd]|!{[ld];[dr]}\hole^{\substack{f\otimes g\\ {\ }\\{\ }\\{\ }\\ {\ }}}\\
&\co pGA\times \co qGB\ar@{->}[rr]_{\hspace{-20pt}\Cup}
\ar@{->}[dd]^{\substack{(\bar f,\bar g)\\ {\ }\\{\ }\\{\ }\\ {\ }}}
&&\co {{p+q}}G{A\otimes B}\ar@{->}[dd]^{f\otimes g}\\
\co 0G{(A')^p}\times \co 0G{(B')^q}\ar@{->}[rr]|!{[rd];[ur]}\hole_{\hspace{-20pt}
\Cup}\ar@{->}[rd]_{(\mu^p,\mu^q)}
&&\co 0G{(A')^p\otimes (B')^q}\ar@{->}[rd]^{\mu^{p+q}}\\
&\co pG{A'}\times \co qG{B'}\ar@{-}[rr]^{\Cup}&&\co {{p+q}}G{A'\otimes B'}
}
\end{gather*}
\end{tiny}
el cual es conmutativo. Los isomorfismos se siguen del hecho de que
$A^p\otimes B^q\cong (A\otimes B^q)^p$, etc.
$\fin$
\end{proof}

\begin{proposicion}\label{CClaseP1.5.14-1}
Sean $A$, $B$ dos $G$-m\'odulos y sea $H<G$ un subgrupo de $G$. Si 
$\bar a\in\co pGA$ y $\bar b\in \co qGB$, entonces
\begin{gather*}
\res(\bar a\Cup\bar  b)=\res(\bar a)\Cup \res(\bar b)\in \co {{p+q}}H{A\otimes B}
\intertext{y}
\cores(\res(\bar a)\Cup \bar b)=\bar a\Cup \cores(\bar b)\in \co {{p+q}}G{A\otimes B}.
\end{gather*}
\end{proposicion}

\begin{proof}
Para $p=q=0$, la primera f\'ormula es inmediata. Para la segunda, sean $a\in
A^G$, $b\in B^H$, $0$-cociclos representando $\bar a$ y $\bar b$. Se tiene
\begin{align*}
\cores(\res(\bar a)\Cup \bar b)&=\cores(a\otimes b +\N_H(A\otimes B))
=\sum_{\sigma\in G/H}\sigma(a\otimes b)+\N_G(A\otimes B)\\
&=\sum_{\sigma\in G/H}(a\otimes \sigma(b))+\N_G(A\otimes B)\quad \text{($a$ es $G$-invariante)}\\
&=a\otimes \big(\sum_{\sigma\in G/H}\sigma(b)\big)+\N_G(A\otimes B)
=\bar a\Cup \cores(\bar b).
\end{align*}
El caso general se sigue por cambio de dimensi\'on.
$\fin$
\end{proof}

Similarmente, el siguiente resultado es inmediato para $p=q=0$ y el
caso general se sigue por cambio de dimensi\'on.

\begin{teorema}\label{CClaseT1.5.14-2}
Sean $\bar a\in \co pGA$, $\bar b\in \co qGB$, $\bar c\in \co rGC$. Entonces
\begin{gather*}
\bar a\Cup \bar b=(-1)^{pq}(\bar b\Cup \bar a)\in \co {{p+q}}G{A\otimes B},
\intertext{y}
(\bar a\Cup \bar b)\Cup c=\bar a\Cup (\bar b\Cup \bar c)\in \co {{p+q+r}}G{
A\otimes(B\otimes C)}. \tag*{$\fin$}
\end{gather*}
\end{teorema}

Sean $a_p$ y $b_q$ un $p$-cociclo de $A$ y un $q$-cociclo de $B$ respectivamente
y sean $\bar a_p$, $\bar b_q$ sus clases de cohomolog\'ia en $\co pGA$ y $\co qGB$
respectivamente.

\begin{lema}\label{CClaseL1.5.14-3}
Si $\bar x_0=\bar a_1\Cup \bar b_{-1}\in \co 0G{A\otimes B}$, entonces
\[
x_0=\sum_{\sigma \in G}a_1(\sigma)\otimes \sigma(b_{-1}), \quad \bar b_{-1}
\in \co {{-1}}GB=\frac{\ker \N_B}{I_G B}.
\]
\end{lema}

\begin{proof}
Sean $A':=\zg\otimes A$ el $G$-m\'odulo inducido, $A''=J_G\otimes A$ y
las sucesiones exactas ($A\cong {\ma Z}\otimes A$ y si un $G$-m\'odulo
es ${\ma Z}$-libre, preserva exactitud como ${\ma Z}$-m\'odulos y por
tanto como $G$-m\'odulos)
\begin{gather*}
0\lra A\lra A'\lra A''\lra 0,\\
0\lra A\otimes B\lra A'\otimes B\lra A''\otimes B\lra 0.
\end{gather*}

Ahora bien $\co 1G{A'}=0=\frac{Z^1(G,A')}{B^1(G,A')}$ y $A\subseteq A'$,
as\'i que existe una $0$-cadena $a'_0\in A'$ con $a_1=\partial a'_0$ y $
a_1(\sigma)=\sigma a'_0-a'_0$ para toda $\sigma \in G$. Sea $a''_0
\in (A'')^G$ la imagen de $a'_0$ en $A''$. El mapeo de conexi\'on
$\delta$ satisface $\bar a_1=\delta(a''_0)$ y obtenemos por el
Teorema \ref{CClaseT1.5.13} y por (\ref{d_0})
\begin{align*}
\bar a_1\Cup \bar b_{-1}&=\delta (\overline{a''_0})\Cup \bar b_{-1}=
\delta(\overline{a''_0}\Cup \bar b_{-1})=\delta(\overline{a''_0\otimes
b_{-1}})=\partial (\overline{a'_0\otimes b_{-1}})\\
&=\overline{\N_G(a'_0\otimes b_{-1})}=
\overline{\sum_{\sigma\in G} \sigma (a'_0)\otimes \sigma(b_{-1})}
\igual_{\substack{\uparrow\\ a_1(\sigma)=\sigma(a'_0)-a'_0}}\\
&=\overline{\sum_{\sigma\in G}(a_1(\sigma)+a'_0)\otimes \sigma(b_{-1})}
=\overline{\sum_{\sigma\in G}(a_1(\sigma)\otimes \sigma(b_{-1})}+
\overline{a'_0\otimes \N_G b_{-1}}\\
&=\overline{\sum_{\sigma\in G}(a_1(\sigma)\otimes \sigma(b_{-1})},
\end{align*}
puesto que $\N_G b_{-1}=0$.
$\fin$
\end{proof}

Tomemos ahora $B={\ma Z}$ e identificamos $A\otimes {\ma Z}
\stackrel{\cong}{\lra} A$, $a\otimes n\longmapsto na$. Recordemos
que $\co {{-2}}G{\ma Z}\cong G/G'\cong \abe G$. Si $\sigma\in G$,
sea $\bar \sigma$ el elemento de $\co {{-2}}G{\ma Z}$ que corresponde
a $\sigma G'\in \abe G$.

\begin{lema}\label{CClaseL1.5.14-4}
$\bar a_1\Cup \bar \sigma=\overline{a_1(\sigma)}\in \co {{-1}}GA$.
\end{lema}

\begin{proof}
De la sucesi\'on exacta $0\lra A\otimes I_G\lra A\otimes \zg\lra A\lra 0$
y de que $A\otimes \zg$ es cohomol\'ogicamente trivial, se sigue que
$\co {{-1}}GA\xrightarrow[\cong]{\ \delta\ }\co 0G{A\otimes I_G}$. Por tanto
basta probar que $\delta(\bar a_1\Cup \bar \sigma)=\delta(\overline{
a_1(\sigma)})$.

Usando la definici\'on de $\delta$ (ver Corolario \ref{CoC17.5.2.1}), se tiene
$\delta(\overline{a_1(\sigma)})=\bar x_0$ con $x_0=\sum_{\theta\in G}
\theta a_1(\sigma)\otimes \theta$.

Por otro lado, bajo el isomorfismo $\co {{-2}}G{\ma Z}\stackrel{\delta}{\lra}
\co {{-1}}G{I_G}$, el elemento $\bar \sigma$ va a $\delta(\bar \sigma)=
\overline{\sigma-1}$, por lo que
\[
\delta(\bar a_1\Cup \bar \sigma)=-(\bar a_1\Cup \delta(\bar \sigma)=-
\bar a_1\Cup (\overline{\sigma-1})=\bar y_0.
\]
Por el Lema \ref{CClaseL1.5.14-3}, se tiene
\begin{gather*}
y_0=-\sum_{\theta\in G}a_1(\theta)\otimes \theta(\sigma -1)=
\sum_{\theta\in G}a_1(\theta)\otimes \theta-\sum_{\theta\in G}
a_1(\theta)\otimes \theta \sigma.
\end{gather*}

El $1$-cociclo $a_1(\theta)$ satisface $a_1(\theta)=a_1(\theta\sigma)-
\theta a_1(\sigma)$, por lo que $y_0=\sum_{\theta\in G}\theta a_1(\sigma)
\otimes \theta\sigma$.

Por tanto $y_0-x_0=\sum_{\theta\in G}\theta a_1(\sigma)\otimes \theta(
\sigma-1)=\N_G(a_1(\sigma)\otimes (\sigma-1))$, esto es, $\bar x_0
=\bar y_0$ y el resultado se sigue.
$\fin$
\end{proof}

Sea $\bar a_2\in \co 2GA$. Se obtiene un homomorfismo
\[
\varphi=\bar a_2\Cup\underline{\ \ }\colon \co {{-2}}G{\ma Z}\lra
\co 0GA(=\co {{2-2=0}}G{A\otimes {\ma Z}}),
\]
dada por 
$\varphi(\bar \sigma)=\bar a_2\Cup \bar \sigma\in \co 0GA$,
\[
\varphi\colon G/G'\cong \abe G\lra A^G/\N_G A.
\]

Por teor\'ia de campos de clase, se seleccionar\'an m\'odulos $A$
tales que $\varphi$ es un isomorfismo de grupos. Esto es, el teorema
principal de teor\'ia de campos de clase es $\abe G\cong A^G/\N_G A$.

\begin{proposicion}\label{CClaseP1.5.14-5}
Se tiene $\bar a_2\Cup \bar \sigma=\overline{\sum_{\sigma \in G}
a_2(\theta,\sigma)}\in \co 0GA$, donde $\bar \sigma\in G/G'$.
\end{proposicion}

\begin{proof}
Sea $B=\zg\otimes A$ y $C=J_G\otimes A$. Entonces se tiene la
sucesi\'on $G$-exacta $0\lra A\lra B\lra C\lra 0$. Ahora bien, por ser
$B$ un $G$-m\'odulo inducido y como $a\in A\subseteq B$, existe
una $1$-cocadena $b\in B$ con $a_2=\partial b$, es decir, se tiene
\begin{gather}\label{2cocadena}
a_2(\theta,\sigma)=\theta b(\sigma)-b(\theta \sigma)+b(\theta).
\end{gather}

La imagen $c$ de $b$ es un $1$-cociclo de $C$ tal que $\bar a_2=
\delta(c)$. Por tanto
\begin{align*}
\bar a_2\Cup \bar\sigma&\isomo_{\substack{\uparrow\\ \text{definici\'on de}\\
\text{producto copa}}} \delta(\bar c)\Cup \bar \sigma =
\delta\big(\bar c\Cup \bar \sigma\big)\isomo_{
\substack{\uparrow\\ \text{Lema \ref{CClaseL1.5.14-4}}}}\delta\big(\overline{
c(\sigma)}\big)=\overline{a_2(\sigma)}=\overline{\partial(b(\sigma))}\\
&\isomo_{\substack{\uparrow\\ \text{(\ref{d_1})}}}
\overline{\sum_{\theta\in G}\theta b(\sigma)}
\isomo_{\substack{\uparrow\\ \text{(\ref{2cocadena})}}}
\overline{\sum_{\theta\in G} a_2(\theta,\sigma)}+\overline{
\sum_{\theta\in G}b(\theta\sigma)-\sum_{\theta\in G}b(\theta)}\\
&=\overline{\sum_{\theta\in G}a_2(\theta,\sigma)}. \tag*{$\fin$}
\end{align*}
\end{proof}

Definimos el cociente de Herbrand para la cohomolog\'ia de grupos
c\'iclicos finitos. Esta definici\'on se puede extender a endomorfismos.

\begin{definicion}\label{cocientedeHerbrand}
Sea $A$ un grupo abeliano y sean $f$ y $g$ dos endomorfismos de
$A$ tales que $f\circ g=g\circ f=0$, esto es, $\im g\subseteq \ker f$ y
$\im f\subseteq \ker g$. Entonces se define el {\em cociente de
Herbarand\index{cociente de Herbrand}\index{Herbrand!cociente de
$\sim$}}
\[
q_{f,g}(A)=\frac{[\ker f:\im g]}{[\ker g:\im f]}
\]
si ambos \'indices son finitos.
\end{definicion}

Se tiene que $q_{f,g}(A)=q_{g,f}(A)^{-1}$. Esta definici\'on
generaliza la definici\'on
anterior de cociente de Herbrand
pues si $G$ es c\'iclico de orden $n$, generado por $\sigma$ y
si $f=\sigma-1=D$ y $g=1+\sigma+\cdots+\cdots+\sigma^{n-1}=\N$, $\N\circ D=
D\circ N=\sigma^n-1=0$ y $\ker f=A^G$, $\im g=\N_G A$, $\ker g={_{
\N_G}}A$ y $\im D=I_G A$, por lo que
\[
h(A)=q_{D,N}(A)=\frac{|\co 0GA|}{|\co {{-1}}GA|}=\frac{|\co 2GA|}
{|\co 1GA|}.
\]
Tambi\'en se tiene que si $G$ es un grupo c\'iclico de orden $n$ y $G$
act\'ua trivialmente en el grupos abeliano $A$, entonces
\[
h(A)=q_{0.n}(A).
\]
Si $A$ es finito, se tiene $q_{f,g}(A)=1$.

\begin{lema}\label{CClaseL1.5.14-6}
Si $f$, $g$ son dos endomorfismos de un grupo abeliano $A$ con $f
\circ g=g\circ f=0$, entonces $q_{0,fg}(A)=q_{0,f}(A) q_{0,g}(A)$, donde
$fg$ denota multiplicaci\'on.
\end{lema}

\begin{proof}
El diagrama
\begin{gather*}
\xymatrix{
0\ar@{->}[r]&g(A)\cap \ker f\ar@{->}[r]^{\ \ \ i}\ar@{->}[d]^i&g(A)\ar@{->}
[r]^f\ar@{->}[d]^i&(fg)(A)\ar@{->}[r]\ar@{->}[d]^i&0\\
0\ar@{->}[r]&\ker f\ar@{->}[r]^i&A\ar@{->}[r]^f&f(A)\ar@{->}[r]&0
}
\intertext{es conmutativo con filas exactas. Por el Lema de la Serpiente,
la sucesi\'on}
0\lra \frac{\ker f}{g(A)\cap \ker f}\lra \frac{A}{g(A)}\lra 
\frac{f(A)}{(fg)(A)}\lra 0
\intertext{es exacta y por tanto}
\frac{[A:(fg)(A)]}{[A:f(A)]}=\frac{[A:g(A)]|g(A)\cap \ker f|}{|\ker f|}.
\intertext{Ahora,}
\frac{\ker fg}{\ker g}=\frac{g^{-1}(g(A)\cap \ker f)}{g^{-1}(0)}\cong
g(A)\cap \ker f
\intertext{por lo que}
\frac{[A:(gf)(A)]}{|\ker gf|}=\frac{[A:g(A)]}{|\ker g|}\cdot \frac{[A:f(A)]}{|\ker f|}.
\tag*{$\fin$}
\end{gather*}
\end{proof}

\begin{teorema}[Tate\index{Tate!teorema de $\sim$}\index{teorema de Tate}]
Sea $G$ un grupo c\'iclico de orden $p$, $p$ un n\'umero primo y sea $A$
un $G$-m\'odulo. Si $q_{0,p}(A)$ est\'a definido, entonces $q_{0,p}(A^G)$
y $h(A)$ est\'an tambi\'en definidos y
\[
h(A)^{p-1}=\frac{q_{0,p}(A^G)^p}{q_{0,p}(A)}.
\]
\end{teorema}

\begin{proof}
Sean $\sigma$ un generador de $G$ y $D=\sigma-1$. Consideremos la
sucesi\'on exacta $0\lra A^G\lra A\stackrel{D}{\lra} I_GA\lra 0$.

Puesto que $I_GA$ es un subgrupo y tambi\'en un grupo cociente de $A$,
se sigue que si $q_{0,p}(A)$ est\'a definido, entonces $q_{0,p}(I_G A)$
tambi\'en est\'a definido. Por tanto $q_{0,p}(A^G)=\frac{q_{0,p}(A)}{
q_{0,p}(I_GA)}$ est\'a definido.

Ahora bien, $G$ act\'ua trivialmente en $A^G$ por lo que $q_{0,p}(A^G)=
h(A^G)$. Se tiene ${\ma Z}\N_G={\ma Z}\big(\sum_{i=0}^{p-1}\sigma^i\big)$
anula a $I_G A$ y por tanto $I_G A$ es un $\g/{\ma Z}\N_G$-m\'odulo.
Como anillos tenemos el isomorfismo
\[
\zg/{\ma Z}\N_G\isomo_{\substack{\uparrow\\ \sigma\leftrightarrow x}}
{\ma Z}[x]/\langle 1+x+\cdots+x^{p-1}\rangle\cong {\ma Z}[\zeta_p],
\]
el anillo de enteros de $\cic p{}$, $\zeta_p=e^{2\pi i/p}$. El isomorfismo est\'a
inducido por el mapeo $\sigma\mapsto \zeta_p$. Se tiene que en ${\ma Z}
[\zeta_p]$, $p=(\zeta_p-1)^{p-1}u$ con $u$ una unidad de ${\ma Z}[\zeta_p]$.
Por tanto $p=(\sigma -1)^{p-1}v$, donde $v$ es una unidad de $\zg/{\ma Z}
\N_G$.

Ahora $v\colon I_G A\lra I_G A$, $x\longmapsto vx$ es un automorfismo, por
lo que $q_{0,v}(I_G A)=1$. Por tanto
\begin{gather*}
q_{0,p}(I_GA)=q_{0,D^{p-1}}(I_GA)q_{0,v}(I_G A)=q_{0,D}(I_G A)^{p-1}=
\frac {1}{q_{D,0}(I_G A)^{p-1}}.\\
\intertext{Ahora, se tiene $\N_G(I_G A)=0$, por lo que}
q_{0,p}(I_GA)=\frac{1}{q_{D,0}(I_G A)^{p-1}}=\frac{1}{q_{D,\N}(I_G A)^{p-1}}
=\frac{1}{h(I_G A)^{p-1}}.
\intertext{Como $q_{0,p}(A)=q_{0,p}(A^G)q_{0,p}(I_G A)$, se sigue que}
q_{0.p}(A^G)=h(A^G),\quad q_{0,p}(I_GA)=\frac {1}{h(I_GA)^{p-1}},
\quad q_{0,p}(A)=\frac{q_{0,p}(A^G)}{h(I_GA)^{p-1}}
\intertext{y tambi\'en}
h(A)^{p-1}=h(A^G)^{p-1}h(I_G A)^{p-1},
\intertext{por lo que}
h(A)^{p-1}=\frac{q_{0,p}(A^G)^{p-1} q_{0,p}(A^G)}{q_{0,p}(A)}=
\frac{q_{0,p}(A^G)^p}{q_{0,p}(A)}. \tag*{$\fin$}
\end{gather*}
\end{proof}

\subsubsection{Cohomolog\'ia trivial y Teorema de Tate}\label{Scohomologiatrivial}

Recordemos que $A$ es {\em cohomol\'ogicamente trivial} si $\co qHA=0$
para todo $q\in {\ma Z}$ y para todo subgrupo $H$ de $G$.

\begin{teorema}\label{dosindices}
Sean $G$ un grupo finito y $A$ un $G$--m\'odulo. Si existen dos
enteros consecutivos $q_0$, $q_0+1$ tales que $\co {{q_0}}HA=
\co {{q_0+1}}HA=0$
para todo subgrupo $H$ de $G$, entonces $A$ es cohomol\'ogicamente
trivial.
\end{teorema}

\begin{proof}
Primero veamos que $\co {{q_0-1}}HA=\co{{q_0+2}}HA=\{0\}$ para
todo subgrupo de $G$. Una vez probado esto, el resultado es
inmediato.

Para la demostraci\'on de lo anterior, aplicamos el cambio de dimensi\'on,
es decir, si tenemos lo anterior para $q_0=1$, esto es, $\co 0HA=\co 3HA
=0$ para todo subgrupo $H$ de $G$, veamos que se cumple para $q_0$
arbitrario.

Por cambio de dimensi\'on tenemos que $\co 1H{A^{q_0-1}}\cong \co {{
q_0}}HA=0$ y $\co 2H{A^{q_0-1}}\cong \co {{q_0+1}}HA=0$, por lo que
$\co {{q_0-(q_0-1)}}H{A^{q_0-1}}
\cong \co {{q_0}}HA=0$ para toda $q_0\in{\ma Z}$.

As\'i, supongamos que $\co 1HA=\co 2HA=0$ para todo subgrupo $H$
de $G$ y queremos probar que 
\begin{gather}\label{co1y3}
\co 0HA=\co 3HA=0 \quad \text{para todo $H<G$.}
\end{gather}

Probaremos (\ref{co1y3}) por inducci\'on en $|G|$, siendo el
caso $G=1$ trivial. Supongamos probado (\ref{co1y3}) para
todo subgrupo $H$ de $G$ tal que $H\neq G$. Por tanto, basta probar
que $\co 0GA=\co 3GA=0$. Si $G$ no es un $p$-grupo, todos los
subgrupos de Sylow de $G$ son propios y satisfacen (\ref{co1y3}),
por lo que (\ref{co1y3}) se sigue para $H=G$ (Corolario \ref{C17.5.6.10-2}).

Supongamos entonces que $G$ es un $p$-grupo. Sea $H\normal G$ tal que
$G/H$ es c\'iclico de orden $p$. Por hip\'otesis de inducci\'on, se tiene
$\co 0HA=\co 3HA=0$ y $\co 1HA=\co 2HA=0$.

Para $q=1,2,3$, $\infla\colon \co q{G/H}{A^H}\lra 
\co qGA$ es un isomorfismo
pues $0\lra\co 1{G/H}{A^H}\stackrel{\infla}{\lra}\co 1GA\stackrel{\res}{\lra}
\co 1HA=0$ es exacta 
y como $\co iHA=0$ para $i=1,\ldots, q-1$, $0\lra \co q{G/H}
{A^H}\stackrel{\infla}{\lra}\co qGA\stackrel{\res}{\lra}\co qHA=0$.

De esta forma, obtenemos que $\co 1GA=0$ implica $\co 1{G/H}{A^H}=0$
y $\co 3{G/H}{A^H}=0$ por ser $G/H$ un grupo c\'iclico. Ahora bien, $\co 
2GA=0$ implica $\co 2{G/H}{A^H}=0$, por lo que 
\[
\co 0{G/H}{A^H}=0=\frac{(A^H)^{G/H}}{\N_{G/H} A^H}=\frac{A^G}{
\N_{G/H} A^H},
\]
esto es, $A^G=\N_{G/H}A^H
\isomo\limits_{\substack{\uparrow\\ \co 0HA=0\\ A^H=
\N_H A}} \N_{G/H}(\N_H A)=\N_G A$. De esta forma se obtiene
que $\co 0GA=0$ y el resultado se sigue. $\fin$
\end{proof}

\begin{observacion}\label{CClaseOT1.5.14Extra}
Pueden existir grupos de cohomolog\'ia $\co {{q_0}}GA=\co
{{q_0+1}}GA=0$ pero $A$ no ser cohomol\'ogicamente
trivial.
\end{observacion}

\begin{ejemplo}\label{CCclaseE1.5.14Extra}
Sean $G={\ma Z}/6{\ma Z}$ y $A={\ma Z}/3{\ma Z}$ con
acci\'on $g\circ a=2^g a$ donde consideramos a $a$
y a $g\in {\ma Z}$ (en particular $3a=0$, etc. y $0^0=1$).

Entonces, $\co {{-1}}GA=\co 0GA=0$ pero $\co 0HA=A\neq 0$,
donde $H=\{0,2,4\}$.
\end{ejemplo}

\begin{teorema}\label{CClaseT1.5.14} 
Si $G$ es un $p$--grupo finito, donde $p$ es un n\'umero primo,
y $A$ es un $G$--m\'odulo sin $p$--torsi\'on, entonces lo siguiente es
equivalente:
\las
\item $H^i(G,A)=H^{i+1}(G,A)=0$ para dos enteros consecutivos $i,i+1$.

\item $A$ es cohomol\'ogicamente trivial.

\item El ${\ma F}_p[G]$--m\'odulo $A/pA$ es libre.
\end{list}
\end{teorema}
\begin{proof} \cite[Theorem 7.1, Part I]{Neu69}, \cite[Theorem 6, Ch. IX]{Ser}. $\fin$
\end{proof}

Consideremos el producto copa $\co pGA\times \co qGB \stackrel{\Cup}{\lra}
\co {{p+q}}G{A\otimes B}$ para dos $G$-m\'odulos, $A$ y $B$ y, como siempre,
$A\otimes B=A\otimes_{\ma Z} B$. Fijemos un elemento $a\in \co pGA$ y
consideremos el mapeo $a\Cup{\underline{\ \ }}\colon \co qGB\lra \co {{p+q}}
G{A\otimes B}$, $b\longmapsto a\Cup b$ donde $b\in \co qGB$.

\begin{teorema}[axioma de la teor\'ia de campos
de clase\index{axioma de la teor\'ia de campos de clase}]\label{CClaseT1.5.15-1}
Sea $G$ un grupo finito y sea $A$ un $G$--m\'odulo tal que para todo
subgrupo $H<G$ se tiene
\las
\item $H^{-1}(H,A)=0$,
\item $H^0(H, A)$ es un grupo c\'iclico de orden $|H|$.
\end{list}

Entonces, si $a$ genera $H^0(H,A)$, el producto copa
\[
a\Cup{\underline{\ \ }} \colon H^q(G,{\ma Z})\lra H^q(G,A)
\]
es un isomorfismo para toda $q\in {\ma Z}$ ($A\cong{\ma Z}\otimes A$).
\end{teorema}

\begin{proof}
Sea $B=A\oplus {\ma Z}[G]$. Puesto que ${\ma Z}[G]$ es cohomol\'ogicamente
trivial, se tiene que $H^q(H,B)=H^q(H,A)$ para todo subgrupo $H<G$
y para todo $q\in{\ma Z}$.

Puesto que $H^0(G,A)=A^G/\N_G A\cong C_{|H|}$, seleccionamos
$a_0\in A^G$ tal que $a=a_0+\N_G A$ es un generador de $H^0
(G,A)$. El mapeo $f\colon {\ma Z}\lra B$ dado por $m\longmapsto
(m a_0, \N_G \cdot m)$ es inyectivo pues el segundo t\'ermino
satisface $\N_G\cdot m=0\iff m=0$.

Sea $\tilde{f}\colon H^q(H,{\ma Z})\lra H^q(G,B)$ el homomorfismo
inducido por $f$. Sea  $\imath \colon A\lra B$ el encaje 
$\imath(a)=(a,0)$. Entonces, por ser $\zg$ cohomol\'ogicamente
trivial, $\imath$ induce el isomorfismo $\bar \imath\colon \co qHA
\lra \co qGB$. Recordemos que para $x_p$, $y_q$ un $p$-cociclo
y un $q$-cociclo respespectivamente, $\bar x_0\Cup \bar y_q=
\overline{x_0\otimes y_q}$ y $\bar x_p\Cup \bar y_0=\overline{
x_p\otimes y_0}$, por lo que el diagrama
\[
\xymatrix{
H^q(G,{\ma Z})\ar[rr]^{a\Cup\underline{\ \ }}\ar[drr]_{\bar{f}}&&
H^q(G,A)\ar[d]^{\bar{\imath}}\\&&H^q(G,B)
}
\]
es conmutativo. As\'i, es suficiente probar que $\bar{f}$ es biyectivo.

Se tiene que $f\colon{\ma Z}\lra B$ es inyectivo, lo cual da lugar
a una sucesi\'on exacta
\begin{gather}\label{CClaseE1.5.15-2}
0\lra{\ma Z}\xrightarrow{\ f\ }B\xrightarrow{\ \varphi\ } C\lra 0.
\end{gather}

Se tiene $H^{-1}(H,B)=H^{-1}(H,A)=0$ y $H^1(H,{\ma Z})=0$
para todo subgrupo $H<G$. Usando el Teorema \ref{CClaseT1.5.3} y
la sucesi\'on exacta (\ref{CClaseE1.5.15-2}) obtenemos
\begin{align*}
H^{-1}(H,B)=0&\lra H^{-1}(H,C)\lra H^0(H,{\ma Z})\xrightarrow{
\ \bar{f}\ }H^0(H,B)\\
&\xrightarrow{\ \bar{\varphi}\ } H^0(H,C)\lra H^1(H,{\ma Z})=0.
\end{align*}

Para $q=0$, $H^0(H,{\ma Z})\xrightarrow{\bar{f}} H^0(H,B)$ est\'a 
dado por
\begin{eqnarray*}
\bar{f}\colon {\ma Z}^H/\N_H {\ma Z}\cong {\ma Z}/|H|{\ma Z}&
\lra & H^0(H,A)\\
m&\longmapsto&ma_0
\end{eqnarray*}
por lo que $\bar{f}$ es un isomorfismo para $q=0$. 
Se sigue que $\ker \bar f=\co{{-1}}HC=0$ e
$\im \bar f=\co 0HB=\ker \varphi$, por lo que $\varphi=0$

Por tanto,
$H^{-1}(H,C)=H^0(H,C)=\{0\}$ para todo subgrupo $H<G$.
Por el Teorema \ref{dosindices}, se sigue que $C$ es cohomol\'ogicamente
trivial. Se sigue de (\ref{CClaseE1.5.15-2}) que $\bar{f}\colon H^q(G,{\ma Z})
\lra H^q(G,B)$ es biyectivo para toda $q\in{\ma Z}$ de donde se
sigue el resultado.
$\fin$
\end{proof}

\begin{teorema}[Teorema de Tate-Nakayama\index{Tate-Nakayama!teorema 
de $\sim$}\index{teorema de Tate-Nakayama}]\label{CClaseT1.5.15}
Sea $G$ un grupo finito y $A$ un $G$--m\'odulo tal que para todo
subgrupo $H<G$ se tiene
\las
\item $H^1(H,A)=0$ y

\item\label{numeral2} $H^2(H,A)$ es c{\'\i}clico de orden $|H|$.
\end{list}

Entonces si $a$ es un generador de $H^2(G,A)$, el mapeo
\[
a\Cup {\underline{\ \ }} \colon H^q(G,{\ma Z})\longrightarrow H^{q+2}(G,A)
\]
es un isomorfismo donde ${\ma Z}$ es un $G$--m\'odulo con
acci\'on trivial.

Adem\'as $\res a\in \co 2HA$ genera al grupo $\co 2HA$ por lo que
tenemos el isomorfismo $\res a\Cup \underline{\ \ }\colon \co qH{\ma Z}
\lra \co {{q+2}}HA$.
\end{teorema}

\begin{proof} 
Consideremos el isomorfismo $\mu^2\colon H^q(H,A^2)
\lra H^{q+2}(H,A)$. Se tiene $\co {{-1}}H{A^2}=\co 1HA=0$ 
y $H^0(H,A^2)$ es c\'iclico de orden $|H|$. Adem\'as, el generador
$a\in \co 2GA$ es la imagen del generador $\mu^{-2}a\in \co
0G{A^2}$. 

Por el Teorema \ref{CClaseT1.5.13} se tiene el diagrama conmutativo
\[
\xymatrix{
H^q(G,{\ma Z})\ar[rr]^{\mu^{-2}a\Cup{\underline{\ \ }}}\ar[d]_{\Id}&&
H^q(G,A^2)\ar[d]^{\mu^2}\\
H^q(G,{\ma Z})\ar[rr]^{a\Cup{\underline{\ \ }}} && H^{q+2}(G,A).
}
\]

Puesto que $\mu^{-2} a\Cup{\underline{\ \ }}$ es una biyecci\'on,
se sigue del Teorema \ref{CClaseT1.5.15-1} que $a\Cup{\underline{\ \ }}$
es un isomorfismo.

Ahora bien, puesto que $(\cores\circ \res)a=[G:H]a$, el orden del 
elmento $\res a\in \co 2HA$ es divisible por $|H|$, por lo que $\res a$
genera $\co 2HA$ por (\ref{numeral2}).
$\fin$
\end{proof}

\begin{observacion}\label{sobreTateNakayama}
Para la teor\'ia de campos de clase, el caso $q=-2$ en el Teorema de
Tate-Nakayama es particularmente importante. 

Aplicamos el Teorema de Tate-Nakaya para $q=-2$ y obtenemos en este caso
que $H^{-2}(G,{\ma Z})\cong G/G^{\prime}=\abe G$
de $G$ y el grupo residual n\'ormico $H^0(G,A)=A^G/\N A$
y un isomorfismo 
\[
\abe G\xrightarrow{\ \ \cong\ \ } A^G/\N A.
\]

Este isomorfismo es la formulaci\'on abstracta del teorema principal en
la teor{\'\i}a de campos de clase y se llama la {\em Ley de 
Reciprocidad\index{ley de reciprocidad}}. El resultado se aplica a los
siguientes casos: $L/K$ es una extensi\'on finita de Galois y $G=
\Gal(L/K)$
\las
\item $A=\*  L$ donde $L/K$ es una extensi\'on finita de campos locales.

\item $A=J_L/\*L$ donde $J_L/\*L$ es el {\em grupo de clases
id\`eles\index{grupo de id\`eles}}
de un campo global $L$.

\item $A= I_L$ donde $I_L$ es el grupo de 
{\em clases de ideales o de divisores\index{grupo de clases
se ideales}\index{grupo de clases de divisores}} de
un campo de n\'umeros o de un campo de funciones congruente.
\end{list}
\end{observacion}

\section{Campos de clase locales}\label{CClaseC3}

\subsection{Cohomolog\'ia de grupos profinitos}\label{S17.5.9}

Recordemos, Definici\'on \ref{D5.8}, que
un grupo topol\'ogico $G$ se llama {\em profinito\index{grupo profinito}}
si $G$ es Hausdorff compacto y tiene una base de vecindades 
abiertas de la identidad $1\in G$ consistente de subgrupos normales.

Sean $x\in G$ y $V$ un conjunto abierto con $x\in V$, entonces
existe $H\normal G$, $H$ abierto en $G$,
tal que $xH\subseteq V$. Ahora si $\{\xi_i\}_{i
\in I}$ es un conjunto de representantes de $G/H$, $G=\cupdot_{i\in
I}\xi_i H$ y $\xi_iH$ es un conjunto abierto con $\xi H\cap \xi_jH=
\emptyset$, para $i\neq j$ y $G$ compacto, por lo que $[G:H]<
\infty$ y $H=G\setminus \cupdot_{\substack{i\in I\\ i\neq 1}} \xi_i H$
es cerrado.

Ahora si $x\neq 1$, existen conjuntos abiertos $U$ y $V$ tales 
que $1\in U$, $x\in V$ y $U\cap V=\emptyset$. Sea $N\normal G$,
$N\subseteq U$, donde $N$ es un conjunto abierto. Entonces $1\in N$
y $G=\big(\bigcup_{g\neq N} gN\big)\cupdot N$. Por tanto $x\in
W=\bigcup_{g\neq N}g N$ y $W\cap N=\emptyset$. Se sigue que
$W$ y $N$ es una disconexi\'on de $G$ y por tanto la componente
conexa de $1$ es $\{1\}$ y por tanto $G$ es totalmente
disconexo. El rec\'iproco tambi\'en se cumple.

En resumen, se tiene $G$ es profinito $\iff$ $G$ es compacto,
Hausdorff y totalmente disconexo. Adem\'as, 
\[
G=\lim_{\substack{\longleftarrow\\ N\normal G\\ N \text{\ abierto}}} G/N
\qquad \text{(Teorema \ref{CClaseT1.4.4})}
\]
tanto topol\'ogica como algebraicamente. Notemos que $G/N$
es finito y que $N$ es abierto.

Si $L/K$ es una extensi\'on de Galois de campos no necesariamente
finita, $L=\bigcup_{\alpha}K(\alpha)$, donde $K(\alpha)/K$ es normal.
De hecho, $L=\bigcup_{\alpha\in L}K(\alpha)$, $\alpha\in L$
por lo que $\alpha$ es algebraico y $K(\alpha)/K$ es una extensi\'on
finita. Si $\widetilde{K(\alpha)}$ es la cerradura normal de $K(\alpha)/
K$, $\widetilde{K(\alpha)}=K(\tilde{\alpha})$ es una extensi\'on
finita y $L=\bigcup_{\tilde\alpha} K(\tilde \alpha)$.

Ahora $K(\alpha)/K$ normal, $K(\alpha)=L^H$ con $H\normal \Gal(
L/K)$. Damos a $G$ la topolog\'ia de Krull: una base de vecindades
abiertos de $g\in G$ es $\{gH\}_{\substack{H\normal G\\ [G:H]<\infty}}$.
Entonces $G$ es compacto y Hausdorff por lo que $G$ es profinito.

En particular, los grupos de Galois son grupos profinitos:
\[
G=\Gal(L/K)=\lim_{\substack{\longleftarrow\\ N\normal G\\ N\text{\ abierto}}}
G/N=\lim_{\substack{\longleftarrow\\ K(\alpha)/K
\\ \text{\ normal}}}
\Gal(K(\alpha)/K)
\]
y $L=\bigcup_{\alpha}K(\alpha)=\lim\limits_{\substack{\longrightarrow\\ 
\alpha}} K(\alpha)$.

En otras palabras $\Gal\Big(\lim\limits_{\substack{\longrightarrow\\ \alpha}}
K(\alpha)/K\Big)=\lim\limits_{\substack{\longleftarrow\\ \alpha}}
\Gal(K(\alpha)/K)$.

Ahora bien, los subgrupos abiertos $H$ de $G$ son de \'indice finito
y por tanto son cerrados y rec\'iprocamente, los subgrupos cerrados
de \'indice finito son abiertos. Se tiene que existen conjuntos
subgrupos cerrados que no son de \'indice finito y por tanto no
son abiertos.

Si $G$ es un grupo pro-finito y $A$ es un $G$-m\'odulo discreto,
esto es, $A$ tiene la topolog\'ia discreta, $A$ es un $G$-m\'odulo
si la acci\'on de $G$-m\'odulo usual $G\times A\lra A$ es
{\underline{continua}}. Esto es equivalente a que $A=\bigcup_N
A^N$ donde $N$ recorre los subgrupos abiertos de $G$.

\begin{definicion}\label{cohomologiainfinita}
Sea $q\geq 0$ y sean $N_1, N_2\normal G$ subgrupos normales y
abiertos de $G$ con $N_1\subseteq N_2$ y $G/N_1\lra
G/N_2$ el mapeo natural, $G\cong \lim\limits_{\substack{\longleftarrow
\\ N}}G/N$, $A=\lim\limits_{\substack{\longrightarrow \\ N}} A^N$. Sea
$\infla^{N_1}_{N_2}\colon\co q{G/N_2}{A^{N_2}}\lra \co q{G/N_1}{A^{N_1}}$
el mapeo inflaci\'on.
Entonces se define el $q$-\'esimos grupo de cohomolog\'ia de $G$
en $A$ por
\[
\co qGA:=\lim_{\substack{\longrightarrow \\ N}}\co q{G/N}{A^N},
\]
con respecto a los mapeos inflaci\'on.
\end{definicion}

Usaremos la siguiente notaci\'on. Si $L/K$ es una extensi\'on
finita de Galois y $G_{L|K}=\Gal(L/K)$, $\co q{G_{L|K}}{\*L}=
\co q{\Gal(L/K)}{\*L}=:H^q(L|K)$\label{CClaseHql/K}.
Como $H^1(L|K)=\{1\}$, se tiene la sucesi\'on exacta:
\begin{gather*}
1\lra H^2(L|K)\xrightarrow[]{\infla_M}H^2(M|K)\xrightarrow[]{\res_L}
H^2(M|L),\quad {\text{(Teorema \ref{inflares})}},
\end{gather*}
con $K\subseteq L\subseteq M$. Puesto que $\infla_M$ es inyectiva,
siempre consideraremos $H^2(L|K)\subseteq H^2(M|K)$. Si $L$
recorre todas las extensiones normales y separables de $K$,
$\bigcup_L L=\sep K$ es una cerradura separable de $K$ y ponemos
$G_K:=\Gal(\sep K/K)$, entonces
\[
\co 2{G_K}{\*{(\sep K)}}=: H^2(\ |K)=\lim_{\substack{\longrightarrow\\ L}}
H^2(L|K)=\bigcup_L H^2(L|K)\label{CClaseH2Brauer}.
\]

\begin{definicion}\label{definicionBrauer}
El {\em grupo de Brauer\index{grupo de Brauer}\index{Brauer!grupo de
$\sim$}} $\Br K$ de un campo $K$ es $\Br K=H^2(\ |K)=\bigcup\limits_{
\substack{L/K\\ \text{finita Galois}}} H^2(L|K)$.
\end{definicion}

Sea $K'/K$ una extensi\'on finita de Galois. Sea $L/K$ una extensi\'on
finita de Galois y $K'\subseteq L$. La restricci\'on $H^2(\ |K)
\xrightarrow[]{\res_{K'}} H^2(\ |K')$ se define as\'i:
\[
\res_{K'}(c)\in H^2(L|K')\subseteq H^2(\ |K').
\]

Puesto que la restricci\'on conmuta con inflaci\'on, la cual es una
inyecci\'on, se sigue que $\res_K c$ no depende de la elecci\'on
del campo $L$.

\begin{proposicion}\label{BrauerExacta}
La sucesi\'on 
\[
1\lra H^2(K'|K)\stackrel{i}{\hooklongrightarrow} H^2(\ |K)
\xrightarrow[]{\res_{K'}} H^2(\ |K')
\]
 es exacta. $\fin$
\end{proposicion}

El isomorfismo de Tate-Nakayama: $\abe {G_{L|K}}\cong
A^G/\N_{L/K} A$ es en cierta forma arbitrario pues depende
del generador $a$ de $\co 2GA$ seleccionado. Para hacer
el mapeo homog\'eneo con respecto a varias extensiones, la
segunda condici\'on del Teorema de Tate-Nakayama, esto es,
$\co 2HA$ es c\'iclico de orden $|H|$, se especializa a que 
hay un isomorfismo entre $H^2(L|K)$ y el grupo c\'iclico
$\Big(\frac {1}{[L:K]} {\ma Z}\Big)/{\ma Z}$ el llamado ``{\em
mapeo invariante\index{mapeo invariante}}'' $\inv_{L|K}$ que
determina un \'unico elemento $u_{L|K}\in H^2(L|K)$ con imagen
$\frac {1}{[L:K]}+{\ma Z}$, esto es
\[
\inv_{L|K}\colon H^2(L|K)\lra \Big(\frac {1}{[L:K]}{\ma Z}\Big)/
{\ma Z}, \quad u_{L|K}\longmapsto \frac{1}{[L:K]}+{\ma Z}.
\]

El elemento $u_{L|K}$ se llama la {\em clase
fundamental de $L/K$\index{clase fundamental}\label{CClaseinvariante}}.

Antes de continuar, veamos algunos ejemplos del grupo de
Brauer.

\begin{ejemplo}\label{Brauer1}
Si $K=\F$ es un campo finito, entonces toda extensi\'on $L/K$ es
c\'iclica, $G=\Gal(L/K)$. Adem\'as $H^1(G|\*L)=\{1\}$ y como $\*L$
es finito, el cociente de Herbrand $h(\*L)=1$, por lo que $\co
2G{\*L}\cong \co 0G{\*L}=\{1\}$, en particular $\N_{L/K}\*L=\*K$.
Por tanto $\Br \F=\bigcup_L H^2(L|K)=\bigcup_L \{0\}=\{0\}$.
De esta forma obtenemos $\Br \F=\{0\}$.
\end{ejemplo}

\begin{ejemplo}\label{Brauer2}
Si $K={\ma R}$, $K$ \'unicamente tiene dos extensiones algebraicas,
${\ma R}$ y ${\ma C}$. Por tanto 
\begin{align*}
\Br{\ma R}&=H^2({\ma C}|{\ma R})=
\co 2{G_{L/K}}{\*{\ma C}}\isomo\limits_{\substack{\uparrow\\ J \text{conjugaci\'on}\\
\text{compleja}}}\co 2{\{1,J\}}{\*{\ma C}}=\frac{(\*{\ma C})^{\{1,J\}}}{\{z\bar z\mid
z\in{\ma C}\}}\\
&=\frac{\*{\ma R}}{{\ma R}^+}\cong {\ma Z}/2{\ma Z}\cong
\frac{\frac{1}{2}{\ma Z}}{\ma Z}\subseteq {\ma Q}/{\ma Z}.
\end{align*}
\end{ejemplo}

\begin{ejemplo}\label{Brauer3}
El campo de los n\'umeros complejos no tienen extensiones de Galois
propias, por lo tanto $\Br{\ma C}=\{0\}$. Esto mismo se aplica a cualquier
campo algebraicamente cerrado.
\end{ejemplo}

\begin{ejemplo}\label{Brauer4}
Veremos m\'as adelante que si $K$ es un campo local, entonces
$\Br K={\ma Q}/{\ma Z}$ (Corolario \ref{CCLC17.6.18}).
\end{ejemplo}

\begin{observacion}\label{usodeBrauer}
Usaremos el grupo de Brauer de un campo local para hallar un mapeo
$\inv_{L|K}\colon H^2(L|K)\lra \big(\frac 1{[L:K]}{\ma Z}\big)/{\ma Z}$
y el elemento $u_{L|K}$ con $\inv_{L|K}(u_{L|K})=\frac 1{[L:K]}
+{\ma Z}\in {\ma Q}/{\ma Z}$. Este elemento nos dar\'a el mapeo
de reciprocidad v\'ia el Teorema de Tate-Nakayam:
\[
u_{L|K}\Cup\underline{\ \ }\colon \co q{G_{L|K}}{\ma Z}
\stackrel{\cong}{\lra} H^{q+2}(L|K).
\]
\end{observacion}

Originalmente, el grupo de Brauer fue definido para
codificar anillos de divisi\'on
sobre campos, es decir, objetivos no conmutativos, y sirve para
hallar la funci\'on que nos da la ley de reciprocidad.
Este enfoque
fue desarrollado por Brauer, Hasse y Noether. De hecho, el
grupo de Brauer tambi\'en describe la correspondencia de los
campos de clase para campos globales. El estudio sistem\'atico
puede ser consultado en los libros de Serre \cite{Ser} y de
Kato, Kurokawa y Saito \cite{KaKuSa2011}.

Antes de usar la teor\'ia de cohomolog\'ia, la teor\'ia de \'algebras
fue usada para describir la teor\'ia de campos de clase, tanto local
como global. Con el uso de la cohomolog\'ia de grupos, tenemos 
los mismos resultados de manera mucho m\'as simple. En esta 
parte describimos la teor\'ia de los grupos de Brauer, pero
la descripci\'on del s\'imbolo residual de la norma lo haremos por
medio del uso de cohomolog\'ia. La obtenci\'on del s\'imbolo residual
de la norma lo haremos en el Teorema \ref{CCLT17.6.24}.

No daremos detalles de la siguiente discusi\'on. Consideremos
$E$ un campo cualquiera. Se define el {\em grupo de
Brauer} $\Br E\label{CClasegrupoBrauer}$
de $E$ como el conjunto de clases de 
$E$--isomorfismos de anillos de divisi\'on finito dimensionales
sobre $E$ y tales que $E$ es el centro del anillo de divisi\'on.

En general $\Br E$ tiene una estructura de grupo abeliano definido
por medio del producto tensorial de \'algebras sobre $E$ y es precisamente
con esta estructura que se llama {\em grupo de Brauer}.

\bigskip

Regresamos a nuestro estudio. Consideremos una extensi\'on $L/K$ no
ramificada de campos locales. Se tiene $f=[L:K]=[\bar L:\bar K]$, $e=1$
y $\Gal(L/K)\cong\Gal(\bar L/\bar K)$. Se tiene un generador de 
$\Gal(\bar L/\bar K)$ es el {\em automorfismo de
Frobenius\index{automorfismo de Frobenius}\index{Frobenius!automorfismo
de $\sim$}} $\varphi_{L|K}\in \Gal(L/K)$ con $\Gal(L/K)\stackrel{\cong}{\lra}
\Gal(\bar L/\bar K)$, $\varphi_{L|K}\longmapsto (x\to x^{|\bar K|}=:
x^{q_K})$, $\bar K={\ma F}_{q_K}$ (ver Subsecci\'on \ref{CClaseS1.3}).

Se tiene que $\varphi_{L|K} x\equiv x^{q_K}\bmod \pK_L$ para toda $x\in
\o_L$ y se verifica que $\varphi_{L|K}=\varphi_{M|K}|_L=\varphi_{M|K}
G_{M|L}\in G_{L|K}$ con $K\subseteq L\subseteq M$ y $\varphi_{M|L}
=\varphi_{M|K}^{[L:K]}$.

\begin{proposicion}\label{CCLP17.6.1}
Sea $L/K$ una extensi\'on finita de Galois de campos locales con
grupo de Galois $G$. Entonces existe un $G$-subm\'odulo $V$ de
$U_L^{(1)}$ de \'indice finito (por tanto $V$ tambi\'en es de
\'indice finito en $U_L$) tal que $V$ es cohomol\'ogicamente
trivial.
\end{proposicion}

\begin{proof}
Sea $\alpha\in L$ tal que $\{\sigma \alpha\}_{\sigma\in G}$ es una
base normal de $L/K$. Sea $a\in\*K$ tal que $a\sigma\alpha\in\o_L$
para toda $\sigma\in G$, esto es, seleccionamos $a\in\*K$ tal que
$v_K(a)\geq -v_L(\sigma\alpha)$ para toda $\sigma \in G$. Sea
$M:=\bigoplus_{\sigma\in G}(a\sigma\alpha)\o_K$. Entonces
$M\cong\o_K[G]\cong \o_k\otimes \zg$ es un $G$-m\'odulo inducido
y por tanto cohomol\'ogicamente trivial.

Se tiene que $M=\{(a\sigma\alpha\xi_{\sigma})_{\sigma\in G}\mid
\xi_{\sigma}\in \o_K\}=\{(m_{\sigma})_{\sigma\in G}\mid v_L(m_{\sigma})
\geq -v_L(a\sigma\alpha)\}=\{(m_{\sigma})_{\sigma\in G}\mid v_L(m_{\sigma})
> -v_L(a\sigma\alpha)-1\}=\{(m_{\sigma})_{\sigma\in G}\mid |m_{\sigma}|_L
< q^{v_L(a\sigma\alpha)+1}\}$. Se sigue que $M$ es abierto en $\o_L$.
Puesto que $\o_L$ es compacto, se sigue que $[\o_L:M]<\infty$ y
puesto que $\bigcap_{n=1}^{\infty}\pi_K^n\o_L=\{0\}$, existe $m\in
{\ma N}$ tal que $\pi_K^m\o_L\subseteq M$.

Sea $V_i:=1+\pi_K^{m+i}M$ el cual es un subm\'odulo de $U_L^{(1)}$
pues $M\subseteq \o_L$. El mapeo $V_i\ni v_i\stackrel{\mu}{\lra}
\pi_K^{-m-i}(v_i-1)\bmod \pi_K M$ define un mapeo biyectivo $V_i/V_{i+1}
\lra M/\pi_K M$ pues $\ker \mu=V_{i+1}$. Adem\'as $V_i/V_{i+1}\cong
M/\pi_K M\cong \big(\o_K/\pi_K\big)[G]\cong (\o_K/\pi_K) \otimes
\zg$, por lo que tanto $V_i/V_{i+1}$ como $M/\pi_K M$ son finitos
y cohomol\'ogicamente triviales.

Se tienen sucesiones exactas
\[
0\lra V_{i-1}/V_i\lra V_{i-2}/V_i\lra V_{i-2}/V_{i-1}\lra 0,
\]
donde $V_{i-1}/V_i$ y $V_{i-2}/V_{i-1}$ son finitos y cohomol\'ogicamente
triviales. Se sigue que $V_{i-2}/V_i$ es finito y cohomol\'ogicamente
trivial. Por inducci\'on obtenemos que $V_1/V_i$ es finito y
cohomol\'ogicamente trivial.

Se tiene que $V_1=\lim\limits_{\substack{\leftarrow\\ i}}V_1/V_i$. Por
tanto $\co qH{V_1}=\co qH{\lim\limits_{\substack{\leftarrow\\i}}V_1/V_i}=
\lim\limits_{\substack{\leftarrow\\ i}}\co qH{V_1/V_i}=0$ para toda
$q\in{\ma Z}$ y todo subgrupo $H$ de $G$. De esta froma obtenemos que
$V_1$ es cohomol\'ogicamente trivial, $V_1=1+\pi_K^{m+1} M$. Finalmente
$[U_L^{(1)}:V_1]\leq [\pK_L:\pK_L^{e(m+1)}][\o_L:M]<\infty$, esto es,
$V_1$ es de \'indice finito en $U_L^{(1)}$, y por tanto tambi\'en en $U_L$,
y $V=V_1$ es cohomol\'ogicamente trivial.
$\fin$
\end{proof}

Como corolario, obtenemos el llamado ``{\em axioma de la teor\'ia de
campos de clase locales\index{axioma de la teor\'ia de campos de 
clase locales}}.

\begin{teorema}[axioma de la teor\'ia de campos de clase locales\index{axioma
de la teor\'ia de campos de clase locales}]\label{axiomacamposlocales}
Sea $L/K$ una extensi\'on c\'iclica finita de campos locales con grupo
de Galois $G$. Entonces $\co 1G{\*L}=\{1\}$ y $\co 0G{\*L}$ tiene
cardinalidad $[L:K]=|G|$.
\end{teorema}

\begin{proof}
Por el Teorema 90 de Hilbert se tiene $\co 1G{\*L}=\{1\}$ para $G$
arbitrario, no necesariamente c\'iclico. Sea $V=V_1$ el subm\'odulo
de $U_L^{(1)}$ obtenido en la Proposici\'on \ref{CCLP17.6.1}.
Entonces, por ser $V$ cohmol\'ogicamente trivial, $h(V)=1$. Puesto
que $U_L^{(1)}/V$ es finito, $h(U_L^{(1)}/V)=1$ y puesto que
$1\lra V\lra U_L^{(1)}\lra U_L^{(1)}/V\lra 1$ es exacta, se sigue
que $h(U_L^{(1)})=1$. Puesto que $[U_L:U_L^{(n)}]<\infty$ para
toda $n\in{\ma N}$, se sigue que $h(U_L/U_L^{(n)})=1$. 

Ahora de la exactitud de la sucesi\'on $1\lra U_L^{(n)}\lra U_L\lra
U_L/U_L^{(n)}\lra 1$, se obtiene que $h(U_L)=h(U_L^{(n)})=1$
para toda $n\in{\ma N}$.

Por otro lado, $1\lra U_L\lra \*L\stackrel{v_L}{\lra} {\ma Z}\lra 0$
es exacta y por tanto $h(\*L)=h({\ma Z})=\frac{|\co 0G{\ma Z}|}
{|\co 1G{\ma Z}|}=\frac{|{\ma Z}/|G|{\ma Z}|}{|\{0\}|}=|G|=[L:K]$.
Puesto que $\co 1G{\*L}=1$, se sigue que $|\co 0G{\*L}|=|G|=
[L:K]$ y $\co 0G{\*L}=\frac{(\*L)^G}{\N_{L/K}\*L}=\frac{\*K}{
\N_{L/K}\*L}$. El resultado se sigue. Adem\'as, $\big|\co 0G{\*L}
\big|=[\*K\colon \N_{L/K}\*L]=|G|=[L:K]$.
$\fin$
\end{proof}

Un corolario importante, es el siguiente.

\begin{teorema}\label{CCLT17.6.3}
Sea $L/K$ una extensi\'on finita no ramificada de campos locales.
Entonces $\co q{G_{L|K}}{U^{(n)}_L}=\{1\}$ para toda $q\in{\ma Z}$
y para toda $n\in{\ma N}\cup\{0\}$. En particular, $\N_{L/K} U^{(n)}_L
=\unidades n$ para toda $n\in{\ma N}\cup\{0\}$.
\end{teorema}

\begin{proof}
Puesto que $L/K$ es no ramificada, se tiene para un elemento
primo $\pi_K$ de $K$, que
\[
v_L(\pi_K)=e(L|K)v_K(\pi_K)=1\cdot 1=1.
\]
Por tanto $\pi_K$ es un elemento primo de $L$.

Se tiene la sucesi\'on exacta
\[
1\lra U_L\lra \*L\stackrel{v_L}{\lra}{\ma Z}\lra 0.
\]
Por otro lado tenemos $\co {{-1}}G{\ma Z}=\co 1G{\ma Z}=\{0\}$
donde $G=\Gal(L/K)$. Por tanto se obtiene la sucesi\'on en
cohomolog\'ia
\[
\co {{-1}}G{\ma Z}=\{0\}\lra \co 0G{U_L}\lra \co 0G{\*L}
\stackrel{\bar v_L}{\lra} \co 0G{\ma Z}\lra\cdots,
\]
con $\bar v_L\colon\co 0G{\*L}=\*K/\N_{L/K}\*L\lra {\ma Z}/|G|
{\ma Z}=\co 0G{\ma Z}$, y $\bar v_L(\bar \pi_K)=\bar 1$. Por
tanto $\bar v_L$ es suprayectiva.

Ahora bien, como $L/K$ es c\'iclica, por el Teorema
\ref{axiomacamposlocales}, obtenemos que
$|\co 0G{\*L}|=[L:K]=|\co 0G{\ma Z}|$, de donde se sigue
que $\bar v_L$ es un isomorfismo. En particular
$\ker \bar v_L=\co 0G{U_L}=\{1\}$. Puesto que $h(U_L)=1$,
se sigue que $\co qG{U_L}=\{1\}$ para toda $q\in{\ma Z}$.

Se tiene la sucesi\'on exacta $1\lra U^{(1)}_L\lra U_L\lra 
\tilde L^*\lra 1$. Puesto que $\co {{-1}}G{\tilde L^*}=\{1\}$ y
$h(\tilde L^*)=1$ por $\tilde L^*$ finito, se sigue que $\tilde L^*$
es cohomol\'ogicamente trivial y por tanto $\co qG{U^{(1)}_L}
\cong \co qG{U_L}=\{1\}$ para toda $q\in{\ma Z}$.

En general, para cualquier $n\in{\ma N}$, se tiene la 
sucesi\'on exacta $1\lra U_L^{(n+1)}\lra U_L^{(n)}\lra
\tilde L\lra 0$ y $\tilde L$ es cohomol\'ogicamente trivial,
por lo que $\co qG{U_L^{(n+1)}}\cong \co qG{U_L^{(n)}}=
\{1\}$ para toda $q\in{\ma Z}$ y para toda $n\geq 1$.
$\fin$
\end{proof}

\begin{observacion}\label{CCLO17.6.4}
En general probaremos que $[U_K:\N_{L/K} U_L]=
e(L|K)$ donde $L/K$ es una extensi\'on abeliana
finita de campos locales (ver Corolario \ref{C17.2.15'}).
\end{observacion}

Consideremos una extensi\'on finita no ramificada de 
campos locales. Se tiene la sucesi\'on exacta
\[
1\lra U_L\lra \*L\stackrel{v_L}{\lra}{\ma Z}\lra 0,
\]
y puesto que $\co q{G_{L/K}}{U_L}=\{1\}$, $\co q{G_{L/K}}{\*L}
\cong \co q{G_{L/K}}{\ma Z}$. Sea $\bar v\colon\co 2{G_{L/K}}{
\*L}\lra \co 2{G_{L/K}}{\ma Z}$ el isomorfismo inducido por la
valuaci\'on.

Ahora, $0\lra{\ma Z}\lra {\ma Q}\lra {\ma Q}/{\ma Z}\lra 0$ es
exacta y ${\ma Q}$ es cohomol\'ogicamente trivial, de
donde obtenemos que el mapeo de conexi\'on $\delta\colon
\co q{G_{L/K}}{{\ma Q}/{\ma Z}}\lra \co {{q+1}}{G_{L/K}}{\ma Z}$
es un isomorfismo. Sea 
\[
\delta^{-1}\colon\co 2{G_{L/K}}{\ma Z}\lra
\co 1{G_{L/K}}{{\ma Q}/{\ma Z}}=\Hom(G_{L/K},{\ma Q}/{\ma Z})
=\chi(G_{L/K}),
\]
el dual de $G_{L/K}$, es decir, el grupo de caracteres de $G_{L/K}$. 

Si $\xi\in \chi(G_{L/K})$, se tiene
que si $\varphi_{L/K}$ denota al automorfismo de Frobenius
de $L/K$, entonces $\xi(\varphi_{L/K})\in \big(\frac{1}{[L:K]}
{\ma Z}\big)/{\ma Z}\subseteq {\ma Q}/{\ma Z}$.
Como $\varphi_{L/K}$ genera al grupo
$G_{L|K}$, el cual es de orden $[L:K]$, $\chi(G_{L|K})
\isomo\limits_{\varphi} \Big(\frac {1}{[L:K]} {\ma Z}\Big)/
{\ma Z}$. Esto es, $\co 1{G_{L|K}}{{\ma Q}/{\ma Z}}=\chi(
G_{L|K})\stackrel{\varphi}{\lra}\Big(\frac{1}{[L:K]}{\ma Z}\Big)/
{\ma Z}$. Entonces se define $\varphi$ para
$\xi\in\chi(G_{L|K})$, por $\varphi(\xi):=\xi(\varphi_{L|K})$.

La composici\'on de estos tres isomorfismos nos da el 
isomorfismo
\[
\co 2{G_{L|K}}{\*L}\stackrel{\bar v_L}{\lra}\co 2{G_{L|K}}
{\ma Z}\stackrel{\delta^{-1}}{\lra}\co 1{G_{L|K}}{{\ma Q}/{\ma
Z}}\stackrel{\varphi}{\lra}\Big(\frac{1}{[L:K]}{\ma Z}\Big)/{\ma Z}.
\]

\begin{definicion}\label{CCLD17.6.5}
Se define el {\em mapeo invariante\index{mapeo invariante}}
para una extensi\'on no ramificada de campos locales $L/K$
\[
\inv_{L|K}\colon \co 2{G_{L|K}}{\*L}\lra \Big(\frac{1}{[L:K]}{\ma
Z}\Big)/{\ma Z},
\]
por $\inv_{L|K}=\varphi\circ \delta^{-1}\circ \bar v_L$.
\end{definicion}

Sea $K_0$ un campo local y sea $\nr K_0$ la m\'axima extensi\'on
no ramificada de $K_0$, es decir, $\nr K_0=\bigcup_L L$, $L/K_0$
no ramificada.

\begin{definicion}\label{CCLD17.6.6}
El campo $\nr K_0$ se llama el {\em campo de inercia\index{campo
de inercia}} sobre $K_0$.
\end{definicion}

\begin{teorema}\label{CCLT17.6.7}
Sean $K_0\subseteq K\subseteq L\subseteq M\subseteq \nr K_0$.
Entonces
\las
\item\label{inv1} $\inv_{L|K}=\inv_{M|K}|_{H^2(L|K)}$.
\item\label{inv2} $\inv_{M|L}\circ \res_L=[L:K]\inv_{M|K}$.
\end{list}
\end{teorema}

\begin{proof}
Para (\ref{inv1}), debemos probar que el diagrama
\begin{gather*}
\xymatrix{
H^2(L|K)\ar@{->}[rr]^(.4){\inv_{L|K}}\ar@{^{(}->}[d]_i&&\Big(
\frac{1}{[L:K]}{\ma Z}\Big)/{\ma Z}\subseteq {\ma Q}/{\ma Z}
\ar@<-5ex>[d]^i\\
H^2(M|K)\ar@{->}[rr]^(.4){\inv_{M|K}}&&\Big(
\frac{1}{[M:K]}{\ma Z}\Big)/{\ma Z}\subseteq {\ma Q}/{\ma Z}
}
\intertext{es conmutativo. Para ello, se probar\'a que el diagrama}
\xymatrix{
H^2(L|K)\ar@{->}[r]^(.45){\bar v_L}\ar@{^{(}->}[d]_i&\co 2{G_{L|K}}{\ma Z}
\ar@{->}[r]^{\delta^{-1}}\ar@{->}[d]^{\infla}&\co 1{G_{L|K}}{{\ma Q}/{\ma Z}}
\ar@{->}[r]^{\varphi}\ar@{->}[d]^{\infla}&\Big(\frac{1}{[L:K]}{\ma Z}\Big)/{\ma Z}
\ar@{->}[d]^i\\
H^2(M|K)\ar@{->}[r]^(.45){\bar v_M}&\co 2{G_{M|K}}{\ma Z}
\ar@{->}[r]^(.5){\delta^{-1}}&\co 1{G_{M|K}}{{\ma Q}/{\ma Z}}
\ar@{->}[r]^{\varphi}&\Big(\frac{1}{[M:K]}{\ma Z}\Big)/{\ma Z}\
}
\end{gather*}
es conmutativo.

Se verifica que el cuadrado de la izquierda es conmutativo usando la
definici\'on de las valuaciones $v_L$, $v_M$ y del mapeo de inflaci\'on
$\infla$ en los $2$-cociclos. Que el cuadro de enmedio es conmutativo
se sigue de que inflaci\'on y los mapeos de conexi\'on conmutan. Para
el cuadro de la derecha, tenemos $\infla\colon\co 1{G_{L|K}}{{\ma Q}/
{\ma Z}}\lra \co 1{G_{M|K}}{{\ma Q}/{\ma Z}}$ se define en $1$-cociclos
$Z^1(G_{L|K},{\ma Q}/{\ma Z})=\Hom(G_{L|K},{\ma Q}/{\ma Z})$ (por
actuar $G_{L|K}$ en ${\ma Q}/{\ma Z}$ de manera trivial) as\'i:
$\infla(\chi)=\tilde\chi$ donde
\begin{gather*}
\xymatrix{
\frac{G_{M|K}}{G_{M|L}}=G_{L|K}\ar@{<-}[d]<3ex>_{\pi}\ar@{->}[r]^(.6){
\chi}&{\ma Q}/{\ma Z}\ar@{<-}[dl]^{\tilde\chi}\\ 
\phantom{\frac{G_{M|K}}{G_{M|L}}=}G_{M|K}
}\qquad
\tilde \chi(g)=\chi(g G_{M|L})=\chi(\bar g).\\
\intertext{Por tanto}
\begin{align*}
(\varphi\circ\infla)(\chi)&=\varphi(\tilde \chi)=\tilde \chi(\varphi_{M|K}),\\
(i\circ \varphi)(\chi)&=i(\varphi(\varphi_{L|K}))=i(\chi(\varphi_{L|K}))=
\chi(\varphi_{M|K}|_L)=\tilde \chi(\varphi_{M|K}).
\end{align*}
\end{gather*}
de donde se sigue la conmutatividad del diagrama de la derecha
probando (\ref{inv1}).

La afirmaci\'on en (\ref{inv2}) es equivalente a la conmutatividad del
diagrama
\begin{gather*}
\xymatrix{
H^2(M|K)\ar@{->}[rr]^(.4){\inv_{M|K}}\ar@{->}[d]_{\res_L}&&\Big(
\frac{1}{[M:K]}{\ma Z}\Big)/{\ma Z}\subseteq {\ma Q}/{\ma Z}
\ar@<-5ex>[d]^{[L:K]}\\
H^2(M|L)\ar@{->}[rr]^(.4){\inv_{M|L}}&&\Big(
\frac{1}{[M:L]}{\ma Z}\Big)/{\ma Z}\subseteq {\ma Q}/{\ma Z}
}
\intertext{Para esta conmutatividad basta probar la conmutatividad del
siguiente diagrama}
\xymatrix{
H^2(M|K)\ar@{->}[r]^(.45){\bar v_M}\ar@{->}[d]^{\res_L}&\co 2{G_{M|K}}{\ma Z}
\ar@{->}[r]^{\delta^{-1}}\ar@{->}[d]^{\res}&\co 1{G_{L|K}}{{\ma Q}/{\ma Z}}
\ar@{->}[r]^{\varphi}\ar@{->}[d]^{\res}&\Big(\frac{1}{[M:K]}{\ma Z}\Big)/{\ma Z}
\ar@{->}[d]^{[L:K]}\\
H^2(M|L)\ar@{->}[r]^(.45){\bar v_M}&\co 2{G_{M|L}}{\ma Z}
\ar@{->}[r]^(.5){\delta^{-1}}&\co 1{G_{M|L}}{{\ma Q}/{\ma Z}}
\ar@{->}[r]^{\varphi}&\Big(\frac{1}{[M:L]}{\ma Z}\Big)/{\ma Z}\
}
\end{gather*}

La conmutatividad del cuadro izquierdo se sigue de la definici\'on de 
$\bar v_M$ y $\res_L$ en los $2$-cociclos. La conmutatividad del cuadro
de enmedio se sigue que $\res$ conmuta con los mapeos de conexi\'on.
Para el cuadro de la derecha, tenemos
\begin{gather*}
\begin{array}{rcl}
\res\colon \co 1{G_{M|K}}{{\ma Q}/{\ma Z}}&\xrightarrow[]{\qquad}
&\co 1{G_{M|L}}{\ma Q}/{\ma Z}\\
\\
\chi&\xrightarrow[]{\ \res\ \ }&\chi|_{G_{M|L}}\\
\\
(\chi\colon G_{M|K}\lra {\ma Q}/{\ma Z})&\xrightarrow[]{\qquad}
& (\chi|_{G_{M|L}}\colon G_{M|L}\lra {\ma Q}/{\ma Z})
\end{array}\\
\\
\begin{align*}
(\varphi\circ \res)(\chi)&=\varphi\big(\chi|_{G_{M|L}}\big)=\big(\chi|_{G_{M|L}}\big)
(\varphi_{M|L})=\chi(\varphi_{M|L}),\\
[L:K]\varphi(\chi)&=[L:K]\chi(\varphi_{M|K})=\chi\big(\varphi_{M|K}^{[L:K]}\big)
=\chi(\varphi_{M|L}),
\end{align*}
\end{gather*}
de donde se sigue (\ref{inv2}).
$\fin$
\end{proof}

Puesto que $\nr K_0$ es ma m\'axima extensi\'on
no ramificada de $K_0$, $\nr K_0$
es la m\'axima extensi\'on no ramificada de cualquier extensi\'on $K/K_0$
no ramificada, $\nr K=\nr K_0$.

Ahora $\infla\colon H^2(L|K)\lra H^2(M|K)$ es inyectiva para $K_0\subseteq
K\subseteq L\subseteq M\subseteq \nr K_0$, por lo que podemos definir
\[
H^2(\nr K|K)=\lim_{\substack{\longrightarrow\\ L}}H^2(L|K)=\bigcup_L H^2(L|K).
\]

Por el Teorema \ref{CCLT17.6.7}, se obtiene un mapeo inyectivo:
\[
\inv_K\colon H^2(\nr K|K)\lra {\ma Q}/{\ma Z}
\]
 llamado el
{\em morfismo invariante\index{morfismo 
invariante}\index{invariante!morfismo $\sim$}} o {\em invariante
de Hasse\index{Hasse!invariante de $\sim$}\index{invariante
de Hasse}}. Adem\'as, $\inv_K$ es biyectivo pues ${\ma Q}/{\ma Z}
=\bigcup_{n=1}^{\infty}\big(\frac 1n{\ma Z}\big)/{\ma Z}$ y puesto
que para cada $n\in{\ma N}$, existe una \'unica extensi\'on no
ramificada $K_n/K$ de $K$ de grado $n$ y 
\[
H^2(K_n|K)
\xrightarrow[\cong]{\ \inv_{K_n|K}\ }\big(\frac 1n{\ma Z}\big)/{\ma Z}.
\]
Se sigue:

\begin{teorema}\label{CCLT17.6.8}
$H^2(\nr K|K)\cong {\ma Q}/{\ma Z}$. $\fin$
\end{teorema}

\subsection{S\'imbolo de la norma residual local}\label{SCClase3.2}

Sea $L/K$ una extensi\'on no ramificada finita, por lo que $\co 1{G_{L|K}}
{\*L}=\{1\}$ y $\co 2{G_{L|K}}{\*L}\isomo\limits_{\inv_{L|K}}\big(\frac {1}
{[L:K]}{\ma Z}\big)/{\ma Z}\cong {\ma Z}/([L\colon K]{\ma Z})$ es c\'iclica de
orden $[L:K]=|G_{L|K}|$, por lo que del Teorema de Tate-Nakayama, si
$u_{L|K}$ es la clase fundamental de $L/K$, se obtiene el isomorfismo
\begin{align*}
\theta_{L|K}:= u_{L|K}\Cup \underline{\ \ }\colon \co {{-2}}{G_{L|K}}{\ma Z}
&\lra \co 0{G_{M|K}}{\*L}\\
\abe G_{L|K}&\lra \*K/\N_{L/K} \*L,
\end{align*}
llamado el {\em mapeo de Nakayama\index{mapeo
de Nakyama}\index{Nakayama!mapeo de $\sim$}} o {\em isomorfismo de
Nakayama}. 

El isomorfismo
inverso de $\theta_{L|K}$, $\theta^{-1}_{L|K}\colon \*K/
\N_{L/K}\*L\lra \abe G_{L|K}=G_{L|K}$
se llama el {\em isomorfismo de reciprocidad\index{isomorfismo de
reciprocidad}}. Como consecuencia, tenemos un epimorfismo
$(\ , L/K)\colon \*K\lra G_{L|K}$ con n\'ucleo $\N_{L/K}\*L$
el cual se llama el {\em s\'imbolo de la norma residual\index{simbolo
de la norma residual@s\'imbolo de la norma residual}} y se tiene
que la sucesi\'on
\[
1\lra \N_{L/K}\*L\hooklongrightarrow \*K\xrightarrow[]{(\ ,L/K)}
G_{L|K}\lra 1
\]
es exacta. Notemos que $(a,L/K)=1\iff a$ es una norma de $\*L$.

El siguiente resultado ser\'a usado para hallar $(a,L/K)$ expl\'icitamente
en el caso no ramificado.

\begin{proposicion}\label{CCLP17.6.9}
Sean $L/K$ una extensi\'on finita no ramificada de campos locales,
$a\in\*K$, $\bar a=a\N_{L/K}\*L\in \co 0{G_{L|K}}{\*L}=H^0(L|K)$. Si
$\mu\in \chi(G_{L|K})=\co 1{G_{L|K}}{{\ma Q}/{\ma Z}}$ es un caracter
de $G_{L|K}$, entonces, si $\delta$ denota al mapeo de conexi\'on
$\delta\colon \co 1{G_{L|K}}{{\ma Q}/{\ma Z}}\lra \co 2{G_{L|K}}{\ma Z}$
obtenido de la sucesi\'on exacta $0\lra {\ma Z}\lra {\ma Q}\lra {\ma Q}/
{\ma Z}\lra 0$, se tiene
\[
\mu\big((a,L/K)\big)=\inv_{L|K}\big(\bar a\Cup \delta\mu\big)\in
\Big(\frac {1}{[L:K]}{\ma Z}\Big)/{\ma Z}\subseteq {\ma Q}/{\ma Z}.
\]
\end{proposicion}

\begin{proof}
Sea $\sigma_a:=(a,L/K)\in G_{L|K}\cong \co {{-2}}{G_{L|K}}{\ma Z}$ y sea
$\bar \sigma_a\in \co {{-2}}{G_{L|K}}{\ma Z}$ el elemento que corresponde
a $\sigma_a$. Se tiene $\bar a=u_{L|K}\Cup \bar \sigma_a\in \co 0{G_{
L|K}}{\*L}$. Ahora bien, el producto copa es asociativo y conmuta con
$\delta$, cualquier mapeo de conexi\'on, por lo que obtenemos
\[
\bar a\Cup\delta\mu=\big(u_{L|K}\Cup \bar\sigma_a\big)\Cup \delta\mu
=u_{L|K}\Cup \big(\bar\sigma_a\Cup \delta\mu\big)=u_{L|K}\Cup
\delta\big(\bar\sigma_a\Cup \mu\big).
\]

Puesto que para $\bar a_1\in \co 1{G_{L|K}}{{\ma Q}/{\ma Z}}$ y para $\bar\sigma\in
\co {{-2}}{G_{L|K}}{{\ma Z}}$ se tiene que $\bar a_1\Cup \bar\sigma=\overline{
a_1(\sigma)}\in \co {{-1}}{G_{L|K}}{{\ma Q}/{\ma Z}}$, se sigue que
\begin{gather*}
\bar \sigma_a\Cup \mu=\mu(\sigma_a)=\frac rn+{\ma Z}\in \Big(\frac 1n {\ma Z}
\Big)/{\ma Z}=\co {{-1}}{G_{L|K}}{{\ma Q}/{\ma Z}},\quad n=[L:K],
\intertext{para alg\'un $r$. Ahora, para $\delta\colon \co {{-1}}{G_{L|K}}{{\ma Q}/{\ma Z}}
\lra \co 0{G_{L|K}}{\ma Z}$ se obtiene}
\delta(\mu(\sigma_a))=n\big(\frac rn+{\ma Z}\big)=r+n{\ma Z}\in \co 0{G_{
L|K}}{\ma Z}={\ma Z}/n{\ma Z},
\end{gather*}
pues $\delta$ es el mapeo de conexi\'on que manda $\ker \N_n\colon
\Big({\ma Q}/{\ma Z}\lra {\ma Q}/{\ma Z}\Big)=\big\{\frac rn +{\ma Z}\mid
0\leq r\leq n-1\big\}$ a $ \co 0{G_{L|K}}{\ma Z}$ definido por $\frac rn+{\ma Z}
\stackrel{\delta}{\lra}r+n{\ma Z}$.

Ahora como el producto copa est\'a dado por el producto tensorial, tenemos
que $\bar a\Cup\delta \mu=u_{L|K}\Cup (r+n{\ma Z})=u_{L|K}^r$. 
Finalmente, 
\begin{gather*}
\inv_{L|K}\big(\bar a\Cup \delta\mu)=r\cdot\inv_{L|K}(u_{L|K})
=\frac rn+{\ma Z}=\mu(\sigma_a).
\tag*{$\fin$}
\end{gather*}
\end{proof}

El siguiente resultado es la expresi\'on expl\'icita para el s\'imbolo
de la norma residual en el caso no ramificado.

\begin{teorema}\label{CCLT17.6.10}
Sean $L/K$ una extensi\'on finita no ramificada de campos locales y
$a\in\*K$. Entonces $(a,L/K)=\varphi_{L|K}^{v_K(a)}$, donde $\varphi_{
L|K}$ es el automorfismo de Frobenius de $L/K$. En particular, si
$\pi_K$ es un elemento primo de $K$, $(\pi_K,L/K)=\varphi_{L|K}$.
\end{teorema}

\begin{proof}
Si $\mu\in\chi(G_{L|K})$, $\delta\mu\in \co 2{G_{L|K}}{\ma Z}\cong
\abe {G_{L|K}}=G_{L|K}$, $\bar a=a\N_{L/K}\*L\in H^0(L|K)$. Entonces,
por la Proposici\'on \ref{CCLP17.6.9}, se tiene
\begin{align*}
\mu\big((a,L/K)\big)&=\inv_{L|K}\big(\bar a\Cup \delta\mu\big)=(\varphi
\circ \delta^{-1}\circ \bar v_K)\big(\bar a\Cup \delta\mu\big)
\igual_{\substack{\uparrow\\ v_L=v_K\\ {\rm por\ ser\ no\ ramificada}}}\\
&=(\varphi\circ \delta^{-1})\big(v_K(\overline{a\otimes \delta\mu})\big).
\end{align*}

Se tiene que los mapeos de la sucesi\'on exacta $0\lra {\ma Z}
\stackrel{\iota}{\lra} {\ma Q}\stackrel{\pi}{\lra}{\ma Q}/{\ma Z}\lra 0$ son
los naturales, y por el Lema de la Serpiente, para $\mu\in\co 1G{{\ma
Q}/{\ma Z}}=\Hom (G,{\ma Q}/{\ma Z})$, se tiene $\delta\mu=(\iota^{-1}
\circ \*\partial_2\circ \pi^{-1})(\mu)=\*\partial_2(\mu)=\mu\circ \partial_2$.
Si $a\in\*K$, $\bar a=a\N_{L/K}\*L$, $\bar a\Cup \delta\mu=\overline{
a\otimes \delta\mu}$ est\'a representado por el $2$-cociclo $f\colon
G\times G\lra \*L$, $f(\sigma,\tau)=a^{\delta\mu(\sigma,\tau)}$.
Por tanto $v_K\big(\overline{a\otimes\delta\mu}\big)=v_K(a^{\delta
\mu})=v_K(a)\delta\mu$.

Por tanto 
\begin{align*}
\mu\big((a,L/K)\big)&=\varphi\circ \delta^{-1}(v_K(a)\delta\mu)
=\varphi(v_K(a)\mu)\\
&=v_K(a)\mu(\varphi_{L|K})=\mu\big(\varphi_{L|K}^{v_K(a)}\big).
\end{align*}

Por tanto, para toda $\mu\in\chi(G_{L|K})$, se tiene que $\mu\big(
(a,L/K)\big)=\mu\big(\varphi_{L|K}^{v_K(a)}\big)$, de donde se
sigue $(a,L/K)=\varphi_{L|K}^{v_K(a)}$.
$\fin$
\end{proof}

\begin{teorema}\label{CCLT17.6.11}
Sean $K$ un campo local, $\pi$ un elemento primo de $K$
y $L/K$ la extensi\'on no ramificada de grado $f$. Entonces
$\N_{L/K}\*L=\langle \pi^f\rangle\times U_K$.
\end{teorema}

\begin{proof}
Se tiene que $\langle \varphi_{L|K}\rangle=G_{L|K}$, $f=[\tilde L:
\tilde K]=[L:K]$ y $o(\varphi_{L|K})=f$. Se tiene que $a\in\*K$ 
pertenece a $\N_{L/K}\*L\iff (a,L/K)=\varphi_{L|K}^{v_K(a)}=1
\iff f|v_K(a)\iff a=u\pi^{lf}$ para alguna $u\in U_K$ y 
alguna $l\in{\ma Z}
\iff a\in \langle\pi^f\rangle\times U_K$.
$\fin$
\end{proof}

Se tiene $G_{\nr K|K}=\lim\limits_{\substack{\longleftarrow\\ L}}G_{L|K}$
donde $L$ recorre las extensiones finitas no ramificadas de $K$.
Para $a\in\*K$ obtendremos $(a,\nr K/K):=\lim\limits_{\substack{\longleftarrow\\
L}}(a,L/K)$ y obtendremos un homomorfismo $\*K\xrightarrow[]{(\underline{\ },
\nr K/K)} G_{\nr K|K}$.

Para ver esto, consideremos $\Lambda_L\colon G_{\nr K|K}\lra G_{L|K}$
la proyecci\'on can\'onica, esto es, el mapeo de restricci\'on $\sigma
\longmapsto \sigma|_L$. Entonces obtenemos
\[
\Lambda_L(a,\nr K/K)=(a,\nr K/K)|_L=(a,L/K)=\varphi_{L|K}^{v_K(a)}\in G_{L|K}.
\]
Puesto que $\varphi_{L|K}=\varphi_{M|K}|_L$ para $K_0\subseteq K
\subseteq L\subseteq M\subseteq \nr K_0$, el sistema $\big\{\varphi_{L|K}
\big\}_{L}$ forma un sistema coherente (sistema compatible) en el sistema
proyectivo y obtenemos
\[
\varphi_K:=\lim_{\substack{\longleftarrow\\ L}}\varphi_{L|K}\in G_{\nr K|K}.
\]
El homomorfismo $\varphi_K$ se llama el {\em automorfismo de
Frobenius universal de $K$\index{automorfismo de Frobenius
universal}\index{Frobenius!automorfismo de $\sim$ universal}} y
se tiene
\[
\varphi_K|_L=\varphi_{L|K}.
\]

\begin{teorema}\label{CCLTT17.6.12}
Si $a\in \*K$, entonces $(a,\nr K/K)=\varphi_K^{v_K(a)}$ y el n\'ucleo del
homomorfismo $\*K\xrightarrow[]{(\underline{\ },\nr K/K)} G_{\nr K|K}$ es el
grupo de unidades $U_K$ de $K$.
\end{teorema}

\begin{proof}
Para un extensi\'on finita no ramificada $L/K$, se tiene
\[
\Lambda_L(a,\nr K/K)=(a,L/K)=\varphi_{L|K}^{v_K(a)}=\Lambda_L(
\varphi_K^{v_K(a)}).
\]
Por tanto $(a,\nr K/K) = \varphi_K^{v_K(a)}$ y se tiene $(a,\nr K/K)=
\varphi_K^{v_K(a)}=1\iff v_K(a)=0$ pues $\varphi_K$ es de orden
infinito. Finalmente, $v_K(a)=0\iff a\in U_K$.
$\fin$
\end{proof}

\begin{observacion}\label{CCLTO17.6.13}
El mapeo $(\underline{\ },\nr K/K)\colon \*K\lra G_{\nr K|K}$ no es suprayectivo
puesto que $G_{\nr K|K}$ es un grupo profinito y $\im\big((\underline{\ },
\nr K/K)\big)\cong \*K/U_K\cong {\ma Z}$ el cual no es profinito por no ser
compacto. Por otro lado, $\langle \varphi_K\rangle$ es denso en 
$G_{\nr K|K}$ y la completaci\'on $\widehat{\langle\varphi_K \rangle}\cong
\widehat{\ma Z}\cong G_{\nr K|K}$.
\end{observacion}

\subsection{Ley de reciprocidad local en general}\label{SCClase3.3}

Sea $K_0$ un campo local y sea $\Omega:=\bar K_0$ una cerradura
separable de $K_0$. Para cada extensi\'on normal $L/K$ con $K_0
\subseteq K\subseteq L\subseteq \Omega$, y $L/K$ finita, nuevamente
usamos la notaci\'on $H^q(L|K):=\co q{G_{L|K}}{\*L}$ y $\Br K=
H^2(\ |K)=\bigcup_{L/K}H^2(L|K)$ donde $L$ recorre las extensions
normales, separables y finitas de $K$. El grupo $\Br K$ es el {\em grupo
de Brauer\index{grupo de Brauer}\index{Brauer!grupo de $\sim$}} de 
$K$. Sea $G_{K_0}:=\Gal(\Omega/K_0)$. Se tiene que $H^1(L|K)=
\{1\}$. Veamos que se satisface lo an\'alogo al caso no ramificado,
esto es, para cada extensi\'on normal finita $L/K$, $L\subseteq \Omega$,
existe un isomorfismo $\inv_{L|K}\colon H^2(L|K)\lra \Big(\frac 1{[L:K]}
{\ma Z}\Big)/{\ma Z}\subseteq {\ma Q}/{\ma Z}$ con las propiedades
de que si $K_0\subseteq K\subseteq L\subseteq M\subseteq \Omega$
es una torre de extensiones normales, $M/K_0$, entonces:
\las
\item $\inv_{L|K}=\inv_{M|K}|_{H^2(L|K)}$.
\item Si $M/K$ es normal, entonces $\inv_{M|L}\circ \res_L=
[L:K] \inv_{M|K}$.
\end{list}

El mecanismo para hacer lo anterior es extender el mapeo invariante
para extensiones no ramificadas a extensiones ramificadas.

Primero, probamos un resultado general de campos locales
(ver Corolario \ref{mpotencias}).
Sea $K$ un campo local y sea $m\in{\ma N}$. Se define
$\mu_m(K):=\{\xi\in K\mid \xi^m=1\}$ las $m$-ra\'ices de la
unidad contenidas en $K$. 

\begin{proposicion}\label{CCLTP17.6.14}
Sea $m$ tal que $\car K\nmid m$ (si $\car K=0$, $m$ es arbitrario).
El subgrupo abierto $(\*K)^m$ de $\*K$ tiene \'indice finito en $\*K$
y $[\*K:(\*K)^m]=mq^{v_K(m)} |\mu_m(K)|=m |m|^{-1}_{\pK}|\mu_m(K)|$
donde el campo residual de $K$, $\tilde K$,
satisface $\tilde K\cong \F$. Adem\'as
$[U_K:(U_K)^m]=q^{v_K(m)}|\mu_m(K)|$.
\end{proposicion}

\begin{proof}
Se tiene el cociente de Herbrand $q_{0,m}(\*K)=\frac{[\*K:(\*K)^m]}
{[\mu_m(K):1]}$. Puesto que $q_{0,m}$ es multiplicativo, obtenemos
\begin{gather*}
q_{0.m}(\*K)=q_{0,m}\big(\*K/U_K\big)q_{0,m}\big(U_K/\unidades n\big)
q_{0,m}(\unidades n).
\intertext{Ahora bien, $q_{0,m}\big(\*K/U_K\big)=q_{0.m}({\ma Z})=m$;
$q_{0,m}\big(U_K/\unidades n\big)=1$ por ser $U_K/\unidades n$
finito; y, para $n>v_K(m)$,}
q_{0,m}(\unidades n\big)=\frac{\big[\unidades n:\big(\unidades n\big)^m
\big]}{[\{1\}:\{1\}]}=\big[\unidades n:\unidades {n+v_K(m)}\big]=q^{v_K(m)},
\end{gather*}
pues $[\unidades n:\unidades {n+1}]=q$
(si $\car K=p>0$, $v_K(m)=0$ y $n\geq 1$).

Se sigue que $\big[\*K:(\*K)^m\big]=q_{o,m}(\*K)|\mu_m(K)|=
m\cdot 1\cdot q^{v_K(m)}\cdot |\mu_m(K)|$. Adem\'as, puesto
que $\mu_m(U_K)=\mu_m(K)$, y tomando $n>v_K(m)$,
\begin{align*}
q_{0,m}(U_K)&=\frac{\big[U_K:(U_K)^m\big]}{[\mu_m(U_K):\{1\}]}=
\big[U_K:(U_K)^m\big] |\mu_m(K)|^{-1}\\
&=q_{0,m}(U_K/\unidades n)
q_{0,m}(\unidades n)=1\cdot q^{v_K(m)}.
\tag*{$\fin$}
\end{align*}
\end{proof}

Para ver que se satisface la segunda condici\'on del Teorema de
Tate-Nakayama para tener un isomorfismo en el producto copa,
primero vemos una desigualdad.

\begin{proposicion}[Segunda desigualdad fundamental]\label{CCLTP17.6.15}
Sea $L/K$ una extensi\'on normal, separable y finita.
Entonces $|H^2(L|K)||[L:K]$.
\end{proposicion}

\begin{proof}
Primero supongamos que la extensi\'on $L/K$ es c\'iclica de grado primo
$p=[L:K]$. Sabemos que el cociente de Herbrand de $\*L$ es $p$ y como
$H^1(L|K)=\{1\}$ se sigue que $|H^2(L|K)|=p=[L:K]$.

Ahora sea $L/K$ normal, separable y finita. Puesto que $L/K$ es una 
extensi\'on de campos locales, tenemos que $G_{L|K}$ es soluble
(Corolario \ref{CClaseGR.6+1}).
Por tanto, existe una extensi\'on $K\subseteq K'\subseteq L$, $[K':K]$
es de orden primo y $K'/K$ es normal. Puesto que $H^1(L|K')=\{1\}$, la
sucesi\'on
\[
1\lra H^2(K'|K)\stackrel{\infla}{\lra} H^2(L|K)\stackrel{\res}{\lra} H^2(L|K'),
\]
es exacta. En particular $|H^2(L|K)|\big| |H^2(L|K)|\cdot |H^2(K'|K)|$ y por el
primer paso, $|H^2(K'|K)|=[K':K]$. Por inducci\'on, suponemos que
$|H^2(L|K)|\big| [L:K']$, de donde obtenemos que
$|H^2(L|K)|\big| [L:K'][K':K]=[L:K]$.
$\fin$
\end{proof}

\begin{teorema}\label{CCLT17.6.16}
Si $L/K$ es una extensi\'on normal, separable y finita de campos locales,
y si $L'/K$ es la extensi\'on no ramificada de grado $[L:K]=[L':K]$,
entonces $H^2(L|K)=H^2(L'|K)\subseteq H^2(\ |K)$.
\end{teorema}

\begin{proof}
Puesto que $|H^2(L'|K)|=[L':K]=[L:K]$ y $|H^2(L|K)|\big|[L:K]$, basta
probar que $H^2(L'|K)\subseteq H^2(L|K)$.

Sea $M:=LL'$. Entonces $M/K$ es una extensi\'on de Galois y $M/L$
es no ramificada.
\[
\xymatrix{
L'\ar@{-}[rr]\ar@{-}[d]_{\text{no ramificada}}&&
M=LL'\ar@{-}[d]^{\text{no ramificada}}\\
K\ar@{-}[rr]&&L
}
\]
Sea $c\in H^2(L'|K)\subseteq H^2(M|K)$ pues $H^1(M|L')=\{1\}$ y
por tanto $1\lra H^2(L'|K)\stackrel{\infla}{\lra}H^2(M|K)
\xrightarrow{\res_{L'}}H^2(M|L')$ es exacta.

De la sucesi\'on exacta $1\lra H^2(L|K)\stackrel{\infla}{\lra}H^2(M|K)
\xrightarrow{\res_L}H^2(M|L)$, se sigue que $c\in H^2(L|K)\iff
\res_L(c)=0$. Ahora bien, se tiene que $\res_L(c)=1\iff \inv_{M|L}
(\res_L(c))=1$ por ser $\inv_{M|L}$ por ser un isomorfismo. El
resultado se sigue de
\begin{gather}\label{CCLE17.6.17}
\inv_{M|L}(\res_L(c))=[L:K]\inv_{L'|K}(c),
\end{gather}
pues como $\inv_{L'|K}(c)\in\Big(\frac {1}
{[L':K]}{\ma Z}\Big)/{\ma Z}=
\Big(\frac {1}
{[L:K]}{\ma Z}\Big)/{\ma Z}$, se
sigue el resultado.
La prueba de (\ref{CCLE17.6.17}) se da manera m\'as general
en el Lema \ref{CCLL17.6.18}.
$\fin$
\end{proof}

\begin{lema}\label{CCLL17.6.18}
Sea $N|K$ una extensi\'on finita de Galois conteniendo a las dos
extensiones $L/K$ y $L'/K$, con $L'/K$ una extensi\'on no ramificada.
Entonces $M=LL'/L$ es no ramificada. Si $c\in H^2(L'|K)\subseteq
H^2(N|K)$, entonces $\res_L(c)\in H^2(M|L)\subseteq H^2(N|L)$ y
\[
\inv_{M|L}(\res_L(c))=[L:K]\inv_{L'|K}(c).
\]
\end{lema}

\begin{proof}
Se tiene que los $2$-cociclos de la  clase $\res_L(c)$ toman valores
en $\*M$ pues $1\lra H^2(L'|K)\xrightarrow{\infla}H^2(M|K)
\xrightarrow{\res_L}H^2(M|L)$ es exacta y $c\mapsto \infla c=c
\mapsto\res_L(c)\in H^2(M|L)=\co 2{G_{M|L}}{\*M}$.

Por tanto $\res_L(c)\in H^2(M|L)$. Sean $f$ el grado de inercia y $e$
el \'indice de ramificaci\'on de la extensi\'on $L/K$, la cual no
necesariamente es normal. Se extienden las valuaciones $v_K$
y $v_L$ a $N$. Se tiene $v_L=ev_K$. El mapeo $\inv$ es la 
composici\'on de $3$ isomorfismos $\bar v_K$, $\delta^{-1}$ y
$\varphi$. La f\'ormula se seguir\'a si probamos que el siguiente
diagrama es conmutativo:
\begin{scriptsize}
\[
\xymatrix{
H^2(L'|K)\ar@{->}[r]^(.45){\bar v_K}\ar@{^{(}->}[d]_i&\co 2{G_{L'|K}}{\ma Z}
\ar@{->}[r]^{\delta^{-1}}\ar@{->}[d]^{\infla}&\co 1{G_{L'|K}}{{\ma Q}/{\ma Z}}
\ar@{->}[r]^{\varphi}\ar@{->}[d]^{\infla}&\Big(\frac{1}{[L':K]}{\ma Z}\Big)/{\ma Z}\subseteq{\ma Q}/{\ma Z}\ar@{^{(}->}[d]^i\\
H^2(N|K)\ar@{->}[d]_{\res_L}&\co 2{G_{N|K}}{\ma Z}
\ar@{->}[d]^{e\res_L}&\co 1{G_{N|K}}{{\ma Q}/{\ma Z}}
\ar@{->}[d]^{e\res_L}&\Big(\frac{1}{[N:K]}{\ma Z}\Big)/{\ma Z}
\subseteq{\ma Q}/{\ma Z}\ar@{->}[d]^{[L:K]}\\
H^2(M|L)\ar@{->}[r]^(.45){\bar v_L}&\co 2{G_{M|L}}{\ma Z}
\ar@{->}[r]^(.5){\delta^{-1}}&\co 1{G_{M|L}}{{\ma Q}/{\ma Z}}
\ar@{->}[r]^{\varphi}&\Big(\frac{1}{[M:L]}{\ma Z}\Big)/{\ma Z}
\subseteq {\ma Q}/{\ma Z}
}
\]
\end{scriptsize}
donde los mapeos verticales inferiores simplemente mandan
las im\'agenes de los mapeos verticales superiores a los grupos de
cohomolog\'ia inferiores.

El cuadrado de la izquierda es conmutativo simplemente por el
comportamiento de los ciclos bajo los mapeos en cuesti\'on. El
de enmedio se sigue de que la inflaci\'on y la restricci\'on conmutan
con los homomorfismos de conexi\'on $\delta$. Finalmente para
el cuadrado de la derecha, se tiene $\varphi_{M|L}|_{L'}=\varphi_{L'/K}^f$,
donde $f=f(L|K)$ pues si $a\in L'$, entonces $\varphi_{M|L}(a)\equiv
a^{|\bar L|}\bmod \pK_M=a^{|\bar K|^f}\bmod \pK_M=a^{|\bar K|^f}
\bmod \pK_{L'}=\varphi_{L'|K}^f(a)$.

Sea $\mu\in \co 1{G_{L'|K}}{{\ma Q}/{\ma Z}}$, entonces
\begin{align*}
[L:K]\varphi(\mu)&=[L:K]\mu\big(\varphi_{L'|K}\big)=ef\mu\big(
\varphi_{L'|K}\big)=e\mu\big(\varphi_{L'|K}^f\big)=e\big(
\varphi_{M|K}|_{L'}\big)\\
&=e\infla \mu\big(\varphi_{M|L}\big)=e(\res\circ \infla)\mu\big(
\varphi_{M|L}\big)=e(\res\circ\infla)\varphi(\mu).
\end{align*}
Por tanto $[L:K]\varphi(\mu)=e(\res\circ\infla)\varphi(\mu)$
lo cual prueba el lema y el Teorema \ref{CCLT17.6.16}.
$\fin$
\end{proof}

\begin{corolario}\label{CCLC17.6.18}
Sea $K$ un campo local. Entonces el grupo de Brauer satisface
$\Br K\cong {\ma Q}/{\ma Z}$.
\end{corolario}

\begin{proof}
De la Teorema \ref{CCLT17.6.16} se sigue que
\begin{gather*}
\Br K=H^2(\ |K)=H^2(\nr K|K)=\bigcup_{\substack{L/K\\ \text{no
ramificas}}}H^2(L|K)\cong {\ma Q}/{\ma Z}.
\tag*{$\fin$}
\end{gather*}
\end{proof}

\begin{definicion}\label{CCLD17.6.19}
Sea $L/K$ una extensi\'on de Galois de campos locales y sea
$L'/K$ la extensi\'on no ramificada del mismo grado $[L:K]=
[L':K]$, de tal forma que $H^2(L|K)=H^2(L'|K)$. Se define el
{\em invariante\index{invariante}} o {\em{invariante de
Hasse}\index{Hasse!invariante de $\sim$}\index{invariante
de Hasse}} como
\[
\inv_{L|K}\colon H^2(L|K)\lra\Big(\frac 1{[L:K]}{\ma Z}/{\ma Z}
\Big)\subseteq {\ma Q}/{\ma Z}
\]
por el isomorfismo $\inv_{L|K}(c):=\inv_{L'|K}(c)$, $c\in
H^2(L|K)=H^2(L'|K)$.
\end{definicion}

\begin{teorema}\label{CCLT17.6.20}
El invariante de Hasse satisface:
\las
\item\label{primera} Si $K\subseteq L\subseteq M$ es una torre de extensiones 
de Galois, entonces
\[
\inv_{L|K}=\inv_{M|K}|_{H^2(L|K)}.
\]

\item\label{segunda} Si $K\subseteq L\subseteq M$ es una torre de extensiones con
$M/K$ de Galois, entonces
\[
\inv_{M|L}\circ\res_L=[L:K]\inv_{M|K}.
\]
\end{list}
\end{teorema}

\begin{proof}
\las
\item
Sean $M'/K$ y $L'/K$ las extensiones no ramificadas de grados
$[M':K]=[M:K]$ y $[L':K]=[L:K]$ respectivamente. Por la unicidad
de las extensiones no ramificadas, tenemos $K\subseteq L'
\subseteq M'$. Para $c\in H^2(L|K)=H^2(L'|K)$, se tiene
\[
\inv_{M|K}(c)=\inv_{M'|K}(c)=\inv_{L'|K}(c)=\inv_{L|K}(c).
\]

\item Consideremos $L/K$ una extensi\'on finita ($K_0\subseteq K
\subseteq L$). Sea $\res_L\colon H^2(\ |K)\lra H^2(\ |L)$.
Si $c\in H^2(\ |K)$, se puede suponer que $c\in H^2(L'|K)$
donde $L'/$ es no ramificada. As\'i, $M=LL'/L$ es no
ramificada y $\res_L(c)\in H^2(M|L)\subseteq H^2(\ |L)$.
Por el Lema \ref{CCLL17.6.18}, $\inv_L(\res_L(c))=
[L:K]\inv_K(c)$ lo cual prueba (\ref{segunda}). $\fin$
\end{list}
\end{proof}

\begin{definicion}\label{CCLD17.6.21}
Si $L/K$ es una extensi\'on de Galois, se define la {\em clase
fundamental\index{clase fundamental}} $u_{L|K}\in H^2(L|K)$
por
\[
\inv_{L|K}(u_{L|K})=\frac 1{[L:K]}+{\ma Z}\in
\Big(\frac 1{[L:K]}{\ma Z}\Big)/{\ma Z}\subseteq {\ma Q}/{\ma Z}.
\]
\end{definicion}

La parte m\'as importante del teorema principal de la teor\'ia
de campos de clase, es:

\begin{teorema}[Teorema general de reciprocidad]\label{CCLT17.6.22}
Sea $L/K$ es una extensi\'on finita de Galois de campos locales
con grupo de Galois $G$.
Entonces, el homomorfismo
\[
u_{L|K}\Cup\underline{\ }\colon \co q{G_{L|K}}{\ma Z}\lra
H^{q+2}(L|K),
\]
es biyectivo para todo $q\in {\ma Z}$.
\end{teorema}

\begin{proof}
Para que se cumpla el Teorema de Tate-Nakayama, se debe
tener que $\co 1G{\*L}=H^1(L|K)=\{1\}$ y $\co 1H{\*L}=\{1\}$
para todo subgrupo $H<G$ de $G$, lo cual es simplemente el
Teorema 90 de Hilbert.

Veamos la segunda condici\'on. Sea $H<G=\Gal(L/K)$ y sea
$E:=L^H$. Entonces $[L:E]=|H|$. Sea $F'/E$ la extensi\'on
no ramificada de grado $|H|$. Del Teorema \ref{CCLT17.6.16}
se tiene que $H^2(L|E)=H^2(F'|E)$ el cual es c\'iclico
de orden $|H|$.
$\fin$
\end{proof}

Aplicando el Teorema \ref{CCLT17.6.22} para $q=1, 2$, se obtiene:

\begin{corolario}\label{CCLC17.6.23} Sea $L/K$ una
extensi\'on finita de Galois de campos locales. Entonces
$H^3(L|K)=\{1\}$ y $H^4(L|K)=\chi(G_{L|K})=\Hom(G_{L|K},
{\ma Q}/{\ma Z})$.
\end{corolario}

\begin{proof}
Se tiene $H^3(L|K)\cong \co 1{G_{L|K}}{\ma Z}=\Hom(G_{L|K},{\ma Q}/{\ma Z})
=\{0\}$ por ser $G_{L|K}$ un grupo finito. Por otro lado, tenemos
$H^4(L|K)\cong \co 2{G_{L|K}}{\ma Z}\cong \abe G_{L|K}\cong
\co 1{G_{L|K}}{{\ma Q}/{\ma Z}}\cong \Hom(G_{L|K},{\ma Q}/{\ma Z})
=\chi(G_{L|K})$.
$\fin$
\end{proof}

Para $q=-2$, se obtiene

\begin{teorema}[Ley local de reciprocidad de Artin]\label{CCLT17.6.24}
Para una extensi\'on finita de Galois $L/K$ de campos locales, se tiene
el {\em isomorfismo de Nakayama\index{Nakayama!isomorfismo de
$\sim$}\index{isomorfismo de Nakayama}}
\begin{gather*}
\xymatrix{
\abe G_{L|K}\cong \co {{-2}}{G_{L|K}}{\ma Z}\ar@{->}[rr]^{u_{L|K}
\Cup \underline{\ \ }\ \ }\ar@/^2pc/@{->}[rr]^{\theta_{L|K} \text{\ mapeo de
Nakayama}}&&H^0(L|K)=\*K/\N_{L|K}\*L.} \tag*{$\fin$}
\end{gather*}
\end{teorema}

Como consecuencia del Teorema \ref{CCLT17.6.24}, se tiene una
demostraci\'on del Corolario \ref{C17.2.15'}.

\begin{proposicion}\label{CCLP17.6.24+1}
Si $L/K$ es una extensi\'on abeliana finita de campos locales,
se tiene que $e=[U_K:\N_{L/K}U_L]$.
\end{proposicion}

\begin{proof}
Por el Teorema \ref{CCLT17.6.24} y el mismo Corolario \ref{C17.2.15'},
se tiene que $ef=[L:K]=[\*K:\N_{L/K}\*L]=[U_K:\N_{L/K} U_L]f$, de donde
se sigue el resultado.
$\fin$
\end{proof}

\begin{corolario}\label{C17.5.33'}
Sea $L/K$ una extensi\'on finita de Galois de campos locales. Entonces
$L/K$ es no ramificada si y s\'olo si $U_K\subseteq \N_{L/K}\*L$.
\end{corolario}

\begin{proof}
Si $L/K$ es no ramificada, $U_K=\N_{L/K} U_L\subseteq \N_{L/K}
\*L$.

Rec\'iprocamente, supongamos $U_K\subseteq \N_{L/K}\*L$. Para
$x\in\* L$, $\N_{L/K} x\in U_K$ implica $x\in U_L$. Por tanto
$\N_{L/K} U_L=U_K$ y $e=[U_K\colon \N_{L/K}U_L]=1$. $\fin$
\end{proof}

El isomorfismo inverso induce la sucesi\'on exacta:
\[
1\lra \N_{L/K}\*L\hooklongrightarrow \*K\xrightarrow{(\ ,L/K)}{}
\abe G_{L|K}\lra 1,
\]
donde $(\ ,L/K)$ se llama el {\em s\'imbolo de la norma residual
local\index{simbolo de la norma residual local!s\'imbolo de la
norma residual local}}. Se tiene

\begin{teorema}\label{CCLT17.6.25} Sea $K_0\subseteq K
\subseteq L\subseteq M$ una torra de campos locales de tal
forma que $M/K$ es finita y de Galois. Los siguientes diagramas
son conmutativos
\lasa
\item\label{(a)} Si $L/K$ es normal
\[
\xymatrix{
\*K\ar@{->}[rr]^{(\ ,M/K)}\ar@{=}[d]&&\abe G_{M|K}
\ar@{->}[d]^{\pi}\\ \*K\ar@{->}[rr]_{(\ ,L/K)}&&\abe G_{L|K}}
\]
es decir, $(\ ,M/K)|_L=(\ ,M/L)$ y
donde $\pi$ es la proyecci\'on natural $\abe{\Gal(M/K)}
\stackrel{\pi}{\lra}\abe{\Gal(L/K)}$, equivalentemente, $\pi$
es mapeo restricci\'on $\abe{\Gal(M/K)}\xrightarrow[]{\rest|_L}
\abe{\Gal(L/K)}$, $\sigma\longmapsto \sigma|_L$.

\item\label{(b)}
\[
\xymatrix{
\*K\ar@{->}[rr]^{(\ ,M/K)}\ar@{^{(}->}[d]&&\abe G_{M|K}
\ar@{->}[d]^{\Ver}\\ \*L\ar@{->}[rr]_{(\ ,M/L)}&&\abe G_{M|L}}
\]
donde $\Ver$ es el mapeo de transferencia $\Ver\colon
\underbracket[0pt]{G/G'}_{\substack{\uigual\\ \abe G}}
\underbracket[0pt]{\lra}_{\substack{\phantom{\uigual}\\ \\ \lra}}
\underbracket[0pt]{H/H'}_{\substack{\uigual\\ \abe H}}$
que es $\res_{-2}\colon \co{{-2}}G{\ma Z}\lra \co {{-2}}H{\ma Z}$.

\item\label{(c)}
\begin{gather*}
\xymatrix{
\*L\ar@{->}[rr]^{(\ ,M/L)}\ar@{->}[d]_{\N_{L/K}}&&\abe G_{M|L}
\ar@{->}[d]^{\kappa}\\ \*K\ar@{->}[rr]_{(\ ,M/K)}&&\abe G_{M|K}}
\intertext{donde}
\xymatrix{
\Gal(M/L)\ar@{^{(}->}[r]^{i}\ar@{->}[d]_{\pi}&\Gal(M/K)
\ar@{->}[d]^{\pi}\\ \abe {\Gal(M/L)}\ar@{->}[r]_{\kappa}
&\abe {\Gal(M/K)}}
\end{gather*}
$\kappa$ inducido por $\cores_{-2}$.

\item\label{(d)} Si $\sigma\in\Gal(\sep K_0/K_0)=G_{K_0}$,
\[
\xymatrix{
\*K\ar@{->}[rr]^{(\ ,M/K)}\ar@{->}[d]_{\sigma}&&\abe G_{M|K}
\ar@{->}[d]^{\*\sigma}\\ \sigma\*K\ar@{->}[rr]_{(\ ,\sigma M/\sigma K)}&&
\abe G_{\sigma M|\sigma K}},
\]
es decir, $(\sigma a,\sigma M/\sigma K)=\sigma (a,M/K)
\sigma^{-1}$ para $a\in \*K$ y donde $\sigma \in G_{K_0}$, los
mapeos $\*K\stackrel{\sigma}{\lra}\sigma \*K$ y $\abe
G_{M|K}\stackrel{\*\sigma}{\lra}\abe G_{\sigma M|\sigma K}$
son $a\mapsto \sigma a$ y $\tau\mapsto \sigma\tau\sigma^{-1}$.
\end{list}
\end{teorema}

\begin{proof}{\ }

\lasa
\item[(a)] Sea $\mu\in \chi(G_{L|K})=\co 1{G_{L|K}}{{\ma Q}/{\ma Z}}$,
$\infla \mu\in \co 1{G_{M|K}}{{\ma Q}/{\ma Z}}$. Entonces
\begin{align*}
\mu(\pi(a,M/K))&=\infla \mu\big((a,M/K)\big)=\inv_{M|K}\big(
\bar a\Cup \delta(\infla \mu)\big)\\
&=\inv_{M|K}\big(\bar a\Cup \infla(\delta\mu)\big)=
\inv_{M|K}\big(\infla(\bar a\Cup (\delta\mu))\big)\\
&=\inv_{L|K}(\bar a\Cup \delta\mu)=\mu((a,L/K)).
\end{align*}
Puesto que esto se cumple para toda $\mu\in\chi(G_{L|K})$,
se sigue que $\pi(a,M/K)=(a,L/K)$.

\item[(b) y (c)] Se seguir\'an de la conmutatividad de los
diagramas
\begin{gather*}
\xymatrix{
\*K\ar@{->}[r]\ar@{^{(}->}[d]_i&H^0(M|K)
\ar@{<-}[rr]^(0.46){u_{M|K}\Cup\underline{\ }}_(0.46){\cong}
\ar@{->}[d]^{\res_0}&&\co {{-2}}{G_{M|K}}{\ma Z}
\ar@{->}[r]^(.60){\cong}_(0.60){\gamma}\ar@{->}[d]^{\res_{-2}}&\abe G_{M|K}
\ar@{->}[d]^{\Ver}\\
\*L\ar@{->}[r]&H^0(M|L)
\ar@{<-}[rr]^(0.46){u_{M|L}\Cup\underline{\ }}_(0.46){\cong}
&&\co {{-2}}{G_{M|L}}{\ma Z}
\ar@{->}[r]^(.60){\cong}_(.60){\gamma}&\abe G_{M|L}
}
\\
\xymatrix{
\*L\ar@{->}[r]\ar@{->}[d]_{\N_{L/K}}&H^0(M|L)
\ar@{<-}[rr]^(0.46){u_{M|L}\Cup\underline{\ }}_(0.46){\cong}
\ar@{->}[d]^{\cores_0}&&\co {{-2}}{G_{M|L}}{\ma Z}
\ar@{->}[r]^(.60){\cong}_(0.60){\gamma}\ar@{->}[d]^{\cores_{-2}}&\abe G_{M|L}
\ar@{->}[d]^{\kappa}\\
\*K\ar@{->}[r]&H^0(M|L)
\ar@{<-}[rr]^(0.46){u_{M|K}\Cup\underline{\ }}_(0.46){\cong}
&&\co {{-2}}{G_{M|K}}{\ma Z}
\ar@{->}[r]^(.60){\cong}_(0.60){\gamma}&\abe G_{M|K}
}
\end{gather*}

La conmutatividad de los cuadros de la izquierda es de inmediata
verificaci\'on siendo todos los mapeos naturales. La conmutatividad
de los cuadros de la derecha es simplemente la definici\'on de
$\Ver$ y de $\kappa$ como los mapeos correspondientes
a $\res_{-2}$ y $\cores_{-2}$ bajo los isomorfismos $\gamma$.

Para los cuadros intermedios, se tiene que si $z\in \co {{-2}}
{G_{M|K}}{\ma Z}$ (resp. $z'\in \co {{-2}}{G_{M|L}}{\ma Z}$)
por los comportamientos de $\cores$ y $\res$ con respecto al 
producto copa, de que
\begin{align*}
\res_0(u_{M|K}\Cup z)&=\res_0(u_{M|K})\Cup (\res_{-2} z)=
u_{M|L}\Cup (\res_{-2} z)
\intertext{y respectivamente}
\cores_0 (u_{M|L}\Cup z')&=\cores_0(\res_0(u_{M|L})\Cup z')=
u_{M|K}\Cup (\cores_{-2} z').
\end{align*}
\item[(d)] Se quiere ver la conmutatividad del diagrama
\begin{scriptsize}
\[
\xymatrix{
\*K\ar@{->}[r]\ar@{->}[d]_{\sigma}&H^0(M|K)
\ar@{<-}[rr]^(0.46){u_{M|K}\Cup\underline{\ }}
\ar@{->}[d]^{\*\sigma}&&\co {{-2}}{G_{M|K}}{\ma Z}
\ar@{->}[r]^(.60){\cong}\ar@{->}[d]^{\*\sigma}&\abe G_{M|K}
\ar@{->}[d]^{\*\sigma}\\
\sigma\*K\ar@{->}[r]&H^0(\sigma M|\sigma K)
\ar@{<-}[rr]^(0.46){u_{\sigma M|\sigma K}
\Cup\underline{\ }}
&&\co {{-2}}{G_{\sigma M|\sigma K}}{\ma Z}
\ar@{->}[r]^(.60){\cong}&\abe G_{\sigma M|\sigma K}
}
\]
\end{scriptsize}

Se tiene que $G_{\sigma M|\sigma K}=\sigma G_{M|K}
\sigma^{-1}$, por lo que $G_{M|K}\stackrel{\*\sigma}{\lra}
G_{\sigma M|\sigma K}$, $\tau\longmapsto \sigma \tau
\sigma^{-1}$ es un isomorfismo y por tanto los $q$-cociclos
del grupo $\co q{G_{\sigma M|\sigma K}}{\*\sigma A}$, $\*\sigma A
=A$, para cualquier $G$-m\'odulo son de la forma $\*\sigma
a_q=\sigma a_q\sigma^{-1}$ donde $a_q$ son los $q$-cociclos
de $\co q{G_{M|K}}A$. De esto se sigue la conmutatividad
del diagrama. $\fin$
\end{list}
\end{proof}

\begin{proposicion}\label{CCLP17.6.26}
Sean $L/K$ na extensi\'on de Galois finita de campos locales,
$a\in \*K$ y $\bar a=a\N_{L/K}\*L\in H^0(L|K)$. Si $\mu\in \chi\big(
\abe G_{L|K}\big)=\chi\big(G_{L|K}\big)=\co 1{G_{L|K}}{{\ma Q}/
{\ma Z}}$, se tiene $\mu(a,L/K)=\inv_{L|K}\big(\bar a\Cup \delta
\mu\big)\in \Big(\frac 1{[L:K]}{\ma Z}\Big)/{\ma Z}\subseteq
{\ma Q}/{\ma Z}$.
\end{proposicion}

\begin{proof}
Se sigue del caso no ramificado. $\fin$
\end{proof}

Consideremos el s\'imbolo residual de la norma (Teorema de
Reciprocidad) en general. Para cada extensi\'on abeliana finita
$L$ de $K$, se tiene el s\'imbolo residual de la norma $\*K
\xrightarrow{(\ ,L/K)}{} G_{L|K}$. Tomando el l\'imite proyectivo
sobre todas las extensiones abelianas finitas $L/K$,
$\abe G_K:=\lim\limits_{\substack{\longleftarrow \\ L/K}}G_{L|K}$,
se obtiene, para cada $a\in \*K$ el 
elemento $(a,K):=\lim\limits_{\substack{\longleftarrow \\ L/K}}
(a,L/K)\in \abe G_K$ el cual es el grupo de Galois
de la m\'axima extensi\'on abeliana de $K$\label{CClaseab}.

La demostraci\'on del siguiente resultado es pospuesto hasta
que estudiemos los grupos de Lubin-Tate. De momento presentamos
una demostraci\'on parcial para el caso $\car K=0$.

\begin{teorema}\label{CCLT17.6.27} El {\em s\'imbolo de la norma
residual universal\index{simbolo de la norma residual
universal!s\'imbolo de la norma residual universal}} define un
monomorfismo continuo $\*K \xhookrightarrow[]{\rho_K=(\ ,K)}
\abe G_K$.
\end{teorema}

\begin{proof}
Se tiene $\rho_K=(\underline{\ },K)\colon \*K \lra \abe G_K=\abe{\Gal(\sep K/K)}$,
$\rho_K(a)=(a,K):=\lim\limits_{\substack{\longleftarrow\\ L/K\\ \text{abeliana finita}}}
(a,L/K)$.

Se tiene que $\rho_K(a)=1 \iff (a,L/K)=1$ para toda extensi\'on $L/K$
abeliana finita $\iff a\in D_K:=\bigcap_{L/K}\N_{L/K}\*L$. En general
se probar\'a m\'as adelante que $D_K=\{1\}$. Para $K$ de caracter\'istica
$0$, se ver\'a que para cada $m\in{\ma N}$, $(\*K)^m$ es un grupo de
normas, esto es, existe una extensi\'on abeliana finita $L_m/K$ tal que
$\N_{L_m/K} (\*L_m)=(\*K)^m$. Se tiene que $[\*K:(\*K)^m]<\infty$ (Proposici\'on
\ref{CCLTP17.6.14}). Se sigue que $D_K\subseteq \bigcap_{m=1}^{\infty}
(\*K)^m=\{1\}$ (Corolario \ref{mpotencias}).

Si $H$ es un subgrupo abierto de $\abe G_K$, $H$ es \'indice finito, por lo que
$H$ es cerrado y $H$ corresponde a una extensi\'on abeliana finita de $L/K$:
$L=(\abe K)^H$ y $\rho_K^{-1}(H)=\N_{L/K}\*L$ el cual es cerrado de \'indice
finito en $\* K$, por lo que es abierto de $\* K$. Por tanto $\rho_K$ es 
continua.$\fin$
\end{proof}

Consideremos en general $a\in \*K$. Se tiene que la restricci\'on de $(a,K)\in 
\abe G_K$ al campo de inercia $\nr K/K$, $\nr K$ 
la m\'axima extensi\'on no ramificada
de $\nr K$, nos da
\[
(a,K)|_{\nr K}=(a,\nr K/K)=\varphi_K^{v_K(a)}\in G_{\nr K|K}
\]
donde $\varphi_K$ es el automorfismo de Frobenius universal.
En particular, $\rho_K=(\ ,K)$ \underline{no} es suprayectiva pues
$\langle \varphi_K\rangle$ es isomorfo a ${\ma Z}$ el cual no es profinito
pues no es compacto y de hecho, la completaci\'on de ${\ma Z}$ es
el anillo de Pr\"ufer $\hat {\ma Z}\cong \prod_p {\ma Z}_p$ y $\hat{\ma Z}
\subseteq \abe G_K$.

Antes de enunciar el resultado principal de la teor\'ia de campos de clase local,
notamos que cuando estudiamos la teor\'ia de campos de clase global, necesitamos
``una ley de reciprocidad para ${\ma R}$''. Resulta ser ${\ma R}$ tiene dos 
extensiones algebraicas $\{{\ma R},{\ma C}\}$. Se tiene $\co 2{G_{{\ma C}|{\ma R}}}
{\*{\ma C}}=\co 0{G_{{\ma C}|{\ma R}}}{\*{\ma C}}=\*{\ma R}/\N_{{\ma C}/{\ma R}}
\*{\ma C}=\*{\ma R}/{\ma R}^+\cong C_2$.

Se tiene que $a\in\*{\ma R}$ es norma $\iff a>0$ y el mapeo invariante
\begin{align*}
\inv_{{\ma C}|{\ma R}}\colon \co 2{G_{{\ma C}|{\ma R}}}{\*{\ma C}}&
\lra \Big(\frac 12{\ma Z}\Big)/{\ma Z}\\
a{\ma R}^+&\lra \begin{cases} 1+{\ma Z}, & a>0,\\
\frac 12+{\ma Z},&a<0,
\end{cases}
\end{align*}
y $(a,{\ma C}/{\ma R})(\sqrt{-1})=(\sqrt{-1})^{\sgn a}$, $(a,{\ma C}/{\ma R})(i)=
\begin{cases} i&\text{si $a>0$},\\
i^{-1}=i^3=-i&\text{si $a<0$}.\end{cases}$

El resultado principal de la Teor\'ia de Campos de Clase locales es:

\begin{teorema}[TCCL]\label{CClaseT3.2.2}\label{CClaseTCCL}
Sea $K$ un campo local o $K\in\{{\ma R}, {\ma C}\}$. 
Sea $\abe K$ la m\'axima extensi\'on abeliana de $K$. Entonces
\las
\item
Existe un \'unico homomorfismo continuo 
\[
\rho_K\colon \*K\to \Gal(\abe K/K)
\]
tal que
\l
\item
Si $L/K$ es una extensi\'on abeliana finita, $\rho_K$ induce un 
isomorfismo
\[
\*K/\N_{L/K}\*L\xrightarrow[\psi_{L/K}=
\widetilde{(\underline{\ },L/K)}]{\cong} \Gal(L/K),
\]
es decir, $\tilde{\rho}_K=\psi_{L/K}$. Se denota $\psi_{L/K}(a)=(a,L/K)$
por abuso del lenguaje.

\item
(Relaci\'on con los campos finitos). Si el campo residual de $K$ es
$\F$, se tiene el siguiente diagrama conmutativo
\[
\begin{CD}
\*K @>{\rho_K}>> \Gal(\abe K/K)\\
@V{\text{valuaci\'on}}V{v_K=v_{\pK}}V @VV{\mu}V\\
{\ma Z}@>>{\rho_{\F}}>\Gal(\abe \F/\F)\\
\end{CD}
\]
donde $\rho_{\F}$ es el mapeo $n\to\tau^n$ donde
$\tau$ es el automorfismo de Frobenius y $\mu$ es la composici\'on
\begin{eqnarray*}
\Gal(\abe K/K)&\xrightarrow{\rest}&\Gal(K^{\rm{nr}} /K)\cong \Gal(\abe \F/\F)\\
\sigma&\longmapsto& \sigma|_{K^{\rm{nr}} }
\end{eqnarray*}
donde $K^{\rm{nr}} $ es la m\'axima extensi\'on no ramificada de $K$
(la cual necesariamente tiene que ser abeliana por el Teorema 
{\rm{\ref{T17.3.3.1}}}).
\end{list}

\item {\underline{Teorema de Existencia\index{teorema de existencia}}}
La correspondencia $U\longmapsto \rho_K^{-1}(U)$ es una biyecci\'on
entre el conjunto de subgrupos abiertos de $\Gal(\abe K/K)$ y el conjunto
de subgrupos abiertos de {\'\i}ndice finito de $\*K$.

En particular, si $H\subseteq \*K$ es un subgrupo abierto de {\'\i}ndice 
finito, existe una \'unica extensi\'on abeliana finita $L/K$ tal que 
$H=\N_{L/K}\*K$. $\fin$
\end{list}
\end{teorema}

Falta probar el Teorema de Existencia, el cual se har\'a cuando estudiemos
los grupos formales de Lubin-Tate. La ley de reciprocidad ya la hemos obtenido
Teorema \ref{CCLT17.6.24}. Como complemento a la ley de reciprocidad
presentamos una forma alternativa de obtenerla, sin usar cohomolog\'ia
de grupos. Esta forma alternativa se debe a J\"urgen Neukirch \cite{Neu86}.

\subsection{Isomorfismo de Neukirch}\label{IsomorfismoNeukirch}

Hay una forma alternativa de obtener el mapeo de Nakayama, 
\[
\theta_{L/K}\colon \abe {\Gal(L/K)}\lra \*K/\N_{L/K}\*L
\]
 esto es, el inverso del mapeo de
reciprocidad. Esta forma alternativa se debe a J\"urgen Neukirch
\cite{Neu84, Neu86}.

Uno de los problemas en el coraz\'on de la teor\'ia de campos de clase, es la
Ley de Reciprocidad y en ella, una de las ideas centrales, es la estrecha
relaci\'on entre los elementos primos de los campos locales y los automorfismos
de Frobenius para extensiones finitas no ramificadas. M\'as precisamente, si
$L/K$ es una extensi\'on finita no ramificada y $\pi$ es un elemento primo de
$K$, entonces la ley de reciprocidad Artin establece que $(\pi, L/K)=\Fr LK$
donde, en esta parte, denotar\'a el automorfismo de Frobenius de la
extensi\'on $L/K$.

La idea del automorfismo de Neukirch es que dada una extensi\'on finita de
Galois $L/K$, se debe identificar $\sigma\in\Gal(L/K)$ con extensi\'on de un
automorfismo de Frobenius, lo que Neukirch define como ``levantamientos
de Frobenius'', y definir un mapeo ${\mc N}_{L/K}$
entre $\Gal(L/K)$ y $\*K/\N_{L/K}
\*L$ identificando $\sigma$ con la norma $\N_{E/K} \pi_E$ de un elemento
primo $\pi_E$ de cierto campo $E$.

Esta n\'itida idea funciona a la perfecci\'on. El problema es t\'ecnico.
La idea es relativamente simple y la definici\'on 
de ${\mc N}_{L/K}$ tambi\'en.
Sin embargo la verificaci\'on, aunque relativamente elemental, es
t\'ecnicamente dif\'icil. Hay que verificar que el mapeo est\'a bien
definido y que es un homomorfismo. Una vez superado el problema 
t\'ecnico, se obtiene una serie de propiedades funtoriales de este mapeo,
importante de por s\'i para toda la teor\'ia de campos de clase, adem\'as
que permiten extender el hecho de
 que ${\mc N}_{L/K}$ es en realidad un isomorfismo para
extensiones finitas no ramificadas, al caso general de un extensi\'on
finita de Galois $L/K$:
\[
{\mc N}_{L/K}\colon \abe{\Gal(L/K)}\stackrel{\cong}{\lra} \*K/\N_{L/K}\*L.
\]
En otras palabras, si la equivalencia funciona bien para extensiones
no ramificadas, la equivalencia se extiende de manera \'unica a 
extensiones finitas de Galois arbitrarias.

Neukirch desarrolla la teor\'ia de campos de clase a partir de este
isomorfismo (\cite{Neu84, Neu86}) evitando la cohomolog\'ia, como fue
su trabajo de 1969 (\cite{Neu69}). Aqu\'i hacemos notar que ambos 
puntos de vista tienen como base las extensiones no ramificadas y en
ambos casos se extiende a extensiones finitas de Galois.

El desarrollo del isomorfismo de Neukirch, lo presentamos directamente
para campos locales, en lugar de hacerlo para ``formaci\'on de clases''
lo cual es m\'as general.

El desarrollo abstracto del isomorfismo de Neukirch (y de hecho para
toda la teor\'ia de campos de clase usando de cohomolog\'ia de grupos),
es como sigue:

Sea $G$ un grupo profinito y $A$ un $G$-m\'odulo. Los subgrupos
cerrados de $G$ se denotan por $G_K$ y los \'indices se les consideran
campos ($G$ se piensa como $G=\Gal(\sep k/k)$ con $k$ un campo
local fijo). De esta forma $K$ juega el papel del campo fijo de $G_K$:
$K:=A^{G_K}$ y en particular $G_{\sep k}=\{1\}$ y $G_k=G$ y $k$
juega el papel del campo base. Definimos $K\subseteq L\iff G_L
\subseteq G_K$ y el grado de la ``extensi\'on'' $L/K$, como $[L\colon K]=
[G_K\colon G_L]$. Si $G_L\normal G_K$, $L/K$ se llama extensi\'on o
de Galois, etc.

\bigskip

Sea $K$ un campo local y sea $\sep K$ una cerradura separable
fija de $K$. Sea $\nr K$ la m\'axima extensi\'on no ramificada de
$K$ contenida en $\sep K$. Recordemos que para cada $n\in {\ma N}$,
existe una \'unica extensi\'on $K_n$ no ramificada de $K$ de 
grado $n$ y $\Gal(K_n/K)\cong {\ma Z}/n{\ma Z}$. Por tanto $\nr K=
\bigcup_{n=1}^{\infty} K_n$ y $\Gal(\nr K/K)=\Gal\Big(\big(\bigcup_{n=1}^{
\infty}K_n\big)/K\Big)=\Gal\Big(\big(\lim\limits_{\longrightarrow n}K_n\big)/K\Big)=
\lim\limits_{\longleftarrow n}\Gal(K_n/K)=\lim\limits_{\longleftarrow n}\big(
{\ma Z}/n{\ma Z}\big)\cong \hat{\ma Z}$, el anillo de Pr\"ufer.

Se tiene el mapeo de restricci\'on $\rest_K\colon\Gal(\sep K/K)\lra 
\Gal(\nr K/K)$. Cada grupo $\Gal(K_n/K)$ est\'a generado por el
automorfismo de Frobenius. De esta forma, $\Gal(K_n/K)=\langle
\Fr {K_n}K\rangle$ donde, si el campo residual de $K$ es $\F$,
el campo residual de $K_n$ es ${\ma F}_{q^n}$ y el automorfismo
de Frobenius es $\Fr {K_n}K(\alpha)=\alpha^q$. Entonces $\Gal(\nr K/K)$
est\'a generado topol\'ogicamente por $\Fro K\colon \overline{\F}\lra \F$, 
$\alpha\lra \alpha^q$ y $\Fro K|_{K_n}=\Fr {K_n}K$, $\Gal(\nr K/K)=
\overline{\langle \Fro K\rangle}\cong \hat {\ma Z}$.

\begin{proposicion}\label{NeuP1}
Sea $L/K$ una extensi\'on algebraica. Entonces $\nr L=L\nr K$.
\end{proposicion}

\begin{proof}
Se tiene que $L$ y $\nr K\subseteq \nr L$, por tanto $L\nr K\subseteq
\nr L$. Recordemos que cada extensi\'on finita no ramificada $M/E$
de un campo local $E$ est\'a un\'ivocamente determinada por la 
extensi\'on de campos residuales, es decir, si $\tilde M=\tilde E(\tilde 
\beta)$, existe $\beta\in \tilde \beta$ con $M=E(\beta)$ y
$\Irr(x,\beta,E)=\Irr(\bar x,\tilde \beta,\tilde E)$ y donde $\tilde \beta
\in \o_M/\pK_M$, $\tilde \beta =\beta+\pK_M$. Por tanto, si probamos
que $\nr {\tilde L}=\nr {\tilde K}\tilde L$, entonces, dada $F/L$ no
ramificada, para cada $\alpha\in F$, $L(\alpha)/L$ es una extensi\'on
finita no ramificada, $\widetilde{L(\alpha)}=\tilde L(\tilde \alpha)$
y tenemos $\widetilde{L(\alpha)}=\tilde L\widetilde
{K(\beta)}$ pues $\widetilde{L(\alpha)}\subseteq \nr{\tilde L}=\nr{\tilde K}
\tilde L$. Se sigue $K(\beta)\subseteq \nr K$. Se sigue $L(\alpha)=
K(\beta)L\subseteq \nr K L$. Por tanto $F\subseteq \nr K L$.
\[
\xymatrix{
K(\beta)\ar@{--}[rr]\ar@{-}[dd]&&\widetilde{K(\beta)}\ar@{-}[r]&
\widetilde{L(\alpha)}\ar@{--}[r]\ar@{-}[d]&L(\alpha)\ar@{-}[d]\\
&&\tilde H=\widetilde{\sep K}\cap \tilde L\ar@{-}[r]&\tilde L\ar@{--}[r]&L\\
K\ar@{--}[r]&\tilde K\ar@{-}[ru]
}
\]
Ahora, $\widetilde{L\nr K}\supseteq \tilde L\widetilde{\nr K}=\tilde L
\widetilde{\sep K}=\tilde L \widetilde{\sep L}$ puesto que $\tilde L/
\tilde K$ es separable. Se sigue que $\widetilde{\nr K}=\widetilde{
\sep K}=\widetilde{\sep L}=\widetilde{\nr L}$.
Por tanto $\widetilde{L\nr K}\supseteq
\tilde L\widetilde{\sep L}=\widetilde{\sep L}=\widetilde{\nr L}$ de
donde obtenemos que $\widetilde{L\nr K}=\widetilde{\nr L}$
y finalmente $L\nr K=\nr L$.
$\fin$
\end{proof}

Introducimos alguna notaci\'on. Sea $k$ un campo local fijo y sea
$\sep k$ una cerradura separable fija de $k$. Todos los campos
considerados estar\'an contenidas en $\sep k$ y usualmente ser\'an
extensiones finitas de $k$ y por supuesto separables. Sea $G:=
\Gal(\sep k/k)$. Si $K$ es una extensi\'on de $k$, sea $G_K:=
\Gal(\sep k/K)$. Entonces $G_K$ es un subgrupo cerrado de
$G$. Si $L/K$ es una extensi\'on, entonces $[L:K]=[G_K:G_L]$
en caso de ser finita. Se tiene $K=(\sep k)^{G_K}$ y $L=(
\sep k)^{G_L}$.
\[
\xymatrix{
\sep k\ar@{-}[d]\ar@{-}@/_1pc/[ddd]_G\ar@{-}@/^/[d]^{G_L}
\ar@{-}@/^2.5pc/[dd]^{G_K}\\
L\ar@{-}[d]\\K\ar@{-}[d]\\k
}
\]
Consideremos $\Gal(\nr K/K)=\overline{\langle \Fro K\rangle}$. Si
$\sigma \in G$, $\sigma|_{\nr k}\in\Gal(\nr k/k)\isomo\limits_{\mu}
\hat{\ma Z}$.

Se tiene que $\sigma|_{\nr k}=\Fro K^{\alpha_{\sigma}}$ con
$\alpha_{\sigma}\in \hat{\ma Z}$. Sea 
$\deg \colon G\lra\hat{\ma Z}$,
definido por $\deg(\sigma):=\alpha_{\sigma}$ el cual es un 
epimorfismo continuo. De esta forma, si $H$ es un subgrupo
cerrado de $\hat{\ma Z}$, $\mu^{-1}(H)$ es cerrado en $\Gal
(\nr k/k)$ y corresponde a un campo $T$ con $k\subseteq T
\subseteq \nr k$. Entonces $\Gal(\sep k/T)=\deg^{-1}(H)=G_T$
y $\deg^{-1}(H)$ es un subgrupo cerrado. 

Se tiene que
$\ker\deg =\Gal(\sep k/\nr k)$. Se define $\deg_K:=\frac {1}{f_K}
\deg\colon G_K\lra \hat{\ma Z}$ y $\deg_K\colon \Gal(\nr K/K)
\stackrel{\cong}{\lra} \hat{\ma Z}$, donde $f_K=f(K|k)$ es el grado
de inercia.
\[
\xymatrix{
\sep k\ar@{-}@/_1pc/[dd]_{\deg^{-1}(H)}\ar@{-}[d]\\
\nr k\ar@{-}[d]\ar@{-}@/^/[d]^{\mu^{-1}(H)}\\
T\ar@{-}[d]\\K
}\qquad
\xymatrix{
\sep k\ar@{=}[rr]\ar@{-}[d]&&\sep K\ar@{-}[d]\\ \nr k\ar@{-}[rr]\ar@{-}[d]
\ar@{-}@/_2pc/[dd]_{\hat{\ma Z}}
&&\nr K=\nr k K\ar@{-}[d]\ar@{-}@/^/[d]^{\hat{\ma Z}}\\
\nr k\cap K\ar@{-}[rr]^{\text{totalmente}}_{\text{ramificada}}
\ar@{-}[d]^{f_K}&&K\\ k
}
\]

El automorfismo de Frobenius sobre $K$, $\Fro K$, el generador
topol\'ogico de $\Gal(\nr K/K)$ satisface $\deg_K(\Fro K)=1$. Si 
$L/K$ es una extensi\'on finita, se tiene
\[
\xymatrix{
\nr K\ar@{-}[d]\ar@{-}[d]\ar@{-}[r]&\nr KL=\nr L\ar@{-}[d]\\
L_0=\nr K\cap L\ar@{-}[d]^{f(L|K)=\frac{f_L}{f_K}:=f_{L|K}}
\ar@{-}[r]&L\\ K
}
\]
Si $M/N$ es una extensi\'on no ramificada, $\Fr MN$ denota
al automorfismo de Frobenius de la extensi\'on $M/N$.
De esta forma $\Fr {L_0}K=\Fro K|_{L_0}$. Se tiene el diagrama
conmutativo
\begin{gather*}
\xymatrix{
G_L\ar@{->}[r]^{\deg_L}\ar@{^{(}->}[d]_i&\hat{\ma Z}\ar@{->}[d]^{f_{L|K}}\\
G_K\ar@{->}[r]_{\deg_K}&\hat{\ma Z}
}
\qquad
\xymatrix{
\nr K\ar@{-}[r]\ar@{-}[d]&\nr L\ar@{-}[d]\\ L_0=\nr K\cap L
\ar@{-}[r]\ar@{-}[d]&L\\ K
}
\intertext{$f_{L|K}\circ \deg_L=\frac{f_L}{f_K}\frac{1}{f_L}\deg =\frac 1{f_K}\deg=
\deg_K=\deg_K\circ i$. En particular $\Fro L|_K=\Fro K^{f_{L|K}}$.
Se tiene}
\xymatrix{
&L\ar@{-}[r]\ar@{-}[d]&\nr K L=\nr L\ar@{-}[d]\\ K\ar@{-}[r]&
\underbracket[0pt]{L_0}_{\substack{\uigual\\ \nr K\cap L}}\ar@{-}[r]&\nr K
}
\qquad
\xymatrix{
G_K\ar@{->}[d]_{\rest|_{\nr L}}\ar@{-->}[dr]^{\deg_K}\\
\Gal(\nr L/K)\ar@{->}[r]\ar@{->}[d]_{\rest|_{\nr K}}&
\hat{\ma Z}\\ \Gal(\nr K/K)\ar@{->}[ru]_{\deg_K \text{\ suprayectiva}}
}
\end{gather*}
Entonces $\deg_K\colon \Gal(\nr L/K)\lra \hat{\ma Z}$ es suprayectiva.

El diagrama
\begin{gather*}
\xymatrix{
G_L\ar@{^{(}->}[d]\ar@{->}[rr]^{\deg_L}&&\hat{\ma Z}\ar@{->}[d]^{f_{L|K}
=f(L|K)=\frac{f_L}{f_K}}\\ G_K\ar@{->}[rr]^{\deg_K}&& \hat{\ma Z}
}
\intertext{se puede reinterpretar como el diagrama conmutativo}
\xymatrix{
\Gal(\sep L/L)\ar@{^{(}->}[d]\ar@{->}[rr]^{\gamma_L=\rest|_{\nr L}}&&
\Gal(\nr L/L)\cong \hat{\ma Z}\ar@{->}[d]^{f_{L|K}}\\
\Gal(\sep K/K)\ar@{->}[rr]^{\gamma_K=\rest|_{\nr K}}&&
\Gal(\nr K/K)\cong\hat{\ma Z}
}
\end{gather*}

Se define $\Lambda(\nr L/L)=\{\tilde \sigma\in \Gal(\nr L/K)\mid \deg_K(
\tilde \sigma)\in{\ma N}\}$. Hacemos notar que se pide que el grado $\deg_K$
sea estrictamente positivo. Se tiene que si $\tilde\sigma\in\Lambda(\nr L/K)$
y si $n=\deg_K(\tilde\sigma)$, entonces $\tilde\sigma|_{\nr K}=\Fro K^n$.

La idea central del isomorfismo de Neukirch es, en cierta forma, convertir
cualquier elemento en un Frobenius. El resultado siguiente es el que
indica el camino a seguir.

\begin{proposicion}\label{NeuP2} El mapeo restricci\'on $\Lambda(\nr L/K)
\lra \Gal(L/K)$, $\tilde\sigma\longmapsto \tilde\sigma|_L$, es suprayectivo.
Si $\sigma=\tilde\sigma|_L$, $\tilde\sigma$ se llama un {\em levantamiento
de Frobenius\index{levantamiento de Frobenius}\index{Frobenius!levantamiento
de $\sim$}} de $\sigma$.
\end{proposicion}

\begin{proof}
Sea $\sigma\in \Gal(L/K)$, entonces $\sigma|_{L_0}=(\Fr {L_0}K)^n$ para
alguna $n\geq 0$. Notemos que $\Gal(L_0/K)=\langle \Fr{L_0}K\rangle=
\{\Fr{L_0}K,\Fr{L_0}K^2,\ldots,\Fr{L_0}K^{[L_0:K]}\}$, por lo que $n\geq 1$.
Sea $\tilde\mu$ una extensi\'on de $\Fro K$ a $\nr L$. Se tiene $\Gal(\nr L/K)
\stackrel{\pi}{\twoheadrightarrow}\Gal(\nr K/K)$, $\ker \pi=H=\Gal(\nr L/\nr K)$
y $\pi(\tilde\mu)=\Fro K$. $\xymatrix{\nr K\ar@{-}[r]^H\ar@{-}[d]&\nr L
\ar@{-}[dl]\\K}$

Entonces $\tilde\mu|_{\nr K}=\Fro K$, $\tilde\mu|_{L_0}=\Fro K|_{L_0}=
\Fr {L_0}K$. Por tanto $\sigma\tilde\mu^{-n}|_{L_0}=1$, es decir, $\sigma
\tilde\mu^{-n}|_L\in\Gal(L/L_0)\cong \Gal(\nr L/\nr K)$. En particular,
$\sigma\tilde\mu^{-n}|_L=\tau|_L$ para alg\'un $\tau\in\Gal(\nr L/\nr K)$.

Sea $\tilde\sigma=\tau\cdot\tilde\mu^n\in\Gal(\nr L/K)$, $\tilde\sigma|_L=\tau|_L
\cdot\tilde\mu^n|_L=\sigma$ y $\tilde\sigma|_{\nr K}=\underbracket[0pt]{\tau|_{\nr K}}_{
\substack{\uigual\\ \Id}}\cdot\tilde\mu^n|_{\nr K}=\Fro K^n$.
Por tanto $\deg_K(\tilde\sigma)=n\in{\ma N}$ y $\tilde \sigma$ es un 
levantamiento de Frobenius de $\sigma$.
$\fin$
\end{proof}

El punto medular de los levantamientos de Frobenius est\'a en el hecho
de que resulta que todo elemento $\sigma\in \Gal(L/K)$ puede 
considerarse ``casi'' un Frobenius.

\begin{proposicion}\label{NeuP3} Sea $\tilde\sigma\in\Lambda(\nr L/K)$
y sea $E$ el campo fijo de $\tilde\sigma$. Entonces
\las
\item $[E:K]<\infty$.
\item $f_{E|K}=f(E|K)=\deg_K(\tilde\sigma)$.
\item $\nr E=\nr L$.
\item $\tilde\sigma=\Fro E$.
\end{list}
\end{proposicion}

\begin{proof}
\[
\xymatrix{
&&\nr E=\nr L\ar@{-}[dl]^{\tilde\sigma}\\ & E\ar@{-}[dl]\\K}
\qquad
\xymatrix{
&\nr K\ar@{-}[r]\ar@{-}[d]&\nr L\ar@{-}[d]^{\tilde\sigma}\\
&E_0\ar@{-}[r]&E\\ K\ar@{-}[ru]
}
\]
\las
\item[(2)] Sea $E_0=E\cap \nr K$. Entonces $E_0$ es el campo fijo de
$\tilde\sigma|_{\nr K}=\Fro K^{\deg_K(\tilde\sigma)}$. Por tanto $f_{E|K}
=f(E|K)=[E_0:K]=\deg_K(\tilde\sigma)$.

\item[(1)] $[E:E_0]=[\nr K E:\nr K]\leq [\nr L:\nr K]<\infty$. Puesto que
\[
\xymatrix{
\nr K\ar@{-}[r]\ar@{-}[d]&\nr KE\ar@{-}[d]\\E_0\ar@{-}[r]&E}
\qquad [E_0:K]=f_{E_0|K}=f_{E|K}<\infty,
\]
se sigue que $[E:K]<\infty$.

\item[(3)] $\Gal(\nr L/E)$ est\'a generado topol\'ogicamente por $\tilde\sigma$. En
particular $\Gal(\nr L/E)$ es proc\'iclico, esto es, generado topol\'ogicamente
por un elemento. Puesto $E\subseteq \nr L$, $\nr E\subseteq \nr L$,
obtenemos el epimorfismo can\'onico $\Gal(\nr L/E)\stackrel{\rho}
{\twoheadrightarrow}\Gal(\nr E/E)$ dado por restricci\'on y se tiene
$\ker\rho=\Gal(\nr L/\nr E)$.

Por ser $\Gal(\nr L/E)$ proc\'iclico, $\Gal(\nr L/E)/\Gal(\nr L/E)^n$ es un
grupo c\'iclico de orden $n$: $\frac{\Gal(\nr L/E)}{\Gal(\nr L/E)^n}=
\frac{\overline{\langle\tilde\sigma\rangle}}{\overline{\langle\tilde\sigma^n
\rangle}}\cong \frac{\langle\tilde\sigma\rangle}{\langle\tilde\sigma^n\rangle}$.
De hecho sea $\frac{\ma Z}{n{\ma Z}}\stackrel{\varepsilon_n}{\lra}\frac{\Gal(\nr L/E)}
{\Gal(\nr L/E)^n}$ dado por $\varepsilon_n(1\bmod n)=\tilde\sigma \bmod
\Gal(\nr L/E)^n$ es un epimorfismo. Tomando la composici\'on inducida
por $\rho$:
\[
\frac{\ma Z}{n{\ma Z}}\stackrel{\varepsilon_n}{\lra}\frac{\Gal(\nr L/E)}
{\Gal(\nr L/E)^n}\stackrel{\tilde\rho}{\lra}\frac{\Gal(\nr E/E)}
{\Gal(\nr E/E)^n}\xrightarrow[]{(\varepsilon'_n)^{-1}}\frac{\ma Z}{n{\ma Z}}
\]
es un isomorfismo y $\rho$ es un isomorfismo (al tomar l\'imites inversos).
Por tanto $\Gal(\nr L/\nr E)=\ker\rho =\{1\}$ y $\nr L=\nr E$.

\item[(4)] Se tiene $f_{E|K}\deg_E(\tilde\sigma)=\deg_K(\tilde\sigma)\igual\limits_{
\substack{\uparrow\\(2)}}f_{E|K}$. Por tanto $\deg_K(\tilde\sigma)=1$ y
$\tilde\sigma=\Fro E$. $\fin$
\end{list}
\end{proof}

Consideremos ahora $L/K$ una extensi\'on finita, $L=(\sep k)^{G_L}$, $K(
\sep k)^{G_K}$, $[L:K]=[G_K:G_L]$. Entonces $\N_{L/K}\colon L\lra K$ dado por
$\N a=\prod_{\sigma\in G_K/G_L} a^{\sigma}$, donde $\{\sigma\}$ recorre
un conjunto de representantes de las clases derechas de $G_K/G_L$ (pues
estamos usando la acci\'on por la derecha: $a^{\sigma}$ en lugar de la
acci\'on por la izquierda: $\sigma a$). Se tiene que si $K\subseteq L
\subseteq M$, $\N_{M/K}=\N_{L/K}\circ \N_{M/L}$ y si $L/K$ es Galois,
$K=L^{\Gal(L/K)}$.

\begin{definicion}\label{NeuD5}
Una {\em valuaci\'on de Hensel\index{valuacion de Hensel!valuaci\'on
de Hensel}\index{Hensel@valuaci\'on de $\sim$}} con respecto a $\deg$,
es un homomorfismo $v\colon\*k\lra \hat{\ma Z}$ tal que
\lasa
\item $v(k)={\mc Z}\supseteq {\ma Z}$ y ${\mc Z}/n{\mc Z}\cong
{\ma Z}/n{\ma Z}$ para toda $n\in {\ma N}$.
\item $v\big(\N_{K/k}\*K)=f_K{\mc Z}$ para todo campo $K$, 
$k\subseteq K\subseteq \sep k$.
\end{list}
\end{definicion}

Se tiene que una valuaci\'on de Hensel, es simplemente una
valuaci\'on usual de $\sep k$. Dada una valuaci\'on de Hensel,
$v\colon\*k\lra \hat{\ma Z}$, obtenemos el homomorfismo
$v_K=\frac 1{f_K}\big(v\circ \N_{K/k}\big)\colon\*K\lra \hat{\ma Z}$ con
$\im v_K={\mc Z}$.

\begin{proposicion}\label{NeuP6}
\las
\item $v_K=v_{K^{\sigma}}\circ \sigma$ para $\sigma\in G=\Gal(\sep k/
k)=G_k$.
\item\label{inciso2} Para cada extensi\'on finita de $L/K$, se tiene el siguiente
diagrama conmutativo $\xymatrix{\*L\ar@{-}[r]^{v_L}\ar@{-}[d]_{\N_{L/K}}
&\hat{\ma Z}\ar@{-}[d]^{f_{L|K}}\\ \*K\ar@{-}[r]_{v_K}&\hat{\ma Z}}$
\end{list}
\end{proposicion}

\begin{proof} \las
\item Se tiene que si $\{\tau\}$ es un conjunto completo de representantes 
de las clases derechas de $G_k/G_K$, entonces $\{\sigma^{-1}\tau\sigma\}$
es un conjunto completo de representatnes  de las clases derechas de
$G_k/\sigma^{-1}G_K\sigma=G_k/G_{K^{\sigma}}$. Por tanto, si $a\in
\*K$, se tiene
\begin{align*}
v_{K^{\sigma}}(a^{\sigma})&=\frac 1{f_{K^{\sigma}}}v\big(\prod_{\tau}
a^{\sigma\sigma^{-1}\tau\sigma}\big)=\frac 1{f_K}v\big(\big(
\underbrace{\prod_{\tau}a^{\tau}}_{\in k}\big)^{\sigma}\big)\\
&=\frac 1{f_K}v\big(\N_{K/k}(a)\big)=v_K(a).
\end{align*}

\item Si $a\in\*L$, entonces
\begin{align*}
f_{L|K}v_L(a)&=f_{L|K}\cdot \frac{1}{f_L}v\big(\N_{L/k}(a)\big)=
\frac 1{f_K}v\big(\N_{K/k}\big(\N_{L/K}(a)\big)\big)\\
&=v_K\big(\N_{L/K}(a)\big).
\tag*{$\fin$}
\end{align*}
\end{list}
\end{proof}

\begin{observacion}\label{NeuO7}
Si $f_{L|K}=[L:K]$, es decir, $L_0=L=\nr K\cap L$, esto es, $L\subseteq
\nr K$, entonces por (\ref{inciso2}), $v_L|_{\*K}=v_K$ pues si $a\in\*K$,
$f_{L|K}v_L(a)=[L:K]v_L(a)=v_K\big(\N_{L/K}(a)\big)=v_K\big(a^{[L:K]}
\big)=[L:K]v_K(a)$.

En particular, un elemento primo de $\*K$ es un elemento primo de
$\*L$ (esto ya nos era conocido y, en este contexto, un elemento primo
$\pi_K$ de $K$ es un elemento con $v_K(\pi_K)=1$). Por otro lado,
si $f_{L|K}=1$, esto es, $L_0=K$, si $\pi_L$ es un elemnto primo de
$A_L$, entonces $\pi_K=\N_{L/K}(\pi_L)$ es un elemento primo de 
$A_K$.
\end{observacion}

\subsection{El isomorfismo de Neukirch es el mapeo de Nakayama}\label{SN2}

\begin{teorema}\label{NeuT8}
Sea $L/K$ una extensi\'on finita no ramificada. Entonces 
\begin{gather*}
\co 0{\Gal(L/K)}{\*L} \cong {\ma Z}/[L\colon K]{\ma Z},\\
|\co 0{\Gal(L/K)}{\*L}|=[L\colon K]\quad  \text{y}\\
\co {{-1}}{\Gal(L/K)}{\*L}=\{1\}.
\end{gather*}
\end{teorema}

\begin{proof}
Se tiene que $\co {{-1}}{\Gal(L/K)}{\*L}=\{1\}$ por el Teorema 90 de Hilbert.
Por otro lado, la sucesi\'on $1\lra U_L\lra \*L\stackrel{v_L}{\lra}{\ma Z}\lra 0$
es una sucesi\'on exacta y $\co n{\Gal(L/K)}{U_L}=\{1\}$ para toda $n\in
{\ma Z}$, por tanto $\co n{\Gal(L/K)}{\*L}\cong\co n{\Gal(L/K)}{\ma Z}$ para
toda $n\in{\ma Z}$. En particular, $\co 0{\Gal(L/K)}{\*L}=\co 0{\Gal(L/K)}
{\ma Z}={\ma Z}/[L:K]{\ma Z}$. 
$\fin$
\end{proof}

\begin{definicion}[Mapeo de Neukirch\index{mapeo de 
Neukirch}\index{Neukirch!mapeo de $\sim$}]\label{NeuD9}
Sea $L/K$ una extensi\'o finita de Galois. El {\em mapeo de 
Neukirch\index{mapeo de Neukirch}\index{Neukirch!mapeo de $\sim$}}
\begin{gather*}
{\mc N}_{L/K}\colon \Gal(L/K)\lra \*K/\N_{L/K}\*L
\intertext{est\'a dado por}
{\mc N}_{L/K}(\sigma):=\N_{E/K}(\pi_E)\bmod \N_{L/K} \*L,
\end{gather*}
donde $E$ es el campo fijo de un levantamiento de Frobenius $\tilde
\sigma\in\Lambda(\nr L/K)$ de $\sigma\in\Gal(L/K)$ y $\pi_E$ es un
elemento primo de $\*E$, $v_E(\pi_E)=1$.
\end{definicion}

\begin{observacion}\label{NeuO10}
Se tiene que probar que la definici\'on de ${\mc N}_{L/K}(\sigma)$ es
independiente de la selecci\'on del levantamiento $\tilde\sigma$ y
del elemento $\pi_E$.
\end{observacion}

Para probar lo anterior, sea $L_0=\nr K\cap L$. Sea $\N:=\N_{L_0/L}$.
\[
\xymatrix{
&\nr K\ar@{-}[r]\ar@{-}[d]&\nr L=\nr K L\ar@{-}[d]\\
&\nr K\cap L=L_0\ar@{-}[r]\ar@{-}[dl]&L\ar@{-}[dl]\\
K\ar@{-}[r]&F
}
\]

Sea $K\subseteq F\subseteq L$ tal que $FL_0=L$, por ejemplo, $F=L$.
Entonces $F_0=\nr K\cap F=\nr K\cap L\cap F=L_0\cap F$ y 
$\N|_F=\N_{F/F_0}$. Fijemos una extensi\'on $\varphi\in\Lambda(\tilde L/K)$
de $\Fro K$ a $\nr L$ y se tiene $\N_{F/K}=\N\circ {\eu N}_F$, donde el
homomorfismo ${\eu N}_F\colon F\lra L$ est\'a definido por
\[
{\eu N}_F(a)=\prod_{\eta=0}^{f-1}a^{\varphi^{\eta}},\quad f=[F_0:K].
\]

De hecho, $\Gal(F_0/K)$ es un grupo c\'iclico de orden $f$ y generado
por $\Fr{F_0}{K}=\Fro K|_{F_0}$, por lo que si $a\in F$, se tiene
\[
\N\big({\eu N}_F(a)\big)=\prod_{\eta=0}^{f-1}\N_{F/F_0}(a)^{\varphi^{
\eta}}=\N_{F_0/K}\big(\N_{F/F_0}(a)\big)=\N_{F/K}(a).
\]

La independencia buscada es consecuencia del siguiente resultado.

\begin{lema}\label{NeuL11}
Sean $\tilde \sigma_1,\tilde\sigma_2,\tilde\sigma_3\in\Lambda(\nr L/K)$
tales que $\tilde\sigma_3=\tilde\sigma_1\tilde\sigma_2$. Si $E_i:=
(\nr L)^{\langle \tilde\sigma_i\rangle}$ es el campo fijo de $\tilde\sigma_i$,
$1\leq i\leq 3$, y si $\pi_i\in E_i$ es un elemento primo de $E_i$, entonces
\[
\N_{E_3/K}(\pi_3)\equiv \N_{E_1/K}(\pi_1)\N_{E_2/K}(\pi_2)
\bmod \N_{L/K}\*L.
\]
\end{lema}

\begin{proof}
Sea $F$ el campo fijo de $\varphi$, el levantamiento de $\Fro K$ a $\nr L$.
\[
\xymatrix{
\nr K\ar@{-}[rr]\ar@{-}[d]_{\Fro K}&&\nr L\ar@{-}[dl]^{\varphi}\ar@{-}[dll]\\
K\ar@{-}[r]&F}
\qquad 
\xymatrix{
\nr K\ar@{-}[d]\ar@{-}[r]\ar@{-}@/_4pc/[ddr]_{\Fro K}&\nr L\ar@{-}[d]^{\varphi}\\
\underbracket[0pt]{F_0}_{\substack{\uigual\\ \nr K\cap F}}
\ar@{-}[r]\ar@{-}[rd]&F\\ &K
}
\qquad
\xymatrix{
\nr K\ar@{-}[r]\ar@{-}[d]&\nr L\ar@{-}[d]_{\varphi}\\ F_0=K\ar@{-}[r]& F
}
\]

Puesto que $\varphi|_{\nr K}=\Fro K$, $F_0=\nr K\cap F=K$. Podemos suponer
que $F$ y que $E_i$, $1\leq i\leq 3$, est\'an contenidos en $L$, pues si no
fuese as\'i, podemos considerar una extensi\'on finita de Galois $L'/K$ contenida
en $\nr L$ y conteniendo tanto a $L$ como a $F,E_1,E_2,E_3$ (por ejemplo, la
cerradura de Galois $L'$ de $F':=LFE_1E_2E_3/K$). Debido a que 
$L\subseteq L'\subseteq \nr L$, se tiene que $\nr {(L')}=\nr L$ y
$\tilde\sigma_1,\tilde\sigma_2,\tilde\sigma_3\in\Lambda(\nr{(L')}/K)=\Lambda
(\nr L/K)$ y si probamos que 
\begin{gather*}
\N_{E_3/K}(\pi_3)\equiv \N_{E_1/K}(\pi_1)\N_{E_2/K}(\pi_2)\bmod \N_{L'/K}\*{(L')},
\intertext{entonces}
N_{E_3/K}(\pi_3)\equiv \N_{E_1/K}(\pi_1)\N_{E_2/K}(\pi_2)\bmod \N_{L/K}\*L,
\end{gather*}
pues $\N_{L'/K}\*{(L')}\subseteq \N_{L/K}(\*L)$.

En resumen, podemos suponer que $F,E_i\subseteq L$, $1\leq i\leq 3$. Sea $n=
[L:K]$ y sea $M/L$ la subextensi\'on de $\nr L/L$ de grado $n$, esto es, $[M:L]=n$.
\[
\xymatrix{
&&&M\ar@{-}[r]\ar@{=}[d]\ar@{-}[dl]_n&\nr L\\
E\ar@{-}[rr]\ar@{-}[dd]&&L\ar@{-}[r]\ar@{-}[dd]&M_0L\ar@{-}[dd]\\
&E_i\ar@{-}[ru]\ar@{-}[d]&&&&M_0=M\cap \nr K\\
K\ar@{-}@/^1pc/[uurr]^n\ar@{-}[r]&(E_i)_0\ar@{-}[r]&L_0\ar@{-}[r]&M_0
}
\]

Ahora bien, $M_0\cap L=M\cap \nr K\cap L\igual\limits_{\substack{\uparrow\\
L\subseteq M}}L\cap \nr K=L_0$ y $[M:M_0]=e(M|K)=e(M|L)e(L|K)=e(L|K)=
[L:L_0]$, por tanto $[LM_0:M_0]=[L:L_0]=[M:M_0]$, por lo que $M=LM_0$.
Adem\'as $[\nr L:\nr K]=[M:M_0]$.
\[
\xymatrix{
L\ar@{-}[r]\ar@{-}[d]&LM_0\ar@{-}[d]\\ L\cap M_0=L_0\ar@{-}[r]&M_0}\qquad
\xymatrix{
\nr K\ar@{-}[r]\ar@{-}[d]&\nr L=M\nr K\ar@{-}[d]\\ M_0=M\cap \nr K\ar@{-}[r]&M}
\]

Por tanto $\N=\N_{L/L_0}$ se extiende a $\N_{M/M_0}$. Sean $n_i=\deg_K(
\tilde\sigma_i)$. Puesto que $\tilde\sigma_3=\tilde\sigma_1\tilde\sigma_2$,
$n_3=n_1+n_2$. Sea $\tilde\sigma_4=\varphi^{-n_2}\tilde\sigma_1\varphi^{n_2}$
el conjugado de $\tilde\sigma_1$ (veremos en un momento para que se hace 
esta definici\'on), $\deg_K(\tilde\sigma_4)=n_4=n_1=\deg_K(\tilde\sigma_1)$.
El campo fijo $L^{\tilde\sigma_4}=E_4$ de $\tilde\sigma_4$ es el conjugado
$E_4=E_1^{\varphi^{n_2}}$ y $\pi_4:=\pi_1^{\varphi^{n_2}}$ es un elemento
primo de $E_4$. Se tiene que $\N_{E_1/K}(\pi_1)=\N_{E_4/K}(\pi_4)$ y por
tanto debemos probar la congruencia
\[
N_{E_3/K}(\pi_3)\equiv \N_{E_4/K}(\pi_4)\N_{E_2/K}(\pi_2)\bmod \N_{L/K}\*L.
\]

Sea $\tau_i:=\tilde\sigma_i^{-1}\varphi^{n_i}\in \Gal(\nr L/\nr K)\cong
\Gal(M/M_0)$, se tiene
\begin{align*}
\tau_3&=\tilde\sigma_3^{-1}\varphi^{n_3}=\big(\varphi^{-n_3}\tilde\sigma_3\big)^{-1}
\igual_{\substack{\uparrow\\ n_3=n_1+n_2}}\big(\varphi^{-n_1-n_2}\tilde\sigma_1
\tilde\sigma_2\big)^{-1}\\
&=\tilde\sigma_2^{-1}\tilde\sigma_1^{-1}\varphi^{n_1+n_2}
=\tilde\sigma_2^{-1}\varphi^{n_2}\underbrace{\varphi^{-n_2}\tilde\sigma_1^{-1}
\varphi^{n_2}}_{\tilde\sigma_4^{-1}}\varphi^{n_1}\igual_{\substack{\uparrow\\ 
n_1=n_4}}\tau_2\tilde\sigma_4^{-1}\varphi^{n_4}\\
&=\tau_2\tau_4,
\intertext{por tanto, $\tau_3=\tau_2\tau_4$ y}
{\eu N}_{E_i}(\pi_i)^{\varphi-1}&=\Big(\prod_{\gamma=0}^{[(E_i)_0:K]-1}
\pi_i^{\varphi^{\gamma}}\Big)^{
\varphi-1}=\pi_i^{(\varphi-1)(1+\varphi+\cdots+\varphi^{[(E_i)_0:K]-1})}\\
&=\pi_i^{\big(\varphi^{[(E_i)_0:K]}-1\big)}=\pi^{\varphi^{n_i}-1}\igual_{
\substack{\uparrow\\ 
\pi_i\in E_i\\ =(\nr L)^{\langle\tilde\sigma_i\rangle}}}\pi_i^{\tilde\sigma_i^{-1}
\varphi^{n_i}-1}=\pi_i^{\tau_i-1},
\end{align*}
es decir ${\eu N}(\pi_i)^{\varphi-1}=\pi_i^{\tau_i-1}$.

Sea $A:={\eu N}_{E_3}(\pi_3){\eu N}_{E_2}(\pi_2)^{-1}{\eu N}_{E_4}
(\pi_4)^{-1}\in U_L$, unidad de $L$ pues
\[
v_L\big({\eu N}_{E_i}(\pi_i)\big)=v_L\big(\prod_{\gamma=0}^{n_i-1}
\pi_i^{\varphi^{\gamma}}\big)=\sum_{\gamma=0}^{n_i-1}1=n_i
\]
y $v_L(A)=n_3-n_2-n_4=n_3-n_2-n_1=0$, esto es, $A\in U_L$.

Adem\'as $A^{\varphi-1}=\pi_3^{\tau_3-1}\pi_2^{1-\tau_2}\pi_4^{1-\tau_4}$.
Ahora bien, por la Proposici\'on \ref{NeuP3} se tiene $\nr L=\nr E_i$ y
$E_i\subseteq L\subseteq \nr L=\nr E_i$, esto es, $L/E_i$ es no
ramificada, en particular $\pi_i$ es un elemento primo de $L$. Se sigue
que existen $\varepsilon_2,\varepsilon_3\in U_L$ con $\pi_2=
\varepsilon_2^{-1}\pi_4$ y $\pi_3=\varepsilon_3\pi_4$. Sea $\varepsilon_4
=\pi_4^{\tau_2-1}\in U_L$. Se obtiene
\begin{gather*}
(\tau_3-1)+(1-\tau_2)+(1-\tau_4)\igual_{\substack{\uparrow\\ \tau_3=\tau_2
\tau_4}}(\tau_2-1)(\tau_4-1)
\intertext{y por tanto}
\begin{align*}
A^{\varphi-1}&=\pi_3^{\tau_3-1}\pi_2^{1-\tau_2}\pi_4^{1-\tau_4}=
\varepsilon_3^{\tau_3-1}\pi_4^{\tau_3-1}\varepsilon_2^{\tau_2-1}
\pi_4^{1-\tau_2}\pi_4^{1-\tau_4}\\
&= \varepsilon_3^{\tau_3-1}\varepsilon_2^{\tau_2-1}\pi_4^{(\tau_3
-1)+(1-\tau_2)+(1-\tau_4)}\\
&=\varepsilon_3^{\tau_3-1}\varepsilon_2^{\tau_2-1}\pi_4^{(\tau_2-1)
(\tau_4-1)}=\varepsilon_3^{\tau_3-1}\varepsilon_2^{\tau_2-1}
\varepsilon_4^{\tau_4-1}.
\end{align*}
\end{gather*}
Es decir, $A^{\varphi-1}=\prod_{i=2}^4 \varepsilon_i^{\tau_i-1}$.

Ahora bien, puesto que $F\subseteq L$, se sigue de la
Proposici\'on \ref{NeuP3} que $\nr F=\nr L$. Por tanto $F\subseteq
L\subseteq \nr L= \nr F$ y $L/F$ es no ramificada. Por tanto
$f=[L:F]=f_{L|F}$ y $\varphi^f=\Fro L$ (recordemos que $\varphi
\in\Lambda(\nr L/K)$ es una extensi\'on de $\Fro K$ a $\nr L$ y
$\xymatrix{\nr K\ar@{-}[rr]\ar@{-}[d]&&\nr L\ar@{-}[d]^{\Fro L}
\ar@{-}[dl]_{\varphi}\\K\ar@{-}[r]&F\ar@{-}[r]_f&L}$,
$\Gal(\nr L/F)=\overline{\langle \varphi\rangle}$, $\Gal(\nr L/L)=
\overline{\langle\Fro L\rangle}$).

Ahora $M/L$ es no ramificada por construcci\'on y $\co i{\Gal(M/L)}
{U_M}=0$ para toda $i\in{\ma Z}$ y como $A,\varepsilon_i\in U_L$ y
$\N_{M/L}U_M= U_L$ pues se tiene $\co 0{\Gal(M/L)}{U_M}=0$, por lo que
existen $\tilde A_i,\tilde\varepsilon_i\in U_M$ con $\N_{M/L}\tilde A=
A$ y $\N_{M/L}\tilde\varepsilon_i=\varepsilon_i$. Por tanto 
$\N_{M/L}\big(\tilde A^{\varphi-1}\big)=\N_{M/L}\big(\prod_{i=2}^4
\tilde\varepsilon_i^{\tau_i-1}\big)$, por lo que existe $x\in U_M$,
$\N_{M/L} x=1$ y $\tilde A^{\varphi-1}=x\prod_{i=2}^4\tilde\varepsilon_i^{
\tau_i-1}$. Por otro lado $\co {{-1}}{\Gal(M/K)}{U_M}=1$ por lo que
$\ker\N_{M/L}=I_{\Gal(M/K)}U_M$, es decir, $x=\tilde u^{\Fro L-1}$
con $\tilde u\in U_M$ y $\Fro L-1=\varphi^f-1=(\varphi-1)(\varphi^{f-1}+
\cdots+\varphi+1)$. Por tanto
\[
\tilde A^{\varphi-1}=\tilde u^{\Fro L-1}\cdot\prod_{i=2}^4\tilde\varepsilon_i^{
\tau_i-1}=\Big(\prod_{\gamma=0}^{f-1}\tilde u^{\varphi^{\gamma}}\Big)^{
\varphi-1}\cdot \prod_{i=2}^{4}\tilde\varepsilon_i^{\tau_i-1}.
\]

Ahora, $\tau_i\in\Gal(M/M_0)$ ($\tau_i=\tilde\sigma_i\varphi^{n_i}\in
\Gal(M/M_0)$). Por tanto, considerando que $\N=\N_{L/L_0}$ se
extiende a $\N_{M/M_0}$, $M_0\subseteq \nr K$ y por tanto
$\varphi-1|_{M_0}=\Fro K-1$, se sigue que
\begin{align*}
\N(\tilde A^{\varphi-1})&=\N(\tilde A)^{\varphi-1}
=\N(\tilde A)^{\Fro K-1}\igual_{\substack{\uparrow\\
\N(\tilde A)\in M_0}}\\
&=\N\Big(\prod_{\gamma=0}^{f-1}\tilde u^{\varphi^{\gamma}}
\Big)^{\Fro K-1}\cdot \underbrace{\prod_{i=2}^4 \overbrace{
\N(\tilde \varepsilon_i^{\tau_i})}^{\substack{\N(\tilde\varepsilon_i)\\ \uigual}}
\N(\tilde\varepsilon_i)^{-1}}_{\substack{\uigual\\ 1}}
=\N\Big(\prod_{\gamma=0}^{f-1}\tilde u^{\varphi^{\gamma}}
\Big)^{\Fro K-1},
\end{align*}
es decir, $\N(\tilde A^{\varphi-1})=\N\Big(\prod_{\gamma=0}^{f-1}\tilde 
u^{\varphi^{\gamma}}\Big)^{\Fro K-1}$.

Se sigue que $\N(\tilde A)=\N\Big(\prod_{\gamma=0}^{f-1}\tilde u^{\varphi^{\gamma}}
\Big)\cdot y$ con $y\in U_{M_0}$ y $y^{\Fro K-1}=1$, esto es, $y\in U_K$.

Sea $u:=\N_{M/L}(\tilde u)$ y recordando que 
\[
\begin{minipage}{5cm}
\xymatrix{
L\ar@{<-}[rr]^{\N_{M/L}}\ar@{->}[d]_{\N}&&M\ar@{->}[d]^{\N}
\ar@{->}[dll]\\ L_0\ar@{<-}[rr]_{\N_{M_0/L_0}}&&M_0}
\end{minipage}
\qquad\qquad
\begin{minipage}{5cm}
$\N_{M_0/L_0}\circ \N=
\N\circ \N_{M/L}$ y que
$\N_{F/K}=\N\circ {\eu N}_F$, se obtiene que
\end{minipage}
\]
\begin{align*}
\N(A)&=\N\big({\eu N}_{E_3}(\pi_3){\eu N}_{E_2}(\pi_2)^{-1}
{\eu N}_{E_4}(\pi_4)^{-1}\big)\\
&=\N_{E_3/K}(\pi_3)\N_{E_2/K}(\pi_2)^{-1}\N_{E_4/K}(\pi_4)^{-1}
\igual_{\substack{\uparrow\\ \N=\N_{L/L_0}}}\N\big(\N_{M/L}(
\tilde A)\big)\\
&\igual_{\substack{\uparrow\\
\N=\N_{M/M_0}}}\N_{M/M_0}\big(
\N_{M/L}(\tilde A)\big)=\N_{M_0/L_0}(\N(\tilde A))\\
&=\N_{M_0/L_0}\Big(\N\big(\prod_{\gamma=0}^{f-1}\tilde u^{\varphi^{\gamma}}
\big)\cdot y\Big)\igual_{\substack{\uparrow\\ \N_{M_0/L_0}\circ \N\\
=\N\circ \N_{M/L}}}\N\Big(\N_{M/L}\big(\prod_{\gamma=0}^{f-1}\tilde u^{
\varphi^{\gamma}}\big)\Big)\underbrace{\N_{M_0/L_0} y}_{
\substack{\uigual\\ y^{[M_0:L_0]}}}\\
&\igual_{\substack{\uparrow\\ n=[M_0:L_0]}}\N\Big(\prod_{\gamma=0}^{
f-1}\N_{M/L}\big(\tilde\mu\big)^{\varphi^{\gamma}}\Big)\cdot y^n\igual_{
\substack{\uparrow\\ u=\N_{L/K}(\tilde u)}}\prod_{\gamma=0}^{f-1}
\N(u)^{\varphi^{\gamma}}\cdot y^n\\
&=\prod_{\gamma=0}^{f-1}\big(\N\circ \N_{M/L}\big)(\tilde u)^{\varphi^{
\gamma}}y^n=\N\Big(\prod_{\gamma=0}^{f-1}u^{\varphi^{\gamma}}\Big)\cdot y^n
=\N\circ {\eu N}_L(u)\cdot y^n\\
&=\N_{L/K}(u)\N_{L/K} (y)\in \N_{L/K}\*L,
\end{align*}
lo cual prueba el resultado. $\fin$
\end{proof}

\begin{corolario}\label{NeuC12}
El mapeo de Neukirch, 
\[
{\mc N}_{L|K}\colon \Gal(L/K)\lra \*K/\N_{L/K}\*L,\quad
{\mc N}_{L|K}(\sigma):=\N_{F/K}(\pi_F)\bmod \N_{L/K}\*L,
\]
est\'a bien definido y es un homomorfismo.
\end{corolario}

\begin{proof}
Sean $\tilde \sigma$, $\tilde\sigma'$ dos levantamientos de Frobenius de
$\sigma$. Sean $F$ y $F'$ los campos fijos de $\tilde\sigma$ y $\tilde\sigma'$
respectivamente y sean $\pi\in\*F$, $\pi'\in \*{(F')}$ elementos primos. Ordenamos
$\tilde\sigma$, $\tilde\sigma'$ tales que $m=\deg_K(\tilde\sigma')-\deg_K(
\tilde\sigma)\geq 0$. Si $m=0$, $\tilde \sigma'|_{\nr K}=\tilde \sigma|_{\nr K}$
y $\tilde\sigma'|_L=\tilde\sigma|_L$ de donde obtenemos que $\tilde\sigma'=
\tilde\sigma$ y $\pi'=\pi u$ con $u\in U_F$. Sean $M$ una extensi\'on finita de
Galois $M/K$, $M\subseteq \nr L$ y tal que $L, F\subseteq M$. Puesto que 
$\co 0{\Gal(M/F)}{U_M}=\{0\}$, existe $\tilde u\in U_M$ con $\N_{M/F}(\tilde u)
=u$. Por tanto
\[
\N_{F'/K}(\pi')=\N_{F/K}(\pi)\N_{M/K}(\tilde u)\equiv \N_{F/K}(\pi)\bmod 
\N_{L/K}\*L
\]
puesto que $\N_{M/K}(\tilde u)\in \N_{M/K}\*M\subseteq \N_{L/K}\*L$.

Supongamos ahora que $m>0$. Entonces $\tilde\tau=\tilde\sigma^{-1}
\tilde\sigma\in\Gal(\nr L/K)$ es un levantamiento de Frobenius de $1\in
\Gal(L/K)$ con $\deg_K(\tilde\tau)=\deg_K(\tilde\sigma')-\deg_K(\tilde
\sigma)=m>0$. Sea $M$ el campo fijo de $\tilde\tau$ ($\tilde\tau|_L=1$,
por lo que $L\subseteq M$) y sea $\pi_M\in\*M$ un elemento primo.
Se sigue del Lema \ref{NeuL11}
\[
\N_{F'/K}(\pi')\equivalente_{\substack{\uparrow\\ \tilde\sigma'=
\tilde\sigma \tilde\tau}}\N_{F/K}(\pi)\N_{M/K}(\pi_M)\equiv \N_{F/K}(\pi)
\bmod \N_{L/K}\*L.
\]

Esto prueba la independencia de ${\mc N}_{L/K}(\pi)(\sigma)$ de la
selecci\'on del levantamiento $\tilde\sigma\in\Lambda(\nr L/K)$ y del
elemento primo $\pi_F\in\*F$. Que ${\mc N}_{L/K}$ es un homomorfismo
se sigue del Lema \ref{NeuL11} puesto que si $\tilde\sigma_1,
\tilde\sigma_2$ son levantamientos de Frobenius de $\sigma_1,
\sigma_2\in\Gal(L/K)$, entonces $\tilde\sigma_3:=\tilde\sigma_1
\tilde\sigma_2$ es un levantamiento de Frobenius de $\sigma_3=
\sigma_1\sigma_2$. $\fin$
\end{proof}

El Teorema de Reciprocidad, el cual asocia a automorfismos de
Frobenius elementos primos, es el siguiente resultado.

\begin{teorema}\label{NeuT13}
Si $L/K$ es una subsextensi\'on finita de $\nr K/K$, entonces el
inverso del mapeo de reciprocidad
\[
{\mc N}_{L|K}\colon\Gal(L/K)\lra \*K/\N_{L/K}\*L
\]
est\'a dado por ${\mc N}_{L|K}(\Fr LK)=\pi_K\bmod \N_{L/K}\*L$
y es un isomorfismo.
\end{teorema}

\begin{proof}
Se tiene que $\Fro K\in\Lambda(\nr K/K)=\Lambda(\nr L/K)$ es un
levantamiento de Frobenius de $\Fr LK\in\Gal(L/K)$ y el campo
fijo de $\Fro K$ es $K$. Por tanto ${\mc N}_{L|K}(\Fr LK)
\igual\limits_{\substack{\uparrow\\ F=K}}\pi_K\bmod \N_{L/K}\*L$.

Ahora bien, $\co 0{\Gal(L/K)}{\*L}=\frac{\*K}{\N_{L/K}\*L}$ tiene
orden $n=[L\colon K]$ el cual tambi\'en es el orden de $\Gal(L/K)$ y 
$\pi_K\bmod \N_{L/K}\*L$ es un generador de $\*K/\N_{L/K}\*L$
pues $\pi_K^m=\N_{L/K}(a)\in\N_{L/K}\*L$ implica que $m=v_K\big(
\N_{L/K}(a)\big)=nv_L(a)\equiv 0\bmod n$. Por tanto ${\mc N}_{L|K}$
es un isomorfismo.
$\fin$
\end{proof}

\begin{observacion}\label{NeuO14}
M\'as adelante veremos que ${\mc N}_{L|K}$ es un isomorfismo para
cualquier extensi\'on abeliana finita $L/K$, no necesariamente $L
\subseteq \nr K$, para lo cual debemos pedir el axioma de campos de
clase, es decir lo del Teorema \ref{NeuT8} para cualquier extensi\'on
$L/K$ es c\'iclica.
\end{observacion}

Una de las propiedades fundamentales de ${\mc N}_{L|K}$ es que
es funtorial.

\begin{teorema}\label{NeuT15}
Sean $L/K$ y $L'/K'$ dos extensiones finitas de Galois de campos
locales tales que $K\subseteq K'$ y $L\subseteq L'$. Entonces
el diagrama
\[
\xymatrix{
\Gal(L'/K')\ar@{->}[rr]^{{\mc N}_{L'|K'}}\ar@{->}[d]_{\rest|_L}&&
\*{(K')}/\N_{L'/K'}\*{(L')}\ar@{->}[d]^{\N_{K'/K}}\\
\Gal(L/K)\ar@{->}[rr]^{{\mc N}_{L|K}}&&\*K/\N_{L/K}\*L}
\]
es conmutativo.
\end{teorema}

\begin{proof}
\[
\xymatrix{
&&L'\ar@{-}[d]\\ &L\ar@{-}[r]\ar@{-}[d]&LK'\ar@{-}[d]\\ &L\cap K'\ar@{-}[r]
\ar@{-}[ld]&K'\\ K
}
\]
Sean $\sigma'\in\Gal(L'/K')$ y $\sigma=\sigma'|_L$. Si $\tilde \sigma'\in
\Lambda(\nr{(L')}/K')$ es un levantamiento de Frobenius de $\sigma'$,
entonces $\tilde \sigma=\tilde\sigma'|_L\in \Gal(\nr L/K)$ es un levantamiento
de Frobenius de $\sigma$ puesto que $\deg_K(\tilde\sigma)=f_{K'|K}\deg_K
(\tilde\sigma')\in{\ma N}$.

Sea $F'$ el campo fijo de $\tilde\sigma'$. Entonces, $F=F'\cap \nr L=F'\cap
\nr F$ es el campo fijo de $\tilde \sigma$ y puesto que $F'/F$ es totalmente
ramificada, $f_{K'|K}=1$. Por tanto, si $\pi_{F'}$ es un elemento primo de
$F'$ y $\pi_F:=\N_{F'/F}(\pi_{F'})$ es un elemento primo de $F$. El teorema
se sigue de
\begin{align*}
{\mc N}_{L|K}(\rest|_L\sigma')&={\mc N}_{L|K}\big(\sigma'|_L\big)=\N_{F/K}(\pi_F)
=\N_{F/K}\big(\N_{F'/F}(\pi'_F)\big)\\
&=\N_{F'/K}(\pi'_F)=\N_{K'/K}\big(\N_{F'/K'}(\pi'_F)\big)\\
&=\N_{K'/K}\big({\mc N}_{F'|K'}(\sigma')\big). \tag*{$\fin$}
\end{align*}
\end{proof}

\begin{proposicion}\label{NeuP16}
Si $L/K$ es una extensi\'on finita de Galois y $\sigma\in G$, entonces el
diagrama
\[
\xymatrix{
\Gal(L/K)\ar@{->}[rr]^{{\mc N}_{L|K}}\ar@{->}[d]_{\*\sigma}&&
\*K/\N_{L/K}\*L\ar@{->}[d]^{\sigma}\\
\Gal(L^{\sigma}/K^{\sigma})\ar@{->}[rr]^(.40){{\mc N}_{L^{\sigma}
|K^{\sigma}}}&&\*{(K^{\sigma})}/\N_{L^{\sigma}/K^{\sigma}}(\*{L^{\sigma})}}
\]
es conmutativa y la flecha izquierda est\'a inducida por conjugaci\'on
$\tau\longmapsto \*{\sigma}(\tau)=\sigma^{-1}\tau\sigma$.
\end{proposicion}

\begin{proof}
Sea $\tilde\tau\in\Lambda(\nr L/K)$ un levantamiento de Frobenius de
$\tau\in\Gal(L/K)$. Sea $\hat\tau\in\Gal(\sep k/K)$ una extensi\'on de
$\tilde\tau$ a $\sep k$. Entonces $\sigma^{-1}\hat\tau \sigma$ es
un levantamiento de Frobenius de $\*\sigma(\tau)$ puesto que $\deg_{
K^{\sigma}}(\sigma^{-1}\hat\tau\sigma)=\deg_K(\hat\tau)=\deg_K(
\tilde\tau)\in{\ma N}$.

Sea $F$ el campo fijo de $\tilde\tau$. Entonces $F^{\sigma}$ es el campo
fijo de $\sigma^{-1}\hat\tau\sigma|_{\nr{(L^{\sigma})}}$ y si $\pi_F$ es un
elemento primo de $F$, entonces $\pi_F^{\sigma}$ es un elemento primo
de $F^{\sigma}$ pues debido a la Proposici\'on \ref{NeuP6} se tiene que
$v_K=v_{K^{\sigma}}\circ \sigma$. Finalmente
\[
{\mc N}_{F^{\sigma}|K^{\sigma}}(\*\sigma(\tau))=\N_{F^{\sigma}/K^{\sigma}}
(\pi^{\sigma})=\N_{F/K}(\pi)^{\sigma}=\sigma\circ {\mc N}_{L|K}(\tau),
\]
es decir, ${\mc N}_{F^{\sigma}|K^{\sigma}}\circ \*\sigma=\sigma\circ
{\mc N}_{L|K}$. $\fin$
\end{proof}

Otra propiedad funtorial de ${\mc N}_{L|K}$ est\'a relacionada con el
mapeo de transferencia $\Ver$ que es $\res_{-2}$. Recordemos como
se define el mapeo de transferencia $\Ver\colon G/G'\lra H/H'$ donde
$H<G$ es un subgrupo de \'indice finito. Sea $R$ un conjunto de
representantes de las clases izquierdas de $H$ en $G$, $G=R\cdot H$,
$1\in R$. Si $\sigma\in G$, para cada $\rho\in R$, se tiene $\sigma\rho=
\rho'\sigma_{\rho}$ con $\sigma_{\rho}\in H$, $\rho'\in R$. Entonces
$\Ver(\sigma\bmod G')=\prod_{\rho\in R}\sigma_{\rho}\bmod H'$.

Otra descripci\'on de $\Ver$, est\'a dada de la siguiente forma.
Sea $\sigma\in G$ y sea $S:=\langle\sigma\rangle$ el subgrupo
generado por $\sigma$. Sea $\tau$ recorriendo un sistema de 
representantes de las clases dobles $SxH$, esto es, $G=\cupdot_{\tau}
S\tau H$. Sea $S_{\tau}:=\tau^{-1}S\tau\cap H$ y sea $n_{\tau}$
el m\'inimo n\'umero natural tal que $\sigma_{\tau}=\tau^{-1}\sigma^{
n_{\tau}}\in H$. Entonces $\sigma_{\tau}$ genera $S_{\tau}$ y
\[
\Ver(\sigma\bmod G')=\prod_{\tau}\sigma_{\tau}\bmod H'.
\]
Esta \'ultima f\'ormula se obtiene de la anterior tomando como $R$
al conjunto $\{\sigma^i\tau\mid i=1,\ldots,n_{\tau}\text{\ para toda 
$n_{\tau}$}\}$.

\begin{proposicion}\label{NeuP17}
Sea $L/K$ una extensi\'on finita de Galois y sea $K\subseteq K'\subseteq L$.
Se tiene un diagrama conmutativo
\[
\xymatrix{
\abe{\Gal(L/K')}\ar@{->}[rr]^{{\mc N}_{L|K'}}\ar@{<-}[d]_{\Ver}&&
\*{(K')}/\N_{L/K'}(\*L)\ar@{<-}[d]^i\\
\abe{\Gal(L/K)}\ar@{->}[rr]^{{\mc N}_{L|K}}&&\*K/\N_{L/K}(\*L)}
\]
donde $i$ es inducido por la contenci\'on $\*K\subseteq\*{(K')}$.
\end{proposicion}

\begin{proof}
Sean ${\mc G}:=\Gal(\nr L/K)$ y ${\mc H}:=\Gal(\nr L/K')\subseteq{\mc G}$.
Sea $\tilde\sigma$ un levantamiento de Frobenius de $\sigma\in\Gal(L/K)$,
$F$ el campo fijo de $\tilde\sigma$ y $S=\Gal(\nr L/F)$. Consideremos la
descomposici\'on de las clases dobles ${\mc G}=\cupdot_{\tau} S\tau {\mc H}$.
Sean $S_{\tau}=\tau^{-1}S\tau\cap {\mc H}$ y $\tilde\sigma_{\tau}=\tau^{-1}
\tilde\sigma^{n_{\tau}}\tau\in{\mc H}$ como antes.

Sean $G:=\Gal(L/K)$, $H:=\Gal(L/K')<G$, $\bar S=\langle\sigma\rangle$, 
$\bar\tau=\tau|_L$ y $\sigma_{\tau}=\tilde\sigma_{\tau}|_L$. Entonces
$G=\cupdot_{\tau} \bar S\bar\tau H$ y por tanto $\Ver(\sigma\bmod G')=
\prod_{\tau}\sigma_{\tau}(\bmod H')$.

Para cada $\tau$, sea $\omega_{\tau}$ un conjunto completo de 
representantes de las clases derechas de ${\mc H}/S_{\tau}$. Entonces
${\mc H}=\cupdot_{\tau}S_{\tau}\omega_{\tau}$ y ${\mc G}=\cupdot_{\tau,
w_{\tau}}S\tau\omega_{\tau}$. Sea $F_{\tau}$ el campo fijo de $\tilde
\sigma_{\tau}$, esto es, el campo fijo de $S_{\tau}$. Se tiene que 
$F^{\tau}=\tau(F)$ es el campo fijo de $\tau^{-1}\tilde\sigma\tau$ de
tal forma que $F_{\tau}/F^{\tau}$ es la subextensi\'on de $\nr L/F^{\tau}$
de grado $n_{\tau}$. Sea $\pi_F$ un elemento primo de $F$. Por
tanto $\pi_{F^{\tau}}:=\pi_F^{\tau}$ es un elemento primo de $F^{\tau}$
y en consecuencia de $F_{\tau}$ por ser $F_{\tau}/F^{\tau}$ no
ramificada. Con la descomposici\'on de clases dobles obtenemos
\[
\N_{F/K}(\pi_F)=\prod_{\tau,\omega_{\tau}}\pi_F^{\tau\omega_{\tau}}=
\prod_{\tau}\big(\prod_{\omega_{\tau}}(\pi_F^{\tau})^{\omega_{\tau}}\big)
=\prod_{\tau}\N_{F_{\tau}/K'}(\pi_F^{\tau}).
\]
Puesto que $\tilde\sigma_{\tau}\in\Lambda(\nr L/K')$ es un 
levantamiento de Frobenius de $\sigma_{\tau}\in H=\Gal(L/K')$ se 
sigue que
\begin{align*}
(i\circ{\mc N}_{L|K})(\sigma)&={\mc N}_{L|K}(\sigma)=\prod_{\tau}
{\mc N}_{L|K'}(\sigma_{\tau})\equiv {\mc N}_{L|K'}\big(\prod_{\tau}
\sigma_{\tau}\big)\\
&\equiv {\mc N}_{L|K'}(\Ver(\sigma \bmod G'))=
({\mc N}_{L|K}\circ \Ver)(\sigma \bmod G').
\tag*{$\fin$}
\end{align*}
\end{proof}

Recordemos que si $L/K$ es una extensi\'on c\'iclica de campos locales.
Entonces $\co 0{\Gal(L/K)}{\*L}\cong {\ma Z}/[L\colon K]{\ma Z}$ y $\co {{-1}}{
\Gal(L/K)}{\*L}=\{1\}$ (Teorema \ref{axiomacamposlocales}, 
Teorema \ref{NeuT8}).

\begin{teorema}[Isomorfismo de Neukirch\index{isomorfismo de
Neukirch}\index{Neukirch!isomorfismo de $\sim$}]\label{NeuT18}
Sea $L/K$ una extensi\'on finita de Galois de campos locales.
Entonces 
\[
{\mc N}_{L|K}\colon \abe{\Gal(L/K)}\lra \*K/\N_{L|K}\*L
\]
es un isomorfismo.
\end{teorema}

\begin{proof}
Si $M/K$ es una subextensi\'on de Galois de $L/K$, entonces por
el Teorema \ref{NeuT15}, se tiene el diagrama conmutativo
\begin{gather}\label{diagramaNeukirch}
\xymatrix{
1\ar@{->}[r]&\Gal(L/M)\ar@{^{(}->}[r]\ar@{->}[d]^{{\mc N}_{L|M}}&
\Gal(L/K)\ar@{->>}[r]\ar@{->}[d]^{{\mc N}_{L|K}}&
\Gal(M/K)\ar@{->}[d]^{{\mc N}_{M|K}}\ar@{->}[r]&1\\
&\*M/\N_{L/M}\*L\ar@{->}[r]_{\N_{M/K}}&\*K/\N_{L/K}\*L\ar@{->}[r]
&\*K/\N_{M/K}\*M\ar@{->}[r]&1
}
\end{gather}

Este diagrama lo usamos para los tres pasos siguientes para la
demostraci\'on del resultado.

\bigskip

\noindent
{\bf Paso 1:} Podemos suponer que $\Gal(L/K)$ es abeliano. Si
probamos el resultado para este caso, y si $M:=\abe L$ es la 
m\'axima subextensi\'on abeliana de $L/K$, entonces $\abe{
\Gal(L/K)}=\Gal(M/K)$ y el subgrupo conmutador $\Gal(L/M)$
de $\Gal(L/K)$ es precisamente el n\'ucleo de ${\mc N}_{L|K}$,
esto es, $\abe{\Gal(L/K)}\lra \*K/\N_{L/K}\*L$ es inyectivo.
\[
\xymatrix{
L\ar@{-}[d]\ar@{-}@/_2pc/[dd]_G\ar@{-}@/^1pc/[d]^{G'=\Gal(L/M)}\\
M=\abe L\ar@{-}[d]\ar@{-}@/^1pc/[d]^{G/G'=\Gal(M/K)}\\ K
}
\]
Para ver la suprayectividad, veamos que esta se sigue para 
extensiones solubles por inducci\'on en el grado del campo.
Esto es, en el caso soluble, ya sea $M=L$ o $[L:M]<[L:K]$ y si
${\mc N}_{M|K}$ y ${\mc N}_{L|M}$ son suprayectivas, 
entonces por el Lema de la Serpiente aplicada al diagrama
(\ref{diagramaNeukirch}), tambi\'en lo es  ${\mc N}_{L|K}$.
Para el caso general, sea $M$ el campo fijo de un 
$p$-subgrupo de Sylow de $\Gal(L/K)$. En general $M/K$ no es
Galois, pero usando la parte izquierda del diagrama 
(\ref{diagramaNeukirch}), en el cual ${\mc N}_{L|K}$ es
suprayectiva. Por tanto es suficiente probar que la imagen de
$\N_{M/K}$ es el $p$-subgrupo de Sylow $S_p$ de $\*K/\N_{
L/K}\*L$, puesto que entonces la imagen de ${\mc N}_{L|K}$
contendr\'ia los $p$-subgrupos de Sylow $S_p$ para toda
$p$ y por tanto ${\mc N}_{L|K}$ ser\'a suprayectiva.

Ahora bien, se tiene que
el encaje $\*K\subseteq \*M$ induce el homomorfismo 
natural $i\colon \*K/\N_{L/K}\*L\lra \*M/\N_{L/M}\*L$ para el 
cual $\N_{M/K}\circ i=[M\colon K]$ ($\N_{M/K}\circ i (x)=\N_{M/K} x
=x^{[M\colon K]}$). Puesto que $\mcd\big([M\colon K],p\big)=1$, $S_p
\xrightarrow[]{[M:K]}S_p$ es suprayectiva, por lo que $S_p$
est\'a en la imagen de $\N_{M/K}$.

\bigskip

\noindent
{\bf Paso 2:} Veamos que podemos suponer que $L/K$ es c\'iclica.
A saber, si $M/K$ recorre las subextensiones c\'iclicas de $L/K$, 
entonces (\ref{diagramaNeukirch}) muestra que el n\'ucleo de 
${\mc N}_{L|K}$ est\'a contenido en el n\'ucleo del mapeo
$\Gal(L/K)\stackrel{\psi}{\lra}\prod_M\Gal(M/K)$. Ahora bien, como
$L/K$ es abeliana, se tiene que $\psi$ es inyectiva y por tanto
${\mc N}_{L|K}$ es inyectiva.

Para ver la suprayectividad, seleccionamos una subextensi\'on
c\'iclica propia $M/K$ de $L/K$, y la suprayectividad ser\'a
obtenida por inducci\'on sobre el grado del campo de la
misma manera que lo fue en el Paso 1 en el caso soluble.

\bigskip

\noindent
{\bf Paso 3:} Sea $L/K$ c\'iclica. Veamos que podemos suponer
que $L/K$ es totalmente ramificada, esto es, $f_{L|K}=1$. Para
ver esta reducci\'on, sea $M=L_0=\nr K\cap L$ y $L/M$ es
totalmente ramificada. Ahora, $M/K$ es no ramificada, por lo
que ${\mc N}_{M|K}$ es un isomorfismo, Por el Lema de 
la Serpiente aplicado al diagrama (\ref{diagramaNeukirch}),
$\N_{M/K}$ es inyectivo (tambi\'en esto puede ser deducido por
el Teorema \ref{axiomacamposlocales}, o el Teorema 
\ref{NeuT8}), los grupos de abajo tienen \'ordenes $[L\colon M],
[L\colon K]$ y $[M\colon K]$.

Por tanto, si ${\mc N}_{L|K}$ es un isomorfismo, ${\mc N}_{L|K}$
lo es. Sea pues $L/K$ c\'iclica totalmente ramificada, $f_{L|K}=1$.
Sea $\sigma$ un generador de $\Gal(L/K)$. Puesto que $\Gal(
\nr L/\nr K)\cong \Gal(L/K)$.
$
\xymatrix{
\nr K\ar@{-}[r]\ar@{-}[d]&\nr L\ar@{-}[d]\\ K=L\cap \nr K\ar@{-}[r]&L}
$
Consideremos $\sigma$ como elemento de $\Gal(\nr L/\nr K)$. 
Por tanto $\tilde\sigma:=\sigma \Fro L\in \Lambda(\nr L/K)$ es un
levantamiento de Frobenius de $\sigma$ con $\deg_K(\tilde\sigma)
=1$. El campo fijo $F/K$ satisface que $F\cap \nr K=K$ y por
tanto $f_{F|K}=1$. Sea $M\subseteq \nr L$ una extensi\'on finita
de Galois de $K$ que contiene a $F$ y a $L$. Sea $\N=\N_{M/M_0}$.
Puesto que $f_{F|K}=f_{L|K}=1$.
\begin{small}
\[
\xymatrix{
F\ar@{-}@/^1pc/[rr]\ar@{-}[rd]&L\ar@{-}[r]\ar@{-}[d]&M\ar@{-}[d]^{\N}
&L\cap M_0=L\cap M\cap \nr K=L\cap \nr  K=L_0=K,\\
&\underbracket[0pt]{K}_{\substack{\uigual\\ L\cap M_0\\ 
\uigual\\ F\cap M_0}}\ar@{-}[r]& M_0
&F\cap M_0=F\cap M\cap \nr K=F\cap \nr K=F_0=K.
}
\]
\end{small}
Por tanto $\N|_F=\N_{F/K}$, $\N|_L=\N_{L/K}$.

Para la inyectividad de ${\mc N}_{L|K}$ debemos probar que si
${\mc N}_{L|K}(\sigma^k)=1$, $0\leq k<n=[L\colon K]$, entonces $k=0$.
Para este fin consideremos $\pi_F\in F$, $\pi_L\in L$ elementos
primos. Puesto que $F,L\subseteq M\subseteq \nr L$, tanto $\pi_F$
como $\pi_L$ son elementos primos de $M$. Sean $\pi_F^k=u\pi_L^k$
con $u\in U_M$. Tenemos
\[
{\mc N}_{L|K}(\sigma^k)\equiv \N\big(\pi_F^k\big)\equiv \N(u)
\N\big(\pi_L^k)\equiv \N(u)\equiv 1\bmod \N_{L/K}\*L.
\]
Por tanto $\N(u)=\N(v)$ con $v\in U_L$. Por tanto $\N(u^{-1}v)=1$ y
por el Teorema 90 de Hilbert existe $a\in M$ con $u^{-1}v=a^{\sigma-1}$.
En $M$, obtenemos $\big(\pi_L^kv\big)^{\sigma-1}=\big(\pi_F^k u^{-1} v
\big)^{\tilde \sigma-1}\igual\limits_{\substack{\uparrow\\ \tilde\sigma(
\pi_F)=\pi_F}} \big(a^{\sigma-1}\big)^{\tilde\sigma-1}=\big(a^{\tilde
\sigma-1}\big)^{\sigma-1}$.

Se obtiene que si $x:=\pi_L^kv a^{1-\tilde\sigma}$, $x^{\sigma-1}=1$
y $x\in\nr K$. Por tanto $x\in \nr K\cap M=M_0$. Ahora bien, 
$v_{M_0}(x)\in \hat{\ma Z}$ y $nv_{M_0}
(x)\igual\limits_{\substack{\uparrow\\ M/M_0\\ {\text{totalmente
ramificada}}}}v_M(x)=k$, se sigue que $k=0$. De esta forma hemos
obtenido que ${\mc N}_{L|K}$ es un mapeo inyectivo.

Finalmente, se tiene que
$|\co 0{\Gal(L/K)}{\*L}|=|\*K/\N_{L/K}\*L|=[L\colon K]=
|\Gal(L/K)|$, por el axioma de la teor\'ia de campos de clase. Por
tanto ${\mc N}_{L|K}$ es suprayectiva y es un isomorfismo.
$\fin$
\end{proof}

Tenemos que ${\mc N}_{L|K}$ es el inverso del mapeo de reciprocidad
\[
(\underline{\ },L/K)\colon \*K\lra \abe{\Gal(L/K)}
\]
con n\'ucleo $\N_{L/K}\*L$, el cual es tambi\'en llamado el s\'imbolo
de la norma residual. En particular ${\mc N}_{L|K}$ coincide con el
mapeo de Nakayama: si $u_{L|K}$ es la clase fundamental de $L/K$:
\[
\theta_{L|K}=u_{L|K}\Cup\underline{\ }\colon \co {{-2}}{\Gal(L/K)}{\ma Z}
\lra \co 0{\Gal(L/K)}{\*L}.
\]

Los diagramas conmutativos del Teoremas \ref{NeuT15} y
de las Proposiciones \ref{NeuP16} y \ref{NeuP17}, se traducen en
general de manera inmediata por el hecho de que ${\mc N}_{L|K}$
es un isomorfismo y de que $(\underline{\ },L/K)/\N_{L/K}\*L=
{\mc N}_{L|K}^{-1}$, de la siguiente forma.

\begin{teorema}\label{NeuT19}
Sean $L/K$ y $L'/K'$ extensiones finitas de Galois tales que $K
\subseteq K'$ y $L\subseteq L'$ y sea $\sigma\in Gal(\sep k/k)$.
Se tienen los diagramas conmutativos
\begin{gather*}
\xymatrix{
\*{(K')}\ar@{->}[rr]^(0.4){(\underline{\ },L'/K')}\ar@{->}[d]_{\N_{K'/K}}
&&\abe{\Gal(L'/K')}\ar@{->}[d]^{\rest|_L}\\
\*K\ar@{->}[rr]^(0.4){(\underline{\ },L/K)}&&\abe{\Gal(L/K)}
}\\
\xymatrix{
\*K\ar@{->}[rr]^(0.4){(\underline{\ },L/K)}\ar@{->}[d]_{\sigma}
&&\abe{\Gal(L/K)}\ar@{->}[d]^{\*\sigma}\\
\*{(K^{\sigma})}\ar@{->}[rr]^(0.4){(\underline{\ },
L^{\sigma}/K^{\sigma})}&&\abe{\Gal(L^{\sigma}/K^{\sigma})}
}
\intertext{donde $\*\sigma(\tau)=\sigma^{-1}\tau\sigma$. 
Finalmente, si $K'\subseteq L$, entonces}
\xymatrix{
\*{(K')}\ar@{->}[rr]^(0.4){(\underline{\ },L/K')}\ar@{<-^{)}}[d]
&&\abe{\Gal(L/K')}\ar@{<-}[d]^{\Ver}\\
\*K\ar@{->}[rr]^(0.4){(\underline{\ },L/K)}&&\abe{\Gal(L/K)}
}
\end{gather*}
$\fin$
\end{teorema}

Pasando al l\'imite proyectivo, el s\'imbolo de la norma residual se
extiende a todas las extensiones de Galois $L/K$, no necesariamente
finitas. M\'as precisamente, si $\{L_i/K\}$ recorre todas las
subextensiones finitas de Galois de $L/K$, entonces $\abe{
\Gal(L/K)}=\lim\limits_{\longleftarrow i} \abe{\Gal(L_i/K)}$ y 
entonces el s\'imbolo de la norma residual $(a,L_i/K)$, con
$a\in\*K$, determina un elemento $(a,L/K)\in \Gal(L/K)$ con
$(a,L/K)|_{L_i}=(a,L_i/K)$ puesto que tenemos que
$(a,L_j/K)|_{L_i}=(a,L_i/K)$ para $L_i\subseteq L_j$
(Teorema \ref{NeuT19}). En particular, para la extensi\'on
$\nr K/K$, el siguiente resultado describe a $(a,\nr K/K)$.

\begin{teorema}\label{NeuT20}
Se tiene que $\deg\circ (\underline{\ },\nr K/K)=v_K$, es decir,
en particular
\[
(a,\nr K/K)=\Fro K^{v_K(a)}.
\]
\end{teorema}

\begin{proof}
Se tiene que $(\pi_K,\nr K/K)=\Fro K$ pues $(\pi_K,\nr K/K)|_L=
(\pi_K,L/K)=\Fr LK=\Fro K|_L$ para toda subextensi\'on finita
$L/K$ de $\nr K/K$ y $(u,\nr K/K)=1$ para $u\in U_K$. Por tanto,
si $a=\pi_K^m u$ con $u\in U_K$, se tiene
\begin{gather*}
\deg_K(a,\nr K/K)=\deg_K(\Fro K^m)=m=v_K(a).
\tag*{$\fin$}
\end{gather*}
\end{proof}

El Teorema \ref{NeuT20} muestra que la valuaci\'on $v$ y, por
tanto $v_K$, es equivalente con la prescripci\'on del s\'imbolo
de la norma residual $(\underline{\ },\nr K/K)$. Por tanto, la
teor\'ia desarrollada por Neukirch se puede interpretar de la
siguiente forma. Si una teor\'ia de campos de clase se desarrolla
para $\hat{\ma Z}$-extensiones $\nr K/K$, entonces, con el
axioma de teor\'ia de campos de clase, la teor\'ia se extiende 
autom\'aticamente, y de manera \'unica, a todas las
extensiones abelianas $L/K$.

\subsection{Aplicaciones del TCCL}\label{Sub17.6.5}

Regresamos al TCCL. Tenemos las siguientes correspondencias biyectivas:
\begin{gather*}
\text{\{extensiones abelianas finitas de $K$\}} \longleftrightarrow\\
\text{\{subgrupos abiertos de $\Gal(\abe K/K)$\}} \longleftrightarrow\\
\text{\{subgrupos abiertos de {\'\i}ndice finito en $\*K$\}}.
\end{gather*}

Las correspondencias est\'an dadas por
\begin{gather*}
L\longleftrightarrow\Big(\ker:\Gal(\abe K/K)\xrightarrow
{\rest|_L} \Gal(L/K)\Big)\cong
\Gal(\abe K/L) \longleftrightarrow\\
\longleftrightarrow \ker\Big(\*K\xrightarrow{\psi_{L|K}} \Gal(L/K)\Big)=\N_{L/K}\*L
\end{gather*}
donde la \'ultima igualdad proviene del isomorfismo 
\[
\*K/\N_{L/K}\*L\xrightarrow[\psi_{L|K}]{\cong}\Gal(L/K).
\]

Notemos que en la correspondencia de Galois, los subgrupos abiertos
de $\Gal(\abe K/K)$ son de \'indice finito y por tanto cerrados (ver
Observaci\'on \ref{CClaseO3.2.3}). En la correspondencia de Galois
obtenemos todas las extensiones abelianas de $K$ por medio de todos
los subgrupos cerrados.

\begin{observacion}\label{CClaseO3.2.3} Si $N$ es un subgrupo abierto de 
$\Gal(\abe K/K)$, entonces $\Gal(\abe K/K)= G=\bigcup_{x\in G} Nx$.
Como $\{Nx\}_{x\in G}$ es una cubierta abierta de $G$ y
$G$ es compacto, existe una subcubierta finita y $G=\bigcup_{i=1}^m
Nx_i$ por lo que $[G:N]<\infty$, es decir, todo subgrupo abierto
de $G$ es de {\'\i}ndice finito.

Tambi\'en tenemos que $N$ es cerrado pues si $G=N\bigcup \big(
\bigcup_{j=1}^n Ny_j\big)$ uni\'on disjunta, entonces 
$N=G\setminus\big(\bigcup_{j=1}^n Ny_j\big)$ el cual es cerrado.

Un subgrupo cerrado de {\'\i}ndice finito es abierto, lo cual se
demuestra de la misma forma. Sin embargo existen subgrupos
cerrados que no son de {\'\i}ndice finito y por tanto no son abiertos.
Por ejemplo, existe una extensi\'on $L$ de ${\ma Q}(\zeta_3)$
tal que $\Gal(L/{\ma Q}(\zeta_3))\cong {\ma Z}_3\times {\ma Z}_3
=\langle \sigma,\theta\rangle$ y $\langle\sigma\rangle\cong{\ma
Z}_3$ es cerrado pero no abierto en $G$ y
se tiene $[G:\langle\sigma\rangle]=\infty$.
\end{observacion}

\begin{corolario}\label{CClaseC3.3.4} Sea $K$ un campo local. Entonces 
hay una correspondencia biyectiva
\begin{multline*}
\text{$\{$extensiones abelianas finitas de $K\}$}\longleftrightarrow\\
\text{$\{$subgrupos abiertos de {\'\i}ndice finito de $\*K\}$}
\end{multline*}
la cual est\'a dada de la siguiente forma: a la extensi\'on abeliana
finita $L/K$ le corresponde el subgrupo $\N_{L/K}\*L$ de $\*K$
($L\longleftrightarrow \N_{L/K}\*L$).

M\'as a\'un esta correspondencia satisface que si a $L$ le corresponde
$H$ ($L\longleftrightarrow H$), entonces $[L:K]=[\*K:H]$ y adem\'as
si a $L^{\prime}$ le corresponde $H^{\prime}$ ($L^{\prime}
\longleftrightarrow H^{\prime}$), se tiene que $L\supseteq L^{\prime}
\iff H\subseteq H^{\prime}$. $\fin$
\end{corolario}

\begin{observacion}\label{CClaseO3.3.6} Se tiene que $U_K$ es un conjunto
abierto de $\*K$ y adem\'as $\*K/U_K\cong {\ma Z}$ por lo que $U_K$ no
puede corresponder a ninguna extensi\'on abeliana de $K$ pues
${\ma Z}$ no es un grupo profinito y por lo tanto no puede ser el
grupo de Galois de ninguna extensi\'on.
\end{observacion}

\begin{notacion y definicion}\label{CClaseND3.3.8} El mapeo
$\rest_L\circ\rho_K\colon\*K\to \Gal(L/K)$ se llama el {\em s{\'\i}mbolo
residual de la norma\index{simbolo de la norma residual@s\'imbolo de la norma residual}}
o {\em s{\'\i}mbolo residual n\'ormico\index{simbolo residual normico@s\'imbolo residual
n\'ormico}} (norm residue symbol en ingl\'es) o {\em mapeo local
de Artin\index{mapeo local de Artin}\index{Artin!mapeo local
de $\sim$}\label{CClasesimboloresidual}}
y se denota $\*K\xrightarrow{(\ ,L/K)}\Gal(L/K)$. Se puede considerar
a $(\ ,L/K)$ o a $\psi_{L/K}$
como el s{\'\i}mbolo de Artin local o mapeo local de Artin.
\end{notacion y definicion}

\begin{definicion}\label{CClaseD3.2.1.1}
Sea $K={\ma R}$. Se define $\unidades 0=\*{\ma R}$ y $\unidades 1
={\ma R}^+$. Para $K={\ma C}$, se define $\unidades 0=\*{\ma C}$.
\end{definicion}

\subsubsection{Ley de Reciprocidad para $K={\ma R}$ y para $K={\ma C}$.}

\begin{proposicion}\label{CClaseP3.2.1.2}
La ley de reciprocidad se cumple para $K={\ma R}$ y para $K={\ma C}$.
\end{proposicion}

\begin{proof} Si $K={\ma R}$, se tiene $\abe {\ma R}={\ma C}$ y
${\ma R}$ y ${\ma C}$ son las \'unicas dos extensiones 
abelianas (y de hecho algebraicas) de ${\ma R}$. Se tiene que
$\*{\ma R}$ tiene \'unicamente dos subgrupos de {\'\i}ndice finito
los cuales son ${\ma R}^+$ y $\*{\ma R}$ Adem\'as se cumple que
$\N_{{\ma R}/{\ma R}}\*{\ma R}=\*{\ma R}$ y $\N_{{\ma C}/{\ma R}}
\*{\ma R}={\ma R}^+=(\*{\ma R})^2$.

Sea $\rho_{\ma R}\colon \*{\ma R}\longrightarrow \Gal(\abe {\ma R}/{\ma R})
=\Gal({\ma C}/{\ma R})=\{1, J\}$ dada por
\[
\rho_{\ma R}(x)=\sgn(x)=\begin{cases}
1&\text{si $x>0$},\\ J=-1&\text{si $x<0$.}
\end{cases}
\]
Entonces $\rho_{\ma R}$ cumple las condiciones del Teorema TCCL
para $K={\ma R}$.

Si ahora consideramos $K={\ma C}$, ${\ma C}$ es la \'unica
extensi\'on algebraica de ${\ma C}$ y el \'unico subgrupo abierto
de {\'\i}ndice finito en $\*{\ma C}$ es $\*{\ma C}$. Por tanto
$\rho_{\ma C}\colon \*{\ma C}\longrightarrow \Gal({\ma C}/{\ma C})=
\{1\}$, $z\mapsto 1$ satisface las condiciones del Teorema TCCL
para $K={\ma C}$. $\fin$
\end{proof}

\begin{proposicion}\label{CClaseP3.2.1.3}
Sea $K$ un campo local con campo residual $\F$
y valuaci\'on $v=v_\pK$. Sea $\o_K$ el anillo de valuaci\'on de $K$.
Entonces
\[
\xymatrix{
*\txt{{\rm{\{extensiones abelinas fini-}}\\ {\rm{tas no ramificadas de $K$\}}}}
\ar@{<->}[rr]_-{\txt{\rm{\scriptsize{Corolario}} \\ 
{\rm{\scriptsize\ref{C17.5.33'}}}}}
\ar@{<->}[dd]_-{\txt{\rm{\scriptsize{teor\'ia
de cam-}}\\ {\rm{\scriptsize{pos locales}}}}}^-{\txt{\rm{\scriptsize{Teorema}}\\ 
{\rm{\scriptsize{\ref{T17.3.3.1}}}}}}
&&*\txt{{\rm{\{subgrupos abiertos de {\'\i}ndice fi-}}\\ {\rm{nito de
  $\*K$ que contienen a $U_K$\}}}}
\ar@{<->}[dd]^{\*K/U_K\cong {\ma Z}}\\
\\
*\txt{{\rm{\{extensiones abelianas}}\\
{\rm{finitas de $\F$\}}}}\ar@{<->}[rr]&&*\txt{{\rm{\{subgrupos abiertos}}\\ 
{\rm{de {\'\i}ndice finito de ${\ma Z}$\}}}}
}
\]
(recordemos (Corolario {\rm{\ref{C17.5.33'}}}) que cuando $L/K$ es 
una extensi\'on abeliana finita no ramificada,
$U_K\subseteq \N_{L/K}\* L$, por eso se tiene la correspondencia dada por
la flecha vertical del diagrama). $\fin$
\end{proposicion}

\begin{observacion}\label{CClaseO3.2.1.4} Se tiene que $\rho_K\colon\*K\to \Gal
(\abe K/K)$ nos proporciona una biyecci\'on entre las extensiones
abelianas finitas de $L/K$ y subgrupos abiertos de {\'\i}ndice
finito en $\*K$: $L\longleftrightarrow \N_{L/K} \*L (U\longleftrightarrow
\rho_K^{-1}(U))$, es decir $\*K/\N_{L/K}\*L\cong\Gal(L/K)$.
\end{observacion}

Si $L$ est\'a dado, es ``f\'acil'' calcular $\N_{L/K}\*L$, pero dado
$H<\*K$ con $H$ subgrupo abierto de {\'\i}ndice finito, ?`como
calcular $L$ tal que $H=\N_{L/K}\*L$?

Ese es el problema que no nos permite dar una descripci\'on
expl{\'\i}cita de todas las extensiones abelianas finitas de $K$.
Por supuesto, si $\rho_K$ se da expl{\'\i}citamente, resolvemos
parcialmente este problema.

Resulta ser que $\rho_K$ es bastante
expl{\'\i}cito si $L/K$ es no ramificada. De hecho se tiene que si
$L/K$ es no ramificada, recordando que estamos en campos 
locales, entonces $\Gal(L/K)\cong \Gal(\tilde{L}/\tilde{K})$
donde $\tilde{L}$ y $\tilde{K}$ son los campos residuales. M\'as a\'un, si
$\tilde{K}=\F$ y $\tilde{L}={\ma F}_{q^n}$ con $n=[\tilde{L}:\tilde{K}]
=[L:K]$, entonces $\rho_K\colon \*K\to \Gal(L/K)$ satisface que
$\rho_K(a)=\tau^{v_\pK(a)}$ (Teorema \ref{CCLT17.6.10}),
donde $\tau$ es el automorfismo
de Frobenius de $L/K$, el cual es el generador de $\Gal(L/K)$
inducido por el automorfismo de Frobenius de $\tilde{L}/\tilde{K}$,
es decir, $\tau\colon\tilde{L}\to\tilde{L}$, $\bar{x}\mapsto 
\bar{x}^q$. Notemos que bajo $\rho_K$ tenemos que si $\pi$
es un elemento primo de $K$ entonces $\pi\longleftrightarrow \tau$.

De hecho, por TCCL (II). que es el contenido del
Teorema \ref{CCLT17.6.10}, se tiene,
$\rho_K(\pi)=\rho_{\F}(v_K(\pi))=\tau$, de donde
\[
\mu\rho_K(a)=\rho_K(a)|_{\abe \F}=\rho_{\F}(v_K(a))=\tau^{v_K(a)}.
\]

Para extensiones ramificadas
la historia es muy diferente y se requiere cohomolog{\'\i}a de
grupos para obtenerla. El Teorema de Tate y el grupo de Brauer nos
dan el mapeo de Nakayama, esto es, $\rho_K^{-1}$. Sin embargo,
$\psi_{L|K}(a)=(a,L/K)$, el s\'imbolo de la norma residual local o
mapeo de Artin local, no es expl\'icito en lo absoluto.
Otras aproximaciones para el estudio de $\rho_K$ son por medio de los grupos formales
de Lubin-Tate (1965) y tambi\'en por medio de \'algebras 
c{\'\i}clicas (Hasse et. at.). Haremos algo expl\'icito el mapeo local
de Artin y obtendremos el Teorema
de Existencia en completa generalidad, por medio
de los grupos formales de Lubin--Tate en la Secci\'on \ref{CClaseS3.3}.

M\'as a\'un, los grupos de Lubin--Tate nos permiten 
encontrar expl\'icitamente la m\'axima extensi\'on abeliana de
un campo local $K$ (Teorema \ref{CClaseT3.2.5.35}). Podemos considerar
este resultado como el an\'alogo al Teorema de Kronecker--Weber
para campos locales.

\subsection{Teorema de Existencia}\label{STeoremadeExistencia}

De la teor\'ia general de campos locales, tenemos que si $L/K$ es una extensi\'on
finita y separable, $\N_{L/K}\colon \*L\lra \*K$ es continua y $\N_{L/K}\*L$ es un subgrupo
cerrado en $\*K$. Por el Teorema de Reciprocidad, se tiene que $L/K$ es
una extensi\'on finita de Galois, entonces $\*K/\N_{L/K}\*L\cong \abe{G_{
L|K}}$ y como $G_{L|K}$ es un grupo finito, $[\*K\colon\N_{L/K}\*L]<\infty$
y por tanto $\N_{L/K}\*L$ es un subgrupo abierto ($\N_{L/K}\*L=\*K\setminus
\underbrace{\big(\bigcup_{x\notin \N_{L/K}\*L}x\N_{L/K}\*L\big)}_{\text{finita}}$).
De hecho, si $A$ es un subgrupo de $\*K$ de \'indice finito, $A$ es
abierto $\iff A$ es cerrado pues $A=\*K\setminus\big(\bigcup_{x\notin A}xA\big)$.

\begin{proposicion}\label{TEP1}
Si $\car K=0$, entonces todo subgrupo de \'indice finito, es abierto y, por
tanto, cerrado.
\end{proposicion}

\begin{proof}
Tenemos que para $m\in{\ma N}$ tal que $\car K\nmid m$, en nuestro caso,
$m\in{\ma N}$ es arbitrario, $(\*K)^m$ es abierto y de \'indice finito en $\*K$
(Proposici\'on \ref{CCLTP17.6.14}). Su $H<\*K$, $[\*K:H]=m<\infty$, entonces
$(\*K)^m\subseteq H$ y $(\*K)^m$ es un subgrupo abierto. Por tanto $H=
\bigcup_{x\in H}x(\*K)^m$ es abierto.
$\fin$
\end{proof}

\begin{observacion}\label{TEO2}
Si $\car K=p>0$, $(\*K)^m$ es abierto para $m$ 
tal que $\car K\nmid m$. Sin embargo el resultado es falso para $m=p$.
Esto es, $(\*K)^p$ no es abierto en $\*K$ pues si
lo fuese, existir\'ia $m\in{\ma N}$ tal que $\unidades m\subseteq (\*K)^p$.
En particular se tendr\'ia que si $\pi$ es un elemento primo de $K$, 
entonces $1+\pi^m=x^p$ para alg\'un $x\in\*K$, por lo que
$\pi^m=x^p-1=(x-1)^p$ y $m=v_K(\pi^m)=pv_K(x-1)$ y en particular
se deber\'ia tener que $p|m$ para toda $m\geq m_0$ pues $\unidades
m\subseteq \unidades {m_0}\subseteq (\*K)^p$, lo cual es
absurdo.

Adem\'as, existen subgrupos de \'indice finito en $\*K$ que no son cerrados
y por tanto no son abiertos. Por ejemplo, tenemos que
$\unidades 1\cong {\ma Z}_p^{\infty}$. Sea $D=\{(\xi_n)_{n=1}^{
\infty}\in {\ma Z}_P^{\infty}\mid \xi_n=0 \text{\ para casi todo $n$}\}$. 
Se tiene que $D$ es denso en $\unidades 1$. Entonces, existe un 
morfismo no trivial de grupos abelianos $\varphi\colon G\lra {\ma Z}/
p{\ma Z}$ tal que $D\subseteq \ker \varphi$ y de aqu\'i se deduce
que $\*K$ tiene subgrupos no cerrados de \'indice $p$.

Por otro lado, se tiene que
$(\*K)^p$ es cerrado pues $\*K ={\ma Z}\times U_K$ por lo que
$(\*K)^p=p{\ma Z}\times (U_K)^p$ y $(\*K)^p/(U_K)^p$ es un subgrupo
discreto de $\*K/U_K$ el cual es Hausdorff por ser $U_K$ un
subgrupo cerrado. Se sigue que $(\*K)^p$ es cerrado en $\*K$
(Lema\ref{discretocerrado}). Tambi\'en se sigue que $(\*K)^p$ no es
de \'indice finito en $\*K$ al ser cerrado pero no abierto.
\end{observacion}

En el Corolario \ref{CClaseC3.2.21'+3} probaremos que 
$D_K=\{1\}$ para $K$ de caracter\'istica $0$.
Para ver que $D_K=\ker\rho_K=\{1\}$ 
en el caso $\car K=p>0$, y que de hecho sirva para cualquier
caracter\'istica, usaremos los
grupos de Lubin-Tate para probar que $\langle\pi_K\rangle\times
\unidades n$ es un grupo de normas
para toda $n\in{\ma N}\cup\{0\}$. Puesto que $\langle \pi_K^f
\rangle\times U_K$ es el grupo de normas de la extensi\'on no
ramificada de $K$ de grado $f$, se tiene que
\begin{gather*}
D_K\subseteq \big(\langle\pi_K^f\rangle\times U_K\big)\cap
\big(\langle \pi_K\rangle\times \unidades n\big)=
\langle\pi_K^f\rangle\cup \unidades n
\intertext{para toda
$f\in{\ma N}$ y para toda $n\in{\ma N}\cup\{0\}$. Por lo tanto}
D_K\subseteq\bigcap_{\substack{f\in{\ma N}\\ n\in{\ma N}\cap\{0\}}}
\big(\langle\pi_K^f\rangle\times \unidades n\big)=\{1\}.
\end{gather*}

En resumen, $\rho_K\colon \*K\lra \abe G_K$ es un monomorfismo.

\bigskip

El Teorema de Existencia establece que si $H<\*K$ es un subgrupo
abierto de \'indice finito, equivalentemente, un subgrupo cerrado de
\'indice finito, existe una extensi\'on abeliana finita $L/K$ tal que
$H=\N_{L/K}\*L$, es decir, $H$ es un {\em grupo de normas}.

\begin{definicion}\label{TED2}
Un {\em grupo de normas\index{normas de normas}} de $\*K$ es un
subgrupo $H<\*K$ de \'indice finito tal que $H=\N_{L/K}\*L$ para
alguna extensi\'on abeliana finita $L$ de $K$.
\end{definicion}

\begin{teorema}\label{GRC2} Si $H$ es un subgrupo de $\*K$ que contiene a
un subgrupo de normas, $H\supseteq\N_{L/K}\*L$,
entonces $H$ es tambi\'en un grupo de normas
y si $H=\N_{M/K}\*M$ entonces $M\subseteq L$.
\end{teorema}

\begin{proof}
Sea $L/K$ una extensi\'on abeliana tal que $\N_{L/K}\*L\subseteq H$,
entonces $\rho_K|_L=\psi_{L|K}\colon \*K\lra G_{L|K}$ es tal que
$\ker \psi_{L|K}=\N_{L/K}\*L$. Como $H$ es abierto y $\psi_{L|K}$
es una biyecci\'on entre los grupos $H$ tales que $\N_{L/K}\*L
\subseteq H\subseteq \*K$ y los subgrupos de $G_{L|K}$, digamos que
$\psi_{L|K}(H)={\mc G}<G_{L|K}$, por lo que ${\mc G}=G_{L|M}$
para alg\'un campo $M$ (de hecho, $M=L^{\mc G}$)
 tal que $K\subseteq M\subseteq L$. Se tiene
\begin{gather*}
\psi_{L|K}|_H\colon H\lra G_{L|K} \quad{\text{por tanto}}\quad
\tilde\psi_{L|K}\colon \*K/H\lra \frac{G_{L|K}}{G_{L|M}}\cong G_{M|K}
\intertext{y}
H=\ker\tilde\psi_{L|K}=\N_{M|K}\*M.
\intertext{Aqu\'i hemos usado la conmutatividad del diagrama}
\xymatrix{
\*K\ar@{->}[rr]^{\psi_{L|K}}_{\cong}\ar@{=}[d]&&G_{L|K}
\ar@{->}[d]^{\rest|_M}&\sigma\ar@{->}[d]\\
\*K\ar@{->}[rr]^{\psi_{M|K}}_{\cong}&&G_{M|K}&\sigma|_M
}
\tag*{$\fin$}
\end{gather*}
\end{proof}

\begin{proposicion}\label{CClaseP3.2.21'-1}
Sea $K$ un campo local de caracter\'istica $0$. Sea $m\in{\ma N}$. Entonces
$(\*K)^m$ es un grupo de normas de $\*K$.
\end{proposicion}

\begin{proof}
Sea $\zeta_m$ una ra\'iz $m$-primitiva de la unidad. Supongamos primero que
$\zeta_m\in\*K$. Para $a\in \*K$, sea $L_a:=K(\sqrt[m]{a})$. Sea $L:=
\bigcup_{a\in\*K}L_a=K\big(\sqrt[m]{\*K}\big)$. Se tiene que $L$ es la
m\'axima extensi\'on abeliana de exponente $m$ (Teorema \ref{CClaseT1.6.2}).
Adem\'as, por la Teor\'ia de Kummer, tenemos que
\begin{align*}
\chi\big(\Gal(L/K)\big)&=\Hom\big(\Gal(L/K),{\ma Q}/{\ma Z}\big)\\
&=\Hom\big(\Gal(L/K),\langle\zeta_m\rangle\big)\cong \*K/(\*K)^m.
\end{align*}

Puesto que $\*K/(\*K)^m$ es finito, $\Gal(L/K)$ es finito y $L/K$ es una
extensi\'on finita. Ahora $\*K/\N_{L/K}\*L\cong \Gal(L/K)$ y $\*K/\N_{L/K}
\*L$ tiene exponente $m$ lo cual implica que $(\*K)^m\subseteq \N_{L/K}
\*L$ y, por otro lado, 
\[
[\*K:(\*K)^m]=\big|\chi\big(\Gal(L/K)\big)\big|=\big|\Gal(L/K)\big|=[L\colon K]
=[\*K\colon\N_{L/K}\*L].
\]
Se sigue que $\N_{L/K}\*L=(\*K)^m$ y por tanto $(\*K)^m$ es un grupo
de normas.

Ahora, si $\zeta_m\notin\*K$, sea $K_1:=K(\zeta_m)$. Por el caso anterior,
sea $L_1:=K_1\big(\sqrt[m]{\*K_1}\big)$ y $\N_{L_1/K_1}\*L_1=(\*K_1)^m$.
Sea $L_0$ cualquier extensi\'on finita de Galois de $K$ con $L_1
\subseteq L_0$, entonces
\begin{align*}
\N_{L_0/K}\*L_0&=\N_{K_1/K}\big(\N_{L_0/K_1}\*L_0\big)\underbracket[0pt]{
\subseteq}_{\substack{\uparrow\\ L_1\subseteq L_0}}\N_{K_1/K}\big(\N_{L_1/
K_1}\*L_1\big)\\
&=\N_{K_1/K}\big((\*K_1)^m\big)\subseteq (\*K)^m.
\end{align*}

Puesto que $(\*K)^m$ contiene a un grupo de normas, por el
Teorema \ref{GRC2}, $(\*K)^m$ es un grupo de normas.
$\fin$
\end{proof}

\begin{corolario}\label{CClaseC3.2.21'+3}
Sea $K$ un campo local de caracter\'istica $0$. Entonces el isomorfismo
de reciprocidad $\rho_K\colon \*K\lra \Gal(\abe K/K)$ es un monomorfismo.
Esto es, $\ker \rho_K=D_K=\bigcap\limits_{\substack{L/K\text{\ abeliana}\\
\text{finita}}}\N_{L/K}\*L=\{1\}$.
\end{corolario}

\begin{proof}
Se tiene $D_K\subseteq \bigcap_{m=1}^{\infty}(\*K)^m=\{1\}$ (Corolario
\ref{mpotencias}).
$\fin$
\end{proof}

\begin{teorema}[Teorema de Existencia en caracter\'istica 
$0$]\label{CClaseT3.2.21'+4}
Sea $K$ un campo local de caracter\'istica $0$. Los grupos de
normas de $\*K$ son precisamente los subgrupos abiertos (y por tanto
cerrados) de \'indice finito en $\*K$.
\end{teorema}

\begin{proof}
Si $H$ es un grupo de normas, $H=\N_{L/K}\*L$ con $L/K$ una extensi\'on
finita abeliana. Se tiene $[\*K\colon H]=[\*K\colon \N_{L/K}\*L]=[L\colon K]
<\infty$ y es cerrado (Teorema \ref{T17.3.2.7}), por tanto abierto.

Rec\'iprocamente, si $H$ es un subgrupo abierto de \'indice finito en $\*K$,
$[\*K\colon H]=m$, entonces $(\*K)^m\subseteq H$ y como $(\*K)^m$
es un grupo de normas, $H$ tambi\'en lo es (Teorema \ref{GRC2}).
$\fin$
\end{proof}

Notemos que no necesitamos que $H$ sea abierto. De hecho, en la 
Proposici\'on \ref{TEP1}, vimos que en caracter\'istica $0$, todo
subgrupo de \'indice finito, es abierto y por tanto cerrado.

\begin{teorema}\label{CClaseT3.2.21'+6}
Sea $K$ un campo local de caracter\'istica $0$. Entonces los grupos
de normas son precisamente los grupos que contienen a alg\'un
$\langle\pi^f\rangle\times \unidades n$ para $n\in{\ma N}\cup\{0\}$,
$f\in{\ma N}$.
\end{teorema}

\begin{proof}
Se tiene que $\langle\pi^f\rangle\times \unidades n$ tiene \'indice
finito en $\*K$: 
\[
\*K=\langle\pi\rangle\times \unidades 0, \quad [\*K\colon
\langle\pi^f\rangle\times\unidades n]=\begin{cases}
f(q-1)q^{n-1}&(n\geq 1),\\ f& (n=0)\end{cases}
\]
y por tanto  $\langle\pi^f\rangle\times \unidades n$ es un grupo
de normas. Se sigue que
si $H$ contiene a $\langle\pi^f\rangle\times \unidades n$, $H$
es un grupo de normas.

Rec\'iprocamente, si $H$ es un grupo de normas, $H$ es abierto y como
$\{\unidades n\}_{n\in{\ma N}\cup \{0\}}$ es un sistema fundamental de
vecindades de $1$, existe $n\in{\ma N}\cup\{0\}$ con $\unidades n
\subseteq H$. Ahora bien, si $f=[\*K\colon H]$, entonces $\pi^f\in H$
por lo que $\langle\pi^f\rangle\times \unidades n\subseteq H$.
$\fin$
\end{proof}

\begin{corolario}\label{CClaseC3.2.21'+5}
Sea $H<\*K$, $K$ un campo local de caracter\'istica $0$. Entonces
lo siguiente es equivalente:
\las
\item $H$ es un grupo de normas;
\item $H$ es un subgrupo abierto de \'indice finito;
\item $H$ es un subgrupo cerrado de \'indice finito;
\item $H$ es un subgrupo de \'indice finito. 
\item $H$ contiene a $(\*K)^m$ para alg\'un $m\in {\ma N}$.
\item $H$ contiene a $\langle\pi\rangle\times\unidades n$ para
algunos $f\in{\ma N}$ y $n\in{\ma N}\cup\{0\}$.
$\fin$
\end{list}
\end{corolario}

\begin{definicion}[Conductor local\index{conductor
local}]\label{CClaseD3.2.23} Dada una extensi\'on abeliana finita
$L/K$ de campos locales, se tiene que $\N_{L/K}\*L$ es un subgrupo
abierto pues es cerrado de \'indice finito (Teorema \ref{T17.3.2.7}) 
y que contiene a $1$. Por tanto $\unidades n\subseteq
\N_{L/K}\*L$ para alguna $n\geq 0$ pues $\{\unidades n\}_{
n\in{\ma N}}$ es un sistema fundamental de vecindades de
$1$ (ver Subsecci\'on \ref{CClaseS1.2.1}). Sea $n_0$ el m{\'\i}nimo
entero no negativo tal que $\unidades {n_0}\subseteq \N_{L/K}\*L$.

Entonces denotamos $n_0:=c_\pK$ con $\pK=\pK_K$ el lugar de $K$.
Se define el {\em conductor local\index{conductor local}\label{CClaseconductorlocal}}
de $L/K$ por
\[
\f{L/K}=\f{\pK}=\f{}=\pK_K^{n_0}=\pK^{n_0}=\pK^{c_\pK}.
\]
\end{definicion}

Si $L/K$ es una extensi\'on de campos globales y consideramos
las completaciones $L_\pL/K_\pK$, entonces denotamos
$\f {L_\pL/K_\pK}=\f{\pK}$.

\begin{teorema}\label{CClaseT3.2.24} Una extensi\'on abeliana finita de
campos locales $L/K$ es no ramificada $\iff \f{L/K}=\f{}=1$, esto
es, $\iff c_\pK=n_0=0$.
\end{teorema}

\begin{proof}
Por el Corolario \ref{C17.5.33'}, $L/K$ es no ramificada $\iff
U_K\subseteq \N_{L/K}\*L \iff n_0=c_{\pK}=0$. $\fin$
\end{proof}

\subsubsection{Red de normas y de subcampos}

\begin{teorema}\label{CClaseT3.2.21'} Sean $L/K$ y $L'/K'$ dos extensiones
de Galois tales que $K\subseteq K'$ y $L\subseteq L'$. Entonces el
siguiente diagrama es conmutativo
\[
\xymatrix{
\Gal^{\ab}(L'/K')\ar[rr]_{\psi_{L'/K'}^{-1}\phantom{xxx}}\ar[d]_{\rest_L}&&
\*{(K')}/\N_{L'/K'}(\*{(L')})\ar[d]^{\N_{K'/K}}\\
\Gal^{\ab}(L/K)\ar[rr]_{\psi_{L/K}^{-1}\phantom{xxx}}&& \*K/\N_{L/K} (\*L)
}
\]

Equivalentemente, tenemos el diagrama conmutativo
\[
\xymatrix{
\*{(K')}\ar[rr]_{\psi_{L'/K'}\phantom{xxx}}\ar[d]_{\N_{K'/K}}&&
\Gal^{\ab}(L'/K')\ar[d]^{\rest_L}\\
 \*K\ar[rr]_{\psi_{L/K}\phantom{xxx}}&&\Gal^{\ab}(L/K)
 }
\]
\end{teorema}

\begin{proof}
Es el Teorema \ref{NeuT19}.
$\fin$
\end{proof}

Como consecuencia del Teorema \ref{CClaseT3.2.21'} tenemos el
siguiente resultado (para el caso global, ver el Teorema\ref{CClaseT4.6.9-1}).

\begin{teorema}\label{CClaseT3.2.21'+1}
Sea $E/F$ una extensi\'on abeliana finita de campos locales y sea 
$E$ el campo de clase de $\Lambda\subseteq \*F$,
es decir, $\N_{E/F}\*E=\Lambda$. Sea $L/F$
una extensi\'on finita y separable. Entonces $LE/L$ es una extensi\'on
finita y el grupo de normas correspondiente es $\N^{-1}_{L/F}(\Lambda)$.
\[
\xymatrix{
L\ar@{-}[r]\ar@{-}[d]&LE\ar@{-}[d]\\ F\ar@{-}[r]_{\Lambda} & E}
\]
\end{teorema}

\begin{proof}
Sea $\psi_{EL/L}\colon \*L\lra\Gal(LE/L)$ el mapeo de Artin. El grupo de
normas correspondiente a $LE/E$ es $\ker \psi_{LE/L}$, esto es,
$\*L/\ker \psi_{LE/E}\cong \Gal(LE/L)$. Por el Teorema \ref{CClaseT3.2.21'}
tenemos $\rest\circ
\psi_{LE/L}=\psi_{E/F}\circ\N_{L/F}$. Por tanto
\begin{gather*}
x\in\ker\psi_{LE/E}\iff \psi_{LE/E}(x)=1\iff \\
\rest\circ \psi_{LE/E}(x)=1=\psi_{E/F}
\circ\N_{L/F}(x)\iff\\
\N_{L/F}(x)=\ker\psi_{E/F}=\Lambda\iff x\in \N^{-1}_{L/F}
(\Lambda). \tag*{$\fin$}
\end{gather*}
\end{proof}

Para probar el Teorema de Existencia en general, primero demostramos el 
siguiente importante resultado.

\begin{teorema}\label{CClaseT3.2.29} Sea $K$ un campo local. Los
grupos de normas $H_L:=\N_{L/K}\*L$ de $\*K$
forman una red y el mapeo $L\longrightarrow H_L$ es una
correspondencia biyectiva que cambia contenciones y que
de hecho es un isomorfismo de redes entre la red de las
extensiones abelianas finitas de $K$ y la red de los grupos
de normas de $\*K$. Por tanto
\las
\item $H_{L_1}\supseteq H_{L_2}\iff L_1\subseteq L_2$;
\item $H_{L_1L_2}=H_{L_1}\cap H_{L_2}$;
\item $H_{L_1\cap L_2}=H_{L_1}H_{L_2}$.
\end{list}
\end{teorema}

\begin{proof} Si $L_1$ y $L_2$ son dos extensiones 
abelianas finitas de $K$
entonces se tiene $\N_{L_1L_2/K}=\N_{L_i/K} \N_{L_1L_2/L_i}$,
$i=1,2$, de donde obtenemos que $H_{L_1L_2}\subseteq
H_{L_1}\cap H_{L_2}$.

Ahora si $a\in H_{L_1}\cap H_{L_2}$, entonces $(a,L_1L_2/K)|_{L_i}
=(a,L_i/K)=1$, para $i=1,2$.

Se tiene el diagrama con filas exactas
\[
\begin{CD}
1@>>>\N_{L_1L_2/K}\big(\*{(L_1L_2)}\big)@>>> \*K @>>{(a,L_1L_2/K)}>
G_{L_1L_2/K}@>>>1\\
&&&&&&@VV{\pi_i=\rest|_{L_i}}V\\
1@>>> \N_{L_i/K}\*L@>>> \*K@>>{(a,L_i.K)}> G_{L_i/K}@>>> 1.
\end{CD}
\]
y tambi\'en tenemos el diagrama conmutativo
\[
\xymatrix{
G_{L_1L_2/K}\ \ar@{^{(}->}[rr]^(.4){\pi_1\times \pi_2}&&G_{L_1/K}\times G_{L_2/K}\\
\\
\*K\ar@{>}[uu]^{(\underline{\ \ },L_1L_2/K)}\ar@{>}[rruu]_{\quad \big(({\underline{\ \ }}, L_1/K),
({\underline{\ \ }}, L_2/K)\big)}}
\]

Por lo tanto $(a,L_1L_2/K)=1$. As{\'\i}, $a\in H_{L_1L_2}$.

De esta manera se obtiene que
$H_{L_1}\cap H_{L_2}=H_{L_1L_2}$.

Ahora $H_{L_1}\supseteq H_{L_2}\iff H_{L_1}\cap H_{L_2}=H_{L_2}=
H_{L_1L_2}\iff [L_1L_2:K]=[L_2:K]\iff L_1L_2=
L_2\iff L_1\subseteq L_2$.

De esta forma tenemos que $L\longrightarrow H_L$ es una
biyecci\'on que cambia contenciones pues si $H_{L_1}=H_{L_2}$
entonces $H_{L_1L_2}=H_{L_1}\cap H_{L_2}=H_{L_1}=H_{L_2}$
por lo que $L_1L_2=L_1=L_2$. 

Finalmente, tenemos que
$H_{L_i}\subseteq H_{L_1\cap L_2}$, $i=1,2$. Por tanto se
sigue que $H_{L_1}H_{L_2}\subseteq H_{L_1\cap L_2}$.
\[
\xymatrix{
&L_1\ar@{-}[r]\ar@{-}[d]&L_1L_2\ar@{-}[d]\\
&L_1\cap L_2\ar@{-}[ld]\ar@{-}[r]&L_2\\ K}
\]

Ahora $H_{L_i}\subseteq H_{L_1}H_{L_2}$ y $H_{L_i}$
es un subgrupo abierto de $\*K$ de {\'\i}ndice finito ya que
$\*K/H_{L_i}\cong \Gal(L_i/K)$. Por tanto se sigue que
$H_{L_1}H_{L_2}$ es un subgrupo abierto de $\*K$
ya que $H_{L_1}H_{L_2}=\bigcup_{a\in 
H_{L_2}} aH_{L_1}$ es una uni\'on
de conjuntos abiertos (ver tambi\'en la 
Observaci\'on \ref{CClaseO3.2.3}).

Ahora consideremos la imagen de $H_{L_1}H_{L_2}$ bajo el s{\'\i}mbolo
residual de la norma del subgrupo correspondiente a las extensiones
$L_i/K$, $i=1,2$, las cuales son
$(H_{L_1}H_{L_2},L_i/K)$, $i=1,2$, $({\underline{\ \ }},L_i/K)\colon
\*K\to \Gal(L_i/K)$, la cual corresponde a alg\'un subcampo
$K\subseteq T\subseteq L_i$. Notemos que $T$ corresponde al mismo campo
tanto de $L_1$ como de $L_2$ pues en ambos casos
$(\underline{\ \ },L_i/K)$ es
la restricci\'on de $\rho_K$, en otras palabras, $T$ corresponde
a $H_{L_1}H_{L_2}$ o, m\'as precisamente, $\rho_K^{-1}\big(\Gal(\abe K/T)
\big)=H_{L_1}H_{L_2}$.

Lo anterior implica que $H_{L_1}
H_{L_2}=\ker ({\underline{\ \ }},T/K)=H_T$. Por otro lado
$T\subseteq L_1\cap L_2$, $H_{L_1\cap L_2}\subseteq
H_T=H_{L_1}H_{L_2}$. Se sigue que 
$H_{L_1\cap L_2}=H_{L_1}H_{L_2}$. $\fin$
\end{proof}

\begin{teorema}[limitaci\'on de normas local]\label{T17.5.46'}
Sea $L/K$ una extensi\'on finita y separable de campos locales.
Sea $M=L^{\ab}$ la m\'axima extensi\'on abeliana de $K$
contenida en $L$. Entonces $\N_{L/K}\*L=\N_{M/K} \*M$.
\end{teorema}

\begin{proof}
Se tiene que $\N_{L/K}\*L=\N_{M/K}\N_{L/M}\*L\subseteq
\N_{M/K}\*M$. 

En caso de que $L/K$ sea una extensi\'on de
Galois, se tiene que $\Gal(M/K)=\Gal(L/K)^{\ab}=G/[G,G]=
G/G'$, donde $G=\Gal(L/K)$. Entonces, por el Teorema
de Tate (Teorema \ref{CClaseT1.5.15}), se tiene
\[
H^{-2}(G,{\ma Z})\cong G/G'\xrightarrow[\ \theta\ ]{\cong} (\*L)^G/
\N_{L/K}\*L=\*K/\N_{L/K}\*L\cong H^0(G,\*L)
\]
y $[\*K:\N_{L/K}\*L]=|G/G'|=[M:K]=[\*K:\N_{M/K}\*M]$ lo que 
implica que $\N_{L/K}\*L=\N_{M/K}\*M$ y se tiene que
$\mu_{L/K}\colon\*K\lra G/G'$ es el mapeo inducido por el inverso
de $\theta$. De hecho $\mu_{L/K}$ es el mapeo de Artin en
el caso abeliano, es decir $\mu_{L/K}=\psi_{L^{\ab}/K}$.

Para el caso general, sea $\tilde L$ una extensi\'on finita de Galois
de $K$ que contiene a $L$. Sean $G=\Gal(\tilde L/K)$ y $H=
\Gal(\tilde L/L)$.
\[
\xymatrix{
\tilde L\ar@{-}[d]\ar@/^1pc/@{-}[d]^H\ar@/_2pc/@{-}[ddd]_G\\
L\ar@{-}[d]\\ M\ar@{-}[d]\\ K
}
\]
Si $R$ es la m\'axima subextensi\'on abeliana de $K$ contenida
en $\tilde L$, se tiene $R=\tilde L^{G'}$ y $M=R\cap L=\tilde L^{G'}
\cap \tilde L^H=\tilde L^{G'H}$. Si $a\in \N_{M/K}\*M$ se tiene
el diagrama conmutativo (Teorema \ref{CClaseT3.2.21'})
\[
\xymatrix{
\*L\ar@{>}[r]^{\psi_{\tilde L/L}}\ar@{>}[d]_{\N_{L/K}}&H/H'\ar@{>}[d]^i\\
\*K\ar@{=}[d]\ar@{>}[r]^{\psi_{\tilde L/K}}&G/G'\ar@{>}[d]^{\pi}\\
\*K\ar@{>}[r]^{\psi_{M/K}\ \ }&G/G'H
}
\]
donde $i$ es la inclusi\'on natural, $\pi$ la proyecci\'on y $\psi_{
\tilde L/L}, \psi_{\tilde L/K}, \psi_{M/K}$ los mapeos de Artin, los
cuales tienen sentido pues $\tilde L/L$ y $\tilde L/K$ son
extensiones de Galois.

Ahora, puesto que $a\in\N_{M/K}\*M$, se tiene $\psi_{M/K}(a)=1
=\pi\circ \psi_{\tilde L/K}(a)$. Por otro lado $\psi_{\tilde L/L}$ es
suprayectivo y $\psi_{\tilde L/K}(a)\in \ker \pi=\im i$ por lo que
existe $b\in \*L$ tal que $\psi_{\tilde L/K}\circ (\N_{L/K} b)=
\psi_{\tilde L/K}(a)$. Se sigue que $a/\N_{L/K} b\in \ker \psi_{\tilde
L/K}=\N_{\tilde L/K}\*{\tilde L}$, esto es, existe $c\in \*{\tilde L}$ 
con $a/\N_{L/K} b=\N_{\tilde L/K} c$ de donde
\[
a=\N_{L/K}b \cdot \N_{\tilde L/K} c=\N_{L/K} b 
\cdot \N_{L/K}(\N_{\tilde L/L} c)=
\N_{L/K}(b \N_{\tilde L/L} c)\in \N_{L/K}\*L
\]
probando que $\N_{M/K}\*M=\N_{L/K}\* L$.
$\fin$
\end{proof}

\subsection{Grupos de ramificaci\'on superior\index{grupos
de ramificaci\'on superior}
y grupos de normas\index{grupos de normas}}\label{CClaseS3.2.4}

Sea $L/K$ una extensi\'on de Galois finita de campos locales
con grupo de Galois $G=\Gal(L/K)$. Se tiene
que $\o_L=\o_K[x]$ para alg\'un $x\in\o_L$
(ver \cite[Ch. III, Section 6, Proposition 12]{Ser}). Entonces se definen
los {\em grupos de ramificaci\'on\index{grupos de ramificaci\'on}} por
(ver Definici\'on \ref{D1.3.10bis}):
\[
G_i:=\{\sigma\in G=\Gal(L/K)\mid v_L(\sigma x-x)\geq i+1\}, \quad i\geq -1.
\]

Se tiene que $G_{-1}=G$ el cual, en el caso global, corresponde al
grupo de descomposici\'on, $G_0$ es el grupo de inercia y
$G=G_{-1}\supseteq G_0\supseteq G_1\supseteq \ldots\supseteq
G_r\supseteq \ldots$, $G_i\lhd G$ para toda $i\geq -1$ y $G_r=\{1\}$
para $r$ suficientemente grande pues para $\sigma\neq 1$, 
$\sigma x\neq x$ y $v_K(\sigma x-x)<\infty$. 

Se define
\[
\imath_G\colon G\lra {\ma Z}\cup\{\infty\},
\]
de la siguiente forma: para $\sigma\neq 1$, 
$\imath_G(\sigma)=v_K(\sigma x-x)\neq \infty$
y $\imath_G(1)=\infty$. Se tiene
\[
\imath_G(\sigma)\geq i+1\iff \sigma\in G_i,
\]
lo cual prueba que la definici\'on de $G_i$ no depende de $x$.
Adem\'as, puesto que para $\tau\in G$ se tiene $\o_L=\o_K[\tau^{-1}x]$,
entonces
\begin{align*}
\imath_G(\tau\sigma\tau^{-1})&=v_K(\tau\sigma\tau^{-1}x-x)
=v_K(\tau(\sigma\tau^{-1}x-\tau^{-1}x))\\
&=v_K(\sigma(\tau^{-1}x)-(\tau^{-1}x))
=\imath_G(\sigma).
\end{align*}

Para $\sigma,\tau\in G$ se tiene $(\sigma\tau)(x)-x=\sigma(\tau x)
-\tau x+\tau x-x$, de donde
\begin{align*}
\imath_G(\sigma\tau)&=v_K((\sigma\tau)x-x)\geq \min\{v_K(\sigma(\tau x)
-(\tau x)),v_K(\tau x -x)\}\\
&=\min\{\imath_G(\sigma),\imath_g(\tau)\}.
\end{align*}

Para un subgrupo $H<G$, sea $E=L^H$. Entonces $\Gal(L/E)=H$.

\begin{proposicion}\label{CClaseGR.1} Para $\sigma\in H$ se tiene
$\imath_H(\sigma)=\imath_G(\sigma)$ y $H_i=G_i\cap H$ para
toda $i\geq -1$.
\end{proposicion}

\begin{proof} Es inmediato pues $v_K(\sigma x-x)$ no depende de $H$. $\fin$
\end{proof}

\begin{corolario}\label{CClaseGR.2}
Sea $E$ la m\'axima subextensi\'on de $L$ no ramificada sobre $K$,
$K\subseteq E\subseteq L$ y sea $H$ el subgrupo correspondiente a
$E$, es decir, $E=L^H$. Entonces $H=G_0$ y los grupos de ramificaci\'on
de $G$ de \'indice mayores o iguales a $0$ son iguales a aqu\'ellos de 
$H$. $\fin$
\end{corolario}

Notemos que $L/E$ es totalmente ramificada.

\begin{proposicion}\label{CClaseGR.3}
Sea $H\lhd G$. Para $\bar{\sigma}\in G/H$ se tiene
\[
\imath_{G/H}(\bar{\sigma})=\frac{1}{e_{L/K}}\sum_{g\in\bar{\sigma}}
\imath_G(g).
\]
\end{proposicion}

\begin{proof} Si $\bar{\sigma}=1$, $1\in \bar{\sigma}$ y ambos lados de la igualdad
es $\infty$.

Sea $\bar{\sigma}\neq 1$. Sean $\o_L=\o_K[x]$ y $\o_E=\o_K[y]$. Se tiene
\[
e_{L/E}\imath_{G/H}(\bar{\sigma})=e_{L/K}v_E(\bar{\sigma} y-y)=v_L(\sigma y
-y) \text{\ \ y\ \ } e_G(\sigma)=v_L(\sigma x-x).
\]

Se tiene que si $\sigma$ es un representante de $\bar{\sigma}$, entonces
$\bar{\sigma}=\{\sigma\tau\mid\tau\in H\}$. Sean $a:=\sigma y-y$ y
$b:=\prod_{\tau\in H}(\sigma\tau(x)-x)$.

Veamos que $\o_L a=\o_L b$. Una vez probado esto se tendr\'a
$v_L(a)=v_L(b)$ y
\begin{align*}
e_{L/E}\imath_{G/H}(\bar{\sigma})&=v_L(a)=
v_L(b)=\prod_{\tau\in H} v_L(\sigma(\tau(x))-x)\\
&=\sum_{\tau\in H}
\imath_G(\sigma\tau)=\sum_{g\in\bar{\sigma}}\imath_G(g).
\end{align*}

Sea $f(T):=\Irr(x,T,E)\in E[T]$. Entonces $f(T)=\prod_{\tau\in H}
(T-\tau x)$. Entonces $\sigma(f)(T)=\sigma(f(T))=\prod_{\tau\in H}
(T-(\sigma\tau)(x))$. 

Puesto que todos los coeficientes de $\sigma f-f$ son
divisibles por $\sigma y-y$ debido a que $\o_E=
\o_K[y]$ y que por tanto $y$ divide a todos los coeficientes de $f$,
se sigue que $a=\sigma y-y$ divide a $\sigma(f)(x)-f(x)=
\sigma(f)(x)=\prod_{\tau\in H}(x-(\sigma\tau)(x))=\pm b$.

Falta ver que $b$ divide a $a$. Puesto que $\o_L=\o_K[x]$,
escribimos $y=g(x)$ con $g(T)\in\o_K[T]$. Se tiene que 
el polinomio $g(T)
-y\in \o_E[T]$ y $x$ es ra\'iz de $g(T)-y$ por lo que $f(T)\mid 
g(T)-y$. 

Escribamos $g(T)-y=f(T)h(T)$ con $h(T)\in \o_E[T]$.
Por tanto $\sigma(g)(x)-\sigma y=\sigma(f)(x)\sigma(h)(x)$.

Ahora bien, $g(T)\in\o_K[T]$, por lo que $\sigma(g)(T)=g(T)$
y $\sigma(g)(x)=g(x)=y$. Se tiene $b=\pm \sigma(f)(x)$. Por tanto
$-a=y-\sigma y=\pm b \sigma(h)(x)$, esto es $b\mid a$
de donde se sigue el resultado. $\fin$
\end{proof}

\begin{corolario}\label{CClaseGR.4}
Si $H=G_j$ para alg\'un $j\geq 0$, entonces
\[
(G/H)_i=G_i/H \quad \text{para}\quad i\leq j\quad\text{y}
\quad (G/H)_i=\{1\}\quad \text{para}\quad i\geq j.
\]
\end{corolario}

\begin{proof} $\{G_i/H\}_{i\leq j}$ es una filtraci\'on decreciente de
subgrupos de $G/H$. Para $\bar{\sigma}\in G/H$, $\bar{\sigma}
\neq 1$, existe un \'indice $i<j$ tal que $\bar{\sigma}\in G_i/H$ y
$\sigma\notin G_{i+1}/H$. Si $\sigma\in G$ representa a la clase
$\bar{\sigma}$, se tiene que $\sigma\in G_i$ y $\sigma\notin 
G_{i+1}$ de donde se sigue que $\imath_G(\sigma)=i+1$.

Ahora bien, $H=G_j\subseteq G_0$, lo cual implica que $L/E$
es totalmente ramificada, donde $E=L^H$ y $e_{L/E}=[L:E]=
|H|$. Entonces 
\[
\imath_{G/H}(\bar{\sigma})=\frac{1}{e_{L/E}}\sum_{g\in\bar{\sigma}}
\imath_G(g)\igual_{\substack{\uparrow\\ \bar{g}=\bar{\sigma}}}
\frac{1}{e_{L/E}}\sum_{g\in\bar{\sigma}}\imath_G(\sigma)=
\imath_G(\sigma)=i+1.
\]

Esto prueba que las filtraciones $\{G_i/H\}_{i\leq j}$ y
$\{(G/H)_i\}_{i\leq j}$ coinciden. Finalmente, tenemos que
$(G/H)_j=G_j/H=H/H=\{1\}$ de donde se sigue que 
$(G/H)_i=\{1\}$ para $i\geq j$. $\fin$
\end{proof}

\begin{proposicion}[Ver la Proposici\'on \ref{P1.3.14}]\label{CClaseGR.5}
Sea $\pi\in\o_L$, $v_L(\pi)=1$
cualquier elemento primo de $L$. Sea $\sigma\in G$ y consideremos
el mapeo $\sigma\longmapsto \sigma\pi/\pi$. Es mapeo induce,
pasando al cociente, un monomorfismo de grupos
\[
G_i/G_{i+1}\xhookrightarrow{\ \theta_i\ } U_L^{(i)}/U_L^{(i+1)}
\cong \pK_L^i/\pK_L^{i+1}.
\]
Se tiene que $\theta_i$ es independiente de $\pi$.
\end{proposicion}

\begin{proof} Sea $\pi'$ otro elemento primo de $L$. Entonces $\pi'=
a\pi$ para alg\'un $a\in U_L$. Sea $\sigma\in G$. Entonces
\[
\frac{\sigma\pi'}{\pi}=\frac{\sigma\pi}{\pi}\cdot \frac{\sigma u}{u}.
\]
Para $\sigma\in G_i$ tenemos $\sigma u-u\in \pK_L^{i+1}$, esto es,
$\sigma u/u\equiv 1\bmod U_L^{(i+1)}$ por lo que $\theta_i$ es
independiente de $\pi$.

Sean $\sigma,\tau\in G_i$, entonces
\[
\frac{(\sigma\tau)(\pi)}{\pi}=\frac{\sigma(\tau(\pi))}{\pi}=
\frac{\sigma\pi}{\pi}\cdot\frac{\tau\pi}{\pi}\cdot\sigma\Big(
\frac{\tau \pi}{\pi}\Big)\cdot\frac{\pi}{\tau \pi}=
\frac{\sigma\pi}{\pi}\cdot\frac{\tau\pi}{\pi}\cdot\frac{\sigma v}
{v}\]
donde $v=\tau\pi/\pi\in U_L$. Se sigue de lo anterior que
$\sigma v/v\equiv 1\bmod U_L^{(i+1)}$ y por tanto
\[
\frac{(\sigma\tau)(\pi)}{\pi}=\frac{\sigma\pi}{\pi}\cdot
\frac{\tau\pi}{\pi}\bmod U_L^{(i+1)}
\]
lo cual implica que $\theta_i$ es un homomorfismo de grupos.
Finalmente, si $\sigma\in\ker\theta_i$, $\sigma\mapsto \sigma\pi/
\pi\in U_L^{(i+1)}$ por lo que $\sigma\in G_{i+1}$ de donde
se sigue que $\bar{\sigma}=1$. $\fin$
\end{proof}

\begin{corolario}[Ver Corolario \ref{C1.3.14'}]\label{CClaseGR.6} 
Se tiene que $G_0/G_1$ es un
grupo c\'iclico de orden un divisor de $q-1$ y $G_i/G_{i+1}$ es un
$p$--grupo elemental abeliano para $i\geq -1$. En particular
$G_1$ es un $p$--grupo. Aqu\'i $q=p^r$ y el campo residual de
$L$  es $\F$.
\end{corolario}

\begin{proof} Se sigue del hecho de
que $G_0/G_1\subseteq U_L/U_L^{(1)}\cong \*\F$
y $G_i/G_{i+1}\subseteq U_L^{(i)}/U_L^{(i+1)}\cong \F$ para
$i\geq 1$. $\fin$
\end{proof}

\begin{corolario}\label{CClaseGR.6+1}
Sea $L/K$ una extensi\'on finita de Galois con grupo de Galois $G$.
Entonces $G$ es un grupo soluble.
\end{corolario}

\begin{proof}
Sea $G_0$ el grupo de inercia de $G$, $G_0\normal G$ y $G/G_0$
es un grupo c\'iclico pues corresponde a la m\'axima extensi\'on no
ramificada de $K$ contenida en $L$ (Teorema \ref{T17.3.3.1}). Del
Corolario \ref{CClaseGR.6} se sigue que el grupo $G_0$ es soluble,
de donde se sigue que $G$ es soluble.
$\fin$
\end{proof}

\begin{definicion}\label{CClaseGR.7} Si $t\in[-1,\infty)$ definimos $G_t:=
G_{\lceil t\rceil}$ donde $\lceil t\rceil$ es la funci\'on techo, esto es,
$\lceil t\rceil$ es el entero m\'as peque\~no mayor o igual a $t$.
\end{definicion}

Se tiene para $t=-1$, $[G_0:G_{-1}]:=[G_{-1}:G_0]^{-1}=[G:G_0]^{-1}$
y para $-1<t\leq 0$, $[G_0:G_t]=1$.

Sea $g_i=|G_i|$, $i\in{\ma Z}$, $i\geq -1$.

\begin{definicion}\label{CClaseGR.8} 
La {\em funci\'on de Herbrand\index{funci\'on
de Herbrand}} 
\begin{gather*}
\varphi=\varphi_{L/K}\colon [-1,\infty)\lra [-1,\infty)
\intertext{se define por}
\varphi(u)=\int_0^u\frac{dt}{[G_0:G_t]}=\frac{1}{g_0}(g_1+
\cdots +g_m+(u-m)g_{m+1})
\end{gather*}
donde $m\leq u\leq m+1$, $m\in {\ma N}$.
\end{definicion}

En particular se tiene $\varphi(m)+1=\frac{1}{g_0}\sum_{i=0}^m g_i$.

Notemos que para $u\geq -1$, $u\notin{\ma Z}$,
\begin{gather}\label{CClaseGR.7'}
\varphi'_{L/K}(u)=\frac{g_{m+1}}{g_0} \quad\text{donde}\quad
m<u<m+1.
\end{gather}

Se tiene que $\varphi$ es continua, lineal, lineal por tramos,
creciente y c\'oncava y por tanto $\varphi$ es una funci\'on
biyectiva y continua.

\begin{definicion}\label{CClaseGR.9}
Sea $\eta=\eta_{L/K}\colon[-1,\infty)\lra[-1,\infty)$ la inversa
de $\varphi$: $\eta=\varphi^{-1}$. Se define el
{\em n\'umero de ramificaci\'on superior
$v$\index{numeros de ramificacion superiores@n\'umeros
de ramificaci\'on superiores}} por
\[
G^v:=G_{\eta(v)}\quad\text{o, equivalentemente,}\quad 
G^{\varphi(u)}=G_u
\]
y $\varphi(u)$ es el n\'umero de ramificaci\'on superior.
\end{definicion}

Se tiene que $\eta$ es continua, lineal por tramos, creciente y convexa.
Adem\'as, $\eta(0)=0$. Si $v=\varphi(u)$ es un entero, entonces 
$u=\eta(v)$ es tambi\'en un entero. En efecto, si $m\in{\ma Z}$ es tal que
$m\leq u\leq m+1$, entonces
\[
g_0 v=g_1+\cdots +g_m+(u-m)g_{m+1}.
\]

Puesto que $G_{m+1}\subseteq G_i$, $0\leq i\leq m$, $g_{m+1}
| g_i$, $0\leq i\leq m$. Por tanto, puesto que $v\in{\ma Z}$,
$u-m\in{\ma Z}$ y $u\in{\ma Z}$. Adem\'as
\[
\eta(v)=\int_0^v [G^0:G^w] dw.
\]

Una de las razones principales para estudiar los n\'umero de
ramificaci\'on superiores, es que tenemos el siguiente resultado.

\begin{teorema}\label{CClaseGR.10} Si $H\lhd G$, entonces, para toda 
$v\in[-1,\infty)$ se tiene
\[
(G/H)^v=G^v H/H.
\]
\end{teorema}

Para probar el Teorema \ref{CClaseGR.10}, primero probamos

\begin{proposicion}\label{CClaseGR.11} Se tiene
\[
\varphi_{L/K}(t)=\frac{1}{g_0}\sum_{\sigma\in G}\min\{
\imath_G(\sigma),t+1\}-1.
\]
\end{proposicion}

\begin{proof} Sea $\theta(t)=\frac{1}{g_0}\sum_{\sigma\in G}\min\{
\imath_G(\sigma),t+1\}-1$. Entonces $\theta$ es una
funci\'on continua, lineal por tramos, $\theta(0)=\varphi(0)=0$.
Si $m\geq -1$ es un entero y $m<t < m+1$, entonces
se tiene que $t+1\notin {\ma Z}$ y
\[
\inf\{\imath_G(\sigma),t+1\}=\imath_G(\sigma)\iff \imath_G(\sigma)<
t+1\iff \imath_G(\sigma)\leq m+1.
\]

Sea $\imath_G(\sigma)=i\leq m+1\ \big(\iff\sigma\notin G_{m+1}\ \big)$,
por lo tanto 
\begin{gather*}
\begin{align*}
\inf\{\imath_G(\sigma),t+1\}=t+1&\iff \imath_G(\sigma)\geq t+1\in
(m+1,m+2) \\
&\iff \imath_G(\sigma)\geq m+2\iff \sigma\in G_{m+1}.
\end{align*}
\intertext{Por tanto}
\begin{align*}
\theta(t)&\igual_{\substack{\uparrow\\ m<t<m+1}}\frac{1}{g_0}
\sum_{\sigma\in G} \inf\{\imath_G(\sigma),t+1\}-1\\
&=\frac{1}{g_0}
\sum_{\sigma\in G\setminus G_{m+1}}\inf\{\imath_G(\sigma),t+1\}
+\frac{1}{g_0}\sum_{\sigma\in G_{m+1}}\inf\{\imath_G(\sigma),
t+1\}-1\\
&=\frac{1}{g_0}\sum_{\sigma\in G\setminus G_{t+1}}
\imath_G(\sigma)+\frac{|G_{m+1}|}{g_0}(t+1)-1.
\end{align*}
\end{gather*}
Se sigue que $\theta'(t)=\frac{g_{m+1}}{g_0}$.

Ahora $\varphi(u)=\int_0^u\frac{dt}{[G:G_t]}$, lo cual implica que
\[
\varphi'(t)=\frac{1}{[G:G_t]}=\frac{1}{[G:G_{\lceil t\rceil}]}=
\frac{1}{[G:G_{m+1}]}=\frac{g_{m+1}}{g_0}=\theta'(t).
\]

De esta forma tenemos que $(\varphi-\theta)'(t)=0$ para toda
$t\in[-1,\infty)\setminus{\ma Z}$, $\varphi(0)=\theta(0)$ y
$(\varphi-\theta)$ es una funci\'on continua. Se sigue que
$\varphi(t)=\theta(t)$ para toda $t\in[-1,\infty)$. $\fin$
\end{proof}

\begin{teorema}[Herbrand\index{teorema de Herbrand}]\label{CClaseGR.12}
Sea $L/K$ una extensi\'on finita de Galois de campos locales con
grupo $G=\Gal(L/K)$. Sea $E/K$ una subextensi\'on de Galois de
$L/K$ con grupo $\Lambda=G/H=\Gal(E/K)$ donde $H=\Gal(L/E)$.
Entonces 
\[
\frac{G_s H}{H}=\Lambda_t=\Big(\frac{G}{H}\Big)_t\quad
\text{con}\quad t=\varphi_{L/E}(s).\quad
\xymatrix{
L\ar@{-}[d]_H\ar@{-}@/^2pc/[dd]_G\\ 
E\ar@{-}[d]_{G/H=\Lambda}\\ K
}
\]
\end{teorema}

\begin{proof} Sea $\bar{\sigma}\in\Lambda$ y seleccionaremos una preimagen
$\sigma \in G$ que toma el m\'aximo valor $\imath_G(\sigma)$. Esto es,
se selecciona $\sigma\in G$ tal que $\imath_G(\sigma)=\sup\{
\imath_G(g)\mid g\in\bar{\sigma}\}$. Probaremos que
\begin{gather}\label{CClaseGR.13}
\imath_{\Lambda}(\bar{\sigma})=\varphi_{L/E}(\imath_G(\sigma)-1).
\end{gather}

Sea $m=\imath_G(\sigma)$, por lo tanto $\sigma\in
H_{m-1}$. Si $\tau\in H$ est\'a en $H_{m-1}$,
entonces $\imath_G(\tau)\geq m$ pues $H_{i-1}=G_{i-1}\cap H$ para
toda $i$ (Proposici\'on \ref{CClaseGR.1}).

Por tanto $m=\imath_G(\sigma)\geq
\imath_G(\tau\sigma)\geq \min\{\imath_G(\tau),\imath_G(
\sigma)\}=m$. Se sigue que $\imath_G(\tau\sigma)=m$. Ahora si
$\tau\in H$ y $\tau\notin H_{m-1}$, se tiene $\imath_G(\tau)<m$ y
$\imath_G(\tau\sigma)=\imath_G(\tau)$. En ambos casos se
tiene $\imath_G(\tau\sigma)=\min\{\imath_G(\tau),m\}$.

Aplicando la Proposici\'on \ref{CClaseGR.3} obtenemos
\begin{gather}\label{CClaseGR.14}
\imath_{\Lambda}(\bar{\sigma})=\frac{1}{e_{L/E}}\sum_{g\in\bar{\sigma}}
\imath_G(g)=\frac{1}{e_{L/E}}\sum_{\tau\in H}\imath_G(\sigma\tau)=
\frac{1}{e_{L/E}}\sum_{\tau\in H}\min\{\imath_G(\tau),m\}.
\end{gather}
Por otro lado, por la Proposici\'on \ref{CClaseGR.1}, se tiene 
$\imath_G(\tau)=\imath_H(\tau)$ y adem\'as $e_{L/E}=|H_0|$.
De la Proposici\'on \ref{CClaseGR.11} y de la Ecuaci\'on (\ref{CClaseGR.14})
obtenemos
\begin{align*}
\imath_{\Lambda}(\bar{\sigma})&=\frac{1}{h_0}\sum_{\tau\in H}
\min\{\imath_G(\tau),m\}=\frac{1}{h_0}
\sum_{\tau\in H}\min\{\imath_H(\tau),m\}\\
&=\varphi_{L/E}(m-1)+1=\varphi_{L/E}(\imath_G(\sigma)-1)+1
\end{align*}
la cual es la Ecuaci\'on (\ref{CClaseGR.13}).

Se tiene que si $\bar{\sigma}\in G_sH/H$,
entonces existe $\sigma'\in G_s$ tal que 
$\overline{\sigma'}=\bar{\sigma}$. De la Ecuaci\'on (\ref{CClaseGR.13}) 
obtenemos
\begin{gather*}
\bar{\sigma}\in G_sH/H\iff \imath_G(\sigma)\geq s+1\iff \imath_G(\sigma)
-1\geq s\\
\iff \varphi_{L/E}(\imath_G(\sigma)-1)\geq\varphi_{L/E}(s)
\iff \imath_{\Lambda}(\sigma')-1\geq \varphi_{L/E}(s)\\
\iff \sigma'\in
\Lambda_{\varphi_{L/E}(s)}=\Lambda_t. \tag*{$\fin$}
\end{gather*}
\end{proof}

\begin{proposicion}\label{CClaseGR.15}
Sea $E/K$ una subextensi\'on de Galois de $L/K$. Entonces si
$\varphi_{L/K}$ denota la funci\'on de Herbrand y $\eta_{L/K}=
\varphi_{L/K}^{-1}$, entonces
\[
\varphi_{L/K}=\varphi_{E/K}\circ \varphi_{L/E}\quad\text{y}\quad
\eta_{L/K}=\eta_{L/E}\circ \eta_{E/K}.
\]
\end{proposicion}

\begin{proof}
Tenemos $e_{L/K}=e_{E/K} e_{L/E}$. Del Teorema de Herbrand
\ref{CClaseGR.12} se tiene $G_s H/H=G_s/G_s\cap H=G_s/H_s=(G/H)_t$
con $t=\varphi_{L/E}(s)$. Por tanto, puesto que $e_{L/E}=|G_0|$,
$e_{E/K}=|(G/H)_0|$ y $e_{L/E}=|H_0|$,
\[
\frac{1}{e_{L/K}}|G_s|=\frac{1}{e_{E/K}}\Big|\Big(\frac{G}{H}\Big)_t\Big|
\cdot \frac{1}{e_{L/E}}|H_s|.
\]

Entonces, de la Ecuaci\'on (\ref{CClaseGR.7'}) se obtiene
\begin{align*}
\varphi_{L/K}'(s)&=\varphi_{E/K}'(t)\cdot \varphi_{L/E}'(s)\igual_{
\substack{\uparrow\\t=\varphi_{L/E}(s)}}\varphi_{E/K}'(\varphi_{L/E}(s))
\varphi_{L/E}'(s)\\
&=(\varphi_{E/K}\circ \varphi_{L/E})'(s).
\end{align*}

Puesto que $\varphi_{L/K}(0)=(\varphi_{E/K}\circ \varphi_{L/E})(0)=0$,
se sigue que 
\[
\varphi_{L/K}=\varphi_{E/K}\circ \varphi_{L/E}.
\]
Tomando
las funciones inversas, se sigue que
\begin{gather*}
\eta_{L/K}=\varphi_{L/K}^{-1}=(\varphi_{E/K}\circ \varphi_{L/E})^{-1}=
\varphi_{L/E}^{-1}\circ \varphi_{E/K}^{-1}=\eta_{L/E}\circ \eta_{E/K}.
\tag*{$\fin$}
\end{gather*}
\end{proof}

Con estos resultados, estamos en condiciones de probar el Teorema
\ref{CClaseGR.10}: si $H\lhd G$, entonces $(G/H)^v=G^v H/H$ para todo
$v\geq -1$.

\begin{proof} (Teorema \ref{CClaseGR.10}): Sea $s=\eta_{E/K}(t)$. Usando el
Teorema de Herbrand \ref{CClaseGR.12} y la Proposici\'on \ref{CClaseGR.15},
obtenemos
\begin{align*}
\frac{G^t H}{H}&=\frac{G_{\eta_{L/K}(t)}H}{H}\igual_{\substack{\uparrow\\
\text{Herbrand}}}\Big(\frac{G}{H}\Big)_{\varphi_{L/E}(\eta_{L/K}(t))}\\
&=\Big(\frac{G}{H}\Big)_{(\varphi_{L/E}\circ \eta_{L/E}\circ \eta_{E/K})(t)}=
\Big(\frac{G}{H}\Big)_{\eta_{E/K}(t)}=\Big(\frac{G}{H}\Big)^{t}. \tag*{$\fin$}
\end{align*}
\end{proof}

\begin{definicion}\label{CClaseD3.2.26} $t$ se llama {\em salto
superior\index{salto superior}} si $G^t(L/K)\neq G^{t+\epsilon}
(L/K)$ para toda $\epsilon>0$.
\end{definicion}

\subsection[Grupos de Lubin-Tate. S\'imbolo residual
de la norma]{Grupos Formales de Lubin-Tate. C\'alculo del
s\'imbolo residual de la norma\index{grupos de 
Lubin--Tate}\index{Lubin-Tate!grupos de $\sim$}}\label{CClaseS3.3}

Nuestra exposici\'on sobre los grupos de Lubin--Tate sigue muy
de cerca a \cite{Neu69} el cual es a su vez, una exposici\'on
detallada de \cite{LubTat65}.

En 1965, J. Lubin y J. Tate \cite{LubTat65} motivados por la analog\'ia
con la teor\'ia de multiplicaci\'on compleja en curvas el\'ipticas, mostraron
como pueden ser usados los grupos formales sobre campos locales para
probar resultados centrales en campos de clase locales. En esta
secci\'on introducimos los grupos formales y por medio de ellos daremos
las demostraciones de los resultados centrales de la teor\'ia de campos
locales. Como ya mencionamos, esta secci\'on est\'a basada en
\cite{Neu69,LubTat65} y tambi\'en en \cite{Iwa86}.

Los grupos formales son los an\'alogos a las extensiones ciclot\'omicas
del campo ${\ma Q}_p$ de los n\'umeros $p$--\'adicos sobre cualquier
campo local. En lugar de las ra\'ices como el n\'ucleo del mapeo
$\*K\xrightarrow{\ n\ }\*K$, se presenta otra acci\'on y el n\'ucleo son los
llamados {\em puntos de divisi\'on\index{puntos de divisi\'on}} los cuales
son tambi\'en las ra\'ices de cierta $n$--potencia de un mapeo.

Aqu\'i queremos mencionar la gran similitud de los grupos formales
con los {\em m\'odulos de Drinfeld\index{m\'odulos de Drinfeld}}, 
y m\'as espec\'ificamente con el {\em m\'odulo de Carlitz\index{m\'odulo
de Carlitz}} lo cual ser\'a evidente a lo largo de esta secci\'on.

Sea $K$ un campo local con valuaci\'on $v$ y sea $\o_K$ el anillo
de valuaci\'on de $K$, es decir $\o_K=\{x\in K\mid |x|_v\leq 1\}=
\{x\in K\mid v(x)\geq 0\}=\bar{B}(1,0)$.

En general, si $A$ es un anillo conmutativo con unidad, el
anillo conmutativo de {\em series formales\index{anillo de series
formales}} en las variables $X_1,\ldots, X_n$ es 
\begin{gather*}
R=A[[X_1,\ldots,X_n]]=\{f(X_1,\ldots,X_n)=\sum_i
a_{i_1,\ldots,i_n}X_1^{i_1}\cdots X_n^{i_n},\\
a_{i_1,\ldots,i_n}\in A \text{\ y $i=(i_1,\ldots,i_n)$ var\'ia en
todos las $n$--tuplas de enteros $n\geq0$}\}.
\end{gather*}

Sean $f,g\in R$, $d\in{\ma Z}$, $d\geq 0$. Se pone $f\equiv g
\bmod \deg d$ si $f-g$ no tiene t\'erminos de grado total menores a $d$.

Si $f\in R$ y $g_1,\ldots,g_n\in A[[Y_1,\ldots,Y_m]]$, se define
\[
f\circ(g_1,\ldots,g_n)=
f(g_1(Y_1,\ldots,Y_m),\ldots,g_n(Y_1,\ldots,Y_m))\in
A[[Y_1,\ldots,Y_m]].
\]

Si $A$ es un anillo topol\'ogico, se considera $R$ como anillo topol\'ogico de
tal forma que el mapeo
\[
f=\sum_i a_{i_1,\ldots,i_n}X_1^{i_1}\cdots X_n^{i_n}\longmapsto
\{a_{i_1,\ldots,i_n}\}_{\{i_1,\cdots,i_n\geq 0\}}\in \prod_i A
\]
define un homomorfismo continuo de $R$ sobre el producto de una cantidad
numerable de $A$ indexadas por $i=(i_1,\ldots,i_n)$.

Sea ahora $K$ un campo local y tomemos $\o_K$. Sea $\pi$ un elemento
primo. Sea $\tau$ es automorfismo de Frobenius, es decir, el mapeo
inducido por $\tau(x)=x^q$ donde $\F\cong \o_K/\pK\o_K=\o_K/\pi\o_K$.

Para $f,g\in \o_K[[T]]$ decimos que $f\equiv g\bmod \pi$ si $f-g=
\sum_{n=0}^{\infty}c_nT^n$ satisface que $\pi|c_n$ para toda
$n\in{\ma N}\cup\{0\}$, esto es $v(c_n)\geq 1$ para toda $n$.

\begin{definicion}\label{CClaseD3.2.5.-1}
Se define 
\[
{\mc F}_{\pi}=\{f\in\o_K[[T]]\mid f(T)\equiv \pi T\bmod \deg 2
\text{\ y\ } f(T)\equiv T^q\bmod \pi\}\label{CClaseFpi}.
\]
\end{definicion}

En general, un elemento de $f\in{\mc F}_{\pi}$ es de la forma
\[
f(T)=\pi T+\pi a_2 T^2+\cdots+\pi a_{q-1}T^{q-1}+T^q+
\pi T^{q+1}g(T)\text{\ con\ } g(T)\in\o_K[[T]].
\]
El elemento m\'as simple de ${\mc F}_{\pi}$ es $f(T)=\pi T+T^q$.
N\'otese la similitud con el m\'odulo de Carlitz.

Sea $f\in{\mc F}_{\pi}$ arbitrario. Se define
\[
f^n(T)=f^{(n)}(T):=(\underbrace{f\circ\cdots\circ f}_{
n})(T)=f(f(\cdots(f(T)\cdots)))\in\o_K[[T]]
\]
y se define $f^{(0)}(T)=T$.

\begin{definicion}\label{CClaseD3.2.5.0}
Se define $\Lambda_{f,n}=\{\lambda\in\Omega\mid v(\lambda)>0\text{\ y\ }
f^{(n)}(\lambda)=0\}\label{CClaselambda-fn}$ donde $\Omega$ es un cerradura
algebraica fija de $K$\label{CClasetorsion}.
\end{definicion}

\begin{definicion}\label{CClaseD3.2.5.1} Se define el campo $L_{f,n}=K(\Lambda_{f,n})$
para $n=1,2,\ldots$ El campo $L_{f,n}$ se llama {\em el campo de los
$\pi^n$ puntos de divisi\'on del m\'odulo de Lubin--Tate $F_f$\index{campos
de Lubin--Tate}} el cual definiremos m\'as adelante.
\end{definicion}

\begin{observacion}\label{CClaseO3.2.5.2}
Se tiene que 
\[
f^{(n)}(T)=f(f^{(n-1)}(T))=f^{(n-1)}(T)g_n(T)
\]
 para alg\'un $g_n(T)\in \o_K[[T]]$.
Por tanto $\Lambda_{f,n-1}\subseteq
\Lambda_{f,n}$ de donde $L_{f,n-1}\subseteq L_{f,n}$, $n=1,
2,\ldots$
\end{observacion}

\begin{definicion}\label{CClaseD3.2.5.3} Sea $f\in{\mc F}_{\pi}$. Se define
$\Lambda_f:=\bigcup_{n=1}^{\infty}\Lambda_{f,n}$ y $L_f:=
K(\Lambda_f)=\bigcup_{n=1}^{\infty}L_{f,n}$.
\end{definicion}

El objetivo inmediato es probar que $L_{f,n}$ son una especie de campos
ciclot\'omicos, o m\'as precisamente, $L_{f,n}$ es una extensi\'on de tipo
ciclot\'omica de $K$. Esto es, queremos probar que $L_{f,n}/K$ es
una extensi\'on abeliana finita y totalmente ramificada. M\'as a\'un, veremos
que $\N_{L_{f,n}/K}(\*{L_{f,n}})=(\pi)\times \unidades n$.

El Teorema \ref{CClaseT3.2.21'+6} tambi\'en 
se cumple en caracter\'istica positivo,
por lo que los grupos de normas de $\*K$ son exactamente los
grupos conteniendo a alg\'un $(\pi^f)\times\unidades n$ para 
$n=0,1,\ldots,$ y $f=1,2,\ldots$. Como consecuencia, se obtiene
el Teorema de Existencia.

La forma en que lo haremos es como sigue. Se usar\'a una serie
de potencias para hacer de $\Lambda_{f,n}$ un $\o_K$--m\'odulo
de tal forma que la multiplicaci\'on de $\Lambda_{f,n}$ por una unidad
$u\in\*{\o_K}=U_K$ produce una permutaci\'on de
$\Lambda_{f,n}$ la cual induce un
automorfismo de $L_{f,n}$ sobre $K$ que resultar\'a ser $(u^{-1},
L_{f,n}/K)$, el s\'imbolo residual de la norma o mapeo de Artin.

Un resultado central para hacer de $\Lambda_{f,n}$ un $\o_K$--m\'odulo,
es el siguiente teorema.

\begin{teorema}\label{CClaseT3.2.5.4} Sean $f,g\in{\mc F}_{\pi}$ y $L(X_1,\ldots,
X_n)=\sum_{i=1}^n a_i X_i$ una forma lineal con coeficientes en $\o_K$.
Entonces existe una \'unica serie de potencias $F(X_1,\ldots,X_n)\in
\o_K[[X_1,\ldots,X_n]]$ tal que
\begin{gather*}
F(X_1,\ldots,X_n)\equiv L(X_1,\ldots,X_n)\bmod \deg 2,\\
f(F(X_1,\ldots,X_n))=F(g(X_1),\ldots,g(X_n)).
\end{gather*}
\end{teorema}

\begin{proof} Pongamos $X:=(X_1,\ldots,X_n)$ y $g(X)=(g(X_1),\ldots,g(X_n))$.
Sea $F_r(X)\in\o_K[X]$ la serie $F(X_1,\ldots,X_n)$ quitando todos
los t\'erminos de grado total mayor o igual a $r$.

Se tiene que si $f(F(X_1,\ldots,X_n))=F(g_1(X_1),\ldots,g(X_n))$
entonces
\[
f(F_r(X))\equiv F_r(g(X))\bmod \deg (r+1)\quad \text{para toda $r$}
\]
y rec\'iprocamente. De esta forma, $F(X_1,\ldots,X_r)$ es soluci\'on
si y s\'olo si
\begin{gather}
F(X)\equiv L(X)\bmod\deg 2\quad\text{y}\nonumber\\
f(F_r(X))\equiv F_r(g(X)) \bmod \deg (r+1) \quad\text{para toda $r$}.\label{CClaseEc3.2.5.1}
\end{gather}

En otras palabras debemos solucionar la Ecuaci\'on (\ref{CClaseEc3.2.5.1})
para toda $r$. Para $r=1$, definimos $F_1(X)=L(X)$ y se 
cumple la Ecuaci\'on (\ref{CClaseEc3.2.5.1}). Notemos que $F_1(X)$ es
\'unico.

Supongamos que hemos hallado un \'unico 
polinomio $F_r(X)$ satisfaciendo
la Ecuaci\'on (\ref{CClaseEc3.2.5.1}). Esto es, 
\[
f(F_r(X))\equiv F_r(g(X))\bmod \deg (r+1).
\]
Notemos que si existe $F_{r+1}(X)$ satisfaciendo la Ecuaci\'on
(\ref{CClaseEc3.2.5.1}) $\bmod \deg(r+2)$ necesariamente se debe tener
$F_{r+1}(X)=F_r(X)+\varphi_{r+1}(X)$ donde $\varphi_{r+1}(X)$
es un polinomio homog\'eneo de grado $r+1$. Veamos que tal
$\varphi_{r+1}(X)$ existe y es \'unico.

Si $F_{r+1}(X)$ es soluci\'on de la Ecuaci\'on (\ref{CClaseEc3.2.5.1}), entonces
\begin{gather*}
f(F_{r+1}(X))=f(F_r(X)+\varphi_{r+1}(X)).
\intertext{Sea}
f(T)=\pi T+\pi a_2 T^2+\cdots+\pi a_{q-1} T^{q-1}+T^q +
\pi T^{q+1} h(T)
\end{gather*}
con $h(T)\in\o_K[[T]]$. Escribamos $f(T)=\sum_{i=1}^{\infty}
a_iT^i$ con $a_1=\pi$, entonces
\begin{align*}
f(F_r(X)&+\varphi_{r+1}(X))=\sum_{i=1}^{\infty}a_i\big(F_r(X)+\varphi_{r+1}(X)\big)^i\\
&=\sum_{i=1}^{\infty}a_i\Big(F_r(X)^i+\sum_{j=1}^i\binom{i}{j}F_r(X)^{i-j}
\varphi_{r+1}(X)^j\Big)\\
&=\sum_{i=1}^{\infty}a_iF_r(X)^i+\varphi_{r+1}(X)\Big(\sum_{i=1}^{\infty}
a_i\sum_{j=1}^i\binom{i}{j}F_r(X)^{i-j}\varphi_{r+1}(X)^{j-1}\Big)\\
&=f(F_r(X))+\varphi_{r+1}(X)\Big(\sum_{i=1}^{\infty} a_i iF_r(X)^{i-1}
+\varphi_{r+1}(T)l(T)\Big).
\end{align*}

Por tanto
\[
f(F_{r+1}(X))\equivalente_{\substack{\uparrow\\ a_1=\pi}} f(F_r(X))+
\pi \varphi_{r+1}(X)\bmod \deg (r+2).
\]
Adem\'as
$F_{r+1}(g(X))=F_r(g(X))+\varphi_{r+1}(g(X))$. Puesto que $g(X)\in
{\mc F}_{\pi}$, se tiene 
\begin{align*}
\varphi_{r+1}(g(x))&=\varphi_{r+1}(g(X_1),\ldots,g(X_n))
\igual\limits_{\substack{\uparrow\\ g(X_i)=\pi X_i+\cdots}}\\
&=\pi^{r+1} \varphi_{r+1}(X)+\text{t\'erminos de grado $\geq r+2$}.
\end{align*}
 Por 
tanto 
\begin{gather*}
F_{r+1}(g(X))\equiv F_r(g(X))+\pi^{r+1}\varphi_{r+1}(X)\bmod \deg (r+2).
\intertext{De las congruencias}
f(F_{r+1}(X))\equiv f(F_r(X))+\pi \varphi_{r+1}(X)\bmod \deg (r+2),\\
F_{r+1}(g(X))\equiv F_r(g(X))+\pi^{r+1}\varphi_{r+1}(X)\bmod \deg (r+2),
\intertext{se sigue que si queremos}
f(F_{r+1}(X))\equiv F_{r+1}(g(X))\bmod\deg (r+2),
\intertext{entonces}
f(F_r(X))+\pi\varphi(X)\equiv F_r(g(X))+\pi^{r+1}\varphi_{r+1}(X)
\bmod\deg (r+2),
\intertext{lo cual equivale a}
\varphi_{r+1}(X)\equiv\frac{f(F_r(X))-F_r(g(X))}{\pi^{r+1}-\pi}\bmod\deg(r+2).
\end{gather*}

De esta forma, $\varphi_{r+1}(X)$ es \'unico pues es la serie
$\frac{f(F_r(X))-F_r(g(X))}{\pi^{r+1}-\pi}$ quitando los t\'erminos de grado
mayor o igual a $(r+2)$ y la serie $f(F_r(X))-F_r(g(X))$ no tiene t\'erminos
de grado diferente a $(r+1)$ por la Ecuaci\'on (\ref{CClaseEc3.2.5.1}).
Ahora bien
\[
f(F_r(X))-F_r(g(X))\equiv F_r(X)^q-F_r(X^q)\equiv 0\bmod \pi.
\]
Por tanto $\varphi_{r+1}(X)\in \o_K[X]$ y $\varphi_{r+1}(X)$
es homog\'eneo de grado $(r+1)$. Esto muestra la existencia y la
unicidad de $F(X)=\lim_{r\to\infty}F_r(X)$. $\fin$
\end{proof}

\begin{observacion}\label{CClaseO3.2.5.5}
La demostraci\'on del Teorema \ref{CClaseT3.2.5.4} de hecho prueba que $F$
es la \'unica serie de potencias en cualquier campo que contenga a
$\o_K$ y que satisface las condiciones del Teorema \ref{CClaseT3.2.5.4}.
\end{observacion}

Definimos en general {\em grupos formales\index{grupos formales}}
y {\em m\'odulos formales\index{m\'odulos formales}}.

Sea $R$ es un anillo conmutativo
con unidad. Sea $\mK=\langle X\rangle=XR[[X]]$ el ideal generado
por $X$. Entonces $\mK=\{f(X)\in R[[X]]\mid f\equiv 0\bmod \deg 1\}$.
Se tiene para $f,g\in\mK$, $(f\circ g)(X)=f(g(X))\in \mK$. En particular
$\mK$ es un semigrupo con la operaci\'on $\circ$. Ahora bien,
para $f\in\mK$ se tiene
$X\circ f=f\circ X=f$, por lo que la identidad del semigrupo es la serie $X$.

Si para $f,g\in\mK$ satisfacen $f\circ g=g\circ f=X$, escribimos $f=g^{-1}$
y $g=f^{-1}$ y se dice que $f$ es {\em invertible\index{elemento invertible
en las series $XR[[X]]$}}.

\begin{proposicion}\label{CClaseP3.2.5.5(2)} Se tiene que $f\in\mK$ es invertible
si y solamente si $a_1\in\* R$ donde $f(X)=a_1X+\sum_{i=2}^{\infty}
a_iX^i$, esto es, si existe $c\in R$ tal que $a_1c=1$.
\end{proposicion}

\begin{proof}
Si $g=f^{-1}$, $(f\circ g)(X)=a_1 g(X)+\sum_{j=2}^{\infty} d_jX^j$ con
$g(X)=cX+\sum_{j=2}^{\infty}c_jX^j$. Por tanto
\[
f(g(X))=a_1 cX+\sum_{j=2}^{\infty} e_jX^j=X,
\]
de donde obtenemos que $a_1c=1$.

Rec\'iprocamente, sea $a_1\in\*R$ y $c_1=c\in R$ con $a_1c_1=1$.
Buscamos $g(X)=\sum_{i=1}^{\infty}c_iX^i$ con $f(g(X))=X$. Se tiene
$f(g(X))=\sum_{i=1}^{\infty}a_i(g(X))^i$, la cual podemos resolverla de
manera recursivo satisfaciendo $f\circ g=X$.

Por la misma raz\'on, existe $h$ tal que $h\circ f=X$. Por tanto
\[
h=h\circ X=h\circ f\circ g=X\circ g=g.
\]
Obtenemos que $f\circ g=g\circ f=X$. $\fin$
\end{proof}

\begin{observacion}\label{CClaseO3.2.5.5(3)}
Si $f\circ g=g\circ f=X$ y $f\circ h=h\circ f=X$, entonces
$h=h\circ X=h\circ(g\circ f)=(h\circ f)\circ g=X\circ g=g$.
\end{observacion}

\begin{definicion}\label{CClaseD3.2.2.5'} Sea $R$ un anillo conmutativo
con unidad. Una serie formal $F(X,Y)\in R[[X,Y]]$ se llama un
{\em grupo formal\index{grupo formal}} sobre $R$ si
\lasa
\item $F(X,Y)\equiv X+Y\bmod \deg 2$.
\item $F(F(X,Y),Z)=F(X,F(Y,Z))$ (asociatividad).
\item $F(X,Y)=F(Y,X)$ (conmutatividad).
\end{list}

En particular se tiene $F(0,0)=0$ lo cual implica que las dos series
en (b) est\'an bien definidas (pues de otra forma podr\'iamos tener
una serie de constantes en $R$ no convergente).
\end{definicion}

\begin{ejemplos}\label{CClaseE3.2.5.5(2')}
\las
\item $F(X,Y)=X+Y$ es un grupo formal llamado el {\em grupo formal
aditivo\index{grupo formal aditivo}} y es denotado por ${\ma G}_a$.

\item $F(X,Y)=X+Y+XY=(1+X)(1+Y)-1$ es un grupo formal llamado
{\em grupo formal multiplicativo\index{grupo formal multiplicativo}}
y es denotado por ${\ma G}_m$.
\end{list}
\end{ejemplos}

De la definici\'on de grupo formal, obtenemos $F(X,0)\equiv X\bmod\deg 2$,
$F(F(X,0),0)=F(X,0)$. De la condici\'on $F(X,0)
\equiv X\bmod \deg 2$ obtenemos
que $f(X):=F(X,0)$ es invertible, $f^{-1}\in\mK=XR[[X]]$.

Se tiene
\[
f(X)=F(X,0)=F(F(X,0),0)=F(f(X),0)=(f\circ f)(X).
\]
Se sigue que $X=(f^{-1}\circ f)(X)=(f^{-1}\circ f\circ f)(X)=f(X)$,
esto es, $f(X)=X$. As\'i, $F(X,0)=X$. Similarmente se obtiene que
$F(0,Y)=Y$. Por tanto
\[
F(X,Y)=X+Y+\sum_{i,j=1}^{\infty} c_{i,j}X^iY^j,
\]
es decir, no existen t\'erminos de la forma $X^i$ o $Y^j$ con 
$i,j\geq 2$.

La ecuaci\'on $F(X,Y)=0$ puede resolverse para $Y$ en $\mK$, esto es,
existe una \'unica serie $i_F(X)=-X+\sum_{i=2}^{\infty}b_iX^i$ con 
$b_i\in R$ tal que $F(X,i_F(X))=F(i_F(X),X)=0$ (esto lo podemos hacer
por recursi\'on).

Notemos que $i_F(X)\equiv -X\bmod \deg 2$. La serie $i_F(X)$
recibe el nombre de {\em inversa formal\index{inversa formal}}.

\begin{ejemplos}\label{CClaseEj3.2.5.5(2')+1}
\las
\item La inversa formal de ${\ma G}_a$ es $i_F(X)=-X$.
\item La inversa formal de ${\ma G}_m$ es $i_F(X)=
(1+X)^{-1}-1=-X+X^2-\cdots =\sum_{i=1}^{\infty}(-1)^i X^i$.
\end{list}
\end{ejemplos}

\begin{definicion}\label{CClaseD3.2.5.5(2')+2} Sea $F(X,Y)\in R[[X,Y]]$
un grupo formal.
Para $f,g\in\mK$, definimos
\begin{gather*}
f\lubintatemas F g:=F(f(X),g(X)).
\end{gather*}
Entonces $f\lubintatemas F g\in \mK$ y $\mK$ es un grupo abeliano con
respecto a esta operaci\'on y el inverso de $f$ es $i_F(f)$.
Este grupo abeliano se denota por $\mK_F$.
\end{definicion}

\begin{ejemplos}\label{CClaseO3.2.5.5(2')+3}
\las
\item $\mK_{{\ma G}_a}=\mK$ con la suma.
\item $\mK_{{\ma G}_m}\cong1+\mK$ con la multiplicaci\'on.
\end{list}
\end{ejemplos}

Sea $G(X,Y)$ otro grupo formal sobre $R$ y sea $f\in\mK$ tal que
\begin{gather}\label{CClaseEc3.2.5.5(5)}
f(F(X,Y))=G(f(X),f(Y))
\end{gather}

\begin{definicion}\label{CClaseD3.2.5.5(4)} Si $f$ satisface la Ecuaci\'on
(\ref{CClaseEc3.2.5.5(5)}), $f$ se llama {\em morfismo\index{morfismo de
grupos formales}} de $F$ en $G$ y se escribe
\[
f\colon F\lra G.
\]
Si $F=G$, $f$ se llama {\em endomorfismo\index{endomorfismo de
grupos formales}}.

Si $f$ tiene inversa $f^{-1}\in\mK$, $f^{-1}\colon G\lra F$ es un morfismo
de $G$ en $F$ y $f$ se llama {\em isomorfismo\index{isomorfismo
de grupos formales}}, $f\colon F\xrightarrow{\ \cong\ }G$.
\end{definicion}

La Ecuaci\'on (\ref{CClaseEc3.2.5.5(5)}) se escribe
\begin{gather}\label{CClaseEq3.8'}
f\circ F=G\circ f.
\end{gather}
N\'otese la similitud con morfismos de m\'odulos de Drinfeld.

\begin{ejemplo}\label{CClaseEj3.2.5.5(4')}
Sea $F$ es un grupo formal. Se define el endomorfismo de $F$,
$[m]\colon F\lra F$
definido por $[0](X)=0$, $[m+1](X)=F([m]X,X)$ si $m\geq 0$ y
$[m-1](X)=F([m](X),i_F(X))$ si $m\leq 0$.

Este endomorfismo se llama {\em multiplicaci\'on por $m$\index{endomorfismo
multiplicaci\'on por $m$ de grupos formales}}. Si $m\in U_K$, entonces
$[m]$ es un isomorfismo.
\end{ejemplo}

\begin{definicion}\label{CClaseD3.2.5.5(6)} Si $F(X_1,\ldots,X_m)\in R[[X_1,\ldots,X_m]]$
y si $f\in\mK=XR[[X]]$ es invertible en $\mK$: $f^{-1}\in \mK$, se define la
serie $F^f(X_1,\ldots,X_m)\in R[[X_1,\ldots,X_m]]$ por
\begin{gather}\label{CClaseEc2.2.5.5(6)+1}
F^f(X_1,\ldots,X_m)=f\circ F\circ f^{-1}=f(F(f^{-1}(X_1),\ldots,f^{-1}(X_m))).
\end{gather}
\end{definicion}

Si $F(X,Y)$ es un grupo formal sobre $R$, entonces $G=F^f$ es nuevamente
un grupo formal y $f\colon F\xrightarrow{\ \cong\ }G$ es un isomorfismo.

Sea
\begin{gather*}
\Hom_R(F,G)=\Hom(F,G)=\{f\mid f\colon F\lra G\text{\ es un 
morfismo}\}\index{grupo de homomorfismos de grupos formales},\\
\End_R(F)=\End(F)=\Hom_R(F,F)\index{grupo de endomorfismos de
grupos formales}.
\end{gather*}

\begin{proposicion}\label{CClaseP3.2.5.5(7)} Se tiene $\Hom(F,G)$ es un
subgrupo de $\mK_G$. Adem\'as $\End(F)$ es un anillo con respecto
a la suma $f\lubintatemas F g$ y multiplicaci\'on $f\circ g$ cuya 
identidad es $X$.
\end{proposicion}

\begin{proof} \cite[Lemma 4.1, p\'agina 51]{Iwa86}. $\fin$
\end{proof}

Para definir m\'odulos formales, 
aplicamos el Teorema \ref{CClaseT3.2.5.4} a
los casos suma y multiplicaci\'on por 
escalar: $L(X,Y)=X+Y$ y $L(X)=aX$ con $a\in\o_K$.

Sea $f\in{\mc F}_{\pi}$ y sea $F_f(X,Y)$ la {\underline{\'unica}} soluci\'on
de
\begin{gather}
F_f(X,Y)\equiv X+Y\bmod \deg 2,\nonumber\\
f(F_f(X,Y))=F_f(f(X),f(Y)) \label{CClaseEc3.5.2.1''}
\end{gather}

Notemos que la Ecuaci\'on (\ref{CClaseEc3.5.2.1''}) nos dice que $f$ es
un endomorfismo de $F_f$: $f\circ F_f=F_f\circ f$ (ver 
Ecuaci\'on (\ref{CClaseEq3.8'})).

Para cada $a\in\o_K$ y $f,g\in {\mc F}_{\pi}$, sea la serie $a_{f,g}(T)
\in \o_K[[T]]$ la \'unica soluci\'on de
\begin{gather}
a_{f,g}(T)\equiv aT\bmod\deg 2,\nonumber\\
f(a_{f,g}(T))=a_{f,g}(g(T))\label{CClaseEc3.5.2.1'''}.
\end{gather}
Por notaci\'on, escribiremos $a_f=a_{f,f}$.

El siguiente teorema probar\'a, entre otras cosas, 
que $F_f$ es un grupo formal.

\begin{teorema}\label{CClaseT3.2.5.6}
Sean $f,g,h\in{\mc F}_{\pi}$ y $a,b\in\o_K$. Entonces
\las
\item $F_f(X,Y)=F_f(Y,X)$.
\item $F_f(F_f(X,Y),Z)=F_f(X,F_f(Y,Z))$.
\item $a_{f,g}(F_g(X,Y))=F_f(a_{f,g}(X),a_{f,g}(Y))$.
\item $a_{f,g}(b_{g,h}(Z))=(a\cdot b)_{f,h}(Z)$.
\item $(a+b)_{f,g}(Z)=F_f(a_{f,g}(Z),b_{f,g}(Z))$.
\item $(\pi^n)_f(Z)=f^{(n)}(Z)$, $n=0,1,2,\ldots $.
\end{list}
\end{teorema}

\begin{proof} Todos estas propiedades se siguen de resolver alg\'un problema
espec\'ifico aplicando el Teorema \ref{CClaseT3.2.5.4}.

\las
\item Sea $F(X,Y)=F_f(X,Y)$. 
Se tiene que 
\[
F_f(X,Y)\equiv X+Y\bmod\deg 2\equiv 
Y+X\bmod \deg 2\equiv F_f(Y,X).
\]
Adem\'as se tiene
que $f(F_f(X,Y))=F_f(f(X),f(Y))$,
por tanto 
\[
f(F_f(Y,X))\igual\limits_{\substack{\uparrow\\ X\leftrightarrow Y}}
F_f(f(Y),f(X)).
\]

Sea $G(X,Y)=F_f(Y,X)$. Entonces
\[
f(G(X,Y))=f(F_f(Y,X))=F_f(f(Y),f(X))=G(f(X),f(Y)).
\]
Por unicidad, se tiene $F=G$, esto es, $F_f(X,Y)=G(X,Y)=F_f(Y,X)$.

\item $F_f(F_f(X,Y),Z)\equiv F_f(X,Y)+Z\bmod\deg 2\equiv
X+Y+Z\bmod\deg 2\equiv X+F_f(Y,Z)\bmod\deg 2=F_f(X,F_f(Y,Z))
\bmod\deg 2$.

Sean $F_1(X,Y,Z)=F_f(F_f(X,Y),Z)$, $F_2(X,Y,Z)=F_f(X,F_f(Y,Z))$. Entonces
\begin{align*}
f(F_1(X,Y,Z))&=f(F_f(F_f(X,Y),Z))=F_f(f(F_f(X,Y),f(Z)))\\
&=F_f(F_f(f(X),f(Y)),f(Z))\\
&=F_f(F_f(f(X),f(Y)),f(Z))\\
&=F_1(f(X),f(Y),f(Z)),\\
f(F_2(X,Y,Z))&=f(F_f(X,F_f(Y,Z)))\\
&=F_f(f(X),f(F_f(Y,Z)))\\
&=F_f(f(X),F_f(f(Y),f(Z)))=F_2(f(X),f(Y),f(Z)).
\end{align*}

Por tanto
\begin{gather*}
F_1(X,Y,Z)\equiv X+Y+Z\bmod\deg 2\equiv F_2(X,Y,Z),\\
f\circ F_1=F_1\circ f\quad\text{y}\quad f\circ F_2=F_2\circ f.
\end{gather*}
Por unicidad tenemos que $F_1=F_2$.

\item Sea $H(X,Y)$ la \'unica soluci\'on de
\begin{gather}
H(X,Y)\equiv aX+aY\bmod \deg 2,\nonumber\\
f(H(X,Y))=H(g(X),g(Y)).\label{CClaseEc3.5.2.2}
\end{gather}

Sea $H_1(X,Y)=F_f(a_{f,g}(X),a_{f,g}(Y))$. Entonces
\begin{align*}
H_1(X,Y)&=F_f(a_{f,g}(X),a_{f,g}(Y))\equiv a_{f,g}(X)+a_{f,g}(Y)\\
&\equiv aX+aY\bmod\deg 2,\\
f(H_1(X,Y))&=f(F_f(a_{f,g}(X),a_{f,g}(Y)))=F_f(f(a_{f,g}(X)),f(a_{f,g}(Y)))\\
&=F_f(a_{f,g}(g(X)),a_{f,g}(g(Y)))=H_1(g(X),g(Y)).
\end{align*}

Se sigue que $H_1(X,Y,Z)$ es soluci\'on de la ecuaci\'on (\ref{CClaseEc3.5.2.2}).

Sea $H_2(X,Y)=a_{f,g}(F_g(X,Y))$. Entonces
\begin{align*}
H_2(X,Y)&=a_{f,g}(F_g(X,Y))\equiv aF_g(X,Y)\equiv a(X+Y)\bmod\deg 2,\\
f(H_2(X,Y))&=f(a_{f,g}(F(X,Y)))=a_{f,g}(g(F_g(X,Y)))\\
&=a_{f,g}(F_g(g(X),g(Y)))=H_2(g(X),g(Y)).
\end{align*}

Se sigue que $H_2(X,Y,Z)$ tambi\'en es soluci\'on de la ecuaci\'on
(\ref{CClaseEc3.5.2.2}). Por la unicidad dada en el Teorema \ref{CClaseT3.2.5.4}, se
sigue que $H_1(X,Y)=H_2(X,Y)$.

\item Sea $H(Z)$ la soluci\'on a 
\begin{gather}\label{CClaseEc3.5.2.3}
H(Z)\equiv ab Z\bmod\deg 2\quad\text{y}\quad F(H(Z))=H(h(Z)).
\end{gather}

Sea $H_1(Z)=a_{f,g}(b_{g,h}(Z))$. Entonces
\begin{align*}
H_1(Z)&=a_{f,g}(b_{g,h}(Z))\equiv a b_{g,h}(Z)\equiv ab Z\bmod\deg 2,\\
f(H_1(Z))&=f(a_{f,g}(b_{g,h}(Z)))=a_{f,g}(g(b_{g,h}(Z)))\\
&=a_{f,g}(b_{g,h}(h(Z)))=H_1(h(Z)).
\end{align*}

Sea ahora $H_2(Z)=(ab)_{f,h}(Z)$. Entonces
\begin{align*}
H_2(Z)&=(ab)_{f,h}(Z)\equiv ab Z\bmod \deg 2,\\
f(H_2(Z))&=f((ab)_{f,h}(Z))=(ab)_{f,h}(h(T))=H_2(h(T)).
\end{align*}

Por tanto $H_1$ y $H_2$ son soluciones de la ecuaci\'on (\ref{CClaseEc3.5.2.3}),
de donde se sigue que $H_1(Z)=H_2(Z)=H(Z)$.

\item Sea $H(Z)$ soluci\'on a
\begin{gather}\label{CClaseEc3.5.2.4}
H(Z)\equiv(a+b) Z\bmod\deg 2\quad\text{y}\quad f(H(Z))=H(g(Z)).
\end{gather}

Sea $H_1(Z)=(a+b)_{f,g}(Z)$. Entonces
\begin{align*}
H_1(Z)&=(a+b)_{f,g}(Z)\equiv (a+b)Z\bmod \deg 2,\\
f(H_1(Z))&=f((a+b)_{f,g}(Z))=(a+b)_{f,g}(g(Z))=H_1(g(Z)).
\end{align*}

Sea ahora $H_2(Z)=F_f(a_{f,g}(Z),b_{f,g}(Z))$.  Entonces
\begin{align*}
H_2(Z)&=F_f(a_{f,g}(Z), b_{f,g}(Z))\equiv a_{f,g}(Z)+b_{f,g}(Z)\bmod\deg 2\\
&\equiv aZ+bZ\bmod\deg 2,\\
f(H_2(Z))&=f(F_f(a_{f,g}(Z),b_{f,g}(Z)))=F_f(f(a_{f,g}(Z)),f(b_{f,g}(Z)))\\
&=F_f(a_{f,g}(g(Z)),b_{f,g}(g(Z)))=H_2(g(Z)).
\end{align*}

Se sigue que $H_1(Z)=H_2(Z)$.

\item Sea $n\geq 0$ y sea $H^{(n)}(Z)$ la soluci\'on a
\begin{gather*}
H^{(n)}(Z)\equiv \pi^n Z\bmod\deg 2 \quad \text{y}\quad
f(H^{(n)}(Z))=H^{(n)}(f(Z)).
\end{gather*}

Sea $H_1^{(n)}(Z)=(\pi^n)_f(Z)$. Entonces
\begin{align*}
H_1^{(n)}(Z)&=(\pi^n)_f(Z)\equiv \pi^n Z\bmod\deg 2,\\
f(H_1^{(n)}(Z))&=f((\pi^n)_f(Z))=(\pi^n)_f(f(Z))=H_1^{(n)}(f(Z)).
\end{align*}

Sea ahora $H_2^{(n)}(Z)=f^{(n)}(Z)$. Entonces
\begin{align*}
H_2^{(n)}(Z)&=f^{(n)}(Z)=f(f^{(n-1)}(Z))\equiv \pi f^{(n-1)}(Z)\bmod\deg 2\\
&\equiv\pi \pi^{n-1} Z\bmod\deg 2,\\
f(H_2^{(n)}(Z))&=f(f^{(n)}(Z))=f^{(n)}(f(Z))=H^{(n)}(f(Z)).
\end{align*}

Por tanto se tiene que $H_1(Z)=H_2(Z)$.
$\fin$
\end{list}
\end{proof}

Sea $L$ una extensi\'on algebraica del campo local $K$. Sea 
$\pK_L=\{x\in L\mid v_L(x)>0\}$. Si $x_1,\ldots,x_n\in\pK_L$ y
$G(X_1,\ldots,X_n)\in\o_K[[X_1,\ldots,X_n]]$, entonces $G(x_1,\ldots,
x_n)$ converge en $K(x_1,\ldots,x_n)$. Si adem\'as el t\'ermino
constante de $G$ es $0$, $G(x_1,\ldots,x_n)\in\pK_L$.

Una demostraci\'on de lo anterior, es como sigue. Como $K$ es completo,
$K[x_1,\ldots,x_n]=k(x_1,\ldots,x_n)$ tambi\'en es completo y si $G_d$
denota a $G$ quitando los t\'erminos de grado total $\geq d+1$, entonces
\begin{align*}
|G(x_1,\ldots,x_n)-G_d(x_1,\ldots,x_n)|&=\Big|\sum_i a_{i_1,\ldots,i_n}
x_1^{i_1}\cdots x_n^{i_n}\Big|\\
&\leq |x_1^{i_1}|\cdots |x_n^{i_n}|\leq
c^{i_1+\cdots+i_n}\leq c^{d+1}\xrightarrow[d\to\infty]{} 0,
\end{align*}
donde $|x_1|,\ldots,|x_n|\leq c<1$.

Entonces $\lim_{d\to\infty}G_d(x_1,\ldots,x_n)=G(x_1,\ldots,x_n)$ pues
$\{G_d(x_1,\ldots,x_n)\}_d$ es una sucesi\'on de Cauchy. Si $a_{0,\ldots,0}
=0$, claramente $|G_d(x_1,\ldots,x_n)|<1$ y $G(x_1,\ldots,x_n)\in
\pK_L$.

Por esta raz\'on, para $f=g=h$ las propiedades del Teorema \ref{CClaseT3.2.5.6},
hacen de $F_f$ un {\em $\o_K$--m\'odulo (de Lie) formal\index{m\'odulo
formal}}.

Si en el Teorema \ref{CClaseT3.2.5.6} tomamos $f=g=h$ y pensamos
$F_f(x,y)$ como la suma y $a_{f,f}$, $a\in\o_K$ como la multiplicaci\'on
por escalar, entonces tendr\'iamos un $\o_K$--m\'odulo haciendo que las
variables $X,Y$ y $Z$ tomen valores de un dominio donde las series
convergen (por ejemplo $L$, una extensi\'on algebraica de $K$).

\begin{proposicion}\label{CClaseP3.2.5.7} Si $f\in {\mc F}_{\pi}$ y $L$ es una 
extensi\'on algebraica de $K$, $\pK_L$ es un $\o_K$--m\'odulo con la
suma y multiplicaci\'on por escalar definidas por:
\begin{gather*}
x\lubintatemas {F_f} y:=F_f(x,y)\quad\text{y}\quad a \lubintatepor 
{F_f} x=a_f(x),\quad
x,y\in\pK_L, a\in\o_K.
\end{gather*}

Se denota $\pK_L^{(f)}$ a este $\o_K$--m\'odulo.
\end{proposicion}

\begin{proof} Es consecuencia inmediata de Teorema \ref{CClaseT3.2.5.6}. $\fin$
\end{proof}

\begin{observacion}\label{CClaseO3.3.20'} Las propiedades (1) y (2) del
Teorema \ref{CClaseT3.2.5.6} prueban que $\pK_L^{(f)}$ es un grupo
aditivo. Las propiedades (3), (4) y (5) con $f=g=h$, prueban que
$\pK_L^{(f)}$ es un $\o_K$--m\'odulo.
\end{observacion}

El inverso aditivo de $x$ es $(-1)_f(x)=(-1)\lubintatepor {F_f} x$.
Es importante distinguir $\pK_L^{(f)}$
de $\pK_L$: ambos son el mismo conjunto y ambos son $\o_K$--m\'odulos
pero con diferente acci\'on.

\begin{teorema}\label{CClaseT3.2.5.8} El conjunto de ceros $\Lambda_{f,n}$ de
$f^{(n)}(x)$ es un $\o_K$--subm\'odulo de $\pK_{L_{f,n}}^{(f)}$, $n\geq 1$,
donde recordemos que $L_{f,n}=K(\Lambda_{f,n})$.

M\'as a\'un, $\Lambda_{f,n}$ son los punto de $\pi^n$--torsi\'on del
$\o_K$--m\'odulo $\pK^{(f)}_{L_{f,n}}$. Esto es
\[
\Lambda_{f,n}=\{\lambda\in \pK_{L_{f,n}}\mid \pi^n\lubintatepor {F_f}
\lambda=0\}=\{\lambda\in\bar{K}\mid \pi^n\lubintatepor {F_f}\
\lambda=0\}.
\]
\end{teorema}

\begin{proof} Por definici\'on, tenemos que si $x\in \Lambda_{f,n}$, entonces
$x\in \pK_{L_{f,n}}$. Se tiene
\begin{align*}
\Lambda_{f,n}&=\{\lambda\in\pK_{L_{f,n}}\mid f^{(n)}(\lambda)=(\pi^n)_f(\lambda)=0\}\\
&=\{\lambda\in\pK_{L_{f,n}}\mid \pi^n\lubintatepor {F_f}
\lambda=0\}=\ker \pi^n.
\end{align*}
Por tanto $\Lambda_{f,n}$ es un $\o_K$--m\'odulo. $\fin$
\end{proof}

Nuevamente hacemos notar la similitud de $\Lambda_{f,n}$ como $\o_K$--m\'odulo
y de $\Lambda_M$, $M\in R_T=\F[T]$ como $R_T$--m\'odulo (Carlitz).

\begin{ejemplo}\label{CClaseE3.3.21'}
Sean $K={\ma Q}_p$ y $f(T)=(T+1)^p-1\in{\mc F}_p$. Entonces
\[
\Lambda_{f,n}=\{\lambda\in{\ma Q}_p\mid p^n\lubintatepor {F_f}\lambda
=0\}=\{\lambda\in\bar{\ma Q}_p\mid f^{(n)}(\lambda)=0\}.
\]

Se tiene $f^{(2)}(T)=f(f(T))=(f(T)+1)^p-1=((T+1)^p-1+1)^p-1=(T+1)^{p^2}-1$
y en general $f^{(n)}(T)=(T+1)^{p^n}-1$. Se sigue que 
\[
\Lambda_{f,n}=\{\lambda\in\bar{\ma Q}_p\mid (\lambda+1)^{p^n}-1=0\}=
\{\zeta_{p^n}^i-1\}_{i=0}^{p^n-1}\cong \langle \zeta_{p^n}\rangle
\cong W_{p^n}.
\]
\end{ejemplo}

\begin{proposicion}\label{CClaseP3.2.5.9} Sean $f,g\in{\mc F}_{\pi}$ y $a\in\o_K$.
El mapeo $\lambda\longmapsto a_{g,f}(\lambda)$ da lugar a un homomorfismo
de $\o_K$--m\'odulos de $\Lambda_{f,n}$ en $\Lambda_{g,n}$. Este homomorfismo
es un isomorfismo si $a\in U_K$. M\'as precisamente, se tiene 
$a_{g,f}\colon F_f\lra F_g$.
\end{proposicion}

\begin{proof} Por el Teorema \ref{CClaseT3.2.5.6} (3) y (4), se tiene
\begin{align*}
\lambda\lubintatemas {F_f} \mu=F_f(\lambda,\mu)\longmapsto a_{g,f}(F_f(\lambda,
\mu))&=F_g(a_{g,f}(\lambda),a_{g,f}(\mu))\\
&=a_{g,f}(\lambda)
\lubintatemas {G_f} a_{g,f}(\mu),
\intertext{y}
b\lubintatepor {F_f} \lambda=b_f(\lambda)\longmapsto a_{g,f}(b_f(\lambda))&=(ab)_{g,f}
(\lambda)=(ba)_{g,f}(\lambda)\\
&=b_{g,g}(a_{g,f}(\lambda))=b\lubintatepor {G_f} a_{g,f}
(\lambda).
\end{align*}

Por tanto el mapeo es un homomorfismo de $\o_K$--m\'odulos. Adem\'as, si
$\lambda\in\Lambda_{f,n}$, 
\begin{align*}
g^{(n)}(a_{g,f}(\lambda))&=(\pi^n)_g(a_{g,f}(\lambda))
=(\pi^n a)_{g,f}(\lambda)=(a\pi^n)_{g,f}(\lambda)\\
&=a_{g,f}((\pi^n)_f (\lambda))=
a_{g,f}(f^{(n)}(\lambda))=a_{g,f}(0)=0.
\end{align*}
Por tanto el homomorfismo manda $\Lambda_{f,n}$ en $\Lambda_{g,n}$.

Ahora si, $a\in U_K$, para $\lambda\in\Lambda_{f,n}$ se tiene
\begin{gather*}
(a^{-1})_{f,g}(a_{g,f}(\lambda))=(a^{-1}a)_{f,f}(\lambda)=1_f(\lambda)
=1\lubintatepor {F_f} \lambda=\lambda
\intertext{y para $\mu\in\Lambda_{g,n}$ se tiene}
(a_{g,f})((a^{-1})_{f,g})(\mu)=(aa^{-1})_{g,g}(\mu)=1_g(\mu)=
1\lubintatepor {G_f} \mu=\mu.
\end{gather*}

Se sigue que el homomorfismo $a_{g,f}\colon \Lambda_{f,n}\lra \Lambda_{g,n}$
tiene como inverso a $(a^{-1})_{f,g}\colon \Lambda_{g,n}\lra \Lambda_{f,n}$. $\fin$
\end{proof}

\begin{teorema}\label{CClaseT3.2.5.10} Para cualquier $f\in{\mc F}_{\pi}$, se tiene
\[
\Lambda_{f,n}\cong \o_K/\pi^n\o_K
\]
como $\o_K$--m\'odulos.
\end{teorema}

\begin{proof} Para cualesquiera $f,g\in{\mc F}_{\pi}$, $1_{g,f}\colon \Lambda_{f,n}
\lra \Lambda_{g,n}$ es un isomorfismo de $\o_K$--m\'odulos. Basta considerar
$f(Z)=\pi Z+Z^q\in{\mc F}_{\pi}$. 

Se har\'a por inducci\'on en $n$. Para
$n=1$, $\Lambda_{f,1}$  es el conjunto de ceros de $f(\lambda)=\pi\lambda
+\lambda^q=0$ el cual es un polinomio separable pues $f'(\lambda)=\pi\neq 0$.
Esto es, $\Lambda_{f,1}$ tiene $q$ elementos y por tanto es un espacio 
vectorial sobre $\F$ de dimensi\'on $1$. En este caso $\F\cong\o_K/\pi\o_K$,
de donde se sigue que $\Lambda_{f,1}\cong \o_K/\pi\o_K$.

Supongamos que $\Lambda_{f,n}\cong \o_K/\pi^n\o_K$ como $\o_K$--m\'odulos.
Consideremos $\pi_f\colon \Lambda_{f,n+1}\lra \Lambda_{f,n}$ dada por
$\pi_f(\lambda)=\pi\lubintatepor {F_f}
\lambda\in \Lambda_{f,n}$ para $\lambda
\in \Lambda_{f,n+1}$ pues $f^{(n)}(\pi_f(\lambda))=
f^{(n)}(f(\lambda))=f^{(n+1)}(\lambda)=0$.

Puesto que $\ker\pi_f=\Lambda_{f,1}$ y $|\Lambda_{f,n+1}|=q^{n+1}$,
$|\Lambda_{f,n}|=q^n$, se tiene la sucesi\'on exacta
\[
0\lra \Lambda_{f,1}\lra \Lambda_{f,n+1}\xrightarrow{\ \pi_f\ }\Lambda_{f,n}
\lra 0.
\]

Otra forma de verificar la suprayectividad de $\pi_f$ es como sigue.
Si $\lambda\in\Lambda_{f,n}$ y $\mu$ es una ra\'iz de $f(Z)-\lambda=
Z^q+\pi Z-\lambda$, entonces $\lambda=f(\mu)$ y 
$f^{(n+1)}(\mu)=f^{(n)}(f(\mu))
=0$ y por tanto $\pi_f(\mu)=f(\mu)=\lambda$ y $\pi_f$ es suprayectiva.

Sea $\lambda\in \Lambda_{f,n+1}\setminus\Lambda_{f,n}$, entonces 
$(\pi^n)_f(\lambda)\neq 0$ y $(\pi^{n+1})_f(\lambda)=0$ por lo que el anulador
de $\lambda$ es $\pi^{n+1}\o_K$. El mapeo 
$a\longmapsto a\lubintatepor {F_f}
\lambda$ da el isomorfismo.
\[
\o_K\lambda\cong\o_K/\pi^{n+1}\o_K
\]
y $\o_K\lambda\subseteq \Lambda_{f,n+1}$. Puesto que $\big|\o_K/\pi^{n+1}\o_K
\big|=|\Lambda_{f,n+1}|=q^{n+1}$, se sigue que $\o_K\lambda=\Lambda_{f,n+1}
\cong \o_K/\pi^{n+1}\o_K$. $\fin$
\end{proof}

\begin{teorema}\label{CClaseT3.2.5.11} Todo automorfismo del $\o_K$--m\'odulo
$\Lambda_{f,n}$ es de la forma $u_f\colon \Lambda_{f,n}\lra \Lambda_{f,n}$
con $u\in U_K$. Se tiene que $u_f=\Id_{\Lambda_{f,n}}$ si y s\'olo si
$u\in \unidades n$. Por tanto
\[
\Aut_{\o_K}(\Lambda_{f,n})\cong U_K/\unidades n.
\]
\end{teorema}

\begin{proof} Si $\sigma\in\Aut_{\o_K}(\Lambda_{f,n})$ se tiene que, usando el isomorfismo
de $\o_K$--m\'odulos $\Lambda_{f,n}\cong\o_K/\pi^n\o_K$, $\sigma\in
\Aut_{\o_K}(\o_K/\pi^n \o_K)$. Sea $\tau\colon \o_K\lra \o_K/\pi^n\o_K$ el 
epimorfismo natural. Por tanto $\sigma\circ \tau\colon\o_K\lra\o_K/\pi^n\o_K$
es un epimorfismo de $\o_K$--m\'odulos y para $\xi\in\o_K$ se tiene
\[
(\sigma\circ \tau)(\xi)=\xi(\sigma\circ \tau)(1)=\xi\sigma(1\bmod\pi^n).
\]
Se sigue que $\sigma(1\bmod \pi^n)=a$ genera a $\o_K/\pi^n\o_K$ y por
tanto $a$ es unidad. 

Se tiene el siguiente diagrama conmutativo
\[
\xymatrix{
\o_K/\pi^n\o_K\ar[r]^{a}\ar[d]_{\cong}^{\varphi}&
\o_K/\pi^n\o_K\ar[d]_{\varphi}^{\cong}\\
\Lambda_{f,n}\ar[r]^{\sigma}&\Lambda_{f,n}
}
\]
donde $\varphi(1\bmod \pi^n)=\lambda$ es un generador $\Lambda_{f,n}$.
Entonces $\sigma(\lambda)=(\varphi a\varphi^{-1})
(\lambda)=a\lubintatepor {F_f}
\lambda=a_{f}(\lambda)$.

Ahora $\sigma=\Id\iff \sigma(1\bmod \pi^n)=a\bmod \pi^n=1\bmod \pi^n
\iff a-1\equiv 0\bmod \pi^n\iff a\in\unidades n$. $\fin$
\end{proof}

\begin{teorema}\label{CClaseT3.2.5.12} El campo $L_{f,n}$ depende \'unicamente de
$\pi$ y no de la elecci\'on de $f\in{\mc F}_{\pi}$. Esto es, para toda $n\geq 1$ y
para cualesquiera $f,g\in{\mc F}_{\pi}$, se tiene
$K(\Lambda_{f,n})=K(\Lambda_{g,n})$.
\end{teorema}

\begin{proof} Sean $f,g\in{\mc F}_{\pi}$ y $\lambda\in\Lambda_{f,n}$. Se tiene que
$1_{g,f}(Z)\in\o_K[[Z]]$, por tanto $1_{g,f}(\lambda)\in K(\lambda)\subseteq
L_{f,n}$. Puesto que $1_{g,f}\colon\Lambda_{f,n}\lra \Lambda_{g,n}$ es
biyectiva, $\Lambda_{g,n}\subseteq L_{f,n}$. Se sigue que $L_{g,n}=K(
\Lambda_{g,n})\subseteq L_{f,n}$. Por simetr\'ia tenemos que
$L_{f,n}\subseteq L_{g,n}$. $\fin$
\end{proof}

\begin{definicion}\label{CClaseD3.2.5.13} Para cualquier $f\in{\mc F}_{\pi}$, se
define $L_{\pi,n}:=L_{f,n}$ y $L_{\pi}:=\bigcup_{n=1}^{\infty} L_{f,n}$.
\end{definicion}

Debido al Teorema \ref{CClaseT3.2.5.12}, se puede suponer que $L_{\pi,n}$
est\'a generado por las ra\'ices de $f^{(n)}(Z)$ donde $f(Z)=\pi Z+Z^q\in
{\mc F}_\pi$. La extensi\'on $L_{\pi,n}/K$ es una extensi\'on de Galois.

\begin{definicion}\label{CClaseD3.2.5.14} Se define $G_{\pi,n}:=\Gal(L_{\pi,n}/K)$
y $G_{\pi}=\Gal(L_{\pi}/K)=\Gal(\bigcup_{n=1}^{\infty}L_{\pi,n}/K)=
\lim\limits_{\leftarrow n}G_{\pi,n}$.
\end{definicion}

Consideremos $\sigma\in G_{\pi,n}$ y $\lambda\in\Lambda_{f,n}$. Sea
$f(Z)=\pi Z+\sum_{i=1}^{\infty}a_i Z^i$ con $a_i\in\o_K$. Se tiene que
$f^{(n)}(Z)=\sum_{j=2}^{\infty}b_j Z^j\in\o_K[[Z]]$
para $n\geq 2$ y $f^{(n)}(\lambda)$
es convergente. Como $\sigma$ act\'ua de manera continua en $L_{\pi,n}$,
\[
0=\sigma(0)=\sigma(f^{(n)}(\lambda))=\sum_{j=2}^{\infty}\sigma(b_j\lambda^j)
\igual_{\substack{\uparrow\\ b_j\in K}} \sum_{j=2}^{\infty} b_j(\sigma\lambda)^j
=f^{(n)}(\sigma\lambda).
\]
Por tanto $f^{(n)}(\sigma\lambda)=0$ y $\sigma\lambda\in\Lambda_{f,n}$.

Puesto que $\Aut_{\o_K}(\Lambda_{f,n})\cong U_K/\unidades n$, cada
clase $u\unidades n\in U_K/\unidades n$ da lugar al automorfismo $\mu_f
\colon \Lambda_{f,n}\lra \Lambda_{f,n}$.

\begin{teorema}\label{CClaseT3.2.5.15}
Para cada $\sigma\in G_{\pi,n}$ existe una \'unica clase $u_{\sigma}
\unidades n=u\unidades n
\in U_K/\unidades n$ tal que $\sigma(\lambda)=(u_\sigma)_f(\lambda)$,
$\lambda\in\Lambda_{f,n}$. El mapeo $\sigma\mapsto u_{\sigma}\unidades n$
da lugar a un isomorfismo $G_{\pi,n}\cong U_K/\unidades n$. Adem\'as
$L_{\pi,n}=K(\lambda)$ donde $\Irr(Z,\lambda,K)=\frac{f^{(n)}(Z)}{f^{(n-1)}(Z)}$.
\end{teorema}

\begin{proof} Notemos que si $a\in\o_K$ y $\lambda\in \Lambda_{f,n}$ entonces
$a\lubintatepor {F_f}\lambda=a_f(\lambda)=\sum_{i=1}^{\infty} c_i\lambda^i
\in \o_K[[\lambda]]$ y por tanto
\[
\sigma(a\lubintatepor {F_f}
\lambda)=\sigma(a_f(\lambda))=\sigma\Big(\sum_{
i=1}^{\infty}c_i\lambda^i\Big)
=\sum_{i=1}^{\infty}c_i(\sigma\lambda)^i=a_f(\sigma\lambda)=
a\lubintatepor {F_f}(\sigma\lambda).
\]

En otras palabras, las acciones de $\o_K$ y de $G_{\pi,n}$ sobre $\Lambda_{
f,n}$ conmutan.

Cada $\sigma\in G_{\pi,n}$ induce un automorfismo del $\o_K$--m\'odulo 
$\Lambda_{f,n}$, esto es, existe $u\unidades n\in U_K/\unidades n$ tal que
$\sigma \lambda=u_f(\lambda)=u\lubintatepor {F_f} \lambda$ para toda $\lambda\in
\Lambda_{f,n}$. 

Se tiene el mapeo $G_{\pi,n}\xrightarrow{\ \varphi\ }U_K/\unidades n$,
$\sigma\longmapsto u_{\sigma}\unidades n$, donde $\sigma \lambda=u_{\sigma}
\lubintatepor {F_f} \lambda$ para toda $\lambda\in \Lambda_{f,n}$.

Como $\Lambda_{f,n}$ genera a $L_{\pi,n}$ sobre $K$, si $\sigma\in\ker \varphi$,
entonces $\sigma\lambda=u_f(\lambda)=\lambda$ para toda $\lambda
\in\Lambda_{f,n}$, por lo tanto $\sigma=\Id$.

Se tiene que $\big|U_K/\unidades n|=q^{n-1}(q-1)$.
Veamos que $|G_{\pi,n}|
\geq q^{n-1}(q-1)$. Se tiene $f^{(n)}(Z)=f(f^{(n-1)}(Z))=
f^{(n-1)}(Z)\phi_n(Z)$ con $\phi_n(Z)=(f^{(n-1)}(Z))^{q-1}
+\pi \in\o_K[Z]$ (recordemos que $f(Z)=\pi Z+Z^q$).

Ahora $f^{(2)}(Z)=f(f(Z))=f(Z)^q+\pi f(Z)=(Z^q+\pi Z)^q+\pi (Z^q+\pi Z)=
Z^{q^2}+\pi^qZ^q+\pi Z^q +\pi^2 Z$. En general se tiene $f^{(n-1)}(Z)=
Z^{q^{n-1}}+\pi (b_{n-2} Z^{q^{n-2}}+\cdots+b_2 Z^{q^2}) +\pi^{n-1} Z$.

Puesto que $\phi_n(Z)=(f^{(n-1)}(Z))^{q-1}+\pi\in\o_K[Z]$, $\phi_n(Z)$ es un 
polinomio de Eisenstein y por tanto irreducible sobre $K$. Si $\lambda$ es 
una ra\'iz de $\phi_n(Z)$, y por tanto ra\'iz de $f^{(n)}(Z)$,
entonces $K(\lambda)$
es una extensi\'on totalmente ramificada de $K$. Adem\'as 
\[
|G_{\pi,n}|\geq [K(\lambda):K]=q^{n-1}(q-1)=\deg\phi_n(Z)=\deg
\frac{f^{(n)}(Z)}{f^{(n-1)}(Z)}=\big|U_K/\unidades n\big|,
\]
por tanto $G_{\pi,n}\cong U_K/\unidades n$ y 
adem\'as $L_{\pi,n}=K(\lambda)$
donde $\lambda$ es una ra\'iz del polinomio de Eisenstein $\phi_n(Z)$ y
$\Irr(Z,\lambda,K)=\phi_n(Z)=\frac{f^{(n)}(Z)}{f^{(n-1)}(Z)}$. $\fin$
\end{proof}

Hemos demostrado m\'as de lo enunciado:

\begin{teorema}\label{CClaseT3.2.5.16} La extensi\'on $L_{\pi,n}/K$ es una 
extensi\'on abeliana y totalmente ramificada de grado $q^{n-1}(q-1)$ y es
generada por una ra\'iz de
\begin{gather*}
\phi_n(Z)=(f^{(n-1)}(Z))^{q-1}+\pi =\frac{f^{(n)}(Z)}{f^{(n-1)}(Z)}. \tag*{$\fin$}
\end{gather*}
\end{teorema}

\begin{corolario}\label{CClaseC3.2.5.16'} $\pi$ es una norma de $L_{\pi,n}$ a $K$.
\end{corolario}

\begin{proof} Se tiene $\phi_n(Z)=\prod_{\sigma\in G_{\pi,n}}(Z-\lambda^{\sigma})=
(f^{(n-1)}(Z))^{q-1}+\pi$. Por lo tanto
$\pi=\prod_{\sigma\in G_{\pi,n}}(-\lambda)^{\sigma}=\N_{L_{\pi,n}/K} (-\lambda)$
por lo que $\pi$ es una norma. $\fin$
\end{proof}

\begin{observacion}\label{CClaseO3.2.5.17} Todo este desarrollo se basa en un
elemento primo $\pi\in\o_K$ fijo. Esto corresponde al caso primo infinito
$\pK$ en el caso de campos de funciones, $\F(z)$, $\pK$ el polo de $x$. As\'i
que, debemos estudiar la misma situaci\'on cuando tomamos otro elemento
primo $\pi'$.
\end{observacion}

\begin{notacion}\label{CClaseN3.2.5.17'} Dado un campo local $K$,  $T$ denotar\'a
la m\'axima extensi\'on\label{CClasemaximaextensionnoramificadalocal}
no ramificada de $K$ contenida en una cerradura algebraica $\bar{K}$
de $K$ fija. Es decir $T=K^{\rm{nr}}$.
\end{notacion}

Se tiene $\Gal(T/K)\cong \Gal(\overline{\F}/\F)\cong \hat{{\ma Z}}$.

\begin{teorema}\label{CClaseT3.2.5.19} Se tiene $G_{\pi}\cong U_K$.
\end{teorema}

\begin{proof}
Del isomorfismo $G_{\pi,n}\cong
U_K/\unidades n$, se sigue que
\begin{gather*}
G_{\pi}=\Gal(L_{\pi}/K)\cong \lim_{\leftarrow n} G_{\pi,n}\cong
\lim_{\leftarrow n} U_K/\unidades n\cong U_K. \tag*{$\fin$}
\end{gather*}
\end{proof}

Sea $\tau_K\in \Gal(T/K)$ el automorfismo de Frobenius de $T/K$, 
y sea $\bar{\tau}_K$ la \'unica extensi\'on continua de $\tau_K$
en la completaci\'on $T_{\pK}$ de $T$  en $\bar{K}_{\pK}$,
una cerradura algebraica de $K_{\pK}$.
Usaremos la notaci\'on $\tau=\tau_K=\bar{\tau}_K$. Por 
definici\'on, $\tau_K$ induce el automorfismo $\alpha\mapsto \alpha^q$
en el campo residual $T(\pK_T)$
\[
\alpha^{\tau}=\tau(\alpha)=\alpha^q\bmod \pK_T \quad\text{para toda}
\quad \alpha\in \o_T,
\]
donde $\o_T$ denota el anillo de enteros de $T$ y $\pK_T$
denota el ideal m\'aximo de $\o_T$.

Usaremos la notaci\'on $\bar{T}$ para la completaci\'on $T_{\pK_T}$
de $T$. Se tiene la igualdad de campos residuales $\bar{T}(\pK_T)=
\bar{T}_{\pK_{\bar{T}}}(\pK_{\bar{T}})$ el cual denotaremos 
simplemente por $\bar{T}(\pK)$.

Consideremos los endomorfismos
\begin{align*}
\tau-1\colon \o_{\bar{T}}&\lra \o_{\bar{T}}\\
\alpha&\longmapsto (\tau-1)(\alpha)=\tau(\alpha)-\alpha,\\
\tau-1\colon U_{\bar{T}}&\lra U_{\bar{T}}\\
u&\longmapsto u^{\tau-1}=\tau(u)/u,
\end{align*}
donde $U_{\bar{T}}=\{x\in\bar{T}\mid |x|=1\}$. Se tiene el siguiente resultado.

\begin{lema}\label{CClaseL3.2.5.20}
Las siguientes sucesiones
\begin{gather*}
0\lra \o_K\lra \o_{\bar{T}}\xrightarrow{\ \tau-1\ }\o_{\bar{T}}\lra 0,\\
1\lra U_K\lra U_{\bar{T}}\xrightarrow{\ \tau-1\ }U_{\bar{T}}\lra 1,
\end{gather*}
son exactas. En particular $(\tau-1)\o_{\bar{T}}=\o_{\bar{T}}$ y 
$U_{\bar{T}}^{\tau -1}= U_{\bar{T}}$.
\end{lema}

\begin{proof} Las demostraciones son totalmente paralelas y por tanto la haremos
principalmente para $U_{\bar{T}}$.

Puesto que el campo residual $S:=\bar{T}(\pK)=\o_{\bar{T}}/
\bar{T}_{\pK_{\bar{T}}}\cong \o_T/{\pK}_T=T(\pK)
\cong \bar{\ma F}_q$ es algebraicamente cerrado,
los mapeos  $\alpha\longmapsto \tau(\alpha)-\alpha$ y $\alpha\longmapsto
\alpha^{\tau-1}=\tau(\alpha)/\alpha$ son suprayectivos de $S\lra S$ y de
$\*S\lra\* S$ respectivamente.

Se tiene 
\begin{gather}
U_{\bar{T}}/U_{\bar{T}}^{(1)}\cong \*S=\*{\bar{T}(\pK)},\quad
U_{\bar{T}}^{(n)}/U_{\bar{T}}^{(n+1)}\cong S^+=\bar{T}(\pK)^+
\quad\text{y}\label{CClaseEc3.2.5.22} \\
\o_{\bar{T}}/\pK_{\bar{T}}\cong \pK_{\bar{T}}^n/
\pK_{\bar{T}}^{n+1}\cong
S^+=\bar{T}(\pK)^+.\nonumber
\end{gather}

Ahora bien, si $x\in U_{\bar{T}}$ (resp. $x\in\o_{\bar{T}}$), existe $y_1\in
\*{T}(\pK)$ (resp. $y_1\in T(\pK)^+$) tal que $\bar{x}=\tau 
\bar{y}_1/\bar{y}_1\in\*{\bar{T}(\pK)}$ (rest. $\bar{x}=
\tau\bar{y}_1-\bar{y}_1\in \bar{T}(\pK)^+$) donde $\bar{z}=
z\bmod \pK$, por lo que $x=\frac{\tau y_1}{y_1} a_1$ con $y_i\in U_{\bar{T}}$
y $a_1\in U_{\bar{T}}^{(1)}$ (resp. $x=\tau y_1-y_1+a_1$ con
$y_1\in \o_{\bar{T}}$, $a_1\in \pK_{T}$).

De la Ecuaci\'on (\ref{CClaseEc3.2.5.22}), podemos continuar el proceso y tenemos
$a_1=\frac{\tau y_2}{y_2}a_2$, con $y_2\in U_{\bar{T}}^{(1)}$ y $a_2
\in U_{\bar{T}}^{(2)}$ (resp. $a_1=\tau y_2-y_2+a_2$ con $y_2\in \pK_{T_{
\pK}}$ y $a_2\in \pK_{\bar{T}}^2$). Se sigue que $x=\frac{\tau(y_1y_2)}{
y_1y_2}a_2$ (resp. $x=\tau(y_1+y_2)-(y_1+y_2)+a_2$).

En general obtendremos
\begin{gather*}
x=\frac{\tau(y_1\cdots y_n)}{y_1\cdots y_n} a_n, \quad y_n\in U_{\bar{T}}^{(n-1)},
\quad a_n\in U_{\bar{T}}^{(n)},\\
\text{(resp.\ }x=\tau(y_1+\cdots+y_n)-(y_1+\cdots+y_n)+a_n, \quad 
y_n\in \pK_{\bar{T}}^{n-1}, \quad a_n\in \pK_{\bar{T}}^n\text{)}.
\end{gather*}

Pasando al l\'imite y puesto que $\bar{T}$ es completo, se obtiene
\[
x=\frac{\tau y}{y}, y=\prod_{n=1}^{\infty}y_n\in U_{\bar{T}} 
\quad \text{(resp.\ }
x=\tau y-y,\quad y=\sum_{n=1}^{\infty}y_n\in \o_{\bar{T}}\text{)}
\]
lo cual prueba en ambos caso la suprayectividad de $\tau -1$.

Puesto que $\tau\in\Gal(\bar{T}/K)$ y $U_K\subseteq K$, 
se tiene que $U_K$ est\'a
contenido en el $\ker(\tau-1)\colon U_{\bar{T}}\lra U_{\bar{T}}$.

Ahora sea $\xi^{\tau-1}=1$, esto es, $\xi^{\tau}=\xi$, con $\xi\in U_{\bar{T}}$.
Se tiene que $\{0\}\bigcup \bigcup_{m=1}^{\infty}V_m=\{0\}\bigcup V_{\infty}$
es un conjunto de representantes de $\bar{T}(\pK)$, donde $V_m=
\{\zeta_{q^{m-1}}^i\mid 0\leq i\leq q^m-2\}$, por lo que $\xi=\sum_{n=0}^{
\infty}a_n\pi^n$ con $a_n\in\{0\}\bigcup V_{\infty}$.

Se tiene que
$\tau(a_n)=a_n^q$ y $\tau(\xi) =\sum_{n=0}^{\infty}a_n^q\pi^n=\sum_{
n=0}^{\infty}a_n \pi^n=\xi$ lo cual implica que $a_n=0$ o $a_n^{q-1}=1$
por lo que $a_n\in V_1=\{\zeta_{q-1}^i\mid 0\leq i\leq q-2\}$. Se sigue
que $\xi\in K\cap U_{\bar{T}}=U_K$ y por tanto
\begin{gather*}
1\lra U_K\lra U_{\bar{T}}\xrightarrow{\ \tau-1\ }U_{\bar{T}}\lra 1
\end{gather*}
es exacta. La demostraci\'on de la exactitud de la otra sucesi\'on es 
similar. $\fin$
\end{proof}

El Lema \ref{CClaseL3.2.5.20} nos permite probar un resultado similar
al Teorema \ref{CClaseT3.2.5.4} pero ahora con respecto a un cambio en 
el elemento primo.

\begin{teorema}\label{CClaseT3.2.5.23} Sean $\pi$ y $\pi'=a\pi$ dos elementos
primos de $K$, $a\in U_K$. Sean $f\in{\mc F}_{\pi}$ y $f'\in{\mc F}_{\pi'}$.
Entonces existe una serie de potencias $\theta(Z)\in \o_{\bar{T}}[[Z]]$,
tal que
\las
\item $\theta(Z)\equiv \varepsilon Z\bmod \deg 2$, $\varepsilon \in U_{\bar{T}}$,
\item $\theta^{\tau}(Z)=\theta(a_f(Z))$,
\item $\theta(F_f(X,Y))=F_{f'}(\theta(X),\theta(Y))$,
\item $\theta(b_f(Z))=b_{f'}(\theta(Z))$ para toda $b\in\o_{K}$,
\end{list}
donde $\theta^{\tau}$ denota a la serie obtenida a partir de la de $\theta$
aplicando el automorfismo de Frobenius $\tau$ a los coeficientes de $\theta$.
\end{teorema}

Antes de probar el teorema, hagamos la siguiente observaci\'on. En general
si $\alpha$ y $\beta$ son dos series en $\o_{\bar{T}}[[Z]]$ y si $\sigma$
es un automorfismo de $\bar{T}$, entonces $\alpha^{\sigma},\beta^{\sigma}$
representan a las series cuyos coeficientes son obtenidos a partir de $\alpha$
y $\beta$ aplicando $\sigma$ a cada uno de los coeficientes, esto es, si
$\alpha(Z)=\sum_{i=0}^{\infty}a_iZ^i$, entonces $\alpha^{\sigma}(Z)=
\sum_{i=0}^{\infty} a_i^{\sigma}Z^i=\sum_{i=0}^{\infty} \sigma(a_i)Z^i$.

Entonces, con c\'alculos directos, se puede demostrar que
 $(\alpha\circ\beta)^{\sigma}=\alpha^{\sigma}\circ \beta^{\sigma}$.

Aqu\'i, $\circ$ denota la composici\'on de series. M\'as precisamente,
para dos series $\alpha, \beta\in\o_{\bar{T}}[[Z]]$, $\alpha\circ \beta$ denota
$(\alpha\circ\beta)(Z)=\alpha(\beta(Z))$. 

\bigskip

\begin{proof} (Teorema \ref{CClaseT3.2.5.23}). Por el Lema \ref{CClaseL3.2.5.20}, se tiene
que existe $\varepsilon\in U_{\bar{T}}$ tal que $a=\frac{\tau(\varepsilon)}{
\varepsilon}$. Sea $\theta_1(Z):=\varepsilon Z$. Supongamos que se ha
construido un polinomio $\theta_n(Z)$ de grado $n$ tal que
\[
\theta_n^{\tau}(Z)\equiv \theta_n(a_f(Z))\bmod \deg (n+1),
\]

Se quiere construir un polinomio $\theta_{n+1}(Z)=\theta_n(Z)+bZ^{n+1}$
que satisfaga 
\[
\theta_{n+1}^{\tau}(Z)\equiv \theta_{n+1}(a_f(Z))\bmod \deg (n+2).
\]

Sea $b=\gamma \varepsilon^{n+1}$ para alg\'un $\gamma$. Si se tiene el polinomio buscado
$\theta_{n+1}(Z)$, entonces
\begin{align*}
\theta_n^{\tau}(Z)-\theta_n(a_f(Z))&=c Z^{n+1} + \text{\ t\'erminos de grado mayor},\\
\theta_{n+1}^{\tau}(Z)-\theta_{n+1}(a_f(Z))&=\theta_n^{\tau}(Z)+\tau(b)Z^{n+1}-
\theta_n(a_f(Z))-b(a_f(Z))^{n+1}\\
&\igual_{\substack{\uparrow\\ a_f(Z)=aZ+\cdots}}
(c+\tau(b)-ba^{n+1})Z^{n+1}\\
&\hspace{2cm}+\text{\ t\'erminos de grado mayor}
\end{align*}

Puesto que queremos $\theta_{n+1}^{\tau}(Z)\equiv \theta_{n+1}(a_f(Z))\mod\deg 
(n+2)$, se tiene que se debe satisfacer $c+\tau(b)-ba^{n+1}=0 =c
+\tau(\gamma)\tau(\varepsilon)^{n+1}-\gamma\frac{\varepsilon^{n+1}\tau(
\varepsilon)^{n+1}}{\varepsilon^{n+1}}$, por lo que se debe tener
\[
c+(\tau(\gamma)-\gamma)\tau(\varepsilon)^{n+1}=0 \quad\text{o, equivalentemente,}
\quad \gamma-\tau(\gamma)=c/\tau(\varepsilon)^{n+1}.
\]

Tal $\gamma$ existe como consecuencia del Lema \ref{CClaseL3.2.5.20}, de
donde obtenemos $\theta_{n+1}(Z)$ y la serie $\theta(Z)=\lim_{n\to\infty}
\theta_n(Z)$ la cual satisface la condici\'on
\begin{gather}\label{CClaseEc8.X1}
\theta^{\tau}(Z)=\theta(a_f(Z)).
\end{gather}

Necesitamos modificar $\theta(Z)$ para satisfacer las condiciones (3) y (4)
del teorema.

Para una serie $\psi\in \o_{\bar{T}}$, $\psi^{-1}$ denota la serie inversa de
$\psi$: $\psi\circ\psi^{-1}=\psi^{-1}\circ\psi=\Id$, donde $\Id(Z)=Z$,
en caso de existir. Si $\psi(Z)=\delta_1Z+\sum_{i=2}^{\infty}\delta_iZ^i$ con
$\delta_i\in \bar{T}$, $\delta_1\neq 0$, $\psi^{-1}$ se encuentra sustituyendo
directamente.

Ahora bien, $\theta(Z)\equiv \varepsilon Z\bmod \deg 2$, $\theta(Z)\in\o_{\bar{T}}
[[Z]]$ con $\varepsilon \in U_{\bar{T}}$. Entonces $\theta^{-1}(Z)=
\sum_{i=1}^{\infty}c_iZ^i$ debe satisfacer 
\[
\theta(\theta^{-1}(Z))=\varepsilon\Big(\sum_{i=1}^{\infty}c_iZ^i\Big)+
\sum_{j=2}^{\infty}d_j\Big(\sum_{i=1}^{\infty}c_iZ^i\Big)^j = (\varepsilon c_1) Z+
(\varepsilon c_2+d_2 c_1)Z^2+\cdots
\]
la cual se resuelve para las $c_i$'s de manera recursiva. En particular, puesto que
$\varepsilon\in U_{\bar{T}}$, $\theta^{-1}(Z)\in \o_{\bar{T}}[[Z]]$.

Consideremos la serie
\begin{gather}\label{CClaseEc8.X2}
h=\theta^{\tau}\circ f\circ \theta^{-1}\in\o_{\bar{T}}[[Z]]
\end{gather}
en donde $\circ$ significa composici\'on o evaluaci\'on, es decir,
$g\circ l(Z)=g(l(Z))$.

De la Ecuaci\'on (\ref{CClaseEc8.X1}) se obtiene $\theta^{\tau}=
\theta\circ a_f$. Por tanto, de la Ecuaci\'on (\ref{CClaseEc8.X2}) obtenemos
\[
h=\theta^{\tau}\circ f\circ \theta^{-1}=\theta\circ a_f\circ f\circ \theta^{-1}.
\]

Del Teorema \ref{CClaseT3.2.5.6} (6) se tiene $f(Z)=(\pi)_f(Z)$. Por tanto
\begin{gather*}
h=\theta\circ (a)_f\circ (\pi)_f\circ \theta^{-1}\igual_{\substack{\uparrow\\
\text{Teorema \ref{CClaseT3.2.5.6} (4)}}} \theta\circ (a\pi)_f\circ \theta^{-1}
\igual_{\substack{\uparrow\\ a\pi=\pi'}}\theta\circ(\pi')_f\circ \theta^{-1}.
\intertext{Se sigue que}
h^{\tau}=\theta^{\tau}\circ(\pi')_f^{\tau}\circ\theta^{-\tau}
\igual_{\substack{\uparrow\\ \text{(\ref{CClaseEc8.X2})}}}\theta^{\tau}
\circ (\pi')_f^{\tau}\circ(a)_f^{-1}\circ\theta^{-1}.
\end{gather*}

Puesto que $\tau\in\Gal(T/K)$, se tiene
\[
(\pi')_f^{\tau}=(\pi')_f=(a\pi)_f=(a)_f\circ(\pi)_f\igual_{\substack{\uparrow\\
\text{Teorema \ref{CClaseT3.2.5.6} (6)}}} (a)_f\circ f.
\]

Se sigue que 
\begin{gather*}
h^{\tau}=\theta^{\tau}\circ (\pi')_f\circ (a)_f^{-1}\circ\theta^{-1}=
\theta^{\tau}\circ (a\pi a^{-1})_f\circ \theta^{-1}=\theta^{\tau}
\circ(\pi)_f\circ \theta^{-1},
\intertext{esto es,}
h^{\tau}=\theta^{\tau}\circ f\circ \theta^{-1}=h
\end{gather*}
por lo que $h\in \o_K$
pues cualquier elemento de $\bar{T}$ fijado por $\tau$
pertenece a $K$ debido a que $\Gal(T/K)\cong \Gal(\bar{\ma F}_q/\F)=
\overline{\langle\tau\rangle}$.

Por otro lado, puesto que $h=\theta\circ \pi'_f\circ \theta^{-1}$ y que
$\theta(Z)\equiv \varepsilon Z\bmod \deg 2$, se obtiene
\begin{gather*}
h(Z)=\theta((\pi')_f(\theta^{-1}(Z)))\equiv \varepsilon \pi'\varepsilon^{-1}Z
\equiv \pi' Z\bmod\deg 2,
\intertext{y}
h(Z)=\theta^{\tau}(f(\theta^{-1}(Z)))\equivalente_{\substack{\uparrow\\ (\ast)}}
\theta^{\tau}(\theta^{-1}(Z)^q)\equivalente_{\substack{\uparrow\\ \tau=q}}
\theta^{\tau}(\theta^{-\tau}(Z^q))\equiv Z^q \bmod \pi' \text{\ y $\pi$}
\end{gather*}
donde ($\ast$): $f(X)\equiv X^q\bmod \pi$, $f(X)\equiv \bmod \pi'$, $f(X)=
\sum_{i=1}^{\infty}\alpha_iZ^i$, $\pi|\alpha_i$ para $i\neq q$ y $\pi'=a^{-1}
\pi$ por lo que $\pi'|\alpha_i$ para $i\neq q$.

Por tanto $h\in{\mc F}_{\pi'}$. Sea $1_{f',h}\colon F_h\lra F_{f'}$ y consideramos
$\theta_1=1_{f',h}\circ \theta$. Entonces
\begin{gather*}
\theta_1(Z)=(1_{f',h}\circ \theta)(Z)\equiv \theta(Z)\bmod\deg 2\equiv \varepsilon
Z \bmod \deg 2,\\
\theta_1^{\tau}=1_{f',h}^{\tau}\circ (\theta^{\tau}(Z))=1_{f',h}^{\tau}(\theta(a_f)
(Z))=(1_{f',h}\circ\theta)(a_f(Z))=\theta_1(a_f(Z)).
\end{gather*}

Se sigue que $\theta_1$ tambi\'en satisface (1) 
y (2) del teorema. Ahora, si $h_1
=\theta_1^{\tau}\circ f\circ \theta_1^{-1}=1_{f',h}\circ h=f'$, se tiene
\[
f'=\theta_1^{\tau}\circ f\circ \theta_1^{-1}=\theta_1\circ \pi'_f\circ \theta^{-1}.
\]

Para probar que $\theta_1(F_f(X,Y))=F_{f'}(\theta(X),\theta(Y))$ basta probar que
la serie $F(X,Y):=\theta_1(F_f(\theta_1^{-1}(X),\theta_1^{-1}(Y)))$ satisface la
caracterizaci\'on de $F_{f'}$, esto es:
\begin{gather*}
F(X,Y)\equiv X+Y\bmod \deg 2\quad\text{y}\\
f'(F(X,Y))=F(f'(X),f'(Y)).
\end{gather*}

Se tiene que 
\begin{align*}
F(X,Y)&=\theta_1(F_f(\theta^{-1}(X),\theta_1^{-1}(Y)))\equiv \varepsilon
F_f(\theta_1^{-1}(X),\theta_1^{-1}(Y))\\
&\equiv \varepsilon(\theta_1^{-1}(X)+
\theta_1^{-1}(Y))\equiv \varepsilon(\varepsilon^{-1} X+\varepsilon^{-1} Y)\\
&\equiv X+Y\bmod\deg 2,\\
F(f'(X),f'(Y))&=\theta_1(F_f(\theta_1^{-1}(f'(X)),\theta^{-1}_1(f'(Y))))
\igual_{\substack{\uparrow\\ f'=\theta_1\circ \pi'_f\circ \theta_1^{-1}\\
\theta_1^{-1}\circ f'=\pi'_f\circ \theta_1^{-1}}}\\
&= \theta_1(F_f(\pi'_f(\theta_1^{-1}(X)),
\pi'_f(\theta_1^{-1}(Y))))\\
&= \theta_1(\pi'_f(F_f(\theta_1^{-1}(X),\theta_1^{-1}(Y))))
\igual_{\substack{\uparrow\\ \theta_1\circ \pi'_f=f'\circ\theta_1}}\\
&=f' \theta_1(F_f(\theta_1^{-1}(X),\theta_1^{-1}(Y)))=(f'\theta_1F_f\theta^{-1})(X,Y)\\
&=f'F(X,Y),
\end{align*}
por lo tanto $F=F_{f'}$.

Para probar (4) se necesita $\theta_1 b_f=b_{f'}\theta_1$, $b\in\o_K$.

Sea $H=\theta_1 b_f\theta_1^{-1}$. Se quiere probar que $H=b_{f'}$ y para ello
es suficiente probar que $H$ satisface la caracterizaci\'on de $b_{f'}$, esto es:
\begin{gather*}
H(X)\equiv b X\bmod\deg 2\quad\text{y}\quad
f(H(X))=H(f'(X)).
\end{gather*}

Ahora 
\begin{align*}
H(X)&=\theta_1 b_f \theta_1^{-1}(X)\equiv (\theta_1 b\theta_1^{-1})(X)
\equiv \varepsilon b\varepsilon^{-1} X\equiv b X\bmod \deg 2\\
H(f'(X))&=\theta_1 b_f\theta_1^{-1}f'(X)\igual_{\substack{\uparrow\\ f'=\theta_1
\pi'_f\theta_1^{-1}}}\theta_1 b_f\theta_1^{-1}\cdot \theta_1\pi'_f\theta_1^{-1}(X)\\
&=(\theta_1b_f\pi'_f\theta_1^{-1})(X)=\theta_1((b\pi')_f\theta_1^{-1})(X)=
(\theta_1\pi'_fb_f\theta_1^{-1})(X)\\
&=\theta_1\pi'_f\theta_1^{-1}\theta_1b_f\theta_1^{-1}(X)\igual_{\substack{\uparrow\\
\theta_1\pi'_f\theta_1^{-1}=f'}} f'(H(X)).
\end{align*}
Por tanto $H=b_{f'}$ y se tiene (4). $\fin$
\end{proof}

La importancia del Teorema \ref{CClaseT3.2.5.23} es el siguiente teorema.

\begin{teorema}\label{CClaseT3.2.5.24} Sean 
$\pi$ y $\pi'=a\pi$ dos elementos primos en 
$K$, $a\in U_K$ y sean $f\in{\mc F}_{\pi}$, $f'\in{\mc F}_{\pi'}$. Entonces
el mapeo $\lambda\longmapsto \theta(\lambda)$ da lugar a un isomorfismo de
$\o_K$--m\'odulos
\[
\Lambda_{f,n}\cong \Lambda_{f',n}
\]
para toda $n\in{\ma N}$.

\end{teorema}

Por el  Teorema \ref{CClaseT3.2.5.10} ya sab\'iamos que 
\[
\Lambda_{f,n}\cong
\o_K/\pi^n\o_K=\o_K/(\pi')^n\o_K\cong \Lambda_{f',n}.
\]
La afirmaci\'on del
teorema es que el mapeo dado por $\theta$ es un isomorfismo.

\begin{proof} Si $\lambda\in\Lambda_{f,n}$, entonces
\[
(f')^{(n)}(\theta(\lambda))=(\pi')_f^{(n)}(\theta(\lambda))=\theta((a^n\pi^n)_f)(\lambda)=
\theta(a^nf^{(n)}(\lambda))=\theta(0)=0,
\]
por tanto $\lambda\in \Lambda_{f',n}$. Se sigue que $\lambda\longmapsto \theta(
\lambda)$ manda $\Lambda_{f,n}$ en $\Lambda_{f',n}$.

Por (3) y (4) del Teorema \ref{CClaseT3.2.5.23}, se tiene
\begin{gather*}
\theta(\lambda\lubintatemas {F_f}\mu)=\theta(F_f(\lambda,\mu))=F_{f'}(\theta(
\lambda),\theta(\mu))=\theta(\lambda)\lubintatemas {F_{f'}} \theta(\mu),\\
\theta (a\lubintatepor {F_f} \lambda)=\theta(a_f(\lambda))=a_{f'}\theta(\lambda)
=a\lubintatepor {F_{f'}}\theta(\lambda).
\end{gather*}
Por tanto $\theta$ es un homomorfismo de $\o_K$--m\'odulos de $\Lambda_{f,n}$
en $\Lambda_{f',n}$. 

Veamos que $\theta$ es 1--1. Sea $\theta(\lambda)=0$, 
$\lambda\in\Lambda_{f,n}$. Ahora 
\[
\theta(\lambda)=\varepsilon \lambda+\sum_{
i=0}^{\infty}\alpha_i\lambda^i=0. 
\]
Si $\lambda\neq 0$ esto implicar\'ia que $0=
\varepsilon+\lambda\big(\sum_{i=2}^{\infty}\alpha_i\lambda^{i-2}\big)$ lo
cual implicar\'ia que $\theta(\varepsilon)=0=\theta(\lambda)\theta\big(
-\sum_{i=2}^{\infty}\alpha_i\lambda^{i-2}\big)$.

Lo anterior no es posible pues $\varepsilon$ es una unidad y 
$\theta$ es una serie invertible.

Puesto que $\Lambda_{f,n}\cong \o_K/\pi^n\o_K\cong \o_K/(\pi')^n\o_K\cong
\Lambda_{f',n}$ tienen la misma cardinalidad, $\theta$ es suprayectiva y
por tanto $\theta$ es un isomorfismo. $\fin$
\end{proof}

\begin{observacion}\label{CClaseO3.2.5.25} Para dos elementos primos $\pi$ y $\pi'$
de $K$, los campos $L_{\pi,n}$ y $L_{\pi',n}$ 
pueden ser diferentes. Sin embargo, se tiene
que $TL_{\pi,n}=TL_{\pi',n}$. Esto es el contenido 
de la Proposici\'on \ref{CClaseP3.2.5.26}.
\end{observacion}

Primero probamos el siguiente lema.

\begin{lema}\label{CClaseL3.3.37'}
Sea $K\subseteq F\subseteq \sep K$, $\sep K$ una cerradura separable
de $K$. Entonces $F$ es un conjunto cerrado en $\sep K$.
\end{lema}

\begin{proof} Sea $H:=\Gal(\sep K/F)$. Entonces $H$ fija a todo elemento
de $F$ y por tanto fija a todo elemento de la cerradura $\bar{F}$
de $F$ en $\sep K$ por continuidad. Por tanto
\begin{gather*}
F\subseteq \bar{F}\subseteq (\sep K)^H=F. \tag*{$\fin$}
\end{gather*}
\end{proof}

\begin{proposicion}\label{CClaseP3.2.5.26} Si $\pi$ y 
$\pi'$ son dos elementos primos
de $K$ y $T$ es la m\'axima extensi\'on 
abeliana no ramificada de $K$, entonces
\[
TL_{\pi,n}=TL_{\pi',n}.
\]
\end{proposicion}

\begin{proof} Por el Teorema \ref{CClaseT3.2.5.24} se tiene que si
$\lambda\in\Lambda_{f,n}$, entonces $\theta(\lambda)=
\varepsilon \lambda+\sum_{i=2}^{\infty} \alpha_i\lambda^i\in \bar{T}
(\lambda)$. Por lo tanto 
\[
\bar{T}(\Lambda_{f',n})=\bar{T}(\theta(\Lambda_{f,n}))\subseteq
\bar{T}(\Lambda_{f,n})=\bar{T}(\theta^{-1}(\Lambda_{f',n}))
\subseteq \bar{T}(\Lambda_{f',n}).
\]

Por lo tanto $\overline{TL}_{\pi,n}=\overline{TL}_{\pi',n}$.

Por el Lema \ref{CClaseL3.3.37'}, se sigue que $T(\Lambda_{\pi,n})=
T(\Lambda_{\pi',n})$, de donde obtenemos que $TL_{\pi,n}
=TL_{\pi',n}$.  $\fin$
\end{proof}

Por otro lado, puesto que $T/K$ es no ramificada y $L_{\pi,n}/K$ es totalmente
ramificada, $T$ y $L_{\pi,n}$ son linealmente disjuntos sobre $K$ y
\begin{align}
\Gal(TL_{\pi,n}/K)&\cong \Gal(T/K)\times \Gal(L_{\pi,n}/K)=\Gal(T/K)\times 
G_{\pi,n}\nonumber\\
&\cong \Gal(\abe {\F}/\F)\times G_{\pi,n}\cong \hat{{\ma Z}}\times
(U_K/\unidades n).\label{CClaseEc3.2.5.27}
\end{align}

Sea $\rho_{\pi}\colon\* K\lra \Gal(TL_{\pi,n}/K)$ el siguiente homomorfismo.
Para $a=u\pi^m\in\* K$ con $u\in U_K$, $m\in{\ma Z}$. Entonces
\begin{align*}
\rho_{\pi}(a)|_T:&= \tau_K^m\in \Gal(T/K), \quad \text{$\tau_K$ el automorfismo de Frobenius},\\
\rho_{\pi}(a)|_{L_{\pi,n}}:&=\sigma_u\in G_{\pi,n},
\end{align*}
donde $\sigma_u$ es el automorfismo de $L_{\pi,n}$ que bajo el isomorfismo
$G_{\pi,n}\cong U_K/\unidades n$, $\sigma_u$ corresponde a la clase
$u^{-1}\unidades n$. En otras palabras $\rho_{\pi}(a)|_{L_{\pi,n}}$ est\'a determinado por
\[
\rho_{\pi}(a)(\lambda)=(u^{-1})_f(\lambda)
=u^{-1}\lubintatepor {F_f}\lambda,\quad \lambda\in \Lambda_{f,n}.
\]

El objetivo central de los grupos formales para nosotros, es probar que $\rho_{\pi}$ es 
precisamente el s\'imbolo residual de la norma $(\underline{\ },K)$, esto es, $\rho_{\pi}$
es el mapeo de reciprocidad.

\begin{teorema}\label{CClaseT3.2.5.28} Para $a\in \*K$ se tiene
\[
\rho_{\pi}(a)=(a,K)|_{TL_{\pi,n}}=
\rho_K|_{TL_{\pi,n}}(a).
\]
\end{teorema}

\begin{proof} Se tiene que los elementos primos de $\*K$ generan $\*K$,
por lo que basta probar el
teorema para elementos primos. Sea $a=\pi$ un elemento primo de $K$. Entonces
\[
\rho_{\pi}(\pi)|_T=\tau=(\pi,T/K)=(\pi,K)|_T
\]
como consecuencia del Teorema \ref{CCLT17.6.10} pues $T/K$ es no ramificada.

Ahora, por el Corolario \ref{CClaseC3.2.5.16'}, $\pi$ es una norma 
de $L_{\pi,n}$ a $K$.
Por el Teorema \ref{CClaseT3.2.2}, $(\pi,L_{\pi,n}
/K)=(\pi,K)|_{L_{\pi,n}}=1$. Se tiene $\pi=1\cdot \pi$, por lo que
\[
\rho_{\pi}(\pi)|_{L_{\pi,n}}=\sigma_1=
\Id_{L_{\pi,n}}=(\pi,K)|_{L_{\pi,n}}=1.
\]
Se sigue que $\rho_{\pi}(\pi)=(\pi,K)|_{TL_{\pi,n}}$.

Ahora sea $\pi'=u\pi$ con $u\in U_K$ otro elemento primo de $K$. Se tiene
$TL_{\pi,n}=TL_{\pi',n}$ y por el Teorema \ref{CCLT17.6.10}
\[
\rho_{\pi}(\pi')|_T=\tau=(\pi',T/K)=(\pi',K)|_T.
\]

Puesto que $\pi'$ es una norma de $L_{\pi',n}$ a $K$, se tiene que $(\pi',K)|_{
L_{\pi',n}}=(\pi',L_{\pi',n}/K)=\Id_{L_{\pi',n}}$. As\'i debemos verificar que $\rho_{\pi}(
\pi')|_{L_{\pi',n}}=\Id_{L_{\pi',n}}$. Esto equivale a ver que $\rho_{\pi}(\pi')(\mu)=\mu$
para toda $\mu\in\Lambda_{f',n}$ con $f'\in{\mc F}_{\pi'}$.

Sabemos del Teorema \ref{CClaseT3.2.5.24} que $\Lambda_{f',n}=\theta(\Lambda_{f,n})$.
Por tanto debemos probar que $\rho_{\pi}(\pi')(\theta(\lambda))=\theta(\lambda)$
para toda $\lambda\in\Lambda_{f,n}$.

Se tiene que $\rho_{\pi}(\pi')=\rho_{\pi}(u\pi)=\rho_{\pi}(u)\circ\rho_u(\pi)$.
Para $\lambda\in\Lambda_{f,n}$, $\rho_{\pi}(\pi)(\lambda)=\lambda$, $\rho_{\pi}
(\pi)|_T=\tau$, $\rho_{\pi}(u)|_T=\Id_T$.

Extendemos $\rho_{\pi}(\pi)|_T$ y $\rho_{\pi}(u)|_T$ continuamente a $\bar{T}$ y
usando el Teorema \ref{CClaseT3.2.5.23}
\begin{align*}
\rho_{\pi}(\pi')(\theta(\lambda))&= (\rho_{\pi}(u)\circ\rho_{\pi}(\pi))(\theta(\lambda))=
\rho_{\pi}(u)(\rho_{\pi}(\pi)(\theta(\lambda)))\\
&=\rho_{\pi}(u)(\theta^{\tau}(\lambda))=(\rho_{\pi}(u)\circ\theta^{\tau})(\lambda)=
\theta^{\tau}(\rho_{\pi}(u))(\lambda)\\
&=\theta^{\tau}((u^{-1})_f(\lambda))=\theta((u)_f(u^{-1})_f(\lambda))=\theta(\lambda).
\tag*{$\fin$}
\end{align*}
\end{proof}

De esta forma se puede describir el s\'imbolo de la norma residual $(\underline{\ },
L_{\pi,n}/K)$ de la extensi\'on abeliana y totalmente ramificada $L_{\pi,n}$:

\begin{teorema}\label{CClaseT3.2.5.29} Sea $a=u\pi^m\in\*K$, $u\in U_K$ y $m\in {\ma Z}$.
Entonces
\[
(a,L_{\pi,n}/K)(\lambda)=(u^{-1})_f(\lambda)\quad\text{para toda $\lambda\in 
\Lambda_{f,n}\subseteq L_{\pi,n}$}.
\]

El grupo de normas de la extensi\'on $L_{\pi,n}/K$ es el grupo $(\pi)\times \unidades n$.
\end{teorema}

\begin{proof} Se tiene del Teorema \ref{CClaseT3.2.5.28} 
que $(a,L_{\pi,n}/K)(\lambda)=\rho_{\pi}(a)(\lambda)$.
Por tanto $(a,K)(\lambda)=(a,L_{\pi,n}/K)(\lambda)=\rho_{\pi}(a)(\lambda)=
(u^{-1})_f(\lambda)$.

En consecuencia, $a\in \N_{L_{\pi,n}/K}(\*{L_{\pi,n}})\iff (a,L_{\pi,n}/K)(\lambda)=
(u^{-1})_f(\lambda)=\lambda$ para toda $\lambda\in \Lambda_{f,n}$. Por el
Teorema \ref{CClaseT3.2.5.11}, $(u^{-1})_f=\Id\iff u^{-1}\in\unidades n\iff u\in\unidades n
\iff a\in (\pi)\times \unidades n$. $\fin$
\end{proof}

El siguiente ejemplo es el origen del uso de los grupos formales en teor\'ia local
de campos locales.

\begin{ejemplo}[Ver el Ejemplo \ref{CClaseE3.3.21'}]\label{CClaseEj3.2.5.30}
Sea $K={\ma Q}_p$ el campo de los n\'umeros $p$--\'adicos. Entonces $p$ es un elemento
primo en $K$. Sea $f\in{\mc F}_p$ definido por
\[
f(Z)=(1+Z)^p-1=pZ+\binom{p}{2} Z^2+\binom{p}{3} Z^3+\cdots+p Z^{p-1}+Z^p.
\]
En nuestro caso, $q=p$. Se tiene
\[
f^{(2)}(Z)=f(f(Z))=(1+f(Z))^p-1=(1+(1+Z)^p-1)^p-1=(1+Z)^{p^2}-1
\]
y en general $f^{(n)}(Z)=(1+Z)^{p^n}-1$. El conjunto de ceros de $f^{(n)}(Z)$ son
$\{\lambda=\xi-1\mid \xi\text{\ es $p^n$--ra\'iz de $1$}\}$. Por tanto $L_{p,n}=
{\ma Q}_p(\zeta_{p^n})$, es decir, $\lambda=\zeta_{p^n}^j-1$, $0\leq j\leq p^n-1$.

Sea $a=up^m\in \*{{\ma Q}_p}$, $u$ una unidad y $\zeta_{p^n}$ una $p^n$--ra\'iz
primitiva de la unidad. Entonces
\[
(a,{\ma Q}_p(\zeta_{p^n})/{\ma Q}_p)(\zeta_{p^n})=\zeta_{p^n}^r
\]
donde $r\in{\ma N}$ tal que $r\equiv u^{-1}\bmod p^n$.

En efecto, si $\lambda=\zeta_{p^n}-1\in \Lambda_{f,n}$, se tiene que
$ru\equiv 1\bmod p^n$ y por los Teoremas \ref{CClaseT3.2.5.29} y \ref{CClaseT3.2.5.11} se
tiene 
\[
(a,{\ma Q}_p(\zeta_{p^n})/{\ma Q}_p)(\lambda)=(u^{-1})_f(\lambda)=r_f(\lambda).
\]

Por otro lado si $g(Z)=(1+Z)^r-1$ se tiene
\begin{align*}
g(Z)&=(1+Z)^r-1=rZ+\cdots+Z^r\equiv rZ\bmod \deg 2\quad\text{y}\\
f(g(Z))&=f((1+Z)^r-1)=((1+Z)^r)^p-1=(1+Z)^{rp}-1\\
&=(1+Z)^{pr}-1=((1+Z)^p)^r-1
=g((1+Z)^p-1)=g(f(Z)).
\end{align*}

Por tanto $g(Z)$ satisface las condiciones que definen $r_f(Z)$, esto es $r_f(Z)=(1+Z)^r-1$.

Se sigue que
\begin{align*}
(a,{\ma Q}_p(\zeta_{p^n})/{\ma Q}_p)(\zeta_{p^n})&=(a,{\ma Q}_p(\zeta_{p^n})/{\ma Q}_p)
(\lambda+1)=r_f(\lambda)+1\\
&=r_f(\zeta_{p^n}-1)+1=(\zeta_{p^n}-1+1)^r-1+1=\zeta_{p^n}^r.
\end{align*}
\end{ejemplo}

Por ser el origen de lo presentado en esta secci\'on, enunciamos el Ejemplo \ref{CClaseEj3.2.5.30}
como un teorema.

\begin{teorema}\label{CClaseT3.2.5.31} Sean $K={\ma Q}_p$, $a=up^m\in\*{{\ma Q}_p}$ con $u$
una unidad, $\zeta_{p^n}$ una ra\'iz $p^n$--primitiva de $1$. Entonces el s\'imbolo
de la norma residual de ${\ma Q}_p(\zeta_{p^n})/{\ma Q}_p$ est\'a dado por
\[
(a,{\ma Q}_p(\zeta_{p^n})/{\ma Q}_p)(\zeta_{p^n})=\zeta_{p^n}^r
\]
donde $r\equiv u^{-1}\bmod p^n$, $r\in{\ma N}$. $\fin$
\end{teorema}

\begin{observacion}\label{O17.5.46'+1}
Usando el teorema de existencia podemos dar otra demostraci\'on
del Teorema \ref{T17.5.46'}. Sea $H=\N_{L/K}\*L$. Por el teorema
de existencia, existe una extensi\'on abeliana finita $M/K$ tal que
$\N_{M/K}\*M=H$. Por el Teorema \ref{CClaseT3.2.21'+1}, el grupo
de normas correspondiente a la extensi\'on abeliana $LM/L$ es
$\N^{-1}_{L/K} H=\N^{-1}_{L/K} \N_{L/K} \*L=\*L$. Por tanto 
$\N_{LM/L}\*{(LM)}=\*L=\N_{L/L}\*L$ de donde se sigue que 
$LM=L$ lo cual equivale a que $M\subseteq L$. Finalmente
si $N/K$ es abeliana y $K\subseteq M\subseteq N\subseteq L$,
entonces $H=\N_{L/K}\*L\subseteq \N_{N/K}\*N\subseteq
\N_{M/K}\*M=H$, de donde se sique que $N=M$ y por tanto
$M$ es la m\'axima extensi\'on abeliana de $K$ contenida en $L$.
\end{observacion}

Para finalizar la teor{\'\i}a local de campos de clase, \'unicamente
resta probar que todo subgrupo abierto de {\'\i}ndice
finito en $\*K$ es un subgrupo de normas.

El Teorema \ref{CCLT17.6.11} prueba que los grupos de las extensiones no ramificadas
son $(\pi^f)\times U_K$. Para las extensiones abelianas totalmente ramificadas
se tiene

\begin{teorema}\label{CClaseT3.2.5.32} Los grupos de normas de las extensiones
abelianas finitas $L/K$ totalmente ramificadas son exactamente los grupos
que contienen a alg\'un grupo
\[
(\pi)\times \unidades n
\]
con $\pi$ un elemento primo de $K$ y $n\in {\ma N}\cup \{0\}$.
\end{teorema}

\begin{proof} Como el grupo de normas de $L_{\pi,n}$ es $(\pi)\times
 \unidades n$, un
subgrupo que contiene a $(\pi)\times \unidades n$ 
corresponde a un subcampo
$L$ con $K\subseteq L\subseteq L_{\pi,n}$ y por tanto $L/K$ 
es abeliana y totalmente ramificada.

Ahora sea $L/K$ totalmente ramificada.
Entonces $L=K(\lambda)$ con $\lambda$ una ra\'iz de un polinomio Eisenstein
\[
X^e+\cdots+\pi=0,
\]
donde $\pi$ es un elemento primo de $K$ el cual es la norma del elemento
$\pm \lambda$. Se sigue que $(\pi)\subseteq \N_{L/K}\*L$. Ahora bien, puesto que
$\N_{L/K}(\*L)$ es un subgrupo abierto en $\*K$, existe $n\in{\ma N}\cup\{0\}$ tal
que $\unidades n\subseteq \N_{L/K}\*L$. Se sigue que 
$(\pi)\times\unidades n\subseteq \N_{L/K} \*L$. $\fin$
\end{proof}

\begin{corolario}\label{CClaseC3.2.5.33} Si $L/K$ es una 
extensi\'on abeliana finita totalmente
ramificada, entonces $L/K$ est\'a contenida en alg\'un $L_{\pi,n}$. 
$\fin$
\end{corolario}

\begin{teorema}\label{CClaseT3.2.5.34} 
El grupo $(\pi^f)\times \unidades n$ es el grupo de normas
del campo $K_f L_{\pi,n}$ donde $K_f/K$ es la extensi\'on no ramificada de $K$ de
grado $f$.
\end{teorema}

\begin{proof} Se tiene 
\begin{align*}
(\pi^f)\times \unidades n&=((\pi^f)\times U_K)\cap ((\pi)\times \unidades n)=
(\N_{K_f/K}\*{K_f})\cap (\N_{L_{\pi,n}/K}\*{L_{\pi,n}})\\
&=\N_{K_fL_{\pi,n}/K} \*{(K_f L_{\pi,n})}
\end{align*}
por los Teorema \ref{CCLT17.6.11} y \ref{CClaseT3.2.29}. $\fin$
\end{proof}

\begin{teorema}\label{CClaseT3.2.23''}
Los grupos de normas de $\* K$ son precisamente los
grupos conteniendo a alg\'un $(\pi^f)\times \unidades n$
con $n=0,1,2,\ldots, f=1,2,\ldots$, donde $\pi$ es un elemento
primo de $K$.
\end{teorema}

\begin{proof} Cada grupo $(\pi^f)\times \unidades n$ tiene \'indice finito
en $\*K=(\pi)\times U_K$ pues
\[
\frac{\*K}{(\pi^f)\times \unidades n}\cong\frac{(\pi)\times U_K}{(\pi^f)\times
\unidades n}\cong\frac{(\pi)}{(\pi^f)}\times \frac{U_K}{\unidades n}
\]
el cual es cardinalidad $f\cdot q\cdot(q^{n-1}-1)$. Adem\'as $(\pi^f)\times
\unidades n$ es un grupo abierto en $\*K$. Por otro lado, por el Teorema
\ref{CClaseT3.2.5.34}, $(\pi^f)\times\unidades n$ es un grupo de normas.

Por tanto, si $H$ es un subgrupo que contiene a uno de estos grupos,
$H$ es un grupo de normas (ver la demostraci\'on del Teorema \ref{CClaseT3.2.29}).

Rec\'iprocamente, si $H$ es un grupo de normas, entonces $H$ es un
subgrupo abierto de $\*K$ de \'indice finito. Entonces $1\in H$ y $\{
\unidades n\}_{n=0}^{\infty}$ es un sistema fundamental de vecindades
de $1$, existe $n\in{\ma N}\cup\{0\}$ con $\unidades n\subseteq H$.
Finalmente, al ser $H$ de \'indice finito en $\*K$, existe $f$ tal que
$\pi^f\in H$. Se sigue que $(\pi^f)\times \unidades n\subseteq H$. 
$\fin$
\end{proof}

\begin{corolario}\label{CClaseC3.2.5.35-1}
Se tiene que
\[
D_K:=\bigcap_{\substack{L/K\\ \text{finita y abeliana}}}\N_{L/K}\*L=
\{1\}\]
y $\rho_K\colon\*K\lra \abe G_K$ es un monomorfismo.
\end{corolario}

\begin{proof}
Se tiene $\{1\}\subseteq D_K\subseteq \bigcap_{f,n}\big(\langle\pi_K^f
\rangle\times \unidades n\big)=\{1\}$.
$\fin$
\end{proof}

\begin{teorema}\label{CClaseT3.2.5.35} 
Sea $T$ la m\'axima extensi\'on no ramificada de $K$. 
Sean $\pi$ un elemento primo de $K$,
$f\in{\mc F}_{\pi}$, $\Lambda_f=\bigcup_{n=1}^{\infty}\Lambda_{f,n}$, 
$L_{\pi}=\bigcup_{n=1}^{\infty}L_{\pi,n}=K(\Lambda_f)$ y $G_{\pi}=
\Gal(L_{\pi}/K)$. Entonces
el campo $TL_{\pi}$ es independiente de $\pi$ y es la m\'axima extensi\'on
abeliana de $K$, $\abe K=TL_{\pi}$. En particular
\[
\abe{G_K}=\Gal(\abe K/K)=G_{T/K}\times G_{\pi}\cong \hat{{\ma Z}}\times U_K.
\]

Si $a=u\pi^m\in \*K$ con $u\in U_K$, entonces el s\'imbolo de la norma residual
$(a,K_=\rho_K(a)$ est\'a dado por
\[
(a,K)|_T=\tau^m, \quad (a,K)(\lambda)=(u^{-1})_f(\lambda)\quad\text{para $\lambda
\in\Lambda_f$}.
\]
\end{teorema}

\begin{proof} Por el Teorema \ref{CClaseT3.2.23''} los grupos de normas son los grupos conteniendo
a alg\'un $(\pi^f)\times \unidades n$, el cual es, por el Teorema \ref{CClaseT3.2.5.34},
el grupo de normas de $L_{\pi,n} K_f$.

As\'i, si $L/K$ es una extensi\'on abeliana finita, existen $n\in{\ma N}\cup\{0\}$,
$f\in{\ma N}$ tales que $(\pi^f)\times \unidades n\subseteq 
\N_{L/K}\*L$ lo cual
implica que $\*L\subseteq L_{\pi,n}K_f$ con $K_f\subseteq T$. Se sigue que
$\* L\subseteq L_{\pi}T$ y $L_{\pi}T$ es la m\'axima 
extensi\'on abeliana de $K$.

Puesto que $\abe{G_K}=G_{T/K}\times G_{\pi}$, el s\'imbolo de la norma residual
est\'a determinad por
\[
(a,K)|_T=\tau^m\quad\text{y}\quad (a,K)(\lambda)=(u^{-1})_f(\lambda)
\quad \text{para $\lambda\in\Lambda_f$}.
\]

Finalmente, por el Teorema \ref{CClaseT3.2.5.19} tenemos que $G_{\pi}
\cong U_K$. $\fin$
\end{proof}

\begin{corolario}[Teorema de Kronecker--Weber local\index{teorema
de Kronecker--Weber local}\index{Kronecker--Weber!teorema local
de $\sim$}]\label{CClaseC3.2.5.36}
La m\'axima extensi\'on abeliana de ${\ma Q}_p$ es
\[
\abe{{\ma Q}_p}=\bigcup_{n=1}^{\infty}{\ma Q}_p(\zeta_n)={\ma Q}_p(\zeta_{\infty}).
\]
\end{corolario}

\begin{proof} Las extensiones abelianas no ramificadas de ${\ma Q}_p$ corresponden
a las extensiones finitas de ${\ma F}_p$ y estas con ${\ma F}_p(\zeta_n)$ con
$\mcd(n,p)=1$. Esto es
\[
T=\bigcup_{\substack{n=1\\ \mcd(n,p)=1}}^{\infty}{\ma Q}_p(\zeta_n).
\]
Por ora lado $L_{p,n}={\ma Q}(\zeta_{p^n})$ por lo cual
$L_{\pi}=\bigcup_{n=1}^{\infty}{\ma Q}_p(\zeta_{p^n})={\ma Q}_p(\zeta_{p^{\infty}})$
de donde se sigue el resultado. $\fin$
\end{proof}

\begin{teorema}[Teorema de Existencia\index{teorema de 
existencia}]\label{CClaseT3.2.30}
Sea $K$ un campo local. Entonces la correspondencia
$L\longrightarrow H_L:=\N_{L/K}\*L\subseteq \*K$ nos da un
isomorfismo de redes que cambia contenciones entre la
red de extensiones abelianas finitas $L/K$ y la red de subgrupos
abiertos de {\'\i}ndice finito de $\*K$. Todo subgrupo que
contenga a un grupo de normas es a su vez un grupo de normas.

En particular, si $K$ es un campo local, entonces las siguientes condiciones 
sobre un subgrupo $H$ de $\*K$, son equivalentes:
\las
\item $H$ es abierto y de \'indice finito en $\*K$.
\item $H$ es un grupo de normas.
\item $H$ contiene a alg\'un grupo $\langle\pi_K^f\rangle\times
\unidades n$ para algunos $f\in{\ma N}$ y $n\in{\ma N}\cup\{0\}$.
\end{list}
\end{teorema}

\begin{observacion}\label{CClaseO3.2.31} El Teorema \ref{CClaseT3.2.30}
es de existencia pues dado una subgrupo abierto de {\'\i}ndice finito
$V$ de $\*K$, existe una \'unica extensi\'on abeliana finita
$L/K$ tal que $\*K/\N_{L/K}\*L\xrightarrow[({\underline{\ \ }},L/K)]{\cong}
\Gal(L/K)$, esto es, $\N_{L/K}\*L=V$.
\end{observacion}

\begin{proof}(Teorema \ref{CClaseT3.2.30}). 
Teoremas \ref{CClaseT3.2.29}, \ref{CClaseT3.2.5.34} y \ref{CClaseT3.2.23''}. 
$\fin$
\end{proof}

\begin{observacion}\label{CClaseO3.2.32} Para otra demostraci\'on del
Teorema de Existencia en caracter{\'\i}stica $0$ se puede consultar
Teorema \ref{CClaseT3.2.21'+4}, 
\cite[Theorem 3.1, Ch. III Section 3, p\'agina 43]{Neu86}; \cite[Theorem
6.2, Ch. II, Section 6, p\'agina]{Neu69}. Para caracter{\'\i}stica $p>0$,
se puede consultar \cite[Ch. XIV, Section 6, p\'agina 218]{Ser}.
\end{observacion}

\begin{observacion}\label{CClaseO3.3.33}
Notemos que $L/K$ es no ramificada $\iff U_K\subseteq \N_{L|K}\*L$.
En efecto, si $L/K$ es no ramificada, entonces $\N_{L/K}\*L=\langle
\pi_K^f\rangle \times U_K$ con $f=[L:K]$. Rec\'iprocamente, si
$U_K\subseteq \N_{L/K}\*L$ y si $f=[\*K\colon \N_{L/K}\*L]$,
entonces $\pi_K^f\in\N_{L/K}\*L$ y por tanto $\langle\pi_K^f\rangle
\times U_K\subseteq \N_{L/K}\*L$. Por tanto $L$ est\'a contenida
en la extensi\'on no ramificada de $K$ de grado $f$ (y de hecho,
es igual).
\end{observacion}

\subsubsection{Grupos de ramificaci\'on, grupos
de descomposici\'on y grupos formales}\label{CClaseS3.6}

Todos los resultados presentados aqu\'i sin demostraci\'on, as\'i como
la teor\'ia general de grupos de ramificaci\'on pueden ser consultados
en \cite[Ch. IV]{Ser}.

En esta parte desarrollamos lo expuesto en la Subsecci\'on \ref{CClaseS3.2.4}.
Sea $L/K$ una extensi\'on finita de Galois con grupo $G=\Gal(L/K)$.
Recordemos que $\varphi\colon[-1,\infty)\lra [-1,\infty)$ dada por 
$\varphi(u)=\int_0^u \frac{dt}{[G:G_t]}$ donde $G_t:=G_{\lceil t\rceil}$ es
una funci\'on biyectiva y sea $\eta$ la funci\'on inversa. Entonces se
define para $v\in[-1,\infty)$
\[
G^v:=G_{\eta(v)}\quad\text{o}\quad G^{\varphi(u)}:=G_u.
\]

\begin{proposicion}\label{CClaseP3.2.5.6.1} Sea $K$ un campo local y sea
$L_{\pi,n}=K(\Lambda_{f,n})$ el campo de los puntos de $\pi^n$--divisi\'on
de un m\'odulo de Lubin--Tate para $\pi$. Entonces
\[
G_i(L_{\pi,n}/K)=\Gal(L_{\pi,n}/L_{\pi,m})\quad\text{para $q^{m-1}\leq i
\leq q^m-1$}.
\]
\end{proposicion}

\begin{proof} Se tiene que el s\'imbolo de la norma residual proporciona un
isomorfismo $U_K/\unidades n\lra \Gal(L_{\pi,n}/K)$ (Teorema \ref{CClaseT3.2.5.15}).
Por tanto se tiene
\begin{gather*}
\Gal(L_{\pi,n}/L_{\pi,m})\cong \frac{\Gal(L_{\pi,n}/K)}{\Gal(L_{\pi,m}/K)}\cong
\frac{U_K/\unidades n}{U_K/\unidades m}\cong \frac{\unidades m}{\unidades n},
\intertext{y por tanto}
\Gal(L_{\pi,n}/L_{\pi,m})=(\unidades m,L_{\pi,n}/K).
\end{gather*}

Lo que se quiere probar es 
que para $q^{m-1}\leq i\leq q^m-1$, $G_i(L_{\pi,n}/K)
=(\unidades m,L_{\pi,n}/K)$.

Sea $\sigma\in G_1(L_{\pi,n}/K)$ y sea $\sigma=(u^{-1},L_{\pi,n}/K)$ para
alg\'un $u$. El mapeo $(\underline{\ \ },L_{\pi,n}/K)\colon U_K/\unidades n
\xrightarrow{\ \cong\ }\Gal(L_{\pi,n}/K)$ manda el $p$--subgrupo de Sylow
$\unidades 1/\unidades n$ de $U_K/\unidades n$ (recordemos que $\big|
U_K/\unidades n\big|=q^{n-1}(q-1)$ y $\big|\unidades 1/\unidades n\big|=q^{n-1}$)
sobre el $p$--subgrupo de Sylow de $\Gal(L_{\pi,n}/K)$, el cual es
$G_1(L_{\pi,n}/K)$.
En particular $u\in \unidades 1$.

Escribamos $u=1+\varepsilon \pi^l$ con $\varepsilon\in U_K$ y $l\geq 1$.
Sea $\lambda\in \Lambda_{f,n}$ un generador como $\o_K$--m\'odulos, esto es,
$\lambda\in\Lambda_{f,n}\setminus\Lambda_{f,n-1}$. Entonces
\begin{align*}
\lambda^{\sigma}&=(u^{-1},L_{\pi,n}/K)(\lambda)\igual_{\substack{\uparrow\\
\text{Teorema \ref{CClaseT3.2.5.35}}}} (u)_f(\lambda)=(1\lubintatemas {F_f}\varepsilon \pi^l)_f
(\lambda)\\
&=1=F_f(\lambda,[\varepsilon\pi^l]_f(\lambda)).
\end{align*}

Si $l\geq n$, $\sigma =1$ por lo que $\lambda^{\sigma}-\lambda=0$ y $v_{L_{\pi,n}}
(\lambda)=\infty$. Si $l<n$ entonces $\lambda_{n-l}:=(\pi^m)_f(\lambda)$ es un
elemento primo $L_{\pi,n-l}$.

Se tiene $(\varepsilon \pi^l)_f(\lambda)=(\varepsilon)_f(\pi^l)_f(\lambda)=
(\varepsilon)_f(\lambda_{n-l})$.

Puesto que $L_{\pi,n}/L_{\pi,n-l}$ es totalmente ramificada de grado $q^l$, se tiene
que $v_{L_{\pi,n}}(\lambda^{q^l})=q^lv_{L_{\pi,n}}(\lambda)=q^l\cdot 1=q^l$ y
$v_{L_{\pi,n}}(\lambda_{n-l})=e(L_n/L_{n-l})v_{L_{\pi,n}}(\lambda_{n-l})=q^l\cdot
1=q^l$, por tanto existe $\varepsilon_0\in U_{L_n}$ tal que
\[
[\varepsilon]_{F_f}(\lambda_{n-l})=[\varepsilon\pi^l]_{F_f}(\lambda)=\varepsilon_0
\lambda^{q^l}.
\]

De las identidades $F_f(X,0)=X$, $F_f(0,Y)=Y$ (ver despu\'es de la Definici\'on 
\ref{CClaseD3.2.2.5'}), se obtiene que
\begin{gather*}
F_f(X,Y)=X+Y+XYG(X,Y)\quad\text{con $G(X,Y)\in\o)K[[X,Y]]$}.
\intertext{De esta forma obtenemos}
\lambda^{\sigma}-\lambda=F_f(\lambda,\varepsilon_0\lambda^{q^l})-\lambda=
\varepsilon_0\lambda^{q^l}+a\lambda^{q^l+1}\quad \text{con $a\in\o_{L_{\pi,n}}$}.
\intertext{Por lo tanto}
\imath_{L_{\pi,n}/K}(\sigma)=v_{L_{\pi,n}}(\lambda^{\sigma}-\lambda)=
\begin{cases}
q^l&\text{si $l<n$}\\
\infty&\text{si $l\geq n$}.
\end{cases}
\end{gather*}

Consideremos $q^{m-1}\leq i\leq q^m-1$ y $u\in\unidades m$. Entonces
$l\geq m$, esto es, $\imath_{L_{\pi,n}/K}(\sigma)\geq q^l\geq i+1$ y de
esta forma $\sigma\in G_i(L_{\pi,n}/K)$.

Esto demuestra que $(\unidades m, L_{\pi,n}/K)\subseteq G_i(L_{\pi,n}/K)$.

Rec\'iprocamente, si $\sigma\in G_i(L_{\pi,n}/K)$ y $\sigma\neq \Id$, entonces
$\imath_{L_{\pi,n}/K}(\sigma)=q^m>i\geq q^{m-1}$, esto es, $l\geq m$ por lo 
que $u\in \unidades m$ lo que demuestra que $G_i(L_{\pi,n}/K)\subseteq
(\unidades m,L_{\pi,n}/K)$. $\fin$
\end{proof}

Un resultado interesante es ver que tipo de saltos superiores
tienen las extensiones abelianas.

\begin{teorema}[Hasse--Arf\index{teorema de 
Hasse--Arf}\index{Hasse--Arf!teorema de $\sim$}]\label{CClaseT3.2.27}
Sea $L/K$ una extensi\'on abeliana finita de campos locales.
Entonces los saltos de la filtraci\'on superior $\{G^t(L/K)\}_{t\geq -1}$ son enteros.
\end{teorema}

\begin{proof} Sea $E/K$ la m\'axima subextensi\'on $K\subseteq E\subseteq L$
no ramificada. Entonces $G^t(L/E)=G^t(L/K)$ para $t>-1$ debido
a que $\eta_{E/K}(s)=s$ donde $\eta_{E/K}$ es la inversa de la funci\'on
de Herbrand. Por tanto
\[
\eta_{L/K}(s)=\eta_{L/E}(\eta_{E/K}(s))=\eta_{L/E}(s).
\]

Se sigue que podemos suponer que $L/K$ es totalmente ramificada.

Sea $\pi_L$ un elemento primo de $L$. Entonces $\pi:=\N_{L/K}(\pi_L)$
es un elemento primo de $K$. Por otro lado, existe $m\in{\ma N}$ tal que
$(\pi)\times \unidades m\subseteq \N_{L/K}\*L$. Se sigue que $L$
est\'a contenida en el campo de clase de $(\pi)\times \unidades m$ el cual
corresponde a $L_{\pi,m}$. Si $t_0$ es un salto de $\{G^t(L/K)\}_t$ entonces,
por el Teorema \ref{CClaseGR.10}, $t_0$ es un salto de
$\{G^t(L_{\pi,m}/K)_t$. Por lo tanto podemos suponer $L=L_{\pi,m}$.

De la Proposici\'on \ref{CClaseP3.2.5.6.1} los saltos de $\{G_s(L_{\pi,m}/K)\}_s$
son los n\'umeros $q^l-1$, $0\leq l\leq m-1$ con la excepci\'on de
que cuando $q=2$, $0$ no es un salto.

Para calcular los saltos $\{G^t(L_{\pi,m}/K)\}_t$ calculamos
$\varphi_{L_{\pi,m}/K}(q^l-1)=l$, $l=0,1,\ldots, m-1$ lo cual prueba el 
teorema. $\fin$
\end{proof}

El siguiente es el resultado central de los grupos 
de ramificaci\'on superior
con la teor\'ia de clase de campos locales.

\begin{teorema}\label{CClaseT3.2.5.6.3} Sea $L/K$ 
una extensi\'on abeliana finita de campos
locales. Entonces el s\'imbolo residual de la norma
\begin{gather*}
(\ \ ,L/K)\colon \* K\lra \Gal(L/K)
\intertext{manda el grupo $\unidades n$, $n\geq 0$, sobre el 
$n$--\'esimo grupo de ramificaci\'on 
superior $G^n(L/K)$:}
(\unidades n,L/K)=G^n(L/K), n\geq 0.
\end{gather*}
\end{teorema}

\begin{proof} Si $E$ es la m\'axima extensi\'on no ramificada de $K$ contenida en
$L$: $K\subseteq E\subseteq L$, entonces $G^n(L/K)=G^n(L/E)$. Por
el Teorema \ref{CCLT17.6.3} y el Teorema \ref{CClaseT3.2.21'} se tiene
\[
(U_E^{(n)},L/E)=(\N_{E/K} U_E^{(n)},L/K)=(\unidades n, L/K).
\]

Por tanto, podemos, sustituir $L/K$ por $L/E$ y por ende podemos suponer
que $L/K$ es totalmente ramificada. Procediendo como en el Teorema de
Hasse--Arf, Teorema \ref{CClaseT3.2.27}, se tiene que $L\subseteq L_{\pi,m}$
para alguna $m\in{\ma N}$ y sin p\'erdida de generalidad podemos suponer
$L=L_{\pi,m}$. 

Ahora, por el Teorema \ref{CClaseT3.2.5.32} y la Proposici\'on 
\ref{CClaseP3.2.5.6.1}, se tiene 
\[
\Gal(L_{\pi,m}/L_{\pi,l})=G_i(L_{\pi,m}/K)
\]
para $q^{l-1}\leq i\leq q^l-1$. Puesto que la funci\'on de Herbrand $\varphi$
satisface $\varphi_{L_{\pi,m}/K}(q^l-1)=l$, obtenemos
\begin{gather*}
(\unidades l,L_{\pi,m}/K)=G_{q^l-1}(L_{\pi,m}/K)=G^l(L_{\pi,m}/K).
\tag*{$\fin$}
\end{gather*}
\end{proof}

El Teorema \ref{CClaseT3.2.5.6.3} nos proporciona una forma de calcular el
conductor local.

\begin{corolario}\label{CClaseC3.2.5.6.4} Sea $L/K$ una extensi\'on abeliana
finita de campos locales y sea $\f{L/K}=\f{}$ el conductor de la
extensi\'on $L/K$. Sea $n$
tal que $G_n(L/K)\neq \{1\}$ y $G_{n+1}(L/K)=\{1\}$, $n\geq 1$. Entonces
$\f{}=\pK^{c_{\pK}}$ donde 
\[
c_{\pK}=1+\varphi_{L/K}(n)=\frac{1}{g_0}(g_0+g_1+\cdots+g_n).
\]

En otras palabras, el conductor local pudo haber sido definido
como 
\[
\f{}=\pK^{1+\varphi_{L/K}(n)}
\]
 donde $G_n(L/K)\neq \{1\}$ y $G_{n+1}(L/K)=\{1\}$.
\end{corolario}

\begin{proof} Por el Teorema \ref{CClaseT3.2.5.6.3} se tiene
\begin{align*}
\unidades m\subseteq \N_{L/K}\*L&\iff (\unidades m,L/K)=1
\iff G^m(L/K)=\{1\}\\
&\iff G_{\eta_{L/K}(m)}(L/K)=\{1\}.
\end{align*}

Se sigue que $G_n(L/K)\neq \{1\}$ y $G_{n+1}(L/K)=\{1\}$ si y solamente si
$\unidades {\varphi_{L/K}(n)}\nsubseteq \N_{L/K}\*L$ y $\unidades {
1+\varphi_{L/K}(n)}\subseteq \N_{L/K}\*L$. $\fin$
\end{proof}

Terminamos esta secci\'on estudiando la descomposici\'on en extensiones
abelianas de campos locales.

\begin{teorema}\label{T17.5.60}
Sea $L/K$ una extensi\'on abeliana finita de campos locales. Sean $e$
y $f$ el \'indice de ramificaci\'on y el grado de inercia. Entonces
\[
e=[U_K:\N_{L/K} U_L]\quad \text{y}\quad f=o(\pi_K\bmod U_K \N_{L/K}
\*L),
\]
donde $\pi_K$ es un elemento uniformizador de $K$.
\end{teorema}

\begin{proof}
El Corolario \ref{CCLP17.6.24+1} prueba que $e=[U_K:\N_{L/K} U_L]$. Otra
demostraci\'on es consecuencia del Teorema \ref{CClaseT3.2.5.6.3}: 
\[
(U_K,L/K)=G_0(L/K)\cong \frac{U_K\N_{L/K}\*L}{\N_{L/K}\*L}\cong
\frac{U_K}{U_K\cap \N_{L/K}\*L}=\frac{U_K}{\N_{L/K}U_L}.
\]
Ahora
bien, $ef=[L:K]=[\*K:\N_{L/K}\*L]$. Notemos que
\[
\frac{U_K\N_{L/K}\*L}{\N_{L/K}\*L}\cong \frac{U_K}{\N_{L/K}\*L\cap
U_K}=\frac{U_K}{\N_{L/K} U_L}.
\]
Se sigue que $f=\frac{[L:K]}{e}=\frac{[\*K:\N_{L/K}\*L]}{[U_K\N_{L/K}\*L:
\N_{L/K}\*L]}=[\*K:U_K\N_{L/K}\*L]$. Puesto que 
$\*K=(\pi_K) \times U_K$, obtenemos
\[
\frac{\*K}{U_K\N_{L/K}\*L}=\frac{(\pi_K) U_K\N_{L/K}\*L}{U_K\N_{L/K}\*L}
=\langle\pi_K \bmod U_K\N_{L/K}\*L\rangle.
\]
Se sigue $f=[\*K:U_K\N_{L/K}\*L]=o(\pi_K\bmod U_K\N_{L/K}\*L)$. $\fin$
\end{proof}

\section{Campos de clase globales}\label{CClaseC4}

\begin{definicion}\label{D17.6.1N}
Un {\em campo global\index{campo global}} $K$ es o bien una extensi\'on
finita de ${\ma Q}$ o bien un campo de funciones con campo de
constantes $\F$, el campo finito de $q=p^r$ elementos.

Las extensi\'on finitas de ${\ma Q}$ son el caso de caracter\'istica $0$
y los llamaremos {\em campos num\'ericos\index{campos num\'ericos}}.
Los campos de funciones son el caso de caracter\'istica $p>0$ y los
llamaremos {\em campos de funciones\index{campos de funciones}}.
\end{definicion}

\begin{observacion}\label{O17.6.2N}
En el caso de campos de funciones {\underline{no}} hay lugares
infinitos.
\end{observacion}

Sea $K$ un campo global. Sea $\pK$ un lugar de $K$. Usaremos
la notaci\'on $\pK|\infty$ para $\pK$ un lugar infinito (arquimediano)
y $\pK\nmid \infty$ si $\pK$ es un lugar finito.

Si $\pK$ es un lugar finito, sea $\o_{\pK}=\{x\in K\mid v_{\pK}(x)
\geq 0\}$ el anillo de valuaci\'on de $\pK$ y $v_{\pK}$ es la valuaci\'on
asociada. La completaci\'on de $K$ en un lugar finito o infinito se
denotar\'a por $K_{\pK}$.

Si $\pK$ es finito, $\N(\pK)$ va a denotar la {\em norma
absoluta\index{norma absoluta}} de $\pK$, esto es $\N(\pK)=
\big|\o_{\pK}/\pK\big|$.

Se tiene que si ${\ma F}_p$ es el campo primo de $\o_{\pK}/
\pK$ y $f_{\pK}=\big[\o_{\pK}/{\pK}:{\ma F}_p\big]$, entonces
$\N(\pK)=p^{f_{\pK}}=q_{\pK}$.

Para el valor absoluto $|\ |_{\pK}$ asociado a $\pK$, seleccionamos
el n\'umero $0<c<1$ como $c=p^{-f_{\pK}}=\N(\pK)^{-1}=
\frac{1}{\big|\o_{\pK}/\pK\big|}=q_{\pK}^{-1}$. De esta foma, tenemos
que si $a\in K$, $a\neq 0$, 
\[
|a|_{\pK}=p^{-f_{\pK}v_{\pK}(a)}=
q_{\pK}^{-v_{\pK}(a)}\quad\text{y}\quad |0|_{\pK}=0.
\]

Si $\pK$ es un lugar finito, denotamos por $U_{\pK}=U_{K_{\pK}}$
al grupo de unidades de $K_{\pK}$. Si $\pK$ es un lugar
infinito, entonces $K_{\pK}={\ma R}$ o ${\ma C}$. Si $\pK$ es
real, es decir si $K_{\pK}={\ma R}$ y si $\sigma\colon K\to 
{\ma R}=K_{\pK}$ es el encaje asociado a $\pK$, entonces
$|a|_{\pK}:=|\sigma a|$, $a\in K$, donde $|\ |$ denota el valor 
absoluto usual de ${\ma R}$.

Si $\pK$ es complejo, $K_{\pK}={\ma C}$, si $\sigma$ es uno de
los dos encajes conjugados $\sigma\colon K\to {\ma C}=K_{\pK}$,
ponemos $|a|_{\pK}:=\|\sigma a\|^2$ para $a\in K$, donde $\|\ \|$
es el valor absoluto usual de ${\ma C}$.

\subsection{Repaso de resultados b\'asicos de campos globales}\label{S17.6.2N}

\subsubsection{Composici\'on de campos}\label{S17.6.1N}

Sean $E$ y $L$ dos campos cualesquiera de la misma 
caracter\'istica, esto es, el campo primo, ${\ma F}_p$
o ${\ma Q}$, es el mismo para ambos campos. Nos preguntamos,
?`qu\'e significa la composici\'on de $E$ y $L$?

Supongamos primero que existe un campo $\Omega$ tal que
$E,L\subseteq \Omega$. Entonces la composici\'on es
simplemente el m\'inimo subcampo de $\Omega$ que contiene
tanto a $E$ como a $L$. Sin embargo, si no tenemos un tal
campo $\Omega$, debemos construir un campo (m\'inimo)
que contenga tanto a $E$ como a $L$.

\begin{definicion}\label{D17.6.3N} Sea $K$ un campo arbitrario
y sean $E/K$ y $L/K$ dos extensiones de $K$ (por ejemplo, $K$
podr\'ia ser el campo primo). Una {\em 
composici\'on\index{composici\'on de campos}} de $E$
y $L$ sobre $K$, es una terna
$(M,\tau,\sigma)$ donde $M$ es un campo que contiene a $K$,
$\tau\colon E\to M$ y $\sigma\colon L\to M$, $\tau$ y $\sigma$
son monomorfismos de campos tales que $\tau|_K=\sigma|_K=
\Id_K$ y $M$ est\'a generado por $\tau(E)$ y $\sigma(L)$.
\end{definicion}

\begin{definicion}\label{D17.6.4N} Dos composiciones $(M,
\tau, \sigma)$, $(M',\tau',\sigma')$ de $E/K$ y $L/K$ se llaman
{\em equivalentes\index{composiciones equivalentes}} si
existe un isomorfismo $\lambda\colon M\lra M'$ tal que
$\lambda\circ \tau=\tau'$ y $\lambda\circ \sigma=\sigma'$.
\[
\xymatrix{
&E\ar@{->}[ld]_{\tau}\ar@{->}[rd]^{\tau'}\ar@{}[d]|{\circlearrowleft}\\
M\ar@{->}[rr]_{\lambda}&& M'
}\qquad
\xymatrix{
&L\ar@{->}[ld]_{\sigma}\ar@{->}[rd]^{\sigma'}\ar@{}[d]|{\circlearrowleft}\\
M\ar@{->}[rr]_{\lambda}&& M'
}
\]
\end{definicion}

Ser equivalentes, es una relaci\'on de equivalencia. Estudiaremos
en el caso en que $L/K$ es finita, $[L:K]<\infty$ y $E/K$ 
es arbitraria, las diferentes clases de equivalencia.

Sean $[L:K]=n<\infty$ y $(M,\tau,\sigma)$ una composici\'on
de $E$ y $L$. Sean $E'=\tau(E)$, $L'=\sigma(L)$ y 
\[
E'L'=\Big\{\sum_{i=1}^r e_il_i\mid e_i\in E', l_i\in L', r
\in {\ma N}\Big\}\subseteq M,
\]
entonces $E'L'$ es una sub\'algebra sobre $K$. Ahora bien,
$E'L'$ es un dominio entero. Sea $\{\alpha_1,\ldots,\alpha_n\}$
una base de $L/K$. Entonces $\{\sigma(\alpha_1),\ldots,\sigma(
\alpha_n)\}$ genera a $E'L'$ sobre $E'$. Puesto que $E'$ es
un campo, $E'L'$ es un campo (recordemos que si $F$ es un
campo y $D$ es un dominio entero con $\dim_F D<\infty$,
entonces $D$ es un campo). Se sigue que $E'L'=M$.

Sea $\theta\colon E\otimes_K L\lra M$ definido por $\theta(e
\otimes_K l)=\tau(e)\sigma(l)$. Entonces $\theta$ es un 
$K$-epimorfismo y $M\cong (E\otimes_K L)/\ker\theta$. Como
$M$ es un campo, $\ker \theta={\eu M}$ es un ideal maximal.
Ahora bien $K\cong K\otimes_K K$ y $\theta|_K=\Id_K$ por 
lo que ${\eu M}\cap K=\{0\}$. 

Los homomorfismos
\begin{gather*}
E\stackrel{i}{\lra} E\otimes_K L, \qquad i(e)=\tau(e)\otimes_K 1,\\
L\stackrel{j}{\lra} E\otimes_K L, \qquad j(l)=1\otimes_K \sigma(l),
\end{gather*}
son monomorfismos pues $(\theta\circ i)(E)=\tau(E)$ y $(\theta
\circ j)(L)=\sigma(L)$. Por tanto ${\eu M}\cap E={\eu M}\cap
L=\{0\}$.

\begin{teorema}\label{T17.6.5N} Las clases de equivalencia de
las composiciones $E$ con $L$ sobre $K$, est\'an en 
correspondencia biyectiva con los ideales maximales de la
$K$-\'algebra $E\otimes_K L$. En particular, la composici\'on
de campos siempre existe.
\end{teorema}

\begin{proof}
Ya vimos que cada composici\'on corresponde a un ideal maximal.
Reciprocamente, si ${\eu M}$ es un ideal maximal de $E\otimes_K
L$, definimos $M=(E\otimes_K L)/{\eu M}$ el cual es un campo.
Sean 
\begin{align*}
i\colon & E\lra (E\otimes_K L)/{\eu M}, \quad i(e)=(e\otimes_K 1)
+{\eu M}\quad \text{y}\\
j\colon & L\lra (E\otimes_K L)/{\eu M}, \quad j(l)=(1\otimes_K l)
+{\eu M}.
\end{align*}
Puesto que ${\eu M}$ no contiene unidades, $i$ y $j$ son inyectivas
y $M$ est\'a generado por $i(E)$ y $j(L)$. Adem\'as $i|_K=j|_K=
\Id_K$. Se sigue que $M$ es una composici\'on de $E$ y $L$.

Ahora sean $(M,\tau,\sigma)$ y $(M',\tau',\sigma')$ dos 
composiciones donde
\[
M\cong (E\otimes_K L)/{\eu M}\quad\text{y}\quad 
M'\cong (E\otimes_K L)/{\eu M}'.
\]
Si $M$ y $M'$ son
equivalentes, existe un isomorfismo $\lambda\colon M\to M'$
con $\lambda\circ\tau=\tau'$ y $\lambda\circ\sigma=\sigma'$.
Sea $\sum_{i=1}^r e_i\otimes_K l_i\in{\eu M}$. Entonces se tiene
que si $\sum_{i=1}^r\tau(e_i)\sigma(l_i)=0$ en $M$ entonces
\begin{gather*}
\lambda\Big(\sum_{i=1}^r\tau(e_i)\sigma(l_i)\Big)
=\sum_{i=1}^r
(\lambda\circ\tau)(e_i)(\lambda\circ \sigma)(l_i)
=\sum_{i=1}^r\tau'(e_i)\sigma'(l_i)=0
\end{gather*}
en $M'$ por lo que $\sum_{i=1}^r e_i\otimes_K l_i\in {\eu M}'$.
Se sigue que ${\eu M}\subseteq {\eu M}'$ y como ambos son 
maximales, se tiene ${\eu M}={\eu M}'$.

Rec\'iprocamente, sean ${\eu M}={\eu M}'$. Si $\theta$ y $\theta'$ son los
isomorfismos de $(E\otimes_K L)/{\eu M}$ y de $(E\otimes_K L)/
{\eu M}'$ a $M$ y $M'$ respectivamente, entonces
\[
\xymatrix{
(E\otimes_K L)/{\eu M}\ar@{=}[r]\ar@{->}^{\theta}[d]&
(E\otimes_K L)/{\eu M}'\ar@{->}[d]^{\theta'}\\
M\ar@{->}[r]_{\theta'\circ \theta^{-1}}& {\eu M}'
}
\]
y $\lambda:=\theta'\circ\theta^{-1}$ es un isomorfismo de $M$ 
y $M'$, y
\begin{align*}
\tau=\theta\circ i,&\quad \tau'=\theta'\circ i=\theta'\circ\theta^{-1}
\circ\theta\circ i=\lambda\circ \tau,\\
\sigma=\theta\circ j,&\quad \sigma'=\theta'\circ j=\theta'\circ\theta^{-1}
\circ\theta\circ j=\lambda\circ \sigma,
\end{align*}
por lo que $M$ y $M'$ son extensiones equivalentes.
$\fin$
\end{proof}

\begin{teorema}\label{T17.6.6N} Sean $T$ un campo y $A$ un
\'algebra sobre $T$ tal que $A$ es de dimemsi\'on finita sobre
$T$ y $A$ tiene unidad. Entonces $A$ tiene un n\'umero finito
de ideales maximales.
\end{teorema}

\begin{proof}
Sea $\dim_T A=n<\infty$ y sean ${\eu M}_1,\ldots, {\eu M}_r$ 
$r$ ideales maximales de $A$ distintos. Sea ${\eu N}=\bigcap_{
i=1}^r {\eu M}_i$. Entonces $A/{\eu N}\cong \bigoplus_{i=1}^r
A/{\eu M}_i$. Ahora bien, $A/{\eu N}$ y $A/{\eu M}_i$, $1\leq
i\leq r$, son $T$-\'algebras y
\[
n=\dim_T A\geq \dim_T A/{\eu N}=\sum_{i=1}^r A/{\eu M}_i\geq r,
\]
de donde se sigue que $r\leq n$.
$\fin$
\end{proof}

\begin{corolario}\label{C17.6.7N} Si $K$ es cualquier campo
y $E/K$ y $L/K$ son dos extensiones de $K$ con $[L:K]=n
<\infty$, entonces el n\'umero de composiciones de $E$ y 
$L$ sobre $K$ es menor o igual $n$ y en particular este
n\'umero es finito.
\end{corolario}

\begin{proof}
Si $\{\alpha_1,\ldots,\alpha_n\}$ es una base de $L/K$,
$\{1\otimes_K\alpha_1,\ldots, 1\otimes_K \alpha_n\}$ genera
a $E\otimes_K L$ sobre $E$ y $\dim_E(E\otimes_K L)
\leq \dim_K L=n$. El resultado se sigue pues $A=E
\otimes_K L$ es una $E$-\'algebra.
$\fin$
\end{proof}

Consideremos ahora $L=K(\theta)$ una extensi\'on simple y
sea $f(x)=\Irr(\theta,x, K)\in K[x]$. Sea $f(x)=p_1(x)^{e_1}
\cdots p_r(x)^{e_r}\in E[x]$ la descomposici\'on de $f(x)$ como
producto de polinomios irreducibles de $E[x]$. Entonces
\begin{align*}
E\otimes_K L&=E\otimes_K K(\theta)\cong E\otimes_K\big(
K[x]/\langle f(x)\rangle\big)\cong \big(E\otimes_K K[x]\big)/\langle f(x)\rangle\\
&\cong (E[x])/\langle f(x)\rangle \cong \bigoplus_{i=1}^r
\big(E[x]/\langle p_i^{e_i}(x)\rangle\big).
\end{align*}

Los ideales maximales de $E\otimes_K L$ est\'an dados por
\[
\langle p_1(x)\rangle/\langle f(x)\rangle,\ldots,
\langle p_r(x)\rangle/\langle f(x)\rangle
\]
y por tanto todas las composiciones inequivalentes de $E$
y $L$ sobre $K$ son $\Big\{E[x]/\langle p_i(x)\rangle\Big\}_{
i=1}^r$.

Se tiene $\dim_E(E\otimes_K L)=\deg f(x)=\sum_{i=1}^r e_i
\deg p_i(x)=[L:K]$. Si $L/K$ es una extensi\'on separable,
se tiene que $e_i$, $1\leq i\leq r$ y
\[
E\otimes_K L\cong \bigoplus_{i=1}^r E[x]/\langle p_i(x)\rangle
\cong \bigoplus_{i=1}^r E(\theta_i)
\]
donde $\theta_i$ es una ra\'iz de $p_i(x)$ y $E\otimes_K L$
es la suma directa de las composiciones de $E$ y $L$ sobre
$K$ y estas composiciones son $\big\{E(\theta_i)\big\}_{i=1}^r$.

De esta forma, si suponemos que $f(x)$ es separable y
enumeramos las ra\'ices de $f(x)$ en $E$, entonces
\[
f(x)=\prod_{i=1}^r p_i(x)=\prod_{i=1}^r\big(\prod_{j=1}^{m_i}(x
-\theta_{ij})\big)
\]
con $p_i(x)=\prod_{j=1}^{m_i}(x-\theta_{ij})$. Se tiene que
$E(\theta_{ij})$ y $E(\theta_{ij'})$ son equivalentes pues
$\theta_{ij}$ y $\theta_{ij'}$ son conjugados sobre $E$.

En resumen, $E$ y $L$ tienen $n=[L:K]$ composiciones con
$r$ clases inequivalentes. Las $r$ clases inequivalentes tienen
$m_1,\ldots, m_r$ elementos respectivamente.

\begin{ejemplo}\label{E17.6.8N} Sea $K={\ma Q}$, $E={\ma R}$
y sea $L=K(\theta)$ con $\theta^3=d$ con $d\in {\ma Q}$ que no
es cubo de ${\ma Q}$. Entonces
\[
x^3-d=(x-\theta)(x^2+\theta x+\theta^2)
\]
con $\theta^3=d$ y $\theta\in{\ma R}$. Las ra\'ices de $x^2+
\theta x+\theta^2\in{\ma R}[x]$ son $\theta_1=\zeta_3\theta$
y $\theta_2=\zeta_3^{-1}\theta$. Hay dos composiciones 
inequivalentes de $L$ con ${\ma R}$, a saber ${\ma R}(
\theta)={\ma R}$ y los encajes conjugados $\{{\ma R}(
\zeta_3\theta)={\ma C}, {\ma R}(\zeta_3^{-1}\theta)={\ma C}\}$.
\end{ejemplo}

\begin{ejemplo}\label{E17.6.9N} Sean $K={\ma Q}_p$ con $p$
un n\'umero primo tal que $p\equiv 1\bmod n$ con $n\in
{\ma N}$, $n\geq 3$. Sea $\zeta_n$ una ra\'iz $n$-\'esima primitiva
de la unidad y sea $L=\cic n{}$. Entonces $p$ se descompone
totalmente en $L/K$. Se tiene que $\zeta_n\in {\ma Q}_p$:
$\zeta_n^p=\zeta_n$ pues $p\equiv 1
\bmod n$ por lo que $\zeta_n\in {\ma F}_p$. Por tanto
$x^n-1\in{\ma F}_p[x]$ se factoriza en factores lineales y por
el Lema de Hensel, $\zeta_n$ tiene un levantamiento a 
${\ma Q}_p$

De esta forma tenemos que
si $\psi_n(x)=\prod_{\mcd(a,n)=1}(x-\zeta_n^a)
=\Irr(\zeta_n,x,{\ma Q})$ y $\prod_{\mcd(a,n)=1}(x-\zeta_n^a)$
es la factorizaci\'on de irreducibles de $\psi_n(x)$ en ${\ma
Q}_p[x]$, se sigue que las composiciones inequivalentes de
$\cic n{}$ con ${\ma Q}_p$ son
precisamente $\big\{{\ma Q}_p(\zeta_n^a)
\big\}_{\mcd(a,n)=1}$, esto es, ${\ma Q}_p(\zeta_n^a)=
{\ma Q}_p(\zeta_n)={\ma Q}_p$ pero son inequivalentes.

En resumen, hay $\varphi(n)$ composiciones de $\cic n{}$
y ${\ma Q}_p$ sobre ${\ma Q}$ inequivalentes y todos ellos
son iguales ${\ma Q}_p$.
\end{ejemplo}

Sea $K$ un campo global y sea $v$ un valor absoluto de $K$
y $K_v$ la completaci\'on de $K$ en $v$. El valor absoluto 
$v$ se extiende de manera \'unica a $K_v$: Si $\{x_n\}_{n=1}^{
\infty}$ es una sucesi\'on
de Cauchy en $K$, $\bar{v}\big(\{x_n\}\big):=\lim\limits_{
n\to\infty}v(x_n)$, $\bar{v}|_K=v$. Equivalentemente,
$\big|\{x_n\}\big|_v=\lim\limits_{n\to\infty}|x_n|_v$.
Se denota nuevamente por $v$ la extensi\'on a $K_v$.

Sea ahora $L/K$ una extensi\'on finita y separable de campos
globales y sea $v$ un valor absoluto de $K$. Sea $K_v$
la completaci\'on de $K$ en $v$. Entonces si $M$ es una
composici\'on de $L$ con $K_v$ sobre $K$, denotada de
momento como $L K_v$, se tiene que $LK_v$ es un campo
que es una extensi\'on finita de $K_v$ y por tanto $v$ se
extiende de manera \'unica a un valor absoluto $\omega$
en $L K_v$. Denotemos $L_{\omega}=LK_v$ con 
$\omega|_K=v$ y $L_{\omega}$ es completo.

De esta forma,
$LK_v\subseteq \bar K_v$ pues $LK_v$ es una extensi\'on
finita de $K_v$. Esta composici\'on es $(L_{\omega},\tau,\iota)$,
$\tau\colon L\lra L_{\omega}\subseteq \bar K_v$ y $\iota$
es el encaje natural $\iota:
K_v\lra L_{\omega}\subseteq \bar K_v$, $\iota|_K=\Id_K$, 
esto es $\iota(K_v)=K_v$ dentro de $\bar K_v$. Se tiene
$L_{\omega}=\tau(L)K_v$ con $\tau$ un encaje de $L$
en $\bar K_v$.

Rec\'iprocamente, si $\omega$ es un valor absoluto de $L$
con $\omega|_K\sim v$, podemos normalizar $\omega$ de
tal forma que $\omega|_K=v$. Entonces, si $L_{\omega}$
es la completaci\'on de $L$ con respecto a $\omega$ y si
$K_{\omega}$ es la cerradura de $K$ en $L_{\omega}$,
$K_{\omega}$ es completo y por tanto $K_{\omega}=K_v$
puesto que $K$ es denso en $K_{\omega}$. Adem\'as,
$L, K_v\subseteq L_{\omega}$ por lo que $LK_v$, la
composici\'on de $L$ y $K_v$, est\'a contenida en $L_{
\omega}$ y $L$ es denso en $L_{\omega}$ por lo que
$LK_v$ es denso y cerrado en $L_{\omega}$ (por ser
completo). Se sigue que $LK_v=L_{\omega}$.

De esta forma, hemos obtenido, todas las completaciones
$L_{\omega}$ con $\omega|_K=v$ son las composiciones
de $L$ con $K_v$ sobre $K$, es decir:
\[
\big\{\tau(L)K_v\big\}_{\tau:L\to \bar K_v}.
\]
Adem\'as, $\big|\{\sigma:L\to \bar K_v\}\big|=[L:K]_s
=[L:K]$, donde $[L:K]_s$ denota el grado de separabilidad
de la extensi\'on $L/K$ y
$|\ |_v$ se extiende de manera \'unica $|\ |_{
\tilde v}$ a $K_v$ y como todas las composiciones de 
$L$ con $K_v$ sobre $K$ est\'an contenidas en una 
extensi\'on finita $C$ de $K_v$, la extensi\'on de $|\ 
|_{\tilde v}$ a $C$ es \'unica, denotada nuevamente por
$|\ |_v$. De esta forma, los valores absolutos $\omega$ 
de $L$ con $\omega|_K=v$ pueden definirse como $|
\alpha|_{\omega}=|\sigma \alpha|_v$ 
donde $\alpha\in L$ y con $\sigma: L\hooklongrightarrow
\bar K_v$. En este caso se tiene $L_{\omega}=\sigma(
L) K_v$.

Dados dos encajes $\tau,\sigma:L\lra \bar K_v$, veremos
cuando $\tau(L)K_v$ y $\sigma(L)K_v$ tienen el mismo
valor absoluto. 

Recordemos que dos valores absolutos $|\ |_1$ y $|\ |_2$
en un campo $E$ se llaman {\em equivalentes\index{valores
absolutos equivalentes}} si para $x\in E$, se tiene
\[
|x|_1<1 \iff |x|_2<1.
\]
Equivalentemente, $|a|_1<|b|_1\iff |a|_2<|b|_2$.

Si $|\ |$ es cualquier valor absoluto en $E$ y $s\in{\ma R}$,
$s>0$, entonces $|\ |$ y $|\ |^s$ son equivalentes. Rec\'iprocamente,
se tiene

\begin{teorema}\label{T17.6.10N} Si $|\ |_1$ y $|\ |_2$ son
equivalentes en $E$, entonces existe $s\in{\ma R}$, $s>0$
con $|\ |_1^s=|\ |_2$.
\end{teorema}

\begin{proof}
Si $|\ |_1$ es trivial, $|\ |_2$ es tambi\'en trivial y rec\'iprocamente
(recordemos que $|\ |$ es {\em trivial\index{valor absoluto trivial}}
significa que $|x|=1$ para todo $x\neq 0$).  

Sea $0<|a_0|_2<1$. Sea $a\in E$, $a\neq 0$ y sea ${\mc A}=
\{(n,m)\in {\ma N}^2\mid |a_0|_2^n<|a|_2^m\}$. Por tanto, si
$(n,m)\in {\mc A}$, $|a_0|_2^n<|a|_2^m$ de donde se sigue que $n\log 
|a_0|_2<m\log |a|_2$. Como $\log |a_0|_2<0$, $(n,m)\in
{\mc A}\iff \frac nm>\frac{\log |a|_2}{\log |a_0|_2}$.

Como $|\ |_1$ es equivalente $|\ |_2$, se tiene $(n,m)\in {\mc A}
\iff \frac n m >\frac{\log|a|_1}{\log|a_0|_1}$. Se sigue que
$\frac{\log|a|_2}{\log|a_0|_2}=\frac{\log|a|_1}{\log|a_0|_1}$. Por
tanto $\frac{\log|a|_2}{\log|a|_1}=\frac {\log|a_0|_2}{\log|a_0|_1}$
y $\frac{\log|a|_2}{\log|a|_1}$ es independiente de $a$.

Sea $s:=\frac{\log|a|_2}{\log|a|_1}$. Entonces $|a|_1^s=|a|_2$
para toda $a\in E$, $a\neq 0$ y $s\in{\ma R}$, $s>0$ pues como
$|\ |_1$ y $|\ |_2$ son equivalentes, $\log |a_0|_1<0$ y
$\log|a_0|_2<0$, por lo que $s=\frac{\log|a_0|_2}{\log|a_0|_1}>0$.
$\fin$
\end{proof}

\begin{teorema}\label{T17.6.11N} Sea $v$ un valor absoluto
no trivial en un campo $E$ y sea $F/E$ una extensi\'on. Si
$\omega$ y $\omega'$ son dos extensiones de $v$ a $F$,
esto es, $\omega|_F=\omega'|_F=v$, y $\omega$ y $\omega'$
son equivalentes, entonces $\omega=\omega'$.
\end{teorema}

\begin{proof}
Se tiene que $|\ |_{\omega}^s=|\ |_{\omega'}$, $s\in {\ma R}$,
$s>0$. Puesto que si $a_0\in E$, con $|a_0|_v\neq 1$, entonces
$|a_0|_{\omega}=|a_0|_v=|a_0|_{\omega'}=|a_0|_{\omega}^s$,
de donde se sigue que $s=1$.
$\fin$
\end{proof}

En otras palabras, si $\omega$ y $\omega'$ son dos valores
absolutos de $F$ no iguales que extienden a $v$, entonces
$\omega$ y $\omega'$ no son equivalentes.

\s

Regresando a $|\ |_v$, un valor absoluto de $K$, y $|\ |_v$ se
extiende de manera \'unica $|\ |_{\tilde v}$ a $K_v$. Entonces,
si $L_{\omega}=\sigma(L)K_v$ y $L_{\omega'}=\tau(L)K_v$,
para $\sigma,\tau:L\lra \bar K_v$ y para $\alpha\in L$, $|\alpha
|_{\omega}=|\sigma\alpha|_v$ y $|\alpha|_{\omega'}=
|\tau\alpha|_v$. Por tanto $\omega$ y $\omega'$
son equivalentes si y s\'olo si $|\sigma\alpha|_v=|\tau\alpha|_v$,
esto es, dan el mismo valor absoluto en $L$.

Sea ahora $\sigma : L\lra \bar K_v$, $\bar K_v$ una cerradura
separable de $K_v$, un encaje. Entonces $\sigma(L)\cong L$ y
$\sigma(L)K_v/K_v$ es una extensi\'on finita. Entonces $\sigma
(L)K_v$ es una completaci\'on de $L$, digamos $L_{\sigma}$.
Esto es, $\sigma$ se identifica con una completaci\'on de $L$,
$L_{\sigma}``="L_{\omega}$ con $\omega|_{K_v}=v$.

\begin{definicion}\label{D17.6.12N} Sea $K$ un campo global.
Los {\em valores absolutos can\'onicos\index{valor absoluto 
can\'onico}} de $K$ estan dados de la siguiente forma:
Si ${\eu p}$ es finito, $\o_{\pK}/\pK=K(\pK)\cong {\ma F}_{
q_{\pK}}$, entonces seleccionamos $c=q_{\pK}^{-1}$ para
el valor absoluto:
\begin{gather*}
|x|_{\pK}=q_{\pK}^{-v_{\pK}(x)},\quad \N(\pK)=q_{\pK}.
\intertext{Si $\pK$ es infinito, $K_{\pK}\cong
{\ma R}$ o ${\ma C}$, si $\sigma: K\lra K_{\pK}$, entonces}
|x|_{\pK}=|\sigma x|_{\ma R}=|\sigma x|\text{\ usual si $\pK=
{\ma R}$ y $|x|_{\pK}=\|\sigma x\|_{\ma C}^2=\|\sigma x\|_{
\ma C}^2$ en ${\ma C}$ si $\pK={\ma C}$}
\end{gather*}
pues en un lugar
complejo consiste de 2 encajes $\sigma$ y $\bar \sigma$ con
$\sigma\neq \bar\sigma$.
\end{definicion}

El conjunto de los valores absolutos can\'onicos se dontar\'a por 
${\mc M}_K$. Se tiene que $\pK, \pK' \in {\mc M}_K$,
$\pK\neq \pK'$, entonces $|\ |_{\pK}$ y $|\ |_{\pK'}$ dan
diferentes topolog\'ias en $K$.

Dado un valor absoluto $\omega$ de $L$, $\omega|_K=v$,
$L_{\omega}/K_v$ es una extensi\'on finita. Ahora si $\sigma
:L_{\omega}\lra\bar K_v$ es uno de los $[L_{\omega}:K_v]=
[L_{\omega}:K_v]_s$ encajes, $\sigma|_{K_v}=\Id_{K_v}$.
Entonces $\sigma(L)\subseteq \sigma(L_{\omega})$ y $
\sigma(L)K_v=\sigma(L_{\omega})$ por el argumento
dado anteriormente.

\begin{teorema}\label{T17.6.13N} Sea $K$ un campo global,
$v\in{\mc M}_K$ y $L/K$ una extensi\'on finita y separable.
Se tiene que dos encajes $\sigma, \tau:L\lra \bar K_v$ sobre
$K$, $\sigma|_K=\tau|_K=\Id_K$, dan lugar al mismo valor
absoluto en $L$ si y s\'olo si son conjugados sobre $K_v$,
lo cual significa que existe un isomorfismo $\lambda:\sigma
(L)K_v\lra \tau(L)K_v$ tal que $\lambda|_{K_v}=\Id_{K_v}$.
En otras palabras, los encajes dan lugar al mismo valor
absoluto si y s\'olo si las composiciones $\sigma(L)K_v$ y
$\tau(L)K_v$ de $L$ y $K_v$ sobre ya sea $K$ o sobre
$K_v$, son equivalentes.
\end{teorema}

\begin{proof}
Dado cualquier encaje $\theta:L\lra \bar K_v$, entonces para
$\beta\in K$, $\theta(\beta)\in\theta(L)$ y $\theta(L)K_v/K_v$
es una extensi\'on finita.

La extensi\'on de $v$ a $\theta(L)K_v$ es \'unica y si ponemos
$\theta(L)K_v=L_{\omega}$, $\omega$ la extensi\'on de $v$ a
$L_{\omega}$, entonces
\[
|\beta|_{\theta}=|\theta \beta|_{\omega}=|\N_{L_{\omega}/K_v}
(\theta\beta)|^{1/n_{\omega}},
\]
donde $n_{\omega}=[L_{\omega}:K_v]$. Sean $\sigma$ y $\tau$
conjugados, $\sigma(L)K_v=L_{\omega}$, $\tau(L)K_v=L_{
\omega'}$, $\lambda:L_{\omega}\lra L_{\omega'}$ isomorfismo
y $\lambda|_{K_v}=\Id_{K_v}$, entonces
\begin{gather*}
\N_{L_{\omega}/K_v}(\sigma\beta)=\N_{L_{\omega'}/K_v}(
\lambda\sigma\beta)=\N_{L_{\omega'}/K_v}(\tau\beta),\\
|\beta|_{\omega}=|\N_{L_{\omega}/K_v}(\sigma
\beta)|^{1/n_{\omega}}=|\N_{L_{\omega'}/K_v}(\tau\beta)|^{
1/n_{\omega'}}=|\beta|_{\omega'}\quad \text{y}\\
n_{\omega}=[L_{\omega}:K_v]=[L_{\omega'}:K_v]=n_{\omega'}.
\end{gather*}
Se sigue que $|\beta|_{\omega}=|\beta|_{\omega'}$.

Rec\'iprocamente, supongamos que los valores absolutos son 
los mismos. Ahora $\tau(L)$ y $\sigma(L)$ son isomorfos a
$L$ sobre $K$. Sea $\lambda:\tau(L)\lra\sigma(L)$ un isomorfismo
sobre $K$. Veremos que $\lambda$ se puede extender de $
\tau(L)K_v$ a $\sigma(L)K_v$ sobre $K_v$.

Se tiene que $\tau(L)$ es denso en $\tau(L)K_v$, por lo que 
dado $x\in \tau(L)K_v$ puede ser escrito como $x=\lim\limits_{
n\to\infty}\tau x_n$ con $x_n\in L$. Ahora bien, como los
valores absolutos inducidos por $\tau$ y por $\sigma$ son el
mismo, se sigue que $\{\lambda\tau x_n\}=\{\sigma x_n\}$
converge a un elemento de $\sigma(L)K_v$. Denotemos este
elemento por $\lambda x$. Se verifica que $\lambda x$ es
independiente de $\{\tau x_n\}$ y que $\lambda:\tau(L)K_v
\lra \sigma(L)K_v$ es un isomorfismo que deja fija a $K_v$.
$\fin$
\end{proof}

Con este resultado, podemos enumerar las extensiones de
$v$ a $L$, una f\'ormula ya muy conocida.

\begin{corolario}\label{C17.6.14N} Sea $L/K$ una extensi\'on
finita y separable de grado $n$ de campos globales y $v\in
{\mc M}_K$. Para cada valor absoluto $\omega$ de $L$, sea
$n_{\omega}=[L_{\omega}:K_v]$ el grado de campos locales
(o de ${\ma R}$ o de ${\ma C}$). Entonces $\sum_{\omega|v}
n_{\omega}=n$.
\end{corolario}

\begin{proof}
Se tiene que $n=[L:K]=[L:K]_s$ es el n\'umero de encajes
$\sigma:L\lra\bar K_v$, $\sigma|_K=\Id_K$. Adem\'as
$[L_{\omega}:K_v]=n_{\omega}$ es el n\'umero de encajes
$L_{\omega}\lra \bar K_v$ y todos ellos son conjugados 
sobre $K_v$ y por tanto dan lugar al mismo valor absoluto.
Se sigue que $n=\sum_{\omega|v}n_{\omega}=\sum_{
\omega|v}e_{\omega}f_{\omega}$.
$\fin$
\end{proof}

\begin{corolario}\label{C17.6.15N} Sean $K$ un campo global
y $L/K$ una extensi\'on finita y separable. Sea $\alpha\in L$.
Entonces, si $v\in{\mc M}_K$, se tiene
\[
\prod_{\omega|v}|\alpha|_{\omega}^{n_{\omega}}=
|\N_{L/K}(\alpha)|_v,
\]
con $n_{\omega}=[L_{\omega}:K_v]$.
\end{corolario}

\begin{proof}
La extensi\'on del valor absoluto de $K_v$ a $L_{\omega}$ 
est\'a dado por $|\alpha|_{\omega}=|\N_{L_{\omega}/K_v}(\alpha)|^{
1/n_{\omega}}$. Por tanto, si $A_{\omega}$ es la clase de
conjugaci\'on de $\omega$, se tiene
\begin{gather*}
|\alpha|_{\omega}^{n_{\omega}}=|\N_{L_{\omega}/K_v}(\alpha)|
=\big|\prod_{\sigma\in A_{\omega}}\sigma\alpha\big|_v,\\
\prod_{\omega|v}|\alpha|_{\omega}^{n_{\omega}}=\prod_{
\omega|v}\big|\prod_{\sigma\in A_{\omega}}\sigma \alpha\big|
=\big|\prod_{\sigma} \sigma \alpha\big|_v,
\end{gather*}
donde $\sigma$ recorre los $[L:K]$ encajes $\sigma:L\lra\bar K$.
Se sigue que $\prod_{\omega|v}|\alpha|_{\omega}^{n_{\omega}}=
\big|\N_{L/K} \alpha\big|_v$.
$\fin$
\end{proof}

\begin{corolario}\label{C17.6.16N} Sea $L/K$ una extensi\'on 
finita y separable de campos globales. Entonces $v\in{\mc M}_K$
se tiene
\[
\N_{L/K}(\alpha)=\prod_{\omega|v}\N_{L_{\omega}/K_v}(\alpha)
\quad \text{y}\quad
\Tr_{L/K}(\alpha)=\sum_{\omega|v}\Tr_{L_{\omega}/K_v}(\alpha)
\]
para $\alpha\in L$.
\end{corolario}

\begin{proof}
Se tiene que si ${\mc A}_{\omega}=\{\sigma:L_{\omega}\lra\bar K_v\mid
\sigma|_{K_v}=\Id_{K_v}\}$ y ${\mc A}'_{\omega}=\{\sigma|_L\mid\sigma
\in {\mc A}_{\omega}\}$, entonces $\bigcup_{\omega|v}{\mc A}'_{\omega}
=\{\sigma:L\lra \bar K\mid \sigma|_K=\Id_K\}$ de donde se sigue el
resultado.
$\fin$
\end{proof}

\begin{proposicion}\label{P17.6.17N}
Sea $K\in\{{\ma Q}, \F(T)\}$. Entonces
\[
\prod_{v\in{\mc M}_K}|\alpha|_v=1.
\]
para toda $\alpha\in \*K$.
\end{proposicion}

\begin{proof}
Sea $\alpha=\gamma p_1^{a_1}\cdots p_r^{a_r}$ con $\{p_1,\ldots,
p_r\}$ ya sea n\'umeros primos distintos si $K={\ma Q}$ o polinomios
m\'onicos irreducibles distintos si $K=\F(T)$ y $\gamma=\pm 1$ en el caso
num\'erico y $\gamma\in \*\F$ en el caso de campos de funciones.

Entonces si $K={\ma Q}$,
\begin{gather*}
|\alpha|_v=\begin{cases}
p_i^{-a_i}&\text{si $v=p_i$},\\
1& \text{si $v\notin\{p_1,\ldots,p_r\}$ y $v$ es finito},\\
|\alpha|&\text{el valor absoluto usual si $v=\infty$},
\end{cases}
\intertext{y si $K=\F$, entonces}
|\alpha|_v=\begin{cases}
q^{-a_i\deg p_i}&\text{si $v=p_i$},\\
1& \text{si $v\notin\{p_1,\ldots,p_r\}$ y $v$ es finito},\\
q^{\deg \alpha}&\text{si $v=\p$, $v=1/T$ es el primo infinito}.
\end{cases}
\intertext{Se sigue que}
\prod_{v\in{\mc M}_K}|\alpha|_v=\begin{cases}
\prod_{i=1}^r p_i^{-a_i}\cdot (p_1^{a_1}\cdots p_r^{a_r})=1&\text{si $K={\ma Q}$},\\
\\
q^{-\sum_{i=1}^r a_i\deg p_i}\cdot q^{\sum_{i=1}^r a_i\deg p_i}=1&
\text{si $K=\F$}.
\end{cases}\tag*{$\fin$}
\end{gather*}
\end{proof}

\begin{teorema}[Producto de los valores absolutos\index{formula de
producto@f\'ormula del producto}]\label{T17.6.18N}
Sea $L$ cualquier campo global. Entonces para $\alpha\in \*L$, se tiene 
\[
\prod_{\omega\in{\mc M}_L}|\alpha|_{\omega}=1.
\]
\end{teorema}

\begin{proof}
Si $L$ es un campo de funciones, sabemos que para $\alpha\in\*L$
$\deg\alpha=\deg(\alpha)_L=0$ donde $(\alpha)_L$ denota al divisor principal
de $\alpha$, $\deg \alpha=\sum_{\pK\in{\mc M}_L}\deg \pK\cdot
v_{\pK}(\alpha)$ y $|\alpha|_{\pK}=q^{-\deg \pK \cdot v_{\pK}(\alpha)}$. Por
tanto $\prod_{\pK\in{\mc M}_L}|\alpha|_{\pK}=q^{-\deg\alpha}=q^0=1$.

Ahora damos otra demostraci\'on que funciona para cualquier campo global.
Si $L$ es campo de funciones, por ser $\F$ un campo perfecto, $L$
es una extensi\'on finita y separable de alg\'un $\F(T)$. 

Para $L$ cualquier campo global, denotamos $n_{\omega}=[L_{\omega}:
K_v]$ donde $\omega\in {\ma P}_L$ y $\omega|_K=v\in {\ma P}_K$ y donde
$K\in\{{\ma Q}, \F(T)\}$ seg\'un sea el caso. Sabemos que $\prod_{\omega|v}
|\alpha|_{\omega}^{n_{\omega}}=\big|\N_{L/K}(\alpha)|_v$, Corolario 
\ref{C17.6.15N}.

Ahora bien $\omega$ es la \'unica extensi\'on de $v$ a $L_{\omega}$, sin 
embargo $|\alpha|_{\omega}=\big|\N_{L_{\omega}/K_v}(\alpha)\big|_v^{1/
n_{\omega}}$ no est\'a normalizada pues $\deg \omega=f_{\omega}\deg v$
donde $f_{\omega}=f(\omega|v)$.

Si $v$ es un lugar finito (es decir, no arquimediano), 
se tiene que el campo residual de $K_v$ es ${\ma F}_{q^{\deg v}}$ y el de
$L_{\omega}$ es ${\ma F}_{q^{\deg \omega}}$. Por otro lado $\deg \omega=
f(\omega|v)\deg v$ donde $n_{\omega}=[L_{\omega}:K_v]=f(\omega|v)
e(\omega|v)$. Adem\'as $v_{L_{\omega}}(\alpha)=\frac {1}{f(\omega|v)}
v_{K_v}\big(\N_{L_{\omega}/K_v}(\alpha)\big)$ (Corolario \ref{C17.3.2.7}),
por lo que 
\begin{align*}
\big|\N_{L_{\omega}/K_v}(\alpha)\big|_v&=q^{-\deg v\cdot v_{K_v}(
\N_{L_{\omega}/K_v}(\alpha))}=q^{-\deg
v\cdot f(\omega|v)\cdot v_{L_{\omega}}(\alpha)}\\
&=q^{-\deg \omega \cdot v_{L_{\omega}}(\alpha)}=\|\alpha\|_{\omega}
\end{align*}
normalizada. Esto es, el valor absoluto normalizado de $\alpha$ en
$L_{\omega}$ es $\|\alpha\|_{\omega}=\big|\N_{L_{\omega}/K_v}
(\alpha)\big|_v$.

Se sigue que, para $v$ un valor absoluto no arquimediano,
\[
\prod_{\omega|v}\|\alpha\|_{\omega}=\big|\prod_{\omega|v}
\N_{L_{\omega}/K_v}(\alpha)\big|_v=\big|\N_{L/K}(\alpha)\big|_v.
\]
Para el caso arquimediano, tenemos que si $\omega, v$ son ambos
reales o ambos complejos, $L_{\omega}=K_v$ y $\N_{L_{\omega}/K_v}
(\alpha)=\alpha$ y 
\[
\big|\N_{L_{\omega}/K_v}(\alpha)\big|_{\omega}=
|\alpha|_v=\begin{cases}
|\alpha|^2&\text{si $w$ y $v$ son complejos},\\
|\alpha|&\text{si $w$ y $v$ son reales}.
\end{cases}
\]
Si $v$ es real y $\omega$ es complejo, $\N_{L_{\omega}/K_v}(\alpha)=
\alpha\bar{\alpha}=|\alpha|^2$ y $\big|\N_{L_{\omega}/K_v}(\alpha)
\big|_{\omega}=\big|\N_{L_{\omega}/K_v}(\alpha)\big|^2$. Entonces,
en el caso arquimediano, tambi\'en tenemos que $|\alpha|_{\omega}=
\big|\N_{L_{\omega}/K_v}(\alpha)\big|_v$.

Se sigue que para cualquier valor absoluto $v$ se tiene
\begin{gather*}
\prod_{\omega\in{\mc M}_L}|\alpha|_{\omega}=\prod_{v\in
{\mc M}_K}\big|\N_{L/K}(\alpha)\big|=1.\tag*{$\fin$}
\end{gather*}
\end{proof}

Sea $S$ un conjunto finito no vac\'io de lugares de un campo global
$K$, $S$ conteniendo a los primos arquimedianos en caracter\'istica $0$.
Sea $K^S:=\{a\in\*K\mid v_{\pK}(a)=0 \text{\ (es decir, $|a|_{\pK}=1$) para
todo ${\pK}\notin S$}\}$. $K^S$ se llama {\em el grupo de las 
$S$-unidades\index{grupo de las $S$-unidades}}. Si $K$ es un campo
num\'erico y $S$ consiste exactamente de los lugares arquimedianos de
$K$, entonces $K^S$ son las unidades de $K$, esto es, $K^S=
\*{\o_K}$ donde $\o_K$ es el anillo de enteros de $K$.

Se tiene el siguiente resultado, cuya demostraci\'on presentamos
m\'as adelante.

\begin{teorema}[Teorema de las unidades de Dirichlet\index{teorema
de las unidades de Dirichlet}\index{Dirichlet!teorema de las
unidades de $\sim$}]\label{T17.6.19N}
Se tiene que, como grupo,
\[
K^S\cong\begin{cases}
W_K\times {\ma Z}^{|S|-1}&\text{si $\car K=0$},\\
\*\F\times {\ma Z}^{|S|-1}&\text{si $\car K=p>0$},
\end{cases}
\]
donde $W_K$ son las ra\'ices de unidad en $K$.
\end{teorema}

El grupo de divisores (caracter\'istica $p>0$) o de ideales fraccionarios
(caracter\'istica $0$), seg\'un sea el caso, se denotar\'a por $D_K$ y
en caracter\'istica $p$\label{CClaseDK}
tenemos $D_{K,0}=\{{\eu a}\in D_K\mid
\deg {\eu a}=0\}$\label{CClaseDK0}. Si $P_K=\{(\alpha)\mid \alpha\in\*K\}$ denota a los
divisores o ideales fraccionarios principales, tenemos los grupos
de clases de divisores o de ideales fraccionarios $D_K/P_K=:I_K
\label{CClaseIK}$
y $D_{K,0}/P_K=I_{K,0}\label{CClaseIK0}$.
Se tiene que $I_K$ es finito cuando $K$
es campo num\'erico y es infinito en el caso de campos de funciones.
Por otro lado $I_{K,0}$ es finito ($K$ un campo de funciones).

Usaremos la siguiente notaci\'on. Sea $\pK$ un lugar,
\[
U_{\pK}=\begin{cases}
\text{grupo de unidades de $K_{\pK}$ si $\pK$ es finito},\\
\*{K_{\pK}}\text{\ si $\pK$ es infinito}.
\end{cases}
\]
Como vimos anteriormente, $\sum_{\pL|\pK}[L_{\pL}:K_{\pK}]=[L:K]$.

Si $L/K$ es una extensi\'on finita de Galois y $\sigma\in G_{L|K}$,
entonces para $\pL|\pK$, entonces $\sigma\pL|\pK$ y $
L_{\pL}\xrightarrow[\ \cong\ ]{\sigma} L_{\sigma\pL}$ es un 
isomorfismo sobre $K_{\pK}$. De hecho, si $\alpha\in L_{\pL}$,
existe una sucesi\'on en $L$, $\alpha=\lim\limits_{n\to\infty} x_n$,
$x_n\in L$, $\sigma x_n$ converge $L_{\sigma\pL}$ a $\sigma
\alpha$ pues el valor absoluto (no normalizado) satisface
\begin{align*}
|\alpha-x_n|_{\pL}&=\big|\N_{L_{\pL}/K_{\pK}}(\alpha-x_n)\big|^{
1/[L_{\pL}:K_{\pK}]}\\
&=\big|\N_{L_{\pL}/K_{\pK}}(\sigma\alpha-
\sigma x_n)\big|^{1/[L_{\pL}:K_{\pK}]}=|\sigma\alpha-\sigma
x_n|_{\sigma\pL}.
\end{align*}

Si $\pL=\sigma\pL$, $L_{\pL}\stackrel{\sigma}{\lra} L_{\pL}$ es un 
automorfismo y $\sigma\in G_{L_{\pL}|K_{\pK}}=D_{L/K}(\pL|\pK)$,
$|D_{L/K}(\pL|\pK)|=e(\pL|\pK)f(\pL|\pK)$ y $h(L/K)=\big|G_{L|K}/D_{L/K}(
\pL|\pK)\big|$.

\s

Con respecto a la Teor\'ia de Kummer (Teorema \ref{CClaseT1.6.2}),
obtenemos la siguiente informaci\'on aritm\'etica

\begin{teorema}\label{T17.6.20N}
Sea $K$ un campo global que contiene una $n$-\'esima
ra\'iz primitiva de la unidad $\zeta_n$, donde $p\nmid n$, y
$p$ es la caracter\'istica de $K$. Sean $L=K\big(\sqrt[n]{\Delta})$,
una extensi\'on finita de Kummer de $K$, y $\pK$ un primo de
$K$. Entonces:
\las
\item $\pK$ se descompone totalmente en $L\iff \Delta\subseteq
(\*K_{\pK})^n$.
\item Si $\pK$ es finito y $\pK\nmid n$, entonces $\pK$ es no 
ramificado en $L \iff$ existe un conjunto de generadores de 
$\Delta/(\*K)^n$ que son unidades en $\pK$. En otras palabras,
$\Delta\subseteq U_{\pK}\cdot (\*K)^n$.
\end{list}
\end{teorema}

\begin{proof}
(1): Se tiene que $\pK$ se descompone totalmente en $L/K \iff
L_{\pL}=K_{\pK}$ para todo $\pL|\pK$, $L_{\pL}=K_{\pK}\big(
\sqrt[n]{\Delta})=K_{\pK}\iff \sqrt[n]{\Delta}\subseteq \*K_{\pK}\iff
\Delta\subseteq (\*K_{\pK})^n$. 

\s

(2): Sea $\pL|\pK$. Supongamos $\pK$ no ramificado. Se tiene
que $L_{\pL}/K_{\pK}$ es no ramificada si el grupo de valores de
la valuaci\'on de $L_{\pL}$ es el mismo que el de $K_{\pK}$. Sea
$\delta$ uno de los generadores de $\Delta$. Puesto que $\delta=
(\delta^{1/n})^n$ en $L_{\pL}$, su valuaci\'on es un $n$-m\'ultiplo
de un valor: $v_{\pL}(\delta)=nv_{\pL}(\sqrt[n]{\delta})$ y $v_{\pL}
(\delta)=e(\pL|\pK)v_{\pK}(\delta)=v_{\pK}(\delta)$. Por tanto se
tiene que $v_{\pK}(\delta)=rn$ para alguna $r\in{\ma Z}$. Sea
$\pi\in K$, $v_{\pK}(\pi)=1$. Entonces $v_{\pK}\big(\delta(\pi^{-r})^n
\big)=0$ y $(\pi^{-r})^n$ es una $n$-potencia, por lo que $\delta
(\pi^{-r})^n$ es tamb\'ien un generador y es una $\pK$-unidad.

Rec\'iprocamente, supongamos que $\Delta/(\*K)^n$ tiene 
generadores $\delta$ que son unidades en $\pK$. Es suficiente
probar que $K(\sqrt[n]{\delta})/K$ es no ramificada para $\delta
\in \Delta$. Se tiene que $\alpha=\sqrt[n]{\delta}$ satisface
$f(x)=x^n-\delta$ y $f'(\alpha)=n\alpha^{n-1}$ es primo relativo
a $\pK$. Entonces $\pK$ no divide al diferente local y por 
tanto es no ramificado.
$\fin$
\end{proof}

\subsubsection{Teor\'ia de Kummer aditiva o de Artin-Schreier-Witt}
\label{Artin-Schreier}

La teor\'ia de Kummer aditiva es necesaria en nuestro enfoque para
probar la segunda desigualdad fundamental. \'Unicamente requerimos
estudiar las extensiones abelianas de un campo global de funciones
de exponente $p$. El tratamiento de exponente $p^n$ se puede
hacer cambiando las extensiones de Artin-Schreier por extensiones
usando vectores de Witt. Los vectores de Witt fueron desarrollados
en el Cap\'itulo \ref{Ch9'}. Aqu\'i presentamos \'unicamente las
extensiones de exponente $p$ para simplificar la presentaci\'on,
en el entendido que la teor\'ia es inmediatamente generalizable
para extensiones de exponente $p^n$ para cualquier $n\in{\ma N}$
cambiando la operaci\'on $\wp(x)=x^p-x$ o $\wp(x)=x^q-x$ por
$\wp(\vec x)=\vec x^p\Witt - \vec x$ o $\wp (\vec x)= \vec x^q
\Witt - \vec x$, ${\ma F}_p$ por $W_n({\ma F}_p)\cong {\ma Z}/p^n
{\ma Z}$, etc.

\s

Sea $K$ un campo de caracter\'istica $p>0$ y sea $L/K$ una 
extensi\'on abeliana ya sea finita o infinita, de exponente $p$,
esto es, $\sigma^p=1$ para toda $\sigma \in G=\Gal(L/K)$. Podemos
usar a ${\ma F}_p$ como el grupo de valores para los caracteres
de $G$, esto es, $\chi(G)=\Hom(G,{\ma F}_p)$. Un caracter $\mu
\in\chi(G)$ satisface
\[
\mu(\sigma\tau)=\mu(\sigma)+\mu(\tau)=\mu(\sigma)+\sigma\mu(\tau),
\]
por lo que $\mu$ es un cociclo aditivo $Z^1(G,K)$. Ahora bien, $\co 1GL
=\{0\}$, por lo que existe $\alpha\in L$ tal que $\mu(\sigma)=(1-\sigma)
\alpha$ para toda $\sigma\in G$.

Sea $\wp:L\lra L$ el operador de Artin-Schreier $\wp(\beta)=\beta^p-\beta$.
Se tiene que $\ker \wp=\{\beta\in L\mid \beta^p-\beta=0\}={\ma F}_p$ y
$\wp(\sigma \beta)=\sigma \wp(\beta)$ para toda $\beta \in L$. Se tiene que
\[
0\lra {\ma F}_p\lra L\stackrel {\wp}{\lra} \wp(L)\lra 0,
\]
es una sucesi\'on $G$-exacta. En cohomolog\'ia obtenemos la sucesi\'on
exacta
\begin{gather*}
L^G=K\stackrel{\wp}{\lra}\wp(L)^G=\wp(L)\cap K\stackrel{\delta}{\lra}
\co 1G{{\ma F}_p}\lra \co 1GL=\{0\}.
\intertext{Ahora bien, $\co 1G{{\ma F}_p}=\Hom (G,{\ma F}_p)\cong \chi(G)$.
Por tanto}
\frac{\wp(L)\cap K}{\wp(K)}\cong \chi(G).
\end{gather*}

El isomorfismo anterior est\'a dado por el mapeo de conexi\'on $\delta$.
Si $a+\wp(K)\in \frac{\wp(L)\cap K}{\wp(K)}$, formalmente, se tiene
\[
\xymatrix{
&{\begin{array}{c} \co 0GL\\ \wp^{-1}(a)\end{array}}
\ar@<1.5ex>@{<-}[r]\ar@<-1.5ex>@{<-}[r]\ar@{->}[d]^{\partial}
&{\begin{array}{c}\co 0G{\wp(L)}\\ a \end{array}}\\
\co 1G{{\ma F}_p}\ar@{<-}[r]& \co 1GL
}
\]
Por tanto $\delta(a+\wp(K))=\mu_a\in\chi(G)$ donde $\mu_a(\sigma)=
\wp^{-1}a-\sigma(\wp^{-1}a)$. Ahora bien, $\wp^{-1}a=b$ significa
$\wp b=a=b^p-b$ por lo que $\sigma(b)=b-\mu_a(\sigma)$ con
$\mu_a(\sigma)\in {\ma F}_p$.

Sea $N$ la m\'axima extensi\'on de Artin-Schreier sobre $K$. Entonces
$\wp(L)\subseteq \wp(N)\cap K\subseteq K$. Si tuvi\'esemos $\wp(N)
\cap K\neq K$, existir\'ia $a\in K$ con $\wp^{-1}a\notin N$ y $N(\wp^{
-1}a)$ ser\'ia una extensi\'on de Artin-Schreier conteniendo propiamente
a $N$. Por tanto $\wp(N)\cap K=K$.

Sea $G=\Gal(N/K)$. Entonces $K/\wp(K)\cong \co 1G{{\ma F}_p}\cong
\chi(G)$. Tenemos un isomorfismo de redes que preserva contenciones
entre la red de extensiones de Artin-Schreier y la red de subgrupos de
$K$ que contienen a $\wp(K)$.

Si $L$ corresponde al subgrupo $\Lambda$, entonces
\lasa
\item $\Lambda=\wp(L)\cap K$: Si $a\in K$, $a\in\Lambda\iff \mu_a(\sigma)
=0 \ (\mu_a\in\chi\big(\Gal(N/K)\big)$ para $\sigma\in\Gal(N/L)\iff \sigma(
\wp^{-1}a)=\wp^{-1}a$ para $\sigma\in\Gal(N/L)\iff \wp^{-1} a\in L\iff a\in
\wp(L)$.
\item $L=K(\wp^{-1}(\Lambda))$: Si $\sigma\in G=\Gal(N/K)$ entonces
$\sigma\in\Gal(N/L)\iff \mu_a(\sigma)=0$ para $a\in\Lambda\iff
\sigma(\wp^{-1}a)=\wp^{-1}a$ para $a\in\Lambda \iff \sigma|_{K(\wp^{-1}
(\Lambda))}=\Id_{K(\wp^{-1}(\Lambda))}\iff \sigma\in \Gal(N/K(\wp^{-1}(
\Lambda)))$.
\end{list}

De esta forma, hemos probado:

\begin{teorema}\label{T17.6.21N}
Sea $K$ un campo de caracter\'istica $p>0$. Existe una correspondencia
entre los subgrupos aditivos $\Lambda$ de $K$ que contienen a $\wp^{-1}
(K)$ y las extensiones abelianas $L/K$ de exponente $p$. La 
correspondencia est\'a dada por $\Lambda\longleftrightarrow L_{\Lambda}
=K(\wp^{-1}(\Lambda))$. Adem\'as $\chi(\Gal(N/K))\cong \frac{\wp(L)
\cap K}{\wp(K)}$. $\fin$
\end{teorema}

\subsection{Anillo de ad\`eles y grupo de id\`eles}\label{S17.6.3N}

Tanto los {\em ad\`eles\index{ad\`eles}\label{CClaseadeles}}
como los {\em id\`eles\index{id\`eles}}
son un caso especial de la siguiente construcci\'on de los llamados
``productos directos restringidos''.

\begin{definicion}[Productos directos restringidos]\label{D17.6.22N}
Sea $\{v\}$ una familia de \'indices y para cada $v$ sea $G_v$
un grupo abeliano localmente compacto, esto es, cada $x\in G_v$
tiene una vecindad compacta. Para casi todo \'indice $v$, sea $H_v$
un subgrupo abierto compacto de $G_v$. Entonces el {\em producto
directo restringido\index{producto directo restringido}} de los $G_v$ con
respecto a los $H_v$ es el subgrupo $G$ de ${\mc G}:=\prod_{v}G_v$
consistente de los elementos de ${\mc G}$ cuyas todas las componentes,
salvo un n\'umero finito, pertenecen a $H_v$:
\[
G=\{(x_v)_v\in{\mc G}\mid x_v\in H_v\text{\ para casi toda $v$}\}.
\]
\end{definicion}

Si $S$ es un conjunto finito de \'indices que incluye a todos los \'indices
$v$ tales que $H_v$ no est\'a definido, entonces se define $G_S:=
\prod_{v\in S}G_v\times\prod_{v\notin S}H_v$.

Entonces $G_S$ es el producto directo de grupos localmente compactos,
por lo que $G_S$ es localmente compacto en la topolog\'ia producto.
Se tiene que $G_S$ es un conjunto abierto de $G$ y $G$ es localmente
compacto.

Se tiene $G_v \xhookrightarrow{\phantom{xx}} G$,
$g_v\xrightarrow{\phantom{xx}} (\ldots,1,\ldots,1,g_v,1\ldots,1,\ldots)$
el cual es un monomorfismo de grupos.

\begin{definicion}\label{D17.6.23N}
Los {\em ad\`eles\index{ad\`eles}} de un campo global $K$, ${\ma A}_K$,
es el producto restringido de $\{K_v\}_{v\in{\ma P}_K}$ con respecto
a $\{\o_v\}_{v\text{\ finito}}$.

Los {\em id\`eles\index{id\`eles}\label{CClaseideles}}
de un campo global $K$, $J_K$, es 
el producto directo restringido de $\{\*K_v\}_{v\in{\ma P}_K}$ con 
respecto a $\{\*\o_v=U_v\}_{v\text{\ finito}}$.
\end{definicion}

Como veremos m\'as adelante, la topolog\'ia de ${\ma A}_K$ restringida
a $J_K$ no es la misma topolog\'ia de $J_K$.

Veamos expl\'icitamente los ad\`eles o reparticiones en un campo global
$K$. Un ad\`ele es una familia $\vec\alpha=(\alpha_{\pK})_{\pK\in{\ma
P}_K}=(\alpha_v)_{v\in{\ma P}_K}$, donde ${\ma P}_K$ denota el
conjunto de lugares (valores absolutos) de $K$, $\alpha_{\pK}\in
K_{\pK}$ para toda $\pK$ y $\alpha_{\pK}\in \o_{\pK}
$ para casi todo $\pK$ finito.
Indistintamente usaremos ${\ma P}_K$ y ${\mc M}_K$, esto es,
siempre supondremos, a menos que se indique lo contrario, que
todos los valores absolutos estar\'an normalizados (de manera 
can\'onica).

De esta forma, ${\ma A}_K$ denota al anillo de los ad\`eles con las
operaciones entrada por entrada y la topolog\'ia de ${\ma A}_K$ 
est\'a dada por la base de vecindades abiertas de $\vec 0$ 
consistente de los conjuntos $\prod_{v\in S}V_v\times \prod_{v
\notin S}\o_v$ donde $S\subseteq {\ma P}_K$ es un conjunto
finito y $V_v$ es un conjunto abierto de $K_v$
y donde $0\in V_v$ para $v\in S$.

\s

Ahora bien el {\em grupo de id\`eles\index{grupo de id\`eles}} de un
campo global $K$ est\'a formado por los vectores $\vec\alpha=
(\alpha_{\pK})_{\pK\in{\ma P}_K}=(\alpha_v)_{v\in{\ma P}_K}$ donde
$\alpha_{\pK}\in\*K_{\pK}$ para toda $\pK\in{\ma P}_K$ y $|\alpha_{
\pK}|_{\pK}=1$, equivalentemente, $v_{\pK}(\alpha_{\pK})=0$, esto es,
$\alpha_{\pK}\in \*\o_{\pK}=U_{\pK}$, para casi toda $\pK\in{\ma P}_K$.

El conjunto de id\`eles $J_K$ forma un grupo con la multiplicaci\'on
entrada por entrada y la topolog\'ia de $J_K$ est\'a dada por el
sistema de vecindades abiertas de $\vec 1$ consiste de los subconjuntos
de la forma del tipo $\prod_{v\in S}W_v\times \prod_{v\notin S} U_v$ 
donde $S\subseteq {\ma P}_K$ es un conjunto finito, $W_v$ es un
conjunto abierto de $\*K_v$ donde $1\in W_v$ para $v\in S$ y
$U_v=\*\o_v$ son las unidades. Recordemos que definimos
$U_v=\begin{cases} \*\o_{K_{\pK}}=\o_{\pK}^*&\text{si $v$ es finito},\\
\*K_v&\text{si $v$ es infinito}
\end{cases}$.

Notemos en particular que $U_K:=\prod_{\pK\in{\ma P}_K}
U_{\pK}$ es un conjunto abierto de $J_K$.

\begin{observacion}\label{O17.6.24N}
El grupo de unidades de ${\ma A}_K$ es $J_K$ pues $\vec\alpha
\in \*{\ma A}_K$ si para toda $v\in{\ma P}_K$, existe $\beta_v\in
K_v$ con $\alpha_v\beta_v=1$. Como $\alpha_v\in\o_v$ para
casi toda $v$, $\alpha_v\in\*\o_v$ para toda $v$. Esto es
$\*{\ma A}_K=J_K$.
\end{observacion}

\begin{definicion}\label{D17.6.25N}
Los {\em id\`eles principales\index{ideles principales@id\`eles
principales}} son los elementos de $J_K$ de la forma $(\ldots,x,
\ldots,x,\ldots)$ con $x\in\*K$.

Los {\em ad\`eles principales\index{adeles principales@ad\`ees
principales}} son los ad\`eles de la forma $(\ldots,x,\ldots,x,\ldots)$
con $x\in K$.
\end{definicion}

La siguiente definici\'on es el concepto central de la teor\'ia global
de campos de clase.

\begin{definicion}\label{D17.6.26N}
El {\em grupo de clases de id\`eles\index{grupo de clases de
id\`eles}\index{clases de id\`eles}\index{ideles@id\`eles!grupo
de clases de $\sim$}}\label{CClaseCK} $C_K$ se define por $J_K/\*K$ donde
$\*K\xhookrightarrow{\ \ \mu\ \ }J_K$, $\vec {\mu(x)}=(\ldots,x,\ldots,x.
\ldots)$, esto es, $\vec{\mu(x)}_v=x$ para toda $v\in{\ma P}_K$.
\end{definicion}

Los elementos de $C_K$ los denotaremos por $\idel \alpha$
donde $\vec\alpha\in J_K$ es el representante del elemento de $C_K$.

\begin{notacion}\label{N17.6.27N} Dado $\pK\in {\ma P}_K$,
donde $K$ es un campo global, $\pK|\infty$ significa que
$\pK$ es arquimediano y $\pK\nmid \infty$ significa que 
$\pK$ es un primo finito.
\end{notacion}

\begin{definicion}[Valor absoluto de id\`eles\index{valor absoluto
de id\'eles}]\label{D17.6.28N}
Se define el {\em valor absoluto} de $J_K$ por $\|\ \|:J_K\lra {\ma R}^+
=\{x\in{\ma R}\mid x>0\}$, por
\[
\|\vec\alpha\|:=\prod_{v\in{\ma P}_K}|\alpha_v|_v.
\]
El valor absoluto est\'a bien definido pues $|\alpha_v|_v=1$ para casi
toda $v$. Adem\'as es f\'acil verificar que el valor absoluto es un mapeo
continuo, donde ${\ma R}^+$ tiene la topolog\'ia usual de los n\'umeros
reales y la im\'agen de $\|\ \|$ tiene la topolog\'ia inducida.
\end{definicion}

Se tiene que si $\vec\alpha$ y $\vec\beta$ est\'an en la misma clase
de equivalencia en $C_K$, entonces existe $x\in\*K$ tal que
$\vec\alpha=x\vec\beta$ y como $\|x\|=1$, se sigue que $\|
\vec\alpha\|=\|\vec\beta\|$ por lo que se puede definir $\|\idel \alpha\|$
por
\[
\|\idel \alpha\|=\|\vec \alpha\|.
\]

Ahora bien, $\ker \|\ \|=\{\vec\alpha\in J_K\mid \|\vec \alpha\|=
\prod_{v\in{\ma P}_K}|\alpha_v|_v=1\}=:J_{K,0}$. 
\label{CClaseidelesgrado0} Se tiene que
$\*K\subseteq J_{K,0}$. 

\begin{definicion}\label{D17.6.29N} El grupo $J_{K,0}$ recibe el nombre
de {\em id\`eles de grado $0$\index{grupo de id\`eles de grado $0$}}
y el grupo $C_{K,0}:=J_{K,0}/\*K$ recibe el nombre de {\em grupo de
clases de id\`eles de grado $0$\index{grupo de clases de id\`eles de
grado $0$}\label{CClaseCK0}}.
\end{definicion}

Si $K$ es un campo de funciones, entonces $\|\vec\alpha\|=q^{-\deg \vec\alpha}$,
es decir, 
\[
J_{K,0}=\{\vec\alpha\in J_K\mid\deg\vec\alpha
\label{CClasegrado global}=\sum_{\pK\in{\ma P}_K}
\deg_{\pK}\alpha_{\pK}=0\},
\]
donde $\deg_{\pK}\alpha_{\pK}=\deg\pK\cdot v_{\pK}
(\alpha_{\pK})\label{CClasegradolocal}$,
$\deg \pK:=\big[\o_{\pK}/\pK:\F]$, $\F$ el campo de constantes de $K$.

Para $K$ un campo global, se tiene 
\begin{align*}
C_K/C_{K,0}&\cong J_K/J_{K,0}\cong \im \|\ \|=:\Delta\\
&\cong
\begin{cases} {\ma R}^+&\text{si $K$ es num\'erico},\\
q^{\ma Z}=\{q^n\mid n\in{\ma Z}\}\cong {\ma Z}&\text{si
$K$ es campo de funciones}
\end{cases}.
\end{align*}

Veremos m\'as adelante con detalle que las sucesiones exactas
\begin{gather*}
1\lra J_{K,0}\lra J_K\lra \Delta\lra 1,\qquad 1\lra C_{K,0}\lra
C_K\lra\Delta\lra 1,
\intertext{se escinden tanto algebraica como topol\'ogicamente, es decir,
tenemos que}
C_K\cong C_{K,0}\times \Delta\qquad\text{y}\qquad 
J_K\cong J_{K,0}\times \Delta
\end{gather*}
tanto algebraica como topol\'ogicamente, donde $\Delta={\ma R}^+$
si $K$ es num\'erico y $\Delta \cong{\ma Z}$ si $K$ es de funciones.
La topolog\'ia inducida de $\im\|\ \|$ en ${\ma Z}$ es la topolog\'ia
discreta (Teorema \ref{T17.6.152N}).

La conexi\'on de los id\`eles con los divisores y de
varios grupos relacionados, nos la da la siguiente definici\'on.

\begin{definicion}\label{D17.6.30N}
Sea $K$ un campo global, $J_K$ el grupo de id\`eles y $D_K$
el grupo de divisores o de ideales fraccionarios. Se define
el epimorfismo\label{CClaseidealdeidele}
\[
\Lambda:J_K\lra D_K, \qquad \Lambda(\vec\alpha)={\eu a}_{\vec\alpha}=
\prod_{\substack{\pK\in{\ma P}_K\\ \pK\nmid \infty}}
\pK^{v_{\pK}(\alpha_{\pK})}.
\]
\end{definicion}

Se tiene que 
\begin{align*}
\ker\Lambda&=U_K:=\prod_{\pK\in{\ma P}_K}U_{\pK}=
\{\vec\alpha\in J_K\mid v_{\pK}(\alpha_{\pK})=0\text{\ para todo $\pK$
finito}\}\\
&=\begin{cases}
\prod_{\pK|\infty}\*K_{\pK}\times \prod_{\pK\nmid \infty}U_{\pK}&
\text{si $K$ es num\'erico},\\ \\
\prod_{\pK} U_{\pK}&\text{si $K$ es de funciones}
\end{cases}.
\end{align*}

Cuando $K$ es campo de funciones, se tiene $\Lambda(J_{K,0})
=D_{K,0}$ y $U_K\subseteq J_{K,0}$. Entonces se tiene
\[
\frac{J_K}{U_K}\cong D_K\text{\ ($K$ campo global)}, \qquad \frac{
J_{K,0}}{U_K} \cong D_{K,0}\text{\ ($K$ campo de funciones)}.
\]
Si $K$ es campo num\'erico, se tiene $\Lambda(J_{K,0})=D_K=
D_{K,0}$. Tambi\'en tenemos para cualquier $K$ campo global que
$\Lambda(\*K)=P_K=\{(x)_K\mid x\in\*K\}$. Se sigue que
\begin{gather*}
\Lambda:J_K/\*K\lra D_K/P_K\qquad\text{y}\qquad \ker \Lambda=
U_K\*K/\*K,\\
\Lambda:J_{K,0}/\*K\lra D_{K,0}/P_K\qquad\text{y}\qquad \ker \Lambda=
U_K\*K/\*K.
\intertext{En particular, tenemos}
J_K/U_K\*K\cong D_K/P_K=I_K \text{\ grupo de clases de divisores},\\
J_{K,0}/U_K\*K\cong D_{K,0}/P_K=I_{K,0} \text{\ grupo de clases de divisores
de grado $0$}.
\end{gather*}
Si $K$ es num\'erico, $D_K/P_K=D_{K,0}/P_K=I_K=I_{K,0}$ es el grupo de
clases y $h_K=|I_K|<\infty$. En el caso de $K$ campo de funciones
$h_K=|I_{K,0}|$ es finito y $I_K\cong I_{K,0}\oplus {\ma Z}$ es infinito.

\subsubsection{Algo de topolog\'ia para los id\`eles}\label{S17.6.4N}

Sea $K$ un campo global. Se tiene que $J_K$ es Hausdorff: si
$\vec\alpha\neq\vec\beta$, entonces $\alpha_v\neq\beta_v$
para alg\'un $v$ y $\*K_v$ es Hausdorff (simplemente por ser espacio
m\'etrico). Sean $U$ y $V$ conjuntos abiertos de $\*K_v$, $\alpha_v
\in U, \beta_v\in V$ y $U\cap V=\emptyset$. Sean $U_0$ y $V_0$ dos
abiertos de $J_K$ con la componente de $U_0$ en $v$ igual a $U$
y la componente en $v$ en $V_0$ igual a $V$ y $\vec\alpha\in U_0$ y
$\vec\beta\in V_0$ y $U_0\cap V_0=\emptyset$.

Sea $S_{\infty}=\{\pK\in{\ma P}_K\mid\pK|\infty\}$ el conjunto de los
lugares arquimedianos de $K$ y sea $S=S_{\infty}\cup \{v_0\}$ un
conjunto no vac\'io finito de lugares. Ahora, cada $\*K_v$ es localmente
compacto. Sea $V_v$ un abierto de $\*K_v$ con $\bar V_v$ compacto y
tal que $1\in V_v$. Entonces $\vec 1\in W=\prod_{v\in S}
V_v\times \prod_{v\notin
S} U_v$ y $U_v$ es compacto. Por tanto $W$ es una vecindad de
$\vec 1$ con cerradura compacta y por tanto $J_K$ es localmente compacto.

Veamos a continuaci\'on una propiedad de $J_K$ que nos es indispensable
para probar el teorema de existencia. La propiedad en discusi\'on es
que si $K$ es un campo de funciones, $J_K$ es totalmente disconexo
y en general si $K$ es un campo global $J_K^{\text{fin}}:=
\{\vec\alpha\in J_K\mid \alpha_v=1\text{\ para\ } v|\infty\}$ es
totalmente disconexo.

En efecto, sean $\vec\alpha, \vec\beta\in J_K^{\text{fin}}$ con $\vec\alpha
\neq \vec\beta$. Existe $v_0$ finito con $\alpha_{v_0}\neq \beta_{v_0}$.
Entonces $\alpha_{v_0},\beta_{v_0}\in \*K_{v_0}$ y como $\*K_{v_0}$
es totalmente disconexo, existen $V_0, W_0$ abiertos de $\*K_{v_0}$,
$\*K_{v_0}=V_0\cup W_0$, $V_0\cap W_0=\emptyset$ y $\alpha_{v_0}
\in V_0$, $\beta_{v_0}\in W_0$. Por tanto $\big(V_0\times
\prod_{v\neq v_0}\*K_v\big)\cap \big(W_0\times \prod_{v\neq v_0}
\*K_v\big)=\emptyset$ y el resultado se sigue.

En resumen, tenemos

\begin{proposicion}\label{P17.6.31N}
Sea $K$ un campo global, entonces $J_K$ es Hausdorff, localmente
compacto y $J_K^{\text{fin}}$ es totalmente disconexo, donde
$J_K^{\text{fin}}=\{\vec\alpha\mid \alpha_v=1 \text{\ para\ } v|\infty\}$. $\fin$
\end{proposicion}

\begin{corolario}\label{C17.6.32N}
Si $K$ es un campo de funciones, entonces $C_K$
y $C_{K,0}$ son totalmente disconexos.
\end{corolario}

\begin{proof}
Se tiene que la proyecci\'on natural $J_K\twoheadrightarrow  C_K$ es
un mapeo continuo, suprayectivo y $J_K$ es totalmente disconexo.
$\fin$
\end{proof}

\begin{observacion}\label{O17.6.33N}
Se tendr\'a que el mapeo global de reciprocidad $\rho_K:C_K
\lra \Gal(\abe K/K)=\abe G_K$, tiene n\'ucleo ${\eu N}_K$ que es
la componente conexa de $\vec 1$ en $C_K$. Por tanto $\rho_K$
resultar\'a ser inyectiva en el caso de campos de funciones. 
En el caso de campos num\'ericos la estructura de ${\eu N}_K$
fue hallada por Tate. En el caso num\'erico, $\rho_K$ es 
suprayectiva y no inyectiva.
\end{observacion}

\begin{proposicion}\label{P17.6.34-1N}
Para cualquier campo global $K$, el subgrupo $C_{K,0}$ es
cerrado de $C_K$. Si $K$ es de funciones, entonces $C_{K,0}$
tambi\'en es abierto. Si $K$ es num\'erico, $C_{K,0}$ no es
abierto en $C_K$.
\end{proposicion}

\begin{proof}
Primero supongamos que $K$ es de funciones. Entonces $U:=
\prod_{v\in{\ma P}_K} U_v$ es abierto en $J_K$ por lo que $\tilde
U:= U \*K/\*K$ es abierto en $C_K$ y adem\'as $\tilde U\subseteq
C_{K,0}$. Si $\idel \alpha\in
C_{K,0}$, entonces $\idel \alpha \tilde U\subseteq C_{K,0}$ probando
que $C_{K,0}$ es abierto. El mismo argumento prueba que si $\idel
\beta\notin \tilde U$, entonces $\idel \beta \tilde U\cap C_{K,0}=
\emptyset$ probando que $C_{K,0}$ es cerrado.

De lo anterior obtenemos que $\deg:C_K\lra {\ma Z}$, donde 
${\ma Z}$ es considerado con la topolog\'ia discreta, es una funci\'on
continua pues $\deg^{-1}(\{0\})=C_{K,0}$ es tanto abierto como
cerrado.

Si $K$ es num\'erico, tenemos que si $\varphi:=\|\ \|$, entonces
$\varphi$ es continua y $\{1\}\in{\ma R}^+$ es cerrado, por lo que
$C_{K,0}=\varphi^{-1}(\{1\})$ es cerrado en $C_K$.

Es claro en el caso num\'erico que toda vecindad de $\idel 1$ 
contiene elementos de valor absoluto diferente de $1$ por lo que
$C_{K,0}$ no es abierto.
$\fin$
\end{proof}

\begin{corolario}\label{C17.6.34-2N}
Si $K$ es cualquier campo global, se tiene que $C_K$ no es compacto.
\end{corolario}

\begin{proof}
Primero consideremos $K$ num\'erico.
Sea $\varphi=\|\ \|$ el valor absoluto de $C_K$. Sea $\{U_i\}_{i\in I}$ una
cubierta abierta de ${\ma R}^+$ que no tenga subcubiertas finitas. 
Entonces $\big\{\varphi^{-1}(U_i)\big\}_{i\in I}$ es una cubierta abierta
de $C_K$ sin subcubiertas finitas.

Si $K$ es de funciones, sea $\idel \alpha_1$ de grado $1$. Entonces
$\bigcup\limits_{n\in{\ma Z}}\idel \alpha_1^n C_{K,0}$ es una cubierta abierta
de $C_K$ sin subcubiertas finitas.
$\fin$
\end{proof}

\begin{proposicion}\label{P17.6.34N}
Se tiene para un campo global $K$ que $\*K\subseteq J_K$ es un
subgrupo discreto de $J_K$. En particular, $\*K$ es cerrado en $J_K$.
\end{proposicion}

\begin{proof}
Sea $S\subseteq {\ma P}_K$ un conjunto finito y no vac\'io de lugares
conteniendo los lugares arquimedianos. Sea $V$ la vecindad de $\vec
1$ en $J_K$ definido por $\vec\alpha\in V\iff |\alpha_v-1|_v<1$ para
toda $v\in S$ y $|\alpha_v|_v=1$ para toda $v\notin S$. Si $x\in\*K$,
$x\neq 1$ se tiene que $\prod_{v\in{\ma P}_K}|x-1|_v=1$. Si $x\in V$
se tendr\'ia que $|x-1|_v<1$ para $v\in S$ y para $v\notin S$,
$\max\{|x|_v,1\}=\max\{1,1\}=1$ y $|x-1|_v\leq \max\{|x|_v,1\}=1$.
Se seguir\'ia que 
\[
\prod_{v\in{\ma P}_K}|x-1|_v=\prod_{v\in S}|x-1|_v\cdot 
\prod_{v\notin S}|x-1|_v\leq \prod_{v\in S}|x-1|_v<1.
\]
Esta contradicci\'on a la f\'ormula del producto 
prueba que $x\notin V$. Por tanto $V\cap \*K=
\{1\}$ y por tanto $\*K$ es discreto.

Por el Lema \ref{discretocerrado}, $\*K$ es cerrado pues $J_K$ es
Hausdorff.
$\fin$
\end{proof}

\begin{observacion}\label{O17.6.35N}
De manera similar, se prueba que $K$ es un subgrupo discreto
de ${\ma A}_K$
\end{observacion}

\begin{observacion}\label{O17.6.36N}
Se tiene que $\*{\ma A}_K=J_K$ (Observaci\'on \ref{O17.6.24N}),
pero la topolog\'ia de $J_K$ no es la inducida por ${\ma A}_K$.
\end{observacion}

\begin{ejemplo}\label{E17.6.37N}
Sea $K={\ma Q}$. Sea $\vec x_p$ el ad\`ele cuya $p$-componente 
es $p$ y cuya $v$-componente es $1$ para $v\neq p$. Veamos
que $\vec x_p\xrightarrow[p\to\infty]{} \vec 1$ en la topolog\'ia de
${\ma A}_{\ma Q}$. 

Sea $V$ un abierto b\'asico conteniendo a $\vec 1$. Puesto que 
$V=\prod_{v\in S}W_v\times \prod_{v\notin S}\o_v$, $W_v$ abierto
de $K_v$ para $v\in S$. Entonces, existe $p_0$ un n\'umero primo
tal que para toda $p\geq p_0$ se tiene que $p\notin S$.
De esta forma, tenemos que $p\in \o_p={\ma Z}_p$. Se sigue que
$\vec x_p\in V$ para toda $p\geq p_0$ y por tanto $\lim\limits_{
p\to\infty} \vec x_p=\vec 1$.

Si la inversi\'on fuese continua en ${\ma A}_{\ma Q}$, entonces la
sucesi\'on $\vec x_p^{-1}$ converger\'ia a $\vec 1^{-1}=\vec 1$, 
sin embargo, si $V$ es un conjunto abierto conteniendo a $\vec 1$,
entonces para toda $p\geq p_0$, $p^{-1}\notin \o_p={\ma Z}_p$.
Por tanto $\vec x_p^{-1}\notin V$ para toda $p\geq p_0$.

Esto prueba que el mapeo $\vec x\lra \vec x^{-1}$ no es continua en la
topolog\'ia de ${\ma A}_{\ma Q}$ restringida a $\*{\ma A}_{\ma Q}=
J_{\ma Q}$. La raz\'on es que $p^{-1}$ no es unidad en ${\ma Z}_p$
y que $p\to\infty$.

Ahora bien, en la topolog\'ia de $J_K$, donde $K$ es un campo global
arbitrario, sea $\varphi:J_K\lra J_K$ dada por $\varphi(\vec \alpha)=
\vec \alpha^{-1}$. Sea ${\mc W}=\prod_{v\in S}W_v\times \prod_{v\notin S}
U_v$ un abierto b\'asico de $J_K$, por tanto $\varphi^{-1}({\mc W})=
\prod_{v\in S}W_v^{-1}\times \prod_{v\notin S}U_v^{-1}$ y si
$W_v$ es abierto en $\*K_v$, entonces $W_v^{-1}$
es abierto en $\*K_v$ por ser $\*K_v$ un grupo topol\'ogico y adem\'as
$U_v^{-1}=U_v$, por lo que $\varphi^{-1}$ es un abierto en $J_K$
y $\varphi$ es continua.

Por tanto $J_K$ es un grupo topol\'ogico. Se puede probar que la
topolog\'ia dada a $J_K$ est\'a dada de la siguiente forma:
Sea $\mu:J_K\hooklongrightarrow {\ma A}_K\times {\ma A}_K$ dada
por $\mu(\vec\alpha)=\big(\vec\alpha,\frac {1}{\vec\alpha}\big)$.
Entonces los abiertos de $J_K$ son $\mu(J_K)\cap (V\times W)$
donde $V$ y $W$ son conjuntos abiertos de ${\ma A}_K$.
\end{ejemplo}

\begin{lema}\label{L17.6.38N}
Sea $K$ un campo global y sea $L/K$ una extensi\'on finita y separable.
Entonces
\[
{\ma A}_K\otimes_K L\cong {\ma A}_L
\]
algebraica y topol\'ogicamente. En esta correspondencia, tenemos que
$K\otimes_K L=L\subseteq {\ma A}_K \otimes_K
L$ donde $K\subseteq {\ma A}_K$
se inyecta en $L\subseteq {\ma A}_L$.
\end{lema}

\begin{proof}
Primero estableceremos el isomorfismo como espacios topol\'ogicos. Sea
$\{\alpha_1,\ldots,\alpha_n\}$ base de $L/K$ y sea $v$ que recorre los 
valores absolutos finitos normalizados de $K$. Se tiene que, ${\ma A}_K
\otimes_K L$ con la topolog\'ia del producto tensorial, es el producto
restringido de
\begin{gather}
K_v\otimes_K L=K_v\alpha_1\oplus\cdots\oplus K_v\alpha_n  \label{Eq17.6.1N}
\intertext{con respecto a}
\o_K\alpha_1\oplus\cdots\oplus \o_v\alpha_n\label{Eq17.6.2N}.
\intertext{Ahora bien, se tiene}
K_v\otimes_K L=L_{\omega_1}\oplus\cdots\oplus L_{\omega_r}, \label{Eq17.6.3N}
\end{gather}
donde $\omega_1,\ldots,\omega_r$ son los valores absolutos
que son extensiones de $v$ y bajo esta identificaci\'on de (\ref{Eq17.6.1N}) con
(\ref{Eq17.6.3N}), se identifica (\ref{Eq17.6.2N}) con $\o_{\omega_1}\oplus\cdots
\oplus \o_{\omega_r}$. As\'i, el producto restringido de ${\ma A}_K\otimes_K L$ con
respecto a $K_v\otimes_K L$ es isomorfo topol\'ogicamente al producto restringido
de los $L_{\omega}$ con respecto a $\o_{\omega}$. Este isomorfismo tambi\'en
es algebraico.
$\fin$
\end{proof}

\begin{corolario}\label{C17.6.39N}
Sea ${\ma A}_K^+$ la parte aditiva de ${\ma A}_K$. Entonces
\[
{\ma A}_L^+\cong \underbrace{{\ma A}_K^+\oplus\cdots\oplus {\ma A}_K^+}_{
n=[L:K]}.
\]
En este isomorfismo se tiene $L\subseteq {\ma A}_L^+$, el subgrupo de los ad\`eles
principales, se manda en $K\oplus\cdots\oplus K$.
\end{corolario}

\begin{proof}
Para $x\in L$, $x\neq 0$, $x{\ma A}_K^+\subseteq {\ma A}_L^+$ y $x{\ma A}_K
\cong {\ma A}_K^+$ como grupos topol\'ogicos. Por tanto
\begin{gather*}
{\ma A}_L^+\cong{\ma A}_K^+\otimes_K L\cong\alpha_1{\ma A}_K^+\oplus\cdots
\oplus \alpha_n{\ma A}_K^+\cong{\ma A}_K^+\oplus\cdots
\oplus {\ma A}_K^+. \tag*{$\fin$}
\end{gather*}
\end{proof}

Del Corolario \ref{C17.6.34-2N}, tenemos que $J_K/\*K=C_K$
no es compacto. Sin embargo, para los ad\`eles, tenemos:

\begin{teorema}\label{T17.6.40N}
$K$ es discreto en ${\ma A}_K$ y ${\ma A}_K^+/K$ es compacto en la
topolog\'ia cociente.
\end{teorema}

\begin{proof}
Por el Corolario \ref{C17.6.39N}, basta probar el teorema para $K={\ma Q}$
en el caso num\'erico y para $K=\F(T)$ en el caso de campos de funciones.

Denotamos por $\infty$ al valor absoluto usual si $K={\ma Q}$ y por el
valor absoluto a $\frac 1T$ si $K=\F(T)$. Sea $W\subseteq {\ma A}_K^+$
el conjunto definido por $W=\{\vec \alpha\in{\ma A}_K^+\mid |\alpha_{\infty}
|_{\infty}\leq \frac 12, |\alpha_v|_v\leq 1\text{\ para todo $v\neq \infty$}\}$.
Veamos que ${\ma A}_K^+=K+W$.

Sea $\vec \alpha\in {\ma A}_K^+$. Para cada $v\neq \infty$, sea $v=p\in
{\ma Z}$ un n\'umero primo o $p\in R_T^+$, donde $R_T^+$ denota
al conjunto de los polinomios m\'onicos irreducibles de $R_T=\F[T]$.
Podemos hallar $r_v=r_p =z_p/p^{x_p}$ con $x_p\in {\ma Z}$ o $R_T$ y
y $x_p\in{\ma Z}$, $x_p\geq 0$ tal que
\[
|\alpha_v-r_v|_v=|\alpha_p-r_p|_p\leq 1.
\]
Puesto que $\vec\alpha\in {\ma A}_K^+$, se tiene $r_v=0$ para casi
toda $v$. Sea $r:=\sum_{v\neq \infty} r_v$. Se tiene $|\alpha_v-r|_v
\leq 1$ para toda $v\neq \infty$.

Sea $s\in{\ma Z}$ o $s\in R_T$ con $|\alpha_{\infty}-r-s|_{\infty}\leq 
\frac 12$. Sea $\vec \beta=\vec \alpha-(r+s)$, $\vec\alpha=\vec \beta+
(r+s)$. Se tiene $\beta_v=\alpha_v-(r+s)$, $|\beta_v|_v=|\alpha_v-(r+s
)|_v\leq \max\{|\alpha_v-r|_v, |s|_v\}\leq 1$ para toda $v\neq \infty$.
Adem\'as $|\beta|_{\infty}=|\alpha_{\infty}-r-s|_{\infty}\leq \frac 12$.

Se tiene que $W\twoheadrightarrow {\ma A}_K^+/K$ es un 
epimorfismo continuo y $W$ es compacto,
de donde se sigue que ${\ma A}_K^+/K$ es compacto.
$\fin$
\end{proof}

\begin{corolario}\label{C17.6.43N}
Sea $K$ un campo global. Existe $W_0\subseteq {\ma A}_K$ 
compacto con
$W_0:=\{\vec\xi\in{\ma A}_K\mid |\xi_v|_v\leq \delta_v\text{\ para toda
$v\in {\ma P}_K$ y $\delta_v=1$ para casi toda $v$}\}$ y tal que
${\ma A}_K=W_0+K$.
\end{corolario}

\begin{proof}
Se sigue inmediatamente de la demostraci\'on del Teorema
\ref{T17.6.40N}.
$\fin$
\end{proof}

\subsubsection{Medida de Haar\index{Haar!medida de 
$\sim$}\index{medida de Haar}}\label{medida de Haar}

Sea $G$ un grupo topol\'ogico localmente compacto. Entonces,
salvo una constante positiva, existe una \'unica medida no
trivial $\mu$ en los subgrupos de Borel, esto es, los elementos
del \'algebra de Borel que es la generada por los conjuntos
abiertos de $G$, la cual es numerablemente aditiva y tal que
satisface lo siguiente:
\las
\item $\mu$ es invariante bajo translaciones por la izquierda:
$\mu(gS)=\mu(S)$ para toda $g\in G$ y para todo 
conjunto de Borel $S$.

\item $\mu$ es finita en conjuntos compactos $T\subseteq G$:
$\mu(T)<\infty$ para $T$ compacto.

\item $\mu$ es regular por fuera para todo subgrupo de Borel:
$\mu(S)=\inf\{\mu(U)\mid S\subseteq U, U \text{\ abierto}\}$.

\item $\mu$ es regular por dentro para conjuntos abiertos $U$:
$\mu(U)=\sup\{\mu(T)\mid T\subseteq U, T\text{\ compacto}\}$.
\end{list}

Se puede probar que $\mu(U)>0$ para todo conjunto abierto
$U$ de $G$ con $U$ no vac\'io. En particular, si $G$ es
compacto, $\mu(G)<\infty$ y se puede normalizar
para que $\mu(G)=1$.

De igual manera, se puede definir un \'unica medida (salvo
constantes positivas) por la derecha. Sin embargo estas dos
medidas no coinciden en general.

\begin{proposicion}\label{P17.6.44N}
Sea $K$ un campo global. Existe una constante $C>0$ que depende
\'unicamente de $K$ con la siguiente propiedad: si
$\vec\alpha\in{\ma A}_K$ satisface
\begin{gather}\label{Eq17.6.44-1N}
\prod_{v} |\alpha_v|_v>C,
\end{gather}
entonces existe $\beta\in K$, $\beta\neq 0$ tal que $|\beta|_v\leq |
\alpha|_v$ para toda $v\in {\ma P}_K$.
\end{proposicion}

\begin{proof}
Primero consideremos cualquier constante $C>0$ y un id\`ele $\vec
\alpha$ que satisface (\ref{Eq17.6.44-1N}). Puesto que para casi toda
$v$ se tiene $|\alpha_v|_v\leq 1$, entonces se sigue que $|\alpha_v|_v
=1$ para casi toda $v$ pues de lo contrario, si para $v$ finito se tiene
$|\alpha_v|_v<1$ entonces $|\alpha_v|_v\leq q_v^{-1}$ y por tanto si
tuvi\'esemos esta situaci\'on para una infinidad de primos finitos $v$,
entonces $\prod_v|\alpha_v|_v=0$.

Sean $c_0$ la medida de Haar $\mu$ de ${\ma A}_K^+/K$
y $c_1$ es la medida de Haar del conjunto $A=\{\vec\gamma\mid
|\gamma_v|_v\leq \frac {1}{10}\text{\ si $v|\infty$ y $|\gamma_v|_v\leq
1$ si $v\nmid\infty$}\}$. Se tiene que por ser ${\ma A}_K^+/K$
compacto, lo mismo que $A$ pues el n\'umero de $v|\infty$ es
finito, entonces $0<c_0<\infty$ y $0<c_1<\infty$. Veamos que
$C=\frac {c_0}{c_1}$ satisfacen las condiciones de la propisici\'on.

Sea $\vec\alpha$ que satisface (\ref{Eq17.6.44-1N}) con $C=\frac{
c_0}{c_1}$. El conjunto $T:=\{\vec\gamma\mid |\gamma_v|_v\leq
\frac 1{10}|\alpha_v|_v, V|\mid, |\tau_v|_v\leq |\alpha_v|_v, v
\nmid \infty\}$ tiene medida igual a $c_1\cdot\prod_v |\alpha_v|_v>
c_1 C=C_0=\mu\big({\ma A}_K^+/K\big)$. Se sigue que el
epimorfismo natural ${\ma A}_K^+\twoheadrightarrow {\ma A}_K^+/K$
debe tener dos elementos distintos $\vec\tau', \vec\tau''$ en $T$
con la misma im\'agen debido a que $\mu(T)>\mu
\big({\ma A}_K^+/K\big)$, y por tanto $\vec\tau'-\vec\tau''=\beta
\in K$. Por tanto, 
\begin{gather*} 
|\beta|_v=|\tau'_v-\tau''_v|_v\leq |\alpha_v|_v
\end{gather*}
para toda $v$.
$\fin$
\end{proof}

\begin{corolario}\label{C17.6.45N}
Sea $v_0$ un valor absoluto normalizado y sean $\delta_v>0$ 
dados para $v\neq v_0$ tales que $\delta_v=1$ para casi toda
$v$. Entonces existe $\beta\in K$ con $\beta\neq 0$ con $|\beta_v
|_v\leq \delta_v$ para todo $v\neq v_0$.
\end{corolario}

\begin{proof}
Seleccionamos $\alpha_v\in K_v$ con $0<|\alpha_v|_v\leq
\delta_v$ y $|\alpha_v|_v=1$ si $\delta_v=1$. Sea $\alpha_{v_0}
\in K_{v_0}$ tal que $\prod_{v\in{\ma P}_K}|\alpha_v|_v>C$.
Por la Proposici\'on \ref{P17.6.44N} tenemos que existe $\beta
\in K$, $\beta\neq 0$ que satisface lo requerido.
$\fin$
\end{proof}

A continuaci\'on presentamos el teorema de aproximaci\'on fuerte.
Una versi\'on para campos de funciones fue demostrada en
el Teorema \ref{T5.3.8} usando el Teorema de
Riemann-Roch, donde el campo de funciones es 
arbitrario y no \'unicamente global. El resultado que ahora
presentamos es v\'alido para cualquier campo global.

\begin{teorema}[Teorema de aproximaci\'on fuerte\index{teorema!de
aproximaci\'on fuerte}\index{aproximaci\'on
fuerte!teorema de $\sim$}]\label{T17.6.46N}
Sea $K$ un campo global y sea $v_0$ cualquier valor absoluto
de $K$. Sea $V$ el producto restringido de los $K_v$ con respecto
a $\o_v$ para toda $v\neq v_0$. Entonces $K$ es denso en $V$.

En otras palabras, dados $S=\{v_1,\ldots, v_n\}\subseteq {\ma P}_K$ con
$v_0\notin S$, $\epsilon_1,\ldots,\epsilon_n>0$ y $\vec\alpha\in
{\ma A}_K$, entonces existe $\beta\in K$ con
\[
|\beta-\alpha_{v_i}|_{v_i}<\epsilon_i,\quad 1\leq i\leq n\quad\text{y}
\quad |\beta|_v\leq 1\quad\text{para toda\ } v\notin S\cup\{v_0\}.
\]
\end{teorema}

\begin{proof}
Existe $W\subseteq {\ma A}_K$ con $\{\vec\xi\mid |\xi_v|_v\leq
\delta_v\text{\ y $\delta_v=1$ para casi toda $v$}\}$ y con ${\ma A}_K
=K+W$ (Corolario \ref{C17.6.43N}). Del Corolario \ref{C17.6.45N} 
tenemos que existe $\lambda\in K$, $\lambda\neq 0$ tal que
\begin{align*}
|\lambda|_v&<\delta_v^{-1}\epsilon \quad v\in S,\\
|\lambda|_v&\leq \delta_v^{-1},\quad v\notin S, v\neq v_0,
\end{align*}
con $\epsilon=\min\{\epsilon_1,\ldots, \epsilon_n\}$.

Dada $\vec \alpha\in{\ma A}_K$, si $\vec\gamma:=\lambda^{-1}\vec
\alpha$, existe $\vec \xi\in W$ tal que $\vec\gamma=\vec\xi+\nu$
para alg\'un $\nu\in K$. Se sigue que $\vec\alpha=\lambda\vec\xi+
\lambda\nu$. Por tanto ${\ma A}_K=\lambda W+K$.

Sea $\vec\mu\in {\ma A}_K$ definido por:
 $(\vec\mu)_{v_i}=\alpha_{v_i}$, $1\leq i\leq n$ y $(\vec\mu)_v
=0$ para $v\notin S$. Sean $\vec\xi\in W$ y $\beta\in K$ tales que
$\vec\mu=\lambda\vec\xi+\beta$. Se tiene
\begin{align*}
|\beta-\alpha_{v_i}|_{v_i}&=|\lambda\xi_{v_i}|_{v_i}\leq |\lambda|_{v_i}
\delta_{v_i}<\epsilon\quad\text{y}\\
|\beta-\alpha_v|_v&=|\beta|_v=|\lambda \xi_v|_v=|\lambda|_v
|\xi_v|_v\leq \delta_v^{-1}\delta_v=1, v\notin S\cup\{v_0\}. \tag*{$\fin$}
\end{align*}
\end{proof}

\begin{lema}\label{L17.6.47N} El subgrupo $J_{K,0}$ es cerrado
en ${\ma A}_K$ en la topolog\'ia de ${\ma A}_K$. Por otro lado, las
topolog\'ias inducidas, tanto de ${\ma A}_K$ como de $J_K$ en
$J_{K.0}$.
\end{lema}

\begin{proof}
Sea $\vec\alpha\in{\ma A}_K$ tal que $\vec\alpha\notin J_{K,0}$.
Debemos hallar una vecindad $W$ de $\vec\alpha$ con $W\cap
J_{K,0}=\emptyset$. Se tiene $|\alpha_v|_v\leq 1$ para casi toda
$v\in {\ma P}_K$ y por tanto el producto infinito $C:=\prod_{v\in
{\ma P}_K}|\alpha_v|_v\in \{x\in{\ma R}\mid x\geq 0\}$ puesto que
en la serie $\sum_v\log |\alpha_v|_v$ se tiene que $\log|\alpha_v
|_v\leq 0$ para casi toda $v$, esto es, $\sum_{v}\log
|\alpha_v|_v\neq \infty$.

(a) Caso en que $C<1$. Puesto que $\vec\alpha$ es un ad\`ele,
\'unicamente un n\'umero finito de lugares de $K$ satisfacen que
$|\alpha_v|_v\geq 1$. Puesto que $C<1$, existe un conjunto 
finito $S$ de lugares de $K$ tal que si $|\alpha_v|_v\geq 1$,
entonces $v\in S$ y se tiene adem\'as $\prod_{v\in S}|\alpha_v
|_v<1$. Sea $W:=\{\vec\xi\mid |\xi_v-\alpha_v|_v<\epsilon, 
v\in S, |\xi_v|_v\leq 1, v\notin S\}$ con $\epsilon$ suficientemente
peque\~no. Se tiene que si $\vec\xi\in W$, entonces
\[
\prod_{v\in{\ma P}_K}|\xi_v|_v=\prod_{v\in S}|\xi_v|_v\cdot
\prod_{v\notin S}|\xi_v|_v\leq \prod_{v\in S}\big(|\alpha_v|_v
+\epsilon\big)<1.
\]
Se sigue que $W\cap J_{K,0}=\emptyset$ y $\vec\alpha \in W$.

\s

(b) Caso en que $C\geq 1$. Puesto que $\vec\alpha\notin J_{K,0}$,
$C>1$. Probemos primero que, en este caso, $\vec\alpha\in J_K$.
Para $v$ finito, si $|x|_v<1$, entonces $|x|_v\leq \frac {1}{q_v}
=q_v^{-1}$, $q_v$ la cardinalidad del campo residual de $K_v$.

De esta forma, existe un conjunto finito $S$ de lugares, tales que 
si $v\notin S$, entonces $|\alpha_v|_v=1$ o $|\alpha_v|_v\leq \frac 12$
pues como $\vec\alpha$ es un ad\`ele, $|\alpha_v|_v\leq 1$ para
$v\notin S$ ($\alpha_v\in\o_v$).

Puesto que $\prod_{v\in{\ma P}_K}|\alpha_v|_v$ converge a $C>0$
($C\geq 1$), $|\alpha_v|_v\xrightarrow[v\to\infty]{}1$.  Por tanto, para un
n\'umero finito de lugares $v$, se tiene $|\alpha_v|_v\leq \frac 12$. As\'i,
agregando un n\'umero finito de lugares, podemos tomar $S$ tal que
$|\alpha_v|_v=1$ para toda $v\notin S$ y por tanto $\vec\alpha\in J_K$
pues $\alpha_v\neq 0$ para todo $v\in{\ma P}_K$ y $C>1$.

Ampliando nuevamente $S$ en caso de ser necesario, podemos
suponer que si $v\notin S$ y $|\alpha_v|_v<1$ entonces $|\alpha_v
|_v<(2C)^{-1}$ pues el n\'umero de lugares de grado menor o igual
a $n$ es finito y $q_v=q^{n_v}$, por lo que tomamos $\frac 1{q_v}
<\frac 1{2C}$.
Nuevamente, ampliando $S$, podemos suponer que
\[
1<\prod_{v\in S}|\alpha_v|_v<2C.
\]

Para $\epsilon >0$ peque\~no, si $|\xi_v-\alpha_v|_v<\epsilon$
para $v\in S$ se tiene
\begin{gather*}
-\epsilon+|\alpha_v|_v<|\xi_v|_v<\epsilon +|\alpha_v|_v\quad \text{y}\\
\prod_{v\in S}\big(|\alpha_v|_v-\epsilon\big)<\prod_{v\in S} |\xi_v|_v<
\prod_{v\in S}\big(\epsilon+|\alpha_v|_v\big).
\intertext{Tomamos $\epsilon$ suficientemente peque\~no de tal forma que}
\prod_{v\in S}\big(|\alpha_v|_v-\epsilon\big)>1 \quad\text{y}\quad
\prod_{v\in S}\big(\epsilon+|\alpha_v|_v\big)<2C.
\end{gather*}

As\'i, sea $W:=\{\vec\xi\mid |\xi_v-\alpha_v|_v<\epsilon, v\in S \text{\ y\ }
|\xi_v|_v\leq 1, v\notin S\}$. Si $\vec \xi\in W$, $\|\vec\xi\|\neq 1$ pues
si $v\notin S$, $|\xi_v|_v=1$ o $|\xi_v|_v<(2C)^{-1}$. As\'i, si existe
$v_0\notin S$ con $|\xi_{v_0}|_{v_0}<(2C)^{-1}$, entonces
\[
\|\vec\xi\|=\prod_{v\in S}|\xi_v|_v\cdot |\xi_{v_0}|_{v_0}\cdot 
\prod_{v\notin S\cup\{v_0\}}|\xi_v|_v< 2C\cdot \frac 1{2C}\cdot 1=1,
\]
esto es $\|\vec\xi\|<1$.

Si para toda $v\notin S$, $|\xi_v|_v=1$, $\|\vec\xi\|=\prod_{v\in S}
|\xi_v|_v>1$. Por tanto $\|\vec\xi\|\neq 1$ y $\vec\xi
\notin J_{K,0}$, $\vec\alpha\in W$
abierto en ${\ma A}_K$ y $W\cap J_{K,0}=\emptyset$. Por tanto
$J_{K,0}$ es cerrado en ${\ma A}_K$.

Sea $\vec\alpha\in J_{K,0}$, en particular $\alpha_v\in \*\o_v
\subseteq \o_v$ para casi toda $v$. Una vecindad de $\vec\alpha$
en ${\ma A}_K$
es una  vecindad de la forma $\vec\alpha+W$ con $W=\prod_{v\in
S}V_v\times \prod_{v\notin S}\o_v$ con $S$ un conjunto finito y
$V_v$ un conjunto abierto de $K_v$. Queremos hallar un conjunto
$W_0=\prod_{v\in T}X_v\times \prod_{v\notin T}\*\o_v$ con $X_v$
un conjunto abierto de $\*K_v$, $T$ un conjunto finito, y tal que
se tenga $\vec\alpha W_0\subseteq \vec\alpha+W$, $\vec\alpha
W_0$ es una vecindad de $\vec\alpha$ en $J_K$.

Pongamos $W=\prod_v Y_v$. Sea $W_0=\prod_v X_v$. Se 
necesita $\alpha_v x_v=\alpha_v+y_v$, lo cual es equivalente a
que $x_v=1+\alpha_v^{-1}y_v\in (1+\alpha_v^{-1} Y_v)\cap \*K_v$.
Esto es, $X_v\subseteq (1+\alpha_v^{-1}Y_v)\cap \*K_v$ para toda $v$.
Como $1+\alpha_v^{-1}Y_v$ es abierto en $K_v$, se tiene que
$(1+\alpha_v^{-1}Y_v)\cap \*K_v$ es abierto en $\*K_v$ 
para toda $v$ y adem\'as claramente
$(1+\alpha_v^{-1}Y_v)\cap \*K_v\neq \emptyset$. 

Por otro lado, se tiene $\alpha_v\in\*\o_v$
y $Y_v=\o_v$ para casi toda $v$. Se sigue que $(1+\alpha_v^{-1}
Y_v)=1+\o_v=\o_v$ para casi toda $v$. Podemos definir
$X_v=(1+\alpha_v^{-1}Y_v)\cap \*\o_v=\*\o_v$ para casi toda $v$ y
$X_v=(1+\alpha_v^{-1}Y_v)\cap \*K_v$ para el resto.
De esto se sigue que la ${\ma A}_K$-topolog\'ia de $J_{K,0}$
contiene a la $J_K$-topolog\'ia.

Rec\'iprocamente, sea ahora $H$ una $J_K$-vecindad de $\vec
\alpha$. Entonces $H$ contiene una $J_K$-vecindad de la forma
$|\xi_v-\alpha_v|_v<\epsilon$ para $v\in S$ y $|\xi_v|_v=1$ para
$v\notin S$ y donde $S$ contiene a todos los $v$ arquimedianos:
$v|\infty$ y a todos los $v$ tales que $|\alpha_v|_v\neq 1$. Puesto
que $\vec\alpha\in J_{K,0}$, se tiene $\|\vec\alpha\|=\prod_{v}
|\alpha_v|_v=1$. Se tiene que para $v\notin S$, si $|\xi_v|_v<1$,
entonces $|\xi_v|_v\leq \frac 12$ por ser $v$ un lugar finito.
Sea $\epsilon$ tal que para $\vec\xi\in J_{K,0}\cap H$ se tiene $|\xi_v-
\alpha_v|_v<\epsilon$ para $v\in S$ y $|\xi_v|_v\leq 1$ para
$v\notin S$ de tal forma que se satisface
\[
\prod_{v\in S}|\xi_v|_v<\prod_{v\in S} \big(|\alpha_v|_v+
\epsilon\big)<2.
\]
Entonces $|\xi_v|_v=1$ para $v\notin S$ pues de lo contrario
existir\'ia $v\notin S$ con $|\xi_v|_v\leq \frac 12$ y por tanto
\[
1=\|\vec\xi\|=\prod_v |\xi_v|_v=\prod_{v\in S}|\xi_v|_v\cdot
\prod_{v\notin S} |\xi_v|_v<2\cdot\frac 12=1,
\]
lo cual es absurdo. Se sigue que $\xi_v\in\*\o_v$.
Por tanto $H\cap J_{K,0}$ contiene una ${\ma A}_K$-vecindad
de $\vec\alpha$ de $J_{K,0}$, esto es, contiene a $\vec\alpha
+H_0$, $H_0:=\prod_{v\in S}V_v\times \prod_{v\notin S}\o_v$.
$\fin$
\end{proof}

\begin{teorema}\label{T17.6.48N}
Sea $K$ un campo global. Se tiene que el 
grupo $C_{K,0}=J_{K,0}/\*K$ es compacto con respecto
a la topolog\'ia cociente.
\end{teorema}

\begin{proof}
Por el Lema \ref{L17.6.47N} es suficiente hallar un conjunto
$W\subseteq {\ma A}_K$ compacto en la topolog\'ia de ${\ma A}_K$
tal que la proyecci\'on natural $W\cap J_{K,0}\twoheadrightarrow
J_{K,0}/*K$ sea suprayectiva.

Sea $W=\{\vec\xi\in{\ma A}_K\mid |\xi_v|_v\leq |\alpha_v|_v\}$ donde
$\vec\alpha$ es cualquier id\`ele de norma mayor a la constante $C$
de la Proposici\'on \ref{P17.6.44N}. Se tiene que $W$ es compacto.
Sea $\vec\beta\in J_{K,0}$. Se tiene que
\begin{gather*}
\|\vec\beta^{-1}\vec\alpha\|=\|\vec\beta^{-1}\|\|\vec\alpha\|>1\cdot C=C,
\intertext{por lo que existe $\mu\in\*K$ tal que}
|\mu|_v\leq |\beta_v^{-1}\alpha_v|_v\quad \text{para toda $v$}.
\end{gather*}

Se sigue que $\mu\vec\beta\in W$ y $\mu\vec\beta\longmapsto \mu
\vec\beta \bmod \*K=\vec\beta\bmod \*K$.
$\fin$
\end{proof}

Sea $K$ un campo num\'erico y sea $\Lambda: J_K\lra D_K$
dada por $\Lambda
(\vec\alpha)={\eu a}_{\vec\alpha}=\prod_{\pK\nmid \infty} \pK^{v_{
\pK}(\alpha_{\pK})}$, donde $D_K$, el grupo de ideales fraccionarios,
tiene la topolog\'ia discreta. Se tiene $\Lambda(\*K)=P_K$ y
$D_K/P_K=I_K$ es el grupo de clases de $K$. Puesto que
hay al menos un lugar arquimediano, $\Lambda$ induce un
epimorfismo continuo $J_{K,0}/\*K\stackrel{\tilde\Lambda}{\lra}
I_K$.

Similarmente, si $K$ es un campo de funciones y $I_{K,0}$ denota
el grupo de clases de divisores de grado $0$,
\[
J_{K,0}/\*K\stackrel{\tilde\Lambda}{\lra} I_{K,0}
\]
es un epimorfismo continuo donde $D_{K,0}$ tiene la topolog\'ia
discreta. Se sigue que $I_K$ en el caso num\'erico y $I_{K,0}$
en el caso de campos de funciones son compactos y discretos,
por lo tanto finitos.

\begin{teorema}\label{T17.6.49N} 
Si $K$ es un campo num\'erico
y $I_K$ es el grupo de clases de $K$, $I_K=D_K/P_K$ es finito.

Si $K$ es un campo de funciones y $I_{K,0}=D_{K,0}/P_K$ es el
grupo de clases de divisores de grado $0$ de $K$, entonces
$I_{K,0}$ es finito. $\fin$
\end{teorema}

\begin{definicion}\label{D17.6.50N} Sea $S$ un conjunto finito de
lugares de un campo global que contiene a los primos infinitos. Se
define el {\em grupo de los $S$-id\`eles\index{grado de los 
$S$-id\`eles}} por
\[
J_{K,S}=\prod_{v\in S}\*K_v\times \prod_{v\notin S}\*\o_v=
\prod_{v\in S}\*K_v\times \prod_{v\notin S} U_v.
\]

Se tiene $J_K=\bigcup_{S\text{\ finito}} J_{K,S}$.
\end{definicion}

\begin{proposicion}\label{P17.6.51N}
Si $K$ es un campo num\'erico y $S=S_{\infty}=\{v\in{\ma P}_K\mid 
v|\infty\}$, entonces $J_K/J_{K,S_{\infty}}=D_K$ y $J_K/J_{K,S_{\infty}}
\*K\cong D_K/P_K\cong I_K$.
\end{proposicion}

\begin{proof}
Nuevamente usamos el epimorfismo $\Lambda:J_K\lra D_K$, 
$\Lambda(\vec\alpha)={\eu a}_{\vec\alpha}=\prod_{\pK\nmid \infty}
\pK^{v_{\pK}(\alpha_{\pK})}$ y se tiene que $\ker \Lambda=
J_{K,S_{\infty}}$ de donde $D_K\cong J_K/J_{K,S_{\infty}}$.

Similarmente, $\tilde\Lambda:J_K\stackrel{\Lambda}{\lra} D_K
\stackrel{\pi}{\lra}D_K/P_K=I_K$ es un epimorfismo y
\begin{gather*}
\vec\alpha\in \ker \tilde\Lambda\iff \prod_{\pK\nmid \infty} \pK^{v_{
\pK}(\alpha_{\pK})}\in P_K,\\
\text{esto es,}\quad \prod_{\pK\nmid\infty}
\pK^{v_{\pK}(\alpha_{\pK})}=(x)=\prod_{\pK\nmid\infty}
\pK^{v_{\pK}(x)}, \quad x\in\*K\\
\iff v_{\pK}(\alpha_{\pK})=v_{\pK}(x) \quad\text{para toda}
\quad \pK\nmid \infty\\
\iff v_{\pK}(\alpha_{\pK}x^{-1})=0
\quad\text{para toda}\quad \pK\nmid\infty\\
\iff \vec\alpha x^{-1}\in J_{K,S_{\infty}}\iff \vec\alpha\in
J_{K,S_{\infty}}\*K.
\end{gather*}
Por tanto $\ker\tilde\Lambda=J_{K,S_{\infty}} \*K$.
$\fin$
\end{proof}

\begin{teorema}\label{T17.6.52N}
Sea $K$ un campo global. Existe $S$, un conjunto finito suficientemente
grande tal que $J_K=J_{K,S}\*K$ y 
\[
C_K=J_K/\*K=J_{K,S}\*K/\*K
\cong J_{K,S}/(J_{K,S}\cap \*K).
\]
\end{teorema}

\begin{proof}
Primero consideremos el caso de $K$ de un campo num\'erico.
Se tiene que $I_K=D_K/P_K$ es finito. Sean ${\eu a}_1,\ldots,
{\eu a}_n$ ideales que representen a todos los elementos de $I_K$.
Tomando la descomposici\'on de los ideales ${\eu a}_1,\ldots, {\eu a}_n$
en producto de primos, se obtiene un n\'umero finito de ideales primos
$\pK_1,\ldots,\pK_s$ que son los diversos divisores primos de los
ideales ${\eu a}_i$'s.

Sea $S$ cualquier conjunto finito que contenga a $\{\pK_1,\ldots,
\pK_s\}$ y a los primos $\pK|\infty$. Sea $\Lambda_0:J_K/J_{K,S_{
\infty}}\stackrel{\cong}{\lra} D_K$. Sea $\vec\alpha\in J_K$, 
$\Lambda_0(\bar{\vec\alpha})={\eu a}=\prod_{\pK\nmid \infty}\pK^{v_{
\pK}(\alpha_{\pK})}\in {\eu a}_i\cdot P_K$ para alg\'un $i$. Esto es,
${\eu a}={\eu a}_i \cdot (x)$ con $(x)\in P_K$, $x\in \*K$. Se tiene que
\[
\vec{\alpha'}=\vec\alpha\cdot x^{-1}\stackrel{\Lambda_0}{\lra} {\eu a}'
=\prod_{\pK\nmid\infty}\pK^{v_{\pK}(\vec{\alpha'})}=\prod_{\pK\nmid
\infty}\pK^{v_{\pK}(\alpha_{\pK})-v_{\pK}(x)}={\eu a}\cdot (x)^{-1}=
{\eu a}_i.
\]

Ahora, las componentes primas de ${\eu a}_i$ pertenecen a $S$, 
por lo que
$v_{\pK}({\eu a}')=v_{\pK}(\alpha'_{\pK})=0$ para $\pK\notin S$.
Se sigue que
$\vec{\alpha'}\in J_{K,S}$, de donde $\vec\alpha=\vec{\alpha'}
\cdot x\in J_{K,S}\*K$ y $J_K=J_{K,S} \*K$.

Ahora consideremos el caso de un campo de funciones con
campo de constantes $\F$. Sea $v_0$ un primo cualquiera de
${\ma P}_K$. Sea $\o_K:=\bigcap_{\pK\neq v_0}\o_{\pK}=\{x\in
K\mid v_{\pK}(x)\geq 0\text{\ para toda $\pK\neq v_0$}\}$. 
Entonces, por el Teorema de Riemann-Roch, existe $T\in K$
cuyo \'unico polo de $T$ es $v_0$. Se sigue que $\o_K$ es la 
cerradura entera de $\F[T]$ en $K$ y $\o_K$ es un dominio
Dedekind.

Sean ${\mc D}_K$ el grupo de ideales fraccionarios de $\o_K$
y ${\mc P}_K$ el subgrupo de ideales principales de ${\mc D}_K$.
Entonces ${\mc D}_K/{\mc P}_K$ es un grupo finito (Corolario
\ref{CRamDed1.2.6}). Sea 
\[
\Lambda_{v_0}:J_K/J_{K,\{v_0\}}\lra {\mc D}_K,\qquad \Lambda_{
v_0}(\bar{\vec\alpha})=\prod_{\pK\neq v_0} \pK^{v_{\pK}(\alpha_{\pK})},
\]
el cual es un isomorfismo y $\Lambda_{v_0}(\bar{\*K})={\mc P}_K$.
Por tanto $J_K/J_{K,\{v_0\}}\*K\cong {\mc D}_K/{\mc P}_K$ es finito.
A partir de este punto, se procede como en el caso num\'erico.
$\fin$
\end{proof}

Ahora consideremos nuevamente un campo global de funciones.
Recordemos que el epimorfismo $\Lambda$ nos induce un
epimorfismo continuo y $C_{K,0} =J_{K,0}/\*K\twoheadrightarrow 
J_{K,0}/U_K \*K\cong I_{K,0}$ y puesto que $C_{K,0}$ es
compacto, $I_{K,0}$ es compacto y discreto, por tanto finito.
Esto ya lo hab\'iamos obtenido en el Teorema \ref{T10.1.1.5}
por otros m\'etodos.

\begin{definicion}\label{D17.6.53N} Sea $S$ un conjunto finito
no vac\'io de lugares de un campo global $K$ conteniendo a los
lugares infinitos en el caso de ser $K$ num\'erico. Sea $K^S
=\{x\in\*K\mid v_{\pK}(x)=0\text{\ para toda $\pK\notin S$}\}=
\*K\cap J_{K,S}$ el grupo de las $S$-unidades de $K$.
\end{definicion}

\begin{proposicion}\label{P17.6.53N+1}
Sea $K$ un campo global de funciones. Existe $S$, un conjunto 
de lugares suficientemente grande, tal que $J_{K,0}=J_{K,S,0}\*K$
donde 
\begin{align*}
J_{K,S,0}&=J_{K,S}\cap J_{K,S}\\
&=\{\vec\alpha\in J_K\mid
\deg\vec\alpha=0\text{\ y\ }v_{\pK}(\alpha_{\pK})=0\text{\ para
casi toda $\pK\in S$}\}.
\end{align*}
\end{proposicion}

\begin{proof}
An\'aloga a la del Teorema \ref{T17.6.52N}.
$\fin$
\end{proof}

Cuando $K$ es num\'erico y $S=S_{\infty}$ es el conjunto de los 
lugares arquimedianos, $K^{S_{\infty}}=E_K=\*\o_K$ las
unidades del anillo de enteros $\o_K$ de $K$.

Se tiene que $K^S$ es el grupo de elementos invertibles
de los $S$-enteros de $K$: $\o_{K,S}=\{x\in K\mid v_{\pK}(x)
\geq 0\text{\ para toda $\pK\notin S$}\}=\bigcap_{\pK\notin S}
\o_{\pK}$, el cual es un dominio Dedekind.

Sea $V^s$ un espacio real de dimensi\'on finita $s$:
$V^s=\bigoplus_{i=1}^s {\ma R} v_i$. Una {\em red\index{red}}
es un subgrupo abeliano libre de $V^s$ de rango $s$
y tal que una ${\ma Z}$-base de este subgrupo es una
${\ma R}$-base de $V^s$. 

\begin{lema}\label{L17.6.54N}
Sea $S$ un conjunto finito no vac\'io de lugares de un campo
global $K$ conteniendo a los lugares arquimedianos cuando $K$
es un campo num\'erico. Sea $V:=\{f:S\lra{\ma R}\}\cong {\ma R}^s$
que es un ${\ma R}$-espacio vectorial $s$ dimensional, en donde
$s=|S|$, Sea $\lambda:K^s\lra V$ dada por $\lambda(a)=f_a$ 
donde $f_a(v)=\log |a|_v$, esto es, $\lambda: a\longmapsto
\big(\log |a|_{v_1},\ldots, \log |a|_{v_s}\big)\in {\ma R}^s$, donde
$S=\{v_1,\ldots,v_s\}$.

Entonces $\ker\lambda$ es un subgrupo finito de $K^S$ y su imagen
es una red generando el ${\ma R}$-espacio vectorial $V_0$ que
consiste de las $f\in V$ tales que $\sum_{v\in S}f(v)=0$.
\end{lema}

\begin{proof}
Si $c$ y $C$ son constantes con $0<c<C$, entonces $\{\eta\in K^S\mid
c\leq |\eta|_v\leq C\text{\ para toda $v\in S$}\}$ es
 la intersecci\'on de un subconjunto compacto de
$J_K$ con el grupo discreta $\*K$, esto es, es compacto y discreto
por lo que es finito.

En particular, $\{x\in \*K\mid |x|_v=1\text{\ para toda $v\in{\ma R}$}\}$
es finito y forma un grupo multiplicativo por lo que es el grupo de
las ra\'ices de unidad de $K$.

Se sigue que $\ker \lambda$ es finito.

Para ver la imagen de $\lambda$, notemos que se definir $\lambda$
por una f\'ormula an\'aloga en los $S$-id\`eles $J_{K,S}$, esto es
$J_{K,S}\lra {\ma R}^s$, $\vec\alpha\longmapsto \big(\log |\alpha_{v_1}|_{v_1},
\ldots,\log |\alpha_{v_s}|_{v_s}\big)$ y la imagen de $J_{K,S,0}=J_{K,S}
\cap J_{K,0}$ genera el ${\ma R}$-espacio vectorial consistente de las
$f$ satisfaciendo $\sum_{v\in S}f(v)=0$ que es de dimensi\'on $s-1$, pues
para un id\`ele $\vec\alpha$ 
en $J_{K,S,0}$ se tiene $\|\vec\alpha\|=\sum_{v\in S}|\alpha_v|_v=1$.

El grupo $\lambda(K^S)$ es discreto puesto que hay un n\'umero finito de
elementos  $\eta\in K^S$ con $1/2<|\eta|_v<2$ para toda $v\in S$. Ahora,
$J_{K,S,0}/K^S$ es un subgrupo abierto, por tanto cerrado, en $J_{K,0}/
\*K=C_{K,0}$ el cual es compacto. Se sigue que $J_{K,S,0}/K^S$ es
compacto y por lo tanto el epimorfismo
$J_{K,S,0}/K^S\stackrel{\lambda}{\longtwoheadrightarrow} \lambda(J_{
K,S,0})/\lambda(K^S)$ implica que $\lambda(J_{K,S,0})/\lambda(K^S)$
es compacto. Se sigue que $\lambda(K^S)$ genera a $V_0$.
$\fin$
\end{proof}

Como corolario, obtenemos el teorema de las unidades de Dirichlet,
mismo que ya hab\'iamos obtenido en el caso de campos de funciones
en el Corolario \ref{CRamDed1.2.5(0)} de otra manera.

\begin{teorema}[Teorema de las Unidades de
Dirichlet\index{teorema de las unidades de
Dirichlet}\index{Dirichlet!teorema de las unidades de
$\sim$}]\label{T17.6.55N}
Sea $K$ un campo global y sea $S$ un conjunto finito no vac\'io
de lugares de $K$ conteniendo a los primos infinitos. Entonces,
como grupos, se tiene $K^S\cong W_K\times {\ma Z}^{s-1}$ donde
$W_K$ son las ra\'ices de unidad en $K$ y $s=|S|$. $\fin$
\end{teorema}

Nuestro objeto de estudio es el grupo de clases de id\`eles $C_K=
J_K/\*K$ donde $K$ es un campo global. Si $L/K$ es una extensi\'on
finita y separable, se define el encaje $J_K\stackrel{\mu}{\lra} J_L$,
definido por $\mu(\vec\alpha)_{\pL}:=\alpha_{\pK}$ donde $\pL|\pK$,
esto es, si $\mu(\vec\alpha)=\vec\beta\in J_L$, las componentes de
$\vec \beta$ est\'an dadas por los provenientes de $K$:
$\beta_{\pL}=\alpha_{\pL|_{K}}$. En particular, si $\beta_{\pL}=
\beta_{\pL'}$ si $\pL|_K=\pL'|_K=\pK$.

Si $L/K$ es una extensi\'on de Galois y $G=\Gal(L/K)=G_{L|K}$,
entonces $G$ act\'ua en $J_L$ de manera natural: si $\sigma\in
G$, $\sigma$ define un isomorfismo $\sigma:L_{\sigma^{-1}\pL}
\lra L_{\pL}$. Por tanto, si $\vec\beta\in J_L$, $\sigma\vec\beta
\in J_L$ est\'a definido por $(\sigma\vec\beta)_{\pL}:=\sigma
\beta_{\sigma^{-1}\pL}\in L_{\pL}, \beta_{\sigma^{-1}\pL}\in L_{
\sigma^{-1}\pL}$, esto es, $\beta_{\sigma^{-1}\pL}\in L_{\sigma^{-1}
\pL}$ es la $\sigma^{-1}\pL$-componente de $\vec\beta$.
Notemos qe $(\sigma\vec\beta)_{\pL}$ es no unidad en $L_{\pL}
\iff (\vec\beta)_{\sigma^{-1}\pL}$ es no unidad.
Por tanto, cuando pasamos a ideales o divisores bajo el mapeo
$\Lambda:J_L\lra D_L$, $\Lambda(\vec\beta)=\prod_{\pL\nmid
\infty}\pL^{v_{\pL}(\beta_{\pL})}$, el mapeo $\sigma:\vec\beta
\lra \sigma\beta$ es el mapeo
\begin{align*}
\Lambda(\sigma\vec\beta)&={\eu a}_{\sigma\vec\beta}=
\prod_{\pL\nmid\infty}\pL^{v_{\pL}((\sigma\vec\beta)_{\pL})}=
\prod_{\pL\nmid \infty}\pL^{v_{\pL}(\sigma\beta_{\sigma^{-1}\pL})}\\
&\igual_{\substack{\uparrow\\ \pL'=\sigma^{-1}\pL\\ \sigma\pL'=\pL}}
\prod_{\pL'\nmid \infty}(\sigma\pL')^{v_{\sigma\pL'}(\sigma\beta_{
\pL'})}=\prod_{\pL'\nmid \infty}(\sigma\pL')^{v_{\pL'}(\beta_{\pL'})}\\
&=\sigma{\eu a}_{\vec\beta}=\sigma\Lambda(\vec\beta).
\intertext{Esto es}
\Lambda(\sigma\vec\beta)&={\eu a}_{\sigma\vec\beta}=\sigma
{\eu a}_{\vec\beta}=\sigma\Lambda(\vec\beta).
\end{align*}

\begin{teorema}\label{T17.6.56N}
Sea $L/K$ una extensi\'on finita de Galois de campos globales
con $G=\Gal(L/K)$. Entonces $J_L^G=J_K$.
\end{teorema}

\begin{proof}
Al considerar $J_K\subseteq J_L$, todas las componentes de 
$\vec\alpha$ en $J_L$ son las mismas en todo $\{\pL\}_{\pL|\pK}$
y $G$ \'unicamente las permuta. Mas precisamente, sea $\sigma
\in G$, $\sigma:L_{\sigma^{-1}\pL}\lra L_{\pL}$ es un 
$K_{\pK}$-isomorfismo para $\pL|\pK$ y si $\vec\alpha\in J_K$
considerado en $J_L$, entonces
\[
(\sigma\vec\alpha)_{\pL}=\sigma\alpha_{\sigma^{-1}\pL}=\sigma
\alpha_{\pL}=\alpha_{\pL}\in K_{\pK}.
\]
Por tanto $\sigma\vec\alpha=\vec\alpha$.

Ahora sea $\vec\beta\in J_L$ con $\sigma\vec\beta=\vec\beta$
para toda $\sigma\in G$. Por tanto, $(\sigma\vec\beta)_{\pL}=
\sigma\beta_{\sigma^{-1}\pL}=\beta_{\pL}$ para todos los primos
$\pL\in{\ma P}_L$. Ahora, el grupo de descomposici\'on $D(\pL|\pK)$
de $\pL/\pK$ de $K$, es isomorfo al grupo de Galois de la
extensi\'on $L_{\pL}/K_{\pK}$. Para $\sigma\in D(\pL|\pK)$ se tiene
$\sigma^{-1}\pL=\pL$ y puesto que $\beta_{\pL}=\sigma\beta_{
\sigma^{-1}\pL}=\sigma \beta_{\pL}$, se sigue que $\beta_{\pL}
\in K_{\pK}$ para todo primo $\pL|\pK$.

Para $\sigma\in G$ arbitrario, $(\sigma\vec\beta)_{\pL}=
\beta_{\pL}=(\sigma\vec\beta)_{\sigma^{-1}\pL}=\beta_{
\sigma^{-1}\pL}\in K_{\pK}$, es decir, los dos lugares $\pL$ y
$\sigma^{-1}\pL$ sobre $\pK$ tienen la misma componente
$\beta_{\pL}=\beta_{\sigma^{-1}\pL}\in K_{\pK}$ por lo que
$\vec\beta\in J_K$. Por tanto $J_L^G\subseteq J_K$ y
$J_K=J_L^G$.
$\fin$
\end{proof}

En general el mapeo $I_{K,0}\lra I_{L,0}$ no es inyectivo, es decir,
un ideal o divisor que no es principal en $D_{K,0}$ puede ser
principal en $D_{L,0}$. Este fen\'omeno no sucede para los
id\`eles.

\begin{proposicion}\label{P17.6.57N} Sea $L/K$ un extensi\'on
finita y separable de campos globales. Entonces $\*L\cap J_K=
\*K$.

En particular, si $\vec\alpha\in J_K$ es un id\`ele de $K$ que es
principal en $J_L$, entonces $\vec\alpha$ es principal en $K$.
\end{proposicion}

\begin{proof}
Es inmediato que $\*K\subseteq \*L\cap J_K$. Sea $\tilde L$
una extensi\'on finita de Galois que contiene a $L$ y sea
$\tilde G=\Gal(\tilde L/K)$. Entonces $J_K, J_L\subseteq J_{
\tilde L}$. Si $\vec\alpha\in \*{\tilde L}\cap J_K$, entonces 
$\vec\alpha\in J_{\tilde L}^{\tilde G}$, es decir, $\sigma\vec
\alpha=\vec\alpha$ para toda $\sigma\in\tilde G$. Ahora, puesto
que $\vec\alpha\in\*{\tilde L}$, se tiene que $\vec\alpha\in(\*{\tilde
L})^{\tilde G}=\*K$, por lo tanto, $\*{\tilde L}\cap J_K=\*K$. Se sigue
que $\*L\cap J_K\subseteq \*{\tilde L}\cap J_K=\* K$.
$\fin$
\end{proof}

\begin{corolario}\label{C17.6.58N}
El encaje $\mu:J_K\lra J_L$ induce un monomorfismo 
$\tilde\mu:C_K=J_K/\*K\lra J_L/\*L=C_L$ donde $L/K$ es
cualquier extensi\'on finita y separable de campos globales.
\end{corolario}

\begin{proof}
Se tiene que $\tilde\mu(\vec\alpha\*K)=\vec\alpha\*L$ es un
homomorfismo. Si $\tilde\mu(\vec\alpha\*K)=\vec\alpha\*L=1$,
es decir, $\vec\alpha\in\*L$, por lo que $\vec\alpha\in\*L\cap J_K
=\*K$. Se sigue que $\vec\alpha\*K=1$ y $\tilde\mu$ es inyectiva.
$\fin$
\end{proof}

En adelante, si $K\subseteq L$, entonces consideramos $C_K
\subseteq C_L$ de manera natural.

\begin{teorema}\label{T17.6.59N}
Sea $L/K$ un extensi\'on finita de Galois de campos globales
con $G=\Gal(L/K)$. Entonces $C_L$ es de manera natural un
$G$-m\'odulo y $C_L^G=C_K$.
\end{teorema}

\begin{proof}
Si $\vec\beta\*L\in C_L$, $\vec\beta\in J_L$, se define $\sigma(
\vec\beta\*L):=\sigma\vec\beta\cdot\*L$. La definici\'on es 
independiente del representante $\vec\beta\in J_L$ de $\vec\beta
\*L\in C_L$ y $C_L$ es un $G$-m\'odulo.

Se tiene la sucesi\'on $G$-exacta
\begin{gather*}
1\lra \*L\lra J_L\lra C_L\lra 1,
\intertext{de donde se sigue la sucesi\'on exacta en cohomolog\'ia}
1\lra (\*L)^G=\*K\lra J_L^G=J_K\lra C_L^G\lra \co 1G{\*L}=\{1\}.
\end{gather*}
Por tanto $C_L^G\cong J_K/\*K=C_K$.
$\fin$
\end{proof}

\subsection{Cohomolog\'ia de $J_L$}\label{S17.6.5N}

Sea $L/K$ una extensi\'on finita de Galois de campos globales
con grupo de Galois $G=\Gal(L/K)=G_{L|K}$. Sean $S$ un
conjunto finito de primos de $K$ y $\bar S=\{\pL\in {\ma P}_L
\mid \pL|_K\in S\}$. Se denota $J_{L,\bar S}=J_{L,S}$ y se
hablar\'a de los $S$-id\`eles de $L$ en lugar de los 
$\bar S$-id\`eles. As\'i
\begin{gather*}
J_{L,S}=\prod_{\substack{\pL|\pK\\ \pK\in S}} \*L_{\pL}
\times \prod_{\substack{\pL|\pK\\ \pK\notin S}} U_{\pL}=
\prod_{\pK\in S}\prod_{\pL|\pK}\*L_{\pL}\times \prod_{\pK\notin S}
\prod_{\pL|\pK}U_{\pL}.
\intertext{Sean}
J_L^{\pK}:=\prod_{\pL|\pK}\*L_{\pL}\quad\text{y}\quad
U_L^{\pK}:=\prod_{\pL|\pK}U_{\pL}
\end{gather*}
subgrupos de $J_L$ en donde las componentes de $J_L^{\pK}$
y de $U_L^{\pK}$ en primos $\pL\nmid \pK$ son iguales a $1$.

Como $G$ permuta los lugares $\pL$ con $\pL|\pK$, se tiene que
$J_L^{\pK}$ y $U_L^{\pK}$ son $G$-m\'odulos. Se tiene
$J_{L,S}=\prod_{\pK\in S}J_L^{\pK}\times \prod_{\pK\notin S}
U_L^{\pK}$.

\subsubsection{Norma de id\'eles}

Sea $L/K$ una extensi\'on finita y separable de campos globales.

\begin{definicion}\label{D17.6.60N} La {\em norma\index{norma
de id\`eles}} $\N_{L/K}$ 
de $J_L$ a $J_K$ se define de la siguiente forma:
Sea $M$ la cerradra normal de $L/K$ y sean $G=\Gal(M/K)$
y $H=\Gal(M/L)$. Entonces $J_L=J_M^H$. Sea $G/H$ el conjunto
de las clases izquierdas de $H$ en $G$. Entonces
\[
\N_{L/K}\vec\alpha=\prod_{\sigma\in G/H}\sigma\vec\alpha,\quad
\vec\alpha\in J_L.
\]
Se tiene que $\N_{L/K}$ no depende de los representantes $\sigma
\in G/H$ pues $J_L$ esta fijo bajo $H$.
\end{definicion}

Para $\tau\in G$, se tiene
\[
\tau\N_{L/K} \vec\alpha=\prod_{\sigma\in G/H}\tau\sigma
\vec\alpha=\prod_{\sigma\in G/H}\sigma\vec\alpha=\N_{L/K}
\vec\alpha,
\]
por tanto $\N_{L/K}\vec\alpha\in J_K$. 

Se tiene que $\N_{L/K}\*L\subseteq \*K$, es decir, $\N_{L/K}$
manda id\`eles principales en id\`eles principales. En particular
$\N_{L/K}$ induce $\N_{L/K}:C_L\lra C_K$.

\begin{teorema}\label{T17.6.61N}
Si $L/K$ es una extensi\'on finita y separable de campos globales,
las componentes locales de $\N_{L/K}\vec\alpha\in J_L$ est\'an
dadas por:
\[
\big(\N_{L/K}\vec\alpha\big)_{\pK}=\prod_{\pL|\pK} \N_{L_{\pL}/
K_{\pK}} \alpha_{\pL}.
\]

En particular si $L/K$ es una extensi\'on finita y abeliana, definimos
para $\vec \alpha \in J_L$
\[
\N_{L/K}\vec\alpha=\vec\beta\in J_K\quad\text{donde}\quad
\beta_{\pK}=\prod_{\pL|\pK}\N_{L_{\pL}/K_{\pK}}\alpha_{\pL}.
\]

Si $L/K$ es una extensi\'on abeliana finita, entonces la norma
$\N_{L/K}:J_L\lra J_K$ es un mapeo abierto y por tanto, la norma
inducida en las clases de id\`eles $\N_{L/K}:C_L\lra C_K$ es un mapeo
abierto.
\end{teorema}

\begin{proof}
Sea $M/K$ la cerradura normal de $L/K$ y sea $G=\Gal(M/K)$ y
$H=\Gal(M/L)$. Para cada primo $\pK$ de $K$, sea ${\eu q}_{\pK}$
un primo fijo de $M$ sobre $\pK$ y sea ${\pL}_{\pK}={\eu 
q}_{\pK}|_L$. Sea $G_{{\eu q}_{\pK}}=\Gal(M_{{\eu q}_{\pK}}/K_{
\pK})\subseteq G$ y sea $\{\sigma_i\}$ un conjunto completo de
representantes de las clases dobles $H\tau G_{{\eu q}_{\pK}}$,
es decir, $G=\cupdot_i H\sigma_i G_{{\eu q}_{\pK}}$.

Entonces ${\pL}_i:={\eu q}_i|_L=(\sigma_i {\eu q}_{\pK})|_L$ es
un conjunto completo de primos de $L$ sobre $\pK$. Para cada
$i$, sea $\{\tau_{ij}\}$ un conjunto completo de representantes 
derechos de las clases de $G_{{\eu q}_{\pK}}\cap H^{\sigma_i}$
en $G_{{\eu q}_{\pK}}$, donde $H^{\sigma_i}:=\sigma_i^{-1}H
\sigma_i$. Entonces
\[
G_{{\eu q}_{\pK}}=\cupdot_j\big(G_{{\eu q}_{\pK}}\cap H^{\sigma_i}\big)\tau_{ij}
\quad\text{y}\quad G=\cupdot_{i,j}H\sigma_i\tau_{ij}=\cupdot_{
i,j}\tau^{-1}_{ij}\sigma^{-1}_i H,
\]
y $G_{{\eu q}_{\pK}}\cap H_i$ es el grupo de descomposici\'on
de ${\eu q}_{\pK}$ sobre $\sigma_i^{-1} L$.

Sea $\vec\alpha\in J_L$, considered como un id\`ele de $M$.
Entonces
\begin{align*}
\N_{L/K}(\vec\alpha)_{\pK}&=\N_{L/K}(\vec\alpha)_{{\eu q}_{\pK}}
=\prod_{i,j}(\tau^{-1}_{ij}\sigma^{-1}_{i}\vec\alpha)_{{\eu q}_{\pK}}
=\prod_{i,j}\tau_{ij}^{-1}(\sigma_i^{-1}\vec\alpha)_{{\eu q}_{\pK}}\\
&=\prod_i \N_{(\sigma_i^{-1}L)_{{\eu q}_{\pK}}/K_{\pK}}\big((
\sigma_i^{-1}\vec\alpha)_{{\eu q}_{\pK}}\big).
\end{align*}

El isomorfismo $\sigma_i:\sigma_i^{-1}L\lra L$ induce
a su vez un $K_{
\pK}$-isomorfismo $\sigma_i:(\sigma_i^{-1}L)_{{\eu q}_{\pK}}
\lra L_{{\eu q}_i}=L_{\pL_i}$ que manda $(\sigma_i^{-1}\vec\alpha
)_{{\eu q}_{\pK}}$ a $\alpha_{{\eu q}_i}=\alpha_{\pL_i}$.
Por tanto
\begin{gather*}
\big(\N_{L/K}(\vec\alpha)\big)_{\pK}=\prod_i \N_{L_{{\pL}_i}/K_{
\pK}}(\vec\alpha)_{\pL_i}=\prod_{\pL|\pK}\N_{L_{\pL}/K_{\pK}}
(\alpha_{\pL}).
\end{gather*}

Veamos que la norma es abierta.
Sean $V$ un conjunto abierto de $J_L$ y $\vec\alpha\in V$. 
Sea $W=\prod_{\pL} W_{\pL}$ una vecindad de $\vec 1$
con $\vec\alpha \cdot W\subseteq V$. Podemos suponer que
$W=\prod_{\pL|\infty}W_{\pL}\times\prod_{\pL\nmid\infty} W_{\pL}$
con $W_{\pL}\subseteq \*L_{\pL}$ un abierto para $\pL$ un primo infinito
y $W_{\pL}=U_{\pL}^{(m_{\pL})}$ si $\pL$ es un primo finito
y donde $m_{\pL}=0$ para casi toda $\pL$. Ahora bien, $U_{\pL}^{
(m_{\pL})}$ es compacto, abierto y cerrado (Proposiciones
\ref{CCUnidades} y \ref{CCCompactos}) 
y como $\N_{L_{\pL}/K{\pK}}$ es continua, se tiene
que $\N_{L_{\pL}/K{\pK}}(U_{\pL}^{(m)})$ es compacto y de \'indice
finito en $U_{\pK}^{(m)}$ por lo que 
$\N_{L_{\pL}/K{\pK}}(U_{\pL}^{(m)})$
es abierto y cerrado. Por tanto,
para cada $m$, existe un $n_m$ tal que $\N_{L_{\pL}/K{\pK}}
U_{\pL}^{(m)}\supseteq U_{\pK}^{(n_m)}$ por ser $\{U_{\pK}^{(n)}
\}_{n=0}^{\infty}$ un sistema fundamental de vecindades de $1\in
\*K_{\pK}$ (Corolario \ref{CCdisconexo}).

Si $\pK$ es infinito, se tiene $\N_{L_{\pL}/K{\pK}}\*L_{\pL}\supseteq
{\ma R}^+$. Sea $\vec\beta=\N_{L/K}\vec\alpha$ y sea $W'=\prod_{
\pK} W'_{\pK}$ tal que $W'_{\pK}\subseteq {\ma R}^+$ y
tal que $\N_{L_{\pL}/K{\pK}} W_{\pL}\supseteq W'_{\pK}$ si $\pK$ es
infinito y $W_{\pK}=U_{\pK}^{(n_{\pK})}$ con $\N_{L_{\pL}/K{\pK}}
U_{\pL}^{(m_{\pL})}\supseteq U_{\pK}^{(n_{\pK})}$ para $\pK$
finito. Entonces
\[
\vec\beta\cdot W'\subseteq \N_{L/K}\big(\vec\alpha\cdot W\big)
\subseteq \N_{L/K} V.
\]
Por tanto $\N_{L/K}$ es un mapeo abierto.
$\fin$
\end{proof}

\begin{proposicion}\label{P17.6.62N} Sea $L/K$ una extensi\'on
abeliana finita de campos globales. Sea $\pL$ un primo de $L$
sobre el primo $\pK$. Entonces
\[
\co qG{J_L^{\pK}}\cong \co q{G_{\pL}}{\*L_{\pL}}\quad\text{y}
\quad \co qG{U_L^{\pK}}\cong \co q{G_{\pL}}{U_{\pL}}\quad
\text{para toda\ } q\in{\ma Z},
\]
donde $G_{\pL}$ es el grupo de descomposici\'on de $\pL$
sobre $K$ y se tiene $G_{\pL}\cong \Gal(L_{\pL}/K_{\pK})$.

Si $\pK$ es un primo finito de $K$ no ramificado en $L$, entonces
$\co qG{U_L^{\pK}}=1$ para toda $q\in{\ma Z}$.
\end{proposicion}

\begin{proof}
Los isomorfismos se siguen del Lema de Shapiro, Teorema
\ref{CClaseT1.5.9''}: $\*L_{\pL}$ y $U_{\pL}$ son
$G_{\pL}$-m\'odulos y se tiene
\[
J_L^{\pK}=\prod_{\pL|\pK}\*L_{\pL}\cong \bigoplus_{\sigma\in
G/G_{\pL}}\sigma \*L_{\pL}\quad\text{y}\quad
U_L^{\pK}=\prod_{\pL|\pK}U_{\pL}\cong \bigoplus_{\sigma\in
G/G_{\pL}}\sigma U_{\pL}.
\]
Cuando $\pK$ es no ramificado en $L/K$, $\pK$ es no ramificado
en $L_{\pL}/K_{\pK}$, por lo tanto $\co q{G_{\pL}}{U_{\pL}}=
\{1\}$ para toda $q\in{\ma Z}$ (Teorema \ref{CCLT17.6.3}).
$\fin$
\end{proof}

\begin{observacion}\label{O17.6.63N}
Si $\pi:J_L^{\pK}\lra \*L_{\pL}$ es la proyecci\'on $\pi(\vec\alpha)=
\alpha_{\pL}$, el isomorfismo por $\co qG{J_L^{\pK}}\cong
\co q{G_{\pL}}{\*L_{\pL}}$ est\'a dado por
\[
\co qG{J_L^{\pK}}\stackrel{\res}{\lra}\co q{G_{\pL}}{J_L^{\pK}}
\stackrel{\bar\pi}{\lra} \co q{G_{\pL}}{\*L_{\pL}}.
\]
\end{observacion}

Si $S$ es un conjunto finito de lugares de $K$, se tiene 
\[
\co qG{J_{L,S}}\cong\prod_{\pK\in S}\co qG{J_L^{\pK}}\times
\prod_{\pK\notin S}\co qG{U_L^{\pK}}.
\]

Cuando $S$ contiene a todos los primos ramificados en $L$, entonces por
la Proposici\'on \ref{P17.6.62N} se obtiene 
$\co qG{J_{L,S}}=\prod_{\pK\in S}\co
q{G_{\pL}}{\*L_{\pL}}$ donde $\pL$ es cualquier primo de $L$ dividiendo
a $\pK$ primo de $K$ (pues para $\pK\notin S$, $\pK$ es no ramificado
y $\co qG{U_L^{\pK}}=\{1\}$).

Ahora $J_L=\bigcup\limits_S J_{L,S}=\lim\limits_{\substack{\lra\\ S}} J_{L,S}$ y
\[
\co qG{J_L}\cong\co qG{\lim_{\substack{\lra\\ S}}J_{L,S}}
\cong\lim_{\substack{\lra\\ S}}
\co qG{J_{L,S}}\cong \bigoplus_{\pK}\co q{G_{\pL}}{\*L_{\pL}},
\]
y donde $S$ var\'ia sobre todos los conjuntos finitos de $K$ que
contienen a todos los primos ramificados. Hemos probado

\begin{teorema}\label{T17.6.64N}
Sea $L/K$ una extensi\'on finita de Galois de campos globales. Para cada
primo de $\pK$ de $K$, sea $\pL$ un primo de $L$ sobre $\pK$. Sea $S$
cualquier conjunto finito de primos de $K$ que contiene a todos los primos
ramificados en $L$. Entonces
\begin{gather*}
\co qG{J_{L,S}}\cong \prod_{\pK\in S}\co q{G_{\pL}}{\*L_{\pL}},\\
\co qG{J_L}\cong \bigoplus_{\pK\in {\ma P}_K}\co q{G_{\pL}}{\*L_{\pL}}.
\tag*{$\fin$}
\end{gather*}
\end{teorema}

\begin{observacion}\label{O17.6.65N} Del resultado del Teorema 
\ref{T17.6.64N} se obtiene que la cohomolog\'ia de $J_L$ se puede
calcular por medio de los grupos de cohomolog\'ia local $\co q{G_{\pL}}
{\*L_{\pL}}$.
\end{observacion}

\begin{observacion}\label{O17.6.66N} El isomorfismo $\co qG{J_L}\cong
\bigoplus_{\pK}\co q{G_{\pL}}{\*L_{\pL}}$ est\'a dado por las proyecciones
$\co qG{J_L}\lra \co q{G_{\pL}}{\*L_{\pL}}$, es decir, la composici\'on
de los mapeos
\[
\co qG{J_L}\stackrel{\res}{\lra} \co q{G_{\pL}}{J_L}\stackrel{\bar\pi}{\lra}
\co q{G_{\pL}}{\*L_{\pL}},
\]
donde $\bar\pi$ est\'a inducido por la proyecci\'on natural $J_L
\stackrel{\pi}{\lra}\*L_{\pL}$, $\vec\alpha\longmapsto \alpha_{\pL}$.
\end{observacion}

Las proyecciones anteriores mapean cada $c\in \co qG{J_L}$ a su
$\pK$-componente $c_{\pK}\in \co q{G_{\pL}}{\*L_{\pL}}$. Por tanto
el Teorema \ref{T17.6.64N} dice que $c$ est\'a un\'ivocamente determinada
por sus componentes locales $c_{\pK}$ y como es suma directa,
$c_{\pK}=1$ para casi toda $\pK$. Para dimensiones $q>0$ el mapeo
$c\lra c_{\pK}$ puede ser descrito como sigue. Dado $c\in \co qG{
J_L}$, seleccionamos un cociclo $\beta(\sigma_1,\ldots,\sigma_q)$
representando a $c$. Esta es una funci\'on del grupo a $G\times
\ldots\times G$ con valores en $J_L$. Restringimos la funci\'on a
$G_{\pL}\times\cdots\times G_{\pL}$ y tomamos la $\pL$-componente
$\beta_{\pL}(\sigma_1,\ldots, \sigma_q)$ del id\`ele $\beta(\sigma_1,
\ldots,\sigma_q)$. La funci\'on resultante de $G_{\pL}\times
\cdots G_{\pL}$ a $\*L_{\pL}$
es un cociclo y su clase de cohomolog\'ia $c_{\pK}\in \co q{G_{\pL}}{
\*L_{\pL}}$ es la $\pK$-componente de $c$.

\begin{proposicion}\label{P17.6.67N}
Sea $K\subseteq L\subseteq M$ una torre de campos globales con
$M/K$ y $L/K$ extensiones finitas y de Galois. Sea ${\eu q}|\pL|\pK$
primos de $M$, $L$ y de $K$ respectivamente. Entonces
\las
\item $(\infla_M c)_{\pK}=\infla_{M_{\eu q}}(c_{\pK}),\quad
c\in \co q{G_{L|K}}{J_L}, \quad q\geq 1$,

\item\label{17.6.2NN} $(\res_L c)_{\pL}=\res_{L_{\pL}}(c_{\pK}), \quad 
c\in \co q{G_{M|K}}{J_M},\quad q\in{\ma Z}$,

\item\label{17.6.3NN} $(\cores_K c)_{\pK}=\sum_{\pL|\pK}\cores_{K_{\pK}}(c_{\pL}),
\quad c\in \co q{G_{M|L}}{J_M},\quad q\in{\ma Z}$,
\end{list}
donde en {\rm{(\ref{17.6.2NN})}} y en {\rm{(\ref{17.6.3NN})}} es suficiente suponer
que $M/K$ es de Galois.
\end{proposicion}

\begin{proof}
En la f\'ormula (\ref{17.6.3NN}), para cada $\pL|\pK$, podemos seleccionar
${\eu q}$ en $M$ con ${\eu q}|\pL$, pero entonces las corestricciones
$\cores_{K_{\pK}}(c_{\pL})$ pertenecen a distintos grupos de cohomolog\'ia
$\co q{G_{M_{\eu q}|K_{\pK}}}{\*M_{\eu q}}$. Sin embargo podemos 
identificar todos estos elementos como sigue. Si ${\eu q}$ y ${\eu q}'$
son dos lugares de $M$ sobre $\pK$, existe $\sigma\in G_{M|K}$ tal
que ${\eu q}^{\sigma}={\eu q}'$ y $\*M_{\eu q}\stackrel{\sigma}{\lra}
\*M_{{\eu q}^\sigma}$ induce un isomorfismo natural $\co q
{G_{M_{\eu q}|K_{\pK}}}{\*M_{\eu q}}\cong \co q
{G_{M_{\eu q}|K_{\pK}}}{\*{M_{{\eu q}^{\sigma}}}}$.

De esta forma seleccionamos, para cada $\pL|\pK$ un ${\eu q}|\pL$
y cada $\cores_{K_{\pK}}(c_{\pL})\in \co q
{G_{M_{\eu q}|K_{\pK}}}{\*M_{\eu q}}$ y la suma se realiza en $\co q
{G_{M_{\eu q}|K_{\pK}}}{\*M_{\eu q}}$.

Para la demostraci\'on de la proposici\'on, se tiene que el mapeo de
restricci\'on que se obtiene cuando pasamos a las componentes locales,
conmuta con los mapeos $\infla$, $\res$ y $\cores$. Esta afirmaci\'on
puede verificarse a nivel de cociclos para $\infla$ y $\res$ si $q\geq 1$,
y para $\cores$ si $q=-1,0$. El caso general se sigue por cambio de
dimensi\'on.
$\fin$
\end{proof}

\begin{teorema}[Teorema de la norma para id\`eles]\label{T17.6.68N}
Sea $L/K$ una extensi\'on finita de Galois de campos globales. Entonces
un id\`ele $\vec\alpha\in J_K$ es la norma de un id\`ele $\vec\beta\in
J_L$ si y s\'olo si cada componente $\alpha_{\pK}\in \*K_{\pK}$ es norma
de un elemento $\beta_{\pL}\in \*L_{\pL}$ con $\pL|\pK$, es decir, si y s\'olo
si es una norma local para toda completaci\'on.
\end{teorema}

\begin{proof}
Se tiene $\co 0G{J_L}=\frac{J_L^G}{\N_{L/K} J_L}=\frac{J_K}{\N_{L/K}J_L}$
y $\co 0{G_{\pL}}{\*L_{\pL}}=\frac{\*K_{\pK}}{\N_{L_{\pL}/K_{\pK}}\*L_{\pL}}$.
Por el Teorema \ref{T17.6.64N}, 
\[
\co 0G{J_L}\cong\frac{J_K}{\N_{L/K}J_L}\cong\bigoplus_{\pK}\frac{\*K_{\pK}}
{\N_{L_{\pL}/K_{\pK}}L_{\pL}}=\bigoplus_{\pK}\co 0{G_{\pL}}{\*L_{\pL}},
\]
donde  para cada $\pK$ se selecciona un $\pL$ con $\pL|\pK$.

Si $\vec\alpha\in J_K$, el isomorfismo anterior manda la clase $\vec\alpha
\N_{L/K}J_L=\idel \alpha$ a las componentes $\tilde\alpha_{\pK}=\alpha_{
\pK}\N_{L_{\pL}/K_{\pK}}\*L_{\pL}$. Se sigue $\idel \alpha=1\iff \tilde\alpha_{
\pK}=1$ para toda $\pK$, es decir, $\vec\alpha\in \N_{L/K}J_L\iff \alpha_{
\pK}\in \N_{L_{\pL}/K_{\pK}}\*L_{\pL}$ para toda $\pK\in{\ma P}_K$.
$\fin$
\end{proof}

El teorema de la norma para id\`eles es un an\'alogo al teorema de la norma
de Hasse que establece que si $L/K$ es una extensi\'on c\'iclica finita de
campos globales y si $x\in\*K$, entonces $x$ es norma de $\*L$ si y s\'olo
si $x$ es norma local para toda completaci\'on (ver Corolario
\ref{C17.6.99N}).

Se tiene que si $x\in\*K\subseteq J_K$, entonces $x\in\N_{L/K}J_L\iff 
x\in\N_{L_{\pL}/K_{\pK}}\*L_{\pL}$ para toda $\pK$, sin embargo no
sabemos que $\vec\beta$ sea principal, es decir, si $\vec\beta =y\in
\*L$. En otras palabras, si $x\in K$ es norma para todas las completaciones,
$x$ es norma de id\`ele $\vec \beta$ pero no necesariamente de
un id\`ele principal.

\begin{corolario}\label{C17.6.69N}
Sea $L/K$ una extensi\'on finita y de Galois de campos globales.
Entonces se tiene $\co 1G{J_L}=\co 3G{J_L}=\{1\}$. Adem\'as, si
$S$ es un conjunto finito de lugares que contiene a los primos
ramificados y a los primos infinitos, entonces $\co 1G{J_{L,S}}=
\{1\}$.
\end{corolario}

\begin{proof}
Se sigue de que
$\co 1{G_{L_{\pL}|K_{\pK}}}{\*L_{\pL}}=
\co 3{G_{L_{\pL}|K_{\pK}}}{\*L_{\pL}}=\{1\}$
 (Corolario \ref{CCLC17.6.23}).
$\fin$
\end{proof}

\begin{corolario}\label{C17.6.70N}
Se tiene que 
\[
\co 2G{J_L}\cong \bigoplus_{\pK\in{\ma P}_K}
\Big(\big(\frac{1}{[L_{\pL}:K_{\pK}]}\ {\ma Z}\big)/{\ma Z}\Big).
\]
\end{corolario}

\begin{proof}
Del Teorema de Tate-Nakayama (Teorema \ref{CClaseT1.5.15}),
se tiene 
\begin{gather*}
\co 2{G_{\pL}}{\*L_{\pL}}\cong \co 0{G_{\pL}}{\ma Z}\cong
\frac{\ma Z}{|G_{\pL}|{\ma Z}}\cong \frac{\Big(\frac{1}{[L_{\pL}
:K_{\pK}]}\ {\ma Z}\Big)}{{\ma Z}},
\intertext{y}
\co 2G{J_L}\cong\bigoplus_{\pK\in{\ma P}_K} \co 2{G_{\pL}}
{\*L_{\pL}}. \tag*{$\fin$}
\end{gather*}
\end{proof}

\subsection{Grupo de Brauer\index{grupo de Brauer}\index{Brauer!grupo
de $\sim$} para campos globales}\label{S17.6.6N}

Puesto que para cualquier extensi\'on finita de Galois $L/K$ de campos
globales se tiene $\co 1{G_{L|K}}{J_L}=\{1\}$, se sigue que $\co 2{G_{
L|K}}{J_L}\subseteq \co 2{G_{M|K}}{J_M}$ para $K\subseteq L
\subseteq M$ una torre de extensiones finitas de Galois, la cual se
obtiene como consecuencia de la sucesi\'on exacta inflaci\'on-restricci\'on.

Sea $\Omega=\sep K$. Puesto que $J_{\Omega}=\lim\limits_{\substack{
\lra\\ L}}J_L=\bigcup\limits_{L/K}J_L$ se obtiene que
\begin{gather*}
\co 2{G_{\Omega|K}}{J_{\Omega}}=\bigcup_L\co 2{G_{L|K}}{J_L}
\intertext{y que para $K\subseteq L\subseteq M$ extensiones finitas de Galois,}
\co 2{G_{L|K}}{J_L}\subseteq \co 2{G_{M|K}}{J_M}\subseteq \co
2{G_{\Omega|K}}{J_{\Omega}}.
\end{gather*}

Por la teor\'ia local, se tiene que el grupo de Brauer $\Br E$ de un
campo local $E$ satisface que $\Br E=\bigcup_F\co 2{G_{F|E}}{\*F}$
donde $F/E$ recorre las extensiones no ramificadas
(Teorema \ref{CCLT17.6.16}). En el caso global
veremos que $\Br K=\bigcup_L\co 2{G_{L|K}}{\*L}$ donde $L$ recorre
las extensiones c\'iclicas ciclot\'omicas de $K$ si $K$ es un campo
num\'erico y las extensiones de constantes de $K$ si $K$ es un
campo de funciones. En otras palabras, $L\in\{\text{subcampos de
$K(\zeta_n)\mid K(\zeta_n)/K$ es c\'iclica}\}$ o 
$L\in\{K{\ma F}_{q^n}\mid n\in{\ma N}\}$.

\begin{lema}\label{L17.6.71N}
Sean $K$ un campo global, $S$ un conjunto finito de primos de $K$ y
$m\in {\ma N}$. Entonces existe una extensi\'on c\'iclica ciclot\'omica 
o una extensi\'on de constantes $L/K$ con la propiedad:
\las
\item[$\bullet$] $m|[L_{\pL}:K_{\pK}]$ para toda $\pK\in S$, $\pK$ finito,
\item[$\bullet$] $[L_{\pL}:K_{\pK}]=2$ para toda $\pK\in S$, $\pK$ real.
\end{list}
\end{lema}

\begin{proof}
Supongamos que hemos probado el lema para los casos especiales
$K_0={\ma Q}$ en el caso
num\'erico y para $K_0=K$ para campos de funciones. Entonces
aplicamos el lema para los casos especiales, encontrando
$M/K_0$ una extensi\'on totalmente imaginaria (caso num\'erico)
c\'iclica ciclot\'omica o una extensi\'on
de constantes en el caso de campos de funciones, tal que para
cada primo $\pK$ de $K_0$ en $S$, $\pK$ finito, 
hay un primo $\pL$ de $M$ para el
cual el grado $[M_{\pL}:(K_0)_{\pK}]$ es divisible por $m\cdot
[K:K_0]$.

Entonces $L:=KM$ satisface las propiedades requeridas para $L/K$.

Ahora, para $K=K_0$, el lema se cumple pues podemos realizar lo
anterior debido a la descomposici\'on de primos en campos 
ciclot\'omicos en el caso num\'erico, o debido a que en el caso
de campos de funciones, en las extensiones de constantes, todo
primo es eventualmente inerte.
$\fin$
\end{proof}

\begin{teorema}\label{T17.6.72N}
Sea $K$ un campo global. Entonces
\[
\Br K=\bigcup_L \co 2{G_{L|K}}{\*L}\quad\text{y}\quad
\co 2{G_{\Omega}|K}{J_{\Omega}}=\bigcup_L
\co 2{G_{L|K}}{J_L},
\]
donde $L/K$ recorre las extensiones c\'iclicas ciclot\'omicas
finitas en el caso num\'erico y las extensiones de constantes
finitas en caso de campos de funciones.
\end{teorema}

\begin{proof}
Sea $c\in \co 2{G_{\Omega|K}}{J_{\Omega}}$, digamos $c\in
\co 2{G_{L|K}}{J_M}$. Sea $m$ el orden de $c$ y sea $S$ un 
conjunto finito de primos de $K$ para los cuales las componentes
$c_{\pK}$ de $c$ son distintas de $1$. Sea $L/K$ una extensi\'on
como en el Lema \ref{L17.6.71N} con $m|[L_{\pL}:K_{\pK}]$ si
$\pK\in S$, $\pK$ finito y $[L_{\pL}:K_{\pK}]=2$ para $\pK\in S$
real.

Sea $N:=ML$. Entonces
\[
\co 2{G_{M|K}}{J_M}, \co 2{G_{L|K}}{J_{L}}\subseteq 
\co 2{G_{N|K}}{J_{N}}.
\]
Veremos que $c\in \co 2{G_{L|K}}{J_{L}}$. Puesto que la sucesi\'on
\[
1\lra \co 2{G_{L|K}}{J_{L}}\stackrel{\infla}{\lra}
\co 2{G_{N|K}}{J_N}\stackrel{\res_L}{\lra}
\co 2{G_{N|K}}{J_{N}}
\]
es exacta, es suficiente probar que $\res_L c=1$.

Por teor\'ia local, junto el Teorema \ref{T17.6.64N} y la
Proposici\'on \ref{P17.6.67N}, se tiene que $\res_L c=1\iff
(\res_L c)_{\pL}=\res_{L_{\pL}} c_{\pK}=1$ para toda
$\pL\in {\ma P}_L\iff \inv_{N_{\eu q}|L_{\pL}}(\res_{L_{\pL}}
c_{\pK})=[L_{\pL}:K_{\pK}]\cdot \inv_{N_{\eu q}|K_{\pK}} c_{\pK}
=\inv_{N_{\eu q}|K_{\pK}} c_{\pK}^{[L_{\pL}:K_{\pK}]}=0$
para toda $\pK\in{\ma P}_K\iff c_{\pK}^{[L_{\pL}:K_{\pK}]}=1$
para toda $\pK\in S$.

Esta \'ultima condici\'on se cumple puesto que $c_{\pK}^m=1$
y $m|[L_{\pL}:K_{\pK}]$ para toda $\pK\in S$ finito y $[L_{\pL}:
K_{\pK}]=2$ para todo $\pK\in S$ real.

Se sigue que $\co 2{G_{\Omega|K}}{J_{\Omega}}=\bigcup_L
\co 2{G_{L|K}}{J_{L}}$.

Para probar que $\Br K =\bigcup_L \co 2{G_{L|K}}{\*L}$ se hace
exactamente lo mismo, sustituyendo, para un campo $E$, $J_E$
por $\*E$.
$\fin$
\end{proof}

\subsection{Primera desigualdad fundamental\index{primera
desigualdad fundamental}}\label{S17.6.7N}

Esta primera desigualdad establece que si $L/K$ es una
extensi\'on c\'iclica de grado $n$ de campos globales, entonces
\[
[C_K:\N_{L/K}C_L]\geq n.
\]

Es importante mencionar, que la primera desigualdad \'unicamente
ser\'a v\'alida para extensiones abelianas y no para otro tipo de
extensiones, incluyendo las de Galois no abelianas. Esto se puede
ver en el Teorema \ref{T17.6.192N}. La segunda desigualdad es
v\'alida para extensiones finitas de Galois arbitrarias.

\begin{teorema}\label{T17.6.73N}
Sea $L/K$ una extensi\'on c\'iclica de grado $n$ de campos
globales. Entonces el cociente de Herbrand de $C_L$, 
satisface
\[
h(C_L)=\frac{|\co 0G{C_L}|}{|\co 1G{C_L}|}=n,
\]
donde $G=\Gal(L/K)$ es el grupo de Galois de la extensi\'on
$L/K$.

En particular, 
\[
|\co 0G{C_L}|=|\co 2G{C_L}|=[C_K:\N_{L/K} C_L]
= n\cdot |\co 1G{C_L}|\geq n.
\]
\end{teorema}

\begin{proof}
Sea $L/K$ una extensi\'on c\'iclica de campos de funciones. Sea $U_L:=
\prod_{\pL}U_{\pL}$ y consideremos el grupo $J_{L,0}$. Entonces
$U_L$, $\*L\subseteq J_{L,0}$, por lo cual $U_L\*L\subseteq J_{L,0}$.
Se tiene
\[
h(C_L)=h(J_L/\*L)=h(J_L/J_{L,0})h(J_{L,0}/U_L\*L)h(U_L\*L/\*L).
\]

Consideremos la sucesi\'on exacta $1\lra J_{L,0}\lra J_L\stackrel{
\deg}{\lra} {\ma Z}\lra 0$, por lo que $J_L/J_{L,0}\cong {\ma Z}$ y 
$h(J_L/J_{L,0})=h({\ma Z})=n$.

Ahora $\Lambda:J_{L,0}\lra D_{L,0}\stackrel{\pi}{\twoheadrightarrow}
I_{L,0}$ y $\ker \Lambda = U_L$. Por tanto $D_{L,0}\cong J_{L,0}/U_L$.
Adem\'as $\ker \pi=P_L=\Lambda(\*L)$, de donde 
obtenemos $J_{L,0}/U_L\*L
\cong I_{L,0}$ el cual es finito por lo que $h(J_{L,0}/U_L\*L)=h(I_{L,0})
=1$.

Ahora, $U_L\*L/\*L\cong U_L/U_L\cap \*L$. Se tiene que $U_L\cap
\*L=\*\F$, donde $\F$ es el campo de constantes de $L$ y por ende
$h(U_L\cap \*L)=1$.

Por otro lado $U=\prod_{\pK\in{\ma P}_K}\prod_{\pL|\pK}U_{\pL}$
y $\sigma\big(\prod_{\pL|\pK}U_{\pL}\big)=\prod_{\pL|\pK}U_{\pL}$
para toda $\sigma\in G$.

Del Lema de Shapiro, obtenemos $\co rG{\prod_{\pL|\pK}U_{\pL}}
\cong \co r{G_{\pL}}{U_{\pL}}$. Por teor\'ia de campos de clase locales
tenemos que $h(G_{\pL},U_{\pL})=1$ (ver la demostraci\'on
del Teorema \ref{axiomacamposlocales}). 
De hecho $|\co 0{G_{\pL}}{U_{\pL}}|=|\co 1{G_{\pL}}{U_{\pL}}|=
e_{\pL}$ y $e_{\pL}=1$ para casi toda $\pL$ (Corolario \ref{C17.2.15'} y
Teoremas \ref{CCLT17.6.3} y \ref{T17.5.60}).

Por tanto $|\co 0G{U_L}|=|\co 1G{U_L}|=\prod_{\pL}e_{\pL}$ y por
tanto $h(U_L)=1$ y $h(U_L/U_L\cap \*L)=h(U_L)h(U_L\cap\*L)^{-1}
=1$. De esta forma obtenemos $h(C_L)=n\cdot 1\cdot 1=n$.

\s

Ahora consideremos una extensi\'on de campos num\'ericos. Sea
$S$ un conjunto finito de primos de $K$ tal que 
\las
\item[1.-] $S$ contiene a los primos infinitos y a los primos
ramificados.
\item[2.-] $J_L=J_{L,S}\*L$.
\item[3.-] $J_K=J_{K,S}\*K$.
\end{list}
Entonces $C_L=J_{L,S}\*L/\*L\cong J_{L,S}/L^S$, donde $L^S=
J_{L,S}\cap \*L$ es el grupo de las $S$-unidades, $L^S=\{x\in\*L\mid
v_{\pL}(x)=0\text{\ para toda $\pL\notin S$}\}$ y $h(C_L)=h(J_{L,S})
h(L^S)^{-1}$. Podemos considerar que $h(J_{L,S})$ es la parte
local de $h(C_L)$ y que $h(L^S)$ es la parte global de $h(C_L)$.

Se va a probar que $h(L^S)=\frac 1n\prod_{\pK\in S} n_{\pK}$ donde
$n_{\pK}=[L_{\pL}:K_{\pK}]$ y donde $\pL$ es un primo fijo de $L$
que divide a $\pK$, $n_{\pK}=|D_{L/K}(\pL|\pK)|=|G_{\pL}|=
|\Gal(L_{\pL}/K_{\pK})|$ donde $D_{L/K}(
\pL|\pK)$ es el grupo de descomposici\'on.

Sea $s=|\bar S|=|\{\pL\in{\ma P}_L\mid \pL|_K\in S\}|$ y para cada
$\pL\in \bar S$ sea $e_{\pL}$ un s\'imbolo formal. Definimos
$E^s:=\bigoplus_{\pL\in \bar S}{\ma R} e_{\pL}$ como el espacio
vectorial sobre los reales de dimensi\'on $s$ con base $\{e_{\pL}\mid
\pL\in\bar S\}$. El grupo $G=\Gal(L/K)$ act\'ua en $E^s$ por
$\sigma (e_{\pL})=e_{\sigma\pL}$ para $\sigma\in G$.

Sea $M$ una red en $E^s$, esto es, $M$ es un subgrupo abeliano
libre de $E^s$ de rango $s$ y tal que una ${\ma Z}$-base de este
subgrupo es una ${\ma R}$-base de $E^s$. Supongamos que $M$
es invariante bajo $G$, es decir, $\sigma M\subseteq M$ para toda
$\sigma\in G$. Veremos que existe una subred $M_0$ de
$M$ de \'indice finito que tiene una ${\ma Z}$-base $\{y_{\pL}\mid
\pL\in\bar S\}$ tal que $\sigma y_{\pL}=y_{\sigma\pL}$.

Para esto tomemos la norma infinita $\|\ \|_{\infty}$ de $E^s$ con
respecto a las coordenadas relativas a la base $\{e_{\pL}\mid \pL
\in \bar S\}$. Esto es, $\|\sum_{\pL\in \bar S} d_{\pL}e_{\pL}\|_{\infty}
=\max_{\pL\in \bar S}\{|d_{\pL}|\}$.
Puesto que $M$ es una red, existe $N$ tal que para
cada $x\in E^s$, existe $y\in M$ con $\|x-y\|_{\infty}<N$. Para cada
$\pK\in S$, sea ${\pL}_0\in \bar S$ un primo fijo de $L$ sobre
$\pK$. Sean $t\in{\ma R}$, $t>s\cdot n\cdot N$ y
$z_{\pL_0}\in M$ tal que $\|te_{\pL}-z_{\pL_0}\|_{\infty}<N$.

Para $\pL|\pK$, sea $y_{\pL}:=\sum_{\sigma\pL_0=\pL}\sigma
z_{\pL_0}$. Sea $M_0=\bigoplus_{\pL\in \bar S} {\ma Z} y_{\pL}$.
Para $\tau\in G$, se tiene
\[
\tau y_{\pL}=\sum_{\sigma \pL_0=\pL}\tau \sigma z_{\pL_0}
\igual_{\substack{\uparrow\\ \rho=\tau\sigma\\ \rho\pL_0=
\tau\sigma\pL_0=\tau \pL}} \sum_{\rho\pL_0=\tau \pL}\rho 
z_{\pL_0}=y_{\tau\pL}.
\]

Debemos ver que $\{y_{\pL}\mid\pL\in\bar S\}$ es 
linealmente independiente sobre
${\ma R}$. Sea $\sum_{\pL\in\bar S}c_{\pL}y_{\pL}=0$ con $c_{\pL}\in
{\ma R}$. Si alg\'un $c_{\pL}\neq 0$, definimos $r=\max\{|c_{\pL}|\mid \pL
\in\bar S\}>0$. Por tanto tenemos $\sum_{\pL\in\bar S}\frac{c_{\pL}}{
\pm r}y_{\pL}=0$. Podemos suponer que $|c_{\pL}|\leq 1$ para toda
$\pL\in\bar S$ y $c_{\pL'}=1$ para alg\'un $\pL'$.

Sea $z_{\pL_0}=te_{\pL_0}+b_{\pL_0}$ con $b_{\pL_0}$ un vector tal
que $\|b_{\pL_0}\|_{\infty}<N$. Entonces
\begin{gather*}
y_{\pL}=\sum_{\sigma\pL_0=\pL}\sigma z_{\pL_0}=t\sum_{\sigma\pL_0
=\pL} e_{\sigma \pL_0}+b'_{\pL}
\intertext{con $\|b'_{\pL}\|_{\infty}\leq n\cdot N$. Por tanto}
y_{\pL}=tm_{\pL}e_{\pL}+b'_{\pL}
\quad\text{con}\quad m_{\pL}=|\{\sigma\in G\mid \sigma\pL_0=\pL\}|
=|D_{L/K}(\pL_0|\pK)|=n_{\pK}.
\end{gather*}
Se sigue que $
0=\sum_{\pL}c_{\pL}y_{\pL}=t\sum_{\pL} c_{\pL}m_{\pL}e_{\pL}+b'$
con $\|b'\|_{\infty}\leq s \cdot n\cdot N$. De esta forma, se tiene
\begin{gather*}
\sum_{\pL}c_{\pL}m_{\pL}e_{\pL}=-\frac {b'}{t},\\
\Big\|\sum_{\pL}c_{\pL}m_{\pL}e_{\pL}\Big\|_{\infty}=\max_{\pL}|c_{\pL}
m_{\pL}|\geq c_{\pL'}m_{\pL'}=m_{\pL'}\geq 1\quad\text{y}\\
\Big\|-\frac{b'}{t}\Big\|_{\infty}=\frac 1t\|b'\|_{\infty}<1, \quad
\text{de donde obtenemos}\quad
m_{\pL'}\leq \frac 1t\|b'\|_{\infty}<1,
\end{gather*}
lo cual es una contradicci\'on, probando que $M_0$ satisface lo
requerido.

Sea $\mu: L^S\lra E^s$ dada por $\mu(x)=\sum_{\pL\in\bar S} \log |x|_{
\pL}\cdot e_{\pL}$. La imagen es una red $s-1$ dimensional con n\'ucleo
finito (Teorema \ref{T17.6.19N}).

Para $\sigma\in G$, se tiene
\begin{align*}
\mu(\sigma x)&=\sum_{\pL}\log|\sigma x|_{\pL} e_{\pL}=
\sum_{\pL}\log |x|_{\sigma^{-1}\pL}\underbrace{\sigma e_{\sigma^{-1}\pL}}_{
\substack{\uigual\\ e_{\pL}}}\\
&=\sigma\big(\sum_{\pL}\log |x|_{\sigma^{-1}\pL}
e_{\sigma^{-1}\pL}\big)=\sigma\mu(x),
\end{align*} 
por lo que $\mu$ es un $G$-homomorfismo.

Sea $e_0:=\sum_{\pL\in \bar S}e_{\pL}$. Entonces $e_0$ y $\mu(L^S)$
generan una red $s$-dimensional $M$ en $E^s$. Puesto que ${\ma Z}
e_0$ es $G$-isomorfo a ${\ma Z}$, de la sucesi\'on exacta
\begin{gather*}
0\lra {\ma Z}e_0\lra M\lra M/{\ma Z}e_0\lra 0
\intertext{y de que $M/{\ma Z}e_0\cong \mu(L^s)$ deducimos que}
h(L^s)=h(\mu(L^s))=h({\ma Z})^{-1}\cdot h(M)=\frac 1n h(M).
\end{gather*}

Sea $M_0$ la subred de $M$ de \'indice finito en $M$ como antes.
Entonces 
\begin{gather*}
M_0\cong\bigoplus_{\pK\in S}\bigoplus_{\pL|\pK}
{\ma Z}y_{\pL}\cong \bigoplus_{\pK \in S} M_0^{\pK}
\quad\text{y}\quad
M_0^{\pK}
=\bigoplus_{\pL|\pK}{\ma Z}y_{\pL}=\bigoplus_{\sigma\in G/G_{\pL}}
\sigma \big({\ma Z} y_{\pL_0}\big).
\intertext{Por el Lema de Shapiro, se tiene}
h(M_0^{\pK})=h(G_{\pL},{\ma Z}y_{\pL_0})=|G_{\pL}|=n_{\pK}
\intertext{y como $M/M_0$ es finito, $h(M)=h(M_0)$, lo que implica que}
h(L^S)=\frac 1n h(M)=\frac 1n \prod_{\pK\in S}h(G_{\pL}, {\ma Z}y_{
\pL_0})=\frac 1n \prod_{\pK\in S}n_{\pK}. 
\end{gather*}

Ahora
\begin{gather*}
\co qG{J_{L,S}}=\prod_{\pK\in S}\co q{G_{\pL}}{\*L_{\pL}}, \quad\text{por tanto}\\
h(G, J_{L,S})=\prod_{\pK\in S}h(\underbracket[0pt]{G_{\pL}}_{
\substack{\uigual\\ G_{\pK}}},\*L_{\pL})=\prod_{\pK\in S}|G_{\pK}|
=\prod_{\pK\in S} n_{\pK}.
\intertext{Se sigue que}
h(C_L)=h(J_{L,S})h(L^s)^{-1}=\Big(\prod_{\pK\in S}n_{\pK}\Big)
\cdot \frac n{\Big(\prod_{\pK\in S}n_{\pK}\Big)}=n.
\tag*{$\fin$}
\end{gather*}
\end{proof}

\begin{corolario}[Caso especial del teorema de densidad
de Chebotarev]\label{C17.6.74N}
Sea $L/K$ una extensi\'on c\'iclica de campos globales de
grado una potencia de un primo $n=p^m$. Entonces existe una
infinidad de primos de $K$ que no se descomponen totalmente
en $L$. Como consecuencia tenemos que hay una
infinidad de primos totalmente inertes en $L/K$.
\end{corolario}

\begin{proof}
Supongamos primero que $L/K$ es c\'iclica de grado $p$. Sea $S$ el
conjunto de primos de $K$ que no se descomponen en $L/K$. Si $S$
fuese finito se probar\'a que  $\N_{L/K} C_L=C_K$ lo que contradice
que $[C_K:\N_{L/K} C_L]\geq p$.

Sea $\vec\alpha\in J_K$. Por el teorema de aproximaci\'on, tenemos
que existe $a\in\*K$ tal que $(\vec\alpha a^{-1})_{\pK}
=\alpha_{\pK} a^{-1}$ est\'a contenido en el
subgrupo abierto $\N_{L_{\pL}/K_{\pK}}\*L_{\pL}$ de $\*K_{\pK}$
(es abierto por el teorema de existencia para campos locales)
para toda $\pK\in S$. 

Para $\pK\notin S$, $\pK$ se descompone totalmente
en $L/K$ y por tanto $L_{\pL}=K_{\pK}$ y por
tanto $(\vec\alpha a^{-1})_{\pK}
\in \N_{L_{\pL}/K_{\pK}}\*L_{\pL}=\*K_{\pL}$.
Entonces $(\vec\alpha a^{-1})_{\pK}=\alpha_{\pK}a^{-1}$
es una norma local para toda $\pK\in{\ma P}_K$.
Por tanto $\vec\alpha a^{-1}$ es una norma de un id\`ele $\vec\beta$
de $L$ (Teorema \ref{T17.6.68N}). Por tanto
$\vec\alpha=\big(\N_{L/K}\vec\beta\big)\cdot a\in \big(\N_{L/K}
J_L\big)\cdot \*K$. Se sigue que $\idel \alpha\in \N_{L/K}C_L$ y
$C_K=\N_{L/K}C_L$.

Ahora sea $L/K$ c\'iclica de grado $p^m$ y sea $M$ la \'unica subextensi\'on
de grado $p$: $[M:K]=p$. Hay una infinidad de primos inertes en $M/K$,
lo que implica que hay una infinidad de primos no totalmente descompuestos
en $L$.

Finalmente, si $\pK$ es inerte en $M/K$, entonces $\pK$ es totalmente
inerte en $L/K$ pues si $G_{\pL}$ es el grupo de descomposici\'on de $\pL|
\pK$ y como la extensi\'on $L/K$
consiste de una \'unica torre de extensiones de grado $p$ (por ser la
extensi\'on c\'icica), se tiene que $\Gal(L/M)\subsetneqq G_{\pL}$ puesto 
que $\pK$ es inerte en $M/K$. Por tanto $G_{\pL}=\Gal(L/K)$ y $\pK$ es
totalmente inerte en $L/K$.
$\fin$
\end{proof}

\begin{corolario}\label{C17.6.75N}
Sea $L/K$ es cualquier extensi\'on finita y separable de campos globales
con $L\neq K$. Entonces hay una infinidad de primos $\pK$ de $K$
que no se descomponen totalmente en $L/K$. En otras palabras, si
$L_{\pL}=K_{\pK}$ para casi toda $\pK$, entonces $L=K$.
\end{corolario}

\begin{proof}
Si $\tilde L$ es la cerradura de Galois de $L/K$, entonces un primo
de $K$ se descompone totalmente en $L$ si y solamente si
se descompone totalmente en $\tilde L$. Por tanto, podemos suponer
que $L/K$ es una extensi\'on de Galois.

Sea $H<G$ un subgrupo c\'iclico de orden una potencia de un
n\'umero primo, $|H|=p^m$ y sea $M:=L^H$. Entonces hay una
infinidad de primos de $M$ que no se descomponen totalmente
en $L$. Se sigue el resultado.
$\fin$
\end{proof}

\begin{observacion}\label{CClaseOF.7} 
El Corolario \ref{C17.6.75N} es un caso particular
del Teorema de densidad de
Chebotarev el cual establece que en cualquier
extensi\'on finita de Galois de campos globales, todas las posibles
descomposiciones posibles de primos, tienen densidad positiva.
\end{observacion}

\begin{teorema}\label{T17.6.76N}
Sea $L=K_1\cdots K_r/K$ una extensi\'on abeliana $p$-elemental
de campos globales, donde $p$ es un n\'umero primo. Se tiene
$\Gal(L/K)\cong C_p^r=C_p\times\cdots\times C_p$. Entonces
existe una infinidad de primos $\pK$ de $K$ que son inertes en
$K_1/K$ y son totalmente descompuestos en $K_i/K$ para toda
$i\geq 2$.
\end{teorema}

\begin{proof}
Se tiene que $L/K_2\cdots K_r$ es c\'iclica de grado $p$. Sea
${\eu q}$ un primo de $K_2\cdots K_r$ que permanece primo
en $L$ y tal que ${\eu q}|_K=\pK$ es no ramificado. De esta
forma se tiene $L_{\eu q}/(K_2\cdots K_r)_{\eu q}$ es c\'iclica
de grado $p$. Por otro lado, $L_{\eu q}/K_{\pK}$ es un grupo
c\'iclico puesto que $\pK$ es no ramificado y $\Gal(L_{\eu q}/
K_{\pK})$ es isomorfo al grupo de Galois de los campos residuales
los cuales son campos finitos.

Puesto que $\Gal(L/K)\cong C_p^r$, 
se tiene que $[L_{\eu q}: K_{\pK}]\leq p$
y, por otro lado,
como $[L_{\eu q}:K_{\pK}]\geq [L_{\eu q}:(K_2\cdots K_r)_{\eu
q}]=p$, se sigue que $(K_2\cdots K_r)_{\eu q}=K_{\pK}$ por lo
que $\pK$ se descompone totalmente en $K_2\cdots K_r/K$.
Adem\'as $\pK$ es inerte en $K_1/K$ pues $f_{\pK}(L|K)=p$.
$\fin$
\end{proof}

Otra consecuencia de la primera desigualdad es el teorema de
F.K. Schmidt.

\begin{teorema}[F. K. Schmidt]\label{T17.6.77N}
Sea $K$ un campo global de funciones. Entonces existe un
divisor de $K$ de grado $1$.
\end{teorema}

\begin{proof}
Sea $\deg:D_K\lra {\ma Z}$ y sea $\deg(D_K)=\rho{\ma Z}$ para
alg\'un $\rho\in{\ma N}$. Sea $k=\F$ el campo de constantes
de $K$. Sea $k_1:={\ma F}_{q^{\rho}}$ la extensi\'on de $k$
de grado $\rho$ y sea $K_1=Kk_1$ la extensi\'on de constantes.
Para cada primo $\pK$ de $K$, el campo residual $K(\pK)
\cong {\ma F}_{q^{\deg \pK}}$ 
contiene a $k_1$ pues $\rho|\deg\pK$. Por tanto, si $\pL|\pK$,
con $\pL$ primo en $K_1$, entonces $(K_1)_{\pL}=K_{\pK}$
puesto que $(K_1)_{\pL}=K_{\pK}\cdot k_1=K_{\pK}$, esto es,
todo primo de $K$ se descompone totalmente en $K_1$
de donde se sigue que ${\ma F}_{q^{\rho}}=\F$ y que $\rho=1$.
$\fin$
\end{proof}

\subsection{Segunda desigualdad fundamental\index{segunda
desigualdad fundamental}}\label{S17.6.8N}

La segunda desigualdad ser\'a v\'alida para cualquier extensi\'on 
finita de Galois de campos globales. Por otro lado, la primera
desigualdad \'unicamente  es v\'alida para extensiones abelianas.
Esto dir\'a m\'as adelante que si $L/K$ es una extensi\'on finita
y separable de campos globales, entonces
\[
[C_K:\N_{L/K}C_L]=[L:K] \iff \text{$L/K$ es abeliana.}
\]

La segunda desigualdad establece que si $L/K$ es una extensi\'on
finita de Galois de campos globales, entonces 
$[C_K:\N_{L/K} C_L]\big|[L:K]$.

Reduciremos este problema para probar que es suficiente
probar la desigualdad para una extensi\'on $L/K$ c\'iclica de
grado primo $l$.

\begin{lema}\label{L17.6.78N} Si $E/K$ es cualquier extensi\'on
finita y separable, entonces $[C_K:\N_{E/K}C_E]$ es finita y divide
a una potencia de $m=[E:K]$.
\end{lema}

\begin{proof}
Sea $L/K$ la cerradura normal de $E/K$. Entonces se tiene $\N_{L/K}C_L
\subseteq \N_{E/K}C_E$. Por tanto 
\[
[C_K:\N_{L/K}C_L]=[C_K:\N_{
E/K}C_E][\N_{E/K}C_E:\N_{L/K}C_L].
\]
Si probamos que $[C_K:\N_{L/K}C_L]$ es finita, entonces $[C_K:
\N_{E/K}C_E]$ es finita, por lo tanto, para la finitud, basta considerar
el caso cuando $E/K$ es de Galois.

Sea pues
$E/K$ Galois. Sea $S$ un conjunto finito de primos incluyendo a los
arquimedianos, a los primos ramificados y suficientemente grande tal
que $J_K=J_{K,S}\*K$ y $J_E=J_{E,S}\*E$ (por ser $E/K$ normal). Entonces 
\begin{gather*}
\*K\N_{E/K}
J_E=\*K\N_{E/K}\*E\cdot \N_{E/K} J_{E,S}=\*K \N_{E/K}J_{E,S}.
\intertext{Por tanto} 
\begin{align*}
[J_K:\*K\N_{E/K}J_E]&=[\*KJ_{K,S}:\*K\N_{E/K}J_{E,S}]\leq
[J_{K,S}:\N_{E/K}J_{E,S}]\\
&= \Big[\prod_{v\in S}\*K_v\times \prod_{v\notin S}U_v:
\prod_{v\in S}\prod_{\omega|v}\N_{E_{\omega}/K_v}\*E_{\omega}
\times \prod_{v\notin S}\prod_{\omega|v} U_{\omega}\Big]\\
&\hspace{-18pt}
\igual_{\substack{\uparrow\\ v\notin S \Rightarrow v\\ \text{no ramificado}}}
\Big[\prod_{v\in S}\*K_v\times \prod_{v\in S} U_v: 
\prod_{v\in S}\prod_{\omega|v}\N_{E_{\omega}/K_v}\*E_{\omega}
\times \prod_{v\notin S} U_v\Big]\\
&\hspace{-15pt}
\des_{\substack{\uparrow\\ E/K\text{\ Galois}}} 
\prod_{v\in S}[\*K_v: \N_{E_{\omega}/K_v}\*E_{\omega}]
\igual_{\substack{\uparrow\\ 
\text{Teorema de}\\ \text{reciprocidad local}}}
\prod_{v\in S}[E_{\omega}:K_v]\\
&=\prod_{v\in S} n_v<\infty.
\end{align*}
\end{gather*}
Esto prueba la finitud.

Ahora, el \'indice de la norma divide a una potencia de grado $m=[E:K]$,
la extensi\'on siendo normal o no, pues si
$\vec\alpha\in C_K$, $\vec\alpha^m=\N_{E/K}\vec\alpha\in
\N_{E/K} C_E$, lo que implica $C_K/\N_{E/K}C_E$ es subgrupo de
$\big({\ma Z}/m{\ma Z}\big)^r$ para alg\'un $r\in{\ma N}$.
$\fin$
\end{proof}

\begin{lema}\label{L17.6.79N}
Sea $K\subseteq F\subseteq E$ dos extensiones finitas. Entonces
\lasa
\item $[C_K:\N_{F/K}C_F]\big|[C_K:\N_{E/K}C_E]$.
\item $[C_K:\N_{E/K}C_E]\big|[C_K:\N_{F/K}C_F]\cdot [C_F:\N_{E/F}C_E]$.
\end{list}

Por tanto, si la desigualdad se presenta en los pasos de la torre, se tiene
la desigualdad en la torre misma.
\end{lema}

\begin{proof}
Se tiene
\[
[C_K:\N_{E/K}C_E]=[C_K:\N_{F/K}C_F][\N_{F/K}C_F:\N_{F/K}(\N_{E/F}
C_E)]
\]
lo cual prueba (a).

Ahora, el mapeo $C_F\lra \N_{F/K}C_F$, $\vec\alpha\longmapsto \N_{
F/K}\vec\alpha$, es un homomorfismo. 
Se tiene el epimorfismo inducido
\[
\frac{C_F}{\N_{E/F}C_E}\xrightarrow{\N_{F/K}} \frac{\N_{F/K}C_F}{\N_{F/K}\big(
\N_{E/F} C_E\big)}=\frac{\N_{F/K}C_F}{\N_{E/K}C_E},
\]
por lo que $[\N_{F/K}C_F:
\N_{F/K}(\N_{E/F}C_E)]$ divide a $[C_F:\N_{E/F}C_E]$ de donde se
sigue (b).
$\fin$
\end{proof}

\begin{corolario}\label{C17.6.80N}
Si la segunda desigualdad se cumple para todas las extensiones
c\'iclicas de grado primo, entonces se cumple para todas las extensiones
de Galois.
\end{corolario}

\begin{proof}
Sea $L/K$ una extensi\'on de Galois y sea $l$ un n\'umero primo. Sea
$E$ el campo fijo de $L$ de un $l$-subgrupo de Sylow del grupo $G=
\Gal(L/K)$. Se tiene que $L/E$ es una torre de extensiones c\'iclicas de 
grado $l$ y la desigualdad se cumple en $L/E$ por el
Lema \ref{L17.6.79N}, esto es
$[C_E:\N_{L/E}C_L]\big|[L:E]$.

Por otro lado $[C_K:\N_{E/K}C_E]$ divide a una potencia $[E:K]$ y por
tanto es primo relativo a $l$. Ahora, por el Lema \ref{L17.6.79N}, se tiene
que 
\[
[C_K:\N_{L/K} C_L]\big|[C_K:\N_{E/K}C_E][C_E:\N_{L/E}C_L],
\]
por lo
que para cada n\'umero primo $l$, la $l$-contribuci\'on de $[C_K:\N_{L/K}
C_L]$ divide a $[L:E]$ y por tanto a $[L:K]$. La desigualdad se sigue
para $L/K$.
$\fin$
\end{proof}

\begin{corolario}\label{C17.6.81N}
Si $l$ un n\'umero primo con
$l\neq p$ donde $p$ es la caracter\'istica de $K$, es suficiente
probar la segunda desigualdad para extensiones c\'iclicas de $K$
de orden $l$ suponiendo que una $l$-ra\'iz primitiva $\zeta_l$
de la unidad est\'a en $K$.
\end{corolario}

\begin{proof}
Sea $L/K$ una extensi\'on c\'iclica de grado $l$. Entonces
\begin{gather*}
[C_K:\N_{L/K}C_L]\big|[C_K:\N_{L(\zeta_l)/K}C_{L(\zeta_l)}]\quad
\text{y}\\
[C_K:\N_{L(\zeta_l)/K}C_{L(\zeta_l)}]\big|[C_K:\N_{K(\zeta_l)/K}C_{
K(\zeta_l)}]\cdot [C_{K(\zeta_l)}:\N_{L(\zeta_l)/K(\zeta_l)}C_{
L(\zeta_l)}].
\end{gather*}

Adem\'as $[C_K:\N_{K(\zeta_l)/K}C_{
K(\zeta_l)}]$ es primo relativo a $l$ pues $[K(\zeta_l):K]|l-1$ y
$[C_K:\N_{L/K}C_L]$ es una potencia de $l$ y divide a 
$ [C_{K(\zeta_l)}:\N_{L(\zeta_l)/K(\zeta_l)}C_{
L(\zeta_l)}]$, una extensi\'on c\'iclica de grado primo de un campo
que contiene a $\zeta_l$.
$\fin$
\end{proof}

\subsubsection{Segunda desigualdad para extensiones de Kummer}

Sea $K$ un campo global que contiene a las $n$-ra\'ices de la unidad,
$\zeta_n\in K$ y donde $n$ es una potencia de un n\'umero primo diferente
a la caracter\'istica $p$ de $K$. Sea $L/K$ una extensi\'on de Galois
con grupo de Galois $\Gal(L/K)\cong ({\ma Z}/n{\ma Z})^r$.

Sea $S$ un conjunto finito no vac\'io de lugares de $K$ conteniendo
a los primos infinitos, a los primos que dividen a $n$, a los primos
ramificados en $L$ y suficientemente grande tal que $J_K=
J_{K,S}\*K$. Sea $s=|S|$.

\begin{proposicion}\label{P17.6.82N}
Se tiene que $s\geq r$ y que existe un conjunto $T$ de $s-r$
primos de $K$, disjunto de $S$, tal que $L=K(\sqrt[n]{\Delta})$
donde $\Delta$ es el n\'ucleo del homomorfismo 
$K^S\lra \prod_{\pK\in T}\*K_{\pK}/(\*K_{\pK})^n$.
\end{proposicion}

\begin{proof}
Se tiene que $L=K(\sqrt[n]{D})$ donde $D=(\*L)^n\cap \*K$ por
la Teor\'ia de Kummer, Teorema \ref{CClaseT1.6.2}. Si $x\in D$,
entonces $K_{\pK}(\sqrt[n]{x})/K_{\pK}$ es no ramificado para 
todo $\pK\notin S$. Por tanto $x=u_{\pK}y_{\pK}^n$ con $u_{\pK}
\in U_{\pK}$ y $y_{\pK}\in\*K_{\pK}$.

Definimos $y_{\pK}=1$ para $\pK\in S$ de tal forma que
se define un id\`ele
$\vec y$ que se puede escribir $\vec y=\vec \alpha\cdot z$ con
$\vec\alpha\in J_{K,S}$, $z\in\*K$. Ahora bien
\[
xz^{-n}=u_{\pK}y_{\pK}^nz^{-n}=u_{\pK}\alpha_{\pK}^nz^nz^{-n}=
u_{\pK}\alpha_{\pK}^n\in U_{\pK}
\]
para toda $\pK\notin S$, de tal forma que $xz^{-n}\in J_{K,S}\cap
\*K=K^S$, por lo que $xz^{-n}\in\Delta:=(\*L)^n\cap K^S$. Por tanto
$D=\Delta (\*K)^n$ (pues claramente, $\Delta
\subseteq D$), de donde $L=K(\sqrt[n]{\Delta})$.

El campo $M:=K(\sqrt[n]{K^S})$ contiene a $L$ puesto que $\Delta
\subseteq K^S$. Nuevamente, por Teor\'ia de Kummer, tenemos
$\Gal(M/K)\cong\Hom\big(K^S/(K^S)^n,\langle \zeta_n\rangle\big)$.
Puesto que $K^S$ tiene rango finito $s-1$ (Teorema \ref{T17.6.19N})
y contiene a $\langle\zeta_n\rangle$, se tiene que $K^S/(K^S)^n
\cong ({\ma Z}/n{\ma Z})^s$. Por tanto $\Gal(M/K)\cong ({\ma Z}/
n{\ma Z})^n$ y como $\frac{\Gal(M/K)}{\Gal(M/L)}\cong \Gal(L/K)
\cong ({\ma Z}/n{\ma Z})^r$, se tiene que $r\leq s$ y $\Gal(M/L)
\cong ({\ma Z}/n{\ma Z})^{s-r}$.

Sean $\{\sigma_1,\ldots,\sigma_{s-r}\}$ una $({\ma Z}/n{\ma Z})$-base
de $\Gal(M/L)$, $M_i$ el campo fijo de $\sigma_i$, $M_i=
M^{\langle\sigma_i\rangle}$ y $L=\bigcap_{i=1}^{s-r} M_i$. Para
cada $1\leq i\leq s-r$, sea $\pL_i$ un primo de $M_i$ totalmente
inerte en $M/M_i$ y tal que $\pK_i:=\pL_i|_K\notin S$. Adem\'as,
seleccionamos estos primos $\pK_i$ distintos. Esto
se puede hacer por el Corolario \ref{C17.6.74N}.

Puesto que $\pK_i\notin S$, se tiene que $\pK_i$ es no ramificado
en $M/K$ pues $M=K(\sqrt[n]{K^S})$ y por tanto $\pK_i$ es no
ramficado en cada $K(\sqrt[n]{x})$, $x\in K^S$.
Sea $Z_i$ el campo de 
descomposici\'on del \'unico primo ${\eu q}_i$ de $M$
sobre $\pL_i$, $1\leq i\leq s-r$ en $M/K$. Como $\pL_i$ es totalmente
inerte en $M/M_i$, se tiene $Z_i\supseteq M_i$. Por otro
lado, el grupo de descomposici\'on $D_i=D_{M/K}({\eu q}_i|{\pK_i})$
es c\'iclico pues $\pK_i$ es no ramificado. Puesto que $D_i
\subseteq \Gal(M/K)=({\ma Z}/n{\ma Z}\big)^s$, $D_i$ tiene 
exponente $n$ de donde se sigue que $D_i$ es de orden
un divisor de $n$. Por tanto $Z_i\subseteq M_i$.
En resumen, $M_i$ es el campo de descomposici\'on de 
${\eu q}_i$ en $M/K$.

Sea $T:=\{\pK_1,\ldots,\pK_{s-r}\}$.
Puesto que $L=\bigcap_{i=1}^{s-r} M_i$, se tiene que $L/K$ es la
m\'axima subextensi\'on de $M/K$ en el cual todos los primos $\pK_1,
\ldots, \pK_{s-r}$ son totalmente descompuestos. Por tanto, si
$x\in K^S$, se tiene 
\begin{gather*}
x\in\Delta\iff K(\sqrt[n]{x})\subseteq L\iff
K_{\pK_i}(\sqrt[n]{x})=K_{\pK_i}, \quad i=1,\ldots, s-r\\
\iff x\in (\*K_{
\pK_i})^n, \quad i=1,\ldots, s-r,
\end{gather*}
 probando que $\Delta$ es el n\'ucleo
del mapeo natural 
\begin{gather*}
K^S\lra \prod_{i=1}^{s-r}
\*K_{\pK_i}/(\*K_{\pK_i})^n, \quad x\longmapsto 
\prod_{i=1}^{s-r} x\bmod (\*K_{\pK_i})^n.
\tag*{$\fin$}
\end{gather*}
\end{proof}

\begin{lema}\label{L17.6.83N}
Sean $S$ y $T$ como en la Proposici\'on {\rm{\ref{P17.6.82N}}}
y sean
\begin{gather*}
J_K(S,T):=\prod_{\pK\in S}(\*K_{\pK})^n\times \prod_{\pK\in T}
\*K_{\pK}\times \prod_{\pK\notin (S\cup T)} U_{\pK},
\intertext{y}
C_K(S,T):=J_K(S,T)\*K/\*K.
\end{gather*}

Entonces $J_K(S,T)\cap \*K=(K^{S\cup T})^n$.
\end{lema}

\begin{proof}
Es inmediato de $(K^{(S\cup T)})^n\subseteq J_K(S,T)\cap
\*K$. Sea $y\in J_K(S,T)\cap \*K$ y sea $N=K(\sqrt[n]{y})$.
Para probar que $N=K$ basta probar que $C_K=\N_{N/K}
C_N$ pues por la primera desigualdad, se tiene $[C_K:
\N_{N/K}C_N]\geq [N:K]$.

Sea $\idel \alpha\in C_K =J_{K,S}\*K/\*K$ y sea $\vec\alpha$
un representante de $\idel \alpha$. El mapeo $K^S\lra \prod_{
\pK\in T} U_{\pK}/(U_{\pK})^n$ es sobre pues $\Delta$ es su
n\'ucleo y se tiene $|K^S/\Delta|=\frac{|K^S/(K^S)^n|}{|\Delta/(K^S)^n|}=
\frac{n^s}{|\Gal(L/K)|}=n^{s-r}$, $|U_{\pK}/(U_{\pK})^n|=
\frac{n}{n_\pK}=n$, donde $n_{\pK}=q^{-v_{\pK}(n)}=1$,
(Proposici\'on \ref{CCLTP17.6.14}) y por
tanto $\big|\prod_{\pK\in T}U_{\pK}/(U_{\pK})^n\big|
=n^{s-r}$.

Por tanto existe $x\in K^S$ con $\alpha_{\pK}=xu_{\pK}^n$ con
$u_{\pK}\in U_{\pK}$ para $\pK\in T$. El id\`ele $\vec\beta=
\vec\alpha x^{-1}$ pertenece a la misma clase $\vec\alpha$ y
veremos que $\vec\beta\in \N_{N/K}J_N$. Para esto, por el
teorema de la norma de id\`eles, Teorema \ref{T17.6.68N},
necesitamos que cada componente $\beta_{\pK}$ sea una
norma de $N_{\pL}/K_{\pK}$. Para $\pK\in S$ esto se cumple
puesto que $y\in(\*K_{\pK})^n$ por lo que $N_{\pL}=K_{\pK}$.
Para $\pK\in T$ esto se cumple pues $\beta_{\pK}=u_{\pK}^n$
es una $n$-potencia. 

Para $\pK\in S\cup T$, $\beta_{\pK}$
es una unidad y $N_{\pL}/K_{\pK}$ es no ramificada, por lo que
$\N_{N_{\pL}/K_{\pK}} U_{\pL}=U_{\pK}$ (Teoremas \ref{CCLT17.6.3} y
\ref{T17.6.20N}). Por tanto $C_K=\N_{N/K}C_N$ de donde
se sigue que $N=K$ y por tanto $y\in(\*K)^n\cap J_K(S,T)
\subseteq \big(K^{(SUT)}\big)^n$. De esta forma obtenemos
que $J_K(S,T)\cap \*K=\big(K^{(S\cup T)}\big)^n$.
$\fin$
\end{proof}

\begin{teorema}\label{T17.6.85N}
Con las notaciones anteriores, se tiene que
$\N_{L/K} C_L\supseteq C_K(S,T)$ y $[C_K:C_K(S,T)]=[L:K]$.

En particular, $[C_K:\N_{L/K}C_L]\big|[C_K:C_K(S,T)]=[L:K]$.

Cuando $L/K$ es c\'iclica, esto es, $r=1$,
entonces $\N_{L/K}C_L=C_K(S,T)$ y $[C_K:\N_{L/K}C_L]=[L:K]$.
\end{teorema}

\begin{proof}
Se tiene la sucesi\'on exacta
\[
1\lra \frac{J_{K,S\cup T}\cap \*K}{J_K(S,T)\cap \*K}\lra \frac{
J_{K,S\cup T}}{J_{K(S,T)}}\twoheadrightarrow \frac{J_{K,S\cup T}
\*K}{J_K(S,T)\*K}\lra 1.
\]

Puesto que $J_{K,S\cup T}\*K=J_K$, el orden del \'ultimo grupo es
\begin{gather*}
[J_{K,S\cup T}\*K: J_K(S,T)\*K]=\big[J_K/\*K:J_K(S,T)\*K/\*K\big]
=[C_K:C_K(S,T)].
\intertext{El orden del primer grupo es}
[J_{K,S\cup T}\cap \*K: J_K(S,T)\cap \*K]=[K^{S\cup T}:(K^{S\cup 
T})^n]=n^{s+s-r}=n^{2s-r},
\intertext{pues $\langle \zeta_n\rangle \subseteq K^{S\cup T}$.
El orden del grupo intermedio es}
\begin{align*}
[J_{K,S\cup T}:J_K(S,T)]&=\prod_{\pK\in S}[\*K_{\pK}:(\*K_{\pK})^n]
\igual_{\substack{\uparrow\\ {\text{Proposici\'on \ref{CCLTP17.6.14}}}}}
\prod_{\pK\in S}\frac{n^2}{|n|_{\pK}}\\
&=n^{2s}\prod_{\pK\in S} |n|_{\pK}^{-1}\igual_{\substack{\uparrow\\
\pK\in S\Rightarrow \pK\nmid n\Rightarrow\\ |n|_{\pK}=q^{-v_{
K}(n)}=q^0=1}} n^{2s}.
\end{align*}
\end{gather*}
Se sigue que $[C_K:C_K(S,T)]=\frac{n^{2s}}{n^{2s-r}}=n^r=[L:K]$.

Ahora veremos que $C_K(S,T)\subseteq \N_{L/K}C_L$. Sea $\vec
\alpha\in J_K(S,T)$. Por el Teorema \ref{T17.6.68N}, $\vec \alpha
\in\N_{L/K}J_L\iff \alpha_{\pK}\in \N_{L_{\pL}/K_{\pK}}\*L_{\pL}$
para toda $\pK\in{\ma P}_K$. Si $\pK\in S$, puesto que $\alpha_{
\pK}\in (\*K_{\pK})^n$, se tiene que es norma $K_{\pK}(\sqrt[n]{
\*K_{\pK}})$ por Teor\'ia de Kummer para campos locales
(Proposici\'on \ref{CClaseP3.2.21'-1}), y por
tanto de $K_{\pK}(\sqrt[n]{\Delta})$, $\Delta=\ker\big(K^S\lra
\prod_{\pK\in T}\*K_{\pK}/(\*K_{\pK})^n\big)$ pues $K_{\pK}(\sqrt[n]{
\Delta})\subseteq K_{\pK}(\sqrt[n]{\*K_{\pK}})$.

Para $\pK\in T$, por el Lema \ref{L17.6.83N}, $\Delta\subseteq
(\*K_{\pK})^n$ de tal forma que $\pK$ se descompone totalmente
y $L_{\pL}=K_{\pK}$.

Para $\pK\notin S\cup T$, $\alpha_{\pK}$ es una unidad y $L_{
\pL}=K_{\pK}(\sqrt[n]{\Delta})$ es no ramificada sobre $K_{\pK}$ y
$\N_{L_{\pL}/K_{\pK}} U_{\pL}=U_{\pK}$.

Por tanto $J_K(S,T)\subseteq \N_{L/K} J_L$ de donde $C_K(S,T)
\subseteq \N_{L/K} C_L$. Por tanto $[C_K:\N_{L/K}C_L]\leq [C_K:
C_K(S,T)]=n^r$.

Cuando $L/K$ es c\'iclica, es decir, cuando $r=1$, entonces, como
consecuencia de la primera desigualdad
\[
n=[L:K]\leq [C_K:\N_{L/K} C_L]\leq [C_K:C_K(S,T)]=[L:K]=n,
\]
de donde se sigue que $\N_{L/K}C_L=C_K(S,T)$ y $[C_K:\N_{
L/K}C_L]=[L:K]$.
$\fin$
\end{proof}

Esto prueba la segunda desigualdad para campos num\'ericos.

\subsubsection{Segunda desigualdad para campos de funciones}

Ahora, sea $K$ un campo global de funciones de caracter\'istica $p>0$.

\begin{definicion}\label{D17.6.86N}
Sea $E$ cualquier campo. Una {\em derivaci\'on\index{derivaci\'on}}
en $E$ es un mapeo $D:E\lra E$, tal que
\las
\item[$\bullet$] $D(x+y)=D(x)+D(y)$ para cualesquiera $x,y\in E$,
\item[$\bullet$] $D(xy)=xD(y)+yD(x)$ para cualesquiera $x,y\in E$.
\end{list}
\end{definicion}

Con $x=y=1$ obtenemos $D(1)=0$ y si $y=x^{-1}$ con $x\neq 0$, 
se obtiene $D(x^{-1})=-x^{-2}D(x)$.

Sea $P_n=\{x\in E\mid D^n(x)=0\}$. $P_n$
es un subgrupo aditivo de $E$.
Se define, $D^0=\Id_E$, $P_0=\{0\}$.

\begin{notacion}\label{N17.6.87N}
El operador aditivo $E\lra E$, $x\longmapsto xy$ se denota ya
sea simplemente por $y$ o por $\mu_y$: $\mu_y(x)=xy$.
De esta forma $Dy$ significa el operador producto $y$ seguido
de $D$ y $D(y)$ significa el efecto que $D$ tiene en el elemento
$y$.
\end{notacion}

\begin{definicion}\label{D17.6.88N}
Un elemento $x\in E$ es una {\em derivada logar\'itmica\index{derivada
logar\'itmica}} en $E$ si existe $y\in E$ tal que $x=\frac{D(y)}{y}$ y
$x$ es {\em la derivada logar\'itmica} de $y$.
\end{definicion}

\begin{lema}\label{D17.6.89N}
Un elemento $x\in E$ es la derivada logar\'itmica de un elemento
$y\in P_n\setminus P_{n-1}$, $n>0$, si y solamente si la $n$-potencia
del operador $D+x$ aplicada a $1$ es igual $0$ y la $(n-1)$-potencia
aplicada a $1$, es diferente a $0$:
\begin{gather*}
(D+x)^n(1)=0\quad\text{y}\quad (D+x)^{n-1}(1)\neq 0,
\intertext{o, equivalentemente,}
(D+\mu_x)^n(1)=0\quad\text{y}\quad (D+\mu_x)^{n-1}(1)\neq 0.
\end{gather*}
\end{lema}

\begin{proof}
Sea $x=\frac{D(y)}{y}$ para alguna $y\in\*E$. Se tiene
\begin{align*}
(D+\mu_x)(z)&=D(z)+xz=D(z)+\frac{D(y)}{y}z\\
&=y^{-1}D(zy)=
y^{-1}D\mu_y(z)=y^{-1}Dy(z).
\end{align*}
Esto es, $D+\mu_x=y^{-1}D\mu_y$ y por tanto, para $n\geq 0$, tenemos
\[
(D+x)^n=(D+\mu_x)^n=y^{-1}D^n\mu_y.
\]

Sea $y\in P_n\setminus P_{n-1}$. Entonces $(D+\mu_x)^n(1)=y^{-1}D^n
\mu_y(1)=y^{-1}D^n(y)=0$ y de manera an\'aloga, tenemos $(D+\mu_x)^{
n-1}(1)\neq 0$.

Rec\'iprocamente, supongamos que $(D+\mu_x)^n(1)=0$ y que $(D+
\mu_x)^{n-1}(1)\neq 0$. Definimos $y^{-1}:=(D+\mu_x)^{n-1}(1)$. Entonces
$(D+\mu_y)(y^{-1})=0=-y^{-2}D(y)+xy^{-1}$, por lo que $x=\frac{D(y)}{y}$.
Ahora bien, $D+\mu_x=y^{-1}Dy=y^{-1}D\mu_y$ y $(D+\mu_x)^{n-1}(1)=
y^{-1}D^{n-1}(y)\neq 0$ y similarmente $y^{-1}D^n(y)=0$. Por tanto
$y\in P_n\setminus P_{n-1}$.
$\fin$
\end{proof}

Consideremos un campo $E$ con la propiedad de
que para cada $x\in E$, exista
$n_x\in {\ma N}$ con $D^{n_x}(x)=0$, esto es, $E=\bigcup_{n=0}^{\infty}
P_n$. Sea $F\subseteq E$ un subcampo tal que
$D(F)\subseteq F$. Entonces, del Lema \ref{D17.6.89N} se sigue que
para $x\in F$, se 
tiene que $x$ es una derivada logar\'itmica en $F\iff x$ es un derivada
logar\'itmica en $E$. De hecho, si $x\in F$ es la derivada logar\'itmica
de un elemento $y\in P_n\setminus P_{n-1}$, entonces $x$ es la 
derivada logar\'itmica de $z$ donde $z^{-1}=(D+\mu_x)^{n-1}(1)$, el
cual es un elemento de $F$ pues $D(F)\subseteq F$.

La situaci\'on anterior se presenta en el siguiente caso. Sea $K$ un
campo de caracter\'istica $p>0$ y $E=K((t))$ 
el campo de las series formales
de Laurent con la derivada usual con respecto a 
$t$: $D\big(\sum_{n=m}^{\infty}a_n t^n\big)=\sum_{n=m}^{\infty}n a_n
t^{n-1}=\sum_{n=m-1}^{\infty}(n+1) a_{n+1}t^n$. Entonces se tiene
que $D^p=0$.

\begin{corolario}\label{C17.6.90N}
Sea $E=K((t))$ las series de Laurent de caracter\'istica $p>0$ y sea
$F$ un subcampo de $E$ estable bajo la derivaci\'on $D:=\frac{d}{dt}$,
$D(F)\subseteq F$. Un elemento $x\in F$ es la derivada logar\'itmica
en $F$ si y solamente si es una derivada logar\'itmica en $E$. $\fin$
\end{corolario}

Sea $K$ un campo de funciones global de caracter\'istica $p>0$. 
Sean $k_0=\F$ el campo de constantes de $K$ y ${\ma F}_p$ es el
campo primo. Sea $K_{\pK}$ la completaci\'on de $K$ con respecto
al primo $\pK$. Sea $K(\pK)$ el campo residual de $K_{\pK}$, 
$\F\subseteq K(\pK)$ por ser $\F$ perfecto y $K_{\pK}$ es de
manera natural isomorfo a $K({\pK})((t))$ donde $t$ es un elemento
uniformizador, $v_{\pK}(t)=1$.

Puesto que $\F$ es perfecto, se puede tomar $t\in K$ que sea un
elemento separador, esto es $K/\F(t)$ es separable (\cite[Corollary 8.2.11
y Remark 8.4.9]{Vil2006}). Sea $\Tr$ la traza de $K(\pK)$ a ${\ma F}_p$.
Sea $\Res(xdy)$ el residuo de la diferencial local $xdy$ calculado con
respecto a un par\'ametro uniformizador.El residuo no depende del
par\'ametro (\cite[Theorem 9.3.9 y Definition 9.3.10]{Vil2006}).

Definimos un pareo local de la siguiente forma. Para $x,y\in K_{\pK}$, sea
\begin{gather}\label{Eq17.6.90-1N}
\int_{\pK} xdy:= \Tr(\Res xdy).
\end{gather}
Si $x\in K_{\pK}$ y $y\in\*K_{\pK}$, se define
\begin{align*}
\varphi_{\pK}(x,y)&=\int_{\pK}x\frac{dy}{y}=\int_{\pK}x\frac{\frac{dy}{dt}}{y}dt.
\intertext{Se verifica que}
\varphi_{\pK}(x+x',y)&=\varphi_{\pK}(x,y)+\varphi_{\pK}(x',y),\\
\varphi_{\pK}(x,yy')&=\varphi_{\pK}(x,y)+\varphi_{\pK}(x,y'),
\end{align*}
por lo que se tiene un pareo de los grupos $K_{\pK}$ (aditivo) y
$\*K_{\pK}$ (multiplicativo) con valores en el grupo aditivo ${\ma F}_p$.
Si $x$ es entero en $K_{\pK}$, esto es, $v_{\pK}(x)\geq 0$ y $y$ es una
unidad en $\*K_{\pK}$, es decir, $v_{\pK}(y)=0$, entonces $x\frac{dy}{y}$
tiene residuo $0$. Esto muestra la continuidad del pareo: $\varphi_{\pK}^{
-1}(0)\supseteq \o_{\pK}\times U_{\pK}$.

\begin{lema}\label{L17.6.91N}
Sea $x\in\o_{\pK}$. Entonces $\varphi_{\pK}(x,y)=0\iff p|v_{\pK}(y)$ o
$x=\wp(z)$ para alg\'un $z\in K_{\pK}$, donde $\wp$ es el operador
de Artin-Schreier: $\wp(z)=z^p-z$.
\end{lema}

\begin{proof}
Sea $x=a_0+a_1t+\cdots$, $a_i\in K(\pK)$, $y=t^n\varepsilon$, 
$\varepsilon\in U_{\pK}$. Entonces $\varphi_{\pK}(x,\varepsilon)=0$.
Por tanto
\[
\varphi_{\pK}(x,y)=\varphi_{\pK}(x,t^n)=\int_{\pK} x\frac{dt^n}{t}=
n\int_{\pK}x\frac{dt}{t}=n\Tr (a_0).
\]

Se sigue que tener $\varphi_{\pK}(x,y)=0$ es equivalente con $p|n$
o $\Tr(a_0)=0$ y debemos probar que $\Tr(a_0)=0$ es equivalente
con la existencia de un $z$ tal que $x=\wp(z)=z^p-z$, $z\in K_{\pK}$.

Primero supongamos que $x=z^p-z$. Si la serie de Laurent para $z$
tiene polos, estos polos ser\'ian dominantes en $z^p$ y entonces $x$
no podr\'ia ser entero: Sea $b$ el t\'ermino constante de $z$, por lo
que $a_0=b^p-b$. Sea $\sigma$ el automorfismo de Frobenius de
$K(\pK)/{\ma F}_p$: $\sigma c=c^p$. De esta forma, $a_0=(\sigma
-1)b$ lo cual es equivalente $\Tr a_0=0$ pues $
\co n{\Gal(K(\pK)/{\ma F}_p)}{K(\pK)} =0$,
en particular, $\co 1{\Gal(K(\pK)/{\ma F}_p)}{K(\pK)}=0$, o 
simplemente, $\Tr a_0=\Tr \sigma b-\Tr b=0$ pues $\sigma b$
y $b$ son conjugado.

Rec\'iprocamente, si $\Tr(a_0)=0$, entonces por la misma raz\'on
anterior,
$a_0=(\sigma-1)b=\wp(b)$ para alg\'un $b\in K(\pK)$. Sea
$x-a_0=x_1=a_1t+a_2t^2+\cdots$ y sea $z_1=-(x_1+x_1^p+x_1^{p^2}+
\cdots )$. La serie converge en la topolog\'ia de $K_{\pK}$ y se obtiene
$x_1=z_1^p-z_1=\wp(z_1)$ y por tanto $x=a_0+x_1=\wp(b)+\wp(z_1)
=\wp(b+z_1)$ y el resultado se sigue.
$\fin$
\end{proof}

\subsubsection{Un pareo global}

La traza $\Tr$ definida en el pareo local (\ref{Eq17.6.90-1N}), se
puede descomponer como 
$\Tr=\Tr_{K(\pK)/{\ma F}_p}=\Tr_{\F/{\ma F}_p}\circ
\Tr_{K(\pK)/\F}$. Sea ${\ma A}_K$ el anillo de ad\`eles. Sea
$\pK$ un primo y sea $t$ un elemento separador para $K$ y
$v_{\pK}(t)=1$. Sea
\[
\lambda(\vec\xi)=\sum_{\pK} \Tr_{\pK}(\Res_{\pK}(\xi_{\pK}dt)),
\]
con $\Tr_{\pK}=\Tr_{K(\pK)/\F}$, $\vec\xi
\in {\ma A}_K$ y $\lambda(\vec\xi)\in\F$. Si
$|\xi_{\pK}|_{\pK}$ es suficientemente peque\~no en todos
los primos donde $t$ tiene polos y $\xi_{\pK}$ es una unidad
en los dem\'as primos, se tiene $\lambda(\vec\xi)=0$ pues
$\Res_{\pK}(\xi_{\pK}dt)=0$ para toda $\pK$. En particular
$\lambda$ es una funci\'on continua sobre ${\ma A}_K$. Por
el Teorema de los residuos (\cite[Theorem 9.3.14]{Vil2006}),
$\lambda(x)=0$ para toda $x\in K$, donde $x$ es el ad\`ele
principal. En otras palabras, $\lambda$ es un diferencial
sobre $K$. Sea
\[
\oint\vec\xi dt =\Tr_{\F/{\ma F}_p}(\lambda(\vec\xi))=\sum_{\pK}
\int_{\pK}\xi_{\pK}dt.
\]

\begin{lema}\label{L17.6.92N}
Si $\vec\xi\in{\ma A}_K$ es tal que $\oint \xi xdt=\Tr_{\F/{\ma F}_p}(
\lambda(\vec\xi x))=0$ para toda $x\in K$, entonces $\vec\xi\in K$.
\end{lema}

\begin{proof}
Podemos reemplazar $x$ por $\alpha x$ por cualquier $\alpha\in
\F$ y obtenemos $\Tr_{\F/{\ma F}_p}(\alpha\cdot\lambda(\vec\xi
x))=0$. Si $\lambda(\vec\xi x)\neq 0$, tendr\'iamos $\Tr_{\F/{\ma
F}_p}(\F)=0$ lo cual no se cumple pues $\F/{\ma F}_p$ es 
separable. Por tanto $\lambda(\vec\xi x)=0$ para toda $x\in K$.

Se tiene que, para $x\in K$, si $\mu$ es una diferencial, $x\mu$ es
una diferencial donde $(x\mu)(\vec\alpha):=\mu(\vec\alpha x)$. Si
seleccionamos $\mu_0\neq 0$, puesto que las diferenciales es un
espacio de dimensi\'on $1$ sobre $K$, toda diferencial es de la
forma $x\mu_0$. Ahora $\lambda\neq 0$ y $(x\lambda)(\vec\xi)=
\lambda(\vec\xi x)=0$, es decir $\mu(\vec\xi)=0$ para toda diferencial
$\mu$. Veamos que necesariamente $\vec\xi\in K$, esto es,
\[
\bigcap_{\mu\text{\ diferencial}}\ker \mu=\bigcap_{{\eu a}\text{\ 
divisor}}(K+{\eu X}_K({\eu a}^{-1}))=K.
\]

Sea $\vec\alpha\xrightarrow[\ \ \delta\ \ ]{} \lambda(\vec\xi \vec\alpha)$.
Se verifica que $\delta$ es una diferencial ($\lambda(\vec\xi x)=\delta
(x)=0$ para toda $x\in K$). Por tanto existe $\beta\in K$ tal que $
\delta =\beta \lambda$, es decir, $\lambda(\vec\xi\vec \alpha)=
\delta(\vec\alpha)=\beta\lambda(\vec
\alpha)=\lambda(\beta\vec\alpha)$.
Por tanto $\lambda_{\pL}(\vec\xi \vec\alpha)=\lambda_{\pL}(\beta
\vec\alpha)$ para toda $\pL\in{\ma P}_K$
lo cual implica $\lambda_{\pL}((\xi_{\pL}
-\beta)\alpha_{\pL})=0$ para toda $\alpha_{\pL}$. Por tanto
$\xi_{\pL}-\beta=0$ para toda $\pL\in {\ma P}_K$, 
de donde se sigue $\vec\xi=\beta\in K$.
$\fin$
\end{proof}

Sea $t$ un elemento separador de $K/\F$, es decir, $K/\F(t)$ es
separable, el cual existe por ser $\F$ perfecto (\cite[Corollary 8.2.11]
{Vil2006}).

\begin{lema}\label{L17.6.93N}
Para cada $\pK\in{\ma P}_K$, sea $t_{\pK}$ un elemento primo de
$K_{\pK}$. Entonces $\frac{dt}{dt_{\pK}}$ es unidad para casi toda
$\pK$.
\end{lema}

\begin{proof}
Tenemos que el divisor de la diferencial $dt$ en $K$ es $(dt)_K
={\eu D}_{K/k(t)}\con_{\F(t)/K}{\mc P}_{\infty}^{-2}$, donde ${\mc P}_{
\infty}$ es el polo de $x$ en $\F(t)$. En particular $v_{\pK}((dt)_K)=0$
para casi toda $\pK$.

Veamos que si $t_{\pK}$ es un elemento primo en $\pK$, entonces si en
$K_{\pK}$, $\frac{dt}{dt_{\pK}}=\sum_{m=n} a_m t_{\pK}^m$ con $a_n\neq
0$, se tiene que $v_{\pK}((dt)_K)=n$.
Una vez probado lo anterior, se seguir\'a que $\frac{dt}{dt_{\pK}}$ es
una unidad para casi toda $\pK$.

Nos referimos a \cite[Subsection 9.3]{Vil2006} donde se ve la relaci\'on entre
las diferenciales $\alpha d\beta$ (de Hasse) y las diferenciales $\omega$ de
Weil. Esto es \cite[Theorem 9.3.15]{Vil2006}: Dada $\alpha d\beta$,
se define $\omega:{\eu X}_K(={\ma A}_K)\lra \F$ por
\begin{gather*}
\omega(\vec \xi)=\sum_{\pK\in{\ma P}_K} \Res_{\pK}(\xi_{\pK}\alpha
d\beta)\quad\text{y se tiene}\quad \omega^{\pK}(\vec \xi) =\Res_{\pK}
\xi_{\pK}\alpha d\beta.
\end{gather*}

En la demostraci\'on del Theorem 9.3.15 de \cite{Vil2006}
se obtiene que el exponente local del divisor $(\omega)_K$, es 
precisamente $v_{\pK}(\alpha d\beta)$. En nuestro caso, $dt=
\frac{dt}{dt_{\pK}}dt_{\pK}$ y $v_{\pK}(dt)=v_{\pK}\big(\frac{dt}{dt_{\pK}}\big)$.
 $\fin$ 
\end{proof}

Ahora sea $\vec\alpha\in J_K$. Sea $\vec\beta$ dada por
$\beta_{\pK}=\frac{1}{\alpha_{\pK}}\frac{d\alpha_{\pK}}{dt}=\frac{1}{
\alpha_{\pK}}\frac{d\alpha_{\pK}}{dt_{\pK}}\big(\frac{dt}{dt_{\pK}}\big)^{-1}$.
Entonces $\beta_{\pK}\in \o_K$ para casi toda $\pK$ y por tanto
$\beta\in{\ma A}_K$. Usamos la notaci\'on $\vec\beta=\frac{1}
{\vec\alpha}\frac{d\vec\alpha}{dt}$.

Sea $\varphi: K\times J_K\lra {\ma F}_p$ el pareo dado por
\begin{align}
\varphi(x,\vec\alpha)&=\oint x\frac{1}{\vec\alpha}\frac{d\vec\alpha}{dt}dt=
\sum_{\pK}\int_{\pK}x\frac 1{\alpha_{\pK}}\frac {d\alpha_{\pK}}{dt}dt=
\sum_{\pK}\int_{\pK}x\frac{d\alpha_{\pK}}{\alpha_{\pK}}\nonumber \\
&=\sum_{\pK}\varphi_{\pK}(x,\alpha_{\pK})=\sum_{\pK}\Tr(\Res
\frac x{\alpha_{\pK}}d\alpha_{\pK}).\label{Ec17.6.93-1N}
\end{align}
Considerando $K$ y ${\ma F}_p$ con la topolog\'ia discreta, se tiene que
si $\alpha_{\pK}$ est\'a muy cerca de $1$, entonces $\frac
{d\alpha_{\pK}}{\alpha_{\pK}}$ est\'a muy cerca de $0$ y viendo las
componentes de los integrandos locales, se sigue la continuidad en 
$J_K$. Por tanto el pareo en (\ref{Ec17.6.93-1N}) es un pareo continuo.

\begin{lema}\label{L17.6.94N}
El n\'ucleo en $J_K$ del pareo {\rm{(\ref{Ec17.6.93-1N})}} es $J_K^p\*K$
y el n\'ucleo en $K$ es $\wp(K)=\{z^p-z\mid z\in K\}$.
\end{lema}

Notemos que cada elemento de $K$ tiene per\'iodo $p$  y que
$J_K/\*KJ_K^p$ es compacto y de hecho, isomorfo a $C_K/C_K^p$.
M\'as precisamente,
se tiene $C_K\cong C_{K,0}\times A$, $A\subseteq {\ma R}^+$,
$C_K^p\cong C_{K,0}^p\times A^p$ y $C_{K,0}$ es compacto,
$A\cong {\ma Z}$ por lo que $C_K/C_K^p\cong C_{K,0}/C_{K,0}^p
\times {\ma Z}/p{\ma Z}$ es compacto.

\begin{proof}
Si $\vec\alpha$ est\'a en el n\'ucleo en $J_K$, 
entonces $\vec\xi =\frac {1}{\vec
\alpha}\frac{d\vec\alpha}{dt}\in {\ma A}_K$
satisface $\oint \vec\xi xdt=0$ para toda
$x\in K$, por lo que $\vec \xi= y\in K$ pues $\omega(\vec\xi)=0$
para toda diferencial de $K$ (Lema \ref{L17.6.92N}).

Tomando la componente $\pK$ de $y$ ($y_{\pK}=y$), que tenga a
$t\in K$ como uniformizador local, entonces $y=\frac{1}{\alpha_{\pK}}
\frac{d\alpha_{\pK}}{dt}$ y por tanto $y$ es una derivada logar\'itmica
en $K_{\pK}$. Puesto que $K\subseteq K_{\pK}$, y $K$ es estable bajo
diferenciaci\'on puesto que $t\in K$, y ya que $y\in K$, el Corolario
\ref{C17.6.90N} prueba que $y$ es derivada logar\'itmica en $K$:
\[
y=\frac{1}{z}\frac{dz}{dt},\qquad z\in K.
\]

Obtenemos 
\[
\vec\xi=\frac{1}{\vec\alpha}\frac{d\vec\alpha}{dt}=
\frac{1}{z}\frac{dz}{dt}=y.
\]
Sea $\vec\gamma=z^{-1}\vec\alpha$. Entonces
$\frac{d\vec\gamma}{dt}=-z^{-2}\frac{dz}{dt}\vec\alpha+z^{-1}\frac{d\vec\alpha}
{dt}=0$. Por tanto $\frac{d\gamma_{\pK}}{dt}=0$ para toda $\pK$
lo cual implica que cada componente es una $p$-potencia por
lo que $\vec\gamma\in J_K^p$ y por tanto $\vec\alpha\in J_K^p\*K$.

Por otro lado, tanto $J_K^p$ como $\* K$ claramente
pertenecen al n\'ucleo del mapeo $\varphi$, $\*K$ puesto que
$\lambda(\vec\xi)=\sum_{\pK}\Tr_{\pK}(\Res_{\pK}(\xi_{\pK}dt))$
es una diferencial. Se sigue que el n\'ucleo de $\varphi(x,\vec\alpha)$
en $J_K$ es $J_K^p\*K$.

Ahora sea $x=z^p-z$. Sea $\vec\xi\in J_K$.
Por el teorema de aproximaci\'on, existe
$y\in K$ tal que $y\vec\xi$ est\'a muy cercano a $1$ en todos los
lugares en donde $x$ tiene un polo. Sea $\pK$ un polo de $x$. Entonces
\begin{gather*}
\varphi_{\pK}(x,y\vec\xi)=\int_{\pK} x\frac{d(y\xi_{\pK})}{y\xi_{\pK}}.
\intertext{Sea $\pi$ un primo en $\pK$. Se tiene}
x\frac{d(y\xi_{\pK})}{y\xi_{\pK}}=x\frac{d(y\xi_{\pK})}{\pi}\Big(
\frac{y\xi_{\pK}}{\pi}\Big)^{-1},\\
y\xi_{\pK}=1+\pi^nb,\quad d(y\xi_{\pK})=
n\pi^{n-1}bd\pi,
\intertext{as\'i que}
\frac{d(y\xi_{\pK})}{\pi}=n\pi^{n-2}bd\pi,\quad \Big(\frac{y\xi_{\pK}}{\pi}
\Big)^{-1}=\pi^{-1}(1-\pi^nb+\cdots)=\pi^{-1}-\pi^{n-1}b+\cdots,\\
\frac{d(y\xi_{\pK})}{\pi}\Big(\frac{y\xi_{\pK}}{\pi}\Big)^{-1}=n\pi^{n-2}+
\cdots,\quad x\frac{d(y\xi_{\pK})}{\pi}\Big(\frac{y\xi_{\pK}}{\pi}\Big)^{-1}
=n\pi^{n-2-m}+\cdots
\end{gather*}
donde $x=\frac{x_0}{\pi^m}$. Por tanto, con $n\geq 2+m$, se tiene
$\Res x\frac{d(y\xi_{\pK})}{y\xi_{\pK}}=0$. Se sigue que $\varphi_{\pK}(
x,y\vec\xi)=0$ para todo lugar $\pK$ que es polo de $x$.

Ahora, como $y$ est\'a en el n\'ucleo de $J_K$. Entonces
\[
\varphi(x,\vec\xi)=\sum_{\pK}\varphi_{\pK}(x,y\vec\xi)=\sum_{\eu q}
\varphi_{\eu q}(x,y\xi_{\eu q}),
\]
donde $\eu q$ var\'ia sobre los primos donde $x$ no tiene polos y
por tanto $x$ es un entero local. Por el Lema \ref{L17.6.91N} se tiene que
$\varphi_{\eu q}(x,y\xi_{\eu q})=0$ pues $x=z^p-z$. Se sigue que
$x$ est\'a en el n\'ucleo de $K$ bajo el pareo $\varphi$.

Rec\'iprocamente, supongamos que $x$ est\'a en el n\'ucleo de $K$.
Entonces, $\varphi(x,\vec\xi)=0$ para toda $\vec\xi\in J_K$, en particular
para los id\`eles $\vec \xi_{\pK}$ que \'unicamente una componente
$\neq 1$, esto es, $(\vec\xi_{\pK})_{\eu q}=\begin{cases}
1&\text{si ${\eu q}\neq \pK$},\\
\xi_{\pK}&\text{si ${\eu q}=\pK$}\end{cases}$, se tiene
 $\varphi(x,\vec\xi_{\pK})=\varphi_{\pK}\big(x,\big(\vec
\xi_{\pK}\big)_{\pK}\big)=0$. Por tanto $\varphi_{\pK}(x,\xi_{\pK})=0$
para toda $\xi_{\pK}\in K_{\pK}$. 

Sea $\pK$ un primo que no es
polo de $x$. En ese primo, se puede aplicar el Lema \ref{L17.6.91N}
y se tiene $x=z_{\pK}^p-z_{\pK}$ con alg\'un $z_{\pK}\in K_{\pK}$.

Se tiene que $\pK$ se descompone en $K\big(\wp^{-1} x\big)/K$ pues
$\big(K(\wp^{-1} x)\big)_{\pK}=\big(K(z_{\pK})\big)_{\pK}=K_{\pK}$
o, simplemente, $Y^p-Y-x=\prod_{i=0}^{p-1}(Y-z_{\pK}-i)$ se 
descompone totalmente m\'odulo $\pK$, por lo tanto $\pK$ se
descompone totalmente (Teorema \ref{T6.3.Ram6-1}).

Como esto se cumple para toda $\pK$, como consecuencia del
Corolario \ref{C17.6.75N}, se sigue que $K(\wp^{-1}x)=K$ por lo
que $x=z^p-z$, $z\in K$.
$\fin$
\end{proof}

\subsubsection{Grupo de Galois de la m\'axima $p$-extensi\'on elemental
abeliana de $K$}

Esta extensi\'on es $L=K(\wp^{-1} (K))$. Sea ${\mc G}=\Gal(L/K)$. El grupo
asociado a $L$ es $A=\wp^{-1}(K)$. El pareo es $(\alpha,\sigma)=(\sigma
-1)\alpha$. $A$ se mapea en $K$ por el mapeo $\alpha\mapsto \wp
\alpha$. La imagen de $A$ es $K$. Se puede ver como el pareo $K
\times {\mc G}\xrightarrow{[\underline{\ },\underline{\ }]}{}
{\ma F}_p$ donde $[x,\sigma]=(\sigma-1)(\wp^{-1}(x))$,
$x\in K$, $\sigma\in {\mc G}$.

El n\'ucleo en $K$ es $\wp (K)$, el mismo que el del pareo $\varphi(x,
\vec\xi)$. El n\'ucleo en ${\mc G}$ es $\{1\}$. Con este pareo, se tiene
${\mc G}\cong\widehat{\big(K/\wp(K)\big)}=\Hom(K/\wp(K),{\ma F}_p)$
donde $K/\wp(K)$ tiene la topolog\'ia discreta. Con el pareo
$\varphi(x,\vec\xi)$, el grupo $J_K/J_K^p\*K$ es naturalmente
isomorfo al mismo grupo $\widehat{K/\wp(K)}$.

Se tiene un isomorfismo natural 
\[
{\mc G}\cong J_K/J_K^p\*K=C_K/C_K^p.
\]
Si $\sigma\in
{\mc G}$ es la imagen del id\`ele $\vec\alpha$ bajo este isomorfismo,
entonces debemos tener $\varphi(x,\vec\alpha)=[x,\sigma]=(\sigma
-1)(\wp^{-1}(x))$.

Denotamos este elemento $\sigma\in{\mc G}$
por $\sigma_{\vec\alpha}$. Esto es, el isomorfismo ${\mc G}\cong
J_L/J_K^p\*K$ lo describimos por $\mu:J_L\stackrel{\mu}{\lra}
{\mc G}$, $\mu(\vec\alpha):=\sigma_{\vec\alpha}$ y donde $\ker
\mu=J_K^p\*K$. La ecuaci\'on
anterior describe el efecto de $\sigma_{\vec\alpha}$ en los generadores
de $L$, una descripci\'on que determina $\sigma_{\vec\alpha}$
completamente:
\begin{gather}\label{Ec17.6.94+1N}
\sigma_{\vec\alpha}:\wp^{-1}(x)\lra \wp^{-1}(x) +\varphi(x,\vec\alpha),
\end{gather}
donde $\varphi$ es el pareo global $\varphi:{\mc G}\times L\lra 
{\ma F}_p$, y si $\vec\alpha\in J_K$ se mapea a 
$\sigma_{\vec\alpha}\in{\mc G}$
bajo $\mu$, tenemos que (\ref{Ec17.6.94+1N}) se lee $\sigma_{
\vec\alpha}(z)=z+\varphi(x,\vec\alpha)$ con $\varphi(x,\vec\alpha)\in
{\ma F}_p$ y donde $z^p-z=x$.

De esta forma, hemos probado

\begin{teorema}\label{T17.6.95N}
El mapeo {\rm{(\ref{Ec17.6.94+1N})}} es un isomorfismo bicontinuo
entre $J_K/J_K^p\*K\cong C_K/C_K^p$ y el grupo de Galois de la
m\'axima extensi\'on abeliana $L$ de $K$ de exponente $p$.
$\fin$
\end{teorema}

\subsubsection{Demostraci\'on de la 2nda. desigualdad en caracter\'istica
$p>0$}

Sea $M/K$ una extensi\'on c\'iclica de grado $p$ de campos globales
de funciones. Como $M$ es subcampo de la m\'axima $p$-extensi\'on 
elemental abeliana de $K$, $L=K(\wp^{-1}(K))$, $M$ 
est\'a determinado por un
subgrupo abierto ${\mc H}$ de ${\mc G}=\Gal(L/K)$ de \'indice $p$,
${\mc H}=\Gal(L/M)$, $[{\mc G}:{\mc H}]=p$.

El grupo de normas $\*K\N_{M/K} J_M$ es un subgrupo abierto
de $J_K$ (Teorema \ref{T17.6.61N}): todos los id\`eles en
un vecindad suficientemente peque\~na de $1$
son normas, y es de
\'indice finito en $J_K$, pues $J_K/\*K\N_{M/K}J_M
\cong C_K/\N_{M/K}C_M$ y si $\Lambda=C_{K,0}\cap \N_{M/K}
C_M$, entonces $\Lambda$ es abierto en $C_{K,0}$.
Como $C_{K,0}$ es compacto, se sigue que
$[C_{K,0}:\Lambda]<\infty$. Puesto que el campo de constantes
de $M$ es ${\ma F}_{q^{\nu}}$ con $\nu\in\{1,p\}$, se tiene que
$\Delta=\N_{M/K}C_M\cap {\ma Z}\neq \{0\}$ de donde se sigue
que ${\ma Z}/\Delta$ es un grupo finito y por tanto $[C_K:\N_{M/K}C_M]
<\infty$.

 Por la primera desigualdad, se tiene que 
 \[
 [C_K:
\N_{M/K}C_M]=[J_K: \*K \N_{M/K}J_M]\geq p=[M:K].
\]
La imagen de este 
grupo bajo el isomorfismo del Teorema \ref{T17.6.95N} 
(el epimorfismo $\mu$) es cierto 
subgrupo abierto ${\mc H}'$ de ${\mc G}$. 
Esquem\'aticamente, tenemos:
\[
\xymatrix{
\*K\N_{M/K}J_M\ar@{~>}[r]&{\mc H}',
}
\qquad
\xymatrix{
M\ar@{~>}[r]&{\mc H}.
}
\]

El \'indice de ${\mc H}'$ en
${\mc G}$ es finito y mayor o igual $p$. Se tiene que ${\mc H}'$ 
determina cierto subcampo $E/K$ de $L/K$ el cual es la composici\'on
de campos c\'iclicos de grado $p$ y cada uno de estos campos se 
quedan  fijos bajo ${\mc H}'$, $E=L^{{\mc H}'}$.

A cada subcampo c\'iclico $M'/K$ de 
grado $p$ que sea diferente a $M$, daremos un elemento de ${\mc H}'$
que no deja fijo a todos los elementos de $M'$. Esto probar\'a que
${\mc H}'={\mc H}$ y probar\'a completamente la 
segunda desigualdad puesto
que el \'indice de ${\mc H}$ en ${\mc G}$ es $p$.

M\'as precisamente, si $[M':K]=p$ con $M'/K$ abeliana y $M'\neq M$,
se probar\'a que $M\not\subseteq L^{{\mc H}'}$ y en consecuencia,
como $L^{{\mc H}'}\neq K$, entonces $L^{{\mc H}'}\subseteq M=
L^{\mc H}$ por lo que ${\mc H}\subseteq {\mc H}'$ y $p=[{\mc G}:
{\mc H}]=[{\mc G}:{\mc H}'][{\mc H}':{\mc H}]\geq [{\mc G}:{\mc H}']\geq p$.
Por tanto obtendremos que $p=[{\mc G}:{\mc H}]=[{\mc G}:{\mc H}']$
y ${\mc H}={\mc H}'$.

Por el Teorema \ref{T17.6.76N}, existe un primo (de hecho, una
infinidad) $\pK$ en $K$ que se
descompone en $M/K$ y es inerte en $M'/K$. Si $K'=K(\wp^{-1} (x))$
se puede seleccionar $\pK$ que no sea polo de $x$. Sea $\vec
\alpha_{\pK}\in J_K$ que tiene en la componente $\pK$ un elemento
primo de $K_{\pK}$ y las dem\'as componentes son $1$. 

Se tiene que $(\vec\alpha_{\pK})_{\pK}\in K_{\pK}$ es norma de
de $M_{\pL}$, $\pL|\pK$, pues $\pK$ se descompone en $M/K$
y por tanto $M_{\pL}=K_{\pK}$. Las dem\'as componentes son $1$
y por tanto son normas locales. Del Teorema \ref{T17.6.68N}, se
tiene que $\vec\alpha_{\pK}$ es una norma de alg\'un elemento de $J_M$.
Se sigue que $\sigma_{\vec\alpha_{\pK}}\in{\mc H}'$. Ahora, calculamos
$\varphi(x,\vec\alpha_{\pK})=\varphi_{\pK}(x,(\vec\alpha_{\pK})_{\pK})$.

Como $\pK$ no se descompone en $M'$, se tiene que $x\notin \wp(
K_{\pK})$. 
Ahora $(\vec\alpha_{\pK})_{\pK}$ es un elemento primo,
por tanto de valuaci\'on $1$ y en particular
no es divisible por $p$. Por el
Lema \ref{L17.6.94N}, se tiene que $\varphi(x,\vec\alpha_{\pK})\neq
0$. Esto significa, de acuerdo a (\ref{Ec17.6.94+1N}), que $\sigma_{
\vec\alpha_{\pK}}$ mueve $\wp^{-1}(x)$:
\[
\sigma_{\vec\alpha_{\pK}}(\wp^{-1}(x))=\wp^{-1}(x)+\varphi_{\pK}(
x,(\vec\alpha_{\pK})_{\pK})\neq \wp^{-1}(x),
\]
y por tanto no es la identidad en $M'$.

\begin{teorema}[Segunda desigualdad fundamental\index{Segunda
desigualdad fundamental}]\label{T17.6.96N}
Para cualquier extensi\'on finita de Galois $L/K$ de campos globales,
se tiene que 
\begin{gather*}
[C_K:\N_{L/K}C_L]\big| [L:K]. \tag*{$\fin$}
\end{gather*}
\end{teorema}

Como consecuencia de la primera y de la segunda desigualdad 
obtenemos.

\begin{teorema}\label{T17.6.97N}
Si $L/K$ es una extensi\'on finita c\'iclica 
de grado primo $p$ de campos
globales, $\co 0G{C_L}\cong \co 2G{C_L}\cong G=\Gal(L/K)$ y
$\co 1G{C_L}=1$.
\end{teorema}

\begin{proof}
Se obtuvo de la primera desigualdad que $h(G,C_L)=p=[L:K]$,
por tanto $|\co 0G{C_L}|\geq p$ y por la segunda desigualdad
que $|\co 0G{C_L}|=[C_K:\N_{L/K}C_L]\leq p$. Se sigue que
$|\co 0G{C_L}|=p$ y $|\co 1G{C_L}|=1$.

Ahora bien, puesto que $\co 0G{C_L}$ tiene orden primo $p$, es
c\'iclico lo mismo que $G$, por lo que $\co 0G{C_L}\cong \co 2G{
C_L}\cong G$.
$\fin$
\end{proof}

\begin{teorema}\label{T17.6.98N}
Si $L/K$ es una extensi\'on finita de Galois de campos globales con
grupo de Galois $G=\Gal(L/K)$, entonces se tiene $\co 1G{C_L}=\{1\}$.
\end{teorema}

\begin{proof}
Si $L/K$ es c\'iclica, entonces por la primera desigualdad, se tiene
$|\co 0G{C_L}|\geq [L:K]$ y por la segunda desigualdad $|\co 0G{C_L}|
\leq [L:K]$, de donde $|\co 0G{C_L}|=[L:K]=h(G,C_L)=\frac{|\co 0G{
C_L}|}{|\co 1G{C_L}|}$ por lo que $|\co 1G{C_L}|=1$. Se sigue que
$\co 1G{C_L}=\{1\}$.

Sea $G$ arbitrario de orden $n=[L:K]$. Por inducci\'on suponemos que
$\co 1H{C_L}=\{1\}$ para toda extensi\'on $L/K$ de Galois con grupo
$H$ y de orden menor a $n$.

Si $G$ no es una $p$-potencia, para todo $p$-subgrupo de Sylow
$S_p$ de $G$, se tiene $\co 1{S_p}{C_L}=\{1\}$. Del Corolario
\ref{C17.5.6.10-2} se sigue que $\co 1G{C_L}=\{1\}$.

Sea ahora $|G|=p^m$ con $p$ un n\'umero primo y $m\geq 1$. 
Si $m=1$ ya se tiene el resultado. Sea $m\geq 2$. Sea $H\normal
G$ de \'indice $p$ y sea $K\subseteq M\subseteq L$ con $H=
\Gal(L/M)$. Por hip\'otesis de inducci\'on tenemos $\co 1H{C_L}
=\{1\}=\co 1{G/H}{C_M}$. De la sucesi\'on inflaci\'on-restricci\'on
\[
1\lra \co 1{G/H}{C_M}\stackrel{\infla}{\lra} \co 1G{C_L}\stackrel{
\res}{\lra} \co 1H{C_L},
\]
se sigue que $\co 1G{C_L}=\{1\}$.
$\fin$
\end{proof}

\begin{corolario}[Teorema de la norma de Hasse\index{teorema
de la norma de Hasse}\index{Hasse!teorema de la norma de
Hasse}]\label{C17.6.99N}
Sea $L/K$ una extensi\'on c\'iclica finita de campos globales.
Entonces un elemento $x\in\*K$ es norma de un elemento de $\*L$
si y solamente si $x$ en un norma local para toda completaci\'on
$L_{\pL}/K_{\pK}$, $\pK\in{\ma P}_K$, $\pL|\pK$.
\end{corolario}

\begin{proof}
La sucesi\'on exacta de $G$-m\'odulos 
\begin{gather*}
1\lra \*L\lra J_L\lra C_L\lra 1,
\intertext{da en cohomolog\'ia la sucesi\'on exacta}
\co {{-1}}G{C_L}=\{1\}\lra \co 0G{\*L}\stackrel{\tilde\theta}{\lra}
 \co 0G{J_L}\cong
\bigoplus_{\pK}\co 0{G_{\pL}}{L_{\pL}},
\end{gather*}
por lo que $\tilde\theta$ es inyectivo.

Como $\co 0G{\*L}\cong \*K/\N_{L/K}\*L$ y $\co 0G{J_L}
\cong \bigoplus_{\pK}\*K_{\pK}/\N_{L_{\pL}/K_{\pK}}\*L_{\pL}$,
se tiene que el homomorfismo
\begin{eqnarray*}
\*K&\stackrel{\theta}{\lra} & \bigoplus_{\pK}\*K_{\pK}/
\N_{L_{\pL}/K_{\pK}}\*L_{\pL},\\
x&\longmapsto& \bigoplus_{\pK}x\bmod \N_{L_{\pL}/K_{\pK}}\*L_{\pL},
\end{eqnarray*}
tiene n\'ucleo $\*K\bigcap\big(\bigcap_{\pK}
\N_{L_{\pL}/K_{\pK}}\*L_{\pL} \big)=\N_{L/K}\*L$, de donde
se sigue el resultado.
$\fin$
\end{proof}

\begin{observacion}\label{O17.6.100-1N}
El teorema de Hasse \'unicamente es v\'alido para extensiones c\'iclicas
de campos globales. La demostraci\'on del Corolario \ref{C17.6.99N}
es \'unicamente v\'alida para extensiones c\'iclicas pues estamos usando
que $\co {{-1}}G{C_L}=\co 1G{C_L}=\{1\}$ lo cual es falso en general
para extensiones no c\'iclicas.
\end{observacion}

\begin{teorema}\label{T17.6.100N}
Sea $L/K$ una extensi\'on finita de Galois de campos globales
con grupo $G=G_{L|K}$. Entonces el orden de $\co 2G{C_L}$ es un
divisor de $[L:K]$.
\end{teorema}

\begin{proof}
Lo probamos por inducci\'on en $n=[L:K]$. Si $n=1$ o $p$, con $p$ un
n\'umero primo ya se tiene (Teorema \ref{T17.6.97N}). Sea $n$ un
n\'umero natural que no es $1$ ni primo. Supongamos el resultado
para toda extensi\'on de grado menor a $n$. Si $n$ no es potencia de un
primo, entonces para cada subgrupo de Sylow $S_p$ de $G$, se tiene
$|\co 2{S_p}{C_L}|\big|n_p=|S_p|$. Sea $\co 2G{C_L}_p$ el $p$-subgrupo de
Sylow de $\co 2G{C_L}$, entonces la restricci\'on $\co 2G{C_L}_p
\stackrel{\res}{\lra} \co 2{S_p}{C_L}$ es inyectiva (Proposici\'on 
\ref{P17.5.6.10-1}). Por tanto, $|\co 2G{C_L}_p|\big||\co 
2{S_p}{C_L}|\big|n_p$ y
como $\co 2G{C_L}\cong \bigoplus_{\pK}\co 2{S_p}{C_L}$
se sigue que $|\co 2G{C_L}|\big|\prod_p n_p=n$.

Sea $G$ un $p$-grupo, $|G|=p^m$, $m\geq 2$. Sea $H\normal G$
de \'indice $p$. Entonces $|H|=\frac np$ y tomando en 
cuenta que $\co 1G{C_L}=\{1\}$, tenemos la sucesi\'on
inflaci\'on-restricci\'on
\[
1\lra \co 2{G/H}{C_L^H}\stackrel{\infla}{\lra}\co 2G{C_L}\stackrel{
\res}{\lra} \co 2H{C_L},
\]
la cual es exacta.

Ahora puesto que $|G/H|=p$ se tiene que 
$G/H$ es c\'iclico, y se tiene $\co 2{G/H}{C_L^H}\cong
\co 0{G/H}{C_L^H}\cong G/H$ ($C_L^H=C_M$ donde $M=L^H$). Por
tanto $\frac{|\co 2G{C_L}|}{p}|\co 2H{C_L}|\big|\frac np$. Se sigue que
$|\co 2G{C_L}|\big|n=[L:K]$.
$\fin$
\end{proof}

\begin{observacion}\label{O17.6.101N}
Del Teorema \ref{T17.6.100N} no se obtiene nuestro objetivo de probar
que $\co 2G{C_L}$ es c\'iclico de orden $[L:K]$.
\end{observacion}

\begin{teorema}[Axioma de la teor\'ia de campos de clases
globales\index{axioma de la teor\'ia de campos de clase}]\label{T17.6.102N}
Sea $L/K$ una extensi\'on c\'iclica finita de campos globales, entonces
\[
|\co i{G_{L|K}}{C_L}|=
\begin{cases}
[L:K],& i\equiv 0\bmod 2,\\
1, &i\equiv 1\bmod 2.
\end{cases}
\]
\end{teorema}

\begin{proof}
El resultado se sigue de que $h(G_{L|K},C_L)=n=[L:K]=\frac{|
\co 0{G_{L|K}}{C_L}|}{|\co 1{G_{L|K}}{C_L}|}$ y $\co 1{G_{L|K}}{C_L}=
\{1\}$ por tanto $|\co 0{G_{L|K}}{C_L}|=[L:K]$.
$\fin$
\end{proof}

Para tener el teorema de reciprocidad, el objetivo ser\'a probar que si $L/K$ es 
una extensi\'on abeliana finita de campos globales, entonces $C_L$ satisface
las condiciones del teorema de Tate-Nakayama:
\l
\item $\co 1G{C_L}=\{1\}$,
\item $\co 2G{C_L}$ es c\'iclica de orden $|G|=[L:K]$, 
\end{list}
donde $G=\Gal(L/K)=G_{L|K}$.

Empezaremos por definir el mapeo invariante
\[
\co 2{G_{L|K}}{C_L}\lra \Big(\frac{1}{[L:K]}{\ma Z}\Big)/{\ma Z}.
\]

Como de costumbre, para cada $\pK\in{\ma P}_K$ seleccionaremos 
\'unicamente un primo $\pL\in{\ma P}_L$ con $\pL|\pK$.

Para cada extensi\'on $L_{\pL}/K_{\pK}$, $\Gal(L_{\pL}/K_{\pK})\cong
G_{\pL}\subseteq G_{L|K}$ donde $G_{\pL}$ es el grupo de descomposici\'on
$D(\pL|\pK)$, se tiene el isomorfismo invariante
\[
\inv_{L_{\pL}|K_{\pK}}:\co 2{G_{L_{\pL}|K_{\pK}}}{\*L_{\pL}}\lra \Big(
\frac{1}{[L_{\pL}:K_{\pK}]}{\ma Z}\Big)/{\ma Z}\subseteq \Big(
\frac{1}{[L:K]}{\ma Z}\Big)/{\ma Z},
\]
y cada $\inv_{L_{\pL}|K_{\pK}}$ es la composici\'on de tres isomorfismos,
ver Definici\'on \ref{CCLD17.6.5}. Se tiene $\co 2{G_{L|K}}{J_L}\cong
\bigoplus_{\pK} \co 2{G_{L_{\pL}|K_{\pK}}}{\*L_{\pL}}$, el cual es un
grupo infinito.

\begin{definicion}\label{D17.6.103N}
Sea $c\in \co 2{G_{L|K}}{J_L}$. Descomponemos $c$ as\'i: 
$c=\oplus_{\pK} c_{\pK}$, $c_{\pK}\in
\co 2{G_{L_{\pL}|K_{\pK}}}{\*L_{\pL}}$. Se define
\begin{gather*}
\inv_{L|K}: \co 2{G_{L|K}}{J_L}\lra \Big(\frac{1}{[L:K]}{\ma Z}\Big)/{\ma Z},
\intertext{por}
\inv_{L|K} c:= \sum_{\pK}\inv_{L_{\pL}|K_{\pK}} c_{\pK},
\end{gather*}
donde $c_{\pK}$ es la componente en $\pK$
de $c\in \co 2{G_{L|K}}{J_L}\cong \bigoplus\limits_{\pK\in
{\ma P}_K}\co 2{G_{L_{\pL}|K_{\pK}}}{\*L_{\pL}}$.
\end{definicion}

Notemos que $c_{\pK}=1$ para casi toda $\pK$ y por tanto 
$\inv_{L_{\pL}|K_{\pK}} c_{\pK}=0$ para casi toda $\pK$.

\begin{proposicion}\label{P17.6.104N}
Si $M\supseteq L\supseteq K$ son extensiones de Galois de $K$ de
campos globales, entonces
\lasa
\item $\inv_{M|K} c=\inv_{L|K} c,\quad c\in \co 2{G_{L|K}}{J_L}\subseteq
\co 2{G_{M|K}}{J_M}$ bajo el mapeo de inflaci\'on puesto que 
$\co 1{G_{M|K}}{J_M}=\{1\}$.

\item $\inv_{M|L}(\res_L c)=[L:K]\cdot \inv_{M|K} c,\quad c\in \co 2{G_{
M|K}}{J_M}$.

\item $\inv_{M|K}(\cores_K c)=\inv_{M|L} c, \quad c\in \co 2{G_{M|L}}{J_M}$.
\end{list}

Para {\rm{(b)}} y {\rm{(c)}} \'unicamente se requiere que $M/K$ sea de
Galois.

En otras palabras, los siguientes diagramas son conmutativos:
\lasa
\setcounter{3bean}{1}
\item $
\xymatrix{
\co 2{G_{M|L}}{J_M}\ar@{->}[rr]^{\inv_{M|L}}\ar@{<-}[d]_{\res_L}&&
\Big(\frac{1}{[M:L]}{\ma Z}\Big)/{\ma Z}\ar@{<-}[d]^{[L:K]}\\
\co 2{G_{M|K}}{J_M}\ar@{->}[rr]^{\inv_{M|K}}&&\Big(\frac{1}{[M:K]}
{\ma Z}\Big)/{\ma Z}
}
$

\item $
\xymatrix{
\co 2{G_{M|L}}{J_M}\ar@{->}[rr]^{\inv_{M|L}}\ar@{->}[d]_{\cores_K}&&
\Big(\frac{1}{[M:L]}{\ma Z}\Big)/{\ma Z}\ar@{->}[d]^{\text{$j$ encaje}}\\
\co 2{G_{M|K}}{J_M}\ar@{->}[rr]^{\inv_{M|K}}&&\Big(\frac{1}{[M:K]}
{\ma Z}\Big)/{\ma Z}
}
$
\end{list}

Si $J=J_{\sep K}=\lim\limits_{\substack{\lra\\ L}} J_L$ y 
$C=C_{\sep K}=\lim\limits_{\substack{\lra\\ L}} C_L=
J/(\sep K)^*$ donde el l\'imite se toma sobre todas las extensiones 
finitas de Galois $L/K$, entonces obtenemos un mapeo
\[
\inv:\co 2KJ\lra {\ma Q}/{\ma Z},
\]
con $\co 2KJ=\co 2{\Gal(\sep K/K)}J$, con la
inclusi\'on $\co 2{G_{L|K}}{J_L}\lra \co 2{G_{M|K}}{J_M}$,
para $K\subseteq L\subseteq M$.
\end{proposicion}

\begin{proof}
La proposici\'on se sigue del comportamiento de los invariantes locales con
respecto a los mapeos de inclusi\'on, restricci\'on y corestricci\'on.
\lasa
\item Sea $c\in \co 2{G_{L|K}}{J_L}$, entonces
\[
\inv_{M|K} c=\sum_{\pK} \inv_{M_{\eu q}|K_{\pK}} c_{\pK}=\sum_{\pK}
\inv_{L_{\pL}|K_{\pK}} c_{\pK}=\inv_{L|K} c,
\]
donde ${\eu q}$ es cualquier primo sobre $\pL$.

\item Sea $c\in \co 2{G_{M|K}}{J_M}$ y $\pL$ recorriendo los primos
de $L$. Entonces
\begin{align*}
\inv_{M|L}(\res_L c)&=\sum_{\pL}\inv_{M_{\eu q}|L_{\pL}}(\res_L c)_{\pL}=
\sum_{\pL}\inv_{M_{\eu q}|L_{\pL}}(\res_{L_{\pL}} c_{\pK})\\
&=\sum_{\pL}[L_{\pL}:K_{\pK}]\cdot \inv_{M_{\eu q}|K_{\pK}} c_{\pK}\\
&=\sum_{\pK}\sum_{\pL|\pK}[L_{\pL}:K_{\pK}]\cdot \inv_{M_{\eu q}|K_{\pK}}
c_{\pK},
\end{align*}
con ${\eu q}|_L=\pL$ y $\pL|_K=\pK$.

Se tiene que $\inv_{M_{\eu q}|K_{\pK}} c_{\pK}$ son independientes del
primo ${\eu q}$ sobre el $\pL$ seleccionado, pues si ${\pL}_1$ es otro
primo de $L$, el $K_{\pK}$-isomorfismo can\'onica $L_{\pL}\stackrel{\cong}
{\lra} L_{{\pL}_1}$ da un isomorfismo can\'onico $\co 2{G_{L_{\pL}|K_{\pK}}}{
\*L_{\pL}}\stackrel{\cong}{\lra} \co 2{G_{L_{{\pL}_1}|K_{\pK}}}{\*L_{{\pL}_1}}$
y el cual preserva el invariante. Por otro lado, tenemos $\sum_{\pL|\pK} [L_{
\pL}:K_{\pK}]=[L:K]$, por lo que
\begin{align*}
\inv_{M|L}(\res_L c)&=\sum_{\pK}\Big(\sum_{\pL|\pK} [L_{\pL}:K_{\pK}]\Big)
\inv_{M_{\eu q}|K_{\pK}} c_{\pK}\\
&=[L:K]\sum_{\pK}\inv_{M_{\eu q}|K_{\pK}}
c_{\pK}=[L:K]\inv_{M|L} c.
\end{align*}

\item Para $c\in \co 2{G_{M|L}}{J_M}$, se tiene que
\begin{align*}
\inv_{M|K}(\cores_K c)&=\sum_{\pK}\inv_{M_{\eu q}|K_{\pK}}
(\cores_K c)_{\pK}\\
&=\sum_{\pK}\sum_{\pL|\pK}\inv_{M_{\eu q}|K_{\pK}}
(\cores_{K_{\pK}} c_{\pL})\\
&=\sum_{\pK}\sum_{\pL|\pK}\inv_{M_{\eu q}|L_{\pL}} c_{\pL}=\inv_{M|L} (c).
\tag*{$\fin$}
\end{align*}
\end{list}
\end{proof}

\begin{observacion}\label{O17.6.105N}
Se tiene todo para poder aplicar el teorema de Tate-Nakayama, excepto 
que $\inv_{L|K}:\co 2{G_{L|K}}{J_L}\lra \Big(\frac{1}{[L:K]}{\ma Z}\Big)/{\ma Z}$
no es un isomorfismo (el primer grupo es infinito y el
segundo finito). Para hacer de esto un isomorfismo, necesitamos
pasar de $J_L$ a $J_L/\*L=C_L$.
\end{observacion}

\begin{definicion}\label{D17.6.106N}
Sea $L/K$ una extensi\'on abeliana finita de campos globales. Para cada
$\pK\in{\ma P}_K$, fijamos un $\pL\in{\ma P}_L$ con $\pL|_K=\pK$. Sea
$\vec\alpha\in J_K$, entonces definimos el mapeo
\[
[\vec\alpha,L/K]:=\prod_{\pK\in{\ma P}_K}(\alpha_{\pK},L_{\pL}/K_{\pK})
\in G_{L|K}.
\]
Aqu\'i, $(\alpha_{\pK},L_{\pL}/K_{\pK})\in G_{\pL}$ para toda $\pK\in
{\ma P}_K$ y $(\alpha_{\pK},L_{\pL}/K_{\pK})=1$ para casi toda
$\pK$ y consideramos $G_{\pL}=G_{L_{\pL}|K_{\pK}}\subseteq G_{
L|K}$. Se tiene que $\alpha_{\pK}$ es una unidad para casi todo $\pK$
y el n\'umero de primos ramificados es finito, por lo que $(\alpha_{\pK},
L_{\pL}/K_{\pK})=1$, para casi toda $\pK$.
\end{definicion}

El producto es independiente del orden de los factores pues el grupo
$G_{L|K}$ es abeliano. A la extensi\'on $L_{\pL}/K_{\pK}$ la denotamos
$L_{\pK}/K_{\pK}$ para enfatizar que para cada $\pK\in{\ma P}_K$
seleccionamos \'unicamente un primo $\pL\in{\ma P}_L$ con $\pL|\pK$.

\begin{proposicion}\label{P17.6.107N}
Sea $L/K$ una extensi\'on abeliana finita de campos globales, $\vec
\alpha\in J_K$ y denotamos $(\vec\alpha):=\vec\alpha\N_{L/K} J_L\in
\co 0{G_{L|K}}{J_L}$. Si $\mu\in\chi(G_{L|K})=\co 1{G_{L|K}}{{\ma Q}/
{\ma Z}}$, entonces
\[
\mu\big([\vec\alpha, L/K]\big)=\inv_{L|K}\big((\vec\alpha)\Cup \delta\mu\big)
\in \Big(\frac {1}{[L:K]}{\ma Z}\Big)/{\ma Z},
\]
$(\vec\alpha)\in \co 0{G_{L|K}}{J_L}$, $\delta\mu\in \co 2{G_{L|K}}{\ma Z}$,
$\big((\vec\alpha)\Cup \delta\mu\big)\in \co 2{G_{L|K}}{J_L}$.
\end{proposicion}

\begin{proof}
El resultado es consecuencia de la f\'ormula an\'aloga que relaciona el
s\'imbolo de la norma residual local $(\underline{\ },L_{\pL}/K_{\pK})$ con el
invariante local $\inv_{L_{\pL}|K_{\pK}}$.

Sea $\mu_{\pK}$ la restricci\'on de $\mu$ a $G_{L_{\pL}|K_{\pK}}$ y sea
$(\alpha_{\pK})=\alpha_{\pK}\N_{L_{\pL}/K_{\pK}}\*L_{\pL}$. Entonces
\[
\mu\big([\vec\alpha,L/K]\big)=\sum_{\pK}\mu_{\pK}(\alpha_{\pK},L_{\pK}
|K_{\pK})=\sum_{\pK}\inv_{L_{\pL}|K_{\pK}}\big((\alpha_{\pK})\Cup
\delta\mu_{\pK}\big).
\]

Adem\'as, se tiene que $\big((\alpha_{\pK})\Cup \delta\mu_{\pK}\big)\in
\co 2{G_{L_{\pL}|K_{\pK}}}{\*L_{\pL}}$ son las componentes locales de
$(\vec\alpha)\Cup \delta\mu\in \co 2{G_{L|K}}{J_L}$, \'unicamente 
necesitamos notar que $\alpha_{\pK}\cdot \delta\mu_{\pK}(\sigma,\tau)$
(resp. $\vec\alpha\cdot \delta\mu(\sigma,\tau)$) es un $2$-cociclos de la
clase $\big((\alpha_{\pK})\Cup \delta\mu_{\pK}\big)$ (resp. $\big((\vec\alpha
)\Cup\delta\mu\big)$). Por tanto
\begin{gather*}
\mu\big([\vec\alpha,L/K]\big)=\inv_{L|K}\big((\vec\alpha)\Cup \delta\mu\big).
\tag*{$\fin$}
\end{gather*}
\end{proof}

Cuando queramos definir $\inv_{L|K}$ en $C_L$, el siguiente resultado
es de capital importancia. De la sucesi\'on exacta
\[
1\lra \*L\lra J_L\lra C_L\lra 1,
\]
obtenemos del hecho de que $\co 1{G_{L|K}}{C_L}=\{1\}$ que el 
homomorfismo $\co 2{G_{L|K}}{\*L}\lra \co 2{G_{L|K}}{J_L}$ es
inyectivo. Con esto en mente, usaremos que $\co 2{G_{L|K}}{\*L}
\subseteq \co 2{G_{L|K}}{J_L}$.

\begin{teorema}\label{T17.6.108N}
Sea $L/K$ una extensi\'on finita de Galois de campos globales.
Si $c\in\co 2{G_{L|K}}{\*L}$, entonces $\inv_{L|K} c=0$.
\end{teorema}

\begin{proof}
La demostraci\'on est\'a basada en la descripci\'on expl\'icita del s\'imbolo
de la norma residual y de la f\'ormula del producto.

Sea $L/K$ una extensi\'on finita de Galois de campos globales. Sea $K_0
={\ma Q}$ o $K_0=K$, dependiendo de si $K$ es num\'erico o de funciones.
Sea $M$ una extensi\'on de Galois de $K_0$ que contiene a $L$. Entonces
\begin{gather*}
c\in\co 2{G_{L|K}}{\*L}\subseteq \co 2{G_{M|L}}{\*M}\subseteq \co 2{G_{
M|K}}{J_M},\\
\cores_{K_0} c\in \co 2{G_{M|K_0}}{\*M} \text{\ y\ }
\inv_{L|K} c=\inv_{M|K} c=\inv_{M|K_0}(\cores_{K_0} c).
\end{gather*}
Por tanto, para probar que $\inv_{L|K} c=0$, es suficiente considerar el caso
$K=K_0$.

Por la estructura del grupo de Brauer, existe una extensi\'on $L_0/K_0$ con
$L_0$ c\'iclica ciclot\'omica o de constantes, con $c\in \co 2{G_{L_0|K_0}}{\*L_0}$.
En el caso num\'erico podemos suponer que $L/K_0$ es c\'iclica ciclot\'omica
y que $L/K_0$ es extensi\'on de constantes en el caso de campos de 
funciones.

Sea $\mu$ un generador del grupo c\'iclico $\chi(G_{L|K_0})=\co 1{G_{L|
K_0}}{{\ma Q}/{\ma Z}}$. Entonces, $\delta\mu$ es un generaodr de $\co
2{G_{L|K_0}}{\ma Z}$. Por el teorema de Tate, se tiene que
\[
{\underline{\ }}\Cup\delta\mu:\co 0{G_{L|K_0}}{\*L}\lra \co 2{G_{L|K_0}}{
\*L}
\]
es una biyecci\'on. Por tanto todo elemento $c\in\co 2{G_{L|K_0}}{\*L}$ 
tiene la forma $c=(a)\Cup \delta\mu$ con $(a)=a\N_{L|K_0}\*L\in \co 2{G_{
L|K_0}}{\*L}$, $a\in \*K_0$. De la Proposici\'on \ref{P17.6.107N} se tiene
\[
\inv_{L|K_0} c=\inv_{L|K_0} \big((a)\Cup \delta\mu\big)=\mu\big([a,L/K_0]\big).
\]
Por tanto necesitamos probar $[a,L/K_0]=\prod_{\pKK}(a_{\pK},L_{\pL}/(K_0
)_{\pK})=1$.

Como $L$ es ciclot\'omica $L\subseteq K_0(\zeta_n)$ para alguna $n$
(ambos casos). El automorfismo $[a,L/K_0]$ es precisamente la restricci\'on de $[a,K_0
(\zeta_n)/K_0]$ a $L$ (Teorema \ref{CCLT17.6.25}).

Por tanto, es suficiente probar que $[a,K_0(\zeta_n)/K_0]=1$ para $a\in
\*K_0$. Puesto que $K_0(\zeta_n)$ est\'a generado por ra\'ices de unidad
de \'ordenes una potencia de un n\'umero primo, es suficiente probar que
$[a,K_0(\zeta_n)/K_0]=1$ para estos generadores. Por tanto, podemos 
suponer que $n=l^m$ con $l$ un n\'umero primo.

Sea $\zeta$ una ra\'iz $l^m$-primitiva de la unidad. Si $l=2$ y $K_0={\ma Q}$,
podemos suponer $m\geq 2$. Las completaciones de $K_0(\zeta_n)/K_0$
son $(K_0)_{\pK}(\zeta)/(K_0)_{\pK}$ 
la cual es no ramificada para $l\neq p=\pK$
y totalmente ramificada para $l=p=\pK$ caso $K_0={\ma Q}$ y $l\neq p$ y no
ramificada siempre para los campos de funciones. Para $\pK$, el primo
infinito real, $(K_0)_{\pK}(\zeta)/(K_0)_{\pK}$ significa ${\ma C}/{\ma R}$.

Se debe probar que para $a\in\*K_0$, 
\[
[a,K_0(\zeta)/K_0]=\prod_{\pK}(a,(K_0)_{\pK}(\zeta)/(K_0)_{\pK})=1.
\]
Se tiene que $K_0(\zeta)/K_0$ es no ramificada excepto si $K_0={\ma Q}$,
$\pK=l$ o $\pK={\mc P}_{\infty}$ el valor absoluto usual de ${\ma Q}$.
Para cualquier otro $\pK$, por ser la extensi\'on no ramificada,
\begin{gather*}
(a,(K_0)_{\pK}(\zeta)/(K_0)_{\pK})=\varphi_{\pK}^{v_{\pK}(a)},
\intertext{donde $\varphi_{\pK}$ es el Frobenius y por tanto}
(a,(K_0)_{\pK}(\zeta)/(K_0)_{\pK})(\zeta)=\varphi_{\pK}^{v_{\pK}(a)}(\zeta).
\end{gather*}

El campo residual tiene $q_{\pK}$-elementos, donde 
\begin{gather*}
q_{\pK}=
\begin{cases} p&\text{si $K_0$ es num\'erico},\\
q^{\deg\pK}&\text{si $K_0$ es de funciones}.
\end{cases}
\intertext{Por tanto}
(a,(K_0)_{\pK}(\zeta)/(K_0)_{\pK})(\zeta)=\zeta^{q_{\pK}^{v_{\pK}(a)}}.
\intertext{Consideremos $K_0$ campo de funciones. Entonces}
(a,(K_0)_{\pK}(\zeta)/(K_0)_{\pK})(\zeta)=\zeta^{q^{v_{\pK}(a)\deg\pK}}.
\intertext{Se sigue}
\prod_{\pK}(a,(K_0)_{\pK}(\zeta)/(K_0)_{\pK})(\zeta)=\prod_{\pK}
\zeta^{q^{v_{\pK}(a)\deg\pK}}=\zeta^{q^s},
\intertext{donde}
s=\sum_{\pK} v_{\pK}(a)\deg \pK=\deg a=0.
\intertext{De esta forma}
[a,K_0(\zeta)/K_0](\zeta)=\zeta^{q^0}=\zeta.
\end{gather*}
Esto prueba que $[a,K_0(\zeta)/K_0]=1$ para $K_0$ un
campo de funciones.

Ahora sea $K_0={\ma Q}$. Para $p\neq l,{\mc P}_{\infty}$, se tiene
\[
(a,{\ma Q}_p
(\zeta)/{\ma Q}_p)=\varphi^{v_p(a)}\quad\text{y}\quad (a,{\ma Q}_p(\zeta)/{\ma Q}_p)
(\zeta)=\zeta^{\varphi^{v_p(a)}}.
\]

Usando los grupos de Lubin-Tate, Ejemplos \ref{CClaseE3.3.21'}
y \ref{CClaseEj3.2.5.30}, se tiene, para $p=l$, $\zeta=\zeta_{p^m}=
\zeta_{l^m}$, que
\[
(a,{\ma Q}_p(\zeta)/{\ma Q}_p)(\zeta)=\zeta^r,
\]
donde $r\equiv u^{-1}\bmod p^m$ y donde $a=up^n$, $u\in U_p$. 
Entonces $r\equiv u^{-1}\equiv a^{-1}p^{v_p(a)}\bmod p^m$.

Para $p={\mc P}_{\infty}$, el automorfismo $(a,{\ma C}/{\ma R})$ es bien
la identidad o bien conjugaci\'on compleja, dependiendo de si $a>0$
o $a<0$ respectivamente. Por tanto 
\begin{gather*}
(a,{\ma Q}_p(\zeta)/{\ma Q}_p)(\zeta)=\zeta^{\sg a}.
\intertext{donde $\sg (a)=\begin{cases}
1&\text{si $a>0$},\\ -1&\text{si $a<0$}. \end{cases}$\quad
Se sigue que}
[a,{\ma Q}(\zeta)/{\ma Q}] (\zeta)=\prod_p (a,{\ma Q}_p(\zeta)/{\ma Q}_p)
(\zeta)=\zeta^{\sg a\cdot \prod_{p\neq l}p^{v_p(a)}\cdot r},
\intertext{Ahora}
\sg a\cdot\prod_{p\neq l} p^{v_p(a)}\cdot r\equiv \sg a\cdot\prod_{p\neq l}
p^{v_p(a)}\cdot l^{v_l(a)}\cdot a^{-1}=\frac{1}{\prod_p|a|_p}=1\bmod l^m.
\end{gather*}

Se sigue que $[a,{\ma Q}(\zeta)/{\ma Q}](\zeta)=\zeta$ lo cual implica que
$[a,{\ma Q}(\zeta)/{\ma Q}]=1$.
$\fin$
\end{proof}

\begin{corolario}\label{C17.6.109N}
Sea $L/K$ una extensi\'on abeliana finita de campos globales y sea 
$\Omega=\sep K$. Entonces, si $a\in\*K$ es un id\`ele principal, se tiene
$[a,L/K]=1$ y $[a,\Omega/K]=1$. $\fin$
\end{corolario}

Con este resultado tenemos que $\co 2{G_{L|K}}{\*L}\subseteq \ker
\inv_{L|K}$ donde $\inv_{L|K}:\co 2{G_{L|K}}{J_L}\lra \big(\frac{1}{[L:K]}{\ma Z}
\big)/{\ma Z}$ y donde $L/K$ es una extensi\'on abeliana
finita de campos globales.
Todav\'ia queda pendiente averiguar exactamente cual es el n\'ucleo de
$\inv_{L|K}$ y si $\inv_{L|K}$ es suprayectiva.

Para el caso c\'iclico, se tiene

\begin{proposicion}\label{P17.6.110N}
Si $L/K$ es una extensi\'on c\'iclica finita de campos globales, entonces
\[
1\lra \co 2{G_{L|K}}{\*L}\xhookrightarrow[]{\ \infla\ }\co 2{G_{L|K}}{J_L}
\xrightarrow[]{\inv_{L|K}} \Big(\frac{1}{[L:K]}{\ma Z}\Big)/{\ma Z}\lra 0
\]
es una sucesi\'on exacta.

En particular $\inv_{L|K}$ es suprayectiva y $\ker\inv_{L|K}=\co 2{G_{
L|K}}{\*L}$.
\end{proposicion}

\begin{proof}
Para ver la suprayectividad de $\inv_{L|K}$, supongamos primero que
$[L:K]=l^r$ con $l$ un n\'umero primo. Se tiene que $\frac 1{[L:K]}+{\ma Z}$
genera $\Big(\frac 1{[L:K]}{\ma Z}\Big)/{\ma Z}$ por lo que es suficiente
hallar $c\in\co 2{G_{L|K}}{J_L}$ tal que $\inv_{L|K} c=\frac 1{[L:K]}+{\ma Z}$.

Como $L/K$ es una extensi\'on c\'iclica de orden $l^r$, existe $\pK_0\in K$
totalmente inerte en $L$. Por tanto, para ${\eu P}_0\in {\ma P}_L$ sobre
$\pK_0$, se tiene que $[L_{\pL_0}:K_{\pK_0}]=[L:K]$. Por el teorema de
reciprocidad local, existe $c_{\pK_0}\in \co 2{G_{L_{\pL_0}|K_{\pK_0}}}
{\*L_{\pL_0}}$ con $\inv_{L_{\pL_0}|K_{\pK_0}}c_{\pK_0}
=\frac 1{[L_{\pL_0}:K_{\pK_0}]}
+{\ma Z}=\frac 1{[L:K]}+{\ma Z}$.

Sea $c_{\pK}=1$ para todo $\pK\neq \pK_0$. De $\co 2{G_{L|K}}{J_L}
\cong \bigoplus_{\pK}\co 2{G_{L_{\pL}|K_{\pK}}}{\*L_{\pL}}$, consideremos
$c=(\ldots,1,c_{\pK_0},1,\ldots)\in \co 2{G_{L|K}}{J_L}$. Entonces
\[
\inv_{L|K}c=\sum_{\pK}\inv_{L_{\pL}|K_{\pK}} c_{\pK}=\inv_{
L_{\pL_0}|K_{\pK_0}} c_{\pK_0}=\frac 1{[L:K]}+{\ma Z}.
\]

Ahora sea $[L:K]=n=p_1^{r_1}\cdots p_s^{r_s}$ su descomposici\'on
en primos. Para cada $1\leq i\leq s$,
existe $K\subseteq L_i\subseteq L$, $[L_i:K]=p_i^{r_i}$. Sea
\[
\frac 1n=\frac{n_1}{p_1^{r_1}}+\cdots+\frac {n_s}{p_s^{r_s}}.
\]
Existe $c_i\in \co 2{G_{L_i|K}}{J_{L_i}}$ con $\inv_{L_i|K} c_i=\inv_{L|K} c_i
=\frac {n_i}{p_i^{r_i}}+{\ma Z}$.

Sea $c=c_1\cdots c_s\in \co 2{G_{L|K}}{J_L}$. Entonces
\[
\inv_{L|K} c=\sum_{i=1}^s \inv_{L|K} c_i=\sum_{i=1}^s \frac {n_i}{p_i^{r_i}}+
{\ma Z}=\frac 1n + {\ma Z}.
\]
Se sigue que $\inv_{L|K}$ es suprayectiva para toda extensi\'on c\'iclica.

Puesto que $\co 2{G_{L|K}}{\*L}\subseteq \ker \inv_{L|K}$, para probar la
igualdad, basta probar que $\co 2{G_{L|K}}{J_L}/\co 2{G_{L|K}}{\*L}$ es 
de orden menor o igual a $\frac 1{[L:K]}=\big|\big(\frac {1}{[L:K]} {\ma Z}\big)
/{\ma Z}\big|$.

Ahora bien, de la sucesi\'on exacta $1\lra \*L\lra J_L\lra C_L\lra 1$ y puesto
que $\co 1{G_{L|K}}{C_L}=1$, obtenemos en cohomolog\'ia la sucesi\'on exacta
\begin{gather*}
1\lra
\co 2{G_{L|K}}{\*L}\lra \co 2{G_{L|K}}{J_L}\lra \co 2{G_{L|K}}{C_L}.
\intertext{Por tanto}
\Big|\frac{\co 2{G_{L|K}}{J_L}}{\co 2{G_{L|K}}{\*L}}\Big|\Big|
\big|\co 2{G_{L|K}}{C_L}\big| \Big| [L:K]. \tag*{$\fin$}
\end{gather*}
\end{proof}

\begin{observacion}\label{O17.6.111N}
Ser\'ia afortunado que tuvi\'esemos que $\inv_{L|K}$ fuese suprayectiva en
general, pero esto es falso. Para que cada elemento de $\big(\frac{1}{[L:K]}
{\ma Z}\big)/{\ma Z}$ est\'e en la imagen del mapeo $\inv_{L|K}$, es
necesario agregar a $L$ una extensi\'on c\'iclica.
\end{observacion}

Consideramos $\co 2{G_{\Omega|K}}{J_{\Omega}}=\bigcup_L\co 2{G_{L|K}}
{J_L}$, donde $L$ recorre a las extensiones finitas de Galois de $K$.
Para extensiones de Galois, $K\subseteq L\subseteq M$, se tiene
$\co 2{G_{L|K}}{J_L}\underbracket[0pt]{\subseteq}_{\substack{\uparrow\\
\infla}} \co 2{G_{M|K}}{J_M}$.

Ahora, el mapeo invariante $\inv$ se puede extender de $\co 2{G_{L|K}}
{J_L}$ a $\co 2{G_{M|K}}{J_M}$ y por ende se obtiene el homomorfismo
\[
\inv_K:\co 2{G_{\Omega|K}}{J_{\Omega}}\lra {\ma Q}/{\ma Z}
\]
y donde $\inv_K|_L=\inv_{L|K}:\co 2{G_{L|K}}{J_L}\lra \big(\frac{1}{[L:K]}
{\ma Z}\big)/{\ma Z}\subseteq {\ma Q}/{\ma Z}$.

Para cada $m\in{\ma N}$, existe $L_m/K$ c\'iclica de grado $m$ y tal que
$\big(\frac 1m {\ma Z}\big)/{\ma Z}\subseteq \im(\inv_K)$, 
lo cual implica que ${\ma Q}/{\ma Z}\subseteq \im(\inv_K)$.

\begin{teorema}\label{T17.6.112N}
El homomorfismo $\inv_K:\co 2{G_{\Omega|K}}{J_{\Omega}}\lra {\ma Q}/
{\ma Z}$ es suprayectivo.
$\fin$
\end{teorema}

Ahora bien, el grupo de Brauer de $K$, $\Br K$, satisface que
$\Br K=\bigcup_L\co 2{G_{L|K}}{\*L}$ donde $L$ recorre las
extensiones finitas de Galois de $K$. Para cada primo $\pK$ de $K$, 
seleccionamos un valor absoluto fijo en $\Omega$, 
entonces este valor absoluto
determina a su vez un primo $\pL$ de cada extensi\'on de Galois $L/K$
sobre $\pK$ y se tiene que el grupo de Brauer del campo local
$K_{\pK}$ es tal que
\begin{gather*}
\Br {K_{\pK}}=\bigcup_L \co 2{G_{L_{\pL}|K_{\pK}}}{\*L_{\pL}}.
\intertext{Se tiene}
\begin{align*}
\co 2{G_{L|K}}{\*L} & \lra \co 2{G_{L|K}}{J_L}\cong \bigoplus_{\pK}
\co 2{G_{L_{\pL}|K_{\pK}}}{\*L_{\pL}}\\
&\xrightarrow[]{\inv_{L|K}}
\Big(\frac 1{[L:K]}{\ma Z}\Big)/{\ma Z}
\end{align*}
\intertext{y tomando el l\'imite directo, en este caso la uni\'on, obtenemos el
homomorfismo can\'onico}
\Br K\lra \co 2{G_{\Omega|K}}{J_{\Omega}}\cong \bigoplus_{\pK} \Br 
{K_{\pK}}\xrightarrow[]{\inv_{L|K}} {\ma Q}/{\ma Z},
\end{gather*}
donde $\inv_K=\sum_{\pK}\inv_{K_{\pK}}$, $\inv_{K_{\pK}}:\Br {K_{\pK}}
\lra {\ma Q}/{\ma Z}$.

\begin{teorema}[Brauer-Hasse-Noether\index{teorema de
Brauer-Hasse-Noether}]\label{T17.6.113N}
Sea $K$ un campo global. Se tiene
\las
\item Sea $\Br K\lra \Br {K_{\pK}}$, $\alpha\longmapsto \alpha_{\pK}$, el
mapeo can\'onico para cada lugar $\pK$ de $K$. Entonces $\alpha_{\pK}
=0$ para casi toda $\pK$ y en particular $\Br K\lra \prod\limits_{\pK} \Br {K_{\pK}}$,
$\alpha\longmapsto (\alpha_{\pK})_{\pK}$ pertenece a $\bigoplus\limits_{\pK}
\Br {K_{\pK}}$.

\s

\item {\bf Teorema principal de Hasse en la teor\'ia de \'algebras}

La sucesi\'on $0\lra \Br K\lra \bigoplus\limits_{\pK} \Br {K_{\pK}}\stackrel{\xi}{\lra}
{\ma Q}/{\ma Z}\lra 0$ es exacta, donde $\xi\big((\alpha_{\pK})_{\pK}\big)=
\sum_{\pK}\inv_{\pK}(\alpha_{\pK})$.

Notemos que $\bigoplus\limits_{\pK}\Br{K_{\pK}}=\bigoplus\limits_{{\pK}\text{\ real}}
\frac{\big(\frac 12 {\ma Z}\big)}{\ma Z}\bigoplus\bigoplus\limits_{\pK\nmid \infty}
{\ma Q}/{\ma Z}$.
\end{list}
\end{teorema}

\begin{proof}
Los grupos $\Br K$, $\bigoplus_{\pK}\Br{K_{\pK}}\cong \co 2{G_{\Omega|K}}
{J_{\Omega}}$ y ${\ma Q}/{\ma Z}$ son las respectivas uniones de $\co 2{G_{
L|K}}{\*L}$, $\bigoplus_{\pK}\co 2{G_{L_{\pL}|L_{\pK}}}{\*L_{\pL}}\cong \co 2{G_{
L|K}}{J_L}$ y $\big(\frac {1}{[L:K]}{\ma Z}\big)/{\ma Z}$ y donde $L/K$ recorre las
extensiones c\'iclicas finitas de $K$. Para cada $L/K$ c\'iclica finita,
se tiene que la sucesi\'on 
\[
1\lra \co 2{G_{L|K}}{\*L}\lra \co 2{G_{L|K}}{J_L}\xrightarrow[]{\inv_{L|K}}
\Big(\frac {1}{[L:K]}{\ma Z}\Big)/{\ma Z}\lra 0
\]
es exacta. El resultado se sigue tomando los l\'imites directos (en nuestro caso
uniones).
$\fin$
\end{proof}

\begin{corolario}[Ley de reciprocidad de Hasse\index{ley de reciprocidad de 
Hasse}\index{Hasse!ley de reciprocidad de $\sim$}]\label{C17.6.114N}
Si $\Br K\stackrel {i}{\lra}\bigoplus\limits_{\pK} \Br {K_{\pK}}$ y $\bigoplus\limits_{
\pK} \Br {K_{\pK}}\stackrel {\xi}{\lra} {\ma Q}/{\ma Z}$, entonces
$\psi=\xi\circ i=0$, esto es, $\psi:\Br K\lra {\ma Q}/{\ma Z}$, $\psi(\alpha)=
\sum\limits_{\pK}\inv_{K_{\pK}} \alpha =0$.
$\fin$
\end{corolario}

\subsection{Ley de reciprocidad}\label{S17.6.9N}

Consideremos la sucesi\'on exacta $1\lra \*L\lra J_L\stackrel{\pi}{\lra} C_L
\lra 1$.

\begin{proposicion}\label{P17.6.115N}
Sea $L/K$ una extensi\'on c\'iclica finita de campos globales. Entonces 
$\tilde\pi: \co 2{G_{L|K}}{J_L}\lra \co 2{G_{L|K}}{C_L}$ es
suprayectiva.
\end{proposicion}

\begin{proof}
De la sucesi\'on exacta, se obtiene en cohomolog\'ia
\begin{gather*}
\co 1{G_{L|K}}{C_L}=1\lra \co 2{G_{L|K}}{\*L}\lra \co 2{G_{L|K}}{J_L}
\stackrel{\tilde\pi}{\lra}\\
\stackrel{\tilde\pi}{\lra}
\co 2{G_{L|K}}{C_L}\lra \co 3{G_{L|K}}{\*L}\cong
\co 1{G_{L|K}}{\*L}=1.
\end{gather*}
Por tanto $\tilde\pi$ es suprayectiva.
$\fin$
\end{proof}

Nos gustar\'ia definir en general para $\bar c\in \co 2{G_{L|K}}{C_L}$
y $c\in \co 2{G_{L|K}}{J_L}$ donde $\bar c=\tilde\pi c$,
\[
\inv_{L|K} \bar c=\inv_{L|K} c\in \Big(\frac {1}{[L:K]} {\ma Z}\Big)/{\ma Z}.
\]

Tenemos que la definici\'on es independiente de la preimagen $c$ pues
2 preim\'agenes cualesquiera 
difieren \'unicamente de un elemento en $\co 2{G_{L|K}}
{\*L}$ pero cualquier elemento de este \'ultimo grupo tiene invariante $0$.
Notemos que en general que se tiene
\begin{gather*}
1\lra \co 2{G_{L|K}}{\*L}\lra \co2{G_{L|K}}{J_L}\stackrel{\tilde\pi}{\lra}
\co 2{G_{L|K}}{C_L}\lra\\
\lra \co 3{G_{L|K}}{\*L}\lra \co 3{G_{L|K}}{J_L}=\{1\}.
\end{gather*}
Se sigue que $\tilde\pi$ es suprayectiva si y solamente si $\co 3{G_{L|K}}
{\*L}=\{1\}$ pero esto no se cumple en general. De hecho, si $L/K$ es una
extensi\'on de Galois de grado $n$ y $m$ es el m\'inimo com\'un m\'ultiplo
de todos los grados locales $n_{\pK}=[L_{\pL}:K_{\pK}]$, entonces $\co 3
{G_{L|K}}{\*L}$ es un grupo c\'iclico de orden $n/m$ (\cite[Chap. 7,
Theorem 12]{ArtTat61}).

Por tanto, no todo elemento $\bar c\in \co 2{G_{L|K}}{C_L}$ proviene de
un elemento $c\in \co 2{G_{L|K}}{J_L}$ y no podemos definir $\inv_{L|K}
\bar c$ como nos hubiera gustado.

Para poder definir un mapeo invariante para una extensi\'on finita de
Galois $L/K$, se procede de manera parecida para el caso de los
grupos $\co 2{G_{L|K}}{J_L}$. Se tiene que el homomorfismo
\[
\co 2{G_{L|K}}{J_L}\stackrel{\tilde\pi}{\lra} \co 2{G_{L|K}}{C_L}
\]
conmuta con los homomorfismos $\infla$ y $\res$, esto es, si $K\subseteq
L\subseteq M$ son extensiones finitas de Galois, entonces se tiene
\[
\tilde\pi\circ\infla_M=\infla_M\circ \tilde\pi\quad\text{y}\quad \tilde\pi\circ
\res_L=\res_L\circ\tilde\pi
\]
y en la segunda f\'ormula, \'unicamente se necesita que $M/K$ sea de 
Galois.

\begin{notacion}\label{N17.6.116N}
Sea $L/K$ una extensi\'on finita de Galois de campos globales. 
Entonces denotamos $ \cob qLK=\co q{G_{L|K}}{C_L}$. El grupo $C_L$
jugar\'a el papel global de $\*L$ en el caso local.
\end{notacion}

Puesto que $\cob 1LK=\{1\}$, las extensiones $L/K$ forman una familia
como en el caso local con $C_L$ en lugar de $\*L$. Como de costumbre,
cuando $K\subseteq L\subseteq M$, la inflaci\'on es inyectiva: $\cob 2LK
\stackrel{\infla}{\lra} \cob 2MK$ y supondremos que $\cob 2LK\subseteq
\cob 2MK$. Si $\Omega=\sep K$ tenemos el l\'imite directo
\[
\cob 2{\Omega}K=\lim_{\substack{\lra\\ L}}\cob 2LK,
\]
donde $L$ var\'ia sobre las extensiones finitas de Galois de $K$ y 
usaremos que $\cob 2LK\subseteq \cob 2{\Omega}K$ v\'ia el mapeo
inflaci\'on. De esta forma se tiene
\[
\cob 2{\Omega}K=\bigcup_L \cob 2LK,
\]
y si $K\subseteq L\subseteq M$ son extensiones de Galois, entonces
$\cob 2LK\subseteq \cob 2MK\subseteq \cob 2{\Omega}K$.

De manera similar como en el caso local, tenemos:

\begin{teorema}\label{T17.6.117N}
Sea $L/K$ una extensi\'on finita de Galois de campos globales de grado $n$.
Sea $N/K$ una extensi\'on c\'iclica de grado $n$. Entonces
\[
\cob 2NK=\cob 2LK\subseteq \cob 2{\Omega}K.
\]
\end{teorema}

\begin{proof}
Sea $M:=NL$. Como $N/K$ es c\'iclica, $M/L$ es c\'iclica de orden un
divisor de $n$.
\[
\xymatrix{
N\ar@{-}[r]\ar@{-}[d]_n&M=NL\ar@{-}[d]\\K\ar@{-}[r]&L}
\]
Sea $\bar c\in \cob 2NK\subseteq \cob 2MK$. De la sucesi\'on exacta
\[
1\lra \cob 2LK\stackrel{\infla}{\lra}\cob 2MK\stackrel{\res}{\lra}\co 2ML,
\]
se tiene que un elemento $\bar c\in \cob 2MK$ es un elemento de
$\cob 2LK$ si y solamente si $\res_L\bar c=1$. Puesto que $N/K$
es c\'iclica, el homomorfismo $\co 2{G_{N|K}}{J_L}\stackrel{\tilde \pi}
{\lra}\cob 2NK$ es suprayectivo. Por tanto, existe $c\in \co 2{G_{N|K}}
{J_L}\subseteq \co 2{G_{M|K}}{J_M}$ con $\tilde \pi c=\bar c$. Como
$\tilde \pi$ conmuta con inflaci\'on, la cual es de hecho un inclusi\'on,
y tambi\'en con restricci\'on, se tiene
\[
\res_L \bar c=\res_L (\tilde\pi c)=\tilde \pi(\res_L c).
\]
Por tanto, $\res_L\bar c=1\iff \res_L c\in \ker \tilde\pi=\co 2{G_{M|L}}{\*M}$.
Puesto que $M/L$ es c\'iclica, $\res_L c\in \co 2{G_{M|L}}{\*M}\iff
\inv_{M|L}(\res_L c)=0$. Ahora bien,
\[
\inv_{M|L}(\res_L c)=[L:K]\cdot\inv_{M|K} c=[N:K]\cdot \inv_{N|K} c=0.
\]
Se sigue que $\cob 2NK\subseteq \co 2LK$.

Puesto que $\cob 1NK=\{1\}$ y $\co 3{G_{N|K}}{\*N}\cong\co 1{G_{N|K}}
{\*N}=\{1\}$, se tiene la sucesi\'on exacta en cohomolog\'ia
\[
1\lra \co 2{G_{N|K}}{\*N}\lra \co 2{G_{N|K}}{J_N}\lra \cob 2NK\lra 1,
\]
y donde $|\cob 2NK|=[N:K]=[L:K]$ pues $N/K$ es una extensi\'on c\'iclica
(Teorema \ref{T17.6.102N}). Por otro lado, $|\cob 2LK|$ divide a $[L:K]$
(Teorema \ref{T17.6.96N}). Por tanto $\cob 2NK=\cob 2LK$.
$\fin$
\end{proof}

\begin{corolario}\label{C17.6.118N}
Se tiene $\cob 2{\Omega}K=\bigcup\limits_{L/K\text{\ c\'iclica}} \cob 2LK$.
\end{corolario}

\begin{proof}
Es consecuencia de Teorema \ref{T17.6.117N} y de que para toda $n\in
{\ma N}$, existe una extensi\'on c\'iclica $L/K$ de grado $n$.
$\fin$
\end{proof}

Consideremos $K\subseteq L\subseteq M$ extensiones finitas de Galois.
Ahora bien, al conmutar $\tilde\pi:\co 2{G_{L|K}}{J_L}\lra \cob 2LK$ con
restricci\'on y con inflaci\'on, usando esto \'ultimo, $\tilde\pi$ se extiende
de manera can\'onica a $\tilde\pi: \co 2{G_{M|K}}{J_M}\lra \cob 2MK$ y
de esta forma obtenemos un homomorfismo
\[
\tilde\pi:\co 2{G_{\Omega|K}}{J_{\Omega}}\lra \cob 2{\Omega}K,
\]
cuya restricci\'on a los grupos $\co 2{G_{L|K}}{J_L}$ son los homorfismos
$\tilde\pi$ originales: $\tilde \pi:\co 2{G_{L|K}}{J_L}\lra \cob 2LK$. Estos
\'ultimos no son suprayectivos en general pero, 
se tiene de cualquier forma el siguiente resultado:

\begin{teorema}\label{T17.6.119}
El homomorfismo
\[
\tilde\pi: \co 2{G_{\Omega|K}}{J_{\Omega}}\lra \cob 2{\Omega}K
\]
es suprayectivo.
\end{teorema}

\begin{proof}
Si $\bar c\in \cob 2{\Omega}K$, entonces existe una extensi\'on 
$L/K$ c\'iclica finita tal que $\bar c\in \cob 2LK$. 
Puesto que para extensiones c\'iclicas finitas
\[
\tilde \pi:\co 2{G_{L|K}}{J_L}\lra \cob 2LK
\]
es suprayectiva, $\bar c=\tilde \pi c$ para alguna $c\in \co 2{G_{L|K}}
{J_L}\subseteq \co 2{G_{\Omega|K}}{J_{\Omega}}$.
$\fin$
\end{proof}

Con este resultado, podemos obtener clases invariantes para los
elementos de $\co 2{\Omega}K=\bigcup_L \cob 2LK$ a partir del mapeo
invariante de las clases de cohomolog\'ia del grupo de id\`eles. De
hecho, del homomorfismo
\[
\inv_K:\co 2{G_{\Omega|K}}{J_{\Omega}}\lra {\ma Q}/{\ma Z},
\]
el cual es suprayectivo, obtenemos:

\begin{definicion}\label{D17.6.120N}
Si $\bar c\in \cob 2{\Omega}K$ y $\bar c=\tilde\pi c$, $c\in \co 2{G_{\Omega|K}}
{J_{\Omega}}$, se define $\inv_K \bar c:=\inv_K c\in {\ma Q}/{\ma Z}$.
\end{definicion}

Veamos que la definici\'on es independiente de la preimagen $c$ seleccionada.
Si $d$ es otra preimagen de $\bar c$ en $\co 2{G_{\Omega|K}}{J_{\Omega}}$,
$\bar c=\tilde \pi d$, entonces podemos escoger $L/K$ una extensi\'on de
Galois con $c,d\in \co 2{G_{L|K}}{J_L}\subseteq \co 2{G_{\Omega|K}}
{J_{\Omega}}$ y podemos agrandar $L$ en caso de ser necesario de tal
forma que $\bar c\in \cob 2LK$. Ahora bien, $\bar c=\tilde \pi c=\tilde \pi d$
por lo que $c$ y $d$ difieren por elemento de $\ker\tilde\pi$, $\tilde \pi: \co
2{G_{L|K}}{J_L}\lra \cob 2LK$ y por tanto de un elemento de $\co 2{G_{
L|K}}{\*L}=\ker \tilde\pi$, pero este elemento tiene invariante $0$. Por tanto
$\inv_K \bar c$ est\'a bien definido.

Con esta definici\'on de $\inv \bar c$, obtenemos un epimorfismo 
\[
\inv_K: \cob 2{\Omega}K\lra {\ma Q}/{\ma Z}.
\]

La restricci\'on de $\inv_K$ a $L$, 
$\inv_K|_L$, donde $L/K$ es una extensi\'on finita de
Galois, $\inv_K|_L=\inv_{L|K}: \cob 2LK\lra \big(\frac{1}{[L:K]}{\ma Z}\big)/
{\ma Z}$ puesto que los \'ordenes de los elementos de $\cob 2LK$ dividen
a $[L:K]$, y por lo tanto son mapeados al \'unico subgrupo $\big(\frac{1}{
[L:K]}{\ma Z}\big)/{\ma Z}$ de ${\ma Q}/{\ma Z}$ de orden $[L:K]$.

\begin{teorema}\label{T17.6.121N}
Cuando $\bar c=\tilde\pi c$, $\bar c\in \cob 2LK$, $c\in \co
2{G_{L|K}}{J_L}$, entonces $\inv_{L|K} \bar c=\inv_{L|K} c$.

Para el caso general, tenemos:
Sean $L/K$ una extensi\'on finita de Galois de campos globales y $\Omega
=\sep K$. Entonces los mapeos invariantes
\begin{align*}
\inv_K&: \cob 2{\Omega}K\lra {\ma Q}/{\ma Z},\\
\inv_{L|K}&: \cob 2LK\lra \Big(\frac 1{[L:K]}{\ma Z}\Big)/{\ma Z},
\end{align*}
son isomorfismos.
\end{teorema}

\begin{proof}
Se tiene que es suficiente verificar que $\inv_{L|K}$ es biyectiva pues como
$\cob 2LK=\cob 2NK$ para $N/K$ c\'iclica con $[N:K]=[L:K]$, se tiene que si
$\{N_n\}_{n=1}^{\infty}$ es una colecci\'on tal que $N_n/K$ es c\'iclica de grado
$n$, $\cob 2{\Omega}K=\bigcup_{n=1}^{\infty}\cob 2{N_n}K$ y $\inv_K\big(
\co 2{N_n}K\big)=\big(\frac 1n {\ma Z}\big)/{\ma Z}$, ${\ma Q}/{\ma Z}=
\bigcup_{n=1}^{\infty}\big(\frac 1n {\ma Z}\big)/{\ma Z}$.

Sea $L/K$ dada y $N/K$ c\'iclica con $[N:K]=[L:K]$. Por tanto $\cob 2LK=
\cob 2NK$. Si $\alpha\in \big(\frac 1{[L:K]}{\ma Z}\big)/{\ma Z}$, existe
$c\in \co 2{G_{N|K}}{J_N}$ con $\inv_{N|K} c=\alpha$. Sea $\bar c=
\tilde \pi c\in \cob 2NK=\cob 2LK$ Entonces $\inv_{L|K}\bar c=\inv_{N|K}
\bar c=\inv_{N|K}c=\alpha$. Por tanto $\inv_{L|K}$ es suprayectiva.
Ahora bien, $|\cob 2LK|\big|[L:K]=\big|\big(\frac{1}{[L:k]}{\ma Z}\big)/{\ma Z}
\big|$. Por tanto $\inv_{L|K}$ es biyectiva.
$\fin$
\end{proof}

Para enunciar el resultado principal de la teor\'ia global de campos de
clase, damos una definici\'on que pudimos haber dado antes de desarrollar
la teor\'ia local de campos de clase, pero, en ese punto, por razones de
claridad, preferimos dar directamente los resultados para campos locales.

\subsubsection{Formaci\'on de clases}

Sea $G$ un grupo profinito, es decir, $G$ es compacto con la topolog\'ia
de los subgrupos normales (topolog\'ia de Krull), es decir, una base de
vecindades abiertas de la identidad $1$ de $G$ consiste de los subgrupos
de $G$ que son normales y de \'indice finito o, equivalentemente, 
$G$ es compacto, Hausdorff y totalmente disconexo.

Sea $\{G_K\mid K\in X\}$ la familia de todos los subgrupos cerrados de
$G$. Los subgrupos abiertos son los subgrupos cerrados de \'indice finito.

Enumeramos cada uno de los subgrupos cerrados $G_K$ con el 
\'indice $K$ que llamaremos {\em campo} y formalmente nos
referimos a $K$ como el ``{\em campo fijo de $G_K$}''. El campo 
$K_0$, con $G_{K_0}=G$ se llama el {\em campo base} y denotamos
por $\Omega=\sep K$ el campo tal que $G_{\Omega}=\{1\}$.

Escribimos formalmene $K\subseteq L$ o $L/K$ si $G_L\subseteq
G_K$ y al par $L/K$ lo llamamos {\em extensi\'on de campos}. La
extensi\'on $L/K$ se llama {\em finita} si $G_L$ es abierto en
$G_K$ (es de \'indice finito) y ponemos $[L:K]:=[G_K:G_L]$ y lo
llamamos el {\em grado de la extensi\'on}. Si $G_L$ es normal en
$G_K$, la extensi\'on $L/K$ se llama {\em normal} o {\em de
Galois} . En este caso se define {\em el grupo de Galois} 
de $L/K$ por $\Gal(L/K)=G_{L|K}:=G_K/G_L$.

Si $K\subseteq L\subseteq M$ son dos extensiones de Galois, se
define la {\em restricci\'on a $L$} de $\sigma\in \Gal(M/K)$ por
\[
\sigma|_L=\sigma\bmod \Gal(M/L)\in \Gal(L/K).
\]
La extensi\'on se llama {\em c\'iclica, abeliana, nilpotente, soluble,
simple}, etc. si $L/K$ es de Galois y $\Gal(L/K)$ es c\'iclica, abeliana,
nilpotente, soluble, simple, etc.

Se define $K=\bigcap_i K_i$ (intersecci\'on) si $G_K$ es topol\'ogicamente
generado por los subgrupos $G_{K_i}$ y $K=\prod_i K_i$ (composici\'on)
si $G_K=\bigcap_i G_{K_i}$. Si $G_{L'}=\sigma^{-1} G_i \sigma$, 
$\sigma\in G$, escribimos $L'=L^{\sigma}=\sigma(L)$.

De esta forma, de cada grupo profinito $G$ obtenemos una teor\'ia
de Galois formal.

\begin{definicion}\label{D17.6.122N}
Sean $G$ un grupo profinito y $A$ un $G$-m\'odulo con $A$ 
dotado de la
topolog\'ia discreta y tal que satisface cualquiera de las tres
siguientes condiciones equivalentes:
\lasa
\item La acci\'on $G\times A\lra A$, $(\sigma,a)\longmapsto
\sigma a$ es continua.
\item Para cada $a\in A$, el estabilizador $\{\sigma\in G\mid \sigma
a=a\}$ es abierto en $G$.
\item $A=\bigcup_U A^U$ donde $U$ recorre todos los subgrupos
abiertos de $G$.
\end{list}

Entonces el par $(G,A)$ se llama una {\em formaci\'on\index{formaci\'on}}.
$A$ puede tener otra topolog\'ia que no sea la discreta, pero para (a)
se considerar\'a $A$ con la topolog\'ia discreta.
\end{definicion}

\begin{ejemplo}\label{E17.6.123N}
Si $G$ es el grupo de Galois de la extensi\'on de campos $E/F$, 
entonces $G$ act\'ua en $\*E$ y el par $(G,\*E)$ es una formaci\'on.
Aqu\'i, para $F\subseteq N\subseteq E$, $G_N=\Gal(E/N)$ (si $G$
es finito, $G$ tiene la topolog\'ia discreta).
\end{ejemplo}

Sea $(G,A)$ una formaci\'on, $A$ escrito multiplicativamente. Se
define $A_K=A^{G_K}=\{a\in A\mid \sigma a=a\text{\ para toda $\sigma
\in G_K$}\}$.

\begin{ejemplo}\label{E17.6.124N}
Sean $K_0$ un campo local y $\Omega= \sep K_0$. Entonces $A_K=
\* K$.
\end{ejemplo}

\begin{ejemplo}\label{E17.6.125N}
Sean $K_0$ un campo global y $\Omega=\sep K_0$. Entonces
$A_K=C_K$.
\end{ejemplo}

\bigskip

Dada una formaci\'on $(G,A)$ consideramos para una extensi\'on de
Galois $L/K$ los grupos de cohomolog\'ia del $G_{L|K}$-m\'odulo 
$A_L$.

\begin{definicion}\label{D17.6.126N}
Una formaci\'on $(G,A)$ se llama una {\em formaci\'on de 
clases\index{formaci\'on de clases}} si satisface las siguientes condiciones:
\las
\item[\bf{Axioma I}:] $\cob 1LK=\co 1{G_{L|K}}{A_L}=\{1\}$ para cada extensi\'on
finita de Galois $L/K$ (formaci\'on de campos).

\item[\bf{Axioma II:}] Para toda extensi\'on finita de Galois $L/K$, existe un
isomorfismo 
\[
\inv_{L|K}:\cob 2LK=\co 2{G_{L|K}}{A_L}\lra
\big(\frac{1}{[L:K]}{\ma Z}\big)/{\ma Z}
\]
llamado el {\em mapeo
invariante\index{mapeo invariante}\index{invariante!mapeo}} con las
siguientes propiedades:
\lasa
\item Si $K\subseteq L\subseteq M$ es una torre de extensiones
finitas de Galois, entonces
\[
\inv_{L|K}=\inv_{M|K}|_{\cob 2LK}.
\]

\item Si $K\subseteq L\subseteq M$ es una torre 
de extensiones con $M/K$ finita
de Galois, entonces
\begin{gather*}
\inv_{M|L}\circ \res_L=[L:K]\cdot \inv_{M|K}\\
\xymatrix{
\cob 2MK\ar@{->}[rr]^{\inv_{M|K}}\ar@{->}[d]_{\res_L}&&
\Big(\frac {1}{[M:K]}{\ma Z}\Big)/{\ma Z}\ar@{->}[d]^{\cdot [L:K]}\\
\cob 2ML\ar@{->}[rr]_{\inv_{M|L}}&&\Big(\frac {1}{[M:L]}{\ma Z}\Big)/{\ma Z}
}
\end{gather*}
\end{list}
\end{list}
\end{definicion}

Si $(G,A)$ es una formaci\'on de clases, se tienen

\begin{teorema}[Principal]\label{T17.6.127N}
Sea $L/K$ una extensi\'on finita de Galois. Entonces el mapeo $u_{L|K}
\Cup {\underline{\ }}:\co q{G_{L|K}}{\ma Z}\lra \cob {{q+2}}LK$ con
$u_{L|K}\in \cob 2LK$ la {\em clase fundamental\index{clase fundamental}}
($\inv_{L|K}(u_{L|K})=\frac 1{[L:K]}+{\ma Z}\in \big(\frac 1{[L:K]}{\ma Z}\big)
/{\ma Z}\subseteq {\ma Q}/{\ma Z}$) es un isomorfismo para toda
$q\in{\ma Z}$.
$\fin$
\end{teorema}

\begin{teorema}[Ley general de reciprocidad\index{ley general de
reciprocidad}\index{reciprocidad!ley general de $\sim$}]\label{T17.6.128N}
Sea $L/K$ una extensi\'on finita de Galois, $u_{L|K}\Cup {\underline{\ }}:
\co {{-2}}{G_{L|K}}{\ma Z}\lra \cob 0LK$ nos da un isomorfismo can\'onico
(mapeo de Nakayama):
\begin{gather*}
\theta_{L|K}:\abe G_{L|K}\lra A_K/\N_{L/K} A_L.
\tag*{$\fin$}
\end{gather*}
\end{teorema}

Se obtienen todos los resultados derivados de cohomolog\'ia que obtuvimos
para los campos locales pues si $K_0$ es un campo local, $\Omega=
\sep K_0$, $G=\Gal(\Omega/K_0)$ y $A=\*{\Omega}$, entonces $(G,A)$
es una formaci\'on de clases.

Para obtener estos resultados para cualquier otra formaci\'on de clases,
las demostraciones se obtienen sustituyendo $\*\Omega$ por el otro $A$
y $\Gal(\Omega/K_0)$ por el otro $G$.

\s

Lo que queremos obtener es que si $K$ es un campo global y $\Omega
=\sep K$, $G=\Gal(\Omega/K)$, entonces $(G, C_{\Omega})$, $C_{
\Omega}$ el grupo de clases de id\`eles de $\Omega$: $C_{\Omega}
=\bigcup_L C_L$, donde $L$ recorre todas las extensiones finitas y 
separables de $K$, es una formaci\'on de clases.

Se tiene que $(G,C_{\Omega})$ es una formaci\'on.

\begin{teorema}\label{T17.6.129N}
La formaci\'on $(G,C_{\Omega})$ es una formaci\'on de clases con
respecto al mapeo invariante $\inv_K:\cob 2{\Omega}K\lra {\ma Q}/
{\ma Z}$ definido en la Definici\'on {\rm{\ref{D17.6.120N}}}.
\end{teorema}

\begin{proof}
Debemos verificar los Axiomas I y II de la Definici\'on \ref{D17.6.126N}:

\s

\noindent
{\underline{\bf{Axioma I:}}} $\cob 1LK=\co 1{G_{L|K}}{C_L}=\{1\}$ para toda
extensi\'on de Galois finita (Teorema \ref{T17.6.98N}).

\s

\noindent
{\underline{\bf{Axioma II:}}} Para cada extensi\'on de Galois finita tenemos el
isomorfismo $\inv_{L|K}: \cob 2LK\lra \big(\frac 1{[L:K]}{\ma Z}\big)/{\ma Z}$
(Teorema \ref{T17.6.121N}).

\s

\noindent
{\underline{\bf {(a):}}} Si $K\subseteq L\subset M$ son dos extensiones finitas de
Galois de $K$ y $\bar c\in \cob 2LK$, entonces $\bar c\in \cob 2MK$ y
$\inv_{M|K} \bar c=\inv_{L|K} \bar c$ pues $\inv_{M|K}$ y $\inv_{L|K}$ est\'an
definidas por restricci\'on de $\inv_K$ a $\co 2MK$ y $\co 2LK\subseteq
\cob 2MK$, respectivamente.

\s

\noindent
{\underline{\bf{(b):}}} Sean $K\subseteq L\subseteq M$ dos extensiones
de $K$ con $M/K$ extensi\'on finita de Galois. Si $\bar c\in \cob 2MK$, entonces
$\res_L\bar c\in \cob 2ML$. Para la demostraci\'on de la f\'ormula
\[
\inv_{M|K}(\res_L \bar c)=[L:K] \inv_{M|L} \bar c
\]
usamos la f\'ormula an\'aloga para el grupo de id\`eles: existe $c\in
\co 2{G_{\Omega|K}}{J_{\Omega}}$ con $\tilde \pi c=\bar c$ donde
podemos suponer que hay una extensi\'on finita de Galois $N/K$ que
contiene a $M$, $K\subseteq L\subseteq M\subseteq N$, tal que 
$c\in \co 2{G_{N|K}}{J_N}$. De la f\'ormulas de la Proposici\'on
\ref{P17.6.104N} y usando al mapeo inflaci\'on como inclusi\'on 
tenemos
\begin{align*}
\inv_{M|L}(\res_L c)&=\inv_{N|L}(\res_L \tilde\pi c)=\inv_{N|L}(
\tilde \pi \res_L c)=\inv_{N|L}(\res_L c)\\
&=[L:K]\cdot \inv_{N|K} c=[L:K]\cdot \inv_{N|K}\tilde\pi c\\
&=[L:K]
\cdot \inv_{M|K}\bar c. \tag*{$\fin$}
\end{align*}
\end{proof}

Con este resultado podemos aplicar los resultados obtenidos para
campos locales a los campos globales.

Sea $u_{L|K}\in \cob 2LK$ la {\em clase fundamental\index{clase
fundamental}} de la extensi\'on de Galois $L/K$ que est\'a un\'ivocamente
por la f\'ormula $\inv_{L|K}u_{L|K}=\frac 1{[L:K]}+{\ma Z}$. Obtenemos
el resultado general:

\begin{teorema}\label{T17.6.130N}
Sea $L/K$ una extensi\'on finita de Galois de campos globales. El 
homomorfismo del producto copa con la clase fundamental
\[
u_{L|K}\Cup {\underline{\ }}:\co q{G_{L|K}}{\ma Z}\lra \cob {{q+2}}LK
\]
es biyectiva para toda $q\in{\ma Z}$.
$\fin$
\end{teorema}

De esto obtenemos:

\begin{corolario}\label{C17.6.131N}
Sea $L/K$ una extensi\'on finita de Galois de campos globales.
Entonces $\cob 3LK=\{1\}$ y $\cob 4LK\cong\chi(G_{L|K})$.
$\fin$
\end{corolario}

Para $q=-2$, obtenemos

\begin{teorema}[Ley de reciprocidad de Artin\index{ley de reciprocidad
de Artin}\index{Artin!ley de reciprocidad de $\sim$}]\label{T17.6.132N}
Sea $L/K$ una extensi\'on finita de Galois de campos globales. El mapeo
del producto copa con la clase fundamental
\[
\abe G_{L|K}\cong \cob {{-2}}{G_{L|K}}{\ma Z}\xrightarrow[]{\ u_{L|K}\Cup
{\underline{\ \ }}\ } \cob 0LK=C_K/\N_{L/K} C_L,
\]
nos da un isomorfismo can\'onico, es decir el {\em mapeo de 
reciprocidad\index{reciprocidad!mapeo de $\sim$}} entre la abelianizaci\'on
$\abe G_{L|K}$ del grupo de Galois $G_{L|K}$ de $L/K$ y el grupo
residual de las normas $C_K/\N_{L/K} C_L$ del grupo de clases de
id\`eles $C_K$ (mapeo de Nakayama\index{mapeo de 
Nakayama}\index{Nakayama!mapeo de $\sim$}):
\[
\theta_{L|K}\colon \abe G_{L|K}\lra C_K/\N_{L/K} C_L.
\]

El inverso del mapeo $\theta_{L|K}$ es el inducido por el homomorfismo
\[
\psi_{L|K}=(\ , L/K):C_K\lra \abe G_{L|K}
\]
con n\'ucleo $\N_{L/K} C_L$ y recibe el nombre del {\em s\'imbolo de la norma
residual (global)\index{simbolo de la norma residual global@s\'imbolo de la norma
residual global}} o {\em mapeo de reciprocidad (global)\index{mapeo
de reciprocidad global}}.

La sucesi\'on 
\[
1\lra \N_{L/K}C_L\lra C_K\xrightarrow[]{(\ ,L/K)} \abe G_{L|K}\lra 1,
\]
es exacta.
$\fin$
\end{teorema}

Puesto que el mapeo invariante es compatible con el mapeo inflaci\'on (que
es una inclusi\'on) y con el mapeo restricci\'on, el s\'imbolo de la norma
residual, se comporta como en el caso local:

\begin{teorema}\label{T17.6.133N}
Sean $K\subseteq L\subseteq M$ dos extensiones finitas de campos 
globales con $M/K$ de Galois. Entonces, los siguientes diagramas son
conmutativos
\lasa
\item Si $L/K$ es de Galois, entonces
\[
\begin{minipage}{5cm}
\xymatrix{
C_K\ar@{->}[rr]^{({\underline{\ }},M/K)}\ar@{->}[d]_{\Id}&&
\abe G_{M|K}\ar@{->}[d]^{\pi}\\
C_K\ar@{->}[rr]^{({\underline{\ }},L/K)}&&\abe G_{L|K}
}
\end{minipage}
\qquad
\begin{minipage}{5cm}
\begin{gather*}
\pi:\Gal(M/K)\lra \Gal(L/K)\cong\\
\cong\frac{\Gal(M/K)}{\Gal(M/L)},\\
\sigma \longmapsto \sigma|_L
\end{gather*}
\end{minipage}
\]
Aqu\'i se tiene $(\idel \alpha, L/K)=\pi(\idel \alpha,M/K)\in\abe G_{L|K}$,
$\idel\alpha\in C_K$.

\item
\begin{gather*}
\xymatrix{
C_K\ar@{->}[rr]^{({\underline{\ }},M/K)}\ar@{^{(}->}[d]_{\iota}&&
\abe G_{M|K}\ar@{->}[d]^{\Ver}\\
C_L\ar@{->}[rr]^{({\underline{\ }},M/L)}&&\abe G_{M|L}
}
\intertext{Esto es, $(\idel\alpha,M/L)=\Ver(\idel\alpha,M/K)\in \abe G_{M|L}$,
$\idel \alpha\in C_K$, donde el mapeo de transferencia $\Ver$ est\'a
inducido por restricci\'on:}
\abe G_{M|K}\cong \co {{-2}}{G_{M|K}}{\ma Z}\xrightarrow{\ \res\ }
\co {{-2}}{G_{M|L}}{\ma Z}\cong \abe G_{M|L}
\end{gather*}
y $\iota$ es el encaje natural.

\item
\[
\xymatrix{
C_L\ar@{->}[rr]^{({\underline{\ }},M/L)}\ar@{->}[d]_{\N_{L|K}}&&
\abe G_{M|L}\ar@{->}[d]^{\kappa}\\
C_K\ar@{->}[rr]^{({\underline{\ }},M/K)}&&\abe G_{M|K}
}
\]
Es decir, $(\N_{L|K}\idel \alpha,M/K)=\kappa\big((\idel \alpha, M/L)\big)\in \abe
G_{M|K}$ para $\idel \alpha\in C_L$ y donde $\kappa$ es el encaje
natural $\kappa:\abe G_{M|L}\lra \abe G_{M|K}$.

\item
\[
\xymatrix{
C_K\ar@{->}[rr]^{({\underline{\ }},M/K)}\ar@{->}[d]_{\sigma}&&
\abe G_{M|K}\ar@{->}[d]^{\*\sigma}\\
C_{\sigma K}\ar@{->}[rr]^{({\underline{\ }},\sigma M/\sigma K)}
&&\abe G_{\sigma M|\sigma K}
}
\]
Es decir, $\big(\sigma \idel\alpha,\sigma M/\sigma K\big)=\sigma
\big(\idel\alpha,M/K\big)\sigma^{-1}$ para $\idel\alpha\in C_K$,
donde para $\sigma\in G:=G_{\Omega|K_0}$ los mapeos 
$C_K\stackrel{\sigma}{\lra} C_{\sigma K}$ y $\abe G_{M|K}
\stackrel{\*\sigma}{\lra}\abe G_{\sigma M|\sigma K}$ est\'an inducidas
por $\vec\alpha\longmapsto \sigma \vec\alpha$ y $\tau\longmapsto
\sigma\tau\sigma^{-1}$. $\fin$
\end{list}
\end{teorema}

\begin{definicion}\label{D17.6.134N}
Sea $K$ un campo global. Un subgrupo $H$ de $C_K$ se llama
{\em grupo de normas\index{grupo de normas}} si existe una extensi\'on
finita de Galois $L/K$ tal que $H=\N_{L/K} C_L$. En este caso
denotamos $H={\mc N}_L=\N_{L/K}C_L$.
\end{definicion}

\subsubsection{Ley de reciprocidad v\'ia el isomorfismo de Neukirch}\label{S17.6.10N}

Aqu\'i presentamos otra forma de obtener el mapeo de Nakayama, usando
ahora el isomorfismo de Neukirch.

Hicimos el desarrollo para campos locales, sin embargo, el desarrollo vale
para cualquier formaci\'on, junto con una funci\'on grado y un axioma sobre
los grupos de cohomolog\'ia. Veremos que $(G_{\Omega},C_{\Omega})$
satisface esta condici\'on, donde $K$ es un campo global, $\Omega$ una
cerradura separable, $\Omega=\sep K$ y $C_{\Omega}=\bigcup_L C_L$,
$L/K$ variando en las extensiones finitas de Galois. Entonces $(G_{\Omega},
C_{\Omega})$ es una formaci\'on.

Necesitamos un epimorfismo continuo llamado {\em grado}: $\deg: G=G_{
\Omega}\lra \hat{\ma Z}$.

Empezamos por usar la Definici\'on \ref{D17.6.106N}. Sea $\vec\alpha
\in J_L$, $L/K$ una extensi\'on finita de Galois de campos globales:
\[
[\vec\alpha,L/K]=\langle \vec\alpha,L/K\rangle:=\prod_{\pK\in{\ma P}_K}
(\alpha_{\pK},L_{\pK}/K_{\pK}),
\]
donde $(\underline{\ },L_{\pK}/K_{\pK})$ son los mapeos de reciprocidad
local y donde para cada $\pK\in{\ma P}_K$ seleccionamos un \'unico $\pL
\in{\ma P}_L$ con $\pL|\pK$ y $L_{\pK}=L_{\pL}$.

\begin{teorema}\label{T17.6.135N}
Si $L/K$ y $L'/K'$ son dos extensiones finitas abelianas de campos
globales, tales que $K\subseteq K'$ y $L\subseteq L'$, entonces
el diagrama
\[
\xymatrix{
J_{K'}\ar@{->}[rr]^{[{\underline{\ }},L'/K']}\ar@{->}[d]_{\N_{K'/K}}&&
\abe G_{L'|K'}\ar@{->}[d]^{\rest}&\sigma\ar@{->}[d]\\
J_K\ar@{->}[rr]_{[{\underline{\ }},L/K]}&&\abe G_{L|K}&\sigma|_L
}
\]
es conmutativo.

En otras palabras, si $\vec\alpha\in J_{K'}$,
\[
[\N_{K'/K}\vec\alpha,L/K]=[\vec\alpha, L'/K']|_L.
\]
\end{teorema}

\begin{proof}
Se tiene para $\vec\alpha\in J_{K'}$,
\[
(\alpha_{\eu q},L'_{\eu q}/K'_{\eu q})|_{L_{\pK}}=\big(\N_{K'_{\eu q}/
K_{\pK}}(\alpha_{\eu q}),L_{\pK}/K_{\pK}\big),
\]
con ${\eu q}\in {\ma P}_{K'}$, ${\eu q}|\pK$, $\pK\in {\ma P}_K$.

Ahora $(\N_{K'/K}\vec\alpha)_{\pK}=\prod_{{\eu q}|\pK}
\N_{K'_{\eu q}/K_{\pK}}(\alpha_{\eu q})$ por lo que
\begin{align*}
[\N_{K'/K}\vec\alpha, L/K]&=\prod_{\pK}\big((\N_{K'/K}\vec\alpha)_{\pK},
L_{\pK}/K_{\pK}\big)=\prod_{\pK}\prod_{{\eu q}|\pK}\big(\N_{K'_{\eu q}/K_{
\pK}}(\alpha_{\eu q}), L_{\pK}/K_{\pK}\big)\\
&=\prod_{\eu q}(\alpha_{\eu q}, L'_{\eu q}/K'_{\eu q})|_L=[\vec\alpha,
L'/K']|_L. \tag*{$\fin$}
\end{align*}
\end{proof}

Se tiene el homomorfismo
\[
[\underline{\ \ },L/K]:J_K\lra \Gal(L/K)
\]
para una extensi\'on abeliana arbitraria $L/K$, la cual podr\'ia ser
de grado infinito, por medio de las restricciones $[\underline{\ \ },
L/K]|_{L'}:=[\underline{\ \ }, L'/K]$ donde $L'$ recorre las subextensiones
finitas de $L/K$.

En otras palabras, si $\vec\alpha\in J_K$, consideramos $[\vec\alpha,
L'/K]\in G_{L'|K}$ que forma parte del l\'imite proyectivo $\lim\limits_{
\substack{\longleftarrow\\ L'}} G_{L'|K}=G_{L|K}$ y $[\vec\alpha,L/K]=
\lim\limits_{\substack{\longleftarrow\\ L'}}[\vec\alpha,L'/K]$ es este
elemento despu\'es de identificar $G_{L|K}$ con $\lim\limits_{
\substack{\longleftarrow\\ L'}}G_{L'|K}$. La igualdad
\[
[\vec\alpha,L/K]=\prod_{\pK}(\alpha_{\pK},L_{\pK}/K_{\pK})
\]
permanece v\'alida en el sentido de que el producto infinito del lado 
derecho converge a $[\vec\alpha,L/K]$ en el grupo topol\'ogico $G_{
L|K}$. M\'as precisamene, si $L'/K$ es finita, el conjunto $S_{L'}=\{\pK
\mid (\alpha_{\pK},L'_{\pK}/K_{\pK})\neq 1\}$ es finito y se pueden
considerar a los elementos
\[
\xi_{L'}:=\prod_{\pK\in S_{L'}} (\alpha_{\pK},L_{\pK}/K_{\pK})\in G_{L|K}.
\]

Sea $[\vec\alpha,L/K]\cdot G_{L|M}$ una vecindad abierta b\'asica de
$[\vec\alpha,L/K]$, esto es, $M/K$ es una subextensi\'on finita de $L/K$.
Entonces $\sigma_{L'}\in [\vec\alpha,L/K]\cdot G_{L/M}$ para toda
$L'\supseteq M$ pues
\[
\sigma_{L'}|_M=\prod_{\pK}(\alpha_{\pK}, M_{\pK}/K_{\pK})=[\vec\alpha,
M/K]=[\vec\alpha,L/M]|_M,
\]
lo cual prueba que $[\vec\alpha, L/K]$ es el \'unico punto de acumulaci\'on
de $\{\sigma_{L'}\}$. Se sigue que el teorema anterior sigue compli\'endose
para extensiones infinitas $L$ y $L'$ de $K$ y $K'$ respectivamente.

El siguiente resultado est\'a relacionado con la teor\'ia de Iwasawa.

\begin{teorema}\label{T17.6.136N}
Sea $\Omega$ la m\'axima extensiones abeliana de ${\ma Q}$, es decir,
(Kronecker-Weber) $\Omega=\bigcup_{n=1}^{\infty} \cic n{}$ es el campo 
generado por las ra\'ices de unidad. Sea $T$ el subgrupo de torsi\'on
de $\Gal(\Omega/{\ma Q})$, es decir, todos los elementos de orden 
finito de $\Gal(\Omega/{\ma Q})$. Sea $\tilde{\ma Q}=\Omega^T$
el campo fijo de $T$. Entonces $\Gal(\tilde{\ma Q}/{\ma Q})\cong
\hat {\ma Z}$.
\end{teorema}

\begin{proof}
Se tiene
\[
\Gal(\Omega/{\ma Q})=\Gal\big(\lim_{\substack{\lra\\ n}} \cic n{}/{\ma Q}
\big)\cong \lim_{\substack{\longleftarrow\\ n}}\Gal(\cic n{}/{\ma Q})\cong
\lim_{\substack{\longleftarrow\\ n}}
\big({\ma Z}/n{\ma Z}\big)^*=\hat {\ma Z}^*.
\]

Ahora, se tiene $\hat{\ma Z}\cong \prod_{\text{$p$ primo}} {\ma Z}_p$
y $\*{\ma Z}_p\cong {\ma Z}_p\times {\ma Z}/(p-1){\ma Z}$ para $p\neq 2$
y $\*{\ma Z}_2\cong {\ma Z}_2\times {\ma Z}/2{\ma Z}$. Por tanto
\[
\Gal(\Omega/{\ma Q})\cong \hat{\ma Z}^*\cong \hat{\ma Z}\times\hat T,
\quad\text{donde}\quad \hat T\cong \prod_{p\neq 2}{\ma Z}/(p-1){\ma Z}
\times {\ma Z}/2{\ma Z}.
\]

Por tanto el grupo de torsi\'on $T$ de $\Gal(\Omega/{\ma Q})$ es isomorfo
al grupo de torsi\'on $\hat T$. Se tiene que $\bigoplus_{p\neq 2}{\ma Z}/
(p-1){\ma Z}\bigoplus {\ma Z}/2{\ma Z}\subseteq \hat T$, y que la
cerradura $\bar T$ de $T$ es isomorfo a $\hat T$.

Se sigue que $\tilde{\ma Q}={\ma Q}^T={\ma Q}^{\bar T}$, $\Gal(\tilde
{\ma Q}/{\ma Q})\cong \Gal(\Omega/{\ma Q})/\bar T\cong \hat{\ma Z}$.
$\fin$
\end{proof}

Procediendo como en el caso de campos locales, fijemos un isomorfismo
$\Gal(\tilde{\ma Q}/{\ma Q})\cong \hat{\ma Z}$. En el caso de campos de
funciones, consideremos $K_0=\F(T)$ y para $n\in{\ma N}$, sea $k_n=
K_0{\ma F}_{q^n}$ la extensi\'on de constantes. Entonces si $R:=\bigcup_{
n=1}^{\infty} k_n$, se tiene
\begin{align*}
\Gal(R/K_0)&=\Gal\Big(\bigcup_{n=1}^{\infty} k_n/K_0\Big)\cong \Gal\Big(
\lim_{\substack{\lra\\ n}} k_n/K_0\Big)\\
&\cong \lim_{\substack{\longleftarrow\\ n}}\Gal(k_n/K_0)\cong \lim_{
\substack{\longleftarrow\\ n}}{\ma Z}/n{\ma Z}\cong \hat{\ma Z}.
\end{align*}

Para considerar ambos casos, sea $K_0\in\{{\ma Q},\F(T)\}$ y $\tilde
K_0\in \{\tilde{\ma Q},R\}$. Tenemos un epimorfismo continuo
\[
\deg_{K_0}:\Gal\nolimits_{\sep K_0}\lra \hat{\ma Z}\cong \Gal(\tilde K_0/K_0),
\]
donde $G_{\sep K_0}={\mc G}=\Gal(\sep K_0/K_0)$. Para $K/K_0$ una
extensi\'on finita y separable, sea $f_K:=[K\cap \tilde K_0:K_0]$ y
obtenemos un epimorfismo
\[
\deg_K:=\frac 1{f_K}\deg_{K_0}:G_K\lra \hat{\ma Z},
\]
que determina la $\hat{\ma Z}$-extensi\'on $\tilde K:=K\cdot \tilde K_0$ de
$K$.
\[
\begin{minipage}{5cm}
\xymatrix{
& \tilde K_0\cdot K\ar@{-}[d]^{\hat{\ma Z}}\ar@{-}[dl]\\
\tilde K_0\ar@{-}[d]&K\ar@{-}[dl]\\ K}
\end{minipage}
\qquad\qquad
\begin{minipage}{4cm}
Podemos llamar a $\tilde K/K$ la $\hat{\ma Z}$-extensi\'on
ciclot\'omica de $K$.
\end{minipage}
\]

El elemento de $\Gal(\tilde K/K)$ que se mapea a $1\in \hat{\ma Z}$
bajo el isomorfismo $\deg_K$ lo denotamos, como en el caso local,
por $\Fro K$, el ``{\em Frobenius}'', esto es, $\deg_K(\Fro K)=1$ y por
$\Fr LK=\Fro K|_L$ en caso de que $L/K$ es una subextensi\'on finita de
$\tilde K/K$. Este elemento $\Fr LK$ no debe confundirse con el 
automorfismo de Frobenius de primos de $L$.

Para el $G_{K_0}$-m\'odulo $A$ tomamos $C_{\Omega}=\bigcup_K C_K$
donde $K$ recorre las extensiones finitas y separables de $K_0$ y se 
tiene $A_K=C_K$. Definimos el homomorfismo
\[
v_K:J_K\lra \tilde{\ma Z}
\]
como $v_K:=\deg_K\circ [{\underline{\ \ }},\tilde K/K]$, $
J_K\xrightarrow{[{\underline{\ \ }},\tilde K/K]} \Gal(\tilde K/K)
\xrightarrow{\deg_K}\hat{\ma Z}$.

\begin{teorema}\label{T17.6.137N}
Para cada id\`ele principal $a\in\*K$, se tiene $[a,\tilde K/K]=1$.
\end{teorema}

\begin{proof}
La prueba est\'a contenida en la demostraci\'on del Teorema
\ref{T17.6.108N}.
$\fin$
\end{proof}

\begin{corolario}\label{C17.6.138N}
El grado $v_K:J_K\lra \hat{\ma Z}$ induce $v_K: C_K\lra \hat{\ma Z}$.
$\fin$
\end{corolario}

\begin{proposicion}\label{P17.6.139N}
El homomorfismo $v_K:C_K\lra \hat{\ma Z}$ es suprayectivo en el
caso num\'erico y $v_K(C_K)={\ma Z}$ en el caso de campos de
funciones y es una {\em valuaci\'on henseliana\index{valuacion
henseliana@valuaci\'on henselian}} con respecto a $\deg_K$, es
decir
\lasa
\item $v_K(C_{\Omega})={\mc Z}\supseteq {\ma Z}$ y ${\mc Z}/
n{\mc Z}\cong {\ma Z}/n{\ma Z}$ para toda $n\in{\ma N}$.

\item $v_K(\N_{L/K_0}C_L)=f_L{\mc Z}$ para toda extensi\'on finita
y separable $L/K_0$.
\end{list}
\end{proposicion}

\begin{proof}
Primero veamos que $v_K$ es suprayectivo. Sea $L/K$ una subextensi\'on
finita de $\tilde K/K$, entoces el mapeo
\[
[\ ,L/K]=\prod_{\pK}(\ , L_{\pK}/K_{\pK}):J_K\lra G_{L|K}
\]
es suprayectivo pues de hecho, debido a que $(\ ,L_{\pK}/K_{\pK}):
\*K_{\pK}\lra \Gal(L_{\pK}/K_{\pK})$ es suprayectivo, $[J_K, L/K]$ contiene a
todos los subgrupos de descomposici\'on $\Gal(L_{\pK}/K_{\pK})$. Por tanto,
todos los primos $\pK$ de $K$ son totalmente descompuestos en el campo
fijo $M$ de $[J_K,L/K]$. Se sigue que $M=K$ (Corolario \ref{C17.6.74N}).
Obtenemos que $[J_K,L/K]=\Gal(L/K)$.

Se sigue $[J_K,\tilde K/K]=[C_K,\tilde K/K]$ es denso en $\Gal(\tilde K/K)$.

En el caso num\'erico, $C_K\cong C_{K,0}\times {\ma R}^+$. Si $x\in 
{\ma R}^+$, entonces $[x,\tilde K/K]|_L=[x,L/K]$. Sea $[L:K]=n$. Existe
$y\in {\ma R}^+$ con $y^n=x$ por lo que 
\begin{gather*}
[x,L/K]=[y^n,L/K]=[y,L/K]^n=1,
\intertext{por tanto
$[x,\tilde K/K]=1$ y $[{\ma R}^+,\tilde K/K]=1$.
Se sigue que $[C_K,\tilde K/K]=[C_{K,0},\tilde K/K]$ es denso y como 
$C_{K,0}$ es compacto, se tiene que}
[C_K,\tilde K/K]=[C_{K,0},\tilde K/K]=\Gal(\tilde K/K),
\end{gather*}
por lo que $v_K=\deg_K\circ [\ ,\tilde K/K]$ es suprayectiva en el caso num\'erico
y denso en el caso de campos de funciones.

En el caso de campos de funciones, $\tilde K/K$ son las extensiones de
constantes por lo que $(\alpha_{\pK},\tilde K_{\pK}/K_{\pK})=\Fr {\tilde K_{\pK}}
{K_{\pK}}^{v_{\pK}(\alpha_{\pK})}$. Por tanto, $[C_K,\tilde K/K]=\langle
\Fr K{K_0}\rangle\cong {\ma Z}$.

Ahora $v_K(C_K)=\begin{cases}
\hat{\ma Z},& \car K=0,\\
{\ma Z},& \car K=p
\end{cases},$ por lo que se satisface (a) (${\mc Z}=\hat {\ma Z}$ o ${\ma Z}$).

\s

Para (b), se tiene que
\begin{align*}
v_K(\N_{L/K} C_L)&=v_K(\N_{L/K} J_L)=f_{L|K}\deg_L[J_L,\tilde L/L]\\
&=f_{L|K}v_L(C_L)=f_{L|K}{\mc Z}.
\tag*{$\fin$}
\end{align*}
\end{proof}

\begin{observacion}\label{O17.6.140N}
Los grupos de clases de id\`eles $C_K$ satisfacen el axioma de los campos
de clase y se tiene que el par
\[
(\deg_{K_0}:G_{K_0}\lra \hat{\ma Z},\qquad v_{K_0}: C_{K_0}\lra \hat{\ma Z})
\]
satisface todas las condiciones de teor\'ia de campos de clase. Si $K/K_0$ es
una extensi\'on finita y separable, entonces por el Teorema \ref{T17.6.135N}
se tiene que el homomorfismo $v_K=\deg_K\circ [\ ,\tilde K/K]:C_K\lra \hat
{\ma Z}$ satisface
\[
v_K=\frac 1{f_K}\deg\circ [\ ,\tilde K_0/K_0]\circ \N_{K/K_0}=\frac 1{f_K}
v_{K_0}\circ \N_{K/K_0},
\]
que es precisamente el mapeo inducido por la valuaci\'on henseliana $v_{K_0}$
en el sentido de la teor\'ia de campos de clase locales.
\end{observacion}

Por tanto, hemos obtenido el teorema de reciprocidad global de Artin.

\begin{teorema}[Ley global de reciprocidad de Artin\index{ley global de reciprocidad
de Artin}\index{Artin!ley de reciprocidad de $\sim$}]\label{T17.6.141N}
Para toda extensi\'on finita de Galois $L/K$ de campos globales, tenemos un
isomorfismo can\'onico ({\em isomorfismo global de Neukirch\index{isomorfismo
global de Neukirch}\index{Neukirch!isomorfismo global de $\sim$}}):
\begin{gather*}
{\eu N}_{L/K}:\abe {\Gal(L/K)}\stackrel{\cong}{\lra} C_K/\N_{L/K} C_L.
\tag*{$\fin$}
\end{gather*}
\end{teorema}

El invesrso de ${\eu N}_{L/K}$ da un epimorfismo
\[
(\ ,L/K):C_K\lra\abe{\Gal(L/K)}
\]
con n\'ucleo $\N_{L/K} C_L$. El mapeo $(\ ,L/K)$ es el s\'imbolo de la norma
residual global.

\subsubsection{Teorema principal de la teor\'ia de campos de clase globales}

Regresamos a nuestra exposici\'on principal. A continuaci\'on enunciamos el
teorema principal de la teor\'ia global de campos de funciones.

\begin{teorema}[Teorema principal de la
teor{\'\i}a global de campos
de clase, TCCG\index{teorema
principal de teor{\'\i}a global de campo de
clase}\label{CClaseTCCG}]\label{CClaseT4.2.1}\label{T17.6.142N}
Sea $K$ un campo global. Entonces
\las
\item Existe un \'unico homomorfismo continuo
\begin{gather*}
\rho_K=(\ ,K)\colon C_K\longrightarrow\Gal(\abe K/K)\quad \text{o}\\
\rho_K=(\ ,L)\colon J_K\longrightarrow\Gal(\abe K/K), 
\quad \*K\subseteq \ker \rho_K)
\end{gather*}
tal que para todo lugar $\pK$ de $K$ el siguiente diagrama es conmutativo:
\begin{gather*}
\xymatrix{
\*{K_\pK}\ar@{->}[rr]^{\rho_{K_\pK}\phantom{xxx}}\ar@{->}[d]_{\theta}
\ar@{}[drr]|{\circlearrowright}&&
\Gal(\abe{K_\pK}/K_\pK)\ar@{->}[d]^{\rest}\\
C_K\ar@{->}[rr]_{\rho_K\phantom{xxx}}&&\Gal(\abe K/K)
}\\
\intertext{donde}
\theta=\pi\circ \enc {\ }{\pK}\colon \*{K_\pK}\xrightarrow{\enc
{\ }{\pK}} J_K\xrightarrow{\pi}J_K/\*K=C_K, \\
x_\pK\mapsto (\ldots,1,1
x_\pK,1,1,\ldots)\mapsto (\ldots, 1,1,x_\pK,1,1,\ldots)\bmod K^{\ast}
\end{gather*}
 y $\rest$ es la restricci\'on
$\sigma\longrightarrow \sigma|_{\abe K}$.

El mapeo $\rho_K=(\ ,K)$ se llama el {\em mapeo de 
reciprocidad\index{mapeo de reciprocidad}}.

\item Si $K\subseteq E\subseteq L$ con $L/K$ una extensi\'on abeliana
finita, se tiene el siguiente diagrama conmutativo
\begin{gather*}
\begin{CD}
C_E@>{\rho_E}>>\Gal(L/E)\\
@V{\N_{E/K}}VV@VV{\mu}V\\
C_K@>>{\rho_K}>\Gal(L/K)\\
\end{CD}
\end{gather*}
donde $\mu$ es el encaje natural.

\item Para cualquier extensi\'on abeliana finita $L$ de $K$, $\rho_K$
induce un isomorfismo
\[
C_K/\N_{L/K} C_L\xrightarrow[{\underbracket[0pt]{
_{\psi_{L/K}=(\ ,L/K)}}_{\substack{\uparrow 
\\ \text{Artin}}}}]{\cong} \Gal(L/K)
\]
donde $\psi_{L/K}=(\ ,L/K)$ es el mapeo de Artin\index{mapeo de 
Artin}\label{CClasemapeoartin} o s\'imbolo de la norma
residual global\label{simbolonormaglobal}: $\psi_{L/K}=
\tilde{\rho}_K\colon J_K\twoheadrightarrow\Gal(L/K)$, $\ker
\psi_{L/K}= \*K\N_{L/K} J_L$ y donde $\N_{L/K}\colon C_L\longrightarrow
C_K$ es la norma\index{norma de grupo de clases de id\`eles}
 inducida por la norma de los grupos de id\`eles\index{norma de
id\`eles}:
\[
\vec y=(y_\pL)_{\pL\in{\ma P}_L}\in J_L\xrightarrow[]{\N_{L/K}}
\Big(\prod_{\pL|\pK}
\N_{L_\pL/K_\pK}y_\pL\Big)_{\pK\in{\ma P}_K}\in J_K.
\]
Esto es, si $\vec\alpha\bmod\*L\in C_L$, entonces 
se define $\N_{L/K}(\vec\alpha
\bmod \*L)=\N_{L/K}(\vec\alpha)\bmod \*K$ y se tiene $\N_{L/K} (C_L)=
\N_{L/K}(J_L)\*K/\*K$.

Adem\'as $\N_{L/K}$ es una funci\'on abierta.

\item\label{CClaseexistencia} $H\longrightarrow \rho_K^{-1}(H)$ 
es una biyecci\'on entre el conjunto
de todos los subgrupos abiertos de $\Gal(\abe K/K)$ y el conjunto
de todos los subgrupos abiertos de {\'\i}ndice finito de $C_K$.
En particular $\rho_K$ es una funci\'on densa.

\item {\rm{$=$ (\ref{CClaseexistencia}) {\bf{Teorema de Existencia}}}} Para cada
subgrupo abierto $H$ de {\'\i}ndice finito en $C_K$, existe una \'unica
extensi\'on abeliana finita $L/K$ tal que $\N_{L/K}C_L = H$. 

\item Si $L/K$ es una extensi\'on abeliana finita y $S$ es el conjunto de
lugares ramificados en $L/K$ m\'as los primos infinitos, entonces para
$\vec x\in J_K$, sea $(\vec x)^S:=\prod_{\pK\notin S}\pK^{v_\pK(x_\pK)}
\in D_K^S$\index{grupo de divisores primos relativos a un
conjunto $S$}\label{CClaseDKS}
 el grupo de los divisores (o ideales fraccionarios)
 primos relativos a $S$, es decir,
$D_K^S=D_K/\langle S\rangle$, entonces $\tilde{\rho}_K(\vec x)=
\psi_{L/K}((\vec x)^S)$, donde $\psi_{L/K}$ es el mapeo de Artin,
para $\vec x\in J_K^S\label{CClaseJKS}
:=\{\vec y\in J_K\mid y_{\eu q}=1 \text{\ para\ }
{\eu q}\in S\}$ (es decir, $\rho_K$ evaluado en id\`eles con componente
$1$ en los primos ramificados y en los primos infinitos, coincide
con el mapeo usual de Artin (o s{\'\i}mbolo de Artin)) en ideales. 
M\'as precisamente, se tiene
\begin{gather*}
\tilde\rho(\vec x)=\psi_{L/K}((\vec x)^S)=\psi_{L/K}\big(
\prod_{\pK\notin S}\pK^{
v_{\pK}(x_{\pK})}\big)=\prod_{\pK\notin S}\artinp{L/K}{\pK}^{v_{\pK}
(x_{\pK})}. \tag*{$\fin$}
\end{gather*}
\end{list}
\end{teorema}

Como corolario, se tiene

\begin{corolario}\label{CClaseC4.2.2}\label{T17.6.143N} 
Sea $K$ un campo global. Entonces
existe una correspondencia biyectiva
\begin{multline*}
\text{$\{$extensiones abelianas finitas de $K\} \longleftrightarrow$} \\
\{\text{subgrupos abiertos de {\'\i}ndice finito de $C_K\}$}
\end{multline*}
donde la correspondencia est\'a dada por $L\longrightarrow \N_{L/K}C_L
\subseteq C_K$. Si $L\longleftrightarrow H$, entonces $[L:K]=[C_K:H]$
y si $L^{\prime}\longleftrightarrow H^{\prime}$, entonces $L\supseteq L^{\prime}
\iff H\subseteq H^{\prime}$. $\fin$
\end{corolario}

\begin{teorema}\label{T17.6.144N}
Sea $K$ un campo global y sea $H$ un subgrupo de $C_K$. Si $H$ contiene
a un grupo de normas, entonces $H$ mismo es un grupo de normas.
\end{teorema}

\begin{proof}
Sea $L/K$ una extensi\'on abeliana finita tal que $\N_{L/K} C_L\subseteq H$.
Entonces $H=\bigcup_{\text{finita}} x\N_{L/K} C_L$ y como $\N_{L/K} C_L$
es abierto, $H$ es abierto. Sea $\psi_{L/K}$ el isomorfismo de Artin
\begin{gather*}
\psi_{L/K}: C_K/\N_{L/K} C_L\lra \Gal(L/K).
\intertext{Sea $\psi_{L/K}\big(H/\N_{L/K} C_L\big)=\Gal(L/M)$. Por tanto}
\tilde\psi_{L/K}:C_K/H\cong \frac{C_K/\N_{L/K}C_L}{H/\N_{L/K}C_L}
\lra \frac{\Gal(L/K)}{\Gal(L/M)}\cong \Gal(M/K).
\end{gather*}

Por tanto $\tilde\psi_{L/K}=\psi_{M/K}:C_K/\N_{M/K}C_M\lra \Gal(M/K)$
es el isomorfismo de Artin y $H=\N_{M/K}C_M$ es un grupo de normas.
$\fin$
\end{proof}

\begin{definicion}\label{D17.6.144+1N}
Si $H=\N_{L/K}C_L$, donde $L/K$ es una extensi\'on abeliana finita de
campos globales, se dice que $L$ es el {\em campo de clase asociado
a $H$\index{campo de clase asociado}}
\end{definicion}

Veamos la versi\'on global de la correspondencia $L\longleftrightarrow 
\N_{L/K} C_L$.

\begin{teorema}\label{CClaseTC.1}\label{T17.6.145N}
Sea $K$ un campo global. Para una extensi\'on abeliana finita $L/K$,
denotamos ${\mc N}_L=\N_{L/K} C_L$.

Sean $L_1, L_2$ dos extensiones abelianas finitas de $K$. Entonces
\las
\item $L_1\subseteq L_2 \iff {\mc N}_{L_1}\supseteq {\mc N}_{L_2}$.
\item ${\mc N}_{L_1L_2}={\mc N}_{L_1}\cap {\mc N}_{L_2}$
\item $ {\mc N}_{L_1\cap L_2}
={\mc N}_{L_1}{\mc N}_{L_2}$.
\end{list}
\end{teorema}

\begin{proof} 
\noindent
\underline{(1) $\Longrightarrow)$} 
Si $L_1 \subseteq L_2$, entonces 
\[
{\mc N}_{L_2}=
\N_{L_2/K}C_{L_2}=\N_{L_1/K}\N_{L_2/L_1}C_{L_2}
\subseteq \N_{L_1/K}C_{L_1}={\mc N}_{L_1}.
\]

\noindent
\underline{(2)}  Se tiene para $i=1,2$, $L_i\subseteq L_1L_2$.
Por la parte (1), ${\mc N}_{L_1L_2}\subseteq {\mc N}_{L_1}
\cap {\mc N}_{L_2}$.

Rec\'iprocamente, si
$\idel \alpha\in {\mc N}_{L_1}
\cap {\mc N}_{L_2}$ se tiene que 
\begin{gather*}
\psi_{L_i/K}(\idel \alpha)=(\idel \alpha,L_i/K)=
\rest_{L_i}\circ \rho_K(\idel\alpha)\in\Gal(L_i/K)\cong
C_K/{\mc N}_{L_i}. 
\intertext{Por tanto $\psi_{L_i/K}(\vec\alpha)=1$, $i=1,2$. Se tiene el
monomorfismo}
\theta\colon \Gal(L_1L_2/K)\hooklongrightarrow
\Gal(L_1/K)\times \Gal(L_2/K), \quad \sigma\longmapsto 
(\sigma|_{L_1},\sigma|_{L_2}).
\end{gather*}

Sea 
\[
\sigma=(\idel\alpha,L_1L_2/K)\xrightarrow{\ \theta\ }(\sigma|_{L_1},
\sigma|_{L_2})=(\psi_{L_1/K}(\idel\alpha),\psi_{L_2/K}(\idel\alpha))=(1,1)
\]
por lo que $\psi_{L_1L_2/K}(\idel\alpha)=1$. Se sigue que
$\idel\alpha\in{\mc N}_{L_1L_2}$ lo cual implica que ${\mc N}_{L_1}
\cap{\mc N}_{L_2}\subseteq {\mc N}_{L_1L_2}$ y por tanto
${\mc N}_{L_1}
\cap{\mc N}_{L_2}= {\mc N}_{L_1L_2}$.

\noindent
\underline{(1) $\Longleftarrow)$} Sea ahora ${\mc N}_{L_2}\subseteq
{\mc N}_{L_1}$. Entonces ${\mc N}_{L_1}\cap{\mc N}_{L_2}={\mc N}_{
L_1L_2}={\mc N}_{L_2}$. Por tanto
\[
[L_1L_2:K]=|C_K/{\mc N}_{L_1L_2}|=|C_K/{\mc N}_{L_2}|=[L_2:K].
\]
Se sigue que $L_1L_2=L_2$ de donde se sigue que $L_1\subseteq L_2$.

\noindent
\underline{(3)} Se tiene que $L_1\cap L_2\subseteq L_i$, $i=1,2$. Por
tanto ${\mc N}_{L_i}\subseteq {\mc N}_{L_1\cap L_2}$, $i=1,2$.
Se sigue que ${\mc N}_{L_1}{\mc N}_{L_2}\subseteq {\mc N}_{
L_1\cap L_2}$.
 
\[
\xymatrix{
&L_1\ar@{-}[rr]\ar@{-}[d]&& L_1L_2\ar@{-}[d]\\
&L_1\cap L_2\ar@{-}[rr]\ar@{-}[dl]&&L_2\\ K
}
\]
Ahora bien, ${\mc N}_{L_i}\subseteq {\mc N}_{L_1}{\mc N}_{L_2}
=:H$, $i=1,2$, donde $H$ es un subgrupo de \'indice finito en $C_K$
pues $H=\bigcup_{\text{finito}}{\mc N}_{L_i} x$.

Sea $T$ el campo que corresponde a $H$. Entonces $T\subseteq
L_1$ y $T\subseteq L_2$ por lo que $T\subseteq L_1\cap L_2$ lo cual
implica ${\mc N}_{L_1\cap L_2}\subseteq \N_{T/K} C_T=H={\mc N}_{L_1}
{\mc N}_{L_2}$ con lo cual se sigue que ${\mc N}_{L_1\cap L_2}=
{\mc N}_{L_1}{\mc N}_{L_2}$. $\fin$
\end{proof}

\begin{teorema}[Hasse]\label{T17.6.146N}
Sea $L/K$ una extensi\'on abeliana finita de 
campos globales. Entonces,
para $\vec\alpha\in J_L$, $\idel \alpha\in C_L=J_L/\*L$ se tiene
\begin{align*}
[\vec\alpha,L/K]&=(\vec\alpha,L/K)=\prod_{\pK}(\alpha_{\pK},L_{\pK}/K_{\pK})
\quad\text{y}\\
[\idel\alpha,L/K]&=(\idel \alpha,L/K)=\psi_{L/K}(\idel \alpha)=\prod_{\pK}(
\alpha_{\pK},L_{\pK}/K_{\pK}).
\end{align*}
\end{teorema}

\begin{proof}
Puesto que $(\ ,L/K)$ es el s\'imbolo de la norma residual, aplicamos la
Proposici\'on \ref{CCLP17.6.9}, que fue demostrada para campos locales,
pero como hicimos notar, vale para cualquier formaci\'on de clase,
y se tiene que si $\bar{\idel\alpha}:=\idel\alpha\N_{L/K} C_L\in
\cob 0LK$, entonces para todo caracter $\mu\in\chi(G_{L|K})=\co 1{G_{L|K}}
{{\ma Q}/{\ma Z}}$, se tiene 
\[
\mu\big((\bar{\idel\alpha},L/K)\big)=\inv_{L|K}\big(
(\bar{\idel\alpha})\Cup \delta\mu\big).
\]

Por otro lado, por la Proposici\'on \ref{P17.6.107N}, se tiene
\begin{gather*}
\mu\big([\vec\alpha,L/K)]\big)=\inv_{L|K}\big(
(\vec\alpha)\Cup \delta\mu\big),
\intertext{donde $(\vec\alpha):=\vec\alpha\N_{L/K}J_L\in \co 0{G_{L|K}}{J_L}$. El
homomorfismo}
\co q{G_{L|K}}{J_L}\stackrel{\tilde\pi}{\lra} \co q{G_{L|K}}{C_L}
\end{gather*}
mapea $(\vec\alpha)\in \co q{G_{L|K}}{J_L}$ a $\bar{\idel\alpha}\in \co
0{G_{L|K}}{C_L}$, y por tanto mapea $(\vec\alpha)\Cup \delta\mu\in
\co 2{G_{L|K}}{J_L}$ a $(\bar{\idel\alpha})\Cup \delta\mu\in \co 2{G_{L|K}}
{C_L}=\cob 2LK$.

Por tanto $\tilde\pi\big((\bar{\idel\alpha})\Cup \delta\mu\big)=(\bar{\idel\alpha})
\Cup \delta\mu$.

Por el Teorema \ref{T17.6.121N}, obtenemos
\begin{gather*}
\mu\big((\bar{\idel\alpha},L/K)\big)=\inv_{L|K}\big((\vec\alpha)\Cup\delta\mu\big)
=\inv_{L|K}\big((\vec\alpha)\Cup\delta\mu\big)=\mu\big([\vec\alpha,L/K]\big),
\intertext{y puesto que esto es para todo caracter 
$\mu\in\chi(G_{L|K})$, se sigue que}
\big(\bar{\idel\alpha},L/K\big)=[\vec\alpha,L/K]=\prod_{\pK}(\alpha_{\pK},
L_{\pK}/K_{\pK}). \tag*{$\fin$}
\end{gather*}
\end{proof}

El teorema de Hasse provee de varios resultados de gran importancia para
la teor\'ia global de campos de clase.

\begin{corolario}\label{C17.6.147N}
Sea $K$ un campo global. Para $\alpha_{\pK}\in\*K_{\pK}$, sea 
$\ence \alpha\pK=(\ldots, 1,1,\alpha_{\pK},1,1,\ldots)$ el id\`ele cuya
$\pK$-componente es $\alpha_{\pK}$ y las dem\'as componentes son $1$.
Entonces $\big(\ence\alpha\pK,L/K\big)=(\alpha_{\pK},L_{\pK}/K_{\pK})$.

En particular el diagrama
\[
\xymatrix{
\*K_{\pK}\ar@{->}[rrr]^{(\ ,L_{\pK}/K_{\pK})}\ar@{->}[d]_{\enc{}\pK}&&&
\Gal(L_{\pK}/K_{\pK})\ar@{^{(}->}[d]^{\iota}\\
C_K\ar@{->}[rrr]_{(\ ,L/K)}&&&\Gal(L/K)
}
\]
donde $\iota$ es el encaje natural,
es conmutativo. Pasando al l\'imite, se tiene que
\[
\xymatrix{
\*K_{\pK}\ar@{->}[rr]^{\rho_{K_{\pK}}}\ar@{->}[d]_{\pi\circ \enc{}{\pK}}
&&\Gal(\abe K_{\pK}/K_{\pK})\ar@{^{(}->}[d]^{\iota}\\
C_K\ar@{->}[rr]_{\rho_K}&&\Gal(\abe K/K)
}
\]
es un diagrama conmutativo. Esto es el contenido de {\rm{(1)}} del Teorema
{\rm{\ref{T17.6.142N}}}.
\end{corolario}

\begin{proof}
Se tiene
\begin{gather*}
\big(\ence\alpha\pK,L/K\big)=\prod_{\eu q}\big(\big(\ence\alpha\pK\big)_{\eu q},
L_{\eu q}/K_{\eu q}\big)=(\alpha_{\pK},L_{\pK}/K_{\pK}).
\tag*{$\fin$}
\end{gather*}
\end{proof}

\begin{corolario}\label{C17.6.148-1N}
Sea $L/K$ una extensi\'on abeliana finita de campos globales.
Para $\vec\alpha\in J_K$, si $S$ es el conjunto es el conjunto de los
primos de $K$ ramificados en $L$ y de los primos infinitos, se tiene
que
\[
\tilde\rho_K(\vec\alpha)=\psi_{L/K}\big((\vec\alpha)^S\big)=\psi_{L/K}
\Big(\prod_{\pK\notin S}\pK^{v_{\pK}(\alpha_{\pK})}\Big)=\prod_{\pK\notin
S}\artinp{L|K}{\pK}^{v_{\pK}(\alpha_{\pK})}.
\]
Este es el contenido de {\rm{(6)}} del Teorema {\rm{\ref{CClaseTCCG}}}.
\end{corolario}

\begin{proof}
Es inmediato del hecho de que $\artin{L|K}{\pK}=\varphi_{\pK}
=(\pi_{\pK},L_{\pK}/K_{\pK})$ donde $\pK$ es no ramificado ni
infinito, $\varphi_{\pK}$ es el Frobenis y $\pi_{\pK}$ es un
elemento primo en $\pK$.
$\fin$
\end{proof}

\begin{corolario}\label{C17.6.148N}
Si $x\in\*K$ es un id\`ele principal de un campo global, entonces
$(x,L/K)=1$.
\end{corolario}

\begin{proof}
Se tiene $(x,L/K)=[x,L/K]=1$.
$\fin$
\end{proof}

\begin{corolario}\label{C17.6.149N}
Sea $L/K$ una extensi\'on abeliana finita de campos globales. Entonces
\[
\N_{L/K} C_L\cap \*K_{\pK}=\N_{L_{\pK}/K_{\pK}}\*L_{\pK},
\]
donde consideramos $\*{K_{\pK}}\subseteq C_K$ bajo el monomorfismo
$\*{K_{\pK}}\xhookrightarrow{\lceil\ \rceil_{\pK}}{} J_K\xrightarrow{\pi}{} J_K/\*K=C_K$.

En particular, las extensiones locales derivadas
de una extensi\'on global, corresponden
a las componentes locales del grupo de normas global.

\end{corolario}

\begin{proof}
Si $\alpha_{\pK}\in \N_{L_{\pK}/K_{\pK}}\*L_{\pK}$ entonces $\big(
\ence \alpha\pK,L/K\big)=(\alpha_{\pK},L_{\pK}/K_{\pK})=1$ por lo que
$\ence\alpha\pK\in \N_{L/K}C_L$, esto es, $\N_{L_{\pK}/K_{\pK}}\*L_{\pK}
\subseteq \N_{L/K} C_L$ y obviamente $\N_{L_{\pK}/K_{\pK}}\*L_{\pK}
\subseteq \*K_{\pK}$ por lo que $\N_{L_{\pK}/K_{\pK}}\*L_{\pK}\subseteq
\N_{L/K} C_L\cap \*K_{\pK}$.

Rec\'iprocamente, sea $\idel\alpha\in \N_{L/K}C_L\cap \*K_{\pK}$. Entonces
$\idel \alpha$ est\'a representado por un id\`ele $\vec\alpha=\N_{L/K}
\vec\beta$ con $\vec\beta\in J_L$ y tambi\'en por un id\`ele $\ence a{\pK}$,
$a_{\pK}\in\*K_{\pK}$. Por tanto existe $x\in\*K$ tal que $\ence a\pK\cdot
x=\N_{L/K}\vec\beta$. Para cualquier ${\eu q}\neq {\eu q}$, $x$ es una norma
de $L_{\eu q}/K_{\eu q}$ pues si ${\eu Q}$ y ${\eu Q}'$ son primos en $L$ sobre
${\eu q}$, entonces $\N_{L_{{\eu Q}}/K_{{\eu q}}}L_{{\eu Q}}^* 
=\N_{L_{{\eu Q}'}/K_{{\eu q}}}L_{{\eu Q}'}^*$
por lo que $x\in \prod_{{\eu Q}|{\eu q}}\N_{L_{{\eu Q}}/K_{{\eu q}}}L_{{\eu Q}}^*=
\N_{L_{{\eu q}}/K_{{\eu q}}}L_{{\eu q}}^*$. Puesto que 
\[
(x,L/K)=\prod_{\eu q}(x,L_{\eu q}/K_{\eu q})=(x,L_{\pK}/K_{\pK})=1,
\]
se sigue que $x$ es una norma de $L_{\pK}/K_{\pK}$. Puesto que $\ence a\pK
\cdot x$ es norma de un id\`ele, por el Teorema \ref{T17.6.68N}, $\ence a\pK
x$ es norma local para toda completaci\'on $L_{\pK}/K_{\pK}$. Se sigue que
$\ence a\pK\in \N_{L_{\pK}/K_{\pK}}\*L_{\pK}$, por tanto, $\idel \alpha=
\ence a\pK\in \N_{L_{\pK}/K_{\pK}}\*L_{\pK}$ probando que $\N_{L/K} C_L
\cap \*K_{\pK}\subseteq \N_{L_{\pK}/K_{\pK}}\*L_{\pK}$.
$\fin$
\end{proof}

\begin{corolario}\label{C17.6.150N}
Sea $L/K$ una extensi\'on abeliana finita de campos globales. Entonces
$\N_{L/K}C_L$ es cerrado en $C_K$ y por tanto es abierto.
\end{corolario}

\begin{proof}
Sea $\idel \alpha\notin \N_{L/K} C_L$. Por tanto existe
$\pK\in{\ma P}_K$ tal que $\alpha_{\pK}\notin
\N_{L_{\pK}/K_{\pK}}\*L_{\pK}$ pues $\N_{L/K} C_L\cap \*K_{\pK}
=\N_{L_{\pK}/K_{\pK}}\*L_{\pK}$. Como
$\N_{L_{\pK}/K_{\pK}}\*L_{\pK}$ es cerrado, existe $U$ conjunto abierto de
$\*K_{\pK}$ tal que $\alpha_{\pK}\in U$ y $U\cap \N_{L_{\pK}/K_{\pK}}\*L_{
\pK}=\emptyset$.

Sea $W$ cualquier abierto de $C_K$ con $\idel\alpha\in W$ de la forma
$W=U\times \prod_{{\eu q}\neq \pK} V_{\eu q}$. Entonces $\N_{L/K} C_K\cap
W=\emptyset$ pues $\N_{L_{\pK}/K_{\pK}} L_{\pK}^*
\cap U=\emptyset$. Se sigue que
$\N_{L/K} C_L$ es cerrado en $C_K$.
$\fin$
\end{proof}

Podemos usar l\'imites y obtener el {\em s\'imbolo universal de la norma
residual\index{simbolo universal de la norma residual@s\'imbolo universal
de la norma residual}}. Para un campo global $K$ y una extensi\'on
abeliana finita, tenemos el homomorfismo $C_K\xrightarrow{(\ ,L/K)}
\Gal(L/K)$. Se tiene $\bigcup_{\substack{\text{$L$ abeliana}\\ \text{finita
de $K$}}} L=\abe K$ es la m\'axima extensi\'on abeliana de $K$. Se tiene
\[
\abe G_K=\Gal(\abe K/K)=\Gal(\lim_{\substack{\lra\\ L}}L/K)\cong
\lim_{\substack{\longleftarrow\\ L}}\Gal(L/K)
\]
donde $L$ recorre todas las extensiones abelianas finitas de $K$.

Para $\idel \alpha\in C_K$, obtenemos 
\begin{gather*}
\big(\idel\alpha,K\big):=\lim_{\substack{\longleftarrow\\ L}}\big(\idel
\alpha, L/K\big)\in \abe G_K
\intertext{formado por los elementos compatibles 
$\big\{\big(\idel\alpha,L/K\big)
\in G_{L|K}\big\}_L$. Por tanto obtenemos}
C_K\xrightarrow {(\ \ ,K)}\abe G_K
\end{gather*}
y el n\'ucleo de $(\ ,K)$ es $\bigcap\limits_{\substack{\text{$L/K$ abeliana}\\
\text{finita}}}\N_{L/K}C_L={\eu N}_K=
\ker(\underline{\ \ }, K)$ y la imagen es densa en
$\abe G_K$. De hecho, se tiene que $(\ \ ,K)$ es suprayectivo en el caso 
num\'erico y ${\eu N}_K$ es la componente conexa de $\idel 1\in C_K$. 
En el caso de campos de funciones, ${\eu N}_K=\{1\}$, es decir $(\ , K)$
es inyectiva pero no suprayectiva. Esto lo veremos m\'as adelante.

En el caso local tenemos $\*K_{\pK}\xrightarrow{(\ ,K_{\pK})}\abe G_{K_{\pK}}
\subseteq \abe G_K$, donde $(\ ,K_{\pK})$ es el s\'imbolo universal de la norma
residual local y pasando a los l\'imites, obtenemos para $\idel \alpha\in C_K$,
$\big(\idel \alpha, K\big)=\prod_{\pK}(\alpha_{\pK}, K_{\pK})$.

\begin{teorema}\label{T17.6.151N}
Sea $K$ un campo global. El s\'imbolo 
universal de la norma residual $\rho_K:C_K\lra
\abe G_K=\Gal(\abe K/K)$ es un mapeo continuo. Este es el contenido
de {\rm{(1)}} del Teorema {\rm{\ref{CClaseTCCG}}}.
\end{teorema}

\begin{proof}
Sea $U$ una vecindad abierta de $1\in\abe G_K$, esto es, $U$ es un
subgrupo abierto de $\abe G_K$ de \'indice finito. Sea $L$ el campo fijo de
$U$: $L=(\abe K)^U$, $L/K$ es una extensi\'on abeliana finita. Sea $\vec
\alpha\in J_K$. Si $\pK\in{\ma P}_K$ es no ramificado en $L$ y $\alpha_{
\pK}$ es una unidad, $(\alpha_{\pK},L_{\pK}/K_{\pK})=1$. Por tanto, si 
$\pK$ es no ramificado, $(U_{\pK},K)\subseteq U$. Ahora, si $\pK\in
{\ma P}_K$ es ramificado, se sigue de la teor\'ia local,
Teorema \ref{CCLT17.6.27}, que existe una vecindad
$V_{\pK}$ de $1\in\*K_{\pK}$ tal que $(V_{\pK},L_{\pK}/K_{\pK})=1$.
De esta forma obtenemos que $(V_{\pK},K)\subseteq U$.

Sea $V=\prod_{\substack{\text{$\pK$ no}\\ \text{ramificado}}}U_{\pK}\times
\prod_{\substack{\text{$\pK$}\\ \text{ramificado}}} V_{\pK}$, entonces $V$
es una vecindad abierta
de $\vec 1\in J_K$ y $(V, K)\subseteq U$. Puesto que $C_K$ 
tiene la topolog\'ia cociente y $(\*K,K)=1$, se sigue la continuidad.
$\fin$
\end{proof}

\subsection{Teorema de existencia}\label{S17.6.11N}

Hemos obtenido una biyecci\'on entre las extensiones abelianas finitas
de $K$, $K$ un campo global, y los grupos de normas de $C_K$. El
teorema de existencia, nos caracteriza estos grupos.

Recordemos las estructuras topol\'ogicas que hemos obtenido:
\l
\item[$\bullet$] $J_K$ es un grupo Hausdorff localmente compacto
(Proposici\'on \ref{P17.6.31N}).

\item[$\bullet$] $C_K$ es un grupo Hausdorff localmente compacto
(Proposici\'on \ref{P17.6.31N}).

\item[$\bullet$] $\*K$ es discreto y por tanto cerrado en $J_K$
(Proposici\'on \ref{P17.6.34N}).

\item[$\bullet$] $C_{K,0}$ es cerrado en $C_K$
 para $K$ global. Si $K$ es de funciones,
$C_{K,0}$ tambi\'en es abierto en $C_K$ (Proposici\'on \ref{P17.6.34-1N}).

\item[$\bullet$] $C_{K,0}$ es compacto (Teorema \ref{T17.6.48N}).

\item[$\bullet$] El epimorfismo can\'onico $J_K\lra C_K$ es continuo y manda
conjuntos abiertos en conjuntos abiertos.

\item[$\bullet$] $\enc{\ }\pK:\*K_{\pK}\lra C_K$ es un monomorfismo continuo.

\item[$\bullet$] Si $K$ es de funciones, $C_K$ y $C_{K,0}$ son totalmente
disconexos (Corolario \ref{C17.6.32N}).

\end{list}

Ahora veamos con m\'as detalle los mapeos norma y grado en campos globales.
Ver Definicion \ref{D17.6.28N} para las propiedades fundamentales de estos
mapeos.

Sea $\|\ \|: J_K\lra {\ma R}^+$ el mapeo norma
\[
\|\vec\alpha\|=\prod_{\pK\in{\ma P}_K}|\alpha_{\pK}|_{\pK},
\]
con todos los valores absolutos $|\ |_{\pK}$ normalizados.
Notemos que la imagen de $\|\ \|$ es ${\ma R}^+$ en el caso num\'erico
y $\{q^n\mid n\in{\ma Z}\}$ en el caso de campos de funciones, donde
$\F$ es el campo de constantes de $K$. Esto se debe al Teorema de
Schmidt, Teorema \ref{T17.6.77N}, que establece la existencia
de un divisor de grado $1$, como detallamos un poco m\'as 
adelante.

Ahora bien, $\{q^n\mid n\in{\ma Z}\}\xrightarrow[\log_q]{\cong}{\ma Z}$
es un isomorfismo de grupos y de espacios topol\'ogicos, ambos con la 
topolog\'ia discreta. 

Para $K$ un campo de funciones, la funci\'on grado $\deg$ es simplemente
$\deg:=\log_q \circ \|\ \|$, esto es,
\[
\deg: J_L\lra {\ma Z},\quad \deg\vec\alpha =\log_q(\|\vec\alpha\|)=
\sum_{\pK}\deg \alpha_{\pK}=\sum_{\pK}\deg \pK \cdot v_{\pK}(\alpha_{\pK}).
\]
Se tiene que $\ker\|\ \|=\ker \deg=J_{K,0}$.

Recordemos que tanto $\|\ \|$ como $\deg$ son funciones continuas.
Puesto que $\|x\|=1$ para $x\in \*K$ (y adem\'as $\deg x=0$ si $K$
es campo de funciones), se tienen los mapeos inducidos
\begin{align*}
\|\ \|&: C_K\lra {\ma R}^+,\quad \text{$K$ un campo global},\\
\deg &: C_K\lra {\ma Z}, \quad \text{$K$ un campo de funciones},
\end{align*}
y $\ker \|\ \|=\ker \deg = C_{K,0}$.
Notemos que, por el Teorema de Schmidt, Teorema \ref{T17.6.77N},
que $\deg$ es suprayectiva.

Las sucesiones exactas de homomorfismos continuos
\begin{gather*}
1\lra J_{K,0}\xhookrightarrow{\ \ \imath\ \ } J_K\xrightarrow{\ \|\ \ \|\ } \Lambda\lra 1,\\
1\lra J_{K,0}\xhookrightarrow{\ \ \imath\ \ } J_K\xrightarrow{\ \deg\ } \Lambda\lra 1,\\
1\lra C_{K,0}\xhookrightarrow{\ \ \jmath\ \ } C_K\xrightarrow{\ \|\ \ \|\ } \Lambda\lra 1,\\
1\lra C_{K,0}\xhookrightarrow{\ \ \jmath\ \ } C_K\xrightarrow{\ \deg\ } \Lambda\lra 1,
\end{gather*}
donde $\Lambda=\begin{cases} 
{\ma R}^+&\text{si $K$ es num\'erico},\\
{\ma Z}&\text{si $K$ es campo de funciones},
\end{cases}$
\quad y $\deg$ est\'a definido \'unicamente si $K$ es campo de funciones,
se escinden. M\'as precisamente, primero consideremos $K$ campo num\'erico,
$\Lambda={\ma R}^+$. Sea $\pK$ un lugar infinito de $K$ y consideremos
${\ma R}^+\subseteq K_{\pK}$. Sea $\varphi:{\ma R}^+\lra J_K$ o $C_K$ donde
$\varphi(x)=(\ldots,1,1,x,1,1,\ldots)$ es el id\`ele o la clase de id\`ele, donde la
entrada en $\pK$ es $x$ y las dem\'as entradas son $1$. Entonces $\varphi$
es continuo, donde consideramos a ${\ma R}^+$ con la topolog\'ia
inducida por la topolog\'ia usual de ${\ma R}$
y se tiene $\|\ \|\circ \varphi=\Id_{{\ma R}^+}$. Tambi\'en, sea
$\psi:J_K\lra J_{K,0}$ (o $\psi:C_K\lra C_{K,0}$), donde $\psi(\vec\alpha)$
es el id\`ele donde todas las entradas ${\eu q}$ distintas de $\pK$ son
$\alpha_{\eu q}$ y su entrada en $\pK$ es $\frac{\alpha_{\pK}}{\|\vec\alpha\|}$.
Entonces $\psi$ es continua y $\psi\circ \imath=\Id_{J_{K,0}}$ (o $
\Id_{C_{K,0}}$).

Ahora, si $K$ es un campo de funciones, $\Lambda={\ma Z}$. Veamos que
debemos considerar ${\ma Z}$ con la topolog\'ia discreta. Sea
$\vec\beta$ un id\`ele fijo de grado $1$, el cual existe por el Teorema de
Schmidt. La topolog\'ia de ${\ma Z}$ debe corresponder a la topolog\'ia de
$\{\vec\beta^n\}_{n\in{\ma Z}}$ con la topolog\'ia 
inducida por la topolog\'ia usual de $J_K$. 
Esto se sigue del hecho de que $U=
\prod_{v\in{\ma P}_K} U_v$ es abierto en $J_K$ y de que
$\vec\beta U\cap \{\vec\beta^n\}_{n\in{\ma Z}}=\{\vec\beta\}$.

Sea $\mu:{\ma Z}\lra J_K$ dada por $\mu(n)=\vec\beta^n$.
Entonces $\mu$ es continua y $\deg\circ\mu=\Id_{\ma Z}$. Tambi\'en
tenemos $\delta:J_K\lra J_{K,0}$ dada por $\delta(\vec\alpha)=
\vec\beta^{-\deg \vec\alpha}\vec\alpha$. Entonces $\delta$ es
continua y $\delta\circ \jmath=\Id_{J_{K,0}}$. Similarmente para $C_K$.

\begin{teorema}\label{T17.6.152N}
Se tienen isomorfismos tanto algebraicos como topol\'ogicos
\begin{align*}
C_K&\xrightarrow[\ \ \nu\ \ ]{\cong} C_{K,0}\times {\ma R}^+, \quad\text{si
$K$ es num\'erico},\\
C_K&\xrightarrow[\ \ \eta\ \ ]{\cong} C_{K,0}\times {\ma Z}, \quad\text{si
$K$ es de funciones},
\end{align*}
donde
\begin{gather*}
\nu(\idel\alpha)=(\psi(\idel\alpha),\|\idel\alpha\|)\quad\text{y}\quad
\eta(\idel\alpha)=(\idel\beta^{-\deg\idel\alpha}\cdot \idel\alpha,\deg\idel\alpha),
\end{gather*}
y donde ${\ma R}^+$ tiene la topolog\'ia inducida por la topolog\'ia usual
de ${\ma R}$ y donde ${\ma Z}$ tiene la topolog\'ia discreta.
$\fin$
\end{teorema}

\begin{observacion}\label{O17.6.153N}
No existe un subgrupo distinguido de representantes de $C_K$ en 
$C_K/C_{K,0}$. En el caso num\'erico, es un lugar infinito $\pK$,
y en el caso de campos de funciones, es una clase de id\`ele
de grado $1$, $\vec\beta$. Esto lo veremos de manera m\'as
clara en el caso de campos de funciones.
\end{observacion}

Los grupos $\prod_{\pK\in S} W_{\pK}\times \prod_{\pK\notin S} U_{\pK}
\subseteq J_K$ es un sistema fundamental de vecindades, donde $S$
es un conjunto finito de lugares y $W_{\pK}$ 
pertenece a un sistema de vecindades
de la identidad de $\*K_{\pK}$.

\subsubsection{Extensiones de constantes}

Sea $K$ un campo global de funciones con campo de constantes $\F$.
Sean $\pK\in{\ma P}_K$ y $K_{\pK}$ la completaci\'on de $K$ en $\pK$.
Sea ${\ma F}_{q_{\pK}}$ el campo residual de $K_{\pK}$ el cual es el
mismo que el de $K$ en $\pK$: ${\ma F}_{q_{\pK}}\cong \o_{\pK}/\pK$.
Se tiene que $q_{\pK}=q^{\deg\pK}$.

\begin{teorema}\label{T17.6.154N}
Sea $L/K$ la extensi\'on de constantes de grado $n$, esto es, $L=K
{\ma F}_{q^n}$. Entonces si escribimos $C_K\cong C_{K,0}\times 
{\ma Z}$, se tiene que $\N_{L/K}C_L=C_{K.0}\times n{\ma Z}$.
\end{teorema}

\begin{proof}
Del Teorema de Schmidt tenemos que 
$\deg:C_K\lra {\ma Z}$ es suprayectiva.
Se tiene que la extensi\'on $L/K$ es no ramificada por lo que $L_{\pK}/
K_{\pK}$ es no ramificada. Se sigue que $\alpha_{\pK}\in\*K_{\pK}$
se tiene
\[
(\alpha_{\pK}, L_{\pK}/K_{\pK})=\Fr {L_{\pK}}{K_{\pK}}^{v_{\pK}(\alpha_{\pK})},
\]
donde $\Fr {L_{\pK}}{K_{\pK}}$ es el automorfismo de Frobenius correspondiente
a $L_{\pK}/K_{\pK}$.

Puesto que $L/K$ es c\'iclica, digamos $\Gal(L/K)=\langle\Fr LK\rangle$, donde
se tiene $\langle\Fr LK^{\deg\pK}\rangle\cong \Gal(L(\pL)/K(\pK))$. Se sigue que 
\[
\Fr {L_{\pK}}{K_{\pK}}=\Fr LK^{\deg \pK}\quad\text{y}\quad
\Fr {L_{\pK}}{K_{\pK}}^{v_{\pK}(\alpha_{\pK})}=\Fr LK^{\deg\pK\cdot v_{\pK}(\alpha_{
\pK})}.
\]

Se sigue que si $\vec\alpha\in J_K$, entonces
\[
(\vec\alpha,L/K)=\prod_{\pK\in{\ma P}_K}(\alpha_{\pK}, L_{\pK}/K_{\pK})=
\prod_{\pK}\Fr LK^{\deg \pK\cdot v_{\pK}(\alpha_{\pK})}=\Fr LK^{\sum_{\pK}
\deg\pK\cdot v_{\pK}(\alpha_{\pK})}=\Fr LK^{\deg \vec\alpha}.
\]

Por tanto $(\vec\alpha, L/K)=1\iff n|\deg\vec\alpha$, 
de donde obtenemos $\N_{L/K}
C_L=\ker (\underline{\ \ },L/K)=C_{K,0}\times n{\ma Z}$.
Adem\'as, se verifica que $\Gal(L/K)\cong\frac{C_K}{\N_{L/K}C_L}
=\frac{C_{K,0}\times {\ma Z}}
{C_{K,0}\times n{\ma Z}}\cong \frac {\ma Z}{n {\ma Z}}$.
$\fin$
\end{proof}

M\'as adelante probaremos (Teorema \ref{CClaseT4.8.1-1.1})

\s

\noindent{\bf{Teorema}} {\it{
Sea $L/K$ una extensi\'on abeliana finita de campos de funciones globales.
Sea $\F$ el campo de constantes de $K$. Sea 
\[
d:=\min\{n\in{\ma N}\mid\text{existe $\vec\alpha\in\N_{L/K}J_L$ 
con $\deg\vec\alpha=n$}\}.
\] 
Entonces ${\ma F}_{q^d}$ es el campo de constantes de $L$.
}}

\subsubsection{Teorema de existencia en caracter\'istica $0$}

\begin{definicion}\label{D17.6.155N}
Sea $S$ un conjunto finito de primos de un campo global $K$. Sean
$U_{K,S}:=\{\vec\alpha\in J_L\mid\text{$\alpha_{\pK}=1$ para $\pK\in S$,
$\alpha_{\pK}\in U_{\pK}$ para $\pK\notin S$}\}\subseteq J_{K,S}=
\prod_{\pK\notin S} U_{\pK}\times \prod_{\pK \in S} \*K_{\pK}$ y 
$\bar U_{K,S}:=U_{K,S}\*K/\*K\subseteq C_{K,S}$.
\end{definicion}

Se tiene que $\bar U_{K,S}$ no es abierto pero, puesto que $U_{\pK}$
es cerrado en $\*K_{\pK}$, se tiene que $\bar U_{K,S}$ es cerrado en $C_K$.
Las vecindades fundamentales 
$\prod_{\pK\notin S} U_{\pK}\times \prod_{\pK\in 
S}W_{\pK}$ de $1$, contienen a $U_{K,S}$.

Primero probamos un resultado de teor\'ia de 
grupos que se usa frecuentemente
en el c\'alculo de diversos de \'indices en subgrupos de id\`eles y de
subgrupos de clases de id\`eles.

\begin{proposicion}\label{P17.6.156N}
Sea $G$ un grupo abeliano y sean $X,Y,Z$ subgrupos de $G$ tales que
$Y\subseteq X$. Entonces el epimorfismo natural $\mu:X/Y\lra XZ/YZ$
tiene n\'ucleo $\frac{X\cap YZ}{Y}\cong \frac{X\cap Z}{Y\cap Z}$.
\end{proposicion}

\begin{proof}
Es inmediato que $\mu$ es un epimorfismo y que $\ker \mu=\frac{X\cap YZ}{Y}$.
Consideremos ahora $X\cap Z\xhookrightarrow{\ \ \imath\ \ }X\cap YZ\stackrel
{\ \ \pi\ \ }{\longtwoheadrightarrow}\frac{X\cap YZ}{Y}$, $\varphi:=\pi\circ \imath$.

Si $\alpha\in X\cap Z$, $\varphi(\alpha)=\alpha\bmod Y$. Sea $\beta\in X\cap YZ$,
es decir, $\beta=yz$ con $\beta\in X$, $y\in Y$, $z\in Z$. Entonces $z=y^{-1}\beta
\in Z\cap X$. Se sigue que $\varphi(z)=z\bmod Y=y^{-1}\beta\bmod Y=\beta
\bmod Y$, por lo que $\varphi$ es un epimorfismo y $\ker\varphi= (X\cap Z)\cap Y=
Y\cap Z$, de donde se obtiene el resultado.
$\fin$
\end{proof}

\begin{corolario}\label{C17.6.157N}
Si $G$ un grupo abeliano y si $X,Y,Z$ son subgrupos de $G$ tales que
$Y\subseteq X$, entonces
\begin{gather*}
[XZ:YZ]=\frac{[X:Y]}{[X\cap Z:Y\cap Z]}. \tag*{$\fin$}
\end{gather*}
\end{corolario}

Un resultado fundamental para el teorema de existencia para campos num\'ericos,
y cierta familia de campos de funciones,
es el siguiente teorema que tiene como base fundamental la teor\'ia de Kummer.

\begin{teorema}\label{T17.6.158N}
Sea $K$ un campo num\'erico que contiene a las $n$-ra\'ices de la unidad. Si
$K$ es un campo de funciones de caracter\'istica $p>0$, se supone que $p\nmid n$.
Sea $S$ un conjunto finito de primos de $K$ tal que
\las
\item[$\bullet$] $S$ contiene a todos los primos infinitos y a todos los primos
encima de los primos que dividen a $n$ (esta \'ultima condici\'on es vac\'ia
en el caso de campos de funciones).
\item[$\bullet$] $J_K=J_{K,S}\*K$.
\end{list}

Entonces $C_K^n\bar U_{K,S}$ es el grupo de normas de 
la extensi\'on de Kummer $L=K\big(\sqrt[n]{K^S}
\big)/K$.
\end{teorema}

\begin{proof}
Por teor\'ia de Kummer y por el teorema de las unidades de Dirichlet, tenemos que
$\chi(G_{L|K})\cong \big((\*L)^n\cap \*K\big)
/(\*K)^n=(K^S)(\*K)^n/(\*K)^n=K^S/(K^S)^n$ y que
$K^S$ es finitamente generado de rango $s-1=|S|-1$ y contiene a las $n$-ra\'ices
de unidad, por lo que $K^S/(K^S)^n\cong C_n^s$.

Por tanto $G_{L|K}\cong C_n^s$.

Para $\idel\alpha^n\in C_K^n$, se tiene $\big(\idel\alpha^n,L/K\big)=\big(\idel\alpha,
L/K)^n=1$. Por tanto $\idel\alpha^n\in\N_{L/K}C_L=\ker(\underline{\ \ },L/K)$.

Ahora sea $\idel\alpha \in U_{K,S}$. Para probar que $\vec\alpha\in \N_{L/K}
J_L$, se debe obtener que para cada $\alpha_{\pK}$ es norma local de $L_{
\pK}/K_{\pK}$, $L_{\pK}=K_{\pK}\big(\sqrt[n]{K^S}\big)/K_{\pK}$ para toda
$\pK\in{\ma P}_K$.

Para $\pK\in S$, $\alpha_{\pK}=1$ por lo que $\alpha_{\pK}$ es norma local.
Para $\pK\notin S$, $\alpha_{\pK}\in U_{\pK}$. Se seguir\'a que $\alpha_{\pK}$
ser\'a norma si $\pK$ es no ramificado. Ahora bien, todo $a\in K^S$ es una 
unidad para $\pK\notin S$ y puesto que $n$ es primo relativo a la caracter\'istica
$\car K(\pK)$ del campo residual para $\pK\notin S$, la ecuaci\'on $x^n-a=0$
es separable sobre $K(\pK)$ por lo que $K_{\pK}(\sqrt[n]{a})/K_{\pK}$ es no
ramificada, Teorema \ref{T17.6.20N}.

Se sigue que $\bar U_{K,S}\subseteq \N_{L/K}C_L$ y $C_K^n\bar U_{K,S}
\subseteq \N_{L/K}C_L$. Ahora bien, por la ley de reciprocidad, se obtiene
$[C_K:\N_{L/K}C_L]=|\Gal(L/K)|=[K^S:(K^S)^n]=n^s$.

Por otro lado, se tiene:
\begin{align}\label{Ec17.6.158-1N}
[C_K:C_K^n\bar U_{K,S}]&=[J_{K,S}\*K:J_{K,S}^n U_{K,S}\*K]
\igual_{\substack{\uparrow\\ \text{Corolario \ref{C17.6.157N}}}}\nonumber\\
&=\frac{[J_{K,S}:J_{K,S}^n U_{K,S}]}{[J_{K,S}\cap\*K:J_{K,S}^n U_{K,S}\cap
\*K]}=\frac AB,
\end{align}
donde $A=[J_{K,S}:J_{K,S}^n U_{K,S}]=\prod_{\pK\in S}[\*K_{\pK}:(\*K_{\pK})^n]$.
La \'ultima igualdad se debe a que el mapeo $J_{K,S}\lra \prod_{\pK\in S}
\*K_{\pK}/(\*K_{\pK})^n$, $\vec\alpha\longmapsto \prod_{\pK\in S} \alpha_{\pK}
(\*K_{\pK})^n$ es un epimorfismo y el n\'ucleo consiste de los id\`eles 
$\vec\alpha\in J_{K,S}$ tales que $\alpha_{\pK}\in (\*K_{\pK}
)^n$ para $\pK\in S$. Estos
\'ultimos son precisamente los id\`eles $J_{K,S}^n U_{K,S}$.

Ahora bien, se tiene que $[\*K_{\pK}:(\*K_{\pK})^n]=nq^{v_{\pK}(n)}|\mu_n
(K_{\pK})|=n |n|_{\pK}^{-1} |\mu_n(K_{\pK})|$ con $q=|K(\pK)|$ (Proposici\'on
\ref{CCLTP17.6.14}). Puesto que $|n|_{\pK}=1$ para $\pK\notin S$, se tiene
\begin{align*}
A&=[J_{K,S}:J_{K,S}^n U_{K,S}]=\prod_{\pK\in S}[\*K_{\pK}:(\*K_{\pK})^n]
=\prod_{\pK\in S}n\cdot |n|_{\pK}^{-1}\cdot n\\
&=\prod_{\pK\in S} n^2\cdot
|n|_{\pK}^{-1}=n^{2s}\prod_{\pK}|n|_{\pK}^{-1}=n^{2s}\cdot 1=n^{2s}.
\end{align*}

Para $B=[J_{K,S}\cap\*K:J_{K,S}^n U_{K,S}\cap\*K]$, se tiene que $J_{K,S}
\cap \*K=K^S$ y $J_{K,S}^nU_{K,S}\cap \*K=(K^S)^n$. Para verificar esta 
\'ultima desigualdad, notemos en primer lugar que se tiene $(K^S)^n
\subseteq J_{K,S}^nU_{K,S}\cap \*K$ y para la otra contenci\'on, si $x\in
J_{K,S}^nU_{K,S}\cap \*K$, entonces $x=\vec\alpha^n u$ con $\vec\alpha
\in J_{K,S}$ y $u\in U_{K,S}$. Formamos el campo $K(\sqrt[n]{x})$ y
probaremos que $K(\sqrt[n]{x})=K$.

Si $\vec\beta\in J_{K,S}$, $\vec\beta$ es una norma de alg\'un id\`ele de
$K(\sqrt[n]{x})$: de hecho, si $\pK\in S$, entonces $\beta_{\pK}\in \*{K_{\pK}}$
es una norma pues $K_{\pK}(\sqrt[n]{x})=K(\sqrt[n]{\alpha_{\pK}^n})=K_{\pK}$.
Si $\pK\notin S$, $\beta_{\pK}\in U_{\pK}$ y $K_{\pK}(\sqrt[n]{x})= K_{\pK}(
\sqrt[n]{u_{\pK}})/K_{\pK}$ pues $\pK\nmid n$, y por
tanto $\pK$ es no ramificado. Puesto que
$J_K=J_{K,S} \*K$, hemos probado que 
\[
\N_{K(\sqrt[n]{x})/K}C_{K(\sqrt[n]{x})}=C_K,
\]
lo que implica, por la ley de reciprocidad, que $K(\sqrt[n]{x})=K$.

De esta forma $\sqrt[n]{x}=y\in\*K$, $x=y^n\in (\*K)^n\cap K^S=(K^S)^n$.
Por tanto $B=[J_{K,S}\cap\*K:J_{K,S}^n U_{K,S}\cap\*K]=[K^S:(K^S)^n]=
n^s$.

Se sigue que 
\[
[C_K:C_K^n\bar U_{K,S}]=\frac AB=\frac {n^{2s}}{n^s}=n^s=[C_K:
\N_{L/K}C_L]
\]
y por tanto $\N_{L/K}C_L= C_K^n \bar U_{K,S}$.
$\fin$
\end{proof}

\begin{corolario}\label{C17.6.159N}
Sea $K$, ya sea un campo num\'erico o un campo global de funciones
con $\car K=p\nmid n$. Sean $S$ y $n$ como antes y no suponemos
que $K$ contenga a las $n$-ra\'ices de unidad. Entonces $C_K^n\bar U_{
K,S}$ es un grupo de normas.
\end{corolario}

\begin{proof}
Sea $K'=K(\zeta_n)$ con $\zeta_n$ una $n$-ra\'iz primitiva de la
unidad. Sea $S'$ el conjunto de primos de $K'$ que contienen a todos
los primos sobre los de $S$ y suficientemente grande tal que $J_{K'}=
J_{K',S'}\*{(K')}$. Entonces $C_{K'}^n\bar U_{K',S'}$ es el grupo de
normas de una extensi\'on de Galois $L'$
de $K'$. Sea $L$ la m\'inima extensi\'on de Galois 
de $K$ conteniendo a $L'$. Se sigue que
\begin{align*}
\N_{L/K}C_L&=\N_{K'/K}(\N_{L'/K'}(\N_{L/L'} (C_L)))\subseteq 
\N_{K'/K}(\N_{L'/K}(C_{L'}))\\
&=\N_{K'/K}(C_{K'}^n\bar U_{K',S'})
=(\N_{K'/K} C_{K'})^n\cdot \N_{K'/K}
\bar U_{K',S'}\subseteq C_K^n \bar U_{K,S}.
\end{align*}
Esto es, $C_K^n\bar U_{K,S}$  contiene a un subgrupo de normas, por lo
que $C_K^n\bar U_{K,S}$ tambi\'en es a su vez. un grupo de normas.
$\fin$
\end{proof}

A continuaci\'on probamos el teorema de existencia, esencialmente en
caracter\'istica $0$. M\'as adelante daremos otra demostraci\'on.

\begin{teorema}[Teorema de existencia en caracter\'istica $0$\index{teorema de
existencia en caracter\'istica $0$}]\label{T17.6.160N}
Sea $K$ un campo num\'erico. Los grupos de normas de $C_K$ son
exactamente los subgrupos cerrados de \'indice finito en $C_K$.
\end{teorema}

\begin{proof}
Sea ${\mc N}_L=\N_{L/K}C_L\subseteq C_K$ un grupo de normas de una
extensi\'on finita y normal, es decir, de Galois, $L/K$. Se tiene, por la ley
de reciprocidad que $[C_K:{\mc N}_L]=|\abe G_{L|K}|<\infty$.

Ahora bien, la norma $\N_{L/K}:C_L\lra C_K$ es una funci\'on continua y
$C_K\cong C_{K,0}\times \Gamma_K$, $C_L=C_{L,0}\times \Gamma_L$
con $\Gamma_K$ y $\Gamma_L\cong {\ma R}^+$. Consideremos el mapeo
de escisi\'on $\varphi:{\ma R}^+\lra C_K$ definido como sigue. Sea $\pK$
cualquier primo infinito de $K$ y sea $\enc {\ }\pK:\*K_{\pK}\lra C_K$ el encaje 
natural. Se tiene que ${\ma R}^+\subseteq \*K_{\pK}$ es un subgrupo.

Si $\pK$ es real $\enc{\ }\pK:{\ma R}^+\lra C_K$ es la inyecci\'on deseada
pues $\|\enc{x}\pK\|=|x|_{\pK}=x\in{\ma R}^+$. Si $\pK$ es complejo, entonces
$\|\enc x\pK\|=|x|_{\pK}=x^2$ y tomamos el mapeo $x\longmapsto \enc{\sqrt x}
\pK$.

Se tiene que la inyecci\'on $\varphi:{\ma R}^+\lra C_K$ nos da un conjunto
completo de representantes de $C_K/C_{K,0}$. Por tanto, podemos 
suponer que $\Gamma_K=\Gamma_L\subseteq C_L$. Se sigue que
\[
\N_{L/K} C_L=\N_{L/K} C_{L,0}\times \N_{L/K} \Gamma_L=\N_{L/K}
C_{L,0}\times \Gamma_K^{[L:K]}=\N_{L/K}C_{L,0}\times \Gamma_K.
\]
La imagen del grupo compacto $C_{L,0}$ en $C_K$ es compacto, por
tanto cerrado y, puesto que $\Gamma_K\subseteq C_K$ es tambi\'en
cerrado, $\N_{L/K}C_L={\mc N}_L$ es cerrado en $C_K$. Esta misma
conclusi\'on est\'a dada en el Corolario \ref{C17.6.150N}.

Rec\'iprocamente, sea ${\mc N}\subseteq C_K$ un subgrupo cerrado
de \'indice finito: $[C_K:{\mc N}]=n$. Se sigue que $C_K^n\subseteq
{\mc N}$. Por otro lado, ${\mc N}$ es tambi\'en un subgrupo abierto
y por ende contiene a alg\'un subgrupo de la forma $\bar U_{K,S}$.
Puesto que $C_K^n\bar U_{K,S}$ es un grupo de normas para $S$
suficientemente grande y como $C_K^n\bar U_{K,S}\subseteq {\mc N}$,
${\mc N}$ es a su vez, un grupo de normas.
$\fin$
\end{proof}

\begin{observacion}\label{O17.6.161N}
La demostraci\'on del teorema de existencia \ref{T17.6.160N} sigue
cumpli\'endose para campos de funciones tales que $\car K=p\nmid n
=[C_K:{\mc N}]$. El caso general de campos de funciones lo daremos
m\'as adelante y requiere otro tratamiento.
\end{observacion}

\subsubsection{Normas universales}

\begin{definicion}\label{D17.6.162N}
Sea $K$ un campo global. Se define el grupo de {\em normas
universales\index{normas universales}} como
\[
{\eu N}_K:=\bigcap_L \N_{L/K} C_L,
\]
donde $L$ recorre todas las extensiones abelianas finitas de $K$.
\end{definicion}

Sea $\rho_K:C_K\lra \abe G_K=\Gal(\abe K/K)$ el s\'imbolo universal de la
norma residual, esto es, $\rho_K=(\underline{\ \ }, K)$.

\begin{teorema}\label{T17.6.163N}
Se tiene 
\[
\ker\rho_K={\eu N}_K=\bigcap\limits_{\substack{\text{$L/K$ abeliana}\\
\text{finita}}}{\mc N}_L=\bigcap\limits_{\substack{\text{$L/K$ abeliana}\\
\text{finita}}}\N_{L/K}C_L.
\]
\end{teorema}

\begin{proof}
Para $\idel \alpha\in C_K$, se tiene $\rho_K\big(\idel\alpha\big)=\lim\limits_{
\substack{\longleftarrow\\ L}}\rho_K|_L\big(\idel\alpha\big)=\lim\limits_{\substack{
\longleftarrow\\ L}}\big(\idel\alpha,L/K)$. Por tanto
\begin{gather*}
\idel\alpha\in\ker\rho_L\iff \big(\idel\alpha,L/K\big)=1\text{\ para toda extensi\'on
$L/K$ abeliana finita\ }\\
\iff \idel\alpha\in{\mc N}_L=\N_{L/K}C_L\text{\ para toda
extensi\'on $L/K$ abeliana finita\ }\\
\iff \idel\alpha\in \bigcap_L\N_{L/K}C_L
\text{\ para $L/K$ extensi\'on abeliana finita}.
\tag*{$\fin$}
\end{gather*}
\end{proof}

Recordemos que para cualquier extensi\'on abeliana finita $E/F$ de campos
globales ($K\subseteq F\subseteq E\subseteq \abe K$), 
con $E/K$ finita, se tiene que $\N_{E/F}
C_E$ es cerrado y abierto en $C_F$ (Corolario \ref{C17.6.150N}).

\begin{lema}\label{L17.6.164N}
Sea $K$ un campo global, $E/K$ una extensi\'on finita con $E\subseteq \abe K$.
Sea $K\subseteq F\subseteq E$. Entonces $\ker\N_{E/F}=\N_{E/F}^{-1}(1)$ es
un subgrupo compacto de $C_E$.
\end{lema}

\begin{proof}
Se tiene $C_E\cong C_{E,0}\times \Lambda$ donde 
\[
\Lambda=
\begin{cases}
{\ma R}^+&\text{si $E$ es num\'erico},\\ {\ma Z}&\text{si $E$ es de 
funciones}.\end{cases}
\]

Sea $[E:F]=n$, $\N_{E/F}\Lambda=n\Lambda\cong \Lambda$ y
por tanto $\ker\N_{E/F}|_{\Lambda}=\{1\}$.

Se tiene el diagrama conmutativo con filas exactas
\[
\xymatrix{
1\ar@{->}[r]&C_{E,0}\ar@{->}[r]^{\imath}\ar@{->}[d]_{\N_0}
\ar@{}[dr]|{\circlearrowleft}&
C_E\ar@{->}[r]^{\|\ \|}\ar@{->}[d]^{\N}\ar@{}[dr]|{\circlearrowleft}
&\Lambda\ar@{->}[r]
\ar@{->}[d]^{\N_1}&1\\
1\ar@{->}[r]&C_{F,0}\ar@{->}[r]^{\imath}&C_F\ar@{->}[r]^{
\|\ \|}&\Lambda\ar@{-}[r]&1
}
\]
donde $\N_0=\N_{E/F}|_{C_{E,0}}$, $\N=\N_{E/F}$, $\N_1=
\N_{E/F}|_{\Lambda}$. Del Lema de la Serpiente, obtenemos
la sucesi\'on exacta
\begin{gather*}
1\lra \ker \N_0\lra \ker\N\lra 1\lra C_{F,0}/\N_{E/F}C_{E,0}\lra\\
\lra C_F/\N_{E/F}C_E\lra \Lambda/n\Lambda\lra 1.
\end{gather*}

Se sigue que $\ker\N=\ker\N_0$ y que $C_F/\N_{E/F}C_E\cong
C_{F,0}/\N_{E/F}C_{E,0}\times \Lambda/n\Lambda$. Adem\'as
$\Lambda/n\Lambda\cong \begin{cases}
1&\text{si $E$ es num\'erico}\\ {\ma Z}/n{\ma Z}&\text{si $E$
es de funciones}.\end{cases}$

Ahora bien, $\ker\N=\ker\N_0<C_{E,0}$ el cual es compacto y
$\ker\N_0$ es cerrado, por lo que $\ker\N$ es compacto.
$\fin$
\end{proof}

\begin{teorema}\label{T17.6.165N}
Sea $E/F$ una extensi\'on de campos globales tal que $[E:K]<\infty$
y $E\subseteq \abe K$. Entonces ${\eu N}_F=\N_{E/F}{\eu N}_E$.
\end{teorema}

\begin{proof}
Primero veamos que $\N_{E/F}{\eu N}_E\subseteq {\eu N}_F$.
Consideremos $\idel\alpha\in{\eu N}_E=\bigcap_L \N_{L/E}C_L$, donde $L$
recorre a las extensiones abelianas finitas de $E$. Para cada tal $L$,
existe $\idel \beta_L$ tal que $\idel \alpha=\N_{L/E}\idel \beta_L$.

Ahora consideremos $M$ cualquier extensi\'on abeliana finita
de $F$. Se tiene
\begin{gather*}
\N_{ME/F}C_{ME}=\N_{M/F}\N_{ME/M}C_{ME}\subseteq \N_{M/F}C_M.
\intertext{Por tanto $\N_{ME/F}C_{ME}\cap \N_{M/F}C_M=\N_{ME/F}
C_{ME}$ y $ME\supseteq E$. Esto es,}
{\eu N}_F=\bigcap_{L\supseteq F}\N_{L/F}C_L=\bigcap_{L\supseteq E}
\N_{L/F}C_L.
\end{gather*}

Se sigue que para $L\supseteq E$, $\N_{E/F}\idel \alpha=\N_{E/F}\N_{L/E}
\idel \beta_L=\N_{L/F}\idel\beta_L\in \N_{L/F}C_L$. Por tanto
\[
\N_{E/F}\idel \alpha\in\bigcap_{L\supseteq E}\N_{L/F} C_L={\eu N}_F
\quad\text{y}\quad \N_{E/F} {\eu N}_E\subseteq {\eu N}_F.
\]

Ahora veamos que ${\eu N}_F\subseteq \N_{E/F}{\eu N}_E$. Sea $\idel
\alpha\in{\eu N}_F$. Para cada $L\supseteq E$, sea ${\mc T}_L:=
\N_{L/E}C_L\cap \N_{E/F}^{-1}(\idel\alpha)=\{\idel \beta\in C_E\mid
\idel \beta\in \N_{L/E}C_L\text{\ y\ }\N_{E/F}\idel\beta=\idel \alpha\}$.

Verificaremos que $\bigcap_{L\supseteq E}{\mc T}_L\neq \emptyset$ y
en este caso, si $\idel\gamma \in\bigcap_{L\supseteq E} {\mc T}_L$,
entonces $\idel \gamma\in \bigcap_{L\supseteq E}\N_{L/E}
C_L={\eu N}_E$ y
$\N_{E/F}\idel\gamma=\idel\alpha$ y por tanto $\idel\alpha\in \N_{E/F}
{\eu N}_E$ mostrando que ${\eu N}_F\subseteq \N_{E/F}{\eu N}_E$.

Tenemos que cada ${\mc T}_L\neq \emptyset$ pues ya que $\idel\alpha\in
{\eu N}_F$, se tiene que para cada $L\supseteq F$, $\idel\alpha\in
\N_{L/F}C_L=\N_{E/F}(\N_{L/E}C_L)$. Por otro lado el conjunto $\{{\mc T}_L
\}_{L\supseteq E}$ tiene la propiedad de la intersecci\'on finita puesto
que $\bigcap_{i=1}^s{\mc T}_{L_i}\supseteq {\mc T}_L$ para $L\supseteq
L_i$, $1\leq i\leq s$. Se sigue que para probar que $\bigcap_{L\supseteq
E}{\mc T}_L\neq \emptyset$, es suficiente probar que cada ${\mc T}_L$ 
es compacto.

Se tiene que $\N_{L/E}C_L$ es cerrado y $\N_{E/F}^{-1}(\idel\alpha)$
es compacto, por lo que ${\mc T}_L$ es cerrado en $\N_{E/F}^{-1}(
\idel\alpha)$ y por tanto es compacto. El resultado se sigue.
$\fin$
\end{proof}

\begin{lema}\label{L17.6.166N}
Para todo $K\subseteq F\subseteq \abe K$ tal que $[F:K]<\infty$, se tiene
que $\ker \varphi_{F,l}$, donde $\varphi_{F,l}:C_F\lra C_F$ est\'a dada por
$\varphi_{F,l}(\idel\alpha)=\idel\alpha^l$, es un conjunto compacto.
\end{lema}

\begin{proof}
Se tiene que $\ker\varphi_{F,l}=\{\idel\alpha\in C_F\mid \idel \alpha^l=1\}$.
Ahora bien, si $\idel\alpha\in\ker\varphi_{F,l}$, entonces
$1=\|\idel\alpha^l\|=\|\idel\alpha\|^l$ y $\|\idel\alpha\|\in 
{\ma R}^+$, por lo que $\|\idel\alpha\|=1$. Se sigue que $\idel\alpha\in 
C_{F,0}$. Puesto que $\varphi_{F,l}$ es continua, $\ker\varphi_{F,l}$ es 
cerrado en $C_{F,0}$ y como $C_{F,0}$ es compacto, se sigue que 
$\ker\varphi_{F,l}$ es compacto.
$\fin$
\end{proof}

\begin{teorema}\label{T17.6.167N}
Para cada n\'umero primo $l$, se tiene que ${\eu N}_E\subseteq C_E^l$
para cada extensi\'on abeliana finita $E/K$ con $\zeta_l\in E$.
\end{teorema}

\begin{proof} Sea $\idel \alpha \in {\eu N}_E$.

\s

\noindent
\underline{\bf {Si $l=p=\car K$:}} Entonces, como $\idel \alpha\in {\eu N}_E$,
se tiene que $(\idel\alpha, M/E)=1$ donde $M$ es la m\'axima extensi\'on
abeliana de $E$ de exponente $p$. Del Teorema \ref{T17.6.95N} se obtiene
que $\idel\alpha\in C_E^p$ (en este caso, $\zeta_p=1\in E$ para todo $E$).

\s

\noindent
\underline{\bf {Sea $l\neq p=\car K$:}} Sea $E$ tal que $\zeta_l\in E$. Sea
$S$ un conjunto finito no vac\'io de primos de $E$ que contiene a todos los
divisores de $l$, $\pK|l$, los primos arquimediano y suficientemente grande
tal que $J_E=J_{E,S}\*E$. 

Sean $D=\prod_{\pK\in S}(\*E_{\pK})^l\times \prod_{\pK\notin S} U_{\pK}$ y
$E_1=E\big((E^S)^{1/l}\big)$. De la teor\'ia de Kummer, obtenemos
\begin{gather*}
[E_1:E]=[E^S:(E^S)^l]=l^s,\quad\text{donde}\quad s=|S|.
\intertext{Se tiene}
\begin{align*}
[J_E:D\*E]&=[J_{E,S}\*E:D\*E]\igual_{\substack{\uparrow\\ \text{Corolario
\ref{C17.6.157N}}}}=\frac{[J_{E,S}:D]}{[J_{E,S}\cap \*E:D\cap \*E]}\\
&=\frac{[J_{E,S}:D]}{[E^S:(E^S)^l]}=\frac{\prod_{\pK\in S}
[\*E_{\pK}:(\*E_{\pK})^l]}{[E^S:(E^S)^l]}.
\end{align*}
\end{gather*}

De la Proposici\'on \ref{CCLTP17.6.14}, se tiene
\begin{gather*}
[\*E_{\pK}:(\*E_{\pK})^l]=\frac{l^2}{|l|_{\pK}}.
\intertext{Para ${\pK}\notin S$, $|l|_{\pK}=1$ por lo que}
\prod_{\pK\in S}[\*E_{\pK}:(\*E_{\pK})^l]=\prod_{\pK\in S} \frac{l^2}{|l|_{\pK}}=
\frac{l^{2s}}{\prod_{\pK}|l|_{\pK}}=l^{2s}.
\intertext{Por tanto}
[J_E:D\*E]=\frac{l^{2s}}{l^s}=l^s=[E_1:E].
\end{gather*}

De esta forma obtenemos
\begin{align*}
[J_E:D\*E]&=[E_1:E]=[C_E:\N_{E_1/E}C_{E_1}]\\
&=\big[J_E/\*E:(\N_{E_1/E}J_{E_1})
\*E/\*E\big]=[J_E:(\N_{E_1/E}J_{E_1})\*E].
\end{align*}
Se sigue que $D\*E=(\N_{E_1/E}J_{E_1})\*E$, puesto que, por el
Teorema \ref{T17.6.85N}, tenemos $D\*E\subseteq (\N_{E_1/E}J_{E_1})\*E$.

Sea $\idel\alpha\in{\eu N}_E$. En particular $\idel \alpha$ es norma desde $E_1$,
$\vec\alpha\in (\N_{E_1/D}J_{E_1})\*E=D\*E$. Por tanto, $\idel \alpha$ tiene
un representante $\vec\beta\in D$. Cualesquiera dos representantes de $\idel
\alpha$ difieren por un elemento $\delta\in D\cap \*E=(E^S)^l$.
Ahora $\vec\beta$ es una $l$-potencia para todos los primos de $S$. Veremos que
tambi\'en es una $l$-potencia para todo $\pK\notin S$.

Para $S'\supseteq S$, definimos $D'=\prod_{\pK\in S'}(\*E_{\pK})^l\times
\prod_{\pK\in S'} U_{\pK}$ y se tiene $D'\*E=(\N_{E_1'/E}J_{
E_1'})\*E$, $E_1'=E\big((E^{s'})^{1/l}\big)$.

Se tiene que $\idel\alpha$ tiene un representante $\vec\gamma\in D'$. Se
puede escribir $\vec\gamma =\vec\eta^l\cdot \vec\gamma'$ donde $\vec\gamma'$
es una unidad fuera de $S'$ y $1$ en los primos
de $S'$, $\vec\eta$ es $1$ fuera de
$S'$. Sea $x\in E$ tal que $\vec\eta =x\cdot\vec\xi$ donde $\vec\xi\in J_{E,S}$. 
Entonces $\vec\eta^l=x^l\cdot \vec\xi^l$ y $\vec\gamma=x^l\cdot\vec\xi^l\cdot
\vec\gamma'$. Se sigue $\vec\gamma$ es una $l$-potencia para todos los
primos $l$ en $S'$.
Por tanto $\vec\gamma$ difiere de $\vec\beta$ por una $l$-potencia
y por tanto $\vec\beta$ es una $l$-potencia en $S'$. Puesto que podemos
incluir cualquier primo en $S'$, se sigue que $\idel \alpha\in C_E^l$.
$\fin$
\end{proof}

\begin{teorema}\label{T17.6.168N}
Para cualquier $K\subseteq F\subseteq \abe K$ tal que $[F:K]<\infty$, se tiene
que ${\eu N}_F$ es (infinitamente) divisible, es decir, para toda $m\in{\ma N}$,
${\eu N}_F^m={\eu N}_F$.
\end{teorema}

\begin{proof}
Es suficiente probar que ${\eu N}_F^l={\eu N}_F$ para cualquier n\'umero 
primo $l$.

Se tiene que ${\eu N}_F=\N_{E/F}{\eu N}_E$ para cada extensi\'on finita
$E/F$. Tomando $E$ suficientemente grande para $l$ en el sentido de que
${\eu N}_E\subseteq C_E^l$, se tiene
\begin{gather}\label{Ec17.6.168-1N}
{\eu N}_F=\N_{E/F}{\eu N}_E\subseteq \N_{E/F}C_E^l=\big(\N_{E/F}
C_E\big)^l.
\end{gather}

Sea $\idel\alpha\in {\eu N}_F$. El s\'imbolo $\big(\idel\alpha\big)^{1/l}$
denota al conjunto de elementos de $C_F$ cuya $l$-potencia es $\idel
\alpha$: $\idel\alpha^{1/l}:=\{\idel\beta\in C_F\mid \idel \beta^l=\idel\alpha\}$.

De (\ref{Ec17.6.168-1N}), se obtiene que los conjuntos
\begin{gather*}
X_E:=\big(\N_{E/F}C_E\big)\cap \big(\idel \alpha\big)^{1/l}
\intertext{son no vac\'ios. Por tanto estos conjuntos satisfacen la propiedad de intersecci\'on
finita:}
\bigcap_{i=1}^r X_{E_i}\supseteq X_E\neq \emptyset\quad\text{para $E$ tal que
$E_i\subseteq E$ para toda $1\leq i\leq r$}.
\end{gather*}

Adem\'as $\N_{E/F}C_E$ es cerrado (Corolario \ref{C17.6.150N}) y $\idel\alpha^{1/l}$
es compacto (Lema \ref{L17.6.166N}) por lo que $X_E$ es un conjunto compacto.
Se sigue que $\bigcap_E X_E\neq \emptyset$. Si $\idel\beta\in\bigcap_E X_E$,
entonces $\idel\beta \in {\eu N}_F\bigcap \big(\idel\alpha\big)^{1/l}$, por tanto $\idel
\beta\in {\eu N}_F$ y $\idel\beta^l=\idel\alpha$ lo que prueba el teorema.
$\fin$
\end{proof}

\begin{proposicion}\label{P17.6.169N}
Sea $K$ cualquier campo global y sea $H\subseteq C_K$ cualquier subgrupo
divisible. Entonces $H\subseteq {\eu N}_K$.
\end{proposicion}

\begin{proof}
Sea $h\in H$ y sea $L/K$ una extensi\'on abeliana finita. 
Se tiene que $(h,L/K)\in\Gal(L/K)$. Sea $|\Gal(L/K)|
=n$. Sea $h_1\in H$ con $h_1^n=h$.
Se sigue que
\begin{gather*}
\rho_K(h)|_L=(h,L/K)=(h_1^n,L/K)=(h_1,L/K)^n=1,
\end{gather*}
lo que implica que
$h\in\N_{L/K}C_L$ y $h\in\bigcap_L \N_{L/K}C_L={\eu N}_K$.
$\fin$
\end{proof}

\begin{corolario}\label{C17.6.170N}
Se tiene que ${\eu N}_K$ es el m\'aximo subgrupo divisible de $C_K$. $\fin$
\end{corolario}

\begin{corolario}\label{C17.6.171N}
Si $K$ es un campo num\'erico, $\rho_K$ es suprayectiva y no inyectiva.
\end{corolario}

\begin{proof}
Se tiene que $C_K\cong C_{K,0}\times {\ma R}^+$ y ${\ma R}^+$ es
divisible, por lo que ${\ma R}^+\subseteq {\eu N}_K=\ker \rho_K$. Por tanto
$\rho_K$ no es inyectiva.

Por otro lado, $\rho_K(C_K)=\rho_K(C_{K,0})\subseteq \abe G_K$ y $\rho_K
(C_K)$ es denso en $\abe G_K$ y $\rho(C_{K,0})$ es compacto, de donde se
sigue que  $\rho_K(C_K)=\abe G_K$.
$\fin$
\end{proof}

\begin{lema}\label{L17.6.172N}
Sea $A$ cualquier subgrupo de \'indice finito en $C_K$, donde $K$ es un
campo global. Entonces, si $H<C_K$ es un subgrupo divisible,
se tiene $H\subseteq A$.
\end{lema}

\begin{proof}
Sea $h\in H$. Sea $[C_K:A]=n<\infty$. Existe $h_1\in H<C_K$ con $h_1^n=h
\in A$.
$\fin$
\end{proof}

Con todos estos elementos, damos otra demostraci\'on del teorema de
existencia para campos num\'ericos (Teorema \ref{T17.6.160N}).

\begin{teorema}[Teorema de existencia en caracter\'istica $0$\index{teorema de
existencia en caracter\'istica $0$}]\label{T17.6.173N}
Sea $K$ un campo num\'erico. Sea $H$ un subgrupo abierto de $C_K$ de
\'indice finito. Entonces $H$ es un grupo de normas.
\end{teorema}

\begin{proof}
Se tiene que $\ker \rho_K={\eu N}_K$ es divisible. Por tanto $\ker\rho_K
\subseteq H$. Se tiene que $H$ es cerrado en $C_K$. Sea $H_0:=H\cap
C_{K,0}$. Entonces $H_0$ es cerrado en $C_{K,0}$, por tanto compacto
y $[C_{K,0}:H_0]=[C_K:H]$ pues ${\ma R}^+\subseteq H$ por ser divisible.
Sea $C_{K,0}\stackrel{\mu}{\lra} C_K/H$, $\idel\xi\longmapsto \idel\xi
\bmod H$, el mapeo natural.

Sea $\idel\xi\cdot x\in C_K$, $\idel\xi\in C_{K,0}$, $x\in{\ma R}^+$.
Entonces $\mu(\idel\xi)=\idel\xi\bmod H=\idel \xi x\bmod H$. Por tanto
$\mu$ es suprayectiva y $\ker\mu=H\cap C_{K,0}=H_0$.

Sea ${\mc H}=\rho_K(H_0)$ el cual es un subgrupo cerrado de $\abe G_K$
por ser compacto, y como ${\eu N}_K\subseteq H$, se tiene que $H_0=
\rho_K^{-1}({\mc H})\cap C_{K,0}$. De hecho, se tiene $\rho_K^{-1}({\mc H})=
\rho_K^{-1}(\rho_K(H_0))\supseteq H_0$, de donde $\rho_K^{-1}({\mc H})
\cap C_{K,0}\supseteq H_0$. Rec\'iprocamente, si $\idel\xi\in \rho_K^{-1}({\mc H})
\cap C_{K,0}$, entonces $\rho_K(\idel\xi)\in{\mc H}=\rho_K(H_0)$. Por tanto
existe $\idel \alpha\in H_0$ tal que $\idel\xi=\idel\alpha\cdot\idel\beta$ con
$\idel\beta\in \ker\rho_K\subseteq H$ y $\idel \beta=\idel\xi\idel\alpha^{-1}\in
C_{K,0}$. Se sigue que $\idel\xi\in H\cap C_{K,0}=H_0$.

Como consecuencia de lo anterior, tenemos que $[\abe G_K:{\mc H}]=
[C_{K,0}:H_0]=[C_K:H]$. Sea $L:=(\abe K)^{\mc H}$. Entonces, si
$\idel \alpha\in \N_{L/K}C_L$, se tiene $\rho_K(\idel\alpha)\big|_L=1$,
lo cual implica que $\rho_K(\idel\alpha)\in {\mc H}$. Se sigue que
$\N_{L/K} C_L\subseteq \rho_K^{-1}({\mc H})$.
\[
\begin{minipage}{5cm}
\xymatrix{
{\eu N}_K\ar@{<->}[r]\ar@{->}[d]\ar@{->}@/^1pc/[rr]^{\rho_K}
&\abe K\ar@{<-}[d]^{\mc H}\ar@{<->}[r]&\{1\} \ar@{->}[d]\\ 
\N_{L/K}C_L\ar@{<->}[r]\ar@{->}[d]\ar@{-->}@/^1pc/[rr]_{\rho_K}&L\ar@{<-}[d]
\ar@{<->}[r]&{\mc H}\ar@{->}[d]\\ 
C_K\ar@{<->}[r]\ar@{->}@/_1pc/[rr]_{\rho_K} &K\ar@{<->}[r]&\abe G_K
}
\end{minipage}
\qquad\qquad
\begin{minipage}{5cm}
Se tiene que $[C_K:\N_{L/K}C_L]=[L:K]=[\abe G_K:{\mc H}]=[C_K:H]$.
\end{minipage}
\]
Se sigue que $\N_{L/K}C_L=H$ y esto concluye el teorema de existencia.
$\fin$
\end{proof}

Con respecto al n\'ucleo de $\rho_K$ para un campo num\'erico $K$, Tate
prob\'o el siguiente resultado.

\begin{teorema}\label{T17.6.174N}
Si $K$ es un campo num\'erico y $\rho_K:C_K\lra \abe G_K$ es el mapeo de
reciprocidad, $\ker\rho_K={\eu N}_K$ es la componente conexa de $1$ y 
se tiene
\[
{\eu N}_K\cong\Big(\frac{\hat{\ma Z}\times {\ma R}}{\ma Z}\Big)^t
\oplus \Big(\frac{\ma R}{\ma Z}\Big)^s\oplus {\ma R},
\]
donde $s$ es el n\'umero de lugares complejos de $K$ y $t$ es el rango de las
unidades totalmente positivas de $K$ (${\ma R}^+\stackrel{\log}{\cong} {\ma R}$).
\end{teorema}

\begin{proof}
Ver \cite[Chapter IX, p\'aginas 65--70]{ArtTat61}, \cite[Chapter VI, Section 1,
Exercises 1--10, p\'aginas 367--368]{Neu99}.
$\fin$
\end{proof}

\subsubsection{Teorema de existencia en caracter\'istica $p>0$}

Ahora estudiemos el caso de campos de funciones. Sean $K$ un campo global
de funciones y $\rho_K:C_K\lra \abe G_K$ el s\'imbolo de la norma residual.

Se tiene que $\Gal(K\abe\F/K)\cong \Gal(\abe \F/\F)\cong \hat{\ma Z}$ y 
$\Gal(K\abe \F/K)=\overline{\langle \Fro K\rangle}$ 
donde $\Fro K$ es el automorfismo de
Frobenius, $\abe \F\xrightarrow{\ \Fro K\ } \abe \F$, $x\longmapsto x^q$. 

Sea $\sigma\in \abe G_K$ y $\sigma|_{\abe\F}=\Fro K^{\nu}$. Se define
$\nu:=\ord (\sigma)\in \hat{\ma Z}$.

Si $\idel\alpha\in C_K$, escribimos $\idel\alpha=\idel\alpha_1^n\idel\alpha_0$
con $\idel\alpha_1$ un elemento de $C_K$ de grado $1$ fijado de antemano
y $\idel\alpha_0\in C_{K,0}$. Se tiene $\deg\idel\alpha=n$. En t\'ermino de
valores absolutos, tenemos $\|\idel\alpha_1\|=q$, $\|\idel\alpha_0\|=1$ y
$\deg\idel\alpha_1=1$, $\deg\idel\alpha_0=0$.

Se tiene que $\rho_K(\idel\alpha)=(\idel\alpha,K)
=\prod_{\pK} (\alpha_{\pK},L_{\pK}/K_{\pK})$.

\begin{teorema}\label{T17.6.175N}
Para $\idel\alpha\in C_K$, $\rho_K(\idel\alpha)=(\idel\alpha,K)$, se tiene
\[
\ord\big((\idel\alpha,K)\big)=\deg \big(\idel\alpha\big)\in {\ma Z}.
\]
En particular, si $\idel\alpha\in \ker\rho_K$, necesariamente $\idel\alpha\in
C_{K,0}$, esto es, $\ker\rho_K\subseteq C_{K,0}$.
\end{teorema}

\begin{proof}
Recordemos que $\deg\idel\alpha=\deg\vec\alpha=\sum_{\pK}\deg_{\pK}
\alpha_{\pK}$, donde $\deg_{\pK}=\deg\pK \cdot v_{\pK}(\alpha_{\pK})$.

Ahora bien, $|\alpha_{\pK}|_{\pK}=q_{\pK}^{-v_{\pK}(\alpha_{\pK})}$ donde
$q_{\pK}=|K(\pK)|$ es la cardinalidad del campo residual y $[K(\pK):\F]=
\deg\pK$. Por tanto $q_{\pK}=|K(\pK)|=q^{\deg\pK}$. Esto es,
\begin{gather*}
|\alpha_{\pK}|_{\pK}=q_{\pK}^{-v_{\pK}(\alpha_{\pK})}=q^{-\deg\pK\cdot v_{\pK}
(\alpha_{\pK})}=q^{-\deg_{\pK}\alpha_{\pK}}.
\intertext{De esta forma, tenemos:}
\|\idel\alpha\|=\|\vec\alpha\|=\prod_{\pK}|\alpha_{\pK}|_{\pK}=q^{-\sum_{\pK}
\deg_{\pK} \alpha_{\pK}}=q^{-\deg\vec\alpha}=q^{-\deg\idel\alpha},
\intertext{y en particular}
\|\idel\alpha\|=1\iff \deg\idel\alpha=0.
\end{gather*}

Sean $K_{\pK}$ es cualquier completaci\'on y $\sigma
\in \abe G_{K_{\pK}}=\Gal(\abe K_{\pK}/K_{\pK})$, $K_{\pK}(\pK)={\ma F}_{q_{\pK}}
={\ma F}_{q^{\deg\pK}}$.
\[
\begin{minipage}{5cm}
\xymatrix{
&\abe K_{\pK}\ar@{-}[d]\ar@{--}[dl]\\ K_{\pK}\ar@{-}[r]&K_{\pK}\abe {\ma F}_{
q_{\pK}}
}
\end{minipage}
\hspace{2cm}
\begin{minipage}{5.5cm}
Sea $\sigma|_{K_{\pK}\abe {\ma F}_{q_{\pK}}}\in \Gal(K_{\pK}\abe {\ma F}_{q_{\pK}}
/K_{\pK})\cong \Gal(\abe {\ma F}_{q_{\pK}}/{\ma F}_{q_{\pK}})$.
\end{minipage}
\]

Se tiene que si $\Fro \pK$ es el Frobenius respectivo a $\abe{\ma F}_{q_{\pK}}/
{\ma F}_{q_{\pK}}$, entonces se tiene que $\overline{\langle \Fro \pK\rangle}=\Gal(
\abe{\ma  F}_{q_{\pK}}/{\ma F}_{q_{\pK}})$. Esto es, $\Fro\pK:\abe{\ma  F}_{q_{
\pK}}\lra \abe{\ma  F}_{q_{\pK}}$ est\'a dada por $\Fro \pK(x)=x^{q_{\pK}}=x^{q^{
\deg \pK}}=\Fro K^{\deg\pK}$ puesto que $\abe{\ma  F}_{q_{\pK}}=\abe\F$ y
$\overline{\langle \Fro\pK\rangle}=\Gal(\abe{\ma  F}_{q_{\pK}}/{\ma  F}_{q_{\pK}})
\subseteq \Gal(\abe{\ma  F}_{q_{\pK}}/\F)=\overline{\langle\Fro K\rangle}$.

Si $\sigma|_{K_{\pK}\abe{\ma  F}_{q_{\pK}}}=\Fro \pK^n$, se define $\ord_{\pK}
(\sigma)
=n$. Por tanto, $\ord_{\pK}(\sigma)=n$ implica $\ord_K(\sigma)=
n\cdot\deg \pK$ y
as\'i, $\ord_K(\sigma)=\deg\pK\cdot\ord_{\pK}(\sigma)$ y, en particular,
\[
\ord_K(\alpha_{\pK},K_{\pK})=\deg\pK\cdot \ord_{\pK}(\alpha_{\pK}, K_{\pK}).
\]

Ahora veamos que $\ord_{\pK}(\alpha_{\pK},K_{\pK})=v_{\pK}(\alpha_{\pK})$.
Se tiene $(\alpha_{\pK},K_{\pK})|_{K_{\pK}\abe{\ma  F}_{q_{\pK}}}=
(\alpha_{\pK},K_{\pK}\abe{\ma  F}_{q_{\pK}}/K_{\pK})=\Fro\pK^{v_{\pK}(
\alpha_{\pK})}$ puesto que la extensi\'on $K_{\pK}\abe{\ma  F}_{q_{\pK}}
/K_{\pK}$ es no ramificada.

Se sigue que $\ord_{\pK}(\alpha_{\pK},K_{\pK})=v_{\pK}(\alpha_{\pK})\in {\ma Z}$.

Ahora bien, puesto que $(\idel\alpha,K)=\prod_{\pK\in{\ma Z}_K}(\alpha_{\pK},
K_{\pK})$, se obtiene
\begin{align*}
\ord_K\big(\idel\alpha,K\big)&=\sum_{\pK\in{\ma P}_K}\ord_{\pK}(\alpha_{\pK},K_{\pK})=
\sum_{\pK\in{\ma P}_K}\deg\pK\cdot v_{\pK}(\alpha_{\pK})\\
&=\sum_{\pK\in{\ma P}_{\pK}}\deg_{\pK}\alpha_{\pK}=\deg(\idel\alpha)\in{\ma Z}.
\end{align*}

En particular, si $\idel\alpha\in\ker\rho_K$, se tiene $\big(\idel\alpha,K\big)|_{K
\abe{\ma  F}_{q_{\pK}}}=\Id=\Fro K^0$ y $\ord_K(\idel\alpha,K)=\deg\big(\idel
\alpha\big)=0$ y $\idel\alpha\in C_{K,0}$.
$\fin$
\end{proof}

\begin{observacion}\label{O17.6.176N}
Si  $A$ es un grupo localmente compacto, entonces si $A_0$ es la componente
conexa de la identidad, $A_0$ es un subgrupo abierto y cerrado de $A$ y $A/A_0$
es totalmente disconexo.
\end{observacion}

\begin{proposicion}\label{P17.6.177N}
Se tiene que $A_0=\bigcap_H H$, donde $H$ recorre todos los subgrupos 
normales abiertos de $A$.
\end{proposicion}

\begin{proof}
Se tiene que para $x\in A$, $x^{-1}A_0 x$ es conexo y contiene a la
identidad, por lo que $x^{-1}A_0 x=A_0$ y $A_0$ es un subgrupo normal.
Ahora, si $H$ es cualquier subgrupo abierto de $A$, escribimos $A=
H\cupdot H_1$ con $H_1=\bigcup\limits_{\substack{x\in G\\ x\notin H}}
xH$, por lo que $H_1$ es abierto. Esta es una disconexi\'on pues
$1\in H$. Por tanto $A_0\subseteq H$ y $A_0\subseteq \bigcap_H H$.

Rec\'iprocamente, $H:=A_0$ es normal y abierto, por lo que
$\bigcap_H H\subseteq A_0$.
$\fin$
\end{proof}

Sea $H$ un subgrupo normal cerrado tal que $A/H$ es totalmente
disconexo. Sea $\pi:A\lra A/H$ la proyecci\'on natural, la cual es continua.
Los abiertos de $A/H$ son los conjuntos $U/H$ con $U$ abierto en $A$
y $H\subseteq U$. Puesto que $A/H$ es totalmente disconexo, $A_0\subseteq
\ker \pi=H$ y por tanto la proyecci\'on 
natural $A/A_0\longtwoheadrightarrow
A/H$ es suprayectiva. En otras palabra, $A/A_0$ es el grupo cociente
m\'aximo totalmente disconexo.

\s

Regresando a $C_K$ y $C_{K,0}$ donde $K$ es un campo global
de funciones, recordemos que tanto $C_K$ como $C_{K,0}$ son
totalmente disconexos (Corolario \ref{C17.6.32N}). Por tanto
$\{1\}=\bigcap\limits_{\substack{H\normal C_K\\ \text{$H$ abierto}}}H=
\bigcap\limits_{\substack{H\normal C_{K,0}\\ \text{$H$ abierto de
$C_{K,0}$}}} H$.

Ahora, si $H$ es cualquier subgrupo, ya sea de $C_K$ o de $C_{K,0}$,
de \'indice finito, digamos $[C_K:H]=n<\infty$ o $[C_{K,0}:H]=n<\infty$,
se tiene $\ker\rho_K\subseteq H$. As\'i,
\[
\ker\rho_K\subseteq
\bigcap_{[C_{K,0}:H]<\infty} H\subseteq \bigcap_{\substack{
H\normal C_{K,0}\\ \text{$H$ abierto}}}=D=\{1\},
\]
donde $D$ es la componente conexa de la identidad y $D=\{1\}$.
Por tanto $\ker\rho_K=\{1\}$, esto es, $\rho_K$ es inyectiva.

\begin{teorema}\label{T17.6.178N}
Si $K$ es un campo global de funciones y si $\rho_K:C_K\lra \abe G_K$
es el s\'imbolo de la norma residual, entonces $\rho_K$ es inyectiva
y no suprayectiva.
\end{teorema}

\begin{proof}
Se tiene que $\deg\big(\big\{\rho_K\big(\idel\alpha\big)|\idel\alpha\in C_K
\big\}\big)={\ma Z}\neq \hat{\ma Z}$, por lo que $\rho_K$ no es
suprayectiva.
$\fin$
\end{proof}

Sea $K/\F$ un campo global de funciones congruente sobre $\F$.
Se tiene $\overline \F=\sep \F=\abe \F=\bigcup_{n=1}^{\infty} {\ma F}_{
q^n}$ una cerradura algebraica de $\F$. Se tiene que $K\abe \F$
es la m\'axima extensi\'on de constantes de $K$ y adem\'as $K
\abe\F\subseteq \abe K$. Adem\'as
\begin{align*}
\Gal(K\abe \F/K)&\cong \Gal(\abe\F/\F)=\Gal\big(\bigcup_{n=1}^{\infty}
{\ma F}_{q^n}/\F\big)=\Gal\big(\lim_{\substack{\lra\\ n}}{\ma F}_{q^n}/\F\big)\\
&\cong \lim_{\substack{\longleftarrow\\ n}}\Gal\big({\ma F}_{q^n}/\F)\cong
 \lim_{\substack{\longleftarrow\\ n}}\big({\ma Z}/n{\ma Z}\big)\cong
 \hat{\ma Z},
\end{align*}
la completaci\'on ${\ma Z}$.

Adem\'as, $C_K\cong C_{K,0}\times {\ma Z}\stackrel{\rho_K}{\lra} \Gal(\abe K/K)
\stackrel{\rest}{\lra}\Gal(K{\ma F}_{q^n}/K)$, $\sigma\longmapsto
\sigma|_{K{\ma F}_{q^n}}$ y se tiene $\ker \big(\rest\circ \rho_K\big)=
C_{K,0}\times n{\ma Z}$ (Teorema \ref{T17.6.154N}).

Tenemos el isomorfismo inducido,
\begin{gather*}
\theta:\frac{C_K}{C_{K,0}\times n{\ma Z}}\stackrel{\cong}{\lra}
\Gal(K{\ma F}_{q^n}/K)\cong \Gal({\ma F}_{q^n}/\F)\cong \frac{\ma Z}{n{\ma Z}}.
\intertext{Se sigue que}
\tilde\rho_K:\frac{C_K}{C_{K,0}}\lra \Gal(K\abe\F/K)\cong \hat{\ma Z}.
\end{gather*}
En particular obtenemos nuevamente que $\rho_K$ no es suprayectiva y que
$C_{K,0}=\ker\rho_K|_{\ma Z}$. M\'as precisamente, si $\idel\alpha\in C_{K,0}$,
se tiene $\rho_K\big(\idel\alpha\big)|_{K\abe\F}=\big(\idel\alpha,K\big)|_{K\abe\F}
=1$, por lo que $\big(\idel\alpha,K\big)\in\Gal(\abe K/K\abe\F)$.

Ahora bien, $\rho_K(C_K)$ es denso en $\abe 
G_K=\Gal(\abe K/K)$ y $C_{K,0}$ es compacto por lo que $\rho_K
(C_{K,0})$ es denso y compacto en $\Gal(\abe K/K\abe\F)$, por tanto es
suprayectiva y como $\rho_K$ es inyectiva obtenemos que
\[
\rho_K|_{C_{K,0}}:C_{K,0}\lra \Gal(\abe K/K\abe\F)
\]
es un isomorfismo algebraico y topol\'ogico.

Sea $\pi:J_K\lra C_K=J_K/\*K$ el epimorfismo natural. Se tiene el siguiente
diagrama conmutativo con filas exactas
\[
\xymatrix{
1\ar@{->}[r]&J_{K,0}\ar@{^{(}->}[r]\ar@{->}[d]_{\pi}
\ar@{}[dr]|{\circlearrowright}&J_K\ar@{->>}[r]^{\deg}
\ar@{->}[d]^{\pi}\ar@{}[dr]|{\circlearrowright}
&{\ma Z}\ar@{->}[r]\ar@{->}[d]^{\Id}&0\\
1\ar@{->}[r]&C_{K,0}\ar@{^{(}->}[r]&C_K
\ar@{->>}[r]^{\deg}&{\ma Z}\ar@{->}[r]&0
}
\]

Sea $\vec\alpha_1\in J_K$ es de grado $1$, entonces para $\vec\alpha\in J_K$
arbitrario, el mapeo $\vec\alpha\lra \big(\vec\alpha_1^{-\deg\vec\alpha}\cdot 
\vec\alpha, \vec\alpha_1^{\deg\vec\alpha}\big)$ induce isomorfismos
\begin{gather*}
J_K\cong J_{K,0}\times \langle\vec\alpha_1\rangle\quad\text{y}\quad
C_K\cong C_{K,0}\times \langle\idel\alpha_1\rangle.
\end{gather*}

Como consecuencia del Teorema \ref{T17.6.175N}, se obtiene que
$\rho_K$ nos da un isomorfismo topol\'ogico:

\begin{teorema}\label{CClaseT4.8.1}\label{T17.6.179N}
Sean ${\mathcal G}=\Gal(\abe K/K)=\abe G_K$,
${\mathcal H}_0=\Gal(\abe K/K\abe\F)$ y ${\mathcal G}_0=\Gal(
K\abe\F/K)\cong \hat{\ma Z}$. Entonces
${\mc G}/{\mc H}_0\cong {\mc G}_0$ y $\rho_K$ es un isomorfismo
algebraico y topol\'ogico de $C_{K,0}$ sobre ${\mathcal H}_0$, es decir,
$C_{K,0}\isomorfo\limits^{\rho_K}\Gal(\abe K/K\abe \F)$ 
y ${\mathcal G}$ es el producto directo de ${\mathcal H}_0$
y ${\mathcal G}_0$,
\begin{gather*}
{\mathcal G}\cong {\mathcal H}_0\times {\mathcal G}_0,
\quad \Gal(\abe K/K)\cong \Gal(\abe K/K\abe\F)\times \Gal(K\abe\F/K).
\\
\xymatrix{
\hat{K}\ar@{-}[rr]^{\hat{\ma Z}}_{{\mathcal G}_0}
\ar@{-}[d]_{\cong C_{K,0}}&&\abe K\ar@{-}[dll]^{\mathcal G}
\ar@{-}[d]\ar@{-}@/^1pc/[d]^{{\mathcal H}_0\cong C_{K,0}}\\
K\ar@{-}[rr]_{\hat{\ma Z}}&&K\abe\F}
\end{gather*}

El campo de constantes de $\hat{K}:=\big(\abe K\big)^{{\mc G}_0}$
es $\F$.
\end{teorema}

\begin{proof}
Sea ${\ma F}_{q^n}$ el campo de
constantes del campo $\hat{K}$. Se tiene que
$K=\hat{K}\cap K\abe\F\supseteq K{\ma F}_{q^n}\cap K\abe\F
\supseteq {\ma F}_{q^n}$, por lo que ${\ma F}_{q^n}\subseteq K$
de donde se sigue que $n=1$ y que $\F$ es el campo de constantes
de $\hat{K}$.

Se tiene que ${\mc H}_0\cong C_{K,0}$
tanto algebraica como topol\'ogicamente.
Ahora bien, puesto que ${\mc G}_0\cong \hat{\ma Z}$, la sucesi\'on
\[
1\lra {\mc H}_0\lra {\mc G}\stackrel{\pi}{\lra} {\mc G}_0\lra 1
\]
se escinde pues $\hat{\ma Z}=\overline{\langle\Fro K\rangle}$ 
es la cerradura de ${\ma Z}$ en ${\mc G}_0$.
M\'as precisamente, si $x\in {\mc G}$ es tal que $\pi(x)=
\Fro K$, entonces $\varphi(\Fro K)=x$ es el mapeo de escisi\'on,
extendi\'endolo de manera continua a todo $\hat{\ma Z}$.
$\fin$
\end{proof}

Detallando m\'as el Teorema \ref{T17.6.179N},
el siguiente resultado nos da la informaci\'on
relevante acerca de las extensiones de constantes.

\begin{teorema}\label{CClaseT4.2.10} Existen dos diagramas conmutativos
de sucesiones exactas. En cada diagrama, la flecha vertical de la izquierda es 
un isomorfismo de grupos topol\'ogicos.
\begin{gather*}
\xymatrix{
1\ar[r]&C_{K,0}\ar[r]^i\ar[d]_{\cong}^{\rho_K^{\prime}=
\rho_K|_{C_{K,0}}}&C_K\ar[r]^{\deg}\ar[d]_{\rho_K}&
{\ma Z}\ar[r]\ar[d]_{\rho_{\F}}&0\\
1\ar[r]& \Gal(\abe K/K \abe\F)\ar[r]_i&\Gal(\abe K/K)\ar[r]_{\rest}
&{\underbracket[0pt]{\Gal(K\abe\F/K)}_{\substack{\uigual\\
 \Gal(\abe\F/\F)\cong \hat{\ma Z}}}}\ar[r]&1}\\
\xymatrix{
1\ar[r]&I_{K,0}\ar[r]^i\ar[d]_{\cong}&I_K\ar[r]^{\deg}\ar[d]&{\ma Z}
\ar[r]\ar[d]^{\rho_{\F}}&0\\
1\ar[r]&\Gal(K^{\rm{nr}} /K \abe\F)\ar[r]_i&\Gal(K^{\rm{nr}} /K)\ar[r]_{\rest}&
\Gal(K\abe\F/K)\ar[r]&1}
\end{gather*}

Aqu\'i se tiene $\rho_{\F}\colon{\ma Z}\lra \Gal(\abe K/K)$ definida
por $\rho_{\F}(n)=\Fro K^n$, donde $\Fro K$ denota al automorfismo
de Frobenius.
\end{teorema}

\begin{proof}  En el primer diagrama, 
$(\rho_{\F}\circ \deg)(\idel x)=\rho_{\F}(
\deg \idel x)=\Fro K^{\deg \idel x}$.
Por el otro lado, $\rest \rho_K(\idel x)
=\rho_{K\abe\F}(\vec x)=\Fro K^{\deg\idel x}$.

El lado izquierdo del diagrama es conmutativo, es decir $i\circ \rho_K^{\prime}=
\rho_K\circ i$ y $\ker (i\circ \rho_K^{\prime})=\ker \rho_K^{\prime}=\ker (
\rho_K\circ i)=0$ pues se tiene que $\rho_K$ es inyectiva 
en el caso de campos de 
funciones. Por lo tanto $\rho_K^{\prime}$ es una inyecci\'on y tenemos que
$\im \rho_K^{\prime}=\ker \rest=\Gal(\abe K/K\abe \F)$.

Ahora bien, por el Lema de la Serpiente (Teorema \ref{CClaseT1.5.2}), se tiene la sucesi\'on
exacta 
\begin{gather*}
1\to\ker \rho_K^{\prime}\to\ker \rho_K\to \ker \rho_{\F}\to\\
\to \coker \rho_K^{\prime}\to \coker \rho_K\to\coker \rho_{\F}\to 1.
\end{gather*}

Usando que $\rho_K$ es inyectiva en el caso de campos de funciones, que
$\coker \rho_K=\hat{\ma Z}/{\ma Z}$, que $\rho_{\F}$ es inyectiva y
que $\coker \rho_{\F}=\hat{\ma Z}/{\ma Z}$, se sigue
\[
1\to\ker \rho_K^{\prime}\to 1\to 1\to \coker \rho_K^{\prime}\to \hat{\ma Z}/{\ma Z}
\to \hat{\ma Z}/{\ma Z}\to 0
\]
es exacta por lo $\ker \rho_K^{\prime}=1$ y $\coker \rho_K^{\prime}=1$ de donde
se sigue que $\rho_K^{\prime}$ es un isomorfismo.

El segundo diagrama, proviene del primero al dividir la primera sucesi\'on entre
$\bar{U}=U\*K/\*K$, con $U=\prod_{\pK\in{\ma P}_K} U_{\pK}$:
$C_{K,0}/\bar{U}=I_{K,0}$, $C_K/\bar{U}=I_K$ y usando que
$\rho_K(\bar{U})=\Gal(\abe K/K^{\rm{nr}})$
pues
$C_K/{\bar{U}}\cong \Gal(K^{\rm{nr}}/K)$ seg\'un el TCCG. $\fin$
\end{proof}

Notemos que, en particular, se tiene que $\bar{U}$ corresponde a la 
m\'axima extensi\'on
no ramificada abeliana de $K$. 

\begin{corolario}\label{CClaseC4.2.11} Se tiene
\[
C_{K,0}\cong \Gal(\abe K/K\abe \F),\quad I_{K,0}\cong \Gal(K^{\rm{nr}}/K\abe\F).
\quad 
\xymatrix{
\abe K\ar@{-}[d]\ar@{-}@/^1pc/[d]^{\bar{U}}\ar@{-}@/^4pc/[dd]^{C_{K,0}}\\
K^{\rm{nr}}\ar@{-}[d]_{\text{finita}}\ar@{-}@/^1pc/[d]^{I_{K,0}}\\
K\abe\F\ar@{-}[d]\\ K^{\ast}}
\]
En particular $K^{\rm{nr}}/K\abe\F$ es una extensi\'on finita. $\fin$
\end{corolario}

En adici\'on tenemos $C_K=\langle\idel\alpha_1\rangle\times
C_{K,0}$, con $\vec\alpha\in J_K$ tal que
$|\idel\alpha_1|=q^{-1}$ o $\deg \idel\alpha_1=1$ y
$C_K$ es isomorfo a ${\ma Z}\times C_{K,0}$ pero no de manera
can\'onica, esto es, $C_{K,0}$ es \'unico no as{\'\i} la copia de
${\ma Z}$ que podemos seleccionar dentro de $C_K$ para obtener
el producto directo. Notemos que tenemos una copia de ${\ma Z}$
en $C_K$ para cada $\idel \alpha_1\in C_K$ con 
$\deg \idel \alpha_1=1$.

\begin{definicion}\label{D17.6.180N}
Se define otra topolog{\'\i}a en $C_K$: las vecindades de $1$ son los
subgrupos abiertos de {\'\i}ndice finito en la topolog{\'\i}a usual de 
$C_K$. 
Esta nueva topolog{\'\i}a se llama la {\em topolog{\'\i}a
de clase\index{topolog{\'\i}a de clase.}} 
\end{definicion}

La topolog\'ia de clase
es estrictamente menos fina que la topolog\'ia usual, como veremos
en el Teorema \ref{CClaseT4.8.2}.

\begin{teorema}\label{CClaseT4.8.2}
La topolog{\'\i}a de clase coincide con la original en $C_{K,0}$, induce
la {\em topolog{\'\i}a de ideales\index{topolog{\'\i}a de ideales}} sobre
${\ma Z}$ y la topolog{\'\i}a producto sobre ${\ma Z}\times C_{K,0}$. 
La topolog{\'\i}a de ideales de ${\ma Z}$ consiste de la topolog{\'\i}a generada
por las vecindades de $0$ por los ideales $\{n{\ma Z}\}_{n\neq 0}$. Es
decir, $V$ es una vecindad de $0$ si contiene a alg\'un ideal $n{\ma Z}$ con
$n\neq 0$. La completaci\'on de ${\ma Z}$ con respecto a esta 
topolog{\'\i}a es el anillo de Pr\"ufer:
\[
\hat{\ma Z}=\lim_{\substack{\longleftarrow\\ n}}{\ma Z}/n{\ma Z}.
\]
\end{teorema}

\begin{proof} 
Sea $B$ un subgrupo abierto de \'indice finito en $C_K$. Sea $\idel\alpha\in
B$ con $\deg\idel\alpha=\min\{n\in{\ma N}\mid \text{existe $\idel\beta\in B$
con $\deg\idel\beta=n$}\}$. Sea $B_0=B\cap C_{K,0}$ y $B=\bigcup_{n\in{\ma Z}}
\idel\alpha^nB_0$. Ahora bien, $B_0$ es un conjunto abierto de $C_K$ pues
$C_{K,0}$ es abierto en $C_K$ (Proposici\'on \ref{P17.6.34-1N}).

Veremos que $B_0$ es de \'indice finito en $C_{K,0}$. Sea ${\eu m}=\prod_{\pK}
\pK^{\gamma_{\pK}}$ con $\gamma_{\pK}\geq 0$ para toda $\pK$ y $\gamma_{
\pK}=0$ para casi toda $\pK$. ${\eu m}$ recibe el nombre de {\em
modulus\index{modulus}}. Un sistema fundamental de vecindades de $1$ en
$J_K$ est\'a dado por los grupos
\[
J_K^{\eu m}=\prod_{\pK\nmid {\eu m}}U_{\pK}\times \prod_{\pK|{\eu m}}U_{
\pK}^{(\gamma_{\pK})},
\]
donde $U_{\pK}^{(\gamma_{\pK})}$ son los grupos consistentes de los elementos 
$\alpha_{\pK}\in\*K_{\pK}$ tales que $\alpha_{\pK}\equiv 1\bmod\pK^{\gamma_{\pK}}$
para $\pK|{\eu m}$. Se tiene que $U_{\pK}^{(\gamma_{\pK})}\subseteq U_{\pK}$
y que si $\vec\alpha\in J_K^{\eu m}$, entonces $\deg\vec\alpha=0$. Sea 
$C_L^{\eu m}:=J_K^{\eu m}\*K/\*K$, los cuales forman un sistema fundamental de
vecindades de $1\in C_K$ y $C_K^{\eu m}\subseteq C_{K,0}$. Se tiene
\[
[C_{K,0}:C_K^{\eu m}]=[C_{K,0}:U_K\*K/\*K][U_K\*K/\*K:C_K^{\eu m}],
\]
donde $U_K=\prod_{\pK}U_{\pK}$. Ahora bien, $C_{K,0}/\big(U_K\*K/\*K\big)
\cong I_{K,0}$, el cual es un grupo finito.

Por otro lado, $\big[U_{\pK}:U_{\pK}^{(\gamma_{\pK})}\big]<\infty$, por lo que
\[
[U_K\*K/\*K:C_K^{\eu m}]\leq \prod_{\pK|{\eu m}}\big[U_{\pK}:U_{\pK}^{
(\gamma_{\pK})}\big]<\infty.
\]

El subgrupo abierto $B_0$ debe contener a alg\'un subgrupo abierto $C_K^{\eu m}$
por lo que $[C_{K,0}:B_0]<\infty$. Por otro lado, $B\cap {\ma Z}$ es un subgrupo
de \'indice finito en ${\ma Z}$ por lo que $B\cap {\ma Z}=n{\ma Z}$ para alg\'un
$n\in{\ma N}$. Por tanto, la topolog\'ia inducida en ${\ma Z}$ es la topolog\'ia de
ideales.
$\fin$
\end{proof}

\begin{definicion}\label{D17.6.182N}
Se define $\hat C_K:=C_{K,0}\times \big(\idel\alpha_1\big)^{\hat{\ma Z}}
\cong C_{K,0}\times \hat {\ma Z}$, es decir, $\hat C_K\cong
C_{K,0}\times \big\{\idel\alpha_1^{\mu}\big\}_{\mu\in\hat{\ma Z}}$.
\end{definicion}

\begin{observacion}\label{D17.6.183N-1}
Se tiene que
\[
\hat C_K\cong\lim_{\substack{\longleftarrow\\ B}}C_K/B,
\]
donde $B$ recorre los subgrupos abiertos de \'indice finito de $C_K$.

Es decir, $\hat C_K$ es la completaci\'on en la topolog\'ia de Krull
de $C_K$.
\end{observacion}

\begin{observacion}\label{D17.6.183N}
Puesto que tanto $C_{K,0}$ como $\hat{\ma Z}$ son compactos, se sigue
que $\hat C_K$ es compacto.
\end{observacion}

Se tiene que $\rho_K:C_K\lra \Gal(\abe K/K)$ es uniformemente continua
en la nueva topolog\'ia, de hecho, si $U$ es un abierto de $G:=\Gal(\abe K/K)$,
existe un subgrupo abierto $H$ de $G$ con $H\subseteq U$
pues los subgrupos abiertos definen la topolog\'ia de Krull. Entonces
$H$ es de \'indice finito. Sea $L:=(\abe K)^H$ el campo fijo bajo $H$ y se
tienen $\Gal(L/K)\cong G/H$. De la ley de reciprocidad, tenemos
\[
\rho_K|_L=\psi_{L/K}: C_K/\N_{L/K}(C_L)\stackrel{\cong}{\lra} G/H.
\]
Ahora bien, $\N_{L/K}C_L$ es abierto en $C_K$ y si $\idel\alpha\idel\beta^{-1}
\in \N_{L/K}C_L$, entonces $\rho_K\big(\idel\alpha\idel\beta^{-1}\big)=
\rho_K(\idel\alpha)\rho_K(\idel\beta)^{-1}\in H$, de donde obtenemos que $\rho_K$
es uniformemente continua.

Por tanto, $\rho_K$ se extiende de manera \'unica a la completaci\'on,
en la nueva topolog\'ia, $\hat C_K$
de $C_K$ y $\hat C_K\cong C_{K,0}\times \hat{\ma Z}$ el cual es compacto.
Por tanto $\rho_K\big(\hat C_K\big)$ es denso y compacto, de donde se sigue
que $\rho_K\big(\hat C_K\big)\cong G$. Puesto que $\rho_K$ es inyectiva, se
tiene $\rho_K: \hat C_K\lra G$ es un isomorfismo tanto algebraico como topol\'ogico.

Con todos estos elementos, se demuestra el Teorema de Existencia
para campos de funciones. De hecho se hace m\'as, se da una
correspondencia biyectiva entre subgrupos cerrados de 
$\hat{C}_K$ y los subgrupos cerrados de ${\mc G}=\Gal(\abe K/K)$.

\begin{teorema}[Teorema de Existencia para campo de 
funciones\index{teorema de existencia para campos de 
funciones}]\label{CClaseT4.8.3}\label{T17.6.184N}
Sea $K$ un campo global de funciones. Sea $H$ un subgrupo abierto
de {\'\i}ndice finito en $C_K$. Entonces existe una extensi\'on
abeliana finita $L$ de $K$ tal que $H={\mc N}_L
=\N_{L/K}C_L$ M\'as a\'un, $L$
es el campo fijo de $\rho_K(H)$: $L=(\abe K)^{\rho_K(H)}$ donde
$\rho_K=(\underline{\ \ }, K)\colon C_K\longrightarrow\Gal(\abe K/K)$
es el s{\'\i}mbolo de la norma residual u homomorfismo de Artin.
\end{teorema}

\begin{proof}  
Sea $H$ un subgrupo abierto de {\'\i}ndice finito de $C_K$.
Sea $H_0=H\cap C_{K,0}$. Puesto que
\[
\frac{H}{H_0}\cong \frac{H}{(H\cap C_{K,0})}
\cong \frac{HC_{K,0}}{C_{K,0}}
\subseteq \frac{C_K}{C_{K,0}}\cong {\ma Z},
\]
se tiene que $H=\cup_{n\in{\ma Z}} h^n H_0$ para
alg\'un $h\in H$.
Por otro lado, debido a que tenemos el isomorfismo
$C_K=C_{K,0}\times {\ma Z}=C_{K,0}\times
\langle\tilde{\vec\alpha}_1\rangle$ con $|\tilde{\vec \alpha}_1|=q^{-1}$,
podemos tomar $h=\tilde{\vec \alpha}_1^d$ para alguna $d$.
Se tiene que $\hat{\ma Z}/d\hat{\ma Z}\cong {\ma Z}/d{\ma Z}$ de
manera natural ($\xymatrix{
{\ma Z}\ar@{^{(}->}[r]\ar@/_1pc/[rr]^{\varphi}&\hat{\ma Z}
\ar@{>>}[r]^{\bmod d\quad\ }&\hat{\ma Z}/d\hat{\ma Z}}$, $\ker \varphi =d{\ma Z}$
y se prueba que $\varphi$ es un epimorfismo). 

Sea $\hat H:=H_0\times d\hat{\ma Z}$, esto es, $\hat H$
es la completaci\'on de $H$ en la nueva topolog\'ia. Entonces 
$\hat C_K=C_{K,0}\times \hat{\ma Z}$, por lo que
\[
\frac{\hat C_K}{\hat H}=\frac{C_{K,0}\times \hat{\ma Z}}
{H_0\times d\hat{\ma Z}}\cong\frac{C_{K,0}}{H_0}\times
 \frac{\ma Z}{d{\ma Z}}\cong
\frac{C_{K,0}\times {\ma Z}}{H_0\times d{\ma Z}}\cong\frac{C_K}{H}
\]
lo cual implica $[C_K:H]=d[C_{K,0}:H_0]=
[\hat C_K:\hat{H}]$.

Sea ${\mathcal H}=\rho_K(\hat{H})\subseteq \Gal(\abe K/K)
= {\mc G}$. Se
tiene que ${\mathcal H}$ es cerrado en ${\mathcal G}$ y 
pues que $\rho_K$ es un isomorfismo
\[
[{\mathcal G}:{\mathcal H}]=[\rho_K(\hat{C}_K):\rho_K(\hat{H})]
\underbracket[0pt]{=}_{\substack{\uparrow\\ \rho_K\text{\ es un}\\ {\text{isomorfismo}}}}
[\hat{C}_K:\hat{H}]=[C_K:H],
\]
es decir, $[C_K:H]=[{\mathcal G}:{\mathcal H}]$. 

Sea $L$ el campo fijo de $\abe K$ por
${\mathcal H}=\rho_K(\hat{H})$, es decir, $L:=(\abe K)^{\mathcal H}$.
Entonces $[L:K]=[{\mathcal G}:{\mathcal H}]=
[C_K:H]$. 

Adem\'as, $\Gal(L/K)\cong \frac{\Gal(\abe K/K)}{\Gal(\abe K/L)}\cong
\frac{{\mathcal G}}{{\mathcal H}}\cong \frac{C_K}{\N_{L/K} C_L}$
y $\rho_K(C_K)\subseteq {\mathcal G}$, por lo que
 $\rho_K(\N_{L/K} C_L)\subseteq
{\mathcal H}$ lo cual implica que $\N_{L/K}C_L\subseteq H$.
Finalmente, puesto que $[C_K:\N_{L/K}C_L]=[L:K]=[C_K:H]$,
se sigue que ${\mc N}_L=\N_{L/K}C_L=H$.
$\fin$
\end{proof}

\subsubsection{Extensiones geom\'etricas y de constantes
de campos de funciones}\label{CClaseS4.8.0}

Analicemos con m\'as detalle parte de la discusi\'on anterior.
Sea $\nr K$ la m\'axima extensi\'on no ramificada de $K$ contenida en 
$\abe K$, donde $K$ es un campo global de funciones.
Se tiene que $K\abe\F\subseteq \nr K$.
Consideremos el siguiente diagrama
\[
\begin{minipage}{5cm}
\xymatrix{
&\abe K\ar@{-}[ddl]_{\mathcal G}\ar@{-}[d]\ar@{-}@/^2pc/[dd]^{{\mathcal H}_0}\\
&K^{\rm{nr}}\ar@{-}[d]\\ K\ar@{-}[r]_{{\mathcal G}_0}&K\abe \F}
\end{minipage}
\qquad
\begin{minipage}{6cm}
donde ${\mathcal H}_0=\Gal(\abe K/K\abe \F)\cong C_{K,0}$,
${\mc G}_0\cong {\mc G}/{\mc H}_0
\cong \Fro K^{\hat{\ma Z}}$, ${\mathcal
G}=\Gal(\abe K/K)$.
\end{minipage}
\]

Sea $\sigma\in {\mathcal G}$ con $\sigma|_{K\abe \F}=
\Fro K^{\ord \sigma}$. En general tenemos $\sigma|_{K\abe \F}
\in{\mc G}_0\cong\overline{\langle\Fro K\rangle}=\Fro K^{\hat{\ma Z}}$.

Sea $\rho_K\colon C_K \longrightarrow {\mathcal G}$ el mapeo de reciprocidad.
Se tiene $\ord \rho_K\big(\idel\alpha\big)=\deg \idel\alpha$.
En particular se tiene que $\ord\rho_K\big(\idel\alpha\big )\in{\ma Z}$
y no en $\hat{\ma Z}$. Adem\'as ${\mathcal H}_0\cong C_{K,0}$ bajo el
mapeo $\rho_K$, ${\mathcal G}\cong {\mathcal H}_0\times {\mathcal G}_0$,
${\mathcal G}_0=\Gal(K\abe\F/K)$ y la sucesi\'on exacta
\[
1\longrightarrow{\mathcal H}_0\longrightarrow {\mathcal G}\longrightarrow
{\mathcal G}_0\longrightarrow 1
\]
se escinde con $\Fro K^{\hat{{\ma Z}}}\cong {\mc G}_0\lra{\mc G}$,
 $\Fro K\mapsto \sigma$ con $\ord \sigma=1$, esto es,
$\sigma|_{K\abe \F}=\Fro K$ y cada uno de estos se puede obtener mediante
$\vec \alpha\in J_K$ con $\deg\vec \alpha=1$, $\rho_K\big(\idel\alpha\big)
=\sigma$.

Para cada $\vec\alpha\in J_K$ tal que $\deg \vec\alpha=1$
definimos un campo
$\tilde{K}(\vec\alpha)$ con $\tilde{K}(\vec\alpha):=(\abe K)^{{\mathcal G}_0
(\vec \alpha)}$ donde ${\mathcal G}_0(\vec\alpha)\subseteq {\mathcal G}$
es la escisi\'on $\psi\colon{\mathcal G}_0\longrightarrow {\mathcal G}_0
(\vec\alpha)\subseteq {\mathcal G}$, $\Fro K\longmapsto \rho_K(\vec\alpha)$.

Se tiene que $\tilde{K}(\vec\alpha)\cap K\abe\F=K$ y por tanto el campo
de constantes de $\tilde{K}(\vec\alpha)$ es $\F$. Con el isomorfismo
$\hat C_K=C_{K,0}\times \idel\alpha^{\hat{\ma Z}}$ con
$\deg\vec \alpha=1$, $\rho_K\colon \hat C_K\longrightarrow {\mathcal G}=
\Gal(\abe K/K)$, $\rho_K(C_{K,0})={\mathcal H}_0$, $\rho_K(\idel\alpha)\in
{\mathcal G}$, tenemos los diagramas
\begin{gather*}
\xymatrix{
\tilde{K}(\vec\alpha)\ar@{-}[rr]^{{\mathcal G}_0
(\vec\alpha)\cong\hat{\ma Z}}\ar@{-}[d]_{
{\mathcal G}/{\mathcal G}_0(\vec \alpha)\cong C_{K,0}}&&
\abe K\ar@{-}[dll]_{\mathcal G}
\ar@{-}[d]^{C_{K,0}\cong {\mathcal H}_0}\\
K\ar@{-}[rr]_{\hat{\ma Z}\cong \Fro K^{\hat{\ma Z}}}&&K\abe \F}\\
\xymatrix{
\tilde{K}(\vec\alpha)\ar@{-}[rr]\ar@{-}[d]_{\bar U=}^{C_K^1}&&
\abe K\ar@{-}[d]^{\Gal(\abe K/\nr K)}_{\bar U=C_K^1\cong}
\ar@{-}@/^7pc/[dd]^{C_{K,0}}\\
K^{\rm{nr}}\cap \tilde{K}(\vec\alpha)\ar@{-}[rr]\ar@{-}[d]^{I_{K,0}}&&
K^{\rm{nr}}\ar@{-}[d]_{\frac{C_{K,0}}{\bar U}\cong}^{I_{K,0}}\\
K\ar@{-}[rr]&& K\abe \F
}
\end{gather*}
donde $\bar U=\frac{U\*K}{\*K}$ y donde el grupo de normas de $\nr K
\cap \tilde K(\vec\alpha)$ es $\bar U\times {\ma Z}$ y 
\[
\frac{C_K}{\bar U
\times {\ma Z}}\cong \frac{C_{K,0}\times {\ma Z}}{\bar U\times {\ma Z}}
\cong \frac{C_{K,0}}{\bar U}\cong I_{K,0}.
\]

Todas las extensiones de $K$ contenidas en $\tilde{K}(\vec\alpha)$ tienen
campo de constantes $\F$ y la m\'axima extensi\'on de $K$ no ramificada
dentro de $\tilde{K}(\vec\alpha)$ es $K^{\rm{nr}}\cap \tilde{K}(\vec\alpha)$ y
\[
\Gal(K^{\rm{nr}} \cap \tilde{K}(\vec\alpha)/K)\cong \Gal(K^{\rm{nr}}/K\abe \F)\cong I_{K,0}.
\]

Notemos que ${\mc G}_0(\vec\alpha)\cong \vec \alpha^{\hat{\ma Z}}\cong
\hat{\ma Z}$, por lo que $K^{\rm{nr}}\cap \tilde{K}(\vec\alpha)
=\tilde{K}(\vec\alpha)^{\text{nr}}$
corresponde a $ C_K^1\times {\mc G}_0(\vec\alpha)\cong 
C_K^1\times \vec\alpha^{
\hat{\ma Z}}$, donde $C_K^1=\Gal(\abe K/\nr K)$.
Es decir, $\tilde{K}(\vec\alpha)^{\text{nr}}=
(\abe K)^{C_K^1\times {\mathcal G}_0(\vec\alpha)}$.

Ahora bien, $\tilde{K}(\vec\alpha)\abe\F=\abe K$, 
$\Gal(\abe K/\tilde{K}(\vec\alpha))
\cong \hat{\ma Z}$ y $\abe K/\tilde{K}(\vec\alpha)$ es una extensi\'on de
constantes. Si $K{\ma F}_{q^d}$ es la extensi\'on de constantes de 
grado $d$ sobre $K$, entonces $\Gal(K{\ma F}_{q^d}/K)=\langle
\bar\Frobenius_K\rangle$ donde $\bar\Frobenius_K:=
\Fro K \bmod\ d\hat{\ma Z}$, es decir
\begin{gather*}
\xymatrix{
K\ar@{-}[rr]_{\hat{\ma Z}/d\hat{\ma Z}\cong \langle\bar\Frobenius_K\rangle}
&&K{\ma F}_{q^d}\ar@{-}[rr]_{d\hat{\ma Z}}
&&K\abe\F\ar@{-}@/_2pc/[llll]_{\hat{\ma Z}}}
\\
\xymatrix{
\tilde{K}(\idel\alpha)\ar@{-}[r]\ar@{-}[d]&\tilde{K}(\idel\alpha){\ma F}_{q^d}
\ar@{-}[rrr]^{\langle\idel\alpha^d\rangle}\ar@{-}[d]&&&\abe K\ar@{-}[d]^{C_{K,0}}
\ar@{-}[dlll]_{C_{K,0}\times\langle \idel \alpha^d\rangle}\\
K\ar@{-}[r]_{\bar\Frobenius_K}&K{\ma F}_{q^d}
\ar@{-}[rrr]_{\overline{\langle\Fro K^d\rangle}\subseteq
\hat{\ma Z}}&&&K\abe\F}
\end{gather*}
esto es, la extensi\'on de constantes de grado $d$ de $K$ corresponde
al grupo $C_{K,0}\times\langle \tilde{\vec\alpha}^d
\rangle\subseteq C_K$ (lo cual ya se hab\'ia obtenido en el
Teorema \ref{T17.6.154N}) y la extensi\'on de
constantes de grado $d$ de $\tilde{K}(\vec\alpha)$ corresponde al grupo
$\langle\tilde{\vec\alpha}^d\rangle$. Es decir $C_{K,0}
\times \langle\tilde{\vec \alpha}^d
\rangle$ es el grupo de normas de $K{\ma F}_{q^d}$.
Aqu\'i estamos identificando, bajo el mapeo de reciprocidad $\rho_K$,
$\idel \alpha$ con $\Fro K$ y $C_{K,0}$ con ${\mc H}_0$.

Sea ahora $L/K$ una extensi\'on abeliana finita y geom\'etrica, es decir,
el campo de constantes de $L$ es $\F$. Entonces $L\cap K\abe\F=K$.
\[
\xymatrix{
&\abe K\ar@{-}[dl]_{\Gal(\abe K/L)}\ar@{-}[d]^{{\mathcal H}}\\
L\ar@{-}[r]^{\hat{\ma Z}}\ar@{-}[d]&L\abe \F\ar@{-}[d]\\
K\ar@{-}[r]_{\hat{\ma Z}}&K\abe \F}
\]
Sea ${\mathcal H}=\Gal(\abe K/L\abe\F)$ y $\Gal(\abe K/L)\cong
{\mathcal H}\times\hat{\ma Z}$ como antes y $\hat{\ma Z}$ se obtiene
como cualquier escisi\'on:
\[
\xymatrix{
1\ar[r] &\underbracket[0pt]{\mathcal H}_{\substack{\uigual\\ \Gal(\abe K/L\abe\F)}}
\ar[r]&\underbracket[0pt]{\Gal(\abe K/L)}_{\vec\alpha}\ar[r]&\underbracket[0pt]{\Gal(
L\abe \F/L)}_{\Fro K}\cong\hat{\ma Z}\ar[r]\ar@/^2pc/[l]\ar@/_1.5pc/[l]&0}
\]
y $\deg\vec\alpha=1$. 
Por tanto $\hat{\ma Z}\cong\langle \vec\alpha^{\hat{\ma Z}}
\rangle$ y $(\abe K)^{\hat{\ma Z}}=\tilde{K}(\vec\alpha)$ de antes. Es decir,
$L\subseteq \tilde{K}(\vec\alpha)$ para alg\'un $\vec\alpha$.
\[
\xymatrix{
\tilde{K}(\idel\alpha)\ar@{-}[r]^{\idel\alpha^{\hat{\ma Z}}}\ar@{-}[d]_{\mathcal H}&
\abe K\ar@{-}[d]^{\mathcal H}\\
L\ar@{-}[d]\ar@{-}[r]^{\idel\alpha^{\hat{\ma Z}}}&L\abe\F\ar@{-}[d]\\
K\ar@{-}[r]_{\Fro K^{\hat{\ma Z}}}&K\abe \F}
\]

Es decir, toda extensi\'on abeliana finita geom\'etrica est\'a contenida
en alguna de las extensiones $\tilde{K}(\vec\alpha)$.

\subsubsection{Sobre los campos de constantes}\label{CClaseS4.8.1-1}

\begin{teorema}\label{CClaseT4.8.1-1.1} Sea $B<C_K$ un subgrupo
abierto de \'indice finito. Sea $d:=\min\{n\in{\ma N}\mid \text{existe
$\tilde{\vec b}\in B$ con $\deg\tilde{\vec b}=n$}\}$. Entonces si
$E$ es el campo asociado a $B$, es decir $B=\N_{E/K} C_E$, 
el campo de constantes de $E$ es ${\ma F}_{q^d}$.
\end{teorema}

\begin{proof}
Primero notemos que para cualquier campo $F$, el campo de 
constantes de $F$ es ${\ma F}_{q^r}$ donde $F\cap K\abe \F=K
{\ma F}_{q^r}$.

Consideremos $\tilde{\vec b}\in B$ con $\deg\tilde{\vec b}=d$. Sea
$B_0:=B\cap C_{K,0}$. Se tiene
\[
\frac{B}{B_0}=\frac{B}{B\cap C_{K,0}}\cong \frac{BC_{K,0}}{C_{K,0}}
\subseteq \frac{C_K}{C_{K.0}}\cong {\ma Z}.
\]

Veamos que $[B:B_0]=\infty$. Esto es claro de lo anterior, pues de 
otra forma, al tener $B/B_0\subseteq {\ma Z}$, necesariamente
$B=B_0\subseteq C_{K,0}$ pero entonces $[C_K:B]\geq
[C_K:C_{K,0}]=\infty$. En resumen $B/B_0\cong {\ma Z}$.
Este isomorfismo viene dado por la sucesi\'on exacta de grupos
\[
\xymatrix{
1\ar[r]& B_0\ar[r]& B\ar[r]^{\deg}& d{\ma Z}
\ar@/^1pc/[l]^{\varphi}\ar[r]& 0,
}
\]
donde $\varphi$ es el mapeo de escisi\'on dado por
$\varphi(d)=\tilde{\vec b}$ y $B\cong \langle\tilde{\vec b}\rangle
\times B_0$.
\[
\xymatrix{
&&&\abe K\ar@{-}[dd]^{B_0}\ar@{-}[lldd]_B
\ar@{-}[dddll]|!{[lldd];[dd]}\hole^{\hspace{-4pt}BC_{K,0}}
\ar@{-}@/^3pc/[ddd]^{C_{K,0}}\\ \\
&E\ar@{-}[d]\ar@{-}[dl]\ar@{-}[rr]_{
\langle\tilde{\vec b}\rangle}
&&E\abe\F\ar@{-}[d]^{\frac{C_{K,0}}{B_0}}\\
K\ar@{-}[r]^d&M\ar@{-}[rr]_{\langle\tilde{\vec b}\rangle}
&&K\abe \F
}
\]

Sea $M$ el campo correspondiente a $BC_{K,0}$. Se tiene
\begin{gather*}
BC_{K,0}=\langle\tilde{\vec b}\rangle B_0 C_{K,0}=
\langle\tilde{\vec b}\rangle C_{K,0}\cong \langle
\tilde{\vec b}\rangle\times C_{K,0}.
\intertext{Se sigue que $M=E\cap K\abe \F$. Puesto que $M\subseteq
K\abe \F$, $M/K$ es una extensi\'on de constantes de grado}
[M:K]=[C_K:BC_{K,0}]=[{\ma Z}\times C_{K,0}:d{\ma Z}
\times C_{K,0}]=[{\ma Z}:d{\ma Z}]=d.
\end{gather*}

Por tanto el campo de constantes de $E$ es ${\ma F}_{q^d}$.
$\fin$
\end{proof}

\subsection{Leyes de descomposici\'on de primos en campos 
globales\index{descomposici\'on en campos globales}}\label{CClaseS4.6}

Consideremos $K$ un campo global arbitrario.

\begin{teorema}\label{CClaseT4.6.1}\label{T17.6.186N}
 Sea $L/K$ una extensi\'on abeliana finita de
campos globales de grado $n$ y sea $\pK$ un primo de $K$
no ramificado en $L$. Sea $\pi\in K_\pK$ un elemento primo,
$v_\pK(\pi)=1$. Sea 
\[
\enc{\pi}{\pK}=(\ldots,1,1,\pi,1,1,\ldots)\in J_K.
\]
Sea $f$ el m{\'\i}nimo entero positivo tal que
$\tilde{\lceil \pi\rceil}_{\pK}^f
\in \N_{L/K}C_L$. Entonces el primo $\pK$ se factoriza en la 
extensi\'on $L$ en $h=n/f$ primos distintos $\pL_1,\ldots,\pL_h$ de
grado relativo $f$, es decir 
\[
\con_{K/L}\pK=\pL_1\cdots\pL_h,\quad [\o_{\pL_i}/\pL_i:\o_\pK/\pK]=f.
\]

Esto es, si conocemos $\N_{L/K} C_L$ podemos determinar
la descomposici\'on de $\pK$ en $L$.
\end{teorema}

\begin{proof} Puesto que $\pK$ es no ramificado, $\pK$ se escribe como
$\pK=\pL_1\cdots \pL_h$ con cada $\pL_i$, $1\leq i\leq h$ de grado
$f^\prime$. Se tiene que $C_K/\N_{L/K}C_L\cong \Gal(L/K)$. Se tiene
$f=o\big(\overline{ 
\tilde{\lceil \pi\rceil}}_{\pK}\big)$ en $C_K/\N_{L/K}C_L$ el cual
es el orden de 
\[
\psi_{L/K}\big(\overline{
\tilde{\lceil \pi\rceil}}_{\pK}\big)=\prod_{{\eu q}\in{\ma P}_K}
\big(\big(\overline{
\tilde{\lceil \pi\rceil}}_{\pK}\big)_{\eu q}, L_{\eu q}/K_{\eu q}\big)=
\big(\pi,L_\pK/K_\pK\big)=\Fro \pK
\]
el automorfismo de Frobenius en $L_\pL/K_\pK$, la cual es 
una extensi\'on no ramificada. Ahora bien $\langle\Fro \pK\rangle
=\Gal(L_\pK/K_\pK)\subseteq \Gal(L/K)$. Por tanto 
$f=o(\Fro {\pK})=[L_\pK:K_\pK]=f^{\prime}$. Se sigue que
$f=f^\prime=n/h$. $\fin$
\end{proof}

\begin{corolario}\label{CClaseT4.3.3.0}\label{C17.6.187N}
Sea $\pK$ finito y no ramificado en $L$. Si $\pi_\pK$ es un
elemento primo de $K_\pK$, entonces $\theta(\pi_\pK)\in C_K/H$
se mapea al automorfismo de Frobenius $\artinp{L/K}{\pK}\in\Gal(L/K)$ bajo
el isomorfismo $C_K/H\isomorfo\limits^{\rho_K} \Gal(L/K)$.
\end{corolario}
\begin{proof}
Es inmediato de los Teoremas \ref{T17.6.186N} y \ref{CCLT17.6.10}.
$\fin$
\end{proof}

\begin{teorema}\label{CClaseTD.1}\label{T17.6.188N}
Sea $L/K$ una extensi\'on abeliana finita de campos globales. Sea
$\pK$ un lugar de $K$. Sea
\[
\theta=\psi_{L/K}\circ \enc{\ }{\pK}\colon \*{K_{\pK}}\xrightarrow{
\enc{\ }{\pK}}J_K\xrightarrow{\psi_{L/K}}\Gal(L/K).
\]

Entonces
\las
\item Para $n\geq 0$, $\theta(U_{\pK}^{(n)})$ es el grupo de ramificaci\'on
superior $G^n(L/K)$ donde $G=\Gal(L/K)$:
\[
\theta(U_{\pK}^{(n)})=G^n(L/K),\quad n\geq 0.
\]

En particular, si $n=0$, 
\[
\theta(U_{\pK}^{(0)})=\theta(U_{\pK})=
G^0(L/K)=G_0(L/K)=I(\pL|\pK)
\] 
es el grupo de inercia.

\item $\theta(\*{K_{\pK}})$ es el grupo de descomposici\'on $D=
D(\pL|\pK)$.
\end{list}
\end{teorema}

\begin{proof} Por el Teorema \ref{T17.6.146N}, se tiene que $\psi\circ \enc{\ }{\pK}=
\theta=(\underline{\ \ },L_{\pK}/K_{\pK})$. Por tanto, por el Teorema
\ref{CClaseT3.2.5.6.3}, se tiene 
\[
\theta(U_{\pK}^{(n)})=(U_{\pK}^{(n)},L_{\pK}
/K_{\pK})=G^n(L/K).
\]
Esto es (1).

Para probar (2), tenemos, del Teorema TCCL, que 
\begin{gather*}
(\underline{\ \ },
L_{\pK}/K_{\pK})\colon \*{K_{\pK}}\lra \Gal(L_{\pK}/K_{\pK})\cong
D(\pL|\pK)
\intertext{es suprayectivo. Por tanto}
\theta(\*{K_{\pK}})=(\*{K_{\pK}},L_{\pK}/K_{\pK})=D(\pL|\pK). \tag*{$\fin$}
\end{gather*}
\end{proof}

\begin{corolario}\label{CClaseCD.2}\label{C17.6.189N}
Sea $K$ un campo global, $L$ una extensi\'on
abeliana finita de $K$ y $H$ el subgrupo abierto de {\'\i}ndice finito de $C_K$
que corresponde a $L$, es decir, $H={\mc N}_L=\N_{L/K} C_L$, $C_K/H\cong
\Gal(L/K)$.

Para un lugar $\pK$ de $K$, consideremos la composici\'on
\[
\mu\colon \*{K_\pK}\xrightarrow{\enc {\ }{\pK}}C_K\xrightarrow{\pi}
C_K/H.
\]
\las
\item $\pK$ se descompone totalmente en $L\iff \mu(\*{K_\pK})=1
\iff \*{K_{\pK}}\subseteq H$.

\item Si $\pK$ es finito, $\pK$ es no ramificado $\iff \mu(U_\pK)=1
\iff U_{\pK}\subseteq H$.
\end{list}

Equivalentemente, si $(\underline{\ \ },L/K)$ denota al s\'imbolo
de Artin, entonces
\las
\item $\pK$ se descompone totalmente $\iff (\*K_{\pK},L/K)=1$.

\item Si $\pK$ es finito, $\pK$ es no ramificado $\iff (U_{\pK},L/K)=1$.
\end{list}
\end{corolario}

\begin{proof} 
\las
\item Se tiene que $\pK$ se descompone totalmente en $L\iff D(
\pL|\pK)=\{1\}\iff \mu(\*{K_{\pK}})=\{1\}$.

\item Si $\pK$ es finito, $\pK$ es no ramificado $\iff I(\pL|\pK)=
\{1\}\iff \mu(U_{\pK})=\{1\}$. $\fin$
\end{list}
\end{proof}

El siguiente teorema es la versi\'on global del Teorema \ref{CClaseT3.2.21'+1}.

\begin{teorema}\label{CClaseT4.6.9-1}\label{T17.6.190N}
Sea $E/F$ una extensi\'on abeliana finita de campos globales y sea 
$E$ el campo de clase de $\Lambda\subseteq C_F$, esto es,
$\N_{E/F} C_E=\Lambda$. Sea $L/F$
una extensi\'on finita y separable. Entonces $LE/L$ es una extensi\'on
finita y el grupo de normas correspondiente es $\N^{-1}_{L/F}(\Lambda)$.
\[
\xymatrix{
L\ar@{-}[rr]^{\N_{L/F}^{-1}(\Lambda)}
\ar@{-}[d]&&LE\ar@{-}[d]\\ F\ar@{-}[rr]_{\Lambda} && E}
\]
\end{teorema}

\begin{proof}
Sea $\psi_{EL/L}\colon C_L\lra\Gal(LE/L)$ el mapeo de Artin. El grupo de
normas correspondiente a $LE/L$ es $\ker \psi_{LE/L}$, esto es,
$C_L/\ker \psi_{LE/L}\cong \Gal(LE/L)$. Por el Teorema
\ref{T17.6.135N} tenemos el diagrama conmutativo
\[
\xymatrix{
C_L\ar@{->}[rr]^{\psi_{LE/L}\phantom{xxx}}
\ar@{->}[d]_{\N_{L/F}}&&\Gal(LE/L)\ar@{->}[d]^{\rest}\\
C_F\ar@{->}[rr]_{\psi_{E/F}\phantom{xxx}}&&\Gal(E/F)
}
\]
Adem\'as el $\rest\colon \Gal(LE/L)\lra \Gal(E/F)$ es inyectivo y $\rest\circ
\psi_{LE/L}=\psi_{E/F}\circ\N_{L/F}$. Por tanto
\begin{gather*}
\idel x\in\ker\psi_{LE/L}\iff \psi_{LE/L}(\idel x)=1\iff \\
\rest\circ \psi_{LE/L}(\idel x)=1=\psi_{E/F}
\circ\N_{L/F}(\idel x)\iff\\
\N_{L/F}(\idel x)\in\ker\psi_{E/F}=\Lambda\iff \idel x\in \N^{-1}_{L/F}
(\Lambda). \tag*{$\fin$}
\end{gather*}
\end{proof}

El siguiente resultado es importante para la teor\'ia de
campos de g\'eneros.

\begin{teorema}\label{CClaseP4.6.9}\label{T17.6.191N}
Sea $L/K$ una extensi\'on finita y separable
de campos globales.
Sea $H$ un subgrupo abierto de \'indice finito en $C_L$ y sea
$L(H)$ su campo de clase. Sea $K_0$ la m\'axima extensi\'on
abeliana de $K$ contenida en $L(H)$. Entonces el grupo de
normas de $K_0$ es $\N_{L/K}(H)$, es decir, $\N_{K_0/K} C_{K_0}=
\N_{L/K} H$.
\end{teorema}

\begin{proof} La norma $\N_{L/K}\colon C_L\lra C_K$ es un mapeo
abierto (Teorema \ref{T17.6.61N}) por lo que $\N_{L/K}(H)$
es abierto en $C_K$. Por otro
lado, del hecho de que $H$ es de \'indice finito en $C_L$,
se sigue que 
$\N_{L/K}(H)$ es de \'indice finito en $\N_{L/K} C_L$.

\[
\xymatrix{
&&L(H)\ar@{<->}[r]\ar@{-}[d]& H\subseteq C_L\\
\N_{L/K}(H)\ar@{<->}[r]&K_0\ar@{-}[ru]\ar@{-}[d]&
LK'\ar@{<->}[r]\ar@{-}[d]&H'\subseteq C_L\\
&K\ar@{-}[r]&L
}
\]

Ahora bien, $L/K$ es una extensi\'on finita por lo que $\N_{L/K}
C_L$ es de \'indice finito en $C_K$. Por tanto $\N_{L/K} (H)$ es
de \'indice finito en $C_K$ y se sigue que
$\N_{L/K}(H)$ es tanto cerrado como abierto en $C_K$.

Sea $K'$ el campo asociado a $\N_{L/K}(H)$ y sea $\Lambda$
el grupo de normas asociado a $K_0$. Por el Teorema
\ref{CClaseT4.6.9-1} se tiene que $LK'/L$ corresponde al
grupo de normas $\N_{L/K}^{-1}(\N_{L/K}(H))\supseteq H$
lo cual implica que $LK'\subseteq L(H)$, de donde obtenemos
que $K'\subseteq K_0$. Por tanto 
$\Lambda\subseteq \N_{L/K}(H)$.
Puesto que $LK_0\subseteq L(H)$, del
Teorema \ref{CClaseT4.6.9-1}, obtenemos que $\N^{-1}_{L/K}
(\Lambda)\supseteq H$. Se sigue que $\N_{L/K}(H)\subseteq
\Lambda$. Por tanto $\Lambda=\N_{L/K}(H)$ probando
el resultado.
$\fin$
\end{proof}

Lo que probamos a continuaci\'on refuerza nuevamente que la teor\'ia que
estamos desarrollando \'unicamente es v\'alida para extensiones abelianas.
Este es la raz\'on del porqu\'e la primera desigualdad \'unicamente es
v\'alida para extensiones abelianas. La demostraci\'on, tanto para el caso
local como para el global, es consecuencia del teorema de existencia.

\begin{teorema}[Limitaci\'on de normas\index{limitacion de normas@limitaci\\'on
de normas}]\label{T17.6.192N}
Sea $L/K$ una extensi\'on finita y separable, ya sea de campos locales o de
campos globales. Sea $L'/K$ la m\'axima subextensi\'on abeliana de $L/K$.
Enconces
\[
\N_{L/K}A_L=\N_{L'/K}A_{L'}\quad\text{donde}\quad A_L=
\begin{cases}
\*L&\text{si $L$ es local,}\\ C_L&\text{si $L$ es global.}
\end{cases}
\]
\end{teorema}

\begin{proof}
Sea $H:=\N_{L/K}A_L$. Entonces $H$ es un subgrupo
abierto y de \'indice finito de $A_K$. Por el teorema de existencia, existe
una extensi\'on abeliana finita $F$ de $K$ tal que $\N_{F/K}A_F=H$.
Puesto que 
\[
H=\N_{F/K}A_F\subseteq \N_{L'/K}A_{L'},
\]
se sigue que $L'\subseteq F$. Puesto que $L'$ es la m\'axima subsextensi\'on
abeliana de $L/K$ y $F/K$ es abeliana, $F\subseteq L'$ y por tanto $F=L'$.
$\fin$
\end{proof}

\section{Grupos de congruencias\index{grupos de
congruencias}}\label{CClaseS4.4}

Con el objetivo de hacer m\'as transparente el Teorema de Existencia, vamos
a estudiar los llamados {\em grupos de congruencias\index{grupos
de congruencias}} para un campo global
$K$, particularmente los campos 
num\'ericos. Para el caso de campos de funciones, debemos posponer un
poco la discusi\'on pues necesitaremos una condici\'on extra que 
autom\'aticamente se cumple para campos num\'ericos.

\begin{definicion}\label{CClaseD4.4.1}\label{D17.6.193N}
Un {\em
m\'odulus\index{m\'odulus}} $\mK$ es un producto formal $\me$ donde
$n_{\pK}\geq 0$ para todo $\pK\in{\ma P}_K$, $n_\pK=0$ para casi
todo $\pK$, $n_\pK=0$ si $\pK$ es complejo y $n_\pK=0$ o $1$ si
$\pK$ es real.
\end{definicion}

Extendemos la definici\'on de unidades a todos los lugares.

\begin{definicion}\label{CClaseD4.4.2}\label{D17.6.194N}
Dado un campo global $K$ y
$\pK$ un lugar arbitrario de $K$ se define el {\em grupo de
las $n_\pK$ unidades de $K$\index{grupo de $n$--unidades}}, donde
$n_\pK\geq 0$ por
\begin{gather}\label{CClaseEq4.3}
U_{\pK}^{(n_\pK)}=
\begin{cases}
1+\pK^{n_\pK} & \text{si $\pK\nmid \infty$ y $n_\pK\geq 1$},\\
U_\pK^{(0)}=U_\pK & \text{si $\pK\nmid\infty$ y $n_\pK=0$},\\
{\ma R}^+\subseteq \*{K_\pK} & \text{si $\pK$ es real y $n_\pK=1$},\\
\*{\ma R}=\*{K_\pK} & \text{si $\pK$ es real y $n_\pK=0$},\\
\*{\ma C}=\*{K_\pK} & \text{si $\pK$ es complejo}.
\end{cases}
\end{gather}
\end{definicion}

De esta forma tenemos que si $\pK$ es real, entonces 
$U_{\pK}^{(1)}={\ma R}^+$ y $[\*{K_{\pK}}:U_{\pK}]=
[\*{\ma R}:{\ma R}^+]=2$.

\begin{definicion}\label{CClaseD4.4.3}\label{D17.6.195N}
Para un lugar $\pK$ y un elemento $\alpha_\pK$
se define 
\[
\alpha_\pK\equiv 1\bmod \pK^{n_\pK}\iff \alpha_\pK
\in \unidadess n.
\]
En particular si $\pK$ es complejo o es real y $n_\pK=0$ la 
congruencia no impone ninguna restricci\'on sobre $\alpha_\pK$.
Si $\pK$ es real y $n_\pK=1$, 
$\alpha_\pK\equiv 1\bmod \pK^{n_\pK}\iff \alpha_\pK
\in \unidadess n={\ma R}^+\iff \alpha_\pK>0$.
\end{definicion}

\begin{definicion}\label{CClaseD4.4.4}\label{D17.6.196N}
Sea $\me$ un m\'odulus en un campo global.
Para un id\`ele $\vec \alpha=(
\alpha_\pK)_{\pK\in{\ma P}_K}\in J_K$ definimos
\[
\vec \alpha\equiv 1\bmod \mK\iff \alpha_\pK\equiv 1\bmod
\pK^{n_\pK}\ \forall\ \pK\in{\ma P}_K\iff \alpha_\pK\in \unidadess n
\ \forall\ \pK\in {\ma P}_K.
\]

Se define el {\em grupo de los id\`eles congruentes a $1$ m\'odulo 
$\mK$\index{id\`eles congruentes a $1$}}\label{CClaseJKm1} por
 \[
 \JKm:=\{\vec \alpha\in J_K\mid \vec\alpha\equiv 1\bmod\mK\}
=\prod_{\pK\in{\ma P}_K}\unidadess n.
\]
\end{definicion}

\begin{definicion}\label{CClaseD4.5.4'}\label{D17.6.197N}
Para un campo global $K$, los elementos
{\em totalmente positivos\index{elemento totalmente positivo}}
$\alpha\in K$ son los elementos tales que ${\eu m}$ involucra
a todos los lugares reales y todos estos lugares tienen exponente
$1$, esto es, $\alpha\equiv 1\bmod U_{\pK}^{(1)}$ para todo 
$\pK$ real.

Equivalentemente $\alpha\in K$ es totalmente positivo
si $\sigma\alpha>0$ para todo $\sigma$ lugar real (es decir
para todo encaje $\sigma \colon K\to {\ma C}$ tal que
$\sigma(K)\subseteq {\ma R}$).
\end{definicion}

De esta forma si $\mK$ es un m\'odulus tal que para todo
lugar real $\pK$, el exponente en $\mK$ de ${\pK}$
es $1$, entonces los elementos
que satisfacen $\alpha\equiv 1\bmod \mK$ son elementos
totalmente positivos.

\begin{observacion}\label{CClaseO3.1.5} El elemento $1+\sqrt{2}$ es positivo pero
no totalmente positivo pues $\varphi\colon {\ma Q}(\sqrt{2})\to{\ma R}$,
$\sqrt{2}\mapsto -\sqrt{2}$, satisface que $\varphi(1+\sqrt{2})=1-\sqrt{2}<0$.

Tambi\'en notemos que el ser totalmente positivo depende del campo
ambiente. Por ejemplo, $-1$ es totalmente positivo en ${\ma Q}(i)$ pero
no en ${\ma Q}$.
\end{observacion}

\begin{definicion}\label{CClaseD4.4.5}\label{D17.6198N}
El grupo $\CKm:=\JKm\*K/\*K
\subseteq C_K$ se llama el {\em subgrupo de congruencias
m\'odulo $\mK$ de $C_K$\index{subgrupo de congruencias
de un grupo de clases de id\`eles}}\label{CClaseCKm1}.
El cociente $C_K/\CKm$
se llama el {\em grupo de clases de rayos m\'odulo $\mK$\index{grupo
de clases de rayos}}.
\end{definicion}

\begin{ejemplo}\label{E17.6.199N} Si $\mK=1$,
esto es, $n_{\pK}=0$ para toda $\pK$. Tenemos 
\begin{gather*}
\JKm=J_K^1
=\prod_{\pK|\infty}\*{K_\pK}
\times \prod_{\pK\nmid \infty}U_\pK=\prod_{\pK\in{\ma P}_K} U_{\pK}^{(0)}=U_K,\\
J_K/J_K^1\*K\cong C_K/C_K^1\cong C_K/U_K\*K \cong I_K\cong D_K/P_K.
\end{gather*}
 (ver despu\'es de la Definici\'on \ref{D17.6.30N}).

Con esta nueva terminolog{\'\i}a, se tiene que
\begin{gather*}
J_K/J_K^1=\bigoplus_{\pK\in{\ma P}_K} \*{K_\pK}/U^{(0)}_\pK\cong
\bigoplus_{\substack{\pK\in{\ma P}_K\\ \pK\text{\ finito}}} {\ma Z}\cong D_K.
\intertext{Por lo tanto}
C_K/C_K^1=\coker (\*K\to J_K/J_K^1)\cong
D_K/(\*K)=D_K/P_K=I_K. 
\end{gather*}
Tambi\'en obtenemos que $C_K^1\subseteq C_{K,0}$ 
y que $C_{K,0}/C_K^1\cong
I_{K,0}$ en el caso de campos de funciones.
\end{ejemplo}

Volvemos a enunciar lo anterior con esta nueva terminolog{\'\i}a
para referencia futura.

\begin{teorema}\label{CClaseT4.4.6}\label{T17.6.200N}
Sea $K$ un campo global. Entonces
\begin{gather*}
C_K/C_K^1\cong I_K\quad \text{$K$ cualquier campo global},\\
C_{K,0}/C_K^1\cong I_{K,0}\quad \text{$K$ un campo de funciones}. \tag*{$\fin$}
\end{gather*}
\end{teorema}

\subsection{Campos num\'ericos}\label{CClaseS4.5}

\begin{observacion}\label{O17.7.1N}
Cuando $K$ es un campo num\'erico, se probar\'a que $[C_K:\CKm]<\infty$,
pero si $K$ es un campo de funciones, entonces $[C_K:\CKm]=\infty$. De
hecho, en este caso, $C_K/C_K^1\cong I_K$ y $[C_K:C_K^1]=|I_K|=\infty$.
Este problema lo solucionaremos m\'as adelante. En 
cualquier caso, esto es, $K$ un campo global arbitrario,
 $\CKm$ es un subgrupo abierto de $C_K$.
\end{observacion}

Debido a lo anterior, en esta subsecci\'on \'unicamente consideraremos
$K$ un campo num\'erico.

\begin{teorema}\label{CClaseT4.5.1}\label{T17.7.2N}
Sea $K$ un campo num\'erico. Los grupos de normas de $C_K$, es decir
los grupos $\N_{L/K} C_L$ donde $L/K$ es una extensi\'on abeliana finita,
son precisamente los subgrupos de $C_K$ que contienen a alg\'un subgrupo
de congruencias $\CKm$.
\end{teorema}

\begin{proof} Se tiene que, por definici\'on, 
$\JKm=\prod_{\pK\in{\ma P}_K}\unidadess n$,
donde $\me$, es un subgrupo abierto de $J_K$. Puesto que $C_K$ tiene la 
topolog{\'\i}a cociente, $\CKm$ es un subgrupo abierto de $C_K$ el cual tiene
{\'\i}ndice
\begin{align*}
[C_K:\CKm]&=[C_K:C_K^1][C_K^1:\CKm]=|I_K|[C_K^1:\CKm]
=h_K[J_K^1\*K:\JKm \*K]\\
&\leq h_K[J_K^1:\JKm]=h_K\prod_{\pKK}
[U_\pK:\unidadess n]<\infty
\end{align*}
pues $n_\pK=0$ para casi toda $\pKK$.

De esta forma obtenemos que $\CKm$ es un subgrupo abierto (y por tanto
cerrado) de {\'\i}ndice finito en $C_K$, por lo que $\CKm$ es un grupo de 
normas (Teorema TCCG). Sea ahora cualquier $H$ subgrupo de $C_K$ con
$H\supseteq \CKm$ para alg\'un m\'odulus ${\eu m}$. Puesto
que $[C_K:\CKm]<\infty$, se tiene $[H:\CKm]<\infty$. Por tanto
$H=\bigcup_{\text{finita}} x\CKm$ de donde se sigue que $H$ es un subgrupo
abierto de $C_K$ y por tanto es un grupo de normas.

Rec{\'\i}procamente, sea $H$ un grupo de normas de $C_K$, es decir,
$H$ es un subgrupo abierto y cerrado de {\'\i}ndice finito en $C_K$. Sean 
$J_K\stackrel{_\theta} {\twoheadrightarrow}C_K$ la proyecci\'on natural y
${\mathcal H}:=\theta^{-1}(H)$. Entonces ${\mathcal H}$ es un subgrupo
abierto de $J_K$, por lo que ${\mathcal H}$ contiene a un subconjunto
de la forma 
\[
W=\prod_{\pK\in S}W_\pK\times\prod_{\pK\notin S}U_\pK,
\]
donde $S$ es un conjunto finito y cada $W_\pK$ es una vecindad abierta
de $1\in \*{K_\pK}$. 

Si $\pK\in S$ es finito, podemos tomar $W_\pK=
\unidadess n$ para alg\'un $n_\pK\geq 0$ debido a que los subgrupos
$\big\{U_\pK^{(n_{\pK})}\big\}_{n_{\pK}\in {\ma N}\cup\{0\}}$
forman una base de vecindades de $1$ en $\*{K_\pK}$.
Si $\pK\in S$ es infinito, $W_\pK$ genera ya sea a ${\ma R}^+$ o a todo
$\*{K_\pK}\in\{\*{\ma R},\*{\ma C}\}$. De esta forma, el subgrupo generado
por $W$ es de la forma $\JKm$ para el m\'odulus $\me$. Se sigue
que ${\mc H}=\theta^{-1}(H)$ contiene a $J_K^{\eu m}$ y por ende
$H$ contiene a $C_K^{\eu m}$.
$\fin$
\end{proof}

\begin{observacion}\label{CClaseO4.5.2}\label{O17.7.3N}
El Teorema \ref{CClaseT4.5.1} no es aplicable
a campos de funciones pues para todo m\'odulus $\mK$, $\CKm
\subseteq C_K^1$ de donde obtenemos
\[
[C_K:\CKm]\geq [C_K:C_K^1]=|I_K|=\infty.
\]

Lo que si se tiene es que $C_{K,0}\cong \Gal(\abe K/K\abe \F)$ y 
$\CKm\subseteq C_{K,0}$ para cualquier m\'odulus $\mK$. M\'as a\'un, si
$H\subseteq C_{K,0}$ es un subgrupo
abierto y de {\'\i}ndice finito en $C_{K,0}$,
entonces $H\supseteq \CKm$ para alg\'un $\mK$ y $H$ es cerrado en
$C_{K,0}$. Por teor{\'\i}a de Galois (Teorema \ref{CClaseT1.4.2}), 
si $(\abe K)^H=L$ y $(\abe K)^{\CKm}=\Km$ entonces $L\subseteq \Km$
exactamente como en el caso num\'erico y podr{\'\i}amos hacer una
discusi\'on totalmente an\'aloga que lo discutido en esta secci\'on
para campos num\'ericos a esta situaci\'on de campos de funciones. La
diferencia es que para campos de funciones $L$ y $\Km$ no corresponden con
grupos de normas. Para un remedio a esta situaci\'on, ver Secci\'on
\ref{CClaseS4.9}.
\end{observacion}

Volvemos a nuestra situaci\'on en que $K$ es un campo num\'erico.

\begin{definicion}\label{CClaseD4.5.3}\label{D17.7.4N}
 El campo de clase $K^\mK/K$ correspondiente
al grupo de congruencias $\CKm$, es decir, $\N_{K^\mK/K}C_{K^\mK}=\CKm$,
se llama el {\em campo de clases de rayos m\'odulo $\mK$\index{campo de 
clases de rayos}}\label{CClaseKm}.
\end{definicion}

Se tiene que, para campos num\'ericos, 
\fbox{$\Gal(K^\mK/K)\cong C_K/\CKm$}.

\begin{observacion}\label{CClaseO4.5.4}\label{O17.7.5N}
 El Teorema \ref{CClaseT4.5.1} prueba que toda extensi\'on
abeliana finita $L/K$ de campos num\'ericos est\'a contenida
en alg\'un campo de clases de rayos $K^\mK$. Esto es, si $L/K$ es una
extensi\'on abeliana finita, existe $K^\mK$ tal que $L\subseteq K^\mK$.
Esto se sigue del Corolario \ref{CClaseC4.2.2} pues si $L$ corresponde a $H$,
y $K^\mK$ corresponde a $\CKm$
entonces $H\supseteq \CKm\iff L\subseteq K^\mK$.

Tambi\'en tenemos que si $\mK|{\eu n}$ entonces $C_L^{\eu n}\subseteq
\CKm$ y por consiguiente $K^\mK\subseteq K^{\eu n}$.
\end{observacion}

\begin{definicion}\label{CClaseD4.5.5}\label{D17.7.6N}
Sean $K$ un campo num\'erico y $L/K$ una extensi\'on abeliana
finita, y sea $\NN L=\N_{L/K}C_L\subseteq C_K$.
Se define el {\em conductor\index{conductor global}} $\f{}=\f {L/K}=\f{}
(L/K)$ de $L/K$ (o de $\NN L$) como el m\'aximo com\'un divisor de todos
los moduli $\mK$ tales que $L\subseteq K^\mK$, 
esto es, $\CKm\subseteq \NN L$.

En otras palabras, $K^{\f{}}/K$ es el m{\'\i}nimo campo de clases de
rayos que contiene a $L/K$ y si $L\subseteq K^{\mK}$, entonces $\f{}|\mK$.
\end{definicion}

Se tiene que dado un campo 
num\'erico $K$ y ${\eu m}$ un m\'odulus, se tiene

\begin{teorema}\label{CClaseT3.1.6} 
{\ }

\las
\item Existe una \'unica extensi\'on $K^{\eu m}$ de $K$ que tiene la
siguiente propiedad: si $\pK$ es un ideal primo no cero de $\o_K$
que no divide a ${\eu m}$ entonces $\pK$ es no ramificado y se tiene
\begin{gather*}
\text{$\pK$ es totalmente descompuesto en $K^{\eu m} \iff$}\\
\text{existe un elemento totalmente positivo $\alpha\in\o_K$ tal que}\\
\text{$\pK=(\alpha), \quad\alpha\equiv 1\bmod {\eu m}$}.
\end{gather*}
(ver la Observaci\'on {\rm{\ref{CClaseO4.5.4}}}).

\item $K^{\eu m}/K$ es una extensi\'on abeliana finita  
de $K$ y toda extensi\'on
abeliana finita de $K$ est\'a contenida en alg\'un $K^{\eu m}$.

\item Si ${\eu n}\subseteq {\eu m}$ 
entonces $K^{\eu m}\subseteq K^{\eu n}$.

\item Si $L/K$ es una extensi\'on abeliana finita entonces existe un
ideal no cero ${\eu f}$ de $\o_K$ m\'aximo tal que 
$L\subseteq K^{\eu m}$.
Adem\'as para todo ideal primo $\pK$ no cero de $\o_K$, $\pK$ es
ramificado en $L\iff \pK|{\eu f}$.
\end{list}

El ideal ${\eu f}$ dado en {\rm (4)} 
es el {\em conductor\index{conductor}} (ver Definici\'on 
{\rm{\ref{CClaseD4.5.5}}}).

Los campos $K^{\eu m}$ son los {\em campos de clases de 
rayos\index{campos de clases de rayos}} (Definici\'on
{\rm{\ref{CClaseD4.5.3}}}). $\fin$
\end{teorema}

\begin{ejemplos}\label{CClaseE3.1.7} 
\las
\item
Sea $K={\ma Q}$, ${\eu m}=(n)$. Para cada primo $p$ de ${\ma Z}$, se
tiene que $p$ es totalmente positivo pero $-p$ no lo es. Sea $\pK$ 
generado por $\pm p$. As{\'\i} que decir ``existe un entero totalmente
positivo $\alpha$ tal que $\pK=(\alpha)$ con $\alpha\equiv 1\bmod n$''
no es equivalente a \'unicamente decir ``$\pK=(p)$ con un n\'umero
primo positivo tal que $p\equiv 1\bmod n$''. Esto es, ${\ma Q}(\zeta_n)$
tiene la propiedad $K^{\eu m}$, m\'as precisamente, ${\ma Q}(\zeta_n)=
{\ma Q}^{(n)}$. Por la unicidad podemos concluir 
que ${\ma Q}(\zeta_n)={\ma Q}^{(n)}$. Se tiene ${\eu f}={\eu m}=(n)$.

\item Para $K={\ma Q}(\zeta_3)$, ${\eu f}=(6)$, $K^{\eu f}=
{\ma Q}(\zeta_3,\sqrt[3]{2})=K(\sqrt[3]{2})$ (ver \cite[table 5.7, p\'agina
10]{KaKuSa2011}).

\item $K={\ma Q}(\sqrt{-5})$, $K^{\o_K}={\ma Q}(\sqrt{-5},\sqrt{-1})=
K(\sqrt{-1})$ (campo de clase de Hilbert, ver Definici\'on \ref{D17.7.16N}).

\item $K={\ma Q}(\sqrt{-6})$, $K^{\o_K}={\ma Q}(\sqrt{-6},\zeta_3)=
K(\zeta_3)=K(\sqrt{-3})$ (campo de clase de Hilbert).
\end{list}

Es f\'acil ver que ${\ma Q}(\sqrt{-5},\sqrt{-1})/{\ma Q}(\sqrt{-5})$ y
que ${\ma Q}(\sqrt{-6},\zeta_3)/{\ma Q}(\sqrt{-6})$ son extensiones
no ramificadas. Ver la Subsecci\'on \ref{S12.4.0}.
\end{ejemplos}

A continuaci\'on presentamos la relaci\'on entre $K^{\mK}$ y $K^{\eu n}$
para dos m\'oduli ${\mK}$, ${\eu n}$ de un campo num\'erico.

Se tiene que $K^{\mK}$ es el campo de clase de $C_K^{\mK}$ y $K^{\eu n}$
es el campo de clase de $C_K^{\eu n}$:
\begin{gather*}
K^{\mK}\longleftrightarrow C_K^{\mK},\quad K^{\eu n}\longleftrightarrow
C_K^{\eu n}.
\intertext{Por el Teorema \ref{CClaseTC.1}, se tiene}
K^{\mK}\subseteq K^{\eu n}\iff C_K^{\eu n}\subseteq C_K^{\mK}\iff
C_K^{{\eu f}_{K^{\eu n}}}\subseteq C_K^{{\eu f}_{K^{\eu m}}} \iff
{\eu f}_{K^{\eu m}}| {\eu f}_{K^{\eu n}}.
\end{gather*}

Escribimos $C_K^{\mK}={\mc N}_{K^{\mK}}$ y $C_K^{\eu n}={\mc N}_{K^{\eu n}}$.
Entonces ${\mc N}_{K^{\mK}\cap K^{\eu n}}={\mc N}_{K^{\mK}}{\mc N}_{K^{\eu n}}
=C_K^{\mK}C_K^{\eu n}$.

Sean ${\eu c}=\mcd({\eu m},{\eu n})$ y ${\eu d}=\mcm ({\eu m},{\eu n})$. Entonces
${\eu c}|{\eu m}$ y ${\eu c}|{\eu n}$ por lo que $K^{\eu c}\subseteq K^{\eu m}\cap
K^{\eu n}$.
Similarmente $K^{\eu m}K^{\eu n}\subseteq K^{\eu d}$.

Ahora, veamos que $C_K^{\eu c}\subseteq C_K^{\eu m}C_K^{\eu n}$. Sean
${\eu m}=\prod_{\pK}\pK^{m_{\pK}}$ y ${\eu n}=\prod_{\pK}\pK^{n_{\pK}}$.
Entonces
${\eu c}=\prod_{\pK}\pK^{\min\{m_{\pK},n_{\pK}\}}=\prod_{\pK}\pK^{c_{\pK}}$
y ${\eu d}=\prod_{\pK}\pK^{\max\{m_{\pK},n_{\pK}\}}=\prod_{\pK}\pK^{d_{\pK}}$.

Sea $\vec\alpha\equiv 1\bmod {\eu c}$. Entonces $\alpha_{\pK}\in U_{\pK}^{
(c_{\pK})}$. Definimos id\`eles $\vec \beta$ y $\vec \gamma$ dados por
\begin{gather*}
\beta_{\pK}= \begin{cases} 
1&\text{si $n_{\pK}<m_{\pK}$}\\ \alpha_{\pK}&\text{si $n_{\pK}\geq m_{\pK}$}
\end{cases},
\qquad
\gamma_{\pK}=\begin{cases}
\alpha_{\pK}&\text{si $n_{\pK}<m_{\pK}$}\\ 1&\text{si $n_{\pK}\geq m_{\pK}$}
\end{cases}.
\end{gather*}
Se sigue que $\vec\beta\equiv 1\bmod {\eu n}$ y que $\vec\gamma\equiv
1\bmod \mK$. Por tanto
$\vec\alpha=\vec \beta\vec\gamma\in C_K^{\eu n}C_K^{\mK}$ lo
cual implica $C_K^{\eu c}\subseteq C_K^{\eu n}C_K^{\mK}$.

Ahora, tenemos $C_K^{\eu m}\subseteq C_K^{\eu c}$ y $C_K^{\eu n}\subseteq
C_K^{\eu c}$ por lo que $C_K^{\mK}C_K^{\eu n}\subseteq C_K^{\eu c}$. Se
siguen las igualdades
\[
C_K^{\eu c}=C_K^{\eu n}C_K^{\mK}=C_K^{\mcd({\eu n}, \mK)}\quad\text{y}
\quad K^{\eu c}=K^{\mcd({\eu n},{\eu m})}=K^{\eu n}\cap K^{\eu m}.
\]

Por otro lado ${\eu n}|{\eu d}$ y ${\eu m}|{\eu d}$. Por tanto $C_K^{\eu d}
\subseteq C_K^{\eu n}\cap C_K^{\eu m}$. Rec\'iprocamente, consideremos
$\vec\alpha\in C_K^{\eu n}\cap C_K^{\eu m}$. Entonces $\alpha_{\pK}\in
U_{\pK}^{(n_{\pK})}\cap U_{\pK}^{(m_{\pK})}=U_{\pK}^{(\max\{n_{\pK},m_{\pK}\})}
=U_{\pK}^{(d_{\pK})}$. Se sigue que $\vec\alpha\in C_K^{\eu d}$ y que
$C_K^{\eu n}\cap C_K^{\eu m}\subseteq C_K^{\eu d}$. Hemos obtenido
las igualdades
\[
C_K^{\eu n}\cap C_K^{\eu m}=C_K^{\eu d}\quad \text{y}\quad
K^{\eu d}=K^{\eu n}K^{\eu m}.
\]

Resumimos nuestra discusi\'on anterior en la siguiente proposici\'on.

\begin{proposicion}\label{CClasePC.1'}\label{P17.7.7N}
Si ${\eu n}$ y ${\eu m}$ son dos moduli, ${\eu c}=\mcd({\eu n},{\eu m})$
y ${\eu d}=\mcm({\eu n},{\eu m})$, entonces
\begin{gather*}
K^{\eu c}=K^{\eu n}\cap K^{\eu m} \quad\text{y}\quad
K^{\eu d}=K^{\eu n}K^{\eu m}. \tag*{$\fin$}
\end{gather*}
\end{proposicion}

\begin{observacion}\label{CClaseO4.5.6}\label{O17.7.8N}
En general no se cumple que $\mK$
sea el conductor $\f{}$ de $K^\mK/K$ aunque por supuesto $\f{}|\mK$.
Esto se debe a que se puede tener $C_K^{\f{}}=\CKm$ con $\f{}|{\eu m}$
y $\f{}\neq {\eu m}$.
M\'as adelante (Consecuencia \ref{CClaseC4.5.9})
veremos un ejemplo de este fen\'omeno.
El lector puede pensar que lo esencial radica en que para $n$ impar
tenemos ${\ma Q}(\zeta_{2n})={\ma Q}(\zeta_n)$ y el conductor
de $\cic {2n}{}$ es $\f{}=n\cdot\infty$ y no $2\cdot n\cdot \infty$.
\end{observacion}

\begin{ejemplo}\label{CClaseE4.5.7}\label{E17.7.9N}
Ver \cite[Theorem 7.7, p\'agina 100]{Neu86} y
\cite[Theorems 7.10 y 7.11, p\'agina 165]{Neu69}.
Consideremos $K={\ma Q}$. Los lugares de ${\ma Q}$ son $\infty$ y $\pK
=(p)$ con $p$ un n\'umero primo.

Se tiene $J_{\ma Q}=\{\vec x\in (\*{\ma R}\times \prod_{p\text{\ primo}}
\*{{\ma Q}_p})\mid x_p\in \*{{\ma Z}_p}=U_{{\ma Q}_p} \text{\ para casi toda
$p$}\}$.

Sea $\mK$ un m\'odulus, $\mK=\mK_0\cdot \infty^{\varepsilon}$, donde
$\mK_0=\prod_{i=1}^r \pK_i^{\alpha_i}$, $\alpha_i\geq 0$,
$1\leq i\leq r$ y $\varepsilon\in\{0,1\}$
con cada $\pK_i$ un primo finito. Por abuso del lenguaje, si $\pK_i=(p_i)$
con $p_i$ un primo racional positivo, pondremos $\mK_0=m=\prod_{i=1}^r
p_i^{\alpha_i}$.

Se tiene 
\begin{align*}
J_{\ma Q}^{\mK_0}=J_{\ma Q}^m&=\*{\ma R}
\times \prod_{p\text{\ finito}} U_p^{
(n_p)}=\*{\ma R}\times \prod_{i=1}^r U_{p_i}^{(\alpha_i)}\times
\prod_{p\notin\{p_1,\ldots, p_r\}}\*{{\ma Z}_p}\\
&=\*{\ma R}\times\prod_{i=1}^r(1+p_i^
{\alpha_i} {\ma Z}_{p_i}) \times \prod_{p\notin\{p_1,\ldots, p_r\}}
\*{{\ma Z}_p}
\intertext{y}
J_{\ma Q}^{\mK_0\infty}&={\ma R}^+\times\prod_{i=1}^r(1+p_i^
{\alpha_i} {\ma Z}_{p_i}) \times \prod_{p\notin\{p_1,\ldots,
p_r\}}\*{{\ma Z}_p}.
\end{align*}

Sea $(m_0)={\eu m}_0$, es decir, $m_0$ es el generador de
${\eu m}_0$ en ${\ma Q}^+$. Por supuesto, $m_0=
\prod_{i=1}^r p_i^{\alpha_i}\in {\ma N}$.

Ahora veamos que $U_p^{(n_p)}$ consiste de normas provenientes de
${\ma Q}_p(\zeta_m)/{\ma Q}_p$ (=$({\ma Q}
(\zeta_m))_{\pL}/{\ma Q}_\pK$, donde $\pL$ es un primo de 
$\cic {m_0}{}$ sobre $\pK=\langle p\rangle$).
Lo anterior se obtiene de la siguiente forma. Escribiendo $m=np^r$ con
$\mcd(n,p)=1$ se tiene ${\ma Q}_p(\zeta_m)={\ma Q}_p(\zeta_n)
{\ma Q}_p(\zeta_{p^r})$ y ${\ma Q}_p(\zeta_n)/{\ma Q}_p$ es no ramificada
pues $p\nmid n$. Por el Teorema \ref{CCLT17.6.3}, se tiene
$\N_{{\ma Q}_p(\zeta_n)/
{\ma Q}_p}U_\pL= U_p$ donde $\pL$ es el primo de ${\ma Q}_p(\zeta_n)$,
lo cual se sigue de que $\co 0G{U_{\pL}}=\{1\}$.

Por otro lado ${\ma Q}_p(\zeta_{p^r})/{\ma Q}_p$
es una extensi\'on totalmente remificada y
se tiene que $\N_{{\ma Q}_p(\zeta_{p^r})/{\ma Q}_p}\big(
\*{{\ma Q}_p(\zeta_{p^r})}\big)=(p)\times U_p^{(r)}$ 
(Teorema \ref{CClaseT3.2.5.29}). Del Teorema \ref{CClaseT3.2.29}, si ponemos
$L_1={\ma Q}_p(\zeta_{n})$ y $L_2={\ma Q}_p(\zeta_{p^r})$, entonces
$L={\ma Q}_p(\zeta_{m})=L_1 L_2$ y se sigue,
con $r=n_p$, que $U_p^{(n_p)}\subseteq 
\NN L=\NN {{L_1 L_2}}=\NN {{L_1}}\cap \NN {{L_2}}=\N_{{\ma Q}_p(\zeta_{m})/{\ma Q}_p}
\big(\*{{\ma Q}_p(\zeta_{m})\big)}$.

Ahora se tiene
que un id\`ele $\vec \alpha\in J_K$ es la norma de un id\`ele $\vec \beta\in
J_L$ si y solamente si cada componente $\alpha_\pK\in\*{K_\pK}$ es la norma de
un elemento $\beta_\pL\in \*{L_\pL}$ para $\pL|\pK$ (ver Teorema \ref{T17.6.68N}).
Por tanto se sigue $C_{\ma Q}^{\mK_0\infty}\subseteq
\N_{{\ma Q}(\zeta_{m})/{\ma Q}}(C_{{\ma Q}(\zeta_{m})})$.

Consideremos $\mK=\mK_0 \infty$. Se tiene que 
\begin{align*}
[C_{\ma Q}:C_{\ma Q}^\mK]&=[C_{\ma Q}:C_{\ma Q}^1][C_{\ma Q}^1:
C_{\ma Q}^\mK]=h_{\ma Q}[J_{\ma Q}^1\*{\ma Q}:J_{\ma Q}^\mK\*{\ma Q}]\\
&= 1\cdot \frac{[J_{\ma Q}^1:J_{\ma Q}^\mK]}{[(J_{\ma Q}^1\cap \*{\ma Q}):
(J_{\ma Q}^\mK\cap \*{\ma Q})]},
\end{align*}
donde hemos utilizado el Corolario \ref{C17.6.157N}.

Ahora bien se tiene que $J_{\ma Q}^1=\*{\ma R}\times\prod_p U_p$
y que $J_{\ma Q}^\mK={\ma R}^+\times\prod_p U_p^{(n_p)}$, por lo que
$J_{\ma Q}^1\cap \*{\ma Q}=\{1,-1\}$ y $J_{\ma Q}^\mK\cap\*{\ma Q}
=\{1\}$. Por tanto obtenemos
\begin{align*}
[C_{\ma Q}:C_{\ma Q}^\mK]&=\frac{1}{2} \prod_p[U_p:U_p^{(n_p)}][\*{\ma R}:
{\ma R}^+]=\prod_p[U_p:U_p^{(n_p)}]=\prod_{p|m}p^{n_p-1}(p-1)\\
&=\varphi(m)=[{\ma Q}(\zeta_m):{\ma Q}]=[C_{\ma Q}:
\N_{{\ma Q}(\zeta_{m})/{\ma Q}}(C_{{\ma Q}(\zeta_{m})})]
\end{align*}
y ya que tenemos $C_{\ma Q}^\mK\subseteq
\N_{{\ma Q}(\zeta_{m})/{\ma Q}}(C_{{\ma Q}(\zeta_{m})})$
entonces $C_{\ma Q}^\mK=
\N_{{\ma Q}(\zeta_{m})/{\ma Q}}(C_{{\ma Q}(\zeta_{m})})$.

En particular obtenemos que ${\ma Q}^\mK={\ma Q}(\zeta_m)$.

Si consideramos ahora $\mK=\mK_0$, entonces se tiene
que $C_{\ma Q}^{\mK_0}\subseteq C_{\ma Q}^{\mK}$, por lo que
${\ma Q}^{\mK_0}\subseteq {\ma Q}^{\mK}$. Por otro lado
$[J_{\ma Q}^{\mK}:J_{\ma Q}^{\mK_0}]=2$ de donde
$[{\ma Q}^\mK:{\ma Q}^{\mK_0}]=2$. Finalmente puesto que
las normas de ${\ma Q}^{\mK_0}$ no est\'an contenidas en
${\ma R}^+$, el campo ${\ma Q}^{\mK_0}$ necesariamente es
real. Se sigue que ${\ma Q}^{\mK_0}={\ma Q}(\zeta_m)^+=
{\ma Q}(\zeta_m+\zeta_m^{-1})$.
\end{ejemplo}

\begin{corolario}[Teorema de Kronecker--Weber\index{teorema
de Kronecker--Weber}\index{Kronecker--Weber!teorema de $\sim$}]
\label{CClaseT.4.5.8}\label{T17.7.10N}
Toda extensi\'on abeliana finita $L/{\ma Q}$ est\'a contenida
en alg\'un campo ciclot\'omico.

M\'as a\'un, si $L/{\ma Q}$ es una extensi\'on abeliana finita
con $L\subseteq {\ma R}$, entonces existe $m\in{\ma N}$ tal que
$L\subseteq {\ma Q}(\zeta_m)^+$.

Se tiene que si $\mK=\mK_0\infty=
m\infty$ es un m\'odulus
arbitrario de ${\ma Q}$, entonces ${\ma Q}^{\mK}={\ma Q}(\zeta_m)$
y ${\ma Q}^{\mK_0}={\ma Q}(\zeta_m)^+$.
\end{corolario}

\begin{proof} Los campos de clases de rayos de ${\ma Q}$ son
los campos ${\ma Q}(\zeta_m)$ y ${\ma Q}(\zeta_m)^+\subseteq
{\ma Q}(\zeta_m)$ y toda extensi\'on abeliana finita est\'a contenida
en alg\'un campo de clases de rayos. Adem\'as $\mK_0|\mK_0
\infty=\mK$ por lo que $K^{\mK_0}\subseteq K^{\mK}={\ma Q}^{\mK}
={\ma Q}(\zeta_m)$.
$\fin$
\end{proof}

\begin{consecuencia}\label{CClaseC4.5.9}\label{C17.7.11N}
El conductor de $K^\mK$  puede ser $\f{}\neq\mK$, $\f{}|\mK$.
Por ejemplo, si $m\in{\ma N}$ es impar y $\mK=2m\infty$ se tiene
${\ma Q}^\mK={\ma Q}(\zeta_{2m})={\ma Q}(\zeta_m)=
{\ma Q}^{\f{}}$ donde $\f{}=m\infty\neq \mK$.
\end{consecuencia}

\begin{observacion}\label{O17.7.12N}
Los campos de clases de rayos deben ser considerados como los
an\'alogos a los campos ciclot\'omicos pues ${\ma Q}^{\mK}=
\cic m{}$ donde $\mK=m\cdot \infty$ y si $\mK_0=m$, entonces
${\ma Q}^{\mK_0}=\cic m{}^+={\ma Q}(\zeta_m+\zeta_m^{-1})$.
\end{observacion}

\begin{observacion}\label{O17.7.13N}
Si consideramos a ${\ma C}/{\ma R}$ como una extensi\'on de campos
locales, entonces $\N_{{\ma C}/{\ma R}}\*{\ma C}=\big(\*{\ma R})^2=
{\ma R}^+=U^{(1)}_{\ma R}$. Por tanto el conductor en este caso es
$\f{}=\pK$ donde $\pK$ es el primo real.
\end{observacion}

\begin{teorema}\label{CClaseT4.5.10}\label{T17.7.14N}
Si $\f{}$ es el conductor de una 
extensi\'on abeliana finita $L/K$ de campos num\'ericos y si
si $\f\pK$ es conductor local de $L_\pK/K_\pK$ para $\pKK$,
entonces
\[
\f{}=\f{L/K}=\prod_{\pKK}\f{\pK}.
\]
\end{teorema}

\begin{proof} 
 Sean ${\mc N}=\N_{L/K} C_L$ y
${\eu n}:=\prod_{\pK\in{\ma P}_K} {\eu f}_{\pK}=\prod_{\pK\in{\ma P}_K}
\pK^{n_{\pK}}$. Por definici\'on se tiene que para un m\'odulus 
$\mK=\prod_{\pK\in{\ma P}_K}\pK^{m_{\pK}}$ tal que $C_K^{\mK}\subseteq
{\mc N}=\N_{L/K}C_L$, entonces $\CKm\subseteq C_K^{\f{}}\subseteq
{\mc N}\iff {\eu f}|{\eu m}$ y $C_K^{\f{}}\subseteq {\mc N}$.
Que $C_K^{\eu n}\subseteq {\mc N}$ se prueba con el argumento
al final de la demostraci\'on.
Por tanto debemos probar que
\[
C_K^{\mK}\subseteq {\mc N}\iff {\eu n}|{\eu m}\iff n_{\pK}\leq
m_{\pK} \text{\ para toda\ } \pK \iff f_{\pK}|\pK^{m_{\pK}}.
\]

Se tiene: 
\begin{align*}
C_K^{\mK} \subseteq \NN{} & \iff \big(\vec\alpha\equiv 1\bmod
\mK \Longrightarrow \idel \alpha\in \NN{} \big) \text{\ para $
\vec\alpha\in J_K$}\\
&\iff \big(\vec\alpha\equiv 1\bmod \mK\Longrightarrow \enc 
{\alpha_\pK}{\pK}=(\ldots,1,1,\alpha_\pK,1,1,\ldots)\\
&\hspace{2cm}\in\NN{}\cap
\*{K_\pK}=\N_\pK\*{L_\pK} \text{\ para toda $\pK$}\big)
\quad\text{(Teorema \ref{T17.6.68N})}\\
&\iff \big(\alpha_\pK\in \unidadess m\Longrightarrow\alpha_\pK
\in\N_\pK\*{L_\pK}\text{\ para toda $\pK$}\big)\\
&\iff \unidadess m\subseteq \N_\pK \*{L_\pK} \text{\ para toda $\pK$}
\iff \f{\pK}|\pK^{m_\pK}\text{\ para toda $\pK$}.
\end{align*}
Esto mismo prueba que $C_K^{\eu n}\subseteq {\mc N}$ y por
tanto ${\eu n}=\f{}$.
$\fin$
\end{proof}

\begin{corolario}\label{CClaseC4.5.11}\label{C17.7.15N}
Sea $L/K$ una extensi\'on abeliana finita de campos num\'e\-ricos.
Un primo $\pK$ de $K$ es ramificado en $L\iff \pK|\f{}$ donde
$\f{}$ es el conductor de $L/K$.
\end{corolario}
\begin{proof} $\pK$ es ramificado en $L\iff\pK$ es ramificado en $L_\pL/
K_\pK\iff \pK|\f\pK$ (Teorema \ref{CClaseT3.2.24}). $\fin$
\end{proof}

\begin{definicion}\label{CClaseD4.5.12}\label{D17.7.16N}
Sea $K$ un campo num\'erico.
El {\em campo de clase de Hilbert\index{campo de clase de
Hilbert}\index{Hilbert!campo de clase de $\sim$}} $K_H$ es la
m\'axima extensi\'on abeliana no ramificada de $K$.
\end{definicion}

\begin{corolario}\label{CClaseC4.5.13}\label{C17.7.17N}
El campo de clase de Hilbert
$K_H$ es el campo de clases de rayos m\'odulo $1$, es decir,
$K_H=K^1$ ($1$ el m\'odulus trivial) y 
\begin{gather*}
\Gal(K_H/K)=\Gal(K^1/K)\cong C_K/C_K^1\cong I_K. \tag*{$\fin$}
\end{gather*}
\end{corolario}

En otras palabras, el campo de clase de Hilbert $K_H$ de $K$ corresponde
a $J_K^{(1)}=\prod_{\pK} U_{\pK}$ o $C_K^{(1)}=J_K^{(1)}\*K/\*K=
(\prod_{\pK}U_{\pK})\*K/\*K$.

\begin{definicion}\label{CClaseD4.5.14}\label{D17.7.18N} 
El {\em campo de clase de Hilbert
extendido\index{campo de clase de Hilbert extendido}\index{Hilbert!campo
de clase extendido de $\sim$}} $K_{H^+}$ de un campo num\'erico $K$
es la m\'axima extensi\'on abeliana de $K$ no ramificada en
los primos finitos.
\end{definicion}

El campo $K_{H^+}$ es el campo de clases de rayos m\'odulo $1_+$
donde $1_+$ es el m\'odulus $1_+:=\prod_{\pK\text{\ real}}
\pK$ (esto permite la ramificaci\'on de los primos reales).
Los subgrupo de id\`eles correspondientes a los grupos de congruencias
m\'odulo $1_+$ y m\'odulo $1$, y los grupos de congruencias mismos son:
\begin{gather*}
J_K^{1_+}=\prod_{\pK\text{\ real}} {\ma R}^+\times
\prod_{\pK\text{\ complejo}}\*{\ma C}
\times \prod_{\pK\nmid \infty}U_\pK,\quad
C_K^{1_+}=J_K^{1_+}\*K/\*K,\\
J_K^{1}=\prod_{\pK\text{\ real}} {\ma R}^*\times
\prod_{\pK\text{\ complejo}}\*{\ma C}
\times \prod_{\pK\nmid \infty}U_\pK, \quad
C_K^{1}=J_K^{1}\*K/\*K,
\intertext{respectivamente. Se sigue que}
\frac{J_K}{J_K^{1_+}}\cong\bigoplus_{\pK{\text\ real}}\frac{\*{\ma R}}{{\ma R}^+}
\bigoplus \bigoplus_{\pK\nmid\infty}\frac{\*{K_\pK}}{U_\pK}
\cong \Big(\bigoplus_{\pK\text{\ real}}
C_2\Big)\bigoplus \Big(\bigoplus_{\pK\nmid\infty}{\ma Z}
\Big)\cong C_2^r\oplus D_K\quad\text{y}\\
C_K/C_K^{1_+}=\coker\big(\*K\longrightarrow J_K/J_1^{1_+}\big)
=J_K/J_K^{1_+}K^{\ast}
\end{gather*}
donde $r$ es el n\'umero de lugares reales en $K$ y $C_n$ denota
el grupo c{\'\i}clico de $n$ elementos.

\begin{ejemplo}\label{CClaseE4.5.15}\label{E17.7.19N}
Por el teorema del discriminante
de Minkowski, esto es,
en toda extensi\'on propia de ${\ma Q}$ hay primos finitos
ramificados, entonces
\[
{\ma Q}_H={\ma Q}_{H^+}={\ma Q}
\]
y $r=1$ en este caso. Notemos que $J^1_{\ma Q}/J^{1_+}_{\ma Q}
\cong C_2$ pero $C^1_{\ma Q}/C^{1_+}_{\ma Q}=1$.
\end{ejemplo}

Para $K$ un campo num\'erico, se tiene
\begin{gather*}
\frac{J_K}{J_K^1}\cong\bigoplus_{\pK\nmid \infty}\frac{\*K_{\pK}}{U_{\pK}}\cong
D_K\quad\text{y}\quad \frac{J_K^1}{J_K^{1_+}}\cong C_2^r.
\intertext{Adem\'as}
\Gal(K_{H^+}/K_H)\cong\frac{\Gal(\abe K/K_H)}{\Gal(\abe K/K_{H^+})}\cong
\frac{C_K^1}{C_K^{1_+}}=\frac{J_K^1\*K/\*K}{J_K^{1_+}\*K}\cong
\frac{J_K^1\*K}{J_K^{1_+}\*K}
\intertext{y}
\frac{J_K^1}{J_K^{1_+}}\cong C_2^r
\longtwoheadrightarrow \frac{J_K^1\*K}{J_K^{1_+}\*K}
\cong \Gal(K_{H^+}/K_H),
\end{gather*}
de donde $\Gal(K_{H^+}/K_H)\cong C_2^s$ con $s\leq r=$ n\'umero de lugares
reales de $K$.

\subsection{Campos de clases de rayos en campos de 
\index{campos de rayos en campos de 
funciones}funciones}\label{CClaseS4.9}

Ahora consideramos $K$ un campo global de funciones.
La m\'axima extensi\'on abeliana no ramificada $\nr K$ de $K$
es el campo de clase correspondiente al
grupo $U_K=\prod_{\pKK} U_\pK$, el cual es abierto
y cerrado en $J_K$ puesto que $U_{\pK}$ es abierto y cerrado en
$\*K_{\pK}$ para toda $\pK$, por lo que, por un lado es abierto por
definici\'on de la topolog\'ia de $J_K$ y por otro lado es compacto
pues $U_{\pK}$ es compacto (ver Proposiciones \ref{CCUnidades} y
\ref{CCCompactos}) y $U_K$ tiene la topolog\'ia producto.
M\'as precisamente $\nr K$ es el campo correspondiente
de $C_K^1=K^{\ast} U_K/K^{\ast}=C_K^1
\subseteq C_{K,0} \subseteq C_K$ y
$[C_K:C_K^1]=\infty =|\Gal(K^{\rm{nr}}/K)|$
(ver Observaci\'on \ref{CClaseO4.5.2}).
$K^{\rm{nr}}$ no es un campo de clase en el sentido de normas.
Adem\'as $K\abe\F\subseteq \nr K$.

Se tiene $U_K\cap \*K=\{x\in\*K\mid v_{\pK}(x)=0\text{\ para toda $\pK\in
{\ma P}_K$}\}=\*\F$. Adem\'as $U_K\subseteq J_{K,0}$ y $J_{K,0}/U_K\cong
D_{K,0}$, el grupo de divisores de grado $0$. Por el isomorfismo $\rho_K:
\hat C_K\lra \abe G_K$, se sigue que $C_{K,0}\cong \Gal(\abe K/K\abe\F)$.
Se tiene
\begin{gather*}
\frac{C_{K,0}}{C_K^1}\cong \frac{C_{K,0}}{\big(U_K\*K/\*K\big)}\cong
\Gal(\nr K/K\abe\F)\quad\text{y}\\
\Lambda:J_{K,0}\lra D_{K,0},
\quad \Lambda(\vec\alpha)={\eu a}_{\vec\alpha}=\prod_{\pK\in{\ma P}_K}
\pK^{v_{\pK}(\alpha_{\pK})},\quad \ker\Lambda=U_K\*K\quad\text{y}\\
J_{K,0}\stackrel{\Lambda}{\longtwoheadrightarrow}D_{K,0}\stackrel{\pi}{
\longtwoheadrightarrow}D_{K,0}/P_K\cong I_{K,0},\quad \ker 
(\pi\circ\Lambda)=U_K\*K.
\intertext{Esto es,}
I_{K,0}\cong \frac{J_{K,0}}{U_K\*K}\cong \frac{J_{K,0}/\*K}{U_K\*K/\*K}
\cong \frac{C_{K,0}}{C_K^1}\cong \Gal(\nr K/K\abe \F)\\
(J_K^1=\prod_{\pK}U_{\pK}^{(0)}
=\prod_{\pK} U_{\pK}=U_K,\quad C_K^1=\frac{J_K^1\*K}{
\*K}=\frac{U_K\*K}{\*K}).
\end{gather*}

\begin{teorema}\label{17.7.28N}
Bajo el mapeo de reciprocidad, se tiene el isomorfismo
\begin{gather*}
I_{K,0}\cong \Gal(\nr K/K\abe\F)\cong\frac{C_{K,0}}{C_K^1}.
\tag*{$\fin$}
\end{gather*}
\end{teorema}
\[
\xymatrix{
&&\abe K\ar@{-}[d]\ar@/^1pc/@{-}[d]^{C_K^1}\ar@/^4pc/@{-}[dd]^{C_{K,0}}
\ar@{-}[ddll]\\&&\nr K\ar@{-}[d]\ar@/^1pc/@{-}[d]^{I_{K,0}}\\
K\ar@{-}[rr]&&K\abe \F
}
\]

Ahora consideremos los grupos de clases de rayos. Sea $\me=
\prod_{i=1}^r\pK_i^{\alpha_i}$ un m\'odulus. Sea
\[
\JKm=\prod_{\pKK} \unidadess n=\prod_{\pK\nmid \mK}U_\pK\times
\prod_{i=1}^r U_{\pK_i}^{(\alpha_i)}.\label{CClaseJKm}
\]

Sea $\CKm=\JKm\*K/\*K\label{CClaseCKm}$. Se tiene $\CKm\subseteq
C_{K,0}$ pues si $\vec x\in\JKm$, $v_\pK(x_\pK)=0$ para toda
$\pKK$. Adem\'as
\begin{align*}
[C_{K,0}:\CKm]&=[C_{K,0}:C_K^1][C_K^1:\CKm],\\
[C_{K,0}:C_K^1]&=|\Gal(K^{\rm{nr}}/K\abe \F)|=|I_{K,0}|=h_K<\infty \quad\text{y}\\
[C_K^1:\CKm]&=[\*K J_K^1/\*K:\*K\JKm/\*K]\leq [J_K^1:\JKm]\\
&:= \Psi(\mK):=\prod_{i=1}^r[U_{\pK_i} :
U_{\pK_i}^{(\alpha_i)}]\\
&=\prod_{i=1}^r q^{(\alpha_i-1)\deg \pK_i}\big(q^{\deg \pK_i}
-1\big)<\infty. \label{CClasePsim}
\end{align*}

Por lo tanto,

\begin{proposicion}\label{CClaseP4.9.0}\label{P17.7.29N}
Tenemos $[C_{K,0}:\CKm]\leq h_K\Psi(\mK)<\infty$, donde 
$\Psi(\mK)=\prod_{\pK\in
\sop(\mK)}(q^{\deg \pK}-1)q^{(n_\pK-1)\deg \pK}$,
$\sop(\mK)$ denota al soporte de $\mK$, es decir
$\sop(\mK)=\{\pK\in{\ma P}_K\mid \pK\big|\mK\}$.
\end{proposicion}

\begin{proof} Para $\mK$ un m\'odulus arbitrario tenemos
\begin{align*}
\Psi(\mK)&=\prod_{\pKK}[U_\pK:\unidadess n]=
[U_\pK:\unidadess n]=\prod_{i=0}^{n_\pK-1}[U_\pK^{(i)}:
U_\pK^{(i+1)}]\\
&=[U_\pK:U_\pK^{(1)}]\prod_{i=1}^{n_\pK-1}[
U_\pK^{(i)}:U_\pK^{(i+1)}]=|{\ma F}_{q^{\deg \pK}}^{\ast}|
\prod_{i=1}^{n_\pK-1}|{\ma F}_{q^{\deg \pK}}|\\
&=(q^{\deg \pK}-1)q^{(n_\pK-1)\deg \pK}. \tag*{$\fin$}
\end{align*}
\end{proof}

\subsubsection{Otro tipo de clases de rayos en campos de funciones}

Sea $B$ un subgrupo abierto de $C_K$ de {\'\i}ndice finito. Sea
$\tilde{\vec b}\in B$ tal que $\deg(\tilde{\vec b}):=
\min\{\deg(\tilde{\vec\alpha})\mid \tilde{\vec\alpha}\in B,
\deg(\tilde{\vec\alpha})>0\}$. 

Sea $B_0:=B\cap C_{K.0}$. Entonces $B=\cup_{n=1}^\infty
\tilde{\vec b}^n B_0$ pues 
\[
B/B_0=B/(B\cap C_{K,0})\cong
BC_{K,0}/C_{K,0}\subseteq C_K/C_{K,0}\cong {\ma Z},
\]
y tenemos que $B_0$ es abierto en $C_{K,0}$ pues
$B$ es abierto y puesto que $C_{K,0}$ tambi\'en es abierto,
se tiene que $B_0$ es abierto.

Ahora bien, se tiene que los conjuntos $\CKm$ forman un sistema
fundamental de vecindades de $1\in C_K$
debido a que los $\JKm$ forman un sistema fundamental de 
vecindades de $1\in J_K$. Como $B_0$ es abierto, entonces
$B_0$ debe contener a alg\'un $\CKm\subseteq B_0$ y
\[
h_K\Psi(\mK)=[C_{K,0}:\CKm]=[C_{K,0}:B_0][B_0:\CKm]<\infty,
\]
se sigue que $[C_{K,0}:B_0]<\infty$.
Sea $L_B$ el campo de clase de $B$, $C_K/B=C_K/\N_{
L_B/K}C_{L_B}\cong \Gal(L_B/K)$ con $B=\N_{L_B/K}C_{L_B}$.
\[
\xymatrix{
&&\abe K\ar@{-}[d]\ar@{-}[ddl]\\
&&{\phantom{L_{B_0}=}}L_{B_0}=\tilde{L}_B\ar@{-}[d]\ar@{-}[dl]\\
&L_B\ar@{-}[dl]_{C_K/\N_{L_B/K}C_{L_B}}&K^{\rm{nr}}\ar@{-}[d]\\
K\ar@{-}[rr]&&K\abe \F}
\]

De esta forma, $B_0$ corresponde al campo $L_{B_0}
=\tilde{L}_B:=
L_B \abe \F$ ({\underline{no}} podemos decir que $C_K/
\N_{\tilde{L}_B/K}C_{\tilde{L}_B}\cong \Gal(\tilde{L}_B/K)$
pues $\tilde{L}_B/K$ es infinita y no tenemos definida la norma).

Las clases de rayos que hemos definido, tienen el inconveniente de no
ser de \'indice finito en $C_K$.
Definimos otros grupos de clases de rayos de tal manera que 
ser\'an abiertos y de {\'\i}ndice finito en $C_K$. Se puede hacer el
estudio para $K$ un campo global arbitrario. Solo haremos el 
desarrollo para campos de funciones en donde $S\neq \emptyset$.
Para el caso num\'erico, se puede tomar $S=\emptyset$.

Sea $K$ un campo de funciones y
sea $S$ un conjunto de lugares de $K$, el cual, eventualmente,
pediremos que sea finito y no vac{\'\i}o.

\begin{definicion}\label{CClaseD4.9.1}\label{D17.7.30N}
Un {\em $S$--m\'odulus\index{S--modulus@$S$--m\'odulus}} es un m\'odulus
$\me$ tal que si $n_\pK>0$ entonces $\pK\notin S$, es decir,
el soporte de $\mK$ es disjunto de $S$, esto es, si $\mK=\prod_{i=1}^r
\pK_i^{\alpha_i}$, entonces $S\cap \{\pK_1,\ldots,\pK_r\}=\emptyset$.
\end{definicion}

\begin{definicion}\label{CClaseD4.9.2}\label{D17.7.31N}
Sea $\me$ un $S$--m\'odulus. Se definen los
{\em subgrupos de $S$--congruencias m\'odulo $\mK$\index{subgrupos
de $S$--congruencias m\'odulo un 
m\'odulus}\index{S--congruencias@$S$--congruencias m\'odulo un m\'odulus
subgrupos de $\sim$}\label{CClaseCKSm}} por
\begin{align*}
\JKSm:&=\Big(\prod_{\pK \in S}\*K_{\pK}\times \prod_{\pKK\setminus S}
\unidadess n\Big)\bigcap J_K\\
&=\Big(\prod_{\pK\in S}\*K_{\pK}\times
\prod_{i=1}^r U_{\pK_i}^{(\alpha_i)}\times \prod_{\substack{
\pK\notin S\cup\\\{\pK_1,\ldots,\pK_r\}}}U_{\pK}\Big)
\bigcap J_K \label{CClaseJKSM}\quad\text{y}\\
\CKSm:&=\frac{\JKSm\*K}{\*K}
\end{align*}
de $J_K$ y de $C_K$ respectivamente, y donde $\me$.
\end{definicion}

\begin{observacion}\label{CClaseO4.9.3}\label{D17.7.32N}
La intersecci\'on con $J_K$ \'unicamente tiene significado
cuando $|S|=\infty$. Cuando $S$ es un conjunto finito, tenemos que
$\big(\prod_{\pK \in S}\*{K_\pK}\times \prod_{\pKK\setminus S}
\unidadess n\big)\subseteq J_K$.
\end{observacion}

\begin{proposicion}\label{CClaseP4.9.4}\label{P17.7.33N}
Sea $S\subseteq {\ma P}_K$ y
$\me$ un $S$--m\'odulus.
\lasa
\item Si $T\subseteq S$ y $\mK|{\eu n}$ donde ${\eu n}$ es un
$T$--m\'odulus, entonces
\[
J_{K,T}^{\eu n}\subseteq \JKSm\quad\text{y}\quad C_{K,T}^{\eu n}
\subseteq \CKSm.
\]

\item Si $T\subseteq {\ma P}_K$ y ${\eu n}$ es un $T$--m\'odulus,
\[
\JKSm J_{K,T}^{\eu n}=J_{K,S\cup T}^{\mcd(\mK,{\eu n})}.
\]

\item Los subgrupos de congruencias $J_{K,\emptyset}^{\eu m}=\JKm$ y
$C_{K,\emptyset}^{\mK}=C_K^{\mK}$ con $\mK$ recorriendo todos
los m\'odulus, forman una base de vecindades abiertas de $1$ en
$J_K$ y $C_K$ respectivamente.

\item $\JKSm$ y $\CKSm$ son abiertos.

\item $J_{K,S}^{1}/\JKSm\cong \prod_\pK \big(U_\pK/\unidadess m\big)$
y $[J_{K,S}^{1}:\JKSm]=\Psi(\mK)$.

\item $[C_K:C_{K,\emptyset}^{\eu m}]=[C_K:\CKm]=\infty$ para
toda $\mK$.

\item $K^{\ast} J_{K,{\ma P}_K\setminus \sop(\mK)}^{\eu m}=J_K$
para toda $\mK$.

\item $\JKSm$ es topol\'ogicamente generado por $\enc
{\*{K_\pK}}{\pK}$, $\pK\in S$ y $\enc {\unidadess n}{\pK}$,
$\pKK \setminus S$, donde $\enc {\ }{\pK}\colon \*{K_\pK}
\longrightarrow J_K$ est\'a definido por $\enc {x}{\pK}=
(\ldots, 1,1,\underbracket[0pt]{x}_{\substack{\uparrow\\ \pK}},1,1,\ldots)$.
\end{list}
\end{proposicion}
\begin{proof} \cite{Aue99, Aue2000}.

\noindent
(f) Sea $\pKK\setminus \sop(\mK)$, $\*{K_\pK}\xrightarrow[
\enc{\ }{\pK}]{} J_K/K^{\ast} J_{K,\emptyset}^{\eu m}$
tiene n\'ucleo $U_\pK$ y $\*{K_\pK}/U_\pK\cong {\ma Z}$.

\noindent
(g) Se sigue del Teorema de Aproximaci\'on de Artin. $\fin$
\end{proof}

Ahora bien, puesto que $[C_K:C_{K,\emptyset}^{\eu m}]=[C_K:
\CKm]=\infty$, se pedir\'a en lo futuro que $S\neq
\emptyset$. En este caso veremos que $[C_K:\CKSm]<\infty$,
el cual es un dominio Dedekind.

Sea $S$ un conjunto no vac{\'\i}o de lugares, $S\neq \emptyset$
y sea $\o_S=\cap_{\pK\notin S} \o_\pK=\{x\in K\mid v_\pK(x)\geq
0\ \forall\ \pK\notin S\}$.

\begin{definicion}\label{CClaseD4.9.5}\label{D17.7.34N}
Sea $\mK$ cualquier $S$--m\'odulus
con $S\neq \emptyset$. Definimos el {\em $S$--grupo de clases 
m\'odulo $\mK$\index{S--grupo de clases modulo un 
modulus@$S$--grupo de clases m\'odulo un
m\'odulus}}, $\ClKSm$\label{CClaseClKSm} como el cociente 
de los ideales fraccionarios de $\o_S$ primos relativos a $\mK$
m\'odulo los ideales principales $z\o_S$ con $z\in K^{\ast}\cap
J_{K,{\ma P}_K\setminus \sop(\mK)}^\mK$ (si $\mK=1$, 
$J_{K,{\ma P}_K\setminus \sop(\mK)}^\mK=J_{K,{\ma P}_K\setminus 
\emptyset}^1=J_K$).
\end{definicion}

En particular $Cl_K^1(\o_S)=Cl_K(\o_S)$ es el grupo de $S$--clases
usuales, es decir, el grupo de clases del dominio Dedekind
$\o_S$. En este caso sabemos que $h_S=|Cl(\o_S)|<\infty$.
De hecho se tiene la sucesi\'on exacta
\[
0\longrightarrow \frac{D_{K,0}(S)}{P_K(S)}\longrightarrow
\frac{D_{K,0}}{P_K}\cong I_{K,0}\longrightarrow
Cl(\o_S)\xrightarrow{\deg}
{\ma Z}/d{\ma Z}\longrightarrow 0,\label{CClaseDKS0}
\]
donde $D_K(S)$ son los divisores con soporte en $S$, $D_{K,
0}(S)=D_K(S)\cap D_{K,0}$ y $d=\deg (D_K^S)=\mcd\{\deg\pK\mid
\pK\in S\}$ y $P_K(S)=P_K\cap D_K(S)$.

\begin{teorema}\label{CClaseT4.9.6}\label{T17.7.35N}
 Sea $S\subseteq {\ma P}_K$, $S\neq
\emptyset$ y sea $\mK$ un $S$--m\'odulus. Entonces
\lasa
\item $\ClKSm\cong C_K/\CKSm$.

\item La sucesi\'on
\[
\*{\o_S}\longrightarrow\frac{J_{K,S}^1}{\JKSm}
\longrightarrow \frac{J_K}{K^{\ast}\JKSm}\longrightarrow \frac{J_K}
{K^{\ast} J_{K,S}^1}\longrightarrow 1
\]
es exacta,

\item $[C_K:\CKSm]$ es finito.
\end{list}
\end{teorema}
\begin{proof}

\noindent
(a) La funci\'on $\theta\colon J_{K,\pKK\setminus \sop(\mK)}^{\mK}
\longrightarrow \ClKSm$ dado por
\[
\theta(\vec \alpha)=\prod_{\pKK\setminus S}(\pK\cap \o_S)^{
v_\pK(\alpha_\pK)}\bmod \PKSm
\]
es suprayectiva y $\ker \theta=\JKSm\big(K^{\ast}\cap 
J_{K,{\ma P}_K \setminus \sop(\mK)}^{\mK}\big)=K^{\ast}\JKSm\cap
J_{K,{\ma P}_K\setminus \sop(\mK)}^{\mK}$.

Por tanto
\begin{align*}
\ClKSm&=\frac{J_{K,{\ma P}_K\setminus \sop(\mK)}^{\mK}}{
K^{\ast}\JKSm\cap J_{K,{\ma P}_K\setminus \sop(\mK)}^{\mK}}=
\frac{K^{\ast} J_{K,{\ma P}_K\setminus \sop(\mK)}^{\mK}}{K^{\ast} \JKSm}\\
&=\frac{J_K}{K^{\ast}\JKSm}\cong \frac{C_K}{\CKSm}.
\end{align*}

\noindent
(b) Se deja al cuidado del lector.

\noindent 
(c)  Se sigue de (b) pues tanto $J_{K,S}^1/\JKSm$ como
$J_K/K^{\ast} J_{K,S}^1 \cong
Cl(\o_S)$ son grupos finitos. $\fin$
\end{proof}

Repetimos la demostraci\'on de lo fundamental del
Teorema \ref{T17.7.35N} con una notaci\'on m\'as simple.
Sea $S$ un conjunto finito, $\o_S=\bigcap_{\pK\notin S}\o_{\pK}=
\{x\in K\mid v_{\pK}(x)\geq 0\text{\ para toda $\pK\notin S$}\}$ el
cual es un dominio Dedekind y definimos $I_{K,S}=\frac{D_K^S}
{P_K^S}$ donde $D_K^S=\frac{D_K}{\langle S\rangle}$ es el
grupo de divisores de $K$ con soporte en el complemento de 
$S$. Se tiene $D_{K}=\langle \pK\mid \pK\in{\ma P}_K\setminus S\rangle$,
y $P_K^S$ son los divisores principales de $\o_S$.

Por el teorema de F.K. Schmidt, $h_S=|I_{K,S}|<\infty$ 
(Corolario \ref{CRamDed1.2.6}). Sean
\begin{gather*}
\Lambda: J_K\longtwoheadrightarrow D_K^S,\quad
\Lambda(\vec\alpha)=\prod_{\pK\notin S}\pK^{v_{\pK}(\alpha_{\pK})}
\quad\text{y}\quad \pi:D_K^S\lra I_{K,S}
\intertext{la proyecci\'on natural. Sea}
\psi=\Lambda\circ\pi:J_K\longtwoheadrightarrow I_{K,S},\quad
\psi(\vec\alpha)=\prod_{\pK\notin S}\pK^{v_{\pK}(\alpha_{\pK})} \bmod P_K^S.
\end{gather*}

Sea $\vec\alpha\in\ker\psi$. Entonces $(\pi\circ\Lambda)(\vec\alpha)=1$,
por lo que $\Lambda(\vec\alpha)\in P_K^S$. Por tanto existe $x
\in\*K$ con $\prod_{\pK\notin S}\pK^{v_{\pK}(\alpha_{\pK})}=
\prod_{\pK\notin S}\pK^{v_{\pK}(x)}$, donde se sigue que
$v_{\pK}(\alpha_{\pK}x^{-1})=0$ para toda $\pK\notin S$. De esta forma
obtenemos $\vec\alpha x^{-1}\in J_{K,S}$. Por tanto $\vec\alpha\in J_{K,S}
\*K$. El rec\'iproco tambi\'en se cumple. Se sigue que $\ker\psi=J_{K,
S}\*K$.

Por tanto
\[
I_{K,S}\cong \frac{J_K}{J_{K,S}\*K}\cong \frac{J_K/\*K}{J_{K,S}\*K/\*K}\cong
\frac{C_K}{C_{K,S}}\quad\text{y por tanto}\quad I_{K,S}\cong\frac{C_K}{C_{K,S}}.
\]

Ahora $C_{K,S}=C_{K,S}^1$, por lo que $[C_K:C_{K,S}^1]<\infty$. Adem\'as
\begin{gather*}
[C_K:C_{K,S}^{\mK}]=[C_K:C_{K,S}^1][C_{K,S}^1:C_{K,S}^{\mK}]=
|I_{K,S}|[C_{K,S}^1:C_{K,S}^{\mK}].
\intertext{Adem\'as}
[J_{K,S}^1:J_{K,S}^{\mK}]=\prod_{\pK\notin S}[U_{\pK}:\unidadess m]
=\prod_{i=1}^r[U_{\pK_i}:U_{\pK_i}^{(\alpha_i)}]<\infty\quad\text{y}\\
[C_{K,S}^1:C_{K,S}^{\mK}]\leq[J_{K,S}^1\*K/\*K:J_{K,S}^{\mK}\*K/\*K]
\leq [J_{K,S}^1:J_{K,S}^{\mK}]<\infty.
\end{gather*}

Se sigue que $[C_K:C_{K,S}^{\mK}]<\infty$.

Esto es, nuevamente hemos obtenido

\begin{teorema}\label{T17.7.36N}
Sean $K$ un campo de funciones, $S$ un 
conjunto finito no vac\'io de primos
de $K$ y $\mK=\prod_{i=1}^r \pK_i^{\alpha_i}$ 
un m\'oduls tal que $S\cap
\{\pK_1,\ldots,\pK_r\}=\emptyset$. Sea $C_{K,S}^{\mK}
=J_{K,S}^{\mK}\*K/\*K$,
donde
\[
J_{K,S}^{\mK}=\prod_{\pK\in S}\*K_{\pK}\times
\prod_{\substack{\pK
\notin S\cup\\ \{\pK_1,\ldots,\pK_r\}}}U_{\pK}
\times \prod_{i=1}^r U_{\pK_i}^{(\alpha_i)}.
\]
Entonces $[C_K:C_{K,S}^{\mK}]<\infty$.
$\fin$
\end{teorema}

\begin{definicion}\label{D17.7.37N}
Por el teorema de existencia, se define $\KSm$\label{CClaseKSm}
como {\em campo de clase asociado al grupo de congruencia
$\CKSm$\index{campo de clase asociado a grupos de
congruencias}} y se llama el {\em campo de clases de $S$--rayos
m\'odulo $\mK$\index{campo de clases de $S$--rayos}}. Se tiene
\[
\Gal(\KSm/K)\cong J_K/K^{\ast} \JKSm\cong C_K/\CKSm\cong
Cl^\mK_K (\o_S).
\]
\end{definicion}

\begin{teorema}\label{T17.7.38N}
En $K_S^{\mK}$, $\pK_1,\ldots,\pK_r$ son los primos ramificados y si 
$\pK\in S$, entonces $\pK$ se descompone totalmente en $K_S^{\mK}/K$.
\end{teorema}

\begin{proof}
Se sigue de que $J_{K,S}^{\mK}=\prod_{\pK\in S}\*K_{\pK}\times\prod_{\substack{\pK
\notin S\cup\\ \{\pK_1,\ldots,\pK_r\}}} U_{\pK}
\times \prod_{i=1}^r U_{\pK_i}^{(\alpha_i)}$.
$\fin$
\end{proof}

\begin{observacion}\label{CClaseO4.9.7}\label{O17.7.39N}
$K_S^1$ es el campo de 
clase de Hilbert de $\o_S$ en el sentido de Rosen, con la 
peque\~na generalizaci\'on que aqu{\'\i}, de momento, $S$
puede ser infinito. Es decir, $K_S^1/K$ es la m\'axima extensi\'on
abeliana no ramificada de $K$ tal que todos los primos de
$S$ se descomponen totalmente. Este es el tercer an\'alogo al
campo de clase de Hilbert en campos de funciones.
\end{observacion}

Cuando $S=\{\p\}$ consiste de un solo primo, $\KSm$ se
puede construir en t\'ermino de m\'odulos de Drinfeld de rango
$1$ (ver Hayes \cite{Hay79} y el Cap\'itulo \ref{DrinfeldCh15}).

\begin{proposicion}\label{CClaseP4.9.11}\label{P17.7.40N}
Si $S$ y $T$ son dos subconjuntos
no vac{\'\i}os de ${\ma P}_K$, $\mK$ es un $S$--m\'odulos
y ${\eu n}$ es un $T$--m\'odulus, entonces
\lasa
\item Si $S\supseteq T$ y $\mK|{\eu n}$, entonces $\KSm
\subseteq K_T^{\eu n}$.

\item $\KSm\cap K_T^{\eu n}=K_{S\cup T}^{\mcd(\mK,{\eu n})}$.

\item Si $T\subseteq S$ y $\mK|{\eu n}$ y ${\eu n}$ tambi\'en
es un $S$--m\'odulus, entonces 
\[
[K_S^{\eu n}:\KSm]\leq
[K_T^{\eu n}:K_T^\mK].
\]
\end{list}
\end{proposicion}
\begin{proof}

\noindent
(a) Se tiene que $C_{K,T}^{\eu n}\subseteq \CKSm$ por lo que
$\KSm\subseteq K_T^{\eu n}$.

\noindent
(b) $\CKSm C_{K,T}^{\eu n}=C_{K,S\cup T}^{\mcd({\eu m},{\eu n})}$
por lo tanto $\KSm\cap K_T^{\eu n}=K_{S\cup T}^{\mcd({\eu m},{\eu n})}$.

\noindent
(c) Los mapeos de Artin $(\underline{\ \ },\KSm/K)$ y $(\underline{\ \ },
K_T^{\eu n}/K)$ inducen isomorfismos $K^{\ast} \JKSm/K^{\ast} J_{K,S}^{\eu n}
\cong \Gal(K_S^{\eu n}/\KSm)$ y $K^{\ast} J_{K,T}^\mK/K^{\ast} J_{K,T}^{\eu n}
\cong \Gal(K_T^{\eu n}/K_T^\mK)$ y el mapeo
\[
K^{\ast} J_{K,T}^\mK/K^{\ast} J_{K,T}^{\eu n} \longrightarrow K^{\ast} \JKSm/K^{\ast} J_{K,S}^{\eu n}
\]
es suprayectivo. $\fin$
\end{proof}

Debido al Teorema \ref{CClaseT4.5.10}, tenemos la siguiente definici\'on.

\begin{definicion}\label{CClaseD4.9.9-1}\label{D17.7.41N}
Sea $L/K$ una extensi\'on abeliana finita de campos de funciones globales.
Entonces definimos el {\em conductor} de la extensi\'on
$L/K$ como
\[
{\eu f}={\eu f}_{L/K}={\eu f}(L/K)=\prod_{\pKK}{\eu f}_{\pK},
\]
donde ${\eu f}_{\pK}$ denota al conductor de la 
extensi\'on de campos locales $L_{\pK}/K_{\pK}$.
\end{definicion}

\begin{teorema}[Teorema del Conductor\index{conductor!teorema
del $\sim$}\index{teorema del conductor}]\label{CClaseT4.9.9}\label{T17.7.42N}
Sea $L/K$ una extensi\'on abeliana finita de campos de funciones
y sea $S$ un subconjunto no vac{\'\i}o de 
\[
\Spl(L/K)=\{\pKK\mid \pK\text{\ se descompone totalmente en \ }
L/K\}.
\]
Entonces el conductor $\f{}=\f{} (L/K)$ es el $S$--m\'odulus m{\'\i}nimo
$\mK$ tal que $L\subseteq K_S^{\mK}$.
\end{teorema}

\begin{proof} Sea $\me$ un $S$--m\'odulus. Por  definici\'on de conductor
local $\f\pK(L/K)|\pK^{n_\pK}\iff G^{n_\pK}(L/K)=1\iff \unidadess n
\subseteq \N_{L_{\pK}/K_{\pK}} \*L_{\pK}\iff \enc {\unidadess n}{\pK}
\subseteq (\N_{L/K} J_L)\*K$.

Adem\'as $\pK\in S$ implica que $\pK$ se descompone totalmente o,
equivalentemente, $\enc {\*K_{\pK}}{\pK}\subseteq (\N_{L/K} J_L)\*K$.
Puesto que $(\N_{L/K} J_L)\*K$ es un subgrupo abierto de {\'\i}ndice
finito en $J_K$ y por lo tanto cerrado, y puesto que $J_{K,S}^{\mK}$ es
topol\'ogicamente generado por sus subgrupos $\{\enc {\*K_{\pK}}{\pK}
\mid \pK\in S\}$ y $\{\enc {\unidadess n}{\pK}\mid \pK\notin S\}$,
se sigue que $\f{}(L/K)|\mK\iff
\JKSm\subseteq (\N_{L/K} J_L)\*K\iff \CKSm\subseteq \N_{L/K} C_L\iff
L\subseteq \KSm$.
$\fin$
\end{proof}

\begin{teorema}\label{CClaseT4.9.10}\label{T17.7.43N}
Sean $S$ un conjunto finito no vac{\'\i}o de lugares de
$K$ y $\mK$ un $S$--m\'odulus.
Entonces $\KSm$ es la m\'axima extensi\'on abeliana $L/K$ tal que
todos los primos de $S$ son totalmente descompuestos y $\f{}(L/K)|\mK$.

Adem\'as, si $d=\deg S:=\mcd \{\deg \pK\mid \pK\in S\}$, ${\ma F}_{q^d}$
es el campo de constantes de $\KSm$.
\end{teorema}

\begin{proof} Sea $L/K$ una extensi\'on abeliana. 
Si todos los primo de $S$ se
descomponen en $L/K$ y $\f{}(L/K)|\mK$, entonces, como
consecuencia del Teorema \ref{T17.7.42N}, $L\subseteq \KSm$.

Rec{\'\i}procamente, si $L\subseteq \KSm$, entonces $\f{}(L/K)|\mK$ y por
definici\'on de $\JKSm$, todos los primos de $S$ se descomponen
totalmente en $\KSm$. Ahora, si $\pK\in S$
se descompone la extensi\'on de constantes de grado $r$, entonces
$r\mid \deg \pK$. Por lo tanto todos los
primos de $S$ se descomponen totalmente en 
$K{\ma F}_{q^r} \iff r | \mcd\{\deg \pK\mid \pK\in S\}$. Se sigue
que ${\ma F}_{q^d}$ es el campo de constantes de $K_S^{\mK}$.
$\fin$
\end{proof}

\begin{corolario}\label{CClaseC4.9.11}\label{C17.7.44N}
(Ver el Teorema {\rm{\ref{T12*.2.2.A}}}). Sea $K$ un campo global de funciones
con campo de constantes $\F$
y sea $S
\subseteq {\ma P}_K$ no vac{\'\i}o. Sea $d=\deg S=\mcd \{\deg \pK\mid
\pK\in S\}$. Entonces el campo de constantes del campo de clase de
Hilbert $K_S^1$ es ${\ma F}_{q^d}$.

En particular, si $K=\F(T)$, entonces $K_S^1={\ma F}_{q^d}(T)
=K{\ma F}_{q^d}$.
\end{corolario}

\begin{proof} La \'ultima parte se sigue de que $I_{\F(T),0}=\{1\}$, lo cual
implica que $\nr {\F(T)}=\F(T)\abe \F$ y puesto que $\F(T)_S^1/\F(T)$
es no ramificado, se sigue que $\F(T)_S^1\subseteq \F(T)\abe \F$,
es decir $\F(T)_S^1/\F(T)$ es una extensi\'on de constantes.
$\fin$
\end{proof}

\begin{observacion}\label{O17.7.45N}
El campo $K_S^1$ es lo que com\'unmente se llama el {\em campo de clase
de Hilbert} con respecto a $S$.
\end{observacion}

\subsubsection{An\'alogos al campo de clase de Hilbert para campos
de funciones}\label{CClaseS4.8.1}

Como hemos mencionado anteriormente, si $K^{\rm{nr}}$ es la m\'axima
extensi\'on abeliana no ramificada de $K$, se tiene que $K\abe\F
\subseteq K^{\rm{nr}}$. Sea ${\mathcal T}=\Gal(\abe K/K^{\rm{nr}})$, $K^{\rm{nr}}=
(\abe K)^{\mathcal T}\supseteq (\abe K)^{{\mathcal H}_0}=K\abe\F$, es decir,
${\mathcal H}_0\supseteq {\mathcal T}$, donde ${\mc H}_0=\Gal(\abe K/K\abe\F)$.
Se tienen los diagramas conmutativos (ver Teorema \ref{CClaseT4.2.10}):
\begin{gather*}
\xymatrix{
1\ar[r]& C_{K,0}\ar@{^{(}->}[r]\ar[d]_{\cong}^{\rho_K}& C_K
\ar@{-}[r]^{\deg}\ar[d]^{\rho_K}&{\ma Z}\ar[r]\ar@
{^{(}->}[d]^{\rho_{\F}}&0\\
1\ar[r]&{\underbracket[0pt]{\Gal(\abe K/K \abe \F)}_{\substack{\ucong\\
 {\mathcal H}_0}}}\ar[r]&{\underbracket[0pt]{\Gal(\abe K/K)}_{\substack{\ucong\\
 {\mathcal G}=\abe G_K}}}\ar[r]_{\rest}&
 {\underbracket[0pt]{\Gal(K\abe \F/K)}_{\substack{\ucong\\ 
 \Gal(\abe \F/\F)}}}\ar[r]&1}\\
\xymatrix{
1\ar[r]& I_{K,0}\ar@{^{(}->}[r]\ar[d]^{\cong}&I_K\ar[r]^{\deg}
\ar[d]^{\mu}&{\ma Z}\ar@{^{(}->}[d]^{\rho_{\F}}\ar[r]&0\\
1\ar[r]&\Gal(K^{\rm{nr}} /K \abe\F)\ar@{^{(}->}[r]&
{\underbracket[0pt]{\Gal(K^{\rm{nr}} /K)}_{\substack{\ucong\\
{\mc G}/{\mc T}}}}\ar[r]_{\rest}&
{\underbracket[0pt]{\Gal(K\abe\F/K)}_{\substack{\ucong\\ \Gal(
\abe\F/\F)}}}\ar[r]&1}
\end{gather*}
pues $\frac{C_{K,0}}{\*KU_K/\*K}\cong I_{K,0}\cong\frac{\Gal(\abe K/K\abe\F)}
{\Gal(\abe K/\nr K)}$ y $\frac{C_K}{U_K\*K/\*K}\cong I_K\cong\frac{\Gal(
\abe K/K)}{\Gal(\abe K/\nr K)}\cong \Gal(\nr K/K)$.

Ahora, si $\mu(\alpha)=1$, $(\rest \mu)(\alpha)=\rho_{\abe\F}(
\deg \alpha)=1$ por lo que $\deg \alpha=0$ y por tanto
$\alpha\in I_{K,0}$. Se sigue que $\mu$ es inyectiva.

En particular $K^{\rm{nr}}/K\abe\F$ es una extensi\'on finita y 
$\Gal(K^{\rm{nr}} /K\abe\F)\cong I_{K,0}$.

Se tiene 
\[
\xymatrix{
&K^{\rm{nr}}\ar@{<-}[d]^{{\mathcal T}_0=\Gal(K^{\rm{nr}}/K\abe\F)\cong I_{K,0}}\\
K\ar[r]_{\substack{\uparrow\\ \Gal(K\abe \F/K)\cong \\
\Gal(\abe \F/\F)\cong \hat{\ma Z}=\overline{\langle\Fro K\rangle}}}\ar[ru]&
K\abe\F}
\]

Se sigue que la sucesi\'on
\[
1\longrightarrow {\mathcal T}_0\xrightarrow{\mu}{\mathcal R}=
\Gal(K^{\rm{nr}}/K)\xrightarrow[\rest]{\varphi}\hat{\ma Z}\longrightarrow 0
\]
es exacta.

\subsubsection{Primer an\'alogo al campo de clase de Hilbert}
\label{CClaseS4.8.1.1}

Se tiene que
\[
1\longrightarrow \Gal(K^{\rm{nr}}/K\abe \F)\longrightarrow\Gal(K^{\rm{nr}}/K)
\longrightarrow \Gal(K\abe\F/K)\longrightarrow 1
\]
es exacta, que $\Gal(K^{\rm{nr}}/K\abe \F)\cong I_{K,0}$ es un grupo
finito y adem\'as $I_K\stackrel{\mu}{\longrightarrow}
\Gal(K^{\rm{nr}}/K)$ es un inyecci\'on. Ahora bien, $\mu$ no es
un isomorfismo pues $I_K$ es discreto y no finito y $\Gal(\nr K/K)$
es profinito, compacto. De hecho $I_K\cong I_{K,0}\times{\ma Z}$.
Se tiene $\hat I_K\cong I_{K,0}\times \hat{\ma Z}\cong \Gal(\nr K/K)$.

Podemos considerar a $K^{\rm{nr}}$ como el primer an\'alogo
al campo de clase de Hilbert.

\subsubsection{Segundo an\'alogo al campo de clase
de Hilbert}\label{CClaseS4.8.1.2}

Consideremos extensiones abelianas no ramificadas de 
$K$ que tienen como campo 
de constantes al mismo campo $\F$.

Sean $U_K=\prod_{\pKK} U_\pK$ y $J_{K,\emptyset}^1$ ($\mK=1$, $S=
\emptyset$), el grupo de rayos m\'odulus $1$ con respecto $S=\emptyset$
y sea
\[
C_{K,\emptyset}^1=C_K^1=U_K\*K/\*K\subseteq C_{K,0}.
\]

Se tiene que $\nr K$ corresponde
a $U_K$ o a $U_K\*K/\* K$ (ver despu\'es 
del Teorema \ref{CClaseT4.2.10}), esto es
$C_K^1\cong \Gal(\abe K/\nr K)$. Adem\'as $\hat C_K=\abe
G_K=\Gal(\abe K/K)$.

Sea $\tilde{\vec b}\in C_K$ tal
que $\deg (\tilde{\vec b})=1$. Entonces el subgrupo 
$B:=C_K^1\times \langle \tilde{\vec b}\rangle$ es abierto en $C_K$ y 
puesto que $C_K\cong  C_{K,0}\times\langle \tilde{\vec b}\rangle$, 
entonces 
\[
\frac{C_K}{B}\cong\frac{C_{K,0}\times
\langle \tilde{\vec b}\rangle}{C_K^1\times
\langle \tilde{\vec b}\rangle}\cong\frac{C_{K.0}}
{C_K^1}=\frac{C_{K,0}}{C_K^1}\cong I_{K,0}
\]
y se sigue que $[C_K:B]=
[C_{K,0}:C_K^1]=h_K=|I_{K,0}|$ es el n\'umero de clase
de $K$.
\[
\xymatrix{
\bullet\ar@{--}[rr]^{\langle\tilde{\vec b}\rangle}
\ar@{--}[d]&&\abe K\ar@{-}[dll]_B\ar@{-}[d]^{C_K^1}\\
K_B\ar@{-}[rr]_{\langle\tilde{\vec b}\rangle} 
\ar@{-}[d]_{I_{K,0}}&&K^{\rm{nr}}\ar@{-}[d]^{I_{K,0}}\\
K\ar@{-}[rr]&&K\abe\F
}
\]

El campo de clase $K_B$ de $B$ es el segundo an\'alogo al campo
de clase de Hilbert. Puesto que $B$ contiene una copia de
${\ma Z}$ ($\longleftrightarrow \langle\Fro K\rangle=\langle 
\tilde{\vec b}\rangle$)
y tales copias de ${\ma Z}$ est\'an indexadas por $C_{K,0}$
puesto que $C_K\cong C_{K,0} \times {\ma Z}$ y $C_K^1
\longleftrightarrow \prod_{\pKK} U_\pK$ es la extensi\'on no
ramificada $\nr K$ de $K$,
se sigue que este campo $K_B$, denotado por $K_H^{(1)}$,
es una extensi\'on abeliana finita no ramificada de $K$ que tiene
como campo de constantes a $\F$ (ver Teorema \ref{CClaseT4.8.1-1.1}).

Veamos que $K_H^{(1)}$ es maximal
en el sentido de ser no ramificada con campo de constantes $\F$.
Sea $K_H^{(1)}\subseteq
\hat{K}$ (alg\'un $\hat{K}= \tilde{K}(\vec \alpha)$, 
ver principio de esta subsecci\'on) y
se tiene que $H=\N_{K_H^{(1)}/K} C_{K_H^{(1)}}\subseteq
C_K$. Ahora bien, las extensiones no ramificadas deben
satisfacer que sus grupos de normas est\'an contenidas en $U$ para
que su conductor local sea $1$ para todo lugar $\pK\in{\ma P}_K$
y entonces su conductor global sea $1$. 

Por otro lado se tiene el diagrama
\[
\xymatrix{
\hat{K}\ar@{-}[rr]^{\langle\tilde{\vec b}_1\rangle}
\ar@{-}[d]\ar@{-}@/^1pc/[d]^{C_K^1}&&\abe K
\ar@{-}[d]\ar@{-}@/^1pc/[d]^{C_K^1}\\
E\ar@{-}[rr]\ar@{-}[d]&&K^{\rm{nr}}\ar@{-}[d]\ar@{-}@/^1pc/[d]^{
I_{K,0}}\\
K\ar@{-}[rr]&&K\abe\F
}
\]
Por la correspondencia de Galois, $K^{\rm{nr}}$ corresponde a $E
\subseteq \hat{K}$ con $E=\hat{K}\cap K^{\rm{nr}}$ y $[E:K]=[K^{\rm{nr}}:
K\abe \F]=h=|I_{K,0}|$.

Ahora bien, $E$ tiene como campo de constantes $\F$ lo mismo
que $K_H^{(1)}$ pues $K_H^{(1)}$ corresponde a $B$, esto es,
$K_H^{(1)}\longleftrightarrow B=C_K^1\times\langle \tilde{\vec b}
\rangle$ con $C_K^1\cong\Gal(\abe K/K^{\rm{nr}})$
y como $\hat{K}$ es fijado por $\tilde{\vec b}_1$ y por tanto por 
$\langle\tilde{\vec b}_1\rangle^{\hat{\ma Z}}$,
se tiene que $\hat{K}=\big(\abe K\big)^{\langle{\tilde{\vec b}_1}
\rangle^{\hat{\ma Z}}}$, de donde, $K_H^{(1)}=E$ y por lo tanto
$K_H^{(1)}$ es maximal con respecto a ser no ramificada sobre
$K$ y tener a $\F$ como campo de constantes. De hecho
deber\'iamos haber escrito $K_H^{(\tilde{\vec b}_1)}\longleftrightarrow
\hat{K}_{\tilde{\vec b}_1}$ pues tanto $K_H^{(1)}$ como $\hat{K}$ est\'an
determinados por $\tilde{\vec b}_1$.
\[
\begin{minipage}{4cm}
\xymatrix{
\bullet\ar@{-}[rr]^{\langle\idel b_1\rangle}\ar@{-}[d]&&\abe K\ar@{-}[d]\\
L\ar@{--}[rr]\ar@{-}[d]&& L K\abe\F\ar@{-}[d]\\
K_B=K_H^{(1)}\ar@{-}[uurr]_B\ar@{-}[rr]\ar@{-}[d]&& \nr K\ar@{-}[d]\\
K\ar@{-}[rr]&&K\abe \F
}
\end{minipage}
\qquad\qquad
\begin{minipage}{4cm}
M\'as precisamente, se tiene $K_H^{(1)}\cap K\abe\F=K$ y $K_H^{(1)}
K\abe \F=K_H^{(1)}\abe\F=\nr K$. 

Veamos que $K_H^{(1)}$ es maximal con respecto a estas propiedades.
\end{minipage}
\]
Sea $L/K$ una extensi\'on abeliana, donde 
el campo de constantes de $L$ es $\F$
y tal que $K_H^{(1)}\subseteq L$.
Entonces $L\cap K\abe \F=K$ y $LK\abe\F\supseteq K_H^{(1)}K\abe\F=\nr K$.
Como $L/K$ es no ramificada, $L K\abe \F=\nr K$ y por tanto $LK\abe \F=
K_H^{(1)}K\abe\F$ y de la teor\'ia de Galois se sigue que $L=K_H^{(1)}$.

Pongamos $K_H^{(1)}=K_H^{(\idel b_1)}$.
Veamos que hay {\underline{exactamente}} $h=h_K=|I_{K,0}|$ de estos
campos: extensiones abelianas no ramificados de $K$
con campo de constantes $\F$ maximales, 
y todos ellos tienen a su grupo de Galois isomorfo a
\[
C_K/B\cong C_{K,0}/C_K^1\cong I_{K,0}.
\]

Por el Observaci\'on \ref{CClaseO1.4.2'},
tenemos que si $L/K$ es una extensi\'on de Galois, finita
o infinita, con grupo de Galois $G$, se tiene que si $H<G$ y si $\bar H$ es la
cerradura de $H$, se tiene $L^H=L^{\bar H}$.

Sean $\tilde{\vec u}_1,\ldots, \tilde{\vec u}_h$, 
$h$ representantes de $C_{K,0}/
C_K^1\cong I_{K,0}$ 
y sean $\tilde{\vec b}_i:=\tilde{\vec b}\tilde{\vec u}_i$, 
$1\leq i\leq h$. Los grupos
$B_i:= C_K^1\times\langle \tilde{\vec b}_i\rangle$ son todos
distintos y $[C_K:B_i]=h$. Los grupos $B_i$ dan lugar a 
$h$ campos que son maximales en el sentido de ser extensiones
abelianas no ramificadas de $K$ y tener como campo de constantes
$\F$ (son maximales puesto que si $K_{H}^{(i)}$ es el campo
asociado a $B_i$, entonces $K_{H}^{(i)}\abe \F=K^{\rm{nr}}$).

\[
\xymatrix{
K_{{\tilde{\vec b}}_h}\ar@{-}[rrrrrd]^{{\tilde{\vec b}}_h\cong 
{\hat{\ma Z}}}\ar@{--}[rd]\ar@{-}[ddd]_{C_K^1}\\
&K_{{\tilde{\vec b}}_2}\ar@{-}[rrrr]^{{\tilde{\vec b}}_2\cong
 {\hat{\ma Z}}}\ar@{-}[rd]\ar@{-}[ddd]_{C_K^1}&&&&\abe K
\ar@{-}[ddd]^{C_K^1}\ar@{-}@/^3pc/[ddddd]^{C_{K,0}}\\
&&K_{{\tilde{\vec b}}_1}\ar@{-}[rrru]^{{\tilde{\vec b}}_1\cong {\hat{\ma Z}}}\ar@{-}[ddd]
|!{[lld];[ddrrr]}\hole_{C_K^1}\\
K_H^{(h)}\ar@{-}[rrrrrd]|!{[ruu];[dr]}\hole^{\hat{\ma Z}}
 \ar@{--}[rd]\ar@{-}[dddr]_{I_{K,0}}\\
&K_H^{(2)}\ar@{-}[dd]^{I_{K,0}}\ar@{-}[dr]
\ar@{-}[rrrr]|!{[ruu];[dr]}\hole_{\hat{\ma Z}}
 &&&&K^{\rm{nr}}\ar@{-}[dd]^{I_{K,0}}\\
&&K_H^{(1)}\ar@{-}[rrru]_{\hat{\ma Z}}\ar@{-}[dl]^{I_{K,0}}\\
&K\ar@{-}[rrrr]^{\hat{\ma Z}}&&&&K\abe \F
}
\]

Sea $\tilde{\vec a}\in C_K$ con $\deg \tilde{\vec a}=1$. 
Por tanto $\deg(\tilde{\vec a}\tilde{\vec b}^{-1})=0$
por lo que $\idel a$ pertenece a alguna clase 
$\tilde{\vec b}\tilde{\vec u}_i\bmod C_K^1$.
As{\'\i}, el grupo $A= C_K^1\times\langle \tilde{\vec a}\rangle$ da
lugar a un campo de clase de Hilbert como los anteriores y se tiene
\[
\frac{C_K}{A}=\frac{C_{K,0}\times\langle \tilde{\vec b}\rangle}{
 C_K^1\times \langle \tilde{\vec a} \rangle}
=\frac{C_{K,0}\times\langle \tilde{\vec b}\rangle}
{C_K^1\times\langle \tilde{\vec b}\tilde{\vec u}_i\rangle}
=\frac{C_K}{B_i}.
\]
Por tanto $A=B_i$.

Los $h$ campos asociados a $B_1,\ldots, B_h$ son todos los
campos con campo de constantes $\F$, que son extensiones
abelianas no ramificadas maximales de $K$. Sean
$K_H^{(1)}, \ldots, K_H^{(h)}$ estos $h$ campos. El campo de
clase correspondiente a $C_{K,0}\times
\langle \tilde{\vec b}^h\rangle$
es la extensi\'on de constantes de grado $h$ de $K$
(ver el principio de esta subsecci\'on). Esto es,
$L:=K{\ma F}_{q^h}$ es el campo correspondiente a
$C_{K,0}\times\langle\tilde{\vec b}^h\rangle$. 
\begin{gather*}
\xymatrix{
&&&\abe K\ar@{-}[d]^{C_{K,0}}\ar@{-}[dll]_{ C_{K,0}
\times\langle\tilde{\vec b}^h\rangle}\\
K\ar@{-}[r]^h&K{\ma F}_{q^h}\ar@{-}[rr]^{\langle\tilde{\vec b}^h
\rangle}&&K\abe \F}
\end{gather*}
\begin{gather*}
\xymatrix{
\tilde{K}\ar@{-}[rrr]^h_{\langle\overline{\tilde{\vec b}}\rangle
=\langle \tilde{\vec b}\bmod \tilde{\vec b}^h\rangle}\ar@{-}[d]_{C_{K,0}}
\ar@{-}@/^2pc/[rrrrrr]^{\langle\tilde{\vec b}\rangle}&&&
\tilde{K}{\ma F}_{q^h}\ar@{-}[rrr]_{\langle\tilde{\vec b}^h\rangle}\ar@{-}[d]
&&&\abe K\ar@{-}[d]^{C_{K,0}}\\
K\ar@{-}[rrr]_h^{\langle\overline{\tilde{\vec b}}\rangle
=\langle \tilde{\vec b}\bmod \tilde{\vec b}^h\rangle}
\ar@{-}@/_2pc/[rrrrrr]_{\langle\tilde{\vec b}\rangle}
&&&K{\ma F}_{q^h}\ar@{-}[rrr]^{\langle\tilde{\vec b}^h\rangle}
&&&K\abe \F
}
\end{gather*}

Afirmamos que
\[
LK_H^{(i)}=LK_H^{(j)}=LK_H^{(1)}\cdots K_H^{(h)}\quad
\text{para todos}\quad 1\leq i, j\leq h.
\]

En efecto, puesto que $\tilde{\vec b}^h=
(\tilde{\vec b}_i\tilde{\vec u}_i^{-1})^h=
\tilde{\vec b}_i^h \tilde{\vec u}_i^{-h}
\equiv \tilde{\vec b}_i^h\bmod C_K^1$, entonces
$\big(C_{K,0}\times\langle \tilde{\vec b}^h\rangle \big)\cap
\big(C_K^1\times\langle \tilde{\vec b}_i\rangle\big)=C_K^1\times
\langle\tilde{\vec b}_i^h\rangle$ y
$\tilde{\vec b}_i^h, \tilde{\vec b}_j^h$
est\'an en la misma clase
$\bmod\ C_K^1$.

Se sigue que $\big( C_{K,0}\times\langle \tilde{\vec b}^h\rangle\big)\cap
\big(C_K^1\times \langle \tilde{\vec b}_i\rangle\big)=
\big(C_{K,0}\times \langle \tilde{\vec b}^h\rangle\big)\cap
\big(C_K^1\times\langle \tilde{\vec b}_j\rangle\big)$ para todo
$1\leq i,j\leq h$.

Se tiene que $LK_H^{(i)}=LK_H^{(j)}$ y por tanto 
\begin{gather*}
\underbracket[0.5pt]{LK_H^{(1)}}K_H^{(2)}\cdots
K_H^{(h)}=\underbracket[0.5pt]{LK_H^{(2)}}K_H^{(2)}K_H^{(3)}
\cdots K_H^{(h)}=\underbracket[0.5pt]{LK_H^{(2)}}K_H^{(3)}\cdots K_H^{(h)}\\
=\underbracket[0.5pt]{LK_H^{(3)}}K_H^{(3)}\cdots K_H^{(h)}=LK_H^{(3)}\cdots
K_H^{(h)}=\cdots =LK_H^{(h)}.
\end{gather*}

Ahora bien, puesto que $L\cap K_H^{(i)}=K$ para toda $i$
por ser $L/K$ extensi\'on de constantes y $K_H^{(i)}/K$
geom\'etrica, se tiene
\begin{gather*}
\xymatrix{
K_H^{(i)}\ar@{-}[ddd]_h\ar@{-}[rrr]^h&&&LK_H^{(i)}
=LK_H^{(j)}\ar@{-}[ddd]^h\\ \\
&&K_H^{(j)}\ar@{-}[ruu]^h\ar@{--}@/_1pc/[uull]\ar@{--}@/^1pc/[dr]\\
K\ar@{-}[rrr]_h\ar@{-}[rru]^h&&&L=K{\ma F}_{q^h}
}
\intertext{$[LK_H^{(1)}\cdots K_H^{(h)}:K]=h^2$. Se tiene}
\Gal(LK_H^{(1)}\cdots K_H^{(h)}/K)\cong \Gal(L/K)\times
\Gal(K_H^{(1)}/K)\cong C_h\times I_{K,0}.
\end{gather*}

Este es el segundo an\'alogo al campo de clase de Hilbert.

\begin{teorema}\label{CClaseT4.8.4} Sea $h=h_K=|I_{K,0}|$ el n\'umero
de clase de $K$. Entonces hay exactamente $h$ campos 
$K_H^{(i)}$, $1\leq i\leq h$ con campo de constantes $\F$ y tales
que $K_H^{(i)}/K$ son extensiones abelianas no ramificadas de $K$
maximales con respecto a estas propiedades. Se tiene que
\[
\Gal(K_H^{(i)}/K)\cong I_{K,0}, \quad 1\leq i\leq h.
\]

Finalmente, si $L$ es la extensi\'on de constantes de $K$
de grado $h$, se tiene
\begin{align*}
\Gal(LK_H^{(i)}/K)&=\Gal(LK_H^{(j)}/K)=\Gal(LK_H^{(1)}\cdots
K_H^{(h)}/K)\\
&\cong \Gal(L/K)\times \Gal(K_H^{(i)}/K)
\cong C_h\oplus I_{K,0},
\end{align*}
para todas $1\leq i, j\leq h$. $\fin$
\end{teorema}

\begin{ejemplo}\label{E17.7.47N}
La \'unica extensi\'on abeliana geom\'etrica no ramificada de $\F(T)$
es $\F(T)$ misma.
\end{ejemplo}

\subsubsection{Tercer an\'alogo al campo de clase de 
Hilbert}\label{CClaseS4.8.1.3}

Una discusi\'on previa a lo que discutimos aqu\'i ya se present\'o
parcialmente en las Secciones \ref{S12*.1} y \ref{S12*.1-1}
del Cap\'itulo \ref{Ch12*} en donde estudiamos los
campos de g\'eneros. Ver particularmente
la Definici\'on \ref{D3.1}. En el Cap\'itulo \ref{Ch12*}
estuvimos interesados en el caso particular de
que el conjunto de primos eran los primos sobre
el divisor de polos de $T$. En la siguiente subsecci\'on, discutiremos
con m\'as detalle este tercer an\'alogo al campo de clase de Hilbert.

En el caso de un campo num\'erico $K$, $K_H$ es la m\'axima
extensi\'on abeliana de $K$ no ramificada y $K_H$ es el campo
de clase del grupo $J_{K,S}\*K/\*K$, donde $S$
es el conjunto de primos infinitos. En general se define
\[
J_{K,S}:= \prod_{\pK\notin S} U_\pK\times \prod_{\pK\in 
S}\*{K_\pK}\label{CClaseJKS2}.
\]

Es decir, al ser $K_H/K$ no ramificada en $S$, los
primos infinitos se descomponen totalmente en $K_H/K$.
Esto se puede copiar en el caso de campos de funciones
tomando como $S$ {\underline{cualquier conjunto
no vac{\'\i}o finito}} de lugares de $K$ y definiendo como $K_H$
a la m\'axima extensi\'on abeliana no ramificada de $K$, tal
que todos los lugares de $S$ se descomponen 
totalmente en $K_H/K$. Resulta ser que $K_H=K_S^1$ el campo
de clases de rayos del m\'odulus $1$ con respecto a $S$.
Esto es, $K_H$ corresponde a $J_{K,S}^1\*K/\*K=C_{K,S}^1$.

Veremos que $K_H/K$ es una extensi\'on finita y que 
$\Gal(K_H/K)\cong Cl(\o_S)$, el grupo de clases del anillo Dedekind
\[
\o_S=\{x\in K\mid v_{\pK}(x)\geq 0\text{\ para toda $\pK\notin S$}\}=
\bigcap_{\pK\notin S}\o_{\pK},
\]
donde para todo $S\neq \emptyset$, $Cl(\o_S)=\frac{
D_{K,S}}{P_{K,S}}=\frac{\text{divisores primos relativos a $S$}}
{\{(\alpha)_S\mid \alpha\in K^{\ast}\}}$ y donde para $\alpha\in K^{\ast}$, se
define $(\alpha)_S=\prod_{\pK\notin S}\pK^{v_\pK(\alpha)}$
(Teorema \ref{CClaseT4.9.6}).

M\'as precisamente, el grupo de ideales fraccionarios de $\o_S$
es $D_{K,S}=\frac{D_K}{D_S}=\frac{D_K}{\langle \pK\mid \pK\in S\rangle}
=\langle \pK\in{\ma P}_K\mid \pK\notin S\rangle$ y se tiene el epimorfismo
\begin{gather*}
\Lambda : J_K\lra D_K/D_S=\langle\pK\in{\ma P}_K\mid
\pK\notin S\rangle,\qquad \Lambda(\vec\alpha)=
\prod_{\pK\notin S}\pK^{v_{\pK}(\alpha_{\pK})}
\intertext{con n\'ucleo $J_{K,S}$, $Cl(\o_S)=\frac{D_K/D_S}{P_S}$ y el
epimorfismo inducido}
\tilde \Lambda:J_K\lra Cl(\o_S)
\intertext{satisface que $\ker\tilde\Lambda=J_{K,S}\*K$. Por tanto}
\frac{J_K}{J_{K,S}\*K}\cong \frac{J_K/\*K}{J_{K,S}/\*K}\cong
\frac{C_K}{C_{K,S}}=\frac{C_K}{C_{K,S}^1}\cong \Gal(K_S^1/K)\cong
Cl(\o_S).
\end{gather*}

En resumen, si $S$ es un conjunto finito no vac\'io de lugares y sea
\[
J_{K,S}:=\prod_{\pK\notin S}U_{\pK}\times \prod_{\pK\in S}\*{K_{\pK}},
\]
entonces la m\'axima extensi\'on abeliana de $K$ no ramificada
tal que los primos de $S$ se descomponen totalmente corresponde
a $J_{K,S}$ o, equivalentemente, a $J_{K,S}\*K/\*K:=C_{K,S}$.

Este ser\'a el tercer an\'alogo al campo de clase de Hilbert $K_H$
de $K$, el cual es igual a $K_S^1$ y es la m\'axima extensi\'on
abeliana no ramificada de $K$ tal que todos los primos de $S$
se descomponen totalmente en $K_S^1/K$ y se tiene
\[
\Gal(K_S^1/K)\cong C_K/C_{K,S}^1\cong Cl(\o_S).
\]
En la siguiente subsecci\'on retomamos este an\'alogo y es el
que consideraremos como el campo de clase de Hilbert para
campos de funciones.

\subsection{Campo de clase de Hilbert en
campos de funciones globales}\label{CClaseS4.9-1}

Como mencionamos antes, el estudio de esta subsecci\'on puede
ser consultado en \cite{Ros87} y en \cite{AngJau2000}.

Sea $K$ un campo global de funciones con campo de constantes
$\F$. Sea $S$ un conjunto finito no vac\'io de lugares de $K$.

\begin{definicion}[Ver Definici\'on \ref{D3.1}]\label{DClaseS4.9-1.1}
El {\em campo de clase de Hilbert relativo a $S$\index{campo de
clase de Hilbert}}, $K_{H,S}$, se define como la m\'axima 
extensi\'on abeliana de $K$ no ramificada y tal que los primos
en $S$ se descomponen totalmente en $K_{H,S}$.
\end{definicion}

Del Teorema \ref{T12*.2.2.A} tenemos que el campo de constantes
de $K_{H,S}$ es ${\ma F}_{q^d}$ donde $d=\mcd(\gr \pK\mid \pK\in S)$,
que $K_{H,S}/K$ es una extensi\'on abeliana finita y que
\[
\Gal(K_{H,S}/K)\cong Cl_S
\]
con $Cl_S=Cl(\AE S)$ donde $\AE S=\{x\in K\mid v_{\eu q}(x)\geq
0\text{\ para todo ${\eu q}\notin S$}\}$.

El isomorfismo $\Gal(K_{H,S}/K)\cong Cl_S$ est\'a dado por
el mapeo de reciprocidad de Artin: $({\underline{\phantom{xx} }},
K_{H,S}/K)=\artinp{K_{H,S}/K}{\ }\colon
Cl_S\lra\Gal(K_{H,S}/K)$.

Para cualquier campo $K$ y cualquier conjunto $S$ finito no vac\'io
de lugares de $K$, denotemos por $h_S$ al n\'umero de clases de
ideales $|Cl_S|$. Sea $L/K$ una extensi\'on finita y separable. Sea
$T:=\{\pL\in{\ma P}_L\mid \pL|_K\in S\}$ y sea $\AE T$ la cerradura
entera de $\AE S$ en $L$.

\begin{proposicion}\label{PClaseS4.9-1.2}
Se tiene $K_{H,S} L\subseteq L_{H,T}$.
\end{proposicion}

\begin{proof}
Se sigue del hecho de que $L=KL\subseteq K_{H,S}L$ y de que
$K_{H,S}/K$ es abeliana, no ramificada y los primos de $S$ se
descomponen totalmente.
$\fin$
\end{proof}

\begin{observacion}\label{OClaseS4.9-1.4}
Si $K_{H,S}\cap L=K$, entonces $h_S|h_T$. En efecto, se tiene
$\Gal(K_{H,S}L/L)\cong\Gal(K_{H,S}/K)$ por lo que $h_S=|
\Gal(K_{H,S}/K)|=|\Gal(K_{H,S}L/L)||h_T$.
\[
\xymatrix{
&L_{H,T}\ar@{-}[d]\ar@/^2pc/@{-}[dd]^{h_T}\\
K_{H,S}\ar@{-}[r]\ar@{-}[d]^{h_S}\ar@/_1pc/@{-}[d]_{
\Gal(K_{H,S}/K)}&K_{H,S}L\ar@{-}[d]_{h_S}\ar@/^1pc/@{-}[d]^{\substack{
\phantom{x}\frac{\Gal(L_{H,T}/L)}{\Gal(L_{H,T}/K_{H,S}L)}
\cong\Gal(K_{H,S}L/L)\\ \phantom{xxxxxxxxxxxxxx}\cong\Gal(K_{H,S}/K)}}\\
K_{H,S}\cap L=K\ar@{-}[r]& L
}
\]
\end{observacion}

\begin{proposicion}\label{PClaseS4.9-1.3}
Si alg\'un primo $\pK\in{\ma P}_K$ es totalmente ramificado en
$\AE T$ o si alg\'un primo $\pK\in S$ es totalmente inerte en
$L$, entonces la norma $\N:=\N_{L_{H,T}/K_{H,S}}
\colon Cl_T\lra Cl_S$ es 
suprayectiva. En particular se tiene que $h_S|h_T$.
\end{proposicion}

\begin{proof}
Se tiene que, en cualquiera de los dos casos mencionados,
$K_{H,S}\cap L=K$, por lo que el mapeo de restricci\'on
$\rest\colon\Gal(L_{H,T}/L)\lra\Gal(K_{H,S}/K)$
es suprayectivo y $h_S|h_T$.

Del Teorema \ref{T7.8.5},
se tiene el siguiente diagrama conmutativo
\[
\xymatrix{
Cl_T\ar@{->}[rrr]^{(\ ,L_{H,T}/L)}_{\cong}\ar@{->}_{\N}[d]
&&&\Gal(L_{H,T}/L)\ar@{->>}[d]^{\rest|_{K_{H,S}}}\\
Cl_S\ar@{->}[rrr]^{(\ , K_{H,S}/K)}_{\cong}&&&\Gal(K_{H,S}/K)
}
\]
Puesto que $\rest|_{K_{H,S}}$, $(\ ,L_{H,T}/L)$ y $(\ ,K_{H,S}/K)$
son suprayectivas, se sigue que $\N\colon Cl_T\lra Cl_S$ es
suprayectiva.
$\fin$
\end{proof}

Aqu\'i recordamos la relaci\'on que existe entre el n\'umero de
clases de divisores $h_K$ de un campo de funciones global $K$
y el n\'umero de clases de ideales, dada por la f\'ormual de
Schmidt (Corolario \ref{CRamDed1.2.7}):
\[
h_{K} d = \gamma h_S,
\]
donde $\gamma=\big|\frac{D_{K}(S)_0}{P_{K}(S)}\big|$
y $d=\mcd(\deg \pK\mid \pK\in S)$.

\begin{observacion}\label{OClaseS4.9-1.5}
Si la extensi\'on $L/K$ es una extensi\'on finita de Galois,
entonces $L_{H,T}/K$ tambi\'en es Galois. La demostraci\'on
es la misma que la de en la Observaci\'on \ref{OClaseS4.8}.
\end{observacion}

A continuaci\'on definimos el {\em campo de clase de Hilbert
extendido} de un campo de funciones global.

En el caso de campos num\'ericos, el campo de clase de 
Hilbert extendido corresponde al grupo de id\`eles $J_K^{1_+}$,
esto es, a elementos totalmente positivos. Debemos desarrollar
un concepto de ``totalmente positivo'' en un campo de funciones
global. Este concepto que usamos aqu\'i 
fue propuesto por Angl\`es y Jaulent en
\cite{AngJau2000}.

Sean $K=\F(T)$, $\p$ el primo infinito, es decir, el polo de $(T)_K$,
y $\Kii\cong\F\big(\big(\frac{1}{{T}}\big)\big)$ la completaci\'on 
de $K$ en $\p$. Sea $x\in\*{\Kii}$. Entonces $x$ se escribe de
manera \'unica como
\[
x=\Big(\frac{1}{T}\Big)^{n_x}\lambda_x\varepsilon_x,\quad n_x\in
{\ma Z},\quad \lambda_x\in \ff\quad\text{y}\quad\varepsilon_x\in 
U_{\Kii}^{(1)}.
\]
Ponemos $\pi_{\infty}=1/T$\label{piinfty}.

\begin{definicion}\label{DClaseS4.9-1.6}
Se define la {\em funci\'on signo\index{funci\'on signo} de $\*{\Kii}$}
como $\phi_{\infty}\colon \*{\Kii}\lra \ff$ la cual est\'a dada por
$\phi_{\infty}(x)=\lambda_x$ para $x\in\*{\Kii}$. El valor $\phi_{
\infty}(x)$ se llama el {\em signo\index{signo} de $x$}.
Tambi\'en ponemos $\sgn (x)=\phi_{\infty}(x)$.
\end{definicion}

Se tiene que $\phi_{\infty}$ es un epimorfismo y $\ker\phi_{\infty}=
\langle\pi_{\infty}\rangle\times U_{\*{\Kii}}^{(1)}$.

\begin{observacion}\label{OClaseS4.9-1.7}
Sea $M\in R_T=\F[T]$, $M\neq 0$, digamos que $M$ es de grado
$d$ y est\'a dado por $M(T)=a_dT^d+a_{d-1}T^{d-1}
\cdots+a_1T+a_0=a_d T^d\big(1+\frac{a_{d-1}}{a_d}\big(\frac{1}{T}
\big)+\cdots+\frac{a_1}{a_d}\big(\frac{1}{T}\big)^{d-1}+
\frac{a_0}{a_d}\big(\frac{1}{T}\big)^{d}\big)=a_d\big(\frac{1}{T}\big)^{-d}
\mu$ con $\mu\in U_{\*{\Kii}}^{(1)}$. Se sigue que el signo de
$M$ es el coeficiente l\'ider de $M$: $\sgn(M)=a_d$.
\end{observacion}

\begin{definicion}\label{DClaseS4.9-1.8}
Sea $L$ una extensi\'on finita y separable de $\Kii$. Se define
el {\em signo de $\*L$\index{signo de un campo local}}
por el morfismo $\phi_L=\phi_{\infty}\circ\N_{L/{\Kii}}\colon
\*L\lra \ff$.
\end{definicion}

\begin{proposicion}\label{PClaseS4.9-1.9}
Sea una torre de campos $\Kii\subseteq E\subseteq L$. Entonces
se tiene
\las
\item $\phi_L=\phi_E\circ\N_{L/E}$,

\item $\N_{L/E}(\ker \phi_L)\subseteq \ker \phi_E$.
\end{list}
\end{proposicion}

\begin{proof}
(1): Se tiene $\phi_L=\phi_{\infty}\circ \N_{
L/\Kii}=\phi_{\infty}\circ (\N_{E/\Kii}\circ\N_{L/E})=
(\phi_E\circ\N_{L/E})$. Por tanto $\phi_L=\phi_E\circ\N_{L/E}$.

(2): Si $x\in\ker \phi_L$, entonces $1=\phi_L(x)=\phi_E(
\N_{L/E}(x))$. Por lo tanto $\N_{L/E}(x)\in \ker\phi_E$.
$\fin$
\end{proof}

\begin{observacion}\label{O17.9.5}
La completaci\'on de ${\ma Q}$ en $\infty$, el lugar arquimediano de
${\ma Q}$, es ${\ma Q}_{\infty}={\ma R}$ y la funci\'on signo
$\phi_{\infty}\colon \*{\ma R}=\*{{\ma Q}_{\infty}}\lra\{\pm 1\}$ es
$\phi_{\infty}(x)=1$ si $x>0$, $\phi_{\infty}(x)=-1$ si $x<0$ y $\ker
\phi_{\infty}={\ma R}^+$.

Ahora si $L_v={\ma C}$, entonces para toda $z\in\*{L_v}$,
$\N_{L_v/{\ma Q}_{\infty}}z=\N_{{\ma C}/{\ma R}}z=z\overline{z}=
|z|^2>0$ y por tanto $\phi_{L_v}(z)=1$ para toda $z\in L_v$.

Tambi\'en notemos que la definici\'on de signo concuerda con
la Definici\'on \ref{CClaseD4.5.4'}.
\end{observacion}

\begin{teorema}\label{TClaseS4.9-1.10}
Sean $E/\Kii$ una extensi\'on finita y separable y $\mo E/\Kii$
la m\'axima subextensi\'on abeliana moderadamente ramificada
m\'axima tal que $\Kii\subseteq \mo E\subseteq E$. Sean $f$
el grado de inercia y $e$ el \'indice de ramificaci\'on de $\mo E/
\Kii$. Entonces ${\ma F}_E={\ma F}_{q^f}$ es el campo residual
de $\mo E$ y
\las
\item existe $\xi\in \*{{\ma F}_E}$ tal que $\mo E={\ma F}_E
\Big(\Big(\sqrt[e]{\frac{-\xi}{T}}\Big)\Big)$ y adem\'as
\[
\phi_E(\*E)=\langle (\*\F)^e, \N_{{\ma F}_E/\F}(\xi)\rangle,
\]

\item si adem\'as tenemos que $\frac{1}{T}\in \N_{\mo E/\Kii}
((\mo E)^*)$, entonces $\mo E=\F\Big(\Big(\sqrt[e]{\frac{-1}{T}}
\Big)\Big)$ y $\phi_E(\*E)=(\*{\F})^e$.
\end{list}
\end{teorema}

\begin{proof}
Sea $\Kii\subseteq \abe E\subseteq E$ la m\'axima extensi\'on abeliana
de $\Kii$ contenida en $E$. Entonces por el Teorema \ref{T17.6.192N}
obtenemos que $\N_{E/\Kii}(\*E)=\N_{\abe E/\Kii}((\abe E)^*)$. Puesto
que $\phi_E=\phi_{\infty}\circ \N_{E/\Kii}$ se sigue que $\phi_E(\* E)=
\phi_{\abe E}((\abe E)^*)$.

Sea $\nr E={\ma F}_E\Big(\Big(\frac 1T\Big)\Big)$ la m\'axima extensi\'on 
no ramificada de $\Kii$ contenida en $E$. Se tiene $v_{\nr E}\big(\frac 1T
\big)=e(\nr E|\Kii)v_{\Kii}\big(\frac 1T\big)=1\cdot 1=1$, $[{\ma F}_E:\F]=
f=[\nr E:\Kii]$ y $\mo E/\nr E$ es totalmente ramificada, digamos 
de grado $e$, con $\mcd(p,e)=1$, esto es $p\nmid e$.
Sea $\rho$ un elemento uniformizador de $\mo E$. Entonces
$v_{\mo E}\big(\frac 1T\big)=e(\mo E|\nr E)v_{\nr E}\big(\frac 1T\big)
=e\cdot 1=e$.

Se tiene $v_{\mo E}\big(\rho^e\big(\frac 1T\big)^{-1}\big)=0$, esto es,
$\rho^e\big(\frac 1T\big)^{-1}\in U_{\mo E}={\ma F}_E^* U_{\mo E}^{(1)}$.
Podemos escribir $\rho^e\big(\frac 1T\big)^{-1}=\xi u$ con $\xi\in {\ma F}_E^*$
y $u\in U_{\mo E}^{(1)}$. Por tanto $\rho^e=\xi u\big(\frac 1T\big)$. Puesto
que $p\nmid e$ se tiene que $(U_{\mo E}^{(1)})^e=U_{\mo E}^{(1)}$ y por 
tanto existe $v\in U_{\mo E}^{(1)}$ con $v^e=u$. De esta forma tenemos
$\rho^e=\xi v^e \big(\frac 1T\big)$ y $v_{\mo E}(v)=0$. Se sigue que
$\big(\frac {\rho}v\big)^e=\xi\big(\frac 1T\big)$ y $v_{\mo E}\big(\frac {\rho}
v\big)=1$. De esta forma, podemos suponer, tomando $\frac {\rho}v$
en lugar de $\rho$, que $\rho^e =\xi \big(\frac
1T\big)$.

Por lo tanto tenemos $\N_{\mo E/\Kii}(\rho)=\N_{\nr E/\Kii}\circ \N_{
\mo E/\nr E}(\rho)$. Se tiene 
\begin{align*}
\Irr(\rho,X,\nr E)&=X^e-\xi\Big(\frac 1T\Big)=
\prod_{\sigma\in\Gal(\mo E/\nr E)}(x-\sigma(\rho))\\
&= X^e-\alpha_{e-1}+
\cdots+(-1)^e\N_{\mo E/\nr E}(\rho).
\end{align*}
Se sigue que $-\xi \big(\frac 1T\big)=(-1)^e
\N_{\mo E/\nr E}(\rho)$ y por tanto $\N_{\mo E/\nr E}(\rho)=(-1)^{e-1}
\xi \big(\frac 1T\big)$ lo cual implica que
\begin{align*}
\N_{\mo E/\Kii}(\rho)&=\N_{\nr E/\Kii}\Big((-1)^{e-1}\xi\Big(\frac 1T\Big)\Big)\\
&=(-1)^{(e-1)f}\big(\N_{\nr E/\Kii}(\xi)\big)\Big(\frac 1T\Big)^f.
\end{align*}
Se sigue $\phi_{\mo E}(\rho)=\phi_{\infty}\circ \N_{\mo E/\Kii}(\rho)=
\phi_{\infty}\big((-1)^{(e-1)f}\big(\N_{{\ma F}_E/\F}(\xi)\big)\frac 1{T^f}\big)=
(-1)^{(e-1)f}\N_{{\ma F}_E/\F}(\xi)$.

Adem\'as, puesto que $\mo E/\nr E$ es totalmente ramificada de grado $e$
y la norma en campos finitos es suprayectiva, se tiene
\begin{align*}
\phi_{\mo E}({\ma F}_E^*)&=\phi_{\nr E}(\N_{\mo E\nr E}({\ma F}_E^*))=
\phi_{\nr E}(({\ma F}_E^*)^e)=\phi_{\infty}((\N_{\nr E/\Kii} ({\ma F}_E^*))^e)\\
&=\phi_{\infty}((\N_{\nr E/\Kii}({\ma F}_{\nr E}^*))^e)=\phi_{\infty}((\*\F)^e)=
(\*\F)^e.
\end{align*}

Ahora bien $(\mo E)^*=\langle \rho\rangle \times {\ma F}_E^*\times
U_{\mo E}^{(1)}$ y
\[
\phi_{\mo E}(U_{\mo E}^{(1)})=\phi_{\infty}(
\N_{\mo E/\Kii} (U_{\mo E}^{(1)}))=\{1\}
\]
 pues $\N_{\mo E/\Kii} (U_{\mo E}^{(1)})
\subseteq U_{\infty}^{(1)}$.

Se sigue que 
\begin{align*}
\phi_{\infty}((\mo E)^*)&=\langle \phi_{\infty}(\rho),\phi_{\infty}
({\ma F}_E^*)\rangle=\langle (-1)^{(e-1)f}\N_{{\ma F}_E/\F}(\xi),(\*\F)^e\rangle\\
&=\langle (-1)^{ef}(-1)^f\N_{{\ma F}_E/\F}(\xi),(\*\F)^e\rangle\\
&=\langle (-1)^{ef}\N_{{\ma F}_E/\F}(-\xi), (\*\F)^e\rangle.
\end{align*}

Se tiene $(-1)^{ef}=\pm 1$. Si $(-1)^{ef}=-1$, entonces en caso de que
$p=2$, $-1=1$ y en caso de que $p>2$, $ef$ es impar por lo que tanto
$e$ como $f$ son impares. Por tanto $(-1)^e=-1\in (\*\F)^*$.
En resumen, tenemos, en cualquier caso,
\[
\phi_{\infty}((\mo E)^*)=\langle \N_{{\ma F}_E/\F}(-\xi),(\*\F)^e\rangle.
\]

Ahora bien, puesto que $\abe E/\mo E$ es de grado una potencia de $p$ y
$\phi_{\infty}(\*E)\subseteq \*\F$ con $\mcd(p,q-1)=1$, se sigue que
\[
\phi_{\infty}(\*E)=\phi_{\infty}(\*{(\abe E)})=\phi_{\infty}(\*{(\mo E)})=
\langle (\*\F)^e,\N_{{\ma F}_E/\F}(-\xi)\rangle.
\]

Por otro lado $\mo E/\Kii$ es de grado $ef$ y $\nr E={\ma F}_E\Big(\Big(
\frac 1T\Big)\Big)={\ma F}_{q^f}\Big(\Big(\frac 1T\Big)\Big)$ y $\mo E=
\nr E(\rho)$ pues $e=\deg\Irr(\rho,X,\mo E)$ y $\rho^e=\xi\frac 1T$ por
lo que $\rho=\sqrt[e]{\xi\frac 1T}=\sqrt[e]{\frac {\xi}T}$. Por tanto
\[
\mo E=\nr E \Big(\sqrt[e]{\frac {\xi} T}\Big)= {\ma F}_E\Big(\Big(\frac
1T\Big)\Big)\Big(\sqrt[e]{\frac {\xi}T}\Big)={\ma F}_E\Big(\Big(\sqrt[e]{\frac
{\xi}T}\Big)\Big).
\]
Esto termina la demostraci\'on de (1).

Para (2), tenemos que si $\frac 1T\in \N_{\mo E/\Kii}(\*{(\mo E)})$, entonces
existe $\mu\in \mo E$ con $\N_{\mo E/\Kii}(\mu)=\frac 1T$ de donde
$-\xi\Big(\frac 1T\big)=-\xi \N_{\mo E/\Kii}(\mu)=(-1)^e\N_{\mo E/\nr E}(\rho)$
y por tanto $-\xi(-1)^e\in \N_{\mo E/\Kii}(\*{(\mo E)})$.

Adem\'as, como $\frac 1T$ es norma, se tiene que $f=1$. De hecho se tiene
que $\frac 1T$ es uniformizador de $\nr E$ y $1=v_{\Kii}\Big(\frac 1T\Big)
=v_{\Kii}(\N_{\mo E/\Kii}(\mu))=
f=f(\mo E|\Kii)v_{\mo E}(\mu)=1$, esto es, $f=f(\mo E|\Kii)=1$ y $v_{\mo E}
(\mu)=1$. Por otro lado $(-1)^e=\N_{\mo E/\Kii}(-1)=(-1)^{[\mo E:\Kii]}$. Se
sigue que $-\xi \in \N_{\mo E/\Kii}(\*{(\mo E)})$.

Sea $-\xi=\N_{\mo E/\Kii}(\delta)$. Como $\rho$ es un elemento 
uniformizador de $\mo E$, tenemos
\[
\*{(\mo E)}=\langle \rho\rangle\times {\ma F}_{\mo E}^*\times U_{\mo E}^{(1)}=
\langle \rho\rangle\times \*\F\times U_{\mo E}^{(1)}.
\]
Ahora bien $v_{\Kii}(\N_{\mo E/\Kii}(\rho^s))=es$ el cual es $0$ si y solamente si
$s=0$ lo cual equivale a que $\rho^s=1$. Adem\'as $\N_{\mo E/\Kii}(U_{\mo E}^{
(1)})\subseteq U_{\Kii}^{(1)}$. Por tanto $-\xi = \N_{\mo E/\Kii}(\delta)$ con
$\delta\in\*\F$ por lo que $-\xi=\delta^e\in (\*\F)^e$. Se sigue que
\begin{align*}
\mo E&=\Kii\Big(\sqrt[e]{\frac {\xi}T}\Big)=\Kii\Big(\sqrt[e]{\frac{-(-\xi)}T}\Big)=
\Kii\Big(\sqrt[e]{\frac {-\delta^e}T}\Big)\\
&=\Kii\Big(\delta \sqrt[e]{-\frac 1T}\Big)=
\Kii\Big(\sqrt[e]{-\frac 1T}\Big)
\end{align*}
lo cual demuestra (2). $\fin$
\end{proof}

Para el siguiente corolario ver las Proposiciones
\ref{P12*.2.2.G} y \ref{P12*.2.2.G'}.

\begin{corolario}\label{CClaseS4.9-1.11}
La extensi\'on abeliana de $\Kii$ asociada a $\ker \phi_{\infty}$
es el campo $\Kii^{\phi}:=\Kii\big[\sqrt[q-1]{-1/T}\big]=\F\big(\big(
\sqrt[q-1]{-1/T}\big)\big)$ la cual es la m\'axima extensi\'on abeliana
moderadamente ramificada de $\Kii$ tal que $1/T$ es una norma.
En particular, si $\K$ es una extensi\'on separable de $\Kii$, la
extensi\'on abeliana de $\K$ asociada a $\ker \phi_{\K}$ por la 
teor\'ia local de campos de clase, es la composici\'on 
$\K^{\phi}:=\K\Kii^{\phi}=\K\big[\sqrt[q-1]{-1/T}\big]=
\K[\sqrt[q-1]{-\pi_{\infty}}]$.
\end{corolario}

\begin{proof}
Sea $E/\Kii$ la extensi\'on asociada a 
$\ker\phi_{\infty}$ dada por la teor\'ia local de
campos de clase, la cual, por construcci\'on, verifica
la identidad n\'ormica: $\N_{E/\Kii}(\*E)=\ker\phi_{\infty}=
U_{\Kii}^{(1)}\Big(\frac{1}{T}\Big)^{\ma Z}$ (ver Proposici\'on
\ref{CClaseP3.2.1.3} y Teoremas \ref{CClaseT3.2.5.34} 
y \ref{CClaseT3.2.5.6.3}). De hecho, como $\frac{1}{T}$
es norma, $f=1$ y $\*\F$ corrsponde a la ramificaci\'on
moderada y $U_{\Kii}^{(1)}$ a la ramificaci\'on salvaje,
de donde $\*\F U_{\Kii}^{(1)}=U_{\Kii}$ corresponde a la
ramificaci\'on. Por tanto a primera afirmaci\'on se sigue 
inmediatamente del Teorema \ref{TClaseS4.9-1.10}
pues $\mo E=\F\Big(\Big(\sqrt[e]{-\frac{1}{T}}\Big)\Big)$
y la $e$ m\'axima es $e=q-1$, esto es,
la m\'axima extensi\'on corresponde
a $\Kii\Big(\sqrt[q-1]{-\frac{1}{T}}
\Big)=\F\Big(\Big(\sqrt[q-1]{-\frac{1}{T}}\Big)\Big)$.

Ahora (ver Teoremas \ref{CClaseT3.2.21'+1} y \ref{CClaseT4.6.9-1}), 
tenemos que $\phi_{\K}=\phi_{\infty}\circ \N_{\K/\Kii}$
lo cual implica $\ker \phi_{\K}=\N^{-1}_{\K/\Kii}(\ker\phi_{\infty})$.
Por tanto, si $H_{\K}$ es el grupo de normas de $\K/\Kii$,
esto es, $\N_{\K/\Kii}\*\K=H_K$. De esta forma tenemos
que $\ker \phi_{\K}=H_K\cap H_{\ker \phi_{\infty}}=
H_{\K}\cap H_{\Kii^{\phi}}=H_{\K\Kii^{\phi}}$.
Por tanto $\ker\phi_{\K}=H_{\K^{\phi}}=H_{\K\Kii^{\phi}}$
y en consecuencia, $\K^{\phi}=\K\Kii^{\phi}=\K\big[
\sqrt[q-1]{-1/T}\big]$ (ver Teorema \ref{CClaseT3.2.29}).
$\fin$
\end{proof}

Consideremos $L/K$ una extensi\'on finita y separable.
Sea $\Pi:=\{\pK\in{\ma P}_L\mid \pK|\p\}$. Sea $L_v$ la completaci\'on
de $L$ en $v\in \Pi$. Se tiene que $\phi_{L_v}(x)$ est\'a
bien definido para $x\in\*L$.

\begin{definicion}\label{DClaseS4.9-1.12}
El elemento $x\in\*L$ se llama {\em totalmente 
positivo\index{elemento totalmente 
positivo}\index{elemento!totalmente positivo}} si
$\phi_{L_v}(x)=1$ para toda $v\in\Pi$. 

Se define
$L^+=\{x\in L\mid x\text{\ es totalmente 
positivo}\}.$\index{totalmente positivos!conjunto de 
elementos $\sim$}\label{totalmentepositivos}
\end{definicion}

\begin{definicion}\label{DClaseS4.9-1.13}
El {\em grupo de signos\index{grupo de signos de un campo
de funciones global}} de un campo global $L$ se define por
$\Sg_L:=\*L/L^+$\label{gruposignoscampoglobal}.
\end{definicion}

Si $R$ es un subgrupo de $\*L$, la imagen de $R$ bajo
$\Sg_L$ se denota por $\sg_L(R)$. En particular, $\Sg_L
=\sg_L(\*L)$.

\begin{proposicion}\label{PClaseS4.9-1.14}
Se tiene $\Sg_L=\sg_L(\*L)\cong\prod_{v\in\Pi}\phi_{
L_v}(\*{L_v})$.
\end{proposicion}

\begin{proof}
Consideremos el mapeo $\theta\colon\*L\lra 
\prod_{v\in\Pi} \phi_{L_v}(\*{L_v})$ definido por
$\theta(x)=\big(\phi_{L_{v_1}}(x),\ldots,\phi_{L_{v_s}}(x)\big)$. Entonces
$\theta$ es un homomorfismo de grupos y $x\in\ker\theta\iff \theta(x)=
\big(\phi_{L_{v_1}}(x),\ldots,\phi_{L_{v_s}}(x)\big)=(1,\ldots,1)\iff
\phi_{L_{v}}(x)=1\text{\ para toda\ } v\in\Pi\iff x\in\*L \text{\ y\ }
\tilde{\theta}\colon \*L/L^+\lra \prod_{v\in\Pi} \phi_{L_v}(\*{L_v})$
es un monomorfismo.

Veamos la suprayectividad de $\tilde{\theta}$. Sea $(\alpha_1,\ldots,
\alpha_s)\in \prod_{v\in\Pi} \phi_{L_v}(\*{L_v})$. Sean $y_i\in\*{L_v}$
tal que $\phi_{L_v}(y_i)=\alpha_i$. Digamos que $y_i=\sum_{j=m}^{
\infty}\beta_j\pi^j$, $m\in {\ma Z}$, $\beta_m\neq 0$. Entonces
$y_i=\beta_m\pi^m\xi$ con $\xi\in U_{L_v}^{(1)}$ y 
$\phi_{L_{v_i}}=\phi_{\infty}(\N_{L_{v_i}/\Kii}(y_i))=\phi_{\infty}
(\N_{L_{v_i}/\Kii}(\beta_m))$.

Sea $x\in\*L$ con $v_{L_{v_i}}(x-y_i)>m$. Ahora bien, $v_{
L_{v_i}}(y_i)=m$ por lo que $v_{L_{v_i}}(x)=m$ y $
\phi_{L_{v_i}}(y_i)=\phi_{L_{v_i}}(x)$. De esta forma, usando
el Teorema de Aproximaci\'on de Artin, tomamos $x\in\*L$ tal
que $v_{L_{v_i}}(x-y_i)>m_i$ donde $v_{L_{v_i}}(y)=m_i$ y
$\big(\phi_{L_{v_1}}(y_1),\ldots,\phi_{L_{v_s}}(y_s)\big)=
\big(\phi_{L_{v_1}}(x),\ldots,\phi_{L_{v_s}}(x)\big)=\theta(x)$.
Por tanto $\theta$ es suprayectiva y $\tilde{\theta}$ es un
isomorfismo.
$\fin$
\end{proof}

\begin{observacion}\label{O17.9.10}
Sea $L$ un campo num\'erico. Entonces $L_v\in\{{\ma R},{\ma C}\}$
para $v|\infty$. Entonces $\N_{L_v/{\ma R}} \*{L_v}\in
\{{\ma R}^+,\*{\ma R}\}$ y $\phi_{L_v}\*{L_v}\subseteq \{\pm 1\}$.
Entonces $\*L/L^+\cong C_2^t$ donde $t$ es e n\'umero de 
lugares reales de $L$ tales que $\N_{L_v/{\ma R}}\*{\ma R}=
\*{\ma R}$ pues en este caso $\phi_{L_v}(\*{\ma R})=\{\pm 1\}$.
\end{observacion}

\begin{definicion}\label{D17.8.21'}
Sea $\AE L=\{x\in L\mid v_{\pK}(x)\geq 0\text{\ para todo\ }
\pK\notin \Pi\}$. En otras palabras, $\AE L=\AE S$ donde
$S=\Pi$. Sean $P_L^+=\{x \AE L\mid x\in L^+\}$
el grupo de ideales principales de $\AE L$ que son
generados por un elemento totalmente positivo
\index{ideales principales totalmente positivos}\index{totalmente
positivos!ideales principales $\sim$} y definimos
el {\em grupo de clases de ideales extendidos}\index{grupo de
ideales extendidos}\label{idealesextendidos}
\[
Cl_L^{\rm{ext}}=Cl_L^+=I_L/P_L^+=C_S^{\rm{ext}}=C_S^+,
\]
donde $I_L$ es el grupo de ideales fraccionarios de 
$\AE L$. Recordemos que $Cl_L=I_L/P_L$ donde $P_L=\{(x)=
x\AE L\mid x\in\*L\}$.
\end{definicion}

Se tiene $|Cl_L|<\infty$ (Corolario \ref{CRamDed1.2.6}).

\begin{definicion}\label{D17.9.13}
Definimos los siguientes subgrupos del grupo de id\`eles
$J_L$ de $L$
\begin{gather*}
U_L=\prod_{v\in \Pi} \*{L_v}\times \prod_{v\notin \Pi} U_{L_v},\\
U_L^+=\prod_{v\in \Pi}\ker\phi_{L_v}\times
\prod_{v\notin\Pi}U_{L_v}.
\end{gather*}
\end{definicion}

Se tiene que $U_L$ y $U_L^+$ son subgrupos abiertos
del grupo de id\`eles de $L$, $J_L$ (ver la Observaci\'on
\ref{O17.6.35N} y usar el Proposici\'on \ref{CCUnidades}).
Los grupos $U_L$ y $U_L^+$ son los correspondientes a $J_L^{1}$
y $J_L^{1_+}$ respectivamente en el caso de campos num\'ericos.

\begin{proposicion}\label{PClaseS4.9-1.15}
Los grupos de clases de ideales en el sentido ordinario y en el
sentido extendido (o restringido) est\'an dados por los isomorfismos
can\'onicos
\las
\item $Cl_L\cong J_L/U_L\*L$,
\item $Cl_L^+=Cl_L^{\rm{ext}}\cong J_L/U_L^+\*L$.
\end{list}
\end{proposicion}

\begin{proof}
Es consecuencia del isomorfimso $I_L\cong J_L/U_L$.
$\fin$
\end{proof}

\begin{corolario}\label{CClaseS4.9-1.16}
Las clases en el sentido ordinario y en el sentido extendido est\'an
ligadas por la sucesi\'on exacta
\[
1\lra \sg_L(E_L)\lra \Sg_L\lra Cl_L^{{\rm ext}}\lra Cl_L\lra 1,
\]
donde $E_L$ son las unidades de $\AE L$.\end{corolario}

\begin{proof}
Sea $\theta\colon Cl_L^{\rm{ext}}\lra Cl_L$ dada por 
$\vec\alpha\bmod U_L^+ \*L\stackrel{\theta}{\lra}
\vec\alpha\bmod U_L\*L$. Entonces $\theta$ es un epimorfismo
con n\'ucleo $\ker \theta =U_L\*L/U_L^+\*L$.
 Por tanto se tiene
la siguiente sucesi\'on exacta
\[
1\lra U_L\*L/U_L^+\*L\lra Cl_L^{\rm{ext}}\lra Cl_L\lra 1.
\]

Sea $\psi\colon U_L/U_L^+\lra U_L\*L/U_L^+\*L$ dada por
$\psi(\vec\alpha\bmod U_L^+)=\vec \alpha \bmod U_L^+\*L$.
Entonces es un epimorfismo bien definido. Ahora, se tiene
\begin{gather*}
\frac{U_L}{U_L^+}\cong \prod_{v\in\Pi}\frac{\*{L_v}}{\ker \phi_{L_v}}
\cong \prod_{v\in\Pi}\phi_{L_v}(\*{L_v})\cong \Sg_L.
\intertext{Por tanto se tiene la sucesi\'on exacta}
1\lra \ker\psi=\frac{U_L^+\*L\cap U_L}{U_L^+}
\lra \frac{U_L}{U_L^+}
\cong \Sg_L\lra \frac{U_L\*L}{U_L^+\*L}\lra 1.
\intertext{Se sigue que la sucesi\'on}
1\lra \frac{U_L^+\*L\cap U_L}{U_L^+}\lra \Sg_L\lra Cl_L^{\rm{
ext}}\lra Cl_S\lra 1
\end{gather*}
es exacta.

Ahora sea $\mu\colon\*L\cap U_L\hooklongrightarrow U_L^+\*L
\cap U_L\lra \frac{U_L^+\*L\cap U_L}{U_L^+}$ dada por
$\mu(\xi)=\xi\bmod U_L^+$. Verificamos que $\mu$ es epimorfismo.
Sea $\xi\in U_L^+\*L\cap U_L$, $\xi=ab\in U_L$ y $a\in
U_L^+$, $b\in \*L$. Esto es, $\xi\equiv b\bmod U_L^+$ y
$\xi\in U_L$, $b\in \*L$. Por tanto $\mu(b)=\bar{\xi}$ y
$b\in \*L$. Entonces se tiene $b=\xi a^{-1}\in U_LU_L^+=U_L$.
Por tanto $b\in U_L\cap \*L$. De esta forma $\mu$ es
suprayectiva.

Por otro lado, $\ker \mu =(\*L\cap U_L)\cap U_L^+=\*L\cap
U_L^+$ por lo que $\frac{U_L^+\*L\cap U_L}{U_L^+}\cong
\frac{\*L\cap U_L}{\*L\cap U_L^+}$.

Puesto que $U_L=\prod_{v\in\Pi}\*{L_v}\times \prod_{v\notin
\Pi}U_{L_v}$, $\*L\cap U_L=E_L$ y $\*L\cap U_L^+=
\*L\cap U_L^+\cap U_L=E_L\cap U_L^+=:E_L^+$. Por tanto
\[
\frac{U_L^+\*L\cap U_L}{U_L^+}\cong \frac{E_L}{E_L^+}\cong
\frac{E_L\*L}{\*L}\cong\sg (E_L),
\]
como consecuencia de la Proposici\'on \ref{PClaseS4.9-1.14}.
$\fin$
\end{proof}

Sea $L_H$ la extensi\'on abeliana de $L$ correspondiente
al subgrupo de id\`eles $\*L U_L$. Entonces $L_H$ es
la m\'axima extensi\'on abeliana de $L$ no ramificada y
donde los elementos de $\Pi$ se descomponen totalmente.
En otras palabras, $L_H$ es el correspondiente
al tercer an\'alogo al campo de clase de Hilbert
con $L_H=L_{H,S}$ y donde
$S=\Pi$. Adem\'as
\[
\Gal(L_H/L)\cong Cl(\AE L)=Cl_L\cong J_L/U_L\*L.
\]

Estamos en condiciones de definir el an\'alogo al campo
de clase de Hilbert extendido en el caso de campos de funciones
globales.

\begin{definicion}\label{DClaseS4.9-1.17}
Sea $L_{H^+}=L_H^{\rm{ext}}$ \label{Hilbertextendido}
la extensi\'on abeliana de $L$ correspondiente
al subgrupo de id\`eles $\*L U_L^+$ donde 
$U_L^+=\prod_{v\in \Pi}\ker\phi_{L_v}\times
\prod_{v\notin\Pi}U_{L_v}$.

El campo $L_{H^+}$ recibe el nombre de {\em el campo de clase
de Hilbert extendido\index{campo de clase de Hilbert
extendido}\index{Hilbert!campo de clase de $\sim$ extendido}}
correspondiente a $\AE L$.
\end{definicion}

Se tiene que $L_{H^+}/L$ es una extensi\'on no ramificada
en ning\'un primo finito, $L_H\subseteq L_{H^+}$ y 
\[
\Gal(L_{H^+}/L)\cong J_L/\*L U_L^+\cong I_L/P_L^+
\cong Cl_L^+\cong Cl_L^{\rm{ext}}.
\]

Por otro lado, se tiene
\[
\frac{U_L}{U_L^+}\cong \prod_{v\in\Pi}\frac{\*{L_v}}{\ker\phi_{L_v}}
\cong \prod_{v\in\Pi}\phi(\*{L_v})\cong \frac{\*L}{L^+}.
\]

\begin{observacion}\label{O17.8.24'}
Como consecuencia
del Corolario \ref{CClaseS4.9-1.16} se tiene que $[L_{H^+}:
L_H]||\Sg_L|=\prod_{v\in S}|\phi_{L_v}(\*{L_v})|| (q-1)^{|S|}$.

Por tanto, $L_{H^+}/L$ es no ramificada en los primos
finitos y moderadamente ramificada en los primos infinitos.

Adem\'as 
\[
\Gal(L_{H^+}/L_H)\cong\frac{\Gal(L_{H^+}/L)}{\Gal(L_H/L)}
\cong \frac{U_L\*L}{U_L^+\*L}\cong \frac{\Sg_L}{\Big(
\frac{U_L\cap U_L^+\*L}{U_L^+}\Big)}.
\]

Todo esto es aplicable a campos num\'ericos sustituyendo
$q-1$ por $2$.
\end{observacion}

\begin{ejemplos}\label{EjClaseS4.9-1.18}{\ }

\las
\item Para $L=K=\F(T)$, $\sg (L)=\sg (\*\F)=\*\F$. Esto es,
$\Sg_L=\*\F=E_LL^+/L^+$, 
$\sg (E_L)=\*\F$. Por tanto, $L_H=L_{H^+}=L$.

\item Sea $M\in R_T=\F[T]$ un polinomio m\'onico de grado 
mayor o igual a $1$. Sea $L=\cicl M{}$. 

Tenemos que para toda
$v\in\Pi$, $L_v=\Kii\big(\sqrt[q-1]{-1/T}\big)$, $\phi_{L_v}=
\phi_{\infty}\circ\N_{L_v/\Kii}$. Por el Teorema 
\ref{CClaseT3.2.5.32}, $\N_{L_v
/\Kii}\*{L_v}=(\pi_{\infty})\times U_{\infty}^{(1)}$ por lo que $\ker 
\phi_{L_v}=\*{L_v}$ y $\Sg_K=\{1\}$. Se sigue que $L_H
=L_{H^+}$.
\end{list}
\end{ejemplos}

\begin{observacion}\label{OClaseS4.9-1.19}
Sea $L/E$ una extensi\'on finita y separable de campos globales.
Entonces $E_{H^+}\subseteq L_{H^+}$. Si $L/E$ es una 
extensi\'on de Galois, entonces $L_{H^+}/E$ tambi\'en es de Galois.

De hecho se tiene
\[
\xymatrix{
&E_{H^+}\ar@{-}[r]\ar@{-}[d]& LE_{H^+}\ar@{-}[d]\\
&E_{H^+}\cap L\ar@{-}[r]\ar@{-}[dl]&L\\ E
}
\]

Ahora bien $E_{H^+}/E_{H^+}\cap L$ es una extensi\'on 
abeliana no ramificada en los primos finitos lo cual implica
que $LE_{H^+}/L$ tambi\'en lo es.

Consideremos la norma $\N_{L/E}\colon C_L=J_L/\*L\lra C_E=
J_E/\*E$ de grupos de id\`eles y puesto que $N_{L_w/
E_v}(\ker\phi_{L_w})\subseteq \ker \phi_{E_v}$,
$\N_{L/E}(U_L^+)\subseteq U_E^+$ y $\N_{L/E}(\*L)
\subseteq \*E$, se sigue que $\N_{L/E}(U_L^+\*L)
\subseteq U_E^+\*E$. Ahora bien $U_E^+\*E$ corresponde
a $E_{H^+}$. Por el Teorema \ref{CClaseT4.6.9-1} obtenemos que
 $E_{H^+}\subseteq E_{H^+} L \subseteq L_{H^+}$.

Ahora si $L/E$ es una extensi\'on de Galois, consideremos
$\sigma\colon L_{H^+}\lra \bar{L}_{H^+}$ un monomorfismo
tal que $\sigma|_E=\Id_E$. Puesto que $L/E$ es Galois,
$\sigma|_L\colon L\lra\bar{L}$ con $\sigma|_E=\Id_E$,
obtenemos que $\sigma(L)=L$. Por otro lado como
$\sigma(\ker\phi_{L_w})=\ker \phi_{L_{\sigma(w)}}$
se tiene que $\sigma(U_L^+)=U_L^+$ y $\sigma(\*L)=
\*L$ lo cual implica que $\sigma(U_L^+ \*L)=U_L^+\*L$.
Por tanto $\sigma(L_{H^+})=L_{H^+}$, esto es,
$\sigma\in\Aut_E(L_{H^+})$ por lo que $L_{H^+}/E$
es Galois.
\end{observacion}

\begin{observacion}\label{OClaseS4.9-1.19'}
De hecho, para
$\star\in\{{\mathfrak{ge}}, {\mathfrak{gex}}, H, H^+$\} y para
una extensi\'on finita y separable $L/E$ de campos globales,
se tiene que $E_{\star}\subseteq L_{\star}$.
\end{observacion}

Para finalizar esta subsecci\'on, veamos cuales son los campos 
de constantes de $L_H$ para un campo
de funciones global $L$ con campo de constantes $\F$.
Hacemos notar que ya hemos calculado el campo de constantes
de $L_H$ en la Proposici\'on \ref{P3.2.B}. Ver
tambi\'en el Teorema \ref{T12*.2.2.A} y el Corolario \ref{CClaseC4.9.11}
m\'as adelante. Aqu\'i daremos una demostraci\'on diferente.

Sea pues $L/K$ una extensi\'on finita y separable donde $K=\F(T)$. 
Recordemos que si tenemos una extensi\'on abeliana
finita $E/F$ tal que el campo de constantes de $F$
es $\F$ y tal que $E$ corresponde al grupo $\Lambda$ de id\`eles de
$J_F$, esto es, $\Gal(E/F)\cong J_F/\Lambda\*F$, entonces 
si $d=\min\{n\in{\ma N}\mid \text{existe $\vec\alpha\in \Lambda$
con $\deg\vec\alpha=n$}\}$, se tiene que el campo de constantes
de $E$ es ${\ma F}_{q^d}$ (Teorema \ref{CClaseT4.8.1-1.1}).

En el Teorema \ref{T12*.2.2.A} se da el campo de constantes
del campo de clases de Hilbert. A continuaci\'on volvemos
a obtener este resultado por otros m\'etodos.

Sea $\con_{K/L}\p=\pK_1^{e_1}\cdots \pK_r^{e_r}$ y sean
$t_i:=\deg \pK_i$, $1\leq i\leq r$, $t_0=\mcd(t_1,\ldots,t_r)$.

\begin{teorema}\label{T17.9.26}
Con las notaciones anteriores, el campo de constantes
de $L_H$ es ${\ma F}_{q^{t_0}}$.
\end{teorema}

\begin{proof}
Sean $a_1,\ldots,a_r\in{\ma Z}$ tales que $t_0=\sum_{i=1}^r
a_it_i$. Sea $\vec\alpha\in J_L$ el id\`ele dado por componentes
como $\alpha_{\pK_i}=\pi_{\pK_i}^{a_i}$, $1\leq i\leq r$ y $\alpha_{\pK}
=1$ para $\pK\nmid \p$. Entonces $\vec\alpha\in U_L$ y
\[
\deg\vec\alpha=\sum_{i=1}^r
a_i\deg\alpha_{\pK_i}=\sum_{i=1}^r a_i\deg \pK_i v_{\pK_i}(\pi_{\pK_i})
=\sum_{i=1}^r a_it_i=t_0.
\]

Ahora sea $\vec\beta\in U_L$. Entonces las componentes de $\vec\beta$
son de la forma $\beta_{\pK_i}=\xi_i\pi_{\pK_i}^{b_i}$, $1\leq i\leq r$ con
$\xi_i\in U_{\pK_i}$ y $b_i\in{\ma Z}$ y $\beta_{\pK}\in U_{\pK}$ para
$\pK\nmid\p$. Entonces $\deg\vec\beta=\sum_{i=1}^r b_it_i$ y puesto
que $t_0|t_i$, $1\leq i\leq r$, $t_0|\deg \vec\beta$. 
Se sigue que $t_0=\min\{n\in{\ma N}\mid \text{existe $\vec\alpha\in U_L$
tal que $\deg\vec\alpha=n$}\}$. Por el Teorema
\ref{CClaseT4.8.1-1.1} obtenemos que el campo de constantes de $L_H$
es ${\ma F}_{q^{t_0}}$. $\fin$
\end{proof}

\begin{teorema}\label{T17.8.33'}
Sea $K/k$ una extensi\'on finita y separable con $k=\F(T)$. Sea
${\mc O}_K$ la cerradura entera de $\F[T]$ en $K$. Entonces un
ideal fraccionario $\pK$ de ${\mc O}_K$ se descompone totalmente en
\las
\item $K_H/K$ si y solamente si $\pK$ es principal;

\item $K_{H^+}/K$ si y solamente si $\pK$ es principal generado
por un elemento totalmente positivo de $K$, es decir 
un elemento de $K^+=\{x\in K\mid \phi_{K_{\pK}}(x)=1 \text{\ 
para toda\ }\pK|\p\}$.
\end{list}
\end{teorema}

\begin{proof}
Del Corolario \ref{C17.6.189N} se tiene que $\pK$ se descompone
totalmente en $K_H/K$
(resp. en $K_{H^+}/K$) si y solamente si $\theta\lceil
\*{K_{\pK}}\rceil_{\pK}:=\{(\ldots, 1, x, 1,\ldots)\mid x\in\*{K_{\pK}}\}
\subseteq \*K U_K$ (resp. $\subseteq \*K U_K^+$) si y solamente si
para toda $x\in\*{K_{\pK}}$ existen
$\beta_x\in\*K$ y $\vec\alpha_x\in U_K$ (resp.
$\vec\alpha_x\in U_K^+$) tal que
$\theta\lceil x\rceil_{\pK}=\beta_x\vec\alpha_x$.
Por tanto $\vec\alpha_x=(\ldots, \beta_x^{-1},\ldots, \beta_x^{-1},
\underbracket[0pt]{\beta_x^{-1}x}_{\substack{\uparrow\\ \pK}},
\beta_x^{-1},\ldots,\beta_x^{-1},\ldots)$. Esto es equivalente a $v_{\mathfrak q}
(\beta_x)=0$ para toda ${\mathfrak q}\neq \pK$ y ${\mathfrak q}\nmid\p$
y $\big(\beta_x^{-1}\big)_{{\mathfrak q}|\p}\in \prod_{{\mathfrak
q}|\p}\*{K_{\mathfrak q}}$ (resp. $\big(\beta_x^{-1}\big)_{{\mathfrak 
q}|\p}\in \prod_{{\mathfrak q}|\p}\ker\phi_{K_{\mathfrak q}}$).

Sea $x\in\*{K_{\pK}}$. El ideal principal en ${\mc O}_K$ generado por
$\beta_x$ satisface $\langle \beta_x
\rangle=\pK^{n_x}=\beta_x{\mc O}_K$. En particular para $x\in\*{K_{\pK}}$ con
$v_{\pK}(\beta_x)=1$ tenemos $\langle \beta_x\rangle=\pK$, $\beta_x\in
\*K$. Por tanto $\pK$ es un ideal principal $\pK=\langle \beta_x\rangle$ (resp.
$\pK$ es un ideal principal $\pK=\langle \beta_x\rangle$
y $\beta_x\in\ker\phi_{K_{\mathfrak q}}$ para toda ${\mathfrak q}|\p$).
$\fin$
\end{proof}

\subsection{Teor{\'\i}a global de campos de clase v{\'\i}a ideales o
divisores}\label{CClaseS4.7}

En la formulaci\'on v{\'\i}a id\`eles de la ley de reciprocidad, las
extensiones abelianas finitas corresponden un{\'\i}vocamente con
los grupos de normas $\N_{L/K}C_L$. En la formulaci\'on v{\'\i}a
ideales (campos num\'ericos), hay una correspondencia (no biyectiva) 
similar, pero m\'as complicada. Aqu{\'\i} tambi\'en las extensiones 
abelianas finitas corresponden a ciertos grupos de normas
(m\'as de uno) en el grupo 
de ideales fraccionarios $D_K$ de $K$. El s{\'\i}mbolo de la norma
residual $(\underline{\ \ }, L/K)\colon C_K\to \Gal(L/K)$
es reemplazado por el s{\'\i}mbolo de Artin, pero este \'ultimo no
est\'a definido en los primos ramificados.

Para un m\'odulus $\me$, $\DKm$\label{CClaseDKm} denotar\'a al grupo de ideales
fraccionarios (divisores) primos relativos a la parte finita de 
$\mK$, $\PKm\label{CClasePKm}$ es el grupo de ideales (divisores) principales
$(a)\in P_K$ con $a\equiv 1\bmod\mK$. Esto \'ultimo significa
que si $\pK|\mK$ es finito, entonces $a\equiv 1\bmod \mK$
significa $a\equiv 1\bmod \pK^{n_\pK}$ en el sentido
usual. Si $\pK$ es real y $n_\pK=1$, $a\equiv 1\bmod\mK$ significa
$a>0$ en el lugar $\pK$. Si $n_\pK=0$ o $\pK$ es complejo
$a\equiv 1\bmod \mK$ no establece ninguna restricci\'on para $a$.

\begin{definicion}\label{CClaseD4.7.1} El grupo $\PKm$ se llama el
{\em rayo m\'odulo $\mK$\index{rayo m\'odulo un m\'odulus}}
y todo grupo $H_K^\mK$ tal que $\PKm\subseteq H_K^\mK
\subseteq \DKm$ se llama el {\em grupo ideal (divisor) m\'odulo
$\mK$\index{grupo ideal}\index{grupo divisor}} y $\DKm/\PKm$
se llama el {\em grupo de rayos m\'odulo $\mK$\index{grupo
de rayos}}.
\end{definicion}

Si $\mK=1$, $\DKm=D_K$ y $\PKm=P_K$, es decir,
$D_K^1/P_K^1=D_K/P_K=I_K$ es el grupo de rayos m\'odulus $1$.

Por ejemplo, si denotamos $\infty$ un lugar real de $K$ (caso
num\'erico) y $\mK=\infty$, se tiene $\DKm=D_K^\infty=D_K=D_K^1$
y $P_K^\infty=\{(a)\in P_K\mid a>0 \text{\ en\ }\infty\}\subseteq P_K
\subseteq D_K$ y $\frac{D_K/P_K^\infty}{P_K/P_K^\infty}\cong
D_K/P_K$, es decir,
\[
1\to \frac{P_K}{P_K^\infty}\to \frac{D_K}{P_K^\infty}\to
\frac{D_K}{P_K}\to 1
\]
es exacta. Se puede pensar que $D_K/P_K$ corresponde al campo 
de clase de Hilbert y que $D_K/P_K^\infty$ corresponde al campo
de clase de Hilbert extendido, de hecho, con m\'as precisi\'on, el
campo de clase de Hilbert extendido, corresponde al m\'odulus 
${\eu m}:=\prod_{\pK\text{\ real}}\pK$ y $P_K^{\eu m}=\{(a)\in P_K
\mid a>0\text{\ para todo $\pK$ real}\}$.

\begin{definicion}\label{CClaseD4.7.2} Sea $K$ un campo num\'erico.
Si $L/K$ es una extensi\'on abeliana finita y $\mK$ es un m\'odulus,
$\mK$ se llama {\em m\'odulus de definici\'on\index{m\'odulus
de definici\'on}} o {\em m\'odulus admisible\index{m\'odulus
admisible}} para $L/K$ si $L$ est\'a contenido en el campo de clases
de rayos m\'odulus $\mK$, es decir si $\CKm\subseteq \N_{L/K}C_L$
o, equivalentemente, $L\subseteq K^\mK$.
\end{definicion}

Si $\mK|\mK^\prime$ y si $\mK$ es un m\'odulus de definici\'on
para $L$, tambi\'en lo es $\mK^\prime$ pues $C_K^{{\eu m}'}
\subseteq \CKm$.

El {\em conductor\index{conductor}} $\f{}=\f{L/K}=\f{}(L/K)$ es el 
m\'aximo com\'un divisor de los moduli de definici\'on de $L/K$
(Definici\'on \ref{CClaseD4.7.2}). 

\begin{definicion}\label{CClaseD4.7.2'}
Sea $K$ un campo num\'erico. Si $L/K$
es una extensi\'on abeliana finita y $\mK$ es un m\'odulus de
definici\'on de $L/K$,
entonces $H_K^\mK\label{CClaseHKm}:=(\N_{L/K}\DKm)\PKm$ se llama el
{\em grupo ideal definido m\'odulus $\mK$\index{grupo ideal
definido un m\'odulus}} asociado a $L/K$.
\end{definicion}

Sea $\artinp{L/K}{\pK}$ el s{\'\i}mbolo de Artin, 
esto es $\artinp{L/K}{\pK}=
\Fro\pK\in \Gal(L/K)$ el automorfismo de Frobenius
para $\pK\in {\ma P}_K$ finito y no ramificado en $L$.

\begin{teorema}\label{CClaseT4.7.3} Sea $K$ un campo global. Sea
$\Lambda\colon J_K\to D_K$ dado por 
\begin{gather*}
\Lambda(\vec \alpha)={\eu a}_{\vec\alpha}=\prod_{\substack{
\pKK\\ \pK\nmid\infty}}\pK^{v_\pK(\alpha_\pK)}.\\
\intertext{Entonces si $\me$ es cualquier m\'odulus, $\Lambda$ induce un
isomorfismo}
\Lambda_\mK\colon C_K/\CKm\longrightarrow \DKm/\PKm.
\end{gather*}
\end{teorema}

\begin{proof} Se tiene $C_K/\CKm\cong J_K/\JKm\*K$,
donde $\JKm=\{\vec\alpha\in J_K\mid \vec\alpha\equiv 
1\bmod \mK\}$. Sean
$S=\{\pKK\mid \pK|\mK\}$ y $\Jm=J_K^S=\{\vec\alpha\in J_K
\mid \alpha_\pK=1\text{\ para todo\ } \pK|\mK\}$.
Notemos que si $S=\{\pK\mid\pK\big|\mK\}$, entonces
$\Jm=J_K^S$.

Primero probemos que $J_K=\Jm\JKm K^{\ast}$. Sea $\vec\alpha\in J_K$.
Por el Teorema de Aproximaci\'on de Artin, existe $a\in K^{\ast}$ tal que
$\alpha_\pK a\equiv 1\bmod \pK^{n_\pK}$ para $\pK|\mK$.
Escribimos $\alpha_\pK a=\beta_\pK \gamma_\pK$ con
\begin{gather*}
\beta_\pK=1\quad {\text{para}}\quad \pK|\mK\quad\text{y}\quad
\beta_\pK=\alpha_\pK a\quad \text{para}\quad \pK\nmid \mK,\\
\gamma_\pK=\alpha_\pK a \quad {\text{para}}\quad \pK|\mK\quad\text{y}\quad
\gamma_\pK=1\quad \text{para}\quad \pK\nmid \mK.
\end{gather*}

Entonces $\vec\beta=(\beta_\pK)_\pK\in \Jm$ y 
$\vec\gamma=(\gamma_{\pK})_{\pK}\in\JKm$.
Adem\'as, se tiene, $\vec\alpha=\vec\beta\vec\gamma a^{-1}\in
\Jm\JKm K^{\ast}$.

Ahora bien $\frac{C_K}{\CKm}\cong \frac{J_K}{\JKm K^{\ast}}
=\frac{\Jm\JKm K^{\ast}}{\JKm K^{\ast}}
\cong \frac{\Jm}{\JKm K^{\ast}\cap \Jm}$. De esta forma se tiene que
$\Lambda\colon J_K\to D_K$ induce un epimorfismo
\[
\Jm\xrightarrow{\mu} \DKm/\PKm, \qquad \vec\alpha
\stackrel{\mu}{\longmapsto}\mu(\vec \alpha)=
\Lambda(\vec\alpha)\bmod \PKm={\eu a}_{\vec\alpha}\bmod \PKm
\]
y $\ker\mu=\JKm K^{\ast}\cap \Jm$ puesto que de la expresi\'on
$\Lambda(\vec \alpha)=\prod_{\pK\nmid\infty} \pK^{v_\pK(\alpha_\pK)}=
(a)=\prod_{\pK\nmid \infty}\pK^{v_\pK(a)}\in \PKm$ se sigue que
$\alpha_\pK a^{-1}\equiv 1\bmod \pK^{n_\pK}$ para toda $\pK$ y
por lo tanto tenemos el isomorfismo buscado. $\fin$
\end{proof}

Sea ahora $K$ un campo num\'erico. Sea
$\mK$ un m\'odulus de $K$ tal que $L$ est\'a contenido en el
campo de clases m\'odulus $\mK$, $L\subseteq K^{\mK}$
esto es, $\CKm\subseteq {\mc N}_L= \N_{
L/K}C_L$. 

Se tiene que si $\pK\nmid \mK$ entonces $\pK$ es
no ramificado y por tanto se obtiene un homomorfismo
\[
\artinp{L/K}{\ }\colon \DKm\longrightarrow \Gal(L/K)
\]
dado por
$\artinp{L/K}{\eu a}=\prod\limits_{\substack{\pK\nmid{\eu m}\\
\pK\nmid \infty}} \artinp{L/K}{\pK}^{a_\pK}$ donde
${\eu a}=\prod_\pK \pK^{a_\pK}$ y $\artinp{L/K}{\pK}$ es el
s{\'\i}mbolo de Artin en $\pK$.

Se tiene que si $\pK$ es un ideal primo y $\pi_\pK$ es un elemento
primo de $K_\pK$, entonces $\artinp{L/K}{\pK}=\big(\enc{\pi_\pK}
{\pK},L/K\big)=\psi_{L/K}\big(\enc{\pi_\pK}{\pK}\big)$ puesto que
ambos son el automorfismo de Frobenius correspondiente a $\pK$.

\begin{teorema}[Reciprocidad v\'ia ideales]\label{CClaseT4.7.4}
Sea $L/K$ una extensi\'on abeliana
finita de campos num\'ericos y se $\mK$ un m\'odulus de
de definici\'on de $L/K$, esto es, $L\subseteq K^\mK$. Entonces
el s{\'\i}mbolo de Artin induce un epimorfismo
\[
\artinp{L/K}{\ }\colon \DKm/\PKm\longrightarrow \Gal(L/K)
\]
el cual tiene n\'ucleo $H_K^\mK/\PKm$ donde $H_K^\mK=(\N_{L/K}
D_L^\mK) \PKm$ donde $D_L^\mK$
son los ideales fraccionarios de $D_L$
primos relativos a $\mK$. Por tanto $D_K^{\mK}/H_K^{\mK}
\cong \Gal(L/K)$.

M\'as a\'un, se tiene un diagrama conmutativo cuyas filas
son exactas
\[
\begin{CD}
1@>>>\N_{L/K}C_L@>>>  C_K@>{(\ , L/K)}>> \Gal(L/K)@>>>1\\
&&@VV{\Lambda_\mK}V@VV{\Lambda_\mK}V @VV{\Id}V\\
1@>>>H_K^\mK/\PKm@>>>\DKm/\PKm@>{\artin{L/K}{\ }}>>
\Gal(L/K)@>>>1
\end{CD}
\]
donde $\Lambda_\mK$ es el mapeo inducido por $\Lambda$ en el Teorema
{\rm{\ref{CClaseT4.7.3}}}.
\end{teorema}

\begin{proof} El isomorfismo $\Lambda_\mK\colon C_K/\CKm\longrightarrow
\DKm/\PKm$ da lugar a un diagrama conmutativo
\[
\begin{CD}
C_K/\CKm@>{(\underline{\ }, L/K)}>> \Gal(L/K)\\
@V{\Lambda_\mK}V\cong V@VV{\Id}V\\
\DKm/\PKm@>>{f}>\Gal(L/K)
\end{CD}
\]
donde $f$ es el homomorfismo que hace conmutativo el diagrama,
es decir $f=(\underline{\ \ }, L/K)\circ \Lambda_m^{-1}$.

Veremos que $f$ est\'a dado por el s{\'\i}mbolo de Artin.

Como se vio anteriormente, cada clase $C_K/\CKm$ est\'a
representado por un id\`ele en $\Jm=\{\vec\alpha\in J_K\mid
\alpha_\pK=1 \text{\ para\ } \pK|\mK\}$. En particular el grupo
$C_K/\CKm$ est\'a generado por las clases de los id\`eles
$\enc{\pi_\pK}{\pK}$ donde $\pK$
es un primo finito que no divide
a $\mK$ y $\pi_\pK$ es un elemento primo de $K_\pK$.

Sea $c\in C_K/\CKm$ la clase de $\enc{\pi_\pK}{\pK}$.
Entonces $\Lambda_\mK(c)=\pK\bmod \PKm$ y
$f(\Lambda_\mK(c))=(c,L/K)=\big(\enc{\pi_\pK}{\pK}, L/K\big)=
\artinp{L/K}{\pK}$. Por lo tanto $f$ es inducido por el s{\'\i}mbolo
de Artin y $\artinp{L/K}{\ }\colon \DKm\to\Gal(L/K)$ es un epimorfismo.

Queda por probar que la imagen de $\N_{L/K} C_L$ bajo
$\Lambda_\mK\colon C_K\to \DKm/\PKm$ es el grupo $H_K^\mK/\PKm$
pues esta imagen corresponde a $\ker f=
\ker \artinp{L/K}{\ }$. M\'as precisamente,
\[
\begin{CD}
1@>>>\N_{L/K}C_L@>i>>C_K@>{(\underline{\ },L/K)}>>
\Gal(L/K)@>>>1\\
&&@V{\Lambda_{\mK}}VV@VV{\Lambda_{\mK}}V@VV{\Id}V\\
1@>>>B@>{\psi}>>D_K^{\mK}/P_K^{\mK}@>{\varphi}>>
\Gal(L/K)@>>>1
\end{CD}
\]
por lo que $\ker\varphi=B=\im \psi=\im\Lambda_{\mK}$.

Definimos $J_L^{\langle\mK\rangle}=\{\vec\alpha\in J_L\mid
\alpha_\pL=1 \text{\ para\ } \pL|\mK\}$ y se tiene que
$J_L=J_L^{\langle \mK\rangle}J_L^{\mK}\L$. Por lo tanto,
ya que $\N_{L/K} J_L^{\mK}\subseteq \JKm$ y que $\N_{L/K}\*L
\subseteq \*K$, se tiene
\[
\frac{\N_{L/K}C_L}{\CKm}=\frac{\N_{L/K}(J_L)K^{\ast}}{\JKm K^{\ast}}=\frac{(\N_{L/K}(
J_L^{\langle\mK\rangle}))\JKm K^{\ast}}{\JKm K^{\ast}}.
\]

El isomorfismo
\[
\Lambda_\mK\colon\frac{C_K}{\CKm}\cong\frac{\Jm\JKm K^{\ast}}{\JKm K^{\ast}}
\xrightarrow{\ \cong\ }\frac{\DKm}{\PKm}
\]
asigna a la clase $\vec\alpha\in \Jm$ la clase del ideal $\Lambda(\vec
\alpha)={\eu a}_{\vec\alpha}=\prod_{\pK\nmid \infty}\pK^{v_\pK(
\alpha_\pK)}\in \DKm$.

Los elementos de $\frac{\N_{L/K}C_L}{\CKm}$ son las clases
representadas por la norma de id\`eles $\N_{L/K}J_L^{\langle
\mK\rangle}$ en $\Jm$. Por tanto son mapeados precisamente
sobre las clases de normas de ideales $\N_{L/K} D_L^\mK$
en $\DKm$ por lo que
\begin{gather*}
\Lambda_\mK\Big(\frac{\N_{L/K}C_L}{\CKm}
\Big)=\frac{(\N_{L/K} D_L^\mK)\PKm}{\PKm}=\frac{H_K^\mK}
{\PKm}. \tag*{$\fin$}
\end{gather*}
\end{proof}

\begin{corolario}\label{CClaseC4.7.5} El s{\'\i}mbolo 
de Artin $\artinp{L/K}{\eu a}$,
${\eu a}\in \DKm$ depende \'unicamente de la clase ${\eu a}\bmod
\PKm$ y da lugar a un isomorfismo en campos num\'ericos
\begin{gather*}
\artinp{L/K}{\ }\colon \frac{\DKm}{H_K^\mK}\xrightarrow{\ \cong\ }\Gal(L/K).
\tag*{$\fin$}
\end{gather*}
\end{corolario}

\begin{observacion}\label{CClaseO4.7.6} El Teorema \ref{CClaseT4.7.4} pone en
evidencia que, a diferencia con el grupo de clases de id\`eles, en
donde para cada extensi\'on abeliana finita $L/K$ correspond{\'\i}a
un \'unico subgrupo de $C_K$, a saber, $\N_{L/K}C_L$, cuando
usamos ideales, para cada m\'odulus $\mK$ tal que $L\subseteq
K^\mK$, nos corresponde un grupo $H_K^\mK$ y por tanto no tenemos
unicidad.
\end{observacion}

\begin{definicion}\label{CClaseD4.7.7} 
El grupo $H_K^\mK\label{CClaseHm1}=\N_{L/K}(\DKm)
\cdot \PKm$ se llama el {\em grupo de ideales declarado m\'odulo
$\mK$\index{grupo declarado por un m\'odulus}} correspondiente
a $L/K$.
\end{definicion}

Como consecuencia del Teorema de existencia del TCCG ($L
\longleftrightarrow \NN L=\N_{L/K}C_L$), se sigue que el mapeo
$L\longmapsto H_K^\mK$ da lugar a una correspondencia biyectiva
entre los campos de clases de rayos m\'odulo $\mK$ y los subgrupos
de $\DKm$ que contienen a $\PKm$: $\PKm\subseteq H_K^\mK
\subseteq \DKm$ en campos num\'ericos.

El siguiente teorema establece que la descomposici\'on de un
primo $\pK$ no ramificado en $L$, se puede leer directamente del 
grupo de ideales $H_K^\mK$ que determinan a $L$.

\begin{teorema}[Ley de descomposici\'on de primos\index{ley de
descomposici\'on de primos}\index{descomposici\'on de 
primos}]\label{CClaseT4.7.8}
Sea $L/K$ una extensi\'on abeliana de grado $n$ de campos 
num\'ericos y sea $\pKK$ un ideal primo no ramificado en $L$.
Sea $\mK$ un m\'odulus de definici\'on de $L$, esto es, $L
\subseteq K^\mK$, el cual no es divisible por $\pK$ (por ejemplo,
se puede tomar como $\mK$ al conductor) y sea $H_K^\mK$ el grupo
de ideales correspondiente a $L$.

Si $f$ es el orden de $\pK\bmod H_K^\mK$ en $\DKm/H_K^\mK$, esto
es, $f$ es el m{\'\i}nimo n\'umero natural tal que $\pK^f\in H_K^\mK$,
entonces $\pK$ se descompone en un producto
\[
\pK=\pL_1\cdots \pL_h
\]
de $h=n/f$ primos distintos $\pL_1,\ldots, \pL_h$ de grado
$f$ sobre $\pK$, donde se tiene $f=[\o_L/\pL_i:\o_K/\pK]$.
\end{teorema}

\begin{proof} Sea $\pK=\pL_1\cdots \pL_h$ la descomposici\'on de $\pK$
en $L$. Puesto que $\pK$ es no ramificada, los $\pL_i$ son
distintos y de grado igual a $f^\prime$, donde $f^\prime$ es el orden
del automorfismo de Frobenius $\Fro\pK=\artinp{L/K}{\pK}$.
Puesto que $\DKm/H_K^\mK\cong \Gal(L/K)$, $f^\prime$ es el orden
de $\pK\bmod H_K^\mK$ en $\DKm/H_K^\mK$.
$\fin$
\end{proof}

\begin{corolario}\label{CClaseC4.7.9} Si $\pK$ es no ramificado en
$L/K$, entonces $\pK$ se descompone
totalmente en $L/K \iff \pK\in H_K^\mK$ (es decir, $\iff f=1$). $\fin$
\end{corolario}

\begin{corolario}\label{CClaseC4.7.9+1} Sea $L/K$ es una extensi\'on
abeliana finita de campos num\'ericos. Entonces hay una infinidad de
lugares de $K$ totalmente descompuestos en $L$. $\fin$
\end{corolario}

\begin{proposicion}\label{CClaseP4.7.10} Sea $K$ un campo num\'erico.
Sea $K_H$ el campo de clase de Hilbert. Entonces $K_H=K^1$
es el campo de clases de rayos m\'odulo $1$. Entonces $\pK$
se descompone totalmente en $K_H/K\iff \pK$ es principal.
\end{proposicion}

\begin{proof} En este caso tenemos $C_K/C_K^1\cong I_K=D_K^1/P_K^1
\isomorfo\limits_{\substack{\uparrow\\ {\mK}=1}} D_K^{\mK}/H_K^{
\mK}$, por lo que $H^1=P_K$,
es decir $D_K=D_K^1$ y $P_K=P_K^1$. Por tanto $\pK$ se
descompone totalmente $\iff \pK\in P_K\iff \pK$ es principal. $\fin$
\end{proof}

Finalmente tenemos el siguiente resultado, conjeturado por Hilbert 
y cuya demostraci\'on es complicada. Se reduce a probar que el
mapeo de transferencia $\Ver\colon G/G^\prime\longrightarrow
G^\prime/G^{\prime\prime}$ es el mapeo trivial.

\begin{teorema}[Teorema del ideal principal\index{teorema del
ideal principal}\index{ideal principal!teorema
del $\sim$}]\label{CClaseT4.7.11}
Todo ideal ${\eu a}$ de $K$ se hace principal en el campo de
clase de Hilbert. 
\end{teorema}
\begin{proof} \cite{Iya34}. $\fin$
\end{proof}

\begin{definicion}[Definici\'on \ref{CClaseD4.5.14}]\label{CClaseD4.7.12}
Sea $K$ un campo num\'erico.
El {\em campo de clase de Hilbert extendido\index{Hilbert!campo
de clase extendido de $\sim$}\index{campo de clase de Hilbert extendido}}
$K_{H^+}$ es la m\'axima extensi\'on abeliana de $K$ no ramificada
en ning\'un primo finito.
\end{definicion}

Se tiene que $K_{H^+}$ es el campo de clase de rayos de $\mK=
\prod_{\substack{\pK\text{\ es real}\\ \pKK}}\pK$. En particular
$\mK_f=1$, $\mK_{\infty}=\infty =\prod_{\pK \text{\ real}}\pK$.

Ahora, $\DKm=D_K^1=D_K$ y $\PKm=\{(\alpha)\in P_K\mid
\alpha \text{\ es totalmente positivo}\}= P_K^\infty=P_K^+$.
As{\'\i} $P_K^+\subseteq P_K\subseteq D_K$ y
\begin{gather*}
\Gal(K_{H^+}/K)\cong D_K/P_K^+\twoheadrightarrow \Gal(K_H/K)
\cong D_K/P_K= I_K.\\
\intertext{y la siguiente sucesi\'on es exacta}
1\to \frac{P_K}{P_K^+}\to \frac{D_K}{P_K^+}\to
\frac{D_K}{P_K}\to 1\qquad\qquad
\xymatrix{
K_{H^+}\ar@{-}[d]_{\frac{P_K}{P_K^+}}^{\text{$2$ grupo elemental}}
\ar@/_3pc/@{-}[dd]_{\frac{D_K}{P_K^+}}\\
K_H\ar@{-}[d]^{\frac{D_K}{P_K}}\\ K}
\end{gather*}

Notemos que si $r$ es el n\'umero de lugares reales en $K$,
entonces $P_K/P_K^+$ es un subgrupo del $2$--grupo elemental
$C_2^r$ y, por supuesto, la contenci\'on puede ser propia.

\begin{proposicion}\label{CClaseP4.7.10'} Para $K$ un campo num\'erico,
sea $K_{H^+}$ el campo de clase de Hilbert extendido. Entonces
$\pKK$, $\pK$ finito, se descompone
totalmente en $K_{H^+}$ si y solamente
si $\pK$ es principal generado por un elemento totalmente positivo.
\end{proposicion}

\begin{proof} El campo $K_{H^+}$ corresponde a 
$C_K^{1^+}$, donde $1^+=\prod_{\pK\text{\ real}}\pK$
y se tiene
$C_K/C_K^{1^+}\cong D_K/P_K^+$. Esto implica que $H^{1^+}=
P_K^{1^+}=P_K^+$. Por tanto $\pK$ se descompone totalmente $\iff
\pK\in P_K^+\iff \pK$ es principal generado por un
elemento totalmente positivo. $\fin$
\end{proof}

\begin{observacion}\label{CClaseO4.7.10'(1)}
Notemos que $K_{H^+}/K$ es una extensi\'on finita a pesar 
de que permitimos a los primos infinitos ramificarse.
Esto no se cumple para los primos finitos. Por ejemplo,
se tiene que
\[
\bigcup_{n=1}^{\infty}\cic pn^+=\cic p{\infty}^+
\]
es una extensi\'on abeliana infinita de ${\ma Q}$
en la cual solamente el primo $p$ se ramifica.
\end{observacion}

\subsection{Sobre la extensi\'on $K_{H^+}/K_H$}\label{S17.7.5}

En esta subsecci\'on profundizamos m\'as sobre la
extensi\'on $K_{H^+}/K_H$.

\begin{observacion}\label{OClaseE4.7.14}
Existen campos num\'ericos $K$ tales que $[K_{H^+}:K_H]=2^s$
con $s>1$. Dos trabajos en donde se discuten problemas de
densidad relacionados con este fen\'omeno son \cite{DumVoi2018,
Hil2006}.
\end{observacion}

\begin{ejemplos}\label{OClaseE4.714'}
Con respecto a la Observaci\'on \ref{OClaseE4.7.14}, usamos el
programa Sage para hacer los c\'alculos de los ejemplos
proporcionados en \cite{DumVoi2018, Hil2006}.
\las
\item Sea $K={\ma Q}(\alpha)$, donde $\Irr(x,\alpha,{\ma Q})=
x^5-2 x^4-32x^3+41x^2+220x-289$. El discriminante de $K$ es
$\delta_K=405,673,292,473$. El grupo de clases extendido de $K$
es $I_K^+\cong C_4\times C_4\times C_2\times C_2$ y por tanto
$\Gal(K_{H^+}/K)\cong C_4\times C_4\times C_2\times C_2$.
El grupo de clase de $K$ es $I_K\cong C_2\times C_2$ y en
particular $\Gal(K_H/K)\cong C_2\times C_2$. Se sigue que
\[
\Gal(K_{H^+}/K_H)\cong C_2\times C_2\times C_2\times C_2=C_2^4.
\]

\item Sea $K={\ma Q}(\beta)$ donde $\Irr(x,\beta,{\ma Q})=x^5-x^4
-21x^3-7x^2+68x+60$. El discriminante de $K$ es $\delta_K=
52,315,684$. El grupo de clase extendido de $K$ es $I_K^+\cong
C_4\times C_2\times C_2$ y el grupo de clase de $K$ es
$I_K\cong C_2$. Entonces
\begin{gather*}
\Gal(K_{H^+}/K)\cong C_4\times C_2\times C_2;\quad
\Gal(K_H/K)\cong C_2;\\
\Gal(K_{H^+}/K_H)\cong C_2\times C_2\times C_2=C_2^3.
\end{gather*}

\item Sea $K={\ma Q}(\gamma)$ donde $\Irr(x,\gamma,{\ma Q})=
x^5-2x^4-6x^3+8x^2+8x+1$. Entonces
\begin{gather*}
\delta_K=638,597;\quad \Gal(K_{H^+}/K)\cong C_2\times C_2;
\quad \Gal(K_H/K)\cong\{1\};\\
\Gal(K_{H^+}/K_H)\cong C_2\times C_2=C_2^2.
\end{gather*}

Los siguientes dos ejemplos, satisfacen, en ambos casos, que
$K$ es una extensi\'on c\'ubica abeliana de ${\ma Q}$ (por supuesto
con grupo de Galois, $\Gal(K/{\ma Q})\cong C_3$).

\item  Sea $f(x)=x^3-237x+316\in {\ma Z}[x]$. Entonces,
si $K$ es el campo de descomposici\'on de $f(x)$ sobre ${\ma Q}$,
entonces $K/{\ma Q}$ es abeliana y $\Gal(K/{\ma Q})\cong
C_3$. Se tiene
\begin{gather*}
\Gal(K_H/K)\cong C_2\oplus C_6
\intertext{y}
\Gal(K_{H^+}/K)\cong C_4\oplus C_{12}.
\intertext{En particular,}
\Gal(K_{H^+}/K_H)\cong C_2 \oplus C_2.
\end{gather*}

\item Sea $g(x)=x^3-x^2-234x+729\in {\ma Z}[x]$. Entonces,
si $K$ es el campo de descomposici\'on de $g(x)$ sobre ${\ma Q}$,
entonces $K/{\ma Q}$ es abeliana y $\Gal(K/{\ma Q})\cong
C_3$. Se tiene
\begin{gather*}
\Gal(K_H/K)\cong C_2\oplus C_6
\intertext{y}
\Gal(K_{H^+}/K)\cong C_4\oplus C_{12}.
\intertext{En particular,}
\Gal(K_{H^+}/K_H)\cong C_2 \oplus C_2.
\end{gather*}

\end{list}
\end{ejemplos}

Sea $K$ un campo num\'erico y sea $J$ la conjugaci\'on compleja.
Si $K/{\ma Q}$ es una extensi\'on de Galois, se tiene que $J\in
\Gal(K/{\ma Q})$ puesto que $J|_{\ma Q}=\Id$. Se tiene que $o(
J|_K)=1$ \'o $2$. M\'as a\'un, se tiene se tiene $J|_K=\Id\iff
K\subseteq {\ma R}$.

Si $o(J|_K)=2$, $\{1,J\}$ no necesariamente es normal en $G=
\Gal(K/{\ma Q})$. Por ejemplo, si $K={\ma Q}\big(\sqrt[3]{2},\zeta_3
\big)$, entonces $K/{\ma Q}$ es una extensi\'on de Galois, $\Gal(
K/{\ma Q})=\langle\sigma,\tau\rangle\cong S_3$ donde
se tiene que $\sigma
\big(\sqrt[3]{2})=\sqrt[3]{2}$ y $\sigma(\zeta_3)=\bar{\zeta}_3=
\zeta_3^2$; $\tau\big(\sqrt[3]{2}\big)=\zeta_3\sqrt[3]{2}$ y $\tau(\zeta_3)=
\zeta_3$. Entonces $K^{\langle\sigma\rangle}=K\cap {\ma R}=
{\ma Q}\big(\sqrt[3]{2}\big)$ y donde tenemos $\sigma =J|_K$.
Adem\'as tenemos $K^{\langle\tau\rangle}=\cic 3{}$ y por tanto
$\langle\tau\rangle\triangleleft G$, $\langle\sigma\rangle
\ntriangleleft G$, $o(\tau)=3$ y $o(\sigma)=2$.
\[
\xymatrix{
{\ma Q}(\sqrt[3]{2})\ar@{-}[r]^{\sigma}\ar@{-}[d]&
{\ma Q}(\sqrt[3]{2},\zeta_3)\ar@{-}[d]^{\tau}\\
{\ma Q}\ar@{-}[r]&{\ma Q}(\zeta_3)
}
\]

En general, dado $K$, consideremos $L:=K J(K)$. Entonces
$J$ act\'ua en $L$ pues $J(L)=J(K)J^2(K)=J(K) K=L$. Si $K\subseteq
{\ma R}$ entonces $L=K$. M\'as generalmente, $L=K\iff J$ act\'ua
en $K$, es decir si $J(K)=K$.

De esta forma, si $K/{\ma Q}$ es Galois, $J\in\Gal(K/{\ma Q})$ y
$[K:K^J]=1$ \'o $2$, $\Gal(K/K^J)=\langle J|_K\rangle$ y $K^J=
K\cap {\ma R}$ pero $K^J$ no necesariamente es normal sobre
${\ma Q}$.

En el caso en que $K/{\ma Q}$ no sea Galois, tenemos que
$K^J=K\cap {\ma R}$ y $[K:K^J]$ es
en general arbitrario. Por ejemplo, si
$K={\ma Q}\big(\zeta_n\sqrt[n]{2}\big)$, $n\in{\ma N}$ 
arbitrario, entonces $K^J=K\cap {\ma R}={\ma Q}$ y $[K:
K^J]=[K:{\ma Q}]=n$.

Como lo establecimos anteriormente, el campo de clase
de Hilbert $K_H$ de un campo num\'erico, corresponde al 
grupo de id\`eles 
\begin{gather*}
J_K^1=\prod_{\pK} U_{\pK}, \quad C_K^1=
J_K^1\*K/\*K\quad\text{y}\quad \Gal(K_H/K)\cong
C_K/C_K^1\cong J_K/J_K^1\*K
\intertext{y el campo de clase de Hilbert extendido
$K_{H^+}$ corresponde al grupo de id\`eles}
J_K^{1_+}=\prod_{\pK\text{\ real}} U_{\pK}^{(1)}\times
\prod_{\pK\text{\ no real}}U_\pK,\quad
C_K^{1_+}=J_K^{1_+}\*K/\*K\quad\text{y}\\
\Gal(K_{H^+}/K)\cong C_K/C_K^{1_+}\cong
J_K/J_K^{1_+}\*K.
\intertext{Se sigue que}
\Gal(K_{H^+}/K_H)\cong (J_K^1\*K)/(J_K^{1_+}\*K).
\end{gather*}

Se tiene la siguiente sucesi\'on exacta
\begin{gather*}
1\lra \underbracket[0pt]{\Gal(K_{H^+}/K_H)}_{
\substack{\ucong\\ (J_K^1\*K)/(J_K^{1_+}\*K)}}
\lra \underbracket[0pt]{\Gal(K_{H^+}/K)
}_{\substack{\ucong\\ J_K/J_K^{1_+}\*K}}
\lra \underbracket[0pt]{\Gal(K_H/K)}_{
\substack{\ucong\\ J_K/J_K^1\*K}}\lra 1.\\
\xymatrix{
K_{H^+}\ar@{-}[d]\ar@/^1pc/@{-}[d]^{\Gal(K_{H^+}/K_H)}
\ar@/_2pc/@{-}[dd]_{\Gal(K_{H^+}/K)}\\
K_H\ar@{-}[d]\ar@/^1pc/@{-}[d]^{\Gal(K_H/K)}\\ K
}
\end{gather*}

En resumen, tenemos
\begin{gather*}
J_K^{1}=\prod_{\pK} U_{\pK}= \prod_{\pK\text{\ real}} {\ma R}^*\times
\prod_{\pK\text{\ complejo}}\*{\ma C}
\times \prod_{\pK\nmid \infty}U_\pK,
\\
J_K^{1_+}=\prod_{\pK\text{\ real}}U_{\pK}^{(1)}
\times \prod_{\pK\text{\ no real}}
U_{\pK}=\prod_{\pK\text{\ real}} {\ma R}^+\times
\prod_{\pK\text{\ complejo}}\*{\ma C}
\times \prod_{\pK\nmid \infty}U_\pK,
\\
\frac{J_K^1}{J_K^{1_+}}\cong\prod_{\pK\text{\ real}}
\frac{{\ma R}^*}{{\ma R}^+} \cong C_2^r,
\end{gather*}
donde $r$ es el n\'umero de lugares reales de $K$.

Sea $\psi\colon J_K^1/J_K^{1_+}\lra (J_K^1\*K)/(J_K^{1_+}\*K)
\cong\Gal(K_{H^+}/K_H)$,
el mapeo natural definido por 
$x\bmod J_K^{1_+}\stackrel{\psi}{\lra} x\bmod J_K^{1_+}\*K$.
Entonces $\psi$ es un epimorfismo y $\ker \psi =(J_K^1\cap 
J_K^{1_+}\*K)/J_K^{1_+}$. Por tanto tenemos la siguiente
sucesi\'on exacta
\[
1\lra (J_K^1\cap J_K^{1_+}\*K)/J_K^{1_+}\lra (J_K^1)/(J_K^{1_+})
\cong C_2^r\lra \Gal(K_{H^+}/K_H)\lra 1.
\]

Se sigue que $\Gal(K_{H^+}/K_H)$ es un $2$--grupo elemental
abeliano de rango menor o igual al n\'umero de lugares reales $r$ de
$K$.

\begin{ejemplo}\label{CClaseE4.7.13}
Primero notemos que ${\ma Q}_{H^+}={\ma Q}_H={\ma Q}$
debido al Teorema de Minkowski el cual establece que en
cualquier extensi\'on propia $K/{\ma Q}$ hay un primo finito
ramific\'andose (ver Ejemplo \ref{CClaseE4.5.15}).

Ahorca consideremos cualquier campo ciclot\'omico 
$K_n={\ma Q}(\zeta_n)$,
$n\in{\ma N}$ y sea $K_n^+={\ma Q}(\zeta_n)^+$ el subcampo real.
Sean $h_n$ y $h_n^+$ los n\'umeros de clase de $K_n$ y $K_n^+$
respectivamente. Entonces $h_n^+|h_n$ (ver 
Teorema \ref{T7.8.35}).

Si $h_n=1$ entonces $h_n^+=1$. En este caso 
se tiene que el campo de clase
de Hilbert de $K_n$ es $K_n$ y el de $K_n^+$ es $K_n^+$.
Consideremos cualquier $n$ tal que $h_n=1$ y que $n$ 
tenga dos factores primos distintos (en ese caso $K_n/
K_n^+$ es no ramificada en todos los primos finitos
(Teorema \ref{T10.2})).

Ahora bien, si $(K_{n})_{H^+}$ y $(K_n^+)_{H^+}$ son los campos
de clase de Hilbert extendidos de $K_n$ y $K_n^+$
respectivamente, entonces,
puesto que $(K_n^+)_{H^+}/K_n^+$ es abeliana y no ramificada
en los primos finitos,
se tiene que $(K_n^+)_{H^+}\subseteq (K_n)_{H^+}=(K_n)_H
=K_n$ pues $h_n=1$ y
todos los lugares de $K_n$ son complejos.
Se sigue que $(K_n^+)_{H^+}=K_n$. En este caso, $r=\varphi(n)/2$
y $P_{K_n^+}/P_{K_n^+}^+\cong C_2$ que es un subgrupo propio
$C_2^r$ para $r>1$.

Hay $30$ campos $K_n$ con $h_n=1$ (Ejemplo \ref{Ej8.9}).
Por ejemplo podemos tomar $n=45$ y en este caso $r=12$.
\end{ejemplo}

\begin{observacion}\label{OClaseS4.8}
Sea $K/{\ma Q}$ una extensi\'on de Galois. Se tiene que tanto $K_H$
como $K_{H^+}$ son extensiones de Galois sobre ${\ma Q}$. En
efecto, consideremos $\sigma$
un monomorfismo de $K_H$ en una de sus cerraduras
algebraicas, $\sigma\colon K_H\lra \bar{K}_H$. Se tiene que
$\sigma|_K\colon K\lra \bar{K}$ y puesto que $K/{\ma Q}$
es una extensi\'on normal, se tiene que $\sigma(K)=K$. Ahora
bien, $\sigma(K_H)/\sigma(K)=K$ es 
una extensi\'on no ramificada en ning\'un
primo y es una extensi\'on abeliana, se sigue que $K_H\sigma
(K_H)/K$ es una extensi\'on abeliana no ramificada, 
lo cual implica 
que $\sigma(K_H) K_H\subseteq K_H$, es decir, $\sigma(
K_H)=K_H$ y $K_H/{\ma Q}$ es una extensi\'on normal.

Un argumento totalmente an\'alogo se aplica a $K_{H^+}/{\ma Q}$.
\end{observacion}

\begin{observacion}\label{O17.7.40} En un campo real $K/{\ma Q}$,
$\p=\infty$ puede ser ramificado en el caso de que $K$ no sea
totalmente real.
\end{observacion}

\begin{ejemplo}\label{E17.7.41}
Sea $K={\ma Q}(\sqrt[3]{2})\subseteq {\ma R}$. Se tiene que
$K$ tiene $2$ encajes complejos:
\[
\sigma_1\colon\sqrt[3]{2}\lra \sqrt[3]{2};\quad
\sigma_2\colon\sqrt[3]{2}\lra \zeta_3\sqrt[3]{2};\quad
\sigma_3\colon\sqrt[3]{2}\lra \zeta_3^2\sqrt[3]{2}.
\]
Los encajes $\sigma_2$ y $\sigma_3$ son complejos y $\sigma_1$
es real y $\p={\mc P}_{\infty,1} {\mc P}_{\infty,2}^2$ donde 
${\mc P}_{\infty,2}$ corresponde a la pareja de encajes
$\{\sigma_2,\sigma_3\}$, $\bar{\sigma}_2=\sigma_3$.
\end{ejemplo}

\begin{observacion}\label{O17.7.42} A\'un en el caso de que la
extensi\'on $K/{\ma Q}$ sea de Galois, no necesariamente
se tiene que $K_H=K_{H^+}^J$.
\end{observacion}

\begin{ejemplo}[Ver el Ejemplo \ref{Ej8.9}]\label{E17.7.43}
Sea $K={\ma Q}(\sqrt{-23})$. Se tiene que $h_K=3$ y los
campos de clase de Hilbert y de Hilbert extendido coinciden
pues $K\nsubseteq {\ma R}$. Se tiene $K_H=K_{H^+}=
{\ma Q}(\sqrt{-23},\alpha)$ donde $\alpha=
\sqrt[3]{(25+3\sqrt{69})/2}+\sqrt[3]{(25-3\sqrt{69})/2}$.
Sea $L={\ma Q}(\alpha)\subseteq {\ma R}$. Se sigue que
$L=K_{H^+}^J\neq K_H$.
\[
\xymatrix{
&K_H=K_{H^+}={\ma Q}(\sqrt{-23},\alpha)\ar@{-}[dl]_3\ar@{-}[dr]^2\\
{\ma Q}(\sqrt{-23})\ar@{-}[dr]_2&&{\ma Q}(\alpha)\ar@{-}[dl]^3\\&{\ma Q} 
}
\]
\end{ejemplo}

\section{Campos de g\'eneros v{\'\i}a campos de clase}\label{CClaseS4.10}

Los campos de g\'eneros fueron tratados en la Subsecci\'on
\ref{S12.4.0} y el Cap\'itulo \ref{Ch12*}.

Sea $K$ un campo global y sea $\K/K$ una extensi\'on finita y separable. 
Sea $\K_H$
el campo de clase de Hilbert de $\K$ (en alguna de sus versiones en el
caso de campos de funciones).

El {\em campo de g\'eneros\index{campo de 
g\'eneros}}\label{CClasecamposdegeneros}
$\g \K$ de $\K/K$ es la m\'axima extensi\'on $\K\subseteq \g \K\subseteq
\K_H$ tal que $\g \K$ es la composici\'on $\K K_1$ con $K_1/K$ abeliana. Se 
tiene que $\g \K$ depende de la definici\'on de $\K_H$ y de $K$. Se tiene
que $K_1$ es la m\'axima extensi\'on abeliana de $K$ contenida en 
$\K_H$.

Equivalentemente, $\g \K=\K K_1$ donde $K_1$ es la m\'axima extensi\'on
abeliana de $K$ tal que $\K K_1/K$ es no ramificada y satisface las condiciones
que se hayan considerado para $\K_H$ (por ejemplo, que los primos
de un subconjunto no vac{\'\i}o $S\subseteq {\ma P}_K$, se descompongan
totalmente en $\K K_1/K$).

En particular, si $\K/{\ma Q}=K$ es una extensi\'on abeliana finita, 
entonces $K_1=\g \K$ y $\g \K$
es la m\'axima extensi\'on abeliana de ${\ma Q}$ tal que $\g \K/\K$ es
no ramificada.

Si $K=\F(T)$, $\K/K$ es una extensi\'on abeliana finita y $\K_H$ se define
como la m\'axima extensi\'on abeliana no ramificada de $\K$ y tal que los
primos de $S_{\infty}(\K)$ se descomponen totalmente, donde
$S_{\infty}(\K):=\{\pK\in{\ma P}_{\K}\mid\pK\cap \F(T)=\p, \p 
\text{\ el polo de\ }$T$\}$,
entonces $\g \K$ es la m\'axima extensi\'on abeliana de $K=\F(T)$ 
tal que $\g \K/\K$ es no ramificada y $S_{\infty}(\K)$ se descompone
totalmente en $\g \K$.

\begin{teorema}\label{T17.8.1N}
Sean $\K$ un campo global y $L/\K$ una extensi\'on abeliana finita. Sea
$L_H$ el campo de clase con grupo de normas $\norm\*\K/\*\K$, donde,
\begin{align*}
\norm&=\prod_{\text{$\pK$ real}}\*L_{\pK}\times \prod_{\pK\nmid\infty}U_{\pK}
=L^1\quad \text{si $\K$ es num\'erico y}\\
\norm&=J_{L,S}^1=\prod_{\pK\in S}\*L_{\pK}\times \prod_{\pK\notin S}U_{\pK}
=L_S^1\quad \text{si $\K$ es es de funciones}.
\end{align*}

Entonces $\g L$ es la extensi\'on abeliana de $\K$ con grupo de normas
igual a $\N_{L/\K} \big(\norm\*\K/\*\K\big)$.
\end{teorema}

\begin{proof}
Puesto que $L/\K$ es abeliana, $\g L$ es la m\'axima extensi\'on abeliana
de $\K$ contenida en $L_H$. El resultado es el contenido del Teorema
\ref{T17.6.191N}.
$\fin$
\end{proof}

M\'as adelante, daremos una generalizaci\'on del Teorema \ref{T17.8.1N}.

Durante este secci\'on usaremos las siguientes definiciones y
notaciones. Usaremos $K$ para denotar el campo de funciones
racionales $K=\F(T)$ o el campo de los n\'umeros racionales ${\ma Q}$.
El \'enfasis ser\'a en campos de funciones. El campo de clase
de Hilbert $L_H$ de una extensi\'on $L/K$ 
finita y separable, ser\'a la m\'axima extensi\'on
abeliana de $L$ no ramificada y tal que los primos infinitos
se descomponen totalmente. Entendemos por los primos infinitos
los primos arquimedianos en el caso num\'erico y los lugares sobre
$\p$, el polo de $T$, en el caso de campos de funciones.
En el caso num\'erico la \'ultima
condici\'on se satisface autom\'aticamente al pedir que la 
extensi\'on sea no ramificada.

Usaremos las notaciones de la Subsecci\'on \ref{CClaseS4.9-1}.

\begin{definicion}\label{D17.9.18}
Sea $L/\K$ una extensi\'on finita y separable de campos globales.
Se define el {\em campo de g\'eneros 
$\ge L{\K}$\label{camposdegeneros} de $L$ con respecto a
$\K$\index{campo de g\'eneros}} 
como la m\'axima extensi\'on de $L$ contenida en $L_H$
tal que es de la forma $L \K_1$ donde $\K_1/\K$ es una extensi\'on
abeliana. Al m\'aximo campo $\K_1$ que satisface $\ge L{\K}=L
\K_1$ lo denotaremos por $\K_L$. En otras palabras, $\K_L$
es la m\'axima extensi\'on abeliana de $\K$ tal que $\ge L{\K}=
L\K_L$.
\end{definicion}

\begin{definicion}\label{D17.9.19}
Sea $L/\K$ una extensi\'on finita y separable de campos globales.
El {\em campo de g\'eneros extendido $\ggex L{\K}$ de $L$
con respecto a $\K$\index{campo de g\'eneros extendido}} se
define como la m\'axima extensi\'on de $L$ contenida en $L_{H^+}$
(ver Definici\'on \ref{DClaseS4.9-1.17}) 
que es de la forma $L \K_{1,+}$ donde $\K_{1,+}/\K$ es una
extensi\'on abeliana. Al m\'aximo campo $\K_{1,+}$
que satisface $\ggex L{\K}
=L\K_{1,+}$ lo denotamos por $\K_{L,+}$. De esta forma, $\K_{L,+}$
es la m\'axima extensi\'on abeliana de $\K$ tal que $\ggex L{\K}
=L \K_{L,+}$.
\end{definicion}

\[
\xymatrix{
L\ar@{-}[r]\ar@{-}[d]&L\K_L=\ge L{\K}\ar@{-}[r]\ar@{-}[d]&L_H\\
\K\ar@{-}[r]&\K_L
}\quad
\xymatrix{
L\ar@{-}[r]\ar@{-}[d]&L\K_{L,+}=\ggex L{\K}\ar@{-}[r]\ar@{-}[d]&L_{H^+}\\
\K\ar@{-}[r]&\K_{L,+}
}
\]

\begin{observacion}\label{O17.9.20}
Cuando $L/\K$ es una extensi\'on abeliana, $\ge L{\K}$ es la m\'axima
extensi\'on abeliana de $\K$ contenida en $L_H$ y $\ggex L{\K}$ es
la m\'axima extensi\'on abeliana de $\K$ contenida en $L_{H^+}$.
\end{observacion}

\begin{observacion}\label{O17.9.20'}
Se tiene que $\K_L=\g{(\K_L)}$ y $\K_{L,+}=\gex{(\K_{L,+})}$.
\end{observacion}

\begin{observacion}\label{O17.9.21}
Cuando $\K=K={\ma Q}$ o $\K=K=\F(T)$ y $L/K$ es una extensi\'on
finita y separable, usaremos la notaci\'on $\g L=\ge LK$ y 
$\gex L=\ggex LK$.
\end{observacion}

El campo $\g L$ lo hemos estudiado extensivamente en el Cap\'itulo
\ref{Ch12*} para el caso de campos de funciones y fue introducido
en la Subsecci\'on \ref{S12.4.0} en el caso de campos num\'ericos.

\begin{teorema}\label{T17.9.22}
Sea $L/\K$ una extensi\'on finita y separable de campos globales.
\las
\item El subgrupo de clases de id\`eles de $\K$ que corresponde a la
extensi\'on $\K_L$ de $\K$ es la imagen 
$\N_{L/\K}(U_L\*L/\*L)$ del
subgrupo de clases de id\`eles $U_L\*L/\*L$ 
asociado a $L_H$.

\item El subgrupo de clases de id\`eles de $L$ que corresponde a la
extensi\'on abeliana $\ge L{\K}$ de $L$  es el
subgrupo $\N_{L/\K}^{-1}(\N_{L/\K}(U_L\*L/\*L))\supseteq
U_L\*L/\*L$  
del subgrupo $U_L\*L/\*L$
 asociado a $L_H$.

\item  El subgrupo de clases de id\`eles de $\K$
que corresponde a la extensi\'on $\K_{L,+}$ de $\K$
es la imagen $\N_{L/\K}(U_L^+\*L/\*L)$ del
subgrupo de clases de id\`eles $U_L^+\*L/\*L$
asociado a $L_{H^+}$.

\item El subgrupo de clases de id\`eles de $L$ que corresponde a la
extensi\'on abeliana $\ggex L{\K}$ de $L$ es el
subgrupo $\N_{L/\K}^{-1}(\N_{L/\K}(U_L^+\*L/\*L))\supseteq
U_L^+\*L/\*L$ del subgrupo 
$U_L^+\*L/\*L$ asociado a $L_{H^+}$.

\end{list}
\end{teorema}

\begin{proof}
(1) y (3). 
\[
\xymatrix{
L\ar@{-}[r]\ar@{-}[d]&\ge L{\K}\ar@{-}[r]\ar@{-}[d]&L_H\\
\K\ar@{-}[r]&\K_L
}\qquad
\xymatrix{
L\ar@{-}[r]\ar@{-}[d]&\ggex L{\K}\ar@{-}[r]\ar@{-}[d]&L_{H^+}\\
\K\ar@{-}[r]&\K_{L,+}
}
\]

Se tiene que $U_L\*L/\*L$ es el subgrupo de clases de id\`eles que
corresponde a $L_H/L$ y $U_L^+\*L/\*L$ a $L_{H^+}/L$. Entonces $\K_L$
es la m\'axima extensi\'on de $\K$ contenida en $L_H$ y $\K_{L,+}$
la respectiva en $L_{H^+}$. Por la Teorema \ref {CClaseP4.6.9}
se sigue que el grupo de normas de $\K_L$ es $\N_{L/\K}(U_L\*L/\*L)$
y el de $\K_{L,+}$ es $\N_{L/\K}(U_L^+\*L/\*L)$. Este es el contenido
de (1) y (3) a nivel de subgrupos de id\`eles.

(2) y (4). Se tiene
\[
\xymatrix{
L\ar@{-}[r]\ar@{-}[d]&L\K_L=\ge L{\K}\ar@{-}[d]\\
\K\ar@{-}[r]&\K_L
}\quad
\xymatrix{
L\ar@{-}[r]\ar@{-}[d]&L\K_{L,+}=\ggex L{\K}\ar@{-}[d]\\
\K\ar@{-}[r]&\K_{L,+}
}
\]

Por el Teorema \ref{T17.6.135N} se tiene que los diagramas
\[
\xymatrix{
C_L\ar@{>}[rr]^{\psi_{\g L/L}\phantom{xx}}\ar@{>}[d]^{\N_{L/\K}}&&
\Gal(\g L/L)\ar@{>}[d]^{\rest}\\
C_K\ar@{>}[rr]_{\psi_{\K_L/\K}\phantom{xx}}&&\Gal(\K_L/\K)
}\qquad
\xymatrix{
C_L\ar@{>}[rr]^{\psi_{\gex L/L}\phantom{xx}}\ar@{>}[d]^{\N_{L/\K}}&&
\Gal(\gex L/L)\ar@{>}[d]^{\rest}\\
C_K\ar@{>}[rr]_{\psi_{\K_{L,+}/\K}\phantom{xx}}&&\Gal(\K_{L,+}/\K)
}
\]
son conmutativos. Entonces
\[
\rest|_{\K_L}\circ \psi_{\g L/L}=\psi_{\K_L/\K}\circ \N_{L/\K}
\text{\ y\ }\rest|_{\K_{L,+}}\circ \psi_{\gex L/L}
=\psi_{\K_{L,+}/\K}\circ\N_{L/\K}.
\]

Como $\Gal(\g L/L)\cong \Gal(\K_L/\K_L\cap L)\subseteq \Gal
(\K_L/\K)$ y $\Gal(\gex L/L)\cong \Gal(\K_{L,+}/\K_{L,+}\cap L)
\subseteq \Gal(\K_{L,+}/\K)$, se tiene que $\rest|_{\K_L}$ y $\rest|_{
\K_{L,+}}$ son mapeos inyectivos. Se sigue que
\begin{gather*}
\vec\alpha\in\ker \psi_{\g L/L} \iff \psi_{\g L/L}(\vec \alpha)=\vec 1 \iff \\
\rest|_{\K_L} \circ \psi_{\g L/L}(\vec \alpha)=\vec 1=
\psi_{\K_L/\K}\circ \N_{L/\K}(\vec \alpha)\iff \\
\N_{L/\K}(\vec \alpha)\in
\ker \psi_{\K_L/\K}\iff \vec\alpha\in \N_{L/\K}^{-1}(\ker \psi_{\K_L/\K})
\intertext{esto es,}
\ker \psi_{\g L/L}=\N_{L/\K}^{-1}(\ker \psi_{\K_L/\K}).
\end{gather*}
Similarmente $\ker \psi_{\gex L/L}=\N^{-1}_{L/\K}(\ker \psi_{\K_{L,+}/\K})$.

De esta forma obtenemos que el subgrupo de clases de id\`eles correspondiente
a $\g L/L$ es $\ker \psi_{\g L/L}=\N_{L/\K}^{-1}(\ker\psi_{\K_L/\K})=
\N_{L/\K}^{-1}(\N_{L/\K}(U_L\*L/\*L))$. Similarmente para $\gex L/L$.
$\fin$
\end{proof}

Se tiene la generalizaci\'on del Corolario \ref{C12*.2.2.D}.

\begin{corolario}\label{C17.9.23}
Sea $L/\K$ una extensi\'on c\'iclica de campos globales con
$\Gal(L/\K)=G=\langle \sigma\rangle$. Sean $Cl_L\cong
\Gal(L_H/L)$, $Cl_L^+\cong\Gal(L_{H^+}/L)$. Sean los cocientes
${\mc G}_{L/\K}:=\Gal(\g L/L)$ y ${\mc G}_{L/\K}^+:=\Gal(\gex L/L)$ los
grupos de g\'eneros y de g\'eneros extendidos respectivamente.
Entonces
\[
{\mc G}_{L/\K}\cong Cl_L/Cl_L^{\langle\sigma-1\rangle}\quad\text{y}\quad
{\mc G}_{L/\K}^+\cong Cl_L^+/(Cl_L^+)^{\langle\sigma-1\rangle}.
\]
\end{corolario}

\begin{proof}
Sea $\vec \alpha\in J_L$ tal que $\N_{L/\K}\vec\alpha\in\N_{L/\K}U_L\*{\K}$.
Pongamos $\N_{L/\K}(\vec\alpha)=\N_{L/\K}(\vec u)\cdot x$ con $\vec u\in
U_L$ y $x\in\*{\K}$. Se sigue que $x=\N_{L/\K}(\vec\alpha/\vec u)\in
\*{\K}\cap\N_{L/\K}(J_L)$. Puesto que $L/\K$ es una extensi\'on c\'iclica, por
el principio local--global de Hasse, se sigue que $x$ es una norma de 
$\*L$: $x=\N_{L/\K} (y)$. Entonces
$\vec x=\N_{L/\K}(\vec y)=\N_{L/\K}(\vec \alpha/\vec u)$, esto es,
$\N_{L/\K}(\vec \alpha/\vec y\vec u)=\vec 1$. Por el Teorema 90 de Hilbert
aplicado a id\`eles (Teorema \ref{T17.6.64N}), existe $\vec \beta\in J_L$
tal que $\vec\beta^{(\sigma-1)}=\vec\alpha/\vec y\vec u$, es decir,
$\vec \alpha=\vec\beta^{(\sigma-1)}\vec y\vec u\in J_L^{\langle\sigma-
1\rangle}\*L U_L$.
Se sigue que
\begin{align*}
{\mc G}_{L/\K}&\cong\frac{J_L/\*L}{\N_{L/\K}^{-1}(\N_{L/\K}((U_L\*L)/\*L))}\cong
\frac{J_L}{J_L^{\langle\sigma-1\rangle}\*LU_L}\\
&\cong \frac{J_L/(U_L\*L)}{(J_L^{\langle\sigma-1\rangle}
U_L\*L)/(U_L\*L)}\cong \frac{J_L/U_L\*L}{(J_L/U_L\*L)^{\langle\sigma-1\rangle}}=
\cong \frac{Cl_L}{Cl_L^{\langle\sigma-1\rangle}}.
\end{align*}

Similarmente se tiene ${\mc G}_{L/\K}^+\cong 
Cl_L^+/(Cl_L^+)^{\langle\sigma-1\rangle}$.
$\fin$
\end{proof}

Como vimos en el Cap\'itulo \ref{Ch12*}, si $E\subseteq \cicl M{}$ para alg\'un
$M\in R_T$ en el caso de campos de funciones, o $E/{\ma Q}$ una extensi\'on
abeliana en el caso num\'erico, se tiene que $\g E$ corresponde a $L^+E$ 
donde, si $X$ es el grupo de caracteres de Dirichlet asociado al campo $E$,
$L$ es el campo asociado a $Y=\prod_{P\in R_T^+}X_P$ y $L^+=L\cap \cicl M{}^+$
y similarmente en el caso de campos num\'ericos.

La pregunta natural es, ?`cual es campo $\gex E$?. Veremos que $\gex E=L$,
es decir, $\gex E$ es el campo asociado a $Y$.

Para resolver este problema, primero determinaremos cual es el subgrupo
de id\`eles de $J_K$ que corresponde al campo de funciones ciclot\'omico
$\cicl M{}$ con $M\in R_T$. Sea ${\mc X}_M$ el subgrupo de id\`eles
que corresponde a $\cicl M{}$. Sea $R_T'=R_T
\setminus\{P_1,\ldots,P_r\}$, $\pi=1/T=\pi_{\infty}$
y $U_{\infty}=U_{\p}$.

Definimos
\begin{gather}\label{idelesKM}
{\mc X}_M=\prod_{i=1}^r U_{P_i}^{(\alpha_i)}\times \prod_{P\in R_T'} U_P
\times [(\pi)\times U_{\infty}^{(1)}]
\end{gather}
donde $M=P_1^{\alpha_1}\cdots P_r^{\alpha_r}\in R_T$.

\begin{observacion}\label{O17.9.23-2}
Se tiene que ${\mc X}_N\*K/\*K\cong {\mc X}_N$.

En efecto, si $x\in{\mc X}_N\cap \*K$, entonces $x\in U_P$ para toda
$P\in R_T$ por lo que $v_P(x)=0$ para toda $P\in R_T$. Puesto que
$\deg(x)_K=0$, se sigue que $v_{\infty}(x)=0$ lo cual implica que $x\in
\*\F$. Puesto que $x\in (\pi)\times U_{\infty}^{(1)}=\ker \phi_{\infty}$,
$x=1$. Por tanto, $\frac{{\mc X}_N\*K}{\*K}\cong \frac{{\mc X}_N}{
{\ma X}_N\cap\*K}\cong {\mc X}_N$.
\end{observacion}

D. Hayes prob\'o (ver \cite[Section 12.8.4,
Theorem 12.8.5]{Vil2006})
que si $U_T=\{\vec\alpha\in J_K\mid \alpha_{\p}=1 \text{\ y 
$\alpha_P\in U_P$ para toda $P\in R_T^+$}\}$, 
entonces $U_T\cong G_T=\Gal(
K_T/K)$ donde $K_T:=\bigcup_{M\in R_T}\cicl M{}$.

\begin{proposicion}\label{P17.9.23-1}
Sea $U':=\prod_{P\in R_T}U_P\times [(\pi)\times U_{\infty}^{(1)}]$.
Entonces existe un epimorfismo $\psi_M\colon U'\lra
\Gal(\cicl M{}/K):=G_M$ con $\ker \psi_M={\mc X}_M$
y por tanto, $U'/{\mc X}_M\cong G_M$.
\end{proposicion}

\begin{proof}
Sea $\vec\xi\in U'$. Entonces 
$\xi_{P_i}\in U_{P_i}=\{\sum_{j=0}^{\infty} a_jP_i^j\mid a_j\in R_T/
(P_i), a_0\neq 0\}$, $1\leq i\leq r$. Puesto que $K\subseteq K_{P_i}$ 
es denso, existe $Q_i\in R_T$ con
$Q_i\equiv \xi_{P_i}\bmod P_i^{\alpha_i}$. Por el Teorema
Chino del Residuo existe $C\in R_T$ tal que 
$C\equiv Q_i\bmod P_i^{\alpha_i}$, $1\leq i\leq r$
y por tanto $C\equiv \xi_{P_i}\bmod P_i^{\alpha_i}$, $1\leq i\leq r$.

Ahora bien, si $C_1\in R_T$ satisface $C_1\equiv \xi_{P_i}
\bmod P_i^{\alpha_i}$, $1\leq i\leq r$, entonces $P_i^{\alpha_i}
|C-C_1$ para $1\leq i\leq r$. Se sigue que $M|C-C_1$ y por
tanto $C\in R_T$ es \'unico m\'odulo $M$. Por otro lado,
$v_{P_i}(\xi_{P_i})=0$, entonces $P_i\nmid \xi_{P_i}$ 
de donde se tiene $\mcd(C,M)=1$. De esta forma tenemos
que $C\bmod M$ define un elemento de $G_M$.

Dado $\sigma\in G_M$, existe $C\in R_T$ tal que $\sigma
\lambda_M=\lambda_M^C$ donde $\lambda_M$
es un generador de $\Lambda_M$. Sea $\vec\xi\in U'$
con $\xi_{P_i}=C$, $1\leq i\leq r$ y $\xi_P=1=\xi_{\infty}$
para toda $P\in R_T'$. Por tanto $\vec\xi\mapsto C\bmod
M$ y $\psi_M$ es suprayectiva. Finalmente, $\ker\psi_M
=\{\vec\xi\in U'\mid\xi_{P_i}\equiv 1 \bmod P_i^{\alpha_i},
1\leq i\leq r\}={\mc X}_M$ de donde se sigue
el resultado.
$\fin$
\end{proof}

Probaremos que $U'/{\mc X}_M\cong J_K/{\mc X}_M\*K$.
Tenemos la composici\'on
\[
\xymatrix{
U'\ar@{^{(}->}[r]\ar@/_1pc/@{>}_{\mu}[rr]&
J_K\ar@{>>}[r]&J_K/{\mc X}_M \*K,
}
\]
con $\im \mu=U'{\mc X}_M\*K/{\mc X}_M\*K$ y $\ker \mu =
U'\cap {\mc X}_M\*K$.

Ahora bien, ${\mc X}_M\subseteq U'$ por lo que ${\mc X}_M
\subseteq U'\cap {\mc X}_M\*K$. Rec\'iprocamente, si $\vec\xi
\in U'\cap {\mc X}_M\*K$, entonces las componentes de
$\vec\xi$ est\'an dadas por
\begin{align*}
\xi_P&=a\cdot \beta_P, \quad P\in R_T, \quad
\vec\beta \in {\mc X}_M, \quad a\in \*K,\\
\xi_{\infty}&=a\cdot\beta_{\infty}, \quad
\beta_{\infty}\in (\pi)\times U_{\infty}^{(1)}.
\end{align*}

Puesto que $\xi_P,\beta_P\in U_P$ se tiene $v_P(\xi_P)=
v_P(\beta_P)=0$ para toda $P\in R_T$. Se sigue que 
$v_P(a)=0$ para toda $P\in R_T$. Adem\'as, como 
$\deg a=0$ entonces $v_{\infty}(a)=0$ y por tanto $a\in\*\F$.

Ahora $\xi_{\infty}, \beta_{\infty}
\in(\pi)\times U_{\infty}^{(1)}=\ker \phi_{\infty}$
por lo que $1=\phi_{\infty}(\xi_{\infty})=
\phi_{\infty} (a) \phi_{\infty} (\beta_{\infty})=
\phi_{\infty} (a)$ de donde
$a=1$. Se sigue que $\vec\xi\in{\mc X}_M$.
Por tanto $\ker \mu={\mc X}_M$ y obtenemos una inyecci\'on
$U'/{\mc X}_M\stackrel{\theta}
{\hooklongrightarrow} J_K/{\mc X}_M \*K$.

Falta verificar la suprayectividad de $\theta$,
debemos probar $J_K=U'{\mc X}_M \*K
=U'\*K$. Se tiene que $U'$ corresponde
a la m\'axima extensi\'on no ramificada en ning\'un
primo finito. Sea $L/K$ esta extensi\'on. Como $U_{\infty}^{
(1)}$ corresponde al primer grupo de ramificaci\'on y por 
tanto corresponde a la ramificaci\'on salvaje de $\p$,
se sigue que en $L/K$ hay a lo m\'as un primo ramificado,
este es moderadamente ramificado y es de grado $1$ ($\p$).
Por la Proposici\'on \ref{Palestine3.3}, $L/K$ es una
extensi\'on de constantes.

Finalmente, puesto que $d=\min\{n\in{\ma N}\mid \deg\vec
\alpha=n, \vec\alpha\in U'\}$, entonces el campo de constantes
de $L$ es $\F$ y por tanto $L=K$. Se sigue que $C_K
\cong U'$, esto es, $J_K/\*K\cong U'$ de donde obtenemos
$J_K=\*K U'$.

Hemos probado

\begin{proposicion}\label{P17.9.23-2}
Sea $M=P_1^{\alpha_1}\cdots P_r^{\alpha_r}\in R_T$. Entonces el subgrupo
de id\`eles 
de $J_K$ que corresponde al campo de funciones ciclot\'omico
$\cicl M{}$ es el subgrupo ${\mc X}_N\*K$, donde
\begin{gather*}
{\mc X}_M=\prod_{i=1}^r U_{P_i}^{(\alpha_i)}\times \prod_{P\in R_T'} U_P
\times [(\pi)\times U_{\infty}^{(1)}]
\end{gather*}
y por tanto el subgrupo de clases de id\`eles de $C_K$ que corresponde a
$\cicl M{}$ es ${\mc X}_M\*K/\*K\cong {\mc X}_M$. $\fin$
\end{proposicion}

\begin{corolario}\label{C17.9.23-3}
El subgrupo de id\`eles correspondiente
al campo $\cicl M{}^+$ es ${\mc X}_M^+\*K$ donde
\[
{\mc X}_M^+=\prod_{i=1}^r U_{P_i}^{(\alpha_i)}\times \prod_{P\in R_T'} U_P
\times \*{\Kii}.
\]
\end{corolario}

\begin{proof}
Si $L_1$ el campo asociado a ${\mc X}_M^+$, el campo de constantes
de $L_1$ es $\F$ por el argumento anterior y se tiene
\[
\frac{{\mc X}_M^+}{{\mc X}_M}\cong\frac{\*{\Kii}}{(\pi)\times
U_{\infty}^{(1)}}\cong \*\F,
\]
lo que implica que $[L:L_1]\leq q-1$. Adem\'as $\p$ se descompone
totalmente en $L_1$ y tiene \'indice de ramificaci\'on $q-1$ en $L$.
Se sigue que $[L:L_1]= q-1$ y $L_1=\cicl M{}^+$. 
$\fin$
\end{proof}

\begin{proposicion}\label{P17.9.23-3(1)}
Sea $E\subseteq \cicl M{}$. Si $\pK$ es un primo infinito
en $E$, entonces $1/T$ es una norm de $E_{\pK}$, esto es,
existe $x\in E_{\pK}$ tal que $\N_{E_{\pK}/k_{\infty}} x=1/T$.
\end{proposicion}

\begin{proof}
Puesto que $E\subseteq \cicl M{}$ se tiene que $E\subseteq \Kii
(\sqrt[q-1]{-1/T})$ (Proposici\'on \ref{P12*.2.2.G'}). Por tanto
$E_{\pK}=k_{\infty}(\sqrt[e]{-1/T})$ con $e|q-1$. Se sigue el
resultado (Teorema \ref{TClaseS4.9-1.10}).
\end{proof}

Regresamos a nuestra discusi\'on sobre $\gex E$ para $E\subseteq\cicl M{}$.
Sea $\Lambda\subseteq J_E$ el subgrupo de id\`eles de $E$ que corresponde
a $E_{H^+}$, esto es, $\N_{E_{H^+}/E}C_{E_{H^+}}=\*E \Lambda/\*E$. Entonces,
por definici\'on, podemos tomar
$\Lambda=U_E^+$. Por el Teorema \ref{T17.9.22}, el subgrupo de clases
de id\`eles de $C_K$ que corresponde $K_{E,+}=\gex E/K$ (esto es,
como $E/K$ es abeliana, $K_{E,+}$ corresponde a $\gex E$) es precisamente
$\N_{E/K}(U_E^+)\*K/\*K\cong \N_{E/K}U_E^+$ (pues $\N_{E/K}(U_E^+)
\subseteq \prod_{P\in R_T}\times \big((\pi)\times U_{\infty}^{(1)}\big)$ y $\*K\cap
\big(\prod_{P\in R_T}\times \big((\pi)\times U_{\infty}^{(1)}\big)\big)=\{1\}$).

Se tiene $U_E^+=\prod_{\pK|\infty}\ker \phi_{E_{\pK}}\times \prod_{\pK\nmid 
\infty} U_{\pK}$, donde denotamos $\infty=\p$.

Se tiene que si $\pK$ no es ramificado, por lo que $\pK\nmid \infty$ y $\pK\nmid
P_i$, $1\leq i\leq r$, donde $M=P_1^{\alpha_1}\cdots P_r^{\alpha_r}$, entonces
se tiene $\N_{E_{\pK}/K_P}(U_{\pK})=U_P$, donde $\pK\cap K=P\in R_T^+$
(Teorema \ref{CClaseT1.2.1.4}). Si $\pK|P_i$ para alg\'un $1\leq i\leq r$, se tiene
$[U_{P_i}:\N_{E_{\pK}/K_{P_i}} U_{\pK}]=e_{E_{\pK}/K_{P_i}}(\pK|P_i)=
\Phi(P^{\alpha_i}_i)=q^{(\alpha_i-1)d_i}(q-1)$ donde $d_i=\deg P_i$.

Ahora bien, puesto que $E\subseteq \cicl M{}$ obtenemos que $\g E\subseteq
\cicl M{}$ y el grupo de $J_E$ correspondiente a $\g E$ es $\*EU_E$ y por tanto
el grupo de $C_K$ que corresponde a $\g E$ es $\N_{E/K}(U_E)\*K/\*K$. 
Se tiene
\begin{gather*}
\N_{E/K} U_E=
\prod_{\pK|\p}\N_{E_{\pK}/K_{\infty}} \*{E_{\pK}}\times \prod_{P\in R_T}
\prod_{\pK|P}\N_{E_{\pK}/K_P} U_{\pK}.
\intertext{Tambi\'en se tiene}
N_{E/K} U_E^+=
\prod_{\pK|\p}\N_{E_{\pK}/K_{\infty}}(\ker \phi_{E_{\pK}})\times \prod_{P\in R_T}
\prod_{\pK|P}\N_{E_{\pK}/K_P} U_{\pK}.
\end{gather*}

Si ponemos $R_T':=R_T\setminus \{P_1,\ldots,P_r\}$, entonces
\begin{gather}\label{Ec17.9.13-1}
\Big(\prod_{P\in R_T'}U_P\times \prod_{i=1}^r U_{P_i}^{(\alpha_i)}
\Big)\*K \subseteq \Big(
\prod_{P\in R_T}\prod_{\pK|P}\N_{E_{\pK}/K_P} U_{\pK}\Big)\*K.
\end{gather}

Sea $\pK$ un primo infinito en $E$.
Como $E_{\pK}/K_{\infty}$ es totalmente ramificado, se tiene
$\N_{E_{\pK}/K_{\infty}} \*{E_{\pK}}\supseteq (\pi_{\infty})\times
U_P^{(n_0)}$ para alg\'un $n_0\in{\ma N}\cup \{0\}$ (Teorema 
\ref{CClaseT3.2.5.32}). Notemos que, como est\'a contenido en
un campo de funciones ciclot\'omico, $E_{\pK}\subseteq
\Kii(\sqrt[q-1]{-1/T})$ y por tanto $1/T$ es una norma de $E_{\pK}$
(Corolario \ref{CClaseS4.9-1.11}).

Se tiene $\*{E_{\pK}}=(\pi_{E_{\pK}})\times \*\F\times U_{\pK}^{(1)}$
y $\phi_{\pK}:=\phi_{E_{\pK}}=\phi_{\infty}\circ \N_{E_{\pK}/\Kii}$.
Si $x\in\ker\phi_{E_{\pK}}$ entonces $\phi_{\infty}(\N_{E_{\pK}/\Kii}(x))=
\phi_{\pK}(x)=1$ por lo que $\N_{E_{\pK}/\Kii} (x)\in\ker\phi_{\infty}$ lo que a su
vez implica que $x\in\N_{E_{\pK}/\Kii}^{-1}(\ker\phi_{\infty})$.
Rec\'iprocamente, si $x\in\N_{E_{\pK}/\Kii}^{-1}(\ker\phi_{\infty})$ entonces
$\N_{E_{\pK}/\Kii}(x)\in\ker\phi_{\infty}$. Se sigue que $\phi_{\infty}(\N_{E_{\pK}/\Kii}x)
=1=\phi_{\pK}(x)$ de donde se obtiene que $x\in \ker \phi_{\pK}$.
En resumen, tenemos $\ker\phi_{\pK}=\N_{E_{\pK}/\Kii}^{-1}(\ker\phi_{\infty})=
\N_{E_{\pK}/\Kii}^{-1}\big((\pi_{\infty})\times U_{\infty}^{(1)}\big)$.

Puesto que $k\subseteq E\subseteq \cicl M{}$, 
se sigue que $\Kii\subseteq E_{\pK}
\subseteq \cicl N{}_{\pL}=\Kii(\sqrt[q-1]{-1/T})$ donde $\pL$ es un primo
en $\cicl M{}$ tal que $\pL|\infty$ y $\pK=\pL\cap E$. El campo asociado a
$\ker\phi_{\infty}$ es $\Kii(\sqrt[q-1]{-1/T})$ y el campo asociado a
$\ker\phi_{\pK}=\N^{-1}_{E_{\pK}/\Kii}(\ker\phi_{\infty})$ es $E_{\pK}
(\sqrt[q-1]{-1/T})=k_{\infty}(\sqrt[q-1]{-1/T})$ (Corolario 
\ref{CClaseS4.9-1.11}). Se sigue que $\N_{E_{\pK}/\Kii}
(\ker \phi_{\pK})=\ker \phi_{\infty}$. Por tanto
\begin{align*}
\N_{E/k}U_E^+&=\prod_{P\in R_T}\prod_{\pK|P}\N_{E_{\pK}/k_{\infty}}U_{\pK}\times
\prod_{\pK|\infty}\N_{E_{\pK}/k_{\infty}}(\ker \phi_{\pK})\\
&=\prod_{P\in R_T}\prod_{\pK|P}\N_{E_{\pK}/k_{\infty}}U_{\pK}\times
(\ker \phi_{\infty}).
\end{align*}

De (\ref{Ec17.9.13-1}) obtenemos que
\begin{multline*}
\Big(\prod_{P\in R_T'}U_P\times \prod_{i=1}^r U_{P_i}^{(\alpha_i)}\times
(\langle\pi_{\infty}\rangle\times U_{\infty}^{(1)}\Big)\*K\\
\subseteq
\Big(\prod_{P\in R_T}\prod_{\pK|P}\N_{E_{\pK}/k_P} U_{\pK}\times
\prod_{\pK|\infty}\N_{E_{\pK}/k_{\infty}}(\ker \phi_{\pK})\Big)\*K.
\end{multline*}

Por tanto ${\mc X}_M\*K\subseteq \N_{E/K}(U_E^+)\*K$ de donde se sigue
que $\gex E\subseteq \cicl M{}$.

Ahora bien, $\gex E/E$ es no ramificada en los primos finitos y como $L$
es el campo asociado a $Y=\prod_{P\in R_T}X_P$, $L$ es la m\'axima extensi\'on
abeliana de $E$ contenida en $\cicl M{}$ no ramificada en los primos finitos,
de donde se sigue que $\gex E\subseteq L$.

Para probar que $L\subseteq \gex E$, basta probar que $L\subseteq E_{H^+}$
pues $L/E$ es abeliana y $\gex E\subseteq E_{H^+}$ es la m\'axima
extensi\'on abeliana de $E$ contenida en $E_{H^+}$.

Ahora bien, para probar que $L\subseteq E_{H^+}$, hay que probar que
$\N_{L/E}C_L\supseteq \N_{E_{H^+}/E} C_{E_{H^+}}=U_E^+\*E/\*E$, donde
$U_E^+=\prod_{\pK|\infty}\ker\phi_{\pK}\times \prod_{\pK\nmid\infty}U_{\pK}$.

Sea $\N_{L/E}C_L=\Lambda \*E/\*E$ donde $\Lambda =\N_{L/E}J_L$. Basta
probar que $U_E^+\subseteq \Lambda$. Puesto que $L/E$ es no ramificada
en todos los primos finitos, si $\pK$ es un primo finito de $E$ y $\pL$ es
un primo de $L$ sobre $\pK$, $\N_{L_{\pL}/E_{\pK}} U_{\pL}=U_{\pK}$
donde denotamos $U_{\pL}=U_{L_{\pL}}$ y $U_{\pK}=U_{E_{\pK}}$
(Teorema \ref{CClaseT1.2.1.4}).
En particular $\N_{L_{\pL}/E_{\pK}}\*{L_{\pL}}\supseteq U_{\pK}$ y
$\N_{L/E}J_L\supseteq \prod_{\pK\nmid \infty} U_{\pK}$.

Por otro lado, $L/\g E$ es totalmente ramificada en los primos infinitos
y los primos infinitos de $E$ son totalmente descompuestos en $\g E$,
por lo que $(\g E)_{\pK}=E_{\pK}$ para $\pK|\infty$ y elemento uniformizador
de $\pK$ en $E$ lo es para $\pK$ en $\g E$.

Por el Teorema \ref{CClaseT3.2.5.32}, $\pi_{E,\infty}:=\pi_{E_{\pK}}$ es
norma de $L_{\pL}$ con $\pL|\pK$. Si $(\underline{\ \ },L_{\pL}/E_{\pK})$
representa el mapeo local de Artin, entonces $(U_{E_{\pK}}^{(1)},L_{\pL}/
E_{\pK})=G^1(L_{\pL}/E_{\pK})$, el primer grupo de ramificaci\'on de
$L_{\pL}/E_{\pK}$. Puesto que los primos infinitos son moderadamente
ramificados en $L_{\pL}/E_{\pK}$, se sigue que $G^1(L_{\pL}/E_{\pK})=
\{1\}$ y $U_{E_{\pK}}^{(1)}\subseteq \N_{L_{\pL}/E_{\pK}}\*L_{\pL}$.

El campo de constantes de $E$ y de $K$ es $\F$ y todo primo 
infinito tiene grado $1$ en $E$, por lo que si $\pK$ es un primo infinito
de $E$, $\*{E_{\pK}}=(\pi_{E,\infty})\times \*\F\times U_{E_{\pK}}^{(1)}$.
Notemos que $[E_{\pK}:\Kii]=ef=e=[\*{\Kii}:\N_{E_{\pK}/\Kii}\*{E_{\pK}}]$
donde $e$ denota al \'indice de ramificaci\'on de $E_{\pK}/\Kii$. Adem\'as
como $E_{\pK}/\Kii$ es total y moderadamente ramificada,
$(\pi_{\infty})\times U_{E_{\pK}}^{(1)}\subseteq \N_{E_{\pK}/\Kii}(\*{E_{\pK}})$
y $\N_{E_{\pK}/\Kii}\*\F=(\*\F)^e$. Se sigue que
$(\pi_{\infty})\times (\*\F)^e\times U_{\infty}^{(1)}\subseteq \N_{E_{\pK}/
\Kii} \*{E_{\pK}}$. De $[\*{\Kii}:\N_{E_{\pK}/\Kii}\*{E_{\pK}}]=e$
obtenemos la igualdad
\[
\N_{E_{\pK}/\Kii} \*{E_{\pK}}=
(\pi_{\infty})\times (\*\F)^e\times U_{\infty}^{(1)}.
\]

Por tanto
$\phi_{\pK}(\*{E_{\pK}})=\phi_{\infty}(\N_{E_{\pK}/\Kii} \*{E_{\pK}})=
\phi_{\infty}\big((\pi_{\infty})\times (\*\F)^e\times U_{\infty}^{(1)})=
(\*\F)^e$. Se sigue que $\ker \phi_{\pK}=(\pi_{E,\infty})\times R\times
U_{E_{\pK}}^{(1)}$ donde $R=\{\lambda\in\*\F\mid \lambda^e=1\}=
(\*\F)^{(q-1)/e}$. Aqu\'i $e$ es el \'indice de ramificaci\'on de $\p$ en
$E/K$.

Notemos que si $\pL$ es un primo infinito de $L$, $\N_{L_{\pL}/E_{\pK}}
\*\F=(\*\F)^{e'}$ donde $e':=e_{\infty}(L/\g E)=e_{\infty}(L/E)$.
Entonces $e'e=e_{\infty}(L/K)$ y en particular 
$e'e|q-1$ y $e'\big|\frac{q-1}{e}$. 
Sea $\*\F=\langle\beta\rangle$ con $o(\beta)=q-1$.
Sea $\lambda\in R$, $\lambda^e=1$. Se tiene $\lambda=\beta^s$
para alg\'un $s$. Por tanto $\lambda^e=\beta^{es}=1$ y $q-1|es$.
Puesto que $e'e|q-1$ se sigue que $e'e|se$ y $e'|s$. De esta
forma se  obtiene que $\lambda=\beta^s=(\beta^{s/e'})^{e'}\in
(\*\F)^{e'}\subseteq \N_{L_{\pL}/E_{\pK}} \*{L_{\pL}}$.

Se sigue que $\ker\phi_{\pK}\subseteq \N_{L_{\pL}/E_{\pK}}\*{L_{\pL}}$
y que $U_E^+=\prod_{\pK|\infty}\ker\phi_{\pK}\times \prod_{
\pK\nmid \infty}U_{\pK}\subseteq \N_{L/K}J_L$ de donde se
obtiene $L\subseteq E_{H^+}$ y por tanto $L=\gex E$.

Hemos probado

\begin{teorema}\label{T17.9.24} Sea $E\subseteq \cicl M{}$. Entonces
el campo de g\'eneros extendido $\gex E$ de $E$ relativo a $K$ es el
campo asociado al grupo de caracteres de Dirichlet $Y=\prod_{P\in
R_T}X_P$, donde $X$ es el grupo de caracteres de Dirichlet
asociado al campo $E$. $\fin$
\end{teorema}

\subsection{Extensiones abelianas finitas}\label{S17.9.1}

Usaremos las notaciones del Cap\'itulo \ref{Ch12*}.
Consideremos ${\K}/K={\ma F}_q(T)$ una extensi\'on abeliana
finita. Sean $n\in{\ma N}\cup\{0\}$, $m\in{\ma N}$ y
$N\in R_T$ tales que ${\K}\subseteq {_n\cicl N{}_m}$. Sea 
$E:={\K}M\cap \cicl N{}$, donde $M=L_nK_m$. 
Tenemos que $E{\K}/{\K}$ es una extensi\'on de constantes,
en particular una extensi\'on no ramificada (ver Teorema \ref{T2.1.A}). 
Sea $H$ el grupo de descomposici\'on de los pirmos
infinitos de ${\K}$ en $E{\K}/{\K}$. Se tiene que $H$
es can\'onicamente isomorfo con el grupo de descomposici\'on
de los primos infinitos de $\g E {\K}/{\K}$. Tenemos que 
$|H|$ es igual al grado de inercia de los primos infinitos de ${\K}$ en
cualquiera de $E{\K}/{\K}$ o de $\g E{\K}/{\K}$. 
Tambi\'en se tiene que $\g {\K}=\g E^H {\K}$.

\begin{teorema}\label{T17.9.1.1}
Con las notaciones anteriores, tenemos que
$\gex {\K} = D{\K}$ con $\gex{(\g{E^H})}
\subseteq D\subseteq \gex E$. En particular, cuando
$H=\{1\}$, tenemos $\gex {\K}=\gex E {\K}$.
\end{teorema}

\begin{proof}
Sabemos que $\g EM=\g {\K}M$ (ver la demostraci\'on
del Teorema \ref{T2.1.A}). Sea 
$C:=\gex {\K} M\cap \cicl N{}\supseteq \g {\K} M\cap \cicl N{}=\g E M\cap
\cicl N{}\supseteq \g E\cap \cicl N{}=\g E$.
\[
\xymatrix{
\cicl N{}\ar@{-}[d]\\ C\ar@{-}[rr]\ar@{-}[d]&&CM=\gex {\K} M\ar@{-}[d]\\
\g E\ar@{-}[rr]&&\g EM=\g {\K} M
}
\]

Notemos que $\gex {\K} M/\g {\K} M$ es no ramificada en los primos
finitos pues esto mismo se cumple en $\gex {\K}/\g {\K}$. Tambi\'en se tiene
que $\g E M/\g E$ es no ramificada en los primos finitos.
En particular $C/\g E$ es no ramificada en los primos finitos
y $C\subseteq \cicl N{}$. Puesto que $\gex E/K$ es la m\'axima extensi\'on
abeliana con $\gex E/E$ no ramificada en los
primos finitos y contenida en un campo de funciones ciclot\'omico,
se sigue que $C\subseteq \gex E$.

Por lo tanto $\g E\subseteq C\subseteq \gex E$ y
$\begin{array}{c}\g E M\\ \|\\ \g {\K} M\end{array} \subseteq \begin{array}{c}
CM\\ \| \\ \gex {\K} M\end{array}\subseteq \begin{array}{c}\gex EM\\\|\\ 
\gex E M\end{array}$. En particular $\gex {\K} M\subseteq \gex E M$.

\[
\xymatrix{
&\cicl N{}\ar@{-}[d]\\
&\gex E\ar@{-}[rr]\ar@{-}[d]\ar@{-}[ldd]&&\gex EM\ar@{-}[d]\\
&C\ar@{-}[rr]\ar@{-}[dd]&&\gex {\K} M\ar@{-}[dl]\\
D\ar@{-}[rr]\ar@{-}[dd]&&\gex {\K}=D{\K}\ar@{-}[dd]\\
&\g E\ar@{-}[dl]\ar@{-}[rr]&&\g E M=\g {\K} M\ar@{-}[uu]\ar@{-}[dl]\\
\g {E^H}\ar@{-}[rr]&& \g {\K}=\g{E^H}{\K}
}
\]

Sea $D:=\gex {\K}\cap \cicl N{}=\g{E^H} {\K}\cap \cicl N{}\supseteq \g{E^H}
\cap \cicl N{} =\g{E^H}$. Por tanto $\g{E^H}\subseteq D$ y $\gex {\K}=D{\K}$.

Se tiene que la extensi\'on $EM/E$ es no ramificada en los primos
finitos. Por lo tanto $DEM=
D{\K}M/EM={\K}M$ es no ramificada en los primos finitos puesto que
$\gex {\K}/{\K}$ satisface esto mismo. Se sigue que $DE/E$ 
es no ramificada en los primos finitos por lo que $D
\subseteq \gex E$. De esta forma, obtenemos que
$\gex {\K}=D{\K}\subseteq \gex E {\K}$.

En el caso particular de que $H=\{1\}$ se tiene que $E\subseteq \g {\K}
\subseteq \gex {\K}$ de tal forma que
$\gex E\subseteq \gex {\K}$ y $\gex {\K}=\gex E {\K}$.

En el caso general, $\g{E^H}\subseteq D\subseteq \gex E$ y $\gex E/
\g{E^H}$ es totalmente ramificada en los primos infinitos. Puesto que
$\g{E^H}\subseteq \g {\K}$, se sigue que
$\gex{(\g{E^H})}\subseteq \gex {\K}$. Por tanto $\gex{(\g{
E^H})}\subseteq D$. $\fin$
\end{proof}

\begin{ejemplo}\label{Ej17.9.1.3}
Sean $\K:=K(\sqrt[l]{\gamma D})$ donde $l$ es un n\'umero primo tal que
$l|q-1$, $D\in R_T$, $D$ m\'onico con $l\nmid\deg D$, 
y $\gamma\not\equiv (-1)^{\deg D}
\bmod (\*\F)^l$. Entonces $E=M\K\cap \cicl D{}=K(\sqrt[l]{\*D})$ donde $M=K_l$
y $D^*=(-1)^{\deg D} D$. Entonces
$H\cong C_l$, $\gex E=\g E= E$, $\g{E^H}=E^H=K$, $\g \K=\K$.

Sea $n:=\deg D$, $D=T^n+a_{n-1}T^{n-1}+\cdots +a_1T+a_0$, $n=lm-r$ con
$0<r<l$. El primo infinito $\p$ se ramifica in $\K/K$,
por lo que solamente hay un \'unico primo infinito 
$\pL_{\infty}$ en $\K$. Por tanto, tenemos que 
$U_\K=\*{\K_{\infty}}\times \prod_{\pK\nmid
\infty} U_{\pK}$ y $U_\K^+=\ker \phi_{\K_{\infty}}\times \prod_{\pK\nmid
\infty} U_{\pK}$, donde $\K_{\infty}=\K_{\pL_{\infty}}$. Ahora
\[
\K_{\infty}=K_{\infty}(\sqrt[l]{\gamma D})=K_{\infty}(\sqrt[l]{\gamma T^n u})=
K_{\infty}\big(\sqrt[l]{\gamma T^{-r} (T^m)^l u}\big),
\]
donde $u=1+a_{n-1}(1/T)+\cdots+a_1(1/T)^{n-1}+a_0(1/T)^n\in U_{\K_{
\infty}}^{(1)}$. Puesto que $p\neq l$, existe $v\in U_{\K_{\infty}}^{(1)}$ 
tal que $u=v^l$. Se sigue que $\K_{\infty}=K_{\infty}(\sqrt[l]{\gamma T^{-r}})
=K_{\infty}(\sqrt[l]{\delta/T})$ para alg\'un $\delta\in\*\F$, $\delta\not\equiv
(-1)\bmod (\*\F)^l$.

Un elemento primo en $\K_{\infty}$ es $\pi^{\star}:=
\pi_{\K_{\infty}}=\sqrt[l]{\delta/T}$,
$(\pi^{\star})^l=\delta/T=\delta\pi_{\infty}$. De aqu\'i que para un
elemento arbitrario $x\in \*{\K_{\infty}}$, digamos
$x=(\pi^{\star})^m \xi w$, $m\in {\mathbb Z}$, $\xi\in\*\F$,
$w\in U_{\K_{\infty}}^{(1)}$, tenemos que
$N_{\K_{\infty}/K_{\infty}}(x)=((-1)^{l-1}\delta\pi_{\infty})^m
\xi^l w'$, con $w'\in U_{\infty}^{(1)}$, de tal forma que
$\phi_{\K_{\infty}}(x)=
((-1)^{l-1}\delta)^m \xi^l$. Se sigue que
si $x\in\ker\phi_{\K_{\infty}}$, entonces
$l|m$. Tambi\'en se tiene que
$\delta^{-1}(\pi^{\star})^l\in \ker\phi_{\K_{\infty}}$. 

Por tanto
$\min\{m\in{\mathbb N}\mid \text{existe $\vec\alpha\in U_K^{+}$ con\ }
\deg\vec\alpha=m\}=l$. Se sigue que el campo de constantes 
de $H_\K^+$ y de $\gex \K$ es ${\mathbb F}_{q^l}$.

Puesto que el campo de constantes de $\g\K$ es 
$\F$, se obtiene que $\gex \K= \K_{{\eu{ge}}, l}= \gex E K$.
\end{ejemplo}

\begin{ejemplo}\label{Ej17.9.1.4}
Sean $l$ un n\'umero primo tal que $l^2|q-1$ y ${\K}=K(\sqrt[l^2]{\gamma D})$
donde $\gamma \in \*\F$, $D=P_1^{\alpha_1}\cdots P_r^{\alpha_r}\in R_T$,
$P_1,\ldots, P_r\in R_T^+$, $1\leq \alpha_i\leq l^2-1$, $1\leq i\leq r$, $r\geq 1$.
Suponemos que $\deg D=ld$ con $l\nmid d$, que $l\nmid \deg P_1$, que
$\mcd(\alpha_i,l)=1$, $1\leq i\leq s$, $s\geq 1$ y que $l|
\alpha_j$, $s+1\leq j\leq r$. Entonces $e_{P_i}({\K}/K)=l^2$ para $1\leq i\leq s$,
$e_{P_j}({\K}/K)=l$ para $s+1\leq j\leq r$, y $e_{\infty}({\K}/K)=l$.
Puesto que $e_{P_1}({\K}/K)=l^2$, el campo de constantes de
${\K}$ es $\F$. Suponemos que $(-1)^{\deg D}\gamma \notin
(\*\F)^l$. Por lo tanto $\F(\sqrt[l^2]{\varepsilon})={\mathbb F}_{q^{l^2}}$ donde
$\varepsilon=(-1)^{\deg D}\gamma$.

Sean $M={\mathbb F}_{q^{l^2}}(T)=K_{l^2}$ y $E={\K}M\cap \cicl D{}$. Entonces
$E=K(\sqrt[l^2]{(-1)^{\deg D}D})=K(\sqrt[l^2]{\*D})\subseteq \cicl D{}$. Se tiene
\[
E{\K}=E(\sqrt[l^2]{\varepsilon})={\K}(\sqrt[l^2]{\varepsilon})=E_{l^2}={\K}_{l^2}.
\]
De esta forma obtenemos que
$f_{\infty}({\K}/K)=[\F(\sqrt[l]{\varepsilon}):\F]=l$ y por tanto
$f_{\infty}(E{\K}/{\K})=\frac{f_{\infty}(E{\K}/K)}
{f_{\infty}({\K}/K)}=\frac{l^2}{l}=l$.
Se sigue que $H \cong C_l$ y por tanto $|H|=l$.

Tambi\'en se tiene que
$\gex E=K\big(\sqrt[l^2]{\*{P_1}},\ldots,\sqrt[l^2]{\*{P_s}},
\sqrt[l]{\*{P_{s+1}}},\ldots \sqrt[l]{\*{P_r}}\big)$. Puesto que 
$l\nmid \deg P_1$, se sigue que
$e_{\infty}\big(K\big(\sqrt[l^2]{\*{P_1}}\big)/K\big)=l^2$ 
y por tanto $e_{\infty}(\gex E/
K)=l^2$. Puesto que $\deg D=ld$ con
$l\nmid d$, tenemos que $e_{\infty}(\g E/K)=l$.
Ahora
\[
\deg D=ld =\sum_{i=1}^r \alpha_i\deg P_i =\alpha_i \deg P_1+\sum_{i=2}^s
\alpha_i\deg P_i+\sum_{j=s+1}^r \alpha_j \deg P_j,
\]
y puesto que $l|\sum_{j=s+1}^r\alpha_j\deg P_j$ y $\mcd(\alpha_i,l)=1$
para $1\leq i\leq s$, se sigue que existe $2\leq i\leq s$ con $l\nmid
\deg P_i$. Digamos que $l\nmid \deg P_2$.

Sean $a,b\in{\ma Z}$ tales que 
$a\deg P_1+bl^2=1$ (en particular, $\mcd (a,l)=1$).
Por tanto, para $i\geq 2$ tenemos que 
$\deg P_i-(a\deg P_i)\deg P_1=
b(\deg P_i)l^2$. Sea $Q_i=P_i P_1^{z_i}$ con 
$z_i:=-a\deg P_i$ para $2
\leq i\leq r$. Sea 
\[
L:=K\big(\sqrt[l]{\*{P_1}},\sqrt[l^2]{Q_2},\ldots, 
\sqrt[l^2]{Q_s},\sqrt[l]{Q_{s+1}},\ldots, \sqrt[l]{Q_r}\big).
\]
 Entonces
$e_{\infty}(L/K)=l$, $e_{P_i}(L/K)=l^2$ para $2\leq i\leq s$ y
$e_{P_j}(L/K)=l$ para $s+1\leq j\leq r$. Para
el primo $P_1$ tenemos que, puesto que
$l\nmid \deg P_2$, que $\sqrt[l^2]{Q_2}=\sqrt[l^2]{P_2P_1^{-a\deg P_2}}$
y que $\mcd (l,-a\deg P_2)=1$, entonces $e_{P_1}(K(\sqrt[l^2]{Q_2})/K)=l^2$.
De esta forma obtenemos que $e_{P_1}(L/K)=l^2$. 
Por tanto $E\subseteq L$ y $[\gex E:L]=l$. Se sigue que
\begin{gather*}
L=\g E=K\big(\sqrt[l]{\*{P_1}},\sqrt[l^2]{Q_2},\ldots, 
\sqrt[l^2]{Q_s},\sqrt[l]{Q_{s+1}},\ldots, \sqrt[l]{Q_r}\big),\\
\gex E= \gex E=K\big(\sqrt[l^2]{\*{P_1}},\ldots,\sqrt[l^2]{\*{P_s}},
\sqrt[l]{\*{P_{s+1}}},\ldots \sqrt[l]{\*{P_r}}\big),\\
[\gex E:\g E]=l=e_{\infty}(\gex E/\g E).
\end{gather*}

Ahora, $E{\K}={\K}_{l^2}$, de tal forma que
$H\cong D_{\infty}(E{\K}/{\K})\cong \Gal(
{\K}_{l^2}/{\K}_{l})$, donde $D_{\infty}$ 
denota el grupo de descomposici\'on de
los primos infinitos.
\[
\xymatrix{
E{\K}={\K}_{l^2}\ar@{-}[d]^{f_{\infty}(E{\K}/{\K}_l)=l}
\\ {\K}_l\ar@{-}[d]^{f_{\infty}({\K}_l/{\K})=1}\\ {\K}
}
\]

Tenemos
\begin{gather*}
E^H=K(\sqrt[l]{\*D})=K(\sqrt[l]{D})\quad \text{y}\\
 \g{E^H}
=K\big(\sqrt[l^2]{Q_2},\ldots, 
\sqrt[l^2]{Q_s},\sqrt[l]{Q_{s+1}},\ldots, \sqrt[l]{Q_r}\big).
\end{gather*}

Tambi\'en se obtiene que $\gex{(\g{E^H})}=\gex E$
debido a que $[\gex{(\g{E^H})}:K]=
\prod_{j=1}^re_{P_j}(\g{E^H}/K)=[\gex E:K]$ y $\g{E^H}
\subseteq \gex E$. Se sigue que 
\begin{align*}
\gex {\K}&=\gex E {\K}=K\big(
\sqrt[l^2]{\*{P_1}},\ldots,\sqrt[l^2]{\*{P_s}},
\sqrt[l]{\*{P_{s+1}}},\ldots \sqrt[l]{\*{P_r}}\sqrt[l^2]{\gamma D}\big)\\
&=\gex E(\sqrt[l^2]{\varepsilon}).
\end{align*}
\end{ejemplo}

\begin{observacion}\label{O17.9.1.5}
\las
\item Cuando $H=\{1\}$ se tiene que
$\gex {\K}={\K}^{\exte}:=\gex E {\K}$.

\item En general tenemos $f_{\infty}(\gex {\K}|\g {\K})>1$.

\item El campo de constantes ${\K}_{H^+}$ puede ser diferente al de
${\K}_H$ (Ejemplos \ref{Ej17.9.1.3} y \ref{Ej17.9.1.4}). En el Ejemplo
\ref{Ej17.9.1.4}
se tiene que ${\ma F}_{q^l}$ es el campo de constantes de
${\K}_H$ y ${\ma F}_{q^{l^2}}$ es el campo de constantes
${\K}_{H^+}$.

\item En general la extensi\'on $\gex {\K}/\g {\K}$ puede contener
subextensiones no ramificadas. En el Ejemplo \ref{Ej17.9.1.4} se
tiene que $\gex {\K}=\g {\K}(\sqrt[l^2]{P_1})$ y
$\g {\K}\subsetneqq ({\g {\K}})_{l^2}=\g {\K}(\sqrt[l]{P_1})
\subsetneqq \gex {\K}$.
Tenemos que $e_{\infty}(\gex {\K}|\g {\K})=f_{\infty}(\gex {\K}|\g {\K})=l$.
\end{list}
\end{observacion}

Para ver una descripci\'on de $\gex{\K}$ para una extensi\'on
finita y separable $\K/K$, en la mayor\'ia de los casos, ver
\cite{RaRzVi2020}.

Terminamos enunciando dos teoremas cl\'asicos de campos
de g\'eneros num\'ericos desde el punto de vista de campos
de clase globales.

Desde el punto de vista de campos de clase, los campos de g\'eneros tienen
las siguientes propiedades. En la construcci\'on de Furtw\"angler del campo
de clase de Hilbert, el siguiente teorema es muy importante.
Ver el Teorema \ref{T12*.2.2.C}.

\begin{teorema}[Furtw\"angler, 1907\index{teorema de 
Furtw\"angler}]\label{CClaseT4.10.1}
Sea $L/K$ una extensi\'on c{\'\i}clica no ramificada de campos num\'ericos
y sea $\sigma$ un generador de $G=\Gal(L/K)$.
Sea $\N\colon I_L\longrightarrow
I_K$ la norma (relativa) de clases de ideales, esto es, $\N {\eu a}={\eu b}\in
I_K$ donde ${\eu a}\in I_L$ y 
$\con_{K/L}{\eu b}=\prod_{\theta\in G}{\eu a}^{\theta}$
(la norma de ${\eu a}$).
Entonces $\ker \N=I_L^{(1-\sigma)}$.
\end{teorema}

\begin{proof}
Se tiene que $I_L\cong \Gal(L_H/L)$ y como $L/K$
es abeliana, el campo de g\'eneros
$\g L$ de $L$ sobre $K$, es la m\'axima extensi\'on
abeliana $K$ contenida en $L_H$.
\[
\begin{minipage}{6cm}
\xymatrix{
&\g L\ar@{-}[dl]_{I_L/I_L^{(1-\sigma)}}\ar@{-}[r]^{I_L^{(1-\sigma)}}
&L_H\ar@{-}@/^1pc/[dll]^{I_L}\ar@{-}@/^2pc/[ddll]^{\mc G}\\
L\ar@{-}[d]_{G=\langle\sigma\rangle}\\ K
}
\end{minipage}
\qquad\qquad
\begin{minipage}{5.5cm}
Sea ${\mc G}=\Gal(L_H/K)$ ($L_H/K$ es Galois por la
maximalidad de $L_H$).

Entonces, como $\g L$ es la m\'axima extensi\'on abeliana de 
$K$ contenida en $L_H$, $\Gal(L_H/\g L)={\mc G}'$ el
subgrupo conmutador de ${\mc G}$.

Veamos que ${\mc G}'=I_L^{(1-\sigma)}$.
\end{minipage}
\]

Se tiene que la sucesi\'on $1\lra I_L\lra {\mc G}\lra G\lra 1$ es exactta
por lo que ${\mc G}=GI_L=\{\tau\bar{\eu a}\mid \tau\in G, 
\bar{\eu a}\in I_L\}$ y $G$ act\'ua en $I_L$ por conjugaci\'on

Para $\bar{\eu a}\in I_L\cong \Gal(L_H/L)$ y para $\tau\in G$, se tiene
\[
\bar{\eu a}^{1-\tau}=(1-\tau)\bar{\eu a}=\bar{\eu a}\circ \underbrace{(\tau^{-1}
\circ \bar{\eu a})}_{\substack{\text{acci\'on}\\ \text{conjugaci\'on}}}=
\bar{\eu a}\circ\tau^{-1}\circ\bar{\eu a}\circ \tau\in {\mc G}',
\]
esto es, $I_L^{(1-\sigma)}\subseteq {\mc G}'$.

Ahora, ${\mc G}/I_L\cong G=\langle \sigma\rangle$ es abeliano, por
lo que ${\mc G}'\subseteq I_L$.

Sean $x,y\in {\mc G}$, digamos $x=\tau\bar{\eu a}$, $y=\gamma \bar{\eu b}$,
entonces
\begin{gather*}
x^{-1}=\bar{\eu a}^{-1}\tau^{-1}=\tau^{-1}\tau\bar{\eu a}^{-1}\tau^{-1}=
\tau^{-1}\bar{\eu a}^{-\tau^{-1}}\quad \text{y}\quad
y^{-1}=\gamma^{-1}\bar{\eu b}^{-\gamma^{-1}}.
\intertext{Por tanto}
xyx^{-1}y^{-1}=(\tau\gamma\bar{\eu a}^{\gamma}\bar{\eu b})(
\tau^{-1}\gamma^{-1}\bar{\eu a}^{-\tau^{-1}\gamma^{-1}}\bar{\eu b}^{
-\gamma^{-1}}).
\intertext{Notemos que si $\bar{\eu c}\in I_L$ y $\delta\in G$, $\bar{\eu c}
\delta=\delta\delta^{-1}\bar{\eu c}\delta=\delta\bar{\eu c}^{\delta}$. Se sigue
que}
\begin{align*}
xyx^{-1}y^{-1}&=\tau\gamma\tau^{-1}\gamma^{-1}\bar{\eu a}^{\gamma\tau^{-1}
\gamma^{-1}}\bar{\eu b}^{\tau^{-1}\gamma^{-1}}\bar{\eu a}^{-\tau^{-1}
\gamma^{-1}}\bar{\eu b}^{-\gamma^{-1}}\\
&=\bar{\eu a}^{\tau^{-1}(1-\gamma^{-1})}\bar{\eu b}^{\gamma^{-1}(\tau^{
-1}-1)}\in I_G I_L=I_L^{\sigma -1}.
\end{align*}
\end{gather*}
Por lo tanto ${\mc G}'\subseteq I_L^{\sigma -1}$ y ${\mc G}'=I_L^{\sigma-1}$.

Por otro lado, $\g L$ corresponde a $\N_{L/K} I_L$ lo que implica que
$\Gal(\g L/K)\cong I_K/\N_{L/K} I_L\cong {\mc G}/{\mc G}'\cong 
{\mc G}/I_L^{(1-\sigma)}$.

En nuestro caso tenemos que, como $L/K$ es no ramificada y abeliana,
entonces $L\subseteq K_H\subseteq L_H$ y $K_H/K$ es maximal, 
abeliana no ramificada, y por tanto $K_H=\g L$ y se tiene el diagrama
\[
\xymatrix{
\ker\N\ar@{->}[rr]^(.3){\cong}\ar@{->}[d]&&
\ker\rest\cong\Gal(L_H/\g L)\phantom{xxx}\ar@<-10ex>@{->}[d]\\
I_L\ar@{->}[rr]^(.3){\cong}\ar@{->}[d]_{\N}&&
\Gal(L_H/L)\phantom{xxxxxxxxxxxx}
\ar@<-10ex>@{->}[d]^{\rest}\\
I_K\ar@{->}[rr]^(.3){\cong}&&\Gal(K_H/K)\cong \Gal(\g L/K)
}
\]

Se sigue que $\ker\N\cong\Gal(L_H/\g L)
\cong {\mc G}'\cong I_L^{(1-\sigma)}$.
$\fin$
\end{proof}

\begin{observacion}\label{CClaseO4.10.1+1}
De la demostraci\'on del Teorema \ref{CClaseT4.10.1},
se sigue que el resultado sigue siendo v\'alido
para $G=\Gal(L/K)$ un grupo abelinao y no \'unicamente c\'iclico
(ver Teoremas \ref{T12*.2.2.B} y \ref{T12*.2.2.C} y Corolario
\ref{C12*.2.2.D}), cambiando $I_L^{(1-\sigma)}$ por $I_G I_L$ 
donde $I_G=\langle \tau-1\mid \tau\in G\rangle$.
\end{observacion}

En \cite{Has67-1}, Hasse prueba que si la norma
satisface que $\N_{L/K}\colon I_L\longrightarrow
I_K^{\f{}}$, en donde $L/K$ es una extensi\'on c{\'\i}clica de grado primo
$l$, y $I_K^{\f{}}$ es un grupo de clases de ideales de rayos m\'odulo el
conductor $\f{}$, se tiene la sucesi\'on exacta
\[
1\longrightarrow \g I\longrightarrow I_L
\stackrel{\N}{\longrightarrow} \frac{H^{\f{}}_{L/K}}{
P_K^{\f{}}}\longrightarrow 1
\qquad (P_K^{\f{}}\subseteq H_{L/K}^{\f{}}\subseteq I_K^{\f{}}),
\]
donde $\g I$ es el {\em g\'enero principal\index{g\'enero principal}}:
\[
\xymatrix{
K_H\ar@{-}@/_2pc/[dd]_{I_K}\ar@{-}@/^1pc/[d]^{\g I\text{\ g\'enero principal}}
\ar@{-}[d]\\
\g K\ar@{-}[d]\\ K^{\ast}}
\]
M\'as precisamente, $I_K\cong\Gal(K_H/K)$ e $\g I\cong \Gal(K_H/\g K)$.
En particular $\g I=I_L^{(1-\sigma)}$.

Hasse tambi\'en prob\'o el {\em Teorema del g\'enero principal\index{teorema
del g\'enero principal}} general:

\begin{proposicion}[Hasse, 1927]\label{CClaseT4.10.2} Sea $L/K$ una extensi\'on
c{\'\i}clica de grado primo de campos num\'ericos 
con generador $\sigma$ y sea $\mK$ un m\'odulus
en $K$. Entonces existe un m\'odulus ${\eu n}$ de $L$ tal que
\l
\item ${\eu n}|\mK\o_L$.
\item ${\eu n}^{\sigma}={\eu n}$.
\item Para $\beta\in L$ primo relativo a ${\eu n}$, se tiene
\[
\N_{L/K}\beta\equiv 1\bmod \mK\iff \beta=\alpha^{1-\sigma}\bmod {\eu n}
\]
para alg\'un $\alpha \in L$. $\fin$
\end{list}
\end{proposicion}

Con este resultado, Hasse define el {\em g\'enero principal\index{g\'enero
principal}} $H_1\bmod \mK$ en $L$ como el grupo de clase de rayos 
m\'odulo ${\eu n}$ cuyas normas relativas caen en el grupo de clases de
rayos m\'odulo $\mK$ en $K$.

\begin{teorema}[Hasse]\label{CClaseT4.10.3} Sean $L/K$, $\mK$, ${\eu n}$ como
antes. Entonces el g\'enero principal $\bar{H}_1$ coincide con el grupo
de potencias $1-\sigma$ de las clases de rayos m\'odulo ${\eu n}$
en $L$. $\fin$
\end{teorema}

\begin{teorema}[Teorema cl\'asico del g\'enero principal\index{teorema
cl\'asico del g\'enero principal}\index{g\'enero principal!teorema 
cl\'asico del $\sim$}]\label{CClaseT4.10.4}
Sea $L/K$ un extensi\'on c{\'\i}clica de campos num\'ericos y sea $\sigma$
un generador de $\Gal(L/K)$. Entonces $[{\eu a}]\in \g I\iff \N_{L/K}
{\eu a}=(\alpha)$ con $\alpha\in\*K$ un residuo n\'ormico en todos
los primos ramificados de $L/K$, esto es,
\[
\artinp{\alpha, L/K}{\pK}=1 \text{\ para todo $\pK$ ramificado en $L/K$},
\]
equivalentemente, $\artinp{\N_{L/K}{\eu a},L/K}{\pK}=1\ \forall\ \pK|
\f{}(L/K)$. $\fin$
\end{teorema}

%% file: Notaciones.tex
\preface[Notaciones y convenciones]\label{Notaciones y convenciones}
\addcontentsline{toc}{chapter}{Notaciones}

\begin{dinglist}{42}

\item $A=\{x\in K\mid v_{\pK}(x)\geq 0, \pK\neq \p\}$

\item ${\ma A}_K$ es el anillo de ad\`eles o reparticiones de $K$, \pageref{CClaseadeles}.

\item ${\eu a}_{\vec{\alpha}} =\prod_{
\pK\in {\ma P}_K}\pK^{v_\pK(\alpha_\pK)}\in D_K$
para $\vec{\alpha}\in J_K$, \pageref{CClaseidealdeidele}

\item $\artin {L/K}{\pK}=(L/K,\pK)=(\_,L/K,\pK)=(\_, L/K)=\psi_{L/K}(\pK)$ 
denota el s{\'\i}mbolo de Artin, \pageref{CClaseartin}

\item $\Br E$ denota al grupo de Brauer del campo $E$, \pageref{CClasegrupoBrauer}

\item $C\colon A\to \aditivo k$ m\'odulo de Carlitz, \pageref{DrinfeldmoduloCarlitz}

\item $\Ci={\ma C}_p$ completaci\'on de la cerradura algebraica de $\Ki$

\item $\ca K$ denota a la caracter{\'\i}stica de $K$

\item $\ca (\rho)$ denota la caracter\'istica de un m\'odulo de Drinfeld

\item $c_i(\rho,a)$ denotan a los coeficientes de $\rho_a$

\item $Cl_A=C_A$ grupo de clases de ideales de un dominio Dedekind

\item $Cl_S=C_S$ grupo de $S$--clases ideales

\item $\coc D$ el campo de cocientes de un dominio entero $D$

\item $C_K=J_K/\* K$ grupo de clases de id\`eles del campo global $K$,
\pageref{CClaseCK}

\item $C_{K,0}= J_{K,0}/\* K$ grupo de clases de id\`eles de grado $0$,
\pageref{CClaseCK0}

\item $C_K^{\eu m}=J_K^{\eu m}\* K/\* K
\subseteq C_K$ para un m\'odulus ${\eu m}$, \pageref{CClaseCKm1}

\item $C_{K,S}^{\eu m}=\* K J_{K,S}^{\eu m}/\* K$
para un $S$--m\'odulus ${\eu m}$, \pageref{CClaseCKSm}

\item $Cl_K^{\eu m}(\o_S)$ es el cociente de los ideales fraccionarios (divisores)
primos relativos a ${\eu m}$ m\'odulo los ideales principales $z\o_S$ con
$z\in \* K\cap J_{K, {\ma P}_K\setminus \sop({\eu m})}^{\eu m}$
para un un m\'odulus arbitrario ${\eu m}$, $\emptyset\neq S\subseteq
{\ma P}_K$, \pageref{CClaseClKSm}

\item $Cl^+_L=Cl^{\rm{ext}}_L$ grupo de clases extendidos de un campo
global, \pageref{idealesextendidos}

\item $D_A$ ideales fraccionarios de un dominio Dedekind $A$

\item $D_K$ grupo de ideales fraccionarios (si $K$ es un 
campo num\'erico) o
grupo de divisores (si $K$ es un campo global de funciones), \pageref{CClaseDK}

\item $D_{K,0}$ grupo de divisores de grado $0$, $K$ un campo global
de funciones, \pageref{CClaseDK0}

\item $D_K^{\eu m}$ es el grupo de ideales fraccionales
(divisores) primos relativos a la parte finita de un
m\'odulus ${\eu m}$, \pageref{CClaseDKm}

\item $D_K^S$ es el grupo
de divisores con soporte en el complemento de $S$, donde
 $S$ es un conjunto de divisores primos, \pageref{CClaseDKS}
 
\item $D_{K,0}^S=D_K^S\cap D_{K,0}$, \pageref{CClaseDKS0}

\item $D_K(S)$ es el grupo
de divisores con soporte en $S$, donde
 $S$ es un conjunto de divisores primos
 
\item $D_{K,0}(S)$ grupo de divisores de grado $0$ con soporte en $S$

\item $D_{\eu m,A}$, \pageref{Drinfeldidealesm}

\item $\dif_K$ el espacio de diferenciales del campo de funciones $K$

\item $\Drin_A(F)$ es el conjunto de $A$--m\'odulos de Drinfeld sobre
$F$, \pageref{DrinfeldDrinf}

\item $e_{L/K}=e_{L/K}(\pL|\pK)=e(L|K)=e(\pL|\pK)$ denota al
\'indice de ramificaci\'on de $\pL|\pK$ en la extensi\'on $L/K$

\item $\bar{E}$ cerradura algebraica de un campo $E$

\item $e_X(x):=\exponencial x{\gamma}X$ funci\'on exponencial,
\pageref{Drinfeldfuncionexponencial}

\item $\End_{\F} B$ grupo de endomorfismos  de la $\F$--\'algebra $B$,
\pageref{Drinfeldendomorfismos}

\item $e_{L/K}=e_{L/K}(\pL|\pK)=e(L|K)=e(\pL|\pK)$ denota
el \'indice de ramificaci\'on de $\pL|\pK$ en la extensi\'on $L/K$

\item $\cond{\pK}=\cond{L/K}=\pK^{c_\pK}=\pK_K^{c_\pK}$, es el conductor
local, \pageref{CClaseconductorlocal}

\item ${\mc F}_{\pi}=\{f\in\o_K[[T]]\mid f\equiv \pi T\bmod\deg 2, f\equiv T^q\bmod
\pi\}$, $\o_K$ el anillo de enteros de un anillo local, \pageref{CClaseFpi}

\item $\frobenius {L/K}{\pL}$ denota al automorfismo de
Frobenius, \pageref{CClasefrobenius}

\item $\aditivo F$ polinomios torcidos

\item $\Gal(L/K)=G_{L/K}$ denota el grupo de Galois de la extensi\'on $L/K$

\item $\deg_\pK \alpha_\pK=[K(\pK):\F]
v_\pK(\alpha_\pK)=[\o_\pK/\pK\colon \F]=v_\pK(\alpha_\pK)=f_\pK v_\pK
(\alpha_\pK)=\deg \pK v_\pK(\alpha_\pK)$
para $\vec{\alpha}\in J_K$, \pageref{CClasegradolocal}

\item $\deg \vec(\alpha)=\sum_{\pK\in {\ma P}_K}
\deg_\pK \alpha_\pK=\deg {\eu A}_{\vec{\alpha}}$
para $\vec(\alpha)\in J_K$, \pageref{CClasegrado global}

\item $H^2(\ |K)=H^2(K^{\text{sep}}|K)=H^2(\Gal(K^{\text{sep}}/K),\K)$,
$K$ un campo local, \pageref{CClaseH2Brauer}

\item $H^q(L/K)=H^q(\Gal(L/K), \K)$, $K$ un campo local, \pageref{CClaseHql/K}

\item $H^{\eu m}=(\N_{L/K} D_L^{\eu m}) P_K^{\eu m}$
para un m\'odulus $\mK$, \pageref{CClaseH^m}

\item $\mathcal H$ denota al conjunto de m\'odulos de Hayes, \pageref{DrinfeldHayes}

\item $h_A=|\pic(A)|$ es el n\'umero de clase de ideales de $A$

\item $h_K$ es el n\'umero de clase de un campo global $K$,
\pageref{numero de clase de un campo}

\item $[i]=T^{q^i}-T$, $[0]=0$, \pageref{DrinfeldD2.3.1}

\item $\Isog(\rho,\rho^{\prime})$ grupo de isogen\'ias entre
dos m\'odulos de Drinfeld, \pageref{Drinfeldisogenias}

\item $I_K=D_K/P_K$ grupo de clases de ideales (divisores) del
campo global $K$, \pageref{CClaseIK}

\item $I_{K,0}$ grupo de clases de
divisores de grado $0$ del campo global $K$ de funciones, \pageref{CClaseIK0}

\item $I_K^{\eu m}=D_{K}^{\eu m}/P_K^{\eu m}$

\item $J$ denota la conjugaci\'on compleja

\item $J_K$ es el grupo de id\`eles del campo global $K$, \pageref{CClaseid'eles},
\pageref{CClaseideles}

\item $J_{K,0}=\{\vec{\alpha}\in J_K\mid \deg \vec{\alpha}=\sum_{\pK}
\deg_{\pK}\alpha_{\pK}=0\} = \{\vec{\alpha}\in J_K\mid \prod_{\pK}|\alpha_{\pK}
|_{\pK} =1\} \supseteq \* K$ grupo de id\`eles 
de grado $0$, \pageref{CClaseidelesgrado0}

\item $J_K^{\eu m}=\{\alpha\in J_K\mid
\alpha\equiv 1\bmod {\eu m}\}=\prod_{\pK\in {\ma P}_K} U_{\pK}^{(n_\pK)}$
donde ${\eu m}$ es un m\'odulus

\item $J_K^S=\{\alpha\in J_K
\mid \alpha_\pK=1 \text{\ para toda\ } \pK\in S\}$
donde $S$ es un conjunto de lugares de $K$

\item $J_{K,S}=\prod_{\pK\notin S}
U_{\pK}\times \prod_{\pK\in S} \* {K_\pK}$ ($J_{K,S}\neq J_K^S$),
donde $S$ es un conjunto finito de lugares

\item $J_{K,S}^{\eu m}=\Big(\prod_{\pK \in S}
\* {K_\pK}\times \prod_{\pK\notin S} U_{\pK}^{(n_{\pK})}\Big)\cap J_K$
para ${\eu m}$ un $S$--m\'odulus

\item $\Lambda_{f,n}=\{\lambda\in\Omega\mid v(\lambda)>0, f^n(\lambda)=0\}$, \pageref{CClasetorsion}

\item $K$, $\K$ un campo global de funciones sobre $\F$

\item $K_H$ campo de clase de Hilbert de un campo global $K$

\item $K_{H^+}$ campo de clase de Hilbert de un campo global $K$

\item $K^+$ conjunto de elementos totalmente positivos en un
campo global $K$, \pageref{totalmentepositivos}

\item $K({\eu p})=\tilde{K}={\mathcal O}_{\pK}/\pK=\F$ el campo residual de $K$
un campo local, \pageref{CClasecamporesidual}

\item $K_\pK$ es la completaci\'on de un campo global $K$ en $\pK\in{\ma P}_K$

\item $\Ki=K_{\p}$

\item $\Kii=\F(T)_{\p}$

\item $K^{\ast}=(\pi)\times \F^{\ast}\times U_{\pK}^{(1)}$, 
$K$ campo local, \pageref{CClasecampolocal}

\item $\bar{K}$ una cerradura algebraica de un campo $K$

\item $K^{\rm{ab}}$ la m\'axima extensi\'on abeliana del campo $K$, 
\pageref{CClaseab}

\item $K^{\rm{nr}}$ es la m\'axima extensi\'on no ramificada de $K$ contenida en
$K^{\rm{ab}}$

\item $K^{\eu m}$ campo de clase de $C_K^{\eu m}$, \pageref{CClaseKm}

\item $K_S^{\eu m}$ es el campo de clase de $C_{K,S}^{\eu m}$

\item $\* K\hookrightarrow J_K$ de manera diagonal: $x\mapsto (\ldots, x,x,x.\ldots)$,
$K$ campo global

\item $K=\F(T)$, \pageref{Drinfeldcampofuncionesracionales}

\item $K_{\eu m}:=K(\rho[\eu m])$, \pageref{DrinfeldKm}

\item $\bar{K}$ una cerradura algebraica de $K$

\item $[L:K]_s$ denota el grado de separabilidad de la extensi\'on $L/K$

\item $L^+=\{x\in L= \mid x\text{\ es totalmente positivo}\}$, 
\pageref{DClaseS4.9-1.12}

\item $\g L$ campo de g\'eneros de $L$

\item $\gex L$ campo de g\'eneros extendido de $L$

\item ${\g L}_{\K}$ campo de g\'eneros de $L$ con respecto a $\K$

\item ${\gex L}_{\K}$ campo de g\'eneros extendido de $L$
con respecto a $\K$

\item $L_i=\prod_{j=1}^i[j]$, \pageref{DrinfeldLi}

\item $L_{\pK}$ denota uno de los campos $\{L_{\pL}\}_{\pL|
\pK}$, para $\pK\in {\ma P}_K$, cualquiera pero s\'olo uno. Es decir, $L_{\pK}/K_{\pK}$ denota $L_{\pL}/
K_{\pK}$ con $L_{\pL}$ seleccionado arbitrariamente

\item Un m\'odulus es ${\eu m}=\prod_{\pK\in {\ma P}_K}\pK^{n_\pK}$
donde $n_\pK\geq 0$ para todo
$\pK\in {\ma P}_K$, $n_\pK=0$ para casi todo $\pK$, $n_\pK=0$ si $\pK$ es complejo
y $\pK\in\{0, 1\}$ si $\pK$ es real

\item $M[a]:=\{m\in M\mid am=0\}$, \pageref{Drinfeldtorsion}

\item $M[I]:=\{m\in M\mid am=0 \text{\ para toda\ } a\in I\}=\cap_{a\in I} M[a]$,
\pageref{Drinfeldtorsion}

\item $\mu_{L/K}\in H^2(L/K)$ la clase fundamental, $\inv_{L/K}(\mu_{L/K})=
\frac{1}{[L:K]}+{\ma Z}\in \big(\frac{1}{[L:K]}{\ma Z}\big)/{\ma Z}$, 
\pageref{CClaseinvariante}

\item $\N_{L/K}=\N$ denota la norma de $L$ a $K$, \pageref{CClasenorma}

\item $\N(\pK)$ denota la norma absoluta, \pageref{CClasenormaabsoluta}

\item $\o_S=\cap_{\pK\notin S}
\o_\pK=\{x\in\* K\mid v_\pK(x)\geq 0 \text{\ para toda \ } \pK\notin S\}$
para $\emptyset\neq S\subseteq {\ma P}_K$

\item ${\mathcal O}_{\eu p}={\mathcal O}_K=\{x\in K\mid v_{\eu p}(x)\geq 0\}=
\{x\in K\mid |x|_{\eu p}\leq 1\}=\bar{B}(0,1)$, $K$ un campo local, \pageref{CClaseanillovaluacion}

\item $P_A$ ideales principales de un dominio Dedekind $A$

\item $P_K$ grupo de ideales fraccionarios 
(divisores) principales de un campo global $K$

\item $P_K^{\eu m}=\{(a)\in P_K\mid a\equiv 1\bmod {\eu m}\}$
para un m\'odulus $\mK$

\item $P^+_{K,\mK}=\{\big(\alpha/\beta\big)\mid \alpha, \beta\in \o_K \text{\ con\ }
\mcd(\alpha\beta,\mK)=1, \alpha\equiv\beta \bmod \mK\text{\ y\ }\alpha/\beta$
es totalmente positivo$\}$, \pageref{CClasePKm+}

\item $P_{\eu m,A}^+$, \pageref{Drinfeldprincipalesm}

\item $P_{K, S}^{\eu m}=\* K\cap
J_{K, {\ma P}_K\setminus \sop({\eu m})}^{\eu m}$

\item $\pic(A)=\frac{D_A}{P_A}$ el grupo de Picard de un dominio Dedekind

\item $\pic_{\eu m}^+(A)=\frac{D_{\eu m,A}}{P_{\eu m,A}^+}$, \pageref{Drinfeldpicm}

\item $P(\Gamma^{\prime}/\Gamma, u)=\exponencial u{\lambda}
{e_{\Gamma}(\Gamma^{\prime})}$, \pageref{DrinfeldpolinomioDrinfeld}

\item ${\ma P}_K=\{\pK\mid \pK \text{\ es un lugar de\ } K\}$
es el conjunto de lugares de $K$, $K$ un campo global

\item ${\mathcal P}(f):=-\limsup\limits_{i\to\infty} \dps \frac{v(a_i)}{i}$ radio de
convergencia de la serie de potencias $\sum_{i=0}^{\infty} a_i u^i$,
\pageref{Drinfeldradioconvergencia}

\item ${\mathcal P}(F)=\{f(x)\in F[x]\mid \text{$f$ es aditivo}\}$, \pageref{Drinfeldaditivos}

\item ${\eu p}=\{x\in K\mid v_{\eu p}(x)>0\}=\{x\in K\mid |x|_{\eu p}<1\}
= B(0,1)$, $K$ es un campo local, \pageref{CClaseidealmaximo}

\item $\p$ un lugar fijo de un campo de funciones global $K$ de grado $\di$

\item $\pi_{\infty}=1/T$ elemento uniformizador del primo infinito $\p$ en 
$K=\F(T)$, \pageref{piinfty} 

\item $\pi_\pK=\pi$ un elemento primo de $K$, $v_\pK=1$, $\pi\in 
\pK\setminus \pK^2$, \pageref{Drinfeldelementoprimo}, 
\pageref{CClaseelementoprimo}

\item $\psi_{L/K,{\eu m}}\colon I_{K}^{\eu m}\to \Gal(L/K)$ es el mapeo de Artin,
\pageref{CClaseartinconm}

\item $\psi_{L/K}\colon K^{\ast}/N_{L/K}L^{\ast}\xrightarrow[\psi_{L/K}]{\iso}
\Gal(L/K)$, $\psi_{L/K}(a)=(a,L/K)$ mapeo local de Artin o s{\'\i}mbolo de la norma
residual o s{\'\i}mbolo residual n\'ormico, $L/K$ una extensi\'on abeliana finita
de campos locales, \pageref{CClasesimboloresidual}

\item $\Psi({\eu m})=\prod_{\pK\in {\ma P}_K}
[U_\pK\colon U_\pK^{(n_\pK)}]=\prod_{\pK\in {\ma P}_K}(q^{\deg \pK}-1)(
q^{(n_{\pK}-1)\deg \pK})$, para un m\'odulus ${\eu m}$, 
\pageref{CClasePsim}

\item $\Psi({\eu A})=\big|A/{\eu A}\big|$ donde $\eu A$ es un ideal no cero de $A$

\item ${\ma Q}_p$ denota al campo de los n\'umeros $p$--\'adicos

\item $\Red_A(\Ci)$ denota al conjunto de las $A$--redes de $\Ci$,
\pageref{Drinfeldredes}

\item $\Red_{A,r}(\Ci)$ es el conjunto de redes de 
rango $r$ en $\Ci$, \pageref{Drinfeldredes1}

\item ${\mathcal R}_{A,r}$ denota al conjunto de clases de 
isomorfismos en $\Red_{A,r}(\Ci)$, \pageref{Drinfeldredes2}

\item $R_T=\F[T]$, \pageref{Drinfeldanillopolinomios}
\item $\rho\colon A\lra \aditivo F$ un $A$--m\'odulo de Drinfeld sobre $F$

\item $\rho_K\colon F\to \Gal(\ab K/K)$ es el
mapeo de reciprocidad o el mapeo local de Artin y $\psi_{L/K}\colon F\to
\Gal(L/K)$ es $\psi_{L/K}=\rest\circ \rho_K$, $\rest\colon \Gal(\ab K/K)\to
\Gal(L/K)$, $\sigma\mapsto \sigma|_L$
para un campo $K$ local o global, $F= J_K$ si $K$ es global y
$F=\* K$, si $K$ es local

\item $\Sg_L:=\*L/L^+$, grupo de signos de un campo global 
$L$, \pageref{gruposignoscampoglobal}

\item El soporte de ${\eu m}$ es $\{\pK\in
{\ma P}_K\mid n_{\pK}>0\}$, donde ${\eu m}$ es un m\'odulus.
Si $S$ es un conjunto de divisores primos, el
soporte de $S$ es $S$ mismo

\item $\Spl (L/K)=\{\pK\in {\ma P}_K\mid \pK \text{\ se descompone totalmente
en \ } L\}$, \pageref{CClaseSpl}

\item $\sgn\colon \*{\Ki} \lra \Fi$ denota una funci\'on signo, 
\pageref{DrinfeldD3.1.19}

\item $T$ denota la m\'axima extensi\'on no ramificada de un campo local
$K$, $T=K^{\rm{nr}}$, \pageref{CClasemaximaextensionnoramificadalocal}

\item TCCG = teorema principal de la teor{\'\i}a de campos de clase globales,
\pageref{CClaseTCCG}

\item TCCL = Teorema de la Teor{\'\i}a de Campos de Clase Locales,
\pageref{CClaseTCCL}

\item $\Tr_{L/K}=\Tr$ denota la traza de $L$ a $K$, \pageref{CClasetraza}

\item $\tau$ el automorfismo de Frobenius sobre $\F$: $\tau u=u^q$

\item Un $S$--m\'odulus ${\eu m}$ es tal que
$S\cap \sop({\eu m})=\emptyset$, donde $S\subseteq {\ma P}_K$

\item $U_\pK=U_K=U_\pK^{(0)}=U_K^{(0)}
=\{x\in K\mid v_\pK(x)=0\}={\mathcal O}_{\pK}^{\ast}$ el 
grupo de unidades de ${\mathcal O}_{\pK}$, $K$ campo local, \pageref{CClaseunidadeslocales}

\item $U_\pK^{(n)}=\unidades n=1+\pK^n=\{x\in K\mid v_\pK(x-1)\geq n\}=\{
x\in K\mid |x-1|_\pK\leq q^{-n}\}$, $K$ campo local,
\pageref{CClasenunidadeslocales}

\item Para $n\in{\ma N}$, $W_n$ denota el grupo de las $n$--ra\'ices de $1$

\item $\vi=v_{\p}$

\item $v_K$ denota a la valuaci\'on de un campo local $K$

\item $|x|_\pK=q^{-v_\pK(x)}$, $x_\pK$ en un campo local con
campo residual de $\F$ elementos \pageref{CClasevalorabsoluto}

\item $\lceil x_\pK\rceil=(\ldots,1,1,\underbracket[0pt]{x_{\pK}}_{
\substack{\uparrow\\ \pK}},1,1,\ldots)\in J_K$, 
para $x_\pK\in K_\pK$, $K$ campo global

\item ${\ma Z}_p$ es el anillo de los enteros $p$--\'adicos

\item $\qed =$ terminaci\'on de una demostraci\'on

\item $\emptyset$ denota al conjunto vac{\'\i}o

\end{dinglist}